\documentclass{amsbook}%
\usepackage{amsfonts}
\usepackage{amscd}
\usepackage{amsmath}
\usepackage{amssymb}
\usepackage{amstext}
\usepackage{answers}
\usepackage{appendix}
\usepackage{accents}
\usepackage{graphicx}
\usepackage{titletoc}%
\providecommand{\U}[1]{\protect\rule{.1in}{.1in}}
\theoremstyle{plain}

\newtheorem{claim}{Claim}

\newtheorem{corollary}{Corollary}

\newtheorem{definition}{Definition}

\newtheorem{lemma}{Lemma}

\newtheorem{proposition}{Proposition}
\newtheorem{remark}{Remark}

\newtheorem{summary}{Summary}
\newtheorem{theorem}{Theorem}
\numberwithin{equation}{chapter}
\ifx\pdfoutput\relax\let\pdfoutput=\undefined\fi
\newcount\msipdfoutput
\ifx\pdfoutput\undefined\else
\ifcase\pdfoutput\else
\ifx\paperwidth\undefined\else
\ifdim\paperheight=0pt\relax\else\pdfpageheight\paperheight\fi
\ifdim\paperwidth=0pt\relax\else\pdfpagewidth\paperwidth\fi
\fi\fi\fi
\begin{document}
\frontmatter
\title[GENERALIZED FEYNMAN-KAC FORMULA]{GENERALIZED FEYNMAN-KAC FORMULA AND ASSOCIATED HEAT KERNEL IN VECTOR BUNDLES}
\author[M. N. Ndumu]{Martin N. Ndumu}
\address{Department of Mathematics and Computer Science, University of Maryland Eastern
Shore, Princess Anne, MD 21853, United States of America}
\curraddr{Department of Mathematics and Computer Science, University of Bamenda,
Bambili, North West Region, Cameroon}
\email{martin.ndumu@gmail.com}
\urladdr{http://www.author.institute}
\thanks{I wish to take this opportunity to thank Professor David Elworthy of Warwick
University, England, for having initiated me to the ideas of Stochastic
Differential Geometry.}
\dedicatory{This work is dedicated to the memory of the late Professor James Eells,
formerly of Warwick University, who suggested the topic to me several years ago.}
\begin{abstract}
Let M be a smooth closed (compact without boundary) Riemannian manifold of
dimension n and P a q-dimensional smooth submanifold of M. M$_{0}$ will denote
the tubular neighborhood of P in M. Let E be a smooth vector bundle over M.
Let L = $\frac{1}{2}\Delta+$ X +V \ be a differential operator on the set of
smooth sections $\Gamma(E)$ of the vector bundle E, where $\Delta$ is the
generalizeed Laplacian (or Laplace-Type operator) on the vector bundle E, X a
smooth vector field on M and V a smooth potential term on M.

Let (x$^{\text{t}}$(s)) 0$\leq$ s$\leq$t be the semi-classical Brownian
Riemannian bridge process from x$\in$M$_{0}$ to P in time t and let
$\tau_{\text{0,s}}^{\text{t}}=$ u$_{\text{0}}^{\text{t}}$(u$_{\text{s}%
}^{\text{t}}$)$^{-1}:$E$_{\text{x}^{\text{t}}\text{(s)}}\longrightarrow
$E$_{\text{x}}$ be the parallel translation along the reversed semi-classical
Brownian bridge process with first exit time $\zeta$ from the tubular
neighborhood M$_{0}$ of P. Let W be the Weitzenb\H{o}ck term and
let,\ e$_{\text{s}}^{\text{t}}\ =$ $\tau_{\text{s,0}}^{\text{t}}+\frac{1}%
{2}\int_{0}^{\text{s}}\tau_{\text{s,r}}^{\text{t}}$e$_{\text{r}}^{\text{t}}%
$W$_{\text{x}_{\text{r}}^{\text{t}}}$dr .

We will show here that for 0$\leq$ s$\leq$t$\Lambda\zeta,$

(Q(t,t-s)$\phi$)(x) $=$ \textbf{E}$_{\text{x}}\left[  \tau_{\text{0,s}%
}^{\text{t}}\text{e}_{\text{s}}^{\text{t}}\phi\text{(x}^{\text{t}}%
\text{(s))}\exp\left\{  \int_{0}^{\text{s}}\frac{\text{L}^{0}\Psi}{\Psi
}\text{(x}^{\text{t}}\text{(u))du}\right\}  \right]  $

is a \textbf{Generalized Feynman-Kac} formula from which we shall deduce the
usual \textbf{Feynman-Kac formula} as well as a stochastic representation of
the \textbf{Generalized Elworthy-Truman heat kernel formula}, and ultimately
the \textbf{heat kernel formula. }

The generalized Feynman-Kac formula shall be expanded and from this expansion
we shall deduce both the generalized \textbf{heat trace} and \textbf{heat
content} expansions.

We then deduce the expansion of the generalized heat kernel and then compute
at the centre of Fermi coordinates (reduced to normal coordinates)\ the first
few coefficients of the expansion.

\end{abstract}
\maketitle

\chapter*{Preface}

This work is in \textbf{Stochastic Differential Geometry. }A description
borrowed from \textbf{Elworthy} $\left[  15\right]  $ states that
\textbf{Differential Geometry} deals with the \textbf{triangle} of
inter-relationships between the \textbf{curvature}, the \textbf{spectrum} and
the \textbf{topology} of a manifold:\ \ \ \ \ \ \ \ \ \ \ \ \ \ \ \ \ \ \ \ \ \ \ \ \ \ \ \ \ \ \ \ 

\begin{center}
\textbf{curvature}\qquad\textbf{spectrum}

\qquad\qquad$\ \ \ \ \ \qquad$

$\qquad\qquad\qquad\ \ \ \ \ \ \ \ \ \qquad\ $\textbf{topology}$\qquad
\ \ \ \ $\qquad\qquad\qquad\qquad\qquad\qquad\qquad
\end{center}

\textbf{Stochastic Differential Geometry} completes the \textbf{square }by
adding \textbf{Brownian motion} into these inter-relationships:\qquad
\qquad\qquad\qquad\qquad\ \ \ \ \ \ 

\begin{center}
\textbf{curvature \ \ \ \ \ \ \ spectrum}

$\qquad$

\qquad\ \ \ \ \ \textbf{topology \ }\ \ \ \ \textbf{\ \-\ \ Brownian motion}
\end{center}

Brownian motion comes in many ways, including its density (with respect to the
Riemannian volume measure on the manifold) called \textbf{the heat kernel
}p$_{\text{t}}$(x$,$y)\textbf{.} The expansion (using normal coordinates), in
powers of t, of the heat kernel generates the spectral invariants of the manifold.

Let M be a \textbf{complete connected} n-dimensional Riemannian manifold and P
a $q-$ dimensional smooth submanifold of M, such that $0$ $\leq q\leq n.$

This paper is a follow-up of \textbf{Ndumu} $\left[  42\right]  ,$
\textbf{Ndumu} $\left[  43\right]  $ and \textbf{Ndumu} $\left[  44\right]  .$
It is a dirct follow-up of \textbf{Ndumu} $\left[  44\right]  $ in which we
defined the generalized \textbf{scalar} heat kernel p$_{\text{t}}$(x$,$P)
relative to the differential operator L = $\frac{1}{2}\Delta^{0}+$ X + V
$=\frac{1}{2}\Delta^{0}+$ $\nabla_{\text{X}}$ $+$ V and the submanifold P,
where $\Delta^{0}$ is the Laplace-Beltrami operator on functions defined on M,
X is a smooth vector field on M and V is a smooth potential term on M. Then
using Fermi coordinates we derived an integral formula for p$_{\text{t}}$%
(x$,$P) given by p$_{\text{t}}$(x$,$P) = $\int_{\text{P}}$p$_{t}^{\text{M}%
_{0}}$(x,y)f(y)$\upsilon_{\text{P}}$(dy)$,$ where f : M$\rightarrow$ R is a
smooth function of compact support in M. p$_{t}^{\text{M}_{0}}(-,-)$ is the
usual Dirichlet heat kernel (relative to the operator L$^{0}$ = $\frac{1}%
{2}\Delta^{0}+$ X + V) of the tubular neighborhood M$_{0}$ of P and
$\upsilon_{\text{P}}$ is the Riemannian volume measure on P. We then derived
an exact and an asymptotic expansion for p$_{\text{t}}$(x,P) and then computed
the leading coefficients of the expansion.

In this paper we generalize \textbf{Ndumu} $\left[  44\right]  $ to the case
of heat kernels of an elliptic operator of the form: L = $\frac{1}{2}\Delta+$
X + V $=$ $\frac{1}{2}\Delta+$ $\nabla_{\text{X}}$ $+$ V on a vector bundle E
over a compact Riemannian manifold M where $\Delta$ is a Laplace -Type
operator, $\nabla$ is a metric connection on the vector bundle E, X is a
smooth vector field on M and V is a smooth potential term on M.

Here we describe a \textbf{generalization of heat flow} on a vector bundle E
relative to the differential operator defined above and derive a
\textbf{Generalized Feynman-Kac formula} from which we deduce the
\textbf{usual Feynman-Kac formula }and a \textbf{Generalized Elworthy-Truman
Heat Kernel Formula} in vector bundles.

We will next obtain the \textbf{Expansion of the Generalized Feynman-Kac
Formula} and then show that it is a double generalization of the \textbf{heat
kernel} and the \textbf{heat content} expansions.

For the existence of heat kernels of such elliptic differential operators see
$\left(  8.2\right)  $ of Chapter 8 of \textbf{Duistermaat }$\left[
13\right]  $ or \textbf{Theorem }$\left(  2.26\right)  $ of \textbf{Berline,
Getzler and Vergne }$\left[  7\right]  .$ See also the discussion in \S 3 of
\textbf{Baudoin} $\left[  5\right]  $ or \textbf{Theorem }$\left(
1.3.5\right)  $ of \textbf{Gilkey} $\left[  21\right]  $ and the discussion in
\S 6 of Chapter III of \textbf{Lawson and Michelsohn} $\left[  35\right]  .$

We will show that the integral formula above generalizes to:

\begin{center}
k$_{\text{t}}$(x,P) = $\int_{\text{P}}$k$_{t}$(x,y)$\phi$(y)$\upsilon
_{\text{P}}$(dy)
\end{center}

We will write L = $\frac{1}{2}\Delta+$ X $+$ V to mean L = $\frac{1}{2}%
\Delta+$ $\nabla_{\text{X}}$ $+$ V$.$

To the best of my knowledge, no author orther than Baudoin, has considered
heat kernel expansions of a vector bundle heat kernel in the presence of a
vector field and/or a potential term. Baudoin $\left[  5\right]  $ included a
vector field (in a more general way). In fact he considered an operator of the type:

\begin{center}
\qquad L = $\frac{1}{2}\underset{i=1}{\overset{n}{\sum}\nabla_{i}^{2}}+$
$\nabla_{0},$
\end{center}

where $\nabla_{i}=\nabla_{X_{i}}+%
\mathcal{F}%
_{i}$ for $i=0,1,...,n$ and where $%
\mathcal{F}%
_{i}$ are smooth potentials, i.e. Weitzenb\H{o}ckians in our sense here and
$X_{i}$ are vector fields on M.

He showed that, under an ellipticiy condition, the associated partial
differential equation:

\begin{center}
\qquad$\frac{\partial\phi_{\text{t}}}{\partial\text{t}}=$ L$\phi_{\text{t}}$
\qquad\qquad\qquad\ \ \ (evolution equation)

$\phi_{0}=$ $\phi$\qquad\qquad\ \ \ \ \ \ (initial condition)
\end{center}

for $\phi_{0},\phi\in\Gamma(E),$ where $\Gamma(E)$ is the space of smooth
sections of the vector bundle E, has for solution:

\begin{center}
\qquad$\phi_{t}(x)=(e^{tL}\phi)(x)=$ P$_{t}\phi(x)=%
{\textstyle\int\limits_{M}}
$k$_{t}$(x,y)$\phi(y)\upsilon_{\text{M}}(dy)$
\end{center}

where k$_{t}(-,-)\in$Hom(E,E). More precisely, the heat kernel has the
property: k$_{t}$(x,y)$\in$Hom(E$_{y}$,E$_{x}$).

\mainmatter

\part{ GENERALIZATION OF HEAT FLOW}

\chapter{Fermi Coordinates}

One of the main geometric tools here is \textbf{Fermi coordinates}. Fermi
coordinates were discovered by Enrico Fermi (see \textbf{Fermi} $\left[
18\right]  ,$ $\left[  19\right]  $). Following \textbf{Gray} $\left[
23\right]  ,\left[  25\right]  $, we define them below in a more modern and
general setting:

Let P be a submanifold of M and let T$_{\text{y}}$P be the tangent space of P
at y$\in$P. Then T$_{\text{y}}$P is a subspace of T$_{\text{y}}$M for each
y$\in$P. Let (T$_{\text{y}}$P)$^{\perp}$ be the orthogonal complement of
T$_{\text{y}}$P in T$_{\text{y}}$M: T$_{\text{y}}$M = T$_{\text{y}}$P$\oplus
$(T$_{\text{y}}$P)$^{\perp}$ for y$\in$P.

Set B = $\underset{\text{y}\in\text{P}}{\cup}$(T$_{\text{y}}$P)$^{\perp}.$
Then B\ is a vector bundle over P with fibers (T$_{\text{y}}$P)$^{\perp}$
called the \textbf{normal bundle} over the submanifold P. We note that it is
usual to denote B = $\underset{\text{y}\in\text{P}}{\cup}$(T$_{\text{y}}%
$P)$^{\perp}$ by (TN)$^{\perp}$ so that:

\qquad TM$\mid_{\text{N}}=$ TN$\oplus$(TN)$^{\perp}.$

Let $\Pi:$B$\longrightarrow$ P be the projection from the normal bundle to the
submanifold P.

Let exp$_{\Pi}:$ B$\rightarrow$ M be the exponential map of the normal bundle.
Let (y$_{\text{1}}$,...,y$_{\text{q}}$) be a local coordinate system on P at a
point y$_{0}\in$P. Let E$_{\text{q+1}},...,$E$_{\text{n}}$ be orthonormal
sections of the normal bundle B defined in a neighborhood U$\subset$ P of
y$_{0}$ . Then the Fermi coordinates (x$_{\text{1}}$,...,x$_{\text{q}}%
$,...,x$_{\text{n}}$) of P at y$_{0}$ relative to (y$_{\text{1}}%
$,...,y$_{\text{q}}$) and E$_{\text{q+1}},...,$E$_{\text{n}}$ are defined for
y$\in$U by:\qquad21`\qquad$\qquad\ \ \ \ \ \ \ \ \ \ \ \ \ \ \ \ \ $

$\left(  1.1\right)  $\qquad x$_{\text{a}}$(exp$_{\Pi}$(y,
$\underset{i=q+1}{\overset{n}{\sum}}$t$_{\text{i}}$E$_{\text{i}}$(y))) =
x$_{\text{a}}$(exp$_{\text{y}}$($\underset{i=q+1}{\overset{n}{\sum}}%
$t$_{\text{i}}$E$_{\text{i}}$(y))) = y$_{\text{a}}$(y) for a =1,...,q$\qquad$

$\left(  1.2\right)  $\qquad x$_{j}$(exp$_{\Pi}$(y,
$\underset{i=q+1}{\overset{n}{\sum}}$t$_{\text{i}}$E$_{\text{i}}$(y))) =
x$_{j}$(exp$_{\text{y}}$($\underset{i=q+1}{\overset{n}{\sum}}$t$_{\text{i}}%
$E$_{\text{i}}$(y))) = $t_{j}$ \ for $j=q+1,...,n$ \qquad

The constants t$_{q+1},...,$t$_{n}$ here are chosen small enough so that
(y,$\underset{i=q+1}{\overset{n}{\sum}}$t$_{\text{i}}$E$_{\text{i}}$%
(y))$\in\ $B$_{0},$ \ where B$_{0}$ is defined below.

In particular if t$_{i}=0$ for $i=q+1,...,n,$ then,

$\left(  1.1\right)  ^{\ast}$\qquad x$_{\text{a}}$(y) = y$_{\text{a}}$(y) for
a $=1,...,q\qquad$

$\left(  1.2\right)  ^{\ast}$\qquad x$_{j}$(y) = 0 \ for $j=q+1,...,n$

\qquad\qquad\qquad\qquad\qquad\qquad\qquad\qquad\qquad\qquad\qquad\qquad
\qquad$\qquad\qquad\blacksquare$

Strictly speaking, the coordinates we have defined above are called Cartesian
Fermi coordinates. Polar Fermi coordinates defined in \textbf{Gray, Karp and
Pinsky} $\left[  29\right]  $ will not be used here, and so without any
ambiguity, Cartesian Fermi coordinates here will simply be called Fermi coordinates.

The \textbf{zero section} of B is defined by \textbf{Zero}(B) = $\left\{
\text{(y,0)}\in\text{B: y}\in\text{P }\right\}  $\qquad

The exponential map exp$_{\pi}:$B$\rightarrow$M of the normal bundle $\Pi
:$B$\rightarrow$P maps \textbf{Zero}(B) diffeomorphically onto P and a
neighborhood B$_{0}$ of \textbf{Zero}(B) diffeomorphically onto a
neighbourhood M$_{0}$ of P in M. To be more specific, we will follow
\textbf{Gray} $\left[  23\right]  ,\left[  25\right]  $ for the definition of
B$_{0}$ and M$_{0}.$M$:$ Let,

\qquad\qquad S(N) = $\left\{  (y,\xi)\in B:\text{ }\left\Vert \xi\right\Vert
=1\right\}  $

be the sphere sub(bundle) of B and let c:S(N)$\rightarrow$R$_{+}$ be the
positive function defined by:

$\left(  1.3\right)  \qquad$c$(y,\xi)$ = sup$\left\{  \rho\geq0:d(\exp_{\Pi
}(y,\rho\xi),P)=\text{ }\rho\right\}  $

where d is the distance on M compatible with the Riemannian metric on M.

Then B$_{0}$ is defined by:

\ $\left(  1.4\right)  $ (i)\qquad B$_{0}$ = $\left\{  (y,\rho\xi)\in
B:0\leq\rho<c(y,\xi)\right\}  $

Then

$\left(  1.4\right)  $ (ii)\qquad M$_{0}=$ exp$_{\Pi}($B$_{0})$

is called a \textbf{tubular neighbourhood }of P in M (see Lemma 2.3 of
\textbf{Gray} $\left[  25\right]  $ where our B$_{0}$ here is $\Theta_{P}$
defined there in $\left(  2.1\right)  $).

M$_{0}$ is called a \textbf{tube} if there exist a constant c$_{0}>0$ such
that $c(y,\xi)$ $=c_{0}.$

The tubular neighbourhood M$_{0}$ of P can also be characterized as follows
(see \textbf{Gray} $\left[  23\right]  $):

M$_{0}=\left\{  \text{x}\in\text{M: there exists a unique unit speed geodesic
}\gamma\text{ from x to P that meets P orthogonally}\right\}  :$

There exists $\gamma:\left[  0.1\right]  \longrightarrow$M$_{0}$ such that:
$\gamma(0)$ = x; $\gamma(1)=$ y$\in$P; $\overset{.}{\gamma}(1)\in$(T$_{y}%
$P)$^{\perp}$ and $\left\Vert \overset{.}{\gamma}(s)\right\Vert =1,$

where:

$\gamma(s)=\exp_{\text{y}}((1-s)$v$)$ where $\gamma$(0) $=\exp_{\text{y}}$(v)
$=$ x and $\gamma$(1)

= $\exp_{\text{y}}$(0) $=$ y.$\qquad\qquad\qquad\qquad\qquad\qquad\qquad$

Next define $\Phi_{P}$:M$_{\text{0}}\rightarrow$R$_{+}$ \ by:

$\left(  1.5\right)  \qquad\Phi_{P}($x) = exp$\left\{  \int_{0}^{1}%
<\text{X(}\gamma\text{(s)) , }\dot{\gamma}\text{(s)%
$>$%
ds}\right\}  $

where X is a vector field on M and $\gamma$ is the unique minimal unit speed
geodesic from x$\in$M$_{\text{0}}$ meeting P orthogonally at a point y$\in$P
in time 1. More generally we have:

$\left(  1.6\right)  $\qquad$\Phi_{P}(\gamma(s))$ $=$ $\exp\left\{  \int%
_{s}^{1}<X(\gamma(u)),\text{ }\dot{\gamma}(u)>du\right\}  .$

Let $g_{ij}(x)=$ $<$ $\frac{\partial}{\partial x_{i}},\frac{\partial}{\partial
x_{j}}>_{x},$ $i,j=1,...,n$ be the components of the metric tensor field
defined by the Fermi coodinates $x_{1},...,x_{q},x_{q+1},...,x_{n}$ relative
to P.

We follow \textbf{Definition 1.12} of Roe $\left[  48\right]  $ for a general
definition and \textbf{Berlne, Getzler, Vergne }$\left[  7\right]  ,$ p. 36
for the \textbf{normal coordinates }version of the definition.

The (infinitesimal) volume (change) function $\theta_{P}:$M$_{0}\rightarrow
$R$_{+}$ \ \ of the exponential map:

exp$_{\Pi}:$ B$\rightarrow$ M where $x=$ exp$_{\pi}(y,v)=$ exp$_{y}(v)$ \ and
exp$_{\text{y}}$:T$_{\text{y}}$M$\rightarrow$M \ is the usual exponential map
of the tangent bundle at y is given by:

$\left(  1.6\right)  \qquad\theta_{P}(x)=$ $\sqrt{\text{det}g(x)}$ where $g=$
$\left(  g_{ij}(x)\right)  $ $i,j=1,...,q,...,n$ is the matrix defined above.

Set,

$\left(  1.7\right)  \qquad\Psi(x)\ \ \ =\theta_{P}(x)^{-\frac{1}{2}}%
\ \ \Phi_{P}(x)\qquad\qquad$

$\left(  1.8\right)  \qquad q_{t}(x,P)=(2\pi t)^{-\frac{n-q}{2}}\Psi
(x)\exp\left\{  -\frac{d(x,P)^{2}}{2t}\right\}  \qquad\ $

We see from the definition above that:

$\left(  1.9\right)  \qquad\theta_{P}(y)=1\forall y\in U\subset$P where $U$ is
the small neighbourhood of $y_{0}\in P.$

Next we note that if x = y$\in$P, then by the definition of $\Phi_{P}$ above,
$\gamma$ is the unique minimal unit speed geodesic from x = y$\in$P meeting P
orthogonally in time 1. The geodesic is thus constant and so
$\overset{.}{\gamma}(s)=0$ for all $s\in\left[  0,1\right]  .$ We conclude
from $\left(  1.5\right)  $ that in this case,

$\left(  1.10\right)  \qquad\Phi_{P}(y)=1\forall y\in$P.

We conclude from the definition of q$_{\text{t}}$(x,P) in $\left(  1.8\right)
$ that for $y\in$P:

$\left(  1.11\right)  \qquad$q$_{\text{t}}$(y,P) = $(2\pi t)^{-\frac
{\text{n-q}}{2}}\Psi(y)\exp\left\{  -\frac{\text{d(y, P)}^{2}}{2t}\right\}
=(2\pi t)^{-\frac{n-q}{2}}$

The definition of Fermi coordinates in $\left(  1.1\right)  $ and $\left(
1.2\right)  $ implies that:

$\qquad\qquad\dim\exp_{\Pi}(B_{0})=\dim$M$_{0}=n$

Consequently,

$\left(  1.12\right)  \qquad$q$_{\text{t}}$(y,M$_{0}$) $=(2\pi t)^{-\frac
{\text{n-q}}{2}}=1$ since $q=\dim M_{0}=n$

We will see in \textbf{Chapter 5} that q$_{\text{t}}$(x,P) is the
\textbf{Euclidean part} of the generalized heat kernel.

When the submanifold P reduces to the point y$_{0},$ then the Fermi
coordinates reduce to the usual \textbf{normal} \textbf{coordinates}
(x$_{\text{1}}$,...,x$_{\text{n}}$) centered at the point y$_{0}\in$M.

Then $\sqrt{\text{det}(g_{\alpha\beta}(x))}=$ $\theta_{y_{0}}(x)$ reduces to
the Jacobian determinant of the exponential map exp$_{\text{y}_{0}}%
$:T$_{\text{y}_{0}}$M$\rightarrow$ M and then q$_{\text{t}}$(x,P) in $\left(
1.12\right)  $ reduces to:

$\left(  1.13\right)  \qquad$q$_{\text{t}}$(x,y$_{0}$) $=(2\pi t)^{-\frac
{\text{n}}{2}}\Psi(x)$ exp$\left\{  -\frac{\text{d(x, y}_{0}\text{)}^{2}%
}{2\text{t}}\right\}  .$

\chapter{Generalized Heat Flow}

\section{Heat Flow in Compact Manifolds}

Here we will derive a generalization of the notion of heat flow. The
description of heat flow is given in many texts (see for example
\textbf{Chavel }$\left[  10\right]  $ in the case of the scalar heat kernel).
The case of the heat equation in vector bundles is given in several texts
cited here. For example: \textbf{Avramidi} $\left[  2\right]  ,$
\textbf{Baudoin} $\left[  5\right]  ,$ \textbf{Driver and Thalmaier} $\left[
12\right]  $, \textbf{Gilkey} $\left[  20\right]  ,$ $\left[  21\right]  ,$
\textbf{Hsu} $\left[  30\right]  $ and \textbf{Lawson and Michelsohn} $\left[
35\right]  .$

For simplicity we will assume that M is a closed (compact without boundary)
n-dimensional Riemannian smooth manifold. We will state the heat equation in a
more general way:

Let $\Delta^{0}$ be the scalar Laplacian on smooth functions on M. Let X be a
smooth vector field on M and V a smooth potential term. W set:

\qquad\qquad\qquad\ L$^{0}=$ $\Delta^{0}+$ X + V

Cauchy's \textbf{evolution equation} for heat flow in M relative to L$^{0}$ is
given by:

$\left(  2.1\right)  $ $\qquad\qquad\frac{\partial f_{t}}{\partial t}=$
L$^{0}f_{t}\qquad\qquad$(evolution equation)

$\left(  2.2\right)  \qquad\qquad\ \underset{\text{t}\rightarrow0}{\lim
}\ \ f_{t}=f$ $\qquad\ \ \ \ $(initial condition)$\ \ $

The \textbf{fundamental solution} (or minimal \textbf{heat kernel})
p$_{\text{t}}(-,-)$ of the \textbf{heat equation} in $\left(  2.1\right)
-\left(  2.2\right)  $ above is a function p:M$\times$M$\times\left(
0,\infty\right)  \longrightarrow$R which is C$^{2}$ on M$\times$M and C$^{1}$
on $\left(  0,\infty\right)  $ which satisfies:

$\left(  2.3\right)  $ $\qquad\qquad\frac{\partial p_{t}}{\partial t}=$
L$^{0}$p$_{t}$ in the first variable x

$\left(  2.4\right)  \ \ \ \ \underset{\text{t}\rightarrow0}{\lim\text{ }}%
$p$_{\text{t}}(-,y)$ $=\delta_{y}$ in the distribution sense

The relation in $\left(  2.4\right)  $ means that:

$\left(  2.5\right)  \qquad\qquad\underset{\text{t}\rightarrow0}{\lim}%
\int_{\text{M}}$f(x)p$_{\text{t}}$(x,y)$\not \upsilon _{M}$(dx) = f(y) =
$\int_{M}$f(x)$\delta_{y}$(dx)

where $\upsilon_{M}$ is the Riemannian volume measure on M and $\delta_{y}$ is
the Dirac measure at y.

We follow \textbf{Chavel} $\left[  11\right]  ,$ pp 134-135):

The \textbf{physical interpretation} of $\left(  2.3\right)  -\left(
2.4\right)  $ is as follows:

p$_{\text{t}}(-,y)$ is the solution of the heat equation with initial (i.e. at
time t = 0) temperature distribution equal to 1 (a spark) completely
concentrated at y$\in$M.$\ $

The \textbf{physical interpretation} of $\left(  2.1\right)  -\left(
2.2\right)  :$

If, in general, we are given an initial tempreture distribution defined by the
function f on the manifold (rather than a spark at the point y), then it can
be proved that the tempreture of the manifold at a point x at time t%
$>$%
0 is given by:

$\left(  2.6\right)  $\ $\qquad\qquad$f$_{t}$(x) $=\int_{\text{M}}%
$f(y)p$_{\text{t}}$(x,y)$\not \upsilon _{M}$(dy)

This means that the equation in $\left(  2.1\right)  $ has for solution:

$\left(  2.7\right)  $\ $\qquad\qquad$f$_{t}$(x) $=\int_{\text{M}}%
$f(y)p$_{\text{t}}$(x,y)$\not \upsilon _{M}$(dy)

where the limit in the initial condition in $\left(  2.2\right)  $ is to be
understood in the weak sense i.e. in the distribution sense: for any
continuous function g:M$\longrightarrow$ R, we have:

$\left(  2.8\right)  $\qquad$\qquad\underset{\text{t}\rightarrow0}{\lim}%
\int_{\text{M}}$f$_{t}$(x)g(x)$\not \upsilon _{M}$(dx) = $\int_{\text{M}}%
$f(x)g(x)$\not \upsilon _{M}$(dx)

Let (exp$_{y}^{-1},$M$_{0}\mathbf{)}$ be the domain of a geodesic chart based
at y. The proposition below gives an expression for the heat kernel of the
compact Riemannian manifold M. This uses the following Cauchy equation for the
heat flow.

\qquad\qquad\qquad\qquad\qquad\qquad\qquad\qquad\qquad\qquad\qquad\qquad
\qquad\qquad\qquad\qquad\qquad$\blacksquare$

\begin{proposition}
Let the Cauchy equation be given as follows:
\end{proposition}

$\qquad\qquad\frac{\partial f_{t}^{\lambda}}{\partial t}=L^{0}f^{\lambda}$

$\qquad\qquad\ \ f_{0}^{\lambda}=(2\pi\lambda)^{-\frac{n}{2}}f\exp\left\{
-\dfrac{d(-,y)^{2}}{2\lambda}\right\}  $

where f is a smooth function on the compact Riemannian manifold M. Then,

$\qquad\qquad\underset{\lambda\rightarrow0}{lim}$f$_{t}^{\lambda}(x)=p_{t}%
^{M}(x,y)f(y)$

\begin{proof}
We recall that the set denoted here by $B_{0}$ here is denoted by $\Theta_{P}$
in Gray $\left[  23\right]  $ and Gray $\left[  25\right]  .$ The statement of
\textbf{Theorem }$8.40$ of \textbf{Gray }$\left[  25\right]  $ has a slight
error in the notation. We correct here:
\end{proof}

When Fermi coordinates reduce to normal coordinates, we have by $\left(
1.4\right)  $ (ii) above:

M$_{0}=$ exp$_{y}($B$_{0})$ and since $\upsilon_{\text{M}}(M)=$ $\upsilon
_{\text{M}}(\exp_{y}(B_{0}))$ by $\left(  3.45\right)  $ of \textbf{Gray
}$\left[  25\right]  ,$ we have:

$\upsilon_{\text{M}}(M)=$ $\upsilon_{\text{M}}(\exp_{y}(B_{0}))=$
$\upsilon_{\text{M}}(M_{0})$

Therefore the solution f$_{\text{t}}^{\lambda}($x) of the above Cauchy Problem
is given by:

f$_{\text{t}}^{\lambda}($x) = $\int_{\text{M}}$f$_{0}^{\lambda}$%
(z)p$_{\text{t}}^{\text{M}}$(x,z)$\upsilon_{\text{M}}$(dz)\ = $\int%
_{\text{M}_{0}}$f$_{0}^{\lambda}$(z)p$_{\text{t}}^{\text{M}}$(x,z)$\upsilon
_{\text{M}}$(dz)

\qquad\ = $(2\pi\lambda)^{-\frac{{\ n}}{2}}$ $\int_{\text{M}_{0}}$%
f(z)$\exp\left\{  -\dfrac{\text{d(z,y)}^{2}}{2\lambda}\right\}  $p$_{\text{t}%
}^{\text{M}}$(x,z)$\upsilon_{\text{M}}$(dz)

Since the integration is over M$_{0}=$ exp$_{y}($B$_{0})$, we make the change
of variable: z = exp$_{y}$v and have:

f$_{\text{t}}^{\lambda}($x) = $(2\pi\lambda)^{-\frac{{\ n}}{2}}$
$\int_{\text{T}_{\text{y}}\text{M}}$f($\exp_{\text{y}}$v)$\exp\left\{
-\dfrac{\text{d(}\exp_{\text{y}}\text{v,y)}^{2}}{2\lambda}\right\}
$p$_{\text{t}}^{\text{M}}$(x,$\exp_{\text{y}}$v)$\theta_{\text{y}}$(v)dv

= $(2\pi\lambda)^{-\frac{{\ n}}{2}}$ $\int_{\text{T}_{\text{y}}\text{M}}%
$f($\exp_{\text{y}}$v)$\exp\left\{  -\frac{\left\Vert \text{v}\right\Vert
^{2}}{2\lambda}\right\}  $p$_{\text{t}}^{\text{M}}$(x,$\exp_{\text{y}}%
$v)$\theta_{\text{y}}$(v)dv

where

$\left(  2.9\right)  \qquad\theta_{\text{y}}(v)$ = $\sqrt{\text{det(g}%
_{\alpha\beta}\text{(x))}}=$ $\left\Vert \det\text{T}_{\text{v}}\exp
_{\text{y}}\right\Vert $ is the Jacobian determinant of the exponential map.

We then make another slight change of variable: v = $\sqrt{\lambda}$w and have:

$\left(  2.10\right)  \qquad\qquad$f$_{\text{t}}^{\lambda}$(x) = (2$\pi
$)$^{-\frac{{n}}{2}}$ $\int_{\text{T}_{\text{y}}\text{M}}$f($\exp_{\text{y}%
}\sqrt{\lambda}$w)$\exp\left\{  -\frac{\left\Vert \text{{\ w}}\right\Vert
^{2}}{2}\right\}  $p$_{\text{t}}^{\text{M}}$(x,$\exp_{\text{y}}\sqrt{\lambda}%
$w)$\theta_{\text{y}}$($\sqrt{\lambda}$w)dw

Then taking limits, as $\lambda\downarrow0,$ on both sides of $\left(
2.8\right)  ,$ we have:

$\underset{\lambda\rightarrow0}{\text{lim}}$f$_{\text{t}}^{\lambda}$(x) =
(2$\pi$)$^{-\frac{{\ n}}{2}}$ $\int_{\text{T}_{\text{y}}\text{M}}$%
f($\exp_{\text{y}}$0)$\exp\left\{  -\frac{\left\Vert \text{w}\right\Vert ^{2}%
}{2}\right\}  $p$_{\text{t}}^{\text{M}}$(x,$\exp_{\text{y}}$0)$\theta
_{\text{y}}$(0)dw

Since $\int_{\text{T}_{\text{y}}\text{M}}$f(y)$\exp\left\{  -\frac{\left\Vert
\text{{\ w}}\right\Vert ^{2}}{2}\right\}  $dw = (2$\pi$)$^{\frac{{n}}{2}}$f(y)
and $\exp_{\text{y}}$0) = y$,$ we have:

$\left(  2.11\right)  \qquad\underset{\lambda\rightarrow0}{\text{lim}}%
$f$_{\text{t}}^{\lambda}$(x) = p$_{\text{t}}^{\text{M}}$(x,y)f(y)

where p$_{\text{t}}^{\text{M}}$(x,y) is the heat kernel of the compact
Riemannian manifold M.

\qquad\qquad\qquad\qquad\qquad\qquad\qquad\qquad\qquad\qquad\qquad\qquad
\qquad\qquad\qquad\qquad\qquad\qquad\qquad$\blacksquare$

\begin{proposition}
Suppose we have:

$\qquad\frac{\partial\text{f}_{\text{t}}}{\partial\text{t}}^{\lambda}=$
L$^{0}$f$_{\text{t}}^{\lambda}$

\qquad f$_{0}^{\lambda}=$ $(2\pi\lambda)^{-\frac{({n-q)}}{2}}$f$\exp\left\{
-\dfrac{d(-,P)^{2}}{2\lambda}\right\}  $

where f is a smooth function on the compact Riemannian manifold M and P is a
compact submanifold of M. Then,

$\underset{\lambda\rightarrow0}{\text{lim}}$f$_{\text{t}}^{\lambda}($x) =
$\int_{\text{P}}$p$_{\text{t}}^{\text{M}}$(x,y)f(y)$\upsilon_{\text{P}}$(dy)
\end{proposition}

\begin{proof}
We recall that the set denoted here by $B_{0}$ is denoted by $\Theta_{P}$ in
Gray $\left[  23\right]  $ and Gray $\left[  25\right]  .$ The statement of
\textbf{Theorem }$8.40$ of \textbf{Gray }$\left[  25\right]  $ has a slight
error in the notation:
\end{proof}

The first line of Equation in $\left(  8.69\right)  $ should read:

\begin{center}
$\upsilon_{\text{M}}(M)=\upsilon_{\text{M}}(\exp_{\Pi}(B_{0}))$ (and not
$\upsilon_{\text{M}}(M)=\upsilon_{\text{M}}(B_{0})):$
\end{center}

We could also, have used Lemma $\left(  2.2\right)  $, Chapter 2 of
\textbf{Ndumu} $\left[  40\right]  $ and assume that the exponential map of
the normal bundle is \ a diffeomorphism of M onto B$_{0}.$

By $\left(  1.4\right)  $ (ii) of Chapter 1,

\begin{center}
M$_{0}=$ exp$_{\Pi}($B$_{0})$
\end{center}

Therefore by \textbf{Theorem }$8.40$ of \textbf{Gray }$\left[  26\right]  ,$

\begin{center}
$\upsilon_{\text{M}}(M)=$ $\upsilon_{\text{M}}(\exp_{\Pi}(B_{0}))=$
$\upsilon_{\text{M}}(M_{0})$
\end{center}

Consequently,

\begin{center}
$f_{\text{t}}^{\lambda}(x)=$ $\int_{\text{M}}$f$_{0}^{\lambda}$(z)p$_{\text{t}%
}^{\text{M}}$(x,z)$\upsilon_{\text{M}}$(dz)\ = $\int_{\text{M}_{0}}$%
f$_{0}^{\lambda}$(z)p$_{\text{t}}^{\text{M}}$(x,z)$\upsilon_{\text{M}}$(dz)
\end{center}

Since,

\begin{center}
f$_{0}^{\lambda}=$ $(2\pi\lambda)^{-\frac{({n-q)}}{2}}$f$\exp\left\{
-\dfrac{d(-,P)^{2}}{2\lambda}\right\}  $

\qquad$f_{\text{t}}^{\lambda}(x)=$ $(2\pi\lambda)^{-\frac{({n-q)}}{2}}%
\int_{\text{M}_{0}}$f(z)$\exp\left\{  -\dfrac{d(z,P)^{2}}{2\lambda}\right\}
$p$_{\text{t}}^{\text{M}}$(x,z)$\upsilon_{\text{M}}$(dz)
\end{center}

Since the integration is over M$_{0}=$ exp$_{\Pi}($B$_{0})$, we make the
change of variable:

\qquad\ z = exp$_{\Pi}(y,$v) = exp$_{y}($v) and have:

\begin{center}
$f_{\text{t}}^{\lambda}(x)=$ $(2\pi\lambda)^{-\frac{{n-q}}{2}}$ $\int%
_{\text{P}\times\text{R}^{n-q}}\exp\left\{  -\dfrac{d(\exp_{\text{y}}%
\text{v},P)^{2}}{2\lambda}\right\}  $p$_{\text{t}}^{\text{M}}$(x,$\exp
_{\text{y}}$v)f($\exp_{\text{y}}$v)$\theta_{\text{y}}$(v)dydv

$\bigskip f_{\text{t}}^{\lambda}(x)=$ $(2\pi\lambda)^{-\frac{{n-q}}{2}}$
$\int_{\text{P}\times\text{R}^{n-q}}\exp\left\{  -\dfrac{\left\Vert
\text{v}\right\Vert ^{2}}{2\lambda}\right\}  $p$_{\text{t}}^{\text{M}}%
$(x,$\exp_{\text{y}}$v)f($\exp_{\text{y}}$v)$\theta_{\text{y}}$(v)dydv
\end{center}

Next, we set: v = $\sqrt{\lambda}$w = ($\sqrt{\lambda}$w$_{q+1},...,\sqrt
{\lambda}$w$_{n})$ and hence dv = $(\sqrt{\lambda})^{\frac{n-q}{2}}$dw. Consequently,

\begin{center}
$f_{\text{t}}^{\lambda}(x)=$ $(2\pi\lambda)^{-\frac{{n-q}}{2}}$ $\int%
_{\text{P}\times\text{R}^{n-q}}\exp\left\{  -\frac{\lambda\left\Vert
\text{w}\right\Vert ^{2}}{2\lambda}\right\}  $p$_{\text{t}}^{\text{M}}%
(x,\exp_{\text{y}}\sqrt{\lambda}w)f(\exp_{\text{y}}\sqrt{\lambda}%
w)\theta_{\text{y}}(\sqrt{\lambda}w)(\lambda)^{\frac{n-q}{2}}$dydw

$=$ $(2\pi)^{-\frac{{n-q}}{2}}$ $\int_{\text{P}\times\text{R}^{n-q}}%
\exp\left\{  -\frac{\left\Vert \text{w}\right\Vert ^{2}}{2}\right\}
$p$_{\text{t}}^{\text{M}}(x,\exp_{\text{y}}\sqrt{\lambda}w)f(\exp_{\text{y}%
}\sqrt{\lambda}w)\theta_{\text{y}}(\sqrt{\lambda}w)$dydw
\end{center}

Taking limits on both sides of the last equation we have:

$\qquad\qquad\underset{\lambda\rightarrow0}{\text{lim}}$f$_{\text{t}}%
^{\lambda}($x) \ = (2$\pi$)$^{-\frac{{\ n-q}}{2}}$ $\int_{\text{R}%
^{\text{n-q}}}\exp\left\{  -\frac{\left\Vert \text{w}\right\Vert ^{2}}%
{2}\right\}  $dw$\int_{\text{P}}$p$_{\text{t}}^{\text{M}}$(x,y)f(y)$\upsilon
_{\text{P}}$(dy)

Since,

\qquad\qquad\ $\int_{\text{R}^{\text{n-q}}}\exp\left\{  -\frac{\left\Vert
\text{w}\right\Vert ^{2}}{2}\right\}  $dw = (2$\pi$)$^{\frac{{n-q}}{2}},$ we have:

$\left(  2.12\right)  $\qquad$\qquad\underset{\lambda\rightarrow0}{\text{lim}%
}$f$_{\text{t}}^{\lambda}($x) = $\int_{\text{P}}$p$_{\text{t}}^{\text{M}}%
$(x,y)f(y)$\upsilon_{\text{P}}$(dy) $\doteq$ p$_{\text{t}}^{\text{M}}$(x,P)
\qquad\qquad\qquad\qquad\qquad\qquad\qquad\qquad\qquad\qquad\qquad\qquad
\qquad\qquad\qquad\qquad\qquad$\qquad\qquad\qquad\qquad\qquad\qquad
\qquad\qquad\qquad\qquad\qquad\qquad\qquad\qquad\qquad\blacksquare$

\begin{proposition}
(i)$\qquad\frac{\partial\text{p}_{\text{t}}^{\text{M}}\text{(}-\text{,P)}%
}{\partial\text{t}}$\ = L$^{0}$p$_{\text{t}}^{\text{M}}$($-$,P)

(ii)\qquad\ $\underset{\text{t}\rightarrow0}{lim}$ p$_{\text{t}}^{\text{M}}%
$($-$,P) $=\delta_{\text{P}}$ in the distribution sense

where $\delta_{\text{P}}$ is a measure on P with density f with respect to the
Riemannian volume measure $\upsilon_{\text{P}}$ on P.$\qquad$
\end{proposition}

\begin{proof}
$\frac{\partial\text{p}_{\text{t}}^{\text{M}}\text{(x,P)}}{\partial\text{t}}$
$\ =$ $\ $ $\frac{\partial}{\partial\text{t}}\int_{\text{P}}$p$_{\text{t}%
}^{\text{M}}$(x,y)f(y)$\upsilon_{\text{P}}$(dy)
\end{proof}

\qquad\qquad$=$ $\int_{\text{P}}\frac{\partial\text{p}_{\text{t}}^{\text{M}%
}\text{(x,y)}}{\partial\text{t}}$f(y)$\upsilon_{\text{P}}$(dy) by
differentiating under the integral sign.

\qquad\qquad\ $=$ $\int_{\text{P}}$L$^{0}$p$_{\text{t}}^{\text{M}}%
$(x,y)f(y)$\upsilon_{\text{P}}$(dy) by $\left(  2.3\right)  $

$\qquad\qquad=$ L$^{0}\int_{\text{P}}$p$_{\text{t}}^{\text{M}}$%
(x,y)f(y)$\upsilon_{\text{P}}$(dy)

The last equality above is due to the fact that the operator L$_{0}$ is
applied to the variable x which is independent of the integration. We have:

$\qquad\qquad\frac{\partial\text{p}_{\text{t}}^{\text{M}}\text{(x,P)}%
}{\partial\text{t}}=$ L$^{0}$p$_{\text{t}}^{\text{M}}$(x,P) by the definition
of p$_{\text{t}}^{\text{M}}$(x,P) in $\left(  2.12\right)  .$

(ii) Using the definition of p$_{\text{t}}^{\text{M}}$(x,P) in $\left(
2.13\right)  ,$ we have for a smooth function f$_{0}$:M$\longrightarrow$R,

$\qquad\qquad\underset{t\longmapsto0}{\lim}\int_{\text{M}}$p$_{t}^{\text{M}}%
$(x,P)f$_{0}$(x)$\upsilon_{\text{M}}$(dx) = $\underset{t\longmapsto0}{\lim
}\int_{\text{M}}\left\{  \int_{\text{P}}\text{p}_{t}^{\text{M}}%
\text{(x,y)f(y)}\upsilon_{\text{P}}\text{(dy)}\right\}  $f$_{0}$%
(x)$\upsilon_{\text{M}}$(dx)

\qquad\qquad\ \ $\ \underset{t\longmapsto0}{=\lim}\int_{\text{P}}\left\{
\int_{\text{M}}\text{p}_{t}^{\text{M}}\text{(x,y)f}_{0}\text{(x)}%
\upsilon_{\text{M}}\text{(dx)}\right\}  $f(y)$\upsilon_{P}$(dy) by the Fubini Theorem.

\qquad\qquad\ \ \ \ $=$ $\int_{\text{P}}\left\{  \underset{t\longmapsto
0}{\lim}\int_{\text{M}}\text{p}_{t}^{\text{M}}\text{(x,y)f}_{0}\text{(x)}%
\upsilon_{\text{M}}\text{(dx)}\right\}  $f(y)$\upsilon_{P}$(dy)

\qquad\qquad\ \ \ \ $=$ $\int_{\text{P}}$f$_{0}$(y)f(y)$\upsilon_{P}$(dy) by
$\left(  2.5\right)  $

\qquad\qquad\ $\ \ =$ $\int_{\text{P}}$f$_{0}$(y)$\delta_{P}$(dy)

The first and last equalities give the important relation:

$\left(  2.13\right)  $\qquad$\underset{t\longmapsto0}{\lim}\int_{\text{M}}%
$p$_{t}^{\text{M}}$(x,P)f$_{0}$(x)$\upsilon_{\text{M}}$(dx) = $\int_{\text{P}%
}$f$_{0}$(y)$\delta_{P}$(dy)

This means that:

$\qquad\qquad\qquad$p$_{t}^{\text{M}}$(-,P) $\longmapsto\delta_{\text{P}}$ as
$t\rightarrow0$ in the distribution sense,

where $\delta_{\text{P}}$ is a measure on P defined by:

$\qquad\qquad\qquad\delta_{\text{P}}$(dy) = f(y)$\upsilon_{\text{P}}$(dy).

It is a measure on P with density f with respect to the Riemannian volume
measure $\upsilon_{\text{P}}$ on P.

We have thus proved (ii).

\qquad\qquad\qquad\qquad\qquad\qquad\qquad\qquad\qquad\qquad\qquad\qquad
\qquad\qquad\qquad\qquad\qquad$\blacksquare$

\begin{definition}
The Generalized Heat Kernel relative to a Compact Riemannian manifold.
\end{definition}

Comparing the equations in $\left(  2.11\right)  $ and\textbf{ }$\left(
2.12\right)  $, and given the properties proved in the last Proposition above,
we are logically led to define the \textbf{generalized heat kernel}
p$_{\text{t}}^{\text{M}}$(x,P) by:

$\left(  2.14\right)  $\qquad p$_{\text{t}}^{\text{M}}$(x,P) = $\int%
_{\text{P}}$p$_{t}^{\text{M}}$(x,y)f(y)$\upsilon_{P}$(dy)

\qquad\qquad\qquad\qquad\qquad\qquad\qquad\qquad\qquad\qquad\qquad\qquad
\qquad\qquad\qquad\qquad\qquad$\blacksquare$

The \textbf{physical interpretation} of the Proposition is as follows:

p$_{t}^{\text{M}}$(-,P) is a solution of (i) under the intial tempreture
distribution equal to f completely concentrated on the submanifold P.

We see that (i) and (ii) of the above Proposition generalize $\left(
2.3\right)  $ and $\left(  2.4\right)  $ above.

If f $\equiv1$ on P, then $\delta_{\text{P}}$= $\upsilon_{\text{P}}$ is the
Riemannian volume measure on P. If further P reduces to the center of Fermi

coordinates y$_{0}$, then $\delta_{\text{P}}$ reduces to the Dirac measure
$\delta_{\text{y}_{0}}$ at y$_{0}$ and we have by $\left(  2.14\right)  :$

$\left(  2.15\right)  $\qquad\ $\underset{t\longmapsto0}{\lim}\int_{\text{M}}%
$p$_{t}^{\text{M}_{0}}$(x,y$_{0}$)f$_{0}$(x)$\upsilon_{\text{M}}$(dx) =
$\int_{\text{M}}$f$_{0}$(x)$\delta_{\text{y}_{0}}$(dx) = f$_{0}$(y$_{0}$)

We see that we recover $\left(  2.5\right)  $ above.

\qquad\qquad\qquad\qquad\qquad\qquad\qquad\qquad\qquad\qquad\qquad\qquad
\qquad\qquad\qquad\qquad\qquad\qquad\qquad$\blacksquare$

\section{Heat Flow in Vector Bundles}

Let M be a compact n-dimensional Riemannian smooth manifold and let E be a
\textbf{vector bundle} over M.

We can go a step further and consider the heat kernel k$_{\text{t}}$(x,y):
E$_{\text{y}}\longrightarrow$ E$_{\text{x}}$ of a vector bundle E over the
compact Riemannian manifold M relative to the operator L = $\frac{1}{2}%
\Delta+$ X $+$ V (where we write X for $\nabla_{\text{X}}$). It is a
homomorphism of the fibers of the vector bundle E. It is the fundamental
solution of the heat equation on E:

For $\phi_{t},\phi\in\Gamma(E),$ where $\Gamma(E)$ is the space of
\textbf{smooth sections} over the vector bundle E,

$\left(  2.16\right)  $ $\qquad\frac{\partial\phi_{\text{t}}}{\partial
\text{t}}=$ L$\phi_{\text{t}}$ \qquad\qquad(evolution equation)

$\qquad\ \ \ \ \qquad\ \ \phi_{\text{0}}=$ $\phi$\qquad\qquad
\ \ \ \ \ \ (initial condition)

where L = $\frac{1}{2}\Delta+$ X $+$ V and $\Delta$ is the Laplacian on vector
bundles defined in $\left(  3.1\right)  $ below.

The above initial condition means that $\phi_{\text{0}}$ $\circeq
$\ $\underset{\text{t}\rightarrow0}{lim}\phi_{\text{t}}$ = $\phi$ in the
distribution sense (see p. 33 of \textbf{Gilkey} $\left[  21\right]  $):

For any specific heat $\rho\in\Gamma(E^{\ast}),$ we have:

\qquad\qquad\ $\underset{\text{t}\rightarrow0}{lim}$ $\int_{\text{M}}%
<\phi(t,x),\rho(x)>\upsilon_{\text{M}}$(dx) = $\int_{\text{M}}<\phi
(x),\rho(x)>\upsilon_{\text{M}}$(dx)

The solution $\phi_{\text{t}}$ of the above evolution equation in $\left(
2.16\right)  $ is given by:

$\left(  2.17\right)  \qquad\phi_{\text{t}}($x$)=$ $\int_{\text{M}}%
$k$_{\text{t}}$(x,y)$\phi$(y)$\upsilon_{\text{M}}$(dy)

The map: k$_{\text{t}}$(x,y):E$_{\text{y}}\longrightarrow$ E$_{\text{x}}$ is
the vector bundle \textbf{heat kernel} (see, for example Definition $2.4$ of
Bismut $\left[  8\right]  $ or section 3 of \textbf{Baudoin} $\left[
5\right]  $).

\qquad\qquad\qquad\qquad\qquad\qquad\qquad\qquad\qquad\qquad\qquad\qquad
\qquad\qquad\qquad\qquad$\blacksquare$

The following is the vector bundle version of $\left(  2.4\right)  -\left(
2.5\right)  $ above:

The vector bundle heat kernel satisfies the heat equation in the first variable:

$\left(  2.18\right)  $ $\qquad\frac{\partial\text{k}_{\text{t}}(-,\text{y}%
)}{\partial\text{t}}=$ Lk$_{\text{t}}$($-$,y)\qquad\qquad\qquad(evolution equation)

$\qquad\ \ \ $k$_{\text{t}}(-,y)\longmapsto\delta_{y}$ as t$\longmapsto0$ in
the distribution sense (initial condition)

\qquad\qquad\qquad\qquad\qquad\qquad\qquad\qquad\qquad\qquad\qquad\qquad
\qquad\qquad\qquad\qquad$\blacksquare$

The above initial condition means that for $\phi\in\Gamma(E),$

\qquad$\underset{\text{t}\rightarrow0}{\lim}\int_{\text{M}}\phi$%
(y)k$_{\text{t}}$(x,y)$\not \upsilon _{M}$(dy) $\longrightarrow$ $\phi$(x)
uniformly in x as t $\longrightarrow0.$

Equivalently,

\qquad$\underset{\text{t}\rightarrow0}{\lim}\int_{\text{M}}\phi$%
(x)k$_{\text{t}}$(x,y)$\not \upsilon _{M}$(dx) $\longrightarrow$ $\phi$(y) =
$\int_{\text{M}}\phi$(x)$\delta_{y}$(dx)

The following is a generalization of \textbf{Proposition 1} here to the case
of vector bundles.

\qquad\qquad\qquad\qquad\qquad\qquad\qquad\qquad\qquad\qquad\qquad\qquad
\qquad\qquad\qquad\qquad$\blacksquare$

\begin{proposition}
\qquad\qquad\qquad\qquad\qquad\qquad\qquad\qquad\qquad\qquad\qquad\qquad
\qquad\qquad
\end{proposition}

Let $\phi_{\text{t}}^{\lambda}\in\Gamma(E)$ be a family of smooth sections of
the vector bundle E. Consider the

following heat equation on E:

$\frac{\partial\phi_{\text{t}}^{\lambda}}{\partial\text{t}}=$ L$\phi
_{\text{t}}^{\lambda}$

$\phi_{0}^{\lambda}=$ $($2$\pi\lambda)^{-\frac{({n-q)}}{2}}\phi\exp\left\{
-\frac{\text{{\ d(-,P)}}^{2}}{2\lambda}\right\}  $

where $\phi$ is a smooth section of E and P is a \textbf{compact submanifold}
of the \textbf{compact manifold} M. Then,

$\underset{\lambda\rightarrow0}{\text{lim}}\phi_{\text{t}}^{\lambda}($x) =
$\int_{\text{P}}$k$_{\text{t}}($x,y$)\phi(y)\upsilon_{\text{P}}($dy)

This is the vector bundle version of \textbf{Proposition }$2.2$ above:

\begin{proof}
As previously given: $\upsilon_{\text{M}}(M)=$ $\upsilon_{\text{M}}(\exp_{\Pi
}(B_{0}))=$ $\upsilon_{\text{M}}(M_{0})$
\end{proof}

and so we have:

$\left(  2.19\right)  \qquad\phi_{\text{t}}^{\lambda}=\int_{\text{M }}%
$k$_{\text{t}}$(x,z)$\phi_{0}^{\lambda}(z)\upsilon_{\text{M}}$(dz) =
$\int_{\text{M}_{0}\text{ }}$k$_{\text{t}}$(x,z)$\phi_{0}^{\lambda}%
(z)\upsilon_{\text{M}}$(dz)

The computations are the same as in \textbf{Proposition 2.2}:

$\phi_{\text{t}}^{\lambda}$ = (2$\pi\lambda)^{-\frac{{n}}{2}}$ $\int%
_{\text{M}_{0}}\exp\left\{  -\frac{\text{{\ d(z,P)}}^{2}}{2\lambda}\right\}
$k$_{\text{t}}$(x,z)$\phi_{0}^{\lambda}(z)\upsilon_{\text{M}}$(dz)

= (2$\pi\lambda)^{-\frac{{n}}{2}}$ $\int_{\text{B}_{0}}\exp\left\{
-\frac{\text{{\ d(exp}}_{\text{y}}\text{v{\ ,P)}}^{2}}{2\lambda}\right\}
$k$_{\text{t}}$(x,$\exp_{\text{y}}$v)$\phi(\exp_{\text{y}}v)\theta_{\text{y}}$(v)dv

= (2$\pi\lambda)^{-\frac{{n}}{2}}$ $\int_{\text{B}_{0}}\exp\left\{
-\frac{\left\Vert \text{v}\right\Vert ^{2}}{2\lambda}\right\}  $k$_{\text{t}}%
$(x,$\exp_{\text{y}}$v)$\phi(\exp_{\text{y}}$v$)\theta_{\text{y}}$(v)dv

= (2$\pi)^{-\frac{{n-q}}{2}}$ $\int_{\text{P}\times\text{R}^{n-q}}\exp\left\{
-\frac{\left\Vert \text{w}\right\Vert ^{2}}{2}\right\}  $k$_{\text{t}}%
$(x,$\exp_{\text{y}}\sqrt{\lambda}$w)$\phi(\exp_{\text{y}}\sqrt{\lambda}%
$w$)\theta_{\text{y}}$($\sqrt{\lambda}$w)dydw

where we have set v = $\sqrt{\lambda}$w. Since $\int_{\text{R}^{\text{n-q}}%
}\exp\left\{  -\frac{\left\Vert \text{w}\right\Vert ^{2}}{2}\right\}  $dw =
(2$\pi)^{\frac{{n-q}}{2}},$ we have:

$\underset{\lambda\rightarrow0}{\text{lim}}$f$_{\text{t}}^{\lambda}($x) $\ =$
(2$\pi)^{-\frac{{n-q}}{2}}$ $\int_{\text{R}^{\text{n-q}}}\exp\left\{
-\frac{\left\Vert \text{w}\right\Vert ^{2}}{2}\right\}  $dw.$\int_{\text{P}}$k
$_{\text{t}}$(x,y)$\phi(y)\upsilon_{\text{P}}$(dy)

\qquad\qquad$\ =$ $\int_{\text{P}}$k$_{\text{t}}$(x,y)$\phi(y)\upsilon
_{\text{P}}$(dy)

The results of \textbf{Proposition 1 }and\textbf{ Proposition 2 }lead us to
define the \textbf{generalized vector bundle heat kernel }by:

$\left(  2.20\right)  $\qquad k$_{\text{t}}$(x,P,$\phi$) \ = $\int_{\text{P}}%
$k$_{\text{t}}$(x,y)$\phi(y)\upsilon_{\text{P}}$(dy)

\qquad\qquad\qquad\qquad\qquad\qquad\qquad\qquad\qquad\qquad\qquad\qquad
\qquad\qquad\qquad\qquad\qquad$\blacksquare$

The vector bundle heat kernel k$_{\text{t}}$(x,y) is a homomorphism
k$_{\text{t}}$(x,y)$:$ E$_{\text{y}}\rightarrow$E$_{\text{x}}$ of the fibers
of the vector bundle E. Later we shall give an explicit stochastic
representation of k$_{\text{t}}$(x,P,$\phi$) from which we deduce that of
k$_{\text{t}}$(x,y)$\phi$(y).

\qquad\qquad\qquad\qquad\qquad\qquad\qquad\qquad\qquad\qquad\qquad\qquad
\qquad\qquad\qquad\qquad\qquad$\blacksquare$

(i)$\qquad\frac{\partial\text{k}_{\text{t}}\text{(x,P,}\phi\text{)}}%
{\partial\text{t}}$ $\ =$ $\ $ Lk$_{\text{t}}$(x,P,$\phi$)

(ii) \qquad k$_{\text{t}}$(x,P,$\phi$)$\longrightarrow\delta_{P}$ in the
distribution sense.

\begin{proof}
We differentiate under and out of the integral sign and have:
\end{proof}

$\left(  2.21\right)  \qquad$ $\frac{\partial\text{k}_{\text{t}}%
\text{(x,P,}\phi\text{)}}{\partial\text{t}}$ $\ =$ $\ $ $\frac{\partial
}{\partial\text{t}}\int_{\text{P}}$k$_{\text{t}}$(x,y)$\phi(y)\upsilon
_{\text{P}}$(dy) = $\int_{\text{P}}\frac{\partial}{\partial\text{t}}%
$k$_{\text{t}}$(x,y)$\phi(y)\upsilon_{\text{P}}$(dz)

\qquad\qquad\qquad\qquad\ \ $\ \ =$ $\int_{\text{P}}$Lk$_{\text{t}}$%
(x,y)$\phi(y)\upsilon_{\text{P}}$(dy)

$\qquad\qquad\qquad\qquad\ \ \ \ =$ \ \ \ \ L$\int_{\text{P}}$k$_{\text{t}}%
$(x,y)$\phi(y)\upsilon_{\text{P}}$(dy) $=$ Lk$_{\text{t}}$(x,P,$\phi$)

(ii)\qquad Using the definition of k$_{\text{t}}$(x,P,$\phi$) in $\left(
2.20\right)  ,$ we have for a smooth function $\phi_{0}\in\Gamma(E)$

$\qquad\qquad\underset{t\longmapsto0}{\lim}\int_{\text{M}}<$k$_{\text{t}}%
$(x,P,$\phi$),$\phi_{0}$(x)$>\upsilon_{\text{M}}$(dx)

\qquad\qquad\qquad= $\underset{t\longmapsto0}{\lim}\int_{\text{M}}\left\{
<\int_{\text{P}}\text{k}_{\text{t}}\text{(x,y)}\phi\text{(y)}\upsilon
_{\text{P}}\text{(dy)}\right\}  ,\phi_{0}$(x)$>\upsilon_{\text{M}}$(dx)

\qquad\qquad\ \ $\ \underset{t\longmapsto0}{=\lim}\int_{\text{P}}\left\{
<\int_{\text{M}}\text{k}_{\text{t}}\text{(x,y)}\phi_{0}\text{(x)}%
\upsilon_{\text{M}}\text{(dx)}\right\}  ,\phi$(y)$>\upsilon_{P}$(dy) by the
Fubini Theorem.

\qquad\qquad\ \ \ \ = $\int_{\text{P}}\left\{  <\underset{t\longmapsto0}{\lim
}\int_{\text{M}}\text{k}_{\text{t}}\text{(x,y)}\phi_{0}\text{(x)}%
\upsilon_{\text{M}}\text{(dx)}\right\}  ,\phi$(y)$>\upsilon_{P}$(dy)

\qquad\qquad\ \ \ \ = $\int_{\text{P}}<\phi_{0}$(y)$,\phi$(y)$>\upsilon_{P}%
$(dy) by $\left(  2.8\right)  $

This confirms that k$_{\text{t}}^{M}$(x,P,$\phi$) is a generalization of the
heat kernel k$_{\text{t}}$(x,y) of the vector bundle E over the compact
manifold M:

\begin{definition}
The map:
\end{definition}

$\left(  2.22\right)  $\qquad k$_{\text{t}}$(x,P,$\phi$) \ = $\int_{\text{P}}%
$k$_{\text{t}}$(x,y)$\phi(y)\upsilon_{\text{P}}$(dy)

is the \textbf{generalized vector bundle heat kernel.}

\qquad\qquad\qquad\qquad\qquad\qquad\qquad\qquad\qquad\qquad\qquad\qquad
\qquad\qquad\qquad\qquad\qquad\qquad$\blacksquare$

\part{SOME VECTOR BUNDLE CALCULUS}

\chapter{Connections on Vector Bundles}

We will first give some definitions and useful properties for vector bundles
and connections on them.

We will follow the following authors for definitions:

\textbf{Definition }$\left(  1.1.4\right)  $ of \textbf{Berline, Getzler and
Vergne} $\left[  7\right]  $, \textbf{Definition }$\left(  4.3\right)  $ of
\textbf{Lawson and Michelsohn }$\left[  35\right]  ,$ $\left(  5.3\right)  $
and $\left(  5.4\right)  $ of \textbf{Morita }$\left[  39\right]  $ and
\textbf{Definition }$\left(  1.1\right)  $ of \textbf{Roe }$\left[  48\right]
.$

\begin{definition}
(Connection, Metric Connection, Extended Connection)
\end{definition}

Let $\Gamma(E)$ denote the space of \textbf{smooth sections} of the
\textbf{vector bundle E} over a closed (compact without boundary) Riemannian
manifold M of dimension n, and let $\Gamma(T^{\ast}M\otimes E)$ denote the
space of \textbf{smooth sections }of the \textbf{vector bundle of 1-forms}
with values in the \textbf{vector bundle E.}

(i) A \textbf{linear connection} is a \textbf{covariant derivative}
$\nabla:\Gamma(E)\rightarrow\Gamma(T^{\ast}M\otimes E)$ which takes smooth
sections of E to smooth sections of the (tensor product) vector bundle
$T^{\ast}M\otimes E$ such that:

\begin{center}
$\nabla\phi\in\Gamma(T^{\ast}M\otimes E)$ and $\nabla_{X}\phi\in\Gamma(E)$ for
$X\in TM$ and $\phi\in\Gamma(E)$
\end{center}

A connection must satisfy the Leibnitz rule: for f$\in$C$^{\infty}($M) and
$\phi\in\Gamma(E),$ we have:

\begin{center}
$\nabla(f\phi)=df\otimes\phi+f\nabla\phi.$
\end{center}

(ii) $\nabla$ is a \textbf{metric connection }if for all $\phi_{1},\phi_{2}%
\in\Gamma(E),$we have: $\qquad$

\begin{center}
$d<\phi_{1},\phi_{2}>$ $=$ $<\nabla\phi_{1},\phi_{2}>+<\phi_{1},\nabla\phi
_{2}>$
\end{center}

where $d$ is the exterior derivative and $<,>$ on the LHS of the last equality
above is the inner product on $\Gamma(E)$ and $<,>$ on the RHS is the the
pairing between $\Gamma(E)$ and $\Gamma(T^{\ast}M\otimes E).$

Equivalently (see $\left(  3.2.2\right)  $ of \textbf{Berline, Getzler and
Vergne} $\left[  1\right]  )$, for any vector field X$\in TM,$

\begin{center}
$X<\phi_{1},\phi_{2}>$ $=$ $<\nabla_{X}\phi_{1},\phi_{2}>+<\phi_{1},\nabla
_{X}\phi_{2}>$
\end{center}

where $<,>$ on the RHS of the last equality above is the inner product on the
sections $\Gamma(E).$ Such a connection is said to be \textbf{compatible} with
the \textbf{Riemannian metric}.

(iii) The notion of the connection $\nabla^{E}:\Gamma(E)\rightarrow$
$\Gamma(T^{\ast}M\otimes E$) can be \textbf{extended} to an operator
$\nabla^{k}$ as follows:

We denote by $\Lambda^{k}(M)$ $=\Lambda^{k}T^{\ast}M$ the \textbf{vector
bundle of k-form}s on M and $\Lambda^{k}(M)\otimes E$ the \textbf{vecctor
bundle of k-forms on M with values in the vector bundle E}. Then, for
$\alpha\in\Gamma(\Lambda^{k}M)$ and $\phi\in$ $\Gamma(E)$ and hence
$\alpha\otimes\phi\in\Gamma(\Lambda^{k}M\otimes E),$

$\nabla^{k}:\Gamma(\Lambda^{k}M\otimes E)\longrightarrow\Gamma(\Lambda
^{k+1}M\otimes E)$

is defined by:

$\nabla^{k}(\alpha\otimes\phi)=$ $d\alpha\otimes\phi$ $+$ $(-1)^{k}%
\alpha\wedge\nabla\phi$ for $k=0,1,...,n,...,$ where $\nabla^{0}=\nabla^{E}.$

(iv) We have the natural identification:

\begin{center}
\qquad$\Gamma(\Lambda^{k}M\otimes E)=\left\{  TM\times TM\times....\times
TM\longrightarrow\Gamma(E)\right\}  $
\end{center}

where the map is \textbf{alternating} and \textbf{multilinear} relative to
C$^{\infty}(M)-$ modules $TM\times TM\times....\times TM$ and $\Gamma(E).$

An arbitrary element $\Theta\in\Gamma(\Lambda^{k}M\otimes E)$ can be written
as $\Theta=\alpha\otimes\phi$ for $\alpha\in\Gamma(\Lambda^{k}M)$ and $\phi
\in\Gamma(E).$

(v) We note that the above definition can be further generalized: see
\textbf{Definition 1.14} on p. 21 of \textbf{Berline} et al:

We take $\alpha\in\Gamma(\Lambda^{k}M)$ and $\phi\in\Gamma(\Lambda^{k}M\otimes
E)=$ space of k-forms on M with values in E and have:

$\nabla^{k}(\theta\wedge\phi)=$ $d\theta\wedge\phi$ $+$ $(-1)^{k}\theta
\wedge\nabla\phi$ for $k=0,1,...,n,...$

\qquad\qquad\qquad\qquad\qquad\qquad\qquad\qquad\qquad\qquad\qquad\qquad
\qquad\qquad\qquad\qquad\qquad\qquad\qquad\qquad$\blacksquare$

We will follow \textbf{Definition} $\left(  2.4\right)  $ of \textbf{Berline,
Getzler and Vergne} $\left[  7\right]  $ or $\left(  1.2.2\right)  $ of
\textbf{Gilkey} $\left[  21\right]  $ for the definition of the
\textbf{Generalized Laplacian} or \textbf{Laplace-Type operator:}

Let $\partial_{i}\circeq\frac{\partial}{\partial x_{i}}$ be a coordinate frame
field on TM. Let $\nabla^{E}$ be the connection on E and $\nabla^{T^{\ast
}M\otimes E}$ the

connection on the vector bundle $T^{\ast}M\otimes E:$

$\left(  3.1\right)  \qquad\Gamma(E)\overset{\nabla^{E}}{\longrightarrow
}\Gamma(T^{\ast}M\otimes E)\overset{\nabla^{T^{\ast}M\otimes E}%
}{\longrightarrow}\Gamma(T^{\ast}M\otimes T^{\ast}M\otimes E)$

\begin{definition}
(Connection Laplacian; Generalized Laplacian)
\end{definition}

(i) The \textbf{connection} \textbf{Laplacian} $\Delta$ on $\Gamma(E)$ is the
second-order differential operator (with the opposite sign convention),
defined for $\phi\in\Gamma(E)$ and vector fields X,Y$\in TM$ \ by:

\begin{center}
$\Delta_{0}\phi=trace\nabla^{T^{\ast}M\otimes E}\nabla^{E}\phi$
\end{center}

(ii) The \textbf{Generalized Laplacian} is defined by:

\begin{center}
$\Delta\phi=\Delta_{0}\phi+$ $W\phi$
\end{center}

where $W\in\Gamma$(End(E)) is a Weitzenb\H{o}ck term.

\qquad\qquad\qquad\qquad\qquad\qquad\qquad\qquad\qquad\qquad\qquad\qquad
\qquad\qquad$\qquad\qquad\blacksquare$

Let $\phi\in\Gamma(E)$ and $X,Y\in TM.$ Then it is known that (see
\textbf{Definition 2.4} of \textbf{Berline, Getzler and Vergne} $\left[
7\right]  )$

\begin{center}
$(\nabla^{T^{\ast}M\otimes E}\nabla^{E}\phi)(X,Y)=(\nabla_{X}^{E}\nabla
_{Y}^{E}-\nabla_{\nabla_{X}Y}^{E})\phi$
\end{center}

\qquad\qquad\qquad\qquad\qquad\qquad\qquad\qquad\qquad\qquad\qquad\qquad
\qquad\qquad\qquad\qquad$\blacksquare$

\begin{proposition}
(The Expression of the Connection Laplacian and the Generalized Laplacian)
\end{proposition}

(i) $\Delta_{0}\phi=$ trace$\nabla^{\text{T}^{\ast}\text{M}\otimes\text{E}%
}\nabla^{\text{E}}\phi=g^{ij}(\nabla_{\partial_{i}}^{E}\nabla_{\partial_{j}%
}^{E}-\Gamma_{ij}^{k}\nabla_{\partial_{k}}^{E})\phi$

(ii) $\Delta\phi$ $=\Delta_{0}\phi+W\phi$\ $=$ $g^{ij}(\nabla_{\partial_{i}%
}^{E}\nabla_{\partial_{j}}^{E}-\Gamma_{ij}^{k}\nabla_{\partial_{k}}^{E})\phi+$
$W\phi$

where $W\in\Gamma(End(E))$ is a Weitzenb\H{o}ck term.

\begin{proof}
(i) Since,
\end{proof}

\begin{center}
$(\nabla^{T^{\ast}M\otimes E}\nabla^{E}\phi)(X,Y)=(\nabla_{X}^{E}\nabla
_{Y}^{E}-\nabla_{\nabla_{X}Y}^{E})\phi,$
\end{center}

we have by definition,

$\qquad\qquad\Delta_{0}\phi=$ trace$\nabla^{\text{T}^{\ast}\text{M}%
\otimes\text{E}}\nabla^{\text{E}}\phi=g^{ij}(\nabla_{\partial_{i}}^{E}%
\nabla_{\partial_{j}}^{E}-\nabla_{\nabla_{\partial_{i}}\partial_{j}}^{E})\phi$

where the basis $\left(  \partial_{1},...,\partial_{d}\right)  =\left(
\frac{\partial}{\partial x_{1}},...,\frac{\partial}{\partial x_{d}}\right)  $
is\textbf{ not} necessarily orthonormal.

We see from p. 11 of \textbf{Roe }$\left[  48\right]  $ that for $X,Y\in TM,$
$X=X^{i}\partial_{i}$ and $Y=Y^{j}\partial_{j},$ the Levi-Cevita connection gives:

\begin{center}
$\qquad\nabla_{X}Y=X^{i}[\partial_{i}Y^{j}+\underset{k=1}{\overset{d}{\sum}%
}\Gamma_{ik}^{j}Y^{k}]\partial_{j}\qquad$
\end{center}

Consequently for $X=\partial_{i}$ and $Y=\partial_{j},$ we have $X^{i}=1;$
$Y^{k}=1$ and so $\partial_{i}Y^{j}=0.$ The last formula above gives:

$\qquad\qquad\nabla_{\partial_{i}}\partial_{k}=$
$\underset{j=1}{\overset{d}{\sum}}\Gamma_{ik}^{j}\partial_{j};$ Equivalently,
$\nabla_{\partial_{i}}\partial_{j}=$ $\underset{k=1}{\overset{d}{\sum}}%
\Gamma_{ij}^{k}\partial_{k}$

where $\Gamma_{ij}^{k}$ are tha \textbf{Christoffel symbols }defined by
\textbf{Levi-Cevita connection}. Therefore,

$\left(  3.2\right)  \qquad\nabla_{\nabla_{\partial i}\partial_{j}}^{E}%
=\nabla_{\Gamma_{ij}^{k}\partial_{k}}^{E}=\Gamma_{ij}^{k}\nabla_{\partial_{k}%
}^{E}$

\qquad\qquad\qquad\qquad\qquad\qquad\qquad\qquad\qquad\qquad\qquad\qquad
\qquad\qquad\qquad\qquad\qquad

\begin{center}
$\Delta_{0}\phi=$ trace$\nabla^{\text{T}^{\ast}\text{M}\otimes\text{E}}%
\nabla^{\text{E}}\phi=g^{ij}(\nabla_{\partial_{i}}^{E}\nabla_{\partial_{j}%
}^{E}-\nabla_{\nabla_{\partial_{i}}\partial_{j}}^{E})\phi=g^{ij}%
(\nabla_{\partial_{i}}^{E}\nabla_{\partial_{j}}^{E}-\Gamma_{ij}^{k}%
\nabla_{\partial_{k}}^{E})\phi$
\end{center}

where $(\Gamma_{ij}^{k})i,j,k=1,...,d$ are the Christoffel symbols of the
Levi-Civita connection on $TM.$

(ii) The \textbf{Generalized Laplacian} (or \textbf{Laplace-Type} operator)
$\Delta$ is defined by adding a Weitzenb\"{o}ck term to the Connection Laplacian:

$\left(  3.3\right)  \qquad$\qquad$\Delta\phi$ $=$ $\Delta_{0}\phi+$ $W\phi
$\ $=$ $g^{ij}(\nabla_{\partial_{i}}^{E}\nabla_{\partial_{j}}^{E}-\Gamma
_{ij}^{k}\nabla_{\partial_{k}}^{E})\phi+$ $W\phi$

where $W\in\Gamma$(End(E)) is a Weitzenb\H{o}ckian. The Weitzenb\H{o}ckian is
a section of the vector bundle $End(E):$ $W_{\text{x}}\in End(E_{\text{x}})$
for each x$\in$M.$\qquad\ \ \ $

$\qquad\qquad\qquad\qquad\qquad\qquad\qquad\qquad\qquad\qquad\qquad
\qquad\qquad\qquad\qquad\qquad\qquad\blacksquare\ \ $

The definition of the \textbf{curvature tensor} $R^{E}$ of the connection
$\nabla^{E}$ is given in several book. See for example: p. 22 of
\textbf{Berline, Getzler and Vergne} $\left[  7\right]  ;$ \textbf{Definition
}$\left(  1.4\right)  $ of \textbf{Roe} $\left[  48\right]  ;$%
\textbf{Definition} $\left(  6.4\right)  $ on p. 54 of \textbf{Duitermaat}
$\left[  13\right]  ;$ \textbf{Theorem} $\left(  3.1.2\right)  $ of
\textbf{Jost} $\left[  33\right]  ;$ \textbf{Proposition} $\left(  4.6\right)
$ of \textbf{Lawson and Michelsohn} $\left[  35\right]  $ and $\left(
5.3\right)  $\textbf{ }of \textbf{Morita }$\left[  39\right]  .$

\qquad\qquad\qquad\qquad\qquad\qquad\qquad\qquad\qquad\qquad\qquad\qquad
\qquad\qquad\qquad\qquad\qquad$\blacksquare$

\begin{definition}
(Curvature)
\end{definition}

For a \textbf{connection} $\nabla^{E}:\Gamma(E)\rightarrow\Gamma(T^{\ast
}M\otimes E)$ and any two vector fields $X,Y\in$TM, the \textbf{curvature
tensor }R$^{E}$ of the \textbf{vector bundle} E is defined by:

$\left(  3.4\right)  \qquad R^{E}(X,Y)\phi$ $=\left(  \nabla_{X}^{E}\nabla
_{Y}^{E}-\nabla_{Y}^{E}\nabla_{X}^{E}-\nabla_{\left[  X,Y\right]  }%
^{E}\right)  \phi=\left[  \nabla_{X}^{E},\nabla_{Y}^{E}\right]  \phi
-\nabla_{\left[  X,Y\right]  }^{E}\phi$ for all $\phi\in\Gamma(E)$

\qquad\qquad\qquad\qquad\qquad\qquad\qquad\qquad\qquad\qquad\qquad\qquad
\qquad\qquad\qquad\qquad\qquad$\blacksquare$

The composition $\nabla^{1}\circ\nabla^{E}:\Gamma(E)\longrightarrow
\Gamma(\Lambda^{2}M\otimes E)$ is a smooth section of the vector bundle
$Hom(E,\Lambda^{2}M\otimes E)$ $\cong$ $\Lambda^{2}M\otimes End(E)$ and
\textbf{coincides} with the \textbf{curvature tensor} $R^{E}$ of the vector
bundle E defined above.

See \textbf{Definition }$\left(  3.1.5\right)  $ of \textbf{Jost }$\left[
33\right]  ,$ $\left(  4.20\right)  $ of \textbf{Proposition }$\left(
4.6\right)  $ in \textbf{Lawson and Michelsohn} $\left[  35\right]  ,$
Definition on p. 293 of \textbf{Milnor and Stasheff }or \textbf{Proposition}
$\left(  5.24\right)  $ of \textbf{Morita} $\left[  39\right]  .$

\qquad\qquad\qquad\qquad\qquad\qquad\qquad\qquad\qquad\qquad\qquad\qquad
\qquad\qquad\qquad\qquad\qquad$\ \ \blacksquare$

\section{Local Expression of the Curvature Tensor}

We have not yet formally defined the notion of a vector bundle. We need the
definition of a local frame field of a vector bundle but before we do that we
need to formally define the notion of a vector bundle.

A (real) \textbf{vector bundle} of \textbf{rank d} over a (compact)
\textbf{n-dimensional} C$^{\infty}-$ Riemannian manifold M is a triple (E,
$\pi,$ M) where E is called the \textbf{total space} and $\pi:E\longrightarrow
M$ is called the \textbf{projection} such that:

(i) $\pi^{-1}(x)=$ $E_{x}$ has the structure of an \textbf{d-dimensional}
vector space.

(ii) E has the property of \textbf{local triviality: }for each x$_{0}\in$M,
there exists an open \textbf{neighbourhood} U$\subset$M of x$_{0}$ and a
\textbf{diffeomorphism} $\varphi_{U}:$ $\pi^{-1}(U)\longrightarrow U\times
R^{d}$ such that for each x$\in U,$ the restriction $\varphi_{x}:$ $\pi
^{-1}(x)\longrightarrow\left\{  x\right\}  \times R^{d}$ is a linear
isomorphism. From the above local trivialization we can choose sections
$\mu_{i}:U\longrightarrow$ $E$ defined by:

$\qquad\mu_{i}(x)=\varphi_{U}^{-1}(x,e_{i})$ where $\left(  e_{1}%
,...,e_{d}\right)  $ is the standard basis of $R^{d}.$

Then, $\left(  \mu_{1}(x),...,\mu_{d}(x)\right)  $ is a basis of $E_{x}$ for
any point x$\in U$ (see, for \textbf{Example, 7.1} of \textbf{Darling
}$\left[  11\right]  $ or \textbf{Theorem 1.5.3 }of\textbf{ Jost} $\left[
33\right]  ).$ The set of sections $\left\{  \mu_{1},...,\mu_{d}\right\}  $ is
called a (local) \textbf{frame field }of the vector bundle E over the chart U
based at the point x$_{0}\in$M. Then any section $\phi\in\Gamma(E)$ can
locally be written as: $\phi=\phi^{j}\mu_{j}$ (summation over repeated indices
is understood), where $\phi^{j}:U\subset M\longrightarrow R$ are smooth functions.

\qquad\qquad\qquad\qquad\qquad\qquad\qquad\qquad\qquad\qquad\qquad\qquad
\qquad\qquad$\blacksquare$

We see that for each point x$_{0}\in M,$ there is an open neighbourhood U of
x$_{0}$ which serves both as a \textbf{manifold chart }of the Riemannian
manifold M and as a \textbf{vector bundle} chart for E.

From now hence we drop the superscript $E$ on $\nabla^{E}:\Gamma
(E)\rightarrow\Gamma(T^{\ast}M\otimes E)$ and have:

$\nabla:\Gamma(E)\rightarrow\Gamma(T^{\ast}M\otimes E)$

\qquad\qquad\qquad\qquad\qquad\qquad\qquad\qquad\qquad\qquad\qquad\qquad
\qquad\qquad\qquad\qquad\qquad\qquad\qquad$\blacksquare$\qquad\qquad
\qquad\qquad\qquad\qquad\qquad\qquad\qquad\qquad\qquad\qquad\qquad\qquad
\qquad\qquad\qquad\qquad\qquad\qquad\qquad\qquad\qquad\qquad\qquad\qquad\qquad

\begin{proposition}
\qquad\qquad\qquad\qquad\qquad\qquad\qquad\qquad\qquad\qquad\qquad\qquad
\qquad\qquad\qquad\qquad
\end{proposition}

Let $\mu_{1},...,\mu_{d}$ be an orthonormal local frame field for the vector
bundle E on a chart U$\subset$M based at y$_{0}\in$U and let $\left(
\text{x}_{1}\text{,...,x}_{q},\text{x}_{q+1},...,\text{x}_{n}\right)  $ be a
coordinate system centred at y$_{0}$. Then (locally), we have for
$i,j=1,...,q,q+1,...,n:$

(i)$\qquad\nabla=\left(  \frac{\partial}{\partial x_{k}}+\Lambda_{k}\right)
\otimes dx_{k}$ and $\nabla_{X}\phi=\left(  \frac{\partial\phi}{\partial
x_{i}}+\Lambda_{i}\phi\right)  X^{i}$

(ii)$\qquad\nabla_{\partial_{i}}=\frac{\partial}{\partial x_{i}}+\Lambda_{i}$

(iii)\qquad$\nabla_{\partial_{i}}\nabla_{\partial_{j}}=\frac{\partial^{2}%
}{\partial x_{i}\partial x_{j}}+\frac{\partial\Lambda_{j}}{\partial x_{i}%
}+\Lambda_{j}\frac{\partial}{\partial x_{i}}+\Lambda_{i}\frac{\partial
}{\partial x_{j}}+\Lambda_{i}\Lambda_{j}$

(iv)\qquad$R^{E}(\partial_{i},\partial_{j})=$\ $\nabla_{\partial_{i}}%
\nabla_{\partial_{j}}-\nabla_{\partial_{j}}\nabla_{\partial_{i}}%
=\frac{\partial\Lambda_{j}}{\partial x_{i}}-\frac{\partial\Lambda_{i}%
}{\partial x_{j}}+\Lambda_{i}\Lambda_{j}-\Lambda_{j}\Lambda_{i}$

$\qquad\qquad\qquad\qquad=\frac{\partial\Lambda_{j}}{\partial x_{i}}%
-\frac{\partial\Lambda_{i}}{\partial x_{j}}+\left[  \Lambda_{i},\Lambda
_{j}\right]  $

(v)\qquad$\ \Omega_{ij}=\ \Omega(\partial_{i},\partial_{j})=\left(
\frac{\partial\Lambda_{j}}{\partial x_{i}}-\frac{\partial\Lambda_{i}}{\partial
x_{j}}+\Lambda_{i}\Lambda_{j}-\Lambda_{j}\Lambda_{i}\right)  =R^{E}%
(\partial_{i},\partial_{j}):=R_{ij}^{E}$\ 

The rest of results below are given in \textbf{normal coordinates }$\left(
x_{1},...,x_{n}\right)  :$

(vi)$\qquad\ \Lambda_{j}(y_{0})=0$ (zero matrix) for $i=1,...,q,q+1,...,n:$

(vii)$\qquad\frac{\partial\Lambda_{j}}{\partial x_{i}}(y_{0})=\frac{1}%
{2}\Omega_{ij}(y_{0})$

(viii)$\qquad\ \frac{\partial^{2}\Lambda_{k}}{\partial\text{x}_{i}%
\partial\text{x}_{j}}(y_{0})=\frac{1}{6}\frac{\partial\Omega_{jk}}%
{\partial\text{x}_{i}}(y_{0});\ \frac{\partial^{2}\Lambda_{k}}{\partial
\text{x}_{i}^{2}}(y_{0})=\frac{1}{6}\frac{\partial\Omega_{ik}}{\partial
\text{x}_{i}}(y_{0})$

(ix)$\qquad\ \frac{\partial^{2}\Lambda_{j}^{2}}{\partial\text{x}_{i}%
\partial\text{x}_{k}}(y_{0})=\frac{1}{4}\left[  \Omega_{kj}\Omega_{ij}%
+\Omega_{ij}\Omega_{kj}\right]  (y_{0})+\frac{1}{3}\left[  \frac
{\partial\Omega_{ij}}{\partial\text{x}_{k}}\Lambda_{j}+\Lambda_{j}%
\frac{\partial\Omega_{kj}}{\partial\text{x}_{i}}\right]  (y_{0})$

\qquad\qquad\qquad\qquad$=\frac{1}{4}\left[  \Omega_{kj}\Omega_{ij}%
+\Omega_{ij}\Omega_{kj}\right]  (y_{0})$ 

Since $\Lambda_{j}(y_{0})=0$

$\qquad\ \frac{\partial^{2}\Lambda_{j}^{2}}{\partial\text{x}_{i}^{2}}%
(y_{0})=\frac{1}{2}\left[  \Omega_{ij}\Omega_{ij}\right]  (y_{0})+\frac{1}%
{3}\left[  \frac{\partial\Omega_{ij}}{\partial\text{x}_{i}}\Lambda_{j}%
+\Lambda_{j}\frac{\partial\Omega_{ij}}{\partial\text{x}_{i}}\right]
(y_{0})=\frac{1}{2}\left[  \Omega_{ij}\Omega_{ij}\right]  (y_{0})$

(x)$\qquad\ \frac{\partial^{3}\Lambda_{l}}{\partial\text{x}_{i}\partial
\text{x}_{j}\partial\text{x}_{k}}(y_{0})=\ \frac{1}{4}\frac{\partial^{2}%
\Omega_{kl}}{\partial\text{x}_{i}\partial\text{x}_{j}}(y_{0})\ $

(xi) Generalization of the Yang-Mills equation:$\qquad$

$\qquad\left[  \frac{\partial\Omega_{ij}}{\partial x_{k}}-\frac{\partial
\Omega_{ik}}{\partial x_{j}}+\frac{\partial\Omega_{jk}}{\partial x_{i}%
}\right]  =[\Omega_{ij},\Lambda_{k}]-[\Omega_{ik},\Lambda_{j}]+[\Omega
_{jk},\Lambda_{i}]=0$ in \textbf{normal coordinates}

\begin{proof}
All computations are carried out locally in a chart $\left(  U;x_{1}%
,...,x_{n}\right)  $ of the Riemannian manifold M.
\end{proof}

We first remark that (i) - (v) are true for \textbf{any coordinate system} and
that is why we use a general chart $\left(  U;x_{1},...,x_{n}\right)  $ of a
Riemannian manifold. In particular they are true for a \textbf{Fermi
coordinate system} $\left(  M_{0};x_{1},...,x_{q},x_{q+1},...,x_{n}\right)  .$

Then, (vi) - (xi) are true for the particular case of a \textbf{normal
coordinate system} only.

(i) By \textbf{Berline, Getzeler and Vergne} $\left[  7\right]  $ p. 22 or
$\left(  3.1.13\right)  $ of \textbf{Jost} $\left[  33\right]  $, p.104 the
connection $\nabla$ can locally (on a local chart) be decomposed as:

\begin{center}
$\nabla=$ $d$ $+$ $\Lambda$
\end{center}

where $\Lambda$ is the End(E)$-$valued $1-$form defined above.

The convention of summation over repeated indices is assumed below.

For $\phi=\phi^{j}\mu_{j}$ (where $\mu_{1},...,\mu_{d}$ is the local frame
field given above), and $\Lambda=\Lambda_{k}dx^{k}$ for $k=1,...,n,$ we have:

\begin{center}
$\nabla\phi=\nabla(\phi^{j}\mu_{j})=d(\phi^{j}\mu_{j})+\Lambda_{k}%
dx^{k}\otimes(\phi^{j}\mu_{j})$

$\ \ \ =d\phi^{j}\otimes\mu_{j}+\phi^{j}d\mu_{j}+\Lambda_{k}\phi^{j}%
dx^{k}\otimes\mu_{j}$
\end{center}

Since $d\mu_{j}=0$ by the fact that $\nabla$ is a metric connection (see, for
example, the proof of \textbf{Lemma }$\left(  3.2.2\right)  $ p. 111 of
\textbf{Jost} $\left[  33\right]  )$ and $d\phi^{j}=\frac{\partial\phi^{j}%
}{\partial x_{k}}dx^{k},$ we have:

$\nabla\phi=\frac{\partial\phi^{j}}{\partial x_{k}}dx^{k}\otimes\mu
_{j}+\Lambda_{k}\phi^{j}dx^{k}\otimes\mu_{j}=\left(  \frac{\partial\phi^{j}%
}{\partial x_{k}}+\Lambda_{k}\phi^{j}\right)  dx^{k}\otimes\mu_{j}$

$\qquad\ =\left(  \frac{\partial}{\partial x_{k}}+\Lambda_{k}\right)  \phi
^{j}\mu_{j}\otimes dx^{k}=\left(  \frac{\partial}{\partial x_{k}}+\Lambda
_{k}\right)  \phi\otimes dx^{k}$

We have:

$\qquad\nabla\phi=\left(  \frac{\partial}{\partial x_{k}}+\Lambda_{k}\right)
\phi\otimes dx^{k}=\left(  \frac{\partial\phi}{\partial x_{k}}+\Lambda_{k}%
\phi\right)  \otimes dx^{k}$

We conclude that locally,

$\qquad\qquad\nabla\phi=\left(  \frac{\partial\phi}{\partial x_{k}}%
+\Lambda_{k}\phi\right)  \otimes dx^{k}$

Equivalently,

$\left(  3.5\right)  \qquad\qquad\nabla=\left(  \frac{\partial}{\partial
x_{k}}+\Lambda_{k}\right)  \otimes dx^{k}\qquad$

Consequently for $X=X^{i}\frac{\partial}{\partial x_{i}},$ we have:

\begin{center}

\end{center}

$\qquad\nabla_{X}\phi:=\nabla\phi(X)=\left(  \frac{\partial\phi}{\partial
x_{k}}+\Lambda_{k}\phi\right)  dx^{k}(X)=\left(  \frac{\partial\phi}{\partial
x_{k}}+\Lambda_{k}\phi\right)  dx^{k}(X^{i}\frac{\partial}{\partial x_{i}})$

$\qquad=\left(  \frac{\partial\phi}{\partial x_{k}}+\Lambda_{k}\phi\right)
X^{i}dx^{k}(\frac{\partial}{\partial x_{i}})=\left(  \frac{\partial\phi
}{\partial x_{k}}+\Lambda_{k}\phi\right)  X^{i}\delta_{i}^{k}=\left(
\frac{\partial\phi}{\partial x_{i}}+\Lambda_{i}\phi\right)  X^{i}$

and so (i) is proved.

The last formula above is given on p.16 of \textbf{Gilkey }$\left[  22\right]
.$ We have proved it here in detail.

(ii)$\qquad\qquad\nabla\phi=\left(  \frac{\partial\phi}{\partial x_{k}%
}+\Lambda_{k}\phi\right)  \otimes dx_{k}$

By definition, $\nabla_{X}\phi=\nabla\phi(X)$ for a vector field $X\in
\Gamma(TM),$ and so,$\qquad\qquad$

$\qquad\nabla_{\partial_{i}}\phi=\nabla\phi(\partial_{i})=\left(
\frac{\partial\phi}{\partial x_{k}}+\Lambda_{k}\phi\right)  dx_{k}%
(\frac{\partial}{\partial x_{i}})=\left(  \frac{\partial\phi}{\partial x_{k}%
}+\Lambda_{k}\phi\right)  \delta_{i}^{k}$

$\qquad\qquad=\left(  \frac{\partial\phi}{\partial x_{i}}+\Lambda_{i}%
\phi\right)  =\left(  \frac{\partial}{\partial x_{i}}+\Lambda_{i}\right)
\phi$

We conclude that:

$\left(  3.6\right)  \qquad\qquad\qquad\qquad\nabla_{\partial_{i}}=\left(
\frac{\partial}{\partial x_{i}}+\Lambda_{i}\right)  $

and (ii) is proved.

This is proof of \textbf{Remark 1.2} of p.10 in\ \textbf{Roe} $\left[
48\right]  .$

(iii) Since,

$\qquad\nabla_{\partial_{j}}\phi=\left(  \frac{\partial}{\partial x_{j}%
}+\Lambda_{j}\right)  (\phi^{k}\mu_{k})=\left(  \frac{\partial\phi^{k}%
}{\partial x_{j}}+\Lambda_{j}\phi^{k}\right)  \mu_{k}=\left(  \frac
{\partial\phi^{k}}{\partial x_{j}}+\Lambda_{j}\phi^{k}\right)  \mu_{k}$

we have:

\qquad\qquad we have:

$\qquad\nabla_{\partial_{i}}\nabla_{\partial_{j}}\phi=\left(  \frac{\partial
}{\partial x_{i}}+\Lambda_{i}\right)  \left(  \frac{\partial\phi^{k}}{\partial
x_{j}}+\Lambda_{j}\phi^{k}\right)  \mu_{k}$

$\qquad\qquad\qquad=\left(  \frac{\partial^{2}\phi^{k}}{\partial x_{i}\partial
x_{j}}+\frac{\partial\Lambda_{j}}{\partial x_{i}}\phi^{k}+\Lambda_{j}%
\frac{\partial\phi^{k}}{\partial x_{i}}+\Lambda_{i}\frac{\partial\phi^{k}%
}{\partial x_{j}}+\Lambda_{i}\Lambda_{j}\phi^{k}\right)  \mu_{k}$

$\qquad\qquad\qquad=\left(  \frac{\partial^{2}}{\partial x_{i}\partial x_{j}%
}+\frac{\partial\Lambda_{j}}{\partial x_{i}}+\Lambda_{j}\frac{\partial
}{\partial x_{i}}+\Lambda_{i}\frac{\partial}{\partial x_{j}}+\Lambda
_{i}\Lambda_{j}\right)  \phi^{k}\mu_{k}$

$\qquad\qquad\qquad=\left(  \frac{\partial^{2}}{\partial x_{i}\partial x_{j}%
}+\frac{\partial\Lambda_{j}}{\partial x_{i}}+\Lambda_{j}\frac{\partial
}{\partial x_{i}}+\Lambda_{i}\frac{\partial}{\partial x_{j}}+\Lambda
_{i}\Lambda_{j}\right)  \phi$

Consequently for $i,j=1,...,n,$ we have:

$\left(  3.7\right)  \qquad\qquad\nabla_{\partial_{i}}\nabla_{\partial_{j}%
}\phi=\left(  \frac{\partial^{2}\phi}{\partial x_{i}\partial x_{j}}%
+\frac{\partial\Lambda_{j}}{\partial x_{i}}\phi+\Lambda_{j}\frac{\partial\phi
}{\partial x_{i}}+\Lambda_{i}\frac{\partial\phi}{\partial x_{j}}+\Lambda
_{i}\Lambda_{j}\phi\right)  $

(iv) Since $\frac{\partial^{2}\phi}{\partial x_{i}\partial x_{j}}%
=\frac{\partial^{2}\phi}{\partial x_{j}\partial x_{i}}$ for any smooth section
$\phi\in\Gamma(E),$ $\left(  3.7\right)  $ gives:

$\left(  3.8\right)  \qquad\nabla_{\partial_{i}}\nabla_{\partial_{j}}%
\phi-\nabla_{\partial_{j}}\nabla_{\partial_{i}}=\frac{\partial\Lambda_{j}%
}{\partial x_{i}}-\frac{\partial\Lambda_{i}}{\partial x_{j}}+\Lambda
_{i}\Lambda_{j}-\Lambda_{j}\Lambda_{i}=\frac{\partial\Lambda_{j}}{\partial
x_{i}}-\frac{\partial\Lambda_{i}}{\partial x_{j}}+\left[  \Lambda_{i}%
,\Lambda_{j}\right]  $

(v) Let $\Lambda$ be the End(E)$-$valued \textbf{connection }$1-$form and
$\Omega$ the End(E)-valued \textbf{curvature} $2-$form of the vector bundle E.
Then $\Lambda$ and $\Omega$ are related by the E. Cartan \textbf{second
structure equation} given in a local chart by:

$\left(  3.9\right)  \qquad\qquad\qquad\qquad\qquad\Omega=$\ $d\Lambda
+\Lambda\wedge\Lambda\qquad\qquad$

\qquad See, for example, \textbf{Theorem} $\left(  5.21\right)  $ of
\textbf{Morita} $\left[  39\right]  .$

Let $\left(  U;x_{1},...,x_{n}\right)  $ be a local chart of the Riemannian
manifold M.

In all that follow, we will write $\partial_{i}$ for $\frac{\partial}{\partial
x_{i}}$ and carry out \ all computations in \textbf{local coordinates}.

\qquad Since $\Lambda$ is an End(E)-valued 1-form and $\Omega$ is an
End(E)-valued 2-form, we follow \textbf{1.19} of \textbf{Roe }$\left[
48\right]  $ and set:

\qquad\qquad$\Lambda=$ $\overset{n}{\underset{k=1}{\sum}}\Lambda_{k}dx_{k}$
and so $d\Lambda=\overset{n}{\underset{l=1}{\sum}}d\Lambda_{l}dx_{l}%
=\overset{n}{\underset{k,l=1}{\sum}}\frac{\partial\Lambda_{l}}{\partial x_{k}%
}dx_{k}\Lambda dx_{l}$\-

\qquad On the other hand,

\qquad$\qquad\Omega=\frac{1}{2!}\overset{n}{\underset{k,l=1}{\sum}}\Omega
_{kl}dx_{k}\wedge dx_{l}$

From the \textbf{Structure Equation} in $\left(  3.9\right)  $ and the last
equations above, we have:

$\left(  3.10\right)  \qquad\Omega=\frac{1}{2!}%
\overset{n}{\underset{k,l=1}{\sum}}\Omega_{kl}dx_{k}\wedge dx_{l}%
=\overset{n}{\underset{k,l=1}{\sum}}[\frac{\partial\Lambda_{l}}{\partial
x_{k}}$ $+$ $\overset{n}{\underset{k,l=1}{\sum}}\Lambda_{k}\Lambda_{l}%
]dx_{k}\Lambda dx_{l}=d\Lambda+\Lambda\wedge\Lambda\qquad\qquad$

where $\Lambda_{k}(x)\in End(E_{x})$ and $\Omega_{kl}(x)\in End(E_{x}).$

\qquad By $\left(  7.1.2\right)  $ of \textbf{Hsu }$\left[  30\right]  ,$
$\Omega_{kl}$ are smooth functions \textbf{alternating} in the indices
$\left(  k,l\right)  .$

Here we will adopt the convention (see, for example, p. 14 of \textbf{Lee
}$\left[  36\right]  $) that for 1-forms $\omega_{1},...,\omega_{k},$ and
vectors $X_{1},...,X_{k},$ the wedge product $\omega_{1}\wedge,...,\wedge
\omega_{k}$ is defined by:

$\qquad\qquad\left(  \omega_{1}\wedge,...,\wedge\omega_{k}\right)  \left(
X_{1},...,X_{k}\right)  =\det\left(  \omega_{i}\left(  X_{j}\right)  \right)
$ $i,j=1,...,k.$

Then the components $\Lambda_{i}$ and $\Omega_{ij}$ are computed as follows:

$\left(  3.11\right)  $ $\qquad\Lambda(\partial_{i})=\Lambda(\frac{\partial
}{\partial x_{i}})=\overset{n}{\underset{k=1}{\sum}}\Lambda_{k}dx_{k}%
(\frac{\partial}{\partial x_{i}})=\overset{n}{\underset{k=1}{\sum}}\Lambda
_{k}\delta_{ik}=\Lambda_{i}$

Next we have:

$\qquad\Omega(\partial_{i},\partial_{j})=\frac{1}{2}%
\overset{n}{\underset{k,l=1}{\sum}}\Omega_{kl}(dx_{k}\wedge dx_{l}%
)(\frac{\partial}{\partial x_{i}},\frac{\partial}{\partial x_{j}})$

$\qquad=\frac{1}{2}$ $\overset{n}{\underset{k,l=1}{\sum}}\Omega_{kl}%
\det\left(
\begin{array}
[c]{cc}%
dx_{k}(\frac{\partial}{\partial x_{i}}) & dx_{k}(\frac{\partial}{\partial
x_{j}})\\
dx_{l}(\frac{\partial}{\partial x_{i}}) & dx_{l}(\frac{\partial}{\partial
x_{j}})
\end{array}
\right)  =\frac{1}{2}$ $\overset{n}{\underset{k,l=1}{\sum}}\Omega_{kl}%
\det\left(
\begin{array}
[c]{cc}%
\delta_{ik} & \text{ }\delta_{jk}\\
\delta_{il} & \delta_{jl}%
\end{array}
\right)  $

$\qquad=\frac{1}{2}$ $\overset{n}{\underset{k,l=1}{\sum}}\left(  \Omega
_{kl}\delta_{ik}\delta_{jl}-\Omega_{kl}\delta_{jk}\delta_{il}\right)
=\frac{1}{2}\left(  \Omega_{ij}-\Omega_{ji}\right)  $

Since $\Omega$ is a 2-form, the coefficients $\Omega_{ij}$ are smooth
functions \textbf{alternating} (skew-symmetric) in the indices $\left(
i,j\right)  .$ We see that $\frac{1}{2}\left(  \Omega_{ij}-\Omega_{ji}\right)
=\Omega_{ij}$ and so we have:

$\left(  3.12\right)  \qquad\Omega(\partial_{i},\partial_{j})=\Omega_{ij}$

From the \textbf{Second Structure Equation} in $\left(  3.9\right)  $, we have:

$\qquad\qquad\Omega(\partial_{i},\partial_{j})=$\ $d\Lambda(\partial
_{i},\partial_{j})+\Lambda\wedge\Lambda(\partial_{i},\partial_{j})$

$\qquad\qquad=\overset{n}{\underset{k,l=1}{\sum}}\frac{\partial\Lambda_{l}%
}{\partial x_{k}}dx_{k}\Lambda dx_{l}(\partial_{i},\partial_{j}%
)+\overset{n}{\underset{k,l=1}{\sum}}\Lambda_{k}\Lambda_{l}dx_{k}\Lambda
dx_{l}(\partial_{i},\partial_{j})$

We have:

$\overset{n}{\underset{k,l=1}{\sum}}\frac{\partial\Lambda_{l}}{\partial x_{k}%
}dx_{k}\Lambda dx_{l}(\partial_{i},\partial_{j})=$
$\overset{n}{\underset{k,l=1}{\sum}}\frac{\partial\Lambda_{l}}{\partial x_{k}%
}\det\left(
\begin{array}
[c]{cc}%
dx_{k}(\frac{\partial}{\partial x_{i}}) & dx_{k}(\frac{\partial}{\partial
x_{j}})\\
dx_{l}(\frac{\partial}{\partial x_{i}}) & dx_{l}(\frac{\partial}{\partial
x_{j}})
\end{array}
\right)  =$ $\overset{n}{\underset{k,l=1}{\sum}}\frac{\partial\Lambda_{l}%
}{\partial x_{k}}\det\left(
\begin{array}
[c]{cc}%
\delta_{ik} & \text{ }\delta_{jk}\\
\delta_{il} & \delta_{jl}%
\end{array}
\right)  $

$=$ $\overset{n}{\underset{k,l=1}{\sum}}\left(  \frac{\partial\Lambda_{l}%
}{\partial x_{k}}\delta_{ik}\delta_{jl}-\frac{\partial\Lambda_{l}}{\partial
x_{k}}\delta_{jk}\delta_{il}\right)  =\frac{\partial\Lambda_{j}}{\partial
x_{i}}-\frac{\partial\Lambda_{i}}{\partial x_{j}}$

Similarly,

$\overset{n}{\underset{k,l=1}{\sum}}\Lambda_{k}\Lambda_{l}dx_{k}\Lambda
dx_{l}(\partial_{i},\partial_{j})=$ $\overset{n}{\underset{k,l=1}{\sum}%
}\Lambda_{k}\Lambda_{l}\det\left(
\begin{array}
[c]{cc}%
\delta_{ik} & \text{ }\delta_{jk}\\
\delta_{il} & \delta_{jl}%
\end{array}
\right)  =\left(  \Lambda_{i}\Lambda_{j}-\Lambda_{j}\Lambda_{i}\right)  $

We conclude that the \textbf{Second Structure Equation} in $\left(
3.9\right)  $ gives:

$\qquad\qquad\Omega(\partial_{i},\partial_{j})=\left(  \frac{\partial
\Lambda_{j}}{\partial x_{i}}-\frac{\partial\Lambda_{i}}{\partial x_{j}%
}+\Lambda_{i}\Lambda_{j}-\Lambda_{j}\Lambda_{i}\right)  $

We see from the equations in $\left(  3.12\right)  $ and $\left(  3.13\right)
$ that:

$\left(  3.13\right)  \qquad\Omega_{ij}=\Omega(\partial_{i},\partial
_{j})=\left(  \frac{\partial\Lambda_{j}}{\partial x_{i}}-\frac{\partial
\Lambda_{i}}{\partial x_{j}}+\Lambda_{i}\Lambda_{j}-\Lambda_{j}\Lambda
_{i}\right)  $

We then use the equalites in $\left(  3.13\right)  ,$ the definition of
$R^{E}$ in $\left(  3.4\right)  $ and the fact that $\left[  \partial
_{i},\partial_{j}\right]  =0,$ to have:

$\left(  3.14\right)  \qquad\Omega_{ij}=\Omega(\partial_{i},\partial
_{j})=\left(  \frac{\partial\Lambda_{j}}{\partial x_{i}}-\frac{\partial
\Lambda_{i}}{\partial x_{j}}+\Lambda_{i}\Lambda_{j}-\Lambda_{j}\Lambda
_{i}\right)  =R^{E}(\partial_{i},\partial_{j}):=R_{ij}^{E}$

The equalites in the last expression above show that the components
$\Omega_{ij}$ of the \textbf{curvature 2-form} $\Omega$ and the components
$R_{ij}^{E}$ of the \textbf{curvature tensor }R$^{E}$ coincide
\textbf{locally} and have for common value equal to:

\qquad\qquad\qquad\ $\left(  \frac{\partial\Lambda_{j}}{\partial x_{i}}%
-\frac{\partial\Lambda_{i}}{\partial x_{j}}+\Lambda_{i}\Lambda_{j}-\Lambda
_{j}\Lambda_{i}\right)  $

To prove (vi) - (xi) here we will use (iii) \textbf{Proposition 13} in
\textbf{Chapter 9} below. We re-state it in full here:

In \textbf{normal coordinates, }we have for $i,j,k,l,p,r=1,...,q,q+1,...,n,$
$x\in M_{0}:$

$\left(  3.15\right)  \ \qquad\Lambda_{l}(x)=\frac{1}{2}\Omega_{il}%
(y_{0})x_{i}+\frac{1}{6}\frac{\partial\Omega_{jl}}{\partial\text{x}_{i}}%
(y_{0})(y_{0})x_{i}x_{j}-\frac{1}{36}[\overset{n}{\underset{r=1}{\sum}}\left(
R_{ijpr}+R_{ipjr}\right)  \Omega_{kr}](y_{0})x_{i}x_{j}x$

$\qquad\qquad+\frac{1}{24}[\frac{\partial^{2}\Omega_{kl}}{\partial
x_{i}\partial x_{j}}+\Omega_{ik}(y_{0})\Omega_{jl}(y_{0})+\Omega_{ij}%
(y_{0})\Omega_{kl}(y_{0})]x_{i}x_{j}x_{k}+$ higher order terms.

Then we have:

(vi) In \textbf{Fermi coordinates}, $\Lambda_{j}(y_{0})=\left(  \Gamma
_{ij}^{k}(y_{0})\right)  $ is \textbf{not} necessarily equal to zero (matirx)
for $j=1,...,q,q+1,...,n.$

\qquad However, when Fermi coordinates reduce to \textbf{normal coordinates,
}we have from the expansion in Proposition 13 here given above:

$\left(  3.16\right)  \qquad\Lambda_{j}(y_{0})=\left(  \Gamma_{ij}^{k}%
(y_{0})\right)  =0$ (zero matrix) for all $i,j,k=1,...,n.$

The above result is in \textbf{Proposition 1.18} of \textbf{Berline, Getzler
and Vergne }$\left[  7\right]  $ or from \textbf{Question 2.33 }of Roe
$\left[  48\right]  .$

(vii) It is immediate from the expansion of $\ \Lambda_{l}(x)$ given above that:

$\left(  3.17\right)  \qquad\frac{\partial\Lambda_{j}}{\partial x_{i}}%
(y_{0})=\frac{1}{2}\Omega_{ij}(y_{0})$

(viii) From the expansion we have:

$\left(  3.18\right)  $\qquad$\frac{\partial^{2}\Lambda_{k}}{\partial
\text{x}_{i}\partial\text{x}_{j}}(y_{0})=\frac{1}{6}\frac{\partial\Omega_{jk}%
}{\partial\text{x}_{i}}(y_{0})$

In particular, when $j=i,$ we have:

$\left(  3.19\right)  \qquad\ \frac{\partial^{2}\Lambda_{k}}{\partial
\text{x}_{i}^{2}}(y_{0})=\frac{1}{6}\frac{\partial\Omega_{ik}}{\partial
\text{x}_{i}}(y_{0})$

(ix) We have:$\qquad$

$\qquad\qquad\frac{\partial^{2}}{\partial\text{x}_{i}\partial\text{x}_{j}%
}(\Lambda_{k}^{2})=\frac{\partial^{2}}{\partial\text{x}_{i}\partial
\text{x}_{j}}(\Lambda_{k}\Lambda_{k})=\frac{\partial}{\partial\text{x}_{i}%
}[\frac{\partial}{\partial\text{x}_{j}}(\Lambda_{k}\Lambda_{k})]$

$\qquad=\frac{\partial}{\partial\text{x}_{i}}[\frac{\partial\Lambda_{k}%
}{\partial\text{x}_{j}}\Lambda_{k}+\Lambda_{k}\frac{\partial\Lambda_{k}%
}{\partial\text{x}_{j}}]=\left[  \frac{\partial^{2}\Lambda_{k}}{\partial
\text{x}_{i}\partial\text{x}_{j}}\Lambda_{k}+\frac{\partial\Lambda_{k}%
}{\partial\text{x}_{j}}\frac{\partial\Lambda_{k}}{\partial\text{x}_{i}}%
+\frac{\partial\Lambda_{k}}{\partial\text{x}_{i}}\frac{\partial\Lambda_{k}%
}{\partial\text{x}_{j}}+\Lambda_{k}\frac{\partial^{2}\Lambda_{k}}%
{\partial\text{x}_{i}\partial\text{x}_{j}}\right]  $

By $\left(  3.17\right)  $ and $\left(  3.18\right)  $ above we have:

$\left(  3.20\right)  $ $\qquad\ \frac{\partial^{2}\Lambda_{k}^{2}}%
{\partial\text{x}_{i}\partial\text{x}_{j}}(y_{0})=\frac{1}{6}\frac
{\partial\Omega_{jk}}{\partial\text{x}_{i}}(y_{0})\Lambda_{k}(y_{0})+\frac
{1}{4}\Omega_{jk}(y_{0})\Omega_{ik}(y_{0})$

$\qquad\qquad\qquad\qquad\qquad+\frac{1}{4}\Omega_{ik}(y_{0})\Omega_{jk}%
(y_{0})+\frac{1}{6}\Lambda_{k}(y_{0})\frac{\partial\Omega_{jk}}{\partial
\text{x}_{i}}(y_{0})$

$\left(  3.21\right)  $ $\qquad\ \frac{\partial^{2}\Lambda_{k}^{2}}%
{\partial\text{x}_{i}\partial\text{x}_{j}}(y_{0})=+\frac{1}{4}[\Omega
_{jk}\Omega_{ik}+\Omega_{ik}\Omega_{jk}](y_{0})+\frac{1}{6}\left[
\frac{\partial\Omega_{jk}}{\partial\text{x}_{i}}\Lambda_{k}+\Lambda_{k}%
\frac{\partial\Omega_{jk}}{\partial\text{x}_{i}}\right]  (y_{0})$

Since we are working in normal coordinates, $\Lambda_{k}(y_{0})=0$ for
$=1,...,n$ and so, we have:

$\left(  3.22\right)  $ \qquad$\ \frac{\partial^{2}\Lambda_{k}^{2}}%
{\partial\text{x}_{i}\partial\text{x}_{j}}(y_{0})=\frac{1}{4}[\Omega
_{jk}\Omega_{ik}+\Omega_{ik}\Omega_{jk}](y_{0})$

In particular, for $k=i,$ we have:\qquad

$\left(  3.23\right)  \qquad\ \frac{\partial^{2}\Lambda_{j}^{2}}%
{\partial\text{x}_{i}^{2}}(y_{0})=\frac{1}{2}\left[  \Omega_{ij}\Omega
_{ij}\right]  (y_{0})$\qquad

(x) We lastly consider: $(\partial^{\alpha}\Lambda_{l})(y_{0})=\ \frac
{\partial^{3}\Lambda_{l}}{\partial\text{x}_{i}\partial\text{x}_{j}%
\partial\text{x}_{k}}(y_{0}):$

\qquad We use the formula in the proof of \textbf{Proposition 1.18} of
\textbf{Berline, Getzler and Vergne }$\left[  7\right]  :$

\qquad We take $\alpha=(1,2)$ on the LHS and hence $\alpha=\left(  0,2\right)
$ with on the RHS

Therefore, on the LHS we have: $\alpha!=1!2!=2$ and $\left\vert \alpha
\right\vert =3$

On the RHS we have: $\alpha!=0!2!$ and $\left\vert \alpha\right\vert =2:$

Consequently,

\qquad$\qquad\qquad\qquad\frac{1}{1!2!}\left(  3+1\right)  \ \frac
{\partial^{3}\Lambda_{l}}{\partial\text{x}_{i}\partial\text{x}_{j}%
\partial\text{x}_{k}}(y_{0})=\frac{1}{0!2!}\ \frac{\partial^{2}\Omega_{kl}%
}{\partial\text{x}_{i}\partial\text{x}_{j}}(y_{0})$

Consequently,

$\left(  3.24\right)  $\qquad\qquad\qquad\qquad$\qquad\ \frac{\partial
^{3}\Lambda_{l}}{\partial\text{x}_{i}\partial\text{x}_{j}\partial\text{x}_{k}%
}(y_{0})=\ \frac{1}{4}\frac{\partial^{2}\Omega_{kl}}{\partial\text{x}%
_{i}\partial\text{x}_{j}}(y_{0})\ $

(xi) We use the \textbf{Second Structure Equation} in $\left(  3.9\right)
:\qquad\qquad\qquad\qquad\qquad$

$\qquad\qquad\Omega=$\ $d\Lambda+\Lambda\wedge\Lambda$

Therefore,

$\qquad\qquad d\Omega=$\ $d(d\Lambda)+d(\Lambda\wedge\Lambda)=$\ $d^{2}%
\Lambda+d\Lambda\wedge\Lambda-\Lambda\wedge d\Lambda=d\Lambda\wedge
\Lambda-\Lambda\wedge d\Lambda$

\qquad$\qquad\ \ \ \ =d\Lambda\wedge\Lambda+\Lambda\wedge\Lambda\wedge
\Lambda-\Lambda\wedge\Lambda\wedge\Lambda-\Lambda\wedge d\Lambda$

$\qquad\ \ \ \ \ \ \ \ \ \ \ =(d\Lambda+\Lambda\wedge\Lambda)\wedge
\Lambda-\Lambda\wedge(\Lambda\wedge\Lambda+d\Lambda)=\Omega\wedge
\Lambda-\Lambda\wedge\Omega$

We thus have:

$\qquad\qquad d\Omega=\Omega\wedge\Lambda-\Lambda\wedge\Omega=\left[
\Omega,\Lambda\right]  $

This is the \textbf{Second Bianchi Identity} given, for example, in $\left(
5.11\right)  ,$ p. 196 of \textbf{Morita }$\left[  39\right]  $ and
\textbf{Theorem 3.1.1 }of\textbf{ Jost }$\left[  33\right]  .$

Now let $\left(  \text{x}_{1}\text{,...,x}_{q},\text{x}_{q+1},...,\text{x}%
_{n}\right)  $ be Fermi coordinates centred at y$_{0}\in$U. We set:

$\qquad\qquad\Lambda=\overset{n}{\underset{r=1}{\sum}}\Lambda_{r}dx_{r}$ and
$\Omega=\frac{1}{2}\overset{n}{\underset{p,q=1}{\sum}}\Omega_{pq}dx_{p}\wedge
dx_{q}$

$\left(  3.25\right)  $ $\ \qquad d\Omega=\frac{1}{2}%
\overset{n}{\underset{p,q=1}{\sum}}d\Omega_{pq}dx_{p}\wedge dx_{q}=\frac{1}%
{2}$ $\overset{n}{\underset{p,q,r=1}{\sum}}\frac{\partial\Omega_{pq}}{\partial
x_{r}}dx_{r}\wedge dx_{p}\wedge dx_{q}$

$\qquad\qquad\qquad\ \ \ =\frac{1}{2}\overset{n}{\underset{p,q,r=1}{\sum}%
}\frac{\partial\Omega_{pq}}{\partial x_{r}}dx_{p}\wedge dx_{q}\wedge
dx_{r}\qquad\qquad$

$\left(  3.26\right)  \qquad\qquad\Omega\wedge\Lambda-\Lambda\wedge
\Omega=\frac{1}{2}\overset{n}{\underset{p,q,r=1}{\sum}}\Omega_{pq}\Lambda
_{r}dx_{p}\wedge dx_{q}\wedge dx_{r}$

$\qquad\qquad-\frac{1}{2}\overset{n}{\underset{p,q,r=1}{\sum}}\Lambda
_{r}\Omega_{pq}dx_{r}\wedge dx_{p}\wedge dx_{q}$

$\qquad\ =\frac{1}{2}\overset{n}{\underset{p,q,r=1}{\sum}}\Omega_{pq}%
\Lambda_{r}dx_{p}\wedge dx_{q}\wedge dx_{r}-\frac{1}{2}%
\overset{n}{\underset{p,q,r=1}{\sum}}\Lambda_{r}\Omega_{pq}dx_{p}\wedge
dx_{q}\wedge dx_{r}\qquad\qquad\qquad\qquad\qquad\ \ \ \qquad\qquad
\qquad\qquad\qquad$

Therefore from the \textbf{Second Bianchi Identity} we equate the final
expressions in $\left(  3.25\right)  $ and $\left(  3.26\right)  :$

$\left(  3.27\right)  \qquad\overset{n}{\underset{p,q,r=1}{\sum}}%
\frac{\partial\Omega_{pq}}{\partial x_{r}}dx_{p}\wedge dx_{q}\wedge
dx_{r}=\overset{n}{\underset{p,q,r=1}{\sum}}\Omega_{pq}\Lambda_{r}dx_{p}\wedge
dx_{q}\wedge dx_{r}$

$\qquad\qquad\qquad-\overset{n}{\underset{p,q,r=1}{\sum}}\Lambda_{r}%
\Omega_{pq}dx_{p}\wedge dx_{q}\wedge dx_{r}$

Therefore,

$\qquad\qquad\overset{n}{\underset{p,q,r=1}{\sum}}\left(  \frac{\partial
\Omega_{pq}}{\partial x_{r}}-\left(  \Omega_{pq}\Lambda_{r}-\Lambda_{r}%
\Omega_{pq}\right)  \right)  dx_{p}\wedge dx_{q}\wedge dx_{r}=0$

Consequently we have,

$\qquad\overset{n}{\underset{p,q,r=1}{\sum}}\left[  \frac{\partial\Omega_{pq}%
}{\partial x_{r}}-\left(  \Omega_{pq}\Lambda_{r}-\Lambda_{r}\Omega
_{pq}\right)  \right]  dx_{p}\wedge dx_{q}\wedge dx_{r}\left(  \frac{\partial
}{\partial x_{i}},\frac{\partial}{\partial x_{j}},\frac{\partial}{\partial
x_{k}}\right)  =0$

We have:

$\qquad=\overset{n}{\underset{p,q,r=1}{\sum}}\left[  \frac{\partial\Omega
_{pq}}{\partial x_{r}}-\left(  \Omega_{pq}\Lambda_{r}-\Lambda_{r}\Omega
_{pq}\right)  \right]  \det%
\begin{array}
[c]{ccc}%
dx_{p}\left(  \frac{\partial}{\partial x_{i}}\right)  & dx_{p}\left(
\frac{\partial}{\partial x_{j}}\right)  & dx_{p}\left(  \frac{\partial
}{\partial x_{k}}\right) \\
dx_{q}\left(  \frac{\partial}{\partial x_{i}}\right)  & dx_{q}\left(
\frac{\partial}{\partial x_{j}}\right)  & dx_{q}\left(  \frac{\partial
}{\partial x_{k}}\right) \\
dx_{r}\left(  \frac{\partial}{\partial x_{i}}\right)  & dx_{r}\left(
\frac{\partial}{\partial x_{j}}\right)  & dx_{r}\left(  \frac{\partial
}{\partial x_{k}}\right)
\end{array}
=0$

Computing the determinant, we have:

$\overset{n}{\underset{p,q,r=1}{\sum}}\left[  \frac{\partial\Omega_{pq}%
}{\partial x_{r}}-\left(  \Omega_{pq}\Lambda_{r}-\Lambda_{r}\Omega
_{pq}\right)  \right]  \left[  \delta_{ip}(\delta_{jq}\delta_{kr}-\delta
_{jr}\delta_{kq})-\delta_{iq}(\delta_{jp}\delta_{kr}-\delta_{jr}\delta
_{kp})+\delta_{ir}(\delta_{jp}\delta_{kq}-\delta_{jq}\delta_{kp})\right]  =0$

We simplify and have for all $i,j=1,...q,q+1,...,n:$

$\qquad\ \ \ \ \qquad\qquad\qquad\qquad\left[  \frac{\partial\Omega_{ij}%
}{\partial x_{k}}-\left(  \Omega_{ij}\Lambda_{k}-\Lambda_{k}\Omega
_{ij}\right)  \right]  -\left[  \frac{\partial\Omega_{ik}}{\partial x_{j}%
}-\left(  \Omega_{ik}\Lambda_{j}-\Lambda_{j}\Omega_{ik}\right)  \right]  $

$\qquad\qquad\qquad\qquad\qquad-\left[  \frac{\partial\Omega_{ji}}{\partial
x_{k}}-\left(  \Omega_{ji}\Lambda_{k}-\Lambda_{k}\Omega_{ji}\right)  \right]
+\left[  \frac{\partial\Omega_{ki}}{\partial x_{j}}-\left(  \Omega_{ki}%
\Lambda_{j}-\Lambda_{j}\Omega_{ki}\right)  \right]  $

$\qquad\qquad\qquad\qquad\qquad+\left[  \frac{\partial\Omega_{jk}}{\partial
x_{i}}-\left(  \Omega_{jk}\Lambda_{i}-\Lambda_{i}\Omega_{jk}\right)  \right]
-\left[  \frac{\partial\Omega_{kj}}{\partial x_{i}}-\left(  \Omega_{kj}%
\Lambda_{i}-\Lambda_{i}\Omega_{kj}\right)  \right]  =0$

Due to the skew-symmetry of $\Omega_{ij}$ in the indices $\left(  i,j\right)
,$ the equation in $\left(  3.19\right)  $ becomes:

$\qquad\qquad\qquad\qquad\qquad\ \ \ \left[  \frac{\partial\Omega_{ij}%
}{\partial x_{k}}-\Omega_{ij}\Lambda_{k}+\Lambda_{k}\Omega_{ij}\right]
-\left[  \frac{\partial\Omega_{ik}}{\partial x_{j}}-\Omega_{ik}\Lambda
_{j}+\Lambda_{j}\Omega_{ik}\right]  $

$\qquad\qquad\qquad\qquad\qquad-\left[  -\frac{\partial\Omega_{ij}}{\partial
x_{k}}+\Omega_{ij}\Lambda_{k}-\Lambda_{k}\Omega_{ij}\right]  +\left[
-\frac{\partial\Omega_{ik}}{\partial x_{j}}+\Omega_{ik}\Lambda_{j}-\Lambda
_{j}\Omega_{ik}\right]  $

$\qquad\qquad\qquad\qquad\qquad+\left[  \frac{\partial\Omega_{jk}}{\partial
x_{i}}-\Omega_{jk}\Lambda_{i}+\Lambda_{i}\Omega_{jk}\right]  -\left[
-\frac{\partial\Omega_{jk}}{\partial x_{i}}+\Omega_{jk}\Lambda_{i}-\Lambda
_{i}\Omega_{jk}\right]  =0$

We simplify and have

$\ \qquad2\left(  \frac{\partial\Omega_{ij}}{\partial x_{k}}-\Omega
_{ij}\Lambda_{k}+\Lambda_{k}\Omega_{ij}\right)  -2\left(  \frac{\partial
\Omega_{ik}}{\partial x_{j}}-\Omega_{ik}\Lambda_{j}+\Lambda_{j}\Omega
_{ik}\right)  +2\left(  \frac{\partial\Omega_{jk}}{\partial x_{i}}-\Omega
_{jk}\Lambda_{i}+\Lambda_{i}\Omega_{jk}\right)  =0$

We have:

$\qquad\left(  \ \frac{\partial\Omega_{ij}}{\partial x_{k}}-\frac
{\partial\Omega_{ik}}{\partial x_{j}}+\frac{\partial\Omega_{jk}}{\partial
x_{i}}\right)  -\Omega_{ij}\Lambda_{k}+\Lambda_{k}\Omega_{ij}+\Omega
_{ik}\Lambda_{j}-\Lambda_{j}\Omega_{ik}-\Omega_{jk}\Lambda_{i}+\Lambda
_{i}\Omega_{jk}=0$

In a shorter notation, we have for all $i,j,k=1,...,q,q+1,...,n:$

$\qquad\left(  \frac{\partial\Omega_{ij}}{\partial x_{k}}-\frac{\partial
\Omega_{ik}}{\partial x_{j}}+\frac{\partial\Omega_{jk}}{\partial x_{i}%
}\right)  +[\Lambda_{i},\Omega_{jk}]-[\Lambda_{j},\Omega_{ik}]+[\Lambda
_{k},\Omega_{ij}]=0$Equivalently,

$\left(  3.28\right)  \qquad\left(  \frac{\partial\Omega_{ij}}{\partial x_{k}%
}-\frac{\partial\Omega_{ik}}{\partial x_{j}}+\frac{\partial\Omega_{jk}%
}{\partial x_{i}}\right)  =[\Omega_{ij},\Lambda_{k}]-[\Omega_{ik},\Lambda
_{j}]+[\Omega_{jk},\Lambda_{i}].$

The last equation above will reduce to the \textbf{Yang-Mill }equation in
$\left(  3.2.16\right)  $ of \textbf{Jost }$\left[  33\right]  $\textbf{ }if
there were some form of linear independence in our equation. However, our
equation can still be regarded as some form of \textbf{generalization of the
Yang-Mills equation}.

\qquad\qquad\qquad\qquad\qquad\qquad\qquad\qquad\qquad\qquad\qquad\qquad
\qquad\qquad\qquad\qquad\qquad\qquad$\blacksquare$\qquad\qquad\qquad
\qquad\qquad\qquad\qquad\qquad\qquad\qquad\qquad\qquad\qquad\qquad\qquad
\qquad\qquad\qquad\qquad\qquad\qquad

\begin{proposition}
Let $\Delta^{0}$ be usual (scalar) Laplacian on C$^{\infty}$(M) and  $\Delta$
the Laplace-Type operator$.$Then we have the following local expressions and
let $\left(  \text{x}_{1}\text{,...,x}_{q},\text{x}_{q+1},...,\text{x}%
_{n}\right)  $ be Fermi coordinates centred at y$_{0}$. Let Let $\mu
_{1},...,\mu_{d}$ be an orthonormal local frame field for the vector bundle E
on a chart U$\subset$M based at y$_{0}\in$U. Then we have:
\end{proposition}

(i)\qquad$\nabla_{X}\phi$ $=$ $X^{j}\left(  \frac{\partial}{\partial x_{j}%
}+\Lambda_{j}\right)  \phi=X^{j}\nabla_{\partial_{j}}\phi$

(ii)\qquad$\nabla_{X}^{0}f$ \ $=$ $<\nabla^{0}f,X>$ $=\frac{\partial
f}{\partial x_{j}}X^{j}$

(iii)\qquad\ $<\nabla\phi,X>$ $=X^{j}\left(  \frac{\partial\phi}{\partial
x_{j}}+\Lambda_{j}\phi\right)  =X^{j}\nabla_{\partial_{j}}\phi=X^{j}\left(
\frac{\partial}{\partial x_{j}}+\Lambda_{j}\right)  \phi$

(iv) $\qquad<\nabla\phi,\nabla^{0}f>$ $=g^{ij}\frac{\partial f}{\partial
x_{i}}\nabla_{\partial_{j}}\phi$ $=$ $g^{ij}\frac{\partial f}{\partial x_{i}%
}\left(  \frac{\partial}{\partial x_{j}}+\Lambda_{j}\right)  \phi$

(v) $\qquad\Delta=g^{ij}\left\{  \frac{\partial^{2}}{\partial x_{i}%
\partial_{j}}+\frac{\partial\Lambda_{j}}{\partial x_{i}}+\Lambda_{j}%
\frac{\partial}{\partial x_{i}}+\Lambda_{i}\frac{\partial}{\partial x_{j}%
}+\Lambda_{i}\Lambda_{j}-\Gamma_{ij}^{k}\left(  \frac{\partial}{\partial
x_{k}}+\Lambda_{k}\right)  \right\}  $

$\qquad\qquad\qquad+W$

\qquad$\qquad\Delta^{0}=g^{ij}\left\{  \frac{\partial^{2}}{\partial
x_{i}\partial_{j}}-\Gamma_{ij}^{k}\frac{\partial}{\partial x_{k}}\right\}  $

(vi) $\ \ L(f\phi)=\frac{1}{2}g^{ij}f\left\{  \frac{\partial^{2}\phi}{\partial
x_{i}\partial_{j}}+\frac{\partial\Lambda_{j}}{\partial x_{i}}\phi+\Lambda
_{j}\frac{\partial\phi}{\partial x_{i}}+\Lambda_{i}\frac{\partial\phi
}{\partial x_{j}}+\Lambda_{i}\Lambda_{j}\phi-\Gamma_{ij}^{k}\left(
\frac{\partial\phi}{\partial x_{k}}+\Lambda_{k}\phi\right)  \right\}  $

$\qquad\qquad+fW\phi+\frac{1}{2}\phi g^{ij}\left\{  \frac{\partial^{2}%
f}{\partial x_{i}\partial_{j}}-\Gamma_{ij}^{k}\left(  \frac{\partial
f}{\partial x_{k}}+\Lambda_{k}f\right)  \right\}  $ $+$ $g^{ij}\frac{\partial
f}{\partial x_{i}}\left(  \frac{\partial\phi}{\partial x_{i}}+\Lambda_{i}%
\phi\right)  \qquad\qquad\qquad\qquad\qquad$

\qquad\qquad\ $+$ $f$ $X^{i}\left(  \frac{\partial\phi}{\partial x_{i}%
}+\Lambda_{i}\phi\right)  +\phi X^{i}\frac{\partial f}{\partial x_{i}}$ $+$
$V(f\phi)$

for L =$\frac{1}{2}\Delta+\nabla_{X}+V$ and L$^{0}=\frac{1}{2}\Delta
^{0}+\nabla_{X}^{0}+V$, where f$\in$C$^{\infty}$(M) and $\phi\in\Gamma(E)$.

\begin{proof}
(i) Let $X$ be a vector field on M.
\end{proof}

Let $X=X^{i}\frac{\partial}{\partial x_{i}}$ and $\phi=\phi^{j}\mu_{j}$

Then,

$\qquad\qquad\nabla_{X}\phi$ $=$ $\nabla_{X^{i}\partial_{i}}(\phi^{j}\mu
_{j})=$ $X^{i}[\nabla_{\partial_{i}}(\phi^{j}\mu_{j})]=$ $X^{i}[\phi
^{j}(\nabla_{\partial_{i}}\mu_{j})+(\partial_{i}\phi^{j})\mu_{j}]$

$\qquad\qquad=$ $X^{i}[\phi^{j}\left(  \frac{\partial}{\partial x_{i}}%
+\Lambda_{i}\right)  \mu_{j}+\frac{\partial\phi^{j}}{\partial x_{i}}\mu_{j}]=$
$X^{i}[\phi^{j}\left(  \frac{\partial\mu_{j}}{\partial x_{i}}+\Lambda_{i}%
\mu_{j}\right)  +\frac{\partial\phi^{j}}{\partial x_{i}}\mu_{j}]$

The connection is metric and hence, $\frac{\partial\mu_{j}}{\partial x_{i}%
}dx_{i}=d\mu_{j}=0$ (see \textbf{Morita }$\left[  39\right]  ,$ p. 111) and so,

$\qquad\qquad\frac{\partial\mu_{j}}{\partial x_{i}}=0.$

We conclude that:

$\qquad\nabla_{X}\phi=X^{i}[\phi^{j}\Lambda_{i}\mu_{j}+\frac{\partial\phi^{j}%
}{\partial x_{i}}\mu_{j}]=$ $X^{i}[\phi^{j}\Lambda_{i}+\frac{\partial\phi^{j}%
}{\partial x_{i}}]\mu_{j}=$ $X^{i}[\Lambda_{i}+\frac{\partial}{\partial x_{i}%
}]\phi^{j}\mu_{j}=$ $X^{i}[\frac{\partial}{\partial x_{i}}+\Lambda_{i}]\phi$

Another very short procedure is to use (ii) of the Proposition here where
$\nabla_{\partial_{i}}$ $=$ $\frac{\partial}{\partial x_{i}}+\Lambda_{i}$ and
the fact that $X=X^{i}\frac{\partial}{\partial x_{i}}$ to have:

$\left(  3.31\right)  \qquad\qquad\qquad\qquad\qquad\nabla_{X}\phi$
$=X^{i}\nabla_{\partial_{i}}\phi$ $=$ $X^{i}\left(  \frac{\partial}{\partial
x_{i}}+\Lambda_{i}\right)  \phi$

(ii) $\nabla_{X}^{0}f=$ $<\nabla^{0}f,X>$ $=$ $<g^{ij}\frac{\partial
f}{\partial x_{i}}\frac{\partial}{\partial x_{j}},X^{k}\frac{\partial
}{\partial x_{k}}>$ $=$ $X^{k}g^{ij}\frac{\partial f}{\partial x_{i}}%
<\frac{\partial}{\partial x_{j}},\frac{\partial}{\partial x_{k}}>$

\qquad\qquad\ $=$ $X^{k}g^{ij}\frac{\partial f}{\partial x_{i}}g_{jk}=$
$X^{k}\frac{\partial f}{\partial x_{i}}\delta_{ik}=X^{i}\frac{\partial
f}{\partial x_{i}}$

(iii) We can pair $\nabla\phi\in\Gamma(T^{\ast}M\otimes E)$ and $X\in TM$ for
$\phi\in\Gamma(E)$ as follows:

From (i) of the last Proposition above, we have:

$\qquad\qquad\nabla\phi=$ $\left(  \frac{\partial\phi}{\partial x_{k}}%
+\Lambda_{k}\phi\right)  \otimes dx_{k}$

Consequently,

$\qquad<\nabla\phi,X>$ $=$ $<\left(  \frac{\partial\phi}{\partial x_{k}%
}+\Lambda_{k}\phi\right)  \otimes dx_{k},X^{i}\frac{\partial}{\partial x_{i}%
}>$ $=$ $X^{i}\left(  \frac{\partial\phi}{\partial x_{k}}+\Lambda_{k}%
\phi\right)  X^{i}<dx_{k},\frac{\partial}{\partial x_{i}}>$

$\qquad=X^{i}\left(  \frac{\partial\phi}{\partial x_{k}}+\Lambda_{k}%
\phi\right)  dx_{k}(\frac{\partial}{\partial x_{i}})=X^{i}\left(
\frac{\partial\phi}{\partial x_{k}}+\Lambda_{k}\phi\right)  \delta_{ik}$

$\qquad=X^{i}\left(  \frac{\partial\phi}{\partial x_{i}}+\Lambda_{i}%
\phi\right)  =X^{i}\left(  \frac{\partial}{\partial x_{i}}+\Lambda_{i}\right)
\phi=X^{j}\nabla_{\partial_{j}}\phi$

See p. 65 of \textbf{Berline}, \textbf{Getzler} and \textbf{Vergne }$\left[
7\right]  $ for such a pairing:

We see from this computation and (i) here and (iii) above that:

$\qquad\nabla_{X}\phi=$ $<\nabla\phi,X>$

(iv) In (i) we take $X=\nabla^{0}f=g^{ij}\frac{\partial f}{\partial x_{i}%
}\frac{\partial}{\partial x_{j}}$ and have:

$\qquad<\nabla\phi,\nabla^{0}f>$ $=g^{ij}\frac{\partial f}{\partial x_{i}%
}\nabla_{\partial_{i}}\phi$ $=$ $g^{ij}\frac{\partial f}{\partial x_{i}%
}\left(  \frac{\partial\phi}{\partial x_{i}}+\Lambda_{i}\phi\right)
=g^{ij}\frac{\partial f}{\partial x_{i}}\nabla_{\partial_{j}}\phi$

(v) By (iii) of \textbf{Proposition }$\left(  3.3\right)  ,$

$\qquad\nabla_{\partial_{i}}\nabla_{\partial_{j}}=\frac{\partial^{2}}{\partial
x_{i}\partial x_{j}}+\frac{\partial\Lambda_{j}}{\partial x_{i}}+\Lambda
_{j}\frac{\partial}{\partial x_{i}}+\Lambda_{i}\frac{\partial}{\partial x_{j}%
}+\Lambda_{i}\Lambda_{j}$

Therefore,

$\qquad\qquad\Delta$ \ $=$ $g^{ij}(\nabla_{\partial_{i}}\nabla_{\partial_{j}%
}-\Gamma_{ij}^{k}\nabla_{\partial_{k}})+$ $W$

$\qquad\qquad\ \ \ =g^{ij}\left\{  \frac{\partial^{2}}{\partial x_{i}%
\partial_{j}}+\frac{\partial\Lambda_{j}}{\partial x_{i}}+\Lambda_{j}%
\frac{\partial}{\partial x_{i}}+\Lambda_{i}\frac{\partial}{\partial x_{j}%
}+\Lambda_{i}\Lambda_{j}-\Gamma_{ij}^{k}\left(  \frac{\partial}{\partial
x_{k}}+\Lambda_{k}\right)  \right\}  +W$

This is the local expression of the \textbf{Generalized Laplacian} (or the
\textbf{Laplace-Type Operator}).

(vi)\qquad Let $L=\frac{1}{2}\Delta+$ $\nabla_{X}$ $+$ $V$

Let $\nabla^{0}$ be the gradient operator on functons $f:M\longrightarrow R$
and let $\Delta^{0}$ be the scalar Laplacians on functions.

Then from (v) and (vi), we have:

$L(f\phi)=\frac{1}{2}\Delta(f\phi)+$ $\nabla_{X}(f\phi)$ $+$ $V(f\phi)$

\qquad$=\frac{1}{2}f\Delta\phi+\frac{1}{2}\phi\Delta^{0}f+$ $<\nabla
\phi,\nabla^{0}f>+$ $f$ $\nabla_{X}\phi+\phi\nabla_{X}^{0}f$ $+$ $V(f\phi)$

$\qquad=fL(\phi)+\phi L^{0}(f)+$ $<\nabla\phi,\nabla^{0}f>-$ $V(f\phi)$

$\qquad=\frac{1}{2}g^{ij}f\left\{  \frac{\partial^{2}\phi}{\partial
x_{i}\partial_{j}}+\frac{\partial\Lambda_{j}}{\partial x_{i}}\phi+\Lambda
_{j}\frac{\partial\phi}{\partial x_{i}}+\Lambda_{i}\frac{\partial\phi
}{\partial x_{j}}+\Lambda_{i}\Lambda_{j}\phi-\Gamma_{ij}^{k}\left(
\frac{\partial\phi}{\partial x_{k}}+\Lambda_{k}\phi\right)  \right\}  +fW\phi$

$\qquad+\frac{1}{2}g^{ij}\left\{  \frac{\partial^{2}f}{\partial x_{i}%
\partial_{j}}-\Gamma_{ij}^{k}\frac{\partial f}{\partial x_{k}}\right\}  \phi$
$+$ $g^{ij}\frac{\partial f}{\partial x_{i}}\left(  \frac{\partial\phi
}{\partial x_{i}}+\Lambda_{i}\phi\right)  $

$\qquad$ $+$ $f$ $X^{i}\left(  \frac{\partial\phi}{\partial x_{i}}+\Lambda
_{i}\phi\right)  +\phi X^{i}\frac{\partial f}{\partial x_{i}}$ $+$ $V(f\phi
)$\qquad\qquad\qquad\qquad\qquad\qquad\qquad\qquad\qquad\qquad\qquad
\qquad\qquad\qquad\qquad\qquad\qquad$\qquad\qquad\qquad\qquad\qquad
\qquad\qquad\qquad\qquad\qquad\qquad\qquad\qquad\qquad\qquad\qquad\qquad
\qquad\blacksquare$\qquad

\chapter{Semigroups of Operators and the Generalized Feynman-Kac Formula}

\section{Semi-Classical Semigroups of Operators}

Here we will use the same techniques as in \textbf{Ndumu} $\left[  42\right]
,\left[  43\right]  ,\left[  44\right]  $ and obtain heat kernel formulae and
heat kernel expansions in the more general context of the Laplace-Type
operator on sections of vector bundles.

Let M be a closed (compact without boundary) Riemannian manifold. Let $\Delta$
be the generalized Laplacian or Laplace-Type operator on \textbf{smooth
sections} $\Gamma(E)$ of a vector bundle E. As defined in $\left(  3.2\right)
,$ $\Delta$ is related to the connection $\nabla$ on the vector bundle E by:

$\left(  4.1\right)  \qquad\Delta=$ $g^{ij}(\nabla_{\partial_{i}}%
\nabla_{\partial_{j}}-\Gamma_{ij}^{k}\nabla_{\partial_{k}})+$ $W$

where the \textbf{Weitzenb\"{o}ck term} W is a smooth section of End(E). Let X
be a smooth vector field on M and V a smooth potential term. We will be
dealing with the operator $L=$ $\frac{1}{2}\Delta+$ $\nabla_{X}$ $+$ $V$. This
will usually be written as $L=$ $\frac{1}{2}\Delta+$ $X+V.$

The heat kernel $k_{t}(-,-)$ is a smooth section of the vector bundle
$E\otimes E^{\ast}$ over the product manifold $M\times M:$

$k_{t}(-,-):M\times M\longrightarrow E\otimes E^{\ast}.$ We thus have:
$k_{t}(x,y)\in E_{x}\otimes E_{y}^{\ast}\cong$ Hom( E$_{y}$,E$_{x}$).

Let $k_{t}(x,y):E_{y}\longrightarrow E_{x}$ be the heat kernel on the vector
bundle E relative to the operator $L=$ $\frac{1}{2}\Delta+$ $\nabla_{\text{X}%
}$ $+$ $V.$ The operator $L=\frac{1}{2}\Delta+$ $\nabla_{\text{X}}$ $+$ V here
is a special case of the more general operators considered in \S \ 3, p. 6 of
\textbf{Baudoin} $\left[  5\right]  .$ See also the heat kernel defined in
\S 6 of \textbf{Chapter III }of \textbf{Lawson and Michelsohn} $\left[
35\right]  .$

The heat kernel is the fundamental solution to the heat equation:

$\left(  4.2\right)  \qquad\frac{\partial\phi_{\text{t}}}{\partial\text{t}}=$
$L\phi_{\text{t}}$ \qquad(evolution equation)

\qquad\ \ $\ \ \qquad\phi_{0}=\phi$ \ \qquad\ (initial condition)

where $\phi_{\text{t}},$ $\phi\in\Gamma(E)$ and $\Gamma(E)$ is the set of
smooth sections of the vector bundle E. As was seen in $\left(  2.16\right)
-\left(  2.17\right)  $ of \textbf{Chapter 2} above, the solution of the above
heat equation is given by:

$\left(  4.3\right)  \qquad\phi_{\text{t}}(x)$\ $=$ $\int_{\text{M}}%
$k$_{\text{t}}(x,z)\phi(z)\upsilon_{\text{M}}(dz)$

\qquad The \textbf{generalized heat kernel} on a vector bundle over a compact
Riemannian manifold M relative to a submanifold P was given in $\left(
2.20\right)  $ above by:

$\left(  4.4\right)  $\qquad k$_{\text{t}}$(x,P,$\phi$) = $\int_{\text{P}}%
$k$_{\text{t}}(x,z)\phi(z)\upsilon_{\text{P}}(dz)$

where $\phi\in\Gamma(E)$.

For t$\geq$s$\geq0,$ we define the two semigroups of operators P$_{\text{t}%
}^{\text{M}}$ and Q(t,s) on $\Gamma(E)$ as follows:

$\left(  4.5\right)  $\qquad$\phi_{\text{t}}$(x) = P$_{\text{t}}^{\text{M}%
}\phi($x) = $\int_{\text{M}}$k$_{\text{t}}$(x,z)$\phi(z)\upsilon_{\text{M}}$(dz)

$\left(  4.6\right)  $\qquad Q(t,s)$\phi$(x) = q$_{\text{t}}$(x,P)$^{-1}%
$P$_{\text{t-s}}^{\text{M}}$(q$_{\text{s}}$(-,P)$\phi$)\ 

\qquad\qquad\qquad\qquad\ \ \ = q$_{\text{t}}$(x,P)$^{-1}\int_{\text{M}}%
$q$_{\text{s}}$(z,P)k$_{\text{t-s}}$(x,z)$\phi$(z)$\upsilon_{\text{M}}$(dz)

where $q_{t}(x,P)=(2\pi t)^{-\frac{n-q}{2}}\Psi(x)\exp\left\{  -\frac
{d(x,P)^{2}}{2t}\right\}  $ as defined in $\left(  1.8\right)  $ above.

Using the Chapman-Kolmogorov equation we easily see that for t$\geq$s$\geq$r%
$>$%
0, we have:

$\left(  4.7\right)  $\qquad\qquad Q(t,s)Q(s,r) = Q(t,r)

and hence (Q(t,s))t$\geq$s is a two-parameter semigroup of operators on
$\Gamma(E)$ which we call the \textbf{semi-classical semigroup of operators
}because it related to the semi-classical Brownian Riemannian bridge process
(x$^{\text{t}}$(s)) 0$\leq$ s$\leq$t from x$\in$M$_{0}$ to P in time t as
given by the Generalized Feynman-Kac formula below.

Then, replacing s by t-s in $\left(  4.6\right)  ,$ we have:

$\left(  4.8\right)  \qquad\qquad$Q(t,t-s)$\phi$(x) = q$_{\text{t}}%
$(x,P)$^{-1}$P$_{\text{s}}^{\text{M}}$(q$_{\text{t-s}}$(-,P)$\phi$)(x$)$

\qquad\qquad\qquad\qquad= q$_{\text{t}}$(x,P$)^{-1}\int_{\text{M}}%
$q$_{\text{t-s}}$(z,P)k$_{\text{s}}$(x,z)$\phi$(z)$\upsilon_{\text{M}}$(dz)

The use of the operators Q(t,s) on smooth sections of the vector bundle E will
prove to be very useful in obtaining both a generalized\textbf{\ Feynman-Kac}
formula and the \textbf{exact} and \textbf{asymptotic} expansion formulae for
the generalized heat kern\qquad ut in \textbf{Ndumu} $\left[  42\right]  .$

First we give some preliminary lemmas which are needed to prove the theorems here:

\qquad\qquad\qquad\qquad\qquad\qquad\qquad\qquad\qquad\qquad\qquad\qquad
\qquad\qquad\qquad\qquad\qquad\qquad$\blacksquare$

\begin{lemma}
Let P$\ $be a q-dimensional compact submanifold of a compact manifold M. Then,
we have:
\end{lemma}

$\underset{\text{s}\uparrow\text{t}}{lim}($Q(t,t-s)$\phi)$(x) = q$_{\text{t}}%
$(x,P)$^{-1}$ $\int_{\text{P}}$k$_{\text{t}}$(x,y)$\phi$(y)$\upsilon
_{\text{P}}$(dy)

\begin{proof}
We have by $\left(  4.8\right)  ,$

\qquad(Q(t,t-s) $\phi$)(x) = q$_{\text{t}}$(x,P$)^{-1}$P$_{\text{s}}%
^{\text{M}}$(q$_{\text{t-s}}$(-,P)$\phi$)(x)
\end{proof}

\qquad\qquad\qquad\ \ \ \ \ \ \ \ \ \qquad

\qquad\ = q$_{\text{t}}$(x,P$)^{-1}\int_{\text{M}}$q$_{\text{t-s}}%
$(z,P)k$_{\text{s}}$(x,z)$\phi$(z)$\upsilon_{\text{M}}$(dz)

We use \textbf{Lemma} $\left(  2.2\right)  $ of Chapter 2 in \textbf{Ndumu}
$\left[  40\right]  $ to assume that the exponential map of the normal bundle
is a global diffeomorphism or we use \textbf{Theorem 8.40} of \textbf{Gray
}$\left[  25\right]  $ (which is a generalization of \textbf{Lemma 8.3 of Gray
}$\left[  25\right]  $) which states (with a slight error there)\ that for M
compact we have:

$(4.9)$\qquad volume(M) = volume(exp$_{\upsilon}(B_{0}))=$ volume(M$_{0}):$

Making the change of variable z = exp$_{\Pi}$(y,v) = exp$_{\text{y}}$v, and
noting that $\phi$(z) depends smoothly on z we have:

\qquad(Q(t,t-s) $\phi$)(x) \ = q$_{\text{t}}$(x,P)$^{-1}\int_{\text{M}}%
$q$_{\text{t-s}}$(z,P)k$_{\text{s}}$(x,z)$\phi$(z)$\upsilon_{\text{M}}$(dz)

\ \qquad\qquad= q$_{\text{t}}$(x,P)$^{-1}\int_{\text{M}_{0}}$q$_{\text{t-s}}%
$(z,P)k$_{\text{s}}$(x,z)$\phi$(z)$\upsilon_{\text{M}}$(dz)

\qquad(The first equation is by definition and the last equation above is by
$\left(  4.9\right)  )$

\qquad\qquad= q$_{\text{t}}$(x,P)$^{-1}\int_{\text{B}_{0}}$q$_{\text{t-s}}%
$(exp$_{\text{y}}$v,P)k$_{\text{s}}$(x,exp$_{\text{y}}$v)$\phi$(exp$_{\text{y}%
}$v)$\theta_{\text{P}}$(v)$\upsilon_{\text{P}}$(dy)dv

\qquad\qquad= q$_{\text{t}}$(x,P)$^{-1}\int_{\text{P}}\int_{R^{\text{n-q}}}%
$q$_{\text{t-s}}$(exp$_{\text{y}}$v,P)k$_{\text{s}}$(x,exp$_{\text{y}}$%
v)$\phi$(exp$_{\text{y}}$v)$\theta_{\text{P}}$(v)$\upsilon_{\text{P}}$(dy)dv

Setting r = t-s and v = $\sqrt{\text{r}}$w, then the first and last equations give:

$(4.10)$\qquad q$_{\text{t}}$(x,P)(Q(t,t-s) $\phi$)(x)

\qquad= $\int_{\text{P}}\int_{R^{\text{n-q}}}$q$_{\text{r}}$(exp$_{\text{y}%
}\sqrt{\text{r}}$w,P)k$_{\text{s}}$(x,exp$_{\text{y}}\sqrt{\text{r}}$w)$\phi
$(exp$_{\text{y}}\sqrt{\text{r}}$w)$\theta_{\text{P}}$($\sqrt{\text{r}}%
$w)$\upsilon_{\text{P}}$r$^{\frac{\text{n-q}}{2}}$(dy)dw

We then take limits as r$\downarrow$0, which is equivalent to s$\uparrow$t, on
both sides of (4.10).

Noting that the limit on the RHS of $\left(  4.10\right)  $ is similar to the
limit in Proposition 2 in Chapter 2. We then use the following elementary facts:

exp$_{\text{y}}$0 = y and $\theta_{\text{P}}$(0)=1 to have:

\qquad q$_{\text{r}}$(exp$_{\text{y}}\sqrt{\text{r}}$w,P) = (2$\pi$%
r)$^{-\frac{\text{n-q}}{2}}\Psi$(exp$_{\text{y}}\sqrt{\text{r}}$w)
exp$\left\{  -\frac{\text{d(exp w}\sqrt{\text{r}}_{y}\text{, P)}^{2}}%
{2r}\right\}  $

\qquad\ = (2$\pi$r)$^{-\frac{\text{n-q}}{2}}\Psi$(exp$_{\text{y}}%
\sqrt{\text{r}}$w) exp$\left\{  -\frac{\left\Vert \text{w}\right\Vert ^{2}}%
{2}\right\}  $

and so,

q$_{\text{t}}$(x,P)(Q(t,t-s) $\phi$)(x) = $\int_{\text{P}}\int_{R^{\text{n-q}%
}}$q$_{\text{r}}$(exp$_{\text{y}}\sqrt{\text{r}}$w,P)k$_{\text{s}}%
$(x,exp$_{\text{y}}\sqrt{\text{r}}$w)$\phi$(exp$_{\text{y}}\sqrt{\text{r}}%
$w)$\theta_{\text{P}}$($\sqrt{\text{r}}$w)$\upsilon_{\text{P}}$r$^{\frac
{\text{n-q}}{2}}$(dy)dw

An easy simplification (r$^{-\frac{\text{n-q}}{2}}$ and r$^{\frac{\text{n-q}%
}{2}}$ cancel out) gives:

q$_{\text{t}}$(x,P)(Q(t,t-s) $\phi$)(x)

\qquad= (2$\pi$)$^{-\frac{\text{n-q}}{2}}\int_{\text{P}}\int_{R^{\text{n-q}}%
}\Psi$(exp$_{\text{y}}\sqrt{\text{r}}$w) exp$\left\{  -\frac{\left\Vert
\text{w}\right\Vert ^{2}}{2}\right\}  $k$_{\text{s}}$(x,exp$_{\text{y}}%
\sqrt{\text{r}}$w)$\phi$(exp$_{\text{y}}\sqrt{\text{r}}$w)$\theta_{\text{P}}%
$($\sqrt{\text{r}}$w)$\upsilon_{\text{P}}$(dy)dw

It is clear that since r = t-s the limit s$\uparrow$t on the RHS is equivalent
to the limit r$\downarrow$0 on the LHS:

We now take the limits on both sides of the last equality above to have:

$(4.11)$\qquad q$_{\text{t}}$(x,P) $\underset{\text{s}\uparrow\text{t}}{lim}%
$(Q(t,t-s) $\phi$)(x) \qquad\qquad\qquad

= (2$\pi$)$^{-\frac{\text{n-q}}{2}}$ $\underset{r\downarrow0}{lim}%
\int_{\text{P}}\int_{R^{\text{n-q}}}\Psi$(exp$_{\text{y}}\sqrt{\text{r}}$w)
exp$\left\{  -\frac{\left\Vert \text{w}\right\Vert ^{2}}{2}\right\}
$k$_{\text{s}}$(x,exp$_{\text{y}}\sqrt{\text{r}}$w)$\phi$(exp$_{\text{y}}%
\sqrt{\text{r}}$w)$\theta_{\text{P}}$($\sqrt{\text{r}}$w)$\upsilon_{\text{P}}$(dy)dw

\qquad= (2$\pi$)$^{-\frac{\text{n-q}}{2}}\int_{\text{P}}\int_{R^{\text{n-q}}%
}\Psi$(0) exp$\left\{  -\frac{\left\Vert \text{w}\right\Vert ^{2}}{2}\right\}
$k$_{\text{s}}$(x,exp$_{\text{y}}$0)$\phi$(exp$_{\text{y}}$0)$\theta
_{\text{P}}$(0)$\upsilon_{\text{P}}$(dy)dw

Since exp$_{\text{y}}$0=y and $\theta_{\text{P}}$(0)=1= $\Phi$(0) and so,
$\Psi$(0)= $\theta_{\text{P}}$(0)$\Phi$(0)=1.

We see that for all x$\in$M$_{0}$:

$\left(  4.12\right)  $\qquad q$_{\text{t}}$(x,P) $\underset{\text{s}%
\uparrow\text{t}}{lim}$(Q(t,t-s) $\phi$)(x) = (2$\pi$)$^{-\frac{\text{n-q}}%
{2}}\int_{\text{P}}\int_{R^{\text{n-q}}}$exp$\left\{  -\frac{\left\Vert
\text{w}\right\Vert ^{2}}{2}\right\}  $k$_{\text{s}}$(x,y)$\phi$%
(y)$\upsilon_{\text{P}}$(dy)dw

Finally, since $\int_{R^{\text{n-q}}}$exp$\left\{  -\frac{\left\Vert
\text{w}\right\Vert ^{2}}{2}\right\}  dw=$ (2$\pi$)$^{\frac{\text{n-q}}{2}}$,
we get the result.

\qquad\qquad\qquad\qquad\qquad\qquad\qquad\qquad\qquad\qquad\qquad\qquad
\qquad\qquad\qquad\qquad\qquad\qquad$\blacksquare$

\begin{lemma}
{\small For }$\phi\in\Gamma${\small (E),\qquad}
\end{lemma}

$\frac{\partial}{\partial\text{t}}${\small (Q(t,s)}$\phi${\small ) =}$\left(
\text{L + }\nabla^{0}\log\text{q}_{\text{t}}\text{(-,P)}+\frac{\text{L}%
^{0}\Psi}{\Psi}-V\right)  ${\small (Q(t,s)}$\phi${\small )\qquad}

{\small where L = }$\frac{1}{2}\Delta${\small +X+V and }$\Delta$ {\small is
the Laplace-Type operator.}

\begin{proof}
$\frac{\partial}{\partial\text{t}}${\small (Q(t,s)}$\phi${\small ) = }%
$\frac{\partial}{\partial\text{t}}\left(  \text{q}_{\text{t}}\text{(-,P)}%
^{-1}\text{P}_{\text{t-s}}^{\text{M}}\text{(q}_{\text{s}}\text{(-,P)}%
\phi\right)  $
\end{proof}

{\small \qquad= - q}$_{\text{t}}${\small (-,P)}$^{-2}\frac{\partial}%
{\partial\text{t}}${\small q}$_{\text{t}}${\small (-,P)}$\left(
\text{P}_{\text{t-s}}^{\text{M}}\text{(q}_{\text{s}}\text{(-,P)}\phi\right)
${\small + q}$_{\text{t}}${\small (-,P)}$^{-1}${\small .}$\frac{\partial
}{\partial\text{t}}\left(  \text{P}_{\text{t-s}}^{\text{M}}\text{(q}%
_{\text{s}}\text{(-,P)}\phi\right)  $

{\small \qquad= I}$_{1}${\small + I}$_{2}${\small , where,}

{\small \ \ I}$_{1}${\small = - q}$_{\text{t}}${\small (-,P)}$^{-2}%
\frac{\partial}{\partial\text{t}}${\small q}$_{\text{t}}${\small (-,P)}%
$\left(  \text{P}_{\text{t-s}}^{\text{M}}\text{(q}_{\text{s}}\text{(-,P)}%
\phi\right)  $

{\small \qquad=- q}$_{\text{t}}${\small (-,P)}$^{-1}\frac{\partial}%
{\partial\text{t}}${\small q}$_{\text{t}}${\small (-,P)}$\left(
\text{P}_{\text{t-s}}^{\text{M}}\text{(q}_{\text{s}}\text{(-,P)}\phi\right)
${\small q}$_{\text{t}}${\small (-,P)}$^{-1}$

{\small \ \ I}$_{2}${\small = q}$_{\text{t}}${\small (-,P)}$^{-1}$%
{\small .}$\frac{\partial}{\partial\text{t}}\left(  \text{P}_{\text{t-s}%
}^{\text{M}}\text{(q}_{\text{s}}\text{(-,P)}\phi\right)  ${\small .}

{\small By Lemma 1 above we have:}

$\left(  4.13\right)  ${\small \qquad I}$_{1}${\small \ = - q}$_{\text{t}}%
${\small (-,P)}$^{-1}${\small \ }$\left\{  \text{L}^{0}\text{q}_{\text{t}%
}\text{(-,P)}-\frac{\text{L}^{0}\Psi}{\Psi}.\text{q}_{\text{t}}\text{(-,P)}%
\right\}  ${\small .q}$_{\text{t}}${\small (-,P)}$^{-1}\left(  P_{\text{t-s}%
}^{\text{M}}(q_{\text{s}}(-,P)\phi\right)  $

{\small Then by }$\left(  4.8\right)  ${\small ,}

{\small \qquad\qquad\qquad I}$_{1}${\small = \ }$\left\{  -\frac{\text{L}%
^{0}\text{q}_{\text{t}}\text{(-,P)}}{\text{q}_{\text{t}}(-,\text{P})}%
+\frac{\text{L}^{0}\Psi}{\Psi}\right\}  ${\small .(Q(t,s)}$\phi${\small )}

{\small Next we have:}

$\left(  4.14\right)  ${\small \qquad I}$_{2}${\small = q}$_{\text{t}}%
${\small (-,P)}$^{-1}${\small .}$\frac{\partial}{\partial\text{t}}\left(
\text{P}_{\text{t-s}}^{\text{M}}\text{(q}_{\text{s}}\text{(-,P)}\phi
\text{)}\right)  $

{\small \ \ \ \qquad\qquad\ = q}$_{\text{t}}${\small (-,P)}$^{-1}$%
{\small L}$\left[  \text{P}_{\text{t-s}}^{\text{M}}\text{(q}_{\text{s}%
}\text{(-,P)}\phi\text{)}\right]  ${\small by Theorem }$\left(  1,6\right)
${\small of Azencott }$\left[  4\right]  $

{\small \ \ \ \qquad\qquad\ = q}$_{\text{t}}${\small (-,P)}$^{-1}$%
{\small L}$\left[  \text{(q}_{\text{t}}\text{(-,P)}^{-1}\text{.P}_{\text{t-s}%
}^{\text{M}}\text{(q}_{\text{s}}\text{(-,P)}\phi\text{).q}_{\text{t}%
}\text{(-,P)}\right]  $

{\small \ \ \qquad\qquad\ \ = q}$_{\text{t}}${\small (-,P)}$^{-1}$%
{\small L}$\left[  \text{(Q(t,s)}\phi\text{).q}_{\text{t}}\text{(-,P)}\right]
${\small by the definition of Q(t,s) in }$\left(  4.6\right)  ${\small .}

{\small \ = q}$_{\text{t}}${\small (-,P)}$^{-1}${\small [L((Q(t,s)}$\phi
${\small )).q}$_{\text{t}}${\small (-,P) + (Q(t,s)}$\phi${\small ).L}$^{0}%
${\small (q}$_{\text{t}}${\small (-,P)) +
$<$%
}$\nabla^{0}${\small q}$_{\text{t}}${\small (-,P),}$\nabla${\small (Q(t,s)}%
$\phi${\small )%
$>$%
}

{\small \qquad- V((Q(t,s)}$\phi${\small ).q}$_{\text{t}}${\small (-,P))]\qquad
(by the property of the differential operator L}$^{0}${\small =}$\frac{1}%
{2}\Delta^{0}${\small +X+V)}

{\small = L}$\left(  \text{(Q(t,s)}\phi\text{)}\right)  ${\small +
(Q(t,s)}$\phi${\small )}$\frac{\text{L}^{0}\text{q}_{\text{t}}\text{(-,P)}%
}{\text{q}_{\text{t}}\text{(-,P)}}${\small +
$<$%
}$\nabla^{0}\log${\small q}$_{\text{t}}${\small (-,P),}$\nabla$%
{\small (Q(t,s)}$\phi${\small )%
$>$
- V((Q(t,s)}$\phi${\small )}

{\small \ = }$\left[  \text{L}+\nabla^{0}\log\text{q}_{\text{t}}%
\text{(-,P)}+\frac{\text{L}^{0}\text{q}_{\text{t}}\text{(-,P)}}{\text{q}%
_{\text{t}}\text{(-,P)}}-\text{V}\right]  ${\small (Q(t,s)}$\phi${\small )}

Then adding the final expressions of I$_{1}$and I$_{2}$in $\left(
4.13\right)  $and $\left(  4.14\right)  $, we have:

$\frac{\partial}{\partial\text{t}}${\small (Q(t,s)}$\phi${\small ) = I}$_{1}%
${\small + I}$_{2}${\small \ \ \ \ =}$\left\{  -\frac{\text{L}^{0}%
\text{q}_{\text{t}}(-,\text{P})}{\text{q}_{\text{t}}(-,\text{P})}%
+\frac{\text{L}^{0}\Psi}{\Psi}\right\}  ${\small .(Q(t,s)}$\phi${\small )}

{\small \qquad\qquad\qquad\qquad\qquad\ \ \ \ +}$\left\{  \text{L}+\nabla
^{0}\log\text{q}_{\text{t}}(-,\text{P})+\frac{\text{L}^{0}\text{q}_{\text{t}%
}(-,\text{P})}{\text{q}_{\text{t}}(-,\text{P})}-\text{V}\right\}
${\small (Q(t,s)}$\phi${\small )}

$\left(  4.15\right)  ${\small \ \ \ \qquad\qquad\ \qquad\ \ =}$\left(
\text{L + }\nabla^{0}\log\text{q}_{\text{t}}\text{(-,P)}+\frac{\text{L}%
^{0}\Psi}{\Psi}-\text{V}\right)  ${\small (Q(t,s)}$\phi${\small )}

{\small \qquad\qquad\qquad\qquad\qquad\qquad\qquad\qquad\qquad\qquad
\qquad\qquad\qquad\qquad\qquad\qquad\qquad}$\blacksquare$

\section{The Generalized Feynman-Kac Formula}

The next objective is to obtain an importnat Generalized Feynman-Kac formula
from which we shall deduce the generalized heat kernel formula for vector bundles.

The definitions and formulations here are adapted to obtaining almost
simultaneously both the Feynman-Kac formula and the heat kernel formula rather
than just the Feynman-Kac formula as given in \textbf{Bismut }$\left[
8\right]  $; \textbf{Driver and Thalmaier }$\left[  12\right]  $ and
\textbf{Norris }$\left[  46\right]  $.

For the definition of the bridge process we will follow \textbf{Elworthy
}$\left[  14\right]  $, Theorem 4B of chapter 5 . For a fixed t%
$>$%
0 let 0$\leq$s%
$<$%
t. Set:

$\left(  4.15\right)  $\qquad L$_{\text{s,t}}^{0}$= L$^{0}$+ $\nabla^{0}\log
$q$_{\text{t-s}}$(-,P) and L$_{\text{s,t}}$= L + $\nabla^{0}\log
$q$_{\text{t-s}}$(-,P)

where L$^{0}$ = $\frac{1}{2}\Delta^{0}$+ X + V and L =$\frac{1}{2}\Delta$+ X +
V for a smooth vector field X on M and a smooth scalar potential term V which
represent heat sources and sinks on the Riemannian manifold (without boundary)
M. The differential operators L$_{\text{s,t}}^{0}$and L$_{\text{s,t}}$depend
on the time variable s and the fixed time parameter t%
$>$%
0.\ \ \ 

\qquad\qquad\qquad\qquad\qquad\qquad\qquad\qquad\qquad\qquad\qquad\qquad
\qquad\qquad\qquad$\blacksquare$

\begin{lemma}
(Existence)
\end{lemma}

Given a curve x$_{s}$0$\leq$s$\leq$t, and a frame u$_{0}$at x$_{0}$i.e. an
isomorphism u$_{0}$:R$^{n}\longrightarrow$E$_{x_{0}}$, there exists a unique

curve u$_{s}$: R$^{n}\longrightarrow$E$_{\text{x}_{s}}$0$\leq$s$\leq$t such
that \U{e31a}(u$_{s}$)= x$_{s}$

See Hsu $\left[  31\right]  $, p. 38 in the special case where E =
TM{\small .}

\begin{definition}
\textsc{\qquad\qquad\qquad\qquad\qquad\qquad\qquad\qquad\qquad\qquad\qquad
}\qquad\qquad$\ \ \ \blacksquare$
\end{definition}

The curve (u$_{s})$ 0$\leq s\leq t$ is called the \textbf{horizontal lift} of
(x$_{s})$ 0$\leq s\leq t$ from u$_{0}$ to the vector bundle E$.$

\qquad\qquad\qquad\qquad\qquad\qquad\qquad\qquad\qquad\qquad\qquad\qquad
\qquad\qquad\qquad$\blacksquare$\qquad\ \ \ 

Let (x$_{s}^{\text{t}}$) $0\leq$s$\leq$t$\wedge\zeta$ be the diffusion process
in M with differential generator L$_{\text{s,t}}^{0}$ = $\frac{1}{2}\Delta
^{0}$ $+$ X $+$ $\nabla^{0}\log$q$_{\text{t-s}}$(-,P) +V applied to smooth
functions on M. The process (x$_{s}^{\text{t}}$) $0\leq$s$\leq$t$\wedge\zeta$
is the \textbf{semi-classical Brownian Riemannian bridge process} from x$\in
$M$_{0}$ to the submanifold P in time t with exit time $\zeta$ from the
tubular neighborhood M$_{0}$ of P$.$

The generator L$_{\text{s,t}}$ = $\frac{1}{2}\Delta$ $+$ X $+$ $\nabla^{0}%
\log$q$_{\text{t-s}}$(-,P) + V here is applied to sections of the vector
bundle E where $\Delta$ is the Laplace-Type operator defined in $\left(
3.2\right)  .$

It was shown in \textbf{Lemma 3.3} of \textbf{Ndumu} $\left[  40\right]  $ and
in \textbf{Ndumu }$\left[  42\right]  $ that x$^{\text{t}}$(t) $=$ y(t) a. s.
where y(t) 0$\leq$t%
$<$%
+$\infty$\ is a process in the small neighborhood U$\subset$ P of the centre
of Fermi coordinates y$_{0}\in P.$

The \textbf{stochastic covariant equation} along the paths of the above bridge
process is similar to the covariant equations in \textbf{Elworthy} $\left[
16\right]  .$ We will follow of \textbf{Elworthy} $\left[  16\right]  $ for
the definition of the stochastic covariant equation along the paths of
Brownian motion:

$\left(  4.16\right)  \qquad$De$_{\text{s}}^{\text{t}}=\frac{1}{2}%
$e$_{\text{s}}^{\text{t}}$W$_{\text{x}_{\text{s}}^{\text{t}}}$

\qquad\qquad\ \ \ \ \ e$_{\text{0}}^{\text{t}}=$ identity

where W$\in\Gamma$(End(E)) is a Weitzenb\H{o}ckian and (e$_{\text{s}%
}^{\text{t}}$ ) 0$\leq$s%
$<$%
t$\Lambda\zeta$ is defined in $\left(  4.20\right)  $ below.

Following $\left(  4.6\right)  -\left(  47\right)  $ of \textbf{Elworthy}
$\left[  16\right]  ,$ the equation in $\left(  4.16\right)  $ above is a
short-hand for the equation:

$\left(  4.17\right)  \qquad$u$_{\text{r}}^{\text{t}}$(u$_{0}^{\text{t}}%
$)$^{-1}\frac{\text{d}}{\text{dr}}\left\{  \text{u}_{\text{0}}^{\text{t}%
}\text{(u}_{\text{r}}^{\text{t}}\text{)}^{-1}\text{e}_{\text{r}}^{\text{t}%
}\right\}  =\frac{1}{2}$e$_{\text{r}}^{\text{t}}$W$_{\text{x}_{s}^{\text{t}}}%
$\qquad

where (u$_{\text{s}}^{\text{t}}$ ) 0$\leq$s$<$t$\Lambda\zeta$ is the unique
\textbf{horizontal lift} of the bridge process (x$_{s}^{\text{t}}$) 0$\leq
$s$<$t$\Lambda\zeta$ starting at u$_{0}$ on the vector bundle GL(R$^{n}$,E) =
$\underset{x\in M}{\cup}$GL$_{x}$(R$^{n}$,E), where GL$_{x}$(R$^{n}$,E) =
$\left\{  u:R^{n}\longrightarrow E_{x}\text{ is a linear isomorphism}\right\}
.$\qquad\qquad\qquad

Therefore, for 0$\leq r\leq s<t,$ we have the map: u$_{\text{0}}^{\text{t}}%
$(u$_{\text{s}}^{\text{t}}$)$^{-1}:$E$_{\text{x}^{\text{t}}\text{(s)}%
}\longrightarrow$ E$_{\text{x}_{0}}$ which has for inverse:

u$_{\text{s}}^{\text{t}}$(u$_{\text{0}}^{\text{t}}$)$^{-1}:$E$_{\text{x}_{0}%
}\longrightarrow$ E$_{_{\text{x}^{\text{t}}\text{(s)}}}$

We take the composition: u$_{\text{r}}^{\text{t}}$(u$_{\text{0}}^{\text{t}}%
$)$^{-1}\circ$ u$_{\text{0}}^{\text{t}}$(u$_{\text{s}}^{\text{t}}$)$^{-1}=$
u$_{\text{r}}^{\text{t}}$(u$_{\text{s}}^{\text{t}}$)$^{-1}$ to have:

$\left(  4.18\right)  \qquad\tau_{\text{r,s}}^{\text{t}}:=$ u$_{\text{r}%
}^{\text{t}}$(u$_{\text{s}}^{\text{t}}$)$^{-1}:$E$_{\text{x}^{\text{t}%
}\text{(s)}}\longrightarrow$ E$_{_{\text{x}^{\text{t}}\text{(r)}}}$

with inverse

$\qquad\qquad\tau_{\text{s,r}}^{\text{t}}:=$ u$_{\text{s}}^{\text{t}}%
$(u$_{\text{r}}^{\text{t}}$)$^{-1}:$E$_{\text{x}^{\text{t}}\text{(r)}%
}\longrightarrow$ E$_{_{\text{x}^{\text{t}}\text{(s)}}}$

These maps are parallel translations on fibers of the vector bundle E along
the semi-classical Brownian Riemannian bridge process (x$_{s}^{\text{t}}$)
0$\leq$s$\leq$t$\Lambda\zeta.$

The equation in $\left(  4.17\right)  $ can be re-written as:

$\qquad\frac{\text{d}}{\text{dr}}\left\{  \text{u}_{\text{0}}^{\text{t}%
}\text{(u}_{\text{r}}^{\text{t}}\text{)}^{-1}\text{e}_{\text{r}}^{\text{t}%
}\right\}  =\frac{1}{2}($u$_{\text{0}}^{\text{t}}$(u$_{\text{r}}^{\text{t}}%
$)$^{-1}$e$_{\text{r}}^{\text{t}})($W$_{\text{x}^{\text{t}}\text{(r)}})$%
\qquad\qquad\qquad\qquad

Integrating, we have:

$\left(  4.19\right)  \qquad$u$_{\text{0}}^{\text{t}}$(u$_{\text{s}}%
^{\text{t}}$)$^{-1}$e$_{\text{s}}^{\text{t}}-$ e$_{\text{0}}^{\text{t}}%
=\frac{1}{2}\int_{0}^{\text{s}}($u$_{\text{0}}^{\text{t}}$(u$_{\text{r}%
}^{\text{t}}$)$^{-1}$e$_{\text{r}}^{\text{t}})($W$_{\text{x}_{\text{r}%
}^{\text{t}}}$)dr

We then take left compositions on each side of the equality in $\left(
4.19\right)  $ by u$_{\text{s}}^{\text{t}}($u$_{0}^{\text{t}}$)$^{-1}$ and get
the expression for the covariant process (e$_{\text{s}}^{\text{t}}$ ) 0$\leq$s%
$<$%
t$\Lambda\zeta$ in terms of the horizontal lift process (u$_{s}^{\text{t}}$)
0$\leq$s$<$t$\Lambda\zeta$ and the Weitzenb\H{o}ck term W:\ \ \ \ \ \ \ \ \ \ \ \ \ \ \ \ \ \ \ \ \ \ \ \ \ \ \ \ \ \ 

\ $\qquad$e$_{\text{s}}^{\text{t}}=$ u$_{\text{s}}^{\text{t}}$(u$_{0}%
^{\text{t}}$)$^{-1}$e$_{\text{0}}^{\text{t}}+\frac{1}{2}\int_{0}^{\text{s}}%
$u$_{\text{s}}^{\text{t}}$(u$_{\text{r}}^{\text{t}}$)$^{-1}$e$_{\text{r}%
}^{\text{t}}$W$_{\text{x}_{r}^{\text{t}}}$)dr \ \ \ 

Since e$_{\text{0}}^{\text{t}}$ is an identity we have for 0$\leq$r$\leq$s%
$<$%
t$\Lambda\zeta,$\ \ \ 

$\left(  4.20\right)  $\qquad e$_{\text{s}}^{\text{t}}\ =$ $\tau_{\text{s,0}%
}^{\text{t}}+\frac{1}{2}\int_{0}^{\text{s}}\tau_{\text{s,r}}^{\text{t}}%
$e$_{\text{r}}^{\text{t}}$W$_{\text{x}_{\text{r}}^{\text{t}}}$dr \ \ \ 

The solution e$_{\text{s}}^{\text{t}}:$E$_{\text{x}}\longrightarrow$
E$_{_{\text{x}^{\text{t}}\text{(s)}}}$ of the equation in $\left(
4.20\right)  $ gives rise to a process (e$_{\text{s}}^{\text{t}})$ 0$\leq
s<$t$\Lambda\zeta$ over the semi-classical Brownian Riemannian bridge process
(x$_{\text{s}}^{\text{t}}$) 0$\leq$s$<$t$\Lambda\zeta.$

\qquad\qquad\qquad\qquad\qquad\qquad\qquad\qquad\qquad\qquad\qquad\qquad
\qquad\qquad\qquad\qquad\qquad\qquad$\blacksquare$

We now come to one of the \textbf{central theorems} of this work.

\begin{theorem}
(Generalized Feynman-Kac Formula)
\end{theorem}

Let $\phi$ be a smooth section of the vector bundle E over the compact
Riemannian manifold M and let $\varsigma$ be the first exit time of the
semi-classical Brownian Riemannian bridge process from the tubular
neighborhood M$_{0}.$ Let \textbf{E}$_{\text{x}}$ be the expectation
corresponding to the Wiener measure \textbf{P}$_{\text{x}}$ on $(\Omega,$
$\digamma,$ $P_{\text{x}})$ of Brownian motion (Wiener process) starting from
x$\in$M$_{0}.$ Then for 0$\leq$ s$\leq$t$,$

\qquad(Q(t,t-s)$\phi$)(x) $=$ \textbf{E}$_{\text{x}}\left[  \tau_{\text{0,s}%
}^{\text{t}}\text{e}_{\text{s}}^{\text{t}}\phi\text{(x}^{\text{t}}%
\text{(s))}\exp\left\{  \int_{0}^{\text{s}}\frac{\text{L}^{0}\Psi}{\Psi
}\text{(x}^{\text{t}}\text{(u))du}\right\}  \right]  $\qquad\qquad\ \qquad\ 

where (x$^{\text{t}}$(s)) 0$\leq$ s$\leq$t is the semi-classical Brownian
Riemannian bridge process from a point x$\in$M$_{0}$ to (a point y$_{0}$) in
the submanifold P in time t and $\tau_{\text{0,s}}^{\text{t}}=$ u$_{\text{0}%
}^{\text{t}}$(u$_{\text{s}}^{\text{t}}$)$^{-1}:$E$_{\text{x}^{\text{t}%
}\text{(s)}}\longrightarrow$E$_{\text{x}}$ is parallel translation along the
reversed semi-classical Brownian bridge process and e$_{\text{s}}^{\text{t}%
}\ \ $is defined in $\left(  4.20\right)  $ above.

\begin{proof}
\qquad\ \ \ \ \ \ 
\end{proof}

For $\lambda\geq$t$\geq$s$\geq r\geq0$ set,

\qquad\qquad\ \ \ $\phi_{\lambda}$(x) $=$ (Q($\lambda$,t-s)$\phi$)(x) and
h($\lambda,$x,w) = $\phi_{\lambda}$(x)w

$\qquad$h(y(r)) = h(t-r,x$^{\text{t}}$(r),w(r)) = $\phi_{\text{t}-\text{r}}%
$(x$^{\text{t}}$(r))w(r) where y(r) $=$ $\left(  \text{t-r,x}^{\text{t}%
}\text{(r),w(r)}\right)  $

and\ w(r) = exp$\left\{  \int_{0}^{\text{r}}\frac{\text{L}^{0}\Psi}{\Psi
}\text{(x}^{\text{t}}\text{(u))du}\right\}  $ and so dw(r) = w(r)$\frac
{\text{L}^{0}\Psi}{\Psi}$(x$^{\text{t}}$(r))dr

We note that since $\phi$ is a smooth section of the vector bundle E, then
$\phi_{\text{t}}=$ Q(t,t-s)$\phi$ is a time-dependent smooth section of the
vector bundle E and hence h is also a time-dependent smooth section of the
vector bundle E.

Then \textbf{It\^{o}'s formula} (differential version)\textbf{ }for a product
of two C$^{2}-$ functions: f$_{s},$g$_{s}:R\longrightarrow R$ gives:

\begin{center}
d(f$_{s}$(x$_{s}$)g$_{s}(x_{s}))=$ df$_{s}$(x$_{s}$)g$_{s}$(x$_{s}$) $+$
f$_{s}$(x$_{s}$)dg$_{s}(x_{s})+\frac{1}{2}$df$_{s}$(x$_{s}$).dg$_{s}(x_{s})$

$=$ df$_{s}$(x$_{s}$)g$_{s}$(x$_{s}$) $+$ f$_{s}$(x$_{s}$)dg$_{s}(x_{s}%
)+\frac{1}{2}$d$<$f$,$g$>_{s}$
\end{center}

d$\left(  \tau_{\text{0,r}}^{\text{t}}\text{e}_{\text{r}}^{\text{t}%
}\text{h(y(r))}\right)  $ = d$\left(  \tau_{\text{0,r}}^{\text{t}}%
\text{e}_{\text{r}}^{\text{t}}\phi_{\text{t-r}}\text{(x}^{\text{t}}%
\text{(r))}\right)  $w(r) + $\tau_{\text{0,r}}^{\text{t}}$e$_{\text{r}%
}^{\text{t}}\phi_{\text{t-r}}$(x$^{\text{t}}$(r))dw(r) + $\frac{1}{2}%
$d$\left(  \tau_{\text{0,r}}^{\text{t}}\text{e}_{\text{r}}^{\text{t}}%
\phi_{\text{t-r}}\text{(x}^{\text{t}}\text{(r))}\right)  $dw(r)

Since dw(r) = w(r)$\frac{\text{L}^{0}\Psi}{\Psi}$(x$^{\text{t}}$(r))dr, we have:

\qquad d$\left(  \tau_{\text{0,r}}^{\text{t}}\text{e}_{\text{r}}^{\text{t}%
}\phi_{\text{t-r}}\text{(x}^{\text{t}}\text{(r))}\right)  $dw(r) =
w(r)$\frac{\text{L}^{0}\Psi}{\Psi}$(x$^{\text{t}}$(r))d$\left(  \tau
_{\text{0,r}}^{\text{t}}\text{e}_{\text{r}}^{\text{t}}\phi_{\text{t-r}%
}\text{(x}^{\text{t}}\text{(r))}\right)  $dr = 0

and so we have:

$\left(  4.22\right)  \qquad$d$\left(  \tau_{\text{0,r}}^{\text{t}}%
\text{e}_{\text{r}}^{\text{t}}\text{h(y(r))}\right)  $ \ \ \ \ = d$\left(
\tau_{\text{0,r}}^{\text{t}}\text{e}_{\text{r}}^{\text{t}}\phi_{\text{t-r}%
}\text{(x}^{\text{t}}\text{(r))}\right)  $w(r) + $\tau_{\text{0,r}}^{\text{t}%
}$e$_{\text{r}}^{\text{t}}\phi_{\text{t-r}}$(x$^{\text{t}}$(r))$\frac
{\text{L}^{0}\Psi}{\Psi}$(x$^{\text{t}}$(r))w(r)dr

We now expand d$\left(  \tau_{\text{0,r}}^{\text{t}}\text{e}_{\text{r}%
}^{\text{t}}\phi_{\text{t-r}}\text{(x}^{\text{t}}\text{(r))}\right)  :$

\qquad\qquad d$\left(  \tau_{\text{0,r}}^{\text{t}}\text{e}_{\text{r}%
}^{\text{t}}\phi_{\text{t-r}}\text{(x}^{\text{t}}\text{(r))}\right)  =$
d$\left(  \tau_{\text{0,r}}^{\text{t}}e_{\text{r}}^{\text{t}}\right)
\phi_{\text{t-r}}$(x$^{\text{t}}$(r)) $+$ $\left(  \tau_{\text{0,r}}%
^{\text{t}}e_{\text{r}}^{\text{t}}\right)  $d$\phi_{\text{t-r}}$(x$^{\text{t}%
}$(r)) $+\frac{1}{2}$d$\left(  \tau_{\text{0,r}}^{\text{t}}e_{\text{r}%
}^{\text{t}}\right)  $d$\phi_{\text{t-r}}$(x$^{\text{t}}$(r))

By $\left(  4.17\right)  $ we have: \qquad d$\left\{  \tau_{\text{0,r}%
}^{\text{t}}\text{e}_{\text{r}}^{\text{t}}\right\}  =\frac{1}{2}%
(\tau_{\text{0,r}}^{\text{t}}$e$_{\text{r}}^{\text{t}})($W$_{\text{x}%
^{\text{t}}\text{(r)}})$drTherefore,

\qquad\qquad$\frac{1}{2}$d$\left(  \tau_{\text{0,r}}^{\text{t}}e_{\text{r}%
}^{\text{t}}\right)  $d$\phi_{\text{t-r}}$(x$^{\text{t}}$(r)) $=\frac{1}%
{4}(\tau_{\text{0,r}}^{\text{t}}e_{\text{r}}^{\text{t}})($W$_{\text{x}%
^{\text{t}}\text{(r)}})$d$\phi_{\text{t-r}}$(x$^{\text{t}}$(r))dr = 0

and hence,

\qquad\qquad d$\left(  \tau_{\text{0,r}}^{\text{t}}\text{e}_{\text{r}%
}^{\text{t}}\phi_{\text{t-r}}\text{(x}^{\text{t}}\text{(r))}\right)  =$
$\frac{1}{2}(\tau_{\text{0,r}}^{\text{t}}$e$_{\text{r}}^{\text{t}}%
)($W$_{\text{x}^{\text{t}}\text{(r)}})\phi_{\text{t-r}}$(x$^{\text{t}}$(r))dr
$+$ $\left(  \tau_{\text{0,r}}^{\text{t}}e_{\text{r}}^{\text{t}}\right)
$d$\phi_{\text{t-r}}$(x$^{\text{t}}$(r))

We insert the RHS of the last equation above on the RHS of $\left(
4.22\right)  $ and have:

$\left(  4.23\right)  \qquad$d$\left(  \tau_{\text{0,r}}^{\text{t}}%
\text{e}_{\text{r}}^{\text{t}}\text{h(y(r))}\right)  $ $=$ $\frac{1}{2}%
(\tau_{\text{0,r}}^{\text{t}}e_{\text{r}}^{\text{t}})($W$_{\text{x}^{\text{t}%
}\text{(r)}})\phi_{\text{t-r}}$(x$^{\text{t}}$(r))w(r)dr $+$ $\left(
\tau_{\text{0,r}}^{\text{t}}e_{\text{r}}^{\text{t}}\right)  $d$\phi
_{\text{t-r}}$(x$^{\text{t}}$(r))w(r)

\qquad\qquad\qquad\qquad\qquad\qquad\qquad$+$ $(\tau_{\text{0,r}}^{\text{t}}%
$e$_{\text{r}}^{\text{t}})\phi_{\text{t-r}}$(x$^{\text{t}}$(r))$\frac
{\text{L}^{0}\Psi}{\Psi}$(x$^{\text{t}}$(r))w(r)dr

We now compute d$\phi_{\text{t-r}}$(x$^{\text{t}}$(r)):

Since the Brownian Riemannian bridge process (x$^{\text{t}}$(s)) 0$\leq$
s$\leq$t$\Lambda\zeta$ has for associated differential generator

L$_{t,s}=\frac{1}{2}\Delta+$ X + $\nabla^{0}\log$q$_{\text{t-s}}$(-,P) + V, we
have by It\^{o}'s formula (differential version):

\qquad d$\phi_{\text{t-r}}$(x$^{\text{t}}$(r)) = $\left[  \frac{\partial
\phi_{\text{t-r}}}{\partial\text{r}}\text{\ +}\frac{1}{2}\Delta\phi
_{\text{t-s}}\text{ +
$<$%
X,}\nabla\phi_{\text{t-r}}\text{%
$>$
+
$<$%
}\nabla^{0}\log\text{q}_{\text{t-r}}\text{(-,P),}\nabla\phi_{\text{t-r}}\text{%
$>$%
}\right]  $(x$^{\text{t}}$(r))dr

\qquad\qquad\qquad\qquad$+$ $<\nabla\phi_{\text{t-r}}$(x$^{\text{t}}%
$(r)),u$_{\text{r}}$dB$_{\text{r}}>$

The expression in $\left(  4.23\right)  $ becomes:

\qquad d$\left(  \tau_{\text{0,r}}^{\text{t}}\text{e}_{\text{r}}^{\text{t}%
}\text{h(y(r))}\right)  $

$=$ $\left(  \tau_{\text{0,r}}^{\text{t}}\text{e}_{\text{r}}^{\text{t}%
}\right)  \left[  \frac{\partial\phi_{\text{t-r}}}{\partial\text{r}}\text{\ +
}\frac{1}{2}\Delta\phi_{\text{t-s}}+\text{ }\frac{1}{2}(\text{W}%
_{\text{x}^{\text{t}}\text{(r)}})\phi_{\text{t-r}}\text{+
$<$%
X,}\nabla\phi_{\text{t-r}}\text{%
$>$
+
$<$%
}\nabla^{0}\log\text{q}_{\text{t-r}}\text{(-,P),}\nabla\phi_{\text{t-r}}\text{%
$>$%
}\right]  $(x$^{\text{t}}$(r))w(r)dr

$+\left(  \tau_{\text{0,r}}^{\text{t}}\text{e}_{\text{r}}^{\text{t}}\right)  $
$<\nabla\phi_{\text{t-r}}$(x$^{\text{t}}$(r)),u$_{\text{r}}$dB$_{\text{r}}%
>$w(r) + $\tau_{\text{0,r}}^{\text{t}}$e$_{\text{r}}^{\text{t}}\phi
_{\text{t-r}}$(x$^{\text{t}}$(r))$\frac{\text{L}^{0}\Psi}{\Psi}$(x$^{\text{t}%
}$(r))w(r)dr

$=$ $\left(  \tau_{\text{0,r}}^{\text{t}}e_{\text{r}}^{\text{t}}\right)
\left[  \frac{\partial\phi_{\text{t-r}}}{\partial\text{r}}\text{\ + }\frac
{1}{2}\Delta\phi_{\text{t-s}}+\text{
$<$%
X,}\nabla\phi_{\text{t-r}}\text{%
$>$
+
$<$%
}\nabla^{0}\log\text{q}_{\text{t-r}}\text{(-,P),}\nabla\phi_{\text{t-r}}\text{%
$>$
+ }\frac{\text{L}^{0}\Psi}{\Psi}\phi_{\text{t-r}}\right]  $(x$^{\text{t}}$(r))w(r)dr

$+\left(  \tau_{\text{0,r}}^{\text{t}}e_{\text{r}}^{\text{t}}\right)  $
$<\nabla\phi_{\text{t-r}}$(x$^{\text{t}}$(r)),u$_{\text{r}}$dB$_{\text{r}}>$w(r)

where $\Delta=\Delta_{0}+$ $W$ is the \textbf{Laplace-Type }operator on vector bundles.

Since L = $\frac{1}{2}\Delta+X+V,$ the equality in $\left(  4.22\right)  $ becomes:

$\left(  4.24\right)  $\qquad d$\left(  \tau_{\text{0,r}}^{\text{t}}%
\text{e}_{\text{r}}^{\text{t}}\text{h(y(r))}\right)  =$ $\left(
\tau_{\text{0,r}}^{\text{t}}e_{\text{r}}^{\text{t}}\right)  \left[
\frac{\partial\phi_{\text{t-r}}}{\partial\text{r}}\text{\ + L + }\nabla
^{0}\log\text{q}_{\text{t-r}}\text{(-,P) + }\frac{\text{L}^{0}\Psi}{\Psi
}-V\right]  \phi_{\text{t-r}}$(x$^{\text{t}}$(r))w(r)dr

\qquad\qquad\qquad\qquad$\qquad\qquad+$ $\left(  \tau_{\text{0,r}}^{\text{t}%
}e_{\text{r}}^{\text{t}}\right)  $w(r) $<\nabla\phi_{\text{t-r}}$%
(x$^{\text{t}}$(r)),u$_{\text{r}}$dB$_{\text{r}}>$

Since $\phi_{\text{t-r}}=$ Q( t-r,t-s)$\phi,$ we set $\lambda=$ t-r and have
by \textbf{Lemma }$2$ above:

$\qquad\frac{\partial\phi_{\text{t-r}}}{\partial\text{r}}=\frac{\partial
\phi_{\lambda}}{\partial r}=\frac{\partial}{\partial\lambda}$(Q($\lambda
$,t-s)$\phi$)$\frac{\partial\lambda}{\partial r}$ $=-\frac{\partial}%
{\partial\lambda}($Q($\lambda$,t-s)$\phi)$

\qquad\qquad\ \ = $-$ $\left[  \text{L + }\nabla^{0}\log\text{q}_{\lambda
}\text{(-,P) + }\frac{\text{L}^{0}\Psi}{\Psi}-\text{V}\right]  $(Q($\lambda
$,t-s)$\phi$)

\qquad\qquad\ \ = $-$ $\left[  \text{L + }\nabla^{0}\log\text{q}_{\text{t-r}%
}\text{(-,P) + }\frac{\text{L}^{0}\Psi}{\Psi}-\text{V}\right]  $%
(Q(t-r,t-s)$\phi$)

Since $\phi_{\text{t-r}}=$ (Q(t-r,t-s)$\phi$), we have:

$\left(  4.25\right)  \qquad\frac{\partial\phi_{\text{t-r}}}{\partial\text{r}%
}=-$ $\left(  \text{L}+\nabla^{0}\log\text{q}_{\text{t-r}}\text{(-,P)}%
+\frac{\text{L}^{0}\Psi}{\Psi}-\text{V}\right)  \phi_{\text{t-r}}$

Inserting the RHS of $\left(  4.25\right)  $ in $\left(  4.24\right)  $ we see
that the RHS of $\left(  4.24\right)  $ is almost all wiped off and we have:

$\qquad\qquad\qquad\qquad\qquad$d$\left(  \tau_{\text{0,r}}^{\text{t}}%
\text{e}_{\text{r}}^{\text{t}}\text{h(y(r))}\right)  =\left(  \tau
_{\text{0,r}}^{\text{t}}e_{\text{r}}^{\text{t}}\right)  $w(r) $<\nabla
\phi_{\text{t-r}}$(x$^{\text{t}}$(r)),u$_{\text{r}}$dB$_{\text{r}}>$

Integrating both sides of the last equation above we have:

$\left(  4.26\right)  $\qquad$\tau_{\text{0,s}}^{\text{t}}$e$_{\text{s}%
}^{\text{t}}$h(y(s)) $-$ $\tau_{\text{0,0}}^{\text{t}}$e$_{\text{0}}%
^{\text{t}}$h(y(0)) $=\int_{0}^{\text{s}}\left(  \tau_{\text{0,r}}^{\text{t}%
}e_{\text{r}}^{\text{t}}\right)  $w(r) $<\nabla\phi_{\text{t-r}}$%
(x$^{\text{t}}$(r)),u$_{\text{r}}$dB$_{\text{r}}>$

Since h(y(s)) = h(t-s,x$^{\text{t}}$(s),w(s)) = $\phi_{\text{t}-\text{s}}%
$(x$^{\text{t}}$(s))w(s), the equation in $\left(  4.26\right)  $ becomes:

$\left(  4.27\right)  \qquad\tau_{\text{0,s}}^{\text{t}}$e$_{\text{s}%
}^{\text{t}}\phi_{\text{t-s}}$(x$^{\text{t}}$(s))w(s) $=\phi_{\text{t}}%
$(x$^{\text{t}}$(0))w(0) + M$_{\text{s}}$ $=\phi_{\text{t}}$(x) +
M$_{\text{s}}$

where,

\qquad\qquad M$_{\text{s}}$ = $\int_{0}^{\text{s}}\left(  \tau_{\text{0,r}%
}^{\text{t}}e_{\text{r}}^{\text{t}}\right)  $w(r) $<\nabla\phi_{\text{t-r}}%
$(x$^{\text{t}}$(r)),u$_{\text{r}}$dB$_{\text{r}}>$

is a \textbf{local martingale} and,

$\qquad\qquad\phi_{\text{t}}($x$^{\text{t}}($0)).w(0) $=$ (Q( t,t-s)$\phi
$)(x$^{\text{t}}$(0))w$(0)=$ (Q( t,t-s)$\phi$)(x)

We recall that by definition,

w(s) = exp$\left\{  \int_{0}^{\text{s}}\frac{\text{L}^{0}\Psi}{\Psi}%
\text{(x}^{\text{t}}\text{(r))dr}\right\}  $ and since since Q(t-s,t-s) is an
identity operator $\phi_{\text{t-s}}=$ Q(t-s,t-s)$\phi=\phi.$ Further since we
have by definition $\phi_{\text{t}}=$ Q( t,t-s)$\phi,$ we can then re-write
$\left(  4.27\right)  $ above, for $0\leq r\leq s\leq$t$\Lambda\zeta$ as follows:

$\left(  4.28\right)  \qquad\tau_{\text{0,s}}^{\text{t}}$e$_{\text{s}%
}^{\text{t}}\phi$(x$^{\text{t}}$(s))exp$\left\{  \int_{0}^{\text{s}}%
\frac{\text{L}^{0}\Psi}{\Psi}\text{(x}^{\text{t}}\text{(r))dr}\right\}  =($Q(
t,t-s)$\phi)$(x) + M$_{\text{s}}$

We then take expectations on both sides of $\left(  4.28\right)  $ to have:

$\left(  4.29\right)  $ \qquad\textbf{E}$_{\text{x}}\left[  \chi
_{\zeta>\text{s}}\tau_{\text{0,s}}^{\text{t}}\text{e}_{\text{s}}^{\text{t}%
}\phi\text{(x}^{\text{t}}\text{(s))exp}\left\{  \int_{0}^{\text{s}}%
\frac{\text{L}^{0}\Psi}{\Psi}\text{(x}^{\text{t}}\text{(r))dr}\right\}
\right]  =$ (Q( t,t-s)$\phi)($x)\ $+$ \textbf{E}$_{\text{x}}($M$_{\text{s}}).$

We need to show that \textbf{E}$_{\text{x}}($M$_{\text{s}})=0:$

By \textbf{Proposition 1.1 }of\textbf{ Ikeda }and\textbf{ Watanabe }$\left[
32\right]  ,$ M$_{\text{s}}=$ $\int_{0}^{\text{s}}\left(  \tau_{\text{0,r}%
}^{\text{t}}e_{\text{r}}^{\text{t}}\right)  $w(r) $<\nabla\phi_{\text{t-r}}%
$(x$^{\text{t}}$(r)),u$_{\text{r}}$dB$_{\text{r}}>$ is a \textbf{martingale}.

Consequently E$_{\text{x}}$(M$_{\text{s}}$) = 0 and $\left(  4.29\right)  $
becomes for 0$\leq$ s$\leq$t$:$

$\left(  4.30\right)  $ \qquad(Q( t,t-s)$\phi)($x)\ $=$ \textbf{E}$_{\text{x}%
}\left[  \chi_{\zeta>\text{s}}\tau_{\text{0,s}}^{\text{t}}\text{e}_{\text{s}%
}^{\text{t}}\phi\text{(x}^{\text{t}}\text{(s))exp}\left\{  \int_{0}^{\text{s}%
}\frac{\text{L}^{0}\Psi}{\Psi}\text{(x}^{\text{t}}\text{(r))dr}\right\}
\right]  $

where $\zeta$ is the first exit time from M$_{0}$ of the bridge process
x$^{\text{t}}$(s) 0$\leq$ s$\leq$t$\wedge\zeta.$

Since vol(M$_{0}$)\ = vol(M), $\zeta$ is the first exit time from the compact
Riemannian manifold M and so $\zeta=+\infty.$

Consequently, we have finally here:

$\bigskip\left(  4.31\right)  \qquad$(Q( t,t-s)$\phi)($x) = \textbf{E}%
$_{\text{x}}\left[  \tau_{\text{0,s}}^{\text{t}}\text{e}_{\text{s}}^{\text{t}%
}\phi\text{(x}^{\text{t}}\text{(s))exp}\left\{  \int_{0}^{\text{s}}%
\frac{\text{L}^{0}\Psi}{\Psi}\text{(x}^{\text{t}}\text{(r))dr}\right\}
\right]  $

The expressions in $\left(  4.30\right)  $ and $\left(  4.31\right)  $ give
the more \textbf{Generalized} \textbf{Feynman-Kac} formula in a vector bundle.
The Corollaries below make this more explicit.

\qquad\qquad\qquad\qquad\qquad\qquad\qquad\qquad\qquad\qquad\qquad\qquad
\qquad\qquad\qquad\qquad\qquad$\qquad\blacksquare$

(i) The expression: M$_{\text{s}}$ = $\int_{0}^{\text{s}}\left(
\tau_{\text{0,r}}^{\text{t}}e_{\text{r}}^{\text{t}}\right)  \exp\left\{
\int_{0}^{\text{s}}\frac{\text{L}^{0}\Psi}{\Psi}\text{(x}^{\text{t}%
}\text{(r))dr}\right\}  <\nabla\phi_{\text{t-r}}$(x$^{\text{t}}$%
(r)),u$_{\text{r}}$dB$_{\text{r}}>$ obtained in $\left(  4.27\right)  $ is
similar to the expression for R$_{s}$V$_{s}$(x$_{s}$) in $\left(  2.21\right)
$ of \textbf{Bismut} $\left[  8\right]  .$

(ii) A similar expression also showed up in the proof of \textbf{Theorem
}$\left(  34\right)  $ in \textbf{Norris }$\left[  45\right]  $ and was proved
to be a martingale.

(iii) Again a similar expression showed up in proof of \textbf{Theorem
}$\left(  7.2.1\right)  $ of \textbf{Hsu }$\left[  30\right]  $ and was
assumed a martingale.

\qquad\qquad\qquad\qquad\qquad\qquad\qquad\qquad\qquad\qquad\qquad\qquad
\qquad\qquad\qquad\qquad\qquad$\ \ \ \ \ \ \blacksquare$

Our proof of the last theorem above is to be compared to proofs of similar
theorems in the following papers: \textbf{Theorem }$\left(  2.5\right)  $ of
\textbf{Bismut }$\left[  8\right]  ,$ \textbf{Proposition }$\left(
4.5\right)  $ of \textbf{Driver and Thalmaier} $\left[  12\right]  ,$
\textbf{Theorem }$(7.2.1)$ of\ \textbf{Hsu} $\left[  30\right]  $ and $\left(
34\right)  $ of \textbf{Norris} $\left[  45\right]  .$ Their theorems are more
adapted to obtaining the Feynman-Kac formula directly, which will be obtained
here as a special case of our theorem here. Even the "generalized Feynma-Kac
formula" obtained in $\left(  34\right)  $ of \textbf{Norris} $\left[
46\right]  $ is a special case of our theorem here as we shall see. The
theorem here is thus the ultimate generalization of the Feynman-Kac formula.

\qquad\qquad\qquad\qquad\qquad\qquad\qquad\qquad\qquad\qquad\qquad\qquad
\qquad\qquad\qquad\qquad\qquad$\ \ \ \ \ \blacksquare$

An analogue of the theorem below was proved in the case of the \textbf{scalar
heat kernel} in \textbf{Theorem }$\left(  4.4\right)  $ of \textbf{Ndumu}
$\left[  42\right]  .$ It is a\textbf{ Generalized Feynman-Kac}
\textbf{formula} from which we shall deduce the usual \textbf{Feynman-Kac
formula} as well as a stochastic representation of the \textbf{Generalized
Elworthy-Truman heat kernel formula}, and ultimately the \textbf{heat kernel
expansion formula.}

We now deduce a more explicit form of the generalized Feynman-Kac Formula and
the usual Feynman-Kac formula:\qquad\qquad\qquad\qquad\qquad\qquad\qquad
\qquad\qquad\qquad\qquad\qquad\qquad\qquad\qquad\qquad\qquad

\begin{corollary}
(The Generalized Feynman-Kac Formula as a Solution of the Heat Equation)
\end{corollary}

For $\phi\in\Gamma(E)$ and 0$\leq$ t$<+\infty,$ when P = M$_{0}$, we have the
\textbf{Generalized Feynman-Kac Formula:}

For almost all x$\in M,$

(i)\qquad\ P$_{\text{s}}^{M}\phi$(x) =$\ \ ($Q(t,t-s)$\phi)$(x) $=$
$\int_{\text{M}}$k$_{\text{s}}$(x,z)$\phi$(z)$\upsilon_{\text{M}}$(dz) =
\textbf{E}$_{\text{x}}\left[  \tau_{\text{0,s}}\text{e}_{\text{s}}%
\phi\text{(x(s))}\phi\text{(x(s))exp}\left\{  \int_{0}^{\text{s}%
}V\text{(x(r))dr}\right\}  \right]  $

where P$_{\text{t}}^{M}$ is the semigroup operator on $\Gamma(E)$ and x(s)
0$\leq$ s$<+\infty$ is \textbf{Brownian motion} with drift X and potential
term V on the compact Riemannian manifold M. It has for associated generator:
$L=\frac{1}{2}\Delta+X+V.$

(ii)\qquad We deduce the usual \textbf{Feynman-Kac Formula }in vector bundles:

\qquad\qquad P$_{\text{t}}^{M}\phi(x)=$ $\int_{\text{M}}$k$_{\text{t}}%
$(x,z)$\phi$(z)$\upsilon_{\text{M}}$(dz) $=$ \textbf{E}$_{\text{x}}\left[
\tau_{\text{0,t}}\text{e}_{\text{t}}\phi\text{(x(t))}\exp\left\{  \int%
_{0}^{\text{t}}\text{V(x(r))dr}\right\}  \right]  $

(iii)\qquad Define $\phi_{s}=\ ($Q(t,t-s)$\phi)$ and let $L=\frac{1}{2}%
\Delta+X+V$

\qquad Then $\phi_{s}$ is a solution of the \textbf{Cauchy Problem for the
Heat Equation:}

\qquad For almost all $x\in M,$ we have:

\qquad$\frac{\partial\phi_{s}}{\partial s}(x)=L\phi_{s}(x)$

$\qquad\ \ \phi_{0}(x)=\phi(x)$

\begin{proof}
(i) Recall that, by definition,
\end{proof}

$\left(  4.32\right)  \qquad($Q(t,t-s)$\phi)$(x) = q$_{\text{t}}$(x,P)$^{-1}%
$P$_{\text{s}}$(q$_{\text{t-s}}$(-,P)$\phi$)(x$_{0}$)

\qquad\qquad\qquad= q$_{\text{t}}$(x,P$)^{-1}\int_{\text{M}}$q$_{\text{t-s}}%
$(z,P)k$_{\text{s}}$(x,z)$\phi$(z)$\upsilon_{\text{M}}$(dz)

where,

\qquad\qquad\qquad q$_{t}$(x,P) $=$ (2$\pi$t)$^{-\frac{n-q}{2}}\Psi
(x)\exp\left\{  -\frac{\text{d(x,P)}^{2}}{2t}\right\}  $

Let P = M$_{0}.$ This is equivalent to the \textbf{Fermi coordinates} become
\textbf{local coordinates} based at y$_{0}\in$P = M$_{0}.$

Consequently,

\qquad\qquad\qquad\qquad\qquad$\Psi(x)=1$ and q = $\dim$P = $\dim$M$_{0}=n$

\qquad We conclude that:

$\left(  4.33\right)  $\qquad\qquad\qquad\qquad\qquad\ \ q$_{\text{t-s}}$(z,P)
= q$_{\text{t-s}}$(z,M$_{0}$) = 1 = q$_{\text{t}}$(x,P$)$

and,

$\left(  4.34\right)  $\qquad\qquad\qquad\qquad$\qquad\frac{\text{L}^{0}\Psi
}{\Psi}(x)=V(x)$\qquad\ 

Consequently for all z$\in$M and 0$\leq s\leq t<+\infty,$

$\left(  4.35\right)  \qquad\qquad\qquad\qquad\qquad\qquad\qquad\nabla^{0}%
$q$_{\text{t-s}}$(z,M$_{0}$) = 0

By $\left(  4.33\right)  $ above, the local coordinates (y$_{\text{1}}%
$,...,y$_{\text{q}}$) at the point y$_{0}\in$P defined in $\left(  1.1\right)
$ of Chapter 1 are now extended to the \textbf{local coordinates}
(y$_{\text{1}}$,...,y$_{\text{n}}$) at y$_{0}\in$P = M$_{0}.$

In this special case, we have by $\left(  4.34\right)  :$

$\left(  4.36\right)  $\qquad\qquad(Q(t,t-s) $\phi$)(x) \ = q$_{\text{t}}%
$(x,M$_{0}$)$^{-1}\int_{\text{M}}$q$_{\text{t-s}}$(z,M$_{0}$)k$_{\text{s}}%
$(x,z)$\phi$(z)$\upsilon_{\text{M}}$(dz)

\ \qquad\qquad\qquad\qquad\qquad\qquad\ \ \ \ \ = $\int_{\text{M}_{0}}%
$k$_{\text{s}}$(x,z)$\phi$(z)$\upsilon_{\text{M}}$(dz) = $\int_{\text{M}}%
$k$_{\text{s}}$(x,z)$\phi$(z)$\upsilon_{\text{M}}$(dz)

The last equality above is due to the fact that volume(M$_{0}$)\ = volume(M)
by \textbf{Theorem }$\left(  8.40\right)  $ of \textbf{Gray }$\left[
25\right]  .$ We remark here that there is a slight error in \textbf{Gray
}$\left[  25\right]  :$ In it, volume(M) = volume($\Theta_{P})$ instead of
volume(M) = volume(exp$_{\upsilon}(\Theta_{P}))=$ volume(M$_{0}$)$\ .$

By $\left(  4.35\right)  ,$ the differential generator L$_{\text{s,t}}$ =
$\frac{1}{2}\Delta$ $+$ X $+$ $\nabla^{0}\log$q$_{\text{t-s}}$(-,P) + V of the
\textbf{semi-classical Brownian Riemannian bridge process} (x$_{s}^{\text{t}}%
$) $0\leq$s$\leq$t$\wedge\zeta$ now reduces to L = $\frac{1}{2}\Delta$ $+$ X +
V. Consequently, the bridge process x$^{\text{t}}$(s) reduces to Brownian
motion x(s) 0$\leq$ s$<+\infty$ having for differential generator L $=\frac
{1}{2}\Delta$ $+$ X $+$ V and with life-time $\zeta=+\infty$ on the compact
Riemannian manifold M. We note that the compactness of M gives $\zeta
=+\infty.$ In this case, e$_{\text{s}}^{\text{t}}$ no longer depends on t%
$>$%
0 and so we denote it by e$_{\text{s}}$ and see that the expression for
$($Q(t,t-s)$\phi)$(x) is independent of t and so we set:

$\left(  4.37\right)  $\qquad$\qquad\qquad\qquad\qquad$P$_{\text{s}}^{M}%
\phi(x)=($Q(t,t-s)$\phi)$(x) = $\int_{\text{M}}$k$_{\text{s}}$(x,z)$\phi
$(z)$\upsilon_{\text{M}}$(dz)

In the present situation, the parallel translations on fibres of the vector
bundle E in $\left(  4.18\right)  $ become the usual translations along the
Brownian motion $(x(s))$ $0\leq r\leq s\leq t<+\infty$ with differential
generator L $=\frac{1}{2}\Delta$ $+$ X $+$ V :

\begin{center}
$\qquad\tau_{\text{r,s}}=$ u$_{\text{r}}$(u$_{\text{s}}$)$^{-1}:$%
E$_{\text{x(s)}}\longrightarrow$ E$_{_{\text{x(r)}}}$ with inverse
$\tau_{\text{s,r}}=$ u$_{\text{s}}$(u$_{\text{r}}$)$^{-1}:$E$_{\text{x(r)}%
}\longrightarrow$ E$_{_{\text{x(s)}}}$\ 
\end{center}

By $\left(  4.20\right)  $ above, we have for $0\leq r\leq s<t<+\infty:$\ \ \ 

$\left(  4.38\right)  $\qquad\qquad\qquad\qquad e$_{\text{s}}\ =$
$\tau_{\text{s,0}}+\frac{1}{2}\int_{0}^{\text{s}}\tau_{\text{s,r}}%
$e$_{\text{r}}$W$_{\text{x(r)}}$dr\ 

We conclude from \textbf{Theorem 1} above and $\left(  4.33\right)  ,$
$\left(  3.34\right)  ,$ $\left(  4.36\right)  ,$ $\left(  4.37\right)  $ and
$\left(  4.38\right)  $ that:

$\left(  4.39\right)  $\qquad\ \textbf{E}$_{\text{x}}\left[  \tau_{\text{0,s}%
}\text{e}_{\text{s}}\phi\text{(x(s))}\phi\text{(x(s))exp}\left\{  \int%
_{0}^{\text{s}}V\text{(x(r))dr}\right\}  \right]  $ = $($Q(t,t-s)$\phi)$(x)
$=$ $\int_{\text{M}}$k$_{\text{s}}$(x,z)$\phi$(z)$\upsilon_{\text{M}}$(dz)

\qquad\qquad\qquad= P$_{\text{s}}^{M}\phi$(x)

The last equation is just notation.

(ii) We take limits as s$\uparrow$t on all sides of $\left(  4.39\right)  $
for 0$\leq$ s$\leq$t$:$

Since the manifold M is compact we can take limits under the expectation sign
and have for $0\leq s\leq t<+\infty:$

$\left(  4.40\right)  $\qquad P$_{\text{t}}^{M}\phi$(x)\ $=$ $\int_{\text{M}}%
$k$_{\text{t}}$(x,z)$\phi$(z)$\upsilon_{\text{M}}$(dz) $\qquad$

\qquad$\qquad=$ $\underset{\text{s}\uparrow\text{t}}{lim}($Q(t,t-s)$\phi)$(x)
= $\underset{\text{s}\uparrow\text{t}}{lim}$\textbf{E}$_{\text{x}}\left[
\tau_{\text{0,s}}\text{e}_{\text{s}}\phi\text{(x(s))}\exp\left\{  \int%
_{0}^{\text{s}}V\text{(x(r))dr}\right\}  \right]  $

\qquad\qquad\qquad\qquad\qquad\qquad\ \ \ \ \ \ \ = \textbf{E}$_{\text{x}%
}\left[  \tau_{\text{0,t}}\text{e}_{\text{t}}\phi\text{(x(t))}\exp\left\{
\int_{0}^{\text{t}}V\text{(x(r))dr}\right\}  \right]  $

We conclude from the last equations that:

$\left(  4.41\right)  \qquad$P$_{\text{t}}^{M}\phi$(x) = $\int_{\text{M}}%
$k$_{\text{t}}$(x,z)$\phi$(z)$\upsilon_{\text{M}}$(dz) $=$ \textbf{E}%
$_{\text{x}}\left[  \tau_{\text{0,t}}\text{e}_{\text{t}}\phi\text{(x(t))}%
\exp\left\{  \int_{0}^{\text{t}}V\text{(x(r))dr}\right\}  \right]  $

(iii) Setting $\phi_{s}=\ ($Q(t,t-s)$\phi),$ we have from $\left(
4.37\right)  $ above:

\qquad$\frac{\partial\phi_{s}}{\partial s}(x)=\frac{\partial}{\partial s}($
P$_{\text{s}}^{M}\phi)(x)=$ $\frac{\partial}{\partial s}\int_{\text{M}}%
$k$_{\text{s}}$(x,z)$\phi$(z)$\upsilon_{\text{M}}$(dz)

Since the manifold M is compact, the differentiation sign can go over the
integral sign and we have:

$\left(  4.42\right)  $\qquad$\frac{\partial\phi_{s}}{\partial s}%
(x)=\frac{\partial}{\partial s}($ P$_{\text{s}}^{M}\phi)(x)=$ $\frac{\partial
}{\partial s}\int_{\text{M}}$k$_{\text{s}}$(x,z)$\phi$(z)$\upsilon_{\text{M}}%
$(dz) $=$ $\int_{\text{M}}\frac{\partial}{\partial s}$k$_{\text{s}}$%
(x,z)$\phi$(z)$\upsilon_{\text{M}}$(dz)

The heat kernel is the fundamental solution of the heat equation and hence,
for L = $\frac{1}{2}\Delta$ $+$ X + V, we have:

$\qquad\frac{\partial}{\partial s}$k$_{\text{s}}$(x,-) = $L$k$_{\text{s}}$(x,-)

Therefore the equations in $\left(  4.42\right)  $ can be further extended:

\qquad$\frac{\partial\phi_{s}}{\partial s}(x)=\frac{\partial}{\partial s}($
P$_{\text{s}}^{M}\phi)(x)=$ $\frac{\partial}{\partial s}\int_{\text{M}}%
$k$_{\text{s}}$(x,z)$\phi$(z)$\upsilon_{\text{M}}$(dz) $=$ $\int_{\text{M}%
}\frac{\partial}{\partial s}$k$_{\text{s}}$(x,z)$\phi$(z)$\upsilon_{\text{M}}$(dz)

\qquad\ $=\int_{\text{M}}L$k$_{\text{s}}$(x,y)$\phi$(z)$\upsilon_{\text{M}}%
$(dz) $=L\int_{\text{M}}$k$_{\text{s}}$(x,y)$\phi$(z)$\upsilon_{\text{M}}$(dz)

The last equation is due to $\left(  2.5\right)  $ and $\left(  2.6\right)  $
of \textbf{Azencott }$\left[  4\right]  :$ The differential operator $L$ and
the integral sign $\int$ can be inter-changed. We equate the first and last
expressions above and have:

$\qquad\frac{\partial\phi_{s}}{\partial s}(x)=L\int_{\text{M}}$k$_{\text{s}}%
$(x,y)$\phi$(z)$\upsilon_{\text{M}}$(dz) $=L(P_{\text{s}}^{M}\phi)(x)$

By notation: $(P_{\text{s}}^{M}\phi)(x)=\int_{\text{M}}$k$_{\text{s}}%
$(x,y)$\phi$(z)$\upsilon_{\text{M}}$(dz)

Consequently, the last two equalities give:

$\qquad\frac{\partial\phi_{s}}{\partial s}(x)=$ $L(P_{\text{s}}^{M}%
\phi)(x)=L\phi_{s}(x)$

The partial differential equation is thus obtained, where:

$\qquad\phi_{s}=\ ($Q(t,t-s)$\phi)$

and since $\ (Q(t,t)$ is the identity operator, we see that:$\qquad$,

\qquad\ $\phi_{0}=\ (Q(t,t)\phi)=\phi$ $\qquad\qquad\qquad\qquad\qquad
\qquad\qquad\qquad\qquad\qquad\qquad\qquad$

So the initial condition is obtained and (iii) is proved.$\qquad\qquad
\qquad\qquad\qquad\qquad\qquad\qquad\qquad\qquad\qquad\qquad\qquad\qquad
\qquad\qquad\qquad\qquad\qquad$

$\qquad\qquad\qquad\qquad\qquad\qquad\qquad\qquad\qquad\qquad\qquad
\qquad\qquad\qquad\qquad\qquad\qquad\qquad\blacksquare$

\begin{remark}
The entire Corollary and the equalites in $\left(  4.41\right)  $ above give a
further justification to the fact that the expression for Q(t,t-s)$\phi$(x) in
\textbf{Theorem 1} above is a generalized Feynman-Kac formula.
\end{remark}

\qquad\qquad\qquad\qquad\qquad\qquad\qquad\qquad\qquad\qquad\qquad\qquad
\qquad\qquad\qquad\qquad\qquad\qquad$\blacksquare$

Compare the result in the last Corollary above with $\left(  4.14\right)  $ of
\textbf{Driver and Thalmaier }$\left[  12\right]  ,$ \textbf{Theorem 7.2.1} of
\textbf{Hsu }$\left[  30\right]  $ and $\left(  34\right)  $ of
\ \textbf{Norris }$\left[  45\right]  .$

\qquad\qquad\qquad\qquad\qquad\qquad\qquad\qquad\qquad\qquad\qquad\qquad
\qquad\qquad\qquad\qquad\qquad\qquad$\blacksquare$

We have thus obtained the usual \textbf{Feynman-Kac Formula} here as a special
case of a more general theorem. This is the reason why we called the (more
general) theorem above, the \textbf{Generalized Feyman-Kac Formula}%
.\qquad\qquad\qquad\qquad\qquad\qquad\qquad\qquad\qquad

The formula in $\left(  4.41\right)  $ above is the well known Feynman-Kac
formula. We compare this with the following:

1. The formula in Theorem $\left(  2.5\right)  $ of \textbf{Bismut} $\left[
8\right]  ,$ where the techniques of proof seem more aligned with those used
in proving the general formula in $\left(  4.31\right)  ,$ is similar to our
Feynman-Kac Formula here in (ii) of the last Corollary or given here in
$\left(  4.41\right)  $ above.

2. The formula for the Feynman-Kac formula of \textbf{Proposition 4.5} in
$\left(  4.14\right)  $ of \textbf{Driver and Thalmaier} $\left[  12\right]  $
(without the potential term)\ is similar to our formula here in (ii) of the
Corollary or $\left(  4.41\right)  $ above.

3. \textbf{Theorem }$\left(  7.2.1\right)  $ of \textbf{Hsu }$\left[
30\right]  ,$ on differential forms, is the same as our formula here (but the
potential term is absent).

4. The Feynman-Kac formula in $\left(  34\right)  $ of \textbf{Norris
}$\left[  45\right]  $ is the same as the formula here except that the
potential term is absent there.

\qquad\qquad\qquad\qquad\qquad\qquad\qquad\qquad\qquad\qquad\qquad\qquad
\qquad\qquad\qquad\qquad\qquad\qquad\qquad$\blacksquare$

\begin{corollary}
(The Generalized Elworthy-Truman Heat Kernel Formula)
\end{corollary}

For $\phi\in\Gamma(E),$

The vector bundle Generalized Heat Kernel k$_{\text{t}}$(x,P,$\phi$) has a
\textbf{deterministic} and a \textbf{stochastic} representations as follows:

(i) \qquad\qquad\qquad k$_{\text{t}}$(x,P,$\phi$) $=$ $\int_{\text{P}}%
$k$_{\text{t}}$(x,y)$\phi($y)$\upsilon_{\text{P}}$(dy) \qquad(deterministic)

(ii) \qquad\qquad\qquad k$_{\text{t}}$(x,P,$\phi$) $=$ q$_{\text{t}}$%
(x,P$)$\textbf{E}$_{\text{x}}\left(  \tau_{\text{0,t}}^{\text{t}}%
\text{e}_{\text{t}}^{\text{t}}\phi\text{(y(t))}\exp\left\{  \int_{0}%
^{\text{t}}\frac{\text{L}^{0}\Psi}{\Psi}\text{(x}^{\text{t}}\text{(r))dr}%
\right\}  \right)  $ \qquad(stochastic)

\begin{proof}
(i) The deterministic representation is given in $\left(  2.22\right)  $ above.
\end{proof}

(ii) For the stochastic representation we use \textbf{Lemma }$1$ and
\textbf{Theorem }$1$ to have:

\begin{center}
q$_{\text{t}}$(x,P)$^{-1}$ $\int_{\text{P}}$k$_{\text{t}}$(x,y)$\phi
(y)\upsilon_{\text{P}}$(dy) = $\underset{\text{s}\uparrow\text{t}}{lim}%
$(Q(t,t-s) $\phi$)(x) $\qquad$

$\qquad=$ $\underset{\text{s}\uparrow\text{t}}{lim}$\textbf{E}$_{\text{x}%
}\left(  \tau_{\text{0,s}}^{\text{t}}\text{e}_{\text{s}}^{\text{t}}%
\phi\text{(x}^{\text{t}}\text{(s))}\exp\left\{  \int_{0}^{\text{s}}%
\frac{\text{L}^{0}\Psi}{\Psi}\text{(x}^{\text{t}}\text{(r))dr}\right\}
\right)  $

\qquad$\ =$ \textbf{E}$_{\text{x}}\left(  \tau_{\text{0,t}}^{\text{t}}%
\text{e}_{\text{t}}^{\text{t}}\phi\text{(y(t))}\exp\left\{  \int_{0}%
^{\text{t}}\frac{\text{L}^{0}\Psi}{\Psi}\text{(x}^{\text{t}}\text{(r))dr}%
\right\}  \right)  $
\end{center}

The second equality is obvious by \textbf{Theorem 1.}

The last equality is due the fact that we can take limits under the
expectation sign since the manifold M is compact.

\qquad\qquad\qquad\qquad\qquad\qquad\qquad\qquad\qquad\qquad\qquad\qquad
\qquad\qquad\qquad\qquad\qquad$\blacksquare$

The results here in (i) generalize \textbf{Theorem} $\left(  4.8\right)  $ of
\textbf{Ndumu} $\left[  42\right]  .$

From the Generalized Heat Kernel formula above, we deduce the vector bundle
version of \ the \textbf{Elworthy-Truman heat kernel formula} relative to the
Generalized Laplacian.

\qquad\qquad\qquad\qquad\qquad\qquad\qquad\qquad\qquad\qquad\qquad\qquad
\qquad\qquad\qquad\qquad$\qquad\blacksquare$

\begin{corollary}
(The Elworthy-Truman Heat Kernel Formula for Vector Bundles)
\end{corollary}

\qquad If P = $\left\{  \text{y}_{0}\right\}  $ (this means that the
\textbf{Fermi coordinates} reduce to \textbf{normal coordinates} centered at
y$_{0}$), we have:

\qquad k$_{\text{t}}$(x,y$_{0},\phi$) = q$_{\text{t}}$(x,y$_{0})$%
\textbf{E}$_{\text{x}}\left(  \tau_{\text{y}_{0}\text{,x}}^{t}\text{e}%
_{\text{t}}^{\text{t}}\phi(\text{y}_{0})\exp\left\{  \int_{0}^{\text{t}}%
\frac{\text{L}^{0}\Psi}{\Psi}\text{(x}^{\text{t}}\text{(s))ds}\right\}
\right)  $

where (x$^{\text{t}}$(s)) 0$\leq$ s$\leq$t$\Lambda\zeta$ \ is the
semi-classical Brownian Riemannian bridge process starting from the point
x$\in$M$_{0}$ and reaching the center of the normal neighborhood y$_{0}\in
$M$_{0}$ in time t and $\tau_{\text{y}_{0}\text{,x}}^{t}=\tau_{\text{0,t}%
}^{\text{t}}:$E$_{\text{y}_{0}}\longrightarrow$ E$_{\text{x}}$ is the parallel
transport along the reversed semi-classical Riemannan Brownian bridge process
from y$_{0}$ to the point x$\in$M$_{0}$ in time t.

Here $\zeta$ is now the first exit time of the bridge process from the normal
neighbourhood M$_{0}$ of y$_{0}.$\qquad\qquad\qquad\qquad\qquad\qquad
\qquad\qquad\qquad\qquad\qquad\qquad\qquad\qquad$\qquad\qquad\qquad
\qquad\qquad\qquad\qquad\qquad\qquad\qquad\qquad\qquad\qquad\qquad\qquad
\qquad\qquad\blacksquare$

The above formula is a generalization of the usual Elworthy-Truman heat kernel
formula to the case of vector bundles. We obtain the usual Elworthy-Truman
heat kernel formula when the vector bundle is the trivial bundle E = M$\times$R.

\qquad\qquad\qquad\qquad\qquad\qquad\qquad\qquad\qquad\qquad\qquad\qquad
\qquad\qquad\qquad\qquad$\ \ \ \blacksquare$

\chapter{Partial Differential Equations Associated to the Feynman-Kac
Formula\qquad\qquad\qquad\qquad\qquad\qquad\qquad}

Here we begin the process of generalized heat kernel expansions in vector
bundles. First we derive a partial differential equation below (\textbf{Lemma
}$4$) which will be refined in \textbf{Theorem }$2$. It will play a central
role in deriving the expansion theorem (\textbf{Theorem }$\mathbf{3}$) below.

In order to avoid a plethora of the superscripts $^{"}0^{"}$ on $\nabla$ and
$\Delta,$ we shall from now henceforth drop the "0" on $\nabla$ and $\Delta$
when applied to C$^{\infty}$(M)-functions and write $\nabla f$ \ and $\Delta
f$ instead of $\nabla^{0}f$ \ and $\Delta^{0}f$ for f$\in$C$^{\infty}$(M).

\qquad\qquad\qquad\qquad\qquad\qquad\qquad\qquad\qquad\qquad\qquad\qquad
\qquad\qquad\qquad\qquad$\blacksquare$

\begin{theorem}
$\frac{\partial}{\partial\text{s}}($Q(t,t-s)$\phi$) = (Q(t,t-s)$\left[  \text{
L}\phi+\text{ }<\nabla^{0}\log\text{q}_{\text{t-s}}\text{(-,P),}\nabla
\phi\text{%
$>$%
) +}\frac{\text{L}^{0}\Psi}{\Psi}\phi\text{ }-\text{V}\phi\right]  $(x)
\end{theorem}

\qquad\qquad\qquad= (Q(t,t-s)$\left[  \text{ L}+\text{ }\nabla^{0}\log
\text{q}_{\text{t-s}}\text{(-,P) +}\frac{\text{L}^{0}\Psi}{\Psi}%
-\text{V}\right]  \phi$(x)

where $\phi$ is a section of the vector bundle E.

\begin{proof}
From the definition of $\frac{\partial}{\partial\text{s}}($Q(t,t-s)$\phi$) in
$\left(  4.8\right)  ,$ we have:\qquad\qquad\qquad
\end{proof}

\ $\frac{\partial}{\partial\text{s}}($Q(t,t-s)$\phi$)\ \ = q$_{\text{t}}%
$(x$_{0}$,P$)^{-1}\frac{\partial}{\partial\text{s}}\left[  \text{P}_{\text{s}%
}\text{(q}_{\text{t-s}}\text{(-,P)}\phi\text{))}\right]  \qquad$

Since M is compact, we can differentiate under the integral sign and have:

$\left(  5.2\right)  \qquad\frac{\partial}{\partial\text{s}}\left[
\text{P}_{\text{s}}\text{(q}_{\text{t-s}}\text{(-,P)}\phi\text{))}\right]
=\frac{\partial}{\partial\text{s}}\int_{\text{M}}\left[  \text{k}_{\text{s}%
}\text{(x,z)q}_{\text{t-s}}\text{(z,P)}\phi\text{(z)}\right]  \upsilon
_{\text{M}}$(dz)

\qquad\qquad$\qquad\qquad\qquad\qquad\ \ \ \ =\int_{\text{M}}\frac{\partial
}{\partial\text{s}}\left[  \text{k}_{\text{s}}\text{(x,z)q}_{\text{t-s}%
}\text{(z,P)}\phi\text{(z)}\right]  \upsilon_{\text{M}}$(dz)

\qquad\qquad\qquad\ = $\int_{\text{M}}\frac{\partial}{\partial\text{s}}%
$(q$_{\text{t-s}}$($z,$P))k$_{\text{s}}$(x,z)$\phi$(z)$\upsilon_{\text{M}}%
$(dz) + $\int_{\text{M}}$q$_{\text{t-s}}$($z,$P)$\frac{\partial}%
{\partial\text{s}}$k$_{\text{s}}$(x,z)$\phi$(z)$\upsilon_{\text{M}}$(dz)

By \textbf{Lemma 2 }and the fact that k$_{\text{s}}$(x$_{\text{0}}$,z) is the
fundamental solution of the heat equation on the vector bundle E, we have:

\qquad\ = $\int_{\text{M}}[-$L$^{0}$ q$_{\text{t-s}}(z,$P) + $\frac
{\text{L}^{0}\Psi}{\Psi}q_{\text{t-s}}(z,$P$)]$k$_{\text{s}}$(x$_{0}$,z)$\phi
$(z)$\upsilon_{\text{M}}$(dz) $\ $

$\qquad+\int_{\text{M}}$q$_{\text{t-s}}(z,$P)$.$Lk$_{\text{s}}$(x$_{0}%
$,z)$\phi$(z)$\upsilon_{\text{M}}$(dz)

\qquad\qquad\qquad= I$_{1}$\ + I$_{2}$

We set:

$\left(  5.3\right)  \qquad$I$_{1}$\ = $\int_{\text{M}}[-$L$^{0}$
q$_{\text{t-s}}$($z,$P) + $\frac{\text{L}^{0}\Psi}{\Psi}q_{\text{t-s}}%
$(z,P)$]\phi$(z)k$_{\text{s}}$(x,z)$\upsilon_{\text{M}}$(dz) $\ $

$\left(  5.4\right)  $\qquad I$_{2}=\int_{\text{M}}$q$_{\text{t-s}}$($z,$%
P)$.$Lk$_{\text{s}}$(x,z)$\phi$(z)$\upsilon_{\text{M}}$(dz)

By the definition of the operator P$_{\text{s}}$ we have:

$\left(  5.5\right)  \qquad$I$_{1}=$ P$_{\text{s}}[-$L$^{0}$(q$_{\text{t-s}%
}(-,$P))$\phi+\frac{\text{L}^{0}\Psi}{\Psi}$q$_{\text{t-s}}$($-$,P)$\phi
]$(x$_{0}$)

We next consider I$_{2}:$ The differential operator L appearing in I$_{2}$ is
taken with respect to the variable x and since q$_{\text{t-s}}(z,$P) is
independent of x$,$ we can use the compactness of M to differentiate outside
the integral sign and have:\qquad\qquad\ 

$\left(  5.6\right)  \qquad$I$_{2}=$ L$\int_{\text{M}}$q$_{\text{t-s}}$%
($z,$P)$\phi$(z)k$_{\text{s}}$(x,z)$\upsilon_{\text{M}}$(dz)

\qquad\qquad\ \ \ \ \ \ = L$\left[  \text{P}_{\text{s}}\text{(q}_{\text{t-s}%
}\text{(z,P))}\phi\right]  $(x)

By $\left(  2.5\right)  $ and $\left(  2.6\right)  $ of \textbf{Azencott
}$\left[  4\right]  $\textbf{ }we can interchange L and P$_{\text{s}}$ and have:

$\left(  5.7\right)  \qquad$I$_{2}=$ P$_{\text{s}}\left[  \text{L}%
(\text{q}_{\text{t-s}}(-,\text{P})\phi))\right]  $(x)

\qquad$=$ P$_{\text{s}}\left[  \text{L}^{0}\text{(q}_{\text{t-s}}%
\text{(-,P))}\phi\text{ + q}_{\text{t-s}}\text{(-,P)L}\phi\text{ +
$<$%
}\nabla^{0}\text{q}_{\text{t-s}}\text{(-,P),}\nabla\phi\text{%
$>$%
}-V\text{q}_{\text{t-s}}(-,\text{P})\phi\right]  $(x)

By $\left(  5.2\right)  ,$ $\left(  5.5\right)  $ and $\left(  5.7\right)  $
we have:

$\left(  5.8\right)  \qquad\frac{\partial}{\partial\text{s}}\int_{\text{M}}%
$q$_{\text{t-s}}$($z,$P)k$_{\text{s}}$(x,z)$\phi$(z)$\upsilon_{\text{M}}$(dz)

\qquad\qquad\ = I$_{1}$ + I$_{2}$

\qquad\qquad\ = P$_{\text{s}}[-$L$^{0}$(q$_{\text{t-s}}$($-,$P))$\phi
+\frac{\text{L}^{0}\Psi}{\Psi}.$q$_{\text{t-s}}$(z,P)$\phi]($x$_{0})$

$\qquad\qquad+$ P$_{\text{s}}\left[  \text{L}^{0}\text{(q}_{\text{t-s}%
}\text{(}-\text{,P))}\phi\text{ + q}_{\text{t-s}}\text{(-,P)L}\phi\text{ +
$<$%
}\nabla^{0}\text{q}_{\text{t-s}}\text{(-,P),}\nabla\phi\text{%
$>$%
}-\text{Vq}_{\text{t-s}}(-,\text{P})\phi\right]  $(x)

\qquad\qquad= P$_{\text{s}}\left[  \frac{\text{L}^{0}\Psi}{\Psi}%
\text{.q}_{\text{t-s}}\text{(-,P)}\phi\text{ + q}_{\text{t-s}}\text{(-,P)L}%
\phi\text{ +
$<$%
}\nabla^{0}\text{q}_{\text{t-s}}\text{(-,P),}\nabla\phi\text{%
$>$%
}-\text{Vq}_{\text{t-s}}(-,\text{P})\phi\right]  ($x$)$\-

Hence by $\left(  5.1\right)  ,$ $\left(  5.2\right)  $ and $\left(
5.8\right)  $ we have:

$\left(  5.9\right)  \qquad\frac{\partial}{\partial\text{s}}$(Q(t,t-s) $\phi)$

= q$_{\text{t}}$(x,P$)^{-1}$ P$_{\text{s}}\left[  \frac{\text{L}^{0}\Psi}%
{\Psi}\text{.q}_{\text{t-s}}\text{(-,P)}\phi\text{ + q}_{\text{t-s}%
}\text{(-,P)L}\phi\text{ +
$<$%
}\nabla^{0}\text{q}_{\text{t-s}}\text{(-,P),}\nabla\phi\text{%
$>$%
}-\text{Vq}_{\text{t-s}}(-,\text{P})\phi\right]  ($x$)$

= q$_{\text{t}}$(x,P$)^{-1}$ P$_{\text{s}}\left[  \frac{\text{L}^{0}\Psi}%
{\Psi}\text{.q}_{\text{t-s}}\text{(-,P)}\phi\text{ + q}_{\text{t-s}%
}\text{(-,P)L}\phi\text{ + q}_{\text{t-s}}\text{(-,P)%
$<$%
}\nabla^{0}\log\text{q}_{\text{t-s}}\text{(-,P),}\nabla\phi\text{%
$>$%
}-\text{Vq}_{\text{t-s}}(-,\text{P})\phi\right]  $(x)

= q$_{\text{t}}$(x,P$)^{-1}$ P$_{\text{s}}\left[  \text{q}_{\text{t-s}%
}\text{(-,P)}\left(  \text{ L}\phi\text{ + }\frac{\text{L}^{0}\Psi}{\Psi}%
\phi\text{ +
$<$%
}\nabla^{0}\log\text{q}_{\text{t-s}}\text{(-,P),}\nabla\phi\text{%
$>$%
}-\text{V}\phi\right)  \right]  $(x)

= (Q(t,t-s)$\left[  \text{ L}\phi+\text{ }<\nabla^{0}\log\text{q}_{\text{t-s}%
}\text{(-,P),}\nabla\phi\text{%
$>$%
) +}\frac{\text{L}^{0}\Psi}{\Psi}\phi\text{ }-\text{V}\phi\right]  $(x) by the
definition of (Q(t,t-s) in $\left(  4.8\right)  $

= (Q(t,t-s)$\left[  \text{ }\left(  \text{L}+\text{ }\nabla^{0}\log
\text{q}_{\text{t-s}}\text{(-,P) + }\frac{\text{L}^{0}\Psi}{\Psi}\text{
}-\text{V}\right)  \phi\right]  $(x) = (Q(t,t-s)(L$_{s,t}\phi)$(x)

\qquad\qquad(where L$_{s,t}=$ L$+$ $\nabla^{0}\log$q$_{\text{t-s}}$(-,P) +
$\frac{\text{L}^{0}\Psi}{\Psi}$ $-$ V)

= $($Q(t,t-s)$\left[  \frac{\text{L(}\Psi\phi\text{)}}{\Psi}-\frac{1}{\Psi
}<\nabla^{0}\Psi,\nabla\phi>+\text{
$<$%
}\nabla\log\text{q}_{\text{t-s}}\text{(-,P),}\nabla\phi\text{%
$>$%
)}\right]  $(x)

Since $\nabla^{0}\log$q$_{\text{t-s}}(-,$P) = $\nabla^{0}\log\Psi-\frac
{1}{\text{2(t-s)}}\nabla^{0}$d(-,P)$^{2}=\frac{1}{\Psi}\nabla^{0}\Psi-\frac
{1}{\text{2(t-s)}}\nabla^{0}$d(-,P)$^{2},$

we have the last equality.

\qquad\qquad$\qquad\qquad\qquad\qquad\qquad\qquad\qquad\qquad\qquad
\qquad\qquad\qquad\qquad\qquad\blacksquare$

\begin{corollary}
When P = M$_{0},$ we have the following important properties:
\end{corollary}

(i) Define $\phi_{s}=\ ($Q(t,t-s)$\phi).$

\qquad Then $\phi_{s}$ is a solution of the \textbf{Cauchy Problem for the
Heat Equation:}

\qquad For almost all $x\in M,$ we have:

\qquad$\frac{\partial\phi_{s}}{\partial s}(x)=L\phi_{s}(x)$

$\qquad\ \phi_{0}(x)=\phi(x)$

(ii) The above solution $\phi_{s}=\ ($Q(t,t-s)$\phi)=P_{s}\phi$ has the property:

$L(P_{s}\phi)(x)=P_{s}(L\phi)(x)$ for almost all $x\in M.$

\begin{proof}
(i) This is (iii) of \textbf{Corollary 1} (of Theorem 1) above. It has been
repeated here for emphasis,
\end{proof}

(ii) By\textbf{Theorem 2} above, we have for P = M$_{0}:$

$\Psi=1$ and so, $\nabla^{0}\log$q$_{\text{t-s}}$(-,M$_{0}$) = 0 and
$\frac{\text{L}^{0}\Psi}{\Psi}=$ V

Consequently, we have the more simplified expression for L$_{s,t}:$

$L_{s,t}=L$ $+$ $\nabla^{0}\log$q$_{\text{t-s}}$(-,P) + $\frac{\text{L}%
^{0}\Psi}{\Psi}$ $-$ $V=L$ $+$ $\nabla^{0}\log$q$_{\text{t-s}}$(-,M$_{0}$) +
$\frac{\text{L}^{0}\Psi}{\Psi}$ $-$ $V=L=$ $\frac{1}{2}\Delta+X+V$

Consequently, the statement of the \textbf{Theorem} here gives for P =
M$_{0}:$

$\frac{\partial}{\partial\text{s}}($Q(t,t-s)$\phi$) = (Q(t,t-s)( L$\phi)$(x)
\qquad\qquad\qquad\qquad\qquad\qquad\qquad\ \ $\left(  1\right)  $

On the other hand, from (i) above, we have for $\phi_{s}=\ ($Q(t,t-s)$\phi),$

\qquad$\frac{\partial\phi_{s}}{\partial s}(x)=L\phi_{s}(x)$ for all x$\in
M_{0}$

That is, for x$\in M_{0},$ we have:

\qquad$\frac{\partial}{\partial s}(Q(t,t-s)\phi)(x)=L(Q(t,t-s)\phi
)(x)\qquad\qquad\qquad\qquad\qquad\left(  2\right)  \qquad\qquad$

By $\left(  4.39\right)  $ above, we have:

$\qquad\phi_{s}=\ ($Q(t,t-s)$\phi)=$ P$_{s}^{M}\phi\qquad\qquad\qquad
\qquad\qquad\qquad\qquad\qquad\ \ \left(  3\right)  $

Then, $\left(  2\right)  $ and $\left(  3\right)  $ combine to give for x$\in
M_{0}:$

\qquad$\frac{\partial}{\partial s}(Q(t,t-s)\phi)(x)=L(Q(t,t-s)\phi
)(x)=L(P_{s}^{M}\phi)(x)\qquad\ \ \ \ \ \ \left(  4\right)  $

On the other hand, the equation of the \textbf{Theorem} here gives

$\qquad\frac{\partial}{\partial\text{s}}(Q(t,t-s)\phi)(x)=(Q(t,t-s)(L\phi
))(x)=P_{s}^{M}(L\phi)(x)\qquad\ \ \ \ \ \left(  5\right)  $

Equations $\left(  4\right)  $ and $\left(  5\right)  $ above combine to give
the desired result:

\qquad$L(P_{s}^{M}\phi)(x)=P_{s}^{M}(L\phi)(x)$ for almost all $x\in M.$

\qquad\qquad\qquad\qquad\qquad\qquad\qquad\qquad\qquad\qquad\qquad\qquad
\qquad\qquad\qquad$\ \ \ \blacksquare$

We can refine the partial differential equation in the above theorem and
remove the singularity at s = t. This is done by using the operators F(r,s) on
smooth sections of the vector bundle E over M.

F(r,s) is defined on $\Gamma(E)$ as follows: if $\gamma$ is the unique minimal
geodesic from the point x $\in$M$_{0}$ to the submanifold P (meeting the
submanifold orthogonally at a point y$\in$P in time r, define F(r,s) as follows:

$\left(  5.10\right)  \qquad$(F(r,s)$\phi)($x$)=\phi(\gamma($r-s)) =
$\phi\circ\gamma($r-s) for s$\in\lbrack0,r],$

This will enable us to remove the singularity at s = t. The theorem below is
one of the most important theorems of this work. Iis the key to the expansion theorem.

\qquad\qquad\qquad\qquad\qquad\qquad\qquad\qquad\qquad\qquad\qquad\qquad
\qquad\qquad\qquad$\blacksquare$

\begin{theorem}
For t$\geq$r$\geq$s%
$>$%
$0,$
\end{theorem}

For $\phi\in\Gamma(E),$ we have:

$\frac{\partial}{\partial\text{s}}($Q(t,t-s)F(t-s,t-r)$\phi$) =
(Q(t,t-s)($\frac{\text{L(}\Psi\text{F}(\text{t}-\text{s},\text{t}%
-\text{r})\phi\text{)}}{\Psi}$)

\begin{proof}
The proof follows the same lines as in the proof of the above theorem, except
that\ the vector bundle \qquad\ 
\end{proof}

section $\phi$ is replaced by F(t-s,t-r)$\phi$ and hence $\left(  5.2\right)
$ becomes:

$\left(  5.11\right)  \qquad\frac{\partial}{\partial\text{s}}\left[
\text{P}_{\text{s}}\text{(q}_{\text{t-s}}\text{(-,P)F(t-s,t-r)}\phi
\text{))}\right]  $

$=\frac{\partial}{\partial\text{s}}\int_{\text{M}}\left[  \text{k}_{\text{s}%
}\text{(x}_{0}\text{,z)(q}_{\text{t-s}}\text{(-,P)F(t-s,t-r)}\phi
\text{)(z)}\right]  \upsilon_{\text{M}}$(dz) = $\int_{\text{M}}\frac{\partial
}{\partial\text{s}}\left[  \text{k}_{\text{s}}\text{(x}_{0}\text{,z)(q}%
_{\text{t-s}}\text{(-,P)F(t-s,t-r)}\phi\text{)(z)}\right]  \upsilon_{\text{M}%
}$(dz)

= $\int_{\text{M}}$k$_{\text{s}}$(x,z)$\frac{\partial}{\partial\text{s}%
}\left[  \text{(q}_{\text{t-s}}\text{(-,P)F(t-s,t-r)}\phi\text{)(z)}\right]
\upsilon_{\text{M}}$(dz) + $\int_{\text{M}}\frac{\partial}{\partial\text{s}}%
$k$_{\text{s}}$(x$_{\text{0}}$,z)$\phi$(z)$\left[  \text{q}_{\text{t-s}%
}\text{(-,P)F(t-s,t-r)}\phi\right]  $(z)$\upsilon_{\text{M}}$(dz)

= $\int_{\text{M}}$k$_{\text{s}}$(x,z)$\frac{\partial}{\partial\text{s}%
}\left[  \text{(q}_{\text{t-s}}\text{(-,P)F(t-s,t-r)}\phi\text{)(z)}\right]
\upsilon_{\text{M}}$(dz) + $\int_{\text{M}}$Lk$_{\text{s}}$(x$_{0}$,z)$\left[
\text{q}_{\text{t-s}}\text{(-,P)F(t-s,t-r)}\phi\right]  $(z)$\upsilon
_{\text{M}}$(dz)

= $\int_{\text{M}}$k$_{\text{s}}$(x,z)$\frac{\partial}{\partial\text{s}%
}\left[  \text{(q}_{\text{t-s}}\text{(-,P)F(t-s,t-r)}\phi\text{)(z)}\right]
\upsilon_{\text{M}}$(dz) + L$\int_{\text{M}}$k$_{\text{s}}$(x$_{0}$,z)$\left[
\text{q}_{\text{t-s}}\text{(-,P)F(t-s,t-r)}\phi\right]  $(z)$\upsilon
_{\text{M}}$(dz)

Setting u = t-s,

$\left(  5.12\right)  \qquad\frac{\partial}{\partial\text{s}}($q$_{\text{t-s}%
}(-,$P$)$F(t-s,t-r)$\phi$)\ =\ $-$ $\frac{\partial}{\partial\text{u}}%
($q$_{\text{u}}(-,$P$).$F(u,t-r)$\phi-$\ \ q$_{\text{u}}$(-,P)$\frac{\partial
}{\partial\text{u}}$(F(u,t-r)$\phi$)

\ = $\left(  -\text{Lq}_{\text{u}}\text{(-,P)\ \ +}\frac{\text{L}\Psi}{\Psi
}\text{.q}_{\text{u}}\text{(-,P)}\right)  $.F(u,t-r)$\phi-$ q$_{\text{u}}%
$(-,P).$\frac{\partial}{\partial\text{u}}$(F(u,t-r)$\phi$) by \textbf{Lemma 2}

= $\left(  -\text{Lq}_{\text{t-s}}\text{(-,P)\ \ +}\frac{\text{L}\Psi}{\Psi
}\text{.q}_{\text{t-s}}\text{(-,P)}\right)  $.F(t-s,t-r)$\phi-$ q$_{\text{t-s}%
}$(-,P).$\frac{\partial}{\partial\text{u}}$(F(u,t-r)$\phi$)\ $\qquad$

Consequently by the last equality in $\left(  5.11\right)  $ and the last
equality in $\left(  5.12\right)  ,$ we have:

$\left(  5.13\right)  \qquad\frac{\partial}{\partial\text{s}}\left[
\text{P}_{\text{s}}\text{(q}_{\text{t-s}}\text{(-,P)F(t-s,t-r)}\phi
\text{))}\right]  $

$=$ $\int_{\text{M}}$k$_{\text{s}}$(x,z)$\left[  \left(  -\text{Lq}%
_{\text{t-s}}\text{(-,P)\ \ +}\frac{\text{L}\Psi}{\Psi}\text{.q}_{\text{t-s}%
}\text{(-,P)}\right)  .\text{F(t-s,t-r)}\phi\text{ }-\text{q}_{\text{t-s}%
}\text{(-,P).}\frac{\partial}{\partial\text{u}}\text{(F(u,t-r)}\phi
\text{)}\right]  $(z)$\upsilon_{\text{M}}$(dz)

+ L$\int_{\text{M}}$k$_{\text{s}}$(x,z)$\left[  \text{q}_{\text{t-s}%
}\text{(-,P)F(t-s,t-r)}\phi\right]  $(z)$\upsilon_{\text{M}}$(dz)

$=$ P$_{\text{s}}\left[  \left(  -\text{Lq}_{\text{t-s}}\text{(-,P)\ \ +}%
\frac{\text{L}\Psi}{\Psi}\text{.q}_{\text{t-s}}\text{(-,P)}\right)
.\text{F(t-s,t-r)}\phi\text{ }-\text{q}_{\text{t-s}}\text{(-,P).}%
\frac{\partial}{\partial\text{u}}\text{(F(u,t-r)}\phi\text{)}\right]  $

+ LP$_{\text{s}}\left[  \text{q}_{\text{t-s}}\text{(-,P)F(t-s,t-r)}%
\phi\right]  $

$=$ P$_{\text{s}}\left[  \left(  -\text{Lq}_{\text{t-s}}\text{(-,P)\ \ +}%
\frac{\text{L}\Psi}{\Psi}\text{.q}_{\text{t-s}}\text{(-,P)}\right)
.\text{F(t-s,t-r)}\phi\text{ }-\text{q}_{\text{t-s}}\text{(-,P).}%
\frac{\partial}{\partial\text{u}}\text{(F(u,t-r)}\phi\text{)}\right]  $

+ P$_{\text{s}}\left[  \text{L}\left(  \text{q}_{\text{t-s}}%
\text{(-,P)F(t-s,t-r)}\phi\right)  \right]  $

( where we have inter-changed L and P$_{\text{s}}$ in the last term above).

$=$ P$_{\text{s}}\left[  \left(  -\text{Lq}_{\text{t-s}}\text{(-,P)\ \ +}%
\frac{\text{L}\Psi}{\Psi}\text{.q}_{\text{t-s}}\text{(-,P)}\right)
.\text{F(t-s,t-r)}\phi\text{ }-\text{q}_{\text{t-s}}\text{(-,P).}%
\frac{\partial}{\partial\text{u}}\text{(F(u,t-r)}\phi\text{)}\right]  $

+ P$_{\text{s}}\left[  \text{L(q}_{\text{t-s}}\text{(-,P)).F(t-s,t-r)}%
\phi\text{ + (q}_{\text{t-s}}\text{(-,P)L(F(t-s,t-r)}\phi\text{) + }%
<\nabla\text{q}_{\text{t-s}}\text{(-,P),}\nabla\text{F(t-s,t-r)}\phi\text{
$>$%
}-V\text{q}_{\text{t-s}}\text{(-,P)F(t-s,t-r)}\phi\right]  $\ \ 

$=$ P$_{\text{s}}\left[  \left(  \frac{\text{L}\Psi}{\Psi}.\text{F(t-s,t-r)}%
\phi\text{ + L}\left(  \text{F(t-s,t-r)}\phi\right)  -\frac{\partial}%
{\partial\text{u}}\text{(F(u,t-r)}\phi\text{)}\right)  \text{.q}_{\text{t-s}%
}\text{(-,P) +
$<$%
}\nabla\text{q}_{\text{t-s}}\text{(-,P),}\nabla\text{F(t-s,t-r)}\phi\text{%
$>$%
}-V\text{q}_{\text{t-s}}\text{(-,P)F(t-s,t-r)}\phi\right]  $

$=$ P$_{\text{s}}\left[  \left(  \frac{\text{L(}\Psi\text{F(t-s,t-r)}%
\phi\text{)}}{\Psi}-\text{%
$<$%
}\frac{\nabla\Psi}{\Psi},\nabla\text{(F(u,t-r)}\phi\text{)%
$>$
}-\frac{\partial}{\partial\text{u}}\text{(F(u,t-r)}\phi\text{)}\right)
\text{.q}_{\text{t-s}}\text{(-,P) +
$<$%
}\nabla\text{q}_{\text{t-s}}\text{(-,P),}\nabla\text{F(t-s,t-r)}\phi>\right]
$

It is an elementary fact that $\nabla$q$_{\text{t-s}}$(-,P) = q$_{\text{t-s}}%
$(-,P)$\nabla\log$q$_{\text{t-s}}(-,$P).

Since $\nabla\log$q$_{\text{t-s}}(-,$P) = $\nabla\log\Psi-\frac{1}%
{\text{2(t-s)}}\nabla$d(-,P)$^{2}=\frac{1}{\Psi}\nabla\Psi-\frac
{1}{\text{2(t-s)}}\nabla$d(-,P)$^{2},$

$\left(  5.13\right)  $ becomes:

$\left(  5.14\right)  \qquad\frac{\partial}{\partial\text{s}}\left[
\text{P}_{\text{s}}\text{(q}_{\text{t-s}}\text{(-,P)F(t-s,t-r)}\phi
\text{))}\right]  $

$=$ P$_{\text{s}}\left[  \left(  \frac{\text{L(}\Psi\text{F(t-s,t-r)}%
\phi\text{)}}{\Psi}\text{ }-\frac{\partial}{\partial\text{u}}\text{(F(u,t-r)}%
\phi\text{)}\right)  \text{.q}_{\text{t-s}}\text{(-,P) }-\frac{1}%
{\text{2(t-s)}}\text{q}_{\text{t-s}}\text{(-,P)}<\nabla\text{d(-,P)}%
^{2}\text{,}\nabla\text{F(t-s,t-r)}\phi>\right]  $

Therefore by the definition of Q(t,t-s) in $\left(  4.6\right)  $ and the last
equality in $\left(  5.14\right)  $ above,

$\left(  5.15\right)  \qquad\ \ \ \frac{\partial}{\partial\text{s}}%
$(Q$_{\text{P}}^{\text{M}_{0}}($t,t-s)F(t-s,t-r)$\phi$)(x$_{0}$) =
q$_{\text{t}}$(x$_{0}$,P$)^{-1}\frac{\partial}{\partial\text{s}}\left[
\text{P}_{\text{s}}\text{(q}_{\text{t-s}}\text{(-,P)F(t-s,t-r)}\phi
\text{)}\right]  $

= q$_{\text{t}}$(x,P$)^{-1}$P$_{\text{s}}\left[  \left(  \frac{\text{L(}%
\Psi\text{F(t-s,t-r)}\phi\text{)}}{\Psi}\text{ -}\frac{\partial}%
{\partial\text{u}}\text{(F(u,t-r)}\phi\text{)}\right)  \text{.q}_{\text{t-s}%
}\text{(-,P) }-\frac{\text{q}_{\text{t-s}}\text{(-,P)}}{\text{2(t-s)}}%
<\nabla\text{d(-,P)}^{2}\text{,}\nabla\text{F(t-s,t-r)}\phi\right]  $

= Q$_{\text{P}}^{\text{M}_{0}}($t,t-s)$\left[  \left(  \frac{\text{L(}%
\Psi\text{F(t-s,t-r)}\phi\text{)}}{\Psi}\text{ }-\frac{\partial}%
{\partial\text{u}}\text{(F(u,t-r)}\phi\text{)}\right)  .\text{q}_{\text{t-s}%
}\text{(-,P)}-\frac{\text{q}_{\text{t-s}}\text{(-,P)}}{\text{2(t-s)}}\text{%
$<$%
}\nabla\text{d(-,P)}^{2}\text{,}\nabla\text{F(t-s,t-r)}\phi\text{%
$>$%
}\right]  $

We next show that $\frac{\partial}{\partial\text{u}}$(F(u,t-r)$\phi$) =
$-\frac{1}{2\text{(t-s)}}$%
$<$%
\ $\nabla$d(-,P)$^{2}$ $,$ $\nabla$F(t-s,t-r)$\phi$%
$>$%

The relation in the last equality above is of crucial importance in this work
and so we prove it in detail.

This will enable us to remove the singularity at s = t.

By the definition of the operator (F(r,s) in $\left(  5.10\right)  $ above,

$\left(  5.16\right)  \qquad\frac{\partial}{\partial\text{u}}$F(u,t-r)$\phi
$)(x)$\mid_{\text{u = t-s}}=$ $\frac{\partial}{\partial\text{u}}(\phi
\circ\gamma)$(u-t+r))$\mid_{\text{u = t-s}}$

where $\gamma$ is the unique minimal geodesic from x$\in$M$_{0}$ to y$\in$P in
time u. We have:

$\left(  5.17\right)  \qquad$ $\frac{\partial}{\partial\text{u}}(\phi
\circ\gamma)$(u-t+r)$\mid_{\text{u= t-s}}=$ $<\nabla\phi$($\gamma$(u-t+r) ,
$\dot{\gamma}($u-t+r))$>\mid_{\text{u=t-s}}$

We write:

$\left(  5.18\right)  \qquad\gamma$($\lambda$) = y $+$ $(1-\frac{\lambda
}{\text{u}})($x$-$ y$)$ (in vector form) to mean:

\qquad\qquad\qquad\ \ \ \ \ = $\left(  \text{x}_{1}\text{,...,x}_{q}\text{,
(1-}\frac{\lambda}{\text{u}}\text{)x}_{q+1}\text{,.., (1-}\frac{\lambda
}{\text{u}}\text{)x}_{n}\right)  $ in Fermi coordinates.

$\qquad\qquad\qquad\ \ \ =\left(  \text{y}_{1}\text{,...,y}_{q}\text{,
(1-}\frac{\lambda}{\text{u}}\text{)x}_{q+1}\text{,.., (1-}\frac{\lambda
}{\text{u}}\text{)x}_{n}\right)  $ by $\left[  2.2\right]  $ in Gray $\left[
4\right]  .$

and hence,

$\left(  5.19\right)  \qquad\gamma$(u-t+r)\ = (y$_{1},..,$y$_{q}%
,\frac{\text{t}-\text{r}}{\text{u}}$x$_{q+1},..,\frac{\text{t}-\text{r}%
}{\text{u}}$x$_{n})$

Therefore, taking derivatives with respect to u at u = t-s:

$\left(  5.20\right)  \qquad\dot{\gamma}($u-t+r)$\mid_{\text{u=t-s}}=$
(0$,..,$0$,-\frac{\text{t}-\text{r}}{\text{(t-s)}^{2}}$x$_{q+1},..,-\frac
{\text{t}-\text{r}}{\text{(t-s)}^{2}}$x$_{n})=$ $-\frac{\text{t-r}%
}{\text{(t-s)}^{2}}\underset{\text{j=q+1}}{\overset{\text{n}}{\sum}}$%
x$_{j}\frac{\partial}{\partial\text{x}_{j}}$

and so the derivative in $\left(  5.16\right)  -\left(  5.17\right)  $ becomes:

$\left(  5.21\right)  \qquad\frac{\partial}{\partial\text{u}}$F(u,t-r)$\phi
$)$\mid_{\text{u = t-s}}=\frac{\partial}{\partial\text{u}}(\phi\circ\gamma
)$(u-t+r)$\mid_{\text{u= t-s}}$

$\qquad\qquad\qquad=$ $-$ $\frac{\text{t-r}}{\text{(t-s)}^{2}}$%
$<$%
$\underset{\text{j=q+1}}{\overset{\text{n}}{\sum}}$x$_{j}\frac{\partial
}{\partial\text{x}_{j}}$ $,\nabla\phi$($\gamma$(r-s))%
$>$%

We next show that $\nabla\phi$($\gamma$(r-s))$\frac{\text{t-r}}{\text{t-s}}$ =
$\nabla$F(t-s,t-r)$\phi:$

Consider $\nabla\phi\circ\gamma$(u-t+r)$\mid_{\text{u=t-s}}$where $\gamma$ is
now the unique minimal geodesic from x to y in time t-s. In local Fermi
coordinates $\gamma$ can be written as:

$\gamma(\lambda)=\left(  y_{1},..,y_{q},(1-\frac{\lambda}{\text{t-s}}%
)x_{q+1},..,(1-\frac{\lambda}{\text{t-s}})x_{n}\right)  $ and so,

$\left(  5.22\right)  \qquad\frac{\partial}{\partial\text{x}_{\text{a}}}%
\gamma(\lambda)=1$ for a = 1,...,q and $\frac{\partial}{\partial\text{x}_{j}%
}\gamma(\lambda)=(1-\frac{\lambda}{\text{t-s}})$ for $j=q+1,...,n.$

Now by $\left(  5.22\right)  $,

$\left(  5.23\right)  \qquad\nabla\phi\circ\gamma(\lambda)=\nabla\phi
(\gamma(\lambda))\frac{\partial}{\partial\text{x}_{j}}\gamma(\lambda)$

$\qquad\qquad\qquad\qquad=\left\{
\begin{array}
[c]{c}%
1\text{ for }j=\text{ }1,...,q\\
\nabla\phi(\gamma(\lambda))(1-\frac{\lambda}{\text{t-s}})\text{ for }j\text{
}=q+1,...,n
\end{array}
\right.  $

Given the expression in $\left(  5.21\right)  ,$ we will use the expression in
$\left(  5.23\right)  $ only for

$j=q+1,...,n$. Consequently,

$\left(  5.24\right)  \qquad\nabla\phi\circ\gamma$(r-s) $=\nabla\phi(\gamma
($r-s$))(1-\frac{\text{r-s}}{\text{t-s}})=\nabla\phi(\gamma($r-s$))(\frac
{\text{t-r}}{\text{t-s}})$

By the definition of the operator F(t-s,t-r) in $\left(  5.10\right)  $

This will enable us to remove the singularity at s = t. The theore$,$

F(t-s,t-r)$\phi=\phi\circ\gamma$(r-s) and so by $\left(  5.24\right)  ,$

$\left(  5.25\right)  \qquad\nabla($F(t-s,t-r)$\phi)=\nabla(\phi\circ\gamma
)$(r-s) $=\nabla\phi(\gamma($r-s$))\frac{\text{t-r}}{\text{t-s}},\qquad$

Hence by $\left(  5.21\right)  $ and $\left(  5.24\right)  $ we have:

$\left(  5.26\right)  \qquad\frac{\partial}{\partial\text{u}}\phi$($\gamma
$(u-t+r))$\mid_{\text{u = t-s}}=-\frac{1}{\text{(t-s)}}$%
$<$%
$\underset{\text{j=q+1}}{\overset{\text{n}}{\sum}}$x$_{j}\frac{\partial
}{\partial\text{x}_{j}}$ $,\nabla$F(t-s,t-r)$\phi$%
$>$%

Following \textbf{Gray }$\left[  23\right]  $ or \textbf{Gray} $\left[
25\right]  $, \textbf{Definition }$\left(  2.19\right)  $, we set:

$\left(  5.27\right)  \qquad\rho^{\text{2}}$ $\ =$ $\underset{\text{j=q+1}%
}{\overset{\text{n}}{\sum}}$x$_{j}^{2}$

Then by \textbf{Lemma} $\left(  2.7\right)  $ of \textbf{Gray} $\left[
25\right]  $, we have:

$\left(  5.28\right)  \qquad$ $\rho$(x) $=$ d(x,P).

By $\left(  5.27\right)  $ above and \textbf{Lemma} $(3.1)$ of \textbf{Gray}
$\left[  23\right]  $ or \textbf{Lemma} $\left(  2.11\right)  $ of
\textbf{Gray} $\left[  25\right]  ,$

we have:

$\left(  5.28\right)  \qquad\nabla\frac{\text{d(x,P)}^{2}}{2}=\nabla\frac
{\rho^{2}}{2}(x)=\underset{j=q+1}{\overset{n}{\sum}}$x$_{j}(x)\frac{\partial
}{\partial\text{x}_{j}}$ \ 

and so we have finally:

$\left(  5.29\right)  \qquad\frac{\partial}{\partial\text{u}}\phi$($\gamma
$(u-t+r))$\mid_{\text{u =t-s}}$ $=-\frac{1}{\text{(t-s)}}$%
$<$%
\ $\nabla\frac{\text{d(-,P)}^{2}}{2}$ $,$ $\nabla$F(t-s,t-r)$\phi$%
$>$%

The result then follows from the last equality in $\left(  5.15\right)  $
above and $\left(  5.29\right)  $ here:

$\left(  5.30\right)  \qquad\frac{\partial}{\partial\text{s}}$(Q$_{\text{P}%
}^{\text{M}_{0}}($t,t-s)F(t-s,t-r)$\phi$)(x) = Q$_{\text{P}}^{\text{M}_{0}}%
($t,t-s)$\left[  \frac{\text{L(}\Psi\text{F(t-s,t-r)}\phi)}{\Psi}\right]
$(x).$\ \ \ \ \ $

\vspace{1pt}The theorem is proved and the singularity at s = t is thus removed
from the

last expression in $\left(  5.15\right)  $.

\qquad\qquad\qquad\qquad\qquad\qquad\qquad\qquad\qquad\qquad\qquad\qquad
\qquad\qquad\qquad$\qquad\ \blacksquare$

The partial differential equation of the above theorem is of crucial importance

in the expansion theorem below.

We recall that $\Psi$(x) $=\theta_{P}$(x)$^{-\frac{1}{2}}\ \ \Phi_{P}($x$)$
was defined in $\left(  1.7\right)  $ of \textbf{Chapter 1}$.$

\qquad\qquad\qquad\qquad\qquad\qquad\qquad\qquad\qquad\qquad\qquad\qquad
\qquad\qquad\qquad\qquad$\blacksquare$

\part{GENERALIZED EXPANSIONS}

\chapter{Expansion of Generalized Feynman-Kac Formula}

We will first prove a preliminary expansion theorem which is an expansion of
the Right Hand Side of the Generalized Feynman-Kac Formula. We shall see that
this expansion generalizes the Heat Kernel Expansion.

By the definition of L$_{\Psi}$, we set:

$\left(  6.1\right)  \qquad$L$_{\Psi}\phi$ = $\frac{\text{L(}\Psi\phi\text{)}%
}{\Psi}$

\begin{theorem}
(Expansion of the Generalized Feynman-Kac Formula)
\end{theorem}

Let $\gamma$ be the unique minimal geodesic from x$\in$M$_{0}$ to y$\in$P in
time t. Then for

t$\geq s\geq0,$ 1$\leq n\leq$N, and smooth section $\phi$ we have:

$\left(  6.2\right)  $ \ (Q(t,t-s)$\phi$)(x) = \textbf{E}$_{\text{x}}\left[
\tau_{\text{0,s}}^{\text{t}}\text{e}_{\text{s}}^{\text{t}}\phi\text{(x}%
^{\text{t}}\text{(s))exp}\left\{  \int_{0}^{\text{s}}\frac{\text{L}^{0}\Psi
}{\Psi}\text{(x}^{\text{t}}\text{(r))dr}\right\}  \right]  $

\qquad\qquad\qquad\qquad\qquad\ \ = $\phi$($\gamma$(s))
\ +\ $\underset{n=1}{\;\overset{\text{N}}{\sum}}$a$_{n}$(s,x,P,$\phi
$)\ +\ F$_{\text{N+1}}$(s,x,P,$\phi$)\qquad\qquad\qquad\qquad\qquad
\ \ \ \qquad\qquad\qquad\qquad\qquad\qquad

where, for 0$\leq$s$_{n}\leq$s$_{n-1}\leq$...$\leq$s$_{1}\leq$s$\leq$t$,$

$\left(  6.3\right)  $ \ a$_{n}$(s,x,P,$\phi$) $=\int_{0}^{\text{s}}\int%
_{0}^{\text{s}_{1}}...\int_{0}^{\text{s}_{n-1}}($F(t,t-s$_{n}$)L$_{\Psi}%
$F(t-s$_{n}$,t-s$_{n-1}$)L$_{\Psi}$F(t-s$_{n-1}$,t-s$_{n-2}$)

\qquad\qquad\qquad\qquad\qquad... L$_{\Psi}$F(t-s$_{2}$,t-s$_{1}$)L$_{\Psi}%
$F(t-s$_{1}$,t-s)$\phi$)$($x$)$ds$_{1}...$ds$_{n}$

\ \ \ \ \ \ \ \ \ and, for 0$\leq$s$_{N+1}\leq$s$_{N}\leq$...$\leq$s$_{1}\leq$s

$\left(  6.4\right)  \qquad$F$_{\text{N+1}}(s,$x$_{0},$P)

\ \qquad\qquad= $\int_{0}^{\text{s}}\int_{0}^{\text{s}_{1}}...\int%
_{0}^{\text{s}_{\text{N}}}($Q$($t,t-s$_{\text{N+1}}$)L$_{\Psi}$%
F(t-s$_{\text{N+1}}$,t-s$_{\text{N}}$)L$_{\Psi}$F(t-s$_{\text{N}}%
$,t-s$_{\text{N-1}})...$

\qquad\qquad... L$_{\Psi}$F(t-s$_{2}$,t-s$_{1}$)L$_{\Psi}$F(t-s$_{1}%
$,t-s)$\phi$)$($x$)$ds$_{1}...$ds$_{N+1}$

\begin{proof}
$\qquad$\qquad
\end{proof}

The first equation in $\left(  6.2\right)  $ above is the \textbf{Generalized
Feynman-Kac Formula}

of \textbf{Theorem 1} above and the second equation gives the generalized expansion.

First we note that using the notation in $\left(  6.1\right)  $, the partial
differential equation

in \textbf{Theorem 2} can be re-written as:

$\left(  6.5\right)  \qquad\frac{\partial}{\partial\text{s}}($%
Q(t,t-s)F(t-s,t-r)$\phi)$ = $($Q(t,t-s)(L$_{\Psi}$F(t-s,t-r)$\phi)$

We integrate each side of the equation in $\left(  6.5\right)  $ above and
have for 0$\leq$s$_{1}\leq$s$\leq$t:

$\left(  6.6\right)  \qquad\int_{0}^{\text{s}}\frac{\partial}{\partial
\text{s}_{1}}$(Q$_{\text{P}}($t,t-s$_{1}$)F(t-s$_{1}$,t-s)$\phi$)(x$_{0}%
$)ds$_{1}$ = $\int_{0}^{\text{s}}$Q$_{\text{P}}($t,t-s$_{1}$)$\left[
\text{L}_{\Psi}\text{F(t-s}_{1}\text{,t-s)}\phi)\right]  $(x)ds$_{1}%
\qquad\ \ \ \ \ \ \ \ $

The Left Hand Side of $\left(  6.6\right)  $ is easily integrated to give:

$\left(  6.7\right)  $\qquad(Q$_{\text{P}}$(t,t-s)F(t-s,t-s)$\phi$)(x) \ -
\ (Q$_{\text{P}}$(t,t)F(t,t-s)$\phi$)(x).\qquad\qquad\qquad\qquad\qquad
\qquad\qquad

Since F(t-s,t-s) and Q$_{\text{P}}$(t,t) are identity operators, $\left(
6.7\right)  $ becomes:

$\left(  6.8\right)  $\qquad(Q$_{\text{P}}$(t,t-s)$\phi$)(x) - F(t,t-s)$\phi
$)(x$_{0}$) = (Q$_{\text{P}}($t,t-s)$\phi$)(x) - $\phi$($\gamma($s))

and so $\left(  6.6\right)  $ becomes:

$\left(  6.9\right)  $\qquad(Q$_{\text{P}}$(t,t-s)$\phi$)(x) = $\phi$%
($\gamma($s)) + $\int_{0}^{\text{s}}$(Q$_{\text{P}}$(t,t-s$_{1}$)L$_{\Psi}%
$F(t-s$_{1}$,t-s)$\phi)($x$)$ds$_{1}$

where $\gamma$ is the unique minimal geodesic from x$\in M_{0}$ to the centre
of Fermi

coordinates y in time 1.$\qquad\qquad\ \ \ \ \ \ $

\qquad\ \ \ \ \ Set $\phi_{1}$ = L$_{\Psi}$F(t-s$_{1}$,t-s)$\phi$. Then
$\phi_{1}$ is a smooth section of E. We then

re-write $\left(  6.9\right)  $ where we replace $\phi$ by $\phi_{1}:$

We have by $\left(  6.9\right)  $ for $\ $0$\leq$s$_{2}\leq$s$_{1}:$

$\left(  6.10\right)  $\qquad(Q$_{\text{P}}$(t,t-s$_{1}$)$\phi_{1})$(x) =
$\phi_{1}$($\gamma$(s$_{1}$)) +$\int_{0}^{\text{s}_{1}}$(Q$_{\text{P}}%
$(t,t-s$_{2}$)L$_{\Psi}$F(t-s$_{2}$,t-s$_{1}$)$\phi_{1})($x)ds$_{2}%
\qquad\ \ \ \ \ \ \ \ \ \ $

\qquad= (L$_{\Psi}$F(t-s$_{1}$,t-s)$\phi)$($\gamma$(s$_{1}$)) $+\int%
_{0}^{\text{s}_{1}}$(Q$_{\text{P}}$(t,t-s$_{2}$)L$_{\Psi}$F(t-s$_{2}$%
,t-s$_{1}$)L$_{\Psi}$F(t-s$_{1}$,t-s)$\phi$)(x)ds$_{2}$\ \ 

\qquad= F(t,t-s$_{1})$(L$_{\Psi}$F(t-s$_{1}$,t-s)$\phi)$($x$) $+\int%
_{0}^{\text{s}_{1}}$(Q$_{\text{P}}$(t,t-s$_{2}$)L$_{\Psi}$F(t-s$_{2}$%
,t-s$_{1}$)L$_{\Psi}$F(t-s$_{1}$,t-s)$\phi$)(x)ds$_{2}$\ 

\qquad Now, since $\phi_{1}$ = L$_{\Psi}$F(t-s$_{1}$,t-s)$\phi,$ the equality
in $\left(  6.9\right)  $ can be re-written as:

$\left(  6.11\right)  \qquad$(Q$_{\text{P}}$(t,t-s)$\phi$)(x) = $\phi$%
($\gamma($s)) + $\int_{0}^{\text{s}}$(Q$_{\text{P}}$(t,t-s$_{1}$)$\phi_{1}%
($x$)$ds$_{1}$

\qquad We now replace (Q$_{\text{P}}$(t,t-s$_{1}$)$\phi_{1}($x$)$ in $\left(
6.11\right)  $ by the \textbf{last} expression on the RHS of $\left(
6.10\right)  $ and have:

$\left(  6.12\right)  $\qquad(Q$_{\text{P}}$(t,t-s)$\phi$)(x) = $\phi$%
($\gamma($s)) \ +\ \ $\int_{0}^{\text{s}}$F(t,t-s$_{1})$L$_{\Psi}$F(t-s$_{1}%
$,t-s)$\phi$)(x)ds$_{1}$\qquad\qquad\qquad\qquad\qquad

\qquad$\ \ \ \qquad\qquad\qquad+\int_{0}^{\text{s}}\int_{0}^{\text{s}_{1}}%
$(Q$_{\text{P}}$(t,t-s$_{2}$)L$_{\Psi}$F(t-s$_{2}$,t-s$_{1}$)L$_{\Psi}%
$F(t-s$_{1}$,t-s)$\mathbf{\phi})($x)ds$_{1}$ds$_{2}\ \ \ \ $

\qquad\qquad\qquad\qquad\qquad= $\phi$($\gamma($s)) \ + a$_{1}$(s,t,x,P,$\phi
$) $+$ F$_{\text{2}}$(s,t$,$x$_{0},$P)

\qquad The formulae of the theorem are thus proved for n = 1. We then use
induction on n and the method above to obtain the general formula noting that
in general,\qquad\qquad\qquad\qquad\ 

$\left(  6.13\right)  $\qquad F$_{\text{n}}$(s,t,x,P,$\phi$) = a$_{\text{n}}%
($s,t,x,P,$\phi$) + F$_{\text{n +1}}($s,t,x,P,$\phi$).$\qquad\qquad
\qquad\qquad\qquad\qquad\qquad\qquad\qquad\qquad\qquad\qquad\qquad\qquad
\qquad\qquad$

$\qquad\qquad\qquad\qquad\qquad\qquad\qquad\qquad\qquad\qquad\qquad
\qquad\qquad\qquad\qquad\qquad\qquad\blacksquare$

\chapter{Generalized Heat Kernel Expansions}

We now come to one of the \textbf{key theorems} of this work. Given its
importance, the proof will be given in full detail.$\qquad\qquad\qquad\ $

\begin{theorem}
(Exact Expansion of the Generalized Heat Kernel)
\end{theorem}

Let $\gamma$ be the unique minimal geodesic from a point x$\in$M$_{0}$ to P in
time t, meeting P at the centre of Fermi coordinates y$_{0}$: $\gamma$(t) =
y$_{0}\in$P. Then we have the \textbf{exact expansion of the Generalized Heat
Kernel}:

(i)$\ \qquad$k$_{\text{t}}$(x,P,$\phi$) $=\int_{\text{P}}$k$_{\text{t}}%
$(x,y)$\phi($y)$\upsilon_{\text{P}}$(dy)\ \ \ \ \ \ \ \qquad

\qquad\ \ \ \ \ \ \ \ \ \ \qquad= q$_{\text{t}}$(x,P)$\left[  \text{b}%
_{0}\text{(x,P,}\phi\text{)\ + \ }\underset{n=1}{\overset{\text{N}}{\sum}%
}\text{b}_{n}\text{(x,P,}\phi\text{)t}^{n}\text{ + R}_{\text{N+1}%
}(\text{t,x,P,}\phi\text{)t}^{\text{N+1}}\right]  $

where the expansion coefficients are given as follows:

(ii)\qquad b$_{0}($x,P,$\phi$) = $\phi$($\gamma$(t)) $=$ $\phi$($y_{0}$) =
$\tau_{\text{x,y}_{0}}\phi($x$)$

where $\gamma:[0,$t$]\rightarrow$ M$_{0}$ is the unique minimal geodesic from
x$\in$M$_{0}$ to y$_{0}\in$P in time t.

For 1$\geq$r$_{0}\geq$r$_{1}\geq$r$_{2}\geq$...$\geq$r$_{\text{N}}\geq
$r$_{\text{N+1}}$ and hence for, (1-r$_{\text{N+1}}$)$\geq$(1-r$_{\text{N}}%
$)$\geq$(1-r$_{\text{N-1}}$)

$\geq$(1-r$_{\text{N-2}}$)...$\geq$(1-r$_{\text{2}}$)$\geq$(1-r$_{\text{1}}%
$)$\geq$0, we have:

(iii)$\ \qquad$b$_{1}($x,P,$\phi$) = $\int_{0}^{1}$F(1,1-r$_{1})$[L$_{\Psi
}\phi\circ\pi_{\text{P}}$](x)dr$_{1}$

(iv)\qquad b$_{2}($x,P$,\phi)$ = $\int_{0}^{1}\int_{0}^{r_{1}}$F(1,1-r$_{2}%
$)[L$_{\Psi}$F(1-r$_{2}$,1-r$_{1}$)L$_{\Psi}\phi\circ\pi_{\text{P}}$%
](x$)$dr$_{1}$dr$_{2}$

\qquad\ \ \ \ \ and for,\ \ 3$\leq n\leq$N, we have the general formula:

(v)\qquad\ b$_{n}($x,P$,\phi)$ = $\int_{0}^{1}\int_{0}^{r_{1}}...\int%
_{0}^{r_{n-1}}$F(1,1-r$_{n}$)[L$_{\Psi}$F(1-r$_{n}$,1-r$_{n-1}$)

\qquad\qquad L$_{\Psi}$F(1-r$_{n-1}$,1-r$_{n-2}$)$...\ $L$_{\Psi}$F(1-r$_{2}%
$,1-r$_{1}$)L$_{\Psi}\phi\circ\pi_{\text{P}}$](x$)$dr$_{1}...$dr$_{n}$

where $\pi_{\text{P}}$ is the projection: $\pi_{\text{P}}$: M$_{0}\rightarrow$
P viewed in Fermi coordinates.

(vi) R$_{\text{N+1}}$(t,x,P,$\phi$) $=\int_{0}^{1}\int_{0}^{\text{r}_{1}%
}...\int_{0}^{\text{r}_{\text{N}}}($Q(t,t-tr$_{\text{N+1}})[$L$_{\Psi}%
$F(1-r$_{\text{N+1}}$,1-r$_{\text{N}}$)L$_{\Psi}$F(1-r$_{\text{N}}%
$,1-r$_{\text{N-1}}$)

\qquad\ \ \qquad\ ... L$_{\Psi}$F(1-r$_{2}$,1-r$_{1}$)L$_{\Psi}\phi\circ
\pi_{\text{P}}]$(x)dr$_{1}...$dr$_{\text{N+1}}$

where $\zeta$ is the first exit time of the bridge process from the tubular
neighborhood M$_{0}$.

\begin{proof}
(i) We take limits on both sides of $\left(  6.2\right)  $ above and by
\textbf{Lemma 1}:
\end{proof}

$\left(  7.1\right)  $\qquad$\underset{\text{s}\uparrow\text{t}}{lim}%
($Q(t,t-s)$\phi)$(x) = q$_{\text{t}}$(x,P)$^{-1}$ $\int_{\text{P}}%
$k$_{\text{t}}$(x,z)$\phi$(x)$\upsilon_{\text{P}}$(dz)

\qquad and by $\left(  6.2\right)  $ of the last theorem above:

$\left(  7.2\right)  \qquad$ $\underset{\text{s}\uparrow\text{t}}{\text{lim}}%
$(Q(t,t-s)$\phi$)(x) $=\mathbf{E}_{\text{x}}\left[  \tau_{\text{0,t}}%
\text{e}_{\text{t}}\phi\text{(x(t))}\exp\left\{  \int_{0}^{\text{t}}%
\frac{\text{L}^{0}\Psi}{\Psi}\text{(x(r))dr}\right\}  \right]  $

\qquad\qquad\qquad\ = $\phi$($\gamma$(t))
\ +\ $\underset{n=1}{\;\overset{\text{N}}{\sum}}$ $\underset{\text{s}%
\uparrow\text{t}}{\text{lim}}$ a$_{n}($s,x,P$,\phi)$\ +\ $\underset{\text{s}%
\uparrow\text{t}}{\text{lim}}$F$_{\text{N+1}}$(s,x,P,$\phi$)

We show below that:

$\phi$($\gamma$(t)) = b$_{0}($x,P,$\phi$); $\underset{\text{s}\uparrow
\text{t}}{\text{lim}}$ a$_{n}($s,x,P$,\phi)=$ b$_{n}$(x,P,$\phi$)t$^{n}$;
$\underset{\text{s}\uparrow\text{t}}{\text{lim}}$F$_{\text{N+1}}$(s,x,P,$\phi
$) = R$_{\text{N+1}}($t,x,P,$\phi$)t$^{N+1}$

where b$_{0}($x,P,$\phi$) ; b$_{n}$(x,P,$\phi$); R$_{\text{N+1}}($t,x,P,$\phi
$) are defined above in the statement of the theorem.

(ii)\qquad b$_{0}($x,P,$\phi$) = $\phi$($\gamma$(t)) by definition

Since $\gamma$ is the (unique minimal) geodesic from x$\in$M$_{0}$ to the
centre of Fermi coordinates y$_{0}\in$P in time t, $\gamma($t$)=y_{0}$ and for
any section $\phi$ of the vector bundle E, we have: $\phi$($\gamma$(t)) =
$\phi(y_{0}).$

Since $\phi(x)\in$E$_{\text{x}}$ and $\phi(y_{0})\in$E$_{\text{y}_{0}}$ , we
see that $\phi(y_{0})=$ $\tau_{\text{x,y}_{0}}\phi($x$)$ by the definition of
the parallel propagator $\tau_{\text{x,y}_{0}}$:E$_{\text{x}}\rightarrow
$E$_{\text{y}_{0}}$ and so (ii) proved.

(iii) and (iv) are special cases of (v) and so the proofs are left out and we
prove only (v):

(v)$\qquad$a$_{n}$(s,x,P,$\phi$) $=\int_{0}^{\text{s}}\int_{0}^{\text{s}_{1}%
}...\int_{0}^{\text{s}_{n-1}}($F(t,t-s$_{n}$)L$_{\Psi}$F(t-s$_{n}$,t-s$_{n-1}%
$)L$_{\Psi}$F(t-s$_{n-1}$,t-s$_{n-2}$)

\qquad\qquad\qquad\qquad\qquad... L$_{\Psi}$F(t-s$_{2}$,t-s$_{1}$)L$_{\Psi}%
$F(t-s$_{1}$,t-s)$\phi$)$($x$)$ds$_{1}...$ds$_{n}$

\qquad(F(r,s)$\phi)($x$)=\phi(\gamma($r-s)) = $\phi\circ\gamma($r-s) for
s$\in\lbrack0,r],$

\qquad First we show that:

$\left(  7.3\right)  $ \qquad F(t, s) = F(1,$\frac{\text{s}}{\text{t}}$):

\qquad By the definition of F(t,s) in $\left(  5.10\right)  $, we have:

$\left(  7.4\right)  \qquad$(F(t,s)$\phi)($x$)=\phi(\gamma($t-s)) = $\phi
\circ\gamma($t-s) for s$\in\lbrack0,t],$

where $\gamma$ is the unique minimal geodesic from x$\in$M to the
\textbf{center of Fermi coordinates} y$\in$P in time t. It is the geodesic
$\gamma(s)=\exp_{\text{y}}((1-\frac{s}{t})$v$)$ where x = $\gamma(0)$
$=\exp_{\text{y}}($v$).$

Hence in Fermi coordinates we can write: $\gamma$(s) = x+$\frac{\text{s}%
}{\text{t}}$(y-x) (vector form) and so,

$\left(  7.5\right)  \qquad\gamma$(t-s) = x+$\frac{\text{t-s}}{\text{t}}$(y-x)
= x + (1-$\frac{\text{s}}{\text{t}}$)(y-x)

By definition, (F(1,$\frac{\text{s}}{\text{t}})\phi$)(x) = $\phi(\eta
(1-\frac{\text{s}}{\text{t}}))$ where $\eta$ is the unique minimal geodesic
from x to y in time 1 and hence in Fermi coordinates, $\eta$(s) = x + s(y-x)
(vector form) and so.

$\left(  7.6\right)  \qquad\eta(1-\frac{\text{s}}{\text{t}})$ = x +
(1-$\frac{\text{s}}{\text{t}}$)(y-x)

We see that the RHS,s of $\left(  7.5\right)  $ and $\left(  7.6\right)  $ are
the same. Consequently we have the important relation:

$\left(  7.7\right)  $\qquad$\gamma($t-s$)=\eta($1$-\frac{\text{s}}{\text{t}%
}),$

and as a result:

$\left(  7.8\right)  $\qquad(F(t,s)$\phi$)(x) = $\phi(\gamma$(t-s) $=$ $\phi
$($\eta($1$-\frac{\text{s}}{\text{t}})=$ (F(1,$\frac{\text{s}}{\text{t}})\phi
$)(x$_{0}$)

We thus have the important relation:

$\left(  7.9\right)  $\qquad F(t,s) = F(1,$\frac{\text{s}}{\text{t}})$ and so,
$\left(  7.3\right)  $ above is proved.

We now set: s = tr$_{0\text{, \ \ \ \ \ \ }}$s$_{\text{n}}$ = tr$_{\text{n}}$,
for 1$\leq$n$\leq$N+1 and we have:

$\left(  7.10\right)  $\qquad F(t,t-s$_{n})=$ F(t,t-tr$_{n})=$ F(1,$\frac
{\text{t-tr}_{n}}{\text{t}})=$ F(1,1-r$_{n}$ )Consequently,

$\left(  7.11\right)  $\qquad F(t-s$_{n}$,t-s$_{n-1}$) = F(1,$\frac
{\text{t-tr}_{n-1}}{\text{t-tr}_{\alpha}})=$ F(1,$\frac{\text{1-r}_{n-1}%
}{\text{1-r}_{\alpha}})=$ F($1$-r$_{n},1-$r$_{n-1})\qquad\qquad\qquad
\qquad\qquad\qquad\qquad\qquad\qquad\qquad\ \ \ $

Consequently, noting that ds$_{\text{n}}$ = tdr$_{\text{n}}$, for 1$\leq
$n$\leq$N+1 and 1$\geq$r$_{0}\geq$r$_{1}\geq$r$_{2}$...r$_{\text{N}}\geq
$r$_{\text{N+1}},$ we have:

$\left(  7,12\right)  $\qquad a$_{n}($s,x,P$,\phi)=\int_{0}^{r_{0}}\int%
_{0}^{r_{1}}...\int_{0}^{r_{n-1}}$F(1,1-r$_{n}$)[L$_{\Psi}$F(1-r$_{n}%
$,1-r$_{n-1}$)L$_{\Psi}$F(1-r$_{n-1}$,1-r$_{n-2}$)...

$\ \ \qquad\qquad\qquad\qquad\qquad\ \ $...\ L$_{\Psi}$F(1-r$_{2}$,1-r$_{1}%
$)L$_{\Psi}$F(1-r$_{1}$,1-r$_{0}$))$\phi$](x)t$^{n}$dr$_{1}$...dr$_{n}$

We see that in the expression for a$_{n}($s,x,P$,\phi)$ given in $\left(
6.25\right)  $ above, we have eliminated the variable s on the Left Hand Side
of the equation. This will be very important when we take limits as
s$\uparrow$t on both sides of $\left(  7.12\right)  .$

Give the change of variable:\ s = tr$_{0},$ the limits s$\uparrow$t and
r$_{0}\uparrow$1 are equivalent.

On the RHS we have by the expression for a$_{n}($s,x,P$,\phi)$ given in
$\left(  7.12\right)  :$\ 

$\left(  7.13\right)  \qquad$ $\underset{\text{s}\uparrow\text{t}}{\text{lim}%
}$ a$_{n}($s,x,P$,\phi)$

\qquad\qquad= $\underset{\text{r}_{0}\uparrow1}{\lim}\int_{0}^{r_{0}}\int%
_{0}^{r_{1}}...\int_{0}^{r_{n-1}}$F(1,1-r$_{n}$)L$_{\Psi}[$F(1-r$_{n}%
$,1-r$_{n-1}$)L$_{\Psi}$F(1-r$_{n-1}$,1-r$_{n-2}$)

$\ \ \ \ \qquad\qquad...\ $L$_{\Psi}$F(1-r$_{2}$,1-r$_{1}$)L$_{\Psi}%
$F(1-r$_{1}$,1-r$_{0}$)$\phi$](x$)$t$^{n}$dr$_{1}...$dr$_{n}$

\qquad\qquad= $\int_{0}^{1}\int_{0}^{r_{1}}...\int_{0}^{r_{n-1}}%
\underset{\text{r}_{0}\uparrow1}{\lim}$F(1,1-r$_{n}$)L$_{\Psi}[$F(1-r$_{n}%
$,1-r$_{n-1}$)L$_{\Psi}$F(1-r$_{n-1}$,1-r$_{n-2}$)

$\ \ \ \ \qquad\qquad...\ $L$_{\Psi}$F(1-r$_{2}$,1-r$_{1}$)L$_{\Psi}%
$F(1-r$_{1}$,1-r$_{0}$)$\phi$](x$)$t$^{n}$dr$_{1}...$dr$_{n}$\qquad\qquad

The only factor in the integrand to examine closely is:

L$_{\Psi}[$F(1-r$_{1}$,1-r$_{0}$)$\phi$] and note that L$_{\Psi}[$F(1-r$_{1}%
$,1-r$_{0}$)$\phi]=$ $\Psi^{-1}$L$[\Psi$F(1-r$_{1}$,1-r$_{0}$)$\phi]=\Psi
^{-1}$L$[\Psi\phi\circ\gamma_{0,1}($r$_{0}-$ r$_{1})]$

Since M is compact and $\phi$ and $\Psi$ are smooth, then L$_{\Psi}\phi$ and
the "horde" of derivatives are continuous.

In particular, L$[\Psi\phi\circ\gamma_{0,1}]$ is continuous.

Consequently,

$\left(  7.14\right)  \qquad\underset{\text{r}_{0}\uparrow1}{\lim}\Psi^{-1}%
$L$[\Psi\phi\circ\gamma_{0,1}($r$_{0}-$ r$_{1})]=\Psi^{-1}$L$[\Psi\phi
\circ\gamma_{0,1}($1$-$ r$_{1})]$

where $\gamma_{0,1}$ is the unique minimal geodesic from x of M$_{0}$ to the
point y$\in$P in time 1-r$_{1}.$

Recall that by definition, $\gamma_{0.1}($s$)=$ y + $\left(  1-\frac{\text{s}%
}{\text{1-r}_{1}}\right)  $ (x$-$y), when expressed as a vector

and,

\qquad\qquad\qquad\ \ = $\left(  x_{1},...,x_{q},\text{ }(1-\frac{\text{s}%
}{\text{1-r}_{1}})x_{q+1},...,(1-\frac{\text{s}}{\text{1-r}_{1}})x_{n}\right)
$, when expressed in local coordinates

Therefore,

$\ \gamma_{0,1}($r$_{0}-$ r$_{1})$\ = $\left(  x_{1},...,x_{q},\text{
}(1-\frac{\text{r}_{0}\text{-r}_{1}}{\text{1-r}_{1}})x_{q+1},...,(1-\frac
{\text{r}_{0}\text{-r}_{1}}{\text{1-r}_{1}})x_{n}\right)  $

and so $\underset{\text{r}_{0}\uparrow1}{\lim}\gamma_{0.1}$(r$_{0}$-r$_{1}$) =
$\left(  x_{1},...,x_{q},0,...,0\right)  =\gamma_{0.1}(1$-r$_{1})$ $=$
$\pi_{\text{P}}$(x)

where $\pi_{\text{P}}$ is the projection: $\pi_{\text{P}}$: M$_{0}\rightarrow$
P on P viewed in Fermi coordinates.

Consequently,

$\left(  7.15\right)  \qquad\underset{\text{r}_{0}\uparrow1}{\lim}$L$_{\Psi}%
[$F(1-r$_{1}$,1-r$_{0}$)$\phi]=$ $\underset{\text{r}_{0}\uparrow1}{\lim}%
\Psi^{-1}$L$[\Psi$F(1-r$_{1}$,1-r$_{0}$)$\phi]=\Psi^{-1}$L$[\Psi\phi\circ
\pi_{\text{P}}]=$ L$_{\Psi}[\phi\circ\pi_{\text{P}}]$

Consequently the expression for the general coefficient is given by:

$\left(  7.16\right)  \qquad$b$_{n}($x$_{0}$,P$,\phi)$ = $\int_{0}^{1}\int%
_{0}^{r_{1}}...\int_{0}^{r_{n-1}}$F(1,1-r$_{n}$)L$_{\Psi}[$F(1-r$_{n}%
$,1-r$_{n-1}$)L$_{\Psi}$F(1-r$_{n-1}$,1-r$_{n-2}$)

$\ \ \ \ \qquad\qquad...\ $L$_{\Psi}$F(1-r$_{2}$,1-r$_{1}$)L$_{\Psi}\phi
\circ\pi_{\text{P}}$](x$)$dr$_{1}...$dr$_{n}$

The \textbf{same change of variables} applies to the remainder term
R$_{\text{N+1}}$(t,x,P) and it is similarly computed to give:$\qquad
\qquad\qquad\qquad\qquad\qquad\qquad\qquad\qquad\qquad\qquad\qquad\qquad
\qquad\qquad\qquad\qquad\qquad$

$\left(  7.17\right)  $\qquad R$_{\text{N+1}}$(t,x,P,$\phi$) $=\int_{0}%
^{1}\int_{0}^{\text{r}_{1}}...\int_{0}^{\text{r}_{\text{N}}}($%
Q(t,t-tr$_{\text{N+1}})[$L$_{\Psi}$F(1-r$_{\text{N+1}}$,1-r$_{\text{N}}%
$)L$_{\Psi}$F(1-r$_{\text{N}}$,1-r$_{\text{N-1}}$)

\qquad\ \ \qquad\ ... L$_{\Psi}$F(1-r$_{2}$,1-r$_{1}$)L$_{\Psi}\phi\circ
\pi_{\text{P}}]$(x)dr$_{1}...$dr$_{\text{N+1}}$

and by the definition of $($Q(t,t-tr$_{\text{N+1}}),$ we have the equivalent
version of the formula above:

$\left(  7.18\right)  $\qquad R$_{\text{N+1}}$(t,x,P,$\phi$) = $\int_{0}%
^{1}\int_{0}^{\text{r}_{1}}....\int_{0}^{\text{r}_{\text{N}}}$E$_{\text{x}%
_{0},\text{y}_{0}}$[$\chi_{\zeta>\text{tr}_{\text{N+1}}}\tau_{\text{0,tr}%
_{\text{N+1}}}^{\text{t}}$e$_{\text{tr}_{\text{N+1}}}^{\text{t}}$L$_{\Psi}%
$F(1-r$_{\text{N+1}}$,1-r$_{\text{N}}$)L$_{\Psi}$F(1-r$_{\text{N}}%
$,1-r$_{\text{N-1}}$)

\ \ \ \ \ \ \ \ ... L$_{\Psi}$F(1-r$_{2}$,1-r$_{1}$)L$_{\Psi}\phi\circ
\pi_{\text{P}}$)((x$^{\text{t}}$(tr$_{\text{N+1}}$))exp\{$\int_{0}%
^{\text{tr}_{\text{N+1}}}\frac{\text{L}\Psi}{\Psi}$(x$^{\text{t}}%
$(s))ds\}]dr$_{1}...$dr$_{\text{N+1}}$

\qquad\qquad\qquad\qquad\qquad\qquad\qquad\qquad\qquad\qquad\qquad\qquad
\qquad\qquad\qquad$\qquad\ \blacksquare$

\begin{theorem}
(Asymptotic Expansion Theorem of the Generalized Heat Kernel)
\end{theorem}

\ \ k$_{\text{t}}$(x,P,$\phi$) $=\ \int_{\text{P}}$k$_{\text{t}}$(x,y)$\phi
$(y)$\upsilon_{\text{P}}$(dy)

\qquad\qquad= q$_{\text{t}}$(x,P)$\left[  \phi\text{(}\gamma\text{(t)) \ +
\ }\underset{n=1}{\overset{\text{N}}{\sum}}\text{b}_{n}\text{(x,P,}%
\phi\text{)t}^{n}\text{ + o(t}^{\text{N}}\text{)}\right]  $

where $\phi$($\gamma$(t)) $=$ b$_{0}($x,P,$\phi$) $=$ $\tau_{\text{x,y}_{0}%
}\phi($x$)$

Since M is compact the integrand of the remainder term is bounded. This gives
the asymptotic expansion.\qquad\qquad\qquad\qquad\qquad\qquad\qquad
\qquad\qquad\qquad\qquad\qquad\qquad\qquad\qquad\qquad\qquad\qquad
$\qquad\qquad\qquad\qquad\qquad\qquad\qquad\qquad\qquad\qquad\qquad
\qquad\qquad\qquad\qquad\qquad\qquad\blacksquare$

The usual \textbf{Heat Kernel Expansion} in vector bundles can then be deduced
from the \textbf{Generalized Heat Kernel Expansion:}

\begin{corollary}
(The Heat Kernel Expansion)
\end{corollary}

\ \ k$_{\text{t}}$(x,y$_{0}$,$\phi$) $=$ q$_{\text{t}}$(x,y$_{0}$)$\left[
\text{b}_{0}\text{(x,y}_{0}\text{,}\phi\text{) \ + \ }%
\underset{n=1}{\overset{\text{N}}{\sum}}\text{b}_{n}\text{(x,y}_{0}%
\text{,}\phi\text{)t}^{n}\text{ + o(t}^{\text{N}}\text{)}\right]  $

where b$_{0}($x,y$_{0}$,$\phi$) $=$ $\tau_{\text{x,y}_{0}}\phi($x$)$ and
b$_{n}$(x,y$_{0}$,$\phi$) = b$_{n}$(x,y$_{0}$)$\phi($y$_{0})$

\begin{proof}
.
\end{proof}

In the theorem we take P = $\left\{  \text{y}_{0}\right\}  $ and recover the
ordinary vector bundle heat kernel on the Left Hand Side and its expansion on
the Right Hand Side.\qquad\qquad\qquad\qquad\qquad\qquad\qquad\qquad
\qquad\qquad\qquad\qquad\qquad\qquad\qquad\qquad\ \qquad\ 

\begin{remark}
Berline, Getzler and Vergne, Theorem 2.26 :
\end{remark}

A similar expansion was carried out in the above cited book.

\qquad\qquad\qquad\qquad\qquad\qquad\qquad\qquad\qquad\qquad\qquad\qquad
\qquad\qquad\qquad\qquad\qquad$\blacksquare$

\chapter{Heat Content Expansion}

We will follow \textbf{Gilkey} $\left[  21\right]  ,$ chapter $2$ for the
definition of the \textbf{heat content} and its asymptotics:

\begin{center}
Let E be a smooth real vector bundle over a closed (compact without boundary)
Riemannian manifold M. Let $\phi\in\Gamma(E)$
\end{center}

be the initial (t = 0) tempreture distribution. Then the subsequent tempreture
distribution at time t%
$>$%
0 is given by:

$\bigskip\left(  8.1\right)  \qquad\qquad\qquad\qquad\phi_{\text{t}}$(x)
$=\phi$(x,t) = $\int_{\text{M}}$ k$_{\text{t}}$(x,y)$\phi$(y)$\upsilon
_{\text{M}}$(dy)

where $\upsilon_{\text{M}}$ is the Riemannian volume measure on M and where
$\phi$(x,t) is the unique solution of:

$\qquad\qquad\qquad\qquad\qquad\ \qquad\ \frac{\partial\phi_{\text{t}}%
}{\partial\text{t}}(x)=$ L$\phi_{\text{t}}(x)$\qquad(evolution equation)

$\left(  8.2\right)  $

$\qquad\qquad\qquad\qquad\qquad\qquad\underset{\text{t}\rightarrow0^{+}}{lim}$
$\phi$(x,t) = $\phi$(x) \qquad(initial condition)

for L = $\frac{1}{2}\Delta+X+V$

Let E$^{\ast}$ the dual vector bundle of E and let $\rho\in\Gamma(E^{\ast})$
be the specific heat. The total \textbf{heat energy content} of the Riemannian
manifold M is defined by:

$\left(  8.3\right)  \qquad\qquad\qquad\beta(\phi,\rho,L)($t$)$ =
$\int_{\text{M}}$
$<$%
$\phi$(x,t),$\rho$(x)%
$>$%
$\upsilon_{\text{M}}$(dx)

\qquad There is an asymptotic expansion:$\qquad\qquad\qquad\qquad$

$\left(  8.4\right)  \qquad\ \qquad\qquad\beta(\phi,\rho,L)($t$)=$
$\underset{m=0}{\overset{+\infty}{\sum}}\beta_{m}(\phi,\rho,L)$t$^{\frac{m}%
{2}}$

where the heat \textbf{content asymptotic expansion coefficients} $\beta
_{m}(\phi,\rho,L)$ are locally computable.

By Theorem $\left(  1.4.7\right)  $ of \textbf{Gilkey} $\left[  21\right]  $
they are given by:

$\left(  8.5\right)  \qquad\qquad\beta_{m}(\phi,\rho,L)=\int_{M}\beta_{m}%
^{M}(\phi,\rho,L)(x)\upsilon_{\text{M}}$(dx) $+\int_{\partial M}\beta
_{m}^{\partial M}(\phi,\rho,L)\upsilon_{\partial\text{M}}$(dx)

We have assumed that the compact Riemannian manifold is without boundary i.e.

\begin{center}
$\partial M=\varnothing$
\end{center}

We have in the case here:

$\left(  8.6\right)  \qquad\qquad\beta_{m}(\phi,\rho,L)=\int_{M}\beta_{m}%
^{M}(\phi,\rho,L)(x)\upsilon_{\text{M}}$(dx)

Recall that the generalized heat kernel expansion is given in \textbf{Theorem
4} by:

$\left(  8.7\right)  \qquad\int_{\text{P}}$k$_{\text{t}}$(x,y)$\phi
($y)$\upsilon_{\text{P}}$(dy) \ \ = q$_{\text{t}}$(x,P)$\left[  \tau
_{\text{x,y}_{0}}\phi\text{(x) \ + \ }\underset{n=1}{\overset{\text{N}}{\sum}%
}\text{b}_{n}\text{(x,P,}\phi\text{)t}^{n}\text{ + R}_{\text{N+1}%
}(\text{t,x,P,}\phi\text{)t}^{\text{N+1}}\right]  $

where we recall that q$_{\text{t}}$(x,P) = $($2$\pi t)^{-\frac{\text{n - q}%
}{2}}\Psi($x)exp$\left\{  -\frac{\text{d(x,P)}^{2}}{2\text{t}}\right\}  .$

Comparing equations $\left(  8.1\right)  $ and $\left(  8.7\right)  $ we see
that in order to make comparison with the heat content expansion we must
assume that the submanifold P = M$_{0}$. Equivalently the \textbf{Fermi
coordinates} based at y$_{0}\in$P all become\textbf{ local coordinates }based
at y$_{0}.$ Equivalently, dimP = q = n = dimM$_{0}$. The following conclusions
then follow:

Since x$\in$ P = M$_{0}$ , we have, d(x,M$_{0}$) = 0 and so:

\bigskip$\left(  8.8\right)  $\qquad\qquad q$_{\text{t}}$(x,P) = q$_{\text{t}%
}$(x,M$_{0}$) = $($2$\pi t)^{-\frac{\text{n - n}}{2}}\Psi($x)exp$\left\{
-\frac{\text{d(x,M}_{0}\text{)}^{2}}{2\text{t}}\right\}  =$ $\Psi($x),

Since P and M$_{0}$ coincide, the unique minimal geodesic from the point
x$\in$M$_{0}$ to P = M$_{0}$ in time t is the constant geodesic $\gamma:$
$\gamma($s) = x for 0$\leq$s$\leq$t. Therefore y$_{0}$ and x coincide.

The computations become very simple:

We recall that: L$^{0}=\frac{1}{2}\Delta^{0}+$ \ X $+$ V\ .

Recall that the function $\Psi:$M$\longrightarrow R_{+}$ is defined by:

$\left(  8.9\right)  \qquad\qquad\qquad\Psi($x$)=\theta^{-\frac{1}{2}}%
($x$)\Phi($x$),$

where $\Phi$ and $\theta$ are defined in $\left(  1.5\right)  $ and $\left(
1.6\right)  $ of \textbf{Chapter 1} respectively by:

\qquad\qquad\qquad\qquad$\Phi_{P}($x) = exp$\left\{  \int_{0}^{1}%
<\text{X(}\gamma\text{(s)) , }\dot{\gamma}\text{(s)%
$>$%
ds}\right\}  $

and,

$\qquad\qquad\qquad\qquad\theta($x$)=$ det(dexp$_{\text{x}_{0}}$) =
$\sqrt{\text{det(g}_{\alpha\beta}\text{(x))}}$\qquad

There is no \textbf{normal} part of the Fermi coordinates (now reduced to
local coordinates), equivalently all normal coordinates are zero. We deduce
from the expansion of $\theta($x$)$ in \textbf{Proposition }$10$ that:

$\left(  8.10\right)  \qquad\qquad\qquad\qquad\qquad\theta($x$)=1.$

Alternatively, we can use \textbf{Proposition }$6$ to see that g$_{\text{ab}%
}($x$)=\delta_{\text{ab}}$ and so

$\qquad\ \qquad\qquad\qquad\qquad\qquad\ \theta($x$)=$ $\sqrt{\text{det(g}%
_{\text{ab}}\text{(x))}}=$ $\sqrt{\text{det(}\delta_{\text{ab}}\text{)}}=1.$

Consequently by $\left(  8.9\right)  $,

$\left(  8.11\right)  \qquad\qquad\qquad\qquad\ \ \Psi($x$)=\Phi($x$)$

On the other hand $\Phi_{P}($x) is defined by:

$\qquad\qquad\qquad\ \ \ \ \ \ \ \ \ \ \ \ \ \ \ \ \ \ \Phi_{P}($x) =
exp$\left\{  \int_{0}^{1}<\text{X(}\gamma\text{(s)) , }\dot{\gamma}\text{(s)%
$>$%
ds}\right\}  $

where $\gamma:[0,1]\longrightarrow M$ $_{0}$ is the unique minimal geodesic
from x$\in$M$_{\text{0}}=$ M to a point y$\in$P = M in time 1.

However, as already pointed out above, Since P and M coincide, the unique
minimal geodesic from a point x of M to P = M is the constant geodesic
$\gamma:$ $\gamma(s)=x$ for 0$\leq$s$\leq$1. Consequently $\overset{\cdot
}{\gamma}(s)=0$ for all s$\in\left[  0,\text{1}\right]  $. We conclude that,

$\qquad\qquad\qquad\qquad\ \ \ \Phi_{P}($x) = exp$\left\{  \int_{0}%
^{1}<\text{X(}\gamma\text{(s)) , }\dot{\gamma}\text{(s)%
$>$%
ds}\right\}  =1$

and so by $\left(  8.9\right)  ,$ $\left(  8.10\right)  $, $\left(
8.11\right)  $ and the last equality above, we have:

$\left(  8.12\right)  \qquad\qquad\qquad\qquad\qquad\Psi($x$)=\theta
^{-\frac{1}{2}}($x$)\Phi($x$)=1$ for all x$\in$M.

Consequently by $\left(  8.8\right)  $ and $\left(  8.12\right)  $:

$\left(  8.13\right)  \qquad\ \ \ \ \qquad$ \qquad\ \ \ \ \ \qquad
q$_{\text{t}}$(x,M$_{0}$) = 1

As a consequence, $\left(  8.7\right)  $ becomes (with M$_{0}$ replacing P):

Since vol(M$_{0}$) = vol(M), we have:

$\left(  8.14\right)  \qquad\qquad\qquad\phi$(x,t) $=\int_{\text{M}}%
$k$_{\text{t}}$(x,y)$\phi$(y)$\upsilon_{\text{M}}$(dy) $=\int_{\text{M}_{0}}%
$k$_{\text{t}}$(x,y)$\phi($y)$\upsilon_{\text{M}}$(dy)

\qquad\qquad\ \qquad\qquad= $\left[  \tau_{\text{x,y}}\phi\text{(x) \ +
\ }\underset{n=1}{\overset{\text{N}}{\sum}}\text{b}_{n}\text{(x,M}_{0}%
\text{,}\phi\text{)t}^{n}\text{ + R}_{\text{N+1}}\text{(t,x,M}_{0}%
\text{)t}^{\text{N+1}}\right]  $

From $\left(  8.1\right)  $ and $\left(  8.3\right)  ,$ we deduce that the
total \textbf{heat energy content} of the Riemannian manifold M is given by:

$\left(  8.15\right)  \qquad\qquad\qquad\beta(\phi,\rho,L)($t$)$ =
$\int_{\text{M}}$
$<$%
$\phi$(x,t),$\rho$(x)%
$>$%
$\upsilon_{\text{M}}$(dx)

\qquad\qquad\qquad= $\int_{\text{M}}$ $\left[  \tau_{\text{x,y}}\phi\text{(x)
\ + \ }\underset{n=1}{\overset{\text{N}}{\sum}}\text{b}_{n}\text{(x,M}%
_{0}\text{,}\phi\text{)t}^{n}\text{ + R}_{\text{N+1}}(\text{t,x,M}_{0}%
\text{,}\phi\text{)t}^{\text{N+1}}\right]  $,$\rho$(x)%
$>$%
$\upsilon_{\text{M}}$(dx)

Since the geodesic $\gamma$ is constant $\gamma(0)=x$ and $\gamma(1)=y$
coincide and so, $\tau_{\text{x,y}}=\tau_{\text{x,x}}$:E$_{\text{x}%
}\rightarrow$E$_{\text{x}}$ is the identity operator and so by $\left(
8.15\right)  ,$

$\left(  8.16\right)  \qquad\qquad\qquad\beta(\phi,\rho,L)($t$)=$
$\underset{n=0}{\overset{+\infty}{\sum}}\beta_{2n}(\phi,\rho,L)$t$^{n}$

\qquad\qquad\qquad= $\int_{\text{M}}$
$<$%
$\left[  \tau_{\text{x,x}}\phi\text{(x) \ + \ }%
\underset{n=1}{\overset{\text{N}}{\sum}}\text{b}_{n}\text{(x,M}_{0}%
\text{,}\phi\text{)t}^{n}\text{ + R}_{\text{N+1}}(\text{t,x,M}_{0}\text{,}%
\phi\text{)t}^{\text{N+1}}\right]  $,$\rho$(x)%
$>$%
$\upsilon_{\text{M}}$(dx)

We see that the expansion coefficients b$_{n}$(x,M$_{0}$,$\phi$) are
\textbf{independent }of t.

Comparing the sum in $\left(  8.4\right)  $ and the sum in $\left(
8.16\right)  ,$ we see that the odd coefficients vanish and $\left(
8.16\right)  $ is the sum of even (m = 2n) coefficients:

\begin{theorem}
(Heat Content Expansion)
\end{theorem}

For $n\geq0,$the general \textbf{heat content expansion general coefficients}
$\beta_{n}(\phi,\rho,L)$ are given by:

(i) $L^{0}=$ $\tau_{\text{x,x}}$

where $L^{0}=$ $\tau_{\text{x,x}}$ is the Identiy parallel propagator along
the fiber E$_{\text{x}}$ of the vector bundle E

(ii)$\qquad\qquad\qquad\qquad\beta_{2n}(\phi,\rho,L)=\frac{1}{n!}%
\int_{\text{M}}$ $<(L^{n}\phi)($x$),\rho($x$)>\upsilon_{\text{M}}$(dx)

for where the operator $L^{n}=LL....L$ (n-times).

(iii)\qquad$\qquad\qquad\qquad\beta_{2n+1}(\phi,\rho,L)=0$

(iv) The \textbf{Remainder Term }is given by:

\qquad R$_{\text{N+1}}$(t,x,M$_{0}$,$\phi$) $=$ R$_{\text{N+1}}$%
(t,x,M,$\phi)=\frac{1}{N!}\int_{0}^{\text{r}_{\text{N}}}$\textbf{E}%
$_{\text{x}}$[$\tau_{\text{0,t}}$e$_{\text{t}}$(L$^{N+1}\phi$%
)(x(tr$_{\text{N+1}}$))exp\{$\int_{0}^{\text{tr}_{\text{N+1}}}$%
V(x(s))ds\}]dr$_{\text{N+1}}$

$\qquad\qquad\qquad\qquad\ \ \ =\frac{1}{N!}\int_{0}^{\text{r}_{\text{N}}}%
($Q(t,t-tr$_{\text{N+1}})$[L$^{N}\phi$]$($x$)$dr$_{\text{N+1}}$

where the operator $L^{N}=LL....L$ (N-times).

\begin{proof}
(i) From the expansion in $\left(  8.16\right)  $ we see that
\end{proof}

(i) The formula for $n=0$ has already been shown in the computation of the
operator $L^{0}=\tau_{\text{x,y}}=\tau_{\text{x,x}}$:E$_{\text{x}}\rightarrow
$E$_{\text{x}}$

(ii) We now prove that it is true for all $n\geq1:$

As pointed out earlier, the points x and y$_{0}$ coincide because P = M$_{0}$.
Consequently all the geodesics from x to y$_{0}$ are constant geodesics in
their various time intervals: \ The geodesic: $\gamma_{n,n-1}:\left[
0,\text{1-r}_{n}\right]  \longrightarrow$M$_{0}$ is constant in the sense that
for 1$\geq$r$_{1}\geq$r$_{2}\geq...\geq$r$_{n-1}\geq$r$_{n}\geq...$%
r$_{\text{N}}\geq$r$_{\text{N+1}}$:

$\left(  8.17\right)  \qquad\qquad\qquad\gamma_{n,n-1}(s)=$ x$\in M_{0}$ for
all $s\in\left[  0,\text{1-r}_{n}\right]  $.

By the definition of the operators F(r,s) given in $\left(  5.10\right)  $ we
conclude that the corresponding operators F(1-r$_{n}$,1-r$_{n-1}$) are
identity operators for 1$\leq n\leq N$ and so the integrand is independent of
r$_{1},$ r$_{2}$,$...,$r$_{n-1},$ r$_{n},...,$r$_{\text{N}}$ and hence by the
definition of the operators F(1-r$_{N}$,1-r$_{N-1}$), we have:

$\left(  8.18\right)  \qquad$F(1,1-r$_{N}$)[L$_{\Psi}$F(1-r$_{N}$,1-r$_{N-1}%
$)L$_{\Psi}$F(1-r$_{N-1}$,1-r$_{N-2}$)$...$F(1-r$_{n}$,1-r$_{n-1}$)

$\ \ \ \ \ \ \ \qquad\qquad\qquad\ ...$L$_{\Psi}$F(1-r$_{2}$,1-r$_{1}%
$)L$_{\Psi}\phi\circ\gamma_{\text{0,1}}$(1-r$_{1}$)]$($x$)$

$\qquad\qquad\qquad\qquad=$ [L$_{\Psi}$L$_{\Psi}\ ...$L$_{\Psi}$L$_{\Psi}\phi
$]$($x$)$

Integrating fron Right to Left it is obvious that:

$\left(  8.19\right)  \qquad\qquad\qquad\qquad$ $\int_{0}^{1}\int_{0}^{r_{1}%
}...\int_{0}^{r_{n-1}}$dr$_{1}...$dr$_{n}=\frac{1}{n!}$

In the light of $\left(  8.18\right)  $ and $\left(  8.19\right)  $ the
general coefficient simply becomes:

$\left(  8.20\right)  \qquad$b$_{n}($x,M$_{0}$,$\phi)=\int_{0}^{1}\int%
_{0}^{r_{1}}...\int_{0}^{r_{n-1}}$[L$_{\Psi}$L$_{\Psi}\ ...$L$_{\Psi}$%
L$_{\Psi}\phi$]$($x$)$dr$_{1}...$dr$_{n}$

\qquad\qquad\qquad\qquad\ \ \qquad$\qquad=\frac{1}{n!}$[L$_{\Psi}$L$_{\Psi
}\ ...$L$_{\Psi}$L$_{\Psi}\phi$]$($x$)$\qquad\qquad\qquad\qquad\qquad

Recall that by definition the operator L$_{\Psi}$ on $\Gamma(E)$ is given by
L$_{\Psi}\phi$ = $\frac{\text{L(}\Psi\phi\text{)}}{\Psi}$

Since by $\left(  8.13\right)  $, $\Psi($x$)=1$ for all x$\in$M$_{0}$ we have:

\qquad\qquad\qquad L$_{\Psi}\phi$ $=$ L$\phi$ and so L$_{\Psi}$ $=$ L and so
$\Psi$ is eliminated from L$_{\Psi}.$

Consequently by $\left(  8.20\right)  ,$

$\left(  8.21\right)  \qquad\qquad\qquad$b$_{n}($x,M$_{0}$,$\phi)$ $=\frac
{1}{n!}$[L$_{\Psi}$L$_{\Psi}$\ ...L$_{\Psi}$L$_{\Psi}\phi$](x)

$\qquad\qquad\qquad\qquad\qquad\qquad\qquad\ =\frac{1}{n!}$[LL$\ ...$%
LL]$\phi($x$)$

\qquad\qquad$\qquad\qquad\ \ \qquad\qquad\ \ \ \ \ \ =\frac{1}{n!}($L$^{n}%
\phi)($x$)$

where L$^{n}=$ L$...$L ($n-$ times)

From $\left(  8.16\right)  $ the general coefficient in the \textbf{heat
content expansion} is given for \textbf{even} coefficients by:

$\left(  8.22\right)  \qquad\qquad\qquad\beta_{2n}(\phi,\rho,L)=\int%
_{\text{M}}$
$<$%
b$_{n}$(x,M$_{0}$,$\phi$),$\rho$(x)%
$>$%
$\upsilon_{\text{M}}$(dx)

We conclude from $\left(  8.21\right)  $ and $\left(  8.22\right)  $ that the
general \textbf{even} coefficient term in the \textbf{heat content expansion}
is given by:

$\left(  8.23\right)  \qquad\qquad\qquad\beta_{2n}(\phi,\rho,L)=\frac{1}%
{n!}\int_{\text{M}}$
$<$%
$(L^{n}\phi)($x$)$,$\rho$(x)%
$>$%
$\upsilon_{\text{M}}$(dx)

where we recall that the differential operator L is given by $L=$ $\frac{1}%
{2}\Delta+$ X $+$ V and $\Delta$ is the Laplace-Type operator (or generalized Laplacian).

We see from the relation in $\left(  8.16\right)  $ and $\left(  8.22\right)
$ that the general \textbf{odd }coefficient term is given by:

$\left(  8.24\right)  \qquad\qquad\qquad\beta_{2n+1}(\phi,\rho,L)=0$

The theorem is thus proved for 1$\leq n\leq N$ and hence for 0$\leq n\leq N$.

Since $N$ is arbitrary, the theorem is proved for all integers $n\geq0.$

(iii) This is obvious from $\left(  8.14\right)  $ and $\left(  8.16\right)
.$

(iv) The corresponding Remainder Term is given by:

\qquad R$_{\text{N+1}}$(t,x,M$_{0}$,$\phi)=\int_{0}^{1}\int_{0}^{\text{r}_{1}%
}....\int_{0}^{\text{r}_{\text{N}}}$\textbf{E}$_{\text{x}}$[$\chi
_{\zeta>\text{tr}_{\text{N+1}}}\tau_{\text{0,t}}^{\text{t}}$e$_{\text{t}%
}^{\text{t}}$L$_{\Psi}$F(1-r$_{\text{N+1}}$,1-r$_{\text{N}}$)L$_{\Psi}%
$F(1-r$_{\text{N}}$,1-r$_{\text{N-1}}$)

\ \ \ \ \ \ \ \ ... L$_{\Psi}$F(1-r$_{2}$,1-r$_{1}$)L$_{\Psi}\phi(\gamma
_{0.1}$(1-r$_{1}$)$)($(x$^{\text{t}}$(tr$_{\text{N+1}}$))exp\{$\int%
_{0}^{\text{tr}_{\text{N+1}}}\frac{\text{L}^{0}\Psi}{\Psi}$(x$^{\text{t}}%
$(s))ds\}]dr$_{1}...$dr$_{\text{N+1}}$

Since $\Psi=$ 1, we have L$_{\Psi}$ $=$ L = $\frac{1}{2}\Delta+X+V$ and
L$^{0}=\frac{1}{2}\Delta^{0}+$ X $+$ V and so,

\begin{center}
$\frac{\text{L}^{0}\Psi}{\Psi}=$ V.
\end{center}

Consequently,

$\left(  8.25\right)  \qquad$ R$_{\text{N+1}}$(t,x$_{0}$,M$_{0}$,$\phi
)=\int_{0}^{1}\int_{0}^{\text{r}_{1}}...\int_{0}^{\text{r}_{\text{N}}}%
($Q(t,t-tr$_{\text{N+1}})[$L$_{\Psi}$L$_{\Psi}...$L$_{\Psi}$L$_{\Psi}\phi$%
]$($x$)$dr$_{1}...$dr$_{\text{N+1}}$

\begin{center}
$=\int_{0}^{1}\int_{0}^{\text{r}_{1}}...\int_{0}^{\text{r}_{N-1}}\int%
_{0}^{\text{r}_{\text{N}}}($Q(t,t-tr$_{\text{N+1}})$[L$^{N+1}\phi$]$($%
x$)$dr$_{1}...$dr$_{\text{N}}$dr$_{\text{N+1}}$
\end{center}

The integrand is independent of r$_{1},...,$r$_{N}$ and so using $\left(
8.19\right)  $,

\qquad\qquad R$_{\text{N+1}}$(t,x,M$_{0}$,$\phi)=\int_{0}^{1}\int%
_{0}^{\text{r}_{1}}...\int_{0}^{\text{r}_{N-1}}$dr$_{1}$...dr$_{\text{N}}%
\int_{0}^{\text{r}_{\text{N}}}($Q(t,t-tr$_{\text{N+1}})$[L$^{N+1}\phi$]$($%
x$)$dr$_{\text{N+1}}$

Since,

\begin{center}
$\int_{0}^{1}\int_{0}^{\text{r}_{1}}...\int_{0}^{\text{r}_{N-1}}$dr$_{1}%
$...dr$_{\text{N}}=\frac{1}{N!},$
\end{center}

and using the definition of $($Q(t,t-s$),$ we have:

$\left(  8.26\right)  \qquad\qquad\qquad$R$_{\text{N+1}}$(t,x,M$_{0}$%
,$\phi)=\frac{1}{N!}\int_{0}^{\text{r}_{\text{N}}}($Q(t,t-tr$_{\text{N+1}}%
)$[L$^{N+1}\phi$]$($x$)$dr$_{\text{N+1}}$

$\qquad\qquad=\frac{1}{N!}\int_{0}^{\text{r}_{\text{N}}}$\textbf{E}%
$_{\text{x}}$[$\chi_{\zeta>\text{tr}_{\text{N+1}}}\tau_{\text{0,t}}%
$e$_{\text{t}}$[L$^{N+1}\phi$](x(tr$_{\text{N+1}}$))exp\{$\int_{0}%
^{\text{tr}_{\text{N+1}}}$V(x(s))ds\}]dr$_{\text{N+1}}$

where $\zeta$ is the first exit time from M$_{0}$ of Brownian motion x(s)
0$\leq$ s$\leq$t$\Lambda\varsigma$ with differential generator L = $\frac
{1}{2}\Delta+X+V.$

The last equality above is due to (i) of \textbf{Corollary 1} of
\textbf{Theorem 1} here.

Since vol(M$_{0}$) = vol(M), we can replace M$_{0}$ by M under the integral
sign in $\left(  8.26\right)  .$

Since the Riemannian manifold M is compact, the explosion time of Brownian
motion on M is $\zeta=+\infty.$ Consequently,

\qquad R$_{\text{N+1}}$(t,x,M,$\phi)=\frac{1}{N!}\int_{0}^{\text{r}_{\text{N}%
}}$\textbf{E}$_{\text{x}}$[$\tau_{\text{0,t}}$e$_{\text{t}}$[L$^{N+1}\phi
$](x(tr$_{\text{N+1}}$))exp\{$\int_{0}^{\text{tr}_{\text{N+1}}}$%
V(x(s))ds\}]dr$_{\text{N+1}}$

\begin{center}
\qquad\qquad\qquad\qquad\qquad\qquad\qquad\qquad\qquad\qquad\qquad
$\blacksquare$
\end{center}

The Heat Content Expansion of \textbf{Theorem 6 } above is \textbf{Theorem
1.3.12} of \textbf{Gilkey} $\left[  21\right]  $ in the case of the second
order operator L = $\frac{1}{2}\Delta+$ X $+$ V where $\Delta$ is the
Laplace-Type operator, X a vector field on M and V a potenial term.

The extra coefficient of $(-1)^{n}$ in (2) of Theorem (1.3.12) of
\textbf{Gilkey} $\left[  21\right]  $ is absent in our formula of $\left(
8.23\right)  $ above. This is due to the fact that we are dealing with the
forward heat equation: $\frac{\partial}{\partial t}=D=\frac{1}{2}\Delta+$ X
$+$ V here, whereas in some litterature, and \textbf{Gilkey }$\left[
1\right]  ,$ $\left[  2\right]  $ in particular, the authors deal with the
backward heat equation $\frac{\partial}{\partial t}=-$ D where D is a
differential operator of some order d$\geq2$.

Consequently in our case D would have been given by D $=-$ $(\frac{1}{2}%
\Delta+$ X $+$ V$)=-L.$

This would then imply that: $D^{n}=(-1)^{n}L^{n}.$

We conclude that our \textbf{generalized heat kernel expansion} is a double
generalization of the \textbf{heat trace} and the \textbf{heat content
expansions}, at least, in the case of closed Riemannian manifolds.

\qquad\qquad\qquad\qquad\qquad\qquad\qquad\qquad\qquad\qquad\qquad\qquad
\qquad\qquad\qquad\qquad\qquad\qquad\qquad\qquad\qquad$\blacksquare$

\begin{theorem}
(Expansion of the Solution of the Feynman-Kac Formula)
\end{theorem}

For $\phi\in\Gamma(E),$ we have for 0$\leq$ t$<+\infty$

\begin{center}
\bigskip\qquad\ $\phi_{\text{t}}$(x) = P$_{\text{t}}\phi$(x) = $\int%
_{\text{M}}$k$_{\text{t}}$(x,y)$\phi($y)$\upsilon_{\text{M}}$(dy) =
$\int_{\text{M}_{0}}$k$_{\text{t}}$(x,y)$\phi($y)$\upsilon_{\text{M}}$(dy)

= \textbf{E}$_{\text{x}}\left[  \tau_{\text{0,t}}^{\text{t}}\text{e}%
_{\text{t}}^{\text{t}}\phi\text{(x(t))exp}\left\{  \int_{0}^{\text{t}%
}\text{V(x}^{\text{t}}\text{(r))dr}\right\}  \right]  $\qquad

\qquad\ \ \ \ \ = $\tau_{\text{x,x}}\phi$(x) \ +
\ $\underset{n=1}{\overset{\text{N}}{\sum}}\frac{1}{n!}$(t$^{n}$L$^{n}$)$\phi
$(x) + R$_{\text{N+1}}$(t,x,M,$\phi$)t$^{\text{N+1}}$
\end{center}

\begin{proof}
The first equality is just notation given in $\left(  8.1\right)  .$ The
second equality is due to the definition of P$_{\text{t}}\phi$(x) in $\left(
4.5\right)  .$
\end{proof}

The third is due to the fact that vol(M) = vol(M$_{0}$)

The\ fourth equality is due to $\left(  4.39\right)  $ above. Then $\left(
8.14\right)  $ and $\left(  8.21\right)  $ combine to give the last
equality$.$ \qquad\qquad

\begin{center}
\qquad\qquad\qquad\qquad\qquad\qquad\qquad\qquad\qquad\qquad\qquad
$\blacksquare$
\end{center}

We note that the last equality is symbolic because L$^{n}$ is the n$^{th}$
derivative and not the n$^{th}$ power!

By $\left(  5.47\right)  $ and $\left(  8.25\right)  ,$ the Remainder term is
given by:

\qquad R$_{\text{N+1}}$(t,x,M,$\phi$) $=\frac{1}{N!}\int_{0}^{\text{r}%
_{\text{N}}}$\textbf{E}$_{\text{x}}$[$\tau_{\text{0,s}}^{\text{t}}%
$e$_{\text{s}}^{\text{t}}$[L$^{N+1}\phi$](x$^{\text{t}}$(tr$_{\text{N+1}}%
$))exp\{$\int_{0}^{\text{tr}_{\text{N+1}}}$V(x$^{\text{t}}$%
(s))\}]dr$_{\text{N+1}}$

Since M is compact, R$_{\text{N+1}}$(t,x,M) is uniformly bounded in (t,x) and hence,

$\qquad\qquad\underset{\text{N}\rightarrow+\infty}{lim}$R$_{\text{N+1}}%
$(t,x,M,$\phi$)t$^{\text{N+1}}=0.$

We then conclude that:

\qquad\qquad P$_{\text{t}}\phi(x)$ = $\underset{n=0}{\overset{\text{+}%
\infty}{\sum}}\frac{1}{n!}(tL)^{n}\phi(x)=$ e$^{\text{tL}}\phi(x)$

and so we can write: P$_{\text{t}}$ = $\underset{n=0}{\overset{\text{+}%
\infty}{\sum}}\frac{1}{n!}(tL)^{n}=$ e$^{\text{tL}}$ where L$^{0}%
=\tau_{\text{x,x}}$

\qquad\qquad\qquad\qquad\qquad\qquad\qquad\qquad\qquad\qquad$\qquad
\qquad\qquad\qquad\qquad\qquad\blacksquare$

We summarize the above as:

$\left(  8.27\right)  \qquad\qquad\phi_{\text{t}}(x)=P_{\text{t}}%
\phi(x)+\ \underset{n=1}{\overset{+\infty}{\sum}}\frac{1}{n!}(L^{n}%
\phi)(x)t^{n}=$ $\underset{n=0}{\overset{+\infty}{\sum}}\frac{1}{n!}%
(tL)^{n}\phi(x)=$ e$^{\text{tL}}\phi(x)$

where L$^{0}=\tau_{\text{x,x}}$ is the Identiy parallel propagator along the
fiber E$_{\text{x}}$ of the vector bundle E.\qquad

The equation: P$_{\text{t}}$ = $\underset{n=0}{\overset{\text{+}\infty}{\sum}%
}\frac{1}{n!}(tL)^{n}=$ e$^{\text{tL}}$ is the reason why the
\textbf{solution} $\phi_{t}$ of the forward heat equation:

$\left(  8.28\right)  $\qquad\ $\ \frac{\partial\phi_{\text{t}}}%
{\partial\text{t}}(x)=$ L$\phi_{\text{t}}(x)$\qquad(evolution equation)

\qquad$\ \ \qquad\qquad\phi_{0}(x)$ = $\phi$(x) \qquad(initial condition)

is sometimes written as:

$\left(  8.29\right)  \qquad\phi_{\text{t}}(x)=P_{\text{t}}\phi
(x)=e^{\text{tL}}\phi(x)=$ $\int_{\text{M}}$ k$_{\text{t}}$(x,y)$\phi
$(y)$\upsilon_{\text{M}}$(dy)

where $\upsilon_{\text{M}}$ is the Riemannian volume measure on M.$\qquad
\qquad\ \qquad$\qquad\qquad\qquad\qquad\qquad\qquad\qquad\qquad\qquad\qquad

Compare our formula in $\left(  8.27\right)  $ with the formula of
\textbf{Berline, Getzeler and Vergne }$\left[  7\right]  $ of
\textbf{Proposition }$\left(  2.13\right)  $ and the expansion formula in
$\left(  2.8\right)  $ of \textbf{Theorem }$\left(  2.30\right)  $ in
\textbf{Berline, Getzeler and Vergne }$\left[  7\right]  .$

The slight difference between our formula here and theirs is that there is the
factor $(-1)^{n}$ in their formula which is due to their use of the
\textbf{backward heat equation}, in contrast to our \textbf{forward heat
equation} in $\left(  8.28\right)  $ above.

Our formula in $\left(  8.27\right)  $ above is the same as the formula on p.6
in \textbf{Baudoin} $\left[  5\right]  .$ This is because Baudoin is using the
forward heat equation as we have done in here.

Beautiful as these results look, the downside is that the expansion
coefficients here do not exhibit the geometric invariants of the underlying
Riemannian manifold M and the vector bundle E.

\qquad\qquad\qquad\qquad\qquad\qquad\qquad\qquad\qquad\qquad\qquad\qquad
\qquad\qquad\qquad\qquad\qquad$\blacksquare$

\begin{remark}
It is possible to compute $\frac{1}{2}\Delta=\frac{1}{2}(\nabla_{\partial_{i}%
}^{E}\nabla_{\partial_{j}}^{E}-\Gamma_{ij}^{k}\nabla_{\partial_{k}}^{E})+W$ in
terms of the geometric invariants of the Riemannian manifold M and the vector
bundle E and hence, compute
\end{remark}

\qquad\qquad\ $L^{n}=\left(  \frac{1}{2}\Delta+X+V\right)  ^{n}$ and
P$_{\text{t}}\phi=$ $\underset{n=0}{\overset{\text{+}\infty}{\sum}}\frac
{1}{n!}(tL)^{n}\phi=$ e$^{tL}\phi$ in terms of geometric invariants.

We note that we have to use the\textbf{ local coordinates }and not
\textbf{normal coordinates}.

We also note that the exponential function above is not the usual one since by
$\left(  8.21\right)  $ above, the $n$ in L$^{n}$ is a differential order and
\textbf{not} a multiplication order in the usual sense.

\qquad\qquad\qquad\qquad\qquad\qquad\qquad\qquad\qquad\qquad\qquad\qquad
\qquad\qquad\qquad\qquad\qquad$\blacksquare$

It is important here to highlight the roles the expression for (Q(t,t-s)$\phi
$)(x) has played in this work:

1. It enabled us to establish the \textbf{Generalized Feynman-Kac Formula} in
\textbf{Theorem 1} which, in turn, enabled us to derive the following:

2. (i) It enabled us to use the \textbf{Generalized Feynman-Kac Formula} as
solution of the Heat Equation in \textbf{Corollary 1}.

\ \ \ (ii) We then deduced the usual \textbf{Feynman-Kac \ Formula} in vector
bundles given in \textbf{Corollary 2.}

3 The \textbf{Generalized Elworthy-Truman Heat Kernel Formula} given in
\textbf{Corollary 3}.

Then come the expansions: Using \textbf{Lemma 4}, we obtained the following expansions:

4. Expansion of the \textbf{Generalized Feynman-Kac Formula }of\textbf{
Theorem 3}.

5. The \textbf{Generalized Heat Kernel Expansion} of \textbf{Theorem 4}

6. The \textbf{Heat Content Expansion} of \textbf{Theorem 6}

7. The \textbf{Expansion of the Solution Feynman-Kac Formula} of
\textbf{Theorem 7}

\qquad\qquad\qquad\qquad\qquad\qquad\qquad\qquad\qquad\qquad\qquad\qquad
\qquad\qquad\qquad$\qquad\blacksquare$

\part{ DIFFERENTIAL GEOMETRIC BACKGROUND}

\chapter{Expansions in Fermi Coordinates}

\begin{definition}
(normal Fermi coordinates, tangential Fermi coordinates)
\end{definition}

Let (x$_{1},..,$x$_{q},$x$_{q+1},$..., x$_{n})$ be the Fermi coordinates
defined in section 2. Then x$_{1},..,$x$_{q}$ are called \textbf{tangential
Fermi coordinates} and x$_{q+1},$..., x$_{n}$ are called \textbf{normal Fermi
coordinates}.

In this section we discuss expansions in normal Fermi coordinates relative to
the submanifold P of M. We defined $\theta_{P}=$ $\sqrt{\text{det(g}%
_{ij}\text{)}}$ where g$_{ij}$ are components of the Riemannian metric tensor
in Fermi coordinates$.$ Power series expansions in\ (normal) Fermi coordinates
have been discussed and used by \ A. Gray and L. Vanhecke in \textbf{Gray and
Vanhecke} $\left[  2\right]  $. In their paper an explicit expansion of
$\theta_{P}$ was given (only up to the fourth term). However here we will need
the fifth term. Also we need the expansions of the various components of the
metric tensor field defined in Fermi coordinates. We follow \textbf{Gray and
Vanhecke }$\left[  2\right]  $ for the following definitions:

\begin{definition}
(normal Fermi vector \ fields, \ tangential Fermi vector fields).
\end{definition}

\ Let (x$_{1},$...,x$_{q},$...,x$_{n})$ be a Fermi coodinate in a (small)
neighbourhood of y$_{0}\in$P.

(a) A vector field \ X given by:

\qquad\qquad\qquad\qquad\qquad X = $\underset{j=q+1}{\overset{n}{\sum}}$%
a$_{j}\frac{\partial}{\partial x_{j}}$

\qquad where the a$_{j},s$ are constant, is called a \textbf{normal Fermi}
vector field.

(b) \ A vector field given by:

\qquad\qquad\qquad\qquad\qquad A = $\underset{\text{a}=1}{\overset{\text{q}%
}{\sum}}$b$_{\text{a}}\frac{\partial}{\partial x_{\text{a}}}$

\qquad where the b$_{\alpha},s$ are constants, is called a \textbf{tangential
Fermi} vector field.

\section{Some Lemmas of Gray and Vanhecke}

\qquad We will denote the Riemannian (Levi-Civita) connection on M by
$\nabla.$ This is not to be confused with the same notation used for the
connection on the vector bundle E. The context in each case will make the distinction.

Following \textbf{Gray} $\left[  23\right]  ,\left[  25\right]  ,$ the
\textbf{second fundamental form} operator T and the \textbf{torsion operator}
$\bot$ are defined as follows:

Let A be a tangential Fermi field and X a normal \ Fermi field. Then,
$\perp_{\text{A}}$X is normal component of $\nabla_{\text{A}}$X and
T$_{\text{A}}$X is the tangential component of $\nabla_{\text{A}}$X. T and
$\perp$ are related by:

\qquad\qquad\qquad\qquad\qquad\qquad\qquad\qquad\qquad\qquad\qquad\qquad
\qquad\qquad\qquad\qquad\qquad$\qquad\blacksquare$

\begin{lemma}
(A. Gray and L. Vanhecke $\left[  27\right]  ,$ Lemma (3.5)).
\end{lemma}

\vspace{1pt} \ \ \ ($\nabla_{\text{A}}$X)$(y)=$ (T$_{\text{A}}$X \ +
$\perp_{\text{A}}$X)$(y)$\qquad

where T$_{\text{A}}$X is the tangential component of $\nabla_{\text{A}}$X and
$\perp_{\text{A}}$X is its normal component.

\qquad\qquad\qquad\qquad\qquad\qquad\qquad\qquad\qquad\qquad\qquad\qquad
\qquad\qquad\qquad\qquad\qquad\qquad$\blacksquare$

A. Gray and L. Vanhecke use the letter T to denote the \textbf{second
fundamental form }operator. This is in contrast to generally accepted notation
that uses the letter T to denote the \textbf{torsion} operator. More
unfortunate is the fact that what is called the torsion operator $\bot$ here
is \textbf{not the usual torsion} operator in differential geometry.

We use the notation as given above because \textbf{Gray} $\left[  24\right]
,\left[  26\right]  $ are our main references here.

The formulae defined in the lemmas below will be used to compute certain
geometric invariants, which in turn will be used to compute the expansion
coefficients of the generalized heat kernel.

\qquad\qquad\qquad\qquad\qquad\qquad\qquad\qquad\qquad\qquad\qquad\qquad
\qquad\qquad\qquad\qquad$\qquad\qquad\blacksquare$\qquad\qquad\qquad
$\ \ \ \ $\qquad\qquad\ \ \ \ \ \qquad\qquad\qquad\ \ \ \qquad\qquad
\qquad\qquad\qquad\qquad

If X is a normal Fermi vector field and A and B are tangential Fermi vector
fields, then:

\begin{lemma}
\qquad\qquad$\nabla_{\text{X}}$A = $\nabla_{\text{A}}$X\qquad and,$\qquad
\nabla_{\text{B}}$A = $\nabla_{\text{A}}$B\qquad\qquad\qquad\qquad\qquad
\qquad\qquad\qquad

where $\nabla$ is the Levi-Civita connection.
\end{lemma}

\qquad$\qquad\qquad\qquad\nabla_{\text{X}}$A -$\nabla_{\text{A}}$X =$\left[
\text{X , A}\right]  \qquad$and,$\qquad\nabla_{\text{B}}$A -$\nabla_{\text{A}%
}$B =$\left[  \text{A , B}\right]  $

By \textbf{Lemma (3.1)} of \textbf{A. Gray} and \textbf{L. Vanhecke }$\left[
27\right]  $

\qquad\qquad\qquad\qquad$\left[  \text{X , A}\right]  =0=\left[  \text{A ,
B}\right]  ,$

and so the results follow.

$\qquad\qquad\qquad\qquad\qquad\qquad\qquad\qquad\qquad\qquad\qquad
\qquad\ \ \ \ \ \ \qquad\qquad\blacksquare\qquad\qquad\qquad\qquad\qquad
\qquad\qquad\ \ $\qquad\qquad\qquad\qquad\qquad\qquad\qquad\qquad\qquad
\qquad\qquad

\vspace{1pt}The p$^{th}$ covariant derivative $\nabla^{\text{p}}$ is defined
inductively as follows:\qquad\qquad\qquad\qquad\qquad\qquad\qquad

$\ \ \ \ \ \ \ \ \ \ \ \ \ \ \ \ \ \ \ \ \ \ \ \ \nabla_{\text{X}%
_{\text{i}_{1}}}^{\text{p}},...,_{\text{X}_{\text{i}_{\text{p}}}}%
=\nabla_{\text{X}_{\text{i}_{1}}}$($\nabla_{\text{X}_{\text{i}_{2}}%
}^{\text{p-1}},...,_{\text{X}_{\text{i}_{\text{p}}}}$)

Following standard notation$,$we set:

\qquad\qquad$\ $\ $\ \ \ \ \ \ \ \ \ \ \ \ \ \nabla_{\text{i}_{1}}^{\text{p}%
},...,_{\text{i}_{\text{p}}}\ \ =\ \nabla_{\text{X}_{\text{i}_{1}}}^{\text{p}%
},...,_{\text{X}_{\text{i}_{\text{p}}}}$

\vspace{1pt}We have the following important lemma on covariant derivatives:

\qquad\qquad\qquad\qquad\qquad\qquad\qquad\qquad\qquad\qquad\qquad\qquad
\qquad\qquad\qquad$\blacksquare$

\begin{lemma}
(A. Gray and L. Vanhecke $\left[  27\right]  ,$ Lemma $\left(  3.3\right)  $
and Lemma $\left(  3.7\right)  $)

Let X,Y, Z be normal Fermi vector fields, A a tangential Fermi vector field
and R the Riemannian curvature tensor of M relative to the connection
$\nabla.$ Then,\qquad\ \ 
\end{lemma}

(i) \qquad\ \ ($\nabla_{\text{XY}}^{2}$A)(y$_{0}$) \ \ $=$ $-$ (R$_{\text{XA}%
}$Y)(y$_{0}$)

(ii)\qquad\ \ ($\nabla_{\text{XXX}}^{3}$A)(y$_{0}$) $=$ ($-\nabla_{\text{X}}%
$(R)$_{\text{XA}}$X $-$ R$_{\text{XT}_{\text{A}}\text{X}}$X $-$ R$_{\text{X}%
\perp_{\text{A}}\text{X}}$X)(y$_{0}$)

(iii)\qquad\ \ ($\nabla_{\text{XXXX}}^{4}$A)(y$_{0}$) $=$ ($-\nabla
_{\text{XX}}^{2}($R)$_{\text{XA}}$X $+$ R$_{\text{XR}_{\text{XA}}\text{X}}$X
$-$ 2$\nabla_{\text{X}}$(R)$_{\text{XT}_{\text{A}}\text{X}}$X

\qquad\qquad\qquad\qquad\qquad$\ \ -$ 2$\nabla_{\text{X}}$(R)$_{\text{X}%
\perp_{\text{A}}\text{X}}$X)(y$_{0}$)

(iv)\qquad\ \ ($\nabla_{\text{X}}$Y)(y$_{0}$) $\ =0$

(v)\qquad\ \ \ ($\nabla_{\text{XY}}^{2}$Z)(y$_{0}$) $\ =-\frac{1}{3}%
($R$_{\text{XY}}$Z + R$_{\text{XZ}}$Y)(y$_{0}$)

(vi)\qquad\ \ ($\nabla_{\text{XXX}}^{3}$Y)(y$_{0}$) $=-\frac{1}{2}%
(\nabla_{\text{X}}($R)$_{\text{XY}}$X)(y$_{0}$)

(vii)\qquad($\nabla_{\text{XXXX}}^{4}$Y)(y$_{0}$) $=(-\frac{3}{5}%
\nabla_{\text{XX}}^{2}($R)$_{\text{XY}}$X + $\frac{1}{5}$R$_{\text{XR}%
_{\text{XY}}\text{X}}$X)(y$_{0}$)

\qquad\qquad\qquad\qquad\qquad\qquad\qquad\qquad\qquad\qquad\qquad\qquad
\qquad\qquad$\blacksquare$

\subsection{Notation}

We introduce the following notation with the following convention for indices:
Tangential Fermi coordinates and tangential Fermi vector fields will \ be
indexed by a,b,c,d (which run from 1 to q). Normal Fermi coordinates \ and
normal Fermi vector fields will be indexed by $i,j,k,l,r,s,t$ (which run from
$q+1$ to $n$).

We will follow the following notation:

(i)$\qquad\perp_{\text{X}_{\text{a}}}$x$_{i}=$ $\perp_{\text{a}i}$

(ii)$\qquad$%
$<$%
$\bot_{\text{X}_{\text{a}}}$X$_{i}$ , X$>$ = $\bot_{\text{a}i\text{X}}$ \ \ 

(iii)\qquad\
$<$%
$\bot_{\text{X}_{\text{a}}}$X$_{i}$ , X$_{j}>$ = $\bot_{\text{a}ij}$

(iv)\qquad T$_{\text{X}_{\text{a}}}$X$_{i}=$ T$_{\text{a}i}$ \ \ \ \ 

(v)\qquad T$_{\text{X}_{\text{a}}}$X$_{\text{b}}=$ T$_{\text{ab}}\qquad$

(vi)$\qquad$%
$<$%
T$_{\text{X}_{\text{a}}}$X$_{\text{b}}$ , X$_{i}>$ = T$_{\text{ab}i}$

(vii)\qquad%
$<$%
T$_{\text{X}_{\text{a}}}$X$_{\text{b}}$ , X$>$ = T$_{\text{abX}}$

(viii\qquad%
$<$%
T$_{\text{X}_{\text{a}}}$X $_{\text{b}}$,T$_{\text{X}_{\text{c}}}$
X$_{\text{d}}>$ $=$
$<$%
T$_{\text{ab}},$T$_{\text{cd}}>$ $=\overset{n}{\underset{i=q+1}{\sum}}%
$T$_{\text{ab}i}$T$_{\text{cd}i}$

(viii)$^{\ast}\qquad$%
$<$%
T$_{\text{X}_{\text{a}}}$X $_{\text{b}}$,$\bot_{\text{X}_{\text{c}}}$X$_{i}>$
$=$ $\underset{j=q+1}{\overset{n}{\sum}}$T$_{\text{ab}j}\perp_{\text{c}ij}$

\qquad We note that by p. 37 of \textbf{Tubes}, \textbf{Gray }$\left[
26\right]  $\textbf{ }T$_{\text{X}_{\text{a}}}$X $_{\text{b}}=$ T$_{\text{ab}%
}=$ T$_{\text{A}}$B is normal and $\bot_{\text{X}_{\text{c}}}$X$_{i}%
=\perp_{\text{c}i}$ is normal.

\qquad Consequently the RHS of (viii)$^{\ast}$ makes sense.\qquad\qquad\qquad

(ix) (a)$\qquad$%
$<$%
R$_{\text{X}_{i}\text{X}_{j}}$X$_{k}$ , X$_{l}>$ = R$_{ijkl}$

$\qquad\qquad$( components of the Riemannian curvature tensor of M)

(ix) (b)$\qquad$%
$<$%
R$_{\text{XX}_{\text{a}}}$X , X$_{\text{b}}>$ = R$_{\text{XX}_{\text{a}%
}\text{XX}_{\text{b}}}$

(x)$\qquad\ \underset{i=1}{\overset{n}{\sum}}$ R$_{\text{X}_{i}\text{XX}%
_{i}\text{X}}=\varrho_{\text{XX}}^{M}$ = components of the Ricci curvature tensor.

(xi)\qquad H = $\underset{\text{a}=1}{\overset{q}{\sum}}$T$_{\text{E}%
_{\text{a}}}$E$_{\text{a}}=$ $\underset{\text{a}=1}{\overset{q}{\sum}}%
$T$_{\text{aa}}=$ the mean curvature vector field,

where, E$_{\text{1}}$,......,E$_{\text{q}}$ is a local frame field on P.\qquad

(xii) \qquad%
$<$%
H , X$_{\text{i}}$%
$>$
=
$<$%
H,i%
$>$
\ 

(xiii) (a) \
$<$%
$\nabla_{\text{X}}$(R)$_{\text{XX}_{\text{a}}}$X , X$_{\text{b}}>$ =
$\nabla_{\text{X}}$R$_{\text{XaX}_{\text{b}}}$

(xiii) (b)\qquad%
$<$%
$\nabla_{\text{XX}}^{2}$(R)$_{\text{XX}_{\text{a}}}$X , X$_{\text{b}}$%
$>$
= $\nabla_{\text{XX}}^{2}$R$_{\text{X}_{\text{a}}\text{X}_{\text{b}}}$

\qquad Note that for a, b = 1,...,q we have:

(xiv)\qquad T$_{\text{ab}X}=$ $-$T$_{\text{a}X\text{b}}$ \ by (3.15) of A.
\textbf{Gray and L. Vanhecke} $\left[  27\right]  $\ 

(xv)\qquad T$_{\text{ab}X}=$ T$_{\text{ba}X}$ \ \ \ by (3.14) of A.
\textbf{Gray and L. Vanhecke} $\left[  27\right]  $

(xvi)\qquad\
$<$%
$\bot_{\text{X}_{\text{a}}}$X , Y$>$ $=$ $-$
$<$%
$\bot_{\text{X}_{\text{a}}}$Y , X$>$ by lemma (3.4) of \textbf{Gray and
Vanhecke }$\left[  27\right]  \qquad$\ 

\qquad We will adopt the following familiar notation:

$\ $\ R$_{ijkl}^{\text{M}}$ = components of the curvature tensor of the
Riemannian manifold M.

$\ $R$_{\text{abcd}}^{\text{P}}$ = components of the curvature tensor of the
submanifold P.

(i)$\qquad\varrho_{ij}^{\text{M}}=$ $\underset{k=1}{\overset{n}{\sum}}%
$R$_{ikjk}^{\text{M}}$ = components of the Ricci curvature tensor of M.

(ii)$\qquad\varrho_{\text{ab}}^{\text{P}}=\underset{\text{c}%
=1}{\overset{q}{\sum}}$R$_{\text{acbc}}^{\text{P}}$ = components of the Ricci
curvature tensor of P

(iii)$\qquad\tau^{\text{M}}=\underset{i,j=1}{\overset{n}{\sum}}$%
R$_{ijij}^{\text{M}}=\underset{i=1}{\overset{n}{\sum}}\varrho_{ii}^{\text{M}}$
denotes the scalar curvature of M.

(iv)$\qquad\tau^{\text{P}}=\underset{\text{a},\text{b}=1}{\overset{q}{\sum}}%
$R$_{\text{abab}}^{\text{P}}=\underset{\text{a}=1}{\overset{q}{\sum}}%
\varrho_{\text{aa}}^{\text{P}}$ denotes the scalar curvature of P.

(v)$\qquad\left\Vert \varrho^{\text{M}}\right\Vert ^{2}%
=\overset{n}{\underset{i,j=1}{\sum}}(\varrho_{ij}^{\text{M}})^{2}$

(vi) \ $\left\Vert \varrho^{\text{P}}\right\Vert ^{2}%
=\overset{q}{\underset{\text{a,b}=1}{\sum}}(\varrho_{\text{ab}}^{\text{P}%
})^{2}$\qquad

\vspace{1pt}(vii)$\qquad\left\Vert \text{R}^{\text{M}}\right\Vert
^{2}\underset{i,j,k,l=1}{=\overset{n}{\sum}}($R$_{i,j,k,l}^{\text{M}}%
)^{2}\qquad$

(viii) $\ \left\Vert \text{R}^{\text{P}}\right\Vert ^{2}%
\underset{\text{a,b,c,d}=1}{=\overset{q}{\sum}}($R$_{\text{abcd}}^{\text{P}%
})^{2}$

We note however that when there is an \textbf{absence of a superscript} on a
geometric invariant (the Riemannian curvature tensor R, the Ricci curvature
$\varrho,$ the scalar curvature $\tau),$ then such an invariant is taken
relative to the Riemannian manifold M.

\qquad\qquad\qquad\qquad\qquad\qquad\qquad\qquad\qquad\qquad\qquad\qquad
\qquad\qquad\qquad\qquad\qquad$\blacksquare$

\section{Preliminary Geometric Expansion Formulae}

We need the expansion of the infinitesimal change of volume function
$\theta_{P}$ as well as those of the components of the metric tensor
g$_{\alpha\beta}$ defined by Fermi \ coordinates. The general computations of
these expansions are givne in \textbf{Theorem (9.21)} \textbf{Gray} $\left[
25\right]  $ which is the same as \textbf{Theorem (4.2)} of \textbf{Gray and
Vanhecke} $\left[  27\right]  $ which states that: if (x$_{1},...,x_{q}%
,x_{q+1},...,$x$_{n}$) are Fermi coordinates and W is a covariant tensor
field, then W($\frac{\partial}{\partial\text{x}_{\alpha_{1}}},...\frac
{\partial}{\partial\text{x}_{\alpha_{r}}}$) can be expanded in the coordinates
x$_{1},...,$x$_{n}.$

In practice we have followed the use of this Theorem as applied to
\textbf{Theorem }$\left(  4.3\right)  $ of \textbf{Gray and Vanhecke} $\left[
26\right]  $ and Theorem $\left(  9.22\right)  $ of \textbf{Gray} $\left[
25\right]  .$

We will need the expansion only in the (normal) coordinates x$_{q+1}%
,...,$x$_{n}.$ First we will give the expansions of the components of the
metric tensor g$_{\alpha\beta}.$ We distinguish \ three cases:\ g$_{\text{ab}%
},$g$_{\text{a}j}$ and g$_{ij}$ for a,b $=1,...,$q and $i,j=q+1,...,n$\ 

The special case of M = R$^{\text{n}}$ was given in \textbf{Gray, Karp and
Pinsky} $\left[  29\right]  $ (with a slight error on the expansion of
g$_{\text{ab}}$).

We now give the expansions of g$_{ij}(x)$ for $x\in$M$_{0}:$

\qquad\qquad\qquad\qquad\qquad\qquad\qquad\qquad\qquad\qquad\qquad\qquad
\qquad\qquad\qquad\qquad\qquad$\blacksquare$\qquad

\begin{proposition}
For a,b = 1,...,q and $x\in$M$_{0},$ we have:

\vspace{1pt}\qquad For a,b = 1,...,q; $r,s,t=q+1,...,n$ and $x\in$M$_{0},$ we have:
\end{proposition}

g$_{\text{ab}}(x)$ $=$ $\delta_{\text{ab}}$ $-$
$2\overset{n}{\underset{r=q+1}{\sum}}$ T$_{\text{ab}r}\,\,\,($y$_{0})$x$_{r}+$
$\underset{r,s\text{=q+1}}{\overset{\text{n}}{\sum}}\{-$R$_{\text{a}%
r\text{b}s}$ + $\underset{\text{c=1}}{\overset{\text{q}}{\sum}}$%
T$_{\text{ac}r}$T$_{\text{bc}s}+\underset{t\text{=q+1}}{\overset{\text{n}%
}{\sum}}$ $\perp_{\text{a}rt}\perp_{\text{b}st}\}$(y$_{0}$)x$_{r}$x$_{s}$

$\qquad\qquad\qquad-\overset{n}{\underset{r,s,t=q+1}{\frac{1}{6}\sum}%
\{}2\nabla_{r}($R)$_{s\text{a}t\text{b}}+$ R$_{\text{b}ki\text{T}_{\text{a}s}%
}+$ R$_{\text{bki}\perp_{\text{a}s}}+$ $3$ (R$_{\text{a}rs\text{T}_{\text{b}%
t}}$ $+\;$R$_{\text{a}rs\perp_{\text{b}t}})$\ 

$\qquad\qquad\qquad+$ $3($R$_{\text{b}rs\text{T}_{\text{a}t}}$ $+$
R$_{\text{b}rs\perp_{\text{a}t}})+$R$_{\text{a}tr\text{T}_{\text{b}s}}+$
R$_{\text{a}tr\perp_{\text{b}s}}\}(y_{0})$x$_{r}$x$_{s}$x$_{t}$

$\qquad+$ higher order terms

\begin{proof}
Let X \ be a normal Fermi vector field . Then by \textbf{Theorem (4.2)} of
\textbf{Gray and Vanhecke} $[27]$:
\end{proof}

g$_{\text{ab}}($x$_{0}$) = g$_{\text{ab}}$(y$_{0}$) + $\frac{1}{1!}%
\underset{i\text{=q+1}}{\overset{\text{n}}{\sum}}($Xg$_{\text{ab}})($y$_{0}%
)$x$_{i}+\frac{1}{2!}\underset{i,j=q+1}{\overset{\text{n}}{\sum\text{ }}}%
($X$^{\text{2}}$g$_{\text{ab}})($y$_{0})$x$_{i}$x$_{j}+\frac{1}{3!}%
\underset{i,j,l\text{=q+1}}{\overset{\text{n}}{\sum}}($X$^{3}$g$_{\text{ab}%
})($y$_{0})$x$_{i}$x$_{j}$x$_{k}$

$\qquad\qquad+$ higher order terms.\qquad

Since the vectors $\frac{\partial}{\partial\text{x}_{1}},..,\frac{\partial
}{\partial\text{x}_{q}},\frac{\partial}{\partial\text{x}_{q+1}},...,\frac
{\partial}{\partial\text{x}_{n}}$ are orthonormal at y$_{0}$, it is clear that:

g$_{\text{ab}}($y$_{0})=$ $<$X$_{\text{a}}$,X$_{\text{b}}>($y$_{0})$ $=$
$<\frac{\partial}{\partial\text{x}_{\text{a}}},\frac{\partial}{\partial
\text{x}_{\text{b}}}>($y$_{0})$ $=\delta_{\text{ab}}$ for a,b =
1,...,q.$\qquad\qquad\qquad\qquad\qquad\qquad\qquad\qquad\qquad\qquad
\qquad\qquad\qquad$

Next let X be a normal Fermi vector field.

\qquad\ \ Xg$_{\text{ab}}=$ X%
$<$%
X$_{\text{a}}$ ,X$_{\text{b}}>$ =
$<$%
$\nabla_{\text{X}}$X$_{\text{a}}$ , X$_{\text{b}}>+$
$<$%
X$_{\text{a}}$ ,$\nabla_{\text{X}}$ X$_{\text{b}}>$

\qquad\qquad\ \ \ \ \ =
$<$%
$\nabla_{\text{X}_{\text{a}}}$X , X$_{\text{b}}>+$
$<$%
X$_{\text{a}}$ ,$\nabla_{\text{X}_{\text{b}}}$ X$> $\qquad

\thinspace\thinspace\thinspace\thinspace\thinspace\thinspace\thinspace
\thinspace\thinspace\thinspace\thinspace\thinspace\thinspace\thinspace
\thinspace\thinspace\thinspace\thinspace The last equality is due to
\textbf{Lemma 6.2}. Then \textbf{Lemma 6.1} gives:

\qquad(\ Xg$_{\text{ab}})(y_{0})=$ $\left\{  \text{%
$<$%
T}_{\text{X}_{\text{a}}}\text{X +}\perp_{\text{X}_{\text{a}}}\text{X,X}%
_{\text{b}}\text{%
$>$
+ \thinspace%
$<$%
T}_{\text{X}_{\text{b}}}\text{X +}\perp_{\text{X}_{\text{b}}}\text{X,X}%
_{\text{a}}\text{%
$>$%
}\right\}  (y_{0})$

Since $\perp_{\text{X}_{\text{a}}}$X is normal and X$_{\text{b}}$ is
tangential, we have:

(\ Xg$_{\text{ab}})(y_{0})$ = $\left\{  \text{%
$<$%
T}_{\text{X}_{\text{a}}}\text{X,X}_{\text{b}}\text{%
$>$
+\thinspace\
$<$%
T}_{\text{X}_{\text{b}}}\text{X,X}_{\text{a}}\text{%
$>$%
}\right\}  (y_{0})$\thinspace

= $-$\thinspace\thinspace$\left\{  \text{%
$<$%
T}_{\text{X}_{\text{a}}}\text{X}_{\text{b}}\text{,X%
$>$
+\thinspace\
$<$%
T}_{\text{X}_{\text{b}}}\text{X}_{\text{a}}\text{,X%
$>$%
}\right\}  (y_{0})$ = $-\left\{  \text{T}_{\text{abX}}\text{\thinspace
\thinspace\thinspace+\thinspace T}_{\text{baX}}\right\}  (y_{0})=$
$-2$T$_{\text{abX}}$(y$_{0}$)\thinspace\thinspace\thinspace\thinspace
\thinspace\thinspace\thinspace\thinspace\thinspace\thinspace\thinspace
\thinspace\thinspace\thinspace\thinspace\thinspace\thinspace\thinspace
\thinspace\thinspace\thinspace\thinspace\thinspace\thinspace\thinspace
\thinspace\thinspace\thinspace\thinspace\thinspace\thinspace\thinspace
\thinspace\thinspace\thinspace\thinspace\thinspace\thinspace\thinspace
\thinspace\thinspace\thinspace\thinspace\thinspace\thinspace\thinspace
\thinspace\thinspace\thinspace\thinspace\thinspace\thinspace\thinspace
\thinspace\thinspace\thinspace\thinspace\thinspace\thinspace\thinspace
\thinspace\thinspace\thinspace\thinspace\thinspace\thinspace\thinspace

Next,

\qquad\thinspace\thinspace\thinspace\thinspace X$^{2}$g$_{\text{ab}}%
$\thinspace\thinspace\thinspace\ = X%
$<$%
$\nabla_{\text{X}}$X$_{\text{a}}$ , X$_{\text{b}}>+$ X%
$<$%
X$_{\text{a}}$ ,$\nabla_{\text{X}}$ X$_{\text{b}}>$

\qquad=
$<$%
$\nabla_{\text{XX}}^{2}$X$_{\text{a}}$ , X$_{\text{b}}>+$
$<$%
$\nabla_{\text{X}}$X$_{\text{a}}$ , $\nabla_{\text{X}}$X$_{\text{b}}>$ +
$<$%
$\nabla_{\text{X}}$X$_{\text{a}}$ ,$\nabla_{\text{X}}$ X$_{\text{b}}>+$
$<$%
X$_{\text{a}}$ ,$\nabla_{\text{XX}}^{2}$ X$_{\text{b}}>$

\qquad\ =
$<$%
$\nabla_{\text{XX}}^{2}$X$_{\text{a}}$ , X$_{\text{b}}>+$ 2%
$<$%
$\nabla_{\text{X}}$X$_{\text{a}}$ , $\nabla_{\text{X}}$X$_{\text{b}}>$ $+$
$<$%
X$_{\text{a}}$ ,$\nabla_{\text{XX}}^{2}$ X$_{\text{b}}>\qquad\qquad\qquad$

Then using the lemmas above, we have:

\ \ \ (X$^{2}$g$_{\text{ab}}$)(y$_{0}$) = (%
$<$%
$-$ R$_{\text{XX}_{\text{a}}}$X ,X$_{\text{b}}>$)(y$_{0}$)

$+$ $2(<$T$_{\text{X}_{\text{a}}}$X +$\perp_{\text{X}_{\text{a}}}$X
,T$_{\text{X}_{\text{b}}}$X +$\perp_{\text{X}_{\text{b}}}$X $>)$(y$_{0}$) $+$
( $<-$R$_{\text{XX}_{\text{b}}}$X ,X$_{\text{a}}>)$(y$_{0}$)

Since T$_{\text{X}_{\text{a}}}$X and T$_{\text{X}_{\text{b}}}$X are
tangential, and $\perp_{\text{X}_{\text{a}}}$X $\ $and $\perp_{\text{X}%
_{\text{b}}}$X are normal, we have at y$_{0}:$

\qquad\qquad\ $<$ T$_{\text{X}_{\text{a}}}$X
$>$%
, $\perp_{\text{X}_{\text{b}}}$X%
$>$
= 0 = $<$ $\perp_{\text{X}_{\text{a}}}$X ,T$_{\text{X}_{\text{b}}}$X
$>$%

and so,

(X$^{2}$g$_{\text{ab}}$)(y$_{0}$) = $\left\{  -\text{R}_{\text{XX}_{\text{a}%
}\text{XX}_{\text{b}}}\text{+2%
$<$%
T}_{\text{X}_{\text{a}}}\text{X,T}_{\text{X}_{\text{b}}}\text{X%
$>$%
+2%
$<$%
}\perp_{\text{X}_{\text{a}}}\text{,}\perp_{\text{X}_{\text{b}}}\text{%
$>$%
}-\text{R}_{\text{XX}_{\text{b}}\text{XX}_{\text{a}}}\right\}  $(y$_{0}$)

\ \qquad\ \ \ \ = $\left\{  -\text{2R}_{\text{XX}_{\text{a}}\text{XX}%
_{\text{b}}}\text{+2%
$<$%
T}_{\text{X}_{\text{a}}}\text{X,T}_{\text{X}_{\text{b}}}\text{X%
$>$%
+2%
$<$%
}\perp_{\text{X}_{\text{a}}}\text{X,}\perp_{\text{X}_{\text{b}}}\text{X%
$>$%
}\right\}  $(y$_{0}$)

Hence,

\qquad\qquad(X$^{2}$g$_{\text{ab}}$)(y$_{0}$)\ = $\left\{  -\text{2R}%
_{\text{XX}_{\text{a}}\text{XX}_{\text{b}}}\text{+2%
$<$%
T}_{\text{aX}}\text{,T}_{\text{bX}}\text{%
$>$%
+2%
$<$%
}\perp_{\text{aX}}\text{,}\perp_{\text{bX}}\text{%
$>$%
}\right\}  $(y$_{0}$)

\qquad\qquad\qquad\qquad\ \ \ \ = $2\underset{\text{i,j=q+1}%
}{\overset{\text{n}}{\sum}}\left\{  (-\text{R}_{\text{iajb}}\text{+%
$<$%
T}_{\text{ai}}\text{,T}_{\text{bj}}\text{%
$>$%
+%
$<$%
}\perp_{\text{ai}}\text{,}\perp_{\text{bj}}\text{%
$>$%
})\right\}  $(y$_{0}$)

\qquad\ = $2\underset{\text{i,j=q+1}}{\overset{\text{n}}{\sum}}\{-$%
R$_{\text{iajb}}$ + $\underset{\text{c,d=1}}{\overset{\text{q}}{\sum}}%
$T$_{\text{aic}}$T$_{\text{bjd}}$
$<$%
$\frac{\partial}{\partial\text{x}_{\text{c}}},\frac{\partial}{\partial
\text{x}_{\text{d}}}>+\underset{\text{k,l=q+1}}{\overset{\text{n}}{\sum}}$
$\perp_{\text{aik}}\perp_{\text{bjl}}<\frac{\partial}{\partial\text{x}%
_{\text{k}}},\frac{\partial}{\partial\text{x}_{\text{l}}}>\}$(y$_{0}$)

Since the system of Cartesian Fermi coordinates is orthonormal at y$_{0},$ we have:%

$<$%
$\frac{\partial}{\partial\text{x}_{\text{c}}},\frac{\partial}{\partial
\text{x}_{\text{d}}}>(y_{0})=\delta_{\text{cd}}$ and $<\frac{\partial
}{\partial\text{x}_{\text{k}}},\frac{\partial}{\partial\text{x}_{\text{l}}}%
>$(y$_{0}$) = $\delta_{\text{kl}}$

and so since T$_{\text{a}i\text{c}}=-$T$_{\text{ac}i}$ by $\left(
3.15\right)  $ of \textbf{Gray and Vanhecke} $\left[  27\right]  ,$ we have:

\qquad(X$^{2}$g$_{\text{ab}}$)(y$_{0}$)\ =
$2\underset{r,s=q+1}{\overset{n}{\sum}}\{-$R$_{r\text{a}s\text{b}}$ +
$\underset{\text{c=1}}{\overset{\text{q}}{\sum}}$T$_{\text{a}r\text{c}}%
$T$_{\text{b}s\text{c}}$ $+\underset{t=q+1}{\overset{n}{\sum}}$ $\perp
_{\text{a}rt}\perp_{\text{b}st}\}$(y$_{0}$)

\qquad\qquad\qquad= $2\underset{r,s=q+1}{\overset{n}{\sum}}\{-$R$_{\text{a}%
r\text{b}s}$ + $\underset{\text{c=1}}{\overset{\text{q}}{\sum}}(-$%
T$_{\text{ac}r})(-$T$_{\text{bc}s}$ ) $+\underset{\text{k=q+1}%
}{\overset{\text{n}}{\sum}}$ $\perp_{\text{a}r\text{k}}\perp_{\text{b}%
s\text{k}}\}$(y$_{0}$)

\qquad\qquad\qquad= $2\underset{r,s=q+1}{\overset{n}{\sum}}\{-$R$_{\text{a}%
r\text{b}s}$ + $\underset{\text{c=1}}{\overset{\text{q}}{\sum}}$%
T$_{\text{ac}r}$T$_{\text{bc}s}$ $+\underset{t=q+1}{\overset{n}{\sum}}$
$\perp_{\text{a}rt}\perp_{\text{b}st}\}$(y$_{0}$)

We note that since we are dealing with the Levi-Cevita (torsion-free)
connection on M, the torsion terms above must disappear.

\qquad\qquad\qquad\qquad\qquad\qquad\qquad\qquad\qquad\qquad\qquad\qquad
\qquad\qquad\qquad$\qquad\blacksquare$

\begin{proposition}
\qquad\qquad\qquad\ \ \qquad\thinspace\thinspace\thinspace\thinspace
\thinspace\thinspace\thinspace\thinspace\thinspace\thinspace\thinspace
\thinspace\thinspace\thinspace\thinspace\thinspace\thinspace\thinspace
\thinspace\thinspace\thinspace\thinspace\thinspace\thinspace\thinspace
\thinspace\thinspace\thinspace\thinspace\thinspace\thinspace\thinspace
\thinspace\thinspace\thinspace\thinspace\thinspace\thinspace\thinspace
\thinspace\thinspace\thinspace\thinspace\thinspace\thinspace\thinspace
\thinspace\thinspace\thinspace\thinspace\thinspace\thinspace\thinspace
\thinspace\thinspace\thinspace\thinspace\thinspace\thinspace\thinspace
\thinspace\thinspace\thinspace\thinspace\thinspace\thinspace\thinspace
\thinspace\thinspace\thinspace\thinspace\thinspace\thinspace\thinspace
\thinspace\thinspace\thinspace\thinspace\thinspace\thinspace\thinspace
\thinspace\qquad\thinspace\thinspace\thinspace\thinspace\thinspace
\thinspace\thinspace\thinspace\thinspace\thinspace\thinspace\thinspace
\thinspace\thinspace\thinspace\thinspace\thinspace\thinspace\thinspace
\thinspace\thinspace\thinspace\thinspace\thinspace\thinspace\thinspace
\thinspace\thinspace\thinspace\thinspace\thinspace\thinspace\thinspace
\thinspace\thinspace\thinspace\thinspace\thinspace\thinspace\thinspace
\thinspace\thinspace\thinspace\thinspace\thinspace\thinspace\thinspace
\thinspace\thinspace\thinspace\thinspace\thinspace\thinspace\thinspace
\thinspace\thinspace\thinspace\thinspace\thinspace\thinspace\thinspace
\thinspace\thinspace\thinspace\thinspace\thinspace\thinspace\thinspace
\thinspace\thinspace\thinspace\thinspace\thinspace\thinspace\thinspace
\thinspace\thinspace\thinspace\thinspace\thinspace\thinspace\thinspace
\thinspace\thinspace\thinspace\thinspace\thinspace\thinspace\thinspace
\thinspace\thinspace\thinspace\thinspace\thinspace\thinspace\thinspace
$\qquad\qquad\qquad\qquad\qquad\qquad\qquad\ \ \ \ \ \ \ \qquad\qquad$%
\qquad\qquad
\end{proposition}

\ For a $=1,...,q$ and $i=q+1,...,n,$ we have for $x\in$M$_{0}:$

\ g$_{\text{a}i}(x)=\delta_{\text{a}i}$\ $-\overset{n}{\underset{r=q+1}{\sum
(}}\perp_{\text{a}ir})(y_{0})$x$_{r}-\frac{4}{3}%
\underset{r,s=q+1}{\overset{\text{n}}{\sum\text{ }}}$R$_{r\text{a}si}(y_{0}%
)$x$_{r}$x$_{s}$

$\qquad\qquad-\frac{1}{6}\underset{r,s,t=q+1}{\overset{\text{n}}{\sum}}%
\{\frac{3}{2}\nabla_{r}($R)$_{s\text{a}ti}+2$R$_{ris\text{T}_{\text{a}t}}+$
$2$R$_{ris\perp_{\text{a}t}}\}(y_{0})$x$_{r}$x$_{s}$x$_{t}+$ higher order terms.

where $\delta_{\text{a}i}=0$

\begin{proof}
\qquad We use the same techniques as in the previous proposition.
\end{proof}

\qquad\qquad\qquad\qquad\qquad\qquad\qquad\qquad\qquad\qquad\qquad\qquad
\qquad\qquad\qquad\qquad\qquad$\blacksquare$

\qquad The proposition below has been proved in the simpler case of
\textbf{normal coordinates} in several papers and books: (\textbf{Gray}%
$\left[  25\right]  ,$\textbf{ Ii}$\left[  31\right]  $\textbf{,}
\textbf{McKean-Singer}$\left[  37\right]  ,$\textbf{ Sakai}$\left[  49\right]
).$ We give an expansion below in the more general case of \textbf{Fermi
coordinates}.\qquad\qquad\qquad\qquad\qquad\qquad\qquad\qquad\qquad
\qquad\qquad\qquad\qquad\qquad\qquad

\begin{proposition}
For $k,l=q+1,...,n$ and $x\in$M$_{0},$ we have:
\end{proposition}

$\qquad g_{kl}(x)=$ $\delta_{kl}-$ $\frac{1}{3}%
\underset{r,s=1}{\overset{n}{\sum}}($ R$_{rksl})(y_{0})$x$_{r}$x$_{s}-\frac
{1}{6}\underset{r,s,t=1}{\overset{n}{\sum}}\nabla_{_{r}}$R$_{sktl}($y$_{0}%
)$x$_{r}$x$_{s}$x$_{t}$

$\ \ +\frac{1}{360}\overset{n}{\underset{r,s,t,u=1}{\sum}}(-18\nabla_{rs}^{2}%
$R$_{tkul}+16\underset{w=\text{1}}{\overset{n}{\sum}}$R$_{rksw}$R$_{tluw}%
)($y$_{0})$x$_{r}$x$_{s}$x$_{t}$x$_{u}$

$\ +$ $\frac{1}{90}\overset{n}{\underset{r,s,t,u,v=1}{\sum\{-\nabla_{rst}^{3}%
}\text{R}_{ukvl}}+$ $2\overset{n}{\underset{w=1}{\sum}}(\nabla_{r}$R$_{sktw}%
$R$_{ulvw}+\nabla_{r}$R$_{sltw}$R$_{ukvw})\}(y_{0})$x$_{r}$x$_{s}$x$_{t}%
$x$_{u}$x$_{v}\qquad\qquad\qquad$

$\ +$ \ higher order terms.

\begin{proof}
The proposition has already been proved in the papers cited above.$\qquad
\qquad\qquad\qquad\qquad\qquad\qquad\qquad\qquad\qquad\qquad\qquad\qquad
\qquad\qquad\qquad\qquad\qquad\qquad\qquad\qquad\qquad\qquad\qquad\qquad
\qquad\qquad\qquad\qquad\qquad$
\end{proof}

\thinspace\thinspace\thinspace\thinspace X$^{2}$g$_{kl}$\thinspace
\thinspace\thinspace$=$ $[$X$<\nabla_{\text{X}}$X$_{k}$ , X$_{l}>+$ X$<$%
X$_{k}$ ,$\nabla_{\text{X}}$ X$_{l}>]$

\qquad$\ \ \ \ \ =$ $[<\nabla_{\text{XX}}^{2}$X$_{k}$ , X$_{l}>+<\nabla
_{\text{X}}$X$_{k}$ , $\nabla_{\text{X}}$X$_{l}>$ + $<\nabla_{\text{X}}$%
X$_{k}$ ,$\nabla_{\text{X}}$ X$_{l}>$

$\qquad\qquad+<$X$_{k}$ ,$\nabla_{\text{XX}}^{2}$ X$_{l}>]$

\qquad$\ \ \ \ \ =$ $[<\nabla_{\text{XX}}^{2}$X$_{k}$ , X$_{l}>+$
2$<\nabla_{\text{X}}$X$_{k}$ , $\nabla_{\text{X}}$X$_{l}>$ $+$ $<$X$_{k}$
,$\nabla_{\text{XX}}^{2}$ X$_{l}>]$

\qquad\qquad$=$ $[<\nabla_{\text{XX}}^{2}$X$_{k}$ , X$_{l}>+$ $2<\nabla
_{\text{X}}$X$_{k}$ , $\nabla_{\text{X}}$X$_{l}>$ $+<\nabla_{\text{XX}}%
^{2}X_{l},X_{k}>]$

By (v) of \textbf{Lemma 6.3} above, we have for $k,l=q+1,...,n$:

\thinspace X$^{2}$g$_{kl}$\thinspace\thinspace\thinspace$=$ $[<-\frac{1}{3}%
($R$_{\text{XX}}$X$_{k}+$ R$_{\text{XX}_{k}}$X$)$ , X$_{l}>$ $+$ $<-\frac
{1}{3}($R$_{\text{XX}}$X$_{l}+$ R$_{\text{XX}_{l}}$X$)$ , X$_{k}>]$

$\qquad\ \ \ \ \ +$ $[-\frac{1}{3}($R$_{\text{XXX}_{k}\text{X}_{l}}+$
R$_{\text{XX}_{k}\text{XX}_{l}})$ $-\frac{1}{3}($R$_{\text{XXX}_{l}%
\text{X}_{k}}+$ R$_{\text{XX}_{l}\text{XX}_{k}})]$

Since R$_{\text{XXX}_{k}\text{X}_{l}}=0=$ R$_{\text{XXX}_{l}\text{X}_{k}}$ and
R$_{\text{XX}_{l}\text{XX}_{k}}=$ R$_{\text{XX}_{k}\text{XX}_{l}},$ we have:

\qquad\qquad X$^{2}$g$_{kl}$\thinspace$\ =$ $-\frac{2}{3}($ R$_{\text{XX}%
_{k}\text{XX}_{l}})$ $=$ $\underset{r,s=q+1}{\overset{n}{\sum}}[-\frac{2}{3}($
R$_{rksl})]$

Finally,

\qquad$\left(  \text{X}^{2}\text{g}_{kl}\right)  (y_{0})=$ $-\frac{2}%
{3}\underset{r,s=q+1}{\overset{n}{\sum}}($ R$_{rksl})(y_{0})$

We next compute the third coefficient:

From the computation of the second coefficient we have:

X$^{2}$g$_{kl}$\thinspace\thinspace\thinspace$=$ $[<\nabla_{\text{XX}}^{2}%
$X$_{k}$ , X$_{l}>+$ 2$<\nabla_{\text{X}}$X$_{k}$ , $\nabla_{\text{X}}$%
X$_{l}>$ $+$ $<\nabla_{\text{XX}}^{2}$ X$_{l},$X$_{k},>]$

Therefore,

X$^{3}$g$_{\alpha\beta}$\thinspace\thinspace\thinspace$=$ X$[<\nabla
_{\text{XX}}^{2}$X$_{k}$ , X$_{l}>+$ 2$<\nabla_{\text{X}}$X$_{k}$ ,
$\nabla_{\text{X}}$X$_{l}>$ $+$ $<\nabla_{\text{XX}}^{2}$ X$_{l},$X$_{k}>]$

\qquad$\ \ \ =$ $[<\nabla_{\text{XXX}}^{3}$X$_{k}$ , X$_{l}>+<\nabla
_{\text{XX}}^{2}$X$_{k}$,$\nabla_{\text{X}}$X$_{l}>+$ 2$<\nabla_{\text{XX}%
}^{2}$X$_{k}$ , $\nabla_{\text{X}}$X$_{l}>$

$\qquad\qquad+$ 2$<\nabla_{\text{X}}$X$_{k}$ , $\nabla_{\text{XX}}^{2}$%
X$_{l}>$

\qquad$\qquad+$ $<\nabla_{\text{XX}}^{3}$ X$_{l},$X$_{k}>]+$ $<\nabla
_{\text{XX}}^{2}$ X$_{l},\nabla_{\text{X}}$X$_{k}>]$

\qquad\ $=$ $[<\nabla_{\text{XXX}}^{3}$X$_{k}$ , X$_{l}>+3<\nabla_{\text{XX}%
}^{2}$X$_{k}$,$\nabla_{\text{X}}$X$_{l}>+$ $3<\nabla_{\text{XX}}^{2}$
X$_{l},\nabla_{\text{X}}$X$_{k}>]$

\qquad\ $+$ $<\nabla_{\text{XXX}}^{3}$ X$_{l},$X$_{k}>]$

We set:

\qquad\ A = $<\nabla_{\text{XXX}}^{3}$X$_{k}$ , X$_{l}>(y_{0});$ B =
$3<\nabla_{\text{XX}}^{2}$X$_{k}$,$\nabla_{\text{X}}$X$_{l}>(y_{0});$

\qquad\ C = $3<\nabla_{\text{XX}}^{2}$ X$_{l},\nabla_{\text{X}}$X$_{k}%
>(y_{0});$ D = $<\nabla_{\text{XXX}}^{3}$ X$_{l},$X$_{k}>(y_{0})$

We compute each of these using the Lemma above:

\qquad($\nabla_{\text{XXX}}^{3}$Y)(y$_{0}$) $=-\frac{1}{2}(\nabla_{\text{X}}%
($R)$_{\text{XY}}$X)(y$_{0}$).

For $k,l=q+1,...,n,$ we have:

\qquad A $=$ $<\nabla_{\text{XXX}}^{3}$X$_{k}$ , X$_{l}>(y_{0})=$ $-\frac
{1}{2}<(\nabla_{\text{X}}($R$)_{\text{XX}_{k}}$X$),$X$_{l}>(y_{0})$

$\qquad\ \ =-\frac{1}{2}(\nabla_{\text{X}}($R$)_{\text{XX}_{k}\text{XX}_{l}%
})(y_{0})$

Similarly,

\qquad D $=$ $<\nabla_{\text{XXX}}^{3}$ X$_{l},$X$_{k}>(y_{0})=-\frac{1}%
{2}(\nabla_{\text{X}}($R$)_{\text{XX}_{l}\text{XX}_{k}})(y_{0})=-\frac{1}%
{2}(\nabla_{\text{X}}($R$)_{\text{XX}_{k}\text{XX}_{l}})(y_{0})$

Since $(\nabla_{\text{X}}$X$_{k}>)(y_{0})=0=(\nabla_{\text{X}}$X$_{l}%
>)(y_{0})$, we have:

\qquad$\left(  \text{X}^{3}\text{g}_{kl}\right)  (y_{0})=-\frac{1}{2}%
(\nabla_{\text{X}}($R$)_{\text{XX}_{k}\text{XX}_{l}})(y_{0})-\frac{1}%
{2}(\nabla_{\text{X}}($R$)_{\text{XX}_{k}\text{XX}_{l}})(y_{0})=-$
$(\nabla_{\text{X}}($R$)_{\text{XX}_{k}\text{XX}_{l}})(y_{0})$

\qquad\qquad\qquad\qquad\qquad\qquad\qquad\qquad\qquad\qquad\qquad\qquad
\qquad\qquad\qquad\qquad\qquad$\blacksquare$

\begin{proposition}
For $k,l=q+1,...,n,$ we have:
\end{proposition}

g$^{kl}(x)=$ $\delta^{kl}+\frac{1}{3}\underset{r,s=q+1}{\overset{n}{\sum}}%
$R$_{\text{r}k\text{s}l}($y$_{0})$x$_{\text{r}}$x$_{\text{s}}+\frac{1}%
{6}\underset{r,s,t=q+1}{\overset{n}{\sum}}\nabla_{\text{r}}$R$_{\text{s}%
k\text{t}l}($y$_{0})x_{\text{r}}x_{\text{s}}x_{\text{t}}$

$\qquad\ \ \ -\frac{1}{360}\overset{n}{\underset{r,s,t,u=q+1}{\sum}}%
(-18\nabla_{\text{rs}}^{2}$R$_{\text{t}k\text{u}l}%
+16\underset{p=q+1}{\overset{n}{\sum}}$R$_{rksp}$R$_{tlup})($y$_{0}%
)x_{\text{r}}x_{\text{s}}x_{\text{t}}x_{\text{u}}\qquad\qquad$

\ + \ higher order terms.\qquad\qquad\qquad\qquad\qquad\qquad\qquad
\qquad\qquad$\qquad\qquad$

\begin{proof}
The expansion is easily proved using the fact that:
\end{proof}

g$_{\alpha\beta}$g$^{\beta\gamma}=\delta_{\alpha\gamma}\qquad$

The details of the proof are given in \textbf{Ndumu} $\left[  40\right]  ,$
p.137$.\qquad\qquad\qquad\qquad\qquad\qquad\qquad\qquad\qquad\qquad
\qquad\qquad\qquad$

$\qquad\qquad\qquad\qquad\qquad\qquad\qquad\qquad\qquad\qquad\qquad
\qquad\qquad\qquad\qquad\qquad\qquad\blacksquare$

\ Lastly we expand the volume change factor defined in $\left(  1.6\right)  :$

\qquad$\theta_{P}($x) =$\sqrt{\text{det(g}_{\alpha\beta}\text{(x)})}$ \ for
$\alpha,\beta=1,...,q,q+1,...,n.$\qquad

We will need at least five terms of the expansion.

\qquad\qquad\qquad\qquad\qquad\qquad\qquad\qquad\qquad\qquad\qquad\qquad
\qquad\qquad\qquad$\qquad\qquad\blacksquare$

\begin{proposition}
$\theta_{P}(x)=1-\underset{r\text{=q+1}}{\overset{n}{\sum}}<H,r>(y_{0})x_{r}$
\end{proposition}

$\qquad\qquad-\frac{1}{6}\underset{r,s=q+1}{\overset{n}{\sum}}[\varrho
_{rs}+\overset{q}{\underset{\text{a}=1}{2\sum}}R_{r\text{a}s\text{a}%
}-3\overset{q}{\underset{\text{a,b=1}}{\sum}}(T_{\text{aa}r}T_{\text{bb}%
s}-T_{\text{ab}r}T_{\text{ab}s})](y_{0})x_{r}x_{s}$

$-\frac{1}{12}\underset{r,s,t=q+1}{\overset{n}{\sum}}[\nabla_{r}\varrho
_{st}-2\varrho_{rs}<H,t>+\overset{q}{\underset{\text{a=1}}{\sum}}(\nabla
_{r}R_{\text{a}s\text{a}t}-4R_{r\text{a}s\text{a}}<H,t>)$

$+4\overset{q}{\underset{\text{a},\text{b}=1}{\sum}}R_{r\text{a}s\text{b}%
}T_{\text{ab}t}\ \ +2\overset{q}{\underset{\text{a},\text{b,c}=1}{\sum}%
}(T_{\text{aa}r}T_{\text{bb}s}T_{\text{cc}t}-3T_{\text{aa}r}TT_{\text{bc}%
s}T_{\text{bc}t}+2T_{\text{ab}r}T_{\text{bc}s}T_{\text{ca}t}](y_{0})x_{r}%
x_{s}x_{t}$

$+\ \ \frac{1}{24}\overset{n}{\underset{r,s,t,u=q+1}{\sum}}[$
$\overset{q}{\underset{\text{a=1}}{\sum}}\{-\nabla_{rs}^{2}(R)_{t\text{a}%
u\text{a}}+\overset{n}{\underset{p=q+1}{\sum}}\overset{q}{\underset{\text{a=1}%
}{\sum}}R_{\text{a}rsp}R_{\text{a}tup}+2\overset{q}{\underset{\text{a,b=1}%
}{\sum}}\nabla_{r}(R)_{\text{a}s\text{b}t}T_{\text{ab}u}\qquad A$

$+\overset{n}{\underset{p=q+1}{\sum}}(-\frac{3}{5}\nabla_{rs}^{2}%
(R)_{tpup}+\frac{1}{5}\overset{n}{\underset{m=q+1}{%
{\textstyle\sum}
}}R_{rpsm}R_{tpum})\}(y_{0})$

$+4\overset{q}{\underset{\text{a,b=1}}{\sum}}\{(\nabla_{r}(R)_{s\text{a}%
t\text{a}}-\overset{q}{\underset{\text{c=1}}{%
{\textstyle\sum}
}}R_{\text{a}r\text{c}s}T_{\text{ac}t})$ $T_{\text{bb}u}%
-4\overset{q}{\underset{\text{a,b=1}}{\sum}}(\nabla_{r}(R)_{s\text{a}%
t\text{b}}-\overset{q}{\underset{\text{c=1}}{%
{\textstyle\sum}
}}R_{\text{b}r\text{c}s}T_{\text{ac}t})T_{\text{ab}u}\}\qquad4B$

$+\frac{4}{3}\overset{q}{\underset{\text{a,b}=1}{\sum}}(R_{r\text{a}s\text{a}%
})(R_{t\text{b}u\text{b}})+\frac{1}{3}\varrho_{rs}\varrho_{tu}+\frac{2}%
{3}\overset{q}{\underset{\text{a}=1}{\sum}}R_{r\text{a}s\text{a}}\varrho
_{tu}\ +\frac{2}{3}\overset{q}{\underset{\text{b=1}}{\sum}}R_{r\text{b}%
s\text{b}}\varrho_{tu}=3C$

\vspace{1pt}$-3\overset{q}{\underset{\text{a,b}=1}{\sum}}R_{r\text{a}%
s\text{b}}R_{t\text{a}u\text{b}}-\frac{1}{3}%
\overset{n}{\underset{p,m=q+1}{\sum}}R_{rpsm}R_{tpum}%
\ -\overset{q}{\underset{\text{a}=1}{\sum}}\overset{n}{\underset{p=q+1}{\sum}%
}R_{r\text{a}sp}R_{t\text{a}up}-\overset{q}{\underset{\text{b}=1}{\sum}%
}\overset{n}{\underset{p=q+1}{\sum}}R_{r\text{b}sp}R_{t\text{b}up}$

$+$ $\overset{q}{\underset{\text{a,b,c=1}}{6\sum}}\{$ $-R_{r\text{a}s\text{a}%
}(T_{\text{bb}t}T_{\text{cc}u}$ $-T_{\text{bc}t}T_{\text{bc}u})\}$%
\ $+6\{$\ $R_{r\text{a}s\text{b}}(T_{\text{ab}t}T_{\text{cc}u}-T_{\text{bc}%
t}T_{\text{ac}u})\}\qquad6D$

$+6\{-R_{r\text{a}s\text{c}}(T_{\text{ba}t}T_{\text{bc}u}-T_{\text{ac}%
t}T_{\text{bb}u})\}+6\overset{q}{\underset{\text{b,c=1}}{\sum}}%
\underset{p=q+1}{\overset{n}{\sum}}$\ $\{-\frac{1}{3}R_{rpsp}(T_{\text{bb}%
t}T_{\text{cc}u}-T_{\text{bc}t}T_{\text{bc}u})\}\qquad$

$+\overset{n}{\underset{r,s,t,u=q+1}{\sum}}%
\overset{q}{\underset{\text{a,b,c,d}=1}{\sum}}T_{\text{aa}r}\{T_{\text{bb}%
s}(T_{\text{cc}t}T_{\text{dd}u}-T_{\text{cd}t}T_{\text{dc}u})-T_{\text{bc}%
s}(T_{\text{bc}t}T_{\text{dd}u}-T_{\text{bd}t}T_{\text{cd}u})$

$+T_{\text{bd}s}(T_{\text{bc}t}T_{\text{cd}u}-T_{\text{bd}t}T_{\text{cc}%
u})\}=E$

$-T_{\text{ab}r}\{T_{\text{ab}s}(T_{\text{cc}t}T_{\text{dd}u}-T_{\text{cd}%
t}T_{\text{dc}u})-T_{\text{bc}s}(T_{\text{ac}t}T_{\text{dd}u}-T_{\text{ad}%
t}T_{\text{cd}u})+T_{\text{bd}s}(T_{\text{ac}t}T_{\text{cd}u}-T_{\text{ad}%
t}T_{\text{cc}u})\}$

\vspace{1pt}$+T_{\text{ac}r}\{T_{\text{ab}s}(T_{\text{bc}t}T_{\text{dd}%
u}-T_{\text{bd}t}T_{\text{dc}u})-T_{\text{bb}s}(T_{\text{ac}t}T_{\text{dd}%
u}-T_{\text{ad}t}T_{\text{cd}u})+T_{\text{bd}s}(T_{\text{ac}t}T_{\text{bd}%
u}-T_{\text{ad}t}T_{\text{bc}u})\}$

\vspace{1pt}$-T_{\text{ad}r}\{T_{\text{ab}s}(T_{\text{bc}t}T_{\text{cd}%
u}-T_{\text{bd}t}T_{\text{cc}u})-T_{\text{bb}s}(T_{\text{ac}t}T_{\text{cd}%
u}-T_{\text{ad}t}T_{\text{cc}u})$

$+T_{\text{bc}s}(T_{\text{ac}t}T_{\text{bd}u}-T_{\text{ad}t}T_{\text{bc}%
u})\}]$(y$_{0}$)x$_{r}$x$_{s}$x$_{t}$x$_{u}$\ $+$ higher order terms.$\ $

\begin{proof}
We drop the subscript P from $\theta_{P}$ and write $\theta.$
\end{proof}

$\theta(x)=\theta(y_{0})+\frac{1}{1!}\overset{n}{\underset{i=q+1}{\sum}}%
($X$\theta)(y_{0})x_{r}+\frac{1}{2!}\overset{n}{\underset{r,s=q+1}{\sum}}%
($X$^{2}\theta)(y_{0})x_{r}x_{s}+\frac{1}{3!}%
\overset{n}{\underset{r,s,t=q+1}{\sum}}($X$^{3}\theta)(y_{0})x_{r}x_{s}x_{t}$

$\qquad+\frac{1}{4!}\overset{n}{\underset{r,s,t,u=q+1}{\sum}}($X$^{4}%
\theta)(y_{0})x_{r}x_{s}x_{t}x_{u}$ \ \ $+$ higher order terns.\ \ \ \ \ \ \ \ \ \ \ \ \ \ \ \ \ \ \ \ \ \ \ \ \ \ \ \ \ \ \ \ \ \ \ \ \ \ \ \ \ \ \ \ \ \ \ \ \ \ \ \ \ \ \ \ \ \ \ \ \ \ \ \ \ \ \ \ \ \ \ \ \ \ \ \ \ \ \ \ \ \ \ \ \ \ \ \ \ \ \ \ \ \ \ \ \ \ \ \ \ \ \ \ \ \ \ \ \ \ \ \ \ \ \ \ \ \ \ \ \ \ \ \ \ \ \ \ \ \ \ \ \ \ \ \ \ \ \ \ \ \ \ \ \ \ \ \ \ \ \ \ \ \ \ \ \ \ \ \ \ \ \ \ \ \ \ \ \ \ \ \ \ \ \ \ \ \ \ \ \ \ \ \ \ \ \ \ \ \ \ \ \ \ \ \ \ \ \ \ \ \ \ \ \ \ \ \ \ \ \ \ \ \ \ \ \ \ \ \ \ \ \ \ \ \ \ \ \ \ \ \ \ \ \ \ \ \ \ \ \ \ \ \ \ \ \ \ \ \ \ \ \ \ \ \ \ \ \ \ \ \ \ \ \ \ \ \ \ \ \ \ \ \ \ \ \ \ \ \ \ \ \ \ \ \ \ \ \ \ \ \ \ \ \ \ \ \ \ \ \ \ \ \ \ \ \ \ \ \ \ \ \ \ \ \ \ \ \ \ \ \ \ \ \ \ \ \ \ \ \ \ \ \ \ \ \ \ \ \ \ \ \ \ \ \ \ \ \ \ \ \ \ \ \ \ \ \ \ \ \ \ \ \ \ \ \ \ \ \ \ \ \ \ \ \ \ \ \ \ \ \ \ \ \ \ \ \ \ \ \ \ \ \ \ \ \ \ \ \ \ \ \ \ \ \ \ \ \ \ \ \ \ \ \ \ \ \ \ \ \ \ \ \ \ \ \ \ \ \ \ \ \ \ \ \ \ \ \ \ \ \ \ \ \ \ \ \ \ \ \ \ \ \ \ \ \ \ \ \ \ \ \ \ \ \ \ \ \ \ \ \ \ \ \ \ \ \ \ \ \ \ \ \ \ \ \ \ \ \ \ \ \ \ \ \ \ \ \ \ \ \ \ \ \ \ \ \ \ \ \ \ \ \ \ \ \ \ \ \ \ \ \ \ \ \ \ \ \ \ \ \ \ \ \ \ \ \ \ \ \ \ \ \ \ \ \ \ \ \ \ \ \ \ \ \ \ \ \ \ \ \ \ \ \ \ \ \ \ \ \ \ \ \ \ \ \ \ \ \ \ \ \ \ \ \ \ \ \ \ \ \ \ \ \ \ \ \ \ \ \ \ \ \ \ \ \ \ \ \ \ \ \ \ \ \ \ \ \ \ \ \ \ \ \ \ \ \ \ \ \ \ \ \ \ \ \ \ \ \ \ \ \ \ \ \ \ \ \ \ \ \ \ \ \ \ \ \ \ \ \ \ \ \ \ \ \ \ \ \ \ \ \ \ \ \ \ \ \ \ \ \ \ \ \ \ \ \ \ \ \ \ \ \ \ \ \ \ \ \ \ \ \ \ \ \ \ \ \ \ \ \ \ \ \ \ \ \ \ \ \ \ \ \ \ \ \ \ \ \ \ \ \ \ \ \ \ \ \ \ \ \ \ \ \ \ \ \ \ \ \ \ \ \ \ \ \ \ \ \ \ \ \ \ \ \ \ \ \ \ \ \ \ \ \ \ \ \ \ \ \ \ \ \ \ \ \ \ \ \ \ \ \ \ \ \ \ \ \ \ \ \ \ \ \ \ \ \ \ \ \ \ \ \ \ \ \ \ \ \ \ \ \ \ \ \ \ \ \ \ \ \ \ \ \ \ \ \ \ \ \ \ \ \ \ \ \ \ \ \ \ \ \ \ \ \ \ \ \ \ \ \ \ \ \ \ \ \ \ \ \ \ \ \ \ \ \ \ \ \ \ \ \ \ \ \ \ \ \ \ \ \ \ \ \ \ \ \ \ \ \ \ \ \ \ \ \ \ \ \ \ \ \ \ \ \ \ \ \ \ \ \ \ \ \ \ \ \ \ \ \ \ \ \ \ \ \ \ \ \ \ \ \ \ \ \ \ \ \ \ \ \ \ \ \ \ \ \ \ \ \ \ \ \ \ \ \ \ \ \ \ \ \ \ \ \ \ \ \ \ \ \ \ \ \ \ \ \ \ \ \ \ \ \ \ \ \ \ \ \ \ \ \ \ \ \ \ \ \ \ \ \ \ \ \ \ \ \ \ \ \ \ \ \ \ \ \ \ \ \ \ \ \ \ \ \ \ \ \ \ \ \ \ \ \ \ \ \ \ \ \ \ \ \ \ \ \ \ \ \ \ \ \ \ \ \ \ \ \ \ \ \ \ \ \ \ \ \ \ \ \ \ \ \ \ \ \ \ \ \ \ \ \ \ \ \ \ \ \ \ \ \ \ \ \ \ \ \ \ \ \ \ \ \ \ \ \ \ \ \ \ \ \ \ \ \ \ \ \ \ \ \ \ \ \ \ \ \ \ \ \ \ \ \ \ \ \ \ \ \ \ \ \ \ \ \ \ \ \ \ \ \ \ \ \ \ \ \ \ \ \ \ \ \ \ \ \ \ \ \ \ \ \ \ \ \ \ \ \ \ \ \ \ \ \ \ \ \ \ \ \ \ \ \ \ \ \ \ \ \ \ \ \ \ \ \ \ \ \ \ \ \ \ \ \ \ \ \ \ \ \ \ \ \ \ \ \ \ \ \ \ \ \ \ \ \ \ \ \ \ \ \ \ \ \ \ \ \ \ \ \ \ \ \ \ \ \ \ \ \ \ \ \ \ \ \ \ \ \ \ \ \ \ \ \ \ \ \ \ \ \ \ \ \ \ \ \ \ \ \ \ \ \ \ \ \ \ \ \ \ \ \ \ \ \ \ \ \ \ \ \ \ \ \ \ \ \ \ \ \ \ \ \ \ \ \ \ \ \ \ \ \ \ \ \ \ \ \ \ \ \ \ \ \ \ \ \ \ \ \ \ \ \ \ \ \ \ \ \ \ \ \ \ \ \ \ \ \ \ \ \ \ \ \ \ \ \ \ \ \ \ \ \ \ \ \ \ \ \ \ \ \ \ \ \ \ \ \ \ \ \ \ \ \ \ \ \ \ \ \ \ \ \ \ \ \ \ \ \ \ \ \ \ \ \ \ \ \ \ \ \ \ \ \ \ \ \ \ \ \ \ \ \ \ \ \ \ \ \ \ \ \ \ \ \ \ \ \ \ \ \ \ \ \ \ \ \ \ \ \ \ \ \ \ \ \ \ \ \ \ \ \ \ \ \ \ \ \ \ \ \ \ \ \ \ \ \ \ \ \ \ \ \ \ \ \ \ \ \ \ \ \ \ \ \ \ \ \ \ \ \ \ \ \ \ \ \ \ \ \ \ \ \ \ \ \ \ \ \ \ \ \ \ \ \ \ \ \ \ \ \ \ \ \ \ \ \ \ \ \ \ \ \ \ \ \ \ \ \ \ \ \ \ \ \ \ \ \ \ \ \ \ \ \ \ \ \ \ \ \ \ \ \ \ \ \ \ \ \ \ \ \ \ \ \ \ \ \ \ \ \ \ \ \ \ \ \ \ \ \ \ \ \ \ \ \ \ \ \ \ \ \ \ \ \ \ \ \ \ \ \ \ \ \ \ \ \ \ \ \ \ \ \ \ \ \ \ \ \ \ \ \ \ \ \ \ \ \ \ \ \ \ \ \ \ \ \ \ \ \ \ \ \ \ \ \ \ \ \ \ \ \ \ \ \ \ \ \ \ \ \ \ \ \ \ \ \ \ \ \ \ \ \ \ \ \ \ \ \ \ \ \ \ \ \ \ \ \ \ \ \ \ \ \ \ \ \ \ \ \ \ \ \ \ \ \ \ \ \ \ \ \ \ \ \ \ \ \ \ \ \ \ \ \ \ \ \ \ \ \ \ \ \ \ \ \ \ \ \ \ \ \ \ \ \ \ \ \ \ \ \ \ \ \ \ \ \ \ \ \ \ \ \ \ \ \ \ \ \ \ \ \ \ \ \ \ \ \ \ \ \ \ \ \ \ \ \ \ \ \ \ \ \ \ \ \ \ \ \ \ \ \ \ \ \ \ \ \ \ \ \ \ \ \ \ \ \ \ \ \ 

All terms up to order 3 (starting with the zeroth order term) have been given
in \textbf{Theorem }$4.3$ of \textbf{Gray} and \textbf{Vanhecke} $\left[
27\right]  .$ See also \textbf{Theorem} $9.22$ of \textbf{Gray} $\left[
25\right]  $ and \textbf{Problem 9.1} on p. 223 of \textbf{Gray }$\left[
25\right]  $\textbf{ }for the \textbf{third orde}r term.

I have \textbf{not seen the 4}$^{th}$\textbf{ order term anywhere} and so I
compute it here using the method of \textbf{Theorem }$9.22$ of \textbf{Gray}
$\left[  25\right]  .$

Let X \ be a normal Fermi vector field . Then,

X$^{4}\theta=\overset{n}{\underset{a=1}{\sum}}<\nabla_{XXXX}^{4}X_{a},X_{a}>$

$+$ $4\overset{n}{\underset{\alpha,\beta=1}{\sum}}\det\left(
\begin{array}
[c]{ll}%
<\nabla_{\text{XXX}}^{3}X_{a},X_{a}> & <\nabla_{\text{XXX}}^{3}X_{\alpha
},X_{\beta}>\\
<\nabla_{\text{X}}X_{\beta},X_{\alpha}> & <\nabla_{\text{X}}X_{\beta}%
,X_{\beta}>
\end{array}
\right)  $

\qquad$+3\overset{n}{\underset{\alpha,\beta=1}{\sum}}\det\left(
\begin{array}
[c]{ll}%
<\nabla_{\text{XX}}^{2}X_{a},X_{a}> & <\nabla_{\text{XXX}}^{2}X_{\alpha
},X_{\beta}>\\
<\nabla_{\text{XX}}^{2}X_{\beta},X_{a}> & <\nabla_{\text{XX}}^{2}X_{\beta
},X_{\beta}>
\end{array}
\right)  $

$+$ $6\overset{n}{\underset{\alpha,\beta,\gamma=1}{\sum}}\det\left(
\begin{array}
[c]{lll}%
<\nabla_{\text{XX}}^{2}X_{\alpha},X_{\alpha}> & <\nabla_{\text{XX}}%
^{2}X_{\alpha},X_{\beta}> & <\nabla_{\text{XX}}^{2}X_{\alpha},X_{\gamma}>\\
<\nabla_{\text{X}}X_{\beta},X_{\alpha}> & <\nabla_{\text{X}}X_{\beta}%
,X_{\beta}> & <\nabla_{\text{X}}X_{\beta},X_{\gamma}>\\
<\nabla_{\text{X}}X_{\gamma},X_{\alpha}> & <\nabla_{\text{X}}X_{\gamma
},X_{\beta}> & <\nabla_{\text{X}}X_{\gamma},X_{\gamma}>
\end{array}
\right)  $

$+\overset{n}{\underset{\alpha,\beta,\gamma,\delta=1}{\sum}}\det\left(
\begin{array}
[c]{llll}%
<\nabla_{\text{X}}X_{\alpha},X_{\alpha}> & <\nabla_{\text{X}}X_{\alpha
},X_{\beta}> & <\nabla_{\text{X}}X_{\alpha},X_{\gamma}> & <\nabla_{\text{X}%
}X_{\alpha},X_{\delta}>\\
<\nabla_{\text{X}}X_{\beta},X_{\alpha}> & <\nabla_{\text{X}}X_{\beta}%
,X_{\beta}> & <\nabla_{\text{X}}X_{\beta},X_{\gamma}> & <\nabla_{\text{X}%
}X_{\beta},X_{\delta}>\\
<\nabla_{\text{X}}X_{\gamma},X_{\alpha}> & <\nabla_{\text{X}}X_{\gamma
},X_{\beta}> & <\nabla_{\text{X}}X_{\gamma},X_{\gamma}> & <\nabla_{\text{X}%
}X_{\gamma},X_{\delta}>\\
<\nabla_{\text{X}}X_{\delta},X_{\alpha}> & <\nabla_{\text{X}}X_{\delta
},X_{\beta}> & <\nabla_{\text{X}}X_{\delta},X_{\gamma}> & <\nabla_{\text{X}%
}X_{\delta},X_{\delta}>
\end{array}
\right)  $

\qquad= A +4B + 3C + 6D + E

where A, B, C, D and E are the appropriate expressions above.

\qquad$A=$ $\overset{q}{\underset{\text{a}=1}{\sum}}<\nabla_{XXXX}%
^{4}X_{\text{a}},X_{\text{a}}>(y_{0})+\overset{n}{\underset{p=q+1}{\sum}%
}<\nabla_{XXXX}^{4}X_{p},X_{p}>(y_{0})$

\qquad By (iii) of \textbf{Lemma 6},

\qquad\ $\overset{q}{\underset{\text{a}=1}{\sum}}<\nabla_{XXXX}^{4}%
X_{\text{a}},X_{\text{a}}>$ $=$ $\overset{q}{\underset{a=1}{\sum}}%
\{-\nabla_{XX}^{2}(R)_{XX_{\text{a}}XX_{\text{a}}}+R_{XX_{\text{a}%
}XR_{XX_{\text{a}}}X}-2\nabla_{X}(R)_{XX_{\text{a}}XT_{X_{\text{a}}}X}%
-2\nabla_{X}(R)_{XX_{\text{a}}X\perp_{X_{\text{a}}}X}\}$

Since $\theta(x)$ \ is independent of the torsion operator by \textbf{Theorem
(8.1)} of \textbf{Gray and Vanhecke} $\left[  27\right]  $, the terms
containing the torsion operator $\perp$ in the expansion $\ $of $\theta(x)$
must disappear and we have:

\qquad By (vii) of \textbf{Lemma 6} above,

\qquad$\overset{n}{\underset{p=q+1}{\sum}}<\nabla_{XXXX}^{4}X_{p},X_{p}%
>(y_{0})$ $=\overset{n}{\underset{p=q+1}{\sum}}\{-\frac{3}{5}\nabla_{XX}%
^{2}(R)_{XX_{p}XX_{p}}+\frac{1}{5}R_{XX_{p}XR_{XX_{p}}X}\}(y_{0})$

\qquad Therefore,

$A=$ $\overset{q}{\underset{\text{a=1}}{\sum}}\{-\nabla_{XX}^{2}%
(R)_{XX_{\text{a}}XX_{\text{a}}}+R_{XX_{\text{a}}XR_{XX_{\text{a}}}X}%
-2\nabla_{X}(R)_{XX_{\text{a}}XT_{X_{\text{a}}}X}%
+\overset{n}{\underset{p=q+1}{\sum}}(-\frac{3}{5}\nabla_{XX}^{2}%
(R)_{XX_{p}XX_{p}}+\frac{1}{5}R_{XX_{p}XR_{XX_{p}}X})\}(y_{0})$

Here we will use the relations: R$_{\text{a}ij\text{T}_{\text{b}k}%
}=\overset{q}{\underset{\text{c=1}}{%
{\textstyle\sum}
}}R_{\text{a}i\text{c}j}T_{\text{bc}k}$ and R$_{i\text{a}jR_{k\text{b}l}}$ $=$
$\overset{n}{\underset{m=q+1}{%
{\textstyle\sum}
}}R_{\text{a}ijm}R_{\text{b}klm}:$

$A=$ $\overset{q}{\underset{\text{a=1}}{\sum}}\{-\nabla_{rs}^{2}%
(R)_{t\text{a}u\text{a}}+\overset{n}{\underset{p=q+1}{\sum}}%
\overset{q}{\underset{\text{a=1}}{\sum}}R_{\text{a}rsp}R_{\text{a}%
tup}+2\overset{q}{\underset{\text{a,b=1}}{\sum}}\nabla_{r}(R)_{\text{a}%
s\text{b}t}T_{\text{ab}u}$

$+\overset{n}{\underset{p=q+1}{\sum}}(-\frac{3}{5}\nabla_{rs}^{2}%
(R)_{tpup}+\frac{1}{5}\overset{n}{\underset{m=q+1}{%
{\textstyle\sum}
}}R_{rpsm}R_{tpum})\}(y_{0})$

$B=\overset{n}{\underset{\alpha,\beta=1}{\sum(}}<\nabla_{XXX}^{3}X_{\alpha
},X_{\alpha}><\nabla_{X}X_{\beta},X_{\beta}>-<\nabla_{XXX}^{3}X_{\alpha
},X_{\beta}><\nabla_{X}X_{\beta},X_{\alpha}>)$

\qquad By (ii) of \textbf{Lemma 6} (where we assume computations at y$_{0}),$

$\ B=\overset{q}{\underset{\text{a,b}=1}{\sum(}}<\nabla_{XXX}^{3}X_{\text{a}%
},X_{\text{a}}><\nabla_{X}X_{\text{b}},X_{\text{b}}>-<\nabla_{XXX}%
^{3}X_{\text{a}},X_{\text{a}}><\nabla_{X}X_{\text{b}},X_{\text{a}}>)$

$\ \ +\overset{q}{\underset{\text{b}=1}{\sum}}$
$\overset{n}{\underset{p=q+1}{\sum(}}<\nabla_{XXX}^{3}X_{p},X_{p}><\nabla
_{X}X_{\text{b}},X_{\text{b}}>-<\nabla_{XXX}^{3}X_{p},X_{\text{b}}><\nabla
_{X}X_{\text{b}},X_{p}>)$

\qquad By (iv) and (vi) of \textbf{Lemma 6}

$B=\overset{q}{\underset{\text{a,b}=1}{\sum}}(<-\nabla_{X}(R)_{XX_{\text{a}}%
}X-R_{XT_{X_{\text{a}}}X}X-R_{X\perp_{X_{\text{a}}}X}$ , $X_{\text{a}}%
>\times<T_{X_{\text{b}}}X+\perp_{X_{\text{b}}}X$ , $X_{\text{b}}>$

$\ \ \ \ -<-\nabla_{X}(R)_{XX_{\text{a}}}X-R_{XT_{X_{\text{a}}}X}%
X-R_{X\perp_{X_{\text{a}}}X}X$ , $X_{\text{b}}>\times<T_{X_{\text{b}}}%
X+\perp_{X_{\text{b}}}X$ , $X_{\text{a}}>)$

\ $\ +\overset{q}{\underset{b=1}{\sum}}$ $\overset{n}{\underset{p=q+1}{\sum(}%
}<-\frac{1}{2}\nabla_{X}(R)_{XX_{p}}X$ , $X_{p}><T_{X_{\text{b}}}%
X+\perp_{X_{\text{b}}}X$ , $X_{\text{b}}>$

\ $\ \ \ \ -$ $<-\frac{1}{2}\nabla_{X}(R)_{XX_{p}}X$ , $X_{\text{b}%
}><T_{X_{\text{b}}}X+\perp_{X_{\text{b}}}X$ , $X_{p}>$

$B\ =\overset{q}{\underset{\text{a,b=1}}{\sum}}(\nabla_{X}(R)_{XX_{\text{a}%
}XX_{\text{a}}}+R_{XX_{\text{a}}XT_{X_{\text{a}}}X}+R_{XX_{\text{a}}%
X\perp_{X_{\text{a}}}X})$ $T_{\text{bb}X}$

$\ \ \ \ -\overset{q}{\underset{\text{a,b=1}}{\sum}}(\nabla_{X}%
(R)_{XX_{\text{a}}XX_{\text{b}}}+R_{XT_{X_{\text{a}}}XXX_{b}}+R_{X\perp
_{X_{\text{a}}}XXX_{\text{b}}}$ $)T_{\text{ab}X}$

\ \ $\ +\overset{q}{\underset{\text{b=1}}{\sum}}$
$\overset{n}{\underset{p=q+1}{\sum}}(\frac{1}{2}\nabla_{X}(R)_{XX_{p}%
XX_{\text{b}}}$ $T_{\text{b}XXp}-\frac{1}{2}\nabla_{X}(R)_{XX_{p}XX_{\text{b}%
}})$ $\perp_{p\text{b}X}$

Since $\theta(x)$ \ is independent of the torsion operator by Theorem $(8.1)$
of $\ $\textbf{Gray and Vanhecke}$\ \left[  27\right]  $, the terms containing
the torsion operator $\perp$ in the expansion $\ $of $\theta(x)$ must
disappear. Further, $T_{\text{b}XXp}=0.$ We recall the following relations we
saw in the computation of \textbf{(vii)} of \textbf{Table A}$_{7}$:

R$_{\text{a}ij\text{T}_{\text{b}k}}=\overset{q}{\underset{\text{c=1}}{%
{\textstyle\sum}
}}R_{\text{a}i\text{c}j}T_{\text{bc}k};\qquad\overset{q}{\underset{\text{c=1}%
}{%
{\textstyle\sum}
}}R_{XT_{X_{\text{a}}}XXX_{\text{b}}}=\overset{q}{\underset{\text{c=1}}{%
{\textstyle\sum}
}}R_{XX_{\text{b}}XT_{X_{\text{a}}}X}=-\overset{q}{\underset{\text{c=1}}{%
{\textstyle\sum}
}}R_{\text{b}r\text{c}s}T_{\text{ac}t}$

$R_{XT_{X_{\text{a}}}XXX_{\text{b}}}=R_{XX_{\text{b}}XT_{X_{\text{a}}}%
X}=-\overset{q}{\underset{\text{c=1}}{%
{\textstyle\sum}
}}R_{\text{b}r\text{c}s}T_{\text{ac}t}$

Consequently,

$B\ =\overset{q}{\underset{\text{a,b=1}}{\sum}}(\nabla_{X}(R)_{XX_{\text{a}%
}XX_{\text{a}}}+R_{XX_{\text{a}}XT_{X_{\text{a}}}X})$ $T_{\text{bb}%
X}-\overset{q}{\underset{\text{a,b=1}}{\sum}}(\nabla_{X}(R)_{XX_{\text{a}%
}XX_{\text{b}}}+R_{XT_{X_{\text{a}}}XXX_{\text{b}}})T_{\text{ab}X}$

$B\ =\overset{q}{\underset{\text{a,b=1}}{\sum}}(\nabla_{r}(R)_{s\text{a}%
t\text{a}}-\overset{q}{\underset{\text{c=1}}{%
{\textstyle\sum}
}}R_{\text{a}r\text{c}s}T_{\text{ac}t})$ $T_{\text{bb}X}%
-\overset{q}{\underset{\text{a,b=1}}{\sum}}(\nabla_{r}(R)_{s\text{a}t\text{b}%
}-\overset{q}{\underset{\text{c=1}}{%
{\textstyle\sum}
}}R_{\text{b}r\text{c}s}T_{\text{ac}t})T_{\text{ab}X}$

\qquad$C\ =\ \overset{n}{\underset{\alpha,\beta=1}{\sum}}\det\left(
\begin{array}
[c]{ll}%
<\nabla_{\text{XX}}^{2}X_{\alpha},X_{\alpha}> & <\nabla_{\text{XXX}}%
^{2}X_{\alpha},X_{\beta}>\\
<\nabla_{\text{XX}}^{2}X_{\beta},X_{\alpha}> & <\nabla_{\text{XX}}^{2}%
X_{\beta},X_{\beta}>
\end{array}
\right)  $

$=\overset{n}{\underset{\alpha,\beta=1}{\sum}}\{<\nabla_{\text{XX}}%
^{2}X_{\alpha},X_{\alpha}><\nabla_{\text{XX}}^{2}X_{\beta},X_{\beta}%
>-<\nabla_{\text{XXX}}^{2}X_{\alpha},X_{\beta}>\times<\nabla_{\text{XX}}%
^{2}X_{\beta},X_{\alpha}>\}$

$=\overset{q}{\underset{\text{a,b}=1}{\sum}}\{<\nabla_{\text{XX}}%
^{2}X_{\text{a}},X_{\text{a}}><\nabla_{\text{XX}}^{2}X_{\text{b}},X_{\text{b}%
}>-<\nabla_{\text{XXX}}^{2}X_{\text{a}},X_{\text{b}}>\times<\nabla_{\text{XX}%
}^{2}X_{\text{b}},X_{\text{a}}>\}$

$+\overset{q}{\underset{\text{a}=1}{\sum}}\overset{n}{\underset{p=q+1}{\sum}%
}\{<\nabla_{\text{XX}}^{2}X_{\text{a}},X_{\text{a}}><\nabla_{\text{XX}}%
^{2}X_{p},X_{p}>-<\nabla_{\text{XXX}}^{2}X_{\text{a}},X_{p}>\times
<\nabla_{\text{XX}}^{2}X_{p},X_{\text{a}}>\}$

$+\overset{q}{\underset{b=1}{\sum}}\overset{n}{\underset{p=q+1}{\sum}%
}\{<\nabla_{\text{XX}}^{2}X_{p},X_{p}><\nabla_{\text{XX}}^{2}X_{\text{b}%
},X_{\text{b}}>-<\nabla_{\text{XXX}}^{2}X_{p},X_{\text{b}}>\times
<\nabla_{\text{XX}}^{2}X_{\text{b}},X_{p}>\}$

$+\overset{n}{\underset{p,m=q+1}{\sum}}\{<\nabla_{\text{XX}}^{2}X_{p}%
,X_{p}><\nabla_{\text{XX}}^{2}X_{m},X_{m}>-<\nabla_{\text{XXX}}^{2}X_{p}%
,X_{m}>\times<\nabla_{\text{XX}}^{2}X_{m},X_{p}>\}$

\qquad By (i) and (iv) of \textbf{Lemma 6},

$C=\overset{q}{\underset{\text{a,b=1}}{\sum\{}}<-R_{XX_{\text{a}}%
}X,X_{\text{a}}><-R_{XX_{\text{b}}}X,X_{\text{b}}>-<-R_{XX_{\text{a}}%
}X,X_{\text{b}}>\times<-R_{XX_{\text{b}}}X,X_{\text{a}}>\}$

$+\overset{q}{\underset{\text{a}=1}{\sum}}\overset{n}{\underset{p=q+1}{\sum
\{}}<-R_{XX_{\text{a}}}X,X_{\text{a}}><-\frac{1}{3}R_{XX_{p}}X,X_{p}%
>-<-R_{XX_{\text{a}}}X,X_{p}>\times<-\frac{1}{3}R_{XX_{p}}X,X_{\text{a}}>\}$

$+\overset{q}{\underset{\text{b}=1}{\sum}}\overset{n}{\underset{p=q+1}{\sum
\{}}<-\frac{1}{3}R_{XX_{p}}X,X_{p}><-R_{XX_{\text{b}}}X,X_{\text{b}}%
>-<-\frac{1}{3}R_{XX_{p}}X,X_{\text{b}}>\times<-R_{XX_{\text{b}}}X,X_{p}>\}$

$+\overset{n}{\underset{p,m=q+1}{\sum\{}}<-\frac{1}{3}R_{XX_{p}}%
X,X_{p}><-\frac{1}{3}R_{XX_{m}}X,X_{m}>-<-\frac{1}{3}R_{XXp}X,X_{m}%
>\times<-\frac{1}{3}R_{XX_{m}}X,X_{p}>\}$

\qquad$=$ $\overset{q}{\underset{\text{a,b}=1}{\sum}}\{(-R_{XX_{\text{a}%
}XX_{\text{a}}})(-R_{XX_{\text{b}}XX_{\text{b}}})-(-R_{XX_{\text{a}%
}XX_{\text{b}}})(-R_{XX_{\text{b}}XX_{\text{a}}})\}$

$\qquad+\overset{q}{\underset{\text{a}=1}{\sum}}%
\overset{n}{\underset{p=q+1}{\sum}}\{(-R_{XX_{\text{a}}XX_{\text{a}}}%
)(-\frac{1}{3}R_{XX_{p}XX_{p}})-(-R_{XX_{\text{a}}XX_{p}})(-\frac{1}%
{3}R_{XX_{p}XX_{\text{a}}})\}$

$\qquad\ +\overset{q}{\underset{\text{b}=1}{\sum}}%
\overset{n}{\underset{p=q+1}{\sum}}\{(-\frac{1}{3}R_{XX_{p}XX_{p}%
})(-R_{XX_{\text{b}}XX_{\text{b}}})-(-\frac{1}{3}R_{XX_{p}XX_{\text{b}}%
})(-R_{XX_{\text{b}}XX_{p}})\}$

$\qquad\ +\overset{n}{\underset{p,m=q+1}{\sum}}\{(-\frac{1}{3}R_{XX_{p}XX_{p}%
})(-\frac{1}{3}R_{XX_{m}XX_{m}})-(-\frac{1}{3}R_{XX_{p}XX_{m}})(-\frac{1}%
{3}R_{XX_{m}XX_{p}})\}$

\qquad$=$ $\ \overset{q}{\underset{\text{a,b}=1}{\sum}}\{R_{XX_{\text{a}%
}XX_{\text{a}}}R_{XX_{\text{b}}XX_{\text{b}}}-R_{XX_{\text{a}}XX_{\text{b}}%
}R_{XX_{\text{b}}XX_{\text{a}}}\}$

$\qquad+\frac{1}{3}\overset{q}{\underset{\text{a}=1}{\sum}}%
\overset{n}{\underset{p=q+1}{\sum}}\{R_{XX_{\text{a}}XX_{\text{a}}}%
R_{XX_{p}XX_{p}})-R_{XX_{\text{a}}XX_{p}}R_{XX_{p}XX_{\text{a}}}\}$

$\qquad\ +\frac{1}{3}\overset{q}{\underset{\text{b}=1}{\sum}}%
\overset{n}{\underset{p=q+1}{\sum}}\{R_{XX_{p}XX_{p}}R_{XX_{\text{b}%
}XX_{\text{b}}}-R_{XX_{p}XX_{\text{b}}}R_{XX_{\text{b}}XX_{p}}\}$

$\qquad+\frac{1}{9}\overset{n}{\underset{p,m=q+1}{\sum}}\{R_{XX_{p}XX_{p}%
}R_{XX_{m}XX_{m}}-R_{XX_{p}XX_{m}}R_{XX_{m}XX_{p}}\}$

$\qquad=$ $\ \overset{q}{\underset{\text{a,b}=1}{\sum}}\{R_{XX_{\text{a}%
}XX_{\text{a}}}R_{XX_{\text{b}}XX_{\text{b}}}-R_{XX_{\text{a}}XX_{\text{b}}%
}R_{XX_{\text{b}}XX_{\text{a}}}\}$

$\qquad+\frac{1}{3}\overset{q}{\underset{\text{a}=1}{\sum}}\{R_{XX_{\text{a}%
}XX_{\text{a}}}\overset{n}{(\underset{r=1}{\sum}}R_{XX_{p}XX_{p}%
}-\overset{q}{\underset{\text{b}=1}{\sum}}R_{XX_{\text{b}}XX_{\text{b}}%
})-R_{XX_{\text{b}}XX_{p}}R_{XX_{p}XX_{\text{b}}}\}$

$\qquad\ +\frac{1}{3}\overset{q}{\underset{\text{b}=1}{\sum\{}}%
\overset{n}{\underset{p=1}{(\sum}}R_{XX_{p}XX_{p}}%
-\overset{q}{\underset{\text{a}=1}{\sum}}R_{XX_{\text{a}}XX_{\text{a}}%
})R_{XX_{\text{b}}XX_{\text{b}}}-R_{XX_{p}XX_{\text{b}}}R_{XX_{\text{b}}%
XX_{p}}\}$

$\qquad+\frac{1}{9}\{\overset{n}{\underset{p=1}{(\sum}}R_{XX_{p}XX_{p}%
}-\overset{q}{\underset{\text{a}=1}{\sum}}R_{XX_{\text{a}}XX_{\text{a}}%
})(\overset{n}{\underset{m=1}{\sum}}R_{XX_{m}XX_{m}}%
-\overset{q}{\underset{\text{b}=1}{\sum}}R_{XX_{\text{b}}XX_{\text{b}}%
})-\overset{n}{\underset{p,m=q+1}{\sum}}R_{XX_{p}XX_{m}}R_{XX_{m}XX_{p}}\}$

$\qquad=$ $\ \overset{q}{\underset{\text{a,b}=1}{\sum}}\{R_{XX_{\text{a}%
}XX_{\text{a}}}R_{XX_{\text{b}}XX_{\text{b}}}-R_{XX_{\text{a}}XX_{\text{b}}%
}R_{XX_{\text{b}}XX_{\text{a}}}\}$

$\qquad+\frac{1}{3}\overset{q}{\underset{\text{a}=1}{\sum}}\{R_{XX_{\text{a}%
}XX_{\text{a}}}\varrho(X,X)-\overset{q}{\underset{\text{b}=1}{\sum}%
}R_{XX_{\text{b}}XX_{\text{b}}})-R_{XX_{\text{b}}XX_{p}}R_{XX_{p}XX_{\text{b}%
}}\}$

$\qquad\ +\frac{1}{3}\overset{q}{\underset{\text{b}=1}{\sum\{}(}%
\varrho(X,X)-\overset{q}{\underset{\text{a}=1}{\sum}}R_{XX_{\text{a}%
}XX_{\text{a}}})R_{XX_{\text{b}}XX_{\text{b}}}-R_{XX_{p}XX_{\text{b}}%
}R_{XX_{\text{b}}XX_{p}}\}$

$\qquad+\frac{1}{9}\{(\varrho(X,X)-\overset{q}{\underset{\text{a}=1}{\sum}%
}R_{XX_{\text{a}}XX_{\text{a}}})(\varrho(X,X)-\overset{q}{\underset{\text{b}%
=1}{\sum}}R_{XX_{\text{b}}XX_{\text{b}}})-\overset{n}{\underset{p,m=q+1}{\sum
}}R_{XX_{p}XX_{m}}R_{XX_{m}XX_{p}}\}$

\qquad$=$ $\ \overset{q}{\underset{\text{a,b}=1}{\sum}}\{R_{XX_{\text{a}%
}XX_{\text{a}}}R_{XX_{\text{b}}XX_{\text{b}}}-R_{XX_{\text{a}}XX_{\text{b}}%
}R_{XX_{\text{a}}XX_{\text{b}}}\}$

$\qquad+\frac{1}{3}\overset{q}{\underset{\text{a}=1}{\sum}}\{\varrho
(X,X)R_{XX_{\text{a}}XX_{\text{a}}}-\overset{q}{\underset{\text{a,b}=1}{\sum}%
}R_{XX_{\text{a}}XX_{\text{a}}}R_{XX_{\text{b}}XX_{\text{b}}})-R_{XX_{\text{b}%
}XX_{p}}R_{XX_{p}XX_{\text{b}}}\}$

$\qquad\ +\frac{1}{3}\overset{q}{\underset{\text{b}=1}{\sum\{}}\varrho
(X,X)R_{XX_{\text{b}}XX_{\text{b}}}-\overset{q}{\underset{\text{a,b}=1}{\sum}%
}R_{XX_{\text{a}}XX_{\text{a}}}R_{XX_{\text{b}}XX_{\text{b}}}-R_{XX_{p}%
XX_{\text{b}}}R_{XX_{\text{b}}XX_{p}}\}$

$\qquad+\frac{1}{9}\{(\varrho(X,X)\varrho(X,X)-\varrho
(X,X)\overset{q}{\underset{\text{a}=1}{\sum}}R_{XX_{\text{a}}XX_{\text{a}}%
}-\varrho(X,X)\overset{q}{\underset{\text{b}=1}{\sum}}R_{XX_{\text{b}%
}XX_{\text{b}}}+\overset{q}{\underset{\text{a,b}=1}{\sum}}\{R_{XX_{\text{a}%
}XX_{\text{a}}}R_{XX_{\text{b}}XX_{\text{b}}}$

$\qquad-\overset{n}{\underset{p,m=q+1}{\sum}}R_{XX_{p}XX_{m}}R_{XX_{m}XX_{p}%
}\}$

$\qquad C=\frac{4}{9}\overset{q}{\underset{\text{a,b}=1}{\sum}}%
(R_{XX_{\text{a}}XX_{\text{a}}})(R_{XX_{\text{b}}XX_{\text{b}}})+\frac{1}%
{9}\varrho(X,X)\varrho(X,X)$

\qquad$\ +\frac{2}{9}\overset{q}{\underset{\text{a}=1}{\sum}}\{\varrho
(X,X)R_{XX_{\text{a}}XX_{\text{a}}}\}\ +\frac{2}{9}%
\overset{q}{\underset{b=1}{\sum}}\{\varrho(X,X)R_{XX_{\text{b}}XX_{\text{b}}%
}\}$

\vspace{1pt}\qquad$\ -\overset{q}{\underset{\text{a,b}=1}{\sum}}%
R_{XX_{\text{a}}XX_{\text{b}}}R_{XX_{\text{a}}XX_{\text{b}}}-\frac{1}%
{9}\overset{n}{\underset{p,m=q+1}{\sum}}R_{XX_{p}XX_{m}}R_{XX_{p}XX_{m}}$

\vspace{1pt}\qquad$\ -\frac{1}{3}\overset{q}{\underset{\text{a}=1}{\sum}%
}\overset{n}{\underset{p=q+1}{\sum}}R_{XX_{\text{a}}XX_{p}}R_{XX_{\text{a}%
}XX_{p}}-\frac{1}{3}\overset{q}{\underset{\text{b}=1}{\sum}}%
\overset{n}{\underset{p=q+1}{\sum}}R_{XX_{\text{b}}XX_{p}}R_{XX_{\text{b}%
}XX_{p}}$

\qquad

\vspace{1pt}$D=\overset{n}{\underset{\alpha,\beta,\gamma=1}{\sum}}\det\left(
\begin{array}
[c]{lll}%
<\nabla_{\text{XX}}^{2}X_{\alpha},X_{\alpha}> & <\nabla_{\text{XX}}%
^{2}X_{\alpha},X_{\beta}> & <\nabla_{\text{XX}}^{2}X_{\alpha},X_{\gamma}>\\
<\nabla_{\text{X}}X_{\beta},X_{\alpha}> & <\nabla_{\text{X}}X_{\beta}%
,X_{\beta}> & <\nabla_{\text{X}}X_{\beta},X_{\gamma}>\\
<\nabla_{\text{X}}X_{\gamma},X_{\alpha}> & <\nabla_{\text{X}}X_{\gamma
},X_{\beta}> & <\nabla_{\text{X}}X_{\gamma},X_{\gamma}>
\end{array}
\right)  $

\vspace{1pt}$=\overset{n}{\underset{\alpha,\beta,\gamma=1}{\sum}}%
\{<\nabla_{\text{XX}}^{2}X_{\alpha},X_{\alpha}>(<\nabla_{\text{X}}X_{\beta
},X_{\beta}><\nabla_{\text{X}}X_{\gamma},X_{\gamma}>$

$-<\nabla_{\text{X}}X_{\beta},X_{\gamma}><\nabla_{\text{X}}X_{\gamma}%
,X_{\beta}>)\}-\{<\nabla_{\text{XX}}^{2}X_{\alpha},X_{\beta}>(<\nabla
_{\text{X}}X_{\beta},X_{\alpha}><\nabla_{\text{X}}X_{\gamma},X_{\gamma}>$

$-$ $<\nabla_{\text{X}}X_{\beta},X_{\gamma}><\nabla_{\text{X}}X_{\gamma
},X_{\alpha}>)\}+\{<\nabla_{\text{XX}}^{2}X_{\alpha},X_{\gamma}>(<\nabla
_{\text{X}}X_{\beta},X_{\alpha}><\nabla_{\text{X}}X_{\gamma},X_{\beta}>$

$-<\nabla_{\text{X}}X_{\gamma},X_{\alpha}><\nabla_{\text{X}}X_{\beta}%
,X_{\beta}>)\}$

\vspace{1pt}\qquad Using the fact that: $\nabla_{\text{X}}X_{s}=0=\nabla
_{\text{X}}X_{t}$ \ at y$_{0}$ for $s,t=q+1,...,n$ and any norml Fermi field
X, we have,

$D=\overset{q}{\underset{\text{a,b,c}=1}{\sum}}[\{<\nabla_{\text{XX}}%
^{2}X_{\text{a}},X_{\text{a}}>(<\nabla_{\text{X}}X_{\text{b}},X_{\text{b}%
}><\nabla_{\text{X}}X_{\text{c}},X_{\text{c}}>$

$-<\nabla_{\text{X}}X_{\text{b}},X_{\text{c}}><\nabla_{\text{X}}X_{\text{c}%
},X_{\text{b}}>)\}-\{<\nabla_{\text{XX}}^{2}X_{\text{a}},X_{\text{b}}%
>(<\nabla_{\text{X}}X_{\text{b}},X_{\text{a}}><\nabla_{\text{X}}X_{\text{c}%
},X_{\text{c}}>$

$-<\nabla_{\text{X}}X_{\text{b}},X_{\text{c}}><\nabla_{\text{X}}X_{\text{c}%
},X_{\text{a}}>)\}+$ $\{<\nabla_{\text{XX}}^{2}X_{\text{a}},X_{\text{c}%
}>(<\nabla_{\text{X}}X_{\text{b}},X_{\text{a}}><\nabla_{\text{X}}X_{\text{c}%
},X_{\text{b}}>$

$\qquad-<\nabla_{\text{X}}X_{\text{c}},X_{\text{a}}><\nabla_{\text{X}%
}X_{\text{b}},X_{\text{b}}>)\}]$

\vspace{1pt}$+\overset{q}{\underset{\text{b,c}=1}{\sum}}%
\underset{p=q+1}{\overset{n}{\sum}}[\{<\nabla_{\text{XX}}^{2}X_{p}%
,X_{p}>(<\nabla_{\text{X}}X_{\text{b}},X_{\text{b}}><\nabla_{\text{X}%
}X_{\text{c}},X_{\text{c}}>$

$-<\nabla_{\text{X}}X_{\text{b}},X_{\text{b}}><\nabla_{\text{X}}X_{\text{c}%
},X_{\text{b}}>)\}-\{<\nabla_{\text{XX}}^{2}X_{p},X_{\text{b}}>(<\nabla
_{\text{X}}X_{\text{b}},X_{p}><\nabla_{\text{X}}X_{\text{c}},X_{\text{c}}>$

$-<\nabla_{\text{X}}X_{\text{b}},X_{\text{c}}><\nabla_{\text{X}}X_{\text{c}%
},X_{p}>)\}$ $+\{<\nabla_{\text{XX}}^{2}X_{p},X_{\text{c}}>(<\nabla_{\text{X}%
}X_{\text{b}},X_{p}><\nabla_{\text{X}}X_{\text{c}},X_{\text{b}}>$

$-<\nabla_{\text{X}}X_{\text{c}},X_{p}><\nabla_{\text{X}}X_{\text{b}%
},X_{\text{b}}>)\}]$

\qquad By (i) and (v) of \textbf{Lemma 6},\qquad

$D=$ $\overset{q}{\underset{\text{a,b,c=1}}{\sum}}\{$ $-R_{XX_{\text{a}%
}XX_{\text{a}}}(T_{\text{bb}X}T_{\text{cc}X}$ $-T_{\text{bc}X}T_{\text{cb}%
X})\}$\ $-\{$\ $-R_{XX_{\text{a}}XX_{\text{b}}}(T_{\text{ba}X}T_{\text{cc}%
X}-T_{\text{bc}X}T_{\text{ca}X})\}$

\qquad\ \ \ $+\{$\ $-R_{XX_{\text{a}}XX_{\text{c}}}(T_{\text{ba}X}%
T_{\text{cb}X}-T_{\text{ca}X}T_{\text{bb}X})\}$

\qquad\ \ $\ +\overset{q}{\underset{\text{b,c=1}}{\sum}}%
\underset{p=q+1}{\overset{n}{\sum}}$\ $\{-\frac{1}{3}R_{X_{p}XX_{p}%
}(T_{\text{bb}X}T_{\text{cc}X}-T_{\text{bc}X}T_{\text{cb}X})\}$

\qquad$\{-R_{XX_{p}XX_{\text{b}}}(T_{\text{b}pX}T_{\text{cc}X}-T_{\text{bc}%
X}T_{\text{c}rX})\}\{-R_{XX_{p}XX_{\text{c}}}(T_{\text{b}pX}T_{\text{cb}%
X}-T_{\text{c}rX}T_{\text{bb}X})\}$

$T_{\text{b}pX}=T_{X_{\text{b}}\text{X}_{p}X}=0$ because $T_{X_{\text{b}%
}\text{X}_{p}}$ is tangential and $X$ is normal. Similarly, $T_{\text{c}%
rX}=0.$ Therefore, the last line is wiped off:

$\{-R_{XX_{p}XX_{\text{b}}}(T_{\text{b}pX}T_{\text{cc}X}-T_{\text{bc}%
X}T_{\text{c}rX})\}\{-R_{XX_{p}XX_{\text{c}}}(T_{\text{b}pX}T_{\text{cb}%
X}-T_{\text{c}rX}T_{\text{bb}X})\}=0$ and so,

$D=$ $\overset{q}{\underset{\text{a,b,c=1}}{\sum}}\{$ $-R_{XX_{\text{a}%
}XX_{\text{a}}}(T_{\text{bb}X}T_{\text{cc}X}$ $-T_{\text{bc}X}T_{\text{bc}%
X})\}$\ $+\{$\ $R_{XX_{\text{a}}XX_{\text{b}}}(T_{\text{ab}X}T_{\text{cc}%
X}-T_{\text{bc}X}T_{\text{ac}X})\}$

\qquad\ \ \ \ $\{-R_{XX_{\text{a}}XX_{\text{c}}}(T_{\text{bc}X}T_{\text{bc}%
X}-T_{\text{ac}X}T_{\text{bb}X})\}+\overset{q}{\underset{\text{b,c=1}}{\sum}%
}\underset{p=q+1}{\overset{n}{\sum}}$\ $\{-\frac{1}{3}R_{X_{p}XX_{p}%
}(T_{\text{bb}X}T_{\text{cc}X}-T_{\text{bc}X}T_{\text{bc}X})\}$

$D=$ $\overset{q}{\underset{\text{a,b,c=1}}{\sum}}\{$ $-R_{r\text{a}s\text{a}%
}(T_{\text{bb}t}T_{\text{cc}u}$ $-T_{\text{bc}t}T_{\text{bc}u})\}$%
\ $+\{$\ $R_{r\text{a}s\text{b}}(T_{\text{ab}t}T_{\text{cc}u}-T_{\text{bc}%
t}T_{\text{ac}u})\}$

\qquad\ \ \ \ $\{-R_{r\text{a}s\text{c}}(T_{\text{bc}t}T_{\text{bc}%
u}-T_{\text{ac}t}T_{\text{bb}u})\}+\overset{q}{\underset{\text{b,c=1}}{\sum}%
}\underset{p=q+1}{\overset{n}{\sum}}$\ $\{-\frac{1}{3}R_{rpsp}(T_{\text{bb}%
t}T_{\text{cc}u}-T_{\text{bc}t}T_{\text{bc}u})\}$

\vspace{1pt}$E=\overset{n}{\underset{\alpha,\beta,\gamma,\delta=1}{\sum}}%
\det\left(
\begin{array}
[c]{llll}%
<\nabla_{\text{X}}X_{\alpha},X_{\alpha}> & <\nabla_{\text{X}}X_{\alpha
},X_{\beta}> & <\nabla_{\text{X}}X_{\alpha},X_{\gamma}> & <\nabla_{\text{X}%
}X_{\alpha},X_{\delta}>\\
<\nabla_{\text{X}}X_{\beta},X_{\alpha}> & <\nabla_{\text{X}}X_{\beta}%
,X_{\beta}> & <\nabla_{\text{X}}X_{\beta},X_{\gamma}> & <\nabla_{\text{X}%
}X_{\beta},X_{\delta}>\\
<\nabla_{\text{X}}X_{\gamma},X_{\alpha}> & <\nabla_{\text{X}}X_{\gamma
},X_{\beta}> & <\nabla_{\text{X}}X_{\gamma},X_{\gamma}> & <\nabla_{\text{X}%
}X_{\gamma},X_{\delta}>\\
<\nabla_{\text{X}}X_{\delta},X_{\alpha}> & <\nabla_{\text{X}}X_{\delta
},X_{\beta}> & <\nabla_{\text{X}}X_{\delta},X_{\gamma}> & <\nabla_{\text{X}%
}X_{\delta},X_{\delta}>
\end{array}
\right)  $

\qquad\qquad

\vspace{1pt}\qquad Because $\nabla_{\text{X}}Y=0$ at y$_{0}$ for all normal
Fermi fields X and Y, we have:

$E=$ $-\overset{q}{\underset{\text{a,b,c,d}=1}{\sum}}T_{\text{aa}X}\left[
\begin{array}
[c]{c}%
-T_{\text{bb}X}-T_{\text{bc}X}-T_{\text{bd}X}\\
-T_{\text{cb}X}-T_{\text{cc}X}-T_{\text{cd}X}\\
-T_{\text{db}X}-T_{\text{dc}X}-T_{\text{dd}X}%
\end{array}
\right]  +\overset{q}{\underset{\text{a,b,c,d}=1}{\sum}}T_{\text{ab}X}\left[
\begin{array}
[c]{c}%
-T_{\text{ba}X}-T_{\text{bc}X}-T_{\text{bd}X}\\
-T_{\text{ca}X}-T_{\text{cc}X}-T_{\text{cd}X}\\
-T_{\text{da}X}-T_{\text{dc}X}-T_{\text{dd}X}%
\end{array}
\right]  $

\vspace{1pt}

\vspace{1pt}$-$ $\overset{q}{\underset{\text{a,b,c,d}=1}{\sum}}T_{\text{ac}%
X}\left[
\begin{array}
[c]{c}%
-T_{\text{ba}X}-T_{\text{bb}X}-T_{\text{bd}X}\\
-T_{\text{ca}X}-T_{\text{cb}X}-T_{\text{cd}X}\\
-T_{\text{da}X}-T_{\text{db}X}-T_{\text{dd}X}%
\end{array}
\right]  +$ $\overset{q}{\underset{\text{a,b,c,d}=1}{\sum}}T_{\text{ad}%
X}\left[
\begin{array}
[c]{c}%
-T_{\text{ba}X}-T_{\text{bb}X}-T_{\text{bc}X}\\
-T_{\text{ca}X}-T_{\text{cb}X}-T_{\text{cc}X}\\
-T_{\text{da}X}-T_{\text{db}X}-T_{\text{dc}X}%
\end{array}
\right]  $

\vspace{1pt}

\vspace{1pt}$=$\qquad$\overset{q}{\underset{\text{a,b,c,d}=1}{\sum}%
}T_{\text{aa}X}\left[
\begin{array}
[c]{c}%
T_{\text{bb}X}\text{ }T_{\text{bc}X}\text{ }T_{\text{bd}X}\\
T_{\text{cb}X}\text{ }T_{\text{cc}X}\text{ }T_{\text{cd}X}\\
\text{ }T_{\text{db}X}\text{ }T_{\text{dc}X}\text{ }T_{\text{dd}X}%
\end{array}
\right]  -\overset{q}{\underset{\text{a,b,c,d}=1}{\sum}}T_{\text{ab}X}\left[
\begin{array}
[c]{c}%
T_{\text{ba}X}\text{ }T_{\text{bc}X}\text{ }T_{\text{bd}X}\\
\text{ }T_{\text{ca}X}\text{ }T_{\text{cc}X}\text{ }T_{\text{cd}X}\\
\text{ }T_{\text{da}X}\text{ }T_{\text{dc}X}\text{ }T_{\text{dd}X}%
\end{array}
\right]  $

\vspace{1pt}

\vspace{1pt}\qquad$+\overset{q}{\underset{\text{a,b,c,d}=1}{\sum}}$
$T_{\text{ac}X}\left[
\begin{array}
[c]{c}%
T_{\text{ba}X}\text{ }T_{\text{bb}X}\text{ }T_{\text{bd}X}\\
\text{ }T_{\text{ca}X}\text{ }T_{\text{cb}X}\text{ }T_{\text{cd}X}\\
\text{ }T_{\text{da}X}\text{ }T_{\text{db}X}\text{ }T_{\text{dd}X}%
\end{array}
\right]  -\overset{q}{\underset{\text{a,b,c,d}=1}{\sum}}T_{\text{ad}X}\left[
\begin{array}
[c]{c}%
\text{ }T_{\text{ba}X}\text{ }T_{\text{bb}X}\text{ }T_{\text{bc}X}\\
\text{ }T_{\text{ca}X}\text{ }T_{\text{cb}X}\text{ }T_{\text{cc}X}\\
\text{ }T_{\text{da}X}\text{ }T_{\text{db}X}\text{ }T_{\text{dc}X}%
\end{array}
\right]  $

$=\overset{q}{\underset{\text{a,b,c,d}=1}{\sum}}[T_{\text{aa}X}\{T_{\text{bb}%
X}(T_{\text{cc}X}T_{\text{dd}X}-T_{\text{cd}X}T_{\text{dc}X})-T_{\text{bc}%
X}(T_{\text{bc}X}T_{\text{dd}X}-T_{\text{bd}X}T_{\text{cd}X})$

$\qquad+T_{\text{bd}X}(T_{\text{bc}X}T_{\text{cd}X}-T_{\text{bd}X}%
T_{\text{cc}X})\}$

\vspace{1pt}$-T_{\text{ab}X}\{T_{\text{ab}X}(T_{\text{cc}X}T_{\text{dd}%
X}-T_{\text{cd}X}T_{\text{dc}X})-T_{\text{bc}X}(T_{\text{ac}X}T_{\text{dd}%
X}-T_{\text{ad}X}T_{\text{cd}X})$

$\qquad+T_{\text{bd}X}(T_{\text{ac}X}T_{\text{cd}X}-T_{\text{ad}X}%
T_{\text{cc}X})\}$

\vspace{1pt}$+T_{\text{ac}X}\{T_{\text{ab}X}(T_{\text{bc}X}T_{ddX}%
-T_{\text{bd}X}T_{\text{dc}X})-T_{\text{bb}X}(T_{\text{ac}X}T_{\text{dd}%
X}-T_{\text{ad}X}T_{\text{cd}X})$

$\qquad+T_{\text{bd}X}(T_{\text{ac}X}T_{\text{bd}X}-T_{\text{ad}X}%
T_{\text{bc}X})\}$

\vspace{1pt}$-T_{\text{ad}X}\{T_{\text{ab}X}(T_{\text{bc}X}T_{\text{cd}%
X}-T_{\text{bd}X}T_{\text{cc}X})-T_{\text{bb}X}(T_{\text{ac}X}T_{\text{cd}%
X}-T_{\text{ad}X}T_{\text{cc}X})$

$\qquad+T_{\text{bc}X}(T_{\text{ac}X}T_{\text{bd}X}-T_{\text{ad}X}%
T_{\text{bc}X})\}]$

We now gather the final expressions for A, B, C, D and E to get the expression
for the 5$^{\text{th}}$ term = A + 4B + 3C + 6D + E.\qquad\qquad\qquad
\qquad\qquad\qquad\qquad\qquad\qquad\qquad\qquad\qquad\qquad\qquad\qquad
\qquad$\qquad\qquad\qquad\qquad\qquad\qquad\qquad\qquad\qquad\qquad
\qquad\qquad\qquad\qquad\qquad\qquad\qquad\qquad\qquad\qquad\blacksquare$

When the submanifold reduces to the centre of Fermi coordinates $\left\{
y_{0}\right\}  ,$ then the Fermi coordinates reduce to normal coordinates and
we have:

\begin{corollary}
$\theta_{P}(x)=1-$ $\frac{1}{6}\underset{r,s=1}{\overset{n}{\sum}}\varrho
_{rs}(y_{0})x_{r}x_{s}-\frac{1}{12}\underset{r,s,t\text{=1}}{\overset{n}{\sum
}}\nabla_{r}\varrho_{st}(y_{0})x_{r}x_{s}x_{t}$
\end{corollary}

\ $+$ $\ \ \frac{1}{24}\overset{n}{\underset{r,s,t,p=1}{\sum}}[$
$\overset{n}{\underset{u=1}{\sum}}(-\frac{3}{5}\nabla_{rs}^{2}R_{tupu}%
+\frac{1}{5}\overset{n}{\underset{p,h=1}{\sum}}R_{rush}R_{tuph})$%
\ $\ +\frac{1}{3}\varrho_{rs}\varrho_{tp}-\frac{1}{3}%
\overset{n}{\underset{p,h=1}{\sum}}R_{rush}R_{tuph}](y_{0})x_{r}x_{s}%
x_{k}x_{p}\ $

\ + higher order terms.$\ \ $

\begin{proof}
We simplify:
\end{proof}

$\qquad\frac{1}{5}\overset{n}{\underset{p,h=1}{\sum}}R_{rush}R_{tuph}-\frac
{1}{3}\overset{n}{\underset{p,h=1}{\sum}}R_{rush}R_{tuph}=-\frac{2}%
{15}\overset{n}{\underset{p,h=1}{\sum}}R_{rush}R_{tuph}$

and have:

$\qquad\theta_{P}(x)=1-$ $\frac{1}{6}\underset{r,s=1}{\overset{n}{\sum}%
}\varrho_{rs}(y_{0})x_{r}x_{s}-\frac{1}{12}\underset{r,s,t=1}{\overset{n}{\sum
}}\nabla_{r}\varrho_{st}(y_{0})x_{r}x_{s}x_{t}$

\ $+$ $\ \ \frac{1}{24}\overset{n}{\underset{r,s,t,p=1}{\sum}}[$
$\overset{n}{\underset{u=1}{\sum}}(-\frac{3}{5}\nabla_{rs}^{2}\varrho
_{tp}-\frac{2}{15}\overset{n}{\underset{p,h=1}{\sum}}R_{rush}R_{tuph}%
)\ \ +\frac{1}{3}\varrho_{rs}\varrho_{tp}](y_{0})x_{r}x_{s}x_{k}x_{p}$\ 

\qquad\qquad\qquad\qquad\qquad\qquad\qquad\qquad\qquad\qquad\qquad\qquad
\qquad\qquad\qquad\qquad\qquad$\blacksquare$

This ties up with the expression in \textbf{Corollary }$9.9$ of \textbf{Gray
}$\left[  25\right]  $

\qquad\qquad\qquad\qquad\qquad\qquad\qquad\qquad\qquad\qquad\qquad\qquad
\qquad\qquad\qquad\qquad\qquad$\blacksquare$\qquad\qquad\qquad\qquad
\qquad\qquad\qquad\qquad\qquad\qquad\qquad\qquad\qquad\qquad\qquad\qquad

Let $\left(  \mu_{1}(x),...,\mu_{d}(x)\right)  $ the basis of the fiber
E$_{x}$\textbf{ }of the vector bundle E over the chart U based at the point
$x\in M_{0}$ as defined in \textbf{Chapter 3} and let $\left(  \text{x}%
_{1}\text{,...,x}_{q},\text{x}_{q+1},...,\text{x}_{n}\right)  $ be
\textbf{Fermi coordinates} in the neighborhood of y$_{0}\in$M$_{0}$ as defined
in \textbf{chapter 1}.

\qquad\qquad\qquad\qquad\qquad\qquad\qquad\qquad\qquad\qquad\qquad\qquad
\qquad\qquad\qquad\qquad\qquad$\blacksquare$

\begin{remark}
Recall that by (ii) of \textbf{Proposition 5} above,
\end{remark}

$\qquad\qquad\nabla_{\partial_{i}}\mu_{j}=\frac{\partial\mu_{j}}{\partial
x_{i}}+\Lambda_{i}\mu_{j}$

Since $\frac{\partial\mu_{j}}{\partial x_{i}}=0,$ we have:

$\qquad\qquad\nabla_{\partial_{i}}\mu_{j}=\Lambda_{i}\mu_{j}=\Gamma_{ij}%
^{k}\mu_{k}$

where $\Lambda_{i}\left(  y_{0}\right)  =\left(  \Gamma_{ij}^{k}\left(
y_{0}\right)  \right)  j,k=1,...,q,q+1,...,n$ is the associated matrix. In the
special case of the \textbf{Levi-Civita} connection $\nabla^{0}$ on TM, the
associated matrix $\Lambda_{i}\left(  y_{0}\right)  =\left(  \Gamma_{ij}%
^{k}\left(  y_{0}\right)  \right)  j,k=1,...,q,q+1,...,n$ is similarly defined:

$\qquad\qquad\nabla_{\partial_{i}}\partial_{j}=\overset{n}{\underset{k=1}{\sum
}}\Gamma_{ij}^{k}\partial_{k}$

The matrix $\Lambda_{i}\left(  y_{0}\right)  =\left(  \Gamma_{ij}^{k}\left(
y_{0}\right)  \right)  $ is equal to zero for \textbf{all} $i,j,k=q+1,...,n.$

In particular, this will be the case if the dimension of the submanifold $q=0$
(which is equivalent to the submanifold P reducing to the singleton $\left\{
y_{0}\right\}  )$ and we have \textbf{normal coordinates.}

However $\Lambda_{i}\left(  y_{0}\right)  $ is not always equal to zero in the
case of expansions in \textbf{normal Fermi coordinates}. For example, for a,b
=1,...,q and $i=q+1,...,n,$ we have:

$\Lambda_{i}\left(  y_{0}\right)  =\Gamma_{i\text{a}}^{\text{b}}\left(
y_{0}\right)  =\Gamma_{\text{a}i}^{\text{b}}\left(  y_{0}\right)
=-\Gamma_{\text{ab}}^{i}\left(  y_{0}\right)  =-T_{\text{ab}i}\left(
y_{0}\right)  \neq0$ (except for a totally geodesic submanifold P) where $T$
is the \textbf{Second Fundamental Form} operator.

We conclude that we have $\Lambda_{i}\left(  y_{0}\right)  \neq0$ for
$i=q+1,...,n$ in the more general case of expansion in \textbf{normal Fermi
coordinates}.

Since the expansion of $\Lambda_{j}$ will be carried out in \textbf{normal
Fermi coordinates}, all derivatives with respect to \textbf{tangential Fermi
coordinates} vanish and hence,

$\frac{\partial\Lambda_{j}}{\partial x_{\text{a}}}(x)=0=$ $\frac
{\partial\Lambda_{\text{b}}}{\partial x_{\text{a}}}(x)$ for all $x\in M_{0}$
and all a,b = 1,...,q ; $j=q+1,...,n.$

\qquad\qquad\qquad\qquad\qquad\qquad\qquad\qquad\qquad\qquad\qquad\qquad
\qquad\qquad\qquad\qquad\qquad\qquad$\blacksquare$

\begin{proposition}
We have the following \textbf{Taylor} \textbf{expansion} formulae: for the
End(E)-valued connection form $\Lambda:$
\end{proposition}

\qquad(i) For a,b = 1,...,q, we have in \textbf{normal Fermi coordinates}
$i,j=q+1,...,n,$ and $x\in M_{0}:$

$\qquad\Lambda_{\text{a}}(x)=\Lambda_{\text{a}}(y_{0}%
)+[\overset{n}{\underset{i=q+1}{\text{ }\sum}}(\nabla_{i}\Lambda_{\text{a}%
})-\ \overset{\text{q}}{\underset{\text{b}=1}{\sum}}\ T_{\text{ab}i}%
)(y_{0})\Lambda_{\text{b}}(y_{0})$

$+\overset{n}{\underset{i,j=q+1}{\sum}}\ \perp_{\text{a}ij}(y_{0})\Lambda
_{j}](y_{0})x_{i}$

$+\frac{1}{2}[\left(  \frac{\partial\Omega_{j\text{a}}}{\partial x_{i}%
}+\Lambda_{i}\frac{\partial\Lambda_{\text{a}}}{\partial x_{j}}+\frac
{\partial\Lambda_{\text{a}}}{\partial x_{i}}\Lambda_{j}+\Lambda_{\text{a}%
}\frac{\partial\Lambda_{j}}{\partial x_{i}}+\Lambda_{i}\Lambda_{j}%
\Lambda_{\text{a}}\right)  $\qquad

$\qquad-R_{i\text{a}j\text{a}}\Lambda_{\text{a}}+T_{\text{aa}i}\left(
\Omega_{\text{a}j}-\Lambda_{\text{a}}\Lambda_{j}\right)  ](y_{0})x_{i}x_{j}+$
higher order terms.

(ii) For a = 1,...q, we have in \textbf{normal Fermi coordinates}
$i,j,k,l=q+1,...,n,$

and $x\in M_{0}:$

$\ \ \Lambda_{l}(x)=\Lambda_{l}(y_{0})+[\frac{\partial\Lambda_{l}}{\partial
x_{i}}+\Lambda_{l}](y_{0})x_{i}$

$+\frac{1}{2}\overset{n}{\underset{i,j=q+1}{\text{ }\sum}}\left[
\frac{\partial^{2}\Lambda_{l}}{\partial x_{i}\partial x_{j}}+\frac
{\partial\Lambda_{j}}{\partial x_{i}}\Lambda_{l}+\Lambda_{j}\frac
{\partial\Lambda_{l}}{\partial x_{i}}+\Lambda_{i}\frac{\partial\Lambda_{l}%
}{\partial x_{j}}+\Lambda_{i}\Lambda_{j}\Lambda_{l}\right]  (y_{0})x_{i}x_{j}$

$-\frac{1}{6}\overset{q}{\underset{\text{a}=1}{\sum}}%
\overset{n}{\underset{i,j,k=q+1}{\text{ }\sum}}[\left(  R_{ijk\text{a}%
}+R_{ikj\text{a}}\right)  \Lambda_{\text{a}}](y_{0})x_{i}x_{j}-\frac{1}%
{6}\overset{n}{\underset{r=q+1}{\sum}}[\left(  R_{ijlr}+R_{iljr}\right)
\Lambda_{r}](y_{0})x_{i}x_{j}$

$+\frac{1}{6}\overset{n}{\underset{i,j,k=q+1}{\text{ }\sum}}[(\nabla
_{ijk}\Lambda_{l})-\frac{1}{2}\overset{n}{\underset{r=1}{\sum}}%
\bigtriangledown_{i}(R)_{jlkr}\Lambda_{r}-\frac{1}{3}%
\overset{n}{\underset{r=1}{\sum}}\left(  R_{ijpr}+R_{ipjr}\right)
\nabla_{X_{k}}\Lambda_{r}](y_{0})x_{i}x_{j}x_{k}$

\begin{center}
$+$ higher order terms
\end{center}

where $(\nabla_{ijk}\Lambda_{l})$ is given in $\left(  9.26\right)  $ below.

The expression is take from $\left(  2.14\right)  ,$ $\left(  9.20\right)  ,$
$\left(  9.25\right)  ,$ $\left(  9.26\right)  $ and $\left(  9.29\right)  $ below.

(iii) In all \textbf{normal coordinates, }we have for
$i,j,k,l,p,r=1,...,q,q+1,...,n,$ $x\in M_{0}:$

$\ \Lambda_{l}(x)=\frac{1}{2}\Omega_{il}(y_{0})x_{i}+\frac{1}{6}\frac
{\partial\Omega_{jl}}{\partial\text{x}_{i}}(y_{0})(y_{0})x_{i}x_{j}$

$\qquad\qquad-\frac{1}{36}[\overset{n}{\underset{r=1}{\sum}}\left(
R_{ijpr}+R_{ipjr}\right)  \Omega_{kr}](y_{0})x_{i}x_{j}x$

$\qquad\qquad+\frac{1}{24}[\frac{\partial^{2}\Omega_{kl}}{\partial
x_{i}\partial x_{j}}+\Omega_{ik}(y_{0})\Omega_{jl}(y_{0})+\Omega_{ij}%
(y_{0})\Omega_{kl}(y_{0})]x_{i}x_{j}x_{k}+$ higher order terms

\begin{proof}
We recall that $\Lambda$ is the End(E)$-$valued \textbf{connection }$1-$form
and $\Omega$ the End(E)-valued \textbf{curvature} $2-$form of the vector
bundle E:
\end{proof}

\qquad Recall that $\Lambda_{j}$ is defined in $\left(  3.11\right)  $ and
$\Omega_{ij}$ is defined in $\left(  3.12\right)  $ above. We have:

$\qquad\qquad\qquad\Lambda=\overset{n}{\underset{i=1}{\sum}}\Lambda_{j}dx_{j}$
and $\Omega=\frac{1}{2!}\overset{n}{\underset{k,l=1}{\sum}}\Omega_{ij}%
(dx_{i}\wedge dx_{j}).$

and,

$\Lambda(\partial_{j})(x)=$ $\Lambda_{j}(x)\in End(E_{x});$ $\Omega
(\partial_{i},\partial_{j})(x)=\Omega_{ij}(x)\in End(E_{x}).$

For $j=1,...,q,q+1,...,n,$ we expand $\Lambda_{j}$ in \textbf{normal Fermi}
\textbf{coordinates} x$_{q+1},...,$x$_{n}$ $.$

For the expansions of $\Lambda_{\text{a}}$ for a = 1,...,q and $\Lambda_{j}$
for $j=q+1,...,n$ here we use the \textbf{notation} and \textbf{definitions}
in $\left(  9.16\right)  $ and $\left(  9.17\right)  $ of \textbf{Gray
}$\left[  25\right]  $ and combine techniques used in proving \textbf{Theorem
9.6 }in \textbf{Gray} $\left[  25\right]  $ and \textbf{Theorem 9.22 }in
\textbf{Gray} $\left[  4\right]  .$

By $\left(  9.16\right)  $ of Gray $\left[  25\right]  ,$ the \textbf{General
Expansion Formula} is given by:

$\qquad\Lambda_{\alpha}(x)=\overset{\infty}{\underset{k=0}{\sum}%
}\overset{n}{\underset{i_{1},...,i_{k}=q+1}{\sum}}\frac{1}{k!}(X_{i_{1}%
}...X_{i_{k}}\Lambda_{\alpha})(y_{0})x_{i_{1}}...x_{i_{k}}$

for $\alpha=1,...,q,q+1,...,n.$

We will use the above formula to compute the first few terms of the expansion
of $\Lambda_{\alpha}(x)$ in Fermi coordinates:

(i) The first coefficient (or zeroth order term) in the expansion of
$\Lambda_{\text{a}}(x)$ at $y_{0}$ is obviously $\Lambda_{\text{a}}(y_{0}).$

\qquad Let X$_{i}$ be a \textbf{normal Fermi field.} The second coefficient
(or the \textbf{first order term})\ is given by:

$\left(  9.0\right)  \qquad\frac{1}{1!}(X_{i}\Lambda_{\text{a}})(y_{0}%
)=(\nabla_{i}\Lambda_{\text{a}})(y_{0})+$\ \ $\overset{n}{\underset{s=1}{\sum
}}<\nabla_{X_{i}}X_{\text{a}},X_{s}>(y_{0})\Lambda_{s}(y_{0})$

\qquad\qquad\qquad\qquad\qquad$=\frac{\partial\Lambda_{\text{a}}}{\partial
x_{i}}(y_{0})+\Lambda_{i}\Lambda_{\text{a}}+$%
\ \ $\overset{n}{\underset{s=1}{\sum}}<\nabla_{X_{i}}X_{\text{a}},X_{s}%
>(y_{0})\Lambda_{s}(y_{0})$

We have from (v) of \textbf{Proposition 6} above,

\qquad\qquad$\Omega_{i\text{a}}=\frac{\partial\Lambda_{\text{a}}}{\partial
x_{i}}-\frac{\partial\Lambda_{i}}{\partial x_{\text{a}}}+\Lambda_{i}%
\Lambda_{\text{a}}-\Lambda_{\text{a}}\Lambda_{i}$

\qquad Since all expansions are made in normal Fermi coordinates and so,
$\frac{\partial\Lambda_{i}}{\partial x_{\text{a}}}=0.$ We have:

\qquad$\left(  9.1\right)  \qquad\left(  X_{i}\Lambda_{\text{a}}\right)
(y_{0})=(\nabla_{i}\Lambda_{\text{a}})(y_{0})=\frac{\partial\Lambda_{\text{a}%
}}{\partial x_{i}}(y_{0})+\Lambda_{i}\Lambda_{\text{a}}=\ \Omega_{\text{a}%
i}+\Lambda_{\text{a}}\Lambda_{i}$

By \textbf{Lemma 4} above to have:

$\qquad\left(  \nabla_{X_{i}}X_{\text{a}}\right)  (y_{0})=\left(
\nabla_{X_{\text{a}}}X_{i}\right)  (y_{0})=\left(  T_{X_{\text{a}}}%
X_{i}\ \right)  (y_{0})+$ $\left(  \perp_{X_{\text{a}}}X_{i}\right)  (y_{0})$

Consequently,

$\qquad\ \overset{n}{\underset{s=1}{\sum}}<\nabla_{X_{i}}X_{\text{a}}%
,X_{s}>(y_{0})\Lambda_{s}(y_{0})$

$\qquad=$\ $\overset{n}{\underset{s=1}{\sum}}$\ $<T_{X_{\text{a}}}X_{i}%
,X_{s}>(y_{0})\Lambda_{s}(y_{0})+$\ $\overset{n}{\underset{s=1}{\sum}}%
$\ $<\perp_{X_{\text{a}}}X_{i},X_{s}>(y_{0})\Lambda_{s}(y_{0})$\ 

Since $T_{X_{\text{a}}}X_{i}$ is tangential and $\perp_{X_{\text{a}}}X_{i}$ is
normal, we have:

$\left(  9.2\right)  \ \overset{n}{\underset{s=1}{\sum}}<\nabla_{X_{i}%
}X_{\text{a}},X_{s}>(y_{0})\Lambda_{s}(y_{0})=$\ $\overset{\text{q}%
}{\underset{\text{b}=1}{\sum}}$\ $<T_{X_{\text{a}}}X_{i},X_{\text{b}}%
>(y_{0})\Lambda_{\text{b}}(y_{0})$

$\qquad+$\ $\overset{n}{\underset{j=q+1}{\sum}}$\ $<\perp_{X_{\text{a}}}%
X_{i},X_{j}>(y_{0})\Lambda_{j}(y_{0})$

By (vi) and (xiv) of the \textbf{Notation} at the beginning of this
Chapterthere, we have:

$\left(  9.3\right)  \qquad\overset{n}{\underset{s=1}{\sum}}<\nabla_{X_{i}%
}X_{\text{a}},X_{s}>(y_{0})\Lambda_{s}(y_{0})=-$\ $\overset{\text{q}%
}{\underset{\text{b}=1}{\sum}}$\ $T_{\text{ab}i}(y_{0})\Lambda_{\text{b}%
}(y_{0})$

$\qquad\qquad+$\ $\overset{n}{\underset{j=q+1}{\sum}}$\ $\perp_{\text{a}%
ij}(y_{0})\Lambda_{j}(y_{0})$

By $\left(  9.0\right)  ,\left(  9.1\right)  ,\left(  9.2\right)  $ and
$\left(  9.3\right)  $ the \textbf{first order term} in the expansion of
$\Lambda_{\text{a}}$ is given by:

$\left(  9.4\right)  \qquad\left(  X_{i}\Lambda_{\text{a}}\right)
(y_{0})=\ \Omega_{\text{a}i}+\Lambda_{\text{a}}\Lambda_{i}-$%
\ $\overset{\text{q}}{\underset{\text{b}=1}{\sum}}$\ $T_{\text{ab}i}%
(y_{0})\Lambda_{\text{b}}(y_{0})+$\ $\overset{n}{\underset{j=q+1}{\sum}}%
$\ $\perp_{\text{a}ij}(y_{0})\Lambda_{j}(y_{0})\qquad\qquad\qquad\qquad\qquad$

\qquad The \textbf{second order term} in the expansion of $\Lambda_{\text{a}%
}(x)$ in \textbf{normal Fermi coordinates }is given by:$\qquad\qquad$

$\left(  9.5\right)  \qquad\frac{1}{2!}\left(  X_{i}X_{j}\Lambda_{\text{a}%
}\right)  (y_{0})=\frac{1}{2}(\nabla_{X_{i}}\nabla_{X_{j}}\Lambda_{\text{a}%
})(y_{0})+\frac{1}{2}$ $<\nabla_{X_{i}X_{j}}X_{\text{a}},X_{\text{a}}%
>(y_{0})\Lambda_{\text{a}}(y_{0})$

$\qquad\qquad\qquad\qquad\qquad\qquad\qquad+\frac{1}{2}<\nabla_{X_{i}%
}X_{\text{a}},X_{\text{a}}>(y_{0})\nabla_{X_{j}}\Lambda_{\text{a}}(y_{0})$

By (iii) of \textbf{Proposition }$6$ above, we have for a = 1,...,q and
$i,j=q+1,...,n:$

$\left(  9.6\right)  \qquad\left(  \nabla_{X_{i}}\nabla_{X_{j}}\Lambda
_{\text{a}}\right)  =\left(  \nabla_{\partial_{i}}\nabla_{\partial_{j}}%
\Lambda_{\text{a}}\right)  $

$\qquad\qquad=\frac{\partial^{2}\Lambda_{\text{a}}}{\partial x_{i}\partial
x_{j}}+\frac{\partial\Lambda_{j}}{\partial x_{i}}\Lambda_{\text{a}}%
+\Lambda_{j}\frac{\partial\Lambda_{\text{a}}}{\partial x_{i}}+\Lambda_{i}%
\frac{\partial\Lambda_{\text{a}}}{\partial x_{j}}+\Lambda_{i}\Lambda
_{j}\Lambda_{\text{a}}$

By (viii) of \textbf{Proposition }$6,$ we have:

$\left(  9.7\right)  $\qquad$\ \frac{\partial^{2}\Lambda_{\text{a}}}{\partial
x_{i}\partial x_{j}}=\frac{\partial\Omega_{j\text{a}}}{\partial x_{i}}%
+\frac{\partial\Lambda_{\text{a}}}{\partial x_{i}}\Lambda_{j}-\Lambda_{j}%
\frac{\partial\Lambda_{\text{a}}}{\partial x_{i}}+\Lambda_{\text{a}}%
\frac{\partial\Lambda_{j}}{\partial x_{i}}-\frac{\partial\Lambda_{j}}{\partial
x_{i}}\Lambda_{\text{a}}$

In $\left(  9.6\right)  $ we replace $\frac{\partial^{2}\Lambda_{\text{a}}%
}{\partial x_{i}\partial x_{j}}$ by the Right Hand Side of $\left(
9.7\right)  $ and have:

$\left(  9.8\right)  \qquad\left(  \nabla_{X_{i}}\nabla_{X_{j}}\Lambda
_{\text{a}}\right)  =\frac{\partial\Omega_{j\text{a}}}{\partial x_{i}}%
+\Lambda_{i}\frac{\partial\Lambda_{\text{a}}}{\partial x_{j}}+\frac
{\partial\Lambda_{\text{a}}}{\partial x_{i}}\Lambda_{j}+\Lambda_{\text{a}%
}\frac{\partial\Lambda_{j}}{\partial x_{i}}+\Lambda_{i}\Lambda_{j}%
\Lambda_{\text{a}}\qquad$

Next we have by $\left(  9.51\right)  $ of \textbf{Gray }$\left[  25\right]
,$

$\left(  9.9\right)  \qquad<\nabla_{X_{i}X_{j}}X_{\text{a}},X_{\text{a}%
}>(y_{0})\Lambda_{\text{a}}(y_{0})=-\left(  R_{i\text{a}j\text{a}}\right)
(y_{0})\Lambda_{\text{a}}(y_{0})$

We compute the last term $<\nabla_{X_{i}}X_{\text{a}},X_{\text{a}}%
>(y_{0})\nabla_{X_{j}}\Lambda_{\text{a}}(y_{0}):$

By \textbf{Lemma 4} here, $\nabla_{X_{i}}X_{\text{a}}=\nabla_{X_{\text{a}}%
}X_{i}=T_{X_{\text{a}}}X_{i}\ +\perp_{X_{\text{a}}}X_{i}$

Therefore,

$<\nabla_{X_{i}}X_{\text{a}},X_{\text{a}}>(y_{0})\nabla_{X_{j}}\Lambda
_{\text{a}}(y_{0})=$ $\overset{n}{\underset{s=1}{\text{ }\sum}}<$%
T$_{X_{\text{a}}}X_{i}$ \ + $\perp_{X_{\text{a}}}X_{i},X_{\text{a}}%
>(y_{0})\nabla_{X_{j}}\Lambda_{\text{a}}(y_{0})$

$=$ $<T_{X_{\text{a}}}X_{i}\ +\perp_{X_{\text{a}}}X_{i},X_{\text{a}}%
>(y_{0})\left[  \frac{\partial\Lambda_{\text{a}}}{\partial x_{j}}+\Lambda
_{j}\Lambda_{\text{a}}\right]  (y_{0})$

Since $T_{X_{\text{a}}}X_{i}$ and $X_{\text{a}}$ are tangential and
$\perp_{X_{\text{a}}}X_{i}$ and is normal, we have:

$<\perp_{X_{\text{a}}}X_{i},X_{\text{a}}>(y_{0})=0$ and so,

$<\nabla_{X_{i}}X_{\text{a}},X_{\text{a}}>(y_{0})\nabla_{X_{j}}\Lambda
_{\text{a}}(y_{0})$

$=$ $<T_{X_{\text{a}}}X_{i}\ ,X_{\text{a}}>(y_{0})\left[  \frac{\partial
\Lambda_{\text{a}}}{\partial x_{j}}+\Lambda_{j}\Lambda_{\text{a}}\right]
(y_{0})=-T_{\text{aa}i}(y_{0})\left[  \frac{\partial\Lambda_{\text{a}}%
}{\partial x_{j}}+\Lambda_{j}\Lambda_{\text{a}}\right]  (y_{0})$

By (v) of \textbf{Proposition }$6$ above,

$\ \frac{\partial\Lambda_{\text{a}}}{\partial x_{j}}=-\Omega_{\text{a}%
j}+[\Lambda_{\text{a}},\Lambda_{j}]=-\Omega_{\text{a}j}+\Lambda_{\text{a}%
}\Lambda_{j}-\Lambda_{j}\Lambda_{\text{a}}$

Consequently,

$\left(  9.10\right)  \qquad<\nabla_{X_{i}}X_{\text{a}},X_{\text{a}}%
>(y_{0})\nabla_{X_{j}}\Lambda_{\text{a}}(y_{0})=-T_{\text{aa}i}(y_{0})\left[
-\Omega_{\text{a}j}+\Lambda_{\text{a}}\Lambda_{j}\right]  (y_{0})$

$\qquad\qquad=T_{\text{aa}i}(y_{0})\left[  \Omega_{\text{a}j}-\Lambda
_{\text{a}}\Lambda_{j}\right]  (y_{0})$

Using $\left(  9.5\right)  ,\left(  8.8\right)  $ and $\left(  9.10\right)  ,$
we see that the final expression for the \textbf{second order term} in
$\left(  9.5\right)  $ is given by:

$\left(  9.11\right)  \qquad\frac{1}{2!}\left(  X_{i}X_{j}\Lambda_{\text{a}%
}\right)  (y_{0})=\frac{1}{2}[\left(  \frac{\partial\Omega_{j\text{a}}%
}{\partial x_{i}}+\Lambda_{i}\frac{\partial\Lambda_{\text{a}}}{\partial x_{j}%
}+\frac{\partial\Lambda_{\text{a}}}{\partial x_{i}}\Lambda_{j}+\Lambda
_{\text{a}}\frac{\partial\Lambda_{j}}{\partial x_{i}}+\Lambda_{i}\Lambda
_{j}\Lambda_{\text{a}}\right)  $\qquad

$\qquad\qquad\qquad\qquad\qquad\qquad\qquad\qquad-R_{i\text{a}j\text{a}%
}\Lambda_{\text{a}}+T_{\text{aa}i}\left(  \Omega_{\text{a}j}-\Lambda
_{\text{a}}\Lambda_{j}\right)  ](y_{0})$

(ii) We follow the same procedure as in (i):

Repalacing W$_{\alpha_{1}....\alpha_{r}}$ defined in $\left(  9.16\right)  $
of \textbf{Gray }$\left[  25\right]  $ by $\Lambda_{\alpha},$ we see that for
a general point $x\in M_{0}$ and setting:

$\left(  9.12\right)  $\qquad$\qquad\qquad\Lambda_{\alpha}(x)=\overset{\infty
}{\underset{k=0}{\sum}}\overset{n}{\underset{i_{1},...,i_{k}=1}{\sum}}\frac
{1}{k!}(X_{i_{1}}...X_{i_{k}}\Lambda_{\alpha})(y_{0})x_{i_{1}}...x_{i_{k}}$

where the Right Hand Side is defined in $\left(  9.17\right)  $ of
\textbf{Gray }$\left[  25\right]  .$

We follow the techiques used in \textbf{Theorem 9.6} and \textbf{Theorem 9.22
}in \textbf{Gray} $\left[  25\right]  $ adapted to our case here where
$\Lambda$ is an End(E)$-$valued $1-$form:

We take $k=0$ and the first coefficient in the expansion of $\Lambda_{k}(x)$
in normal Fermi coordinates centred at $y_{0}$ is obviously $\Lambda_{k}%
(y_{0}).$

We will repeatedly use the definition given in $\left(  9.16\right)  $ of
\textbf{Gray} $\left[  25\right]  $ to compute the coefficients. We note that
the above expansion formula is a form of the usual \textbf{Taylor Expansion}
where the terms of the expansion in $\left(  9.16\right)  $ of \textbf{Gray
}$\left[  25\right]  $ are defined in $\left(  9.17\right)  $ of \textbf{Gray
}$\left[  25\right]  .$

The \textbf{first order term} is given for $i,l=q+1,...,n:$

$\left(  9.13\right)  \qquad\frac{1}{1!}\left(  X_{i}\Lambda_{l}\right)
(y_{0})=\frac{1}{1!}\left(  \nabla_{_{X_{i}}}\Lambda_{l}\right)  (y_{0}%
)+\frac{1}{1!}$\ $\overset{n}{\underset{s=1}{\sum}}<\nabla_{X_{i}}X_{l}%
,X_{s}>(y_{0})\Lambda_{s}(y_{0})$

By $\left(  9.3\right)  $ of \textbf{Lemma 9.1} of \textbf{Gray} $\left[
25\right]  $, we have: $\left(  \nabla_{X_{i}}X_{l}\right)  (y_{0})=0$

for $i,l=q+1,...,n.$

The second term (the \textbf{first order term}) is thus given by:

$\left(  9.14\right)  \qquad\frac{1}{1!}\left(  X_{i}\Lambda_{l}\right)
(y_{0})=\left(  \nabla_{_{X_{i}}}\Lambda_{l}\right)  (y_{0})=\left(
\nabla_{_{\partial_{i}}}\Lambda_{l}\right)  (y_{0})$

$\qquad\qquad=\nabla_{i}\Lambda_{l}(y_{0})=[\frac{\partial\Lambda_{l}%
}{\partial x_{i}}+\Lambda_{i}\Lambda_{l}](y_{0})$

\qquad\qquad\qquad\qquad\qquad\qquad\qquad\qquad\qquad\qquad\qquad\qquad
\qquad$\blacksquare$

We follow the method in \textbf{Theorem 9.22 }in \textbf{Gray} $\left[
25\right]  $ and compute the \textbf{second order term}

in the expansion of $\Lambda_{k}(x)$ given by.

We have for $i,j,l=q+1,...,n:$

$\left(  9.15\right)  \qquad\frac{1}{2!}\left(  X_{i}X_{j}\Lambda_{l}\right)
(y_{0})=\frac{1}{2!}(\nabla_{X_{i}}\nabla_{X_{j}}\Lambda_{l})(y_{0})$

$\qquad\qquad+\frac{1}{2!}$ $\overset{n}{\underset{r=1}{\sum}}$ $<\nabla
_{X_{i}X_{j}}^{2}X_{l},X_{r}>(y_{0})\Lambda_{r}(y_{0})$

$\qquad\qquad+\frac{1}{2!}.2$ $\overset{n}{\underset{r=1}{\sum}}$
$<\nabla_{X_{i}}X_{j},X_{r}>(y_{0})\nabla_{X_{l}}\Lambda_{r}(y_{0})$\ \ \ \ \ \ \ \ \ \ \ \ \ \ \ \ \ \ 

We compute each term on the RHS of $\left(  9.15\right)  :$

By (iii) of \textbf{Proposition }$6$ above, we have for $i,j,k,l=q+1,...,n:$

$\left(  9.16\right)  \qquad(\nabla_{X_{i}}\nabla_{X_{j}}\Lambda_{l}%
)(y_{0})=(\nabla_{\partial_{i}}\nabla_{\partial_{j}}\Lambda_{l})(y_{0})$

$\qquad\qquad=[\frac{\partial^{2}\Lambda_{l}}{\partial x_{i}\partial x_{j}%
}+\frac{\partial\Lambda_{j}}{\partial x_{i}}\Lambda_{l}+\Lambda_{j}%
\frac{\partial\Lambda_{l}}{\partial x_{i}}+\Lambda_{i}\frac{\partial
\Lambda_{l}}{\partial x_{j}}+\Lambda_{i}\Lambda_{j}\Lambda_{l}](y_{0})$

We next compute:

$\qquad\overset{n}{\underset{r=1}{\sum}}$ $<\nabla_{X_{i}X_{j}}^{2}X_{l}%
,X_{r}>(y_{0})\Lambda_{r}(y_{0})=$ $\overset{q}{\underset{\text{a}=1}{\sum}}$
$<\nabla_{X_{i}X_{j}}^{2}X_{l},X_{\text{a}}>(y_{0})\Lambda_{\text{a}}(y_{0})$

$\qquad+\overset{n}{\underset{r=q+1}{\sum}}$ $<\nabla_{X_{i}X_{j}}^{2}%
X_{l},X_{r}>(y_{0})\Lambda_{r}(y_{0})$

By $\left(  9.11\right)  $ of \textbf{Gray} $\left[  26\right]  ,$

$\left(  9.17\right)  \qquad\overset{q}{\underset{\text{a}=1}{\sum}}%
<\nabla_{X_{i}X_{j}}^{2}X_{l},X_{\text{a}}>(y_{0})\Lambda_{\text{a}}(y_{0})$

$\qquad\qquad=-\frac{1}{3}\overset{q}{\underset{\text{a}=1}{\sum}}\left(
<R_{X_{i}X_{j}}X_{l}+R_{X_{i}X_{l}}X_{j},X_{\text{a}}>\right)  (y_{0}%
)\Lambda_{\text{a}}(y_{0})$

$\qquad\qquad=-\frac{1}{3}\overset{q}{\underset{\text{a}=1}{\sum}}\left(
R_{ijl\text{a}}+R_{ilj\text{a}}\right)  (y_{0})\Lambda_{\text{a}}(y_{0})$

Next we have:

$\qquad\overset{n}{\underset{r=q+1}{\sum}}$ $<\nabla_{X_{i}X_{j}}^{2}%
X_{l},X_{r}>(y_{0})\Lambda_{r}(y_{0})$

$\qquad=-\frac{1}{3}\overset{n}{\underset{r=q+1}{\sum}}\left(  <R_{X_{i}X_{j}%
}X_{l}+R_{X_{i}X_{l}}X_{j},X_{r}>\right)  \Lambda_{r}(y_{0})$

$\qquad=-\frac{1}{3}\overset{n}{\underset{r=q+1}{\sum}}\left(  R_{ijlr}%
+R_{iljr}\right)  \Lambda_{r}(y_{0})$

We conclude that,

$\left(  9.18\right)  \qquad\overset{n}{\underset{r=1}{\sum}}$ $<\nabla
_{X_{i}X_{j}}X_{l},X_{r}>(y_{0})\Lambda_{r}(y_{0})$

$=-\frac{1}{3}\overset{q}{\underset{\text{a}=1}{\sum}}\left(  R_{ijl\text{a}%
}+R_{ilj\text{a}}\right)  (y_{0})\Lambda_{\text{a}}(y_{0})-\frac{1}%
{3}\overset{n}{\underset{r=q+1}{\sum}}\left(  R_{ijlr}+R_{iljr}\right)
\Lambda_{r}(y_{0})$

To compute the last expression here, we note that $\nabla_{X_{i}}X_{k}(y_{0})$
$=0$ for $i,k=q+1,...,n,$ and so we have $<\nabla_{X_{i}}X_{k},X_{s}%
>(y_{0})=0$ and so the (the determinant) disappears.

We then conclude by $\left(  9.15\right)  ,\left(  9.16\right)  ,\left(
9.17\right)  $ and $\left(  9.18\right)  $ that the \textbf{second order term}
is given by:

$\left(  9,19\right)  \qquad\frac{1}{2!}\left(  X_{i}X_{j}\Lambda_{l}\right)
(y_{0})=\frac{1}{2}\left[  \frac{\partial^{2}\Lambda_{l}}{\partial
x_{i}\partial x_{j}}+\frac{\partial\Lambda_{j}}{\partial x_{i}}\Lambda
_{l}+\Lambda_{j}\frac{\partial\Lambda_{l}}{\partial x_{i}}+\Lambda_{i}%
\frac{\partial\Lambda_{l}}{\partial x_{j}}+\Lambda_{i}\Lambda_{j}\Lambda
_{l}\right]  (y_{0})$

$\qquad\qquad-\frac{1}{6}\overset{q}{\underset{\text{a}=1}{\sum}}[\left(
R_{ijl\text{a}}+R_{ilj\text{a}}\right)  \Lambda_{\text{a}}](y_{0})-\frac{1}%
{6}\overset{n}{\underset{r=q+1}{\sum}}[\left(  R_{ijlr}+R_{iljr}\right)
\Lambda_{r}](y_{0})$

By (viii) of \textbf{Proposition 5}, we have \textbf{normal coordinates},
$\qquad\qquad\qquad\qquad\qquad\qquad$

$\ \qquad\qquad\qquad\qquad\ \ \frac{\partial^{2}\Lambda_{l}}{\partial
\text{x}_{i}\partial\text{x}_{j}}(y_{0})=\frac{1}{3}\frac{\partial\Omega_{il}%
}{\partial\text{x}_{j}}(y_{0})$

and,

$\qquad\qquad\qquad\qquad\Lambda_{l}(y_{0})=0$

Consequently, $\left(  9,19\right)  $ becomes,

$\left(  9.20\right)  \qquad\frac{1}{2!}\left(  X_{i}X_{j}\Lambda_{l}\right)
(y_{0})=\frac{1}{6}\frac{\partial\Omega_{il}}{\partial\text{x}_{j}}%
(y_{0})(y_{0})$

We next compute the \textbf{third order term}: For $i,j,k,l=q+1,...,n,$ we have:

$\left(  9.21\right)  \qquad\frac{1}{3!}\left(  X_{i}X_{j}X_{k}\Lambda
_{l}\right)  (y_{0})=\frac{1}{3!}(\nabla_{X_{i}}\nabla_{X_{j}}\nabla_{X_{k}%
}\Lambda_{l})(y_{0})$

$\qquad\qquad+\frac{1}{3!}$ $\overset{n}{\underset{r=1}{\sum}}$ $<\nabla
_{X_{i}X_{j}X_{k}}^{3}X_{l},X_{r}>(y_{0})\Lambda_{r}(y_{0})$

$\qquad\qquad+\frac{1}{3!}$ $\overset{n}{\underset{r=1}{\sum}}$ $<\nabla
_{X_{i}X_{j}}^{2}X_{l},X_{r}>(y_{0})\nabla_{X_{k}}\Lambda_{r}(y_{0})$

$\qquad\qquad+\frac{1}{3!}$ $\overset{n}{\underset{r=1}{\sum}}$ $<\nabla
_{X_{i}}X_{l},X_{r}>(y_{0})\nabla_{X_{j}X_{k}}^{2}\Lambda_{r}(y_{0})$

We compute each terms above expression:

First, we set:

$\left(  9.22\right)  \qquad\qquad\qquad\qquad\frac{1}{3!}(\nabla_{X_{i}%
}\nabla_{X_{j}}\nabla_{X_{k}}\Lambda_{l})(y_{0})=\frac{1}{6}(\nabla
_{ijk}\Lambda_{l})(y_{0})$

Details of computation of this term will be given below.

Then next, we have by $\left(  9.12\right)  $ of \textbf{Lemma }$\left(
9.3\right)  $ of \textbf{Gray }$\left[  25\right]  ,$ we have:

$\qquad<\nabla_{X_{i}X_{j}X_{k}}^{3}X_{l},X_{r}>(y_{0})\Lambda_{r}(y_{0})=$
$\ <-\frac{1}{2}(\bigtriangledown_{X_{i}}(R)_{X_{j}X_{l}}X_{k},X_{r}%
>)(y_{0})\Lambda_{r}(y_{0})$

$\qquad<\nabla_{X_{i}X_{j}X_{k}}^{3}X_{l},X_{r}>(y_{0})\Lambda_{r}%
(y_{0})=-\frac{1}{2}\left(  \bigtriangledown_{X_{i}}(R)_{X_{j}X_{l}X_{k}X_{r}%
}\right)  (y_{0})\Lambda_{r}(y_{0})$

$\qquad\qquad\qquad\qquad\qquad\qquad\qquad\qquad\ \ =-\frac{1}{2}%
(\bigtriangledown_{i}(R)_{jlkr})(y_{0})\Lambda_{r}(y_{0})$

Since in \textbf{normal coordinates}, $\Lambda_{r}(y_{0})=0$ for
$r=1,...,n,$we have:

$\left(  9.23\right)  \qquad<\nabla_{X_{i}X_{j}X_{k}}^{3}X_{l},X_{r}%
>(y_{0})\Lambda_{r}(y_{0})=0$

Next we have:

$\qquad<\nabla_{X_{i}X_{j}}^{2}X_{l},X_{r}>(y_{0})\nabla_{X_{k}}\Lambda
_{r}(y_{0})=-\frac{1}{3}\left(  R_{X_{i}X_{j}}X_{l}+R_{X_{i}X_{l}}X_{j}%
,X_{r}\right)  (y_{0})\nabla_{X_{k}}\Lambda_{r}(y_{0})$

$=-\frac{1}{3}\left(  R_{X_{i}X_{j}X_{p}X_{r}}+R_{X_{i}X_{l}X_{j}X_{r}%
}\right)  (y_{0})\nabla_{X_{k}}\Lambda_{r}(y_{0})=-\frac{1}{3}\left(
R_{ijlr}+R_{iljr}\right)  (y_{0})\nabla_{X_{k}}\Lambda_{r}(y_{0}).$Consequently,

$\left(  9.24\right)  \qquad<\nabla_{X_{i}X_{j}}^{2}X_{l},X_{r}>(y_{0}%
)\nabla_{X_{k}}\Lambda_{r}(y_{0})=-\frac{1}{3}\left(  R_{ijlr}+R_{iljr}%
\right)  (y_{0})\nabla_{X_{k}}\Lambda_{r}(y_{0})$

\qquad Since $(\nabla_{X_{i}}X_{l})(y_{0})=0$ for $i,l=q+1,...,n,$ the last
term in $\left(  9.21\right)  $ disappears and consequently, we have by
$\left(  9.21\right)  ,$ $\left(  9.22\right)  ,$ $\left(  9.23\right)  $ and
$\left(  9.24\right)  $ that the expression for the \textbf{third order term}
is given by:

$\left(  9.25\right)  \qquad\frac{1}{3!}\left(  X_{i}X_{j}X_{k}\Lambda
_{l}\right)  (y_{0})=\frac{1}{6}[\nabla_{ijk}\Lambda_{l}-\frac{1}%
{3}\overset{n}{\underset{r=1}{\sum}}\left(  R_{ijpr}+R_{ipjr}\right)
\nabla_{X_{k}}\Lambda_{r}](y_{0})$

We give it more detaily here: Since $\nabla_{i}=\frac{\partial}{\partial
x_{i}}+\Lambda_{i}$

$\qquad\qquad\nabla_{ijk}\Lambda_{l}=\nabla_{X_{i}}\nabla_{X_{j}}\nabla
_{X_{k}}\Lambda_{l}$

$\qquad\qquad\ \nabla_{X_{j}}\nabla_{X_{k}}\Lambda_{l}=\frac{\partial
^{2}\Lambda_{l}}{\partial x_{j}\partial x_{k}}+\frac{\partial\Lambda_{k}%
}{\partial x_{j}}\Lambda_{l}+\Lambda_{k}\frac{\partial\Lambda_{l}}{\partial
x_{j}}+\Lambda_{j}\frac{\partial\Lambda_{l}}{\partial x_{k}}+\Lambda
_{j}\Lambda_{k}\Lambda_{l}$

$\nabla_{ijk}=\nabla_{X_{i}}\nabla_{X_{j}}\nabla_{X_{k}}=\frac{\partial
}{\partial x_{i}}\left(  \frac{\partial^{2}}{\partial x_{j}\partial x_{k}%
}+\frac{\partial\Lambda_{k}}{\partial x_{j}}+\Lambda_{k}\frac{\partial
}{\partial x_{j}}+\Lambda_{j}\frac{\partial}{\partial x_{k}}+\Lambda
_{j}\Lambda_{k}\right)  $

$\qquad\qquad\qquad\qquad+\Lambda_{i}\left(  \frac{\partial^{2}}{\partial
x_{j}\partial x_{k}}+\frac{\partial\Lambda_{k}}{\partial x_{j}}+\Lambda
_{k}\frac{\partial}{\partial x_{j}}+\Lambda_{j}\frac{\partial}{\partial x_{k}%
}+\Lambda_{j}\Lambda_{k}\right)  $

Therefore,

$\left(  9.26\right)  \qquad\nabla_{ijk}\Lambda_{l}=\frac{\partial^{3}%
\Lambda_{l}}{\partial x_{i}\partial x_{j}\partial x_{k}}+\frac{\partial
^{2}\Lambda_{k}}{\partial x_{i}\partial x_{j}}\Lambda_{l}+\frac{\partial
\Lambda_{k}}{\partial x_{i}}\frac{\partial\Lambda_{l}}{\partial x_{j}}%
+\Lambda_{k}\frac{\partial^{2}\Lambda_{l}}{\partial x_{i}\partial x_{j}}%
+\frac{\partial\Lambda_{j}}{\partial x_{i}}\frac{\partial\Lambda_{l}}{\partial
x_{k}}+\Lambda_{j}\frac{\partial^{2}\Lambda_{l}}{\partial x_{i}\partial x_{k}%
}$

$\qquad+\frac{\partial\Lambda_{j}}{\partial x_{i}}\Lambda_{k}\Lambda
_{l}+\Lambda_{j}\frac{\partial\Lambda_{k}}{\partial x_{j}}\Lambda_{l}%
+\Lambda_{i}\left(  \frac{\partial^{2}\Lambda_{l}}{\partial x_{j}\partial
x_{k}}+\frac{\partial\Lambda_{k}}{\partial x_{j}}\Lambda_{l}+\Lambda_{k}%
\frac{\partial\Lambda_{l}}{\partial x_{j}}+\Lambda_{j}\frac{\partial
\Lambda_{l}}{\partial x_{k}}+\Lambda_{j}\Lambda_{k}\Lambda_{l}\right)  $

(iii)\qquad The expansion of $\Lambda_{l}(x)$ in Fermi coordinates, up to the
third order term, is thus given from by:

\qquad$\ \ \Lambda_{l}(x)=\Lambda_{l}(y_{0})+[\frac{\partial\Lambda_{l}%
}{\partial x_{i}}+\Lambda_{l}](y_{0})x_{i}$

$+\frac{1}{2}\overset{n}{\underset{i,j=q+1}{\text{ }\sum}}\left[
\frac{\partial^{2}\Lambda_{l}}{\partial x_{i}\partial x_{j}}+\frac
{\partial\Lambda_{j}}{\partial x_{i}}\Lambda_{l}+\Lambda_{j}\frac
{\partial\Lambda_{l}}{\partial x_{i}}+\Lambda_{i}\frac{\partial\Lambda_{l}%
}{\partial x_{j}}+\Lambda_{i}\Lambda_{j}\Lambda_{l}\right]  (y_{0})x_{i}x_{j}$

$-\frac{1}{6}\overset{q}{\underset{\text{a}=1}{\sum}}%
\overset{n}{\underset{i,j,k=q+1}{\text{ }\sum}}[\left(  R_{ijk\text{a}%
}+R_{ikj\text{a}}\right)  \Lambda_{\text{a}}](y_{0})x_{i}x_{j}-\frac{1}%
{6}\overset{n}{\underset{r=q+1}{\sum}}[\left(  R_{ijlr}+R_{iljr}\right)
\Lambda_{r}](y_{0})x_{i}x_{j}$

$+\frac{1}{6}\overset{n}{\underset{i,j,k=q+1}{\text{ }\sum}}[(\nabla
_{ijk}\Lambda_{l})-\frac{1}{2}\overset{n}{\underset{r=1}{\sum}}%
\bigtriangledown_{i}R_{jlkr}\Lambda_{r}-\frac{1}{3}%
\overset{n}{\underset{r=1}{\sum}}\left(  R_{ijpr}+R_{ipjr}\right)
\nabla_{X_{k}}\Lambda_{r}](y_{0})x_{i}x_{j}x_{k}$

where $(\nabla_{ijk}\Lambda_{l})$ is given in $\left(  9.26\right)  $
\textbf{above} and finally in $\left(  9.29\right)  $ \textbf{below}.

\ In \textbf{normal coordinates}, we have $q=0$ and $\Lambda_{j}%
(y_{0})=0=\Lambda_{r}(y_{0})$ and so the expansion above becomes:

\qquad$\ \ \Lambda_{l}(x)=+\frac{\partial\Lambda_{l}}{\partial x_{i}}%
(y_{0})x_{i}+\frac{1}{2}\overset{n}{\underset{i,j=1}{\text{ }\sum}}\left[
\frac{\partial^{2}\Lambda_{l}}{\partial x_{i}\partial x_{j}}\right]
(y_{0})x_{i}x_{j}$

$\qquad+\frac{1}{6}\overset{n}{\underset{i,j,k=1}{\text{ }\sum}}[(\nabla
_{ijk}\Lambda_{l})-\frac{1}{3}\overset{n}{\underset{r=1}{\sum}}\left(
R_{ijpr}+R_{ipjr}\right)  \nabla_{X_{k}}\Lambda_{r}](y_{0})x_{i}x_{j}x_{k}$

\qquad From $\left(  2.14\right)  $ the first order term is given by
$[\nabla_{i}\Lambda_{l}](y_{0})$ and by (viii) of \textbf{Proposition 4
}above, we have:

\qquad Since $\Lambda_{j}(y_{0})=0,$ we have:

$\left(  9.27\right)  \qquad\lbrack\nabla_{i}\Lambda_{l}](y_{0})=\frac
{\partial\Lambda_{l}}{\partial x_{i}}(y_{0})+\Lambda_{i}(y_{0})\Lambda
_{l}(y_{0})=\frac{\partial\Lambda_{l}}{\partial x_{i}}(y_{0})=\frac{1}%
{2}\Omega_{il}(y_{0})$

\qquad$\qquad\nabla_{X_{k}}\Lambda_{r}](y_{0})=[\nabla_{k}\Lambda_{r}%
](y_{0})=\frac{1}{2}\Omega_{kr}(y_{0})$

\qquad By \textbf{Proposition 1.18 }of \textbf{Berline, Getzeler and Vergne
}$\left[  7\right]  :$

$\left(  9.28\right)  \qquad\frac{1}{1!}\left(  X_{i}\Lambda_{l}\right)
(y_{0})=\left(  \nabla_{_{X_{i}}}\Lambda_{l}\right)  (y_{0})=\left(
\nabla_{_{\partial_{i}}}\Lambda_{l}\right)  (y_{0})=\nabla_{i}\Lambda
_{l}(y_{0})=\frac{\partial\Lambda_{l}}{\partial x_{i}}(y_{0})=\frac{1}%
{2}\Omega_{il}(y_{0})$

$\qquad\qquad\ \ \frac{\partial^{2}\Lambda_{l}}{\partial\text{x}_{i}%
\partial\text{x}_{k}}(y_{0})=\frac{1}{3}\frac{\partial\Omega_{kl}}%
{\partial\text{x}_{i}}(y_{0})$ by (viii) of \textbf{Proposition 5} above.

$\qquad\qquad\frac{\partial^{3}\Lambda_{l}}{\partial x_{i}\partial
x_{j}\partial x_{k}}(y_{0})=\frac{1}{4}\frac{\partial^{2}\Omega_{kl}}{\partial
x_{i}\partial x_{j}}(y_{0})$ by using the formula in \textbf{Proposition 1.18}
of \textbf{Berline, Getzler, Vergne }$\left[  7\right]  .$

Using the fact that $\Lambda_{l}(y_{0})=0$ and the expression in $\left(
9.26\right)  ,$ $\nabla_{ijk}\Lambda_{l}(y_{0})$ simplifies to:

$\left(  9.29\right)  \qquad\nabla_{ijk}\Lambda_{l}(y_{0})=\frac{\partial
^{3}\Lambda_{l}}{\partial x_{i}\partial x_{j}\partial x_{k}}(y_{0}%
)+[\frac{\partial\Lambda_{k}}{\partial x_{i}}\frac{\partial\Lambda_{l}%
}{\partial x_{j}}+\frac{\partial\Lambda_{j}}{\partial x_{i}}\frac
{\partial\Lambda_{l}}{\partial x_{k}}](y_{0})$

$\qquad\qquad\qquad\qquad\ \ \ \ \ =\frac{1}{4}\frac{\partial^{2}\Omega_{kl}%
}{\partial x_{i}\partial x_{j}}(y_{0})+\frac{1}{4}\Omega_{ik}(y_{0}%
)\Omega_{jl}(y_{0})+\frac{1}{4}\Omega_{ij}(y_{0})\Omega_{kl}(y_{0})$

We use the equalites in $\left(  9.28\right)  $ and $\left(  9.29\right)  $ to
have the expansion in (iii) in terms of \textbf{geometric invariants}:

$\qquad\ \Lambda_{l}(x)=\frac{1}{2}\Omega_{il}(y_{0})x_{i}+\frac{1}{6}%
\frac{\partial\Omega_{jl}}{\partial\text{x}_{i}}(y_{0})(y_{0})x_{i}x_{j}$

$\qquad\qquad-\frac{1}{36}[\overset{n}{\underset{r=1}{\sum}}\left(
R_{ijpr}+R_{ipjr}\right)  \Omega_{kr}](y_{0})x_{i}x_{j}x_{k}$

$\qquad\qquad+\frac{1}{24}[\frac{\partial^{2}\Omega_{kl}}{\partial
x_{i}\partial x_{j}}+\Omega_{ik}\Omega_{jl}+\Omega_{ij}\Omega_{kl}%
](y_{0})x_{i}x_{j}x_{k}+$ higher order terms

\qquad\qquad\qquad\qquad\qquad\qquad\qquad\qquad\qquad\qquad\qquad\qquad
\qquad\qquad\qquad\qquad$\blacksquare$

\part{LEADING COEFFICIENTS OF THE EXPANSIONS}

\chapter{Computations of the First Two Coefficients}

In this Chapter we give the expressions for the coefficients b$_{0}$(x$,$P)
and b$_{1}$($y_{0},$P) in terms of the \textbf{geometry} of the
\textbf{Riemannian manifold} M and the \textbf{submanifold} P. Given the
volume of computations involved, the third and fourth coefficients b$_{2}%
$($y_{0},$P) and b$_{3}$($y_{0},$P) will deserve full Chapters of their own.

The definitions of \textbf{b}$_{0}($x$,$P,$\phi)$ and \textbf{b}$_{1}($%
x$,$P,$\phi)$ are taken from Chapter 5.

\section{The First Coefficient}

\begin{theorem}
(\ Computation of b$_{0})\qquad$
\end{theorem}

\begin{center}
\textbf{b}$_{0}($x$,$P,$\phi)=$ $\tau_{\text{x,y}}\phi($x$)$
\end{center}

where $\tau_{\text{x,y}}:$ E$_{\text{x}}\longrightarrow$ E$_{\text{y}}$ is the
parallel translation on fibers along the (unique minimal) geodesic $\gamma$
from x to y$\in$P in time t.

The first coefficient is thus given by:

\begin{center}
\textbf{b}$_{0}($x$,$P,$\phi)=$ \textbf{b}$_{0}($x$,$P$)\phi($x$)=$
$\tau_{\text{x,y}}\phi($x$)$
\end{center}

In \textbf{Theorem 3} it was shown that:

\begin{center}
\textbf{b}$_{0}($x$,$P,$\phi)=\phi(\gamma(t))$
\end{center}

where $\gamma:[0,$t$]\rightarrow$ M$_{0}$ is the unique minimal geodesic from
x$\in$M$_{0}$ to the submanifold P at the point y$\in$P in time:

$\gamma(0)=x$ and $\gamma(t)=y\in$P and so $\phi(\gamma(t))=\phi(y).$

Let $\tau_{\text{x,y}}:$E$_{\text{x}}\longrightarrow$ E$_{\text{y}}$ be the
\textbf{parallel translation} on fibers along the (unique minimal) geodesic
$\gamma:[0,$t$]\rightarrow$ M$_{0}$ from x to y$\in$P in time t. Then from the
(expansion theorem) \textbf{Theorem 3 }and the defintion of $\tau_{\text{x,y}%
}$, we have:

\begin{center}
\textbf{b}$_{0}($x$,$P,$\phi)=\phi(\gamma(t))=\phi(y)=$ $\tau_{\text{x,y}}%
\phi($x$)$
\end{center}

Since \textbf{b}$_{0}($x$,$P,$\phi)=\tau_{\text{x,y}}\phi($x$),$ we can set,

\begin{center}
\textbf{b}$_{0}$(x,P) = $\tau_{\text{x,y}}$
\end{center}

and have:

\begin{center}
\textbf{b}$_{0}($x$,$P,$\phi)=$ $\tau_{\text{x,y}}\phi($x$)$
\end{center}

where \textbf{b}$_{0}($x$,$P$)=\tau_{\text{x,y}}\in$Hom$($E$_{\text{x}}%
,$E$_{\text{y}})$ is the parallel translation on fibers along the (unique
minimal) geodesic $\gamma.$

We remark that this is the only coefficient we obtain at a general point
x$\in$M$_{0}.$

In particular, when the Fermi coordinates reduce to normal coordinates, or
equivalently, P = $\left\{  \text{y}_{0}\right\}  ,$ then

\begin{center}
\qquad\textbf{b}$_{0}$(x,P) = \textbf{b}$_{0}$(x,y$_{0}$) $=\tau
_{\text{x,y}_{0}}\in$Hom$($E$_{\text{x}},$E$_{\text{y}_{0}})$
\end{center}

This is to be compared to A$_{0}$(x,y$_{0}$) in Theorem $\left(  8.1\right)  $
in \textbf{Duistermaat} $\left[  1\right]  $and to $\Theta$(x,y$_{0}$) in
\textbf{Theorem }$\left(  7.15\right)  $ of \textbf{Roe }$\left[  1\right]  .$

If further, x = y$_{0}$ then, the geodesic $\gamma$ is a constant geodesic and,

\begin{center}
\qquad\ \textbf{b}$_{0}($y$_{0},$y$_{0}$) $=\tau_{\text{y}_{0}\text{,y}_{0}%
}\in$End(E$_{\text{y}_{0}}$),
\end{center}

\qquad\qquad\qquad\qquad\qquad\qquad\qquad\qquad\qquad\qquad\qquad\qquad
\qquad\qquad\qquad\qquad\qquad$\blacksquare$

\section{The Second Coefficient$\qquad$}

The computation for b$_{1}$($x,$P) will be carried out first at the general
point x$\in$M$_{0}$ but this will not reveal the geometry of the underlying
geometric objects: the submanifold P, the manifold M and the vector bundle E.
Computing the second coefficient at the centre of Fermi coordinates
$x=y_{0}\in$P gives a more elegant expression revealing the geometry of P, M
and E.$\qquad\qquad\qquad\qquad\qquad\qquad\qquad\qquad\qquad\qquad
\qquad\qquad\qquad\qquad\qquad\qquad\qquad\qquad\qquad\qquad\qquad\qquad
\qquad\qquad\qquad\qquad\qquad\qquad$

In order to avoid too many of the superscript $"0"$ on $\nabla$ and $\Delta,$
we will, as stated earlier, write $\Delta$ instead of $\Delta^{0}$ and
$\nabla$ instead of $\nabla^{0}$ when applied to functions.

The expression for \textbf{b}$_{1}($x$,$P,$\phi)$ is given in $\left(
C_{17}\right)  $ of \textbf{Appendix C}. If we replace the general point x by
the center of Fermi coordinates y$_{0},$then the second coefficient
\textbf{b}$_{1}($y$_{0},$P,$\phi)$ will fully \textbf{exibit the geometry of
}the underlying geometric entities: the Riemannian manifold M, the submanifold
P and the vector bundle E.\textbf{ }Computing \textbf{b}$_{1}($y$_{0},$%
P,$\phi)$ is analagous to computing the corresponding Minakshisundaram-Pleijel
\textbf{heat kernel} expansion coefficient along the diagonal of the
Riemannian manifold.

\qquad\qquad\qquad\qquad\qquad\qquad\qquad\qquad\qquad\qquad\qquad\qquad
\qquad\qquad\qquad\qquad$\blacksquare$

\subsection{The Raw Expression for the Second Term}

\begin{proposition}
b$_{1}$(x,P,$\phi)=\int_{0}^{1}$F(1,1-r$_{1}$)$[$L$_{\Psi}[\phi\circ
\pi_{\text{P}}](x_{0})$dr$_{1}=\int_{0}^{1}$L$_{\Psi}[\phi\circ\pi_{\text{P}%
}](z_{0})$dr\qquad\qquad\qquad\qquad$\ \ \ $
\end{proposition}

$\qquad\qquad\ =\int_{0}^{1}[\frac{\text{L}\Psi}{\Psi}\left(  z_{0}\right)
.\phi(y)$\qquad\qquad

\qquad\qquad$+$ \ $\frac{1}{2}\underset{\text{a,b=1}}{\overset{\text{q}}{\sum
}}$g$^{\text{ab}}(z_{0})\frac{\partial^{2}\phi}{\partial\text{x}_{\text{a}%
}\partial\text{x}_{\text{b}}}(y)$ $+\frac{1}{2}%
\underset{i,j=1}{\overset{n}{\sum}}$g$^{ij}(z_{0})\frac{\partial\Lambda_{j}%
}{\partial\text{x}_{i}}(z_{0})\phi(y)$

$\qquad\qquad+\frac{1}{2}\underset{j=1}{\overset{n}{\sum}}\underset{\text{a=1}%
}{\overset{\text{q}}{\sum}}$g$^{\text{a}j}(z_{0})[$ $\Lambda_{j}(z_{0}%
)\frac{\partial\phi}{\partial\text{x}_{\text{a}}}(y)]+\frac{1}{2}%
\underset{i=1}{\overset{n}{\sum}}\underset{\text{b=1}}{\overset{\text{q}%
}{\sum}}$g$^{i\text{b}}(z_{0})[\Lambda_{i}(z_{0})\frac{\partial\phi}%
{\partial\text{x}_{\text{b}}}(y)]$

\qquad\qquad$-\frac{1}{2}\underset{i,j=1}{\overset{n}{\sum}}$g$^{ij}%
(z_{0})[\Gamma_{ij}^{\text{c}}(z_{0})\frac{\partial\phi}{\partial
\text{x}_{\text{c}}}(y)$ + $\Gamma_{ij}^{k}(z_{0})\Lambda_{k}(z_{0}%
)\phi(y)]\ +\frac{1}{2}$W$(z_{0})\phi(y)$

\qquad\ \ \qquad+ $\underset{\text{a=1}}{\overset{\text{q}}{\sum}}(\nabla
^{0}\log\Psi)_{\text{a}}(z_{0})\frac{\partial\phi}{\partial\text{x}_{\text{a}%
}}(y)$ + $\underset{j=1}{\overset{n}{\sum}}$ $(\nabla^{0}\log\Psi)_{j}%
(z_{0})\Lambda_{j}(z_{0})\phi(y)$

\qquad\qquad\ + $\underset{\text{a=1}}{\overset{\text{q}}{\sum}}$X$_{\text{a}%
}(z_{0})\frac{\partial\phi}{\partial\text{x}_{\text{a}}}(y)$ +
$\underset{j=1}{\overset{n}{\sum}}$ X$_{j}(z_{0})\Lambda_{j}(z_{0})\phi
(y)]$dr$_{1}$

\begin{proof}
\qquad
\end{proof}

\qquad The expression of \textbf{b}$_{1}($x$,$P,$\phi)$ above is taken from
$\left(  C_{17}\right)  $ of \textbf{Appendix C}

\qquad\qquad\qquad\qquad\qquad\qquad\qquad\qquad\qquad\qquad\qquad\qquad
\qquad\qquad\qquad\qquad\qquad$\blacksquare$\qquad

We call for \textbf{b}$_{1}($x$,$P,$\phi)$ the Raw Expression for the second
term in the expansion of the Generalized Heat Kernel. As we can see, computing
the second term at a general point x$_{0}\in M_{0}$ does not, unfortunately,
reveal the geometry of the underlying spaces. In order to do so, we must
compute it at the \textbf{centre} of Fermi coordinates y$_{0}\in P.$ Computing
the second term at this point will yield an expression given in geometric
invariants of Riemmanian manifold M, the submanifold P and the vector bundle E.

For the purpose of clarity, we proceed by first giving some Computational
Lemmas. The expression in the Lemma below is a preliminary version of the
second term of the heat kernel expansion.

\begin{lemma}
$\qquad\qquad\qquad$
\end{lemma}

\textbf{b}$_{1}($y$_{0},$P,$\phi)=\frac{\text{L}\Psi}{\Psi}(y_{0})\phi(y_{0})$
+ \ $\frac{1}{2}\underset{\text{a=1}}{\overset{\text{q}}{\sum}}\frac
{\partial^{2}\phi}{\partial\text{x}_{\text{a}}^{2}}(y_{0})+$
$\underset{\text{a=1}}{\overset{\text{q}}{\sum}}\Lambda_{\text{a}}(y_{0}%
)\frac{\partial\phi}{\partial x_{\text{a}}}\left(  y_{0}\right)  $

$\qquad+\frac{1}{2}$ $\underset{\text{a=1}}{\overset{\text{q}}{\sum}}%
\Lambda_{\text{a}}(y_{0})\Lambda_{\text{a}}(y_{0})\phi\left(  y_{0}\right)
+\frac{1}{2}$W$_{y_{0}}(\phi\left(  y_{0}\right)  )+$ $\underset{\text{a=1}%
}{\overset{\text{q}}{\sum}}X_{\text{a}}(y_{0})\frac{\partial\phi}%
{\partial\text{x}_{\text{a}}}(y_{0})$

\qquad\ + $\underset{\text{a=1}}{\overset{\text{q}}{\sum}}$ $X_{\text{a}%
}(y_{0})\Lambda_{\text{a}}(y_{0})\phi(y_{0})+$
$\underset{j=q+1}{\overset{n}{\sum}}$ X$_{j}(y_{0})\Lambda_{j}(y_{0}%
)\phi(y_{0})$

\qquad\qquad\qquad\qquad\qquad\qquad\qquad\qquad\qquad\qquad\qquad\qquad
\qquad\qquad$\qquad\qquad\blacksquare$

The expression of the Lemma is taken from $\left(  10.20\right)  $ below. The
details of computation are as follows:

The proof is purely computational. It is however very long because we want to
give all details of computation.

\textbf{b}$_{1}$\textbf{(x}$_{0}$\textbf{,P,}$\phi)=\int_{0}^{1}$F(1,1-r$_{1}%
$)$[$L$_{\Psi}[\phi\circ\pi_{\text{P}}](x_{0})$dr$_{1}=\int_{0}^{1}$L$_{\Psi
}[\phi\circ\pi_{\text{P}}](z_{0})$dr

The integrand L$_{\Psi}$[$\phi\circ\pi_{\text{P}}$]$(z_{0})$ is \textbf{not
}independent of r$_{1}$ since z$_{0}=\gamma_{0,1}(r_{1})$ with r$_{1}%
\in\left[  0,1\right]  $ where $\gamma_{0,1}:\left[  0,1\right]
\longrightarrow$ M$_{0}$ is the unique minimal geodesic from a general point
$x_{0}\in$M$_{0}$ to a general point $y\in U\subset$P in time 1.

Since z$_{0}$ depends on r$_{1}\in\left[  0,1\right]  ,$ the integral above
depends on r$_{1}\in\left[  0,1\right]  .$

By the definition of $\gamma_{0,1}$ above, we have in Fermi coordinates:

$\left(  10.1\right)  \qquad\qquad\gamma_{0,1}(s)=\left(  x_{1},...,x_{q}%
,(1-s)x_{q+1},...,(1-\frac{s}{1-r_{2}})x_{n}\right)  =$ General Point

\qquad\qquad\qquad\qquad\qquad$=y_{0}+(1-s)(x_{0}-y_{0})$

In particular, we see that: $\gamma_{0,1}(0)=\left(  x_{1},...,x_{q}%
,x_{q+1},...,x_{n}\right)  =x_{0}\in M_{0}=$ Starting Point

$\left(  10.2\right)  \qquad\qquad z_{0}=\gamma_{0,1}(r_{1})=\left(
x_{1},...,x_{q},(1-r_{1})x_{q+1},...,(1-r_{1})x_{n}\right)  =$ Mid Point

\qquad\qquad\qquad\qquad$\gamma_{0,1}(1)=\left(  x_{1},...,x_{q}%
,0,...,0\right)  \approxeq\left(  x_{1},...,x_{q}\right)  =y_{0}\in P=$ End
Point\qquad\qquad\qquad$\ \qquad\qquad\qquad$

From the definition of $\gamma_{0,1}(s)$ above, we have:

$\left(  10.3\right)  $\qquad\qquad$\frac{\partial}{\partial x_{i}}%
\gamma_{0,1}(s)=\left\{
\begin{array}
[c]{c}%
1\text{ for }i=1,...,q\\
(1-s)\text{ for }i=q+1,...,n
\end{array}
\right.  $

In particular, we have:

$\left(  10.4\right)  $ \qquad$\frac{\partial z_{0}}{\partial x_{i}}%
=\frac{\partial}{\partial x_{i}}\gamma_{0,1}(r_{1})=\left\{
\begin{array}
[c]{c}%
1\text{ for }i=1,...,q\\
1-r_{1}\text{ for }i=q+1,...,n
\end{array}
\right.  $

From the definition of $z_{0},$ we have:

$\left(  10.5\right)  \qquad\qquad\qquad\frac{\partial}{\partial\text{x}_{i}%
}\pi_{\text{P}}$(z$_{0}$) $=\frac{\partial}{\partial\text{x}_{i}}\pi
_{\text{P}}(x_{0})=\left\{
\begin{array}
[c]{c}%
1\text{ for }i=1,...,q\\
0\text{ for }i=q+1,...,n
\end{array}
\right.  $

$\left(  10.6\right)  $\qquad\qquad$\qquad\frac{\partial^{2}}{\partial
x_{i}\partial x_{j}}\gamma_{0,1}(r_{1})=0$ for all $i,j=1,...,q,q+1,...,n.$

We also have:

$\qquad$

$\left(  10.7\right)  \qquad\qquad\qquad\frac{\partial^{2}}{\partial
\text{x}_{i}\partial\text{x}_{j}}\pi_{\text{P}}$(z$_{0}$) $=0=\frac
{\partial^{2}}{\partial\text{x}_{j}\partial\text{x}_{i}}\pi_{\text{P}}$(x) for
all $i,j=1,...,q,q+1,...,n$

\qquad\qquad\qquad\qquad\qquad\qquad\qquad\qquad\qquad\qquad\qquad\qquad
\qquad\qquad\qquad\qquad$\blacksquare$

To compute the second coefficient at the centre of Fermi coordinates $y_{0}%
\in$P, we assume:

$\left(  10.8\right)  \qquad\qquad\qquad$ $x_{0}=y_{0}$

and compute the second coefficient at the centre of Fermi coordinates
$y_{0}\in$P. In this case, $\gamma_{0,1}$ is the constant geodesic:
$\gamma_{0,1}(s)=y_{0}\forall s\in\left[  0,1\right]  $ and so, $z_{0}%
=\gamma_{0,1}(r_{1})=y_{0}.$ The integrand is thus independent of r$_{1}$ and
so the integration in $\left(  C_{17}\right)  $ is trivial and becomes:

and so by $\left(  10.8\right)  ,$

$\left(  10.9\right)  \qquad$ $z_{0}=\gamma(r_{1})=y_{0}=x_{0}.$

\qquad The integrand is thus independent of r$_{1}$ and so by
\textbf{Proposition 11} above,

we have:

$\left(  10.10\right)  $\qquad\textbf{b}$_{1}$\textbf{(}$y_{0}$\textbf{,P,}%
$\phi)=\int_{0}^{1}$F(1,1-r$_{1}$)$[$L$_{\Psi}[\phi\circ\pi_{\text{P}}%
](y_{0})$dr$_{1}=$ L$_{\Psi}[\phi\circ\pi_{\text{P}}](y_{0})$

\qquad$=[\frac{\text{L}\Psi}{\Psi}\left(  y_{0}\right)  .\phi(y_{0})$%
\qquad\qquad

\qquad$+$ \ $\frac{1}{2}\underset{\text{a,b=1}}{\overset{\text{q}}{\sum}}%
$g$^{\text{ab}}(y_{0})\frac{\partial^{2}\phi}{\partial\text{x}_{\text{a}%
}\partial\text{x}_{\text{b}}}(y_{0})$ $+\frac{1}{2}%
\underset{i,j=1}{\overset{n}{\sum}}$g$^{ij}(y_{0})\frac{\partial\Lambda_{j}%
}{\partial\text{x}_{i}}(y_{0})\phi(y_{0})$

$\qquad+\frac{1}{2}\underset{j=1}{\overset{n}{\sum}}\underset{\text{a=1}%
}{\overset{\text{q}}{\sum}}$g$^{\text{a}j}(y_{0})[$ $\Lambda_{j}(y_{0}%
)\frac{\partial\phi}{\partial\text{x}_{\text{a}}}(y_{0})]+\frac{1}%
{2}\underset{i=1}{\overset{n}{\sum}}\underset{\text{b=1}}{\overset{\text{q}%
}{\sum}}$g$^{i\text{b}}(y_{0})[\Lambda_{i}(y_{0})\frac{\partial\phi}%
{\partial\text{x}_{\text{b}}}(y_{0})]$

\qquad$-\frac{1}{2}\underset{i,j=1}{\overset{n}{\sum}}$g$^{ij}(y_{0}%
)[\Gamma_{ij}^{\text{c}}(y_{0})\frac{\partial\phi}{\partial\text{x}_{\text{c}%
}}(y_{0})$ + $\Gamma_{ij}^{k}(y_{0})\Lambda_{k}(y_{0})\phi(y_{0})]\ +\frac
{1}{2}$W$(y_{0})\phi(y_{0})$

\qquad\ \ $+$ $\underset{\text{a=1}}{\overset{\text{q}}{\sum}}(\nabla^{0}%
\log\Psi)_{\text{a}}(y_{0})\frac{\partial\phi}{\partial\text{x}_{\text{a}}%
}(y)$ + $\underset{j=1}{\overset{n}{\sum}}$ $(\nabla^{0}\log\Psi)_{j}%
(y_{0})\Lambda_{j}(y_{0})\phi(y)$

\qquad$+$ $\underset{\text{a=1}}{\overset{\text{q}}{\sum}}$X$_{\text{a}}%
(y_{0})\frac{\partial\phi}{\partial\text{x}_{\text{a}}}(y_{0})$ +
$\underset{j=1}{\overset{n}{\sum}}$ X$_{j}(y_{0})\Lambda_{j}(y_{0})\phi
(y_{0})$\qquad\qquad

\qquad\qquad\qquad\qquad\qquad\qquad\qquad\qquad\qquad\qquad\qquad\qquad
\qquad\qquad\qquad\qquad$\blacksquare$\qquad

We simplify the above expression:

$g^{ij}(y_{0})=\delta^{ij}$ and $\Gamma_{ii}^{\text{c}}(y_{0})=0=$
$\Gamma_{\text{aa}}^{\text{c}}(y_{0})=0=\Gamma_{ij}^{k}(y_{0})$

for a,c = 1,....,q; $i,j,k=q+1,...,n,$ and so we have:

$\qquad-$ $\frac{1}{2}\underset{i,j=1}{\overset{n}{\sum}}g^{ij}(y_{0}%
)\Gamma_{ij}^{\text{c}}(y_{0})\frac{\partial\phi}{\partial x_{\text{c}}%
}\left(  y_{0}\right)  =-$ $\frac{1}{2}\underset{i=1}{\overset{n}{\sum}}%
\Gamma_{ii}^{\text{c}}(y_{0})\frac{\partial\phi}{\partial x_{\text{c}}}\left(
y_{0}\right)  =0$

Next we have:

\qquad g$^{ij}(y_{0})\frac{\partial\Lambda_{j}}{\partial\text{x}_{i}}%
(y_{0})=\frac{\partial\Lambda i}{\partial\text{x}_{i}}(y_{0})=\frac{1}%
{2}\Omega_{ii}(y_{0})=0$

The first equality is obvious by the fact that g$^{ij}(y_{0})=\delta^{ij}.$
The second is by $\left(  3.23\right)  $ above

and the last is due to the fact that $\Omega_{ij}(y_{0})$ is skew-symmetric in
the indices by $\left(  3.14\right)  $ above.

We have:

$\left(  10.11\right)  \qquad-\frac{1}{2}\underset{i,j=1}{\overset{n}{\sum}%
}g^{ij}(y_{0})$ $\underset{k=1}{\overset{n}{\sum}}\Gamma_{ij}^{k}%
(y_{0})\Lambda_{k}(y_{0})\phi\left(  y_{0}\right)  =-$ $\frac{1}{2}$
$\underset{i,k=1}{\overset{n}{\sum}}\Gamma_{ii}^{k}(y_{0})\Lambda_{k}%
(y_{0})\phi\left(  y_{0}\right)  $

$\ =-$ $\frac{1}{2}$ $\underset{k=1}{\overset{n}{\sum}}\underset{\text{a}%
=1}{\overset{\text{q}}{\sum}}\Gamma_{\text{aa}}^{k}(y_{0})\Lambda_{k}%
(y_{0})\phi\left(  y_{0}\right)  -$ $\frac{1}{2}$
$\underset{k=1}{\overset{n}{\sum}}\underset{i=q+1}{\overset{n}{\sum}}%
\Gamma_{ii}^{k}(y_{0})\Lambda_{k}(y_{0})\phi\left(  y_{0}\right)  $

$=-$ $\frac{1}{2}$ $\underset{\text{b}=1}{\overset{\text{q}}{\sum}%
}\underset{\text{a}=1}{\overset{\text{q}}{\sum}}\Gamma_{\text{aa}}^{\text{b}%
}(y_{0})\Lambda_{\text{b}}(y_{0})\phi\left(  y_{0}\right)  -\frac{1}%
{2}\underset{k=q+1}{\overset{n}{\sum}}\underset{\text{a}=1}{\overset{\text{q}%
}{\sum}}\Gamma_{\text{aa}}^{k}(y_{0})\Lambda_{k}(y_{0})\phi\left(
y_{0}\right)  $

$-$ $\frac{1}{2}$ $\underset{\text{a}=1}{\overset{\text{q}}{\sum}%
}\underset{i=q+1}{\overset{n}{\sum}}\Gamma_{ii}^{\text{a}}(y_{0})\Lambda
_{k}(y_{0})\phi\left(  y_{0}\right)  -$ $\frac{1}{2}$
$\underset{k=q+1}{\overset{n}{\sum}}\underset{i=q+1}{\overset{n}{\sum}}%
\Gamma_{ii}^{k}(y_{0})\Lambda_{k}(y_{0})\phi\left(  y_{0}\right)  $

$=-\frac{1}{2}\underset{k=q+1}{\overset{n}{\sum}}\underset{\text{a}%
=1}{\overset{\text{q}}{\sum}}\Gamma_{\text{aa}}^{k}(y_{0})\Lambda_{k}%
(y_{0})\phi\left(  y_{0}\right)  =-\frac{1}{2}%
\underset{k=q+1}{\overset{n}{\sum}}\underset{\text{a}=1}{\overset{\text{q}%
}{\sum}}T_{\text{aa}k}(y_{0})\Lambda_{k}(y_{0})\phi\left(  y_{0}\right)  $

Next, we have:

$\underset{j=1}{\overset{n}{\sum}}(\nabla\log\Psi)^{j}\Lambda_{j}\left(
y_{0}\right)  \phi\left(  y_{0}\right)  =$ $\underset{\text{a=1}%
}{\overset{\text{q}}{\sum}}(\nabla\log\Psi)^{\text{a}}\left(  y_{0}\right)
\Lambda_{\text{a}}\left(  y_{0}\right)  \phi\left(  y_{0}\right)  +$
$\underset{j=q+1}{\overset{n}{\sum}}(\nabla\log\Psi)^{j}\Lambda_{j}\left(
y_{0}\right)  \phi\left(  y_{0}\right)  $

From $\left(  10.10\right)  $ \ and above simplifications, we have the
expression in the Lemma:

$\left(  10.12\right)  \qquad\qquad$\textbf{b}$_{1}($y$_{0},$P,$\phi
)=\frac{\text{L}\Psi}{\Psi}(y_{0})\phi(y_{0})$ $+$\ $\frac{1}{2}%
\underset{\text{a=1}}{\overset{\text{q}}{\sum}}\frac{\partial^{2}\phi
}{\partial\text{x}_{\text{a}}^{2}}(y_{0})$

$\qquad\qquad+\frac{1}{2}[$ $\underset{\text{a=1}}{\overset{\text{q}}{\sum}%
}\Lambda_{\text{a}}(y_{0})\frac{\partial\phi}{\partial x_{\text{a}}}\left(
y_{0}\right)  +$ $\underset{\text{b=1}}{\overset{\text{q}}{\sum}}%
\Lambda_{\text{b}}(y_{0})\frac{\partial\phi}{\partial x_{\text{b}}}\left(
y_{0}\right)  ]$

$\qquad\qquad+\frac{1}{2}[$ $\underset{\text{a=1}}{\overset{\text{q}}{\sum}%
}\Lambda_{\text{a}}(y_{0})\Lambda_{\text{a}}(y_{0})\phi\left(  y_{0}\right)
]-$ $\frac{1}{2}$ $\underset{k=q+1}{\overset{n}{\sum}}\underset{\text{a}%
=1}{\overset{\text{q}}{\sum}}T_{\text{aa}k}(y_{0})\Lambda_{k}(y_{0}%
)\phi\left(  y_{0}\right)  $

$\qquad\qquad+\frac{1}{2}$W$_{y_{0}}(\phi(y_{0}))$

\qquad$\qquad+$ $\underset{\text{a=1}}{\overset{\text{q}}{\sum}}(\nabla
\log\Psi)_{\text{a}}(y_{0})\frac{\partial\phi}{\partial x_{\text{a}}}\left(
y_{0}\right)  +$ $\underset{\text{a=1}}{\overset{\text{q}}{\sum}}(\nabla
\log\Psi)_{\text{a}}\left(  y_{0}\right)  \Lambda_{\text{a}}\left(
y_{0}\right)  \phi\left(  y_{0}\right)  $

$\qquad\qquad+$ $\underset{j=q+1}{\overset{n}{\sum}}(\nabla\log\Psi
)_{j}\Lambda_{j})\left(  y_{0}\right)  \phi\left(  y_{0}\right)  $

\qquad\qquad$+$ $\underset{\text{a=1}}{\overset{\text{q}}{\sum}}$X$_{\text{a}%
}(y_{0})\frac{\partial\phi}{\partial\text{x}_{\text{a}}}(y_{0})$ +
$\underset{\text{a=1}}{\overset{\text{q}}{\sum}}$ X$_{\text{a}}(y_{0}%
)\Lambda_{\text{a}}(y_{0})\phi(y_{0})+$ $\underset{j=q+1}{\overset{\text{n}%
}{\sum}}$ X$_{j}(y_{0})\Lambda_{j}(y_{0})\phi(y_{0})\mathbf{.}$

\qquad\qquad\qquad\qquad\qquad\qquad\qquad\qquad\qquad\qquad\qquad\qquad
\qquad\qquad\qquad\qquad$\qquad\qquad\blacksquare$

We need to express \textbf{b}$_{1}($y$_{0},$P,$\phi)$ in terms of the geometry
of the Riemannian manifold M, the submanifold P and the vector bundle E and so
we make the following computations below:

$\nabla\log\Psi=\nabla\log\theta^{-\frac{1}{2}}+\nabla\log\Phi$

where $\Phi$, $\theta$ and $\Psi$ are defined respectively in $\left(
1.5\right)  $, $\left(  1.6\right)  $ and $\left(  1.7\right)  $ of

Chapter 1. Since,

\qquad\qquad\qquad\qquad\qquad\qquad\qquad\ $\theta(y_{0})=1=\Phi(y_{0}),$

By (iii)$^{\ast}$ of Table A$_{9,}$

$\qquad\qquad\qquad\qquad\qquad\left(  \nabla\log\theta^{-\frac{1}{2}}\right)
(y_{0})_{\text{a}}=0$ for a = 1,...,q

and by (iv)$^{\ast}$ of \textbf{Table A}$_{9},$

$\left(  10.13\right)  \qquad\qquad\qquad\left(  \nabla\log\theta^{-\frac
{1}{2}}\right)  _{i}(y_{0})=\frac{1}{2}<H,i>(y_{0})$ for $i=q+1,...,n;$

By (xi) of \textbf{Table B}$_{1},$

$\left(  10.14\right)  \qquad\qquad\qquad\left(  \nabla\log\Phi\right)
_{\text{a}}(y_{0})=0$ for a = 1,...,q

and by (vi) of Table \textbf{B}$_{1}.$

$\left(  10.15\right)  \qquad\qquad\qquad\left(  \nabla\log\Phi\right)
_{j}(y_{0})=-X_{j}(y_{0})$ for $j=q+1,...,n$

We conclude that for for a = 1,...,q and $j,k=q+1,...,n,$ we have:

$\left(  10.16\right)  \qquad<\nabla\log\theta^{-\frac{1}{2}},\nabla\log
\Phi>(y_{0})=-\frac{1}{2}<H,j>(y_{0})X_{j}(y_{0})$

$\left(  10.17\right)  \qquad\left(  \nabla\log\Psi\right)  _{\text{a}}%
(y_{0})=\left(  \nabla\log\theta^{-\frac{1}{2}}\right)  _{\text{a}}%
(y_{0})+\left(  \nabla\log\Phi\right)  _{\text{a}}(y_{0})=0$

$\left(  10.18\right)  \qquad(\nabla\log\Psi)_{j}(y_{0})=\left(  \nabla
\log\theta^{-\frac{1}{2}}\right)  _{j}(y_{0})+\left(  \nabla\log\Phi\right)
_{j}(y_{0})$

$\qquad\qquad\qquad=\frac{1}{2}<H,j>-$ $X_{j}$

$\left(  10.19\right)  \qquad<\nabla\log\Psi,X>_{j}(y_{0})=\frac{1}%
{2}<H,j>X_{j}-$ $X_{j}^{2}$

We use the fact that:

$(\nabla\log\Psi)_{\text{a}}\left(  y_{0}\right)  =0,$ $(\nabla\log\Psi
)_{j}(y_{0})=$ $\frac{1}{2}<H,j>-$ $X_{j}$ and

$\underset{\text{a}=1}{\overset{\text{q}}{\sum}}T_{\text{aa}k}(y_{0})=$
$<H,k>$

to simplify the expression in $\left(  10.12\right)  $ and have:

$\qquad$\textbf{b}$_{1}($y$_{0},$P,$\phi)=\frac{\text{L}\Psi}{\Psi}(y_{0}%
)\phi(y_{0})$ + \ $\frac{1}{2}\underset{\text{a=1}}{\overset{\text{q}}{\sum}%
}\frac{\partial^{2}\phi}{\partial\text{x}_{\text{a}}^{2}}(y_{0})+$
$\underset{\text{a=1}}{\overset{\text{q}}{\sum}}\Lambda_{\text{a}}(y_{0}%
)\frac{\partial\phi}{\partial x_{\text{a}}}\left(  y_{0}\right)  $

$\qquad+\frac{1}{2}$ $\underset{\text{a=1}}{\overset{\text{q}}{\sum}}%
\Lambda_{\text{a}}(y_{0})\Lambda_{\text{a}}(y_{0})\phi\left(  y_{0}\right)
-\frac{1}{2}\underset{k=q+1}{\overset{n}{\sum}}<H,k>(y_{0})\Lambda_{k}%
(y_{0})\phi\left(  y_{0}\right)  +\frac{1}{2}$W$_{y_{0}}(\phi\left(
y_{0}\right)  )$

$\qquad+\frac{1}{2}$ $\underset{j=q+1}{\overset{n}{\sum}}[<H,j>-X_{j}%
]\Lambda_{j}\left(  y_{0}\right)  \phi\left(  y_{0}\right)  +$
$\underset{\text{a=1}}{\overset{\text{q}}{\sum}}X_{\text{a}}(y_{0}%
)\frac{\partial\phi}{\partial\text{x}_{\text{a}}}(y_{0})$ +
$\underset{\text{a=1}}{\overset{\text{q}}{\sum}}$ $X_{\text{a}}(y_{0}%
)\Lambda_{\text{a}}(y_{0})\phi(y_{0})$

$\qquad\qquad\qquad\qquad+$ $\underset{j=q+1}{\overset{\text{n}}{\sum}}$
X$_{j}(y_{0})\Lambda_{j}(y_{0})\phi(y_{0})$

There is an obvious cancellation above and so we have the final expression:

$\left(  10.20\right)  \qquad$\textbf{b}$_{1}($y$_{0},$P,$\phi)=\frac
{\text{L}\Psi}{\Psi}(y_{0})\phi(y_{0})$ + \ $\frac{1}{2}\underset{\text{a=1}%
}{\overset{\text{q}}{\sum}}\frac{\partial^{2}\phi}{\partial\text{x}_{\text{a}%
}^{2}}(y_{0})+$ $\underset{\text{a=1}}{\overset{\text{q}}{\sum}}%
\Lambda_{\text{a}}(y_{0})\frac{\partial\phi}{\partial x_{\text{a}}}\left(
y_{0}\right)  $

$\qquad\qquad\qquad+\frac{1}{2}$ $\underset{\text{a=1}}{\overset{\text{q}%
}{\sum}}\Lambda_{\text{a}}(y_{0})\Lambda_{\text{a}}(y_{0})\phi\left(
y_{0}\right)  +\frac{1}{2}$W$_{y_{0}}(\phi\left(  y_{0}\right)  )+$
$\underset{\text{a=1}}{\overset{\text{q}}{\sum}}X_{\text{a}}(y_{0}%
)\frac{\partial\phi}{\partial\text{x}_{\text{a}}}(y_{0})$

\qquad\qquad\qquad\ + $\underset{\text{a=1}}{\overset{\text{q}}{\sum}}$
$X_{\text{a}}(y_{0})\Lambda_{\text{a}}(y_{0})\phi(y_{0})+$
$\underset{j=q+1}{\overset{\text{n}}{\sum}}$ X$_{j}(y_{0})\Lambda_{j}%
(y_{0})\phi(y_{0})$

The Lemma above is thus proved.

\qquad\qquad\qquad\qquad\qquad\qquad\qquad\qquad\qquad\qquad\qquad\qquad
\qquad\qquad\qquad\qquad$\blacksquare$

We next give a detailed computation of $\frac{\text{L}\Psi}{\Psi}(y_{0}%
):$\qquad

The computation of $\frac{\text{L}\Psi}{\Psi}(y_{0})$ will \textbf{reveal role
the geometry} of the \textbf{Riemannian manifold M}, the \textbf{submanifold
P} and the \textbf{vector bundle E.} This will give us a final expression of
the second term of the heat kernel expansion.

\qquad\qquad\qquad\qquad\qquad\qquad\qquad\qquad\qquad\qquad\qquad\qquad
\qquad\qquad\qquad\qquad$\blacksquare$

We start by recalling some properties of geometric invariants:

Tables in Appendix A give for a,b = 1,...,q and $i,j=q+1,...,n:\qquad$

$g^{\text{ab}}(y_{0})=\delta_{\text{ab}}$; g$^{ij}(y_{0})=\delta_{ij}$
$\qquad$

$g^{\text{a}j}(y_{0})=0=g^{i\text{b}}(y_{0})$ for a,b = 1,...,q and
$i,j=q+1,...,n$.

$\Lambda_{i}(y_{0})\Lambda_{j}(y_{0})=0$ for $i,j=q+1,...,n$ by $\left(
6.13\right)  $ above.

$\frac{\partial\Lambda_{i}}{\partial x_{i}}(y_{0})=0$ since $\frac
{\partial\Lambda_{j}}{\partial x_{i}}(y_{0})=\frac{1}{2}\Omega_{ij}(y_{0})$ is
skew-symmetric in $\left(  i,j\right)  $

$i,j=q+1,...,n$

$\Gamma_{\text{ab}}^{i}(y_{0})=T_{\text{ab}i}(y_{0})$ by (i) of Table A$_{7}$

$\Gamma_{\text{ab}}^{\text{c}}(y_{0})=0$ for a,b,c=1,...,q by (ii) of Table
A$_{7}$

$\Gamma_{ij}^{k}(y_{0})=0$ for $i,j,k=q+1,...,n$ by (i) of Table A$_{8}$

$\Gamma_{ij}^{\text{c}}(y_{0})=0$ for c = 1,...,q and $i,j=q+1,...,n$ by (ii)
of Table A$_{8}$

$\Gamma_{\text{a}j}^{\text{b}}(y_{0})=-\Gamma_{\text{ab}}^{j}(y_{0}%
)=-$T$_{\text{ab}j}(y_{0})$ by (iii) of Table A$_{8}$

$\Gamma_{\text{a}j}^{k}(y_{0})=$ $\perp_{\text{a}jk}(y_{0})$ by (iv)\ of Table
A$_{8}$

and the fact that $\underset{\text{a}=1}{\overset{\text{q}}{\sum}}%
T_{\text{aa}k}(y_{0})=$ $<H,k>,$ we have the simpler expression:

For the computations below we will often use the formula:

L(f$\phi$) = (Lf)$\phi$ + f(L$\phi$) +
$<$%
$\nabla^{0}$f,$\nabla\phi$%
$>$
$-$ V$($f$\phi).$

The objective below is to express $\frac{\text{L}\Psi}{\Psi}(y_{0})$ in terms
of the geometric invariants of the Riemannian manifold M and the submanifold
P. We recall that:

\begin{center}
L = $\frac{1}{2}\Delta+X+V=\frac{1}{2}\Delta+\nabla_{X}+V$
\end{center}

Since \ $\theta(y_{0})=1=\Phi(y_{0}),$ we have:

$\left(  10.21\right)  \qquad\frac{\text{L}\Psi}{\Psi}(y_{0})=$ $\frac{1}%
{2}\Delta\Psi(y_{0})+$ $<\nabla\Psi,X>(y_{0})+$ V$(y_{0})$

\qquad\qquad$=\frac{1}{2}\Delta\theta^{-\frac{1}{2}}(y_{0})+$ $\frac{1}%
{2}\Delta\Phi(y_{0})+$ $<\nabla\theta^{-\frac{1}{2}},\nabla\Phi>(y_{0})$

$\qquad\qquad+$ $<\nabla\Psi,X>(y_{0})+$ V$(y_{0})$

\qquad\qquad$=\frac{1}{2}\Delta\theta^{-\frac{1}{2}}(y_{0})+$ $\frac{1}%
{2}\Delta\Phi(y_{0})+$ $<\nabla\log\theta^{-\frac{1}{2}},\nabla\log\Phi
>(y_{0})$

$\qquad\qquad+$ $<\nabla\log\Psi,X>(y_{0})+$ V$(y_{0})$

$\left(  \nabla\log\Phi\right)  _{j}(y_{0})=-X_{j}(y_{0})$ for $j=q+1,...,n$

We conclude that for for a = 1,...,q and $j,k=q+1,...,n,$ we have:

$\left(  10.22\right)  \qquad<\nabla\log\theta^{-\frac{1}{2}},\nabla\log
\Phi>(y_{0})=-\frac{1}{2}<H,j>(y_{0})X_{j}(y_{0})$

$\left(  10.23\right)  \qquad\left(  \nabla\log\Psi\right)  _{\text{a}}%
(y_{0})=\left(  \nabla\log\theta^{-\frac{1}{2}}\right)  _{\text{a}}%
(y_{0})+\left(  \nabla\log\Phi\right)  _{\text{a}}(y_{0})=0$

$\left(  10.24\right)  \qquad(\nabla\log\Psi)_{j}(y_{0})=\left(  \nabla
\log\theta^{-\frac{1}{2}}\right)  _{j}(y_{0})+\left(  \nabla\log\Phi\right)
_{j}(y_{0})=\frac{1}{2}<H,j>-$ $X_{j}$

$\left(  10.25\right)  \qquad<\nabla\log\Psi,X>_{j}(y_{0})=\frac{1}%
{2}[<H,j>X_{j}](y_{0})-$ $X_{j}^{2}(y_{0})$

Recall that $\left(  10.21\right)  $ gives:

$\qquad\frac{\text{L}\Psi}{\Psi}(y_{0})=\frac{1}{2}\Delta\theta^{-\frac{1}{2}%
}(y_{0})+$ $\frac{1}{2}\Delta\Phi(y_{0})+$ $<\nabla\log\theta^{-\frac{1}{2}%
},\nabla\log\Phi>(y_{0})$

$\qquad\qquad\qquad\qquad\qquad\qquad+$ $<\nabla\log\Psi,X>(y_{0})+$
V$(y_{0})$

Therefore by $\left(  10.22\right)  $ and $\left(  10.25\right)  ,$ we have:

$\left(  10.26\right)  \qquad\frac{\text{L}\Psi}{\Psi}(y_{0})=\frac{1}%
{2}\Delta\theta^{-\frac{1}{2}}(y_{0})+$ $\frac{1}{2}\Delta\Phi(y_{0})-$
$\underset{j=q+1}{\overset{n}{\sum}}($X$_{j})^{2}(y_{0})+$ V$(y_{0})$

We re-write the last expression above as:

$\left(  10.27\right)  \qquad\frac{\text{L}\Psi}{\Psi}(y_{0})=\frac{1}%
{2}\Delta\theta^{-\frac{1}{2}}(y_{0})+\frac{1}{2}\Delta\Phi(y_{0})-$
$\left\Vert \text{X}\right\Vert ^{2}(y_{0})+$ $\underset{\text{a}%
=1}{\overset{q}{\sum}}X_{\text{a}}^{2}(y_{0})+$ V$(y_{0})$

By (ii) of \textbf{TableA}$_{10},\qquad$

$\left(  10.28\right)  \qquad\frac{1}{2}\Delta\theta^{-\frac{1}{2}}%
(y_{0})=\frac{1}{24}[\underset{i=q+1}{\overset{n}{\sum}}3<H,i>^{2}+2(\tau
^{M}-3\tau^{P}\ +\overset{q}{\underset{\text{a=1}}{\sum}}\varrho_{\text{aa}%
}^{M}+\overset{q}{\underset{\text{a,b}=1}{\sum}}R_{\text{abab}}^{M})](y_{0})$

Next by (iii) of \textbf{Table B}$_{3}.$ This is also given by $\left(
B_{26}\right)  :$

$\left(  10.29\right)  \qquad\frac{1}{2}\Delta\Phi(y_{0})=\frac{1}{2}$
$\left\Vert \text{X}\right\Vert _{M}^{2}(y_{0})-\frac{1}{2}$
$\operatorname{div}X_{M}(y_{0})-\frac{1}{2}$ $\left\Vert \text{X}\right\Vert
_{P}^{2}(y_{0})$ $+\frac{1}{2}$ $\operatorname{div}X_{P}(y_{0})\qquad$

\qquad\qquad\qquad\qquad\qquad\qquad\qquad\qquad\qquad\qquad\qquad\qquad
\qquad\qquad\qquad\qquad\qquad$\blacksquare$

We conclude from $\left(  10.27\right)  ,$ $\left(  10.28\right)  $ and
$\left(  10.29\right)  $ that:

$\qquad\qquad\frac{\text{L}\Psi}{\Psi}(y_{0})=\frac{1}{24}%
[\underset{i=q+1}{\overset{n}{\sum}}3<H,i>^{2}+2(\tau^{M}-3\tau^{P}%
\ +\overset{q}{\underset{\text{a=1}}{\sum}}\varrho_{\text{aa}}^{M}%
+\overset{q}{\underset{\text{a,b}=1}{\sum}}R_{\text{abab}}^{M})](y_{0})$

\qquad\qquad\qquad$\qquad+\frac{1}{2}$ $\left\Vert \text{X}\right\Vert
_{M}^{2}(y_{0})-\frac{1}{2}$ $\operatorname{div}X_{M}(y_{0})-\frac{1}{2}$
$\left\Vert \text{X}\right\Vert _{P}^{2}(y_{0})$ $+$ $\frac{1}{2}%
\operatorname{div}X_{P}(y_{0})$

\qquad\qquad\qquad\qquad$-$ $\left\Vert \text{X}\right\Vert ^{2}(y_{0})+$
$\underset{\text{a}=1}{\overset{q}{\sum}}X_{\text{a}}^{2}(y_{0})+$ V$(y_{0})$

Since $\underset{\text{a}=1}{\overset{q}{\sum}}X_{\text{a}}^{2}(y_{0}%
)=\left\Vert \text{X}\right\Vert _{P}^{2}(y_{0}),$ we simplify the last
expression above and have:

$\left(  10.30\right)  $\qquad$\frac{\text{L}\Psi}{\Psi}(y_{0})=\frac{1}%
{24}[\underset{i=q+1}{\overset{n}{\sum}}3<H,i>^{2}+2(\tau^{M}-3\tau
^{P}\ +\overset{q}{\underset{\text{a=1}}{\sum}}\varrho_{\text{aa}}%
^{M}+\overset{q}{\underset{\text{a,b}=1}{\sum}}R_{\text{abab}}^{M})](y_{0})$

\qquad\qquad$-\frac{1}{2}$ $\left\Vert \text{X}\right\Vert _{M}^{2}%
(y_{0})-\frac{1}{2}$ $\operatorname{div}X_{M}(y_{0})+\frac{1}{2}$ $\left\Vert
\text{X}\right\Vert _{P}^{2}(y_{0})$ $+$ $\frac{1}{2}\operatorname{div}%
X_{P}(y_{0})+$ V$(y_{0})$

\qquad\qquad\qquad\qquad\qquad\qquad\qquad\qquad\qquad\qquad\qquad\qquad
\qquad\qquad\qquad\qquad\qquad\qquad$\blacksquare$

We insert $\left(  10.30\right)  $ in the expression for \textbf{b}$_{1}%
($y$_{0},$P,$\phi)$ in $\left(  10.20\right)  $ and obtain the:

\begin{theorem}
\qquad\textbf{b}$_{1}($y$_{0},$P,$\phi)=\Theta(y_{0})\phi\left(  y_{0}\right)
$
\end{theorem}

$\qquad=\frac{1}{24}[\underset{i=q+1}{\overset{n}{\sum}}3<H,i>^{2}+2(\tau
^{M}-3\tau^{P}\ +\overset{q}{\underset{\text{a=1}}{\sum}}\varrho_{\text{aa}%
}^{M}+\overset{q}{\underset{\text{a,b}=1}{\sum}}R_{\text{abab}}^{M}%
)](y_{0})\phi\left(  y_{0}\right)  $

\qquad\qquad$-\frac{1}{2}[$ $\left\Vert \text{X}\right\Vert _{M}^{2}-$
$\left\Vert \text{X}\right\Vert _{P}^{2}$ $+\operatorname{div}X_{M}-$
$\operatorname{div}X_{P}](y_{0})\phi\left(  y_{0}\right)  +$ V$(y_{0}%
)\phi\left(  y_{0}\right)  $

$\qquad\qquad+$ $\frac{1}{2}\underset{\text{a=1}}{\overset{\text{q}}{\sum}%
}\frac{\partial^{2}\phi}{\partial\text{x}_{\text{a}}^{2}}(y_{0})$ $+$
$\underset{\text{a=1}}{\overset{\text{q}}{\sum}}\Lambda_{\text{a}}(y_{0}%
)\frac{\partial\phi}{\partial x_{\text{a}}}\left(  y_{0}\right)  \ +\frac
{1}{2}$ $\underset{\text{a=1}}{\overset{\text{q}}{\sum}}\Lambda_{\text{a}%
}(y_{0})\Lambda_{\text{a}}(y_{0})\phi\left(  y_{0}\right)  $

\qquad\qquad$+$ $\underset{\text{a=1}}{\overset{\text{q}}{\sum}}X_{\text{a}%
}(y_{0})\frac{\partial\phi}{\partial\text{x}_{\text{a}}}(y_{0})$ +
$\underset{\text{a=1}}{\overset{\text{q}}{\sum}}$ $X_{\text{a}}(y_{0}%
)\Lambda_{\text{a}}(y_{0})\phi(y_{0})+\frac{1}{2}$W$\left(  y_{0}\right)
\phi\left(  y_{0}\right)  $

\qquad\qquad\qquad\qquad\qquad\qquad\qquad\qquad\qquad\qquad\qquad\qquad
\qquad\qquad\qquad\qquad\qquad$\blacksquare$\qquad

\begin{corollary}
When Fermi coordinates reduce to local coordinates, we have:
\end{corollary}

$\qquad$\textbf{b}$_{1}($y$_{0},$P,$\phi)=\Theta(y_{0})\phi\left(
y_{0}\right)  =$ \textbf{b}$_{1}($y$_{0},$P$)\phi\left(  y_{0}\right)
=\Theta(y_{0})\phi\left(  y_{0}\right)  $

$\qquad=\frac{1}{12}[\tau^{M}](y_{0})\phi\left(  y_{0}\right)  +\frac{1}%
{2}W\left(  y_{0}\right)  \phi\left(  y_{0}\right)  -\frac{1}{2}[$ $\left\Vert
\text{X}\right\Vert _{M}^{2}+\operatorname{div}X_{M}](y_{0})\phi\left(
y_{0}\right)  +$ V$(y_{0})\phi\left(  y_{0}\right)  $

\qquad\qquad\qquad\qquad\qquad\qquad\qquad\qquad\qquad\qquad\qquad\qquad
\qquad\qquad\qquad\qquad\qquad

We see that we recover the usual \textbf{first order coeffiicient} of the heat
kernel expansion when we take:

X = 0 and V = 0.

\qquad\qquad\qquad\qquad\qquad\qquad\qquad\qquad\qquad\qquad\qquad\qquad
\qquad\qquad\qquad\qquad\qquad$\blacksquare$

\chapter{Computation of the Third Coefficient\qquad\qquad\qquad\qquad}

We will not attempt to give the expression of the third coefficient b$_{2}%
$(x$,$P) at a general point x$\in$M$_{0}.$ The computation at a general point
x$\in$M proves intractable just as in the case of the ordinary heat kernel
expansion. We will compute it at the center y$_{0}\in P$ of Fermi coordinates
of the submanifold P (which is analogous to expanding along the diagonal in
the case of the Minakshisundaram-Pleijel heat kernel expansion). However
b$_{2}$($y_{0},$P,$\phi$) will still be weirdly too long.

The third coefficient b$_{2}(y_{0},P,\phi)$ defined at the centre of Fermi
coordinates is given in \textbf{Theorem }$4$ (\textbf{Exact Expansion of the
Generalized Heat Kernel)} by:

$\left(  11.1\right)  \qquad$b$_{2}($y$_{0}$,P$,\phi)=\int_{0}^{1}\int%
_{0}^{r_{1}}$F(1,1-r$_{2}$)[L$_{\Psi}$F(1-r$_{2}$,1-r$_{1}$)L$_{\Psi}%
[\phi\circ\pi_{P}$](y$_{0}$)dr$_{1}$dr$_{2}$

We will present it here, expressed in \textbf{geometric invariants}. It is one
of the most \textbf{significant achievements} of this work.$\qquad$

The expression for the third coefficient is long and unwieldly$.$ The
computation at a general point x$\in$M$_{0}$ will be too long and cumbersome
and so to simplify computations we shall compute the coefficient at the center
of Fermi coordinates y$_{0}\in$P. This is analagous to computing the
Minakshisundaram-Pleijel \textbf{heat kernel} expansion coefficients along the
diagonal of the manifold M.

We start by giving some purely \textbf{computational lemmas.}

\qquad\qquad\qquad\qquad\qquad\qquad\qquad\qquad\qquad\qquad\qquad\qquad
\qquad\qquad\qquad\qquad\qquad$\blacksquare$

\begin{lemma}
(Preliminary Version of b$_{2}$(y$_{0},P,\phi)$
\end{lemma}

Setting $\Theta=$ L$_{\Psi}[\phi\circ\pi_{\text{P}}],$ we have:

b$_{2}$(y$_{0},P,\phi)=$ I$_{1}+$ I$_{2}+$ I$_{3}+$ I$_{4}$

$=\frac{1}{2}\frac{\text{L}\Psi}{\Psi}(y_{0})\Theta(y_{0})\qquad$I$_{1}$

$+\frac{1}{4}\underset{\text{a=1}}{\overset{\text{q}}{\sum}}\frac{\partial
^{2}\Theta}{\partial x_{\text{a}}^{2}}(y_{0})+\frac{1}{12}%
\underset{i=q+1}{\overset{n}{\sum}}\frac{\partial^{2}\Theta}{\partial
x_{i}^{2}}(y_{0})+\frac{1}{2}\underset{\text{a=1}}{\overset{\text{q}}{\sum}%
}\Lambda_{\text{a}}(y_{0})\frac{\partial\Theta}{\partial x_{\text{a}}}%
(y_{0})+\frac{1}{4}\underset{i=q+1}{\overset{n}{\sum}}\Lambda_{i}(y_{0}%
)\frac{\partial\Theta}{\partial x_{i}}(y_{0})$

$+$ $\frac{1}{4}\underset{\text{a=1}}{\overset{\text{q}}{\sum}}\Lambda
_{\text{a}}^{2}(y_{0})\Theta(y_{0})-\frac{1}{4}%
\underset{k=q+1}{\overset{n}{\sum}}<H,k>(y_{0})\Lambda_{k}(y_{0})+\frac{1}{4}%
$W$(y_{0})\Theta(y_{0})$

$+\frac{1}{2}$ $\underset{\text{a=1}}{\overset{\text{q}}{\sum}}$X$_{\text{a}%
}(y_{0})\frac{\partial\Theta}{\partial x_{\text{a}}}(y_{0})+\frac{1}{2}$
$\underset{\text{a=1}}{\overset{\text{q}}{\sum}}$X$_{\text{a}}(y_{0}%
)\Lambda_{\text{a}}(y_{0})\Theta(y_{0})+\frac{1}{2}$
$\underset{j=q+1}{\overset{n}{\sum}}$X$_{j}(y_{0})\Lambda_{j}(y_{0}%
)\Theta(y_{0})$

The third expansion coefficient for the generalized vector bundle heat kernel
is given for 1$\geq r_{1}\geq r_{2}$ by:

b$_{2}($y$_{0}$,P$,\phi)=\int_{0}^{1}\int_{0}^{r_{1}}$F(1,1-r$_{2}$)[L$_{\Psi
}$F(1-r$_{2}$,1-r$_{1}$)L$_{\Psi}\phi\circ\pi_{\text{P}}$](y$_{0}$)dr$_{1}%
$dr$_{2}$

\begin{center}
\qquad\qquad\ \ \ \ \ = I$_{1}$ + I$_{2}$ + I$_{3}$
\end{center}

where,

I$_{1}=\int_{0}^{1}\int_{0}^{r_{1}}\frac{\text{L}\Psi}{\Psi}(y_{0})$%
F(1-r$_{2},$1-r$_{1}$)L$_{\Psi}[\phi\circ\pi_{\text{P}}]$(y$_{0}$)dr$_{1}%
$dr$_{2}$

I$_{2}=$ $\int_{0}^{1}\int_{0}^{r_{1}}<\nabla\log\Psi,\nabla\lbrack$%
F(1-r$_{2},$1-r$_{1}$)L$_{\Psi}\phi\circ\pi_{\text{P}}]>(y_{0})$dr$_{1}%
$dr$_{2}$

I$_{3}=$ $\int_{0}^{1}\int_{0}^{r_{1}}$L$[$F(1-r$_{2}$,1-r$_{1}$)L$_{\Psi}%
\phi\circ\pi_{\text{P}}](y_{0})$dr$_{1}$dr$_{2}$

The proofs are purely computational. Recall that F(r,s) is defined on
$\Gamma(E)$ by the formula:

$[$F(r,s)$\phi](x_{0})=\phi\circ\gamma(r-s)$ where $\gamma$ is the unique
minimal geodesic from x$_{0}$ to the point y$_{0}\in$P in time r.

$\left(  11.2\right)  $ \qquad F(1,1-r$_{2}$)[L$_{\Psi}$F(1-r$_{2}$,1-r$_{1}%
$)L$_{\Psi}[\phi\circ\pi_{\text{P}}$](y$_{0}$) = [L$_{\Psi}$F(1-r$_{2}%
$,1-r$_{1}$)L$_{\Psi}\phi\circ\pi_{\text{P}}$]$(\gamma(r_{2}))$

where $\gamma$ is now the unique minimal geodesic from y$_{0}$ to y$_{0}$ time
1. In this case $\gamma$ is a constant geodesic: $\gamma(s)=y_{0}$ for all
s$\in\lbrack0,1].$ Since $r_{2}\in\lbrack0,1]$ we have : $(\gamma
(r_{2}))=y_{0}.$

The integrand becomes:

$\left(  11.3\right)  \qquad$[L$_{\Psi}$F(1-r$_{2}$,1-r$_{1}$)L$_{\Psi}%
[\phi\circ\pi_{\text{P}}$](y$_{0}$)

The definition of L$_{\Psi}$ in $\left(  5.31\right)  $ gives L$_{\Psi}\phi=$
$\frac{\text{L(}\Psi\phi)}{\Psi}$ and since,

\qquad\qquad L$($f$\phi)=\phi$ L(f) + fL($\phi$) + $<\nabla f,\nabla\phi>-$
V$(f\phi)$

for any twice differentiable function f:M$\longrightarrow R$ and any smooth
section $\phi\in\Gamma(E),$ we have:

$\left(  11.4\right)  \qquad$[L$_{\Psi}$F(1-r$_{2}$,1-r$_{1}$)L$_{\Psi}%
\phi\circ\pi_{\text{P}}$](y$_{0}$) = $\frac{1}{\Psi(y_{0})}$[L$\Psi$%
F(1-r$_{2}$,1-r$_{1}$)L$_{\Psi}\phi\circ\pi_{\text{P}}$](y$_{0}$)

\qquad\qquad= $\frac{\text{L}\Psi}{\Psi}(y_{0}).$F(1-r$_{2},$1-r$_{1}%
$)L$_{\Psi}[\phi\circ\pi_{\text{P}}]$(y$_{0}$) + L[F(1-r$_{2}$,1-r$_{1}%
$)L$_{\Psi}\phi\circ\pi_{\text{P}}$](y$_{0}$)

\qquad+ $<\nabla\log\Psi,\nabla\lbrack$F(1-r$_{2},$1-r$_{1}$)L$_{\Psi}%
\phi\circ\pi_{\text{P}}]>(y_{0})-$ V$[$F(1-r$_{2}$,1-r$_{1}$)L$_{\Psi}%
\phi\circ\pi_{\text{P}}](y_{0})$

We label the above expressions as follows:

\qquad\lbrack L$_{\Psi}$F(1-r$_{2}$,1-r$_{1}$)L$_{\Psi}\phi\circ\pi_{\text{P}%
}$](y$_{0}$) = $J_{1}+J_{2}+J_{3}+J_{4}$

where,

$\left(  11.5\right)  \qquad J_{1}=\frac{\text{L}\Psi}{\Psi}(y_{0}%
).$F(1-r$_{2},$1-r$_{1}$)L$_{\Psi}[\phi\circ\pi_{\text{P}}](y_{0})$

\qquad\qquad$\ \ \ J_{2}=$ $<\nabla\log\Psi,\nabla\lbrack$F(1-r$_{2},$%
1-r$_{1}$)L$_{\Psi}\phi\circ\pi_{\text{P}}]>(y_{0})$

$\ \qquad\qquad\ J_{3}=$ L[F(1-r$_{2}$,1-r$_{1}$)L$_{\Psi}\phi\circ
\pi_{\text{P}}$]$(y_{0})$

\qquad\qquad$\ \ J_{4}=-$ V(y$_{0}$)[F(1-r$_{2}$,1-r$_{1}$)L$_{\Psi}\phi
\circ\pi_{\text{P}}$]$(y_{0})$\qquad\qquad

\underline{\textbf{Computation of I}$_{1}$}\qquad\qquad\qquad\qquad

From the previous lemma,

\qquad I$_{1}=\int_{0}^{1}\int_{0}^{r_{1}}J_{1}$dr$_{1}$dr$_{2}=\int_{0}%
^{1}\int_{0}^{r_{1}}\frac{\text{L}\Psi}{\Psi}(y_{0})$F(1-r$_{2},$1-r$_{1}%
$)L$_{\Psi}[\phi\circ\pi_{\text{P}}](y_{0})$dr$_{1}$dr$_{2}$

$\qquad\ \ =\int_{0}^{1}\int_{0}^{r_{1}}\frac{\text{L}\Psi}{\Psi}(y_{0}%
)$L$_{\Psi}[\phi\circ\pi_{\text{P}}](\gamma_{1,2}(r_{1}-r_{2}))$dr$_{1}%
$dr$_{2}$

where by the definition of F(1-r$_{2},$1-r$_{1}$), $\gamma_{1,2}%
:[0,1-r_{2}]\longrightarrow$M$_{0}$ is the unique minimal geodesic from
y$_{0}$ to y$_{0}$ in time 1-r$_{2}$ and so $\gamma_{1,2}$ is a (constant)
geodesic $\gamma_{1,2}(s)=y_{0}$ for all s$\in\lbrack0,1-r_{2}].$ Since
$(r_{1}-r_{2})\in\lbrack0,1-r_{2}],$ we have $\gamma_{1,2}(r_{1}-r_{2}%
)=y_{0}.$ Consequently,

\qquad I$_{1}=\int_{0}^{1}\int_{0}^{r_{1}}\frac{\text{L}\Psi}{\Psi}(y_{0}%
)$L$_{\Psi}[\phi\circ\pi_{\text{P}}](y_{0})$dr$_{1}$dr$_{2}$

The integrand is independent of r$_{1}$ and r$_{2}$ and so integration gives:

$\left(  11.6\right)  \qquad$I$_{1}=\frac{1}{2}\frac{\text{L}\Psi}{\Psi}%
(y_{0})$L$_{\Psi}[\phi\circ\pi_{\text{P}}](y_{0})$ $=\frac{1}{2}\frac
{\text{L}\Psi}{\Psi}(y_{0})\Theta(y_{0})$

where we set $\Theta=$ L$_{\Psi}[\phi\circ\pi_{\text{P}}]$

\underline{\textbf{Computation of I}$_{2}$}

\qquad I$_{2}=$ $\int_{0}^{1}\int_{0}^{r_{1}}J_{2}$dr$_{1}$dr$_{2}=$ $\int%
_{0}^{1}\int_{0}^{r_{1}}<\nabla\log\Psi,\nabla\lbrack$F(1-r$_{2},$1-r$_{1}%
$)L$_{\Psi}\phi\circ\pi_{\text{P}}]>(y_{0})$dr$_{1}$dr$_{2}$

By the definition of the of the operator F(1-r$_{2},$1-r$_{1}$),

F(1-r$_{2},$1-r$_{1}$)L$_{\Psi}[\phi\circ\pi_{\text{P}}](x)=[$L$_{\Psi}%
\phi\circ\pi_{\text{P}}$]$\circ\gamma_{1,2}(r_{1}-r_{2})$

F(1-r$_{2},$1-r$_{1}$)L$_{\Psi}[\phi\circ\pi_{\text{P}}](x)=[$L$_{\Psi}%
\phi\circ\pi_{\text{P}}$]$\circ\gamma_{1,2}(r_{1}-r_{2})$

where $\gamma_{1,2}$ is the unique minimal geodesic from some point $x\in
M_{0}$ to the submanifold P in time $1-r_{2}.$ Assume the geodesic
$\gamma_{1,2}$ meets the submanifold at a point y$\in$ P. Then,

$\gamma_{1,2}(0)=$ $x;$ $\gamma_{1,2}(1-r_{2})=y\in$P.

By the definition of $\gamma_{1,2}$ above, we have in Fermi coordinates:

$\left(  11.7\right)  \qquad z=\gamma_{1,2}(s)=\left(  x_{1},...,x_{q}%
,(1-\frac{s}{1-r_{2}})x_{q+1},...,(1-\frac{s}{1-r_{2}})x_{n}\right)  $

\qquad\qquad\qquad$=y+(1-\frac{s}{1-r_{2}})(x-y)$

$\left(  11.8\right)  \qquad z_{1}=\gamma_{1,2}(r_{1}-r_{2})=\left(
x_{1},...,x_{q},\frac{1-r_{1}}{1-r_{2}}x_{q+1},...,\frac{1-r_{1}}{1-r_{2}%
}x_{n}\right)  $

From the definition of $\gamma_{1,2}(s)$ above, we have:

$\left(  11.9\right)  $\qquad$\frac{\partial}{\partial x_{j}}\gamma
_{1,2}(s)=\left\{
\begin{array}
[c]{c}%
1\text{ for }j=1,...,q\\
(1-\frac{s}{1-r_{2}})\text{ for }j=q+1,...,n
\end{array}
\right.  $

In particular, we have from $\left(  7.20\right)  :$

$\left(  11.10\right)  $\qquad$\frac{\partial}{\partial x_{j}}\gamma
_{1,2}(r_{1}-r_{2})=\left\{
\begin{array}
[c]{c}%
1\text{ for }j=1,...,q\\
\frac{1-r_{1}}{1-r_{2}}\text{ for }j=q+1,...,n
\end{array}
\right.  $

$\left(  11.11\right)  $\qquad$\frac{\partial^{2}}{\partial x_{j}\partial
x_{i}}\gamma_{1,2}(r_{1}-r_{2})=0$ for all $i,j=1,...,q,q+1,...,n.$

From $\left(  11.8\right)  ,$

\qquad$\qquad\qquad\pi_{\text{P}}(z_{1})=\pi_{\text{P}}\circ\gamma_{1,2}%
(r_{1}-r_{2})=\left(  x_{1},...,x_{q}\right)  =y\in P$

Therefore,

$\left(  11.12\right)  $ (i)\qquad$\frac{\partial}{\partial\text{x}_{\text{a}%
}}\pi_{\text{P}}$(z$_{1}$) $=1$ for a=1,...,q

$\left(  11.13\right)  $ (ii)\qquad$\frac{\partial}{\partial\text{x}_{i}}%
\pi_{\text{P}}$(z$_{1}$) $=0$ for $i=q+1,...,n$

$\left(  11.14\right)  $ (ii)\qquad$\frac{\partial^{2}}{\partial\text{x}%
_{i}\partial\text{x}_{j}}\pi_{\text{P}}$(z$_{1}$) $=0$ for $i,j=1,...,n$

From $\left(  11.5\right)  :J_{1}=\frac{\text{L}\Psi}{\Psi}(y_{0}).$%
F(1-r$_{2},$1-r$_{1}$)L$_{\Psi}[\phi\circ\pi_{\text{P}}](y_{0})=\frac
{\text{L}\Psi}{\Psi}(y_{0})\Theta(y_{0})$

We then consider:

$\left(  11.15\right)  \qquad J_{2}=$ $<\nabla\log\Psi,\nabla\lbrack
$F(1-r$_{2},$1-r$_{1}$)L$_{\Psi}\phi\circ\pi_{\text{P}}]>(y_{0})$

By (iii) of \textbf{Theorem 1}, we have for a vector field Y the formula:

(The Einstein convention for summation over repeated indices is used below)

\qquad$<\nabla\phi,Y>$ $=$ $Y_{j}\left(  \frac{\partial}{\partial x_{j}%
}+\Lambda_{j}\right)  \phi=Y_{j}\nabla_{\partial_{j}}\phi$

We now take $Y=$ $\nabla\log\Psi$ and use $\left(  7.3\right)  $ above to have:

$\qquad\left(  \nabla\log\Psi\right)  _{j}(y_{0})=\left(  \nabla\log
\theta^{-\frac{1}{2}}\right)  _{j}(y_{0})+\left(  \nabla\log\Phi\right)
_{j}(y_{0})$

$\qquad\qquad\qquad\ \ \qquad=\frac{1}{2}<H,j>(y_{0})-$ $X_{j}(y_{0})$ for
$j=q+1,...,n.$

By (ix)$^{\ast}$ of \textbf{Table A}$_{9}$ in \textbf{Appendix A},

$\frac{\partial}{\partial x_{i}}(\nabla\log\theta^{-\frac{1}{2}})_{\text{a}%
}(y_{0})=\frac{1}{2}\perp_{\text{a}ki}(y_{0})<H,k>(y_{0})=-\frac{1}{2}%
\perp_{\text{a}ik}(y_{0})<H,k>(y_{0}).$Consequently,

$\ J_{2}=$ $\frac{1}{2}<H,j>(y_{0})\left(  \frac{\partial}{\partial x_{j}%
}+\Lambda_{j}\right)  (y_{0})$F(1-r$_{2},$1-r$_{1}$)$[$L$_{\Psi}\phi\circ
\pi_{\text{P}}](y_{0})$

\qquad$\ -$ $X_{j}(y_{0})\left(  \frac{\partial}{\partial x_{j}}+\Lambda
_{j}\right)  (y_{0})$F(1-r$_{2},$1-r$_{1}$)$[$L$_{\Psi}\phi\circ\pi_{\text{P}%
}](y_{0})$

We set:

$\ J_{21}=$ $\frac{1}{2}<H,j>(y_{0})\Lambda_{j}(y_{0})$F(1-r$_{2},$1-r$_{1}%
$)$[$L$_{\Psi}\phi\circ\pi_{\text{P}}](y_{0})$

$\ J_{22}=$ $\frac{1}{2}<H,j>(y_{0})\frac{\partial}{\partial x_{j}}$%
F(1-r$_{2},$1-r$_{1}$)$[$L$_{\Psi}\phi\circ\pi_{\text{P}}](y_{0})$

$J_{23}=$ $-$ $X_{j}(y_{0})\Lambda_{j}(y_{0})$F(1-r$_{2},$1-r$_{1}$%
)$[$L$_{\Psi}\phi\circ\pi_{\text{P}}](y_{0})$

$J_{24}=$ $-$ $X_{j}(y_{0})\frac{\partial}{\partial x_{j}}$F(1-r$_{2}%
,$1-r$_{1}$)$[$L$_{\Psi}\phi\circ\pi_{\text{P}}](y_{0})$

We compute each of the above expressions:

$\ J_{21}=$ $\frac{1}{2}<H,j>(y_{0})\Lambda_{j}(y_{0})$F(1-r$_{2},$1-r$_{1}%
$)$[$L$_{\Psi}\phi\circ\pi_{\text{P}}](y_{0})$

$=$ $\frac{1}{2}<H,j>(y_{0})\Lambda_{j}(y_{0})$ $[L_{\Psi}\phi\circ
\pi_{\text{P}}](\gamma_{1,2}(r_{1}-r_{2}))$

where, by the definition above, $\gamma_{1,2}$ is now the unique minimal
geodesic from the centre of Fermi coordinates point $y_{0}\in$P to the
submanifold P in time $1-r_{2}.$ It is immediate that, in this case,
$\gamma_{1,2}$ is the constant geodesic: $\gamma_{1,2}(s)=y_{0}$ for all
$s\in\left[  0,1-r_{2}\right]  $ and so, in particular, $(\gamma_{1,2}%
(r_{1}-r_{2}))=y_{0}.$ Consequantly,

$\ J_{21}=$ $\frac{1}{2}<H,j>(y_{0})\Lambda_{j}(y_{0})[$L$_{\Psi}\phi\circ
\pi_{\text{P}}](y_{0})$

We have:

$\left(  11.16\right)  \qquad J_{21}+J_{23}=\frac{1}{2}<H,j>(y_{0})\Lambda
_{j}(y_{0})[$L$_{\Psi}\phi\circ\pi_{\text{P}}](y_{0})$ $-$ $X_{j}%
(y_{0})\Lambda_{j}(y_{0})[$L$_{\Psi}\phi\circ\pi_{\text{P}}](y_{0})$

We next compute the other terms of the expression:

$\ \qquad\qquad J_{22}+J_{24}=$ $\frac{1}{2}<H,j>(y_{0})\frac{\partial
}{\partial x_{j}}$F(1-r$_{2},$1-r$_{1}$)$[$L$_{\Psi}\phi\circ\pi_{\text{P}%
}](y_{0})$

$\qquad\qquad\qquad\ \qquad-$ $X_{j}(y_{0})\frac{\partial}{\partial x_{j}}%
$F(1-r$_{2},$1-r$_{1}$)$[$L$_{\Psi}\phi\circ\pi_{\text{P}}](y_{0})$

Now,

\qquad$\frac{\partial}{\partial x_{j}}$F(1-r$_{2},$1-r$_{1}$)$[$L$_{\Psi}%
\phi\circ\pi_{\text{P}}](y_{0})$ = $\frac{\partial}{\partial x_{j}}[$L$_{\Psi
}\phi\circ\pi_{\text{P}}\circ\gamma_{1,2}$(r$_{1}$-r$_{2}$)$](y_{0})$

\ $\ =\frac{\partial}{\partial x_{j}}[$L$_{\Psi}\phi\circ\pi_{\text{P}%
}](\gamma_{1,2}$(r$_{1}$-r$_{2}$)$)\frac{\partial}{\partial x_{j}}\gamma
_{1,2}$(r$_{1}$-r$_{2}$)

Since $\gamma_{1,2}$(r$_{1}$-r$_{2}$) = y$_{0}$ and $\frac{\partial}{\partial
x_{j}}\gamma_{1,2}$(r$_{1}$-r$_{2}$) = $\frac{1-r_{1}}{1-r_{2}}$ for
$j=q+1,...,n$ by $\left(  7.21\right)  $ above, we have:

$\left(  11.17\right)  $ $\ J_{22}+J_{24}=$ $\frac{1}{2}<H,j>(y_{0}%
)\frac{\partial}{\partial x_{j}}[$L$_{\Psi}\phi\circ\pi_{\text{P}}%
](y_{0})\frac{1-r_{1}}{1-r_{2}}\ -$ $X_{j}(y_{0})\frac{\partial}{\partial
x_{j}}[$L$_{\Psi}\phi\circ\pi_{\text{P}}](y_{0})\frac{1-r_{1}}{1-r_{2}}$

An elementary computation gives:

$\left(  11.18\right)  $ $\int_{0}^{1}\int_{0}^{r_{1}}\frac{1-r_{1}}{1-r_{2}}%
$dr$_{1}$dr$_{2}=\int_{0}^{1}\int_{0}^{r_{1}}\frac{1-r_{1}}{1-r_{2}}$dr$_{2}%
$dr$_{1}=\frac{1}{4}$

Consequently,

$\left(  11.19\right)  $\qquad I$_{2}=$ $\int_{0}^{1}\int_{0}^{r_{1}}%
<\nabla\log\Psi,\nabla\lbrack$F(1-r$_{2},$1-r$_{1}$)L$_{\Psi}[\phi\circ
\pi_{\text{P}}]>(y_{0})$dr$_{1}$dr$_{2}$

$=$ $\frac{1}{2}<H,j>(y_{0})\frac{\partial}{\partial x_{j}}[$L$_{\Psi}%
\phi\circ\pi_{\text{P}}](y_{0})$ $\int_{0}^{1}\int_{0}^{r_{1}}\frac{1-r_{1}%
}{1-r_{2}}$dr$_{1}$dr$_{2}$

$-X_{j}(y_{0})\frac{\partial}{\partial x_{j}}[$L$_{\Psi}\phi\circ\pi
_{\text{P}}](y_{0})\int_{0}^{1}\int_{0}^{r_{1}}\frac{1-r_{1}}{1-r_{2}}$%
dr$_{1}$dr$_{2}$

$=$ $\frac{1}{8}<H,j>(y_{0})\frac{\partial}{\partial x_{j}}[$L$_{\Psi}%
\phi\circ\pi_{\text{P}}](y_{0})$ $-\frac{1}{4}X_{j}(y_{0})\frac{\partial
}{\partial x_{j}}[$L$_{\Psi}\phi\circ\pi_{\text{P}}](y_{0})$

We see that for $j=q+1,...,n,$

$\left(  11.20\right)  $ I$_{2}=$ $\frac{1}{8}<H,j>(y_{0})\frac{\partial
}{\partial x_{j}}[$L$_{\Psi}\phi\circ\pi_{\text{P}}](y_{0})$ $-\frac{1}%
{4}X_{j}(y_{0})\frac{\partial}{\partial x_{j}}[$L$_{\Psi}\phi\circ
\pi_{\text{P}}](y_{0})$

Recalling that L$_{\Psi}[\phi\circ\pi_{\text{P}}]=\Theta,$ we have:

$\left(  11.21\right)  $ \ I$_{2}=$ $\frac{1}{8}%
\underset{j=q+1}{\overset{n}{\sum}}<H,j>$ $\frac{\partial\Theta}{\partial
x_{j}}(y_{0})-\frac{1}{4}\underset{j=q+1}{\overset{n}{\sum}}X_{j}(y_{0}%
)\frac{\partial\Theta}{\partial x_{j}}(y_{0})=$ $\frac{1}{8}%
\underset{j=q+1}{\overset{n}{\sum}}\left[  <H,j>-2X_{j}\right]  (y_{0}%
)\frac{\partial\Theta}{\partial x_{j}}(y_{0})$

\underline{\textbf{Computation of I}$_{3}$}\qquad\qquad\qquad\qquad
\qquad\qquad\qquad\qquad

I$_{3}=$ $\int_{0}^{1}\int_{0}^{r_{1}}J_{3}$dr$_{1}$dr$_{2}=$ $\int_{0}%
^{1}\int_{0}^{r_{1}}$L$[$F(1-r$_{2}$,1-r$_{1}$)L$_{\Psi}\phi\circ\pi
_{\text{P}}](y_{0})$dr$_{1}$dr$_{2}$

We simplify the integrand L$[$F(1-r$_{2}$,1-r$_{1}$)L$_{\Psi}\phi\circ
\pi_{\text{P}}](y_{0}):$

We first consider the expression:

$\left(  11.22\right)  \qquad$F(1-r$_{2}$,1-r$_{1}$)L$_{\Psi}\phi\circ
\pi_{\text{P}}\left(  x\right)  $ $=$ L$_{\Psi}[\phi\circ\pi_{\text{P}}%
]$($\gamma_{1,2}$(r$_{1}$-r$_{2}$)) $=$ L$_{\Psi}[\phi\circ\pi_{\text{P}%
}](z_{1})$ $=\Theta(z_{1})$

where z$_{1}=\gamma_{1,2}$(r$_{1}$-r$_{2}$) and $\gamma_{1,2}:$ $[0,$%
1-r$_{2}]\longrightarrow$ M$_{0}$ \ is the unique minimal geodesic from a
point $x\in$M$_{0}$ to the submanifold P in time 1-r$_{2}.$ We assume it meets
P at some point y$\in U\subset$P in time 1-r$_{2}$, where $U$ is a small
neighbourhood on P of the centre of Ferm coordinates (see definition of Fermi coordinates).

We see that L$_{\Psi}[\phi\circ\pi_{\text{P}}](z_{1})$ $=\Theta(z_{1})$ here
is the same as $\Theta(z_{0})$ of $\left(  C_{16}\right)  $ in Appendix C,
with $z_{0}=\gamma($r$_{1})$ there replaced by z$_{1}=\gamma_{1,2}$(r$_{1}%
$-r$_{2}$) here. From $\left(  C_{16}\right)  $ of Appendix C, we have:

$\left(  11.23\right)  \qquad\Theta(z_{1})=$ L$_{\Psi}$[$\phi\circ
\pi_{\text{P}}$]$(z_{1})$

$\qquad=$ $\frac{\text{L}\Psi}{\Psi}(z_{1})\phi\circ\pi_{\text{P}}(z_{1})+$
L[$\phi\circ\pi_{\text{P}}$]$(z_{1})+$ $<\nabla\log\Psi,\nabla\phi\circ
\pi_{\text{P}}>(z_{1})-$ V$[\phi\circ\pi_{\text{P}}](z_{1})$

\qquad\textbf{ }$=\frac{\text{L}\Psi}{\Psi}\left(  z_{1}\right)  \phi\circ
\pi_{\text{P}}(z_{1})+$ $\frac{1}{2}\Delta$[$\phi\circ\pi_{\text{P}}$%
]$(z_{1})+$ $<\triangledown$L[$\phi\circ\pi_{\text{P}}$]$(z_{1})+$
L[$\phi\circ\pi_{\text{P}}$]$(z_{1})$\qquad\qquad

\qquad+ \ $\frac{1}{2}\underset{\text{a,b=1}}{\overset{\text{q}}{\sum}}%
$g$^{\text{ab}}(z_{1})[\frac{\partial^{2}\phi}{\partial\text{x}_{\text{a}%
}\partial\text{x}_{\text{b}}}\circ\pi_{\text{P}}(z_{1})$ $]+\frac{1}%
{2}\underset{i,j=1}{\overset{n}{\sum}}$g$^{ij}(z_{1})[\frac{\partial
\Lambda_{j}}{\partial\text{x}_{i}}(z_{1})\phi\circ\pi_{\text{P}}(z_{1})]$

$\qquad+\frac{1}{2}\underset{j=1}{\overset{n}{\sum}}\underset{\text{a=1}%
}{\overset{\text{q}}{\sum}}$g$^{\text{a}j}(z_{1})[$ $\Lambda_{j}(z_{1}%
)\frac{\partial\phi}{\partial\text{x}_{\text{a}}}\circ\pi_{\text{P}}%
(z_{1})]+\frac{1}{2}\underset{i=1}{\overset{n}{\sum}}\underset{\text{b=1}%
}{\overset{\text{q}}{\sum}}$g$^{i\text{b}}(z_{1})[\Lambda_{i}$(z$_{1}$%
)$\frac{\partial\phi}{\partial\text{x}_{\text{b}}}\circ\pi_{\text{P}}(z_{1})]$

$\qquad+\frac{1}{2}\underset{i,j=1}{\overset{n}{\sum}}$g$^{ij}(z_{1}%
)\Lambda_{i}$(z$_{1}$)$\Lambda_{j}$(z$_{1}$)$\phi\circ\pi_{\text{P}}(z_{1})$

$\qquad-\frac{1}{2}\underset{i,j=1}{\overset{n}{\sum}}$g$^{ij}(z_{1}%
)[\underset{\text{c=1}}{\overset{\text{q}}{\sum}}\Gamma_{ij}^{\text{c}}%
(z_{1})\frac{\partial\phi}{\partial\text{x}_{\text{c}}}\circ\pi_{\text{P}%
}(z_{1})$ + $\underset{k=1}{\overset{n}{\sum}}\Gamma_{ij}^{k}(z_{1}%
)\Lambda_{k}(z_{1})\phi\circ\pi_{\text{P}}(z_{1})]$

$\qquad\ +\frac{1}{2}$W$(z_{1})\phi\circ\pi_{\text{P}}(z_{1})$

\qquad\ \ + $\underset{\text{a=1}}{\overset{\text{q}}{\sum}}(\nabla\log
\Psi)_{\text{a}}(z_{1})\frac{\partial\phi}{\partial\text{x}_{\text{a}}}%
\circ\pi_{\text{P}}(z_{1})$ + $\underset{\text{j=1}}{\overset{\text{n}}{\sum}%
}$ $(\nabla\log\Psi)_{j}(z_{1})\Lambda_{j}(z_{1})\phi\circ\pi_{\text{P}}%
(z_{1}).$

\qquad+ $\underset{\text{a=1}}{\overset{\text{q}}{\sum}}$X$_{\text{a}}%
(z_{1})\frac{\partial\phi}{\partial\text{x}_{\text{a}}}\circ\pi_{\text{P}%
}(z_{1})$ + $\underset{j=1}{\overset{n}{\sum}}$ X$_{j}(z_{1})\Lambda_{j}%
(z_{1})\phi\circ\pi_{\text{P}}(z_{1})$

\qquad\qquad\qquad\qquad\qquad\qquad\qquad\qquad\qquad\qquad\qquad\qquad
\qquad\qquad\qquad\qquad$\blacksquare$

Since $\pi_{\text{P}}(z_{1})=\pi_{\text{P}}\circ\gamma_{1,2}(r_{1}%
-r_{2})=\left(  x_{1},...,x_{q}\right)  =y$ by $\left(  11.8\right)  $ above,
we have:

$\left(  11.24\right)  \qquad\qquad\Theta(z_{0})=$ L$_{\Psi}[\phi\circ
\pi_{\text{P}}](z_{1})$

\qquad\qquad\textbf{ }$=\frac{\text{L}\Psi}{\Psi}\left(  z_{1}\right)
\phi(y)$\qquad\qquad

\qquad\qquad$+$\ $\frac{1}{2}\underset{\text{a,b=1}}{\overset{\text{q}}{\sum}%
}$g$^{\text{ab}}(z_{1})\left\{  \frac{\partial^{2}\phi}{\partial
\text{x}_{\text{a}}\partial\text{x}_{\text{b}}}(y)\text{ }\right\}  +\frac
{1}{2}\underset{i,j=1}{\overset{n}{\sum}}$g$^{ij}(z_{0})\left\{
\frac{\partial\Lambda_{j}}{\partial\text{x}_{i}}(z_{1})\phi(y)\right\}  $

\qquad\qquad$+\frac{1}{2}\underset{j=1}{\overset{n}{\sum}}\underset{\text{a=1}%
}{\overset{\text{q}}{\sum}}$g$^{\text{a}j}(z_{1})\left\{  \text{ }\Lambda
_{j}(z_{1})\frac{\partial\phi}{\partial\text{x}_{\text{a}}}(y)\right\}
+\frac{1}{2}\underset{i=1}{\overset{n}{\sum}}\underset{\text{b=1}%
}{\overset{\text{q}}{\sum}}$g$^{i\text{b}}(z_{1})\left\{  \Lambda_{i}%
\text{(z}_{1}\text{)}\frac{\partial\phi}{\partial\text{x}_{\text{b}}}(y)\text{
}\right\}  $

$\qquad\qquad+\frac{1}{2}\underset{i,j=1}{\overset{n}{\sum}}$g$^{ij}%
(z_{1})\Lambda_{i}$(z$_{1}$)$\Lambda_{j}$(z$_{1}$)$\phi(y)$

\qquad$\qquad-\frac{1}{2}\underset{i,j=1}{\overset{n}{\sum}}$g$^{ij}%
(z_{1})\left\{  \Gamma_{ij}^{\text{c}}(z_{1})\frac{\partial\phi}%
{\partial\text{x}_{\text{c}}}(y)\text{ + }\Gamma_{ij}^{k}(z_{0})\Lambda
_{k}(z_{1})\phi(y)\right\}  \ +\frac{1}{2}$W$(z_{1})\phi(y)$

\qquad\qquad$+$ $\underset{\text{a=1}}{\overset{\text{q}}{\sum}}(\nabla
^{0}\log\Psi)_{\text{a}}(z_{1})\frac{\partial\phi}{\partial\text{x}_{\text{a}%
}}(y)$ + $\underset{j=1}{\overset{n}{\sum}}$ $(\nabla^{0}\log\Psi)_{j}%
(z_{1})\Lambda_{j}(z_{1})\phi(y)$

\qquad\qquad$+$ $\underset{\text{a=1}}{\overset{\text{q}}{\sum}}$X$_{\text{a}%
}(z_{1})\frac{\partial\phi}{\partial\text{x}_{\text{a}}}(y)+$
$\underset{\text{j=1}}{\overset{\text{n}}{\sum}}$ X$_{j}(z_{1})\Lambda
_{j}(z_{1})\phi(y)+$ V$(z_{1})\phi(y)-$ V(z$_{1}$)$\phi(y)$

Recall that $\Theta=$ L$_{\Psi}$[$\phi\circ\pi_{P}$] and that z$_{1}%
=\gamma_{1,2}$(r$_{1}$-r$_{2}$) where $\gamma_{1,2}:[0,$1-r$_{2}%
]\longrightarrow$ M$_{0}$ is the unique minimal geodesic from some point
x$\in$M$_{0}$ to y$_{0}\in$P in time 1-r$_{2}.$ With this notation in mind,

$\left(  11.25\right)  \qquad$I$_{3}=$ $\int_{0}^{1}\int_{0}^{r_{1}}$%
L$[$F(1-r$_{2}$,1-r$_{1}$)L$_{\Psi}\phi\circ\pi_{P}](y_{0})$dr$_{1}$dr$_{2}$

$\ \ \qquad\qquad=$ $\int_{0}^{1}\int_{0}^{r_{1}}$L$[$L$_{\Psi}[\phi\circ
\pi_{P}\circ\gamma_{1,2}($r$_{1}$-r$_{2}$)$](y_{0})$dr$_{1}$dr$_{2}$

$\ \qquad\qquad=\int_{0}^{1}\int_{0}^{r_{1}}$L$[\Theta\circ\gamma_{1,2}%
($r$_{1}$-r$_{2}$)$](y_{0})$dr$_{1}$dr$_{2}=\int_{0}^{1}\int_{0}^{r_{1}}%
$L$[\Theta(z_{1})](y_{0})$dr$_{1}$dr$_{2}$

where z$_{1}=\gamma_{1,2}$(r$_{1}$-r$_{2}$). We compute L$[\Theta(z_{1})]$
with full details:

The operator L = $\frac{1}{2}\Delta+$ X $+$ V where $\Delta$ is the
generalized Laplacian (Laplace-Type operator) on sections of the vector bundle
E is given by (v) of \textbf{Proposition 7:} $\qquad\qquad\frac{1}{2}%
\Delta=\frac{1}{2}g^{ij}\left\{  \frac{\partial^{2}}{\partial x_{i}%
\partial_{j}}+\frac{\partial\Lambda_{j}}{\partial x_{i}}+\Lambda_{j}%
\frac{\partial}{\partial x_{i}}+\Lambda_{i}\frac{\partial}{\partial x_{j}%
}+\Lambda_{i}\Lambda_{j}-\Gamma_{ij}^{k}\left(  \frac{\partial}{\partial
x_{k}}+\Lambda_{k}\right)  \right\}  +\frac{1}{2}$W

The integrand in $\left(  11.25\right)  $ is then expressed as:

$\left(  11.26\right)  \qquad$L$\Theta(z_{1})=\frac{1}{2}\Delta\Theta(z_{1})+$
$<$X$,\nabla\Theta>(z_{1})+$ [V$\Theta](z_{1})$

We compute each term on the RHS side of $\left(  11.26\right)  $ above. We
start with:

$\left(  11.27\right)  \qquad\frac{1}{2}\Delta\Theta(z_{1})=\frac{1}{2}%
g^{ij}(z_{1})[\left\{  \frac{\partial^{2}}{\partial x_{i}\partial_{j}}%
+\frac{\partial\Lambda_{j}}{\partial x_{i}}+\Lambda_{j}\frac{\partial
}{\partial x_{i}}+\Lambda_{i}\frac{\partial}{\partial x_{j}}+\Lambda
_{i}\Lambda_{j}-\Gamma_{ij}^{k}\left(  \frac{\partial}{\partial x_{k}}%
+\Lambda_{k}\right)  \right\}  \Theta(z_{1})$

$\qquad\qquad\ \ \qquad+\frac{1}{2}$W$(z_{1})\Theta(z_{1})$

where $z_{1}=\gamma_{1,2}(r_{1}-r_{2})$

We have from $\left(  11.10\right)  :$

$\left(  11.28\right)  \qquad\frac{\partial}{\partial x_{j}}\left[
\Theta(z_{1})\right]  =$ $\left\{
\begin{array}
[c]{c}%
\frac{\partial\Theta}{\partial x_{\text{b}}}(z_{1})\text{ for b}=1,...,q\\
\frac{\partial\Theta}{\partial x_{j}}(z_{1})\frac{1-r_{1}}{1-r_{2}}\text{ for
}j=q+1,...,n
\end{array}
\right.  $

From the above, we have the following second derivatives: for a,b =1,...,q and
$i,j=q+1,...,n:$

$\left(  11.29\right)  \qquad\frac{\partial^{2}}{\partial x_{\text{a}}\partial
x_{\text{b}}}\left[  \Theta\circ\gamma_{1,2}(r_{1}-r_{2})\right]  =$
$\frac{\partial^{2}\Theta}{\partial x_{\text{a}}\partial x_{\text{b}}}\left[
\gamma_{1,2}(r_{1}-r_{2})\right]  =$ $\frac{\partial^{2}\Theta}{\partial
x_{\text{a}}\partial x_{\text{b}}}(z_{1})$

\qquad$\frac{\partial^{2}}{\partial x_{i}\partial x_{\text{b}}}\left[
\Theta\circ\gamma_{1,2}(r_{1}-r_{2})\right]  =$ $\frac{\partial^{2}\Theta
}{\partial x_{i}\partial x_{\text{b}}}\left[  \gamma_{1,2}(r_{1}%
-r_{2})\right]  (\frac{1-r_{1}}{1-r_{2}})=$ $\frac{\partial^{2}\Theta
}{\partial x_{i}\partial x_{\text{b}}}(z_{1})(\frac{1-r_{1}}{1-r_{2}})$

\qquad$\frac{\partial^{2}}{\partial x_{\text{a}}\partial x_{j}}\left[
\Theta\circ\gamma_{1,2}(r_{1}-r_{2})\right]  =$ $\frac{\partial^{2}\Theta
}{\partial x_{\text{a}}\partial x_{j}}\left[  \gamma_{1,2}(r_{1}%
-r_{2})\right]  (\frac{1-r_{1}}{1-r_{2}})=$ $\frac{\partial^{2}\Theta
}{\partial x_{\text{a}}\partial x_{j}}(z_{1})(\frac{1-r_{1}}{1-r_{2}})$

\qquad$\frac{\partial^{2}}{\partial x_{i}\partial x_{j}}\left[  \Theta
\circ\gamma_{1,2}(r_{1}-r_{2})\right]  =$ $\frac{\partial^{2}\Theta}{\partial
x_{i}\partial x_{j}}\left[  \gamma_{1,2}(r_{1}-r_{2})\right]  (\frac{1-r_{1}%
}{1-r_{2}})^{2}=$ $\frac{\partial^{2}\Theta}{\partial x_{i}\partial x_{j}%
}(z_{1})(\frac{1-r_{1}}{1-r_{2}})^{2}$

Therefore by $\left(  11.27\right)  ,$ $\left(  11.28\right)  $ and $\left(
11.29\right)  $ we have for a,b = 1,...,q and $i,j=q+1,...,n:$

$\left(  11.30\right)  $ $\frac{1}{2}\Delta\lbrack\Theta(z_{1})]=\frac{1}%
{2}\underset{\text{a,b=1}}{\overset{\text{q}}{\sum}}g^{\text{ab}}(z_{1}%
)\frac{\partial^{2}\Theta}{\partial x_{\text{a}}\partial x_{\text{b}}}%
(z_{1})+\frac{1}{2}\underset{i=q+1}{\overset{n}{\sum}}\underset{\text{b=1}%
}{\overset{\text{q}}{\sum}}g^{i\text{b}}(z_{1})\frac{\partial^{2}\Theta
}{\partial x_{i}\partial x_{\text{b}}}(z_{1})(\frac{1-r_{1}}{1-r_{2}})$

$+\frac{1}{2}\underset{j=q+1}{\overset{n}{\sum}}\underset{\text{a=1}%
}{\overset{\text{q}}{\sum}}g^{\text{a}j}(z_{1})\frac{\partial^{2}\Theta
}{\partial x_{\text{a}}\partial x_{j}}(z_{1})(\frac{1-r_{1}}{1-r_{2}}%
)+\frac{1}{2}\underset{i,j=q+1}{\overset{n}{\sum}}$g$^{ij}(z_{1}%
)\frac{\partial^{2}\Theta}{\partial x_{i}\partial x_{j}}(z_{1})(\frac{1-r_{1}%
}{1-r_{2}})^{2}$

$+\frac{1}{2}\underset{i,j=1}{\overset{n}{\sum}}g^{ij}(z_{1})(\frac
{\partial\Lambda_{j}}{\partial x_{i}}\Theta)(z_{1})+\frac{1}{2}$
$\underset{j=1}{\overset{n}{\sum}}\underset{\text{a=1}}{\overset{\text{q}%
}{\sum}}g^{\text{a}j}(z_{1})(\Lambda_{j}\frac{\partial\Theta}{\partial
x_{\text{a}}})(z_{1})$

$+\frac{1}{2}\underset{j=1}{\overset{n}{\sum}}\underset{\text{i}%
=q+1}{\overset{n}{\sum}}g^{ij}(z_{1})(\Lambda_{j}\frac{\partial\Theta
}{\partial x_{i}})(z_{1})\frac{1-r_{1}}{1-r_{2}}+\frac{1}{2}$
$\underset{j=1}{\overset{n}{\sum}}\underset{\text{a=1}}{\overset{\text{q}%
}{\sum}}g^{i\text{b}}(z_{1})(\Lambda_{i}\frac{\partial\Theta}{\partial
x_{\text{b}}})(z_{1})$

$+\frac{1}{2}\underset{i=1}{\overset{n}{\sum}}%
\underset{j=q+1}{\overset{n}{\sum}}g^{ij}(z_{1})(\Lambda_{i}\frac
{\partial\Theta}{\partial x_{j}})(z_{1})\frac{1-r_{1}}{1-r_{2}}+\frac{1}%
{2}\underset{i,j=1}{\overset{n}{\sum}}g^{ij}(z_{1})(\Lambda_{i}\Lambda
_{j})(z_{1})\Theta(z_{1})$

$-\frac{1}{2}\underset{i,j=1}{\overset{n}{\sum}}g^{ij}(z_{1}%
)\underset{\text{a}=1}{\overset{\text{q}}{\sum}}\left\{  \Gamma_{ij}%
^{\text{a}}(z_{1})\left(  \frac{\partial\Theta}{\partial x_{\text{a}}}%
(z_{1})+\Lambda_{\text{a}}(z_{1})\Theta(z_{1})\right)  \right\}  $

$-\frac{1}{2}\underset{i,j=1}{\overset{n}{\sum}}g^{ij}(z_{1}%
)\underset{k=q+1}{\overset{n}{\sum}}\left\{  \Gamma_{ij}^{k}(z_{1})\left(
\frac{\partial\Theta}{\partial x_{k}}(z_{1})\frac{1-r_{1}}{1-r_{2}}%
+\Lambda_{k}(z_{1})\Theta(z_{1})\right)  \right\}  +\frac{1}{2}$%
W$(z_{1})\Theta(z_{1})$

Since $g^{\text{a}j}(z_{1})=g^{j\text{a}}(z_{1})$ and $\frac{\partial
^{2}\Theta}{\partial x_{\text{a}}\partial x_{j}}=\frac{\partial^{2}\Theta
}{\partial x_{j}\partial x_{\text{a}}},$ we have the following equalities:

$\underset{i=q+1}{\overset{n}{\sum}}\underset{\text{b=1}}{\overset{\text{q}%
}{\sum}}g^{i\text{b}}(z_{1})\frac{\partial^{2}\Theta}{\partial x_{i}\partial
x_{\text{b}}}(z_{1})(\frac{1-r_{1}}{1-r_{2}}%
)=\underset{j=q+1}{\overset{n}{\sum}}\underset{\text{a=1}}{\overset{\text{q}%
}{\sum}}g^{\text{a}j}(z_{1})\frac{\partial^{2}\Theta}{\partial x_{\text{a}%
}\partial x_{j}}(z_{1})(\frac{1-r_{1}}{1-r_{2}})$

$\underset{j=1}{\overset{n}{\sum}}\underset{\text{a=1}}{\overset{\text{q}%
}{\sum}}g^{i\text{b}}(z_{1})(\Lambda_{i}\frac{\partial\Theta}{\partial
x_{\text{b}}})(z_{1})=$ $\underset{j=1}{\overset{n}{\sum}}\underset{\text{a=1}%
}{\overset{\text{q}}{\sum}}g^{\text{a}j}(z_{1})(\Lambda_{j}\frac
{\partial\Theta}{\partial x_{\text{a}}})(z_{1})$

$\underset{j=1}{\overset{n}{\sum}}\underset{i=q+1}{\overset{n}{\sum}}%
g^{ij}(z_{1})(\Lambda_{j}\frac{\partial\Theta}{\partial x_{i}})(z_{1}%
)\frac{1-r_{1}}{1-r_{2}}=$ $\underset{i=1}{\overset{n}{\sum}}%
\underset{j=q+1}{\overset{n}{\sum}}g^{ij}(z_{1})(\Lambda_{i}\frac
{\partial\Theta}{\partial x_{j}})(z_{1})\frac{1-r_{1}}{1-r_{2}}$

Then the expression in $\left(  11.30\right)  $ becomes:

$\left(  11.31\right)  \qquad\frac{1}{2}\Delta\lbrack\Theta(z_{1})]=\frac
{1}{2}\underset{\text{a,b=1}}{\overset{\text{q}}{\sum}}g^{\text{ab}}%
(z_{1})\frac{\partial^{2}\Theta}{\partial x_{\text{a}}\partial x_{\text{b}}%
}(z_{1})+\underset{j=q+1}{\overset{n}{\sum}}\underset{\text{a=1}%
}{\overset{\text{q}}{\sum}}g^{\text{a}j}(z_{1})\frac{\partial^{2}\Theta
}{\partial x_{\text{a}}\partial x_{j}}(z_{1})(\frac{1-r_{1}}{1-r_{2}})$

$+\frac{1}{2}\underset{i,j=q+1}{\overset{n}{\sum}}$g$^{ij}(z_{1}%
)\frac{\partial^{2}\Theta}{\partial x_{i}\partial x_{j}}(z_{1})(\frac{1-r_{1}%
}{1-r_{2}})^{2}\qquad$

$+\frac{1}{2}\underset{i,j=1}{\overset{n}{\sum}}g^{ij}(z_{1})(\frac
{\partial\Lambda_{j}}{\partial x_{i}}\Theta)(z_{1})+$
$\underset{j=1}{\overset{n}{\sum}}\underset{\text{a=1}}{\overset{\text{q}%
}{\sum}}g^{\text{a}j}(z_{1})(\Lambda_{j}\frac{\partial\Theta}{\partial
x_{\text{a}}})(z_{1})$

$+\frac{1}{2}\underset{i=1}{\overset{n}{\sum}}%
\underset{j=q+1}{\overset{n}{\sum}}g^{ij}(z_{1})(\Lambda_{i}\frac
{\partial\Theta}{\partial x_{j}})(z_{1})\frac{1-r_{1}}{1-r_{2}}+\frac{1}%
{2}\underset{i,j=1}{\overset{n}{\sum}}g^{ij}(z_{1})(\Lambda_{i}\Lambda
_{j})(z_{1})\Theta(z_{1})$

$-\frac{1}{2}\underset{i,j=1}{\overset{n}{\sum}}g^{ij}(z_{1}%
)\underset{\text{a}=1}{\overset{\text{q}}{\sum}}\left\{  \Gamma_{ij}%
^{\text{a}}(z_{1})\left(  \frac{\partial\Theta}{\partial x_{\text{a}}}%
(z_{1})+\Lambda_{\text{a}}(z_{1})\Theta(z_{1})\right)  \right\}  $

$-\frac{1}{2}\underset{i,j=1}{\overset{n}{\sum}}g^{ij}(z_{1}%
)\underset{k=q+1}{\overset{n}{\sum}}\left\{  \Gamma_{ij}^{k}(z_{1})\left(
\frac{\partial\Theta}{\partial x_{k}}(z_{1})\frac{1-r_{1}}{1-r_{2}}%
+\Lambda_{k}(z_{1})\Theta(z_{1})\right)  \right\}  $

\begin{center}
$+\frac{1}{2}$W$(z_{1})\Theta(z_{1})$
\end{center}

We use the fact that z$_{1}=\gamma_{1,2}(r_{1}-r_{2})$ and have:

$\left(  11.32\right)  \qquad\frac{1}{2}\Delta\lbrack\Theta(z_{1})](y_{0})$

$=$ F(1-r$_{2}$,1-r$_{1}$)$[\frac{1}{2}\underset{\text{a,b=1}%
}{\overset{\text{q}}{\sum}}g^{\text{ab}}\frac{\partial^{2}\Theta}{\partial
x_{\text{a}}\partial x_{\text{b}}}+\underset{\text{j}=q+1}{\overset{n}{\sum}%
}\underset{\text{a=1}}{\overset{\text{q}}{\sum}}g^{\text{a}j}\frac
{\partial^{2}\Theta}{\partial x_{\text{a}}\partial x_{j}}(\frac{1-r_{1}%
}{1-r_{2}})+\frac{1}{2}\underset{i,j=q+1}{\overset{n}{\sum}}$g$^{ij}%
\frac{\partial^{2}\Theta}{\partial x_{i}\partial x_{j}}(\frac{1-r_{1}}%
{1-r_{2}})^{2}\qquad$

$+\frac{1}{2}\underset{i,j=1}{\overset{n}{\sum}}g^{ij}(\frac{\partial
\Lambda_{j}}{\partial x_{i}}\Theta)+$ $\underset{j=1}{\overset{n}{\sum}%
}\underset{\text{a=1}}{\overset{\text{q}}{\sum}}g^{\text{a}j}(\Lambda_{j}%
\frac{\partial\Theta}{\partial x_{\text{a}}})+$
$\underset{i=1}{\overset{n}{\sum}}\underset{j=q+1}{\overset{n}{\sum}}%
g^{ij}(\Lambda_{i}\frac{\partial\Theta}{\partial x_{j}})\frac{1-r_{1}}%
{1-r_{2}}+\frac{1}{2}\underset{i,j=1}{\overset{n}{\sum}}g^{ij}(\Lambda
_{i}\Lambda_{j})\Theta$

$-\frac{1}{2}\underset{i,j=1}{\overset{n}{\sum}}g^{ij}\underset{\text{a}%
=1}{\overset{\text{q}}{\sum}}\left\{  \Gamma_{ij}^{\text{a}}\left(
\frac{\partial\Theta}{\partial x_{\text{a}}}+\Lambda_{\text{a}}\right)
\right\}  -\frac{1}{2}\underset{i,j=1}{\overset{n}{\sum}}g^{ij}%
\underset{k=q+1}{\overset{n}{\sum}}\left\{  \Gamma_{ij}^{k}\left(
\frac{\partial\Theta}{\partial x_{k}}\frac{1-r_{1}}{1-r_{2}}+\Lambda_{k}%
\Theta\right)  \right\}  +\frac{1}{2}W\Theta](y_{0})$

$=$ $[\frac{1}{2}\underset{\text{a,b=1}}{\overset{\text{q}}{\sum}}%
g^{\text{ab}}\frac{\partial^{2}\Theta}{\partial x_{\text{a}}\partial
x_{\text{b}}}+\underset{j=q+1}{\overset{n}{\sum}}\underset{\text{a=1}%
}{\overset{\text{q}}{\sum}}g^{\text{aj}}\frac{\partial^{2}\Theta}{\partial
x_{\text{a}}\partial x_{j}}(\frac{1-r_{1}}{1-r_{2}})+\frac{1}{2}%
\underset{i,j=q+1}{\overset{n}{\sum}}$g$^{ij}\frac{\partial^{2}\Theta
}{\partial x_{i}\partial x_{j}}(\frac{1-r_{1}}{1-r_{2}})^{2}\qquad$

$+\frac{1}{2}\underset{i,j=1}{\overset{n}{\sum}}g^{ij}(\frac{\partial
\Lambda_{j}}{\partial x_{i}}\Theta)+$ $\underset{j=1}{\overset{n}{\sum}%
}\underset{\text{a=1}}{\overset{\text{q}}{\sum}}g^{\text{a}j}(\Lambda_{j}%
\frac{\partial\Theta}{\partial x_{\text{a}}})+$
$\underset{i=1}{\overset{n}{\sum}}\underset{j=q+1}{\overset{n}{\sum}}%
g^{ij}(\Lambda_{j}\frac{\partial\Theta}{\partial x_{i}})\frac{1-r_{1}}%
{1-r_{2}}+\frac{1}{2}\underset{i,j=1}{\overset{n}{\sum}}g^{ij}(\Lambda
_{i}\Lambda_{j})\Theta$

$-\frac{1}{2}\underset{i,j=1}{\overset{n}{\sum}}g^{ij}\underset{\text{a}%
=1}{\overset{\text{q}}{\sum}}\left\{  \Gamma_{ij}^{\text{a}}\left(
\frac{\partial\Theta}{\partial x_{\text{a}}}+\Lambda_{\text{a}}\Theta
(z_{1})\right)  \right\}  -\frac{1}{2}\underset{i,j=1}{\overset{n}{\sum}%
}g^{ij}\underset{k=q+1}{\overset{n}{\sum}}\left\{  \Gamma_{ij}^{k}\left(
\frac{\partial\Theta}{\partial x_{k}}\frac{1-r_{1}}{1-r_{2}}+\Lambda_{k}%
\Theta\right)  \right\}  $

$+\frac{1}{2}$W$\Theta](\gamma_{1,2}(r_{1}-r_{2}))$

where $\gamma_{1,2}:[0,1-r_{2}]\longrightarrow$M$_{0}$ is now the unique
minimal geodesic from y$_{0}$ to y$_{0}$ in time $1-r_{2}.$ Therefore it is
the constant geodesic $\gamma_{1,2}(s)=y_{0}$ for all $s\in\lbrack0,1-r_{2}].$
Since $(r_{1}-r_{2})\in\lbrack0,1-r_{2}],$we have: $\gamma_{1,2}(r_{1}%
-r_{2})=y_{0}$

Since $g^{\text{ab}}(y_{0})=\delta^{\text{ab}}$ for a,b = 1,...,q$;g^{ij}%
(y_{0})=\delta^{ij};g^{\text{aj}}(y_{0})=$ $\delta^{\text{a}j}=0$ for a =
1,...,q and $j=q+1,...,n,$ the last expression in $\left(  11.32\right)  $ becomes:

$\left(  11.33\right)  \qquad\frac{1}{2}\Delta\lbrack\Theta(z_{1}%
)](y_{0})=\frac{1}{2}\underset{\text{a=1}}{\overset{\text{q}}{\sum}}%
\frac{\partial^{2}\Theta}{\partial x_{\text{a}}^{2}}(y_{0})+\frac{1}%
{2}\underset{i=q+1}{\overset{n}{\sum}}\frac{\partial^{2}\Theta}{\partial
x_{i}^{2}}(y_{0})(\frac{1-r_{1}}{1-r_{2}})^{2}+\frac{1}{2}%
\underset{i=1}{\overset{n}{\sum}}\frac{\partial\Lambda_{i}}{\partial x_{i}%
}(y_{0}).\Theta(y_{0})$

$+\underset{\text{a=1}}{\overset{\text{q}}{\sum}}\Lambda_{\text{a}}%
(y_{0})\frac{\partial\Theta}{\partial x_{\text{a}}}(y_{0})+$
$\underset{i=q+1}{\overset{n}{\sum}}\Lambda_{i}(y_{0})\frac{\partial\Theta
}{\partial x_{i}}(y_{0})\frac{1-r_{1}}{1-r_{2}}+$ $\frac{1}{2}%
\underset{i=1}{\overset{n}{\sum}}\Lambda_{i}^{2}(y_{0})\Theta(y_{0})$

$-\frac{1}{2}\underset{i=1}{\overset{n}{\sum}}\underset{\text{a}%
=1}{\overset{\text{q}}{\sum}}\left\{  \Gamma_{ii}^{\text{a}}(y_{0})\left(
\frac{\partial\Theta}{\partial x_{\text{a}}}(y_{0})+\Lambda_{\text{a}}%
(y_{0})\Theta(y_{0})\right)  \right\}  $

$-\frac{1}{2}\underset{i=1}{\overset{n}{\sum}}%
\underset{k=q+1}{\overset{n}{\sum}}\left\{  \Gamma_{ii}^{k}(y_{0})\left(
\frac{\partial\Theta}{\partial x_{k}}(y_{0})\frac{1-r_{1}}{1-r_{2}}%
+\Lambda_{k}(y_{0})\Theta(y_{0})\right)  \right\}  +\frac{1}{2}$%
W$(y_{0})\Theta(y_{0})$

We note that:

1. $\frac{\partial\Lambda_{i}}{\partial x_{i}}(y_{0})=0=\Omega_{ii}(y_{0})$
since $\Omega_{ij}$ is skew-symmetric in the indices $\left(  i,j\right)  .$

2.$\Lambda_{i}^{2}(y_{0})=0$ for $i=q+1,...,n$

3. $\Gamma_{ij}^{\text{a}}(y_{0})=0$ for a = 1,...,q and
$i,j=1,...,q,q+1,...,n$ by (ii) of Table A$_{8}.$

4. $\Gamma_{ij}^{k}(y_{0})=0$ for $i,j,k=q+1,...,n$ by (i) of Table of
A$_{8}.$

5. $\Gamma_{\text{ab}}^{k}(y_{0})=T_{\text{ab}}^{k}(y_{0})$ for a,b = 1,...,q
and $k=q+1,...,n$ by (i) of Table A$_{7}.$

Consequently the expression in $\left(  11.33\right)  $ simplifies to:

$\left(  11.34\right)  \qquad\frac{1}{2}\Delta\lbrack\Theta(z_{1}%
)](y_{0})=\frac{1}{2}\underset{\text{a=1}}{\overset{\text{q}}{\sum}}%
\frac{\partial^{2}\Theta}{\partial x_{\text{a}}^{2}}(y_{0})+\frac{1}%
{2}\underset{i=q+1}{\overset{n}{\sum}}\frac{\partial^{2}\Theta}{\partial
x_{i}^{2}}(y_{0})(\frac{1-r_{1}}{1-r_{2}})^{2}$

$+\underset{\text{a=1}}{\overset{\text{q}}{\sum}}\Lambda_{\text{a}}%
(y_{0})\frac{\partial\Theta}{\partial x_{\text{a}}}(y_{0})+$
$\underset{i=q+1}{\overset{n}{\sum}}\Lambda_{i}(y_{0})\frac{\partial\Theta
}{\partial x_{i}}(y_{0})\frac{1-r_{1}}{1-r_{2}}+$ $\frac{1}{2}%
\underset{\text{a=1}}{\overset{\text{q}}{\sum}}\Lambda_{\text{a}}^{2}%
(y_{0})\Theta(y_{0})$

$-\frac{1}{2}\underset{\text{a=1}}{\overset{\text{q}}{\sum}}%
\underset{k=q+1}{\overset{n}{\sum}}\left\{  T_{\text{aa}}^{k}(y_{0})\left(
\frac{\partial\Theta}{\partial x_{k}}(y_{0})\frac{1-r_{1}}{1-r_{2}}%
+\Lambda_{k}(y_{0})\Theta(y_{0})\right)  \right\}  +\frac{1}{2}$%
W$(y_{0})\Theta(y_{0}).$

Elementary computations give:

\begin{center}
$\qquad\int_{0}^{1}\int_{0}^{r_{1}}(\frac{1-r_{1}}{1-r_{2}})^{2}$dr$_{1}%
$dr$_{2}=\frac{1}{6};\qquad\int_{0}^{1}\int_{0}^{r_{1}}(\frac{1-r_{1}}%
{1-r_{2}})$dr$_{1}$dr$_{2}=\frac{1}{4};\qquad\int_{0}^{1}\int_{0}^{r_{1}}%
$dr$_{1}$dr$_{2}=\frac{1}{2}$
\end{center}

Therefore we have:

$\left(  11.35\right)  \qquad\frac{1}{2}\int_{0}^{1}\int_{0}^{r_{1}}%
\Delta\lbrack\Theta(z_{1})](y_{0})$dr$_{1}$dr$_{2}\ =\frac{1}{4}%
\underset{\text{a=1}}{\overset{\text{q}}{\sum}}\frac{\partial^{2}\Theta
}{\partial x_{\text{a}}^{2}}(y_{0})+\frac{1}{12}%
\underset{i=q+1}{\overset{n}{\sum}}\frac{\partial^{2}\Theta}{\partial
x_{i}^{2}}(y_{0})$

$\qquad\qquad+\frac{1}{2}\underset{\text{a=1}}{\overset{\text{q}}{\sum}%
}\Lambda_{\text{a}}(y_{0})\frac{\partial\Theta}{\partial x_{\text{a}}}%
(y_{0})+\frac{1}{4}$ $\underset{i=q+1}{\overset{n}{\sum}}\Lambda_{i}%
(y_{0})\frac{\partial\Theta}{\partial x_{i}}(y_{0})+$ $\frac{1}{4}%
\underset{\text{a=1}}{\overset{\text{q}}{\sum}}\Lambda_{\text{a}}^{2}%
(y_{0})\Theta(y_{0})$

$\qquad\qquad-\frac{1}{8}\underset{\text{a=1}}{\overset{\text{q}}{\sum}%
}\underset{k=q+1}{\overset{n}{\sum}}T_{\text{aa}}^{k}(y_{0})\frac
{\partial\Theta}{\partial x_{k}}(y_{0})-\frac{1}{4}\underset{\text{a=1}%
}{\overset{\text{q}}{\sum}}\underset{k=q+1}{\overset{n}{\sum}}T_{\text{aa}%
}^{k}(y_{0})\Lambda_{k}(y_{0})\Theta(y_{0})+\frac{1}{4}$W$(y_{0})\Theta
(y_{0})$

By (i) of Table 7,

\qquad\qquad\qquad\ $\underset{\text{a=1}}{\overset{\text{q}}{\sum}%
}T_{\text{aak}}=$ $<H,\frac{\partial}{\partial x_{k}}>$ $\circeq$ $<H,k>$

where $T$ is \textbf{the second fundamental form} operator for the submanifold
P and $H$ is the mean curvature vector field. Consequently we have:

$\left(  11.36\right)  \qquad\frac{1}{2}\int_{0}^{1}\int_{0}^{r_{1}}%
\Delta\lbrack\Theta(z_{1})](y_{0})$dr$_{1}$dr$_{2}\ =\frac{1}{4}%
\underset{\text{a=1}}{\overset{\text{q}}{\sum}}\frac{\partial^{2}\Theta
}{\partial x_{\text{a}}^{2}}(y_{0})+\frac{1}{12}%
\underset{i=q+1}{\overset{n}{\sum}}\frac{\partial^{2}\Theta}{\partial
x_{i}^{2}}(y_{0})$

$+\frac{1}{2}\underset{\text{a=1}}{\overset{\text{q}}{\sum}}\Lambda_{\text{a}%
}(y_{0})\frac{\partial\Theta}{\partial x_{\text{a}}}(y_{0})+\frac{1}{4}$
$\underset{i=q+1}{\overset{n}{\sum}}\Lambda_{i}(y_{0})\frac{\partial\Theta
}{\partial x_{i}}(y_{0})+$ $\frac{1}{4}\underset{\text{a=1}}{\overset{\text{q}%
}{\sum}}\Lambda_{\text{a}}^{2}(y_{0})\Theta(y_{0})$

$-\frac{1}{8}\underset{k=q+1}{\overset{n}{\sum}}<H,k>(y_{0})\frac
{\partial\Theta}{\partial x_{k}}(y_{0})-\frac{1}{4}%
\underset{k=q+1}{\overset{n}{\sum}}<H,k>(y_{0})\Lambda_{k}(y_{0})+\frac{1}{4}%
$W$(y_{0})\Theta(y_{0})$

We next compute $<X,\nabla\lbrack\Theta(z_{1})]>(y_{0}):$

By (iii) of \textbf{Theorem 1},

$<$X$,\nabla\lbrack\Theta(z_{1})]>$ $=$ X$_{j}\left(  \frac{\partial}{\partial
x_{j}}+\Lambda_{j}\right)  \Theta(z_{1})=$ X$_{j}\frac{\partial}{\partial
x_{j}}\Theta(z_{1})+$ X$^{j}\Lambda_{j}\Theta(z_{1})$

From $\left(  11.28\right)  ,$ we have:

$\left(  11.37\right)  $\qquad$\qquad<$X$,\nabla\lbrack\Theta(z_{1})]>$ $=$
$\underset{\text{a=1}}{\overset{\text{q}}{\sum}}$X$_{\text{a}}\frac
{\partial\Theta}{\partial x_{\text{a}}}(z_{1})+$
$\underset{i=q+1}{\overset{n}{\sum}}$X$_{j}\frac{\partial\Theta}{\partial
x_{j}}(z_{1})\frac{1-r_{1}}{1-r_{2}}+$ $\underset{i=1}{\overset{n}{\sum}}%
$X$_{j}\Lambda_{j}\Theta(z_{1})$

We re-write the last equation above at the point $y_{0}:$

$<$X$,\nabla\lbrack$F(1-r$_{2}$,1-r$_{1}\Theta>(y_{0})$

$=$ $\underset{\text{a=1}}{\overset{\text{q}}{\sum}}$X$_{\text{a}}(y_{0}%
)$F(1-r$_{2}$,1-r$_{1}\frac{\partial\Theta}{\partial x_{\text{a}}}(y_{0})+$
$\underset{i=q+1}{\overset{n}{\sum}}$X$_{j}(y_{0})$F(1-r$_{2}$,1-r$_{1}%
)\frac{\partial\Theta}{\partial x_{j}}(y_{0})\frac{1-r_{1}}{1-r_{2}}+$
$\underset{i=1}{\overset{n}{\sum}}$X$^{j}(y_{0})\Lambda_{j}(y_{0})$F(1-r$_{2}%
$,1-r$_{1})\Theta(y_{0})$

Since by definition we we now have $\gamma_{1,2}(r_{1}-r_{2})=y_{0},$

$<$X$,\nabla\lbrack$F(1-r$_{2}$,1-r$_{1}\Theta>(y_{0})$

$=$ $\underset{\text{a=1}}{\overset{\text{q}}{\sum}}$X$_{\text{a}}(y_{0}%
)\frac{\partial\Theta}{\partial x_{\text{a}}}(\gamma_{1,2}(r_{1}-r_{2}))+$
$\underset{i=q+1}{\overset{n}{\sum}}$X$_{j}(y_{0})\frac{\partial\Theta
}{\partial x_{j}}(\gamma_{1,2}(r_{1}-r_{2}))\frac{1-r_{1}}{1-r_{2}}+$
$\underset{i=1}{\overset{n}{\sum}}$X$^{j}(y_{0})\Lambda_{j}(y_{0}%
)\Theta(\gamma_{1,2}(r_{1}-r_{2}))$

$<X,\nabla\lbrack$F(1-r$_{2}$,1-r$_{1}\Theta>(y_{0})=$ $\underset{\text{a=1}%
}{\overset{\text{q}}{\sum}}$X$_{\text{a}}(y_{0})\frac{\partial\Theta}{\partial
x_{\text{a}}}(y_{0})+$ $\underset{j=q+1}{\overset{n}{\sum}}$X$_{j}(y_{0}%
)\frac{\partial\Theta}{\partial x_{j}}(y_{0})\frac{1-r_{1}}{1-r_{2}}+$
$\underset{j=1}{\overset{n}{\sum}}$X$^{j}(y_{0})\Lambda_{j}(y_{0})\Theta
(y_{0})$

Since $\int_{0}^{1}\int_{0}^{r_{1}}(\frac{1-r_{1}}{1-r_{2}})$dr$_{1}$%
dr$_{2}=\frac{1}{4}$ and $\int_{0}^{1}\int_{0}^{r_{1}}$dr$_{1}$dr$_{2}%
=\frac{1}{2},$ we have:

$\left(  11.38\right)  \qquad\int_{0}^{1}\int_{0}^{r_{1}}<X,\nabla
\lbrack\Theta(z_{1})]>(y_{0})$dr$_{1}$dr$_{2}$

$=\frac{1}{2}$ $\underset{\text{a=1}}{\overset{\text{q}}{\sum}}$X$_{\text{a}%
}(y_{0})\frac{\partial\Theta}{\partial x_{\text{a}}}(y_{0})+\frac{1}{4}$
$\underset{j=q+1}{\overset{n}{\sum}}$X$_{j}(y_{0})\frac{\partial\Theta
}{\partial x_{j}}(y_{0})+\frac{1}{2}$ $\underset{j=1}{\overset{n}{\sum}}%
$X$_{j}(y_{0})\Lambda_{j}(y_{0})\Theta(y_{0})$

$=\frac{1}{2}$ $\underset{\text{a=1}}{\overset{\text{q}}{\sum}}$X$_{\text{a}%
}(y_{0})\frac{\partial\Theta}{\partial x_{\text{a}}}(y_{0})+\frac{1}{4}$
$\underset{j=q+1}{\overset{n}{\sum}}$X$_{j}(y_{0})\frac{\partial\Theta
}{\partial x_{j}}(y_{0})+\frac{1}{2}$ $\underset{\text{a=1}}{\overset{\text{q}%
}{\sum}}$X$_{\text{a}}(y_{0})\Lambda_{\text{a}}(y_{0})\Theta(y_{0})$

$+\frac{1}{2}$ $\underset{j=q+1}{\overset{n}{\sum}}$X$_{j}(y_{0})\Lambda
_{j}(y_{0})\Theta(y_{0})$

The computation of the last integral is simple:

$\left(  11.39\right)  \qquad\int_{0}^{1}\int_{0}^{r_{1}}[$V$\Theta
(z_{1})](y_{0})$dr$_{1}$dr$_{2}=\int_{0}^{1}\int_{0}^{r_{1}}$V$(y_{0}%
)\Theta(y_{0})$dr$_{1}$dr$_{2}=\frac{1}{2}$V$(y_{0})\Theta(y_{0})$

We now gather all terms of I$_{3}$ and have by $\left(  11.36\right)  ,\left(
11.38\right)  $ and $\left(  11.39\right)  :$

$\left(  11.40\right)  \qquad$I$_{3}=$ $\int_{0}^{1}\int_{0}^{r_{1}}$%
L$[$F(1-r$_{2}$,1-r$_{1}$)L$_{\Psi}\phi\circ\pi_{P}](y_{0})$dr$_{1}$%
dr$_{2}=\int_{0}^{1}\int_{0}^{r_{1}}$L$[\Theta(z_{1})](y_{0})$dr$_{1}$dr$_{2}$

$=\frac{1}{4}\underset{\text{a=1}}{\overset{\text{q}}{\sum}}\frac{\partial
^{2}\Theta}{\partial x_{\text{a}}^{2}}(y_{0})+\frac{1}{12}%
\underset{i=q+1}{\overset{n}{\sum}}\frac{\partial^{2}\Theta}{\partial
x_{i}^{2}}(y_{0})$

$+\frac{1}{2}\underset{\text{a=1}}{\overset{\text{q}}{\sum}}\Lambda_{\text{a}%
}(y_{0})\frac{\partial\Theta}{\partial x_{\text{a}}}(y_{0})+\frac{1}{4}$
$\underset{i=q+1}{\overset{n}{\sum}}\Lambda_{i}(y_{0})\frac{\partial\Theta
}{\partial x_{i}}(y_{0})+$ $\frac{1}{4}\underset{\text{a=1}}{\overset{\text{q}%
}{\sum}}\Lambda_{\text{a}}^{2}(y_{0})\Theta(y_{0})$

$-\frac{1}{8}\underset{k=q+1}{\overset{n}{\sum}}<H,k>(y_{0})\frac
{\partial\Theta}{\partial x_{k}}(y_{0})-\frac{1}{4}%
\underset{j=q+1}{\overset{n}{\sum}}<H,j>(y_{0})\Lambda_{j}(y_{0})\Theta
(y_{0})+\frac{1}{4}$W$(y_{0})\Theta(y_{0})$

$+\frac{1}{2}$ $\underset{\text{a=1}}{\overset{\text{q}}{\sum}}$X$_{\text{a}%
}(y_{0})\frac{\partial\Theta}{\partial x_{\text{a}}}(y_{0})+\frac{1}{4}$
$\underset{j=q+1}{\overset{n}{\sum}}$X$_{j}(y_{0})\frac{\partial\Theta
}{\partial x_{j}}(y_{0})+\frac{1}{2}$ $\underset{\text{a=1}}{\overset{\text{q}%
}{\sum}}$X$_{\text{a}}(y_{0})\Lambda_{\text{a}}(y_{0})\Theta(y_{0})+\frac
{1}{2}$ $\underset{j=q+1}{\overset{n}{\sum}}$X$_{j}(y_{0})\Lambda_{j}%
(y_{0})\Theta(y_{0})$

$+\frac{1}{2}$V$(y_{0})\Theta(y_{0})$

Next we have:

\qquad\qquad\qquad\qquad\qquad\qquad\qquad\qquad\qquad\qquad\qquad\qquad
\qquad$\qquad$

$\ $I$_{4}=-\int_{0}^{1}\int_{0}^{r_{1}}\ $J$_{4}$(y$_{0}$)dr$_{1}$dr$_{2}=-$
$\int_{0}^{1}\int_{0}^{r_{1}}$V(y$_{0}$)[F(1-r$_{2}$,1-r$_{1}$)L$_{\Psi}%
\phi\circ\pi_{\text{P}}$](y$_{0}$)dr$_{1}$dr$_{2}$

$\ \ \ \ \ =-\int_{0}^{1}\int_{0}^{r_{1}}$V(y$_{0}$)L$_{\Psi}\phi\circ
\pi_{\text{P}}](y_{0})$dr$_{1}$dr$_{2}=-\int_{0}^{1}\int_{0}^{r_{1}}$V(y$_{0}%
$)$\Theta(y_{0})$dr$_{1}$dr$_{2}$

$\ \ \ \ =-\frac{1}{2}$V$(y_{0})\Theta(y_{0})$

$\left(  11.41\right)  \qquad$I$_{4}=-\frac{1}{2}$V$(y_{0})\Theta(y_{0})$

\qquad\qquad\qquad\qquad\qquad\qquad\qquad\qquad\qquad\qquad\qquad\qquad
\qquad\qquad$\blacksquare$\qquad\qquad\qquad\qquad\qquad\qquad\qquad
\qquad\qquad\qquad\qquad\qquad\qquad\qquad\qquad

We gather all terms of b$_{2}($y$_{0}$,P$,\phi)$ in $\left(  11.6\right)  ,$
$\left(  11.21\right)  ,$ $\left(  11.40\right)  $ and $\left(  11.41\right)
$ respectively and have the formula:

$\left(  11.42\right)  \qquad$b$_{2}($y$_{0}$,P$,\phi)=$ I$_{1}+$ I$_{2}+$
I$_{3}+$ I$_{4}$

$=\frac{1}{2}\frac{\text{L}\Psi}{\Psi}(y_{0})\Theta(y_{0})\qquad$I$_{1}$

$+$ $\frac{1}{8}\underset{j=q+1}{\overset{n}{\sum}}<H,j>(y_{0})\frac
{\partial\Theta}{\partial x_{j}}(y_{0})-\frac{1}{4}%
\underset{j=q+1}{\overset{n}{\sum}}X_{j}(y_{0})\frac{\partial\Theta}{\partial
x_{j}}(y_{0})\qquad$I$_{2}$

$+\frac{1}{4}\underset{\text{a=1}}{\overset{\text{q}}{\sum}}\frac{\partial
^{2}\Theta}{\partial x_{\text{a}}^{2}}(y_{0})+\frac{1}{12}%
\underset{i=q+1}{\overset{n}{\sum}}\frac{\partial^{2}\Theta}{\partial
x_{i}^{2}}(y_{0})\qquad$I$_{3}$ starts

$+\frac{1}{2}\underset{\text{a=1}}{\overset{\text{q}}{\sum}}\Lambda_{\text{a}%
}(y_{0})\frac{\partial\Theta}{\partial x_{\text{a}}}(y_{0})+\frac{1}{4}$
$\underset{i=q+1}{\overset{n}{\sum}}\Lambda_{i}(y_{0})\frac{\partial\Theta
}{\partial x_{i}}(y_{0})+$ $\frac{1}{4}\underset{\text{a=1}}{\overset{\text{q}%
}{\sum}}\Lambda_{\text{a}}^{2}(y_{0})\Theta(y_{0})$

$-\frac{1}{8}\underset{j=q+1}{\overset{n}{\sum}}<H,j>(y_{0})\frac
{\partial\Theta}{\partial x_{j}}(y_{0})-\frac{1}{4}%
\underset{j=q+1}{\overset{n}{\sum}}<H,j>(y_{0})\Lambda_{j}(y_{0})\Theta
(y_{0})+\frac{1}{4}$W$(y_{0})\Theta(y_{0})$

$+\frac{1}{2}$ $\underset{\text{a=1}}{\overset{\text{q}}{\sum}}X_{\text{a}%
}(y_{0})\frac{\partial\Theta}{\partial x_{\text{a}}}(y_{0})+\frac{1}{4}$
$\underset{j=q+1}{\overset{n}{\sum}}X_{j}(y_{0})\frac{\partial\Theta}{\partial
x_{j}}(y_{0})+\frac{1}{2}$ $\underset{\text{a=1}}{\overset{\text{q}}{\sum}%
}X_{\text{a}}(y_{0})\Lambda_{\text{a}}(y_{0})\Theta(y_{0})$

$+\frac{1}{2}$ $\underset{j=q+1}{\overset{n}{\sum}}X_{j}(y_{0})\Lambda
_{j}(y_{0})\Theta(y_{0})+\frac{1}{2}$V$(y_{0})\Theta(y_{0})\qquad$I$_{3}$ ends

$-\frac{1}{2}$V$(y_{0})\Theta(y_{0})\qquad$I$_{4}$

There are obvious cancellations above: We notice that I$_{2}$ and I$_{4}$ have
been wiped off by some parts of I$_{3}$ in the above expression. In
particular, we see that:

$+$ $\frac{1}{8}\underset{j=q+1}{\overset{n}{\sum}}<H,j>(y_{0})\frac
{\partial\Theta}{\partial x_{j}}(y_{0})-\frac{1}{4}%
\underset{j=q+1}{\overset{n}{\sum}}X_{j}(y_{0})\frac{\partial\Theta}{\partial
x_{j}}(y_{0})\qquad$I$_{2}$

$-\frac{1}{8}\underset{j=q+1}{\overset{n}{\sum}}<H,j>(y_{0})\frac
{\partial\Theta}{\partial x_{j}}(y_{0})+\frac{1}{4}$
$\underset{j=q+1}{\overset{n}{\sum}}$X$_{j}(y_{0})\frac{\partial\Theta
}{\partial x_{j}}(y_{0})=0$

\qquad\qquad\qquad\qquad\qquad\qquad\qquad\qquad\qquad\qquad\qquad\qquad
\qquad\qquad\qquad\qquad$\blacksquare$

We set:

\qquad I$_{1}=\frac{1}{2}\frac{\text{L}\Psi}{\Psi}(y_{0})\Theta(y_{0})$

\qquad I$_{31}=\frac{1}{12}\underset{\text{a=1}}{\overset{\text{q}}{\sum}%
}\frac{\partial^{2}\Theta}{\partial x_{\text{a}}^{2}}(y_{0});$ \qquad
I$_{32}=\frac{1}{12}\underset{i=q+1}{\overset{n}{\sum}}\frac{\partial
^{2}\Theta}{\partial x_{i}^{2}}(y_{0})$

\qquad I$_{33}=\frac{1}{2}\underset{\text{a=1}}{\overset{\text{q}}{\sum}%
}\Lambda_{\text{a}}(y_{0})\frac{\partial\Theta}{\partial x_{\text{a}}}%
(y_{0});\qquad$I$_{34}=\frac{1}{4}\underset{i=q+1}{\overset{n}{\sum}}%
\Lambda_{i}(y_{0})\frac{\partial\Theta}{\partial x_{i}}(y_{0});\qquad$

\qquad I$_{35}=\frac{1}{4}\underset{\text{a=1}}{\overset{\text{q}}{\sum}%
}\Lambda_{\text{a}}^{2}(y_{0})\Theta(y_{0});\qquad$

\qquad I$_{36}=-\frac{1}{4}\underset{j=q+1}{\overset{n}{\sum}}<H,j>(y_{0}%
)\Lambda_{j}(y_{0})\Theta(y_{0});\qquad$I$_{37}=$ $\frac{1}{4}$W$(y_{0}%
)\Theta(y_{0});\qquad$

\qquad I$_{38}=$ $\frac{1}{2}$ $\underset{\text{a=1}}{\overset{\text{q}}{\sum
}}X_{\text{a}}(y_{0})\frac{\partial\Theta}{\partial x_{\text{a}}}%
(y_{0});\qquad$I$_{39}=\frac{1}{2}$ $\underset{j=1}{\overset{n}{\sum}}$%
X$_{j}(y_{0})\Lambda_{j}(y_{0})\Theta(y_{0})$Then,

\qquad\qquad\qquad\qquad I$_{3}=$ I$_{31}+$ I$_{32}+$ I$_{33}+$I$_{34}+$
I$_{35}+$ I$_{36}+$ I$_{37}+$ I$_{38}+$ I$_{39}$

\qquad\qquad\qquad\qquad\qquad\qquad\qquad\qquad\qquad\qquad\qquad\qquad
\qquad\qquad\qquad\qquad\qquad$\blacksquare$

Consequently, we have the\textbf{ "Raw" Expression} for the third term as a:

\begin{lemma}
b$_{2}($y$_{0}$,P$,\phi)=$ I$_{1}+$ I$_{3}=$ I$_{1}+$ I$_{31}+$ I$_{32}+$
I$_{33}+$I$_{34}+$ I$_{35}+$ I$_{36}+$ I$_{37}+$ I$_{38}+$ I$_{39}$
\end{lemma}

$\qquad=\frac{1}{2}\frac{\text{L}\Psi}{\Psi}(y_{0})\Theta(y_{0})\qquad$I$_{1}$

$+\frac{1}{4}\underset{\text{a=1}}{\overset{\text{q}}{\sum}}\frac{\partial
^{2}\Theta}{\partial x_{\text{a}}^{2}}(y_{0})+\frac{1}{12}%
\underset{i=q+1}{\overset{n}{\sum}}\frac{\partial^{2}\Theta}{\partial
x_{i}^{2}}(y_{0})+\frac{1}{2}\underset{\text{a=1}}{\overset{\text{q}}{\sum}%
}\Lambda_{\text{a}}(y_{0})\frac{\partial\Theta}{\partial x_{\text{a}}}%
(y_{0})+\frac{1}{4}\underset{i=q+1}{\overset{n}{\sum}}\Lambda_{i}(y_{0}%
)\frac{\partial\Theta}{\partial x_{i}}(y_{0})\qquad$I$_{3}$ starts

$+$ $\frac{1}{4}\underset{\text{a=1}}{\overset{\text{q}}{\sum}}\Lambda
_{\text{a}}^{2}(y_{0})\Theta(y_{0})-\frac{1}{4}%
\underset{j=q+1}{\overset{n}{\sum}}<H,j>(y_{0})\Lambda_{j}(y_{0})\Theta
(y_{0})+\frac{1}{4}\Theta(y_{0})$W$(y_{0})$

$+\frac{1}{2}$ $\underset{\text{a=1}}{\overset{\text{q}}{\sum}}$X$_{\text{a}%
}(y_{0})\frac{\partial\Theta}{\partial x_{\text{a}}}(y_{0})+\frac{1}{2}$
$\underset{\text{a=1}}{\overset{\text{q}}{\sum}}$X$_{\text{a}}(y_{0}%
)\Lambda_{\text{a}}(y_{0})\Theta(y_{0})+\frac{1}{2}$
$\underset{j=q+1}{\overset{n}{\sum}}$X$_{j}(y_{0})\Lambda_{j}(y_{0}%
)\Theta(y_{0})\qquad\qquad$I$_{3}$ ends

\qquad\qquad\qquad\qquad\qquad\qquad\qquad\qquad\qquad\qquad\qquad\qquad
\qquad\qquad\qquad$\qquad\qquad\qquad\blacksquare$

The computations that express b$_{2}($y$_{0}$,P$,\phi)$ in terms of
\textbf{geometric invariants} have been done in \textbf{Appendix D:}

\qquad\qquad\qquad\qquad\qquad\qquad\qquad\qquad\qquad\qquad\qquad\qquad
\qquad\qquad\qquad\qquad\qquad\qquad$\blacksquare$

We now come to one of the most important theorems of this work. After lengthy
computations, the expression for the third coefficient given in
\textbf{geometric invariants} of the Riemannian manifold M, of the submanifold
P and of the vector bundle E, is given in the theorem below.

Our work here can be regarded as the \textbf{ultimate generalization} of
\textbf{heat kernel expansions} in the following sense:

Firstly we are working in the more general context of a \textbf{vector bundle}
E over a Riemannian manifold M.

Secondly we are working with \textbf{Fermi coordinates} which generalize
normal coordinates (equivalently the center of normal coordinates y$_{0}$ is
generalized to a submanifold P).

\qquad\qquad\qquad\qquad\qquad\qquad\qquad\qquad\qquad\qquad\qquad\qquad
\qquad\qquad\qquad\qquad$\blacksquare$

\begin{theorem}
b$_{2}($y$_{0}$,P$,\phi)=$ I$_{1}+$ I$_{31}+$ I$_{32}+$ I$_{33}+$ I$_{34}+$
I$_{35}+$ I$_{36}+$ I$_{37}$
\end{theorem}

\qquad$=\frac{1}{2}[\frac{1}{24}(\underset{i=q+1}{\overset{n}{\sum}}%
3<H,i>^{2}+2(\tau^{M}-3\tau^{P}\ +\overset{q}{\underset{\text{a=1}}{\sum}%
}\varrho_{\text{aa}}^{M}+\overset{q}{\underset{\text{a,b}=1}{\sum}%
}R_{\text{abab}}^{M}))\qquad$I$_{1}$

\qquad$-\frac{1}{2}($ $\left\Vert \text{X}\right\Vert _{M}^{2}+$
$\operatorname{div}X_{M}-\left\Vert \text{X}\right\Vert _{P}^{2}$ $-$
$\operatorname{div}X_{P})+$ $V](y_{0})$

$\times\lbrack\frac{1}{24}(\underset{i=q+1}{\overset{n}{\sum}}3<H,i>^{2}%
+2(\tau^{M}-3\tau^{P}\ +\overset{q}{\underset{\text{a=1}}{\sum}}%
\varrho_{\text{aa}}^{M}+\overset{q}{\underset{\text{a,b}=1}{\sum}%
}R_{\text{abab}}^{M})\phi\left(  y_{0}\right)  $

$-\frac{1}{2}($ $\left\Vert \text{X}\right\Vert _{M}^{2}+$ $\operatorname{div}%
X_{M}-\left\Vert \text{X}\right\Vert _{P}^{2}$ $-\operatorname{div}X_{P})+$
$V+\frac{1}{2}W)\phi\left(  y_{0}\right)  $

$+$ $(\frac{1}{2}\underset{\text{a=1}}{\overset{\text{q}}{\sum}}\frac
{\partial^{2}\phi}{\partial\text{x}_{\text{a}}^{2}}$ $+$ $\underset{\text{a=1}%
}{\overset{\text{q}}{\sum}}\Lambda_{\text{a}}(y_{0})\frac{\partial\phi
}{\partial x_{\text{a}}}\ +\frac{1}{2}$ $\underset{\text{a=1}%
}{\overset{\text{q}}{\sum}}\Lambda_{\text{a}}\Lambda_{\text{a}})+$
$\underset{\text{a=1}}{\overset{\text{q}}{\sum}}X_{\text{a}}\frac{\partial
\phi}{\partial\text{x}_{\text{a}}}+$ $\underset{\text{a=1}}{\overset{\text{q}%
}{\sum}}$ $X_{\text{a}}\Lambda_{\text{a}})\phi]\left(  y_{0}\right)  $

$+\frac{1}{96}[\underset{i=q+1}{\overset{n}{\sum}}3<H,i>^{2}+2(\tau^{M}%
-3\tau^{P}\ +\overset{q}{\underset{\text{a=1}}{\sum}}\varrho_{\text{aa}}%
^{M}+\overset{q}{\underset{\text{a,b}=1}{\sum}}R_{\text{abab}}^{M}%
)](y_{0})\frac{\partial^{2}\phi}{\partial x_{\text{c}}^{2}}(y_{0})\qquad
$I$_{31}\qquad$I$_{311}$

$-\frac{1}{4}[\left\Vert \text{X}(y_{0})\right\Vert ^{2}+$ divX$(y_{0})-$
$\underset{\text{a}=1}{\overset{q}{\sum}}($X$_{\text{a}})^{2}(y_{0})$ $-$
$\underset{\text{a}=1}{\overset{q}{\sum}}\frac{\partial X_{\text{a}}}{\partial
x_{\text{a}}}(y_{0})]\frac{\partial^{2}\phi}{\partial x_{\text{c}}^{2}}%
(y_{0})+$ $\frac{1}{4}$V(y$_{0}$)$\frac{\partial^{2}\phi}{\partial
x_{\text{c}}^{2}}(y_{0})$

$-\frac{1}{2}$ $[X_{j}\frac{\partial X_{j}}{\partial x_{\text{c}}}+\frac{1}%
{2}\frac{\partial^{2}X_{j}}{\partial x_{\text{c}}\partial x_{j}}](y_{0}%
).\frac{\partial\phi}{\partial x_{\text{c}}}(y_{0})+\frac{1}{4}[<H,j>\frac
{\partial X_{j}}{\partial x_{\text{c}}}](y_{0}).\frac{\partial\phi}{\partial
x_{\text{c}}}(y_{0})+$ $\frac{1}{2}\frac{\partial\text{V}}{\partial
x_{\text{c}}}(y_{0}).\frac{\partial\phi}{\partial x_{\text{c}}}(y_{0})$

$+\frac{1}{4}[(\frac{\partial X_{j}}{\partial x_{\text{c}}})^{2}-X_{j}%
\frac{\partial^{2}X_{j}}{\partial x_{\text{c}}^{2}}](y_{0})\phi(y_{0})$
$-\frac{1}{8}\frac{\partial^{3}X_{j}}{\partial x_{\text{c}}^{2}\partial x_{j}%
}(y_{0})\phi(y_{0})-\frac{1}{2}\frac{\partial X_{i}}{\partial x_{\text{c}}%
}(y_{0}))\frac{\partial X_{i}}{\partial x_{\text{c}}}(y_{0})\phi(y_{0}%
)+\frac{1}{4}\frac{\partial^{2}\text{V}}{\partial x_{\text{c}}^{2}}(y_{0}%
)\phi(y_{0})$

$+\frac{1}{8}\underset{\text{a=1}}{\overset{\text{q}}{\sum}}\frac{\partial
^{4}\phi}{\partial x_{\text{a}}^{2}\partial x_{\text{c}}^{2}}(y_{0}%
)\qquad\qquad\qquad$I$_{312}$

$+\frac{1}{4}\underset{\text{a=1}}{\overset{\text{q}}{\sum}}[\Lambda
_{\text{a}}\frac{\partial^{3}\phi}{\partial\text{x}_{\text{a}}\partial
x_{\text{c}}^{2}}](y_{0})\qquad\qquad$I$_{314}$

$+$ $\frac{1}{8}\underset{\text{a=1}}{\overset{\text{q}}{\sum}}[\Lambda
_{\text{a}}^{2}\frac{\partial^{2}\phi}{\partial x_{\text{c}}^{2}}%
](y_{0})\qquad\ \ \ \ $I$_{315}\qquad$

$+\frac{1}{8}\frac{\partial^{2}\text{W}}{\partial x_{\text{a}}^{2}}(y_{0}%
)\phi(y_{0})+\frac{1}{4}\frac{\partial\text{W}}{\partial x_{\text{a}}}%
(y_{0})\frac{\partial\phi}{\partial x_{\text{a}}}(y_{0})+\frac{1}{8}$%
W$(y_{0})\frac{\partial^{2}\phi}{\partial x_{\text{a}}^{2}}(y_{0})\qquad
$I$_{319}$

$+\frac{1}{4}$ $\underset{\text{a=1}}{\overset{\text{q}}{\sum}}\frac
{\partial^{2}\text{X}_{\text{a}}}{\partial x_{\text{a}}^{2}}(y_{0}%
)\frac{\partial\phi}{\partial\text{x}_{\text{a}}}(y_{0})+\frac{1}{4}$
$\underset{\text{a=1}}{\overset{\text{q}}{\sum}}$X$_{\text{a}}(y_{0}%
)\frac{\partial^{3}\phi}{\partial x_{\text{a}}^{3}}(y_{0})+\frac{1}{2}$
$\underset{\text{a=1}}{\overset{\text{q}}{\sum}}\frac{\partial\text{X}%
_{\text{a}}}{\partial x_{\text{a}}}(y_{0})\frac{\partial^{2}\phi}%
{\partial\text{x}_{\text{a}}^{2}}(y_{0})\qquad$L$_{1}$

$+\frac{1}{4}[$ $\underset{\text{b=1}}{\overset{\text{q}}{\sum}}$
$\frac{\partial^{2}X_{\text{b}}}{\partial x_{\text{a}}^{2}}\Lambda_{\text{b}%
}(y_{0})\phi(y_{0})+\frac{1}{4}[$ $\underset{\text{b=1}}{\overset{\text{q}%
}{\sum}}$ X$_{\text{b}}(y_{0})\Lambda_{\text{b}}(y_{0})\frac{\partial^{2}\phi
}{\partial x_{\text{a}}^{2}}(y_{0})+\frac{1}{2}[$ $\underset{\text{b=1}%
}{\overset{\text{q}}{\sum}}$ $\frac{\partial X_{\text{b}}}{\partial
x_{\text{a}}}(y_{0})\Lambda_{\text{b}}(y_{0})\frac{\partial\phi}{\partial
x_{\text{a}}}(y_{0})\qquad$L$_{2}$

$-\frac{1}{3456}[3<H,i>^{2}\ +2(\tau^{M}-3\tau^{P}%
\ +\overset{q}{\underset{\text{a=1}}{\sum}}\varrho_{\text{aa}}^{M}%
+\overset{q}{\underset{\text{a,b}=1}{\sum}}R_{\text{abab}}^{M})]^{2}%
(y_{0})\phi(y_{0})\phi(y_{0})$\ \qquad\textbf{I}$_{32}\qquad$I$_{321}$

$+\frac{1}{24}[2<H,i>^{2}(y_{0})+\frac{1}{3}(\tau^{M}-3\tau^{P}%
+\overset{q}{\underset{\text{a}=1}{\sum}}\varrho_{\text{aa}}%
+\overset{q}{\underset{\text{a,b}=1}{\sum}}R_{\text{abab}})](y_{0})\qquad
$I$_{3212}=\frac{1}{24}(L_{1}+L_{2}+L_{3})$

$\times\lbrack\frac{1}{4}<H,j>^{2}(y_{0})+\frac{1}{6}(\tau^{M}-3\tau
^{P}+\overset{q}{\underset{\text{a}=1}{\sum}}\varrho_{\text{aa}}%
^{M}+\overset{q}{\underset{\text{a,b}=1}{\sum}}R_{\text{abab}}^{M}%
)](y_{0})\phi(y_{0})$

$-\frac{1}{96}[<H,i><H,j>](y_{0})$

$\times\lbrack2\varrho_{ij}+$ $\overset{q}{\underset{\text{a}=1}{4\sum}%
}R_{i\text{a}j\text{a}}-3\overset{q}{\underset{\text{a,b=1}}{\sum}%
}(T_{\text{aa}i}T_{\text{bb}j}-T_{\text{ab}i}T_{\text{ab}j}%
)-3\overset{q}{\underset{\text{a,b=1}}{\sum}}(T_{\text{aa}j}T_{\text{bb}%
i}-T_{\text{ab}j}T_{\text{ab}i}](y_{0})\phi(y_{0})\qquad L_{2}\qquad
L_{21}\qquad\ \ \ \ \ $

$-\frac{1}{864}[2\varrho_{ij}+$ $\overset{q}{\underset{\text{a}=1}{4\sum}%
}R_{i\text{a}j\text{a}}-3\overset{q}{\underset{\text{a,b=1}}{\sum}%
}(T_{\text{aa}i}T_{\text{bb}j}-T_{\text{ab}i}T_{\text{ab}j}%
)-3\overset{q}{\underset{\text{a,b=1}}{\sum}}(T_{\text{aa}j}T_{\text{bb}%
i}-T_{\text{ab}j}T_{\text{ab}i}]^{2}(y_{0})\phi(y_{0})$

$-\frac{1}{288}[<H,j>](y_{0})\times\lbrack\{\nabla_{i}\varrho_{ij}%
-2\varrho_{ij}<H,i>+\overset{q}{\underset{\text{a}=1}{\sum}}(\nabla
_{i}R_{\text{a}i\text{a}j}-4R_{i\text{a}j\text{a}}<H,i>)\qquad L_{212}$

$+4\overset{q}{\underset{\text{a,b=1}}{\sum}}R_{i\text{a}j\text{b}%
}T_{\text{ab}i}+2\overset{q}{\underset{\text{a,b,c=1}}{\sum}}(T_{\text{aa}%
i}T_{\text{bb}j}T_{\text{cc}i}-3T_{\text{aa}i}T_{\text{bc}j}T_{\text{bc}%
i}+2T_{\text{ab}i}T_{\text{bc}j}T_{\text{ac}i})](y_{0})\phi(y_{0})$%
\qquad\qquad\qquad\qquad\qquad\ \ 

$-\frac{1}{288}[<H,j>](y_{0})\times\lbrack\nabla_{j}\varrho_{ii}-2\varrho
_{ij}<H,i>+\overset{q}{\underset{\text{a}=1}{\sum}}(\nabla_{j}R_{\text{a}%
i\text{a}i}-4R_{i\text{a}j\text{a}}<H,i>)$

$+4\overset{q}{\underset{\text{a,b=1}}{\sum}}R_{j\text{a}i\text{b}%
}T_{\text{ab}i}+2\overset{q}{\underset{\text{a,b,c=1}}{\sum}}(T_{\text{aa}%
j}T_{\text{bb}i}T_{\text{cc}i}-3T_{\text{aa}j}T_{\text{bc}i}T_{\text{bc}%
i}+2T_{\text{ab}j}T_{\text{bc}i}T_{\text{ac}i})](y_{0})\phi(y_{0})$

$-\frac{1}{288}[<H,j>](y_{0})\times\lbrack\nabla_{i}\varrho_{ij}-2\varrho
_{ii}<H,j>+\overset{q}{\underset{\text{a}=1}{\sum}}(\nabla_{i}R_{\text{a}%
i\text{a}j}-4R_{i\text{a}i\text{a}}<H,j>)$

$+4\overset{q}{\underset{\text{a,b=1}}{\sum}}R_{i\text{a}i\text{b}%
}T_{\text{ab}j}+2\overset{q}{\underset{\text{a,b,c}=1}{\sum}}(T_{\text{aa}%
i}T_{\text{bb}i}T_{\text{cc}j}-3T_{\text{aa}i}T_{\text{bc}i}T_{\text{bc}%
j}+2T_{\text{ab}i}T_{\text{bc}i}T_{\text{ac}j})](y_{0})\phi(y_{0})$

$-\frac{1}{3}[<H,j><H,k>](y_{0})R_{ijik}(y_{0})\phi(y_{0})-$ $\frac{5}%
{64}[<H,i>^{2}<H,j>^{2}](y_{0})\phi(y_{0})\qquad L_{213}$

$-\frac{1}{96}<H,i><H,j>$

$\times\lbrack2\varrho_{ij}+\overset{q}{\underset{\text{a}=1}{4\sum}%
}R_{i\text{a}j\text{a}}-3\overset{q}{\underset{\text{a,b=1}}{\sum}%
}(T_{\text{aa}i}T_{\text{bb}j}-T_{\text{ab}i}T_{\text{ab}j}%
)-3\overset{q}{\underset{\text{a,b=1}}{\sum}}(T_{\text{aa}j}T_{\text{bb}%
i}-T_{\text{ab}j}T_{\text{ab}i}](y_{0})\phi(y_{0})$

$-\frac{1}{96}<H,j>^{2}[\tau^{M}\ -3\tau^{P}+\ \underset{\text{a}%
=1}{\overset{\text{q}}{\sum}}\varrho_{\text{aa}}^{M}+$
$\overset{q}{\underset{\text{a},\text{b}=1}{\sum}}R_{\text{abab}}^{M}$
$](y_{0})\phi(y_{0})$

$+\frac{1}{288}<H,j>[\nabla_{i}\varrho_{ij}-2\varrho_{ij}%
<H,i>+\overset{q}{\underset{\text{a}=1}{\sum}}(\nabla_{i}R_{\text{a}%
i\text{a}j}-4R_{i\text{a}j\text{a}}<H,i>)$

$+4\overset{q}{\underset{\text{a,b=1}}{\sum}}R_{i\text{a}j\text{b}%
}T_{\text{ab}i}+2\overset{q}{\underset{\text{a,b,c=1}}{\sum}}(T_{\text{aa}%
i}T_{\text{bb}j}T_{\text{cc}i}-3T_{\text{aa}i}T_{\text{bc}j}T_{\text{bc}%
i}+2T_{\text{ab}i}T_{\text{bc}j}T_{\text{ac}i})](y_{0})\phi(y_{0})$%
\qquad\qquad\qquad\qquad\qquad\ \ 

$+\frac{1}{12}<H,j>[\nabla_{j}\varrho_{ii}-2\varrho_{ij}%
<H,i>+\overset{q}{\underset{\text{a}=1}{\sum}}(\nabla_{j}R_{\text{a}%
i\text{a}i}-4R_{i\text{a}j\text{a}}<H,i>)$

$+4\overset{q}{\underset{\text{a,b=1}}{\sum}}R_{j\text{a}i\text{b}%
}T_{\text{ab}i}+2\overset{q}{\underset{\text{a,b,c=1}}{\sum}}(T_{\text{aa}%
j}T_{\text{bb}i}T_{\text{cc}i}-3T_{\text{aa}j}T_{\text{bc}i}T_{\text{bc}%
i}+2T_{\text{ab}j}T_{\text{bc}i}T_{\text{ac}i})](y_{0})\phi(y_{0})$

$+\frac{1}{288}<H,j>[\nabla_{i}\varrho_{ij}-2\varrho_{ii}%
<H,j>+\overset{q}{\underset{\text{a}=1}{\sum}}(\nabla_{i}R_{\text{a}%
i\text{a}j}-4R_{i\text{a}i\text{a}}<H,j>)+4\overset{q}{\underset{\text{a,b=1}%
}{\sum}}R_{i\text{a}i\text{b}}T_{\text{ab}j}$

$+2\overset{q}{\underset{\text{a,b,c}=1}{\sum}}(T_{\text{aa}i}T_{\text{bb}%
i}T_{\text{cc}j}-3T_{\text{aa}i}T_{\text{bc}i}T_{\text{bc}j}+2T_{\text{ab}%
i}T_{\text{bc}i}T_{\text{ac}j})](y_{0})\phi(y_{0})$

$-\frac{1}{144}R_{jijk}(y_{0})$ \ $[<H,i><H,k>](y_{0})\phi(y_{0})\qquad\qquad
L_{22}$

$-\frac{1}{432}R_{jijk}(y_{0})[2\varrho_{ik}+$ $\overset{q}{\underset{\text{a}%
=1}{4\sum}}R_{i\text{a}k\text{a}}-3\overset{q}{\underset{\text{a,b=1}}{\sum}%
}(T_{\text{aa}i}T_{\text{bb}k}-T_{\text{ab}i}T_{\text{ab}k}%
)-3\overset{q}{\underset{\text{a,b=1}}{\sum}}(T_{\text{aa}k}T_{\text{bb}%
i}-T_{\text{ab}k}T_{\text{ab}i}](y_{0})\phi(y_{0})$

$+\frac{1}{144}<H,k>(y_{0})[\nabla_{j}$R$_{ijik}(y_{0})-\nabla_{i}$%
R$_{jijk}](y_{0})\phi(y_{0})$

$-\frac{5}{32}<H,i>^{2}(y_{0})<H,j>^{2}(y_{0})\phi(y_{0})\qquad\qquad
L_{23}\qquad L_{231}$

$-\frac{1}{48}<H,i>(y_{0})<H,j>(y_{0})$

$\times\lbrack2\varrho_{ij}+$ $\overset{q}{\underset{\text{a}=1}{4\sum}%
}R_{i\text{a}j\text{a}}-3\overset{q}{\underset{\text{a,b=1}}{\sum}%
}(T_{\text{aa}i}T_{\text{bb}j}-T_{\text{ab}i}T_{\text{ab}j}%
)-3\overset{q}{\underset{\text{a,b=1}}{\sum}}(T_{\text{aa}j}T_{\text{bb}%
i}-T_{\text{ab}j}T_{\text{ab}i}](y_{0})\phi(y_{0})$

$-\frac{1}{48}<H,i>^{2}(y_{0})[\varrho_{jj}+$ $\overset{q}{\underset{\text{a}%
=1}{2\sum}}R_{j\text{a}j\text{a}}-3\overset{q}{\underset{\text{a,b=1}}{\sum}%
}(T_{\text{aa}j}T_{\text{bb}j}-T_{\text{ab}j}T_{\text{ab}j})](y_{0})\phi
(y_{0})$

$-\frac{1}{144}<H,i>(y_{0})[\nabla_{i}\varrho_{jj}-2\varrho_{ij}%
<H,j>+\overset{q}{\underset{\text{a}=1}{\sum}}(\nabla_{i}R_{\text{a}%
j\text{a}j}-4R_{i\text{a}j\text{a}}<H,j>)+4\overset{q}{\underset{\text{a,b=1}%
}{\sum}}R_{i\text{a}j\text{b}}T_{\text{ab}j}$

$+2\overset{q}{\underset{\text{a,b,c=1}}{\sum}}(T_{\text{aa}i}T_{\text{bb}%
j}T_{\text{cc}j}-3T_{\text{aa}i}T_{\text{bc}j}T_{\text{bc}j}+2T_{\text{ab}%
i}T_{\text{bc}j}T_{\text{ca}j})](y_{0})\phi(y_{0})$\qquad\qquad\qquad
\qquad\qquad\ \ 

$-\frac{1}{24}\times\frac{1}{6}<H,i>(y_{0})[\nabla_{j}\varrho_{ij}%
-2\varrho_{ij}<H,j>+\overset{q}{\underset{\text{a}=1}{\sum}}(\nabla
_{j}R_{\text{a}i\text{a}j}-4R_{j\text{a}i\text{a}}%
<H,j>)+4\overset{q}{\underset{\text{a,b=1}}{\sum}}R_{j\text{a}i\text{b}%
}T_{\text{ab}j}$

$+2\overset{q}{\underset{\text{a,b,c=1}}{\sum}}(T_{\text{aa}j}T_{\text{bb}%
i}T_{\text{cc}j}-3T_{\text{aa}j}T_{\text{bc}i}T_{\text{bc}j}+2T_{\text{ab}%
j}T_{\text{bc}i}T_{\text{ac}j})](y_{0})\phi(y_{0})$

$-\frac{1}{144}<H,i>(y_{0})[\nabla_{j}\varrho_{ij}-2\varrho_{jj}%
<H,i>+\overset{q}{\underset{\text{a}=1}{\sum}}(\nabla_{j}R_{\text{a}%
i\text{a}j}-4R_{j\text{a}j\text{a}}<H,i>)+4\overset{q}{\underset{\text{a,b=1}%
}{\sum}}R_{j\text{a}j\text{b}}T_{\text{ab}i}$

$+2\overset{q}{\underset{\text{a,b,c=1}}{\sum}}(T_{\text{aa}j}T_{\text{bb}%
j}T_{\text{cc}i}-3T_{\text{aa}j}T_{\text{bc}j}T_{\text{bc}i}+2T_{\text{ab}%
j}T_{\text{bc}j}T_{\text{ac}i})](y_{0})\phi(y_{0})$

$-\frac{1}{96}<H,j>^{2}(y_{0})[\varrho_{ii}+$ $\overset{q}{\underset{\text{a}%
=1}{2\sum}}R_{i\text{a}i\text{a}}-3\overset{q}{\underset{\text{a,b=1}}{\sum}%
}(T_{\text{aa}i}T_{\text{bb}i}-T_{\text{ab}i}T_{\text{ab}i})](y_{0})\phi
(y_{0})\qquad L_{232}$

$-\frac{1}{432}[\varrho_{ii}+$ $\overset{q}{\underset{\text{a}=1}{2\sum}%
}R_{i\text{a}i\text{a}}-3\overset{q}{\underset{\text{a,b=1}}{\sum}%
}(T_{\text{aa}i}T_{\text{bb}i}-T_{\text{ab}i}T_{\text{ab}i})](y_{0})\phi
(y_{0})$

\ $\times\lbrack\varrho_{jj}+$ $\overset{q}{\underset{\text{a}=1}{2\sum}%
}R_{j\text{a}j\text{a}}-3\overset{q}{\underset{\text{a,b=1}}{\sum}%
}(T_{\text{aa}j}T_{\text{bb}j}-T_{\text{ab}j}T_{\text{ab}j})]\}(y_{0}%
)\phi(y_{0})$

$+\frac{1}{48}R_{ijik}(y_{0})$ \ $[<H,j><H,k>](y_{0})\phi(y_{0})\qquad
\qquad\qquad\qquad$\ $L_{233}$

$+\frac{1}{432}R_{ijik}(y_{0})[2\varrho_{jk}+$ $\overset{q}{\underset{\text{a}%
=1}{4\sum}}R_{j\text{a}k\text{a}}-3\overset{q}{\underset{\text{a,b=1}}{\sum}%
}(T_{\text{aa}j}T_{\text{bb}k}-T_{\text{ab}j}T_{\text{ab}k}%
)-3\overset{q}{\underset{\text{a,b=1}}{\sum}}(T_{\text{aa}k}T_{\text{bb}%
j}-T_{\text{ab}k}T_{\text{ab}j}](y_{0})\phi(y_{0})$

$+\overset{n}{\underset{i,j=q+1}{\sum}}\frac{35}{128}<H,i>^{2}(y_{0}%
)<H,j>^{2}(y_{0})\phi(y_{0})\qquad\qquad\ \frac{1}{24}\frac{\partial^{4}%
\theta^{-\frac{1}{2}}}{\partial x_{i}^{2}\partial x_{j}^{2}}(y_{0})$

$+\frac{5}{192}\overset{n}{\underset{j=q+1}{\sum}}<H,j>^{2}(y_{0})[\tau
^{M}\ -3\tau^{P}+\ \underset{\text{a}=1}{\overset{\text{q}}{\sum}}%
\varrho_{\text{aa}}^{M}+\overset{q}{\underset{\text{a},\text{b}=1}{\sum}%
}R_{\text{abab}}^{M}](y_{0})\phi(y_{0})\qquad\ \ \ \ \ \ \ \ $

$+\frac{5}{192}\overset{n}{\underset{i=q+1}{\sum}}<H,i>^{2}(y_{0})[\tau
^{M}\ -3\tau^{P}+\ \underset{\text{a}=1}{\overset{\text{q}}{\sum}}%
\varrho_{\text{aa}}^{M}+\overset{q}{\underset{\text{a},\text{b}=1}{\sum}%
}R_{\text{abab}}^{M}](y_{0})\phi(y_{0})\qquad\qquad$

$+\frac{5}{192}\overset{n}{\underset{i,j=q+1}{\sum}}[<H,i><H,j>](y_{0}%
)\qquad\qquad\qquad\qquad\qquad\qquad\qquad\qquad$

$\times\lbrack2\varrho_{ij}+4\overset{q}{\underset{\text{a}=1}{\sum}%
}R_{i\text{a}j\text{a}}-3\overset{q}{\underset{\text{a,b=1}}{\sum}%
}(T_{\text{aa}i}T_{\text{bb}j}-T_{\text{ab}i}T_{\text{ab}j}%
)-3\overset{q}{\underset{\text{a,b=1}}{\sum}}(T_{\text{aa}j}T_{\text{bb}%
i}-T_{\text{ab}j}T_{\text{ab}i})](y_{0})\phi(y_{0})$

$+\frac{1}{96}\overset{n}{\underset{i,j=q+1}{\sum}}<H,j>(y_{0})[\{\nabla
_{i}\varrho_{ij}-2\varrho_{ij}<H,i>+\overset{q}{\underset{\text{a}=1}{\sum}%
}(\nabla_{i}R_{\text{a}i\text{a}j}-4R_{i\text{a}j\text{a}}<H,i>)\qquad$

$+4\overset{q}{\underset{\text{a,b=1}}{\sum}}R_{i\text{a}j\text{b}%
}T_{\text{ab}i}+2\overset{q}{\underset{\text{a,b,c=1}}{\sum}}(T_{\text{aa}%
i}T_{\text{bb}j}T_{\text{cc}i}-T_{\text{aa}i}T_{\text{bc}j}T_{\text{bc}%
i}-2T_{\text{bc}j}(T_{\text{aa}i}T_{\text{bc}i}-T_{\text{ab}i}T_{\text{ac}%
i}))\}$\qquad\qquad\qquad\ \ 

$+\{\nabla_{j}\varrho_{ii}-2\varrho_{ij}<H,i>+\overset{q}{\underset{\text{a}%
=1}{\sum}}(\nabla_{j}R_{\text{a}i\text{a}i}-4R_{i\text{a}j\text{a}}<H,i>)$

$+4\overset{q}{\underset{\text{a,b=1}}{\sum}}R_{j\text{a}i\text{b}%
}T_{\text{ab}i}+2\overset{q}{\underset{\text{a,b,c=1}}{\sum}}(T_{\text{aa}%
j}(T_{\text{bb}i}T_{\text{cc}i}-T_{\text{bc}i}T_{\text{bc}i})-2T_{\text{aa}%
j}T_{\text{bc}i}T_{\text{bc}i}+2T_{\text{ab}j}T_{\text{bc}i}T_{\text{ac}%
i})\}\qquad$

$+\{\nabla_{i}\varrho_{ij}-2\varrho_{ii}<H,j>+\overset{q}{\underset{\text{a}%
=1}{\sum}}(\nabla_{i}R_{\text{a}i\text{a}j}-4R_{i\text{a}i\text{a}%
}<H,j>)+4\overset{q}{\underset{\text{a,b=1}}{\sum}}R_{i\text{a}i\text{b}%
}T_{\text{ab}j}$

$+2\overset{q}{\underset{\text{a,b,c}=1}{\sum}}(T_{\text{aa}i}T_{\text{bb}%
i}T_{\text{cc}j}-3T_{\text{aa}i}T_{\text{bc}i}T_{\text{bc}j}+2T_{\text{ab}%
i}T_{\text{bc}i}T_{\text{ac}j})\}](y_{0})\phi(y_{0})$

$+\frac{1}{96}\overset{n}{\underset{i,j=q+1}{\sum}}<H,i>(y_{0})[\{\nabla
_{i}\varrho_{jj}-2\varrho_{ij}<H,j>+\overset{q}{\underset{\text{a}=1}{\sum}%
}(\nabla_{i}R_{\text{a}j\text{a}j}-4R_{i\text{a}j\text{a}}<H,j>)\qquad$

$+4\overset{q}{\underset{\text{a,b=1}}{\sum}}R_{i\text{a}j\text{b}%
}T_{\text{ab}j}+2\overset{q}{\underset{\text{a,b,c=1}}{\sum}}T_{\text{aa}%
i}(T_{\text{bb}j}T_{\text{cc}j}-T_{\text{bc}j}T_{\text{bc}j})-2T_{\text{aa}%
i}T_{\text{bc}j}T_{\text{bc}j}+2T_{\text{ab}i}T_{\text{bc}j}T_{\text{ac}%
j})\}(y_{0})\qquad$\qquad\qquad\qquad\qquad\qquad\ \ 

$+\{\nabla_{j}\varrho_{ij}-2\varrho_{ij}<H,j>+\overset{q}{\underset{\text{a}%
=1}{\sum}}(\nabla_{j}R_{\text{a}i\text{a}j}-4R_{j\text{a}i\text{a}}<H,j>)$

$+4\overset{q}{\underset{\text{a,b=1}}{\sum}}R_{j\text{a}i\text{b}%
}T_{\text{ab}j}+2\overset{q}{\underset{\text{a,b,c=1}}{\sum}}(T_{\text{aa}%
j}T_{\text{bb}i}T_{\text{cc}j}-T_{\text{ab}j}T_{\text{bc}i}T_{\text{ac}%
j}-2T_{\text{bc}i}(T_{\text{aa}j}T_{\text{bc}j}-T_{\text{ab}j}T_{\text{ac}%
j}))\}(y_{0})$

$+\{\nabla_{j}\varrho_{ij}-2\varrho_{jj}<H,i>+\overset{q}{\underset{\text{a}%
=1}{\sum}}(\nabla_{j}R_{\text{a}i\text{a}j}-4R_{j\text{a}j\text{a}%
}<H,i>)+4\overset{q}{\underset{\text{a,b=1}}{\sum}}R_{j\text{a}j\text{b}%
}T_{\text{ab}i}$

$+2\overset{q}{\underset{\text{a,b,c=1}}{\sum}}(T_{\text{aa}j}T_{\text{bb}%
j}T_{\text{cc}i}-3T_{\text{aa}j}T_{\text{bc}j}T_{\text{bc}i}+2T_{\text{ab}%
j}T_{\text{bc}j}T_{\text{ac}i})\}](y_{0})\phi(y_{0})$

$+\frac{1}{576}\overset{n}{\underset{i,j=q+1}{\sum}}[2\varrho_{ij}%
+4\overset{q}{\underset{\text{a}=1}{\sum}}R_{i\text{a}j\text{a}}%
-3\overset{q}{\underset{\text{a,b=1}}{\sum}}(T_{\text{aa}i}T_{\text{bb}%
j}-T_{\text{ab}i}T_{\text{ab}j})-3\overset{q}{\underset{\text{a,b=1}}{\sum}%
}(T_{\text{aa}j}T_{\text{bb}i}-T_{\text{ab}j}T_{\text{ab}i})]^{2}(y_{0}%
)\phi(y_{0})$

$+\frac{1}{288}[\tau^{M}\ -3\tau^{P}+\ \underset{\text{a}=1}{\overset{\text{q}%
}{\sum}}\varrho_{\text{aa}}^{M}+\overset{q}{\underset{\text{a},\text{b}%
=1}{\sum}}R_{\text{abab}}^{M}]^{2}(y_{0})\phi(y_{0})$

$-\ \frac{1}{288}\overset{n}{\underset{i,j=q+1}{\sum}}[$
$\overset{q}{\underset{\text{a=1}}{\sum}}\{-(\nabla_{ii}^{2}R_{j\text{a}%
j\text{a}}+\nabla_{jj}^{2}R_{i\text{a}i\text{a}}+4\nabla_{ij}^{2}%
R_{i\text{a}j\text{a}}+2R_{ij}R_{i\text{a}j\text{a}})\qquad A$

$+\overset{n}{\underset{p=q+1}{\sum}}\overset{q}{\underset{\text{a=1}}{\sum}%
}(R_{\text{a}iip}R_{\text{a}jjp}+R_{\text{a}jjp}R_{\text{a}iip}+R_{\text{a}%
ijp}R_{\text{a}ijp}+R_{\text{a}ijp}R_{\text{a}jip}+R_{\text{a}jip}%
R_{\text{a}ijp}+R_{\text{a}jip}R_{\text{a}jip})$

$+2\overset{q}{\underset{\text{a,b=1}}{\sum}}\nabla_{i}(R)_{\text{a}%
i\text{b}j}T_{\text{ab}j}+2\overset{q}{\underset{\text{a,b=1}}{\sum}}%
\nabla_{j}(R)_{\text{a}j\text{b}i}T_{\text{ab}i}%
+2\overset{q}{\underset{\text{a,b=1}}{\sum}}\nabla_{i}(R)_{\text{a}j\text{b}%
i}T_{\text{ab}j}+2\overset{q}{\underset{\text{a,b=1}}{\sum}}\nabla
_{i}(R)_{\text{a}j\text{b}j}T_{\text{ab}i}$

$+2\overset{q}{\underset{\text{a,b=1}}{\sum}}\nabla_{j}(R)_{\text{a}%
i\text{b}i}T_{\text{ab}j}+2\overset{q}{\underset{\text{a,b=1}}{\sum}}%
\nabla_{j}(R)_{\text{a}i\text{b}j}T_{\text{ab}i}$

$+\overset{n}{\underset{p=q+1}{\sum}}(-\frac{3}{5}\nabla_{ii}^{2}%
(R)_{jpjp}+\overset{n}{\underset{p=q+1}{\sum}}(-\frac{3}{5}\nabla_{jj}%
^{2}(R)_{ipip}+\overset{n}{\underset{p=q+1}{\sum}}(-\frac{3}{5}\nabla_{ij}%
^{2}(R)_{ipjp}+\overset{n}{\underset{p=q+1}{\sum}}(-\frac{3}{5}\nabla_{ij}%
^{2}(R)_{jpip}$

$+\overset{n}{\underset{p=q+1}{\sum}}(-\frac{3}{5}\nabla_{ji}^{2}%
(R)_{ipjp}+\overset{n}{\underset{p=q+1}{\sum}}(-\frac{3}{5}\nabla_{ji}%
^{2}(R)_{jpip}$

$+\frac{1}{5}\overset{n}{\underset{m,p=q+1}{%
{\textstyle\sum}
}}R_{ipim}R_{jpjm}+\frac{1}{5}\overset{n}{\underset{m,p=q+1}{%
{\textstyle\sum}
}}R_{jpjm}R_{ipim}+\frac{1}{5}\overset{n}{\underset{m,p=q+1}{%
{\textstyle\sum}
}}R_{ipjm}R_{ipjm}+\frac{1}{5}\overset{n}{\underset{m,p=q+1}{%
{\textstyle\sum}
}}R_{ipjm}R_{jpim}$

$+\frac{1}{5}\overset{n}{\underset{m,p=q+1}{%
{\textstyle\sum}
}}R_{jpim}R_{ipjm}+\frac{1}{5}\overset{n}{\underset{m,p=q+1}{%
{\textstyle\sum}
}}R_{jpim}R_{jpim}\}(y_{0})$

$+4\overset{q}{\underset{\text{a,b=1}}{\sum}}\{(\nabla_{i}(R)_{i\text{a}%
j\text{a}}-\overset{q}{\underset{\text{c=1}}{%
{\textstyle\sum}
}}R_{\text{a}i\text{c}i}T_{\text{ac}j})$ $T_{\text{bb}j}+4(\nabla
_{j}(R)_{j\text{a}i\text{a}}-\overset{q}{\underset{\text{c=1}}{%
{\textstyle\sum}
}}R_{\text{a}j\text{c}j}T_{\text{ac}i})$ $T_{\text{bb}i}+4(\nabla
_{i}(R)_{j\text{a}i\text{a}}-\overset{q}{\underset{\text{c=1}}{%
{\textstyle\sum}
}}R_{\text{a}i\text{c}j}T_{\text{ac}i})$ $T_{\text{bb}j}$ $4B\ $

$+4(\nabla_{i}(R)_{j\text{a}j\text{a}}-\overset{q}{\underset{\text{c=1}}{%
{\textstyle\sum}
}}R_{\text{a}i\text{c}j}T_{\text{ac}j})$ $T_{\text{bb}i}+4(\nabla
_{j}(R)_{i\text{a}i\text{a}}-\overset{q}{\underset{\text{c=1}}{%
{\textstyle\sum}
}}R_{\text{a}j\text{c}i}T_{\text{ac}i})$ $T_{\text{bb}j}+4(\nabla
_{j}(R)_{i\text{a}j\text{a}}-\overset{q}{\underset{\text{c=1}}{%
{\textstyle\sum}
}}R_{\text{a}j\text{c}i}T_{\text{ac}j})$ $T_{\text{bb}i}$

$-4\overset{q}{\underset{\text{a,b=1}}{\sum}}(\nabla_{i}(R)_{i\text{a}%
j\text{b}}-\overset{q}{\underset{\text{c=1}}{%
{\textstyle\sum}
}}R_{\text{b}r\text{c}s}T_{\text{ac}t})T_{\text{ab}j}%
-4\overset{q}{\underset{\text{a,b=1}}{\sum}}(\nabla_{j}(R)_{j\text{a}%
i\text{b}}-\overset{q}{\underset{\text{c=1}}{%
{\textstyle\sum}
}}R_{\text{b}j\text{c}j}T_{\text{ac}i})T_{\text{ab}i}$

$-4\overset{q}{\underset{\text{a,b=1}}{\sum}}(\nabla_{i}(R)_{j\text{a}%
i\text{b}}-\overset{q}{\underset{\text{c=1}}{%
{\textstyle\sum}
}}R_{\text{b}i\text{c}j}T_{\text{ac}i})T_{\text{ab}j}%
-4\overset{q}{\underset{\text{a,b=1}}{\sum}}(\nabla_{i}(R)_{j\text{a}%
j\text{b}}-\overset{q}{\underset{\text{c=1}}{%
{\textstyle\sum}
}}R_{\text{b}i\text{c}j}T_{\text{ac}j})T_{\text{ab}i}$

$-4\overset{q}{\underset{\text{a,b=1}}{\sum}}(\nabla_{j}(R)_{i\text{a}%
i\text{b}}-\overset{q}{\underset{\text{c=1}}{%
{\textstyle\sum}
}}R_{\text{b}j\text{c}i}T_{\text{ac}i})T_{\text{ab}j}%
-4\overset{q}{\underset{\text{a,b=1}}{\sum}}(\nabla_{j}(R)_{i\text{a}%
j\text{b}}-\overset{q}{\underset{\text{c=1}}{%
{\textstyle\sum}
}}R_{\text{b}j\text{c}i}T_{\text{ac}j})T_{\text{ab}i}\}](y_{0})$

$-\frac{1}{48}$ $[\frac{4}{9}\overset{q}{\underset{\text{a,b=1}}{\sum}%
}(\varrho_{\text{aa}}-\overset{q}{\underset{\text{c}=1}{\sum}}R_{\text{acac}%
})(\varrho_{\text{bb}}-\overset{q}{\underset{\text{d}=1}{\sum}}R_{\text{bdbd}%
})+\frac{8}{9}\overset{n}{\underset{i,j=q+1}{\sum}}%
\overset{q}{\underset{\text{a,b}=1}{\sum}}(R_{i\text{a}j\text{a}}%
R_{i\text{b}j\text{b}})\qquad3C$

$+\frac{2}{9}\overset{q}{\underset{\text{a}=1}{\sum}}(\varrho_{\text{aa}}%
^{M}-\varrho_{\text{aa}}^{P})(\tau^{M}-\overset{q}{\underset{\text{c}=1}{\sum
}}\varrho_{\text{cc}}^{M})+\frac{4}{9}\overset{n}{\underset{i,j=q+1}{\sum}%
}\overset{q}{\underset{\text{a}=1}{\sum}}R_{i\text{a}j\text{a}}\varrho_{ij}\ $

$\ +\frac{2}{9}\overset{q}{\underset{\text{b}=1}{\sum}}(\varrho_{\text{bb}%
}^{M}-\varrho_{\text{bb}}^{P})(\tau^{M}-\overset{q}{\underset{\text{c}%
=1}{\sum}}\varrho_{\text{cc}}^{M})+\frac{4}{9}%
\overset{n}{\underset{i,j=q+1}{\sum}}\overset{q}{\underset{\text{b}=1}{\sum}%
}R_{i\text{b}j\text{b}}\varrho_{ij}\ $

$+\frac{1}{9}(\tau^{M}-\overset{q}{\underset{\text{a=1}}{\sum}}\varrho
_{\text{aa}})(\tau^{M}-\overset{q}{\underset{\text{b=1}}{\sum}}\varrho
_{\text{bb}})+\frac{2}{9}(\left\Vert \varrho^{M}\right\Vert ^{2}%
-\overset{q}{\underset{\text{a,b}=1}{\sum}}\varrho_{\text{ab}})$

$-\overset{n}{\underset{i,j=q+1}{\sum}}\overset{q}{\underset{\text{a,b}%
=1}{\sum}}R_{i\text{a}i\text{b}}R_{j\text{a}j\text{b}}\ -\frac{1}%
{2}\overset{n}{\underset{i,j=q+1}{\sum}}\overset{q}{\underset{\text{a,b}%
=1}{\sum}}R_{i\text{a}j\text{b}}^{2}-\overset{n}{\underset{i,j=q+1}{\sum}%
}\overset{q}{\underset{\text{a,b}=1}{\sum}}R_{i\text{a}j\text{b}}%
R_{j\text{a}i\text{b}}-\frac{1}{2}\overset{n}{\underset{i,j=q+1}{\sum}%
}\overset{q}{\underset{\text{a,b}=1}{\sum}}R_{j\text{a}i\text{b}}^{2}$

$-\frac{1}{9}\overset{n}{\underset{i,j,p,m=q+1}{\sum}}R_{ipim}R_{jpjm}%
\ -\frac{1}{18}\overset{n}{\underset{i,j,p,m=q+1}{\sum}}R_{ipjm}^{2}-\frac
{1}{9}\overset{n}{\underset{i,j,p,m=q+1}{\sum}}R_{ipjm}R_{jpim}-\frac{1}%
{18}\overset{n}{\underset{i,j,p,m=q+1}{\sum}}R_{jpim}^{2}$

$-\frac{1}{3}\overset{q}{\underset{\text{a}=1}{\sum}}%
\overset{n}{\underset{i,j,p=q+1}{\sum}}R_{i\text{a}ip}R_{j\text{a}jp}-\frac
{1}{6}\overset{q}{\underset{\text{a}=1}{\sum}}%
\overset{n}{\underset{i,j,p=q+1}{\sum}}R_{i\text{a}jp}^{2}-\frac{1}%
{3}\overset{q}{\underset{\text{a}=1i,j,}{\sum}}%
\overset{n}{\underset{p=q+1}{\sum}}R_{i\text{a}jp}R_{j\text{a}ip}-\frac{1}%
{6}\overset{q}{\underset{\text{a}=1}{\sum}}%
\overset{n}{\underset{i,j,p=q+1}{\sum}}R_{j\text{a}ip}^{2}$

$-\frac{1}{3}\overset{q}{\underset{\text{b}=1i,j,}{\sum}}%
\overset{n}{\underset{p=q+1}{\sum}}R_{i\text{b}ip}R_{j\text{b}jp}-\frac{1}%
{6}\overset{q}{\underset{\text{b}=1}{\sum}}%
\overset{n}{\underset{i,j,p=q+1}{\sum}}R_{i\text{b}jp}^{2}-\frac{1}%
{3}\overset{q}{\underset{\text{b}=1}{\sum}}%
\overset{n}{\underset{i.j,p=q+1}{\sum}}R_{i\text{b}jp}R_{j\text{b}ip}-\frac
{1}{6}\overset{q}{\underset{\text{b}=1}{\sum}}%
\overset{n}{\underset{i,j,p=q+1}{\sum}}R_{j\text{b}ip}^{2}](y_{0})\phi(y_{0})$

$-\frac{1}{48}$ $\overset{q}{\underset{\text{a,b,c=1}}{\sum}}[$
$-\overset{n}{\underset{i=q+1}{\sum}}R_{i\text{a}i\text{a}}(R_{\text{bcbc}%
}^{P}-R_{\text{bcbc}}^{M})$ $-\overset{n}{\underset{j=q+1}{\sum}}%
R_{j\text{a}j\text{a}}(R_{\text{bcbc}}^{P}-R_{\text{bcbc}}^{M})\qquad\qquad6D$

$+\overset{n}{\underset{i=q+1}{\sum}}R_{i\text{a}i\text{b}}(R_{\text{acbc}%
}^{P}-R_{\text{acbc}}^{M})\ -\overset{n}{\underset{i=q+1}{\sum}}%
R_{i\text{a}i\text{c}}(R_{\text{abbc}}^{P}-R_{\text{abbc}}^{M})$

$+\overset{n}{\underset{j=q+1}{\sum}}R_{j\text{a}j\text{b}}(R_{\text{acbc}%
}^{P}-R_{\text{acbc}}^{M})$\ $-\overset{n}{\underset{j=q+1}{\sum}}%
R_{j\text{a}j\text{c}}(R_{\text{abbc}}^{P}-R_{\text{abbc}}^{M})$

$+\underset{i,j=q+1}{\overset{n}{\sum}}$ $-R_{i\text{a}j\text{a}}%
(T_{\text{bb}i}T_{\text{cc}j}$ $-T_{\text{bc}i}T_{\text{bc}j})$
$-\underset{i,j=q+1}{\overset{n}{\sum}}R_{i\text{a}j\text{a}}(T_{\text{bb}%
j}T_{\text{cc}i}$ $-T_{\text{bc}j}T_{\text{bc}i})$

$+$ $\underset{i,j=q+1}{\overset{n}{\sum}}$ $-R_{j\text{a}i\text{a}%
}(T_{\text{bb}i}T_{\text{cc}j}$ $-T_{\text{bc}i}T_{\text{bc}j})$
$-\underset{i,j=q+1}{\overset{n}{\sum}}R_{j\text{a}i\text{a}}(T_{\text{bb}%
j}T_{\text{cc}i}$ $-T_{\text{bc}j}T_{\text{bc}i})$

$+\underset{i,j=q+1}{\overset{n}{\sum}}\ R_{i\text{a}j\text{b}}(T_{\text{ab}%
i}T_{\text{cc}j}-T_{\text{bc}i}T_{\text{ac}j}%
)\ +\underset{i,j=q+1}{\overset{n}{\sum}}\ R_{i\text{a}j\text{b}}%
(T_{\text{ab}j}T_{\text{cc}i}-T_{\text{bc}j}T_{\text{ac}i})$

$+\underset{i,j=q+1}{\overset{n}{\sum}}\ R_{j\text{a}i\text{ib}}%
(T_{\text{ab}i}T_{\text{cc}j}-T_{\text{bc}i}T_{\text{ac}j}%
)\ +\underset{i,j=q+1}{\overset{n}{\sum}}\ R_{j\text{a}i\text{b}}%
(T_{\text{ab}j}T_{\text{cc}i}-T_{\text{bc}j}T_{\text{ac}i})\qquad$

$+\underset{i,j=q+1}{\overset{n}{\sum}}-R_{i\text{a}j\text{c}}(T_{\text{ab}%
i}T_{\text{bc}j}-T_{\text{ac}i}T_{\text{bb}j}%
)-\underset{i,j=q+1}{\overset{n}{\sum}}R_{i\text{a}j\text{c}}(T_{\text{ba}%
j}T_{\text{bc}i}-T_{\text{ac}j}T_{\text{bb}i})$

$+\underset{i,j=q+1}{\overset{n}{\sum}}-R_{j\text{a}i\text{c}}(T_{\text{ba}%
i}T_{\text{bc}j}-T_{\text{ac}i}T_{\text{bb}j}%
)-\underset{i,j=q+1}{\overset{n}{\sum}}R_{j\text{a}i\text{c}}(T_{\text{ba}%
j}T_{\text{bc}i}-T_{\text{ac}j}T_{\text{bb}i})](y_{0})\phi(y_{0})$

$+\frac{1}{144}\underset{p=q+1}{\overset{n}{\sum}}%
[\underset{i=q+1}{\overset{n}{\sum}}\overset{q}{\underset{\text{b,c=1}}{\sum}%
}R_{ipip}(R_{\text{bcbc}}^{P}-R_{\text{bcbc}}^{M}%
)+\underset{j=q+1}{\overset{n}{\sum}}$ $\overset{q}{\underset{\text{b,c=1}%
}{\sum}}R_{jpjp}(R_{\text{bcbc}}^{P}-R_{\text{bcbc}}^{M})](y_{0})\phi(y_{0})$

$+\frac{1}{72}\underset{i,j,p=q+1}{\overset{n}{\sum}}%
\overset{q}{\underset{\text{b,c=1}}{\sum}}[R_{ipjp}(T_{\text{bb}i}%
T_{\text{cc}j}-T_{\text{bc}i}T_{\text{bc}j})+R_{ipjp}(T_{\text{bb}%
j}T_{\text{cc}i}-T_{\text{bc}j}T_{\text{bc}i})](y_{0})\phi(y_{0})\qquad$

$-\frac{1}{288}\underset{i,j=q+1}{\overset{n}{\sum}}[T_{\text{aa}%
i}T_{\text{bb}j}(T_{\text{cc}i}T_{\text{dd}j}-T_{\text{cd}i}T_{\text{dc}%
j})+T_{\text{aa}i}T_{\text{bb}j}(T_{\text{cc}j}T_{\text{dd}i}-T_{\text{cd}%
j}T_{\text{dc}i})\qquad E$

$+T_{\text{aa}j}T_{\text{bb}i}(T_{\text{cc}i}T_{\text{dd}j}-T_{\text{cd}%
i}T_{\text{dc}j})+T_{\text{aa}j}T_{\text{bb}i}(T_{\text{cc}j}T_{\text{dd}%
i}-T_{\text{cd}j}T_{\text{dc}i})](y_{0})\phi(y_{0})$

$+\frac{1}{288}\underset{i,j=q+1}{\overset{n}{\sum}}[T_{\text{aa}%
i}T_{\text{bc}j}(T_{\text{bc}i}T_{\text{dd}j}-T_{\text{bd}i}T_{\text{cd}%
j})+T_{\text{aa}i}T_{\text{bc}j}(T_{\text{bc}j}T_{\text{dd}i}-T_{\text{bd}%
j}T_{\text{cd}i})$

$+T_{\text{aa}j}T_{\text{bc}i}(T_{\text{bc}i}T_{\text{dd}j}-T_{\text{bd}%
i}T_{\text{cd}j})+T_{\text{aa}j}T_{\text{bc}i}(T_{\text{bc}j}T_{\text{dd}%
i}-T_{\text{bd}j}T_{\text{cd}i})](y_{0})\phi(y_{0})$

$-\frac{1}{288}\underset{i,j=q+1}{\overset{n}{\sum}}[T_{\text{aa}%
i}T_{\text{bd}j}(T_{\text{bc}i}T_{\text{cd}j}-T_{\text{bd}i}T_{\text{cc}%
j})+T_{\text{aa}i}T_{\text{bd}j}(T_{\text{bc}j}T_{\text{cd}i}-T_{\text{bd}%
j}T_{\text{cc}i})$

$+T_{\text{aa}j}T_{\text{bd}i}(T_{\text{bc}i}T_{\text{cd}j}-T_{\text{bd}%
i}T_{\text{cc}j})+T_{\text{aa}j}T_{\text{bd}i}(T_{\text{bc}j}T_{\text{cd}%
i}-T_{\text{bd}j}T_{\text{cc}i})](y_{0})\phi(y_{0})\qquad$

$+\frac{1}{288}\underset{i,j=q+1}{\overset{n}{\sum}}[T_{\text{ab}%
i}T_{\text{ab}j}(T_{\text{cc}i}T_{\text{dd}j}-T_{\text{cd}i}T_{\text{dc}%
j})+T_{\text{ab}i}T_{\text{ab}j}(T_{\text{cc}j}T_{\text{dd}i}-T_{\text{cd}%
j}T_{\text{dc}i})$

$+T_{\text{ab}j}T_{\text{ab}i}(T_{\text{cc}i}T_{\text{dd}j}-T_{\text{cd}%
i}T_{\text{dc}j})+T_{\text{ab}j}T_{\text{ab}i}(T_{\text{cc}j}T_{\text{dd}%
i}-T_{\text{cd}j}T_{\text{dc}i})](y_{0})\phi(y_{0})$

$-\frac{1}{288}\underset{i,j=q+1}{\overset{n}{\sum}}[T_{\text{ab}%
i}T_{\text{bc}j}(T_{\text{ac}i}T_{\text{dd}j}-T_{\text{ad}i}T_{\text{cd}%
j})+T_{\text{ab}i}T_{\text{bc}j}(T_{\text{ac}j}T_{\text{dd}i}-T_{\text{ad}%
j}T_{\text{cd}i})$

$+T_{\text{ab}j}T_{\text{bc}i}(T_{\text{ac}i}T_{\text{dd}j}-T_{\text{ad}%
i}T_{\text{cd}j})+T_{\text{ab}j}T_{\text{bc}i}(T_{\text{ac}j}T_{\text{dd}%
i}-T_{\text{ad}j}T_{\text{cd}i})](y_{0})\phi(y_{0})$

$+\frac{1}{288}\underset{i,j=q+1}{\overset{n}{\sum}}[T_{\text{ab}%
i}T_{\text{bd}j}(T_{\text{ac}i}T_{\text{cd}j}-T_{\text{ad}i}T_{\text{cc}%
j})+T_{\text{ab}i}T_{\text{bd}j}(T_{\text{ac}j}T_{\text{cd}i}-T_{\text{ad}%
j}T_{\text{cc}i})$

$+T_{\text{ab}i}T_{\text{bd}j}(T_{\text{ac}j}T_{\text{cd}i}-T_{\text{ad}%
j}T_{\text{cc}i})+T_{\text{ab}j}T_{\text{bd}i}(T_{\text{ac}j}T_{\text{cd}%
i}-T_{\text{ad}j}T_{\text{cc}i})](y_{0})\phi(y_{0})$

$-\ \frac{1}{288}\underset{i,j=q+1}{\overset{n}{\sum}}[T_{\text{ac}%
i}T_{\text{ab}j}(T_{\text{bc}i}T_{\text{dd}j}-T_{\text{bd}i}T_{\text{dc}%
j})+T_{\text{ac}i}T_{\text{ab}j}(T_{\text{bc}j}T_{\text{dd}i}-T_{\text{bd}%
j}T_{\text{dc}i})$

$+T_{\text{ac}j}T_{\text{ab}i}(T_{\text{bc}i}T_{\text{dd}j}-T_{\text{bd}%
i}T_{\text{dc}j})+T_{\text{ac}j}T_{\text{ab}i}(T_{\text{bc}j}T_{\text{dd}%
i}-T_{\text{bd}j}T_{\text{dc}i})](y_{0})\phi(y_{0})$

$+\ \frac{1}{288}\underset{i,j=q+1}{\overset{n}{\sum}}[T_{\text{ac}%
i}T_{\text{bb}j}(T_{\text{ac}i}T_{\text{dd}j}-T_{\text{ad}i}T_{\text{cd}%
j})+T_{\text{ac}i}T_{\text{bb}j}(T_{\text{ac}j}T_{\text{dd}i}-T_{\text{ad}%
j}T_{\text{cd}i})$

$+T_{\text{ac}j}T_{\text{bb}i}(T_{\text{ac}i}T_{\text{dd}j}-T_{\text{ad}%
i}T_{\text{cd}i})+T_{\text{ac}j}T_{\text{bb}i}(T_{\text{ac}j}T_{\text{dd}%
i}-T_{\text{ad}j}T_{\text{cd}i})](y_{0})\phi(y_{0})$

$-\ \frac{1}{288}\underset{i,j=q+1}{\overset{n}{\sum}}[T_{\text{ac}%
i}T_{\text{bd}j}(T_{\text{ac}i}T_{\text{bd}j}-T_{\text{ad}i}T_{\text{bc}%
j})+T_{\text{ac}i}T_{\text{bd}j}(T_{\text{ac}j}T_{\text{bd}i}-T_{\text{ad}%
j}T_{\text{bc}i})$

$+T_{\text{ac}j}T_{\text{bd}i}(T_{\text{ac}i}T_{\text{bd}j}-T_{\text{ad}%
i}T_{\text{bc}j})+T_{\text{ac}j}T_{\text{bd}i}(T_{\text{ac}j}T_{\text{bd}%
i}-T_{\text{ad}j}T_{\text{bc}i})](y_{0})\phi(y_{0})$

$+\frac{1}{288}\underset{i,j=q+1}{\overset{n}{\sum}}[T_{\text{ad}%
i}T_{\text{ab}j}(T_{\text{bc}i}T_{\text{cd}j}-T_{\text{bd}i}T_{\text{cc}%
j})+T_{\text{ad}i}T_{\text{ab}j}(T_{\text{bc}j}T_{\text{cd}i}-T_{\text{bd}%
j}T_{\text{cc}i})$

$+T_{\text{ad}j}T_{\text{ab}i}(T_{\text{bc}i}T_{\text{cd}j}-T_{\text{bd}%
i}T_{\text{cc}j})+T_{\text{ad}j}T_{\text{ab}i}(T_{\text{bc}j}T_{\text{cd}%
i}-T_{\text{bd}j}T_{\text{cc}i})](y_{0})\phi(y_{0})$

$-\ \frac{1}{288}\underset{i,j=q+1}{\overset{n}{\sum}}[T_{\text{ad}%
i}T_{\text{bb}j}(T_{\text{ac}i}T_{\text{cd}j}-T_{\text{ad}i}T_{\text{cc}%
j})+T_{\text{ad}i}T_{\text{bb}j}(T_{\text{ac}j}T_{\text{cd}i}-T_{\text{ad}%
j}T_{\text{cc}i})$

$+T_{\text{ad}j}T_{\text{bb}i}(T_{\text{ac}i}T_{\text{cd}j}-T_{\text{ad}%
i}T_{\text{cc}j})+T_{\text{ad}j}T_{\text{bb}i}(T_{\text{ac}j}T_{\text{cd}%
i}-T_{\text{ad}j}T_{\text{cc}i})](y_{0})\phi(y_{0})$

$+\ \frac{1}{288}\underset{i,j=q+1}{\overset{n}{\sum}}[T_{\text{ad}%
i}T_{\text{bc}j}(T_{\text{ac}i}T_{\text{bd}j}-T_{\text{ad}i}T_{\text{bc}%
j})+T_{\text{ad}i}T_{\text{bc}j}(T_{\text{ac}j}T_{\text{bd}i}-T_{\text{ad}%
j}T_{\text{bc}i})$

$+T_{\text{ad}j}T_{\text{bc}i}(T_{\text{ac}i}T_{\text{bd}j}-T_{\text{ad}%
i}T_{\text{bc}j})+T_{\text{ad}j}T_{\text{bc}i}(T_{\text{ac}j}T_{\text{bd}%
i}-T_{\text{ad}j}T_{\text{bc}i})](y_{0})\phi(y_{0})$

$-\ \frac{1}{144}[(R_{\text{cdcd}}^{P}-R_{\text{cdcd}}^{M})(R_{\text{abab}%
}^{P}-R_{\text{abab}}^{M})](y_{0})\phi(y_{0})$

$+\frac{1}{144}[(R_{\text{bdcd}}^{P}-R_{\text{bdcd}}^{M})(R_{\text{abac}}%
^{P}-R_{\text{abac}}^{M})](y_{0})\phi(y_{0})$

$\ +\ \frac{1}{144}[(R_{\text{bcdc}}^{P}-R_{\text{bcdc}}^{M})(R_{\text{abad}%
}^{P}-R_{\text{abad}}^{M})](y_{0})\phi(y_{0})$

$\ -\ \frac{1}{144}[(R_{\text{adcd}}^{P}-R_{\text{adcd}}^{M})(R_{\text{abbc}%
}^{P}-R_{\text{abbc}}^{M})](y_{0})\phi(y_{0})\qquad$

$\ +\ \frac{1}{144}[(R_{\text{acdc}}^{P}-R_{\text{acdc}}^{M})(R_{\text{abdb}%
}^{P}-R_{\text{abdb}}^{M})](y_{0})\phi(y_{0})$

$\ -\ \frac{1}{576}[(R_{\text{abcd}}^{P}-R_{\text{abcd}}^{M})]^{2}(y_{0}%
)\phi(y_{0})$

$-\frac{1}{144}<H,i>(y_{0})<H,j>(y_{0})\qquad\qquad\qquad\qquad\qquad
\qquad\qquad L_{3}$

$\times\lbrack2\varrho_{ij}+$ $\overset{q}{\underset{\text{a}=1}{4\sum}%
}R_{i\text{a}j\text{a}}-3\overset{q}{\underset{\text{a,b=1}}{\sum}%
}(T_{\text{aa}i}T_{\text{bb}j}-T_{\text{ab}i}T_{\text{ab}j}%
)-3\overset{q}{\underset{\text{a,b=1}}{\sum}}(T_{\text{aa}j}T_{\text{bb}%
i}-T_{\text{ab}j}T_{\text{ab}i}](y_{0})\phi(y_{0})$

$-\frac{1}{16}[<H,i>^{2}(y_{0})<H,j>^{2}](y_{0})\phi(y_{0})$

$-\frac{1}{144}<H,i>(y_{0})<H,j>(y_{0})$

$\times\lbrack2\varrho_{ij}+$ $\overset{q}{\underset{\text{a}=1}{4\sum}%
}R_{i\text{a}j\text{a}}-3\overset{q}{\underset{\text{a,b=1}}{\sum}%
}(T_{\text{aa}i}T_{\text{bb}j}-T_{\text{ab}i}T_{\text{ab}j}%
)-3\overset{q}{\underset{\text{a,b=1}}{\sum}}(T_{\text{aa}j}T_{\text{bb}%
i}-T_{\text{ab}j}T_{\text{ab}i}](y_{0})\phi(y_{0})$

$-\frac{1}{72}<H,i>(y_{0})<H,k>(y_{0})R_{jijk}(y_{0})\phi(y_{0})$

$-\frac{1}{16}<H,i>^{2}(y_{0})<H,j>^{2}(y_{0})\phi(y_{0})$

$-\frac{1}{72}<H,i>^{2}(y_{0})[\varrho_{jj}+$ $\overset{q}{\underset{\text{a}%
=1}{2\sum}}R_{j\text{a}j\text{a}}-3\overset{q}{\underset{\text{a,b=1}}{\sum}%
}(T_{\text{aa}j}T_{\text{bb}j}-T_{\text{ab}j}T_{\text{ab}j})](y_{0})\phi
(y_{0})$

$+\frac{5}{32}[<H,i>^{2}<H,j>^{2}](y_{0})\phi(y_{0})$

$+\frac{1}{48}<H,i>(y_{0})<H,j>$

$\times\lbrack2\varrho_{ij}+$ $\overset{q}{\underset{\text{a}=1}{4\sum}%
}R_{i\text{a}j\text{a}}-3\overset{q}{\underset{\text{a,b=1}}{\sum}%
}(T_{\text{aa}i}T_{\text{bb}j}-T_{\text{ab}i}T_{\text{ab}j}%
)-3\overset{q}{\underset{\text{a,b=1}}{\sum}}(T_{\text{aa}j}T_{\text{bb}%
i}-T_{\text{ab}j}T_{\text{ab}i}](y_{0})\phi(y_{0})$

$+\frac{1}{48}<H,i>^{2}(y_{0})[\tau^{M}\ -3\tau^{P}+\ \underset{\text{a}%
=1}{\overset{\text{q}}{\sum}}\varrho_{\text{aa}}^{M}+$
$\overset{q}{\underset{\text{a},\text{b}=1}{\sum}}R_{\text{abab}}^{M}$
$](y_{0})\phi(y_{0})$

$+\frac{1}{144}<H,i>(y_{0})[\nabla_{i}\varrho_{jj}-2\varrho_{ij}%
<H,j>+\overset{q}{\underset{\text{a}=1}{\sum}}(\nabla_{i}R_{\text{a}%
j\text{a}j}-4R_{i\text{a}j\text{a}}<H,j>)$

$+4\overset{q}{\underset{\text{a,b=1}}{\sum}}R_{i\text{a}j\text{b}%
}T_{\text{ab}j}+2\overset{q}{\underset{\text{a,b,c=1}}{\sum}}(T_{\text{aa}%
i}T_{\text{bb}j}T_{\text{cc}j}-3T_{\text{aa}i}T_{\text{bc}j}T_{\text{bc}%
j}+2T_{\text{ab}i}T_{\text{bc}j}T_{\text{ca}j})](y_{0})\phi(y_{0})$%
\qquad\qquad\qquad\qquad\qquad\ \ 

$+\frac{1}{144}<H,i>(y_{0})[\nabla_{j}\varrho_{ij}-2\varrho_{ij}%
<H,j>+\overset{q}{\underset{\text{a}=1}{\sum}}(\nabla_{j}R_{\text{a}%
i\text{a}j}-4R_{j\text{a}i\text{a}}<H,j>)$

$+4\overset{q}{\underset{\text{a,b=1}}{\sum}}R_{j\text{a}i\text{b}%
}T_{\text{ab}j}+2\overset{q}{\underset{\text{a,b,c=1}}{\sum}}(T_{\text{aa}%
j}T_{\text{bb}i}T_{\text{cc}j}-3T_{\text{aa}j}T_{\text{bc}i}T_{\text{bc}%
j}+2T_{\text{ab}j}T_{\text{bc}i}T_{\text{ac}j})](y_{0})\phi(y_{0})$

$+\frac{1}{144}<H,i>(y_{0})[\nabla_{j}\varrho_{ij}-2\varrho_{jj}%
<H,i>+\overset{q}{\underset{\text{a}=1}{\sum}}(\nabla_{j}R_{\text{a}%
i\text{a}j}-4R_{j\text{a}j\text{a}}<H,i>)$

$+4\overset{q}{\underset{\text{a,b=1}}{\sum}}R_{j\text{a}j\text{b}%
}T_{\text{ab}i}+2\overset{q}{\underset{\text{a,b,c=1}}{\sum}}(T_{\text{aa}%
j}T_{\text{bb}j}T_{\text{cc}i}-3T_{\text{aa}j}T_{\text{bc}j}T_{\text{bc}%
i}+2T_{\text{ab}j}T_{\text{bc}j}T_{\text{ac}i})](y_{0})\phi(y_{0})$

$-\frac{1}{192}<H,i>^{2}<H,j>^{2}(y_{0})\phi(y_{0})\qquad\qquad\qquad
\qquad\qquad\qquad\qquad$I$_{3213}\ $

$-\frac{1}{288}<H,i>^{2}(y_{0})[\tau^{M}\ -3\tau^{P}+\ \underset{\text{a}%
=1}{\overset{\text{q}}{\sum}}\varrho_{\text{aa}}^{M}+$
$\overset{q}{\underset{\text{a},\text{b}=1}{\sum}}R_{\text{abab}}^{M}$
$](y_{0})\phi(y_{0})$

$-\frac{1}{288}<H,i>(y_{0})<H,j>(y_{0})$

$\times\lbrack2\varrho_{ij}+$ $\overset{q}{\underset{\text{a}=1}{4\sum}%
}R_{i\text{a}j\text{a}}-3\overset{q}{\underset{\text{a,b=1}}{\sum}%
}(T_{\text{aa}i}T_{\text{bb}j}-T_{\text{ab}i}T_{\text{ab}j}%
)-3\overset{q}{\underset{\text{a,b=1}}{\sum}}(T_{\text{aa}j}T_{\text{bb}%
i}-T_{\text{ab}j}T_{\text{ab}i}](y_{0})\phi(y_{0})$

$+\frac{1}{144}<H,i>(y_{0})<H,k>(y_{0})R_{jijk}(y_{0})\phi(y_{0})$

$+\frac{1}{144}<H,i>^{2}(y_{0})[\varrho_{jj}+$ $\overset{q}{\underset{\text{a}%
=1}{2\sum}}R_{j\text{a}j\text{a}}-3\overset{q}{\underset{\text{a,b=1}}{\sum}%
}(T_{\text{aa}j}T_{\text{bb}j}-T_{\text{ab}j}T_{\text{ab}j})](y_{0})\phi
(y_{0})$

$-\frac{1}{288}<H,i>(y_{0})[\nabla_{i}\varrho_{jj}-2\varrho_{ij}%
<H,j>+\overset{q}{\underset{\text{a}=1}{\sum}}(\nabla_{i}R_{\text{a}%
j\text{a}j}-4R_{i\text{a}j\text{a}}<H,j>)$

$+4\overset{q}{\underset{\text{a,b=1}}{\sum}}R_{i\text{a}j\text{b}%
}T_{\text{ab}j}+2\overset{q}{\underset{\text{a,b,c=1}}{\sum}}(T_{\text{aa}%
i}T_{\text{bb}j}T_{\text{cc}j}-3T_{\text{aa}i}T_{\text{bc}j}T_{\text{bc}%
j}+2T_{\text{ab}i}T_{\text{bc}j}T_{\text{ca}j})](y_{0})\phi(y_{0})$%
\qquad\qquad\qquad\qquad\qquad\ \ 

$-\frac{1}{288}<H,i>(y_{0})[\nabla_{j}\varrho_{ij}-2\varrho_{ij}%
<H,j>+\overset{q}{\underset{\text{a}=1}{\sum}}(\nabla_{j}R_{\text{a}%
i\text{a}j}-4R_{j\text{a}i\text{a}}<H,j>)$

$+4\overset{q}{\underset{\text{a,b=1}}{\sum}}R_{j\text{a}i\text{b}%
}T_{\text{ab}j}+2\overset{q}{\underset{\text{a,b,c=1}}{\sum}}(T_{\text{aa}%
j}T_{\text{bb}i}T_{\text{cc}j}-3T_{\text{aa}j}T_{\text{bc}i}T_{\text{bc}%
j}+2T_{\text{ab}j}T_{\text{bc}i}T_{\text{ac}j})](y_{0})\phi(y_{0})$

$-\frac{1}{288}<H,i>(y_{0})[\nabla_{j}\varrho_{ij}-2\varrho_{jj}%
<H,i>+\overset{q}{\underset{\text{a}=1}{\sum}}(\nabla_{j}R_{\text{a}%
i\text{a}j}-4R_{j\text{a}j\text{a}}<H,i>)+4\overset{q}{\underset{\text{a,b=1}%
}{\sum}}R_{j\text{a}j\text{b}}T_{\text{ab}i}$

$+2\overset{q}{\underset{\text{a,b,c=1}}{\sum}}(T_{\text{aa}j}T_{\text{bb}%
j}T_{\text{cc}i}-3T_{\text{aa}j}T_{\text{bc}j}T_{\text{bc}i}+2T_{\text{ab}%
j}T_{\text{bc}j}T_{\text{ac}i})](y_{0})\phi(y_{0})$

$+\frac{1}{24}[\left\Vert \text{X}\right\Vert _{M}^{2}+\operatorname{div}%
$X$_{M}-\left\Vert \text{X}\right\Vert _{P}^{2}-\operatorname{div}X_{P}%
](y_{0})[\left\Vert \text{X}\right\Vert _{M}^{2}-\operatorname{div}$%
X$_{M}-\left\Vert \text{X}\right\Vert _{P}^{2}+\operatorname{div}$%
X$_{P}](y_{0})\phi(y_{0})\qquad$\ I$_{3212}$

$+\frac{1}{6}X_{i}(y_{0})T_{\text{ab}i}(y_{0})T_{\text{ab}j}(y_{0})X_{j}%
(y_{0})+\frac{1}{3}\perp_{\text{a}ij}(y_{0})X_{i}(y_{0})[\frac{\partial X_{j}%
}{\partial x_{\text{a}}}-\perp_{\text{a}jk}X_{k}](y_{0})\phi(y_{0})\qquad
$I$_{32122}\qquad Q_{1}$

$+$ $\frac{2}{3}X_{i}(y_{0})X_{j}(y_{0})\frac{\partial X_{j}}{\partial
x_{\text{a}}}(y_{0})-\frac{1}{6}X_{i}(y_{0})\frac{\partial^{2}X_{j}}{\partial
x_{\text{a}}\partial x_{j}}(y_{0})\phi(y_{0})\qquad\qquad Q_{2}$

$-\frac{1}{12}X_{i}(y_{0})\frac{\partial^{2}X_{i}}{\partial x_{\text{a}}^{2}%
}(y_{0})+\frac{1}{12}X_{i}^{2}(y_{0})[\operatorname{div}X_{M}-\left\Vert
X\right\Vert _{M}^{2}+\left\Vert X\right\Vert _{P}^{2}-\operatorname{div}%
X_{P}-$ $<H,j>X_{j}](y_{0})\phi(y_{0})$

$+\frac{1}{6}X_{i}(y_{0})X_{j}(y_{0})\frac{\partial X_{i}}{\partial x_{j}%
}(y_{0})\phi(y_{0})+$ $\frac{1}{18}X_{i}(y_{0})X_{k}(y_{0})R_{jijk}(y_{0}%
)\phi(y_{0})-\frac{1}{12}X_{i}(y_{0})\frac{\partial^{2}X_{i}}{\partial
x_{j}^{2}}(y_{0})\phi(y_{0})$

$+\frac{1}{12}[R_{\text{a}i\text{a}k}-\underset{\text{c=1}}{\overset{\text{q}%
}{\sum}}T_{\text{ac}i}T_{\text{ac}k}-\perp_{\text{a}ik}\perp_{\text{a}%
jk}](y_{0})X_{k}(y_{0})\phi(y_{0})+\frac{1}{18}R_{ijkj}(y_{0})X_{i}%
(y_{0})X_{k}(y_{0})\phi(y_{0})$

$+\frac{1}{12}<H,j>(y_{0})X_{i}(y_{0})[X_{i}X_{j}-\frac{1}{2}\left(
\frac{\partial X_{j}}{\partial x_{i}}+\frac{\partial X_{i}}{\partial x_{j}%
}\right)  ](y_{0})\phi(y_{0})\qquad\qquad\qquad\qquad\qquad\qquad\qquad
\qquad\qquad\qquad\ \qquad\qquad\qquad\qquad\qquad\qquad\qquad\qquad
\qquad\qquad\qquad$

$-\frac{1}{6}[-R_{\text{a}i\text{b}i}+5\overset{q}{\underset{\text{c}=1}{\sum
}}T_{\text{ac}i}T_{\text{bc}i}+2\overset{n}{\underset{j=q+1}{\sum}}%
\perp_{\text{a}ij}\perp_{\text{b}ij}](y_{0})\underset{k=q+1}{\overset{n}{\sum
}}T_{\text{ab}k}(y_{0})X_{k}(y_{0})\phi(y_{0})\qquad\qquad$I$_{32123}\qquad
S_{1}\qquad$

$-\frac{2}{9}\underset{j=q+1}{\overset{n}{\sum}}R_{i\text{a}ij}(y_{0}%
)[\frac{\partial X_{j}}{\partial x_{\text{a}}}%
-\underset{k=q+1}{\overset{n}{\sum}}\perp_{\text{a}jk}X_{k}](y_{0})\phi
(y_{0})$

$+\frac{1}{12}\times\frac{2}{3}\underset{j,k=q+1}{\overset{n}{\sum}}%
R_{ijik}(y_{0})[X_{j}X_{k}-\frac{1}{2}(\frac{\partial X_{j}}{\partial x_{k}%
}+\frac{\partial X_{k}}{\partial x_{j}})](y_{0})\phi(y_{0})$

$-\frac{1}{6}T_{\text{ab}i}(y_{0})\frac{\partial^{2}X_{i}}{\partial
x_{\text{a}}\partial x_{\text{b}}}(y_{0})\phi(y_{0})\qquad\qquad\qquad
S_{2}\qquad\qquad S_{21}\qquad\qquad\qquad\qquad\qquad\qquad$

$+$ $\frac{1}{12}T_{\text{ab}i}(y_{0})$

$\times\lbrack$ $(R_{\text{a}i\text{b}j}+R_{\text{a}j\text{b}i})$
$-\underset{\text{c=1}}{\overset{\text{q}}{\sum}}(T_{\text{ac}i}T_{\text{bc}%
j}+T_{\text{ac}j}T_{\text{bc}i})-\overset{n}{\underset{k=q+1}{\sum}}%
(\perp_{\text{a}ik}\perp_{\text{b}jk}+$ $\perp_{\text{a}jk}\perp_{\text{b}%
ik})](y_{0})X_{j}(y_{0})\phi(y_{0})$

$-$ $\frac{1}{6}T_{\text{ab}i}(y_{0})T_{\text{ab}j}(y_{0})[X_{i}X_{j}-\frac
{1}{2}\left(  \frac{\partial X_{i}}{\partial x_{j}}+\frac{\partial X_{j}%
}{\partial x_{i}}\right)  ](y_{0})\phi(y_{0})$

$\qquad-\frac{1}{3}\perp_{\text{a}ij}(y_{0})[(X_{i}\frac{\partial X_{j}%
}{\partial x_{\text{a}}}+X_{j}\frac{\partial X_{i}}{\partial x_{\text{a}}%
})-\frac{1}{4}\left(  \frac{\partial^{2}X_{i}}{\partial x_{\text{a}}\partial
x_{j}}+\frac{\partial^{2}X_{j}}{\partial x_{\text{a}}\partial x_{i}}\right)
](y_{0})\phi(y_{0})\qquad S_{22}$

$\qquad-\frac{1}{6}\perp_{\text{a}ij}(y_{0})[T_{\text{ab}j}\frac{\partial
X_{i}}{\partial x_{\text{b}}}](y_{0})\phi(y_{0})$

$\qquad+\frac{1}{6}\perp_{\text{a}ij}(y_{0})[(\perp_{\text{b}ik}T_{\text{ab}%
j})+\frac{2}{3}(2R_{\text{a}ijk}+R_{\text{a}jik}+R_{\text{a}kji})](y_{0}%
)X_{k}(y_{0})\phi(y_{0})$

$\qquad-\frac{1}{6}\perp_{\text{a}ij}(y_{0})\perp_{\text{a}jk}(y_{0}%
)[X_{i}X_{k}-\frac{1}{2}\left(  \frac{\partial X_{i}}{\partial x_{k}}%
+\frac{\partial X_{k}}{\partial x_{i}}\right)  ](y_{0})\phi(y_{0})\qquad
\qquad\qquad\qquad\qquad\qquad\qquad\qquad$

$\qquad+\frac{1}{12}[(\frac{\partial X_{j}}{\partial x_{\text{a}}})^{2}%
+X_{j}\frac{\partial^{2}X_{j}}{\partial x_{\text{a}}^{2}}-\frac{1}{2}%
\frac{\partial^{3}X_{j}}{\partial x_{\text{a}}^{2}\partial x_{j}}](y_{0}%
)\phi(y_{0})-\frac{1}{6}\overset{n}{\underset{k=q+1}{\sum}}[\perp_{\text{b}%
ik}$T$_{\text{aa}k}\frac{\partial X_{i}}{\partial x_{\text{b}}^{2}}%
](y_{0})\phi(y_{0})\qquad\qquad S_{3}\qquad S_{31}$

$\qquad+\frac{1}{144}[\{4\nabla_{i}R_{i\text{a}j\text{a}}+2\nabla
_{j}R_{i\text{a}i\text{a}}+$ $8(\overset{q}{\underset{\text{c=1}}{%
{\textstyle\sum}
}}R_{\text{a}i\text{c}i}^{{}}T_{\text{ac}j}+\;\overset{n}{\underset{k=q+1}{%
{\textstyle\sum}
}}R_{\text{a}iik}\perp_{\text{a}jk})$

$\qquad+8(\overset{q}{\underset{\text{c=1}}{%
{\textstyle\sum}
}}R_{\text{a}i\text{c}j}^{{}}T_{\text{ac}i}+\;\overset{n}{\underset{k=q+1}{%
{\textstyle\sum}
}}R_{\text{a}ijk}\perp_{\text{a}ik})+8(\overset{q}{\underset{\text{c=1}}{%
{\textstyle\sum}
}}R_{\text{a}j\text{c}i}^{{}}T_{\text{ac}i}+\;\overset{n}{\underset{k=q+1}{%
{\textstyle\sum}
}}R_{\text{a}jik}\perp_{\text{a}ik})\}$\ 

$\qquad+\frac{2}{3}\underset{k=q+1}{\overset{n}{\sum}}\{T_{\text{aa}%
k}(R_{ijik}+3\overset{q}{\underset{\text{c}=1}{\sum}}\perp_{\text{c}ij}%
\perp_{\text{c}ik})\}](y_{0})X_{k}(y_{0})\phi(y_{0})$

$-\frac{1}{12}[$ R$_{\text{a}i\text{a}k}$ $-\underset{\text{c=1}%
}{\overset{\text{q}}{\sum}}T_{\text{ac}i}T_{\text{ac}k}%
-\overset{n}{\underset{l=q+1}{\sum}}(\perp_{\text{a}il}\perp_{\text{a}%
kl}](y_{0})\times\lbrack X_{i}X_{k}-\frac{1}{2}\left(  \frac{\partial X_{i}%
}{\partial x_{k}}+\frac{\partial X_{k}}{\partial x_{i}}\right)  ](y_{0}%
)\phi(y_{0})$

$-\frac{1}{24}T_{\text{aa}k}(y_{0})[-X_{i}^{2}X_{k}+X_{k}\frac{\partial X_{i}%
}{\partial x_{i}}\ +X_{i}\left(  \frac{\partial X_{k}}{\partial x_{i}}%
+\frac{\partial X_{i}}{\partial x_{k}}\right)  -\frac{1}{3}\left(
\frac{\partial^{2}X_{k}}{\partial x_{i}^{2}}+2\frac{\partial^{2}X_{i}%
}{\partial x_{i}\partial x_{k}}\right)  ](y_{0})\phi(y_{0})$

$+\frac{1}{18}[R_{\text{a}jij}\frac{\partial X_{i}}{\partial x_{\text{a}}^{2}%
}](y_{0})\qquad\qquad\qquad\qquad\qquad\qquad\qquad\qquad\qquad S_{32}$

$+\frac{1}{24}[\frac{4}{3}\overset{q}{\underset{\text{a}=1}{\sum}}%
\perp_{\text{a}ki}R_{ij\text{a}j}-\frac{1}{3}(\nabla_{i}R_{kjij}+\nabla
_{j}R_{ijik}+\nabla_{k}R_{ijij})](y_{0})X_{k}(y_{0})\phi(y_{0})$

\ $-\frac{1}{18}R_{ijkj}(y_{0})[X_{i}X_{k}-\frac{1}{2}\left(  \frac{\partial
X_{i}}{\partial x_{k}}+\frac{\partial X_{k}}{\partial x_{i}}\right)
](y_{0})\phi(y_{0})$

$+\frac{1}{24}[X_{i}^{2}X_{j}^{2}-2X_{i}X_{j}\left(  \frac{\partial X_{j}%
}{\partial x_{i}}+\frac{\partial X_{i}}{\partial x_{j}}\right)  -X_{i}%
^{2}\frac{\partial X_{j}}{\partial x_{j}}-X_{j}^{2}\frac{\partial X_{i}%
}{\partial x_{i}}](y_{0})\phi(y_{0})$

$+\frac{1}{48}\left(  \frac{\partial X_{j}}{\partial x_{i}}+\frac{\partial
X_{i}}{\partial x_{j}}\right)  ^{2}(y_{0})\phi(y_{0})+\frac{1}{24}\left(
\frac{\partial X_{i}}{\partial x_{i}}\frac{\partial X_{j}}{\partial x_{j}%
}\right)  (y_{0})\phi(y_{0})\qquad$

$+\frac{1}{36}X_{i}(y_{0})\left(  2\frac{\partial^{2}X_{j}}{\partial
x_{i}\partial x_{j}}+\frac{\partial^{2}X_{i}}{\partial x_{j}^{2}}\right)
(y_{0})\phi(y_{0})+\frac{1}{36}X_{j}(y_{0})\left(  \frac{\partial^{2}X_{j}%
}{\partial x_{i}^{2}}+2\frac{\partial^{2}X_{i}}{\partial x_{i}\partial x_{j}%
}\right)  (y_{0})\phi(y_{0})$

$-\frac{1}{48}\left(  \frac{\partial^{3}X_{i}}{\partial x_{i}\partial
x_{j}^{2}}+\frac{\partial^{3}X_{j}}{\partial x_{i}^{2}\partial x_{j}}\right)
(y_{0})\phi(y_{0})$

$+$ $\frac{2}{3}<H,j>(y_{0})\left(  \frac{\partial^{2}X_{i}}{\partial
x_{i}\partial x_{j}}+2\frac{\partial^{2}X_{j}}{\partial x_{i}^{2}}\right)
(y_{0})\phi(y_{0})+$ $\frac{2}{3}<H,j>(y_{0})R_{ijik}(y_{0})X_{k}(y_{0}%
)\phi(y_{0})\qquad$I$_{3213}$

$+\frac{1}{12}[<H,i><H,j>\ +\frac{1}{6}(2\varrho_{ij}%
+4\overset{q}{\underset{\text{a}=1}{\sum}}R_{i\text{a}j\text{a}}%
-6\overset{q}{\underset{\text{a,b}=1}{\sum}}T_{\text{aa}i}T_{\text{bb}%
j}-T_{\text{ab}i}T_{\text{ab}j})](y_{0})\phi(y_{0})$

$\times\frac{1}{2}[\left(  \frac{\partial X_{j}}{\partial x_{i}}%
-\frac{\partial X_{i}}{\partial x_{j}}\right)  ](y_{0})\phi(y_{0})$

$-\frac{1}{12}\perp_{\text{a}ij}(y_{0})<H,i>(y_{0})[(X_{j}\perp_{\text{a}%
ij}-\frac{\partial X_{i}}{\partial x_{\text{a}}})+\frac{\partial X_{\text{a}}%
}{\partial x_{i}}](y_{0})\phi(y_{0})$

$-\frac{1}{18}[X_{j}\left(  2\frac{\partial^{2}X_{j}}{\partial x_{i}^{2}%
}+\frac{\partial^{2}X_{i}}{\partial x_{i}\partial x_{j}}\right)  ](y_{0}%
)\phi(y_{0})-\frac{1}{12}[\left(  \frac{\partial X_{i}}{\partial x_{j}}%
+\frac{\partial X_{j}}{\partial x_{i}}\right)  ]\frac{\partial X_{j}}{\partial
x_{i}}(y_{0})\phi(y_{0})\qquad\ $I$_{3214}\qquad\qquad\qquad$\qquad
\qquad\qquad\qquad\qquad\qquad\qquad\qquad\qquad\qquad\qquad\qquad\qquad
\qquad\qquad\qquad

$+\frac{1}{12}\frac{\partial^{2}\text{V}}{\partial x_{i}^{2}}(y_{0})\phi
(y_{0})\qquad\qquad$I$_{3215}$

$+\frac{1}{12}\underset{i=q+\text{1}}{\overset{\text{n}}{\sum}}%
\underset{\text{a,b=1}}{\overset{\text{q}}{\sum}}[-R_{\text{a}i\text{b}%
i}+5\overset{q}{\underset{\text{c}=1}{\sum}}T_{\text{ac}i}T_{\text{bc}%
i}+2\overset{n}{\underset{j=q+1}{\sum}}\perp_{\text{a}ij}\perp_{\text{b}%
ij}](y_{0})\times\frac{\partial^{2}\phi}{\partial\text{x}_{\text{a}}%
\partial\text{x}_{\text{b}}}(y_{0})\qquad$\textbf{I}$_{322}$

$+\frac{1}{72}\underset{i,j,k=q+1}{\overset{n}{\sum}}R_{ijik}(y_{0}%
)\Omega_{jk}(y_{0})\phi(y_{0})\qquad\qquad\qquad\qquad$\ I$_{323}$

$+\frac{1}{24}\underset{\text{a=1}}{\overset{\text{q}}{\sum}}%
\underset{i,j=q+1}{\overset{n}{\sum}}\left\{  \frac{8}{3}R_{i\text{a}%
ij}+4\underset{\text{b=1}}{\overset{\text{q}}{\sum}}T_{\text{ab}i}%
\perp_{\text{b}ji}\right\}  (y_{0})\left\{  -\Omega_{\text{a}j}+[\Lambda
_{\text{a}},\Lambda_{j}]\right\}  (y_{0})\phi(y_{0})$

$+\frac{1}{12}\underset{i\text{=}q+1}{\overset{\text{n}}{\sum}}%
\underset{\text{a,b=1}}{\overset{\text{q}}{\sum}}$ $[-R_{\text{a}i\text{b}%
i}+5\overset{q}{\underset{\text{c}=1}{\sum}}T_{\text{ac}i}T_{\text{bc}%
i}+2\overset{n}{\underset{\text{k}=q+1}{\sum}}\perp_{\text{a}i\text{k}}%
\perp_{\text{b}i\text{k}}](y_{0})\times\lbrack\Lambda_{\text{a}}(y_{0}%
)\Lambda_{\text{b}}(y_{0})\phi(y_{0})]\qquad$\ I$_{324}$

$+\frac{1}{12}[\frac{8}{3}R_{i\text{a}ij}-4\underset{\text{b=1}%
}{\overset{\text{q}}{\sum}}T_{\text{ab}i}(y_{0})\perp_{\text{b}ij}%
](y_{0})[\Lambda_{\text{a}}\Lambda_{j}\phi](y_{0})$

$+\frac{1}{12}\underset{i\text{=}q+1}{\overset{\text{n}}{\sum}}%
\underset{\text{a,b=1}}{\overset{\text{q}}{\sum}}$ $[-R_{\text{a}i\text{b}%
i}+5\overset{q}{\underset{\text{c}=1}{\sum}}T_{\text{ac}i}T_{\text{bc}%
i}+2\overset{n}{\underset{\text{k}=q+1}{\sum}}\perp_{\text{a}i\text{k}}%
\perp_{\text{b}i\text{k}}](y_{0})\qquad$I$_{325}\qquad$I$_{3251}$

$\times\lbrack\Lambda_{\text{a}}(y_{0})\Lambda_{\text{b}}(y_{0})\phi(y_{0})]$

$+\frac{1}{12}[\frac{8}{3}R_{i\text{a}ij}-4\underset{\text{b=1}%
}{\overset{\text{q}}{\sum}}T_{\text{ab}i}(y_{0})\perp_{\text{b}ij}%
](y_{0})[\Lambda_{\text{a}}\Lambda_{j}\phi](y_{0})$

$+\frac{1}{24}\underset{i=q+1}{\overset{n}{\sum}}\underset{\text{a=1}%
}{\overset{\text{q}}{\sum}}\left(  \frac{\partial\Omega_{i\text{a}}}{\partial
x_{i}}\Lambda_{\text{a}}+[\Omega_{i\text{a}}+[\Lambda_{\text{a}},\Lambda
_{i}],\Lambda_{i}]\right)  \Lambda_{\text{a}}(y_{0})\phi(y_{0})\qquad
$I$_{3252}$

$+\frac{1}{24}\underset{i=q+1}{\overset{n}{\sum}}\underset{\text{a=1}%
}{\overset{\text{q}}{\sum}}\Lambda_{\text{a}}(y_{0})\left(  \frac
{\partial\Omega_{i\text{a}}}{\partial x_{i}}\Lambda_{\text{a}}+[\Omega
_{i\text{a}}+[\Lambda_{\text{a}},\Lambda_{i}],\Lambda_{i}]\right)  (y_{0}%
)\phi(y_{0})$

$+\frac{1}{12}\underset{i=q+1}{\overset{n}{\sum}}\underset{\text{a=1}%
}{\overset{\text{q}}{\sum}}\left(  \Omega_{i\text{a}}+[\Lambda_{\text{a}%
},\Lambda_{i}]\right)  ^{2}(y_{0})\phi(y_{0})$

$+\frac{1}{48}\underset{i,j=q+1}{\overset{n}{\sum}}\left(  \Omega_{ij}%
\Omega_{ij}\right)  (y_{0})\phi(y_{0})+\frac{1}{72}%
\underset{i,j=q+1}{\overset{n}{\sum}}\left(  \frac{\partial\Omega_{ij}%
}{\partial\text{x}_{i}}\Lambda_{j}+\Lambda_{j}\frac{\partial\Omega_{ij}%
}{\partial\text{x}_{i}}\right)  (y_{0})\phi(y_{0})$

$+\frac{1}{12}[\underset{i=q+1}{\overset{n}{\sum}}\underset{\text{a,b=1}%
}{\overset{\text{q}}{\sum}}2T_{\text{ab}i}(y_{0})\left\{  (\Omega_{i\text{a}%
}+[\Lambda_{\text{a}},\Lambda_{i}])\Lambda_{\text{b}}+\Lambda_{\text{a}%
}(\Omega_{i\text{a}}+[\Lambda_{\text{a}},\Lambda_{i}])\right\}  (y_{0}%
)\phi(y_{0})\qquad$I$_{3253}$

$-\frac{1}{12}[\underset{i,j=q+1}{\overset{n}{\sum}}\underset{\text{a=1}%
}{\overset{\text{q}}{\sum}}\perp_{\text{a}ij}(y_{0})\left\{  (\Omega
_{i\text{a}}+[\Lambda_{\text{a}},\Lambda_{i}])\Lambda_{j}+\frac{1}{2}%
\Lambda_{\text{a}}\Omega_{ij}\right\}  ](y_{0})\phi(y_{0})$

$-\frac{1}{12}[\underset{i,j=q+1}{\overset{n}{\sum}}\underset{\text{b=1}%
}{\overset{\text{q}}{\sum}}\perp_{\text{b}ij}(y_{0})\left\{  \frac{1}{2}%
\Omega_{ij}\Lambda_{\text{b}}+\Lambda_{j}(\Omega_{i\text{b}}+[\Lambda
_{\text{b}},\Lambda_{i}])\right\}  ](y_{0})\phi(y_{0})$

$\mathbf{-}\frac{1}{12}\underset{\text{a,b=1}}{\overset{\text{q}}{\sum}%
}\underset{i,j=q+1}{\overset{n}{\sum}}T_{\text{ab}i}(y_{0})[-R_{\text{a}%
i\text{b}i}+5\overset{q}{\underset{\text{c}=1}{\sum}}T_{\text{ac}%
i}T_{\text{bc}i}+2\overset{n}{\underset{k=q+1}{\sum}}\perp_{\text{a}ik}%
\perp_{\text{b}ik}](y_{0})\Lambda_{j}(y_{0})\phi(y_{0})\ \ \qquad$%
\textbf{I}$_{326}\qquad$\textbf{I}$_{3261}\qquad$

$+\frac{1}{12}\underset{i\text{=q+1}}{\overset{\text{n}}{\sum}}%
\underset{j\text{=q+1}}{\overset{\text{n}}{\sum}}\underset{\text{a=1}%
}{\overset{\text{q}}{\sum}}[4\underset{\text{c=1}}{\overset{q}{\sum}%
}(T_{\text{ac}i})(\perp_{j\text{c}i})+\frac{8}{3}R_{i\text{a}ij}](y_{0}%
)\qquad\qquad$I$_{32613}$

$\times\lbrack\underset{\text{c=1}}{\overset{\text{q}}{\sum}}T_{\text{ac}%
j}\frac{\partial\phi}{\partial\text{x}_{\text{c}}}+\underset{\text{b=1}%
}{\overset{\text{q}}{\sum}}T_{\text{ab}j}\Lambda_{\text{b}}%
-\underset{k=q+1}{\overset{n}{\sum}}\perp_{\text{a}jk}\Lambda_{k}](y_{0}%
)\phi(y_{0})$

$-\frac{1}{24}\underset{i=q+1}{\overset{n}{\sum}}\underset{\text{a,b=1}%
}{\overset{\text{q}}{\sum}}$ $\overset{n}{\underset{k=q+1}{\sum}}%
$T$_{\text{aa}k}[\frac{8}{3}R_{i\text{c}ik}+4\underset{\text{d}%
=1}{\overset{q}{\sum}}(T_{\text{db}k})(\perp_{\text{d}ik})]\frac{\partial\phi
}{\partial\text{x}_{\text{b}}}$ $(y_{0})\qquad$ I$_{32621}\qquad$I$_{3262}$

$-\frac{1}{12}\underset{i=q+1}{\overset{n}{\sum}}\underset{\text{a,b=1}%
}{\overset{\text{q}}{\sum}}[$ $\overset{n}{\underset{k,l=q+1}{\sum}}%
\perp_{\text{b}ik}(-R_{\text{a}k\text{a}l}+\underset{\text{d=1}%
}{\overset{\text{q}}{\sum}}T_{\text{ad}k}T_{\text{ad}l}%
))-\overset{n}{\underset{k,l=q+1}{\sum}}\perp_{\text{b}ik}%
(\underset{r=q+1}{\overset{n}{\sum}}\perp_{\text{a}kr}\perp_{\text{a}%
lr})](y_{0})\frac{\partial\phi}{\partial\text{x}_{\text{b}}}$ $(y_{0})$

$-\frac{1}{12}\underset{i=q+1}{\overset{n}{\sum}}\underset{\text{b=1}%
}{\overset{\text{q}}{\sum}}[\frac{8}{3}\overset{q}{\underset{\text{c}=1}{\sum
}}($T$_{\text{bc}i}R_{ij\text{c}j})+\frac{2}{3}%
\overset{n}{\underset{k=q+1}{\sum}}(\perp_{\text{b}ik}R_{ijjk})](y_{0}%
)\frac{\partial\phi}{\partial\text{x}_{\text{b}}}(y_{0})$

$-\frac{1}{6}\underset{i=q+1}{\overset{n}{\sum}}\underset{\text{a,b=1}%
}{\overset{\text{q}}{\sum}}[4\overset{q}{\underset{\text{c=1}}{%
{\textstyle\sum}
}}$R$_{ij\text{c}i}^{{}}T_{\text{bc}j}+$ $4\overset{n}{\underset{k=q+1}{%
{\textstyle\sum}
}}$R$_{ijik}\perp_{\text{b}jk}+3\nabla_{i}$R$_{j\text{b}ij}%
+4\overset{q}{\underset{\text{c=1}}{%
{\textstyle\sum}
}}$R$_{ij\text{c}j}^{{}}T_{\text{bc}i}+$ $4$R$_{ijjk}\perp_{\text{b}ik}%
](y_{0})\frac{\partial\phi}{\partial\text{x}_{\text{b}}}(y_{0})$

$-\frac{1}{24}\overset{n}{\underset{k=q+1}{\sum}}$T$_{\text{aa}k}[\frac{8}%
{3}R_{i\text{c}ik}+4\underset{\text{d}=1}{\overset{q}{\sum}}(T_{\text{db}%
k})(\perp_{\text{d}ik})]\Lambda_{\text{b}}(y_{0})\phi(y_{0})\qquad$%
I$_{326221}\qquad$I$_{32622}$

$-\frac{1}{12}$ $\overset{n}{\underset{k,l=q+1}{\sum}}\perp_{\text{b}ik}%
(y_{0})[-R_{\text{a}k\text{a}l}+\underset{\text{d=1}}{\overset{\text{q}}{\sum
}}T_{\text{ad}k}T_{\text{ad}l}](y_{0})\Lambda_{\text{b}}(y_{0})\phi(y_{0})$

$-\frac{1}{12}\overset{n}{\underset{k,l=q+1}{\sum}}\perp_{\text{b}ik}%
(y_{0})[\underset{r=q+1}{\overset{n}{\sum}}\perp_{\text{a}kr}\perp
_{\text{a}lr}](y_{0})\Lambda_{\text{b}}(y_{0})\phi(y_{0})$

$-\frac{1}{144}[\{4\nabla_{i}$R$_{i\text{a}j\text{a}}$ $+2\nabla_{j}%
$R$_{i\text{a}i\text{a}}+$ $8$ $(\overset{q}{\underset{\text{c=1}}{%
{\textstyle\sum}
}}R_{\text{a}i\text{c}i}^{{}}T_{\text{ac}j}+\;\overset{n}{\underset{k=q+1}{%
{\textstyle\sum}
}}R_{\text{a}iik}\perp_{\text{a}jk})$

\ $+8(\overset{q}{\underset{\text{c=1}}{%
{\textstyle\sum}
}}R_{\text{a}i\text{c}j}^{{}}T_{\text{ac}i}+\overset{n}{\underset{l=q+1}{%
{\textstyle\sum}
}}R_{\text{a}ijl}\perp_{\text{a}il})+8(\overset{q}{\underset{\text{c=1}}{%
{\textstyle\sum}
}}R_{\text{a}j\text{c}i}^{{}}T_{\text{ac}i}+\overset{q}{\underset{\text{c=1}}{%
{\textstyle\sum}
}}R_{\text{a}j\text{c}i}^{{}}T_{\text{ac}i})\}$\ \qquad$\qquad\ \ \ \ \ $

$+\frac{2}{3}\underset{k=q+1}{\overset{n}{\sum}}\{T_{\text{aa}k}%
(R_{ijik}+3\overset{q}{\underset{\text{c}=1}{\sum}}\perp_{\text{c}ij}%
\perp_{\text{c}ik})\}](y_{0})\Lambda_{k}(y_{0})\phi(y_{0})$

\ $+\frac{1}{24}[\frac{4}{3}\overset{q}{\underset{\text{a}=1}{\sum}}%
\perp_{\text{a}ik}R_{ij\text{a}j}+\frac{1}{3}(\nabla_{i}R_{kjij}+\nabla
_{j}R_{ijik}+\nabla_{k}R_{ijij})](y_{0})\Lambda_{k}(y_{0})\phi(y_{0})$

\ $\mathbf{-}$ $\frac{1}{72}\underset{i,j=q+1}{\overset{n}{\sum}%
}\underset{\text{a}=1}{\overset{q}{\sum}}T_{\text{aa}j}(y_{0})\frac
{\partial\Omega_{ij}}{\partial\text{x}_{i}}(y_{0})\phi(y_{0})\qquad
$I$_{326222}$

$+$ $\frac{1}{12}$ $\overset{n}{\underset{j=q+1}{\sum}}(\perp_{\text{b}ij}%
$T$_{\text{aa}j})(y_{0})\left(  \Omega_{i\text{b}}(y_{0})+[\Lambda_{\text{b}%
},\Lambda_{i}]\right)  (y_{0})\phi(y_{0})$ \qquad I$_{326223}$

$\mathbf{-}\frac{1}{18}$ $\underset{i,j=q+1}{\overset{n}{\sum}}%
\underset{\text{b=1}}{\overset{q}{\sum}}R_{\text{b}jij}(y_{0})\left(
\Omega_{i\text{a}}(y_{0})+[\Lambda_{\text{a}},\Lambda_{i}]\right)  (y_{0}%
)\phi(y_{0})$

$\mathbf{-}$ $\frac{1}{24}\underset{i,j=q+1}{\overset{n}{\sum}}%
\underset{\text{a=1}}{\overset{q}{\sum}}[$ R$_{\text{a}i\text{a}j}$
$-\underset{\text{c=1}}{\overset{\text{q}}{\sum}}T_{\text{ac}i}T_{\text{ac}%
j}-\overset{n}{\underset{k=q+1}{\sum}}(\perp_{\text{a}ik}\perp_{\text{a}%
jk}](y_{0})\Omega_{ij}(y_{0})\phi(y_{0})$

\ $\mathbf{-}$ $\frac{1}{36}\underset{i,j,k=q+1}{\overset{n}{\sum}}%
R_{ijkj}(y_{0})\Omega_{ik}(y_{0})(y_{0})\phi(y_{0})$

$+\frac{1}{6}\underset{i,k=q+1}{\overset{n}{\sum}}\underset{\text{a,b,c}%
=1}{\overset{q}{\sum}}T_{\text{ab}i}(y_{0})[(\perp_{\text{c}ik}T_{\text{ab}%
k})(\frac{\partial\phi}{\partial\text{x}_{\text{c}}}$ $+$ $\Lambda_{\text{c}%
}\phi)](y_{0})\qquad$\textbf{I}$_{32631}\qquad$\textbf{I}$_{3263}$

$-\frac{1}{12}[$ $(R_{\text{a}i\text{b}j}+R_{\text{a}j\text{b}i})$
$-\underset{\text{c=1}}{\overset{\text{q}}{\sum}}(T_{\text{ac}i}T_{\text{bc}%
j}+T_{\text{ac}j}T_{\text{bc}i})$

$-\overset{n}{\underset{k=q+1}{\sum}}(\perp_{\text{a}ik}\perp_{\text{b}jk}+$
$\perp_{\text{a}jk}\perp_{\text{b}ik})](y_{0})T_{\text{ab}i}(y_{0})\Lambda
_{j}(y_{0})\phi(y_{0})-\frac{1}{12}T_{\text{ab}i}^{2}(y_{0})\Omega_{ij}%
(y_{0})\phi(y_{0}$

$+\frac{1}{12}\underset{i,j=q+1}{\overset{n}{\sum}}\underset{\text{a,b=1}%
}{\overset{q}{\sum}}[\perp_{\text{a}ij}(\frac{\partial\phi}{\partial
\text{x}_{\text{b}}}+\Lambda_{\text{b}})](y_{0})\qquad\qquad$\ \textbf{I}%
$_{32632}$

$\times\lbrack-R_{\text{a}i\text{b}j}-R_{\text{a}j\text{b}i}%
+\underset{\text{c=1}}{\overset{\text{q}}{\sum}}T_{\text{ac}i}T_{\text{bc}%
j}-3\underset{\text{c=1}}{\overset{\text{q}}{\sum}}T_{\text{ac}j}%
T_{\text{bc}i}+\underset{\text{k=q+1}}{\overset{\text{n}}{\sum}}%
\perp_{\text{a}i\text{k}}\perp_{\text{b}j\text{k}}-\underset{\text{k=q+1}%
}{\overset{\text{n}}{\sum}}\perp_{\text{a}j\text{k}}\perp_{\text{b}i\text{k}%
}\ ](y_{0})\phi(y_{0})$

$-\frac{1}{6}\underset{i,j=q+1}{\overset{n}{\sum}}\underset{\text{a,b=1}%
}{\overset{q}{\sum}}T_{\text{ab}j}(y_{0})\perp_{\text{a}ij}(y_{0}%
)\frac{\partial\Lambda_{\text{b}}}{\partial x_{i}}(y_{0})\phi(y_{0})$

$-\frac{1}{6}\underset{i,j,k=q+1}{\overset{n}{\sum}}\underset{\text{a}%
=1}{\overset{q}{\sum}}\perp_{\text{a}ij}(y_{0})[\overset{q}{\underset{\text{b}%
=1}{\sum(}}\perp_{\text{b}ik}T_{\text{ab}j})(y_{0})+\frac{2}{3}(2R_{\text{a}%
ijk}+R_{\text{a}jik}+R_{\text{a}kji})](y_{0})\Lambda_{k}(y_{0})\phi(y_{0})$

$+\frac{1}{6}\underset{i,j,k=q+1}{\overset{n}{\sum}}\underset{\text{a}%
=1}{\overset{q}{\sum}}\perp_{\text{a}ij}(y_{0})\perp_{\text{a}jk}(y_{0}%
)\Omega_{ik}(y_{0})\phi(y_{0})$

$+\ \frac{1}{24}\underset{i\text{=q+1}}{\overset{\text{n}}{\sum}}%
\frac{\partial^{2}\text{W}}{\partial x_{i}^{2}}(y_{0})\phi(y_{0})\qquad\qquad
$\ \textbf{I}$_{327}$

$+\frac{1}{24}\underset{i,j\text{=q+1}}{\overset{\text{n}}{\sum}%
}\underset{\text{a=1}}{\overset{\text{q}}{\sum}}<H,j>[4\underset{\text{c}%
=1}{\overset{q}{\sum}}(T_{\text{ac}i})(\perp_{j\text{c}i})+\frac{8}%
{3}R_{i\text{a}ij}](y_{0})\frac{\partial\phi}{\partial\text{x}_{\text{a}}%
}(y_{0})\qquad$I$_{328}$

$+\frac{1}{24}\underset{i,j\text{=q+1}}{\overset{\text{n}}{\sum}%
}\underset{\text{a=1}}{\overset{\text{q}}{\sum}}\perp_{\text{a}ji}%
(y_{0})[<H,i><H,j>](y_{0})\frac{\partial\phi}{\partial\text{x}_{\text{a}}%
}(y_{0})$\qquad$\ \ \ +\frac{1}{72}\underset{i,j=q+1}{\overset{n}{\sum}%
}\underset{\text{a=1}}{\overset{\text{q}}{\sum}}\perp_{\text{a}ji}%
(y_{0})[2\varrho_{ij}+4\overset{q}{\underset{\text{a}=1}{\sum}}R_{i\text{a}%
j\text{a}}-6\overset{q}{\underset{\text{b,c}=1}{\sum}}T_{\text{cc}%
i}T_{\text{bb}j}-T_{\text{bc}i}T_{\text{bc}j}](y_{0})\frac{\partial\phi
}{\partial\text{x}_{\text{a}}}(y_{0})$

$+\frac{8}{3}R_{j\text{a}ji}(y_{0})X_{i}(y_{0})+[2X_{j}\frac{\partial X_{j}%
}{\partial x_{\text{a}}}-\frac{\partial^{2}X_{j}}{\partial x_{\text{a}%
}\partial x_{j}}](y_{0})\frac{\partial\phi}{\partial\text{x}_{\text{a}}}%
(y_{0})$

$+\frac{1}{24}\underset{i,j=q+1}{\overset{n}{\sum}}\underset{\text{a=1}%
}{\overset{\text{q}}{\sum}}<H,j>(y_{0})[4\underset{\text{c}%
=1}{\overset{q}{\sum}}(T_{\text{ac}i})(\perp_{j\text{c}i})+\frac{8}%
{3}R_{i\text{a}ij}](y_{0})\Lambda_{\text{a}}(y_{0})\phi(y_{0})\qquad$%
I$_{329}\qquad$I$_{3291}\qquad$I$_{32911}\qquad$

$+\frac{1}{12}\underset{i,j=q+1}{\overset{n}{\sum}}\underset{\text{a=1}%
}{\overset{\text{q}}{\sum}}\perp_{\text{a}ji}(y_{0})[<H,i><H,j>](y_{0}%
)\Lambda_{\text{a}}(y_{0})\phi(y_{0})$

$+\frac{1}{72}\underset{i,j=q+1}{\overset{n}{\sum}}\underset{\text{a=1}%
}{\overset{\text{q}}{\sum}}\perp_{\text{a}ji}(y_{0})[2\varrho_{ij}%
+4\overset{q}{\underset{\text{a}=1}{\sum}}R_{i\text{a}j\text{a}}%
-6\overset{q}{\underset{\text{b,c}=1}{\sum}}T_{\text{cc}i}T_{\text{bb}%
j}-T_{\text{bc}i}T_{\text{bc}j}](y_{0})\Lambda_{\text{a}}(y_{0})\phi(y_{0})$

$+\frac{1}{36}\underset{i,j=q+1}{\overset{n}{\sum}}%
\underset{k=q+1}{\overset{n}{\sum}}<H,k>(y_{0})$R$_{ijik}(y_{0})\Lambda
_{j}(y_{0})\phi(y_{0})$

$-\frac{1}{288}\underset{i,j=q+1}{\overset{n}{\sum}}<H,j>(y_{0})[3<H,i>^{2}%
+2(\tau^{M}-3\tau^{P}+\overset{q}{\underset{\text{a}=1}{\sum}}\varrho
_{\text{aa}}+\overset{q}{\underset{\text{a,b}=1}{\sum}}R_{\text{abab}}%
)](y_{0})\Lambda_{j}(y_{0})\phi(y_{0})$

$-\frac{1}{12}\underset{i,j=q+1}{\overset{n}{\sum}}<H,i>(y_{0})$

$\times\lbrack\frac{3}{4}<H,i><H,j>$\ $+\frac{1}{6}(\varrho_{ij}%
+2\overset{q}{\underset{\text{a}=1}{\sum}}R_{i\text{a}j\text{a}}%
-3\overset{q}{\underset{\text{a,b=1}}{\sum}}T_{\text{aa}i}T_{\text{bb}%
j}-T_{\text{ab}i}T_{\text{ab}j})](y_{0})\Lambda_{j}(y_{0})\phi(y_{0})$

$+\frac{5}{32}\underset{i,j=q+1}{\overset{n}{\sum}}<H,i>^{2}<H,j>\Lambda
_{j}(y_{0})\phi(y_{0})$

$+\frac{1}{48}\underset{i,j=q+1}{\overset{n}{\sum}}$%
$<$%
H,$i$%
$>$%
(y$_{0}$)$[$(2$\varrho_{ij}$+4$\overset{q}{\underset{\text{a}=1}{\sum}}%
$R$_{i\text{a}j\text{a}}$-3$\overset{q}{\underset{\text{a,b=1}}{\sum}}%
$T$_{\text{aa}i}$T$_{\text{bb}j}$-T$_{\text{ab}i}$T$_{\text{ab}j}%
$-3$\overset{q}{\underset{\text{a,b=1}}{\sum}}$T$_{\text{aa}j}$T$_{\text{bb}%
i}$-T$_{\text{ab}j}$T$_{\text{ab}i}$)$](y_{0})\Lambda_{j}(y_{0})\phi(y_{0})$

$+\frac{1}{48}\underset{i,j=q+1}{\overset{n}{\sum}}<H,j>(y_{0})[$ $\tau
^{M}-3\tau^{P}+\overset{q}{\underset{\text{a}=1}{\sum}}\varrho_{\text{aa}%
}+\overset{q}{\underset{\text{a,b}=1}{\sum}}R_{\text{abab}}](y_{0})\Lambda
_{j}(y_{0})\phi(y_{0})$

$+\frac{1}{144}\underset{i,j=q+1}{\overset{n}{\sum}}[\nabla_{i}\varrho
_{ij}-2\varrho_{ij}<H,i>+\overset{q}{\underset{\text{a}=1}{\sum}}(\nabla
_{i}R_{\text{a}i\text{a}j}-4R_{i\text{a}j\text{a}}<H,i>)$

$+4\overset{q}{\underset{\text{a,b=1}}{\sum}}R_{i\text{a}j\text{b}%
}T_{\text{ab}i}+2\overset{q}{\underset{\text{a,b,c=1}}{\sum}}(T_{\text{aa}%
i}T_{\text{bb}j}T_{\text{cc}i}-3T_{\text{aa}i}T_{\text{bc}j}T_{\text{bc}%
i}+2T_{\text{ab}i}T_{\text{bc}j}T_{\text{ca}i})](y_{0})\Lambda_{j}(y_{0}%
)\phi(y_{0})$\qquad\qquad\qquad\qquad\qquad\ \ 

$+\frac{1}{144}\underset{i,j=q+1}{\overset{n}{\sum}}[\nabla_{j}\varrho
_{ii}-2\varrho_{ji}<H,i>+\overset{q}{\underset{\text{a}=1}{\sum}}(\nabla
_{j}R_{\text{a}i\text{a}i}-4R_{j\text{a}i\text{a}}<H,i>)$

$+4\overset{q}{\underset{\text{a,b=1}}{\sum}}R_{j\text{a}i\text{b}%
}T_{\text{ab}i}$+2$\overset{q}{\underset{\text{a,b,c=1}}{\sum}}(T_{\text{aa}%
j}T_{\text{bb}i}T_{\text{cc}i}-3T_{\text{aa}j}T_{\text{bc}i}T_{\text{bc}%
i}+2T_{\text{ab}j}T_{\text{bc}i}T_{\text{ca}i})](y_{0})\Lambda_{j}(y_{0}%
)\phi(y_{0})$

$+\frac{1}{144}\underset{i,j=q+1}{\overset{n}{\sum}}[\nabla_{i}\varrho
_{ij}-2\varrho_{ii}<H,j>+\overset{q}{\underset{\text{a}=1}{\sum}}(\nabla
_{i}R_{\text{a}i\text{a}j}-4R_{i\text{a}i\text{a}}<H,j>)$

$+4\overset{q}{\underset{\text{a,b=1}}{\sum}}R_{i\text{a}i\text{b}%
}T_{\text{ab}j}$+2$\overset{q}{\underset{\text{a,b,c=1}}{\sum}}(T_{\text{aa}%
i}T_{\text{bb}i}T_{\text{cc}j}-3T_{\text{aa}i}T_{\text{bc}i}T_{\text{bc}%
j}+2T_{\text{ab}i}T_{\text{bc}i}T_{\text{ac}j})](y_{0})\Lambda_{j}(y_{0}%
)\phi(y_{0})$

$+\frac{1}{72}\underset{i,j=q+1}{\overset{n}{\sum}}$ $<H,j>(y_{0}%
)\frac{\partial\Omega_{ij}}{\partial\text{x}_{i}}(y_{0})\phi(y_{0}%
)\qquad\qquad$I$_{32912}$

$-\frac{1}{12}\underset{i=q+1}{\overset{n}{\sum}}\underset{\text{a=1}%
}{\overset{\text{q}}{\sum}}\perp_{\text{a}ij}(y_{0})<H,j>(y_{0})[\Omega
_{i\text{a}}+[\Lambda_{\text{a}},\Lambda_{i}]](y_{0})\phi(y_{0})\qquad\qquad
$I$_{32913}$

$+$ $\frac{1}{72}\underset{i,j=q+1}{\overset{n}{\sum}}[3<H,i><H,j>\ +(\varrho
_{ij}+2\overset{q}{\underset{\text{a}=1}{\sum}}R_{i\text{a}j\text{a}%
}-3\overset{q}{\underset{\text{a,b}=1}{\sum}}T_{\text{aa}i}T_{\text{bb}%
j}-T_{\text{ab}i}T_{\text{ab}j})](y_{0})\Omega_{ij}(y_{0})\phi(y_{0})$

$+\frac{1}{12}\underset{i=q+1}{\overset{n}{\sum}}$ $\underset{\text{a=1}%
}{\overset{\text{q}}{\sum}}$ $[-4\underset{\text{b=1}}{\overset{\text{q}}{%
{\textstyle\sum}
}}T_{\text{ab}i}\frac{\partial X_{i}}{\partial x_{\text{b}}}%
+\underset{j=q+1}{\overset{n}{%
{\textstyle\sum}
}}\perp_{\text{a}ij}\left(  \frac{\partial X_{i}}{\partial x_{j}}%
+\frac{\partial X_{j}}{\partial x_{i}}\right)  \qquad\qquad$I$_{3292}\qquad
$I$_{32921}$

$+\frac{8}{3}\underset{j=q+1}{\overset{n}{%
{\textstyle\sum}
}}R_{i\text{a}ij}X_{j}+\left(  2X_{i}\frac{\partial X_{i}}{\partial
x_{\text{a}}}-\frac{\partial^{2}X_{i}}{\partial x_{\text{a}}\partial x_{i}%
}\right)  ](y_{0})\Lambda_{\text{a}}(y_{0})\phi(y_{0}$

$+\frac{1}{6}$ $\underset{i=q+1}{\overset{n}{\sum}}\underset{\text{a=1}%
}{\overset{\text{q}}{\sum}}$ $[\underset{j=q+1}{\overset{n}{%
{\textstyle\sum}
}}X_{j}\perp_{\text{a}ij}+\frac{\partial X_{i}}{\partial x_{\text{a}}}]$
$(y_{0})[\Omega_{i\text{a}}+[\Lambda_{\text{a}},\Lambda_{i}](y_{0})\phi
(y_{0})$

$-\frac{1}{36}[\left(  2\frac{\partial^{2}X_{i}}{\partial x_{i}\partial x_{j}%
}+\frac{\partial^{2}X_{j}}{\partial x_{i}^{2}}\right)
+2\underset{k=q+1}{\overset{n}{\sum}}R_{ijik}X_{k}](y_{0})\Lambda_{j}%
(y_{0})\phi(y_{0})\qquad$I$_{32922}$

$-\frac{1}{36}X_{j}(y_{0})$ $\frac{\partial\Omega_{ij}}{\partial x_{i}}%
(y_{0})\phi(y_{0})-\frac{1}{12}\frac{\partial X_{j}}{\partial x_{i}}%
(y_{0})\Omega_{ij}(y_{0})$ $\phi(y_{0})$

$+\frac{1}{12}$ $\underset{\text{a=1}}{\overset{\text{q}}{\sum}}\frac
{\partial^{2}X_{\text{a}}}{\partial x_{i}^{2}}(y_{0})\frac{\partial\phi
}{\partial\text{x}_{\text{a}}}(y_{0})\qquad$\ \ \textbf{L}$_{1}$\qquad\qquad

$+\frac{1}{12}$ $\underset{\text{a}=1}{\overset{\text{q}}{\sum}}$
$\frac{\partial^{2}X_{\text{a}}}{\partial x_{i}^{2}}(y_{0})\Lambda_{\text{a}%
}(y_{0})\phi(y_{0})+\frac{1}{12}\underset{j=q+1}{\overset{n}{\sum}}%
\frac{\partial^{2}X_{j}}{\partial x_{i}^{2}}(y_{0})\Lambda_{j}(y_{0}%
)\phi(y_{0})\qquad$\textbf{L}$_{2}\qquad$\textbf{L}$_{21}$

$+\frac{1}{36}$ $\underset{j=q+1}{\overset{n}{\sum}}$ $X_{j}(y_{0}%
)\frac{\partial\Omega_{ij}}{\partial x_{i}}(y_{0})\phi(y_{0})\qquad$%
\textbf{L}$_{22}$

$+\frac{1}{12}$ $\underset{j=q+1}{\overset{\text{n}}{\sum}}$ $\frac{\partial
X_{j}}{\partial x_{i}}(y_{0})\Omega_{ij}(y_{0})\phi(y_{0})\qquad$%
\textbf{L}$_{23}$\qquad$\qquad$

$+\frac{1}{48}$ $\underset{\text{c=1}}{\overset{\text{q}}{\sum}}%
\Lambda_{\text{c}}(y_{0})[\underset{\alpha=q+1}{\overset{n}{\sum}}%
3<H,i>^{2}+2(\tau^{M}-3\tau^{P}\ +\overset{q}{\underset{\text{a=1}}{\sum}%
}\varrho_{\text{aa}}^{M}+\overset{q}{\underset{\text{a,b}=1}{\sum}%
}R_{\text{abab}}^{M})](y_{0})\frac{\partial\phi}{\partial x_{\text{c}}}%
(y_{0})\qquad$I$_{33}\qquad$I$_{331}$

$-\frac{1}{4}\underset{\text{c=1}}{\overset{\text{q}}{\sum}}\Lambda_{\text{c}%
}(y_{0})[\left\Vert \text{X}\right\Vert ^{2}+$ divX $-$ $\underset{\text{a}%
=1}{\overset{q}{\sum}}($X$_{\text{a}})^{2}$ $-$ $\underset{\text{a}%
=1}{\overset{q}{\sum}}\frac{\partial X_{\text{a}}}{\partial x_{\text{a}}%
}](y_{0})\frac{\partial\phi}{\partial x_{\text{c}}}(y_{0})+$ $\frac{1}%
{2}\underset{\text{c=1}}{\overset{\text{q}}{\sum}}\Lambda_{\text{c}}(y_{0}%
)$V$(y_{0})\frac{\partial\phi}{\partial x_{\text{c}}}(y_{0})$

$+\frac{1}{2}\underset{\text{c=1}}{\overset{\text{q}}{\sum}}\Lambda_{\text{c}%
}(y_{0})[-(X_{j}\frac{\partial X_{j}}{\partial x_{\text{c}}}+\frac{1}{2}%
\frac{\partial^{2}X_{j}}{\partial x_{\text{c}}\partial x_{j}})(y_{0})+\frac
{1}{2}(<H,j>\frac{\partial X_{j}}{\partial x_{\text{c}}})(y_{0})+\frac
{\partial\text{V}}{\partial x_{\text{c}}}(y_{0})]\phi(y_{0})$

$+\frac{1}{4}\underset{\text{a,c=1}}{\overset{\text{q}}{\sum}}\Lambda
_{\text{c}}(y_{0})\left\{  \frac{\partial^{3}\phi}{\partial\text{x}_{\text{a}%
}^{2}\partial x_{\text{c}}}\text{ }\right\}  (y_{0})\qquad$I$_{332}$

$+\frac{1}{2}\underset{\text{a,b=1}}{\overset{\text{q}}{\sum}}\left\{
\Lambda_{\text{a}}(y_{0})\Lambda_{\text{b}}(y_{0})\frac{\partial^{2}\phi
}{\partial\text{x}_{\text{a}}\partial x_{\text{b}}}\text{ }\right\}
(y_{0})\qquad$I$_{334}$

$+\frac{1}{4}\underset{\text{c=1}}{\overset{\text{q}}{\sum}}%
[\underset{i=1}{\overset{n}{\sum}}\Lambda_{\text{c}}(y_{0})\Lambda_{i}%
^{2}(y_{0})\frac{\partial\phi}{\partial x_{\text{c}}}(y_{0})\qquad$I$_{335}$

$-\frac{1}{4}\underset{j=q+1}{\overset{n}{\sum}}\underset{\text{a,c=1}%
}{\overset{\text{q}}{\sum}}\Lambda_{\text{c}}(y_{0})T_{\text{aa}j}%
(y_{0})\left\{  \Lambda_{j}\frac{\partial\phi}{\partial x_{\text{c}}}\right\}
(y_{0})\qquad$I$_{336}$

$-$ $\frac{1}{2}\underset{j=q+1}{\overset{n}{\sum}}\underset{\text{c=1}%
}{\overset{\text{q}}{\sum}}\Lambda_{\text{c}}(y_{0})\frac{\partial X_{j}%
}{\partial x_{\text{a}}}(y_{0})\Lambda_{\text{c}}(y_{0})\Lambda_{j}(y_{0}%
)\phi(y_{0})\qquad$I$_{338}$

$+\frac{1}{4}\underset{\text{c=1}}{\overset{\text{q}}{\sum}}\left[
\Lambda_{\text{c}}\frac{\partial\text{W}}{\partial x_{\text{c}}}\phi
+\Lambda_{\text{c}}\text{W}\frac{\partial\phi}{\partial x_{\text{c}}}\right]
(y_{0})\qquad$I$_{339}$

$+$ $\underset{\text{a,c=1}}{\overset{\text{q}}{\sum}}[\Lambda_{\text{c}}%
\frac{\partial\text{X}_{\text{a}}}{\partial x_{\text{a}}}](y_{0}%
)[\frac{\partial\phi}{\partial\text{x}_{\text{a}}}+$ X$_{\text{a}}%
\frac{\partial^{2}\phi}{\partial\text{x}_{\text{a}}^{2}}](y_{0})\qquad$E$_{1}$

$+$ $\underset{\text{b=1}}{\overset{\text{q}}{\sum}}$ $[\frac{\partial
\text{X}_{\text{b}}}{\partial x_{\text{a}}}\Lambda_{\text{b}}\Lambda
_{\text{c}}](y_{0})\phi(y_{0})+$ $\underset{\text{b=1}}{\overset{\text{q}%
}{\sum}}$ $[$X$_{\text{b}}\Lambda_{\text{c}}\Lambda_{\text{b}}\frac
{\partial\phi}{\partial x_{\text{a}}}](y_{0})\qquad$E$_{2}$

$+$ $\underset{j=q+1}{\overset{\text{n}}{\sum}}$ $[\frac{\partial X_{j}%
}{\partial x_{\text{a}}}\Lambda_{\text{c}}\Lambda_{j}](y_{0})\phi(y_{0})+$
$\underset{j=q+1}{\overset{\text{n}}{\sum}}$ $\Lambda_{\text{c}}%
(y_{0})[\Lambda_{j}\frac{\partial\phi}{\partial x_{\text{a}}}](y_{0})$

$+\frac{1}{96}\underset{\text{c=1}}{\overset{\text{q}}{\sum}}\Lambda
_{\text{c}}^{2}(y_{0})[\underset{i=q+1}{\overset{n}{\sum}}3<H,i>^{2}%
+2(\tau^{M}-3\tau^{P}\ +\overset{q}{\underset{\text{a=1}}{\sum}}%
\varrho_{\text{aa}}^{M}+\overset{q}{\underset{\text{a,b}=1}{\sum}%
}R_{\text{abab}}^{M})](y_{0})\phi\left(  y_{0}\right)  \qquad$I$_{34}$

$-\frac{1}{8}\underset{\text{c=1}}{\overset{\text{q}}{\sum}}\Lambda_{\text{c}%
}^{2}(y_{0})[$ $\left\Vert \text{X}\right\Vert _{M}^{2}+\frac{1}{2}$
$\operatorname{div}X_{M}-\frac{1}{2}$ $\left\Vert \text{X}\right\Vert _{P}%
^{2}$ $-$ $\frac{1}{2}\operatorname{div}X_{P}](y_{0})\phi\left(  y_{0}\right)
$

$+$ $\frac{1}{8}\underset{\text{c=1}}{\overset{\text{q}}{\sum}}\Lambda
_{\text{c}}^{2}(y_{0})[\underset{\text{a=1}}{\overset{\text{q}}{\sum}}%
\frac{\partial^{2}\phi}{\partial\text{x}_{\text{a}}^{2}}$ $+$ $2$
$\underset{\text{a=1}}{\overset{\text{q}}{\sum}}\Lambda_{\text{a}}%
\frac{\partial\phi}{\partial x_{\text{a}}}\ +$ $\underset{\text{a=1}%
}{\overset{\text{q}}{\sum}}\Lambda_{\text{a}}^{2}](y_{0})\phi\left(
y_{0}\right)  $

$+\frac{1}{4}\underset{\text{c=1}}{\overset{\text{q}}{\sum}}\Lambda_{\text{c}%
}^{2}(y_{0})[$ $\underset{\text{a=1}}{\overset{\text{q}}{\sum}}X_{\text{a}%
}\frac{\partial\phi}{\partial\text{x}_{\text{a}}}+$ $\underset{\text{a=1}%
}{\overset{\text{q}}{\sum}}$ $X_{\text{a}}\Lambda_{\text{a}}+\frac{1}{2}$W $+$
V$](y_{0})\phi(y_{0})$

$+\frac{1}{96}[3\underset{j=q+1}{\overset{n}{\sum}}<H,j>^{2}+2(\tau^{M}%
-3\tau^{P}\ +\overset{q}{\underset{\text{a=1}}{\sum}}\varrho_{\text{aa}}%
^{M}+\overset{q}{\underset{\text{a,b}=1}{\sum}}R_{\text{abab}}^{M})](y_{0}%
)$W$(y_{0})\phi\left(  y_{0}\right)  \qquad$I$_{35}$

$-\frac{1}{8}[$ $\left\Vert \text{X}\right\Vert _{M}^{2}+$ $\operatorname{div}%
X_{M}-$ $\left\Vert \text{X}\right\Vert _{P}^{2}$ $-\operatorname{div}%
X_{P}](y_{0})$W$(y_{0})\phi\left(  y_{0}\right)  $

$+$ $\frac{1}{8}\underset{\text{a=1}}{\overset{\text{q}}{\sum}}\frac
{\partial^{2}\phi}{\partial\text{x}_{\text{a}}^{2}}(y_{0})$W$(y_{0})$
$+\frac{1}{4}$ $\underset{\text{a=1}}{\overset{\text{q}}{\sum}}\Lambda
_{\text{a}}(y_{0})\frac{\partial\phi}{\partial x_{\text{a}}}\left(
y_{0}\right)  $W$(y_{0})\ +\frac{1}{8}$ $\underset{\text{a=1}%
}{\overset{\text{q}}{\sum}}\Lambda_{\text{a}}^{2}(y_{0})$W$(y_{0})\phi\left(
y_{0}\right)  $

$+\frac{1}{4}$ $\underset{\text{a=1}}{\overset{\text{q}}{\sum}}$X$_{\text{a}%
}(y_{0})\frac{\partial\phi}{\partial\text{x}_{\text{a}}}(y_{0})$%
W$(y_{0})+\frac{1}{4}$ $\underset{\text{a=1}}{\overset{\text{q}}{\sum}}$
X$_{\text{a}}(y_{0})\Lambda_{\text{a}}(y_{0})$W$(y_{0})\phi(y_{0})+\frac{1}%
{8}$W$^{2}\left(  y_{0}\right)  \phi(y_{0})+\frac{1}{4}$ V$(y_{0})$%
W$(y_{0})\phi\left(  y_{0}\right)  $

$+\frac{1}{48}$ $\underset{\text{c=1}}{\overset{\text{q}}{\sum}}$X$_{\text{c}%
}(y_{0})[\underset{\alpha=q+1}{\overset{n}{\sum}}3<H,i>^{2}+2(\tau^{M}%
-3\tau^{P}\ +\overset{q}{\underset{\text{a=1}}{\sum}}\varrho_{\text{aa}}%
^{M}+\overset{q}{\underset{\text{a,b}=1}{\sum}}R_{\text{abab}}^{M}%
)](y_{0})\frac{\partial\phi}{\partial x_{\text{c}}}(y_{0})\qquad$I$_{36}%
\qquad$I$_{361}$

$-\frac{1}{4}\underset{\text{c=1}}{\overset{\text{q}}{\sum}}$X$_{\text{c}%
}(y_{0})[\left\Vert \text{X}\right\Vert ^{2}+$ divX $-$ $\underset{\text{a}%
=1}{\overset{q}{\sum}}($X$_{\text{a}})^{2}$ $-$ $\underset{\text{a}%
=1}{\overset{q}{\sum}}\frac{\partial X_{\text{a}}}{\partial x_{\text{a}}%
}](y_{0})\frac{\partial\phi}{\partial x_{\text{c}}}(y_{0})+$ $\frac{1}%
{2}\underset{\text{c=1}}{\overset{\text{q}}{\sum}}\Lambda_{\text{c}}(y_{0}%
)$V$(y_{0})\frac{\partial\phi}{\partial x_{\text{c}}}(y_{0})$

$+\frac{1}{2}\underset{\text{c=1}}{\overset{\text{q}}{\sum}}$X$_{\text{c}%
}(y_{0})[-(X_{j}\frac{\partial X_{j}}{\partial x_{\text{c}}}+\frac{1}{2}%
\frac{\partial^{2}X_{j}}{\partial x_{\text{c}}\partial x_{j}})(y_{0})+\frac
{1}{2}(<H,j>\frac{\partial X_{j}}{\partial x_{\text{c}}})(y_{0})+\frac
{\partial\text{V}}{\partial x_{\text{c}}}(y_{0})]\phi(y_{0})$

$+\frac{1}{4}\underset{\text{a,c=1}}{\overset{\text{q}}{\sum}}$X$_{\text{c}%
}(y_{0})\frac{\partial^{3}\phi}{\partial\text{x}_{\text{a}}^{2}\partial
x_{\text{c}}}$ $(y_{0})$ \qquad\qquad I$_{362}$

$+\frac{1}{2}\underset{\text{a,c=1}}{\overset{\text{q}}{\sum}}[$X$_{\text{c}%
}\Lambda_{\text{a}}\frac{\partial^{2}\phi}{\partial\text{x}_{\text{a}}\partial
x_{\text{c}}}$ $](y_{0})\qquad\qquad$I$_{364}$

$+\frac{1}{4}\underset{\text{b,c=1}}{\overset{\text{q}}{\sum}}[$X$_{\text{c}%
}\Lambda_{\text{b}}^{2}](y_{0})\frac{\partial\phi}{\partial x_{\text{c}}%
}(y_{0})\qquad$I$_{365}$

$+\frac{1}{4}\underset{\text{a,c=1}}{\overset{\text{q}}{\sum}}\left[
\text{X}_{\text{c}}\frac{\partial\text{W}}{\partial x_{\text{a}}}\phi
+\text{X}_{\text{c}}\text{W}\frac{\partial\phi}{\partial x_{\text{c}}}\right]
(y_{0})\qquad$I$_{369}\qquad$

$+$ $\underset{\text{a,c=1}}{\overset{\text{q}}{\sum}}[$X$_{\text{c}}%
\frac{\partial\text{X}_{\text{a}}}{\partial x_{\text{a}}}[\frac{\partial\phi
}{\partial\text{x}_{\text{a}}}+$ X$_{\text{a}}\frac{\partial^{2}\phi}%
{\partial\text{x}_{\text{a}}^{2}}](y_{0})\qquad\qquad$E$_{1}$

$+\underset{\text{a,b,c=1}}{\overset{\text{q}}{\sum}}$ $[$X$_{\text{c}}%
\frac{\partial\text{X}_{\text{b}}}{\partial x_{\text{a}}}\Lambda_{\text{b}%
}](y_{0})\phi(y_{0})+$ $\underset{\text{a,b,c=1}}{\overset{\text{q}}{\sum}}$
$[$X$_{\text{c}}$X$_{\text{b}}\Lambda_{\text{b}}](y_{0})\frac{\partial\phi
}{\partial x_{\text{a}}}(y_{0})\qquad$E$_{2}$

$+\frac{1}{48}\underset{\text{a=1}}{\overset{\text{q}}{\sum}}X_{\text{a}%
}(y_{0})\Lambda_{\text{a}}(y_{0})[\underset{i=q+1}{\overset{n}{\sum}%
}3<H,i>^{2}+2(\tau^{M}-3\tau^{P}\ +\overset{q}{\underset{\text{a=1}}{\sum}%
}\varrho_{\text{aa}}^{M}+\overset{q}{\underset{\text{a,b}=1}{\sum}%
}R_{\text{abab}}^{M})](y_{0})\phi\left(  y_{0}\right)  \qquad$I$_{37}$

$-\frac{1}{4}\underset{\text{a=1}}{\overset{\text{q}}{\sum}}X_{\text{a}}%
(y_{0})\Lambda_{\text{a}}(y_{0})[$ $\left\Vert \text{X}\right\Vert _{M}%
^{2}+\frac{1}{2}$ $\operatorname{div}X_{M}-\frac{1}{2}$ $\left\Vert
\text{X}\right\Vert _{P}^{2}$ $-$ $\frac{1}{2}\operatorname{div}X_{P}%
](y_{0})\phi\left(  y_{0}\right)  $

$+$ $\frac{1}{4}\underset{\text{a,b=1}}{\overset{\text{q}}{\sum}}X_{\text{a}%
}(y_{0})\Lambda_{\text{a}}(y_{0})[\frac{\partial^{2}\phi}{\partial
\text{x}_{\text{b}}^{2}}$ $+$ $\Lambda_{\text{b}}\frac{\partial\phi}{\partial
x_{\text{b}}}]\left(  y_{0}\right)  \ +\frac{1}{4}$ $\underset{\text{a,b=1}%
}{\overset{\text{q}}{\sum}}[X_{\text{a}}\Lambda_{\text{a}}\Lambda_{\text{b}%
}^{2}](y_{0})\phi\left(  y_{0}\right)  $

$+\frac{1}{4}$ $\underset{\text{a,b=1}}{\overset{\text{q}}{\sum}}%
[$X$_{\text{a}}\Lambda_{\text{a}}$ $X_{\text{b}}](y_{0})\frac{\partial\phi
}{\partial\text{x}_{\text{b}}}(y_{0})$

$+\frac{1}{2}$ $\underset{\text{a,b=1}}{\overset{\text{q}}{\sum}}$
$[X_{\text{a}}X_{\text{b}}\Lambda_{\text{a}}\Lambda_{\text{b}}](y_{0}%
)\phi(y_{0})+\frac{1}{4}$ $\underset{\text{a=1}}{\overset{\text{q}}{\sum}}$
$[X_{\text{a}}\Lambda_{\text{a}}$W$]\left(  y_{0}\right)  \phi\left(
y_{0}\right)  $

$+\frac{1}{2}$ $\underset{\text{a=1}}{\overset{\text{q}}{\sum}}$
$[X_{\text{a}}\Lambda_{\text{a}}$V$](y_{0})\phi(y_{0})$

$+\frac{1}{48}$ $\underset{j=q+1}{\overset{n}{\sum}}[$X$_{j}\Lambda_{j}%
](y_{0})[\underset{i=q+1}{\overset{n}{\sum}}3<H,i>^{2}+2(\tau^{M}-3\tau
^{P}\ +\overset{q}{\underset{\text{a=1}}{\sum}}\varrho_{\text{aa}}%
^{M}+\overset{q}{\underset{\text{a,b}=1}{\sum}}R_{\text{abab}}^{M}%
)](y_{0})\phi\left(  y_{0}\right)  $

$-\frac{1}{4}\underset{j=q+1}{\overset{n}{\sum}}[X_{j}\Lambda_{j}](y_{0})[$
$\left\Vert \text{X}\right\Vert _{M}^{2}+\frac{1}{2}$ $\operatorname{div}%
X_{M}-\frac{1}{2}$ $\left\Vert \text{X}\right\Vert _{P}^{2}$ $-$ $\frac{1}%
{2}\operatorname{div}X_{P}](y_{0})\phi\left(  y_{0}\right)  $

$+$ $\frac{1}{4}\underset{j=q+1}{\overset{n}{\sum}}\underset{\text{a=1}%
}{\overset{\text{q}}{\sum}}[X_{j}\Lambda_{j}](y_{0})\frac{\partial^{2}\phi
}{\partial\text{x}_{\text{a}}^{2}}(y_{0})$ $+\frac{1}{2}%
\underset{j=q+1}{\overset{n}{\sum}}\underset{\text{a=1}}{\overset{\text{q}%
}{\sum}}[X_{j}\Lambda_{j}\Lambda_{\text{a}}](y_{0})\frac{\partial\phi
}{\partial x_{\text{a}}}\left(  y_{0}\right)  $

$+\frac{1}{4}$ $\underset{j=q+1}{\overset{n}{\sum}}\underset{\text{a=1}%
}{\overset{\text{q}}{\sum}}[X_{j}\Lambda_{j}\Lambda_{\text{a}}^{2}](y_{0}%
)\phi\left(  y_{0}\right)  +$ $\frac{1}{4}$ $\underset{j=q+1}{\overset{n}{\sum
}}\underset{\text{a=1}}{\overset{\text{q}}{\sum}}[X_{j}\Lambda_{j}X_{\text{a}%
}](y_{0})\frac{\partial\phi}{\partial\text{x}_{\text{a}}}(y_{0})$

$+$ $\frac{1}{2}$ $\underset{j=q+1}{\overset{n}{\sum}}\underset{\text{a=1}%
}{\overset{\text{q}}{\sum}}[X_{\text{a}}\Lambda_{\text{a}}X_{j}\Lambda
_{j}](y_{0})$ $\phi(y_{0})+\frac{1}{4}$ $\underset{j=q+1}{\overset{n}{\sum}%
}[X_{j}\Lambda_{j}$W$]\left(  y_{0}\right)  \phi\left(  y_{0}\right)  $

$+$ $\frac{1}{2}$ $\underset{j=q+1}{\overset{n}{\sum}}[X_{j}\Lambda_{j}%
$V$](y_{0})\phi\left(  y_{0}\right)  $

\begin{proof}
From Appendix $\left(  D_{57}\right)  $\qquad\qquad\qquad\qquad\qquad
\qquad\qquad\qquad\qquad\qquad\qquad\qquad
\end{proof}

\qquad\qquad\qquad\qquad\qquad\qquad\qquad\qquad\qquad\qquad\qquad\qquad
\qquad\qquad\qquad\qquad$\blacksquare$

The generalized heat kernel expansion coefficients exhibit the geometric
invariants of the Riemannian manifold M, the submanifold P and the vector
bundle E. There are possible simplifications. However the expression is too
long and unwieldly. The expression will look more elegant if we assume that
the submanifold P is \textbf{totally geodesic. }In fact some light will be
thrown into the somewhat "dark jungle" of terms in the theorem above if we
make the assumption that the submanifold P is \textbf{totally geodesic. }We
will look for simplifications in this case.

A submanifold P of a Riemannian manifold M is totally geodesic if the
geodesics in P are also geodesics in M. An important consequence is that the
Second Fundamental Form \textbf{T of P vanishes}. The Mean Curvature H also
vanishes since it is defined by:

\qquad\qquad H = $\underset{\text{a=1}}{\overset{\text{q}}{\sum}}%
$T$_{\text{aa}}$

We recall here the Gauss equation: for a,b,c,d = 1,...,q,

$\qquad\ $ $\overset{n}{\underset{i=q+1}{\sum}}$(T$_{\text{ac}i}$%
T$_{\text{bd}i}$ - T$_{\text{ad}i}$T$_{\text{bc}i}$) = R$_{\text{abcd}%
}^{\text{P}}$ - R$_{\text{abcd}}^{\text{M}}$

From the Gauss equation above, another consequence of the vaniahing of the
second fundamental form is that:

$\qquad$R$_{\text{abcd}}^{\text{M}}$ = R$_{\text{abcd}}^{\text{P}}$

\qquad$\overset{q}{\underset{\text{a=1}}{\sum}}R_{\text{abac}}^{M}%
=\overset{q}{\underset{\text{a=1}}{\sum}}R_{\text{abac}}^{P}=\varrho
_{\text{bc}}^{P}$

\ \ \ $\overset{q}{\underset{\text{a,b=1}}{\sum}}R_{\text{abab}}%
^{M}=\overset{q}{\underset{\text{a,b=1}}{\sum}}R_{\text{abab}}^{P}=\tau^{P}$

We assume from now henceforth that the submanifold P is totally geodesic in M.
Hence, the Second Fundamental Form operator T and the related Mean Curvature
Field H vanish identically: T = 0 = H. Further R$_{\text{abcd}}^{\text{M}}$ =
R$_{\text{abcd}}^{\text{P}}.$

As a conequence the expressions for b$_{1}$($y_{0}$,P,$\phi$)(y$_{0})$ and of
b$_{2}$($y_{0}$,P,$\phi$)(y$_{0})$ become shorter and slightly more elegant.

\qquad\qquad\qquad\qquad\qquad\qquad\qquad\qquad\qquad\qquad\qquad\qquad
\qquad\qquad\qquad\qquad\qquad\qquad\qquad$\blacksquare$

\begin{corollary}
Reduced Version: Totally Geodesic Submanifold$\qquad$
\end{corollary}

b$_{2}($y$_{0}$,P$,\phi)=$ I$_{1}+$ I$_{31}+$ I$_{32}+$ I$_{33}+$ I$_{34}+$
I$_{35}+$ I$_{36}+$ I$_{37}$

$=\frac{1}{2}[\frac{1}{24}\left(  2(\tau^{M}-3\tau^{P}%
\ +\overset{q}{\underset{\text{a=1}}{\sum}}\varrho_{\text{aa}}^{M}%
+\overset{q}{\underset{\text{a,b}=1}{\sum}}R_{\text{abab}}^{M})\right)
\qquad\qquad$I$_{1}$

$-\frac{1}{2}\left(  \left\Vert \text{X}\right\Vert ^{2}+\text{divX}%
-\underset{\text{a}=1}{\overset{q}{\sum}}\text{X}_{\text{a}}^{2}%
-\underset{\text{a}=1}{\overset{q}{\sum}}\frac{\partial X_{\text{a}}}{\partial
x_{\text{a}}}-2\text{\textbf{V}}\right)  ]^{2}(y_{0})\phi\left(  y_{0}\right)
$

$+\frac{1}{2}[\frac{1}{24}\left(  2(\tau^{M}-3\tau^{P}%
\ +\overset{q}{\underset{\text{a=1}}{\sum}}\varrho_{\text{aa}}^{M}%
+\overset{q}{\underset{\text{a,b}=1}{\sum}}R_{\text{abab}}^{M})\right)  $

$-\frac{1}{2}\left(  \left\Vert \text{X}\right\Vert ^{2}+\text{divX}%
-\underset{\text{a}=1}{\overset{q}{\sum}}\text{X}_{\text{a}}^{2}%
-\underset{\text{a}=1}{\overset{q}{\sum}}\frac{\partial X_{\text{a}}}{\partial
x_{\text{a}}}-2\text{\textbf{V}}\right)  ](y_{0})$

$\times\lbrack$ $\frac{1}{2}\underset{\text{a=1}}{\overset{\text{q}}{\sum}%
}\frac{\partial^{2}\phi}{\partial\text{x}_{\text{a}}^{2}}(y_{0})$ $+$
$\underset{\text{a=1}}{\overset{\text{q}}{\sum}}\Lambda_{\text{a}}(y_{0}%
)\frac{\partial\phi}{\partial x_{\text{a}}}\left(  y_{0}\right)  \ +\frac
{1}{2}$ $\underset{\text{a=1}}{\overset{\text{q}}{\sum}}(\Lambda_{\text{a}%
}\Lambda_{\text{a}})(y_{0})\phi\left(  y_{0}\right)  $

$+$ $\underset{\text{a=1}}{\overset{\text{q}}{\sum}}$X$_{\text{a}}(y_{0}%
)\frac{\partial\phi}{\partial\text{x}_{\text{a}}}(y_{0})$ +
$\underset{\text{a=1}}{\overset{\text{q}}{\sum}}$ X$_{\text{a}}(y_{0}%
)\Lambda_{\text{a}}(y_{0})\phi(y_{0})+\frac{1}{2}$W$\left(  y_{0}\right)
\phi\left(  y_{0}\right)  ]$

$+\frac{1}{144}[(\tau^{M}-3\tau^{P}\ +\overset{q}{\underset{\text{a=1}}{\sum}%
}\varrho_{\text{aa}}^{M}+\overset{q}{\underset{\text{a,b}=1}{\sum}%
}R_{\text{abab}}^{M})](y_{0}).\frac{\partial^{2}\phi}{\partial x_{\text{c}%
}^{2}}(y_{0})\qquad$I$_{311}\qquad$I$_{31}$

$-\frac{1}{12}[\left\Vert \text{X}(y_{0})\right\Vert ^{2}+$ divX$(y_{0})-$
$\underset{\text{a}=1}{\overset{q}{\sum}}($X$_{\text{a}})^{2}(y_{0})$ $-$
$\underset{\text{a}=1}{\overset{q}{\sum}}\frac{\partial X_{\text{a}}}{\partial
x_{\text{a}}}(y_{0})]\frac{\partial^{2}\phi}{\partial x_{\text{c}}^{2}}%
(y_{0})+$ $\frac{1}{6}$V(y$_{0}$)$\frac{\partial^{2}\phi}{\partial
x_{\text{c}}^{2}}(y_{0})$

$-\frac{1}{6}$ $[X_{j}\frac{\partial X_{j}}{\partial x_{\text{c}}}+\frac{1}%
{2}\frac{\partial^{2}X_{j}}{\partial x_{\text{c}}\partial x_{j}}](y_{0}%
).\frac{\partial\phi}{\partial x_{\text{c}}}(y_{0})+$ $\frac{1}{6}%
\frac{\partial\text{V}}{\partial x_{\text{c}}}(y_{0}).\frac{\partial\phi
}{\partial x_{\text{c}}}(y_{0})$

$+\frac{1}{12}[(\frac{\partial X_{j}}{\partial x_{\text{c}}})^{2}-X_{j}%
\frac{\partial^{2}X_{j}}{\partial x_{\text{c}}^{2}}](y_{0})\phi(y_{0})$
$-\frac{1}{24}\frac{\partial^{3}X_{j}}{\partial x_{\text{c}}^{2}\partial
x_{j}}(y_{0})\phi(y_{0})-\frac{1}{6}\frac{\partial X_{i}}{\partial
x_{\text{c}}}(y_{0}))\frac{\partial X_{i}}{\partial x_{\text{c}}}(y_{0}%
)\phi(y_{0})$

$+\frac{1}{12}\frac{\partial^{2}\text{V}}{\partial x_{\text{c}}^{2}}%
(y_{0})\phi(y_{0})$

$+\frac{1}{24}\underset{\text{a=1}}{\overset{\text{q}}{\sum}}\frac
{\partial^{4}\phi}{\partial x_{\text{a}}^{2}\partial x_{\text{c}}^{2}}%
(y_{0})\qquad\qquad\qquad$I$_{312}$

$+\frac{1}{12}\underset{\text{a=1}}{\overset{\text{q}}{\sum}}[\Lambda
_{\text{a}}\frac{\partial^{3}\phi}{\partial\text{x}_{\text{a}}\partial
x_{\text{c}}^{2}}](y_{0})\qquad\qquad$I$_{314}$

$+$ $\frac{1}{24}\underset{\text{a=1}}{\overset{\text{q}}{\sum}}%
[\Lambda_{\text{a}}^{2}(y_{0})\frac{\partial^{2}\phi}{\partial x_{\text{c}%
}^{2}}](y_{0})\qquad\ \ \ \ $I$_{315}$

$+\frac{1}{24}\frac{\partial^{2}\text{W}}{\partial x_{\text{a}}^{2}}%
(y_{0})\phi(y_{0})+\frac{1}{12}\frac{\partial\text{W}}{\partial x_{\text{a}}%
}(y_{0})\frac{\partial\phi}{\partial x_{\text{a}}}(y_{0})+\frac{1}{24}%
$W$(y_{0})\frac{\partial^{2}\phi}{\partial x_{\text{a}}^{2}}(y_{0})\qquad
$I$_{319}$

$+\frac{1}{12}$ $\underset{\text{a=1}}{\overset{\text{q}}{\sum}}\frac
{\partial^{2}\text{X}_{\text{a}}}{\partial x_{\text{a}}^{2}}(y_{0}%
)\frac{\partial\phi}{\partial\text{x}_{\text{a}}}(y_{0})+\frac{1}{12}$
$\underset{\text{a=1}}{\overset{\text{q}}{\sum}}$X$_{\text{a}}(y_{0}%
)\frac{\partial^{3}\phi}{\partial x_{\text{a}}^{3}}(y_{0})+\frac{1}{6}$
$\underset{\text{a=1}}{\overset{\text{q}}{\sum}}\frac{\partial\text{X}%
_{\text{a}}}{\partial x_{\text{a}}}(y_{0})\frac{\partial^{2}\phi}%
{\partial\text{x}_{\text{a}}^{2}}(y_{0})\qquad$L$_{1}$

$+\frac{1}{12}$ $\underset{\text{a,b=1}}{\overset{\text{q}}{\sum}}$
$\frac{\partial^{2}X_{\text{b}}}{\partial x_{\text{a}}^{2}}\Lambda_{\text{b}%
}(y_{0})\phi(y_{0})+\frac{1}{12}[$ $\underset{\text{a,b=1}}{\overset{\text{q}%
}{\sum}}$ X$_{\text{b}}(y_{0})\Lambda_{\text{b}}(y_{0})\frac{\partial^{2}\phi
}{\partial x_{\text{a}}^{2}}(y_{0})\qquad\qquad\qquad$L$_{2}$

$+\frac{1}{6}[$ $\underset{\text{a,b=1}}{\overset{\text{q}}{\sum}}$
$\frac{\partial X_{\text{b}}}{\partial x_{\text{a}}}(y_{0})\Lambda_{\text{b}%
}(y_{0})\frac{\partial\phi}{\partial x_{\text{a}}}(y_{0})\qquad\qquad\qquad$

$-\frac{1}{1728}\left\{  (\tau^{M}-3\tau^{P}%
\ +\overset{q}{\underset{\text{a=1}}{\sum}}\varrho_{\text{aa}}^{M}%
+\overset{q}{\underset{\text{a,b}=1}{\sum}}R_{\text{abab}}^{M})\right\}  ^{2}%
$(y$_{0}$)$\phi(y_{0})\qquad$I$_{32}\qquad$I$_{321}\qquad$J$_{1}\qquad
$I$_{3211}\qquad$

$+\frac{1}{216}[\overset{n}{\underset{\alpha=q+1}{\sum}}\underset{\beta
,\gamma=q+1}{\sum}$R$_{\alpha\beta\alpha\gamma}]\times\lbrack(\varrho
_{\beta\gamma}+2\overset{q}{\underset{\text{a}=1}{\sum}}$R$_{\beta
\text{a}\gamma\text{a}}$)$](y_{0})\phi(y_{0})$

$-\frac{1}{72}[($ R$_{\text{a}\alpha\text{a}\lambda}](y_{0})\times
\lbrack(\varrho_{\alpha\lambda}+2\overset{q}{\underset{\text{a}=1}{\sum}%
}R_{\alpha\text{a}\lambda\text{a}}](y_{0})\phi(y_{0})$\qquad\qquad\qquad
\qquad\qquad\ \ 

$-\frac{1}{108}\overset{n}{\underset{\alpha,\lambda\text{=q+1}}{\sum}%
}\overset{q}{\underset{\text{a=1}}{\sum}}$ R$_{\alpha\beta\lambda\beta}%
(y_{0})\times\lbrack(\varrho_{\alpha\lambda}+2\overset{q}{\underset{\text{a}%
=1}{\sum}}R_{\alpha\text{a}\lambda\text{a}})](y_{0})\phi(y_{0})\qquad
\qquad\qquad\qquad$

$+\frac{1}{144}[\varrho_{\alpha\beta}+2\overset{q}{\underset{\text{a}=1}{\sum
}}R_{\alpha\text{a}\beta\text{a}}]^{2}(y_{0})\phi(y_{0})\qquad$(k)

$\bigskip+\frac{1}{288}\left\{  \tau-3\tau^{P}+\overset{q}{\underset{\text{a}%
=1}{\sum}}\varrho_{\text{aa}}+\overset{q}{\underset{\text{a,b}=1}{\sum}%
}R_{\text{abab}}\right\}  ^{2}(y_{0})\phi(y_{0})\qquad\qquad\qquad\qquad
\qquad\ \ \ $(l)

$-\frac{1}{288}[-$\ $2\Delta\tau+2$ $\underset{\text{b}=1}{\overset{q}{\sum}%
}\nabla_{\text{bb}}^{2}\tau+\underset{\alpha=q+1}{\overset{n}{\sum}%
}\underset{\text{a}=1}{\overset{q}{\sum}}\nabla_{\alpha\alpha}^{2}%
\varrho_{\text{aa}}+\underset{\beta=q+1}{\overset{n}{\sum}}\underset{\text{a}%
=1}{\overset{q}{\sum}}\nabla_{\beta\beta}^{2}\varrho_{\text{aa}}\qquad
$(m$_{11})+($ m$_{31}$) start

$+\frac{2}{5}\overset{n}{\underset{\alpha,\beta,r=q+1}{\sum}}\nabla
_{\beta\beta}^{2}R_{\alpha r\alpha r}+\frac{2}{5}\overset{n}{\underset{\alpha
,\beta,r=q+1}{\sum}}\nabla_{\alpha\alpha}^{2}R_{\beta r\beta r}\qquad
\qquad\qquad$(m$_{11})+($ m$_{31}$) ends

$-2\Delta\tau$ $+$ $4\underset{\text{a}=1}{\overset{q}{\sum}}%
\underset{\text{j}=q+1}{\overset{n}{\sum}}\nabla_{\text{aj}}^{2}%
\varrho_{\text{aj}}+4\underset{\text{b}=1}{\overset{q}{\sum}}%
\underset{\text{i}=q+1}{\overset{n}{\sum}}\nabla_{\text{ib}}^{2}%
\varrho_{\text{ib}}+4\underset{\text{a,b=1}}{\overset{n}{\sum}}\nabla
_{\text{ab}}^{2}\varrho_{\text{ab}}\qquad\qquad$(m$_{12})+$ (m$_{32})$ starts

$+\frac{4}{5}\underset{\alpha,\beta=q+1}{\overset{n}{\sum}}%
\underset{r=q+1}{\overset{n}{\sum}}\nabla_{\alpha\beta}^{2}R_{\alpha r\beta
r}$ $+\frac{4}{5}\underset{\alpha,\beta=q+1}{\overset{n}{\sum}}%
\underset{r=q+1}{\overset{n}{\sum}}\nabla_{\beta\alpha}^{2}R_{\alpha r\beta
r}\qquad\qquad\qquad\qquad$(m$_{12}$) + (m$_{32}$) ends$\qquad\qquad$

$+$ $\{R_{\alpha\text{a}\alpha s}R_{\beta\text{a}\beta s}+R_{\beta
\text{a}\beta s}R_{\alpha\text{a}\alpha s}+R_{\alpha\text{a}\beta s}%
R_{\alpha\text{a}\beta s}+R_{\alpha\text{a}\beta s}R_{\beta\text{a}\alpha
s}+R_{\beta\text{a}\alpha s}R_{\alpha\text{a}\beta s}\qquad$(m$_{2})$ starts

$+R_{\beta\text{a}\alpha s}R_{\beta\text{a}\alpha s}\}\qquad$\qquad(m$_{2})$ ends

$+\frac{2}{5}\varrho_{rs}^{2}\qquad\qquad\qquad\qquad\qquad\qquad\qquad
\qquad\qquad\qquad\qquad\qquad\ \ \ ($m$_{41})$

$-\frac{2}{5}\left\Vert R^{M}\right\Vert ^{2}+\frac{2}{15}%
\{\overset{n}{\underset{\alpha,\beta=\text{1}}{\sum}}%
\overset{q}{\underset{\text{a,b=1}}{\sum}}(2R_{\alpha\text{a}\beta\text{b}%
}^{2}+R_{\alpha\beta\text{ba}}^{2})\qquad\qquad($m$_{42})+($o$_{2})$ starts

+$\overset{n}{\underset{\alpha,\beta\text{=1a}}{\sum}}%
\overset{q}{\underset{\text{=1}}{\sum}}\overset{n}{\underset{\text{s=q+1}%
}{\sum}}(2R_{\alpha\text{a}\beta\text{s}}^{2}+R_{\alpha\beta\text{sa}}%
^{2})+\overset{n}{\underset{\alpha,\beta=\text{1}}{\sum}}%
\overset{q}{\underset{\text{b=1}}{\sum}}\overset{n}{\underset{\text{r=q+1}%
}{\sum}}(2R_{\alpha\text{r}\beta\text{b}}^{2}+R_{\alpha\beta\text{ra}}^{2})$

$+\overset{q}{\underset{\text{a,b=1}}{\sum}}%
\overset{n}{\underset{r,s\text{=q+1}}{\sum}}(2R_{\text{a}r\text{b}s}%
^{2}+R_{\text{ab}sr}^{2})+\overset{q}{\underset{\text{a=1}}{\sum}%
}\overset{n}{\underset{\beta,r,s\text{=q+1}}{\sum}}(2R_{\text{a}r\beta s}%
^{2}+R_{\text{a}\beta sr}^{2})$

$\bigskip\overset{q}{+\underset{\text{b=1}}{\sum}}%
\overset{n}{\underset{r,s\text{=q+1}}{\sum}}(2R_{\alpha r\text{b}s}%
^{2}+R_{\alpha\text{b}sr}^{2})\}\qquad\qquad\qquad\qquad\qquad\qquad($%
m$_{42})+($o$_{2})$ ends

\ $+\frac{8}{3}\overset{q}{\underset{\text{a},\text{b}=1}{\sum}}%
(\varrho_{\text{aa}}\varrho_{\text{bb}}-\varrho_{\text{aa}}%
\overset{q}{\underset{\text{d}=1}{\sum}}R_{\text{bdbd}}-\varrho_{\text{bb}%
}\overset{q}{\underset{\text{c}=1}{\sum}}R_{\text{acac}}%
+\overset{q}{\underset{\text{c,d}=1}{\sum}}R_{\text{acac}}R_{\text{bdbd}%
})\qquad$(n) starts

\ +$\frac{2}{3}(\tau^{2}-\tau\overset{q}{\underset{\text{a=}1}{\sum}}%
\varrho_{\text{aa}}-\tau\overset{q}{\underset{\text{b=}1}{\sum}}%
\varrho_{\text{bb}}+\overset{q}{\underset{\text{a,b=}1}{\sum}}\varrho
_{\text{aa}}\varrho_{\text{bb}})$

\ + $\frac{4}{3}\{$($\tau\varrho_{\text{aa}}-\varrho_{\text{aa}}%
\overset{q}{\underset{\text{c}=1}{\sum}}\varrho_{\text{cc}}-\tau
\overset{q}{\underset{\text{b}=1}{\sum}}R_{\text{baba}}%
+\overset{q}{\underset{\text{b,c}=1}{\sum}}\varrho_{\text{cc}}R_{\text{baba}%
})$

+ ($\tau\varrho_{\text{bb}}-\varrho_{\text{bb}}\overset{q}{\underset{\text{c}%
=1}{\sum}}\varrho_{\text{cc}}-\tau\overset{q}{\underset{\text{a}=1}{\sum}%
}R_{\text{abab}}+\overset{q}{\underset{\text{a,c}=1}{\sum}}\varrho_{\text{cc}%
}R_{\text{abab}})\}$

$-6(\varrho_{\text{ab}}^{2}-\varrho_{\text{ab}}\overset{q}{\underset{\text{c}%
=1}{\sum}}R_{\text{acbc}}-\varrho_{\text{ab}}\underset{\text{d}%
=1}{\overset{q}{\sum}}R_{\text{adbd}}+\underset{\text{a,b,c,d}%
=1}{\overset{q}{\sum}}R_{\text{acbc}}R_{\text{adbd}})$

\vspace{1pt}$-2(\varrho_{\text{ar}}^{2}-\varrho_{\text{ar}}%
\overset{q}{\underset{\text{c}=1}{\sum}}R_{\text{acrc}}-\varrho_{\text{ar}%
}\underset{\text{d}=1}{\overset{q}{\sum}}R_{\text{adrd}}+\underset{\text{c,d}%
=1}{\overset{q}{\sum}}R_{\text{acrc}}R_{\text{adrd}})$

$-2(\varrho_{\text{br}}^{2}-\varrho_{\text{br}}\overset{q}{\underset{\text{c}%
=1}{\sum}}R_{\text{bcrc}}-\varrho_{\text{br}}\underset{\text{d}%
=1}{\overset{q}{\sum}}R_{\text{bdrd}}+\underset{\text{c,d}=1}{\overset{q}{\sum
}}R_{\text{bcrc}}R_{\text{bdrd}})\qquad$

$+\frac{4}{3}(\varrho_{\alpha\beta}^{2})+$ $\frac{16}{3}%
\overset{q}{\underset{\text{a},\text{b}=1}{\sum}}(R_{\alpha\text{a}%
\beta\text{a}}R_{\alpha\text{b}\beta\text{b}})+\frac{16}{3}%
\overset{q}{\underset{\text{a}=1}{\sum}}(R_{\alpha\text{a}\beta\text{a}%
}\varrho_{\alpha\beta})$

$-3\overset{q}{\underset{\text{a}=1}{\sum}}(R_{\alpha\text{a}\beta r}%
+R_{\beta\text{a}\alpha r})^{2}-2\overset{q}{\underset{\text{a}=1}{\sum}%
}(R_{\alpha\text{a}\beta r}+R_{\beta\text{a}\alpha r})^{2}\qquad\qquad$(n) ends

$\ -\frac{2}{3}(\varrho_{\text{rs}}^{2}-\varrho_{\text{rs}}%
\overset{q}{\underset{\text{b}=1}{\sum}}R_{\text{brbs}}-\varrho_{\text{rs}%
}\underset{\text{d}=1}{\overset{q}{\sum}}R_{\text{aras}}%
+\overset{q}{\underset{\text{a,b=1}}{\sum}}R_{\text{aras}}R_{\text{brbs}%
})\qquad$(o$_{1}$)

$+$ $6\{$ $2\tau\overset{q}{(\underset{b,c=1}{\sum}}R_{\text{bcbc}}%
^{M}-R_{\text{bcbc}}^{P})-2$ $\underset{\text{a=1}}{\overset{n}{\sum}}%
\varrho_{\text{aa}}\overset{q}{(\underset{b,c=1}{\sum}}R_{\text{bcbc}}%
^{M}-R_{\text{bcbc}}^{P})\qquad$(p) starts

$+\underset{\alpha=q+1}{\overset{n}{\sum}}\overset{q}{\underset{\text{a}%
,b,c=1}{\sum}}R_{\alpha\text{a}\alpha b}(R_{\text{acbc}}^{P}-R_{\text{acbc}%
}^{M})+\underset{\beta=q+1}{\overset{n}{\sum}}\overset{q}{\underset{\text{a}%
,b,c=1}{\sum}}R_{\beta\text{a}\beta b}(R_{\text{acbc}}^{P}-R_{\text{acbc}}%
^{M})$

$+\underset{\alpha=q+1}{\overset{n}{\sum}}\overset{q}{\underset{\text{a}%
,b,c=1}{\sum}}R_{\alpha\text{a}\alpha c}(R_{\text{abcb}}^{P}-R_{\text{abcb}%
}^{M})+$ $\underset{\beta=q+1}{\overset{n}{\sum}}%
\overset{q}{\underset{\text{a},b,c=1}{\sum}}R_{\beta\text{a}\beta
c}(R_{\text{abcb}}^{P}-R_{\text{abcb}}^{M})\}$

\vspace{1pt}$+2\underset{\alpha,r=q+1}{\overset{n}{\sum}}%
\overset{q}{\underset{b,c=1}{\sum}}R_{\alpha r\alpha r}(R_{\text{bcbc}}%
^{M}-R_{\text{bcbc}}^{P})+2\underset{\beta,r=q+1}{\overset{n}{\sum}%
}\overset{q}{\underset{b,c=1}{\sum}}R_{\beta r\beta r}(R_{\text{bcbc}}%
^{M}-R_{\text{bcbc}}^{P})$

$+6\underset{\alpha,r=q+1}{\overset{n}{\sum}}\overset{q}{\underset{b,c=1}{\sum
}}R_{\alpha r\alpha b}(R_{\text{bcrc}}^{M}-R_{\text{bcrc}}^{P}%
)+6\underset{\beta,r=q+1}{\overset{n}{\sum}}\overset{q}{\underset{b,c=1}{\sum
}}R_{\beta r\beta b}(R_{\text{bcrc}}^{M}-R_{\text{bcrc}}^{P})$

\vspace{1pt}+6$\underset{\alpha,r=q+1}{\overset{n}{\sum}}%
\overset{q}{\underset{b,c=1}{\sum}}R_{\alpha r\alpha c}(R_{\text{bcbr}}%
^{P}-R_{\text{bcbr}}^{M})+6\underset{\beta,r=q+1}{\overset{n}{\sum}%
}\overset{q}{\underset{b,c=1}{\sum}}R_{\beta r\beta c}(R_{\text{bcbr}}%
^{P}-R_{\text{bcbr}}^{M})\qquad$(p) ends\qquad\qquad\qquad\qquad\qquad\ \ 

$+\frac{1}{24}[\left\Vert \text{X}\right\Vert ^{2}+$ divX $-\underset{\text{a}%
=1}{\overset{q}{\sum}}$X$_{\text{a}}^{2}-\underset{\text{a}%
=1}{\overset{q}{\sum}}\frac{\partial X_{\text{a}}}{\partial x_{\text{a}}%
}](y_{0})[\left\Vert \text{X}\right\Vert ^{2}-$ divX$-\underset{\text{a}%
=1}{\overset{q}{\sum}}X_{\text{a}}^{2}+\underset{\text{a}=1}{\overset{q}{\sum
}}\frac{\partial X_{\text{a}}}{\partial x_{\text{a}}}](y_{0})\phi(y_{0}%
)\qquad$J$_{2}\qquad$J$_{21}\qquad\qquad\qquad\qquad$

$-$ $\frac{2}{9}R_{i\text{a}ij}(y_{0})\frac{\partial X_{\beta}}{\partial
x_{\text{a}}}(y_{0})\phi\left(  y_{0}\right)  +\frac{1}{18}$R$_{ijik}%
(y_{0})[X_{j}(y_{0})X_{k}(y_{0})-\frac{\partial X_{j}}{\partial x_{k}}%
(y_{0})]\phi\left(  y_{0}\right)  $\ \ J$_{22}$\ 

$+$ $\frac{1}{12}[(\frac{\partial X_{i}}{\partial x_{\text{c}}})^{2}%
+X_{i}\frac{\partial^{2}X_{i}}{\partial x_{\text{c}}^{2}}](y_{0})\phi\left(
y_{0}\right)  $ $-\frac{1}{24}\frac{\partial^{3}X_{i}}{\partial x_{i}\partial
x_{\text{c}}^{2}}(y_{0})\phi\left(  y_{0}\right)  \qquad\qquad\qquad
\qquad\qquad\qquad$

$+\frac{1}{24}[\left\Vert \text{X}\right\Vert ^{2}$ $-\underset{\text{a}%
=1}{\overset{q}{\sum}}$X$_{\text{a}}^{2}][\left\Vert \text{X}\right\Vert
^{2}-\underset{\text{a}=1}{\overset{q}{\sum}}X_{\text{a}}^{2}%
-\operatorname{div}$X $+\underset{\text{a}=1}{\overset{q}{\sum}}\frac{\partial
X_{\text{a}}}{\partial x_{\text{a}}}](y_{0})\phi\left(  y_{0}\right)  $

$-\frac{1}{12}X_{i}(y_{0})X_{j}(y_{0})\frac{\partial X_{i}}{\partial x_{j}%
}(y_{0})\phi\left(  y_{0}\right)  \qquad\qquad\qquad\qquad$

$-$ $\frac{1}{36}\underset{j=q+1}{\overset{n}{\sum}}X_{i}(y_{0})X_{k}%
(y_{0})R_{jijk}(y_{0})\phi\left(  y_{0}\right)  +\frac{1}{24}X_{i}(y_{0}%
)\frac{\partial^{2}X_{i}}{\partial x_{j}^{2}}(y_{0})\phi\left(  y_{0}\right)
$

$-\frac{1}{18}[\left\Vert \text{X}\right\Vert ^{2}-\underset{\text{a}%
=1}{\overset{q}{\sum}}X_{\text{a}}^{2}-\operatorname{div}X+\underset{\text{a}%
=1}{\overset{q}{\sum}}\frac{\partial X_{\text{a}}}{\partial x_{\text{a}}%
}](y_{0})\operatorname{div}$X$(y_{0})\phi\left(  y_{0}\right)  $

$-$ $\frac{1}{24}[X_{i}X_{j}-\frac{\partial X_{i}}{\partial x_{j}}%
](y_{0})\frac{\partial X_{i}}{\partial x_{j}}(y_{0})\phi\left(  y_{0}\right)
$

$-\frac{1}{18}\underset{\text{a=1}}{\overset{\text{q}}{\sum}}%
\underset{i=q+1}{\overset{n}{\sum}}X_{j}(y_{0})X_{k}(y_{0})R_{i\text{a}%
ji}(y_{0})\phi\left(  y_{0}\right)  $

$+\frac{1}{24}[\frac{1}{3}\underset{i=q+1}{\overset{n}{\sum}}X_{j}(y_{0}%
)X_{k}(y_{0})R_{ijik}(y_{0}))+X_{j}(y_{0})\frac{\partial^{2}X_{i}}{\partial
x_{i}\partial x_{j}}(y_{0})]\phi\left(  y_{0}\right)  $

$+\frac{1}{24}X_{j}(y_{0})\underset{i=q+1}{\overset{n}{\sum}}X_{k}%
(y_{0})[-\frac{4}{3}\underset{\text{a=1}}{\overset{\text{q}}{\sum}%
}(R_{i\text{a}jk}+R_{j\text{a}ik})+\frac{1}{3}R_{ijik}(y_{0})]\phi\left(
y_{0}\right)  \qquad\qquad\qquad\qquad\qquad\ \ \qquad$

$-$ $\frac{1}{24}X_{i}(y_{0})\underset{j=q+1}{\overset{n}{\sum}}X_{k}%
(y_{0})[\frac{8}{3}\underset{\text{a=1}}{\overset{\text{q}}{\sum}}%
R_{j\text{a}ji}+\frac{2}{3}R_{jijk}](y_{0})\phi\left(  y_{0}\right)  \qquad$

$+\frac{1}{144}\underset{i=q+1}{\overset{n}{\sum}}X_{k}(y_{0}%
)[\underset{\text{a=1}}{\overset{\text{q}}{\sum}}(3\nabla_{j}$R$_{i\text{a}%
ji}+$ $3\nabla_{i}$R$_{j\text{a}ji})](y_{0})\phi\left(  y_{0}\right)
\qquad\qquad\qquad$

$+\frac{1}{72}\underset{i.j=q+1}{\overset{n}{\sum}}[\nabla_{j}R_{jiik}%
+\nabla_{i}R_{jijk}](y_{0})X_{k}(y_{0})\phi\left(  y_{0}\right)  +\frac{1}%
{36}\underset{j=q+1}{\overset{n}{\sum}}R_{jijk}(y_{0})\frac{\partial X_{k}%
}{\partial x_{i}}(y_{0})\phi\left(  y_{0}\right)  $

$-\frac{1}{24}\underset{i=q+1}{\overset{n}{\sum}}[-\frac{4}{3}%
\underset{\text{a=1}}{\overset{\text{q}}{\sum}}(R_{i\text{a}jk}+R_{j\text{a}%
ik})+\frac{1}{3}R_{ijik}(y_{0})]\frac{\partial X_{k}}{\partial x_{j}}%
(y_{0})\phi\left(  y_{0}\right)  $

$+\frac{1}{24}$ $[\frac{\partial X_{i}}{\partial x_{j}}-X_{i}X_{j}%
](y_{0})\frac{\partial X_{i}}{\partial x_{j}}(y_{0})+\frac{1}{24}X_{j}%
(y_{0})\frac{\partial^{2}X_{i}}{\partial x_{i}\partial x_{j}}(y_{0}%
)\phi\left(  y_{0}\right)  \qquad$

$-\frac{1}{24}\underset{i=q+1}{\overset{n}{\sum}}$ $[-\frac{4}{3}%
\underset{\text{a=1}}{\overset{\text{q}}{\sum}}(R_{i\text{a}jk}+R_{j\text{a}%
ik})+\frac{1}{3}R_{ijik}](y_{0})\frac{\partial X_{k}}{\partial x_{j}}%
(y_{0})\phi\left(  y_{0}\right)  $

$+\frac{1}{24}X_{i}(y_{0})\frac{\partial^{2}X_{i}}{\partial x_{j}^{2}}%
(y_{0})-\frac{1}{24}\frac{\partial^{3}X_{i}}{\partial x_{i}\partial x_{j}^{2}%
}(y_{0})\phi\left(  y_{0}\right)  $

$-\frac{1}{144}\underset{\text{a}=1}{\overset{q}{\sum}}%
\underset{i=q+1}{\overset{n}{\sum}}[\{4\nabla_{i}$R$_{i\text{a}l\text{a}}$ $+$
$2\nabla_{l}$R$_{i\text{a}i\text{a}})$\ $\ +\frac{1}{18}%
\underset{k=q+1}{\overset{n}{\sum}}T_{\text{aa}k}R_{ilik}](y_{0})X_{l}%
(y_{0})\phi\left(  y_{0}\right)  $

$+\frac{1}{72}\underset{i,j=q+1}{\overset{n}{\sum}}[(\nabla_{i}R_{ljij}%
+\nabla_{j}R_{ijil}+\nabla_{l}R_{ijij})](y_{0})X_{l}(y_{0})\phi\left(
y_{0}\right)  $

$+\frac{1}{18}\underset{\text{a}=1}{\overset{q}{\sum}}%
\underset{j=q+1}{\overset{n}{\sum}}[R_{\text{a}jij}\frac{\partial X_{i}%
}{\partial x_{\text{a}}}](y_{0})\phi\left(  y_{0}\right)  -\frac{1}{12}$
$\underset{\text{b}=1}{\overset{q}{\sum}}[R_{\text{a}i\text{a}l}](y_{0}%
)[X_{i}X_{l}-\frac{\partial X_{i}}{\partial x_{l}}](y_{0})\phi\left(
y_{0}\right)  $

$-\frac{1}{12}\underset{j=q+1}{\overset{n}{\sum}}\frac{2}{3}R_{ijlj}%
(y_{0})[X_{i}X_{l}-\frac{\partial X_{i}}{\partial x_{l}}](y_{0})\phi\left(
y_{0}\right)  -\frac{1}{24}X_{i}(y_{0})\frac{\partial^{2}X_{i}}{\partial
x_{\text{a}}^{2}}(y_{0})\phi\left(  y_{0}\right)  $

$-\frac{1}{24}[\left\Vert \text{X}\right\Vert ^{2}$ $-\underset{\text{a}%
=1}{\overset{q}{\sum}}$X$_{\text{a}}^{2}][\left\Vert \text{X}\right\Vert
^{2}-\underset{\text{a}=1}{\overset{q}{\sum}}X_{\text{a}}^{2}%
-\operatorname{div}$X $+\underset{\text{a}=1}{\overset{q}{\sum}}\frac{\partial
X_{\text{a}}}{\partial x_{\text{a}}}](y_{0})\phi\left(  y_{0}\right)  $

$+\frac{1}{12}X_{i}(y_{0})X_{j}(y_{0})\frac{\partial X_{i}}{\partial x_{j}%
}(y_{0})\phi\left(  y_{0}\right)  \qquad\ \qquad\ \qquad\qquad\qquad$

$+$ $\frac{1}{18}X_{i}(y_{0})X_{k}(y_{0})[R_{jijk}](y_{0})\phi\left(
y_{0}\right)  -\frac{1}{12}X_{i}(y_{0})\frac{\partial^{2}X_{i}}{\partial
x_{j}^{2}}(y_{0})\phi\left(  y_{0}\right)  $

$+\frac{1}{12}[R_{\text{a}i\text{a}k}](y_{0})X_{i}(y_{0})X_{k}(y_{0}%
)\phi\left(  y_{0}\right)  +\frac{1}{18}R_{ijkj}(y_{0})X_{i}(y_{0})X_{k}%
(y_{0})\phi\left(  y_{0}\right)  $

$-\frac{1}{12}[\left\Vert \text{X}\right\Vert ^{2}(y_{0})+$ divX$(y_{0})-$
$\underset{\text{a}=1}{\overset{q}{\sum}}$X$_{\text{a}}^{2}(y_{0})-$
$\underset{\text{a}=1}{\overset{q}{\sum}}\frac{\partial X_{\text{a}}}{\partial
x_{\text{a}}}(y_{0})][\left\Vert \text{X}\right\Vert ^{2}(y_{0})-$
$\underset{\text{a}=1}{\overset{q}{\sum}}$X$_{\text{a}}^{2}(y_{0})]\phi\left(
y_{0}\right)  \qquad$J$_{4}\qquad$J$_{41}$

$+\frac{1}{12}\left[  X_{i}X_{j}\frac{\partial X_{j}}{\partial x_{i}}%
-X_{i}^{2}X_{j}^{2}\right]  (y_{0})\phi(y_{0})+\frac{1}{12}X_{j}^{2}%
(y_{0})\frac{\partial X_{i}}{\partial x_{i}}(y_{0})\phi(y_{0})$ $\qquad
\ \qquad$J$_{42}\qquad\qquad\qquad$

$-$ $\frac{1}{36}X_{l}(y_{0})X_{j}\varrho_{jl}(y_{0})\phi(y_{0})+\frac{1}%
{12}X_{i}(y_{0})X_{j}(y_{0})\frac{\partial X_{i}}{\partial x_{j}}(y_{0}%
)\phi(y_{0})-\frac{1}{12}X_{j}(y_{0})\frac{\partial^{2}X_{i}}{\partial
x_{i}\partial x_{j}}(y_{0})\phi(y_{0})$

$+\frac{1}{12}[-X_{j}\frac{\partial^{2}X_{j}}{\partial x_{i}^{2}}](y_{0}%
).\phi(y_{0})+\frac{1}{12}[X_{i}X_{j}-\frac{\partial X_{i}}{\partial x_{j}%
}](y_{0})\frac{\partial X_{j}}{\partial x_{i}}(y_{0}).\phi(y_{0})$

$+$ $\frac{1}{6}[\left\Vert \text{X}\right\Vert ^{2}(y_{0})-$
$\underset{\text{a}=1}{\overset{q}{\sum}}$X$_{\text{a}}^{2}(y_{0})]^{2}%
-\frac{1}{6}X_{i}(y_{0})X_{j}\frac{\partial X_{i}}{\partial x_{j}}(y_{0}%
)\phi(y_{0})-\frac{1}{6}X_{i}(y_{0})X_{j}(y_{0})\frac{\partial X_{j}}{\partial
x_{i}}(y_{0}).\phi(y_{0})\qquad$J$_{43}$

$+\frac{1}{12}\frac{\partial^{2}\text{V}}{\partial x_{i}^{2}}(y_{0}%
).\phi(y_{0})\qquad\qquad$J$_{5}$

$-\frac{1}{12}\underset{i\text{=q+1}}{\overset{\text{n}}{\sum}}%
\underset{\text{a,b=1}}{\overset{\text{q}}{\sum}}R_{\text{a}i\text{b}i}%
(y_{0})\frac{\partial^{2}\phi}{\partial\text{x}_{\text{a}}\partial
\text{x}_{\text{b}}}(y_{0})\qquad$\textbf{I}$_{322}$

$+\frac{1}{72}\overset{\text{n}}{\underset{i=q+1}{\sum}}R_{ijik}(y_{0}%
)\Omega_{jk}(y_{0})\phi(y_{0})\qquad$\ I$_{323}$

$-\frac{1}{12}\underset{i\text{=q+1}}{\overset{\text{n}}{\sum}}$
$\underset{\text{a,b=1}}{\overset{\text{q}}{\sum}}R_{\text{a}i\text{b}i}%
(y_{0})\times\Lambda_{\text{b}}(y_{0})\frac{\partial\phi}{\partial
\text{x}_{\text{a}}}(y_{0})\qquad$I$_{324}\qquad$I$_{3241}$

$-\frac{1}{12}\underset{i\text{=}q+1}{\overset{\text{n}}{\sum}}%
\underset{\text{a,b=1}}{\overset{\text{q}}{\sum}}$ $R_{\text{a}i\text{b}%
i}(y_{0})\times\lbrack\Lambda_{\text{a}}(y_{0})\Lambda_{\text{b}}(y_{0}%
)\phi(y_{0})]\qquad$I$_{3251}\qquad$I$_{325}\qquad\qquad\qquad\qquad\qquad$

$+$ $\frac{1}{48}\underset{i,j\text{=q+1}}{\overset{\text{n}}{\sum}}%
(\Omega_{ij}\Omega_{ij})(y_{0})\phi(y_{0})\qquad\qquad\qquad\qquad\qquad
\qquad$I$_{3252}$

$\mathbf{+}\frac{1}{54}\underset{i,j\text{=q+1}}{\overset{\text{n}}{\sum}%
}\underset{\text{b=1}}{\overset{\text{q}}{\sum}}%
[\overset{q}{\underset{\text{c}=1}{\sum}}\{3\nabla_{i}$R$_{j\text{b}%
ij}+3\nabla_{j}$R$_{i\text{b}ij}\}](y_{0})\frac{\partial\phi}{\partial
\text{x}_{\text{c}}}(y_{0})\qquad$

$\mathbf{+}\frac{1}{144}\underset{i,j\text{=q+1}}{\overset{\text{n}}{\sum}%
}\underset{\text{c=1}}{\overset{\text{q}}{\sum}}[\{3\nabla_{i}$R$_{j\text{c}%
ij}+3\nabla_{j}$R$_{i\text{c}ij}\}](y_{0})\times\Lambda_{\text{c}}(y_{0}%
)\phi(y_{0})\qquad$I$_{32622}\qquad\qquad\qquad$

$\mathbf{-}$ $\frac{1}{24}\underset{\text{a=1}}{\overset{\text{q}}{\sum}%
}R_{\text{a}i\text{a}j}(y_{0})\Omega_{il}(y_{0})\phi(y_{0})$

$\mathbf{-}$ $\frac{1}{36}\underset{j=q+1}{\overset{n}{\sum}}R_{ijlj}%
(y_{0})\Omega_{il}(y_{0})\phi(y_{0})$

$+\ \frac{1}{24}\underset{i=q+1}{\overset{n}{\sum}}\frac{\partial^{2}\text{W}%
}{\partial x_{i}^{2}}(y_{0})\phi(y_{0})\qquad$\textbf{I}$_{327}$

$+\frac{8}{3}\underset{j=q+1}{\overset{n}{\sum}}R_{j\text{a}ji}(y_{0}%
)X_{i}(y_{0})+[2X_{j}\frac{\partial X_{j}}{\partial x_{\text{a}}}%
-\frac{\partial^{2}X_{j}}{\partial x_{\text{a}}\partial x_{j}}](y_{0}%
)\qquad\qquad$

$+$ $\frac{1}{72}[(\varrho_{ij}+2\overset{q}{\underset{\text{a}=1}{\sum}%
}R_{i\text{a}j\text{a}}](y_{0})\Omega_{ij}(y_{0})\phi(y_{0})\qquad$%
I$_{32913}\qquad\qquad$\qquad

$+\frac{2}{9}\underset{j=q+1}{\overset{n}{\sum}}R_{j\text{a}ji}(y_{0}%
)X_{i}(y_{0})\Lambda_{\text{a}}(y_{0})\phi(y_{0})+\frac{1}{12}[2X_{j}%
\frac{\partial X_{j}}{\partial x_{\text{a}}}-\frac{\partial^{2}X_{j}}{\partial
x_{\text{a}}\partial x_{j}}](y_{0})\Lambda_{\text{a}}(y_{0})\phi(y_{0})\qquad
$I$_{3292}\qquad$I$_{32921}$

$-\frac{1}{36}X_{j}(y_{0})$ $\frac{\partial\Omega_{ij}}{\partial x_{i}}%
(y_{0})\phi(y_{0})-\frac{1}{12}\frac{\partial X_{j}}{\partial x_{i}}%
(y_{0})\Omega_{ij}(y_{0})$ $\phi(y_{0})\qquad$I$_{32922}$

$+\frac{1}{12}$ $\underset{\text{a=1}}{\overset{\text{q}}{\sum}}\frac
{\partial^{2}X_{\text{a}}}{\partial x_{i}^{2}}(y_{0})\frac{\partial\phi
}{\partial\text{x}_{\text{a}}}(y_{0})\qquad$\textbf{L}$_{1}$

$+\frac{1}{12}$ $\underset{\text{a}=1}{\overset{\text{q}}{\sum}}$
$\frac{\partial^{2}X_{\text{a}}}{\partial x_{i}^{2}}(y_{0})\Lambda_{\text{a}%
}(y_{0})\phi(y_{0})+\frac{1}{12}\underset{j=q+1}{\overset{n}{\sum}}%
\frac{\partial^{2}X_{j}}{\partial x_{i}^{2}}(y_{0})\Lambda_{j}(y_{0}%
)\phi(y_{0})\qquad$\textbf{L}$_{2}\qquad$\textbf{L}$_{21}$

$+\frac{1}{36}$ $\underset{j=q+1}{\overset{n}{\sum}}$ $X_{j}(y_{0}%
)\frac{\partial\Omega_{ij}}{\partial x_{i}}(y_{0})\phi(y_{0})\qquad$%
\textbf{L}$_{22}$

$+\frac{1}{36}$ $\underset{j=q+1}{\overset{n}{\sum}}$ X$_{j}(y_{0}%
)\frac{\partial\Omega_{ij}}{\partial x_{i}}(y_{0})\phi(y_{0})\qquad$%
\textbf{L}$_{22}$

$+\frac{1}{12}$ $\underset{j=q+1}{\overset{\text{n}}{\sum}}$ $\frac{\partial
X_{j}}{\partial x_{i}}(y_{0})\Omega_{ij}(y_{0})\phi(y_{0})\qquad$%
\textbf{L}$_{23}$

$+\frac{1}{24}$ $\underset{\text{c=1}}{\overset{\text{q}}{\sum}}%
\Lambda_{\text{c}}(y_{0})[(\tau^{M}-3\tau^{P}%
\ +\overset{q}{\underset{\text{a=1}}{\sum}}\varrho_{\text{aa}}^{M}%
+\overset{q}{\underset{\text{a,b}=1}{\sum}}R_{\text{abab}}^{M})](y_{0}%
)\frac{\partial\phi}{\partial x_{\text{c}}}(y_{0})\qquad$I$_{331}\qquad
$I$_{33}$

$-\frac{1}{4}\underset{\text{c=1}}{\overset{\text{q}}{\sum}}\Lambda_{\text{c}%
}(y_{0})[\left\Vert \text{X}\right\Vert ^{2}+$ divX $-$ $\underset{\text{a}%
=1}{\overset{q}{\sum}}($X$_{\text{a}})^{2}$ $-$ $\underset{\text{a}%
=1}{\overset{q}{\sum}}\frac{\partial X_{\text{a}}}{\partial x_{\text{a}}%
}](y_{0})\frac{\partial\phi}{\partial x_{\text{c}}}(y_{0})$

$+$ $\frac{1}{2}\underset{\text{c=1}}{\overset{\text{q}}{\sum}}\Lambda
_{\text{c}}(y_{0})$V$(y_{0})\frac{\partial\phi}{\partial x_{\text{c}}}%
(y_{0})+\frac{1}{2}\underset{\text{c=1}}{\overset{\text{q}}{\sum}}%
\Lambda_{\text{c}}(y_{0})[-(X_{j}\frac{\partial X_{j}}{\partial x_{\text{c}}%
}+\frac{1}{2}\frac{\partial^{2}X_{j}}{\partial x_{\text{c}}\partial x_{j}%
})(y_{0})+\frac{\partial\text{V}}{\partial x_{\text{c}}}(y_{0})]\phi(y_{0})$

$+\frac{1}{4}\underset{\text{a,c=1}}{\overset{\text{q}}{\sum}}\Lambda
_{\text{c}}(y_{0})\left\{  \frac{\partial^{3}\phi}{\partial\text{x}_{\text{a}%
}^{2}\partial x_{\text{c}}}\text{ }\right\}  (y_{0})$ \qquad\qquad I$_{332}$

$+\frac{1}{2}\underset{\text{a,c=1}}{\overset{\text{q}}{\sum}}\left\{
\Lambda_{\text{a}}(y_{0})\Lambda_{\text{c}}(y_{0})\frac{\partial^{2}\phi
}{\partial\text{x}_{\text{a}}\partial x_{\text{c}}}\text{ }\right\}
(y_{0})\qquad$I$_{334}$

$+\frac{1}{4}\underset{\text{a,b,c=1}}{\overset{\text{q}}{\sum}}%
\Lambda_{\text{c}}(y_{0})\Lambda_{\text{b}}^{2}(y_{0})\frac{\partial\phi
}{\partial x_{\text{c}}}(y_{0})\qquad$I$_{335}$

$+\frac{1}{4}\underset{\text{c=1}}{\overset{\text{q}}{\sum}}\left[
\Lambda_{\text{c}}(y_{0})\frac{\partial\text{W}}{\partial x_{\text{c}}}%
(y_{0})\phi(y_{0})+\Lambda_{\text{c}}(y_{0})\text{W}(y_{0})\frac{\partial\phi
}{\partial x_{\text{c}}}(y_{0})\right]  \qquad$I$_{339}$

$+$ $\underset{\text{a,c=1}}{\overset{\text{q}}{\sum}}\Lambda_{\text{c}}%
(y_{0})\frac{\partial\text{X}_{\text{a}}}{\partial x_{\text{a}}}(y_{0}%
)[\frac{\partial\phi}{\partial\text{x}_{\text{a}}}+$ X$_{\text{a}}%
\frac{\partial^{2}\phi}{\partial\text{x}_{\text{a}}^{2}}](y_{0})\qquad\qquad
$E$_{1}$

$+\underset{\text{a,b,c=1}}{\overset{\text{q}}{\sum}}$ $\Lambda_{\text{c}%
}(y_{0})\frac{\partial\text{X}_{\text{b}}}{\partial x_{\text{a}}}%
(y_{0})\Lambda_{\text{b}}(y_{0})\phi(y_{0})+$ $\underset{\text{a,b,c=1}%
}{\overset{\text{q}}{\sum}}$ $\Lambda_{\text{c}}(y_{0})$X$_{\text{b}}%
(y_{0})\Lambda_{\text{b}}(y_{0})\frac{\partial\phi}{\partial x_{\text{a}}%
}(y_{0})\qquad$E$_{2}$

$+\frac{1}{48}\underset{\text{c=1}}{\overset{\text{q}}{\sum}}\Lambda
_{\text{c}}^{2}(y_{0})[(\tau^{M}-3\tau^{P}\ +\overset{q}{\underset{\text{a=1}%
}{\sum}}\varrho_{\text{aa}}^{M}+\overset{q}{\underset{\text{a,b}=1}{\sum}%
}R_{\text{abab}}^{M})](y_{0})\phi\left(  y_{0}\right)  \qquad$I$_{34}$

$-\frac{1}{8}\underset{\text{c=1}}{\overset{\text{q}}{\sum}}\Lambda_{\text{c}%
}^{2}(y_{0})[\left\Vert \text{X}\right\Vert ^{2}+\operatorname{div}$X
$-\underset{\text{a}=1}{\overset{q}{\sum}}$X$_{\text{a}}^{2}%
-\underset{\text{a}=1}{\overset{q}{\sum}}\frac{\partial X_{\text{a}}}{\partial
x_{\text{a}}}]\phi\left(  y_{0}\right)  +$ $\frac{1}{4}\underset{\text{c=1}%
}{\overset{\text{q}}{\sum}}\Lambda_{\text{c}}^{2}(y_{0})$\textbf{V}$\left(
y_{0}\right)  \phi\left(  y_{0}\right)  \ $

$+\frac{1}{8}\underset{\text{c=1}}{\overset{\text{q}}{\sum}}\Lambda_{\text{c}%
}^{2}(y_{0})[$ $\underset{\text{a=1}}{\overset{\text{q}}{\sum}}\frac
{\partial^{2}\phi}{\partial\text{x}_{\text{a}}^{2}}$ $+$ $\underset{\text{a=1}%
}{\overset{\text{q}}{\sum}}\Lambda_{\text{a}}\frac{\partial\phi}{\partial
x_{\text{a}}}\ +\frac{1}{2}$ $\underset{\text{a=1}}{\overset{\text{q}}{\sum}%
}(\Lambda_{\text{a}}\Lambda_{\text{a}})](y_{0})\phi\left(  y_{0}\right)  $

$+\frac{1}{4}\underset{\text{c=1}}{\overset{\text{q}}{\sum}}\Lambda_{\text{c}%
}^{2}(y_{0})[$ $\underset{\text{a=1}}{\overset{\text{q}}{\sum}}$X$_{\text{a}%
}\frac{\partial\phi}{\partial\text{x}_{\text{a}}}$ + $\underset{\text{a=1}%
}{\overset{\text{q}}{\sum}}$ X$_{\text{a}}\Lambda_{\text{a}}\phi(y_{0}%
)+\frac{1}{2}$W$]\left(  y_{0}\right)  \phi\left(  y_{0}\right)  $

$+\frac{1}{48}[(\tau^{M}-3\tau^{P}\ +\overset{q}{\underset{\text{a=1}}{\sum}%
}\varrho_{\text{aa}}^{M}+\overset{q}{\underset{\text{a,b}=1}{\sum}%
}R_{\text{abab}}^{M})](y_{0})\phi\left(  y_{0}\right)  $W$(y_{0})\qquad
$I$_{35}$

$-\frac{1}{8}[\left\Vert \text{X}\right\Vert ^{2}+\operatorname{div}$X
$-\underset{\text{a}=1}{\overset{q}{\sum}}$X$_{\text{a}}^{2}%
-\underset{\text{a}=1}{\overset{q}{\sum}}\frac{\partial X_{\text{a}}}{\partial
x_{\text{a}}}]\phi\left(  y_{0}\right)  $W$(y_{0})\ +$ $\frac{1}{4}$%
\textbf{V}$\left(  y_{0}\right)  $W$(y_{0})\phi\left(  y_{0}\right)  $

$+$ $\frac{1}{8}\underset{\text{a=1}}{\overset{\text{q}}{\sum}}\frac
{\partial^{2}\phi}{\partial\text{x}_{\text{a}}^{2}}(y_{0})$W$(y_{0})$
$+\frac{1}{4}$ $\underset{\text{a=1}}{\overset{\text{q}}{\sum}}\Lambda
_{\text{a}}(y_{0})\frac{\partial\phi}{\partial x_{\text{a}}}\left(
y_{0}\right)  $W$(y_{0})\ +\frac{1}{8}$ $\underset{\text{a=1}%
}{\overset{\text{q}}{\sum}}(\Lambda_{\text{a}}\Lambda_{\text{a}})(y_{0}%
)$W$(y_{0})\phi\left(  y_{0}\right)  $

$+\frac{1}{4}$ $\underset{\text{a=1}}{\overset{\text{q}}{\sum}}$X$_{\text{a}%
}(y_{0})\frac{\partial\phi}{\partial\text{x}_{\text{a}}}(y_{0})$W$(y_{0})$
+$\frac{1}{4}$ $\underset{\text{a=1}}{\overset{\text{q}}{\sum}}$ X$_{\text{a}%
}(y_{0})\Lambda_{\text{a}}(y_{0})$W$(y_{0})\phi(y_{0})+\frac{1}{8}$%
W$^{2}\left(  y_{0}\right)  \phi(y_{0})$

$+\frac{1}{48}$ $\underset{\text{c=1}}{\overset{\text{q}}{\sum}}$X$_{\text{c}%
}(y_{0})[\underset{i=q+1}{\overset{n}{\sum}}3<H,i>^{2}+2(\tau^{M}-3\tau
^{P}\ +\overset{q}{\underset{\text{a=1}}{\sum}}\varrho_{\text{aa}}%
^{M}+\overset{q}{\underset{\text{a,b}=1}{\sum}}R_{\text{abab}}^{M}%
)](y_{0})\frac{\partial\phi}{\partial x_{\text{c}}}(y_{0})\qquad$%
I$_{361}\qquad$I$_{36}$

$-\frac{1}{4}\underset{\text{c=1}}{\overset{\text{q}}{\sum}}$X$_{\text{c}%
}(y_{0})[\left\Vert \text{X}\right\Vert ^{2}+$ divX $-$ $\underset{\text{a}%
=1}{\overset{q}{\sum}}($X$_{\text{a}})^{2}$ $-$ $\underset{\text{a}%
=1}{\overset{q}{\sum}}\frac{\partial X_{\text{a}}}{\partial x_{\text{a}}%
}](y_{0})\frac{\partial\phi}{\partial x_{\text{c}}}(y_{0})+$ $\frac{1}%
{2}\underset{\text{c=1}}{\overset{\text{q}}{\sum}}\Lambda_{\text{c}}(y_{0}%
)$V$(y_{0})\frac{\partial\phi}{\partial x_{\text{c}}}(y_{0})$

$+\frac{1}{2}\underset{\text{c=1}}{\overset{\text{q}}{\sum}}$X$_{\text{c}%
}(y_{0})[-(X_{j}\frac{\partial X_{j}}{\partial x_{\text{c}}}+\frac{1}{2}%
\frac{\partial^{2}X_{j}}{\partial x_{\text{c}}\partial x_{j}})(y_{0}%
)+\frac{\partial\text{V}}{\partial x_{\text{c}}}(y_{0})]\phi(y_{0})$

$+\frac{1}{4}\underset{\text{a,c=1}}{\overset{\text{q}}{\sum}}$X$_{\text{c}%
}(y_{0})\left\{  \frac{\partial^{3}\phi}{\partial\text{x}_{\text{a}}%
^{2}\partial x_{\text{c}}}\text{ }\right\}  (y_{0})$ \qquad\qquad I$_{362}$

$+\frac{1}{2}\underset{\text{a,c=1}}{\overset{\text{q}}{\sum}}\left\{
\text{X}_{\text{c}}(y_{0})\Lambda_{\text{a}}(y_{0})\frac{\partial^{2}\phi
}{\partial\text{x}_{\text{a}}\partial x_{\text{c}}}\text{ }\right\}
(y_{0})\qquad$I$_{364}$

$+\frac{1}{4}\underset{\text{b,c=1}}{\overset{\text{q}}{\sum}}$X$_{\text{c}%
}(y_{0})\Lambda_{\text{b}}^{2}(y_{0})\frac{\partial\phi}{\partial x_{\text{c}%
}}(y_{0})\qquad$I$_{365}$

$+\frac{1}{4}\underset{\text{a,c=1}}{\overset{\text{q}}{\sum}}\left[
\text{X}_{\text{c}}(y_{0})\frac{\partial\text{W}}{\partial x_{\text{a}}}%
(y_{0})\phi(y_{0})+\text{X}_{\text{c}}(y_{0})\text{W}(y_{0})\frac{\partial
\phi}{\partial x_{\text{c}}}(y_{0})\right]  \qquad$I$_{369}\qquad$

$+$ $\underset{\text{a,c=1}}{\overset{\text{q}}{\sum}}$X$_{\text{c}}%
(y_{0})\frac{\partial\text{X}_{\text{a}}}{\partial x_{\text{a}}}(y_{0}%
)[\frac{\partial\phi}{\partial\text{x}_{\text{a}}}+$ X$_{\text{a}}%
\frac{\partial^{2}\phi}{\partial\text{x}_{\text{a}}^{2}}](y_{0})\qquad\qquad
$E$_{1}$

$+\underset{\text{a,b,c=1}}{\overset{\text{q}}{\sum}}$ X$_{\text{c}}%
(y_{0})\frac{\partial\text{X}_{\text{b}}}{\partial x_{\text{a}}}(y_{0}%
)\Lambda_{\text{b}}(y_{0})\phi(y_{0})+$ $\underset{\text{a,b,c=1}%
}{\overset{\text{q}}{\sum}}$ X$_{\text{c}}(y_{0})$X$_{\text{b}}(y_{0}%
)\Lambda_{\text{b}}(y_{0})\frac{\partial\phi}{\partial x_{\text{a}}}%
(y_{0})\qquad$E$_{2}$

$+\frac{1}{48}\underset{\text{a=1}}{\overset{\text{q}}{\sum}}$X$_{\text{a}%
}(y_{0})\Lambda_{\text{a}}(y_{0})[3\underset{j=q+1}{\overset{n}{\sum}%
}<H,j>^{2}+2(\tau^{M}-3\tau^{P}\ +\overset{q}{\underset{\text{a=1}}{\sum}%
}\varrho_{\text{aa}}^{M}+\overset{q}{\underset{\text{a,b}=1}{\sum}%
}R_{\text{abab}}^{M})](y_{0})\phi\left(  y_{0}\right)  \qquad$I$_{37}$

$-\frac{1}{4}\underset{\text{a=1}}{\overset{\text{q}}{\sum}}$X$_{\text{a}%
}(y_{0})\Lambda_{\text{a}}(y_{0})[\left\Vert \text{X}\right\Vert
^{2}+\operatorname{div}$X $-\underset{\text{a}=1}{\overset{q}{\sum}}%
$X$_{\text{a}}^{2}-\underset{\text{a}=1}{\overset{q}{\sum}}\frac{\partial
X_{\text{a}}}{\partial x_{\text{a}}}]\phi\left(  y_{0}\right)  \ $

$+\underset{\text{a=1}}{\overset{\text{q}}{\frac{1}{2}\sum}}$X$_{\text{a}%
}(y_{0})\Lambda_{\text{a}}(y_{0})[$ $\frac{1}{2}\underset{\text{a=1}%
}{\overset{\text{q}}{\sum}}\frac{\partial^{2}\phi}{\partial\text{x}_{\text{a}%
}^{2}}(y_{0})$ $+$ $\underset{\text{a=1}}{\overset{\text{q}}{\sum}}%
\Lambda_{\text{a}}(y_{0})\frac{\partial\phi}{\partial x_{\text{a}}}\left(
y_{0}\right)  \ +\frac{1}{4}$ $\underset{\text{a=1}}{\overset{\text{q}}{\sum}%
}(\Lambda_{\text{a}}\Lambda_{\text{a}})(y_{0})\phi\left(  y_{0}\right)  ]$

$+\frac{1}{2}\underset{\text{a=1}}{\overset{\text{q}}{\sum}}$X$_{\text{a}%
}(y_{0})\Lambda_{\text{a}}(y_{0})[$ $\underset{\text{a=1}}{\overset{\text{q}%
}{\sum}}$X$_{\text{a}}(y_{0})\frac{\partial\phi}{\partial\text{x}_{\text{a}}%
}(y_{0})$ + $\underset{\text{a=1}}{\overset{\text{q}}{\sum}}$ X$_{\text{a}%
}(y_{0})\Lambda_{\text{a}}(y_{0})\phi(y_{0})]$

$+\frac{1}{4}\underset{\text{a=1}}{\overset{\text{q}}{\sum}}$X$_{\text{a}%
}(y_{0})\Lambda_{\text{a}}(y_{0})$W$\left(  y_{0}\right)  \phi\left(
y_{0}\right)  +$ $\underset{\text{a=1}}{\overset{\text{q}}{\frac{1}{2}\sum}}%
$X$_{\text{a}}(y_{0})\Lambda_{\text{a}}(y_{0})$V$\left(  y_{0}\right)
\phi\left(  y_{0}\right)  $

\begin{proof}
\qquad We delete all \textbf{fundamental form-related} items wherever they
occur from the last theorem above.\qquad\qquad\qquad\qquad\qquad\qquad
\qquad\qquad\qquad$\qquad\qquad$
\end{proof}

$\qquad\qquad\qquad\qquad\qquad\qquad\qquad\qquad\qquad\qquad\qquad
\qquad\qquad\qquad\qquad\qquad\qquad\qquad\blacksquare$

In order to compare our work to the work of previous authors, we must give the
expression for b$_{2}($y$_{0}$,P$,\phi)$ in the case that the submanifold P
reduces to the singleton $\left\{  y_{0}\right\}  $ which is the centre of
Fermi coordinates, now reduced to normal coordinates. The formula below
generalizes the third coefficient of the usual vector bundle heat kernel expansion.

\qquad\qquad\qquad\qquad\qquad\qquad\qquad\qquad\qquad\qquad\qquad\qquad
\qquad\qquad\qquad\qquad\qquad\qquad$\blacksquare$

\begin{corollary}
b$_{2}($y$_{0},y_{0})\phi\left(  y_{0}\right)  $
\end{corollary}

$\qquad=\frac{1}{2}[\frac{1}{24}(2(\tau^{M}))-\frac{1}{2}(\left\Vert
\text{X}\right\Vert _{M}^{2}+\operatorname{div}X_{M})+V](y_{0})\phi\left(
y_{0}\right)  \qquad$I$_{1}$

$\qquad\times\lbrack\frac{1}{24}(2(\tau^{M})-\frac{1}{2}($ $\left\Vert
\text{X}\right\Vert _{M}^{2}+$ $\operatorname{div}X_{M})+$ $V+\frac{1}%
{2}W)]\phi\left(  y_{0}\right)  $

$\qquad-\frac{1}{3456}[2(\tau^{M}]^{2}(y_{0})\phi(y_{0})+\frac{1}{24}[\frac
{1}{3}(\tau^{M})](y_{0})\phi(y_{0})\times\lbrack\frac{1}{6}(\tau^{M}%
)](y_{0})\phi(y_{0})$\qquad I$_{321}\qquad\ \ \ \ \ $

\qquad$\ -\frac{1}{864}[2\varrho_{ij}]^{2}(y_{0})\phi(y_{0})\ -\frac{1}%
{432}R_{jijk}(y_{0})[2\varrho_{ik}](y_{0})\phi(y_{0})\qquad L_{2}\qquad
L_{21}\qquad L_{22}\ \ \ \ \ \ \ \ \ \ \ \ \ \ \ \ \ \ \ $

\qquad$\ -\frac{1}{432}[\varrho_{ii}](y_{0})\times\lbrack\varrho_{jj}%
](y_{0})\phi(y_{0})+\frac{1}{432}R_{ijik}(y_{0})[2\varrho_{jk}](y_{0}%
)\phi(y_{0})\qquad L_{23}\qquad L_{233}\qquad\qquad\qquad\ \ \ \ $

$\qquad\ +\frac{1}{576}[2\varrho_{ij})]^{2}(y_{0})\phi(y_{0})+\frac{1}%
{288}[\tau^{M}\ ]^{2}(y_{0})\phi(y_{0})\ \ \ \ \ $

$\qquad-\ \frac{1}{288}\overset{n}{\underset{i,j=1}{\sum}}[-\frac{3}%
{5}\overset{n}{\underset{p=1}{\sum}}\nabla_{ii}^{2}(R)_{jpjp}-\frac{3}%
{5}\overset{n}{\underset{p=1}{\sum}}\nabla_{jj}^{2}(R)_{ipip}-\frac{6}%
{5}\overset{n}{\underset{p=1}{\sum}}\nabla_{ij}^{2}(R)_{ipjp}$

$\qquad-\frac{6}{5}\overset{n}{\underset{p=1}{\sum}}\nabla_{ji}^{2}%
(R)_{jpip}\qquad\qquad$

$\qquad\ \ +\frac{1}{5}\overset{n}{\underset{m,p=q+1}{%
{\textstyle\sum}
}}R_{ipim}R_{jpjm}+\frac{1}{5}\overset{n}{\underset{m,p=q+1}{%
{\textstyle\sum}
}}R_{jpjm}R_{ipim}$ $+\frac{1}{5}\overset{n}{\underset{m,p=q+1}{%
{\textstyle\sum}
}}R_{ipjm}R_{ipjm}$

\qquad$\ +\frac{1}{5}\overset{n}{\underset{m,p=q+1}{%
{\textstyle\sum}
}}R_{ipjm}R_{jpim}$

$\qquad\ +\frac{1}{5}\overset{n}{\underset{m,p=q+1}{%
{\textstyle\sum}
}}R_{jpim}R_{ipjm}+\frac{1}{5}\overset{n}{\underset{m,p=q+1}{%
{\textstyle\sum}
}}R_{jpim}R_{jpim}](y_{0})\phi(y_{0})$

$\qquad-\frac{1}{48}$ $[\frac{1}{9}(\tau^{M})(\tau^{M})+\frac{2}{9}(\left\Vert
\varrho^{M}\right\Vert ^{2})](y_{0})\phi(y_{0})\qquad\qquad3C$

$\qquad-\frac{1}{48}[-\frac{1}{9}\overset{n}{\underset{i,j,p,m=q+1}{\sum}%
}R_{ipim}R_{jpjm}\ -\frac{1}{18}\overset{n}{\underset{i,j,p,m=q+1}{\sum}%
}R_{ipjm}^{2}-\frac{1}{9}\overset{n}{\underset{i,j,p,m=q+1}{\sum}}%
R_{ipjm}R_{jpim}$

$\qquad-\frac{1}{18}\overset{n}{\underset{i,j,p,m=q+1}{\sum}}R_{jpim}%
^{2}](y_{0})\phi(y_{0})\qquad$

\qquad$+\frac{1}{24}[\left\Vert \text{X}\right\Vert _{M}^{2}%
+\operatorname{div}$X$_{M}](y_{0})[\left\Vert \text{X}\right\Vert _{M}%
^{2}-\operatorname{div}$X$_{M}](y_{0})\phi(y_{0})\qquad$\ I$_{3212}$

$\qquad+\frac{1}{12}X_{i}^{2}(y_{0})[\operatorname{div}X_{M}-\left\Vert
X\right\Vert _{M}^{2}](y_{0})\phi(y_{0})\qquad\qquad$I$_{32122}\qquad Q_{2}$

$\qquad+\frac{1}{6}X_{i}(y_{0})X_{j}(y_{0})\frac{\partial X_{i}}{\partial
x_{j}}(y_{0})\phi(y_{0})+$ $\frac{1}{18}X_{i}(y_{0})X_{k}(y_{0})R_{jijk}%
(y_{0})\phi\left(  y_{0}\right)  $

\qquad$-\frac{1}{12}X_{i}(y_{0})\frac{\partial^{2}X_{i}}{\partial x_{j}^{2}%
}(y_{0})\phi(y_{0})+\frac{1}{18}R_{ijkj}(y_{0})X_{i}(y_{0})X_{k}(y_{0}%
)\phi(y_{0})\qquad\qquad\qquad\qquad\qquad\qquad\qquad\qquad\qquad
\ \qquad\qquad\qquad\qquad\qquad\qquad\qquad\qquad\qquad\qquad\qquad
\qquad\qquad\qquad\qquad$

$\qquad+\frac{1}{18}R_{ijik}(y_{0})[X_{j}X_{k}-\frac{1}{2}(\frac{\partial
X_{j}}{\partial x_{k}}+\frac{\partial X_{k}}{\partial x_{j}})](y_{0}%
)\phi(y_{0})\qquad$I$_{32123}\qquad S_{1}\qquad\qquad\qquad\qquad\qquad
\qquad\qquad\qquad\qquad\qquad\qquad\qquad\qquad\qquad$

$\qquad+\frac{1}{24}[-\frac{1}{3}(\nabla_{i}R_{kjij}+\nabla_{j}R_{ijik}%
+\nabla_{k}R_{ijij})](y_{0})X_{k}(y_{0})\phi(y_{0})\qquad S_{32}$

\qquad\ $-\frac{1}{18}R_{ijkj}(y_{0})[X_{i}X_{k}-\frac{1}{2}\left(
\frac{\partial X_{i}}{\partial x_{k}}+\frac{\partial X_{k}}{\partial x_{i}%
}\right)  ](y_{0})\phi(y_{0})$

$\qquad+\frac{1}{24}[X_{i}^{2}X_{j}^{2}-2X_{i}X_{j}\left(  \frac{\partial
X_{j}}{\partial x_{i}}+\frac{\partial X_{i}}{\partial x_{j}}\right)
-X_{i}^{2}\frac{\partial X_{j}}{\partial x_{j}}-X_{j}^{2}\frac{\partial X_{i}%
}{\partial x_{i}}](y_{0})\phi(y_{0})$

$\qquad+\frac{1}{48}\left(  \frac{\partial X_{j}}{\partial x_{i}}%
+\frac{\partial X_{i}}{\partial x_{j}}\right)  ^{2}(y_{0})\phi(y_{0})+\frac
{1}{24}\left(  \frac{\partial X_{i}}{\partial x_{i}}\frac{\partial X_{j}%
}{\partial x_{j}}\right)  (y_{0})\phi(y_{0})$

$\qquad+\frac{1}{36}X_{i}(y_{0})\left(  2\frac{\partial^{2}X_{j}}{\partial
x_{i}\partial x_{j}}+\frac{\partial^{2}X_{i}}{\partial x_{j}^{2}}\right)
(y_{0})+\frac{1}{36}X_{j}(y_{0})\left(  \frac{\partial^{2}X_{j}}{\partial
x_{i}^{2}}+2\frac{\partial^{2}X_{i}}{\partial x_{i}\partial x_{j}}\right)
(y_{0})\phi\left(  y_{0}\right)  $

$\qquad-\frac{1}{48}\left(  \frac{\partial^{3}X_{i}}{\partial x_{i}\partial
x_{j}^{2}}+\frac{\partial^{3}X_{j}}{\partial x_{i}^{2}\partial x_{j}}\right)
(y_{0})\phi\left(  y_{0}\right)  $

\qquad$+\frac{1}{12}[\frac{1}{6}(2\varrho_{ij})](y_{0})\phi(y_{0})\times
\frac{1}{2}[\left(  \frac{\partial X_{j}}{\partial x_{i}}-\frac{\partial
X_{i}}{\partial x_{j}}\right)  ](y_{0})\phi(y_{0})\qquad\qquad$I$_{3213}$

\qquad$-\frac{1}{18}[X_{j}\left(  2\frac{\partial^{2}X_{j}}{\partial x_{i}%
^{2}}+\frac{\partial^{2}X_{i}}{\partial x_{i}\partial x_{j}}\right)
](y_{0})\phi(y_{0})-\frac{1}{12}[\left(  \frac{\partial X_{i}}{\partial x_{j}%
}+\frac{\partial X_{j}}{\partial x_{i}}\right)  ]\frac{\partial X_{j}%
}{\partial x_{i}}(y_{0})\phi(y_{0})\qquad\ $I$_{3214}\qquad$\qquad\qquad
\qquad\qquad\qquad\qquad\qquad\qquad\qquad\qquad\qquad\qquad\qquad\qquad
\qquad\qquad

$\qquad+\frac{1}{12}\frac{\partial^{2}\text{V}}{\partial x_{i}^{2}}(y_{0}%
)\phi(y_{0})\qquad\qquad$I$_{3215}$

\qquad$+\frac{1}{72}\underset{i,j,k=1}{\overset{n}{\sum}}R_{ijik}(y_{0}%
)\Omega_{jk}(y_{0})\phi(y_{0})\qquad\qquad\qquad\qquad$\ I$_{323}$

$\qquad+\frac{1}{48}\underset{i,j=1}{\overset{n}{\sum}}\left(  \Omega
_{ij}\Omega_{ij}\right)  (y_{0})\phi(y_{0})+\frac{1}{72}%
\underset{i,j=1}{\overset{n}{\sum}}\left(  \frac{\partial\Omega_{ij}}%
{\partial\text{x}_{i}}\Lambda_{j}+\Lambda_{j}\frac{\partial\Omega_{ij}%
}{\partial\text{x}_{i}}\right)  (y_{0})\phi(y_{0})\qquad\qquad$I$_{3252}$

\qquad\ $+\frac{1}{24}[\frac{1}{3}(\nabla_{i}R_{kjij}+\nabla_{j}%
R_{ijik}+\nabla_{k}R_{ijij})](y_{0})\Lambda_{k}(y_{0})\phi(y_{0})$

\qquad\ $\mathbf{-}$ $\frac{1}{36}\underset{i,j,k=1}{\overset{n}{\sum}%
}R_{ijkj}(y_{0})\Omega_{ik}(y_{0})(y_{0})\phi(y_{0})\qquad\qquad$I$_{326223}$

\qquad$+\ \frac{1}{24}\underset{i=1}{\overset{n}{\sum}}\frac{\partial
^{2}\text{W}}{\partial x_{i}^{2}}(y_{0})\phi(y_{0})\qquad\qquad$%
\ \textbf{I}$_{327}$

$\qquad+\frac{1}{144}\underset{i,j=1}{\overset{n}{\sum}}[\nabla_{i}%
\varrho_{ij}](y_{0})\Lambda_{j}(y_{0})\phi(y_{0})$\qquad\qquad\qquad
\qquad\qquad\ \ 

$\qquad+\frac{1}{144}\underset{i,j=1}{\overset{n}{\sum}}[\nabla_{j}%
\varrho_{ii}](y_{0})\Lambda_{j}(y_{0})\phi(y_{0})$

$\qquad+\frac{1}{144}\underset{i,j=1}{\overset{n}{\sum}}[\nabla_{i}%
\varrho_{ij}](y_{0})\Lambda_{j}(y_{0})\phi(y_{0})\qquad$

$\qquad+$ $\frac{1}{72}\underset{i,j=1}{\overset{n}{\sum}}\varrho_{ij}%
(y_{0})\Omega_{ij}(y_{0})\phi(y_{0})\qquad\qquad$I$_{32913}$

\qquad$-\frac{1}{36}[\left(  2\frac{\partial^{2}X_{i}}{\partial x_{i}\partial
x_{j}}+\frac{\partial^{2}X_{j}}{\partial x_{i}^{2}}\right)
+2\underset{k=1}{\overset{n}{\sum}}R_{ijik}X_{k}](y_{0})\Lambda_{j}(y_{0}%
)\phi(y_{0})\qquad$I$_{32922}$

$\qquad-\frac{1}{36}X_{j}(y_{0})$ $\frac{\partial\Omega_{ij}}{\partial x_{i}%
}(y_{0})\phi(y_{0})-\frac{1}{12}\frac{\partial X_{j}}{\partial x_{i}}%
(y_{0})\Omega_{ij}(y_{0})$ $\phi(y_{0})$\qquad\qquad

\qquad$+\frac{1}{12}\underset{j=1}{\overset{n}{\sum}}\frac{\partial^{2}X_{j}%
}{\partial x_{i}^{2}}(y_{0})\Lambda_{j}(y_{0})\phi(y_{0})\qquad$%
\textbf{L}$_{2}\qquad$\textbf{L}$_{21}$

\qquad$+\frac{1}{36}$ $\underset{j=1}{\overset{n}{\sum}}$ $X_{j}(y_{0}%
)\frac{\partial\Omega_{ij}}{\partial x_{i}}(y_{0})\phi(y_{0})\qquad$%
\textbf{L}$_{22}$

\qquad$+\frac{1}{12}$ $\underset{j=1}{\overset{n}{\sum}}$ $\frac{\partial
X_{j}}{\partial x_{i}}(y_{0})\Omega_{ij}(y_{0})\phi(y_{0})\qquad$%
\textbf{L}$_{23}$\qquad

$\qquad+\frac{1}{96}[2(\tau^{M})](y_{0})$W$(y_{0})\phi\left(  y_{0}\right)
\qquad$I$_{35}$

\qquad$-\frac{1}{8}[$ $\left\Vert \text{X}\right\Vert _{M}^{2}+$
$\operatorname{div}X_{M}](y_{0})$W(y$_{0}$)$\phi\left(  y_{0}\right)  $

$\qquad+\frac{1}{8}$W$^{2}\left(  y_{0}\right)  \phi(y_{0})+\frac{1}{4}$
V$(y_{0})$W$(y_{0})\phi\left(  y_{0}\right)  $

\qquad$+\frac{1}{48}$ $\underset{j=1}{\overset{n}{\sum}}X_{j}(y_{0}%
)\Lambda_{j}(y_{0})[2(\tau^{M})](y_{0})\phi\left(  y_{0}\right)  $

\qquad$-\frac{1}{4}\underset{j=1}{\overset{n}{\sum}}X_{j}(y_{0})\Lambda
_{j}(y_{0})[$ $\left\Vert \text{X}\right\Vert _{M}^{2}+\frac{1}{2}$
$\operatorname{div}X_{M}](y_{0})\phi\left(  y_{0}\right)  $

$\qquad+\frac{1}{4}$ $\underset{j=1}{\overset{n}{\sum}}X_{j}(y_{0})\Lambda
_{j}(y_{0})$W$\left(  y_{0}\right)  \phi\left(  y_{0}\right)  $

$\qquad+$ $\frac{1}{2}$ $\underset{j=1}{\overset{n}{\sum}}X_{j}(y_{0}%
)\Lambda_{j}(y_{0})$V$(y_{0})\phi\left(  y_{0}\right)  $

\begin{proof}
We \textbf{delete all submanifold-related items} from the Theorem and have the
above expression
\end{proof}

We then carry out \textbf{simplifications}:

$\left(  11.43\right)  $ \qquad\underline{\textbf{ Simplifications of the
above Expression}}

Before we carry out the above simplifications, we note the following:

It is immediate that,

$\qquad\qquad\overset{n}{\underset{i,j=1}{\sum}}\nabla_{ii}^{2}\varrho
_{jj}=\Delta\tau\qquad$

By \textbf{Mckean and Singer }$\left[  1\right]  ,$ p.65,

\qquad\qquad$\underset{i,j=1}{\overset{n}{\sum}}$ $\nabla_{ij}^{2}\varrho
_{ij}=\frac{1}{2}\Delta\tau\qquad$

By $\left(  9.27\right)  $ of \textbf{Lemma }$9.7$ of\textbf{ Gray }$\left[
4\right]  ,$

$\qquad\qquad\underset{i=1}{\overset{n}{\sum}}$ $\nabla_{i}\varrho_{ij}%
=\frac{1}{2}\nabla_{j}\tau$

Since the Ricci curvature $\varrho_{ij}$ is symmetric in the indices $i,j,$ we have,

\qquad$\qquad\underset{i,j=1}{\overset{n}{\sum}}\nabla_{ij}^{2}\varrho
_{ji}=\underset{i,j=1}{\overset{n}{\sum}}\nabla_{ij}^{2}\varrho_{ij}=\frac
{1}{2}\Delta\tau$

By definition,

\qquad\qquad\ $\underset{i,j=1}{\overset{n}{\sum}}\varrho_{ij}^{2}=\left\Vert
\varrho^{M}\right\Vert ^{2}$and $\overset{n}{\underset{i,j,k,l=1}{\sum}%
}R_{ijkl}^{2}=\left\Vert R\right\Vert ^{2}$

By \textbf{Lemma (9.11)} of Gilkey $\left[  2\right]  :$

$\qquad\underset{i,j,k,l=1}{\overset{n}{\sum}}R_{ijkl}R_{kjil}=\frac{1}%
{2}\underset{i,j,k,l=1}{\overset{n}{\sum}}R_{ijkl}R_{ijkl}=\frac{1}%
{2}\underset{i,j,k,l=1}{\overset{n}{\sum}}R_{ijkl}^{2}=\frac{1}{2}\left\Vert
R\right\Vert ^{2}$

We finally note that the \textbf{divergence} of a vector field X is given in
\textbf{Fermi coordinates} by$\left(  B_{22}\right)  :$\qquad$\qquad\qquad$

$\qquad\operatorname{div}_{M}X(y_{0})$ $=-\underset{j=q+1}{\overset{n}{\sum}%
}<H,j>(y_{0})X_{j}(y_{0})+$ $\underset{j=1}{\overset{n}{\sum}}\frac{\partial
X_{j}}{\partial x_{j}}(y_{0})$

When the Fermi coordinates reduce to \textbf{normal coordinates}, then $H=0$
and we have, as is well known:

$\qquad\operatorname{div}_{M}X(y_{0})$ $=$ $\underset{j=1}{\overset{n}{\sum}%
}\frac{\partial X_{j}}{\partial x_{j}}(y_{0})$

With summation over repeated indices understood,

$\qquad X_{i}^{2}=\left\Vert X\right\Vert _{M}^{2}$

$\qquad\Lambda_{j}(y_{0})=0$ for $j=1,...,q,q+1,...,n$

\begin{center}
\qquad\qquad\qquad\qquad\qquad\qquad\qquad\qquad\qquad\qquad$\blacksquare$
\end{center}

We use the above relations and further simplify the expression for
b$_{2}(y_{0},y_{0})\phi\left(  y_{0}\right)  $ in the last Corollary above and
have:$\qquad$

\begin{corollary}
b$_{2}($y$_{0},$y$_{0})\phi($y$_{0})=$ I$_{1}+$ I$_{31}+$ I$_{32}+$ I$_{33}+$
I$_{34}+$ I$_{35}+$ I$_{36}+$ I$_{37}$
\end{corollary}

$=\frac{1}{288}[\tau^{M}-6\left\Vert X\right\Vert ^{2}+6\operatorname{div}$X
$-$ $12$\textbf{V}$]^{2}(y_{0})\phi\left(  y_{0}\right)  \qquad$I$_{1}%
\qquad\qquad$

$\ +\frac{1}{24}[\tau^{M}-6$ $\left\Vert X\right\Vert _{M}^{2}%
-6\operatorname{div}X_{M}+12$V$](y_{0})$W$(y_{0})\phi(y_{0})\qquad\qquad$

$-\frac{1}{864}(\tau^{M})^{2}(y_{0})\phi(y_{0})\qquad\qquad$J$_{1}$

$+\frac{1}{216}\left\Vert \varrho^{M}\right\Vert ^{2}\phi(y_{0})-\frac{1}%
{108}\left\Vert \varrho^{M}\right\Vert ^{2}\phi(y_{0})+\frac{1}{144}\left\Vert
\varrho^{M}\right\Vert ^{2}(y_{0})\phi(y_{0})+\frac{1}{288}(\tau^{M}%
)^{2}(y_{0})\phi(y_{0})\qquad I_{32}\qquad I_{321}\ \qquad\ \ \ $

$-\frac{1}{288}[-$\ $2\Delta\tau^{M}+\frac{4}{5}\Delta\tau^{M}-2\Delta\tau
^{M}+\frac{4}{5}\Delta\tau^{M}+\frac{2}{5}\varrho^{2}-\frac{2}{5}\left\Vert
R^{M}\right\Vert ^{2}\ +\frac{2}{3}\tau^{2}+\frac{4}{3}\varrho^{2}-\frac{2}%
{3}\varrho^{2}]\phi(y_{0})$

$=[+\frac{1}{144}\Delta\tau^{M}-\frac{2}{720}\Delta\tau^{M}+\frac{1}%
{144}\Delta\tau^{M}$

$-\frac{2}{720}\Delta\tau^{M}-\frac{1}{720}\left\Vert \varrho^{M}\right\Vert
^{2}+\frac{1}{720}\left\Vert R^{M}\right\Vert ^{2}-\frac{1}{432}(\tau^{M}%
)^{2}-\frac{2}{432}\left\Vert \varrho^{M}\right\Vert ^{2}+\frac{1}%
{432}\left\Vert \varrho^{M}\right\Vert ^{2}](y_{0})\phi(y_{0})\qquad$

$+\frac{1}{6}X_{i}(y_{0})X_{j}(y_{0})\frac{\partial X_{i}}{\partial x_{j}%
}(y_{0})\phi(y_{0})+$ $\frac{1}{18}X_{i}(y_{0})X_{k}(y_{0})R_{jijk}(y_{0}%
)\phi\left(  y_{0}\right)  \qquad\qquad Q_{2}$

$=+\frac{1}{6}X_{i}(y_{0})X_{j}(y_{0})\frac{\partial X_{i}}{\partial x_{j}%
}(y_{0})\phi(y_{0})+$ $\frac{1}{18}\varrho_{ik}^{M}X_{i}(y_{0})X_{k}%
(y_{0})\phi\left(  y_{0}\right)  $

$=+\frac{1}{6}[X_{i}X_{j}\frac{\partial X_{i}}{\partial x_{j}}](y_{0}%
)\phi(y_{0})+$ $\frac{1}{18}[\varrho_{ij}^{M}X_{i}X_{j}](y_{0})\phi
(y_{0})\qquad$

\qquad$-\frac{1}{12}X_{i}(y_{0})\frac{\partial^{2}X_{i}}{\partial x_{j}^{2}%
}(y_{0})\phi(y_{0})+\frac{1}{18}R_{ijkj}(y_{0})X_{i}(y_{0})X_{k}(y_{0}%
)\phi(y_{0})$

$\qquad=-\frac{1}{12}X_{i}(y_{0})\frac{\partial^{2}X_{i}}{\partial x_{j}^{2}%
}(y_{0})\phi(y_{0})+\frac{1}{18}\varrho_{ik}^{M}(y_{0})X_{i}(y_{0})X_{k}%
(y_{0})\phi(y_{0})$

$\qquad=-\frac{1}{12}[X_{i}\frac{\partial^{2}X_{i}}{\partial x_{j}^{2}}%
](y_{0})\phi(y_{0})+\frac{1}{18}[\varrho_{ij}^{M}X_{i}X_{j}](y_{0})\phi
(y_{0})\qquad\qquad\qquad\qquad\qquad\qquad\qquad\qquad\ \qquad\qquad
\qquad\qquad\qquad\qquad\qquad\qquad\qquad\qquad\qquad\qquad\qquad\qquad
\qquad$

$\qquad+\frac{1}{18}R_{ijik}^{M}(y_{0})[X_{j}X_{k}-\frac{1}{2}(\frac{\partial
X_{j}}{\partial x_{k}}+\frac{\partial X_{k}}{\partial x_{j}})](y_{0}%
)\phi(y_{0})\qquad$I$_{32123}\qquad S_{1}\qquad$

$\qquad=+\frac{1}{18}\varrho_{jk}^{M}(y_{0})[X_{j}X_{k}-\frac{1}{2}%
(\frac{\partial X_{j}}{\partial x_{k}}+\frac{\partial X_{k}}{\partial x_{j}%
})](y_{0})\phi(y_{0})\qquad$

$\qquad=+\frac{1}{18}\varrho_{ji}^{M}(y_{0})[X_{j}X_{i}-\frac{1}{2}%
(\frac{\partial X_{j}}{\partial x_{i}}+\frac{\partial X_{i}}{\partial x_{j}%
})](y_{0})\phi(y_{0})$

$\qquad=+\frac{1}{18}\varrho_{ij}^{M}(y_{0})[X_{i}X_{j}-\frac{1}{2}%
(\frac{\partial X_{j}}{\partial x_{i}}+\frac{\partial X_{i}}{\partial x_{j}%
})](y_{0})\phi(y_{0})$

$=+\frac{1}{18}[\varrho_{ij}^{M}X_{i}X_{j}](y_{0})\phi(y_{0})-\frac{1}%
{36}[\varrho_{ij}^{M}\frac{\partial X_{j}}{\partial x_{i}}](y_{0})\phi
(y_{0})-\frac{1}{36}[\varrho_{ij}^{M}\frac{\partial X_{i}}{\partial x_{j}%
}](y_{0})\phi(y_{0})\qquad\qquad\qquad\qquad\qquad\qquad\qquad\qquad
\qquad\qquad\qquad\qquad\qquad$\qquad\qquad\qquad\qquad\qquad\qquad
\qquad\qquad\qquad$\qquad\qquad\qquad\qquad\qquad\qquad\qquad\qquad
\qquad\qquad\qquad\qquad\qquad\qquad\qquad\qquad$

$+\frac{1}{24}[-\frac{1}{3}(\nabla_{i}R_{kjij}+\nabla_{j}R_{ijik}+\nabla
_{k}R_{ijij})](y_{0})X_{k}(y_{0})\phi(y_{0})\qquad S_{32}$

$=-\frac{1}{72}[(\nabla_{i}\varrho_{ki}^{M}+\nabla_{j}\varrho_{jk}^{M}%
+\nabla_{k}\tau^{M})](y_{0})X_{k}(y_{0})\phi(y_{0})$

$=-\frac{1}{72}[(\nabla_{i}\varrho_{ki}^{M}+\nabla_{i}\varrho_{ik}^{M}%
+\nabla_{k}\tau^{M})](y_{0})X_{k}(y_{0})\phi(y_{0})$

$=-\frac{1}{72}[(\nabla_{i}\varrho_{ik}^{M}+\nabla_{i}\varrho_{ik}^{M}%
+\nabla_{k}\tau^{M})](y_{0})X_{k}(y_{0})\phi(y_{0})$

$=-\frac{1}{72}[(\nabla_{i}\varrho_{ij}^{M}+\nabla_{i}\varrho_{ij}^{M}%
+\nabla_{j}\tau^{M})](y_{0})X_{j}(y_{0})\phi(y_{0})$

$=-\frac{1}{72}[\frac{1}{2}\nabla_{j}\tau^{M}+\frac{1}{2}\nabla_{j}\tau
^{M}+\nabla_{j}\tau^{M})](y_{0})X_{j}(y_{0})\phi(y_{0})$

$=-\frac{1}{36}[\nabla_{j}\tau^{M})](y_{0})X_{j}(y_{0})\phi(y_{0})$

$-\frac{1}{18}R_{ijkj}(y_{0})[X_{i}X_{k}-\frac{1}{2}\left(  \frac{\partial
X_{i}}{\partial x_{k}}+\frac{\partial X_{k}}{\partial x_{i}}\right)
](y_{0})\phi(y_{0})$

$=-\frac{1}{18}\varrho_{ik}^{M}(y_{0})[X_{i}X_{k}-\frac{1}{2}\left(
\frac{\partial X_{i}}{\partial x_{k}}+\frac{\partial X_{k}}{\partial x_{i}%
}\right)  ](y_{0})\phi(y_{0})$

$=-\frac{1}{18}[\varrho_{ij}^{M}X_{i}X_{j}](y_{0})\phi(y_{0})+\frac{1}%
{36}[\varrho_{ij}^{M}\frac{\partial X_{i}}{\partial x_{j}}](y_{0})\phi
(y_{0})+\frac{1}{36}[\varrho_{ij}^{M}\frac{\partial X_{j}}{\partial x_{i}%
}](y_{0})\phi(y_{0})\qquad$

$\ \ -\frac{1}{12}\left(  X_{i}X_{j}\frac{\partial X_{j}}{\partial x_{i}%
}\right)  (y_{0})\phi(y_{0})-\frac{1}{12}\left(  X_{i}X_{j}\frac{\partial
X_{i}}{\partial x_{j}}\right)  (y_{0})\phi(y_{0})$

$+\frac{1}{24}[X_{i}^{2}X_{j}^{2}-2X_{i}X_{j}\left(  \frac{\partial X_{j}%
}{\partial x_{i}}+\frac{\partial X_{i}}{\partial x_{j}}\right)  -X_{i}%
^{2}\frac{\partial X_{j}}{\partial x_{j}}-X_{j}^{2}\frac{\partial X_{i}%
}{\partial x_{i}}](y_{0})\phi(y_{0})$

$=+\frac{1}{24}[\left\Vert \text{X}\right\Vert _{M}^{2}\left\Vert
\text{X}\right\Vert _{M}^{2}-2X_{i}X_{j}\left(  \frac{\partial X_{j}}{\partial
x_{i}}+\frac{\partial X_{i}}{\partial x_{j}}\right)  -\left\Vert
\text{X}\right\Vert _{M}^{2}\operatorname{div}X_{M_{{}}}-\left\Vert
\text{X}\right\Vert _{M}^{2}\operatorname{div}X_{M_{{}}}](y_{0})\phi(y_{0})$

$=+\frac{1}{24}[\left\Vert \text{X}\right\Vert _{M}^{2}\left\Vert
\text{X}\right\Vert _{M}^{2}-2X_{i}X_{j}\left(  \frac{\partial X_{j}}{\partial
x_{i}}+\frac{\partial X_{i}}{\partial x_{j}}\right)  -2\left\Vert
\text{X}\right\Vert _{M}^{2}\operatorname{div}X_{M}](y_{0})\phi(y_{0})$

$=+\frac{1}{24}[\left\Vert \text{X}\right\Vert _{M}^{4}](y_{0})\phi
(y_{0})-\frac{1}{12}X_{i}X_{j}\left(  \frac{\partial X_{j}}{\partial x_{i}%
}+\frac{\partial X_{i}}{\partial x_{j}}\right)  -\frac{1}{12}\left\Vert
\text{X}\right\Vert _{M}^{2}\operatorname{div}X_{M}](y_{0})\phi(y_{0})$

$=+\frac{1}{24}[\left\Vert \text{X}\right\Vert _{M}^{4}](y_{0})\phi
(y_{0})-\frac{1}{12}\left(  X_{i}X_{j}\frac{\partial X_{j}}{\partial x_{i}%
}\right)  -\frac{1}{12}\left(  X_{i}X_{j}\frac{\partial X_{i}}{\partial x_{j}%
}\right)  -\frac{1}{12}\left\Vert \text{X}\right\Vert _{M}^{2}%
\operatorname{div}X_{M}](y_{0})\phi(y_{0})$

$\ +\frac{1}{48}\left(  \frac{\partial X_{j}}{\partial x_{i}}+\frac{\partial
X_{i}}{\partial x_{j}}\right)  ^{2}(y_{0})\phi(y_{0})+\frac{1}{24}\left(
\frac{\partial X_{i}}{\partial x_{i}}\frac{\partial X_{j}}{\partial x_{j}%
}\right)  (y_{0})\phi(y_{0})$

$=+\frac{1}{48}\left(  \frac{\partial X_{j}}{\partial x_{i}}\right)
^{2}(y_{0})\phi(y_{0})+\frac{1}{48}\left(  \frac{\partial X_{i}}{\partial
x_{j}}\right)  ^{2}(y_{0})\phi(y_{0})+\frac{1}{24}\left(  \frac{\partial
X_{j}}{\partial x_{i}}\frac{\partial X_{i}}{\partial x_{j}}\right)
(y_{0})\phi(y_{0})$

$\ +\frac{1}{24}\left(  \frac{\partial X_{i}}{\partial x_{i}}\frac{\partial
X_{j}}{\partial x_{j}}\right)  (y_{0})\phi(y_{0})$

$=+\frac{1}{48}\left(  \frac{\partial X_{j}}{\partial x_{i}}\right)
^{2}(y_{0})\phi(y_{0})+\frac{1}{48}\left(  \frac{\partial X_{i}}{\partial
x_{j}}\right)  ^{2}(y_{0})\phi(y_{0})+\frac{1}{24}\left(  \frac{\partial
X_{j}}{\partial x_{i}}\frac{\partial X_{i}}{\partial x_{j}}\right)
(y_{0})\phi(y_{0})$

$\ +\frac{1}{24}(\operatorname{div}X)^{2}$

$+\frac{1}{36}X_{i}(y_{0})\left(  2\frac{\partial^{2}X_{j}}{\partial
x_{i}\partial x_{j}}+\frac{\partial^{2}X_{i}}{\partial x_{j}^{2}}\right)
(y_{0})\phi(y_{0})+\frac{1}{36}X_{j}(y_{0})\left(  \frac{\partial^{2}X_{j}%
}{\partial x_{i}^{2}}+2\frac{\partial^{2}X_{i}}{\partial x_{i}\partial x_{j}%
}\right)  (y_{0})\phi(y_{0})$

$=+\frac{1}{18}X_{i}(y_{0})\frac{\partial^{2}X_{j}}{\partial x_{i}\partial
x_{j}}+\frac{1}{36}X_{i}(y_{0})\frac{\partial^{2}X_{i}}{\partial x_{j}^{2}%
}(y_{0})\phi(y_{0})+\frac{1}{36}X_{j}\frac{\partial^{2}X_{j}}{\partial
x_{i}^{2}}(y_{0})\phi(y_{0})+\frac{1}{18}X_{j}\frac{\partial^{2}X_{i}%
}{\partial x_{i}\partial x_{j}}(y_{0})\phi(y_{0})$

$\qquad-\frac{1}{48}\left(  \frac{\partial^{3}X_{i}}{\partial x_{i}\partial
x_{j}^{2}}+\frac{\partial^{3}X_{j}}{\partial x_{i}^{2}\partial x_{j}}\right)
(y_{0})\phi(y_{0})$

\qquad$+\frac{1}{12}[(\frac{1}{6}(2\varrho_{ij}^{M}))](y_{0})\times\frac{1}%
{2}[\left(  \frac{\partial X_{j}}{\partial x_{i}}-\frac{\partial X_{i}%
}{\partial x_{j}}\right)  ](y_{0})\phi(y_{0})\qquad\qquad$I$_{3213}$

$=+\frac{1}{12}[\frac{1}{6}(\varrho_{ij}^{M})](y_{0})\times\lbrack\left(
\frac{\partial X_{j}}{\partial x_{i}}-\frac{\partial X_{i}}{\partial x_{j}%
}\right)  ](y_{0})\phi(y_{0})$

$=+\frac{1}{72}\left(  \varrho_{ij}^{M}\frac{\partial X_{j}}{\partial x_{i}%
}\right)  (y_{0})\phi(y_{0})-\frac{1}{72}\left(  \varrho_{ij}^{M}%
\frac{\partial X_{i}}{\partial x_{j}}\right)  (y_{0})\phi(y_{0})$

$-\frac{1}{18}[X_{j}\left(  2\frac{\partial^{2}X_{j}}{\partial x_{i}^{2}%
}+\frac{\partial^{2}X_{i}}{\partial x_{i}\partial x_{j}}\right)  ](y_{0}%
)\phi(y_{0})-\frac{1}{12}[\left(  \frac{\partial X_{i}}{\partial x_{j}}%
+\frac{\partial X_{j}}{\partial x_{i}}\right)  \frac{\partial X_{j}}{\partial
x_{i}}](y_{0})\phi(y_{0})\qquad\ $I$_{3214}\qquad$

$=-\frac{1}{9}\left(  X_{j}\frac{\partial^{2}X_{j}}{\partial x_{i}^{2}%
}\right)  (y_{0})\phi(y_{0})-\frac{1}{18}\left(  X_{j}\frac{\partial^{2}X_{i}%
}{\partial x_{i}\partial x_{j}}\right)  (y_{0})\phi(y_{0})$

$-\frac{1}{12}\left(  \frac{\partial X_{i}}{\partial x_{j}}\frac{\partial
X_{j}}{\partial x_{i}}\right)  (y_{0})\phi(y_{0})-\frac{1}{12}\left(
\frac{\partial X_{j}}{\partial x_{i}}\right)  ^{2}(y_{0})\phi(y_{0})$%
\qquad\qquad\qquad\qquad\qquad\qquad\qquad\qquad\qquad\qquad\qquad\qquad
\qquad\qquad\qquad\qquad

$\qquad+\frac{1}{12}\frac{\partial^{2}\text{V}}{\partial x_{i}^{2}}(y_{0}%
)\phi(y_{0})\qquad\qquad$I$_{3215}$

\qquad$+\frac{1}{72}[R_{ijik}\Omega_{jk}](y_{0})\phi(y_{0})=+\frac{1}%
{72}[\varrho_{jk}^{M}\Omega_{jk}](y_{0})\phi(y_{0})=+\frac{1}{72}[\varrho
_{jk}^{M}\Omega_{jk}](y_{0})\phi(y_{0})\qquad$\ I$_{323}$

\qquad$=+\frac{1}{72}[\varrho_{ij}^{M}\Omega_{ij}](y_{0})\phi(y_{0})$

$\qquad+\frac{1}{48}\underset{i,j=1}{\overset{n}{\sum}}\left(  \Omega
_{ij}\Omega_{ij}\right)  (y_{0})\phi(y_{0})\qquad\qquad\qquad\qquad\qquad
$I$_{3252}$

\qquad\ $+\frac{1}{24}[\frac{1}{3}(\nabla_{i}R_{kjij}+\nabla_{j}%
R_{ijik}+\nabla_{k}R_{ijij})](y_{0})\Lambda_{k}(y_{0})\phi(y_{0})=0$

\qquad\ $\mathbf{-}$ $\frac{1}{36}[R_{ijkj}\Omega_{ik}](y_{0})(y_{0}%
)\phi(y_{0})=\mathbf{-}\frac{1}{36}[\varrho_{ik}^{M}\Omega_{ik}](y_{0}%
)\phi(y_{0})=\mathbf{-}\frac{1}{36}[\varrho_{ij}^{M}\Omega_{ij}](y_{0}%
)\phi(y_{0})\qquad$I$_{326223}$

\qquad$+\ \frac{1}{24}\underset{i=1}{\overset{n}{\sum}}\frac{\partial
^{2}\text{W}}{\partial x_{i}^{2}}(y_{0})\phi(y_{0})\qquad\qquad$%
\ \textbf{I}$_{327}$

$\qquad+\frac{1}{144}\underset{i,j=1}{\overset{n}{\sum}}[\nabla_{i}%
\varrho_{ij}](y_{0})\Lambda_{j}(y_{0})\phi(y_{0})=0$\qquad\qquad\qquad
\qquad\qquad\ \ 

$\qquad+\frac{1}{144}\underset{i,j=1}{\overset{n}{\sum}}[\nabla_{j}%
\varrho_{ii}](y_{0})\Lambda_{j}(y_{0})\phi(y_{0})=0$

$\qquad+\frac{1}{144}\underset{i,j=1}{\overset{n}{\sum}}[\nabla_{i}%
\varrho_{ij}](y_{0})\Lambda_{j}(y_{0})\phi(y_{0})=0\qquad$

$\qquad+$ $\frac{1}{72}\underset{i,j=1}{\overset{n}{\sum}}\varrho_{ij}%
(y_{0})\Omega_{ij}(y_{0})\phi(y_{0})\qquad\qquad$I$_{32913}$

\qquad$-\frac{1}{36}[\left(  2\frac{\partial^{2}X_{i}}{\partial x_{i}\partial
x_{j}}+\frac{\partial^{2}X_{j}}{\partial x_{i}^{2}}\right)
+2\underset{k=1}{\overset{n}{\sum}}R_{ijik}X_{k}](y_{0})\Lambda_{j}(y_{0}%
)\phi(y_{0})=0\qquad$I$_{32922}$

$\qquad-\frac{1}{36}X_{j}(y_{0})$ $\frac{\partial\Omega_{ij}}{\partial x_{i}%
}(y_{0})\phi(y_{0})-\frac{1}{12}\frac{\partial X_{j}}{\partial x_{i}}%
(y_{0})\Omega_{ij}(y_{0})$ $\phi(y_{0})$\qquad\qquad

\qquad$+\frac{1}{12}\underset{j=1}{\overset{n}{\sum}}\frac{\partial^{2}X_{j}%
}{\partial x_{i}^{2}}(y_{0})\Lambda_{j}(y_{0})\phi(y_{0})=0\qquad$%
\textbf{L}$_{2}\qquad$\textbf{L}$_{21}$

\qquad$+\frac{1}{36}$ $\underset{j=1}{\overset{n}{\sum}}$ $X_{j}(y_{0}%
)\frac{\partial\Omega_{ij}}{\partial x_{i}}(y_{0})\phi(y_{0})\qquad$%
\textbf{L}$_{22}$

\qquad$+\frac{1}{12}$ $\underset{j=1}{\overset{n}{\sum}}$ $\frac{\partial
X_{j}}{\partial x_{i}}(y_{0})\Omega_{ij}(y_{0})\phi(y_{0})\qquad$%
\textbf{L}$_{23}$\qquad

$\qquad+\frac{1}{96}[2(\tau^{M})](y_{0})$W$(y_{0})\phi\left(  y_{0}\right)
=+\frac{1}{48}[\tau^{M}$W$](y_{0})\phi\left(  y_{0}\right)  \qquad$I$_{35}$

\qquad$-\frac{1}{8}[$ $\left\Vert \text{X}\right\Vert _{M}^{2}+$
$\operatorname{div}X_{M}]\left(  y_{0}\right)  $W(y$_{0}$)$\phi\left(
y_{0}\right)  $

$\qquad+\frac{1}{8}$W$^{2}\left(  y_{0}\right)  \phi(y_{0})+\frac{1}{4}$
V$(y_{0})$W$(y_{0})\phi\left(  y_{0}\right)  $

\qquad$+\frac{1}{48}$ $\underset{j=1}{\overset{n}{\sum}}X_{j}(y_{0}%
)\Lambda_{j}(y_{0})[2(\tau^{M})](y_{0})\phi\left(  y_{0}\right)  =0$

\qquad$-\frac{1}{4}\underset{j=1}{\overset{n}{\sum}}X_{j}(y_{0})\Lambda
_{j}(y_{0})[$ $\left\Vert \text{X}\right\Vert _{M}^{2}+\frac{1}{2}$
$\operatorname{div}X_{M}](y_{0})\phi\left(  y_{0}\right)  =0$

$\qquad+\frac{1}{4}$ $\underset{j=1}{\overset{n}{\sum}}X_{j}(y_{0})\Lambda
_{j}(y_{0})$W$\left(  y_{0}\right)  \phi\left(  y_{0}\right)  =0$

$\qquad+$ $\frac{1}{2}$ $\underset{j=1}{\overset{n}{\sum}}X_{j}(y_{0}%
)\Lambda_{j}(y_{0})$V$(y_{0})\phi\left(  y_{0}\right)  =0$

\qquad\qquad\qquad\qquad\qquad\qquad\qquad\qquad\qquad\qquad\qquad\qquad
\qquad\qquad\qquad\qquad\qquad$\blacksquare$

We add similar terms to arrive at the final expression for b$_{2}(y_{0}%
,y_{0})\phi\left(  y_{0}\right)  $.

In order to do so, we first make a Table of the Fractions involved in the
expression above:

\begin{summary}
Table of Fractions
\end{summary}

$\qquad\qquad\ \ \ \ [\tau^{M}-6\left(  \left\Vert \text{X}\right\Vert
_{M}^{2}+\operatorname{div}X_{M}\right)  +12V$ $]^{2}(y_{0})\phi\left(
y_{0}\right)  :$ $\frac{1}{288}=+\frac{1}{288}$

$\qquad\qquad\ \ \ \ [\tau^{M}-6$ $\left\Vert X\right\Vert _{M}^{2}%
]-6\operatorname{div}X_{M}+12$V$](y_{0})$W$(y_{0}):\frac{1}{48}+\frac{1}%
{48}=+\frac{1}{24}$

\qquad(i)$\qquad\left\Vert R^{M}\right\Vert ^{2}:+\frac{1}{720}=+\frac{1}%
{720}$

$\qquad$(ii)$\qquad(\tau^{M})^{2}:-\frac{1}{864}+\frac{1}{288}-\frac{1}%
{432}=0$

\qquad(iii)$\qquad\left\Vert \varrho^{M}\right\Vert ^{2}:+\frac{1}{216}%
-\frac{1}{108}+\frac{1}{144}-\frac{1}{720}-\frac{2}{432}+\frac{1}{432}%
=-\frac{1}{720}$

\qquad(iv)$\qquad\Delta\tau^{M}:\frac{1}{144}-\frac{2}{720}+\frac{1}%
{144}=+\frac{6}{720}$

\qquad(v)$\qquad X_{i}X_{j}\frac{\partial X_{i}}{\partial x_{j}}:+\frac{1}%
{6}-\frac{1}{12}=\frac{1}{12}=\frac{3}{36}$

\ \ \ \ \ \ \ (vi)$\qquad X_{i}X_{j}\frac{\partial X_{j}}{\partial x_{i}%
}:-\frac{1}{12}=-\frac{1}{12}=--\frac{3}{36}$

\qquad(vii)$\qquad\varrho_{ij}X_{i}X_{j}:\frac{1}{18}+\frac{1}{18}+\frac
{1}{18}-\frac{1}{18}=\frac{1}{9}=\frac{4}{36}$

\ \ \ \ \ (viii)$\qquad\varrho_{ij}\frac{\partial X_{i}}{\partial x_{j}%
}:-\frac{1}{36}+\frac{1}{36}-\frac{1}{72}=-\frac{1}{72}$

\ \ \ \ (ix)$\qquad\varrho_{ij}\frac{\partial X_{j}}{\partial x_{i}}:-\frac
{1}{36}+\frac{1}{36}+\frac{1}{72}=\frac{1}{72}\qquad$

$\ \ \ \ \ \ $(x)$\qquad\ X_{i}\frac{\partial^{2}X_{i}}{\partial x_{j}^{2}%
}:-\frac{1}{12}+\frac{1}{36}=-\frac{2}{36}$

\ \ \ \ \ \ (xi)$\qquad X_{j}\frac{\partial^{2}X_{j}}{\partial x_{i}^{2}%
}:\frac{1}{36}-\frac{1}{9}=-\frac{3}{36}$

\ \ \ \ \ \ (xii)$\qquad X_{j}\frac{\partial^{2}X_{i}}{\partial x_{i}\partial
x_{j}}:\frac{1}{18}-\frac{1}{18}=0$

\ \ \ \ \ (xiii)\qquad$X_{i}\frac{\partial^{2}X_{j}}{\partial x_{i}\partial
x_{j}}:\frac{1}{18}=\frac{1}{18}$

\ \ \ \ \ \ (xiv)$\qquad(\frac{\partial X_{i}}{\partial x_{j}})^{2}:\frac
{1}{48}=\frac{1}{48}$

\ \ \ \ (xv) $\ \ \ \ \ \ \ \left(  \frac{\partial X_{j}}{\partial x_{i}%
}\right)  ^{2}:\frac{1}{48}-\frac{1}{12}=\frac{1-4}{48}=-\frac{3}{48}$

\ \ \ \ \ \ (xvi)$\qquad\frac{\partial X_{i}}{\partial x_{j}}\frac{\partial
X_{j}}{\partial x_{i}}:\frac{1}{24}-\frac{1}{12}=-\frac{1}{24}$

\ \ \ \ (xvii)$\qquad\ \left\Vert X\right\Vert ^{4}:\frac{1}{24}-\frac{1}%
{12}+\frac{1}{24}=0$

\ \ \ \ \ \ (xviii)\qquad$(\operatorname{div}X)^{2}:-\frac{1}{24}+\frac{1}%
{24}=0$

\ \ \ \ \ \ (xix)$\qquad\ [\left\Vert X\right\Vert ^{2}\operatorname{div}%
X:\frac{1}{12}-\frac{1}{12}=0$

\ \ \ \ \ \ (xx) $\ \ \ \ \ \ \ \ \frac{\partial^{3}X_{j}}{\partial x_{i}%
^{2}\partial x_{j}}:-\frac{1}{48}=-\frac{1}{48}$

\ \ \ \ \ (xxi) $\ \ \ \ \ \ \ \frac{\partial^{3}X_{i}}{\partial x_{i}\partial
x_{j}^{2}}:-\frac{1}{48}=-\frac{1}{48}$

\ \ \ \ \ \ (xxii)$\qquad X_{j}$ $\frac{\partial\Omega_{ij}}{\partial x_{i}%
}:-\frac{1}{36}+\frac{1}{36}=0$

\ \ \ \ \ \ (xxiii) $\ \ \ \ \ \frac{\partial X_{j}}{\partial x_{i}}%
\Omega_{ij}:-\frac{1}{12}+\frac{1}{12}=0$

\ \ \ \ \ \ (xxiv)$\qquad\varrho_{ij}^{M}\Omega_{ij}:\frac{1}{72}-\frac{1}%
{36}+\frac{1}{72}=0$

\ \ \ \ \ \ (xxv) $\ \ \ \ \ \ \ \Omega_{ij}\Omega_{ij}:\frac{1}{48}=\frac
{1}{48}$

\qquad(xxvi) $\ \ \ \ \ \left(  \nabla_{j}\tau^{M}\right)  X_{j}:$
\ $-\frac{1}{36}=-\frac{1}{36}$

\qquad\qquad\qquad\qquad\qquad\qquad\qquad\qquad\qquad\qquad\qquad\qquad
\qquad\qquad\qquad$\qquad\blacksquare$\qquad

We give the result based on the \textbf{Table of Fractions} above:

\begin{corollary}
b$_{2}($y$_{0}$,y$_{0},\phi)=$ b$_{2}($y$_{0},$y$_{0})\phi($y$_{0})=$ I$_{1}+$
I$_{32}+$ I$_{35}$
\end{corollary}

$=\frac{1}{288}[\tau^{M}-6\left\Vert X_{M}\right\Vert ^{2}+6\operatorname{div}%
$X$_{M}$ $-$ $12$\textbf{V}$]^{2}(y_{0})\phi\left(  y_{0}\right)  \qquad$Part
I$_{1}+$ Part I$_{35}\qquad\qquad$

$\ +\frac{1}{24}[\tau^{M}-6$ $\left\Vert X_{M}\right\Vert ^{2}%
]-6\operatorname{div}X_{M}+12$V$](y_{0})$W$(y_{0})\phi(y_{0})\qquad$

$+\frac{1}{720}[\left\Vert R^{M}\right\Vert ^{2}-\ \left\Vert \varrho
^{M}\right\Vert ^{2}+\ 6\Delta\tau^{M}](y_{0})\phi(y_{0})\qquad$I$_{32}\qquad
$I$_{321}\qquad$J$_{1}$

$-\frac{1}{36}[\nabla_{j}\tau^{M}](y_{0})X_{j}(y_{0})\phi(y_{0})\qquad
\qquad\qquad\qquad S_{32}$

$+$ $\frac{8}{72}[\varrho_{ij}^{M}X_{i}X_{j}](y_{0})\phi(y_{0})+\frac{1}%
{72}[\varrho_{ij}^{M}\frac{\partial X_{j}}{\partial x_{i}}-\varrho_{ij}%
^{M}\frac{\partial X_{i}}{\partial x_{j}}](y_{0})\phi(y_{0})\qquad$%
I$_{3213}\qquad$I$_{32123}\qquad S_{1}\qquad$

$+\frac{1}{12}[X_{i}X_{j}\frac{\partial X_{i}}{\partial x_{j}}](y_{0}%
)\phi(y_{0})-\frac{1}{12}\left(  X_{i}X_{j}\frac{\partial X_{j}}{\partial
x_{i}}\right)  (y_{0})\phi(y_{0})$\qquad I$_{32122}\qquad Q_{2}\qquad$

$+\frac{2}{36}X_{i}(y_{0})\frac{\partial^{2}X_{j}}{\partial x_{i}\partial
x_{j}}-\frac{2}{36}X_{i}(y_{0})\frac{\partial^{2}X_{i}}{\partial x_{j}^{2}%
}(y_{0})\phi(y_{0})-\frac{3}{36}X_{j}\frac{\partial^{2}X_{j}}{\partial
x_{i}^{2}}(y_{0})\phi(y_{0})\qquad\qquad\qquad\qquad\qquad\qquad\ \qquad
\qquad\qquad\qquad\qquad\qquad\qquad\qquad\qquad\qquad\qquad\qquad\qquad
\qquad\qquad\qquad\qquad\qquad\qquad\qquad\qquad\qquad\qquad\qquad\qquad
\qquad\qquad\qquad\qquad\qquad$\qquad\qquad\qquad\qquad\qquad$\ \ $

$+\frac{1}{48}\left(  \frac{\partial X_{j}}{\partial x_{i}}\right)  ^{2}%
(y_{0})\phi(y_{0})-\frac{3}{48}\left(  \frac{\partial X_{i}}{\partial x_{j}%
}\right)  ^{2}(y_{0})\phi(y_{0})-\frac{2}{48}\left(  \frac{\partial X_{j}%
}{\partial x_{i}}\frac{\partial X_{i}}{\partial x_{j}}\right)  (y_{0}%
)\phi(y_{0})$

$-\frac{1}{48}\left(  \frac{\partial^{3}X_{i}}{\partial x_{i}\partial
x_{j}^{2}}+\frac{\partial^{3}X_{j}}{\partial x_{i}^{2}\partial x_{j}}\right)
(y_{0})\phi(y_{0})$

$+\frac{1}{48}\underset{i,j=1}{\overset{n}{\sum}}\left(  \Omega_{ij}%
\Omega_{ij}\right)  (y_{0})\phi(y_{0})\qquad\qquad$I$_{3252}\qquad$%
\qquad\qquad\qquad\qquad\qquad\qquad\qquad\qquad\qquad\qquad\qquad\qquad
\qquad\qquad\qquad

$+\frac{1}{24}[2\frac{\partial^{2}\text{V}}{\partial x_{i}^{2}}(y_{0}%
)\phi(y_{0})+\ \frac{\partial^{2}\text{W}}{\partial x_{i}^{2}}(y_{0}%
)\phi(y_{0})+3W^{2}\left(  y_{0}\right)  \phi(y_{0})+6VW](y_{0})\phi\left(
y_{0}\right)  \qquad$I$_{3215}\qquad$I$_{327}\qquad$ Part I$_{35}\qquad$

\begin{proof}
We use the last Corollary above and Table of Fractions above
\end{proof}

The \textbf{reduced expression} of the second order term is particularly
fascinating in the sense that it brings out very clearly the roles played by
the underlying geometric objects: the Riemannian manifold M, the vector bundle
E, the Weitzenbock term W, the vector field X and the potential term
V.\qquad\qquad\qquad\qquad\qquad\qquad\qquad\qquad\qquad\qquad\qquad
\qquad\qquad\qquad$\qquad\blacksquare$

No author (to the best of my knowledge) has computed the second order
coefficient of a vector bundle heat kernel in the presence of a \textbf{vector
field} and/or a \textbf{scalar potential} term. Consequently in order to make
\textbf{comparison} with the work of previous authors, we take X $\equiv0;$ V
$\equiv0$ and \textbf{delete} all the terms containing the \textbf{vector
field X} and the \textbf{potential term V} and have the \textbf{Ultimate
Corollary:}

\begin{center}
\qquad\qquad\qquad\qquad\qquad\qquad$\qquad\qquad\blacksquare$
\end{center}

\begin{corollary}
(Reduced Third Coefficient at P = $\left\{  y_{0}\right\}  $)
\end{corollary}

\qquad b$_{2}($y$_{0}$,y$_{0},\phi)=$ b$_{2}($y$_{0},$y$_{0})\phi($y$_{0})$

\qquad$=\frac{1}{1440}[5(\tau^{M})^{2}$\ $-2(\varrho^{M})^{2}+12\Delta\tau
^{M}+2\left\Vert R^{M}\right\Vert ^{2}$\ $+30(\Omega_{ij}\Omega_{ij}%
)+60\frac{\partial^{2}\text{W}}{\partial x_{i}^{2}}$

$\qquad+60\tau^{M}$W $+180$W$^{2}]\left(  y_{0}\right)  \phi(y_{0})$

\begin{proof}
We delete all terms related to the \textbf{vector field} X and the
\textbf{potential vector} V and get the result.
\end{proof}

$\left(  44\right)  \qquad b_{2}(y_{0},y_{0})\phi(y_{0})=\frac{1}{288}%
[\tau^{M}]^{2}(y_{0})\phi\left(  y_{0}\right)  \ +\frac{1}{24}[\tau
^{M}]W(y_{0})\phi(y_{0})\qquad$

$\qquad+\frac{1}{720}[\left\Vert R^{M}\right\Vert ^{2}-\ \left\Vert
\varrho^{M}\right\Vert ^{2}+\ 6\Delta\tau^{M}](y_{0})\phi(y_{0})\qquad
$I$_{3211}=\frac{1}{24}\frac{\partial^{2}}{\partial\text{x}_{i}^{2}}%
(\theta^{\frac{1}{2}}\Delta\theta^{-\frac{1}{2}})\phi(y_{0})\qquad\qquad
\qquad\ \qquad\qquad\qquad\qquad\qquad\qquad\qquad\qquad\qquad\qquad
\qquad\qquad\qquad\qquad\qquad\qquad\qquad\qquad\qquad\qquad\qquad\qquad
\qquad\qquad\qquad\qquad\qquad\qquad\qquad\qquad$\qquad\qquad\qquad
\qquad\qquad$\ \ $

$\qquad+\frac{1}{48}\underset{i,j=1}{\overset{n}{\sum}}\left(  \Omega
_{ij}\Omega_{ij}\right)  (y_{0})\phi(y_{0})+\ \frac{1}{24}\frac{\partial
^{2}\text{W}}{\partial x_{i}^{2}}(y_{0})\phi(y_{0})+\frac{1}{8}$W$^{2}\left(
y_{0}\right)  \phi(y_{0})$

$\qquad=\frac{1}{1440}[5(\tau^{M})^{2}+2\left\Vert R^{M}\right\Vert
^{2}-\ 2\left\Vert \varrho^{M}\right\Vert ^{2}+\ 12\Delta\tau^{M}+60\tau^{M}$W
$+$ $30\Omega_{ij}\Omega_{ij}$

$\qquad+60\frac{\partial^{2}\text{W}}{\partial x_{i}^{2}}+180$W$^{2}]\left(
y_{0}\right)  \phi(y_{0})\qquad$

\qquad\qquad\qquad\qquad\qquad\qquad\qquad\qquad\qquad\qquad\qquad\qquad
\qquad\qquad\qquad\qquad$\blacksquare$

\begin{remark}
Comparison with Previous Results
\end{remark}

(i) b$_{2}(y_{0},y_{0})$ here is $\left[  \text{a}_{2}\right]  $ of $\left(
34\right)  $ in \textbf{Avramidi }$\left[  1\right]  ,$ and is e$_{4}(x,D)$ in
Theorem $\left(  4.1.6\right)  $ of \textbf{Gilkey }$\left[  1\right]  $ and
a$_{4}($F,D) in Theorem $\left(  3.3.1\right)  $ of \textbf{Gilkey} $\left[
2\right]  .$\qquad$\qquad\qquad\qquad\qquad\qquad\qquad\qquad\qquad
\qquad\qquad\qquad\qquad\qquad\qquad\qquad\qquad$

(ii) We have thus recovered theorems of previous authors. Notice that instead
of the factor $\frac{1}{1440}$ appearing here, \textbf{Gilkey} has $\frac
{1}{360}$ in the above references. The extra factors of $\frac{1}{2}$ in
b$_{1}(y_{0},y_{0})\phi(y_{0})$ and $\frac{1}{4}$ in b$_{2}(y_{0},y_{0}%
)\phi(y_{0})$ here are due to the fact that we are dealing with half the
Laplacian $\frac{1}{2}\Delta$ whereas Gilkey is dealing with the full
Laplacian $\Delta.$

\textbf{Mckean} and \textbf{Singer} in $\left[  1\right]  $ were able to
compute \textbf{b}$_{2}($y$_{0},$y$_{0})$ for only $\frac{1}{2}\Delta.$ Here
we have computed it for $\frac{1}{2}\Delta+$ X + V where X is a smooth vector
field and V is a smooth scalar potential term.

\qquad\qquad\qquad\qquad\qquad\qquad\qquad\qquad\qquad\qquad\qquad\qquad
\qquad\qquad\qquad$\qquad\blacksquare$

\section{The case of a gradient vector field:}

\qquad\qquad\qquad\qquad\qquad\qquad\qquad\qquad\qquad\qquad\qquad\qquad
\qquad\qquad\qquad\qquad$\blacksquare$

$\left(  11.45\right)  \qquad$b$_{2}($y$_{0}$,y$_{0},\phi)=$ b$_{2}(y_{0}%
,$y$_{0})\phi($y$_{0})=[$ I$_{1}+$ I$_{32}+$ I$_{35}]\phi\left(  y_{0}\right)
$where,

\qquad I$_{1}=\frac{1}{2}\frac{\text{L}\Psi}{\Psi}(y_{0})\Theta(y_{0});\qquad
$I$_{32}=\frac{1}{12}\underset{i=1}{\overset{n}{\sum}}\frac{\partial^{2}%
\Theta}{\partial x_{i}^{2}}(y_{0});\qquad$I$_{35}=$ $\frac{1}{4}\Theta(y_{0}%
)$W$(y_{0})\qquad$

From $\left(  11.30\right)  $ of Chapter 7, we have in \textbf{Fermi
coordinates}:

\qquad$\frac{\text{L}\Psi}{\Psi}(y_{0})=\frac{1}{24}%
[\underset{i=q+1}{\overset{n}{\sum}}3<H,i>^{2}+2(\tau^{M}-3\tau^{P}%
\ +\overset{q}{\underset{\text{a=1}}{\sum}}\varrho_{\text{aa}}^{M}%
+\overset{q}{\underset{\text{a,b}=1}{\sum}}R_{\text{abab}}^{M})](y_{0})$

\qquad\qquad$-\frac{1}{2}$ $\left\Vert \text{X}\right\Vert _{M}^{2}%
(y_{0})-\frac{1}{2}$ $\operatorname{div}X_{M}(y_{0})+\frac{1}{2}$ $\left\Vert
\text{X}\right\Vert _{P}^{2}(y_{0})$ $+$ $\frac{1}{2}\operatorname{div}%
X_{P}(y_{0})+$ V$(y_{0})$

Consequently in \textbf{normal coordinates}, we have:

\qquad$\frac{\text{L}\Psi}{\Psi}(y_{0})=\frac{1}{12}(\tau^{M})(y_{0})-\frac
{1}{2}$ $\left\Vert \text{X}\right\Vert _{M}^{2}(y_{0})-\frac{1}{2}$
$\operatorname{div}X_{M}(y_{0})+$ V$(y_{0})$

$\left(  11.46\right)  \qquad\frac{\text{L}\Psi}{\Psi}(y_{0})=\frac{1}%
{12}[\tau^{M}-6$ $\left\Vert \text{X}\right\Vert _{M}^{2}-6\operatorname{div}%
X_{M}+12$V$](y_{0})$

Next we have the general formula in \textbf{Fermi coordinates} from $\left(
11.31\right)  :$

\qquad$\Theta(y_{0})\phi\left(  y_{0}\right)  =$ L$_{\Psi}[\phi\circ
\pi_{\text{P}}](y_{0})$

$\qquad\qquad=\frac{1}{24}[\underset{i=q+1}{\overset{n}{\sum}}3<H,i>^{2}%
+2(\tau^{M}-3\tau^{P}\ +\overset{q}{\underset{\text{a=1}}{\sum}}%
\varrho_{\text{aa}}^{M}+\overset{q}{\underset{\text{a,b}=1}{\sum}%
}R_{\text{abab}}^{M})](y_{0})\phi\left(  y_{0}\right)  $

\qquad\qquad$\ -\frac{1}{2}[$ $\left\Vert \text{X}\right\Vert _{M}^{2}-$
$\left\Vert \text{X}\right\Vert _{P}^{2}$ $+$ $\operatorname{div}X_{M}-$
$\operatorname{div}X_{P}](y_{0})\phi\left(  y_{0}\right)  +$ V$(y_{0}%
)\phi\left(  y_{0}\right)  $

$\qquad\qquad+$ $\frac{1}{2}\underset{\text{a=1}}{\overset{\text{q}}{\sum}%
}\frac{\partial^{2}\phi}{\partial\text{x}_{\text{a}}^{2}}(y_{0})$ $+$
$\underset{\text{a=1}}{\overset{\text{q}}{\sum}}\Lambda_{\text{a}}(y_{0}%
)\frac{\partial\phi}{\partial x_{\text{a}}}\left(  y_{0}\right)  \ +\frac
{1}{2}$ $\underset{\text{a=1}}{\overset{\text{q}}{\sum}}\Lambda_{\text{a}%
}(y_{0})\Lambda_{\text{a}}(y_{0})\phi\left(  y_{0}\right)  $

\qquad\qquad$+$ $\underset{\text{a=1}}{\overset{\text{q}}{\sum}}X_{\text{a}%
}(y_{0})\frac{\partial\phi}{\partial\text{x}_{\text{a}}}(y_{0})$ +
$\underset{\text{a=1}}{\overset{\text{q}}{\sum}}$ $X_{\text{a}}(y_{0}%
)\Lambda_{\text{a}}(y_{0})\phi(y_{0})+\frac{1}{2}$W$\left(  y_{0}\right)
\phi\left(  y_{0}\right)  $

Consequently in \textbf{normal coordinates} we have:

$\qquad\qquad\ \ \ \ \Theta(y_{0})\phi\left(  y_{0}\right)  =\frac{1}{12}%
(\tau^{M}(y_{0})\phi\left(  y_{0}\right)  -\frac{1}{2}[$ $\left\Vert
\text{X}\right\Vert _{M}^{2}+\operatorname{div}X_{M}](y_{0})\phi\left(
y_{0}\right)  $

$\qquad\qquad\qquad\qquad\qquad\qquad+$ V$(y_{0})\phi\left(  y_{0}\right)
+\frac{1}{2}$W$\left(  y_{0}\right)  \phi\left(  y_{0}\right)  $

We have:

$\left(  11.47\right)  \qquad\Theta(y_{0})\phi\left(  y_{0}\right)  =\frac
{1}{12}[\tau^{M}-6$ $\left\Vert \text{X}\right\Vert _{M}^{2}-6$
$\operatorname{div}X_{M}+12V](y_{0})\phi\left(  y_{0}\right)  $

$\qquad\qquad\qquad\qquad\qquad\qquad+\frac{1}{2}$W$\left(  y_{0}\right)
\phi\left(  y_{0}\right)  $

We conclude that in \textbf{normal coordinates,} we have from $\left(
11.46\right)  $ and $\left(  11.47\right)  :$

\bigskip$\qquad\qquad\qquad$I$_{1}=\frac{1}{2}\frac{\text{L}\Psi}{\Psi}%
(y_{0})\Theta(y_{0})\phi\left(  y_{0}\right)  $

$\qquad\qquad\qquad=\frac{1}{2}\times\frac{1}{12}[\tau^{M}-6$ $\left\Vert
\text{X}\right\Vert _{M}^{2}-6$ $\operatorname{div}X_{M}+$ $12$V$](y_{0})$

$\qquad\qquad\qquad\qquad\times\lbrack\frac{1}{12}[\tau^{M}-6[\left\Vert
\text{X}\right\Vert _{M}^{2}+6\operatorname{div}X_{M}+12$V $+$ $6$W$]\left(
y_{0}\right)  \phi\left(  y_{0}\right)  ]$

We have:

$\left(  11.48\right)  \qquad$I$_{1}=\frac{1}{2}\frac{\text{L}\Psi}{\Psi
}(y_{0})\Theta(y_{0})\phi\left(  y_{0}\right)  $

$\qquad\qquad\qquad=\frac{1}{288}[\tau^{M}-6$ $\left\Vert \text{X}\right\Vert
_{M}^{2}-6$ $\operatorname{div}X_{M}+$ $12$V$]^{2}(y_{0})\phi\left(
y_{0}\right)  $

$\qquad\qquad\qquad\ \ +\frac{1}{48}[\tau^{M}-6$ $\left\Vert \text{X}%
\right\Vert _{M}^{2}-6$ $\operatorname{div}X_{M}+$ $6$V$](y_{0})$W$\left(
y_{0}\right)  \phi\left(  y_{0}\right)  $

It is immediate from $\left(  8.11\right)  $ that:

$\left(  11.49\right)  \qquad$I$_{35}=$ $\frac{1}{4}\Theta(y_{0})$%
W$(y_{0})\phi\left(  y_{0}\right)  $

$\qquad\qquad\qquad=\frac{1}{48}[(\tau^{M}-6\left\Vert \text{X}\right\Vert
_{M}^{2}-6$ $\operatorname{div}X_{M}+12$V $+$ $6$W$](y_{0})$W$(y_{0}%
)\phi\left(  y_{0}\right)  $

We next compute the much longer expression:

$\left(  11.50\right)  \qquad\ \ \ $I$_{32}=\frac{1}{12}%
\underset{k=1}{\overset{n}{\sum}}\frac{\partial^{2}\Theta}{\partial x_{k}^{2}%
}(y_{0})=$ I$_{321}+$ I$_{323}+$ I$_{325}+$ I$_{326}+$ I$_{327}+$ I$_{329}+$
\textbf{L}$_{2}$

where the labelling is taken from Appendix D.

\qquad\textbf{I}$_{321}=\frac{1}{12}\underset{i=q+1}{\overset{n}{\sum}}%
\frac{\partial^{2}}{\partial x_{i}^{2}}[\frac{\text{L}\Psi}{\Psi}\phi\circ
\pi_{\text{P}}](y_{0})$

\qquad\textbf{I}$_{323}=\frac{1}{24}\underset{i=q+1}{\overset{n}{\sum}}%
\frac{\partial^{2}}{\partial x_{i}^{2}}[\underset{j,k=1}{\overset{n}{\sum}}%
$g$^{jk}\left\{  \frac{\partial\Lambda_{k}}{\partial\text{x}_{j}}\phi\circ
\pi_{\text{P}}\right\}  ](y_{0})$

\qquad\textbf{I}$_{325}=\frac{1}{24}\underset{i=q+1}{\overset{n}{\sum}}%
\frac{\partial^{2}}{\partial x_{i}^{2}}[\underset{i,j=1}{\overset{n}{\sum}}%
$g$^{jk}\Lambda_{j}\Lambda_{k}\phi\circ\pi_{\text{P}}](y_{0})$

\qquad\textbf{I}$_{326}=-\frac{1}{24}\underset{i=q+1}{\overset{n}{\sum}}%
\frac{\partial^{2}}{\partial x_{i}^{2}}[\underset{j,k=1}{\overset{n}{\sum}}%
$g$^{jk}\left\{  \text{ }\underset{l=1}{\overset{n}{\sum}}\Gamma_{jk}%
^{l}\Lambda_{l}\phi\circ\pi_{\text{P}}\right\}  ](y_{0})$

\qquad\textbf{I}$_{327}=\frac{1}{24}\underset{i=q+1}{\overset{n}{\sum}}%
\frac{\partial^{2}}{\partial x_{i}^{2}}[$W$\phi\circ\pi_{\text{P}}](y_{0})$

\qquad\textbf{I}$_{329}=\frac{1}{12}\underset{i=q+1}{\overset{n}{\sum}}$
$\underset{j=1}{\overset{n}{\sum}}$ $\frac{\partial^{2}}{\partial x_{i}^{2}%
}[(\nabla\log\Psi)_{j}\Lambda_{j}\phi\circ\pi_{\text{P}}](y_{0})$

\qquad\textbf{L}$_{2}=\frac{1}{12}\underset{i=q+1}{\overset{n}{\sum}}$
$\underset{j=1}{\overset{n}{\sum}}$ $\frac{\partial^{2}}{\partial x_{i}^{2}}%
[$X$_{j}\Lambda_{j}\phi\circ\pi_{\text{P}}](y_{0})$

\underline{\textbf{Computations}}

$\left(  11.51\right)  \qquad$I$_{321}$ $=\frac{1}{12}%
\underset{k=1}{\overset{n}{\sum}}\frac{\partial^{2}}{\partial x_{k}^{2}}%
[\frac{\text{L}\Psi}{\Psi}](y_{0}).\phi(y_{0})$

Recall that the differential operator L is given by: L = $\frac{1}{2}%
\Delta+X+V=\frac{1}{2}\Delta+\nabla_{X}+V$

Next recall that for any two smooth functions: $\Phi,f:$M$\longrightarrow R,$
we have:

\qquad$\qquad L(\Phi f)=(L\Phi)f+\Phi(Lf)+<\nabla\Phi,\nabla f>-V(\Phi f)$

Therefore for $\Psi=\Phi\theta^{-\frac{1}{2}},$ we have:

$\qquad\qquad\frac{\text{L}\Psi}{\Psi}=\frac{\text{L}(\Phi\theta^{-\frac{1}%
{2}})}{\Phi\theta^{-\frac{1}{2}}}=\frac{(L\Phi)}{\Phi}+\frac{(L\theta
^{-\frac{1}{2}})}{\theta^{-\frac{1}{2}}}+\frac{1}{\Phi\theta^{-\frac{1}{2}}%
}<\nabla\theta^{-\frac{1}{2}},\bigtriangledown\Phi>-\frac{V(\Phi\theta
^{-\frac{1}{2}})}{\Phi\theta^{-\frac{1}{2}}}$

\qquad\qquad\qquad\qquad\qquad\qquad$=\frac{(L\Phi)}{\Phi}+\frac
{(L\theta^{-\frac{1}{2}})}{\theta^{-\frac{1}{2}}}+$ $<\nabla\log\theta
^{-\frac{1}{2}},\triangledown\log\Phi>-$ $V$

\qquad\qquad\qquad$=\frac{1}{2}\frac{\Delta\Phi}{\Phi}+\frac{1}{2}\frac
{\Delta\theta^{-\frac{1}{2}}}{\theta^{-\frac{1}{2}}}+\frac{1}{\Phi
}<\bigtriangledown\Phi,X>+\frac{1}{\theta^{-\frac{1}{2}}}<\nabla\theta
^{-\frac{1}{2}},X>$

$\qquad\qquad\qquad+$ $<\nabla\log\theta^{-\frac{1}{2}},\Phi\log\nabla>+V+V-V$

We have finally here:

$\left(  11.52\right)  $\qquad\qquad$\frac{\text{L}\Psi}{\Psi}=\frac{1}{2}%
\Phi^{-1}\triangle\Phi+\frac{1}{2}\theta^{\frac{1}{2}}\bigtriangleup
\theta^{-\frac{1}{2}}+$ $<\bigtriangledown\log\Phi,X>-\frac{1}{2}<\nabla
\log\theta,X>$

$\qquad\qquad\qquad\qquad\qquad-\frac{1}{2}$ $<\nabla\log\theta,\triangledown
\log\Phi>+V$

Since $X$ is a gradient vector field, by (i) of Table B$_{2},$ then for any
point $x_{0}$ in the normal neighborhood M$_{0},$ we have:

$\left(  11.53\right)  \qquad\qquad\qquad\qquad\qquad\qquad\qquad
\qquad\triangledown\log\Phi(x_{0})$ $=-X(x_{0})$

Consequently in M$_{0},$

$\qquad\qquad\qquad\qquad\frac{\text{L}\Psi}{\Psi}=\frac{1}{2}\Phi
^{-1}\triangle\Phi+\frac{1}{2}\theta^{\frac{1}{2}}\bigtriangleup\theta
^{-\frac{1}{2}}-$ $<X,X>-\frac{1}{2}<\nabla\log\theta,X>$

$\qquad\qquad\qquad\qquad\qquad+\frac{1}{2}$ $<\nabla\log\theta,X>+$ $V$

$\left(  11.54\right)  $\qquad\qquad\qquad$\frac{\text{L}\Psi}{\Psi}=\frac
{1}{2}\Phi^{-1}\triangle\Phi+\frac{1}{2}\theta^{\frac{1}{2}}\bigtriangleup
\theta^{-\frac{1}{2}}-\left\Vert X\right\Vert ^{2}+V$

We note here that the Laplacian $\triangle$ here is the Laplacian on functions
and hence the Weitzenbock term W is absent.

Since the Christoffel symbols $\Gamma_{ij}^{k}(y_{0})=0$ and $g^{ij}%
(y_{0})=\delta^{ij}$ for $i,j,k=1,...,n,$ where $y_{0}$ is the centre of
\textbf{normal neighborhood}, we have from $\left(  8.17\right)  :$

$\left(  11.55\right)  \qquad\qquad$I$_{321}$ $=\frac{1}{12}%
\underset{k=1}{\overset{n}{\sum}}\frac{\partial^{2}}{\partial x_{k}^{2}}%
[\frac{\text{L}\Psi}{\Psi}](y_{0}).\phi(y_{0})=\frac{1}{12}\Delta\lbrack
\frac{\text{L}\Psi}{\Psi}](y_{0}).\phi(y_{0})$

\qquad$\qquad\frac{1}{12}\Delta\lbrack\frac{\text{L}\Psi}{\Psi}](y_{0}%
)=\frac{1}{24}\Delta\lbrack\Phi^{-1}\triangle\Phi](y_{0})+\frac{1}{24}%
\Delta\lbrack\theta^{\frac{1}{2}}\bigtriangleup\theta^{-\frac{1}{2}}](y_{0})$

$\qquad\qquad\qquad\qquad\qquad-\frac{1}{12}[\Delta\left\Vert X\right\Vert
^{2}](y_{0})+\frac{1}{12}[\Delta V](y_{0})=L_{1}+L_{2}+L_{3}$

where,

\qquad$L_{1}=\frac{1}{24}\Delta\lbrack\Phi^{-1}\triangle\Phi](y_{0});\qquad
L_{2}=\frac{1}{24}\Delta\lbrack\theta^{\frac{1}{2}}\bigtriangleup
\theta^{-\frac{1}{2}}](y_{0});\qquad$

$L_{3}=\frac{1}{12}\left(  -[\Delta\left\Vert X\right\Vert ^{2}]+[\Delta
V]\right)  (y_{0})$

\qquad We compute each of the terms in the last two lines above:

\qquad$L_{1}=\frac{1}{24}\Delta\lbrack\Phi^{-1}\triangle\Phi](y_{0})=\frac
{1}{24}\Delta\lbrack\Phi^{-1}](y_{0})[\triangle\Phi](y_{0})+\frac{1}{24}%
\Phi^{-1}(y_{0})[\triangle\triangle\Phi](y_{0})$

$+\frac{1}{12}$ $<\nabla\Phi^{-1},\nabla\triangle\Phi>(y_{0})$

Since $\Phi^{-1}(y_{0})=1$ and $\nabla\Phi^{-1}(y_{0})=-\Phi^{-1}%
(y_{0})[\nabla\log\nabla\Phi](y_{0})=X(y_{0}),$we have:

$\left(  11.56\right)  \qquad\qquad L_{1}=\frac{1}{24}\Delta\lbrack\Phi
^{-1}](y_{0})[\triangle\Phi](y_{0})+\frac{1}{24}[\triangle^{2}\Phi
](y_{0})+\frac{1}{12}$ $<X,\nabla\triangle\Phi>(y_{0})$

\qquad All terms in the last expression above have been computed except:

\qquad\qquad\ $\Delta\lbrack\Phi^{-1}](y_{0})$ and $[\nabla\triangle
\Phi](y_{0})$

\qquad We use a trick to compute $\Delta\lbrack\Phi^{-1}](y_{0}):$

Since $\Phi(y_{0})=1=\Phi^{-1}(y_{0})$ and $\nabla\log\Phi(y_{0})=-X(y_{0}),$
we have:

$\qquad\qquad0=\Delta\lbrack\Phi^{-1}\Phi](y_{0})=\Delta\lbrack\Phi
^{-1}](y_{0})\Phi(y_{0})+\Phi^{-1}(y_{0})\Delta\Phi+2$ $<\nabla\Phi
^{-1},\nabla\Phi>(y_{0})$

\qquad\qquad$\ =\Delta\lbrack\Phi^{-1}](y_{0})+\Delta\Phi-2$ $<\nabla\log
\Phi,\nabla\log\Phi>(y_{0})$

\qquad We thus have:

$\left(  11.57\right)  $\qquad$\qquad\Delta\lbrack\Phi^{-1}](y_{0}%
)=-\Delta\Phi+2\left\Vert X\right\Vert ^{2}$

\qquad\qquad\qquad We insert the expression on the Right Hand Side of $\left(
11.56\right)  $ and get:

$\left(  11.58\right)  \qquad L_{1}=\frac{1}{24}[-\Delta\Phi+2\left\Vert
X\right\Vert ^{2}](y_{0})[\triangle\Phi](y_{0})+\frac{1}{24}[\triangle^{2}%
\Phi](y_{0})+\frac{1}{12}$ $<X,\nabla\triangle\Phi>(y_{0})$

\qquad We next compute $[\nabla\triangle\Phi](y_{0})$, where,

$\left(  11.59\right)  $\qquad$\Delta\Phi_{P}(x_{0})=\Phi_{P}(x_{0})\left(
\left\Vert \text{X}\right\Vert _{M}^{2}-\operatorname{div}X\right)  (x_{0})$

It is therefore immediate from the expression for in $\left(  11.59\right)  $ that:

\qquad$\qquad\lbrack\nabla\Delta\Phi_{P}](x_{0})=\nabla\Phi_{P}(x_{0})\left(
\left\Vert \text{X}\right\Vert _{M}^{2}-\operatorname{div}X\right)
(x_{0})+\Phi_{P}(x_{0})\left(  \nabla\left\Vert \text{X}\right\Vert _{M}%
^{2}-\nabla\operatorname{div}X\right)  (x_{0})$

Since $\nabla\Phi_{P}(x_{0})=\Phi(x_{0})\nabla\log\Phi_{P}(x_{0}),$ we have:

\qquad\qquad$\lbrack\nabla\Delta\Phi_{P}](x_{0})=\Phi(x_{0})\nabla\log\Phi
_{P}(x_{0})\left(  \left\Vert \text{X}\right\Vert _{M}^{2}-\operatorname{div}%
X\right)  (x_{0})+\left(  \nabla\left\Vert \text{X}\right\Vert _{M}^{2}%
-\nabla\operatorname{div}X\right)  (x_{0})$

We have: $\Phi(y_{0})=1;$ $\nabla\log\Phi_{P}(x_{0})=-X(x_{0})$ and
$\Delta\Phi_{P}(y_{0})=\left(  \left\Vert \text{X}\right\Vert _{M}%
^{2}-\operatorname{div}X\right)  (y_{0});$

Consequently, we have:

$\left(  11.60\right)  $ $\qquad\frac{1}{24}[-\Delta\Phi+2\left\Vert
X\right\Vert ^{2}](y_{0})[\triangle\Phi](y_{0})=\frac{1}{24}[\left(
\left\Vert \text{X}\right\Vert _{M}^{2}+\operatorname{div}X\right)
(y_{0})\left(  \left\Vert \text{X}\right\Vert _{M}^{2}-\operatorname{div}%
X\right)  (y_{0})$

\qquad\qquad\qquad$=\frac{1}{24}\left(  \left\Vert \text{X}\right\Vert
_{M}^{4}-(\operatorname{div}X)^{2}\right)  (y_{0})$\qquad

$\left(  11.61\right)  \qquad\qquad\lbrack\nabla\Delta\Phi_{P}](y_{0}%
)=-X(y_{0})\left(  \left\Vert X\right\Vert _{M}^{2}-\operatorname{div}%
X\right)  (y_{0})+\left(  \nabla\left\Vert X\right\Vert _{M}^{2}%
-\nabla\operatorname{div}X\right)  (y_{0})$

By (vi) of \textbf{Appendix B}$_{3},$ we have:

$\left(  11.62\right)  \qquad\Delta^{2}\Phi(y_{0})=$ $\left(  \left\Vert
\text{X}\right\Vert _{M}^{2}-\operatorname{div}X\right)  ^{2}(y_{0})+$
$\left(  \Delta\left\Vert \text{X}\right\Vert _{M}^{2}-\Delta
\operatorname{div}X\right)  (y_{0})$

$\qquad\qquad\qquad\qquad\qquad\qquad-2\left\langle X,\bigtriangledown\left(
\left\Vert \text{X}\right\Vert _{M}^{2}-\operatorname{div}X\right)
\right\rangle (y_{0})$

We conclude from $\left(  11.58\right)  ,$ $\left(  11.59\right)  ,$ $\left(
11.60\right)  ;$ $\left(  11.61\right)  $ and $\left(  11.62\right)  $ that:

$\qquad\qquad L_{1}=\frac{1}{24}[-\Delta\Phi+2\left\Vert X\right\Vert _{M}%
^{2}][\triangle\Phi](y_{0})+\frac{1}{24}[\triangle^{2}\Phi](y_{0})+\frac
{1}{12}$ $<X,\nabla\triangle\Phi>(y_{0})$

\qquad$=-\frac{1}{24}[\triangle\Phi]^{2}(y_{0})+\frac{1}{12}\left\Vert
X\right\Vert _{M}^{2}[\triangle\Phi](y_{0})+\frac{1}{24}[\triangle^{2}%
\Phi](y_{0})+\frac{1}{12}$ $<X,\nabla\triangle\Phi>(y_{0})$

is given by:

$\qquad L_{1}=\frac{1}{24}\left(  \left\Vert \text{X}\right\Vert _{M}%
^{4}-(\operatorname{div}X)^{2}\right)  (y_{0})+\frac{1}{12}\left\Vert
X\right\Vert _{M}^{2}\left(  \left\Vert X\right\Vert _{M}^{2}%
-\operatorname{div}X_{M}\right)  (y_{0})$

$\qquad\qquad+\frac{1}{24}\left(  \left\Vert X\right\Vert _{M}^{2}%
-\operatorname{div}X_{M}\right)  ^{2}(y_{0})+\frac{1}{24}$ $\left(
\Delta\left\Vert X\right\Vert _{M}^{2}-\Delta\operatorname{div}X_{M}\right)
(y_{0})$

$\qquad\qquad-\frac{1}{12}\left\langle X,\left(  \bigtriangledown\left\Vert
X\right\Vert _{M}^{2}-\bigtriangledown\operatorname{div}X_{M}\right)
\right\rangle (y_{0})$

$\qquad+\frac{1}{12}$ $<X,-X\left(  \left\Vert X\right\Vert _{M}%
^{2}-\operatorname{div}X_{M}\right)  >(y_{0})+\frac{1}{12}$ $<X,\left(
\nabla\left\Vert X\right\Vert _{M}^{2}-\nabla\operatorname{div}X_{M}\right)
>(y_{0})$

Since,\qquad

$\qquad\frac{1}{12}$ $<X,-X\left(  \left\Vert X\right\Vert _{M}^{2}%
-\operatorname{div}X_{M}\right)  >(y_{0})=-\frac{1}{12}\left\Vert X\right\Vert
_{M}^{2}\left(  \left\Vert X\right\Vert _{M}^{2}-\operatorname{div}%
X_{M}\right)  (y_{0})$

$\qquad L_{1}$ simplifies to:

$\left(  11.63\right)  $\qquad$L_{1}=\frac{1}{24}\left(  \left\Vert
\text{X}\right\Vert _{M}^{4}-(\operatorname{div}X)^{2}\right)  (y_{0}%
)+\frac{1}{24}\left(  \left\Vert X\right\Vert _{M}^{2}-\operatorname{div}%
X_{M}\right)  ^{2}(y_{0})$

$\qquad\qquad\qquad\qquad+\frac{1}{24}$ $\left(  \Delta\left\Vert X\right\Vert
_{M}^{2}-\Delta\operatorname{div}X_{M}\right)  (y_{0})$

$\qquad\qquad\qquad\qquad\qquad\qquad\qquad\qquad\qquad\qquad\qquad
\qquad\qquad\qquad\qquad\qquad\qquad\ \ \ \ \ \blacksquare$

We next compute the expression:

$\qquad\qquad\qquad\qquad\ \ \ L_{2}=\frac{1}{24}\Delta\lbrack\theta^{\frac
{1}{2}}\bigtriangleup\theta^{-\frac{1}{2}}](y_{0})=\frac{1}{24}\frac
{\partial^{2}}{\partial x_{k}^{2}}[\theta^{\frac{1}{2}}\bigtriangleup
\theta^{-\frac{1}{2}}](y_{0})$

\qquad The expression is given in Appendix A$_{321}$ and in \textbf{Corollary
9} as I$_{3211}$

$\left(  11.64\right)  \qquad L_{2}=\frac{1}{24}\frac{\partial^{2}}{\partial
x_{k}^{2}}[\theta^{\frac{1}{2}}\bigtriangleup\theta^{-\frac{1}{2}}](y_{0})$

$\qquad\qquad\qquad L_{2}=+\frac{1}{720}[\left\Vert R^{M}\right\Vert
^{2}-\ \left\Vert \varrho^{M}\right\Vert ^{2}+\ 6\Delta\tau^{M}](y_{0}%
)\phi(y_{0})\qquad$I$_{3211}$

\qquad\ \ \qquad\qquad\qquad\qquad\qquad\qquad\qquad\qquad\qquad\qquad
\qquad\qquad\qquad\qquad\qquad$\qquad\blacksquare$

The last term $L_{3}$ is in final form:$\qquad$

$\left(  11.65\right)  \qquad\qquad\qquad L_{3}=+\frac{1}{12}\left(
-[\Delta\left\Vert X\right\Vert _{M}^{2}]+[\Delta V]\right)  (y_{0})$

\qquad\qquad\qquad\qquad\qquad\qquad\qquad\qquad\qquad\qquad\qquad\qquad
\qquad\qquad\qquad\qquad\qquad$\ \ \blacksquare$

We conclude from $\left(  11.55\right)  ,$ $\left(  11.62\right)  ,$ $\left(
11.63\right)  $ and $\left(  11.64\right)  $ that:

$\left(  11.66\right)  \qquad\qquad$I$_{321}$ $=\frac{1}{12}%
\underset{k=1}{\overset{n}{\sum}}\frac{\partial^{2}}{\partial x_{k}^{2}}%
[\frac{\text{L}\Psi}{\Psi}](y_{0}).\phi(y_{0})=\frac{1}{12}\Delta\lbrack
\frac{\text{L}\Psi}{\Psi}](y_{0}).\phi(y_{0})$

$\qquad$I$_{321}$ $=\frac{1}{12}\underset{k=1}{\overset{n}{\sum}}%
\frac{\partial^{2}}{\partial x_{k}^{2}}[\frac{\text{L}\Psi}{\Psi}](y_{0}%
)\phi(y_{0})=L_{1}+L_{2}+L_{3}$

\qquad$=\frac{1}{24}\left(  \left\Vert \text{X}\right\Vert _{M}^{4}%
-(\operatorname{div}X)^{2}\right)  (y_{0})+\frac{1}{24}\left(  \left\Vert
X\right\Vert _{M}^{2}-\operatorname{div}X_{M}\right)  ^{2}(y_{0})\qquad L_{1}$

$\qquad+\frac{1}{24}$ $\left(  \Delta\left\Vert X\right\Vert _{M}^{2}%
-\Delta\operatorname{div}X_{M}\right)  (y_{0})$

\qquad$+\frac{1}{720}[\left\Vert R^{M}\right\Vert ^{2}-\ \left\Vert
\varrho^{M}\right\Vert ^{2}+\ 6\Delta\tau^{M}](y_{0})\phi(y_{0})\qquad
\qquad\qquad\ \ \ L_{2}=$ $I_{3211}$

\qquad$+\frac{1}{12}\left(  -\Delta\left\Vert X\right\Vert _{M}^{2}+\Delta
V\right)  (y_{0})\qquad\qquad\qquad\qquad\qquad\qquad\qquad L_{3}$

$\frac{1}{24}\left(  \left\Vert \text{X}\right\Vert _{M}^{4}%
-(\operatorname{div}X)^{2}\right)  (y_{0})+\frac{1}{24}\left(  \left\Vert
X\right\Vert _{M}^{2}-\operatorname{div}X_{M}\right)  ^{2}(y_{0})=\frac{1}%
{12}\left(  \left\Vert \text{X}\right\Vert _{M}^{4}-\left\Vert X\right\Vert
_{M}^{2}\operatorname{div}X_{M}\right)  (y_{0})$

A last simplification here gives:

\qquad$\left(  11.67\right)  \qquad$I$_{321}$ $=\frac{1}{12}%
\underset{k=1}{\overset{n}{\sum}}\frac{\partial^{2}}{\partial x_{k}^{2}}%
[\frac{\text{L}\Psi}{\Psi}](y_{0})\phi(y_{0})=L_{1}+L_{2}+L_{3}$

\qquad\qquad$=\frac{1}{12}\left(  \left\Vert \text{X}\right\Vert _{M}%
^{4}-\left\Vert X\right\Vert _{M}^{2}\operatorname{div}X_{M}\right)
(y_{0})\qquad L_{1}+L_{3}$

$\qquad\qquad-\frac{1}{24}$ $\left(  \Delta\left\Vert X\right\Vert _{M}%
^{2}+\Delta\operatorname{div}X_{M}\right)  (y_{0})+\frac{1}{12}\frac
{\partial^{2}V}{\partial x_{i}^{2}}(y_{0})$

\qquad$\qquad+\frac{1}{720}[\left\Vert R^{M}\right\Vert ^{2}-\ \left\Vert
\varrho^{M}\right\Vert ^{2}+\ 6\Delta\tau^{M}](y_{0})\phi(y_{0})\qquad L_{2}=$
$I_{3211}$\qquad\qquad\qquad\qquad\qquad\qquad\qquad\qquad\qquad\qquad
\qquad\qquad\qquad\qquad\qquad\qquad\qquad$\qquad\qquad\qquad\qquad
\qquad\qquad\qquad\qquad\qquad\qquad\qquad\qquad\qquad\qquad\qquad\qquad
\qquad\qquad\blacksquare$\qquad\qquad\qquad\qquad\qquad\qquad\qquad
\qquad\qquad\qquad\qquad\qquad\qquad\qquad\qquad\qquad

We next compute:

$\qquad$I$_{323}$ $\mathbf{=}$ $\frac{1}{24}\underset{i=1}{\overset{n}{\sum}%
}\frac{\partial^{2}}{\partial x_{i}^{2}}[\underset{i,j=1}{\overset{n}{\sum}}%
$g$^{jk}\frac{\partial\Lambda_{k}}{\partial\text{x}_{j}}\phi\circ\pi
_{\text{P}}$ $](y_{0})$

$\qquad\mathbf{=}$ $\frac{1}{24}\underset{i=1}{\overset{n}{\sum}%
}\underset{j,k=1}{\overset{n}{\sum}}[\frac{\partial^{2}\text{g}^{jk}}{\partial
x_{i}^{2}}\frac{\partial\Lambda_{k}}{\partial\text{x}_{j}}\phi\circ
\pi_{\text{P}}$ $](y_{0})$ $\mathbf{+}$ $\frac{1}{24}%
\underset{i=1}{\overset{n}{\sum}}[\underset{j,k=1}{\overset{n}{\sum}}$%
g$^{jk}\frac{\partial^{2}}{\partial x_{i}^{2}}(\frac{\partial\Lambda_{k}%
}{\partial\text{x}_{j}}\phi\circ\pi_{\text{P}})$ ]$(y_{0})$

\qquad$\mathbf{+}\frac{1}{12}\underset{i=1}{\overset{n}{\sum}}%
\underset{j,k=1}{\overset{n}{\sum}}[\frac{\partial\text{g}^{jk}}{\partial
x_{i}}\frac{\partial}{\partial x_{i}}(\frac{\partial\Lambda_{k}}%
{\partial\text{x}_{j}}\phi\circ\pi_{\text{P}})$ $](y_{0})$

$\left(  11.68\right)  $\qquad I$_{323}$ $=$ I$_{3231}+$ I$_{3232}$ +
I$_{3233}$where,

\ I$_{3231}\mathbf{=}$ $\frac{1}{24}\underset{i=1}{\overset{n}{\sum}%
}\underset{j,k=1}{\overset{n}{\sum}}[\frac{\partial^{2}\text{g}^{jk}}{\partial
x_{i}^{2}}\frac{\partial\Lambda_{k}}{\partial\text{x}_{j}}\phi\circ
\pi_{\text{P}}](y_{0})$ $\mathbf{=}$ $\frac{1}{24}%
\underset{i=1}{\overset{n}{\sum}}\underset{j,k=1}{\overset{n}{\sum}}%
[\frac{\partial^{2}\text{g}^{jk}}{\partial x_{i}^{2}}\frac{\partial\Lambda
_{k}}{\partial\text{x}_{j}}(y_{0})\phi(y_{0})$

\ I$_{3232}=\frac{1}{24}\underset{i=1}{\overset{n}{\sum}}%
[\underset{i,j=1}{\overset{n}{\sum}}$g$^{jk}\frac{\partial^{2}}{\partial
x_{i}^{2}}(\frac{\partial\Lambda_{k}}{\partial\text{x}_{j}}\phi\circ
\pi_{\text{P}})(y_{0})=\frac{1}{24}\underset{i=1}{\overset{n}{\sum}%
}\underset{i,j=1}{\overset{n}{\sum}}[$g$^{jk}\frac{\partial^{3}\Lambda_{k}%
}{\partial x_{i}^{2}\partial\text{x}_{j}}](y_{0})\phi(y_{0})$

$\ $I$_{3233}=\frac{1}{12}\underset{i=1}{\overset{n}{\sum}}%
\underset{j,k=1}{\overset{n}{\sum}}[\frac{\partial\text{g}^{jk}}{\partial
x_{i}}\frac{\partial}{\partial x_{i}}(\frac{\partial\Lambda_{k}}%
{\partial\text{x}_{j}}\phi\circ\pi_{\text{P}})](y_{0})=\frac{1}{12}%
\underset{i=1}{\overset{n}{\sum}}\underset{j,k=1}{\overset{n}{\sum}}%
[\frac{\partial\text{g}^{jk}}{\partial x_{i}}\frac{\partial^{2}\Lambda_{k}%
}{\partial x_{i}\partial\text{x}_{j}}](y_{0})\phi(y_{0})$

We will compute \textbf{each} of the above expressions in terms of invariants
of the manifold M, the submanifold P and the vector bundle E:

Therefore we have:

\qquad I$_{3231}\mathbf{=}$ $\frac{1}{24}\underset{i,j,k=1}{\overset{n}{\sum}%
}\frac{\partial^{2}\text{g}^{jk}}{\partial x_{i}^{2}}(y_{0})[\frac
{\partial\Lambda_{k}}{\partial\text{x}_{j}}\phi\circ\pi_{\text{P}}](y_{0})$

$\frac{\partial^{2}\text{g}^{jk}}{\partial\text{x}_{i}^{2}}(y_{0})=\frac{2}%
{3}R_{jiki}(y_{0})=\frac{2}{3}R_{ijik}(y_{0})$ by (iii) of Table A$_{2},$
$\frac{\partial\Lambda_{k}}{\partial\text{x}_{j}}(y_{0})=\frac{1}{2}%
\Omega_{jk}(y_{0})$ by (vii) of \textbf{Proposition 5.}

We conclude that:

$\left(  11.69\right)  \qquad$I$_{3231}=\frac{1}{72}%
\underset{i,j,k=1}{\overset{n}{\sum}}R_{ijik}(y_{0})\Omega_{jk}(y_{0}%
)\phi(y_{0})=\frac{1}{72}\underset{j,k=1}{\overset{n}{\sum}}\varrho_{jk}%
^{M}(y_{0})\Omega_{jk}(y_{0})\phi(y_{0})$

\qquad\qquad\qquad\qquad\qquad\qquad\qquad\qquad\qquad\qquad\qquad\qquad
\qquad\qquad\qquad\qquad\qquad\qquad$\blacksquare$

We next consider:

\qquad\qquad I$_{3232}=\frac{1}{24}\underset{i=1}{\overset{n}{\sum}%
}[\underset{j,k=1}{\overset{n}{\sum}}$g$^{jk}\frac{\partial^{2}}{\partial
x_{i}^{2}}(\frac{\partial\Lambda_{k}}{\partial\text{x}_{j}}\phi\circ
\pi_{\text{P}})](y_{0})=\frac{1}{24}\underset{i=1}{\overset{n}{\sum}%
}[\underset{j,k=1}{\overset{n}{\sum}}$g$^{jk}\frac{\partial^{3}\Lambda_{k}%
}{\partial x_{i}^{2}\partial\text{x}_{j}}\phi\circ\pi_{\text{P}}](y_{0})$

Since g$^{jk}(y_{0})=\delta^{jk},$

\ \ \ I$_{3232}=\frac{1}{24}\underset{i=1}{\overset{n}{\sum}}%
[\underset{j=1}{\overset{n}{\sum}}\frac{\partial^{3}\Lambda_{j}}{\partial
x_{i}^{2}\partial\text{x}_{j}}\phi](y_{0})=\frac{1}{24}%
\underset{i,j=1}{\overset{n}{\sum}}\frac{\partial^{3}\Lambda_{j}}{\partial
x_{i}^{2}\partial\text{x}_{j}}(y_{0})\phi(y_{0})=\frac{1}{24}%
\underset{i,j=1}{\overset{n}{\sum}}\frac{\partial^{3}\Lambda_{j}}%
{\partial\text{x}_{j}\partial x_{i}^{2}}(y_{0})\phi(y_{0})$

Recall by using the formula in \textbf{Proposition 1.18} of \textbf{Berline,
Getzler, Vergne }$\left[  7\right]  ,$ we have:

\qquad\qquad\qquad\qquad$\frac{\partial^{3}\Lambda_{l}}{\partial x_{i}\partial
x_{j}\partial x_{k}}(y_{0})=\frac{1}{4}\frac{\partial^{2}\Omega_{kl}}{\partial
x_{i}\partial x_{j}}(y_{0})$

and so,

$\left(  11.70\right)  $\qquad\qquad\qquad$\frac{\partial^{3}\Lambda_{j}%
}{\partial\text{x}_{i}^{2}\partial\text{x}_{j}}(y_{0})=\frac{1}{4}%
\frac{\partial^{2}\Omega_{jj}}{\partial\text{x}_{i}^{2}}(y_{0})=0$

\qquad\qquad\qquad\qquad\qquad\qquad\qquad\qquad\qquad\qquad\qquad\qquad
\qquad\qquad\qquad\qquad\qquad$\blacksquare$

We next consider:

$\qquad\qquad$I$_{3233}=\frac{1}{12}\underset{i=1}{\overset{n}{\sum}%
}\underset{j,k=1}{\overset{n}{\sum}}[\frac{\partial\text{g}^{jk}}{\partial
x_{i}}\frac{\partial^{2}\Lambda_{k}}{\partial x_{i}\partial\text{x}_{j}%
}](y_{0})$

Since $\frac{\partial\text{g}^{jk}}{\partial x_{i}}(y_{0})=0$ for
$i,j,k=q+1,...,n,$ by (ii) of \textbf{Table A}$_{2}$\textbf{ }in
\textbf{Appendix A, }we have:

$\left(  11.71\right)  $\qquad\qquad\qquad\qquad\qquad\qquad\qquad
\ \ I$_{3233}=0$

\qquad\qquad\qquad\qquad\qquad\qquad\qquad$\qquad\qquad\qquad\qquad
\qquad\qquad\qquad\qquad\qquad\qquad\ \ \blacksquare$

We conclude from $\left(  11.68\right)  ,\left(  11.69\right)  ,$ $\left(
11.70\right)  $ and $\left(  11.71\right)  $ that:

\qquad\ I$_{323}=$ I$_{3231}+$ I$_{3232}+$ I$_{3233}$

$\qquad\ $I$_{323}=\frac{1}{72}\underset{j,k=1}{\overset{n}{\sum}}\varrho
_{jk}^{M}(y_{0})\Omega_{jk}(y_{0})\phi(y_{0})$

Changing indices, we have:

$\left(  11.72\right)  \qquad$I$_{323}=\frac{1}{72}%
\underset{i,j=1}{\overset{n}{\sum}}\varrho_{ij}^{M}(y_{0})\Omega_{ij}%
(y_{0})\phi(y_{0})$

\qquad\qquad\qquad\qquad\qquad\qquad\qquad\qquad\qquad\qquad\qquad\qquad
\qquad\qquad\qquad\qquad\qquad$\ \ \ \blacksquare$\qquad\qquad\qquad
\qquad\qquad\qquad\qquad\qquad\qquad\qquad\qquad\qquad\qquad

We next compute the expression for I$_{325}:$

\qquad I$_{325}$ $=\frac{1}{24}\frac{\partial^{2}}{\partial x_{i}^{2}%
}[\underset{i,j=1}{\overset{n}{\sum}}$g$^{jk}\Lambda_{j}\Lambda_{k}\phi
\circ\pi_{\text{P}}](y_{0})$

Again here we will repeatedly use the fact given in $\left(  7.5\right)  $ and
$\left(  7.7\right)  $ of \textbf{Chapter 7}, that:

$\frac{\partial}{\partial\text{x}_{i}}\pi_{\text{P}}$(z$_{0}$) $=\left\{
\begin{array}
[c]{c}%
1\text{ for }i=1,...,q\\
0\text{ for }i=q+1,...,n
\end{array}
\right.  $ and $\frac{\partial^{2}}{\partial\text{x}_{j}\partial\text{x}_{i}%
}\pi_{\text{P}}$(z$_{0}$) $=0$ for all $i,j=1,...,q,q+1,...,n.$

$\left(  11.73\right)  $\qquad\ I$_{325}$ $\mathbf{=}$ $\frac{1}%
{24}\underset{i=1}{\overset{n}{\sum}}\frac{\partial^{2}}{\partial x_{i}^{2}}[$
$\underset{j,k=1}{\overset{n}{\sum}}$g$^{jk}\Lambda_{j}\Lambda_{k}\phi\circ
\pi_{\text{P}}](y_{0})$

$=$ $\frac{1}{24}[\underset{i=1}{\overset{n}{\sum}}%
\underset{j,k=1}{\overset{n}{\sum}}\frac{\partial^{2}\text{g}^{jk}}{\partial
x_{i}^{2}}\left\{  \text{ }\Lambda_{j}\Lambda_{k}\phi\circ\pi_{\text{P}%
}\right\}  ](y_{0})+$ $\frac{1}{24}[\underset{j=1}{\overset{n}{\sum}%
}\underset{j,k=1}{\overset{n}{\sum}}$g$^{jk}\frac{\partial^{2}}{\partial
x_{i}^{2}}\left\{  \text{ }\Lambda_{j}\Lambda_{k}\phi\circ\pi_{\text{P}%
}\right\}  ](y_{0})$

\qquad$\qquad+$ $\frac{1}{12}[\underset{i=1}{\overset{n}{\sum}}%
\underset{j,k=1}{\overset{n}{\sum}}\frac{\partial\text{g}^{jk}}{\partial
x_{i}}\frac{\partial}{\partial x_{i}}\left\{  \text{ }\Lambda_{j}\Lambda
_{k}\phi\circ\pi_{\text{P}}\right\}  ](y_{0})$

\qquad$\qquad=$ I$_{3251}+$ I$_{3252}$ + I$_{3253}$ where,

I$_{3251}=\frac{1}{24}[\underset{i=1}{\overset{n}{\sum}}%
\underset{j,k=1}{\overset{n}{\sum}}\frac{\partial^{2}\text{g}^{jk}}{\partial
x_{i}^{2}}\left\{  \text{ }\Lambda_{j}\Lambda_{k}\phi\circ\pi_{\text{P}%
}\right\}  ](y_{0})$

I$_{3252}=$ $\frac{1}{24}[\underset{i=1}{\overset{n}{\sum}}%
\underset{j,k=1}{\overset{n}{\sum}}$g$^{jk}\frac{\partial^{2}}{\partial
x_{i}^{2}}\left\{  \text{ }\Lambda_{j}\Lambda_{k}\phi\circ\pi_{\text{P}%
}\right\}  ](y_{0})$

I$_{3253}=\frac{1}{12}[\underset{i=1}{\overset{n}{\sum}}%
\underset{j,k=1}{\overset{n}{\sum}}\frac{\partial\text{g}^{jk}}{\partial
x_{i}}\frac{\partial}{\partial x_{i}}\left\{  \text{ }\Lambda_{j}\Lambda
_{k}\phi\circ\pi_{\text{P}}\right\}  ](y_{0})$

We start with:

I$_{3251}=\frac{1}{24}[\underset{i=1}{\overset{n}{\sum}}%
\underset{j,k=1}{\overset{n}{\sum}}\frac{\partial^{2}\text{g}^{jk}}{\partial
x_{i}^{2}}\left\{  \Lambda_{j}\Lambda_{k}\phi\circ\pi_{\text{P}}\right\}
](y_{0})$

Since $\Lambda_{j}(y_{0})\Lambda_{k}(y_{0})=0$ for $j,k=q+1,...,n$ by $\left(
6.13\right)  ,$ we have:

$\left(  11.74\right)  $\qquad\qquad\qquad\qquad\qquad\qquad I$_{3251}\ =0$

Since g$^{jk}(y_{0})=\delta^{jk},$

I$_{3252}=$ $\frac{1}{24}[\underset{i=1}{\overset{n}{\sum}}%
\underset{j,k=1}{\overset{n}{\sum}}$g$^{jk}$ $\frac{\partial^{2}}{\partial
x_{i}^{2}}\left\{  \Lambda_{j}\Lambda_{k}\phi\circ\pi_{\text{P}}\right\}
](y_{0})=$ $\frac{1}{24}[\underset{i,j=1}{\overset{n}{\sum}}\frac{\partial
^{2}}{\partial x_{i}^{2}}(\Lambda_{j}^{2})\phi\circ\pi_{\text{P}}](y_{0})$

I$_{3252}=$ $\frac{1}{24}[\underset{i,j=1}{\overset{n}{\sum}}\frac
{\partial^{2}\Lambda_{j}^{2}}{\partial x_{i}^{2}}(y_{0})\phi(y_{0})\ $

By (ix) of \textbf{Proposition 5} above, we have:

$\qquad\qquad\qquad\ \frac{\partial^{2}\Lambda_{j}^{2}}{\partial\text{x}%
_{i}^{2}}(y_{0})=\frac{1}{2}\left(  \Omega_{ij}\Omega_{ij}\right)
(y_{0})+\frac{1}{3}[\frac{\partial\Omega_{ij}}{\partial\text{x}_{i}}%
\Lambda_{j}+\Lambda_{j}\frac{\partial\Omega_{ij}}{\partial\text{x}_{i}}%
](y_{0})$

Since $\Lambda_{j}(y_{0})=0$ for $j=1,...,q,q+1,...,n$ in \textbf{normal
coordinates,} we have:

$\left(  11.75\right)  $\qquad I$_{3252}=\frac{1}{48}%
\underset{i,j=1}{\overset{n}{\sum}}\left(  \Omega_{ij}\Omega_{ij}\right)
(y_{0})\phi(y_{0})\qquad$\qquad$\qquad$\qquad\qquad\qquad\qquad\qquad
\qquad\qquad\qquad\qquad\qquad\qquad\qquad\qquad\qquad\qquad$\qquad
\qquad\qquad\qquad\qquad\qquad\qquad\qquad\qquad\qquad\qquad\qquad\qquad
\qquad\qquad\qquad\qquad\qquad\blacksquare$

Since in normal coordinates, we have: $\frac{\partial\text{g}^{jk}}{\partial
x_{i}}(y_{0})=0$ for $i,j,k=1,...,n,$ then it is immediate that:

$\left(  11.76\right)  $\qquad I$_{3253}=\frac{1}{12}%
[\underset{i=1}{\overset{n}{\sum}}\underset{j,k=1}{\overset{n}{\sum}}%
\frac{\partial\text{g}^{jk}}{\partial x_{i}}\frac{\partial}{\partial x_{i}%
}\left\{  \text{ }\Lambda_{j}\Lambda_{k}\phi\circ\pi_{\text{P}}\right\}
](y_{0})=0$

Therefore, we have by $\left(  11.74\right)  ,\left(  11.75\right)  $ and
$\left(  11.76\right)  $ that:

$\left(  11.77\right)  \qquad$I$_{325}=$ I$_{3251}+$ I$_{3252}+$
I$_{3253}=\frac{1}{48}\underset{i,j=1}{\overset{n}{\sum}}\left(  \Omega
_{ij}\Omega_{ij}\right)  (y_{0})\phi(y_{0})$

Next we consider:

\qquad\textbf{I}$_{326}=-\frac{1}{24}\underset{i=1}{\overset{n}{\sum}}%
\frac{\partial^{2}}{\partial x_{i}^{2}}[\underset{j,k=1}{\overset{n}{\sum}}%
$g$^{jk}\left\{  \underset{l=1}{\overset{n}{\sum}}\Gamma_{jk}^{l}\Lambda
_{l}\phi\circ\pi_{\text{P}}\right\}  ](y_{0})$

\qquad\qquad$=-\frac{1}{24}\underset{i=1}{\overset{n}{\sum}}%
[\underset{j,k=1}{\overset{n}{\sum}}\frac{\partial^{2}\text{g}^{jk}}{\partial
x_{i}^{2}}\left\{  \underset{l=1}{\overset{n}{\sum}}\Gamma_{jk}^{l}\Lambda
_{l}\phi\circ\pi_{\text{P}}\right\}  ](y_{0})$

\qquad\qquad$-\frac{1}{24}\underset{i=1}{\overset{n}{\sum}}%
[\underset{j,k=1}{\overset{n}{\sum}}$g$^{jk}\left\{
\underset{l=1}{\overset{n}{\sum}}\frac{\partial^{2}}{\partial x_{i}^{2}%
}(\Gamma_{jk}^{l}\Lambda_{l}\phi\circ\pi_{\text{P}})\right\}  ](y_{0})$

\qquad\qquad$-\frac{1}{12}\underset{i=1}{\overset{n}{\sum}}%
[\overset{n}{\underset{j,k=1}{\sum}}\frac{\partial\text{g}^{jk}}{\partial
x_{i}}\left\{  \underset{k=1}{\overset{n}{\sum}}\frac{\partial}{\partial
x_{i}}(\Gamma_{jk}^{l}\Lambda_{l}\phi\circ\pi_{\text{P}})\right\}  ](y_{0})$

Since $\frac{\partial\text{g}^{jk}}{\partial x_{i}}(y_{0})=0;$g$^{jk}%
=\delta^{jk}$ and $\Gamma_{jk}^{l}(y_{0})=0$ by (i) of \textbf{TableA}$_{8}$
in\textbf{ normal coordinates}, we have:

\qquad\qquad I$_{326}=-\frac{1}{24}\underset{i,j,l=1}{\overset{n}{\sum}}%
\frac{\partial^{2}}{\partial x_{i}^{2}}[(\Gamma_{jj}^{l}\Lambda_{l})\phi
\circ\pi_{\text{P}})](y_{0})$

$\mathbf{\qquad\qquad=-}$ $\frac{1}{24}\underset{i,j,l=1}{\overset{n}{\sum}%
}[(\frac{\partial^{2}\Gamma_{jj}^{l}}{\partial x_{i}^{2}}\Lambda_{l}\phi
\circ\pi_{\text{P}}+\Gamma_{jj}^{l}\frac{\partial^{2}\Lambda_{l}}{\partial
x_{i}^{2}}\phi\circ\pi_{\text{P}})](y_{0})$

$\qquad\qquad\qquad\qquad\qquad\mathbf{-}$ $\frac{1}{12}%
\underset{i=q+1}{\overset{n}{\sum}}[\underset{j=1}{\overset{n}{\sum}%
}\underset{l=q+1}{\overset{n}{\sum}}\frac{\partial\Gamma_{jj}^{l}}{\partial
x_{i}}\frac{\partial\Lambda_{l}}{\partial x_{i}}\phi\circ\pi_{\text{P}}%
](y_{0})$

Since $\Lambda_{l}(y_{0})=0$ and $\Gamma_{jj}^{l}(y_{0})=0$ in \textbf{normal
coordinates,}

\qquad\qquad I$_{326}=\mathbf{-}$ $\frac{1}{12}%
\underset{i,j,k=q+1}{\overset{n}{\sum}}[\frac{\partial\Gamma_{jj}^{k}%
}{\partial x_{i}}\frac{\partial\Lambda_{k}}{\partial x_{i}}\phi](y_{0})$

We have by (viii) of Table A$_{8}$:

$\qquad\frac{\partial\Gamma_{jj}^{k}}{\partial\text{x}_{i}}(y_{0})=$ $\frac
{2}{3}R_{ijkj}(y_{0})$

and by (vii) of \textbf{Proposition 5,}

$\qquad\frac{\partial\Lambda_{k}}{\partial x_{i}}(y_{0})=\frac{1}{2}%
\Omega_{ik}(y_{0})$

Consequently we have:

\qquad\qquad I$_{326}=$\ $\mathbf{-}$ $\frac{1}{36}%
\underset{i,j,k=1}{\overset{n}{\sum}}[R_{ijkj}\Omega_{ik}](y_{0})(y_{0}%
)\phi(y_{0})=\mathbf{-}$ $\frac{1}{36}\underset{i,k=1}{\overset{n}{\sum}%
}[\varrho_{ik}^{M}\Omega_{ik}](y_{0})(y_{0})\phi(y_{0})$

Changing indices,

$\left(  11.78\right)  $\qquad I$_{326}=\mathbf{-}$ $\frac{1}{36}%
\underset{i,j=1}{\overset{n}{\sum}}[\varrho_{ij}^{M}\Omega_{ij}](y_{0}%
)(y_{0})\phi(y_{0})$

\qquad\qquad\qquad\qquad

\qquad\qquad\qquad\qquad\qquad\qquad\qquad\qquad$\qquad\qquad\qquad
\qquad\qquad\qquad\qquad\blacksquare$

Next we have:\qquad

$\left(  11.79\right)  \qquad$\textbf{I}$_{327}=\frac{1}{24}%
\underset{i=1}{\overset{n}{\sum}}\frac{\partial^{2}}{\partial x_{i}^{2}}%
[$W$\phi\circ\pi_{\text{P}}](y_{0})=\frac{1}{24}%
\underset{i=1}{\overset{n}{\sum}}\frac{\partial^{2}\text{W}}{\partial
x_{i}^{2}}(y_{0})\phi(y_{0})$

\qquad\qquad\qquad\qquad\qquad\qquad\qquad\qquad\qquad\qquad\qquad\qquad
\qquad\qquad\qquad$\blacksquare$\qquad\qquad\qquad\qquad\qquad\qquad
\qquad\qquad\qquad\qquad\qquad\qquad\qquad\qquad\qquad\qquad\qquad\qquad
\qquad\qquad\qquad$\qquad$

We next compute:

\qquad I$_{329}=\frac{1}{12}\underset{k=1}{\overset{n}{\sum}}\frac
{\partial^{2}}{\partial x_{k}^{2}}[\underset{j=1}{\overset{n}{\sum}}%
(\nabla\log\Psi)_{j}\Lambda_{j}](y_{0})\phi(y_{0})$

$\qquad\qquad=\frac{1}{12}\underset{j,k=1}{\overset{n}{\sum}}[\frac
{\partial^{2}}{\partial x_{k}^{2}}(\nabla\log\Psi)_{j}](y_{0})\Lambda
_{j}(y_{0})\phi(y_{0})+\frac{1}{12}\underset{j,k=1}{\overset{n}{\sum}}%
(\nabla\log\Psi)_{j}\frac{\partial^{2}\Lambda_{j}}{\partial x_{k}^{2}}%
](y_{0})\phi(y_{0})$

$\qquad\qquad+\frac{1}{6}\underset{j,k=1}{\overset{n}{\sum}}\frac{\partial
}{\partial x_{k}}(\nabla\log\Psi)_{j})(y_{0})\frac{\partial\Lambda_{j}%
}{\partial x_{k}}](y_{0})\phi(y_{0})$

Since $\Lambda_{j}(y_{0})=0,$ we have:

\qquad I$_{329}=\frac{1}{12}\underset{j,k=1}{\overset{n}{\sum}}(\nabla\log
\Psi)_{j}\frac{\partial^{2}\Lambda_{j}}{\partial x_{k}^{2}}](y_{0})\phi
(y_{0})+\frac{1}{6}\underset{j,k=1}{\overset{n}{\sum}}\frac{\partial}{\partial
x_{k}}(\nabla\log\Psi)_{j})(y_{0})\frac{\partial\Lambda_{j}}{\partial x_{k}%
}](y_{0})\phi(y_{0})$

$\qquad\qquad=\frac{1}{12}[(\nabla$log$\Phi_{P})_{j}\frac{\partial^{2}%
\Lambda_{j}}{\partial x_{k}^{2}}](y_{0})\phi(y_{0})+\frac{1}{12}[(\nabla
\log\theta^{-\frac{1}{2}})_{i}\frac{\partial^{2}\Lambda_{j}}{\partial
x_{k}^{2}}](y_{0})\phi(y_{0})$

$\qquad+\frac{1}{6}\underset{j,k=1}{\overset{n}{\sum}}\frac{\partial}{\partial
x_{k}}(\nabla\log\Phi)_{j})\frac{\partial\Lambda_{j}}{\partial x_{k}}%
](y_{0})\phi(y_{0})+\frac{1}{6}\underset{j,k=1}{\overset{n}{\sum}}%
[\frac{\partial}{\partial x_{k}}(\nabla\log\theta^{-\frac{1}{2}})_{j}%
)(y_{0})\frac{\partial\Lambda_{j}}{\partial x_{k}}](y_{0})\phi(y_{0})$

From Table A$_{9}$ of \textbf{Appendix A,} we have in \textbf{normal
coordinates}:

$(\nabla\log\theta^{-\frac{1}{2}})_{i}(y_{0})=0$ by (iv)$^{\ast}$ and
$\frac{\partial}{\partial x_{k}}(\nabla\log\theta^{-\frac{1}{2}})_{j}%
(y_{0})=\frac{1}{6}\varrho_{kj}(y_{0})$ by (ix)$^{\ast\ast}$ for
$i,j=1,...,q,q+1,...n.$

From Table B$_{2}$ of \textbf{Appendix B}, we have in \textbf{normal
coordinates}:

$(\nabla$log$\Phi_{P})_{j}(y_{0})=-X_{j}(y_{0})$ by (i) and $\frac{\partial
}{\partial x_{k}}(\nabla$log$\Phi_{P})_{j}(y_{0})=-\frac{\partial X_{j}%
}{\partial x_{k}}(y_{0})$ by (iii) for $i,j=1,...,q,q+1,...,n.$

Further, we have:

$\frac{\partial\Lambda_{j}}{\partial x_{k}}(y_{0})=\frac{1}{2}\Omega
_{kj}(y_{0})$ and $\frac{\partial^{2}\Lambda_{j}}{\partial x_{k}^{2}}%
(y_{0})=\frac{1}{3}$ $\frac{\partial\Omega_{kj}}{\partial x_{k}}(y_{0})$ as
already seen above.

Changing $k$ to $i,$ we have:

\qquad I$_{329}=-\frac{1}{36}[X_{j}\frac{\partial\Omega_{ij}}{\partial x_{i}%
}](y_{0})\phi(y_{0})-\frac{1}{12}[\frac{\partial X_{j}}{\partial x_{i}}%
\Omega_{ij}](y_{0})\phi(y_{0})+\frac{1}{72}[\varrho_{ij}\Omega_{ij}%
](y_{0})\phi(y_{0})$

$\left(  11.80\right)  \qquad$I$_{329}\qquad=\frac{1}{72}%
\underset{i,j=1}{\overset{n}{\sum}}[\varrho_{ij}\Omega_{ij}](y_{0})\phi
(y_{0})-\frac{1}{36}\underset{i,j=1}{\overset{n}{\sum}}[X_{j}\frac
{\partial\Omega_{ij}}{\partial\text{x}_{i}})](y_{0})\phi(y_{0})$

$\qquad\qquad\qquad\qquad\qquad-\frac{1}{12}\underset{i,j=1}{\overset{n}{\sum
}}[\frac{\partial X_{j}}{\partial x_{i}}\Omega_{ij}](y_{0})\phi(y_{0})$

Next we have the direct computation:

\qquad\textbf{L}$_{2}=\frac{1}{12}\underset{k=1}{\overset{n}{\sum}}%
\frac{\partial^{2}}{\partial x_{k}^{2}}[\underset{j=1}{\overset{n}{\sum}}%
X_{j}\Lambda_{j}](y_{0})\phi(y_{0})$

$\qquad=\frac{1}{12}\underset{k=1}{\overset{n}{\sum}}%
[\underset{j=1}{\overset{n}{\sum}}\frac{\partial^{2}X_{j}}{\partial x_{k}^{2}%
}\Lambda_{j}](y_{0})\phi(y_{0})+\frac{1}{12}\underset{k=1}{\overset{n}{\sum}%
}[\underset{j=1}{\overset{n}{\sum}}X_{j}\frac{\partial^{2}\Lambda_{j}%
}{\partial x_{k}^{2}}](y_{0})\phi(y_{0})+\frac{1}{6}%
\underset{k=1}{\overset{n}{\sum}}[\underset{j=1}{\overset{n}{\sum}}%
\frac{\partial X_{j}}{\partial x_{k}}\frac{\partial\Lambda_{j}}{\partial
x_{k}}](y_{0})\phi(y_{0})$

We use $\Lambda_{j}(y_{0})=0$ and change $k\longrightarrow i$ to have:
$\frac{\partial^{2}\Lambda_{j}}{\partial x_{i}^{2}}(y_{0})=\frac{1}{3}%
\frac{\partial\Omega_{ij}}{\partial\text{x}_{i}}(y_{0})$ and $\frac
{\partial\Lambda_{j}}{\partial x_{i}}(y_{0})=\frac{1}{2}\Omega_{ij}(y_{0})$

Consequently,

$\left(  11.81\right)  \qquad\mathbf{L}_{2}=\frac{1}{36}%
\underset{i,j=1}{\overset{n}{\sum}}[X_{j}\frac{\partial\Omega_{ij}}%
{\partial\text{x}_{i}}](y_{0})\phi(y_{0})+\frac{1}{12}%
\underset{i,j=1}{\overset{n}{\sum}}[\frac{\partial X_{j}}{\partial x_{i}%
}\Omega_{ij}](y_{0})\phi(y_{0})$

The terms of I$_{32}$ are given in $\left(  11.50\right)  :$

We collect all \textbf{finally computed} terms of I$_{32}$ from $\left(
11.67\right)  $ for I$_{321},$ $\left(  11.72\right)  $ for I$_{323},$
$\left(  11.77\right)  $ for I$_{325},$ $\left(  11.78\right)  $ for
I$_{326},$ $\left(  11.79\right)  $ for \textbf{I}$_{327},$ $\left(
11.80\right)  $ for I$_{329}$ and $\left(  8.81\right)  $ for $\mathbf{L}_{2}$
and have from $\left(  11.81\right)  $:

\qquad$\qquad$I$_{32}=\frac{1}{12}\underset{k=1}{\overset{n}{\sum}}%
\frac{\partial^{2}\Theta}{\partial x_{k}^{2}}(y_{0})=$ I$_{321}+$ I$_{323}+$
I$_{325}+$ I$_{326}+$ I$_{327}+$ I$_{329}+$ \textbf{L}$_{2}$

\qquad\qquad$=\frac{1}{12}\left(  \left\Vert \text{X}\right\Vert _{M}%
^{4}-\left\Vert X\right\Vert _{M}^{2}\operatorname{div}X_{M}\right)
(y_{0})\qquad\qquad\qquad$I$_{321}\qquad L_{1}+L_{3}$

$\qquad\qquad-\frac{1}{24}$ $\left(  \Delta\left\Vert X\right\Vert _{M}%
^{2}+\Delta\operatorname{div}X_{M}\right)  (y_{0})+\frac{1}{12}\frac
{\partial^{2}V}{\partial x_{i}^{2}}(y_{0})$

\qquad$\qquad+\frac{1}{720}[\left\Vert R^{M}\right\Vert ^{2}-\ \left\Vert
\varrho^{M}\right\Vert ^{2}+\ 6\Delta\tau^{M}](y_{0})\phi(y_{0})\qquad\qquad
L_{2}=$ $I_{3211}$

\qquad$\qquad+\frac{1}{72}\underset{i,j=1}{\overset{n}{\sum}}\varrho_{ij}%
^{M}(y_{0})\Omega_{ij}(y_{0})\phi(y_{0})\qquad$I$_{323}$

\qquad$\qquad+\frac{1}{48}\underset{i,j=1}{\overset{n}{\sum}}\left(
\Omega_{ij}\Omega_{ij}\right)  (y_{0})\phi(y_{0})\qquad$I$_{325}$

\qquad$\qquad\mathbf{-}$ $\frac{1}{36}\underset{i,j=1}{\overset{n}{\sum}%
}[\varrho_{ij}^{M}\Omega_{ij}](y_{0})(y_{0})\phi(y_{0})\qquad$I$_{326}$

\qquad$\qquad+\frac{1}{24}\underset{i=1}{\overset{n}{\sum}}\frac{\partial
^{2}\text{W}}{\partial x_{i}^{2}}(y_{0})\phi(y_{0})\qquad\qquad$%
\textbf{I}$_{327}$

\qquad$\qquad+\frac{1}{72}\underset{i,j=1}{\overset{n}{\sum}}[\varrho
_{ij}\Omega_{ij}](y_{0})\phi(y_{0})-\frac{1}{36}%
\underset{i,j=1}{\overset{n}{\sum}}[X_{j}\frac{\partial\Omega_{ij}}%
{\partial\text{x}_{i}})](y_{0})\phi(y_{0})\qquad$I$_{329}$

$\qquad\qquad-\frac{1}{12}\underset{i,j=1}{\overset{n}{\sum}}[\frac{\partial
X_{j}}{\partial x_{i}}\Omega_{ij}](y_{0})\phi(y_{0})$

\qquad$\qquad\ +\frac{1}{36}\underset{i,j=1}{\overset{n}{\sum}}[X_{j}%
\frac{\partial\Omega_{ij}}{\partial\text{x}_{i}}](y_{0})\phi(y_{0})+\frac
{1}{12}\underset{i,j=1}{\overset{n}{\sum}}[\frac{\partial X_{j}}{\partial
x_{i}}\Omega_{ij}](y_{0})\phi(y_{0})\qquad\mathbf{L}_{2}$

There are many cancellations. We see that:

I$_{323}+$ I$_{326}+$ I$_{329}+\mathbf{L}_{2}=+\frac{1}{72}%
\underset{i,j=1}{\overset{n}{\sum}}\varrho_{ij}^{M}(y_{0})\Omega_{ij}%
(y_{0})\phi(y_{0})$ $\mathbf{-}$ $\frac{1}{36}%
\underset{i,j=1}{\overset{n}{\sum}}[\varrho_{ij}^{M}\Omega_{ij}](y_{0}%
)(y_{0})\phi(y_{0})$

$\qquad\qquad+\frac{1}{72}\underset{i,j=1}{\overset{n}{\sum}}[\varrho
_{ij}\Omega_{ij}](y_{0})\phi(y_{0})-\frac{1}{36}%
\underset{i,j=1}{\overset{n}{\sum}}[X_{j}\frac{\partial\Omega_{ij}}%
{\partial\text{x}_{i}})](y_{0})\phi(y_{0})-\frac{1}{12}%
\underset{i,j=1}{\overset{n}{\sum}}[\frac{\partial X_{j}}{\partial x_{i}%
}\Omega_{ij}](y_{0})\phi(y_{0})\qquad\qquad$

\qquad\qquad$\ +\frac{1}{36}\underset{i,j=1}{\overset{n}{\sum}}[X_{j}%
\frac{\partial\Omega_{ij}}{\partial\text{x}_{i}}](y_{0})\phi(y_{0})+\frac
{1}{12}\underset{i,j=1}{\overset{n}{\sum}}[\frac{\partial X_{j}}{\partial
x_{i}}\Omega_{ij}](y_{0})\phi(y_{0})$

$\qquad\qquad=0$

We have the final expression for I$_{32}:$

$\left(  11.82\right)  \qquad$I$_{32}$ $=\frac{1}{12}%
\underset{k=1}{\overset{n}{\sum}}\frac{\partial^{2}\Theta}{\partial x_{k}^{2}%
}(y_{0})$

\qquad$=\frac{1}{24}[\left(  2\left\Vert \text{X}\right\Vert _{M}%
^{4}-2\left\Vert X\right\Vert _{M}^{2}\operatorname{div}X_{M}\right)  -\left(
\Delta\left\Vert X\right\Vert _{M}^{2}+\Delta\operatorname{div}X_{M}\right)
](y_{0})\qquad$I$_{321}\qquad L_{1}+L_{3}$

$\qquad\qquad+\frac{1}{12}\frac{\partial^{2}V}{\partial x_{i}^{2}}(y_{0})$

\qquad$\qquad+\frac{1}{720}[\left\Vert R^{M}\right\Vert ^{2}-\ \left\Vert
\varrho^{M}\right\Vert ^{2}+\ 6\Delta\tau^{M}](y_{0})\phi(y_{0})\qquad\qquad
L_{2}=$ $I_{3211}$

\qquad$\qquad+\frac{1}{48}\underset{i,j=1}{\overset{n}{\sum}}\left(
\Omega_{ij}\Omega_{ij}\right)  (y_{0})\phi(y_{0})\qquad\qquad\qquad$I$_{325}$

\qquad$\qquad+\frac{1}{24}\underset{i=1}{\overset{n}{\sum}}\frac{\partial
^{2}\text{W}}{\partial x_{i}^{2}}(y_{0})\phi(y_{0})\qquad\qquad\qquad\qquad
$\textbf{I}$_{327}$

\qquad\qquad\qquad\qquad\qquad\qquad\qquad\qquad\qquad\qquad\qquad\qquad
\qquad\qquad\qquad\qquad$\blacksquare$

The expression for b$_{2}($y$_{0}$,y$_{0},\phi)$ is given in $\left(
11.45\right)  $ as:

\qquad b$_{2}($y$_{0}$,y$_{0},\phi)=$ b$_{2}(y_{0},$y$_{0})\phi($y$_{0})=$
I$_{1}+$ I$_{32}+$ I$_{35}$where,

\qquad I$_{1}=\frac{1}{2}\frac{\text{L}\Psi}{\Psi}(y_{0})\Theta(y_{0});\qquad
$I$_{32}=\frac{1}{12}\underset{i=1}{\overset{n}{\sum}}\frac{\partial^{2}%
\Theta}{\partial x_{i}^{2}}(y_{0});\qquad$I$_{35}=$ $\frac{1}{4}\Theta(y_{0}%
)$W$(y_{0})$

The expression for I$_{1}$ is given in $\left(  11.48\right)  $ by:

$\left(  11.83\right)  \qquad$I$_{1}=\frac{1}{288}[\tau^{M}-6$ $\left\Vert
\text{X}\right\Vert _{M}^{2}-6$ $\operatorname{div}X_{M}+$ $12$V$]^{2}%
(y_{0})\phi\left(  y_{0}\right)  $

$\qquad\qquad\qquad\ \ +\frac{1}{48}[\tau^{M}-6$ $\left\Vert \text{X}%
\right\Vert _{M}^{2}-6$ $\operatorname{div}X_{M}+$ $6$V$](y_{0})$W$\left(
y_{0}\right)  \phi\left(  y_{0}\right)  $

The expression for I$_{35}$ is given in $\left(  11.49\right)  $ by:

$\left(  11.84\right)  $\qquad I$_{35}=\frac{1}{48}[(\tau^{M}-6\left\Vert
\text{X}\right\Vert _{M}^{2}-6$ $\operatorname{div}X_{M}+12$V $+$ $6$%
W$](y_{0})$W$(y_{0})\phi\left(  y_{0}\right)  $

\qquad\qquad\qquad\qquad\qquad\qquad\qquad\qquad\qquad\qquad\qquad\qquad
\qquad\qquad\qquad$\blacksquare$

Consequently, from $\left(  11.82\right)  $ for I$_{32},$ $\left(
11.83\right)  $ for I$_{1}$ and $\left(  11.84\right)  $ for I$_{35},$ we have
the final expression for b$_{2}($y$_{0}$,y$_{0},\phi):$

\begin{theorem}
In the case that the vector field X is a gradient vector field, we have:
\end{theorem}

b$_{2}($y$_{0}$,y$_{0},\phi)=$ b$_{2}(y_{0},$y$_{0})\phi($y$_{0})=$ I$_{1}+$
I$_{32}+$ I$_{35}$

$=+\frac{1}{288}[\tau^{M}-6$ $\left\Vert \text{X}\right\Vert _{M}^{2}-6$
$\operatorname{div}X_{M}+$ $12$V$]^{2}(y_{0})\phi\left(  y_{0}\right)
\qquad\qquad$I$_{1}$

$\ +\frac{1}{48}[\tau^{M}-6$ $\left\Vert \text{X}\right\Vert _{M}^{2}-6$
$\operatorname{div}X_{M}+$ $6$V$](y_{0})$W$\left(  y_{0}\right)  \phi\left(
y_{0}\right)  $

$+\frac{1}{48}[(\tau^{M}-6\left\Vert \text{X}\right\Vert _{M}^{2}-6$
$\operatorname{div}X_{M}+12$V $+$ $6$W$](y_{0})$W$(y_{0})\phi\left(
y_{0}\right)  \qquad$I$_{35}$

$+\frac{1}{24}[\left(  2\left\Vert \text{X}\right\Vert _{M}^{4}-2\left\Vert
X\right\Vert _{M}^{2}\operatorname{div}X_{M}\right)  -\left(  \Delta\left\Vert
X\right\Vert _{M}^{2}+\Delta\operatorname{div}X_{M}\right)  ](y_{0})\qquad
$I$_{32}\qquad$I$_{321}\qquad L_{1}+L_{3}$

$+\frac{1}{12}\frac{\partial^{2}\text{V}}{\partial x_{i}^{2}}(y_{0})$

$+\frac{1}{720}[\left\Vert R^{M}\right\Vert ^{2}-\ \left\Vert \varrho
^{M}\right\Vert ^{2}+\ 6\Delta\tau^{M}](y_{0})\phi(y_{0})\qquad\qquad L_{2}=$
$I_{3211}$

$+\frac{1}{48}\underset{i,j=1}{\overset{n}{\sum}}\left(  \Omega_{ij}%
\Omega_{ij}\right)  (y_{0})\phi(y_{0})\qquad\qquad\qquad$I$_{325}$

$+\frac{1}{24}\underset{i=1}{\overset{n}{\sum}}\frac{\partial^{2}\text{W}%
}{\partial x_{i}^{2}}(y_{0})\phi(y_{0})\qquad\qquad\qquad\qquad$%
\textbf{I}$_{327}$

\begin{proof}
We have the above expression from $\left(  11.82\right)  $ for I$_{32},$
$\left(  11.83\right)  $ for I$_{1}$ and $\left(  11.84\right)  $ for I$_{35}$
as stated above.
\end{proof}

The Corollary is thus proved.\qquad

\qquad\qquad\qquad\qquad\qquad\qquad\qquad\qquad\qquad\qquad\qquad\qquad
\qquad\qquad\qquad\qquad\qquad\qquad\qquad\qquad$\blacksquare$\qquad
\qquad\qquad\qquad\qquad\qquad\qquad\qquad

We observe that I$_{322},$ I$_{324},$ I$_{325}$ and I$_{326}$ are "casualties"
as they have been wiped off.

\appendix{}

\part{APPENDICES}

\chapter{Derivatives of Components of the Metric Tensor}

It is important to note that since all expansions were carried out in Fermi
coordinates, all differentiations in tangential Fermi coordinates are zero and
so we will consider all differentiation with respect to normal Fermi
coordinates only.

For all the computations in this \textbf{Appendix Chapter}, we will use the
\textbf{Preliminary Geometric Expansion Formulae} of \textbf{Chapter 10 }in
\textbf{Part 4}. In particular, Table A$_{1}-$ Table A$_{8}$ use the
expansions of the components g$_{ij}$ of the metric tensor and its inverse
g$^{ij}.$ Tables A$_{9}-$ Table A$_{10}$ will use the expansion of $\theta
_{P}^{-\frac{1}{2}}$ where $\ \theta_{P}(x)=\sqrt{\text{det}g(x)}$ is the
determinant of the matrix $\left(  g_{ij}(x)\right)  i,j=1,...,q,...,n$
defined $\left(  1.6\right)  $ in \textbf{Chapter 1} here.\qquad\qquad
\qquad\qquad\qquad\qquad\qquad\qquad\qquad\qquad\qquad\qquad\qquad\qquad
\qquad\qquad\qquad\qquad\qquad\qquad$\qquad\qquad\qquad\qquad\qquad
\qquad\qquad\qquad\qquad\qquad\qquad\qquad\qquad\qquad\qquad\qquad\qquad
\qquad\qquad\blacksquare$

\section{Table A$_{1}$\qquad\qquad\qquad\qquad\qquad\qquad\qquad\qquad
\qquad\qquad\qquad\qquad\qquad\qquad\qquad\qquad\qquad\qquad}

For $i,j,k,l=q+1,...,n$

(i) g$_{kl}(y_{0})=\delta_{kl}$

\vspace{1pt}(ii) $\frac{\partial\text{g}_{kl}}{\partial\text{x}_{i}}(y_{0})=0$

(iii) \ $\frac{\partial^{2}\text{g}_{pl}}{\partial\text{x}_{i}\partial
\text{x}_{j}}(y_{0})=-\frac{1}{3}(R_{ipjl}+R_{jpil})(y_{0})$

In particular,

\ \ \ \ \ \ $\frac{\partial^{2}\text{g}_{kl}}{\partial\text{x}_{i}^{2}}%
(y_{0})=-\frac{1}{3}(R_{ikil}+R_{ikil})(y_{0})=-\frac{2}{3}(R_{ikil})(y_{0})$

When Fermi coordinates reduce to normal coordinates then,

\qquad\ $\frac{\partial^{2}\text{g}_{kl}}{\partial\text{x}_{i}^{2}}%
(y_{0})=-\frac{2}{3}\varrho_{kl}(y_{0})$

(iv) $\frac{\partial^{3}g_{pl}}{\partial x_{i}\partial x_{j}\partial x_{k}%
}(y_{0})=-\frac{1}{6}[\nabla_{_{k}}R_{jkil}+\nabla_{_{k}}R_{ikjl}+\nabla
_{_{j}}R_{kpil}+\nabla_{_{i}}R_{kpjl}+\nabla_{_{j}}R_{ipkl}+\nabla_{_{i}%
}R_{jpkl}](y_{0})$

In particular,

$\qquad\frac{\partial^{3}g_{kl}}{\partial x_{i}\partial x_{j}\partial x_{k}%
}(y_{0})=-\frac{1}{6}[\nabla_{_{k}}R_{jkil}+\nabla_{_{k}}R_{ikjl}-\nabla
_{_{j}}\varrho_{il}-\nabla_{_{i}}\varrho_{jl}](y_{0}$

In particular,

\qquad\ $\frac{\partial^{3}\text{g}_{kl}}{\partial\text{x}_{i}^{2}%
\partial\text{x}_{j}}(y_{0})=-\frac{1}{3}(\nabla_{i}R_{ikjl}+\nabla
_{i}R_{jkil}+\nabla_{j}R_{ikil})(y_{0})$\qquad

In particular,\ when Fermi coordinates reduce to normal coordinates,

$\qquad\frac{\partial^{3}\text{g}_{jl}}{\partial\text{x}_{i}^{2}%
\partial\text{x}_{j}}(y_{0})=-\frac{1}{3}(\nabla_{i}R_{ijjl}+\nabla
_{i}R_{jjil}+\nabla_{j}R_{ijil})(y_{0})=-\frac{1}{3}(\nabla_{i}R_{ijjl}%
+\nabla_{j}R_{ijil})(y_{0})$

$\qquad\qquad\qquad\ \ =-\frac{1}{3}(-\nabla_{i}R_{jijl}+\nabla_{j}%
R_{ijil})(y_{0})=-\frac{1}{3}(-\nabla_{i}\varrho_{il}+\nabla_{j}\varrho
_{jl})(y_{0})=-\frac{1}{3}\frac{1}{2}(-\nabla_{l}\tau+\nabla_{l}\tau
)(y_{0})=0$

In particular,

$\qquad\frac{\partial^{3}\text{g}_{il}}{\partial x_{i}\partial x_{j}^{2}%
}(y_{0})$ $=-\frac{1}{3}(\nabla_{j}R_{jiil}+\nabla_{j}R_{iijl}+\nabla
_{i}R_{jijl})(y_{0})=-\frac{1}{3}(-\nabla_{j}R_{ijil}+\nabla_{i}%
R_{jijl})(y_{0})$

\qquad\qquad\qquad$\ \ \ =-\frac{1}{3}(-\nabla_{j}\varrho_{jq}+\nabla
_{i}\varrho_{iq})(y_{0})=-\frac{1}{3}\frac{1}{2}(-\nabla_{q}\tau+\nabla
_{q}\tau)(y_{0})=0$

We see from the last two equations above that:

\begin{center}
$\frac{\partial^{3}\text{g}_{jl}}{\partial\text{x}_{i}^{2}\partial\text{x}%
_{j}}(y_{0})=0=\frac{\partial^{3}\text{g}_{il}}{\partial x_{i}\partial
x_{j}^{2}}(y_{0})$\qquad
\end{center}

The index $q$ in the computations below should \textbf{not} be confused with
the dimension $q$ of the submanifold P.

(v) $\frac{\partial^{4}g_{pq}}{\partial x_{i}\partial x_{j}\partial
x_{k}\partial x_{l}}(y_{0})$

\begin{center}
$=\frac{1}{360}[(-18\nabla_{lk}^{2}$R$_{jpiq}+16\underset{w=\text{1}%
}{\overset{n}{\sum}}$R$_{lpkw}$R$_{jqiw})+(-18\nabla_{lk}^{2}$R$_{ipjq}%
+16\underset{w=\text{1}}{\overset{n}{\sum}}$R$_{lpkw}$R$_{iqjw})+(-18\nabla
_{lj}^{2}$R$_{kpiq}+16\underset{w=\text{1}}{\overset{n}{\sum}}$R$_{lpjw}%
$R$_{kqiw})$

$+(-18\nabla_{li}^{2}$R$_{kpjq}+16\underset{w=\text{1}}{\overset{n}{\sum}}%
$R$_{lpiw}$R$_{kqjw})+(-18\nabla_{lj}^{2}$R$_{ipkq}+16\underset{w=\text{1}%
}{\overset{n}{\sum}}$R$_{lpjw}$R$_{iqkw})+(-18\nabla_{li}^{2}$R$_{jpkq}%
+16\underset{w=\text{1}}{\overset{n}{\sum}}$R$_{lpiw}$R$_{jqkw})$

$+(-18\nabla_{kl}^{2}$R$_{jpiq}+16\underset{w=\text{1}}{\overset{n}{\sum}}%
$R$_{kplw}$R$_{jqiw})+(-18\nabla_{kl}^{2}$R$_{ipjq}+16\underset{w=\text{1}%
}{\overset{n}{\sum}}$R$_{kplw}$R$_{iqjw})+(-18\nabla_{jl}^{2}$R$_{kpiq}%
+16\underset{w=\text{1}}{\overset{n}{\sum}}$R$_{jplw}$R$_{kqiw})$

$+(-18\nabla_{il}^{2}$R$_{kpjq}+16\underset{w=\text{1}}{\overset{n}{\sum}}%
$R$_{iplw}$R$_{kqjw})+(-18\nabla_{jl}^{2}$R$_{ipkq}+16\underset{w=\text{1}%
}{\overset{n}{\sum}}$R$_{jplw}$R$_{iqkw})+(-18\nabla_{il}^{2}$R$_{jpkq}%
+16\underset{w=\text{1}}{\overset{n}{\sum}}$R$_{iplw}$R$_{jqkw})$

$+(-18\nabla_{kj}^{2}$R$_{lpiq}+16\underset{w=\text{1}}{\overset{n}{\sum}}%
$R$_{kpjw}$R$_{lqiw})+(-18\nabla_{ki}^{2}$R$_{lpjq}+16\underset{w=\text{1}%
}{\overset{n}{\sum}}$R$_{kpiw}$R$_{lqjw})+(-18\nabla_{jk}^{2}$R$_{lpiq}%
+16\underset{w=\text{1}}{\overset{n}{\sum}}$R$_{jpkw}$R$_{lqiw})$

$+(-18\nabla_{ik}^{2}$R$_{lpjq}+16\underset{w=\text{1}}{\overset{n}{\sum}}%
$R$_{ipkw}$R$_{lqjw})+(-18\nabla_{ji}^{2}$R$_{lpkq}+16\underset{w=\text{1}%
}{\overset{n}{\sum}}$R$_{jpiw}$R$_{lqkw})+(-18\nabla_{ij}^{2}$R$_{lpkq}%
+16\underset{w=\text{1}}{\overset{n}{\sum}}$R$_{ipjw}$R$_{lqkw}$

$+(-18\nabla_{kj}^{2}$R$_{iplq}+16\underset{w=\text{1}}{\overset{n}{\sum}}%
$R$_{kpjw}$R$_{iqlw})+(-18\nabla_{ki}^{2}$R$_{jplq}+16\underset{w=\text{1}%
}{\overset{n}{\sum}}$R$_{kpiw}$R$_{jqlw})+(-18\nabla_{jk}^{2}$R$_{iplq}%
+16\underset{w=\text{1}}{\overset{n}{\sum}}$R$_{jpkw}$R$_{iqlw})$

$\bigskip+(-18\nabla_{ik}^{2}$R$_{jplq}+16\underset{w=\text{1}%
}{\overset{n}{\sum}}$R$_{ipkw}$R$_{jqlw})+(-18\nabla_{ji}^{2}$R$_{kplq}%
+16\underset{w=\text{1}}{\overset{n}{\sum}}$R$_{jpiw}$R$_{kqlw})+(-18\nabla
_{ij}^{2}$R$_{kplq}+16\underset{w=\text{1}}{\overset{n}{\sum}}$R$_{ipjw}%
$R$_{kqlw})$
\end{center}

(vi)$\qquad\frac{\partial^{4}g_{jq}}{\partial x_{i}\partial x_{j}\partial
x_{k}^{2}}(y_{0})\qquad$

$=\frac{1}{360}[(16\underset{w=\text{1}}{\overset{n}{\sum}}\varrho_{jw}%
$R$_{jqiw})+(18\nabla_{kk}^{2}\varrho_{iq}+16\underset{w=\text{1}%
}{\overset{n}{\sum}}\varrho_{jw}$R$_{iqjw})-(18\nabla_{kj}^{2}$R$_{kjiq}%
+16\underset{w=\text{1}}{\overset{n}{\sum}}$R$_{kw}$R$_{kqiw})$

$+(18\nabla_{ki}^{2}\varrho_{kq}+16\underset{w=\text{1}}{\overset{n}{\sum}}%
$R$_{kjiw}$R$_{kqjw})+(-18\nabla_{kj}^{2}$R$_{ijkq}-16\underset{w=\text{1}%
}{\overset{n}{\sum}}\varrho_{kw}$R$_{iqkw})+(16\underset{w=\text{1}%
}{\overset{n}{\sum}}$R$_{kjiw}$R$_{jqkw})$

$+(16\underset{w=\text{1}}{\overset{n}{\sum}}\varrho_{jw}$R$_{jqiw}%
)+(18\nabla_{kk}^{2}\varrho_{iq}+16\underset{w=\text{1}}{\overset{n}{\sum}%
}\varrho_{jw}$R$_{iqjw})-18\nabla_{jk}^{2}$R$_{kjiq}$

$+(18\nabla_{ik}^{2}\varrho_{kq}+16\underset{w=\text{1}}{\overset{n}{\sum}}%
$R$_{ijkw}$R$_{kqjw})-18\nabla_{jk}^{2}$R$_{ijkq}+16\underset{w=\text{1}%
}{\overset{n}{\sum}}$R$_{ijkw}$R$_{jqkw}$

$-(18\nabla_{kj}^{2}$R$_{kjiq}+16\underset{w=\text{1}}{\overset{n}{\sum}%
}\varrho_{kw}$R$_{kqiw})+(18\nabla_{ki}^{2}\varrho_{kq}+16\underset{w=\text{1}%
}{\overset{n}{\sum}}$R$_{kjiw}$R$_{kqjw})+18\nabla_{jk}^{2}$R$_{jkiq}$

$+(18\nabla_{ik}^{2}\varrho_{kq}+16\underset{w=\text{1}}{\overset{n}{\sum}}%
$R$_{ijkw}$R$_{kqjw})-18\nabla_{ji}^{2}\varrho_{jq}-(18\nabla_{ij}^{2}%
\varrho_{jq}+16\underset{w=\text{1}}{\overset{n}{\sum}}\varrho_{iw}$%
R$_{kqkw})$

$-(18\nabla_{kj}^{2}$R$_{ijkq}+16\underset{w=\text{1}}{\overset{n}{\sum}%
}\varrho_{kw}$R$_{iqkw})+(16\underset{w=\text{1}}{\overset{n}{\sum}}$%
R$_{kjiw}$R$_{jqkw})-18\nabla_{jk}^{2}$R$_{ijkq}$

$+(16\underset{w=\text{1}}{\overset{n}{\sum}}$R$_{ijkw}$R$_{jqkw}%
)-18\nabla_{ji}^{2}\varrho_{jq}-(18\nabla_{ij}^{2}\varrho_{jq}%
+16\underset{w=\text{1}}{\overset{n}{\sum}}\varrho_{iw}\varrho_{qw})](y_{0})$

\subsection{\protect\underline{\textbf{Computations}}}

\vspace{1pt}This is a direct computation using \textbf{Proposition 7}. For
example (iv) is computed as follows:

(iv) g$_{pq}$(x$_{0}$) = $\delta_{pq}-\frac{1}{3}%
\underset{r,s=q+1}{\overset{n}{\sum}}$R$_{\text{r}p\text{s}q}($y$_{0}%
)$x$_{\text{r}}$x$_{\text{s}}-\frac{1}{6}\underset{r,s,t=q+1}{\overset{n}{\sum
}}\nabla_{\text{r}}$R$_{\text{s}p\text{t}q}($y$_{0})$x$_{\text{r}}%
$x$_{\text{s}}$x$_{\text{t}}$

$\qquad\ \ \ \ +\frac{1}{360}\overset{n}{\underset{r,s,t,u=q+1}{\sum}%
}(-18\nabla_{\text{rs}}^{2}$R$_{\text{t}p\text{u}q}+16\underset{\text{w=q+1}%
}{\overset{\text{n}}{\sum}}$R$_{\text{r}p\text{js}}$R$_{\text{k}q\text{lw}}%
)($y$_{0})$x$_{\text{r}}$x$_{\text{s}}$x$_{\text{t}}$x$_{\text{u}}$

\ + $\frac{1}{90}\overset{n}{\underset{\text{r,s,t,u,v=q+1}}{\sum
\{-\nabla_{\text{rst}}^{3}}\text{R}_{\text{u}p\text{v}q}}%
+2\overset{n}{\underset{\text{w=q+1}}{\sum}}(\nabla_{r}$R$_{\text{s}%
p\text{tw}}$R$_{\text{l}q\text{hw}}+\nabla_{r}$R$_{\text{j}p\text{ks}}%
$R$_{\text{l}q\text{vs}})\}(y_{0})$x$_{\text{r}}$x$_{\text{s}}$x$_{\text{t}}%
$x$_{\text{u}}$x$_{\text{v}}\qquad\qquad\qquad$

\ + \ higher order terms.

$\frac{\partial^{3}g_{pq}}{\partial x_{i}\partial x_{j}\partial x_{k}}%
(y_{0})=\frac{\partial^{3}g_{pq}}{\partial x_{i}\partial x_{j}\partial x_{k}%
}(y_{0})=-\frac{1}{6}\underset{r,s,t=1}{\overset{n}{\sum}}\nabla_{_{r}%
}R_{sk??tq}(y_{0})\frac{\partial^{3}}{\partial x_{i}\partial x_{j}\partial
x_{k}}($x$_{r}$x$_{s}$x$_{t})$

$\frac{\partial^{3}}{\partial x_{i}\partial x_{j}\partial x_{k}}($x$_{r}%
$x$_{s}$x$_{t})=\frac{\partial^{2}}{\partial x_{i}\partial x_{j}}(\delta
_{kr}x_{s}x_{t}+x_{r}\delta_{ks}x_{t}+x_{r}x_{s}\delta_{kt})$

$=\frac{\partial^{2}}{\partial x_{i}\partial x_{j}}(\delta_{kr}x_{s}%
x_{t})+\frac{\partial^{2}}{\partial x_{i}\partial x_{j}}(x_{r}\delta_{ks}%
x_{t})+\frac{\partial^{2}}{\partial x_{i}\partial x_{j}}(x_{r}x_{s}\delta
_{kt})$

$=\delta_{kr}\frac{\partial}{\partial x_{i}}(\delta_{js}x_{t}+x_{s}\delta
_{jt})+\delta_{ks}\frac{\partial}{\partial x_{i}}(\delta_{jr}x_{t}+x_{r}%
\delta_{jt})+\delta_{kt}\frac{\partial}{\partial x_{i}}(\delta_{jr}x_{s}%
+x_{r}\delta_{js})$

$=\delta_{kr}(\delta_{js}\delta_{it}+\delta_{is}\delta_{jt})+\delta
_{ks}(\delta_{jr}\delta_{it}+\delta_{ir}\delta_{jt})+\delta_{kt}(\delta
_{jr}\delta_{is}+\delta_{ir}\delta_{js})$

$=\delta_{kr}\delta_{js}\delta_{it}+\delta_{kr}\delta_{is}\delta_{jt}%
+\delta_{ks}\delta_{jr}\delta_{it}+\delta_{ks}\delta_{ir}\delta_{jt}%
+\delta_{kt}\delta_{jr}\delta_{is}+\delta_{kt}\delta_{ir}\delta_{js}$

Therefore,

$\frac{\partial^{3}g_{pq}}{\partial x_{i}\partial x_{j}\partial x_{k}}%
(y_{0})=-\frac{1}{6}\underset{r,s,t=1}{\overset{n}{\sum}}\nabla_{_{r}}%
R_{sptq}(y_{0})[\delta_{kr}\delta_{js}\delta_{it}+\delta_{kr}\delta_{is}%
\delta_{jt}+\delta_{ks}\delta_{jr}\delta_{it}$

$\qquad\qquad\qquad\qquad+\delta_{ks}\delta_{ir}\delta_{jt}+\delta_{kt}%
\delta_{jr}\delta_{is}+\delta_{kt}\delta_{ir}\delta_{js}]$

$=-\frac{1}{6}[\nabla_{_{k}}R_{jkiq}+\nabla_{_{k}}R_{ikjq}+\nabla_{_{j}%
}R_{kpiq}+\nabla_{_{i}}R_{kpjq}+\nabla_{_{j}}R_{ipkq}+\nabla_{_{i}}%
R_{jpkq}](y_{0})$

(v) We use the expansion formula in \textbf{Proposition 7} above, or in
\textbf{Corollary 9.8} of \textbf{Gray }$\left[  4\right]  .$

\qquad For $p,q=1,...,n$ and $x\in$M$_{0},$ we have:

$\qquad g_{pq}(x)=$ $\delta_{pq}-$ $\frac{1}{3}%
\underset{r,s=1}{\overset{n}{\sum}}($ R$_{rpsq})(y_{0})$x$_{r}$x$_{s}$

$-\frac{1}{6}\underset{r,s,t=1}{\overset{n}{\sum}}\nabla_{_{r}}$R$_{sptq}%
($y$_{0})$x$_{r}$x$_{s}$x$_{t}$

$\ \ +\frac{1}{360}\overset{n}{\underset{r,s,t,u=1}{\sum}}(-18\nabla_{rs}^{2}%
$R$_{tpuq}+16\underset{w=\text{1}}{\overset{n}{\sum}}$R$_{rpsw}$R$_{tquw}%
)($y$_{0})$x$_{r}$x$_{s}$x$_{t}$x$_{u}$

$\ +$ $\frac{1}{90}\overset{n}{\underset{r,s,t,u,v=1}{\sum\{-\nabla_{rst}^{3}%
}\text{R}_{u\alpha v\beta}}+$ $2\overset{n}{\underset{w=1}{\sum}}(\nabla_{r}%
$R$_{s\alpha tw}$R$_{u\beta vw}+\nabla_{r}$R$_{s\beta tw}$R$_{u\alpha
vw})\}(y_{0})$x$_{r}$x$_{s}$x$_{t}$x$_{u}$x$_{v}\qquad\qquad\qquad$

$\ +$ \ higher order terms.

Consequently,

$\qquad\frac{\partial^{4}g_{pq}}{\partial x_{i}\partial x_{j}\partial
x_{k}\partial x_{l}}(y_{0})=\frac{1}{360}\overset{n}{\underset{r,s,t,u=1}{\sum
}}(-18\nabla_{rs}^{2}$R$_{tpuq}+16\underset{w=\text{1}}{\overset{n}{\sum}}%
$R$_{rpsw}$R$_{tquw})(y_{0})\frac{\partial^{4}}{\partial x_{i}\partial
x_{j}\partial x_{k}\partial x_{l}}($x$_{r}$x$_{s}$x$_{t}$x$_{u})$

$\qquad=\frac{1}{360}\overset{n}{\underset{r,s,t,u=1}{\sum}}(-18\nabla
_{rs}^{2}$R$_{tpuq}+16\underset{w=\text{1}}{\overset{n}{\sum}}$R$_{rpsw}%
$R$_{tquw})(y_{0})$

$\qquad\times\frac{\partial^{3}}{\partial x_{i}\partial x_{j}\partial x_{k}%
}[\delta_{lr}$x$_{s}$x$_{t}$x$_{u}+x_{r}\delta_{ls}$x$_{t}$x$_{u}+x_{r}%
x_{s}\delta_{lt}$x$_{u}+x_{r}x_{s}x_{t}\delta_{lu}]$

$\qquad=\frac{1}{360}\overset{n}{\underset{r,s,t,u=1}{\sum}}(-18\nabla
_{rs}^{2}$R$_{tpuq}+16\underset{w=\text{1}}{\overset{n}{\sum}}$R$_{rpsw}%
$R$_{tquw})(y_{0})$

$\qquad\times\frac{\partial^{2}}{\partial x_{i}\partial x_{j}}[\delta
_{lr}\delta_{ks}$x$_{t}$x$_{u}+\delta_{lr}x_{s}\delta_{kt}x_{u}+\delta
_{lr}x_{s}x_{t}\delta_{ku}+\delta_{kr}\delta_{ls}$x$_{t}$x$_{u}+x_{r}%
\delta_{ls}\delta_{kt}x_{u}+x_{r}\delta_{ls}x_{t}\delta_{ku}]$

$\qquad+\frac{\partial^{2}}{\partial x_{i}\partial x_{j}}[\delta_{kr}%
x_{s}\delta_{lt}$x$_{u}+x_{r}\delta_{ks}\delta_{lt}x_{u}+x_{r}x_{s}\delta
_{lt}\delta_{ku}+\delta_{kr}x_{s}x_{t}\delta_{lu}+x_{r}\delta_{ks}x_{t}%
\delta_{lu}+x_{r}x_{s}\delta_{kt}\delta_{lu}]$

\qquad$=\frac{1}{360}\overset{n}{\underset{r,s,t,u=1}{\sum}}(-18\nabla
_{rs}^{2}$R$_{tpuq}+16\underset{w=\text{1}}{\overset{n}{\sum}}$R$_{rpsw}%
$R$_{tquw})(y_{0})$

$\qquad\times\frac{\partial}{\partial x_{i}}[\delta_{lr}\delta_{ks}\delta
_{jt}$x$_{u}+\delta_{lr}\delta_{ks}x_{t}\delta_{ju}+\delta_{lr}\delta
_{js}\delta_{kt}x_{u}+\delta_{lr}x_{s}\delta_{kt}\delta_{ju}+\delta_{lr}%
\delta_{js}x_{t}\delta_{ku}+\delta_{lr}x_{s}\delta_{jt}\delta_{ku}]$

\qquad$\ +\frac{\partial}{\partial x_{i}}[\delta_{kr}\delta_{ls}\delta_{jt}%
$x$_{u}+\delta_{kr}\delta_{ls}x_{t}\delta_{ju}+\delta_{jr}\delta_{ls}%
\delta_{kt}x_{u}+x_{r}\delta_{ls}\delta_{kt}\delta_{ju}+\delta_{jr}\delta
_{ls}x_{t}\delta_{ku}+x_{r}\delta_{ls}\delta_{jt}\delta_{ku}]$

$\qquad+\frac{\partial}{\partial x_{i}}[\delta_{kr}\delta_{js}\delta_{lt}%
$x$_{u}+\delta_{kr}x_{s}\delta_{lt}\delta_{ju}+\delta_{jr}\delta_{ks}%
\delta_{lt}x_{u}+x_{r}\delta_{ks}\delta_{lt}\delta_{ju}+\delta_{jr}x_{s}%
\delta_{lt}\delta_{ku}+x_{r}\delta_{js}\delta_{lt}\delta_{ku}]$

\qquad$+\frac{\partial}{\partial x_{i}}[\delta_{kr}\delta_{js}x_{t}\delta
_{lu}+\delta_{kr}x_{s}\delta_{jt}\delta_{lu}+\delta_{jr}\delta_{ks}x_{t}%
\delta_{lu}+x_{r}\delta_{ks}\delta_{jt}\delta_{lu}+\delta_{jr}x_{s}\delta
_{kt}\delta_{lu}+x_{r}\delta_{js}\delta_{kt}\delta_{lu}]$

\qquad$=\frac{1}{360}\overset{n}{\underset{r,s,t,u=1}{\sum}}(-18\nabla
_{rs}^{2}$R$_{tpuq}+16\underset{w=\text{1}}{\overset{n}{\sum}}$R$_{rpsw}%
$R$_{tquw})(y_{0})$

$\qquad\times\lbrack\delta_{lr}\delta_{ks}\delta_{jt}\delta_{iu}+\delta
_{lr}\delta_{ks}\delta_{it}\delta_{ju}+\delta_{lr}\delta_{js}\delta_{kt}%
\delta_{iu}+\delta_{lr}\delta_{is}\delta_{kt}\delta_{ju}+\delta_{lr}%
\delta_{js}\delta_{it}\delta_{ku}+\delta_{lr}\delta_{is}\delta_{jt}\delta
_{ku}]$

\qquad$\ +[\delta_{kr}\delta_{ls}\delta_{jt}\delta_{iu}+\delta_{kr}\delta
_{ls}\delta_{it}\delta_{ju}+\delta_{jr}\delta_{ls}\delta_{kt}\delta
_{iu}+\delta_{ir}\delta_{ls}\delta_{kt}\delta_{ju}+\delta_{jr}\delta
_{ls}\delta_{it}\delta_{ku}+\delta_{ir}\delta_{ls}\delta_{jt}\delta_{ku}]$

$\qquad+[\delta_{kr}\delta_{js}\delta_{lt}\delta_{iu}+\delta_{kr}\delta
_{is}\delta_{lt}\delta_{ju}+\delta_{jr}\delta_{ks}\delta_{lt}\delta
_{iu}+\delta_{ir}\delta_{ks}\delta_{lt}\delta_{ju}+\delta_{jr}\delta
_{is}\delta_{lt}\delta_{ku}+\delta_{ir}\delta_{js}\delta_{lt}\delta_{ku}]$

\qquad$+[\delta_{kr}\delta_{js}\delta_{it}\delta_{lu}+\delta_{kr}\delta
_{is}\delta_{jt}\delta_{lu}+\delta_{jr}\delta_{ks}\delta_{it}\delta
_{lu}+\delta_{ir}\delta_{ks}\delta_{jt}\delta_{lu}+\delta_{jr}\delta
_{is}\delta_{kt}\delta_{lu}+\delta_{ir}\delta_{js}\delta_{kt}\delta_{lu}]$

$\frac{\partial^{4}g_{pq}}{\partial x_{i}\partial x_{j}\partial x_{k}\partial
x_{l}}(y_{0})\qquad$

\begin{center}
$=\frac{1}{360}[(-18\nabla_{lk}^{2}$R$_{jpiq}+16\underset{w=\text{1}%
}{\overset{n}{\sum}}$R$_{lpkw}$R$_{jqiw})+(-18\nabla_{lk}^{2}$R$_{ipjq}%
+16\underset{w=\text{1}}{\overset{n}{\sum}}$R$_{lpkw}$R$_{iqjw})+(-18\nabla
_{lj}^{2}$R$_{kpiq}+16\underset{w=\text{1}}{\overset{n}{\sum}}$R$_{lpjw}%
$R$_{kqiw})$

$+(-18\nabla_{li}^{2}$R$_{kpjq}+16\underset{w=\text{1}}{\overset{n}{\sum}}%
$R$_{lpiw}$R$_{kqjw})+(-18\nabla_{lj}^{2}$R$_{ipkq}+16\underset{w=\text{1}%
}{\overset{n}{\sum}}$R$_{lpjw}$R$_{iqkw})+(-18\nabla_{li}^{2}$R$_{jpkq}%
+16\underset{w=\text{1}}{\overset{n}{\sum}}$R$_{lpiw}$R$_{jqkw})$

$+(-18\nabla_{kl}^{2}$R$_{jpiq}+16\underset{w=\text{1}}{\overset{n}{\sum}}%
$R$_{kplw}$R$_{jqiw})+(-18\nabla_{kl}^{2}$R$_{ipjq}+16\underset{w=\text{1}%
}{\overset{n}{\sum}}$R$_{kplw}$R$_{iqjw})+(-18\nabla_{jl}^{2}$R$_{kpiq}%
+16\underset{w=\text{1}}{\overset{n}{\sum}}$R$_{jplw}$R$_{kqiw})$

$+(-18\nabla_{il}^{2}$R$_{kpjq}+16\underset{w=\text{1}}{\overset{n}{\sum}}%
$R$_{iplw}$R$_{kqjw})+(-18\nabla_{jl}^{2}$R$_{ipkq}+16\underset{w=\text{1}%
}{\overset{n}{\sum}}$R$_{jplw}$R$_{iqkw})+(-18\nabla_{il}^{2}$R$_{jpkq}%
+16\underset{w=\text{1}}{\overset{n}{\sum}}$R$_{iplw}$R$_{jqkw})$

$+(-18\nabla_{kj}^{2}$R$_{lpiq}+16\underset{w=\text{1}}{\overset{n}{\sum}}%
$R$_{kpjw}$R$_{lqiw})+(-18\nabla_{ki}^{2}$R$_{lpjq}+16\underset{w=\text{1}%
}{\overset{n}{\sum}}$R$_{kpiw}$R$_{lqjw})+(-18\nabla_{jk}^{2}$R$_{lpiq}%
+16\underset{w=\text{1}}{\overset{n}{\sum}}$R$_{jpkw}$R$_{lqiw})$

$+(-18\nabla_{ik}^{2}$R$_{lpjq}+16\underset{w=\text{1}}{\overset{n}{\sum}}%
$R$_{ipkw}$R$_{lqjw})+(-18\nabla_{ji}^{2}$R$_{lpkq}+16\underset{w=\text{1}%
}{\overset{n}{\sum}}$R$_{jpiw}$R$_{lqkw})+(-18\nabla_{ij}^{2}$R$_{lpkq}%
+16\underset{w=\text{1}}{\overset{n}{\sum}}$R$_{ipjw}$R$_{lqkw}$

$+(-18\nabla_{kj}^{2}$R$_{iplq}+16\underset{w=\text{1}}{\overset{n}{\sum}}%
$R$_{kpjw}$R$_{iqlw})+(-18\nabla_{ki}^{2}$R$_{jplq}+16\underset{w=\text{1}%
}{\overset{n}{\sum}}$R$_{kpiw}$R$_{jqlw})+(-18\nabla_{jk}^{2}$R$_{iplq}%
+16\underset{w=\text{1}}{\overset{n}{\sum}}$R$_{jpkw}$R$_{iqlw})$

$\bigskip+(-18\nabla_{ik}^{2}$R$_{jplq}+16\underset{w=\text{1}%
}{\overset{n}{\sum}}$R$_{ipkw}$R$_{jqlw})+(-18\nabla_{ji}^{2}$R$_{kplq}%
+16\underset{w=\text{1}}{\overset{n}{\sum}}$R$_{jpiw}$R$_{kqlw})+(-18\nabla
_{ij}^{2}$R$_{kplq}+16\underset{w=\text{1}}{\overset{n}{\sum}}$R$_{ipjw}%
$R$_{kqlw})$
\end{center}

(vi) We take $l=k$ and $p=j$ and have:

$\qquad\frac{\partial^{4}g_{jq}}{\partial x_{i}\partial x_{j}\partial
x_{k}^{2}}(y_{0})$

$=\frac{1}{360}[(-18\nabla_{kk}^{2}R_{jjiq}+16\underset{w=\text{1}%
}{\overset{n}{\sum}}$R$_{kjkw}$R$_{jqiw})+(-18\nabla_{kk}^{2}$R$_{ijjq}%
+16\underset{w=\text{1}}{\overset{n}{\sum}}$R$_{kjkw}$R$_{iqjw})+(-18\nabla
_{kj}^{2}$R$_{kjiq}+16\underset{w=\text{1}}{\overset{n}{\sum}}$R$_{kjjw}%
$R$_{kqiw})$

$+(-18\nabla_{ki}^{2}$R$_{kjjq}+16\underset{w=\text{1}}{\overset{n}{\sum}}%
$R$_{kjiw}$R$_{kqjw})+(-18\nabla_{kj}^{2}$R$_{ijkq}+16\underset{w=\text{1}%
}{\overset{n}{\sum}}$R$_{kjjw}$R$_{iqkw})+(-18\nabla_{ki}^{2}$R$_{jjkq}%
+16\underset{w=\text{1}}{\overset{n}{\sum}}$R$_{kjiw}$R$_{jqkw})$

$+(-18\nabla_{kk}^{2}$R$_{jjiq}+16\underset{w=\text{1}}{\overset{n}{\sum}}%
$R$_{kjkw}$R$_{jqiw})+(-18\nabla_{kk}^{2}$R$_{ijjq}+16\underset{w=\text{1}%
}{\overset{n}{\sum}}$R$_{kjkw}$R$_{iqjw})+(-18\nabla_{jk}^{2}$R$_{kjiq}%
+16\underset{w=\text{1}}{\overset{n}{\sum}}$R$_{jjlw}$R$_{kqiw})$

$+(-18\nabla_{ik}^{2}$R$_{kjjq}+16\underset{w=\text{1}}{\overset{n}{\sum}}%
$R$_{ijkw}$R$_{kqjw})+(-18\nabla_{jk}^{2}$R$_{ijkq}+16\underset{w=\text{1}%
}{\overset{n}{\sum}}$R$_{jjkw}$R$_{iqkw})+(-18\nabla_{ik}^{2}$R$_{jjkq}%
+16\underset{w=\text{1}}{\overset{n}{\sum}}$R$_{ijkw}$R$_{jqkw})$

$+(-18\nabla_{kj}^{2}$R$_{kjiq}+16\underset{w=\text{1}}{\overset{n}{\sum}}%
$R$_{kjjw}$R$_{kqiw})+(-18\nabla_{ki}^{2}$R$_{kjjq}+16\underset{w=\text{1}%
}{\overset{n}{\sum}}$R$_{kjiw}$R$_{kqjw})+(-18\nabla_{jk}^{2}$R$_{kjiq}%
+16\underset{w=\text{1}}{\overset{n}{\sum}}$R$_{jjkw}$R$_{kqiw})$

$+(-18\nabla_{ik}^{2}$R$_{kjjq}+16\underset{w=\text{1}}{\overset{n}{\sum}}%
$R$_{ijkw}$R$_{kqjw})+(-18\nabla_{ji}^{2}$R$_{kjkq}+16\underset{w=\text{1}%
}{\overset{n}{\sum}}$R$_{jjiw}$R$_{kqkw})+(-18\nabla_{ij}^{2}$R$_{kjkq}%
+16\underset{w=\text{1}}{\overset{n}{\sum}}$R$_{ijjw}$R$_{kqkw})$

$+(-18\nabla_{kj}^{2}$R$_{ijkq}+16\underset{w=\text{1}}{\overset{n}{\sum}}%
$R$_{kjjw}$R$_{iqkw})+(-18\nabla_{ki}^{2}$R$_{jjkq}+16\underset{w=\text{1}%
}{\overset{n}{\sum}}$R$_{kjiw}$R$_{jqkw})+(-18\nabla_{jk}^{2}$R$_{ijkq}%
+16\underset{w=\text{1}}{\overset{n}{\sum}}$R$_{jjkw}$R$_{iqkw})$

$+(-18\nabla_{ik}^{2}$R$_{jjkq}+16\underset{w=\text{1}}{\overset{n}{\sum}}%
$R$_{ijkw}$R$_{jqkw})+(-18\nabla_{ji}^{2}$R$_{kjkq}+16\underset{w=\text{1}%
}{\overset{n}{\sum}}$R$_{jjiw}$R$_{kqkw})+(-18\nabla_{ij}^{2}$R$_{kjkq}%
+16\underset{w=\text{1}}{\overset{n}{\sum}}$R$_{ijjw}$R$_{kqkw})$

The last expression above can be simplified. We have, for example:
R$_{jjiq}=0$ and $\underset{k=\text{1}}{\overset{n}{\sum}}$R$_{kjkq}%
=\varrho_{jq}:$

$\frac{\partial^{4}g_{jq}}{\partial x_{i}\partial x_{j}\partial x_{k}^{2}%
}(y_{0})$

$=\frac{1}{360}[(16\underset{w=\text{1}}{\overset{n}{\sum}}\varrho_{jw}%
$R$_{jqiw})+(18\nabla_{kk}^{2}\varrho_{iq}+16\underset{w=\text{1}%
}{\overset{n}{\sum}}\varrho_{jw}$R$_{iqjw})-(18\nabla_{kj}^{2}$R$_{kjiq}%
+16\underset{w=\text{1}}{\overset{n}{\sum}}$R$_{kw}$R$_{kqiw})$

$+(18\nabla_{ki}^{2}\varrho_{kq}+16\underset{w=\text{1}}{\overset{n}{\sum}}%
$R$_{kjiw}$R$_{kqjw})+(-18\nabla_{kj}^{2}$R$_{ijkq}-16\underset{w=\text{1}%
}{\overset{n}{\sum}}\varrho_{kw}$R$_{iqkw})+(16\underset{w=\text{1}%
}{\overset{n}{\sum}}$R$_{kjiw}$R$_{jqkw})$

$+(16\underset{w=\text{1}}{\overset{n}{\sum}}\varrho_{jw}$R$_{jqiw}%
)+(18\nabla_{kk}^{2}\varrho_{iq}+16\underset{w=\text{1}}{\overset{n}{\sum}%
}\varrho_{jw}$R$_{iqjw})-18\nabla_{jk}^{2}$R$_{kjiq}$

$+(18\nabla_{ik}^{2}\varrho_{kq}+16\underset{w=\text{1}}{\overset{n}{\sum}}%
$R$_{ijkw}$R$_{kqjw})-18\nabla_{jk}^{2}$R$_{ijkq}+16\underset{w=\text{1}%
}{\overset{n}{\sum}}$R$_{ijkw}$R$_{jqkw}$

$-(18\nabla_{kj}^{2}$R$_{kjiq}+16\underset{w=\text{1}}{\overset{n}{\sum}%
}\varrho_{kw}$R$_{kqiw})+(18\nabla_{ki}^{2}\varrho_{kq}+16\underset{w=\text{1}%
}{\overset{n}{\sum}}$R$_{kjiw}$R$_{kqjw})+18\nabla_{jk}^{2}$R$_{jkiq}$

$+(18\nabla_{ik}^{2}\varrho_{kq}+16\underset{w=\text{1}}{\overset{n}{\sum}}%
$R$_{ijkw}$R$_{kqjw})-18\nabla_{ji}^{2}\varrho_{jq}-(18\nabla_{ij}^{2}%
\varrho_{jq}+16\underset{w=\text{1}}{\overset{n}{\sum}}\varrho_{iw}$%
R$_{kqkw})$

$-(18\nabla_{kj}^{2}$R$_{ijkq}+16\underset{w=\text{1}}{\overset{n}{\sum}%
}\varrho_{kw}$R$_{iqkw})+(16\underset{w=\text{1}}{\overset{n}{\sum}}$%
R$_{kjiw}$R$_{jqkw})-18\nabla_{jk}^{2}$R$_{ijkq}$

$+(16\underset{w=\text{1}}{\overset{n}{\sum}}$R$_{ijkw}$R$_{jqkw}%
)-18\nabla_{ji}^{2}\varrho_{jq}-(18\nabla_{ij}^{2}\varrho_{jq}%
+16\underset{w=\text{1}}{\overset{n}{\sum}}\varrho_{iw}\varrho_{qw})](y_{0})$

\bigskip

We set $p=j$ and $l=i$ and have:

$\frac{\partial^{4}\text{g}_{jq}}{\partial x_{i}^{2}\partial x_{j}\partial
x_{k}}(y_{0})$

$=\frac{1}{360}[(16\underset{w=\text{1}}{\overset{n}{\sum}}\varrho_{jw}%
$R$_{jqkw})+(18\nabla_{ii}^{2}\varrho_{kq}+16\underset{w=\text{1}%
}{\overset{n}{\sum}}\varrho_{jw}$R$_{kqjw})-(18\nabla_{ij}^{2}$R$_{ijkq}%
+16\underset{w=\text{1}}{\overset{n}{\sum}}$R$_{iw}$R$_{iqkw})$

$+(18\nabla_{ik}^{2}\varrho_{iq}+16\underset{w=\text{1}}{\overset{n}{\sum}}%
$R$_{ijkw}$R$_{iqjw})+(-18\nabla_{ij}^{2}$R$_{kjiq}-16\underset{w=\text{1}%
}{\overset{n}{\sum}}\varrho_{iw}$R$_{kqiw})+(16\underset{w=\text{1}%
}{\overset{n}{\sum}}$R$_{ijkw}$R$_{jqiw})$

$+(16\underset{w=\text{1}}{\overset{n}{\sum}}\varrho_{jw}$R$_{jqkw}%
)+(18\nabla_{ii}^{2}\varrho_{kq}+16\underset{w=\text{1}}{\overset{n}{\sum}%
}\varrho_{jw}$R$_{kqjw})-18\nabla_{ji}^{2}$R$_{ijkq}$

$+(18\nabla_{ki}^{2}\varrho_{iq}+16\underset{w=\text{1}}{\overset{n}{\sum}}%
$R$_{kjiw}$R$_{iqjw})-18\nabla_{ji}^{2}$R$_{kjiq}+16\underset{w=\text{1}%
}{\overset{n}{\sum}}$R$_{kjiw}$R$_{jqiw}$

$-(18\nabla_{ij}^{2}$R$_{ijkq}+16\underset{w=\text{1}}{\overset{n}{\sum}%
}\varrho_{iw}$R$_{iqkw})+(18\nabla_{ik}^{2}\varrho_{iq}+16\underset{w=\text{1}%
}{\overset{n}{\sum}}$R$_{ijkw}$R$_{iqjw})+18\nabla_{ji}^{2}$R$_{jikq}$

$+(18\nabla_{ii}^{2}\varrho_{iq}+16\underset{w=\text{1}}{\overset{n}{\sum}}%
$R$_{kjiw}$R$_{iqjw})-18\nabla_{jk}^{2}\varrho_{jq}-(18\nabla_{kj}^{2}%
\varrho_{jq}+16\underset{w=\text{1}}{\overset{n}{\sum}}\varrho_{kw}$%
R$_{iqiw})$

$-(18\nabla_{ij}^{2}$R$_{kjiq}+16\underset{w=\text{1}}{\overset{n}{\sum}%
}\varrho_{iw}$R$_{kqiw})+(16\underset{w=\text{1}}{\overset{n}{\sum}}$%
R$_{ijkw}$R$_{jqiw})-18\nabla_{ji}^{2}$R$_{kjiq}$

$+(16\underset{w=\text{1}}{\overset{n}{\sum}}$R$_{kjiw}$R$_{jqkw}%
)-18\nabla_{jk}^{2}\varrho_{jq}-(18\nabla_{kj}^{2}\varrho_{jq}%
+16\underset{w=\text{1}}{\overset{n}{\sum}}\varrho_{kw}\varrho_{qw})](y_{0})$

\qquad\qquad\qquad\qquad\qquad\qquad\qquad\qquad\qquad\qquad\qquad\qquad
\qquad\qquad\qquad\qquad\qquad$\blacksquare$

\section{\textbf{Table A}$_{2}$}

For $i,j,k=q+1,...,n$

(i) g$^{jk}(y_{0})=\delta_{jk}$

(ii) $\frac{\partial\text{g}^{jk}}{\partial\text{x}_{i}}(y_{0})=0$

(iii) $\frac{\partial^{2}\text{g}^{kl}}{\partial\text{x}_{i}\partial
\text{x}_{j}}(y_{0})=\frac{1}{3}(R_{ikjl}+R_{jkil})(y_{0})$

\qquad In particular, in normal coordinates, we have:

\qquad$\frac{\partial^{2}\text{g}^{kl}}{\partial\text{x}_{i}^{2}}(y_{0}%
)=\frac{1}{3}(R_{ikil}+R_{ikil})(y_{0})=\frac{2}{3}(R_{ikil}(y_{0})=\frac
{2}{3}\varrho_{kl}(y_{0})$

(iv) $\frac{\partial^{3}g^{pq}}{\partial x_{i}\partial x_{j}\partial x_{k}%
}(y_{0})$

$=-\frac{1}{6}\underset{r,s,t=1}{\overset{n}{\sum}}\nabla_{_{r}}R_{sptq}%
(y_{0})$

$\times\lbrack\delta_{kr}\delta_{js}\delta_{it}+\delta_{kr}\delta_{is}%
\delta_{jt}+\delta_{ks}\delta_{jr}\delta_{it}+\delta_{ks}\delta_{ir}%
\delta_{jt}+\delta_{kt}\delta_{jr}\delta_{is}+\delta_{kt}\delta_{ir}%
\delta_{js}]$

$=-\frac{1}{6}[\nabla_{_{k}}R_{jkiq}+\nabla_{_{k}}R_{ikjq}+\nabla_{_{j}%
}R_{kpiq}+\nabla_{_{i}}R_{kpjq}+\nabla_{_{j}}R_{ipkq}+\nabla_{_{i}}%
R_{jpkq}](y_{0})$

$-\frac{1}{6}[\nabla_{_{k}}R_{jkim}+\nabla_{_{k}}R_{ikjq}+\nabla_{_{j}%
}R_{kpiq}+\nabla_{_{i}}R_{kpjq}+\nabla_{_{j}}R_{ipkq}+\nabla_{_{i}}%
R_{jpkq}](y_{0})$

(v) In particular,

$\qquad\frac{\partial^{3}\text{g}^{pq}}{\partial\text{x}_{i}^{2}%
\partial\text{x}_{j}}(y_{0})=$ $\frac{1}{3}\nabla_{j}$R$_{ipiq}(y_{0}%
)+\frac{1}{3}\nabla_{i}$R$_{jpiq}(y_{0})+\frac{1}{3}\nabla_{i}$R$_{ipjq}%
(y_{0})$

(vi) $\frac{\partial^{3}\text{g}^{pq}}{\partial\text{x}_{i}\partial
\text{x}_{j}^{2}}(y_{0})=$ $\frac{1}{3}\nabla_{i}$R$_{jpjq}(y_{0})+\frac{1}%
{3}\nabla_{j}$R$_{ipjq}(y_{0})+\frac{1}{3}\nabla_{j}$R$_{jpiq}(y_{0})$

\subsection{\protect\underline{\textbf{Computations}}}

For the computations below, we use the expansion of g$^{\alpha\beta}(x)$ given
in Proposition

For $k,l=q+1,...,n,$ we have:

g$^{kl}(x_{0})=$ $\delta^{kl}+\frac{1}{3}\underset{r,s=q+1}{\overset{n}{\sum}%
}$R$_{\text{r}k\text{s}l}($y$_{0})$x$_{\text{r}}$x$_{\text{s}}+\frac{1}%
{6}\underset{r,s,t=q+1}{\overset{n}{\sum}}\nabla_{\text{r}}$R$_{\text{s}%
k\text{t}l}($y$_{0})x_{\text{r}}x_{\text{s}}x_{\text{t}}$

$\qquad\ \ \ -\frac{1}{360}\overset{n}{\underset{r,s,t,u=q+1}{\sum}}%
(-18\nabla_{\text{rs}}^{2}$R$_{\text{t}k\text{u}l}%
+16\underset{p=q+1}{\overset{n}{\sum}}$R$_{rksp}$R$_{tlup})($y$_{0}%
)x_{\text{r}}x_{\text{s}}x_{\text{t}}x_{\text{u}}\qquad$

(i) and (ii) are obvious.

(iii) $\frac{\partial^{2}\text{g}^{kl}}{\partial\text{x}_{i}\partial
\text{x}_{j}}(x_{0})=\frac{1}{3}\underset{r,s=q+1}{\overset{n}{\sum}}%
$R$_{rksl}$(y$_{0}$)$(\frac{\partial\text{x}_{\text{r}}}{\partial\text{x}_{i}%
}\frac{\partial\text{x}_{\text{s}_{{}}}}{\partial\text{x}_{j}}+\frac
{\partial\text{x}_{\text{r}}}{\partial\text{x}_{j}}\frac{\partial
\text{x}_{\text{s}}}{\partial\text{x}_{i}})+O($x$_{0})$

\qquad\qquad\qquad$\ =\frac{1}{3}\underset{r,s=q+1}{\overset{\text{n}}{\sum}}%
$R$_{rksl}$(y$_{0}$)$(\delta_{ri}\delta_{sj}+\delta_{rj}\delta_{si})+O($%
x$_{0})$

\ \ \ \ $\frac{\partial^{2}\text{g}^{kl}}{\partial\text{x}_{i}\partial
\text{x}_{j}}(y_{0})=\frac{1}{3}$R$_{ikjl}(y_{0})$ $+\frac{1}{3}$%
R$_{jkil}(y_{0})=\frac{1}{3}[$R$_{ikjl}+$ R$_{jkil}](y_{0})$

(iii) We take $\delta=\gamma$ in (iv).

(iv) \ $\frac{\partial^{3}\text{g}^{kl}}{\partial\text{x}_{i}^{2}%
\partial\text{x}_{j}}(x_{0})=$ $\frac{1}{6}%
\underset{r,s,t=q+1}{\overset{n}{\sum}}\nabla_{\text{r}}$R$_{\text{s}%
k\text{t}l}(y_{0})\frac{\partial^{3}}{\partial\text{x}_{i}^{2}\partial
\text{x}_{j}}(x_{\text{r}}x_{\text{s}}x_{\text{t}})+O(x_{0})$

\qquad$\frac{\partial^{3}}{\partial\text{x}_{i}^{2}\partial\text{x}_{j}%
}(x_{\text{r}}x_{\text{s}}x_{\text{t}})=\frac{\partial^{2}}{\partial
\text{x}_{i}^{2}}(\delta_{rj}x_{\text{s}}x_{\text{t}}+x_{\text{r}}%
\frac{\partial}{\partial\text{x}_{j}}(x_{\text{s}}x_{\text{t}})$

$\qquad=\frac{\partial^{2}}{\partial\text{x}_{i}^{2}}(\delta_{rj}x_{\text{s}%
}x_{\text{t}}+x_{\text{r}}\delta_{sj}x_{\text{t}}+x_{\text{r}}x_{\text{s}%
}\delta_{tj})$

\qquad$=\frac{\partial}{\partial\text{x}_{i}}(\delta_{rj}\delta_{si}%
x_{\text{t}}+\delta_{rj}x_{\text{s}}\delta_{it}+\delta_{ri}\delta
_{sj}x_{\text{t}}+x_{\text{r}}\delta_{sj}\delta_{ti}+\delta_{ir}x_{\text{s}%
}\delta_{tj}+x_{\text{r}}\delta_{si}\delta_{tj})$

\qquad$=(\delta_{rj}\delta_{si}\delta_{ti}+\delta_{rj}\delta_{si}\delta
_{it}+\delta_{ri}\delta_{sj}\delta_{ti}+\delta_{ri}\delta_{sj}\delta
_{ti}+\delta_{ir}\delta_{si}\delta_{tj}+\delta_{ri}\delta_{si}\delta_{tj})$

Therefore,

(v) $\frac{\partial^{3}\text{g}^{kl}}{\partial\text{x}_{i}^{2}\partial
\text{x}_{j}}(x_{0})=$ $\frac{1}{6}\underset{r,s,t=q+1}{\overset{n}{\sum}%
}\nabla_{\text{r}}$R$_{\text{s}k\text{t}l}(y_{0})(2\delta_{rj}\delta
_{si}\delta_{ti}+2\delta_{ri}\delta_{sj}\delta_{ti}+2\delta_{ri}\delta
_{si}\delta_{tj})+O(x_{0})$

$\qquad\frac{\partial^{3}\text{g}^{kl}}{\partial\text{x}_{i}^{2}%
\partial\text{x}_{j}}(x_{0})=$ $\frac{1}{3}%
\underset{r,s,t=q+1}{\overset{n}{\sum}}\nabla_{r}$R$_{sktl}(y_{0})(\delta
_{rj}\delta_{si}\delta_{ti}+\delta_{ri}\delta_{sj}\delta_{ti}+\delta
_{ri}\delta_{si}\delta_{tj})+O(x_{0})$

$\qquad\frac{\partial^{3}\text{g}^{kl}}{\partial\text{x}_{i}^{2}%
\partial\text{x}_{j}}(x_{0})\qquad=$ $\frac{1}{3}\nabla_{j}$R$_{ikil}%
(y_{0})+\frac{1}{3}\nabla_{i}$R$_{jkil}(y_{0})+\frac{1}{3}\nabla_{i}$%
R$_{ikjl}(y_{0})+O(x_{0})$

In particular,

$\frac{\partial^{3}\text{g}^{kl}}{\partial\text{x}_{i}^{2}\partial\text{x}%
_{j}}(y_{0})=$ $\frac{1}{3}\nabla_{j}$R$_{ikil}(y_{0})+\frac{1}{3}\nabla_{i}%
$R$_{jkil}(y_{0})+\frac{1}{3}\nabla_{i}$R$_{ikjl}(y_{0})$

\section{Table A$_{3}$}

The expansion of g$_{\text{a}\alpha}($x) in Fermi coordinates is given as follows:

g$_{\text{a}\alpha}($x) = \ $-\overset{n}{\underset{i=q+1}{\sum(}}%
\perp_{\text{a}\alpha i})(y_{0})$x$_{i}-\frac{4}{3}%
\underset{i,j=q+1}{\overset{\text{n}}{\sum\text{ }}}$R$_{i\text{a}j\alpha
}(y_{0})$x$_{i}$x$_{j}$

We re-label it as follows:

g$_{\text{a}k}($x) = \ $-\overset{n}{\underset{l=q+1}{\sum(}}\perp
_{\text{a}kl})(y_{0})$x$_{l}-\frac{4}{3}\underset{l,m=q+1}{\overset{\text{n}%
}{\sum\text{ }}}$R$_{l\text{a}mk}(y_{0})$x$_{l}$x$_{m}$

\vspace{1pt}For a=1,...,q \ and $r=q+1,...,n$ we have:

g$_{\text{a}k}(y_{0})=0$

(i) g$_{\text{a}r}(y_{0})=0$

$\frac{\partial\text{g}_{\text{a}k}}{\partial\text{x}_{i}}(y_{0}%
)=-\overset{n}{\underset{l=q+1}{\sum(}}\perp_{\text{a}kl})(y_{0})\delta
_{il}=-\perp_{\text{a}ki}(y_{0})=\perp_{\text{a}ik}(y_{0})$

(ii) $\frac{\partial\text{g}_{\text{a}r}}{\partial\text{x}_{\alpha}}(y_{0})=$
$\perp_{\text{a}\alpha r}(y_{0})=-\perp_{\text{a}r\alpha}(y_{0})$

\qquad$\frac{\partial\text{g}_{\text{a}k}}{\partial\text{x}_{i}}%
(y_{0})=-\overset{n}{\underset{l=q+1}{\sum(}}\perp_{\text{a}kl})(y_{0}%
)\delta_{il}=-\perp_{\text{a}ki}(y_{0})=\perp_{\text{a}ik}(y_{0})$

(iii) $\frac{\partial^{2}\text{g}_{\text{a}r}}{\partial\text{x}_{\alpha
}\partial\text{x}_{\beta}}(y_{0})=-\frac{4}{3}(R_{\alpha\text{a}\beta
r}+R_{\beta\text{a}\alpha r})(y_{0})$

\qquad$\frac{\partial^{2}\text{g}_{\text{a}k}}{\partial\text{x}_{i}%
\partial\text{x}_{j}}(y_{0})=-\frac{4}{3}(R_{i\text{a}jk}+R_{j\text{a}%
ik})(y_{0})$

\qquad

In particular,

(iii)$^{\ast}$ $\frac{\partial^{2}\text{g}_{\text{a}r}}{\partial
\text{x}_{\alpha}^{2}}(y_{0})=-\frac{4}{3}(R_{\alpha\text{a}\alpha
r}+R_{\alpha\text{a}\alpha r})(y_{0})=-\frac{8}{3}R_{\alpha\text{a}\alpha
r}(y_{0})$

(iv) $\frac{\partial^{3}\text{g}_{\text{a}r}}{\partial\text{x}_{\alpha}%
^{2}\partial\text{x}_{\beta}}(y_{0})=-\frac{1}{6}\underset{\text{i,j,k=q+1}%
}{\overset{\text{n}}{\sum}}\{(3\nabla_{\alpha}$R$_{\alpha\text{a}\beta r}%
+4$R$_{\alpha r\alpha\text{T}_{\text{a}\beta}}+$ $4$R$_{\alpha r\alpha
\perp_{\text{a}\beta}})$

$\qquad\qquad\qquad\qquad+(3\nabla_{\alpha}$R$_{\beta\text{a}\alpha r}%
+4$R$_{\alpha\text{r}\beta\text{T}_{\text{a}\alpha}}+$ $4$R$_{\alpha
\text{r}\beta\perp_{\text{a}\alpha}})$

\qquad\qquad\qquad\qquad\ \ +$(3\nabla_{\beta}$R$_{\alpha\text{a}\alpha r}%
+4$R$_{\beta\text{r}\alpha\text{T}_{\text{a}\alpha}}+$ $4$R$_{\beta
\text{r}\alpha\perp_{\text{a}\alpha}})\}(y_{0})$

\subsection{\protect\underline{\textbf{Computations}}}

We use \textbf{Proposition 8}:

(i) immediate

(ii) $\frac{\partial\text{g}_{\text{a}r}}{\partial\text{x}_{\alpha}}%
(y_{0})=-\underset{i=q+1}{\overset{n}{\sum}}\perp_{\text{a}ri}(y_{0}%
)\delta_{\alpha i}=$ $-\perp_{\text{a}r\alpha}(y_{0})=$ $\perp_{\text{a}\alpha
r}(y_{0})$

(iii) $\frac{\partial^{2}\text{g}_{\text{a}r}}{\partial\text{x}_{\alpha
}\partial\text{x}_{\beta}}(y_{0})=-\frac{4}{3}\underset{\text{i,j}%
=q+1}{\overset{\text{n}}{\sum\text{ }}}$R$_{i\text{a}jr}$($y_{0}%
)(\delta_{\alpha i}\delta_{\beta j}+\delta_{\beta i}\delta_{\alpha j})$

\vspace{1pt}\qquad\qquad\qquad\ \ \ \ $=-\frac{4}{3}($R$_{\alpha\text{a}\beta
r}+$ R$_{\beta\text{a}\alpha r})(y_{0})$

In particular,

\qquad$\frac{\partial^{2}\text{g}_{\text{a}r}}{\partial\text{x}_{\alpha}^{2}%
}(y_{0})=-\frac{4}{3}($R$_{\alpha\text{a}\alpha\text{r}}+$ R$_{\alpha
\text{a}\alpha\text{r}})(y_{0})=-\frac{8}{3}$R$_{\alpha\text{a}\alpha\text{r}%
}(y_{0})$

(iv) $\frac{\partial^{3}\text{g}_{\text{a}r}}{\partial\text{x}_{\alpha}%
^{2}\partial\text{x}_{\beta}}(y_{0})=-\frac{1}{6}\underset{\text{i,j,k=q+1}%
}{\overset{\text{n}}{\sum}}\{\frac{3}{2}\nabla_{\text{i}}$R$_{\text{jakr}}%
+2$R$_{\text{irjT}_{\text{ak}}}+$ $2$R$_{\text{irj}\perp_{\text{ak}}}%
\}(y_{0})\frac{\partial^{3}}{\partial\text{x}_{\alpha}^{2}\partial
\text{x}_{\beta}}($x$_{\text{i}}$x$_{\text{j}}$x$_{\text{k}})$

\qquad\qquad\qquad\qquad

\vspace{1pt}\qquad Now,

$\qquad\frac{\partial^{3}}{\partial\text{x}_{\alpha}^{2}\partial
\text{x}_{\beta}}($x$_{\text{i}}$x$_{\text{j}}$x$_{\text{k}})=\frac
{\partial^{2}}{\partial\text{x}_{\alpha}^{2}}($x$_{\text{i}}$x$_{\text{j}%
}\delta_{\beta k}+$x$_{\text{i}}$x$_{\text{k}}\delta_{\beta j}+$x$_{\text{j}}%
$x$_{\text{k}}\delta_{\beta i})$

\qquad\qquad\qquad\qquad\qquad\ \ \ =2($\delta_{\alpha i}\delta_{\alpha
j}\delta_{\beta k}+\delta_{\alpha i}\delta_{\alpha k}\delta_{\beta j}%
+\delta_{\alpha j}\delta_{\alpha k}\delta_{\beta i})$

\qquad Therefore,

\qquad\ $\frac{\partial^{3}\text{g}_{\text{a}r}}{\partial\text{x}_{\alpha}%
^{2}\partial\text{x}_{\beta}}(y_{0})=-\frac{1}{3}\underset{\text{i,j,k=q+1}%
}{\overset{\text{n}}{\sum}}\{\frac{3}{2}\nabla_{\text{i}}$R$_{\text{jakr}}%
+2$R$_{\text{irjT}_{\text{ak}}}+$ $2$R$_{\text{irj}\perp_{\text{ak}}}%
\}(y_{0})$($\delta_{\alpha i}\delta_{\alpha j}\delta_{\beta k}$

$\qquad\qquad\qquad\qquad\qquad+\delta_{\alpha i}\delta_{\alpha k}%
\delta_{\beta j}+\delta_{\alpha j}\delta_{\alpha k}\delta_{\beta i})$

\qquad

\vspace{1pt}\qquad\qquad\qquad$\ \ =-\frac{1}{3}\underset{\text{i,j,k=q+1}%
}{\overset{\text{n}}{\sum}}[\{\frac{3}{2}\nabla_{\alpha}$R$_{\alpha
\text{a}\beta\text{r}}+2$R$_{\alpha\text{r}\alpha\text{T}_{\text{a}\beta}}+$
$2$R$_{\alpha\text{r}\alpha\perp_{\text{a}\beta}}\}+\{\frac{3}{2}%
\nabla_{\alpha}$R$_{\beta\text{a}\alpha\text{r}}$

$\qquad\qquad\qquad+2$R$_{\alpha\text{r}\beta\text{T}_{\text{a}\alpha}}+$
$2$R$_{\alpha\text{r}\beta\perp_{\text{a}\alpha}}\}$ +$\{\frac{3}{2}%
\nabla_{\beta}$R$_{\alpha\text{a}\alpha\text{r}}+2$R$_{\beta\text{r}%
\alpha\text{T}_{\text{a}\alpha}}+$ $2$R$_{\beta\text{r}\alpha\perp
_{\text{a}\alpha}}\}](y_{0})$

\section{Table A$_{4}$}

For a = 1,...,q \ and \ $\alpha,r=q+1,...,n$

\vspace{1pt}(i) \ g$^{\text{a}r}(y_{0})=0$

(ii) $\frac{\partial\text{g}^{\text{a}r}}{\partial\text{x}_{\alpha}}%
(y_{0})=-\perp_{\text{a}\alpha r}(y_{0})=$ $\perp_{\text{a}r\alpha}(y_{0})$

(iii) $\frac{\partial^{2}\text{g}^{\text{a}r}}{\partial\text{x}_{i}^{2}}%
(y_{0})=\frac{8}{3}R_{i\text{a}ir}+4\underset{\text{b=1}}{\overset{\text{q}%
}{\sum}}T_{\text{ab}i}(y_{0})\perp_{\text{b}ri}(y_{0})$

(iv) $\frac{\partial^{2}g^{\text{a}r}}{\partial\text{x}_{i}\partial
\text{x}_{j}}(y_{0})=\frac{4}{3}(R_{i\text{a}jr}+R_{j\text{a}ir}%
)(y_{0})+4\underset{\text{b=1}}{\overset{\text{q}}{\sum}}T_{\text{ab}i}%
(y_{0})\perp_{\text{b}rj}(y_{0})$

\subsection{\protect\underline{\textbf{Computations}}}

\vspace{1pt}(i) \ g$^{\text{a}r}(y_{0})=\delta^{\text{a}r}=0$ for a = 1,...,q
and $r=q+1,...,n.$

(ii) For a,b =1,...,q ; $j,r=q+1,...,n$ and $\beta=1,...,n,$ we have, with
summation over $\beta,$b and $j$ understood:

$\qquad0=\delta_{\text{a}}^{r}=$ g$_{\text{a}\beta}$g$^{\beta r}=$
g$_{\text{ab}}$g$^{\text{b}r}+$ g$_{\text{a}j}$g$^{jr}$

Hence,

\qquad$0=\frac{\partial\text{g}_{\text{ab}}}{\partial\text{x}_{\alpha}}%
(y_{0})\delta^{\text{b}r}+\delta_{\text{ab}}\frac{\partial\text{g}^{\text{b}%
r}}{\partial\text{x}_{\alpha}}(y_{0})+\frac{\partial\text{g}_{\text{a}j}%
}{\partial\text{x}_{\alpha}}(y_{0})\delta^{j\text{r}}+\delta_{\text{a}j}%
\frac{\partial\text{g}^{jr}}{\partial\text{x}_{\alpha}}(y_{0})$

\qquad$0=\frac{\partial\text{g}^{\text{a}r}}{\partial\text{x}_{\alpha}}%
(y_{0})+\frac{\partial\text{g}_{\text{a}r}}{\partial\text{x}_{\alpha}}(y_{0})$

\qquad(since $\delta^{\text{b}r}=0=\delta_{\text{a}j})$

\qquad Hence,

\qquad$\frac{\partial\text{g}^{\text{a}r}}{\partial\text{x}_{\alpha}}%
(y_{0})=-\frac{\partial\text{g}_{\text{a}r}}{\partial\text{x}_{\alpha}}%
(y_{0})=-\perp_{\text{a}\alpha r}(y_{0})$ \ by (ii) of Table 3.

(iii) Similarly, with summation over $\beta,$c and $j$ understood:

$\qquad0=\frac{\partial^{2}}{\partial\text{x}_{\alpha}^{2}}($g$_{\text{a}%
\beta}$g$^{\beta r})(y_{0})=\frac{\partial^{2}}{\partial\text{x}_{\alpha}^{2}%
}($g$_{\text{ac}}$g$^{\text{c}r})(y_{0})+\frac{\partial^{2}}{\partial
\text{x}_{\alpha}^{2}}($g$_{\text{a}j}$g$^{jr})$

$\qquad=\delta_{\text{ac}}$ $\frac{\partial^{2}\text{g}^{\text{c}r}}%
{\partial\text{x}_{\alpha}^{2}}(y_{0})+2\frac{\partial\text{g}_{\text{ac}}%
}{\partial\text{x}_{\alpha}}(y_{0})\frac{\partial\text{g}^{\text{c}r}%
}{\partial\text{x}_{\alpha}}(y_{0})+$ $\frac{\partial^{2}\text{g}_{\text{ac}}%
}{\partial\text{x}_{\alpha}^{2}}(y_{0})\delta^{\text{c}r}$

\qquad$\ \ \ +\delta_{\text{a}j}$ $\frac{\partial^{2}\text{g}^{jr}}%
{\partial\text{x}_{\alpha}^{2}}(y_{0})+2\frac{\partial\text{g}_{\text{a}j}%
}{\partial\text{x}_{\alpha}}(y_{0})\frac{\partial\text{g}^{jr}}{\partial
\text{x}_{\alpha}}(y_{0})+$ $\frac{\partial^{2}\text{g}_{\text{a}j}}%
{\partial\text{x}_{\alpha}^{2}}(y_{0})\delta^{jr}$

\qquad$=\frac{\partial^{2}\text{g}^{\text{a}r}}{\partial\text{x}_{\alpha}^{2}%
}(y_{0})+2\frac{\partial\text{g}_{\text{ac}}}{\partial\text{x}_{\alpha}}%
(y_{0})\frac{\partial\text{g}^{\text{c}r}}{\partial\text{x}_{\alpha}}(y_{0})+$
$\frac{\partial^{2}\text{g}_{\text{a}r}}{\partial\text{x}_{\alpha}^{2}}%
(y_{0})$

\qquad\qquad(since $\delta_{\text{a}j}=0=\delta^{\text{c}r}$ and
$\frac{\partial\text{g}^{jr}}{\partial\text{x}_{\alpha}}(y_{0})=0$\ by (ii) of
Table 1 $)$

\qquad Therefore,

\ \ \ \ $\frac{\partial^{2}\text{g}^{\text{a}r}}{\partial\text{x}_{\alpha}%
^{2}}(y_{0})=-2\underset{c=1}{\overset{n}{\sum}}\frac{\partial\text{g}%
_{\text{ac}}}{\partial\text{x}_{\alpha}}(y_{0})\frac{\partial\text{g}%
^{\text{c}r}}{\partial\text{x}_{\alpha}}(y_{0})-\frac{\partial^{2}%
\text{g}_{\text{a}r}}{\partial\text{x}_{\alpha}^{2}}(y_{0})$

Sequentially, we have by (ii) of Table 5 below (which was proved using

only Table 3), (ii) of Table 4 above and by (iii) of Table 3:

$\frac{\partial^{2}\text{g}^{\text{a}r}}{\partial\text{x}_{\alpha}^{2}}%
(y_{0})=-2\underset{c=1}{\overset{q}{\sum}}(-2T_{\text{ac}\alpha}%
)(y_{0})(\perp_{\text{c}r\alpha})(y_{0})+\frac{4}{3}(R_{\alpha\text{a}\alpha
r}+R_{\alpha\text{a}\alpha r})(y_{0})$

\qquad\qquad$=$ $4\underset{\text{c=1}}{\overset{\text{q}}{\sum}}%
(T_{\text{ac}\alpha})(y_{0})(\perp_{\text{c}r\alpha})(y_{0})+\frac{8}%
{3}R_{\alpha\text{a}\alpha r}(y_{0})$\qquad

(iv) $0=\delta_{\text{a}r}=\underset{\beta=1}{\overset{n}{\sum}}$
g$_{\text{a}\beta}$g$^{\beta r}=\underset{b=1}{\overset{q}{\sum}}$
g$_{\text{ab}}$g$^{\text{b}r}+$ $\underset{\beta=q+1}{\overset{n}{\sum}}%
$g$_{\text{a}\beta}$g$^{\beta r}$

Therefore for a,b = 1,...,q and $i,j,r=q+1,...,n,$ we have:

$0=\underset{b=1}{\overset{q}{\sum}}$ $\frac{\partial^{2}}{\partial
\text{x}_{i}\partial\text{x}_{j}}($g$_{\text{ab}}$g$^{\text{b}r})+$
$\underset{\beta=q+1}{\overset{n}{\sum}}\frac{\partial^{2}}{\partial
\text{x}_{i}\partial\text{x}_{j}}($g$_{\text{a}\beta}$g$^{\beta r})$

$0=\underset{b=1}{\overset{q}{\sum}}$ $\frac{\partial^{2}g_{\text{ab}}%
}{\partial\text{x}_{i}\partial\text{x}_{j}}(y_{0})$g$^{\text{b}r}%
(y_{0})+2\underset{b=1}{\overset{q}{\sum}}\frac{\partial g_{\text{ab}}%
}{\partial x_{i}}(y_{0})\frac{\partial g^{\text{b}r}}{\partial x_{j}}%
(y_{0})+\underset{b=1}{\overset{q}{\sum}}g_{\text{ab}}(y_{0})\frac
{\partial^{2}g^{\text{b}r}}{\partial\text{x}_{i}\partial\text{x}_{j}}(y_{0})$

$+$ $\underset{\beta=q+1}{\overset{n}{\sum}}\frac{\partial^{2}g_{\text{a}%
\beta}}{\partial\text{x}_{i}\partial\text{x}_{j}}(y_{0})g^{\beta r}%
(y_{0})+2\underset{\beta=q+1}{\overset{n}{\sum}}\frac{\partial g_{\text{a}%
\beta}}{\partial x_{i}}(y_{0})\frac{\partial g^{\beta r}}{\partial x_{j}%
}(y_{0})+\underset{\beta=q+1}{\overset{n}{\sum}}g_{\text{a}\beta}(y_{0}%
)\frac{\partial^{2}g^{\beta r}}{\partial\text{x}_{i}\partial\text{x}_{j}%
}(y_{0})$

Since g$^{\text{b}r}(y_{0})=\delta_{\text{b}r}=0$ and $g_{\text{a}\beta}%
(y_{0})=\delta_{\text{a}\beta}=0,$ (ii) $\frac{\partial\text{g}^{\beta\gamma}%
}{\partial\text{x}_{\alpha}}(y_{0})=0$ by (ii) of Table A$_{2},$ we have:

$0=2\underset{b=1}{\overset{q}{\sum}}\frac{\partial g_{\text{ab}}}{\partial
x_{i}}(y_{0})\frac{\partial g^{\text{b}r}}{\partial x_{j}}(y_{0}%
)+\frac{\partial^{2}g^{\text{a}r}}{\partial\text{x}_{i}\partial\text{x}_{j}%
}(y_{0})+$ $\frac{\partial^{2}g_{\text{a}r}}{\partial\text{x}_{i}%
\partial\text{x}_{j}}(y_{0})$

We have:

$\frac{\partial^{2}g^{\text{a}r}}{\partial\text{x}_{i}\partial\text{x}_{j}%
}(y_{0})=-$ $\frac{\partial^{2}g_{\text{a}r}}{\partial\text{x}_{i}%
\partial\text{x}_{j}}(y_{0})-2\underset{\text{b=1}}{\overset{\text{q}}{\sum}%
}\frac{\partial g_{\text{ab}}}{\partial x_{i}}(y_{0})\frac{\partial
g^{\text{b}r}}{\partial x_{j}}(y_{0})$

Now by $\frac{\partial^{2}\text{g}_{\text{a}r}}{\partial\text{x}_{i}%
\partial\text{x}_{j}}(y_{0})=-\frac{4}{3}(R_{i\text{a}jr}+R_{j\text{a}%
ir})(y_{0})$ (iii) of \textbf{Table A}$_{3};\frac{\partial g_{\text{ab}}%
}{\partial x_{i}}(y_{0})=-2$T$_{\text{ab}i}(y_{0})$ by (ii) of Table A$_{5}$
below and $\frac{\partial\text{g}^{\text{b}r}}{\partial\text{x}_{j}}%
(y_{0})=-\perp_{\text{b}jr}(y_{0})=$ $\perp_{\text{b}rj}(y_{0})$ by (ii) of
\textbf{Table A}$_{4}$ here. Therefore,

$\frac{\partial^{2}g^{\text{a}r}}{\partial\text{x}_{i}\partial\text{x}_{j}%
}(y_{0})=\frac{4}{3}(R_{i\text{a}jr}+R_{j\text{a}ir})(y_{0}%
)+4\underset{\text{b=1}}{\overset{\text{q}}{\sum}}T_{\text{ab}i}(y_{0}%
)\perp_{\text{b}rj}(y_{0})$

\section{Table A$_{5}$}

This is computed using the expansion of g$_{\text{ab}}$(x$_{0}$) in
\textbf{Proposition 6.1}:

For a,b = 1,...,q we have:

(i) g$_{\text{ab}}(y_{0})=\delta_{\text{ab}}$

\vspace{1pt}

(ii) $\frac{\partial\text{g}_{\text{ab}}}{\partial\text{x}_{\alpha}}%
(y_{0})=-2$T$_{\text{ab}\alpha}(y_{0})$

\vspace{1pt}

(iii) \ $\frac{\partial^{2}\text{g}_{\text{ab}}}{\partial\text{x}_{\alpha
}\text{x}_{\lambda}}(y_{0})=\underset{\alpha,\lambda\text{=q+1}%
}{\overset{\text{n}}{\sum}}\{-$R$_{\text{a}\alpha\text{b}\lambda}$
$-$R$_{\text{a}\lambda\text{b}\alpha}$+ $\underset{\text{c=1}%
}{\overset{\text{q}}{\sum}}$T$_{\text{ac}\alpha}$T$_{\text{bc}\lambda}$
$+\underset{\text{c=1}}{\overset{\text{q}}{\sum}}$T$_{\text{ac}\lambda}%
$T$_{\text{bc}\alpha}$

$\qquad+\underset{\text{k=q+1}}{\overset{\text{n}}{\sum}}$ $\perp
_{\text{a}\alpha\text{k}}\perp_{\text{b}\lambda\text{k}}%
+\underset{\text{k=q+1}}{\overset{\text{n}}{\sum}}$ $\perp_{\text{a}%
\lambda\text{k}}\perp_{\text{b}\alpha\text{k}}$\}(y$_{0}$)$\ $

(iv) $\qquad\frac{\partial^{2}\text{g}_{\text{ab}}}{\partial\text{x}_{\alpha
}^{2}}(y_{0})=2\underset{\alpha,\lambda\text{=q+1}}{\overset{\text{n}}{\sum}%
}\{-$R$_{\text{a}\alpha\text{b}\alpha}$ + $\underset{\text{c=1}%
}{\overset{\text{q}}{\sum}}$T$_{\text{ac}\alpha}$T$_{\text{bc}\alpha}$
$+\underset{\text{k=q+1}}{\overset{\text{n}}{\sum}}$ $\perp_{\text{a}%
\alpha\text{k}}\perp_{\text{b}\alpha\text{k}}\}(y_{0})$\qquad\ \ \ 

(v)\qquad\ $\frac{\partial^{3}\text{g}_{\text{ab}}}{\partial\text{x}_{\alpha
}^{2}\partial\text{x}_{\lambda}}(y_{0})=-\frac{1}{3}\underset{\alpha
,\lambda\text{=q+1}}{\overset{\text{n}}{\sum}}\{2\nabla_{\alpha}$%
R$_{\alpha\text{a}\lambda b}+$ R$_{b\lambda\alpha\text{T}_{\text{a}\alpha}}$
$+$R$_{b\lambda\alpha\perp_{\text{a}\alpha}}$

$\qquad\qquad\qquad+3$ (R$_{\text{a}\alpha\alpha\text{T}_{b\lambda}}$
$+\;$R$_{\text{a}\alpha\alpha\perp_{b\lambda}})+3($R$_{b\alpha\alpha
\text{T}_{\text{a}\lambda}}$ $+$R$_{b\alpha\alpha\perp_{\text{a}\lambda}})$
$+$R$_{\text{a}\lambda\alpha\text{T}_{b\alpha}}$ $+$ R$_{\text{a}\lambda
\alpha\perp_{b\alpha}}$

\ \ \ \ \ \qquad\ \ +2$\nabla_{\alpha}$R$_{\lambda\text{a}\alpha b}+$
R$_{b\alpha\alpha\text{T}_{\text{a}\lambda}}$ $+$R$_{b\alpha\alpha
\perp_{\text{a}\lambda}}$ $+3$ (R$_{\text{a}\alpha\lambda\text{T}_{b\alpha}}$
$+\;$R$_{\text{a}\alpha\lambda\perp_{b\alpha}})$

\qquad$\ \qquad+3($R$_{b\alpha\lambda\text{T}_{\text{a}\alpha}}$
$+$R$_{b\alpha\lambda\perp_{\text{a}\alpha}})$ $+$R$_{\text{a}\alpha
\alpha\text{T}_{b\text{j}}}$ $+$ R$_{\text{a}\alpha\alpha\perp_{b\lambda}}$

$\qquad\qquad+2\nabla_{\lambda}$R$_{\alpha\text{a}\alpha b}+$ R$_{b\alpha
\lambda\text{T}_{\text{a}\alpha}}$ $+$R$_{b\alpha\lambda\perp_{\text{a}\alpha
}}$ $+3$ (R$_{\text{a}\lambda\alpha\text{T}_{b\alpha}}$ $+\;$R$_{\text{a}%
\lambda\alpha\perp_{b\alpha}})$

\qquad$\qquad+3($R$_{b\lambda\alpha\text{T}_{\text{a}\alpha}}$ $+$%
R$_{b\lambda\alpha\perp_{\text{a}\alpha}})$ $+$R$_{\text{a}\alpha
\lambda\text{T}_{b\alpha\text{j}}}$ $+$ R$_{\text{a}\alpha\lambda
\perp_{b\alpha}}\}(y_{0})$\ 

\subsection{\protect\underline{\textbf{Computations}}}

(i) immediate

(ii) immediate

$\vspace{1pt}$(iii)$\qquad\frac{\partial^{2}\text{g}_{\text{ab}}}%
{\partial\text{x}_{\alpha}\text{x}_{\lambda}}(y_{0})$%
\ $=\underset{\text{i,j=q+1}}{\overset{\text{n}}{\sum}}\{-$R$_{\text{aibj}}$ +
$\underset{\text{c=1}}{\overset{\text{q}}{\sum}}$T$_{\text{aci}}%
$T$_{\text{bcj}}$ $+\underset{\text{k=q+1}}{\overset{\text{n}}{\sum}}$
$\perp_{\text{aik}}\perp_{\text{bjk}}$\}(y$_{0}$)

$\qquad\qquad\qquad\qquad\times(\delta_{\alpha i}\delta_{\lambda j}%
+\delta_{\alpha j}\delta_{\lambda i})$

$\vspace{1pt}$ \ \ \ \ \ $\qquad=\underset{\text{i,j=q+1}}{\overset{\text{n}%
}{\sum}}\{-$R$_{\text{a}\alpha\text{b}\lambda}$ $-$R$_{\text{a}\lambda
\text{b}\alpha}$+ $\underset{\text{c=1}}{\overset{\text{q}}{\sum}}%
$T$_{\text{ac}\alpha}$T$_{\text{bc}\lambda}$ $+\underset{\text{c=1}%
}{\overset{\text{q}}{\sum}}$T$_{\text{ac}\lambda}$T$_{\text{bc}\alpha}$

$\qquad\qquad+\underset{\text{k=q+1}}{\overset{\text{n}}{\sum}}$
$\perp_{\text{a}\alpha\text{k}}\perp_{\text{b}\lambda\text{k}}%
+\underset{\text{k=q+1}}{\overset{\text{n}}{\sum}}$ $\perp_{\text{a}%
\lambda\text{k}}\perp_{\text{b}\alpha\text{k}}$\}(y$_{0}$)$\ $

$\vspace{1pt}$(iv) Obvious.

\qquad$\ \ \ \ $\qquad\qquad

(v) $\ \ \frac{\partial^{3}\text{g}_{\text{ab}}}{\partial\text{x}_{\alpha}%
^{2}\partial\text{x}_{\lambda}}(y_{0}%
)\ =\overset{n}{\underset{i,j,k=q+1}{\text{ }-\frac{1}{6}\sum}\{}%
2\nabla_{\text{i}}($R)$_{\text{jakb}}+$ R$_{\text{bkiT}_{\text{aj}}}$
$+$R$_{\text{bki}\perp_{\text{aj}}}$

$\qquad+3$ (R$_{\text{aijT}_{\text{bk}}}$ $+\;$R$_{\text{aij}\perp_{\text{bk}%
}})+3($R$_{\text{bijT}_{\text{aK}}}$ $+$R$_{\text{bij}\perp_{\text{aK}}})$
$+$R$_{\text{akiT}_{\text{bj}}}$ $+$ R$_{\text{aki}\perp_{\text{bj}}}\}$

$\qquad\times(2\delta_{\alpha i}\delta_{\alpha j}\delta_{\lambda k}%
+2\delta_{\alpha i}\delta_{\alpha k}\delta_{\lambda j}+2\delta_{\alpha
j}\delta_{\alpha k}\delta_{\lambda i})$

\vspace{1pt}\qquad\ \ $\overset{n}{\underset{i,j,k=q+1}{=-\frac{1}{3}\sum}%
\{}2\nabla_{\alpha}$R$_{\alpha\text{a}\lambda b}+$ R$_{b\lambda\alpha
\text{T}_{\text{a}\alpha}}$ $+$R$_{b\lambda\alpha\perp_{\text{a}\alpha}}$ $+3$
(R$_{\text{a}\alpha\alpha\text{T}_{b\lambda}}$ $+\;$R$_{\text{a}\alpha
\alpha\perp_{b\lambda}})$

\qquad$+3($R$_{b\alpha\alpha\text{T}_{\text{a}\lambda}}$ $+$R$_{b\alpha
\alpha\perp_{\text{a}\lambda}})$ $+$R$_{\text{a}\lambda\alpha\text{T}%
_{b\alpha}}$ $+$ R$_{\text{a}\lambda\alpha\perp_{b\alpha}}$

\ \ \ \ \ \ \ \ +2$\nabla_{\alpha}$R$_{\lambda\text{a}\alpha b}+$
R$_{b\alpha\alpha\text{T}_{\text{a}\lambda}}$ $+$R$_{b\alpha\alpha
\perp_{\text{a}\lambda}}$ $+3$ (R$_{\text{a}\alpha\lambda\text{T}_{b\alpha}}$
$+\;$R$_{\text{a}\alpha\lambda\perp_{b\alpha}})$

\qquad\ $+3($R$_{b\alpha\lambda\text{T}_{\text{a}\alpha}}$ $+$R$_{b\alpha
\lambda\perp_{\text{a}\alpha}})$ $+$R$_{\text{a}\alpha\alpha\text{T}%
_{b\text{j}}}$ $+$ R$_{\text{a}\alpha\alpha\perp_{b\lambda}}$

\qquad$+2\nabla_{\lambda}$R$_{\alpha\text{a}\alpha b}+$ R$_{b\alpha
\lambda\text{T}_{\text{a}\alpha}}$ $+$R$_{b\alpha\lambda\perp_{\text{a}\alpha
}}$ $+3$ (R$_{\text{a}\lambda\alpha\text{T}_{b\alpha}}$ $+\;$R$_{\text{a}%
\lambda\alpha\perp_{b\alpha}})$

\qquad$+3($R$_{b\lambda\alpha\text{T}_{\text{a}\alpha}}$ $+$R$_{b\lambda
\alpha\perp_{\text{a}\alpha}})$ $+$R$_{\text{a}\alpha\lambda\text{T}%
_{b\alpha\text{j}}}$ $+$ R$_{\text{a}\alpha\lambda\perp_{b\alpha}}\}(y_{0})$

\section{{\protect\LARGE Table A}$_{6}$}

\vspace{1pt}For a,b = 1,...,q \ and $\alpha=q+1,...,n$

(i) \qquad g$^{\text{ab}}(y_{0})=\delta_{\text{ab}}$

\vspace{1pt}

(ii)\qquad\ $\frac{\partial\text{g}^{\text{ab}}}{\partial\text{x}_{i}}%
(y_{0})=$ $2$T$_{\text{ab}i}(y_{0})$

\vspace{1pt}

(iii) $\qquad\frac{\partial^{2}\text{g}^{\text{ab}}}{\partial\text{x}_{i}^{2}%
}(y_{0})=2\left\{  -R_{\text{a}i\text{b}i}+5\overset{q}{\underset{\text{c}%
=1}{\sum}}T_{\text{ac}i}T_{\text{bc}i}+2\overset{n}{\underset{k=q+1}{\sum}%
}\perp_{\text{a}ik}\perp_{\text{b}ik}\right\}  (y_{0})$

\subsection{\protect\underline{\textbf{Computations}}}

(i) \qquad From, \ g$_{\text{ac}}$g$^{\text{cb}}=\delta_{\text{ab}},$ \ and
g$_{\text{ac}}$(y$_{0})=\delta_{\text{ac}}$ we have:

\qquad g$_{\text{ac}}$(y$_{0})$g$^{\text{cb}}($y$_{0})=\delta_{\text{ab}}$

\qquad g$_{\text{ac}}$(y$_{0})$g$^{\text{cb}}($y$_{0})=\delta_{ac}%
$g$^{\text{cb}}($y$_{0})=$ g$^{\text{ab}}($y$_{0})$Hence,

\qquad g$^{\text{ab}}($y$_{0})=\delta_{\text{ab}}$

\vspace{1pt}

(ii) Again from g$_{ac}$g$^{cb}=\delta_{ab}$

\qquad we have:

\qquad0 = $\frac{\partial}{\partial\text{x}\alpha}($g$_{\text{cb}}%
$g$^{\text{cb}})($y$_{0})=$ g$_{\text{ac}}($y$_{0})\frac{\partial
\text{g}^{\text{cb}}}{\partial\text{x}_{\alpha}}($y$_{0})+$g$^{\text{cb}}%
($y$_{0})\frac{\partial\text{g}_{\text{ac}}}{\partial\text{x}_{\alpha}}%
($y$_{0})$

$\qquad=\delta_{\text{ac}}($y$_{0})\frac{\partial\text{g}^{\text{cb}}%
}{\partial\text{x}_{\alpha}}(y_{0})+\delta^{\text{cb}}\frac{\partial
\text{g}_{\text{ac}}}{\partial\text{x}_{\alpha}}($y$_{0})$

$\ \qquad=\frac{\partial\text{g}^{\text{ab}}}{\partial\text{x}_{\alpha}}%
(y_{0})+\frac{\partial\text{g}_{\text{ab}}}{\partial\text{x}_{\alpha}}($%
y$_{0})$

Therefore,

$\frac{\partial\text{g}^{\text{ab}}}{\partial\text{x}_{\alpha}}(y_{0}%
)=-\frac{\partial\text{g}_{\text{ab}}}{\partial\text{x}_{\alpha}}($y$_{0}%
)=2$T$_{\text{ab}\alpha}\qquad$

(iii) Then, we have:

0 = $\frac{\partial^{2}}{\partial\text{x}^{2}\alpha}($g$_{\text{ac}}%
$g$^{\text{cb}})($y$_{0})=$ $\delta_{\text{ac}}\frac{\partial^{2}%
\text{g}^{\text{cb}}}{\partial\text{x}^{2}\alpha}($y$_{0})+2\frac
{\partial\text{g}_{\text{ac}}}{\partial\text{x}_{\alpha}}($y$_{0}%
)\frac{\partial\text{g}^{\text{cb}}}{\partial\text{x}_{\alpha}}(y_{0}%
)+\delta^{cb}\frac{\partial^{2}\text{g}_{\text{ac}}}{\partial\text{x}%
^{2}\alpha}($y$_{0})$

$\qquad=\frac{\partial^{2}\text{g}^{\text{ab}}}{\partial\text{x}^{2}\alpha
}(y_{0})+2\frac{\partial\text{g}_{\text{ac}}}{\partial\text{x}_{\alpha}}%
(y_{0})\frac{\partial\text{g}^{\text{cb}}}{\partial\text{x}_{\alpha}}%
(y_{0})+\frac{\partial^{2}\text{g}_{\text{ab}}}{\partial\text{x}^{2}\alpha}%
($y$_{0})$Therefore,

$\frac{\partial^{2}\text{g}^{\text{ab}}}{\partial\text{x}^{2}\alpha}%
(y_{0})=-2\overset{q}{\underset{c=1}{\sum}}\frac{\partial\text{g}_{\text{ac}}%
}{\partial\text{x}_{\alpha}}(y_{0})\frac{\partial\text{g}^{\text{cb}}%
}{\partial\text{x}_{\alpha}}(y_{0})-2\overset{n}{\underset{j=q+1}{\sum}}%
\frac{\partial\text{g}_{\text{a}j}}{\partial\text{x}_{\alpha}}(y_{0}%
)\frac{\partial\text{g}^{jb}}{\partial\text{x}_{\alpha}}(y_{0})-\frac
{\partial^{2}\text{g}_{\text{ab}}}{\partial\text{x}^{2}\alpha}($y$_{0})$

Sequentially, we have by (ii) of Table 5 and (ii) of Table 6, by (ii)

of Table 3 and (ii) of Table 4; then (iii) of Table 5:

$\qquad\frac{\partial^{2}\text{g}^{\text{ab}}}{\partial\text{x}^{2}\alpha
}(y_{0})=-2\overset{q}{\underset{c=1}{\sum}}(-2T_{\text{ac}\alpha
})(2T_{\text{cb}\alpha})-2\overset{n}{\underset{j=q+1}{\sum}}(-\perp
_{\text{a}j\alpha})(\perp_{\text{b}j\alpha})$

\qquad\qquad\qquad$+2(-R_{\text{a}\alpha\text{b}\alpha}+$
$\underset{\text{c=1}}{\overset{\text{q}}{\sum}}T_{\text{ac}\alpha
}T_{\text{bc}\alpha}+\overset{\text{n}}{\underset{\text{k=q+1}}{\sum}}%
\perp_{\text{a}\alpha\text{k}}\perp_{\text{b}\alpha\text{k}})$

\qquad\qquad\ \ \ = $8\overset{q}{\underset{\text{c}=1}{\sum}}(T_{\text{ac}%
\alpha})(T_{\text{bc}\alpha})+2\overset{n}{\underset{\text{k}=q+1}{\sum}%
}(\perp_{\text{a}\alpha\text{k}})(\perp_{\text{b}\alpha\text{k}})$

\qquad\qquad\ \ \ $+2(-R_{\text{a}\alpha\text{b}\alpha}+$
$\underset{\text{c=1}}{\overset{\text{q}}{\sum}}T_{\text{ac}\alpha
}T_{\text{bc}\alpha}\overset{\text{n}}{+\underset{\text{k=q+1}}{\sum}}%
\perp_{\text{a}\alpha\text{k}}\perp_{\text{b}\alpha\text{k}})$

\qquad\qquad\ \ \ $=2(-R_{\text{a}\alpha\text{b}\alpha}+5$
$\underset{\text{c=1}}{\overset{\text{q}}{\sum}}T_{\text{ac}\alpha
}T_{\text{bc}\alpha}\overset{\text{n}}{+2\underset{\text{k=q+1}}{\sum}}%
\perp_{\text{a}\alpha\text{k}}\perp_{\text{b}\alpha\text{k}})$\qquad

\section{\textbf{Table A}$_{7}$}

For a,b,c = 1,...,q and $i,j,k=q+1,...,n$ we have:

(i) $\Gamma_{\text{ab}}^{j}(y_{0})=T_{\text{ab}j}(y_{0})$

(ii) $\Gamma_{\text{ab}}^{\text{c}}(y_{0})=0$

(iii ) $\frac{\partial\Gamma_{\text{ab}}^{j}}{\partial\text{x}_{i}}%
(y_{0})=\frac{1}{2}\{$ $(R_{\text{a}i\text{b}j}+R_{\text{a}j\text{b}i})$
$-\underset{\text{c=1}}{\overset{\text{q}}{\sum}}(T_{\text{ac}i}T_{\text{bc}%
j}+T_{\text{ac}j}T_{\text{bc}i})$

$\qquad-\overset{n}{\underset{k=q+1}{\sum}}(\perp_{\text{a}ki}\perp
_{\text{b}kj}+$ $\perp_{\text{a}kj}\perp_{\text{b}ki})\}(y_{0})$

(iv) $\frac{\partial\Gamma_{\text{aa}}^{j}}{\partial\text{x}_{i}}(y_{0})=\{$
R$_{\text{a}i\text{a}j}$ $-\underset{\text{c=1}}{\overset{\text{q}}{\sum}%
}T_{\text{ac}i}T_{\text{ac}j}-\overset{n}{\underset{k=q+1}{\sum}}%
(\perp_{\text{a}ki}\perp_{\text{a}kj}\}(y_{0})$

(v) $\ \frac{\partial\Gamma_{\text{ab}}^{\text{c}}}{\partial\text{x}_{i}%
}(y_{0})=-\overset{n}{\underset{k=q+1}{\sum}}(\perp_{\text{c}ik}%
)($T$_{\text{ab}k})(y_{0})$

(vi) $\ \frac{\partial^{2}\Gamma_{\text{ab}}^{\lambda}}{\partial
\text{x}_{\alpha}^{2}}(y_{0})=\frac{1}{6}[\{2\nabla_{\alpha}$R$_{\alpha
\text{a}\lambda\text{b}}+$ R$_{\text{b}\lambda\alpha\text{T}_{\text{a}\alpha}%
}$ $+$R$_{\text{b}\lambda\alpha\perp_{\text{a}\alpha}}$ $+3$ (R$_{\text{a}%
\alpha\alpha\text{T}_{\text{b}\lambda}}$ $+\;$R$_{\text{a}\alpha\alpha
\perp_{\text{b}\lambda}})$

\qquad\qquad\ $+3($R$_{\text{b}\alpha\alpha\text{T}_{\text{a}\lambda}}$
$+$R$_{\text{b}\alpha\alpha\perp_{\text{a}\lambda}})$ $+$R$_{\text{a}%
\lambda\alpha\text{T}_{\text{b}\alpha}}$ $+$ R$_{\text{a}\lambda\alpha
\perp_{\text{b}\alpha}}$

\ \ \ \ \ \ \ \ \ \ \ \ \ \ \ \ +2$\nabla_{\alpha}$R$_{\lambda\text{a}%
\alpha\text{b}}+$ R$_{\text{b}\alpha\alpha\text{T}_{\text{a}\lambda}}$
$+$R$_{\text{b}\alpha\alpha\perp_{\text{a}\lambda}}$ $+3$ (R$_{\text{a}%
\alpha\lambda\text{T}_{\text{b}\alpha}}$ $+\;$R$_{\text{a}\alpha\lambda
\perp_{\text{b}\alpha}})$

\qquad\qquad\ $+3($R$_{\text{b}\alpha\lambda\text{T}_{\text{a}\alpha}}$
$+$R$_{\text{b}\alpha\lambda\perp_{\text{a}\alpha}})$ $+$R$_{\text{a}%
\alpha\alpha\text{T}_{\text{b}\lambda}}$ $+$ R$_{\text{a}\alpha\alpha
\perp_{\text{b}\lambda}}$

\qquad\qquad$+2\nabla_{\lambda}$R$_{\alpha\text{a}\alpha\text{b}}+$
R$_{\text{b}\alpha\lambda\text{T}_{\text{a}\alpha}}$ $+$R$_{\text{b}%
\alpha\lambda\perp_{\text{a}\alpha}}$ $+3$ (R$_{\text{a}\lambda\alpha
\text{T}_{\text{b}\alpha}}$ $+\;$R$_{\text{a}\lambda\alpha\perp_{\text{b}%
\alpha}})$

\qquad\ \ \ \ \ \ $+3($R$_{\text{b}\lambda\alpha\text{T}_{\text{a}\alpha}}$
$+$R$_{\text{b}\lambda\alpha\perp_{\text{a}\alpha}})$ $+$R$_{\text{a}%
\alpha\lambda\text{T}_{\text{b}\alpha}}$ $+$ R$_{\text{a}\alpha\lambda
\perp_{\text{b}\alpha}}\}$\ 

\qquad$\ \ \ \ +\frac{2}{3}\underset{\gamma=q+1}{\overset{n}{\sum}%
}\{T_{\text{ab}\gamma}(R_{\alpha\lambda\alpha\gamma}%
-3\overset{q}{\underset{\text{c}=1}{\sum}}\perp_{\text{c}\alpha\lambda}%
\perp_{\text{c}\alpha\gamma})\}](y_{0})$

$\ \frac{\partial^{2}\Gamma_{\text{aa}}^{\lambda}}{\partial\text{x}_{\alpha
}^{2}}(y_{0})=\frac{1}{6}[\{2\nabla_{\alpha}$R$_{\alpha\text{a}\lambda
\text{a}}+$ R$_{\text{a}\lambda\alpha\text{T}_{\text{a}\alpha}}$
$+$R$_{\text{a}\lambda\alpha\perp_{\text{a}\alpha}}$ $+3$ (R$_{\text{a}%
\alpha\alpha\text{T}_{\text{a}\lambda}}$ $+\;$R$_{\text{a}\alpha\alpha
\perp_{\text{a}\lambda}})$

\qquad\qquad\ $+3($R$_{\text{a}\alpha\alpha\text{T}_{\text{a}\lambda}}$
$+$R$_{\text{a}\alpha\alpha\perp_{\text{a}\lambda}})$ $+$R$_{\text{a}%
\lambda\alpha\text{T}_{\text{a}\alpha}}$ $+$ R$_{\text{a}\lambda\alpha
\perp_{\text{a}\alpha}}$

\ \ \ \ \ \ \ \ \ \ \ \ \ \ \ \ $+2\nabla_{\alpha}$R$_{\lambda\text{a}%
\alpha\text{a}}+$ R$_{\text{a}\alpha\alpha\text{T}_{\text{a}\lambda}}$ $+$
R$_{\text{a}\alpha\alpha\perp_{\text{a}\lambda}}$ $+3$ (R$_{\text{a}%
\alpha\lambda\text{T}_{\text{a}\alpha}}$ $+\;$R$_{\text{a}\alpha\lambda
\perp_{\text{a}\alpha}})$

\qquad\qquad\ $+3($R$_{\text{a}\alpha\lambda\text{T}_{\text{a}\alpha}}$ $+$
R$_{\text{a}\alpha\lambda\perp_{\text{a}\alpha}})$ $+$R$_{\text{a}\alpha
\alpha\text{T}_{\text{a}\lambda}}$ $+$ R$_{\text{a}\alpha\alpha\perp
_{\text{a}\lambda}}$

\qquad\qquad$+2\nabla_{\lambda}$R$_{\alpha\text{a}\alpha\text{a}}+$
R$_{\text{a}\alpha\lambda\text{T}_{\text{a}\alpha}}$ $+$ R$_{\text{a}%
\alpha\lambda\perp_{\text{a}\alpha}}$ $+3$(R$_{\text{a}\lambda\alpha
\text{T}_{\text{a}\alpha}}$ $+\;$R$_{\text{a}\lambda\alpha\perp_{\text{a}%
\alpha}})$

\qquad\ \ \ \ \ \ $+3($R$_{\text{a}\lambda\alpha\text{T}_{\text{a}\alpha}}$
$+$ R$_{\text{a}\lambda\alpha\perp_{\text{a}\alpha}})$ $+$ R$_{\text{a}%
\alpha\lambda\text{T}_{\text{a}\alpha}}$ $+$ R$_{\text{a}\alpha\lambda
\perp_{\text{a}\alpha}}\}$\ 

\qquad$\ \ \ \ +\frac{2}{3}\underset{\gamma=q+1}{\overset{n}{\sum}%
}\{T_{\text{aa}\gamma}(R_{\alpha\lambda\alpha\gamma}%
-3\overset{q}{\underset{\text{c}=1}{\sum}}\perp_{\text{c}\alpha\lambda}%
\perp_{\text{c}\alpha\gamma})\}](y_{0})$

(vii) $\frac{\partial^{2}\Gamma_{\text{aa}}^{\lambda}}{\partial\text{x}%
_{\alpha}^{2}}(y_{0})=\frac{1}{6}[\{4\nabla_{\alpha}$R$_{\alpha\text{a}%
\lambda\text{a}}$ $+2\nabla_{\lambda}$R$_{\alpha\text{a}\alpha\text{a}}+$ $8$
(R$_{\text{a}\alpha\alpha\text{T}_{\text{a}\lambda}}$ $+\;$R$_{\text{a}%
\alpha\alpha\perp_{\text{a}\lambda}})$

\ \ \ \ \ \ \ \ \ \ \ \ \qquad\ \ \ \ \ \ \ $+8$ (R$_{\text{a}\alpha
\lambda\text{T}_{\text{a}\alpha}}$ $+\;$R$_{\text{a}\alpha\lambda
\perp_{\text{a}\alpha}})+8$ (R$_{\text{a}\lambda\alpha\text{T}_{\text{a}%
\alpha}}$ $+\;$R$_{\text{a}\lambda\alpha\perp_{\text{a}\alpha}})\}$\ 

\qquad$\qquad\ \ \ \ \ \ \ \ \ +\frac{2}{3}\underset{\gamma
=q+1}{\overset{n}{\sum}}\{T_{\text{aa}\gamma}(R_{\alpha\lambda\alpha\gamma
}+3\overset{q}{\underset{\text{c}=1}{\sum}}\perp_{\text{c}\alpha\lambda}%
\perp_{\text{c}\alpha\gamma})\}](y_{0})\qquad$

(vii) $\frac{\partial^{2}\Gamma_{\text{aa}}^{j}}{\partial\text{x}_{i}^{2}%
}(y_{0})=\frac{1}{6}[\{4\nabla_{i}$R$_{i\text{a}j\text{a}}$ $+2\nabla_{j}%
$R$_{i\text{a}i\text{a}}+$ $8$ $($R$_{\text{a}ii\text{T}_{\text{a}j}}$
$+\overset{n}{\underset{k=q+1}{%
{\textstyle\sum}
}}R_{\text{a}iik}\perp_{\text{a}jk})$

\ \ \ \ \ \ \ \ \ \ \ \ \qquad\ \ \ \ \ \ \ $+8($R$_{\text{a}ij\text{T}%
_{\text{a}i}}$ $+\;$R$_{\text{a}ij\perp_{\text{a}i}})+8($R$_{\text{a}%
ji\text{T}_{\text{a}i}}$ $+\;$R$_{\text{a}ji\perp_{\text{a}i}})\}$\ 

\qquad$\qquad\ \ \ \ \ \ \ \ \ +\frac{2}{3}\underset{k=q+1}{\overset{n}{\sum}%
}\{T_{\text{aa}k}(R_{ijik}+3\overset{q}{\underset{\text{c}=1}{\sum}}%
\perp_{\text{c}ij}\perp_{\text{c}ik})\}](y_{0})$

$\frac{\partial^{2}\Gamma_{\text{aa}}^{j}}{\partial\text{x}_{i}^{2}}%
(y_{0})=\frac{1}{6}[\{4\nabla_{i}$R$_{i\text{a}j\text{a}}$ $+2\nabla_{j}%
$R$_{i\text{a}i\text{a}}+$ $8$ $(\overset{q}{\underset{\text{c=1}}{%
{\textstyle\sum}
}}R_{\text{a}i\text{c}i}^{{}}T_{\text{ac}j}+\;\overset{n}{\underset{k=q+1}{%
{\textstyle\sum}
}}R_{\text{a}iik}\perp_{\text{a}jk})$

\ \ \ \ \ \ \ \ $+8(\overset{q}{\underset{\text{c=1}}{%
{\textstyle\sum}
}}R_{\text{a}i\text{c}j}^{{}}T_{\text{ac}i}+\overset{n}{\underset{l=q+1}{%
{\textstyle\sum}
}}R_{\text{a}ijl}\perp_{\text{a}il})+8(\overset{q}{\underset{\text{c=1}}{%
{\textstyle\sum}
}}R_{\text{a}j\text{c}i}^{{}}T_{\text{ac}i}+\overset{q}{\underset{\text{c=1}}{%
{\textstyle\sum}
}}R_{\text{a}j\text{c}i}^{{}}T_{\text{ac}i})\}$\ 

\qquad$+\frac{2}{3}\underset{k=q+1}{\overset{n}{\sum}}\{T_{\text{aa}%
k}(R_{ijik}+3\overset{q}{\underset{\text{c}=1}{\sum}}\perp_{\text{c}ij}%
\perp_{\text{c}ik})\}](y_{0})$

Note that we have at the center of Fermi coordinates $y_{0}\in P,$

R$_{\text{b}ki\text{T}_{\text{a}j}}=$ $<R_{\text{X}_{\text{b}}X_{k}X_{i}}^{{}%
},T_{\text{X}_{\text{a}}}X_{j}>$ $=$ $<R_{\text{b}ki}^{{}},T_{\text{a}j}>$
$=\overset{q}{\underset{\text{c,d=1}}{%
{\textstyle\sum}
}}R_{\text{b}ki\text{c}}^{{}}T_{\text{a}j\text{d}}<\frac{\partial}%
{\partial\text{x}_{\text{c}}},\frac{\partial}{\partial\text{x}_{\text{d}}}>$

$=\overset{q}{\underset{\text{c,d=1}}{%
{\textstyle\sum}
}}R_{\text{b}ki\text{c}}^{{}}T_{\text{a}j\text{d}}\delta_{\text{cd}%
}=\overset{q}{\underset{\text{c=1}}{%
{\textstyle\sum}
}}R_{\text{b}ki\text{c}}^{{}}T_{\text{a}j\text{c}}%
=\overset{q}{\underset{\text{c=1}}{%
{\textstyle\sum}
}}R_{\text{b}k\text{c}i}^{{}}T_{\text{ac}j}$

A slight change of indices gives:

R$_{\text{a}ij\text{T}_{\text{b}k}}=\overset{q}{\underset{\text{c=1}}{%
{\textstyle\sum}
}}R_{\text{a}i\text{c}j}^{{}}T_{\text{bc}k}$

R$_{\text{a}ij\text{T}_{\text{a}i}}=\overset{q}{\underset{\text{c=1}}{%
{\textstyle\sum}
}}R_{\text{a}i\text{c}j}^{{}}T_{\text{ac}i}\qquad$R$_{\text{a}ii\text{T}%
_{\text{a}j}}=\overset{q}{\underset{\text{c=1}}{%
{\textstyle\sum}
}}R_{\text{a}i\text{c}i}^{{}}T_{\text{ac}j}$

R$_{\text{a}ii\text{T}_{\text{a}j}}=\overset{q}{\underset{\text{c=1}}{%
{\textstyle\sum}
}}R_{\text{a}i\text{c}i}^{{}}T_{\text{ac}j}$

Similarly,

R$_{\text{a}ij\perp_{\text{b}k}}=$ $<R_{\text{X}_{\text{a}}X_{i}X_{j}}^{{}%
},\perp_{\text{X}_{\text{b}}}X_{k}>$ $=$ $<R_{\text{a}ij},\perp_{\text{b}k}>$
$=$ $\overset{n}{\underset{l,m=q+1}{%
{\textstyle\sum}
}}R_{\text{a}ijl}\perp_{\text{b}km}<\frac{\partial}{\partial\text{x}_{l}%
},\frac{\partial}{\partial\text{x}_{m}}>$

$=\overset{n}{\underset{l,m}{%
{\textstyle\sum}
}}R_{\text{a}ijl}\perp_{\text{b}km}\delta_{lm}=\overset{n}{\underset{l=q+1}{%
{\textstyle\sum}
}}R_{\text{a}ijl}\perp_{\text{b}kl}$

R$_{\text{a}ij\perp_{\text{b}k}}=\overset{n}{\underset{l=q+1}{%
{\textstyle\sum}
}}R_{\text{a}ijl}\perp_{\text{b}kl}$

R$_{\text{a}ij\perp_{\text{a}i}}=\overset{n}{\underset{k=q+1}{%
{\textstyle\sum}
}}R_{\text{a}ijk}\perp_{\text{a}ik}$

R$_{\text{a}ii\perp_{\text{a}j}}=\overset{n}{\underset{k=q+1}{%
{\textstyle\sum}
}}R_{\text{a}iik}\perp_{\text{a}jk}$

We also have at $y_{0}\in P,$

R$_{i\text{a}jR_{k\text{b}l}}=$ $<R_{X_{i}\text{X}_{\text{a}}X_{j}}^{{}%
},R_{X_{k}\text{X}_{\text{b}}X_{l}}^{{}}>$ $=$ $\overset{n}{\underset{m,s=1}{%
{\textstyle\sum}
}}R_{i\text{a}jm}R_{k\text{b}ls}<\frac{\partial}{\partial\text{x}_{m}}%
,\frac{\partial}{\partial\text{x}_{s}}>$

$=$ $\overset{n}{\underset{m=1}{%
{\textstyle\sum}
}}R_{i\text{a}jm}R_{k\text{b}lm}=$ $\overset{n}{\underset{m=1}{%
{\textstyle\sum}
}}R_{\text{a}ijm}R_{\text{b}klm}$

We summarize:

R$_{\text{a}ij\text{T}_{\text{b}k}}=\overset{q}{\underset{\text{c=1}}{%
{\textstyle\sum}
}}R_{\text{a}i\text{c}j}^{{}}T_{\text{bc}k}$

R$_{\text{a}ij\perp_{\text{b}k}}=\overset{n}{\underset{l=q+1}{%
{\textstyle\sum}
}}R_{\text{a}ijl}\perp_{\text{b}kl}$

R$_{i\text{a}jR_{k\text{b}l}}=$ $\overset{n}{\underset{m=1}{%
{\textstyle\sum}
}}R_{i\text{a}jm}R_{k\text{b}lm}=$ $\overset{n}{\underset{m=1}{%
{\textstyle\sum}
}}R_{\text{a}ijm}R_{\text{b}klm}$

Therefore,

$\qquad\frac{\partial^{2}\Gamma_{\text{aa}}^{j}}{\partial\text{x}_{i}^{2}%
}(y_{0})=\frac{1}{6}[\{4\nabla_{i}R_{i\text{a}j\text{a}}+2\nabla
_{j}R_{i\text{a}i\text{a}}+$ $8(\overset{q}{\underset{\text{c=1}}{%
{\textstyle\sum}
}}R_{\text{a}i\text{c}i}^{{}}T_{\text{ac}j}+\;\overset{n}{\underset{k=q+1}{%
{\textstyle\sum}
}}R_{\text{a}iik}\perp_{\text{a}jk})$

\ \ \ \ \ \ \ $+8(\overset{q}{\underset{\text{c=1}}{%
{\textstyle\sum}
}}R_{\text{a}i\text{c}j}^{{}}T_{\text{ac}i}+\;\overset{n}{\underset{k=q+1}{%
{\textstyle\sum}
}}R_{\text{a}ijk}\perp_{\text{a}ik})+8(\overset{q}{\underset{\text{c=1}}{%
{\textstyle\sum}
}}R_{\text{a}j\text{c}i}^{{}}T_{\text{ac}i}+\;\overset{n}{\underset{k=q+1}{%
{\textstyle\sum}
}}R_{\text{a}jik}\perp_{\text{a}ik})\}$\ 

\qquad$+\frac{2}{3}\underset{k=q+1}{\overset{n}{\sum}}\{T_{\text{aa}%
k}(R_{ijik}+3\overset{q}{\underset{\text{c}=1}{\sum}}\perp_{\text{c}ij}%
\perp_{\text{c}ik})\}](y_{0})$

(viii) $\qquad\frac{\partial^{2}\Gamma_{\text{ab}}^{\text{c}}}{\partial
\text{x}_{i}^{2}}(y_{0})$\ = $4\overset{n}{\underset{k=q+1}{\sum}}%
$T$_{\text{ab}k}[\underset{\text{d}=1}{\overset{q}{\sum}}(T_{\text{dc}%
k})(\perp_{\text{d}ik})+\frac{2}{3}R_{i\text{c}ik}]$

\qquad\ $\ -$ $\underset{k,l=q+1}{\overset{n}{\sum}}\perp_{\text{ck}i}%
(y_{0})[-(R_{\text{a}k\text{b}l}+R_{\text{a}l\text{b}k})+\underset{\text{d=1}%
}{\overset{\text{q}}{\sum}}(T_{\text{ad}k}T_{\text{bd}l}+T_{\text{ad}%
l}T_{\text{bd}k})](y_{0})$

\qquad$\ -\underset{k,l=q+1}{\overset{n}{\sum}}\perp_{\text{ck}i}%
(y_{0})[\underset{r=q+1}{\overset{n}{\sum}}(\perp_{\text{a}kr}\perp
_{\text{b}lr}+\perp_{\text{a}lr}\perp_{\text{b}kr})](y_{0})$

(ix)\qquad\ $\frac{\partial^{2}\Gamma_{\text{aa}}^{\text{c}}}{\partial
\text{x}_{i}^{2}}(y_{0})=$ $\overset{n}{\underset{k=q+1}{\sum}}$%
T$_{\text{aa}k}(y_{0})[\frac{8}{3}R_{i\text{c}ik}+4\underset{\text{d}%
=1}{\overset{q}{\sum}}(T_{\text{dc}k})(\perp_{\text{d}ik})](y_{0})$

$-$ $2\overset{n}{\underset{k,l=q+1}{\sum}}\perp_{\text{c}ki}(y_{0}%
)[-R_{\text{a}k\text{a}l}+\underset{\text{d=1}}{\overset{\text{q}}{\sum}%
}T_{\text{ad}k}T_{\text{ad}l}](y_{0})-2\overset{n}{\underset{k,l=q+1}{\sum}%
}\perp_{\text{c}ki}(y_{0})[\underset{r=q+1}{\overset{n}{\sum}}\perp
_{\text{a}kr}\perp_{\text{a}lr}](y_{0})$

(x)\qquad\ $\Gamma_{\text{a}\beta}^{\lambda}(y_{0})=$ $\perp_{\text{a}%
\beta\lambda}(y_{0})$

(xi)\qquad\ $\frac{\partial\Gamma_{\text{a}\beta}^{\lambda}}{\partial
\text{x}_{\alpha}}(y_{0})=\ \overset{q}{\underset{\text{b}=1}{\sum(}}%
\perp_{\text{b}\alpha\lambda}T_{\text{ab}\beta})(y_{0})+\frac{2}%
{3}[2R_{\text{a}\alpha\beta\lambda}+R_{\text{a}\beta\alpha\lambda}%
+R_{\text{a}\lambda\beta\alpha}](y_{0})$

(xii)$\ \ \qquad\ \frac{\partial\Gamma_{\beta\text{a}}^{\lambda}}{\partial
x_{\alpha}}(y_{0})=$ $\overset{q}{\underset{\text{b}=1}{\sum}}(\perp
_{\text{b}\alpha\lambda})(T_{\text{ab}\beta})+[\frac{2}{3}(2R_{\text{a}%
\alpha\beta\lambda}+R_{\text{a}\beta\alpha\lambda}+R_{\text{a}\lambda
\beta\alpha})](y_{0})$

\qquad

\subsection{\qquad\protect\underline{\textbf{Computations}}}

(i) Since derivatives with respect to the tangential variables

\qquad\ \ \ x$_{\text{a}}$ and x$_{\text{b}}$ all vanish,

$\ \ \ \Gamma_{\text{ab}}^{\lambda}=\frac{1}{2}\overset{n}{\underset{\gamma
=1}{\sum}}$g$^{\lambda\gamma}(\frac{\partial\text{g}_{\text{b}\gamma}%
}{\partial\text{x}_{\text{a}}}+\frac{\partial\text{g}_{\text{a}\gamma}%
}{\partial\text{x}_{\text{b}}}-\frac{\partial\text{g}_{\text{ab}}}%
{\partial\text{x}_{\gamma}})$\ \ = $-\frac{1}{2}\overset{n}{\underset{\gamma
=1}{\sum}}$g$^{\lambda\gamma}\frac{\partial\text{g}_{\text{ab}}}%
{\partial\text{x}_{\gamma}}$

\ \ \ Since g$^{\lambda\gamma}(y_{0})=\delta^{\lambda\gamma}$,

\qquad\qquad$\Gamma_{\text{ab}}^{\lambda}(y_{0})$\ \ = $-\frac{1}%
{2}\overset{n}{\underset{\gamma=1}{\sum}}$g$^{\lambda\gamma}(y_{0}%
)\frac{\partial\text{g}_{\text{ab}}}{\partial\text{x}_{\gamma}}(y_{0}%
)=-\frac{1}{2}\frac{\partial\text{g}_{\text{ab}}}{\partial\text{x}_{\lambda}%
}(y_{0})$

$\qquad=T_{\text{ab}\lambda}$\ by (ii) of Table 5.

(ii) $\Gamma_{\text{ab}}^{\text{c}}(y_{0})=-\frac{1}{2}%
\overset{n}{\underset{k=1}{\sum}}$g$^{\text{c}k}(y_{0})\frac{\partial
\text{g}_{\text{ab}}}{\partial\text{x}_{k}}(y_{0})$

$\qquad=-\frac{1}{2}\overset{\text{q}}{\underset{\text{d}=1}{\sum}%
}g^{\text{cd}}(y_{0})\frac{\partial\text{g}_{\text{ab}}}{\partial
\text{x}_{\text{d}}}(y_{0})-\frac{1}{2}\overset{n}{\underset{k=q+1}{\sum}%
}g^{\text{c}k}(y_{0})\frac{\partial\text{g}_{\text{ab}}}{\partial\text{x}_{k}%
}(y_{0})$

Since $\frac{\partial\text{g}_{\text{ab}}}{\partial\text{x}_{\text{d}}}%
(y_{0})=0$ and g$^{\text{c}k}(y_{0})=\delta^{\text{c}k}=0$ for a,b,c,d
=1,...,q and k = q+1,...,n, we have:

\qquad\qquad$\Gamma_{\text{ab}}^{\text{c}}(y_{0})=0$

\vspace{1pt}(iii) $\frac{\partial\Gamma_{\text{ab}}^{\lambda}}{\partial
\text{x}_{\alpha}}(y_{0})=-\frac{1}{2}\overset{n}{\underset{\gamma=1}{\sum}%
}\frac{\partial}{\partial x_{\alpha}}($g$^{\lambda\gamma}\frac{\partial
\text{g}_{\text{ab}}}{\partial\text{x}_{\gamma}})(y_{0})$

$\qquad=-\frac{1}{2}\overset{n}{\underset{\gamma=1}{\sum}}$g$^{\lambda\gamma
}(y_{0})\frac{\partial^{2}\text{g}_{\text{ab}}}{\partial x_{\alpha}%
\partial\text{x}_{\gamma}}(y_{0})-\frac{1}{2}\overset{n}{\underset{\gamma
=1}{\sum}}\frac{\partial g^{\lambda\gamma}}{\partial x_{\alpha}}(y_{0}%
)\frac{\partial\text{g}_{\text{ab}}}{\partial\text{x}_{\gamma}}(y_{0})$

\qquad\qquad\qquad\qquad

\qquad$=-\frac{1}{2}\frac{\partial^{2}\text{g}_{\text{ab}}}{\partial
x_{\alpha}\partial\text{x}_{\lambda}}(y_{0})-\frac{1}{2}%
\overset{q}{\underset{\gamma=1}{\sum}}\frac{\partial g^{\lambda\gamma}%
}{\partial x_{\alpha}}(y_{0})\frac{\partial\text{g}_{\text{ab}}}%
{\partial\text{x}_{\text{c}}}(y_{0})-\frac{1}{2}\overset{n}{\underset{\gamma
=q+1}{\sum}}\frac{\partial g^{\lambda\gamma}}{\partial x_{\alpha}}(y_{0}%
)\frac{\partial\text{g}_{\text{ab}}}{\partial\text{x}_{\gamma}}(y_{0})$

\qquad\qquad For $c=1,...,q$ and $\alpha,\gamma,\lambda=q+1,...,n$

\qquad\qquad\qquad\ \ \ \ $\frac{\partial\text{g}_{\text{ab}}}{\partial
\text{x}_{\text{c}}}(y_{0})=0=\frac{\partial g^{\lambda\gamma}}{\partial
x_{\alpha}}(y_{0})$

by the fact that derivatives with respect the tangential variables
x$_{1},...,$x$_{q}$

vanish, and\ by (iii) of Table 5. Hence,

\qquad$\frac{\partial\Gamma_{\text{ab}}^{\lambda}}{\partial\text{x}_{\alpha}%
}(y_{0})=-\frac{1}{2}\frac{\partial^{2}\text{g}_{\text{ab}}}{\partial
x_{\alpha}\partial\text{x}_{\lambda}}(y_{0})$

$\ \ \ \ \ \ \ \ \ \ \ =\frac{1}{2}$ $\{($R$_{\text{a}\alpha b\lambda}+$
R$_{\text{a}\lambda b\alpha})-$ $\overset{\text{q}}{\underset{\text{c=1}%
}{\sum}}$($T_{\text{ac}\alpha}T_{\text{bc}\lambda}+T_{\text{ac}\lambda
}T_{\text{bc}\alpha})$

$\qquad-\underset{\text{k=q+1}}{\overset{\text{n}}{\sum}}(\perp_{\text{ak}%
\alpha}\perp_{\text{bk}\lambda}+$ $\ \perp_{\text{ak}\lambda}\perp
_{\text{bk}\alpha})\}(y_{0})$

(iv) This is easily deduced from (iii)\ \ 

(v) From $\Gamma_{\text{ab}}^{\text{c}}=-\frac{1}{2}%
\overset{n}{\underset{k=1}{\sum}}$g$^{\text{c}k}\frac{\partial\text{g}%
_{\text{ab}}}{\partial\text{x}_{k}},$ we have:

$\frac{\partial\Gamma_{\text{ab}}^{\text{c}}}{\partial\text{x}_{i}}=-\frac
{1}{2}\overset{n}{\underset{k=1}{\sum}}(\frac{\partial\text{g}^{\text{ck}}%
}{\partial\text{x}_{i}}\frac{\partial\text{g}_{\text{ab}}}{\partial
\text{x}_{k}}+$ g$^{\text{ck}}\frac{\partial^{2}\text{g}_{\text{ab}}}%
{\partial\text{x}_{i}\partial\text{x}_{k}})$

and so,

$\frac{\partial\Gamma_{\text{ab}}^{\text{c}}}{\partial\text{x}_{i}}%
(y_{0})=-\frac{1}{2}\overset{n}{\underset{k=1}{\sum}}(\frac{\partial
\text{g}^{\text{ck}}}{\partial\text{x}_{i}}(y_{0})\frac{\partial
\text{g}_{\text{ab}}}{\partial\text{x}_{k}}(y_{0})+$ $\frac{\partial
^{2}\text{g}_{\text{ab}}}{\partial\text{x}_{i}\partial\text{x}_{\text{c}}%
}(y_{0}))$

$=-\frac{1}{2}\overset{n}{\underset{k=1}{\sum}}\frac{\partial\text{g}%
^{\text{ck}}}{\partial\text{x}_{i}}(y_{0})\frac{\partial\text{g}_{\text{ab}}%
}{\partial\text{x}_{k}}(y_{0})$ = $-\frac{1}{2}%
\overset{n}{\underset{k=q+1}{\sum}}(-\perp_{\text{cik}})(-2$T$_{\text{ab}%
k})(y_{0})$

$\frac{\partial\Gamma_{\text{ab}}^{\text{c}}}{\partial\text{x}_{i}}%
(y_{0})=-\overset{n}{\underset{k=q+1}{\sum}}(\perp_{\text{c}ik})($%
T$_{\text{ab}k})(y_{0})$

(vi) \ $\frac{\partial^{2}\Gamma_{\text{ab}}^{\lambda}}{\partial
\text{x}_{\alpha}^{2}}(y_{0})=-\frac{1}{2}\overset{n}{\underset{\gamma
=1}{\sum}}\frac{\partial^{2}}{\partial x_{\alpha}}($g$^{\lambda\gamma}%
\frac{\partial\text{g}_{\text{ab}}}{\partial\text{x}_{\gamma}})(y_{0})$

$\ =-\frac{1}{2}\overset{n}{\underset{\gamma=1}{\sum}}$g$^{\lambda\gamma
}(y_{0})\frac{\partial^{3}\text{g}_{\text{ab}}}{\partial x_{\alpha}%
^{2}\partial\text{x}_{\gamma}}(y_{0})-\overset{n}{\underset{\gamma=1}{\sum}%
}\frac{\partial g^{\lambda\gamma}}{\partial x_{\alpha}}(y_{0})\frac
{\partial^{2}\text{g}_{\text{ab}}}{\partial x_{\alpha}\partial\text{x}%
_{\gamma}}(y_{0})$

$-\frac{1}{2}\underset{\gamma=1}{\overset{n}{\sum}}\frac{\partial
^{2}g^{\lambda\gamma}}{\partial x_{\alpha}^{2}}(y_{0})\frac{\partial
\text{g}_{\text{ab}}}{\partial\text{x}_{\gamma}}(y_{0})$

$=$ $-\frac{1}{2}\frac{\partial^{3}\text{g}_{\text{ab}}}{\partial x_{\alpha
}^{2}\partial\text{x}_{\lambda}}(y_{0})-\overset{n}{\underset{\gamma
=q+1}{\sum}}\frac{\partial g^{\lambda\gamma}}{\partial x_{\alpha}}(y_{0}%
)\frac{\partial^{2}\text{g}_{\text{ab}}}{\partial x_{\alpha}\partial
\text{x}_{\gamma}}(y_{0})$

$-\frac{1}{2}\underset{\gamma=q+1}{\overset{n}{\sum}}\frac{\partial
^{2}g^{\lambda\gamma}}{\partial x_{\alpha}^{2}}(y_{0})\frac{\partial
\text{g}_{\text{ab}}}{\partial\text{x}_{\gamma}}(y_{0})$

Since $\ \frac{\partial g^{\lambda\gamma}}{\partial x_{\alpha}}(y_{0})=0$
\ for $\alpha,\gamma,\lambda=q+1,...,n$ \ by (ii) of Table 2,

\ \ \ \qquad\qquad\qquad

\vspace{1pt}\qquad$\frac{\partial^{2}\Gamma_{\text{ab}}^{\lambda}}%
{\partial\text{x}_{\alpha}^{2}}(y_{0})=-\frac{1}{2}\frac{\partial^{3}%
\text{g}_{\text{ab}}}{\partial x_{\alpha}^{2}\partial\text{x}_{\lambda}}%
(y_{0})-\frac{1}{2}\underset{\gamma=q+1}{\overset{n}{\sum}}\frac{\partial
^{2}g^{\lambda\gamma}}{\partial x_{\alpha}^{2}}(y_{0})\frac{\partial
\text{g}_{\text{ab}}}{\partial\text{x}_{\gamma}}(y_{0})$

By (iv) of Table 5, (ii) of Table 5 and (iii) of Table 2 , the last

\ equation above becomes:

$\frac{\partial^{2}\Gamma_{\text{ab}}^{\lambda}}{\partial\text{x}_{\alpha}%
^{2}}(y_{0})$ \ = $\overset{n}{\underset{i,j,k=q+1}{\frac{1}{6}\sum}\{}%
2\nabla_{\alpha}$R$_{\alpha\text{a}\lambda\text{b}}+$ R$_{\text{b}%
\lambda\alpha\text{T}_{\text{a}\alpha}}$ $+$R$_{\text{b}\lambda\alpha
\perp_{\text{a}\alpha}}$ $+3$ (R$_{\text{a}\alpha\alpha\text{T}_{\text{b}%
\lambda}}$ $+\;$R$_{\text{a}\alpha\alpha\perp_{\text{b}\lambda}})$

\qquad\qquad\ \ \ $+3($R$_{\text{b}\alpha\alpha\text{T}_{\text{a}\lambda}}$
$+$R$_{\text{b}\alpha\alpha\perp_{\text{a}\lambda}})$ $+$R$_{\text{a}%
\lambda\alpha\text{T}_{\text{b}\alpha}}$ $+$ R$_{\text{a}\lambda\alpha
\perp_{\text{b}\alpha}}$

\ \ \ \ \ \ \ \ \ \ \ \ \ \ \ \ +2$\nabla_{\alpha}$R$_{\lambda\text{a}%
\alpha\text{b}}+$ R$_{\text{b}\alpha\alpha\text{T}_{\text{a}\lambda}}$
$+$R$_{\text{b}\alpha\alpha\perp_{\text{a}\lambda}}$ $+3$ (R$_{\text{a}%
\alpha\lambda\text{T}_{\text{b}\alpha}}$ $+\;$R$_{\text{a}\alpha\lambda
\perp_{\text{b}\alpha}})$

\qquad\qquad\ \ \ $+3($R$_{\text{b}\alpha\lambda\text{T}_{\text{a}\alpha}}$
$+$R$_{\text{b}\alpha\lambda\perp_{\text{a}\alpha}})$ $+$R$_{\text{a}%
\alpha\alpha\text{T}_{\text{b}\lambda}}$ $+$ R$_{\text{a}\alpha\alpha
\perp_{\text{b}\lambda}}$

\qquad\qquad$\ +2\nabla_{\lambda}$R$_{\alpha\text{a}\alpha\text{b}}+$
R$_{\text{b}\alpha\lambda\text{T}_{\text{a}\alpha}}$ $+$R$_{\text{b}%
\alpha\lambda\perp_{\text{a}\alpha}}$ $+3$ (R$_{\text{a}\lambda\alpha
\text{T}_{\text{b}\alpha}}$ $+\;$R$_{\text{a}\lambda\alpha\perp_{\text{b}%
\alpha}})$

\qquad\qquad\ $+3($R$_{\text{b}\lambda\alpha\text{T}_{\text{a}\alpha}}$
$+$R$_{\text{b}\lambda\alpha\perp_{\text{a}\alpha}})$ $+$R$_{\text{a}%
\alpha\lambda\text{T}_{\text{b}\alpha}}$ $+$ R$_{\text{a}\alpha\lambda
\perp_{\text{b}\alpha}}\}(y_{0})$\qquad

\ \qquad\ \ \ \ \ \ +$\frac{2}{3}\underset{\gamma=q+1}{\overset{n}{\sum}%
}\{T_{\text{ab}\gamma}($R$_{\alpha\lambda\alpha\gamma}%
+3\overset{q}{\underset{\text{a}=1}{\sum}}\perp_{\text{a}\alpha\lambda}%
\perp_{\text{a}\alpha\gamma})\}(y_{0})$

$\qquad$R$_{\text{b}ki\text{T}_{\text{a}j}}=$ $<R_{\text{X}_{\text{b}}%
X_{k}X_{i}}^{{}},T_{\text{X}_{\text{a}}}X_{j}>$ $=$ $<R_{\text{b}ki}^{{}%
},T_{\text{a}j}>$

$=R_{\text{b}ki\text{c}}^{{}}T_{\text{a}j\text{d}}<\frac{\partial}%
{\partial\text{x}_{\text{c}}},\frac{\partial}{\partial\text{x}_{\text{d}}}>$
$=R_{\text{b}ki\text{c}}^{{}}T_{\text{a}j\text{d}}\delta_{\text{cd}%
}=R_{\text{b}ki\text{c}}^{{}}T_{\text{a}j\text{c}}=R_{\text{b}k\text{c}i}^{{}%
}T_{\text{ac}j}$

(vii) \ This is immediate from (vi)

(viii) \ For a,b,c=1,..,q we have:

$\qquad\frac{\partial^{2}\Gamma_{\text{ab}}^{\text{c}}}{\partial\text{x}%
_{i}^{2}}(y_{0})$

$\qquad=-\frac{1}{2}\overset{n}{\underset{k=q+1}{\sum}}(\frac{\partial
^{2}\text{g}^{\text{ck}}}{\partial^{2}\text{x}_{i}}\frac{\partial
\text{g}_{\text{ab}}}{\partial\text{x}_{k}}+\frac{\partial\text{g}^{\text{ck}%
}}{\partial\text{x}_{i}}\frac{\partial^{2}\text{g}_{\text{ab}}}{\partial
\text{x}_{i}\partial\text{x}_{k}}+\frac{\partial\text{g}^{\text{ck}}}%
{\partial\text{x}_{i}}\frac{\partial^{2}\text{g}_{\text{ab}}}{\partial
\text{x}_{i}\partial\text{x}_{k}})(y_{0})$

\qquad$\ =$ $-\frac{1}{2}\overset{n}{\underset{k,l=q+1}{\sum}}(\frac
{\partial^{2}\text{g}^{\text{c}k}}{\partial^{2}\text{x}_{i}}\frac
{\partial\text{g}_{\text{ab}}}{\partial\text{x}_{k}}+2\frac{\partial
\text{g}^{\text{c}k}}{\partial\text{x}_{i}}\frac{\partial^{2}\text{g}%
_{\text{ab}}}{\partial\text{x}_{i}\partial\text{x}_{k}})(y_{0})$

and so:

\qquad$\frac{\partial^{2}\Gamma_{\text{ab}}^{\text{c}}}{\partial\text{x}%
_{i}^{2}}(y_{0})$\ $=$ $\frac{1}{2}\overset{n}{\underset{k=q+1}{\sum}}%
2$T$_{\text{ab}k}\left\{  \frac{8}{3}R_{i\text{c}ik}+4\underset{\text{d}%
=1}{\overset{q}{\sum}}(T_{\text{dc}k})(\perp_{\text{d}ik})\right\}  $

$-$ $\overset{n}{\underset{k,l=q+1}{\sum}}\perp_{\text{c}ki}(y_{0})\left\{
-(R_{\text{a}k\text{b}l}+R_{\text{a}l\text{b}k})+\underset{\text{d=1}%
}{\overset{\text{q}}{\sum}}(T_{\text{ad}k}T_{\text{bd}l}+T_{\text{ad}%
l}T_{\text{bd}k})\right\}  (y_{0})$

$-\overset{n}{\underset{k,l=q+1}{\sum}}\perp_{\text{c}ki}(y_{0})\left\{
\underset{r=q+1}{\overset{n}{\sum}}(\perp_{\text{a}kr}\perp_{\text{b}lr}%
+\perp_{\text{a}lr}\perp_{\text{b}kr})\right\}  (y_{0})$

(ix) The expression for $\frac{\partial^{2}\Gamma_{\text{aa}}^{\text{c}}%
}{\partial\text{x}_{i}^{2}}(y_{0})$\ is easily deduced from (viii).

$\frac{\partial^{2}\Gamma_{\text{aa}}^{\text{c}}}{\partial\text{x}_{i}^{2}%
}(y_{0})$\ $=$ $\frac{1}{2}\overset{n}{\underset{k=q+1}{\sum}}2$%
T$_{\text{aa}k}\left\{  \frac{8}{3}R_{i\text{c}ik}+4\underset{\text{d}%
=1}{\overset{q}{\sum}}(T_{\text{dc}k})(\perp_{\text{d}ik})\right\}  $

$-$ $\overset{n}{\underset{k,l=q+1}{\sum}}\perp_{\text{c}ki}(y_{0})\left\{
-(R_{\text{a}k\text{a}l}+R_{\text{a}l\text{a}k})+\underset{\text{d=1}%
}{\overset{\text{q}}{\sum}}(T_{\text{ad}k}T_{\text{ad}l}+T_{\text{ad}%
l}T_{\text{ad}k})\right\}  (y_{0})$

$-\overset{n}{\underset{k,l=q+1}{\sum}}\perp_{\text{c}ki}(y_{0})\left\{
\underset{r=q+1}{\overset{n}{\sum}}(\perp_{\text{a}kr}\perp_{\text{a}lr}%
+\perp_{\text{a}lr}\perp_{\text{a}kr})\right\}  (y_{0})$

\qquad$=$ $\overset{n}{\underset{k=q+1}{\sum}}$T$_{\text{aa}k}\left\{
\frac{8}{3}R_{i\text{c}ik}+4\underset{\text{d}=1}{\overset{q}{\sum}%
}(T_{\text{dc}k})(\perp_{\text{d}ik})\right\}  $

$-$ $2\overset{n}{\underset{k,l=q+1}{\sum}}\perp_{\text{c}ki}(y_{0})\left\{
-R_{\text{a}k\text{a}l}+\underset{\text{d=1}}{\overset{\text{q}}{\sum}%
}T_{\text{ad}k}T_{\text{ad}l}\right\}  (y_{0})$

$-2\overset{n}{\underset{k,l=q+1}{\sum}}\perp_{\text{c}ki}(y_{0})\left\{
\underset{r=q+1}{\overset{n}{\sum}}\perp_{\text{a}kr}\perp_{\text{a}%
lr}\right\}  (y_{0})$

(x) $\ \Gamma_{\beta\text{a}}^{\lambda}(y_{0})=\frac{1}{2}%
\overset{n}{\underset{\gamma=1}{\sum}}$g$^{\lambda\gamma}(y_{0})(\frac
{\partial\text{g}_{\text{a}\gamma}}{\partial\text{x}_{\beta}}+\frac
{\partial\text{g}_{\beta\gamma}}{\partial\text{x}_{\text{a}}}-\frac
{\partial\text{g}_{\beta\text{a}}}{\partial\text{x}_{\gamma}})(y_{0})$

\qquad Since differentiation with respect to x$_{\text{a}}$ for a = 1,...,q
vanish and,

\qquad g$^{\lambda\gamma}(y_{0})=\delta^{\lambda\gamma}$

we have:

\ \ \ \ \ \ $\Gamma_{\beta\text{a}}^{\lambda}(y_{0})=\Gamma_{\text{a}\beta
}^{\lambda}(y_{0})$\ = \ \ $\frac{1}{2}(\frac{\partial\text{g}_{\text{a}%
\lambda}}{\partial\text{x}_{\beta}}-\frac{\partial\text{g}_{\text{a}\beta}%
}{\partial\text{x}_{\lambda}})(y_{0})$\qquad

Hence, by (ii) of Table 2,

\ $\Gamma_{\text{a}\beta}^{\lambda}(y_{0})$ \ =$\frac{1}{2}(\perp
_{\text{a}\beta\lambda}-\perp_{\text{a}\lambda\beta})(y_{0})=\frac{1}{2}%
(\perp_{\text{a}\beta\lambda}-\perp_{\text{a}\lambda\beta})(y_{0})$

$\qquad\qquad=\frac{1}{2}(\perp_{\text{a}\beta\lambda}+\perp_{\text{a}%
\beta\lambda})(y_{0})=$ $\perp_{\text{a}\beta\lambda}$

\ 

(xi) $\frac{\partial\Gamma_{\text{a}\beta}^{\lambda}}{\partial x_{\alpha}%
}(y_{0})=\frac{1}{2}\overset{n}{\underset{\gamma=1}{\sum}}\frac{\partial
g^{\lambda\gamma}}{\partial x_{\alpha}}(y_{0})(\frac{\partial\text{g}%
_{\text{a}\gamma}}{\partial\text{x}_{\beta}}+\frac{\partial\text{g}%
_{\beta\gamma}}{\partial\text{x}_{\text{a}}}-\frac{\partial\text{g}%
_{\text{a}\beta}}{\partial\text{x}_{\gamma}})(y_{0})$

$\qquad\qquad\qquad+\frac{1}{2}\overset{n}{\underset{\gamma=1}{\sum}}%
$g$^{\lambda\gamma}(y_{0})(\frac{\partial^{2}\text{g}_{\text{a}\gamma}%
}{\partial\text{x}_{\alpha}\partial\text{x}_{\beta}}+\frac{\partial
^{2}\text{g}_{\beta\gamma}}{\partial\text{x}_{\alpha}\partial\text{x}%
_{\text{a}}}-\frac{\partial^{2}\text{g}_{\text{a}\beta}}{\partial
\text{x}_{\alpha}\partial\text{x}_{\gamma}})(y_{0})$

Since differentiations with respect to x$_{\text{a}}$ vanish and
g$^{\lambda\gamma}(y_{0})=\delta^{\lambda\gamma})$,

we have:

$\ \frac{\partial\Gamma_{\beta\text{a}}^{\lambda}}{\partial x_{\alpha}}%
(y_{0})\ =\frac{1}{2}\overset{n}{\underset{\gamma=1}{\sum}}\frac{\partial
g^{\lambda\gamma}}{\partial x_{\alpha}}(y_{0})(\frac{\partial\text{g}%
_{\text{a}\gamma}}{\partial\text{x}_{\beta}}-\frac{\partial\text{g}%
_{\text{a}\beta}}{\partial\text{x}_{\gamma}})(y_{0})+\frac{1}{2}%
(\frac{\partial^{2}\text{g}_{\text{a}\lambda}}{\partial\text{x}_{\alpha
}\partial\text{x}_{\beta}}-\frac{\partial^{2}\text{g}_{\text{a}\beta}%
}{\partial\text{x}_{\alpha}\partial\text{x}_{\lambda}})(y_{0})$

\ $=\frac{1}{2}\overset{q}{\underset{\gamma=1}{\sum}}\frac{\partial
g^{\lambda\gamma}}{\partial x_{\alpha}}(y_{0})(\frac{\partial\text{g}%
_{\text{a}\gamma}}{\partial\text{x}_{\beta}}-\frac{\partial\text{g}%
_{\text{a}\beta}}{\partial\text{x}_{\gamma}})(y_{0})+$ $\ \frac{1}%
{2}\overset{n}{\underset{\gamma=q+1}{\sum}}\frac{\partial g^{\lambda\gamma}%
}{\partial x_{\alpha}}(y_{0})(\frac{\partial\text{g}_{\text{a}\gamma}%
}{\partial\text{x}_{\beta}}-\frac{\partial\text{g}_{\text{a}\beta}}%
{\partial\text{x}_{\gamma}})(y_{0})$

$\qquad\ +\frac{1}{2}(\frac{\partial^{2}\text{g}_{\text{a}\lambda}}%
{\partial\text{x}_{\alpha}\partial\text{x}_{\beta}}-\frac{\partial^{2}%
\text{g}_{\text{a}\beta}}{\partial\text{x}_{\alpha}\partial\text{x}_{\lambda}%
})(y_{0})$\qquad\qquad\qquad

\qquad\ \ $=$\ $\frac{1}{2}\overset{q}{\underset{\text{b}=1}{\sum}}%
\frac{\partial g^{\lambda\text{b}}}{\partial x_{\alpha}}(y_{0})\frac
{\partial\text{g}_{\text{ab}}}{\partial\text{x}_{\beta}}(y_{0})+\frac{1}%
{2}(\frac{\partial^{2}\text{g}_{\text{a}\lambda}}{\partial\text{x}_{\alpha
}\partial\text{x}_{\beta}}-\frac{\partial^{2}\text{g}_{\text{a}\beta}%
}{\partial\text{x}_{\alpha}\partial\text{x}_{\lambda}})(y_{0})$\qquad
\qquad\qquad

(since $\frac{\partial\text{g}_{\text{a}\beta}}{\partial\text{x}_{\gamma}%
})(y_{0})=0=\frac{\partial g^{\lambda\gamma}}{\partial x_{\alpha}}(y_{0})$
\ for $\gamma=0,...,q$ \ and $\lambda,\gamma=q+1,...,n)$

\qquad By (iii) of Table 4 and (ii) of Table 5,

\qquad$\frac{\partial\Gamma_{\text{a}\beta}^{\lambda}}{\partial x_{\alpha}%
}(y_{0})=$\ $\frac{1}{2}\overset{q}{\underset{\text{b}=1}{\sum}}\left\{
(-\perp_{\text{b}\alpha\lambda})(-2T_{\text{ab}\beta})\right\}  (y_{0})$

$\qquad+\frac{1}{2}\left\{  -\frac{4}{3}(R_{\alpha\text{a}\beta\lambda
}+R_{\beta\text{a}\alpha\lambda})+\frac{4}{3}(R_{\alpha\text{a}\lambda\beta
}+R_{\lambda\text{a}\alpha\beta})\right\}  (y_{0})$Since,

\qquad$R_{\alpha\text{a}\beta\lambda}=-$ $R_{\text{a}\alpha\beta\lambda},$
$\ R_{\beta\text{a}\alpha\lambda}=-R_{\text{a}\beta\alpha\lambda},$
$R_{\alpha\text{a}\lambda\beta}=R_{\text{a}\alpha\beta\lambda},$ \ 

\ and \ \ 

$\qquad R_{\lambda\text{a}\alpha\beta}=R_{\text{a}\lambda\beta\alpha}$, \ 

we have:

$\frac{\partial\Gamma_{\text{a}\beta}^{\lambda}}{\partial x_{\alpha}}%
(y_{0})=\ \overset{q}{\underset{\text{b}=1}{\sum(}}\perp_{\text{b}%
\alpha\lambda}.T_{\text{ab}\beta})(y_{0})+\frac{2}{3}\left\{  2R_{\text{a}%
\alpha\beta\lambda}+R_{\text{a}\beta\alpha\lambda}+R_{\text{a}\lambda
\beta\alpha}\right\}  (y_{0})$

(xii) This is immediate from (xi) since $\Gamma_{\alpha\beta}^{\gamma}%
=\Gamma_{\beta\alpha}^{\gamma}$

\section{\textbf{Table A}$_{8}$}

\vspace{1pt}For a = 1,...,q and $\ \alpha,\beta,\gamma,\lambda=q+1,...,n$ and
$i,j=q+1,...,n,$ we have:

(i) $\Gamma_{\beta\gamma}^{\lambda}(y_{0})=0$

(ii) $\Gamma_{\beta\gamma}^{\text{a}}(y_{0})=0$

(iii) $\Gamma_{\text{a}j}^{\text{b}}(y_{0})=-\Gamma_{\text{ab}}^{j}(y_{0}%
)=-$T$_{\text{ab}j}(y_{0})$

(iv) $\Gamma_{\text{a}\gamma}^{\lambda}(y_{0})=$ $\perp_{\text{a}\gamma
\lambda}(y_{0})$

(v) $\frac{\partial\Gamma_{jj}^{\text{a}}}{\partial\text{x}_{i}}(y_{0}%
)=\frac{2}{3}R_{\text{a}jij}(y_{0})$

(vi)\ $\frac{\partial\Gamma_{\text{a}j}^{\text{b}}}{\partial\text{x}_{i}%
}(y_{0})$

$=\frac{1}{2}[-R_{\text{a}i\text{b}j}-R_{\text{a}j\text{b}i}%
+\underset{\text{c=1}}{\overset{\text{q}}{\sum}}T_{\text{ac}i}T_{\text{bc}%
j}-3\underset{\text{c=1}}{\overset{\text{q}}{\sum}}T_{\text{ac}j}%
T_{\text{bc}i}$

$+\underset{\text{k=q+1}}{\overset{\text{n}}{\sum}}\perp_{\text{a}i\text{k}%
}\perp_{\text{b}j\text{k}}-\underset{\text{k=q+1}}{\overset{\text{n}}{\sum}%
}\perp_{\text{a}j\text{k}}\perp_{\text{b}i\text{k}}\ ](y_{0})$

(vii) $\frac{\partial\Gamma_{jk}^{l}}{\partial\text{x}_{i}}(y_{0})=$ $-$
$\frac{1}{3}(R_{ikjl}+R_{ijkl})(y_{0})$

\vspace{1pt}(viii) $\frac{\partial\Gamma_{jj}^{k}}{\partial\text{x}_{i}}%
(y_{0})=$ $\frac{2}{3}R_{ijkj}(y_{0})$

\vspace{1pt}(ix) $\frac{\partial^{2}\Gamma_{jj}^{k}}{\partial\text{x}_{i}^{2}%
}(y_{0})=-[\frac{4}{3}\overset{q}{\underset{\text{a}=1}{\sum}}\perp
_{\text{a}ik}R_{ij\text{a}j}+\frac{1}{3}(\nabla_{i}R_{kjij}+\nabla_{j}%
R_{ijik}+\nabla_{k}R_{ijij})](y_{0})$

\bigskip(x) $\frac{\partial^{2}\Gamma_{jj}^{\text{b}}}{\partial\text{x}%
_{i}^{2}}(y_{0})=\frac{8}{3}\overset{q}{\underset{\text{c}=1}{\sum}}%
($T$_{\text{bc}i}R_{ij\text{c}j})(y_{0})+\frac{2}{3}%
\overset{n}{\underset{k=q+1}{\sum}}(\perp_{\text{b}ik}R_{ijjk})(y_{0})$

$\qquad-\frac{1}{6}[4\overset{q}{\underset{\text{c=1}}{%
{\textstyle\sum}
}}$R$_{ij\text{c}i}^{{}}T_{\text{bc}j}+$ $4\overset{n}{\underset{k=q+1}{%
{\textstyle\sum}
}}$R$_{ijik}\perp_{\text{b}jk}+3\nabla_{i}$R$_{j\text{b}ij}$

$\qquad+4\overset{q}{\underset{\text{c=1}}{%
{\textstyle\sum}
}}$R$_{ij\text{c}j}^{{}}T_{\text{bc}i}+$ $4$R$_{ijjk}\perp_{\text{b}ik}%
](y_{0})$

(xi) $\frac{\partial^{2}\Gamma_{\text{aa}}^{\text{c}}}{\partial\text{x}%
_{i}^{2}}(y_{0})=$ $\overset{n}{\underset{k=q+1}{\sum}}$T$_{\text{aa}%
k}\left\{  \frac{8}{3}R_{i\text{c}ik}+4\underset{\text{d}=1}{\overset{q}{\sum
}}(T_{\text{dc}k})(\perp_{\text{d}ik})\right\}  $

$\qquad+$ $2\overset{n}{\underset{k,l=q+1}{\sum}}\perp_{\text{c}ik}%
(y_{0})\left\{  -R_{\text{a}k\text{a}l}+\underset{\text{d=1}%
}{\overset{\text{q}}{\sum}}T_{\text{ad}k}T_{\text{ad}l}\right\}  (y_{0})$

$\qquad+2\overset{n}{\underset{k,l=q+1}{\sum}}\perp_{\text{c}ik}%
(y_{0})\left\{  \underset{r=q+1}{\overset{n}{\sum}}\perp_{\text{a}kr}%
\perp_{\text{a}lr}\right\}  (y_{0})$

$\qquad$\ 

\subsection{\protect\underline{\textbf{Computations}}}

(i)\qquad$\Gamma_{\beta\gamma}^{\lambda}=$ $\ \ \frac{1}{2}%
\overset{n}{\underset{\sigma=1}{\sum}}$g$^{\lambda\sigma}(\frac{\partial
\text{g}_{\gamma\sigma}}{\partial\text{x}_{\beta}}+\frac{\partial
\text{g}_{\beta\sigma}}{\partial\text{x}_{\gamma}}-\frac{\partial
\text{g}_{\beta\gamma}}{\partial\text{x}_{\sigma}})$

\qquad\ \ \ $\Gamma_{\beta\gamma}^{\lambda}(y_{0})$ =$\ \ \ \frac{1}%
{2}\overset{q}{\underset{\text{b}=1}{\sum}}$g$^{\lambda\text{b}}(y_{0}%
)(\frac{\partial\text{g}_{\gamma\sigma}}{\partial\text{x}_{\beta}}%
+\frac{\partial\text{g}_{\beta\sigma}}{\partial\text{x}_{\gamma}}%
-\frac{\partial\text{g}_{\beta\gamma}}{\partial\text{x}_{\sigma}})(y_{0})$

$\qquad\qquad\qquad\ \ \ +$ $\frac{1}{2}\overset{n}{\underset{\sigma
=q+1}{\sum}}$g$^{\lambda\sigma}(y_{0})(\frac{\partial\text{g}_{\gamma\sigma}%
}{\partial\text{x}_{\beta}}+\frac{\partial\text{g}_{\beta\sigma}}%
{\partial\text{x}_{\gamma}}-\frac{\partial\text{g}_{\beta\gamma}}%
{\partial\text{x}_{\sigma}})(y_{0})$

\qquad The first sum on the R.H.S. of the equation above is zero because

\qquad g$^{\lambda\text{b}}(y_{0})=\delta_{\lambda\text{b}}=0.$ This is
because $b=1,...,q$ and $\lambda=q+1,...,n.$

The derivatives\ in the second sum are all zero by (ii) of Table 1. Hence,
$\Gamma_{\beta\gamma}^{\lambda}(y_{0})\ =0.$

(ii) $\Gamma_{\beta\gamma}^{\text{a}}(y_{0})=\frac{1}{2}%
\overset{q}{\underset{\text{b}=1}{\sum}}$g$^{\text{ab}}(y_{0})(\frac
{\partial\text{g}_{\gamma\text{b}}}{\partial\text{x}_{\beta}}+\frac
{\partial\text{g}_{\beta\text{b}}}{\partial\text{x}_{\gamma}}-\frac
{\partial\text{g}_{\beta\gamma}}{\partial\text{x}_{\text{b}}})(y_{0})$

$\qquad\qquad\ +\frac{1}{2}\overset{n}{\underset{\sigma=q+1}{\sum}}%
$g$^{\text{a}\sigma}(y_{0})(\frac{\partial\text{g}_{\beta\sigma}}%
{\partial\text{x}_{\beta}}+\frac{\partial\text{g}_{\gamma\sigma}}%
{\partial\text{x}_{\gamma}}-\frac{\partial\text{g}_{\beta\gamma}}%
{\partial\text{x}_{\sigma}})(y_{0})$

\qquad$=\frac{1}{2}\overset{q}{\underset{\text{b}=1}{\sum}}$g$^{\text{ab}%
}(y_{0})(\frac{\partial\text{g}_{\gamma\text{b}}}{\partial\text{x}_{\beta}%
}+\frac{\partial\text{g}_{\beta\text{b}}}{\partial\text{x}_{\gamma}}%
-\frac{\partial\text{g}_{\beta\gamma}}{\partial\text{x}_{\text{b}}}%
)(y_{0})=\frac{1}{2}\overset{q}{\underset{\text{b}=1}{\sum}}(\frac
{\partial\text{g}_{\gamma\text{a}}}{\partial\text{x}_{\beta}}+\frac
{\partial\text{g}_{\beta\text{a}}}{\partial\text{x}_{\gamma}})(y_{0})$

\qquad$=\frac{1}{2}(\perp_{\text{a}\beta\gamma}+$ $\perp_{\text{a}\gamma\beta
})(y_{0})=\frac{1}{2}(\perp_{\text{a}\beta\gamma}-$ $\perp_{\text{a}%
\beta\gamma})(y_{0})=0$

(iii) \ $\Gamma_{\beta\gamma}^{\lambda}(y_{0})=\frac{1}{2}%
\overset{n}{\underset{\sigma=1}{\sum}}$g$^{\lambda\sigma}(y_{0})(\frac
{\partial\text{g}_{\gamma\sigma}}{\partial\text{x}_{\beta}}+\frac
{\partial\text{g}_{\beta\sigma}}{\partial\text{x}_{\gamma}}-\frac
{\partial\text{g}_{\beta\gamma}}{\partial\text{x}_{\sigma}})(y_{0}%
)=(\frac{\partial\text{g}_{\gamma\lambda}}{\partial\text{x}_{\beta}}%
+\frac{\partial\text{g}_{\beta\lambda}}{\partial\text{x}_{\gamma}}%
-\frac{\partial\text{g}_{\beta\gamma}}{\partial\text{x}_{\lambda}})(y_{0})$

Therefore,

$\qquad\Gamma_{\text{a}j}^{\text{b}}(y_{0})=\frac{1}{2}(\frac{\partial
\text{g}_{j\text{b}}}{\partial\text{x}_{\text{a}}}+\frac{\partial
\text{g}_{\text{ab}}}{\partial\text{x}_{j}}-\frac{\partial\text{g}_{\text{a}%
j}}{\partial\text{x}_{\text{b}}})(y_{0})=\frac{1}{2}(\frac{\partial
\text{g}_{\text{ab}}}{\partial\text{x}_{j}})(y_{0})=-$T$_{\text{ab}j}(y_{0})$

(iv) $\Gamma_{\text{a}\gamma}^{\lambda}(y_{0})=$ $\frac{1}{2}%
\overset{n}{\underset{\sigma=1}{\sum}}$g$^{\lambda\sigma}(y_{0})(\frac
{\partial\text{g}_{\gamma\sigma}}{\partial\text{x}_{\text{a}}}+\frac
{\partial\text{g}_{\text{a}\sigma}}{\partial\text{x}_{\gamma}}-\frac
{\partial\text{g}_{\text{a}\gamma}}{\partial\text{x}_{\sigma}})(y_{0})$

\bigskip\qquad\qquad$\qquad=\frac{1}{2}\overset{n}{\underset{\sigma=1}{\sum}}%
$g$^{\lambda\sigma}(y_{0})(\frac{\partial\text{g}_{\text{a}\sigma}}%
{\partial\text{x}_{\gamma}}-\frac{\partial\text{g}_{\text{a}\gamma}}%
{\partial\text{x}_{\sigma}})(y_{0})$

$\qquad=\frac{1}{2}\overset{n}{\underset{j=q+1}{\sum}}g^{\lambda j}%
(y_{0})(\frac{\partial\text{g}_{\text{a}j}}{\partial\text{x}_{\gamma}}%
-\frac{\partial\text{g}_{\text{a}\gamma}}{\partial\text{x}_{j}})(y_{0}%
)=\frac{1}{2}(\frac{\partial\text{g}_{\text{a}\lambda}}{\partial
\text{x}_{\gamma}}-\frac{\partial\text{g}_{\text{a}\gamma}}{\partial
\text{x}_{\lambda}})(y_{0})$ $=\perp_{\text{a}\gamma\lambda}(y_{0})\qquad$

(v) $\ \Gamma_{jj}^{\text{a}}=\frac{1}{2}\overset{n}{\underset{k=1}{\sum}}%
$g$^{\text{a}k}(2\frac{\partial\text{g}_{j\text{k}}}{\partial\text{x}_{j}%
}-\frac{\partial\text{g}_{jj}}{\partial\text{x}_{k}})$

Therefore,

$\qquad\frac{\partial\Gamma_{jj}^{\text{a}}}{\partial\text{x}_{i}}%
(y_{0})=\ \frac{1}{2}\overset{q}{\underset{\text{b}=1}{\sum}}\frac{\partial
g^{\text{ab}}}{\partial x_{i}}(y_{0})(2\frac{\partial\text{g}_{j\text{b}}%
}{\partial\text{x}_{j}}-\frac{\partial\text{g}_{jj}}{\partial\text{x}%
_{\text{b}}})(y_{0})$

$\qquad\qquad\ \ \ \ +\ \frac{1}{2}\overset{q}{\underset{\text{b=}1}{\sum}}%
$g$^{\text{ab}}(y_{0})(\frac{\partial^{2}\text{g}_{j\text{b}}}{\partial
x_{i}\partial\text{x}_{j}}-\frac{\partial^{2}\text{g}_{jj}}{\partial
x_{i}\partial\text{x}_{\text{b}}})(y_{0})$

\qquad\qquad\qquad+$\ \frac{1}{2}\overset{n}{\underset{k=q+1}{\sum}}%
\frac{\partial g^{\text{ak}}}{\partial x_{i}}(y_{0})(2\frac{\partial
\text{g}_{jk}}{\partial\text{x}_{j}}-\frac{\partial\text{g}_{jj}}%
{\partial\text{x}_{k}})(y_{0})$

$\qquad\qquad\ \ \ \ +\ \frac{1}{2}\overset{n}{\underset{k=q+1}{\sum}}%
$g$^{\text{ak}}(y_{0})(\frac{\partial^{2}\text{g}_{jk}}{\partial x_{i}%
\partial\text{x}_{j}}-\frac{\partial^{2}\text{g}_{jj}}{\partial x_{i}%
\partial\text{x}_{k}})(y_{0})$

\qquad\qquad\qquad= $\frac{1}{2}\overset{q}{\underset{\text{b}=1}{\sum}}%
\frac{\partial g^{\text{ab}}}{\partial x_{i}}(y_{0})(2\frac{\partial
\text{g}_{j\text{b}}}{\partial\text{x}_{j}})(y_{0})+\ \frac{1}{2}%
\overset{q}{\underset{\text{b=}1}{\sum}}$g$^{\text{ab}}(y_{0})(\frac
{\partial^{2}\text{g}_{j\text{b}}}{\partial x_{i}\partial\text{x}_{j}}%
)(y_{0})$

\qquad\qquad\qquad+$\ \frac{1}{2}\overset{n}{\underset{k=q+1}{\sum}}%
\frac{\partial g^{\text{ak}}}{\partial x_{i}}(y_{0})(2\frac{\partial
\text{g}_{jk}}{\partial\text{x}_{j}}-\frac{\partial\text{g}_{jj}}%
{\partial\text{x}_{k}})(y_{0})$

\qquad\qquad\qquad= $\frac{1}{2}\overset{q}{\underset{\text{b}=1}{\sum}}%
\frac{\partial g^{\text{ab}}}{\partial x_{i}}(y_{0})(2\frac{\partial
\text{g}_{j\text{b}}}{\partial\text{x}_{j}})(y_{0})+\ \frac{1}{2}%
(\frac{\partial^{2}\text{g}_{j\text{a}}}{\partial x_{i}\partial\text{x}_{j}%
})(y_{0})$

(since $\frac{\partial\text{g}_{jk}}{\partial\text{x}_{j}}(y_{0}%
)=0=\frac{\partial\text{g}_{jj}}{\partial\text{x}_{k}}(y_{0})$ for j, k = q+1,...,n)

\qquad\qquad\qquad= $\overset{q}{\underset{\text{b}=1}{\sum}}$T$_{\text{abi}%
}(y_{0})\perp_{\text{bjj}}(y_{0})-\ \frac{1}{2}\frac{4}{3}\left\{
R_{\text{iajj}}+R_{\text{jaij}}\right\}  (y_{0})$

\qquad\qquad\qquad= $-\ \frac{2}{3}R_{\text{jaij}}(y_{0})=\frac{2}%
{3}R_{\text{ijaj}}(y_{0})$ since $\perp_{\text{bjj}}(y_{0})=0.$

(vi) \ $\Gamma_{\text{a}j}^{\text{b}}$ =$\ \ \ \frac{1}{2}%
\overset{n}{\underset{k=1}{\sum}}$g$^{\text{b}k}(\frac{\partial\text{g}_{jk}%
}{\partial\text{x}_{\text{a}}}+\frac{\partial\text{g}_{\text{a}k}}%
{\partial\text{x}_{j}}-\frac{\partial\text{g}_{\text{a}j}}{\partial
\text{x}_{k}})=\ \frac{1}{2}\overset{n}{\underset{k=1}{\sum}}$g$^{\text{b}%
k}(\frac{\partial\text{g}_{\text{a}k}}{\partial\text{x}_{j}}-\frac
{\partial\text{g}_{\text{a}j}}{\partial\text{x}_{k}})$

Therefore,

$\frac{\partial\Gamma_{\text{a}j}^{\text{b}}}{\partial\text{x}_{i}}(y_{0})$ =
$\frac{1}{2}\overset{n}{\underset{k=1}{\sum}}\frac{\partial g^{\text{b}k}%
}{\partial\text{x}_{i}}(y_{0})(\frac{\partial\text{g}_{\text{a}k}}%
{\partial\text{x}_{j}}-\frac{\partial\text{g}_{\text{a}j}}{\partial
\text{x}_{k}})(y_{0})+\ \frac{1}{2}\overset{n}{\underset{k=1}{\sum}}%
$g$^{\text{b}k}(y_{0})(\frac{\partial^{2}\text{g}_{\text{a}k}}{\partial
\text{x}_{i}\partial\text{x}_{j}}-\frac{\partial^{2}\text{g}_{\text{a}j}%
}{\partial\text{x}_{i}\partial\text{x}_{k}})(y_{0})$

\qquad=$\frac{1}{2}\overset{q}{\underset{\text{c}=1}{\sum}}\frac{\partial
g^{\text{bc}}}{\partial\text{x}_{i}}(y_{0})(\frac{\partial\text{g}_{\text{ac}%
}}{\partial\text{x}_{j}}-\frac{\partial\text{g}_{\text{a}j}}{\partial
\text{x}_{\text{c}}})(y_{0})+\frac{1}{2}\overset{n}{\underset{k=q+1}{\sum}%
}\frac{\partial g^{\text{b}k}}{\partial\text{x}_{i}}(y_{0})(\frac
{\partial\text{g}_{\text{a}k}}{\partial\text{x}_{j}}-\frac{\partial
\text{g}_{\text{a}j}}{\partial\text{x}_{k}})(y_{0})$

$\qquad+\ \frac{1}{2}(\frac{\partial^{2}\text{g}_{\text{ab}}}{\partial
\text{x}_{i}\partial\text{x}_{j}}-\frac{\partial^{2}\text{g}_{\text{a}j}%
}{\partial\text{x}_{i}\partial\text{x}_{\text{b}}})(y_{0})$

= $\frac{1}{2}\overset{q}{\underset{\text{c}=1}{\sum}}\frac{\partial
g^{\text{bc}}}{\partial\text{x}_{i}}(y_{0})\frac{\partial\text{g}_{\text{ac}}%
}{\partial\text{x}_{j}}(y_{0})+\frac{1}{2}\overset{n}{\underset{k=q+1}{\sum}%
}\frac{\partial g^{\text{b}k}}{\partial\text{x}_{i}}(y_{0})(\frac
{\partial\text{g}_{\text{a}k}}{\partial\text{x}_{j}}-\frac{\partial
\text{g}_{\text{a}j}}{\partial\text{x}_{k}})(y_{0})+\ \frac{1}{2}%
\frac{\partial^{2}\text{g}_{\text{ab}}}{\partial\text{x}_{i}\partial
\text{x}_{j}}(y_{0})$

= $-2\overset{q}{\underset{\text{c}=1}{\sum}}T_{\text{bc}i}(y_{0}%
)T_{\text{ac}j}(y_{0})+\frac{1}{2}\overset{n}{\underset{k=q+1}{\sum}}%
\perp_{\text{b}ki}(y_{0})(\perp_{\text{a}jk}-\perp_{\text{a}kj})(y_{0})$

+$\frac{1}{2}\left\{  -R_{\text{a}i\text{b}j}-R_{\text{a}j\text{b}%
i}+\underset{\text{c=1}}{\overset{\text{q}}{\sum}}T_{\text{ac}i}T_{\text{bc}%
j}+\underset{\text{c=1}}{\overset{\text{q}}{\sum}}T_{\text{ac}j}T_{\text{bc}%
i}+\underset{\text{k=q+1}}{\overset{\text{n}}{\sum}}\perp_{\text{a}i\text{k}%
}\perp_{\text{b}j\text{k}}+\underset{\text{k=q+1}}{\overset{\text{n}}{\sum}%
}\perp_{\text{a}j\text{k}}\perp_{\text{b}i\text{k}}\ \right\}  (y_{0})$

= $-2\overset{q}{\underset{\text{c}=1}{\sum}}T_{\text{bc}i}(y_{0}%
)T_{\text{ac}j}(y_{0})-\overset{n}{\underset{k=q+1}{\sum}}\perp_{\text{b}%
i\text{k}}(y_{0})(\perp_{\text{a}j\text{k}})(y_{0})$

+$\frac{1}{2}\left\{  -R_{\text{a}i\text{b}j}-R_{\text{a}j\text{b}%
i}+\underset{\text{c=1}}{\overset{\text{q}}{\sum}}T_{\text{ac}i}T_{\text{bc}%
j}+\underset{\text{c=1}}{\overset{\text{q}}{\sum}}T_{\text{ac}j}T_{\text{bc}%
i}+\underset{\text{k=q+1}}{\overset{\text{n}}{\sum}}\perp_{\text{a}i\text{k}%
}\perp_{\text{b}j\text{k}}+\underset{\text{k=q+1}}{\overset{\text{n}}{\sum}%
}\perp_{\text{a}j\text{k}}\perp_{\text{b}i\text{k}}\ \right\}  (y_{0})$

= $\frac{1}{2}\left\{  -R_{\text{a}i\text{b}j}-R_{\text{a}j\text{b}%
i}+\underset{\text{c=1}}{\overset{\text{q}}{\sum}}T_{\text{ac}i}T_{\text{bc}%
j}-3\underset{\text{c=1}}{\overset{\text{q}}{\sum}}T_{\text{ac}j}%
T_{\text{bc}i}+\underset{\text{k=q+1}}{\overset{\text{n}}{\sum}}%
\perp_{\text{a}i\text{k}}\perp_{\text{b}j\text{k}}-\underset{\text{k=q+1}%
}{\overset{\text{n}}{\sum}}\perp_{\text{a}j\text{k}}\perp_{\text{b}i\text{k}%
}\ \right\}  (y_{0})$

(vii) $\ \frac{\partial\Gamma_{\beta\gamma}^{\lambda}}{\partial\text{x}%
_{\alpha}}=\ \frac{1}{2}\overset{n}{\underset{\sigma=1}{\sum}}\frac{\partial
g^{\lambda\sigma}}{\partial x_{\alpha}}(y_{0})(\frac{\partial\text{g}%
_{\gamma\sigma}}{\partial\text{x}_{\beta}}+\frac{\partial\text{g}_{\beta
\sigma}}{\partial\text{x}_{\gamma}}-\frac{\partial\text{g}_{\beta\gamma}%
}{\partial\text{x}_{\sigma}})(y_{0})$

$\qquad\qquad\ +\ \frac{1}{2}\overset{n}{\underset{\sigma=1}{\sum}}%
$g$^{\lambda\sigma}(y_{0})(\frac{\partial^{2}\text{g}_{\gamma\sigma}}{\partial
x_{\alpha}\partial\text{x}_{\beta}}+\frac{\partial^{2}\text{g}_{\beta\sigma}%
}{\partial x_{\alpha}\partial\text{x}_{\gamma}}-\frac{\partial^{2}%
\text{g}_{\beta\gamma}}{\partial x_{\alpha}\partial\text{x}_{\sigma}})(y_{0})$

\qquad\ \ \ \ \ \ \ \ $\ =\frac{1}{2}\overset{q}{\underset{\sigma=1}{\sum}%
}\frac{\partial g^{\lambda\sigma}}{\partial x_{\alpha}}(y_{0})(\frac
{\partial\text{g}_{\gamma\sigma}}{\partial\text{x}_{\beta}}+\frac
{\partial\text{g}_{\beta\sigma}}{\partial\text{x}_{\gamma}}-\frac
{\partial\text{g}_{\beta\gamma}}{\partial\text{x}_{\sigma}})(y_{0})$

$\qquad\qquad+\ \frac{1}{2}\overset{n}{\underset{\sigma=q+1}{\sum}}%
\frac{\partial g^{\lambda\sigma}}{\partial x_{\alpha}}(y_{0})(\frac
{\partial\text{g}_{\gamma\sigma}}{\partial\text{x}_{\beta}}+\frac
{\partial\text{g}_{\beta\sigma}}{\partial\text{x}_{\gamma}}-\frac
{\partial\text{g}_{\beta\gamma}}{\partial\text{x}_{\sigma}})(y_{0})$

$\qquad\ \ \ \ \ \ \ +\ \frac{1}{2}\overset{q}{\underset{\sigma=1}{\sum}}%
$g$^{\lambda\sigma}(y_{0})(\frac{\partial^{2}\text{g}_{\gamma\sigma}}{\partial
x_{\alpha}\partial\text{x}_{\beta}}+\frac{\partial^{2}\text{g}_{\beta\sigma}%
}{\partial x_{\alpha}\partial\text{x}_{\gamma}}-\frac{\partial^{2}%
\text{g}_{\beta\gamma}}{\partial x_{\alpha}\partial\text{x}_{\sigma}})(y_{0})$

$\qquad\qquad+\ \frac{1}{2}\overset{n}{\underset{\sigma=q+1}{\sum}}%
$g$^{\lambda\sigma}(y_{0})(\frac{\partial^{2}\text{g}_{\gamma\sigma}}{\partial
x_{\alpha}\partial\text{x}_{\beta}}+\frac{\partial^{2}\text{g}_{\beta\sigma}%
}{\partial x_{\alpha}\partial\text{x}_{\gamma}}-\frac{\partial^{2}%
\text{g}_{\beta\gamma}}{\partial x_{\alpha}\partial\text{x}_{\sigma}})(y_{0})$

\qquad\qquad\ =A + B+ C+ D, \ where,

\qquad A = $\frac{1}{2}\overset{q}{\underset{\sigma=1}{\sum}}\frac{\partial
g^{\lambda\sigma}}{\partial x_{\alpha}}(y_{0})(\frac{\partial\text{g}%
_{\gamma\sigma}}{\partial\text{x}_{\beta}}+\frac{\partial\text{g}_{\beta
\sigma}}{\partial\text{x}_{\gamma}}-\frac{\partial\text{g}_{\beta\gamma}%
}{\partial\text{x}_{\sigma}})(y_{0})$

$\qquad=$ $\frac{1}{2}\overset{q}{\underset{\sigma=1}{\sum}}\frac{\partial
g^{\lambda\sigma}}{\partial x_{\alpha}}(y_{0})(\frac{\partial\text{g}%
_{\gamma\sigma}}{\partial\text{x}_{\beta}}+\frac{\partial\text{g}_{\beta
\sigma}}{\partial\text{x}_{\gamma}})(y_{0})$

\qquad\qquad(derivatives with respect to the tangential variable x$_{\sigma}$
are zero)

\qquad By (ii) of Table 4 and (ii) of Table 3,

\qquad A = $\frac{1}{2}(-\perp_{\sigma\alpha\lambda})(\perp_{\sigma\beta
\gamma}+\perp_{\sigma\gamma\beta})=0$ \ because $\perp_{\sigma\gamma\beta
}=-\perp_{\sigma\beta\gamma}$by Lemma\ (3.4) of $\left[  20\right]  .$

\qquad B =$\ \frac{1}{2}\overset{n}{\underset{\sigma=q+1}{\sum}}\frac{\partial
g^{\lambda\sigma}}{\partial x_{\alpha}}(y_{0})(\frac{\partial\text{g}%
_{\gamma\sigma}}{\partial\text{x}_{\beta}}+\frac{\partial\text{g}_{\beta
\sigma}}{\partial\text{x}_{\gamma}}-\frac{\partial\text{g}_{\beta\gamma}%
}{\partial\text{x}_{\sigma}})(y_{0})=0$ \ because $\frac{\partial
g^{\lambda\sigma}}{\partial x_{\alpha}}(y_{0})=0$

\qquad\ by (ii) of Table 1.

\qquad C = $\ \frac{1}{2}\overset{q}{\underset{\sigma=1}{\sum}}$%
g$^{\lambda\sigma}(y_{0})(\frac{\partial^{2}\text{g}_{\gamma\sigma}}{\partial
x_{\alpha}\partial\text{x}_{\beta}}+\frac{\partial^{2}\text{g}_{\beta\sigma}%
}{\partial x_{\alpha}\partial\text{x}_{\gamma}}-\frac{\partial^{2}%
\text{g}_{\beta\gamma}}{\partial x_{\alpha}\partial\text{x}_{\sigma}}%
)(y_{0})\ $

$\qquad\ =\ \frac{1}{2}\overset{q}{\underset{\sigma=1}{\sum}}\delta
_{\lambda\sigma}(\frac{\partial^{2}\text{g}_{\gamma\sigma}}{\partial
x_{\alpha}\partial\text{x}_{\beta}}+\frac{\partial^{2}\text{g}_{\beta\sigma}%
}{\partial x_{\alpha}\partial\text{x}_{\gamma}}-\frac{\partial^{2}%
\text{g}_{\beta\gamma}}{\partial x_{\alpha}\partial\text{x}_{\sigma}}%
)(y_{0})=0$

\qquad\qquad because $\delta_{\lambda\sigma}=0$ \ for $\sigma=1,...,q$ and
$\lambda=q+1,...,n$

\qquad D = $\frac{1}{2}\overset{n}{\underset{\sigma=q+1}{\sum}}$%
g$^{\lambda\sigma}(y_{0})(\frac{\partial^{2}\text{g}_{\gamma\sigma}}{\partial
x_{\alpha}\partial\text{x}_{\beta}}+\frac{\partial^{2}\text{g}_{\beta\sigma}%
}{\partial x_{\alpha}\partial\text{x}_{\gamma}}-\frac{\partial^{2}%
\text{g}_{\beta\gamma}}{\partial x_{\alpha}\partial\text{x}_{\sigma}})(y_{0})$

$\qquad=$ $\frac{1}{2}\overset{n}{\underset{\sigma=q+1}{\sum}}\delta
^{\lambda\sigma}(\frac{\partial^{2}\text{g}_{\gamma\sigma}}{\partial
x_{\alpha}\partial\text{x}_{\beta}}+\frac{\partial^{2}\text{g}_{\beta\sigma}%
}{\partial x_{\alpha}\partial\text{x}_{\gamma}}-\frac{\partial^{2}%
\text{g}_{\beta\gamma}}{\partial x_{\alpha}\partial\text{x}_{\sigma}})(y_{0})$

\qquad= $\frac{1}{2}(\frac{\partial^{2}\text{g}_{\gamma\lambda}}{\partial
x_{\alpha}\partial\text{x}_{\beta}}+\frac{\partial^{2}\text{g}_{\beta\lambda}%
}{\partial x_{\alpha}\partial\text{x}_{\gamma}}-\frac{\partial^{2}%
\text{g}_{\beta\gamma}}{\partial x_{\alpha}\partial\text{x}_{\lambda}}%
)(y_{0})$.

\ \ \ \ \ By (iii) of Table 1,

\qquad D \ = $-$ $\frac{1}{6}\{(R_{\alpha\gamma\beta\lambda}+R_{\beta
\gamma\alpha\lambda})+(R_{\alpha\beta\gamma\lambda}+R_{\gamma\beta
\alpha\lambda})-(R_{\alpha\beta\lambda\gamma}+R_{\lambda\beta\alpha\gamma})\}$

\qquad Since, $-R_{\alpha\beta\lambda\gamma}=R_{\alpha\beta\gamma\lambda
},-R_{\lambda\beta\alpha\gamma}=-R_{\alpha\gamma\lambda\beta}=R_{\alpha
\gamma\beta\lambda}$ and

$\qquad\qquad\qquad R_{\gamma\beta\alpha\lambda}=-R_{\beta\gamma\alpha\lambda
}$

\qquad D = $-$ $\frac{1}{6}(2R_{\alpha\gamma\beta\lambda}+2R_{\alpha
\beta\gamma\lambda})=$ $-$ $\frac{1}{3}(R_{\alpha\gamma\beta\lambda}%
+R_{\alpha\beta\gamma\lambda})$

(viii) is directly deduced from (vii).

(ix) \ \ $\frac{\partial\ \Gamma_{\beta\beta}^{\lambda}}{\partial x_{\alpha}%
}=\frac{1}{2}\overset{n}{\underset{\gamma=1}{\sum}}\frac{\partial
\text{g}^{\lambda\gamma}}{\partial\text{x}_{\alpha}}(2\frac{\partial
\text{g}_{\beta\gamma}}{\partial\text{x}\beta}-\frac{\partial\text{g}%
_{\beta\beta}}{\partial\text{x}_{\gamma}})+\frac{1}{2}%
\overset{n}{\underset{\gamma=1}{\sum}}$g$^{\lambda\gamma}(2\frac{\partial
^{2}\text{g}_{\beta\gamma}}{\partial x_{\alpha}\partial\text{x}_{\beta}}%
-\frac{\partial^{2}\text{g}_{\beta\beta}}{\partial x_{\alpha}\partial
\text{x}_{\gamma}})$

\vspace{1pt}\qquad Hence,

\ \ \ \ \ $\frac{\partial^{2}\ \Gamma_{\beta\beta}^{\lambda}}{\partial
^{2}x_{\alpha}}(y_{0})$

$=\frac{1}{2}\overset{n}{\underset{\gamma=1}{\sum}}\frac{\partial^{2}%
\text{g}^{\lambda\gamma}}{\partial\text{x}_{\alpha}^{2}}(y_{0})(2\frac
{\partial\text{g}_{\beta\gamma}}{\partial\text{x}_{\beta}}-\frac
{\partial\text{g}_{\beta\beta}}{\partial\text{x}_{\gamma}})(y_{0})+\frac{1}%
{2}\overset{n}{\underset{\gamma=1}{\sum}}\frac{\partial\text{g}^{\lambda
\gamma}}{\partial\text{x}_{\alpha}}(y_{0})(2\frac{\partial^{2}\text{g}%
_{\beta\gamma}}{\partial\text{x}_{\alpha}\partial\text{x}_{\beta}}%
-\frac{\partial^{2}\text{g}_{\beta\beta}}{\partial\text{x}_{\alpha}%
\partial\text{x}_{\gamma}})(y_{0})$

$+\frac{1}{2}\overset{n}{\underset{\gamma=1}{\sum}}$g$^{\lambda\gamma}%
(y_{0})(2\frac{\partial^{3}\text{g}_{\beta\gamma}}{\partial x_{\alpha}%
^{2}\partial\text{x}_{\beta}}-\frac{\partial^{3}\text{g}_{\beta\beta}%
}{\partial x_{\alpha}^{2}\partial\text{x}_{\gamma}})(y_{0})=A+B+C,$ where,

$A=\frac{1}{2}\overset{n}{\underset{\gamma=1}{\sum}}\frac{\partial^{2}%
\text{g}^{\lambda\gamma}}{\partial\text{x}_{\alpha}^{2}}(y_{0})(2\frac
{\partial\text{g}_{\beta\gamma}}{\partial\text{x}_{\beta}}-\frac
{\partial\text{g}_{\beta\beta}}{\partial\text{x}_{\gamma}})(y_{0})=$
$\overset{q}{\underset{\text{a}=1}{\sum}}\frac{\partial^{2}\text{g}%
^{\lambda\text{a}}}{\partial\text{x}_{\alpha}^{2}}(y_{0})\frac{\partial
\text{g}_{\beta\text{a}}}{\partial\text{x}_{\beta}}(y_{0})=0$

\ since $\frac{\partial\text{g}_{\beta\text{a}}}{\partial\text{x}_{\beta}%
}(y_{0})=$ $\perp_{\text{a}\beta\beta}=0$

$B=\frac{1}{2}\overset{n}{\underset{\gamma=1}{\sum}}\frac{\partial
\text{g}^{\lambda\gamma}}{\partial\text{x}_{\alpha}}(y_{0})(2\frac
{\partial^{2}\text{g}_{\beta\gamma}}{\partial\text{x}_{\alpha}\partial
\text{x}_{\beta}}-\frac{\partial^{2}\text{g}_{\beta\beta}}{\partial
\text{x}_{\alpha}\partial\text{x}_{\gamma}})(y_{0})$

$=\frac{1}{2}\overset{q}{\underset{\gamma=1}{\sum}}\frac{\partial
\text{g}^{\lambda\text{a}}}{\partial\text{x}_{\alpha}}(y_{0})(2\frac
{\partial^{2}\text{g}_{\beta\text{a}}}{\partial\text{x}_{\alpha}%
\partial\text{x}_{\beta}}-\frac{\partial^{2}\text{g}_{\beta\beta}}%
{\partial\text{x}_{\alpha}\partial\text{x}_{\text{a}}})(y_{0})$
=$\overset{q}{\underset{\text{a}=1}{\sum}}\frac{\partial\text{g}%
^{\lambda\text{a}}}{\partial\text{x}_{\alpha}}(y_{0})\frac{\partial
^{2}\text{g}_{\beta\text{a}}}{\partial x_{\alpha}\partial\text{x}_{\beta}%
}(y_{0})$

= $\overset{q}{\underset{\text{a}=1}{\sum}}$ $(-\perp_{\text{a}\alpha\lambda
})(-\frac{4}{3}(R_{\alpha\text{a}\beta\beta}+R_{\beta\text{a}\alpha\beta
}))=\frac{4}{3}\overset{q}{\underset{\text{a}=1}{\sum}}\perp_{\text{a}%
\alpha\lambda}R_{\beta\text{a}\alpha\beta}$

$=\frac{4}{3}\overset{q}{\underset{\gamma=1}{\sum}}\perp_{\text{a}%
\lambda\alpha}.R_{\alpha\beta\text{a}\beta}$

$\ C=\frac{1}{2}\overset{n}{\underset{\gamma=1}{\sum}}$g$^{\lambda\gamma
}(y_{0})(2\frac{\partial^{3}\text{g}_{\beta\gamma}}{\partial x_{\alpha}%
^{2}\partial\text{x}_{\beta}}-\frac{\partial^{3}\text{g}_{\beta\beta}%
}{\partial x_{\alpha}^{2}\partial\text{x}_{\gamma}})(y_{0})$

\ $\ =\overset{n}{\underset{\gamma=q+1}{\sum}}$g$^{\lambda\gamma}(y_{0}%
)(\frac{\partial^{3}\text{g}_{\beta\gamma}}{\partial x_{\alpha}^{2}%
\partial\text{x}_{\beta}}-\frac{\partial^{3}\text{g}_{\beta\beta}}{\partial
x_{\alpha}^{2}\partial\text{x}_{\gamma}})(y_{0})$

(\ because g$^{\text{a}\gamma}(y_{0})=\delta^{\text{a}\gamma}=0$ \ for
a=1,......,q and $\gamma=q+1,...,n.)$

\vspace{1pt}

\qquad\ =$($ $\frac{\partial^{3}\text{g}_{\beta\lambda}}{\partial x_{\alpha
}^{2}\partial\text{x}_{\beta}}-\frac{\partial^{3}\text{g}_{\beta\beta}%
}{\partial x_{\alpha}^{2}\partial\text{x}_{\gamma}})(y_{0})=(\frac
{\partial^{3}\text{g}_{\beta\lambda}}{\partial x_{\alpha}^{2}\partial
\text{x}_{\beta}}$ $-$\ $\frac{\partial^{3}\text{g}_{\beta\beta}}{\partial
x_{\alpha}^{2}\partial\text{x}_{\lambda}})(y_{0})$

By (iv) of Table 1,

$\ \ \ \frac{\partial^{3}\text{g}_{\lambda\beta}}{\partial\text{x}_{\alpha
}^{2}\partial\text{x}_{\gamma}}(y_{0})=-\frac{1}{3}(\nabla_{\alpha}%
R_{\alpha\lambda\gamma\beta}+\nabla_{\alpha}R_{\gamma\lambda\alpha\beta
}+\nabla_{\gamma}R_{\alpha\lambda\alpha\beta})(y_{0})$

Hence,

$\frac{\partial^{3}\text{g}_{\beta\lambda}}{\partial x_{\alpha}^{2}%
\partial\text{x}_{\beta}}$ $-$\ $\frac{\partial^{3}\text{g}_{\beta\beta}%
}{\partial x_{\alpha}^{2}\partial\text{x}_{\lambda}})(y_{0})$

= $-\frac{1}{3}(\nabla_{\alpha}R_{\alpha\lambda\beta\beta}+\nabla_{\alpha
}R_{\beta\lambda\alpha\beta}+\nabla_{\beta}R_{\alpha\lambda\alpha\beta}%
)(y_{0})$

\ \ $-\frac{1}{3}(\nabla_{\alpha}R_{\alpha\beta\lambda\beta}+\nabla_{\alpha
}R_{\lambda\beta\alpha\beta}+\nabla_{\lambda}R_{\alpha\beta\alpha\beta}%
)(y_{0})$

= $-\frac{1}{3}(\nabla_{\alpha}R_{\alpha\beta\beta\lambda}+\nabla_{\beta
}R_{\alpha\beta\alpha\lambda}+\nabla_{\alpha}R_{\alpha\beta\lambda\beta
}+\nabla_{\alpha}R_{\lambda\beta\alpha\beta}+\nabla_{\lambda}R_{\alpha
\beta\alpha\beta})$

= $-\frac{1}{3}(-\nabla_{\alpha}R_{\alpha\beta\lambda\beta}+\nabla_{\beta
}R_{\alpha\beta\alpha\lambda}+\nabla_{\alpha}R_{\alpha\beta\lambda\beta
}+\nabla_{\alpha}R_{\lambda\beta\alpha\beta}+\nabla_{\lambda}R_{\alpha
\beta\alpha\beta})$

= $-\frac{1}{3}(\nabla_{\alpha}R_{\lambda\beta\alpha\beta}+\nabla_{\beta
}R_{\alpha\beta\alpha\lambda}+\nabla_{\lambda}R_{\alpha\beta\alpha\beta}%
)$\qquad\qquad\qquad\qquad$\ \ \ $\qquad\ 

$\frac{\partial^{2}\ \Gamma_{\beta\beta}^{\lambda}}{\partial x_{\alpha}^{2}%
}(y_{0})=\frac{4}{3}\overset{n}{\underset{\lambda=q+1}{\sum}}%
\overset{q}{\underset{\text{a}=1}{\sum}}[\perp_{\text{a}\lambda\alpha
}.R_{\alpha\beta\text{a}\beta}-\frac{1}{3}(\nabla_{\alpha}R_{\lambda
\beta\alpha\beta}+\nabla_{\beta}R_{\alpha\beta\alpha\lambda}+\nabla_{\lambda
}R_{\alpha\beta\alpha\beta})](y_{0})$

(x) $\frac{\partial^{2}\ \Gamma_{\beta\beta}^{\text{b}}}{\partial x_{\alpha
}^{2}}(y_{0})$

$=\frac{1}{2}\overset{n}{\underset{\gamma=1}{\sum}}\frac{\partial^{2}%
\text{g}^{\text{b}\gamma}}{\partial\text{x}_{\alpha}^{2}}(y_{0})(2\frac
{\partial\text{g}_{\beta\gamma}}{\partial\text{x}_{\beta}}-\frac
{\partial\text{g}_{\beta\beta}}{\partial\text{x}_{\gamma}})(y_{0})$

$+\frac{1}{2}\overset{n}{\underset{\gamma=1}{\sum}}\frac{\partial
\text{g}^{\text{b}\gamma}}{\partial\text{x}_{\alpha}}(y_{0})(2\frac
{\partial^{2}\text{g}_{\beta\gamma}}{\partial\text{x}_{\alpha}\partial
\text{x}_{\beta}}-\frac{\partial^{2}\text{g}_{\beta\beta}}{\partial
\text{x}_{\alpha}\partial\text{x}_{\gamma}})(y_{0})$

$+\frac{1}{2}\overset{n}{\underset{\gamma=1}{\sum}}$g$^{\text{b}\gamma}%
(y_{0})(2\frac{\partial^{3}\text{g}_{\beta\gamma}}{\partial x_{\alpha}%
^{2}\partial\text{x}_{\beta}}-\frac{\partial^{3}\text{g}_{\beta\beta}%
}{\partial x_{\alpha}^{2}\partial\text{x}_{\gamma}})(y_{0})=A+B+C,$

where,

$A=\frac{1}{2}\overset{n}{\underset{\gamma=1}{\sum}}\frac{\partial^{2}%
\text{g}^{\text{b}\gamma}}{\partial\text{x}_{\alpha}^{2}}(y_{0})(2\frac
{\partial\text{g}_{\beta\gamma}}{\partial\text{x}_{\beta}}-\frac
{\partial\text{g}_{\beta\beta}}{\partial\text{x}_{\gamma}})(y_{0})$

$=\overset{n}{\underset{\gamma=q+1}{\sum}}\frac{\partial^{2}\text{g}%
^{\text{b}\gamma}}{\partial\text{x}_{\alpha}^{2}}(y_{0})(\frac{\partial
\text{g}_{\beta\gamma}}{\partial\text{x}_{\beta}}-\frac{\partial
\text{g}_{\beta\beta}}{\partial\text{x}_{\gamma}})(y_{0}%
)+\overset{q}{\underset{\text{a}=1}{\sum}}\frac{\partial^{2}\text{g}%
^{\text{ba}}}{\partial\text{x}_{\alpha}^{2}}(y_{0})\frac{\partial
\text{g}_{\beta\text{a}}}{\partial\text{x}_{\beta}}(y_{0})=0$

\ since $\frac{\partial\text{g}_{\beta\gamma}}{\partial\text{x}_{\beta}}%
(y_{0})=0=\frac{\partial\text{g}_{\beta\beta}}{\partial\text{x}_{\gamma}%
}(y_{0})$ $\ $and $\frac{\partial\text{g}_{\beta\text{a}}}{\partial
\text{x}_{\beta}}(y_{0})=$ $\perp_{\text{a}\beta\beta}=0$

B = $\frac{1}{2}\overset{n}{\underset{\gamma=1}{\sum}}\frac{\partial
\text{g}^{\text{b}\gamma}}{\partial\text{x}_{\alpha}}(y_{0})(2\frac
{\partial^{2}\text{g}_{\beta\gamma}}{\partial\text{x}_{\alpha}\partial
\text{x}_{\beta}}-\frac{\partial^{2}\text{g}_{\beta\beta}}{\partial
\text{x}_{\alpha}\partial\text{x}_{\gamma}})(y_{0})$

=$\frac{1}{2}\overset{q}{\underset{c=1}{\sum}}\frac{\partial\text{g}%
^{\text{bc}}}{\partial\text{x}_{\alpha}}(y_{0})(2\frac{\partial^{2}%
\text{g}_{\beta\text{c}}}{\partial\text{x}_{\alpha}\partial\text{x}_{\beta}%
}-\frac{\partial^{2}\text{g}_{\beta\beta}}{\partial\text{x}_{\alpha}%
\partial\text{x}_{\text{c}}})(y_{0})+\frac{1}{2}\overset{n}{\underset{\gamma
=q+1}{\sum}}\frac{\partial\text{g}^{\text{b}\gamma}}{\partial\text{x}_{\alpha
}}(y_{0})(2\frac{\partial^{2}\text{g}_{\beta\gamma}}{\partial\text{x}_{\alpha
}\partial\text{x}_{\beta}}-\frac{\partial^{2}\text{g}_{\beta\beta}}%
{\partial\text{x}_{\alpha}\partial\text{x}_{\gamma}})(y_{0})$

=$\frac{1}{2}\overset{q}{\underset{c=1}{\sum}}\frac{\partial\text{g}%
^{\text{bc}}}{\partial\text{x}_{\alpha}}(y_{0})(2\frac{\partial^{2}%
\text{g}_{\beta\text{c}}}{\partial\text{x}_{\alpha}\partial\text{x}_{\beta}%
})(y_{0})+\frac{1}{2}\overset{n}{\underset{\gamma=q+1}{\sum}}\frac
{\partial\text{g}^{\text{b}\gamma}}{\partial\text{x}_{\alpha}}(y_{0}%
)(2\frac{\partial^{2}\text{g}_{\beta\gamma}}{\partial\text{x}_{\alpha}%
\partial\text{x}_{\beta}}-\frac{\partial^{2}\text{g}_{\beta\beta}}%
{\partial\text{x}_{\alpha}\partial\text{x}_{\gamma}})(y_{0})$

We use (iii) of Table 1 to have:

B =$\frac{8}{3}\overset{q}{\underset{\text{c}=1}{\sum}}$T$_{\text{bc}\alpha
}(y_{0})R_{\alpha\beta\text{c}\beta}(y_{0})-\frac{1}{2}\frac{1}{3}%
\overset{n}{\underset{\gamma=q+1}{\sum}}\{\perp_{\text{b}\gamma\alpha
}(2R_{\alpha\beta\beta\gamma}+R_{\beta\beta\alpha\gamma}-R_{\alpha\beta
\gamma\beta}-R_{\gamma\beta\alpha\beta})\}(y_{0})$

\qquad=$\frac{8}{3}\overset{q}{\underset{\text{c}=1}{\sum}}$T$_{\text{bc}%
\alpha}(y_{0})R_{\alpha\beta\text{c}\beta}(y_{0})-\frac{1}{2}\frac{1}%
{3}\overset{n}{\underset{\gamma=q+1}{\sum}}\{\perp_{\text{b}\gamma\alpha
}(-2R_{\alpha\beta\gamma\beta}-2R_{\alpha\beta\gamma\beta})\}(y_{0})$

\qquad=$\frac{8}{3}\overset{q}{\underset{\text{c}=1}{\sum}}$T$_{\text{bc}%
\alpha}(y_{0})R_{\alpha\beta\text{c}\beta}(y_{0})+\frac{1}{3}%
\overset{n}{\underset{\gamma=q+1}{\sum}}\{\perp_{\text{b}\gamma\alpha
}(R_{\alpha\beta\gamma\beta}+R_{\alpha\beta\gamma\beta})\}(y_{0})$

\qquad=$\frac{8}{3}\overset{q}{\underset{\text{c}=1}{\sum}}$T$_{\text{bc}%
\alpha}(y_{0})R_{\alpha\beta\text{c}\beta}(y_{0})+\frac{2}{3}%
\overset{n}{\underset{\gamma=q+1}{\sum}}\{\perp_{\text{b}\gamma\alpha
}R_{\alpha\beta\gamma\beta}\}(y_{0})$

C = $\frac{1}{2}\overset{n}{\underset{\gamma=1}{\sum}}$g$^{\text{b}\gamma
}(y_{0})(2\frac{\partial^{3}\text{g}_{\beta\gamma}}{\partial x_{\alpha}%
^{2}\partial\text{x}_{\beta}}-\frac{\partial^{3}\text{g}_{\beta\beta}%
}{\partial x_{\alpha}^{2}\partial\text{x}_{\gamma}})(y_{0})$

\qquad=$\frac{1}{2}\overset{q}{\underset{\text{c}=1}{\sum}}$g$^{\text{bc}%
}(y_{0})(2\frac{\partial^{3}\text{g}_{\beta\text{c}}}{\partial x_{\alpha}%
^{2}\partial\text{x}_{\beta}}-\frac{\partial^{3}\text{g}_{\beta\beta}%
}{\partial x_{\alpha}^{2}\partial\text{x}_{\text{c}}})(y_{0}%
)=\overset{q}{\underset{\text{c}=1}{\sum}}$g$^{\text{bc}}(y_{0})\frac
{\partial^{3}\text{g}_{\beta\text{c}}}{\partial x_{\alpha}^{2}\partial
\text{x}_{\beta}}(y_{0})$

\qquad=$\frac{\partial^{3}\text{g}_{\beta\text{b}}}{\partial x_{\alpha}%
^{2}\partial\text{x}_{\beta}}(y_{0})$

\qquad= $-\frac{1}{6}\{4$R$_{\alpha\beta\alpha\text{T}_{\text{b}\beta}}+$
$4$R$_{\alpha\beta\alpha\perp_{\text{b}\beta}}+3\nabla_{\alpha}$%
R$_{\beta\text{b}\alpha\beta}+4$R$_{\alpha\beta\beta\text{T}_{\text{b}\alpha}%
}+$ $4$R$_{\alpha\beta\beta\perp_{\text{b}\alpha}}$

\qquad\ $\ \ \ \ +3\nabla_{\beta}$R$_{\alpha\text{b}\alpha\beta}\}(y_{0})$

Use the formula:

$\frac{\partial^{2}\Gamma_{jj}^{\text{b}}}{\partial\text{x}_{i}^{2}}%
(y_{0})=\frac{8}{3}\overset{q}{\underset{\text{c}=1}{\sum}}($T$_{\text{bc}%
i}R_{ij\text{c}j})(y_{0})+\frac{2}{3}\overset{n}{\underset{k=q+1}{\sum}}%
(\perp_{\text{b}ki}R_{ijkj})(y_{0})$

$\qquad\qquad-\frac{1}{6}[4$R$_{iji\text{T}_{\text{b}j}}+$ $4$R$_{iji\perp
_{\text{b}j}}+3\nabla_{i}$R$_{j\text{b}ij}+4$R$_{ijj\text{T}_{\text{b}i}}+$
$4$R$_{ijj\perp_{\text{b}i}}](y_{0})$

\qquad\qquad R$_{\text{a}ij\text{T}_{\text{b}k}}%
=\overset{q}{\underset{\text{c=1}}{%
{\textstyle\sum}
}}R_{\text{a}i\text{c}j}^{{}}T_{\text{bc}k}\qquad$R$_{\text{a}ij\perp
_{\text{b}k}}=\overset{n}{\underset{l=q+1}{%
{\textstyle\sum}
}}R_{\text{a}ijl}\perp_{\text{b}kl}$

Therefore,

$\frac{\partial^{2}\Gamma_{jj}^{\text{b}}}{\partial\text{x}_{i}^{2}}%
(y_{0})=\frac{8}{3}\overset{q}{\underset{\text{c}=1}{\sum}}($T$_{\text{bc}%
i}R_{ij\text{c}j})(y_{0})+\frac{2}{3}\overset{n}{\underset{k=q+1}{\sum}}%
(\perp_{\text{b}ik}R_{ijjk})(y_{0})$

$-\frac{1}{6}[4\overset{q}{\underset{\text{c=1}}{%
{\textstyle\sum}
}}$R$_{ij\text{c}i}^{{}}T_{\text{bc}j}+$ $4\overset{n}{\underset{k=q+1}{%
{\textstyle\sum}
}}$R$_{ijik}\perp_{\text{b}jk}+3\nabla_{i}$R$_{j\text{b}ij}%
+4\overset{q}{\underset{\text{c=1}}{%
{\textstyle\sum}
}}$R$_{ij\text{c}j}^{{}}T_{\text{bc}i}+$ $4$R$_{ijjk}\perp_{\text{b}ik}%
](y_{0})\frac{\partial\phi}{\partial\text{x}_{\text{b}}}(y_{0})$

\section{\textbf{Table A}$_{9}$}

The computations in this Table use mostly the expansion formula given by
\textbf{Proposition 6.5.}

We recall that since all expansions are carried out in normal Fermi
coordinates, all derivatives with respect to tangential Fermi coordinates vanish.

For $i,j,k=q+1,...,n,$ we have:

(i)$\ \ \ \ \theta(y_{0})$ \ \ $=1$

(ii)$\ \frac{\partial\theta}{\partial\text{x}_{i}}(y_{0})=-<H,i>(y_{0})$

(iii) $\frac{\partial\theta^{\frac{1}{2}}}{\partial\text{x}_{i}}%
(y_{0})=\ \frac{1}{2}\frac{\partial\theta}{\partial\text{x}_{i}}(y_{0}%
)=-\frac{1}{2}<H,i>(y_{0})$

(iii)$^{\ast}$ $(\nabla\log\theta^{-\frac{1}{2}})_{\text{a}}(y_{0})=0$

(iv) $\frac{\partial\theta^{-\frac{1}{2}}}{\partial\text{x}_{i}}(y_{0}%
)=\frac{1}{2}<H,i>(y_{0})$

\vspace{1pt}(iv)$^{\ast}$ $(\nabla\log\theta^{-\frac{1}{2}})_{i}(y_{0}%
)=\frac{1}{2}<H,i>(y_{0})$

(v) $\frac{\partial^{2}\theta}{\partial\text{x}_{i}\partial\text{x}_{j}}%
(y_{0})$

$=-\frac{1}{6}[2\varrho_{ij}+4\overset{q}{\underset{\text{a}=1}{\sum}%
}R_{i\text{a}j\text{a}}-3\overset{q}{\underset{\text{a,b=1}}{\sum}%
}(T_{\text{aa}i}T_{\text{bb}j}-T_{\text{ab}i}T_{\text{ab}j}%
)-3\overset{q}{\underset{\text{a,b=1}}{\sum}}(T_{\text{aa}j}T_{\text{bb}%
i}-T_{\text{ab}j}T_{\text{ab}i})](y_{0})$

$\frac{\partial^{2}\theta}{\partial\text{x}_{i}^{2}}(y_{0})=-\frac{1}%
{6}[2\varrho_{ii}+4\overset{q}{\underset{\text{a=1}}{\sum}}R_{i\text{a}%
i\text{a}}-6\overset{q}{\underset{\text{a,b=1}}{\sum}}(T_{\text{aa}%
i}T_{\text{bb}i}-T_{\text{ab}i}T_{\text{ab}i})](y_{0})$

\qquad$\frac{\partial^{2}\theta}{\partial\text{x}_{i}^{2}}(y_{0})=-\frac{1}%
{3}[\varrho_{ii}+2\overset{q}{\underset{\text{a=1}}{\sum}}R_{i\text{a}%
i\text{a}}-3\overset{q}{\underset{\text{a,b=1}}{\sum}}(T_{\text{aa}%
i}T_{\text{bb}i}-T_{\text{ab}i}T_{\text{ab}i})](y_{0})$

(vi) $\frac{\partial^{2}\theta}{\partial\text{x}_{i}^{2}}(y_{0})=-\frac{1}%
{3}[\tau^{M}-3\tau^{P}+\overset{q}{\underset{\text{a}=1}{\sum}}\varrho
_{\text{aa}}^{M}+\overset{q}{\underset{\text{a,b}=1}{\sum}}R_{\text{abab}}%
^{M}](y_{0})$

(vii) $\frac{\partial^{2}\theta^{\frac{1}{2}}}{\partial\text{x}_{i}%
\partial\text{x}_{j}}(y_{0})=-\frac{1}{4}\frac{\partial\theta}{\partial
\text{x}_{i}}(y_{0})\frac{\partial\theta}{\partial\text{x}_{j}}(y_{0}%
)+\frac{1}{2}\frac{\partial^{2}\theta}{\partial\text{x}_{i}\partial
\text{x}_{j}}(y_{0})$

\qquad$=-\frac{1}{4}<H,i><H,j>$\ $-\frac{1}{12}[2\varrho_{ij}%
+4\overset{q}{\underset{\text{a}=1}{\sum}}R_{i\text{a}j\text{a}}%
-3\overset{q}{\underset{\text{a,b=1}}{\sum}}(T_{\text{aa}i}T_{\text{bb}%
j}+T_{\text{aa}j}T_{\text{bb}i}-2T_{\text{ab}i}T_{\text{ab}j}](y_{0})$

(viii) $\frac{\partial^{2}\theta^{\frac{1}{2}}}{\partial\text{x}_{i}^{2}%
}(y_{0})=-\frac{1}{4}(\frac{\partial\theta}{\partial\text{x}_{i}})^{2}%
(y_{0})+\frac{1}{2}\frac{\partial^{2}\theta}{\partial\text{x}_{i}^{2}}(y_{0})$

\qquad\qquad$=$ $-\frac{1}{4}<H,i>^{2}(y_{0})-\frac{1}{3}[\varrho
_{ii}+2\overset{q}{\underset{\text{a=1}}{\sum}}R_{i\text{a}i\text{a}%
}-3\overset{q}{\underset{\text{a,b=1}}{\sum}}(T_{\text{aa}i}T_{\text{bb}%
i}-T_{\text{ab}i}T_{\text{ab}i})](y_{0})$\qquad

\bigskip(ix)\qquad\ $\frac{\partial^{2}\theta^{-\frac{1}{2}}}{\partial
\text{x}_{i}\partial\text{x}_{j}}(y_{0})=\frac{3}{4}\frac{\partial\theta
}{\partial\text{x}_{i}}(y_{0})\frac{\partial\theta}{\partial\text{x}_{j}%
}(y_{0})-\frac{1}{2}\frac{\partial^{2}\theta}{\partial\text{x}_{i}%
\partial\text{x}_{j}}(y_{0})$

\ $=\frac{3}{4}<H,i><H,j>$

$+\frac{1}{12}[2\varrho_{ij}+4\overset{q}{\underset{\text{a}=1}{\sum}%
}R_{i\text{a}j\text{a}}-3\overset{q}{\underset{\text{a,b=1}}{\sum}%
}(T_{\text{aa}i}T_{\text{bb}j}-T_{\text{ab}i}T_{\text{ab}j}%
)-3\overset{q}{\underset{\text{a,b=1}}{\sum}}(T_{\text{aa}j}T_{\text{bb}%
i}-T_{\text{ab}j}T_{\text{ab}i})](y_{0})$

(ix)$^{\ast}\qquad\frac{\partial}{\partial x_{i}}(\nabla\log\theta^{-\frac
{1}{2}})_{\text{a}}(y_{0})=-\frac{1}{2}\overset{q}{\underset{j=q+1}{\sum}%
}\perp_{\text{a}ij}(y_{0})<H,j>(y_{0})$

(ix)$^{\ast\ast}$\qquad$\frac{\partial}{\partial x_{i}}(\nabla\log
\theta^{-\frac{1}{2}})_{j}(y_{0})$ for $i,j=q+1,...n$

\qquad$\qquad=\frac{1}{2}<H,i><H,j>\ $

\qquad\qquad$+\frac{1}{12}[2\varrho_{ij}+4\overset{q}{\underset{\text{a}%
=1}{\sum}}R_{i\text{a}j\text{a}}-3\overset{q}{\underset{\text{a,b=1}}{\sum}%
}(T_{\text{aa}i}T_{\text{bb}j}-T_{\text{ab}i}T_{\text{ab}j}%
)-3\overset{q}{\underset{\text{a,b=1}}{\sum}}(T_{\text{aa}j}T_{\text{bb}%
i}-T_{\text{ab}j}T_{\text{ab}i})](y_{0})$

\vspace{1pt}(x)\qquad\ $\frac{\partial^{2}\theta^{-\frac{1}{2}}}%
{\partial\text{x}_{i}^{2}}(y_{0})=\frac{3}{4}(\frac{\partial\theta}%
{\partial\text{x}_{i}})^{2}(y_{0})-\frac{1}{2}\frac{\partial^{2}\theta
}{\partial\text{x}_{i}^{2}}(y_{0})$

\qquad\qquad$=\frac{3}{4}<H,i>^{2}(y_{0})+\frac{1}{6}(\tau^{M}-3\tau
^{P}+\overset{q}{\underset{\text{a}=1}{\sum}}\varrho_{\text{aa}}%
^{M}+\overset{q}{\underset{\text{a,b}=1}{\sum}}R_{\text{abab}}^{M})(y_{0})$

\qquad\qquad$=\frac{1}{12}[9<H,i>^{2}+2(\tau^{M}-3\tau^{P}%
+\overset{q}{\underset{\text{a=1}}{\sum}}\varrho_{\text{aa}}^{M}%
+\overset{q}{\underset{\text{a,b=1}}{\sum}}R_{\text{abab}}^{M})](y_{0})$

(xi)\qquad\ $\frac{\partial^{3}\theta}{\partial\text{x}_{i}\partial
\text{x}_{j}\partial\text{x}_{k}}(y_{0})$

$\qquad=-\frac{1}{12}\{\nabla_{i}\varrho_{jk}-2\varrho_{ij}%
<H,k>+\overset{\text{q}}{\underset{\text{a=1}}{\sum}}(\nabla_{i}%
R_{\text{a}j\text{a}k}-4R_{i\text{a}j\text{a}}<H,k>)$

$\qquad+4\overset{\text{q}}{\underset{\text{a,b=1}}{\sum}}R_{i\text{a}%
j\text{b}}T_{\text{ab}k}$

$+2\overset{q}{\underset{\text{a,b,c=1}}{\sum}}(T_{\text{aa}i}T_{\text{bb}%
j}T_{\text{cc}k}-3T_{\text{aa}i}T_{\text{bc}j}T_{\text{bc}k}+2T_{\text{ab}%
i}T_{\text{bc}j}T_{\text{ca}k})\}(y_{0})$\qquad\qquad\qquad\qquad\qquad\ \ 

$-\frac{1}{12}\{\nabla_{j}\varrho_{ik}-2\varrho_{ji}%
<H,k>+\overset{q}{\underset{\text{a}=1}{\sum}}(\nabla_{j}R_{\text{a}%
i\text{a}k}-4R_{j\text{a}i\text{a}}<H,k>)+4\overset{q}{\underset{\text{a,b=1}%
}{\sum}}R_{j\text{a}i\text{b}}T_{\text{ab}k}$

$+2\overset{q}{\underset{\text{a,b,c=1}}{\sum}}(T_{\text{aa}j}T_{\text{bb}%
i}T_{\text{cc}k}-3T_{\text{aa}j}T_{\text{bc}i}T_{\text{bc}k}+2T_{\text{ab}%
j}T_{\text{bc}i}T_{\text{ca}k})\}(y_{0})$

$-\frac{1}{12}\{\nabla_{i}\varrho_{kj}-2\varrho_{ik}%
<H,j>+\overset{q}{\underset{\text{a}=1}{\sum}}(\nabla_{i}R_{\text{a}%
k\text{a}j}-4R_{i\text{a}k\text{a}}<H,j>)$

$+4\overset{q}{\underset{\text{a,b=1}}{\sum}}R_{i\text{a}k\text{b}%
}T_{\text{ab}j}+2\overset{q}{\underset{\text{a,b,c=1}}{\sum}}(T_{\text{aa}%
i}T_{\text{bb}k}T_{\text{cc}j}-3T_{\text{aa}i}T_{\text{bc}k}T_{\text{bc}%
j}+2T_{\text{ab}i}T_{\text{bc}k}T_{\text{ca}j})\}(y_{0})$

$-\frac{1}{12}\{\nabla_{j}\varrho_{ki}-2\varrho_{jk}%
<H,i>+\overset{q}{\underset{\text{a}=1}{\sum}}(\nabla_{j}R_{\text{a}%
k\text{a}i}-4R_{j\text{a}k\text{a}}<H,i>)$

$+4\overset{q}{\underset{\text{a,b=1}}{\sum}}R_{j\text{a}j\text{b}%
}T_{\text{ab}i}+2\overset{q}{\underset{\text{a,b,c=1}}{\sum}}(T_{\text{aa}%
j}T_{\text{bb}k}T_{\text{cc}i}-3T_{\text{aa}j}T_{\text{bc}k}T_{\text{bc}%
i}+2T_{\text{ab}j}T_{\text{bc}k}T_{\text{ca}i})\}(y_{0})$\qquad\qquad
\qquad\qquad$\ $

$-\frac{1}{12}\{\nabla_{k}\varrho_{ij}-2\varrho_{ki}%
<H,j>+\overset{q}{\underset{\text{a}=1}{\sum}}(\nabla_{k}R_{\text{a}%
i\text{a}j}-4R_{k\text{a}i\text{a}}<H,j>)$

$+4\overset{q}{\underset{\text{a,b=1}}{\sum}}R_{k\text{a}i\text{b}%
}T_{\text{ab}j}+2\overset{q}{\underset{\text{a,b,c=1}}{\sum}}(T_{\text{aa}%
k}T_{\text{bb}i}T_{\text{cc}j}-3T_{\text{aa}k}T_{\text{bc}i}T_{\text{bc}%
j}+2T_{\text{ab}k}T_{\text{bc}i}T_{\text{ca}j})\}(y_{0})$\qquad

$-\frac{1}{12}\{\nabla_{k}\varrho_{ji}-2\varrho_{kj}%
<H,i>+\overset{q}{\underset{\text{a}=1}{\sum}}(\nabla_{k}R_{\text{a}%
j\text{a}i}-4R_{k\text{a}j\text{a}}<H,i>)$

$+4\overset{q}{\underset{\text{a,b=1}}{\sum}}R_{k\text{a}j\text{b}%
}T_{\text{ab}i}$+2$\overset{q}{\underset{\text{a,b,c=1}}{\sum}}(T_{\text{aa}%
k}T_{\text{bb}j}T_{\text{cc}i}-3T_{\text{aa}k}T_{\text{bc}j}T_{\text{bc}%
i}+2T_{\text{ab}k}T_{\text{bc}j}T_{\text{ca}i})\}(y_{0})$

\vspace{1pt}

\vspace{1pt}(xii) $\frac{\partial^{3}\theta}{\partial\text{x}_{i}^{2}%
\partial\text{x}_{j}}(y_{0})$

$=-\frac{1}{6}[\nabla_{i}\varrho_{ij}-2\varrho_{ij}%
<H,i>+\overset{q}{\underset{\text{a}=1}{\sum}}(\nabla_{i}R_{\text{a}%
i\text{a}j}-4R_{i\text{a}j\text{a}}<H,i>)+4\overset{q}{\underset{\text{a,b=1}%
}{\sum}}R_{i\text{a}j\text{b}}T_{\text{ab}i}$

$+2\overset{q}{\underset{\text{a,b,c=1}}{\sum}}(T_{\text{aa}i}T_{\text{bb}%
j}T_{\text{cc}i}-3T_{\text{aa}i}T_{\text{bc}j}T_{\text{bc}i}+2T_{\text{ab}%
i}T_{\text{bc}j}T_{\text{ac}i})](y_{0})$\qquad\qquad\qquad\qquad\qquad\ \ 

$-\frac{1}{6}[\nabla_{j}\varrho_{ii}-2\varrho_{ij}%
<H,i>+\overset{q}{\underset{\text{a}=1}{\sum}}(\nabla_{j}R_{\text{a}%
i\text{a}i}-4R_{i\text{a}j\text{a}}<H,i>)+4\overset{q}{\underset{\text{a,b=1}%
}{\sum}}R_{j\text{a}i\text{b}}T_{\text{ab}i}$

$+2\overset{q}{\underset{\text{a,b,c=1}}{\sum}}(T_{\text{aa}j}T_{\text{bb}%
i}T_{\text{cc}i}-3T_{\text{aa}j}T_{\text{bc}i}T_{\text{bc}i}+2T_{\text{ab}%
j}T_{\text{bc}i}T_{\text{ac}i})](y_{0})$

$-\frac{1}{6}[\nabla_{i}\varrho_{ij}-2\varrho_{ii}%
<H,j>+\overset{q}{\underset{\text{a}=1}{\sum}}(\nabla_{i}R_{\text{a}%
i\text{a}j}-4R_{i\text{a}i\text{a}}<H,j>)+4\overset{q}{\underset{\text{a,b=1}%
}{\sum}}R_{i\text{a}i\text{b}}T_{\text{ab}j}$

+2$\overset{q}{\underset{\text{a,b,c}=1}{\sum}}(T_{\text{aa}i}T_{\text{bb}%
i}T_{\text{cc}j}-3T_{\text{aa}i}T_{\text{bc}i}T_{\text{bc}j}+2T_{\text{ab}%
i}T_{\text{bc}i}T_{\text{ac}j})](y_{0})$\qquad\qquad\qquad\qquad$\ $\qquad

\vspace{1pt}(xiii) $\frac{\partial^{3}\theta}{\partial\text{x}_{i}%
\partial\text{x}_{j}^{2}}(y_{0})$

$=-\frac{1}{6}[\nabla_{i}\varrho_{jj}-2\varrho_{ij}%
<H,j>+\overset{q}{\underset{\text{a}=1}{\sum}}(\nabla_{i}R_{\text{a}%
j\text{a}j}-4R_{i\text{a}j\text{a}}<H,j>)+4\overset{q}{\underset{\text{a,b=1}%
}{\sum}}R_{i\text{a}j\text{b}}T_{\text{ab}j}$

$+2\overset{q}{\underset{\text{a,b,c=1}}{\sum}}(T_{\text{aa}i}T_{\text{bb}%
j}T_{\text{cc}j}-3T_{\text{aa}i}T_{\text{bc}j}T_{\text{bc}j}+2T_{\text{ab}%
i}T_{\text{bc}j}T_{\text{ca}j})](y_{0})$\qquad\qquad\qquad\qquad\qquad\ \ 

$-\frac{1}{6}[\nabla_{j}\varrho_{ij}-2\varrho_{ij}%
<H,j>+\overset{q}{\underset{\text{a}=1}{\sum}}(\nabla_{j}R_{\text{a}%
i\text{a}j}-4R_{j\text{a}i\text{a}}<H,j>)+4\overset{q}{\underset{\text{a,b=1}%
}{\sum}}R_{j\text{a}i\text{b}}T_{\text{ab}j}$

$+2\overset{q}{\underset{\text{a,b,c=1}}{\sum}}(T_{\text{aa}j}T_{\text{bb}%
i}T_{\text{cc}j}-3T_{\text{aa}j}T_{\text{bc}i}T_{\text{bc}j}+2T_{\text{ab}%
j}T_{\text{bc}i}T_{\text{ac}j})](y_{0})$

$-\frac{1}{6}[\nabla_{j}\varrho_{ij}-2\varrho_{jj}%
<H,i>+\overset{q}{\underset{\text{a}=1}{\sum}}(\nabla_{j}R_{\text{a}%
i\text{a}j}-4R_{j\text{a}j\text{a}}<H,i>)+4\overset{q}{\underset{\text{a,b=1}%
}{\sum}}R_{j\text{a}j\text{b}}T_{\text{ab}i}$

$+2\overset{q}{\underset{\text{a,b,c=1}}{\sum}}(T_{\text{aa}j}T_{\text{bb}%
j}T_{\text{cc}i}-3T_{\text{aa}j}T_{\text{bc}j}T_{\text{bc}i}+2T_{\text{ab}%
j}T_{\text{bc}j}T_{\text{ac}i})](y_{0})$\qquad\qquad\qquad\qquad$\ $\qquad

\vspace{1pt}(xiv) $\qquad\frac{\partial^{3}\theta^{\frac{1}{2}}}%
{\partial\text{x}_{i}\partial\text{x}_{j}\partial\text{x}_{k}}(y_{0})=\frac
{1}{8}(\frac{\partial\theta}{\partial\text{x}_{i}}\frac{\partial\theta
}{\partial\text{x}_{j}}\frac{\partial\theta}{\partial\text{x}_{k}}%
)(y_{0})-\frac{1}{4}(\frac{\partial\theta}{\partial\text{x}_{i}}\frac
{\partial^{2}\theta}{\partial\text{x}_{j}\partial\text{x}_{k}})(y_{0})$

$\qquad\qquad+\frac{1}{8}(\frac{\partial\theta}{\partial\text{x}_{j}}%
\frac{\partial\theta}{\partial\text{x}_{i}}\frac{\partial\theta}%
{\partial\text{x}_{k}})(y_{0})-\frac{1}{4}(\frac{\partial\theta}%
{\partial\text{x}_{j}}\frac{\partial^{2}\theta}{\partial\text{x}_{i}%
\partial\text{x}_{k}})(y_{0})$

\qquad\qquad\ $+\frac{1}{8}(\frac{\partial\theta}{\partial\text{x}_{k}}%
\frac{\partial\theta}{\partial\text{x}_{i}}\frac{\partial\theta}%
{\partial\text{x}_{j}})(y_{0})-\frac{1}{4}(\frac{\partial\theta}%
{\partial\text{x}_{k}}\frac{\partial^{2}\theta}{\partial\text{x}_{i}%
\partial\text{x}_{j}})(y_{0})+\frac{1}{2}\frac{\partial^{3}\theta}%
{\partial\text{x}_{i}\partial\text{x}_{j}\partial x_{k}}(y_{0})$

$\qquad=\frac{3}{8}(\frac{\partial\theta}{\partial\text{x}_{i}}\frac
{\partial\theta}{\partial\text{x}_{j}}\frac{\partial\theta}{\partial
\text{x}_{k}})(y_{0})-\frac{1}{4}(\frac{\partial\theta}{\partial\text{x}_{i}%
}\frac{\partial^{2}\theta}{\partial\text{x}_{j}\partial\text{x}_{k}}%
)(y_{0})-\frac{1}{4}(\frac{\partial\theta}{\partial\text{x}_{j}}\frac
{\partial^{2}\theta}{\partial\text{x}_{i}\partial\text{x}_{k}})(y_{0})$

$\qquad-\frac{1}{4}(\frac{\partial\theta}{\partial\text{x}_{k}}\frac
{\partial^{2}\theta}{\partial\text{x}_{i}\partial\text{x}_{j}})(y_{0}%
)+\frac{1}{2}\frac{\partial^{3}\theta}{\partial\text{x}_{i}\partial
\text{x}_{j}\partial x_{k}}(y_{0})$

We use the expressions already computed above.

\vspace{1pt}

\vspace{1pt}(xv) $\frac{\partial^{3}\theta^{-\frac{1}{2}}}{\partial
\text{x}_{i}\partial\text{x}_{j}\partial\text{x}_{k}}(y_{0})=-\frac{15}%
{8}(\frac{\partial\theta}{\partial x_{i}}\frac{\partial\theta}{\partial x_{j}%
}\frac{\partial\theta}{\partial x_{k}})(y_{0})+\frac{3}{4}\frac{\partial
^{2}\theta}{\partial\text{x}_{i}\partial\text{x}_{j}}(y_{0})\frac
{\partial\theta}{\partial x_{k}}(y_{0})$

\qquad\qquad$+\frac{3}{4}\frac{\partial\theta}{\partial x_{i}}(y_{0}%
)\frac{\partial^{2}\theta}{\partial\text{x}_{j}\partial\text{x}_{k}}(y_{0})$
$+\frac{3}{4}\frac{\partial\theta}{\partial x_{j}}(y_{0})\frac{\partial
^{2}\theta}{\partial\text{x}_{i}\partial\text{x}_{k}}(y_{0})-$ $\frac{1}%
{2}\frac{\partial^{3}\theta}{\partial\text{x}_{i}\partial\text{x}_{j}\partial
x_{k}}(y_{0})$

We use the expressions already computed above.\qquad

(xvi) $\frac{\partial^{3}\theta^{-\frac{1}{2}}}{\partial\text{x}_{i}%
^{2}\partial\text{x}_{j}}(y_{0})=-\frac{15}{8}(\frac{\partial\theta}{\partial
x_{i}})^{2}\frac{\partial\theta}{\partial x_{j}})(y_{0})+\frac{3}{2}%
\frac{\partial\theta}{\partial x_{i}}(y_{0})\frac{\partial^{2}\theta}%
{\partial\text{x}_{i}\partial\text{x}_{j}}(y_{0})$

\qquad\qquad\qquad\qquad\qquad\ $\ +\frac{3}{4}\frac{\partial\theta}{\partial
x_{j}}(y_{0})\frac{\partial^{2}\theta}{\partial\text{x}_{i}^{2}}(y_{0})-$
$\frac{1}{2}\frac{\partial^{3}\theta}{\partial\text{x}_{i}^{2}\partial
\text{x}_{j}}(y_{0})$\qquad

(xvii) $\frac{\partial^{3}\theta^{-\frac{1}{2}}}{\partial\text{x}_{i}%
\partial\text{x}_{j}^{2}}(y_{0})=-\frac{15}{8}(\frac{\partial\theta}{\partial
x_{i}}(\frac{\partial\theta}{\partial x_{j}})^{2})(y_{0})+\frac{3}{2}%
\frac{\partial\theta}{\partial x_{j}}(y_{0})\frac{\partial^{2}\theta}%
{\partial\text{x}_{i}\partial\text{x}_{j}}(y_{0})$

\qquad\qquad\qquad\qquad\qquad\ $+\frac{3}{4}\frac{\partial\theta}{\partial
x_{i}}(y_{0})\frac{\partial^{2}\theta}{\partial\text{x}_{j}^{2}}(y_{0})-$
$\frac{1}{2}\frac{\partial^{3}\theta}{\partial\text{x}_{i}\partial\text{x}%
_{j}^{2}}(y_{0})$

We use the expression already computed above

(xviii) We have for a = 1,...,q and $i=q+1,...,n:$

$[\frac{\partial^{2}}{\partial x_{i}^{2}}(\nabla\log\theta^{-\frac{1}{2}%
})_{\text{a}}](y_{0})=\frac{1}{2}\overset{n}{\underset{j=q+1}{\sum}%
}<H,j>(y_{0})[\frac{8}{3}R_{i\text{a}ij}-4\underset{\text{b=1}%
}{\overset{\text{q}}{\sum}}T_{\text{ab}i}\perp_{\text{b}ij}](y_{0})$

$\qquad\qquad\qquad\qquad\qquad+2\overset{n}{\underset{k=q+1}{\sum}}%
\perp_{\text{a}ij}(y_{0})[\frac{1}{4}<H,i><H,j>$

$+\frac{3}{4}<H,i>(y_{0})<H,j>(y_{0})$

$+\frac{1}{12}[2\varrho_{ij}+$ $\overset{q}{\underset{\text{a}=1}{4\sum}%
}R_{i\text{a}j\text{a}}-3\overset{q}{\underset{\text{a,b=1}}{\sum}%
}(T_{\text{aa}i}T_{\text{bb}j}+T_{\text{aa}j}T_{\text{bb}i}-2T_{\text{ab}%
i}T_{\text{ab}j})](y_{0})$

(xix) Then for $i,j=q+1,...,n:$

$\qquad\lbrack\frac{\partial^{2}}{\partial x_{i}^{2}}(\nabla\log\theta
^{-\frac{1}{2}})_{j}](y_{0})=\frac{1}{3}\underset{k=q+1}{\overset{n}{\sum}%
}<H,k>(y_{0})$R$_{ijik}(y_{0})$

$-\frac{1}{24}<H,j>(y_{0})[3<H,i>^{2}+2(\tau^{M}-3\tau^{P}%
+\overset{q}{\underset{\text{a}=1}{\sum}}\varrho_{\text{aa}}%
+\overset{q}{\underset{\text{a,b}=1}{\sum}}R_{\text{abab}})](y_{0})$

$-<H,i>(y_{0})[\frac{3}{4}<H,i><H,j>$\ 

$+\frac{1}{6}(\varrho_{ij}+2\overset{q}{\underset{\text{a}=1}{\sum}%
}R_{i\text{a}j\text{a}}-3\overset{q}{\underset{\text{a,b=1}}{\sum}%
}T_{\text{aa}i}T_{\text{bb}j}-T_{\text{ab}i}T_{\text{ab}j})](y_{0})$

$+\frac{15}{8}$%
$<$%
H,$i$%
$>$%
$^{2}$%
$<$%
H,$j$%
$>$%

$+\frac{1}{4}$%
$<$%
H,$i$%
$>$%
(y$_{0}$)$[$(2$\varrho_{ij}$+4$\overset{q}{\underset{\text{a}=1}{\sum}}%
$R$_{i\text{a}j\text{a}}$-3$\overset{q}{\underset{\text{a,b=1}}{\sum}}%
$T$_{\text{aa}i}$T$_{\text{bb}j}$-T$_{\text{ab}i}$T$_{\text{ab}j}%
$-3$\overset{q}{\underset{\text{a,b=1}}{\sum}}$T$_{\text{aa}j}$T$_{\text{bb}%
i}$-T$_{\text{ab}j}$T$_{\text{ab}i}$)$](y_{0})$

$+\frac{1}{4}<H,j>[$ $\tau^{M}-3\tau^{P}+\overset{q}{\underset{\text{a}%
=1}{\sum}}\varrho_{\text{aa}}+\overset{q}{\underset{\text{a,b}=1}{\sum}%
}R_{\text{abab}}](y_{0})$

$+\frac{1}{12}[\nabla_{i}\varrho_{ij}-2\varrho_{ij}%
<H,i>+\overset{q}{\underset{\text{a}=1}{\sum}}(\nabla_{i}R_{\text{a}%
i\text{a}j}-4R_{i\text{a}j\text{a}}<H,i>)$

$+4\overset{q}{\underset{\text{a,b=1}}{\sum}}R_{i\text{a}j\text{b}%
}T_{\text{ab}i}+2\overset{q}{\underset{\text{a,b,c=1}}{\sum}}(T_{\text{aa}%
i}T_{\text{bb}j}T_{\text{cc}i}-3T_{\text{aa}i}T_{\text{bc}j}T_{\text{bc}%
i}+2T_{\text{ab}i}T_{\text{bc}j}T_{\text{ca}i})](y_{0})$\qquad\qquad
\qquad\qquad\qquad\ \ 

$+\frac{1}{12}[\nabla_{j}\varrho_{ii}-2\varrho_{ji}%
<H,i>+\overset{q}{\underset{\text{a}=1}{\sum}}(\nabla_{j}R_{\text{a}%
i\text{a}i}-4R_{j\text{a}i\text{a}}<H,i>)$

$+4\overset{q}{\underset{\text{a,b=1}}{\sum}}R_{j\text{a}i\text{b}%
}T_{\text{ab}i}+2\overset{q}{\underset{\text{a,b,c=1}}{\sum}}(T_{\text{aa}%
j}T_{\text{bb}i}T_{\text{cc}i}-3T_{\text{aa}j}T_{\text{bc}i}T_{\text{bc}%
i}+2T_{\text{ab}j}T_{\text{bc}i}T_{\text{ca}i})](y_{0})$

$+\frac{1}{12}[\nabla_{i}\varrho_{ij}-2\varrho_{ii}%
<H,j>+\overset{q}{\underset{\text{a}=1}{\sum}}(\nabla_{i}R_{\text{a}%
i\text{a}j}-4R_{i\text{a}i\text{a}}<H,j>)$

$+4\overset{q}{\underset{\text{a,b=1}}{\sum}}R_{i\text{a}i\text{b}%
}T_{\text{ab}j}+2\overset{q}{\underset{\text{a,b,c=1}}{\sum}}(T_{\text{aa}%
i}T_{\text{bb}i}T_{\text{cc}j}-3T_{\text{aa}i}T_{\text{bc}i}T_{\text{bc}%
j}+2T_{\text{ab}i}T_{\text{bc}i}T_{\text{ac}j})](y_{0})$

\qquad\qquad\qquad\qquad\qquad\qquad\qquad\qquad\qquad\qquad$\blacksquare$

$\left(  A_{21}\right)  $\qquad$\frac{\partial^{4}\theta_{p}}{\partial
\text{x}_{i}^{2}\partial\text{x}_{j}^{2}}(y_{0})$

$=4\times\ \frac{1}{24}[$ $\overset{q}{\underset{\text{a=1}}{\sum}}%
\{-(\nabla_{ii}^{2}R_{j\text{a}j\text{a}}+\nabla_{jj}^{2}R_{i\text{a}%
i\text{a}}+\nabla_{ij}^{2}R_{i\text{a}j\text{a}}+\nabla_{ij}^{2}%
R_{j\text{a}i\text{a}}+\nabla_{ji}^{2}R_{i\text{a}j\text{a}}+\nabla_{ji}%
^{2}R_{j\text{a}i\text{a}})$ $A$

$+\overset{n}{\underset{p=q+1}{\sum}}\overset{q}{\underset{\text{a=1}}{\sum}%
}(R_{\text{a}iip}R_{\text{a}jjp}+R_{\text{a}jjp}R_{\text{a}iip}+R_{\text{a}%
ijp}R_{\text{a}ijp}+R_{\text{a}ijp}R_{\text{a}jip}+R_{\text{a}jip}%
R_{\text{a}ijp}+R_{\text{a}jip}R_{\text{a}jip})$

$+2\overset{q}{\underset{\text{a,b=1}}{\sum}}\nabla_{i}(R)_{\text{a}%
i\text{b}j}T_{\text{ab}j}+2\overset{q}{\underset{\text{a,b=1}}{\sum}}%
\nabla_{j}(R)_{\text{a}j\text{b}i}T_{\text{ab}i}%
+2\overset{q}{\underset{\text{a,b=1}}{\sum}}\nabla_{i}(R)_{\text{a}j\text{b}%
i}T_{\text{ab}j}+2\overset{q}{\underset{\text{a,b=1}}{\sum}}\nabla
_{i}(R)_{\text{a}j\text{b}j}T_{\text{ab}i}$

$+2\overset{q}{\underset{\text{a,b=1}}{\sum}}\nabla_{j}(R)_{\text{a}%
i\text{b}i}T_{\text{ab}j}+2\overset{q}{\underset{\text{a,b=1}}{\sum}}%
\nabla_{j}(R)_{\text{a}i\text{b}j}T_{\text{ab}i}%
+\overset{n}{\underset{p=q+1}{\sum}}(-\frac{3}{5}\nabla_{ii}^{2}%
(R)_{jpjp}+\overset{n}{\underset{p=q+1}{\sum}}(-\frac{3}{5}\nabla_{jj}%
^{2}(R)_{ipip}$

$+\overset{n}{\underset{p=q+1}{\sum}}(-\frac{3}{5}\nabla_{ij}^{2}%
(R)_{ipjp}+\overset{n}{\underset{p=q+1}{\sum}}(-\frac{3}{5}\nabla_{ij}%
^{2}(R)_{jpip}+\overset{n}{\underset{p=q+1}{\sum}}(-\frac{3}{5}\nabla_{ji}%
^{2}(R)_{ipjp}+\overset{n}{\underset{p=q+1}{\sum}}(-\frac{3}{5}\nabla_{ji}%
^{2}(R)_{jpip}$

$+\frac{1}{5}\overset{n}{\underset{m,p=q+1}{%
{\textstyle\sum}
}}R_{ipim}R_{jpjm}+\frac{1}{5}\overset{n}{\underset{m,p=q+1}{%
{\textstyle\sum}
}}R_{jpjm}R_{ipim}+\frac{1}{5}\overset{n}{\underset{m,p=q+1}{%
{\textstyle\sum}
}}R_{ipjm}R_{ipjm}+\frac{1}{5}\overset{n}{\underset{m,p=q+1}{%
{\textstyle\sum}
}}R_{ipjm}R_{jpim}$

$+\frac{1}{5}\overset{n}{\underset{m,p=q+1}{%
{\textstyle\sum}
}}R_{jpim}R_{ipjm}+\frac{1}{5}\overset{n}{\underset{m,p=q+1}{%
{\textstyle\sum}
}}R_{jpim}R_{jpim}\}(y_{0})$

$+4\overset{q}{\underset{\text{a,b=1}}{\sum}}\{(\nabla_{i}(R)_{i\text{a}%
j\text{a}}-\overset{q}{\underset{\text{c=1}}{%
{\textstyle\sum}
}}R_{\text{a}i\text{c}i}T_{\text{ac}j})$ $T_{\text{bb}j}+4(\nabla
_{j}(R)_{j\text{a}i\text{a}}-\overset{q}{\underset{\text{c=1}}{%
{\textstyle\sum}
}}R_{\text{a}j\text{c}j}T_{\text{ac}i})$ $T_{\text{bb}i}$

$+4(\nabla_{i}(R)_{j\text{a}i\text{a}}-\overset{q}{\underset{\text{c=1}}{%
{\textstyle\sum}
}}R_{\text{a}i\text{c}j}T_{\text{ac}i})$ $T_{\text{bb}j}$ $4B\ $

$+4(\nabla_{i}(R)_{j\text{a}j\text{a}}-\overset{q}{\underset{\text{c=1}}{%
{\textstyle\sum}
}}R_{\text{a}i\text{c}j}T_{\text{ac}j})$ $T_{\text{bb}i}+4(\nabla
_{j}(R)_{i\text{a}i\text{a}}-\overset{q}{\underset{\text{c=1}}{%
{\textstyle\sum}
}}R_{\text{a}j\text{c}i}T_{\text{ac}i})$ $T_{\text{bb}j}+4(\nabla
_{j}(R)_{i\text{a}j\text{a}}-\overset{q}{\underset{\text{c=1}}{%
{\textstyle\sum}
}}R_{\text{a}j\text{c}i}T_{\text{ac}j})$ $T_{\text{bb}i}$

$-4\overset{q}{\underset{\text{a,b=1}}{\sum}}(\nabla_{i}(R)_{i\text{a}%
j\text{b}}-\overset{q}{\underset{\text{c=1}}{%
{\textstyle\sum}
}}R_{\text{b}r\text{c}s}T_{\text{ac}t})T_{\text{ab}j}%
-4\overset{q}{\underset{\text{a,b=1}}{\sum}}(\nabla_{j}(R)_{j\text{a}%
i\text{b}}-\overset{q}{\underset{\text{c=1}}{%
{\textstyle\sum}
}}R_{\text{b}j\text{c}j}T_{\text{ac}i})T_{\text{ab}i}$

$-4\overset{q}{\underset{\text{a,b=1}}{\sum}}(\nabla_{i}(R)_{j\text{a}%
i\text{b}}-\overset{q}{\underset{\text{c=1}}{%
{\textstyle\sum}
}}R_{\text{b}i\text{c}j}T_{\text{ac}i})T_{\text{ab}j}%
-4\overset{q}{\underset{\text{a,b=1}}{\sum}}(\nabla_{i}(R)_{j\text{a}%
j\text{b}}-\overset{q}{\underset{\text{c=1}}{%
{\textstyle\sum}
}}R_{\text{b}i\text{c}j}T_{\text{ac}j})T_{\text{ab}i}$

$\bigskip-4\overset{q}{\underset{\text{a,b=1}}{\sum}}(\nabla_{j}%
(R)_{i\text{a}i\text{b}}-\overset{q}{\underset{\text{c=1}}{%
{\textstyle\sum}
}}R_{\text{b}j\text{c}i}T_{\text{ac}i})T_{\text{ab}j}%
-4\overset{q}{\underset{\text{a,b=1}}{\sum}}(\nabla_{j}(R)_{i\text{a}%
j\text{b}}-\overset{q}{\underset{\text{c=1}}{%
{\textstyle\sum}
}}R_{\text{b}j\text{c}i}T_{\text{ac}j})T_{\text{ab}i}\}\qquad$

$+\frac{4}{3}\overset{q}{\underset{\text{a,b}=1}{\sum}}(R_{i\text{a}i\text{a}%
}R_{j\text{b}j\text{b}})+\frac{4}{3}\overset{q}{\underset{\text{a,b}=1}{\sum}%
}(R_{j\text{a}j\text{a}}R_{i\text{b}i\text{b}})+\frac{4}{3}%
\overset{q}{\underset{\text{a,b}=1}{\sum}}(R_{i\text{a}j\text{a}}%
R_{i\text{b}j\text{b}})\qquad\qquad3C$

$+\frac{4}{3}\overset{q}{\underset{\text{a,b}=1}{\sum}}(R_{i\text{a}j\text{a}%
}R_{j\text{b}i\text{b}})+\frac{4}{3}\overset{q}{\underset{\text{a,b}=1}{\sum}%
}(R_{j\text{a}i\text{a}}R_{i\text{b}j\text{b}})+\frac{4}{3}%
\overset{q}{\underset{\text{a,b}=1}{\sum}}(R_{j\text{a}i\text{a}}%
R_{j\text{b}i\text{b}})$

$+\frac{1}{3}\varrho_{ii}\varrho_{jj}+\frac{1}{3}\varrho_{jj}\varrho
_{ii}+\frac{1}{3}\varrho_{ij}\varrho_{ij}+\frac{1}{3}\varrho_{ij}\varrho
_{ji}+\frac{1}{3}\varrho_{ji}\varrho_{ij}+\frac{1}{3}\varrho_{ji}\varrho_{ji}$

$+\frac{2}{3}\overset{q}{\underset{\text{a}=1}{\sum}}R_{i\text{a}i\text{a}%
}\varrho_{jj}\ +\frac{2}{3}\overset{q}{\underset{\text{a}=1}{\sum}%
}R_{j\text{a}j\text{a}}\varrho_{ii}\ +\frac{2}{3}%
\overset{q}{\underset{\text{a}=1}{\sum}}R_{i\text{a}j\text{a}}\varrho
_{ij}\ +\frac{2}{3}\overset{q}{\underset{\text{a}=1}{\sum}}R_{i\text{a}%
j\text{a}}\varrho_{ji}\ $

$+\frac{2}{3}\overset{q}{\underset{\text{a}=1}{\sum}}R_{j\text{a}i\text{a}%
}\varrho_{ij}\ +\frac{2}{3}\overset{q}{\underset{\text{a}=1}{\sum}%
}R_{j\text{a}i\text{a}}\varrho_{ji}\ +\frac{2}{3}%
\overset{q}{\underset{\text{b}=1}{\sum}}R_{i\text{b}i\text{b}}\varrho
_{jj}\ +\frac{2}{3}\overset{q}{\underset{\text{b}=1}{\sum}}R_{j\text{b}%
j\text{b}}\varrho_{ii}\ $

$+\frac{2}{3}\overset{q}{\underset{\text{b}=1}{\sum}}R_{i\text{b}j\text{b}%
}\varrho_{ij}\ +\frac{2}{3}\overset{q}{\underset{\text{b}=1}{\sum}%
}R_{i\text{b}j\text{b}}\varrho_{ji}\ +\frac{2}{3}%
\overset{q}{\underset{\text{b}=1}{\sum}}R_{j\text{b}i\text{b}}\varrho
_{ij}\ +\frac{2}{3}\overset{q}{\underset{\text{b}=1}{\sum}}R_{j\text{b}%
i\text{b}}\varrho_{ji}$

$-3\overset{q}{\underset{\text{a,b}=1}{\sum}}R_{i\text{a}i\text{b}%
}R_{j\text{a}j\text{b}}\ -3\overset{q}{\underset{\text{a,b}=1}{\sum}%
}R_{j\text{a}j\text{b}}R_{i\text{a}i\text{b}}%
\ -3\overset{q}{\underset{\text{a,b}=1}{\sum}}R_{i\text{a}j\text{b}%
}R_{i\text{a}j\text{b}}\ -3\overset{q}{\underset{\text{a,b}=1}{\sum}%
}R_{i\text{a}j\text{b}}R_{j\text{a}i\text{b}}$

$-3\overset{q}{\underset{\text{a,b}=1}{\sum}}R_{j\text{a}i\text{b}%
}R_{i\text{a}j\text{b}}\ -3\overset{q}{\underset{\text{a,b}=1}{\sum}%
}R_{j\text{a}i\text{b}}R_{j\text{a}i\text{b}}\ $

$-\frac{1}{3}\overset{n}{\underset{p,m=q+1}{\sum}}R_{ipim}R_{jpjm}\ -\frac
{1}{3}\overset{n}{\underset{p,m=q+1}{\sum}}R_{jpjm}R_{ipim}-\frac{1}%
{3}\overset{n}{\underset{p,m=q+1}{\sum}}R_{ipjm}R_{ipjm}$

$-\frac{1}{3}\overset{n}{\underset{p,m=q+1}{\sum}}R_{ipjm}R_{jpim}-\frac{1}%
{3}\overset{n}{\underset{p,m=q+1}{\sum}}R_{jpim}R_{ipjm}-\frac{1}%
{3}\overset{n}{\underset{p,m=q+1}{\sum}}R_{jpim}R_{jpim}$

$-\overset{q}{\underset{\text{a}=1}{\sum}}\overset{n}{\underset{p=q+1}{\sum}%
}R_{i\text{a}ip}R_{j\text{a}jp}-\overset{q}{\underset{\text{a}=1}{\sum}%
}\overset{n}{\underset{p=q+1}{\sum}}R_{j\text{a}jp}R_{i\text{a}ip}%
-\overset{q}{\underset{\text{a}=1}{\sum}}\overset{n}{\underset{p=q+1}{\sum}%
}R_{i\text{a}jp}R_{i\text{a}jp}$

$-\overset{q}{\underset{\text{a}=1}{\sum}}\overset{n}{\underset{p=q+1}{\sum}%
}R_{i\text{a}jp}R_{j\text{a}ip}-\overset{q}{\underset{\text{a}=1}{\sum}%
}\overset{n}{\underset{p=q+1}{\sum}}R_{j\text{a}ip}R_{i\text{a}jp}%
-\overset{q}{\underset{\text{a}=1}{\sum}}\overset{n}{\underset{p=q+1}{\sum}%
}R_{j\text{a}ip}R_{j\text{a}ip}$

$-\overset{q}{\underset{\text{b}=1}{\sum}}\overset{n}{\underset{p=q+1}{\sum}%
}R_{i\text{b}ip}R_{j\text{b}jp}-\overset{q}{\underset{\text{b}=1}{\sum}%
}\overset{n}{\underset{p=q+1}{\sum}}R_{j\text{b}jp}R_{i\text{b}ip}%
-\overset{q}{\underset{\text{b}=1}{\sum}}\overset{n}{\underset{p=q+1}{\sum}%
}R_{i\text{b}jp}R_{i\text{b}jp}$

$-\overset{q}{\underset{\text{b}=1}{\sum}}\overset{n}{\underset{p=q+1}{\sum}%
}R_{i\text{b}jp}R_{j\text{b}ip}-\overset{q}{\underset{\text{b}=1}{\sum}%
}\overset{n}{\underset{p=q+1}{\sum}}R_{j\text{b}ip}R_{i\text{b}jp}%
-\overset{q}{\underset{\text{b}=1}{\sum}}\overset{n}{\underset{p=q+1}{\sum}%
}R_{j\text{b}ip}R_{j\text{b}ip}$

$+$ $\overset{q}{\underset{\text{a,b,c=1}}{6\sum}}\{$ $-R_{i\text{a}i\text{a}%
}(T_{\text{bb}j}T_{\text{cc}j}$ $-T_{\text{bc}j}T_{\text{bc}j})\}+$
$\overset{q}{\underset{\text{a,b,c=1}}{6\sum}}\{$ $-R_{j\text{a}j\text{a}%
}(T_{\text{bb}i}T_{\text{cc}i}$ $-T_{\text{bc}i}T_{\text{bc}i})\}\qquad6D$

$+$ $\overset{q}{\underset{\text{a,b,c=1}}{6\sum}}\{$ $-R_{i\text{a}j\text{a}%
}(T_{\text{bb}i}T_{\text{cc}j}$ $-T_{\text{bc}i}T_{\text{bc}j})\}+$
$\overset{q}{\underset{\text{a,b,c=1}}{6\sum}}\{$ $-R_{i\text{a}j\text{a}%
}(T_{\text{bb}j}T_{\text{cc}i}$ $-T_{\text{bc}j}T_{\text{bc}i})\}$

$+$ $\overset{q}{\underset{\text{a,b,c=1}}{6\sum}}\{$ $-R_{j\text{a}i\text{a}%
}(T_{\text{bb}i}T_{\text{cc}j}$ $-T_{\text{bc}i}T_{\text{bc}j})\}+$
$\overset{q}{\underset{\text{a,b,c=1}}{6\sum}}\{$ $-R_{j\text{a}i\text{a}%
}(T_{\text{bb}j}T_{\text{cc}i}$ $-T_{\text{bc}j}T_{\text{bc}i})\}$

\ $+6\{$\ $R_{i\text{a}i\text{b}}(T_{\text{ab}j}T_{\text{cc}j}-T_{\text{bc}%
j}T_{\text{ac}j})\}\ +6\{\ R_{j\text{a}j\text{b}}(T_{\text{ab}i}T_{\text{cc}%
i}-T_{\text{bc}i}T_{\text{ac}i})\}$

$+6\{\ R_{i\text{a}j\text{b}}(T_{\text{ab}i}T_{\text{cc}j}-T_{\text{bc}%
i}T_{\text{ac}j})\}\ +6\{\ R_{i\text{a}j\text{b}}(T_{\text{ab}j}T_{\text{cc}%
i}-T_{\text{bc}j}T_{\text{ac}i})\}$

$+6\{\ R_{j\text{a}i\text{ib}}(T_{\text{ab}i}T_{\text{cc}j}-T_{\text{bc}%
i}T_{\text{ac}j})\}\ +6\{\ R_{j\text{a}i\text{b}}(T_{\text{ab}j}T_{\text{cc}%
i}-T_{\text{bc}j}T_{\text{ac}i})\}\qquad$

$+6\{-R_{i\text{a}i\text{c}}(T_{\text{ba}j}T_{\text{bc}j}-T_{\text{ac}%
j}T_{\text{bb}j})\}$\ $+6\{-R_{j\text{a}j\text{c}}(T_{\text{ba}i}%
T_{\text{bc}i}-T_{\text{ac}i}T_{\text{bb}i})\}$

$+6\{-R_{i\text{a}j\text{c}}(T_{\text{ba}i}T_{\text{bc}j}-T_{\text{ac}%
i}T_{\text{bb}j})\}+6\{-R_{i\text{a}j\text{c}}(T_{\text{ba}j}T_{\text{bc}%
i}-T_{\text{ac}j}T_{\text{bb}i})\}$

$+6\{-R_{j\text{a}i\text{c}}(T_{\text{ba}i}T_{\text{bc}j}-T_{\text{ac}%
i}T_{\text{bb}j})\}+6\{-R_{j\text{a}i\text{c}}(T_{\text{ba}j}T_{\text{bc}%
i}-T_{\text{ac}j}T_{\text{bb}i})\}$

$+6\overset{q}{\underset{\text{b,c=1}}{\sum}}\underset{p=q+1}{\overset{n}{\sum
}}$\ $\{-\frac{1}{3}R_{ipip}(T_{\text{bb}j}T_{\text{cc}j}-T_{\text{bc}%
j}T_{\text{bc}j})\}+6\overset{q}{\underset{\text{b,c=1}}{\sum}}%
\underset{p=q+1}{\overset{n}{\sum}}\ \{-\frac{1}{3}R_{jpjp}(T_{\text{bb}%
i}T_{\text{cc}i}-T_{\text{bc}i}T_{\text{bc}i})\}$

$+6\overset{q}{\underset{\text{b,c=1}}{\sum}}\underset{p=q+1}{\overset{n}{\sum
}}$\ $\{-\frac{1}{3}R_{ipjp}(T_{\text{bb}i}T_{\text{cc}j}-T_{\text{bc}%
i}T_{\text{bc}j})\}+6\overset{q}{\underset{\text{b,c=1}}{\sum}}%
\underset{p=q+1}{\overset{n}{\sum}}\ \{-\frac{1}{3}R_{ipjp}(T_{\text{bb}%
j}T_{\text{cc}i}-T_{\text{bc}j}T_{\text{bc}i})\}$

$+6\overset{q}{\underset{\text{b,c=1}}{\sum}}\underset{p=q+1}{\overset{n}{\sum
}}$\ $\{-\frac{1}{3}R_{jpip}(T_{\text{bb}i}T_{\text{cc}j}-T_{\text{bc}%
i}T_{\text{bc}j})\}+6\overset{q}{\underset{\text{b,c=1}}{\sum}}%
\underset{p=q+1}{\overset{n}{\sum}}\ \{-\frac{1}{3}R_{jpip}(T_{\text{bb}%
j}T_{\text{cc}i}-T_{\text{bc}j}T_{\text{bc}i})\}\qquad$

$+\overset{q}{\underset{\text{a,b,c,d}=1}{\sum}}T_{\text{aa}i}T_{\text{bb}%
i}(T_{\text{cc}j}T_{\text{dd}j}-T_{\text{cd}j}T_{\text{dc}j})+T_{\text{aa}%
j}T_{\text{bb}j}(T_{\text{cc}i}T_{\text{dd}i}-T_{\text{cd}i}T_{\text{dc}%
i})\qquad E$

$+T_{\text{aa}i}T_{\text{bb}j}(T_{\text{cc}i}T_{\text{dd}j}-T_{\text{cd}%
i}T_{\text{dc}j})+T_{\text{aa}i}T_{\text{bb}j}(T_{\text{cc}j}T_{\text{dd}%
i}-T_{\text{cd}j}T_{\text{dc}i})$

$+T_{\text{aa}j}T_{\text{bb}i}(T_{\text{cc}i}T_{\text{dd}j}-T_{\text{cd}%
i}T_{\text{dc}j})+T_{\text{aa}j}T_{\text{bb}i}(T_{\text{cc}j}T_{\text{dd}%
i}-T_{\text{cd}j}T_{\text{dc}i})$

$-T_{\text{aa}i}T_{\text{bc}i}(T_{\text{bc}j}T_{\text{dd}j}-T_{\text{bd}%
j}T_{\text{cd}j})-T_{\text{aa}j}T_{\text{bc}j}(T_{\text{bc}i}T_{\text{dd}%
i}-T_{\text{bd}i}T_{\text{cd}i})-T_{\text{aa}i}T_{\text{bc}j}(T_{\text{bc}%
i}T_{\text{dd}j}-T_{\text{bd}i}T_{\text{cd}j})$

$-T_{\text{aa}i}T_{\text{bc}j}(T_{\text{bc}j}T_{\text{dd}i}-T_{\text{bd}%
j}T_{\text{cd}i})-T_{\text{aa}j}T_{\text{bc}i}(T_{\text{bc}i}T_{\text{dd}%
j}-T_{\text{bd}i}T_{\text{cd}j})-T_{\text{aa}j}T_{\text{bc}i}(T_{\text{bc}%
j}T_{\text{dd}i}-T_{\text{bd}j}T_{\text{cd}i})$

$+T_{\text{aa}i}T_{\text{bd}i}(T_{\text{bc}j}T_{\text{cd}j}-T_{\text{bd}%
j}T_{\text{cc}j})+T_{\text{aa}j}T_{\text{bd}j}(T_{\text{bc}i}T_{\text{cd}%
i}-T_{\text{bd}i}T_{\text{cc}i})+T_{\text{aa}i}T_{\text{bd}j}(T_{\text{bc}%
i}T_{\text{cd}j}-T_{\text{bd}i}T_{\text{cc}j})$

$+T_{\text{aa}i}T_{\text{bd}j}(T_{\text{bc}j}T_{\text{cd}i}-T_{\text{bd}%
j}T_{\text{cc}i})+T_{\text{aa}j}T_{\text{bd}i}(T_{\text{bc}i}T_{\text{cd}%
j}-T_{\text{bd}i}T_{\text{cc}j})+T_{\text{aa}j}T_{\text{bd}i}(T_{\text{bc}%
j}T_{\text{cd}i}-T_{\text{bd}j}T_{\text{cc}i})\qquad$

$-T_{\text{ab}i}T_{\text{ab}i}(T_{\text{cc}j}T_{\text{dd}j}-T_{\text{cd}%
j}T_{\text{dc}j})-T_{\text{ab}j}T_{\text{ab}j}(T_{\text{cc}i}T_{\text{dd}%
i}-T_{\text{cd}i}T_{\text{dc}i})-T_{\text{ab}i}T_{\text{ab}j}(T_{\text{cc}%
i}T_{\text{dd}j}-T_{\text{cd}i}T_{\text{dc}j})$

$-T_{\text{ab}i}T_{\text{ab}j}(T_{\text{cc}j}T_{\text{dd}i}-T_{\text{cd}%
j}T_{\text{dc}i})-T_{\text{ab}j}T_{\text{ab}i}(T_{\text{cc}i}T_{\text{dd}%
j}-T_{\text{cd}i}T_{\text{dc}j})-T_{\text{ab}j}T_{\text{ab}i}(T_{\text{cc}%
j}T_{\text{dd}i}-T_{\text{cd}j}T_{\text{dc}i})$

$+T_{\text{ab}i}T_{\text{bc}i}(T_{\text{ac}j}T_{\text{dd}j}-T_{\text{ad}%
j}T_{\text{cd}j})+T_{\text{ab}j}T_{\text{bc}j}(T_{\text{ac}i}T_{\text{dd}%
i}-T_{\text{ad}i}T_{\text{cd}i})+T_{\text{ab}i}T_{\text{bc}j}(T_{\text{ac}%
i}T_{\text{dd}j}-T_{\text{ad}i}T_{\text{cd}j})$

$+T_{\text{ab}i}T_{\text{bc}j}(T_{\text{ac}j}T_{\text{dd}i}-T_{\text{ad}%
j}T_{\text{cd}i})+T_{\text{ab}j}T_{\text{bc}i}(T_{\text{ac}i}T_{\text{dd}%
j}-T_{\text{ad}i}T_{\text{cd}j})+T_{\text{ab}j}T_{\text{bc}i}(T_{\text{ac}%
j}T_{\text{dd}i}-T_{\text{ad}j}T_{\text{cd}i})$

$-T_{\text{ab}i}T_{\text{bd}i}(T_{\text{ac}j}T_{\text{cd}j}-T_{\text{ad}%
j}T_{\text{cc}j})-T_{\text{ab}j}T_{\text{bd}j}(T_{\text{ac}i}T_{\text{cd}%
i}-T_{\text{ad}i}T_{\text{cc}i})-T_{\text{ab}i}T_{\text{bd}j}(T_{\text{ac}%
i}T_{\text{cd}j}-T_{\text{ad}i}T_{\text{cc}j})$

$-T_{\text{ab}i}T_{\text{bd}j}(T_{\text{ac}j}T_{\text{cd}i}-T_{\text{ad}%
j}T_{\text{cc}i})-T_{\text{ab}i}T_{\text{bd}j}(T_{\text{ac}j}T_{\text{cd}%
i}-T_{\text{ad}j}T_{\text{cc}i})-T_{\text{ab}j}T_{\text{bd}i}(T_{\text{ac}%
j}T_{\text{cd}i}-T_{\text{ad}j}T_{\text{cc}i})$

\vspace{1pt}$+T_{\text{ac}i}T_{\text{ab}i}(T_{\text{bc}j}T_{\text{dd}%
j}-T_{\text{bd}j}T_{\text{dc}j})+T_{\text{ac}j}T_{\text{ab}j}(T_{\text{bc}%
i}T_{\text{dd}i}-T_{\text{bd}i}T_{\text{dc}i})+T_{\text{ac}i}T_{\text{ab}%
j}(T_{\text{bc}i}T_{\text{dd}j}-T_{\text{bd}i}T_{\text{dc}j})$

$+T_{\text{ac}i}T_{\text{ab}j}(T_{\text{bc}j}T_{\text{dd}i}-T_{\text{bd}%
j}T_{\text{dc}i})+T_{\text{ac}j}T_{\text{ab}i}(T_{\text{bc}i}T_{\text{dd}%
j}-T_{\text{bd}i}T_{\text{dc}j})+T_{\text{ac}j}T_{\text{ab}i}(T_{\text{bc}%
j}T_{\text{dd}i}-T_{\text{bd}j}T_{\text{dc}i})$

$-T_{\text{ac}i}T_{\text{bb}i}(T_{\text{ac}j}T_{\text{dd}j}-T_{\text{ad}%
j}T_{\text{cd}j})-T_{\text{ac}j}T_{\text{bb}j}(T_{\text{ac}i}T_{\text{dd}%
i}-T_{\text{ad}i}T_{\text{cd}i})-T_{\text{ac}i}T_{\text{bb}j}(T_{\text{ac}%
i}T_{\text{dd}j}-T_{\text{ad}i}T_{\text{cd}j})$

$-T_{\text{ac}i}T_{\text{bb}j}(T_{\text{ac}j}T_{\text{dd}i}-T_{\text{ad}%
j}T_{\text{cd}i})-T_{\text{ac}j}T_{\text{bb}i}(T_{\text{ac}i}T_{\text{dd}%
j}-T_{\text{ad}i}T_{\text{cd}i})-T_{\text{ac}j}T_{\text{bb}i}(T_{\text{ac}%
j}T_{\text{dd}i}-T_{\text{ad}j}T_{\text{cd}i})$

$+T_{\text{ac}i}T_{\text{bd}i}(T_{\text{ac}j}T_{\text{bd}j}-T_{\text{ad}%
j}T_{\text{bc}j})+T_{\text{ac}j}T_{\text{bd}j}(T_{\text{ac}i}T_{\text{bd}%
i}-T_{\text{ad}i}T_{\text{bc}i})+T_{\text{ac}i}T_{\text{bd}j}(T_{\text{ac}%
i}T_{\text{bd}j}-T_{\text{ad}i}T_{\text{bc}j})$

$+T_{\text{ac}i}T_{\text{bd}j}(T_{\text{ac}j}T_{\text{bd}i}-T_{\text{ad}%
j}T_{\text{bc}i})+T_{\text{ac}j}T_{\text{bd}i}(T_{\text{ac}i}T_{\text{bd}%
j}-T_{\text{ad}i}T_{\text{bc}j})+T_{\text{ac}j}T_{\text{bd}i}(T_{\text{ac}%
j}T_{\text{bd}i}-T_{\text{ad}j}T_{\text{bc}i})$

\vspace{1pt}$-T_{\text{ad}i}T_{\text{ab}i}(T_{\text{bc}j}T_{\text{cd}%
j}-T_{\text{bd}j}T_{\text{cc}j})-T_{\text{ad}j}T_{\text{ab}j}(T_{\text{bc}%
i}T_{\text{cd}i}-T_{\text{bd}i}T_{\text{cc}i})-T_{\text{ad}i}T_{\text{ab}%
j}(T_{\text{bc}i}T_{\text{cd}j}-T_{\text{bd}i}T_{\text{cc}j})$

$-T_{\text{ad}i}T_{\text{ab}j}(T_{\text{bc}j}T_{\text{cd}i}-T_{\text{bd}%
j}T_{\text{cc}i})-T_{\text{ad}j}T_{\text{ab}i}(T_{\text{bc}i}T_{\text{cd}%
j}-T_{\text{bd}i}T_{\text{cc}j})-T_{\text{ad}j}T_{\text{ab}i}(T_{\text{bc}%
j}T_{\text{cd}i}-T_{\text{bd}j}T_{\text{cc}i})$

$+T_{\text{ad}i}T_{\text{bb}i}(T_{\text{ac}j}T_{\text{cd}j}-T_{\text{ad}%
j}T_{\text{cc}j})+T_{\text{ad}j}T_{\text{bb}j}(T_{\text{ac}i}T_{\text{cd}%
i}-T_{\text{ad}i}T_{\text{cc}i})+T_{\text{ad}i}T_{\text{bb}j}(T_{\text{ac}%
i}T_{\text{cd}j}-T_{\text{ad}i}T_{\text{cc}j})$

$+T_{\text{ad}i}T_{\text{bb}j}(T_{\text{ac}j}T_{\text{cd}i}-T_{\text{ad}%
j}T_{\text{cc}i})+T_{\text{ad}j}T_{\text{bb}i}(T_{\text{ac}i}T_{\text{cd}%
j}-T_{\text{ad}i}T_{\text{cc}j})+T_{\text{ad}j}T_{\text{bb}i}(T_{\text{ac}%
j}T_{\text{cd}i}-T_{\text{ad}j}T_{\text{cc}i})$

$-T_{\text{ad}i}T_{\text{bc}i}(T_{\text{ac}j}T_{\text{bd}j}-T_{\text{ad}%
j}T_{\text{bc}j})-T_{\text{ad}j}T_{\text{bc}j}(T_{\text{ac}i}T_{\text{bd}%
i}-T_{\text{ad}i}T_{\text{bc}i})-T_{\text{ad}i}T_{\text{bc}j}(T_{\text{ac}%
i}T_{\text{bd}j}-T_{\text{ad}i}T_{\text{bc}j})$

$-T_{\text{ad}i}T_{\text{bc}j}(T_{\text{ac}j}T_{\text{bd}i}-T_{\text{ad}%
j}T_{\text{bc}i})-T_{\text{ad}j}T_{\text{bc}i}(T_{\text{ac}i}T_{\text{bd}%
j}-T_{\text{ad}i}T_{\text{bc}j})-T_{\text{ad}j}T_{\text{bc}i}(T_{\text{ac}%
j}T_{\text{bd}i}-T_{\text{ad}j}T_{\text{bc}i})](y_{0})$

\begin{center}
\qquad\qquad\qquad\qquad\qquad\qquad\qquad$\blacksquare$
\end{center}

Some items of 3C can be can be more elegantly expressed geometrically. We
recall the properties:

$R_{ijkl}=R_{klij}$ and $\varrho_{ij}=\overset{n}{\underset{k=1}{\sum}%
}R_{ikjk}$ which is symmetric in the indices $\left(  i,j\right)  .$The first is:

$+\frac{4}{3}\overset{n}{\underset{i,j=q+1}{\sum}}(R_{i\text{a}i\text{a}%
}R_{j\text{b}j\text{b}})+\frac{4}{3}\overset{q}{\underset{\text{a,b}=1}{\sum}%
}(R_{j\text{a}j\text{a}}R_{i\text{b}i\text{b}})$

$=\frac{4}{3}(\overset{n}{\underset{i=1}{\sum}}R_{i\text{a}i\text{a}%
}-\overset{q}{\underset{\text{c}=1}{\sum}}R_{\text{caca}}%
)(\overset{n}{\underset{j=1}{\sum}}R_{j\text{b}j\text{b}}%
-\overset{q}{\underset{\text{d}=1}{\sum}}R_{\text{dbdb}})+\frac{4}%
{3}(\overset{n}{\underset{j=1}{\sum}}R_{j\text{b}j\text{b}}%
-\overset{q}{\underset{\text{d}=1}{\sum}}R_{\text{dbdb}}%
)(\overset{n}{\underset{i=1}{\sum}}R_{i\text{a}i\text{a}}%
-\overset{q}{\underset{\text{c}=1}{\sum}}R_{\text{caca}})$

$=\frac{4}{3}(\varrho_{\text{aa}}-\overset{q}{\underset{\text{c}=1}{\sum}%
}R_{\text{acac}})(\varrho_{\text{bb}}-\overset{q}{\underset{\text{d}=1}{\sum}%
}R_{\text{bdbd}})+\frac{4}{3}(\varrho_{\text{bb}}%
-\overset{q}{\underset{\text{d}=1}{\sum}}R_{\text{bdbd}})(\varrho_{\text{aa}%
}-\overset{q}{\underset{\text{c}=1}{\sum}}R_{\text{acac}})$

$=\frac{8}{3}(\varrho_{\text{aa}}-\overset{q}{\underset{\text{c}=1}{\sum}%
}R_{\text{acac}})(\varrho_{\text{bb}}-\overset{q}{\underset{\text{d}=1}{\sum}%
}R_{\text{bdbd}})$

The next is:

$\frac{2}{3}\overset{n}{\underset{i,j=q+1}{\sum}}(\varrho_{ii}\varrho
_{jj}+2\varrho_{ij}\varrho_{ij})=\frac{2}{3}[\overset{n}{\underset{i=q+1}{\sum
}}(\varrho_{ii}\overset{n}{\underset{j=q+1}{\sum}}\varrho_{jj}]+\frac{4}%
{3}\overset{n}{\underset{i,j=q+1}{\sum}}\varrho_{ij}^{2}$

$=\frac{2}{3}[\overset{n}{\underset{i=1}{(\sum}}\varrho_{ii}%
-\overset{q}{\underset{\text{a=1}}{\sum}}\varrho_{\text{aa}}%
)(\overset{n}{\underset{j=1}{\sum}}\varrho_{jj}%
-\overset{q}{\underset{\text{b=1}}{\sum}}\varrho_{\text{bb}})]+\frac{4}%
{3}[\overset{n}{\underset{i,j=1}{\sum}}\varrho_{ij}^{2}%
-\overset{q}{\underset{\text{a,b}=1}{\sum}}\varrho_{\text{ab}}^{2}]$

$=\frac{2}{3}[(\tau^{M}-\overset{q}{\underset{\text{a=1}}{\sum}}%
\varrho_{\text{aa}}^{M})(\tau^{M}-\overset{q}{\underset{\text{b=1}}{\sum}%
}\varrho_{\text{bb}}^{M})]+\frac{4}{3}[\left\Vert \varrho^{M}\right\Vert
^{2}-\overset{q}{\underset{\text{a,b}=1}{\sum}}(\varrho_{\text{ab}}^{M})^{2}]$

Then we have,

$+\frac{2}{3}\overset{n}{\underset{i,j=q+1}{\sum}}R_{i\text{a}i\text{a}%
}\varrho_{jj}\ +\frac{2}{3}\overset{n}{\underset{i,j=q+1}{\sum}}%
R_{j\text{a}j\text{a}}\varrho_{ii}\ $

$=\frac{2}{3}(\overset{n}{\underset{i=1}{\sum}}R_{i\text{a}i\text{a}%
}.-\overset{q}{\underset{\text{b}=1}{\sum}}R_{\text{baba}}%
)(\overset{n}{\underset{j=1}{\sum}}\varrho_{jj}-\overset{q}{\underset{\text{c}%
=1}{\sum}}\varrho_{\text{cc}})+\frac{2}{3}(\overset{n}{\underset{i=1}{\sum}%
}R_{j\text{a}j\text{a}}.-\overset{q}{\underset{\text{b}=1}{\sum}%
}R_{\text{baba}})(\overset{n}{\underset{j=1}{\sum}}\varrho_{jj}%
-\overset{q}{\underset{\text{c}=1}{\sum}}\varrho_{\text{cc}})$

$=\frac{2}{3}(\varrho_{\text{aa}}^{M}-\varrho_{\text{aa}}^{P})(\tau
^{M}-\overset{q}{\underset{\text{c}=1}{\sum}}\varrho_{\text{cc}}^{M})+\frac
{2}{3}(\varrho_{\text{aa}}^{M}-\varrho_{\text{aa}}^{P})(\tau^{M}%
-\overset{q}{\underset{\text{c}=1}{\sum}}\varrho_{\text{cc}}^{M})$

$=\frac{4}{3}(\varrho_{\text{aa}}^{M}-\varrho_{\text{aa}}^{P})(\tau
^{M}-\overset{q}{\underset{\text{c}=1}{\sum}}\varrho_{\text{cc}}^{M})$

We will often use the \textbf{Gauss Equation} in $\left(  4.28\right)  $ of
\textbf{Gray }$\left[  4\right]  :$

$\ $ $\overset{n}{\underset{i=q+1}{\sum}}$(T$_{\text{ac}i}$T$_{\text{bd}i}%
-$T$_{\text{ad}i}$T$_{\text{bc}i}$) $=$ R$_{\text{abcd}}^{P}-$ R$_{\text{abcd}%
}^{M}$

Given the above properties, the expression for $\frac{\partial^{4}\theta
}{\partial\text{x}_{i}^{2}\partial\text{x}_{j}^{2}}(y_{0})$ is slightly more
geometrically expressed:

\qquad\qquad\qquad\qquad\qquad\qquad\qquad\qquad\qquad\qquad\qquad\qquad
\qquad\qquad\qquad\qquad\qquad\qquad\qquad$\blacksquare$

$\left(  A_{20}\right)  $

\qquad\qquad\qquad\qquad\qquad\qquad\qquad\qquad\qquad\qquad\qquad\qquad
\qquad\qquad\qquad\qquad\qquad\qquad\qquad$\blacksquare$

The "horde" of Second Fundamental Forms in the last expression above contain
terms having curvature differences which can be paired. Each member of the
same pair is marked with \ the same number. The numbers run from 1 to 6.

$+\ \frac{1}{6}\underset{i=q+1}{\overset{n}{\sum}}[T_{\text{aa}i}%
T_{\text{bb}i}(R_{\text{cdcd}}^{P}-R_{\text{cdcd}}^{M}%
)+\underset{j=q+1}{\overset{n}{\sum}}T_{\text{aa}j}T_{\text{bb}j}%
(R_{\text{cdcd}}^{P}-R_{\text{cdcd}}^{M})](y_{0})\qquad(1)$

$-\ \frac{1}{6}[\underset{i=q+1}{\overset{n}{\sum}}T_{\text{ab}i}%
T_{\text{ab}i}(R_{\text{cdcd}}^{P}-R_{\text{cdcd}}^{M}%
)+\underset{j=q+1}{\overset{n}{\sum}}T_{\text{ab}j}T_{\text{ab}j}%
(R_{\text{cdcd}}^{P}-R_{\text{cdcd}}^{M})](y_{0})\qquad\left(  1\right)  $

$+\frac{1}{6}[\underset{i=q+1}{\overset{n}{\sum}}T_{\text{ac}i}T_{\text{ab}%
i}(R_{\text{bdcd}}^{P}-R_{\text{bdcd}}^{M})+\underset{j=q+1}{\overset{n}{\sum
}}T_{\text{ac}j}T_{\text{ab}j}(R_{\text{bdcd}}^{P}-R_{\text{bdcd}}^{M}%
)](y_{0})\qquad(2)$

$-\ \frac{1}{6}[\underset{i=q+1}{\overset{n}{\sum}}T_{\text{aa}i}%
T_{\text{bc}i}(R_{\text{bdcd}}^{P}-R_{\text{bdcd}}^{M}%
)+\underset{j=q+1}{\overset{n}{\sum}}T_{\text{aa}j}T_{\text{bc}j}%
(R_{\text{bdcd}}^{P}-R_{\text{bdcd}}^{M})](y_{0})\qquad(2)$

$+\ \frac{1}{6}[\underset{i=q+1}{\overset{n}{\sum}}T_{\text{aa}i}%
T_{\text{bd}i}(R_{\text{bccd}}^{P}-R_{\text{bccd}}^{M}%
)+\underset{j=q+1}{\overset{n}{\sum}}T_{\text{aa}j}T_{\text{bd}j}%
(R_{\text{bccd}}^{P}-R_{\text{bccd}}^{M})](y_{0})\qquad(3)$

$-\ \frac{1}{6}[\underset{i=q+1}{\overset{n}{\sum}}T_{\text{ad}i}%
T_{\text{ab}i}(R_{\text{bccd}}^{P}-R_{\text{abcd}}^{M}%
)+\underset{j=q+1}{\overset{n}{\sum}}T_{\text{ad}j}T_{\text{ab}j}%
(R_{\text{bccd}}^{P}-R_{\text{bccd}}^{M})](y_{0})\qquad(3)$

$+\ \frac{1}{6}[\underset{i=q+1}{\overset{n}{\sum}}T_{\text{ab}i}%
T_{\text{bc}i}(R_{\text{adcd}}^{P}-R_{\text{adcd}}^{M}%
)+\underset{j=q+1}{\overset{n}{\sum}}T_{\text{ab}j}T_{\text{bc}j}%
(R_{\text{adcd}}^{P}-R_{\text{adcd}}^{M})](y_{0})\qquad(4)$

$-\ \frac{1}{6}[\underset{i=q+1}{\overset{n}{\sum}}T_{\text{ac}i}%
T_{\text{bb}i}(R_{\text{adcd}}^{P}-R_{\text{adcd}}^{M}%
)+\underset{j=q+1}{\overset{n}{\sum}}T_{\text{ac}j}T_{\text{bb}j}%
(R_{\text{adcd}}^{P}-R_{\text{adcd}}^{M})](y_{0})\qquad(4)$

$+\ \frac{1}{6}[\underset{i=q+1}{\overset{n}{\sum}}T_{\text{ad}i}%
T_{\text{bb}i}(R_{\text{accd}}^{P}-R_{\text{accd}}^{M}%
)+\underset{j=q+1}{\overset{n}{\sum}}T_{\text{ad}j}T_{\text{bb}j}%
(R_{\text{accd}}^{P}-R_{\text{accd}}^{M})](y_{0})\qquad(5)$

$-\ \frac{1}{6}[\underset{i=q+1}{\overset{n}{\sum}}T_{\text{ab}i}%
T_{\text{bd}i}(R_{\text{accd}}^{P}-R_{\text{accd}}^{M}%
)+\underset{j=q+1}{\overset{n}{\sum}}T_{\text{ab}j}T_{\text{bd}j}%
(R_{\text{accd}}^{P}-R_{\text{accd}}^{M})](y_{0})\qquad(5)$

$+\ \frac{1}{6}[\underset{i=q+1}{\overset{n}{\sum}}T_{\text{ac}i}%
T_{\text{bd}i}(R_{\text{abcd}}^{P}-R_{\text{abcd}}^{M}%
)+\underset{j=q+1}{\overset{n}{\sum}}T_{\text{ac}j}T_{\text{bd}j}%
(R_{\text{abcd}}^{P}-R_{\text{abcd}}^{M})](y_{0})\qquad(6)$

$-\ \frac{1}{6}[\underset{i=q+1}{\overset{n}{\sum}}T_{\text{ad}i}%
T_{\text{bc}i}(R_{\text{abcd}}^{P}-R_{\text{abcd}}^{M}%
)+\underset{j=q+1}{\overset{n}{\sum}}T_{\text{ad}j}T_{\text{bc}j}%
(R_{\text{abcd}}^{P}-R_{\text{abcd}}^{M})](y_{0})\qquad(6)$

Each of the above pair is factorable and the Gauss Equation is again
applicable and we have:

$=+\ \frac{1}{3}[(R_{\text{cdcd}}^{P}-R_{\text{cdcd}}^{M})(R_{\text{abab}}%
^{P}-R_{\text{abab}}^{M})](y_{0})\qquad\left(  1\right)  $

$\ -\frac{1}{3}[(R_{\text{bdcd}}^{P}-R_{\text{bdcd}}^{M})(R_{\text{abac}}%
^{P}-R_{\text{abac}}^{M})](y_{0})\qquad(2)$

$\ -\ \frac{1}{3}[(R_{\text{bcdc}}^{P}-R_{\text{bcdc}}^{M})(R_{\text{abad}%
}^{P}-R_{\text{abad}}^{M})](y_{0})\qquad(3)$

$\qquad\ +\ \frac{1}{3}[(R_{\text{adcd}}^{P}-R_{\text{adcd}}^{M}%
)(R_{\text{abbc}}^{P}-R_{\text{abbc}}^{M})](y_{0})\qquad(4)\qquad$

$\ -\ \frac{1}{3}[(R_{\text{acdc}}^{P}-R_{\text{acdc}}^{M})(R_{\text{abdb}%
}^{P}-R_{\text{abdb}}^{M})](y_{0})\qquad(5)$

$\ +\ \frac{1}{3}[(R_{\text{abcd}}^{P}-R_{\text{abcd}}^{M})]^{2}(y_{0}%
)\qquad(6)$

We see that $\overset{n}{\underset{i,j=q+1}{\sum}}\frac{\partial^{4}\theta
_{p}}{\partial\text{x}_{i}^{2}\partial\text{x}_{j}^{2}}(y_{0})$ expressed in
more refined \textbf{geometric invariants} is:

(xx)$\qquad\frac{\partial^{4}\theta_{p}}{\partial\text{x}_{i}^{2}%
\partial\text{x}_{j}^{2}}(y_{0})=\frac{1}{6}%
\overset{n}{\underset{i,j=q+1}{\sum}}[$ $\overset{q}{\underset{\text{a=1}%
}{\sum}}\{-(\nabla_{ii}^{2}R_{j\text{a}j\text{a}}+\nabla_{jj}^{2}%
R_{i\text{a}i\text{a}}+4\nabla_{ij}^{2}R_{i\text{a}j\text{a}}+2R_{ij}%
R_{i\text{a}j\text{a}})\qquad A$

$+\overset{n}{\underset{p=q+1}{\sum}}\overset{q}{\underset{\text{a=1}}{\sum}%
}(R_{\text{a}iip}R_{\text{a}jjp}+R_{\text{a}jjp}R_{\text{a}iip}+R_{\text{a}%
ijp}R_{\text{a}ijp}+R_{\text{a}ijp}R_{\text{a}jip}+R_{\text{a}jip}%
R_{\text{a}ijp}+R_{\text{a}jip}R_{\text{a}jip})$

$+2\overset{q}{\underset{\text{a,b=1}}{\sum}}\nabla_{i}(R)_{\text{a}%
i\text{b}j}T_{\text{ab}j}+2\overset{q}{\underset{\text{a,b=1}}{\sum}}%
\nabla_{j}(R)_{\text{a}j\text{b}i}T_{\text{ab}i}%
+2\overset{q}{\underset{\text{a,b=1}}{\sum}}\nabla_{i}(R)_{\text{a}j\text{b}%
i}T_{\text{ab}j}+2\overset{q}{\underset{\text{a,b=1}}{\sum}}\nabla
_{i}(R)_{\text{a}j\text{b}j}T_{\text{ab}i}$

$+2\overset{q}{\underset{\text{a,b=1}}{\sum}}\nabla_{j}(R)_{\text{a}%
i\text{b}i}T_{\text{ab}j}+2\overset{q}{\underset{\text{a,b=1}}{\sum}}%
\nabla_{j}(R)_{\text{a}i\text{b}j}T_{\text{ab}i}$

$+\overset{n}{\underset{p=q+1}{\sum}}(-\frac{3}{5}\nabla_{ii}^{2}%
(R)_{jpjp}+\overset{n}{\underset{p=q+1}{\sum}}(-\frac{3}{5}\nabla_{jj}%
^{2}(R)_{ipip}$

$+\overset{n}{\underset{p=q+1}{\sum}}(-\frac{3}{5}\nabla_{ij}^{2}%
(R)_{ipjp}+\overset{n}{\underset{p=q+1}{\sum}}(-\frac{3}{5}\nabla_{ij}%
^{2}(R)_{jpip}$

$+\overset{n}{\underset{p=q+1}{\sum}}(-\frac{3}{5}\nabla_{ji}^{2}%
(R)_{ipjp}+\overset{n}{\underset{p=q+1}{\sum}}(-\frac{3}{5}\nabla_{ji}%
^{2}(R)_{jpip}$

$+\frac{1}{5}\overset{n}{\underset{m,p=q+1}{%
{\textstyle\sum}
}}R_{ipim}R_{jpjm}+\frac{1}{5}\overset{n}{\underset{m,p=q+1}{%
{\textstyle\sum}
}}R_{jpjm}R_{ipim}$

$+\frac{1}{5}\overset{n}{\underset{m,p=q+1}{%
{\textstyle\sum}
}}R_{ipjm}R_{ipjm}+\frac{1}{5}\overset{n}{\underset{m,p=q+1}{%
{\textstyle\sum}
}}R_{ipjm}R_{jpim}$

$+\frac{1}{5}\overset{n}{\underset{m,p=q+1}{%
{\textstyle\sum}
}}R_{jpim}R_{ipjm}+\frac{1}{5}\overset{n}{\underset{m,p=q+1}{%
{\textstyle\sum}
}}R_{jpim}R_{jpim}\}(y_{0})$

$+4\overset{q}{\underset{\text{a,b=1}}{\sum}}\{(\nabla_{i}(R)_{i\text{a}%
j\text{a}}-\overset{q}{\underset{\text{c=1}}{%
{\textstyle\sum}
}}R_{\text{a}i\text{c}i}T_{\text{ac}j})$ $T_{\text{bb}j}+4(\nabla
_{j}(R)_{j\text{a}i\text{a}}-\overset{q}{\underset{\text{c=1}}{%
{\textstyle\sum}
}}R_{\text{a}j\text{c}j}T_{\text{ac}i})$ $T_{\text{bb}i}$

$+4(\nabla_{i}(R)_{j\text{a}i\text{a}}-\overset{q}{\underset{\text{c=1}}{%
{\textstyle\sum}
}}R_{\text{a}i\text{c}j}T_{\text{ac}i})$ $T_{\text{bb}j}$ $4B\ $

$+4(\nabla_{i}(R)_{j\text{a}j\text{a}}-\overset{q}{\underset{\text{c=1}}{%
{\textstyle\sum}
}}R_{\text{a}i\text{c}j}T_{\text{ac}j})$ $T_{\text{bb}i}+4(\nabla
_{j}(R)_{i\text{a}i\text{a}}-\overset{q}{\underset{\text{c=1}}{%
{\textstyle\sum}
}}R_{\text{a}j\text{c}i}T_{\text{ac}i})$ $T_{\text{bb}j}$

$+4(\nabla_{j}(R)_{i\text{a}j\text{a}}-\overset{q}{\underset{\text{c=1}}{%
{\textstyle\sum}
}}R_{\text{a}j\text{c}i}T_{\text{ac}j})$ $T_{\text{bb}i}$

$-4\overset{q}{\underset{\text{a,b=1}}{\sum}}(\nabla_{i}(R)_{i\text{a}%
j\text{b}}-\overset{q}{\underset{\text{c=1}}{%
{\textstyle\sum}
}}R_{\text{b}r\text{c}s}T_{\text{ac}t})T_{\text{ab}j}%
-4\overset{q}{\underset{\text{a,b=1}}{\sum}}(\nabla_{j}(R)_{j\text{a}%
i\text{b}}-\overset{q}{\underset{\text{c=1}}{%
{\textstyle\sum}
}}R_{\text{b}j\text{c}j}T_{\text{ac}i})T_{\text{ab}i}$

$-4\overset{q}{\underset{\text{a,b=1}}{\sum}}(\nabla_{i}(R)_{j\text{a}%
i\text{b}}-\overset{q}{\underset{\text{c=1}}{%
{\textstyle\sum}
}}R_{\text{b}i\text{c}j}T_{\text{ac}i})T_{\text{ab}j}%
-4\overset{q}{\underset{\text{a,b=1}}{\sum}}(\nabla_{i}(R)_{j\text{a}%
j\text{b}}-\overset{q}{\underset{\text{c=1}}{%
{\textstyle\sum}
}}R_{\text{b}i\text{c}j}T_{\text{ac}j})T_{\text{ab}i}$

$-4\overset{q}{\underset{\text{a,b=1}}{\sum}}(\nabla_{j}(R)_{i\text{a}%
i\text{b}}-\overset{q}{\underset{\text{c=1}}{%
{\textstyle\sum}
}}R_{\text{b}j\text{c}i}T_{\text{ac}i})T_{\text{ab}j}%
-4\overset{q}{\underset{\text{a,b=1}}{\sum}}(\nabla_{j}(R)_{i\text{a}%
j\text{b}}-\overset{q}{\underset{\text{c=1}}{%
{\textstyle\sum}
}}R_{\text{b}j\text{c}i}T_{\text{ac}j})T_{\text{ab}i}\}](y_{0})$

$+$ $[\frac{4}{9}\overset{q}{\underset{\text{a,b=1}}{\sum}}(\varrho
_{\text{aa}}-\overset{q}{\underset{\text{c}=1}{\sum}}R_{\text{acac}}%
)(\varrho_{\text{bb}}-\overset{q}{\underset{\text{d}=1}{\sum}}R_{\text{bdbd}%
})+\frac{8}{9}\overset{n}{\underset{i,j=q+1}{\sum}}%
\overset{q}{\underset{\text{a,b}=1}{\sum}}(R_{i\text{a}j\text{a}}%
R_{i\text{b}j\text{b}})\qquad3C$

$+\frac{4}{3}\overset{q}{\underset{\text{a}=1}{\sum}}(\varrho_{\text{aa}}%
^{M}-\varrho_{\text{aa}}^{P})(\tau^{M}-\overset{q}{\underset{\text{c}=1}{\sum
}}\varrho_{\text{cc}}^{M})+\frac{4}{9}\overset{n}{\underset{i,j=q+1}{\sum}%
}\overset{q}{\underset{\text{a}=1}{\sum}}R_{i\text{a}j\text{a}}\varrho_{ij}\ $

$\ +\frac{2}{9}\overset{q}{\underset{\text{b}=1}{\sum}}(\varrho_{\text{bb}%
}^{M}-\varrho_{\text{bb}}^{P})(\tau^{M}-\overset{q}{\underset{\text{c}%
=1}{\sum}}\varrho_{\text{cc}}^{M})+\frac{4}{9}%
\overset{n}{\underset{i,j=q+1}{\sum}}\overset{q}{\underset{\text{b}=1}{\sum}%
}R_{i\text{b}j\text{b}}\varrho_{ij}\ $

$+\frac{1}{9}(\tau^{M}-\overset{q}{\underset{\text{a=1}}{\sum}}\varrho
_{\text{aa}})(\tau^{M}-\overset{q}{\underset{\text{b=1}}{\sum}}\varrho
_{\text{bb}})+\frac{2}{9}(\left\Vert \varrho^{M}\right\Vert ^{2}%
-\overset{q}{\underset{\text{a,b}=1}{\sum}}\varrho_{\text{ab}})$

$-\overset{n}{\underset{i,j=q+1}{\sum}}\overset{q}{\underset{\text{a,b}%
=1}{\sum}}R_{i\text{a}i\text{b}}R_{j\text{a}j\text{b}}\ -\frac{1}%
{2}\overset{n}{\underset{i,j=q+1}{\sum}}\overset{q}{\underset{\text{a,b}%
=1}{\sum}}R_{i\text{a}j\text{b}}^{2}-\overset{n}{\underset{i,j=q+1}{\sum}%
}\overset{q}{\underset{\text{a,b}=1}{\sum}}R_{i\text{a}j\text{b}}%
R_{j\text{a}i\text{b}}$

$-\frac{1}{2}\overset{n}{\underset{i,j=q+1}{\sum}}%
\overset{q}{\underset{\text{a,b}=1}{\sum}}R_{j\text{a}i\text{b}}^{2}$

$-\frac{1}{9}\overset{n}{\underset{i,j,p,m=q+1}{\sum}}R_{ipim}R_{jpjm}%
\ -\frac{1}{18}\overset{n}{\underset{i,j,p,m=q+1}{\sum}}R_{ipjm}^{2}-\frac
{1}{9}\overset{n}{\underset{i,j,p,m=q+1}{\sum}}R_{ipjm}R_{jpim}$

$-\frac{1}{18}\overset{n}{\underset{i,j,p,m=q+1}{\sum}}R_{jpim}^{2}$

$-\frac{1}{3}\overset{q}{\underset{\text{a}=1}{\sum}}%
\overset{n}{\underset{i,j,p=q+1}{\sum}}R_{i\text{a}ip}R_{j\text{a}jp}-\frac
{1}{6}\overset{q}{\underset{\text{a}=1}{\sum}}%
\overset{n}{\underset{i,j,p=q+1}{\sum}}R_{i\text{a}jp}^{2}-\frac{1}%
{3}\overset{q}{\underset{\text{a}=1i,j,}{\sum}}%
\overset{n}{\underset{p=q+1}{\sum}}R_{i\text{a}jp}R_{j\text{a}ip}$

$-\frac{1}{6}\overset{q}{\underset{\text{a}=1}{\sum}}%
\overset{n}{\underset{i,j,p=q+1}{\sum}}R_{j\text{a}ip}^{2}$

$-\frac{1}{3}\overset{q}{\underset{\text{b}=1i,j,}{\sum}}%
\overset{n}{\underset{p=q+1}{\sum}}R_{i\text{b}ip}R_{j\text{b}jp}-\frac{1}%
{6}\overset{q}{\underset{\text{b}=1}{\sum}}%
\overset{n}{\underset{i,j,p=q+1}{\sum}}R_{i\text{b}jp}^{2}-\frac{1}%
{3}\overset{q}{\underset{\text{b}=1}{\sum}}%
\overset{n}{\underset{i.j,p=q+1}{\sum}}R_{i\text{b}jp}R_{j\text{b}ip}$

$-\frac{1}{6}\overset{q}{\underset{\text{b}=1}{\sum}}%
\overset{n}{\underset{i,j,p=q+1}{\sum}}R_{j\text{b}ip}^{2}](y_{0})$

$+$ $\overset{q}{\underset{\text{a,b,c=1}}{\sum}}[$
$-\overset{n}{\underset{i=q+1}{\sum}}R_{i\text{a}i\text{a}}(R_{\text{bcbc}%
}^{P}-R_{\text{bcbc}}^{M})$ $-\overset{n}{\underset{j=q+1}{\sum}}%
R_{j\text{a}j\text{a}}(R_{\text{bcbc}}^{P}-R_{\text{bcbc}}^{M})\qquad\qquad6D$

\ $+\overset{n}{\underset{i=q+1}{\sum}}R_{i\text{a}i\text{b}}(R_{\text{acbc}%
}^{P}-R_{\text{acbc}}^{M})\ -\overset{n}{\underset{i=q+1}{\sum}}%
R_{i\text{a}i\text{c}}(R_{\text{abbc}}^{P}-R_{\text{abbc}}^{M})$

$+\overset{n}{\underset{j=q+1}{\sum}}R_{j\text{a}j\text{b}}(R_{\text{acbc}%
}^{P}-R_{\text{acbc}}^{M})$\ $-\overset{n}{\underset{j=q+1}{\sum}}%
R_{j\text{a}j\text{c}}(R_{\text{abbc}}^{P}-R_{\text{abbc}}^{M})$

$+\underset{i,j=q+1}{\overset{n}{\sum}}$ $-R_{i\text{a}j\text{a}}%
(T_{\text{bb}i}T_{\text{cc}j}$ $-T_{\text{bc}i}T_{\text{bc}j})$
$-\underset{i,j=q+1}{\overset{n}{\sum}}R_{i\text{a}j\text{a}}(T_{\text{bb}%
j}T_{\text{cc}i}$ $-T_{\text{bc}j}T_{\text{bc}i})$

$+$ $\underset{i,j=q+1}{\overset{n}{\sum}}$ $-R_{j\text{a}i\text{a}%
}(T_{\text{bb}i}T_{\text{cc}j}$ $-T_{\text{bc}i}T_{\text{bc}j})$
$-\underset{i,j=q+1}{\overset{n}{\sum}}R_{j\text{a}i\text{a}}(T_{\text{bb}%
j}T_{\text{cc}i}$ $-T_{\text{bc}j}T_{\text{bc}i})$

$+\underset{i,j=q+1}{\overset{n}{\sum}}\ R_{i\text{a}j\text{b}}(T_{\text{ab}%
i}T_{\text{cc}j}-T_{\text{bc}i}T_{\text{ac}j}%
)\ +\underset{i,j=q+1}{\overset{n}{\sum}}\ R_{i\text{a}j\text{b}}%
(T_{\text{ab}j}T_{\text{cc}i}-T_{\text{bc}j}T_{\text{ac}i})$

$+\underset{i,j=q+1}{\overset{n}{\sum}}\ R_{j\text{a}i\text{ib}}%
(T_{\text{ab}i}T_{\text{cc}j}-T_{\text{bc}i}T_{\text{ac}j}%
)\ +\underset{i,j=q+1}{\overset{n}{\sum}}\ R_{j\text{a}i\text{b}}%
(T_{\text{ab}j}T_{\text{cc}i}-T_{\text{bc}j}T_{\text{ac}i})\qquad$

$+\underset{i,j=q+1}{\overset{n}{\sum}}-R_{i\text{a}j\text{c}}(T_{\text{ab}%
i}T_{\text{bc}j}-T_{\text{ac}i}T_{\text{bb}j}%
)-\underset{i,j=q+1}{\overset{n}{\sum}}R_{i\text{a}j\text{c}}(T_{\text{ba}%
j}T_{\text{bc}i}-T_{\text{ac}j}T_{\text{bb}i})$

$+\underset{i,j=q+1}{\overset{n}{\sum}}-R_{j\text{a}i\text{c}}(T_{\text{ba}%
i}T_{\text{bc}j}-T_{\text{ac}i}T_{\text{bb}j}%
)-\underset{i,j=q+1}{\overset{n}{\sum}}R_{j\text{a}i\text{c}}(T_{\text{ba}%
j}T_{\text{bc}i}-T_{\text{ac}j}T_{\text{bb}i})](y_{0})$

$-\frac{1}{3}\underset{p=q+1}{\overset{n}{\sum}}%
[\underset{i=q+1}{\overset{n}{\sum}}\overset{q}{\underset{\text{b,c=1}}{\sum}%
}R_{ipip}(R_{\text{bcbc}}^{P}-R_{\text{bcbc}}^{M}%
)+\underset{j=q+1}{\overset{n}{\sum}}$ $\overset{q}{\underset{\text{b,c=1}%
}{\sum}}R_{jpjp}(R_{\text{bcbc}}^{P}-R_{\text{bcbc}}^{M})](y_{0})$

$-\frac{2}{3}\underset{i,j,p=q+1}{\overset{n}{\sum}}%
\overset{q}{\underset{\text{b,c=1}}{\sum}}[R_{ipjp}(T_{\text{bb}i}%
T_{\text{cc}j}-T_{\text{bc}i}T_{\text{bc}j})+R_{ipjp}(T_{\text{bb}%
j}T_{\text{cc}i}-T_{\text{bc}j}T_{\text{bc}i})](y_{0})\qquad$

$+\ \frac{1}{6}\underset{i,j=q+1}{\overset{n}{\sum}}[T_{\text{aa}%
i}T_{\text{bb}j}(T_{\text{cc}i}T_{\text{dd}j}-T_{\text{cd}i}T_{\text{dc}%
j})+T_{\text{aa}i}T_{\text{bb}j}(T_{\text{cc}j}T_{\text{dd}i}-T_{\text{cd}%
j}T_{\text{dc}i})\qquad E$

$+T_{\text{aa}j}T_{\text{bb}i}(T_{\text{cc}i}T_{\text{dd}j}-T_{\text{cd}%
i}T_{\text{dc}j})+T_{\text{aa}j}T_{\text{bb}i}(T_{\text{cc}j}T_{\text{dd}%
i}-T_{\text{cd}j}T_{\text{dc}i})](y_{0})$

$-\frac{1}{6}\underset{i,j=q+1}{\overset{n}{\sum}}[T_{\text{aa}i}%
T_{\text{bc}j}(T_{\text{bc}i}T_{\text{dd}j}-T_{\text{bd}i}T_{\text{cd}%
j})+T_{\text{aa}i}T_{\text{bc}j}(T_{\text{bc}j}T_{\text{dd}i}-T_{\text{bd}%
j}T_{\text{cd}i})$

$+T_{\text{aa}j}T_{\text{bc}i}(T_{\text{bc}i}T_{\text{dd}j}-T_{\text{bd}%
i}T_{\text{cd}j})+T_{\text{aa}j}T_{\text{bc}i}(T_{\text{bc}j}T_{\text{dd}%
i}-T_{\text{bd}j}T_{\text{cd}i})](y_{0})$

$+\ \frac{1}{6}\underset{i,j=q+1}{\overset{n}{\sum}}[T_{\text{aa}%
i}T_{\text{bd}j}(T_{\text{bc}i}T_{\text{cd}j}-T_{\text{bd}i}T_{\text{cc}%
j})+T_{\text{aa}i}T_{\text{bd}j}(T_{\text{bc}j}T_{\text{cd}i}-T_{\text{bd}%
j}T_{\text{cc}i})$

$+T_{\text{aa}j}T_{\text{bd}i}(T_{\text{bc}i}T_{\text{cd}j}-T_{\text{bd}%
i}T_{\text{cc}j})+T_{\text{aa}j}T_{\text{bd}i}(T_{\text{bc}j}T_{\text{cd}%
i}-T_{\text{bd}j}T_{\text{cc}i})](y_{0})\qquad$

$-\ \frac{1}{6}\underset{i,j=q+1}{\overset{n}{\sum}}[T_{\text{ab}%
i}T_{\text{ab}j}(T_{\text{cc}i}T_{\text{dd}j}-T_{\text{cd}i}T_{\text{dc}%
j})+T_{\text{ab}i}T_{\text{ab}j}(T_{\text{cc}j}T_{\text{dd}i}-T_{\text{cd}%
j}T_{\text{dc}i})$

$+T_{\text{ab}j}T_{\text{ab}i}(T_{\text{cc}i}T_{\text{dd}j}-T_{\text{cd}%
i}T_{\text{dc}j})+T_{\text{ab}j}T_{\text{ab}i}(T_{\text{cc}j}T_{\text{dd}%
i}-T_{\text{cd}j}T_{\text{dc}i})](y_{0})$

$+\ \frac{1}{6}\underset{i,j=q+1}{\overset{n}{\sum}}[T_{\text{ab}%
i}T_{\text{bc}j}(T_{\text{ac}i}T_{\text{dd}j}-T_{\text{ad}i}T_{\text{cd}%
j})+T_{\text{ab}i}T_{\text{bc}j}(T_{\text{ac}j}T_{\text{dd}i}-T_{\text{ad}%
j}T_{\text{cd}i})$

$+T_{\text{ab}j}T_{\text{bc}i}(T_{\text{ac}i}T_{\text{dd}j}-T_{\text{ad}%
i}T_{\text{cd}j})+T_{\text{ab}j}T_{\text{bc}i}(T_{\text{ac}j}T_{\text{dd}%
i}-T_{\text{ad}j}T_{\text{cd}i})](y_{0})$

$-\ \frac{1}{6}\underset{i,j=q+1}{\overset{n}{\sum}}[T_{\text{ab}%
i}T_{\text{bd}j}(T_{\text{ac}i}T_{\text{cd}j}-T_{\text{ad}i}T_{\text{cc}%
j})+T_{\text{ab}i}T_{\text{bd}j}(T_{\text{ac}j}T_{\text{cd}i}-T_{\text{ad}%
j}T_{\text{cc}i})$

$+T_{\text{ab}i}T_{\text{bd}j}(T_{\text{ac}j}T_{\text{cd}i}-T_{\text{ad}%
j}T_{\text{cc}i})+T_{\text{ab}j}T_{\text{bd}i}(T_{\text{ac}j}T_{\text{cd}%
i}-T_{\text{ad}j}T_{\text{cc}i})](y_{0})$

$+\ \frac{1}{6}\underset{i,j=q+1}{\overset{n}{\sum}}[T_{\text{ac}%
i}T_{\text{ab}j}(T_{\text{bc}i}T_{\text{dd}j}-T_{\text{bd}i}T_{\text{dc}%
j})+T_{\text{ac}i}T_{\text{ab}j}(T_{\text{bc}j}T_{\text{dd}i}-T_{\text{bd}%
j}T_{\text{dc}i})$

$+T_{\text{ac}j}T_{\text{ab}i}(T_{\text{bc}i}T_{\text{dd}j}-T_{\text{bd}%
i}T_{\text{dc}j})+T_{\text{ac}j}T_{\text{ab}i}(T_{\text{bc}j}T_{\text{dd}%
i}-T_{\text{bd}j}T_{\text{dc}i})](y_{0})$

$-\ \frac{1}{6}\underset{i,j=q+1}{\overset{n}{\sum}}[T_{\text{ac}%
i}T_{\text{bb}j}(T_{\text{ac}i}T_{\text{dd}j}-T_{\text{ad}i}T_{\text{cd}%
j})+T_{\text{ac}i}T_{\text{bb}j}(T_{\text{ac}j}T_{\text{dd}i}-T_{\text{ad}%
j}T_{\text{cd}i})$

$+T_{\text{ac}j}T_{\text{bb}i}(T_{\text{ac}i}T_{\text{dd}j}-T_{\text{ad}%
i}T_{\text{cd}i})+T_{\text{ac}j}T_{\text{bb}i}(T_{\text{ac}j}T_{\text{dd}%
i}-T_{\text{ad}j}T_{\text{cd}i})](y_{0})$

$+\ \frac{1}{6}\underset{i,j=q+1}{\overset{n}{\sum}}[T_{\text{ac}%
i}T_{\text{bd}j}(T_{\text{ac}i}T_{\text{bd}j}-T_{\text{ad}i}T_{\text{bc}%
j})+T_{\text{ac}i}T_{\text{bd}j}(T_{\text{ac}j}T_{\text{bd}i}-T_{\text{ad}%
j}T_{\text{bc}i})$

$+T_{\text{ac}j}T_{\text{bd}i}(T_{\text{ac}i}T_{\text{bd}j}-T_{\text{ad}%
i}T_{\text{bc}j})+T_{\text{ac}j}T_{\text{bd}i}(T_{\text{ac}j}T_{\text{bd}%
i}-T_{\text{ad}j}T_{\text{bc}i})](y_{0})$

$-\frac{1}{6}\underset{i,j=q+1}{\overset{n}{\sum}}[T_{\text{ad}i}%
T_{\text{ab}j}(T_{\text{bc}i}T_{\text{cd}j}-T_{\text{bd}i}T_{\text{cc}%
j})+T_{\text{ad}i}T_{\text{ab}j}(T_{\text{bc}j}T_{\text{cd}i}-T_{\text{bd}%
j}T_{\text{cc}i})$

$+T_{\text{ad}j}T_{\text{ab}i}(T_{\text{bc}i}T_{\text{cd}j}-T_{\text{bd}%
i}T_{\text{cc}j})+T_{\text{ad}j}T_{\text{ab}i}(T_{\text{bc}j}T_{\text{cd}%
i}-T_{\text{bd}j}T_{\text{cc}i})](y_{0})$

$+\ \frac{1}{6}\underset{i,j=q+1}{\overset{n}{\sum}}[T_{\text{ad}%
i}T_{\text{bb}j}(T_{\text{ac}i}T_{\text{cd}j}-T_{\text{ad}i}T_{\text{cc}%
j})+T_{\text{ad}i}T_{\text{bb}j}(T_{\text{ac}j}T_{\text{cd}i}-T_{\text{ad}%
j}T_{\text{cc}i})$

$+T_{\text{ad}j}T_{\text{bb}i}(T_{\text{ac}i}T_{\text{cd}j}-T_{\text{ad}%
i}T_{\text{cc}j})+T_{\text{ad}j}T_{\text{bb}i}(T_{\text{ac}j}T_{\text{cd}%
i}-T_{\text{ad}j}T_{\text{cc}i})](y_{0})$

$-\ \frac{1}{6}\underset{i,j=q+1}{\overset{n}{\sum}}[T_{\text{ad}%
i}T_{\text{bc}j}(T_{\text{ac}i}T_{\text{bd}j}-T_{\text{ad}i}T_{\text{bc}%
j})+T_{\text{ad}i}T_{\text{bc}j}(T_{\text{ac}j}T_{\text{bd}i}-T_{\text{ad}%
j}T_{\text{bc}i})$

$+T_{\text{ad}j}T_{\text{bc}i}(T_{\text{ac}i}T_{\text{bd}j}-T_{\text{ad}%
i}T_{\text{bc}j})+T_{\text{ad}j}T_{\text{bc}i}(T_{\text{ac}j}T_{\text{bd}%
i}-T_{\text{ad}j}T_{\text{bc}i})](y_{0})$

$=+\ \frac{1}{3}[(R_{\text{cdcd}}^{P}-R_{\text{cdcd}}^{M})(R_{\text{abab}}%
^{P}-R_{\text{abab}}^{M})](y_{0})\qquad\left(  1\right)  $

$\ -\frac{1}{3}[(R_{\text{bdcd}}^{P}-R_{\text{bdcd}}^{M})(R_{\text{abac}}%
^{P}-R_{\text{abac}}^{M})](y_{0})\qquad(2)$

$\ -\ \frac{1}{3}[(R_{\text{bcdc}}^{P}-R_{\text{bcdc}}^{M})(R_{\text{abad}%
}^{P}-R_{\text{abad}}^{M})](y_{0})\qquad(3)$

$\qquad\ +\ \frac{1}{3}[(R_{\text{adcd}}^{P}-R_{\text{adcd}}^{M}%
)(R_{\text{abbc}}^{P}-R_{\text{abbc}}^{M})](y_{0})\qquad(4)\qquad$

$\ -\ \frac{1}{3}[(R_{\text{acdc}}^{P}-R_{\text{acdc}}^{M})(R_{\text{abdb}%
}^{P}-R_{\text{abdb}}^{M})](y_{0})\qquad(5)$

$\ +\ \frac{1}{3}[(R_{\text{abcd}}^{P}-R_{\text{abcd}}^{M})]^{2}(y_{0}%
)\qquad(6)$

$\qquad\qquad\qquad\qquad\qquad\qquad\qquad\qquad\qquad\qquad\qquad
\qquad\qquad\qquad\qquad\qquad\blacksquare$

(xix) We compute the expression for $\ \frac{\partial^{4}\theta^{-\frac{1}{2}%
}}{\partial\text{x}_{i}^{2}\partial\text{x}_{j}^{2}}(y_{0}):$

Since,

\ $\frac{\partial^{2}\theta^{-\frac{1}{2}}}{\partial\text{x}_{j}^{2}}%
=\frac{\partial}{\partial\text{x}_{j}}(\frac{\partial\theta^{-\frac{1}{2}}%
}{\partial\text{x}_{j}})=\frac{\partial}{\partial\text{x}_{j}}(-\frac{1}%
{2}\theta^{-\frac{3}{2}}\frac{\partial\theta}{\partial\text{x}_{j}})=-\frac
{1}{2}[\frac{\partial\theta^{-\frac{3}{2}}}{\partial\text{x}_{j}}%
\frac{\partial\theta}{\partial\text{x}_{j}}+\theta^{-\frac{3}{2}}%
\frac{\partial^{2}\theta}{\partial\text{x}_{j}^{2}}],$

we have,

$\frac{\partial^{4}\theta^{-\frac{1}{2}}}{\partial\text{x}_{i}^{2}%
\partial\text{x}_{j}^{2}}=$ $\frac{\partial^{2}}{\partial\text{x}_{i}^{2}%
}(\frac{\partial^{2}\theta^{-\frac{1}{2}}}{\partial\text{x}_{j}^{2}}%
)=-\frac{1}{2}\frac{\partial^{2}}{\partial\text{x}_{i}^{2}}[\frac
{\partial\theta^{-\frac{3}{2}}}{\partial\text{x}_{j}}\frac{\partial\theta
}{\partial\text{x}_{j}}+\theta^{-\frac{3}{2}}\frac{\partial^{2}\theta
}{\partial\text{x}_{j}^{2}}]$

$=-\frac{1}{2}[\frac{\partial^{3}\theta^{-\frac{3}{2}}}{\partial\text{x}%
_{i}^{2}\partial\text{x}_{j}}\frac{\partial\theta}{\partial\text{x}_{j}%
}+2\frac{\partial^{2}\theta^{-\frac{3}{2}}}{\partial\text{x}_{i}%
\partial\text{x}_{j}}\frac{\partial^{2}\theta}{\partial\text{x}_{i}%
\partial\text{x}_{j}}+\frac{\partial\theta^{-\frac{3}{2}}}{\partial
\text{x}_{j}}\frac{\partial^{3}\theta}{\partial\text{x}_{i}^{2}\partial
\text{x}_{j}}+\frac{\partial^{2}\theta^{-\frac{3}{2}}}{\partial\text{x}%
_{i}^{2}}\frac{\partial^{2}\theta}{\partial\text{x}_{j}^{2}}$

$+2\frac{\partial\theta^{-\frac{3}{2}}}{\partial\text{x}_{i}}\frac
{\partial^{3}\theta}{\partial\text{x}_{i}\partial\text{x}_{j}^{2}}%
+\theta^{-\frac{3}{2}}\frac{\partial^{4}\theta}{\partial\text{x}_{i}%
^{2}\partial\text{x}_{j}^{2}}]$ \qquad\qquad\qquad\qquad\qquad

Next, we have,

$\ \ \ \frac{\partial\theta^{-\frac{3}{2}}}{\partial\text{x}_{i}}%
(y_{0})=-\frac{3}{2}\frac{\partial\theta}{\partial\text{x}_{i}}(y_{0});$
$\frac{\partial^{2}\theta^{-\frac{3}{2}}}{\partial\text{x}_{i}\partial
\text{x}_{j}}(y_{0})$

$=\frac{15}{4}\frac{\partial\theta}{\partial\text{x}_{i}}(y_{0}).\frac
{\partial\theta}{\partial\text{x}_{j}}(y_{0})-\frac{3}{2}\frac{\partial
^{2}\theta}{\partial\text{x}_{i}\partial\text{x}_{j}}(y_{0})$\qquad
\qquad\qquad

\vspace{1pt} $\ \frac{\partial^{3}\theta^{-\frac{3}{2}}}{\partial\text{x}%
_{i}^{2}\partial\text{x}_{j}}(y_{0})=-\frac{105}{8}\frac{\partial\theta
}{\partial\text{x}_{j}}(y_{0})(\frac{\partial\theta}{\partial\text{x}_{i}%
})^{2}(y_{0})+\frac{15}{4}\frac{\partial\theta}{\partial\text{x}_{j}}%
(y_{0})\frac{\partial^{2}\theta}{\partial\text{x}_{i}^{2}}(y_{0})$

$\ +\frac{15}{2}\frac{\partial\theta}{\partial\text{x}_{i}}(y_{0}%
)\frac{\partial^{2}\theta}{\partial\text{x}_{i}\partial\text{x}_{j}}%
(y_{0})-\frac{3}{2}\frac{\partial^{3}\theta}{\partial\text{x}_{i}^{2}%
\partial\text{x}_{j}}(y_{0})]$

We conclude that:

$\ \ \frac{\partial^{4}\theta^{-\frac{1}{2}}}{\partial\text{x}_{i}^{2}%
\partial\text{x}_{j}^{2}}(y_{0})=\frac{105}{16}(\frac{\partial\theta}%
{\partial\text{x}_{i}})^{2}(y_{0})(\frac{\partial\theta}{\partial\text{x}_{j}%
})^{2}(y_{0})-\frac{15}{8}(\frac{\partial\theta}{\partial\text{x}_{j}}%
)^{2}(y_{0})\frac{\partial^{2}\theta}{\partial\text{x}_{i}^{2}}(y_{0})$

$-\frac{15}{8}(\frac{\partial\theta}{\partial\text{x}_{i}})^{2}(y_{0}%
)\frac{\partial^{2}\theta}{\partial\text{x}_{j}^{2}}(y_{0})-\frac{15}{2}%
\frac{\partial\theta}{\partial\text{x}_{i}}(y_{0})\frac{\partial\theta
}{\partial\text{x}_{j}}(y_{0})\frac{\partial^{2}\theta}{\partial\text{x}%
_{i}\partial\text{x}_{j}}(y_{0})+\frac{3}{2}\frac{\partial\theta}%
{\partial\text{x}_{j}}\frac{\partial^{3}\theta}{\partial\text{x}_{i}%
^{2}\partial\text{x}_{j}}(y_{0})$

$+\frac{3}{2}\frac{\partial\theta}{\partial\text{x}_{i}}(y_{0})\frac
{\partial^{3}\theta}{\partial\text{x}_{i}\partial\text{x}_{j}^{2}}%
(y_{0})+\frac{3}{2}(\frac{\partial^{2}\theta}{\partial\text{x}_{i}%
\partial\text{x}_{j}})^{2}(y_{0})+\frac{3}{4}\frac{\partial^{2}\theta
}{\partial\text{x}_{i}^{2}}\frac{\partial^{2}\theta}{\partial\text{x}_{j}^{2}%
}-\frac{1}{2}\ \frac{\partial^{4}\theta}{\partial\text{x}_{i}^{2}%
\partial\text{x}_{j}^{2}}(y_{0})$

Consequently, we have:

$\left(  A_{22}\right)  \qquad\ \frac{\partial^{4}\theta^{-\frac{1}{2}}%
}{\partial\text{x}_{i}^{2}\partial\text{x}_{j}^{2}}(y_{0})=\frac{105}%
{16}<H,i>^{2}(y_{0})<H,j>^{2}(y_{0})$

$+\frac{15}{24}<H,j>^{2}(y_{0})[\varrho_{ii}+2\overset{q}{\underset{\text{a=1}%
}{\sum}}R_{i\text{a}i\text{a}}-3\overset{q}{\underset{\text{a,b=1}}{\sum}%
}(T_{\text{aa}i}T_{\text{bb}i}-T_{\text{ab}i}T_{\text{ab}i})](y_{0})$

$+\frac{15}{24}<H,i>^{2}(y_{0})[\varrho_{jj}+2\overset{q}{\underset{\text{a=1}%
}{\sum}}R_{j\text{a}j\text{a}}-3\overset{q}{\underset{\text{a,b=1}}{\sum}%
}(T_{\text{aa}j}T_{\text{bb}j}-T_{\text{ab}j}T_{\text{ab}j})](y_{0})$

$+\frac{15}{12}[<H,i><H,j>](y_{0})$

$\times\lbrack2\varrho_{ij}+4\overset{q}{\underset{\text{a}=1}{\sum}%
}R_{i\text{a}j\text{a}}-3\overset{q}{\underset{\text{a,b=1}}{\sum}%
}(T_{\text{aa}i}T_{\text{bb}j}-T_{\text{ab}i}T_{\text{ab}j}%
)-3\overset{q}{\underset{\text{a,b=1}}{\sum}}(T_{\text{aa}j}T_{\text{bb}%
i}-T_{\text{ab}j}T_{\text{ab}i})](y_{0})$

$+\frac{1}{4}<H,j>(y_{0})[\{\nabla_{i}\varrho_{ij}-2\varrho_{ij}%
<H,i>+\overset{q}{\underset{\text{a}=1}{\sum}}(\nabla_{i}R_{\text{a}%
i\text{a}j}-4R_{i\text{a}j\text{a}}<H,i>)$

$+4\overset{q}{\underset{\text{a,b=1}}{\sum}}R_{i\text{a}j\text{b}%
}T_{\text{ab}i}+2\overset{q}{\underset{\text{a,b,c=1}}{\sum}}(T_{\text{aa}%
i}T_{\text{bb}j}T_{\text{cc}i}-3T_{\text{aa}i}T_{\text{bc}j}T_{\text{bc}%
i}+2T_{\text{ab}i}T_{\text{bc}j}T_{\text{ac}i})\}$\qquad\qquad\qquad
\qquad\qquad\ \ 

$+\{\nabla_{j}\varrho_{ii}-2\varrho_{ij}<H,i>+\overset{q}{\underset{\text{a}%
=1}{\sum}}(\nabla_{j}R_{\text{a}i\text{a}i}-4R_{i\text{a}j\text{a}%
}<H,i>)+4\overset{q}{\underset{\text{a,b=1}}{\sum}}R_{j\text{a}i\text{b}%
}T_{\text{ab}i}$

$+2\overset{q}{\underset{\text{a,b,c=1}}{\sum}}(T_{\text{aa}j}T_{\text{bb}%
i}T_{\text{cc}i}-3T_{\text{aa}j}T_{\text{bc}i}T_{\text{bc}i}+2T_{\text{ab}%
j}T_{\text{bc}i}T_{\text{ac}i})\}$

$+\{\nabla_{i}\varrho_{ij}-2\varrho_{ii}<H,j>+\overset{q}{\underset{\text{a}%
=1}{\sum}}(\nabla_{i}R_{\text{a}i\text{a}j}-4R_{i\text{a}i\text{a}%
}<H,j>)+4\overset{q}{\underset{\text{a,b=1}}{\sum}}R_{i\text{a}i\text{b}%
}T_{\text{ab}j}$

$+2\overset{q}{\underset{\text{a,b,c}=1}{\sum}}(T_{\text{aa}i}T_{\text{bb}%
i}T_{\text{cc}j}-3T_{\text{aa}i}T_{\text{bc}i}T_{\text{bc}j}+2T_{\text{ab}%
i}T_{\text{bc}i}T_{\text{ac}j})\}](y_{0})$

$+\frac{1}{4}<H,i>(y_{0})[\{\nabla_{i}\varrho_{jj}-2\varrho_{ij}%
<H,j>+\overset{q}{\underset{\text{a}=1}{\sum}}(\nabla_{i}R_{\text{a}%
j\text{a}j}-4R_{i\text{a}j\text{a}}<H,j>)$

$+4\overset{q}{\underset{\text{a,b=1}}{\sum}}R_{i\text{a}j\text{b}%
}T_{\text{ab}j}+2\overset{q}{\underset{\text{a,b,c=1}}{\sum}}(T_{\text{aa}%
i}T_{\text{bb}j}T_{\text{cc}j}-3T_{\text{aa}i}T_{\text{bc}j}T_{\text{bc}%
j}+2T_{\text{ab}i}T_{\text{bc}j}T_{\text{ca}j})\}$\qquad\qquad\qquad
\qquad\qquad\ \ 

$+\{\nabla_{j}\varrho_{ij}-2\varrho_{ij}<H,j>+\overset{q}{\underset{\text{a}%
=1}{\sum}}(\nabla_{j}R_{\text{a}i\text{a}j}-4R_{j\text{a}i\text{a}}<H,j>)$

$+4\overset{q}{\underset{\text{a,b=1}}{\sum}}R_{j\text{a}i\text{b}%
}T_{\text{ab}j}+2\overset{q}{\underset{\text{a,b,c=1}}{\sum}}(T_{\text{aa}%
j}T_{\text{bb}i}T_{\text{cc}j}-3T_{\text{aa}j}T_{\text{bc}i}T_{\text{bc}%
j}+2T_{\text{ab}j}T_{\text{bc}i}T_{\text{ac}j})\}$

$+\{\nabla_{j}\varrho_{ij}-2\varrho_{jj}<H,i>+\overset{q}{\underset{\text{a}%
=1}{\sum}}(\nabla_{j}R_{\text{a}i\text{a}j}-4R_{j\text{a}j\text{a}%
}<H,i>)+4\overset{q}{\underset{\text{a,b=1}}{\sum}}R_{j\text{a}j\text{b}%
}T_{\text{ab}i}$

$+2\overset{q}{\underset{\text{a,b,c=1}}{\sum}}(T_{\text{aa}j}T_{\text{bb}%
j}T_{\text{cc}i}-3T_{\text{aa}j}T_{\text{bc}j}T_{\text{bc}i}+2T_{\text{ab}%
j}T_{\text{bc}j}T_{\text{ac}i})\}](y_{0})$

$+\frac{3}{2}\times\frac{1}{36}[2\varrho_{ij}+4\overset{q}{\underset{\text{a}%
=1}{\sum}}R_{i\text{a}j\text{a}}-3\overset{q}{\underset{\text{a,b=1}}{\sum}%
}(T_{\text{aa}i}T_{\text{bb}j}-T_{\text{ab}i}T_{\text{ab}j})$

$-3\overset{q}{\underset{\text{a,b=1}}{\sum}}(T_{\text{aa}j}T_{\text{bb}%
i}-T_{\text{ab}j}T_{\text{ab}i})]^{2}(y_{0})$

$+\frac{1}{12}[\varrho_{ii}+2\overset{q}{\underset{\text{a=1}}{\sum}%
}R_{i\text{a}i\text{a}}-3\overset{q}{\underset{\text{a,b=1}}{\sum}%
}(T_{\text{aa}i}T_{\text{bb}i}-T_{\text{ab}i}T_{\text{ab}i})](y_{0})$

$\times\lbrack\varrho_{jj}+2\overset{q}{\underset{\text{a=1}}{\sum}%
}R_{j\text{a}j\text{a}}-3\overset{q}{\underset{\text{a,b=1}}{\sum}%
}(T_{\text{aa}j}T_{\text{bb}j}-T_{\text{ab}j}T_{\text{ab}j})](y_{0})-\frac
{1}{2}\ \frac{\partial^{4}\theta}{\partial\text{x}_{i}^{2}\partial\text{x}%
_{j}^{2}}(y_{0})$

where the value of $\frac{\partial^{4}\theta}{\partial\text{x}_{i}^{2}%
\partial\text{x}_{j}^{2}}(y_{0})$ is given by (xviii).

\qquad\qquad\qquad\qquad\qquad\qquad\qquad\qquad\qquad\qquad\qquad\qquad
\qquad\qquad\qquad\qquad\qquad$\blacksquare$

(xx)\qquad We express I$_{32124}=\frac{1}{24}\frac{\partial^{4}\theta
^{-\frac{1}{2}}}{\partial x_{i}^{2}\partial x_{j}^{2}}(y_{0})$ in
\textbf{geometric invariants}:

Here we will use (vi) to have:

$\qquad\frac{\partial^{2}\theta}{\partial\text{x}_{i}^{2}}(y_{0})=-\frac{1}%
{3}[\tau^{M}-3\tau^{P}+\overset{q}{\underset{\text{a}=1}{\sum}}\varrho
_{\text{aa}}^{M}+\overset{q}{\underset{\text{a,b}=1}{\sum}}R_{\text{abab}}%
^{M}](y_{0}).$

We have:

$\left(  A_{23}\right)  \qquad$ I$_{32124}=\ \frac{1}{24}\frac{\partial
^{4}\theta^{-\frac{1}{2}}}{\partial x_{i}^{2}\partial x_{j}^{2}}%
(y_{0})=\overset{n}{\underset{i,j=q+1}{\sum}}\frac{1}{24}\times\frac{105}%
{16}<H,i>^{2}(y_{0})<H,j>^{2}(y_{0})\qquad\qquad$

$+\frac{1}{24}\times\frac{15}{24}\overset{n}{\underset{j=q+1}{\sum}}%
<H,j>^{2}(y_{0})[\tau^{M}\ -3\tau^{P}+\ \underset{\text{a}%
=1}{\overset{\text{q}}{\sum}}\varrho_{\text{aa}}^{M}%
+\overset{q}{\underset{\text{a},\text{b}=1}{\sum}}R_{\text{abab}}^{M}%
](y_{0})\qquad\ \ \ \ \ \ \ \ $

$+\frac{1}{24}\times\frac{15}{24}\overset{n}{\underset{i=q+1}{\sum}}%
<H,i>^{2}(y_{0})[\tau^{M}\ -3\tau^{P}+\ \underset{\text{a}%
=1}{\overset{\text{q}}{\sum}}\varrho_{\text{aa}}^{M}%
+\overset{q}{\underset{\text{a},\text{b}=1}{\sum}}R_{\text{abab}}^{M}%
](y_{0})\qquad\qquad$

$+\frac{1}{24}\times\frac{15}{12}\overset{n}{\underset{i,j=q+1}{\sum}%
}[<H,i><H,j>](y_{0})\qquad\qquad\qquad\qquad\qquad\qquad\qquad\qquad$

$\times\lbrack2\varrho_{ij}+4\overset{q}{\underset{\text{a}=1}{\sum}%
}R_{i\text{a}j\text{a}}-3\overset{q}{\underset{\text{a,b=1}}{\sum}%
}(T_{\text{aa}i}T_{\text{bb}j}-T_{\text{ab}i}T_{\text{ab}j}%
)-3\overset{q}{\underset{\text{a,b=1}}{\sum}}(T_{\text{aa}j}T_{\text{bb}%
i}-T_{\text{ab}j}T_{\text{ab}i})](y_{0})$

$+\frac{1}{24}\times\frac{1}{4}\overset{n}{\underset{i,j=q+1}{\sum}%
}<H,j>(y_{0})[\{\nabla_{i}\varrho_{ij}-2\varrho_{ij}%
<H,i>+\overset{q}{\underset{\text{a}=1}{\sum}}(\nabla_{i}R_{\text{a}%
i\text{a}j}-4R_{i\text{a}j\text{a}}<H,i>)\qquad$

$+4\overset{q}{\underset{\text{a,b=1}}{\sum}}R_{i\text{a}j\text{b}%
}T_{\text{ab}i}+2\overset{q}{\underset{\text{a,b,c=1}}{\sum}}(T_{\text{aa}%
i}T_{\text{bb}j}T_{\text{cc}i}-T_{\text{aa}i}T_{\text{bc}j}T_{\text{bc}%
i}-2T_{\text{bc}j}(T_{\text{aa}i}T_{\text{bc}i}-T_{\text{ab}i}T_{\text{ac}%
i}))\}$\qquad\qquad\qquad\ \ 

$+\{\nabla_{j}\varrho_{ii}-2\varrho_{ij}<H,i>+\overset{q}{\underset{\text{a}%
=1}{\sum}}(\nabla_{j}R_{\text{a}i\text{a}i}-4R_{i\text{a}j\text{a}}<H,i>)$

$+4\overset{q}{\underset{\text{a,b=1}}{\sum}}R_{j\text{a}i\text{b}%
}T_{\text{ab}i}+2\overset{q}{\underset{\text{a,b,c=1}}{\sum}}(T_{\text{aa}%
j}(T_{\text{bb}i}T_{\text{cc}i}-T_{\text{bc}i}T_{\text{bc}i})-2T_{\text{aa}%
j}T_{\text{bc}i}T_{\text{bc}i}+2T_{\text{ab}j}T_{\text{bc}i}T_{\text{ac}%
i})\}\qquad$

$+\{\nabla_{i}\varrho_{ij}-2\varrho_{ii}<H,j>+\overset{q}{\underset{\text{a}%
=1}{\sum}}(\nabla_{i}R_{\text{a}i\text{a}j}-4R_{i\text{a}i\text{a}%
}<H,j>)+4\overset{q}{\underset{\text{a,b=1}}{\sum}}R_{i\text{a}i\text{b}%
}T_{\text{ab}j}$

$+2\overset{q}{\underset{\text{a,b,c}=1}{\sum}}(T_{\text{aa}i}T_{\text{bb}%
i}T_{\text{cc}j}-3T_{\text{aa}i}T_{\text{bc}i}T_{\text{bc}j}+2T_{\text{ab}%
i}T_{\text{bc}i}T_{\text{ac}j})\}](y_{0})$

$+\frac{1}{24}\times\frac{1}{4}\overset{n}{\underset{i,j=q+1}{\sum}%
}<H,i>(y_{0})[\{\nabla_{i}\varrho_{jj}-2\varrho_{ij}%
<H,j>+\overset{q}{\underset{\text{a}=1}{\sum}}(\nabla_{i}R_{\text{a}%
j\text{a}j}-4R_{i\text{a}j\text{a}}<H,j>)\qquad$

$\qquad+4\overset{q}{\underset{\text{a,b=1}}{\sum}}R_{i\text{a}j\text{b}%
}T_{\text{ab}j}+2\overset{q}{\underset{\text{a,b,c=1}}{\sum}}T_{\text{aa}%
i}(T_{\text{bb}j}T_{\text{cc}j}-T_{\text{bc}j}T_{\text{bc}j})-2T_{\text{aa}%
i}T_{\text{bc}j}T_{\text{bc}j}+2T_{\text{ab}i}T_{\text{bc}j}T_{\text{ac}%
j})\}\qquad$\qquad\qquad\qquad\qquad\qquad\ \ 

$+\{\nabla_{j}\varrho_{ij}-2\varrho_{ij}<H,j>+\overset{q}{\underset{\text{a}%
=1}{\sum}}(\nabla_{j}R_{\text{a}i\text{a}j}-4R_{j\text{a}i\text{a}}<H,j>)$

$+4\overset{q}{\underset{\text{a,b=1}}{\sum}}R_{j\text{a}i\text{b}%
}T_{\text{ab}j}+2\overset{q}{\underset{\text{a,b,c=1}}{\sum}}(T_{\text{aa}%
j}T_{\text{bb}i}T_{\text{cc}j}-T_{\text{ab}j}T_{\text{bc}i}T_{\text{ac}%
j}-2T_{\text{bc}i}(T_{\text{aa}j}T_{\text{bc}j}-T_{\text{ab}j}T_{\text{ac}%
j})\}$

$+\{\nabla_{j}\varrho_{ij}-2\varrho_{jj}<H,i>+\overset{q}{\underset{\text{a}%
=1}{\sum}}(\nabla_{j}R_{\text{a}i\text{a}j}-4R_{j\text{a}j\text{a}%
}<H,i>)+4\overset{q}{\underset{\text{a,b=1}}{\sum}}R_{j\text{a}j\text{b}%
}T_{\text{ab}i}$

$+2\overset{q}{\underset{\text{a,b,c=1}}{\sum}}(T_{\text{aa}j}T_{\text{bb}%
j}T_{\text{cc}i}-3T_{\text{aa}j}T_{\text{bc}j}T_{\text{bc}i}+2T_{\text{ab}%
j}T_{\text{bc}j}T_{\text{ac}i})\}](y_{0})$

$+\frac{1}{24}\times\frac{3}{2}\times\frac{1}{36}%
\overset{n}{\underset{i,j=q+1}{\sum}}[2\varrho_{ij}%
+4\overset{q}{\underset{\text{a}=1}{\sum}}R_{i\text{a}j\text{a}}%
-3\overset{q}{\underset{\text{a,b=1}}{\sum}}(T_{\text{aa}i}T_{\text{bb}%
j}-T_{\text{ab}i}T_{\text{ab}j})$

$-3\overset{q}{\underset{\text{a,b=1}}{\sum}}(T_{\text{aa}j}T_{\text{bb}%
i}-T_{\text{ab}j}T_{\text{ab}i})]^{2}(y_{0})$

$+\frac{1}{24}\times\frac{1}{12}[\tau^{M}\ -3\tau^{P}+\ \underset{\text{a}%
=1}{\overset{\text{q}}{\sum}}\varrho_{\text{aa}}^{M}%
+\overset{q}{\underset{\text{a},\text{b}=1}{\sum}}R_{\text{abab}}^{M}%
]^{2}(y_{0})$

$-\frac{1}{24}\times\frac{1}{2}$

$\times\frac{1}{6}\overset{n}{\underset{i,j=q+1}{\sum}}[$
$\overset{q}{\underset{\text{a=1}}{\sum}}\{-(\nabla_{ii}^{2}R_{j\text{a}%
j\text{a}}+\nabla_{jj}^{2}R_{i\text{a}i\text{a}}+4\nabla_{ij}^{2}%
R_{i\text{a}j\text{a}}+2R_{ij}R_{i\text{a}j\text{a}})\qquad-\frac{1}{2}%
\ \frac{\partial^{4}\theta}{\partial\text{x}_{i}^{2}\partial\text{x}_{j}^{2}%
}(y_{0})\qquad A$

$+\overset{n}{\underset{p=q+1}{\sum}}\overset{q}{\underset{\text{a=1}}{\sum}%
}(R_{\text{a}iip}R_{\text{a}jjp}+R_{\text{a}jjp}R_{\text{a}iip}+R_{\text{a}%
ijp}R_{\text{a}ijp}+R_{\text{a}ijp}R_{\text{a}jip}+R_{\text{a}jip}%
R_{\text{a}ijp}+R_{\text{a}jip}R_{\text{a}jip})$

$+2\overset{q}{\underset{\text{a,b=1}}{\sum}}\nabla_{i}(R)_{\text{a}%
i\text{b}j}T_{\text{ab}j}+2\overset{q}{\underset{\text{a,b=1}}{\sum}}%
\nabla_{j}(R)_{\text{a}j\text{b}i}T_{\text{ab}i}%
+2\overset{q}{\underset{\text{a,b=1}}{\sum}}\nabla_{i}(R)_{\text{a}j\text{b}%
i}T_{\text{ab}j}$

$+2\overset{q}{\underset{\text{a,b=1}}{\sum}}\nabla_{i}(R)_{\text{a}%
j\text{b}j}T_{\text{ab}i}$

$+2\overset{q}{\underset{\text{a,b=1}}{\sum}}\nabla_{j}(R)_{\text{a}%
i\text{b}i}T_{\text{ab}j}+2\overset{q}{\underset{\text{a,b=1}}{\sum}}%
\nabla_{j}(R)_{\text{a}i\text{b}j}T_{\text{ab}i}$

$+\overset{n}{\underset{p=q+1}{\sum}}(-\frac{3}{5}\nabla_{ii}^{2}%
(R)_{jpjp}+\overset{n}{\underset{p=q+1}{\sum}}(-\frac{3}{5}\nabla_{jj}%
^{2}(R)_{ipip}$

$+\overset{n}{\underset{p=q+1}{\sum}}(-\frac{3}{5}\nabla_{ij}^{2}%
(R)_{ipjp}+\overset{n}{\underset{p=q+1}{\sum}}(-\frac{3}{5}\nabla_{ij}%
^{2}(R)_{jpip}$

$+\overset{n}{\underset{p=q+1}{\sum}}(-\frac{3}{5}\nabla_{ji}^{2}%
(R)_{ipjp}+\overset{n}{\underset{p=q+1}{\sum}}(-\frac{3}{5}\nabla_{ji}%
^{2}(R)_{jpip}$

$+\frac{1}{5}\overset{n}{\underset{m,p=q+1}{%
{\textstyle\sum}
}}R_{ipim}R_{jpjm}+\frac{1}{5}\overset{n}{\underset{m,p=q+1}{%
{\textstyle\sum}
}}R_{jpjm}R_{ipim}$

$+\frac{1}{5}\overset{n}{\underset{m,p=q+1}{%
{\textstyle\sum}
}}R_{ipjm}R_{ipjm}+\frac{1}{5}\overset{n}{\underset{m,p=q+1}{%
{\textstyle\sum}
}}R_{ipjm}R_{jpim}$

$+\frac{1}{5}\overset{n}{\underset{m,p=q+1}{%
{\textstyle\sum}
}}R_{jpim}R_{ipjm}+\frac{1}{5}\overset{n}{\underset{m,p=q+1}{%
{\textstyle\sum}
}}R_{jpim}R_{jpim}\}(y_{0})$

$+4\overset{q}{\underset{\text{a,b=1}}{\sum}}\{(\nabla_{i}(R)_{i\text{a}%
j\text{a}}-\overset{q}{\underset{\text{c=1}}{%
{\textstyle\sum}
}}R_{\text{a}i\text{c}i}T_{\text{ac}j})$ $T_{\text{bb}j}+4(\nabla
_{j}(R)_{j\text{a}i\text{a}}-\overset{q}{\underset{\text{c=1}}{%
{\textstyle\sum}
}}R_{\text{a}j\text{c}j}T_{\text{ac}i})$ $T_{\text{bb}i}$

$+4(\nabla_{i}(R)_{j\text{a}i\text{a}}-\overset{q}{\underset{\text{c=1}}{%
{\textstyle\sum}
}}R_{\text{a}i\text{c}j}T_{\text{ac}i})$ $T_{\text{bb}j}$ $4B\ +4(\nabla
_{i}(R)_{j\text{a}j\text{a}}-\overset{q}{\underset{\text{c=1}}{%
{\textstyle\sum}
}}R_{\text{a}i\text{c}j}T_{\text{ac}j})$ $T_{\text{bb}i}$

$+4(\nabla_{j}(R)_{i\text{a}i\text{a}}-\overset{q}{\underset{\text{c=1}}{%
{\textstyle\sum}
}}R_{\text{a}j\text{c}i}T_{\text{ac}i})$ $T_{\text{bb}j}$

$+4(\nabla_{j}(R)_{i\text{a}j\text{a}}-\overset{q}{\underset{\text{c=1}}{%
{\textstyle\sum}
}}R_{\text{a}j\text{c}i}T_{\text{ac}j})$ $T_{\text{bb}i}%
-4\overset{q}{\underset{\text{a,b=1}}{\sum}}(\nabla_{i}(R)_{i\text{a}%
j\text{b}}-\overset{q}{\underset{\text{c=1}}{%
{\textstyle\sum}
}}R_{\text{b}r\text{c}s}T_{\text{ac}t})T_{\text{ab}j}$

$-4\overset{q}{\underset{\text{a,b=1}}{\sum}}(\nabla_{j}(R)_{j\text{a}%
i\text{b}}-\overset{q}{\underset{\text{c=1}}{%
{\textstyle\sum}
}}R_{\text{b}j\text{c}j}T_{\text{ac}i})T_{\text{ab}i}$

$-4\overset{q}{\underset{\text{a,b=1}}{\sum}}(\nabla_{i}(R)_{j\text{a}%
i\text{b}}-\overset{q}{\underset{\text{c=1}}{%
{\textstyle\sum}
}}R_{\text{b}i\text{c}j}T_{\text{ac}i})T_{\text{ab}j}%
-4\overset{q}{\underset{\text{a,b=1}}{\sum}}(\nabla_{i}(R)_{j\text{a}%
j\text{b}}-\overset{q}{\underset{\text{c=1}}{%
{\textstyle\sum}
}}R_{\text{b}i\text{c}j}T_{\text{ac}j})T_{\text{ab}i}$

$-4\overset{q}{\underset{\text{a,b=1}}{\sum}}(\nabla_{j}(R)_{i\text{a}%
i\text{b}}-\overset{q}{\underset{\text{c=1}}{%
{\textstyle\sum}
}}R_{\text{b}j\text{c}i}T_{\text{ac}i})T_{\text{ab}j}%
-4\overset{q}{\underset{\text{a,b=1}}{\sum}}(\nabla_{j}(R)_{i\text{a}%
j\text{b}}-\overset{q}{\underset{\text{c=1}}{%
{\textstyle\sum}
}}R_{\text{b}j\text{c}i}T_{\text{ac}j})T_{\text{ab}i}\}](y_{0})$

$+$ $[\frac{4}{9}\overset{q}{\underset{\text{a,b=1}}{\sum}}(\varrho
_{\text{aa}}-\overset{q}{\underset{\text{c}=1}{\sum}}R_{\text{acac}}%
)(\varrho_{\text{bb}}-\overset{q}{\underset{\text{d}=1}{\sum}}R_{\text{bdbd}%
})+\frac{8}{9}\overset{n}{\underset{i,j=q+1}{\sum}}%
\overset{q}{\underset{\text{a,b}=1}{\sum}}(R_{i\text{a}j\text{a}}%
R_{i\text{b}j\text{b}})\qquad3C$

$+\frac{4}{3}\overset{q}{\underset{\text{a}=1}{\sum}}(\varrho_{\text{aa}}%
^{M}-\varrho_{\text{aa}}^{P})(\tau^{M}-\overset{q}{\underset{\text{c}=1}{\sum
}}\varrho_{\text{cc}}^{M})+\frac{4}{9}\overset{n}{\underset{i,j=q+1}{\sum}%
}\overset{q}{\underset{\text{a}=1}{\sum}}R_{i\text{a}j\text{a}}\varrho_{ij}\ $

$\ +\frac{2}{9}\overset{q}{\underset{\text{b}=1}{\sum}}(\varrho_{\text{bb}%
}^{M}-\varrho_{\text{bb}}^{P})(\tau^{M}-\overset{q}{\underset{\text{c}%
=1}{\sum}}\varrho_{\text{cc}}^{M})+\frac{4}{9}%
\overset{n}{\underset{i,j=q+1}{\sum}}\overset{q}{\underset{\text{b}=1}{\sum}%
}R_{i\text{b}j\text{b}}\varrho_{ij}\ $

$+\frac{1}{9}(\tau^{M}-\overset{q}{\underset{\text{a=1}}{\sum}}\varrho
_{\text{aa}})(\tau^{M}-\overset{q}{\underset{\text{b=1}}{\sum}}\varrho
_{\text{bb}})+\frac{2}{9}(\left\Vert \varrho^{M}\right\Vert ^{2}%
-\overset{q}{\underset{\text{a,b}=1}{\sum}}\varrho_{\text{ab}})$

$-\overset{n}{\underset{i,j=q+1}{\sum}}\overset{q}{\underset{\text{a,b}%
=1}{\sum}}R_{i\text{a}i\text{b}}R_{j\text{a}j\text{b}}\ -\frac{1}%
{2}\overset{n}{\underset{i,j=q+1}{\sum}}\overset{q}{\underset{\text{a,b}%
=1}{\sum}}R_{i\text{a}j\text{b}}^{2}-\overset{n}{\underset{i,j=q+1}{\sum}%
}\overset{q}{\underset{\text{a,b}=1}{\sum}}R_{i\text{a}j\text{b}}%
R_{j\text{a}i\text{b}}$

$-\frac{1}{2}\overset{n}{\underset{i,j=q+1}{\sum}}%
\overset{q}{\underset{\text{a,b}=1}{\sum}}R_{j\text{a}i\text{b}}^{2}$

$-\frac{1}{9}\overset{n}{\underset{i,j,p,m=q+1}{\sum}}R_{ipim}R_{jpjm}%
\ -\frac{1}{18}\overset{n}{\underset{i,j,p,m=q+1}{\sum}}R_{ipjm}^{2}-\frac
{1}{9}\overset{n}{\underset{i,j,p,m=q+1}{\sum}}R_{ipjm}R_{jpim}$

$-\frac{1}{18}\overset{n}{\underset{i,j,p,m=q+1}{\sum}}R_{jpim}^{2}$

$-\frac{1}{3}\overset{q}{\underset{\text{a}=1}{\sum}}%
\overset{n}{\underset{i,j,p=q+1}{\sum}}R_{i\text{a}ip}R_{j\text{a}jp}-\frac
{1}{6}\overset{q}{\underset{\text{a}=1}{\sum}}%
\overset{n}{\underset{i,j,p=q+1}{\sum}}R_{i\text{a}jp}^{2}-\frac{1}%
{3}\overset{q}{\underset{\text{a}=1i,j,}{\sum}}%
\overset{n}{\underset{p=q+1}{\sum}}R_{i\text{a}jp}R_{j\text{a}ip}$

$-\frac{1}{6}\overset{q}{\underset{\text{a}=1}{\sum}}%
\overset{n}{\underset{i,j,p=q+1}{\sum}}R_{j\text{a}ip}^{2}$

$-\frac{1}{3}\overset{q}{\underset{\text{b}=1i,j,}{\sum}}%
\overset{n}{\underset{p=q+1}{\sum}}R_{i\text{b}ip}R_{j\text{b}jp}-\frac{1}%
{6}\overset{q}{\underset{\text{b}=1}{\sum}}%
\overset{n}{\underset{i,j,p=q+1}{\sum}}R_{i\text{b}jp}^{2}-\frac{1}%
{3}\overset{q}{\underset{\text{b}=1}{\sum}}%
\overset{n}{\underset{i.j,p=q+1}{\sum}}R_{i\text{b}jp}R_{j\text{b}ip}$

$-\frac{1}{6}\overset{q}{\underset{\text{b}=1}{\sum}}%
\overset{n}{\underset{i,j,p=q+1}{\sum}}R_{j\text{b}ip}^{2}](y_{0})$

$+$ $\overset{q}{\underset{\text{a,b,c=1}}{\sum}}[$
$-\overset{n}{\underset{i=q+1}{\sum}}R_{i\text{a}i\text{a}}(R_{\text{bcbc}%
}^{P}-R_{\text{bcbc}}^{M})$ $-\overset{n}{\underset{j=q+1}{\sum}}%
R_{j\text{a}j\text{a}}(R_{\text{bcbc}}^{P}-R_{\text{bcbc}}^{M})\qquad\qquad6D$

\ $+\overset{n}{\underset{i=q+1}{\sum}}R_{i\text{a}i\text{b}}(R_{\text{acbc}%
}^{P}-R_{\text{acbc}}^{M})\ -\overset{n}{\underset{i=q+1}{\sum}}%
R_{i\text{a}i\text{c}}(R_{\text{abbc}}^{P}-R_{\text{abbc}}^{M})$

$+\overset{n}{\underset{j=q+1}{\sum}}R_{j\text{a}j\text{b}}(R_{\text{acbc}%
}^{P}-R_{\text{acbc}}^{M})$\ $-\overset{n}{\underset{j=q+1}{\sum}}%
R_{j\text{a}j\text{c}}(R_{\text{abbc}}^{P}-R_{\text{abbc}}^{M})$

$+\underset{i,j=q+1}{\overset{n}{\sum}}$ $-R_{i\text{a}j\text{a}}%
(T_{\text{bb}i}T_{\text{cc}j}$ $-T_{\text{bc}i}T_{\text{bc}j})$
$-\underset{i,j=q+1}{\overset{n}{\sum}}R_{i\text{a}j\text{a}}(T_{\text{bb}%
j}T_{\text{cc}i}$ $-T_{\text{bc}j}T_{\text{bc}i})$

$+$ $\underset{i,j=q+1}{\overset{n}{\sum}}$ $-R_{j\text{a}i\text{a}%
}(T_{\text{bb}i}T_{\text{cc}j}$ $-T_{\text{bc}i}T_{\text{bc}j})$
$-\underset{i,j=q+1}{\overset{n}{\sum}}R_{j\text{a}i\text{a}}(T_{\text{bb}%
j}T_{\text{cc}i}$ $-T_{\text{bc}j}T_{\text{bc}i})$

$+\underset{i,j=q+1}{\overset{n}{\sum}}\ R_{i\text{a}j\text{b}}(T_{\text{ab}%
i}T_{\text{cc}j}-T_{\text{bc}i}T_{\text{ac}j}%
)\ +\underset{i,j=q+1}{\overset{n}{\sum}}\ R_{i\text{a}j\text{b}}%
(T_{\text{ab}j}T_{\text{cc}i}-T_{\text{bc}j}T_{\text{ac}i})$

$+\underset{i,j=q+1}{\overset{n}{\sum}}\ R_{j\text{a}i\text{ib}}%
(T_{\text{ab}i}T_{\text{cc}j}-T_{\text{bc}i}T_{\text{ac}j}%
)\ +\underset{i,j=q+1}{\overset{n}{\sum}}\ R_{j\text{a}i\text{b}}%
(T_{\text{ab}j}T_{\text{cc}i}-T_{\text{bc}j}T_{\text{ac}i})\qquad$

$+\underset{i,j=q+1}{\overset{n}{\sum}}-R_{i\text{a}j\text{c}}(T_{\text{ab}%
i}T_{\text{bc}j}-T_{\text{ac}i}T_{\text{bb}j}%
)-\underset{i,j=q+1}{\overset{n}{\sum}}R_{i\text{a}j\text{c}}(T_{\text{ba}%
j}T_{\text{bc}i}-T_{\text{ac}j}T_{\text{bb}i})$

$+\underset{i,j=q+1}{\overset{n}{\sum}}-R_{j\text{a}i\text{c}}(T_{\text{ba}%
i}T_{\text{bc}j}-T_{\text{ac}i}T_{\text{bb}j}%
)-\underset{i,j=q+1}{\overset{n}{\sum}}R_{j\text{a}i\text{c}}(T_{\text{ba}%
j}T_{\text{bc}i}-T_{\text{ac}j}T_{\text{bb}i})](y_{0})$

$-\frac{1}{3}\underset{p=q+1}{\overset{n}{\sum}}%
[\underset{i=q+1}{\overset{n}{\sum}}\overset{q}{\underset{\text{b,c=1}}{\sum}%
}R_{ipip}(R_{\text{bcbc}}^{P}-R_{\text{bcbc}}^{M}%
)+\underset{j=q+1}{\overset{n}{\sum}}$ $\overset{q}{\underset{\text{b,c=1}%
}{\sum}}R_{jpjp}(R_{\text{bcbc}}^{P}-R_{\text{bcbc}}^{M})](y_{0})$

$-\frac{2}{3}\underset{i,j,p=q+1}{\overset{n}{\sum}}%
\overset{q}{\underset{\text{b,c=1}}{\sum}}[R_{ipjp}(T_{\text{bb}i}%
T_{\text{cc}j}-T_{\text{bc}i}T_{\text{bc}j})+R_{ipjp}(T_{\text{bb}%
j}T_{\text{cc}i}-T_{\text{bc}j}T_{\text{bc}i})](y_{0})\qquad$

$+\ \frac{1}{6}\underset{i,j=q+1}{\overset{n}{\sum}}[T_{\text{aa}%
i}T_{\text{bb}j}(T_{\text{cc}i}T_{\text{dd}j}-T_{\text{cd}i}T_{\text{dc}%
j})+T_{\text{aa}i}T_{\text{bb}j}(T_{\text{cc}j}T_{\text{dd}i}-T_{\text{cd}%
j}T_{\text{dc}i})\qquad E$

$+T_{\text{aa}j}T_{\text{bb}i}(T_{\text{cc}i}T_{\text{dd}j}-T_{\text{cd}%
i}T_{\text{dc}j})+T_{\text{aa}j}T_{\text{bb}i}(T_{\text{cc}j}T_{\text{dd}%
i}-T_{\text{cd}j}T_{\text{dc}i})](y_{0})$

$-\frac{1}{6}\underset{i,j=q+1}{\overset{n}{\sum}}[T_{\text{aa}i}%
T_{\text{bc}j}(T_{\text{bc}i}T_{\text{dd}j}-T_{\text{bd}i}T_{\text{cd}%
j})+T_{\text{aa}i}T_{\text{bc}j}(T_{\text{bc}j}T_{\text{dd}i}-T_{\text{bd}%
j}T_{\text{cd}i})$

$+T_{\text{aa}j}T_{\text{bc}i}(T_{\text{bc}i}T_{\text{dd}j}-T_{\text{bd}%
i}T_{\text{cd}j})+T_{\text{aa}j}T_{\text{bc}i}(T_{\text{bc}j}T_{\text{dd}%
i}-T_{\text{bd}j}T_{\text{cd}i})](y_{0})$

$+\ \frac{1}{6}\underset{i,j=q+1}{\overset{n}{\sum}}[T_{\text{aa}%
i}T_{\text{bd}j}(T_{\text{bc}i}T_{\text{cd}j}-T_{\text{bd}i}T_{\text{cc}%
j})+T_{\text{aa}i}T_{\text{bd}j}(T_{\text{bc}j}T_{\text{cd}i}-T_{\text{bd}%
j}T_{\text{cc}i})$

$+T_{\text{aa}j}T_{\text{bd}i}(T_{\text{bc}i}T_{\text{cd}j}-T_{\text{bd}%
i}T_{\text{cc}j})+T_{\text{aa}j}T_{\text{bd}i}(T_{\text{bc}j}T_{\text{cd}%
i}-T_{\text{bd}j}T_{\text{cc}i})](y_{0})\qquad$

$-\ \frac{1}{6}\underset{i,j=q+1}{\overset{n}{\sum}}[T_{\text{ab}%
i}T_{\text{ab}j}(T_{\text{cc}i}T_{\text{dd}j}-T_{\text{cd}i}T_{\text{dc}%
j})+T_{\text{ab}i}T_{\text{ab}j}(T_{\text{cc}j}T_{\text{dd}i}-T_{\text{cd}%
j}T_{\text{dc}i})$

$+T_{\text{ab}j}T_{\text{ab}i}(T_{\text{cc}i}T_{\text{dd}j}-T_{\text{cd}%
i}T_{\text{dc}j})+T_{\text{ab}j}T_{\text{ab}i}(T_{\text{cc}j}T_{\text{dd}%
i}-T_{\text{cd}j}T_{\text{dc}i})](y_{0})$

$+\ \frac{1}{6}\underset{i,j=q+1}{\overset{n}{\sum}}[T_{\text{ab}%
i}T_{\text{bc}j}(T_{\text{ac}i}T_{\text{dd}j}-T_{\text{ad}i}T_{\text{cd}%
j})+T_{\text{ab}i}T_{\text{bc}j}(T_{\text{ac}j}T_{\text{dd}i}-T_{\text{ad}%
j}T_{\text{cd}i})$

$+T_{\text{ab}j}T_{\text{bc}i}(T_{\text{ac}i}T_{\text{dd}j}-T_{\text{ad}%
i}T_{\text{cd}j})+T_{\text{ab}j}T_{\text{bc}i}(T_{\text{ac}j}T_{\text{dd}%
i}-T_{\text{ad}j}T_{\text{cd}i})](y_{0})$

$-\ \frac{1}{6}\underset{i,j=q+1}{\overset{n}{\sum}}[T_{\text{ab}%
i}T_{\text{bd}j}(T_{\text{ac}i}T_{\text{cd}j}-T_{\text{ad}i}T_{\text{cc}%
j})+T_{\text{ab}i}T_{\text{bd}j}(T_{\text{ac}j}T_{\text{cd}i}-T_{\text{ad}%
j}T_{\text{cc}i})$

$+T_{\text{ab}i}T_{\text{bd}j}(T_{\text{ac}j}T_{\text{cd}i}-T_{\text{ad}%
j}T_{\text{cc}i})+T_{\text{ab}j}T_{\text{bd}i}(T_{\text{ac}j}T_{\text{cd}%
i}-T_{\text{ad}j}T_{\text{cc}i})](y_{0})$

$+\ \frac{1}{6}\underset{i,j=q+1}{\overset{n}{\sum}}[T_{\text{ac}%
i}T_{\text{ab}j}(T_{\text{bc}i}T_{\text{dd}j}-T_{\text{bd}i}T_{\text{dc}%
j})+T_{\text{ac}i}T_{\text{ab}j}(T_{\text{bc}j}T_{\text{dd}i}-T_{\text{bd}%
j}T_{\text{dc}i})$

$+T_{\text{ac}j}T_{\text{ab}i}(T_{\text{bc}i}T_{\text{dd}j}-T_{\text{bd}%
i}T_{\text{dc}j})+T_{\text{ac}j}T_{\text{ab}i}(T_{\text{bc}j}T_{\text{dd}%
i}-T_{\text{bd}j}T_{\text{dc}i})](y_{0})$

$-\ \frac{1}{6}\underset{i,j=q+1}{\overset{n}{\sum}}[T_{\text{ac}%
i}T_{\text{bb}j}(T_{\text{ac}i}T_{\text{dd}j}-T_{\text{ad}i}T_{\text{cd}%
j})+T_{\text{ac}i}T_{\text{bb}j}(T_{\text{ac}j}T_{\text{dd}i}-T_{\text{ad}%
j}T_{\text{cd}i})$

$+T_{\text{ac}j}T_{\text{bb}i}(T_{\text{ac}i}T_{\text{dd}j}-T_{\text{ad}%
i}T_{\text{cd}i})+T_{\text{ac}j}T_{\text{bb}i}(T_{\text{ac}j}T_{\text{dd}%
i}-T_{\text{ad}j}T_{\text{cd}i})](y_{0})$

$+\ \frac{1}{6}\underset{i,j=q+1}{\overset{n}{\sum}}[T_{\text{ac}%
i}T_{\text{bd}j}(T_{\text{ac}i}T_{\text{bd}j}-T_{\text{ad}i}T_{\text{bc}%
j})+T_{\text{ac}i}T_{\text{bd}j}(T_{\text{ac}j}T_{\text{bd}i}-T_{\text{ad}%
j}T_{\text{bc}i})$

$+T_{\text{ac}j}T_{\text{bd}i}(T_{\text{ac}i}T_{\text{bd}j}-T_{\text{ad}%
i}T_{\text{bc}j})+T_{\text{ac}j}T_{\text{bd}i}(T_{\text{ac}j}T_{\text{bd}%
i}-T_{\text{ad}j}T_{\text{bc}i})](y_{0})$

$-\frac{1}{6}\underset{i,j=q+1}{\overset{n}{\sum}}[T_{\text{ad}i}%
T_{\text{ab}j}(T_{\text{bc}i}T_{\text{cd}j}-T_{\text{bd}i}T_{\text{cc}%
j})+T_{\text{ad}i}T_{\text{ab}j}(T_{\text{bc}j}T_{\text{cd}i}-T_{\text{bd}%
j}T_{\text{cc}i})$

$+T_{\text{ad}j}T_{\text{ab}i}(T_{\text{bc}i}T_{\text{cd}j}-T_{\text{bd}%
i}T_{\text{cc}j})+T_{\text{ad}j}T_{\text{ab}i}(T_{\text{bc}j}T_{\text{cd}%
i}-T_{\text{bd}j}T_{\text{cc}i})](y_{0})$

$+\ \frac{1}{6}\underset{i,j=q+1}{\overset{n}{\sum}}[T_{\text{ad}%
i}T_{\text{bb}j}(T_{\text{ac}i}T_{\text{cd}j}-T_{\text{ad}i}T_{\text{cc}%
j})+T_{\text{ad}i}T_{\text{bb}j}(T_{\text{ac}j}T_{\text{cd}i}-T_{\text{ad}%
j}T_{\text{cc}i})$

$+T_{\text{ad}j}T_{\text{bb}i}(T_{\text{ac}i}T_{\text{cd}j}-T_{\text{ad}%
i}T_{\text{cc}j})+T_{\text{ad}j}T_{\text{bb}i}(T_{\text{ac}j}T_{\text{cd}%
i}-T_{\text{ad}j}T_{\text{cc}i})](y_{0})$

$-\ \frac{1}{6}\underset{i,j=q+1}{\overset{n}{\sum}}[T_{\text{ad}%
i}T_{\text{bc}j}(T_{\text{ac}i}T_{\text{bd}j}-T_{\text{ad}i}T_{\text{bc}%
j})+T_{\text{ad}i}T_{\text{bc}j}(T_{\text{ac}j}T_{\text{bd}i}-T_{\text{ad}%
j}T_{\text{bc}i})$

$+T_{\text{ad}j}T_{\text{bc}i}(T_{\text{ac}i}T_{\text{bd}j}-T_{\text{ad}%
i}T_{\text{bc}j})+T_{\text{ad}j}T_{\text{bc}i}(T_{\text{ac}j}T_{\text{bd}%
i}-T_{\text{ad}j}T_{\text{bc}i})](y_{0})$

$+\ \frac{1}{3}[(R_{\text{cdcd}}^{P}-R_{\text{cdcd}}^{M})(R_{\text{abab}}%
^{P}-R_{\text{abab}}^{M})](y_{0})\qquad\left(  1\right)  $

$-\frac{1}{3}[(R_{\text{bdcd}}^{P}-R_{\text{bdcd}}^{M})(R_{\text{abac}}%
^{P}-R_{\text{abac}}^{M})](y_{0})\qquad(2)$

$-\ \frac{1}{3}[(R_{\text{bcdc}}^{P}-R_{\text{bcdc}}^{M})(R_{\text{abad}}%
^{P}-R_{\text{abad}}^{M})](y_{0})\qquad(3)$

$+\ \frac{1}{3}[(R_{\text{adcd}}^{P}-R_{\text{adcd}}^{M})(R_{\text{abbc}}%
^{P}-R_{\text{abbc}}^{M})](y_{0})\qquad(4)\qquad$

$-\ \frac{1}{3}[(R_{\text{acdc}}^{P}-R_{\text{acdc}}^{M})(R_{\text{abdb}}%
^{P}-R_{\text{abdb}}^{M})](y_{0})\qquad(5)$

$+\ \frac{1}{3}[(R_{\text{abcd}}^{P}-R_{\text{abcd}}^{M})]^{2}(y_{0}%
)\qquad(6)$

\qquad\qquad\qquad\qquad\qquad\qquad\qquad\qquad\qquad\qquad\qquad\qquad
\qquad\qquad\qquad\qquad\qquad$\blacksquare$

We simplify all fractions and give a final expression for:

$\left(  A_{24}\right)  \qquad$I$_{32124}=\ \frac{1}{24}%
\overset{n}{\underset{i,j=q+1}{\sum}}\frac{\partial^{4}\theta^{-\frac{1}{2}}%
}{\partial x_{i}^{2}\partial x_{j}^{2}}(y_{0})$

$=\overset{n}{\underset{i,j=q+1}{\sum}}\frac{35}{128}<H,i>^{2}(y_{0}%
)<H,j>^{2}(y_{0})\qquad\qquad\qquad$

$+\frac{5}{192}\overset{n}{\underset{j=q+1}{\sum}}<H,j>^{2}(y_{0})[\tau
^{M}\ -3\tau^{P}+\ \underset{\text{a}=1}{\overset{\text{q}}{\sum}}%
\varrho_{\text{aa}}^{M}+\overset{q}{\underset{\text{a},\text{b}=1}{\sum}%
}R_{\text{abab}}^{M}](y_{0})\qquad\ \ \ \ \ \ \ \ $

$+\frac{5}{192}\overset{n}{\underset{i=q+1}{\sum}}<H,i>^{2}(y_{0})[\tau
^{M}\ -3\tau^{P}+\ \underset{\text{a}=1}{\overset{\text{q}}{\sum}}%
\varrho_{\text{aa}}^{M}+\overset{q}{\underset{\text{a},\text{b}=1}{\sum}%
}R_{\text{abab}}^{M}](y_{0})\qquad\qquad$

$+\frac{5}{192}\overset{n}{\underset{i,j=q+1}{\sum}}[<H,i><H,j>](y_{0}%
)\qquad\qquad\qquad\qquad\qquad\qquad\qquad\qquad$

$\times\lbrack2\varrho_{ij}+4\overset{q}{\underset{\text{a}=1}{\sum}%
}R_{i\text{a}j\text{a}}-3\overset{q}{\underset{\text{a,b=1}}{\sum}%
}(T_{\text{aa}i}T_{\text{bb}j}-T_{\text{ab}i}T_{\text{ab}j}%
)-3\overset{q}{\underset{\text{a,b=1}}{\sum}}(T_{\text{aa}j}T_{\text{bb}%
i}-T_{\text{ab}j}T_{\text{ab}i})](y_{0})$

$+\frac{1}{96}\overset{n}{\underset{i,j=q+1}{\sum}}<H,j>(y_{0})[\{\nabla
_{i}\varrho_{ij}-2\varrho_{ij}<H,i>+\overset{q}{\underset{\text{a}=1}{\sum}%
}(\nabla_{i}R_{\text{a}i\text{a}j}-4R_{i\text{a}j\text{a}}<H,i>)\qquad$

$+4\overset{q}{\underset{\text{a,b=1}}{\sum}}R_{i\text{a}j\text{b}%
}T_{\text{ab}i}+2\overset{q}{\underset{\text{a,b,c=1}}{\sum}}(T_{\text{aa}%
i}T_{\text{bb}j}T_{\text{cc}i}-T_{\text{aa}i}T_{\text{bc}j}T_{\text{bc}%
i}-2T_{\text{bc}j}(T_{\text{aa}i}T_{\text{bc}i}-T_{\text{ab}i}T_{\text{ac}%
i}))\}$\qquad\qquad\qquad\ \ 

$+\{\nabla_{j}\varrho_{ii}-2\varrho_{ij}<H,i>+\overset{q}{\underset{\text{a}%
=1}{\sum}}(\nabla_{j}R_{\text{a}i\text{a}i}-4R_{i\text{a}j\text{a}}<H,i>)$

$+4\overset{q}{\underset{\text{a,b=1}}{\sum}}R_{j\text{a}i\text{b}%
}T_{\text{ab}i}+2\overset{q}{\underset{\text{a,b,c=1}}{\sum}}(T_{\text{aa}%
j}(T_{\text{bb}i}T_{\text{cc}i}-T_{\text{bc}i}T_{\text{bc}i})-2T_{\text{aa}%
j}T_{\text{bc}i}T_{\text{bc}i}+2T_{\text{ab}j}T_{\text{bc}i}T_{\text{ac}%
i})\}\qquad$

$+\{\nabla_{i}\varrho_{ij}-2\varrho_{ii}<H,j>+\overset{q}{\underset{\text{a}%
=1}{\sum}}(\nabla_{i}R_{\text{a}i\text{a}j}-4R_{i\text{a}i\text{a}%
}<H,j>)+4\overset{q}{\underset{\text{a,b=1}}{\sum}}R_{i\text{a}i\text{b}%
}T_{\text{ab}j}$

$+2\overset{q}{\underset{\text{a,b,c}=1}{\sum}}(T_{\text{aa}i}T_{\text{bb}%
i}T_{\text{cc}j}-3T_{\text{aa}i}T_{\text{bc}i}T_{\text{bc}j}+2T_{\text{ab}%
i}T_{\text{bc}i}T_{\text{ac}j})\}](y_{0})$

$+\frac{1}{96}\overset{n}{\underset{i,j=q+1}{\sum}}<H,i>(y_{0})[\{\nabla
_{i}\varrho_{jj}-2\varrho_{ij}<H,j>+\overset{q}{\underset{\text{a}=1}{\sum}%
}(\nabla_{i}R_{\text{a}j\text{a}j}-4R_{i\text{a}j\text{a}}<H,j>)\qquad$

$+4\overset{q}{\underset{\text{a,b=1}}{\sum}}R_{i\text{a}j\text{b}%
}T_{\text{ab}j}+2\overset{q}{\underset{\text{a,b,c=1}}{\sum}}T_{\text{aa}%
i}(T_{\text{bb}j}T_{\text{cc}j}-T_{\text{bc}j}T_{\text{bc}j})-2T_{\text{aa}%
i}T_{\text{bc}j}T_{\text{bc}j}+2T_{\text{ab}i}T_{\text{bc}j}T_{\text{ac}%
j})\}(y_{0})\qquad$\qquad\qquad\qquad\qquad\qquad\ \ 

$+\{\nabla_{j}\varrho_{ij}-2\varrho_{ij}<H,j>+\overset{q}{\underset{\text{a}%
=1}{\sum}}(\nabla_{j}R_{\text{a}i\text{a}j}-4R_{j\text{a}i\text{a}}<H,j>)$

$+4\overset{q}{\underset{\text{a,b=1}}{\sum}}R_{j\text{a}i\text{b}%
}T_{\text{ab}j}+2\overset{q}{\underset{\text{a,b,c=1}}{\sum}}(T_{\text{aa}%
j}T_{\text{bb}i}T_{\text{cc}j}-T_{\text{ab}j}T_{\text{bc}i}T_{\text{ac}%
j}-2T_{\text{bc}i}(T_{\text{aa}j}T_{\text{bc}j}-T_{\text{ab}j}T_{\text{ac}%
j}))\}(y_{0})$

$+\{\nabla_{j}\varrho_{ij}-2\varrho_{jj}<H,i>+\overset{q}{\underset{\text{a}%
=1}{\sum}}(\nabla_{j}R_{\text{a}i\text{a}j}-4R_{j\text{a}j\text{a}%
}<H,i>)+4\overset{q}{\underset{\text{a,b=1}}{\sum}}R_{j\text{a}j\text{b}%
}T_{\text{ab}i}$

$+2\overset{q}{\underset{\text{a,b,c=1}}{\sum}}(T_{\text{aa}j}T_{\text{bb}%
j}T_{\text{cc}i}-3T_{\text{aa}j}T_{\text{bc}j}T_{\text{bc}i}+2T_{\text{ab}%
j}T_{\text{bc}j}T_{\text{ac}i})\}](y_{0})$

$+\frac{1}{576}\overset{n}{\underset{i,j=q+1}{\sum}}[2\varrho_{ij}%
+4\overset{q}{\underset{\text{a}=1}{\sum}}R_{i\text{a}j\text{a}}%
-3\overset{q}{\underset{\text{a,b=1}}{\sum}}(T_{\text{aa}i}T_{\text{bb}%
j}-T_{\text{ab}i}T_{\text{ab}j})-3\overset{q}{\underset{\text{a,b=1}}{\sum}%
}(T_{\text{aa}j}T_{\text{bb}i}-T_{\text{ab}j}T_{\text{ab}i})]^{2}(y_{0})$

$+\frac{1}{288}[\tau^{M}\ -3\tau^{P}+\ \underset{\text{a}=1}{\overset{\text{q}%
}{\sum}}\varrho_{\text{aa}}^{M}+\overset{q}{\underset{\text{a},\text{b}%
=1}{\sum}}R_{\text{abab}}^{M}]^{2}(y_{0})$

$-\frac{1}{288}\overset{n}{\underset{i,j=q+1}{\sum}}[$
$\overset{q}{\underset{\text{a=1}}{\sum}}\{-(\nabla_{ii}^{2}R_{j\text{a}%
j\text{a}}+\nabla_{jj}^{2}R_{i\text{a}i\text{a}}+4\nabla_{ij}^{2}%
R_{i\text{a}j\text{a}}+2R_{ij}R_{i\text{a}j\text{a}})\qquad-\frac{1}%
{2}\overset{n}{\underset{i,j=q+1}{\sum}}\frac{\partial^{4}\theta_{p}}%
{\partial\text{x}_{i}^{2}\partial\text{x}_{j}^{2}}(y_{0})\qquad A$

$\qquad+\overset{n}{\underset{p=q+1}{\sum}}\overset{q}{\underset{\text{a=1}%
}{\sum}}(R_{\text{a}iip}R_{\text{a}jjp}+R_{\text{a}jjp}R_{\text{a}%
iip}+R_{\text{a}ijp}R_{\text{a}ijp}+R_{\text{a}ijp}R_{\text{a}jip}%
+R_{\text{a}jip}R_{\text{a}ijp}+R_{\text{a}jip}R_{\text{a}jip})$

$\qquad+2\overset{q}{\underset{\text{a,b=1}}{\sum}}\nabla_{i}(R)_{\text{a}%
i\text{b}j}T_{\text{ab}j}+2\overset{q}{\underset{\text{a,b=1}}{\sum}}%
\nabla_{j}(R)_{\text{a}j\text{b}i}T_{\text{ab}i}%
+2\overset{q}{\underset{\text{a,b=1}}{\sum}}\nabla_{i}(R)_{\text{a}j\text{b}%
i}T_{\text{ab}j}$

$\qquad+2\overset{q}{\underset{\text{a,b=1}}{\sum}}\nabla_{i}(R)_{\text{a}%
j\text{b}j}T_{\text{ab}i}$

$\qquad+2\overset{q}{\underset{\text{a,b=1}}{\sum}}\nabla_{j}(R)_{\text{a}%
i\text{b}i}T_{\text{ab}j}+2\overset{q}{\underset{\text{a,b=1}}{\sum}}%
\nabla_{j}(R)_{\text{a}i\text{b}j}T_{\text{ab}i}$

$\qquad+\overset{n}{\underset{p=q+1}{\sum}}(-\frac{3}{5}\nabla_{ii}%
^{2}(R)_{jpjp}+\overset{n}{\underset{p=q+1}{\sum}}(-\frac{3}{5}\nabla_{jj}%
^{2}(R)_{ipip}$

$+\overset{n}{\underset{p=q+1}{\sum}}(-\frac{3}{5}\nabla_{ij}^{2}%
(R)_{ipjp}+\overset{n}{\underset{p=q+1}{\sum}}(-\frac{3}{5}\nabla_{ij}%
^{2}(R)_{jpip}+\overset{n}{\underset{p=q+1}{\sum}}(-\frac{3}{5}\nabla_{ji}%
^{2}(R)_{ipjp}+\overset{n}{\underset{p=q+1}{\sum}}(-\frac{3}{5}\nabla_{ji}%
^{2}(R)_{jpip}$

$+\frac{1}{5}\overset{n}{\underset{m,p=q+1}{%
{\textstyle\sum}
}}R_{ipim}R_{jpjm}+\frac{1}{5}\overset{n}{\underset{m,p=q+1}{%
{\textstyle\sum}
}}R_{jpjm}R_{ipim}+\frac{1}{5}\overset{n}{\underset{m,p=q+1}{%
{\textstyle\sum}
}}R_{ipjm}R_{ipjm}$

$+\frac{1}{5}\overset{n}{\underset{m,p=q+1}{%
{\textstyle\sum}
}}R_{ipjm}R_{jpim}$

$+\frac{1}{5}\overset{n}{\underset{m,p=q+1}{%
{\textstyle\sum}
}}R_{jpim}R_{ipjm}+\frac{1}{5}\overset{n}{\underset{m,p=q+1}{%
{\textstyle\sum}
}}R_{jpim}R_{jpim}\}(y_{0})$

$+4\overset{q}{\underset{\text{a,b=1}}{\sum}}\{(\nabla_{i}(R)_{i\text{a}%
j\text{a}}-\overset{q}{\underset{\text{c=1}}{%
{\textstyle\sum}
}}R_{\text{a}i\text{c}i}T_{\text{ac}j})$ $T_{\text{bb}j}+4(\nabla
_{j}(R)_{j\text{a}i\text{a}}-\overset{q}{\underset{\text{c=1}}{%
{\textstyle\sum}
}}R_{\text{a}j\text{c}j}T_{\text{ac}i})$ $T_{\text{bb}i}$

$+4(\nabla_{i}(R)_{j\text{a}i\text{a}}-\overset{q}{\underset{\text{c=1}}{%
{\textstyle\sum}
}}R_{\text{a}i\text{c}j}T_{\text{ac}i})$ $T_{\text{bb}j}$ $4B\ $

$+4(\nabla_{i}(R)_{j\text{a}j\text{a}}-\overset{q}{\underset{\text{c=1}}{%
{\textstyle\sum}
}}R_{\text{a}i\text{c}j}T_{\text{ac}j})$ $T_{\text{bb}i}+4(\nabla
_{j}(R)_{i\text{a}i\text{a}}-\overset{q}{\underset{\text{c=1}}{%
{\textstyle\sum}
}}R_{\text{a}j\text{c}i}T_{\text{ac}i})$ $T_{\text{bb}j}$

$+4(\nabla_{j}(R)_{i\text{a}j\text{a}}-\overset{q}{\underset{\text{c=1}}{%
{\textstyle\sum}
}}R_{\text{a}j\text{c}i}T_{\text{ac}j})$ $T_{\text{bb}i}$

$-4\overset{q}{\underset{\text{a,b=1}}{\sum}}(\nabla_{i}(R)_{i\text{a}%
j\text{b}}-\overset{q}{\underset{\text{c=1}}{%
{\textstyle\sum}
}}R_{\text{b}r\text{c}s}T_{\text{ac}t})T_{\text{ab}j}%
-4\overset{q}{\underset{\text{a,b=1}}{\sum}}(\nabla_{j}(R)_{j\text{a}%
i\text{b}}-\overset{q}{\underset{\text{c=1}}{%
{\textstyle\sum}
}}R_{\text{b}j\text{c}j}T_{\text{ac}i})T_{\text{ab}i}$

$-4\overset{q}{\underset{\text{a,b=1}}{\sum}}(\nabla_{i}(R)_{j\text{a}%
i\text{b}}-\overset{q}{\underset{\text{c=1}}{%
{\textstyle\sum}
}}R_{\text{b}i\text{c}j}T_{\text{ac}i})T_{\text{ab}j}%
-4\overset{q}{\underset{\text{a,b=1}}{\sum}}(\nabla_{i}(R)_{j\text{a}%
j\text{b}}-\overset{q}{\underset{\text{c=1}}{%
{\textstyle\sum}
}}R_{\text{b}i\text{c}j}T_{\text{ac}j})T_{\text{ab}i}$

$-4\overset{q}{\underset{\text{a,b=1}}{\sum}}(\nabla_{j}(R)_{i\text{a}%
i\text{b}}-\overset{q}{\underset{\text{c=1}}{%
{\textstyle\sum}
}}R_{\text{b}j\text{c}i}T_{\text{ac}i})T_{\text{ab}j}%
-4\overset{q}{\underset{\text{a,b=1}}{\sum}}(\nabla_{j}(R)_{i\text{a}%
j\text{b}}-\overset{q}{\underset{\text{c=1}}{%
{\textstyle\sum}
}}R_{\text{b}j\text{c}i}T_{\text{ac}j})T_{\text{ab}i}\}](y_{0})$

$+$ $[\frac{4}{9}\overset{q}{\underset{\text{a,b=1}}{\sum}}(\varrho
_{\text{aa}}-\overset{q}{\underset{\text{c}=1}{\sum}}R_{\text{acac}}%
)(\varrho_{\text{bb}}-\overset{q}{\underset{\text{d}=1}{\sum}}R_{\text{bdbd}%
})+\frac{8}{9}\overset{n}{\underset{i,j=q+1}{\sum}}%
\overset{q}{\underset{\text{a,b}=1}{\sum}}(R_{i\text{a}j\text{a}}%
R_{i\text{b}j\text{b}})\qquad3C$

$+\frac{4}{3}\overset{q}{\underset{\text{a}=1}{\sum}}(\varrho_{\text{aa}}%
^{M}-\varrho_{\text{aa}}^{P})(\tau^{M}-\overset{q}{\underset{\text{c}=1}{\sum
}}\varrho_{\text{cc}}^{M})+\frac{4}{9}\overset{n}{\underset{i,j=q+1}{\sum}%
}\overset{q}{\underset{\text{a}=1}{\sum}}R_{i\text{a}j\text{a}}\varrho_{ij}\ $

$+\frac{2}{9}\overset{q}{\underset{\text{b}=1}{\sum}}(\varrho_{\text{bb}}%
^{M}-\varrho_{\text{bb}}^{P})(\tau^{M}-\overset{q}{\underset{\text{c}=1}{\sum
}}\varrho_{\text{cc}}^{M})+\frac{4}{9}\overset{n}{\underset{i,j=q+1}{\sum}%
}\overset{q}{\underset{\text{b}=1}{\sum}}R_{i\text{b}j\text{b}}\varrho_{ij}\ $

$+\frac{1}{9}(\tau^{M}-\overset{q}{\underset{\text{a=1}}{\sum}}\varrho
_{\text{aa}})(\tau^{M}-\overset{q}{\underset{\text{b=1}}{\sum}}\varrho
_{\text{bb}})+\frac{2}{9}(\left\Vert \varrho^{M}\right\Vert ^{2}%
-\overset{q}{\underset{\text{a,b}=1}{\sum}}\varrho_{\text{ab}})$

$-\overset{n}{\underset{i,j=q+1}{\sum}}\overset{q}{\underset{\text{a,b}%
=1}{\sum}}R_{i\text{a}i\text{b}}R_{j\text{a}j\text{b}}\ -\frac{1}%
{2}\overset{n}{\underset{i,j=q+1}{\sum}}\overset{q}{\underset{\text{a,b}%
=1}{\sum}}R_{i\text{a}j\text{b}}^{2}-\overset{n}{\underset{i,j=q+1}{\sum}%
}\overset{q}{\underset{\text{a,b}=1}{\sum}}R_{i\text{a}j\text{b}}%
R_{j\text{a}i\text{b}}$

$-\frac{1}{2}\overset{n}{\underset{i,j=q+1}{\sum}}%
\overset{q}{\underset{\text{a,b}=1}{\sum}}R_{j\text{a}i\text{b}}^{2}%
\qquad\qquad$

$-\frac{1}{9}\overset{n}{\underset{i,j,p,m=q+1}{\sum}}R_{ipim}R_{jpjm}%
\ -\frac{1}{18}\overset{n}{\underset{i,j,p,m=q+1}{\sum}}R_{ipjm}^{2}-\frac
{1}{9}\overset{n}{\underset{i,j,p,m=q+1}{\sum}}R_{ipjm}R_{jpim}$

$-\frac{1}{18}\overset{n}{\underset{i,j,p,m=q+1}{\sum}}R_{jpim}^{2}$

$-\frac{1}{3}\overset{q}{\underset{\text{a}=1}{\sum}}%
\overset{n}{\underset{i,j,p=q+1}{\sum}}R_{i\text{a}ip}R_{j\text{a}jp}-\frac
{1}{6}\overset{q}{\underset{\text{a}=1}{\sum}}%
\overset{n}{\underset{i,j,p=q+1}{\sum}}R_{i\text{a}jp}^{2}-\frac{1}%
{3}\overset{q}{\underset{\text{a}=1i,j,}{\sum}}%
\overset{n}{\underset{p=q+1}{\sum}}R_{i\text{a}jp}R_{j\text{a}ip}$

$-\frac{1}{6}\overset{q}{\underset{\text{a}=1}{\sum}}%
\overset{n}{\underset{i,j,p=q+1}{\sum}}R_{j\text{a}ip}^{2}$

$-\frac{1}{3}\overset{q}{\underset{\text{b}=1i,j,}{\sum}}%
\overset{n}{\underset{p=q+1}{\sum}}R_{i\text{b}ip}R_{j\text{b}jp}-\frac{1}%
{6}\overset{q}{\underset{\text{b}=1}{\sum}}%
\overset{n}{\underset{i,j,p=q+1}{\sum}}R_{i\text{b}jp}^{2}-\frac{1}%
{3}\overset{q}{\underset{\text{b}=1}{\sum}}%
\overset{n}{\underset{i.j,p=q+1}{\sum}}R_{i\text{b}jp}R_{j\text{b}ip}$

$-\frac{1}{6}\overset{q}{\underset{\text{b}=1}{\sum}}%
\overset{n}{\underset{i,j,p=q+1}{\sum}}R_{j\text{b}ip}^{2}](y_{0})$

$+$ $\overset{q}{\underset{\text{a,b,c=1}}{\sum}}[$
$-\overset{n}{\underset{i=q+1}{\sum}}R_{i\text{a}i\text{a}}(R_{\text{bcbc}%
}^{P}-R_{\text{bcbc}}^{M})$ $-\overset{n}{\underset{j=q+1}{\sum}}%
R_{j\text{a}j\text{a}}(R_{\text{bcbc}}^{P}-R_{\text{bcbc}}^{M})\qquad\qquad6D$

$+\overset{n}{\underset{i=q+1}{\sum}}R_{i\text{a}i\text{b}}(R_{\text{acbc}%
}^{P}-R_{\text{acbc}}^{M})\ -\overset{n}{\underset{i=q+1}{\sum}}%
R_{i\text{a}i\text{c}}(R_{\text{abbc}}^{P}-R_{\text{abbc}}^{M})$

$+\overset{n}{\underset{j=q+1}{\sum}}R_{j\text{a}j\text{b}}(R_{\text{acbc}%
}^{P}-R_{\text{acbc}}^{M})$\ $-\overset{n}{\underset{j=q+1}{\sum}}%
R_{j\text{a}j\text{c}}(R_{\text{abbc}}^{P}-R_{\text{abbc}}^{M})$

$+\underset{i,j=q+1}{\overset{n}{\sum}}$ $-R_{i\text{a}j\text{a}}%
(T_{\text{bb}i}T_{\text{cc}j}$ $-T_{\text{bc}i}T_{\text{bc}j})$
$-\underset{i,j=q+1}{\overset{n}{\sum}}R_{i\text{a}j\text{a}}(T_{\text{bb}%
j}T_{\text{cc}i}$ $-T_{\text{bc}j}T_{\text{bc}i})$

$+$ $\underset{i,j=q+1}{\overset{n}{\sum}}$ $-R_{j\text{a}i\text{a}%
}(T_{\text{bb}i}T_{\text{cc}j}$ $-T_{\text{bc}i}T_{\text{bc}j})$
$-\underset{i,j=q+1}{\overset{n}{\sum}}R_{j\text{a}i\text{a}}(T_{\text{bb}%
j}T_{\text{cc}i}$ $-T_{\text{bc}j}T_{\text{bc}i})$

$+\underset{i,j=q+1}{\overset{n}{\sum}}\ R_{i\text{a}j\text{b}}(T_{\text{ab}%
i}T_{\text{cc}j}-T_{\text{bc}i}T_{\text{ac}j}%
)\ +\underset{i,j=q+1}{\overset{n}{\sum}}\ R_{i\text{a}j\text{b}}%
(T_{\text{ab}j}T_{\text{cc}i}-T_{\text{bc}j}T_{\text{ac}i})$

$+\underset{i,j=q+1}{\overset{n}{\sum}}\ R_{j\text{a}i\text{ib}}%
(T_{\text{ab}i}T_{\text{cc}j}-T_{\text{bc}i}T_{\text{ac}j}%
)\ +\underset{i,j=q+1}{\overset{n}{\sum}}\ R_{j\text{a}i\text{b}}%
(T_{\text{ab}j}T_{\text{cc}i}-T_{\text{bc}j}T_{\text{ac}i})\qquad$

$+\underset{i,j=q+1}{\overset{n}{\sum}}-R_{i\text{a}j\text{c}}(T_{\text{ab}%
i}T_{\text{bc}j}-T_{\text{ac}i}T_{\text{bb}j}%
)-\underset{i,j=q+1}{\overset{n}{\sum}}R_{i\text{a}j\text{c}}(T_{\text{ba}%
j}T_{\text{bc}i}-T_{\text{ac}j}T_{\text{bb}i})$

$+\underset{i,j=q+1}{\overset{n}{\sum}}-R_{j\text{a}i\text{c}}(T_{\text{ba}%
i}T_{\text{bc}j}-T_{\text{ac}i}T_{\text{bb}j}%
)-\underset{i,j=q+1}{\overset{n}{\sum}}R_{j\text{a}i\text{c}}(T_{\text{ba}%
j}T_{\text{bc}i}-T_{\text{ac}j}T_{\text{bb}i})](y_{0})$

$-\frac{1}{3}\underset{p=q+1}{\overset{n}{\sum}}%
[\underset{i=q+1}{\overset{n}{\sum}}\overset{q}{\underset{\text{b,c=1}}{\sum}%
}R_{ipip}(R_{\text{bcbc}}^{P}-R_{\text{bcbc}}^{M}%
)+\underset{j=q+1}{\overset{n}{\sum}}$ $\overset{q}{\underset{\text{b,c=1}%
}{\sum}}R_{jpjp}(R_{\text{bcbc}}^{P}-R_{\text{bcbc}}^{M})](y_{0})$

$-\frac{2}{3}\underset{i,j,p=q+1}{\overset{n}{\sum}}%
\overset{q}{\underset{\text{b,c=1}}{\sum}}[R_{ipjp}(T_{\text{bb}i}%
T_{\text{cc}j}-T_{\text{bc}i}T_{\text{bc}j})+R_{ipjp}(T_{\text{bb}%
j}T_{\text{cc}i}-T_{\text{bc}j}T_{\text{bc}i})](y_{0})\qquad$

$+\ \frac{1}{6}\underset{i,j=q+1}{\overset{n}{\sum}}[T_{\text{aa}%
i}T_{\text{bb}j}(T_{\text{cc}i}T_{\text{dd}j}-T_{\text{cd}i}T_{\text{dc}%
j})+T_{\text{aa}i}T_{\text{bb}j}(T_{\text{cc}j}T_{\text{dd}i}-T_{\text{cd}%
j}T_{\text{dc}i})\qquad E$

$+T_{\text{aa}j}T_{\text{bb}i}(T_{\text{cc}i}T_{\text{dd}j}-T_{\text{cd}%
i}T_{\text{dc}j})+T_{\text{aa}j}T_{\text{bb}i}(T_{\text{cc}j}T_{\text{dd}%
i}-T_{\text{cd}j}T_{\text{dc}i})](y_{0})$

$-\frac{1}{6}\underset{i,j=q+1}{\overset{n}{\sum}}[T_{\text{aa}i}%
T_{\text{bc}j}(T_{\text{bc}i}T_{\text{dd}j}-T_{\text{bd}i}T_{\text{cd}%
j})+T_{\text{aa}i}T_{\text{bc}j}(T_{\text{bc}j}T_{\text{dd}i}-T_{\text{bd}%
j}T_{\text{cd}i})$

$+T_{\text{aa}j}T_{\text{bc}i}(T_{\text{bc}i}T_{\text{dd}j}-T_{\text{bd}%
i}T_{\text{cd}j})+T_{\text{aa}j}T_{\text{bc}i}(T_{\text{bc}j}T_{\text{dd}%
i}-T_{\text{bd}j}T_{\text{cd}i})](y_{0})$

$+\ \frac{1}{6}\underset{i,j=q+1}{\overset{n}{\sum}}[T_{\text{aa}%
i}T_{\text{bd}j}(T_{\text{bc}i}T_{\text{cd}j}-T_{\text{bd}i}T_{\text{cc}%
j})+T_{\text{aa}i}T_{\text{bd}j}(T_{\text{bc}j}T_{\text{cd}i}-T_{\text{bd}%
j}T_{\text{cc}i})$

$+T_{\text{aa}j}T_{\text{bd}i}(T_{\text{bc}i}T_{\text{cd}j}-T_{\text{bd}%
i}T_{\text{cc}j})+T_{\text{aa}j}T_{\text{bd}i}(T_{\text{bc}j}T_{\text{cd}%
i}-T_{\text{bd}j}T_{\text{cc}i})](y_{0})\qquad$

$-\ \frac{1}{6}\underset{i,j=q+1}{\overset{n}{\sum}}[T_{\text{ab}%
i}T_{\text{ab}j}(T_{\text{cc}i}T_{\text{dd}j}-T_{\text{cd}i}T_{\text{dc}%
j})+T_{\text{ab}i}T_{\text{ab}j}(T_{\text{cc}j}T_{\text{dd}i}-T_{\text{cd}%
j}T_{\text{dc}i})$

$+T_{\text{ab}j}T_{\text{ab}i}(T_{\text{cc}i}T_{\text{dd}j}-T_{\text{cd}%
i}T_{\text{dc}j})+T_{\text{ab}j}T_{\text{ab}i}(T_{\text{cc}j}T_{\text{dd}%
i}-T_{\text{cd}j}T_{\text{dc}i})](y_{0})$

$+\ \frac{1}{6}\underset{i,j=q+1}{\overset{n}{\sum}}[T_{\text{ab}%
i}T_{\text{bc}j}(T_{\text{ac}i}T_{\text{dd}j}-T_{\text{ad}i}T_{\text{cd}%
j})+T_{\text{ab}i}T_{\text{bc}j}(T_{\text{ac}j}T_{\text{dd}i}-T_{\text{ad}%
j}T_{\text{cd}i})$

$+T_{\text{ab}j}T_{\text{bc}i}(T_{\text{ac}i}T_{\text{dd}j}-T_{\text{ad}%
i}T_{\text{cd}j})+T_{\text{ab}j}T_{\text{bc}i}(T_{\text{ac}j}T_{\text{dd}%
i}-T_{\text{ad}j}T_{\text{cd}i})](y_{0})$

$-\ \frac{1}{6}\underset{i,j=q+1}{\overset{n}{\sum}}[T_{\text{ab}%
i}T_{\text{bd}j}(T_{\text{ac}i}T_{\text{cd}j}-T_{\text{ad}i}T_{\text{cc}%
j})+T_{\text{ab}i}T_{\text{bd}j}(T_{\text{ac}j}T_{\text{cd}i}-T_{\text{ad}%
j}T_{\text{cc}i})$

$+T_{\text{ab}i}T_{\text{bd}j}(T_{\text{ac}j}T_{\text{cd}i}-T_{\text{ad}%
j}T_{\text{cc}i})+T_{\text{ab}j}T_{\text{bd}i}(T_{\text{ac}j}T_{\text{cd}%
i}-T_{\text{ad}j}T_{\text{cc}i})](y_{0})$

$+\ \frac{1}{6}\underset{i,j=q+1}{\overset{n}{\sum}}[T_{\text{ac}%
i}T_{\text{ab}j}(T_{\text{bc}i}T_{\text{dd}j}-T_{\text{bd}i}T_{\text{dc}%
j})+T_{\text{ac}i}T_{\text{ab}j}(T_{\text{bc}j}T_{\text{dd}i}-T_{\text{bd}%
j}T_{\text{dc}i})$

$+T_{\text{ac}j}T_{\text{ab}i}(T_{\text{bc}i}T_{\text{dd}j}-T_{\text{bd}%
i}T_{\text{dc}j})+T_{\text{ac}j}T_{\text{ab}i}(T_{\text{bc}j}T_{\text{dd}%
i}-T_{\text{bd}j}T_{\text{dc}i})](y_{0})$

$-\ \frac{1}{6}\underset{i,j=q+1}{\overset{n}{\sum}}[T_{\text{ac}%
i}T_{\text{bb}j}(T_{\text{ac}i}T_{\text{dd}j}-T_{\text{ad}i}T_{\text{cd}%
j})+T_{\text{ac}i}T_{\text{bb}j}(T_{\text{ac}j}T_{\text{dd}i}-T_{\text{ad}%
j}T_{\text{cd}i})$

$+T_{\text{ac}j}T_{\text{bb}i}(T_{\text{ac}i}T_{\text{dd}j}-T_{\text{ad}%
i}T_{\text{cd}i})+T_{\text{ac}j}T_{\text{bb}i}(T_{\text{ac}j}T_{\text{dd}%
i}-T_{\text{ad}j}T_{\text{cd}i})](y_{0})$

$+\ \frac{1}{6}\underset{i,j=q+1}{\overset{n}{\sum}}[T_{\text{ac}%
i}T_{\text{bd}j}(T_{\text{ac}i}T_{\text{bd}j}-T_{\text{ad}i}T_{\text{bc}%
j})+T_{\text{ac}i}T_{\text{bd}j}(T_{\text{ac}j}T_{\text{bd}i}-T_{\text{ad}%
j}T_{\text{bc}i})$

$+T_{\text{ac}j}T_{\text{bd}i}(T_{\text{ac}i}T_{\text{bd}j}-T_{\text{ad}%
i}T_{\text{bc}j})+T_{\text{ac}j}T_{\text{bd}i}(T_{\text{ac}j}T_{\text{bd}%
i}-T_{\text{ad}j}T_{\text{bc}i})](y_{0})$

$-\frac{1}{6}\underset{i,j=q+1}{\overset{n}{\sum}}[T_{\text{ad}i}%
T_{\text{ab}j}(T_{\text{bc}i}T_{\text{cd}j}-T_{\text{bd}i}T_{\text{cc}%
j})+T_{\text{ad}i}T_{\text{ab}j}(T_{\text{bc}j}T_{\text{cd}i}-T_{\text{bd}%
j}T_{\text{cc}i})$

$+T_{\text{ad}j}T_{\text{ab}i}(T_{\text{bc}i}T_{\text{cd}j}-T_{\text{bd}%
i}T_{\text{cc}j})+T_{\text{ad}j}T_{\text{ab}i}(T_{\text{bc}j}T_{\text{cd}%
i}-T_{\text{bd}j}T_{\text{cc}i})](y_{0})$

$+\ \frac{1}{6}\underset{i,j=q+1}{\overset{n}{\sum}}[T_{\text{ad}%
i}T_{\text{bb}j}(T_{\text{ac}i}T_{\text{cd}j}-T_{\text{ad}i}T_{\text{cc}%
j})+T_{\text{ad}i}T_{\text{bb}j}(T_{\text{ac}j}T_{\text{cd}i}-T_{\text{ad}%
j}T_{\text{cc}i})$

$+T_{\text{ad}j}T_{\text{bb}i}(T_{\text{ac}i}T_{\text{cd}j}-T_{\text{ad}%
i}T_{\text{cc}j})+T_{\text{ad}j}T_{\text{bb}i}(T_{\text{ac}j}T_{\text{cd}%
i}-T_{\text{ad}j}T_{\text{cc}i})](y_{0})$

$-\ \frac{1}{6}\underset{i,j=q+1}{\overset{n}{\sum}}[T_{\text{ad}%
i}T_{\text{bc}j}(T_{\text{ac}i}T_{\text{bd}j}-T_{\text{ad}i}T_{\text{bc}%
j})+T_{\text{ad}i}T_{\text{bc}j}(T_{\text{ac}j}T_{\text{bd}i}-T_{\text{ad}%
j}T_{\text{bc}i})$

$+T_{\text{ad}j}T_{\text{bc}i}(T_{\text{ac}i}T_{\text{bd}j}-T_{\text{ad}%
i}T_{\text{bc}j})+T_{\text{ad}j}T_{\text{bc}i}(T_{\text{ac}j}T_{\text{bd}%
i}-T_{\text{ad}j}T_{\text{bc}i})](y_{0})$

$+\ \frac{1}{3}[(R_{\text{cdcd}}^{P}-R_{\text{cdcd}}^{M})(R_{\text{abab}}%
^{P}-R_{\text{abab}}^{M})](y_{0})\qquad\qquad$

$-\frac{1}{3}[(R_{\text{bdcd}}^{P}-R_{\text{bdcd}}^{M})(R_{\text{abac}}%
^{P}-R_{\text{abac}}^{M})](y_{0})\qquad\qquad$

$-\ \frac{1}{3}[(R_{\text{bcdc}}^{P}-R_{\text{bcdc}}^{M})(R_{\text{abad}}%
^{P}-R_{\text{abad}}^{M})](y_{0})\qquad\qquad$

$+\ \frac{1}{3}[(R_{\text{adcd}}^{P}-R_{\text{adcd}}^{M})(R_{\text{abbc}}%
^{P}-R_{\text{abbc}}^{M})](y_{0})\qquad\qquad$

$-\ \frac{1}{3}[(R_{\text{acdc}}^{P}-R_{\text{acdc}}^{M})(R_{\text{abdb}}%
^{P}-R_{\text{abdb}}^{M})](y_{0})$

$+\frac{1}{3}[(R_{\text{abcd}}^{P}-R_{\text{abcd}}^{M})]^{2}(y_{0}%
)\qquad\qquad\qquad$

\qquad\qquad\qquad\qquad\qquad\qquad\qquad\qquad\qquad\qquad\qquad\qquad
\qquad\qquad\qquad\qquad$\blacksquare$

\subsection{Computations}

We use the expansion in \textbf{Proposition 11} of \textbf{Chapter 10} in
\textbf{Part 4} given there as follows:

$\theta_{P}(x)=1-\underset{r=q+1}{\overset{n}{\sum}}<H,r>(y_{0})x_{r}$

$\qquad\qquad-\frac{1}{6}\underset{r,s=q+1}{\overset{n}{\sum}}[\varrho
_{rs}+\overset{q}{\underset{\text{a}=1}{2\sum}}R_{r\text{a}s\text{a}%
}-3\overset{q}{\underset{\text{a,b=1}}{\sum}}(T_{\text{aa}r}T_{\text{bb}%
s}-T_{\text{ab}r}T_{\text{ab}s})](y_{0})x_{r}x_{s}$

$-\frac{1}{12}\underset{r,s,t=q+1}{\overset{n}{\sum}}[\nabla_{r}\varrho
_{st}-2\varrho_{rs}<H,t>+\overset{q}{\underset{\text{a=1}}{\sum}}(\nabla
_{r}R_{\text{a}s\text{a}t}-4R_{r\text{a}s\text{a}}<H,t>)$

$+4\overset{q}{\underset{\text{a},\text{b}=1}{\sum}}R_{r\text{a}s\text{b}%
}T_{\text{ab}t}\ \ +2\overset{q}{\underset{\text{a},\text{b,c}=1}{\sum}%
}(T_{\text{aa}r}T_{\text{bb}s}T_{\text{cc}t}-3T_{\text{aa}r}TT_{\text{bc}%
s}T_{\text{bc}t}+2T_{\text{ab}r}T_{\text{bc}s}T_{\text{ca}t}](y_{0})x_{r}%
x_{s}x_{t}$

$+\ \ \frac{1}{24}\overset{n}{\underset{r,s,t,u=q+1}{\sum}}[$
$\overset{q}{\underset{\text{a=1}}{\sum}}\{-\nabla_{rs}^{2}(R)_{t\text{a}%
u\text{a}}+\overset{n}{\underset{p=q+1}{\sum}}\overset{q}{\underset{\text{a=1}%
}{\sum}}R_{\text{a}rsp}R_{\text{a}tup}+2\overset{q}{\underset{\text{a,b=1}%
}{\sum}}\nabla_{r}(R)_{\text{a}s\text{b}t}T_{\text{ab}u}\qquad A$

$+\overset{n}{\underset{p=q+1}{\sum}}(-\frac{3}{5}\nabla_{rs}^{2}%
(R)_{tpup}+\frac{1}{5}\overset{n}{\underset{m=q+1}{%
{\textstyle\sum}
}}R_{rpsm}R_{tpum})\}(y_{0})$

$+4\overset{q}{\underset{\text{a,b=1}}{\sum}}\{(\nabla_{r}(R)_{s\text{a}%
t\text{a}}-\overset{q}{\underset{\text{c=1}}{%
{\textstyle\sum}
}}R_{\text{a}r\text{c}s}T_{\text{ac}t})$ $T_{\text{bb}u}%
-4\overset{q}{\underset{\text{a,b=1}}{\sum}}(\nabla_{r}(R)_{s\text{a}%
t\text{b}}-\overset{q}{\underset{\text{c=1}}{%
{\textstyle\sum}
}}R_{\text{b}r\text{c}s}T_{\text{ac}t})T_{\text{ab}u}\}\qquad4B$

$+\frac{4}{3}\overset{q}{\underset{\text{a,b}=1}{\sum}}(R_{r\text{a}s\text{a}%
})(R_{t\text{b}u\text{b}})+\frac{1}{3}\varrho_{rs}\varrho_{tu}+\frac{2}%
{3}\overset{q}{\underset{\text{a}=1}{\sum}}R_{r\text{a}s\text{a}}\varrho
_{tu}\ +\frac{2}{3}\overset{q}{\underset{\text{b=1}}{\sum}}R_{r\text{b}%
s\text{b}}\varrho_{tu}=3C$

\vspace{1pt}$-3\overset{q}{\underset{\text{a,b}=1}{\sum}}R_{r\text{a}%
s\text{b}}R_{t\text{a}u\text{b}}-\frac{1}{3}%
\overset{n}{\underset{p,m=q+1}{\sum}}R_{rpsm}R_{tpum}%
\ -\overset{q}{\underset{\text{a}=1}{\sum}}\overset{n}{\underset{p=q+1}{\sum}%
}R_{r\text{a}sp}R_{t\text{a}up}-\overset{q}{\underset{\text{b}=1}{\sum}%
}\overset{n}{\underset{p=q+1}{\sum}}R_{r\text{b}sp}R_{t\text{b}up}$

$+$ $\overset{q}{\underset{\text{a,b,c=1}}{6\sum}}\{$ $-R_{r\text{a}s\text{a}%
}(T_{\text{bb}t}T_{\text{cc}u}$ $-T_{\text{bc}t}T_{\text{bc}u})\}$%
\ $+6\{$\ $R_{r\text{a}s\text{b}}(T_{\text{ab}t}T_{\text{cc}u}-T_{\text{bc}%
t}T_{\text{ac}u})\}\qquad6D$

$+6\{-R_{r\text{a}s\text{c}}(T_{\text{ba}t}T_{\text{bc}u}-T_{\text{ac}%
t}T_{\text{bb}u})\}+6\overset{q}{\underset{\text{b,c=1}}{\sum}}%
\underset{p=q+1}{\overset{n}{\sum}}$\ $\{-\frac{1}{3}R_{rpsp}(T_{\text{bb}%
t}T_{\text{cc}u}-T_{\text{bc}t}T_{\text{bc}u})\}\qquad$

$+\overset{n}{\underset{r,s,t,u=q+1}{\sum}}%
\overset{q}{\underset{\text{a,b,c,d}=1}{\sum}}T_{\text{aa}r}\{T_{\text{bb}%
s}(T_{\text{cc}t}T_{\text{dd}u}-T_{\text{cd}t}T_{\text{dc}u})-T_{\text{bc}%
s}(T_{\text{bc}t}T_{\text{dd}u}-T_{\text{bd}t}T_{\text{cd}u})$

$+T_{\text{bd}s}(T_{\text{bc}t}T_{\text{cd}u}-T_{\text{bd}t}T_{\text{cc}%
u})\}=E$

$-T_{\text{ab}r}\{T_{\text{ab}s}(T_{\text{cc}t}T_{\text{dd}u}-T_{\text{cd}%
t}T_{\text{dc}u})-T_{\text{bc}s}(T_{\text{ac}t}T_{\text{dd}u}-T_{\text{ad}%
t}T_{\text{cd}u})+T_{\text{bd}s}(T_{\text{ac}t}T_{\text{cd}u}-T_{\text{ad}%
t}T_{\text{cc}u})\}$

\vspace{1pt}$+T_{\text{ac}r}\{T_{\text{ab}s}(T_{\text{bc}t}T_{\text{dd}%
u}-T_{\text{bd}t}T_{\text{dc}u})-T_{\text{bb}s}(T_{\text{ac}t}T_{\text{dd}%
u}-T_{\text{ad}t}T_{\text{cd}u})+T_{\text{bd}s}(T_{\text{ac}t}T_{\text{bd}%
u}-T_{\text{ad}t}T_{\text{bc}u})\}$

\vspace{1pt}$-T_{\text{ad}r}\{T_{\text{ab}s}(T_{\text{bc}t}T_{\text{cd}%
u}-T_{\text{bd}t}T_{\text{cc}u})-T_{\text{bb}s}(T_{\text{ac}t}T_{\text{cd}%
u}-T_{\text{ad}t}T_{\text{cc}u})$

$+T_{\text{bc}s}(T_{\text{ac}t}T_{\text{bd}u}-T_{\text{ad}t}T_{\text{bc}%
u})\}]$(y$_{0}$)x$_{r}$x$_{s}$x$_{t}$x$_{u}$\ $+$ higher order terms.

\qquad\qquad\qquad\qquad\qquad\qquad\qquad\qquad\qquad\qquad\qquad\qquad
\qquad\qquad\qquad\qquad\qquad\qquad$\blacksquare$

(i), (ii) and (iii) are immediate from the above expansion:

For $i,j,k=q+1,...,n,$ we have:

(i)$\ \ \ \ \theta(y_{0})$ \ \ $=1$

$(\nabla\log\theta^{-\frac{1}{2}})_{i}(y_{0})=-\frac{1}{2}(\nabla\log
\theta)_{i}(y_{0})=-\frac{1}{2}\frac{1}{\theta(y_{0})}(\nabla\theta)_{i}%
(y_{0})=-\frac{1}{2}(\nabla\theta)_{i}(y_{0})$

\qquad$=-\ \frac{1}{2}\frac{\partial\theta}{\partial\text{x}_{i}}%
(y_{0})=\left\{
\begin{array}
[c]{c}%
0\text{ for 1=1,...,q}\\
\frac{1}{2}<H,i>(y_{0})\text{ for }i=q+1,...,n
\end{array}
\right.  $

The short computation above is proof of (ii), (iii), (iii)$^{\ast},$ (iv) and
(iv)$^{\ast}$

(ii)$\ \frac{\partial\theta}{\partial\text{x}_{i}}(y_{0})=-<H,i>(y_{0})$

(iii) $\frac{\partial\theta^{\frac{1}{2}}}{\partial\text{x}_{i}}%
(y_{0})=\ \frac{1}{2}\frac{\partial\theta}{\partial\text{x}_{i}}(y_{0}%
)=-\frac{1}{2}<H,i>(y_{0})$

(iii)$^{\ast}$ $(\nabla\log\theta^{-\frac{1}{2}})_{\text{a}}(y_{0})=0$

(iv) $\frac{\partial\theta^{-\frac{1}{2}}}{\partial\text{x}_{i}}(y_{0}%
)=\frac{1}{2}<H,i>(y_{0})$

\vspace{1pt}(iv)$^{\ast}$ $(\nabla\log\theta^{-\frac{1}{2}})_{i}(y_{0}%
)=\frac{1}{2}<H,i>(y_{0})$

(v) It is immediate from the expansion given above that:

$\frac{\partial^{2}\theta}{\partial\text{x}_{i}\partial\text{x}_{j}}(y_{0})$

$=-\frac{1}{6}[2\varrho_{ij}+4\overset{q}{\underset{\text{a}=1}{\sum}%
}R_{i\text{a}j\text{a}}-3\overset{q}{\underset{\text{a,b=1}}{\sum}%
}(T_{\text{aa}i}T_{\text{bb}j}-T_{\text{ab}i}T_{\text{ab}j}%
)-3\overset{q}{\underset{\text{a,b=1}}{\sum}}(T_{\text{aa}j}T_{\text{bb}%
i}-T_{\text{ab}j}T_{\text{ab}i})](y_{0})$

We easily deduce from the above equality that:

$\qquad\frac{\partial^{2}\theta}{\partial\text{x}_{i}^{2}}(y_{0})=-\frac{1}%
{6}[2\varrho_{ii}+4\overset{q}{\underset{\text{a=1}}{\sum}}R_{i\text{a}%
i\text{a}}-6\overset{q}{\underset{\text{a,b=1}}{\sum}}(T_{\text{aa}%
i}T_{\text{bb}i}-T_{\text{ab}i}T_{\text{ab}i})](y_{0})$

\qquad$\frac{\partial^{2}\theta}{\partial\text{x}_{i}^{2}}(y_{0})=-\frac{1}%
{3}[\varrho_{ii}+2\overset{q}{\underset{\text{a=1}}{\sum}}R_{i\text{a}%
i\text{a}}-3\overset{q}{\underset{\text{a,b=1}}{\sum}}(T_{\text{aa}%
i}T_{\text{bb}i}-T_{\text{ab}i}T_{\text{ab}i})](y_{0})$

From the \textbf{Gauss Equation} we have:

(vi) $\frac{\partial^{2}\theta}{\partial\text{x}_{i}^{2}}(y_{0})=-\frac{1}%
{3}[\tau^{M}-3\tau^{P}+\overset{q}{\underset{\text{a}=1}{\sum}}\varrho
_{\text{aa}}^{M}+\overset{q}{\underset{\text{a,b}=1}{\sum}}R_{\text{abab}}%
^{M}](y_{0})$

(vii) $\frac{\partial^{2}\theta^{\frac{1}{2}}}{\partial\text{x}_{i}%
\partial\text{x}_{j}}(y_{0})=-\frac{1}{4}\frac{\partial\theta}{\partial
\text{x}_{i}}(y_{0})\frac{\partial\theta}{\partial\text{x}_{j}}(y_{0}%
)+\frac{1}{2}\frac{\partial^{2}\theta}{\partial\text{x}_{i}\partial
\text{x}_{j}}(y_{0})$

\qquad$=-\frac{1}{4}<H,i><H,j>$\ $-\frac{1}{12}[2\varrho_{ij}%
+4\overset{q}{\underset{\text{a}=1}{\sum}}R_{i\text{a}j\text{a}}%
-3\overset{q}{\underset{\text{a,b=1}}{\sum}}(T_{\text{aa}i}T_{\text{bb}%
j}+T_{\text{aa}j}T_{\text{bb}i}-2T_{\text{ab}i}T_{\text{ab}j}](y_{0})$

(viii) $\frac{\partial^{2}\theta^{\frac{1}{2}}}{\partial\text{x}_{i}^{2}%
}(y_{0})=-\frac{1}{4}(\frac{\partial\theta}{\partial\text{x}_{i}})^{2}%
(y_{0})+\frac{1}{2}\frac{\partial^{2}\theta}{\partial\text{x}_{i}^{2}}(y_{0})$

\qquad\qquad$=$ $-\frac{1}{4}<H,i>^{2}(y_{0})-\frac{1}{3}[\varrho
_{ii}+2\overset{q}{\underset{\text{a=1}}{\sum}}R_{i\text{a}i\text{a}%
}-3\overset{q}{\underset{\text{a,b=1}}{\sum}}(T_{\text{aa}i}T_{\text{bb}%
i}-T_{\text{ab}i}T_{\text{ab}i})](y_{0})$\qquad

\bigskip(ix)\qquad\ $\frac{\partial^{2}\theta^{-\frac{1}{2}}}{\partial
\text{x}_{i}\partial\text{x}_{j}}(y_{0})=\frac{3}{4}\frac{\partial\theta
}{\partial\text{x}_{i}}(y_{0})\frac{\partial\theta}{\partial\text{x}_{j}%
}(y_{0})-\frac{1}{2}\frac{\partial^{2}\theta}{\partial\text{x}_{i}%
\partial\text{x}_{j}}(y_{0})$

\ $=\frac{3}{4}<H,i><H,j>$

$+\frac{1}{12}[2\varrho_{ij}+4\overset{q}{\underset{\text{a}=1}{\sum}%
}R_{i\text{a}j\text{a}}-3\overset{q}{\underset{\text{a,b=1}}{\sum}%
}(T_{\text{aa}i}T_{\text{bb}j}-T_{\text{ab}i}T_{\text{ab}j}%
)-3\overset{q}{\underset{\text{a,b=1}}{\sum}}(T_{\text{aa}j}T_{\text{bb}%
i}-T_{\text{ab}j}T_{\text{ab}i})](y_{0})$

(ix)$^{\ast}\qquad\frac{\partial}{\partial x_{i}}(\nabla\log\theta^{-\frac
{1}{2}})_{\text{a}}(y_{0})=-\frac{1}{2}\overset{q}{\underset{j=q+1}{\sum}%
}\perp_{\text{a}ij}(y_{0})<H,j>(y_{0})$

(ix)$^{\ast\ast}$\qquad$\frac{\partial}{\partial x_{i}}(\nabla\log
\theta^{-\frac{1}{2}})_{j}(y_{0})$ for $i,j=q+1,...n$

\qquad$\qquad=\frac{1}{2}<H,i><H,j>\ $

\qquad\qquad$+\frac{1}{12}[2\varrho_{ij}+4\overset{q}{\underset{\text{a}%
=1}{\sum}}R_{i\text{a}j\text{a}}-3\overset{q}{\underset{\text{a,b=1}}{\sum}%
}(T_{\text{aa}i}T_{\text{bb}j}-T_{\text{ab}i}T_{\text{ab}j}%
)-3\overset{q}{\underset{\text{a,b=1}}{\sum}}(T_{\text{aa}j}T_{\text{bb}%
i}-T_{\text{ab}j}T_{\text{ab}i})](y_{0})$

\vspace{1pt}(x)\qquad\ $\frac{\partial^{2}\theta^{-\frac{1}{2}}}%
{\partial\text{x}_{i}^{2}}(y_{0})=\frac{3}{4}(\frac{\partial\theta}%
{\partial\text{x}_{i}})^{2}(y_{0})-\frac{1}{2}\frac{\partial^{2}\theta
}{\partial\text{x}_{i}^{2}}(y_{0})$

\qquad\qquad$=\frac{3}{4}<H,i>^{2}(y_{0})+\frac{1}{6}(\tau^{M}-3\tau
^{P}+\overset{q}{\underset{\text{a}=1}{\sum}}\varrho_{\text{aa}}%
^{M}+\overset{q}{\underset{\text{a,b}=1}{\sum}}R_{\text{abab}}^{M})(y_{0})$

\qquad\qquad$=\frac{1}{12}[9<H,i>^{2}+2(\tau^{M}-3\tau^{P}%
+\overset{q}{\underset{\text{a=1}}{\sum}}\varrho_{\text{aa}}^{M}%
+\overset{q}{\underset{\text{a,b=1}}{\sum}}R_{\text{abab}}^{M})](y_{0})$

(xi)\qquad A direct computation using the expansion formula gives:

\qquad\ $\frac{\partial^{3}\theta}{\partial\text{x}_{i}\partial\text{x}%
_{j}\partial\text{x}_{k}}(y_{0})$

$\qquad=-\frac{1}{12}\{\nabla_{i}\varrho_{jk}-2\varrho_{ij}%
<H,k>+\overset{\text{q}}{\underset{\text{a=1}}{\sum}}(\nabla_{i}%
R_{\text{a}j\text{a}k}-4R_{i\text{a}j\text{a}}<H,k>)$

$\qquad+4\overset{\text{q}}{\underset{\text{a,b=1}}{\sum}}R_{i\text{a}%
j\text{b}}T_{\text{ab}k}$

$+2\overset{q}{\underset{\text{a,b,c=1}}{\sum}}(T_{\text{aa}i}T_{\text{bb}%
j}T_{\text{cc}k}-3T_{\text{aa}i}T_{\text{bc}j}T_{\text{bc}k}+2T_{\text{ab}%
i}T_{\text{bc}j}T_{\text{ca}k})\}(y_{0})$\qquad\qquad\qquad\qquad\qquad\ \ 

$-\frac{1}{12}\{\nabla_{j}\varrho_{ik}-2\varrho_{ji}%
<H,k>+\overset{q}{\underset{\text{a}=1}{\sum}}(\nabla_{j}R_{\text{a}%
i\text{a}k}-4R_{j\text{a}i\text{a}}<H,k>)+4\overset{q}{\underset{\text{a,b=1}%
}{\sum}}R_{j\text{a}i\text{b}}T_{\text{ab}k}$

$+2\overset{q}{\underset{\text{a,b,c=1}}{\sum}}(T_{\text{aa}j}T_{\text{bb}%
i}T_{\text{cc}k}-3T_{\text{aa}j}T_{\text{bc}i}T_{\text{bc}k}+2T_{\text{ab}%
j}T_{\text{bc}i}T_{\text{ca}k})\}(y_{0})$

$-\frac{1}{12}\{\nabla_{i}\varrho_{kj}-2\varrho_{ik}%
<H,j>+\overset{q}{\underset{\text{a}=1}{\sum}}(\nabla_{i}R_{\text{a}%
k\text{a}j}-4R_{i\text{a}k\text{a}}<H,j>)$

$+4\overset{q}{\underset{\text{a,b=1}}{\sum}}R_{i\text{a}k\text{b}%
}T_{\text{ab}j}+2\overset{q}{\underset{\text{a,b,c=1}}{\sum}}(T_{\text{aa}%
i}T_{\text{bb}k}T_{\text{cc}j}-3T_{\text{aa}i}T_{\text{bc}k}T_{\text{bc}%
j}+2T_{\text{ab}i}T_{\text{bc}k}T_{\text{ca}j})\}(y_{0})$

$-\frac{1}{12}\{\nabla_{j}\varrho_{ki}-2\varrho_{jk}%
<H,i>+\overset{q}{\underset{\text{a}=1}{\sum}}(\nabla_{j}R_{\text{a}%
k\text{a}i}-4R_{j\text{a}k\text{a}}<H,i>)$

$+4\overset{q}{\underset{\text{a,b=1}}{\sum}}R_{j\text{a}j\text{b}%
}T_{\text{ab}i}+2\overset{q}{\underset{\text{a,b,c=1}}{\sum}}(T_{\text{aa}%
j}T_{\text{bb}k}T_{\text{cc}i}-3T_{\text{aa}j}T_{\text{bc}k}T_{\text{bc}%
i}+2T_{\text{ab}j}T_{\text{bc}k}T_{\text{ca}i})\}(y_{0})$\qquad\qquad
\qquad\qquad$\ $

$-\frac{1}{12}\{\nabla_{k}\varrho_{ij}-2\varrho_{ki}%
<H,j>+\overset{q}{\underset{\text{a}=1}{\sum}}(\nabla_{k}R_{\text{a}%
i\text{a}j}-4R_{k\text{a}i\text{a}}<H,j>)$

$+4\overset{q}{\underset{\text{a,b=1}}{\sum}}R_{k\text{a}i\text{b}%
}T_{\text{ab}j}+2\overset{q}{\underset{\text{a,b,c=1}}{\sum}}(T_{\text{aa}%
k}T_{\text{bb}i}T_{\text{cc}j}-3T_{\text{aa}k}T_{\text{bc}i}T_{\text{bc}%
j}+2T_{\text{ab}k}T_{\text{bc}i}T_{\text{ca}j})\}(y_{0})$\qquad

$-\frac{1}{12}\{\nabla_{k}\varrho_{ji}-2\varrho_{kj}%
<H,i>+\overset{q}{\underset{\text{a}=1}{\sum}}(\nabla_{k}R_{\text{a}%
j\text{a}i}-4R_{k\text{a}j\text{a}}<H,i>)$

$+4\overset{q}{\underset{\text{a,b=1}}{\sum}}R_{k\text{a}j\text{b}%
}T_{\text{ab}i}$+2$\overset{q}{\underset{\text{a,b,c=1}}{\sum}}(T_{\text{aa}%
k}T_{\text{bb}j}T_{\text{cc}i}-3T_{\text{aa}k}T_{\text{bc}j}T_{\text{bc}%
i}+2T_{\text{ab}k}T_{\text{bc}j}T_{\text{ca}i})\}(y_{0})$

\vspace{1pt}

\vspace{1pt}(xii) $\frac{\partial^{3}\theta}{\partial\text{x}_{i}^{2}%
\partial\text{x}_{j}}(y_{0})$

$=-\frac{1}{6}[\nabla_{i}\varrho_{ij}-2\varrho_{ij}%
<H,i>+\overset{q}{\underset{\text{a}=1}{\sum}}(\nabla_{i}R_{\text{a}%
i\text{a}j}-4R_{i\text{a}j\text{a}}<H,i>)+4\overset{q}{\underset{\text{a,b=1}%
}{\sum}}R_{i\text{a}j\text{b}}T_{\text{ab}i}$

$+2\overset{q}{\underset{\text{a,b,c=1}}{\sum}}(T_{\text{aa}i}T_{\text{bb}%
j}T_{\text{cc}i}-3T_{\text{aa}i}T_{\text{bc}j}T_{\text{bc}i}+2T_{\text{ab}%
i}T_{\text{bc}j}T_{\text{ac}i})](y_{0})$\qquad\qquad\qquad\qquad\qquad\ \ 

$-\frac{1}{6}[\nabla_{j}\varrho_{ii}-2\varrho_{ij}%
<H,i>+\overset{q}{\underset{\text{a}=1}{\sum}}(\nabla_{j}R_{\text{a}%
i\text{a}i}-4R_{i\text{a}j\text{a}}<H,i>)+4\overset{q}{\underset{\text{a,b=1}%
}{\sum}}R_{j\text{a}i\text{b}}T_{\text{ab}i}$

$+2\overset{q}{\underset{\text{a,b,c=1}}{\sum}}(T_{\text{aa}j}T_{\text{bb}%
i}T_{\text{cc}i}-3T_{\text{aa}j}T_{\text{bc}i}T_{\text{bc}i}+2T_{\text{ab}%
j}T_{\text{bc}i}T_{\text{ac}i})](y_{0})$

$-\frac{1}{6}[\nabla_{i}\varrho_{ij}-2\varrho_{ii}%
<H,j>+\overset{q}{\underset{\text{a}=1}{\sum}}(\nabla_{i}R_{\text{a}%
i\text{a}j}-4R_{i\text{a}i\text{a}}<H,j>)+4\overset{q}{\underset{\text{a,b=1}%
}{\sum}}R_{i\text{a}i\text{b}}T_{\text{ab}j}$

+2$\overset{q}{\underset{\text{a,b,c}=1}{\sum}}(T_{\text{aa}i}T_{\text{bb}%
i}T_{\text{cc}j}-3T_{\text{aa}i}T_{\text{bc}i}T_{\text{bc}j}+2T_{\text{ab}%
i}T_{\text{bc}i}T_{\text{ac}j})](y_{0})$\qquad\qquad\qquad\qquad$\ $\qquad

\vspace{1pt}(xiii) $\frac{\partial^{3}\theta}{\partial\text{x}_{i}%
\partial\text{x}_{j}^{2}}(y_{0})$

$=-\frac{1}{6}[\nabla_{i}\varrho_{jj}-2\varrho_{ij}%
<H,j>+\overset{q}{\underset{\text{a}=1}{\sum}}(\nabla_{i}R_{\text{a}%
j\text{a}j}-4R_{i\text{a}j\text{a}}<H,j>)+4\overset{q}{\underset{\text{a,b=1}%
}{\sum}}R_{i\text{a}j\text{b}}T_{\text{ab}j}$

$+2\overset{q}{\underset{\text{a,b,c=1}}{\sum}}(T_{\text{aa}i}T_{\text{bb}%
j}T_{\text{cc}j}-3T_{\text{aa}i}T_{\text{bc}j}T_{\text{bc}j}+2T_{\text{ab}%
i}T_{\text{bc}j}T_{\text{ca}j})](y_{0})$\qquad\qquad\qquad\qquad\qquad\ \ 

$-\frac{1}{6}[\nabla_{j}\varrho_{ij}-2\varrho_{ij}%
<H,j>+\overset{q}{\underset{\text{a}=1}{\sum}}(\nabla_{j}R_{\text{a}%
i\text{a}j}-4R_{j\text{a}i\text{a}}<H,j>)+4\overset{q}{\underset{\text{a,b=1}%
}{\sum}}R_{j\text{a}i\text{b}}T_{\text{ab}j}$

$+2\overset{q}{\underset{\text{a,b,c=1}}{\sum}}(T_{\text{aa}j}T_{\text{bb}%
i}T_{\text{cc}j}-3T_{\text{aa}j}T_{\text{bc}i}T_{\text{bc}j}+2T_{\text{ab}%
j}T_{\text{bc}i}T_{\text{ac}j})](y_{0})$

$-\frac{1}{6}[\nabla_{j}\varrho_{ij}-2\varrho_{jj}%
<H,i>+\overset{q}{\underset{\text{a}=1}{\sum}}(\nabla_{j}R_{\text{a}%
i\text{a}j}-4R_{j\text{a}j\text{a}}<H,i>)+4\overset{q}{\underset{\text{a,b=1}%
}{\sum}}R_{j\text{a}j\text{b}}T_{\text{ab}i}$

$+2\overset{q}{\underset{\text{a,b,c=1}}{\sum}}(T_{\text{aa}j}T_{\text{bb}%
j}T_{\text{cc}i}-3T_{\text{aa}j}T_{\text{bc}j}T_{\text{bc}i}+2T_{\text{ab}%
j}T_{\text{bc}j}T_{\text{ac}i})](y_{0})$\qquad\qquad\qquad\qquad$\ $\qquad

\vspace{1pt}(xiv) $\qquad\frac{\partial^{3}\theta^{\frac{1}{2}}}%
{\partial\text{x}_{i}\partial\text{x}_{j}\partial\text{x}_{k}}(y_{0})=\frac
{1}{8}(\frac{\partial\theta}{\partial\text{x}_{i}}\frac{\partial\theta
}{\partial\text{x}_{j}}\frac{\partial\theta}{\partial\text{x}_{k}}%
)(y_{0})-\frac{1}{4}(\frac{\partial\theta}{\partial\text{x}_{i}}\frac
{\partial^{2}\theta}{\partial\text{x}_{j}\partial\text{x}_{k}})(y_{0})$

$\qquad\qquad+\frac{1}{8}(\frac{\partial\theta}{\partial\text{x}_{j}}%
\frac{\partial\theta}{\partial\text{x}_{i}}\frac{\partial\theta}%
{\partial\text{x}_{k}})(y_{0})-\frac{1}{4}(\frac{\partial\theta}%
{\partial\text{x}_{j}}\frac{\partial^{2}\theta}{\partial\text{x}_{i}%
\partial\text{x}_{k}})(y_{0})$

\qquad\qquad\ $+\frac{1}{8}(\frac{\partial\theta}{\partial\text{x}_{k}}%
\frac{\partial\theta}{\partial\text{x}_{i}}\frac{\partial\theta}%
{\partial\text{x}_{j}})(y_{0})-\frac{1}{4}(\frac{\partial\theta}%
{\partial\text{x}_{k}}\frac{\partial^{2}\theta}{\partial\text{x}_{i}%
\partial\text{x}_{j}})(y_{0})+\frac{1}{2}\frac{\partial^{3}\theta}%
{\partial\text{x}_{i}\partial\text{x}_{j}\partial x_{k}}(y_{0})$

$\qquad=\frac{3}{8}(\frac{\partial\theta}{\partial\text{x}_{i}}\frac
{\partial\theta}{\partial\text{x}_{j}}\frac{\partial\theta}{\partial
\text{x}_{k}})(y_{0})-\frac{1}{4}(\frac{\partial\theta}{\partial\text{x}_{i}%
}\frac{\partial^{2}\theta}{\partial\text{x}_{j}\partial\text{x}_{k}}%
)(y_{0})-\frac{1}{4}(\frac{\partial\theta}{\partial\text{x}_{j}}\frac
{\partial^{2}\theta}{\partial\text{x}_{i}\partial\text{x}_{k}})(y_{0})$

$\qquad-\frac{1}{4}(\frac{\partial\theta}{\partial\text{x}_{k}}\frac
{\partial^{2}\theta}{\partial\text{x}_{i}\partial\text{x}_{j}})(y_{0}%
)+\frac{1}{2}\frac{\partial^{3}\theta}{\partial\text{x}_{i}\partial
\text{x}_{j}\partial x_{k}}(y_{0})$

We use the expressions already computed above.

\vspace{1pt}

\vspace{1pt}(xv) $\frac{\partial^{3}\theta^{-\frac{1}{2}}}{\partial
\text{x}_{i}\partial\text{x}_{j}\partial\text{x}_{k}}(y_{0})=-\frac{15}%
{8}(\frac{\partial\theta}{\partial x_{i}}\frac{\partial\theta}{\partial x_{j}%
}\frac{\partial\theta}{\partial x_{k}})(y_{0})+\frac{3}{4}\frac{\partial
^{2}\theta}{\partial\text{x}_{i}\partial\text{x}_{j}}(y_{0})\frac
{\partial\theta}{\partial x_{k}}(y_{0})$

\qquad\qquad$+\frac{3}{4}\frac{\partial\theta}{\partial x_{i}}(y_{0}%
)\frac{\partial^{2}\theta}{\partial\text{x}_{j}\partial\text{x}_{k}}(y_{0})$
$+\frac{3}{4}\frac{\partial\theta}{\partial x_{j}}(y_{0})\frac{\partial
^{2}\theta}{\partial\text{x}_{i}\partial\text{x}_{k}}(y_{0})-$ $\frac{1}%
{2}\frac{\partial^{3}\theta}{\partial\text{x}_{i}\partial\text{x}_{j}\partial
x_{k}}(y_{0})$

We use the expressions already computed above.\qquad

(xvi) $\frac{\partial^{3}\theta^{-\frac{1}{2}}}{\partial\text{x}_{i}%
^{2}\partial\text{x}_{j}}(y_{0})=-\frac{15}{8}(\frac{\partial\theta}{\partial
x_{i}})^{2}\frac{\partial\theta}{\partial x_{j}})(y_{0})+\frac{3}{2}%
\frac{\partial\theta}{\partial x_{i}}(y_{0})\frac{\partial^{2}\theta}%
{\partial\text{x}_{i}\partial\text{x}_{j}}(y_{0})$

\qquad\qquad\qquad\qquad\qquad\ $\ +\frac{3}{4}\frac{\partial\theta}{\partial
x_{j}}(y_{0})\frac{\partial^{2}\theta}{\partial\text{x}_{i}^{2}}(y_{0})-$
$\frac{1}{2}\frac{\partial^{3}\theta}{\partial\text{x}_{i}^{2}\partial
\text{x}_{j}}(y_{0})$\qquad

(xvii) $\frac{\partial^{3}\theta^{-\frac{1}{2}}}{\partial\text{x}_{i}%
\partial\text{x}_{j}^{2}}(y_{0})=-\frac{15}{8}(\frac{\partial\theta}{\partial
x_{i}}(\frac{\partial\theta}{\partial x_{j}})^{2})(y_{0})+\frac{3}{2}%
\frac{\partial\theta}{\partial x_{j}}(y_{0})\frac{\partial^{2}\theta}%
{\partial\text{x}_{i}\partial\text{x}_{j}}(y_{0})$

\qquad\qquad\qquad\qquad\qquad\ $+\frac{3}{4}\frac{\partial\theta}{\partial
x_{i}}(y_{0})\frac{\partial^{2}\theta}{\partial\text{x}_{j}^{2}}(y_{0})-$
$\frac{1}{2}\frac{\partial^{3}\theta}{\partial\text{x}_{i}\partial\text{x}%
_{j}^{2}}(y_{0})$

\section{Table A$_{10}$}

(i)$\qquad<\nabla\theta,\nabla f>(y_{0})$ $=$ $<\nabla\log\theta,\nabla
f>(y_{0})=-$ $\overset{n}{\underset{i=q+1}{\sum}}<H,i>(y_{0})\frac{\partial
f}{\partial\text{x}_{i}}(y_{0})$

(ii)$\qquad\frac{1}{2}\Delta\theta^{-\frac{1}{2}}(y_{0})=\frac{1}%
{24}[\underset{i=q+1}{\overset{n}{\sum}}3<H,i>^{2}+2(\tau^{M}-3\tau
^{P}\ +\overset{q}{\underset{\text{a=1}}{\sum}}\varrho_{\text{aa}}%
^{M}+\overset{q}{\underset{\text{a,b}=1}{\sum}}R_{\text{abab}}^{M})](y_{0})$

(iii)$\qquad\frac{1}{4}\Delta\theta^{-\frac{1}{2}}(y_{0})(\frac{1}{2}%
\Delta\theta^{-\frac{1}{2}})(y_{0})=\frac{1}{2}(\frac{1}{2}\Delta
\theta^{-\frac{1}{2}})^{2}(y_{0})$

\qquad\qquad$\ =\frac{1}{1152}[\underset{i=q+1}{\overset{n}{\sum}}%
3<H,i>^{2}+2(\tau^{M}-3\tau^{P}\ +\overset{q}{\underset{\text{a=1}}{\sum}%
}\varrho_{\text{aa}}^{M}+\overset{q}{\underset{\text{a,b}=1}{\sum}%
}R_{\text{abab}}^{M})]^{2}(y_{0})$

(iv)$\qquad\frac{1}{48}(\Delta\theta^{-\frac{1}{2}})(y_{0})[6\Delta
\theta^{-\frac{1}{2}}+2\frac{\partial^{2}\theta^{\frac{1}{2}}}{\partial
x_{i}^{2}}+\frac{\partial\theta^{\frac{1}{2}}}{\partial x_{i}}(6\frac
{\partial\theta^{-\frac{1}{2}}}{\partial x_{i}}+4\frac{\partial\theta
^{\frac{1}{2}}}{\partial x_{i}}-3\overset{q}{\underset{\text{a}=1}{\sum}%
}\Gamma_{\text{aa}}^{i})](y_{0})\qquad$

$\qquad\qquad=\frac{1}{1728}[3<H,i>^{2}\ +2(\tau^{M}-3\tau^{P}%
\ +\overset{q}{\underset{\text{a=1}}{\sum}}\varrho_{\text{aa}}^{M}%
+\overset{q}{\underset{\text{a,b=1}}{\sum}}R_{\text{abab}}^{M})]^{2}(y_{0})$\ 

\qquad$+\frac{1}{576}[<H,i>^{2}](y_{0})[\underset{i=q+1}{\overset{n}{\sum}%
}3<H,i>^{2}+2(\tau^{M}-3\tau^{P}\ +\overset{q}{\underset{\text{a=1}}{\sum}%
}\varrho_{\text{aa}}^{M}+\overset{q}{\underset{\text{a,b}=1}{\sum}%
}R_{\text{abab}}^{M})](y_{0})$

(v)\qquad$\frac{1}{4}\frac{\partial^{2}\theta}{\partial\text{x}_{i}^{2}}%
(y_{0})\times\frac{1}{8}\frac{\partial^{2}\theta}{\partial\text{x}_{j}^{2}%
}(y_{0})$

$\qquad=\frac{1}{24}\times\frac{1}{48}[2\varrho_{ii}%
+4\overset{q}{\underset{\text{a}=1}{\sum}}R_{i\text{a}i\text{a}}%
-6\overset{q}{\underset{\text{a,b=}1}{\sum}}(T_{\text{aa}i}T_{\text{bb}%
i}-T_{\text{ab}i}T_{\text{ab}i})](y_{0})$ $\ \ \ ($l$)$

$\qquad\times\lbrack2\varrho_{jj}+4\overset{q}{\underset{\text{a}=1}{\sum}%
}R_{j\text{a}j\text{a}}-6\overset{q}{\underset{\text{a,b}=1}{\sum}%
}(T_{\text{aa}j}T_{\text{bb}j}-T_{\text{ab}j}T_{\text{ab}j})](y_{0}%
)\qquad\qquad\qquad\ \ \ \ $

$\qquad=\frac{1}{288}[\tau^{M}-3\tau^{P}\ +\overset{q}{\underset{\text{a=1}%
}{\sum}}\varrho_{\text{aa}}^{M}+\overset{q}{\underset{\text{a,b}=1}{\sum}%
}R_{\text{abab}}^{M}]^{2}(y_{0})$

(vi)$\qquad$A$_{3211}=\frac{1}{24}(\frac{\partial^{2}\theta^{\frac{1}{2}}%
}{\partial x_{i}^{2}}\Delta\theta^{-\frac{1}{2}})(y_{0})$

\qquad$=-\frac{1}{3456}[3<H,i>^{2}\ +2(\tau^{M}-3\tau^{P}%
\ +\overset{q}{\underset{\text{a=1}}{\sum}}\varrho_{\text{aa}}^{M}%
+\overset{q}{\underset{\text{a,b}=1}{\sum}}R_{\text{abab}}^{M})]^{2}(y_{0})$

(vii)$\qquad$A$_{3212}=$ $\frac{1}{24}\left[  \frac{\partial^{2}}%
{\partial\text{x}_{i}^{2}}(\Delta\theta^{-\frac{1}{2}})\right]  (y_{0}%
)=\frac{1}{24}(L_{1}+L_{2}+$ $L_{3})$

$=\frac{1}{24}[2<H,i>^{2}(y_{0})+\frac{1}{3}(\tau^{M}-3\tau^{P}%
+\overset{q}{\underset{\text{a}=1}{\sum}}\varrho_{\text{aa}}%
+\overset{q}{\underset{\text{a,b}=1}{\sum}}R_{\text{abab}})](y_{0})\qquad
L_{1}$

$\times\lbrack\frac{1}{4}<H,j>^{2}(y_{0})+\frac{1}{6}(\tau^{M}-3\tau
^{P}+\overset{q}{\underset{\text{a}=1}{\sum}}\varrho_{\text{aa}}%
^{M}+\overset{q}{\underset{\text{a,b}=1}{\sum}}R_{\text{abab}}^{M})](y_{0})$

$-\frac{1}{24}\times\frac{1}{8}[<H,i><H,j>](y_{0}\qquad\qquad\qquad\qquad
L_{211}\qquad L_{21}$ \qquad$L_{2}$

$\times\lbrack2\varrho_{ij}+$ $\overset{q}{\underset{\text{a}=1}{4\sum}%
}R_{i\text{a}j\text{a}}-3\overset{q}{\underset{\text{a,b=1}}{\sum}%
}(T_{\text{aa}i}T_{\text{bb}j}-T_{\text{ab}i}T_{\text{ab}j}%
)-3\overset{q}{\underset{\text{a,b=1}}{\sum}}(T_{\text{aa}j}T_{\text{bb}%
i}-T_{\text{ab}j}T_{\text{ab}i}](y_{0})\qquad$

$-\frac{1}{24}\times\frac{1}{72}[2\varrho_{ij}+$
$\overset{q}{\underset{\text{a}=1}{4\sum}}R_{i\text{a}j\text{a}}%
-3\overset{q}{\underset{\text{a,b=1}}{\sum}}(T_{\text{aa}i}T_{\text{bb}%
j}-T_{\text{ab}i}T_{\text{ab}j})-3\overset{q}{\underset{\text{a,b=1}}{\sum}%
}(T_{\text{aa}j}T_{\text{bb}i}-T_{\text{ab}j}T_{\text{ab}i}]^{2}(y_{0})$

$-\frac{1}{24}\times\frac{1}{12}[<H,j>](y_{0})\times\lbrack\{\nabla_{i}%
\varrho_{ij}-2\varrho_{ij}<H,i>+\overset{q}{\underset{\text{a}=1}{\sum}%
}(\nabla_{i}R_{\text{a}i\text{a}j}-4R_{i\text{a}j\text{a}}<H,i>)\qquad
L_{212}$

$+4\overset{q}{\underset{\text{a,b=1}}{\sum}}R_{i\text{a}j\text{b}%
}T_{\text{ab}i}+2\overset{q}{\underset{\text{a,b,c=1}}{\sum}}(T_{\text{aa}%
i}T_{\text{bb}j}T_{\text{cc}i}-3T_{\text{aa}i}T_{\text{bc}j}T_{\text{bc}%
i}+2T_{\text{ab}i}T_{\text{bc}j}T_{\text{ac}i})](y_{0})$\qquad\qquad
\qquad\qquad\qquad\ \ 

$-\frac{1}{24}\times\frac{1}{12}[<H,j>](y_{0})\times\lbrack\nabla_{j}%
\varrho_{ii}-2\varrho_{ij}<H,i>+\overset{q}{\underset{\text{a}=1}{\sum}%
}(\nabla_{j}R_{\text{a}i\text{a}i}-4R_{i\text{a}j\text{a}}<H,i>)$

$+4\overset{q}{\underset{\text{a,b=1}}{\sum}}R_{j\text{a}i\text{b}%
}T_{\text{ab}i}+2\overset{q}{\underset{\text{a,b,c=1}}{\sum}}(T_{\text{aa}%
j}T_{\text{bb}i}T_{\text{cc}i}-3T_{\text{aa}j}T_{\text{bc}i}T_{\text{bc}%
i}+2T_{\text{ab}j}T_{\text{bc}i}T_{\text{ac}i})](y_{0})$

$-\frac{1}{24}\times\frac{1}{12}[<H,j>](y_{0})\times\lbrack\nabla_{i}%
\varrho_{ij}-2\varrho_{ii}<H,j>+\overset{q}{\underset{\text{a}=1}{\sum}%
}(\nabla_{i}R_{\text{a}i\text{a}j}-4R_{i\text{a}i\text{a}}<H,j>)$

$+4\overset{q}{\underset{\text{a,b=1}}{\sum}}R_{i\text{a}i\text{b}%
}T_{\text{ab}j}+2\overset{q}{\underset{\text{a,b,c}=1}{\sum}}(T_{\text{aa}%
i}T_{\text{bb}i}T_{\text{cc}j}-3T_{\text{aa}i}T_{\text{bc}i}T_{\text{bc}%
j}+2T_{\text{ab}i}T_{\text{bc}i}T_{\text{ac}j})](y_{0})$

$-\frac{1}{24}\times\frac{1}{3}[<H,j><H,k>](y_{0})R_{ijik}(y_{0})\qquad\qquad
L_{213}$

$-\frac{1}{24}$ $\times\frac{15}{8}[<H,i>^{2}<H,j>^{2}](y_{0})\qquad
\qquad\qquad$

$-\frac{1}{24}\times\frac{1}{4}<H,i><H,j>$

$\times\lbrack2\varrho_{ij}+\overset{q}{\underset{\text{a}=1}{4\sum}%
}R_{i\text{a}j\text{a}}-3\overset{q}{\underset{\text{a,b=1}}{\sum}%
}(T_{\text{aa}i}T_{\text{bb}j}-T_{\text{ab}i}T_{\text{ab}j}%
)-3\overset{q}{\underset{\text{a,b=1}}{\sum}}(T_{\text{aa}j}T_{\text{bb}%
i}-T_{\text{ab}j}T_{\text{ab}i}](y_{0})$

$-\frac{1}{24}\times\frac{1}{4}<H,j>^{2}[\tau^{M}\ -3\tau^{P}%
+\ \underset{\text{a}=1}{\overset{\text{q}}{\sum}}\varrho_{\text{aa}}^{M}+$
$\overset{q}{\underset{\text{a},\text{b}=1}{\sum}}R_{\text{abab}}^{M}$
$](y_{0})$

$+\frac{1}{24}\times\frac{1}{12}<H,j>[\nabla_{i}\varrho_{ij}-2\varrho
_{ij}<H,i>+\overset{q}{\underset{\text{a}=1}{\sum}}(\nabla_{i}R_{\text{a}%
i\text{a}j}-4R_{i\text{a}j\text{a}}<H,i>)$

$+4\overset{q}{\underset{\text{a,b=1}}{\sum}}R_{i\text{a}j\text{b}%
}T_{\text{ab}i}+2\overset{q}{\underset{\text{a,b,c=1}}{\sum}}(T_{\text{aa}%
i}T_{\text{bb}j}T_{\text{cc}i}-3T_{\text{aa}i}T_{\text{bc}j}T_{\text{bc}%
i}+2T_{\text{ab}i}T_{\text{bc}j}T_{\text{ac}i})](y_{0})$\qquad\qquad
\qquad\qquad\qquad\ \ 

$+\frac{1}{24}\times\frac{1}{12}<H,j>[\nabla_{j}\varrho_{ii}-2\varrho
_{ij}<H,i>+\overset{q}{\underset{\text{a}=1}{\sum}}(\nabla_{j}R_{\text{a}%
i\text{a}i}-4R_{i\text{a}j\text{a}}<H,i>)$

$+4\overset{q}{\underset{\text{a,b=1}}{\sum}}R_{j\text{a}i\text{b}%
}T_{\text{ab}i}+2\overset{q}{\underset{\text{a,b,c=1}}{\sum}}(T_{\text{aa}%
j}T_{\text{bb}i}T_{\text{cc}i}-3T_{\text{aa}j}T_{\text{bc}i}T_{\text{bc}%
i}+2T_{\text{ab}j}T_{\text{bc}i}T_{\text{ac}i})](y_{0})$

$+\frac{1}{24}\times\frac{1}{12}<H,j>[\nabla_{i}\varrho_{ij}-2\varrho
_{ii}<H,j>+\overset{q}{\underset{\text{a}=1}{\sum}}(\nabla_{i}R_{\text{a}%
i\text{a}j}-4R_{i\text{a}i\text{a}}<H,j>)$

$+4\overset{q}{\underset{\text{a,b=1}}{\sum}}R_{i\text{a}i\text{b}%
}T_{\text{ab}j}+2\overset{q}{\underset{\text{a,b,c}=1}{\sum}}(T_{\text{aa}%
i}T_{\text{bb}i}T_{\text{cc}j}-3T_{\text{aa}i}T_{\text{bc}i}T_{\text{bc}%
j}+2T_{\text{ab}i}T_{\text{bc}i}T_{\text{ac}j})](y_{0})$

$+\frac{1}{24}\times\frac{1}{2}R_{ijjk}(y_{0})$ \ $[<H,i><H,k>](y_{0}%
)\qquad\qquad\qquad\qquad\qquad\qquad L_{22}$

$+\frac{1}{24}\times\frac{1}{18}R_{ijjk}(y_{0})$

$\times\lbrack2\varrho_{ik}+$ $\overset{q}{\underset{\text{a}=1}{4\sum}%
}R_{i\text{a}k\text{a}}-3\overset{q}{\underset{\text{a,b=1}}{\sum}%
}(T_{\text{aa}i}T_{\text{bb}k}-T_{\text{ab}i}T_{\text{ab}k}%
)-3\overset{q}{\underset{\text{a,b=1}}{\sum}}(T_{\text{aa}k}T_{\text{bb}%
i}-T_{\text{ab}k}T_{\text{ab}i}](y_{0})$

$\ +\frac{1}{24}\times\frac{1}{6}<H,k>(y_{0})[\nabla_{j}$R$_{ikij}%
(y_{0})+\nabla_{i}$R$_{jkij}](y_{0})$

$+\frac{1}{24}\times\frac{1}{2}R_{ijik}(y_{0})$ \ $[<H,j><H,k>](y_{0}%
)\qquad\qquad\qquad\qquad\qquad\qquad$ $L_{23}$

$+\frac{1}{24}\times\frac{1}{18}R_{ijik}(y_{0})$

$\times\lbrack2\varrho_{jk}+$ $\overset{q}{\underset{\text{a}=1}{4\sum}%
}R_{j\text{a}k\text{a}}-3\overset{q}{\underset{\text{a,b=1}}{\sum}%
}(T_{\text{aa}j}T_{\text{bb}k}-T_{\text{ab}j}T_{\text{ab}k}%
)-3\overset{q}{\underset{\text{a,b=1}}{\sum}}(T_{\text{aa}k}T_{\text{bb}%
j}-T_{\text{ab}k}T_{\text{ab}j}](y_{0})$

$+\overset{n}{\underset{i,j=q+1}{\sum}}\frac{35}{128}<H,i>^{2}(y_{0}%
)<H,j>^{2}(y_{0})\qquad\qquad\ \frac{1}{24}\frac{\partial^{4}\theta^{-\frac
{1}{2}}}{\partial x_{i}^{2}\partial x_{j}^{2}}(y_{0})$

$+\frac{5}{192}\overset{n}{\underset{j=q+1}{\sum}}<H,j>^{2}(y_{0})[\tau
^{M}\ -3\tau^{P}+\ \underset{\text{a}=1}{\overset{\text{q}}{\sum}}%
\varrho_{\text{aa}}^{M}+\overset{q}{\underset{\text{a},\text{b}=1}{\sum}%
}R_{\text{abab}}^{M}](y_{0})\qquad\ \ \ \ \ \ \ \ $

$+\frac{5}{192}\overset{n}{\underset{i=q+1}{\sum}}<H,i>^{2}(y_{0})[\tau
^{M}\ -3\tau^{P}+\ \underset{\text{a}=1}{\overset{\text{q}}{\sum}}%
\varrho_{\text{aa}}^{M}+\overset{q}{\underset{\text{a},\text{b}=1}{\sum}%
}R_{\text{abab}}^{M}](y_{0})\qquad\qquad$

$+\frac{5}{192}\overset{n}{\underset{i,j=q+1}{\sum}}[<H,i><H,j>](y_{0}%
)\qquad\qquad\qquad\qquad\qquad\qquad\qquad\qquad$

$\times\lbrack2\varrho_{ij}+4\overset{q}{\underset{\text{a}=1}{\sum}%
}R_{i\text{a}j\text{a}}-3\overset{q}{\underset{\text{a,b=1}}{\sum}%
}(T_{\text{aa}i}T_{\text{bb}j}-T_{\text{ab}i}T_{\text{ab}j}%
)-3\overset{q}{\underset{\text{a,b=1}}{\sum}}(T_{\text{aa}j}T_{\text{bb}%
i}-T_{\text{ab}j}T_{\text{ab}i})](y_{0})$

$+\frac{1}{96}\overset{n}{\underset{i,j=q+1}{\sum}}<H,j>(y_{0})[\{\nabla
_{i}\varrho_{ij}-2\varrho_{ij}<H,i>+\overset{q}{\underset{\text{a}=1}{\sum}%
}(\nabla_{i}R_{\text{a}i\text{a}j}-4R_{i\text{a}j\text{a}}<H,i>)\qquad$

$+4\overset{q}{\underset{\text{a,b=1}}{\sum}}R_{i\text{a}j\text{b}%
}T_{\text{ab}i}+2\overset{q}{\underset{\text{a,b,c=1}}{\sum}}(T_{\text{aa}%
i}T_{\text{bb}j}T_{\text{cc}i}-T_{\text{aa}i}T_{\text{bc}j}T_{\text{bc}%
i}-2T_{\text{bc}j}(T_{\text{aa}i}T_{\text{bc}i}-T_{\text{ab}i}T_{\text{ac}%
i}))\}$\qquad\qquad\qquad\ \ 

$+\{\nabla_{j}\varrho_{ii}-2\varrho_{ij}<H,i>+\overset{q}{\underset{\text{a}%
=1}{\sum}}(\nabla_{j}R_{\text{a}i\text{a}i}-4R_{i\text{a}j\text{a}}<H,i>)$

$+4\overset{q}{\underset{\text{a,b=1}}{\sum}}R_{j\text{a}i\text{b}%
}T_{\text{ab}i}+2\overset{q}{\underset{\text{a,b,c=1}}{\sum}}(T_{\text{aa}%
j}(T_{\text{bb}i}T_{\text{cc}i}-T_{\text{bc}i}T_{\text{bc}i})-2T_{\text{aa}%
j}T_{\text{bc}i}T_{\text{bc}i}+2T_{\text{ab}j}T_{\text{bc}i}T_{\text{ac}%
i})\}\qquad$

$+\{\nabla_{i}\varrho_{ij}-2\varrho_{ii}<H,j>+\overset{q}{\underset{\text{a}%
=1}{\sum}}(\nabla_{i}R_{\text{a}i\text{a}j}-4R_{i\text{a}i\text{a}%
}<H,j>)+4\overset{q}{\underset{\text{a,b=1}}{\sum}}R_{i\text{a}i\text{b}%
}T_{\text{ab}j}$

$+2\overset{q}{\underset{\text{a,b,c}=1}{\sum}}(T_{\text{aa}i}T_{\text{bb}%
i}T_{\text{cc}j}-3T_{\text{aa}i}T_{\text{bc}i}T_{\text{bc}j}+2T_{\text{ab}%
i}T_{\text{bc}i}T_{\text{ac}j})\}](y_{0})$

$+\frac{1}{96}\overset{n}{\underset{i,j=q+1}{\sum}}<H,i>(y_{0})[\{\nabla
_{i}\varrho_{jj}-2\varrho_{ij}<H,j>+\overset{q}{\underset{\text{a}=1}{\sum}%
}(\nabla_{i}R_{\text{a}j\text{a}j}-4R_{i\text{a}j\text{a}}<H,j>)\qquad$

$+4\overset{q}{\underset{\text{a,b=1}}{\sum}}R_{i\text{a}j\text{b}%
}T_{\text{ab}j}+2\overset{q}{\underset{\text{a,b,c=1}}{\sum}}T_{\text{aa}%
i}(T_{\text{bb}j}T_{\text{cc}j}-T_{\text{bc}j}T_{\text{bc}j})-2T_{\text{aa}%
i}T_{\text{bc}j}T_{\text{bc}j}+2T_{\text{ab}i}T_{\text{bc}j}T_{\text{ac}%
j})\}(y_{0})\qquad$\qquad\qquad\qquad\qquad\qquad\ \ 

$+\{\nabla_{j}\varrho_{ij}-2\varrho_{ij}<H,j>+\overset{q}{\underset{\text{a}%
=1}{\sum}}(\nabla_{j}R_{\text{a}i\text{a}j}-4R_{j\text{a}i\text{a}}<H,j>)$

$+4\overset{q}{\underset{\text{a,b=1}}{\sum}}R_{j\text{a}i\text{b}%
}T_{\text{ab}j}+2\overset{q}{\underset{\text{a,b,c=1}}{\sum}}(T_{\text{aa}%
j}T_{\text{bb}i}T_{\text{cc}j}-T_{\text{ab}j}T_{\text{bc}i}T_{\text{ac}%
j}-2T_{\text{bc}i}(T_{\text{aa}j}T_{\text{bc}j}-T_{\text{ab}j}T_{\text{ac}%
j}))\}(y_{0})$

$+\{\nabla_{j}\varrho_{ij}-2\varrho_{jj}<H,i>+\overset{q}{\underset{\text{a}%
=1}{\sum}}(\nabla_{j}R_{\text{a}i\text{a}j}-4R_{j\text{a}j\text{a}%
}<H,i>)+4\overset{q}{\underset{\text{a,b=1}}{\sum}}R_{j\text{a}j\text{b}%
}T_{\text{ab}i}$

$+2\overset{q}{\underset{\text{a,b,c=1}}{\sum}}(T_{\text{aa}j}T_{\text{bb}%
j}T_{\text{cc}i}-3T_{\text{aa}j}T_{\text{bc}j}T_{\text{bc}i}+2T_{\text{ab}%
j}T_{\text{bc}j}T_{\text{ac}i})\}](y_{0})$

$+\frac{1}{576}\overset{n}{\underset{i,j=q+1}{\sum}}[2\varrho_{ij}%
+4\overset{q}{\underset{\text{a}=1}{\sum}}R_{i\text{a}j\text{a}}%
-3\overset{q}{\underset{\text{a,b=1}}{\sum}}(T_{\text{aa}i}T_{\text{bb}%
j}-T_{\text{ab}i}T_{\text{ab}j})-3\overset{q}{\underset{\text{a,b=1}}{\sum}%
}(T_{\text{aa}j}T_{\text{bb}i}-T_{\text{ab}j}T_{\text{ab}i})]^{2}(y_{0})$

$+\frac{1}{288}[\tau^{M}\ -3\tau^{P}+\ \underset{\text{a}=1}{\overset{\text{q}%
}{\sum}}\varrho_{\text{aa}}^{M}+\overset{q}{\underset{\text{a},\text{b}%
=1}{\sum}}R_{\text{abab}}^{M}]^{2}(y_{0})$

$-\ \frac{1}{288}\overset{n}{\underset{i,j=q+1}{\sum}}[$
$\overset{q}{\underset{\text{a=1}}{\sum}}\{-(\nabla_{ii}^{2}R_{j\text{a}%
j\text{a}}+\nabla_{jj}^{2}R_{i\text{a}i\text{a}}+4\nabla_{ij}^{2}%
R_{i\text{a}j\text{a}}+2R_{ij}R_{i\text{a}j\text{a}})\qquad A$

$+\overset{n}{\underset{p=q+1}{\sum}}\overset{q}{\underset{\text{a=1}}{\sum}%
}(R_{\text{a}iip}R_{\text{a}jjp}+R_{\text{a}jjp}R_{\text{a}iip}+R_{\text{a}%
ijp}R_{\text{a}ijp}+R_{\text{a}ijp}R_{\text{a}jip}+R_{\text{a}jip}%
R_{\text{a}ijp}+R_{\text{a}jip}R_{\text{a}jip})$

$+2\overset{q}{\underset{\text{a,b=1}}{\sum}}\nabla_{i}(R)_{\text{a}%
i\text{b}j}T_{\text{ab}j}+2\overset{q}{\underset{\text{a,b=1}}{\sum}}%
\nabla_{j}(R)_{\text{a}j\text{b}i}T_{\text{ab}i}%
+2\overset{q}{\underset{\text{a,b=1}}{\sum}}\nabla_{i}(R)_{\text{a}j\text{b}%
i}T_{\text{ab}j}+2\overset{q}{\underset{\text{a,b=1}}{\sum}}\nabla
_{i}(R)_{\text{a}j\text{b}j}T_{\text{ab}i}$

$+2\overset{q}{\underset{\text{a,b=1}}{\sum}}\nabla_{j}(R)_{\text{a}%
i\text{b}i}T_{\text{ab}j}+2\overset{q}{\underset{\text{a,b=1}}{\sum}}%
\nabla_{j}(R)_{\text{a}i\text{b}j}T_{\text{ab}i}$

$+\overset{n}{\underset{p=q+1}{\sum}}(-\frac{3}{5}\nabla_{ii}^{2}%
(R)_{jpjp}+\overset{n}{\underset{p=q+1}{\sum}}(-\frac{3}{5}\nabla_{jj}%
^{2}(R)_{ipip}$

$+\overset{n}{\underset{p=q+1}{\sum}}(-\frac{3}{5}\nabla_{ij}^{2}%
(R)_{ipjp}+\overset{n}{\underset{p=q+1}{\sum}}(-\frac{3}{5}\nabla_{ij}%
^{2}(R)_{jpip}+\overset{n}{\underset{p=q+1}{\sum}}(-\frac{3}{5}\nabla_{ji}%
^{2}(R)_{ipjp}+\overset{n}{\underset{p=q+1}{\sum}}(-\frac{3}{5}\nabla_{ji}%
^{2}(R)_{jpip}$

$+\frac{1}{5}\overset{n}{\underset{m,p=q+1}{%
{\textstyle\sum}
}}R_{ipim}R_{jpjm}+\frac{1}{5}\overset{n}{\underset{m,p=q+1}{%
{\textstyle\sum}
}}R_{jpjm}R_{ipim}+\frac{1}{5}\overset{n}{\underset{m,p=q+1}{%
{\textstyle\sum}
}}R_{ipjm}R_{ipjm}+\frac{1}{5}\overset{n}{\underset{m,p=q+1}{%
{\textstyle\sum}
}}R_{ipjm}R_{jpim}$

$+\frac{1}{5}\overset{n}{\underset{m,p=q+1}{%
{\textstyle\sum}
}}R_{jpim}R_{ipjm}+\frac{1}{5}\overset{n}{\underset{m,p=q+1}{%
{\textstyle\sum}
}}R_{jpim}R_{jpim}\}(y_{0})$

$+4\overset{q}{\underset{\text{a,b=1}}{\sum}}\{(\nabla_{i}(R)_{i\text{a}%
j\text{a}}-\overset{q}{\underset{\text{c=1}}{%
{\textstyle\sum}
}}R_{\text{a}i\text{c}i}T_{\text{ac}j})$ $T_{\text{bb}j}+4(\nabla
_{j}(R)_{j\text{a}i\text{a}}-\overset{q}{\underset{\text{c=1}}{%
{\textstyle\sum}
}}R_{\text{a}j\text{c}j}T_{\text{ac}i})$ $T_{\text{bb}i}$

$+4(\nabla_{i}(R)_{j\text{a}i\text{a}}-\overset{q}{\underset{\text{c=1}}{%
{\textstyle\sum}
}}R_{\text{a}i\text{c}j}T_{\text{ac}i})$ $T_{\text{bb}j}$ $4B\ $

$+4(\nabla_{i}(R)_{j\text{a}j\text{a}}-\overset{q}{\underset{\text{c=1}}{%
{\textstyle\sum}
}}R_{\text{a}i\text{c}j}T_{\text{ac}j})$ $T_{\text{bb}i}$

$+4(\nabla_{j}(R)_{i\text{a}i\text{a}}-\overset{q}{\underset{\text{c=1}}{%
{\textstyle\sum}
}}R_{\text{a}j\text{c}i}T_{\text{ac}i})$ $T_{\text{bb}j}+4(\nabla
_{j}(R)_{i\text{a}j\text{a}}-\overset{q}{\underset{\text{c=1}}{%
{\textstyle\sum}
}}R_{\text{a}j\text{c}i}T_{\text{ac}j})$ $T_{\text{bb}i}$

$-4\overset{q}{\underset{\text{a,b=1}}{\sum}}(\nabla_{i}(R)_{i\text{a}%
j\text{b}}-\overset{q}{\underset{\text{c=1}}{%
{\textstyle\sum}
}}R_{\text{b}r\text{c}s}T_{\text{ac}t})T_{\text{ab}j}%
-4\overset{q}{\underset{\text{a,b=1}}{\sum}}(\nabla_{j}(R)_{j\text{a}%
i\text{b}}-\overset{q}{\underset{\text{c=1}}{%
{\textstyle\sum}
}}R_{\text{b}j\text{c}j}T_{\text{ac}i})T_{\text{ab}i}$

$-4\overset{q}{\underset{\text{a,b=1}}{\sum}}(\nabla_{i}(R)_{j\text{a}%
i\text{b}}-\overset{q}{\underset{\text{c=1}}{%
{\textstyle\sum}
}}R_{\text{b}i\text{c}j}T_{\text{ac}i})T_{\text{ab}j}%
-4\overset{q}{\underset{\text{a,b=1}}{\sum}}(\nabla_{i}(R)_{j\text{a}%
j\text{b}}-\overset{q}{\underset{\text{c=1}}{%
{\textstyle\sum}
}}R_{\text{b}i\text{c}j}T_{\text{ac}j})T_{\text{ab}i}$

$-4\overset{q}{\underset{\text{a,b=1}}{\sum}}(\nabla_{j}(R)_{i\text{a}%
i\text{b}}-\overset{q}{\underset{\text{c=1}}{%
{\textstyle\sum}
}}R_{\text{b}j\text{c}i}T_{\text{ac}i})T_{\text{ab}j}%
-4\overset{q}{\underset{\text{a,b=1}}{\sum}}(\nabla_{j}(R)_{i\text{a}%
j\text{b}}-\overset{q}{\underset{\text{c=1}}{%
{\textstyle\sum}
}}R_{\text{b}j\text{c}i}T_{\text{ac}j})T_{\text{ab}i}\}](y_{0})$

$-\frac{1}{48}$ $[\frac{4}{9}\overset{q}{\underset{\text{a,b=1}}{\sum}%
}(\varrho_{\text{aa}}-\overset{q}{\underset{\text{c}=1}{\sum}}R_{\text{acac}%
})(\varrho_{\text{bb}}-\overset{q}{\underset{\text{d}=1}{\sum}}R_{\text{bdbd}%
})+\frac{8}{9}\overset{n}{\underset{i,j=q+1}{\sum}}%
\overset{q}{\underset{\text{a,b}=1}{\sum}}(R_{i\text{a}j\text{a}}%
R_{i\text{b}j\text{b}})\qquad3C$

$+\frac{2}{9}\overset{q}{\underset{\text{a}=1}{\sum}}(\varrho_{\text{aa}}%
^{M}-\varrho_{\text{aa}}^{P})(\tau^{M}-\overset{q}{\underset{\text{c}=1}{\sum
}}\varrho_{\text{cc}}^{M})+\frac{4}{9}\overset{n}{\underset{i,j=q+1}{\sum}%
}\overset{q}{\underset{\text{a}=1}{\sum}}R_{i\text{a}j\text{a}}\varrho_{ij}\ $

$\ +\frac{2}{9}\overset{q}{\underset{\text{b}=1}{\sum}}(\varrho_{\text{bb}%
}^{M}-\varrho_{\text{bb}}^{P})(\tau^{M}-\overset{q}{\underset{\text{c}%
=1}{\sum}}\varrho_{\text{cc}}^{M})+\frac{4}{9}%
\overset{n}{\underset{i,j=q+1}{\sum}}\overset{q}{\underset{\text{b}=1}{\sum}%
}R_{i\text{b}j\text{b}}\varrho_{ij}\ $

$+\frac{1}{9}(\tau^{M}-\overset{q}{\underset{\text{a=1}}{\sum}}\varrho
_{\text{aa}})(\tau^{M}-\overset{q}{\underset{\text{b=1}}{\sum}}\varrho
_{\text{bb}})+\frac{2}{9}(\left\Vert \varrho^{M}\right\Vert ^{2}%
-\overset{q}{\underset{\text{a,b}=1}{\sum}}\varrho_{\text{ab}})$

$-\overset{n}{\underset{i,j=q+1}{\sum}}\overset{q}{\underset{\text{a,b}%
=1}{\sum}}R_{i\text{a}i\text{b}}R_{j\text{a}j\text{b}}\ -\frac{1}%
{2}\overset{n}{\underset{i,j=q+1}{\sum}}\overset{q}{\underset{\text{a,b}%
=1}{\sum}}R_{i\text{a}j\text{b}}^{2}-\overset{n}{\underset{i,j=q+1}{\sum}%
}\overset{q}{\underset{\text{a,b}=1}{\sum}}R_{i\text{a}j\text{b}}%
R_{j\text{a}i\text{b}}-\frac{1}{2}\overset{n}{\underset{i,j=q+1}{\sum}%
}\overset{q}{\underset{\text{a,b}=1}{\sum}}R_{j\text{a}i\text{b}}^{2}\qquad$

$-\frac{1}{9}\overset{n}{\underset{i,j,p,m=q+1}{\sum}}R_{ipim}R_{jpjm}%
\ -\frac{1}{18}\overset{n}{\underset{i,j,p,m=q+1}{\sum}}R_{ipjm}^{2}-\frac
{1}{9}\overset{n}{\underset{i,j,p,m=q+1}{\sum}}R_{ipjm}R_{jpim}-\frac{1}%
{18}\overset{n}{\underset{i,j,p,m=q+1}{\sum}}R_{jpim}^{2}\qquad$

$-\frac{1}{3}\overset{q}{\underset{\text{a}=1}{\sum}}%
\overset{n}{\underset{i,j,p=q+1}{\sum}}R_{i\text{a}ip}R_{j\text{a}jp}-\frac
{1}{6}\overset{q}{\underset{\text{a}=1}{\sum}}%
\overset{n}{\underset{i,j,p=q+1}{\sum}}R_{i\text{a}jp}^{2}-\frac{1}%
{3}\overset{q}{\underset{\text{a}=1i,j,}{\sum}}%
\overset{n}{\underset{p=q+1}{\sum}}R_{i\text{a}jp}R_{j\text{a}ip}-\frac{1}%
{6}\overset{q}{\underset{\text{a}=1}{\sum}}%
\overset{n}{\underset{i,j,p=q+1}{\sum}}R_{j\text{a}ip}^{2}$

$-\frac{1}{3}\overset{q}{\underset{\text{b}=1i,j,}{\sum}}%
\overset{n}{\underset{p=q+1}{\sum}}R_{i\text{b}ip}R_{j\text{b}jp}-\frac{1}%
{6}\overset{q}{\underset{\text{b}=1}{\sum}}%
\overset{n}{\underset{i,j,p=q+1}{\sum}}R_{i\text{b}jp}^{2}-\frac{1}%
{3}\overset{q}{\underset{\text{b}=1}{\sum}}%
\overset{n}{\underset{i.j,p=q+1}{\sum}}R_{i\text{b}jp}R_{j\text{b}ip}-\frac
{1}{6}\overset{q}{\underset{\text{b}=1}{\sum}}%
\overset{n}{\underset{i,j,p=q+1}{\sum}}R_{j\text{b}ip}^{2}](y_{0})$

$-\frac{1}{48}$ $\overset{q}{\underset{\text{a,b,c=1}}{\sum}}[$
$-\overset{n}{\underset{i=q+1}{\sum}}R_{i\text{a}i\text{a}}(R_{\text{bcbc}%
}^{P}-R_{\text{bcbc}}^{M})$ $-\overset{n}{\underset{j=q+1}{\sum}}%
R_{j\text{a}j\text{a}}(R_{\text{bcbc}}^{P}-R_{\text{bcbc}}^{M})\qquad\qquad6D$

\ $+\overset{n}{\underset{i=q+1}{\sum}}R_{i\text{a}i\text{b}}(R_{\text{acbc}%
}^{P}-R_{\text{acbc}}^{M})\ -\overset{n}{\underset{i=q+1}{\sum}}%
R_{i\text{a}i\text{c}}(R_{\text{abbc}}^{P}-R_{\text{abbc}}^{M})$

$+\overset{n}{\underset{j=q+1}{\sum}}R_{j\text{a}j\text{b}}(R_{\text{acbc}%
}^{P}-R_{\text{acbc}}^{M})$\ $-\overset{n}{\underset{j=q+1}{\sum}}%
R_{j\text{a}j\text{c}}(R_{\text{abbc}}^{P}-R_{\text{abbc}}^{M})$

$+\underset{i,j=q+1}{\overset{n}{\sum}}$ $-R_{i\text{a}j\text{a}}%
(T_{\text{bb}i}T_{\text{cc}j}$ $-T_{\text{bc}i}T_{\text{bc}j})$
$-\underset{i,j=q+1}{\overset{n}{\sum}}R_{i\text{a}j\text{a}}(T_{\text{bb}%
j}T_{\text{cc}i}$ $-T_{\text{bc}j}T_{\text{bc}i})$

$+$ $\underset{i,j=q+1}{\overset{n}{\sum}}$ $-R_{j\text{a}i\text{a}%
}(T_{\text{bb}i}T_{\text{cc}j}$ $-T_{\text{bc}i}T_{\text{bc}j})$
$-\underset{i,j=q+1}{\overset{n}{\sum}}R_{j\text{a}i\text{a}}(T_{\text{bb}%
j}T_{\text{cc}i}$ $-T_{\text{bc}j}T_{\text{bc}i})$

$+\underset{i,j=q+1}{\overset{n}{\sum}}\ R_{i\text{a}j\text{b}}(T_{\text{ab}%
i}T_{\text{cc}j}-T_{\text{bc}i}T_{\text{ac}j}%
)\ +\underset{i,j=q+1}{\overset{n}{\sum}}\ R_{i\text{a}j\text{b}}%
(T_{\text{ab}j}T_{\text{cc}i}-T_{\text{bc}j}T_{\text{ac}i})$

$+\underset{i,j=q+1}{\overset{n}{\sum}}\ R_{j\text{a}i\text{ib}}%
(T_{\text{ab}i}T_{\text{cc}j}-T_{\text{bc}i}T_{\text{ac}j}%
)\ +\underset{i,j=q+1}{\overset{n}{\sum}}\ R_{j\text{a}i\text{b}}%
(T_{\text{ab}j}T_{\text{cc}i}-T_{\text{bc}j}T_{\text{ac}i})\qquad$

$+\underset{i,j=q+1}{\overset{n}{\sum}}-R_{i\text{a}j\text{c}}(T_{\text{ab}%
i}T_{\text{bc}j}-T_{\text{ac}i}T_{\text{bb}j}%
)-\underset{i,j=q+1}{\overset{n}{\sum}}R_{i\text{a}j\text{c}}(T_{\text{ba}%
j}T_{\text{bc}i}-T_{\text{ac}j}T_{\text{bb}i})$

$+\underset{i,j=q+1}{\overset{n}{\sum}}-R_{j\text{a}i\text{c}}(T_{\text{ba}%
i}T_{\text{bc}j}-T_{\text{ac}i}T_{\text{bb}j}%
)-\underset{i,j=q+1}{\overset{n}{\sum}}R_{j\text{a}i\text{c}}(T_{\text{ba}%
j}T_{\text{bc}i}-T_{\text{ac}j}T_{\text{bb}i})](y_{0})$

$+\frac{1}{144}\underset{p=q+1}{\overset{n}{\sum}}%
[\underset{i=q+1}{\overset{n}{\sum}}\overset{q}{\underset{\text{b,c=1}}{\sum}%
}R_{ipip}(R_{\text{bcbc}}^{P}-R_{\text{bcbc}}^{M}%
)+\underset{j=q+1}{\overset{n}{\sum}}$ $\overset{q}{\underset{\text{b,c=1}%
}{\sum}}R_{jpjp}(R_{\text{bcbc}}^{P}-R_{\text{bcbc}}^{M})](y_{0})$

$+\frac{1}{72}\underset{i,j,p=q+1}{\overset{n}{\sum}}%
\overset{q}{\underset{\text{b,c=1}}{\sum}}[R_{ipjp}(T_{\text{bb}i}%
T_{\text{cc}j}-T_{\text{bc}i}T_{\text{bc}j})+R_{ipjp}(T_{\text{bb}%
j}T_{\text{cc}i}-T_{\text{bc}j}T_{\text{bc}i})](y_{0})\qquad$

$-\frac{1}{288}\underset{i,j=q+1}{\overset{n}{\sum}}[T_{\text{aa}%
i}T_{\text{bb}j}(T_{\text{cc}i}T_{\text{dd}j}-T_{\text{cd}i}T_{\text{dc}%
j})+T_{\text{aa}i}T_{\text{bb}j}(T_{\text{cc}j}T_{\text{dd}i}-T_{\text{cd}%
j}T_{\text{dc}i})\qquad E$

$+T_{\text{aa}j}T_{\text{bb}i}(T_{\text{cc}i}T_{\text{dd}j}-T_{\text{cd}%
i}T_{\text{dc}j})+T_{\text{aa}j}T_{\text{bb}i}(T_{\text{cc}j}T_{\text{dd}%
i}-T_{\text{cd}j}T_{\text{dc}i})](y_{0})$

$+\frac{1}{288}\underset{i,j=q+1}{\overset{n}{\sum}}[T_{\text{aa}%
i}T_{\text{bc}j}(T_{\text{bc}i}T_{\text{dd}j}-T_{\text{bd}i}T_{\text{cd}%
j})+T_{\text{aa}i}T_{\text{bc}j}(T_{\text{bc}j}T_{\text{dd}i}-T_{\text{bd}%
j}T_{\text{cd}i})$

$+T_{\text{aa}j}T_{\text{bc}i}(T_{\text{bc}i}T_{\text{dd}j}-T_{\text{bd}%
i}T_{\text{cd}j})+T_{\text{aa}j}T_{\text{bc}i}(T_{\text{bc}j}T_{\text{dd}%
i}-T_{\text{bd}j}T_{\text{cd}i})](y_{0})$

$-\frac{1}{288}\underset{i,j=q+1}{\overset{n}{\sum}}[T_{\text{aa}%
i}T_{\text{bd}j}(T_{\text{bc}i}T_{\text{cd}j}-T_{\text{bd}i}T_{\text{cc}%
j})+T_{\text{aa}i}T_{\text{bd}j}(T_{\text{bc}j}T_{\text{cd}i}-T_{\text{bd}%
j}T_{\text{cc}i})$

$+T_{\text{aa}j}T_{\text{bd}i}(T_{\text{bc}i}T_{\text{cd}j}-T_{\text{bd}%
i}T_{\text{cc}j})+T_{\text{aa}j}T_{\text{bd}i}(T_{\text{bc}j}T_{\text{cd}%
i}-T_{\text{bd}j}T_{\text{cc}i})](y_{0})\qquad$

$+\frac{1}{288}\underset{i,j=q+1}{\overset{n}{\sum}}[T_{\text{ab}%
i}T_{\text{ab}j}(T_{\text{cc}i}T_{\text{dd}j}-T_{\text{cd}i}T_{\text{dc}%
j})+T_{\text{ab}i}T_{\text{ab}j}(T_{\text{cc}j}T_{\text{dd}i}-T_{\text{cd}%
j}T_{\text{dc}i})$

$+T_{\text{ab}j}T_{\text{ab}i}(T_{\text{cc}i}T_{\text{dd}j}-T_{\text{cd}%
i}T_{\text{dc}j})+T_{\text{ab}j}T_{\text{ab}i}(T_{\text{cc}j}T_{\text{dd}%
i}-T_{\text{cd}j}T_{\text{dc}i})](y_{0})$

$-\frac{1}{288}\underset{i,j=q+1}{\overset{n}{\sum}}[T_{\text{ab}%
i}T_{\text{bc}j}(T_{\text{ac}i}T_{\text{dd}j}-T_{\text{ad}i}T_{\text{cd}%
j})+T_{\text{ab}i}T_{\text{bc}j}(T_{\text{ac}j}T_{\text{dd}i}-T_{\text{ad}%
j}T_{\text{cd}i})$

$+T_{\text{ab}j}T_{\text{bc}i}(T_{\text{ac}i}T_{\text{dd}j}-T_{\text{ad}%
i}T_{\text{cd}j})+T_{\text{ab}j}T_{\text{bc}i}(T_{\text{ac}j}T_{\text{dd}%
i}-T_{\text{ad}j}T_{\text{cd}i})](y_{0})$

$+\frac{1}{288}\underset{i,j=q+1}{\overset{n}{\sum}}[T_{\text{ab}%
i}T_{\text{bd}j}(T_{\text{ac}i}T_{\text{cd}j}-T_{\text{ad}i}T_{\text{cc}%
j})+T_{\text{ab}i}T_{\text{bd}j}(T_{\text{ac}j}T_{\text{cd}i}-T_{\text{ad}%
j}T_{\text{cc}i})$

$+T_{\text{ab}i}T_{\text{bd}j}(T_{\text{ac}j}T_{\text{cd}i}-T_{\text{ad}%
j}T_{\text{cc}i})+T_{\text{ab}j}T_{\text{bd}i}(T_{\text{ac}j}T_{\text{cd}%
i}-T_{\text{ad}j}T_{\text{cc}i})](y_{0})$

$-\ \frac{1}{288}\underset{i,j=q+1}{\overset{n}{\sum}}[T_{\text{ac}%
i}T_{\text{ab}j}(T_{\text{bc}i}T_{\text{dd}j}-T_{\text{bd}i}T_{\text{dc}%
j})+T_{\text{ac}i}T_{\text{ab}j}(T_{\text{bc}j}T_{\text{dd}i}-T_{\text{bd}%
j}T_{\text{dc}i})$

$+T_{\text{ac}j}T_{\text{ab}i}(T_{\text{bc}i}T_{\text{dd}j}-T_{\text{bd}%
i}T_{\text{dc}j})+T_{\text{ac}j}T_{\text{ab}i}(T_{\text{bc}j}T_{\text{dd}%
i}-T_{\text{bd}j}T_{\text{dc}i})](y_{0})$

$+\ \frac{1}{288}\underset{i,j=q+1}{\overset{n}{\sum}}[T_{\text{ac}%
i}T_{\text{bb}j}(T_{\text{ac}i}T_{\text{dd}j}-T_{\text{ad}i}T_{\text{cd}%
j})+T_{\text{ac}i}T_{\text{bb}j}(T_{\text{ac}j}T_{\text{dd}i}-T_{\text{ad}%
j}T_{\text{cd}i})$

$+T_{\text{ac}j}T_{\text{bb}i}(T_{\text{ac}i}T_{\text{dd}j}-T_{\text{ad}%
i}T_{\text{cd}i})+T_{\text{ac}j}T_{\text{bb}i}(T_{\text{ac}j}T_{\text{dd}%
i}-T_{\text{ad}j}T_{\text{cd}i})](y_{0})$

$-\ \frac{1}{288}\underset{i,j=q+1}{\overset{n}{\sum}}[T_{\text{ac}%
i}T_{\text{bd}j}(T_{\text{ac}i}T_{\text{bd}j}-T_{\text{ad}i}T_{\text{bc}%
j})+T_{\text{ac}i}T_{\text{bd}j}(T_{\text{ac}j}T_{\text{bd}i}-T_{\text{ad}%
j}T_{\text{bc}i})$

$+T_{\text{ac}j}T_{\text{bd}i}(T_{\text{ac}i}T_{\text{bd}j}-T_{\text{ad}%
i}T_{\text{bc}j})+T_{\text{ac}j}T_{\text{bd}i}(T_{\text{ac}j}T_{\text{bd}%
i}-T_{\text{ad}j}T_{\text{bc}i})](y_{0})$

$+\frac{1}{288}\underset{i,j=q+1}{\overset{n}{\sum}}[T_{\text{ad}%
i}T_{\text{ab}j}(T_{\text{bc}i}T_{\text{cd}j}-T_{\text{bd}i}T_{\text{cc}%
j})+T_{\text{ad}i}T_{\text{ab}j}(T_{\text{bc}j}T_{\text{cd}i}-T_{\text{bd}%
j}T_{\text{cc}i})$

$+T_{\text{ad}j}T_{\text{ab}i}(T_{\text{bc}i}T_{\text{cd}j}-T_{\text{bd}%
i}T_{\text{cc}j})+T_{\text{ad}j}T_{\text{ab}i}(T_{\text{bc}j}T_{\text{cd}%
i}-T_{\text{bd}j}T_{\text{cc}i})](y_{0})$

$-\ \frac{1}{288}\underset{i,j=q+1}{\overset{n}{\sum}}[T_{\text{ad}%
i}T_{\text{bb}j}(T_{\text{ac}i}T_{\text{cd}j}-T_{\text{ad}i}T_{\text{cc}%
j})+T_{\text{ad}i}T_{\text{bb}j}(T_{\text{ac}j}T_{\text{cd}i}-T_{\text{ad}%
j}T_{\text{cc}i})$

$+T_{\text{ad}j}T_{\text{bb}i}(T_{\text{ac}i}T_{\text{cd}j}-T_{\text{ad}%
i}T_{\text{cc}j})+T_{\text{ad}j}T_{\text{bb}i}(T_{\text{ac}j}T_{\text{cd}%
i}-T_{\text{ad}j}T_{\text{cc}i})](y_{0})$

$+\ \frac{1}{288}\underset{i,j=q+1}{\overset{n}{\sum}}[T_{\text{ad}%
i}T_{\text{bc}j}(T_{\text{ac}i}T_{\text{bd}j}-T_{\text{ad}i}T_{\text{bc}%
j})+T_{\text{ad}i}T_{\text{bc}j}(T_{\text{ac}j}T_{\text{bd}i}-T_{\text{ad}%
j}T_{\text{bc}i})$

$+T_{\text{ad}j}T_{\text{bc}i}(T_{\text{ac}i}T_{\text{bd}j}-T_{\text{ad}%
i}T_{\text{bc}j})+T_{\text{ad}j}T_{\text{bc}i}(T_{\text{ac}j}T_{\text{bd}%
i}-T_{\text{ad}j}T_{\text{bc}i})](y_{0})$

$-\ \frac{1}{144}[(R_{\text{cdcd}}^{P}-R_{\text{cdcd}}^{M})(R_{\text{abab}%
}^{P}-R_{\text{abab}}^{M})](y_{0})$

$+\frac{1}{144}[(R_{\text{bdcd}}^{P}-R_{\text{bdcd}}^{M})(R_{\text{abac}}%
^{P}-R_{\text{abac}}^{M})](y_{0})$

$+\ \frac{1}{144}[(R_{\text{bcdc}}^{P}-R_{\text{bcdc}}^{M})(R_{\text{abad}%
}^{P}-R_{\text{abad}}^{M})](y_{0})$

$-\ \frac{1}{144}[(R_{\text{adcd}}^{P}-R_{\text{adcd}}^{M})(R_{\text{abbc}%
}^{P}-R_{\text{abbc}}^{M})](y_{0})\qquad$

$+\ \frac{1}{144}[(R_{\text{acdc}}^{P}-R_{\text{acdc}}^{M})(R_{\text{abdb}%
}^{P}-R_{\text{abdb}}^{M})](y_{0})$

$-\ \frac{1}{576}[(R_{\text{abcd}}^{P}-R_{\text{abcd}}^{M})]^{2}(y_{0}%
)\qquad(6)$

$+\frac{1}{24}\times\frac{21}{2}<H,i>^{2}(y_{0})<H,j>^{2}(y_{0})\qquad
\qquad\qquad L_{3}$

$-\frac{1}{24}\times\frac{1}{3}[<H,i><H,j>](y_{0})$

$\times\lbrack2\varrho_{ij}+$ $\overset{q}{\underset{\text{a}=1}{4\sum}%
}R_{i\text{a}j\text{a}}-3\overset{q}{\underset{\text{a,b=1}}{\sum}%
}(T_{\text{aa}i}T_{\text{bb}j}-T_{\text{ab}i}T_{\text{ab}j}%
)-3\overset{q}{\underset{\text{a,b=1}}{\sum}}(T_{\text{aa}j}T_{\text{bb}%
i}-T_{\text{ab}j}T_{\text{ab}i}](y_{0})$

$-\frac{1}{24}\times\frac{1}{3}<H,i>(y_{0})<H,k>(y_{0})R_{ijjk}(y_{0})$

\ $-\frac{1}{24}\times<H,i>^{2}(y_{0})[\tau^{M}\ -3\tau^{P}%
+\ \underset{\text{a}=1}{\overset{\text{q}}{\sum}}\varrho_{\text{aa}}%
^{M}+\overset{q}{\underset{\text{a},\text{b}=1}{\sum}}R_{\text{abab}}%
^{M}](y_{0})$

$-\frac{1}{24}\times\frac{1}{6}<H,i>(y_{0})[\nabla_{i}\varrho_{jj}%
-2\varrho_{ij}<H,j>+\overset{q}{\underset{\text{a}=1}{\sum}}(\nabla
_{i}R_{\text{a}j\text{a}j}-4R_{i\text{a}j\text{a}}<H,j>)$

$+4\overset{q}{\underset{\text{a,b=1}}{\sum}}R_{i\text{a}j\text{b}%
}T_{\text{ab}j}+2\overset{q}{\underset{\text{a,b,c=1}}{\sum}}(T_{\text{aa}%
i}T_{\text{bb}j}T_{\text{cc}j}-3T_{\text{aa}i}T_{\text{bc}j}T_{\text{bc}%
j}+2T_{\text{ab}i}T_{\text{bc}j}T_{\text{ca}j})](y_{0})$\qquad\qquad
\qquad\qquad\qquad\ \ 

$-\frac{1}{24}\times\frac{1}{6}<H,i>(y_{0})[\nabla_{j}\varrho_{ij}%
-2\varrho_{ij}<H,j>+\overset{q}{\underset{\text{a}=1}{\sum}}(\nabla
_{j}R_{\text{a}i\text{a}j}-4R_{j\text{a}i\text{a}}<H,j>)$

$+4\overset{q}{\underset{\text{a,b=1}}{\sum}}R_{j\text{a}i\text{b}%
}T_{\text{ab}j}+2\overset{q}{\underset{\text{a,b,c=1}}{\sum}}(T_{\text{aa}%
j}T_{\text{bb}i}T_{\text{cc}j}-3T_{\text{aa}j}T_{\text{bc}i}T_{\text{bc}%
j}+2T_{\text{ab}j}T_{\text{bc}i}T_{\text{ac}j})](y_{0})$

$-\frac{1}{24}\times\frac{1}{6}<H,i>(y_{0})[\nabla_{j}\varrho_{ij}%
-2\varrho_{jj}<H,i>+\overset{q}{\underset{\text{a}=1}{\sum}}(\nabla
_{j}R_{\text{a}i\text{a}j}-4R_{j\text{a}j\text{a}}<H,i>)$

$+4\overset{q}{\underset{\text{a,b=1}}{\sum}}R_{j\text{a}j\text{b}%
}T_{\text{ab}i}+2\overset{q}{\underset{\text{a,b,c=1}}{\sum}}(T_{\text{aa}%
j}T_{\text{bb}j}T_{\text{cc}i}-3T_{\text{aa}j}T_{\text{bc}j}T_{\text{bc}%
i}+2T_{\text{ab}j}T_{\text{bc}j}T_{\text{ac}i})](y_{0})$

$\qquad-\frac{1}{24}\times\frac{1}{3}[<H,i><H,k>](y_{0})R_{ijjk}(y_{0})$

$-\frac{1}{24}\times\frac{1}{6}<H,i>(y_{0})[\nabla_{i}\varrho_{jj}%
-2\varrho_{ij}<H,j>+\overset{q}{\underset{\text{a}=1}{\sum}}(\nabla
_{i}R_{\text{a}j\text{a}j}-4R_{i\text{a}j\text{a}}<H,j>)$

$+4\overset{q}{\underset{\text{a,b=1}}{\sum}}R_{i\text{a}j\text{b}%
}T_{\text{ab}j}+2\overset{q}{\underset{\text{a,b,c=1}}{\sum}}(T_{\text{aa}%
i}T_{\text{bb}j}T_{\text{cc}j}-3T_{\text{aa}i}T_{\text{bc}j}T_{\text{bc}%
j}+2T_{\text{ab}i}T_{\text{bc}j}T_{\text{ca}j})](y_{0})$\qquad\qquad
\qquad\qquad\qquad\ \ 

$-\frac{1}{24}\times\frac{1}{6}<H,i>(y_{0})[\nabla_{j}\varrho_{ij}%
-2\varrho_{ij}<H,j>+\overset{q}{\underset{\text{a}=1}{\sum}}(\nabla
_{j}R_{\text{a}i\text{a}j}-4R_{j\text{a}i\text{a}}<H,j>)$

$+4\overset{q}{\underset{\text{a,b=1}}{\sum}}R_{j\text{a}i\text{b}%
}T_{\text{ab}j}+2\overset{q}{\underset{\text{a,b,c=1}}{\sum}}(T_{\text{aa}%
j}T_{\text{bb}i}T_{\text{cc}j}-3T_{\text{aa}j}T_{\text{bc}i}T_{\text{bc}%
j}+2T_{\text{ab}j}T_{\text{bc}i}T_{\text{ac}j})](y_{0})$

$-\frac{1}{24}\times\frac{1}{6}<H,i>(y_{0})[\nabla_{j}\varrho_{ij}%
-2\varrho_{jj}<H,i>+\overset{q}{\underset{\text{a}=1}{\sum}}(\nabla
_{j}R_{\text{a}i\text{a}j}-4R_{j\text{a}j\text{a}}<H,i>)$

$+4\overset{q}{\underset{\text{a,b=1}}{\sum}}R_{j\text{a}j\text{b}%
}T_{\text{ab}i}+2\overset{q}{\underset{\text{a,b,c=1}}{\sum}}(T_{\text{aa}%
j}T_{\text{bb}j}T_{\text{cc}i}-3T_{\text{aa}j}T_{\text{bc}j}T_{\text{bc}%
i}+2T_{\text{ab}j}T_{\text{bc}j}T_{\text{ac}i})](y_{0})$

$\qquad\qquad\qquad\qquad\qquad\qquad\qquad\qquad\qquad\qquad\qquad
\qquad\qquad\qquad\qquad\qquad\blacksquare\qquad\qquad\qquad\qquad\qquad
\qquad$

(viii)$\qquad$A$_{3213}$\ $=\frac{1}{12}(\frac{\partial\theta^{\frac{1}{2}}%
}{\partial\text{x}_{i}})\frac{\partial}{\partial\text{x}_{i}}(\Delta
\theta^{-\frac{1}{2}})(y_{0})$

$=\frac{1}{288}<H,i>^{2}(y_{0})[3<H,i>^{2}(y_{0})+2(\tau^{M}-3\tau
^{P}+\overset{q}{\underset{\text{a}=1}{\sum}}\varrho_{\text{aa}}%
^{M}+\overset{q}{\underset{\text{a,b}=1}{\sum}}R_{\text{abab}}^{M})](y_{0})$

$+\frac{1}{32}[<H,i>^{2}<H,j>^{2}](y_{0})$

$+\frac{1}{288}[<H,i><H,j>](y_{0})$

$\times\lbrack2\varrho_{ij}+$ $\overset{q}{\underset{\text{a}=1}{4\sum}%
}R_{i\text{a}j\text{a}}-3\overset{q}{\underset{\text{a,b=1}}{\sum}%
}(T_{\text{aa}i}T_{\text{bb}j}-T_{\text{ab}i}T_{\text{ab}j}%
)-3\overset{q}{\underset{\text{a,b=1}}{\sum}}(T_{\text{aa}j}T_{\text{bb}%
i}-T_{\text{ab}j}T_{\text{ab}i}](y_{0})\qquad$

$+\frac{1}{288}[<H,i><H,j>](y_{0})$

$\times\lbrack2\varrho_{ij}+$ $\overset{q}{\underset{\text{a}=1}{4\sum}%
}R_{i\text{a}j\text{a}}-3\overset{q}{\underset{\text{a,b=1}}{\sum}%
}(T_{\text{aa}i}T_{\text{bb}j}-T_{\text{ab}i}T_{\text{ab}j}%
)-3\overset{q}{\underset{\text{a,b=1}}{\sum}}(T_{\text{aa}j}T_{\text{bb}%
i}-T_{\text{ab}j}T_{\text{ab}i}](y_{0})$

$+\frac{1}{144}[<H,i><H,k>](y_{0})R_{jijk}(y_{0})$

$+\frac{5}{64}[<H,i>^{2}<H,j>^{2}](y_{0})$

$-\frac{1}{96}[<H,i><H,j>](y_{0})$

$\times\lbrack2\varrho_{ij}+\overset{q}{\underset{\text{a}=1}{4\sum}%
}R_{i\text{a}j\text{a}}-3\overset{q}{\underset{\text{a,b=1}}{\sum}%
}(T_{\text{aa}i}T_{\text{bb}j}-T_{\text{ab}i}T_{\text{ab}j}%
)-3\overset{q}{\underset{\text{a,b=1}}{\sum}}(T_{\text{aa}j}T_{\text{bb}%
i}-T_{\text{ab}j}T_{\text{ab}i}](y_{0})$

\ $-\frac{1}{96}<H,i>^{2}(y_{0})[\tau^{M}\ -3\tau^{P}+\ \underset{\text{a}%
=1}{\overset{\text{q}}{\sum}}\varrho_{\text{aa}}^{M}%
+\overset{q}{\underset{\text{a},\text{b}=1}{\sum}}R_{\text{abab}}^{M}](y_{0})$

$-\frac{1}{288}<H,i>(y_{0})[\nabla_{i}\varrho_{jj}-2\varrho_{ij}%
<H,j>+\overset{q}{\underset{\text{a}=1}{\sum}}(\nabla_{i}R_{\text{a}%
j\text{a}j}-4R_{i\text{a}j\text{a}}<H,j>)$

$+4\overset{q}{\underset{\text{a,b=1}}{\sum}}R_{i\text{a}j\text{b}%
}T_{\text{ab}j}+2\overset{q}{\underset{\text{a,b,c=1}}{\sum}}(T_{\text{aa}%
i}T_{\text{bb}j}T_{\text{cc}j}-3T_{\text{aa}i}T_{\text{bc}j}T_{\text{bc}%
j}+2T_{\text{ab}i}T_{\text{bc}j}T_{\text{ca}j})](y_{0})$\qquad\qquad
\qquad\qquad\qquad\ \ 

$-\frac{1}{288}<H,i>(y_{0})[\nabla_{j}\varrho_{ij}-2\varrho_{ij}%
<H,j>+\overset{q}{\underset{\text{a}=1}{\sum}}(\nabla_{j}R_{\text{a}%
i\text{a}j}-4R_{j\text{a}i\text{a}}<H,j>)$

$\qquad+4\overset{q}{\underset{\text{a,b=1}}{\sum}}R_{j\text{a}i\text{b}%
}T_{\text{ab}j}+2\overset{q}{\underset{\text{a,b,c=1}}{\sum}}(T_{\text{aa}%
j}T_{\text{bb}i}T_{\text{cc}j}-3T_{\text{aa}j}T_{\text{bc}i}T_{\text{bc}%
j}+2T_{\text{ab}j}T_{\text{bc}i}T_{\text{ac}j})](y_{0})$

$-\frac{1}{288}<H,i>(y_{0})[\nabla_{j}\varrho_{ij}-2\varrho_{jj}%
<H,i>+\overset{q}{\underset{\text{a}=1}{\sum}}(\nabla_{j}R_{\text{a}%
i\text{a}j}-4R_{j\text{a}j\text{a}}<H,i>)$

$+4\overset{q}{\underset{\text{a,b=1}}{\sum}}R_{j\text{a}j\text{b}%
}T_{\text{ab}i}+2\overset{q}{\underset{\text{a,b,c=1}}{\sum}}(T_{\text{aa}%
j}T_{\text{bb}j}T_{\text{cc}i}-3T_{\text{aa}j}T_{\text{bc}j}T_{\text{bc}%
i}+2T_{\text{ab}j}T_{\text{bc}j}T_{\text{ac}i})](y_{0})$

$+\frac{1}{288}[<H,i><H,j>](y_{0})$

$\times\lbrack2\varrho_{ij}+$ $\overset{q}{\underset{\text{a}=1}{4\sum}%
}R_{i\text{a}j\text{a}}-3\overset{q}{\underset{\text{a,b=1}}{\sum}%
}(T_{\text{aa}i}T_{\text{bb}j}-T_{\text{ab}i}T_{\text{ab}j}%
)-3\overset{q}{\underset{\text{a,b=1}}{\sum}}(T_{\text{aa}j}T_{\text{bb}%
i}-T_{\text{ab}j}T_{\text{ab}i}](y_{0})$

$+\frac{1}{32}[<H,i>^{2}<H,j>^{2}](y_{0})$

$+\frac{1}{288}[<H,i><H,j>](y_{0})$

$\times\lbrack2\varrho_{ij}+$ $\overset{q}{\underset{\text{a}=1}{4\sum}%
}R_{i\text{a}j\text{a}}-3\overset{q}{\underset{\text{a,b=1}}{\sum}%
}(T_{\text{aa}i}T_{\text{bb}j}-T_{\text{ab}i}T_{\text{ab}j}%
)-3\overset{q}{\underset{\text{a,b=1}}{\sum}}(T_{\text{aa}j}T_{\text{bb}%
i}-T_{\text{ab}j}T_{\text{ab}i}](y_{0})$

$-\frac{1}{144}<H,i>(y_{0})<H,k>(y_{0})R_{ijjk}(y_{0})$

$+\frac{5}{64}[<H,i>^{2}<H,j>^{2}](y_{0})$

$-\frac{1}{96}[<H,i><H,j>](y_{0})$

$\times\lbrack2\varrho_{ij}+\overset{q}{\underset{\text{a}=1}{4\sum}%
}R_{i\text{a}j\text{a}}-3\overset{q}{\underset{\text{a,b=1}}{\sum}%
}(T_{\text{aa}i}T_{\text{bb}j}-T_{\text{ab}i}T_{\text{ab}j}%
)-3\overset{q}{\underset{\text{a,b=1}}{\sum}}(T_{\text{aa}j}T_{\text{bb}%
i}-T_{\text{ab}j}T_{\text{ab}i}](y_{0})$

$-\frac{1}{96}<H,i>^{2}(y_{0})[\tau^{M}\ -3\tau^{P}+\ \underset{\text{a}%
=1}{\overset{\text{q}}{\sum}}\varrho_{\text{aa}}^{M}%
+\overset{q}{\underset{\text{a},\text{b}=1}{\sum}}R_{\text{abab}}^{M}](y_{0})$

$-\frac{1}{288}<H,i>(y_{0})[\nabla_{i}\varrho_{jj}-2\varrho_{ij}%
<H,j>+\overset{q}{\underset{\text{a}=1}{\sum}}(\nabla_{i}R_{\text{a}%
j\text{a}j}-4R_{i\text{a}j\text{a}}<H,j>)$

$+4\overset{q}{\underset{\text{a,b=1}}{\sum}}R_{i\text{a}j\text{b}%
}T_{\text{ab}j}+2\overset{q}{\underset{\text{a,b,c=1}}{\sum}}(T_{\text{aa}%
i}T_{\text{bb}j}T_{\text{cc}j}-3T_{\text{aa}i}T_{\text{bc}j}T_{\text{bc}%
j}+2T_{\text{ab}i}T_{\text{bc}j}T_{\text{ca}j})](y_{0})$\qquad\qquad
\qquad\qquad\qquad\ \ 

$-\frac{1}{288}<H,i>(y_{0})[\nabla_{j}\varrho_{ij}-2\varrho_{ij}%
<H,j>+\overset{q}{\underset{\text{a}=1}{\sum}}(\nabla_{j}R_{\text{a}%
i\text{a}j}-4R_{j\text{a}i\text{a}}<H,j>)$

$+4\overset{q}{\underset{\text{a,b=1}}{\sum}}R_{j\text{a}i\text{b}%
}T_{\text{ab}j}+2\overset{q}{\underset{\text{a,b,c=1}}{\sum}}(T_{\text{aa}%
j}T_{\text{bb}i}T_{\text{cc}j}-3T_{\text{aa}j}T_{\text{bc}i}T_{\text{bc}%
j}+2T_{\text{ab}j}T_{\text{bc}i}T_{\text{ac}j})](y_{0})$

$-\frac{1}{288}[\nabla_{j}\varrho_{ij}-2\varrho_{jj}%
<H,i>+\overset{q}{\underset{\text{a}=1}{\sum}}(\nabla_{j}R_{\text{a}%
i\text{a}j}-4R_{j\text{a}j\text{a}}<H,i>)+4\overset{q}{\underset{\text{a,b=1}%
}{\sum}}R_{j\text{a}j\text{b}}T_{\text{ab}i}$

$+2\overset{q}{\underset{\text{a,b,c=1}}{\sum}}(T_{\text{aa}j}T_{\text{bb}%
j}T_{\text{cc}i}-3T_{\text{aa}j}T_{\text{bc}j}T_{\text{bc}i}+2T_{\text{ab}%
j}T_{\text{bc}j}T_{\text{ac}i})](y_{0})$

(ix)$\qquad$A$_{321}=\frac{1}{24}\frac{\partial^{2}}{\partial\text{x}_{i}^{2}%
}(\theta^{\frac{1}{2}}\Delta\theta^{-\frac{1}{2}})(y_{0})\phi(y_{0})$

$\qquad\qquad=\frac{1}{24}\left[  (\frac{\partial^{2}\theta^{\frac{1}{2}}%
}{\partial\text{x}_{i}^{2}})(\Delta\theta^{-\frac{1}{2}})+\frac{\partial^{2}%
}{\partial\text{x}_{i}^{2}}(\Delta\theta^{-\frac{1}{2}})+2(\frac
{\partial\theta^{\frac{1}{2}}}{\partial\text{x}_{i}})\frac{\partial}%
{\partial\text{x}_{i}}(\Delta\theta^{-\frac{1}{2}})\right]  (y_{0})\phi
(y_{0})$

$=\ $A$_{3211}+$ A$_{3212}+$ A$_{3213}$ is given below at the end of
\textbf{Table A}$_{10}$.

\qquad\qquad\qquad\qquad\qquad\qquad\qquad\qquad\qquad\qquad\qquad\qquad
\qquad\qquad\qquad\qquad\qquad$\blacksquare$

\subsection{\protect\underline{\textbf{Computations}}}

(i)$\qquad\ \nabla$ is the gradient operator here:

$\qquad<\nabla\theta,\nabla f>(y_{0})$ $=$ $\theta(y_{0})<\nabla\log
\theta,\nabla f>(y_{0})$

$\qquad=$ $<\nabla\log\theta,\nabla f>(y_{0})=(\nabla\log\theta)_{i}%
(y_{0})\frac{\partial f}{\partial\text{x}_{i}}(y_{0})$

By (i) of Table 10, we have:

$\left(  A_{25}\right)  \qquad<\nabla\theta,\nabla f>(y_{0})=-$
$\underset{i=q+1}{\overset{n}{\sum}}<H,i>(y_{0})\frac{\partial f}%
{\partial\text{x}_{i}}(y_{0})$

(ii$\;\not )  \;$ We can also use the usual definition of the \textbf{scalar}
Laplacian given by:

$\Delta f$ \ $=g^{ij}[\frac{\partial^{2}f}{\partial x_{i}\partial_{j}}%
-\Gamma_{ij}^{k}\frac{\partial f}{\partial x_{k}}]$

$\Delta\theta^{-\frac{1}{2}}(y_{0})$\ $=g^{ij}(y_{0})[\frac{\partial^{2}%
\theta^{-\frac{1}{2}}}{\partial x_{i}\partial_{j}}-\Gamma_{ij}^{k}%
\frac{\partial\theta^{-\frac{1}{2}}}{\partial x_{k}}](y_{0})=\delta^{ij}%
(y_{0})[\frac{\partial^{2}\theta^{-\frac{1}{2}}}{\partial x_{i}\partial_{j}%
}-\Gamma_{ij}^{k}\frac{\partial\theta^{-\frac{1}{2}}}{\partial x_{k}}](y_{0})$

$\qquad\qquad=[\frac{\partial^{2}\theta^{-\frac{1}{2}}}{\partial x_{i}^{2}%
}-\Gamma_{ii}^{k}\frac{\partial\theta^{-\frac{1}{2}}}{\partial x_{k}}](y_{0})$

Since the expansion of $\theta=\theta_{P}$ in \textbf{Proposition 11} is given
\textbf{normal} Fermi coordinates, we have:

$\frac{\partial^{2}\theta^{-\frac{1}{2}}}{\partial x_{\text{a}}^{2}}%
(y_{0})=0=\frac{\partial\theta^{-\frac{1}{2}}}{\partial x_{\text{a}}}(y_{0})$
for a = 1,...,q and since $\Gamma_{\text{ab}}^{\text{c}}(y_{0})=0=\Gamma
_{ii}^{\text{c}}(y_{0})$ for a, b, c = 1,...,q and $i=q+1,...,n.$

$\Delta\theta^{-\frac{1}{2}}(y_{0})=[-\Gamma_{\text{aa}}^{k}\frac
{\partial\theta^{-\frac{1}{2}}}{\partial x_{k}}](y_{0})+[\frac{\partial
^{2}\theta^{-\frac{1}{2}}}{\partial x_{i}^{2}}-\Gamma_{ii}^{k}\frac
{\partial\theta^{-\frac{1}{2}}}{\partial x_{k}}](y_{0})$ for a = 1,...,q and
$i,j,k=q+1,...,n.$

We further have $\Gamma_{ii}^{k}(y_{0})=0$ for $i,j,k=q+1,...,n$ and so the
final expression for the Laplacian is given by:

\qquad$\Delta\theta^{-\frac{1}{2}}(y_{0})=[-\Gamma_{\text{aa}}^{k}%
\frac{\partial\theta^{-\frac{1}{2}}}{\partial x_{k}}](y_{0})+\frac
{\partial^{2}\theta^{-\frac{1}{2}}}{\partial x_{i}^{2}}(y_{0})$ for a =
1,...,q and $i,j,k=q+1,...,n.$

By (x) of \textbf{Table A}$_{9},$ we have:

$\frac{\partial^{2}\theta^{-\frac{1}{2}}}{\partial\text{x}_{i}^{2}}%
(y_{0})=\frac{3}{4}<H,i>^{2}(y_{0})+\frac{1}{6}(\tau^{M}-3\tau^{P}%
+\overset{q}{\underset{\text{a}=1}{\sum}}\varrho_{\text{aa}}^{M}%
+\overset{q}{\underset{\text{a,b}=1}{\sum}}R_{\text{abab}}^{M})(y_{0})$

$\Gamma_{\text{aa}}^{i}(y_{0})=T_{\text{aa}i}(y_{0})=$ $<H,i>(y_{0})$ by (i)
of \textbf{Table A}$_{7}$

$\frac{\partial\theta^{-\frac{1}{2}}}{\partial\text{x}_{i}}(y_{0})=\frac{1}%
{2}<H,i>(y_{0})$ by (iv) of \textbf{Table A}$_{9}.$

Therefore,

$\Delta\theta^{-\frac{1}{2}}(y_{0})=[-\Gamma_{\text{aa}}^{i}\frac
{\partial\theta^{-\frac{1}{2}}}{\partial x_{i}}](y_{0})+\frac{\partial
^{2}\theta^{-\frac{1}{2}}}{\partial x_{i}^{2}}(y_{0})$

$=-\frac{1}{2}<H,i>^{2}(y_{0})+\frac{3}{4}<H,i>^{2}(y_{0})+\frac{1}{6}%
(\tau^{M}-3\tau^{P}+\overset{q}{\underset{\text{a}=1}{\sum}}\varrho
_{\text{aa}}^{M}+\overset{q}{\underset{\text{a,b}=1}{\sum}}R_{\text{abab}}%
^{M})(y_{0})$

$=\frac{1}{4}<H,i>^{2}(y_{0})+\frac{1}{6}(\tau^{M}-3\tau^{P}%
+\overset{q}{\underset{\text{a}=1}{\sum}}\varrho_{\text{aa}}^{M}%
+\overset{q}{\underset{\text{a,b}=1}{\sum}}R_{\text{abab}}^{M})(y_{0})$

For the Laplace-Beltrami operator $\Delta$ we can avoid the use of Christoffel
symbols $\Gamma_{ij}^{k}$ and use the version of the formula given by:

$\Delta f=\theta^{-1}\frac{\partial}{\partial\text{x}_{i}}[\theta g^{ij}%
\frac{\partial f}{\partial\text{x}_{j}}]$ for a smooth function
$f:M\longrightarrow R$

where we assume the \textbf{Einstein convention} of summation over repeated
indices, as usual.

$\Delta\theta^{-\frac{1}{2}}(y_{0})=\theta^{-1}(y_{0})\frac{\partial}%
{\partial\text{x}_{i}}[\theta g^{ij}\frac{\partial\theta^{-\frac{1}{2}}%
}{\partial\text{x}_{j}}](y_{0})$

Since $\theta(y_{0})=1,$ we have:

$\Delta\theta^{-\frac{1}{2}}(y_{0})=\frac{\partial}{\partial\text{x}_{i}%
}[\theta g^{ij}\frac{\partial\theta^{-\frac{1}{2}}}{\partial\text{x}_{j}%
}](y_{0})=\frac{\partial\theta}{\partial\text{x}_{i}}(y_{0})[g^{ij}%
\frac{\partial\theta^{-\frac{1}{2}}}{\partial\text{x}_{j}}](y_{0}%
)+\frac{\partial}{\partial\text{x}_{i}}[g^{ij}\frac{\partial\theta^{-\frac
{1}{2}}}{\partial\text{x}_{j}}](y_{0})$

$=\frac{\partial\theta}{\partial\text{x}_{i}}(y_{0})[g^{ij}\frac
{\partial\theta^{-\frac{1}{2}}}{\partial\text{x}_{j}}](y_{0})+\frac{\partial
g^{ij}}{\partial\text{x}_{i}}(y_{0})\frac{\partial\theta^{-\frac{1}{2}}%
}{\partial\text{x}_{j}}(y_{0})+g^{ij}(y_{0})\frac{\partial^{2}\theta
^{-\frac{1}{2}}}{\partial\text{x}_{i}\partial\text{x}_{j}}](y_{0})$

Since $g^{ij}(y_{0})=\delta^{ij},$ we have:

$\Delta\theta^{-\frac{1}{2}}(y_{0})=\frac{\partial\theta}{\partial\text{x}%
_{i}}(y_{0})\frac{\partial\theta^{-\frac{1}{2}}}{\partial\text{x}_{i}}%
(y_{0})+\frac{\partial g^{ij}}{\partial\text{x}_{i}}(y_{0})\frac
{\partial\theta^{-\frac{1}{2}}}{\partial\text{x}_{j}}(y_{0})+\frac
{\partial^{2}\theta^{-\frac{1}{2}}}{\partial\text{x}_{i}^{2}}(y_{0})$

$g^{ij}$ and $\theta$ are both expanded in normal Fermi coordinates and so
derivatives of $g^{ij}$ and $\theta$ with respect to tangential Fermi
coordinates must vanish. On the other hand,

$\frac{\partial g^{ij}}{\partial\text{x}_{i}}(y_{0})=0$ for $i,j=q+1,...,n$
\ by the expansion of $g^{ij}$ in \textbf{Proposition 6.4.}

Thereore,

$[\frac{\partial g^{ij}}{\partial\text{x}_{j}}\frac{\partial\theta^{-\frac
{1}{2}}}{\partial\text{x}_{j}}](y_{0})=0$ for $i,j=1,...,q,q+1,...,n$

Therefore,

$\Delta\theta^{-\frac{1}{2}}(y_{0})=\frac{\partial\theta}{\partial\text{x}%
_{i}}(y_{0})\frac{\partial\theta^{-\frac{1}{2}}}{\partial\text{x}_{i}}%
(y_{0})+\frac{\partial^{2}\theta^{-\frac{1}{2}}}{\partial\text{x}_{i}^{2}%
}(y_{0})$

$\frac{\partial\theta}{\partial xi}(y_{0})=-<H,i>(y_{0})$ by (ii) of
\textbf{Table A}$_{9}$ and,

$\frac{\partial\theta^{-\frac{1}{2}}}{\partial\text{x}_{i}}(y_{0})=\frac{1}%
{2}<H,i>(y_{0})$ by (iv) of \textbf{Table A}$_{9}.$

By (x) of \textbf{Table A}$_{9},$

$\frac{\partial^{2}\theta^{-\frac{1}{2}}}{\partial\text{x}_{i}^{2}}%
(y_{0})=\frac{3}{4}<H,i>^{2}(y_{0})+\frac{1}{6}(\tau^{M}-3\tau^{P}%
+\overset{q}{\underset{\text{a}=1}{\sum}}\varrho_{\text{aa}}^{M}%
+\overset{q}{\underset{\text{a,b}=1}{\sum}}R_{\text{abab}}^{M})(y_{0})$

By (x) of \textbf{Table A}$_{9},$

$\left(  A_{26}\right)  \qquad\Delta\theta^{-\frac{1}{2}}(y_{0})$

$\qquad=-\frac{1}{2}<H,i>^{2}+\frac{3}{4}<H,i>^{2}(y_{0})+\frac{1}{6}(\tau
^{M}-3\tau^{P}+\overset{q}{\underset{\text{a}=1}{\sum}}\varrho_{\text{aa}}%
^{M}+\overset{q}{\underset{\text{a,b}=1}{\sum}}R_{\text{abab}}^{M})](y_{0})$

We simplify and obtain the same formula:

$\Delta\theta^{-\frac{1}{2}}(y_{0})=\frac{1}{4}<H,i>^{2}(y_{0})+\frac{1}%
{6}(\tau^{M}-3\tau^{P}+\overset{q}{\underset{\text{a}=1}{\sum}}\varrho
_{\text{aa}}^{M}+\overset{q}{\underset{\text{a,b}=1}{\sum}}R_{\text{abab}}%
^{M})](y_{0})$

Consequently, we have:

$\frac{1}{2}\Delta\theta^{-\frac{1}{2}}(y_{0})=\frac{1}{24}%
[\underset{i=q+1}{\overset{n}{\sum}}3<H,i>^{2}+2(\tau^{M}-3\tau^{P}%
\ +\overset{q}{\underset{\text{a=1}}{\sum}}\varrho_{\text{aa}}^{M}%
+\overset{q}{\underset{\text{a,b}=1}{\sum}}R_{\text{abab}}^{M})](y_{0})$

(iii)$\qquad\frac{1}{4}(\Delta\theta^{-\frac{1}{2}})(y_{0})(\frac{1}{2}%
\Delta\theta^{-\frac{1}{2}})(y_{0})=\frac{1}{2}(\frac{1}{2}\Delta
\theta^{-\frac{1}{2}})^{2}(y_{0})$

\qquad$\qquad=\frac{1}{2}(\frac{1}{24})^{2}[\underset{\alpha
=q+1}{\overset{n}{\sum}}3<H,\alpha>^{2}+2(\tau^{M}-3\tau^{P}%
\ +\overset{q}{\underset{\text{a=1}}{\sum}}\varrho_{\text{aa}}^{M}%
+\overset{q}{\underset{\text{a,b}=1}{\sum}}R_{\text{abab}}^{M})]^{2}(y_{0})$

$\left(  A_{27}\right)  $\qquad$\frac{1}{2}(\frac{1}{2}\Delta\theta^{-\frac
{1}{2}})^{2}(y_{0})$

$\qquad\qquad=\frac{1}{1152}[\underset{\alpha=q+1}{\overset{n}{\sum}%
}3<H,i>^{2}+2(\tau^{M}-3\tau^{P}\ +\overset{q}{\underset{\text{a=1}}{\sum}%
}\varrho_{\text{aa}}^{M}+\overset{q}{\underset{\text{a,b}=1}{\sum}%
}R_{\text{abab}}^{M})]^{2}(y_{0})$

(iv)$\qquad\frac{1}{48}(\Delta\theta^{-\frac{1}{2}})(y_{0})[6\Delta
\theta^{-\frac{1}{2}}+2\frac{\partial^{2}\theta^{\frac{1}{2}}}{\partial
x_{i}^{2}}+\frac{\partial\theta^{\frac{1}{2}}}{\partial x_{i}}(6\frac
{\partial\theta^{-\frac{1}{2}}}{\partial x_{i}}+4\frac{\partial\theta
^{\frac{1}{2}}}{\partial x_{i}}-3\overset{q}{\underset{\text{a}=1}{\sum}%
}\Gamma_{\text{aa}}^{i})](y_{0})=K$

$\qquad\qquad K=\frac{1}{2}(\frac{1}{2}\Delta\theta^{-\frac{1}{2}})^{2}%
(y_{0})+\frac{1}{12}(\frac{1}{2}\Delta\theta^{-\frac{1}{2}})(y_{0}%
)(\frac{\partial^{2}\theta^{\frac{1}{2}}}{\partial x_{i}^{2}})(y_{0})$

$\qquad\qquad+\frac{1}{24}(\frac{1}{2}\Delta\theta^{-\frac{1}{2}}%
)(y_{0})[6\frac{\partial\theta^{\frac{1}{2}}}{\partial x_{i}}\frac
{\partial\theta^{-\frac{1}{2}}}{\partial x_{i}}+4\frac{\partial\theta
^{\frac{1}{2}}}{\partial x_{i}}\frac{\partial\theta^{\frac{1}{2}}}{\partial
x_{i}}-3\frac{\partial\theta^{\frac{1}{2}}}{\partial x_{i}}%
\overset{q}{\underset{\text{a}=1}{\sum}}\Gamma_{\text{aa}}^{i})](y_{0}%
)=K_{1}+K_{2}+K_{3}$

where,

\qquad$K_{1}=\frac{1}{2}(\frac{1}{2}\Delta\theta^{-\frac{1}{2}})^{2}(y_{0})$

\qquad$K_{2}=\frac{1}{12}(\frac{1}{2}\Delta\theta^{-\frac{1}{2}})(y_{0}%
)(\frac{\partial^{2}\theta^{\frac{1}{2}}}{\partial x_{i}^{2}})(y_{0})$

\qquad$K_{3}=\frac{1}{24}(\frac{1}{2}\Delta\theta^{-\frac{1}{2}}%
)(y_{0})[6\frac{\partial\theta^{\frac{1}{2}}}{\partial x_{i}}\frac
{\partial\theta^{-\frac{1}{2}}}{\partial x_{i}}+4\frac{\partial\theta
^{\frac{1}{2}}}{\partial x_{i}}\frac{\partial\theta^{\frac{1}{2}}}{\partial
x_{i}}-3\frac{\partial\theta^{\frac{1}{2}}}{\partial x_{i}}%
\overset{q}{\underset{\text{a}=1}{\sum}}\Gamma_{\text{aa}}^{i})](y_{0})$

By (iii) here above,

$\qquad\qquad K_{1}=\frac{1}{2}(\frac{1}{2}\Delta\theta^{-\frac{1}{2}}%
)^{2}(y_{0})=\frac{1}{8}(\Delta\theta^{-\frac{1}{2}})^{2}(y_{0})$ \qquad

\qquad\qquad$\ =\frac{1}{1152}[\underset{i=q+1}{\overset{n}{\sum}}%
3<H,i>^{2}+2(\tau^{M}-3\tau^{P}\ +\overset{q}{\underset{\text{a=1}}{\sum}%
}\varrho_{\text{aa}}^{M}+\overset{q}{\underset{\text{a,b}=1}{\sum}%
}R_{\text{abab}}^{M})]^{2}(y_{0})$

By (ii) of \textbf{Table 10} above here,

$\frac{1}{2}\Delta\theta^{-\frac{1}{2}}(y_{0})=\frac{1}{24}%
[\underset{i=q+1}{\overset{n}{\sum}}3<H,i>^{2}+2(\tau^{M}-3\tau^{P}%
\ +\overset{q}{\underset{\text{a=1}}{\sum}}\varrho_{\text{aa}}^{M}%
+\overset{q}{\underset{\text{a,b}=1}{\sum}}R_{\text{abab}}^{M})](y_{0})$

By (viii) of \textbf{Table A}$_{9},$

$\frac{\partial^{2}\theta^{\frac{1}{2}}}{\partial\text{x}_{i}^{2}}(y_{0}%
)$\ $=-\frac{1}{12}[3<H,i>^{2}+2(\tau^{M}-3\tau^{P}%
+\overset{q}{\underset{\text{a}=1}{\sum}}\varrho_{\text{aa}}%
+\overset{q}{\underset{\text{a,b}=1}{\sum}}R_{\text{abab}})](y_{0})$

Therefore,

$K_{2}=\frac{1}{12}(\frac{1}{2}\Delta\theta^{-\frac{1}{2}})(y_{0}%
)(\frac{\partial^{2}\theta^{\frac{1}{2}}}{\partial x_{i}^{2}})(y_{0})$

$=-\frac{1}{3456}[3<H,i>^{2}\ +2(\tau^{M}-3\tau^{P}%
\ +\overset{q}{\underset{\text{a=1}}{\sum}}\varrho_{\text{aa}}^{M}%
+\overset{q}{\underset{\text{a,b}=1}{\sum}}R_{\text{abab}}^{M})]^{2}(y_{0})$

We have:

\qquad\ $\frac{\partial\theta^{-\frac{1}{2}}}{\partial x_{i}}(y_{0})=\frac
{1}{2}<H,i>(y_{0})$ by (iv) of \textbf{TableA}$_{9};$

\qquad\ $\frac{\partial\theta^{\frac{1}{2}}}{\partial x_{i}}(y_{0})=-\frac
{1}{2}<H,i>(y_{0})$ by (iii) of \textbf{Table A}$_{9};$

\qquad\ $\overset{q}{\underset{\text{a}=1}{\sum}}\Gamma_{\text{aa}}^{i}%
(y_{0})=\overset{q}{\underset{\text{a}=1}{\sum}}T_{\text{aa}i}(y_{0})=$
$<H,i>(y_{0})$ by (i) of \textbf{Table A}$_{7}.$

Consequently,

\qquad$\frac{\partial\theta^{\frac{1}{2}}}{\partial x_{i}}(y_{0}%
)[6\frac{\partial\theta^{-\frac{1}{2}}}{\partial x_{i}}+4\frac{\partial
\theta^{\frac{1}{2}}}{\partial x_{i}}-3\overset{q}{\underset{\text{a}=1}{\sum
}}\Gamma_{\text{aa}}^{i}](y_{0})$

$=-\frac{1}{2}<H,i>(y_{0})[6.\frac{1}{2}<H,i>-$ $4.\frac{1}{2}<H,i>-$ $3$
$<H,i>](y_{0})$

$=-\frac{1}{2}<H,i>(y_{0})[3<H,i>-2<H,i>-3$ $<H,i>](y_{0})$

$=$ $<H,i>^{2}(y_{0})$

From (ii) above,

$\frac{1}{2}\Delta\theta^{-\frac{1}{2}}(y_{0})=\frac{1}{24}%
[\underset{i=q+1}{\overset{n}{\sum}}3<H,i>^{2}+2(\tau^{M}-3\tau^{P}%
\ +\overset{q}{\underset{\text{a=1}}{\sum}}\varrho_{\text{aa}}^{M}%
+\overset{q}{\underset{\text{a,b}=1}{\sum}}R_{\text{abab}}^{M})](y_{0})$

Consequently,

$K_{3}=\frac{1}{24}(\frac{1}{2}\Delta\theta^{-\frac{1}{2}})(y_{0}%
)[6\frac{\partial\theta^{\frac{1}{2}}}{\partial x_{i}}\frac{\partial
\theta^{-\frac{1}{2}}}{\partial x_{i}}+4\frac{\partial\theta^{\frac{1}{2}}%
}{\partial x_{i}}\frac{\partial\theta^{\frac{1}{2}}}{\partial x_{i}}%
-3\frac{\partial\theta^{\frac{1}{2}}}{\partial x_{i}}%
\overset{q}{\underset{\text{a}=1}{\sum}}\Gamma_{\text{aa}}^{i})](y_{0})$

$\qquad=\frac{1}{24}\times\frac{1}{24}[\underset{i=q+1}{\overset{n}{\sum}%
}3<H,i>^{2}+2(\tau^{M}-3\tau^{P}\ +\overset{q}{\underset{\text{a=1}}{\sum}%
}\varrho_{\text{aa}}^{M}+\overset{q}{\underset{\text{a,b}=1}{\sum}%
}R_{\text{abab}}^{M})](y_{0})\times<H,i>^{2}(y_{0})$

Therefore the expression above denoted $K,$ is given by:

$K=K_{1}+K_{2}+K_{3}$

$=\frac{1}{48}(\Delta\theta^{-\frac{1}{2}})(y_{0})[6\Delta\theta^{-\frac{1}%
{2}}+2\frac{\partial^{2}\theta^{\frac{1}{2}}}{\partial x_{i}^{2}}%
+\frac{\partial\theta^{\frac{1}{2}}}{\partial x_{i}}(6\frac{\partial
\theta^{-\frac{1}{2}}}{\partial x_{i}}+4\frac{\partial\theta^{\frac{1}{2}}%
}{\partial x_{i}}-3\overset{q}{\underset{\text{a}=1}{\sum}}\Gamma_{\text{aa}%
}^{i})](y_{0})$

$=\frac{1}{1152}[\underset{i=q+1}{\overset{n}{\sum}}3<H,i>^{2}+2(\tau
^{M}-3\tau^{P}\ +\overset{q}{\underset{\text{a=1}}{\sum}}\varrho_{\text{aa}%
}^{M}+\overset{q}{\underset{\text{a,b}=1}{\sum}}R_{\text{abab}}^{M}%
)]^{2}(y_{0})$

$-\frac{1}{3456}[\underset{i=q+1}{\overset{n}{\sum}}3<H,i>^{2}\ +2(\tau
^{M}-3\tau^{P}\ +\overset{q}{\underset{\text{a=1}}{\sum}}\varrho_{\text{aa}%
}^{M}+\overset{q}{\underset{\text{a,b}=1}{\sum}}R_{\text{abab}}^{M}%
)]^{2}(y_{0})$

$+\frac{1}{576}<H,i>^{2}(y_{0})\times\lbrack\underset{i=q+1}{\overset{n}{\sum
}}3<H,i>^{2}+2(\tau^{M}-3\tau^{P}\ +\overset{q}{\underset{\text{a=1}}{\sum}%
}\varrho_{\text{aa}}^{M}+\overset{q}{\underset{\text{a,b}=1}{\sum}%
}R_{\text{abab}}^{M})](y_{0})$

$=\frac{1}{1728}[\underset{i=q+1}{\overset{n}{\sum}}3<H,i>^{2}+2(\tau
^{M}-3\tau^{P}\ +\overset{q}{\underset{\text{a=1}}{\sum}}\varrho_{\text{aa}%
}^{M}+\overset{q}{\underset{\text{a,b}=1}{\sum}}R_{\text{abab}}^{M}%
)]^{2}(y_{0})$

$+\frac{1}{576}<H,i>^{2}(y_{0})[\underset{i=q+1}{\overset{n}{\sum}}%
3<H,i>^{2}+2(\tau^{M}-3\tau^{P}\ +\overset{q}{\underset{\text{a=1}}{\sum}%
}\varrho_{\text{aa}}^{M}+\overset{q}{\underset{\text{a,b}=1}{\sum}%
}R_{\text{abab}}^{M})](y_{0})$

(v) By (v) of \textbf{Table A}$_{9},$ we have:

\qquad\ $\frac{\partial^{2}\theta}{\partial\text{x}_{i}^{2}}(y_{0})=-$
$\frac{1}{6}\underset{i=q+1}{\overset{n}{\sum}}[2\varrho_{ii}%
+4\overset{q}{\underset{\text{a}=1}{\sum}}R_{i\text{a}i\text{a}}%
-6\overset{q}{\underset{\text{a,b=1}}{\sum}}(T_{\text{aa}i}T_{\text{bb}%
i}-T_{\text{ab}i}T_{\text{ab}i})](y_{0})$

\qquad$=-\frac{1}{3}[\tau^{M}-3\tau^{P}+\overset{q}{\underset{\text{a}%
=1}{\sum}}\varrho_{\text{aa}}+\overset{q}{\underset{\text{a,b}=1}{\sum}%
}R_{\text{abab}}](y_{0})$Therefore,

$\qquad\frac{1}{4}\frac{\partial^{2}\theta}{\partial\text{x}_{i}^{2}}%
(y_{0})\times\frac{1}{8}\frac{\partial^{2}\theta}{\partial\text{x}_{j}^{2}%
}(y_{0})$

$\qquad=\frac{1}{24}\times\frac{1}{48}[2\varrho_{ii}%
+4\overset{q}{\underset{\text{a}=1}{\sum}}R_{i\text{a}i\text{a}}%
-6\overset{q}{\underset{\text{a,b=1}}{\sum}}(T_{\text{aa}i}T_{\text{bb}%
i}-T_{\text{ab}i}T_{\text{ab}i})](y_{0})$ $\ $

$\qquad\qquad\times\lbrack2\varrho_{jj}+4\overset{q}{\underset{\text{a}%
=1}{\sum}}R_{j\text{a}j\text{a}}-6\overset{q}{\underset{\text{a,b=1}}{\sum}%
}(T_{\text{aa}j}T_{\text{bb}j}-T_{\text{ab}j}T_{\text{ab}j})](y_{0})\qquad$

\qquad$=(\frac{1}{4})(\frac{1}{8})(-\frac{1}{3})(-\frac{1}{3})[\tau^{M}%
-3\tau^{P}+\overset{q}{\underset{\text{a}=1}{\sum}}\varrho_{\text{aa}%
}+\overset{q}{\underset{\text{a,b}=1}{\sum}}R_{\text{abab}}]^{2}(y_{0})$

\qquad$=\frac{1}{288}[\tau^{M}-3\tau^{P}+\overset{q}{\underset{\text{a}%
=1}{\sum}}\varrho_{\text{aa}}+\overset{q}{\underset{\text{a,b}=1}{\sum}%
}R_{\text{abab}}]^{2}(y_{0})$\qquad

(vi) The expression for $\frac{1}{24}(\frac{\partial^{2}\theta^{\frac{1}{2}}%
}{\partial x_{i}^{2}}\Delta\theta^{-\frac{1}{2}})(y_{0})$ has already been
given in (v) above and so,$\qquad$

$\left(  A_{28}\right)  $\qquad A$_{3211}=\frac{1}{24}(\frac{\partial
^{2}\theta^{\frac{1}{2}}}{\partial x_{i}^{2}}\Delta\theta^{-\frac{1}{2}%
})(y_{0})$

$\qquad=-\frac{1}{3456}[3<H,i>^{2}\ +2(\tau^{M}-3\tau^{P}%
\ +\overset{q}{\underset{\text{a=1}}{\sum}}\varrho_{\text{aa}}^{M}%
+\overset{q}{\underset{\text{a,b}=1}{\sum}}R_{\text{abab}}^{M})]^{2}(y_{0}%
)$\ $\qquad\qquad$

(vii) We next consider the important but complicated term below: For
$i=q+1,...,n,$

$\qquad$A$_{3212}$ $=\frac{1}{24}[\frac{\partial^{2}}{\partial x_{i}^{2}%
}(\Delta\theta^{-\frac{1}{2}})](y_{0})$

$\Delta\theta^{-\frac{1}{2}}=\theta^{-1}[\frac{\partial}{\partial\text{x}_{j}%
}(\theta g^{jk}\frac{\partial\theta^{-\frac{1}{2}}}{\partial\text{x}_{k}%
})]=\theta^{-1}[\frac{\partial\theta}{\partial\text{x}_{j}}(g^{jk}%
\frac{\partial\theta^{-\frac{1}{2}}}{\partial\text{x}_{k}})+\theta
(\frac{\partial g^{jk}}{\partial\text{x}_{j}}\frac{\partial\theta^{-\frac
{1}{2}}}{\partial\text{x}_{k}}+g^{jk}\frac{\partial^{2}\theta^{-\frac{1}{2}}%
}{\partial\text{x}_{j}\partial\text{x}_{k}})]$

We will use the following formula: for any smooth functions
$f,g:M\longrightarrow R,$

we have: $\frac{\partial^{2}}{\partial x\partial y}(fg)=\frac{\partial^{2}%
f}{\partial x\partial y}g+f\frac{\partial^{2}g}{\partial x\partial y}%
+2\frac{\partial f}{\partial x}\frac{\partial g}{\partial y}$ in any chart of M.

We set: $f=\theta^{-1}$ and $g=[\frac{\partial}{\partial\text{x}_{j}}(\theta
g^{jk}\frac{\partial\theta^{-\frac{1}{2}}}{\partial\text{x}_{k}}%
)]=[\frac{\partial\theta}{\partial\text{x}_{j}}(g^{jk}\frac{\partial
\theta^{-\frac{1}{2}}}{\partial\text{x}_{k}})+\theta(\frac{\partial g^{jk}%
}{\partial\text{x}_{j}}\frac{\partial\theta^{-\frac{1}{2}}}{\partial
\text{x}_{k}}+g^{jk}\frac{\partial^{2}\theta^{-\frac{1}{2}}}{\partial
\text{x}_{j}\partial\text{x}_{k}})],$

and have:

Since $\theta(y_{0})=1,$ we have:

$\frac{\partial^{2}}{\partial x_{i}^{2}}(\Delta\theta^{-\frac{1}{2}}%
)(y_{0})=\frac{\partial^{2}\theta^{-1}}{\partial x_{i}^{2}}(y_{0}%
)[\frac{\partial\theta}{\partial\text{x}_{j}}(g^{jk}\frac{\partial
\theta^{-\frac{1}{2}}}{\partial\text{x}_{k}})+\theta(\frac{\partial g^{jk}%
}{\partial\text{x}_{j}}\frac{\partial\theta^{-\frac{1}{2}}}{\partial
\text{x}_{k}}+g^{jk}\frac{\partial^{2}\theta^{-\frac{1}{2}}}{\partial
\text{x}_{j}\partial\text{x}_{k}})](y_{0})$

$+\frac{\partial^{2}}{\partial x_{i}^{2}}[\frac{\partial\theta}{\partial
\text{x}_{j}}(g^{jk}\frac{\partial\theta^{-\frac{1}{2}}}{\partial\text{x}_{k}%
})+\theta(\frac{\partial g^{jk}}{\partial\text{x}_{j}}\frac{\partial
\theta^{-\frac{1}{2}}}{\partial\text{x}_{k}}+g^{jk}\frac{\partial^{2}%
\theta^{-\frac{1}{2}}}{\partial\text{x}_{j}\partial\text{x}_{k}})](y_{0})$

$+2\frac{\partial\theta^{-1}}{\partial x_{i}}(y_{0})\frac{\partial}{\partial
x_{i}}[\frac{\partial\theta}{\partial\text{x}_{j}}(g^{jk}\frac{\partial
\theta^{-\frac{1}{2}}}{\partial\text{x}_{k}})+\theta(\frac{\partial g^{jk}%
}{\partial\text{x}_{j}}\frac{\partial\theta^{-\frac{1}{2}}}{\partial
\text{x}_{k}}+g^{jk}\frac{\partial^{2}\theta^{-\frac{1}{2}}}{\partial
\text{x}_{j}\partial\text{x}_{k}})](y_{0})=L_{1}+L_{2}+L_{3}$

where,

$L_{1}=\frac{\partial^{2}\theta^{-1}}{\partial x_{i}^{2}}(y_{0})[\frac
{\partial\theta}{\partial\text{x}_{j}}(g^{jk}\frac{\partial\theta^{-\frac
{1}{2}}}{\partial\text{x}_{k}})+\theta(\frac{\partial g^{jk}}{\partial
\text{x}_{j}}\frac{\partial\theta^{-\frac{1}{2}}}{\partial\text{x}_{k}}%
+g^{jk}\frac{\partial^{2}\theta^{-\frac{1}{2}}}{\partial\text{x}_{j}%
\partial\text{x}_{k}})](y_{0})$

$L_{2}=\frac{\partial^{2}}{\partial x_{i}^{2}}[\frac{\partial\theta}%
{\partial\text{x}_{j}}(g^{jk}\frac{\partial\theta^{-\frac{1}{2}}}%
{\partial\text{x}_{k}})+\theta(\frac{\partial g^{jk}}{\partial\text{x}_{j}%
}\frac{\partial\theta^{-\frac{1}{2}}}{\partial\text{x}_{k}}+g^{jk}%
\frac{\partial^{2}\theta^{-\frac{1}{2}}}{\partial\text{x}_{j}\partial
\text{x}_{k}})](y_{0})$

$L_{3}=2\frac{\partial\theta^{-1}}{\partial x_{i}}(y_{0})\frac{\partial
}{\partial x_{i}}[\frac{\partial\theta}{\partial\text{x}_{j}}(g^{jk}%
\frac{\partial\theta^{-\frac{1}{2}}}{\partial\text{x}_{k}})+\theta
(\frac{\partial g^{jk}}{\partial\text{x}_{j}}\frac{\partial\theta^{-\frac
{1}{2}}}{\partial\text{x}_{k}}+g^{jk}\frac{\partial^{2}\theta^{-\frac{1}{2}}%
}{\partial\text{x}_{j}\partial\text{x}_{k}})](y_{0})$

Since $\frac{\partial^{2}\theta^{-1}}{\partial x_{i}^{2}}(y_{0})=2(\frac
{\partial\theta}{\partial\text{x}_{i}}\frac{\partial\theta}{\partial
\text{x}_{i}})(y_{0})-\frac{\partial^{2}\theta}{\partial\text{x}_{i}^{2}%
}(y_{0});$ $\theta(y_{0})=1$ and $g^{jk}(y_{0})=\delta^{jk},$ we have:

$L_{1}=[2(\frac{\partial\theta}{\partial\text{x}_{i}}\frac{\partial\theta
}{\partial\text{x}_{i}})(y_{0})-\frac{\partial^{2}\theta}{\partial\text{x}%
_{i}^{2}}(y_{0})][\frac{\partial\theta}{\partial\text{x}_{j}}\frac
{\partial\theta^{-\frac{1}{2}}}{\partial\text{x}_{j}}+\frac{\partial g^{jk}%
}{\partial\text{x}_{j}}\frac{\partial\theta^{-\frac{1}{2}}}{\partial
\text{x}_{k}}+\frac{\partial^{2}\theta^{-\frac{1}{2}}}{\partial\text{x}%
_{j}^{2}}](y_{0})$

Expansions of $\theta$ and $g^{jk}$ in \textbf{Chapter 6} are carried out in
normal coordinates,

hence differentiation of $\theta$ and $g^{jk}$ with tangential coordinates vanish.

In the particular case of $\frac{\partial g^{jk}}{\partial\text{x}_{j}}%
\frac{\partial\theta^{-\frac{1}{2}}}{\partial\text{x}_{k}},$ the indics $j$
and $k$ must be those of normal

coordinates: $j,k=q+1,...,n.$ In this case the expansion of $g^{jk}$ given in
\textbf{Proposition }$6.4$

shows that $\frac{\partial g^{jk}}{\partial\text{x}_{j}}(y_{0})=0.$ Therefore
$\frac{\partial g^{jk}}{\partial\text{x}_{j}}(y_{0})\frac{\partial
\theta^{-\frac{1}{2}}}{\partial\text{x}_{k}}(y_{0})$ must vanish.

This argument will be valid in other contexts below. We have:

$L_{1}=[2(\frac{\partial\theta}{\partial\text{x}_{i}}\frac{\partial\theta
}{\partial\text{x}_{i}})(y_{0})-\frac{\partial^{2}\theta}{\partial\text{x}%
_{i}^{2}}(y_{0})][\frac{\partial\theta}{\partial\text{x}_{j}}\frac
{\partial\theta^{-\frac{1}{2}}}{\partial\text{x}_{j}}+\frac{\partial^{2}%
\theta^{-\frac{1}{2}}}{\partial\text{x}_{j}^{2}}](y_{0})$

We have:\ $\frac{\partial\theta}{\partial\text{x}_{j}}(y_{0})=-<H,i>(y_{0})$
\ and $\frac{\partial\theta^{-\frac{1}{2}}}{\partial\text{x}_{j}}(y_{0}%
)=\frac{1}{2}<H,i>$

\ $\frac{\partial^{2}\theta}{\partial\text{x}_{i}^{2}}(y_{0})=-\frac{1}%
{3}[\tau^{M}-3\tau^{P}+\overset{q}{\underset{\text{a}=1}{\sum}}\varrho
_{\text{aa}}+\overset{q}{\underset{\text{a,b}=1}{\sum}}R_{\text{abab}}%
](y_{0})$ by (v) of \textbf{Table A}$_{9}$

$\frac{\partial^{2}\theta^{-\frac{1}{2}}}{\partial\text{x}_{i}^{2}}%
(y_{0})=\frac{3}{4}<H,i>^{2}(y_{0})+\frac{1}{6}(\tau^{M}-3\tau^{P}%
+\overset{q}{\underset{\text{a}=1}{\sum}}\varrho_{\text{aa}}^{M}%
+\overset{q}{\underset{\text{a,b}=1}{\sum}}R_{\text{abab}}^{M})(y_{0})$ by (x)
of \textbf{Table A}$_{9}.$

\textbf{ A}$_{9}.$ Since $g^{jk}(y_{0})=\delta^{jk},$ we have:

\qquad$L_{1}=[2<H,i>^{2}(y_{0})+\frac{1}{3}(\tau^{M}-3\tau^{P}%
+\overset{q}{\underset{\text{a}=1}{\sum}}\varrho_{\text{aa}}%
+\overset{q}{\underset{\text{a,b}=1}{\sum}}R_{\text{abab}})](y_{0})$

\qquad$\times\lbrack\frac{1}{4}<H,j>^{2}(y_{0})+\frac{1}{6}(\tau^{M}-3\tau
^{P}+\overset{q}{\underset{\text{a}=1}{\sum}}\varrho_{\text{aa}}%
^{M}+\overset{q}{\underset{\text{a,b}=1}{\sum}}R_{\text{abab}}^{M})](y_{0})$

\begin{center}
\qquad\qquad\qquad\qquad\qquad\qquad\qquad\qquad\qquad\qquad\qquad\qquad
\qquad$\blacksquare$
\end{center}

We next compute:

$L_{2}=\frac{\partial^{2}}{\partial x_{i}^{2}}[\frac{\partial\theta}%
{\partial\text{x}_{j}}(g^{jk}\frac{\partial\theta^{-\frac{1}{2}}}%
{\partial\text{x}_{k}})+\theta(\frac{\partial g^{jk}}{\partial\text{x}_{j}%
}\frac{\partial\theta^{-\frac{1}{2}}}{\partial\text{x}_{k}}+g^{jk}%
\frac{\partial^{2}\theta^{-\frac{1}{2}}}{\partial\text{x}_{j}\partial
\text{x}_{k}})](y_{0})=L_{21}+L_{22}+L_{23}$where,

$L_{21}=\frac{\partial^{2}}{\partial x_{i}^{2}}[\frac{\partial\theta}%
{\partial\text{x}_{j}}(g^{jk}\frac{\partial\theta^{-\frac{1}{2}}}%
{\partial\text{x}_{k}}](y_{0});$ $L_{22}=\frac{\partial^{2}}{\partial
x_{i}^{2}}[\theta\frac{\partial g^{jk}}{\partial\text{x}_{j}}\frac
{\partial\theta^{-\frac{1}{2}}}{\partial\text{x}_{k}}](y_{0});$ $L_{23}%
=\frac{\partial^{2}}{\partial x_{i}^{2}}[\theta g^{jk}\frac{\partial^{2}%
\theta^{-\frac{1}{2}}}{\partial\text{x}_{j}\partial\text{x}_{k}}](y_{0})$

We have,

$L_{21}=\frac{\partial^{2}}{\partial x_{i}^{2}}[\frac{\partial\theta}%
{\partial\text{x}_{j}}(g^{jk}\frac{\partial\theta^{-\frac{1}{2}}}%
{\partial\text{x}_{k}})](y_{0})=2\frac{\partial^{2}\theta}{\partial
x_{i}\partial\text{x}_{j}}(y_{0})\frac{\partial}{\partial\text{x}_{i}}%
(g^{jk}\frac{\partial\theta^{-\frac{1}{2}}}{\partial\text{x}_{k}}%
)(y_{0})+\frac{\partial^{3}\theta}{\partial x_{i}^{2}\partial\text{x}_{j}%
}(y_{0})(g^{jk}\frac{\partial\theta^{-\frac{1}{2}}}{\partial\text{x}_{k}%
})(y_{0})$

$\qquad+\frac{\partial\theta}{\partial\text{x}_{j}}(y_{0})\frac{\partial^{2}%
}{\partial x_{i}^{2}}(g^{jk}\frac{\partial\theta^{-\frac{1}{2}}}%
{\partial\text{x}_{k}})(y_{0})=L_{211}+L_{212}+L_{213}$

where,

$L_{211}=2\frac{\partial^{2}\theta}{\partial x_{i}\partial\text{x}_{j}}%
(y_{0})\frac{\partial}{\partial\text{x}_{i}}(g^{jk}\frac{\partial
\theta^{-\frac{1}{2}}}{\partial\text{x}_{k}})(y_{0})=2\frac{\partial^{2}%
\theta}{\partial x_{i}\partial\text{x}_{j}}(y_{0})[\frac{\partial\text{g}%
^{jk}}{\partial\text{x}_{i}}\frac{\partial\theta^{-\frac{1}{2}}}%
{\partial\text{x}_{k}}+g^{jk}\frac{\partial^{2}\theta^{-\frac{1}{2}}}%
{\partial\text{x}_{i}\partial\text{x}_{k}}](y_{0})$

To ensure that the expression on the RHS\ of the last equation does not vanish,

we take $i,j,k=q+1,...,n.$ In this case the expansion of $g^{jk}$ shows that
$\frac{\partial g^{jk}}{\partial\text{x}_{i}}(y_{0})=0.$

On the other hand, $g^{jk}(y_{0})=\delta^{ij}$ and we have:

$L_{211}=2\frac{\partial^{2}\theta}{\partial x_{i}\partial\text{x}_{j}}%
(y_{0})[g^{jk}\frac{\partial^{2}\theta^{-\frac{1}{2}}}{\partial\text{x}%
_{i}\partial\text{x}_{k}}](y_{0})=2\frac{\partial^{2}\theta}{\partial
x_{i}\partial\text{x}_{j}}(y_{0})\frac{\partial^{2}\theta^{-\frac{1}{2}}%
}{\partial\text{x}_{i}\partial\text{x}_{j}}(y_{0})$

From (v) and (ix) of \textbf{Appendix A}$_{9}:$

$\frac{\partial^{2}\theta}{\partial\text{x}_{i}\partial\text{x}_{j}}%
(y_{0})=-\frac{1}{6}[2\varrho_{ij}+$ $\overset{q}{\underset{\text{a}=1}{4\sum
}}R_{i\text{a}j\text{a}}-3\overset{q}{\underset{\text{a,b=1}}{\sum}%
}(T_{\text{aa}i}T_{\text{bb}j}-T_{\text{ab}i}T_{\text{ab}j})$

$-3\overset{q}{\underset{\text{a,b=1}}{\sum}}(T_{\text{aa}j}T_{\text{bb}%
i}-T_{\text{ab}j}T_{\text{ab}i}](y_{0})$

$\frac{\partial^{2}\theta^{-\frac{1}{2}}}{\partial\text{x}_{i}\partial
\text{x}_{j}}(y_{0})$ \ $=\frac{3}{4}<H,i>(y_{0})<H,j>(y_{0})$

$+\frac{1}{12}[2\varrho_{ij}+$ $\overset{q}{\underset{\text{a}=1}{4\sum}%
}R_{i\text{a}j\text{a}}-3\overset{q}{\underset{\text{a,b=1}}{\sum}%
}(T_{\text{aa}i}T_{\text{bb}j}-T_{\text{ab}i}T_{\text{ab}j}%
)-3\overset{q}{\underset{\text{a,b=1}}{\sum}}(T_{\text{aa}j}T_{\text{bb}%
i}-T_{\text{ab}j}T_{\text{ab}i}](y_{0})$

Consequently,

$L_{211}=-\frac{1}{4}[<H,i><H,j>](y_{0})$

$\times\lbrack2\varrho_{ij}+$ $\overset{q}{\underset{\text{a}=1}{4\sum}%
}R_{i\text{a}j\text{a}}-3\overset{q}{\underset{\text{a,b=1}}{\sum}%
}(T_{\text{aa}i}T_{\text{bb}j}-T_{\text{ab}i}T_{\text{ab}j}%
)-3\overset{q}{\underset{\text{a,b=1}}{\sum}}(T_{\text{aa}j}T_{\text{bb}%
i}-T_{\text{ab}j}T_{\text{ab}i}](y_{0})$

$-\frac{1}{36}[2\varrho_{ij}+$ $\overset{q}{\underset{\text{a}=1}{4\sum}%
}R_{i\text{a}j\text{a}}-3\overset{q}{\underset{\text{a,b=1}}{\sum}%
}(T_{\text{aa}i}T_{\text{bb}j}-T_{\text{ab}i}T_{\text{ab}j}%
)-3\overset{q}{\underset{\text{a,b=1}}{\sum}}(T_{\text{aa}j}T_{\text{bb}%
i}-T_{\text{ab}j}T_{\text{ab}i}]^{2}(y_{0})$

\qquad\qquad\qquad\qquad\qquad\qquad\qquad\qquad\qquad\qquad\qquad\qquad
\qquad\qquad\qquad\qquad\qquad\qquad\qquad\qquad$\blacksquare$

Next we have:

$L_{212}=\frac{\partial^{3}\theta}{\partial x_{i}^{2}\partial\text{x}_{j}%
}(y_{0})(g^{jk}\frac{\partial\theta^{-\frac{1}{2}}}{\partial\text{x}_{k}%
})(y_{0})=\frac{\partial^{3}\theta}{\partial x_{i}^{2}\partial\text{x}_{j}%
}(y_{0})\frac{\partial\theta^{-\frac{1}{2}}}{\partial\text{x}_{j}}(y_{0})$

$\frac{\partial\theta^{-\frac{1}{2}}}{\partial\text{x}_{j}}(y_{0})=\frac{1}%
{2}[<H,j>](y_{0}).$

The expression for $\frac{\partial^{3}\theta}{\partial x_{i}^{2}%
\partial\text{x}_{j}}(y_{0})$ is taken from (xii) of \textbf{Appendix A}%
$_{9}:$

$L_{212}=\frac{\partial\theta^{-\frac{1}{2}}}{\partial\text{x}_{j}}%
(y_{0})\frac{\partial^{3}\theta}{\partial x_{i}^{2}\partial\text{x}_{j}}%
(y_{0})$

$L_{212}=-\frac{1}{12}[<H,j>](y_{0})\times\lbrack\{\nabla_{i}\varrho
_{ij}-2\varrho_{ij}<H,i>+\overset{q}{\underset{\text{a}=1}{\sum}}(\nabla
_{i}R_{\text{a}i\text{a}j}-4R_{i\text{a}j\text{a}}<H,i>)$

$+4\overset{q}{\underset{\text{a,b=1}}{\sum}}R_{i\text{a}j\text{b}%
}T_{\text{ab}i}+2\overset{q}{\underset{\text{a,b,c=1}}{\sum}}(T_{\text{aa}%
i}T_{\text{bb}j}T_{\text{cc}i}-3T_{\text{aa}i}T_{\text{bc}j}T_{\text{bc}%
i}+2T_{\text{ab}i}T_{\text{bc}j}T_{\text{ac}i})](y_{0})$\qquad\qquad
\qquad\qquad\qquad\ \ 

$-\frac{1}{12}[<H,j>](y_{0})\times\lbrack\nabla_{j}\varrho_{ii}-2\varrho
_{ij}<H,i>+\overset{q}{\underset{\text{a}=1}{\sum}}(\nabla_{j}R_{\text{a}%
i\text{a}i}-4R_{i\text{a}j\text{a}}<H,i>)$

$+4\overset{q}{\underset{\text{a,b=1}}{\sum}}R_{j\text{a}i\text{b}%
}T_{\text{ab}i}+2\overset{q}{\underset{\text{a,b,c=1}}{\sum}}(T_{\text{aa}%
j}T_{\text{bb}i}T_{\text{cc}i}-3T_{\text{aa}j}T_{\text{bc}i}T_{\text{bc}%
i}+2T_{\text{ab}j}T_{\text{bc}i}T_{\text{ac}i})](y_{0})$

$-\frac{1}{12}[<H,j>](y_{0})\times\lbrack\nabla_{i}\varrho_{ij}-2\varrho
_{ii}<H,j>+\overset{q}{\underset{\text{a}=1}{\sum}}(\nabla_{i}R_{\text{a}%
i\text{a}j}-4R_{i\text{a}i\text{a}}<H,j>)$

$+4\overset{q}{\underset{\text{a,b=1}}{\sum}}R_{i\text{a}i\text{b}%
}T_{\text{ab}j}+2\overset{q}{\underset{\text{a,b,c}=1}{\sum}}(T_{\text{aa}%
i}T_{\text{bb}i}T_{\text{cc}j}-3T_{\text{aa}i}T_{\text{bc}i}T_{\text{bc}%
j}+2T_{\text{ab}i}T_{\text{bc}i}T_{\text{ac}j})](y_{0})$

\qquad\qquad\qquad\qquad\qquad\qquad\qquad\qquad\qquad\qquad\qquad\qquad
\qquad\qquad\qquad\qquad$\blacksquare$

The last expression in this subset is:

$L_{213}=\frac{\partial\theta}{\partial\text{x}_{j}}(y_{0})\frac{\partial^{2}%
}{\partial x_{i}^{2}}(g^{jk}\frac{\partial\theta^{-\frac{1}{2}}}%
{\partial\text{x}_{k}})(y_{0})=\frac{\partial\theta}{\partial\text{x}_{j}%
}(y_{0})[2\frac{\partial\text{g}^{jk}}{\partial\text{x}_{i}}\frac{\partial
^{2}\theta^{-\frac{1}{2}}}{\partial\text{x}_{i}\partial\text{x}_{k}}%
+\frac{\partial^{2}\text{g}^{jk}}{\partial x_{i}^{2}}\frac{\partial
\theta^{-\frac{1}{2}}}{\partial\text{x}_{k}}$

$\qquad+g^{jk}\frac{\partial^{3}\theta^{-\frac{1}{2}}}{\partial x_{i}%
^{2}\partial\text{x}_{k}}](y_{0})$

$\qquad=\frac{\partial\theta}{\partial\text{x}_{j}}(y_{0})[2\frac
{\partial\text{g}^{jk}}{\partial\text{x}_{i}}\frac{\partial^{2}\theta
^{-\frac{1}{2}}}{\partial\text{x}_{i}\partial\text{x}_{k}}+\frac{\partial
^{2}\text{g}^{jk}}{\partial x_{i}^{2}}\frac{\partial\theta^{-\frac{1}{2}}%
}{\partial\text{x}_{k}}+\frac{\partial^{3}\theta^{-\frac{1}{2}}}{\partial
x_{i}^{2}\partial\text{x}_{j}}](y_{0})$

The same argument above shows that the indices must all be normal:

$\ i,j,k=q+1,...,n.$

In this case, $\frac{\partial\text{g}^{jk}}{\partial\text{x}_{i}}(y_{0}%
)=0$\textbf{ and }$\frac{\partial^{2}\text{g}^{jk}}{\partial x_{i}^{2}}%
=\frac{2}{3}R_{ijik}(y_{0})$ by (iii) of \textbf{Appendix A}$_{2}.$

The expression of $\frac{\partial^{3}\theta^{-\frac{1}{2}}}{\partial x_{i}%
^{2}\partial\text{x}_{j}}(y_{0})$ is given by (xvi) of \textbf{Appendix
A}$_{9}.$

$L_{213}=\frac{\partial\theta}{\partial\text{x}_{j}}(y_{0})[\frac{\partial
^{2}\text{g}^{jk}}{\partial x_{i}^{2}}\frac{\partial\theta^{-\frac{1}{2}}%
}{\partial\text{x}_{k}}+\frac{\partial^{3}\theta^{-\frac{1}{2}}}{\partial
x_{i}^{2}\partial\text{x}_{j}}](y_{0})$

\qquad$=-\frac{1}{3}[<H,j><H,k>](y_{0})R_{ijik}(y_{0})\qquad\qquad L_{213}$

$\qquad-$ $\frac{15}{8}[<H,i>^{2}<H,j>^{2}](y_{0})\qquad$

$\qquad-\frac{1}{4}<H,i><H,j>[2\varrho_{ij}+\overset{q}{\underset{\text{a}%
=1}{4\sum}}R_{i\text{a}j\text{a}}-3\overset{q}{\underset{\text{a,b=1}}{\sum}%
}(T_{\text{aa}i}T_{\text{bb}j}-T_{\text{ab}i}T_{\text{ab}j})$

$\qquad-3\overset{q}{\underset{\text{a,b=1}}{\sum}}(T_{\text{aa}j}%
T_{\text{bb}i}-T_{\text{ab}j}T_{\text{ab}i}](y_{0})$

$\qquad-\frac{1}{4}<H,j>^{2}[\tau^{M}\ -3\tau^{P}+\ \underset{\text{a}%
=1}{\overset{\text{q}}{\sum}}\varrho_{\text{aa}}^{M}+$
$\overset{q}{\underset{\text{a},\text{b}=1}{\sum}}R_{\text{abab}}^{M}$
$](y_{0})$

$\qquad+\frac{1}{12}<H,j>[\nabla_{i}\varrho_{ij}-2\varrho_{ij}%
<H,i>+\overset{q}{\underset{\text{a}=1}{\sum}}(\nabla_{i}R_{\text{a}%
i\text{a}j}-4R_{i\text{a}j\text{a}}<H,i>)$

$\qquad+4\overset{q}{\underset{\text{a,b=1}}{\sum}}R_{i\text{a}j\text{b}%
}T_{\text{ab}i}$

$\qquad+2\overset{q}{\underset{\text{a,b,c=1}}{\sum}}(T_{\text{aa}%
i}T_{\text{bb}j}T_{\text{cc}i}-3T_{\text{aa}i}T_{\text{bc}j}T_{\text{bc}%
i}+2T_{\text{ab}i}T_{\text{bc}j}T_{\text{ac}i})](y_{0})$\qquad\qquad
\qquad\qquad\qquad\ \ 

$\qquad+\frac{1}{12}<H,j>[\nabla_{j}\varrho_{ii}-2\varrho_{ij}%
<H,i>+\overset{q}{\underset{\text{a}=1}{\sum}}(\nabla_{j}R_{\text{a}%
i\text{a}i}-4R_{i\text{a}j\text{a}}<H,i>)$

$\qquad+4\overset{q}{\underset{\text{a,b=1}}{\sum}}R_{j\text{a}i\text{b}%
}T_{\text{ab}i}$

$\qquad+2\overset{q}{\underset{\text{a,b,c=1}}{\sum}}(T_{\text{aa}%
j}T_{\text{bb}i}T_{\text{cc}i}-3T_{\text{aa}j}T_{\text{bc}i}T_{\text{bc}%
i}+2T_{\text{ab}j}T_{\text{bc}i}T_{\text{ac}i})](y_{0})$

$\qquad+\frac{1}{12}<H,j>[\nabla_{i}\varrho_{ij}-2\varrho_{ii}%
<H,j>+\overset{q}{\underset{\text{a}=1}{\sum}}(\nabla_{i}R_{\text{a}%
i\text{a}j}-4R_{i\text{a}i\text{a}}<H,j>)$

$\qquad+4\overset{q}{\underset{\text{a,b=1}}{\sum}}R_{i\text{a}i\text{b}%
}T_{\text{ab}j}$

$\qquad+2\overset{q}{\underset{\text{a,b,c}=1}{\sum}}(T_{\text{aa}%
i}T_{\text{bb}i}T_{\text{cc}j}-3T_{\text{aa}i}T_{\text{bc}i}T_{\text{bc}%
j}+2T_{\text{ab}i}T_{\text{bc}i}T_{\text{ac}j})](y_{0})$

\qquad\qquad\qquad\qquad\qquad\qquad\qquad\qquad\qquad\qquad\qquad\qquad
\qquad\qquad\qquad\qquad\qquad\qquad$\blacksquare$

Therefore,

$\qquad L_{21}=L_{211}+L_{212}+L_{213}$

$=-\frac{1}{4}[2\varrho_{ij}+$ $\overset{q}{\underset{\text{a}=1}{4\sum}%
}R_{i\text{a}j\text{a}}-3\overset{q}{\underset{\text{a,b=1}}{\sum}%
}(T_{\text{aa}i}T_{\text{bb}j}-T_{\text{ab}i}T_{\text{ab}j}%
)-3\overset{q}{\underset{\text{a,b=1}}{\sum}}(T_{\text{aa}j}T_{\text{bb}%
i}-T_{\text{ab}j}T_{\text{ab}i}](y_{0})\qquad L_{211}$

$\qquad\ \ \ \ \ \times\lbrack<H,i><H,j>](y_{0})$

$\qquad-\frac{1}{36}[2\varrho_{ij}+$ $\overset{q}{\underset{\text{a}=1}{4\sum
}}R_{i\text{a}j\text{a}}-3\overset{q}{\underset{\text{a,b=1}}{\sum}%
}(T_{\text{aa}i}T_{\text{bb}j}-T_{\text{ab}i}T_{\text{ab}j}%
)-3\overset{q}{\underset{\text{a,b=1}}{\sum}}(T_{\text{aa}j}T_{\text{bb}%
i}-T_{\text{ab}j}T_{\text{ab}i}]^{2}(y_{0})$

\qquad$-\frac{1}{12}[<H,j>](y_{0})\times\lbrack\{\nabla_{i}\varrho
_{ij}-2\varrho_{ij}<H,i>+\overset{q}{\underset{\text{a}=1}{\sum}}(\nabla
_{i}R_{\text{a}i\text{a}j}-4R_{i\text{a}j\text{a}}<H,i>)\qquad L_{212}$

$\qquad+4\overset{q}{\underset{\text{a,b=1}}{\sum}}R_{i\text{a}j\text{b}%
}T_{\text{ab}i}+2\overset{q}{\underset{\text{a,b,c=1}}{\sum}}(T_{\text{aa}%
i}T_{\text{bb}j}T_{\text{cc}i}-3T_{\text{aa}i}T_{\text{bc}j}T_{\text{bc}%
i}+2T_{\text{ab}i}T_{\text{bc}j}T_{\text{ac}i})](y_{0})$\qquad\qquad
\qquad\qquad\qquad\ \ 

$\qquad-\frac{1}{12}[<H,j>](y_{0})\times\lbrack\nabla_{j}\varrho_{ii}%
-2\varrho_{ij}<H,i>+\overset{q}{\underset{\text{a}=1}{\sum}}(\nabla
_{j}R_{\text{a}i\text{a}i}-4R_{i\text{a}j\text{a}}<H,i>)$

$\qquad+4\overset{q}{\underset{\text{a,b=1}}{\sum}}R_{j\text{a}i\text{b}%
}T_{\text{ab}i}+2\overset{q}{\underset{\text{a,b,c=1}}{\sum}}(T_{\text{aa}%
j}T_{\text{bb}i}T_{\text{cc}i}-3T_{\text{aa}j}T_{\text{bc}i}T_{\text{bc}%
i}+2T_{\text{ab}j}T_{\text{bc}i}T_{\text{ac}i})](y_{0})$

$\qquad-\frac{1}{12}[<H,j>](y_{0})\times\lbrack\nabla_{i}\varrho_{ij}%
-2\varrho_{ii}<H,j>+\overset{q}{\underset{\text{a}=1}{\sum}}(\nabla
_{i}R_{\text{a}i\text{a}j}-4R_{i\text{a}i\text{a}}<H,j>)$

$\qquad+4\overset{q}{\underset{\text{a,b=1}}{\sum}}R_{i\text{a}i\text{b}%
}T_{\text{ab}j}+2\overset{q}{\underset{\text{a,b,c}=1}{\sum}}(T_{\text{aa}%
i}T_{\text{bb}i}T_{\text{cc}j}-3T_{\text{aa}i}T_{\text{bc}i}T_{\text{bc}%
j}+2T_{\text{ab}i}T_{\text{bc}i}T_{\text{ac}j})](y_{0})$

$\qquad-\frac{1}{3}[<H,j><H,k>](y_{0})R_{ijik}(y_{0})-$ $\frac{15}%
{8}[<H,i>^{2}<H,j>^{2}](y_{0})\qquad L_{213}$

$\qquad-\frac{1}{4}<H,i><H,j>[2\varrho_{ij}+\overset{q}{\underset{\text{a}%
=1}{4\sum}}R_{i\text{a}j\text{a}}-3\overset{q}{\underset{\text{a,b=1}}{\sum}%
}(T_{\text{aa}i}T_{\text{bb}j}-T_{\text{ab}i}T_{\text{ab}j})$

$\qquad-3\overset{q}{\underset{\text{a,b=1}}{\sum}}(T_{\text{aa}j}%
T_{\text{bb}i}-T_{\text{ab}j}T_{\text{ab}i}](y_{0})$

$\qquad-\frac{1}{4}<H,j>^{2}[\tau^{M}\ -3\tau^{P}+\ \underset{\text{a}%
=1}{\overset{\text{q}}{\sum}}\varrho_{\text{aa}}^{M}+$
$\overset{q}{\underset{\text{a},\text{b}=1}{\sum}}R_{\text{abab}}^{M}$
$](y_{0})$

$\qquad+\frac{1}{12}<H,j>[\nabla_{i}\varrho_{ij}-2\varrho_{ij}%
<H,i>+\overset{q}{\underset{\text{a}=1}{\sum}}(\nabla_{i}R_{\text{a}%
i\text{a}j}-4R_{i\text{a}j\text{a}}<H,i>)$

$\qquad+4\overset{q}{\underset{\text{a,b=1}}{\sum}}R_{i\text{a}j\text{b}%
}T_{\text{ab}i}$

$\qquad+2\overset{q}{\underset{\text{a,b,c=1}}{\sum}}(T_{\text{aa}%
i}T_{\text{bb}j}T_{\text{cc}i}-3T_{\text{aa}i}T_{\text{bc}j}T_{\text{bc}%
i}+2T_{\text{ab}i}T_{\text{bc}j}T_{\text{ac}i})](y_{0})$\qquad\qquad
\qquad\qquad\qquad\ \ 

$\qquad+\frac{1}{12}<H,j>[\nabla_{j}\varrho_{ii}-2\varrho_{ij}%
<H,i>+\overset{q}{\underset{\text{a}=1}{\sum}}(\nabla_{j}R_{\text{a}%
i\text{a}i}-4R_{i\text{a}j\text{a}}<H,i>)$

$\qquad+4\overset{q}{\underset{\text{a,b=1}}{\sum}}R_{j\text{a}i\text{b}%
}T_{\text{ab}i}$

$\qquad+2\overset{q}{\underset{\text{a,b,c=1}}{\sum}}(T_{\text{aa}%
j}T_{\text{bb}i}T_{\text{cc}i}-3T_{\text{aa}j}T_{\text{bc}i}T_{\text{bc}%
i}+2T_{\text{ab}j}T_{\text{bc}i}T_{\text{ac}i})](y_{0})$

$\qquad+\frac{1}{12}<H,j>[\nabla_{i}\varrho_{ij}-2\varrho_{ii}%
<H,j>+\overset{q}{\underset{\text{a}=1}{\sum}}(\nabla_{i}R_{\text{a}%
i\text{a}j}-4R_{i\text{a}i\text{a}}<H,j>)$

$\qquad+4\overset{q}{\underset{\text{a,b=1}}{\sum}}R_{i\text{a}i\text{b}%
}T_{\text{ab}j}$

$\qquad+2\overset{q}{\underset{\text{a,b,c}=1}{\sum}}(T_{\text{aa}%
i}T_{\text{bb}i}T_{\text{cc}j}-3T_{\text{aa}i}T_{\text{bc}i}T_{\text{bc}%
j}+2T_{\text{ab}i}T_{\text{bc}i}T_{\text{ac}j})](y_{0}$

\qquad\qquad\qquad\qquad\qquad\qquad\qquad\qquad\qquad\qquad\qquad\qquad
\qquad\qquad\qquad\qquad\qquad\qquad\qquad$\blacksquare$

We next compute $L_{22}=\frac{\partial^{2}}{\partial x_{i}^{2}}[\theta
\frac{\partial g^{jk}}{\partial\text{x}_{j}}\frac{\partial\theta^{-\frac{1}%
{2}}}{\partial\text{x}_{k}}](y_{0}):$ Since $\theta(y_{0})=1$

$L_{22}=\frac{\partial^{2}}{\partial x_{i}^{2}}[\theta\frac{\partial
\text{g}^{jk}}{\partial\text{x}_{j}}\frac{\partial\theta^{-\frac{1}{2}}%
}{\partial\text{x}_{k}}](y_{0})=2\frac{\partial\theta}{\partial x_{i}}%
(y_{0})\frac{\partial}{\partial x_{i}}[\frac{\partial\text{g}^{jk}}%
{\partial\text{x}_{j}}\frac{\partial\theta^{-\frac{1}{2}}}{\partial
\text{x}_{k}}](y_{0})$

$\qquad+\frac{\partial^{2}\theta}{\partial x_{i}^{2}}(y_{0})[\frac
{\partial\text{g}^{jk}}{\partial\text{x}_{j}}\frac{\partial\theta^{-\frac
{1}{2}}}{\partial\text{x}_{k}}](y_{0})+\frac{\partial^{2}}{\partial x_{i}^{2}%
}[\frac{\partial\text{g}^{jk}}{\partial\text{x}_{j}}\frac{\partial
\theta^{-\frac{1}{2}}}{\partial\text{x}_{k}}](y_{0})$

$\qquad=L_{221}+L_{222}+L_{223}$

where,

$L_{221}=2\frac{\partial\theta}{\partial x_{i}}(y_{0})\frac{\partial}{\partial
x_{i}}[\frac{\partial\text{g}^{jk}}{\partial\text{x}_{j}}\frac{\partial
\theta^{-\frac{1}{2}}}{\partial\text{x}_{k}}](y_{0});\qquad L_{222}%
=\frac{\partial^{2}\theta}{\partial x_{i}^{2}}(y_{0})[\frac{\partial
\text{g}^{jk}}{\partial\text{x}_{j}}\frac{\partial\theta^{-\frac{1}{2}}%
}{\partial\text{x}_{k}}](y_{0});$

$L_{223}=\frac{\partial^{2}}{\partial x_{i}^{2}}[\frac{\partial\text{g}^{jk}%
}{\partial\text{x}_{j}}\frac{\partial\theta^{-\frac{1}{2}}}{\partial
\text{x}_{k}}](y_{0})=2\frac{\partial^{2}\text{g}^{jk}}{\partial x_{i}%
\partial\text{x}_{j}}(y_{0})\frac{\partial^{2}\theta^{-\frac{1}{2}}}{\partial
x_{i}\partial\text{x}_{k}}(y_{0})+\frac{\partial\theta^{-\frac{1}{2}}%
}{\partial\text{x}_{k}}(y_{0})\frac{\partial^{3}\text{g}^{jk}}{\partial
x_{i}^{2}\partial\text{x}_{j}}(y_{0})$

$\qquad+\frac{\partial\text{g}^{jk}}{\partial\text{x}_{j}}(y_{0}%
)\frac{\partial^{3}\theta^{-\frac{1}{2}}}{\partial x_{i}^{2}\partial
\text{x}_{k}}(y_{0})$

We compute:

\qquad$L_{221}=2\frac{\partial\theta}{\partial x_{i}}(y_{0})\frac{\partial
}{\partial x_{i}}[\frac{\partial\text{g}^{jk}}{\partial\text{x}_{j}}%
\frac{\partial\theta^{-\frac{1}{2}}}{\partial\text{x}_{k}}](y_{0}%
)=2\frac{\partial\theta}{\partial x_{i}}(y_{0})[\frac{\partial^{2}%
\text{g}^{jk}}{\partial x_{i}\partial\text{x}_{j}}\frac{\partial\theta
^{-\frac{1}{2}}}{\partial\text{x}_{k}}+\frac{\partial\text{g}^{jk}}%
{\partial\text{x}_{j}}\frac{\partial^{2}\theta^{-\frac{1}{2}}}{\partial
x_{i}\partial\text{x}_{k}}](y_{0})$

Here again, we must have $i,j,k=q+1,...,n:$ Therefore, $\frac{\partial
\text{g}^{jk}}{\partial\text{x}_{j}}(y_{0})=0.$ Consequently

$L_{221}=2\frac{\partial\theta}{\partial x_{i}}(y_{0})[\frac{\partial
^{2}\text{g}^{jk}}{\partial x_{i}\partial\text{x}_{j}}\frac{\partial
\theta^{-\frac{1}{2}}}{\partial\text{x}_{k}}](y_{0})$

We have by (iv) of \textbf{Appendix A}$_{2},$ $\frac{\partial^{2}\text{g}%
^{kl}}{\partial\text{x}_{i}\partial\text{x}_{j}}(y_{0})=\frac{1}{3}%
(R_{ikjl}+R_{jkil})(y_{0}).$

In particular we have the following:

$\frac{\partial^{2}\text{g}^{jk}}{\partial\text{x}_{i}\partial\text{x}_{j}%
}(y_{0})=\frac{1}{3}(R_{ijjk}+R_{jjik})(y_{0})=\frac{1}{3}R_{ijjk}%
(y_{0})=-\frac{1}{3}R_{jijk}(y_{0})$

$L_{221}=2\frac{\partial\theta}{\partial x_{i}}(y_{0})[\frac{\partial
^{2}\text{g}^{jk}}{\partial x_{i}\partial\text{x}_{j}}\frac{\partial
\theta^{-\frac{1}{2}}}{\partial\text{x}_{k}}](y_{0})=-2.\frac{1}%
{2}<H,i><H,k>[-\frac{1}{3}R_{jijk}](y_{0})$

$L_{221}=\frac{1}{3}<H,i><H,k>R_{jijk}(y_{0})$

\qquad\qquad\qquad\qquad\qquad\qquad\qquad\qquad\qquad\qquad\qquad
\qquad$\blacksquare$

Next, since $\frac{\partial\text{g}^{jk}}{\partial\text{x}_{j}}(y_{0})=0$ for
$i,j,k=q+1,...,n,$ we have:

\bigskip$L_{222}=\frac{\partial^{2}\theta}{\partial x_{i}^{2}}(y_{0}%
)[\frac{\partial\text{g}^{jk}}{\partial\text{x}_{j}}\frac{\partial
\theta^{-\frac{1}{2}}}{\partial\text{x}_{k}}](y_{0})=0$

\qquad\qquad\qquad\qquad\qquad\qquad\qquad\qquad\qquad\qquad\qquad
\qquad$\blacksquare$

Next, we have

$L_{223}=2\frac{\partial^{2}\text{g}^{jk}}{\partial x_{i}\partial\text{x}_{j}%
}(y_{0})\frac{\partial^{2}\theta^{-\frac{1}{2}}}{\partial x_{i}\partial
\text{x}_{k}}(y_{0})+\frac{\partial\theta^{-\frac{1}{2}}}{\partial\text{x}%
_{k}}(y_{0})\frac{\partial^{3}\text{g}^{jk}}{\partial x_{i}^{2}\partial
\text{x}_{j}}(y_{0})+\frac{\partial\text{g}^{jk}}{\partial\text{x}_{j}}%
(y_{0})\frac{\partial^{3}\theta^{-\frac{1}{2}}}{\partial x_{i}^{2}%
\partial\text{x}_{k}}(y_{0})$

Again, $\frac{\partial\text{g}^{jk}}{\partial\text{x}_{j}}(y_{0})=0$ for
$i,j,k=q+1,...,n$ and so,

$L_{223}=2\frac{\partial^{2}\text{g}^{jk}}{\partial x_{i}\partial\text{x}_{j}%
}(y_{0})\frac{\partial^{2}\theta^{-\frac{1}{2}}}{\partial x_{i}\partial
\text{x}_{k}}(y_{0})+\frac{\partial\theta^{-\frac{1}{2}}}{\partial\text{x}%
_{k}}(y_{0})\frac{\partial^{3}\text{g}^{jk}}{\partial x_{i}^{2}\partial
\text{x}_{j}}(y_{0})$\qquad

$\qquad\frac{\partial^{2}\text{g}^{kl}}{\partial\text{x}_{i}\partial
\text{x}_{j}}(y_{0})=\frac{1}{3}(R_{ikjl}+R_{jkil})(y_{0})$ by (iv) of
\textbf{Appendix A}$_{2}.$

In particular we have the following:

\ $\frac{\partial^{2}\text{g}^{jk}}{\partial\text{x}_{i}\partial\text{x}_{j}%
}(y_{0})=\frac{1}{3}(R_{ijjk}+R_{jjik})(y_{0})=\frac{1}{3}R_{ijjk}%
(y_{0})=-\frac{1}{3}R_{jijk}(y_{0})$

By (v) of \textbf{Appendix A}$_{2},$

$\frac{\partial^{3}\text{g}^{kl}}{\partial\text{x}_{i}^{2}\partial\text{x}%
_{j}}(y_{0})=$ $\frac{1}{3}\nabla_{j}$R$_{ikil}(y_{0})+\frac{1}{3}\nabla_{i}%
$R$_{jkil}(y_{0})+\frac{1}{3}\nabla_{i}$R$_{ikjl}(y_{0}).$

Therefore,

$\frac{\partial^{3}\text{g}^{kj}}{\partial\text{x}_{i}^{2}\partial\text{x}%
_{j}}(y_{0})=\frac{1}{3}\nabla_{j}$R$_{ikij}(y_{0})+\frac{1}{3}\nabla_{i}%
$R$_{jkij}(y_{0})+\frac{1}{3}\nabla_{i}$R$_{ikjj}(y_{0})$

$=\frac{1}{3}\nabla_{j}$R$_{ikij}(y_{0})+\frac{1}{3}\nabla_{i}$R$_{jkij}%
(y_{0})$

As before, $\frac{\partial^{2}\theta^{-\frac{1}{2}}}{\partial\text{x}%
_{i}\partial\text{x}_{k}}(y_{0})$ is taken from (ix) of \textbf{Appendix
A}$_{9}:$

$\frac{\partial^{2}\theta^{-\frac{1}{2}}}{\partial\text{x}_{i}\partial
\text{x}_{j}}(y_{0})$ \ $=\frac{3}{4}<H,i>(y_{0})<H,j>(y_{0})$

$+\frac{1}{12}[2\varrho_{ij}+$ $\overset{q}{\underset{\text{a}=1}{4\sum}%
}R_{i\text{a}j\text{a}}-3\overset{q}{\underset{\text{a,b=1}}{\sum}%
}(T_{\text{aa}i}T_{\text{bb}j}-T_{\text{ab}i}T_{\text{ab}j}%
)-3\overset{q}{\underset{\text{a,b=1}}{\sum}}(T_{\text{aa}j}T_{\text{bb}%
i}-T_{\text{ab}j}T_{\text{ab}i}](y_{0})$

We have:

$L_{223}=2(-\frac{1}{3}R_{ijjk})(y_{0})$ \ $[\frac{3}{4}<H,i><H,k>](y_{0})$

$+2(-\frac{1}{3}R_{ijjk})(y_{0}).\frac{1}{12}[2\varrho_{ik}+$
$\overset{q}{\underset{\text{a}=1}{4\sum}}R_{i\text{a}k\text{a}}%
-3\overset{q}{\underset{\text{a,b=1}}{\sum}}(T_{\text{aa}i}T_{\text{bb}%
k}-T_{\text{ab}i}T_{\text{ab}k})$

$-3\overset{q}{\underset{\text{a,b=1}}{\sum}}(T_{\text{aa}k}T_{\text{bb}%
i}-T_{\text{ab}k}T_{\text{ab}i}](y_{0}$

$+\frac{1}{2}\times\frac{1}{3}<H,k>(y_{0})[\nabla_{j}$R$_{ikij}(y_{0}%
)+\nabla_{i}$R$_{jkij}](y_{0})$

$L_{223}=\frac{1}{2}R_{ijjk}(y_{0})$ \ $[<H,i><H,k>](y_{0})$

$-\frac{1}{18}R_{jijk}(y_{0})[2\varrho_{ik}+$ $\overset{q}{\underset{\text{a}%
=1}{4\sum}}R_{i\text{a}k\text{a}}-3\overset{q}{\underset{\text{a,b=1}}{\sum}%
}(T_{\text{aa}i}T_{\text{bb}k}-T_{\text{ab}i}T_{\text{ab}k})$

$-3\overset{q}{\underset{\text{a,b=1}}{\sum}}(T_{\text{aa}k}T_{\text{bb}%
i}-T_{\text{ab}k}T_{\text{ab}i}](y_{0}$

$+\frac{1}{6}<H,k>(y_{0})[\nabla_{j}$R$_{ikij}(y_{0})+\nabla_{i}$%
R$_{jkij}](y_{0})$

\qquad\qquad\qquad\qquad\qquad\qquad\qquad\qquad\qquad\qquad\qquad\qquad
\qquad\qquad\qquad\qquad$\blacksquare$

Noting that $L_{222}=0,$ we conclude that,

$L_{22}=L_{221}+L_{222}+L_{223}$

$=\frac{1}{3}<H,i><H,k>R_{jijk}(y_{0})\qquad\qquad\qquad L_{221}$

$+\frac{2}{3}\times\frac{3}{4}R_{ijjk}(y_{0})$ \ $[<H,i><H,k>](y_{0})$

$=-\frac{1}{2}R_{jijk}(y_{0})\ [<H,i><H,k>](y_{0})\qquad L_{223}$

$-\frac{2}{3}\times\frac{1}{12}R_{jijk}(y_{0})[2\varrho_{ik}+$
$\overset{q}{\underset{\text{a}=1}{4\sum}}R_{i\text{a}k\text{a}}%
-3\overset{q}{\underset{\text{a,b=1}}{\sum}}(T_{\text{aa}i}T_{\text{bb}%
k}-T_{\text{ab}i}T_{\text{ab}k})$

$-3\overset{q}{\underset{\text{a,b=1}}{\sum}}(T_{\text{aa}k}T_{\text{bb}%
i}-T_{\text{ab}k}T_{\text{ab}i}](y_{0}$

$+\frac{1}{2}\times\frac{1}{3}<H,k>(y_{0})[\nabla_{j}$R$_{ikij}(y_{0}%
)+\nabla_{i}$R$_{jkij}](y_{0})$

$L_{22}=-\frac{1}{6}R_{jijk}(y_{0})$ \ $[<H,i><H,k>](y_{0})\qquad\qquad
L_{22}$

$-\frac{1}{18}R_{jijk}(y_{0})[2\varrho_{ik}+$ $\overset{q}{\underset{\text{a}%
=1}{4\sum}}R_{i\text{a}k\text{a}}-3\overset{q}{\underset{\text{a,b=1}}{\sum}%
}(T_{\text{aa}i}T_{\text{bb}k}-T_{\text{ab}i}T_{\text{ab}k})$

$-3\overset{q}{\underset{\text{a,b=1}}{\sum}}(T_{\text{aa}k}T_{\text{bb}%
i}-T_{\text{ab}k}T_{\text{ab}i}](y_{0})$

$+\frac{1}{6}<H,k>(y_{0})[\nabla_{j}$R$_{ijik}(y_{0})-\nabla_{i}$%
R$_{jijk}](y_{0})$

\qquad\qquad\qquad\qquad\qquad\qquad\qquad\qquad\qquad\qquad\qquad\qquad
\qquad\qquad$\blacksquare$

We then compute:

$L_{23}=\frac{\partial^{2}}{\partial x_{i}^{2}}[\theta g^{jk}\frac
{\partial^{2}\theta^{-\frac{1}{2}}}{\partial\text{x}_{j}\partial\text{x}_{k}%
}](y_{0})$

$=2\frac{\partial\theta}{\partial\text{x}_{i}}(y_{0})\frac{\partial}%
{\partial\text{x}_{i}}[g^{jk}\frac{\partial^{2}\theta^{-\frac{1}{2}}}%
{\partial\text{x}_{j}\partial\text{x}_{k}}](y_{0})+\frac{\partial^{2}\theta
}{\partial x_{i}^{2}}(y_{0})[g^{jk}\frac{\partial^{2}\theta^{-\frac{1}{2}}%
}{\partial\text{x}_{j}\partial\text{x}_{k}}](y_{0})$

$+\frac{\partial^{2}}{\partial x_{i}^{2}}[g^{jk}\frac{\partial^{2}%
\theta^{-\frac{1}{2}}}{\partial\text{x}_{j}\partial\text{x}_{k}}%
](y_{0})=L_{231}+L_{232}+L_{233}$

where,

$L_{231}=2\frac{\partial\theta}{\partial\text{x}_{i}}(y_{0})\frac{\partial
}{\partial\text{x}_{i}}[g^{jk}\frac{\partial^{2}\theta^{-\frac{1}{2}}%
}{\partial\text{x}_{j}\partial\text{x}_{k}}](y_{0});\qquad L_{232}%
=\frac{\partial^{2}\theta}{\partial x_{i}^{2}}(y_{0})[g^{jk}\frac{\partial
^{2}\theta^{-\frac{1}{2}}}{\partial\text{x}_{j}\partial\text{x}_{k}}%
](y_{0});\qquad$

$L_{233}=\frac{\partial^{2}}{\partial x_{i}^{2}}[g^{jk}\frac{\partial
^{2}\theta^{-\frac{1}{2}}}{\partial\text{x}_{j}\partial\text{x}_{k}}](y_{0})$

$L_{231}=2\frac{\partial\theta}{\partial\text{x}_{i}}(y_{0})\frac{\partial
}{\partial\text{x}_{i}}[g^{jk}\frac{\partial^{2}\theta^{-\frac{1}{2}}%
}{\partial\text{x}_{j}\partial\text{x}_{k}}](y_{0})=2\frac{\partial\theta
}{\partial\text{x}_{i}}(y_{0})\frac{\partial g^{jk}}{\partial\text{x}_{i}%
}(y_{0})\frac{\partial^{2}\theta^{-\frac{1}{2}}}{\partial\text{x}_{j}%
\partial\text{x}_{k}}(y_{0})$

$+2\frac{\partial\theta}{\partial\text{x}_{i}}(y_{0})g^{jk}(y_{0}%
)\frac{\partial^{3}\theta^{-\frac{1}{2}}}{\partial\text{x}_{i}\partial
\text{x}_{j}\partial\text{x}_{k}}(y_{0})$

Since $g^{jk}(y_{0})=\delta^{ij}$ and $\frac{\partial g^{jk}}{\partial
\text{x}_{i}}(y_{0})=0$ for $i,j,k=q+1,...,n,$ we have:

$L_{231}=2\frac{\partial\theta}{\partial\text{x}_{i}}(y_{0})\frac{\partial
^{3}\theta^{-\frac{1}{2}}}{\partial\text{x}_{i}\partial\text{x}_{j}^{2}}%
(y_{0})$

We use $\left(  A_{17}\right)  $ where,

\ $\frac{\partial^{3}\theta^{-\frac{1}{2}}}{\partial\text{x}_{i}%
\partial\text{x}_{j}^{2}}(y_{0})=\frac{15}{8}[<H,i><H,j>^{2}](y_{0})$\qquad

$+\frac{1}{4}<H,j>[2\varrho_{ij}+$ $\overset{q}{\underset{\text{a}=1}{4\sum}%
}R_{i\text{a}j\text{a}}-3\overset{q}{\underset{\text{a,b=1}}{\sum}%
}(T_{\text{aa}i}T_{\text{bb}j}-T_{\text{ab}i}T_{\text{ab}j})$

$-3\overset{q}{\underset{\text{a,b=1}}{\sum}}(T_{\text{aa}j}T_{\text{bb}%
i}-T_{\text{ab}j}T_{\text{ab}i}](y_{0})$

$+\frac{1}{4}<H,i>(y_{0})[\tau^{M}\ -3\tau^{P}+\ \underset{\text{a}%
=1}{\overset{\text{q}}{\sum}}\varrho_{\text{aa}}^{M}+$
$\overset{q}{\underset{\text{a},\text{b}=1}{\sum}}R_{\text{abab}}^{M}$
$](y_{0})$

$+\frac{1}{12}[\nabla_{i}\varrho_{jj}-2\varrho_{ij}%
<H,j>+\overset{q}{\underset{\text{a}=1}{\sum}}(\nabla_{i}R_{\text{a}%
j\text{a}j}-4R_{i\text{a}j\text{a}}<H,j>)\qquad\frac{\partial^{3}\theta
}{\partial\text{x}_{i}\partial\text{x}_{j}^{2}}(y_{0})$

$+4\overset{q}{\underset{\text{a,b=1}}{\sum}}R_{i\text{a}j\text{b}%
}T_{\text{ab}j}+2\overset{q}{\underset{\text{a,b,c=1}}{\sum}}(T_{\text{aa}%
i}T_{\text{bb}j}T_{\text{cc}j}-3T_{\text{aa}i}T_{\text{bc}j}T_{\text{bc}%
j}+2T_{\text{ab}i}T_{\text{bc}j}T_{\text{ca}j})](y_{0})$\qquad\qquad
\qquad\qquad\qquad\ \ 

$+\frac{1}{12}[\nabla_{j}\varrho_{ij}-2\varrho_{ij}%
<H,j>+\overset{q}{\underset{\text{a}=1}{\sum}}(\nabla_{j}R_{\text{a}%
i\text{a}j}-4R_{j\text{a}i\text{a}}<H,j>)$

$+4\overset{q}{\underset{\text{a,b=1}}{\sum}}R_{j\text{a}i\text{b}%
}T_{\text{ab}j}+2\overset{q}{\underset{\text{a,b,c=1}}{\sum}}(T_{\text{aa}%
j}T_{\text{bb}i}T_{\text{cc}j}-3T_{\text{aa}j}T_{\text{bc}i}T_{\text{bc}%
j}+2T_{\text{ab}j}T_{\text{bc}i}T_{\text{ac}j})](y_{0})$

$+\frac{1}{12}[\nabla_{j}\varrho_{ij}-2\varrho_{jj}%
<H,i>+\overset{q}{\underset{\text{a}=1}{\sum}}(\nabla_{j}R_{\text{a}%
i\text{a}j}-4R_{j\text{a}j\text{a}}<H,i>)+4\overset{q}{\underset{\text{a,b=1}%
}{\sum}}R_{j\text{a}j\text{b}}T_{\text{ab}i}$

$+2\overset{q}{\underset{\text{a,b,c=1}}{\sum}}(T_{\text{aa}j}T_{\text{bb}%
j}T_{\text{cc}i}-3T_{\text{aa}j}T_{\text{bc}j}T_{\text{bc}i}+2T_{\text{ab}%
j}T_{\text{bc}j}T_{\text{ac}i})](y_{0})$

We thus have:

$L_{231}=2\frac{\partial\theta}{\partial\text{x}_{i}}(y_{0})\frac{\partial
^{3}\theta^{-\frac{1}{2}}}{\partial\text{x}_{i}\partial\text{x}_{j}^{2}}%
(y_{0})$

$=-\frac{15}{4}<H,i>^{2}(y_{0})<H,j>^{2}(y_{0})$

$-\frac{1}{2}<H,i>(y_{0})<H,j>(y_{0})$

$\times\lbrack2\varrho_{ij}+$ $\overset{q}{\underset{\text{a}=1}{4\sum}%
}R_{i\text{a}j\text{a}}-3\overset{q}{\underset{\text{a,b=1}}{\sum}%
}(T_{\text{aa}i}T_{\text{bb}j}-T_{\text{ab}i}T_{\text{ab}j}%
)-3\overset{q}{\underset{\text{a,b=1}}{\sum}}(T_{\text{aa}j}T_{\text{bb}%
i}-T_{\text{ab}j}T_{\text{ab}i}](y_{0})$

$-\frac{1}{2}<H,i>^{2}(y_{0})[\varrho_{jj}+$ $\overset{q}{\underset{\text{a}%
=1}{2\sum}}R_{j\text{a}j\text{a}}-3\overset{q}{\underset{\text{a,b=1}}{\sum}%
}(T_{\text{aa}j}T_{\text{bb}j}-T_{\text{ab}j}T_{\text{ab}j})](y_{0})$

$-\frac{1}{6}<H,i>(y_{0})[\nabla_{i}\varrho_{jj}-2\varrho_{ij}%
<H,j>+\overset{q}{\underset{\text{a}=1}{\sum}}(\nabla_{i}R_{\text{a}%
j\text{a}j}-4R_{i\text{a}j\text{a}}<H,j>)$

$+4\overset{q}{\underset{\text{a,b=1}}{\sum}}R_{i\text{a}j\text{b}%
}T_{\text{ab}j}+2\overset{q}{\underset{\text{a,b,c=1}}{\sum}}(T_{\text{aa}%
i}T_{\text{bb}j}T_{\text{cc}j}-3T_{\text{aa}i}T_{\text{bc}j}T_{\text{bc}%
j}+2T_{\text{ab}i}T_{\text{bc}j}T_{\text{ca}j})](y_{0})$\qquad\qquad
\qquad\qquad\qquad\ \ 

$-\frac{1}{6}<H,i>(y_{0})[\nabla_{j}\varrho_{ij}-2\varrho_{ij}%
<H,j>+\overset{q}{\underset{\text{a}=1}{\sum}}(\nabla_{j}R_{\text{a}%
i\text{a}j}-4R_{j\text{a}i\text{a}}<H,j>)$

$+4\overset{q}{\underset{\text{a,b=1}}{\sum}}R_{j\text{a}i\text{b}%
}T_{\text{ab}j}+2\overset{q}{\underset{\text{a,b,c=1}}{\sum}}(T_{\text{aa}%
j}T_{\text{bb}i}T_{\text{cc}j}-3T_{\text{aa}j}T_{\text{bc}i}T_{\text{bc}%
j}+2T_{\text{ab}j}T_{\text{bc}i}T_{\text{ac}j})](y_{0})$

$-\frac{1}{6}<H,i>(y_{0})[\nabla_{j}\varrho_{ij}-2\varrho_{jj}%
<H,i>+\overset{q}{\underset{\text{a}=1}{\sum}}(\nabla_{j}R_{\text{a}%
i\text{a}j}-4R_{j\text{a}j\text{a}}<H,i>)$

$+4\overset{q}{\underset{\text{a,b=1}}{\sum}}R_{j\text{a}j\text{b}%
}T_{\text{ab}i}+2\overset{q}{\underset{\text{a,b,c=1}}{\sum}}(T_{\text{aa}%
j}T_{\text{bb}j}T_{\text{cc}i}-3T_{\text{aa}j}T_{\text{bc}j}T_{\text{bc}%
i}+2T_{\text{ab}j}T_{\text{bc}j}T_{\text{ac}i})](y_{0})$

\qquad\qquad\qquad\qquad\qquad\qquad\qquad\qquad\qquad\qquad\qquad\qquad
\qquad\qquad\qquad\qquad\qquad\qquad\qquad$\blacksquare$

$L_{232}=\frac{\partial^{2}\theta}{\partial x_{i}^{2}}(y_{0})[g^{jk}%
\frac{\partial^{2}\theta^{-\frac{1}{2}}}{\partial\text{x}_{j}\partial
\text{x}_{k}}](y_{0})=\frac{\partial^{2}\theta}{\partial x_{i}^{2}}%
(y_{0})\frac{\partial^{2}\theta^{-\frac{1}{2}}}{\partial\text{x}_{j}^{2}%
}(y_{0})$

From (v) of \textbf{Appendix A}$_{9}:$

$\frac{\partial^{2}\theta}{\partial\text{x}_{i}^{2}}(y_{0})=-\frac{1}%
{3}[\varrho_{ii}+$ $\overset{q}{\underset{\text{a}=1}{2\sum}}R_{i\text{a}%
i\text{a}}-3\overset{q}{\underset{\text{a,b=1}}{\sum}}(T_{\text{aa}%
i}T_{\text{bb}i}-T_{\text{ab}i}T_{\text{ab}i})](y_{0})$

As before, $\frac{\partial^{2}\theta^{-\frac{1}{2}}}{\partial\text{x}%
_{i}\partial\text{x}_{k}}(y_{0})$ (ix) of \textbf{Appendix A}$_{9}:$

$\frac{\partial^{2}\theta^{-\frac{1}{2}}}{\partial\text{x}_{j}^{2}}(y_{0})$
\ $=\frac{3}{4}<H,j>^{2}(y_{0})+\frac{1}{6}[\varrho_{jj}+$
$\overset{q}{\underset{\text{a}=1}{2\sum}}R_{j\text{a}j\text{a}}%
-3\overset{q}{\underset{\text{a,b=1}}{\sum}}(T_{\text{aa}j}T_{\text{bb}%
j}-T_{\text{ab}j}T_{\text{ab}j})](y_{0})$

Therefore,

$L_{232}=-\frac{1}{4}<H,j>^{2}(y_{0})[\varrho_{ii}+$
$\overset{q}{\underset{\text{a}=1}{2\sum}}R_{i\text{a}i\text{a}}%
-3\overset{q}{\underset{\text{a,b=1}}{\sum}}(T_{\text{aa}i}T_{\text{bb}%
i}-T_{\text{ab}i}T_{\text{ab}i})](y_{0})$

$\qquad-\frac{1}{18}[\varrho_{ii}+$ $\overset{q}{\underset{\text{a}=1}{2\sum}%
}R_{i\text{a}i\text{a}}-3\overset{q}{\underset{\text{a,b=1}}{\sum}%
}(T_{\text{aa}i}T_{\text{bb}i}-T_{\text{ab}i}T_{\text{ab}i})](y_{0})$

\qquad\qquad$\times\lbrack\varrho_{jj}+$ $\overset{q}{\underset{\text{a}%
=1}{2\sum}}R_{j\text{a}j\text{a}}-3\overset{q}{\underset{\text{a,b=1}}{\sum}%
}(T_{\text{aa}j}T_{\text{bb}j}-T_{\text{ab}j}T_{\text{ab}j})](y_{0})$

\qquad\qquad\qquad\qquad\qquad\qquad\qquad\qquad\qquad\qquad\qquad\qquad
\qquad\qquad$\blacksquare$

Next we have,

$L_{233}=\frac{\partial^{2}}{\partial x_{i}^{2}}[$g$^{jk}\frac{\partial
^{2}\theta^{-\frac{1}{2}}}{\partial\text{x}_{j}\partial\text{x}_{k}}%
](y_{0})=2\frac{\partial\text{g}^{jk}}{\partial\text{x}_{i}}(y_{0}%
)\frac{\partial^{3}\theta^{-\frac{1}{2}}}{\partial x_{i}\partial\text{x}%
_{j}\partial\text{x}_{k}}(y_{0})$

$+\frac{\partial^{2}\text{g}^{jk}}{\partial x_{i}^{2}}(y_{0})\frac
{\partial^{2}\theta^{-\frac{1}{2}}}{\partial\text{x}_{j}\partial\text{x}_{k}%
}(y_{0})+$ g$^{jk}(y_{0})\frac{\partial^{4}\theta^{-\frac{1}{2}}}{\partial
x_{i}^{2}\partial\text{x}_{j}\partial\text{x}_{k}}(y_{0})$

Since g$^{jk}(y_{0})=\delta^{jk},$ we have:

$L_{233}=2\frac{\partial\text{g}^{jk}}{\partial\text{x}_{i}}(y_{0}%
)\frac{\partial^{3}\theta^{-\frac{1}{2}}}{\partial x_{i}\partial\text{x}%
_{j}\partial\text{x}_{k}}(y_{0})+\frac{\partial^{2}\text{g}^{jk}}{\partial
x_{i}^{2}}(y_{0})\frac{\partial^{2}\theta^{-\frac{1}{2}}}{\partial\text{x}%
_{j}\partial\text{x}_{k}}(y_{0})+\frac{\partial^{4}\theta^{-\frac{1}{2}}%
}{\partial x_{i}^{2}\partial x_{j}^{2}}(y_{0})$

Again $\frac{\partial\text{g}^{jk}}{\partial\text{x}_{i}}=0$ for
$i,j,k=q+1,...,n$ and so,

$L_{233}=\frac{\partial^{2}\text{g}^{jk}}{\partial x_{i}^{2}}(y_{0}%
)\frac{\partial^{2}\theta^{-\frac{1}{2}}}{\partial\text{x}_{j}\partial
\text{x}_{k}}(y_{0})+$ $\frac{\partial^{4}\theta^{-\frac{1}{2}}}{\partial
x_{i}^{2}\partial x_{j}^{2}}(y_{0})$

$\frac{\partial^{2}\text{g}^{jk}}{\partial x_{i}^{2}}(y_{0})=\frac{2}%
{3}R_{ijik}(y_{0})$ by (iii) of \textbf{Appendix A}$_{2}$ and $\frac
{\partial^{2}\theta^{-\frac{1}{2}}}{\partial\text{x}_{j}\partial\text{x}_{k}%
}(y_{0})$

is given by (ix) of \textbf{Appendix A}$_{9}.$

Finally $\frac{\partial^{4}\theta^{-\frac{1}{2}}}{\partial x_{i}^{2}\partial
x_{j}^{2}}(y_{0})$ is given by (xx) of \textbf{Appendix A}$_{9}.$ We have:

$L_{233}=\frac{\partial^{2}\text{g}^{jk}}{\partial x_{i}^{2}}(y_{0}%
)\frac{\partial^{2}\theta^{-\frac{1}{2}}}{\partial\text{x}_{j}\partial
\text{x}_{k}}(y_{0})+$ $\frac{\partial^{4}\theta^{-\frac{1}{2}}}{\partial
x_{i}^{2}\partial x_{j}^{2}}(y_{0})$

We have:

\ $L_{233}=\frac{1}{2}R_{ijik}(y_{0})$ \ $[<H,j><H,k>](y_{0})$

$+\frac{1}{18}R_{ijik}(y_{0})[2\varrho_{jk}+$ $\overset{q}{\underset{\text{a}%
=1}{4\sum}}R_{j\text{a}k\text{a}}-3\overset{q}{\underset{\text{a,b=1}}{\sum}%
}(T_{\text{aa}j}T_{\text{bb}k}-T_{\text{ab}j}T_{\text{ab}k})$

$-3\overset{q}{\underset{\text{a,b=1}}{\sum}}(T_{\text{aa}k}T_{\text{bb}%
j}-T_{\text{ab}k}T_{\text{ab}j}](y_{0})+$ $\frac{\partial^{4}\theta^{-\frac
{1}{2}}}{\partial x_{i}^{2}\partial x_{j}^{2}}(y_{0})$

\qquad\qquad\qquad\qquad\qquad\qquad\qquad\qquad\qquad\qquad\qquad\qquad
\qquad\qquad\qquad$\blacksquare$

Therefore,

$L_{23}=L_{231}+L_{232}+\ L_{233}$

$=-\frac{15}{4}<H,i>^{2}(y_{0})<H,j>^{2}(y_{0})\qquad\qquad L_{231}$

$-\frac{1}{2}<H,i>(y_{0})<H,j>(y_{0})$

$\times\lbrack2\varrho_{ij}+$ $\overset{q}{\underset{\text{a}=1}{4\sum}%
}R_{i\text{a}j\text{a}}-3\overset{q}{\underset{\text{a,b=1}}{\sum}%
}(T_{\text{aa}i}T_{\text{bb}j}-T_{\text{ab}i}T_{\text{ab}j}%
)-3\overset{q}{\underset{\text{a,b=1}}{\sum}}(T_{\text{aa}j}T_{\text{bb}%
i}-T_{\text{ab}j}T_{\text{ab}i}](y_{0})$

$-\frac{1}{2}<H,i>^{2}(y_{0})[\varrho_{jj}+$ $\overset{q}{\underset{\text{a}%
=1}{2\sum}}R_{j\text{a}j\text{a}}-3\overset{q}{\underset{\text{a,b=1}}{\sum}%
}(T_{\text{aa}j}T_{\text{bb}j}-T_{\text{ab}j}T_{\text{ab}j})](y_{0})$

$-\frac{1}{6}<H,i>(y_{0})[\nabla_{i}\varrho_{jj}-2\varrho_{ij}%
<H,j>+\overset{q}{\underset{\text{a}=1}{\sum}}(\nabla_{i}R_{\text{a}%
j\text{a}j}-4R_{i\text{a}j\text{a}}<H,j>)$

$+4\overset{q}{\underset{\text{a,b=1}}{\sum}}R_{i\text{a}j\text{b}%
}T_{\text{ab}j}$

$+2\overset{q}{\underset{\text{a,b,c=1}}{\sum}}(T_{\text{aa}i}T_{\text{bb}%
j}T_{\text{cc}j}-3T_{\text{aa}i}T_{\text{bc}j}T_{\text{bc}j}+2T_{\text{ab}%
i}T_{\text{bc}j}T_{\text{ca}j})](y_{0})$\qquad\qquad\qquad\qquad\qquad\ \ 

$-\frac{1}{6}<H,i>(y_{0})[\nabla_{j}\varrho_{ij}-2\varrho_{ij}%
<H,j>+\overset{q}{\underset{\text{a}=1}{\sum}}(\nabla_{j}R_{\text{a}%
i\text{a}j}-4R_{j\text{a}i\text{a}}<H,j>)$

$+4\overset{q}{\underset{\text{a,b=1}}{\sum}}R_{j\text{a}i\text{b}%
}T_{\text{ab}j}+2\overset{q}{\underset{\text{a,b,c=1}}{\sum}}(T_{\text{aa}%
j}T_{\text{bb}i}T_{\text{cc}j}-3T_{\text{aa}j}T_{\text{bc}i}T_{\text{bc}%
j}+2T_{\text{ab}j}T_{\text{bc}i}T_{\text{ac}j})](y_{0})$

$-\frac{1}{6}<H,i>(y_{0})[\nabla_{j}\varrho_{ij}-2\varrho_{jj}%
<H,i>+\overset{q}{\underset{\text{a}=1}{\sum}}(\nabla_{j}R_{\text{a}%
i\text{a}j}-4R_{j\text{a}j\text{a}}<H,i>)$

$+4\overset{q}{\underset{\text{a,b=1}}{\sum}}R_{j\text{a}j\text{b}%
}T_{\text{ab}i}+2\overset{q}{\underset{\text{a,b,c=1}}{\sum}}(T_{\text{aa}%
j}T_{\text{bb}j}T_{\text{cc}i}-3T_{\text{aa}j}T_{\text{bc}j}T_{\text{bc}%
i}+2T_{\text{ab}j}T_{\text{bc}j}T_{\text{ac}i})](y_{0})$

$-\frac{1}{4}<H,j>^{2}(y_{0})[\varrho_{ii}+$ $\overset{q}{\underset{\text{a}%
=1}{2\sum}}R_{i\text{a}i\text{a}}-3\overset{q}{\underset{\text{a,b=1}}{\sum}%
}(T_{\text{aa}i}T_{\text{bb}i}-T_{\text{ab}i}T_{\text{ab}i})](y_{0})\qquad
L_{232}$

$\qquad-\frac{1}{18}[\varrho_{ii}+$ $\overset{q}{\underset{\text{a}=1}{2\sum}%
}R_{i\text{a}i\text{a}}-3\overset{q}{\underset{\text{a,b=1}}{\sum}%
}(T_{\text{aa}i}T_{\text{bb}i}-T_{\text{ab}i}T_{\text{ab}i})](y_{0})$

\qquad\qquad$\times\lbrack\varrho_{jj}+$ $\overset{q}{\underset{\text{a}%
=1}{2\sum}}R_{j\text{a}j\text{a}}-3\overset{q}{\underset{\text{a,b=1}}{\sum}%
}(T_{\text{aa}j}T_{\text{bb}j}-T_{\text{ab}j}T_{\text{ab}j})]\}(y_{0})$

$\qquad+\frac{1}{2}R_{ijik}(y_{0})$ \ $[<H,j><H,k>](y_{0})\qquad\qquad
\qquad\qquad$\ $L_{233}$

$\qquad+\frac{1}{18}R_{ijik}(y_{0})[2\varrho_{jk}+$
$\overset{q}{\underset{\text{a}=1}{4\sum}}R_{j\text{a}k\text{a}}%
-3\overset{q}{\underset{\text{a,b=1}}{\sum}}(T_{\text{aa}j}T_{\text{bb}%
k}-T_{\text{ab}j}T_{\text{ab}k})$

$\qquad-3\overset{q}{\underset{\text{a,b=1}}{\sum}}(T_{\text{aa}k}%
T_{\text{bb}j}-T_{\text{ab}k}T_{\text{ab}j}](y_{0})+$ $\frac{\partial
^{4}\theta^{-\frac{1}{2}}}{\partial x_{i}^{2}\partial x_{j}^{2}}(y_{0})$

\qquad\qquad\qquad\qquad\qquad\qquad\qquad\qquad\qquad\qquad\qquad\qquad
\qquad\qquad\qquad\qquad\qquad\qquad\qquad\qquad$\blacksquare$

We collect all terms of $\ L_{2}$ and have:

$\qquad L_{2}=L_{21}+L_{22}+L_{23}$

$=-\frac{1}{4}[2\varrho_{ij}+$ $\overset{q}{\underset{\text{a}=1}{4\sum}%
}R_{i\text{a}j\text{a}}-3\overset{q}{\underset{\text{a,b=1}}{\sum}%
}(T_{\text{aa}i}T_{\text{bb}j}-T_{\text{ab}i}T_{\text{ab}j})$

$-3\overset{q}{\underset{\text{a,b=1}}{\sum}}(T_{\text{aa}j}T_{\text{bb}%
i}-T_{\text{ab}j}T_{\text{ab}i}](y_{0})\times\lbrack<H,i><H,j>](y_{0})\qquad
L_{21}\qquad L_{211}$

$-\frac{1}{36}[2\varrho_{ij}+$ $\overset{q}{\underset{\text{a}=1}{4\sum}%
}R_{i\text{a}j\text{a}}-3\overset{q}{\underset{\text{a,b=1}}{\sum}%
}(T_{\text{aa}i}T_{\text{bb}j}-T_{\text{ab}i}T_{\text{ab}j}%
)-3\overset{q}{\underset{\text{a,b=1}}{\sum}}(T_{\text{aa}j}T_{\text{bb}%
i}-T_{\text{ab}j}T_{\text{ab}i}]^{2}(y_{0})$

$-\frac{1}{12}[<H,j>](y_{0})\times\lbrack\{\nabla_{i}\varrho_{ij}%
-2\varrho_{ij}<H,i>+\overset{q}{\underset{\text{a}=1}{\sum}}(\nabla
_{i}R_{\text{a}i\text{a}j}-4R_{i\text{a}j\text{a}}<H,i>)\qquad L_{212}$

$+4\overset{q}{\underset{\text{a,b=1}}{\sum}}R_{i\text{a}j\text{b}%
}T_{\text{ab}i}+2\overset{q}{\underset{\text{a,b,c=1}}{\sum}}(T_{\text{aa}%
i}T_{\text{bb}j}T_{\text{cc}i}-3T_{\text{aa}i}T_{\text{bc}j}T_{\text{bc}%
i}+2T_{\text{ab}i}T_{\text{bc}j}T_{\text{ac}i})](y_{0})$\qquad\qquad
\qquad\qquad\qquad\ \ 

$-\frac{1}{12}[<H,j>](y_{0})\times\lbrack\nabla_{j}\varrho_{ii}-2\varrho
_{ij}<H,i>+\overset{q}{\underset{\text{a}=1}{\sum}}(\nabla_{j}R_{\text{a}%
i\text{a}i}-4R_{i\text{a}j\text{a}}<H,i>)$

$+4\overset{q}{\underset{\text{a,b=1}}{\sum}}R_{j\text{a}i\text{b}%
}T_{\text{ab}i}+2\overset{q}{\underset{\text{a,b,c=1}}{\sum}}(T_{\text{aa}%
j}T_{\text{bb}i}T_{\text{cc}i}-3T_{\text{aa}j}T_{\text{bc}i}T_{\text{bc}%
i}+2T_{\text{ab}j}T_{\text{bc}i}T_{\text{ac}i})](y_{0})$

$-\frac{1}{12}[<H,j>](y_{0})\times\lbrack\nabla_{i}\varrho_{ij}-2\varrho
_{ii}<H,j>+\overset{q}{\underset{\text{a}=1}{\sum}}(\nabla_{i}R_{\text{a}%
i\text{a}j}-4R_{i\text{a}i\text{a}}<H,j>)$

$+4\overset{q}{\underset{\text{a,b=1}}{\sum}}R_{i\text{a}i\text{b}%
}T_{\text{ab}j}+2\overset{q}{\underset{\text{a,b,c}=1}{\sum}}(T_{\text{aa}%
i}T_{\text{bb}i}T_{\text{cc}j}-3T_{\text{aa}i}T_{\text{bc}i}T_{\text{bc}%
j}+2T_{\text{ab}i}T_{\text{bc}i}T_{\text{ac}j})](y_{0})$

$-\frac{1}{3}[<H,j><H,k>](y_{0})R_{ijik}(y_{0})-$ $\frac{15}{8}[<H,i>^{2}%
<H,j>^{2}](y_{0})\qquad L_{213}$

$-\frac{1}{4}<H,i><H,j>[2\varrho_{ij}+\overset{q}{\underset{\text{a}=1}{4\sum
}}R_{i\text{a}j\text{a}}-3\overset{q}{\underset{\text{a,b=1}}{\sum}%
}(T_{\text{aa}i}T_{\text{bb}j}-T_{\text{ab}i}T_{\text{ab}j})$

$-3\overset{q}{\underset{\text{a,b=1}}{\sum}}(T_{\text{aa}j}T_{\text{bb}%
i}-T_{\text{ab}j}T_{\text{ab}i}](y_{0})$

$-\frac{1}{4}<H,j>^{2}[\tau^{M}\ -3\tau^{P}+\ \underset{\text{a}%
=1}{\overset{\text{q}}{\sum}}\varrho_{\text{aa}}^{M}+$
$\overset{q}{\underset{\text{a},\text{b}=1}{\sum}}R_{\text{abab}}^{M}$
$](y_{0})$

$+\frac{1}{12}<H,j>[\nabla_{i}\varrho_{ij}-2\varrho_{ij}%
<H,i>+\overset{q}{\underset{\text{a}=1}{\sum}}(\nabla_{i}R_{\text{a}%
i\text{a}j}-4R_{i\text{a}j\text{a}}<H,i>)$

$+4\overset{q}{\underset{\text{a,b=1}}{\sum}}R_{i\text{a}j\text{b}%
}T_{\text{ab}i}+2\overset{q}{\underset{\text{a,b,c=1}}{\sum}}(T_{\text{aa}%
i}T_{\text{bb}j}T_{\text{cc}i}-3T_{\text{aa}i}T_{\text{bc}j}T_{\text{bc}%
i}+2T_{\text{ab}i}T_{\text{bc}j}T_{\text{ac}i})](y_{0})$\qquad\qquad
\qquad\qquad\qquad\ \ 

$+\frac{1}{12}<H,j>[\nabla_{j}\varrho_{ii}-2\varrho_{ij}%
<H,i>+\overset{q}{\underset{\text{a}=1}{\sum}}(\nabla_{j}R_{\text{a}%
i\text{a}i}-4R_{i\text{a}j\text{a}}<H,i>)$

$+4\overset{q}{\underset{\text{a,b=1}}{\sum}}R_{j\text{a}i\text{b}%
}T_{\text{ab}i}+2\overset{q}{\underset{\text{a,b,c=1}}{\sum}}(T_{\text{aa}%
j}T_{\text{bb}i}T_{\text{cc}i}-3T_{\text{aa}j}T_{\text{bc}i}T_{\text{bc}%
i}+2T_{\text{ab}j}T_{\text{bc}i}T_{\text{ac}i})](y_{0})$

$+\frac{1}{12}<H,j>[\nabla_{i}\varrho_{ij}-2\varrho_{ii}%
<H,j>+\overset{q}{\underset{\text{a}=1}{\sum}}(\nabla_{i}R_{\text{a}%
i\text{a}j}-4R_{i\text{a}i\text{a}}<H,j>)$

$+4\overset{q}{\underset{\text{a,b=1}}{\sum}}R_{i\text{a}i\text{b}%
}T_{\text{ab}j}+2\overset{q}{\underset{\text{a,b,c}=1}{\sum}}(T_{\text{aa}%
i}T_{\text{bb}i}T_{\text{cc}j}-3T_{\text{aa}i}T_{\text{bc}i}T_{\text{bc}%
j}+2T_{\text{ab}i}T_{\text{bc}i}T_{\text{ac}j})](y_{0})$

$-\frac{1}{6}R_{jijk}(y_{0})$ \ $[<H,i><H,k>](y_{0})\qquad\qquad L_{22}$

$-\frac{1}{18}R_{jijk}(y_{0})[2\varrho_{ik}+$ $\overset{q}{\underset{\text{a}%
=1}{4\sum}}R_{i\text{a}k\text{a}}-3\overset{q}{\underset{\text{a,b=1}}{\sum}%
}(T_{\text{aa}i}T_{\text{bb}k}-T_{\text{ab}i}T_{\text{ab}k})$

$-3\overset{q}{\underset{\text{a,b=1}}{\sum}}(T_{\text{aa}k}T_{\text{bb}%
i}-T_{\text{ab}k}T_{\text{ab}i}](y_{0})+\frac{1}{6}<H,k>(y_{0})[\nabla_{j}%
$R$_{ijik}(y_{0})-\nabla_{i}$R$_{jijk}](y_{0})$

$-\frac{15}{4}<H,i>^{2}(y_{0})<H,j>^{2}(y_{0})\qquad\qquad L_{23}\qquad
L_{231}$

$-\frac{1}{2}<H,i>(y_{0})<H,j>(y_{0})$

$\times\lbrack2\varrho_{ij}+$ $\overset{q}{\underset{\text{a}=1}{4\sum}%
}R_{i\text{a}j\text{a}}-3\overset{q}{\underset{\text{a,b=1}}{\sum}%
}(T_{\text{aa}i}T_{\text{bb}j}-T_{\text{ab}i}T_{\text{ab}j}%
)-3\overset{q}{\underset{\text{a,b=1}}{\sum}}(T_{\text{aa}j}T_{\text{bb}%
i}-T_{\text{ab}j}T_{\text{ab}i}](y_{0})$

$-\frac{1}{2}<H,i>^{2}(y_{0})[\varrho_{jj}+$ $\overset{q}{\underset{\text{a}%
=1}{2\sum}}R_{j\text{a}j\text{a}}-3\overset{q}{\underset{\text{a,b=1}}{\sum}%
}(T_{\text{aa}j}T_{\text{bb}j}-T_{\text{ab}j}T_{\text{ab}j})](y_{0})$

$-\frac{1}{6}<H,i>(y_{0})[\nabla_{i}\varrho_{jj}-2\varrho_{ij}%
<H,j>+\overset{q}{\underset{\text{a}=1}{\sum}}(\nabla_{i}R_{\text{a}%
j\text{a}j}-4R_{i\text{a}j\text{a}}<H,j>)$

$+4\overset{q}{\underset{\text{a,b=1}}{\sum}}R_{i\text{a}j\text{b}%
}T_{\text{ab}j}+2\overset{q}{\underset{\text{a,b,c=1}}{\sum}}(T_{\text{aa}%
i}T_{\text{bb}j}T_{\text{cc}j}-3T_{\text{aa}i}T_{\text{bc}j}T_{\text{bc}%
j}+2T_{\text{ab}i}T_{\text{bc}j}T_{\text{ca}j})](y_{0})$\qquad\qquad
\qquad\qquad\qquad\ \ 

$-\frac{1}{6}<H,i>(y_{0})[\nabla_{j}\varrho_{ij}-2\varrho_{ij}%
<H,j>+\overset{q}{\underset{\text{a}=1}{\sum}}(\nabla_{j}R_{\text{a}%
i\text{a}j}-4R_{j\text{a}i\text{a}}<H,j>)$

$+4\overset{q}{\underset{\text{a,b=1}}{\sum}}R_{j\text{a}i\text{b}%
}T_{\text{ab}j}+2\overset{q}{\underset{\text{a,b,c=1}}{\sum}}(T_{\text{aa}%
j}T_{\text{bb}i}T_{\text{cc}j}-3T_{\text{aa}j}T_{\text{bc}i}T_{\text{bc}%
j}+2T_{\text{ab}j}T_{\text{bc}i}T_{\text{ac}j})](y_{0})$

$-\frac{1}{6}<H,i>(y_{0})[\nabla_{j}\varrho_{ij}-2\varrho_{jj}%
<H,i>+\overset{q}{\underset{\text{a}=1}{\sum}}(\nabla_{j}R_{\text{a}%
i\text{a}j}-4R_{j\text{a}j\text{a}}<H,i>)$

$+4\overset{q}{\underset{\text{a,b=1}}{\sum}}R_{j\text{a}j\text{b}%
}T_{\text{ab}i}+2\overset{q}{\underset{\text{a,b,c=1}}{\sum}}(T_{\text{aa}%
j}T_{\text{bb}j}T_{\text{cc}i}-3T_{\text{aa}j}T_{\text{bc}j}T_{\text{bc}%
i}+2T_{\text{ab}j}T_{\text{bc}j}T_{\text{ac}i})](y_{0})$

$-\frac{1}{4}<H,j>^{2}(y_{0})[\varrho_{ii}+$ $\overset{q}{\underset{\text{a}%
=1}{2\sum}}R_{i\text{a}i\text{a}}-3\overset{q}{\underset{\text{a,b=1}}{\sum}%
}(T_{\text{aa}i}T_{\text{bb}i}-T_{\text{ab}i}T_{\text{ab}i})](y_{0})\qquad
L_{232}$

$-\frac{1}{18}[\varrho_{ii}+$ $\overset{q}{\underset{\text{a}=1}{2\sum}%
}R_{i\text{a}i\text{a}}-3\overset{q}{\underset{\text{a,b=1}}{\sum}%
}(T_{\text{aa}i}T_{\text{bb}i}-T_{\text{ab}i}T_{\text{ab}i})](y_{0})$

$\times\lbrack\varrho_{jj}+$ $\overset{q}{\underset{\text{a}=1}{2\sum}%
}R_{j\text{a}j\text{a}}-3\overset{q}{\underset{\text{a,b=1}}{\sum}%
}(T_{\text{aa}j}T_{\text{bb}j}-T_{\text{ab}j}T_{\text{ab}j})]\}(y_{0})$

$+\frac{1}{2}R_{ijik}(y_{0})$ \ $[<H,j><H,k>](y_{0})\qquad\qquad\qquad\qquad
$\ $L_{233}$

$+\frac{1}{18}R_{ijik}(y_{0})[2\varrho_{jk}+$ $\overset{q}{\underset{\text{a}%
=1}{4\sum}}R_{j\text{a}k\text{a}}-3\overset{q}{\underset{\text{a,b=1}}{\sum}%
}(T_{\text{aa}j}T_{\text{bb}k}-T_{\text{ab}j}T_{\text{ab}k})$

$-3\overset{q}{\underset{\text{a,b=1}}{\sum}}(T_{\text{aa}k}T_{\text{bb}%
j}-T_{\text{ab}k}T_{\text{ab}j}](y_{0})+$ $\frac{\partial^{4}\theta^{-\frac
{1}{2}}}{\partial x_{i}^{2}\partial x_{j}^{2}}(y_{0})$

We now come to the computation of:

$L_{3}=2\frac{\partial\theta^{-1}}{\partial x_{i}}(y_{0})\frac{\partial
}{\partial x_{i}}[\frac{\partial\theta}{\partial\text{x}_{j}}(g^{jk}%
\frac{\partial\theta^{-\frac{1}{2}}}{\partial\text{x}_{k}})+\theta
(\frac{\partial g^{jk}}{\partial\text{x}_{j}}\frac{\partial\theta^{-\frac
{1}{2}}}{\partial\text{x}_{k}}+g^{jk}\frac{\partial^{2}\theta^{-\frac{1}{2}}%
}{\partial\text{x}_{j}\partial\text{x}_{k}})](y_{0})$

We set:

$L_{\Delta}=\frac{\partial}{\partial x_{i}}[\frac{\partial\theta}%
{\partial\text{x}_{j}}(g^{jk}\frac{\partial\theta^{-\frac{1}{2}}}%
{\partial\text{x}_{k}})+\theta(\frac{\partial g^{jk}}{\partial\text{x}_{j}%
}\frac{\partial\theta^{-\frac{1}{2}}}{\partial\text{x}_{k}}+g^{jk}%
\frac{\partial^{2}\theta^{-\frac{1}{2}}}{\partial\text{x}_{j}\partial
\text{x}_{k}})](y_{0})$

We now compute,

$L_{\Delta}=\frac{\partial}{\partial x_{i}}[\frac{\partial\theta}%
{\partial\text{x}_{j}}(g^{jk}\frac{\partial\theta^{-\frac{1}{2}}}%
{\partial\text{x}_{k}})+\theta(\frac{\partial g^{jk}}{\partial\text{x}_{j}%
}\frac{\partial\theta^{-\frac{1}{2}}}{\partial\text{x}_{k}}+g^{jk}%
\frac{\partial^{2}\theta^{-\frac{1}{2}}}{\partial\text{x}_{j}\partial
\text{x}_{k}})](y_{0})$

\qquad$=[\frac{\partial^{2}\theta}{\partial x_{i}\partial\text{x}_{j}}%
(g^{jk}\frac{\partial\theta^{-\frac{1}{2}}}{\partial\text{x}_{k}}%
)+\frac{\partial\theta}{\partial\text{x}_{j}}(\frac{\partial\text{g}^{jk}%
}{\partial x_{i}}\frac{\partial\theta^{-\frac{1}{2}}}{\partial\text{x}_{k}%
}+g^{jk}\frac{\partial^{2}\theta^{-\frac{1}{2}}}{\partial x_{i}\partial
\text{x}_{k}})](y_{0})$

$+\theta(y_{0})[\frac{\partial^{2}\text{g}^{jk}}{\partial x_{i}\partial
\text{x}_{j}}\frac{\partial\theta^{-\frac{1}{2}}}{\partial\text{x}_{k}}%
+\frac{\partial g^{jk}}{\partial\text{x}_{j}}\frac{\partial^{2}\theta
^{-\frac{1}{2}}}{\partial x_{i}\partial\text{x}_{k}}](y_{0})+\theta
(y_{0})[\frac{\partial\text{g}^{jk}}{\partial x_{i}}\frac{\partial^{2}%
\theta^{-\frac{1}{2}}}{\partial\text{x}_{j}\partial\text{x}_{k}}$

$+g^{jk}\frac{\partial^{3}\theta^{-\frac{1}{2}}}{\partial x_{i}\partial
\text{x}_{j}\partial\text{x}_{k}}](y_{0})$

$+\frac{\partial\theta}{\partial\text{x}_{i}}(y_{0})[(\frac{\partial g^{jk}%
}{\partial\text{x}_{j}}\frac{\partial\theta^{-\frac{1}{2}}}{\partial
\text{x}_{k}}+g^{jk}\frac{\partial^{2}\theta^{-\frac{1}{2}}}{\partial
\text{x}_{j}\partial\text{x}_{k}})](y_{0})$

Since $g^{jk}(y_{0})=\delta^{jk}$ and $\theta(y_{0})=1,$ we have,

$L_{\Delta}=[\frac{\partial^{2}\theta}{\partial x_{i}\partial\text{x}_{j}%
}\frac{\partial\theta^{-\frac{1}{2}}}{\partial\text{x}_{j}}+\frac
{\partial\theta}{\partial\text{x}_{j}}(\frac{\partial\text{g}^{jk}}{\partial
x_{i}}\frac{\partial\theta^{-\frac{1}{2}}}{\partial\text{x}_{k}}%
+\frac{\partial^{2}\theta^{-\frac{1}{2}}}{\partial x_{i}\partial\text{x}_{j}%
})+\frac{\partial^{2}\text{g}^{jk}}{\partial x_{i}\partial\text{x}_{j}}%
\frac{\partial\theta^{-\frac{1}{2}}}{\partial\text{x}_{k}}$

$+\frac{\partial g^{jk}}{\partial\text{x}_{j}}\frac{\partial^{2}\theta
^{-\frac{1}{2}}}{\partial x_{i}\partial\text{x}_{k}}](y_{0})$

$+[\frac{\partial\text{g}^{jk}}{\partial x_{i}}\frac{\partial^{2}%
\theta^{-\frac{1}{2}}}{\partial\text{x}_{j}\partial\text{x}_{k}}%
+\frac{\partial^{3}\theta^{-\frac{1}{2}}}{\partial x_{i}\partial\text{x}%
_{j}^{2}}](y_{0})+\frac{\partial\theta}{\partial\text{x}_{i}}(y_{0}%
)[\frac{\partial g^{jk}}{\partial\text{x}_{j}}\frac{\partial\theta^{-\frac
{1}{2}}}{\partial\text{x}_{k}}+\frac{\partial^{2}\theta^{-\frac{1}{2}}%
}{\partial\text{x}_{j}^{2}}](y_{0})$

Since $\frac{\partial\text{g}^{jk}}{\partial x_{i}}(y_{0})=0$ for
$i,j,k=q+1,...,n,$ we have:

$[\frac{\partial\theta}{\partial\text{x}_{j}}\frac{\partial\text{g}^{jk}%
}{\partial x_{i}}\frac{\partial\theta^{-\frac{1}{2}}}{\partial\text{x}_{k}%
}](y_{0})=0=[\frac{\partial\text{g}^{jk}}{\partial x_{j}}\frac{\partial
^{2}\theta^{-\frac{1}{2}}}{\partial\text{x}_{i}\partial\text{x}_{k}}%
](y_{0})=0=[\frac{\partial g^{jk}}{\partial\text{x}_{i}}\frac{\partial
^{2}\theta^{-\frac{1}{2}}}{\partial x_{j}\partial\text{x}_{k}}]$

$=0=[\frac{\partial g^{jk}}{\partial\text{x}_{j}}\frac{\partial\theta
^{-\frac{1}{2}}}{\partial\text{x}_{k}}](y_{0})$

Consequently,

$L_{\Delta}=[\frac{\partial^{2}\theta}{\partial x_{i}\partial\text{x}_{j}%
}\frac{\partial\theta^{-\frac{1}{2}}}{\partial\text{x}_{j}}+\frac
{\partial\theta}{\partial\text{x}_{j}}\frac{\partial^{2}\theta^{-\frac{1}{2}}%
}{\partial x_{i}\partial\text{x}_{j}}+\frac{\partial^{2}\text{g}^{jk}%
}{\partial x_{i}\partial\text{x}_{j}}\frac{\partial\theta^{-\frac{1}{2}}%
}{\partial\text{x}_{k}}](y_{0})\ +\frac{\partial^{3}\theta^{-\frac{1}{2}}%
}{\partial x_{i}\partial\text{x}_{j}^{2}}(y_{0})$

$+\frac{\partial\theta}{\partial\text{x}_{i}}(y_{0})\frac{\partial^{2}%
\theta^{-\frac{1}{2}}}{\partial\text{x}_{j}^{2}}(y_{0})$

Since $\frac{\partial\theta}{\partial x_{i}}(y_{0})=-$ $<H,i>$ and
$\theta(y_{0})=1$ and $\frac{\partial\theta^{-1}}{\partial x_{i}}(y_{0})$

$=-$ $\theta^{-2}(y_{0})\frac{\partial\theta}{\partial x_{i}}(y_{0})=$
$<H,i>(y_{0})$

The expression for $\frac{\partial^{2}\theta}{\partial\text{x}_{i}%
\partial\text{x}_{j}}(y_{0})$ is in (v) of \textbf{Appendix A}$_{9}$ and that of

$\frac{\partial^{2}\theta^{-\frac{1}{2}}}{\partial\text{x}_{i}\partial
\text{x}_{j}}(y_{0})$ is in (ix) of \textbf{Appendix A}$_{9}:$

$\frac{\partial^{2}\theta}{\partial\text{x}_{i}\partial\text{x}_{j}}(y_{0})$

$=-\frac{1}{6}[2\varrho_{ij}+$ $\overset{q}{\underset{\text{a}=1}{4\sum}%
}R_{i\text{a}j\text{a}}-3\overset{q}{\underset{\text{a,b=1}}{\sum}%
}(T_{\text{aa}i}T_{\text{bb}j}-T_{\text{ab}i}T_{\text{ab}j}%
)-3\overset{q}{\underset{\text{a,b=1}}{\sum}}(T_{\text{aa}j}T_{\text{bb}%
i}-T_{\text{ab}j}T_{\text{ab}i}](y_{0})$

$\frac{\partial^{2}\theta^{-\frac{1}{2}}}{\partial\text{x}_{i}\partial
\text{x}_{j}}(y_{0})$ \ $=\frac{3}{4}<H,i>(y_{0})<H,j>(y_{0})$

$+\frac{1}{12}[2\varrho_{ij}+$ $\overset{q}{\underset{\text{a}=1}{4\sum}%
}R_{i\text{a}j\text{a}}-3\overset{q}{\underset{\text{a,b=1}}{\sum}%
}(T_{\text{aa}i}T_{\text{bb}j}-T_{\text{ab}i}T_{\text{ab}j}%
)-3\overset{q}{\underset{\text{a,b=1}}{\sum}}(T_{\text{aa}j}T_{\text{bb}%
i}-T_{\text{ab}j}T_{\text{ab}i}](y_{0})$

$\frac{\partial^{2}\text{g}^{jk}}{\partial x_{i}\partial\text{x}_{j}}=\frac
{1}{3}R_{ijjk}(y_{0})$ by (iv) of \textbf{Appendix A}$_{2}.$ The expression of
$\frac{\partial^{3}\theta^{-\frac{1}{2}}}{\partial\text{x}_{i}\partial
\text{x}_{j}^{2}}(y_{0})$

is given by (xvii) of \textbf{Appendix A}$_{9}.$ We have:

$L_{\Delta}=[\frac{\partial^{2}\theta}{\partial x_{i}\partial\text{x}_{j}%
}\frac{\partial\theta^{-\frac{1}{2}}}{\partial\text{x}_{j}}+\frac
{\partial\theta}{\partial\text{x}_{j}}\frac{\partial^{2}\theta^{-\frac{1}{2}}%
}{\partial x_{i}\partial\text{x}_{j}}+\frac{\partial^{2}\text{g}^{jk}%
}{\partial x_{i}\partial\text{x}_{j}}\frac{\partial\theta^{-\frac{1}{2}}%
}{\partial\text{x}_{k}}](y_{0})+[\frac{\partial\theta}{\partial\text{x}_{i}%
}\frac{\partial^{2}\theta^{-\frac{1}{2}}}{\partial\text{x}_{j}^{2}}%
\ +\frac{\partial^{3}\theta^{-\frac{1}{2}}}{\partial x_{i}\partial\text{x}%
_{j}^{2}}](y_{0})\qquad$

$=-\frac{1}{12}<H,j>(y_{0})[2\varrho_{ij}+$ $\overset{q}{\underset{\text{a}%
=1}{4\sum}}R_{i\text{a}j\text{a}}-3\overset{q}{\underset{\text{a,b=1}}{\sum}%
}(T_{\text{aa}i}T_{\text{bb}j}-T_{\text{ab}i}T_{\text{ab}j})$

$-3\overset{q}{\underset{\text{a,b=1}}{\sum}}(T_{\text{aa}j}T_{\text{bb}%
i}-T_{\text{ab}j}T_{\text{ab}i}](y_{0})$

$-\frac{3}{4}[<H,i>(y_{0})<H,j>^{2}](y_{0})$

$-\frac{1}{12}<H,j>(y_{0})[2\varrho_{ij}+$ $\overset{q}{\underset{\text{a}%
=1}{4\sum}}R_{i\text{a}j\text{a}}-3\overset{q}{\underset{\text{a,b=1}}{\sum}%
}(T_{\text{aa}i}T_{\text{bb}j}-T_{\text{ab}i}T_{\text{ab}j})$

$-3\overset{q}{\underset{\text{a,b=1}}{\sum}}(T_{\text{aa}j}T_{\text{bb}%
i}-T_{\text{ab}j}T_{\text{ab}i}](y_{0})$

$\ +\frac{1}{6}<H,k>(y_{0})R_{ijjk}(y_{0})$

$-\frac{3}{4}[<H,i><H,j>^{2}](y_{0})-\frac{1}{6}<H,i>(y_{0})$

$\times\lbrack\varrho_{jj}+$ $\overset{q}{\underset{\text{a}=1}{2\sum}%
}R_{j\text{a}j\text{a}}-3\overset{q}{\underset{\text{a,b=1}}{\sum}%
}(T_{\text{aa}j}T_{\text{bb}j}-T_{\text{ab}j}T_{\text{ab}j})](y_{0})$

$+\frac{15}{8}[<H,i><H,j>^{2}](y_{0})\qquad$\ $\frac{\partial^{3}%
\theta^{-\frac{1}{2}}}{\partial\text{x}_{i}\partial\text{x}_{j}^{2}}(y_{0})$

$+\frac{1}{4}<H,j>$

$\times\lbrack2\varrho_{ij}+$ $\overset{q}{\underset{\text{a}=1}{4\sum}%
}R_{i\text{a}j\text{a}}-3\overset{q}{\underset{\text{a,b=1}}{\sum}%
}(T_{\text{aa}i}T_{\text{bb}j}-T_{\text{ab}i}T_{\text{ab}j}%
)-3\overset{q}{\underset{\text{a,b=1}}{\sum}}(T_{\text{aa}j}T_{\text{bb}%
i}-T_{\text{ab}j}T_{\text{ab}i}](y_{0})$

$+\frac{1}{4}<H,i>(y_{0})[\tau^{M}\ -3\tau^{P}+\ \underset{\text{a}%
=1}{\overset{\text{q}}{\sum}}\varrho_{\text{aa}}^{M}+$
$\overset{q}{\underset{\text{a},\text{b}=1}{\sum}}R_{\text{abab}}^{M}$
$](y_{0})$

$+\frac{1}{12}[\nabla_{i}\varrho_{jj}-2\varrho_{ij}%
<H,j>+\overset{q}{\underset{\text{a}=1}{\sum}}(\nabla_{i}R_{\text{a}%
j\text{a}j}-4R_{i\text{a}j\text{a}}<H,j>)\qquad\frac{\partial^{3}\theta
}{\partial\text{x}_{i}\partial\text{x}_{j}^{2}}(y_{0})$

$+4\overset{q}{\underset{\text{a,b=1}}{\sum}}R_{i\text{a}j\text{b}%
}T_{\text{ab}j}+2\overset{q}{\underset{\text{a,b,c=1}}{\sum}}(T_{\text{aa}%
i}T_{\text{bb}j}T_{\text{cc}j}-3T_{\text{aa}i}T_{\text{bc}j}T_{\text{bc}%
j}+2T_{\text{ab}i}T_{\text{bc}j}T_{\text{ca}j})](y_{0})$\qquad\qquad
\qquad\qquad\qquad\ \ 

$+\frac{1}{12}[\nabla_{j}\varrho_{ij}-2\varrho_{ij}%
<H,j>+\overset{q}{\underset{\text{a}=1}{\sum}}(\nabla_{j}R_{\text{a}%
i\text{a}j}-4R_{j\text{a}i\text{a}}<H,j>)$

$+4\overset{q}{\underset{\text{a,b=1}}{\sum}}R_{j\text{a}i\text{b}%
}T_{\text{ab}j}+2\overset{q}{\underset{\text{a,b,c=1}}{\sum}}(T_{\text{aa}%
j}T_{\text{bb}i}T_{\text{cc}j}-3T_{\text{aa}j}T_{\text{bc}i}T_{\text{bc}%
j}+2T_{\text{ab}j}T_{\text{bc}i}T_{\text{ac}j})](y_{0})$

$+\frac{1}{12}[\nabla_{j}\varrho_{ij}-2\varrho_{jj}%
<H,i>+\overset{q}{\underset{\text{a}=1}{\sum}}(\nabla_{j}R_{\text{a}%
i\text{a}j}-4R_{j\text{a}j\text{a}}<H,i>)+4\overset{q}{\underset{\text{a,b=1}%
}{\sum}}R_{j\text{a}j\text{b}}T_{\text{ab}i}$

$+2\overset{q}{\underset{\text{a,b,c=1}}{\sum}}(T_{\text{aa}j}T_{\text{bb}%
j}T_{\text{cc}i}-3T_{\text{aa}j}T_{\text{bc}j}T_{\text{bc}i}+2T_{\text{ab}%
j}T_{\text{bc}j}T_{\text{ac}i})](y_{0})$

\qquad\qquad\qquad\qquad\qquad\qquad\qquad\qquad\qquad\qquad\qquad\qquad
\qquad\qquad\qquad\qquad$\blacksquare$

$\qquad L_{3}=2\frac{\partial\theta^{-1}}{\partial x_{i}}(y_{0})L_{\Delta
}=2<H,i>(y_{0})L_{\Delta}$

Using the expression of $L_{\Delta}$ given above, we have:

$L_{3}=2<H,i>(y_{0})L_{\Delta}$

$=-\frac{1}{6}<H,i>(y_{0})<H,j>(y_{0})$

$\times\lbrack2\varrho_{ij}+$ $\overset{q}{\underset{\text{a}=1}{4\sum}%
}R_{i\text{a}j\text{a}}-3\overset{q}{\underset{\text{a,b=1}}{\sum}%
}(T_{\text{aa}i}T_{\text{bb}j}-T_{\text{ab}i}T_{\text{ab}j}%
)-3\overset{q}{\underset{\text{a,b=1}}{\sum}}(T_{\text{aa}j}T_{\text{bb}%
i}-T_{\text{ab}j}T_{\text{ab}i}](y_{0})$

$-\frac{3}{2}[<H,i>^{2}(y_{0})<H,j>^{2}](y_{0})$

$-\frac{1}{6}<H,i>(y_{0})<H,j>(y_{0})$

$\times\lbrack2\varrho_{ij}+$ $\overset{q}{\underset{\text{a}=1}{4\sum}%
}R_{i\text{a}j\text{a}}-3\overset{q}{\underset{\text{a,b=1}}{\sum}%
}(T_{\text{aa}i}T_{\text{bb}j}-T_{\text{ab}i}T_{\text{ab}j}%
)-3\overset{q}{\underset{\text{a,b=1}}{\sum}}(T_{\text{aa}j}T_{\text{bb}%
i}-T_{\text{ab}j}T_{\text{ab}i}](y_{0})$

$-\frac{1}{3}<H,i>(y_{0})<H,k>(y_{0})R_{jijk}(y_{0})$

$-\frac{3}{2}<H,i>^{2}(y_{0})<H,j>^{2}(y_{0})$

$-\frac{1}{3}<H,i>^{2}(y_{0})[\varrho_{jj}+$ $\overset{q}{\underset{\text{a}%
=1}{2\sum}}R_{j\text{a}j\text{a}}-3\overset{q}{\underset{\text{a,b=1}}{\sum}%
}(T_{\text{aa}j}T_{\text{bb}j}-T_{\text{ab}j}T_{\text{ab}j})](y_{0})$

$+\frac{15}{4}[<H,i>^{2}<H,j>^{2}](y_{0})\qquad$\ $\frac{\partial^{3}%
\theta^{-\frac{1}{2}}}{\partial\text{x}_{i}\partial\text{x}_{j}^{2}}%
(y_{0})<H,i>(y_{0})$

$+\frac{1}{2}<H,i>(y_{0})<H,j>$

$\times\lbrack2\varrho_{ij}+$ $\overset{q}{\underset{\text{a}=1}{4\sum}%
}R_{i\text{a}j\text{a}}-3\overset{q}{\underset{\text{a,b=1}}{\sum}%
}(T_{\text{aa}i}T_{\text{bb}j}-T_{\text{ab}i}T_{\text{ab}j}%
)-3\overset{q}{\underset{\text{a,b=1}}{\sum}}(T_{\text{aa}j}T_{\text{bb}%
i}-T_{\text{ab}j}T_{\text{ab}i}](y_{0})$

$+\frac{1}{2}<H,i>(y_{0})<H,i>(y_{0})[\tau^{M}\ -3\tau^{P}%
+\ \underset{\text{a}=1}{\overset{\text{q}}{\sum}}\varrho_{\text{aa}}^{M}+$
$\overset{q}{\underset{\text{a},\text{b}=1}{\sum}}R_{\text{abab}}^{M}$
$](y_{0})$

$+\frac{1}{6}<H,i>(y_{0})[\nabla_{i}\varrho_{jj}-2\varrho_{ij}%
<H,j>+\overset{q}{\underset{\text{a}=1}{\sum}}(\nabla_{i}R_{\text{a}%
j\text{a}j}-4R_{i\text{a}j\text{a}}<H,j>)$

$+4\overset{q}{\underset{\text{a,b=1}}{\sum}}R_{i\text{a}j\text{b}%
}T_{\text{ab}j}+2\overset{q}{\underset{\text{a,b,c=1}}{\sum}}(T_{\text{aa}%
i}T_{\text{bb}j}T_{\text{cc}j}-3T_{\text{aa}i}T_{\text{bc}j}T_{\text{bc}%
j}+2T_{\text{ab}i}T_{\text{bc}j}T_{\text{ca}j})](y_{0})$\qquad\qquad
\qquad\qquad\qquad\ \ 

$+\frac{1}{6}<H,i>(y_{0})[\nabla_{j}\varrho_{ij}-2\varrho_{ij}%
<H,j>+\overset{q}{\underset{\text{a}=1}{\sum}}(\nabla_{j}R_{\text{a}%
i\text{a}j}-4R_{j\text{a}i\text{a}}<H,j>)$

$+4\overset{q}{\underset{\text{a,b=1}}{\sum}}R_{j\text{a}i\text{b}%
}T_{\text{ab}j}+2\overset{q}{\underset{\text{a,b,c=1}}{\sum}}(T_{\text{aa}%
j}T_{\text{bb}i}T_{\text{cc}j}-3T_{\text{aa}j}T_{\text{bc}i}T_{\text{bc}%
j}+2T_{\text{ab}j}T_{\text{bc}i}T_{\text{ac}j})](y_{0})$

$+\frac{1}{6}<H,i>(y_{0})[\nabla_{j}\varrho_{ij}-2\varrho_{jj}%
<H,i>+\overset{q}{\underset{\text{a}=1}{\sum}}(\nabla_{j}R_{\text{a}%
i\text{a}j}-4R_{j\text{a}j\text{a}}<H,i>)$

$+4\overset{q}{\underset{\text{a,b=1}}{\sum}}R_{j\text{a}j\text{b}%
}T_{\text{ab}i}+2\overset{q}{\underset{\text{a,b,c=1}}{\sum}}(T_{\text{aa}%
j}T_{\text{bb}j}T_{\text{cc}i}-3T_{\text{aa}j}T_{\text{bc}j}T_{\text{bc}%
i}+2T_{\text{ab}j}T_{\text{bc}j}T_{\text{ac}i})](y_{0})$

$\qquad\qquad\qquad\qquad\qquad\qquad\qquad\qquad\qquad\qquad\qquad
\qquad\qquad\qquad\qquad\qquad\qquad\qquad\qquad\blacksquare\qquad\qquad$

We now give the expression for $\frac{1}{24}\frac{\partial^{2}}{\partial
x_{i}^{2}}(\Delta\theta^{-\frac{1}{2}})(y_{0})$ which includes the expression
of $\frac{1}{24}$ $\frac{\partial^{4}\theta^{-\frac{1}{2}}}{\partial x_{i}%
^{2}\partial x_{j}^{2}}(y_{0})$ obtained earlier:

$\left(  A_{29}\right)  \qquad A_{3212}=\frac{1}{24}\frac{\partial^{2}%
}{\partial x_{i}^{2}}(\Delta\theta^{-\frac{1}{2}})(y_{0})=\frac{1}{24}%
(L_{1}+L_{2}+$ $L_{3})$

\qquad$=\frac{1}{24}[2<H,i>^{2}(y_{0})+\frac{1}{3}(\tau^{M}-3\tau
^{P}+\overset{q}{\underset{\text{a}=1}{\sum}}\varrho_{\text{aa}}%
+\overset{q}{\underset{\text{a,b}=1}{\sum}}R_{\text{abab}})](y_{0}%
)\qquad\qquad L_{1}$

$\qquad\times\lbrack\frac{1}{4}<H,j>^{2}(y_{0})+\frac{1}{6}(\tau^{M}-3\tau
^{P}+\overset{q}{\underset{\text{a}=1}{\sum}}\varrho_{\text{aa}}%
^{M}+\overset{q}{\underset{\text{a,b}=1}{\sum}}R_{\text{abab}}^{M})](y_{0})$

$-\frac{1}{4}\times\frac{1}{24}[2\varrho_{ij}+$
$\overset{q}{\underset{\text{a}=1}{4\sum}}R_{i\text{a}j\text{a}}%
-3\overset{q}{\underset{\text{a,b=1}}{\sum}}(T_{\text{aa}i}T_{\text{bb}%
j}-T_{\text{ab}i}T_{\text{ab}j})-3\overset{q}{\underset{\text{a,b=1}}{\sum}%
}(T_{\text{aa}j}T_{\text{bb}i}-T_{\text{ab}j}T_{\text{ab}i}](y_{0})\qquad
L_{2}\qquad L_{21}\qquad L_{211}$

$\qquad\ \ \ \ \ \times\lbrack<H,i><H,j>](y_{0})$

$\qquad-\frac{1}{24}\times\frac{1}{36}[2\varrho_{ij}+$
$\overset{q}{\underset{\text{a}=1}{4\sum}}R_{i\text{a}j\text{a}}%
-3\overset{q}{\underset{\text{a,b=1}}{\sum}}(T_{\text{aa}i}T_{\text{bb}%
j}-T_{\text{ab}i}T_{\text{ab}j})-3\overset{q}{\underset{\text{a,b=1}}{\sum}%
}(T_{\text{aa}j}T_{\text{bb}i}-T_{\text{ab}j}T_{\text{ab}i}]^{2}(y_{0})$

\qquad$-\frac{1}{24}\times\frac{1}{12}[<H,j>](y_{0})\times\lbrack\{\nabla
_{i}\varrho_{ij}-2\varrho_{ij}<H,i>+\overset{q}{\underset{\text{a}=1}{\sum}%
}(\nabla_{i}R_{\text{a}i\text{a}j}-4R_{i\text{a}j\text{a}}<H,i>)\qquad
L_{212}$

$\qquad+4\overset{q}{\underset{\text{a,b=1}}{\sum}}R_{i\text{a}j\text{b}%
}T_{\text{ab}i}+2\overset{q}{\underset{\text{a,b,c=1}}{\sum}}(T_{\text{aa}%
i}T_{\text{bb}j}T_{\text{cc}i}-3T_{\text{aa}i}T_{\text{bc}j}T_{\text{bc}%
i}+2T_{\text{ab}i}T_{\text{bc}j}T_{\text{ac}i})](y_{0})$\qquad\qquad
\qquad\qquad\qquad\ \ 

$\qquad-\frac{1}{24}\times\frac{1}{12}[<H,j>](y_{0})\times\lbrack\nabla
_{j}\varrho_{ii}-2\varrho_{ij}<H,i>+\overset{q}{\underset{\text{a}=1}{\sum}%
}(\nabla_{j}R_{\text{a}i\text{a}i}-4R_{i\text{a}j\text{a}}<H,i>)$

$\qquad+4\overset{q}{\underset{\text{a,b=1}}{\sum}}R_{j\text{a}i\text{b}%
}T_{\text{ab}i}+2\overset{q}{\underset{\text{a,b,c=1}}{\sum}}(T_{\text{aa}%
j}T_{\text{bb}i}T_{\text{cc}i}-3T_{\text{aa}j}T_{\text{bc}i}T_{\text{bc}%
i}+2T_{\text{ab}j}T_{\text{bc}i}T_{\text{ac}i})](y_{0})$

$\qquad-\frac{1}{24}\times\frac{1}{12}[<H,j>](y_{0})\times\lbrack\nabla
_{i}\varrho_{ij}-2\varrho_{ii}<H,j>+\overset{q}{\underset{\text{a}=1}{\sum}%
}(\nabla_{i}R_{\text{a}i\text{a}j}-4R_{i\text{a}i\text{a}}<H,j>)$

$\qquad+4\overset{q}{\underset{\text{a,b=1}}{\sum}}R_{i\text{a}i\text{b}%
}T_{\text{ab}j}+2\overset{q}{\underset{\text{a,b,c}=1}{\sum}}(T_{\text{aa}%
i}T_{\text{bb}i}T_{\text{cc}j}-3T_{\text{aa}i}T_{\text{bc}i}T_{\text{bc}%
j}+2T_{\text{ab}i}T_{\text{bc}i}T_{\text{ac}j})](y_{0})$

$\qquad-\frac{1}{3}[<H,j><H,k>](y_{0})R_{ijik}(y_{0})-$ $\frac{1}{24}%
\times\frac{15}{8}[<H,i>^{2}<H,j>^{2}](y_{0})\qquad L_{213}$

$\qquad-\frac{1}{24}\times\frac{1}{4}<H,i><H,j>[2\varrho_{ij}%
+\overset{q}{\underset{\text{a}=1}{4\sum}}R_{i\text{a}j\text{a}}%
-3\overset{q}{\underset{\text{a,b=1}}{\sum}}(T_{\text{aa}i}T_{\text{bb}%
j}-T_{\text{ab}i}T_{\text{ab}j})$

$\qquad-3\overset{q}{\underset{\text{a,b=1}}{\sum}}(T_{\text{aa}j}%
T_{\text{bb}i}-T_{\text{ab}j}T_{\text{ab}i}](y_{0})$

$\qquad-\frac{1}{24}\times\frac{1}{4}<H,j>^{2}[\tau^{M}\ -3\tau^{P}%
+\ \underset{\text{a}=1}{\overset{\text{q}}{\sum}}\varrho_{\text{aa}}^{M}+$
$\overset{q}{\underset{\text{a},\text{b}=1}{\sum}}R_{\text{abab}}^{M}$
$](y_{0})$

$\qquad+\frac{1}{24}\times\frac{1}{12}<H,j>[\nabla_{i}\varrho_{ij}%
-2\varrho_{ij}<H,i>+\overset{q}{\underset{\text{a}=1}{\sum}}(\nabla
_{i}R_{\text{a}i\text{a}j}-4R_{i\text{a}j\text{a}}<H,i>)$

$+4\overset{q}{\underset{\text{a,b=1}}{\sum}}R_{i\text{a}j\text{b}%
}T_{\text{ab}i}+2\overset{q}{\underset{\text{a,b,c=1}}{\sum}}(T_{\text{aa}%
i}T_{\text{bb}j}T_{\text{cc}i}-3T_{\text{aa}i}T_{\text{bc}j}T_{\text{bc}%
i}+2T_{\text{ab}i}T_{\text{bc}j}T_{\text{ac}i})](y_{0})$\qquad\qquad
\qquad\qquad\qquad\ \ 

$\qquad+\frac{1}{12}<H,j>[\nabla_{j}\varrho_{ii}-2\varrho_{ij}%
<H,i>+\overset{q}{\underset{\text{a}=1}{\sum}}(\nabla_{j}R_{\text{a}%
i\text{a}i}-4R_{i\text{a}j\text{a}}<H,i>)$

$+4\overset{q}{\underset{\text{a,b=1}}{\sum}}R_{j\text{a}i\text{b}%
}T_{\text{ab}i}+2\overset{q}{\underset{\text{a,b,c=1}}{\sum}}(T_{\text{aa}%
j}T_{\text{bb}i}T_{\text{cc}i}-3T_{\text{aa}j}T_{\text{bc}i}T_{\text{bc}%
i}+2T_{\text{ab}j}T_{\text{bc}i}T_{\text{ac}i})](y_{0})$

$\qquad+\frac{1}{24}\times\frac{1}{12}<H,j>[\nabla_{i}\varrho_{ij}%
-2\varrho_{ii}<H,j>+\overset{q}{\underset{\text{a}=1}{\sum}}(\nabla
_{i}R_{\text{a}i\text{a}j}-4R_{i\text{a}i\text{a}}<H,j>)$

$+4\overset{q}{\underset{\text{a,b=1}}{\sum}}R_{i\text{a}i\text{b}%
}T_{\text{ab}j}+2\overset{q}{\underset{\text{a,b,c}=1}{\sum}}(T_{\text{aa}%
i}T_{\text{bb}i}T_{\text{cc}j}-3T_{\text{aa}i}T_{\text{bc}i}T_{\text{bc}%
j}+2T_{\text{ab}i}T_{\text{bc}i}T_{\text{ac}j})](y_{0})$

\qquad$-\frac{1}{24}\times\frac{1}{6}R_{jijk}(y_{0})$ \ $[<H,i><H,k>](y_{0}%
)\qquad\qquad L_{22}$

$-\frac{1}{24}\times\frac{1}{18}R_{jijk}(y_{0})[2\varrho_{ik}+$
$\overset{q}{\underset{\text{a}=1}{4\sum}}R_{i\text{a}k\text{a}}%
-3\overset{q}{\underset{\text{a,b=1}}{\sum}}(T_{\text{aa}i}T_{\text{bb}%
k}-T_{\text{ab}i}T_{\text{ab}k})$

$-3\overset{q}{\underset{\text{a,b=1}}{\sum}}(T_{\text{aa}k}T_{\text{bb}%
i}-T_{\text{ab}k}T_{\text{ab}i}](y_{0})$

$+\frac{1}{24}\times\frac{1}{6}<H,k>(y_{0})[\nabla_{j}$R$_{ijik}(y_{0}%
)-\nabla_{i}$R$_{jijk}](y_{0})$

$-\frac{1}{24}\times\frac{15}{4}<H,i>^{2}(y_{0})<H,j>^{2}(y_{0})\qquad\qquad
L_{23}\qquad L_{231}$

$-\frac{1}{24}\times\frac{1}{2}<H,i>(y_{0})<H,j>(y_{0})$

$\times\lbrack2\varrho_{ij}+$ $\overset{q}{\underset{\text{a}=1}{4\sum}%
}R_{i\text{a}j\text{a}}-3\overset{q}{\underset{\text{a,b=1}}{\sum}%
}(T_{\text{aa}i}T_{\text{bb}j}-T_{\text{ab}i}T_{\text{ab}j}%
)-3\overset{q}{\underset{\text{a,b=1}}{\sum}}(T_{\text{aa}j}T_{\text{bb}%
i}-T_{\text{ab}j}T_{\text{ab}i}](y_{0})$

$-\frac{1}{24}\times\frac{1}{2}<H,i>^{2}(y_{0})[\varrho_{jj}+$
$\overset{q}{\underset{\text{a}=1}{2\sum}}R_{j\text{a}j\text{a}}%
-3\overset{q}{\underset{\text{a,b=1}}{\sum}}(T_{\text{aa}j}T_{\text{bb}%
j}-T_{\text{ab}j}T_{\text{ab}j})](y_{0})$

$-\frac{1}{24}\times\frac{1}{6}<H,i>(y_{0})[\nabla_{i}\varrho_{jj}%
-2\varrho_{ij}<H,j>+\overset{q}{\underset{\text{a}=1}{\sum}}(\nabla
_{i}R_{\text{a}j\text{a}j}-4R_{i\text{a}j\text{a}}<H,j>)$

$+4\overset{q}{\underset{\text{a,b=1}}{\sum}}R_{i\text{a}j\text{b}%
}T_{\text{ab}j}+2\overset{q}{\underset{\text{a,b,c=1}}{\sum}}(T_{\text{aa}%
i}T_{\text{bb}j}T_{\text{cc}j}-3T_{\text{aa}i}T_{\text{bc}j}T_{\text{bc}%
j}+2T_{\text{ab}i}T_{\text{bc}j}T_{\text{ca}j})](y_{0})$\qquad\qquad
\qquad\qquad\qquad\ \ 

$-\frac{1}{24}\times\frac{1}{6}<H,i>(y_{0})[\nabla_{j}\varrho_{ij}%
-2\varrho_{ij}<H,j>+\overset{q}{\underset{\text{a}=1}{\sum}}(\nabla
_{j}R_{\text{a}i\text{a}j}-4R_{j\text{a}i\text{a}}<H,j>)$

$+4\overset{q}{\underset{\text{a,b=1}}{\sum}}R_{j\text{a}i\text{b}%
}T_{\text{ab}j}+2\overset{q}{\underset{\text{a,b,c=1}}{\sum}}(T_{\text{aa}%
j}T_{\text{bb}i}T_{\text{cc}j}-3T_{\text{aa}j}T_{\text{bc}i}T_{\text{bc}%
j}+2T_{\text{ab}j}T_{\text{bc}i}T_{\text{ac}j})](y_{0})$

$-\frac{1}{24}\times\frac{1}{6}<H,i>(y_{0})[\nabla_{j}\varrho_{ij}%
-2\varrho_{jj}<H,i>+\overset{q}{\underset{\text{a}=1}{\sum}}(\nabla
_{j}R_{\text{a}i\text{a}j}-4R_{j\text{a}j\text{a}}<H,i>)$

$+4\overset{q}{\underset{\text{a,b=1}}{\sum}}R_{j\text{a}j\text{b}%
}T_{\text{ab}i}+2\overset{q}{\underset{\text{a,b,c=1}}{\sum}}(T_{\text{aa}%
j}T_{\text{bb}j}T_{\text{cc}i}-3T_{\text{aa}j}T_{\text{bc}j}T_{\text{bc}%
i}+2T_{\text{ab}j}T_{\text{bc}j}T_{\text{ac}i})](y_{0})$

$-\frac{1}{24}\times\frac{1}{4}<H,j>^{2}(y_{0})[\varrho_{ii}+$
$\overset{q}{\underset{\text{a}=1}{2\sum}}R_{i\text{a}i\text{a}}%
-3\overset{q}{\underset{\text{a,b=1}}{\sum}}(T_{\text{aa}i}T_{\text{bb}%
i}-T_{\text{ab}i}T_{\text{ab}i})](y_{0})\qquad L_{232}$

$\qquad-\frac{1}{24}\times\frac{1}{18}[\varrho_{ii}+$
$\overset{q}{\underset{\text{a}=1}{2\sum}}R_{i\text{a}i\text{a}}%
-3\overset{q}{\underset{\text{a,b=1}}{\sum}}(T_{\text{aa}i}T_{\text{bb}%
i}-T_{\text{ab}i}T_{\text{ab}i})](y_{0})$

\qquad$\times\lbrack\varrho_{jj}+$ $\overset{q}{\underset{\text{a}=1}{2\sum}%
}R_{j\text{a}j\text{a}}-3\overset{q}{\underset{\text{a,b=1}}{\sum}%
}(T_{\text{aa}j}T_{\text{bb}j}-T_{\text{ab}j}T_{\text{ab}j})]\}(y_{0})$

$\qquad+\frac{1}{24}\times\frac{1}{2}R_{ijik}(y_{0})$ \ $[<H,j><H,k>](y_{0}%
)\qquad\qquad\qquad\qquad$\ $L_{233}$

$\qquad+\frac{1}{24}\times\frac{1}{18}R_{ijik}(y_{0})[2\varrho_{jk}+$
$\overset{q}{\underset{\text{a}=1}{4\sum}}R_{j\text{a}k\text{a}}%
-3\overset{q}{\underset{\text{a,b=1}}{\sum}}(T_{\text{aa}j}T_{\text{bb}%
k}-T_{\text{ab}j}T_{\text{ab}k})$

$-3\overset{q}{\underset{\text{a,b=1}}{\sum}}(T_{\text{aa}k}T_{\text{bb}%
j}-T_{\text{ab}k}T_{\text{ab}j}](y_{0})$

$+\overset{n}{\underset{i,j=q+1}{\sum}}\frac{35}{128}<H,i>^{2}(y_{0}%
)<H,j>^{2}(y_{0})\qquad\qquad\ \frac{1}{24}\frac{\partial^{4}\theta^{-\frac
{1}{2}}}{\partial x_{i}^{2}\partial x_{j}^{2}}(y_{0})$

$+\frac{5}{192}\overset{n}{\underset{j=q+1}{\sum}}<H,j>^{2}(y_{0})[\tau
^{M}\ -3\tau^{P}+\ \underset{\text{a}=1}{\overset{\text{q}}{\sum}}%
\varrho_{\text{aa}}^{M}+\overset{q}{\underset{\text{a},\text{b}=1}{\sum}%
}R_{\text{abab}}^{M}](y_{0})\qquad\ \ \ \ \ \ \ \ $

$+\frac{5}{192}\overset{n}{\underset{i=q+1}{\sum}}<H,i>^{2}(y_{0})[\tau
^{M}\ -3\tau^{P}+\ \underset{\text{a}=1}{\overset{\text{q}}{\sum}}%
\varrho_{\text{aa}}^{M}+\overset{q}{\underset{\text{a},\text{b}=1}{\sum}%
}R_{\text{abab}}^{M}](y_{0})\qquad\qquad$

$+\frac{5}{192}\overset{n}{\underset{i,j=q+1}{\sum}}[<H,i><H,j>](y_{0}%
)\qquad\qquad\qquad\qquad\qquad\qquad\qquad\qquad$

$\times\lbrack2\varrho_{ij}+4\overset{q}{\underset{\text{a}=1}{\sum}%
}R_{i\text{a}j\text{a}}-3\overset{q}{\underset{\text{a,b=1}}{\sum}%
}(T_{\text{aa}i}T_{\text{bb}j}-T_{\text{ab}i}T_{\text{ab}j}%
)-3\overset{q}{\underset{\text{a,b=1}}{\sum}}(T_{\text{aa}j}T_{\text{bb}%
i}-T_{\text{ab}j}T_{\text{ab}i})](y_{0})$

$+\frac{1}{96}\overset{n}{\underset{i,j=q+1}{\sum}}<H,j>(y_{0})[\{\nabla
_{i}\varrho_{ij}-2\varrho_{ij}<H,i>+\overset{q}{\underset{\text{a}=1}{\sum}%
}(\nabla_{i}R_{\text{a}i\text{a}j}-4R_{i\text{a}j\text{a}}<H,i>)\qquad$

$+4\overset{q}{\underset{\text{a,b=1}}{\sum}}R_{i\text{a}j\text{b}%
}T_{\text{ab}i}+2\overset{q}{\underset{\text{a,b,c=1}}{\sum}}(T_{\text{aa}%
i}T_{\text{bb}j}T_{\text{cc}i}-T_{\text{aa}i}T_{\text{bc}j}T_{\text{bc}%
i}-2T_{\text{bc}j}(T_{\text{aa}i}T_{\text{bc}i}-T_{\text{ab}i}T_{\text{ac}%
i}))\}$\qquad\qquad\qquad\ \ 

$+\{\nabla_{j}\varrho_{ii}-2\varrho_{ij}<H,i>+\overset{q}{\underset{\text{a}%
=1}{\sum}}(\nabla_{j}R_{\text{a}i\text{a}i}-4R_{i\text{a}j\text{a}}<H,i>)$

$+4\overset{q}{\underset{\text{a,b=1}}{\sum}}R_{j\text{a}i\text{b}%
}T_{\text{ab}i}+2\overset{q}{\underset{\text{a,b,c=1}}{\sum}}(T_{\text{aa}%
j}(T_{\text{bb}i}T_{\text{cc}i}-T_{\text{bc}i}T_{\text{bc}i})-2T_{\text{aa}%
j}T_{\text{bc}i}T_{\text{bc}i}+2T_{\text{ab}j}T_{\text{bc}i}T_{\text{ac}%
i})\}\qquad$

$+\{\nabla_{i}\varrho_{ij}-2\varrho_{ii}<H,j>+\overset{q}{\underset{\text{a}%
=1}{\sum}}(\nabla_{i}R_{\text{a}i\text{a}j}-4R_{i\text{a}i\text{a}%
}<H,j>)+4\overset{q}{\underset{\text{a,b=1}}{\sum}}R_{i\text{a}i\text{b}%
}T_{\text{ab}j}$

$+2\overset{q}{\underset{\text{a,b,c}=1}{\sum}}(T_{\text{aa}i}T_{\text{bb}%
i}T_{\text{cc}j}-3T_{\text{aa}i}T_{\text{bc}i}T_{\text{bc}j}+2T_{\text{ab}%
i}T_{\text{bc}i}T_{\text{ac}j})\}](y_{0})$

$+\frac{1}{96}\overset{n}{\underset{i,j=q+1}{\sum}}<H,i>(y_{0})[\{\nabla
_{i}\varrho_{jj}-2\varrho_{ij}<H,j>+\overset{q}{\underset{\text{a}=1}{\sum}%
}(\nabla_{i}R_{\text{a}j\text{a}j}-4R_{i\text{a}j\text{a}}<H,j>)\qquad$

$+4\overset{q}{\underset{\text{a,b=1}}{\sum}}R_{i\text{a}j\text{b}%
}T_{\text{ab}j}+2\overset{q}{\underset{\text{a,b,c=1}}{\sum}}T_{\text{aa}%
i}(T_{\text{bb}j}T_{\text{cc}j}-T_{\text{bc}j}T_{\text{bc}j})-2T_{\text{aa}%
i}T_{\text{bc}j}T_{\text{bc}j}+2T_{\text{ab}i}T_{\text{bc}j}T_{\text{ac}%
j})\}(y_{0})\qquad$\qquad\qquad\qquad\qquad\qquad\ \ 

$+\{\nabla_{j}\varrho_{ij}-2\varrho_{ij}<H,j>+\overset{q}{\underset{\text{a}%
=1}{\sum}}(\nabla_{j}R_{\text{a}i\text{a}j}-4R_{j\text{a}i\text{a}}<H,j>)$

$+4\overset{q}{\underset{\text{a,b=1}}{\sum}}R_{j\text{a}i\text{b}%
}T_{\text{ab}j}+2\overset{q}{\underset{\text{a,b,c=1}}{\sum}}(T_{\text{aa}%
j}T_{\text{bb}i}T_{\text{cc}j}-T_{\text{ab}j}T_{\text{bc}i}T_{\text{ac}%
j}-2T_{\text{bc}i}(T_{\text{aa}j}T_{\text{bc}j}-T_{\text{ab}j}T_{\text{ac}%
j}))\}(y_{0})$

$+\{\nabla_{j}\varrho_{ij}-2\varrho_{jj}<H,i>+\overset{q}{\underset{\text{a}%
=1}{\sum}}(\nabla_{j}R_{\text{a}i\text{a}j}-4R_{j\text{a}j\text{a}%
}<H,i>)+4\overset{q}{\underset{\text{a,b=1}}{\sum}}R_{j\text{a}j\text{b}%
}T_{\text{ab}i}$

$+2\overset{q}{\underset{\text{a,b,c=1}}{\sum}}(T_{\text{aa}j}T_{\text{bb}%
j}T_{\text{cc}i}-3T_{\text{aa}j}T_{\text{bc}j}T_{\text{bc}i}+2T_{\text{ab}%
j}T_{\text{bc}j}T_{\text{ac}i})\}](y_{0})$

$+\frac{1}{576}\overset{n}{\underset{i,j=q+1}{\sum}}[2\varrho_{ij}%
+4\overset{q}{\underset{\text{a}=1}{\sum}}R_{i\text{a}j\text{a}}%
-3\overset{q}{\underset{\text{a,b=1}}{\sum}}(T_{\text{aa}i}T_{\text{bb}%
j}-T_{\text{ab}i}T_{\text{ab}j})-3\overset{q}{\underset{\text{a,b=1}}{\sum}%
}(T_{\text{aa}j}T_{\text{bb}i}-T_{\text{ab}j}T_{\text{ab}i})]^{2}(y_{0})$

$+\frac{1}{288}[\tau^{M}\ -3\tau^{P}+\ \underset{\text{a}=1}{\overset{\text{q}%
}{\sum}}\varrho_{\text{aa}}^{M}+\overset{q}{\underset{\text{a},\text{b}%
=1}{\sum}}R_{\text{abab}}^{M}]^{2}(y_{0})$

$-\ \frac{1}{288}\overset{n}{\underset{i,j=q+1}{\sum}}[$
$\overset{q}{\underset{\text{a=1}}{\sum}}\{-(\nabla_{ii}^{2}R_{j\text{a}%
j\text{a}}+\nabla_{jj}^{2}R_{i\text{a}i\text{a}}+4\nabla_{ij}^{2}%
R_{i\text{a}j\text{a}}+2R_{ij}R_{i\text{a}j\text{a}})\qquad A$

$\qquad\qquad\qquad+\overset{n}{\underset{p=q+1}{\sum}}%
\overset{q}{\underset{\text{a=1}}{\sum}}(R_{\text{a}iip}R_{\text{a}%
jjp}+R_{\text{a}jjp}R_{\text{a}iip}+R_{\text{a}ijp}R_{\text{a}ijp}%
+R_{\text{a}ijp}R_{\text{a}jip}+R_{\text{a}jip}R_{\text{a}ijp}+R_{\text{a}%
jip}R_{\text{a}jip})$

$+2\overset{q}{\underset{\text{a,b=1}}{\sum}}\nabla_{i}(R)_{\text{a}%
i\text{b}j}T_{\text{ab}j}+2\overset{q}{\underset{\text{a,b=1}}{\sum}}%
\nabla_{j}(R)_{\text{a}j\text{b}i}T_{\text{ab}i}%
+2\overset{q}{\underset{\text{a,b=1}}{\sum}}\nabla_{i}(R)_{\text{a}j\text{b}%
i}T_{\text{ab}j}+2\overset{q}{\underset{\text{a,b=1}}{\sum}}\nabla
_{i}(R)_{\text{a}j\text{b}j}T_{\text{ab}i}$

$+2\overset{q}{\underset{\text{a,b=1}}{\sum}}\nabla_{j}(R)_{\text{a}%
i\text{b}i}T_{\text{ab}j}+2\overset{q}{\underset{\text{a,b=1}}{\sum}}%
\nabla_{j}(R)_{\text{a}i\text{b}j}T_{\text{ab}i}$

$+\overset{n}{\underset{p=q+1}{\sum}}(-\frac{3}{5}\nabla_{ii}^{2}%
(R)_{jpjp}+\overset{n}{\underset{p=q+1}{\sum}}(-\frac{3}{5}\nabla_{jj}%
^{2}(R)_{ipip}$

$+\overset{n}{\underset{p=q+1}{\sum}}(-\frac{3}{5}\nabla_{ij}^{2}%
(R)_{ipjp}+\overset{n}{\underset{p=q+1}{\sum}}(-\frac{3}{5}\nabla_{ij}%
^{2}(R)_{jpip}+\overset{n}{\underset{p=q+1}{\sum}}(-\frac{3}{5}\nabla_{ji}%
^{2}(R)_{ipjp}+\overset{n}{\underset{p=q+1}{\sum}}(-\frac{3}{5}\nabla_{ji}%
^{2}(R)_{jpip}$

$+\frac{1}{5}\overset{n}{\underset{m,p=q+1}{%
{\textstyle\sum}
}}R_{ipim}R_{jpjm}+\frac{1}{5}\overset{n}{\underset{m,p=q+1}{%
{\textstyle\sum}
}}R_{jpjm}R_{ipim}+\frac{1}{5}\overset{n}{\underset{m,p=q+1}{%
{\textstyle\sum}
}}R_{ipjm}R_{ipjm}+\frac{1}{5}\overset{n}{\underset{m,p=q+1}{%
{\textstyle\sum}
}}R_{ipjm}R_{jpim}$

$+\frac{1}{5}\overset{n}{\underset{m,p=q+1}{%
{\textstyle\sum}
}}R_{jpim}R_{ipjm}+\frac{1}{5}\overset{n}{\underset{m,p=q+1}{%
{\textstyle\sum}
}}R_{jpim}R_{jpim}\}(y_{0})$

$+4\overset{q}{\underset{\text{a,b=1}}{\sum}}\{(\nabla_{i}(R)_{i\text{a}%
j\text{a}}-\overset{q}{\underset{\text{c=1}}{%
{\textstyle\sum}
}}R_{\text{a}i\text{c}i}T_{\text{ac}j})$ $T_{\text{bb}j}+4(\nabla
_{j}(R)_{j\text{a}i\text{a}}-\overset{q}{\underset{\text{c=1}}{%
{\textstyle\sum}
}}R_{\text{a}j\text{c}j}T_{\text{ac}i})$ $T_{\text{bb}i}$

$+4(\nabla_{i}(R)_{j\text{a}i\text{a}}-\overset{q}{\underset{\text{c=1}}{%
{\textstyle\sum}
}}R_{\text{a}i\text{c}j}T_{\text{ac}i})$ $T_{\text{bb}j}$ $4B\ $

$+4(\nabla_{i}(R)_{j\text{a}j\text{a}}-\overset{q}{\underset{\text{c=1}}{%
{\textstyle\sum}
}}R_{\text{a}i\text{c}j}T_{\text{ac}j})$ $T_{\text{bb}i}+4(\nabla
_{j}(R)_{i\text{a}i\text{a}}-\overset{q}{\underset{\text{c=1}}{%
{\textstyle\sum}
}}R_{\text{a}j\text{c}i}T_{\text{ac}i})$ $T_{\text{bb}j}+4(\nabla
_{j}(R)_{i\text{a}j\text{a}}$

$-\overset{q}{\underset{\text{c=1}}{%
{\textstyle\sum}
}}R_{\text{a}j\text{c}i}T_{\text{ac}j})$ $T_{\text{bb}i}$

$-4\overset{q}{\underset{\text{a,b=1}}{\sum}}(\nabla_{i}(R)_{i\text{a}%
j\text{b}}-\overset{q}{\underset{\text{c=1}}{%
{\textstyle\sum}
}}R_{\text{b}r\text{c}s}T_{\text{ac}t})T_{\text{ab}j}%
-4\overset{q}{\underset{\text{a,b=1}}{\sum}}(\nabla_{j}(R)_{j\text{a}%
i\text{b}}-\overset{q}{\underset{\text{c=1}}{%
{\textstyle\sum}
}}R_{\text{b}j\text{c}j}T_{\text{ac}i})T_{\text{ab}i}$

$-4\overset{q}{\underset{\text{a,b=1}}{\sum}}(\nabla_{i}(R)_{j\text{a}%
i\text{b}}-\overset{q}{\underset{\text{c=1}}{%
{\textstyle\sum}
}}R_{\text{b}i\text{c}j}T_{\text{ac}i})T_{\text{ab}j}%
-4\overset{q}{\underset{\text{a,b=1}}{\sum}}(\nabla_{i}(R)_{j\text{a}%
j\text{b}}-\overset{q}{\underset{\text{c=1}}{%
{\textstyle\sum}
}}R_{\text{b}i\text{c}j}T_{\text{ac}j})T_{\text{ab}i}$

$-4\overset{q}{\underset{\text{a,b=1}}{\sum}}(\nabla_{j}(R)_{i\text{a}%
i\text{b}}-\overset{q}{\underset{\text{c=1}}{%
{\textstyle\sum}
}}R_{\text{b}j\text{c}i}T_{\text{ac}i})T_{\text{ab}j}%
-4\overset{q}{\underset{\text{a,b=1}}{\sum}}(\nabla_{j}(R)_{i\text{a}%
j\text{b}}-\overset{q}{\underset{\text{c=1}}{%
{\textstyle\sum}
}}R_{\text{b}j\text{c}i}T_{\text{ac}j})T_{\text{ab}i}\}](y_{0})$

$-\frac{1}{48}$ $[\frac{4}{9}\overset{q}{\underset{\text{a,b=1}}{\sum}%
}(\varrho_{\text{aa}}-\overset{q}{\underset{\text{c}=1}{\sum}}R_{\text{acac}%
})(\varrho_{\text{bb}}-\overset{q}{\underset{\text{d}=1}{\sum}}R_{\text{bdbd}%
})+\frac{8}{9}\overset{n}{\underset{i,j=q+1}{\sum}}%
\overset{q}{\underset{\text{a,b}=1}{\sum}}(R_{i\text{a}j\text{a}}%
R_{i\text{b}j\text{b}})\qquad3C$

$+\frac{2}{9}\overset{q}{\underset{\text{a}=1}{\sum}}(\varrho_{\text{aa}}%
^{M}-\varrho_{\text{aa}}^{P})(\tau^{M}-\overset{q}{\underset{\text{c}=1}{\sum
}}\varrho_{\text{cc}}^{M})+\frac{4}{9}\overset{n}{\underset{i,j=q+1}{\sum}%
}\overset{q}{\underset{\text{a}=1}{\sum}}R_{i\text{a}j\text{a}}\varrho_{ij}\ $

$\ +\frac{2}{9}\overset{q}{\underset{\text{b}=1}{\sum}}(\varrho_{\text{bb}%
}^{M}-\varrho_{\text{bb}}^{P})(\tau^{M}-\overset{q}{\underset{\text{c}%
=1}{\sum}}\varrho_{\text{cc}}^{M})+\frac{4}{9}%
\overset{n}{\underset{i,j=q+1}{\sum}}\overset{q}{\underset{\text{b}=1}{\sum}%
}R_{i\text{b}j\text{b}}\varrho_{ij}\ $

$+\frac{1}{9}(\tau^{M}-\overset{q}{\underset{\text{a=1}}{\sum}}\varrho
_{\text{aa}})(\tau^{M}-\overset{q}{\underset{\text{b=1}}{\sum}}\varrho
_{\text{bb}})+\frac{2}{9}(\left\Vert \varrho^{M}\right\Vert ^{2}%
-\overset{q}{\underset{\text{a,b}=1}{\sum}}\varrho_{\text{ab}})$

$\bigskip-\overset{n}{\underset{i,j=q+1}{\sum}}%
\overset{q}{\underset{\text{a,b}=1}{\sum}}R_{i\text{a}i\text{b}}%
R_{j\text{a}j\text{b}}\ -\frac{1}{2}\overset{n}{\underset{i,j=q+1}{\sum}%
}\overset{q}{\underset{\text{a,b}=1}{\sum}}R_{i\text{a}j\text{b}}%
^{2}-\overset{n}{\underset{i,j=q+1}{\sum}}\overset{q}{\underset{\text{a,b}%
=1}{\sum}}R_{i\text{a}j\text{b}}R_{j\text{a}i\text{b}}-\frac{1}{2}%
\overset{n}{\underset{i,j=q+1}{\sum}}\overset{q}{\underset{\text{a,b}=1}{\sum
}}R_{j\text{a}i\text{b}}^{2}$

$-\frac{1}{9}\overset{n}{\underset{i,j,p,m=q+1}{\sum}}R_{ipim}R_{jpjm}%
\ -\frac{1}{18}\overset{n}{\underset{i,j,p,m=q+1}{\sum}}R_{ipjm}^{2}-\frac
{1}{9}\overset{n}{\underset{i,j,p,m=q+1}{\sum}}R_{ipjm}R_{jpim}$

$-\frac{1}{18}\overset{n}{\underset{i,j,p,m=q+1}{\sum}}R_{jpim}^{2}-\frac
{1}{3}\overset{q}{\underset{\text{a}=1}{\sum}}%
\overset{n}{\underset{i,j,p=q+1}{\sum}}R_{i\text{a}ip}R_{j\text{a}jp}-\frac
{1}{6}\overset{q}{\underset{\text{a}=1}{\sum}}%
\overset{n}{\underset{i,j,p=q+1}{\sum}}R_{i\text{a}jp}^{2}$

$-\frac{1}{3}\overset{q}{\underset{\text{a}=1i,j,}{\sum}}%
\overset{n}{\underset{p=q+1}{\sum}}R_{i\text{a}jp}R_{j\text{a}ip}-\frac{1}%
{6}\overset{q}{\underset{\text{a}=1}{\sum}}%
\overset{n}{\underset{i,j,p=q+1}{\sum}}R_{j\text{a}ip}^{2}-\frac{1}%
{3}\overset{q}{\underset{\text{b}=1i,j,}{\sum}}%
\overset{n}{\underset{p=q+1}{\sum}}R_{i\text{b}ip}R_{j\text{b}jp}$

$-\frac{1}{6}\overset{q}{\underset{\text{b}=1}{\sum}}%
\overset{n}{\underset{i,j,p=q+1}{\sum}}R_{i\text{b}jp}^{2}-\frac{1}%
{3}\overset{q}{\underset{\text{b}=1}{\sum}}%
\overset{n}{\underset{i.j,p=q+1}{\sum}}R_{i\text{b}jp}R_{j\text{b}ip}-\frac
{1}{6}\overset{q}{\underset{\text{b}=1}{\sum}}%
\overset{n}{\underset{i,j,p=q+1}{\sum}}R_{j\text{b}ip}^{2}](y_{0})$

$-\frac{1}{48}$ $\overset{q}{\underset{\text{a,b,c=1}}{\sum}}[$
$-\overset{n}{\underset{i=q+1}{\sum}}R_{i\text{a}i\text{a}}(R_{\text{bcbc}%
}^{P}-R_{\text{bcbc}}^{M})$ $-\overset{n}{\underset{j=q+1}{\sum}}%
R_{j\text{a}j\text{a}}(R_{\text{bcbc}}^{P}-R_{\text{bcbc}}^{M})\qquad\qquad6D$

\ $+\overset{n}{\underset{i=q+1}{\sum}}R_{i\text{a}i\text{b}}(R_{\text{acbc}%
}^{P}-R_{\text{acbc}}^{M})\ -\overset{n}{\underset{i=q+1}{\sum}}%
R_{i\text{a}i\text{c}}(R_{\text{abbc}}^{P}-R_{\text{abbc}}^{M})$

$+\overset{n}{\underset{j=q+1}{\sum}}R_{j\text{a}j\text{b}}(R_{\text{acbc}%
}^{P}-R_{\text{acbc}}^{M})$\ $-\overset{n}{\underset{j=q+1}{\sum}}%
R_{j\text{a}j\text{c}}(R_{\text{abbc}}^{P}-R_{\text{abbc}}^{M})$

$+\underset{i,j=q+1}{\overset{n}{\sum}}$ $-R_{i\text{a}j\text{a}}%
(T_{\text{bb}i}T_{\text{cc}j}$ $-T_{\text{bc}i}T_{\text{bc}j})$
$-\underset{i,j=q+1}{\overset{n}{\sum}}R_{i\text{a}j\text{a}}(T_{\text{bb}%
j}T_{\text{cc}i}$ $-T_{\text{bc}j}T_{\text{bc}i})$

$+$ $\underset{i,j=q+1}{\overset{n}{\sum}}$ $-R_{j\text{a}i\text{a}%
}(T_{\text{bb}i}T_{\text{cc}j}$ $-T_{\text{bc}i}T_{\text{bc}j})$
$-\underset{i,j=q+1}{\overset{n}{\sum}}R_{j\text{a}i\text{a}}(T_{\text{bb}%
j}T_{\text{cc}i}$ $-T_{\text{bc}j}T_{\text{bc}i})$

$\ +\underset{i,j=q+1}{\overset{n}{\sum}}\ R_{i\text{a}j\text{b}}%
(T_{\text{ab}i}T_{\text{cc}j}-T_{\text{bc}i}T_{\text{ac}j}%
)\ +\underset{i,j=q+1}{\overset{n}{\sum}}\ R_{i\text{a}j\text{b}}%
(T_{\text{ab}j}T_{\text{cc}i}-T_{\text{bc}j}T_{\text{ac}i})$

$+\underset{i,j=q+1}{\overset{n}{\sum}}\ R_{j\text{a}i\text{ib}}%
(T_{\text{ab}i}T_{\text{cc}j}-T_{\text{bc}i}T_{\text{ac}j}%
)\ +\underset{i,j=q+1}{\overset{n}{\sum}}\ R_{j\text{a}i\text{b}}%
(T_{\text{ab}j}T_{\text{cc}i}-T_{\text{bc}j}T_{\text{ac}i})\qquad$

$+\underset{i,j=q+1}{\overset{n}{\sum}}-R_{i\text{a}j\text{c}}(T_{\text{ab}%
i}T_{\text{bc}j}-T_{\text{ac}i}T_{\text{bb}j}%
)-\underset{i,j=q+1}{\overset{n}{\sum}}R_{i\text{a}j\text{c}}(T_{\text{ba}%
j}T_{\text{bc}i}-T_{\text{ac}j}T_{\text{bb}i})$

$+\underset{i,j=q+1}{\overset{n}{\sum}}-R_{j\text{a}i\text{c}}(T_{\text{ba}%
i}T_{\text{bc}j}-T_{\text{ac}i}T_{\text{bb}j}%
)-\underset{i,j=q+1}{\overset{n}{\sum}}R_{j\text{a}i\text{c}}(T_{\text{ba}%
j}T_{\text{bc}i}-T_{\text{ac}j}T_{\text{bb}i})](y_{0})$

$+\frac{1}{144}\underset{p=q+1}{\overset{n}{\sum}}%
[\underset{i=q+1}{\overset{n}{\sum}}\overset{q}{\underset{\text{b,c=1}}{\sum}%
}R_{ipip}(R_{\text{bcbc}}^{P}-R_{\text{bcbc}}^{M}%
)+\underset{j=q+1}{\overset{n}{\sum}}$ $\overset{q}{\underset{\text{b,c=1}%
}{\sum}}R_{jpjp}(R_{\text{bcbc}}^{P}-R_{\text{bcbc}}^{M})](y_{0})$

$+\frac{1}{72}\underset{i,j,p=q+1}{\overset{n}{\sum}}%
\overset{q}{\underset{\text{b,c=1}}{\sum}}[R_{ipjp}(T_{\text{bb}i}%
T_{\text{cc}j}-T_{\text{bc}i}T_{\text{bc}j})+R_{ipjp}(T_{\text{bb}%
j}T_{\text{cc}i}-T_{\text{bc}j}T_{\text{bc}i})](y_{0})\qquad$

$-\frac{1}{288}\underset{i,j=q+1}{\overset{n}{\sum}}[T_{\text{aa}%
i}T_{\text{bb}j}(T_{\text{cc}i}T_{\text{dd}j}-T_{\text{cd}i}T_{\text{dc}%
j})+T_{\text{aa}i}T_{\text{bb}j}(T_{\text{cc}j}T_{\text{dd}i}-T_{\text{cd}%
j}T_{\text{dc}i})\qquad E$

$\qquad+T_{\text{aa}j}T_{\text{bb}i}(T_{\text{cc}i}T_{\text{dd}j}%
-T_{\text{cd}i}T_{\text{dc}j})+T_{\text{aa}j}T_{\text{bb}i}(T_{\text{cc}%
j}T_{\text{dd}i}-T_{\text{cd}j}T_{\text{dc}i})](y_{0})$

$\qquad+\frac{1}{288}\underset{i,j=q+1}{\overset{n}{\sum}}[T_{\text{aa}%
i}T_{\text{bc}j}(T_{\text{bc}i}T_{\text{dd}j}-T_{\text{bd}i}T_{\text{cd}%
j})+T_{\text{aa}i}T_{\text{bc}j}(T_{\text{bc}j}T_{\text{dd}i}-T_{\text{bd}%
j}T_{\text{cd}i})$

$\qquad+T_{\text{aa}j}T_{\text{bc}i}(T_{\text{bc}i}T_{\text{dd}j}%
-T_{\text{bd}i}T_{\text{cd}j})+T_{\text{aa}j}T_{\text{bc}i}(T_{\text{bc}%
j}T_{\text{dd}i}-T_{\text{bd}j}T_{\text{cd}i})](y_{0})$

$\qquad-\frac{1}{288}\underset{i,j=q+1}{\overset{n}{\sum}}[T_{\text{aa}%
i}T_{\text{bd}j}(T_{\text{bc}i}T_{\text{cd}j}-T_{\text{bd}i}T_{\text{cc}%
j})+T_{\text{aa}i}T_{\text{bd}j}(T_{\text{bc}j}T_{\text{cd}i}-T_{\text{bd}%
j}T_{\text{cc}i})$

$\qquad+T_{\text{aa}j}T_{\text{bd}i}(T_{\text{bc}i}T_{\text{cd}j}%
-T_{\text{bd}i}T_{\text{cc}j})+T_{\text{aa}j}T_{\text{bd}i}(T_{\text{bc}%
j}T_{\text{cd}i}-T_{\text{bd}j}T_{\text{cc}i})](y_{0})\qquad$

$\qquad+\frac{1}{288}\underset{i,j=q+1}{\overset{n}{\sum}}[T_{\text{ab}%
i}T_{\text{ab}j}(T_{\text{cc}i}T_{\text{dd}j}-T_{\text{cd}i}T_{\text{dc}%
j})+T_{\text{ab}i}T_{\text{ab}j}(T_{\text{cc}j}T_{\text{dd}i}-T_{\text{cd}%
j}T_{\text{dc}i})$

$\qquad+T_{\text{ab}j}T_{\text{ab}i}(T_{\text{cc}i}T_{\text{dd}j}%
-T_{\text{cd}i}T_{\text{dc}j})+T_{\text{ab}j}T_{\text{ab}i}(T_{\text{cc}%
j}T_{\text{dd}i}-T_{\text{cd}j}T_{\text{dc}i})](y_{0})$

$\qquad-\frac{1}{288}\underset{i,j=q+1}{\overset{n}{\sum}}[T_{\text{ab}%
i}T_{\text{bc}j}(T_{\text{ac}i}T_{\text{dd}j}-T_{\text{ad}i}T_{\text{cd}%
j})+T_{\text{ab}i}T_{\text{bc}j}(T_{\text{ac}j}T_{\text{dd}i}-T_{\text{ad}%
j}T_{\text{cd}i})$

$\qquad+T_{\text{ab}j}T_{\text{bc}i}(T_{\text{ac}i}T_{\text{dd}j}%
-T_{\text{ad}i}T_{\text{cd}j})+T_{\text{ab}j}T_{\text{bc}i}(T_{\text{ac}%
j}T_{\text{dd}i}-T_{\text{ad}j}T_{\text{cd}i})](y_{0})$

$\qquad+\frac{1}{288}\underset{i,j=q+1}{\overset{n}{\sum}}[T_{\text{ab}%
i}T_{\text{bd}j}(T_{\text{ac}i}T_{\text{cd}j}-T_{\text{ad}i}T_{\text{cc}%
j})+T_{\text{ab}i}T_{\text{bd}j}(T_{\text{ac}j}T_{\text{cd}i}-T_{\text{ad}%
j}T_{\text{cc}i})$

$\qquad+T_{\text{ab}i}T_{\text{bd}j}(T_{\text{ac}j}T_{\text{cd}i}%
-T_{\text{ad}j}T_{\text{cc}i})+T_{\text{ab}j}T_{\text{bd}i}(T_{\text{ac}%
j}T_{\text{cd}i}-T_{\text{ad}j}T_{\text{cc}i})](y_{0})$

$\qquad-\ \frac{1}{288}\underset{i,j=q+1}{\overset{n}{\sum}}[T_{\text{ac}%
i}T_{\text{ab}j}(T_{\text{bc}i}T_{\text{dd}j}-T_{\text{bd}i}T_{\text{dc}%
j})+T_{\text{ac}i}T_{\text{ab}j}(T_{\text{bc}j}T_{\text{dd}i}-T_{\text{bd}%
j}T_{\text{dc}i})$

$\qquad+T_{\text{ac}j}T_{\text{ab}i}(T_{\text{bc}i}T_{\text{dd}j}%
-T_{\text{bd}i}T_{\text{dc}j})+T_{\text{ac}j}T_{\text{ab}i}(T_{\text{bc}%
j}T_{\text{dd}i}-T_{\text{bd}j}T_{\text{dc}i})](y_{0})$

$\qquad+\ \frac{1}{288}\underset{i,j=q+1}{\overset{n}{\sum}}[T_{\text{ac}%
i}T_{\text{bb}j}(T_{\text{ac}i}T_{\text{dd}j}-T_{\text{ad}i}T_{\text{cd}%
j})+T_{\text{ac}i}T_{\text{bb}j}(T_{\text{ac}j}T_{\text{dd}i}-T_{\text{ad}%
j}T_{\text{cd}i})$

$\qquad+T_{\text{ac}j}T_{\text{bb}i}(T_{\text{ac}i}T_{\text{dd}j}%
-T_{\text{ad}i}T_{\text{cd}i})+T_{\text{ac}j}T_{\text{bb}i}(T_{\text{ac}%
j}T_{\text{dd}i}-T_{\text{ad}j}T_{\text{cd}i})](y_{0})$

$\qquad-\ \frac{1}{288}\underset{i,j=q+1}{\overset{n}{\sum}}[T_{\text{ac}%
i}T_{\text{bd}j}(T_{\text{ac}i}T_{\text{bd}j}-T_{\text{ad}i}T_{\text{bc}%
j})+T_{\text{ac}i}T_{\text{bd}j}(T_{\text{ac}j}T_{\text{bd}i}-T_{\text{ad}%
j}T_{\text{bc}i})$

$\qquad+T_{\text{ac}j}T_{\text{bd}i}(T_{\text{ac}i}T_{\text{bd}j}%
-T_{\text{ad}i}T_{\text{bc}j})+T_{\text{ac}j}T_{\text{bd}i}(T_{\text{ac}%
j}T_{\text{bd}i}-T_{\text{ad}j}T_{\text{bc}i})](y_{0})$

$\qquad+\frac{1}{288}\underset{i,j=q+1}{\overset{n}{\sum}}[T_{\text{ad}%
i}T_{\text{ab}j}(T_{\text{bc}i}T_{\text{cd}j}-T_{\text{bd}i}T_{\text{cc}%
j})+T_{\text{ad}i}T_{\text{ab}j}(T_{\text{bc}j}T_{\text{cd}i}-T_{\text{bd}%
j}T_{\text{cc}i})$

$\qquad+T_{\text{ad}j}T_{\text{ab}i}(T_{\text{bc}i}T_{\text{cd}j}%
-T_{\text{bd}i}T_{\text{cc}j})+T_{\text{ad}j}T_{\text{ab}i}(T_{\text{bc}%
j}T_{\text{cd}i}-T_{\text{bd}j}T_{\text{cc}i})](y_{0})$

$\qquad-\ \frac{1}{288}\underset{i,j=q+1}{\overset{n}{\sum}}[T_{\text{ad}%
i}T_{\text{bb}j}(T_{\text{ac}i}T_{\text{cd}j}-T_{\text{ad}i}T_{\text{cc}%
j})+T_{\text{ad}i}T_{\text{bb}j}(T_{\text{ac}j}T_{\text{cd}i}-T_{\text{ad}%
j}T_{\text{cc}i})$

$\qquad+T_{\text{ad}j}T_{\text{bb}i}(T_{\text{ac}i}T_{\text{cd}j}%
-T_{\text{ad}i}T_{\text{cc}j})+T_{\text{ad}j}T_{\text{bb}i}(T_{\text{ac}%
j}T_{\text{cd}i}-T_{\text{ad}j}T_{\text{cc}i})](y_{0})$

$\qquad+\ \frac{1}{288}\underset{i,j=q+1}{\overset{n}{\sum}}[T_{\text{ad}%
i}T_{\text{bc}j}(T_{\text{ac}i}T_{\text{bd}j}-T_{\text{ad}i}T_{\text{bc}%
j})+T_{\text{ad}i}T_{\text{bc}j}(T_{\text{ac}j}T_{\text{bd}i}-T_{\text{ad}%
j}T_{\text{bc}i})$

$\qquad\qquad+T_{\text{ad}j}T_{\text{bc}i}(T_{\text{ac}i}T_{\text{bd}%
j}-T_{\text{ad}i}T_{\text{bc}j})+T_{\text{ad}j}T_{\text{bc}i}(T_{\text{ac}%
j}T_{\text{bd}i}-T_{\text{ad}j}T_{\text{bc}i})](y_{0})$

$\qquad\qquad-\ \frac{1}{144}[(R_{\text{cdcd}}^{P}-R_{\text{cdcd}}%
^{M})(R_{\text{abab}}^{P}-R_{\text{abab}}^{M})](y_{0})\qquad\left(  1\right)
$

$\qquad\qquad\ +\frac{1}{144}[(R_{\text{bdcd}}^{P}-R_{\text{bdcd}}%
^{M})(R_{\text{abac}}^{P}-R_{\text{abac}}^{M})](y_{0})\qquad(2)$

$\ \qquad\qquad+\ \frac{1}{144}[(R_{\text{bcdc}}^{P}-R_{\text{bcdc}}%
^{M})(R_{\text{abad}}^{P}-R_{\text{abad}}^{M})](y_{0})\qquad(3)$

$\qquad\qquad\ -\ \frac{1}{144}[(R_{\text{adcd}}^{P}-R_{\text{adcd}}%
^{M})(R_{\text{abbc}}^{P}-R_{\text{abbc}}^{M})](y_{0})\qquad(4)\qquad$

$\ \qquad\qquad+\ \frac{1}{144}[(R_{\text{acdc}}^{P}-R_{\text{acdc}}%
^{M})(R_{\text{abdb}}^{P}-R_{\text{abdb}}^{M})](y_{0})\qquad(5)$

$\ \qquad\qquad-\ \frac{1}{576}[(R_{\text{abcd}}^{P}-R_{\text{abcd}}^{M}%
)]^{2}(y_{0})\qquad(6)$

\qquad\qquad$-\frac{1}{24}\times\frac{1}{6}<H,i>(y_{0})<H,j>(y_{0}%
)\qquad\qquad\qquad\qquad\qquad\qquad\qquad L_{3}$

$\qquad\qquad\times\lbrack2\varrho_{ij}+$ $\overset{q}{\underset{\text{a}%
=1}{4\sum}}R_{i\text{a}j\text{a}}-3\overset{q}{\underset{\text{a,b=1}}{\sum}%
}(T_{\text{aa}i}T_{\text{bb}j}-T_{\text{ab}i}T_{\text{ab}j})$

$\qquad\qquad-3\overset{q}{\underset{\text{a,b=1}}{\sum}}(T_{\text{aa}%
j}T_{\text{bb}i}-T_{\text{ab}j}T_{\text{ab}i}](y_{0})$

$\qquad\qquad-\frac{1}{24}\times\frac{3}{2}[<H,i>^{2}(y_{0})<H,j>^{2}](y_{0})$

$\qquad\qquad-\frac{1}{24}\times\frac{1}{6}<H,i>(y_{0})<H,j>(y_{0})$

$\qquad\qquad\times\lbrack2\varrho_{ij}+$ $\overset{q}{\underset{\text{a}%
=1}{4\sum}}R_{i\text{a}j\text{a}}-3\overset{q}{\underset{\text{a,b=1}}{\sum}%
}(T_{\text{aa}i}T_{\text{bb}j}-T_{\text{ab}i}T_{\text{ab}j})$

$\qquad\qquad-3\overset{q}{\underset{\text{a,b=1}}{\sum}}(T_{\text{aa}%
j}T_{\text{bb}i}-T_{\text{ab}j}T_{\text{ab}i}](y_{0})$

\qquad\qquad$\ -\frac{1}{24}\times\frac{1}{3}<H,i>(y_{0})<H,k>(y_{0}%
)R_{jijk}(y_{0})$

\qquad\qquad$-\frac{1}{24}\times\frac{3}{2}<H,i>^{2}(y_{0})<H,j>^{2}(y_{0})$

$\qquad\qquad-\frac{1}{24}\times\frac{1}{3}<H,i>^{2}(y_{0})[\varrho_{jj}+$
$\overset{q}{\underset{\text{a}=1}{2\sum}}R_{j\text{a}j\text{a}}$

$\qquad\qquad-3\overset{q}{\underset{\text{a,b=1}}{\sum}}(T_{\text{aa}%
j}T_{\text{bb}j}-T_{\text{ab}j}T_{\text{ab}j})](y_{0})$

\qquad$\qquad+\frac{1}{24}\times\frac{15}{4}[<H,i>^{2}<H,j>^{2}](y_{0}%
)<H,i>(y_{0})\qquad$\ $\frac{\partial^{3}\theta^{-\frac{1}{2}}}{\partial
\text{x}_{i}\partial\text{x}_{j}^{2}}(y_{0})$

\qquad$\qquad+\frac{1}{24}\times\frac{1}{2}<H,i>(y_{0})<H,j>$

$\times\lbrack2\varrho_{ij}+$ $\overset{q}{\underset{\text{a}=1}{4\sum}%
}R_{i\text{a}j\text{a}}-3\overset{q}{\underset{\text{a,b=1}}{\sum}%
}(T_{\text{aa}i}T_{\text{bb}j}-T_{\text{ab}i}T_{\text{ab}j}%
)-3\overset{q}{\underset{\text{a,b=1}}{\sum}}(T_{\text{aa}j}T_{\text{bb}%
i}-T_{\text{ab}j}T_{\text{ab}i}](y_{0})$

$+\frac{1}{24}\times\frac{1}{2}<H,i>(y_{0})<H,i>(y_{0})[\tau^{M}\ -3\tau
^{P}+\ \underset{\text{a}=1}{\overset{\text{q}}{\sum}}\varrho_{\text{aa}}%
^{M}+$ $\overset{q}{\underset{\text{a},\text{b}=1}{\sum}}R_{\text{abab}}^{M}$
$](y_{0})$

$+\frac{1}{24}\times\frac{1}{6}<H,i>(y_{0})[\nabla_{i}\varrho_{jj}%
-2\varrho_{ij}<H,j>+\overset{q}{\underset{\text{a}=1}{\sum}}(\nabla
_{i}R_{\text{a}j\text{a}j}-4R_{i\text{a}j\text{a}}<H,j>)$

$+4\overset{q}{\underset{\text{a,b=1}}{\sum}}R_{i\text{a}j\text{b}%
}T_{\text{ab}j}+2\overset{q}{\underset{\text{a,b,c=1}}{\sum}}(T_{\text{aa}%
i}T_{\text{bb}j}T_{\text{cc}j}-3T_{\text{aa}i}T_{\text{bc}j}T_{\text{bc}%
j}+2T_{\text{ab}i}T_{\text{bc}j}T_{\text{ca}j})](y_{0})$\qquad\qquad
\qquad\qquad\qquad\ \ 

$+\frac{1}{24}\times\frac{1}{6}<H,i>(y_{0})[\nabla_{j}\varrho_{ij}%
-2\varrho_{ij}<H,j>+\overset{q}{\underset{\text{a}=1}{\sum}}(\nabla
_{j}R_{\text{a}i\text{a}j}-4R_{j\text{a}i\text{a}}<H,j>)$

$+4\overset{q}{\underset{\text{a,b=1}}{\sum}}R_{j\text{a}i\text{b}%
}T_{\text{ab}j}+2\overset{q}{\underset{\text{a,b,c=1}}{\sum}}(T_{\text{aa}%
j}T_{\text{bb}i}T_{\text{cc}j}-3T_{\text{aa}j}T_{\text{bc}i}T_{\text{bc}%
j}+2T_{\text{ab}j}T_{\text{bc}i}T_{\text{ac}j})](y_{0})$

$+\frac{1}{24}\times\frac{1}{6}<H,i>(y_{0})[\nabla_{j}\varrho_{ij}%
-2\varrho_{jj}<H,i>+\overset{q}{\underset{\text{a}=1}{\sum}}(\nabla
_{j}R_{\text{a}i\text{a}j}-4R_{j\text{a}j\text{a}}<H,i>)$

$+4\overset{q}{\underset{\text{a,b=1}}{\sum}}R_{j\text{a}j\text{b}%
}T_{\text{ab}i}+2\overset{q}{\underset{\text{a,b,c=1}}{\sum}}(T_{\text{aa}%
j}T_{\text{bb}j}T_{\text{cc}i}-3T_{\text{aa}j}T_{\text{bc}j}T_{\text{bc}%
i}+2T_{\text{ab}j}T_{\text{bc}j}T_{\text{ac}i})](y_{0})$

$\qquad\qquad\qquad\qquad\qquad\qquad\qquad\qquad\qquad\qquad\qquad
\qquad\qquad\qquad\qquad\qquad\qquad\qquad\blacksquare$

(viii)\qquad We next compute: A$_{3213}\ =\frac{1}{12}[\frac{\partial
\theta^{\frac{1}{2}}}{\partial\text{x}_{i}}.\frac{\partial}{\partial
\text{x}_{i}}(\Delta\theta^{-\frac{1}{2}})](y_{0})\phi(y_{0})$

Since $\frac{\partial\theta^{\frac{1}{2}}}{\partial\text{x}_{i}}(y_{0}%
)=-\frac{1}{2}<H,i>(y_{0}),$ we have:

\qquad\qquad\qquad\qquad\ A$_{3213}\ =-\frac{1}{24}<H,i>(y_{0})[$%
\ $\frac{\partial}{\partial\text{x}_{i}}(\Delta\theta^{-\frac{1}{2}}%
)](y_{0})\phi(y_{0})$

We have, by the definition of the scalar Laplacian:

$\Delta\theta^{-\frac{1}{2}}=\theta^{-1}[\frac{\partial\theta}{\partial
\text{x}_{j}}(g^{jk}\frac{\partial\theta^{-\frac{1}{2}}}{\partial\text{x}_{k}%
})+\theta(\frac{\partial g^{jk}}{\partial\text{x}_{j}}\frac{\partial
\theta^{-\frac{1}{2}}}{\partial\text{x}_{k}}+g^{jk}\frac{\partial^{2}%
\theta^{-\frac{1}{2}}}{\partial\text{x}_{j}\partial\text{x}_{k}})]$

We have:

$\frac{\partial}{\partial\text{x}_{i}}(\Delta\theta^{-\frac{1}{2}}%
)(y_{0})=\frac{\partial\theta^{-1}}{\partial x_{i}}(y_{0})[\frac
{\partial\theta}{\partial\text{x}_{j}}(g^{jk}\frac{\partial\theta^{-\frac
{1}{2}}}{\partial\text{x}_{k}})+\theta(\frac{\partial g^{jk}}{\partial
\text{x}_{j}}\frac{\partial\theta^{-\frac{1}{2}}}{\partial\text{x}_{k}}%
+g^{jk}\frac{\partial^{2}\theta^{-\frac{1}{2}}}{\partial\text{x}_{j}%
\partial\text{x}_{k}})](y_{0})$

$+\theta^{-1}(y_{0})\frac{\partial}{\partial x_{i}}[\frac{\partial\theta
}{\partial\text{x}_{j}}(g^{jk}\frac{\partial\theta^{-\frac{1}{2}}}%
{\partial\text{x}_{k}})+\theta(\frac{\partial g^{jk}}{\partial\text{x}_{j}%
}\frac{\partial\theta^{-\frac{1}{2}}}{\partial\text{x}_{k}}+g^{jk}%
\frac{\partial^{2}\theta^{-\frac{1}{2}}}{\partial\text{x}_{j}\partial
\text{x}_{k}})](y_{0})$

Since $\theta^{-1}(y_{0})=1$ and $\frac{\partial\theta^{-1}}{\partial x_{i}%
}(y_{0})=-$ $\theta^{-2}(y_{0})\frac{\partial\theta}{\partial x_{i}}%
(y_{0})=-<H,i>,$ we have:

$\qquad\frac{\partial}{\partial\text{x}_{i}}(\Delta\theta^{-\frac{1}{2}%
})(y_{0})$

$\qquad=-<H,i>(y_{0})[\frac{\partial\theta}{\partial\text{x}_{j}}(g^{jk}%
\frac{\partial\theta^{-\frac{1}{2}}}{\partial\text{x}_{k}})+\theta
(\frac{\partial g^{jk}}{\partial\text{x}_{j}}\frac{\partial\theta^{-\frac
{1}{2}}}{\partial\text{x}_{k}}+g^{jk}\frac{\partial^{2}\theta^{-\frac{1}{2}}%
}{\partial\text{x}_{j}\partial\text{x}_{k}})](y_{0})$

$\qquad\qquad\qquad\qquad+\frac{\partial}{\partial x_{i}}[\frac{\partial
\theta}{\partial\text{x}_{j}}(g^{jk}\frac{\partial\theta^{-\frac{1}{2}}%
}{\partial\text{x}_{k}})+\theta(\frac{\partial g^{jk}}{\partial\text{x}_{j}%
}\frac{\partial\theta^{-\frac{1}{2}}}{\partial\text{x}_{k}}+g^{jk}%
\frac{\partial^{2}\theta^{-\frac{1}{2}}}{\partial\text{x}_{j}\partial
\text{x}_{k}})](y_{0})$

In previous computations, we saw that,

$\Delta\theta^{-\frac{1}{2}}(y_{0})=[\frac{\partial\theta}{\partial
\text{x}_{j}}(g^{jk}\frac{\partial\theta^{-\frac{1}{2}}}{\partial\text{x}_{k}%
})+\theta(\frac{\partial g^{jk}}{\partial\text{x}_{j}}\frac{\partial
\theta^{-\frac{1}{2}}}{\partial\text{x}_{k}}+g^{jk}\frac{\partial^{2}%
\theta^{-\frac{1}{2}}}{\partial\text{x}_{j}\partial\text{x}_{k}})](y_{0})$

and,

$L_{\Delta}=\frac{\partial}{\partial x_{i}}[\frac{\partial\theta}%
{\partial\text{x}_{j}}(g^{jk}\frac{\partial\theta^{-\frac{1}{2}}}%
{\partial\text{x}_{k}})+\theta(\frac{\partial g^{jk}}{\partial\text{x}_{j}%
}\frac{\partial\theta^{-\frac{1}{2}}}{\partial\text{x}_{k}}+g^{jk}%
\frac{\partial^{2}\theta^{-\frac{1}{2}}}{\partial\text{x}_{j}\partial
\text{x}_{k}})](y_{0})$Therefore,

$\frac{\partial}{\partial\text{x}_{i}}(\Delta\theta^{-\frac{1}{2}}%
)(y_{0})=-<H,i>(y_{0})\Delta\theta^{-\frac{1}{2}}(y_{0})+L_{\Delta}$

and so,

A$_{3213}\ =-\frac{1}{24}<H,i>(y_{0})[\frac{\partial}{\partial\text{x}_{i}%
}(\Delta\theta^{-\frac{1}{2}})](y_{0})$

$\qquad=$ $\frac{1}{24}<H,i>^{2}(y_{0})\Delta\theta^{-\frac{1}{2}}%
(y_{0})-\frac{1}{24}<H,i>(y_{0})L_{\Delta}$

We insert the values of $\ \Delta\theta^{-\frac{1}{2}}(y_{0})$ and $L_{\Delta
}$ above and have:

\qquad A$_{3213}\ =\frac{1}{288}<H,i>^{2}(y_{0})[3<H,j>^{2}(y_{0})$

$+2(\tau^{M}-3\tau^{P}+\overset{q}{\underset{\text{a}=1}{\sum}}\varrho
_{\text{aa}}^{M}+\overset{q}{\underset{\text{a,b}=1}{\sum}}R_{\text{abab}}%
^{M})](y_{0})+\frac{1}{12}\times\frac{1}{24}<H,i>(y_{0})<H,j>(y_{0})$

$\times\lbrack2\varrho_{ij}+$ $\overset{q}{\underset{\text{a}=1}{4\sum}%
}R_{i\text{a}j\text{a}}-3\overset{q}{\underset{\text{a,b=1}}{\sum}%
}(T_{\text{aa}i}T_{\text{bb}j}-T_{\text{ab}i}T_{\text{ab}j}%
)-3\overset{q}{\underset{\text{a,b=1}}{\sum}}(T_{\text{aa}j}T_{\text{bb}%
i}-T_{\text{ab}j}T_{\text{ab}i}](y_{0})$

$+\frac{1}{24}\times\frac{3}{4}<H,i>(y_{0})[<H,i>(y_{0})<H,j>^{2}](y_{0})$

$+\frac{1}{24}\times\frac{1}{12}<H,i>(y_{0})<H,j>(y_{0})$

$\times\lbrack2\varrho_{ij}+$ $\overset{q}{\underset{\text{a}=1}{4\sum}%
}R_{i\text{a}j\text{a}}-3\overset{q}{\underset{\text{a,b=1}}{\sum}%
}(T_{\text{aa}i}T_{\text{bb}j}-T_{\text{ab}i}T_{\text{ab}j}%
)-3\overset{q}{\underset{\text{a,b=1}}{\sum}}(T_{\text{aa}j}T_{\text{bb}%
i}-T_{\text{ab}j}T_{\text{ab}i}](y_{0})$

$\ -\frac{1}{24}\times\frac{1}{6}<H,i>(y_{0})<H,k>(y_{0})R_{ijjk}(y_{0})$

$+\frac{1}{24}\times\frac{3}{4}<H,i>^{2}(y_{0})<H,j>^{2}(y_{0})$

$+\frac{1}{24}\times\frac{1}{6}<H,i>^{2}(y_{0})[\varrho_{jj}+$
$\overset{q}{\underset{\text{a}=1}{2\sum}}R_{j\text{a}j\text{a}}%
-3\overset{q}{\underset{\text{a,b=1}}{\sum}}(T_{\text{aa}j}T_{\text{bb}%
j}-T_{\text{ab}j}T_{\text{ab}j})](y_{0})$

$-\frac{1}{24}\times\frac{15}{8}<H,i>(y_{0})[<H,i><H,j>^{2}](y_{0})\qquad
$\ $\frac{\partial^{3}\theta^{-\frac{1}{2}}}{\partial\text{x}_{i}%
\partial\text{x}_{j}^{2}}(y_{0})$

$-\frac{1}{24}\times\frac{1}{4}<H,i>(y_{0})<H,j>$

$\times\lbrack2\varrho_{ij}+$ $\overset{q}{\underset{\text{a}=1}{4\sum}%
}R_{i\text{a}j\text{a}}-3\overset{q}{\underset{\text{a,b=1}}{\sum}%
}(T_{\text{aa}i}T_{\text{bb}j}-T_{\text{ab}i}T_{\text{ab}j}%
)-3\overset{q}{\underset{\text{a,b=1}}{\sum}}(T_{\text{aa}j}T_{\text{bb}%
i}-T_{\text{ab}j}T_{\text{ab}i}](y_{0})$

$-\frac{1}{24}\times\frac{1}{4}<H,i>(y_{0})<H,i>(y_{0})[\tau^{M}\ -3\tau
^{P}+\ \underset{\text{a}=1}{\overset{\text{q}}{\sum}}\varrho_{\text{aa}}%
^{M}+$ $\overset{q}{\underset{\text{a},\text{b}=1}{\sum}}R_{\text{abab}}^{M}$
$](y_{0})$

$-\frac{1}{24}\times\frac{1}{12}<H,i>(y_{0})[\nabla_{i}\varrho_{jj}%
-2\varrho_{ij}<H,j>$

$+\overset{q}{\underset{\text{a}=1}{\sum}}(\nabla_{i}R_{\text{a}j\text{a}%
j}-4R_{i\text{a}j\text{a}}<H,j>)$ $\ \frac{\partial^{3}\theta}{\partial
\text{x}_{i}\partial\text{x}_{j}^{2}}(y_{0})$

$+4\overset{q}{\underset{\text{a,b=1}}{\sum}}R_{i\text{a}j\text{b}%
}T_{\text{ab}j}+2\overset{q}{\underset{\text{a,b,c=1}}{\sum}}(T_{\text{aa}%
i}T_{\text{bb}j}T_{\text{cc}j}-3T_{\text{aa}i}T_{\text{bc}j}T_{\text{bc}%
j}+2T_{\text{ab}i}T_{\text{bc}j}T_{\text{ca}j})](y_{0})$\qquad\qquad
\qquad\qquad\qquad\ \ 

$-\frac{1}{24}\times\frac{1}{12}<H,i>(y_{0})[\nabla_{j}\varrho_{ij}%
-2\varrho_{ij}<H,j>+\overset{q}{\underset{\text{a}=1}{\sum}}(\nabla
_{j}R_{\text{a}i\text{a}j}-4R_{j\text{a}i\text{a}}<H,j>)$

$+4\overset{q}{\underset{\text{a,b=1}}{\sum}}R_{j\text{a}i\text{b}%
}T_{\text{ab}j}+2\overset{q}{\underset{\text{a,b,c=1}}{\sum}}(T_{\text{aa}%
j}T_{\text{bb}i}T_{\text{cc}j}-3T_{\text{aa}j}T_{\text{bc}i}T_{\text{bc}%
j}+2T_{\text{ab}j}T_{\text{bc}i}T_{\text{ac}j})](y_{0})$

$-\frac{1}{24}\times\frac{1}{12}<H,i>(y_{0})[\nabla_{j}\varrho_{ij}%
-2\varrho_{jj}<H,i>+\overset{q}{\underset{\text{a}=1}{\sum}}(\nabla
_{j}R_{\text{a}i\text{a}j}-4R_{j\text{a}j\text{a}}<H,i>)$

$+4\overset{q}{\underset{\text{a,b=1}}{\sum}}R_{j\text{a}j\text{b}%
}T_{\text{ab}i}+2\overset{q}{\underset{\text{a,b,c=1}}{\sum}}(T_{\text{aa}%
j}T_{\text{bb}j}T_{\text{cc}i}-3T_{\text{aa}j}T_{\text{bc}j}T_{\text{bc}%
i}+2T_{\text{ab}j}T_{\text{bc}j}T_{\text{ac}i})](y_{0})$

\qquad\qquad\qquad\qquad\qquad\qquad\qquad\qquad\qquad\qquad\qquad\qquad
\qquad\qquad\qquad$\blacksquare$

We simplify and have:

$\left(  A_{30}\right)  $\qquad A$_{3213}\ =\frac{1}{12}[(\frac{\partial
\theta^{\frac{1}{2}}}{\partial\text{x}_{i}})\frac{\partial}{\partial
\text{x}_{i}}(\Delta\theta^{-\frac{1}{2}})](y_{0})$

$=-\frac{1}{192}<H,i>^{2}<H,j>^{2}(y_{0})$

$-\frac{1}{288}<H,i>^{2}(y_{0})[\tau^{M}\ -3\tau^{P}+\ \underset{\text{a}%
=1}{\overset{\text{q}}{\sum}}\varrho_{\text{aa}}^{M}+$
$\overset{q}{\underset{\text{a},\text{b}=1}{\sum}}R_{\text{abab}}^{M}$
$](y_{0})$

$-\frac{1}{288}<H,i>(y_{0})<H,j>(y_{0})$

$\times\lbrack2\varrho_{ij}+$ $\overset{q}{\underset{\text{a}=1}{4\sum}%
}R_{i\text{a}j\text{a}}-3\overset{q}{\underset{\text{a,b=1}}{\sum}%
}(T_{\text{aa}i}T_{\text{bb}j}-T_{\text{ab}i}T_{\text{ab}j}%
)-3\overset{q}{\underset{\text{a,b=1}}{\sum}}(T_{\text{aa}j}T_{\text{bb}%
i}-T_{\text{ab}j}T_{\text{ab}i}](y_{0})$

$+\frac{1}{144}<H,i>(y_{0})<H,k>(y_{0})R_{jijk}(y_{0})$

$+\frac{1}{144}<H,i>^{2}(y_{0})[\varrho_{jj}+$ $\overset{q}{\underset{\text{a}%
=1}{2\sum}}R_{j\text{a}j\text{a}}-3\overset{q}{\underset{\text{a,b=1}}{\sum}%
}(T_{\text{aa}j}T_{\text{bb}j}-T_{\text{ab}j}T_{\text{ab}j})](y_{0})$

$-\frac{1}{288}<H,i>(y_{0})[\nabla_{i}\varrho_{jj}-2\varrho_{ij}%
<H,j>+\overset{q}{\underset{\text{a}=1}{\sum}}(\nabla_{i}R_{\text{a}%
j\text{a}j}-4R_{i\text{a}j\text{a}}<H,j>)$

$+4\overset{q}{\underset{\text{a,b=1}}{\sum}}R_{i\text{a}j\text{b}%
}T_{\text{ab}j}+2\overset{q}{\underset{\text{a,b,c=1}}{\sum}}(T_{\text{aa}%
i}T_{\text{bb}j}T_{\text{cc}j}-3T_{\text{aa}i}T_{\text{bc}j}T_{\text{bc}%
j}+2T_{\text{ab}i}T_{\text{bc}j}T_{\text{ca}j})](y_{0})$\qquad\qquad
\qquad\qquad\qquad\ \ 

$-\frac{1}{288}<H,i>(y_{0})[\nabla_{j}\varrho_{ij}-2\varrho_{ij}%
<H,j>+\overset{q}{\underset{\text{a}=1}{\sum}}(\nabla_{j}R_{\text{a}%
i\text{a}j}-4R_{j\text{a}i\text{a}}<H,j>)$

$+4\overset{q}{\underset{\text{a,b=1}}{\sum}}R_{j\text{a}i\text{b}%
}T_{\text{ab}j}+2\overset{q}{\underset{\text{a,b,c=1}}{\sum}}(T_{\text{aa}%
j}T_{\text{bb}i}T_{\text{cc}j}-3T_{\text{aa}j}T_{\text{bc}i}T_{\text{bc}%
j}+2T_{\text{ab}j}T_{\text{bc}i}T_{\text{ac}j})](y_{0})$

$-\frac{1}{288}<H,i>(y_{0})[\nabla_{j}\varrho_{ij}-2\varrho_{jj}%
<H,i>+\overset{q}{\underset{\text{a}=1}{\sum}}(\nabla_{j}R_{\text{a}%
i\text{a}j}-4R_{j\text{a}j\text{a}}<H,i>)$

$+4\overset{q}{\underset{\text{a,b=1}}{\sum}}R_{j\text{a}j\text{b}%
}T_{\text{ab}i}+2\overset{q}{\underset{\text{a,b,c=1}}{\sum}}(T_{\text{aa}%
j}T_{\text{bb}j}T_{\text{cc}i}-3T_{\text{aa}j}T_{\text{bc}j}T_{\text{bc}%
i}+2T_{\text{ab}j}T_{\text{bc}j}T_{\text{ac}i})](y_{0})$

\qquad\qquad\qquad\qquad\qquad\qquad\qquad\qquad\qquad\qquad\qquad\qquad
\qquad\qquad\qquad$\blacksquare\qquad$

(ix) We come to the very long expression of:

\qquad A$_{321}=\frac{1}{24}\frac{\partial^{2}}{\partial\text{x}_{i}^{2}%
}(\theta^{\frac{1}{2}}\Delta\theta^{-\frac{1}{2}})(y_{0})\phi(y_{0})$

$\qquad\qquad=\frac{1}{24}\left[  (\frac{\partial^{2}\theta^{\frac{1}{2}}%
}{\partial\text{x}_{i}^{2}})(\Delta\theta^{-\frac{1}{2}})+\frac{\partial^{2}%
}{\partial\text{x}_{i}^{2}}(\Delta\theta^{-\frac{1}{2}})+2(\frac
{\partial\theta^{\frac{1}{2}}}{\partial\text{x}_{i}})\frac{\partial}%
{\partial\text{x}_{i}}(\Delta\theta^{-\frac{1}{2}})\right]  (y_{0})\phi
(y_{0})$

\qquad\qquad\ \ $=$\ \ \ A$_{3211}+$ A$_{3212}+$ A$_{3213}$where,

\bigskip$\qquad$A$_{3211}$ $=\frac{1}{24}\left[  (\frac{\partial^{2}%
\theta^{\frac{1}{2}}}{\partial\text{x}_{i}^{2}})(\Delta\theta^{-\frac{1}{2}%
})\right]  (y_{0})\phi(y_{0})$ is taken from (vi) of \textbf{Table A}$_{10}$ here.

\qquad A$_{3212}$ $=\frac{1}{24}\left[  \frac{\partial^{2}}{\partial
\text{x}_{i}^{2}}(\Delta\theta^{-\frac{1}{2}})\right]  (y_{0})\phi(y_{0})$ is
taken from $\left(  A_{29}\right)  $ or from (vii) of \textbf{Table A}$_{10}$ here.

\qquad A$_{3213}=\frac{1}{12}[(\frac{\partial\theta^{\frac{1}{2}}}%
{\partial\text{x}_{i}})\frac{\partial}{\partial\text{x}_{i}}(\Delta
\theta^{-\frac{1}{2}})](y_{0})\phi(y_{0})$ is taken from $\left(
A_{30}\right)  $ above.

We now gather all terms of A$_{321}=\frac{1}{24}\frac{\partial^{2}}%
{\partial\text{x}_{i}^{2}}(\theta^{\frac{1}{2}}\Delta\theta^{-\frac{1}{2}%
})\phi(y_{0})$ expressed in

\textbf{geometric invariants} of the Riemannian manifold M and and submanifold P.

We use the expressions in $\left(  A_{28}\right)  ,$ $\left(  A_{29}\right)  $
and $\left(  A_{30}\right)  $ to have:

$\left(  A_{31}\right)  \qquad\qquad$A$_{321}=\frac{1}{24}\frac{\partial^{2}%
}{\partial\text{x}_{i}^{2}}(\theta^{\frac{1}{2}}\Delta\theta^{-\frac{1}{2}%
})(y_{0})\phi(y_{0})$\ \ $=$\ \ \ A$_{3211}+$ A$_{3212}+$ A$_{3213}$

$\qquad=-\frac{1}{3456}[3<H,i>^{2}\ +2(\tau^{M}-3\tau^{P}%
\ +\overset{q}{\underset{\text{a=1}}{\sum}}\varrho_{\text{aa}}^{M}%
+\overset{q}{\underset{\text{a,b}=1}{\sum}}R_{\text{abab}}^{M})]^{2}%
(y_{0})\phi(y_{0})\phi(y_{0})$\ A$_{321}$\qquad A$_{3211}$

$+\frac{1}{24}[2<H,i>^{2}(y_{0})+\frac{1}{3}(\tau^{M}-3\tau^{P}%
+\overset{q}{\underset{\text{a}=1}{\sum}}\varrho_{\text{aa}}%
+\overset{q}{\underset{\text{a,b}=1}{\sum}}R_{\text{abab}})](y_{0})\phi
(y_{0})$ A$_{3212}=\frac{1}{24}(L_{1}+L_{2}+L_{3})$

$\times\lbrack\frac{1}{4}<H,j>^{2}(y_{0})+\frac{1}{6}(\tau^{M}-3\tau
^{P}+\overset{q}{\underset{\text{a}=1}{\sum}}\varrho_{\text{aa}}%
^{M}+\overset{q}{\underset{\text{a,b}=1}{\sum}}R_{\text{abab}}^{M}%
)](y_{0})\phi(y_{0})$

$-\frac{1}{96}[2\varrho_{ij}+$ $\overset{q}{\underset{\text{a}=1}{4\sum}%
}R_{i\text{a}j\text{a}}-3\overset{q}{\underset{\text{a,b=1}}{\sum}%
}(T_{\text{aa}i}T_{\text{bb}j}-T_{\text{ab}i}T_{\text{ab}j}%
)-3\overset{q}{\underset{\text{a,b=1}}{\sum}}(T_{\text{aa}j}T_{\text{bb}%
i}-T_{\text{ab}j}T_{\text{ab}i}](y_{0})$ $L_{2}$ $\ L_{21}$

$\times\lbrack<H,i><H,j>](y_{0})\phi(y_{0})$

$-\frac{1}{864}[2\varrho_{ij}+$ $\overset{q}{\underset{\text{a}=1}{4\sum}%
}R_{i\text{a}j\text{a}}-3\overset{q}{\underset{\text{a,b=1}}{\sum}%
}(T_{\text{aa}i}T_{\text{bb}j}-T_{\text{ab}i}T_{\text{ab}j}%
)-3\overset{q}{\underset{\text{a,b=1}}{\sum}}(T_{\text{aa}j}T_{\text{bb}%
i}-T_{\text{ab}j}T_{\text{ab}i}]^{2}(y_{0})\phi(y_{0})$

$-\frac{1}{288}[<H,j>](y_{0})\times\lbrack\{\nabla_{i}\varrho_{ij}%
-2\varrho_{ij}<H,i>+\overset{q}{\underset{\text{a}=1}{\sum}}(\nabla
_{i}R_{\text{a}i\text{a}j}-4R_{i\text{a}j\text{a}}<H,i>)\qquad L_{212}$

$\qquad+4\overset{q}{\underset{\text{a,b=1}}{\sum}}R_{i\text{a}j\text{b}%
}T_{\text{ab}i}+2\overset{q}{\underset{\text{a,b,c=1}}{\sum}}(T_{\text{aa}%
i}T_{\text{bb}j}T_{\text{cc}i}-3T_{\text{aa}i}T_{\text{bc}j}T_{\text{bc}%
i}+2T_{\text{ab}i}T_{\text{bc}j}T_{\text{ac}i})](y_{0})\phi(y_{0})$%
\qquad\qquad\qquad\qquad\qquad\ \ 

$\qquad-\frac{1}{24}\times\frac{1}{12}[<H,j>](y_{0})\times\lbrack\nabla
_{j}\varrho_{ii}-2\varrho_{ij}<H,i>+\overset{q}{\underset{\text{a}=1}{\sum}%
}(\nabla_{j}R_{\text{a}i\text{a}i}-4R_{i\text{a}j\text{a}}<H,i>)$

$\qquad+4\overset{q}{\underset{\text{a,b=1}}{\sum}}R_{j\text{a}i\text{b}%
}T_{\text{ab}i}+2\overset{q}{\underset{\text{a,b,c=1}}{\sum}}(T_{\text{aa}%
j}T_{\text{bb}i}T_{\text{cc}i}-3T_{\text{aa}j}T_{\text{bc}i}T_{\text{bc}%
i}+2T_{\text{ab}j}T_{\text{bc}i}T_{\text{ac}i})](y_{0})\phi(y_{0})$

$\qquad-\frac{1}{24}\times\frac{1}{12}[<H,j>](y_{0})\times\lbrack\nabla
_{i}\varrho_{ij}-2\varrho_{ii}<H,j>+\overset{q}{\underset{\text{a}=1}{\sum}%
}(\nabla_{i}R_{\text{a}i\text{a}j}-4R_{i\text{a}i\text{a}}<H,j>)$

$\qquad+4\overset{q}{\underset{\text{a,b=1}}{\sum}}R_{i\text{a}i\text{b}%
}T_{\text{ab}j}+2\overset{q}{\underset{\text{a,b,c}=1}{\sum}}(T_{\text{aa}%
i}T_{\text{bb}i}T_{\text{cc}j}-3T_{\text{aa}i}T_{\text{bc}i}T_{\text{bc}%
j}+2T_{\text{ab}i}T_{\text{bc}i}T_{\text{ac}j})](y_{0})\phi(y_{0})$

$\qquad-\frac{1}{3}[<H,j><H,k>](y_{0})R_{ijik}(y_{0})-$ $\frac{1}{24}%
\times\frac{15}{8}[<H,i>^{2}<H,j>^{2}](y_{0})\phi(y_{0})\qquad L_{213}$

$-\frac{1}{96}<H,i><H,j>$

$\times\lbrack2\varrho_{ij}+\overset{q}{\underset{\text{a}=1}{4\sum}%
}R_{i\text{a}j\text{a}}-3\overset{q}{\underset{\text{a,b=1}}{\sum}%
}(T_{\text{aa}i}T_{\text{bb}j}-T_{\text{ab}i}T_{\text{ab}j}%
)-3\overset{q}{\underset{\text{a,b=1}}{\sum}}(T_{\text{aa}j}T_{\text{bb}%
i}-T_{\text{ab}j}T_{\text{ab}i}](y_{0})\phi(y_{0})$

$-\frac{1}{96}<H,j>^{2}[\tau^{M}\ -3\tau^{P}+\ \underset{\text{a}%
=1}{\overset{\text{q}}{\sum}}\varrho_{\text{aa}}^{M}+$
$\overset{q}{\underset{\text{a},\text{b}=1}{\sum}}R_{\text{abab}}^{M}$
$](y_{0})\phi(y_{0})$

$+\frac{1}{288}<H,j>[\nabla_{i}\varrho_{ij}-2\varrho_{ij}%
<H,i>+\overset{q}{\underset{\text{a}=1}{\sum}}(\nabla_{i}R_{\text{a}%
i\text{a}j}-4R_{i\text{a}j\text{a}}<H,i>)+4\overset{q}{\underset{\text{a,b=1}%
}{\sum}}R_{i\text{a}j\text{b}}T_{\text{ab}i}$

$+2\overset{q}{\underset{\text{a,b,c=1}}{\sum}}(T_{\text{aa}i}T_{\text{bb}%
j}T_{\text{cc}i}-3T_{\text{aa}i}T_{\text{bc}j}T_{\text{bc}i}+2T_{\text{ab}%
i}T_{\text{bc}j}T_{\text{ac}i})](y_{0})\phi(y_{0})$\qquad\qquad\qquad
\qquad\qquad\ \ 

$+\frac{1}{12}<H,j>[\nabla_{j}\varrho_{ii}-2\varrho_{ij}%
<H,i>+\overset{q}{\underset{\text{a}=1}{\sum}}(\nabla_{j}R_{\text{a}%
i\text{a}i}-4R_{i\text{a}j\text{a}}<H,i>)+4\overset{q}{\underset{\text{a,b=1}%
}{\sum}}R_{j\text{a}i\text{b}}T_{\text{ab}i}$

$+2\overset{q}{\underset{\text{a,b,c=1}}{\sum}}(T_{\text{aa}j}T_{\text{bb}%
i}T_{\text{cc}i}-3T_{\text{aa}j}T_{\text{bc}i}T_{\text{bc}i}+2T_{\text{ab}%
j}T_{\text{bc}i}T_{\text{ac}i})](y_{0})\phi(y_{0})$

$+\frac{1}{288}<H,j>[\nabla_{i}\varrho_{ij}-2\varrho_{ii}%
<H,j>+\overset{q}{\underset{\text{a}=1}{\sum}}(\nabla_{i}R_{\text{a}%
i\text{a}j}-4R_{i\text{a}i\text{a}}<H,j>)$

$+4\overset{q}{\underset{\text{a,b=1}}{\sum}}R_{i\text{a}i\text{b}%
}T_{\text{ab}j}+2\overset{q}{\underset{\text{a,b,c}=1}{\sum}}(T_{\text{aa}%
i}T_{\text{bb}i}T_{\text{cc}j}-3T_{\text{aa}i}T_{\text{bc}i}T_{\text{bc}%
j}+2T_{\text{ab}i}T_{\text{bc}i}T_{\text{ac}j})](y_{0})\phi(y_{0})$

$-\frac{1}{144}R_{jijk}(y_{0})$ \ $[<H,i><H,k>](y_{0})\phi(y_{0})\qquad\qquad
L_{22}$

$-\frac{1}{432}R_{jijk}(y_{0})[2\varrho_{ik}+$ $\overset{q}{\underset{\text{a}%
=1}{4\sum}}R_{i\text{a}k\text{a}}-3\overset{q}{\underset{\text{a,b=1}}{\sum}%
}(T_{\text{aa}i}T_{\text{bb}k}-T_{\text{ab}i}T_{\text{ab}k})$

$-3\overset{q}{\underset{\text{a,b=1}}{\sum}}(T_{\text{aa}k}T_{\text{bb}%
i}-T_{\text{ab}k}T_{\text{ab}i}](y_{0})\phi(y_{0})$

$+\frac{1}{144}<H,k>(y_{0})[\nabla_{j}$R$_{ijik}(y_{0})-\nabla_{i}$%
R$_{jijk}](y_{0})\phi(y_{0})$

$-\frac{5}{32}<H,i>^{2}(y_{0})<H,j>^{2}(y_{0})\phi(y_{0})\qquad\qquad
L_{23}\qquad L_{231}$

$-\frac{1}{48}<H,i>(y_{0})<H,j>(y_{0})$

$\times\lbrack2\varrho_{ij}+$ $\overset{q}{\underset{\text{a}=1}{4\sum}%
}R_{i\text{a}j\text{a}}-3\overset{q}{\underset{\text{a,b=1}}{\sum}%
}(T_{\text{aa}i}T_{\text{bb}j}-T_{\text{ab}i}T_{\text{ab}j}%
)-3\overset{q}{\underset{\text{a,b=1}}{\sum}}(T_{\text{aa}j}T_{\text{bb}%
i}-T_{\text{ab}j}T_{\text{ab}i}](y_{0})\phi(y_{0})$

$-\frac{1}{48}<H,i>^{2}(y_{0})[\varrho_{jj}+$ $\overset{q}{\underset{\text{a}%
=1}{2\sum}}R_{j\text{a}j\text{a}}-3\overset{q}{\underset{\text{a,b=1}}{\sum}%
}(T_{\text{aa}j}T_{\text{bb}j}-T_{\text{ab}j}T_{\text{ab}j})](y_{0})\phi
(y_{0})$

$-\frac{1}{144}<H,i>(y_{0})[\nabla_{i}\varrho_{jj}-2\varrho_{ij}%
<H,j>+\overset{q}{\underset{\text{a}=1}{\sum}}(\nabla_{i}R_{\text{a}%
j\text{a}j}-4R_{i\text{a}j\text{a}}<H,j>)$

$+4\overset{q}{\underset{\text{a,b=1}}{\sum}}R_{i\text{a}j\text{b}%
}T_{\text{ab}j}+2\overset{q}{\underset{\text{a,b,c=1}}{\sum}}(T_{\text{aa}%
i}T_{\text{bb}j}T_{\text{cc}j}-3T_{\text{aa}i}T_{\text{bc}j}T_{\text{bc}%
j}+2T_{\text{ab}i}T_{\text{bc}j}T_{\text{ca}j})](y_{0})\phi(y_{0})$%
\qquad\qquad\qquad\qquad\qquad\ \ 

$-\frac{1}{144}<H,i>(y_{0})[\nabla_{j}\varrho_{ij}-2\varrho_{ij}%
<H,j>+\overset{q}{\underset{\text{a}=1}{\sum}}(\nabla_{j}R_{\text{a}%
i\text{a}j}-4R_{j\text{a}i\text{a}}<H,j>)$

$+4\overset{q}{\underset{\text{a,b=1}}{\sum}}R_{j\text{a}i\text{b}%
}T_{\text{ab}j}+2\overset{q}{\underset{\text{a,b,c=1}}{\sum}}(T_{\text{aa}%
j}T_{\text{bb}i}T_{\text{cc}j}-3T_{\text{aa}j}T_{\text{bc}i}T_{\text{bc}%
j}+2T_{\text{ab}j}T_{\text{bc}i}T_{\text{ac}j})](y_{0})$

$-\frac{1}{144}<H,i>(y_{0})[\nabla_{j}\varrho_{ij}-2\varrho_{jj}%
<H,i>+\overset{q}{\underset{\text{a}=1}{\sum}}(\nabla_{j}R_{\text{a}%
i\text{a}j}-4R_{j\text{a}j\text{a}}<H,i>)$

$+4\overset{q}{\underset{\text{a,b=1}}{\sum}}R_{j\text{a}j\text{b}%
}T_{\text{ab}i}+2\overset{q}{\underset{\text{a,b,c=1}}{\sum}}(T_{\text{aa}%
j}T_{\text{bb}j}T_{\text{cc}i}-3T_{\text{aa}j}T_{\text{bc}j}T_{\text{bc}%
i}+2T_{\text{ab}j}T_{\text{bc}j}T_{\text{ac}i})](y_{0})\phi(y_{0})$

$-\frac{1}{96}<H,j>^{2}(y_{0})[\varrho_{ii}+$ $\overset{q}{\underset{\text{a}%
=1}{2\sum}}R_{i\text{a}i\text{a}}-3\overset{q}{\underset{\text{a,b=1}}{\sum}%
}(T_{\text{aa}i}T_{\text{bb}i}-T_{\text{ab}i}T_{\text{ab}i})](y_{0})\qquad
L_{232}$

$-\frac{1}{432}[\varrho_{ii}+$ $\overset{q}{\underset{\text{a}=1}{2\sum}%
}R_{i\text{a}i\text{a}}-3\overset{q}{\underset{\text{a,b=1}}{\sum}%
}(T_{\text{aa}i}T_{\text{bb}i}-T_{\text{ab}i}T_{\text{ab}i})](y_{0})\phi
(y_{0})$

$\times\lbrack\varrho_{jj}+$ $\overset{q}{\underset{\text{a}=1}{2\sum}%
}R_{j\text{a}j\text{a}}-3\overset{q}{\underset{\text{a,b=1}}{\sum}%
}(T_{\text{aa}j}T_{\text{bb}j}-T_{\text{ab}j}T_{\text{ab}j})](y_{0})\phi
(y_{0})$

$+\frac{1}{48}R_{ijik}(y_{0})$ \ $[<H,j><H,k>](y_{0})\phi(y_{0})\qquad
\qquad\qquad\qquad$\ $L_{233}$

$+\frac{1}{432}R_{ijik}(y_{0})$

$\times\lbrack2\varrho_{jk}+$ $\overset{q}{\underset{\text{a}=1}{4\sum}%
}R_{j\text{a}k\text{a}}-3\overset{q}{\underset{\text{a,b=1}}{\sum}%
}(T_{\text{aa}j}T_{\text{bb}k}-T_{\text{ab}j}T_{\text{ab}k}%
)-3\overset{q}{\underset{\text{a,b=1}}{\sum}}(T_{\text{aa}k}T_{\text{bb}%
j}-T_{\text{ab}k}T_{\text{ab}j}](y_{0})\phi(y_{0})$

$+\frac{35}{128}\overset{n}{\underset{i,j=q+1}{\sum}}<H,i>^{2}(y_{0}%
)<H,j>^{2}(y_{0})\qquad\qquad\ \frac{1}{24}\frac{\partial^{4}\theta^{-\frac
{1}{2}}}{\partial x_{i}^{2}\partial x_{j}^{2}}(y_{0})$

$+\frac{5}{192}\overset{n}{\underset{j=q+1}{\sum}}<H,j>^{2}(y_{0})[\tau
^{M}\ -3\tau^{P}+\ \underset{\text{a}=1}{\overset{\text{q}}{\sum}}%
\varrho_{\text{aa}}^{M}+\overset{q}{\underset{\text{a},\text{b}=1}{\sum}%
}R_{\text{abab}}^{M}](y_{0})\qquad\ \ \ \ \ \ \ \ $

$+\frac{5}{192}\overset{n}{\underset{i=q+1}{\sum}}<H,i>^{2}(y_{0})[\tau
^{M}\ -3\tau^{P}+\ \underset{\text{a}=1}{\overset{\text{q}}{\sum}}%
\varrho_{\text{aa}}^{M}+\overset{q}{\underset{\text{a},\text{b}=1}{\sum}%
}R_{\text{abab}}^{M}](y_{0})\qquad\qquad$

$+\frac{5}{192}\overset{n}{\underset{i,j=q+1}{\sum}}[<H,i><H,j>](y_{0}%
)\qquad\qquad\qquad\qquad\qquad\qquad\qquad\qquad$

$\times\lbrack2\varrho_{ij}+4\overset{q}{\underset{\text{a}=1}{\sum}%
}R_{i\text{a}j\text{a}}-3\overset{q}{\underset{\text{a,b=1}}{\sum}%
}(T_{\text{aa}i}T_{\text{bb}j}-T_{\text{ab}i}T_{\text{ab}j}%
)-3\overset{q}{\underset{\text{a,b=1}}{\sum}}(T_{\text{aa}j}T_{\text{bb}%
i}-T_{\text{ab}j}T_{\text{ab}i})](y_{0})$

$+\frac{1}{96}\overset{n}{\underset{i,j=q+1}{\sum}}<H,j>(y_{0})[\{\nabla
_{i}\varrho_{ij}-2\varrho_{ij}<H,i>+\overset{q}{\underset{\text{a}=1}{\sum}%
}(\nabla_{i}R_{\text{a}i\text{a}j}-4R_{i\text{a}j\text{a}}<H,i>)\qquad$

$+4\overset{q}{\underset{\text{a,b=1}}{\sum}}R_{i\text{a}j\text{b}%
}T_{\text{ab}i}+2\overset{q}{\underset{\text{a,b,c=1}}{\sum}}(T_{\text{aa}%
i}T_{\text{bb}j}T_{\text{cc}i}-T_{\text{aa}i}T_{\text{bc}j}T_{\text{bc}%
i}-2T_{\text{bc}j}(T_{\text{aa}i}T_{\text{bc}i}-T_{\text{ab}i}T_{\text{ac}%
i}))\}$\qquad\qquad\qquad\ \ 

$+\{\nabla_{j}\varrho_{ii}-2\varrho_{ij}<H,i>+\overset{q}{\underset{\text{a}%
=1}{\sum}}(\nabla_{j}R_{\text{a}i\text{a}i}-4R_{i\text{a}j\text{a}}<H,i>)$

$+4\overset{q}{\underset{\text{a,b=1}}{\sum}}R_{j\text{a}i\text{b}%
}T_{\text{ab}i}+2\overset{q}{\underset{\text{a,b,c=1}}{\sum}}(T_{\text{aa}%
j}(T_{\text{bb}i}T_{\text{cc}i}-T_{\text{bc}i}T_{\text{bc}i})-2T_{\text{aa}%
j}T_{\text{bc}i}T_{\text{bc}i}+2T_{\text{ab}j}T_{\text{bc}i}T_{\text{ac}%
i})\}\qquad$

$+\{\nabla_{i}\varrho_{ij}-2\varrho_{ii}<H,j>+\overset{q}{\underset{\text{a}%
=1}{\sum}}(\nabla_{i}R_{\text{a}i\text{a}j}-4R_{i\text{a}i\text{a}%
}<H,j>)+4\overset{q}{\underset{\text{a,b=1}}{\sum}}R_{i\text{a}i\text{b}%
}T_{\text{ab}j}$

$+2\overset{q}{\underset{\text{a,b,c}=1}{\sum}}(T_{\text{aa}i}T_{\text{bb}%
i}T_{\text{cc}j}-3T_{\text{aa}i}T_{\text{bc}i}T_{\text{bc}j}+2T_{\text{ab}%
i}T_{\text{bc}i}T_{\text{ac}j})\}](y_{0})$

$+\frac{1}{96}\overset{n}{\underset{i,j=q+1}{\sum}}<H,i>(y_{0})[\{\nabla
_{i}\varrho_{jj}-2\varrho_{ij}<H,j>+\overset{q}{\underset{\text{a}=1}{\sum}%
}(\nabla_{i}R_{\text{a}j\text{a}j}-4R_{i\text{a}j\text{a}}<H,j>)\qquad$

$+4\overset{q}{\underset{\text{a,b=1}}{\sum}}R_{i\text{a}j\text{b}%
}T_{\text{ab}j}+2\overset{q}{\underset{\text{a,b,c=1}}{\sum}}T_{\text{aa}%
i}(T_{\text{bb}j}T_{\text{cc}j}-T_{\text{bc}j}T_{\text{bc}j})-2T_{\text{aa}%
i}T_{\text{bc}j}T_{\text{bc}j}+2T_{\text{ab}i}T_{\text{bc}j}T_{\text{ac}%
j})\}(y_{0})\qquad$\qquad\qquad\qquad\qquad\qquad\ \ 

$+\{\nabla_{j}\varrho_{ij}-2\varrho_{ij}<H,j>+\overset{q}{\underset{\text{a}%
=1}{\sum}}(\nabla_{j}R_{\text{a}i\text{a}j}-4R_{j\text{a}i\text{a}}<H,j>)$

$+4\overset{q}{\underset{\text{a,b=1}}{\sum}}R_{j\text{a}i\text{b}%
}T_{\text{ab}j}+2\overset{q}{\underset{\text{a,b,c=1}}{\sum}}(T_{\text{aa}%
j}T_{\text{bb}i}T_{\text{cc}j}-T_{\text{ab}j}T_{\text{bc}i}T_{\text{ac}%
j}-2T_{\text{bc}i}(T_{\text{aa}j}T_{\text{bc}j}-T_{\text{ab}j}T_{\text{ac}%
j}))\}(y_{0})$

$+\{\nabla_{j}\varrho_{ij}-2\varrho_{jj}<H,i>+\overset{q}{\underset{\text{a}%
=1}{\sum}}(\nabla_{j}R_{\text{a}i\text{a}j}-4R_{j\text{a}j\text{a}%
}<H,i>)+4\overset{q}{\underset{\text{a,b=1}}{\sum}}R_{j\text{a}j\text{b}%
}T_{\text{ab}i}$

$+2\overset{q}{\underset{\text{a,b,c=1}}{\sum}}(T_{\text{aa}j}T_{\text{bb}%
j}T_{\text{cc}i}-3T_{\text{aa}j}T_{\text{bc}j}T_{\text{bc}i}+2T_{\text{ab}%
j}T_{\text{bc}j}T_{\text{ac}i})\}](y_{0})$

$+\frac{1}{576}\overset{n}{\underset{i,j=q+1}{\sum}}[2\varrho_{ij}%
+4\overset{q}{\underset{\text{a}=1}{\sum}}R_{i\text{a}j\text{a}}%
-3\overset{q}{\underset{\text{a,b=1}}{\sum}}(T_{\text{aa}i}T_{\text{bb}%
j}-T_{\text{ab}i}T_{\text{ab}j})$

$-3\overset{q}{\underset{\text{a,b=1}}{\sum}}(T_{\text{aa}j}T_{\text{bb}%
i}-T_{\text{ab}j}T_{\text{ab}i})]^{2}(y_{0})$

$+\frac{1}{288}[\tau^{M}\ -3\tau^{P}+\ \underset{\text{a}=1}{\overset{\text{q}%
}{\sum}}\varrho_{\text{aa}}^{M}+\overset{q}{\underset{\text{a},\text{b}%
=1}{\sum}}R_{\text{abab}}^{M}]^{2}(y_{0})$

$-\ \frac{1}{288}\overset{n}{\underset{i,j=q+1}{\sum}}[$
$\overset{q}{\underset{\text{a=1}}{\sum}}\{-(\nabla_{ii}^{2}R_{j\text{a}%
j\text{a}}+\nabla_{jj}^{2}R_{i\text{a}i\text{a}}+4\nabla_{ij}^{2}%
R_{i\text{a}j\text{a}}+2R_{ij}R_{i\text{a}j\text{a}})$

$+\overset{n}{\underset{p=q+1}{\sum}}\overset{q}{\underset{\text{a=1}}{\sum}%
}(R_{\text{a}iip}R_{\text{a}jjp}+R_{\text{a}jjp}R_{\text{a}iip}+R_{\text{a}%
ijp}R_{\text{a}ijp}+R_{\text{a}ijp}R_{\text{a}jip}+R_{\text{a}jip}%
R_{\text{a}ijp}+R_{\text{a}jip}R_{\text{a}jip})$

$+2\overset{q}{\underset{\text{a,b=1}}{\sum}}\nabla_{i}(R)_{\text{a}%
i\text{b}j}T_{\text{ab}j}+2\overset{q}{\underset{\text{a,b=1}}{\sum}}%
\nabla_{j}(R)_{\text{a}j\text{b}i}T_{\text{ab}i}%
+2\overset{q}{\underset{\text{a,b=1}}{\sum}}\nabla_{i}(R)_{\text{a}j\text{b}%
i}T_{\text{ab}j}$

$+2\overset{q}{\underset{\text{a,b=1}}{\sum}}\nabla_{i}(R)_{\text{a}%
j\text{b}j}T_{\text{ab}i}+2\overset{q}{\underset{\text{a,b=1}}{\sum}}%
\nabla_{j}(R)_{\text{a}i\text{b}i}T_{\text{ab}j}%
+2\overset{q}{\underset{\text{a,b=1}}{\sum}}\nabla_{j}(R)_{\text{a}i\text{b}%
j}T_{\text{ab}i}$

$+\overset{n}{\underset{p=q+1}{\sum}}(-\frac{3}{5}\nabla_{ii}^{2}%
(R)_{jpjp}+\overset{n}{\underset{p=q+1}{\sum}}(-\frac{3}{5}\nabla_{jj}%
^{2}(R)_{ipip}$

$+\overset{n}{\underset{p=q+1}{\sum}}(-\frac{3}{5}\nabla_{ij}^{2}%
(R)_{ipjp}+\overset{n}{\underset{p=q+1}{\sum}}(-\frac{3}{5}\nabla_{ij}%
^{2}(R)_{jpip}+\overset{n}{\underset{p=q+1}{\sum}}(-\frac{3}{5}\nabla_{ji}%
^{2}(R)_{ipjp}$

$+\overset{n}{\underset{p=q+1}{\sum}}(-\frac{3}{5}\nabla_{ji}^{2}%
(R)_{jpip}+\frac{1}{5}\overset{n}{\underset{m,p=q+1}{%
{\textstyle\sum}
}}R_{ipim}R_{jpjm}+\frac{1}{5}\overset{n}{\underset{m,p=q+1}{%
{\textstyle\sum}
}}R_{jpjm}R_{ipim}$

$+\frac{1}{5}\overset{n}{\underset{m,p=q+1}{%
{\textstyle\sum}
}}R_{ipjm}R_{ipjm}+\frac{1}{5}\overset{n}{\underset{m,p=q+1}{%
{\textstyle\sum}
}}R_{ipjm}R_{jpim}$

$+\frac{1}{5}\overset{n}{\underset{m,p=q+1}{%
{\textstyle\sum}
}}R_{jpim}R_{ipjm}+\frac{1}{5}\overset{n}{\underset{m,p=q+1}{%
{\textstyle\sum}
}}R_{jpim}R_{jpim}\}(y_{0})$

$+4\overset{q}{\underset{\text{a,b=1}}{\sum}}\{(\nabla_{i}(R)_{i\text{a}%
j\text{a}}-\overset{q}{\underset{\text{c=1}}{%
{\textstyle\sum}
}}R_{\text{a}i\text{c}i}T_{\text{ac}j})$ $T_{\text{bb}j}+4(\nabla
_{j}(R)_{j\text{a}i\text{a}}-\overset{q}{\underset{\text{c=1}}{%
{\textstyle\sum}
}}R_{\text{a}j\text{c}j}T_{\text{ac}i})$ $T_{\text{bb}i}$

$+4(\nabla_{i}(R)_{j\text{a}i\text{a}}-\overset{q}{\underset{\text{c=1}}{%
{\textstyle\sum}
}}R_{\text{a}i\text{c}j}T_{\text{ac}i})$ $T_{\text{bb}j}$

$+4(\nabla_{i}(R)_{j\text{a}j\text{a}}-\overset{q}{\underset{\text{c=1}}{%
{\textstyle\sum}
}}R_{\text{a}i\text{c}j}T_{\text{ac}j})$ $T_{\text{bb}i}+4(\nabla
_{j}(R)_{i\text{a}i\text{a}}-\overset{q}{\underset{\text{c=1}}{%
{\textstyle\sum}
}}R_{\text{a}j\text{c}i}T_{\text{ac}i})$ $T_{\text{bb}j}$

$+4(\nabla_{j}(R)_{i\text{a}j\text{a}}-\overset{q}{\underset{\text{c=1}}{%
{\textstyle\sum}
}}R_{\text{a}j\text{c}i}T_{\text{ac}j})$ $T_{\text{bb}i}$

$-4\overset{q}{\underset{\text{a,b=1}}{\sum}}(\nabla_{i}(R)_{i\text{a}%
j\text{b}}-\overset{q}{\underset{\text{c=1}}{%
{\textstyle\sum}
}}R_{\text{b}r\text{c}s}T_{\text{ac}t})T_{\text{ab}j}%
-4\overset{q}{\underset{\text{a,b=1}}{\sum}}(\nabla_{j}(R)_{j\text{a}%
i\text{b}}-\overset{q}{\underset{\text{c=1}}{%
{\textstyle\sum}
}}R_{\text{b}j\text{c}j}T_{\text{ac}i})T_{\text{ab}i}$

$-4\overset{q}{\underset{\text{a,b=1}}{\sum}}(\nabla_{i}(R)_{j\text{a}%
i\text{b}}-\overset{q}{\underset{\text{c=1}}{%
{\textstyle\sum}
}}R_{\text{b}i\text{c}j}T_{\text{ac}i})T_{\text{ab}j}%
-4\overset{q}{\underset{\text{a,b=1}}{\sum}}(\nabla_{i}(R)_{j\text{a}%
j\text{b}}-\overset{q}{\underset{\text{c=1}}{%
{\textstyle\sum}
}}R_{\text{b}i\text{c}j}T_{\text{ac}j})T_{\text{ab}i}$

$-4\overset{q}{\underset{\text{a,b=1}}{\sum}}(\nabla_{j}(R)_{i\text{a}%
i\text{b}}-\overset{q}{\underset{\text{c=1}}{%
{\textstyle\sum}
}}R_{\text{b}j\text{c}i}T_{\text{ac}i})T_{\text{ab}j}%
-4\overset{q}{\underset{\text{a,b=1}}{\sum}}(\nabla_{j}(R)_{i\text{a}%
j\text{b}}-\overset{q}{\underset{\text{c=1}}{%
{\textstyle\sum}
}}R_{\text{b}j\text{c}i}T_{\text{ac}j})T_{\text{ab}i}\}](y_{0})$

$-\frac{1}{48}$ $[\frac{4}{9}\overset{q}{\underset{\text{a,b=1}}{\sum}%
}(\varrho_{\text{aa}}-\overset{q}{\underset{\text{c}=1}{\sum}}R_{\text{acac}%
})(\varrho_{\text{bb}}-\overset{q}{\underset{\text{d}=1}{\sum}}R_{\text{bdbd}%
})+\frac{8}{9}\overset{n}{\underset{i,j=q+1}{\sum}}%
\overset{q}{\underset{\text{a,b}=1}{\sum}}(R_{i\text{a}j\text{a}}%
R_{i\text{b}j\text{b}})$

$+\frac{2}{9}\overset{q}{\underset{\text{a}=1}{\sum}}(\varrho_{\text{aa}}%
^{M}-\varrho_{\text{aa}}^{P})(\tau^{M}-\overset{q}{\underset{\text{c}=1}{\sum
}}\varrho_{\text{cc}}^{M})+\frac{4}{9}\overset{n}{\underset{i,j=q+1}{\sum}%
}\overset{q}{\underset{\text{a}=1}{\sum}}R_{i\text{a}j\text{a}}\varrho_{ij}\ $

$\ +\frac{2}{9}\overset{q}{\underset{\text{b}=1}{\sum}}(\varrho_{\text{bb}%
}^{M}-\varrho_{\text{bb}}^{P})(\tau^{M}-\overset{q}{\underset{\text{c}%
=1}{\sum}}\varrho_{\text{cc}}^{M})+\frac{4}{9}%
\overset{n}{\underset{i,j=q+1}{\sum}}\overset{q}{\underset{\text{b}=1}{\sum}%
}R_{i\text{b}j\text{b}}\varrho_{ij}\ $

$+\frac{1}{9}(\tau^{M}-\overset{q}{\underset{\text{a=1}}{\sum}}\varrho
_{\text{aa}})(\tau^{M}-\overset{q}{\underset{\text{b=1}}{\sum}}\varrho
_{\text{bb}})+\frac{2}{9}(\left\Vert \varrho^{M}\right\Vert ^{2}%
-\overset{q}{\underset{\text{a,b}=1}{\sum}}\varrho_{\text{ab}})$

$-\overset{n}{\underset{i,j=q+1}{\sum}}\overset{q}{\underset{\text{a,b}%
=1}{\sum}}R_{i\text{a}i\text{b}}R_{j\text{a}j\text{b}}\ -\frac{1}%
{2}\overset{n}{\underset{i,j=q+1}{\sum}}\overset{q}{\underset{\text{a,b}%
=1}{\sum}}R_{i\text{a}j\text{b}}^{2}-\overset{n}{\underset{i,j=q+1}{\sum}%
}\overset{q}{\underset{\text{a,b}=1}{\sum}}R_{i\text{a}j\text{b}}%
R_{j\text{a}i\text{b}}-\frac{1}{2}\overset{n}{\underset{i,j=q+1}{\sum}%
}\overset{q}{\underset{\text{a,b}=1}{\sum}}R_{j\text{a}i\text{b}}^{2}$

$-\frac{1}{9}\overset{n}{\underset{i,j,p,m=q+1}{\sum}}R_{ipim}R_{jpjm}%
\ -\frac{1}{18}\overset{n}{\underset{i,j,p,m=q+1}{\sum}}R_{ipjm}^{2}-\frac
{1}{9}\overset{n}{\underset{i,j,p,m=q+1}{\sum}}R_{ipjm}R_{jpim}$

$-\frac{1}{18}\overset{n}{\underset{i,j,p,m=q+1}{\sum}}R_{jpim}^{2}$

$-\frac{1}{3}\overset{q}{\underset{\text{a}=1}{\sum}}%
\overset{n}{\underset{i,j,p=q+1}{\sum}}R_{i\text{a}ip}R_{j\text{a}jp}-\frac
{1}{6}\overset{q}{\underset{\text{a}=1}{\sum}}%
\overset{n}{\underset{i,j,p=q+1}{\sum}}R_{i\text{a}jp}^{2}-\frac{1}%
{3}\overset{q}{\underset{\text{a}=1i,j,}{\sum}}%
\overset{n}{\underset{p=q+1}{\sum}}R_{i\text{a}jp}R_{j\text{a}ip}$

$-\frac{1}{6}\overset{q}{\underset{\text{a}=1}{\sum}}%
\overset{n}{\underset{i,j,p=q+1}{\sum}}R_{j\text{a}ip}^{2}$

$-\frac{1}{3}\overset{q}{\underset{\text{b}=1i,j,}{\sum}}%
\overset{n}{\underset{p=q+1}{\sum}}R_{i\text{b}ip}R_{j\text{b}jp}-\frac{1}%
{6}\overset{q}{\underset{\text{b}=1}{\sum}}%
\overset{n}{\underset{i,j,p=q+1}{\sum}}R_{i\text{b}jp}^{2}-\frac{1}%
{3}\overset{q}{\underset{\text{b}=1}{\sum}}%
\overset{n}{\underset{i.j,p=q+1}{\sum}}R_{i\text{b}jp}R_{j\text{b}ip}$

$-\frac{1}{6}\overset{q}{\underset{\text{b}=1}{\sum}}%
\overset{n}{\underset{i,j,p=q+1}{\sum}}R_{j\text{b}ip}^{2}](y_{0})$

$-\frac{1}{48}$ $\overset{q}{\underset{\text{a,b,c=1}}{\sum}}[$
$-\overset{n}{\underset{i=q+1}{\sum}}R_{i\text{a}i\text{a}}(R_{\text{bcbc}%
}^{P}-R_{\text{bcbc}}^{M})$ $-\overset{n}{\underset{j=q+1}{\sum}}%
R_{j\text{a}j\text{a}}(R_{\text{bcbc}}^{P}-R_{\text{bcbc}}^{M})$

\ $+\overset{n}{\underset{i=q+1}{\sum}}R_{i\text{a}i\text{b}}(R_{\text{acbc}%
}^{P}-R_{\text{acbc}}^{M})\ -\overset{n}{\underset{i=q+1}{\sum}}%
R_{i\text{a}i\text{c}}(R_{\text{abbc}}^{P}-R_{\text{abbc}}^{M})$

$+\overset{n}{\underset{j=q+1}{\sum}}R_{j\text{a}j\text{b}}(R_{\text{acbc}%
}^{P}-R_{\text{acbc}}^{M})$\ $-\overset{n}{\underset{j=q+1}{\sum}}%
R_{j\text{a}j\text{c}}(R_{\text{abbc}}^{P}-R_{\text{abbc}}^{M})$

$+\underset{i,j=q+1}{\overset{n}{\sum}}$ $-R_{i\text{a}j\text{a}}%
(T_{\text{bb}i}T_{\text{cc}j}$ $-T_{\text{bc}i}T_{\text{bc}j})$
$-\underset{i,j=q+1}{\overset{n}{\sum}}R_{i\text{a}j\text{a}}(T_{\text{bb}%
j}T_{\text{cc}i}$ $-T_{\text{bc}j}T_{\text{bc}i})$

$+$ $\underset{i,j=q+1}{\overset{n}{\sum}}$ $-R_{j\text{a}i\text{a}%
}(T_{\text{bb}i}T_{\text{cc}j}$ $-T_{\text{bc}i}T_{\text{bc}j})$
$-\underset{i,j=q+1}{\overset{n}{\sum}}R_{j\text{a}i\text{a}}(T_{\text{bb}%
j}T_{\text{cc}i}$ $-T_{\text{bc}j}T_{\text{bc}i})$

$+\underset{i,j=q+1}{\overset{n}{\sum}}\ R_{i\text{a}j\text{b}}(T_{\text{ab}%
i}T_{\text{cc}j}-T_{\text{bc}i}T_{\text{ac}j}%
)\ +\underset{i,j=q+1}{\overset{n}{\sum}}\ R_{i\text{a}j\text{b}}%
(T_{\text{ab}j}T_{\text{cc}i}-T_{\text{bc}j}T_{\text{ac}i})$

$+\underset{i,j=q+1}{\overset{n}{\sum}}\ R_{j\text{a}i\text{ib}}%
(T_{\text{ab}i}T_{\text{cc}j}-T_{\text{bc}i}T_{\text{ac}j}%
)\ +\underset{i,j=q+1}{\overset{n}{\sum}}\ R_{j\text{a}i\text{b}}%
(T_{\text{ab}j}T_{\text{cc}i}-T_{\text{bc}j}T_{\text{ac}i})\qquad$

$+\underset{i,j=q+1}{\overset{n}{\sum}}-R_{i\text{a}j\text{c}}(T_{\text{ab}%
i}T_{\text{bc}j}-T_{\text{ac}i}T_{\text{bb}j}%
)-\underset{i,j=q+1}{\overset{n}{\sum}}R_{i\text{a}j\text{c}}(T_{\text{ba}%
j}T_{\text{bc}i}-T_{\text{ac}j}T_{\text{bb}i})$

$+\underset{i,j=q+1}{\overset{n}{\sum}}-R_{j\text{a}i\text{c}}(T_{\text{ba}%
i}T_{\text{bc}j}-T_{\text{ac}i}T_{\text{bb}j}%
)-\underset{i,j=q+1}{\overset{n}{\sum}}R_{j\text{a}i\text{c}}(T_{\text{ba}%
j}T_{\text{bc}i}-T_{\text{ac}j}T_{\text{bb}i})](y_{0})$

$+\frac{1}{144}\underset{p=q+1}{\overset{n}{\sum}}%
[\underset{i=q+1}{\overset{n}{\sum}}\overset{q}{\underset{\text{b,c=1}}{\sum}%
}R_{ipip}(R_{\text{bcbc}}^{P}-R_{\text{bcbc}}^{M}%
)+\underset{j=q+1}{\overset{n}{\sum}}$ $\overset{q}{\underset{\text{b,c=1}%
}{\sum}}R_{jpjp}(R_{\text{bcbc}}^{P}-R_{\text{bcbc}}^{M})](y_{0})$

$+\frac{1}{72}\underset{i,j,p=q+1}{\overset{n}{\sum}}%
\overset{q}{\underset{\text{b,c=1}}{\sum}}[R_{ipjp}(T_{\text{bb}i}%
T_{\text{cc}j}-T_{\text{bc}i}T_{\text{bc}j})+R_{ipjp}(T_{\text{bb}%
j}T_{\text{cc}i}-T_{\text{bc}j}T_{\text{bc}i})](y_{0})\qquad$

$-\frac{1}{288}\underset{i,j=q+1}{\overset{n}{\sum}}[T_{\text{aa}%
i}T_{\text{bb}j}(T_{\text{cc}i}T_{\text{dd}j}-T_{\text{cd}i}T_{\text{dc}%
j})+T_{\text{aa}i}T_{\text{bb}j}(T_{\text{cc}j}T_{\text{dd}i}-T_{\text{cd}%
j}T_{\text{dc}i})$

$+T_{\text{aa}j}T_{\text{bb}i}(T_{\text{cc}i}T_{\text{dd}j}-T_{\text{cd}%
i}T_{\text{dc}j})+T_{\text{aa}j}T_{\text{bb}i}(T_{\text{cc}j}T_{\text{dd}%
i}-T_{\text{cd}j}T_{\text{dc}i})](y_{0})$

$+\frac{1}{288}\underset{i,j=q+1}{\overset{n}{\sum}}[T_{\text{aa}%
i}T_{\text{bc}j}(T_{\text{bc}i}T_{\text{dd}j}-T_{\text{bd}i}T_{\text{cd}%
j})+T_{\text{aa}i}T_{\text{bc}j}(T_{\text{bc}j}T_{\text{dd}i}-T_{\text{bd}%
j}T_{\text{cd}i})$

$+T_{\text{aa}j}T_{\text{bc}i}(T_{\text{bc}i}T_{\text{dd}j}-T_{\text{bd}%
i}T_{\text{cd}j})+T_{\text{aa}j}T_{\text{bc}i}(T_{\text{bc}j}T_{\text{dd}%
i}-T_{\text{bd}j}T_{\text{cd}i})](y_{0})$

$-\frac{1}{288}\underset{i,j=q+1}{\overset{n}{\sum}}[T_{\text{aa}%
i}T_{\text{bd}j}(T_{\text{bc}i}T_{\text{cd}j}-T_{\text{bd}i}T_{\text{cc}%
j})+T_{\text{aa}i}T_{\text{bd}j}(T_{\text{bc}j}T_{\text{cd}i}-T_{\text{bd}%
j}T_{\text{cc}i})$

$+T_{\text{aa}j}T_{\text{bd}i}(T_{\text{bc}i}T_{\text{cd}j}-T_{\text{bd}%
i}T_{\text{cc}j})+T_{\text{aa}j}T_{\text{bd}i}(T_{\text{bc}j}T_{\text{cd}%
i}-T_{\text{bd}j}T_{\text{cc}i})](y_{0})\qquad$

$+\frac{1}{288}\underset{i,j=q+1}{\overset{n}{\sum}}[T_{\text{ab}%
i}T_{\text{ab}j}(T_{\text{cc}i}T_{\text{dd}j}-T_{\text{cd}i}T_{\text{dc}%
j})+T_{\text{ab}i}T_{\text{ab}j}(T_{\text{cc}j}T_{\text{dd}i}-T_{\text{cd}%
j}T_{\text{dc}i})$

$+T_{\text{ab}j}T_{\text{ab}i}(T_{\text{cc}i}T_{\text{dd}j}-T_{\text{cd}%
i}T_{\text{dc}j})+T_{\text{ab}j}T_{\text{ab}i}(T_{\text{cc}j}T_{\text{dd}%
i}-T_{\text{cd}j}T_{\text{dc}i})](y_{0})$

$-\frac{1}{288}\underset{i,j=q+1}{\overset{n}{\sum}}[T_{\text{ab}%
i}T_{\text{bc}j}(T_{\text{ac}i}T_{\text{dd}j}-T_{\text{ad}i}T_{\text{cd}%
j})+T_{\text{ab}i}T_{\text{bc}j}(T_{\text{ac}j}T_{\text{dd}i}-T_{\text{ad}%
j}T_{\text{cd}i})$

$+T_{\text{ab}j}T_{\text{bc}i}(T_{\text{ac}i}T_{\text{dd}j}-T_{\text{ad}%
i}T_{\text{cd}j})+T_{\text{ab}j}T_{\text{bc}i}(T_{\text{ac}j}T_{\text{dd}%
i}-T_{\text{ad}j}T_{\text{cd}i})](y_{0})$

$+\frac{1}{288}\underset{i,j=q+1}{\overset{n}{\sum}}[T_{\text{ab}%
i}T_{\text{bd}j}(T_{\text{ac}i}T_{\text{cd}j}-T_{\text{ad}i}T_{\text{cc}%
j})+T_{\text{ab}i}T_{\text{bd}j}(T_{\text{ac}j}T_{\text{cd}i}-T_{\text{ad}%
j}T_{\text{cc}i})$

$+T_{\text{ab}i}T_{\text{bd}j}(T_{\text{ac}j}T_{\text{cd}i}-T_{\text{ad}%
j}T_{\text{cc}i})+T_{\text{ab}j}T_{\text{bd}i}(T_{\text{ac}j}T_{\text{cd}%
i}-T_{\text{ad}j}T_{\text{cc}i})](y_{0})$

$-\ \frac{1}{288}\underset{i,j=q+1}{\overset{n}{\sum}}[T_{\text{ac}%
i}T_{\text{ab}j}(T_{\text{bc}i}T_{\text{dd}j}-T_{\text{bd}i}T_{\text{dc}%
j})+T_{\text{ac}i}T_{\text{ab}j}(T_{\text{bc}j}T_{\text{dd}i}-T_{\text{bd}%
j}T_{\text{dc}i})$

$+T_{\text{ac}j}T_{\text{ab}i}(T_{\text{bc}i}T_{\text{dd}j}-T_{\text{bd}%
i}T_{\text{dc}j})+T_{\text{ac}j}T_{\text{ab}i}(T_{\text{bc}j}T_{\text{dd}%
i}-T_{\text{bd}j}T_{\text{dc}i})](y_{0})$

$+\ \frac{1}{288}\underset{i,j=q+1}{\overset{n}{\sum}}[T_{\text{ac}%
i}T_{\text{bb}j}(T_{\text{ac}i}T_{\text{dd}j}-T_{\text{ad}i}T_{\text{cd}%
j})+T_{\text{ac}i}T_{\text{bb}j}(T_{\text{ac}j}T_{\text{dd}i}-T_{\text{ad}%
j}T_{\text{cd}i})$

$+T_{\text{ac}j}T_{\text{bb}i}(T_{\text{ac}i}T_{\text{dd}j}-T_{\text{ad}%
i}T_{\text{cd}i})+T_{\text{ac}j}T_{\text{bb}i}(T_{\text{ac}j}T_{\text{dd}%
i}-T_{\text{ad}j}T_{\text{cd}i})](y_{0})$

$-\ \frac{1}{288}\underset{i,j=q+1}{\overset{n}{\sum}}[T_{\text{ac}%
i}T_{\text{bd}j}(T_{\text{ac}i}T_{\text{bd}j}-T_{\text{ad}i}T_{\text{bc}%
j})+T_{\text{ac}i}T_{\text{bd}j}(T_{\text{ac}j}T_{\text{bd}i}-T_{\text{ad}%
j}T_{\text{bc}i})$

$+T_{\text{ac}j}T_{\text{bd}i}(T_{\text{ac}i}T_{\text{bd}j}-T_{\text{ad}%
i}T_{\text{bc}j})+T_{\text{ac}j}T_{\text{bd}i}(T_{\text{ac}j}T_{\text{bd}%
i}-T_{\text{ad}j}T_{\text{bc}i})](y_{0})$

$+\frac{1}{288}\underset{i,j=q+1}{\overset{n}{\sum}}[T_{\text{ad}%
i}T_{\text{ab}j}(T_{\text{bc}i}T_{\text{cd}j}-T_{\text{bd}i}T_{\text{cc}%
j})+T_{\text{ad}i}T_{\text{ab}j}(T_{\text{bc}j}T_{\text{cd}i}-T_{\text{bd}%
j}T_{\text{cc}i})$

$+T_{\text{ad}j}T_{\text{ab}i}(T_{\text{bc}i}T_{\text{cd}j}-T_{\text{bd}%
i}T_{\text{cc}j})+T_{\text{ad}j}T_{\text{ab}i}(T_{\text{bc}j}T_{\text{cd}%
i}-T_{\text{bd}j}T_{\text{cc}i})](y_{0})$

$-\ \frac{1}{288}\underset{i,j=q+1}{\overset{n}{\sum}}[T_{\text{ad}%
i}T_{\text{bb}j}(T_{\text{ac}i}T_{\text{cd}j}-T_{\text{ad}i}T_{\text{cc}%
j})+T_{\text{ad}i}T_{\text{bb}j}(T_{\text{ac}j}T_{\text{cd}i}-T_{\text{ad}%
j}T_{\text{cc}i})$

$+T_{\text{ad}j}T_{\text{bb}i}(T_{\text{ac}i}T_{\text{cd}j}-T_{\text{ad}%
i}T_{\text{cc}j})+T_{\text{ad}j}T_{\text{bb}i}(T_{\text{ac}j}T_{\text{cd}%
i}-T_{\text{ad}j}T_{\text{cc}i})](y_{0})$

$+\ \frac{1}{288}\underset{i,j=q+1}{\overset{n}{\sum}}[T_{\text{ad}%
i}T_{\text{bc}j}(T_{\text{ac}i}T_{\text{bd}j}-T_{\text{ad}i}T_{\text{bc}%
j})+T_{\text{ad}i}T_{\text{bc}j}(T_{\text{ac}j}T_{\text{bd}i}-T_{\text{ad}%
j}T_{\text{bc}i})$

$+T_{\text{ad}j}T_{\text{bc}i}(T_{\text{ac}i}T_{\text{bd}j}-T_{\text{ad}%
i}T_{\text{bc}j})+T_{\text{ad}j}T_{\text{bc}i}(T_{\text{ac}j}T_{\text{bd}%
i}-T_{\text{ad}j}T_{\text{bc}i})](y_{0})$

$-\ \frac{1}{144}[(R_{\text{cdcd}}^{P}-R_{\text{cdcd}}^{M})(R_{\text{abab}%
}^{P}-R_{\text{abab}}^{M})](y_{0})$

$\ +\frac{1}{144}[(R_{\text{bdcd}}^{P}-R_{\text{bdcd}}^{M})(R_{\text{abac}%
}^{P}-R_{\text{abac}}^{M})](y_{0})$

$\ +\ \frac{1}{144}[(R_{\text{bcdc}}^{P}-R_{\text{bcdc}}^{M})(R_{\text{abad}%
}^{P}-R_{\text{abad}}^{M})](y_{0})$

$-\ \frac{1}{144}[(R_{\text{adcd}}^{P}-R_{\text{adcd}}^{M})(R_{\text{abbc}%
}^{P}-R_{\text{abbc}}^{M})](y_{0})$

$+\frac{1}{144}[(R_{\text{acdc}}^{P}-R_{\text{acdc}}^{M})(R_{\text{abdb}}%
^{P}-R_{\text{abdb}}^{M})](y_{0})$

$\ -\ \frac{1}{576}[(R_{\text{abcd}}^{P}-R_{\text{abcd}}^{M})]^{2}(y_{0})$

$-\frac{1}{144}<H,i>(y_{0})<H,j>(y_{0})\qquad\qquad\qquad\qquad\qquad
\qquad\qquad L_{3}$

$\times\lbrack2\varrho_{ij}+$ $\overset{q}{\underset{\text{a}=1}{4\sum}%
}R_{i\text{a}j\text{a}}-3\overset{q}{\underset{\text{a,b=1}}{\sum}%
}(T_{\text{aa}i}T_{\text{bb}j}-T_{\text{ab}i}T_{\text{ab}j})$

$-3\overset{q}{\underset{\text{a,b=1}}{\sum}}(T_{\text{aa}j}T_{\text{bb}%
i}-T_{\text{ab}j}T_{\text{ab}i}](y_{0})\phi(y_{0})$

$-\frac{1}{16}[<H,i>^{2}(y_{0})<H,j>^{2}](y_{0})\phi(y_{0})$

$-\frac{1}{144}<H,i>(y_{0})<H,j>(y_{0})$

$\times\lbrack2\varrho_{ij}+$ $\overset{q}{\underset{\text{a}=1}{4\sum}%
}R_{i\text{a}j\text{a}}-3\overset{q}{\underset{\text{a,b=1}}{\sum}%
}(T_{\text{aa}i}T_{\text{bb}j}-T_{\text{ab}i}T_{\text{ab}j})$

$-3\overset{q}{\underset{\text{a,b=1}}{\sum}}(T_{\text{aa}j}T_{\text{bb}%
i}-T_{\text{ab}j}T_{\text{ab}i}](y_{0})\phi(y_{0})$

$-\frac{1}{72}<H,i>(y_{0})<H,k>(y_{0})R_{jijk}(y_{0})\phi(y_{0})$

$-\frac{1}{16}<H,i>^{2}(y_{0})<H,j>^{2}(y_{0})\phi(y_{0})$

$-\frac{1}{72}<H,i>^{2}(y_{0})[\varrho_{jj}+$ $\overset{q}{\underset{\text{a}%
=1}{2\sum}}R_{j\text{a}j\text{a}}-3\overset{q}{\underset{\text{a,b=1}}{\sum}%
}(T_{\text{aa}j}T_{\text{bb}j}-T_{\text{ab}j}T_{\text{ab}j})](y_{0})\phi
(y_{0})$

$+\frac{5}{32}[<H,i>^{2}<H,j>^{2}](y_{0})\phi(y_{0})$

$+\frac{1}{48}<H,i>(y_{0})<H,j>$

$\times\lbrack2\varrho_{ij}+$ $\overset{q}{\underset{\text{a}=1}{4\sum}%
}R_{i\text{a}j\text{a}}-3\overset{q}{\underset{\text{a,b=1}}{\sum}%
}(T_{\text{aa}i}T_{\text{bb}j}-T_{\text{ab}i}T_{\text{ab}j}%
)-3\overset{q}{\underset{\text{a,b=1}}{\sum}}(T_{\text{aa}j}T_{\text{bb}%
i}-T_{\text{ab}j}T_{\text{ab}i}](y_{0})\phi(y_{0})$

$+\frac{1}{48}<H,i>(y_{0})<H,i>(y_{0})[\tau^{M}\ -3\tau^{P}%
+\ \underset{\text{a}=1}{\overset{\text{q}}{\sum}}\varrho_{\text{aa}}^{M}+$
$\overset{q}{\underset{\text{a},\text{b}=1}{\sum}}R_{\text{abab}}^{M}$
$](y_{0})\phi(y_{0})$

$+\frac{1}{144}<H,i>(y_{0})[\nabla_{i}\varrho_{jj}-2\varrho_{ij}%
<H,j>+\overset{q}{\underset{\text{a}=1}{\sum}}(\nabla_{i}R_{\text{a}%
j\text{a}j}-4R_{i\text{a}j\text{a}}<H,j>)$

$+4\overset{q}{\underset{\text{a,b=1}}{\sum}}R_{i\text{a}j\text{b}%
}T_{\text{ab}j}+2\overset{q}{\underset{\text{a,b,c=1}}{\sum}}(T_{\text{aa}%
i}T_{\text{bb}j}T_{\text{cc}j}-3T_{\text{aa}i}T_{\text{bc}j}T_{\text{bc}%
j}+2T_{\text{ab}i}T_{\text{bc}j}T_{\text{ca}j})](y_{0})\phi(y_{0})$%
\qquad\qquad\qquad\qquad\qquad\ \ 

$+\frac{1}{144}<H,i>(y_{0})[\nabla_{j}\varrho_{ij}-2\varrho_{ij}%
<H,j>+\overset{q}{\underset{\text{a}=1}{\sum}}(\nabla_{j}R_{\text{a}%
i\text{a}j}-4R_{j\text{a}i\text{a}}<H,j>)$

$+4\overset{q}{\underset{\text{a,b=1}}{\sum}}R_{j\text{a}i\text{b}%
}T_{\text{ab}j}+2\overset{q}{\underset{\text{a,b,c=1}}{\sum}}(T_{\text{aa}%
j}T_{\text{bb}i}T_{\text{cc}j}-3T_{\text{aa}j}T_{\text{bc}i}T_{\text{bc}%
j}+2T_{\text{ab}j}T_{\text{bc}i}T_{\text{ac}j})](y_{0})\phi(y_{0})$

$+\frac{1}{144}<H,i>(y_{0})[\nabla_{j}\varrho_{ij}-2\varrho_{jj}%
<H,i>+\overset{q}{\underset{\text{a}=1}{\sum}}(\nabla_{j}R_{\text{a}%
i\text{a}j}-4R_{j\text{a}j\text{a}}<H,i>)$

$+4\overset{q}{\underset{\text{a,b=1}}{\sum}}R_{j\text{a}j\text{b}%
}T_{\text{ab}i}+2\overset{q}{\underset{\text{a,b,c=1}}{\sum}}(T_{\text{aa}%
j}T_{\text{bb}j}T_{\text{cc}i}-3T_{\text{aa}j}T_{\text{bc}j}T_{\text{bc}%
i}+2T_{\text{ab}j}T_{\text{bc}j}T_{\text{ac}i})](y_{0})\phi(y_{0})$

$-\frac{1}{192}<H,i>^{2}<H,j>^{2}(y_{0})\phi(y_{0})\qquad\qquad\qquad
\qquad\qquad\qquad\qquad$A$_{3213}\ $

$-\frac{1}{288}<H,i>^{2}(y_{0})[\tau^{M}\ -3\tau^{P}+\ \underset{\text{a}%
=1}{\overset{\text{q}}{\sum}}\varrho_{\text{aa}}^{M}+$
$\overset{q}{\underset{\text{a},\text{b}=1}{\sum}}R_{\text{abab}}^{M}$
$](y_{0})\phi(y_{0})$

$-\frac{1}{288}<H,i>(y_{0})<H,j>(y_{0})$

$\times\lbrack2\varrho_{ij}+$ $\overset{q}{\underset{\text{a}=1}{4\sum}%
}R_{i\text{a}j\text{a}}-3\overset{q}{\underset{\text{a,b=1}}{\sum}%
}(T_{\text{aa}i}T_{\text{bb}j}-T_{\text{ab}i}T_{\text{ab}j}%
)-3\overset{q}{\underset{\text{a,b=1}}{\sum}}(T_{\text{aa}j}T_{\text{bb}%
i}-T_{\text{ab}j}T_{\text{ab}i}](y_{0})\phi(y_{0})$

$+\frac{1}{144}<H,i>(y_{0})<H,k>(y_{0})R_{jijk}(y_{0})$

$+\frac{1}{144}<H,i>^{2}(y_{0})[\varrho_{jj}+$ $\overset{q}{\underset{\text{a}%
=1}{2\sum}}R_{j\text{a}j\text{a}}-3\overset{q}{\underset{\text{a,b=1}}{\sum}%
}(T_{\text{aa}j}T_{\text{bb}j}-T_{\text{ab}j}T_{\text{ab}j})](y_{0})\phi
(y_{0})$

$-\frac{1}{288}<H,i>(y_{0})[\nabla_{i}\varrho_{jj}-2\varrho_{ij}%
<H,j>+\overset{q}{\underset{\text{a}=1}{\sum}}(\nabla_{i}R_{\text{a}%
j\text{a}j}-4R_{i\text{a}j\text{a}}<H,j>)$

$+4\overset{q}{\underset{\text{a,b=1}}{\sum}}R_{i\text{a}j\text{b}%
}T_{\text{ab}j}+2\overset{q}{\underset{\text{a,b,c=1}}{\sum}}(T_{\text{aa}%
i}T_{\text{bb}j}T_{\text{cc}j}-3T_{\text{aa}i}T_{\text{bc}j}T_{\text{bc}%
j}+2T_{\text{ab}i}T_{\text{bc}j}T_{\text{ca}j})](y_{0})\phi(y_{0})$%
\qquad\qquad\qquad\qquad\qquad\ \ 

$-\frac{1}{288}<H,i>(y_{0})[\nabla_{j}\varrho_{ij}-2\varrho_{ij}%
<H,j>+\overset{q}{\underset{\text{a}=1}{\sum}}(\nabla_{j}R_{\text{a}%
i\text{a}j}-4R_{j\text{a}i\text{a}}<H,j>)$

$+4\overset{q}{\underset{\text{a,b=1}}{\sum}}R_{j\text{a}i\text{b}%
}T_{\text{ab}j}+2\overset{q}{\underset{\text{a,b,c=1}}{\sum}}(T_{\text{aa}%
j}T_{\text{bb}i}T_{\text{cc}j}-3T_{\text{aa}j}T_{\text{bc}i}T_{\text{bc}%
j}+2T_{\text{ab}j}T_{\text{bc}i}T_{\text{ac}j})](y_{0})\phi(y_{0})$

$-\frac{1}{288}<H,i>(y_{0})[\nabla_{j}\varrho_{ij}-2\varrho_{jj}%
<H,i>+\overset{q}{\underset{\text{a}=1}{\sum}}(\nabla_{j}R_{\text{a}%
i\text{a}j}-4R_{j\text{a}j\text{a}}<H,i>)$

$+4\overset{q}{\underset{\text{a,b=1}}{\sum}}R_{j\text{a}j\text{b}%
}T_{\text{ab}i}+2\overset{q}{\underset{\text{a,b,c=1}}{\sum}}(T_{\text{aa}%
j}T_{\text{bb}j}T_{\text{cc}i}-3T_{\text{aa}j}T_{\text{bc}j}T_{\text{bc}%
i}+2T_{\text{ab}j}T_{\text{bc}j}T_{\text{ac}i})](y_{0})\phi(y_{0})$

\qquad\qquad\qquad\qquad\qquad\qquad\qquad\qquad\qquad\qquad\qquad\qquad
\qquad\qquad\qquad\qquad\qquad\qquad$\blacksquare$

In \textbf{Normal Coordinates} the Submanifold reduces to the singleton
$\left\{  y_{0}\right\}  $ and we have:

$\left(  A_{31}\right)  \qquad$A$_{321}=\frac{1}{24}\frac{\partial^{2}%
}{\partial\text{x}_{i}^{2}}(\theta^{\frac{1}{2}}\Delta\theta^{-\frac{1}{2}%
})\phi(y_{0})$\ \ $=$\ \ \ A$_{3211}+$ A$_{3212}+$ A$_{3213}$

$\qquad=-\frac{1}{3456}[+2(\tau^{M}]^{2}(y_{0})\phi(y_{0})\phi(y_{0})$\ \qquad
A$_{321}$\qquad A$_{3211}$

\qquad$+\frac{1}{24}[+\frac{1}{3}(\tau^{M}]\times\lbrack+\frac{1}{6}(\tau
^{M})](y_{0})\phi(y_{0})\qquad$A$_{3212}=\frac{1}{24}(L_{1}+L_{2}+L_{3}%
)\qquad$

$\qquad-\frac{1}{24}\times\frac{1}{36}[2\varrho_{ij}]^{2}(y_{0})\phi
(y_{0})\qquad\qquad\qquad L_{2}\qquad L_{21}$\qquad\qquad\qquad\qquad
\ \ $\qquad\qquad$

$\qquad-\frac{1}{24}\times\frac{1}{18}R_{jijk}(y_{0})[2\varrho_{ik}%
](y_{0})\phi(y_{0})L_{22}$\qquad\qquad\qquad\qquad\ \ $\qquad$

$\qquad-\frac{1}{24}\times\frac{1}{18}[\varrho_{ii}](y_{0})\times
\lbrack\varrho_{jj}(y_{0})\phi(y_{0})\qquad L_{232}\qquad$\ $L_{233}$

$+\frac{1}{24}\times\frac{1}{18}R_{ijik}(y_{0})[2\varrho_{jk}](y_{0}%
)\phi(y_{0})$

$+\frac{1}{576}\overset{n}{\underset{i,j=q+1}{\sum}}[2\varrho_{ij}]^{2}%
(y_{0})+\frac{1}{288}[\tau^{M}\ ]^{2}(y_{0})\qquad\ \frac{1}{24}\frac
{\partial^{4}\theta^{-\frac{1}{2}}}{\partial x_{i}^{2}\partial x_{j}^{2}%
}(y_{0})$

$-\ \frac{1}{288}\overset{n}{\underset{i,j=q+1}{\sum}}%
[\overset{n}{\underset{p=q+1}{\sum}}(-\frac{3}{5}\nabla_{ii}^{2}%
(R)_{jpjp}+\overset{n}{\underset{p=q+1}{\sum}}(-\frac{3}{5}\nabla_{jj}%
^{2}(R)_{ipip}\qquad\qquad A$

$+\overset{n}{\underset{p=q+1}{\sum}}(-\frac{3}{5}\nabla_{ij}^{2}%
(R)_{ipjp}+\overset{n}{\underset{p=q+1}{\sum}}(-\frac{3}{5}\nabla_{ij}%
^{2}(R)_{jpip}+\overset{n}{\underset{p=q+1}{\sum}}(-\frac{3}{5}\nabla_{ji}%
^{2}(R)_{ipjp}$

$+\overset{n}{\underset{p=q+1}{\sum}}(-\frac{3}{5}\nabla_{ji}^{2}%
(R)_{jpip})+\frac{1}{5}\overset{n}{\underset{m,p=q+1}{%
{\textstyle\sum}
}}R_{ipim}R_{jpjm}+\frac{1}{5}\overset{n}{\underset{m,p=q+1}{%
{\textstyle\sum}
}}R_{jpjm}R_{ipim}$

$+\frac{1}{5}\overset{n}{\underset{m,p=q+1}{%
{\textstyle\sum}
}}R_{ipjm}R_{ipjm}+\frac{1}{5}\overset{n}{\underset{m,p=q+1}{%
{\textstyle\sum}
}}R_{ipjm}R_{jpim}+\frac{1}{5}\overset{n}{\underset{m,p=q+1}{%
{\textstyle\sum}
}}R_{jpim}R_{ipjm}$

$+\frac{1}{5}\overset{n}{\underset{m,p=q+1}{%
{\textstyle\sum}
}}R_{jpim}R_{jpim}\}(y_{0})+\frac{1}{9}(\tau^{M})(\tau^{M})+\frac{2}%
{9}(\left\Vert \varrho^{M}\right\Vert ^{2})-\frac{1}{9}%
\overset{n}{\underset{i,j,p,m=q+1}{\sum}}R_{ipim}R_{jpjm}\ $

$-\frac{1}{18}\overset{n}{\underset{i,j,p,m=q+1}{\sum}}R_{ipjm}^{2}-\frac
{1}{9}\overset{n}{\underset{i,j,p,m=q+1}{\sum}}R_{ipjm}R_{jpim}-\frac{1}%
{18}\overset{n}{\underset{i,j,p,m=q+1}{\sum}}R_{jpim}^{2}](y_{0})$

\qquad\qquad\qquad$\qquad\qquad\qquad\qquad\qquad\qquad\qquad\qquad
\qquad\qquad\qquad\qquad\qquad\blacksquare$\qquad

\chapter{The Vector Field X and its Derivatives}

We recall that the smooth map $\Phi:M\longrightarrow R$ defined in $\left(
1.5\right)  $ of Chapter 1 using the vector field X. Here we will compute
derivatives of the map $\Phi$ and those of the associated vector field
$\nabla\log\Phi$ and Laplacian $\Delta\Phi.$

\section{Table B$_{1}:$ \textbf{Derivatives of }$\nabla\log\Phi$}

\qquad For a general vector field X on M and a general point of a tubular
neighbourhood x$_{0}\in$M$_{0},$ we have:

Throughout the computations in this Appendix the \textbf{Einstein summation
convention} for repeated indices will often be used. However, in some cases,
the summation symbol will be explicitely written for emphasis.

General formulae are given in (i)-(v) and more precise formulae are given in (vi)-(xvi).

\subsection{Normal Derivatives}

We have the follwing beautiful formulae:

(i) \ For $j=q+1,...,n,$

$\qquad(\nabla$log$\Phi_{P})_{j}(x_{0})+X_{j}(x_{0})=-$
$\underset{k=q+1}{\overset{n}{\sum}}x_{k}(x_{0})\left(  \frac{\partial
}{\partial x_{j}}(\nabla\log\Phi_{P})_{k}+\frac{\partial X_{k}}{\partial
x_{j}}\right)  (x_{0})$

Recall that the definition of Fermi coordinates in $\left(  1.1\right)
-\left(  1.2\right)  $ of Chapter 1,

we have the important property: x$_{j}(y)=0$ \ for $j=q+1,...,n$ for all $y\in
U\subset P$

where U is a (small) neighbourhood of the centre of Fermi coordinates
y$_{0}\in$P.

An important consequence of this property is that:

$(\nabla$log$\Phi_{P})_{j}(y)=-X_{j}(y)$ for all $y\in U\subset P.$

In particular, $(\nabla$log$\Phi_{P})_{j}(y_{0})=-X_{j}(y_{0})$ where $y_{0}$
is the centre of Fermi coordinates.

This property will be consistently used in the second part of the appendices:

(ii) For $i,j=q+1,...,n,$ we have from $\left(  B_{7}\right)  :$

\qquad$\lbrack\frac{\partial}{\partial x_{i}}(\nabla$log$\Phi_{P})_{j}%
+\frac{\partial}{\partial x_{j}}(\nabla\log\Phi_{P})_{i}+\frac{\partial X_{j}%
}{\partial x_{i}}+\frac{\partial X_{i}}{\partial x_{j}}](x_{0})$

$=-\underset{k=q+1}{\overset{n}{\sum}}x_{k}(x_{0})\left(  \frac{\partial^{2}%
}{\partial x_{i}\partial x_{j}}(\nabla\log\Phi_{P})_{k}+\frac{\partial
^{2}X_{k}}{\partial x_{i}\partial x_{j}}\right)  (x_{0})\qquad$

In particular, for all $y\in U\subset P,$

$[\frac{\partial}{\partial x_{i}}(\nabla$log$\Phi_{P})_{j}+\frac{\partial
}{\partial x_{j}}(\nabla\log\Phi_{P})_{i}](y)=-[\frac{\partial X_{j}}{\partial
x_{i}}+\frac{\partial X_{i}}{\partial x_{j}}](y)$

(iii) \ For $i,j,k=q+1,...,n,$ we have from $\left(  B_{8}\right)  :$

$[\frac{\partial^{2}}{\partial x_{i}\partial x_{j}}(\nabla$log$\Phi_{P}%
)_{k}+\frac{\partial^{2}}{\partial x_{k}\partial x_{i}}(\nabla\log\Phi
_{P})_{j}+\frac{\partial^{2}}{\partial x_{j}\partial x_{k}}(\nabla\log\Phi
_{P})_{i}$

$+\frac{\partial^{2}X_{k}}{\partial x_{i}\partial x_{j}}+\frac{\partial
^{2}X_{j}}{\partial x_{k}\partial x_{i}}+\frac{\partial^{2}X_{i}}{\partial
x_{j}\partial x_{k}}](x_{0})\qquad\ \ \ \ \ \ \ \ \ \qquad$

$=-\underset{l=q+1}{\overset{n}{\sum}}x_{l}(x_{0})\left(  \frac{\partial^{3}%
}{\partial x_{i}\partial x_{j}\partial x_{k}}(\nabla\log\Phi_{P})_{l}%
+\frac{\partial^{3}X_{l}}{\partial x_{i}\partial x_{j}\partial x_{k}}\right)
(x_{0})$

In particular, for all $y\in U\subset P,$

$[\frac{\partial^{2}}{\partial x_{i}\partial x_{j}}(\nabla$log$\Phi_{P}%
)_{k}+\frac{\partial^{2}}{\partial x_{k}\partial x_{i}}(\nabla\log\Phi
_{P})_{j}+\frac{\partial^{2}}{\partial x_{j}\partial x_{k}}(\nabla\log\Phi
_{P})_{i}](y)$

$=-[\frac{\partial^{2}X_{k}}{\partial x_{i}\partial x_{j}}+\frac{\partial
^{2}X_{j}}{\partial x_{k}\partial x_{i}}+\frac{\partial^{2}X_{i}}{\partial
x_{j}\partial x_{k}}](y)$

(iii)$^{\ast}$ For $i,j,k,l=q+1,...,n,$ we have from $\left(  B_{9}\right)  :$

\ $[\frac{\partial^{3}}{\partial x_{i}\partial x_{j}\partial x_{k}}(\nabla
$log$\Phi_{P})_{l}+\frac{\partial^{3}}{\partial x_{l}\partial x_{i}\partial
x_{j}}(\nabla\log\Phi_{P})_{k}+\frac{\partial^{3}}{\partial x_{k}\partial
x_{l}\partial x_{i}}(\nabla\log\Phi_{P})_{j}$

$+\frac{\partial^{3}}{\partial x_{j}\partial x_{k}\partial x_{l}}(\nabla
\log\Phi_{P})_{i}\qquad\ \qquad$

$+\frac{\partial^{3}X_{l}}{\partial x_{i}\partial x_{j}\partial x_{k}}%
+\frac{\partial^{3}X_{k}}{\partial x_{l}\partial x_{i}\partial x_{j}}%
+\frac{\partial^{3}X_{j}}{\partial x_{k}\partial x_{l}\partial x_{i}}%
+\frac{\partial^{3}X_{i}}{\partial x_{j}\partial x_{k}\partial x_{l}}](x_{0})$

$=-\underset{m=q+1}{\overset{n}{\sum}}x_{m}(x_{0})\left(  \frac{\partial^{4}%
}{\partial x_{i}\partial x_{j}\partial x_{k}\partial x_{l}}(\nabla\log\Phi
_{P})_{m}+\frac{\partial^{4}X_{m}}{\partial x_{i}\partial x_{j}\partial
x_{k}\partial x_{l}}\right)  (x_{0})$

In particular, for all $y\in U\subset P,$ we have:

\ $[\frac{\partial^{3}}{\partial x_{i}\partial x_{j}\partial x_{k}}(\nabla
$log$\Phi_{P})_{l}+\frac{\partial^{3}}{\partial x_{l}\partial x_{i}\partial
x_{j}}(\nabla\log\Phi_{P})_{k}+\frac{\partial^{3}}{\partial x_{k}\partial
x_{l}\partial x_{i}}(\nabla\log\Phi_{P})_{j}$

$+\frac{\partial^{3}}{\partial x_{j}\partial x_{k}\partial x_{l}}(\nabla
\log\Phi_{P})_{i}](y)\qquad\ \qquad$

$=-[\frac{\partial^{3}X_{l}}{\partial x_{i}\partial x_{j}\partial x_{k}}%
+\frac{\partial^{3}X_{k}}{\partial x_{l}\partial x_{i}\partial x_{j}}%
+\frac{\partial^{3}X_{j}}{\partial x_{k}\partial x_{l}\partial x_{i}}%
+\frac{\partial^{3}X_{i}}{\partial x_{j}\partial x_{k}\partial x_{l}}](y)$
$\qquad$

The formulae for higher derivatives follow.$\qquad$

(iv) For a $=$ 1,...,q and $j=q+1,...,n$

$\frac{\partial}{\partial x_{\text{a}}}(\nabla$log$\Phi_{P})_{j}(x_{0}%
)+\frac{\partial X_{j}}{\partial x_{\text{a}}}(x_{0})=\frac{1}{2}%
\underset{i=q+1}{\overset{n}{\sum}}x_{i}(x_{0})\left(  \frac{\partial
}{\partial x_{\text{a}}\partial x_{j}}(\nabla\log\Phi_{P})_{i}+\frac
{\partial^{2}X_{i}}{\partial x_{\text{a}}\partial x_{j}}\right)  (x_{0})$

\qquad For a,b = 1,...,q and for $j=q+1,...,n,$

(v) $\frac{\partial^{2}}{\partial x_{\text{a}}\partial x_{\text{b}}}(\nabla
$log$\Phi_{P})_{j}(x_{0})+\frac{\partial^{2}X_{j}}{\partial x_{\text{a}%
}\partial x_{\text{b}}}(x_{0})$

$=-\frac{1}{2}\underset{i=q+1}{\overset{n}{\sum}}x_{i}(x_{0})\left(
\frac{\partial^{2}}{\partial x_{\text{a}}\partial x_{\text{b}}\partial x_{j}%
}(\nabla\log\Phi_{P})_{i}+\frac{\partial^{3}X_{i}}{\partial x_{\text{a}%
}\partial x_{\text{b}}\partial x_{j}}\right)  (x_{0})$

For $y\in$U$\subset$P where U is a (small) neighbourhood of the centre of

Fermi coordinates y$_{0}\in$P:

For $j,k=q+1,...,n.$

(vi) $\ (\nabla$log$\Phi_{P})_{j}(y)=-X_{j}(y)$

(vii) $\frac{\partial}{\partial x_{i}}(\nabla$log$\Phi_{P})_{j}(y)+\frac
{\partial}{\partial x_{j}}(\nabla\log\Phi_{P})_{i}(y)=-\frac{\partial X_{j}%
}{\partial x_{i}}(y)-\frac{\partial X_{i}}{\partial x_{j}}(y)$

(vii)$^{\ast}$ We shall see in $\left(  39\right)  $ in \textbf{Appendix
B}$_{4}$ that there is an improved formula:

$\qquad$\ $\frac{\partial}{\partial x_{i}}(\nabla\log\Phi_{P})_{j}%
(y)=-\frac{1}{2}\left(  \frac{\partial X_{j}}{\partial x_{i}}+\frac{\partial
X_{i}}{\partial x_{j}}\right)  (y)=\frac{\partial}{\partial x_{j}}(\nabla
\log\Phi_{P})_{i}(y)\qquad\qquad\qquad\qquad$

(viii) From $\left(  B_{12}\right)  :$

$\qquad\frac{\partial^{2}}{\partial x_{j}\partial x_{k}}(\nabla$log$\Phi
_{P})_{i}(y)+\frac{\partial^{2}}{\partial x_{i}\partial x_{k}}(\nabla\log
\Phi_{P})_{j}(y)+\frac{\partial^{2}}{\partial x_{i}\partial x_{j}}(\nabla
\log\Phi_{P})_{k}(y)$

$\qquad=-\frac{\partial^{2}X_{i}}{\partial x_{j}\partial x_{k}}(y)-\frac
{\partial^{2}X_{j}}{\partial x_{i}\partial x_{k}}(y)-\frac{\partial^{2}X_{k}%
}{\partial x_{i}\partial x_{j}}(y)$

(viii)$^{\ast}$ We shall see in $\left(  59\right)  $ of \textbf{Appendix
B}$_{4}$ that there is an improved but complicated formula:

$[\frac{\partial^{2}}{\partial x_{i}\partial x_{j}}(\nabla\log\Phi_{P}%
)_{k}](y)$

$=-\frac{1}{3}\left(  \frac{\partial^{2}X_{k}}{\partial x_{i}\partial x_{j}%
}+\frac{\partial^{2}X_{j}}{\partial x_{i}\partial x_{k}}+\frac{\partial
^{2}X_{i}}{\partial x_{j}\partial x_{k}}\right)  (y)-\frac{1}{3}%
\underset{l=q+1}{\overset{n}{\sum}}(R_{ikjl}+R_{jkil})(y_{0})X_{l}(y)$

\qquad$+[\perp_{\text{a}ik}\frac{\partial X_{j}}{\partial x_{\text{a}}}%
+\perp_{\text{a}jk}\frac{\partial X_{i}}{\partial x_{\text{a}}}%
](y)+\underset{l=q+1}{\overset{n}{\sum}}[\perp_{\text{a}ik}\perp_{\text{a}%
jl}X_{l}+\perp_{\text{a}jk}\perp_{\text{a}il}X_{l}](y)$

(viii)$^{\ast\ast}$ From $\left(  B_{12}\right)  ^{\ast}:$ For all $y\in P$

$[\frac{\partial^{3}}{\partial x_{i}\partial x_{j}\partial x_{k}}(\nabla
$log$\Phi_{P})_{l}+\frac{\partial^{3}}{\partial x_{l}\partial x_{i}\partial
x_{j}}(\nabla\log\Phi_{P})_{k}+\frac{\partial^{3}}{\partial x_{k}\partial
x_{l}\partial x_{i}}(\nabla\log\Phi_{P})_{j}$

$+\frac{\partial^{3}}{\partial x_{j}\partial x_{k}\partial x_{l}}(\nabla
\log\Phi_{P})_{i}](y)\qquad\ \qquad$

$=-[\frac{\partial^{3}X_{l}}{\partial x_{i}\partial x_{j}\partial x_{k}}%
+\frac{\partial^{3}X_{k}}{\partial x_{l}\partial x_{i}\partial x_{j}}%
+\frac{\partial^{3}X_{j}}{\partial x_{k}\partial x_{l}\partial x_{i}}%
+\frac{\partial^{3}X_{i}}{\partial x_{j}\partial x_{k}\partial x_{l}}](y)$

In \textbf{particular},

$\qquad\lbrack\frac{\partial^{3}}{\partial x_{i}^{2}\partial x_{k}}(\nabla
\log\Phi_{P})_{j}+\frac{\partial^{3}}{\partial x_{j}\partial x_{i}^{2}}%
(\nabla\log\Phi_{P})_{k}+\frac{\partial^{3}}{\partial x_{k}\partial
x_{j}\partial x_{i}}(\nabla\log\Phi_{P})_{i}+\frac{\partial^{3}}{\partial
x_{i}\partial x_{k}\partial x_{j}}(\nabla\log\Phi_{P})_{i}](y)$

\qquad$=-[\frac{\partial^{3}X_{j}}{\partial x_{i}^{2}\partial x_{k}}%
+\frac{\partial^{3}X_{k}}{\partial x_{j}\partial x_{i}^{2}}+\frac{\partial
^{3}X_{i}}{\partial x_{k}\partial x_{j}\partial x_{i}}+\frac{\partial^{3}%
X_{i}}{\partial x_{i}\partial x_{k}\partial x_{j}}](y)=-[\frac{\partial
^{3}X_{j}}{\partial x_{i}^{2}\partial x_{k}}+\frac{\partial^{3}X_{k}}{\partial
x_{i}^{2}\partial x_{j}}+2\frac{\partial^{3}X_{i}}{\partial x_{i}\partial
x_{j}\partial x_{k}}](y)$

In \textbf{particular} for $k=j,$ we have:

$[\frac{\partial^{3}}{\partial x_{i}^{2}\partial x_{j}}(\nabla\log\Phi
_{P})_{j}+\frac{\partial^{3}}{\partial x_{i}\partial x_{j}^{2}}(\nabla
$log$\Phi_{P})_{i}](y)=-\left(  \frac{\partial^{3}X_{j}}{\partial x_{i}%
^{2}\partial x_{j}}+\frac{\partial^{3}X_{i}}{\partial x_{i}\partial x_{j}^{2}%
}\right)  (y)$

(viii)$^{\ast\ast\ast}\qquad\frac{\partial^{4}}{\partial x_{i}\partial
x_{j}\partial x_{k}\partial x_{l}}(\nabla\log\Phi_{P})_{r}+\frac{\partial^{4}%
}{\partial x_{i}\partial x_{j}\partial x_{k}\partial x_{r}}(\nabla$%
log$\Phi_{P})_{l}+\frac{\partial^{4}}{\partial x_{i}\partial x_{j}\partial
x_{l}\partial x_{r}}(\nabla\log\Phi_{P})_{k}$

$\qquad+\frac{\partial^{4}}{\partial x_{i}\partial x_{k}\partial x_{l}\partial
x_{r}}(\nabla\log\Phi_{P})_{j}+\frac{\partial^{4}}{\partial x_{j}\partial
x_{k}\partial x_{l}\partial x_{r}}(\nabla\log\Phi_{P})_{i}\qquad\ \qquad
\qquad\qquad$

$\qquad+\frac{\partial^{4}X_{r}}{\partial x_{i}\partial x_{j}\partial
x_{k}\partial x_{l}}+\frac{\partial^{4}X_{l}}{\partial x_{i}\partial
x_{j}\partial x_{k}\partial x_{r}}+\frac{\partial^{4}X_{k}}{\partial
x_{i}\partial x_{j}\partial x_{l}\partial x_{r}}+\frac{\partial^{4}X_{j}%
}{\partial x_{i}\partial x_{k}\partial x_{l}\partial x_{r}}+\frac{\partial
^{4}X_{i}}{\partial x_{j}\partial x_{k}\partial x_{l}\partial x_{r}}$

$=-\underset{m=q+1}{\overset{n}{\sum}}x_{m}\left(  \frac{\partial^{5}%
}{\partial x_{i}\partial x_{j}\partial x_{k}\partial x_{l}\partial x_{r}%
}(\nabla\log\Phi_{P})_{m}+\frac{\partial^{5}X_{m}}{\partial x_{i}\partial
x_{j}\partial x_{k}\partial x_{l}\partial x_{r}}\right)  $

\qquad In particular, when $k=i$ and $l=j,$ we have:

(viii)$^{\ast\ast\ast\ast}\qquad\frac{\partial^{4}}{\partial x_{i}^{2}\partial
x_{j}^{2}}(\nabla\log\Phi_{P})_{r}+\frac{\partial^{4}}{\partial x_{i}%
^{2}\partial x_{j}\partial x_{r}}(\nabla$log$\Phi_{P})_{j}+\frac{\partial^{4}%
}{\partial x_{i}^{2}\partial x_{j}\partial x_{r}}(\nabla\log\Phi_{P})_{j}$

$\qquad+\frac{\partial^{4}}{\partial x_{i}\partial x_{j}^{2}\partial x_{r}%
}(\nabla\log\Phi_{P})_{i}+\frac{\partial^{4}}{\partial x_{j}^{2}\partial
x_{k}\partial x_{r}}(\nabla\log\Phi_{P})_{i}\qquad\ \qquad\qquad\qquad$

$\qquad+\frac{\partial^{4}X_{r}}{\partial x_{i}^{2}\partial x_{j}^{2}}%
+\frac{\partial^{4}X_{j}}{\partial x_{i}^{2}\partial x_{j}\partial x_{r}%
}+\frac{\partial^{4}X_{j}}{\partial x_{i}^{2}\partial x_{j}\partial x_{r}%
}+\frac{\partial^{4}X_{i}}{\partial x_{i}\partial x_{j}^{2}\partial x_{r}%
}+\frac{\partial^{4}X_{i}}{\partial x_{i}\partial x_{j}^{2}\partial x_{r}}$

$\qquad=-\underset{m=q+1}{\overset{n}{\sum}}x_{m}\left(  \frac{\partial^{5}%
}{\partial x_{i}^{2}\partial x_{j}^{2}\partial x_{r}}(\nabla\log\Phi_{P}%
)_{m}+\frac{\partial^{5}X_{m}}{\partial x_{i}^{2}\partial x_{j}^{2}\partial
x_{r}}\right)  $

(viii)$^{\ast\ast\ast\ast\ast}\qquad$We have at all points of M$_{0},$

$\frac{\partial^{5}}{\partial x_{i}\partial x_{j}\partial x_{k}\partial
x_{l}\partial x_{q}}(\nabla\log\Phi_{P})_{r}+\frac{\partial^{5}}{\partial
x_{i}\partial x_{j}\partial x_{k}\partial x_{l}\partial x_{r}}(\nabla\log
\Phi_{P})_{q}+\frac{\partial^{4}}{\partial x_{i}\partial x_{j}\partial
x_{k}\partial x_{q}\partial x_{r}}(\nabla\log\Phi_{P})_{l}$

$+\frac{\partial^{5}}{\partial x_{i}\partial x_{j}\partial x_{l}\partial
x_{q}\partial x_{r}}(\nabla\log\Phi_{P})_{k}+\frac{\partial^{5}}{\partial
x_{i}\partial x_{k}\partial x_{l}\partial x_{q}\partial x_{r}}(\nabla\log
\Phi_{P})_{j}+\frac{\partial^{5}}{\partial x_{j}\partial x_{k}\partial
x_{l}\partial x_{q}\partial x_{r}}(\nabla\log\Phi_{P})_{i}\qquad\ \qquad
\qquad\qquad$

$+\frac{\partial^{5}X_{r}}{\partial x_{i}\partial x_{j}\partial x_{k}\partial
x_{l}\partial x_{q}}+\frac{\partial^{5}X_{q}}{\partial x_{i}\partial
x_{j}\partial x_{k}\partial x_{l}\partial x_{r}}+\frac{\partial^{5}X_{l}%
}{\partial x_{i}\partial x_{j}\partial x_{k}\partial x_{q}\partial x_{r}%
}+\frac{\partial^{5}X_{k}}{\partial x_{i}\partial x_{j}\partial x_{l}\partial
x_{q}\partial x_{r}}+\frac{\partial^{5}X_{j}}{\partial x_{i}\partial
x_{k}\partial x_{l}\partial x_{q}\partial x_{r}}$

$+\frac{\partial^{5}X_{i}}{\partial x_{j}\partial x_{k}\partial x_{l}\partial
x_{q}\partial x_{r}}=-\underset{m=q+1}{\overset{n}{\sum}}x_{m}\left(
\frac{\partial^{6}}{\partial x_{i}\partial x_{j}\partial x_{k}\partial
x_{l}\partial x_{q}\partial x_{r}}(\nabla\log\Phi_{P})_{m}+\frac{\partial
^{6}X_{m}}{\partial x_{i}\partial x_{j}\partial x_{k}\partial x_{l}\partial
x_{q}\partial x_{r}}\right)  $

In particuar for $l=i;q=j;r=k,$ we have at all points of M$_{0},$

(viii)$^{\ast\ast\ast\ast\ast\ast}\frac{\partial^{5}}{\partial x_{i}%
^{2}\partial x_{j}^{2}\partial x_{k}}(\nabla\log\Phi_{P})_{k}+\frac
{\partial^{5}}{\partial x_{i}^{2}\partial x_{j}\partial x_{k}^{2}}(\nabla
\log\Phi_{P})_{j}+\frac{\partial^{4}}{\partial x_{i}\partial x_{j}^{2}\partial
x_{k}^{2}}(\nabla\log\Phi_{P})_{i}$

$+\frac{\partial^{5}}{\partial x_{i}^{2}\partial x_{j}^{2}\partial x_{k}%
}(\nabla\log\Phi_{P})_{k}+\frac{\partial^{5}}{\partial x_{i}^{2}\partial
x_{j}\partial x_{k}^{2}}(\nabla\log\Phi_{P})_{j}+\frac{\partial^{5}}{\partial
x_{j}^{2}\partial x_{k}^{2}\partial x_{i}}(\nabla\log\Phi_{P})_{i}%
\qquad\ \qquad\qquad\qquad$

$+\frac{\partial^{5}X_{k}}{\partial x_{i}^{2}\partial x_{j}^{2}\partial x_{k}%
}+\frac{\partial^{5}X_{j}}{\partial x_{i}^{2}\partial x_{j}\partial x_{k}^{2}%
}+\frac{\partial^{5}X_{i}}{\partial x_{i}\partial x_{j}^{2}\partial x_{k}^{2}%
}+\frac{\partial^{5}X_{k}}{\partial x_{i}^{2}\partial x_{j}^{2}\partial x_{k}%
}+\frac{\partial^{5}X_{j}}{\partial x_{i}^{2}\partial x_{k}^{2}\partial x_{j}%
}+\frac{\partial^{5}X_{i}}{\partial x_{i}\partial x_{j}^{2}\partial x_{k}^{2}%
}$

$=-\underset{m=q+1}{\overset{n}{\sum}}x_{m}\left(  \frac{\partial^{6}%
}{\partial x_{i}^{2}\partial x_{j}^{2}\partial x_{k}^{2}}(\nabla\log\Phi
_{P})_{m}+\frac{\partial^{6}X_{m}}{\partial x_{i}^{2}\partial x_{j}%
^{2}\partial x_{k}^{2}}\right)  $

Simplifying, we have:

$2\frac{\partial^{5}}{\partial x_{i}^{2}\partial x_{j}^{2}\partial x_{k}%
}(\nabla\log\Phi_{P})_{k}+2\frac{\partial^{5}}{\partial x_{i}^{2}\partial
x_{j}\partial x_{k}^{2}}(\nabla\log\Phi_{P})_{j}+2\frac{\partial^{4}}{\partial
x_{i}\partial x_{j}^{2}\partial x_{k}^{2}}(\nabla\log\Phi_{P})_{i}\qquad
\qquad\qquad$

$+2\frac{\partial^{5}X_{k}}{\partial x_{i}^{2}\partial x_{j}^{2}\partial
x_{k}}+2\frac{\partial^{5}X_{j}}{\partial x_{i}^{2}\partial x_{j}\partial
x_{k}^{2}}+2\frac{\partial^{5}X_{i}}{\partial x_{i}\partial x_{j}^{2}\partial
x_{k}^{2}}$

$=-\underset{m=q+1}{\overset{n}{\sum}}x_{m}\left(  \frac{\partial^{6}%
}{\partial x_{i}^{2}\partial x_{j}^{2}\partial x_{k}^{2}}(\nabla\log\Phi
_{P})_{m}+\frac{\partial^{6}X_{m}}{\partial x_{i}^{2}\partial x_{j}%
^{2}\partial x_{k}^{2}}\right)  $

In particular, we have:

$[\frac{\partial^{5}}{\partial x_{i}^{2}\partial x_{j}^{2}\partial x_{k}%
}(\nabla\log\Phi_{P})_{k}+\frac{\partial^{5}}{\partial x_{i}^{2}\partial
x_{j}\partial x_{k}^{2}}(\nabla\log\Phi_{P})_{j}+\frac{\partial^{5}}{\partial
x_{i}\partial x_{j}^{2}\partial x_{k}^{2}}(\nabla\log\Phi_{P})_{i}%
](y)\qquad\qquad\qquad$

$-[\frac{\partial^{5}X_{k}}{\partial x_{i}^{2}\partial x_{j}^{2}\partial
x_{k}}+\frac{\partial^{5}X_{j}}{\partial x_{i}^{2}\partial x_{j}\partial
x_{k}^{2}}+\frac{\partial^{5}X_{i}}{\partial x_{i}\partial x_{j}^{2}\partial
x_{k}^{2}}](y)$

\subsection{\textbf{Tangential Derivatives:}}

(ix) For a $=$ 1,...,q and for $j=q+1,...,n,$

$\qquad\frac{\partial}{\partial x_{\text{a}}}(\nabla$log$\Phi_{P}%
)_{j}(y)=-\frac{\partial X_{j}}{\partial x_{\text{a}}}(y)\qquad$\qquad

(x) For a,b = 1,...,q,

\qquad$\frac{\partial^{2}}{\partial x_{\text{a}}\partial x_{\text{b}}}(\nabla
$log$\Phi_{P})_{j}(y)=-\frac{\partial^{2}X_{j}}{\partial x_{\text{a}}\partial
x_{\text{b}}}(y)$

Fomulae for higher derivatives follow.

(xi) For a = 1,...,q \ and $y\in U\subset P,$ we have:

\qquad$(\nabla$log$\Phi_{P})_{\text{a}}(y)=0$

(xii) For a, b = 1,...,q$\ $

$\qquad\frac{\partial}{\partial x_{\text{b}}}(\nabla$log$\Phi_{P})_{\text{a}%
}(y)\ =0$

(xiii) For a, b, c = 1,...,q,

\qquad$\ \frac{\partial^{2}}{\partial x_{\text{c}}\partial x_{\text{b}}%
}(\nabla$log$\Phi_{P})_{\text{a}}(y)\ =0\qquad$

\subsection{\textbf{Mixed Derivatives: }}

For a =1,...,q and\textbf{ }$i,j,k=q+1,...,n:$

(xiv) $\frac{\partial^{2}}{\partial x_{\text{a}}\partial x_{k}}\nabla\log
\Phi_{P})_{j}(y)+\frac{\partial^{2}}{\partial x_{\text{a}}\partial x_{j}%
}\nabla\log\Phi_{P})_{k}(y)=-\frac{\partial^{2}X_{j}}{\partial x_{\text{a}%
}\partial x_{k}}(y)-$ $\frac{\partial^{2}X_{k}}{\partial x_{\text{a}}\partial
x_{j}}(y).$

In particular for $k=j,$

\qquad$\frac{\partial^{2}}{\partial x_{\text{a}}\partial x_{j}}\nabla\log
\Phi_{P})_{j}(y)=-\frac{\partial^{2}X_{j}}{\partial x_{\text{a}}\partial
x_{j}}(y)$

(xv) $\frac{\partial}{\partial x_{j}}(\nabla$log$\Phi_{P})_{\text{a}}(y)=$
$\underset{i=q+1}{\overset{n}{%
{\textstyle\sum}
}}X_{i}(y)\perp_{\text{a}ij}(y)-\frac{\partial X_{j}}{\partial x_{\text{a}}%
}(y)$

(xvi) $\frac{\partial^{2}}{\partial x_{i}\partial x_{j}}(\nabla$log$\Phi
_{P})_{\text{a}}(y)$

\qquad$=-$ $2\underset{\text{b=1}}{\overset{\text{q}}{%
{\textstyle\sum}
}}$T$_{\text{ab}j}(y_{0})\frac{\partial X_{i}}{\partial x_{\text{b}}%
}(y)-2\underset{\text{b=1}}{\overset{\text{q}}{%
{\textstyle\sum}
}}T_{\text{ab}i}(y)\frac{\partial X_{j}}{\partial x_{\text{b}}}(y)$

$+\frac{1}{2}\underset{k=q+1}{\overset{n}{%
{\textstyle\sum}
}}\perp_{\text{a}jk}(y)\left[  \left(  \frac{\partial X_{i}}{\partial x_{k}%
}+\frac{\partial X_{k}}{\partial x_{i}}\right)  \right]  (y)+\frac{1}%
{2}\underset{k=q+1}{\overset{n}{%
{\textstyle\sum}
}}\perp_{\text{a}ik}(y)[\left(  \frac{\partial X_{k}}{\partial x_{j}}%
+\frac{\partial X_{j}}{\partial x_{k}}\right)  ](y)$\qquad

$+\frac{4}{3}\underset{k=q+1}{\overset{n}{%
{\textstyle\sum}
}}\left[  R_{i\text{a}jk}+R_{j\text{a}ik}\right]  (y)X_{k}(y)+[X_{i}%
\frac{\partial X_{j}}{\partial x_{\text{a}}}+X_{j}\frac{\partial X_{i}%
}{\partial x_{\text{a}}}-\frac{1}{2}\left(  \frac{\partial^{2}X_{i}}{\partial
x_{\text{a}}\partial x_{j}}+\frac{\partial^{2}X_{j}}{\partial x_{\text{a}%
}\partial x_{i}}\right)  ](y)$\qquad\qquad

In particular, taking $j=i,$ we have:

$\frac{\partial^{2}}{\partial x_{i}^{2}}(\nabla$log$\Phi_{P})_{\text{a}%
}(y)=-4\underset{\text{b=1}}{\overset{\text{q}}{%
{\textstyle\sum}
}}T_{\text{ab}i}(y)\frac{\partial X_{i}}{\partial x_{\text{b}}}(y)$

$+\underset{k=q+1}{\overset{n}{%
{\textstyle\sum}
}}\perp_{\text{a}ik}(y)\left[  \left(  \frac{\partial X_{i}}{\partial x_{k}%
}+\frac{\partial X_{k}}{\partial x_{i}}\right)  \right]  (y)+\frac{8}%
{3}\underset{k=q+1}{\overset{n}{%
{\textstyle\sum}
}}R_{i\text{a}ik}(y)X_{k}(y)$

$+[2X_{i}\frac{\partial X_{i}}{\partial x_{\text{a}}}-\frac{\partial^{2}X_{i}%
}{\partial x_{\text{a}}\partial x_{i}}](y)$

\subsection{\textbf{Computations of B}$_{1}$}

\subsubsection{Normal Derivatives}

(i) Recalling that by the definition in $\left(  1.5\right)  $ in Chapter 1,

$\left(  B_{1}\right)  \qquad\Phi_{P}(x)=$ exp$\left\{  \int_{0}^{1}%
<\text{X(}\gamma\text{(s)) , }\dot{\gamma}\text{(s)%
$>$%
ds}\right\}  \qquad\qquad\qquad\qquad\qquad$

where $\gamma$ is the unique minimal geodesic from x to P in time 1 and
meeting P orthogonally at a point y$\in$P:

The geodesic $\gamma:[0,1]\longrightarrow$ M$_{0}$ is given in Fermi
coordinates as:

$\qquad\gamma(s)=\left(  \text{x}_{1},...,\text{x}_{q},(1-s)\text{x}%
_{q+1},...,(1-s)\text{x}_{n}\right)  $

Consequently.

\qquad$\dot{\gamma}$(s) $=\left(  \text{0},...,\text{0},-\text{x}%
_{q+1},...,-\text{x}_{n}\right)  =-$ $\underset{j=q+1}{\overset{n}{\sum}}%
x_{j}\frac{\partial}{\partial x_{j}}\mid_{\gamma(s)}$

By \textbf{Definition, }p. 22 of \textbf{Gray }$\left[  4\right]  ,$ we set:

\qquad$\sigma^{2}=\underset{j=q+1}{\overset{n}{\sum}}x_{j}^{2}$ and
$X=\underset{j=1}{\overset{n}{\sum}}X_{j}\frac{\partial}{\partial x_{j}}$

Then,

$\ \sigma<\nabla\sigma,X>$ $\ =$ $\sigma X(\sigma)=\frac{1}{2}X(\sigma
^{2})=\frac{1}{2}\underset{i=q+1}{\overset{n}{\sum}}X(x_{i}^{2})$

$=\frac{1}{2}\underset{j=1}{\overset{n}{\sum}}%
\underset{i=q+1}{\overset{n}{\sum}}X_{j}\frac{\partial}{\partial x_{j}}%
(x_{i}^{2})\qquad\left(  1\right)  $

$=$ $\underset{j=1}{\overset{n}{\sum}}\underset{i=q+1}{\overset{n}{\sum}}%
x_{i}\not X  _{j}\frac{\partial x_{i}}{\partial x_{j}}=$
$\underset{j=1}{\overset{n}{\sum}}\underset{i=q+1}{\overset{n}{\sum}}%
x_{i}\not X  _{j}\delta_{ij}=$ $\underset{i=q+1}{\overset{n}{\sum}}%
x_{i}\not X  _{i}$

We thus have the formula:$\qquad\qquad$

$\left(  B_{2}\right)  $\qquad$\ \ \sigma<\nabla\sigma,X>$ $=$
$\underset{i=q+1}{\overset{n}{\sum}}x_{i}X_{i}\qquad\qquad\qquad\left(
2\right)  \qquad\qquad\qquad\qquad\qquad\qquad\qquad\qquad\qquad\qquad
\qquad\qquad\qquad\qquad\qquad\qquad\qquad\qquad\qquad\qquad\qquad\qquad$

\ \ By $\left(  3.21\right)  $ of \textbf{Ndumu }$\left[  3\right]  ,$ we have
for a general smooth vector field X:\qquad\qquad\qquad$\qquad\qquad$

$\left(  B_{3}\right)  $\qquad$\qquad<\nabla\sigma,X>$ $=-<\nabla\sigma
,\nabla$log$\Phi_{P}>\qquad\qquad\qquad\qquad\qquad\qquad\qquad\ \ \ \ \qquad
\qquad\qquad\qquad\qquad\qquad\qquad\qquad\qquad\qquad\qquad\qquad\qquad
\qquad\ \ $

Therefore by $\left(  B_{2}\right)  $ and $\left(  B_{3}\right)  $, we have:

$\underset{i=q+1}{\overset{n}{\sum}}x_{i}X_{i}=\sigma<\nabla\sigma,X>$
$=-\sigma<\nabla\sigma,\nabla$log$\Phi_{P}>$
$=-\underset{i=q+1}{\overset{n}{\sum}}x_{i}(\nabla\log\Phi_{P})_{i}%
\qquad\qquad\qquad\ \ \ \ \qquad\qquad\qquad\qquad\qquad\qquad\ \ \qquad
\qquad$

The first and last equalities above give:$\qquad$

$\left(  B_{4}\right)  \qquad\underset{i=q+1}{\overset{n}{\sum}}x_{i}X_{i}=-$
$\underset{i=q+1}{\overset{n}{\sum}}x_{i}(\nabla$log$\Phi_{P})_{i}\qquad
\qquad\qquad\qquad\qquad\qquad\qquad\qquad\qquad\qquad\qquad\qquad\qquad
\qquad\qquad\qquad\ \ \ \ \ \qquad$

$\qquad\qquad\ \ \ \ \ \qquad\ \ \ \ \ \ \qquad\qquad\qquad\qquad\qquad
\qquad\qquad$

Differentiating both sides of $\left(  B_{4}\right)  $ above, we have for
$j=1,...,q,q+1,...,n:$

\ $\underset{i=q+1}{\overset{n}{\sum}}\frac{\partial x_{i}}{\partial x_{j}%
}X_{i}+\underset{i=q+1}{\overset{n}{\sum}}x_{i}\frac{\partial X_{i}}{\partial
x_{j}}=-$ $\underset{i=q+1}{\overset{n}{\sum}}\frac{\partial x_{i}}{\partial
x_{j}}(\nabla\log\Phi_{P})_{i}-\underset{i=q+1}{\overset{n}{\sum}}x_{i}%
\frac{\partial}{\partial x_{j}}(\nabla\log\Phi_{P})_{i}$\qquad

Re-arranging the above equation, we have for $j=1,...,q,q+1,...,n$

$\underset{i=q+1}{\overset{n}{\sum}}x_{i}\left(  \frac{\partial X_{i}%
}{\partial x_{j}}+\frac{\partial}{\partial x_{j}}\nabla\log\Phi_{P}%
)_{i}\right)  =-$\ $\underset{i=q+1}{\overset{n}{\sum}}\frac{\partial x_{i}%
}{\partial x_{j}}\left(  X_{i}+(\nabla\log\Phi_{P})_{i}\right)  $

\qquad\qquad\qquad\qquad\qquad\qquad\qquad$=-$%
\ $\underset{i=q+1}{\overset{n}{\sum}}\delta_{ij}\left(  X_{i}+(\nabla\log
\Phi_{P})_{i}\right)  $

\qquad$\qquad\qquad\qquad=-\left\{
\begin{array}
[c]{c}%
0\text{ for }j=1,...,q\\
X_{j}+(\nabla\log\Phi_{P})_{j}\text{ for }j=q+1,...,n
\end{array}
\right.  $

We can re-write the above equation in two separate equations (on M$_{0}$) as follows:

For a =1,...,q and for $i,j,k,l=q+1,...,n,$

$\left(  B_{5}\right)  \qquad\underset{i=q+1}{\overset{n}{\sum}}x_{i}\left(
\frac{\partial}{\partial x_{\text{a}}}\nabla\log\Phi_{P})_{i}+\frac{\partial
X_{i}}{\partial x_{\text{a}}}\right)  =0$ $\qquad\qquad\qquad\qquad
\qquad\qquad\qquad\qquad\qquad\ \ \ \ \ \qquad$

$\left(  B_{6}\right)  \qquad(\nabla$log$\Phi_{P})_{j}+X_{j}%
=-\underset{i=q+1}{\text{ }\overset{n}{\sum}}x_{i}\left(  \frac{\partial
}{\partial x_{j}}(\nabla\log\Phi_{P})_{i}+\frac{\partial X_{i}}{\partial
x_{j}}\right)  \qquad\left(  3\right)  \qquad\qquad\qquad$

(ii) Changing indices on the RHS, we can re-write $\left(  B_{6}\right)  $ as:

\qquad\qquad$(\nabla$log$\Phi_{P})_{j}+X_{j}=-\underset{k=q+1}{\text{
}\overset{n}{\sum}}x_{k}\left(  \frac{\partial}{\partial x_{j}}(\nabla\log
\Phi_{P})_{k}+\frac{\partial X_{k}}{\partial x_{j}}\right)  $

Differentiating on both sides of the last equation above gives:

\qquad$\frac{\partial}{\partial x_{i}}(\nabla$log$\Phi_{P})_{j}+\frac{\partial
X_{j}}{\partial x_{i}}=-\underset{k=q+1}{\text{ }\overset{n}{\sum}}%
\frac{\partial x_{k}}{\partial x_{i}}\left(  \frac{\partial}{\partial x_{j}%
}(\nabla\log\Phi_{P})_{k}+\frac{\partial X_{k}}{\partial x_{j}}\right)  $

$\qquad\qquad\qquad\qquad\qquad\qquad-\underset{k=q+1}{\text{ }%
\overset{n}{\sum}}x_{k}\left(  \frac{\partial^{2}}{\partial x_{i}\partial
x_{j}}(\nabla\log\Phi_{P})_{k}+\frac{\partial^{2}X_{k}}{\partial x_{i}\partial
x_{j}}\right)  $

Since $\frac{\partial x_{k}}{\partial x_{i}}=\delta_{ik},$ the last equation
above gives:

$\left(  B_{7}\right)  \qquad\frac{\partial}{\partial x_{i}}(\nabla$%
log$\Phi_{P})_{j}+\frac{\partial}{\partial x_{j}}(\nabla\log\Phi_{P}%
)_{i}+\frac{\partial X_{j}}{\partial x_{i}}+\frac{\partial X_{i}}{\partial
x_{j}}\qquad\qquad\qquad\qquad\left(  4\right)  \ $

$\qquad\qquad=-\underset{k=q+1}{\overset{n}{\sum}}x_{k}\left(  \frac
{\partial^{2}}{\partial x_{i}\partial x_{j}}(\nabla\log\Phi_{P})_{k}%
+\frac{\partial^{2}X_{k}}{\partial x_{i}\partial x_{j}}\right)  \qquad
\ \ \ \ $

(iii) Further differentiating, we have:

$\left(  B_{8}\right)  \qquad\frac{\partial^{2}}{\partial x_{i}\partial x_{j}%
}(\nabla$log$\Phi_{P})_{k}+\frac{\partial^{2}}{\partial x_{k}\partial x_{i}%
}(\nabla\log\Phi_{P})_{j}+\frac{\partial^{2}}{\partial x_{j}\partial x_{k}%
}(\nabla\log\Phi_{P})_{i}\qquad\left(  5\right)  $

$\qquad\qquad+\frac{\partial^{2}X_{k}}{\partial x_{i}\partial x_{j}}%
+\frac{\partial^{2}X_{j}}{\partial x_{k}\partial x_{i}}+\frac{\partial
^{2}X_{i}}{\partial x_{j}\partial x_{k}}\qquad\qquad\qquad
\ \ \ \ \ \ \ \ \ \qquad$

$\qquad\qquad=-\underset{m=q+1}{\overset{n}{\sum}}x_{m}\left(  \frac
{\partial^{3}}{\partial x_{i}\partial x_{j}\partial x_{k}}(\nabla\log\Phi
_{P})_{m}+\frac{\partial^{3}X_{m}}{\partial x_{i}\partial x_{j}\partial x_{k}%
}\right)  \qquad$

(iii)$^{\ast}$ Another further differentiation gives:

$\left(  B_{9}\right)  \qquad\frac{\partial^{3}}{\partial x_{i}\partial
x_{j}\partial x_{k}}(\nabla$log$\Phi_{P})_{l}+\frac{\partial^{3}}{\partial
x_{l}\partial x_{i}\partial x_{j}}(\nabla\log\Phi_{P})_{k}$

$\qquad\qquad+\frac{\partial^{3}}{\partial x_{k}\partial x_{l}\partial x_{i}%
}(\nabla\log\Phi_{P})_{j}+\frac{\partial^{3}}{\partial x_{j}\partial
x_{k}\partial x_{l}}(\nabla\log\Phi_{P})_{i}\qquad\ \qquad$

$\qquad+\frac{\partial^{3}X_{l}}{\partial x_{i}\partial x_{j}\partial x_{k}%
}+\frac{\partial^{3}X_{k}}{\partial x_{l}\partial x_{i}\partial x_{j}}%
+\frac{\partial^{3}X_{j}}{\partial x_{k}\partial x_{l}\partial x_{i}}%
+\frac{\partial^{3}X_{i}}{\partial x_{j}\partial x_{k}\partial x_{l}}$

$\qquad=-\underset{m=q+1}{\overset{n}{\sum}}x_{m}\left(  \frac{\partial^{4}%
}{\partial x_{i}\partial x_{j}\partial x_{k}\partial x_{l}}(\nabla\log\Phi
_{P})_{m}+\frac{\partial^{4}X_{m}}{\partial x_{i}\partial x_{j}\partial
x_{k}\partial x_{l}}\right)  $

\qquad\qquad\qquad\qquad\qquad\qquad\qquad\qquad\qquad\qquad\qquad\qquad
\qquad\qquad\qquad\qquad\qquad\qquad$\blacksquare$

$\left(  B_{9}\right)  ^{\ast}\qquad\frac{\partial^{4}}{\partial x_{i}\partial
x_{j}\partial x_{k}\partial x_{l}}(\nabla\log\Phi_{P})_{r}+\frac{\partial^{4}%
}{\partial x_{i}\partial x_{j}\partial x_{k}\partial x_{r}}(\nabla$%
log$\Phi_{P})_{l}$

$\qquad+\frac{\partial^{4}}{\partial x_{l}\partial x_{i}\partial x_{j}\partial
x_{r}}(\nabla\log\Phi_{P})_{k}+\frac{\partial^{4}}{\partial x_{k}\partial
x_{l}\partial x_{i}\partial x_{r}}(\nabla\log\Phi_{P})_{j}+\frac{\partial^{4}%
}{\partial x_{j}\partial x_{k}\partial x_{l}\partial x_{r}}(\nabla\log\Phi
_{P})_{i}\qquad\ \qquad\qquad\qquad$

$\qquad+\frac{\partial^{4}X_{r}}{\partial x_{j}\partial x_{k}\partial
x_{l}\partial x_{i}}+\frac{\partial^{4}X_{l}}{\partial x_{i}\partial
x_{j}\partial x_{k}\partial x_{r}}+\frac{\partial^{4}X_{k}}{\partial
x_{l}\partial x_{i}\partial x_{j}\partial x_{r}}+\frac{\partial^{4}X_{j}%
}{\partial x_{k}\partial x_{l}\partial x_{i}\partial x_{r}}+\frac{\partial
^{4}X_{i}}{\partial x_{j}\partial x_{k}\partial x_{l}\partial x_{r}}$

$\qquad=-\underset{m=q+1}{\overset{n}{\sum}}x_{m}\left(  \frac{\partial^{5}%
}{\partial x_{i}\partial x_{j}\partial x_{k}\partial x_{l}\partial x_{r}%
}(\nabla\log\Phi_{P})_{m}+\frac{\partial^{5}X_{m}}{\partial x_{i}\partial
x_{j}\partial x_{k}\partial x_{l}\partial x_{r}}\right)  $

\qquad In particular, when $k=i$ and $l=j,$ we have:

$\left(  B_{9}\right)  ^{\ast\ast}\qquad\frac{\partial^{4}}{\partial x_{i}%
^{2}\partial x_{j}^{2}}(\nabla\log\Phi_{P})_{r}+\frac{\partial^{4}}{\partial
x_{i}^{2}\partial x_{j}\partial x_{r}}(\nabla$log$\Phi_{P})_{j}+\frac
{\partial^{4}}{\partial x_{i}\partial x_{j}^{2}\partial x_{r}}(\nabla\log
\Phi_{P})_{i}$

$\qquad\qquad+\frac{\partial^{4}}{\partial x_{j}\partial x_{i}^{2}\partial
x_{r}}(\nabla\log\Phi_{P})_{j}+\frac{\partial^{4}}{\partial x_{j}^{2}\partial
x_{i}\partial x_{r}}(\nabla\log\Phi_{P})_{i}\qquad\ \qquad\qquad\qquad$

$\qquad\qquad+\frac{\partial^{4}X_{r}}{\partial x_{i}^{2}\partial x_{j}^{2}%
}+\frac{\partial^{4}X_{j}}{\partial x_{i}^{2}\partial x_{j}\partial x_{r}%
}+\frac{\partial^{4}X_{i}}{\partial x_{i}\partial^{2}x_{j}\partial x_{r}%
}+\frac{\partial^{4}X_{j}}{\partial x_{i}^{2}\partial x_{j}\partial x_{r}%
}+\frac{\partial^{4}X_{i}}{\partial x_{i}\partial x_{j}^{2}\partial x_{r}}$

$\qquad\qquad=-\underset{m=q+1}{\overset{n}{\sum}}x_{m}\left(  \frac
{\partial^{5}}{\partial x_{i}^{2}\partial x_{j}^{2}\partial x_{r}}(\nabla
\log\Phi_{P})_{m}+\frac{\partial^{5}X_{m}}{\partial x_{i}^{2}\partial
x_{j}^{2}\partial x_{r}}\right)  $

$\qquad\ \ $

Higher derivatives follow.$\ \ \qquad$

The above are general formulae relating a general vector field X and
$\nabla\log\Phi_{P}$ and their

derivatives with respect \textbf{normal Fermi} \textbf{coordinates} in the
tubular neighbourhood M$_{0}$ of P.

\qquad\qquad\qquad\qquad\qquad\qquad\qquad\qquad\qquad\qquad\qquad\qquad
\qquad\qquad\qquad\qquad\qquad$\blacksquare$

\subsubsection{\textbf{Tangential Derivatives}:}

For a,b =1,...,q and $i,j,k=q+1,...,n:$

(iv) We differentiate both sides of $\left(  B_{5}\right)  :$

\qquad\ \ $\underset{i=q+1}{\overset{n}{\sum}}x_{i}\left(  \frac{\partial
}{\partial x_{\text{a}}}\nabla\log\Phi_{P})_{i}+\frac{\partial X_{i}}{\partial
x_{\text{a}}}\right)  =0$

to have for a =1,...,q and $j=1,...,q,q+1,...,n:$

\qquad\ $\underset{i=q+1}{\overset{n}{\sum}}\frac{\partial x_{i}}{\partial
x_{j}}\left(  \frac{\partial}{\partial x_{\text{a}}}\nabla\log\Phi_{P}%
)_{i}+\frac{\partial X_{i}}{\partial x_{\text{a}}}\right)
+\underset{i=q+1}{\overset{n}{\sum}}x_{i}\left(  \frac{\partial^{2}}{\partial
x_{\text{a}}\partial x_{j}}\nabla\log\Phi_{P})_{i}+\frac{\partial^{2}X_{i}%
}{\partial x_{\text{a}}\partial x_{j}}\right)  =0$

Since $\frac{\partial x_{i}}{\partial x_{j}}=\delta_{ij},$ the last equation
above becomes for $j=q+1,...,n,$

\qquad$\left(  \frac{\partial}{\partial x_{\text{a}}}\nabla\log\Phi_{P}%
)_{j}+\frac{\partial X_{j}}{\partial x_{\text{a}}}\right)
=-\underset{i=q+1}{\overset{n}{\sum}}x_{i}\left(  \frac{\partial^{2}}{\partial
x_{\text{a}}\partial x_{j}}\nabla\log\Phi_{P})_{i}+\frac{\partial^{2}X_{i}%
}{\partial x_{\text{a}}\partial x_{j}}\right)  $

\qquad Alternatively, differentiate (with respect to tangential coordinates)
both sides of $\left(  B_{6}\right)  :$

$\qquad\frac{\partial}{\partial x_{\text{a}}}(\nabla$log$\Phi_{P})_{j}%
(x_{0})+\frac{\partial X_{j}}{\partial x_{\text{a}}}(x_{0})$

$\qquad=-\underset{i=q+1}{\overset{n}{\sum}}\frac{\partial x_{i}}{\partial
x_{\text{a}}}(x_{0})\left(  \frac{\partial}{\partial x_{j}}(\nabla\log\Phi
_{P})_{i}+\frac{\partial^{2}X_{i}}{\partial x_{j}}\right)  (x_{0})$

\qquad$-\underset{i=q+1}{\overset{n}{\sum}}x_{i}(x_{0})\left(  \frac
{\partial^{2}}{\partial x_{\text{a}}\partial x_{j}}(\nabla\log\Phi_{P}%
)_{i}+\frac{\partial^{2}X_{i}}{\partial x_{\text{a}}\partial x_{j}}\right)
(x_{0})$

Since $\frac{\partial x_{i}}{\partial x_{\text{a}}}=\delta_{i\text{a}}=0$ for
a = 1,...,q and $i=q+1,...,n,$ we have for $x_{0}\in M_{0}:$

$\frac{\partial}{\partial x_{\text{a}}}(\nabla$log$\Phi_{P})_{j}(x_{0}%
)+\frac{\partial X_{j}}{\partial x_{\text{a}}}(x_{0})=-$
$\underset{i=q+1}{\overset{n}{\sum}}x_{i}(x_{0})\left(  \frac{\partial^{2}%
}{\partial x_{\text{a}}\partial x_{j}}(\nabla\log\Phi_{P})_{i}+\frac
{\partial^{2}X_{i}}{\partial x_{\text{a}}\partial x_{j}}\right)  (x_{0})$

We conclude that for $y\in P,$ we have:

$\frac{\partial}{\partial x_{\text{a}}}(\nabla$log$\Phi_{P})_{j}%
(y)+\frac{\partial X_{j}}{\partial x_{\text{a}}}(y)=0$

(v) For a,b = 1,...,q and for $j=q+1,...,n,$

\qquad We repeat the process in (iv) and obtain:

$\frac{\partial^{2}}{\partial x_{\text{a}}\partial x_{\text{b}}}(\nabla
$log$\Phi_{P})_{j}(x_{0})+\frac{\partial^{2}X_{j}}{\partial x_{\text{a}%
}\partial x_{\text{b}}}(x_{0})$

$=-\underset{i=q+1}{\overset{n}{\sum}}x_{i}(x_{0})\left(  \frac{\partial^{2}%
}{\partial x_{\text{a}}\partial x_{\text{b}}\partial x_{j}}(\nabla\log\Phi
_{P})_{i}+\frac{\partial^{3}X_{i}}{\partial x_{\text{a}}\partial x_{\text{b}%
}\partial x_{j}}\right)  (x_{0})$

Recall that by the definition of Fermi coordinates, $x_{i}(y)=0$ for
$i=q+1,...,n$ for any $y\in$U$\subset$P where U is a small neighbourhood of
the centre of Fermi coordinates y$\in$P.

$\qquad\frac{\partial^{2}}{\partial x_{\text{a}}\partial x_{\text{b}}}(\nabla
$log$\Phi_{P})_{j}(y)+\frac{\partial^{2}X_{j}}{\partial x_{\text{a}}\partial
x_{\text{b}}}(y)=0$

In this case, the expressions in $\left(  B_{6}\right)  ,$ $\left(
B_{7}\right)  $ and $\left(  B_{8}\right)  $ become for $j,k,l=q+1,...,n:$

(vi) $\left(  B_{10}\right)  \qquad(\nabla\log\Phi_{P})_{j}(y)=-X_{j}(y)$ for
$j=q+1,...,n\qquad\qquad\qquad\qquad\qquad\ \left(  6\right)  \qquad
\qquad\qquad\qquad\qquad\ \ \ \qquad\qquad\ $

(vii) $\left(  B_{11}\right)  \qquad\frac{\partial}{\partial x_{i}}(\nabla
$log$\Phi_{P})_{j}(y)+\frac{\partial}{\partial x_{j}}(\nabla\log\Phi_{P}%
)_{i}(y)=-\left(  \frac{\partial X_{j}}{\partial x_{i}}+\frac{\partial X_{i}%
}{\partial x_{j}}\right)  (y)\qquad\ \ \ \ \left(  7\right)  $

We will see in $\left(  39\right)  $ that:

$\qquad\frac{\partial}{\partial x_{i}}(\nabla\log\Phi_{P})_{j}(y)=-\frac{1}%
{2}\left(  \frac{\partial X_{j}}{\partial x_{i}}+\frac{\partial X_{i}%
}{\partial x_{j}}\right)  (y)=\frac{\partial}{\partial x_{j}}(\nabla\log
\Phi_{P})_{i}(y)\qquad\qquad\qquad$\qquad\ 

(viii) $\left(  B_{12}\right)  $ $\frac{\partial^{2}}{\partial x_{i}\partial
x_{j}}(\nabla$log$\Phi_{P})_{k}(y)+\frac{\partial^{2}}{\partial x_{k}\partial
x_{i}}(\nabla\log\Phi_{P})_{j}(y)+\frac{\partial^{2}}{\partial x_{j}\partial
x_{k}}(\nabla\log\Phi_{P})_{i}(y)\qquad\left(  8\right)  $

$\qquad\qquad=-$ $[\frac{\partial^{2}X_{k}}{\partial x_{i}\partial x_{j}%
}(y)+\frac{\partial^{2}X_{j}}{\partial x_{k}\partial x_{i}}(y)+\frac
{\partial^{2}X_{i}}{\partial x_{j}\partial x_{k}}(y)]\qquad$

(viii)$^{\ast}$ Further we have from $\left(  B_{9}\right)  :$

$\left(  B_{12}\right)  ^{\ast}\qquad\lbrack\frac{\partial^{3}}{\partial
x_{i}\partial x_{j}\partial x_{k}}(\nabla$log$\Phi_{P})_{l}+\frac{\partial
^{3}}{\partial x_{l}\partial x_{i}\partial x_{j}}(\nabla\log\Phi_{P})_{k}$

$\qquad\qquad+\frac{\partial^{3}}{\partial x_{k}\partial x_{l}\partial x_{i}%
}(\nabla\log\Phi_{P})_{j}+\frac{\partial^{3}}{\partial x_{j}\partial
x_{k}\partial x_{l}}(\nabla\log\Phi_{P})_{i}](y)\qquad\ \qquad$

$\qquad=-$ $[\frac{\partial^{3}X_{l}}{\partial x_{i}\partial x_{j}\partial
x_{k}}+\frac{\partial^{3}X_{k}}{\partial x_{l}\partial x_{i}\partial x_{j}%
}+\frac{\partial^{3}X_{j}}{\partial x_{k}\partial x_{l}\partial x_{i}}%
+\frac{\partial^{3}X_{i}}{\partial x_{j}\partial x_{k}\partial x_{l}%
}](y)\qquad\qquad\qquad\qquad\qquad\qquad\left(  8^{\ast}\right)  \qquad$

Higher derivatives follow.

(ix) It is immediate from (iv) and (v) that:

\qquad\ $\frac{\partial}{\partial x_{\text{a}}}(\nabla$log$\Phi_{P}%
)_{j}(y)=-\frac{\partial X_{j}}{\partial x_{\text{a}}}(y)$ and $\frac
{\partial^{2}}{\partial x_{\text{a}}\partial x_{\text{b}}}(\nabla$log$\Phi
_{P})_{j}(y)=-\frac{\partial^{2}X_{j}}{\partial x_{\text{a}}\partial
x_{\text{b}}}(y)$

\qquad The property for higher derivatives follow.

Compare the last two formulae with the bizarre formula:

\qquad$\ \frac{\partial}{\partial x_{i}}(\nabla$log$\Phi_{P})_{\text{a}}(y)=$
J$_{1}$ + J$_{2}=$ $\underset{j=q+1}{\overset{n}{%
{\textstyle\sum}
}}X_{j}(y)\perp_{\text{a}ij}(y)-\frac{\partial X_{i}}{\partial x_{\text{a}}%
}(y)$

from (xv) below.

(x) It is immediate from (v) that:

$\qquad\frac{\partial^{2}}{\partial x_{\text{a}}\partial x_{\text{b}}}(\nabla
$log$\Phi_{P})_{j}(y)=-\frac{\partial^{2}X_{j}}{\partial x_{\text{a}}\partial
x_{\text{b}}}(y)$

Higher derivatives follow.

In particular when $y=y_{0}$ is the centre of Fermi coordinates, we have for
a,b = 1,...,q and $j=q+1,...,n$

$\qquad(\nabla$log$\Phi_{P})_{j}(y_{0})=-X_{j}(y_{0})\qquad\qquad\qquad
\qquad\qquad\qquad\qquad\qquad\qquad\qquad\left(  9\right)  \qquad\qquad
\qquad\qquad\qquad\qquad\qquad\ \ \ \ $

$\qquad\frac{\partial}{\partial x_{k}}(\nabla\log\Phi_{P})_{j}(y_{0}%
)+\frac{\partial}{\partial x_{j}}(\nabla\log\Phi_{P})_{k}(y_{0})=-\frac
{\partial X_{j}}{\partial x_{k}}(y_{0})-\frac{\partial X_{k}}{\partial x_{j}%
}(y_{0})\ \ \ \ \ \ \ \ $

\qquad$\frac{\partial^{2}}{\partial x_{k}\partial x_{l}}(\nabla\log\Phi
_{P})_{j}(y_{0})+\frac{\partial^{2}}{\partial x_{j}\partial x_{l}}(\nabla
\log\Phi_{P})_{k}(y_{0})+\frac{\partial^{2}}{\partial x_{k}\partial x_{j}%
}(\nabla\log\Phi_{P})_{l}(y)$

$\qquad=-\frac{\partial^{2}X_{j}}{\partial x_{k}\partial x_{l}}(y_{0}%
)-\frac{\partial^{2}X_{k}}{\partial x_{j}\partial x_{l}}(y_{0})-\frac
{\partial^{2}X_{l}}{\partial x_{k}\partial x_{j}}(y_{0})$

$\qquad\frac{\partial}{\partial x_{\text{a}}}(\nabla$log$\Phi_{P})_{j}%
(y_{0})=-\frac{\partial X_{j}}{\partial x_{\text{a}}}(y_{0})$

$\qquad\frac{\partial^{2}}{\partial x_{\text{a}}\partial x_{\text{b}}}(\nabla
$log$\Phi_{P})_{j}(y_{0})=-\frac{\partial^{2}X_{j}}{\partial x_{\text{a}%
}\partial x_{\text{b}}}(y_{0})$ $\qquad$

Higher derivatives follow.$\ \qquad$

(xi) By $\left(  B_{1}\right)  $ and $\left(  B_{2}\right)  $ above,

$\left(  B_{13}\right)  $\qquad log$\Phi_{P}(x)=$ $\int_{0}^{1}<$X($\gamma
$(s)) , $\dot{\gamma}$(s)%
$>$%
ds $=-\int_{0}^{1}[$ $\sigma(<X,\nabla\sigma>](\gamma(s))ds$ \qquad$\left(
10\right)  $

$\qquad\qquad\qquad\qquad\ \ \ =-\underset{i=q+1}{\overset{n}{\sum}}\int%
_{0}^{1}$ $(x_{i}X_{i})(\gamma(s))ds$\qquad

where $\gamma$ is the unique minimal geodesic from x to P in time 1 and
meeting P orthogonally at a point $y\in P:$

The geodesic $\gamma:[0,1]\longrightarrow$ M$_{0}$ is given in Fermi
coordinates as:

\qquad$\qquad\gamma(s)=\left(  \text{x}_{1},...,\text{x}_{q},(1-s)\text{x}%
_{q+1},...,(1-s)\text{x}_{n}\right)  $\qquad\qquad\qquad\qquad\qquad
\qquad\qquad$\qquad$\qquad\qquad\qquad

Now by the \textbf{definition} of the \textbf{gradient operator} we have for a
= 1,...,q and $j=1,...,q,q+1,...,n:$

\qquad$(\nabla$log$\Phi_{P})_{\text{a}}(x)=$ $\underset{j=1}{\overset{n}{\sum
}}$g$^{j\text{a}}(x)[\frac{\partial}{\partial x_{j}}$log$\Phi_{P}](x)$

$\qquad=$ $\underset{\text{b=1}}{\overset{q}{\sum}}$g$^{\text{ab}}%
(x)[\frac{\partial}{\partial x_{\text{b}}}$log$\Phi_{P}](x)+$
$\underset{j=q+1}{\overset{n}{\sum}}$g$^{\text{a}j}(x)[\frac{\partial
}{\partial x_{j}}$log$\Phi_{P}](x)$

Since g$^{\text{ab}}(y)=\delta^{\text{ab}}$ and g$^{\text{a}j}(y)=\delta
^{\text{a}j}=0$ for a,b = 1,...,q and $j=q+1,...,n,$ we have:

$\qquad(\nabla$log$\Phi_{P})_{\text{a}}(y)=[\frac{\partial}{\partial
x_{\text{a}}}$log$\Phi_{P}](y)$

From the last equation in $\left(  B_{13}\right)  $ above, we have for $y\in
U\subset P,$\qquad

\qquad$(\nabla\log\Phi_{P})_{\text{a}}(y)=-\underset{i=q+1}{\overset{n}{\sum}%
}\int_{\text{0}}^{1}$ $\frac{\partial}{\partial x_{\text{a}}}[(x_{i}%
X_{i})(\gamma(s))]ds$

$\qquad\qquad\qquad\qquad=-\underset{i=q+1}{\overset{n}{\sum}}\int_{\text{0}%
}^{1}$ $[x_{i}\frac{\partial X_{i}}{\partial x_{\text{a}}}](\gamma
(s))\frac{\partial}{\partial x_{\text{a}}}\gamma(s)]ds$

Since $x=y,$ the geodesic $\gamma:[0,1]\longrightarrow M_{0}$ is now the
constant geodesic defined by:

$\gamma(s)=y$ for all $s\in\lbrack0,1].$ From the definition of the \ geodesic
$\gamma$ in normal coordinates given above, we have:

\qquad\ $\frac{\partial}{\partial x_{\text{a}}}\gamma(s)=1$ for a = 1,...,q. Consequently,

\qquad$(\nabla\log\Phi_{P})_{\text{a}}(y)=-\underset{i=q+1}{\overset{n}{\sum}%
}\int_{\text{0}}^{1}$ $[x_{i}\frac{\partial X_{i}}{\partial x_{\text{a}}%
}](y)]ds=-\underset{i=q+1}{\overset{n}{\sum}}x_{i}(y)\frac{\partial X_{i}%
}{\partial x_{\text{a}}}(y)$

From the property of Fermi coordinates in $\left(  1.2\right)  ^{\ast},$ we
have $x_{i}(y)=0$ for $i=q+1,...,n.$ We conclude that:

$(\nabla\log\Phi_{P})_{\text{a}}(y)=0\qquad\qquad\qquad\qquad\left(
11\right)  $

(xii) We note that in \textbf{Chapter 6} we carried out expansions of the the
components of the metric tensor g$_{ij}$ and its inverse g$^{ij}$ in
\textbf{normal} Fermi coordinates. From those expansions we see that the
derivatives of g$_{ij}$ and g$^{ij}$ with respect to tangential coordinates
vanish: $\frac{\partial g^{i\text{a}}}{\partial x_{\text{b}}}(x)\ =0$ for b =
1,...,q$.$ All higher derivatives must vanish as well.

Secondly we note that by (iii) of $\left(  C_{7}\right)  $, we have:
$\frac{\partial^{2}}{\partial\text{x}_{i}\partial\text{x}_{j}}\gamma(s)=0$ for
$i,j=1,...,q,q+1,...,n=0.$

The above properties will be used below without explicit mention:

Therefore from $\left(  B_{14}\right)  ,$ we have:

$\ \ \left(  B_{15}\right)  \ \ \ \frac{\partial}{\partial x_{\text{c}}%
}(\nabla$log$\Phi_{P})_{\text{a}}(x)\ =-\underset{i=q+1}{\overset{n}{\sum}%
}\int_{\text{0}}^{1}$ g$^{i\text{a}}(x)[\frac{\partial X_{i}}{\partial
x_{\text{c}}}(\gamma(s))\frac{\partial}{\partial x_{\text{b}}}\gamma
(s)\frac{\partial}{\partial x_{i}}\gamma(s)]ds$

$\qquad\qquad\qquad\qquad-\underset{j=1}{\overset{n}{\sum}}%
\underset{i=q+1}{\overset{n}{\sum}}\int_{\text{0}}^{1}$ g$^{j\text{a}%
}(x)[\frac{\partial x_{i}}{\partial x_{\text{b}}}+x_{i}\frac{\partial
}{\partial x_{\text{b}}}\frac{\partial X_{i}}{\partial x_{j}})(\gamma
(s))\frac{\partial}{\partial x_{\text{b}}}\gamma(s)\frac{\partial}{\partial
x_{i}}\gamma(s)]ds$

Since $\frac{\partial x_{i}}{\partial x_{\text{b}}}=\delta_{i\text{b}}=0$ $,$
we have:

$\left(  B_{16}\right)  \qquad\ \frac{\partial}{\partial x_{\text{b}}}(\nabla
$log$\Phi_{P})_{\text{a}}(x)\ =-\underset{i=q+1}{\overset{n}{\sum}}%
\int_{\text{0}}^{1}$ g$^{i\text{a}}$(x)[$\frac{\partial X_{i}}{\partial
x_{\text{b}}}$($\gamma$(s))$\frac{\partial}{\partial x_{\text{b}}}\gamma
$(s)$\frac{\partial}{\partial x_{i}}\gamma(s)]ds$

$\qquad\qquad\qquad\qquad\qquad\ -\underset{j=1}{\overset{n}{\sum}%
}\underset{i=q+1}{\overset{n}{\sum}}\int_{\text{0}}^{1}$ g$^{j\text{a}}%
$(x)[$x_{i}\frac{\partial}{\partial x_{\text{b}}}\frac{\partial X_{i}%
}{\partial x_{j}})$($\gamma$(s))$\frac{\partial}{\partial x_{\text{b}}}\gamma
$(s)$\frac{\partial}{\partial x_{i}}\gamma(s)]ds$\qquad

As we saw earlier, for $x=y,$ the geodesic $\gamma:[0,1]\longrightarrow M_{0}$
is the constant geodesic:

$\gamma(s)=y$ for all $s\in\lbrack0,1]$ and\ $\frac{\partial}{\partial
x_{\text{a}}}\gamma(s)=1.$ Therefore,\qquad

$\ \frac{\partial}{\partial x_{\text{b}}}(\nabla$log$\Phi_{P})_{\text{a}%
}(y)\ =-\underset{i=q+1}{\overset{n}{\sum}}\int_{\text{0}}^{1}$ g$^{i\text{a}%
}(y)$[$\frac{\partial X_{i}}{\partial x_{\text{b}}}(y)\frac{\partial}{\partial
x_{\text{b}}}\gamma(s)\frac{\partial}{\partial x_{i}}\gamma(s)]ds$

$\qquad\qquad\qquad\qquad-\underset{j=1}{\overset{n}{\sum}}%
\underset{i=q+1}{\overset{n}{\sum}}\int_{\text{0}}^{1}$ g$^{j\text{a}%
}(y)[x_{i}\frac{\partial}{\partial x_{\text{b}}}\frac{\partial X_{i}}{\partial
x_{j}})(y)\frac{\partial}{\partial x_{\text{b}}}\gamma(s)\frac{\partial
}{\partial x_{i}}\gamma(s)]ds$\qquad\qquad\qquad\qquad\qquad

Now, g$^{i\text{a}}(y)=0$ for a = 1,...,q and $i=q+1,...,n$ and $x_{i}(y)=0$
for $i=q+1,...,n.$

We see that:

\qquad$\ \frac{\partial}{\partial x_{\text{b}}}(\nabla$log$\Phi_{P}%
)_{\text{a}}(y)\ =0$ \qquad\qquad\qquad\qquad\qquad\qquad\qquad\qquad
\qquad\qquad\qquad$\left(  12\right)  $

(xiii) Further differentiating both sides of $\left(  B_{16}\right)  $ gives
for b,c = 1,...,q:

$\left(  B_{17}\right)  $\qquad$\ \frac{\partial^{2}}{\partial x_{\text{c}%
}\partial x_{\text{b}}}(\nabla$log$\Phi_{P})_{\text{a}}(y)\ =0\qquad
\qquad\qquad\qquad\qquad\qquad\qquad\qquad\qquad\left(  13\right)
\qquad\qquad\qquad\qquad\qquad\qquad\qquad\qquad\qquad\qquad\qquad
\qquad\ \qquad\qquad\qquad\qquad\qquad\qquad\qquad\qquad\qquad\qquad
\qquad\qquad$

\subsubsection{M\textbf{ixed Derivatives:}}

(xiv) From the expression in $\left(  B_{7}\right)  ,$ we have:$\qquad$

$\qquad\frac{\partial^{2}}{\partial x_{\text{a}}\partial x_{k}}(\nabla
$log$\Phi_{P})_{j}(x_{0})+\frac{\partial^{2}}{\partial x_{\text{a}}\partial
x_{j}}(\nabla\log\Phi_{P})_{k}(x_{0})+\frac{\partial^{2}X_{j}}{\partial
x_{\text{a}}\partial x_{k}}(x_{0})+\frac{\partial^{2}X_{k}}{\partial
x_{\text{a}}\partial x_{j}}(x_{0})$

$\qquad=-\underset{i=q+1}{\overset{n}{\sum}}x_{i}(x_{0})\left(  \frac
{\partial^{3}}{\partial x_{\text{a}}\partial x_{j}\partial x_{k}}(\nabla
\log\Phi_{P})_{i}+\frac{\partial^{3}X_{i}}{\partial x_{\text{a}}\partial
x_{j}\partial x_{k}}\right)  (x_{0})$

$\qquad$Therefore at $x_{0}=y,$ we have:

$\frac{\partial^{2}}{\partial x_{\text{a}}\partial x_{k}}(\nabla$log$\Phi
_{P})_{j}(y)+\frac{\partial^{2}}{\partial x_{\text{a}}\partial x_{j}}%
(\nabla\log\Phi_{P})_{k}(y)=-\left(  \frac{\partial^{2}X_{j}}{\partial
x_{\text{a}}\partial x_{k}}-\frac{\partial^{2}X_{k}}{\partial x_{\text{a}%
}\partial x_{j}}\right)  (y)$

In particular for $k=j,$

$\frac{\partial^{2}}{\partial x_{\text{a}}\partial x_{j}}(\nabla$log$\Phi
_{P})_{j}(y)=-\frac{\partial^{2}X_{j}}{\partial x_{\text{a}}\partial x_{j}%
}(y)$

As an alternative proof of the above formula, we can also differentiate both
sides of the equation in $\left(  B_{5}\right)  $:$\qquad\qquad\qquad
\qquad\qquad\qquad\ \ $\qquad\qquad\qquad\qquad\qquad\qquad\qquad\qquad
\qquad\qquad\qquad\qquad\qquad

We differentiate both sides with respect to normal Fermi coordinates: For a
=1,...,q and $j=q+1,...,n:$

$\underset{i=q+1}{\overset{n}{\sum}}\frac{\partial x_{i}}{\partial x_{j}%
}\left(  \frac{\partial}{\partial x_{\text{a}}}\nabla\log\Phi_{P})_{i}%
+\frac{\partial X_{i}}{\partial x_{\text{a}}}\right)  +$
$\underset{i=q+1}{\overset{n}{\sum}}x_{i}\left(  \frac{\partial^{2}}{\partial
x_{\text{a}}\partial x_{j}}\nabla\log\Phi_{P})_{i}+\frac{\partial^{2}X_{i}%
}{\partial x_{\text{a}}\partial x_{j}}\right)  =0$

Since $\frac{\partial x_{i}}{\partial x_{j}}=\delta_{j}^{i}$ we have:

$\frac{\partial}{\partial x_{\text{a}}}\nabla\log\Phi_{P})_{j}+\frac{\partial
X_{j}}{\partial x_{\text{a}}}=-$ $\underset{i=q+1}{\overset{n}{\sum}}%
x_{i}\left(  \frac{\partial^{2}}{\partial x_{\text{a}}\partial x_{j}}%
\nabla\log\Phi_{P})_{i}+\frac{\partial^{2}X_{i}}{\partial x_{\text{a}}\partial
x_{j}}\right)  $

We differentiate both sides again with respect to the normal coordinate
$x_{k}$, and have for $k=q+1,...,n:$

$\frac{\partial^{2}}{\partial x_{\text{a}}\partial x_{k}}\nabla\log\Phi
_{P})_{j}+\frac{\partial^{2}X_{j}}{\partial x_{\text{a}}\partial x_{k}}+$
$\underset{i=q+1}{\overset{n}{\sum}}\frac{\partial x_{i}}{\partial x_{k}%
}\left(  \frac{\partial^{2}}{\partial x_{\text{a}}\partial x_{j}}\nabla
\log\Phi_{P})_{i}+\frac{\partial^{2}X_{i}}{\partial x_{\text{a}}\partial
x_{j}}\right)  $

$+\underset{i=q+1}{\overset{n}{\sum}}x_{i}\left(  \frac{\partial^{3}}{\partial
x_{\text{a}}\partial x_{j}\partial x_{k}}\nabla\log\Phi_{P})_{i}%
+\frac{\partial^{3}X_{i}}{\partial x_{\text{a}}\partial x_{j}\partial x_{k}%
}\right)  =0$

Since $\frac{\partial x_{i}}{\partial x_{k}}=\delta_{k}^{i}$ we have:

$\frac{\partial^{2}}{\partial x_{\text{a}}\partial x_{k}}\nabla\log\Phi
_{P})_{j}+\frac{\partial^{2}X_{j}}{\partial x_{\text{a}}\partial x_{k}}+$
$\frac{\partial^{2}}{\partial x_{\text{a}}\partial x_{j}}\nabla\log\Phi
_{P})_{k}+\frac{\partial^{2}X_{k}}{\partial x_{\text{a}}\partial x_{j}}$

$+\underset{i=q+1}{\overset{n}{\sum}}x_{i}\left(  \frac{\partial^{3}}{\partial
x_{\text{a}}\partial x_{j}\partial x_{k}}\nabla\log\Phi_{P})_{i}%
+\frac{\partial^{3}X_{i}}{\partial x_{\text{a}}\partial x_{j}\partial x_{k}%
}\right)  =0$

Since $x_{i}(y)=0$ for $i=q+1,...,n,$ we arrive at the same formula:

$\frac{\partial^{2}}{\partial x_{\text{a}}\partial x_{k}}(\nabla\log\Phi
_{P})_{j}(y)+$ $\frac{\partial^{2}}{\partial x_{\text{a}}\partial x_{j}}%
\nabla\log\Phi_{P})_{k}(y)=-\frac{\partial^{2}X_{j}}{\partial x_{\text{a}%
}\partial x_{k}}(y)-\frac{\partial^{2}X_{k}}{\partial x_{\text{a}}\partial
x_{j}}(y)$

In particular, for $k=j,$

\qquad$\frac{\partial^{2}}{\partial x_{\text{a}}\partial x_{j}}(\nabla\log
\Phi_{P})_{j}(y)=-\frac{\partial^{2}X_{j}}{\partial x_{\text{a}}\partial
x_{j}}(y)$

(xv) Recalling the Einstein convention of summation over repeated indices, we
have:\qquad\qquad\qquad\qquad

$\qquad(\nabla$log$\Phi_{P})_{\text{a}}(x)=g^{j\text{a}}(x)\frac{\partial
}{\partial x_{j}}\log\Phi_{P}(x)=g^{j\text{a}}(x)\frac{\partial}{\partial
x_{j}}\log\Phi_{P}(x)$

$\qquad=[g^{j\text{a}}\frac{1}{\Phi}\frac{\partial\Phi_{P}}{\partial x_{j}%
}](x)$

For $i=q+1,...,n$ and $j=1,....n,$ we have:

J$_{1}$ + J$_{2}=\frac{\partial}{\partial x_{i}}(\nabla$log$\Phi
_{P})_{\text{a}}(x)=\frac{\partial}{\partial x_{i}}[g^{j\text{a}}%
\frac{\partial}{\partial x_{j}}\log\Phi_{P}](x)=\frac{\partial}{\partial
x_{i}}[g^{j\text{a}}\frac{1}{\Phi}\frac{\partial\Phi_{P}}{\partial x_{j}}](x)$

$=\frac{\partial g^{j\text{a}}}{\partial x_{i}}(x)[\frac{1}{\Phi}%
\frac{\partial\Phi_{P}}{\partial x_{j}}](x)+g^{j\text{a}}(x)[\frac{1}{\Phi
^{2}}(\Phi_{P}\frac{\partial^{2}\Phi_{P}}{\partial x_{i}\partial x_{j}}%
-\frac{\partial\Phi_{P}}{\partial x_{i}}\frac{\partial\Phi_{P}}{\partial
x_{j}})](x)$

where,

J$_{1}$ $=\frac{\partial g^{j\text{a}}}{\partial x_{i}}(x)[\frac{1}{\Phi}%
\frac{\partial\Phi_{P}}{\partial x_{j}}](x);$J$_{2}=g^{j\text{a}}(x)[\frac
{1}{\Phi^{2}}(\Phi_{P}\frac{\partial^{2}\Phi_{P}}{\partial x_{i}\partial
x_{j}}-\frac{\partial\Phi_{P}}{\partial x_{i}}\frac{\partial\Phi_{P}}{\partial
x_{j}})](x)$

Therefore at $x=y\in U\subset P$ we have for a = 1,...,q,

J$_{1}=\underset{j=1}{\overset{n}{%
{\textstyle\sum}
}}\frac{\partial g^{j\text{a}}}{\partial x_{i}}(y)[\frac{1}{\Phi}%
\frac{\partial\Phi_{P}}{\partial x_{j}}](y)=$ $\underset{\text{b=1}%
}{\overset{\text{q}}{%
{\textstyle\sum}
}}\frac{\partial g^{\text{ab}}}{\partial x_{i}}(y)[\frac{1}{\Phi}%
\frac{\partial\Phi_{P}}{\partial x_{\text{b}}}%
](y)+\underset{j=q+1}{\overset{n}{%
{\textstyle\sum}
}}\frac{\partial g^{j\text{a}}}{\partial x_{i}}(y)[\frac{1}{\Phi}%
\frac{\partial\Phi_{P}}{\partial x_{j}}](y)$

Since $\Phi_{P}(y)=1;\frac{\partial g^{j\text{a}}}{\partial x_{i}}(y)=$
$\perp_{\text{a}ji}=-\perp_{\text{a}ij}$by (ii) of \textbf{Table A}$_{4}$ for
$i,j=q+1,...,n;$

Next we have: $\frac{\partial\Phi_{P}}{\partial x_{\text{b}}}(y)=0$ for b =
1,...,q by (vii) of \textbf{Table B}$_{4}$ below and $\frac{\partial\Phi_{P}%
}{\partial x_{j}}(y)=-X_{j}(y)$ by (i) of Table B$_{4}$ below. Consequently,
we have for a = 1,...,q and $i,j=q+1,...,n,$

J$_{1}=$ $\underset{j=q+1}{\overset{n}{%
{\textstyle\sum}
}}X_{j}(y)\perp_{\text{a}ij}(y)$

Since $g^{j\text{a}}(y)=\delta^{j\text{a}};\Phi_{P}(y)=1$ and $\frac
{\partial\Phi_{P}}{\partial x_{\text{a}}}(y)=0$ by (vii) of Table B$_{4}$, we
have at $x=y:$

By (xi) of Table B$_{4},$

J$_{2}=\underset{j=1}{\overset{n}{%
{\textstyle\sum}
}}g^{j\text{a}}(y)[\frac{1}{\Phi^{2}}(\Phi_{P}\frac{\partial^{2}\Phi_{P}%
}{\partial x_{i}\partial x_{j}}-\frac{\partial\Phi_{P}}{\partial x_{i}}%
\frac{\partial\Phi_{P}}{\partial x_{j}})](y)=[\frac{\partial^{2}\Phi_{P}%
}{\partial x_{i}\partial x_{\text{a}}}](y)=-\frac{\partial X_{i}}{\partial
x_{\text{a}}}(y)$

We conclude that at $x=y$ we have for a =1,...,q and $i,j=q+1,...,n,$

$\ \frac{\partial}{\partial x_{i}}(\nabla$log$\Phi_{P})_{\text{a}}(y)=$
J$_{1}$ + J$_{2}=$ $\underset{j=q+1}{\overset{n}{%
{\textstyle\sum}
}}\left(  X_{j}\perp_{\text{a}ij}\right)  (y)-\frac{\partial X_{i}}{\partial
x_{\text{a}}}(y)$

(xvi) We next compute for $i,j=q+1,...,n:$\qquad\qquad

$\frac{\partial^{2}}{\partial x_{i}\partial x_{j}}(\nabla$log$\Phi
_{P})_{\text{a}}=\frac{\partial}{\partial x_{i}}$J$_{1}$ $+\frac{\partial
}{\partial x_{i}}$J$_{2}$

where we recall,

J$_{1}$ $=\underset{j=1}{\overset{n}{%
{\textstyle\sum}
}}\frac{\partial g^{j\text{a}}}{\partial x_{i}}(x)[\frac{1}{\Phi}%
\frac{\partial\Phi_{P}}{\partial x_{j}}](x)$ and J$_{2}%
=\underset{j=1}{\overset{n}{%
{\textstyle\sum}
}}g^{j\text{a}}(x)[\frac{1}{\Phi^{2}}(\Phi_{P}\frac{\partial^{2}\Phi_{P}%
}{\partial x_{i}\partial x_{j}}-\frac{\partial\Phi_{P}}{\partial x_{i}}%
\frac{\partial\Phi_{P}}{\partial x_{j}})](x)$

For $k=1,...,q,q+1,...,n,$ we have at a general point $x\in M_{0}:$

J$_{1}$ $=\underset{k=1}{\overset{n}{%
{\textstyle\sum}
}}\frac{\partial g^{k\text{a}}}{\partial x_{j}}(x)[\frac{1}{\Phi}%
\frac{\partial\Phi_{P}}{\partial x_{k}}](x)=$ $\underset{\text{b=1}%
}{\overset{\text{q}}{%
{\textstyle\sum}
}}\frac{\partial g^{\text{ab}}}{\partial x_{j}}(x)[\frac{1}{\Phi}%
\frac{\partial\Phi_{P}}{\partial x_{\text{b}}}%
](x)+\underset{k=q+1}{\overset{n}{%
{\textstyle\sum}
}}\frac{\partial g^{k\text{a}}}{\partial x_{j}}(x)[\frac{1}{\Phi}%
\frac{\partial\Phi_{P}}{\partial x_{k}}](x)$

J$_{2}=\underset{k=1}{\overset{n}{%
{\textstyle\sum}
}}g^{k\text{a}}(x)[\frac{1}{\Phi^{2}}(\Phi_{P}\frac{\partial^{2}\Phi_{P}%
}{\partial x_{j}\partial x_{k}}-\frac{\partial\Phi_{P}}{\partial x_{j}}%
\frac{\partial\Phi_{P}}{\partial x_{k}})](x)=$ $\underset{\text{b=1}%
}{\overset{\text{q}}{%
{\textstyle\sum}
}}g^{\text{ab}}(x)[\frac{1}{\Phi^{2}}(\Phi_{P}\frac{\partial^{2}\Phi_{P}%
}{\partial x_{j}\partial x_{\text{b}}}-\frac{\partial\Phi_{P}}{\partial x_{j}%
}\frac{\partial\Phi_{P}}{\partial x_{\text{b}}})](x)$

$+\underset{k=q+1}{\overset{n}{%
{\textstyle\sum}
}}g^{k\text{a}}(x)[\frac{1}{\Phi^{2}}(\Phi_{P}\frac{\partial^{2}\Phi_{P}%
}{\partial x_{j}\partial x_{k}}-\frac{\partial\Phi_{P}}{\partial x_{j}}%
\frac{\partial\Phi_{P}}{\partial x_{k}})](x)$

We have for $i,j=q+1,...,n,$ (omitting the summation sign over b = 1,...,q and

$k=q+1,...,n):$

\qquad$\frac{\partial}{\partial x_{i}}$J$_{1}$ $=[\frac{\partial
^{2}g^{\text{ab}}}{\partial x_{i}\partial x_{j}}\frac{1}{\Phi}\frac
{\partial\Phi_{P}}{\partial x_{\text{b}}}](y)$ $+\frac{\partial g^{\text{ab}}%
}{\partial x_{j}}(y)[-\frac{1}{\Phi_{P}^{2}}\frac{\partial\Phi_{P}}{\partial
x_{i}}\frac{\partial\Phi_{P}}{\partial x_{\text{b}}}+\frac{1}{\Phi}%
\frac{\partial^{2}\Phi_{P}}{\partial x_{i}\partial x_{\text{b}}}](y)$

$\qquad\qquad$\ $\ +[\frac{\partial^{2}g^{k\text{a}}}{\partial x_{i}\partial
x_{j}}\frac{1}{\Phi}\frac{\partial\Phi_{P}}{\partial x_{k}}](y)$
$+\frac{\partial g^{k\text{a}}}{\partial x_{j}}(y)[-\frac{1}{\Phi_{P}^{2}%
}\frac{\partial\Phi_{P}}{\partial x_{i}}\frac{\partial\Phi_{P}}{\partial
x_{k}}+\frac{1}{\Phi}\frac{\partial^{2}\Phi_{P}}{\partial x_{i}\partial x_{k}%
}](y)$

We have:

$\frac{\partial\Phi_{P}}{\partial x_{\text{b}}}(y)=0$ for b = 1,...,q by (vii)
of \textbf{Table B}$_{4}$ below and $\frac{\partial g^{\text{ab}}}{\partial
x_{j}}(y)$

$=$ $2$T$_{\text{ab}j}(y_{0})$ by (ii) of \textbf{Table A}$_{6}.$

$\frac{\partial\Phi_{P}}{\partial x_{k}}(y)=-X_{k}(y)$ by (i) of \textbf{Table
B}$_{4}$ below, $\frac{\partial^{2}\Phi_{P}}{\partial x_{i}\partial
x_{\text{b}}}(y)=-\frac{\partial X_{i}}{\partial x_{\text{b}}}(y)$ by (xi) of
\textbf{Table B}$_{4}$

$\frac{\partial^{2}\Phi_{P}}{\partial x_{i}\partial x_{k}}(y)=X_{i}%
(y)X_{k}(y)-\frac{1}{2}\left(  \frac{\partial X_{k}}{\partial x_{i}}%
+\frac{\partial X_{i}}{\partial x_{k}}\right)  (y)$ by (ii) of \textbf{Table
B}$_{4}$

$\frac{\partial g^{k\text{a}}}{\partial x_{j}}(y)=-\perp_{\text{a}jk}(y)$ by
(ii) of \textbf{Table A}$_{4}$ and $\frac{\partial^{2}g^{k\text{a}}}{\partial
x_{i}\partial x_{j}}(y)=-\frac{4}{3}\left[  R_{i\text{a}jk}+R_{j\text{a}%
ik}\right]  (y)$

by (iii) of \textbf{Table A}$_{3}$. Therefore we have for a,b = 1,...,q and
$i,j,k=q+1,...,n,$

$\frac{\partial}{\partial x_{i}}$J$_{1}$ $=-$ $2\underset{\text{b=1}%
}{\overset{\text{q}}{%
{\textstyle\sum}
}}$T$_{\text{ab}j}(y_{0})\frac{\partial X_{i}}{\partial x_{\text{b}}}(y)$

$+\frac{4}{3}\underset{k=q+1}{\overset{n}{%
{\textstyle\sum}
}}\left[  R_{i\text{a}jk}+R_{j\text{a}ik}\right]  (y)X_{k}(y)-$
$\underset{k=q+1}{\overset{n}{%
{\textstyle\sum}
}}\perp_{\text{a}jk}(y)\left[  -X_{i}X_{k}+X_{i}X_{k}-\frac{1}{2}\left(
\frac{\partial X_{i}}{\partial x_{k}}+\frac{\partial X_{k}}{\partial x_{i}%
}\right)  \right]  (y)\qquad$

$=-$ $2\underset{\text{b=1}}{\overset{\text{q}}{%
{\textstyle\sum}
}}$T$_{\text{ab}j}(y_{0})\frac{\partial X_{i}}{\partial x_{\text{b}}}%
(y)+\frac{4}{3}\underset{k=q+1}{\overset{n}{%
{\textstyle\sum}
}}\left[  R_{i\text{a}jk}+R_{j\text{a}ik}\right]  (y)X_{k}(y)$

$+\frac{1}{2}\underset{k=q+1}{\overset{n}{%
{\textstyle\sum}
}}\perp_{\text{a}jk}(y)\left[  \left(  \frac{\partial X_{i}}{\partial x_{k}%
}+\frac{\partial X_{k}}{\partial x_{i}}\right)  \right]  (y)\qquad\qquad
\qquad\qquad$

Next we have at $x\in M_{0}:$

\qquad J$_{2}=g^{k\text{a}}(x)[\frac{1}{\Phi^{2}}(\Phi_{P}\frac{\partial
^{2}\Phi_{P}}{\partial x_{j}\partial x_{k}}-\frac{\partial\Phi_{P}}{\partial
x_{j}}\frac{\partial\Phi_{P}}{\partial x_{k}})](x)$

We differentiate at at $x=y\in U\subset P\subset M_{0}:$

$\frac{\partial}{\partial x_{i}}$J$_{2}=\frac{\partial g^{k\text{a}}}{\partial
x_{i}}(y)[\frac{1}{\Phi^{2}}(\Phi_{P}\frac{\partial^{2}\Phi_{P}}{\partial
x_{j}\partial x_{k}}-\frac{\partial\Phi_{P}}{\partial x_{j}}\frac{\partial
\Phi_{P}}{\partial x_{k}})](y)$

$+g^{k\text{a}}(y)\frac{\partial}{\partial x_{i}}[\Phi_{P}^{-1}\frac
{\partial^{2}\Phi_{P}}{\partial x_{j}\partial x_{k}}-\Phi_{P}^{-2}%
\frac{\partial\Phi_{P}}{\partial x_{j}}\frac{\partial\Phi_{P}}{\partial x_{k}%
}](y)$

$=\frac{\partial g^{k\text{a}}}{\partial x_{i}}(y)[\frac{1}{\Phi^{2}}(\Phi
_{P}\frac{\partial^{2}\Phi_{P}}{\partial x_{j}\partial x_{k}}-\frac
{\partial\Phi_{P}}{\partial x_{j}}\frac{\partial\Phi_{P}}{\partial x_{k}%
})](y)$

$+g^{k\text{a}}(y)[-\Phi_{P}^{-2}\frac{\partial\Phi_{P}}{\partial x_{i}}%
\frac{\partial^{2}\Phi_{P}}{\partial x_{j}\partial x_{k}}+\Phi_{P}^{-1}%
\frac{\partial^{3}\Phi_{P}}{\partial x_{i}\partial x_{j}\partial x_{k}}](y)$

$+g^{k\text{a}}(y)[2\Phi_{P}^{-3}\frac{\partial\Phi_{P}}{\partial x_{i}}%
\frac{\partial\Phi_{P}}{\partial x_{j}}\frac{\partial\Phi_{P}}{\partial x_{k}%
}-\Phi_{P}^{-2}(\frac{\partial^{2}\Phi_{P}}{\partial x_{i}\partial x_{j}}%
\frac{\partial\Phi_{P}}{\partial x_{k}}+\frac{\partial\Phi_{P}}{\partial
x_{j}}\frac{\partial^{2}\Phi_{P}}{\partial x_{i}\partial x_{k}})](y)$

$\qquad\qquad\ \ =J_{21}+J_{22}+J_{23}$

For a,b =1,...,q and $i,j,k=q+1,...,n,$ we have:

$J_{21}=\frac{\partial g^{k\text{a}}}{\partial x_{i}}(y)[\frac{1}{\Phi^{2}%
}(\Phi_{P}\frac{\partial^{2}\Phi_{P}}{\partial x_{j}\partial x_{k}}%
-\frac{\partial\Phi_{P}}{\partial x_{j}}\frac{\partial\Phi_{P}}{\partial
x_{k}})](y)$

$=\frac{\partial g^{\text{ab}}}{\partial x_{i}}(y)[\frac{1}{\Phi^{2}}(\Phi
_{P}\frac{\partial^{2}\Phi_{P}}{\partial x_{j}\partial x_{\text{b}}}%
-\frac{\partial\Phi_{P}}{\partial x_{j}}\frac{\partial\Phi_{P}}{\partial
x_{\text{b}}})](y)$

$+\frac{\partial g^{k\text{a}}}{\partial x_{i}}(y)[\frac{1}{\Phi^{2}}(\Phi
_{P}\frac{\partial^{2}\Phi_{P}}{\partial x_{j}\partial x_{k}}-\frac
{\partial\Phi_{P}}{\partial x_{j}}\frac{\partial\Phi_{P}}{\partial x_{k}%
})](y)$

Values of all terms of the above expression are in Table B$_{4}$ and have been cited

above and so we use them here without giving the references:

We have for a = 1,...,q and $i,j=q+1,...,n:$

$J_{21}=-$ $2\underset{\text{b=1}}{\overset{\text{q}}{%
{\textstyle\sum}
}}$T$_{\text{ab}i}(y)\frac{\partial X_{j}}{\partial x_{\text{b}}%
}(y)-\underset{k=q+1}{\overset{n}{%
{\textstyle\sum}
}}\perp_{\text{a}ik}(y)[X_{j}X_{k}-\frac{1}{2}\left(  \frac{\partial X_{k}%
}{\partial x_{j}}+\frac{\partial X_{j}}{\partial x_{k}}\right)  -X_{j}%
X_{k}](y)$

\qquad$=-$ $2\underset{\text{b=1}}{\overset{\text{q}}{%
{\textstyle\sum}
}}$T$_{\text{ab}i}(y)\frac{\partial X_{j}}{\partial x_{\text{b}}}(y)+\frac
{1}{2}\underset{k=q+1}{\overset{n}{%
{\textstyle\sum}
}}\perp_{\text{a}ik}(y)[\left(  \frac{\partial X_{k}}{\partial x_{j}}%
+\frac{\partial X_{j}}{\partial x_{k}}\right)  ](y)$

Next we have for a = 1,...,q, $i,j=q+1,...,n$ and $k=1,...,q,q+1,...,n$

$J_{22}=g^{k\text{a}}(y)[-\Phi_{P}^{-2}\frac{\partial\Phi_{P}}{\partial x_{i}%
}\frac{\partial^{2}\Phi_{P}}{\partial x_{j}\partial x_{k}}+\Phi_{P}^{-1}%
\frac{\partial^{3}\Phi_{P}}{\partial x_{i}\partial x_{j}\partial x_{k}%
}](y)=[-\frac{\partial\Phi_{P}}{\partial x_{i}}\frac{\partial^{2}\Phi_{P}%
}{\partial x_{j}\partial x_{\text{a}}}+\frac{\partial^{3}\Phi_{P}}{\partial
x_{i}\partial x_{j}\partial x_{\text{a}}}](y)$

All terms of the last expression above have already been given except the last

term which is given by (xiv) of

\textbf{Table B}$_{4}:$

$-\frac{\partial\Phi_{P}}{\partial x_{i}}(y)\frac{\partial^{2}\Phi_{P}%
}{\partial x_{\text{a}}\partial x_{j}}(y)=-X_{i}(y)\frac{\partial X_{j}%
}{\partial x_{\text{a}}}(y)$

$\frac{\partial^{3}\Phi_{P}}{\partial x_{\text{a}}\partial x_{i}\partial
x_{j}}(y)=2\left(  X_{i}\frac{\partial X_{j}}{\partial x_{\text{a}}}%
+X_{j}\frac{\partial X_{i}}{\partial x_{\text{a}}}\right)  (y)-\frac{1}%
{2}\left(  \frac{\partial^{2}X_{i}}{\partial x_{\text{a}}\partial x_{j}}%
+\frac{\partial^{2}X_{j}}{\partial x_{\text{a}}\partial x_{i}}\right)  (y)$

Therefore,

$J_{22}=[-X_{i}\frac{\partial X_{j}}{\partial x_{\text{a}}}+2\left(
X_{i}\frac{\partial X_{j}}{\partial x_{\text{a}}}+X_{j}\frac{\partial X_{i}%
}{\partial x_{\text{a}}}\right)  -\frac{1}{2}\left(  \frac{\partial^{2}X_{i}%
}{\partial x_{\text{a}}\partial x_{j}}+\frac{\partial^{2}X_{j}}{\partial
x_{\text{a}}\partial x_{i}}\right)  ](y)$

\qquad$=[X_{i}\frac{\partial X_{j}}{\partial x_{\text{a}}}+2X_{j}%
\frac{\partial X_{i}}{\partial x_{\text{a}}}-\frac{1}{2}\left(  \frac
{\partial^{2}X_{i}}{\partial x_{\text{a}}\partial x_{j}}+\frac{\partial
^{2}X_{j}}{\partial x_{\text{a}}\partial x_{i}}\right)  ](y)$

$J_{23}=g^{k\text{a}}(y)[2\Phi_{P}^{-3}\frac{\partial\Phi_{P}}{\partial x_{i}%
}\frac{\partial\Phi_{P}}{\partial x_{j}}\frac{\partial\Phi_{P}}{\partial
x_{k}}-\Phi_{P}^{-2}(\frac{\partial^{2}\Phi_{P}}{\partial x_{i}\partial x_{j}%
}\frac{\partial\Phi_{P}}{\partial x_{k}}+\frac{\partial\Phi_{P}}{\partial
x_{j}}\frac{\partial^{2}\Phi_{P}}{\partial x_{i}\partial x_{k}})](y)$

\qquad$=[2\frac{\partial\Phi_{P}}{\partial x_{i}}\frac{\partial\Phi_{P}%
}{\partial x_{j}}\frac{\partial\Phi_{P}}{\partial x_{\text{a}}}-(\frac
{\partial^{2}\Phi_{P}}{\partial x_{i}\partial x_{j}}\frac{\partial\Phi_{P}%
}{\partial x_{\text{a}}}+\frac{\partial\Phi_{P}}{\partial x_{j}}\frac
{\partial^{2}\Phi_{P}}{\partial x_{i}\partial x_{\text{a}}})](y)$

Since $\frac{\partial\Phi_{P}}{\partial x_{\text{a}}}(y)=0$ and $\frac
{\partial^{2}\Phi_{P}}{\partial x_{i}\partial x_{\text{a}}}(y)=-\frac{\partial
X_{i}}{\partial x_{\text{a}}}(y)$ we have:

$J_{23}=-[\frac{\partial\Phi_{P}}{\partial x_{j}}\frac{\partial^{2}\Phi_{P}%
}{\partial x_{i}\partial x_{\text{a}}})](y)=-[(-X_{j})(-\frac{\partial X_{i}%
}{\partial x_{\text{a}}})](y)=-[X_{j}\frac{\partial X_{i}}{\partial
x_{\text{a}}}](y)$

Therefore,

$\frac{\partial}{\partial x_{i}}$J$_{2}=J_{21}+J_{22}+J_{23}$

$=-$ $2\underset{\text{b=1}}{\overset{\text{q}}{%
{\textstyle\sum}
}}$T$_{\text{ab}i}(y)\frac{\partial X_{j}}{\partial x_{\text{b}}}(y)+\frac
{1}{2}\underset{k=q+1}{\overset{n}{%
{\textstyle\sum}
}}\perp_{\text{a}ik}(y)[\left(  \frac{\partial X_{k}}{\partial x_{j}}%
+\frac{\partial X_{j}}{\partial x_{k}}\right)  ](y)$

$+[X_{i}\frac{\partial X_{j}}{\partial x_{\text{a}}}+2X_{j}\frac{\partial
X_{i}}{\partial x_{\text{a}}}-\frac{1}{2}\left(  \frac{\partial^{2}X_{i}%
}{\partial x_{\text{a}}\partial x_{j}}+\frac{\partial^{2}X_{j}}{\partial
x_{\text{a}}\partial x_{i}}\right)  ](y)-[X_{j}\frac{\partial X_{i}}{\partial
x_{\text{a}}})](y)$

$=-$ $2\underset{\text{b=1}}{\overset{\text{q}}{%
{\textstyle\sum}
}}$T$_{\text{ab}i}(y)\frac{\partial X_{j}}{\partial x_{\text{b}}}(y)+\frac
{1}{2}\underset{k=q+1}{\overset{n}{%
{\textstyle\sum}
}}\perp_{\text{a}ik}(y)[\left(  \frac{\partial X_{k}}{\partial x_{j}}%
+\frac{\partial X_{j}}{\partial x_{k}}\right)  ](y)$

$+[X_{i}\frac{\partial X_{j}}{\partial x_{\text{a}}}+X_{j}\frac{\partial
X_{i}}{\partial x_{\text{a}}}-\frac{1}{2}\left(  \frac{\partial^{2}X_{i}%
}{\partial x_{\text{a}}\partial x_{j}}+\frac{\partial^{2}X_{j}}{\partial
x_{\text{a}}\partial x_{i}}\right)  ](y)$

We conclude that:

$\frac{\partial^{2}}{\partial x_{i}\partial x_{j}}(\nabla$log$\Phi
_{P})_{\text{a}}(y)=\frac{\partial}{\partial x_{i}}$J$_{1}$ $+\frac{\partial
}{\partial x_{i}}$J$_{2}$

$=-$ $2\underset{\text{b=1}}{\overset{\text{q}}{%
{\textstyle\sum}
}}$T$_{\text{ab}j}(y_{0})\frac{\partial X_{i}}{\partial x_{\text{b}}}%
(y)+\frac{4}{3}\underset{k=q+1}{\overset{n}{%
{\textstyle\sum}
}}\left[  R_{i\text{a}jk}+R_{j\text{a}ik}\right]  (y)X_{k}(y)$

$+\frac{1}{2}\underset{k=q+1}{\overset{n}{%
{\textstyle\sum}
}}\perp_{\text{a}jk}(y)\left[  \left(  \frac{\partial X_{i}}{\partial x_{k}%
}+\frac{\partial X_{k}}{\partial x_{i}}\right)  \right]  (y)$

$-$ $2\underset{\text{b=1}}{\overset{\text{q}}{%
{\textstyle\sum}
}}$T$_{\text{ab}i}(y)\frac{\partial X_{j}}{\partial x_{\text{b}}}(y)+\frac
{1}{2}\underset{k=q+1}{\overset{n}{%
{\textstyle\sum}
}}\perp_{\text{a}ik}(y)[\left(  \frac{\partial X_{k}}{\partial x_{j}}%
+\frac{\partial X_{j}}{\partial x_{k}}\right)  ](y)$

$+[X_{i}\frac{\partial X_{j}}{\partial x_{\text{a}}}+X_{j}\frac{\partial
X_{i}}{\partial x_{\text{a}}}-\frac{1}{2}\left(  \frac{\partial^{2}X_{i}%
}{\partial x_{\text{a}}\partial x_{j}}+\frac{\partial^{2}X_{j}}{\partial
x_{\text{a}}\partial x_{i}}\right)  ](y)$

We have the more elegant expression:

$\frac{\partial^{2}}{\partial x_{i}\partial x_{j}}(\nabla$log$\Phi
_{P})_{\text{a}}(y)=-$ $2\underset{\text{b=1}}{\overset{\text{q}}{%
{\textstyle\sum}
}}$T$_{\text{ab}j}(y_{0})\frac{\partial X_{i}}{\partial x_{\text{b}}%
}(y)-2\underset{\text{b=1}}{\overset{\text{q}}{%
{\textstyle\sum}
}}T_{\text{ab}i}(y)\frac{\partial X_{j}}{\partial x_{\text{b}}}(y)$

$+\frac{1}{2}\underset{k=q+1}{\overset{n}{%
{\textstyle\sum}
}}\perp_{\text{a}jk}(y)\left[  \left(  \frac{\partial X_{i}}{\partial x_{k}%
}+\frac{\partial X_{k}}{\partial x_{i}}\right)  \right]  (y)+\frac{1}%
{2}\underset{k=q+1}{\overset{n}{%
{\textstyle\sum}
}}\perp_{\text{a}ik}(y)[\left(  \frac{\partial X_{k}}{\partial x_{j}}%
+\frac{\partial X_{j}}{\partial x_{k}}\right)  ](y)$\qquad

$+\frac{4}{3}\underset{k=q+1}{\overset{n}{%
{\textstyle\sum}
}}\left[  R_{i\text{a}jk}+R_{j\text{a}ik}\right]  (y)X_{k}(y)+[X_{i}%
\frac{\partial X_{j}}{\partial x_{\text{a}}}+X_{j}\frac{\partial X_{i}%
}{\partial x_{\text{a}}}-\frac{1}{2}\left(  \frac{\partial^{2}X_{i}}{\partial
x_{\text{a}}\partial x_{j}}+\frac{\partial^{2}X_{j}}{\partial x_{\text{a}%
}\partial x_{i}}\right)  ](y)$

\qquad\qquad\qquad\qquad\qquad\qquad\qquad\qquad\qquad\qquad\qquad\qquad
\qquad\qquad\qquad\qquad\qquad$\blacksquare$

In particular,

$\frac{\partial^{2}}{\partial x_{i}^{2}}(\nabla$log$\Phi_{P})_{\text{a}%
}(y)=-4\underset{\text{b=1}}{\overset{\text{q}}{%
{\textstyle\sum}
}}T_{\text{ab}i}(y)\frac{\partial X_{i}}{\partial x_{\text{b}}}(y)$

$+\underset{k=q+1}{\overset{n}{%
{\textstyle\sum}
}}\perp_{\text{a}ik}(y)\left(  \frac{\partial X_{i}}{\partial x_{k}}%
+\frac{\partial X_{k}}{\partial x_{i}}\right)  (y)$

$\bigskip+\frac{8}{3}\underset{k=q+1}{\overset{n}{%
{\textstyle\sum}
}}R_{i\text{a}ik}(y)X_{k}(y)+\left(  2X_{i}\frac{\partial X_{i}}{\partial
x_{\text{a}}}-\frac{\partial^{2}X_{i}}{\partial x_{\text{a}}\partial x_{i}%
}\right)  (y$

$\qquad\qquad\qquad\qquad\qquad\qquad\qquad\qquad\qquad\qquad\qquad
\qquad\qquad\qquad\qquad\blacksquare$

\section{\textbf{Table B}$_{2}:$\textbf{ Gradient Vector Fields}}

We consider \textbf{Fermi coordinates} $x_{1},...,x_{q},x_{q+1},...,x_{n}$
based at based at $y_{0}\in M_{0}$ in the

star-shaped Fermi neighbourhood M$_{0}.$We have the following \textbf{more
general gradient formulae}

when X = gradf for a smooth function f:M$\rightarrow R$ and a general point
x$_{0}\in M_{0}:$

(i)$\qquad\nabla\log\Phi_{P}(x_{0})=\nabla f\circ\pi_{\text{P}}(x_{0})-\nabla
f(x_{0})=X\circ\pi_{\text{P}}(x_{0})-X$

In particular, in normal coordinates,

\qquad$\qquad\nabla\log\Phi_{P}(x_{0})=-X(x_{0})$

Equivalently, we have component-wise:

For a,b,c = 1,...,q and $i,j,k=1,...,q,q+1,...,n$ and for $y\in P$ in the neighbourhood

of $y_{0}\in P,$ and for $x_{0}\in M_{0},$ we have:

$\qquad(\nabla\log\Phi_{P})_{j}(x_{0})$ $=$ $\underset{\text{c=1}%
}{\overset{q}{\sum}}$g$^{\text{c}j}(x_{0})\frac{\partial\text{f}}{\partial
x_{\text{c}}}(y)-X_{j}(x_{0})\qquad$

(ii) $\ \frac{\partial}{\partial x_{\text{a}}}(\nabla\log\Phi_{P})_{j}(x_{0})$
$=$ $\underset{\text{c=1}}{\overset{q}{\sum}}g^{\text{c}j}(x_{0}%
)\frac{\partial^{2}\text{f}}{\partial x_{\text{a}}\partial x_{\text{c}}}(y)-$
$\frac{\partial X_{j}}{\partial x_{\text{a}}}(x_{0})$

\qquad Higher tangential derivatives follow. For normal derivatives, we hava:

(iii) $\frac{\partial}{\partial x_{i}}(\nabla\log\Phi_{P})_{j}(x_{0})$ $=$
$\underset{\text{c=1}}{\overset{q}{\sum}}\frac{\partial g^{\text{c}j}%
}{\partial x_{i}}(x_{0})\frac{\partial\text{f}}{\partial x_{\text{c}}}(y)-$
$\frac{\partial X_{j}}{\partial x_{i}}(x_{0})$

$\ \ \ \frac{\partial^{2}}{\partial x_{i}\partial x_{k}}(\nabla\log\Phi
_{P})_{j}(x_{0})=$ $\underset{\text{c=1}}{\overset{q}{\sum}}\frac{\partial
^{2}g^{\text{c}j}}{\partial x_{i}\partial x_{k}}(x_{0})\frac{\partial\text{f}%
}{\partial x_{\text{c}}}(y)+$ $\underset{\text{c,d=1}}{\overset{q}{\sum}}%
\frac{\partial g^{\text{c}j}}{\partial x_{i}}(x_{0})\frac{\partial^{2}%
\text{f}}{\partial x_{\text{c}}\partial x_{k}}(y)-$ $\frac{\partial^{2}X_{j}%
}{\partial x_{i}\partial x_{k}}(x_{0})$

\qquad Higher normal derivatives follow.

In particular, for $x_{0}=y\in P,$ we have:

(iv) $\frac{\partial}{\partial x_{\text{a}}}(\nabla\log\Phi_{P})_{\text{b}%
}(y)$ $=$ $\underset{\text{b=1}}{\overset{q}{\sum}}\frac{\partial^{2}\text{f}%
}{\partial x_{\text{a}}\partial x_{\text{b}}}(y)-$ $\frac{\partial
X_{\text{b}}}{\partial x_{\text{a}}}(y)$

(v) $\frac{\partial}{\partial x_{\text{a}}}(\nabla\log\Phi_{P})_{j}(y)$ $=$
$-$ $\frac{\partial X_{j}}{\partial x_{\text{a}}}(y)$

(vi) $\frac{\partial}{\partial x_{i}}(\nabla\log\Phi_{P})_{\text{a}}(y)$
$=2\underset{\text{b=1}}{\overset{q}{\sum}}T_{\text{ab}i}(y)\frac
{\partial\text{f}}{\partial x_{\text{b}}}(y)-$ $\frac{\partial X_{\text{a}}%
}{\partial x_{i}}(y)$

(vii) $\frac{\partial}{\partial x_{i}}(\nabla\log\Phi_{P})_{j}(y)=-$
$\underset{\text{a=1}}{\overset{q}{\sum}}\perp_{\text{a}ij}(y)\frac
{\partial\text{f}}{\partial x_{\text{a}}}(y)-\frac{\partial X_{j}}{\partial
x_{i}}(y)$

We specialize to the simpler case of a \textbf{gradient vector field} X and a
\textbf{singleton}:

P = $\left\{  y_{0}\right\}  .$ In this case, the Fermi coordinates
$x_{1},...,x_{q},x_{q+1},...,x_{n}$ become

\textbf{normal coordinates} based at $y_{0}\in M_{0}.$ We have for all
$x_{0}\in M_{0}:$

$i,j,k=1,...,q,q+1,...,n,$

(viii) $(\nabla$log$\Phi_{P})_{j}(x_{0})=-X_{j}(x_{0})\qquad\qquad\qquad
$\qquad\qquad\qquad\qquad\qquad\qquad$\qquad\qquad\qquad\qquad\qquad
\qquad\qquad\qquad\qquad\qquad\ \ \ \qquad\qquad\qquad\qquad\qquad\qquad$

(ix) $\frac{\partial}{\partial x_{i}}(\nabla$log$\Phi_{P})_{j}(x_{0}%
)=-\frac{\partial X_{j}}{\partial x_{i}}(x_{0})$

(x) $\frac{\partial^{2}}{\partial x_{i}\partial x_{j}}(\nabla\log\Phi_{P}%
)_{k}(x_{0})=-\frac{\partial^{2}X_{k}}{\partial x_{i}\partial x_{j}}(x_{0})$
$\qquad$

Formulae for higher derivatives follow.

\begin{center}
\qquad\qquad\qquad\qquad\qquad\qquad\qquad\qquad$\blacksquare$
\end{center}

\subsection{\textbf{Computations of B}$_{2}$}

\qquad(i) Let X be a \textbf{smooth gradient vector field} on M i.e. X =
$\nabla$f for

some smooth function f: M$\longrightarrow$R. Then,

\qquad$\Phi_{P}($x) = exp$\left\{  \int_{\text{0}}^{1}<\text{X(}%
\gamma\text{(s)) , }\dot{\gamma}\text{(s)%
$>$%
ds}\right\}  =$ exp$\left\{  \int_{0}^{1}<\nabla\text{f(}\gamma\text{(s)) ,
}\dot{\gamma}\text{(s)%
$>$%
ds}\right\}  $

$\qquad\qquad\ \ \ \ =$ exp$\left\{  \int_{\text{0}}^{1}\frac{\text{d}%
}{\text{ds}}\text{f(}\gamma\text{(s))ds}\right\}  =$ exp$\left\{
[\text{f(}\gamma\text{(s))}]_{0}^{1}\right\}  $

$\qquad\qquad\ \ \ \ =$ exp$\left\{  [\text{f(}\gamma\text{(1)) }-\text{
f(}\gamma\text{(0))}]\right\}  =$ exp$\left\{  [\text{f(y) }-\text{ f(x}%
_{0}\text{)}]\right\}  .$

Therefore,

\qquad log$\Phi_{P}($x$_{0}$) = f(y) $-$ f(x$_{0}$) = f$\circ\pi_{\text{P}}%
$(x$_{0}$) $-$ f(x$_{0}$)

Since X(x$_{0}$) = $\nabla$f(x$_{0}$), we have,

\qquad$\nabla\log\Phi_{P}(x_{0})=\nabla f\circ\pi_{\text{P}}(x_{0})-\nabla
f(x_{0})=X\circ\pi_{\text{P}}(x_{0})-X(x_{0})$

We wish to re-write the above equation \textbf{component-wise:}

By definition, the gradient operator $\nabla$ has the coordinate expression:

\qquad$\nabla$f $=$ g$^{ij}\partial_{i}$f$\partial_{j}$ where $\partial
_{i}=\frac{\partial}{\partial x_{i}}$Therefore,

\qquad$\nabla($f$\circ\pi_{\text{P}})$(x$_{0}$) $=$ g$^{ij}(x_{0})\partial
_{i}($f$\circ\pi_{P})(x_{0})\partial_{j}=$ g$^{ij}(x_{0})\partial_{i}$%
f$(\pi_{P}(x_{0}))\partial_{i}\pi_{P}(x_{0}))\partial_{j}$

\qquad Recall that for b = 1,...,q and $i=1,...,q,q+1,...,n,$

\qquad$\partial_{i}\pi_{P}(x_{0})=\frac{\partial}{\partial x_{i}}\pi_{P}%
(x_{0})=\left\{
\begin{array}
[c]{c}%
1\text{ for i=1,...,q}\\
0\text{ for }i=q+1,...,n
\end{array}
\right.  $

Therefore for b = 1,...,q and $j=1,...,q,q+1,...,n,$

$\qquad\nabla$f$\circ\pi_{\text{P}}$(x$_{0}$) $=$ g$^{\text{b}j}%
(x_{0})\partial_{\text{b}}$f$(y)\partial_{j}=$ g$^{\text{b}j}(x)\frac
{\partial\text{f}}{\partial x_{\text{b}}}(y)\frac{\partial}{\partial x_{j}}$

Consequently for b = 1,...,q and $j=1,...,q,q+1,...,n,$

$\left(  B_{18}\right)  $\qquad$(\nabla\log\Phi_{P})_{j}(x_{0})$ $=$
$\underset{\text{b=1}}{\overset{q}{\sum}}$g$^{\text{b}j}(x_{0})\frac
{\partial\text{f}}{\partial x_{\text{b}}}(y)-X_{j}(x_{0})\qquad\qquad\left(
14\right)  $

The above is the component-wise version of the above equation.

We note that in normal coordinates, q =0 and so f(y$_{0}$) is a constant.
Consequently,$\qquad\qquad$

\qquad\qquad\ $(\nabla\log\Phi_{P})(x_{0})$ $=$ $\bigtriangledown
f(x_{0})=-X(x_{0})$

(ii) Taking derivatives on both sides of $\left(  B_{18}\right)  $ above, it
is immediate that:

Here, unlike in $\left(  B_{6}\right)  ,$ we are able to define $(\nabla
\log\Phi_{P})_{j}(x_{0})$ for all coordinates

or a, b = 1,...,q ; $j=1,...,q,q+1,...,n.$

$\qquad\frac{\partial}{\partial x_{\text{a}}}(\nabla\log\Phi_{P})_{j}(x_{0})$
$=$ $\underset{\text{b=1}}{\overset{q}{\sum}}g^{\text{b}j}(x_{0}%
)\frac{\partial^{2}\text{f}}{\partial x_{\text{a}}\partial x_{\text{b}}}(y)-$
$\frac{\partial X_{j}}{\partial x_{\text{a}}}(x_{0})$

(iii) For $i,j=1,...,q,q+1,...,n,$ we differentiate both sides of $\left(
B_{18}\right)  $ given by:

$\qquad(\nabla\log\Phi_{P})_{j}(x_{0})$ $=$ $\underset{\text{b}%
=1}{\overset{q}{\sum}}$g$^{\text{b}j}(x_{0})\frac{\partial\text{f}}{\partial
x_{\text{b}}}(y)-X_{j}(x_{0}):$

$\frac{\partial}{\partial x_{i}}(\nabla\log\Phi_{P})_{j}(x_{0})$ $=$
$\underset{\text{b=1}}{\overset{q}{\sum}}\frac{\partial g^{\text{b}j}%
}{\partial x_{i}}(x_{0})\frac{\partial\text{f}}{\partial x_{\text{b}}}(y)+$
$\underset{\text{b=1}}{\overset{q}{\sum}}g^{\text{b}j}(x_{0})\frac
{\partial^{2}\text{f}}{\partial x_{i}\partial x_{\text{b}}}(\pi_{P}%
(x_{0}))\partial_{i}\pi_{P}(x_{0})-$ $\frac{\partial X_{j}}{\partial x_{i}%
}(x_{0})$

Since,

$\pi_{P}(x_{0}))=y$ ; $\partial_{i}\pi_{P}(x_{0})=$ $\left\{
\begin{array}
[c]{c}%
1\text{ for }i=1,...,q\\
0\text{ for }i=q+1,...,n
\end{array}
\right.  ,$

we have the equation for a, b = 1,...,q:

$\frac{\partial}{\partial x_{\text{a}}}(\nabla\log\Phi_{P})_{j}(x_{0})$ $=$
$\underset{\text{b=1}}{\overset{q}{\sum}}\frac{\partial g^{\text{b}j}%
}{\partial x_{\text{a}}}(x_{0})\frac{\partial\text{f}}{\partial x_{\text{b}}%
}(y)+$ $\underset{\text{b=1}}{\overset{q}{\sum}}g^{\text{b}j}(x_{0}%
)\frac{\partial^{2}\text{f}}{\partial x_{\text{a}}\partial x_{\text{b}}}(y)-$
$\frac{\partial X_{j}}{\partial x_{\text{a}}}(x_{0})$

Since expansions are in normal Fermi coordinates, $\frac{\partial
g^{\text{b}j}}{\partial x_{\text{a}}}(x_{0})=0$ and so we have

two \textbf{general equations} which give (iv) and (v) respectively:

$\frac{\partial}{\partial x_{i}}(\nabla\log\Phi_{P})_{j}(x_{0})$ $=$
$\underset{\text{b=1}}{\overset{q}{\sum}}\frac{\partial g^{\text{b}j}%
}{\partial x_{i}}(x_{0})\frac{\partial\text{f}}{\partial x_{\text{b}}}(y)-$
$\frac{\partial X_{j}}{\partial x_{i}}(x_{0})$

(vi) We take $x_{0}=y\in P$ and consider four equations taken from the

general equations above. From the first one, we have for a, b = 1,...,q :

$\frac{\partial}{\partial x_{\text{a}}}(\nabla\log\Phi_{P})_{\text{c}}(y)$ $=$
$\underset{\text{b=1}}{\overset{q}{\sum}}g^{\text{bc}}(y)\frac{\partial
^{2}\text{f}}{\partial x_{\text{a}}\partial x_{\text{b}}}(y)-$ $\frac{\partial
X_{\text{c}}}{\partial x_{\text{a}}}(y)$

Since $g^{\text{bc}}(y)=\delta^{\text{bc}},$ we have (switching the roles of b
and c) for a, b = 1,...,q:

$\qquad\frac{\partial}{\partial x_{\text{a}}}(\nabla\log\Phi_{P})_{\text{b}%
}(y)$ $=$ $\underset{\text{b=1}}{\overset{q}{\sum}}\frac{\partial^{2}\text{f}%
}{\partial x_{\text{a}}\partial x_{\text{b}}}(y)-$ $\frac{\partial
X_{\text{b}}}{\partial x_{\text{a}}}(y)$ \ \ \ \ \ \ \ \ \ \ \ \ \ \ \ \ \ \ 

(vii) Next, still from the first equation, we have for a, b = 1,...,q and
$j=q+1,...,n$:

$\qquad\frac{\partial}{\partial x_{\text{a}}}(\nabla\log\Phi_{P})_{j}(y)$ $=$
$\underset{\text{b=1}}{\overset{q}{\sum}}g^{\text{b}j}(y)\frac{\partial
^{2}\text{f}}{\partial x_{\text{a}}\partial x_{\text{b}}}(y)-$ $\frac{\partial
X_{j}}{\partial x_{\text{a}}}(y)$

It is clear that: $g^{\text{b}j}(y)=\delta^{\text{b}j}=0$ and so,

\qquad$\frac{\partial}{\partial x_{\text{a}}}(\nabla\log\Phi_{P})_{j}(y)$ $=$
$-$ $\frac{\partial X_{j}}{\partial x_{\text{a}}}(y)$

(viii) From the second equation, we have for a,b =1,...,q and $i=q+1,...,n,$

$\qquad\frac{\partial}{\partial x_{i}}(\nabla\log\Phi_{P})_{\text{a}}(y)$ $=$
$\underset{\text{b=1}}{\overset{q}{\sum}}\frac{\partial g^{\text{ba}}%
}{\partial x_{i}}(y)\frac{\partial\text{f}}{\partial x_{\text{b}}}(y)-$
$\frac{\partial X_{\text{a}}}{\partial x_{i}}(y)$

Since $\frac{\partial g^{\text{ba}}}{\partial x_{i}}(y)=2T_{\text{ab}i}(y),$
we have for a,b =1,...,n and $i=q+1,...,n,$

$\qquad\frac{\partial}{\partial x_{i}}(\nabla\log\Phi_{P})_{\text{a}}(y)$
$=2\underset{\text{b=1}}{\overset{q}{\sum}}T_{\text{ab}i}(y)\frac
{\partial\text{f}}{\partial x_{\text{b}}}(y)-$ $\frac{\partial X_{\text{a}}%
}{\partial x_{i}}(y)$

(ix) We have the last equation from the second general equation: for b = 1,...,q

and $i,j=q+1,...,n,$

$\qquad\frac{\partial}{\partial x_{i}}(\nabla\log\Phi_{P})_{j}(x_{0})$ $=$
$\underset{\text{b=1}}{\overset{q}{\sum}}\frac{\partial g^{\text{b}j}%
}{\partial x_{i}}(x_{0})\frac{\partial\text{f}}{\partial x_{\text{b}}}(y)-$
$\frac{\partial X_{j}}{\partial x_{i}}(x_{0})$

Since $\frac{\partial g^{\text{b}j}}{\partial x_{i}}(y)=\perp_{\text{b}%
ji}(y)=-\perp_{\text{b}ij}(y)$ we have (replacing b by a) for a = 1,...,q

and $i,j=q+1,...,n,$

$\qquad\frac{\partial}{\partial x_{i}}(\nabla\log\Phi_{P})_{j}(y)=-$
$\underset{\text{a=1}}{\overset{q}{\sum}}\perp_{\text{a}ij}(y)\frac
{\partial\text{f}}{\partial x_{\text{a}}}(y)-\frac{\partial X_{j}}{\partial
x_{i}}(y)\qquad$

We see from the formula in $\left(  B_{18}\right)  $ that when the submanifold
reduces to the point y$_{0}$

(the centre of Fermi coordinates and so all $y=y_{0}$ and so we have normal coordinates),

equivalently, $q=0$ and so the first expression on the RHS of $\left(
B_{18}\right)  $ vanishes and we have

for all x$_{0}\in$M$_{0}$ and for $j=1,...q,q+1,...,n$ $:$

$\left(  B_{19}\right)  \qquad(\nabla\log\Phi_{P})_{j}(x_{0})=-X_{j}(x_{0})$
\qquad$\left(  15\right)  $

We conclude that if X is a gradient vector field and the submanifold

\textbf{reduces to the singleton} $\left\{  y_{0}\right\}  ,$

\qquad$\qquad\qquad\qquad\qquad\qquad\qquad\nabla$log$\Phi_{P}(x_{0})=-\nabla
f(x_{0})=-X(x_{0})$

\begin{center}
$\ \qquad\qquad\qquad\ \ \ $
\end{center}

(ii) then we have for all x$_{0}\in$M$_{0}$:

$\ \left(  B_{20}\right)  \qquad\frac{\partial}{\partial x_{i}}(\nabla\log
\Phi_{P})_{k}(x_{0})=-\frac{\partial X_{k}}{\partial x_{i}}(x_{0})$ for
$i,k=1,...,n\qquad\qquad\qquad\left(  16\right)  \qquad\qquad\qquad
\qquad\qquad\qquad\qquad\qquad\qquad\qquad\qquad\qquad\qquad$\qquad
\qquad\qquad\qquad\qquad\qquad

(iii)\qquad Further differentiation in this case gives for
$i,j,k=1,...,q,q+1,...,n,$

$\left(  B_{21}\right)  \qquad\frac{\partial^{2}}{\partial x_{i}\partial
x_{j}}(\nabla\log\Phi_{P})_{k}(x_{0})=-\frac{\partial^{2}X_{k}}{\partial
x_{i}\partial x_{j}}(x_{0})$ $\qquad\qquad\qquad\qquad\qquad\ \ \ \left(
17\right)  \qquad\qquad\qquad\qquad\qquad$

Higher derivatives follow.$\qquad$

\begin{center}
\qquad\qquad\qquad\qquad\qquad\qquad\qquad\qquad\qquad\qquad\qquad\qquad
\qquad\qquad\qquad\qquad\qquad\qquad$\blacksquare$ $\qquad$\qquad\qquad
\qquad\qquad\qquad\qquad$\qquad\qquad\qquad\qquad\qquad\qquad\qquad
\qquad\ \qquad\qquad\qquad\ \ \qquad\qquad\qquad\qquad\qquad\qquad\qquad
\qquad$\qquad
\end{center}

\section{Table B$_{3}:$ The Laplacian of $\Phi$\qquad\qquad\qquad\qquad
\qquad$\qquad$}

\qquad\ For a general vector field X on M we have at the centre of Fermi
coordinates y$_{0}\in$P$:$

$\qquad$(i) \qquad div($\nabla$log$\Phi_{P}$)$(y_{0})=-$ divX$(y_{0})+$
$\underset{\text{a}=1}{\overset{q}{\sum}}\frac{\partial X_{\text{a}}}{\partial
x_{\text{a}}}(y_{0})$

When Fermi coordinates reduce to normal coordinates.

\qquad(ii)\qquad div($\nabla$log$\Phi_{P}$)$(y_{0})=-$ divX$(y_{0})$

$\qquad$(iii)\qquad\ $\Delta\Phi_{P}(y_{0})=$ $\left\Vert \text{X}\right\Vert
^{2}(y_{0})-$ divX$(y_{0})-\underset{\text{a}=1}{\overset{q}{\sum}}%
$X$_{\text{a}}^{2}(y_{0})+$ $\underset{\text{a}=1}{\overset{q}{\sum}}%
\frac{\partial X_{\text{a}}}{\partial x_{\text{a}}}(y_{0})$

\qquad\qquad\qquad$\qquad\qquad=$ $\left\Vert \text{X}\right\Vert _{M}%
^{2}(y_{0})-$ divX$_{M}(y_{0})-$ $\left\Vert \text{X}\right\Vert _{P}%
^{2}(y_{0})$ $+$ divX$_{P}(y_{0})$

\qquad(iv) When Fermi coordinates reduce to normal coordinates.

\qquad$\qquad\ \Delta\Phi_{P}(y_{0})=$ $\left\Vert \text{X}\right\Vert
^{2}(y_{0})-$ divX$(y_{0})$

\qquad(v) The \textbf{Laplacian }in the case of \textbf{gradient vector field} X:

In this case we have a more \textbf{general formula:}

The formula in (iv) can be obtained at a more general point. We have:

\qquad$\left(  B_{25}\right)  $\qquad\qquad\qquad$\Delta\Phi_{P}(x_{0}%
)=\Phi_{P}(x_{0})\left(  \left\Vert \text{X}\right\Vert _{M}^{2}%
-\operatorname{div}X\right)  (x_{0})$

\qquad(vi) Repeating the last case above, we have:

\qquad$\Delta^{2}\Phi(x_{0})=$ $\Phi_{P}(x_{0})\left(  \left\Vert
\text{X}\right\Vert _{M}^{2}-\operatorname{div}X\right)  ^{2}(x_{0})+$
$\Phi_{P}(x_{0})[\left(  \Delta\left\Vert \text{X}\right\Vert _{M}^{2}%
-\Delta\operatorname{div}X\right)  ](x_{0})$

$\qquad\qquad\qquad\qquad\qquad\qquad-\Phi(x_{0})\left\langle
X,\bigtriangledown\left(  \left\Vert \text{X}\right\Vert _{M}^{2}%
-\operatorname{div}X\right)  \right\rangle (x_{0})$

\qquad In particular, $\Phi_{P}(y_{0}),$ we have:

\qquad\qquad$\Delta^{2}\Phi(y_{0})=$ $\left(  \left\Vert \text{X}\right\Vert
_{M}^{2}-\operatorname{div}X\right)  ^{2}(y_{0})+$ $[\left(  \Delta\left\Vert
\text{X}\right\Vert _{M}^{2}-\Delta\operatorname{div}X\right)  ](y_{0})$

$\qquad\qquad\qquad\qquad-2\left\langle X,\bigtriangledown\left(  \left\Vert
\text{X}\right\Vert _{M}^{2}-\operatorname{div}X\right)  \right\rangle
(y_{0})$

\qquad\qquad\qquad\qquad\qquad\qquad\qquad\qquad\qquad\qquad\qquad\qquad
\qquad\qquad\qquad\qquad\qquad$\blacksquare$

\begin{remark}
The factor $\Phi_{P}(x_{0})$ on the RHS of \ $\left(  B_{25}\right)  $ above
was absent in my Thesis
\end{remark}

\qquad\qquad defended in Warwick University in 1989. This stands corrected here.

\qquad\qquad\qquad\qquad\qquad\qquad\qquad\qquad\qquad\qquad\qquad\qquad
\qquad\qquad\qquad\qquad\qquad$\blacksquare$

\subsection{\textbf{Computations of B}$_{3}$}

(i) Let X be a vector field on M and let $\left(  x_{1},...,x_{n}\right)  $ be
local coordinates on M.

Then by definition, the divergence operator is defined by (see for example,
\textbf{Hsu} $\left[  1\right]  ,$ p.74):

\qquad$\qquad\qquad\qquad\qquad\qquad\operatorname{div}X=$ $\frac{1}{\theta
}\underset{j=1}{\overset{n}{\sum}}\frac{\partial}{\partial x_{j}}(\theta
X_{j})$

where $\theta$ is defined in 1.6 of Chapter 1 here.

We choose Fermi coordinates and compute divX at the centre of Fermi
coordinates $y_{0}\in P\subset M_{0}:$

$\qquad\frac{1}{\theta(y_{0})}\underset{j=1}{\overset{n}{\sum}}\frac{\partial
}{\partial x_{j}}(\theta X_{j})(y_{0})=\frac{1}{\theta(y_{0})}%
\underset{j=1}{\overset{n}{\sum}}\frac{\partial\theta}{\partial x_{j}}%
(y_{0})X_{j}(y_{0})+$ $\underset{j=1}{\overset{n}{\sum}}\frac{\partial X_{j}%
}{\partial x_{j}}(y_{0})$

Now by (i) and (ii) \textbf{Table A}$_{9},$

$\qquad\theta(y_{0})=1;\frac{\partial\theta}{\partial\text{x}_{j}}%
(y_{0})=\left\{
\begin{array}
[c]{c}%
0\text{ for }j=1,...,q\\
-<H,j>(y_{0})\text{ for }j=q+1,...,n
\end{array}
\right.  ,$

and so we have:

$\left(  B_{22}\right)  $\qquad$\qquad\qquad\operatorname{div}_{M}X(y_{0})$
$=-\underset{j=q+1}{\overset{n}{\sum}}<H,j>(y_{0})X_{j}(y_{0})+$
$\underset{j=1}{\overset{n}{\sum}}\frac{\partial X_{j}}{\partial x_{j}}%
(y_{0})\qquad\qquad\left(  18\right)  $

Since $\frac{\partial\theta}{\partial\text{x}_{\text{a}}}(y_{0})=0$ for
tangntial coordinates a =1,...,q we have:

\qquad\qquad\qquad$\operatorname{div}_{P}X=$ $\underset{\text{a=1}%
}{\overset{\text{q}}{\sum}}\frac{\partial X_{\text{a}}}{\partial x_{\text{a}}%
}$

$\qquad$\ $\underset{j=q+1}{\overset{n}{\sum}}\frac{\partial X_{j}}{\partial
x_{j}}=\operatorname{div}_{M}X-\operatorname{div}_{P}%
X+\underset{j=q+1}{\overset{n}{\sum}}<H,j>(y_{0})X_{j}(y_{0})\qquad\left(
18\right)  ^{\ast}$

and,

\qquad\ $\underset{j=1}{\overset{n}{\sum}}\frac{\partial X_{j}}{\partial
x_{j}}(y_{0})=$ $\operatorname{div}_{M}X+\underset{j=q+1}{\overset{n}{\sum}%
}<H,j>(y_{0})X_{j}(y_{0})\qquad\qquad\ \ \left(  18\right)  ^{\ast\ast}%
\qquad\ \ \qquad\qquad$

We have by $\left(  B_{22}\right)  :$

\qquad div$(\nabla$log$\Phi_{P})(y_{0})=$
$\ -\underset{j=q+1}{\overset{n}{\sum}}<H,j>(y_{0})(\nabla$log$\Phi_{P}%
)_{j}(y_{0})+$ $\underset{j=1}{\overset{n}{\sum}}\frac{\partial}{\partial
x_{j}}(\nabla$log$\Phi_{P})_{j}(y_{0})\qquad$

\qquad$=$ $\ -\underset{j=q+1}{\overset{n}{\sum}}<H,j>(y_{0})(\nabla$%
log$\Phi_{P})_{j}(y_{0})+$ $\underset{j=q+1}{\overset{n}{\sum}}\frac{\partial
}{\partial x_{j}}(\nabla$log$\Phi_{P})_{j}(y_{0})+$ $\underset{\text{a}%
}{\overset{q}{\sum}}\frac{\partial}{\partial x_{\text{a}}}(\nabla$log$\Phi
_{P})_{\text{a}}(y_{0})$

$(\nabla$log$\Phi_{P})_{j}(y_{0})=-X_{j}(y_{0})$ by (vi) \textbf{Table B}%
$_{1}$ and $\frac{\partial}{\partial x_{j}}(\nabla$log$\Phi_{P})_{j}%
(y_{0})=-\frac{\partial X_{j}}{\partial x_{j}}(y_{0})$ by (vii) of
\textbf{Table B}$_{1},$

$\qquad\frac{\partial}{\partial x_{\text{b}}}(\nabla$log$\Phi_{P})_{\text{a}%
}(y)\ =0$ for a = 1,...,q by (xii) of \textbf{Table B}$_{1}$

Consequently we have:

$\left(  B_{23}\right)  $\qquad$\operatorname{div}(\nabla\log\Phi_{P}%
)(y_{0})=$ $\underset{j=q+1}{\overset{n}{\sum}}<H,j>(y_{0})X_{j}(y_{0})-$
$\underset{j=q+1}{\overset{n}{\sum}}\frac{\partial X_{j}}{\partial x_{j}%
}(y_{0})\qquad\qquad\qquad$

\qquad$\ \qquad\qquad=\underset{j=q+1}{\overset{n}{\sum}}<H,j>(y_{0}%
)X_{j}(y_{0})-$ $\underset{j=1}{\overset{n}{\sum}}\frac{\partial X_{j}%
}{\partial x_{j}}(y_{0})+$ $\underset{\text{a}=1}{\overset{q}{\sum}}%
\frac{\partial X_{\text{a}}}{\partial x_{\text{a}}}(y_{0})$

$\qquad\qquad\qquad=-$ $\operatorname{div}X(y_{0})+\operatorname{div}%
_{P}X(y_{0})$ $\qquad\qquad\qquad\qquad\qquad\qquad\qquad\left(  19\right)  $

We have thus shown that:

$\left(  B_{24}\right)  \qquad\qquad\operatorname{div}(\nabla\log\Phi_{P}%
$)$(y_{0})=-$ $\operatorname{div}_{M}X(y_{0})+$ $\operatorname{div}_{P}%
X(y_{0})\qquad\qquad\ \ \ \left(  20\right)  \qquad\qquad\qquad\qquad
\qquad\qquad\qquad\qquad\qquad$

(ii) In particular when Fermi coordinates reduce to \textbf{normal
coordinates}, we have the nice formula:

\qquad$\qquad\qquad\operatorname{div}(\nabla\log\Phi_{P}$)$(y_{0})=-$
$\operatorname{div}X(y_{0})$

$\qquad\qquad\qquad\qquad\qquad\qquad\qquad\qquad\qquad\qquad\qquad
\qquad\qquad\qquad\qquad\qquad\ \ \qquad$

(iii) By definition,$\ \qquad$

$\qquad\qquad\Delta\Phi_{P}$(y$_{0}$) $=$ $\operatorname{div}(\nabla\Phi_{P}%
$)(y$_{0}$) $=$ $\operatorname{div}$($\Phi_{P}\nabla\log\Phi_{P})(y_{0})$

\qquad It is well known that for a smooth vector field X and a smooth function
f:M$\longrightarrow R,$we have:

\qquad\qquad\qquad$\operatorname{div}(fX)=<\nabla f,X>+f\operatorname{div}X$.

Recalling that $\Phi_{P}(y_{0})=1,$ and $(\nabla$log$\Phi_{P})_{j}%
(y_{0})=-X_{j}(y_{0})$ by (vi) \textbf{Table B}$_{1},$ we have:

$\qquad\qquad\Delta\Phi_{P}(y_{0})=$ $<\nabla\log\Phi_{P},\nabla\log\Phi
_{P})>(y_{0})+$ $\operatorname{div}(\nabla\log\Phi_{P})(y_{0})$

We have by $B_{24}$ and $\left(  B_{25}\right)  :$

\qquad$\qquad\Delta\Phi_{P}(y_{0})=$ $\underset{j=q+1}{\overset{n}{\sum}}%
$X$_{j}^{2}(y_{0})-\operatorname{div}_{M}X(y_{0})+$ $\operatorname{div}%
X_{P}(y_{0})$

\qquad\qquad$=$ $\left\Vert \text{X}\right\Vert ^{2}(y_{0})-$
$\operatorname{div}X(y_{0})-$ $\underset{\text{a}=1}{\overset{q}{\sum}}%
$X$_{\text{a}}^{2}(y_{0})+$ $\underset{\text{a}=1}{\overset{q}{\sum}}%
\frac{\partial X_{\text{a}}}{\partial x_{\text{a}}}(y_{0})\qquad\qquad
\qquad\qquad\qquad\qquad\qquad\qquad$

$\left(  B_{26}\right)  $\qquad$\qquad\Delta\Phi_{P}(y_{0})=$ $\left\Vert
\text{X}\right\Vert _{M}^{2}(y_{0})-$ $\operatorname{div}X_{M}(y_{0})-$
$\left\Vert \text{X}\right\Vert _{P}^{2}(y_{0})$ $+$ $\operatorname{div}%
X_{P}(y_{0})\qquad\left(  21\right)  $

(iv) In particular when Fermi coordinates reduce to normal coordinates, we have:

$\qquad\qquad$\qquad$\Delta\Phi_{P}(y_{0})=\left\Vert \text{X}\right\Vert
_{M}^{2}(y_{0})-\operatorname{div}$X$_{M}(y_{0})\qquad\qquad\qquad\qquad
\qquad\qquad\qquad\left(  22\right)  $

(v) By definining equalities in (iii) above, we have:

$\qquad\qquad\qquad\Delta\Phi_{P}(x_{0})=\operatorname{div}(\Phi
\log\bigtriangledown\Phi)(x_{0})$

$\qquad\Delta\Phi_{P}(x_{0})=$
$<$%
$\nabla\Phi_{P},\nabla\log\Phi_{P}>(x_{0})$ $+$ $\Phi_{P}(x_{0}%
)\operatorname{div}(\nabla\log\Phi_{P}$)(x$_{0}$)

\qquad\qquad$\ \ \ \ \ \ \ =$ $\Phi_{P}(x_{0})[$%
$<$%
$\nabla\log\Phi_{P},\nabla\log\Phi_{P}>+$ div($\nabla\log\Phi_{P}$)$]$(x$_{0}$)

Since we have: $(\nabla$log$\Phi_{P})_{j}(x_{0})=-X_{j}(x_{0})$ by (vi)
\textbf{Table B}$_{1},$

$\Delta\Phi_{P}(x_{0})=\Phi_{P}(x_{0})[<-X,-X>+$ $\operatorname{div}%
(-X)](x_{0})=$ $\Phi_{P}(x_{0})\left(  \left\Vert \text{X}\right\Vert _{M}%
^{2}-\operatorname{div}X\right)  (x_{0})$

We have:

$\left(  B_{25}\right)  $\qquad\qquad$\Delta\Phi_{P}(x_{0})=\Phi_{P}%
(x_{0})\left(  \left\Vert \text{X}\right\Vert _{M}^{2}-\operatorname{div}%
X\right)  (x_{0})$

(vi) The last formula above gives:

\qquad\qquad$\Delta^{2}\Phi(x_{0})=\Delta\lbrack\Delta\Phi_{P}](x_{0})=$
$\Delta\lbrack\Phi_{P}\left(  \left\Vert \text{X}\right\Vert _{M}%
^{2}-\operatorname{div}X\right)  ](x_{0})$

$\qquad=$ $\Delta\Phi_{P}(x_{0})\left(  \left\Vert \text{X}\right\Vert
_{M}^{2}-\operatorname{div}X\right)  (x_{0})+$ $\Phi_{P}(x_{0})[\Delta\left(
\left\Vert \text{X}\right\Vert _{M}^{2}-\operatorname{div}X\right)  ](x_{0})$

$\qquad+2\left\langle \bigtriangledown\Phi_{P},\bigtriangledown\left(
\left\Vert \text{X}\right\Vert _{M}^{2}-\operatorname{div}X\right)
\right\rangle (x_{0})$

By $\left(  B_{25}\right)  $ and the fact that $\bigtriangledown\Phi_{P}%
=\Phi\bigtriangledown\log\Phi_{P},$we have:

$\Delta^{2}\Phi(x_{0})=$ $\Phi_{P}(x_{0})\left(  \left\Vert \text{X}%
\right\Vert _{M}^{2}-\operatorname{div}X\right)  ^{2}(x_{0})+$ $\Phi_{P}%
(x_{0})[\left(  \Delta\left\Vert \text{X}\right\Vert _{M}^{2}-\Delta
\operatorname{div}X\right)  ](x_{0})$

$\qquad\qquad+2\Phi\left\langle \bigtriangledown\log\Phi_{P},\bigtriangledown
\left(  \left\Vert \text{X}\right\Vert _{M}^{2}-\operatorname{div}X\right)
\right\rangle (x_{0})$

Since $\bigtriangledown\log\Phi_{P}=-X,$ we have the final expression:

\qquad$\left(  B_{26}\right)  \qquad\Delta^{2}\Phi(x_{0})=$ $\Phi_{P}%
(x_{0})\left(  \left\Vert \text{X}\right\Vert _{M}^{2}-\operatorname{div}%
X\right)  ^{2}(x_{0})$

$\qquad+$ $\Phi_{P}(x_{0})[\left(  \Delta\left\Vert \text{X}\right\Vert
_{M}^{2}-\Delta\operatorname{div}X\right)  ](x_{0})-2\Phi(x_{0})\left\langle
X,\bigtriangledown\left(  \left\Vert \text{X}\right\Vert _{M}^{2}%
-\operatorname{div}X\right)  \right\rangle (x_{0})$

In particular, since $\Phi_{P}(y_{0})=1,$ we have:

\qquad\qquad\qquad$\Delta^{2}\Phi(y_{0})=$ $\left(  \left\Vert \text{X}%
\right\Vert _{M}^{2}-\operatorname{div}X\right)  ^{2}(y_{0})+$ $[\left(
\Delta\left\Vert \text{X}\right\Vert _{M}^{2}-\Delta\operatorname{div}%
X\right)  ](y_{0})$

$\qquad\qquad\qquad\qquad\qquad\qquad-2\left\langle X,\bigtriangledown\left(
\left\Vert \text{X}\right\Vert _{M}^{2}-\operatorname{div}X\right)
\right\rangle (y_{0})$

We can compute \textbf{higher order Laplacians }from the formula in $\left(
B_{26}\right)  .$

\qquad\qquad\qquad\qquad\qquad\qquad\qquad\qquad\qquad\qquad\qquad\qquad
\qquad\qquad\qquad\qquad\qquad$\qquad\blacksquare$

\section{Table B$_{4}:$ \textbf{ Derivatives of }$\Phi$}

\subsection{\textbf{Normal Derivatives}}

We recall here again, that the Einstein convention of summation over repeated
indices is undertood for all that is here and beyond.

For $i,j,k,l=q+1,...,n,$ we have:

\qquad(i) $\frac{\partial\Phi_{P}}{\partial x_{i}}(y_{0})=-X_{i}(y_{0})$

\qquad(ii) $\frac{\partial^{2}\Phi_{P}}{\partial x_{i}\partial x_{j}}%
(y_{0})=X_{i}(y_{0})X_{j}(y_{0})-\frac{1}{2}\left(  \frac{\partial X_{i}%
}{\partial x_{j}}+\frac{\partial X_{j}}{\partial x_{i}}\right)  (y_{0})$

In particular, for a \textbf{totally geodesic} submanifold P (the second
fundamental form vanishes)$,$ we have:

\qquad(ii)$^{\ast}$ $\frac{\partial^{2}\Phi_{P}}{\partial x_{i}^{2}}%
(y_{0})=\overset{n}{\underset{i=1}{\sum}}[X_{i}^{2}-\frac{\partial X_{i}%
}{\partial x_{i}}](y_{0})-\underset{\text{a}=1}{\overset{q}{\sum}}X_{\text{a}%
}^{2}(y_{0})+\underset{\text{a}=1}{\overset{q}{\sum}}\frac{\partial
X_{\text{a}}}{\partial x_{\text{a}}}(y_{0})$

\qquad\qquad$=\left\Vert \text{X}\right\Vert _{M}^{2}(y_{0})-\left\Vert
\text{X}\right\Vert _{P}^{2}(y_{0})-$ $\operatorname{div}$X$_{M}(y_{0})$ $+$
$\operatorname{div}$X$_{P}(y_{0})\qquad\ \ \qquad\ \ \ \ \left(  23\right)  $

\qquad(ii)$^{\ast\ast}\frac{\partial}{\partial x_{j}}(\nabla$log$\Phi_{P}%
)_{i}(y_{0})=-\frac{1}{2}\left(  \frac{\partial X_{j}}{\partial x_{i}}%
+\frac{\partial X_{i}}{\partial x_{j}}\right)  (y_{0})=$ $\frac{\partial
}{\partial x_{i}}(\nabla\log\Phi_{P})_{j}(y_{0})$ $\left(  24\right)  $

\qquad

\qquad(iii) $\frac{\partial\Phi_{P}^{-1}}{\partial x_{i}}(y_{0})=-\frac
{\partial\Phi_{P}}{\partial x_{i}}(y_{0})=X_{i}(y_{0})$

\qquad(iv) $\frac{\partial^{2}\Phi_{P}^{-1}}{\partial x_{i}\partial x_{j}%
}(y_{0})=X_{i}(y_{0})X_{j}(y_{0})+\frac{1}{2}\left(  \frac{\partial X_{i}%
}{\partial x_{j}}+\frac{\partial X_{j}}{\partial x_{i}}\right)  (y_{0})$

\qquad\qquad In particular,

\qquad$\qquad\frac{\partial^{2}\Phi_{P}^{-1}}{\partial x_{i}^{2}}%
(y_{0})=[X_{i}^{2}+\frac{\partial X_{i}}{\partial x_{i}}](y_{0})=\left\Vert
\text{X}\right\Vert _{M}^{2}(y_{0})+\operatorname{div}X_{M}(y_{0})\qquad
\qquad\left(  25\right)  $

$\qquad\qquad\qquad\qquad\qquad-\left\Vert \text{X}\right\Vert _{P}^{2}%
(y_{0})$ $-$ $\operatorname{div}$X$_{P}(y_{0})$ $\qquad$

\qquad\qquad From $\left(  B_{54}\right)  :$

$\qquad$(v) $\qquad\frac{\partial^{3}\Phi_{P}}{\partial x_{i}\partial
x_{j}\partial x_{k}}(y)=-X_{i}(y)X_{j}(y)X_{k}(y)+\frac{1}{2}X_{i}(y)\left(
\frac{\partial X_{k}}{\partial x_{j}}+\frac{\partial X_{j}}{\partial x_{k}%
}\right)  (y)$

$\qquad\qquad\qquad\qquad+\frac{1}{2}X_{j}(y)\left(  \frac{\partial X_{k}%
}{\partial x_{i}}+\frac{\partial X_{i}}{\partial x_{k}}\right)  (y)+\frac
{1}{2}X_{k}(y)\left(  \frac{\partial X_{j}}{\partial x_{i}}+\frac{\partial
X_{i}}{\partial x_{j}}\right)  (y)$

$\qquad\qquad\qquad\qquad-\frac{1}{3}\left(  \frac{\partial^{2}X_{i}}{\partial
x_{j}\partial x_{k}}+\frac{\partial^{2}X_{j}}{\partial x_{i}\partial x_{k}%
}+\frac{\partial^{2}X_{k}}{\partial x_{i}\partial x_{j}}\right)  (y)$

\qquad\qquad In particular,

$\qquad\frac{\partial^{3}\Phi_{P}}{\partial x_{i}^{2}\partial x_{j}%
}(y)=[-X_{i}^{2}X_{j}+X_{j}\frac{\partial X_{i}}{\partial x_{i}}+X_{i}\left(
\frac{\partial X_{i}}{\partial x_{j}}+\frac{\partial X_{j}}{\partial x_{i}%
}\right)  -\frac{1}{3}\left(  \frac{\partial^{2}X_{j}}{\partial x_{i}^{2}%
}+2\frac{\partial^{2}X_{i}}{\partial x_{i}\partial x_{j}}\right)  ](y)$

\qquad(v)$^{\ast}$ $\frac{\partial^{3}\Phi_{P}}{\partial x_{i}^{2}\partial
x_{j}}(y_{0})=X_{j}[\operatorname{div}X-\left\Vert \text{X}\right\Vert
_{M}^{2}-\operatorname{div}X_{P}+\left\Vert \text{X}\right\Vert _{P}^{2}+$
$<H,i>X_{i}](y_{0})$ $\qquad\qquad\qquad\qquad\qquad\qquad\qquad\qquad
\qquad\qquad$

$\qquad\qquad+X_{i}(y_{0})\left(  \frac{\partial X_{j}}{\partial x_{i}}%
+\frac{\partial X_{i}}{\partial x_{j}}\right)  (y_{0})-\frac{1}{3}%
[\frac{\partial^{2}X_{j}}{\partial x_{i}^{2}}+2\frac{\partial^{2}X_{i}%
}{\partial x_{i}\partial x_{j}}](y_{0})$

\qquad(v)$^{\ast\ast}$ $\frac{\partial^{3}\Phi_{P}}{\partial x_{i}\partial
x_{j}^{2}}(y_{0})=\underset{i=q+1}{\overset{n}{\sum}}X_{i}(y_{0}%
)\operatorname{div}X_{M}-\operatorname{div}X_{P}-\left\Vert \text{X}%
\right\Vert _{M}^{2}+\left\Vert \text{X}\right\Vert _{P}^{2}$

$\qquad\qquad+<H,j>X_{j}](y_{0})$ $+X_{j}(y_{0})\left(  \frac{\partial X_{i}%
}{\partial x_{j}}+\frac{\partial X_{j}}{\partial x_{i}}\right)  (y_{0}%
)-\frac{1}{3}\left(  \frac{\partial^{2}X_{i}}{\partial x_{j}^{2}}%
+2\frac{\partial^{2}X_{j}}{\partial x_{i}\partial x_{j}}\right)  (y_{0})$

\qquad(v)$^{\ast\ast\ast}$ $[\frac{\partial^{2}}{\partial x_{i}\partial x_{j}%
}(\nabla\log\Phi_{P})_{k}](y_{0})=-\frac{1}{3}\left(  \frac{\partial^{2}X_{k}%
}{\partial x_{i}\partial x_{j}}+\frac{\partial^{2}X_{j}}{\partial
x_{i}\partial x_{k}}+\frac{\partial^{2}X_{i}}{\partial x_{j}\partial x_{k}%
}\right)  (y_{0})$

$\qquad\qquad\qquad\qquad\qquad-\frac{1}{3}(R_{ikjl}+R_{jkil})(y_{0}%
)X_{l}(y_{0})\qquad$

$\qquad\qquad+[\perp_{\text{a}ik}\frac{\partial X_{j}}{\partial x_{\text{a}}%
}+\perp_{\text{a}jk}\frac{\partial X_{i}}{\partial x_{\text{a}}}%
](y_{0})+[\perp_{\text{a}ik}\perp_{\text{a}jl}X_{l}+\perp_{\text{a}jk}%
\perp_{\text{a}il}X_{l}](y_{0})$

\qquad(vi) $\frac{\partial^{4}\Phi_{P}}{\partial x_{i}^{2}\partial x_{j}^{2}%
}(y_{0})=[X_{i}^{2}X_{j}^{2}-2X_{i}X_{j}\left(  \frac{\partial X_{j}}{\partial
x_{i}}+\frac{\partial X_{i}}{\partial x_{j}}\right)  -X_{i}^{2}\frac{\partial
X_{j}}{\partial x_{j}}-X_{j}^{2}\frac{\partial X_{i}}{\partial x_{i}}](y_{0})$

$\qquad\qquad+\frac{1}{2}\left(  \frac{\partial X_{j}}{\partial x_{i}}%
+\frac{\partial X_{i}}{\partial x_{j}}\right)  ^{2}(y_{0})+\left(
\frac{\partial X_{i}}{\partial x_{i}}\frac{\partial X_{j}}{\partial x_{j}%
}\right)  (y_{0})+\frac{2}{3}[X_{i}\left(  2\frac{\partial^{2}X_{j}}{\partial
x_{i}\partial x_{j}}+\frac{\partial^{2}X_{i}}{\partial x_{j}^{2}}\right)  $

$\qquad\qquad+X_{j}\left(  \frac{\partial^{2}X_{j}}{\partial x_{i}^{2}}%
+2\frac{\partial^{2}X_{i}}{\partial x_{i}\partial x_{j}}\right)
](y_{0})-\frac{1}{2}\left(  \frac{\partial^{3}X_{i}}{\partial x_{i}\partial
x_{j}^{2}}+\frac{\partial^{3}X_{j}}{\partial x_{i}^{2}\partial x_{j}}\right)
(y_{0})$

\qquad\qquad\qquad\qquad\qquad\qquad\qquad\qquad\qquad\qquad\qquad\qquad
\qquad\qquad

\qquad(vi)$^{\ast}$ Role the Riemannian manifold M and the submanifold P in
the above formula:

$\qquad\frac{\partial^{4}\Phi_{P}}{\partial x_{i}^{2}\partial x_{j}^{2}}%
(y_{0})=[\left\Vert X\right\Vert _{M}^{2}-\left\Vert X\right\Vert _{P}%
^{2}]^{2}(y_{0})$

$\qquad-2[\left\Vert X\right\Vert _{M}^{2}-\left\Vert X\right\Vert _{P}%
^{2}](y_{0})[\operatorname{div}_{M}X-\operatorname{div}_{P}%
X+\underset{j=q+1}{\overset{n}{\sum}}<H,j>X_{j}](y_{0})$

$\qquad+$ $[\operatorname{div}_{M}X-\operatorname{div}_{P}X]^{2}%
(y_{0})+[\operatorname{div}_{M}X-\operatorname{div}_{P}X](y_{0}%
)[\underset{i=q+1}{\overset{n}{\sum}}<H,i>X_{i}](y_{0})$

$\qquad+$ $[\operatorname{div}_{M}X-\operatorname{div}_{P}X](y_{0}%
)[\underset{j=q+1}{\overset{n}{\sum}}<H,j>X_{j}](y_{0}%
)+\underset{i,j=q+1}{\overset{n}{\sum}}<H,i><H,j>X_{i}X_{j}](y_{0})$

$\qquad-2\underset{i,j=q+1}{\overset{n}{%
{\textstyle\sum}
}}[X_{i}X_{j}\left(  \frac{\partial X_{j}}{\partial x_{i}}+\frac{\partial
X_{i}}{\partial x_{j}}\right)  ](y_{0})+\frac{1}{2}%
\underset{i,j=q+1}{\overset{n}{%
{\textstyle\sum}
}}\left(  \frac{\partial X_{j}}{\partial x_{i}}+\frac{\partial X_{i}}{\partial
x_{j}}\right)  ^{2}(y_{0})\qquad$

$\qquad+\frac{2}{3}\underset{i,j=q+1}{\overset{n}{%
{\textstyle\sum}
}}[X_{i}\left(  2\frac{\partial^{2}X_{j}}{\partial x_{i}\partial x_{j}}%
+\frac{\partial^{2}X_{i}}{\partial x_{j}^{2}}\right)  +X_{j}\left(
\frac{\partial^{2}X_{j}}{\partial x_{i}^{2}}+2\frac{\partial^{2}X_{i}%
}{\partial x_{i}\partial x_{j}}\right)  ](y_{0})$

$\qquad-\frac{1}{2}\underset{i,j=q+1}{\overset{n}{%
{\textstyle\sum}
}}\left(  \frac{\partial^{3}X_{i}}{\partial x_{i}\partial x_{j}^{2}}%
+\frac{\partial^{3}X_{j}}{\partial x_{i}^{2}\partial x_{j}}\right)  (y_{0})$

\qquad\qquad\qquad\qquad\qquad\qquad\qquad\qquad\qquad\qquad\qquad\qquad
\qquad\qquad\qquad\qquad\qquad\qquad\qquad\qquad$\blacksquare$

The above is a fairly more geometric presentation of the formula in which we
see the roles played by the divergence of the vector field X on the Riemannian
manifold M and the submanifold P as well as the norms on the tangent bundles
of the Riemannian manifold and the submanifold. We also see the role played by
the \textbf{mean curvatur}e of the submanifold P. The mean curvature will
disappear if we assume that the submanifold is \textbf{totally geodesic}.

(vi)$^{\ast\ast}$ We see that if the Fermi coordinates reduce to normal
coordinates, which is equivalent to the submanifold reducing to the centre of
Fermi coordinates $\left\{  y_{0}\right\}  ,$ then we have a simpler formula
in which all the submanifold terms disappear: For $i,j=1,...,q,q+1,...,n:$

$\qquad\frac{\partial^{4}\Phi_{P}}{\partial x_{i}^{2}\partial x_{j}^{2}}%
(y_{0})=[\left\Vert X\right\Vert _{M}^{2}](y_{0})-[\operatorname{div}%
_{M}X]]^{2}(y_{0})$

$\qquad-2[X_{i}X_{j}\left(  \frac{\partial X_{j}}{\partial x_{i}}%
+\frac{\partial X_{i}}{\partial x_{j}}\right)  ](y_{0})+\frac{1}{2}\left(
\frac{\partial X_{j}}{\partial x_{i}}+\frac{\partial X_{i}}{\partial x_{j}%
}\right)  ^{2}(y_{0})\qquad$

$\bigskip+\frac{2}{3}[X_{i}\left(  2\frac{\partial^{2}X_{j}}{\partial
x_{i}\partial x_{j}}+\frac{\partial^{2}X_{i}}{\partial x_{j}^{2}}\right)
+X_{j}\left(  \frac{\partial^{2}X_{j}}{\partial x_{i}^{2}}+2\frac{\partial
^{2}X_{i}}{\partial x_{i}\partial x_{j}}\right)  ](y_{0})-\frac{1}{2}\left(
\frac{\partial^{3}X_{i}}{\partial x_{i}\partial x_{j}^{2}}+\frac{\partial
^{3}X_{j}}{\partial x_{i}^{2}\partial x_{j}}\right)  (y_{0})$

\qquad\qquad\qquad\qquad\qquad\qquad\qquad\qquad\qquad\qquad\qquad\qquad
\qquad\qquad\qquad\qquad\qquad$\blacksquare$

The formula for $\frac{\partial^{4}\Phi_{P}}{\partial x_{i}^{4}}(y_{0})$ is
shorter and more elegant:\qquad

$\qquad\frac{\partial^{4}\Phi_{P}}{\partial x_{i}^{4}}(y_{0})=[\left\Vert
X\right\Vert _{M}^{4}](y_{0})-2[\left\Vert X\right\Vert _{M}^{2}%
\operatorname{div}_{M}X]](y_{0})+3[\operatorname{div}_{M}X]^{2}(y_{0})$

$\qquad\qquad-4[X_{i}^{2}\left(  \frac{\partial X_{i}}{\partial x_{i}}\right)
](y_{0})+4[X_{i}\left(  \frac{\partial^{2}X_{i}}{\partial x_{i}^{2}}\right)
](y_{0})-\left(  \frac{\partial^{3}X_{i}}{\partial x_{i}^{3}}\right)  (y_{0})$

\qquad\qquad\qquad\qquad\qquad\qquad\qquad\qquad\qquad\qquad\qquad\qquad
\qquad\qquad\qquad\qquad\qquad$\blacksquare$

\begin{claim}
The pattern in the formulae below give us a clue as to how higher order
derivatives should look.
\end{claim}

(i) $\frac{\partial\Phi_{P}}{\partial x_{i}}(y_{0})=-X_{i}(y_{0})$

(ii) $\frac{\partial^{2}\Phi_{P}}{\partial x_{i}\partial x_{j}}(y_{0}%
)=X_{i}(y_{0})X_{j}(y_{0})-\frac{1}{2}\left(  \frac{\partial X_{i}}{\partial
x_{j}}+\frac{\partial X_{j}}{\partial x_{i}}\right)  (y_{0})$

(v) $\qquad\frac{\partial^{3}\Phi_{P}}{\partial x_{i}\partial x_{j}\partial
x_{k}}(y)=-X_{i}(y)X_{j}(y)X_{k}(y)+\frac{1}{2}X_{i}(y)\left(  \frac{\partial
X_{k}}{\partial x_{j}}+\frac{\partial X_{j}}{\partial x_{k}}\right)  (y)$

$\qquad\qquad\qquad\qquad+\frac{1}{2}X_{j}(y)\left(  \frac{\partial X_{k}%
}{\partial x_{i}}+\frac{\partial X_{i}}{\partial x_{k}}\right)  (y)+\frac
{1}{2}X_{k}(y)\left(  \frac{\partial X_{j}}{\partial x_{i}}+\frac{\partial
X_{i}}{\partial x_{j}}\right)  (y)$

$\qquad\qquad\qquad\qquad-\frac{1}{3}\left(  \frac{\partial^{2}X_{i}}{\partial
x_{j}\partial x_{k}}+\frac{\partial^{2}X_{j}}{\partial x_{i}\partial x_{k}%
}+\frac{\partial^{2}X_{k}}{\partial x_{i}\partial x_{j}}\right)  (y)$

\qquad\qquad\qquad\qquad\qquad\qquad\qquad\qquad\qquad\qquad\qquad\qquad
\qquad\qquad\qquad\qquad\qquad$\blacksquare$

\subsection{\textbf{Tangential Derivatives:}}

For a, b = 1,...,q we have:

\qquad(vii)\qquad$\frac{\partial\Phi_{P}}{\partial x_{\text{a}}}(y_{0})=0$

\qquad(viii)$\qquad\frac{\partial^{2}\Phi_{P}}{\partial x_{\text{a}}\partial
x_{\text{b}}}(y_{0})=0$

\qquad(ix)\qquad$\ \frac{\partial\Phi_{P}^{-1}}{\partial x_{\text{a}}}%
(y_{0})=0$

\qquad(x)$\qquad\ \ \frac{\partial^{2}\Phi_{P}^{-1}}{\partial x_{\text{a}%
}\partial x_{\text{b}}}(y_{0})=0$

Higher derivatives are all equal to zero.

\subsection{\textbf{Mixed Derivatives }}

For a = 1,...,q and $i=q+1,...,n,$ we have:

\qquad(xi)\qquad$\frac{\partial^{2}\Phi_{P}}{\partial x_{\text{a}}\partial
x_{i}}(y_{0})=$ $\frac{\partial}{\partial x_{\text{a}}}(\nabla\log\Phi
_{P})_{i}(y_{0})=-\frac{\partial X_{i}}{\partial x_{\text{a}}}(y_{0})$

\qquad(xii)\qquad$\frac{\partial^{3}\Phi_{P}}{\partial x_{\text{a}}\partial
x_{\text{b}}\partial x_{i}}(x_{0})=\frac{\partial^{2}}{\partial x_{\text{a}%
}\partial x_{\text{b}}}(\nabla\log\Phi_{P})_{i}(y_{0})=-\frac{\partial
^{2}X_{i}}{\partial x_{\text{a}}\partial x_{\text{b}}}(y_{0})$

\qquad(xiii)$\qquad\frac{\partial^{3}\Phi_{P}}{\partial x_{\text{a}}%
^{2}\partial x_{i}}(y_{0})=\frac{\partial^{2}}{\partial x_{\text{a}}^{2}%
}(\nabla\log\Phi_{P})_{i}(y_{0})=-\frac{\partial^{2}X_{i}}{\partial
x_{\text{a}}^{2}}(y_{0})$

\qquad(xiv)\qquad$\frac{\partial^{3}\Phi_{P}}{\partial x_{\text{c}}\partial
x_{i}\partial x_{j}}(y_{0})$ \qquad\qquad\qquad\qquad\qquad\qquad\qquad from
$\left(  B_{88}\right)  $

$\qquad\qquad=[X_{i}\frac{\partial X_{j}}{\partial x_{\text{c}}}+X_{j}%
\frac{\partial X_{i}}{\partial x_{\text{c}}}](y_{0})+[X_{j}\frac{\partial
X_{i}}{\partial x_{\text{c}}}+X_{i}\frac{\partial X_{j}}{\partial x_{\text{c}%
}}](y_{0})-\frac{1}{2}\left(  \frac{\partial^{2}X_{i}}{\partial x_{\text{c}%
}\partial x_{j}}+\frac{\partial^{2}X_{j}}{\partial x_{\text{c}}\partial x_{i}%
}\right)  (y_{0})$

$\qquad\qquad=2[X_{i}\frac{\partial X_{j}}{\partial x_{\text{c}}}+X_{j}%
\frac{\partial X_{i}}{\partial x_{\text{c}}}](y_{0})-\frac{1}{2}\left(
\frac{\partial^{2}X_{i}}{\partial x_{\text{c}}\partial x_{j}}+\frac
{\partial^{2}X_{j}}{\partial x_{\text{c}}\partial x_{i}}\right)  (y_{0}%
)\qquad$

In particular,

$\qquad$(xiv)$^{\ast}$ $\frac{\partial^{3}\Phi_{P}}{\partial x_{\text{c}%
}\partial x_{i}^{2}}(y_{0})=$ $4X_{i}\frac{\partial X_{i}}{\partial
x_{\text{c}}}(y_{0})-\frac{\partial^{2}X_{i}}{\partial x_{\text{c}}\partial
x_{i}}(y_{0})$\qquad

\qquad(xv) $\frac{\partial^{4}\Phi_{P}}{\partial x_{\text{c}}^{2}\partial
x_{i}\partial x_{j}}(y_{0})=$ $2\frac{\partial X_{i}}{\partial x_{\text{c}}%
}(y_{0})\frac{\partial X_{j}}{\partial x_{\text{c}}}(y_{0})+X_{i}(y_{0}%
)\frac{\partial^{2}X_{j}}{\partial x_{\text{c}}^{2}}(y_{0})$

$\qquad\qquad+X_{j}(y_{0})\frac{\partial^{2}X_{i}}{\partial x_{\text{c}}^{2}%
}(y_{0})$ $-\frac{\partial^{3}X_{i}}{\partial x_{\text{c}}^{2}\partial x_{j}%
}(y_{0})\qquad\qquad\qquad\qquad\qquad$

In particular,

\qquad(xvi)$\qquad\frac{\partial^{4}\Phi_{P}}{\partial x_{\text{c}}%
^{2}\partial x_{i}^{2}}(y_{0})=$ $2[(\frac{\partial X_{i}}{\partial
x_{\text{c}}})^{2}+X_{i}\frac{\partial^{2}X_{i}}{\partial x_{\text{c}}^{2}%
}](y_{0})$ $-\frac{\partial^{3}X_{i}}{\partial x_{\text{c}}^{2}\partial x_{i}%
}(y_{0})\qquad$from $\left(  B_{99}\right)  \qquad\qquad$

\qquad(xvii)\qquad$\frac{1}{2}\Delta\Phi(y_{0}%
)=\underset{i=q+1}{\overset{n}{\sum}}\frac{1}{2}X_{i}^{2}(y_{0})-\frac{1}%
{2}\underset{i=q+1}{\overset{n}{\sum}}\frac{\partial X_{i}}{\partial x_{i}%
}(y_{0})$

\qquad\qquad$=\frac{1}{2}$ $\left\Vert X\right\Vert ^{2}(y_{0})-$ $\frac{1}%
{2}\operatorname{div}X(y_{0})-\frac{1}{2}\underset{\text{a}%
=1}{\overset{q}{\sum}}X_{\text{a}}^{2}(y_{0})+\frac{1}{2}$ $\underset{\text{a}%
=1}{\overset{q}{\sum}}\frac{\partial X_{\text{a}}}{\partial x_{\text{a}}%
}(y_{0})$

\qquad\qquad$=\frac{1}{2}$ $\left\Vert X\right\Vert _{M}^{2}(y_{0})-\frac
{1}{2}\operatorname{div}_{M}X(y_{0})-\frac{1}{2}\left\Vert X\right\Vert
_{P}^{2}(y_{0})+\frac{1}{2}\operatorname{div}_{P}X(y_{0})(y_{0})$

This ties up with the expression of $\left(  B_{26}\right)  $ in $\left(
21\right)  $ which proves (iii) of \textbf{Table B}$_{3}.$

$\qquad$(xviii)$\qquad\frac{\partial^{4}\Phi_{P}}{\partial x_{i}^{2}\partial
x_{j}\partial x_{k}}(x_{0})\qquad$from$\qquad\left(  114\right)  $

$=-X_{j}(y_{0})[-X_{i}^{2}X_{k}+X_{k}\frac{\partial X_{i}}{\partial x_{i}%
}+X_{i}\left(  \frac{\partial X_{i}}{\partial x_{k}}+\frac{\partial X_{k}%
}{\partial x_{i}}\right)  -\frac{1}{3}\left(  \frac{\partial^{2}X_{k}%
}{\partial x_{i}^{2}}+2\frac{\partial^{2}X_{i}}{\partial x_{i}\partial x_{k}%
}\right)  ](y_{0})\qquad\frac{\partial^{2}L_{1}}{\partial x_{i}\partial x_{k}%
}(y_{0})$

$-\frac{1}{2}[X_{i}^{2}-\frac{\partial X_{i}}{\partial x_{i}}](y_{0})\left(
\frac{\partial X_{j}}{\partial x_{k}}+\frac{\partial X_{k}}{\partial x_{j}%
}\right)  (y_{0})\qquad\qquad$

$-\frac{1}{2}[X_{i}X_{k}-\frac{1}{2}\left(  \frac{\partial X_{k}}{\partial
x_{i}}+\frac{\partial X_{i}}{\partial x_{k}}\right)  ](y_{0})[\left(
\frac{\partial X_{j}}{\partial x_{i}}+\frac{\partial X_{i}}{\partial x_{j}%
}\right)  ](y_{0})$

$-\frac{2}{3}[X_{i}X_{q}](y_{0})[(R_{ijkq}+R_{kjiq})](y_{0})\qquad\qquad
\qquad\qquad\ \ \left(  1\right)  $

$+\frac{2}{3}X_{i}(y_{0})[\left(  \frac{\partial^{2}X_{k}}{\partial
x_{i}\partial x_{j}}+\frac{\partial^{2}X_{j}}{\partial x_{i}\partial x_{k}%
}+\frac{\partial^{2}X_{i}}{\partial x_{j}\partial x_{k}}\right)  +$
$\underset{l=1}{\overset{n}{\sum}}(R_{ikjl}+R_{jkil})X_{l}](y_{0}%
)\qquad\left(  2\right)  \qquad$

$-\frac{2}{3}[X_{k}X_{q}](y_{0})\varrho_{jq}(y_{0})\qquad\frac{\partial
^{2}L_{2}}{\partial x_{i}\partial x_{k}}(y_{0})$

$+\frac{1}{3}[\nabla_{k}R_{ijiq}+\nabla_{i}R_{kjiq}+\nabla_{i}R_{ijkq}%
](y_{0})X_{q}(y_{0})+\frac{1}{3}\varrho_{jq}(y_{0})\left(  \frac{\partial
Xq}{\partial x_{k}}+\frac{\partial X_{k}}{\partial x_{q}}\right)  ](y_{0})$

$+\frac{1}{3}[(R_{ijkq}+R_{kjiq})\left(  \frac{\partial X_{q}}{\partial x_{i}%
}+\frac{\partial X_{i}}{\partial x_{q}}\right)  ](y_{0})\qquad\left(
3\right)  $

$-\frac{1}{2}[X_{i}X_{k}-\frac{1}{2}\left(  \frac{\partial X_{k}}{\partial
x_{i}}+\frac{\partial X_{i}}{\partial x_{k}}\right)  \left(  \frac{\partial
X_{j}}{\partial x_{i}}+\frac{\partial X_{i}}{\partial x_{j}}\right)
](y_{0})\qquad\frac{\partial^{2}L_{3}}{\partial x_{i}\partial x_{k}}%
(y_{0})\qquad$

$+\frac{1}{3}X_{k}(y_{0})\left(  \frac{\partial^{2}X_{j}}{\partial x_{i}^{2}%
}+2\frac{\partial^{2}X_{i}}{\partial x_{i}\partial x_{j}}\right)
(y_{0})+\frac{2}{3}\underset{l=1}{\overset{n}{\sum}}\varrho_{jl}(y_{0}%
)X_{k}(y_{0})X_{l}(y_{0})$

$-\frac{1}{4}[\frac{\partial^{3}X_{j}}{\partial x_{i}^{2}\partial x_{k}}%
+\frac{\partial^{3}X_{k}}{\partial x_{i}^{2}\partial x_{j}}+2\frac
{\partial^{3}X_{i}}{\partial x_{i}\partial x_{j}\partial x_{k}}](y_{0})$

(xviii)$^{\ast}\frac{\partial^{4}\Phi_{P}}{\partial x_{i}\partial
x_{j}\partial x_{k}^{2}}(x_{0})=\qquad$from $\left(  115\right)  $

$-X_{j}(y_{0})[-X_{k}^{2}X_{i}+X_{i}\frac{\partial X_{k}}{\partial x_{k}%
}+X_{i}\left(  \frac{\partial X_{k}}{\partial x_{i}}+\frac{\partial X_{i}%
}{\partial x_{k}}\right)  -\frac{1}{3}\left(  \frac{\partial^{2}X_{i}%
}{\partial x_{k}^{2}}+2\frac{\partial^{2}X_{k}}{\partial x_{k}\partial x_{i}%
}\right)  ](y_{0})\qquad$

$-\frac{1}{2}[X_{k}^{2}-\frac{\partial X_{k}}{\partial x_{k}}](y_{0})\left(
\frac{\partial X_{j}}{\partial x_{i}}+\frac{\partial X_{i}}{\partial x_{j}%
}\right)  (y_{0})-\frac{1}{2}[X_{i}X_{k}-\frac{1}{2}\left(  \frac{\partial
X_{k}}{\partial x_{i}}+\frac{\partial X_{i}}{\partial x_{k}}\right)
](y_{0})[\left(  \frac{\partial X_{j}}{\partial x_{k}}+\frac{\partial X_{k}%
}{\partial x_{j}}\right)  ](y_{0})$

$-\frac{2}{3}[X_{k}X_{q}](y_{0})[R_{kjiq}+R_{ijkq}](y_{0})+\frac{2}{3}%
X_{k}(y_{0})[\left(  \frac{\partial^{2}X_{k}}{\partial x_{i}\partial x_{j}%
}+\frac{\partial^{2}X_{j}}{\partial x_{i}\partial x_{k}}+\frac{\partial
^{2}X_{i}}{\partial x_{j}\partial x_{k}}\right)  +$ $(R_{kijq}+R_{jikq}%
)X_{q}](y_{0})$

$-\frac{2}{3}[X_{i}X_{q}](y_{0})\varrho_{jq}(y_{0})+\frac{1}{3}[\nabla
_{i}R_{kjkq}+\nabla_{k}R_{ijkq}+\nabla_{k}R_{kjiq}](y_{0})X_{q}(y_{0}%
)+\frac{1}{3}\varrho_{jq}(y_{0})\left(  \frac{\partial Xq}{\partial x_{i}%
}+\frac{\partial X_{i}}{\partial x_{q}}\right)  ](y_{0})$

$+\frac{1}{3}[(R_{ijkq}+R_{kjiq})\left(  \frac{\partial X_{q}}{\partial x_{k}%
}+\frac{\partial X_{k}}{\partial x_{q}}\right)  ](y_{0})-\frac{1}{2}%
[X_{i}X_{k}-\frac{1}{2}\left(  \frac{\partial X_{k}}{\partial x_{i}}%
+\frac{\partial X_{i}}{\partial x_{k}}\right)  \left(  \frac{\partial X_{j}%
}{\partial x_{k}}+\frac{\partial X_{k}}{\partial x_{j}}\right)  ](y_{0})$

$+\frac{1}{3}X_{i}(y_{0})\left(  \frac{\partial^{2}X_{j}}{\partial x_{k}^{2}%
}+2\frac{\partial^{2}X_{k}}{\partial x_{k}\partial x_{j}}\right)
(y_{0})+\frac{2}{3}\varrho_{jl}(y_{0})X_{i}(y_{0})X_{l}(y_{0})$

$-\frac{1}{4}[\frac{\partial^{3}X_{j}}{\partial x_{k}^{2}\partial x_{i}}%
+\frac{\partial^{3}X_{i}}{\partial x_{k}^{2}\partial x_{j}}+2\frac
{\partial^{3}X_{k}}{\partial x_{i}\partial x_{j}\partial x_{k}}](y_{0})$

\qquad\qquad\qquad\qquad\qquad\qquad\qquad\qquad\qquad\qquad\qquad\qquad
\qquad\qquad$\blacksquare$

\subsection{\textbf{Computations of Appendix B}$_{4}$}

\subsubsection{NORMAL DERIVATIVES\qquad\qquad}

\qquad\qquad\qquad\qquad\qquad\qquad\qquad\qquad\qquad\qquad

(i) We use the definition of the gradient operator at a general point
x$_{0}\in$M$_{0}:$

For $j,k=1,...,q,q+1,...,n,$

$(\nabla\log\Phi_{P})_{k}(x_{0})=$ g$^{jk}(x_{0})\frac{\partial}{\partial
x_{j}}\log\Phi_{P}(x_{0})$

\qquad$(\nabla\log\Phi_{P})_{k}(x_{0})=$ g$^{jk}(x_{0})\frac{\partial
}{\partial x_{j}}\log\Phi_{P}(x_{0})\qquad\qquad\qquad\qquad(27)\qquad
\qquad\qquad\qquad\qquad\qquad\ \ \ \ \ $

and so,

\qquad$(\nabla\log\Phi_{P})_{k}(y_{0})=\frac{\partial}{\partial x_{k}}\log
\Phi_{P}(y_{0})\qquad\qquad\qquad\qquad\qquad\qquad\qquad\ \ $

Therefore for $k=q+1,...,n,$ we have:

\qquad$(\nabla\log\Phi_{P})_{k}(y_{0})=\frac{\partial}{\partial x_{k}}\log
\Phi_{P}(y_{0})=-X_{k}(y_{0})\qquad\qquad\qquad(28)\qquad\ \ \ \ \qquad
\qquad\qquad\ $

From Elementary Calculus, we have:

\qquad$\frac{\partial\Phi_{P}}{\partial x_{i}}(x_{0})=\Phi_{P}(x_{0}%
)\frac{\partial}{\partial x_{i}}\log\Phi_{P}(x_{0})\qquad\qquad\qquad
\qquad\qquad\ \ \ \ \ (29)\qquad$

$\qquad\qquad\qquad\ \ \ \ \ \ \qquad\qquad\qquad\ \ $

Since $\Phi_{P}(y_{0})=1$, we have by $\left(  28\right)  $ and $\left(
29\right)  $ that$:$

\qquad$\frac{\partial\Phi_{P}}{\partial x_{i}}(y_{0})=\frac{\partial}{\partial
x_{i}}\log\Phi_{P}(y_{0})=$ $(\nabla$log$\Phi_{P})_{i}(y_{0})=-X_{i}%
(y_{0})\qquad\ \left(  30\right)  \qquad\qquad\qquad\qquad\ \ \qquad
\ \ \ \ \qquad\qquad\qquad\qquad\qquad\qquad\qquad\qquad\qquad\qquad
\qquad\qquad\qquad\qquad\qquad\qquad$

Concluding, we have:

$\qquad\frac{\partial\Phi_{P}}{\partial x_{i}}(y_{0})=-X_{i}(y_{0}%
)\qquad\qquad\qquad\qquad\qquad\qquad\qquad\qquad\qquad\ (31)$

(ii) By the definition of the gradient operator in $\left(  27\right)  $, we have:

\qquad g$_{ik}(x_{0})(\nabla$log$\Phi_{P})_{k}(x_{0})=$ g$_{ik}(x_{0})$%
g$^{jk}(x_{0})\frac{\partial}{\partial x_{j}}$log$\Phi_{P}(x_{0}%
)=\frac{\partial}{\partial x_{i}}$log$\Phi_{P}(x_{0})$

We have:

\qquad$\frac{\partial}{\partial x_{i}}$log$\Phi_{P}(x_{0})=$ g$_{ik}%
(x_{0})(\nabla$log$\Phi_{P})_{k}(x_{0})$

From the last equation above and the relation in $\left(  29\right)  ,$ we have:

\qquad$\frac{\partial\Phi_{P}}{\partial x_{j}}(x_{0})=\Phi_{P}(x_{0}%
)\frac{\partial}{\partial x_{j}}\log\Phi_{P}(x_{0})=\Phi_{P}(x_{0})$%
g$_{jk}(x_{0})(\nabla$log$\Phi_{P})_{k}(x_{0})$Conequently,

$\qquad\frac{\partial^{2}\Phi_{P}}{\partial x_{i}\partial x_{j}}(x_{0}%
)=\frac{\partial}{\partial x_{i}}[\Phi_{P}$g$_{jk}(\nabla$log$\Phi_{P}%
)_{k}](x_{0})\qquad\qquad\qquad\qquad\qquad\qquad\qquad\qquad$

\qquad$=\frac{\partial\Phi_{P}}{\partial x_{i}}(x_{0})$g$_{jk}(x_{0})(\nabla
$log$\Phi_{P})_{k}(x_{0})\qquad\qquad\qquad\qquad\qquad\qquad(32)$

$\qquad+\Phi_{P}(x_{0})[\frac{\partial\text{g}_{jk}}{\partial x_{i}}%
(x_{0})(\nabla$log$\Phi_{P})_{k}(x_{0})+$ g$_{jk}(x_{0})\frac{\partial
}{\partial x_{i}}(\nabla$log$\Phi_{P})_{k}(x_{0})]$\ \ \ 

Since $\Phi_{P}(y_{0})=1,$ g$_{jk}(y_{0})=\delta_{jk}$ and $\frac{\partial
g_{jk}}{\partial x_{i}}(y_{0})=0$ for $i,j,k=q+1,...,n,$

by (ii) of \textbf{Table 1 in Appendix A}, we have:

\qquad$\frac{\partial^{2}\Phi_{P}}{\partial x_{i}\partial x_{j}}(y_{0}%
)=\frac{\partial\Phi_{P}}{\partial x_{i}}(y_{0})(\nabla$log$\Phi_{P}%
)_{j}(y_{0})+$ $\frac{\partial}{\partial x_{i}}(\nabla$log$\Phi_{P})_{j}%
(y_{0})\qquad\qquad\qquad\ $

\qquad\qquad\qquad$\ =X_{i}(y_{0})X_{j}(y_{0})+\frac{\partial}{\partial x_{i}%
}(\nabla\log\Phi_{P})_{j}(y_{0})\qquad\qquad\qquad(33)$

Similarly we have:

\qquad$\frac{\partial^{2}\Phi_{P}}{\partial x_{j}\partial x_{i}}(y_{0}%
)=X_{j}(y_{0})X_{i}(y_{0})+\frac{\partial}{\partial x_{j}}(\nabla\log\Phi
_{P})_{i}(y_{0})\qquad\qquad\qquad\left(  34\right)  $

Therefore by $\left(  33\right)  $ and $\left(  34\right)  ,$

\qquad$\frac{\partial^{2}\Phi_{P}}{\partial x_{i}\partial x_{j}}(y_{0}%
)+\frac{\partial^{2}\Phi_{P}}{\partial x_{j}\partial x_{i}}(y_{0})=X_{i}%
(y_{0})X_{j}(y_{0})+\frac{\partial}{\partial x_{i}}(\nabla\log\Phi_{P}%
)_{j}(y_{0})$

$\qquad+X_{j}(y_{0})X_{i}(y_{0})+\frac{\partial}{\partial x_{j}}(\nabla
\log\Phi_{P})_{i}(y_{0})$\qquad\qquad\qquad\qquad

Since $\Phi_{P}:M\rightarrow R$ is a smooth function,

$\frac{\partial^{2}\Phi_{P}}{\partial x_{i}\partial x_{j}}(y_{0}%
)=\frac{\partial^{2}\Phi_{P}}{\partial x_{j}\partial x_{i}}(y_{0})$

and since

$X_{i}(y_{0})X_{j}(y_{0})=X_{j}(y_{0})X_{i}(y_{0})$ ,

we have:

$2\frac{\partial^{2}\Phi_{P}}{\partial x_{i}\partial x_{j}}(y_{0}%
)=2X_{i}(y_{0})X_{j}(y_{0})+\frac{\partial}{\partial x_{i}}(\nabla\log\Phi
_{P})_{j}(y_{0})+\frac{\partial}{\partial x_{j}}(\nabla\log\Phi_{P})_{i}%
(y_{0})\qquad$

By (vii) of Table B$_{1}$, we have:

$\frac{\partial}{\partial x_{i}}(\nabla\log\Phi_{P})_{j}(y_{0})+\frac
{\partial}{\partial x_{j}}(\nabla\log\Phi_{P})_{i}(y_{0})=-\left(
\frac{\partial X_{j}}{\partial x_{i}}+\frac{\partial X_{i}}{\partial x_{j}%
}\right)  (y_{0})\qquad\left(  35\right)  $

Therefore from the last two equations,

$2\frac{\partial^{2}\Phi_{P}}{\partial x_{i}\partial x_{j}}(y_{0}%
)=2X_{i}(y_{0})X_{j}(y_{0})-\left(  \frac{\partial X_{j}}{\partial x_{i}%
}+\frac{\partial X_{i}}{\partial x_{j}}\right)  (y_{0})$

Hence, we have:

$\frac{\partial^{2}\Phi_{P}}{\partial x_{i}\partial x_{j}}(y_{0})=X_{i}%
(y_{0})X_{j}(y_{0})-\frac{1}{2}\left(  \frac{\partial X_{j}}{\partial x_{i}%
}+\frac{\partial X_{i}}{\partial x_{j}}\right)  (y_{0})\qquad\qquad
\qquad\ \ \ \ \ \ \left(  36\right)  $

\qquad\qquad\qquad\qquad\qquad\qquad\qquad\qquad\qquad\qquad\qquad\qquad
\qquad\qquad\qquad\qquad\qquad\qquad\qquad

From $(33)$ or $\left(  34\right)  $ and $\left(  36\right)  $ we see that:

$\frac{\partial}{\partial x_{i}}(\nabla\log\Phi_{P})_{j}(y_{0})=-\frac{1}%
{2}\left(  \frac{\partial X_{j}}{\partial x_{i}}+\frac{\partial X_{i}%
}{\partial x_{j}}\right)  (y_{0})\qquad\qquad\qquad\qquad\qquad\left(
37\right)  $

\qquad\qquad\qquad\qquad\qquad\qquad\qquad\qquad\qquad\qquad\qquad\qquad
\qquad\qquad\qquad\qquad\qquad\qquad\qquad

From the last equation above, it is immediate (by inter-changing indices) that:

\qquad$\frac{\partial}{\partial x_{j}}(\nabla\log\Phi_{P})_{i}(y_{0}%
)=-\frac{1}{2}\left(  \frac{\partial X_{i}}{\partial x_{j}}+\frac{\partial
X_{j}}{\partial x_{i}}\right)  (y_{0})\qquad\qquad\qquad\qquad\left(
38\right)  $

We conclude from $\left(  37\right)  $ and $\left(  38\right)  $ that for
$i,j=q+1,...,n,$ we have:

\qquad$\frac{\partial}{\partial x_{i}}(\nabla\log\Phi_{P})_{j}(y_{0}%
)=-\frac{1}{2}\left(  \frac{\partial X_{j}}{\partial x_{i}}+\frac{\partial
X_{i}}{\partial x_{j}}\right)  (y_{0})=\frac{\partial}{\partial x_{j}}%
(\nabla\log\Phi_{P})_{i}(y_{0})$ $\left(  39\right)  $\qquad\qquad\qquad
\qquad\qquad\qquad\qquad\qquad\qquad\qquad\qquad\qquad\qquad\qquad\qquad
\qquad\qquad\qquad\qquad

The very important formula in $\left(  39\right)  $ can also be proved using
simple Calculus:

Since

$\frac{\partial^{2}\Phi_{P}}{\partial x_{j}\partial x_{i}}(y_{0}%
)=\frac{\partial^{2}\Phi_{P}}{\partial x_{i}\partial x_{j}}(y_{0});$
$X_{i}(y_{0})X_{j}(y_{0})=X_{j}(y_{0})X_{i}(y_{0}),$

we see from $\left(  33\right)  $ and $\left(  34\right)  $ that:

\qquad\ $\frac{\partial}{\partial x_{j}}(\nabla$log$\Phi_{P})_{i}(y_{0})=$
$\frac{\partial}{\partial x_{i}}(\nabla$log$\Phi_{P})_{j}(y_{0})\qquad
\qquad\qquad\qquad\qquad\ \ \ \left(  40\right)  $

\qquad\qquad\qquad\qquad\qquad\qquad\qquad\qquad\qquad\qquad\qquad\qquad
\qquad\qquad\qquad\qquad\qquad$\qquad\qquad\blacksquare$

(ii)$^{\ast\ast}$ Then $\left(  35\right)  $ and $\left(  40\right)  $ give:

\qquad$\frac{\partial}{\partial x_{j}}(\nabla$log$\Phi_{P})_{i}(y_{0}%
)=-\frac{1}{2}\left(  \frac{\partial X_{j}}{\partial x_{i}}+\frac{\partial
X_{i}}{\partial x_{j}}\right)  (y_{0})=$ $\frac{\partial}{\partial x_{i}%
}(\nabla$log$\Phi_{P})_{j}(y_{0})$

(iii) Since $\Phi\Phi^{-1}=1,$ we have:

$\qquad\frac{\partial\Phi_{P}}{\partial x_{i}}\Phi^{-1}+\Phi\frac{\partial
\Phi_{P}^{-1}}{\partial x_{i}}=0\qquad\qquad\qquad\qquad\qquad\qquad
\qquad\qquad\qquad\qquad\qquad\qquad$

Therefore by $\left(  31\right)  $ and the fact that $\Phi(y_{0})=1=\Phi
^{-1}(y_{0})$, we have:

$\qquad\frac{\partial\Phi_{P}^{-1}}{\partial x_{i}}(y_{0})=-\frac{\partial
\Phi_{P}}{\partial x_{i}}(y_{0})=X_{i}(y_{0}).$

\bigskip(iv) It is obvious that:

$\qquad\frac{\partial\Phi_{P}^{-1}}{\partial x_{i}}=-\frac{\partial\Phi_{P}%
}{\partial x_{i}}\Phi^{-2}$

$\qquad\frac{\partial}{\partial x_{j}}[\frac{\partial\Phi_{P}^{-1}}{\partial
x_{i}}]=-\frac{\partial}{\partial x_{j}}[\frac{\partial\Phi_{P}}{\partial
x_{i}}\Phi^{-2}]$

$\qquad=-\frac{\partial}{\partial xj}[\frac{\partial\Phi_{P}}{\partial x_{i}%
}]\Phi^{-2}-\frac{\partial\Phi_{P}}{\partial x_{i}}\frac{\partial}{\partial
xj}\Phi^{-2}=-\frac{\partial^{2}\Phi_{P}}{\partial x_{j}\partial x_{i}}%
\Phi^{-2}+2\frac{\partial\Phi_{P}}{\partial x_{i}}\frac{\partial\Phi}{\partial
xj}\Phi^{-3}$

$\qquad\frac{\partial^{2}\Phi_{P}^{-1}}{\partial x_{i}\partial x_{j}}%
=-\frac{\partial^{2}\Phi_{P}}{\partial x_{i}\partial x_{j}}\Phi^{-2}%
+2\frac{\partial\Phi_{P}}{\partial x_{i}}\frac{\partial\Phi}{\partial x_{j}%
}\Phi^{-3}$

Since $\Phi^{-2}(y_{0})=1=\Phi^{-3}(y_{0})$ and $\frac{\partial\Phi_{P}%
}{\partial x_{i}}(y_{0})=-X_{i}(y_{0}),$ we have by $\left(  36\right)  :$

$\qquad\frac{\partial^{2}\Phi_{P}^{-1}}{\partial x_{i}\partial x_{j}}%
(y_{0})=2\frac{\partial\Phi_{P}}{\partial x_{i}}(y_{0})\frac{\partial\Phi
}{\partial x_{j}}(y_{0})-\frac{\partial^{2}\Phi_{P}}{\partial x_{i}\partial
x_{j}}(y_{0})=X_{i}(y_{0})X_{j}(y_{0})+\frac{1}{2}\left(  \frac{\partial
X_{j}}{\partial x_{i}}+\frac{\partial X_{i}}{\partial x_{j}}\right)
(y_{0})\qquad\qquad\qquad\qquad\ \ \ \qquad$

(v) We change $k$ to $l$ on the RHS of $(32)$ and have:\ 

$\qquad\frac{\partial^{2}\Phi_{P}}{\partial x_{j}\partial x_{i}}%
(x_{0})\ =\frac{\partial\Phi_{P}}{\partial x_{j}}(x_{0})$g$_{il}(x_{0}%
)(\nabla$log$\Phi_{P})_{l}(x_{0})\qquad\qquad$

$\qquad\qquad+\Phi_{P}(x_{0})\frac{\partial\text{g}_{il}}{\partial x_{j}%
}(x_{0})(\nabla$log$\Phi_{P})_{l}(x_{0})+$ $\Phi_{P}(x_{0})$g$_{il}%
(x_{0})\frac{\partial}{\partial x_{j}}(\nabla$log$\Phi_{P})_{l}(x_{0})$

$\qquad\qquad=$ L$_{1}$ + L$_{2}$ + L$_{3}$\ 

Therefore,

$\qquad\frac{\partial^{3}\Phi_{P}}{\partial x_{k}\partial x_{j}\partial x_{i}%
}(x_{0})=\frac{\partial}{\partial x_{k}}$L$_{1}(x_{0})+\frac{\partial
}{\partial x_{k}}$L$_{2}(x_{0})+\frac{\partial}{\partial x_{k}}$L$_{3}%
(x_{0})\qquad\qquad\qquad\qquad\left(  41\right)  $

where,

\qquad$\frac{\partial}{\partial x_{k}}$L$_{1}(x_{0})=$ $\frac{\partial
}{\partial x_{k}}[\frac{\partial\Phi_{P}}{\partial x_{j}}$g$_{il}(\nabla
$log$\Phi_{P})_{l}](x_{0})\qquad\qquad\qquad\qquad\qquad\qquad\qquad\left(
42\right)  $

$\ \ \ \ =$ $\frac{\partial^{2}\Phi_{P}}{\partial x_{k}\partial x_{j}}%
(x_{0})[$g$_{il}(\nabla$log$\Phi_{P})_{l}](x_{0})+$ $\frac{\partial\Phi_{P}%
}{\partial x_{j}}(x_{0})[\frac{\partial\text{g}_{il}}{\partial x_{k}}(\nabla
$log$\Phi_{P})_{l}+$ g$_{il}\frac{\partial}{\partial x_{k}}(\nabla\log\Phi
_{P})_{l}](x_{0})\qquad$

Since $(\nabla$log$\Phi_{P})_{\text{a}}(y)=0$ for a = 1,...,q and
$\frac{\partial\text{g}_{il}}{\partial x_{k}}(y)=0$ for $i,k,l=q+1,...,n,$ we have,

$\qquad\frac{\partial\text{g}_{il}}{\partial x_{k}}(y)(\nabla$log$\Phi
_{P})_{l}(y)=0.$

Consequently,

$\qquad\frac{\partial}{\partial x_{k}}$L$_{1}(y_{0})=\frac{\partial^{2}%
\Phi_{P}}{\partial x_{j}\partial x_{k}}(y)[(\nabla$log$\Phi_{P})_{i}](y)+$
$\frac{\partial\Phi_{P}}{\partial x_{j}}(y)[$ $\frac{\partial}{\partial x_{k}%
}(\nabla\log\Phi_{P})_{i}](y)$

By (i) of Table B$_{1},$ and then by $\left(  36\right)  $ and $\left(
37\right)  $ above, we have:

$\qquad(\nabla$log$\Phi_{P})_{i}](y)=-X_{i}(y)$

Therefore for $i,j,k=q+1,...,n,$

$\frac{\partial}{\partial x_{k}}$L$_{1}(y_{0})=-X_{i}(y)\left[  X_{j}%
X_{k}-\frac{1}{2}\left(  \frac{\partial X_{k}}{\partial x_{j}}+\frac{\partial
X_{j}}{\partial x_{k}}\right)  \right]  (y)+\frac{1}{2}X_{j}(y)\left(
\frac{\partial X_{i}}{\partial x_{k}}+\frac{\partial X_{k}}{\partial x_{i}%
}\right)  (y)\qquad\left(  43\right)  $

$\qquad=-X_{i}(y)X_{j}(y)X_{k}(y)+\frac{1}{2}X_{i}(y)\left(  \frac{\partial
X_{k}}{\partial x_{j}}+\frac{\partial X_{j}}{\partial x_{k}}\right)
(y)+\frac{1}{2}X_{j}(y)\left(  \frac{\partial X_{i}}{\partial x_{k}}%
+\frac{\partial X_{k}}{\partial x_{i}}\right)  (y)\ \ $

$\qquad\qquad\qquad\qquad\qquad\qquad\qquad\qquad\qquad\qquad\qquad
\qquad\qquad\qquad\qquad\qquad\qquad\qquad\qquad\blacksquare\ \ \ \ \ \ \qquad
\qquad\qquad\qquad\qquad\qquad\ \ \ \ \ \ $

Next we have:

$\qquad\frac{\partial}{\partial x_{k}}$L$_{2}(x_{0})=\frac{\partial}{\partial
x_{k}}[\Phi_{P}\frac{\partial\text{g}_{il}}{\partial x_{j}}(\nabla\log\Phi
_{P})_{l}](x_{0})\qquad\qquad\qquad\qquad\qquad\qquad\qquad\qquad\ \left(
44\right)  \qquad\qquad\qquad\qquad\qquad\qquad\qquad\ \ \ $

$=\frac{\partial\Phi_{P}}{\partial x_{k}}(x_{0})[\frac{\partial\text{g}_{il}%
}{\partial x_{j}}(\nabla$log$\Phi_{P})_{l}](x_{0})+\Phi_{P}(x_{0}%
)\frac{\partial}{\partial x_{k}}[\frac{\partial\text{g}_{il}}{\partial x_{j}%
}(\nabla$log$\Phi_{P})_{l}](x_{0})$

$=\frac{\partial\Phi_{P}}{\partial x_{k}}(x_{0})[\frac{\partial\text{g}_{il}%
}{\partial x_{j}}(\nabla$log$\Phi_{P})_{l}](x_{0})+\Phi_{P}(x_{0}%
)[\frac{\partial^{2}\text{g}_{il}}{\partial x_{k}\partial x_{j}}(\nabla
\log\Phi_{P})_{l}+\frac{\partial\text{g}_{il}}{\partial x_{j}}\frac{\partial
}{\partial x_{k}}(\nabla\log\Phi_{P})_{l}](x_{0})\qquad$

Summing over repeated indices as usual, we have for a = 1,...,q and
$i,j,k,l=q+1,...,n:\qquad$

$\frac{\partial}{\partial x_{k}}$L$_{2}(y)$

$=\frac{\partial\Phi_{P}}{\partial x_{k}}(y)[\frac{\partial\text{g}%
_{i\text{a}}}{\partial x_{j}}(\nabla$log$\Phi_{P})_{\text{a}}](y)+\Phi
_{P}(y)[\frac{\partial^{2}\text{g}_{i\text{a}}}{\partial x_{k}\partial x_{j}%
}(\nabla\log\Phi_{P})_{\text{a}}+\frac{\partial\text{g}_{i\text{a}}}{\partial
x_{j}}\frac{\partial}{\partial x_{k}}(\nabla\log\Phi_{P})_{\text{a}}%
](y)\qquad$

$+\frac{\partial\Phi_{P}}{\partial x_{k}}(y)[\frac{\partial\text{g}_{il}%
}{\partial x_{j}}(\nabla$log$\Phi_{P})_{l}](y)+\Phi_{P}(y)[\frac{\partial
^{2}\text{g}_{il}}{\partial x_{k}\partial x_{j}}(\nabla\log\Phi_{P})_{l}%
+\frac{\partial\text{g}_{il}}{\partial x_{j}}\frac{\partial}{\partial x_{k}%
}(\nabla\log\Phi_{P})_{l}](y)\qquad\qquad\qquad\qquad\qquad\qquad\qquad$\qquad

Since

$\Phi_{P}(y)=1,\frac{\partial\Phi_{P}}{\partial x_{k}}(y)=-X_{k}(y),$
$(\nabla$log$\Phi_{P})_{\text{a}}](y)$ $=0$ for a=1,...,q,

$\frac{\partial\text{g}_{il}}{\partial x_{j}}(y)=0$ for $i,j,l=q+1,...,n,$

$\frac{\partial\text{g}_{\text{a}i}}{\partial\text{x}_{j}}(y)=$ $\perp
_{\text{a}ji}(y),$ $\frac{\partial^{2}g^{k\text{a}}}{\partial x_{i}\partial
x_{j}}(y)=\frac{4}{3}(R_{i\text{a}jk}+R_{j\text{a}ik}%
)(y)+4\underset{\text{b=1}}{\overset{\text{q}}{\sum}}T_{\text{ab}i}%
(y)\perp_{\text{b}kj}(y)$

For $i,j,l=q+1,...,n$ we have by (iv) of A$_{4}:$

$\frac{\partial}{\partial x_{k}}(\nabla\log\Phi_{P})_{\text{a}}](y)=-X_{l}%
(y)\perp_{\text{a}kl}(y)-\frac{\partial X_{k}}{\partial x_{\text{a}}}(y)$

and by (xv) of B$_{1}:$

$(\nabla\log\Phi_{P})_{l}(y)=-X_{l}(y)$

by (iii) of \textbf{Table A}$_{1}$

$\frac{\partial\text{g}_{il}}{\partial x_{j}}(y)=0;$ $\frac{\partial
^{2}\text{g}_{il}}{\partial x_{k}\partial x_{j}}(y_{0})=\frac{1}{3}%
(R_{ijkl}+R_{ikjl})(y_{0})$ for $i,j,k,l=q+1,...,n$

Therefore,

$\frac{\partial}{\partial x_{k}}$L$_{2}(y)=[-$ $X_{l}\perp_{\text{a}ji}%
\perp_{\text{a}kl}-\perp_{\text{a}ji}\frac{\partial X_{k}}{\partial
x_{\text{a}}}](y)-\frac{1}{3}(R_{ijkl}+R_{ikjl})(y)X_{l}(y)\qquad\left(
45\right)  \qquad$

$\qquad\qquad\qquad=[\perp_{\text{a}ij}\perp_{\text{a}kl}X_{l}+\perp
_{\text{a}ij}\frac{\partial X_{k}}{\partial x_{\text{a}}}](y)-\frac{1}%
{3}(R_{ijkl}+R_{ikjl})(y)X_{l}(y)\qquad$

$\qquad\qquad\qquad\qquad\qquad\qquad\qquad\qquad\qquad\qquad\qquad
\qquad\qquad\qquad\qquad\qquad\qquad\blacksquare\qquad\qquad\qquad\qquad
\qquad\qquad\qquad\qquad\qquad$

Next we have:

$\qquad\frac{\partial}{\partial x_{k}}$L$_{3}(x_{0})=\frac{\partial}{\partial
x_{k}}[\Phi_{P}$g$_{il}\frac{\partial}{\partial x_{j}}(\nabla$log$\Phi
_{P})_{l}](x_{0})$

\qquad\qquad$\ \ \ =\frac{\partial\Phi_{P}}{\partial x_{k}}(x_{0})[$%
g$_{il}\frac{\partial}{\partial x_{j}}(\nabla$log$\Phi_{P})_{l}](x_{0}%
)+\Phi_{P}(x_{0})\frac{\partial}{\partial x_{k}}[$g$_{il}\frac{\partial
}{\partial x_{j}}(\nabla$log$\Phi_{P})_{l}](x_{0})$

$\frac{\partial}{\partial x_{k}}$L$_{3}(x_{0})\ \ \ =\frac{\partial\Phi_{P}%
}{\partial x_{k}}(x_{0})[$g$_{il}\frac{\partial}{\partial x_{j}}(\nabla
$log$\Phi_{P})_{l}](x_{0})\qquad\qquad\qquad\qquad\qquad\qquad\left(
46\right)  $

$+\Phi_{P}(x_{0})\frac{\partial\text{g}_{il}}{\partial x_{k}}(x_{0}%
)[\frac{\partial}{\partial x_{j}}(\nabla$log$\Phi_{P})_{l}](x_{0})+\Phi
_{P}(x_{0})$g$_{il}(x_{0})[\frac{\partial^{2}}{\partial x_{j}\partial x_{k}%
}(\nabla$log$\Phi_{P})_{l}](x_{0})$

\qquad\qquad\qquad\qquad\qquad\qquad\qquad\qquad\qquad\qquad\qquad\qquad
\qquad\qquad\qquad\qquad\qquad\qquad$\blacksquare$

Since $\Phi_{P}(y)=1,$ g$_{il}(y)=\delta_{il}$ and $\frac{\partial
\text{g}_{il}}{\partial x_{k}}(y)=0$ for $i,k,l=q+1,...,n$,

we have for a = 1,...,q:

$\frac{\partial}{\partial x_{k}}$L$_{3}(y)=\frac{\partial\Phi_{P}}{\partial
x_{k}}(y)[\frac{\partial}{\partial x_{j}}(\nabla$log$\Phi_{P})_{i}%
](y)+\frac{\partial\text{g}_{i\text{a}}}{\partial x_{k}}(y)[\frac{\partial
}{\partial x_{j}}(\nabla$log$\Phi_{P})_{\text{a}}](y)$

$+[\frac{\partial^{2}}{\partial x_{j}\partial x_{k}}(\nabla$log$\Phi_{P}%
)_{i}](y)$

We have by (xv) of B$_{1}\qquad\qquad\qquad\ $

\ $\frac{\partial\text{g}_{\text{a}i}}{\partial\text{x}_{k}}(y)=$
$\perp_{\text{a}ki}(y)=-\perp_{\text{a}ik}(y),$ $\frac{\partial}{\partial
x_{j}}(\nabla\log\Phi_{P})_{\text{a}}(y)$

$=-X_{l}(y)\perp_{\text{a}jl}(y)-\frac{\partial X_{j}}{\partial x_{\text{a}}%
}(y)$

We conclude from $\left(  39\right)  $ and the last line above that:

$\frac{\partial}{\partial x_{k}}$L$_{3}(y)=\frac{1}{2}X_{k}(y)\left(
\frac{\partial X_{j}}{\partial x_{i}}+\frac{\partial X_{i}}{\partial x_{j}%
}\right)  (y)+[\perp_{\text{a}ik}\perp_{\text{a}jl}X_{l}](y)+[\perp
_{\text{a}ik}\frac{\partial X_{j}}{\partial x_{\text{a}}}](y)\qquad\left(
47\right)  $

$\qquad\qquad\qquad\qquad+[\frac{\partial^{2}}{\partial x_{j}\partial x_{k}%
}(\nabla$log$\Phi_{P})_{i}](y)\qquad\qquad\qquad\qquad\qquad\qquad\qquad
\qquad\qquad\ \ $

\qquad\qquad\qquad\qquad\qquad\qquad\qquad\qquad\qquad\qquad\qquad\qquad
\qquad\qquad\qquad\qquad\qquad$\blacksquare$

From $\left(  41\right)  ,\left(  43\right)  ,\left(  45\right)  ,\left(
47\right)  ,$ we have:

$\frac{\partial^{3}\Phi_{P}}{\partial x_{k}\partial x_{j}\partial x_{i}%
}(y)=-X_{i}(y)X_{j}(y)X_{k}(y)+\frac{1}{2}X_{i}(y)\left(  \frac{\partial
X_{k}}{\partial x_{j}}+\frac{\partial X_{j}}{\partial x_{k}}\right)
(y)+\frac{1}{2}X_{j}(y)\left(  \frac{\partial X_{i}}{\partial x_{k}}%
+\frac{\partial X_{k}}{\partial x_{i}}\right)  (y)\ \ \ \ \ \ $

\qquad\qquad\qquad$\lbrack\perp_{\text{a}ij}\perp_{\text{a}kl}X_{l}%
+\perp_{\text{a}ij}\frac{\partial X_{k}}{\partial x_{\text{a}}}](y)-\frac
{1}{3}(R_{ijkl}+R_{ikjl})(y)X_{l}(y)\qquad\qquad\qquad\qquad\qquad$

$+\frac{1}{2}X_{k}(y)\left(  \frac{\partial X_{j}}{\partial x_{i}}%
+\frac{\partial X_{i}}{\partial x_{j}}\right)  (y)+[\perp_{\text{a}ik}%
\perp_{\text{a}jl}X_{l}](y)+[\perp_{\text{a}ik}\frac{\partial X_{j}}{\partial
x_{\text{a}}}](y)+[\frac{\partial^{2}}{\partial x_{j}\partial x_{k}}(\nabla
$log$\Phi_{P})_{i}](y)$

We re-order the terms of the above expression as follows:

$\frac{\partial^{3}\Phi_{P}}{\partial x_{k}\partial x_{j}\partial x_{i}%
}(y)=-X_{i}(y)X_{j}(y)X_{k}(y)+\frac{1}{2}X_{i}(y)\left(  \frac{\partial
X_{k}}{\partial x_{j}}+\frac{\partial X_{j}}{\partial x_{k}}\right)
(y)+\frac{1}{2}X_{j}(y)\left(  \frac{\partial X_{i}}{\partial x_{k}}%
+\frac{\partial X_{k}}{\partial x_{i}}\right)  (y)\ \ \ \ \ \ \ \left(
48\right)  $

\qquad\qquad$+\frac{1}{2}X_{k}(y)\left(  \frac{\partial X_{j}}{\partial x_{i}%
}+\frac{\partial X_{i}}{\partial x_{j}}\right)  (y)+[\frac{\partial^{2}%
}{\partial x_{j}\partial x_{k}}(\nabla\log\Phi_{P})_{i}](y)-\frac{1}%
{3}(R_{ijkl}+R_{ikjl})(y)X_{l}(y)$

$\qquad\qquad+[\perp_{\text{a}ij}\perp_{\text{a}kl}X_{l}+\perp_{\text{a}%
ij}\frac{\partial X_{k}}{\partial x_{\text{a}}}](y)+[\perp_{\text{a}ik}%
\perp_{\text{a}jl}X_{l}+\perp_{\text{a}ik}\frac{\partial X_{j}}{\partial
x_{\text{a}}}](y)\qquad$

$\qquad\qquad\qquad\qquad\qquad\qquad\qquad\qquad\qquad\qquad\qquad
\qquad\qquad\qquad\qquad\qquad\qquad\blacksquare\qquad\qquad\qquad\qquad
\qquad\qquad\qquad$\qquad\qquad

In $\left(  48\right)  $ above, we switch the positions of the indices $i$ and
$k$

(the first four terms do not change) and have: \qquad

$\frac{\partial^{3}\Phi_{P}}{\partial x_{i}\partial x_{j}\partial x_{k}}(y)$

$=-X_{i}(y)X_{j}(y)X_{k}(y)+\frac{1}{2}X_{i}(y)\left(  \frac{\partial X_{k}%
}{\partial x_{j}}+\frac{\partial X_{j}}{\partial x_{k}}\right)  (y)+\frac
{1}{2}X_{j}(y)\left(  \frac{\partial X_{i}}{\partial x_{k}}+\frac{\partial
X_{k}}{\partial x_{i}}\right)  (y)\ \ \ \ \ \left(  49\right)  $

$+\frac{1}{2}X_{k}(y)\left(  \frac{\partial X_{j}}{\partial x_{i}}%
+\frac{\partial X_{i}}{\partial x_{j}}\right)  (y)+[\frac{\partial^{2}%
}{\partial x_{j}\partial x_{i}}(\nabla\log\Phi_{P})_{k}](y)-\frac{1}%
{3}(R_{kjil}+R_{kijl})(y)X_{l}(y)$

$\qquad\qquad+[\perp_{\text{a}kj}\perp_{\text{a}il}X_{l}+\perp_{\text{a}%
kj}\frac{\partial X_{i}}{\partial x_{\text{a}}}](y)+[\perp_{\text{a}ki}%
\perp_{\text{a}jl}X_{l}+\perp_{\text{a}ki}\frac{\partial X_{j}}{\partial
x_{\text{a}}}](y)$

\qquad\qquad\qquad\qquad\qquad\qquad\qquad\qquad\qquad\qquad\qquad\qquad
\qquad\qquad\qquad\qquad\qquad$\blacksquare$\qquad\qquad\qquad\qquad
\qquad\qquad\qquad\qquad\qquad

We lastly switch the positions of $j$ and $k$ in $\left(  49\right)  :$

$\frac{\partial^{3}\Phi_{P}}{\partial x_{i}\partial x_{k}\partial x_{j}%
}(y)=-X_{i}(y)X_{j}(y)X_{k}(y)+\frac{1}{2}X_{i}(y)\left(  \frac{\partial
X_{k}}{\partial x_{j}}+\frac{\partial X_{j}}{\partial x_{k}}\right)
(y)\qquad\ \left(  50\right)  $

$\qquad\qquad\qquad+\frac{1}{2}X_{j}(y)\left(  \frac{\partial X_{i}}{\partial
x_{k}}+\frac{\partial X_{k}}{\partial x_{i}}\right)  (y)\ \ \ \ \ \ \ $

$+\frac{1}{2}X_{k}(y)\left(  \frac{\partial X_{j}}{\partial x_{i}}%
+\frac{\partial X_{i}}{\partial x_{j}}\right)  (y)+[\frac{\partial^{2}%
}{\partial x_{k}\partial x_{i}}(\nabla\log\Phi_{P})_{j}](y)-\frac{1}%
{3}(R_{jkil}+R_{jikl})(y)X_{l}(y)$

$+[\perp_{\text{a}jk}\perp_{\text{a}il}X_{l}+\perp_{\text{a}jk}\frac{\partial
X_{i}}{\partial x_{\text{a}}}](y)+[\perp_{\text{a}ji}\perp_{\text{a}kl}%
X_{l}+\perp_{\text{a}ji}\frac{\partial X_{k}}{\partial x_{\text{a}}}](y)$

\qquad\qquad\qquad\qquad\qquad\qquad\qquad\qquad\qquad\qquad\qquad\qquad
\qquad\qquad\qquad\qquad\qquad$\blacksquare$

We notice that the first four terms in each of the three equations are
identical. Only the last two terms are different. Consequently, adding terms
on each side of the equations in $\left(  48\right)  ,\left(  49\right)  $ and
$\left(  50\right)  $ we have:

$\frac{\partial^{3}\Phi_{P}}{\partial x_{k}\partial x_{j}\partial x_{i}%
}(y)+\frac{\partial^{3}\Phi_{P}}{\partial x_{i}\partial x_{j}\partial x_{k}%
}(y)+\frac{\partial^{3}\Phi_{P}}{\partial x_{i}\partial x_{k}\partial x_{j}%
}(y)\qquad\qquad\qquad\qquad\qquad\qquad\qquad\qquad\left(  51\right)  $

$=-3(X_{i}X_{j}X_{k}(y)+\frac{3}{2}X_{i}(y)\left(  \frac{\partial X_{k}%
}{\partial x_{j}}+\frac{\partial X_{j}}{\partial x_{k}}\right)  (y)+\frac
{3}{2}X_{j}(y)\left(  \frac{\partial X_{k}}{\partial x_{i}}+\frac{\partial
X_{i}}{\partial x_{k}}\right)  (y)$

$+\frac{3}{2}X_{k}(y)\left(  \frac{\partial X_{j}}{\partial x_{i}}%
+\frac{\partial X_{i}}{\partial x_{j}}\right)  (y)$

$+[\frac{\partial^{2}}{\partial x_{j}\partial x_{k}}(\nabla\log\Phi_{P}%
)_{i}+\frac{\partial^{2}}{\partial x_{j}\partial x_{i}}(\nabla\log\Phi
_{P})_{k}+\frac{\partial^{2}}{\partial x_{k}\partial x_{i}}(\nabla\log\Phi
_{P})_{j}](y)$

\qquad$\ -\frac{1}{3}\left[  R_{ijkl}+R_{ikjl}+R_{kjil}+R_{kijl}%
+R_{jkil}+R_{jikl}\right]  (y)X_{l}(y)$

\qquad$+[\perp_{\text{a}ij}\perp_{\text{a}kl}X_{l}+\perp_{\text{a}ij}%
\frac{\partial X_{k}}{\partial x_{\text{a}}}](y)+[\perp_{\text{a}ik}%
\perp_{\text{a}jl}X_{l}+\perp_{\text{a}ik}\frac{\partial X_{j}}{\partial
x_{\text{a}}}](y)$

\qquad$+[\perp_{\text{a}kj}\perp_{\text{a}il}X_{l}+\perp_{\text{a}kj}%
\frac{\partial X_{i}}{\partial x_{\text{a}}}](y)+[\perp_{\text{a}ki}%
\perp_{\text{a}jl}X_{l}+\perp_{\text{a}ki}\frac{\partial X_{j}}{\partial
x_{\text{a}}}](y)$

\qquad$+[\perp_{\text{a}jk}\perp_{\text{a}il}X_{l}+\perp_{\text{a}jk}%
\frac{\partial X_{i}}{\partial x_{\text{a}}}](y)+[\perp_{\text{a}ji}%
\perp_{\text{a}kl}X_{l}+\perp_{\text{a}ji}\frac{\partial X_{k}}{\partial
x_{\text{a}}}](y)$

\qquad\qquad\qquad\qquad\qquad\qquad\qquad\qquad\qquad\qquad\qquad\qquad
\qquad\qquad\qquad\qquad\qquad\qquad$\blacksquare$

We observe that:

$R_{ijkl}+R_{ikjl}+R_{kjil}+R_{kijl}+R_{jkil}+R_{jikl}=R_{ijkl}+R_{ikjl}%
-R_{jkil}$

$+R_{jkil}-R_{ikjl}-R_{ijkl}=0$

We next observe that:\qquad

\qquad$+[\perp_{\text{a}ij}\perp_{\text{a}kl}X_{l}+\perp_{\text{a}ij}%
\frac{\partial X_{k}}{\partial x_{\text{a}}}](y)+[\perp_{\text{a}ik}%
\perp_{\text{a}jl}X_{l}+\perp_{\text{a}ik}\frac{\partial X_{j}}{\partial
x_{\text{a}}}](y)$

\qquad$+[\perp_{\text{a}kj}\perp_{\text{a}il}X_{l}+\perp_{\text{a}kj}%
\frac{\partial X_{i}}{\partial x_{\text{a}}}](y)+[\perp_{\text{a}ki}%
\perp_{\text{a}jl}X_{l}+\perp_{\text{a}ki}\frac{\partial X_{j}}{\partial
x_{\text{a}}}](y)$

\qquad$+[\perp_{\text{a}jk}\perp_{\text{a}il}X_{l}+\perp_{\text{a}jk}%
\frac{\partial X_{i}}{\partial x_{\text{a}}}](y)+[\perp_{\text{a}ji}%
\perp_{\text{a}kl}X_{l}+\perp_{\text{a}ji}\frac{\partial X_{k}}{\partial
x_{\text{a}}}](y)$

$=+[\perp_{\text{a}ij}\perp_{\text{a}kl}X_{l}+\perp_{\text{a}ij}\frac{\partial
X_{k}}{\partial x_{\text{a}}}](y)+[\perp_{\text{a}ik}\perp_{\text{a}jl}%
X_{l}+\perp_{\text{a}ik}\frac{\partial X_{j}}{\partial x_{\text{a}}}](y)$

$\ +[-\perp_{\text{a}jk}\perp_{\text{a}il}X_{l}-\perp_{\text{a}jk}%
\frac{\partial X_{i}}{\partial x_{\text{a}}}](y)+[-\perp_{\text{a}ik}%
\perp_{\text{a}jl}X_{l}-\perp_{\text{a}ik}\frac{\partial X_{j}}{\partial
x_{\text{a}}}](y)$

$\ +[\perp_{\text{a}jk}\perp_{\text{a}il}X_{l}+\perp_{\text{a}jk}%
\frac{\partial X_{i}}{\partial x_{\text{a}}}](y)+[-\perp_{\text{a}ij}%
\perp_{\text{a}kl}X_{l}-\perp_{\text{a}ij}\frac{\partial X_{k}}{\partial
x_{\text{a}}}](y)$

$=0\qquad\qquad\qquad$

We then have the beautiful equation: For all $y\in U\subset P\subset M_{0},$\ \ \ \ 

$\frac{\partial^{3}\Phi_{P}}{\partial x_{k}\partial x_{j}\partial x_{i}%
}(y)+\frac{\partial^{3}\Phi_{P}}{\partial x_{i}\partial x_{j}\partial x_{k}%
}(y)+\frac{\partial^{3}\Phi_{P}}{\partial x_{i}\partial x_{k}\partial x_{j}%
}(y)\qquad\qquad\qquad\qquad\qquad\qquad\qquad\qquad\qquad\left(  52\right)
\qquad\qquad\qquad\qquad\qquad\qquad\qquad$

$=-3(X_{i}X_{j}X_{k})(y)+\frac{3}{2}X_{i}(y)\left(  \frac{\partial X_{k}%
}{\partial x_{j}}+\frac{\partial X_{j}}{\partial x_{k}}\right)  (y)+\frac
{3}{2}X_{j}(y)\left(  \frac{\partial X_{k}}{\partial x_{i}}+\frac{\partial
X_{i}}{\partial x_{k}}\right)  (y)$

$\qquad+\frac{3}{2}X_{k}(y)\left(  \frac{\partial X_{j}}{\partial x_{i}}%
+\frac{\partial X_{i}}{\partial x_{j}}\right)  (y)$

$\qquad\ +[\frac{\partial^{2}}{\partial x_{j}\partial x_{k}}(\nabla\log
\Phi_{P})_{i}+\frac{\partial^{2}}{\partial x_{j}\partial x_{i}}(\nabla\log
\Phi_{P})_{k}+\frac{\partial^{2}}{\partial x_{k}\partial x_{i}}(\nabla\log
\Phi_{P})_{j}](y)$

\qquad\qquad\qquad\qquad\qquad\qquad\qquad\qquad\qquad\qquad\qquad\qquad
\qquad\qquad\qquad\qquad\qquad\qquad$\blacksquare$\qquad$\ \ \ \ \ \ \qquad
$\qquad\qquad\qquad$\ \ \ \ \ \ $

Since $\Phi:M_{0}\longrightarrow R$ is a smooth function, we then have from
Calculus that for all $x_{0}\in M_{0}:$

\qquad\qquad$\frac{\partial^{3}\Phi_{P}}{\partial x_{k}\partial x_{j}\partial
x_{i}}(x_{0})=\frac{\partial^{3}\Phi_{P}}{\partial x_{i}\partial x_{j}\partial
x_{k}}(x_{0})=\frac{\partial^{3}\Phi_{P}}{\partial x_{i}\partial x_{k}\partial
x_{j}}(x_{0})$ $\qquad$

We then can re-write $\left(  52\right)  $ as:

$3\frac{\partial^{3}\Phi_{P}}{\partial x_{i}\partial x_{j}\partial x_{k}%
}(y)=-3(X_{i}X_{j}X_{k})(y)$

$+\frac{3}{2}X_{i}(y)\left(  \frac{\partial X_{k}}{\partial x_{j}}%
+\frac{\partial X_{j}}{\partial x_{k}}\right)  (y)+\frac{3}{2}X_{j}(y)\left(
\frac{\partial X_{k}}{\partial x_{i}}+\frac{\partial X_{i}}{\partial x_{k}%
}\right)  (y)+\frac{3}{2}X_{k}(y)\left(  \frac{\partial X_{j}}{\partial x_{i}%
}+\frac{\partial X_{i}}{\partial x_{j}}\right)  (y)$

$+[\frac{\partial^{2}}{\partial x_{j}\partial x_{k}}(\nabla\log\Phi_{P}%
)_{i}+\frac{\partial^{2}}{\partial x_{j}\partial x_{i}}(\nabla\log\Phi
_{P})_{k}+\frac{\partial^{2}}{\partial x_{k}\partial x_{i}}(\nabla\log\Phi
_{P})_{j}](y)$

$\qquad\qquad\qquad\qquad\qquad\qquad\qquad\qquad$

By (viii) of \textbf{Table B}$_{1},$ we have for $y\in U\subset P:$

\qquad$\lbrack\frac{\partial^{2}}{\partial x_{j}\partial x_{k}}(\nabla
$log$\Phi_{P})_{i}+\frac{\partial^{2}}{\partial x_{j}\partial x_{i}}%
(\nabla\log\Phi_{P})_{k}+\frac{\partial^{2}}{\partial x_{k}\partial x_{i}%
}(\nabla\log\Phi_{P})_{j}](y)\qquad\qquad\qquad\qquad\left(  53\right)  $

$\qquad=-\left(  \frac{\partial^{2}X_{i}}{\partial x_{j}\partial x_{k}}%
+\frac{\partial^{2}X_{j}}{\partial x_{i}\partial x_{k}}+\frac{\partial
^{2}X_{k}}{\partial x_{i}\partial x_{j}}\right)  (y)$

\qquad\qquad\qquad\qquad\qquad\qquad\qquad\qquad\qquad\qquad\qquad\qquad
\qquad\qquad\qquad\qquad$\blacksquare$

Therefore by $\left(  53\right)  ,$ we can re-write $\left(  52\right)  $
as:\qquad

$\qquad\frac{\partial^{3}\Phi_{P}}{\partial x_{i}\partial x_{j}\partial x_{k}%
}(y)$

$\qquad=-X_{i}(y)X_{j}(y)X_{k}(y)+\frac{1}{2}X_{i}(y)\left(  \frac{\partial
X_{k}}{\partial x_{j}}+\frac{\partial X_{j}}{\partial x_{k}}\right)
(y)+\frac{1}{2}X_{j}(y)\left(  \frac{\partial X_{k}}{\partial x_{i}}%
+\frac{\partial X_{i}}{\partial x_{k}}\right)  (y)\qquad\left(  54\right)  $

$\qquad+\frac{1}{2}X_{k}(y)\left(  \frac{\partial X_{j}}{\partial x_{i}}%
+\frac{\partial X_{i}}{\partial x_{j}}\right)  (y)-\frac{1}{3}\left(
\frac{\partial^{2}X_{i}}{\partial x_{j}\partial x_{k}}+\frac{\partial^{2}%
X_{j}}{\partial x_{i}\partial x_{k}}+\frac{\partial^{2}X_{k}}{\partial
x_{i}\partial x_{j}}\right)  (y)$

\qquad\qquad\qquad\qquad\qquad\qquad\qquad\qquad\qquad\qquad\qquad\qquad
\qquad\qquad\qquad\qquad\qquad\qquad\qquad$\blacksquare$

(v$^{\ast}$)\ In particular, taking $k=i,$ in $\left(  54\right)  $ above, we have:

$\qquad\frac{\partial^{3}\Phi_{P}}{\partial x_{i}^{2}\partial x_{j}}%
(y)=-X_{i}^{2}(y)X_{j}(y)+\frac{1}{2}X_{i}(y)\left(  \frac{\partial X_{i}%
}{\partial x_{j}}+\frac{\partial X_{j}}{\partial x_{i}}\right)  (y)\qquad$

$\qquad+\frac{1}{2}X_{j}(y)\left(  \frac{\partial X_{i}}{\partial x_{i}}%
+\frac{\partial X_{i}}{\partial x_{i}}\right)  \ (y)+\frac{1}{2}%
X_{i}(y)\left(  \frac{\partial X_{j}}{\partial x_{i}}+\frac{\partial X_{i}%
}{\partial x_{j}}\right)  (y)$

\qquad$-\frac{1}{3}\left(  \frac{\partial^{2}X_{i}}{\partial x_{j}\partial
x_{i}}+\frac{\partial^{2}X_{j}}{\partial x_{i}\partial x_{i}}+\frac
{\partial^{2}X_{i}}{\partial x_{i}\partial x_{j}}\right)  (y)$$\qquad$

$=-X_{i}^{2}(y)X_{j}(y)+X_{i}(y)\left(  \frac{\partial X_{i}}{\partial x_{j}%
}+\frac{\partial X_{j}}{\partial x_{i}}\right)  (y)+\frac{1}{2}X_{j}(y)\left(
\frac{\partial X_{i}}{\partial x_{i}}+\frac{\partial X_{i}}{\partial x_{i}%
}\right)  \ (y)-\frac{1}{3}\left(  \frac{\partial^{2}X_{j}}{\partial x_{i}%
^{2}}+2\frac{\partial^{2}X_{i}}{\partial x_{i}\partial x_{j}}\right)
(y)\qquad$

$\frac{\partial^{3}\Phi_{P}}{\partial x_{i}^{2}\partial x_{j}}(y)=-X_{i}%
^{2}(y)X_{j}(y)+X_{i}(y)\left(  \frac{\partial X_{j}}{\partial x_{i}}%
+\frac{\partial X_{i}}{\partial x_{j}}\right)  (y)+X_{j}(y)\frac{\partial
X_{i}}{\partial x_{i}}\ (y)-\frac{1}{3}\left(  \frac{\partial^{2}X_{j}%
}{\partial x_{i}^{2}}+2\frac{\partial^{2}X_{i}}{\partial x_{i}\partial x_{j}%
}\right)  (y)\qquad\left(  55\right)  $

We recall that the expression defining the divergence of a vector field X on
the Riemannian manifold M was given in $\left(  B_{22}\right)  :$

\qquad$\operatorname{div}X(y)$ $=$ $\underset{i=1}{\overset{n}{\sum}}%
\frac{\partial X_{i}}{\partial x_{i}}(y)-\underset{i=q+1}{\overset{n}{\sum}%
}<H,i>(y_{0})X_{i}(y)$

Therefore,$\qquad\qquad\qquad\qquad\qquad$

$\qquad\underset{i=1}{\overset{n}{\sum}}\frac{\partial X_{i}}{\partial x_{i}%
}(y)=\operatorname{div}X(y)+\underset{i=q+1}{\overset{n}{\sum}}<H,i>(y_{0}%
)X_{i}(y)$

\qquad$\underset{i=q+1}{\overset{n}{\sum}}\frac{\partial X_{i}}{\partial
x_{i}}(y)=\operatorname{div}X(y)+\underset{i=q+1}{\overset{n}{\sum}%
}<H,i>(y)X_{i}(y)-\underset{\text{a}=1}{\overset{q}{\sum}}\frac{\partial
X_{\text{a}}}{\partial x_{\text{a}}}(y)$

\qquad$\underset{i,j=q+1}{\overset{n}{\sum}}\frac{\partial^{3}\Phi_{P}%
}{\partial x_{i}^{2}\partial x_{j}}(y)=-\underset{i,j=q+1}{\overset{n}{\sum}%
}X_{j}(y)X_{i}^{2}(y)+\underset{i,j=q+1}{\overset{n}{\sum}}X_{j}%
(y)\frac{\partial X_{i}}{\partial x_{i}}%
\ (y)+\underset{i,j=q+1}{\overset{n}{\sum}}X_{i}(y)\left(  \frac{\partial
X_{j}}{\partial x_{i}}+\frac{\partial X_{i}}{\partial x_{j}}\right)
(y)\qquad$

\qquad$-\frac{1}{3}\underset{i,j=q+1}{\overset{n}{\sum}}\left(  \frac
{\partial^{2}X_{j}}{\partial x_{i}^{2}}+2\frac{\partial^{2}X_{i}}{\partial
x_{i}\partial x_{j}}\right)  (y)=-\underset{j=q+1}{\overset{n}{\sum}}%
X_{j}\left(  \left\Vert \text{X}\right\Vert _{M}^{2}-\left\Vert \text{X}%
\right\Vert _{P}^{2}\right)  (y)\qquad\qquad\qquad\qquad\qquad\qquad
\qquad\qquad\qquad\qquad\qquad$

$\qquad\qquad\qquad+\underset{j=q+1}{\overset{n}{\sum}}X_{j}(y)\left(
\operatorname{div}X(y)+\underset{i=q+1}{\overset{n}{\sum}}<H,i>(y_{0}%
)X_{i}(y)-\underset{\text{a=1}}{\overset{\text{q}}{\sum}}\frac{\partial
X_{\text{a}}}{\partial x_{\text{a}}}(y)\right)  $

$\qquad\qquad\qquad+\underset{i,j=q+1}{\overset{n}{\sum}}X_{i}(y)\left(
\frac{\partial X_{j}}{\partial x_{i}}+\frac{\partial X_{i}}{\partial x_{j}%
}\right)  (y)-\frac{1}{3}\underset{i,j=q+1}{\overset{n}{\sum}}\left(
\frac{\partial^{2}X_{j}}{\partial x_{i}^{2}}+2\frac{\partial^{2}X_{i}%
}{\partial x_{i}\partial x_{j}}\right)  (y)$

\qquad\qquad\qquad\qquad\qquad\qquad\qquad\qquad\qquad\qquad\qquad\qquad
\qquad\qquad\qquad\qquad\qquad$\blacksquare$\qquad
$\underset{i,j=q+1}{\overset{n}{\sum}}\frac{\partial^{3}\Phi_{P}}{\partial
x_{i}^{2}\partial x_{j}}(y)=\underset{j=q+1}{\overset{n}{\sum}}X_{j}\left(
\operatorname{div}X-\left\Vert \text{X}\right\Vert _{M}^{2}-\operatorname{div}%
X_{P}+\left\Vert \text{X}\right\Vert _{P}^{2}%
+\underset{i=q+1}{\overset{n}{\sum}}<H,i>(y_{0})X_{i}\right)  (y_{0}%
)\qquad\qquad\left(  56\right)  \qquad\qquad\qquad\qquad\qquad\qquad
\qquad\qquad\qquad\qquad\left(  56\right)  $

$\qquad\qquad\qquad+\underset{i,j=q+1}{\overset{n}{\sum}}X_{i}(y)\left(
\frac{\partial X_{j}}{\partial x_{i}}+\frac{\partial X_{i}}{\partial x_{j}%
}\right)  (y)-\frac{1}{3}\underset{i,j=q+1}{\overset{n}{\sum}}\left(
\frac{\partial^{2}X_{j}}{\partial x_{i}^{2}}+2\frac{\partial^{2}X_{i}%
}{\partial x_{i}\partial x_{j}}\right)  (y)$

\qquad\qquad\qquad\qquad\qquad\qquad\qquad\qquad\qquad\qquad\qquad\qquad
\qquad\qquad\qquad\qquad\qquad$\blacksquare$

(v)$^{\ast\ast}$ We obtain $\underset{i,j=q+1}{\overset{n}{\sum}}%
\frac{\partial^{3}\Phi_{P}}{\partial x_{i}\partial x_{j}^{2}}(y)$ by
inter-changing the positions of $i$ and $j$ in $\left(  56\right)  :$

$\frac{\partial^{3}\Phi_{P}}{\partial x_{i}\partial x_{j}^{2}}(y)=-X_{i}%
(y)X_{j}^{2}(y)+X_{i}(y)\frac{\partial X_{j}}{\partial x_{j}}(y)+X_{j}%
(y)\left(  \frac{\partial X_{i}}{\partial x_{j}}+\frac{\partial X_{j}%
}{\partial x_{i}}\right)  \ (y)-\frac{1}{3}\left(  \frac{\partial^{2}X_{i}%
}{\partial x_{j}^{2}}+2\frac{\partial^{2}X_{j}}{\partial x_{i}\partial x_{j}%
}\right)  (y)\qquad\left(  57\right)  \qquad$

We have finally,

$\underset{i,j=q+1}{\overset{n}{\sum}}\frac{\partial^{3}\Phi_{P}}{\partial
x_{i}\partial x_{j}^{2}}(y)$

$\qquad=\underset{i=q+1}{\overset{n}{\sum}}X_{i}(y)\left(  \operatorname{div}%
X_{M}-\left\Vert \text{X}\right\Vert _{M}^{2}-\operatorname{div}%
X_{P}+\left\Vert \text{X}\right\Vert _{P}^{2}%
-\underset{j=q+1}{\overset{n}{\sum}}<H,j>(y_{0})X_{j}\right)  (y)\qquad\left(
58\right)  \qquad\qquad$

\qquad$+\underset{i,j=q+1}{\overset{n}{\sum}}X_{j}(y)\left(  \frac{\partial
X_{i}}{\partial x_{j}}+\frac{\partial X_{j}}{\partial x_{i}}\right)
(y)-\frac{1}{3}\underset{i,j=q+1}{\overset{n}{\sum}}\left(  \frac{\partial
^{2}X_{i}}{\partial x_{j}^{2}}+2\frac{\partial^{2}X_{j}}{\partial
x_{i}\partial x_{j}}\right)  (y)$

\qquad\qquad\qquad\qquad\qquad\qquad\qquad\qquad\qquad\qquad\qquad\qquad
\qquad\qquad\qquad\qquad\qquad\qquad\qquad$\blacksquare$

We shall need the expression for $[\frac{\partial^{2}}{\partial x_{j}\partial
x_{i}}(\nabla\log\Phi_{P})_{k}](y):$

(v)$^{\ast\ast\ast}$ We compare $\left(  49\right)  $ and $\left(  54\right)
$ here in \textbf{Table B}$_{4}$ and see that:$\qquad\qquad\qquad\qquad
\qquad\qquad\qquad$

$\qquad\frac{\partial^{3}\Phi_{P}}{\partial x_{i}\partial x_{j}\partial x_{k}%
}(y)=-X_{i}(y)X_{j}(y)X_{k}(y)+\frac{1}{2}X_{i}(y)\left(  \frac{\partial
X_{k}}{\partial x_{j}}+\frac{\partial X_{j}}{\partial x_{k}}\right)  (y)$

$+\frac{1}{2}X_{j}(y)\left(  \frac{\partial X_{i}}{\partial x_{k}}%
+\frac{\partial X_{k}}{\partial x_{i}}\right)  (y)+\frac{1}{2}X_{k}(y)\left(
\frac{\partial X_{j}}{\partial x_{i}}+\frac{\partial X_{i}}{\partial x_{j}%
}\right)  (y)+[\frac{\partial^{2}}{\partial x_{j}\partial x_{i}}(\nabla
\log\Phi_{P})_{k}](y)$

$-\frac{1}{3}(R_{kjil}+R_{kijl})(y)X_{l}(y)$

$+[\perp_{\text{a}kj}\perp_{\text{a}il}X_{l}+\perp_{\text{a}kj}\frac{\partial
X_{i}}{\partial x_{\text{a}}}](y)+[\perp_{\text{a}ki}\perp_{\text{a}jl}%
X_{l}](y)+[\perp_{\text{a}ki}\frac{\partial X_{j}}{\partial x_{\text{a}}}](y)$

$\qquad$

$=-X_{i}(y)X_{j}(y)X_{k}(y)+\frac{1}{2}X_{i}(y)\left(  \frac{\partial X_{k}%
}{\partial x_{j}}+\frac{\partial X_{j}}{\partial x_{k}}\right)  (y)+\frac
{1}{2}X_{j}(y)\left(  \frac{\partial X_{k}}{\partial x_{i}}+\frac{\partial
X_{i}}{\partial x_{k}}\right)  (y)$

$\qquad+\frac{1}{2}X_{k}(y)\left(  \frac{\partial X_{j}}{\partial x_{i}}%
+\frac{\partial X_{i}}{\partial x_{j}}\right)  (y)-\frac{1}{3}\left(
\frac{\partial^{2}X_{i}}{\partial x_{j}\partial x_{k}}+\frac{\partial^{2}%
X_{j}}{\partial x_{i}\partial x_{k}}+\frac{\partial^{2}X_{k}}{\partial
x_{i}\partial x_{j}}\right)  (y)$

Therefore, at any point $y\in P,$ we have:

$\qquad\lbrack\frac{\partial^{2}}{\partial x_{i}\partial x_{j}}(\nabla\log
\Phi_{P})_{k}](y)$

$\qquad=-\frac{1}{3}\left(  \frac{\partial^{2}X_{i}}{\partial x_{j}\partial
x_{k}}+\frac{\partial^{2}X_{j}}{\partial x_{i}\partial x_{k}}+\frac
{\partial^{2}X_{k}}{\partial x_{i}\partial x_{j}}\right)  (y)+\frac{1}%
{3}(R_{kjil}+R_{kijl})(y)X_{l}(y)$

$\qquad-[\perp_{\text{a}kj}\perp_{\text{a}il}X_{l}+\perp_{\text{a}kj}%
\frac{\partial X_{i}}{\partial x_{\text{a}}}](y)-[\perp_{\text{a}ki}%
\perp_{\text{a}jl}X_{l}+\perp_{\text{a}ki}\frac{\partial X_{j}}{\partial
x_{\text{a}}}](y)$

\qquad$=-\frac{1}{3}\left(  \frac{\partial^{2}X_{k}}{\partial x_{i}\partial
x_{j}}+\frac{\partial^{2}X_{j}}{\partial x_{i}\partial x_{k}}+\frac
{\partial^{2}X_{i}}{\partial x_{j}\partial x_{k}}\right)  (y)-\frac{1}%
{3}(R_{ikjl}+R_{jkil})(y)X_{l}(y)$

$\qquad+[\perp_{\text{a}jk}\perp_{\text{a}il}X_{l}+\perp_{\text{a}jk}%
\frac{\partial X_{i}}{\partial x_{\text{a}}}](y)+[\perp_{\text{a}ik}%
\perp_{\text{a}jl}X_{l}+\perp_{\text{a}ik}\frac{\partial X_{j}}{\partial
x_{\text{a}}}](y)$

\qquad\qquad\qquad\qquad\qquad\qquad\qquad\qquad\qquad\qquad\qquad\qquad
\qquad\qquad\qquad\qquad\qquad$\blacksquare$

We have thus proved the formula in (v)$^{\ast\ast\ast}$ of \textbf{Table
B}$_{4}$ in \textbf{Appendix B}.

In particular, we have at the centre of Fermi coordinates $y_{0}\in P:$

$\qquad\lbrack\frac{\partial^{2}}{\partial x_{i}\partial x_{j}}(\nabla\log
\Phi_{P})_{k}](y_{0})$

$=-\frac{1}{3}\left(  \frac{\partial^{2}X_{k}}{\partial x_{i}\partial x_{j}%
}+\frac{\partial^{2}X_{j}}{\partial x_{i}\partial x_{k}}+\frac{\partial
^{2}X_{i}}{\partial x_{j}\partial x_{k}}\right)  (y_{0})-\frac{1}%
{3}\underset{l=q+1}{\overset{n}{\sum}}(R_{ikjl}+R_{jkil})(y_{0})X_{l}%
(y_{0})\qquad\left(  59\right)  $

$\qquad+[\perp_{\text{a}ik}\frac{\partial X_{j}}{\partial x_{\text{a}}}%
+\perp_{\text{a}jk}\frac{\partial X_{i}}{\partial x_{\text{a}}}](y_{0}%
)+\underset{l=q+1}{\overset{n}{\sum}}[\perp_{\text{a}ik}\perp_{\text{a}%
jl}X_{l}+\perp_{\text{a}jk}\perp_{\text{a}il}X_{l}](y_{0})$

\qquad\qquad\qquad\qquad\qquad\qquad\qquad\qquad\qquad\qquad\qquad\qquad
\qquad\qquad\qquad\qquad\qquad$\blacksquare$

In particular,

$[\frac{\partial^{2}}{\partial x_{i}\partial x_{j}}(\nabla\log\Phi_{P}%
)_{j}](y_{0})$

$=-\frac{1}{3}\left(  \frac{\partial^{2}X_{j}}{\partial x_{i}\partial x_{j}%
}+\frac{\partial^{2}X_{j}}{\partial x_{i}\partial x_{j}}+\frac{\partial
^{2}X_{i}}{\partial x_{j}^{2}}\right)  (y_{0})-\frac{1}{3}%
\underset{l=q+1}{\overset{n}{\sum}}(R_{ijjl}+R_{jjil})(y_{0})X_{l}%
(y_{0})\qquad$

$+[\perp_{\text{a}ij}\frac{\partial X_{j}}{\partial x_{\text{a}}}%
+\perp_{\text{a}jj}\frac{\partial X_{i}}{\partial x_{\text{a}}}](y_{0}%
)+\underset{l=q+1}{\overset{n}{\sum}}[\perp_{\text{a}ij}\perp_{\text{a}%
jl}X_{l}+\perp_{\text{a}jj}\perp_{\text{a}il}X_{l}](y_{0})$

$=-\frac{1}{3}\left(  2\frac{\partial^{2}X_{j}}{\partial x_{i}\partial x_{j}%
}+\frac{\partial^{2}X_{i}}{\partial x_{j}^{2}}\right)  (y_{0})-\frac{1}%
{3}\underset{l=q+1}{\overset{n}{\sum}}R_{ijjl}(y_{0})X_{l}(y_{0})$

$+[\perp_{\text{a}ij}\frac{\partial X_{j}}{\partial x_{\text{a}}}%
](y_{0})+\underset{l=q+1}{\overset{n}{\sum}}[\perp_{\text{a}ij}\perp
_{\text{a}jl}X_{l}](y_{0})$

$[\frac{\partial^{2}}{\partial x_{i}^{2}}(\nabla\log\Phi_{P})_{k}%
](y_{0})=-\frac{1}{3}\left(  \frac{\partial^{2}X_{k}}{\partial x_{i}^{2}%
}+\frac{\partial^{2}X_{i}}{\partial x_{i}\partial x_{k}}+\frac{\partial
^{2}X_{i}}{\partial x_{i}\partial x_{k}}\right)  (y_{0})$

$-\frac{1}{3}\underset{l=q+1}{\overset{n}{\sum}}(R_{ikil}+R_{ikil}%
)(y_{0})X_{l}(y_{0})\qquad$

$+[\perp_{\text{a}ik}\frac{\partial X_{i}}{\partial x_{\text{a}}}%
+\perp_{\text{a}ik}\frac{\partial X_{i}}{\partial x_{\text{a}}}](y_{0}%
)+\underset{l=q+1}{\overset{n}{\sum}}[\perp_{\text{a}ik}\perp_{\text{a}%
il}X_{l}+\perp_{\text{a}ik}\perp_{\text{a}il}X_{l}](y_{0})$

$[\frac{\partial^{2}}{\partial x_{i}^{2}}(\nabla\log\Phi_{P})_{k}%
](y_{0})=-\frac{1}{3}\left(  \frac{\partial^{2}X_{k}}{\partial x_{i}^{2}%
}+2\frac{\partial^{2}X_{i}}{\partial x_{i}\partial x_{k}}\right)
(y_{0})-\frac{2}{3}\underset{l=q+1}{\overset{n}{\sum}}R_{ikil}(y_{0}%
)X_{l}(y_{0})\qquad$

\qquad\qquad\qquad\qquad\qquad$\qquad+[2\perp_{\text{a}ik}\frac{\partial
X_{i}}{\partial x_{\text{a}}}](y_{0})+\underset{l=q+1}{\overset{n}{\sum}%
}[2\perp_{\text{a}ik}\perp_{\text{a}il}X_{l}](y_{0})$

$[\frac{\partial^{2}}{\partial x_{i}^{2}}(\nabla\log\Phi_{P})_{j}%
](y_{0})=-\frac{1}{3}\left(  \frac{\partial^{2}X_{j}}{\partial x_{i}^{2}%
}+2\frac{\partial^{2}X_{i}}{\partial x_{i}\partial x_{j}}\right)
(y_{0})-\frac{2}{3}\underset{l=q+1}{\overset{n}{\sum}}R_{ijil}(y_{0}%
)X_{l}(y_{0})\qquad$

\qquad\qquad\qquad\qquad\qquad$\qquad+[2\perp_{\text{a}ij}\frac{\partial
X_{i}}{\partial x_{\text{a}}}](y_{0})+\underset{l=q+1}{\overset{n}{\sum}%
}[2\perp_{\text{a}ij}\perp_{\text{a}il}X_{l}](y_{0})$

(vi) We next express $\frac{\partial^{4}\Phi_{P}}{\partial x_{i}^{2}\partial
x_{j}^{2}}(y)$ in terms of the vector field X and geometric invariants:

From $(32)$ we have:

$\frac{\partial\Phi_{P}}{\partial x_{j}}(x_{0})=\Phi_{P}(x_{0})\frac{\partial
}{\partial x_{j}}\log\Phi_{P}(x_{0})=\Phi_{P}(x_{0})$g$_{jk}(x_{0})(\nabla
$log$\Phi_{P})_{k}(x_{0})$

Conequently,

$\frac{\partial^{2}\Phi_{P}}{\partial x_{i}\partial x_{j}}(x_{0}%
)=\frac{\partial}{\partial x_{i}}[\Phi_{P}$g$_{jk}(\nabla$log$\Phi_{P}%
)_{k}](x_{0})$

$\frac{\partial^{2}\Phi_{P}}{\partial x_{i}\partial x_{j}}(x_{0}%
)=\frac{\partial}{\partial x_{i}}[\Phi_{P}$g$_{jk}](x_{0})(\nabla$log$\Phi
_{P})_{k}](x_{0})+[\Phi_{P}$g$_{jk}](x_{0})\frac{\partial}{\partial x_{i}%
}[(\nabla$log$\Phi_{P})_{k}](x_{0})$

$\qquad\qquad=\frac{\partial\Phi_{P}}{\partial x_{i}}(x_{0})[$g$_{jk}(\nabla
$log$\Phi_{P})_{k}](x_{0})\qquad\qquad\qquad\qquad\qquad\qquad$

$\qquad\qquad+\Phi_{P}(x_{0})[\frac{\partial\text{g}_{jk}}{\partial x_{i}%
}(\nabla$log$\Phi_{P})_{k}+$ g$_{jk}\frac{\partial}{\partial x_{i}}(\nabla
$log$\Phi_{P})_{k}](x_{0})=L_{1}+L_{2}+L_{3}$\ \ 

where,

$\qquad\qquad L_{1}=\frac{\partial\Phi_{P}}{\partial x_{i}}(x_{0})[$%
g$_{jk}(\nabla$log$\Phi_{P})_{k}](x_{0})\qquad\qquad\qquad\qquad\qquad
\qquad\left(  60\right)  $

$\qquad\qquad L_{2}=$\ $\Phi_{P}(x_{0})[\frac{\partial\text{g}_{jk}}{\partial
x_{i}}(\nabla$log$\Phi_{P})_{k}](x_{0})\qquad\qquad\qquad\qquad\qquad
\qquad\left(  61\right)  $

$\qquad\qquad L_{3}=\Phi_{P}(x_{0})[$ g$_{jk}\frac{\partial}{\partial x_{i}%
}(\nabla$log$\Phi_{P})_{k}](x_{0})\qquad\qquad\ \ \qquad\qquad\qquad\left(
62\right)  $\ \ 

Then,

$\frac{\partial^{4}\Phi_{P}}{\partial x_{i}\partial x_{j}\partial
x_{i}\partial x_{j}}(y_{0})=\frac{\partial^{2}L_{1}}{\partial x_{i}\partial
x_{j}}(y_{0})+\frac{\partial^{2}L_{2}}{\partial x_{i}\partial x_{j}}%
(y_{0})+\frac{\partial^{2}L_{3}}{\partial x_{i}\partial x_{j}}(y_{0}%
)\qquad\qquad\qquad\ \ \ \ \ \left(  63\right)  $

We compute each of the terms of the expression on the RHS of $\left(
63\right)  :$

$L_{1}=[\frac{\partial\Phi_{P}}{\partial x_{i}}(x_{0})$g$_{jk}(\nabla$%
log$\Phi_{P})_{k}](x_{0})$

$\frac{\partial L_{1}}{\partial x_{j}}(x_{0})=\frac{\partial}{\partial x_{j}%
}[\frac{\partial\Phi_{P}}{\partial x_{i}}$g$_{jk}](x_{0})(\nabla$log$\Phi
_{P})_{k}(x_{0})+[\frac{\partial\Phi_{P}}{\partial x_{i}}$g$_{jk}](x_{0}%
)\frac{\partial}{\partial x_{j}}(\nabla$log$\Phi_{P})_{k}](x_{0})$

\qquad$\ \ \ \ =[\frac{\partial^{2}\Phi_{P}}{\partial x_{j}\partial x_{i}}%
$g$_{jk}+\frac{\partial\Phi_{P}}{\partial x_{i}}\frac{\partial\text{g}_{jk}%
}{\partial x_{j}}](x_{0})(\nabla$log$\Phi_{P})_{k}(x_{0})+[\frac{\partial
\Phi_{P}}{\partial x_{i}}$g$_{jk}\frac{\partial}{\partial x_{j}}(\nabla
$log$\Phi_{P})_{k}](x_{0})$

$=[\frac{\partial^{2}\Phi_{P}}{\partial x_{j}\partial x_{i}}[$g$_{jk}(\nabla
$log$\Phi_{P})_{k}](x_{0})+[\frac{\partial\Phi_{P}}{\partial x_{i}}%
\frac{\partial\text{g}_{jk}}{\partial x_{j}}(\nabla$log$\Phi_{P})_{k}%
](x_{0})+[\frac{\partial\Phi_{P}}{\partial x_{i}}$g$_{jk}](x_{0}%
)\frac{\partial}{\partial x_{j}}(\nabla$log$\Phi_{P})_{k}](x_{0})$

$=[\frac{\partial^{2}\Phi_{P}}{\partial x_{j}\partial x_{i}}$g$_{jk}(\nabla
$log$\Phi_{P})_{k}](x_{0})+\frac{\partial\Phi_{P}}{\partial x_{i}}%
[\frac{\partial\text{g}_{jk}}{\partial x_{j}}(\nabla$log$\Phi_{P})_{k}+$
g$_{jk}\frac{\partial}{\partial x_{j}}(\nabla$log$\Phi_{P})_{k}](x_{0})$

Then,

$\frac{\partial^{2}L_{1}}{\partial x_{i}\partial x_{j}}(y_{0})=\frac
{\partial^{3}\Phi_{P}}{\partial x_{j}\partial x_{i}^{2}}(y_{0})[$%
g$_{jk}(\nabla$log$\Phi_{P})_{k}](y_{0})+\frac{\partial^{2}\Phi_{P}}{\partial
x_{j}\partial x_{i}}(y_{0})[\frac{\partial\text{g}_{jk}}{\partial x_{i}%
}(\nabla\log\Phi_{P})_{k}$

$+g_{jk}\frac{\partial}{\partial x_{i}}(\nabla\log\Phi_{P})_{k}](y_{0}%
)\qquad\left(  64\right)  $

$+\frac{\partial^{2}\Phi_{P}}{\partial x_{i}^{2}}(y_{0})[\frac{\partial
g_{jk}}{\partial x_{j}}(\nabla$log$\Phi_{P})_{k}+$ g$_{jk}\frac{\partial
}{\partial x_{j}}(\nabla$log$\Phi_{P})_{k}](y_{0})$

$+\frac{\partial\Phi_{P}}{\partial x_{i}}(y_{0})[\frac{\partial^{2}g_{jk}%
}{\partial x_{i}\partial x_{j}}(\nabla$log$\Phi_{P})_{k}+\frac{\partial
g_{jk}}{\partial x_{j}}\frac{\partial}{\partial x_{i}}(\nabla$log$\Phi
_{P})_{k}$ $](y_{0})$

$+\frac{\partial\Phi_{P}}{\partial x_{i}}(y_{0})[\frac{\partial g_{jk}%
}{\partial x_{i}}\frac{\partial}{\partial x_{j}}(\nabla\log\Phi_{P}%
)_{k}+g_{jk}\frac{\partial^{2}}{\partial x_{i}\partial x_{j}}(\nabla\log
\Phi_{P})_{k}](y_{0})$

$\qquad\qquad\qquad\qquad\qquad=M_{1}+M_{2}+M_{3}+M_{4}$

where,

$M_{1}=\frac{\partial^{3}\Phi_{P}}{\partial x_{j}\partial x_{i}^{2}}(y_{0}%
)[$g$_{jk}(\nabla\log\Phi_{P})_{k}](y_{0})$

$+\frac{\partial^{2}\Phi_{P}}{\partial x_{j}\partial x_{i}}(y_{0}%
)[\frac{\partial g_{jk}}{\partial x_{i}}(\nabla\log\Phi_{P})_{k}+g_{jk}%
\frac{\partial}{\partial x_{i}}(\nabla\log\Phi_{P})_{k}](y_{0})$

$M_{2}=\frac{\partial^{2}\Phi_{P}}{\partial x_{i}^{2}}(y_{0})[\frac{\partial
g_{jk}}{\partial x_{j}}(\nabla\log\Phi_{P})_{k}+$ g$_{jk}\frac{\partial
}{\partial x_{j}}(\nabla$log$\Phi_{P})_{k}](y_{0})$

$M_{3}=\frac{\partial\Phi_{P}}{\partial x_{i}}(y_{0})[\frac{\partial^{2}%
g_{jk}}{\partial x_{i}\partial x_{j}}(\nabla\log\Phi_{P})_{k}+\frac{\partial
g_{jk}}{\partial x_{j}}\frac{\partial}{\partial x_{i}}(\nabla$log$\Phi
_{P})_{k}$ $](y_{0})$

$M_{4}=\frac{\partial\Phi_{P}}{\partial x_{i}}(y_{0})[\frac{\partial g_{jk}%
}{\partial x_{i}}\frac{\partial}{\partial x_{j}}(\nabla\log\Phi_{P}%
)_{k}+g_{jk}\frac{\partial^{2}}{\partial x_{i}\partial x_{j}}(\nabla\log
\Phi_{P})_{k}](y_{0})$

We compute each of the above in terms of the vector field X on M and the
geometric invariants of M.

We recall the range of indices: $i,j=q+1,...,n$ and $k=1,...,q,q+1,...,n:$

Since g$_{jk}(y_{0})=\delta_{jk}$ and so,

$M_{1}=\frac{\partial^{3}\Phi_{P}}{\partial x_{j}\partial x_{i}^{2}}%
(y_{0})[(\nabla$log$\Phi_{P})_{j}](y_{0})+\frac{\partial^{2}\Phi_{P}}{\partial
x_{j}\partial x_{i}}(y_{0})[\frac{\partial g_{jk}}{\partial x_{i}}(\nabla
$log$\Phi_{P})_{k}+\frac{\partial}{\partial x_{i}}(\nabla\log\Phi_{P}%
)_{j}](y_{0})$

Since there is summation over the index $k$ and $\frac{\partial g_{jk}%
}{\partial x_{i}}(y_{0})=0$ for $i,j,k=q+1,...,n,$ we have:

$M_{1}=\frac{\partial^{3}\Phi_{P}}{\partial x_{j}\partial x_{i}^{2}}%
(y_{0})(\nabla$log$\Phi_{P})_{j}(y_{0})+\frac{\partial^{2}\Phi_{P}}{\partial
x_{j}\partial x_{i}}(y_{0})[\underset{\text{a}=1}{\overset{q}{\sum}}%
\frac{\partial g_{j\text{a}}}{\partial x_{i}}(\nabla$log$\Phi_{P})_{\text{a}%
}+\frac{\partial}{\partial x_{i}}(\nabla\log\Phi_{P})_{j}](y_{0})$

We use the following for all computations: for a = 1,...,q by (xi) of
\textbf{Table B}$_{1};$

$(\nabla$log$\Phi_{P})_{\text{a}}(y_{0})=0$

For $j=q+1,...,n$ by (vi) of \ \textbf{Table B}$_{1}$

$(\nabla$log$\Phi_{P})_{j}(y_{0})=-X_{j}(y_{0})$

We have from $\left(  36\right)  $ above:

$\frac{\partial^{2}\Phi_{P}}{\partial x_{i}\partial x_{j}}(y_{0})=[X_{i}%
X_{j}-\frac{1}{2}\left(  \frac{\partial X_{j}}{\partial x_{i}}+\frac{\partial
X_{i}}{\partial x_{j}}\right)  ](y_{0})$

We have from $\left(  39\right)  $ above:

$\frac{\partial}{\partial x_{i}}(\nabla\log\Phi_{P})_{j}(y_{0})=-\frac{1}%
{2}\left(  \frac{\partial X_{j}}{\partial x_{i}}+\frac{\partial X_{i}%
}{\partial x_{j}}\right)  (y_{0})=\frac{\partial}{\partial x_{j}}(\nabla
\log\Phi_{P})_{i}(y_{0})$

From $\left(  54\right)  $ above

$\frac{\partial^{3}\Phi_{P}}{\partial x_{i}^{2}\partial x_{j}}(y_{0}%
)=[-X_{i}^{2}X_{j}+X_{i}\left(  \frac{\partial X_{j}}{\partial x_{i}}%
+\frac{\partial X_{i}}{\partial x_{j}}\right)  +X_{j}\frac{\partial X_{i}%
}{\partial x_{i}}\ -\frac{1}{3}\left(  \frac{\partial^{2}X_{j}}{\partial
x_{i}^{2}}+2\frac{\partial^{2}X_{i}}{\partial x_{i}\partial x_{j}}\right)
](y_{0})$

Therefore,

$M_{1}=[X_{i}^{2}X_{j}^{2}-X_{i}X_{j}\left(  \frac{\partial X_{j}}{\partial
x_{i}}+\frac{\partial X_{i}}{\partial x_{j}}\right)  -X_{j}^{2}\frac{\partial
X_{i}}{\partial x_{i}}\ +\frac{1}{3}X_{j}\left(  \frac{\partial^{2}X_{j}%
}{\partial x_{i}^{2}}+2\frac{\partial^{2}X_{i}}{\partial x_{i}\partial x_{j}%
}\right)  ](y_{0})$

\qquad$\ +[X_{i}X_{j}-\frac{1}{2}\left(  \frac{\partial X_{j}}{\partial x_{i}%
}+\frac{\partial X_{i}}{\partial x_{j}}\right)  ](y_{0})[-\frac{1}{2}\left(
\frac{\partial X_{j}}{\partial x_{i}}+\frac{\partial X_{i}}{\partial x_{j}%
}\right)  ](y_{0})$

We have finally here:

$M_{1}=[X_{i}^{2}X_{j}^{2}-X_{i}X_{j}\left(  \frac{\partial X_{j}}{\partial
x_{i}}+\frac{\partial X_{i}}{\partial x_{j}}\right)  -X_{j}^{2}\frac{\partial
X_{i}}{\partial x_{i}}\ +\frac{1}{3}X_{j}\left(  \frac{\partial^{2}X_{j}%
}{\partial x_{i}^{2}}+2\frac{\partial^{2}X_{i}}{\partial x_{i}\partial x_{j}%
}\right)  ](y_{0})\qquad\left(  65\right)  $

\qquad$\ -\frac{1}{2}[X_{i}X_{j}\left(  \frac{\partial X_{j}}{\partial x_{i}%
}+\frac{\partial X_{i}}{\partial x_{j}}\right)  ](y_{0})+\frac{1}{4}\left(
\frac{\partial X_{j}}{\partial x_{i}}+\frac{\partial X_{i}}{\partial x_{j}%
}\right)  ^{2}(y_{0})$

Next we have for $i,j=q+1,...,n$ and $k=1,...,q,q+1,...,n:$

$M_{2}=\frac{\partial^{2}\Phi_{P}}{\partial x_{i}^{2}}(y_{0})[\frac{\partial
g_{jk}}{\partial x_{j}}(\nabla$log$\Phi_{P})_{k}+$ g$_{jk}\frac{\partial
}{\partial x_{j}}(\nabla$log$\Phi_{P})_{k}](y_{0})$ \ \ \ \ \ \ \ \ \ \ \ \ \ \ 

Since g$_{jk}(y_{0})=\delta_{jk};$ $\frac{\partial g_{jk}}{\partial x_{j}%
}(y_{0})=0$ for $i,j,k=q+1,...,n$

and $(\nabla$log$\Phi_{P})_{\text{a}}(y_{0})=0$ for a = 1,...,q, we have:

$M_{2}=\frac{\partial^{2}\Phi_{P}}{\partial x_{i}^{2}}(y_{0})[\frac{\partial
}{\partial x_{j}}(\nabla$log$\Phi_{P})_{j}](y_{0})=[X_{i}^{2}-\frac{\partial
X_{i}}{\partial x_{i}}](y_{0})[-\frac{\partial X_{j}}{\partial x_{j}}](y_{0})$

$M_{2}=[-X_{i}^{2}\frac{\partial X_{j}}{\partial x_{j}}+\frac{\partial X_{i}%
}{\partial x_{i}}\frac{\partial X_{j}}{\partial x_{j}}](y_{0})\qquad
\qquad\qquad\qquad\qquad\qquad\qquad\qquad\qquad\qquad\left(  66\right)  $

Next we compute, recalling that there is summation over the index $k:$

$M_{3}=\frac{\partial\Phi_{P}}{\partial x_{i}}(y_{0})[\frac{\partial^{2}%
g_{jk}}{\partial x_{i}\partial x_{j}}(\nabla\log\Phi_{P})_{k}+\frac{\partial
g_{jk}}{\partial x_{j}}\frac{\partial}{\partial x_{i}}(\nabla$log$\Phi
_{P})_{k}$ $](y_{0})$

$=\frac{\partial\Phi_{P}}{\partial x_{i}}(y_{0})\underset{\text{a=1}%
}{\overset{q}{\sum}}[\frac{\partial^{2}g_{j\text{a}}}{\partial x_{i}\partial
x_{j}}(\nabla$log$\Phi_{P})_{\text{a}}+\frac{\partial g_{j\text{a}}}{\partial
x_{j}}\frac{\partial}{\partial x_{i}}(\nabla$log$\Phi_{P})_{\text{a}}$
$](y_{0})$

$+\frac{\partial\Phi_{P}}{\partial x_{i}}(y_{0}%
)\underset{k=q+1}{\overset{n}{\sum}}[\frac{\partial^{2}g_{jk}}{\partial
x_{i}\partial x_{j}}(\nabla$log$\Phi_{P})_{k}+\frac{\partial g_{jk}}{\partial
x_{j}}\frac{\partial}{\partial x_{i}}(\nabla$log$\Phi_{P})_{k}$ $](y_{0})$

We recall here that:

$(\nabla$log$\Phi_{P})_{\text{a}}(y_{0})=0$ and $\frac{\partial g_{j\text{a}}%
}{\partial x_{j}}=\perp_{\text{a}jj}=0;$ $\frac{\partial g_{ij}}{\partial
x_{k}}(y_{0})=0$

$(\nabla$log$\Phi_{P})_{j}(y_{0})=-X_{j}(y_{0})$

By (iii) of \textbf{Appendix A}$_{1}:$

\ $\frac{\partial^{2}\text{g}_{kl}}{\partial\text{x}_{i}\partial\text{x}_{j}%
}(y_{0})=-\frac{1}{3}(R_{ikjl}+R_{jkil})(y_{0})$

In particular,

$\frac{\partial^{2}\text{g}_{jk}}{\partial\text{x}_{i}\partial\text{x}_{j}%
}(y_{0})=-\frac{1}{3}(R_{ijjk}+R_{jjik})(y_{0})=-\frac{1}{3}(R_{ijjk}%
)(y_{0})=\frac{1}{3}(R_{jijk})(y_{0})$

Therefore we have:

$M_{3}=\frac{\partial\Phi_{P}}{\partial x_{i}}(y_{0}%
)\underset{k=q+1}{\overset{n}{\sum}}[\frac{\partial^{2}g_{jk}}{\partial
x_{i}\partial x_{j}}(\nabla$log$\Phi_{P})_{k}$ $](y_{0})$

$=-X_{i}(y_{0})\underset{k=q+1}{\overset{n}{\sum}}[\frac{1}{3}(R_{jijk}%
)(y_{0})(-X_{k})$ $](y_{0})$

We see that,

$M_{3}=\frac{1}{3}\underset{k=q+1}{\overset{n}{\sum}}[R_{jijk}X_{i}X_{k}$
$](y_{0})\qquad\qquad\qquad\qquad\qquad\qquad\qquad\qquad\left(  67\right)  $

We compute the last term $M_{4}$ of $\frac{\partial^{2}L_{1}}{\partial
x_{i}\partial x_{j}}(y_{0}):$

$M_{4}=\frac{\partial\Phi_{P}}{\partial x_{i}}(y_{0})[\frac{\partial g_{jk}%
}{\partial x_{i}}\frac{\partial}{\partial x_{j}}(\nabla$log$\Phi_{P}%
)_{k}+g_{jk}\frac{\partial^{2}}{\partial x_{i}\partial x_{j}}(\nabla\log
\Phi_{P})_{k}](y_{0})$

$=\frac{\partial\Phi_{P}}{\partial x_{i}}(y_{0})\underset{\text{a=1}%
}{\overset{q}{\sum}}[\frac{\partial g_{j\text{a}}}{\partial x_{i}}%
\frac{\partial}{\partial x_{j}}(\nabla$log$\Phi_{P})_{\text{a}}+g_{j\text{a}%
}\frac{\partial^{2}}{\partial x_{i}\partial x_{j}}(\nabla\log\Phi
_{P})_{\text{a}}](y_{0})$

$+\frac{\partial\Phi_{P}}{\partial x_{i}}(y_{0}%
)\underset{k=q+1}{\overset{n}{\sum}}[\frac{\partial g_{jk}}{\partial x_{i}%
}\frac{\partial}{\partial x_{j}}(\nabla$log$\Phi_{P})_{k}+g_{jk}\frac
{\partial^{2}}{\partial x_{i}\partial x_{j}}(\nabla\log\Phi_{P})_{k}](y_{0})$

We remind that for a = 1,...,q and $i,j,k=q+1,...,n,$ we have:

$g_{jk}(y_{0})=\delta_{jk};$ $g_{j\text{a}}(y_{0})=\delta_{\text{a}j}=0;$
$\frac{\partial g_{jk}}{\partial x_{i}}(y_{0})=0;$ $\frac{\partial
g_{j\text{a}}}{\partial x_{i}}(y_{0})=\perp_{\text{a}ij}(y_{0})$

By (xv) of \textbf{Appendix B}$_{1},$

$\frac{\partial}{\partial x_{j}}(\nabla$log$\Phi_{P})_{\text{a}}(y_{0})=$
$\underset{k=q+1}{\overset{n}{%
{\textstyle\sum}
}}X_{k}(y_{0})\perp_{\text{a}jk}(y_{0})-\frac{\partial X_{j}}{\partial
x_{\text{a}}}(y_{0})$

Therefore,

$M_{4}=\frac{\partial\Phi_{P}}{\partial x_{i}}(y_{0})\underset{\text{a=1}%
}{\overset{q}{\sum}}[\frac{\partial g_{j\text{a}}}{\partial x_{i}}%
\frac{\partial}{\partial x_{j}}(\nabla$log$\Phi_{P})_{\text{a}}](y_{0}%
)+\frac{\partial\Phi_{P}}{\partial x_{i}}(y_{0})[\frac{\partial^{2}}{\partial
x_{i}\partial x_{j}}(\nabla\log\Phi_{P})_{j}](y_{0})$

$M_{4}=-X_{i}(y_{0})\underset{\text{a=1}}{\overset{q}{\sum}}[\perp
_{\text{a}ij}(y_{0})\underset{k=q+1}{\overset{n}{%
{\textstyle\sum}
}}X_{k}\perp_{\text{a}jk}-\frac{\partial X_{j}}{\partial x_{\text{a}}}%
](y_{0})-X_{i}(y_{0})[\frac{\partial^{2}}{\partial x_{i}\partial x_{j}}%
(\nabla\log\Phi_{P})_{j}](y_{0})\qquad\left(  68\right)  $

We conclude from $\left(  64\right)  ,$ $\left(  65\right)  ,$ $\left(
66\right)  ,$ $\left(  67\right)  $ and $\left(  68\right)  $ that:

$\frac{\partial^{2}L_{1}}{\partial x_{i}\partial x_{j}}(y_{0})=M_{1}%
+M_{2}+M_{3}+M_{4}\qquad\qquad\qquad\qquad\qquad\qquad\qquad\qquad
\qquad\left(  69\right)  $

$=[X_{i}^{2}X_{j}^{2}-X_{i}X_{j}\left(  \frac{\partial X_{j}}{\partial x_{i}%
}+\frac{\partial X_{i}}{\partial x_{j}}\right)  -X_{j}^{2}\frac{\partial
X_{i}}{\partial x_{i}}\ +\frac{1}{3}X_{j}\left(  \frac{\partial^{2}X_{j}%
}{\partial x_{i}^{2}}+2\frac{\partial^{2}X_{i}}{\partial x_{i}\partial x_{j}%
}\right)  ](y_{0})\qquad\qquad\qquad\qquad\qquad M_{1}$

$-\frac{1}{2}[X_{i}X_{j}\left(  \frac{\partial X_{j}}{\partial x_{i}}%
+\frac{\partial X_{i}}{\partial x_{j}}\right)  ](y_{0})+\frac{1}{4}\left(
\frac{\partial X_{j}}{\partial x_{i}}+\frac{\partial X_{i}}{\partial x_{j}%
}\right)  ^{2}(y_{0})$

$+[-X_{i}^{2}\frac{\partial X_{j}}{\partial x_{j}}+\frac{\partial X_{i}%
}{\partial x_{i}}\frac{\partial X_{j}}{\partial x_{j}}](y_{0})\qquad\qquad
M_{2}$

$+\frac{1}{3}\underset{k=q+1}{\overset{n}{\sum}}[R_{jijk}X_{i}X_{k}$
$](y_{0})\qquad M_{3}$

$-X_{i}(y_{0})\underset{\text{a=1}}{\overset{q}{\sum}}[\perp_{\text{a}%
ij}(y_{0})\underset{k=q+1}{\overset{n}{%
{\textstyle\sum}
}}X_{k}\perp_{\text{a}jk}-\frac{\partial X_{j}}{\partial x_{\text{a}}}%
](y_{0})-X_{i}(y_{0})[\frac{\partial^{2}}{\partial x_{i}\partial x_{j}}%
(\nabla\log\Phi_{P})_{j}](y_{0})\qquad M_{4}$

We next compute $\frac{\partial^{2}L_{2}}{\partial x_{i}\partial x_{j}}%
(y_{0})$ where from $\left(  61\right)  $,

$L_{2}=$\ $\Phi_{P}(x_{0})[\frac{\partial\text{g}_{jk}}{\partial x_{i}}%
(\nabla$log$\Phi_{P})_{k}](x_{0})$

Then,

$\frac{\partial^{2}L_{2}}{\partial x_{i}\partial x_{j}}(y_{0})=\frac
{\partial^{2}}{\partial x_{i}\partial x_{j}}[\Phi_{P}\frac{\partial
\text{g}_{jk}}{\partial x_{i}}(\nabla$log$\Phi_{P})_{k}](y_{0})$

$=\frac{\partial}{\partial x_{i}}[\frac{\partial}{\partial x_{j}}\{\Phi
_{P}\frac{\partial\text{g}_{jk}}{\partial x_{i}}(\nabla$log$\Phi_{P}%
)_{k}\}](y_{0})$

$=\frac{\partial}{\partial x_{i}}[\frac{\partial}{\partial x_{j}}\{\Phi
_{P}\frac{\partial\text{g}_{jk}}{\partial x_{i}}\}(\nabla$log$\Phi_{P}%
)_{k}+\Phi_{P}\frac{\partial\text{g}_{jk}}{\partial x_{i}}\frac{\partial
}{\partial x_{j}}(\nabla\log\Phi_{P})_{k}](y_{0})$

$\frac{\partial^{2}L_{2}}{\partial x_{i}\partial x_{j}}(y_{0})=\frac{\partial
}{\partial x_{i}}[\{\frac{\partial\Phi_{P}}{\partial x_{j}}\frac
{\partial\text{g}_{jk}}{\partial x_{i}}+\Phi_{P}\frac{\partial^{2}%
\text{g}_{jk}}{\partial x_{i}\partial x_{j}}\}(\nabla$log$\Phi_{P})_{k}%
+\Phi_{P}\frac{\partial\text{g}_{jk}}{\partial x_{i}}\frac{\partial}{\partial
x_{j}}(\nabla\log\Phi_{P})_{k}](y_{0})\qquad\left(  70\right)  $

$=N_{1}+N_{2}$

where,

$N_{1}=\frac{\partial}{\partial x_{i}}[\{\frac{\partial\Phi_{P}}{\partial
x_{j}}\frac{\partial\text{g}_{jk}}{\partial x_{i}}+\Phi_{P}\frac{\partial
^{2}\text{g}_{jk}}{\partial x_{i}\partial x_{j}}\}(\nabla$log$\Phi_{P}%
)_{k}](y_{0})$

$N_{2}=\frac{\partial}{\partial x_{i}}[\Phi_{P}\frac{\partial\text{g}_{jk}%
}{\partial x_{i}}\frac{\partial}{\partial x_{j}}(\nabla\log\Phi_{P}%
)_{k}](y_{0})$

We compute each of these:

$N_{1}=[\{\frac{\partial^{2}\Phi_{P}}{\partial x_{i}\partial x_{j}}%
\frac{\partial\text{g}_{jk}}{\partial x_{i}}+\frac{\partial\Phi_{P}}{\partial
x_{j}}\frac{\partial^{2}\text{g}_{jk}}{\partial x_{i}^{2}}+\frac{\partial
\Phi_{P}}{\partial x_{i}}\frac{\partial^{2}\text{g}_{jk}}{\partial
x_{i}\partial x_{j}}+\Phi_{P}\frac{\partial^{3}\text{g}_{jk}}{\partial
x_{i}^{2}\partial x_{j}}\}(\nabla$log$\Phi_{P})_{k}](y_{0})$

$+$ $[\{\frac{\partial\Phi_{P}}{\partial x_{j}}\frac{\partial\text{g}_{jk}%
}{\partial x_{i}}+\Phi_{P}\frac{\partial^{2}\text{g}_{jk}}{\partial
x_{i}\partial x_{j}}\}\frac{\partial}{\partial x_{i}}(\nabla$log$\Phi_{P}%
)_{k}](y_{0})$

For a = 1,...,q and $i,j,k=q+1,...,n,$ we have:

$N_{1}=\underset{\text{a=1}}{\overset{q}{\sum}}[\{\frac{\partial^{2}\Phi_{P}%
}{\partial x_{i}\partial x_{j}}\frac{\partial\text{g}_{j\text{a}}}{\partial
x_{i}}+\frac{\partial\Phi_{P}}{\partial x_{j}}\frac{\partial^{2}%
\text{g}_{j\text{a}}}{\partial x_{i}^{2}}+\frac{\partial\Phi_{P}}{\partial
x_{i}}\frac{\partial^{2}\text{g}_{j\text{a}}}{\partial x_{i}\partial x_{j}%
}+\Phi_{P}\frac{\partial^{3}\text{g}_{j\text{a}}}{\partial x_{i}^{2}\partial
x_{j}}\}(\nabla$\textbf{log}$\Phi_{P})_{\text{a}}](y_{0})$

$+\underset{\text{a=1}}{\overset{q}{\sum}}$\textbf{ }$[\{\frac{\partial
\Phi_{P}}{\partial x_{j}}\frac{\partial\text{g}_{j\text{a}}}{\partial x_{i}%
}+\Phi_{P}\frac{\partial^{2}\text{g}_{j\text{a}}}{\partial x_{i}\partial
x_{j}}\}\frac{\partial}{\partial x_{i}}(\nabla$\textbf{log}$\Phi
_{P})_{\text{a}}](y_{0})$

$+\underset{k=q+1}{\overset{n}{\sum}}[\{\frac{\partial^{2}\Phi_{P}}{\partial
x_{i}\partial x_{j}}\frac{\partial\text{g}_{jk}}{\partial x_{i}}%
+\frac{\partial\Phi_{P}}{\partial x_{j}}\frac{\partial^{2}\text{g}_{jk}%
}{\partial x_{i}^{2}}+\frac{\partial\Phi_{P}}{\partial x_{i}}\frac
{\partial^{2}\text{g}_{jk}}{\partial x_{i}\partial x_{j}}+\Phi_{P}%
\frac{\partial^{3}\text{g}_{jk}}{\partial x_{i}^{2}\partial x_{j}}\}(\nabla
$log$\Phi_{P})_{k}](y_{0})$

$+\underset{k=q+1}{\overset{n}{\sum}}$ $[\{\frac{\partial\Phi_{P}}{\partial
x_{j}}\frac{\partial\text{g}_{jk}}{\partial x_{i}}+\Phi_{P}\frac{\partial
^{2}\text{g}_{jk}}{\partial x_{i}\partial x_{j}}\}\frac{\partial}{\partial
x_{i}}(\nabla$log$\Phi_{P})_{k}](y_{0})$

Since $(\nabla$log$\Phi_{P})_{\text{a}}(y_{0})=0$ and $\frac{\partial
\text{g}_{jk}}{\partial x_{i}}(y_{0})=0$ for a = 1,...,q and
$i,j,k=q+1,...,n,$ we have:

$N_{1}=$ $\underset{\text{a=1}}{\overset{q}{\sum}}[\{\frac{\partial\Phi_{P}%
}{\partial x_{j}}\frac{\partial\text{g}_{j\text{a}}}{\partial x_{i}}%
+\frac{\partial^{2}\text{g}_{j\text{a}}}{\partial x_{i}\partial x_{j}}%
\}\frac{\partial}{\partial x_{i}}(\nabla$log$\Phi_{P})_{\text{a}}](y_{0})$

$+\underset{k=q+1}{\overset{n}{\sum}}[\{\frac{\partial\Phi_{P}}{\partial
x_{j}}\frac{\partial^{2}\text{g}_{jk}}{\partial x_{i}^{2}}+\frac{\partial
\Phi_{P}}{\partial x_{i}}\frac{\partial^{2}\text{g}_{jk}}{\partial
x_{i}\partial x_{j}}+\frac{\partial^{3}\text{g}_{jk}}{\partial x_{i}%
^{2}\partial x_{j}}\}(\nabla$log$\Phi_{P})_{k}](y_{0})$

$+\underset{k=q+1}{\overset{n}{\sum}}$ $[\frac{\partial^{2}\text{g}_{jk}%
}{\partial x_{i}\partial x_{j}}\frac{\partial}{\partial x_{i}}(\nabla$%
log$\Phi_{P})_{k}](y_{0})$

From the Tables we have by (i) of \textbf{Table B}$_{4}$ and by (i) of
\textbf{TableB}$_{1};$

$\Phi_{P}(y_{0})=1;$ $\frac{\partial\Phi_{P}}{\partial x_{j}}(y_{0}%
)=-X_{j}(y_{0})$

$[(\nabla$log$\Phi_{P})_{j}](y_{0})$ $=-X_{j}(y_{0})$

Then by$\ \left(  36\right)  $ and by $\left(  39\right)  $ above:

$\frac{\partial^{2}\Phi_{P}}{\partial x_{i}\partial x_{j}}(y_{0})=X_{i}%
(y_{0})X_{j}(y_{0})-\frac{1}{2}\left(  \frac{\partial X_{j}}{\partial x_{i}%
}+\frac{\partial X_{i}}{\partial x_{j}}\right)  (y_{0})$

$\frac{\partial}{\partial x_{i}}(\nabla\log\Phi_{P})_{j}(y_{0})=-\frac{1}%
{2}\left(  \frac{\partial X_{j}}{\partial x_{i}}+\frac{\partial X_{i}%
}{\partial x_{j}}\right)  (y_{0})=\frac{\partial}{\partial x_{j}}(\nabla
\log\Phi_{P})_{i}(y_{0})$

By (iii) of \textbf{Appendix A}$_{1},$

\ $\frac{\partial^{2}\text{g}_{kl}}{\partial\text{x}_{i}\partial\text{x}_{j}%
}(y_{0})=-\frac{1}{3}(R_{ikjl}+R_{jkil})(y_{0})$

Therefore,

$\frac{\partial^{2}\text{g}_{jk}}{\partial x_{i}^{2}}=-\frac{2}{3}%
R_{ijik}(y_{0});$ $\frac{\partial^{2}\text{g}_{jk}}{\partial\text{x}%
_{i}\partial\text{x}_{j}}(y_{0})=\frac{1}{3}R_{jijk}(y_{0})$

\ $\frac{\partial^{2}\text{g}_{kl}}{\partial\text{x}_{i}\partial\text{x}_{j}%
}(y_{0})=-\frac{1}{3}(R_{ikjl}+R_{jkil})(y_{0})$

$\frac{\partial\text{g}_{\text{a}j}}{\partial\text{x}_{i}}(y_{0}%
)=\perp_{\text{a}ij}(y_{0})$ by (ii) of \textbf{Table A}$_{3}$ and so
$\frac{\partial\text{g}_{j\text{a}}}{\partial x_{j}}(y_{0})=\perp_{\text{a}%
jj}(y_{0})=0$

$\frac{\partial^{2}\text{g}_{\text{a}k}}{\partial\text{x}_{i}\partial
\text{x}_{j}}(y_{0})=-\frac{4}{3}(R_{i\text{a}jk}+R_{j\text{a}ik})(y_{0})$ by
(iii) of \textbf{Table A}$_{3};$ $\frac{\partial^{2}\text{g}_{\text{a}j}%
}{\partial\text{x}_{i}^{2}}(y_{0})=-\frac{8}{3}R_{i\text{a}ij}(y_{0})$

$\frac{\partial^{2}\text{g}_{j\text{a}}}{\partial x_{i}\partial x_{j}}%
(y_{0})=-\frac{4}{3}(R_{i\text{a}jj}+R_{j\text{a}ij})(y_{0})=-\frac{4}%
{3}R_{j\text{a}ij}(y_{0})$

$\frac{\partial^{3}\text{g}_{jk}}{\partial x_{i}^{2}\partial x_{j}}%
(y_{0})=-\frac{1}{3}(\nabla_{i}R_{ijjk}+\nabla_{j}R_{ijik})(y_{0})$ by
(iv)$^{\ast}$ of \textbf{Table A}$_{1}$

$\frac{\partial}{\partial x_{i}}(\nabla$log$\Phi_{P})_{\text{a}}(y_{0})=$
$\underset{k=q+1}{\overset{n}{%
{\textstyle\sum}
}}X_{k}(y_{0})\perp_{\text{a}ik}(y_{0})-\frac{\partial X_{i}}{\partial
x_{\text{a}}}(y_{0})$ by (xv) of \textbf{Table B}$_{1}$

Therefore,

$N_{1}=$ $\underset{\text{a=1}}{\overset{q}{\sum}}[\{-X_{j}\perp_{\text{a}%
ij}-\frac{4}{3}R_{j\text{a}ij}\}\{\underset{k=q+1}{\overset{n}{%
{\textstyle\sum}
}}X_{k}\perp_{\text{a}ik}-\frac{\partial X_{i}}{\partial x_{\text{a}}%
}\}](y_{0})$

$+\underset{k=q+1}{\overset{n}{\sum}}[\{(-X_{j})(-\frac{2}{3}R_{ijik}%
)+(-X_{i})(\frac{1}{3}R_{jijk})-\frac{1}{3}(\nabla_{i}R_{ijjk}+\nabla
_{j}R_{ijik})\}(-X_{k})](y_{0})$

$+\underset{k=q+1}{\overset{n}{\sum}}$ $[\frac{1}{3}R_{jijk}(y_{0})(-\frac
{1}{2}\left(  \frac{\partial X_{k}}{\partial x_{i}}+\frac{\partial X_{i}%
}{\partial x_{k}}\right)  )](y_{0})$

Simplifying, we have:

$N_{1}=\frac{1}{3}\underset{k=q+1}{\overset{n}{\sum}}[R_{jijk}X_{i}%
X_{k}-2R_{ijik}X_{j}X_{k}+(\nabla_{i}R_{ijjk}+\nabla_{j}R_{ijik})X_{k}%
](y_{0})\qquad\left(  71\right)  $

$-\frac{1}{6}\underset{k=q+1}{\overset{n}{\sum}}$ $[R_{jijk}\left(
\frac{\partial X_{k}}{\partial x_{i}}+\frac{\partial X_{i}}{\partial x_{k}%
}\right)  ](y_{0})$

$+$ $\underset{\text{a=1}}{\overset{q}{\sum}}[\underset{k=q+1}{\overset{n}{%
{\textstyle\sum}
}}-X_{j}X_{k}\perp_{\text{a}ij}\perp_{\text{a}ik}+X_{j}\frac{\partial X_{i}%
}{\partial x_{\text{a}}}\perp_{\text{a}ij}](y_{0})\qquad$

$+\frac{4}{3}\underset{\text{a=1}}{\overset{q}{\sum}}$
$[\underset{k=q+1}{\overset{n}{%
{\textstyle\sum}
}}-R_{j\text{a}ij}X_{k}\perp_{\text{a}ik}+R_{j\text{a}ij}\frac{\partial X_{i}%
}{\partial x_{\text{a}}}](y_{0})$

\qquad\qquad\qquad\qquad\qquad\qquad\qquad\qquad\qquad\qquad\qquad\qquad
\qquad\qquad\qquad\qquad$\blacksquare$

We next compute:

$N_{2}=\frac{\partial}{\partial x_{i}}[\Phi_{P}\frac{\partial\text{g}_{jk}%
}{\partial x_{i}}\frac{\partial}{\partial x_{j}}(\nabla\log\Phi_{P}%
)_{k}](y_{0})$

$=\frac{\partial}{\partial x_{i}}[\Phi_{P}\frac{\partial\text{g}_{jk}%
}{\partial x_{i}}](y_{0})\frac{\partial}{\partial x_{j}}(\nabla\log\Phi
_{P})_{k}(y_{0})+[\Phi_{P}\frac{\partial\text{g}_{jk}}{\partial x_{i}}%
\frac{\partial^{2}}{\partial x_{i}\partial x_{j}}(\nabla\log\Phi_{P}%
)_{k}](y_{0})$

$=\frac{\partial}{\partial x_{i}}[\Phi_{P}\frac{\partial\text{g}_{jk}%
}{\partial x_{i}}](y_{0})\frac{\partial}{\partial x_{j}}(\nabla\log\Phi
_{P})_{k}(y_{0})+[\Phi_{P}\frac{\partial\text{g}_{jk}}{\partial x_{i}}%
\frac{\partial^{2}}{\partial x_{i}\partial x_{j}}(\nabla\log\Phi_{P}%
)_{k}](y_{0})$

$N_{2}=[\frac{\partial\Phi_{P}}{\partial x_{i}}\frac{\partial\text{g}_{jk}%
}{\partial x_{i}}+\Phi_{P}\frac{\partial^{2}\text{g}_{jk}}{\partial x_{i}^{2}%
}](y_{0})\frac{\partial}{\partial x_{j}}(\nabla\log\Phi_{P})_{k}(y_{0}%
)+[\Phi_{P}\frac{\partial\text{g}_{jk}}{\partial x_{i}}\frac{\partial^{2}%
}{\partial x_{i}\partial x_{j}}(\nabla\log\Phi_{P})_{k}](y_{0})$

Since there is summation over $k=1,...,q,q+1,...,n,$ we have:

$N_{2}=\underset{\text{a=1}}{\overset{q}{\sum}}[\frac{\partial\Phi_{P}%
}{\partial x_{i}}\frac{\partial\text{g}_{j\text{a}}}{\partial x_{i}}+\Phi
_{P}\frac{\partial^{2}\text{g}_{j\text{a}}}{\partial x_{i}^{2}}](y_{0}%
)\frac{\partial}{\partial x_{j}}(\nabla\log\Phi_{P})_{\text{a}}(y_{0}%
)+\underset{\text{a=1}}{\overset{q}{\sum}}[\Phi_{P}\frac{\partial
\text{g}_{j\text{a}}}{\partial x_{i}}\frac{\partial^{2}}{\partial
x_{i}\partial x_{j}}(\nabla\log\Phi_{P})_{\text{a}}](y_{0})$

$+\underset{k=q+1}{\overset{n}{\sum}}[\frac{\partial\Phi_{P}}{\partial x_{i}%
}\frac{\partial\text{g}_{jk}}{\partial x_{i}}+\Phi_{P}\frac{\partial
^{2}\text{g}_{jk}}{\partial x_{i}^{2}}](y_{0})\frac{\partial}{\partial x_{j}%
}(\nabla\log\Phi_{P})_{k}(y_{0})+\underset{k=q+1}{\overset{n}{\sum}}[\Phi
_{P}\frac{\partial\text{g}_{jk}}{\partial x_{i}}\frac{\partial^{2}}{\partial
x_{i}\partial x_{j}}(\nabla\log\Phi_{P})_{k}](y_{0})$

Values from Tables used for computing $N_{1}$ apply here except for
$\frac{\partial^{3}\Phi_{P}}{\partial x_{i}^{2}\partial x_{j}}(y_{0}).$

By (v)$^{\ast}$ of Table B$_{4},$ we have:

$\frac{\partial^{3}\Phi_{P}}{\partial x_{i}^{2}\partial x_{j}}(y_{0}%
)=[-X_{j}X_{i}^{2}+X_{j}\frac{\partial X_{i}}{\partial x_{i}}\ +X_{i}\left(
\frac{\partial X_{j}}{\partial x_{i}}+\frac{\partial X_{i}}{\partial x_{j}%
}\right)  -\frac{1}{3}\left(  \frac{\partial^{2}X_{j}}{\partial x_{i}^{2}%
}+2\frac{\partial^{2}X_{i}}{\partial x_{i}\partial x_{j}}\right)  ](y_{0})$

Therefore,

$N_{2}=\underset{\text{a=1}}{\overset{q}{\sum}}[-X_{i}\perp_{\text{a}ij}%
-\frac{8}{3}R_{i\text{a}ij}](y_{0})[\underset{k=q+1}{\overset{n}{%
{\textstyle\sum}
}}X_{k}\perp_{\text{a}jk}-\frac{\partial X_{j}}{\partial x_{\text{a}}}%
](y_{0})\qquad\qquad\qquad\left(  72\right)  $

$+\underset{\text{a=1}}{\overset{q}{\sum}}\perp_{\text{a}ij}(y_{0}%
)\frac{\partial^{2}}{\partial x_{i}\partial x_{j}}(\nabla\log\Phi
_{P})_{\text{a}}(y_{0})\qquad$

$+\underset{k=q+1}{\overset{n}{\sum}}\frac{1}{3}R_{ijik}(y_{0})\left(
\frac{\partial X_{k}}{\partial x_{j}}+\frac{\partial X_{j}}{\partial x_{k}%
}\right)  (y_{0})$

\qquad\qquad\qquad\qquad\qquad\qquad\qquad$\blacksquare$

By $\left(  69\right)  ,$ $\left(  70\right)  ,$ $\left(  71\right)  $ and
$\left(  72\right)  ,$ we have:

$\frac{\partial^{2}L_{2}}{\partial x_{i}\partial x_{j}}(y_{0})=$
$\underset{\text{a=1}}{\overset{q}{\sum}}[\underset{k=q+1}{\overset{n}{%
{\textstyle\sum}
}}-X_{j}X_{k}\perp_{\text{a}ij}\perp_{\text{a}ik}+X_{j}\frac{\partial X_{i}%
}{\partial x_{\text{a}}}\perp_{\text{a}ij}](y_{0})\qquad N_{1}\qquad\left(
73\right)  $

$+\frac{4}{3}\underset{\text{a=1}}{\overset{q}{\sum}}$
$[\underset{k=q+1}{\overset{n}{%
{\textstyle\sum}
}}-R_{j\text{a}ij}X_{k}\perp_{\text{a}ik}+R_{j\text{a}ij}\frac{\partial X_{i}%
}{\partial x_{\text{a}}}](y_{0})$

$+\underset{k=q+1}{\overset{n}{\sum}}[\frac{1}{3}R_{jijk}X_{i}X_{k}-\frac
{2}{3}R_{ijik}X_{j}X_{k}+\frac{1}{3}(\nabla_{i}R_{ijjk}+\nabla_{j}%
R_{ijik})X_{k}](y_{0})$

$-\frac{1}{6}\underset{k=q+1}{\overset{n}{\sum}}$ $[R_{jijk}\left(
\frac{\partial X_{k}}{\partial x_{i}}+\frac{\partial X_{i}}{\partial x_{k}%
}\right)  ](y_{0})$

$+\underset{\text{a=1}}{\overset{q}{\sum}}[-X_{i}\perp_{\text{a}ij}-\frac
{8}{3}R_{i\text{a}ij}](y_{0})[\underset{k=q+1}{\overset{n}{%
{\textstyle\sum}
}}X_{k}\perp_{\text{a}jk}-\frac{\partial X_{j}}{\partial x_{\text{a}}}%
](y_{0})\qquad N_{2}$

$+\underset{\text{a=1}}{\overset{q}{\sum}}\perp_{\text{a}ij}(y_{0}%
)\frac{\partial^{2}}{\partial x_{i}\partial x_{j}}(\nabla$log$\Phi
_{P})_{\text{a}}(y_{0})$

$+\frac{1}{3}\underset{k=q+1}{\overset{n}{\sum}}R_{ijik}(y_{0})\left(
\frac{\partial X_{k}}{\partial x_{j}}+\frac{\partial X_{j}}{\partial x_{k}%
}\right)  (y_{0})$

\qquad\qquad\qquad\qquad\qquad\qquad\qquad$\blacksquare$

By (xvi) of \textbf{Table B}$_{1}$ of \textbf{Appendix B,}

$\frac{\partial^{2}}{\partial x_{i}\partial x_{j}}(\nabla$log$\Phi
_{P})_{\text{a}}(y_{0})=-$ $2\underset{\text{b=1}}{\overset{\text{q}}{%
{\textstyle\sum}
}}$T$_{\text{ab}j}(y_{0})\frac{\partial X_{i}}{\partial x_{\text{b}}}%
(y_{0})-2\underset{\text{b=1}}{\overset{\text{q}}{%
{\textstyle\sum}
}}T_{\text{ab}i}(y)\frac{\partial X_{j}}{\partial x_{\text{b}}}(y_{0})$

$+\frac{1}{2}\underset{k=q+1}{\overset{n}{%
{\textstyle\sum}
}}\perp_{\text{a}jk}(y_{0})\left[  \left(  \frac{\partial X_{i}}{\partial
x_{k}}+\frac{\partial X_{k}}{\partial x_{i}}\right)  \right]  (y_{0})+\frac
{1}{2}\underset{k=q+1}{\overset{n}{%
{\textstyle\sum}
}}\perp_{\text{a}ik}(y_{0})[\left(  \frac{\partial X_{k}}{\partial x_{j}%
}+\frac{\partial X_{j}}{\partial x_{k}}\right)  ](y_{0})$\qquad

$+\frac{4}{3}\underset{k=q+1}{\overset{n}{%
{\textstyle\sum}
}}\left(  R_{i\text{a}jk}+R_{j\text{a}ik}\right)  (y_{0})X_{k}(y_{0}%
)+[X_{i}\frac{\partial X_{j}}{\partial x_{\text{a}}}+X_{j}\frac{\partial
X_{i}}{\partial x_{\text{a}}}-\frac{1}{2}\left(  \frac{\partial^{2}X_{i}%
}{\partial x_{\text{a}}\partial x_{j}}+\frac{\partial^{2}X_{j}}{\partial
x_{\text{a}}\partial x_{i}}\right)  ](y_{0})$

We insert the expression for $\frac{\partial^{2}}{\partial x_{i}\partial
x_{j}}[(\nabla$log$\Phi_{P})_{\text{a}}](y_{0})$ and have:

$\frac{\partial^{2}L_{2}}{\partial x_{i}\partial x_{j}}(y_{0})=$
$\underset{\text{a=1}}{\overset{q}{\sum}}[\underset{k=q+1}{\overset{n}{%
{\textstyle\sum}
}}-X_{j}X_{k}\perp_{\text{a}ij}\perp_{\text{a}ik}+X_{j}\frac{\partial X_{i}%
}{\partial x_{\text{a}}}\perp_{\text{a}ij}](y_{0})\qquad N_{1}$

$+\frac{4}{3}\underset{\text{a=1}}{\overset{q}{\sum}}$
$[\underset{k=q+1}{\overset{n}{%
{\textstyle\sum}
}}-R_{j\text{a}ij}X_{k}\perp_{\text{a}ik}+R_{j\text{a}ij}\frac{\partial X_{i}%
}{\partial x_{\text{a}}}](y_{0})$

$+\underset{k=q+1}{\overset{n}{\sum}}[\frac{1}{3}R_{jijk}X_{i}X_{k}-\frac
{2}{3}R_{ijik}X_{j}X_{k}+\frac{1}{3}(\nabla_{i}R_{ijjk}+\nabla_{j}%
R_{ijik})X_{k}](y_{0})$

$-\frac{1}{6}\underset{k=q+1}{\overset{n}{\sum}}$ $[R_{jijk}\left(
\frac{\partial X_{k}}{\partial x_{i}}+\frac{\partial X_{i}}{\partial x_{k}%
}\right)  ](y_{0})$

$+\underset{\text{a=1}}{\overset{q}{\sum}}[-X_{i}\perp_{\text{a}ij}-\frac
{8}{3}R_{i\text{a}ij}](y_{0})[\underset{k=q+1}{\overset{n}{%
{\textstyle\sum}
}}X_{k}\perp_{\text{a}jk}-\frac{\partial X_{j}}{\partial x_{\text{a}}}%
](y_{0})\qquad N_{2}$

$-$ $2\underset{\text{a,b=1}}{\overset{\text{q}}{%
{\textstyle\sum}
}}$T$_{\text{ab}j}(y_{0})\perp_{\text{a}ij}(y_{0})\frac{\partial X_{i}%
}{\partial x_{\text{b}}}(y_{0})-2\underset{\text{a,b=1}}{\overset{\text{q}}{%
{\textstyle\sum}
}}T_{\text{ab}i}(y)\perp_{\text{a}ij}(y_{0})\frac{\partial X_{j}}{\partial
x_{\text{b}}}(y_{0})$

$+\frac{1}{2}\underset{k=q+1}{\overset{n}{%
{\textstyle\sum}
}}\underset{\text{a=1}}{\overset{q}{\sum}}\perp_{\text{a}ij}(y_{0}%
)\perp_{\text{a}jk}(y_{0})\left[  \left(  \frac{\partial X_{i}}{\partial
x_{k}}+\frac{\partial X_{k}}{\partial x_{i}}\right)  \right]  (y_{0})$

$+\frac{1}{2}\underset{k=q+1}{\overset{n}{%
{\textstyle\sum}
}}\underset{\text{a=1}}{\overset{q}{\sum}}\perp_{\text{a}ij}(y_{0}%
)\perp_{\text{a}ik}(y_{0})[\left(  \frac{\partial X_{k}}{\partial x_{j}}%
+\frac{\partial X_{j}}{\partial x_{k}}\right)  ](y_{0})$\qquad

$+\frac{4}{3}\underset{k=q+1}{\overset{n}{%
{\textstyle\sum}
}}\underset{\text{a=1}}{\overset{q}{\sum}}R_{i\text{a}jk}(y_{0})\perp
_{\text{a}ij}(y_{0})X_{k}(y_{0})+\frac{4}{3}\underset{k=q+1}{\overset{n}{%
{\textstyle\sum}
}}\underset{\text{a=1}}{\overset{q}{\sum}}R_{j\text{a}ik}(y_{0})\perp
_{\text{a}ij}(y_{0})X_{k}(y_{0})$

$+\underset{\text{a=1}}{\overset{q}{\sum}}\perp_{\text{a}ij}(y_{0})[X_{i}%
\frac{\partial X_{j}}{\partial x_{\text{a}}}+X_{j}\frac{\partial X_{i}%
}{\partial x_{\text{a}}}-\frac{1}{2}\left(  \frac{\partial^{2}X_{i}}{\partial
x_{\text{a}}\partial x_{j}}+\frac{\partial^{2}X_{j}}{\partial x_{\text{a}%
}\partial x_{i}}\right)  ](y_{0})$

$+\frac{1}{3}\underset{k=q+1}{\overset{n}{\sum}}R_{ijik}(y_{0})\left(
\frac{\partial X_{k}}{\partial x_{j}}+\frac{\partial X_{j}}{\partial x_{k}%
}\right)  (y_{0})$

\qquad\qquad\qquad\qquad\qquad\qquad\qquad\qquad\qquad\qquad\qquad\qquad
\qquad\qquad\qquad\qquad\qquad\qquad\qquad$\blacksquare$

We re-write the expression of $\frac{\partial^{2}L_{2}}{\partial x_{i}\partial
x_{j}}(y_{0})$ in which we allign similar terms:

$\frac{\partial^{2}L_{2}}{\partial x_{i}\partial x_{j}}(y_{0})=\frac{1}%
{3}\underset{k=q+1}{\overset{n}{\sum}}[R_{jijk}X_{i}X_{k}-2R_{ijik}X_{j}%
X_{k}+(\nabla_{i}R_{ijjk}+\nabla_{j}R_{ijik})X_{k}](y_{0})\qquad\left(
73\right)  $

$-\frac{1}{6}\underset{k=q+1}{\overset{n}{\sum}}$ $R_{jijk}\left(
\frac{\partial X_{k}}{\partial x_{i}}+\frac{\partial X_{i}}{\partial x_{k}%
}\right)  (y_{0})+\frac{1}{3}\underset{k=q+1}{\overset{n}{\sum}}R_{ijik}%
(y_{0})\left(  \frac{\partial X_{k}}{\partial x_{j}}+\frac{\partial X_{j}%
}{\partial x_{k}}\right)  (y_{0})\qquad$

$-$ $2\underset{\text{a,b=1}}{\overset{\text{q}}{%
{\textstyle\sum}
}}$T$_{\text{ab}j}(y_{0})\perp_{\text{a}ij}(y_{0})\frac{\partial X_{i}%
}{\partial x_{\text{b}}}(y_{0})-2\underset{\text{a,b=1}}{\overset{\text{q}}{%
{\textstyle\sum}
}}T_{\text{ab}i}(y)\perp_{\text{a}ij}(y_{0})\frac{\partial X_{j}}{\partial
x_{\text{b}}}(y_{0})$

$+\frac{1}{2}\underset{k=q+1}{\overset{n}{%
{\textstyle\sum}
}}\underset{\text{a=1}}{\overset{q}{\sum}}\perp_{\text{a}ij}(y_{0}%
)\perp_{\text{a}jk}(y_{0})\left[  \left(  \frac{\partial X_{i}}{\partial
x_{k}}+\frac{\partial X_{k}}{\partial x_{i}}\right)  \right]  (y_{0})$

$+\frac{1}{2}\underset{k=q+1}{\overset{n}{%
{\textstyle\sum}
}}\underset{\text{a=1}}{\overset{q}{\sum}}\perp_{\text{a}ij}(y_{0}%
)\perp_{\text{a}ik}(y_{0})\left(  \frac{\partial X_{k}}{\partial x_{j}}%
+\frac{\partial X_{j}}{\partial x_{k}}\right)  (y_{0})$\qquad

$+\frac{4}{3}\underset{k=q+1}{\overset{n}{%
{\textstyle\sum}
}}\underset{\text{a=1}}{\overset{q}{\sum}}R_{i\text{a}jk}(y_{0})\perp
_{\text{a}ij}(y_{0})X_{k}(y_{0})+\frac{4}{3}\underset{k=q+1}{\overset{n}{%
{\textstyle\sum}
}}\underset{\text{a=1}}{\overset{q}{\sum}}R_{j\text{a}ik}(y_{0})\perp
_{\text{a}ij}(y_{0})X_{k}(y_{0})$

$+\underset{\text{a=1}}{\overset{q}{\sum}}\perp_{\text{a}ij}(y_{0})[X_{i}%
\frac{\partial X_{j}}{\partial x_{\text{a}}}+X_{j}\frac{\partial X_{i}%
}{\partial x_{\text{a}}}-\frac{1}{2}\left(  \frac{\partial^{2}X_{i}}{\partial
x_{\text{a}}\partial x_{j}}+\frac{\partial^{2}X_{j}}{\partial x_{\text{a}%
}\partial x_{i}}\right)  ](y_{0})$

$+$ $\underset{\text{a=1}}{\overset{q}{\sum}}[\underset{k=q+1}{\overset{n}{%
{\textstyle\sum}
}}-X_{j}X_{k}\perp_{\text{a}ij}\perp_{\text{a}ik}+X_{j}\frac{\partial X_{i}%
}{\partial x_{\text{a}}}\perp_{\text{a}ij}](y_{0})\qquad$

$+\underset{k=q+1}{\overset{n}{%
{\textstyle\sum}
}}\underset{\text{a=1}}{\overset{q}{\sum}}[-\perp_{\text{a}ij}\perp
_{\text{a}jk}X_{i}X_{k}](y_{0})+$ $\underset{\text{a=1}}{\overset{q}{\sum}%
}[\perp_{\text{a}ij}X_{i}\frac{\partial X_{j}}{\partial x_{\text{a}}}](y_{0})$

$-\frac{8}{3}\underset{k=q+1}{\overset{n}{%
{\textstyle\sum}
}}$ $\underset{\text{a=1}}{\overset{q}{\sum}}[R_{i\text{a}ij}\perp
_{\text{a}jk}X_{k}](y_{0})+\frac{8}{3}\underset{\text{a=1}}{\overset{q}{\sum}%
}[R_{i\text{a}ij}\frac{\partial X_{j}}{\partial x_{\text{a}}}](y_{0})$ \ \ 

$+\frac{4}{3}\underset{\text{a=1}}{\overset{q}{\sum}}$
$[\underset{k=q+1}{\overset{n}{%
{\textstyle\sum}
}}R_{j\text{a}ji}\perp_{\text{a}ik}X_{k}](y_{0})-\frac{4}{3}%
\underset{\text{a=1}}{\overset{q}{\sum}}[R_{j\text{a}ji}\frac{\partial X_{i}%
}{\partial x_{\text{a}}}](y_{0})$\ 

\qquad\qquad\qquad\qquad\qquad\qquad\qquad\qquad\qquad\qquad\qquad\qquad
\qquad\qquad\qquad\qquad\qquad\qquad$\blacksquare$

We now compute $\frac{\partial^{2}L_{3}}{\partial x_{i}\partial x_{j}}(y_{0})$
where by $\left(  62\right)  ,$

$L_{3}=\Phi_{P}(x_{0})[$ g$_{jk}\frac{\partial}{\partial x_{i}}(\nabla
$log$\Phi_{P})_{k}](x_{0})$

Then we have,

$\frac{\partial^{2}L_{3}}{\partial x_{i}\partial x_{j}}(y_{0})=\frac
{\partial^{2}}{\partial x_{i}\partial x_{j}}[\Phi_{P}$ g$_{jk}\frac{\partial
}{\partial x_{i}}(\nabla$log$\Phi_{P})_{k}](y_{0})$

$=\frac{\partial}{\partial x_{i}}[\frac{\partial}{\partial x_{j}}(\Phi_{P}$
g$_{jk})\frac{\partial}{\partial x_{i}}(\nabla$log$\Phi_{P})_{k}+\Phi_{P}%
$g$_{jk}\frac{\partial^{2}}{\partial x_{j}\partial x_{i}}(\nabla\log\Phi
_{P})_{k}](y_{0})=Q_{1}+Q_{2}$

where,

\ \ $Q_{1}=\frac{\partial}{\partial x_{i}}[\frac{\partial}{\partial x_{j}%
}(\Phi_{P}$ g$_{jk})\frac{\partial}{\partial x_{i}}(\nabla$log$\Phi_{P}%
)_{k}](y_{0})$

$Q_{2}=\frac{\partial}{\partial x_{i}}[\Phi_{P}$g$_{jk}\frac{\partial^{2}%
}{\partial x_{j}\partial x_{i}}(\nabla\log\Phi_{P})_{k}](y_{0})$

We compute each of the above expressions:

$Q_{1}=\frac{\partial}{\partial x_{i}}[\frac{\partial}{\partial x_{j}}%
(\Phi_{P}$ g$_{jk})\frac{\partial}{\partial x_{i}}(\nabla$log$\Phi_{P}%
)_{k}](y_{0})$

$=\frac{\partial}{\partial x_{i}}[\{\frac{\partial\Phi_{P}}{\partial x_{j}}%
$g$_{jk}+\Phi_{P}\frac{\partial\text{g}_{jk}}{\partial x_{j}}\}\frac{\partial
}{\partial x_{i}}(\nabla$log$\Phi_{P})_{k}](y_{0})$

$=[\{\frac{\partial^{2}\Phi_{P}}{\partial x_{i}\partial x_{j}}$g$_{jk}%
+\frac{\partial\Phi_{P}}{\partial x_{j}}\frac{\partial\text{g}_{jk}}{\partial
x_{i}}+\frac{\partial\Phi_{P}}{\partial x_{i}}\frac{\partial\text{g}_{jk}%
}{\partial x_{j}}+\Phi_{P}\frac{\partial^{2}\text{g}_{jk}}{\partial
x_{i}\partial x_{j}}\}\frac{\partial}{\partial x_{i}}(\nabla$log$\Phi_{P}%
)_{k}](y_{0})$

$+[\{\frac{\partial\Phi_{P}}{\partial x_{j}}$g$_{jk}+\Phi_{P}\frac
{\partial\text{g}_{jk}}{\partial x_{j}}\}\frac{\partial^{2}}{\partial
x_{i}^{2}}(\nabla$log$\Phi_{P})_{k}](y_{0})$

Since g$_{jk}(y_{0})=\delta_{jk}$ and $\Phi_{P}(y_{0})=1,$ we have:

$Q_{1}=[\{\frac{\partial\Phi_{P}}{\partial x_{j}}\frac{\partial\text{g}_{jk}%
}{\partial x_{i}}+\frac{\partial\Phi_{P}}{\partial x_{i}}\frac{\partial
\text{g}_{jk}}{\partial x_{j}}+\frac{\partial^{2}\text{g}_{jk}}{\partial
x_{i}\partial x_{j}}\}\frac{\partial}{\partial x_{i}}(\nabla$log$\Phi_{P}%
)_{k}](y_{0})$

$+\frac{\partial^{2}\Phi_{P}}{\partial x_{i}\partial x_{j}}\frac{\partial
}{\partial x_{i}}(\nabla\log\Phi_{P})_{j}(y_{0})$

$+$ $[\frac{\partial\text{g}_{jk}}{\partial x_{j}}\frac{\partial^{2}}{\partial
x_{i}^{2}}(\nabla$log$\Phi_{P})_{k}](y_{0})+[\frac{\partial\Phi_{P}}{\partial
x_{j}}\frac{\partial^{2}}{\partial x_{i}^{2}}(\nabla$log$\Phi_{P})_{j}%
](y_{0})$

We re-write the last expression above as:

$Q_{1}=\frac{\partial^{2}\Phi_{P}}{\partial x_{i}\partial x_{j}}\frac
{\partial}{\partial x_{i}}(\nabla\log\Phi_{P})_{j}(y_{0})+[\frac{\partial
\Phi_{P}}{\partial x_{j}}\frac{\partial^{2}}{\partial x_{i}^{2}}(\nabla
$log$\Phi_{P})_{j}](y_{0})$

$+\frac{\partial^{2}\text{g}_{jk}}{\partial x_{i}\partial x_{j}}\frac
{\partial}{\partial x_{i}}(\nabla$log$\Phi_{P})_{k}](y_{0})$

$+$ $[\{\frac{\partial\Phi_{P}}{\partial x_{j}}\frac{\partial\text{g}_{jk}%
}{\partial x_{i}}+\frac{\partial\Phi_{P}}{\partial x_{i}}\frac{\partial
\text{g}_{jk}}{\partial x_{j}}\}\frac{\partial}{\partial x_{i}}(\nabla
$log$\Phi_{P})_{k}](y_{0})+$ $[\frac{\partial\text{g}_{jk}}{\partial x_{j}%
}\frac{\partial^{2}}{\partial x_{i}^{2}}(\nabla$log$\Phi_{P})_{k}](y_{0})$

Recall that $\frac{\partial\text{g}_{jk}}{\partial x_{i}}(y_{0})=0$ for
$i,j,k=q+1,...,n+1.$

$Q_{1}=\frac{\partial^{2}\Phi_{P}}{\partial x_{i}\partial x_{j}}\frac
{\partial}{\partial x_{i}}(\nabla\log\Phi_{P})_{j}(y_{0})+[\frac{\partial
\Phi_{P}}{\partial x_{j}}\frac{\partial^{2}}{\partial x_{i}^{2}}(\nabla
$log$\Phi_{P})_{j}](y_{0})$

$+\underset{k=q+1}{\overset{n}{%
{\textstyle\sum}
}}[\frac{\partial^{2}\text{g}_{jk}}{\partial x_{i}\partial x_{j}}%
\frac{\partial}{\partial x_{i}}(\nabla$log$\Phi_{P})_{k}](y_{0})$

$+$ $\underset{\text{a=1}}{\overset{q}{\sum}}[\frac{\partial^{2}%
\text{g}_{j\text{a}}}{\partial x_{i}\partial x_{j}}\frac{\partial}{\partial
x_{i}}(\nabla$log$\Phi_{P})_{\text{a}}](y_{0})+$ $\underset{\text{a=1}%
}{\overset{q}{\sum}}$ $[\{\frac{\partial\Phi_{P}}{\partial x_{j}}%
\frac{\partial\text{g}_{j\text{a}}}{\partial x_{i}}+\frac{\partial\Phi_{P}%
}{\partial x_{i}}\frac{\partial\text{g}_{j\text{a}}}{\partial x_{j}}%
\}\frac{\partial}{\partial x_{i}}(\nabla$log$\Phi_{P})_{\text{a}}](y_{0})$

$+$ $[\frac{\partial\text{g}_{j\text{a}}}{\partial x_{j}}\frac{\partial^{2}%
}{\partial x_{i}^{2}}(\nabla$log$\Phi_{P})_{\text{a}}](y_{0})$

The values of all terms of $Q_{1}$ have already been given in prprevious calculations:

$\qquad Q_{1}=\frac{\partial^{2}\Phi_{P}}{\partial x_{i}\partial x_{j}}%
\frac{\partial}{\partial x_{i}}(\nabla\log\Phi_{P})_{j}(y_{0})+[\frac
{\partial\Phi_{P}}{\partial x_{j}}\frac{\partial^{2}}{\partial x_{i}^{2}%
}(\nabla$log$\Phi_{P})_{j}](y_{0})$

$+\underset{k=q+1}{\overset{n}{%
{\textstyle\sum}
}}[\frac{\partial^{2}\text{g}_{jk}}{\partial x_{i}\partial x_{j}}%
\frac{\partial}{\partial x_{i}}(\nabla$log$\Phi_{P})_{k}](y_{0})$

$+$ $\underset{\text{a=1}}{\overset{q}{\sum}}[\frac{\partial^{2}%
\text{g}_{j\text{a}}}{\partial x_{i}\partial x_{j}}\frac{\partial}{\partial
x_{i}}(\nabla$log$\Phi_{P})_{\text{a}}](y_{0})+$ $\underset{\text{a=1}%
}{\overset{q}{\sum}}$ $[\{\frac{\partial\Phi_{P}}{\partial x_{j}}%
\frac{\partial\text{g}_{j\text{a}}}{\partial x_{i}}+\frac{\partial\Phi_{P}%
}{\partial x_{i}}\frac{\partial\text{g}_{j\text{a}}}{\partial x_{j}}%
\}\frac{\partial}{\partial x_{i}}(\nabla$log$\Phi_{P})_{\text{a}}](y_{0})$

$+$ $[\frac{\partial\text{g}_{j\text{a}}}{\partial x_{j}}\frac{\partial^{2}%
}{\partial x_{i}^{2}}(\nabla$log$\Phi_{P})_{\text{a}}](y_{0})$

Therefore we have:

$\qquad Q_{1}=[X_{i}(y_{0})X_{j}(y_{0})-\frac{1}{2}\left(  \frac{\partial
X_{j}}{\partial x_{i}}+\frac{\partial X_{i}}{\partial x_{j}}\right)
](y_{0})[-\frac{1}{2}\left(  \frac{\partial X_{j}}{\partial x_{i}}%
+\frac{\partial X_{i}}{\partial x_{j}}\right)  ](y_{0})\qquad\qquad\left(
74\right)  $

$\qquad-[X_{j}\frac{\partial^{2}}{\partial x_{i}^{2}}(\nabla$log$\Phi_{P}%
)_{j}](y_{0})-\frac{1}{6}\underset{k=q+1}{\overset{n}{%
{\textstyle\sum}
}}R_{jijk}(y_{0})\left(  \frac{\partial X_{k}}{\partial x_{i}}+\frac{\partial
X_{i}}{\partial x_{k}}\right)  (y_{0})$

$\qquad-\frac{4}{3}\underset{\text{a=1}}{\overset{q}{\sum}}R_{j\text{a}%
ij}(y_{0})[\underset{k=q+1}{\overset{n}{%
{\textstyle\sum}
}}X_{k}\perp_{\text{a}ik}-\frac{\partial X_{i}}{\partial x_{\text{a}}}%
](y_{0})$

$-$ $\underset{\text{a=1}}{\overset{q}{\sum}}\perp_{\text{a}ij}X_{j}(y_{0})[$
$\underset{k=q+1}{\overset{n}{%
{\textstyle\sum}
}}X_{k}\perp_{\text{a}ik}-\frac{\partial X_{i}}{\partial x_{\text{a}}}%
](y_{0})$

\qquad\qquad\qquad\qquad\qquad\qquad\qquad\qquad\qquad\qquad\qquad\qquad
\qquad\qquad\qquad\qquad\qquad\qquad$\blacksquare$

$Q_{2}=\frac{\partial}{\partial x_{i}}[\Phi_{P}$g$_{jk}\frac{\partial^{2}%
}{\partial x_{j}\partial x_{i}}(\nabla\log\Phi_{P})_{k}](y_{0})$

$=\frac{\partial}{\partial x_{i}}[\Phi_{P}$g$_{jk}](y_{0})\frac{\partial^{2}%
}{\partial x_{j}\partial x_{i}}(\nabla\log\Phi_{P})_{k}](y_{0})+[\Phi_{P}%
$g$_{jk}](y_{0})\frac{\partial^{3}}{\partial x_{i}^{2}\partial x_{j}}%
[(\nabla\log\Phi_{P})_{k}](y_{0})$

$=[\frac{\partial\Phi_{P}}{\partial x_{i}}$g$_{jk}+\Phi_{P}\frac{\partial
g_{jk}}{\partial x_{i}}](y_{0})\frac{\partial^{2}}{\partial x_{j}\partial
x_{i}}[(\nabla\log\Phi_{P})_{k}](y_{0})+[\Phi_{P}$g$_{jk}](y_{0}%
)\frac{\partial^{3}}{\partial x_{i}^{2}\partial x_{j}}[(\nabla\log\Phi
_{P})_{k}](y_{0})$

Since $\Phi_{P}(y_{0})=1$ and g$_{jk}(y_{0})=\delta_{jk},$ we have:

$Q_{2}=\frac{\partial\Phi_{P}}{\partial x_{i}}(y_{0})\frac{\partial^{2}%
}{\partial x_{j}\partial x_{i}}[(\nabla\log\Phi_{P})_{j}](y_{0})+\frac
{\partial g_{jk}}{\partial x_{i}}(y_{0})\frac{\partial^{2}}{\partial
x_{j}\partial x_{i}}[(\nabla\log\Phi_{P})_{k}](y_{0})$

$\qquad+\frac{\partial^{3}}{\partial x_{i}^{2}\partial x_{j}}[(\nabla\log
\Phi_{P})_{j}](y_{0})$

$Q_{2}=\frac{\partial\Phi_{P}}{\partial x_{i}}(y_{0})\frac{\partial^{2}%
}{\partial x_{j}\partial x_{i}}(\nabla\log\Phi_{P})_{j}](y_{0})+\frac
{\partial^{3}}{\partial x_{i}^{2}\partial x_{j}}(\nabla\log\Phi_{P}%
)_{j}](y_{0})$

$+\underset{\text{a=1}}{\overset{q}{\sum}}\frac{\partial g_{j\text{a}}%
}{\partial x_{i}}(y_{0})\frac{\partial^{2}}{\partial x_{j}\partial x_{i}%
}[(\nabla\log\Phi_{P})_{\text{a}}](y_{0})+\underset{k=q+1}{\overset{n}{%
{\textstyle\sum}
}}\frac{\partial g_{jk}}{\partial x_{i}}(y_{0})\frac{\partial^{2}}{\partial
x_{j}\partial x_{i}}[(\nabla\log\Phi_{P})_{k}](y_{0})$

Since $\frac{\partial g_{jk}}{\partial x_{i}}(y_{0})=0$ for $i,j,k=q+1,...,n,$
we have:

$Q_{2}=\frac{\partial\Phi_{P}}{\partial x_{i}}(y_{0})\frac{\partial^{2}%
}{\partial x_{j}\partial x_{i}}[(\nabla\log\Phi_{P})_{j}](y_{0})+\frac
{\partial^{3}}{\partial x_{i}^{2}\partial x_{j}}[(\nabla\log\Phi_{P}%
)_{j}](y_{0})$

$+\underset{\text{a=1}}{\overset{q}{\sum}}\frac{\partial g_{j\text{a}}%
}{\partial x_{i}}(y_{0})\frac{\partial^{2}}{\partial x_{j}\partial x_{i}%
}[(\nabla\log\Phi_{P})_{\text{a}}](y_{0})$

$Q_{2}=-X_{i}(y_{0})\frac{\partial^{2}}{\partial x_{j}\partial x_{i}}%
[(\nabla\log\Phi_{P})_{j}](y_{0})+\frac{\partial^{3}}{\partial x_{i}%
^{2}\partial x_{j}}[(\nabla\log\Phi_{P})_{j}](y_{0})$

$+\underset{\text{a=1}}{\overset{q}{\sum}}\perp_{\text{a}ij}(y_{0}%
)\frac{\partial^{2}}{\partial x_{j}\partial x_{i}}[(\nabla\log\Phi
_{P})_{\text{a}}](y_{0})\qquad\qquad\qquad\qquad\qquad\qquad\qquad\left(
75\right)  $

From $\left(  74\right)  $ and $\left(  75\right)  ,$ we have:

$\qquad\frac{\partial^{2}L_{3}}{\partial x_{i}\partial x_{j}}(y_{0}%
)=Q_{1}+Q_{2}$

$=[X_{i}(y_{0})X_{j}(y_{0})-\frac{1}{2}\left(  \frac{\partial X_{j}}{\partial
x_{i}}+\frac{\partial X_{i}}{\partial x_{j}}\right)  ](y_{0})[-\frac{1}%
{2}\left(  \frac{\partial X_{j}}{\partial x_{i}}+\frac{\partial X_{i}%
}{\partial x_{j}}\right)  ](y_{0})$ $Q_{1}\qquad\left(  76\right)  $

$-[X_{j}\frac{\partial^{2}}{\partial x_{i}^{2}}(\nabla\log\Phi_{P})_{j}%
](y_{0})-\frac{1}{6}\underset{k=q+1}{\overset{n}{%
{\textstyle\sum}
}}R_{jijk}(y_{0})\left(  \frac{\partial X_{k}}{\partial x_{i}}+\frac{\partial
X_{i}}{\partial x_{k}}\right)  (y_{0})$

$-\frac{4}{3}\underset{\text{a=1}}{\overset{q}{\sum}}R_{j\text{a}ij}%
(y_{0})[\underset{k=q+1}{\overset{n}{%
{\textstyle\sum}
}}X_{k}\perp_{\text{a}ik}-\frac{\partial X_{i}}{\partial x_{\text{a}}}%
](y_{0})$

$-$ $\underset{\text{a=1}}{\overset{q}{\sum}}\perp_{\text{a}ij}X_{j}(y_{0})[$
$\underset{k=q+1}{\overset{n}{%
{\textstyle\sum}
}}X_{k}\perp_{\text{a}ik}-\frac{\partial X_{i}}{\partial x_{\text{a}}}%
](y_{0})$

$-X_{i}(y_{0})\frac{\partial^{2}}{\partial x_{j}\partial x_{i}}[(\nabla
\log\Phi_{P})_{j}](y_{0})+\frac{\partial^{3}}{\partial x_{i}^{2}\partial
x_{j}}[(\nabla\log\Phi_{P})_{j}](y_{0})\qquad\qquad Q_{2}$

$+\underset{\text{a=1}}{\overset{q}{\sum}}\perp_{\text{a}ij}(y_{0}%
)\frac{\partial^{2}}{\partial x_{j}\partial x_{i}}[(\nabla\log\Phi
_{P})_{\text{a}}](y_{0})\qquad$

\qquad\qquad\qquad\qquad\qquad\qquad\qquad\qquad\qquad\qquad\qquad\qquad
\qquad\qquad\qquad$\blacksquare$

By $\left(  63\right)  ,$ $\left(  69\right)  ,$ $\left(  72\right)  $ and
$\left(  76\right)  ,$ we have:\qquad\qquad\qquad\qquad\qquad\qquad\qquad

$\qquad\qquad\frac{\partial^{4}\Phi_{P}}{\partial x_{i}\partial x_{j}\partial
x_{i}\partial x_{j}}(y_{0})=\frac{\partial^{2}L_{1}}{\partial x_{i}\partial
x_{j}}(y_{0})+\frac{\partial^{2}L_{2}}{\partial x_{i}\partial x_{j}}%
(y_{0})+\frac{\partial^{2}L_{3}}{\partial x_{i}\partial x_{j}}(y_{0}%
)\qquad\qquad\qquad\left(  77\right)  $

$=[X_{i}^{2}X_{j}^{2}-X_{i}X_{j}\left(  \frac{\partial X_{j}}{\partial x_{i}%
}+\frac{\partial X_{i}}{\partial x_{j}}\right)  -X_{j}^{2}\frac{\partial
X_{i}}{\partial x_{i}}\ +\frac{1}{3}X_{j}\left(  \frac{\partial^{2}X_{j}%
}{\partial x_{i}^{2}}+2\frac{\partial^{2}X_{i}}{\partial x_{i}\partial x_{j}%
}\right)  ](y_{0})\qquad M_{1}\qquad\frac{\partial^{2}L_{1}}{\partial
x_{i}\partial x_{j}}(y_{0})$

\qquad$\qquad\ -\frac{1}{2}[X_{i}X_{j}\left(  \frac{\partial X_{j}}{\partial
x_{i}}+\frac{\partial X_{i}}{\partial x_{j}}\right)  ](y_{0})+\frac{1}%
{4}\left(  \frac{\partial X_{j}}{\partial x_{i}}+\frac{\partial X_{i}%
}{\partial x_{j}}\right)  ^{2}(y_{0})$

$\qquad\qquad+[-X_{i}^{2}\frac{\partial X_{j}}{\partial x_{j}}+\frac{\partial
X_{i}}{\partial x_{i}}\frac{\partial X_{j}}{\partial x_{j}}](y_{0}%
)\qquad\qquad M_{2}$

$\qquad\qquad+\frac{1}{3}\underset{k=q+1}{\overset{n}{\sum}}[R_{jijk}%
X_{i}X_{k}$ $](y_{0})\qquad M_{3}$

$-X_{i}(y_{0})\underset{\text{a=1}}{\overset{q}{\sum}}[\perp_{\text{a}%
ij}(y_{0})\underset{k=q+1}{\overset{n}{%
{\textstyle\sum}
}}X_{k}\perp_{\text{a}jk}-\frac{\partial X_{j}}{\partial x_{\text{a}}}%
](y_{0})-X_{i}(y_{0})[\frac{\partial^{2}}{\partial x_{i}\partial x_{j}}%
(\nabla\log\Phi_{P})_{j}](y_{0})\qquad M_{4}$

$\qquad\qquad+$ $\underset{\text{a=1}}{\overset{q}{\sum}}%
[\underset{k=q+1}{\overset{n}{%
{\textstyle\sum}
}}-X_{j}X_{k}\perp_{\text{a}ij}\perp_{\text{a}ik}+X_{j}\frac{\partial X_{i}%
}{\partial x_{\text{a}}}\perp_{\text{a}ij}](y_{0})\qquad N_{1}\qquad
\qquad\frac{\partial^{2}L_{2}}{\partial x_{i}\partial x_{j}}(y_{0}%
)\qquad\qquad\qquad\qquad$

$\qquad\qquad+\frac{4}{3}\underset{\text{a=1}}{\overset{q}{\sum}}$
$[\underset{k=q+1}{\overset{n}{%
{\textstyle\sum}
}}-R_{j\text{a}ij}X_{k}\perp_{\text{a}ik}+R_{j\text{a}ij}\frac{\partial X_{i}%
}{\partial x_{\text{a}}}](y_{0})$

$\qquad\qquad+\underset{k=q+1}{\overset{n}{\sum}}[\frac{1}{3}R_{jijk}%
X_{i}X_{k}-\frac{2}{3}R_{ijik}X_{j}X_{k}+\frac{1}{3}(\nabla_{i}R_{ijjk}%
+\nabla_{j}R_{ijik})X_{k}](y_{0})$

$\qquad\qquad-\frac{1}{6}\underset{k=q+1}{\overset{n}{\sum}}$ $[R_{jijk}%
\left(  \frac{\partial X_{k}}{\partial x_{i}}+\frac{\partial X_{i}}{\partial
x_{k}}\right)  ](y_{0})$

$\qquad\qquad+\underset{\text{a=1}}{\overset{q}{\sum}}[-X_{i}\perp
_{\text{a}ij}-\frac{8}{3}R_{i\text{a}ij}](y_{0})[\underset{k=q+1}{\overset{n}{%
{\textstyle\sum}
}}X_{k}\perp_{\text{a}jk}-\frac{\partial X_{j}}{\partial x_{\text{a}}}%
](y_{0})\qquad N_{2}$

$\qquad\qquad+\underset{\text{a=1}}{\overset{q}{\sum}}\perp_{\text{a}ij}%
(y_{0})\frac{\partial^{2}}{\partial x_{i}\partial x_{j}}(\nabla\log\Phi
_{P})_{\text{a}}(y_{0})$

$\qquad\qquad+\frac{1}{3}\underset{k=q+1}{\overset{n}{\sum}}R_{ijik}%
(y_{0})\left(  \frac{\partial X_{k}}{\partial x_{j}}+\frac{\partial X_{j}%
}{\partial x_{k}}\right)  (y_{0})$

$\qquad\qquad+$ $[X_{i}(y_{0})X_{j}(y_{0})-\frac{1}{2}\left(  \frac{\partial
X_{j}}{\partial x_{i}}+\frac{\partial X_{i}}{\partial x_{j}}\right)
](y_{0})[-\frac{1}{2}\left(  \frac{\partial X_{j}}{\partial x_{i}}%
+\frac{\partial X_{i}}{\partial x_{j}}\right)  ](y_{0})\qquad Q_{1}\qquad
\frac{\partial^{2}L_{3}}{\partial x_{i}\partial x_{j}}(y_{0})$

$\qquad\qquad-X_{j}\frac{\partial^{2}}{\partial x_{i}^{2}}[(\nabla\log\Phi
_{P})_{j}](y_{0})-\frac{1}{6}\underset{k=q+1}{\overset{n}{%
{\textstyle\sum}
}}R_{jijk}(y_{0})\left(  \frac{\partial X_{k}}{\partial x_{i}}+\frac{\partial
X_{i}}{\partial x_{k}}\right)  (y_{0})$

$-\frac{4}{3}\underset{\text{a=1}}{\overset{q}{\sum}}R_{j\text{a}ij}%
(y_{0})[\underset{k=q+1}{\overset{n}{%
{\textstyle\sum}
}}X_{k}\perp_{\text{a}ik}-\frac{\partial X_{i}}{\partial x_{\text{a}}}%
](y_{0})-$ $\underset{\text{a=1}}{\overset{q}{\sum}}\perp_{\text{a}ij}%
X_{j}(y_{0})[$ $\underset{k=q+1}{\overset{n}{%
{\textstyle\sum}
}}X_{k}\perp_{\text{a}ik}-\frac{\partial X_{i}}{\partial x_{\text{a}}}%
](y_{0})$

$\qquad-X_{i}(y_{0})\frac{\partial^{2}}{\partial x_{i}\partial x_{j}}%
[(\nabla\log\Phi_{P})_{j}](y_{0})+\frac{\partial^{3}}{\partial x_{i}%
^{2}\partial x_{j}}[(\nabla\log\Phi_{P})_{j}](y_{0})\qquad\qquad Q_{2}$

$+\underset{\text{a=1}}{\overset{q}{\sum}}\perp_{\text{a}ij}(y_{0}%
)\frac{\partial^{2}}{\partial x_{j}\partial x_{i}}[(\nabla\log\Phi
_{P})_{\text{a}}](y_{0})\qquad$

\qquad\qquad\qquad\qquad\qquad\qquad\qquad\qquad\qquad\qquad\qquad\qquad
\qquad\qquad\qquad\qquad\qquad$\blacksquare$

We re-write the last expression above, keeping like-terms together:

$\frac{\partial^{4}\Phi_{P}}{\partial x_{i}\partial x_{j}\partial
x_{i}\partial x_{j}}(y_{0})=\frac{\partial^{2}L_{1}}{\partial x_{i}\partial
x_{j}}(y_{0})+\frac{\partial^{2}L_{2}}{\partial x_{i}\partial x_{j}}%
(y_{0})+\frac{\partial^{2}L_{3}}{\partial x_{i}\partial x_{j}}(y_{0}%
)\qquad\qquad\qquad\qquad\left(  78\right)  $

$=[X_{i}^{2}X_{j}^{2}-X_{i}X_{j}\left(  \frac{\partial X_{j}}{\partial x_{i}%
}+\frac{\partial X_{i}}{\partial x_{j}}\right)  -X_{i}^{2}\frac{\partial
X_{j}}{\partial x_{j}}-X_{j}^{2}\frac{\partial X_{i}}{\partial x_{i}}%
\ +\frac{1}{3}X_{j}\left(  \frac{\partial^{2}X_{j}}{\partial x_{i}^{2}}%
+2\frac{\partial^{2}X_{i}}{\partial x_{i}\partial x_{j}}\right)  ](y_{0})$

\qquad$\qquad\ -\frac{1}{2}[X_{i}X_{j}\left(  \frac{\partial X_{j}}{\partial
x_{i}}+\frac{\partial X_{i}}{\partial x_{j}}\right)  ](y_{0})+\frac{1}%
{4}\left(  \frac{\partial X_{j}}{\partial x_{i}}+\frac{\partial X_{i}%
}{\partial x_{j}}\right)  ^{2}(y_{0})+[\frac{\partial X_{i}}{\partial x_{i}%
}\frac{\partial X_{j}}{\partial x_{j}}](y_{0})\qquad$

$\qquad\qquad+\frac{1}{3}\underset{k=q+1}{\overset{n}{\sum}}[R_{jijk}%
X_{i}X_{k}$ $](y_{0})$

$\qquad\qquad+\underset{k=q+1}{\overset{n}{\sum}}[\frac{1}{3}R_{jijk}%
X_{i}X_{k}-\frac{2}{3}R_{ijik}X_{j}X_{k}+\frac{1}{3}(\nabla_{i}R_{ijjk}%
+\nabla_{j}R_{ijik})X_{k}](y_{0})$

$\qquad\qquad-\frac{1}{6}\underset{k=q+1}{\overset{n}{\sum}}$ $[R_{jijk}%
\left(  \frac{\partial X_{k}}{\partial x_{i}}+\frac{\partial X_{i}}{\partial
x_{k}}\right)  ](y_{0})$

$\qquad\qquad-\frac{1}{6}\underset{k=q+1}{\overset{n}{%
{\textstyle\sum}
}}R_{jijk}(y_{0})\left(  \frac{\partial X_{k}}{\partial x_{i}}+\frac{\partial
X_{i}}{\partial x_{k}}\right)  (y_{0})$

$\qquad\qquad+\frac{1}{3}\underset{k=q+1}{\overset{n}{\sum}}R_{ijik}%
(y_{0})\left(  \frac{\partial X_{k}}{\partial x_{j}}+\frac{\partial X_{j}%
}{\partial x_{k}}\right)  (y_{0})$

$\qquad\qquad+$ $[X_{i}(y_{0})X_{j}(y_{0})-\frac{1}{2}\left(  \frac{\partial
X_{j}}{\partial x_{i}}+\frac{\partial X_{i}}{\partial x_{j}}\right)
](y_{0})[-\frac{1}{2}\left(  \frac{\partial X_{j}}{\partial x_{i}}%
+\frac{\partial X_{i}}{\partial x_{j}}\right)  ](y_{0})$

$\qquad\qquad-X_{i}(y_{0})\frac{\partial^{2}}{\partial x_{i}\partial x_{j}%
}[(\nabla\log\Phi_{P})_{j}](y_{0})-X_{j}\frac{\partial^{2}}{\partial x_{i}%
^{2}}[(\nabla\log\Phi_{P})_{j}](y_{0})$

$\qquad\qquad-X_{i}(y_{0})\frac{\partial^{2}}{\partial x_{i}\partial x_{j}%
}[(\nabla\log\Phi_{P})_{j}](y_{0})+\frac{\partial^{3}}{\partial x_{i}%
^{2}\partial x_{j}}[(\nabla\log\Phi_{P})_{j}](y_{0})$

$\qquad\qquad-X_{i}(y_{0})\underset{\text{a=1}}{\overset{q}{\sum}}%
[\perp_{\text{a}ij}(y_{0})\underset{k=q+1}{\overset{n}{%
{\textstyle\sum}
}}X_{k}\perp_{\text{a}jk}-\frac{\partial X_{j}}{\partial x_{\text{a}}}%
](y_{0})$

$\qquad\qquad+$ $\underset{\text{a=1}}{\overset{q}{\sum}}%
[\underset{k=q+1}{\overset{n}{%
{\textstyle\sum}
}}-X_{j}X_{k}\perp_{\text{a}ij}\perp_{\text{a}ik}+X_{j}\frac{\partial X_{i}%
}{\partial x_{\text{a}}}\perp_{\text{a}ij}](y_{0})\qquad\qquad\qquad$

$\qquad\qquad+\frac{4}{3}\underset{\text{a=1}}{\overset{q}{\sum}}$
$[\underset{k=q+1}{\overset{n}{%
{\textstyle\sum}
}}-R_{j\text{a}ij}X_{k}\perp_{\text{a}ik}+R_{j\text{a}ij}\frac{\partial X_{i}%
}{\partial x_{\text{a}}}](y_{0})$

$\qquad\qquad+\underset{\text{a=1}}{\overset{q}{\sum}}[-X_{i}\perp
_{\text{a}ij}-\frac{8}{3}R_{i\text{a}ij}](y_{0})[\underset{k=q+1}{\overset{n}{%
{\textstyle\sum}
}}X_{k}\perp_{\text{a}jk}-\frac{\partial X_{j}}{\partial x_{\text{a}}}%
](y_{0})$

$\qquad\qquad+\underset{\text{a=1}}{\overset{q}{\sum}}\perp_{\text{a}ij}%
(y_{0})\frac{\partial^{2}}{\partial x_{i}\partial x_{j}}[(\nabla$log$\Phi
_{P})_{\text{a}}](y_{0})+$ $\underset{\text{a=1}}{\overset{q}{\sum}}%
\perp_{\text{a}ij}(y_{0})\frac{\partial^{2}}{\partial x_{j}\partial x_{i}%
}[(\nabla\log\Phi_{P})_{\text{a}}](y_{0})$

$-\frac{4}{3}\underset{\text{a=1}}{\overset{q}{\sum}}R_{j\text{a}ij}%
(y_{0})[\underset{k=q+1}{\overset{n}{%
{\textstyle\sum}
}}X_{k}\perp_{\text{a}ik}-\frac{\partial X_{i}}{\partial x_{\text{a}}}%
](y_{0})-$ $\underset{\text{a=1}}{\overset{q}{\sum}}\perp_{\text{a}ij}%
X_{j}(y_{0})[$ $\underset{k=q+1}{\overset{n}{%
{\textstyle\sum}
}}X_{k}\perp_{\text{a}ik}-\frac{\partial X_{i}}{\partial x_{\text{a}}}%
](y_{0})$

\qquad\qquad\qquad\qquad\qquad\qquad\qquad\qquad\qquad\qquad\qquad\qquad
\qquad\qquad\qquad\qquad$\blacksquare$

We simplify the expression: For the ease of simplifications, we have marked
the same expression with the same number in the main expression for
$\frac{\partial^{4}\Phi_{P}}{\partial x_{i}\partial x_{j}\partial
x_{i}\partial x_{j}}(y_{0})$ above:$\qquad\qquad\qquad\qquad\qquad\qquad$

$\qquad\frac{\partial^{4}\Phi_{P}}{\partial x_{i}\partial x_{j}\partial
x_{i}\partial x_{j}}(y_{0})=[X_{i}^{2}X_{j}^{2}-X_{i}X_{j}\left(
\frac{\partial X_{j}}{\partial x_{i}}+\frac{\partial X_{i}}{\partial x_{j}%
}\right)  -X_{i}^{2}\frac{\partial X_{j}}{\partial x_{j}}-X_{j}^{2}%
\frac{\partial X_{i}}{\partial x_{i}}\ \qquad\left(  79\right)  $

$\qquad+\frac{1}{3}X_{j}\left(  \frac{\partial^{2}X_{j}}{\partial x_{i}^{2}%
}+2\frac{\partial^{2}X_{i}}{\partial x_{i}\partial x_{j}}\right)
](y_{0})\qquad$

\qquad$\ -[X_{i}X_{j}\left(  \frac{\partial X_{j}}{\partial x_{i}}%
+\frac{\partial X_{i}}{\partial x_{j}}\right)  ](y_{0})+\frac{1}{2}\left(
\frac{\partial X_{j}}{\partial x_{i}}+\frac{\partial X_{i}}{\partial x_{j}%
}\right)  ^{2}(y_{0})+[\frac{\partial X_{i}}{\partial x_{i}}\frac{\partial
X_{j}}{\partial x_{j}}](y_{0})\qquad$

$\qquad+\frac{2}{3}\underset{k=q+1}{\overset{n}{\sum}}[R_{jijk}X_{i}%
X_{k}-R_{ijik}X_{j}X_{k}+\frac{1}{2}(\nabla_{i}R_{ijjk}+\nabla_{j}%
R_{ijik})X_{k}](y_{0})$

$\qquad+\frac{1}{3}\underset{k=q+1}{\overset{n}{\sum}}R_{ijik}(y_{0})\left(
\frac{\partial X_{k}}{\partial x_{j}}+\frac{\partial X_{j}}{\partial x_{k}%
}\right)  (y_{0})-\frac{1}{3}\underset{k=q+1}{\overset{n}{\sum}}$
$R_{jijk}(y_{0})\left(  \frac{\partial X_{k}}{\partial x_{i}}+\frac{\partial
X_{i}}{\partial x_{k}}\right)  (y_{0})$

$-2X_{i}(y_{0})\frac{\partial^{2}}{\partial x_{i}\partial x_{j}}[(\nabla
\log\Phi_{P})_{j}](y_{0})-X_{j}\frac{\partial^{2}}{\partial x_{i}^{2}}%
[(\nabla\log\Phi_{P})_{j}](y_{0})+\frac{\partial^{3}}{\partial x_{i}%
^{2}\partial x_{j}}[(\nabla\log\Phi_{P})_{j}](y_{0})$

$\qquad+2\underset{\text{a=1}}{\overset{q}{\sum}}\perp_{\text{a}ij}%
(y_{0})\frac{\partial^{2}}{\partial x_{i}\partial x_{j}}[(\nabla\log\Phi
_{P})_{\text{a}}](y_{0})$

$\qquad-\underset{\text{a=1}}{\overset{q}{\sum}}\underset{k=q+1}{[\overset{n}{%
{\textstyle\sum}
}}\perp_{\text{a}ij}\perp_{\text{a}jk}X_{i}X_{k}-\perp_{\text{a}ij}X_{i}%
\frac{\partial X_{j}}{\partial x_{\text{a}}}](y_{0})\qquad(1)$

$\qquad-$ $\underset{\text{a=1}}{\overset{q}{\sum}}%
[\underset{k=q+1}{\overset{n}{%
{\textstyle\sum}
}}\perp_{\text{a}ij}\perp_{\text{a}ik}X_{j}X_{k}-\perp_{\text{a}ij}X_{j}%
\frac{\partial X_{i}}{\partial x_{\text{a}}}](y_{0})\qquad(2)\qquad\qquad$

$\qquad-\frac{4}{3}\underset{\text{a=1}}{\overset{q}{\sum}}$
$[\underset{k=q+1}{\overset{n}{%
{\textstyle\sum}
}}\perp_{\text{a}ik}R_{j\text{a}ij}X_{k}-R_{j\text{a}ij}\frac{\partial X_{i}%
}{\partial x_{\text{a}}}](y_{0})\qquad(3)$

$\qquad-\underset{\text{a=1}}{\overset{q}{\sum}}[\underset{k=q+1}{\overset{n}{%
{\textstyle\sum}
}}\perp_{\text{a}ij}\perp_{\text{a}jk}X_{i}X_{k}-\perp_{\text{a}ij}X_{i}%
\frac{\partial X_{j}}{\partial x_{\text{a}}}](y_{0})\qquad(1)$

$\qquad-\frac{8}{3}\underset{\text{a=1}}{\overset{q}{\sum}}%
\underset{k=q+1}{[\overset{n}{%
{\textstyle\sum}
}}\perp_{\text{a}jk}R_{i\text{a}ij}X_{k}-R_{i\text{a}ij}\frac{\partial X_{j}%
}{\partial x_{\text{a}}}](y_{0})\qquad(4)$

$\qquad-\frac{4}{3}\underset{\text{a=1}}{\overset{q}{\sum}}%
[\underset{k=q+1}{\overset{n}{%
{\textstyle\sum}
}}\perp_{\text{a}ik}R_{j\text{a}ij}X_{k}-R_{j\text{a}ij}(y_{0})\frac{\partial
X_{i}}{\partial x_{\text{a}}}](y_{0})\qquad(3)$

$\qquad-$ $\underset{\text{a=1}}{\overset{q}{\sum}}[$
$\underset{k=q+1}{\overset{n}{%
{\textstyle\sum}
}}\perp_{\text{a}ij}\perp_{\text{a}ik}X_{j}X_{k}-\perp_{\text{a}ij}X_{j}%
\frac{\partial X_{i}}{\partial x_{\text{a}}}](y_{0})\qquad(2)$

\qquad\qquad\qquad\qquad\qquad\qquad\qquad\qquad\qquad\qquad\qquad
$\qquad\qquad\qquad\blacksquare$

\ Adding the marked items, the expresssion simplifies to:

$\frac{\partial^{4}\Phi_{P}}{\partial x_{i}\partial x_{j}\partial
x_{i}\partial x_{j}}(y_{0})=[X_{i}^{2}X_{j}^{2}-X_{i}X_{j}\left(
\frac{\partial X_{j}}{\partial x_{i}}+\frac{\partial X_{i}}{\partial x_{j}%
}\right)  -X_{i}^{2}\frac{\partial X_{j}}{\partial x_{j}}-X_{j}^{2}%
\frac{\partial X_{i}}{\partial x_{i}}\ \qquad\qquad\left(  80\right)  $

$\qquad\qquad\qquad+\frac{1}{3}X_{j}\left(  \frac{\partial^{2}X_{j}}{\partial
x_{i}^{2}}+2\frac{\partial^{2}X_{i}}{\partial x_{i}\partial x_{j}}\right)
](y_{0})\qquad$

\qquad$\ -[X_{i}X_{j}\left(  \frac{\partial X_{j}}{\partial x_{i}}%
+\frac{\partial X_{i}}{\partial x_{j}}\right)  ](y_{0})+\frac{1}{2}\left(
\frac{\partial X_{j}}{\partial x_{i}}+\frac{\partial X_{i}}{\partial x_{j}%
}\right)  ^{2}(y_{0})+[\frac{\partial X_{i}}{\partial x_{i}}\frac{\partial
X_{j}}{\partial x_{j}}](y_{0})\qquad$

$\qquad+\frac{2}{3}\underset{k=q+1}{\overset{n}{\sum}}[R_{jijk}X_{i}%
X_{k}-R_{ijik}X_{j}X_{k}+\frac{1}{2}(\nabla_{i}R_{ijjk}+\nabla_{j}%
R_{ijik})X_{k}](y_{0})$

$\qquad+\frac{1}{3}\underset{k=q+1}{\overset{n}{\sum}}R_{ijik}(y_{0})\left(
\frac{\partial X_{k}}{\partial x_{j}}+\frac{\partial X_{j}}{\partial x_{k}%
}\right)  (y_{0})-\frac{1}{3}\underset{k=q+1}{\overset{n}{\sum}}$
$R_{jijk}(y_{0})\left(  \frac{\partial X_{k}}{\partial x_{i}}+\frac{\partial
X_{i}}{\partial x_{k}}\right)  (y_{0})$

$\qquad-2X_{i}(y_{0})\frac{\partial^{2}}{\partial x_{i}\partial x_{j}}%
[(\nabla\log\Phi_{P})_{j}](y_{0})-X_{j}(y_{0})\frac{\partial^{2}}{\partial
x_{i}^{2}}[(\nabla\log\Phi_{P})_{j}](y_{0})$

$\qquad+\frac{\partial^{3}}{\partial x_{i}^{2}\partial x_{j}}[(\nabla\log
\Phi_{P})_{j}](y_{0})$

$\qquad+2\underset{\text{a=1}}{\overset{q}{\sum}}\perp_{\text{a}ij}%
(y_{0})\frac{\partial^{2}}{\partial x_{i}\partial x_{j}}[(\nabla\log\Phi
_{P})_{\text{a}}](y_{0})$

$\qquad-2\underset{\text{a=1}}{\overset{q}{\sum}}%
\underset{k=q+1}{[\overset{n}{%
{\textstyle\sum}
}}\perp_{\text{a}ij}\perp_{\text{a}jk}X_{i}X_{k}-\perp_{\text{a}ij}X_{i}%
\frac{\partial X_{j}}{\partial x_{\text{a}}}](y_{0})\qquad(1)$

$\qquad-2$ $\underset{\text{a=1}}{\overset{q}{\sum}}%
[\underset{k=q+1}{\overset{n}{%
{\textstyle\sum}
}}\perp_{\text{a}ij}\perp_{\text{a}ik}X_{j}X_{k}-\perp_{\text{a}ij}X_{j}%
\frac{\partial X_{i}}{\partial x_{\text{a}}}](y_{0})\qquad(2)\qquad\qquad$

$\qquad-\frac{8}{3}\underset{\text{a=1}}{\overset{q}{\sum}}$
$[\underset{k=q+1}{\overset{n}{%
{\textstyle\sum}
}}\perp_{\text{a}ik}R_{j\text{a}ij}X_{k}-R_{j\text{a}ij}\frac{\partial X_{i}%
}{\partial x_{\text{a}}}](y_{0})\qquad(3)$

$\qquad-\frac{8}{3}\underset{\text{a=1}}{\overset{q}{\sum}}%
\underset{k=q+1}{[\overset{n}{%
{\textstyle\sum}
}}\perp_{\text{a}jk}R_{i\text{a}ij}X_{k}-R_{i\text{a}ij}\frac{\partial X_{j}%
}{\partial x_{\text{a}}}](y_{0})\qquad(4)$

\qquad\qquad\qquad\qquad\qquad\qquad\qquad\qquad\qquad\qquad\qquad\qquad
\qquad\qquad\qquad\qquad\qquad\qquad$\blacksquare$

From (v)$^{\ast\ast\ast}$ of \textbf{Table B}$_{4}$ in \textbf{Appendix B} or
$\left(  59\right)  $ above, we have:

$\frac{\partial^{2}}{\partial x_{i}\partial x_{j}}[(\nabla\log\Phi_{P}%
)_{j}](y_{0})=-\frac{1}{3}\left(  2\frac{\partial^{2}X_{j}}{\partial
x_{i}\partial x_{j}}+\frac{\partial^{2}X_{i}}{\partial x_{j}^{2}}\right)
(y_{0})-\frac{1}{3}\underset{k=q+1}{\overset{n}{\sum}}R_{ijjk}(y_{0}%
)X_{k}(y_{0})$

$+\underset{\text{a=1}}{\overset{q}{\sum}}[\perp_{\text{a}ij}\frac{\partial
X_{j}}{\partial x_{\text{a}}}](y_{0})+\underset{k=q+1}{\overset{n}{\sum}%
}\underset{\text{a=1}}{\overset{q}{\sum}}[\perp_{\text{a}ij}\perp_{\text{a}%
jk}X_{k}](y_{0})$

$[\frac{\partial^{2}}{\partial x_{i}^{2}}(\nabla\log\Phi_{P})_{j}%
](y_{0})=-\frac{1}{3}\left(  \frac{\partial^{2}X_{j}}{\partial x_{i}^{2}%
}+2\frac{\partial^{2}X_{i}}{\partial x_{i}\partial x_{j}}\right)
(y_{0})-\frac{2}{3}\underset{k=q+1}{\overset{n}{\sum}}R_{ijik}(y_{0}%
)X_{k}(y_{0})\qquad$

$+\underset{\text{a=1}}{\overset{q}{\sum}}[2\perp_{\text{a}ij}\frac{\partial
X_{i}}{\partial x_{\text{a}}}](y_{0})+\underset{k=q+1}{\overset{n}{\sum}%
}\underset{\text{a=1}}{\overset{q}{\sum}}[2\perp_{\text{a}ij}\perp
_{\text{a}ik}X_{k}](y_{0})$

We insert the expressions of $\frac{\partial^{2}}{\partial x_{i}\partial
x_{j}}[(\nabla\log\Phi_{P})_{j}](y_{0})$ and $\frac{\partial^{2}}{\partial
x_{i}^{2}}[(\nabla\log\Phi_{P})_{j}](y_{0})$

into the expression of $\qquad\frac{\partial^{4}\Phi_{P}}{\partial
x_{i}\partial x_{j}\partial x_{i}\partial x_{j}}(y_{0})$ in $\left(
80\right)  $ above and have:

$\frac{\partial^{4}\Phi_{P}}{\partial x_{i}\partial x_{j}\partial
x_{i}\partial x_{j}}(y_{0})=[X_{i}^{2}X_{j}^{2}-2X_{i}X_{j}\left(
\frac{\partial X_{j}}{\partial x_{i}}+\frac{\partial X_{i}}{\partial x_{j}%
}\right)  -X_{i}^{2}\frac{\partial X_{j}}{\partial x_{j}}-X_{j}^{2}%
\frac{\partial X_{i}}{\partial x_{i}}\ \qquad\left(  81\right)  $

$+\frac{1}{3}X_{j}\left(  \frac{\partial^{2}X_{j}}{\partial x_{i}^{2}}%
+2\frac{\partial^{2}X_{i}}{\partial x_{i}\partial x_{j}}\right)
](y_{0})\ +\frac{1}{2}\left(  \frac{\partial X_{j}}{\partial x_{i}}%
+\frac{\partial X_{i}}{\partial x_{j}}\right)  ^{2}(y_{0})+\left(
\frac{\partial X_{i}}{\partial x_{i}}\frac{\partial X_{j}}{\partial x_{j}%
}\right)  (y_{0})\qquad$

$(3)\qquad+\frac{2}{3}\underset{k=q+1}{\overset{n}{\sum}}[R_{jijk}X_{i}%
X_{k}](y_{0})-\frac{2}{3}\underset{k=q+1}{\overset{n}{\sum}}[R_{ijik}%
X_{j}X_{k}](y_{0})\qquad\qquad\left(  4\right)  $

$\qquad+\frac{1}{3}\underset{k=q+1}{\overset{n}{\sum}}(\nabla_{i}%
R_{ijjk}+\nabla_{j}R_{ijik})X_{k}](y_{0})$

$\qquad+\frac{1}{3}\underset{k=q+1}{\overset{n}{\sum}}R_{ijik}(y_{0})\left(
\frac{\partial X_{k}}{\partial x_{j}}+\frac{\partial X_{j}}{\partial x_{k}%
}\right)  (y_{0})-\frac{1}{3}\underset{k=q+1}{\overset{n}{\sum}}$
$R_{jijk}(y_{0})\left(  \frac{\partial X_{k}}{\partial x_{i}}+\frac{\partial
X_{i}}{\partial x_{k}}\right)  (y_{0})$

$\qquad+\frac{2}{3}X_{i}(y_{0})\left(  2\frac{\partial^{2}X_{j}}{\partial
x_{i}\partial x_{j}}+\frac{\partial^{2}X_{i}}{\partial x_{j}^{2}}\right)
(y_{0})-\frac{2}{3}\underset{k=q+1}{\overset{n}{\sum}}[R_{jijk}X_{i}%
X_{k}](y_{0})\qquad$ $(3)$

$\qquad-2\underset{\text{a=1}}{\overset{q}{\sum}}[\perp_{\text{a}ij}X_{i}%
\frac{\partial X_{j}}{\partial x_{\text{a}}}](y_{0}%
)-2\underset{k=q+1}{\overset{n}{\sum}}\underset{\text{a=1}}{\overset{q}{\sum}%
}[\perp_{\text{a}ij}\perp_{\text{a}jk}X_{i}X_{k}](y_{0})\qquad(1)$

$\qquad+\frac{1}{3}X_{j}(y_{0})\left(  \frac{\partial^{2}X_{j}}{\partial
x_{i}^{2}}+2\frac{\partial^{2}X_{i}}{\partial x_{i}\partial x_{j}}\right)
(y_{0})+\frac{2}{3}\underset{k=q+1}{\overset{n}{\sum}}[R_{ijik}X_{j}%
X_{k}](y_{0})\qquad\left(  4\right)  $

$\qquad+2\underset{\text{a=1}}{\overset{q}{\sum}}\perp_{\text{a}ij}%
(y_{0})\frac{\partial^{2}}{\partial x_{i}\partial x_{j}}[(\nabla\log\Phi
_{P})_{\text{a}}](y_{0})+\frac{\partial^{3}}{\partial x_{i}^{2}\partial x_{j}%
}[(\nabla\log\Phi_{P})_{j}](y_{0})\qquad$

$\qquad-2\underset{\text{a=1}}{\overset{q}{\sum}}[\perp_{\text{a}ij}X_{j}%
\frac{\partial X_{i}}{\partial x_{\text{a}}}](y_{0}%
)-2\underset{k=q+1}{\overset{n}{\sum}}\underset{\text{a=1}}{\overset{q}{\sum}%
}[\perp_{\text{a}ij}\perp_{\text{a}ik}X_{j}X_{k}](y_{0})\qquad(2)$

$\qquad-2\underset{\text{a=1}}{\overset{q}{\sum}}%
\underset{k=q+1}{[\overset{n}{%
{\textstyle\sum}
}}\perp_{\text{a}ij}\perp_{\text{a}jk}X_{i}X_{k}](y_{0})+2\underset{\text{a=1}%
}{\overset{q}{\sum}}[\perp_{\text{a}ij}X_{i}\frac{\partial X_{j}}{\partial
x_{\text{a}}}](y_{0})\qquad(1)$

$\qquad-2$ $\underset{\text{a=1}}{\overset{q}{\sum}}%
[\underset{k=q+1}{\overset{n}{%
{\textstyle\sum}
}}\perp_{\text{a}ij}\perp_{\text{a}ik}X_{j}X_{k}-\perp_{\text{a}ij}X_{j}%
\frac{\partial X_{i}}{\partial x_{\text{a}}}](y_{0})\qquad(2)\qquad\qquad$

$\qquad-\frac{8}{3}\underset{\text{a=1}}{\overset{q}{\sum}}$
$[\underset{k=q+1}{\overset{n}{%
{\textstyle\sum}
}}\perp_{\text{a}ik}R_{j\text{a}ij}X_{k}-R_{j\text{a}ij}\frac{\partial X_{i}%
}{\partial x_{\text{a}}}](y_{0})\qquad$

$-\frac{8}{3}\underset{\text{a=1}}{\overset{q}{\sum}}%
\underset{k=q+1}{[\overset{n}{%
{\textstyle\sum}
}}\perp_{\text{a}jk}R_{i\text{a}ij}X_{k}-R_{i\text{a}ij}\frac{\partial X_{j}%
}{\partial x_{\text{a}}}](y_{0})$

$\qquad\qquad\qquad\qquad\qquad\qquad\qquad\qquad\qquad\qquad\qquad
\qquad\qquad\qquad\qquad\qquad\qquad\qquad\blacksquare\qquad$

In the expression above, we have marked similar items with the same number.

Each pair either adds up to zero or to a simpler expression:

$\left(  1\right)  \qquad-2\underset{\text{a=1}}{\overset{q}{\sum}}%
[\perp_{\text{a}ij}X_{i}\frac{\partial X_{j}}{\partial x_{\text{a}}}%
](y_{0})-2\underset{k=q+1}{\overset{n}{\sum}}\underset{\text{a=1}%
}{\overset{q}{\sum}}[\perp_{\text{a}ij}\perp_{\text{a}jk}X_{i}X_{k}](y_{0})$

$\qquad-2\underset{\text{a=1}}{\overset{q}{\sum}}%
\underset{k=q+1}{[\overset{n}{%
{\textstyle\sum}
}}\perp_{\text{a}ij}\perp_{\text{a}jk}X_{i}X_{k}](y_{0})+2\underset{\text{a=1}%
}{\overset{q}{\sum}}[\perp_{\text{a}ij}X_{i}\frac{\partial X_{j}}{\partial
x_{\text{a}}}](y_{0})$

$\qquad=-4\underset{k=q+1}{\overset{n}{\sum}}\underset{\text{a=1}%
}{\overset{q}{\sum}}[\perp_{\text{a}ij}\perp_{\text{a}jk}X_{i}X_{k}](y_{0})$

$\left(  2\right)  \qquad-2\underset{\text{a=1}}{\overset{q}{\sum}}%
[\perp_{\text{a}ij}X_{j}\frac{\partial X_{i}}{\partial x_{\text{a}}}%
](y_{0})+\underset{k=q+1}{\overset{n}{\sum}}\underset{\text{a=1}%
}{\overset{q}{\sum}}\perp_{\text{a}ij}\perp_{\text{a}ik}X_{j}X_{k}](y_{0})$

\qquad\qquad$-2$ $\underset{\text{a=1}}{\overset{q}{\sum}}%
[\underset{k=q+1}{\overset{n}{%
{\textstyle\sum}
}}\perp_{\text{a}ij}\perp_{\text{a}ik}X_{j}X_{k}-\perp_{\text{a}ij}X_{j}%
\frac{\partial X_{i}}{\partial x_{\text{a}}}](y_{0})$

$\qquad=-4\underset{k=q+1}{\overset{n}{\sum}}\underset{\text{a=1}%
}{\overset{q}{\sum}}[\perp_{\text{a}ij}\perp_{\text{a}ik}X_{j}X_{k}](y_{0})$

$\left(  3\right)  \qquad+\frac{2}{3}\underset{k=q+1}{\overset{n}{\sum}%
}[R_{jijk}X_{i}X_{k}](y_{0})-\frac{2}{3}\underset{k=q+1}{\overset{n}{\sum}%
}[R_{jijk}X_{i}X_{k}](y_{0})=0\qquad$

$(4)\qquad-\frac{2}{3}\underset{k=q+1}{\overset{n}{\sum}}[R_{ijik}X_{j}%
X_{k}](y_{0})+\frac{2}{3}\underset{k=q+1}{\overset{n}{\sum}}[R_{ijik}%
X_{j}X_{k}](y_{0})=0$

Therefore the expression simplifies to:

$\qquad\frac{\partial^{4}\Phi_{P}}{\partial x_{i}\partial x_{j}\partial
x_{i}\partial x_{j}}(y_{0})=[X_{i}^{2}X_{j}^{2}-2X_{i}X_{j}\left(
\frac{\partial X_{j}}{\partial x_{i}}+\frac{\partial X_{i}}{\partial x_{j}%
}\right)  -X_{i}^{2}\frac{\partial X_{j}}{\partial x_{j}}-X_{j}^{2}%
\frac{\partial X_{i}}{\partial x_{i}}\ ](y_{0})\qquad\left(  82\right)  $

$\qquad+\frac{2}{3}X_{i}(y_{0})\left(  2\frac{\partial^{2}X_{j}}{\partial
x_{i}\partial x_{j}}+\frac{\partial^{2}X_{i}}{\partial x_{j}^{2}}\right)
(y_{0})+\frac{2}{3}X_{j}(y_{0})\left(  \frac{\partial^{2}X_{j}}{\partial
x_{i}^{2}}+2\frac{\partial^{2}X_{i}}{\partial x_{i}\partial x_{j}}\right)
(y_{0})$

\qquad$\ +\frac{1}{2}\left(  \frac{\partial X_{j}}{\partial x_{i}}%
+\frac{\partial X_{i}}{\partial x_{j}}\right)  ^{2}(y_{0})+\left(
\frac{\partial X_{i}}{\partial x_{i}}\frac{\partial X_{j}}{\partial x_{j}%
}\right)  (y_{0})+\frac{1}{3}\underset{k=q+1}{\overset{n}{\sum}}(\nabla
_{i}R_{ijjk}+\nabla_{j}R_{ijik})X_{k}](y_{0})$

$\qquad+\frac{1}{3}\underset{k=q+1}{\overset{n}{\sum}}R_{ijik}(y_{0})\left(
\frac{\partial X_{k}}{\partial x_{j}}+\frac{\partial X_{j}}{\partial x_{k}%
}\right)  (y_{0})-\frac{1}{3}\underset{k=q+1}{\overset{n}{\sum}}$
$R_{jijk}(y_{0})\left(  \frac{\partial X_{k}}{\partial x_{i}}+\frac{\partial
X_{i}}{\partial x_{k}}\right)  (y_{0})$

$\qquad+\frac{\partial^{3}}{\partial x_{i}^{2}\partial x_{j}}[(\nabla\log
\Phi_{P})_{j}](y_{0})\qquad$

$\qquad+2\underset{\text{a=1}}{\overset{q}{\sum}}\perp_{\text{a}ij}%
(y_{0})\frac{\partial^{2}}{\partial x_{i}\partial x_{j}}[(\nabla\log\Phi
_{P})_{\text{a}}](y_{0})$

$\qquad-4\underset{k=q+1}{\overset{n}{\sum}}\underset{\text{a=1}%
}{\overset{q}{\sum}}[\perp_{\text{a}ij}\perp_{\text{a}jk}X_{i}X_{k}%
](y_{0})\qquad\left(  1\right)  $

$\qquad-4\underset{k=q+1}{\overset{n}{\sum}}\underset{\text{a=1}%
}{\overset{q}{\sum}}[\perp_{\text{a}ij}\perp_{\text{a}ik}X_{j}X_{k}%
](y_{0})\qquad\left(  2\right)  \qquad$

$\qquad-\frac{8}{3}\underset{\text{a=1}}{\overset{q}{\sum}}$
$[\underset{k=q+1}{\overset{n}{%
{\textstyle\sum}
}}\perp_{\text{a}ik}R_{j\text{a}ij}X_{k}-R_{j\text{a}ij}\frac{\partial X_{i}%
}{\partial x_{\text{a}}}](y_{0})\qquad$

$\qquad-\frac{8}{3}\underset{\text{a=1}}{\overset{q}{\sum}}%
\underset{k=q+1}{[\overset{n}{%
{\textstyle\sum}
}}\perp_{\text{a}jk}R_{i\text{a}ij}X_{k}-R_{i\text{a}ij}\frac{\partial X_{j}%
}{\partial x_{\text{a}}}](y_{0})\qquad$

\qquad\qquad\qquad\qquad\qquad\qquad\qquad\qquad\qquad\qquad\qquad\qquad
\qquad\qquad\qquad$\blacksquare$

We have not yet computed the expression for $\frac{\partial^{3}}{\partial
x_{i}\partial x_{j}^{2}}[(\nabla\log\Phi_{P})_{i}](y_{0}).$ We use a simple
trick to do so:

\ Since the function $\Phi:M\longrightarrow R$ is smooth (we only need it
differentiable to the order four), we can switch the order of differentiation
as we want. In this case it means switching the positions of the indices $i$
and $j$. Therefore,

$\frac{\partial^{4}\Phi_{P}}{\partial x_{i}^{2}\partial x_{j}^{2}}%
(y_{0})=\frac{\partial^{4}\Phi_{P}}{\partial x_{i}\partial x_{j}\partial
x_{i}\partial x_{j}}(y_{0})=\frac{\partial^{4}\Phi_{P}}{\partial x_{j}\partial
x_{i}\partial x_{j}\partial x_{i}}(y_{0})=\frac{\partial^{4}\Phi_{P}}{\partial
x_{j}^{2}\partial x_{i}^{2}}(y_{0})$

We obtain $\frac{\partial^{4}\Phi_{P}}{\partial x_{j}\partial x_{i}\partial
x_{j}\partial x_{i}}(y_{0})$ by switching the positions of the indices $i$ and
$j$ in $\frac{\partial^{4}\Phi_{P}}{\partial x_{i}\partial x_{j}\partial
x_{i}\partial x_{j}}(y_{0})$ above.

However, the first three lines remain unchanged because of the very symmetric
roles of $i$ and $j$ in them.

We thus have:

$\qquad\frac{\partial^{4}\Phi_{P}}{\partial x_{i}^{2}\partial x_{j}^{2}}%
(y_{0})=\frac{\partial^{4}\Phi_{P}}{\partial x_{i}\partial x_{j}\partial
x_{i}\partial x_{j}}(y_{0})=\frac{1}{2}[\frac{\partial^{4}\Phi_{P}}{\partial
x_{i}\partial x_{j}\partial x_{i}\partial x_{j}}+\frac{\partial^{4}\Phi_{P}%
}{\partial x_{j}\partial x_{i}\partial x_{j}\partial x_{i}}](y_{0})$

$\qquad=\frac{\partial^{4}\Phi_{P}}{\partial x_{j}\partial x_{i}\partial
x_{j}\partial x_{i}}(y_{0})=\frac{\partial^{4}\Phi_{P}}{\partial x_{j}%
^{2}\partial x_{i}^{2}}(y_{0})$

$\qquad=[X_{i}^{2}X_{j}^{2}-2X_{i}X_{j}\left(  \frac{\partial X_{j}}{\partial
x_{i}}+\frac{\partial X_{i}}{\partial x_{j}}\right)  -X_{i}^{2}\frac{\partial
X_{j}}{\partial x_{j}}-X_{j}^{2}\frac{\partial X_{i}}{\partial x_{i}}%
\ ](y_{0})\qquad\qquad\qquad\qquad\left(  83\right)  $

$\qquad+\frac{2}{3}X_{i}(y_{0})\left(  2\frac{\partial^{2}X_{j}}{\partial
x_{i}\partial x_{j}}+\frac{\partial^{2}X_{i}}{\partial x_{j}^{2}}\right)
(y_{0})+\frac{2}{3}X_{j}(y_{0})\left(  \frac{\partial^{2}X_{j}}{\partial
x_{i}^{2}}+2\frac{\partial^{2}X_{i}}{\partial x_{i}\partial x_{j}}\right)
(y_{0})$

\qquad$\ +\frac{1}{2}\left(  \frac{\partial X_{j}}{\partial x_{i}}%
+\frac{\partial X_{i}}{\partial x_{j}}\right)  ^{2}(y_{0})+\left(
\frac{\partial X_{i}}{\partial x_{i}}\frac{\partial X_{j}}{\partial x_{j}%
}\right)  (y_{0})\qquad$

$\qquad+\frac{1}{6}\underset{k=q+1}{\overset{n}{\sum}}(\nabla_{i}%
R_{ijjk}+\nabla_{j}R_{ijik})X_{k}](y_{0})+\frac{1}{6}%
\underset{k=q+1}{\overset{n}{\sum}}(\nabla_{j}R_{jiik}+\nabla_{i}%
R_{jijk})X_{k}](y_{0})$

$\qquad+\frac{1}{6}\underset{k=q+1}{\overset{n}{\sum}}R_{ijik}(y_{0})\left(
\frac{\partial X_{k}}{\partial x_{j}}+\frac{\partial X_{j}}{\partial x_{k}%
}\right)  (y_{0})-\frac{1}{6}\underset{k=q+1}{\overset{n}{\sum}}$
$R_{jijk}(y_{0})\left(  \frac{\partial X_{k}}{\partial x_{i}}+\frac{\partial
X_{i}}{\partial x_{k}}\right)  (y_{0})$

$\qquad+\frac{1}{6}\underset{k=q+1}{\overset{n}{\sum}}R_{jijk}(y_{0})\left(
\frac{\partial X_{k}}{\partial x_{i}}+\frac{\partial X_{i}}{\partial x_{k}%
}\right)  (y_{0})-\frac{1}{6}\underset{k=q+1}{\overset{n}{\sum}}$
$R_{ijik}(y_{0})\left(  \frac{\partial X_{k}}{\partial x_{j}}+\frac{\partial
X_{j}}{\partial x_{k}}\right)  (y_{0})$

$\qquad+\frac{1}{2}[\frac{\partial^{3}}{\partial x_{i}^{2}\partial x_{j}%
}(\nabla\log\Phi_{P})_{j}+\frac{\partial^{3}}{\partial x_{j}^{2}\partial
x_{i}}(\nabla\log\Phi_{P})_{i}](y_{0})\qquad$

$\qquad+\underset{\text{a=1}}{\overset{q}{\sum}}\perp_{\text{a}ij}(y_{0}%
)\frac{\partial^{2}}{\partial x_{i}\partial x_{j}}[(\nabla\log\Phi
_{P})_{\text{a}}](y_{0})+$ $\underset{\text{a=1}}{\overset{q}{\sum}}%
\perp_{\text{a}ji}(y_{0})\frac{\partial^{2}}{\partial x_{j}\partial x_{i}%
}[(\nabla\log\Phi_{P})_{\text{a}}](y_{0})\qquad\left(  3\right)  $

$\qquad-2\underset{k=q+1}{\overset{n}{\sum}}\underset{\text{a=1}%
}{\overset{q}{\sum}}[\perp_{\text{a}ij}\perp_{\text{a}jk}X_{i}X_{k}%
](y_{0})+2\underset{k=q+1}{\overset{n}{\sum}}\underset{\text{a=1}%
}{\overset{q}{\sum}}[\perp_{\text{a}ij}\perp_{\text{a}ik}X_{j}X_{k}%
](y_{0})\qquad\left(  1\right)  $

$\qquad-2\underset{k=q+1}{\overset{n}{\sum}}\underset{\text{a=1}%
}{\overset{q}{\sum}}[\perp_{\text{a}ij}\perp_{\text{a}ik}X_{j}X_{k}%
](y_{0})+2\underset{k=q+1}{\overset{n}{\sum}}\underset{\text{a=1}%
}{\overset{q}{\sum}}[\perp_{\text{a}ij}\perp_{\text{a}jk}X_{i}X_{k}%
](y_{0})\qquad\left(  2\right)  \qquad$

$\qquad-\frac{4}{3}\underset{\text{a=1}}{\overset{q}{\sum}}$
$[\underset{k=q+1}{\overset{n}{%
{\textstyle\sum}
}}\perp_{\text{a}ik}R_{j\text{a}ij}X_{k}-R_{j\text{a}ij}\frac{\partial X_{i}%
}{\partial x_{\text{a}}}](y_{0})+\frac{4}{3}\underset{\text{a=1}%
}{\overset{q}{\sum}}[\underset{k=q+1}{\overset{n}{%
{\textstyle\sum}
}}\perp_{\text{a}jk}R_{i\text{a}ij}X_{k}-R_{i\text{a}ij}\frac{\partial X_{j}%
}{\partial x_{\text{a}}}](y_{0})\qquad\left(  4\right)  \qquad$

$\qquad-\frac{4}{3}\underset{\text{a=1}}{\overset{q}{\sum}}%
\underset{k=q+1}{[\overset{n}{%
{\textstyle\sum}
}}\perp_{\text{a}jk}R_{i\text{a}ij}X_{k}-R_{i\text{a}ij}\frac{\partial X_{j}%
}{\partial x_{\text{a}}}](y_{0})+\frac{4}{3}\underset{\text{a=1}%
}{\overset{q}{\sum}}\underset{k=q+1}{[\overset{n}{%
{\textstyle\sum}
}}\perp_{\text{a}ik}R_{j\text{a}ij}X_{k}-R_{j\text{a}ij}\frac{\partial X_{i}%
}{\partial x_{\text{a}}}](y_{0})\qquad\left(  5\right)  $

\qquad\qquad\qquad\qquad\qquad\qquad\qquad\qquad\qquad\qquad\qquad\qquad
\qquad\qquad\qquad\qquad\qquad\qquad$\blacksquare$

We see from the above that,

$\qquad\left(  1\right)  +\left(  2\right)
=-2\underset{k=q+1}{\overset{n}{\sum}}\underset{\text{a=1}}{\overset{q}{\sum}%
}[\perp_{\text{a}ij}\perp_{\text{a}jk}X_{i}X_{k}](y_{0}%
)+2\underset{k=q+1}{\overset{n}{\sum}}\underset{\text{a=1}}{\overset{q}{\sum}%
}[\perp_{\text{a}ij}\perp_{\text{a}ik}X_{j}X_{k}](y_{0})$

$\qquad\qquad\qquad-2\underset{k=q+1}{\overset{n}{\sum}}\underset{\text{a=1}%
}{\overset{q}{\sum}}[\perp_{\text{a}ij}\perp_{\text{a}ik}X_{j}X_{k}%
](y_{0})+2\underset{k=q+1}{\overset{n}{\sum}}\underset{\text{a=1}%
}{\overset{q}{\sum}}[\perp_{\text{a}ij}\perp_{\text{a}jk}X_{i}X_{k}](y_{0})=0$

$\qquad\left(  3\right)  =+\underset{\text{a=1}}{\overset{q}{\sum}}%
\perp_{\text{a}ij}(y_{0})\frac{\partial^{2}}{\partial x_{i}\partial x_{j}%
}[(\nabla\log\Phi_{P})_{\text{a}}](y_{0})+$ $\underset{\text{a=1}%
}{\overset{q}{\sum}}\perp_{\text{a}ji}(y_{0})\frac{\partial^{2}}{\partial
x_{j}\partial x_{i}}[(\nabla\log\Phi_{P})_{\text{a}}](y_{0})=0$

\qquad\qquad The last line vanishes because $\perp_{\text{a}ij}$ is
skew-symmetric in the indics $i,j.$

Next, we have:

$\qquad\left(  4\right)  +\left(  5\right)  $

$=-\frac{4}{3}\underset{\text{a=1}}{\overset{q}{\sum}}$
$[\underset{k=q+1}{\overset{n}{%
{\textstyle\sum}
}}\perp_{\text{a}ik}R_{j\text{a}ij}X_{k}-R_{j\text{a}ij}\frac{\partial X_{i}%
}{\partial x_{\text{a}}}](y_{0})+\frac{4}{3}\underset{\text{a=1}%
}{\overset{q}{\sum}}[\underset{k=q+1}{\overset{n}{%
{\textstyle\sum}
}}\perp_{\text{a}jk}R_{i\text{a}ij}X_{k}-R_{i\text{a}ij}\frac{\partial X_{j}%
}{\partial x_{\text{a}}}](y_{0})\qquad\qquad$

$\qquad-\frac{4}{3}\underset{\text{a=1}}{\overset{q}{\sum}}%
\underset{k=q+1}{[\overset{n}{%
{\textstyle\sum}
}}\perp_{\text{a}jk}R_{i\text{a}ij}X_{k}-R_{i\text{a}ij}\frac{\partial X_{j}%
}{\partial x_{\text{a}}}](y_{0})+\frac{4}{3}\underset{\text{a=1}%
}{\overset{q}{\sum}}\underset{k=q+1}{[\overset{n}{%
{\textstyle\sum}
}}\perp_{\text{a}ik}R_{j\text{a}ij}X_{k}-R_{j\text{a}ij}\frac{\partial X_{i}%
}{\partial x_{\text{a}}}](y_{0})=0$

We see that the expression above simplifies to:

$\qquad\frac{\partial^{4}\Phi_{P}}{\partial x_{i}^{2}\partial x_{j}^{2}}%
(y_{0})=\frac{1}{2}[\frac{\partial^{4}\Phi_{P}}{\partial x_{i}\partial
x_{j}\partial x_{i}\partial x_{j}}+\frac{\partial^{4}\Phi_{P}}{\partial
x_{j}\partial x_{i}\partial x_{j}\partial x_{i}}](y_{0})=\frac{\partial
^{4}\Phi_{P}}{\partial x_{j}^{2}\partial x_{i}^{2}}(y_{0})\qquad\qquad\left(
84\right)  $

$\qquad=[X_{i}^{2}X_{j}^{2}-2X_{i}X_{j}\left(  \frac{\partial X_{j}}{\partial
x_{i}}+\frac{\partial X_{i}}{\partial x_{j}}\right)  -X_{i}^{2}\frac{\partial
X_{j}}{\partial x_{j}}-X_{j}^{2}\frac{\partial X_{i}}{\partial x_{i}}%
](y_{0})\qquad\qquad$

$\qquad+\frac{2}{3}X_{i}(y_{0})\left(  2\frac{\partial^{2}X_{j}}{\partial
x_{i}\partial x_{j}}+\frac{\partial^{2}X_{i}}{\partial x_{j}^{2}}\right)
(y_{0})+\frac{2}{3}X_{j}(y_{0})\left(  \frac{\partial^{2}X_{j}}{\partial
x_{i}^{2}}+2\frac{\partial^{2}X_{i}}{\partial x_{i}\partial x_{j}}\right)
(y_{0}$

$\ +\frac{1}{2}\left(  \frac{\partial X_{j}}{\partial x_{i}}+\frac{\partial
X_{i}}{\partial x_{j}}\right)  ^{2}(y_{0})+\left(  \frac{\partial X_{i}%
}{\partial x_{i}}\frac{\partial X_{j}}{\partial x_{j}}\right)  (y_{0}%
)+\frac{1}{2}[\frac{\partial^{3}}{\partial x_{i}^{2}\partial x_{j}}(\nabla
\log\Phi_{P})_{j}+\frac{\partial^{3}}{\partial x_{j}^{2}\partial x_{i}}%
(\nabla\log\Phi_{P})_{i}](y_{0})\qquad$

$\qquad+\frac{1}{6}\underset{k=q+1}{\overset{n}{\sum}}(\nabla_{i}%
R_{ijjk}+\nabla_{j}R_{ijik})X_{k}](y_{0})+\frac{1}{6}%
\underset{k=q+1}{\overset{n}{\sum}}(\nabla_{j}R_{jiik}+\nabla_{i}%
R_{jijk})X_{k}](y_{0})$

$\qquad+\frac{1}{6}\underset{k=q+1}{\overset{n}{\sum}}R_{ijik}(y_{0})\left(
\frac{\partial X_{k}}{\partial x_{j}}+\frac{\partial X_{j}}{\partial x_{k}%
}\right)  (y_{0})-\frac{1}{6}\underset{k=q+1}{\overset{n}{\sum}}$
$R_{jijk}(y_{0})\left(  \frac{\partial X_{k}}{\partial x_{i}}+\frac{\partial
X_{i}}{\partial x_{k}}\right)  (y_{0})$

$\qquad+\frac{1}{6}\underset{k=q+1}{\overset{n}{\sum}}R_{jijk}(y_{0})\left(
\frac{\partial X_{k}}{\partial x_{i}}+\frac{\partial X_{i}}{\partial x_{k}%
}\right)  (y_{0})-\frac{1}{6}\underset{k=q+1}{\overset{n}{\sum}}$
$R_{ijik}(y_{0})\left(  \frac{\partial X_{k}}{\partial x_{j}}+\frac{\partial
X_{j}}{\partial x_{k}}\right)  (y_{0})$

$\qquad\qquad\qquad\qquad\qquad\qquad\qquad\qquad\qquad\qquad\qquad
\qquad\qquad\qquad\qquad\qquad\qquad\qquad\qquad\blacksquare\qquad$

An important remark here is that we have made use of the skew-symmetry of the
"torsion" $\perp_{\text{a}ij}$in the pair of indices $\left(  i,j\right)  .$
We next have the obvious cancellations:

$+\frac{1}{6}\underset{k=q+1}{\overset{n}{\sum}}R_{ijik}(y_{0})\left(
\frac{\partial X_{k}}{\partial x_{j}}+\frac{\partial X_{j}}{\partial x_{k}%
}\right)  (y_{0})-\frac{1}{6}\underset{k=q+1}{\overset{n}{\sum}}$
$R_{jijk}(y_{0})\left(  \frac{\partial X_{k}}{\partial x_{i}}+\frac{\partial
X_{i}}{\partial x_{k}}\right)  (y_{0})$

$+\frac{1}{6}\underset{k=q+1}{\overset{n}{\sum}}R_{jijk}(y_{0})\left(
\frac{\partial X_{k}}{\partial x_{i}}+\frac{\partial X_{i}}{\partial x_{k}%
}\right)  (y_{0})-\frac{1}{6}\underset{k=q+1}{\overset{n}{\sum}}$
$R_{ijik}(y_{0})\left(  \frac{\partial X_{k}}{\partial x_{j}}+\frac{\partial
X_{j}}{\partial x_{k}}\right)  (y_{0})=0$

We now make use of the skew-symmetry of the Riemannian curvature tensor in its
first two and last two indices as needed:

$+\frac{1}{6}\underset{k=q+1}{\overset{n}{\sum}}[(\nabla_{i}R_{ijjk}%
+\nabla_{j}R_{ijik})X_{k}](y_{0})+\frac{1}{6}\underset{k=q+1}{\overset{n}{\sum
}}[(\nabla_{j}R_{jiik}+\nabla_{i}R_{jijk})X_{k}](y_{0})$

$=+\frac{1}{6}\underset{k=q+1}{\overset{n}{\sum}}[(-\nabla_{i}R_{jijk}%
+\nabla_{j}R_{ijik})X_{k}](y_{0})+\frac{1}{6}\underset{k=q+1}{\overset{n}{\sum
}}[(-\nabla_{j}R_{ijik}+\nabla_{i}R_{jijk})X_{k}](y_{0})=0$

Consequently we have:

$\frac{\partial^{4}\Phi_{P}}{\partial x_{i}^{2}\partial x_{j}^{2}}%
(y_{0})=\frac{1}{2}[\frac{\partial^{4}\Phi_{P}}{\partial x_{i}\partial
x_{j}\partial x_{i}\partial x_{j}}(y_{0})+\frac{\partial^{4}\Phi_{P}}{\partial
x_{j}\partial x_{i}\partial x_{j}\partial x_{i}}(y_{0})=\frac{\partial^{4}%
\Phi_{P}}{\partial x_{j}^{2}\partial x_{i}^{2}}(y_{0})\qquad\left(  85\right)
$

$=[X_{i}^{2}X_{j}^{2}-2X_{i}X_{j}\left(  \frac{\partial X_{j}}{\partial x_{i}%
}+\frac{\partial X_{i}}{\partial x_{j}}\right)  -X_{i}^{2}\frac{\partial
X_{j}}{\partial x_{j}}-X_{j}^{2}\frac{\partial X_{i}}{\partial x_{i}}%
](y_{0})\qquad\qquad$

$+\frac{2}{3}X_{i}(y_{0})\left(  2\frac{\partial^{2}X_{j}}{\partial
x_{i}\partial x_{j}}+\frac{\partial^{2}X_{i}}{\partial x_{j}^{2}}\right)
(y_{0})+\frac{2}{3}X_{j}(y_{0})\left(  \frac{\partial^{2}X_{j}}{\partial
x_{i}^{2}}+2\frac{\partial^{2}X_{i}}{\partial x_{i}\partial x_{j}}\right)
(y_{0}$

\qquad$\ +\frac{1}{2}\left(  \frac{\partial X_{j}}{\partial x_{i}}%
+\frac{\partial X_{i}}{\partial x_{j}}\right)  ^{2}(y_{0})+\left(
\frac{\partial X_{i}}{\partial x_{i}}\frac{\partial X_{j}}{\partial x_{j}%
}\right)  (y_{0})$

$+\frac{1}{2}[\frac{\partial^{3}}{\partial x_{i}^{2}\partial x_{j}}(\nabla
\log\Phi_{P})_{j}+\frac{\partial^{3}}{\partial x_{j}^{2}\partial x_{i}}%
(\nabla\log\Phi_{P})_{i}](y_{0})$

\qquad\qquad\qquad\qquad\qquad\qquad\qquad\qquad\qquad\qquad\qquad\qquad
\qquad\qquad\qquad\qquad\qquad$\blacksquare$

We have from (viii)$^{\ast\ast\ast}$ of \textbf{Table B}$_{1}:$

$[\frac{\partial^{3}}{\partial x_{i}^{2}\partial x_{j}}(\nabla\log\Phi
_{P})_{j}+\frac{\partial^{3}}{\partial x_{i}\partial x_{j}^{2}}(\nabla
$log$\Phi_{P})_{i}](y_{0})=-$ $\left(  \frac{\partial^{3}X_{i}}{\partial
x_{i}\partial x_{j}^{2}}+\frac{\partial^{3}X_{j}}{\partial x_{i}^{2}\partial
x_{j}}\right)  (y_{0}).$

Consequently, we have the very nice \textbf{final }expression for
$\frac{\partial^{4}\Phi_{P}}{\partial x_{i}^{2}\partial x_{j}^{2}}(y_{0}):$

$\frac{\partial^{4}\Phi_{P}}{\partial x_{i}^{2}\partial x_{j}^{2}}%
(y_{0})=[X_{i}^{2}X_{j}^{2}-2X_{i}X_{j}\left(  \frac{\partial X_{j}}{\partial
x_{i}}+\frac{\partial X_{i}}{\partial x_{j}}\right)  -X_{i}^{2}\frac{\partial
X_{j}}{\partial x_{j}}-X_{j}^{2}\frac{\partial X_{i}}{\partial x_{i}}%
](y_{0})\qquad\left(  86\right)  $

$+\frac{1}{2}\left(  \frac{\partial X_{j}}{\partial x_{i}}+\frac{\partial
X_{i}}{\partial x_{j}}\right)  ^{2}(y_{0})+\left(  \frac{\partial X_{i}%
}{\partial x_{i}}\frac{\partial X_{j}}{\partial x_{j}}\right)  (y_{0})\qquad$

$+\frac{2}{3}X_{i}(y_{0})\left(  2\frac{\partial^{2}X_{j}}{\partial
x_{i}\partial x_{j}}+\frac{\partial^{2}X_{i}}{\partial x_{j}^{2}}\right)
(y_{0})+\frac{2}{3}X_{j}(y_{0})\left(  \frac{\partial^{2}X_{j}}{\partial
x_{i}^{2}}+2\frac{\partial^{2}X_{i}}{\partial x_{i}\partial x_{j}}\right)
(y_{0})$

$-\frac{1}{2}\left(  \frac{\partial^{3}X_{i}}{\partial x_{i}\partial x_{j}%
^{2}}+\frac{\partial^{3}X_{j}}{\partial x_{i}^{2}\partial x_{j}}\right)
(y_{0})$

\qquad\qquad\qquad\qquad\qquad\qquad\qquad\qquad\qquad\qquad\qquad\qquad
\qquad\qquad\qquad\qquad$\blacksquare$

We recall the Einstein Convention of summation over repeated indices. As we
have done before, we set:

$\left\Vert X\right\Vert _{M}^{2}=$ $X_{i}^{2}$ where $\left\Vert X\right\Vert
_{M}$ is the norm on the tangent bundle TM and $\left\Vert X\right\Vert
_{P}^{2})=$ $X_{\text{a}}^{2}$ \ where $\left\Vert X\right\Vert _{P}$ is the
norm on the tangent bundle TP and see that:

$[X_{i}^{2}X_{j}^{2}]$ $=(\underset{i=q+1}{\overset{n}{%
{\textstyle\sum}
}}X_{i}^{2})(\underset{j=q+1}{\overset{n}{%
{\textstyle\sum}
}}X_{j}^{2})=(\underset{i=1}{\overset{n}{%
{\textstyle\sum}
}}X_{i}^{2}-\underset{\text{a}=1}{\overset{q}{%
{\textstyle\sum}
}}X_{\text{a}}^{2})(\underset{j=1}{\overset{n}{%
{\textstyle\sum}
}}X_{j}^{2}-\underset{\text{a}=1}{\overset{q}{%
{\textstyle\sum}
}}X_{\text{a}}^{2})$

$\qquad=(\left\Vert X\right\Vert _{M}^{2}-\left\Vert X\right\Vert _{P}%
^{2})(\left\Vert X\right\Vert _{M}^{2}-\left\Vert X\right\Vert _{P}%
^{2})=(\left\Vert X\right\Vert _{M}^{2}-\left\Vert X\right\Vert _{P}^{2})^{2}$

By $\left(  18\right)  ^{\ast}$ above, we have for $j=q+1,...,n,$

\ $\frac{\partial X_{j}}{\partial x_{j}}=[\operatorname{div}_{M}%
X-\operatorname{div}_{P}X+\underset{j=q+1}{\overset{n}{\sum}}<H,j>X_{j}]$

Therefore,

$\underset{i,j=q+1}{\overset{n}{%
{\textstyle\sum}
}}\left(  \frac{\partial X_{i}}{\partial x_{i}}\frac{\partial X_{j}}{\partial
x_{j}}\right)  =(\underset{i=q+1}{\overset{n}{%
{\textstyle\sum}
}}\frac{\partial X_{i}}{\partial x_{i}})(\underset{i,j=q+1}{\overset{n}{%
{\textstyle\sum}
}}\frac{\partial X_{j}}{\partial x_{j}})$

$=[\operatorname{div}_{M}X-\operatorname{div}_{P}%
X+\underset{i=q+1}{\overset{n}{\sum}}<H,i>X_{i}][\operatorname{div}%
_{M}X-\operatorname{div}_{P}X+\underset{j=q+1}{\overset{n}{\sum}}<H,j>X_{j}]$

$=[\operatorname{div}_{M}X-\operatorname{div}_{P}X]^{2}+[\operatorname{div}%
_{M}X-\operatorname{div}_{P}X][\underset{i=q+1}{\overset{n}{\sum}}<H,i>X_{i}]$

$+$ $[\operatorname{div}_{M}X-\operatorname{div}_{P}%
X][\underset{j=q+1}{\overset{n}{\sum}}<H,j>X_{j}%
]+\underset{i,j=q+1}{\overset{n}{\sum}}<H,i><H,j>X_{i}X_{j}]$

We thus have:

$\frac{\partial^{4}\Phi_{P}}{\partial x_{i}^{2}\partial x_{j}^{2}}%
(y_{0})=[\left\Vert X\right\Vert _{M}^{2}-\left\Vert X\right\Vert _{P}%
^{2}]^{2}(y_{0})-2[X_{i}X_{j}\left(  \frac{\partial X_{j}}{\partial x_{i}%
}+\frac{\partial X_{i}}{\partial x_{j}}\right)  ](y_{0})$

$-[\left\Vert X\right\Vert _{M}^{2}-\left\Vert X\right\Vert _{P}^{2}%
](y_{0})[\operatorname{div}_{M}X-\operatorname{div}_{P}%
X+\underset{j=q+1}{\overset{n}{\sum}}<H,j>X_{j}](y_{0})$

$-[\left\Vert X\right\Vert _{M}^{2}-\left\Vert X\right\Vert _{P}^{2}%
](y_{0})[\operatorname{div}_{M}X-\operatorname{div}_{P}%
X+\underset{j=q+1}{\overset{n}{\sum}}<H,j>X_{j}](y_{0}))$

$+\frac{1}{2}\underset{i,j=q+1}{\overset{n}{%
{\textstyle\sum}
}}\left(  \frac{\partial X_{j}}{\partial x_{i}}+\frac{\partial X_{i}}{\partial
x_{j}}\right)  ^{2}(y_{0})$

$+$ $[\operatorname{div}_{M}X-\operatorname{div}_{P}X]^{2}(y_{0}%
)+[\operatorname{div}_{M}X-\operatorname{div}_{P}X](y_{0}%
)[\underset{i=q+1}{\overset{n}{\sum}}<H,i>X_{i}](y_{0})$

$+$ $[\operatorname{div}_{M}X-\operatorname{div}_{P}X](y_{0}%
)[\underset{j=q+1}{\overset{n}{\sum}}<H,j>X_{j}](y_{0}%
)+\underset{i,j=q+1}{\overset{n}{\sum}}[<H,i><H,j>X_{i}X_{j}](y_{0})\qquad$

$+\frac{2}{3}\underset{i,j=q+1}{\overset{n}{%
{\textstyle\sum}
}}[X_{i}\left(  2\frac{\partial^{2}X_{j}}{\partial x_{i}\partial x_{j}}%
+\frac{\partial^{2}X_{i}}{\partial x_{j}^{2}}\right)  +X_{j}\left(
\frac{\partial^{2}X_{j}}{\partial x_{i}^{2}}+2\frac{\partial^{2}X_{i}%
}{\partial x_{i}\partial x_{j}}\right)  ](y_{0})-\frac{1}{2}%
\underset{i,j=q+1}{\overset{n}{%
{\textstyle\sum}
}}\left(  \frac{\partial^{3}X_{i}}{\partial x_{i}\partial x_{j}^{2}}%
+\frac{\partial^{3}X_{j}}{\partial x_{i}^{2}\partial x_{j}}\right)  (y_{0})$

We re-write this in a more elegant way as follows:$\qquad\qquad\qquad
\qquad\qquad\qquad$

$\frac{\partial^{4}\Phi_{P}}{\partial x_{i}^{2}\partial x_{j}^{2}}%
(y_{0})=[\left\Vert X\right\Vert _{M}^{2}-\left\Vert X\right\Vert _{P}%
^{2}]^{2}(y_{0})-2[\left\Vert X\right\Vert _{M}^{2}-\left\Vert X\right\Vert
_{P}^{2}](y_{0})[\operatorname{div}_{M}X\qquad\left(  87\right)  $

$-\operatorname{div}_{P}X+\underset{j=q+1}{\overset{n}{\sum}}<H,j>X_{j}%
](y_{0})\qquad$

$+$ $[\operatorname{div}_{M}X-\operatorname{div}_{P}X]^{2}(y_{0}%
)+[\operatorname{div}_{M}X-\operatorname{div}_{P}X](y_{0}%
)[\underset{i=q+1}{\overset{n}{\sum}}<H,i>X_{i}](y_{0})$

$+$ $[\operatorname{div}_{M}X-\operatorname{div}_{P}X](y_{0}%
)[\underset{j=q+1}{\overset{n}{\sum}}<H,j>X_{j}](y_{0}%
)+\underset{i,j=q+1}{\overset{n}{\sum}}<H,i><H,j>X_{i}X_{j}](y_{0})$

$-2\underset{i,j=q+1}{\overset{n}{%
{\textstyle\sum}
}}[X_{i}X_{j}\left(  \frac{\partial X_{j}}{\partial x_{i}}+\frac{\partial
X_{i}}{\partial x_{j}}\right)  ](y_{0})+\frac{1}{2}%
\underset{i,j=q+1}{\overset{n}{%
{\textstyle\sum}
}}\left(  \frac{\partial X_{j}}{\partial x_{i}}+\frac{\partial X_{i}}{\partial
x_{j}}\right)  ^{2}(y_{0})\qquad$

$+\frac{2}{3}\underset{i,j=q+1}{\overset{n}{%
{\textstyle\sum}
}}[X_{i}\left(  2\frac{\partial^{2}X_{j}}{\partial x_{i}\partial x_{j}}%
+\frac{\partial^{2}X_{i}}{\partial x_{j}^{2}}\right)  +X_{j}\left(
\frac{\partial^{2}X_{j}}{\partial x_{i}^{2}}+2\frac{\partial^{2}X_{i}%
}{\partial x_{i}\partial x_{j}}\right)  ](y_{0})-\frac{1}{2}%
\underset{i,j=q+1}{\overset{n}{%
{\textstyle\sum}
}}\left(  \frac{\partial^{3}X_{i}}{\partial x_{i}\partial x_{j}^{2}}%
+\frac{\partial^{3}X_{j}}{\partial x_{i}^{2}\partial x_{j}}\right)  (y_{0})$

\qquad\qquad\qquad\qquad\qquad\qquad\qquad\qquad\qquad\qquad\qquad\qquad
\qquad\qquad\qquad\qquad\qquad\qquad\qquad\qquad$\blacksquare$

The above is a fairly more geometric presentation of the formula in which we
see the roles played by the divergence of the vector field X on the Riemannian
manifold M and the submanifold P as well as the norms on the tangent bundles
of the Riemannian manifold and the submanifold. We also see the role played by
the \textbf{mean curvatur}e of the submanifold P. The mean curvature will
disappear if we assume that the submanifold is \textbf{totally geodesic}.

We also see that if the Fermi coordinates reduce to normal coordinates, which
is equivalent to the submanifold reducing to the centre of Fermi coordinates
$\left\{  y_{0}\right\}  ,$ then we have a simpler formula in which all the
submanifold terms disappear:

$\frac{\partial^{4}\Phi_{P}}{\partial x_{i}^{2}\partial x_{j}^{2}}%
(y_{0})=[\left\Vert X\right\Vert _{M}^{4}](y_{0})-2[\left\Vert X\right\Vert
_{M}^{2}(y_{0})[\operatorname{div}_{M}X](y_{0})+$ $[\operatorname{div}%
_{M}X]^{2}(y_{0})\qquad\qquad$

$-2\underset{i,j=1}{\overset{n}{%
{\textstyle\sum}
}}[X_{i}X_{j}\left(  \frac{\partial X_{j}}{\partial x_{i}}+\frac{\partial
X_{i}}{\partial x_{j}}\right)  ](y_{0})+\frac{1}{2}%
\underset{i,j=1}{\overset{n}{%
{\textstyle\sum}
}}\left(  \frac{\partial X_{j}}{\partial x_{i}}+\frac{\partial X_{i}}{\partial
x_{j}}\right)  ^{2}(y_{0})\qquad$

$+\frac{2}{3}\underset{i,j=1}{\overset{n}{%
{\textstyle\sum}
}}[X_{i}\left(  2\frac{\partial^{2}X_{j}}{\partial x_{i}\partial x_{j}}%
+\frac{\partial^{2}X_{i}}{\partial x_{j}^{2}}\right)  +X_{j}\left(
\frac{\partial^{2}X_{j}}{\partial x_{i}^{2}}+2\frac{\partial^{2}X_{i}%
}{\partial x_{i}\partial x_{j}}\right)  ](y_{0})-\frac{1}{2}%
\underset{i,j=1}{\overset{n}{%
{\textstyle\sum}
}}\left(  \frac{\partial^{3}X_{i}}{\partial x_{i}\partial x_{j}^{2}}%
+\frac{\partial^{3}X_{j}}{\partial x_{i}^{2}\partial x_{j}}\right)  (y_{0})$

\qquad\qquad\qquad\qquad\qquad\qquad\qquad\qquad\qquad\qquad\qquad\qquad
\qquad\qquad\qquad\qquad\qquad\qquad$\blacksquare$

Simplifying the first line, we have the \textbf{final} expression:

$\frac{\partial^{4}\Phi_{P}}{\partial x_{i}^{2}\partial x_{j}^{2}}%
(y_{0})=[\left\Vert X\right\Vert _{M}^{2}](y_{0})-[\operatorname{div}%
_{M}X]]^{2}(y_{0})\qquad\qquad\qquad\qquad\qquad\qquad\left(  88\right)  $

$-2\underset{i,j=1}{\overset{n}{%
{\textstyle\sum}
}}[X_{i}X_{j}\left(  \frac{\partial X_{j}}{\partial x_{i}}+\frac{\partial
X_{i}}{\partial x_{j}}\right)  ](y_{0})+\frac{1}{2}%
\underset{i,j=1}{\overset{n}{%
{\textstyle\sum}
}}\left(  \frac{\partial X_{j}}{\partial x_{i}}+\frac{\partial X_{i}}{\partial
x_{j}}\right)  ^{2}(y_{0})\qquad$

$\bigskip+\frac{2}{3}[X_{i}\left(  2\frac{\partial^{2}X_{j}}{\partial
x_{i}\partial x_{j}}+\frac{\partial^{2}X_{i}}{\partial x_{j}^{2}}\right)
+X_{j}\left(  \frac{\partial^{2}X_{j}}{\partial x_{i}^{2}}+2\frac{\partial
^{2}X_{i}}{\partial x_{i}\partial x_{j}}\right)  ](y_{0})-\frac{1}{2}\left(
\frac{\partial^{3}X_{i}}{\partial x_{i}\partial x_{j}^{2}}+\frac{\partial
^{3}X_{j}}{\partial x_{i}^{2}\partial x_{j}}\right)  (y_{0})$

\qquad\qquad\qquad\qquad\qquad\qquad\qquad\qquad\qquad\qquad\qquad\qquad
\qquad\qquad\qquad\qquad\qquad\qquad$\blacksquare$

In particular,

$\frac{\partial^{4}\Phi_{P}}{\partial x_{i}^{4}}(y_{0})=[\left\Vert
X\right\Vert _{M}^{2}](y_{0})-\operatorname{div}_{M}X]^{2}(y_{0})$

$-2\underset{i=1}{\overset{n}{%
{\textstyle\sum}
}}2X_{i}^{2}\left(  \frac{\partial X_{i}}{\partial x_{i}}\right)
(y_{0})+\frac{1}{2}\left(  2\frac{\partial X_{i}}{\partial x_{i}}\right)
^{2}(y_{0})\qquad$

$\bigskip+\frac{2}{3}[X_{i}\left(  3\frac{\partial^{2}X_{i}}{\partial
x_{i}^{2}}\right)  +X_{i}\left(  3\frac{\partial^{2}X_{i}}{\partial x_{i}^{2}%
}\right)  ](y_{0})-\frac{1}{2}\left(  2\frac{\partial^{3}X_{i}}{\partial
x_{i}^{3}}\right)  (y_{0})$

$\frac{\partial^{4}\Phi_{P}}{\partial x_{i}^{4}}(y_{0})=[\left\Vert
X\right\Vert _{M}^{2}](y_{0})-\operatorname{div}_{M}X]^{2}(y_{0}%
)-4\underset{i=1}{\overset{n}{%
{\textstyle\sum}
}}X_{i}^{2}\left(  \frac{\partial X_{i}}{\partial x_{i}}\right)  (y_{0})$

$+2\left(  \frac{\partial X_{i}}{\partial x_{i}}\right)  ^{2}(y_{0}%
)+4\underset{i=1}{\overset{n}{%
{\textstyle\sum}
}}[X_{i}\left(  \frac{\partial^{2}X_{i}}{\partial x_{i}^{2}}\right)
](y_{0})-\left(  \frac{\partial^{3}X_{i}}{\partial x_{i}^{3}}\right)  (y_{0})$

$=[\left\Vert X\right\Vert _{M}^{2}](y_{0})-\operatorname{div}_{M}X]^{2}%
(y_{0})$

$-4\underset{i=1}{\overset{n}{%
{\textstyle\sum}
}}X_{i}^{2}\left(  \frac{\partial X_{i}}{\partial x_{i}}\right)
(y_{0})+2[\operatorname{div}_{M}X]^{2}(y_{0})+4\underset{i=1}{\overset{n}{%
{\textstyle\sum}
}}[X_{i}\left(  \frac{\partial^{2}X_{i}}{\partial x_{i}^{2}}\right)
](y_{0})-\underset{i=1}{\overset{n}{%
{\textstyle\sum}
}}\left(  \frac{\partial^{3}X_{i}}{\partial x_{i}^{3}}\right)  (y_{0})$

Finally we have:

$\frac{\partial^{4}\Phi_{P}}{\partial x_{i}^{4}}(y_{0})=[\left\Vert
X\right\Vert _{M}^{4}](y_{0})-2[\left\Vert X\right\Vert _{M}^{2}%
\operatorname{div}_{M}X]](y_{0})+3[\operatorname{div}_{M}X]^{2}(y_{0}%
)\qquad\qquad\left(  89\right)  $

$\qquad\qquad-4[X_{i}^{2}\left(  \frac{\partial X_{i}}{\partial x_{i}}\right)
](y_{0})+4[X_{i}\left(  \frac{\partial^{2}X_{i}}{\partial x_{i}^{2}}\right)
](y_{0})-\left(  \frac{\partial^{3}X_{i}}{\partial x_{i}^{3}}\right)  (y_{0}%
)$\qquad\qquad\qquad\qquad\qquad\qquad\qquad\qquad\qquad\qquad\qquad
\qquad\qquad\qquad\qquad\qquad\qquad\qquad\qquad$\qquad\qquad\qquad
\qquad\qquad\qquad\qquad\qquad\qquad\qquad\qquad\blacksquare$

\subsubsection{TANGENTIAL DERIVATIVES}

(vii) By the definition of the gradient operator we have at a general point
x$_{0}\in$M$_{0}:$

$j,k=1,...,q,q+1,...,n,$

$(\nabla$log$\Phi_{P})_{k}(x_{0})=$ $g^{jk}(x_{0})\frac{\partial}{\partial
x_{j}}$log$\Phi_{P}(x_{0})$

Consequently,

$g_{ik}(x_{0})(\nabla$log$\Phi_{P})_{k}(x_{0})=$ $g_{ik}(x_{0})g^{jk}%
(x_{0})\frac{\partial}{\partial x_{j}}$log$\Phi_{P}(x_{0})=\delta_{i}^{j}%
\frac{\partial}{\partial x_{j}}$log$\Phi_{P}(x_{0})$

$=\frac{\partial}{\partial x_{i}}$log$\Phi_{P}(x_{0})$

From the last equalities above we have:

$\frac{\partial}{\partial x_{i}}$log$\Phi_{P}(x_{0})=$ $g_{ik}(x_{0})(\nabla
$log$\Phi_{P})_{k}(x_{0})$

Hence for $i,k=1,...,q,q+1,...,n$ we have:

$\frac{\partial\Phi_{P}}{\partial x_{i}}(x_{0})=\Phi_{P}(x_{0})\frac{\partial
}{\partial x_{i}}\log\Phi_{P}(x_{0})=\Phi_{P}(x_{0})$g$_{ik}(x_{0})(\nabla
$log$\Phi_{P})_{k}(x_{0})\qquad\left(  90\right)  $

Hence for a = 1,...,q and $k=1,...,q,q+1,...,n,$ we have by (xi) of Table
B$_{1},$

$\frac{\partial\Phi_{P}}{\partial x_{\text{a}}}(y_{0})=\Phi_{P}(y_{0}%
)\frac{\partial}{\partial x_{\text{a}}}\log\Phi_{P}(y_{0})=$ $\Phi_{P}(y_{0}%
)$g$_{\text{a}k}(y_{0})(\nabla$log$\Phi_{P})_{k}(y_{0})$

$\qquad\qquad=(\nabla$log$\Phi_{P})_{\text{a}}(y_{0})=0$

Alternatively, for a = 1,...,q and $k=q+1,...,n$ we have : g$_{\text{a}%
k}(y_{0})$

$=\delta_{\text{a}k}=0$ and so,

$\frac{\partial\Phi_{P}}{\partial x_{\text{a}}}(y_{0})=\Phi_{P}(y_{0}%
)\frac{\partial}{\partial x_{\text{a}}}\log\Phi_{P}(y_{0})=$ $\Phi_{P}(y_{0}%
)$g$_{\text{a}k}(y_{0})(\nabla$log$\Phi_{P})_{k}(y_{0})=0$

(viii) From $\left(  74\right)  $ we have:

$\frac{\partial^{2}\Phi_{P}}{\partial x_{\text{a}}\partial x_{\text{b}}}%
(x_{0})=\frac{\partial}{\partial x_{\text{b}}}[\Phi_{P}(x_{0})$g$_{\text{a}%
k}(x_{0})(\nabla$log$\Phi_{P})_{k}](x_{0})$

$=\frac{\partial\Phi_{P}}{\partial x_{\text{b}}}(x_{0})[$g$_{\text{a}k}%
(x_{0})(\nabla$log$\Phi_{P})_{k}](x_{0})+\Phi_{P}(x_{0})\frac{\partial
}{\partial x_{\text{b}}}[$g$_{\text{a}k}(x_{0})(\nabla$log$\Phi_{P}%
)_{k}](x_{0})$

$=\frac{\partial\Phi_{P}}{\partial x_{\text{b}}}(x_{0})[$g$_{\text{a}k}%
(x_{0})(\nabla$log$\Phi_{P})_{k}](x_{0})+\Phi_{P}(x_{0})\frac{\partial
g_{\text{a}k}}{\partial x_{\text{b}}}(x_{0})(\nabla$log$\Phi_{P})_{k}](x_{0})$

$+$ $\Phi_{P}(x_{0})$g$_{\text{a}k}(x_{0})\frac{\partial}{\partial
x_{\text{b}}}(\nabla$log$\Phi_{P})_{k}](x_{0})$

Therefore,

$\frac{\partial^{2}\Phi_{P}}{\partial x_{\text{a}}\partial x_{\text{b}}}%
(y_{0})=\frac{\partial\Phi_{P}}{\partial x_{\text{b}}}(y_{0})[$g$_{\text{a}%
k}(y_{0})(\nabla$log$\Phi_{P})_{k}](y_{0})$

$+\Phi_{P}(y_{0})\frac{\partial g_{\text{a}k}}{\partial x_{\text{b}}}%
(y_{0})(\nabla$log$\Phi_{P})_{k}](y_{0})+$ $\Phi_{P}(y_{0})$g$_{\text{a}%
k}(y_{0})\frac{\partial}{\partial x_{\text{b}}}(\nabla$log$\Phi_{P}%
)_{k}](y_{0})$

Since $\Phi_{P}(y_{0})=1,$ $\frac{\partial\Phi_{P}}{\partial x_{\text{b}}%
}(y_{0})=0$ by (vii) above, g$_{\text{a}k}(y_{0})=\delta_{\text{a}k}$ and
$\frac{\partial g_{\text{a}k}}{\partial x_{\text{b}}}(y_{0}=0,$

By (xii) of \textbf{Table B}$_{1}.$

$\qquad\frac{\partial^{2}\Phi_{P}}{\partial x_{\text{a}}\partial x_{\text{b}}%
}(y_{0})=$ $\frac{\partial}{\partial x_{\text{b}}}(\nabla$log$\Phi
_{P})_{\text{a}}](y_{0})=0$

(ix) From $\Phi\Phi^{-1}=1$

$\qquad\frac{\partial\Phi_{P}}{\partial x_{\text{a}}}\Phi^{-1}+\Phi
\frac{\partial\Phi_{P}^{-1}}{\partial x_{\text{a}}}=0\qquad$

Since $\Phi(y_{0})=1=\Phi^{-1}(y_{0}),$ we have:

$\qquad\frac{\partial\Phi_{P}^{-1}}{\partial x_{\text{a}}}(y_{0}%
)=-\frac{\partial\Phi_{P}}{\partial x_{\text{a}}}(y_{0})=0$

(x) $\frac{\partial\Phi_{P}^{-1}}{\partial x_{\text{a}}}=-\frac{\partial
\Phi_{P}}{\partial x_{\text{a}}}\Phi^{-2}$

$\qquad\frac{\partial^{2}\Phi_{P}^{-1}}{\partial x_{\text{a}}\partial
x_{\text{b}}}(y_{0})=\frac{\partial}{\partial x_{\text{b}}}[\frac{\partial
\Phi_{P}^{-1}}{\partial x_{\text{a}}}](y_{0})=-\frac{\partial}{\partial
x_{\text{b}}}[\frac{\partial\Phi_{P}}{\partial x_{\text{a}}}\Phi^{-2}](y_{0})$

$\qquad=-\frac{\partial^{2}\Phi_{P}}{\partial x_{\text{a}}\partial
x_{\text{b}}}(y_{0})\Phi^{-2}(y_{0})+2\frac{\partial\Phi_{P}}{\partial
x_{\text{a}}}(y_{0})\frac{\partial\Phi}{\partial x_{\text{b}}}(y_{0})\Phi
^{-3}(y_{0})=0$.

The last equality is due to the fact that:

$\qquad\qquad\frac{\partial^{2}\Phi_{P}}{\partial x_{\text{a}}\partial
x_{\text{b}}}(y_{0})=0=\frac{\partial\Phi_{P}}{\partial x_{\text{a}}}(y_{0})$

\subsubsection{Mixed Derivatives}

(xi) As before we have:

$\qquad\frac{\partial\Phi_{P}}{\partial x_{i}}(x_{0})=\Phi_{P}(x_{0}%
)\frac{\partial}{\partial x_{i}}\log\Phi_{P}(x_{0})=\Phi_{P}(x_{0})$%
g$_{ik}(x_{0})(\nabla$log$\Phi_{P})_{k}(x_{0})\qquad\qquad$

Therefore,

$\frac{\partial^{2}\Phi_{P}}{\partial x_{\text{a}}\partial x_{i}}(x_{0}%
)=\frac{\partial\Phi_{P}}{\partial x_{\text{a}}}(x_{0})$g$_{ik}(x_{0})(\nabla
$log$\Phi_{P})_{k}(x_{0})+\Phi_{P}(x_{0})\frac{\partial}{\partial x_{\text{a}%
}}[$g$_{ik}(\nabla$log$\Phi_{P})_{k}](x_{0})$

$=\frac{\partial\Phi_{P}}{\partial x_{\text{a}}}(x_{0})$g$_{ik}(x_{0})(\nabla
$log$\Phi_{P})_{k}(x_{0})+\Phi_{P}(x_{0})[\frac{\partial\text{g}_{ik}%
}{\partial x_{\text{a}}}(\nabla$log$\Phi_{P})_{k}+$ g$_{ik}\frac{\partial
}{\partial x_{\text{a}}}$($\nabla$log$\Phi_{P}$)$_{k}](x_{0})\qquad\left(
91\right)  $

Therefore,

$\frac{\partial^{2}\Phi_{P}}{\partial x_{\text{a}}\partial x_{i}}(y_{0}%
)=\frac{\partial\Phi_{P}}{\partial x_{\text{a}}}(y_{0})$g$_{ik}(y_{0})(\nabla
$log$\Phi_{P})_{k}(y_{0})+\Phi_{P}(y_{0})[\frac{\partial\text{g}_{ik}%
}{\partial x_{\text{a}}}(\nabla$log$\Phi_{P})_{k}$

$+$ g$_{ik}\frac{\partial}{\partial x_{\text{a}}}$($\nabla$log$\Phi_{P}$%
)$_{k}](y_{0})$

Since $\frac{\partial\Phi_{P}}{\partial x_{\text{a}}}(y_{0})=0,$
$\frac{\partial\text{g}_{ik}}{\partial x_{\text{a}}}(y_{0})=0,$ $\Phi
_{P}(y_{0})=1$ and g$_{ik}(y_{0})=\delta_{ik},$ we have:

By (ix) of \textbf{Table B}$_{1}$.

$\frac{\partial^{2}\Phi_{P}}{\partial x_{\text{a}}\partial x_{i}}(y_{0})=$
$\frac{\partial}{\partial x_{\text{a}}}(\nabla\log\Phi_{P})_{i}(y_{0}%
)=-\frac{\partial X_{i}}{\partial x_{\text{a}}}(y_{0})$

(xii) From $\left(  91\right)  ,$ we have:

$\frac{\partial^{3}\Phi_{P}}{\partial x_{\text{a}}\partial x_{\text{b}%
}\partial x_{i}}(x_{0})=\frac{\partial^{2}\Phi_{P}}{\partial x_{\text{a}%
}\partial x_{\text{b}}}(x_{0})[$g$_{ik}(\nabla$log$\Phi_{P})_{k}](x_{0}%
)+\frac{\partial\Phi_{P}}{\partial x_{\text{a}}}(x_{0})[\frac{\partial
\text{g}_{ik}}{\partial x_{\text{b}}}(\nabla$log$\Phi_{P})_{k}$

$+$ g$_{ik}\frac{\partial}{\partial x_{\text{b}}}(\nabla$log$\Phi_{P}%
)_{k}](x_{0})+\frac{\partial\Phi_{P}}{\partial x_{\text{b}}}(x_{0}%
)[\frac{\partial\text{g}_{ik}}{\partial x_{\text{a}}}(\nabla$log$\Phi_{P}%
)_{k}+$ g$_{ik}\frac{\partial}{\partial x_{\text{a}}}$($\nabla$log$\Phi_{P}%
$)$_{k}](x_{0})$

$\Phi(x_{0})[\frac{\partial^{2}\text{g}_{ik}}{\partial x_{\text{a}}\partial
x_{\text{b}}}(\nabla$log$\Phi_{P})_{k}+\frac{\partial\text{g}_{ik}}{\partial
x_{\text{a}}}\frac{\partial}{\partial x_{\text{b}}}(\nabla$log$\Phi_{P})_{k}+$
$\frac{\partial\text{g}_{ik}}{\partial x_{\text{b}}}\frac{\partial}{\partial
x_{\text{a}}}$($\nabla$log$\Phi_{P}$)$_{k}$

$+$ g$_{ik}\frac{\partial^{2}}{\partial x_{\text{a}}\partial x_{\text{b}}%
}(\nabla\log\Phi_{P})_{k}](x_{0})$

Since $\frac{\partial\Phi_{P}}{\partial x_{\text{a}}}(y_{0})=0=\frac
{\partial^{2}\Phi_{P}}{\partial x_{\text{a}}\partial x_{\text{b}}}(y_{0})$ and
$\frac{\partial\text{g}_{ik}}{\partial x_{\text{a}}}(y_{0})=0=\frac
{\partial^{2}\text{g}_{ik}}{\partial x_{\text{a}}\partial x_{\text{b}}}%
(y_{0}),$ we have:

$\frac{\partial^{3}\Phi_{P}}{\partial x_{\text{a}}\partial x_{\text{b}%
}\partial x_{i}}(y_{0})=\Phi(y_{0})[$ g$_{ik}\frac{\partial^{2}}{\partial
x_{\text{a}}\partial x_{\text{b}}}(\nabla\log\Phi_{P})_{k}](y_{0})$

Since $\Phi(y_{0})=1$ and g$_{ik}(y_{0})=\delta_{ik},$ we have by (x) of
\textbf{Table B}$_{1}$for a,b = 1,...,q:

\qquad$\frac{\partial^{3}\Phi_{P}}{\partial x_{\text{a}}\partial x_{\text{b}%
}\partial x_{i}}(y_{0})=\frac{\partial^{2}}{\partial x_{\text{a}}\partial
x_{\text{b}}}(\nabla\log\Phi_{P})_{i}(y_{0})=-\frac{\partial^{2}X_{i}%
}{\partial x_{\text{a}}\partial x_{\text{b}}}(y_{0})$

(xiii) In particular for a = b, we have:

\qquad$\frac{\partial^{3}\Phi_{P}}{\partial x_{\text{a}}^{2}\partial x_{i}%
}(y_{0})=\frac{\partial^{2}}{\partial x_{\text{a}}^{2}}(\nabla\log\Phi
_{P})_{i}(y_{0})=-\frac{\partial^{2}X_{i}}{\partial x_{\text{a}}^{2}}(y_{0})$

(xiv) By $\left(  32\right)  $ we have:

$\qquad\frac{\partial^{2}\Phi_{P}}{\partial x_{i}\partial x_{j}}(x_{0}%
)=\frac{\partial\Phi_{P}}{\partial x_{j}}(x_{0})$g$_{ik}(x_{0})(\nabla
$log$\Phi_{P})_{k}(x_{0})$

$\qquad+\Phi_{P}(x_{0})\frac{\partial\text{g}_{ik}}{\partial x_{j}}%
(x_{0})(\nabla$log$\Phi_{P})_{k}(x_{0})+$ $\Phi_{P}(x_{0})$g$_{ik}(x_{0}%
)\frac{\partial}{\partial x_{j}}(\nabla$log$\Phi_{P})_{k}(x_{0})$\ 

Therefore,

$\frac{\partial^{3}\Phi_{P}}{\partial x_{\text{c}}\partial x_{i}\partial
x_{j}}(x_{0})$

$=\frac{\partial}{\partial x_{\text{c}}}[\frac{\partial\Phi_{P}}{\partial
x_{j}}$g$_{ik}(\nabla$log$\Phi_{P})_{k}](x_{0})+\frac{\partial}{\partial
x_{\text{c}}}[\Phi_{P}\frac{\partial\text{g}_{ik}}{\partial x_{j}}(\nabla
$log$\Phi_{P})_{k}](x_{0})\qquad\qquad\left(  92\right)  $

$+\frac{\partial}{\partial x_{\text{a}}}[\Phi_{P}$g$_{ik}\frac{\partial
}{\partial x_{j}}$($\nabla$log$\Phi_{P}$)$_{k}](x_{0})$ $=$ R$_{1}(x_{0})$ +
R$_{2}(x_{0})$ + R$_{3}(x_{0})$

where (we omit to explicitely write down

$\frac{\partial\text{g}_{ik}}{\partial x_{\text{c}}}$ and $\frac{\partial
^{2}\text{g}_{ik}}{\partial x_{\text{c}}\partial x_{j}}$ since $\frac
{\partial\text{g}_{ik}}{\partial x_{\text{c}}}(y_{0})=0=\frac{\partial
^{2}\text{g}_{ik}}{\partial x_{\text{c}}\partial x_{j}}(y_{0})$),

R$_{1}(x_{0})=\frac{\partial}{\partial x_{\text{c}}}[\frac{\partial\Phi_{P}%
}{\partial x_{j}}$g$_{ik}(\nabla$log$\Phi_{P})_{k}](x_{0})=[\frac{\partial
^{2}\Phi_{P}}{\partial x_{\text{c}}\partial x_{j}}$g$_{ik}(\nabla$log$\Phi
_{P})_{k}+$ $\frac{\partial\Phi_{P}}{\partial x_{j}}$g$_{ik}\frac{\partial
}{\partial x_{\text{c}}}(\nabla\log\Phi_{P})_{k}](x_{0})$

R$_{2}(x_{0})=\frac{\partial}{\partial x_{\text{c}}}[\Phi_{P}\frac
{\partial\text{g}_{ik}}{\partial x_{j}}(\nabla$log$\Phi_{P})_{k}%
](x_{0})=[\frac{\partial\Phi_{P}}{\partial x_{\text{c}}}\frac{\partial
\text{g}_{ik}}{\partial x_{j}}(\nabla$log$\Phi_{P})_{k}+\Phi_{P}\frac
{\partial\text{g}_{ik}}{\partial x_{j}}\frac{\partial}{\partial x_{\text{c}}%
}(\nabla\log\Phi_{P})_{k}](x_{0})$

R$_{3}(x_{0})=\frac{\partial}{\partial x_{\text{c}}}[\Phi_{P}$g$_{ik}%
\frac{\partial}{\partial x_{j}}$($\nabla$log$\Phi_{P}$)$_{k}](x_{0}%
)=[\frac{\partial\Phi_{P}}{\partial x_{\text{c}}}$g$_{ik}\frac{\partial
}{\partial x_{j}}$($\nabla$log$\Phi_{P}$)$_{k}+\Phi_{P}$g$_{ik}\frac
{\partial^{2}}{\partial x_{\text{c}}\partial x_{j}}(\nabla\log\Phi_{P}%
)_{k}](x_{0})$

Since g$_{ik}(y_{0})=\delta_{ik}$ and $\frac{\partial\text{g}_{ik}}{\partial
x_{j}}(y_{0})=0$ for $i,j=q+1,...,n$ and $k=1,...,q,q+1,...n$, we
have:$\qquad$

$\frac{\partial^{3}\Phi_{P}}{\partial x_{\text{c}}\partial x_{i}\partial
x_{j}}(y_{0})=$ R$_{1}(y_{0})$ + R$_{2}(y_{0})$ + R$_{3}(y_{0})$

$=[\frac{\partial^{2}\Phi_{P}}{\partial x_{\text{c}}\partial x_{j}}(\nabla
$log$\Phi_{P})_{i}+$ $\frac{\partial\Phi_{P}}{\partial x_{j}}\frac{\partial
}{\partial x_{\text{c}}}(\nabla\log\Phi_{P})_{i}](y_{0})$

$+$ $[\frac{\partial\Phi_{P}}{\partial x_{\text{c}}}\frac{\partial}{\partial
x_{j}}(\nabla\log\Phi_{P})_{i}+$ $\Phi_{P}\frac{\partial^{2}}{\partial
x_{\text{c}}\partial x_{j}}(\nabla\log\Phi_{P})_{i}](y_{0})$

Since,

$\Phi_{P}(y_{0})=1;\frac{\partial\Phi_{P}}{\partial x_{\text{c}}}%
(y_{0})=0;\frac{\partial\Phi_{P}}{\partial x_{j}}(y_{0})=-X_{j}(y_{0}%
);(\nabla\log\Phi_{P})_{i}(y_{0})=-X_{i}(y_{0});$

By (ix) of \textbf{Table B}$_{1},$

$\qquad\frac{\partial^{2}\Phi_{P}}{\partial x_{\text{c}}\partial x_{j}}%
(y_{0})=\frac{\partial}{\partial x_{\text{c}}}(\nabla\log\Phi_{P})_{i}%
(y_{0})=-\frac{\partial X_{j}}{\partial x_{\text{c}}}(y_{0})$Therefore,

$\qquad\frac{\partial^{3}\Phi_{P}}{\partial x_{\text{c}}\partial x_{i}\partial
x_{j}}(y_{0})=[(-\frac{\partial X_{j}}{\partial x_{\text{c}}})(-X_{i}%
)+(-X_{j})(-\frac{\partial X_{i}}{\partial x_{\text{c}}})](y_{0})+$
$\frac{\partial^{2}}{\partial x_{\text{c}}\partial x_{j}}(\nabla\log\Phi
_{P})_{i}](y_{0})$

\qquad$\frac{\partial^{3}\Phi_{P}}{\partial x_{\text{c}}\partial x_{i}\partial
x_{j}}(y_{0})=[X_{i}\frac{\partial X_{j}}{\partial x_{\text{c}}}+X_{j}%
\frac{\partial X_{i}}{\partial x_{\text{c}}}](y_{0})+$ $\frac{\partial^{2}%
}{\partial x_{\text{c}}\partial x_{j}}(\nabla\log\Phi_{P})_{i}](y_{0}%
)\qquad\qquad\ \left(  93\right)  $

Similarly,

\qquad$\frac{\partial^{3}\Phi_{P}}{\partial x_{\text{c}}\partial x_{j}\partial
x_{i}}(y_{0})=[X_{j}\frac{\partial X_{i}}{\partial x_{\text{c}}}+X_{i}%
\frac{\partial X_{j}}{\partial x_{\text{c}}}](y_{0})+$ $\frac{\partial^{2}%
}{\partial x_{\text{c}}\partial x_{i}}(\nabla\log\Phi_{P})_{j}](y_{0}%
)\qquad\qquad\left(  94\right)  \ $

Since,

$\qquad\frac{\partial^{3}\Phi_{P}}{\partial x_{\text{c}}\partial x_{i}\partial
x_{j}}(y_{0})=\frac{\partial^{3}\Phi_{P}}{\partial x_{\text{c}}\partial
x_{j}\partial x_{i}}(y_{0}),$

we have by $\left(  93\right)  $ and $\left(  94\right)  :$

\qquad$\frac{\partial^{3}\Phi_{P}}{\partial x_{\text{c}}\partial x_{i}\partial
x_{j}}(y_{0})=\frac{1}{2}[\frac{\partial^{3}\Phi_{P}}{\partial x_{\text{c}%
}\partial x_{i}\partial x_{j}}+\frac{\partial^{3}\Phi_{P}}{\partial
x_{\text{c}}\partial x_{j}\partial x_{i}}](y_{0})\qquad\qquad\qquad
\qquad\qquad\qquad\left(  95\right)  $

\qquad$=\frac{1}{2}[X_{i}\frac{\partial X_{j}}{\partial x_{\text{c}}}%
+X_{j}\frac{\partial X_{i}}{\partial x_{\text{c}}}](y_{0})+\frac{1}{2}%
[X_{j}\frac{\partial X_{i}}{\partial x_{\text{c}}}+X_{i}\frac{\partial X_{j}%
}{\partial x_{\text{c}}}](y_{0})$

\qquad$+\frac{1}{2}[$ $\frac{\partial^{2}}{\partial x_{\text{c}}\partial
x_{j}}(\nabla\log\Phi_{P})_{i}+$ $\frac{\partial^{2}}{\partial x_{\text{c}%
}\partial x_{i}}(\nabla\log\Phi_{P})_{j}](y_{0})$

By (xiv) of \textbf{Table B}$_{1},$

$\frac{\partial^{2}}{\partial x_{\text{a}}\partial x_{j}}(\nabla$log$\Phi
_{P})_{i}(y_{0})+\frac{\partial^{2}}{\partial x_{\text{a}}\partial x_{i}%
}(\nabla\log\Phi_{P})_{j}(y_{0})=-\left(  \frac{\partial^{2}X_{i}}{\partial
x_{\text{a}}\partial x_{j}}-\frac{\partial^{2}X_{j}}{\partial x_{\text{a}%
}\partial x_{i}}\right)  (y_{0})\qquad\left(  96\right)  $

The formulas in $\left(  95\right)  $ and $\left(  96\right)  $ then give:

$\frac{\partial^{3}\Phi_{P}}{\partial x_{\text{c}}\partial x_{i}\partial
x_{j}}(y_{0})=\frac{1}{2}[\frac{\partial^{3}\Phi_{P}}{\partial x_{\text{c}%
}\partial x_{i}\partial x_{j}}+\frac{\partial^{3}\Phi_{P}}{\partial
x_{\text{c}}\partial x_{j}\partial x_{i}}](y_{0})$

\qquad$=[X_{i}\frac{\partial X_{j}}{\partial x_{\text{c}}}+X_{j}\frac{\partial
X_{i}}{\partial x_{\text{c}}}](y_{0})-\frac{1}{2}\left(  \frac{\partial
^{2}X_{i}}{\partial x_{\text{c}}\partial x_{j}}+\frac{\partial^{2}X_{j}%
}{\partial x_{\text{c}}\partial x_{i}}\right)  (y_{0})\qquad\qquad\left(
97\right)  \ \qquad$

From $\left(  92\right)  $ we have:

$\frac{\partial^{3}\Phi_{P}}{\partial x_{\text{c}}\partial x_{i}\partial
x_{j}}(x_{0})=$ R$_{1}(x_{0})$ + R$_{2}(x_{0})$ + R$_{3}(x_{0})\qquad
\qquad\qquad\qquad\qquad\left(  98\right)  \qquad\qquad\qquad\qquad$

where,

\qquad R$_{1}(x_{0})=[\frac{\partial^{2}\Phi_{P}}{\partial x_{\text{c}%
}\partial x_{j}}$g$_{ik}(\nabla$log$\Phi_{P})_{k}+$ $\frac{\partial\Phi_{P}%
}{\partial x_{j}}$g$_{ik}\frac{\partial}{\partial x_{\text{c}}}(\nabla\log
\Phi_{P})_{k}](x_{0})$

\qquad R$_{2}(x_{0})=[\frac{\partial\Phi_{P}}{\partial x_{\text{c}}}%
\frac{\partial\text{g}_{ik}}{\partial x_{j}}(\nabla$log$\Phi_{P})_{k}+\Phi
_{P}\frac{\partial\text{g}_{ik}}{\partial x_{j}}\frac{\partial}{\partial
x_{\text{c}}}(\nabla\log\Phi_{P})_{k}](x_{0})$

\qquad R$_{3}(x_{0})=[\frac{\partial\Phi_{P}}{\partial x_{\text{c}}}$%
g$_{ik}\frac{\partial}{\partial x_{j}}$($\nabla$log$\Phi_{P}$)$_{k}+\Phi_{P}%
$g$_{ik}\frac{\partial^{2}}{\partial x_{\text{c}}\partial x_{j}}(\nabla
\log\Phi_{P})_{k}](x_{0})$

$\frac{\partial^{4}\Phi_{P}}{\partial x_{\text{c}}^{2}\partial x_{i}\partial
x_{j}}(y_{0})=$ $\frac{\partial}{\partial x_{\text{c}}}$R$_{1}(y_{0})+$
$\frac{\partial}{\partial x_{\text{c}}}$R$_{2}(y_{0})+$ $\frac{\partial
}{\partial x_{\text{c}}}$R$_{3}(y_{0})\qquad\left(  99\right)  \qquad
\qquad\qquad\qquad\qquad\qquad\qquad\qquad\qquad\qquad$where,

$\frac{\partial}{\partial x_{\text{c}}}$R$_{1}(y_{0})=\frac{\partial}{\partial
x_{\text{c}}}[\frac{\partial^{2}\Phi_{P}}{\partial x_{\text{c}}\partial x_{j}%
}$g$_{ik}(\nabla$log$\Phi_{P})_{k}](y_{0})+\frac{\partial}{\partial
x_{\text{c}}}[\frac{\partial\Phi_{P}}{\partial x_{j}}$g$_{ik}\frac{\partial
}{\partial x_{\text{c}}}(\nabla$log$\Phi_{P})_{k}](y_{0})\qquad\left(
100\right)  \qquad$

\qquad$\ \ \qquad\qquad\qquad\qquad\ \ \ \ =$ R$_{11}(y_{0})$ + R$_{12}%
(y_{0})$

and,

\qquad R$_{11}(y_{0})=\frac{\partial}{\partial x_{\text{c}}}[\frac
{\partial^{2}\Phi_{P}}{\partial x_{\text{c}}\partial x_{j}}$g$_{ik}(\nabla
$log$\Phi_{P})_{k}](y_{0})$

\qquad R$_{12}(y_{0})=\frac{\partial}{\partial x_{\text{c}}}[\frac
{\partial\Phi_{P}}{\partial x_{j}}$g$_{ik}\frac{\partial}{\partial
x_{\text{c}}}(\nabla$log$\Phi_{P})_{k}](y_{0})$

Recalling that $\frac{\partial g_{ik}}{\partial x_{\text{c}}}(y_{0})=0,$ we have:

R$_{11}(y_{0})=\frac{\partial^{3}\Phi_{P}}{\partial x_{\text{c}}^{2}\partial
x_{j}}(y_{0})[$g$_{ik}(\nabla$log$\Phi_{P})_{k}](y_{0})+$ $\frac{\partial
^{2}\Phi_{P}}{\partial x_{\text{c}}\partial x_{j}}(y_{0})\frac{\partial
}{\partial x_{\text{c}}}[$g$_{ik}(\nabla$log$\Phi_{P})_{k}](y_{0})$

\qquad\qquad$=\frac{\partial^{3}\Phi_{P}}{\partial x_{\text{c}}^{2}\partial
x_{j}}(y_{0})[$g$_{ik}(\nabla$log$\Phi_{P})_{k}](y_{0})+$ $\frac{\partial
^{2}\Phi_{P}}{\partial x_{\text{c}}\partial x_{j}}(y_{0})[$ g$_{ik}%
\frac{\partial}{\partial x_{\text{c}}}(\nabla\log\Phi_{P})_{k}](y_{0})$

\qquad R$_{11}(y_{0})=\frac{\partial^{3}\Phi_{P}}{\partial x_{\text{c}}%
^{2}\partial x_{j}}(y_{0})[(\nabla$log$\Phi_{P})_{i}](y_{0})+$ $\frac
{\partial^{2}\Phi_{P}}{\partial x_{\text{c}}\partial x_{j}}(y_{0}%
)[\frac{\partial}{\partial x_{\text{c}}}(\nabla\log\Phi_{P})_{i}](y_{0})$

By (xi) and (xiii) of \textbf{Table B}$_{4}$\textbf{, and} (i) and (ix) of
\textbf{Table B}$_{1},$ we have:

\qquad R$_{11}(y_{0})=\frac{\partial^{3}\Phi_{P}}{\partial x_{\text{c}}%
^{2}\partial x_{j}}(y_{0})[(\nabla$log$\Phi_{P})_{i}](y_{0})+$ $\frac
{\partial^{2}\Phi_{P}}{\partial x_{\text{c}}\partial x_{j}}(y_{0}%
)[\frac{\partial}{\partial x_{\text{c}}}(\nabla\log\Phi_{P})_{i}](y_{0})$

\qquad\qquad$=-\frac{\partial^{2}X_{j}}{\partial x_{\text{c}}^{2}}%
(y_{0})[-X_{i}(y_{0})]+[-\frac{\partial X_{j}}{\partial x_{\text{c}}}%
(y_{0})][-\frac{\partial X_{i}}{\partial x_{\text{c}}}(y_{0})]$

\qquad\qquad\ \ R$_{11}(y_{0})=X_{i}(y_{0})\frac{\partial^{2}X_{j}}{\partial
x_{\text{c}}^{2}}(y_{0})+\frac{\partial X_{i}}{\partial x_{\text{c}}}%
(y_{0})\frac{\partial X_{j}}{\partial x_{\text{c}}}(y_{0})\qquad\qquad
\qquad\left(  101\right)  \qquad\qquad\qquad\qquad\qquad\qquad\qquad
\qquad\qquad\ $

Next we have:

\qquad R$_{12}(y_{0})=\frac{\partial}{\partial x_{\text{c}}}[\frac
{\partial\Phi_{P}}{\partial x_{j}}$g$_{ik}\frac{\partial}{\partial
x_{\text{c}}}(\nabla$log$\Phi_{P})_{k}](y_{0})$

\qquad$=\frac{\partial^{2}\Phi_{P}}{\partial x_{\text{c}}\partial x_{j}}%
(y_{0})[$g$_{ik}\frac{\partial}{\partial x_{\text{c}}}(\nabla$log$\Phi
_{P})_{k}](y_{0})+$ $\frac{\partial\Phi_{P}}{\partial x_{j}}(y_{0})$%
g$_{ik}[\frac{\partial^{2}}{\partial x_{\text{c}}^{2}}(\nabla\log\Phi_{P}%
)_{k}](y_{0})$

\qquad R$_{12}(y_{0})=\frac{\partial^{2}\Phi_{P}}{\partial x_{\text{c}%
}\partial x_{j}}(y_{0})[\frac{\partial}{\partial x_{\text{c}}}(\nabla$%
log$\Phi_{P})_{i}](y_{0})+\frac{\partial\Phi_{P}}{\partial x_{j}}(y_{0}%
)[\frac{\partial^{2}}{\partial x_{\text{c}}^{2}}(\nabla\log\Phi_{P}%
)_{i}](y_{0})$

By (xi) of \textbf{Table B}$_{4},$ (ix) of \textbf{Table B}$_{1}$, (i) of
\textbf{Table B}$_{1}$ and (x) of \textbf{Table B}$_{1},$ we have:

\qquad R$_{12}=(-\frac{\partial X_{j}}{\partial x_{\text{c}}}(y_{0}%
))[-\frac{\partial X_{i}}{\partial x_{\text{c}}}(y_{0})]+(-X_{j}%
(y_{0}))[-\frac{\partial^{2}X_{i}}{\partial x_{\text{c}}^{2}}(y_{0}%
)]\qquad\qquad\qquad\left(  102\right)  \qquad\qquad$

\qquad\qquad$=\frac{\partial X_{i}}{\partial x_{\text{c}}}(y_{0}%
)\frac{\partial X_{j}}{\partial x_{\text{c}}}(y_{0})+\frac{\partial^{2}X_{i}%
}{\partial x_{\text{c}}^{2}}(y_{0})X_{j}(y_{0})\qquad\qquad\qquad
\qquad\ \ \ \qquad\qquad$

We see from $\left(  100\right)  ,$ $\left(  101\right)  $ and $\left(
102\right)  $ that:

\qquad$\frac{\partial}{\partial x_{\text{c}}}$R$_{1}(y_{0})=$ R$_{11}$(y$_{0}%
$) + R$_{12}$(y$_{0}$)

$=X_{i}(y_{0})\frac{\partial^{2}X_{j}}{\partial x_{\text{c}}^{2}}(y_{0}%
)+\frac{\partial X_{i}}{\partial x_{\text{c}}}(y_{0})\frac{\partial X_{j}%
}{\partial x_{\text{c}}}(y_{0})+\frac{\partial X_{i}}{\partial x_{\text{c}}%
}(y_{0})\frac{\partial X_{j}}{\partial x_{\text{c}}}(y_{0})+\frac{\partial
^{2}X_{i}}{\partial x_{\text{c}}^{2}}(y_{0})X_{j}(y_{0})\qquad\left(
103\right)  $

We next compute $\frac{\partial}{\partial x_{\text{c}}}$R$_{2}(y_{0})$ where,

\qquad R$_{2}(x_{0})=[\frac{\partial\Phi_{P}}{\partial x_{\text{c}}}%
\frac{\partial\text{g}_{ik}}{\partial x_{j}}(\nabla$log$\Phi_{P})_{k}+\Phi
_{P}\frac{\partial\text{g}_{ik}}{\partial x_{j}}\frac{\partial}{\partial
x_{\text{c}}}(\nabla\log\Phi_{P})_{k}](x_{0})\qquad\qquad\qquad$

$\qquad\qquad\ \ \ \ =$ R$_{21}(x_{0})+$ R$_{22}(x_{0})\qquad$

and where for c,...,q; $i,j=q+1,...,n$ and $k=1,...,q,q+1,...,n$ we have,

\qquad R$_{21}(x_{0})=\frac{\partial\Phi_{P}}{\partial x_{\text{c}}}%
(x_{0})[\frac{\partial\text{g}_{ik}}{\partial x_{j}}(\nabla$log$\Phi_{P}%
)_{k}](x_{0})$

\qquad R$_{22}(x_{0})=\Phi_{P}(x_{0})[\frac{\partial\text{g}_{ik}}{\partial
x_{j}}\frac{\partial}{\partial x_{\text{c}}}(\nabla\log\Phi_{P})_{k}](x_{0})$

Now,

$\frac{\partial}{\partial x_{\text{c}}}$R$_{21}(y_{0})=\frac{\partial^{2}%
\Phi_{P}}{\partial x_{\text{c}}^{2}}(y_{0})[\frac{\partial\text{g}_{ik}%
}{\partial x_{j}}(\nabla$log$\Phi_{P})_{k}](y_{0})+\frac{\partial\Phi_{P}%
}{\partial x_{\text{c}}}(y_{0})\frac{\partial}{\partial x_{\text{c}}}%
[\frac{\partial\text{g}_{ik}}{\partial x_{j}}(\nabla$log$\Phi_{P})_{k}%
](y_{0})$

Since $\frac{\partial^{2}\Phi_{P}}{\partial x_{\text{c}}^{2}}(y_{0}%
)=0=\frac{\partial\Phi_{P}}{\partial x_{\text{c}}}(y_{0}),$ we have,

$\frac{\partial}{\partial x_{\text{c}}}$R$_{21}(y_{0})=0\qquad\qquad
\qquad\qquad\qquad\qquad\qquad\qquad\qquad\qquad\qquad$

We next consider:

$\qquad\frac{\partial}{\partial x_{\text{c}}}$R$_{22}=\frac{\partial\Phi_{P}%
}{\partial x_{\text{c}}}(x_{0})[\frac{\partial\text{g}_{ik}}{\partial x_{j}%
}\frac{\partial}{\partial x_{\text{c}}}(\nabla\log\Phi_{P})_{k}](x_{0})$

$\qquad+\Phi_{P}(x_{0})[\frac{\partial^{2}\text{g}_{ik}}{\partial x_{\text{c}%
}\partial x_{j}}\frac{\partial}{\partial x_{\text{c}}}(\nabla\log\Phi_{P}%
)_{k}+\frac{\partial\text{g}_{ik}}{\partial x_{j}}\frac{\partial^{2}}{\partial
x_{\text{c}}^{2}}(\nabla\log\Phi_{P})_{k}](x_{0})$

Since $\frac{\partial^{2}\text{g}_{ik}}{\partial x_{\text{c}}\partial x_{j}%
}(y_{0})=0=\frac{\partial\Phi_{P}}{\partial x_{\text{c}}}(y_{0})$ for c =
1,...,q and $i,j=q++1,...,n,$ we have:

$\frac{\partial}{\partial x_{\text{c}}}$R$_{22}(y_{0})=$
$\underset{k=1}{\overset{n}{%
{\textstyle\sum}
}}[\frac{\partial\text{g}_{ik}}{\partial x_{j}}\frac{\partial^{2}}{\partial
x_{\text{c}}^{2}}(\nabla\log\Phi_{P})_{k}](y_{0})$

$=$ $\underset{k=1}{\overset{n}{%
{\textstyle\sum}
}}[\frac{\partial\text{g}_{i\text{a}}}{\partial x_{j}}\frac{\partial^{2}%
}{\partial x_{\text{c}}^{2}}(\nabla\log\Phi_{P})_{\text{a}}](y_{0})+$
$\underset{k=qa+1}{\overset{n}{%
{\textstyle\sum}
}}[\frac{\partial\text{g}_{ik}}{\partial x_{j}}\frac{\partial^{2}}{\partial
x_{\text{c}}^{2}}(\nabla\log\Phi_{P})_{k}](y_{0})$

By (xiii) of \textbf{Table B}$_{1},$ $\frac{\partial^{2}}{\partial
x_{\text{c}}^{2}}(\nabla\log\Phi_{P})_{\text{a}}](y_{0})=0$ for a,c = 1,...,q
and since

$\frac{\partial\text{g}_{ik}}{\partial x_{j}}(y_{0})=0$ for
$i,j,k=q++1,...,n,$

$\qquad\qquad\qquad\qquad\qquad\frac{\partial}{\partial x_{\text{c}}}$%
R$_{22}(y_{0})=0$

We conclude that,$\qquad\qquad\qquad\qquad\qquad\qquad\qquad\qquad\qquad
\qquad\qquad\qquad\qquad\qquad\qquad\qquad\qquad\qquad\qquad\qquad\qquad
\qquad\qquad\qquad\qquad\qquad\qquad\qquad\qquad\qquad\qquad\qquad\qquad
\qquad\qquad\qquad\qquad\qquad\qquad\qquad\qquad\qquad\qquad\qquad\qquad
\qquad\qquad\qquad\qquad\qquad\qquad\qquad\qquad\qquad\qquad\qquad\qquad$

$\frac{\partial}{\partial x_{\text{c}}}$R$_{2}=\frac{\partial}{\partial
x_{\text{c}}}$R$_{21}+\frac{\partial}{\partial x_{\text{c}}}$R$_{22}%
=0\qquad\qquad\qquad\qquad\qquad\left(  104\right)  $

We then consider:

R$_{3}(x_{0})=[\frac{\partial\Phi_{P}}{\partial x_{\text{c}}}$g$_{ik}%
\frac{\partial}{\partial x_{j}}(\nabla\log\Phi_{P})_{k}+\Phi_{P}$g$_{ik}%
\frac{\partial^{2}}{\partial x_{\text{c}}\partial x_{j}}(\nabla\log\Phi
_{P})_{k}](x_{0})$

\qquad$\qquad\qquad=$ R$_{31}(x_{0})+$ R$_{32}(x_{0})\qquad\qquad\qquad
\qquad\left(  105\right)  $

\qquad R$_{31}(x_{0})=\frac{\partial\Phi_{P}}{\partial x_{\text{c}}}(x_{0}%
)[$g$_{ik}\frac{\partial}{\partial x_{j}}(\nabla\Phi_{P})_{k}](x_{0})$

\qquad R$_{32}(x_{0})=\Phi_{P}(x_{0})[$ g$_{ik}\frac{\partial^{2}}{\partial
x_{\text{c}}\partial x_{j}}(\nabla\log\Phi_{P})_{k}](x_{0})$

\qquad$\frac{\partial}{\partial x_{\text{c}}}$R$_{31}(y_{0})=\frac{\partial
}{\partial x_{\text{c}}}[\frac{\partial\Phi_{P}}{\partial x_{\text{c}}}%
($g$_{ik}\frac{\partial}{\partial x_{j}}(\nabla\log\Phi_{P})_{k})](y_{0})$

\qquad$=\frac{\partial^{2}\Phi_{P}}{\partial x_{\text{c}}^{2}}(y_{0}%
)[($g$_{ik}\frac{\partial}{\partial x_{j}}(\nabla\log\Phi_{P})_{k}%
](y_{0})+\frac{\partial\Phi_{P}}{\partial x_{\text{c}}}(y_{0})\frac{\partial
}{\partial x_{\text{c}}}[$g$_{ik}\frac{\partial}{\partial x_{j}}(\nabla
\log\Phi_{P})_{k}](y_{0})$

Since $\frac{\partial\Phi_{P}}{\partial x_{\text{c}}}(y_{0})=0=\frac
{\partial^{2}\Phi_{P}}{\partial x_{\text{c}}^{2}}(y_{0})$ by (vii) and (viii)
of \textbf{Table B}$_{4},$ we have:

\qquad$\frac{\partial}{\partial x_{\text{c}}}$R$_{31}(y_{0})=0\qquad
\qquad\qquad\qquad\qquad\qquad\qquad\qquad\left(  105\right)  $

Finally we consider:

\qquad$\frac{\partial}{\partial x_{\text{c}}}$R$_{32}=\frac{\partial}{\partial
x_{\text{c}}}[\Phi_{P}($ g$_{ik}\frac{\partial^{2}}{\partial x_{\text{c}%
}\partial x_{j}}(\nabla\log\Phi_{P})_{k})](x_{0})$

$=\frac{\partial\Phi_{P}}{\partial x_{\text{c}}}[$g$_{ik}\frac{\partial^{2}%
}{\partial x_{\text{c}}\partial x_{j}}(\nabla\log\Phi_{P})_{k}](x_{0}%
)+\Phi_{P}(x_{0})\frac{\partial}{\partial x_{\text{c}}}[$g$_{ik}\frac
{\partial^{2}}{\partial x_{\text{c}}\partial x_{j}}(\nabla\log\Phi_{P}%
)_{k}](x_{0})$

$=\frac{\partial\Phi_{P}}{\partial x_{\text{c}}}[$g$_{ik}\frac{\partial^{2}%
}{\partial x_{\text{c}}\partial x_{j}}(\nabla\log\Phi_{P})_{k}](x_{0}%
)+\Phi_{P}(x_{0})[\frac{\partial\text{g}_{ik}}{\partial x_{\text{c}}}%
\frac{\partial^{2}}{\partial x_{\text{c}}\partial x_{j}}(\nabla\log\Phi
_{P})_{k}$

$+$ g$_{ik}\frac{\partial^{3}}{\partial x_{\text{c}}^{2}\partial x_{j}}%
(\nabla\log\Phi_{P})_{k}](x_{0})$

Since $\frac{\partial\Phi_{P}}{\partial x_{\text{c}}}(y_{0})=0,\frac
{\partial\text{g}_{ik}}{\partial x_{\text{c}}}(y_{0})=0$ and g$_{ik}%
(y_{0})=\delta_{ik},$ we have by (xiv) of \textbf{Table B}$_{1}:$

$\frac{\partial}{\partial x_{\text{c}}}$R$_{32}=\frac{\partial^{3}}{\partial
x_{\text{c}}^{2}\partial x_{j}}(\nabla\log\Phi_{P})_{i}(y_{0})=-\frac
{\partial^{3}X_{i}}{\partial x_{\text{c}}^{2}\partial x_{j}}(y_{0}%
)\qquad\qquad\left(  106\right)  \qquad\qquad\qquad$

Consequently by $\left(  104\right)  ,$ $\left(  105\right)  $ and $\left(
106\right)  ,$

$\frac{\partial}{\partial x_{\text{c}}}$R$_{3}=\frac{\partial}{\partial
x_{\text{c}}}$R$_{31}+\frac{\partial}{\partial x_{\text{c}}}$R$_{32}%
=-\frac{\partial^{3}X_{i}}{\partial x_{\text{c}}^{2}\partial x_{j}}%
(y_{0})\qquad\qquad\qquad\left(  107\right)  $

Finally we have by $\left(  103\right)  ,$ $\left(  104\right)  $ and $\left(
107\right)  $ that:

$\frac{\partial^{4}\Phi_{P}}{\partial x_{\text{c}}^{2}\partial x_{i}\partial
x_{j}}(y_{0})=$ $\frac{\partial}{\partial x_{\text{c}}}$R$_{1}(y_{0})+$
$\frac{\partial}{\partial x_{\text{c}}}$R$_{2}(y_{0})+$ $\frac{\partial
}{\partial x_{\text{c}}}$R$_{3}(y_{0})\qquad\left(  108\right)  $

$=X_{i}(y_{0})\frac{\partial^{2}X_{j}}{\partial x_{\text{c}}^{2}}%
(y_{0})+2\frac{\partial X_{i}}{\partial x_{\text{c}}}(y_{0})\frac{\partial
X_{j}}{\partial x_{\text{c}}}(y_{0})+\frac{\partial^{2}X_{i}}{\partial
x_{\text{c}}^{2}}(y_{0})X_{j}(y_{0})-\frac{\partial^{3}X_{i}}{\partial
x_{\text{c}}^{2}\partial x_{j}}(y_{0})\qquad\qquad\ \ \ \ \ \ \qquad
\qquad\qquad\qquad\qquad\ \ \qquad\qquad\qquad\qquad$

In particular,

$\frac{\partial^{4}\Phi_{P}}{\partial x_{\text{c}}^{2}\partial x_{i}^{2}%
}(y_{0})=2X_{i}(y_{0})\frac{\partial^{2}X_{i}}{\partial x_{\text{c}}^{2}%
}(y_{0})+2[\frac{\partial X_{i}}{\partial x_{\text{c}}}]^{2}(y_{0}%
)-\frac{\partial^{3}X_{i}}{\partial x_{\text{c}}^{2}\partial x_{i}}%
(y_{0})\qquad\left(  109\right)  $

(xvi)\qquad$\qquad\qquad\frac{1}{2}\Delta\Phi(y_{0})=\frac{1}{2}%
\underset{i,j=1}{\overset{n}{\sum}}g^{ij}(y_{0})[\frac{\partial^{2}\Phi
}{\partial x_{i}\partial x_{j}}-\underset{k=1}{\overset{n}{\sum}}\Gamma
_{ij}^{k}\frac{\partial\Phi}{\partial x_{k}}](y_{0})$

Since $g^{ij}(y_{0})=\delta^{ij}$ and $\frac{\partial\Phi}{\partial
x_{\text{a}}}(y_{0})=0=\frac{\partial^{2}\Phi}{\partial x_{\text{a}}^{2}}$ for
a = 1,...,q by (vii) and (viii) of \textbf{Table B}$_{4},$ we have:

\qquad$\frac{1}{2}\Delta\Phi(y_{0})=\frac{1}{2}%
\underset{i=q+1}{\overset{n}{\sum}}[\frac{\partial^{2}\Phi}{\partial x_{i}%
^{2}}-\underset{k=q+1}{\overset{n}{\sum}}\Gamma_{ii}^{k}\frac{\partial\Phi
}{\partial x_{k}}](y_{0})$

Further, $\Gamma_{ii}^{k}(y_{0})=0$ for $i,k=q+1,...n$ and so,

$\qquad\frac{1}{2}\Delta\Phi(y_{0})=\frac{1}{2}%
\underset{i=q+1}{\overset{n}{\sum}}[\frac{\partial^{2}\Phi}{\partial x_{i}%
^{2}}-\underset{k=q+1}{\overset{n}{\sum}}\underset{\text{a=1}%
}{\overset{q}{\sum}}\Gamma_{\text{aa}}^{k}(y_{0})\frac{\partial\Phi}{\partial
x_{k}}](y_{0})$

$\frac{\partial\Phi}{\partial x_{k}}(y_{0})=-X_{k}(y_{0})$ by (i) of Table
B$_{4}$ and $\frac{\partial^{2}\Phi_{P}}{\partial x_{i}^{2}}(y_{0})=X_{i}%
^{2}(y_{0})-\frac{\partial X_{i}}{\partial x_{i}}(y_{0})$ by (ii) of Table
B$_{4}.$

Further, by (i) of \textbf{Table A}$_{7}$,

$\underset{\text{a=1}}{\overset{q}{\sum}}\Gamma_{\text{aa}}^{k}(y_{0}%
)=\underset{\text{a=1}}{\overset{q}{\sum}}T_{\text{aa}k}(y_{0})=$
$<H,k>(y_{0}),$

Consequently,$\qquad$

\qquad$\Delta\Phi(y_{0})=\underset{i=q+1}{\overset{n}{\sum}}[X_{i}^{2}%
-\frac{\partial X_{i}}{\partial x_{i}}](y_{0}%
)+\underset{k=q+1}{\overset{n}{\sum}}$ $<H,k>(y_{0})X_{k}(y_{0})$

\qquad$=$ $\underset{i=1}{\overset{n}{\sum}}X_{i}^{2}(y_{0})-$
$\underset{i=1}{\overset{n}{\sum}}\frac{\partial X_{i}}{\partial x_{i}}%
(y_{0})-\underset{\text{a}=1}{\overset{q}{\sum}}X_{\text{a}}^{2}%
(y_{0})+\underset{\text{a}=1}{\overset{q}{\sum}}\frac{\partial X_{i}}{\partial
x_{\text{a}}}(y_{0})+\underset{k=q+1}{\overset{n}{\sum}}<H,k>(y_{0}%
)X_{k}(y_{0})$

$\qquad\Delta\Phi(y_{0})=$ $\left\Vert \text{X}\right\Vert _{M}^{2}(y_{0})-$
divX$_{M}(y_{0})-$ $\left\Vert \text{X}\right\Vert _{P}^{2}(y_{0})$ $+$
divX$_{P}(y_{0})\qquad\left(  110\right)  $

\qquad\qquad\qquad\qquad\qquad\qquad\qquad\qquad\qquad\qquad\qquad\qquad
\qquad\qquad\qquad\qquad$\qquad\qquad\blacksquare$

\section{Table B$_{5}:$ Derivatives of the Scalar Laplacian}

(i)\qquad For $i=q+1,...,n,$ we have$\qquad\qquad\qquad\qquad\qquad
\qquad\qquad$

$\qquad\qquad\frac{1}{24}\frac{\partial}{\partial x_{i}}[\Delta\Phi
](y_{0})=Q_{1}+Q_{2}\qquad$from $\left(  B_{102}\right)  \qquad\qquad
\qquad\qquad\left(  111\right)  \qquad\qquad$

\qquad$\qquad=\frac{1}{12}T_{\text{ab}i}(y_{0})T_{\text{ab}j}(y_{0}%
)X_{j}(y_{0})+\frac{1}{6}\perp_{\text{a}ij}(y_{0})[\frac{\partial X_{j}%
}{\partial x_{\text{a}}}-\perp_{\text{a}jk}X_{k}](y_{0})\qquad\qquad Q_{1}$

\qquad$\qquad+\frac{1}{12}[$ $4X_{j}\frac{\partial X_{j}}{\partial
x_{\text{a}}}-\frac{\partial^{2}X_{j}}{\partial x_{\text{a}}\partial x_{j}%
}](y_{0})\qquad Q_{2}$

$\qquad\qquad+\frac{1}{12}[-\frac{1}{2}\frac{\partial^{2}X_{i}}{\partial
x_{\text{a}}^{2}}+\frac{1}{2}X_{i}\left(  \operatorname{div}X_{M}-\left\Vert
X\right\Vert _{M}^{2}+\left\Vert X\right\Vert _{P}^{2}-\operatorname{div}%
X_{P}-<H,j>X_{j}\right)  ](y_{0})$

$\qquad\qquad+\frac{1}{12}[X_{j}\frac{\partial X_{i}}{\partial x_{j}}+$
$\frac{1}{3}X_{k}R_{jijk}-\frac{1}{2}\frac{\partial^{2}X_{i}}{\partial
x_{j}^{2}}](y_{0})$

\qquad$\qquad+\frac{1}{24}[R_{\text{a}i\text{a}k}-\underset{\text{c=1}%
}{\overset{\text{q}}{\sum}}T_{\text{ac}i}T_{\text{ac}k}-\perp_{\text{a}%
ik}\perp_{\text{a}jk}](y_{0})X_{k}(y_{0})+\frac{1}{36}R_{ijkj}(y_{0}%
)X_{k}(y_{0})$

\qquad$\qquad+\frac{1}{24}<H,j>(y_{0})[X_{i}X_{j}-\frac{1}{2}\left(
\frac{\partial X_{j}}{\partial x_{i}}+\frac{\partial X_{i}}{\partial x_{j}%
}\right)  ](y_{0})$

(ii) For $i=q+1,...,n,$ we have from $\left(  B_{106}\right)  :$ $\qquad$

$\qquad\frac{1}{2}\frac{\partial^{2}}{\partial x_{i}^{2}}[\Delta\Phi
](y_{0})=S_{1}+S_{2}+S_{3}\qquad\qquad\qquad\qquad\qquad\qquad\qquad
\qquad\qquad\qquad\qquad$

$=-2[-R_{\text{a}i\text{b}i}+5\overset{q}{\underset{\text{c}=1}{\sum}%
}T_{\text{ac}i}T_{\text{bc}i}+2\overset{n}{\underset{j=q+1}{\sum}}%
\perp_{\text{a}ij}\perp_{\text{b}ij}](y_{0})\underset{k=q+1}{\overset{n}{\sum
}}T_{\text{ab}k}(y_{0})X_{k}(y_{0})\qquad S_{1}\qquad\left(  112\right)  $

$\ -\frac{8}{3}R_{i\text{a}ij}(y_{0})[\frac{\partial X_{j}}{\partial
x_{\text{a}}}-\underset{k=q+1}{\overset{n}{\sum}}\perp_{\text{a}jk}%
X_{k}](y_{0})+\frac{2}{3}R_{ijik}(y_{0})[X_{j}X_{k}-\frac{1}{2}(\frac{\partial
X_{j}}{\partial x_{k}}+\frac{\partial X_{k}}{\partial x_{j}})](y_{0})$

$-4\perp_{\text{a}ij}(y_{0})[(X_{i}\frac{\partial X_{j}}{\partial x_{\text{a}%
}}+X_{j}\frac{\partial X_{i}}{\partial x_{\text{a}}})-\frac{1}{4}\left(
\frac{\partial^{2}X_{i}}{\partial x_{\text{a}}\partial x_{j}}+\frac
{\partial^{2}X_{j}}{\partial x_{\text{a}}\partial x_{i}}\right)
](y_{0})\qquad\qquad S_{2}$

$\ -2\perp_{\text{a}ij}(y_{0})$ $\overset{q}{\underset{\text{b}=1}{\sum}%
}[\perp_{\text{b}ik}T_{\text{ab}j}+\frac{2}{3}(2R_{\text{a}ijk}+R_{\text{a}%
jik}+R_{\text{a}kji})](y_{0})X_{k}(y_{0}$

$+\frac{1}{2}T_{\text{bb}k}(y_{0})\frac{\partial^{2}X_{k}}{\partial
x_{\text{a}}^{2}}(y_{0})+[(\frac{\partial X_{j}}{\partial x_{\text{a}}}%
)^{2}+X_{j}\frac{\partial^{2}X_{j}}{\partial x_{\text{a}}^{2}}](y_{0})$
$-\frac{1}{2}\frac{\partial^{3}X_{j}}{\partial x_{\text{a}}^{2}\partial x_{j}%
}(y_{0})\qquad S_{3}$

$+\frac{1}{2}[T_{\text{aa}j}\frac{\partial^{2}X_{j}}{\partial x_{\text{b}}%
^{2}}](y_{0})+\frac{1}{2}[2[(\frac{\partial X_{i}}{\partial x_{\text{a}}}%
)^{2}+X_{i}\frac{\partial^{2}X_{i}}{\partial x_{\text{a}}^{2}}]-\frac
{\partial^{3}X_{i}}{\partial x_{\text{a}}^{2}\partial x_{i}}](y_{0})$

$+\frac{1}{2}.\frac{1}{6}[\{4\nabla_{i}R_{i\text{a}j\text{a}}+2\nabla
_{j}R_{i\text{a}i\text{a}}+$ $8(\overset{q}{\underset{\text{c=1}}{%
{\textstyle\sum}
}}R_{\text{a}i\text{c}i}^{{}}T_{\text{ac}j}+\;\overset{n}{\underset{k=q+1}{%
{\textstyle\sum}
}}R_{\text{a}iik}\perp_{\text{a}jk})$

\ $+8(\overset{q}{\underset{\text{c=1}}{%
{\textstyle\sum}
}}R_{\text{a}i\text{c}j}^{{}}T_{\text{ac}i}+\;\overset{n}{\underset{k=q+1}{%
{\textstyle\sum}
}}R_{\text{a}ijk}\perp_{\text{a}ik})+8(\overset{q}{\underset{\text{c=1}}{%
{\textstyle\sum}
}}R_{\text{a}j\text{c}i}^{{}}T_{\text{ac}i}+\;\overset{n}{\underset{k=q+1}{%
{\textstyle\sum}
}}R_{\text{a}jik}\perp_{\text{a}ik})\}$\ 

$\ +\frac{2}{3}\underset{k=q+1}{\overset{n}{\sum}}\{T_{\text{aa}k}%
(R_{ijik}+3\overset{q}{\underset{\text{c}=1}{\sum}}\perp_{\text{c}ij}%
\perp_{\text{c}ik})\}](y_{0})X_{j}(y_{0})$

$-[R_{\text{a}i\text{a}j}-\underset{\text{c=1}}{\overset{\text{q}}{\sum}%
}T_{\text{ac}i}T_{\text{ac}j}-\overset{n}{\underset{k=q+1}{\sum}}%
(\perp_{\text{a}ik}\perp_{\text{a}jk}](y_{0})\times\lbrack X_{i}X_{j}-\frac
{1}{2}\left(  \frac{\partial X_{i}}{\partial x_{j}}+\frac{\partial X_{j}%
}{\partial x_{i}}\right)  ](y_{0})$

$-\frac{1}{2}T_{\text{aa}j}[-X_{i}^{2}X_{j}+X_{j}\frac{\partial X_{i}%
}{\partial x_{i}}\ +X_{i}\left(  \frac{\partial X_{j}}{\partial x_{i}}%
+\frac{\partial X_{i}}{\partial x_{j}}\right)  -\frac{1}{3}\left(
\frac{\partial^{2}X_{j}}{\partial x_{i}^{2}}+2\frac{\partial^{2}X_{i}%
}{\partial x_{i}\partial x_{j}}\right)  ](y_{0})$

$+$ $[X_{i}^{2}X_{j}^{2}-2X_{i}X_{j}\left(  \frac{\partial X_{j}}{\partial
x_{i}}+\frac{\partial X_{i}}{\partial x_{j}}\right)  -X_{i}^{2}\frac{\partial
X_{j}}{\partial x_{j}}-X_{j}^{2}\frac{\partial X_{i}}{\partial x_{i}}%
](y_{0})\qquad$

$+\frac{1}{2}\left(  \frac{\partial X_{j}}{\partial x_{i}}+\frac{\partial
X_{i}}{\partial x_{j}}\right)  ^{2}(y_{0})+\left(  \frac{\partial X_{i}%
}{\partial x_{i}}\frac{\partial X_{j}}{\partial x_{j}}\right)  (y_{0})\qquad$

$+\frac{2}{3}X_{i}(y_{0})\left(  2\frac{\partial^{2}X_{j}}{\partial
x_{i}\partial x_{j}}+\frac{\partial^{2}X_{i}}{\partial x_{j}^{2}}\right)
(y_{0})+\frac{2}{3}X_{j}(y_{0})\left(  \frac{\partial^{2}X_{j}}{\partial
x_{i}^{2}}+2\frac{\partial^{2}X_{i}}{\partial x_{i}\partial x_{j}}\right)
(y_{0})$

$-\frac{1}{2}\left(  \frac{\partial^{3}X_{i}}{\partial x_{i}\partial x_{j}%
^{2}}+\frac{\partial^{3}X_{j}}{\partial x_{i}^{2}\partial x_{j}}\right)
(y_{0})+\frac{1}{2}\underset{\text{a=1}}{\overset{q}{\sum}}[\frac{4}%
{3}R_{\text{a}jij}\frac{\partial X_{i}}{\partial x_{\text{a}}}](y_{0})$

$\ -\frac{1}{6}\overset{n}{\underset{k=q+1}{\sum}}%
[4\overset{q}{\underset{\text{a}=1}{\sum}}\perp_{\text{a}ik}R_{ij\text{a}%
j}+(\nabla_{i}R_{kjij}+\nabla_{j}R_{ijik}+\nabla_{k}R_{ijij})](y_{0}%
)X_{k}(y_{0})$

$-\frac{2}{3}\underset{k=q+1}{\overset{n}{\sum}}R_{ijkj}(y_{0})[X_{i}%
X_{k}-\frac{1}{2}\left(  \frac{\partial X_{i}}{\partial x_{k}}+\frac{\partial
X_{k}}{\partial x_{i}}\right)  ](y_{0})\qquad$

\qquad\qquad\qquad\qquad\qquad\qquad\qquad\qquad\qquad\qquad\qquad\qquad
\qquad\qquad\qquad\qquad$\blacksquare$

(iii) For c = 1,...,q, we have from $\left(  B_{108}\right)  $ below:$\qquad$

$\qquad\frac{1}{2}\frac{\partial}{\partial x_{\text{c}}}[\Delta\Phi
](y_{0})=2[X_{j}\frac{\partial X_{j}}{\partial x_{\text{c}}}](y_{0})-\frac
{1}{2}\left(  \frac{\partial^{2}X_{j}}{\partial x_{\text{c}}\partial x_{j}%
}\right)  (y_{0})+\frac{1}{2}[T_{\text{aa}i}\frac{\partial X_{i}}{\partial
x_{\text{c}}}](y_{0})$ $\qquad\left(  113\right)  $

\qquad\qquad\qquad\qquad$=2[X_{j}\frac{\partial X_{j}}{\partial x_{\text{c}}%
}](y_{0})-\frac{1}{2}\left(  \frac{\partial^{2}X_{j}}{\partial x_{\text{c}%
}\partial x_{j}}\right)  (y_{0})+\frac{1}{2}<H,i>(y_{0})\frac{\partial X_{i}%
}{\partial x_{\text{c}}}(y_{0})$

\qquad\qquad\qquad\qquad\qquad\qquad\qquad\qquad\qquad\qquad\qquad\qquad
\qquad\qquad\qquad\qquad$\blacksquare$

(iv) For c = 1,...,q, we have from $\left(  B_{109}\right)  $ below:$\qquad$

$\qquad\frac{1}{2}\frac{\partial^{2}}{\partial x_{\text{c}}^{2}}[\Delta
\Phi](y_{0})=\ \frac{1}{2}[\frac{\partial^{4}\Phi}{\partial x_{\text{c}}%
^{2}\partial x_{j}^{2}}-\Gamma_{jj}^{k}\frac{\partial^{3}\Phi}{\partial
x_{\text{c}}^{2}\partial x_{k}}](y_{0})$

$\qquad=$ $[(\frac{\partial X_{j}}{\partial x_{\text{c}}})^{2}+X_{j}%
\frac{\partial^{2}X_{j}}{\partial x_{\text{c}}^{2}}](y_{0})$ $-\frac{1}%
{2}\frac{\partial^{3}X_{j}}{\partial x_{\text{c}}^{2}\partial x_{j}}%
(y_{0})+\frac{1}{2}\underset{\text{a}=1}{\overset{q}{\sum}}[T_{\text{aa}%
k}\frac{\partial^{2}X_{k}}{\partial x_{\text{c}}^{2}}](y_{0})\qquad
\qquad\left(  114\right)  $

$\qquad=$ $[(\frac{\partial X_{i}}{\partial x_{\text{c}}})^{2}+X_{i}%
\frac{\partial^{2}X_{i}}{\partial x_{\text{c}}^{2}}](y_{0})$ $-\frac{1}%
{2}\frac{\partial^{3}X_{i}}{\partial x_{\text{c}}^{2}\partial x_{i}}%
(y_{0})+\frac{1}{2}\underset{\text{a}=1}{\overset{q}{\sum}}[T_{\text{aa}%
j}\frac{\partial^{2}X_{j}}{\partial x_{\text{c}}^{2}}](y_{0})$

$\qquad\qquad\qquad\qquad\qquad\qquad\qquad\qquad\qquad\qquad\qquad
\qquad\qquad\qquad\qquad\qquad\blacksquare\qquad\qquad\qquad\qquad\qquad
\qquad\qquad\qquad\qquad$

(v) For c = 1,...,q, we have from $\left(  B_{111}\right)  $ below:\qquad

$\qquad\frac{\partial}{\partial x_{\text{c}}}[\frac{\text{L}\Psi}{\Psi}%
](y_{0})$

$=-[2X_{i}\frac{\partial X_{i}}{\partial x_{\text{c}}}](y_{0})+2[X_{j}%
\frac{\partial X_{j}}{\partial x_{\text{c}}}](y_{0})-\frac{1}{2}\left(
\frac{\partial^{2}X_{j}}{\partial x_{\text{c}}\partial x_{j}}\right)  (y_{0})$
$+\frac{1}{2}[T_{\text{aa}i}\frac{\partial X_{i}}{\partial x_{\text{c}}%
}](y_{0})+$ $\frac{\partial\text{V}}{\partial x_{\text{c}}}(y_{0}%
)\qquad\left(  115\right)  $

$=-[2X_{i}\frac{\partial X_{i}}{\partial x_{\text{c}}}](y_{0})+2[X_{j}%
\frac{\partial X_{j}}{\partial x_{\text{c}}}](y_{0})-\frac{1}{2}\left(
\frac{\partial^{2}X_{j}}{\partial x_{\text{c}}\partial x_{j}}\right)
(y_{0})+\frac{1}{2}[<H,i>\frac{\partial X_{i}}{\partial x_{\text{c}}}%
](y_{0})+$ $\frac{\partial\text{V}}{\partial x_{\text{c}}}(y_{0})$

(vi) For c = 1,...,q, we have from $\left(  B_{112}\right)  $ below:

$\qquad\frac{\partial^{2}}{\partial x_{\text{c}}^{2}}[\frac{\text{L}\Psi}%
{\Psi}](y_{0})$

$=[(\frac{\partial X_{j}}{\partial x_{\text{c}}})^{2}+X_{j}\frac{\partial
^{2}X_{j}}{\partial x_{\text{c}}^{2}}](y_{0})$ $-\frac{1}{2}\frac{\partial
^{3}X_{j}}{\partial x_{\text{c}}^{2}\partial x_{j}}(y_{0})+\frac{1}%
{2}<H,i>(y_{0})\frac{\partial^{2}X_{i}}{\partial x_{\text{c}}^{2}}%
(y_{0})\qquad\qquad\qquad\qquad\left(  116\right)  $

\qquad\qquad\qquad\qquad$\qquad-2[(\frac{\partial X_{i}}{\partial x_{\text{c}%
}})^{2}+X_{i}\frac{\partial^{2}X_{i}}{\partial x_{\text{c}}^{2}})](y_{0})+$
$\frac{\partial^{2}\text{V}}{\partial x_{\text{c}}^{2}}(y_{0})$

\qquad\qquad\qquad\qquad\qquad\qquad\qquad\qquad\qquad\qquad\qquad\qquad
\qquad\qquad\qquad\qquad\qquad\qquad$\blacksquare$

(vii) For $i=q+1,...,n,$ we have from $\left(  B_{113}\right)  $ below:

$R=\frac{1}{12}$ $\frac{\partial}{\partial x_{i}}\left(  \frac{\text{L}\Psi
}{\Psi}\right)  (y_{0})=R_{1}+R_{2}+R_{3}+R_{4}+R_{5}\qquad\qquad\qquad\qquad$

$=-\frac{1}{24}<H,i>(y_{0})\times\frac{1}{24}[3<H,i>^{2}+2(\tau^{M}-3\tau
^{P}\ +\overset{q}{\underset{\text{a=1}}{\sum}}\varrho_{\text{aa}}%
^{M}+\overset{q}{\underset{\text{a,b}=1}{\sum}}R_{\text{abab}}^{M}%
)](y_{0})\qquad R_{1}\ $

$-\frac{1}{24}\underset{k=q+1}{\overset{n}{%
{\textstyle\sum}
}}<H,k>(y_{0})\underset{\text{a,b=}1}{\overset{q}{%
{\textstyle\sum}
}}T_{\text{ab}k}^{2}(y_{0})\phi(y_{0})\qquad$

$+\frac{5}{64}<H,i><H,j>^{2}\phi(y_{0})\qquad\qquad\qquad\qquad\qquad
\qquad\ \ \ \ $\qquad\qquad

$+\frac{1}{96}<H,j>(y_{0})$

$\times\lbrack(2\varrho_{ij}+4\overset{q}{\underset{\text{a}=1}{\sum}%
}R_{i\text{a}j\text{a}}-3\overset{q}{\underset{\text{a},b=1}{\sum}%
}T_{\text{aa}j}T_{\text{bb}i}-T_{\text{ab}j}T_{\text{ab}i}%
-3\overset{q}{\underset{\text{a},\text{b}=1}{\sum}}T_{\text{aa}i}%
T_{\text{bb}j}-T_{\text{ab}i}T_{\text{ab}j})](y_{0})\phi(y_{0})$

$+\frac{1}{96}<H,i>[$ $\tau^{M}-3\tau^{P}+\overset{q}{\underset{\text{a}%
=1}{\sum}}\varrho_{\text{aa}}+\overset{q}{\underset{\text{a,b}=1}{\sum}%
}R_{\text{abab}}](y_{0})\phi(y_{0})$

$+\frac{1}{288}[\nabla_{i}\varrho_{jj}-2\varrho_{ij}%
<H,j>+\overset{q}{\underset{\text{a}=1}{\sum}}(\nabla_{i}R_{\text{a}%
j\text{a}j}-4R_{i\text{a}j\text{a}}<H,j>)+4\overset{q}{\underset{\text{a,b}%
=1}{\sum}}R_{i\text{a}j\text{b}}T_{\text{ab}j}$

$+2\overset{q}{\underset{\text{a,b,c}=1}{\sum}}(T_{\text{aa}i}T_{\text{bb}%
j}T_{\text{cc}j}-3T_{\text{aa}i}T_{\text{bc}j}T_{\text{bc}j}+2T_{\text{ab}%
i}T_{\text{bc}j}T_{\text{ca}j})](y_{0})\phi(y_{0})$\qquad\qquad\qquad
\qquad\qquad\ \ 

$+\frac{1}{288}[\nabla_{j}\varrho_{ij}-2\varrho_{ij}%
<H,j>+\overset{q}{\underset{\text{a}=1}{\sum}}(\nabla_{j}R_{\text{a}%
i\text{a}j}-4R_{j\text{a}i\text{a}}<H,j>)$

$+4\overset{q}{\underset{\text{a,b}=1}{\sum}}R_{j\text{a}i\text{b}%
}T_{\text{ab}j}+2\overset{q}{\underset{\text{a,b,c}=1}{\sum}}(T_{\text{aa}%
j}T_{\text{bb}i}T_{\text{cc}j}-3T_{\text{aa}j}T_{\text{bc}i}T_{\text{bc}%
j}+2T_{\text{ab}j}T_{\text{bc}i}T_{\text{ca}j})](y_{0})\phi(y_{0})$

$+\frac{1}{288}[\nabla_{j}\varrho_{ij}-2\varrho_{jj}%
<H,i>+\overset{q}{\underset{\text{a}=1}{\sum}}(\nabla_{j}R_{\text{a}%
i\text{a}j}-4R_{j\text{a}j\text{a}}<H,i>)+4\overset{q}{\underset{\text{a}%
,b=1}{\sum}}R_{j\text{a}j\text{b}}T_{\text{ab}i}$

$+2\overset{q}{\underset{\text{a,b,c}=1}{\sum}}(T_{\text{aa}j}T_{\text{bb}%
j}T_{\text{cc}i}-3T_{\text{aa}j}T_{\text{bc}j}T_{\text{bc}i}+2T_{\text{ab}%
j}T_{\text{bc}j}T_{\text{ca}i})](y_{0})\phi(y_{0})$

$-\frac{1}{48}$%
$<$%
H,$k$%
$>$%
(y$_{0}$)$[R_{\text{a}i\text{a}k}-\underset{\text{c=1}}{\overset{\text{q}%
}{\sum}}T_{\text{ac}i}.T_{\text{ac}k}+\frac{2}{3}R_{ijkj}](y_{0})\phi(y_{0})$

$-\frac{1}{24}$%
$<$%
H,$k$%
$>$%
(y$_{0}$)$[\frac{3}{4}$%
$<$%
H,$i$%
$>$%
$<$%
H,$k$%
$>$%
\ $+\frac{1}{12}(2\varrho_{ik}+4\overset{q}{\underset{\text{a}=1}{\sum}%
}R_{i\text{a}k\text{a}}-6\overset{q}{\underset{\text{a},b=1}{\sum}%
}T_{\text{aa}i}T_{\text{bb}k}-T_{\text{ab}i}T_{\text{ab}k})]\phi(y_{0}$

$+\frac{1}{24}[\left\Vert \text{X}\right\Vert ^{2}-$ divX $-\left\Vert
\text{X}\right\Vert _{P}^{2}+$ divX$_{P}](y_{0})$X$_{i}(y_{0})+\frac{1}%
{24}\frac{\partial}{\partial x_{i}}[\Delta\Phi](y_{0})\qquad\qquad$

$+\frac{1}{12}T_{\text{ab}i}(y_{0})T_{\text{ab}j}(y_{0})X_{j}(y_{0}%
)\qquad\qquad\qquad\qquad\qquad\qquad\qquad\qquad\qquad\qquad$

$+\frac{1}{12}[X_{j}\frac{\partial X_{j}}{\partial x_{\text{a}}}+X_{j}%
\frac{\partial X_{j}}{\partial x_{\text{a}}}](y_{0})+[X_{j}\frac{\partial
X_{j}}{\partial x_{\text{a}}}+X_{j}\frac{\partial X_{j}}{\partial x_{\text{a}%
}}](y_{0})-\frac{1}{2}\left(  \frac{\partial^{2}X_{j}}{\partial x_{\text{a}%
}\partial x_{j}}+\frac{\partial^{2}X_{j}}{\partial x_{\text{a}}\partial x_{j}%
}\right)  (y_{0})$

$-\frac{1}{24}\frac{\partial^{2}X_{i}}{\partial x_{\text{a}}^{2}}(y_{0}%
)+\frac{1}{2}X_{i}(y_{0})[\frac{\partial X_{j}}{\partial x_{j}}-X_{j}%
X_{j}](y_{0})+X_{j}(y_{0})\frac{\partial X_{i}}{\partial x_{j}}(y_{0}%
)\qquad\ \qquad\ \qquad\qquad\qquad$

$+$ $\frac{1}{36}X_{k}(y_{0})R_{jijk}(y_{0})-\frac{1}{2}\frac{\partial
^{2}X_{i}}{\partial x_{j}^{2}}(y_{0})$

$+\frac{1}{24}[R_{\text{a}i\text{a}k}-\underset{\text{c=1}}{\overset{\text{q}%
}{\sum}}T_{\text{ac}i}T_{\text{ac}k}](y_{0})X_{k}(y_{0})+\frac{1}{3}%
R_{ijkj}(y_{0})X_{k}(y_{0})$

$+\frac{1}{24}<H,j>(y_{0})[X_{i}X_{j}-\frac{1}{2}\left(  \frac{\partial X_{j}%
}{\partial x_{i}}+\frac{\partial X_{i}}{\partial x_{j}}\right)  ](y_{0})$

$-\frac{1}{6}X_{j}(y_{0})\frac{\partial X_{j}}{\partial x_{i}}(y_{0})+\frac
{1}{12}\frac{\partial\text{V}}{\partial x_{i}}(y_{0})\qquad\qquad\qquad
R_{4}\qquad R_{5}$

\qquad\qquad\qquad\qquad\qquad\qquad\qquad\qquad\qquad\qquad\qquad\qquad
\qquad\qquad\qquad\qquad\qquad\qquad$\blacksquare$

(viii)\qquad I$_{321}$ $\mathbf{=}\frac{1}{12}$ $\frac{\partial^{2}}{\partial
x_{i}^{2}}\left\{  \Psi^{-1}L\Psi\right\}  (y_{0})=$ I$_{3211}+$ I$_{3212}+$
I$_{3213}+$I$_{3214}$ $+$ I$_{3215}$

$\qquad$I$_{3211}=\frac{1}{24}$ $\frac{\partial^{2}}{\partial x_{i}^{2}%
}(\theta^{\frac{1}{2}}\Delta\theta^{-\frac{1}{2}})(y_{0})$ is given by (xii)
in \textbf{Appendix A}$_{10}$

$\qquad$I$_{3212}=\frac{1}{24}$ $\frac{\partial^{2}}{\partial x_{i}^{2}}%
(\Phi^{-1}\Delta\Phi)(y_{0})\phi(y_{0})$ is from $\left(  B_{114}\right)  $ below

\qquad\ I$_{3213}=\frac{1}{12}\frac{\partial^{2}}{\partial x_{i}^{2}}%
(<\nabla\log\theta^{-\frac{1}{2}},\nabla\log\Phi+X>)\phi(y_{0})$ is from
$\left(  B_{115}\right)  $ below

\qquad I$_{3214}=\frac{1}{12}\frac{\partial^{2}}{\partial x_{i}^{2}}%
[<\nabla\log\Phi,X>](y_{0})\phi(y_{0})$ is from $\left(  B_{116}\right)  $ below

\qquad I$_{3215}=\frac{1}{12}\frac{\partial^{2}\text{V}}{\partial x_{i}^{2}%
}(y_{0})\phi(y_{0})$ is from $\left(  B_{117}\right)  $ below

We thus have:

I$_{321}$ $\mathbf{=}\frac{1}{12}$ $\frac{\partial^{2}}{\partial x_{i}^{2}%
}\left\{  \frac{\text{L}\Psi}{\Psi}\right\}  (y_{0})\phi(y_{0})$

\qquad$=\frac{1}{24}$ $\frac{\partial^{2}}{\partial x_{i}^{2}}(\theta
^{\frac{1}{2}}\Delta\theta^{-\frac{1}{2}})(y_{0})\phi(y_{0})$ \qquad
I$_{3211}$

$+\frac{1}{24}[\left\Vert \text{X}\right\Vert _{M}^{2}+\operatorname{div}%
$X$_{M}-\left\Vert \text{X}\right\Vert _{P}^{2}-\operatorname{div}X_{P}%
](y_{0})[\left\Vert \text{X}\right\Vert _{M}^{2}-\operatorname{div}$%
X$_{M}-\left\Vert \text{X}\right\Vert _{P}^{2}+\operatorname{div}$%
X$_{P}](y_{0})\phi(y_{0})\qquad$I$_{3212}\mid$I$_{32121}$

$+\frac{1}{6}X_{i}(y_{0})T_{\text{ab}i}(y_{0})T_{\text{ab}j}(y_{0})X_{j}%
(y_{0})\phi(y_{0})+\frac{1}{3}\perp_{\text{a}ij}(y_{0})X_{i}(y_{0}%
)[\frac{\partial X_{j}}{\partial x_{\text{a}}}-\perp_{\text{a}jk}X_{k}%
](y_{0})\phi(y_{0})\qquad$I$_{32122}\qquad Q_{1}$

$+$ $\frac{2}{3}X_{i}(y_{0})X_{j}(y_{0})\frac{\partial X_{j}}{\partial
x_{\text{a}}}(y_{0})\phi(y_{0})-\frac{1}{6}X_{i}(y_{0})\frac{\partial^{2}%
X_{j}}{\partial x_{\text{a}}\partial x_{j}}(y_{0})\phi(y_{0})\qquad
\qquad\qquad Q_{2}$

$-\frac{1}{12}X_{i}(y_{0})\frac{\partial^{2}X_{i}}{\partial x_{\text{a}}^{2}%
}(y_{0})\phi(y_{0})$

$+\frac{1}{12}X_{i}^{2}(y_{0})[\operatorname{div}X_{M}-\left\Vert X\right\Vert
_{M}^{2}+\left\Vert X\right\Vert _{P}^{2}-\operatorname{div}X_{P}%
-<H,j>(y_{0})X_{j}(y_{0})]\phi(y_{0})$

$+\frac{1}{6}X_{i}(y_{0})X_{j}(y_{0})\frac{\partial X_{i}}{\partial x_{j}%
}(y_{0})\phi(y_{0})+$ $\frac{1}{18}X_{i}(y_{0})X_{k}(y_{0})R_{jijk}(y_{0}%
)\phi(y_{0})-\frac{1}{12}X_{i}(y_{0})\frac{\partial^{2}X_{i}}{\partial
x_{j}^{2}}(y_{0})\phi(y_{0})$

$+\frac{1}{12}[R_{\text{a}i\text{a}k}-\underset{\text{c=1}}{\overset{\text{q}%
}{\sum}}T_{\text{ac}i}T_{\text{ac}k}-\perp_{\text{a}ik}\perp_{\text{a}%
jk}](y_{0})X_{k}(y_{0})\phi(y_{0})+\frac{1}{18}R_{ijkj}(y_{0})X_{i}%
(y_{0})X_{k}(y_{0})\phi(y_{0})$

$+\frac{1}{12}<H,j>(y_{0})X_{i}(y_{0})[X_{i}X_{j}-\frac{1}{2}\left(
\frac{\partial X_{j}}{\partial x_{i}}+\frac{\partial X_{i}}{\partial x_{j}%
}\right)  ](y_{0})\phi(y_{0})\qquad\qquad\qquad\qquad\qquad\qquad\qquad
\qquad\qquad\qquad\ \qquad\qquad\qquad\qquad\qquad\qquad\qquad\qquad
\qquad\qquad\qquad$

$-\frac{1}{6}[-R_{\text{a}i\text{b}i}+5\overset{q}{\underset{\text{c}=1}{\sum
}}T_{\text{ac}i}T_{\text{bc}i}+2\overset{n}{\underset{j=q+1}{\sum}}%
\perp_{\text{a}ij}\perp_{\text{b}ij}](y_{0})\underset{k=q+1}{\overset{n}{\sum
}}T_{\text{ab}k}(y_{0})X_{k}(y_{0})\phi(y_{0})\qquad\qquad S_{1}\qquad$

$-\frac{2}{9}\underset{j=q+1}{\overset{n}{\sum}}R_{i\text{a}ij}(y_{0}%
)[\frac{\partial X_{j}}{\partial x_{\text{a}}}%
-\underset{k=q+1}{\overset{n}{\sum}}\perp_{\text{a}jk}X_{k}](y_{0})\phi
(y_{0})$

$+\frac{1}{12}\times\frac{2}{3}\underset{j,k=q+1}{\overset{n}{\sum}}%
R_{ijik}(y_{0})[X_{j}X_{k}-\frac{1}{2}(\frac{\partial X_{j}}{\partial x_{k}%
}+\frac{\partial X_{k}}{\partial x_{j}})](y_{0})\phi(y_{0})$

$-\frac{1}{6}T_{\text{ab}i}(y_{0})\frac{\partial^{2}X_{i}}{\partial
x_{\text{a}}\partial x_{\text{b}}}(y_{0})\phi(y_{0})\qquad\qquad\qquad
S_{2}\qquad\qquad S_{21}\qquad\qquad\qquad\qquad\qquad\qquad$

$+$ $\frac{1}{12}T_{\text{ab}i}(y_{0})[$ $(R_{\text{a}i\text{b}j}%
+R_{\text{a}j\text{b}i})$ $-\underset{\text{c=1}}{\overset{\text{q}}{\sum}%
}(T_{\text{ac}i}T_{\text{bc}j}+T_{\text{ac}j}T_{\text{bc}i})$

$-\overset{n}{\underset{k=q+1}{\sum}}(\perp_{\text{a}ik}\perp_{\text{b}jk}+$
$\perp_{\text{a}jk}\perp_{\text{b}ik})](y_{0})X_{j}(y_{0})\phi(y_{0})$

$-$ $\frac{1}{6}T_{\text{ab}i}(y_{0})T_{\text{ab}j}(y_{0})[X_{i}X_{j}-\frac
{1}{2}\left(  \frac{\partial X_{i}}{\partial x_{j}}+\frac{\partial X_{j}%
}{\partial x_{i}}\right)  ](y_{0})\phi(y_{0})$

$-\frac{1}{3}\perp_{\text{a}ij}(y_{0})[(X_{i}\frac{\partial X_{j}}{\partial
x_{\text{a}}}+X_{j}\frac{\partial X_{i}}{\partial x_{\text{a}}})-\frac{1}%
{4}\left(  \frac{\partial^{2}X_{i}}{\partial x_{\text{a}}\partial x_{j}}%
+\frac{\partial^{2}X_{j}}{\partial x_{\text{a}}\partial x_{i}}\right)
](y_{0})\phi(y_{0})\qquad S_{22}$

$-\frac{1}{6}\perp_{\text{a}ij}(y_{0})[T_{\text{ab}j}\frac{\partial X_{i}%
}{\partial x_{\text{b}}}](y_{0})$

$+\frac{1}{6}\perp_{\text{a}ij}(y_{0})[(\perp_{\text{b}ik}T_{\text{ab}%
j})+\frac{2}{3}(2R_{\text{a}ijk}+R_{\text{a}jik}+R_{\text{a}kji})](y_{0}%
)X_{k}(y_{0})\phi(y_{0})$

$-\frac{1}{6}\perp_{\text{a}ij}(y_{0})\perp_{\text{a}jk}(y_{0})[X_{i}%
X_{k}-\frac{1}{2}\left(  \frac{\partial X_{i}}{\partial x_{k}}+\frac{\partial
X_{k}}{\partial x_{i}}\right)  ](y_{0})\phi(y_{0})\qquad\qquad\qquad
\qquad\qquad\qquad\qquad\qquad$

$+\frac{1}{12}[(\frac{\partial X_{j}}{\partial x_{\text{a}}})^{2}+X_{j}%
\frac{\partial^{2}X_{j}}{\partial x_{\text{a}}^{2}}-\frac{1}{2}\frac
{\partial^{3}X_{j}}{\partial x_{\text{a}}^{2}\partial x_{j}}](y_{0})\phi
(y_{0})-\frac{1}{6}\overset{n}{\underset{k=q+1}{\sum}}[\perp_{\text{b}ik}%
$T$_{\text{aa}k}\frac{\partial X_{i}}{\partial x_{\text{b}}^{2}}](y_{0}%
)\phi(y_{0})\qquad\qquad S_{3}\qquad S_{31}$

$+\frac{1}{144}[\{4\nabla_{i}R_{i\text{a}j\text{a}}+2\nabla_{j}R_{i\text{a}%
i\text{a}}+$ $8(\overset{q}{\underset{\text{c=1}}{%
{\textstyle\sum}
}}R_{\text{a}i\text{c}i}^{{}}T_{\text{ac}j}+\;\overset{n}{\underset{k=q+1}{%
{\textstyle\sum}
}}R_{\text{a}iik}\perp_{\text{a}jk})$

$+8(\overset{q}{\underset{\text{c=1}}{%
{\textstyle\sum}
}}R_{\text{a}i\text{c}j}^{{}}T_{\text{ac}i}+\;\overset{n}{\underset{k=q+1}{%
{\textstyle\sum}
}}R_{\text{a}ijk}\perp_{\text{a}ik})+8(\overset{q}{\underset{\text{c=1}}{%
{\textstyle\sum}
}}R_{\text{a}j\text{c}i}^{{}}T_{\text{ac}i}+\;\overset{n}{\underset{k=q+1}{%
{\textstyle\sum}
}}R_{\text{a}jik}\perp_{\text{a}ik})\}$\ 

$\ +\frac{2}{3}\underset{k=q+1}{\overset{n}{\sum}}\{T_{\text{aa}k}%
(R_{ijik}+3\overset{q}{\underset{\text{c}=1}{\sum}}\perp_{\text{c}ij}%
\perp_{\text{c}ik})\}](y_{0})X_{k}(y_{0})$

$-\frac{1}{12}[$ R$_{\text{a}i\text{a}k}$ $-\underset{\text{c=1}%
}{\overset{\text{q}}{\sum}}T_{\text{ac}i}T_{\text{ac}k}%
-\overset{n}{\underset{l=q+1}{\sum}}(\perp_{\text{a}il}\perp_{\text{a}%
kl}](y_{0})\times\lbrack X_{i}X_{k}-\frac{1}{2}\left(  \frac{\partial X_{i}%
}{\partial x_{k}}+\frac{\partial X_{k}}{\partial x_{i}}\right)  ](y_{0})$

$-\frac{1}{24}T_{\text{aa}k}(y_{0})[-X_{i}^{2}X_{k}+X_{k}\frac{\partial X_{i}%
}{\partial x_{i}}\ +X_{i}\left(  \frac{\partial X_{k}}{\partial x_{i}}%
+\frac{\partial X_{i}}{\partial x_{k}}\right)  -\frac{1}{3}\left(
\frac{\partial^{2}X_{k}}{\partial x_{i}^{2}}+2\frac{\partial^{2}X_{i}%
}{\partial x_{i}\partial x_{k}}\right)  ](y_{0})$

$+\frac{1}{18}[R_{\text{a}jij}\frac{\partial X_{i}}{\partial x_{\text{a}}^{2}%
}](y_{0})\qquad\qquad\qquad\qquad\qquad\qquad\qquad\qquad\qquad S_{32}$

$+\frac{1}{24}[\frac{4}{3}\overset{q}{\underset{\text{a}=1}{\sum}}%
\perp_{\text{a}ki}R_{ij\text{a}j}-\frac{1}{3}(\nabla_{i}R_{kjij}+\nabla
_{j}R_{ijik}+\nabla_{k}R_{ijij})](y_{0})X_{k}(y_{0})\phi(y_{0})$

$-\frac{1}{18}R_{ijkj}(y_{0})[X_{i}X_{k}-\frac{1}{2}\left(  \frac{\partial
X_{i}}{\partial x_{k}}+\frac{\partial X_{k}}{\partial x_{i}}\right)
](y_{0})\phi(y_{0})$

$+\frac{1}{24}[X_{i}^{2}X_{j}^{2}-2X_{i}X_{j}\left(  \frac{\partial X_{j}%
}{\partial x_{i}}+\frac{\partial X_{i}}{\partial x_{j}}\right)  -X_{i}%
^{2}\frac{\partial X_{j}}{\partial x_{j}}-X_{j}^{2}\frac{\partial X_{i}%
}{\partial x_{i}}](y_{0})\phi(y_{0})$

$+\frac{1}{48}\left(  \frac{\partial X_{j}}{\partial x_{i}}+\frac{\partial
X_{i}}{\partial x_{j}}\right)  ^{2}(y_{0})+\frac{1}{24}\left(  \frac{\partial
X_{i}}{\partial x_{i}}\frac{\partial X_{j}}{\partial x_{j}}\right)
(y_{0})\phi(y_{0})\qquad$

$+\frac{1}{36}X_{i}(y_{0})\left(  2\frac{\partial^{2}X_{j}}{\partial
x_{i}\partial x_{j}}+\frac{\partial^{2}X_{i}}{\partial x_{j}^{2}}\right)
(y_{0})+\frac{1}{36}X_{j}(y_{0})\left(  \frac{\partial^{2}X_{j}}{\partial
x_{i}^{2}}+2\frac{\partial^{2}X_{i}}{\partial x_{i}\partial x_{j}}\right)
(y_{0})\phi(y_{0})$

$-\frac{1}{48}\left(  \frac{\partial^{3}X_{i}}{\partial x_{i}\partial
x_{j}^{2}}+\frac{\partial^{3}X_{j}}{\partial x_{i}^{2}\partial x_{j}}\right)
(y_{0})\phi(y_{0})$

$+$ $\frac{2}{3}<H,j>(y_{0})\left(  \frac{\partial^{2}X_{i}}{\partial
x_{i}\partial x_{j}}+2\frac{\partial^{2}X_{j}}{\partial x_{i}^{2}}\right)
(y_{0})\phi(y_{0})+$ $\frac{2}{3}<H,j>(y_{0})R_{ijik}(y_{0})X_{k}(y_{0}%
)\phi(y_{0})\qquad$I$_{32123}$

$+\frac{1}{12}[<H,i><H,j>\ +\frac{1}{6}(2\varrho_{ij}%
+4\overset{q}{\underset{\text{a}=1}{\sum}}R_{i\text{a}j\text{a}}%
-6\overset{q}{\underset{\text{a,b}=1}{\sum}}T_{\text{aa}i}T_{\text{bb}%
j}-T_{\text{ab}i}T_{\text{ab}j})](y_{0})$

$\qquad\times\frac{1}{2}[\left(  \frac{\partial X_{j}}{\partial x_{i}}%
-\frac{\partial X_{i}}{\partial x_{j}}\right)  ](y_{0})\phi(y_{0})$

\qquad$-\frac{1}{12}\perp_{\text{a}ij}(y_{0})<H,i>(y_{0})[(X_{j}%
\perp_{\text{a}ij}-\frac{\partial X_{i}}{\partial x_{\text{a}}})+\frac
{\partial X_{\text{a}}}{\partial x_{i}}](y_{0})\phi(y_{0})$

$-\frac{1}{18}[X_{j}\left(  2\frac{\partial^{2}X_{j}}{\partial x_{i}^{2}%
}+\frac{\partial^{2}X_{i}}{\partial x_{i}\partial x_{j}}\right)  ](y_{0}%
)\phi(y_{0})\qquad\ $I$_{32124}\qquad$

$+\frac{1}{12}[2X_{i}X_{j}-\left(  \frac{\partial X_{i}}{\partial x_{j}}%
+\frac{\partial X_{j}}{\partial x_{i}}\right)  ]\frac{\partial X_{j}}{\partial
x_{i}}(y_{0})\phi(y_{0})-\frac{1}{6}[X_{i}^{2}\frac{\partial X_{j}}{\partial
x_{i}}](y_{0})\phi(y_{0})$\qquad\qquad\qquad\qquad\qquad\qquad\qquad
\qquad\qquad\qquad\qquad\qquad\qquad\qquad\qquad\qquad\qquad

$+\frac{1}{12}\frac{\partial^{2}\text{V}}{\partial x_{i}^{2}}(y_{0})\phi
(y_{0})\qquad$I$_{32125}$

\qquad\qquad\qquad\qquad\qquad\qquad\qquad\qquad\qquad\qquad\qquad\qquad
\qquad\qquad$\blacksquare$

\subsection{\textbf{Computations of Table B}$_{5}$}

\subsubsection{Normal Derivatives$\qquad\qquad\qquad$}

(i) Here we will use the other version of the definition of the \textbf{scalar
Laplacian} given by:

$\qquad\frac{1}{2}\Delta\Phi=\frac{1}{2}g^{jk}[\frac{\partial^{2}\Phi
}{\partial x_{j}\partial x_{k}}-\Gamma_{jk}^{l}\frac{\partial\Phi}{\partial
x_{l}}]$

where computations will be carried out for $j,k,l=1,...,q,q+1,...,n.$

Therefore for $i=q+1,...,n,$ we have:

$\qquad\frac{\partial}{\partial x_{i}}[\frac{1}{2}\Delta\Phi](y_{0})=\frac
{1}{2}\frac{\partial g^{jk}}{\partial x_{i}}(y_{0})[\frac{\partial^{2}\Phi
}{\partial x_{j}\partial x_{k}}-\Gamma_{jk}^{l}\frac{\partial\Phi}{\partial
x_{l}}](y_{0})+\frac{1}{2}g^{jk}(y_{0})\frac{\partial}{\partial x_{i}}%
[\frac{\partial^{2}\Phi}{\partial x_{j}\partial x_{k}}-\Gamma_{jk}^{l}%
\frac{\partial\Phi}{\partial x_{l}}](y_{0})$

\qquad\qquad\qquad\qquad$\ \ =Q_{1}+Q_{2}$\qquad\qquad\qquad$\ $

where, \qquad

$\qquad Q_{1}=\frac{1}{2}\frac{\partial g^{jk}}{\partial x_{i}}(y_{0}%
)[\frac{\partial^{2}\Phi}{\partial x_{j}\partial x_{k}}-\Gamma_{jk}^{l}%
\frac{\partial\Phi}{\partial x_{l}}](y_{0})$

$\qquad Q_{2}=\frac{1}{2}g^{jk}(y_{0})\frac{\partial}{\partial x_{i}}%
[\frac{\partial^{2}\Phi}{\partial x_{j}\partial x_{k}}-\Gamma_{jk}^{l}%
\frac{\partial\Phi}{\partial x_{l}}](y_{0})$

We recall that there is summation over repeated indices. Therefore,

\qquad$Q_{1}=\frac{1}{2}\frac{\partial g^{\text{ab}}}{\partial x_{i}}%
(y_{0})[\frac{\partial^{2}\Phi}{\partial x_{\text{a}}\partial x_{\text{b}}%
}-\Gamma_{\text{ab}}^{l}\frac{\partial\Phi}{\partial x_{l}}](y_{0})+\frac
{1}{2}\frac{\partial g^{\text{a}k}}{\partial x_{i}}(y_{0})[\frac{\partial
^{2}\Phi}{\partial x_{\text{a}}\partial x_{k}}-\Gamma_{\text{a}k}^{l}%
\frac{\partial\Phi}{\partial x_{l}}](y_{0})$

$\qquad\qquad+\frac{1}{2}\frac{\partial g^{j\text{a}}}{\partial x_{i}}%
(y_{0})[\frac{\partial^{2}\Phi}{\partial x_{j}\partial x_{\text{a}}}%
-\Gamma_{j\text{a}}^{l}\frac{\partial\Phi}{\partial x_{l}}](y_{0})+\frac{1}%
{2}\frac{\partial g^{jk}}{\partial x_{i}}(y_{0})[\frac{\partial^{2}\Phi
}{\partial x_{j}\partial x_{k}}-\Gamma_{jk}^{l}\frac{\partial\Phi}{\partial
x_{l}}](y_{0})$

The second and third terms on the RHS of the last equation above are equal.

Then we have: $\frac{\partial^{2}\Phi}{\partial x_{\text{a}}\partial
x_{\text{b}}}(y_{0})=0=\frac{\partial g^{jk}}{\partial x_{i}}(y_{0})$ for a,b
= 1,...,q and $i,j,k=q+1,...,n.$

The equation simplifies to:

$Q_{1}=\frac{\partial g^{\text{ab}}}{\partial x_{i}}(y_{0})[-\Gamma
_{\text{ab}}^{k}\frac{\partial\Phi}{\partial x_{k}}](y_{0})+2\frac{\partial
g^{\text{a}j}}{\partial x_{i}}(y_{0})[\frac{\partial^{2}\Phi}{\partial
x_{\text{a}}\partial x_{j}}-\Gamma_{\text{a}j}^{k}\frac{\partial\Phi}{\partial
x_{k}}](y_{0})$

We have for a,b = 1,...,q and $i,j,k=1,...,q+1,...,n,$

Now $\frac{\partial\text{g}^{\text{ab}}}{\partial\text{x}_{i}}(y_{0})=$
$2$T$_{\text{ab}i}(y_{0})$ by (ii) of \textbf{Table A}$_{6}.$ $\frac{\partial
g^{jk}}{\partial x_{i}}(y_{0})=0;$ $\Gamma_{\text{ab}}^{l}(y_{0}%
)=T_{\text{ab}l}(y_{0})$ by (i) of Table A$_{7};$ $\frac{\partial
g^{\text{a}j}}{\partial x_{i}}(y_{0})=\perp_{\text{a}ji}(y_{0})=-\perp
_{\text{a}ij}(y_{0})$ by (ii) of \textbf{Table A}$_{7};$ $\Gamma_{\text{a}%
j}^{k}(y_{0})=\perp_{\text{a}jk}(y_{0})$ by (iv) of \textbf{Table A}$_{8};$

$\frac{\partial\Phi}{\partial xk}=-X_{k}(y_{0})$ and $\frac{\partial^{2}\Phi
}{\partial x_{\text{a}}\partial x_{j}}(y_{0})=-\frac{\partial X_{j}}{\partial
x_{\text{a}}}(y_{0})$ by (xi) of \textbf{Table B}$_{4},$ we have:

$\qquad Q_{1}=T_{\text{ab}i}(y_{0})T_{\text{ab}k}(y_{0})X_{k}(y_{0}%
)-2\perp_{\text{a}ij}(y_{0})[-\frac{\partial X_{j}}{\partial x_{\text{a}}%
}(y_{0})-\perp_{\text{a}jk}(y_{0})(-X_{k}(y_{0})]\qquad\qquad\qquad\qquad$

$\left(  B_{100}\right)  $\qquad$Q_{1}=T_{\text{ab}i}(y_{0})T_{\text{ab}%
j}(y_{0})X_{j}(y_{0})+2\perp_{\text{a}ij}(y_{0})[\frac{\partial X_{j}%
}{\partial x_{\text{a}}}-\perp_{\text{a}jk}X_{k}](y_{0})\qquad\qquad$

$\qquad\qquad\qquad\qquad\qquad\qquad\qquad\qquad\qquad\qquad\qquad
\qquad\qquad\qquad\qquad\qquad\qquad\qquad\qquad\blacksquare\qquad\qquad
\qquad\qquad\qquad\qquad$

$\qquad Q_{2}=\frac{1}{2}g^{jk}(y_{0})\frac{\partial}{\partial x_{i}}%
[\frac{\partial^{2}\Phi}{\partial x_{j}\partial x_{k}}-\Gamma_{jk}^{l}%
\frac{\partial\Phi}{\partial x_{l}}](y_{0})=\frac{1}{2}\delta^{jk}(y_{0}%
)\frac{\partial}{\partial x_{i}}[\frac{\partial^{2}\Phi}{\partial
x_{j}\partial x_{k}}-\Gamma_{jk}^{l}\frac{\partial\Phi}{\partial x_{l}}%
](y_{0})$

$\ \ \ =\frac{1}{2}\frac{\partial}{\partial x_{i}}[\frac{\partial^{2}\Phi
}{\partial x_{j}^{2}}-\Gamma_{jj}^{l}\frac{\partial\Phi}{\partial x_{l}%
}](y_{0})$

$\ \ \ =\frac{1}{2}[\frac{\partial^{3}\Phi}{\partial x_{i}\partial x_{j}^{2}%
}-\frac{\partial\Gamma_{jj}^{l}}{\partial x_{i}}\frac{\partial\Phi}{\partial
x_{l}}+\Gamma_{jj}^{l}\frac{\partial^{2}\Phi}{\partial x_{i}\partial x_{l}%
}](y_{0})=Q_{21}+Q_{22}+Q_{23}$

where,

$\qquad Q_{21}=\frac{1}{2}\frac{\partial^{3}\Phi}{\partial x_{i}\partial
x_{j}^{2}}(y_{0});Q_{22}=\frac{1}{2}[-\frac{\partial\Gamma_{jj}^{l}}{\partial
x_{i}}\frac{\partial\Phi}{\partial x_{l}}](y_{0});Q_{23}=\frac{1}{2}%
[\Gamma_{jj}^{l}\frac{\partial^{2}\Phi}{\partial x_{i}\partial x_{l}}](y_{0})$

For a = 1,...,q and $i,j=q+1,...,n,$

$\qquad Q_{21}=\frac{1}{2}\frac{\partial^{3}\Phi}{\partial x_{i}\partial
x_{j}^{2}}(y_{0})=\frac{1}{2}\frac{\partial^{3}\Phi}{\partial x_{\text{a}%
}\partial x_{j}^{2}}(y_{0})+\frac{1}{2}\frac{\partial^{3}\Phi}{\partial
x_{i}\partial x_{\text{a}}^{2}}(y_{0})+\frac{1}{2}\frac{\partial^{3}\Phi
}{\partial x_{i}\partial x_{j}^{2}}(y_{0})$

By (xiv) of \textbf{Table B}$_{4}$ (which is B$_{88}$ above), (v) of
\textbf{Table B}$_{4};$ (xiii) of \textbf{Table B}$_{4}$ we have:

\qquad$Q_{21}=[X_{j}\frac{\partial X_{j}}{\partial x_{\text{a}}}+X_{j}%
\frac{\partial X_{j}}{\partial x_{\text{a}}}](y_{0})+[X_{j}\frac{\partial
X_{j}}{\partial x_{\text{a}}}+X_{j}\frac{\partial X_{j}}{\partial x_{\text{a}%
}}](y_{0})-\frac{1}{2}\left(  \frac{\partial^{2}X_{j}}{\partial x_{\text{a}%
}\partial x_{j}}+\frac{\partial^{2}X_{j}}{\partial x_{\text{a}}\partial x_{j}%
}\right)  (y_{0})$

$\qquad\qquad-\frac{1}{2}\frac{\partial^{2}X_{i}}{\partial x_{\text{a}}^{2}%
}(y_{0})+\frac{1}{2}X_{i}(y_{0})\left(  \frac{\partial X_{j}}{\partial x_{j}%
}-X_{j}^{2}\right)  (y_{0})+X_{j}(y_{0})\frac{\partial X_{i}}{\partial x_{j}%
}(y_{0})\qquad\ \qquad\ \qquad\qquad\qquad$

\qquad\qquad$\ +$ $\frac{1}{3}X_{k}(y_{0})[R_{jijk}](y_{0})-\frac{1}{2}%
\frac{\partial^{2}X_{i}}{\partial x_{j}^{2}}(y_{0})$

\qquad$\qquad\ Q_{21}=4X_{j}(y_{0})\frac{\partial X_{j}}{\partial x_{\text{a}%
}}(y_{0})-\frac{\partial^{2}X_{j}}{\partial x_{\text{a}}\partial x_{j}}%
(y_{0})-\frac{1}{2}\frac{\partial^{2}X_{i}}{\partial x_{\text{a}}^{2}}(y_{0})$

$\qquad\qquad+\frac{1}{2}X_{i}(y_{0})\left(  \operatorname{div}X_{M}%
-\left\Vert X\right\Vert _{M}^{2}+\left\Vert X\right\Vert _{P}^{2}%
-\operatorname{div}X_{P}-\underset{j=q+1}{\overset{n}{\sum}}<H,j>(y_{0}%
)X_{j}(y_{0})\right)  (y_{0})$

$\qquad\qquad+X_{j}(y_{0})\frac{\partial X_{i}}{\partial x_{j}}(y_{0})+$
$\frac{1}{3}X_{k}(y_{0})[R_{jijk}](y_{0})-\frac{1}{2}\frac{\partial^{2}X_{i}%
}{\partial x_{j}^{2}}(y_{0})$

$\qquad Q_{22}=\frac{1}{2}[-\frac{\partial\Gamma_{jj}^{l}}{\partial x_{i}%
}\frac{\partial\Phi}{\partial x_{l}}](y_{0})=-\frac{1}{2}\frac{\partial
\Gamma_{jj}^{l}}{\partial x_{i}}(y_{0})\frac{\partial\Phi}{\partial x_{l}%
}(y_{0})$

For a = 1,...,q and $l=q+1,...n,$ we have:

$\qquad Q_{22}=-\frac{1}{2}\frac{\partial\Gamma_{jj}^{\text{a}}}{\partial
x_{i}}(y_{0})\frac{\partial\Phi}{\partial x_{\text{a}}}(y_{0})-\frac{1}%
{2}\frac{\partial\Gamma_{jj}^{k}}{\partial x_{i}}(y_{0})\frac{\partial\Phi
}{\partial x_{l}}(y_{0})$

Since $\frac{\partial\Phi}{\partial x_{\text{a}}}(y_{0})=0,$

\qquad$Q_{22}=-\frac{1}{2}\frac{\partial\Gamma_{jj}^{k}}{\partial x_{i}}%
(y_{0})\frac{\partial\Phi}{\partial x_{k}}(y_{0})=-\frac{1}{2}\frac
{\partial\Gamma_{\text{aa}}^{k}}{\partial x_{i}}(y_{0})\frac{\partial\Phi
}{\partial x_{k}}(y_{0})-\frac{1}{2}\frac{\partial\Gamma_{jj}^{k}}{\partial
x_{i}}(y_{0})\frac{\partial\Phi}{\partial x_{k}}(y_{0})$

By (iv) of \textbf{Table A}$_{7},$ we have, $\frac{\partial\Gamma_{\text{aa}%
}^{k}}{\partial\text{x}_{i}}(y_{0})=[$ R$_{\text{a}i\text{a}k}$
$-\underset{\text{c=1}}{\overset{\text{q}}{\sum}}T_{\text{ac}i}T_{\text{ac}%
}-\perp_{\text{a}ik}\perp_{\text{a}jk}](y_{0})$

By (viii) of \textbf{Table A}$_{8},$ $\frac{\partial\Gamma_{jj}^{k}}%
{\partial\text{x}_{i}}(y_{0})=$ $\frac{2}{3}R_{ijkj}(y_{0})$

Therefore,

\qquad$Q_{22}=-\frac{1}{2}[R_{\text{a}i\text{a}k}-\underset{\text{c=1}%
}{\overset{\text{q}}{\sum}}T_{\text{ac}i}T_{\text{ac}k}-\perp_{\text{a}%
ik}\perp_{\text{a}jk}](y_{0})(-X_{k})(y_{0})-\frac{1}{2}\frac{2}{3}%
R_{ijkj}(y_{0})(-X_{k})(y_{0})$

\qquad$Q_{22}$\qquad$=\frac{1}{2}[R_{\text{a}i\text{a}k}-\underset{\text{c=1}%
}{\overset{\text{q}}{\sum}}T_{\text{ac}i}T_{\text{ac}k}-\perp_{\text{a}%
ik}\perp_{\text{a}jk}](y_{0})X_{k}(y_{0})+\frac{1}{3}R_{ijkj}(y_{0}%
)X_{k}(y_{0})\qquad\qquad$

\qquad$Q_{23}=\frac{1}{2}[\Gamma_{jj}^{l}\frac{\partial^{2}\Phi}{\partial
x_{i}\partial x_{l}}](y_{0})=\frac{1}{2}[\Gamma_{jj}^{\text{a}}\frac
{\partial^{2}\Phi}{\partial x_{i}\partial x_{l}}](y_{0})+\frac{1}{2}%
[\Gamma_{jj}^{l}\frac{\partial^{2}\Phi}{\partial x_{i}\partial x_{l}}](y_{0})$

\qquad\qquad$=\frac{1}{2}[\Gamma_{\text{bb}}^{\text{a}}\frac{\partial^{2}\Phi
}{\partial x_{i}\partial x_{l}}](y_{0})+\frac{1}{2}[\Gamma_{jj}^{\text{a}%
}\frac{\partial^{2}\Phi}{\partial x_{i}\partial x_{l}}](y_{0})+\frac{1}%
{2}[\Gamma_{\text{aa}}^{l}\frac{\partial^{2}\Phi}{\partial x_{i}\partial
x_{l}}](y_{0})+\frac{1}{2}[\Gamma_{jj}^{l}\frac{\partial^{2}\Phi}{\partial
x_{i}\partial x_{l}}](y_{0})$

Form \textbf{Tables A}$_{7}$\textbf{ and A}$_{8},$ we have for a,b = 1,...,q
and $i,j,k=q+1,...,n:$

\qquad$\Gamma_{\text{bb}}^{\text{a}}(y_{0})=0;\Gamma_{jj}^{\text{a}}%
(y_{0})=0;\Gamma_{jj}^{l}(y_{0})=0$ and $\Gamma_{\text{aa}}^{j}(y_{0})=$
T$_{\text{aa}j}(y_{0})=$ $<H,j>(y_{0}).$

Therefore,

\qquad$Q_{23}=\frac{1}{2}<H,j>(y_{0})\frac{\partial^{2}\Phi}{\partial
x_{i}\partial x_{j}}(y_{0})$

By (ii) of \textbf{Table B}$_{4}$, we have:

\qquad$Q_{23}=\frac{1}{2}<H,j>(y_{0})[X_{i}X_{j}-\frac{1}{2}\left(
\frac{\partial X_{j}}{\partial x_{i}}+\frac{\partial X_{i}}{\partial x_{j}%
}\right)  ](y_{0})$

Therefore,

\qquad$Q_{2}=Q_{21}+Q_{22}+Q_{23}\qquad\qquad\qquad\qquad\qquad\qquad
\qquad\qquad$

\qquad\qquad$=2X_{j}(y_{0})[\frac{\partial X_{j}}{\partial x_{\text{a}}}%
+\frac{\partial X_{j}}{\partial x_{\text{a}}}](y_{0})-\frac{1}{2}\left(
\frac{\partial^{2}X_{j}}{\partial x_{\text{a}}\partial x_{j}}+\frac
{\partial^{2}X_{j}}{\partial x_{\text{a}}\partial x_{j}}\right)  (y_{0})\qquad
Q_{21}\qquad\qquad\qquad$

$\qquad\qquad-\frac{1}{2}\frac{\partial^{2}X_{i}}{\partial x_{\text{a}}^{2}%
}(y_{0})+\frac{1}{2}X_{i}(y_{0})\left(  \operatorname{div}X_{M}-\left\Vert
X\right\Vert _{M}^{2}+\left\Vert X\right\Vert _{P}^{2}-\operatorname{div}%
X_{P}-\underset{j=q+1}{\overset{n}{\sum}}<H,j>(y_{0})X_{j}(y_{0})\right)
(y_{0})$

$\qquad\qquad+X_{j}(y_{0})\frac{\partial X_{i}}{\partial x_{j}}(y_{0})+$
$\frac{1}{3}X_{k}(y_{0})R_{jijk}(y_{0})-\frac{1}{2}\frac{\partial^{2}X_{i}%
}{\partial x_{j}^{2}}(y_{0})$

\qquad$\qquad+\frac{1}{2}[R_{\text{a}i\text{a}k}-\underset{\text{c=1}%
}{\overset{\text{q}}{\sum}}T_{\text{ac}i}T_{\text{ac}k}-\perp_{\text{a}%
ik}\perp_{\text{a}jk}](y_{0})X_{k}(y_{0})+\frac{1}{3}R_{ijkj}(y_{0}%
)X_{k}(y_{0})\qquad\qquad Q_{22}$

\qquad\qquad$+\frac{1}{2}<H,j>(y_{0})[X_{i}X_{j}-\frac{1}{2}\left(
\frac{\partial X_{j}}{\partial x_{i}}+\frac{\partial X_{i}}{\partial x_{j}%
}\right)  ](y_{0})\qquad\qquad\qquad\qquad Q_{23}$

We see that $\frac{1}{3}X_{k}(y_{0})R_{jijk}(y_{0})+\frac{1}{3}R_{ijkj}%
(y_{0})X_{k}(y_{0})=\frac{2}{3}X_{k}(y_{0})R_{jijk}(y_{0})$ and so,

$\left(  B_{101}\right)  \qquad Q_{2}=Q_{21}+Q_{22}+Q_{23}\qquad\qquad
\qquad\qquad\qquad\qquad\qquad\qquad$

\qquad\qquad$=2X_{j}(y_{0})[\frac{\partial X_{j}}{\partial x_{\text{a}}}%
+\frac{\partial X_{j}}{\partial x_{\text{a}}}](y_{0})-\frac{1}{2}\left(
\frac{\partial^{2}X_{j}}{\partial x_{\text{a}}\partial x_{j}}+\frac
{\partial^{2}X_{j}}{\partial x_{\text{a}}\partial x_{j}}\right)  (y_{0})\qquad
Q_{21}\qquad\qquad$

$\qquad\qquad-\frac{1}{2}\frac{\partial^{2}X_{i}}{\partial x_{\text{a}}^{2}%
}(y_{0})+\frac{1}{2}X_{i}(y_{0})\left(  \operatorname{div}X_{M}-\left\Vert
X\right\Vert _{M}^{2}+\left\Vert X\right\Vert _{P}^{2}-\operatorname{div}%
X_{P}-\underset{j=q+1}{\overset{n}{\sum}}<H,j>(y_{0})X_{j}(y_{0})\right)
(y_{0})$

$\qquad\qquad+X_{j}(y_{0})\frac{\partial X_{i}}{\partial x_{j}}(y_{0})+$
$\frac{2}{3}X_{k}(y_{0})R_{jijk}(y_{0})-\frac{1}{2}\frac{\partial^{2}X_{i}%
}{\partial x_{j}^{2}}(y_{0})$

\qquad$\qquad+\frac{1}{2}[R_{\text{a}i\text{a}k}-\underset{\text{c=1}%
}{\overset{\text{q}}{\sum}}T_{\text{ac}i}T_{\text{ac}k}-\perp_{\text{a}%
ik}\perp_{\text{a}jk}](y_{0})X_{k}(y_{0})\qquad\qquad\qquad\qquad\qquad
Q_{22}$

\qquad\qquad$+\frac{1}{2}<H,j>(y_{0})[X_{i}X_{j}-\frac{1}{2}\left(
\frac{\partial X_{j}}{\partial x_{i}}+\frac{\partial X_{i}}{\partial x_{j}%
}\right)  ](y_{0})\qquad\qquad\qquad Q_{23}$

\qquad\qquad\qquad\qquad\qquad\qquad\qquad\qquad\qquad\qquad\qquad\qquad
\qquad\qquad\qquad\qquad\qquad\qquad\qquad$\blacksquare$

Finally we have by $\left(  B_{100}\right)  $ and $\left(  B_{101}\right)  ,$

$\left(  B_{102}\right)  \qquad\frac{1}{2}\frac{\partial}{\partial x_{i}%
}[\Delta\Phi](y_{0})=Q_{1}+Q_{2}\qquad\qquad\qquad\qquad\qquad\qquad
\qquad\qquad\qquad$

\qquad$\qquad=T_{\text{ab}i}(y_{0})T_{\text{ab}j}(y_{0})X_{j}(y_{0}%
)+2\perp_{\text{a}ij}(y_{0})[\frac{\partial X_{j}}{\partial x_{\text{a}}%
}-\perp_{\text{a}jk}X_{k}](y_{0})\qquad\qquad Q_{1}$

\qquad$\qquad+$ $4X_{j}(y_{0})\frac{\partial X_{j}}{\partial x_{\text{a}}%
}(y_{0})-\frac{\partial^{2}X_{j}}{\partial x_{\text{a}}\partial x_{j}}%
(y_{0})\qquad Q_{2}$

$\qquad\qquad-\frac{1}{2}\frac{\partial^{2}X_{i}}{\partial x_{\text{a}}^{2}%
}(y_{0})+\frac{1}{2}X_{i}(y_{0})\left(  \operatorname{div}X_{M}-\left\Vert
X\right\Vert _{M}^{2}+\left\Vert X\right\Vert _{P}^{2}-\operatorname{div}%
X_{P}-<H,j>(y_{0})X_{j}\right)  (y_{0})$

$\qquad\qquad+X_{j}(y_{0})\frac{\partial X_{i}}{\partial x_{j}}(y_{0})+$
$\frac{1}{3}X_{k}(y_{0})[R_{jijk}](y_{0})-\frac{1}{2}\frac{\partial^{2}X_{i}%
}{\partial x_{j}^{2}}(y_{0})$

\qquad$\qquad+\frac{1}{2}[R_{\text{a}i\text{a}k}-\underset{\text{c=1}%
}{\overset{\text{q}}{\sum}}T_{\text{ac}i}T_{\text{ac}k}-\perp_{\text{a}%
ik}\perp_{\text{a}jk}](y_{0})X_{k}(y_{0})+\frac{1}{3}R_{ijkj}(y_{0}%
)X_{k}(y_{0})$

\qquad\qquad$+\frac{1}{2}<H,j>(y_{0})[X_{i}X_{j}-\frac{1}{2}\left(
\frac{\partial X_{j}}{\partial x_{i}}+\frac{\partial X_{i}}{\partial x_{j}%
}\right)  ](y_{0})$

\begin{center}
\qquad\qquad\qquad\qquad\qquad\qquad\qquad\qquad\qquad\qquad\qquad\qquad
\qquad\qquad\qquad\qquad\qquad\qquad\qquad$\blacksquare\qquad\qquad$
\end{center}

(ii) Here we will choose the alternative but equavalent definition of the
\textbf{scalar Laplacian} is given by:

\begin{center}
$\qquad\frac{1}{2}\Delta\Phi=\frac{1}{2}g^{jk}[\frac{\partial^{2}\Phi
}{\partial x_{j}\partial x_{k}}-\Gamma_{jk}^{l}\frac{\partial\Phi}{\partial
x_{l}}]$
\end{center}

where computations will be carried out for $j,k,l=1,...,q,q+1,...,n.$

For $i=q+1,...,n,$ we have:

$\qquad\qquad\frac{1}{2}\frac{\partial^{2}}{\partial x_{i}^{2}}[\Delta
\Phi](y_{0})=\frac{1}{2}\frac{\partial^{2}g^{jk}}{\partial x_{i}^{2}}%
(y_{0})[\frac{\partial^{2}\Phi}{\partial x_{j}\partial x_{k}}-\Gamma_{jk}%
^{l}\frac{\partial\Phi}{\partial x_{l}}](y_{0})$

$\qquad\qquad\qquad+\frac{\partial g^{jk}}{\partial x_{i}}(y_{0}%
)\frac{\partial}{\partial x_{i}}[\frac{\partial^{2}\Phi}{\partial
x_{j}\partial x_{k}}-\Gamma_{jk}^{l}\frac{\partial\Phi}{\partial x_{l}}%
](y_{0})\ \ +\frac{1}{2}g^{jk}(y_{0})\frac{\partial^{2}}{\partial x_{i}^{2}%
}[\frac{\partial^{2}\Phi}{\partial x_{j}\partial x_{k}}-\Gamma_{jk}^{l}%
\frac{\partial\Phi}{\partial x_{l}}](y_{0})\qquad\qquad\qquad\qquad
\qquad\qquad\qquad$

$\qquad\qquad\qquad\ =S_{1}+S_{2}++S_{3}$

where we set:

$\qquad S_{1}=\frac{1}{2}\frac{\partial^{2}g^{jk}}{\partial x_{i}^{2}}%
(y_{0})[\frac{\partial^{2}\Phi}{\partial x_{j}\partial x_{k}}-\Gamma_{jk}%
^{l}\frac{\partial\Phi}{\partial x_{l}}](y_{0})$

$\qquad S_{2}=\frac{\partial g^{jk}}{\partial x_{i}}(y_{0})\frac{\partial
}{\partial x_{i}}[\frac{\partial^{2}\Phi}{\partial x_{j}\partial x_{k}}%
-\Gamma_{jk}^{l}\frac{\partial\Phi}{\partial x_{l}}](y_{0})=\frac{\partial
g^{jk}}{\partial x_{i}}(y_{0})[\frac{\partial^{3}\Phi}{\partial x_{i}\partial
x_{j}\partial x_{k}}-\frac{\partial\Gamma_{jk}^{l}}{\partial x_{i}}%
\frac{\partial\Phi}{\partial x_{l}}-\Gamma_{jk}^{l}\frac{\partial^{2}\Phi
}{\partial x_{i}\partial x_{l}}](y_{0})$

$\qquad S_{3}=\frac{1}{2}g^{jk}(y_{0})\frac{\partial^{2}}{\partial x_{i}^{2}%
}[\frac{\partial^{2}\Phi}{\partial x_{j}\partial x_{k}}-\Gamma_{jk}^{l}%
\frac{\partial\Phi}{\partial x_{l}}](y_{0})$

$\ \ \ \qquad=\frac{1}{2}g^{jk}(y_{0})[\frac{\partial^{4}\Phi}{\partial
x_{i}^{2}\partial x_{j}\partial x_{k}}-\frac{\partial^{2}\Gamma_{jk}^{l}%
}{\partial x_{i}^{2}}\frac{\partial\Phi}{\partial x_{l}}-2\frac{\partial
\Gamma_{jk}^{l}}{\partial x_{i}}\frac{\partial^{2}\Phi}{\partial x_{i}\partial
x_{l}}-\Gamma_{jk}^{l}\frac{\partial^{3}\Phi}{\partial x_{i}^{2}\partial
x_{l}}](y_{0})$

$\ \qquad S_{3}=\frac{1}{2}[\frac{\partial^{4}\Phi}{\partial x_{i}^{2}\partial
x_{j}^{2}}-\frac{\partial^{2}\Gamma_{jj}^{l}}{\partial x_{i}^{2}}%
\frac{\partial\Phi}{\partial x_{l}}-2\frac{\partial\Gamma_{jj}^{l}}{\partial
x_{i}}\frac{\partial^{2}\Phi}{\partial x_{i}\partial x_{l}}-\Gamma_{jj}%
^{l}\frac{\partial^{3}\Phi}{\partial x_{i}^{2}\partial x_{l}}](y_{0})$

For a,b = 1,...,q and $i,j,k=q+1,...,n:$

\qquad$S_{1}=\frac{1}{2}\frac{\partial^{2}g^{\text{ab}}}{\partial x_{i}^{2}%
}(y_{0})[\frac{\partial^{2}\Phi}{\partial x_{\text{a}}\partial x_{\text{b}}%
}-\underset{l=1}{\overset{n}{\sum}}\Gamma_{\text{ab}}^{l}\frac{\partial\Phi
}{\partial x_{l}}](y_{0})+\frac{\partial^{2}g^{\text{a}j}}{\partial x_{i}^{2}%
}(y_{0})[\frac{\partial^{2}\Phi}{\partial x_{\text{a}}\partial x_{j}%
}-\underset{l=1}{\overset{n}{\sum}}\Gamma_{\text{a}j}^{l}\frac{\partial\Phi
}{\partial x_{l}}](y_{0})$

\qquad\qquad$+$ $\frac{1}{2}\frac{\partial^{2}g^{jk}}{\partial x_{i}^{2}%
}(y_{0})[\frac{\partial^{2}\Phi}{\partial x_{j}\partial x_{k}}%
-\underset{l=1}{\overset{n}{\sum}}\Gamma_{jk}^{l}\frac{\partial\Phi}{\partial
x_{l}}](y_{0})$

Since $\frac{\partial^{2}\Phi}{\partial x_{\text{a}}\partial x_{\text{b}}%
}(y_{0})=0=\frac{\partial\Phi}{\partial x_{\text{c}}}(y_{0}),$ we have for
a,b,c =1,...,q and $i,j,k,l=q+1,....n,$

\qquad$S_{1}=\frac{1}{2}\frac{\partial^{2}g^{\text{ab}}}{\partial x_{i}^{2}%
}(y_{0})[\underset{l=q+1}{\overset{n}{\sum}}\Gamma_{\text{ab}}^{l}%
\frac{\partial\Phi}{\partial x_{l}}](y_{0})+\frac{\partial^{2}g^{\text{a}j}%
}{\partial x_{i}^{2}}(y_{0})[\frac{\partial^{2}\Phi}{\partial x_{\text{a}%
}\partial x_{j}}-\underset{l=q+1}{\overset{n}{\sum}}\Gamma_{\text{a}j}%
^{l}\frac{\partial\Phi}{\partial x_{l}}](y_{0})$

\qquad\qquad$+$ $\frac{1}{2}\frac{\partial^{2}g^{jk}}{\partial x_{i}^{2}%
}(y_{0})[\frac{\partial^{2}\Phi}{\partial x_{j}\partial x_{k}}%
-\underset{l=q+1}{\overset{n}{\sum}}\Gamma_{jk}^{l}\frac{\partial\Phi
}{\partial x_{l}}](y_{0})$

\qquad$\ \ \ =S_{11}+S_{12}+S_{13}\qquad\qquad\qquad\qquad\qquad\qquad
\qquad\qquad\qquad\qquad$

where,

\qquad$S_{11}=\frac{1}{2}\frac{\partial^{2}g^{\text{ab}}}{\partial x_{i}^{2}%
}(y_{0})[\underset{l=q+1}{\overset{n}{\sum}}\Gamma_{\text{ab}}^{l}%
\frac{\partial\Phi}{\partial x_{l}}](y_{0})$

\qquad\ $S_{12}=\frac{\partial^{2}g^{\text{a}j}}{\partial x_{i}^{2}}%
(y_{0})[\frac{\partial^{2}\Phi}{\partial x_{\text{a}}\partial x_{j}%
}-\underset{l=q+1}{\overset{n}{\sum}}\Gamma_{\text{a}j}^{l}\frac{\partial\Phi
}{\partial x_{l}}](y_{0})$

\qquad$S_{13}=\frac{1}{2}\frac{\partial^{2}g^{jk}}{\partial x_{i}^{2}}%
(y_{0})[\frac{\partial^{2}\Phi}{\partial x_{j}\partial x_{k}}%
-\underset{l=q+1}{\overset{n}{\sum}}\Gamma_{jk}^{l}\frac{\partial\Phi
}{\partial x_{l}}](y_{0})$

By (iii) of \textbf{Table A}$_{6}$, (i) of \textbf{Table A}$_{7}$ and (i) of
\textbf{Table B}$_{4},$

\qquad$S_{11}=-[-R_{\text{a}i\text{b}i}+5\overset{q}{\underset{\text{c}%
=1}{\sum}}T_{\text{ac}i}T_{\text{bc}i}+\overset{n}{\underset{j=q+1}{\sum}%
}\perp_{\text{a}ij}\perp_{\text{b}ij}](y_{0})T_{\text{ab}l}(y_{0})X_{l}%
(y_{0})\qquad\qquad\qquad$

$\Gamma_{\text{a}j}^{k}(y_{0})=\perp_{\text{a}jk}(y_{0})$ by (iv)
\textbf{Table A}$_{8}$

By (iii) of \textbf{Table A}$_{4},$ (xi) of \textbf{Table B}$_{4}$ and (x) of
\textbf{Table A}$_{7}\qquad$

\qquad$S_{12}=$ $\frac{8}{3}R_{i\text{a}ij}(y_{0})[-\frac{\partial X_{j}%
}{\partial x_{\text{a}}}(y_{0})+\perp_{\text{a}jk}X_{k}](y_{0})=-$ $\frac
{8}{3}R_{i\text{a}ij}(y_{0})[\frac{\partial X_{j}}{\partial x_{\text{a}}%
}-\perp_{\text{a}jk}X_{k}](y_{0})$

Since $\Gamma_{jk}^{l}(y_{0})=0$ for $j,k,l=q+1,....n$ by (i) of \textbf{Table
A}$_{8},$ we have

$\qquad S_{13}=\frac{1}{2}\frac{\partial^{2}g^{jk}}{\partial x_{i}^{2}}%
(y_{0})[\frac{\partial^{2}\Phi}{\partial x_{j}\partial x_{k}}-\Gamma_{jk}%
^{l}\frac{\partial\Phi}{\partial x_{l}}](y_{0})=\frac{1}{2}\frac{\partial
^{2}g^{jk}}{\partial x_{i}^{2}}(y_{0})[\frac{\partial^{2}\Phi}{\partial
x_{j}\partial x_{k}}](y_{0})$

By (iii) of \textbf{Table A}$_{2}$ and by (ii) of \textbf{Table B}$_{4}$, we
have for a = 1,...,q and $i,j,k=q+1,....n:$

$\qquad S_{13}=\frac{2}{3}R_{ijik}(y_{0})[X_{j}X_{k}-\frac{1}{2}%
(\frac{\partial X_{j}}{\partial x_{k}}+\frac{\partial X_{k}}{\partial x_{j}%
})](y_{0})\qquad\qquad\qquad\qquad\qquad\qquad\qquad$

Then,$\qquad\qquad\qquad\qquad$

$\left(  B_{103}\right)  $\qquad$S_{1}=S_{11}+S_{12}+S_{13}\qquad\qquad
\qquad\qquad\qquad\qquad\qquad\qquad\qquad\qquad\qquad\qquad\qquad\qquad
\qquad\qquad\qquad\qquad\qquad\qquad\qquad$

$\qquad=-2[-R_{\text{a}i\text{b}i}+5\overset{q}{\underset{\text{c}=1}{\sum}%
}T_{\text{ac}i}T_{\text{bc}i}+2\overset{n}{\underset{j=q+1}{\sum}}%
\perp_{\text{a}ij}\perp_{\text{b}ij}](y_{0})\underset{k=q+1}{\overset{n}{\sum
}}T_{\text{ab}k}(y_{0})X_{k}(y_{0})\qquad$

$\qquad\ -\frac{8}{3}\underset{j=q+1}{\overset{n}{\sum}}R_{i\text{a}ij}%
(y_{0})[\frac{\partial X_{j}}{\partial x_{\text{a}}}%
-\underset{k=q+1}{\overset{n}{\sum}}\perp_{\text{a}jk}X_{k}](y_{0})+\frac
{2}{3}\underset{j,k=q+1}{\overset{n}{\sum}}R_{ijik}(y_{0})[X_{j}X_{k}-\frac
{1}{2}(\frac{\partial X_{j}}{\partial x_{k}}+\frac{\partial X_{k}}{\partial
x_{j}})](y_{0})$

$\qquad\qquad\qquad\qquad\qquad\qquad\qquad\qquad\qquad\qquad\qquad
\qquad\qquad\qquad\qquad\qquad\qquad\blacksquare\qquad\qquad\qquad$

We next compute for $i=q+1,...,n$ and $j,k,l=1,...,q,q+1,...,n$

$\qquad S_{2}=\frac{\partial g^{jk}}{\partial x_{i}}(y_{0})\frac{\partial
}{\partial x_{i}}[\frac{\partial^{2}\Phi}{\partial x_{j}\partial x_{k}}%
-\Gamma_{jk}^{l}\frac{\partial\Phi}{\partial x_{l}}](y_{0})=\frac{\partial
g^{jk}}{\partial x_{i}}(y_{0})[\frac{\partial^{3}\Phi}{\partial x_{i}\partial
x_{j}\partial x_{k}}-\frac{\partial\Gamma_{jk}^{l}}{\partial x_{i}}%
\frac{\partial\Phi}{\partial x_{l}}-\Gamma_{jk}^{l}\frac{\partial^{2}\Phi
}{\partial x_{i}\partial x_{l}}](y_{0})$

Then for a,b = 1,...,q ; $i,j.k=q+1,...,n$ and $l=1,...,q,q+1,...,n,$ we have:

$S_{2}=\frac{\partial g^{\text{ab}}}{\partial x_{i}}(y_{0})[\frac{\partial
^{3}\Phi}{\partial x_{i}\partial x_{\text{a}}\partial x_{\text{b}}}%
-\frac{\partial\Gamma_{\text{ab}}^{l}}{\partial x_{i}}\frac{\partial\Phi
}{\partial x_{l}}-\Gamma_{\text{ab}}^{l}\frac{\partial^{2}\Phi}{\partial
x_{i}\partial x_{l}}](y_{0})+2\frac{\partial g^{\text{a}j}}{\partial x_{i}%
}(y_{0})[\frac{\partial^{3}\Phi}{\partial x_{\text{a}}\partial x_{i}\partial
x_{j}}-\frac{\partial\Gamma_{\text{a}j}^{l}}{\partial x_{i}}\frac{\partial
\Phi}{\partial x_{l}}-\Gamma_{\text{a}j}^{l}\frac{\partial^{2}\Phi}{\partial
x_{i}\partial x_{l}}](y_{0})$

$+\frac{\partial g^{jk}}{\partial x_{i}}(y_{0})[\frac{\partial^{3}\Phi
}{\partial x_{i}\partial x_{j}\partial x_{k}}-\frac{\partial\Gamma_{jk}^{l}%
}{\partial x_{i}}\frac{\partial\Phi}{\partial x_{l}}-\Gamma_{jk}^{l}%
\frac{\partial^{2}\Phi}{\partial x_{i}\partial x_{l}}](y_{0})$

Since $\frac{\partial g^{jk}}{\partial x_{i}}(y_{0})=0$ for $i,j,k=q+1,...,n$
and we have:

$S_{2}=\frac{\partial g^{\text{ab}}}{\partial x_{i}}(y_{0})[\frac{\partial
^{3}\Phi}{\partial x_{i}\partial x_{\text{a}}\partial x_{\text{b}}}%
-\frac{\partial\Gamma_{\text{ab}}^{l}}{\partial x_{i}}\frac{\partial\Phi
}{\partial x_{l}}-\Gamma_{\text{ab}}^{l}\frac{\partial^{2}\Phi}{\partial
x_{i}\partial x_{l}}](y_{0})+2\frac{\partial g^{\text{a}j}}{\partial x_{i}%
}(y_{0})[\frac{\partial^{3}\Phi}{\partial x_{\text{a}}\partial x_{i}\partial
x_{j}}-\frac{\partial\Gamma_{\text{a}j}^{l}}{\partial x_{i}}\frac{\partial
\Phi}{\partial x_{l}}-\Gamma_{\text{a}j}^{l}\frac{\partial^{2}\Phi}{\partial
x_{i}\partial x_{l}}](y_{0})$

$=S_{21}+S_{22}$

where,

$S_{21}=\frac{\partial g^{\text{ab}}}{\partial x_{i}}(y_{0})[\frac
{\partial^{3}\Phi}{\partial x_{i}\partial x_{\text{a}}\partial x_{\text{b}}%
}-\frac{\partial\Gamma_{\text{ab}}^{l}}{\partial x_{i}}\frac{\partial\Phi
}{\partial x_{l}}-\Gamma_{\text{ab}}^{l}\frac{\partial^{2}\Phi}{\partial
x_{i}\partial x_{l}}](y_{0})$

$S_{22}=2\frac{\partial g^{\text{a}j}}{\partial x_{i}}(y_{0})[\frac
{\partial^{3}\Phi}{\partial x_{\text{a}}\partial x_{i}\partial x_{j}}%
-\frac{\partial\Gamma_{\text{a}j}^{l}}{\partial x_{i}}\frac{\partial\Phi
}{\partial x_{l}}-\Gamma_{\text{a}j}^{l}\frac{\partial^{2}\Phi}{\partial
x_{i}\partial x_{l}}](y_{0})$

$S_{21}=\frac{\partial g^{\text{ab}}}{\partial x_{i}}(y_{0})[\frac
{\partial^{3}\Phi}{\partial x_{i}\partial x_{\text{a}}\partial x_{\text{b}}%
}](y_{0})+\frac{\partial g^{\text{ab}}}{\partial x_{i}}(y_{0})[-\frac
{\partial\Gamma_{\text{ab}}^{l}}{\partial x_{i}}\frac{\partial\Phi}{\partial
x_{l}}-\Gamma_{\text{ab}}^{l}\frac{\partial^{2}\Phi}{\partial x_{i}\partial
x_{l}}](y_{0})$

Since $\frac{\partial\Phi}{\partial x_{\text{c}}}(y_{0})=0=\Gamma_{\text{ab}%
}^{\text{c}}(y_{0})$ for a,b,c = 1,...,q, we have for a,b = 1,...,q and
$i,j=q+1,...,n,$

$S_{21}=\frac{\partial g^{\text{ab}}}{\partial x_{i}}(y_{0})[\frac
{\partial^{3}\Phi}{\partial x_{i}\partial x_{\text{a}}\partial x_{\text{b}}%
}](y_{0})+\frac{\partial g^{\text{ab}}}{\partial x_{i}}(y_{0})[-\frac
{\partial\Gamma_{\text{ab}}^{j}}{\partial x_{i}}\frac{\partial\Phi}{\partial
x_{j}}-\Gamma_{\text{ab}}^{j}\frac{\partial^{2}\Phi}{\partial x_{i}\partial
x_{j}}](y_{0})$

$\frac{\partial g^{\text{ab}}}{\partial x_{i}}(y_{0})=2T_{\text{ab}i}(y_{0})$
by(ii) of \textbf{TableA}$_{6}.;\frac{\partial g^{\text{a}j}}{\partial x_{i}%
}(y_{0})=\perp_{\text{a}ji}(y_{0})=-\perp_{\text{a}ij}(y_{0})$ by (ii) of
\textbf{TableA}$_{4};$

$\Gamma_{\text{ab}}^{k}(y_{0})=T_{\text{ab}k}(y_{0})$ is from (i)
\textbf{Table A}$_{7};$ $\Gamma_{\text{a}j}^{\text{b}}(y_{0})=-\Gamma
_{\text{ab}}^{j}(y_{0})$

$\qquad\Gamma_{j\text{a}}^{k}(y_{0})=\Gamma_{\text{a}j}^{k}(y_{0}%
)=\perp_{\text{a}jk}(y_{0})$ by (iv) of \textbf{Table A}$_{8}$

$\qquad\frac{\partial\Gamma_{j\text{a}}^{k}}{\partial\text{x}_{i}}(y_{0})=$
$\overset{q}{\underset{\text{b}=1}{\sum}}[(\perp_{\text{b}ik}T_{\text{ab}%
j})+\frac{2}{3}(2R_{\text{a}ijk}+R_{\text{a}jik}+R_{\text{a}kji})](y_{0})$ by
(xii) of \textbf{Table A}$_{7}$

By (iii) of \textbf{Table A}$_{7},$

$\qquad\frac{\partial\Gamma_{\text{ab}}^{j}}{\partial x_{i}}(y_{0})=\frac
{1}{2}[$ $(R_{\text{a}i\text{b}j}+R_{\text{a}j\text{b}i})$
$-\underset{\text{c=1}}{\overset{\text{q}}{\sum}}(T_{\text{ac}i}T_{\text{bc}%
j}+T_{\text{ac}j}T_{\text{bc}i})$

$\qquad\qquad\qquad-\overset{n}{\underset{k=q+1}{\sum}}(\perp_{\text{a}%
ik}\perp_{\text{b}jk}+$ $\perp_{\text{a}jk}\perp_{\text{b}ik})](y_{0})$

$\qquad\frac{\partial^{2}\Phi}{\partial x_{i}\partial x_{j}}(y_{0}%
)=[X_{i}X_{j}-\frac{1}{2}\left(  \frac{\partial X_{i}}{\partial x_{j}}%
+\frac{\partial X_{j}}{\partial x_{i}}\right)  ](y_{0})$ by (ii) of
\textbf{Table B}$_{4}.$\textbf{ }

$\qquad\frac{\partial^{3}\Phi}{\partial x_{i}\partial x_{j}\partial
x_{\text{a}}}(y_{0})=2(X_{i}\frac{\partial X_{j}}{\partial x_{\text{a}}}%
+X_{j}\frac{\partial X_{i}}{\partial x_{\text{a}}})(y_{0})-\frac{1}{2}\left(
\frac{\partial^{2}X_{i}}{\partial x_{\text{a}}\partial x_{j}}+\frac
{\partial^{2}X_{j}}{\partial x_{\text{a}}\partial x_{i}}\right)  (y_{0})$ is
from (xiv) of \textbf{Table B}$_{4}$\textbf{ }

$\qquad\frac{\partial^{3}\Phi}{\partial x_{i}\partial x_{\text{a}}\partial
x_{\text{b}}}(y_{0})=-\frac{\partial^{2}X_{i}}{\partial x_{\text{a}}\partial
x_{\text{b}}}(y_{0})$ by (xii) of \textbf{Table B}$_{4}$

Therefore we have for a,b = 1,...,q and $i,j=q+1,...,n,$

$\qquad\qquad S_{21}=-2T_{\text{ab}i}(y_{0})\frac{\partial^{2}X_{i}}{\partial
x_{\text{a}}\partial x_{\text{b}}}(y_{0})\qquad\qquad\qquad\qquad\qquad
\qquad\qquad\qquad\left(  117\right)  $

$+$ $T_{\text{ab}i}(y_{0})[$ $(R_{\text{a}i\text{b}j}+R_{\text{a}j\text{b}i})$
$-\underset{\text{c=1}}{\overset{\text{q}}{\sum}}(T_{\text{ac}i}T_{\text{bc}%
j}+T_{\text{ac}j}T_{\text{bc}i})$

$-\overset{n}{\underset{k=q+1}{\sum}}(\perp_{\text{a}ik}\perp_{\text{b}jk}+$
$\perp_{\text{a}jk}\perp_{\text{b}ik})](y_{0})X_{j}(y_{0})$

$-$ $2T_{\text{ab}i}(y_{0})T_{\text{ab}j}(y_{0})[X_{i}X_{j}-\frac{1}{2}\left(
\frac{\partial X_{i}}{\partial x_{j}}+\frac{\partial X_{j}}{\partial x_{i}%
}\right)  ](y_{0})$

\qquad\qquad\qquad\qquad\qquad\qquad\qquad\qquad\qquad\qquad\qquad\qquad
\qquad\qquad\qquad\qquad$\blacksquare$

Next we have for a,b =1,...,q and $i,j,k=q+1,...,n:$

$S_{22}=2\frac{\partial g^{\text{a}j}}{\partial x_{i}}(y_{0})[\frac
{\partial^{3}\Phi}{\partial x_{\text{a}}\partial x_{i}\partial x_{j}}%
-\frac{\partial\Gamma_{\text{a}j}^{l}}{\partial x_{i}}\frac{\partial\Phi
}{\partial x_{l}}-\Gamma_{\text{a}j}^{l}\frac{\partial^{2}\Phi}{\partial
x_{i}\partial x_{l}}](y_{0})$

$=2\frac{\partial g^{\text{a}j}}{\partial x_{i}}(y_{0})[\frac{\partial^{3}%
\Phi}{\partial x_{\text{a}}\partial x_{i}\partial x_{j}}](y_{0})-2\frac
{\partial g^{\text{a}j}}{\partial x_{i}}(y_{0})[\frac{\partial\Gamma
_{\text{a}j}^{\text{b}}}{\partial x_{i}}\frac{\partial\Phi}{\partial
x_{\text{b}}}+\Gamma_{\text{a}j}^{\text{b}}\frac{\partial^{2}\Phi}{\partial
x_{i}\partial x_{\text{b}}}](y_{0})$

$-2\frac{\partial g^{\text{a}j}}{\partial x_{i}}(y_{0})[\frac{\partial
\Gamma_{\text{a}j}^{k}}{\partial x_{i}}\frac{\partial\Phi}{\partial x_{k}%
}+\Gamma_{\text{a}j}^{k}\frac{\partial^{2}\Phi}{\partial x_{i}\partial x_{k}%
}](y_{0})$

We have:

$S_{22}=-4\perp_{\text{a}ij}(y_{0})[(X_{i}\frac{\partial X_{j}}{\partial
x_{\text{a}}}+X_{j}\frac{\partial X_{i}}{\partial x_{\text{a}}})-\frac{1}%
{4}\left(  \frac{\partial^{2}X_{i}}{\partial x_{\text{a}}\partial x_{j}}%
+\frac{\partial^{2}X_{j}}{\partial x_{\text{a}}\partial x_{i}}\right)
](y_{0})\qquad\qquad\left(  118\right)  $

$\qquad-2\perp_{\text{a}ij}(y_{0})[T_{\text{ab}j}\frac{\partial X_{i}%
}{\partial x_{\text{b}}}](y_{0})$

$\qquad+2\perp_{\text{a}ij}(y_{0})[(\perp_{\text{b}ik}T_{\text{ab}j})+\frac
{2}{3}(2R_{\text{a}ijk}+R_{\text{a}jik}+R_{\text{a}kji})](y_{0})X_{k}(y_{0})$

$\qquad-2\perp_{\text{a}ij}(y_{0})\perp_{\text{a}jk}(y_{0})[X_{i}X_{k}%
-\frac{1}{2}\left(  \frac{\partial X_{i}}{\partial x_{k}}+\frac{\partial
X_{k}}{\partial x_{i}}\right)  ](y_{0})$

\qquad\qquad\qquad\qquad\qquad\qquad\qquad\qquad\qquad\qquad\qquad\qquad
\qquad\qquad\qquad\qquad$\blacksquare$

We conclude from $\left(  117\right)  $ and $\left(  118\right)  $ that:

$\left(  B_{104}\right)  \qquad\qquad S_{2}=S_{21}+S_{22}\qquad\qquad
\qquad\qquad\qquad\qquad\qquad$

$=-2T_{\text{ab}i}(y_{0})\frac{\partial^{2}X_{i}}{\partial x_{\text{a}%
}\partial x_{\text{b}}}(y_{0})\qquad\qquad\qquad S_{21}\qquad\qquad
\qquad\qquad\qquad\qquad$

$+$ $T_{\text{ab}i}(y_{0})[$ $(R_{\text{a}i\text{b}j}+R_{\text{a}j\text{b}i})$
$-\underset{\text{c=1}}{\overset{\text{q}}{\sum}}(T_{\text{ac}i}T_{\text{bc}%
j}+T_{\text{ac}j}T_{\text{bc}i})$

$-\overset{n}{\underset{k=q+1}{\sum}}(\perp_{\text{a}ik}\perp_{\text{b}jk}+$
$\perp_{\text{a}jk}\perp_{\text{b}ik})](y_{0})X_{j}(y_{0})$

$-$ $2T_{\text{ab}i}(y_{0})T_{\text{ab}j}(y_{0})[X_{i}X_{j}-\frac{1}{2}\left(
\frac{\partial X_{i}}{\partial x_{j}}+\frac{\partial X_{j}}{\partial x_{i}%
}\right)  ](y_{0})$

$-4\perp_{\text{a}ij}(y_{0})[(X_{i}\frac{\partial X_{j}}{\partial x_{\text{a}%
}}+X_{j}\frac{\partial X_{i}}{\partial x_{\text{a}}})-\frac{1}{4}\left(
\frac{\partial^{2}X_{i}}{\partial x_{\text{a}}\partial x_{j}}+\frac
{\partial^{2}X_{j}}{\partial x_{\text{a}}\partial x_{i}}\right)
](y_{0})\qquad S_{22}$

$-2\perp_{\text{a}ij}(y_{0})[T_{\text{ab}j}\frac{\partial X_{i}}{\partial
x_{\text{b}}}](y_{0})$

$+2\perp_{\text{a}ij}(y_{0})[(\perp_{\text{b}ik}T_{\text{ab}j})+\frac{2}%
{3}(2R_{\text{a}ijk}+R_{\text{a}jik}+R_{\text{a}kji})](y_{0})X_{k}(y_{0})$

$-2\perp_{\text{a}ij}(y_{0})\perp_{\text{a}jk}(y_{0})[X_{i}X_{k}-\frac{1}%
{2}\left(  \frac{\partial X_{i}}{\partial x_{k}}+\frac{\partial X_{k}%
}{\partial x_{i}}\right)  ](y_{0})$

\qquad\qquad\qquad\qquad\qquad\qquad\qquad\qquad\qquad\qquad\qquad\qquad
\qquad\qquad\qquad$\blacksquare$

We now compute the last term in the expression for $\frac{1}{2}\frac
{\partial^{2}}{\partial x_{i}^{2}}[\Delta\Phi](y_{0}):$

For $i=q+1,..,n$ and $j,k=1,...,q,q+1,..,n$

$S_{3}=\frac{1}{2}[\frac{\partial^{4}\Phi}{\partial x_{i}^{2}\partial
x_{j}^{2}}-\frac{\partial^{2}\Gamma_{jj}^{k}}{\partial x_{i}^{2}}%
\frac{\partial\Phi}{\partial x_{k}}-2\frac{\partial\Gamma_{jj}^{k}}{\partial
x_{i}}\frac{\partial^{2}\Phi}{\partial x_{i}\partial x_{k}}-\Gamma_{jj}%
^{k}\frac{\partial^{3}\Phi}{\partial x_{i}^{2}\partial x_{k}}](y_{0}%
)\qquad\qquad$

For a =1,...,q ; $i,j=q+1,...,n$ and $k=1,...,q,q+1,..,n$

$=\frac{1}{2}[\frac{\partial^{4}\Phi}{\partial x_{i}^{2}\partial x_{\text{a}%
}^{2}}-\frac{\partial^{2}\Gamma_{\text{aa}}^{k}}{\partial x_{i}^{2}}%
\frac{\partial\Phi}{\partial x_{k}}-2\frac{\partial\Gamma_{\text{aa}}^{k}%
}{\partial x_{i}}\frac{\partial^{2}\Phi}{\partial x_{i}\partial x_{k}}%
-\Gamma_{\text{aa}}^{k}\frac{\partial^{3}\Phi}{\partial x_{i}^{2}\partial
x_{k}}](y_{0})$ \qquad$\left(  119\right)  $

$+\frac{1}{2}[\frac{\partial^{4}\Phi}{\partial x_{i}^{2}\partial x_{j}^{2}%
}-\frac{\partial^{2}\Gamma_{jj}^{k}}{\partial x_{i}^{2}}\frac{\partial\Phi
}{\partial x_{k}}-2\frac{\partial\Gamma_{jj}^{k}}{\partial x_{i}}%
\frac{\partial^{2}\Phi}{\partial x_{i}\partial x_{k}}-\Gamma_{jj}^{k}%
\frac{\partial^{3}\Phi}{\partial x_{i}^{2}\partial x_{k}}](y_{0})$

$=S_{31}+S_{32}$

where,

$S_{31}=\frac{1}{2}[\frac{\partial^{4}\Phi}{\partial x_{i}^{2}\partial
x_{\text{a}}^{2}}-\frac{\partial^{2}\Gamma_{\text{aa}}^{k}}{\partial x_{i}%
^{2}}\frac{\partial\Phi}{\partial x_{k}}-2\frac{\partial\Gamma_{\text{aa}}%
^{k}}{\partial x_{i}}\frac{\partial^{2}\Phi}{\partial x_{i}\partial x_{k}%
}-\Gamma_{\text{aa}}^{k}\frac{\partial^{3}\Phi}{\partial x_{i}^{2}\partial
x_{k}}](y_{0})$

$S_{32}=\frac{1}{2}[\frac{\partial^{4}\Phi}{\partial x_{i}^{2}\partial
x_{j}^{2}}-\frac{\partial^{2}\Gamma_{jj}^{k}}{\partial x_{i}^{2}}%
\frac{\partial\Phi}{\partial x_{k}}-2\frac{\partial\Gamma_{jj}^{k}}{\partial
x_{i}}\frac{\partial^{2}\Phi}{\partial x_{i}\partial x_{k}}-\Gamma_{jj}%
^{k}\frac{\partial^{3}\Phi}{\partial x_{i}^{2}\partial x_{k}}](y_{0})$

Then for a,b = 1,...,q and $i,k=q+1,...,n,$

$\ \ \ S_{31}=\frac{1}{2}[\frac{\partial^{4}\Phi}{\partial x_{i}^{2}\partial
x_{\text{a}}^{2}}-\frac{\partial^{2}\Gamma_{\text{aa}}^{\text{b}}}{\partial
x_{i}^{2}}\frac{\partial\Phi}{\partial x_{\text{b}}}-2\frac{\partial
\Gamma_{\text{aa}}^{\text{b}}}{\partial x_{i}}\frac{\partial^{2}\Phi}{\partial
x_{i}\partial x_{\text{b}}}-\Gamma_{\text{aa}}^{\text{b}}\frac{\partial
^{3}\Phi}{\partial x_{i}^{2}\partial x_{\text{b}}}](y_{0})$

$\qquad+\frac{1}{2}[-\frac{\partial^{2}\Gamma_{\text{aa}}^{k}}{\partial
x_{i}^{2}}\frac{\partial\Phi}{\partial x_{k}}-2\frac{\partial\Gamma
_{\text{aa}}^{k}}{\partial x_{i}}\frac{\partial^{2}\Phi}{\partial
x_{i}\partial x_{k}}-\Gamma_{\text{aa}}^{k}\frac{\partial^{3}\Phi}{\partial
x_{i}^{2}\partial x_{k}}](y_{0})$

We have:

$\frac{\partial\Phi}{\partial x_{\text{b}}}(y_{0})=0;$ $\frac{\partial\Phi
}{\partial x_{k}}(y_{0})=-X_{k}(y_{0});$ $\frac{\partial^{2}\Phi}{\partial
x_{i}\partial x_{\text{b}}}(y_{0})=-\frac{\partial X_{i}}{\partial
x_{\text{b}}}(y_{0});$

$\frac{\partial^{3}\Phi_{P}}{\partial x_{i}^{2}\partial x_{k}}(y_{0}%
)=[-X_{i}^{2}X_{k}+X_{k}\frac{\partial X_{i}}{\partial x_{i}}\ +X_{i}\left(
\frac{\partial X_{k}}{\partial x_{i}}+\frac{\partial X_{i}}{\partial x_{k}%
}\right)  -\frac{1}{3}\left(  \frac{\partial^{2}X_{k}}{\partial x_{i}^{2}%
}+2\frac{\partial^{2}X_{i}}{\partial x_{i}\partial x_{k}}\right)  ](y_{0})$
from (v) of \textbf{Table B}$_{4}.$

$\frac{\partial^{3}\Phi_{P}}{\partial x_{\text{a}}^{2}\partial x_{k}}%
(y_{0})=-\frac{\partial^{2}X_{k}}{\partial x_{\text{a}}^{2}}(y_{0})$ is from
(xiii) of \textbf{Table B}$_{4};$ $\Gamma_{\text{aa}}^{k}(y_{0})=T_{\text{aa}%
k}(y_{0})$ is from (i) \textbf{Table A}$_{7}.$

\qquad$\frac{\partial^{4}\Phi_{P}}{\partial x_{\text{a}}^{2}\partial x_{j}%
^{2}}(y_{0})=$ $2[(\frac{\partial X_{j}}{\partial x_{\text{a}}})^{2}%
+X_{j}\frac{\partial^{2}X_{j}}{\partial x_{\text{a}}^{2}}](y_{0})$
$-\frac{\partial^{3}X_{j}}{\partial x_{\text{a}}^{2}\partial x_{j}}(y_{0})$ is
from (xvi) of \textbf{Table B}$_{4}.$

$\Gamma_{\text{aa}}^{\text{b}}(y_{0})=0;$ $\Gamma_{\text{aa}}^{j}%
(y_{0})=T_{\text{aa}j}(y_{0});\ \frac{\partial\Gamma_{\text{aa}}^{\text{c}}%
}{\partial\text{x}_{i}}(y_{0})=-\overset{n}{\underset{k=q+1}{\sum}}%
(\perp_{\text{c}ik})($T$_{\text{aa}k})(y_{0})$ by (v) of \ \textbf{Table
A}$_{7}.$

$\frac{\partial\Gamma_{\text{aa}}^{k}}{\partial\text{x}_{i}}(y_{0})=[$
R$_{\text{a}i\text{a}k}$ $-\underset{\text{c=1}}{\overset{\text{q}}{\sum}%
}T_{\text{ac}i}T_{\text{ac}k}-\overset{n}{\underset{l=q+1}{\sum}}%
(\perp_{\text{a}il}\perp_{\text{a}kl}](y_{0})$ by (iv) of \ \textbf{Table
A}$_{7}.$

By (vi) of Table A$_{7},$

$\frac{\partial^{2}\Gamma_{\text{aa}}^{j}}{\partial\text{x}_{i}^{2}}%
(y_{0})=+\frac{1}{6}[\{4\nabla_{i}R_{i\text{a}j\text{a}}+2\nabla
_{j}R_{i\text{a}i\text{a}}+$ $8(\overset{q}{\underset{\text{c=1}}{%
{\textstyle\sum}
}}R_{\text{a}i\text{c}i}^{{}}T_{\text{ac}j}+\;\overset{n}{\underset{k=q+1}{%
{\textstyle\sum}
}}R_{\text{a}iik}\perp_{\text{a}jk})$

\ $+8(\overset{q}{\underset{\text{c=1}}{%
{\textstyle\sum}
}}R_{\text{a}i\text{c}j}^{{}}T_{\text{ac}i}+\;\overset{n}{\underset{k=q+1}{%
{\textstyle\sum}
}}R_{\text{a}ijk}\perp_{\text{a}ik})+8(\overset{q}{\underset{\text{c=1}}{%
{\textstyle\sum}
}}R_{\text{a}j\text{c}i}^{{}}T_{\text{ac}i}+\;\overset{n}{\underset{k=q+1}{%
{\textstyle\sum}
}}R_{\text{a}jik}\perp_{\text{a}ik})\}$\ 

$+\frac{2}{3}\underset{k=q+1}{\overset{n}{\sum}}\{T_{\text{aa}k}%
(R_{ijik}+3\overset{q}{\underset{\text{c}=1}{\sum}}\perp_{\text{c}ij}%
\perp_{\text{c}ik})\}](y_{0})$

Therefore, we have the expression for $S_{31}:$

$S_{31}=\frac{1}{2}[\frac{\partial^{4}\Phi}{\partial x_{i}^{2}\partial
x_{\text{a}}^{2}}-\frac{\partial^{2}\Gamma_{\text{aa}}^{\text{b}}}{\partial
x_{i}^{2}}\frac{\partial\Phi}{\partial x_{\text{b}}}-2\frac{\partial
\Gamma_{\text{aa}}^{\text{b}}}{\partial x_{i}}\frac{\partial^{2}\Phi}{\partial
x_{i}\partial x_{\text{b}}}-\Gamma_{\text{aa}}^{\text{b}}\frac{\partial
^{3}\Phi}{\partial x_{i}^{2}\partial x_{\text{b}}}](y_{0})$

$+\frac{1}{2}[-\frac{\partial^{2}\Gamma_{\text{aa}}^{k}}{\partial x_{i}^{2}%
}\frac{\partial\Phi}{\partial x_{k}}-2\frac{\partial\Gamma_{\text{aa}}^{k}%
}{\partial x_{i}}\frac{\partial^{2}\Phi}{\partial x_{i}\partial x_{k}}%
-\Gamma_{\text{aa}}^{k}\frac{\partial^{3}\Phi}{\partial x_{i}^{2}\partial
x_{k}}](y_{0})$

From the values above,

$S_{31}=[(\frac{\partial X_{j}}{\partial x_{\text{a}}})^{2}+X_{j}%
\frac{\partial^{2}X_{j}}{\partial x_{\text{a}}^{2}}-\frac{1}{2}\frac
{\partial^{3}X_{j}}{\partial x_{\text{a}}^{2}\partial x_{j}}](y_{0}%
)-2\overset{n}{\underset{k=q+1}{\sum}}[\perp_{\text{b}ik}$T$_{\text{aa}k}%
\frac{\partial X_{i}}{\partial x_{\text{b}}^{2}}](y_{0})\qquad\left(
120\right)  $

$+\frac{1}{12}[\{4\nabla_{i}R_{i\text{a}j\text{a}}+2\nabla_{j}R_{i\text{a}%
i\text{a}}+$ $8(\overset{q}{\underset{\text{c=1}}{%
{\textstyle\sum}
}}R_{\text{a}i\text{c}i}^{{}}T_{\text{ac}j}+\;\overset{n}{\underset{k=q+1}{%
{\textstyle\sum}
}}R_{\text{a}iik}\perp_{\text{a}jk})$

$+8(\overset{q}{\underset{\text{c=1}}{%
{\textstyle\sum}
}}R_{\text{a}i\text{c}j}^{{}}T_{\text{ac}i}+\;\overset{n}{\underset{k=q+1}{%
{\textstyle\sum}
}}R_{\text{a}ijk}\perp_{\text{a}ik})+8(\overset{q}{\underset{\text{c=1}}{%
{\textstyle\sum}
}}R_{\text{a}j\text{c}i}^{{}}T_{\text{ac}i}+\;\overset{n}{\underset{k=q+1}{%
{\textstyle\sum}
}}R_{\text{a}jik}\perp_{\text{a}ik})\}$\ 

$\ +\frac{2}{3}\underset{k=q+1}{\overset{n}{\sum}}\{T_{\text{aa}k}%
(R_{ijik}+3\overset{q}{\underset{\text{c}=1}{\sum}}\perp_{\text{c}ij}%
\perp_{\text{c}ik})\}](y_{0})X_{k}(y_{0})$

$-[$ R$_{\text{a}i\text{a}k}$ $-\underset{\text{c=1}}{\overset{\text{q}}{\sum
}}T_{\text{ac}i}T_{\text{ac}k}-\overset{n}{\underset{l=q+1}{\sum}}%
(\perp_{\text{a}il}\perp_{\text{a}kl}](y_{0})\times\lbrack X_{i}X_{k}-\frac
{1}{2}\left(  \frac{\partial X_{i}}{\partial x_{k}}+\frac{\partial X_{k}%
}{\partial x_{i}}\right)  ](y_{0})$

$-\frac{1}{2}T_{\text{aa}k}(y_{0})[-X_{i}^{2}X_{k}+X_{k}\frac{\partial X_{i}%
}{\partial x_{i}}\ +X_{i}\left(  \frac{\partial X_{k}}{\partial x_{i}}%
+\frac{\partial X_{i}}{\partial x_{k}}\right)  -\frac{1}{3}\left(
\frac{\partial^{2}X_{k}}{\partial x_{i}^{2}}+2\frac{\partial^{2}X_{i}%
}{\partial x_{i}\partial x_{k}}\right)  ](y_{0})$

We next consider $S_{32}$ for $i,j=q+1,..,n$ and $k=1,...,q,q+1,..,n:$

$S_{32}=\frac{1}{2}[\frac{\partial^{4}\Phi}{\partial x_{i}^{2}\partial
x_{j}^{2}}-\frac{\partial^{2}\Gamma_{jj}^{k}}{\partial x_{i}^{2}}%
\frac{\partial\Phi}{\partial x_{k}}-2\frac{\partial\Gamma_{jj}^{k}}{\partial
x_{i}}\frac{\partial^{2}\Phi}{\partial x_{i}\partial x_{k}}-\Gamma_{jj}%
^{k}\frac{\partial^{3}\Phi}{\partial x_{i}^{2}\partial x_{k}}](y_{0})$

Then for a = 1,...,q and $i,j,k=q+1,..,n,$ we have:

$S_{32}=\frac{1}{2}[-\frac{\partial^{2}\Gamma_{jj}^{\text{a}}}{\partial
x_{i}^{2}}\frac{\partial\Phi}{\partial x_{\text{a}}}-2\frac{\partial
\Gamma_{jj}^{\text{a}}}{\partial x_{i}}\frac{\partial^{2}\Phi}{\partial
x_{i}\partial x_{\text{a}}}-\Gamma_{jj}^{\text{a}}\frac{\partial^{3}\Phi
}{\partial x_{i}^{2}\partial x_{\text{a}}}](y_{0})\qquad\qquad\qquad\left(
121\right)  $

$\qquad+\frac{1}{2}[-\frac{\partial^{2}\Gamma_{jj}^{k}}{\partial x_{i}^{2}%
}\frac{\partial\Phi}{\partial x_{k}}-2\frac{\partial\Gamma_{jj}^{k}}{\partial
x_{i}}\frac{\partial^{2}\Phi}{\partial x_{i}\partial x_{k}}-\Gamma_{jj}%
^{k}\frac{\partial^{3}\Phi}{\partial x_{i}^{2}\partial x_{k}}+\frac
{\partial^{4}\Phi}{\partial x_{i}^{2}\partial x_{j}^{2}}](y_{0})$

$\Gamma_{jj}^{k}(y_{0})=0;$ $\Gamma_{jj}^{\text{a}}(y_{0})=0$ by (ii) of
\textbf{Table A}$_{8}$ and $\frac{\partial\Gamma_{jj}^{\text{a}}}%
{\partial\text{x}_{i}}(y_{0})=\frac{2}{3}R_{\text{a}jij}(y_{0})$ by (v) of
\textbf{Table A}$_{8}.$

$\frac{\partial\Gamma_{jj}^{k}}{\partial x_{i}}(y_{0})=$ $\frac{2}{3}%
R_{ijkj}(y_{0})$ by (viii) of Table

$\frac{\partial^{2}\Gamma_{jj}^{k}}{\partial\text{x}_{i}^{2}}(y_{0})=[\frac
{4}{3}\overset{q}{\underset{\text{a}=1}{\sum}}\perp_{\text{a}ki}%
R_{ij\text{a}j}-\frac{1}{3}(\nabla_{i}R_{kjij}+\nabla_{j}R_{ijik}+\nabla
_{k}R_{ijij})](y_{0})$ by (ix) of \textbf{Table A}$_{8}$

$\frac{\partial\Phi}{\partial x_{\text{a}}}(y_{0})=0,$ $\frac{\partial\Phi
}{\partial x_{k}}(y_{0})=-X_{k}(y_{0}),$ $\frac{\partial^{2}\Phi}{\partial
x_{i}\partial x_{\text{a}}}(y_{0})=-\frac{\partial X_{i}}{\partial
x_{\text{a}}}(y_{0})$

$\frac{\partial^{2}\Phi}{\partial x_{i}\partial x_{k}}(y_{0})=[X_{i}%
X_{k}-\frac{1}{2}\left(  \frac{\partial X_{i}}{\partial x_{k}}+\frac{\partial
X_{k}}{\partial x_{i}}\right)  ](y_{0})\qquad\qquad\qquad\qquad\qquad
\qquad\qquad\qquad\left(  122\right)  $

$\frac{\partial^{3}\Phi}{\partial x_{i}^{2}\partial x_{k}}(y_{0})=[-X_{i}%
^{2}X_{k}+X_{k}\frac{\partial X_{i}}{\partial x_{i}}\ +X_{i}\left(
\frac{\partial X_{k}}{\partial x_{i}}+\frac{\partial X_{i}}{\partial x_{k}%
}\right)  -\frac{1}{3}\left(  \frac{\partial^{2}X_{k}}{\partial x_{i}^{2}%
}+2\frac{\partial^{2}X_{i}}{\partial x_{i}\partial x_{k}}\right)
](y_{0})\qquad\left(  123\right)  $

are already known. . The expression for $\frac{\partial^{4}\Phi}{\partial
x_{i}^{2}\partial x_{j}^{2}}(y_{0})$ is taken from $\left(  86\right)  $
above. We have:

$S_{32}=+\frac{2}{3}[R_{\text{a}jij}(y_{0})\frac{\partial X_{i}}{\partial
x_{\text{a}}^{2}}(y_{0})](y_{0})\qquad\qquad\qquad\qquad\qquad\qquad
\qquad\qquad\qquad\left(  124\right)  $

\qquad$+\frac{1}{2}[\frac{4}{3}\overset{q}{\underset{\text{a}=1}{\sum}}%
\perp_{\text{a}ki}R_{ij\text{a}j}-\frac{1}{3}(\nabla_{i}R_{kjij}+\nabla
_{j}R_{ijik}+\nabla_{k}R_{ijij})](y_{0})X_{k}(y_{0})$

$\ \ \ -\frac{2}{3}R_{ijkj}(y_{0})[X_{i}X_{k}-\frac{1}{2}\left(
\frac{\partial X_{i}}{\partial x_{k}}+\frac{\partial X_{k}}{\partial x_{i}%
}\right)  ](y_{0})$

\qquad$+\frac{1}{2}[X_{i}^{2}X_{j}^{2}-2X_{i}X_{j}\left(  \frac{\partial
X_{j}}{\partial x_{i}}+\frac{\partial X_{i}}{\partial x_{j}}\right)
-X_{i}^{2}\frac{\partial X_{j}}{\partial x_{j}}-X_{j}^{2}\frac{\partial X_{i}%
}{\partial x_{i}}](y_{0})$

$\qquad+\frac{1}{4}\left(  \frac{\partial X_{j}}{\partial x_{i}}%
+\frac{\partial X_{i}}{\partial x_{j}}\right)  ^{2}(y_{0})+\frac{1}{2}\left(
\frac{\partial X_{i}}{\partial x_{i}}\frac{\partial X_{j}}{\partial x_{j}%
}\right)  (y_{0})\qquad$

$\qquad+\frac{1}{3}X_{i}(y_{0})\left(  2\frac{\partial^{2}X_{j}}{\partial
x_{i}\partial x_{j}}+\frac{\partial^{2}X_{i}}{\partial x_{j}^{2}}\right)
(y_{0})+\frac{1}{3}X_{j}(y_{0})\left(  \frac{\partial^{2}X_{j}}{\partial
x_{i}^{2}}+2\frac{\partial^{2}X_{i}}{\partial x_{i}\partial x_{j}}\right)
(y_{0})$

$\qquad-\frac{1}{4}\left(  \frac{\partial^{3}X_{i}}{\partial x_{i}\partial
x_{j}^{2}}+\frac{\partial^{3}X_{j}}{\partial x_{i}^{2}\partial x_{j}}\right)
(y_{0})$

\qquad\qquad\qquad\qquad\qquad\qquad\qquad\qquad\qquad\qquad\qquad\qquad
\qquad\qquad\qquad\qquad\qquad\qquad$\blacksquare$

We now gather all the expressions that make up $S_{3}$

We have by the final expression for $S_{31}$in $\left(  122\right)  $ and
$S_{32}$ in $\left(  124\right)  :$

$\left(  B_{105}\right)  \qquad S_{3}=S_{31}+S_{32}\qquad\qquad\qquad
\qquad\qquad\qquad\qquad\qquad\qquad\qquad\left(  125\right)  $

\qquad$=+[(\frac{\partial X_{j}}{\partial x_{\text{a}}})^{2}+X_{j}%
\frac{\partial^{2}X_{j}}{\partial x_{\text{a}}^{2}}-\frac{1}{2}\frac
{\partial^{3}X_{j}}{\partial x_{\text{a}}^{2}\partial x_{j}}](y_{0}%
)-2\overset{n}{\underset{k=q+1}{\sum}}[\perp_{\text{b}ik}$T$_{\text{aa}k}%
\frac{\partial X_{i}}{\partial x_{\text{b}}^{2}}](y_{0})\qquad S_{31}$

$+\frac{1}{12}[\{4\nabla_{i}R_{i\text{a}j\text{a}}+2\nabla_{j}R_{i\text{a}%
i\text{a}}+$ $8(\overset{q}{\underset{\text{c=1}}{%
{\textstyle\sum}
}}R_{\text{a}i\text{c}i}^{{}}T_{\text{ac}j}+\;\overset{n}{\underset{k=q+1}{%
{\textstyle\sum}
}}R_{\text{a}iik}\perp_{\text{a}jk})$

$+8(\overset{q}{\underset{\text{c=1}}{%
{\textstyle\sum}
}}R_{\text{a}i\text{c}j}^{{}}T_{\text{ac}i}+\;\overset{n}{\underset{k=q+1}{%
{\textstyle\sum}
}}R_{\text{a}ijk}\perp_{\text{a}ik})+8(\overset{q}{\underset{\text{c=1}}{%
{\textstyle\sum}
}}R_{\text{a}j\text{c}i}^{{}}T_{\text{ac}i}+\;\overset{n}{\underset{k=q+1}{%
{\textstyle\sum}
}}R_{\text{a}jik}\perp_{\text{a}ik})\}$\ 

$\ +\frac{2}{3}\underset{k=q+1}{\overset{n}{\sum}}\{T_{\text{aa}k}%
(R_{ijik}+3\overset{q}{\underset{\text{c}=1}{\sum}}\perp_{\text{c}ij}%
\perp_{\text{c}ik})\}](y_{0})X_{k}(y_{0})$

$-[$ R$_{\text{a}i\text{a}k}$ $-\underset{\text{c=1}}{\overset{\text{q}}{\sum
}}T_{\text{ac}i}T_{\text{ac}k}-\overset{n}{\underset{l=q+1}{\sum}}%
(\perp_{\text{a}il}\perp_{\text{a}kl}](y_{0})\times\lbrack X_{i}X_{k}-\frac
{1}{2}\left(  \frac{\partial X_{i}}{\partial x_{k}}+\frac{\partial X_{k}%
}{\partial x_{i}}\right)  ](y_{0})$

$-\frac{1}{2}T_{\text{aa}k}(y_{0})[-X_{i}^{2}X_{k}+X_{k}\frac{\partial X_{i}%
}{\partial x_{i}}\ +X_{i}\left(  \frac{\partial X_{k}}{\partial x_{i}}%
+\frac{\partial X_{i}}{\partial x_{k}}\right)  -\frac{1}{3}\left(
\frac{\partial^{2}X_{k}}{\partial x_{i}^{2}}+2\frac{\partial^{2}X_{i}%
}{\partial x_{i}\partial x_{k}}\right)  ](y_{0})$

$+\frac{2}{3}[R_{\text{a}jij}(y_{0})\frac{\partial X_{i}}{\partial
x_{\text{a}}^{2}}(y_{0})](y_{0})\qquad\qquad\qquad\qquad\qquad\qquad
\qquad\qquad\qquad S_{32}$

$+\frac{1}{2}[\frac{4}{3}\overset{q}{\underset{\text{a}=1}{\sum}}%
\perp_{\text{a}ki}R_{ij\text{a}j}-\frac{1}{3}(\nabla_{i}R_{kjij}+\nabla
_{j}R_{ijik}+\nabla_{k}R_{ijij})](y_{0})X_{k}(y_{0})$

$-\frac{2}{3}R_{ijkj}(y_{0})[X_{i}X_{k}-\frac{1}{2}\left(  \frac{\partial
X_{i}}{\partial x_{k}}+\frac{\partial X_{k}}{\partial x_{i}}\right)  ](y_{0})$

$+\frac{1}{2}[X_{i}^{2}X_{j}^{2}-2X_{i}X_{j}\left(  \frac{\partial X_{j}%
}{\partial x_{i}}+\frac{\partial X_{i}}{\partial x_{j}}\right)  -X_{i}%
^{2}\frac{\partial X_{j}}{\partial x_{j}}-X_{j}^{2}\frac{\partial X_{i}%
}{\partial x_{i}}](y_{0})$

$+\frac{1}{4}\left(  \frac{\partial X_{j}}{\partial x_{i}}+\frac{\partial
X_{i}}{\partial x_{j}}\right)  ^{2}(y_{0})+\frac{1}{2}\left(  \frac{\partial
X_{i}}{\partial x_{i}}\frac{\partial X_{j}}{\partial x_{j}}\right)
(y_{0})\qquad$

$+\frac{1}{3}X_{i}(y_{0})\left(  2\frac{\partial^{2}X_{j}}{\partial
x_{i}\partial x_{j}}+\frac{\partial^{2}X_{i}}{\partial x_{j}^{2}}\right)
(y_{0})+\frac{1}{3}X_{j}(y_{0})\left(  \frac{\partial^{2}X_{j}}{\partial
x_{i}^{2}}+2\frac{\partial^{2}X_{i}}{\partial x_{i}\partial x_{j}}\right)
(y_{0})$

$-\frac{1}{4}\left(  \frac{\partial^{3}X_{i}}{\partial x_{i}\partial x_{j}%
^{2}}+\frac{\partial^{3}X_{j}}{\partial x_{i}^{2}\partial x_{j}}\right)
(y_{0})$

$\qquad\qquad\qquad\qquad\qquad\qquad\qquad\qquad\qquad\qquad\qquad
\qquad\qquad\qquad\blacksquare$

\begin{center}
We have by the final expressions of $S_{1},$ $S_{2}$ and $S_{3}$ respectively
from $\left(  B_{103}\right)  ,\left(  B_{104}\right)  $ and $\left(
B_{105}\right)  :\qquad\qquad\qquad\qquad\qquad\qquad\qquad\qquad\qquad
\qquad\qquad\qquad\ \ \ \qquad\qquad\qquad\qquad\qquad\qquad\qquad\qquad
\qquad\qquad\qquad$
\end{center}

$\left(  B_{106}\right)  \qquad\frac{1}{2}\frac{\partial^{2}}{\partial
x_{i}^{2}}[\Delta\Phi](y_{0})=S_{1}+S_{2}++S_{3}\qquad\qquad\qquad\qquad
\qquad\qquad\qquad\qquad\left(  126\right)  \qquad\qquad\qquad\qquad
\qquad\qquad\qquad\qquad\qquad\qquad\qquad\ \qquad\qquad\qquad\qquad
\qquad\qquad\qquad\qquad\qquad\qquad\qquad$

$\qquad=-2[-R_{\text{a}i\text{b}i}+5\overset{q}{\underset{\text{c}=1}{\sum}%
}T_{\text{ac}i}T_{\text{bc}i}+2\overset{n}{\underset{j=q+1}{\sum}}%
\perp_{\text{a}ij}\perp_{\text{b}ij}](y_{0})\underset{k=q+1}{\overset{n}{\sum
}}T_{\text{ab}k}(y_{0})X_{k}(y_{0})\qquad S_{1}\qquad$

$-\frac{8}{3}\underset{j=q+1}{\overset{n}{\sum}}R_{i\text{a}ij}(y_{0}%
)[\frac{\partial X_{j}}{\partial x_{\text{a}}}%
-\underset{k=q+1}{\overset{n}{\sum}}\perp_{\text{a}jk}X_{k}](y_{0})+\frac
{2}{3}\underset{j,k=q+1}{\overset{n}{\sum}}R_{ijik}(y_{0})[X_{j}X_{k}-\frac
{1}{2}(\frac{\partial X_{j}}{\partial x_{k}}+\frac{\partial X_{k}}{\partial
x_{j}})](y_{0})$

$-2T_{\text{ab}i}(y_{0})\frac{\partial^{2}X_{i}}{\partial x_{\text{a}}\partial
x_{\text{b}}}(y_{0})\qquad\qquad\qquad S_{2}\qquad\qquad S_{21}\qquad
\qquad\qquad\qquad\qquad\qquad$

$+$ $T_{\text{ab}i}(y_{0})[$ $(R_{\text{a}i\text{b}j}+R_{\text{a}j\text{b}i})$
$-\underset{\text{c=1}}{\overset{\text{q}}{\sum}}(T_{\text{ac}i}T_{\text{bc}%
j}+T_{\text{ac}j}T_{\text{bc}i})-\overset{n}{\underset{k=q+1}{\sum}}%
(\perp_{\text{a}ik}\perp_{\text{b}jk}+$ $\perp_{\text{a}jk}\perp_{\text{b}%
ik})](y_{0})X_{j}(y_{0})$

$-$ $2T_{\text{ab}i}(y_{0})T_{\text{ab}j}(y_{0})[X_{i}X_{j}-\frac{1}{2}\left(
\frac{\partial X_{i}}{\partial x_{j}}+\frac{\partial X_{j}}{\partial x_{i}%
}\right)  ](y_{0})$

$-4\perp_{\text{a}ij}(y_{0})[(X_{i}\frac{\partial X_{j}}{\partial x_{\text{a}%
}}+X_{j}\frac{\partial X_{i}}{\partial x_{\text{a}}})-\frac{1}{4}\left(
\frac{\partial^{2}X_{i}}{\partial x_{\text{a}}\partial x_{j}}+\frac
{\partial^{2}X_{j}}{\partial x_{\text{a}}\partial x_{i}}\right)
](y_{0})\qquad S_{22}$

$-2\perp_{\text{a}ij}(y_{0})[T_{\text{ab}j}\frac{\partial X_{i}}{\partial
x_{\text{b}}}](y_{0})$

$+2\perp_{\text{a}ij}(y_{0})[(\perp_{\text{b}ik}T_{\text{ab}j})+\frac{2}%
{3}(2R_{\text{a}ijk}+R_{\text{a}jik}+R_{\text{a}kji})](y_{0})X_{k}(y_{0})$

$-2\perp_{\text{a}ij}(y_{0})\perp_{\text{a}jk}(y_{0})[X_{i}X_{k}-\frac{1}%
{2}\left(  \frac{\partial X_{i}}{\partial x_{k}}+\frac{\partial X_{k}%
}{\partial x_{i}}\right)  ](y_{0})\qquad\qquad\qquad\qquad\qquad\qquad
\qquad\qquad$

$+[(\frac{\partial X_{j}}{\partial x_{\text{a}}})^{2}+X_{j}\frac{\partial
^{2}X_{j}}{\partial x_{\text{a}}^{2}}-\frac{1}{2}\frac{\partial^{3}X_{j}%
}{\partial x_{\text{a}}^{2}\partial x_{j}}](y_{0}%
)-2\overset{n}{\underset{k=q+1}{\sum}}[\perp_{\text{b}ik}$T$_{\text{aa}k}%
\frac{\partial X_{i}}{\partial x_{\text{b}}^{2}}](y_{0})\qquad\qquad
S_{3}\qquad S_{31}$

$+\frac{1}{12}[\{4\nabla_{i}R_{i\text{a}j\text{a}}+2\nabla_{j}R_{i\text{a}%
i\text{a}}+$ $8(\overset{q}{\underset{\text{c=1}}{%
{\textstyle\sum}
}}R_{\text{a}i\text{c}i}^{{}}T_{\text{ac}j}+\;\overset{n}{\underset{k=q+1}{%
{\textstyle\sum}
}}R_{\text{a}iik}\perp_{\text{a}jk})$

$+8(\overset{q}{\underset{\text{c=1}}{%
{\textstyle\sum}
}}R_{\text{a}i\text{c}j}^{{}}T_{\text{ac}i}+\;\overset{n}{\underset{k=q+1}{%
{\textstyle\sum}
}}R_{\text{a}ijk}\perp_{\text{a}ik})+8(\overset{q}{\underset{\text{c=1}}{%
{\textstyle\sum}
}}R_{\text{a}j\text{c}i}^{{}}T_{\text{ac}i}+\;\overset{n}{\underset{k=q+1}{%
{\textstyle\sum}
}}R_{\text{a}jik}\perp_{\text{a}ik})\}$\ 

$\ +\frac{2}{3}\underset{k=q+1}{\overset{n}{\sum}}\{T_{\text{aa}k}%
(R_{ijik}+3\overset{q}{\underset{\text{c}=1}{\sum}}\perp_{\text{c}ij}%
\perp_{\text{c}ik})\}](y_{0})X_{k}(y_{0})$

$-[$ R$_{\text{a}i\text{a}k}$ $-\underset{\text{c=1}}{\overset{\text{q}}{\sum
}}T_{\text{ac}i}T_{\text{ac}k}-\overset{n}{\underset{l=q+1}{\sum}}%
(\perp_{\text{a}il}\perp_{\text{a}kl}](y_{0})\times\lbrack X_{i}X_{k}-\frac
{1}{2}\left(  \frac{\partial X_{i}}{\partial x_{k}}+\frac{\partial X_{k}%
}{\partial x_{i}}\right)  ](y_{0})$

$-\frac{1}{2}T_{\text{aa}k}(y_{0})[-X_{i}^{2}X_{k}+X_{k}\frac{\partial X_{i}%
}{\partial x_{i}}\ +X_{i}\left(  \frac{\partial X_{k}}{\partial x_{i}}%
+\frac{\partial X_{i}}{\partial x_{k}}\right)  -\frac{1}{3}\left(
\frac{\partial^{2}X_{k}}{\partial x_{i}^{2}}+2\frac{\partial^{2}X_{i}%
}{\partial x_{i}\partial x_{k}}\right)  ](y_{0})$

$+\frac{2}{3}[R_{\text{a}jij}(y_{0})\frac{\partial X_{i}}{\partial
x_{\text{a}}^{2}}(y_{0})](y_{0})\qquad\qquad\qquad\qquad\qquad\qquad
\qquad\qquad\qquad S_{32}$

$+\frac{1}{2}[\frac{4}{3}\overset{q}{\underset{\text{a}=1}{\sum}}%
\perp_{\text{a}ki}R_{ij\text{a}j}-\frac{1}{3}(\nabla_{i}R_{kjij}+\nabla
_{j}R_{ijik}+\nabla_{k}R_{ijij})](y_{0})X_{k}(y_{0})$

$-\frac{2}{3}R_{ijkj}(y_{0})[X_{i}X_{k}-\frac{1}{2}\left(  \frac{\partial
X_{i}}{\partial x_{k}}+\frac{\partial X_{k}}{\partial x_{i}}\right)  ](y_{0})$

$+\frac{1}{2}[X_{i}^{2}X_{j}^{2}-2X_{i}X_{j}\left(  \frac{\partial X_{j}%
}{\partial x_{i}}+\frac{\partial X_{i}}{\partial x_{j}}\right)  -X_{i}%
^{2}\frac{\partial X_{j}}{\partial x_{j}}-X_{j}^{2}\frac{\partial X_{i}%
}{\partial x_{i}}](y_{0})$

$+\frac{1}{4}\left(  \frac{\partial X_{j}}{\partial x_{i}}+\frac{\partial
X_{i}}{\partial x_{j}}\right)  ^{2}(y_{0})+\frac{1}{2}\left(  \frac{\partial
X_{i}}{\partial x_{i}}\frac{\partial X_{j}}{\partial x_{j}}\right)
(y_{0})\qquad$

$+\frac{1}{3}X_{i}(y_{0})\left(  2\frac{\partial^{2}X_{j}}{\partial
x_{i}\partial x_{j}}+\frac{\partial^{2}X_{i}}{\partial x_{j}^{2}}\right)
(y_{0})+\frac{1}{3}X_{j}(y_{0})\left(  \frac{\partial^{2}X_{j}}{\partial
x_{i}^{2}}+2\frac{\partial^{2}X_{i}}{\partial x_{i}\partial x_{j}}\right)
(y_{0})$

$-\frac{1}{4}\left(  \frac{\partial^{3}X_{i}}{\partial x_{i}\partial x_{j}%
^{2}}+\frac{\partial^{3}X_{j}}{\partial x_{i}^{2}\partial x_{j}}\right)
(y_{0})$

\qquad\qquad\qquad\qquad\qquad\qquad\qquad\qquad\qquad\qquad\qquad\qquad
\qquad\qquad\qquad\qquad\qquad\qquad\qquad$\blacksquare$

\subsubsection{Tangential Derivatives}

(iii)$\qquad\frac{\partial}{\partial x_{\text{c}}}[\frac{1}{2}\Delta
\Phi](y_{0})$

$\qquad=\frac{1}{2}\frac{\partial g^{jk}}{\partial x_{\text{c}}}(y_{0}%
)[\frac{\partial^{2}\Phi}{\partial x_{j}\partial x_{k}}-\Gamma_{jk}^{l}%
\frac{\partial\Phi}{\partial x_{l}}](y_{0})+\frac{1}{2}g^{jk}(y_{0}%
)\frac{\partial}{\partial x_{\text{c}}}[\frac{\partial^{2}\Phi}{\partial
x_{j}\partial x_{k}}-\Gamma_{jk}^{l}\frac{\partial\Phi}{\partial x_{l}}%
](y_{0})$

Since $\frac{\partial g^{jk}}{\partial x_{\text{c}}}(y_{0})=0$ and
$g^{jk}(y_{0})=\delta^{jk},$ we have:

$\qquad\frac{\partial}{\partial x_{\text{c}}}[\frac{1}{2}\Delta\Phi
](y_{0})=\ \frac{1}{2}\frac{\partial}{\partial x_{\text{c}}}[\frac
{\partial^{2}\Phi}{\partial x_{j}^{2}}-\Gamma_{jj}^{l}\frac{\partial\Phi
}{\partial x_{l}}](y_{0})=\ \frac{1}{2}\frac{\partial}{\partial x_{\text{c}}%
}[\frac{\partial^{2}\Phi}{\partial x_{j}^{2}}-\Gamma_{jj}^{l}\frac
{\partial\Phi}{\partial x_{l}}](y_{0})$

Since $\frac{\partial\Gamma_{jj}^{i}}{\partial x_{\text{c}}}(y_{0})=0,$ we
have for c = 1,...,q and (changing indices) $i,j=1,...,q,q+1,...,n,$

\qquad$\frac{1}{2}\frac{\partial}{\partial x_{\text{c}}}[\Delta\Phi
](y_{0})=\ \frac{1}{2}[\frac{\partial^{3}\Phi}{\partial x_{\text{c}}\partial
x_{i}^{2}}-\Gamma_{ii}^{j}\frac{\partial^{2}\Phi}{\partial x_{\text{c}%
}\partial x_{j}}](y_{0})=$ J$_{1}+$ J$_{2}$

Since for a,c,d = 1,...,q, $\frac{\partial^{3}\Phi}{\partial x_{\text{c}%
}\partial x_{\text{a}}^{2}}(y_{0})=0=\frac{\partial^{2}\Phi}{\partial
x_{\text{c}}\partial x_{\text{d}}}(y_{0}),$ we can set,

\qquad\ J$_{1}=\frac{1}{2}[\frac{\partial^{3}\Phi}{\partial x_{\text{c}%
}\partial x_{i}^{2}}](y_{0})$ for $i=q+1,...,n$

\qquad\ J$_{2}=\frac{1}{2}[-\Gamma_{ii}^{j}\frac{\partial^{2}\Phi}{\partial
x_{\text{c}}\partial x_{j}}](y_{0})$ for $i=1,...,q,q+1,...,n$ and
$j=q+1,...,n$

By (xiv) of \textbf{Table B}$_{4}$, we have for c = 1,...,q and $i=q+1,...,n:$

$\left(  B_{107}\right)  $\qquad J$_{1}=\frac{1}{2}\frac{\partial^{3}\Phi_{P}%
}{\partial x_{\text{c}}\partial x_{i}^{2}}(x_{0})=2[X_{i}\frac{\partial X_{i}%
}{\partial x_{\text{c}}}](y_{0})-\frac{1}{2}\left(  \frac{\partial^{2}X_{i}%
}{\partial x_{\text{c}}\partial x_{i}}\right)  (y_{0})$

Then since $\Gamma_{ii}^{j}(y_{0})=0$ for $i,j=q+1,...,n,$ we have,

\ J$_{2}=\frac{1}{2}[-\Gamma_{\text{aa}}^{i}\frac{\partial^{2}\Phi}{\partial
x_{\text{c}}\partial x_{i}}](y_{0})$ for a = 1,...,q and $i=q+1,...,n:$

$\left(  B_{107}\right)  ^{\ast}\qquad$J$_{2}=-\frac{1}{2}[\Gamma_{\text{aa}%
}^{i}\frac{\partial^{2}\Phi}{\partial x_{\text{c}}\partial x_{i}}%
](y_{0})=-\frac{1}{2}[T_{\text{aa}i}(-\frac{\partial X_{i}}{\partial
x_{\text{c}}})](y_{0})=\frac{1}{2}[T_{\text{aa}i}\frac{\partial X_{i}%
}{\partial x_{\text{c}}}](y_{0})$

We we have finally here,

$\left(  B_{108}\right)  $\qquad$\frac{1}{2}\frac{\partial}{\partial
x_{\text{c}}}[\Delta\Phi](y_{0})=$ J$_{1}+$ J$_{2}\qquad\qquad\qquad
\qquad\qquad\qquad\qquad\qquad\qquad\qquad\qquad$

$\qquad\qquad=2[X_{i}\frac{\partial X_{i}}{\partial x_{\text{c}}}%
](y_{0})-\frac{1}{2}\left(  \frac{\partial^{2}X_{i}}{\partial x_{\text{c}%
}\partial x_{i}}\right)  (y_{0})+\frac{1}{2}[T_{\text{aa}i}\frac{\partial
X_{i}}{\partial x_{\text{c}}}](y_{0})\qquad\qquad\left(  127\right)  $

\qquad\qquad$=2[X_{i}\frac{\partial X_{i}}{\partial x_{\text{c}}}%
](y_{0})-\frac{1}{2}\left(  \frac{\partial^{2}X_{i}}{\partial x_{\text{c}%
}\partial x_{i}}\right)  (y_{0})+\frac{1}{2}<H,i>(y_{0})\frac{\partial X_{i}%
}{\partial x_{\text{c}}}(y_{0})$

$\qquad$

(iv) We next compute the expression for $\frac{1}{2}\frac{\partial^{2}%
}{\partial x_{\text{c}}^{2}}[\Delta\Phi]$

The procedure will be the same as in (ii) but with the advantage that
tangential differentiation of many items will vanish at the centre of Fermi
coordinates y$_{0}:$

$\qquad\frac{1}{2}\Delta\Phi=\frac{1}{2}g^{jk}[\frac{\partial^{2}\Phi
}{\partial x_{j}\partial x_{k}}-\Gamma_{jk}^{l}\frac{\partial\Phi}{\partial
x_{l}}]$ for $j,k,l=1,...,q,q+1,...,n.$

Therefore for c = 1,...,q and$,$ we have:

$\frac{\partial^{2}}{\partial x_{\text{c}}^{2}}[\frac{1}{2}\Delta\Phi
](y_{0})=\frac{1}{2}\frac{\partial^{2}g^{jk}}{\partial x_{\text{c}}^{2}}%
(y_{0})[\frac{\partial^{2}\Phi}{\partial x_{j}\partial x_{k}}-\Gamma_{jk}%
^{l}\frac{\partial\Phi}{\partial x_{l}}](y_{0})+\frac{\partial g^{jk}%
}{\partial x_{\text{c}}}(y_{0})\frac{\partial}{\partial x_{\text{c}}}%
[\frac{\partial^{2}\Phi}{\partial x_{j}\partial x_{k}}-\Gamma_{jk}^{l}%
\frac{\partial\Phi}{\partial x_{l}}](y_{0})$

$\qquad\qquad\qquad\ \ \ \qquad\ +\frac{1}{2}g^{jk}(y_{0})\frac{\partial^{2}%
}{\partial x_{\text{c}}^{2}}[\frac{\partial^{2}\Phi}{\partial x_{j}\partial
x_{k}}-\Gamma_{jk}^{l}\frac{\partial\Phi}{\partial x_{l}}](y_{0})$

Since $\frac{\partial^{2}g^{jk}}{\partial x_{\text{c}}^{2}}(y_{0}%
)=0=\frac{\partial g^{jk}}{\partial x_{\text{c}}}(y_{0})$ and $g^{jk}%
(y_{0})=\delta^{jk},$ we have:

$\qquad\frac{\partial^{2}}{\partial x_{\text{c}}^{2}}[\frac{1}{2}\Delta
\Phi](y_{0})=\ \frac{1}{2}\frac{\partial^{2}}{\partial x_{\text{c}}^{2}}%
[\frac{\partial^{2}\Phi}{\partial x_{j}^{2}}-\Gamma_{jj}^{l}\frac{\partial
\Phi}{\partial x_{l}}](y_{0})$

Since $\frac{\partial^{2}\Gamma_{jj}^{l}}{\partial x_{\text{c}}^{2}}%
(y_{0})=0=\frac{\partial\Gamma_{jj}^{l}}{\partial x_{\text{c}}}(y_{0}),$ we
have for $j,k=1,...,q,q+1,...,n$

$\frac{\partial^{2}}{\partial x_{\text{c}}^{2}}[\frac{1}{2}\Delta\Phi
](y_{0})=\ \frac{1}{2}[\frac{\partial^{4}\Phi}{\partial x_{\text{c}}%
^{2}\partial x_{j}^{2}}-\Gamma_{jj}^{k}\frac{\partial^{3}\Phi}{\partial
x_{\text{c}}^{2}\partial x_{k}}](y_{0})=$ T$_{1}+$ T$_{2}\qquad\qquad\left(
128\right)  $

where,$\qquad\qquad\qquad$

\qquad\qquad\ T$_{1}\ =\ \frac{1}{2}[\frac{\partial^{4}\Phi}{\partial
x_{\text{c}}^{2}\partial x_{j}^{2}}](y_{0})$ and T$_{2}=\ \frac{1}{2}%
[-\Gamma_{jj}^{k}\frac{\partial^{3}\Phi}{\partial x_{\text{c}}^{2}\partial
x_{k}}](y_{0})$

For a,c =1,...,q and $j=q+1,...,n,$

T$_{1}\ =\ \frac{1}{2}[\frac{\partial^{4}\Phi}{\partial x_{\text{c}}%
^{2}\partial x_{\text{a}}^{2}}](y_{0})+\ \frac{1}{2}[\frac{\partial^{4}\Phi
}{\partial x_{\text{c}}^{2}\partial x_{j}^{2}}](y_{0})$

Then since $\frac{\partial^{4}\Phi}{\partial x_{\text{c}}^{2}\partial
x_{\text{a}}^{2}}(y_{0})=0$ and $\frac{\partial^{4}\Phi}{\partial x_{\text{c}%
}^{2}\partial x_{j}^{2}}(y_{0})$ is given by (xvi) of \textbf{Table B}$_{4},$

T$_{1}\ =\ \ \frac{1}{2}\frac{\partial^{4}\Phi}{\partial x_{\text{c}}%
^{2}\partial x_{j}^{2}}(y_{0})=$ $[(\frac{\partial X_{j}}{\partial
x_{\text{a}}})^{2}+X_{j}\frac{\partial^{2}X_{j}}{\partial x_{\text{a}}^{2}%
}](y_{0})$ $-\frac{1}{2}\frac{\partial^{3}X_{j}}{\partial x_{\text{a}}%
^{2}\partial x_{j}}(y_{0})$ \qquad\qquad$\left(  129\right)  $

\qquad\qquad\qquad\qquad\qquad\qquad\qquad\qquad\qquad\qquad\qquad\qquad
\qquad\qquad\qquad\qquad$\blacksquare$

For $j,k=1,...,q,q+1,...,n,$

T$_{2}=\ \frac{1}{2}[-\Gamma_{jj}^{k}\frac{\partial^{3}\Phi}{\partial
x_{\text{c}}^{2}\partial x_{k}}](y_{0})$

Then for a,c = 1,...,q and $j=1,...,q,q+1,...,n$ and $k=q+1,...,n,$

T$_{2}=\frac{1}{2}[-\Gamma_{jj}^{\text{a}}\frac{\partial^{3}\Phi}{\partial
x_{\text{c}}^{2}\partial x_{\text{a}}}](y_{0})+\frac{1}{2}[-\Gamma_{jj}%
^{k}\frac{\partial^{3}\Phi}{\partial x_{\text{c}}^{2}\partial x_{k}}](y_{0})$

Since $\frac{\partial^{3}\Phi}{\partial x_{\text{c}}^{2}\partial x_{\text{a}}%
}(y_{0})=0,$ we have for $j=1,...,q,q+1,...,n,$

T$_{2}=\frac{1}{2}[-\Gamma_{jj}^{k}\frac{\partial^{3}\Phi}{\partial
x_{\text{c}}^{2}\partial x_{k}}](y_{0})$

Then for a,c = 1,...,q and $j,k=q+1,...,n,$ we have:

T$_{2}=-\frac{1}{2}[\Gamma_{\text{aa}}^{k}\frac{\partial^{3}\Phi}{\partial
x_{\text{c}}^{2}\partial x_{k}}](y_{0})-\frac{1}{2}[\Gamma_{jj}^{k}%
\frac{\partial^{3}\Phi}{\partial x_{\text{c}}^{2}\partial x_{k}}](y_{0})$

Now, $\Gamma_{jj}^{k}(y_{0})=0$ by (i) of Table A$_{8};$ $\Gamma_{\text{aa}%
}^{k}(y_{0})=T_{\text{aa}k}(y_{0})$ by (i) of \textbf{Table A}$_{7}.$

$\frac{\partial^{3}\Phi}{\partial x_{\text{c}}^{2}\partial x_{k}}%
(y_{0})=-\frac{\partial^{2}X_{k}}{\partial x_{\text{a}}^{2}}(y_{0})$ by (xiii)
of \textbf{Table B}$_{4}$

Therefore,

T$_{2}=-\frac{1}{2}[\Gamma_{\text{aa}}^{k}\frac{\partial^{3}\Phi}{\partial
x_{\text{c}}^{2}\partial x_{k}}](y_{0})=\frac{1}{2}[T_{\text{aa}k}%
\frac{\partial^{2}X_{k}}{\partial x_{\text{c}}^{2}}](y_{0})=\frac{1}%
{2}[<H,k>\frac{\partial^{2}X_{k}}{\partial x_{\text{c}}^{2}}](y_{0}%
)\qquad\left(  130\right)  $

Therefore by $\left(  120\right)  $ and $\left(  121\right)  ,$ we have:

$\left(  B_{109}\right)  \qquad\frac{1}{2}\frac{\partial^{2}}{\partial
x_{\text{c}}^{2}}[\Delta\Phi](y_{0})=$ $[(\frac{\partial X_{j}}{\partial
x_{\text{c}}})^{2}+X_{j}\frac{\partial^{2}X_{j}}{\partial x_{\text{c}}^{2}%
}](y_{0})$ $-\frac{1}{2}\frac{\partial^{3}X_{j}}{\partial x_{\text{c}}%
^{2}\partial x_{j}}(y_{0})+\frac{1}{2}\underset{\text{a}=1}{\overset{q}{\sum}%
}[T_{\text{aa}k}\frac{\partial^{2}X_{k}}{\partial x_{\text{c}}^{2}}%
](y_{0})\qquad\left(  131\right)  $

$\qquad\qquad\qquad=$ $[(\frac{\partial X_{i}}{\partial x_{\text{c}}}%
)^{2}+X_{i}\frac{\partial^{2}X_{i}}{\partial x_{\text{c}}^{2}}](y_{0})$
$-\frac{1}{2}\frac{\partial^{3}X_{i}}{\partial x_{\text{c}}^{2}\partial x_{i}%
}(y_{0})+\frac{1}{2}[<H,j>\frac{\partial^{2}X_{j}}{\partial x_{\text{c}}^{2}%
}](y_{0})$

\begin{center}
\qquad\qquad\qquad\qquad\qquad\qquad\qquad\qquad\qquad\qquad\qquad\qquad
\qquad\qquad\qquad$\qquad\qquad\qquad\qquad\blacksquare\qquad\qquad
\qquad\qquad\qquad$\qquad$\qquad\qquad\qquad$
\end{center}

We expand $\frac{\text{L}\Psi}{\Psi}:$

$\qquad\frac{\text{L}\Psi}{\Psi}=$ $\frac{1}{2}\frac{1}{\Psi}\Delta\Psi+$
$\frac{1}{\Psi}<\nabla\Psi,X>+$ V

\qquad$=\frac{1}{2}\theta^{\frac{1}{2}}\Delta\theta^{-\frac{1}{2}}+$ $\frac
{1}{2}\Phi^{-1}\Delta\Phi+$ $\frac{1}{\Psi}<\nabla\theta^{-\frac{1}{2}}%
,\nabla\Phi>+\frac{1}{\Psi}$ $<\nabla\Psi,X>+$ V

\qquad$=\frac{1}{2}\theta^{\frac{1}{2}}\Delta\theta^{-\frac{1}{2}}+$ $\frac
{1}{2}\Phi^{-1}\Delta\Phi+$ $<\nabla\log\theta^{-\frac{1}{2}},\nabla\log
\Phi>+$ $<\nabla\log\Psi,X>+$ V

\qquad$=\frac{1}{2}\theta^{\frac{1}{2}}\Delta\theta^{-\frac{1}{2}}+$ $\frac
{1}{2}\Phi^{-1}\Delta\Phi+$ $<\nabla\log\theta^{-\frac{1}{2}},\nabla\log
\Phi>+$ $<\nabla\log\theta^{-\frac{1}{2}},X>$

$\qquad+$ $<\nabla\log\Phi,X>+$ V

$\ \frac{\text{L}\Psi}{\Psi}=\frac{1}{2}\theta^{\frac{1}{2}}\Delta
\theta^{-\frac{1}{2}}+$ $\frac{1}{2}\Phi^{-1}\Delta\Phi+$ $<\nabla\log
\theta^{-\frac{1}{2}},\nabla\log\Phi+X>+$ $\Phi^{-1}$ $<\nabla\Phi,X>+$ V
$\qquad\left(  132\right)  \qquad$

We note that since the expansion of $\theta$ is in normal Fermi coordinates,
derivatives with respect to tangential Fermi coordinates vanish.
Paradoxically, all \textbf{purely tangential} derivatives of $\Phi$ also
vanish as seen in \textbf{Table B}$_{4},$ even though it is not expanded in
any variable. Then all derivatives, \textbf{mixed} or \textbf{purely
tangential}, of $\theta,$ $g^{ij}$ and $\Gamma_{ij}^{k}$ vanish because the
expansions of Chapter 6 are carried out in normal Fermi coordinates. We will
often use these properties here without mentioning them explicitely.

In particular, we note that $\frac{\partial}{\partial x_{\text{c}}}%
[\Delta\theta^{-\frac{1}{2}}](y_{0})=0.$

It is obvious that for any smooth functon f: M$\rightarrow R$ and any smooth
vector field X on M, we have:

$\qquad<\nabla f,X>$ $=$ $\frac{\partial f}{\partial x_{i}}X_{i}.$

Therefore from $\left(  132\right)  ,$

$\qquad\frac{\partial}{\partial x_{\text{c}}}[\frac{\text{L}\Psi}{\Psi}%
](y_{0})=\frac{1}{2}\frac{\partial}{\partial x_{\text{c}}}[\theta^{\frac{1}%
{2}}\Delta\theta^{-\frac{1}{2}}](y_{0})+$ $\frac{1}{2}\frac{\partial}{\partial
x_{\text{c}}}[\Phi^{-1}\Delta\Phi](y_{0})$

$\qquad\qquad\qquad+$ $\frac{\partial}{\partial x_{\text{c}}}[<\nabla
\log\theta^{-\frac{1}{2}},\nabla\log\Phi+X>](y_{0})+\frac{\partial}{\partial
x_{\text{c}}}[$ $\Phi^{-1}$ $<\nabla\Phi,X>](y_{0})+$ $\frac{\partial\text{V}%
}{\partial x_{\text{c}}}(y_{0})$

From the properties of tangential derivatives above, we have:

$\qquad\frac{\partial}{\partial x_{\text{c}}}[\frac{\text{L}\Psi}{\Psi}%
](y_{0})=$ $\frac{1}{2}\frac{\partial}{\partial x_{\text{c}}}[\Delta
\Phi](y_{0})+$ $\frac{\partial}{\partial x_{\text{c}}}[(\frac{\partial
}{\partial x_{i}}\log\theta^{-\frac{1}{2}})((\nabla\log\Phi)_{i}+X_{i}%
)](y_{0})$

$\qquad\qquad\qquad\qquad+\frac{\partial}{\partial x_{\text{c}}}[$
$\frac{\partial\Phi}{\partial x_{i}}X_{i}](y_{0})+$ $\frac{\partial\text{V}%
}{\partial x_{\text{c}}}(y_{0})$

\qquad\qquad\qquad\qquad$=$ $\frac{1}{2}\frac{\partial}{\partial x_{\text{c}}%
}[\Delta\Phi](y_{0})+[(\frac{\partial}{\partial x_{i}}\log\theta^{-\frac{1}%
{2}})((\frac{\partial}{\partial x_{\text{c}}}(\nabla\log\Phi)_{i}%
+\frac{\partial X_{i}}{\partial x_{\text{c}}})](y_{0})$

$\qquad\qquad\qquad\qquad+[$ $\frac{\partial^{2}\Phi}{\partial x_{\text{c}%
}\partial x_{i}}X_{i}+\frac{\partial\Phi}{\partial x_{i}}\frac{\partial X_{i}%
}{\partial x_{\text{c}}}](y_{0})+$ $\frac{\partial\text{V}}{\partial
x_{\text{c}}}(y_{0})$

Since $\frac{\partial}{\partial x_{\text{c}}}(\nabla\log\Phi)_{i}%
(y_{0})=-\frac{\partial X_{i}}{\partial x_{\text{c}}}(y_{0})$ by (ix) of
\textbf{Appendix B}$_{1},$ we see that,

$\frac{\partial}{\partial x_{\text{c}}}[\frac{\text{L}\Psi}{\Psi}%
](y_{0})=\frac{1}{2}\frac{\partial}{\partial x_{\text{c}}}[\Delta\Phi
](y_{0})+[$ $\frac{\partial^{2}\Phi}{\partial x_{\text{c}}\partial x_{i}}%
X_{i}+\frac{\partial\Phi}{\partial x_{i}}\frac{\partial X_{i}}{\partial
x_{\text{c}}}](y_{0})+$ $\frac{\partial\text{V}}{\partial x_{\text{c}}}%
(y_{0})$

Now $\frac{\partial\Phi}{\partial x_{i}}(y_{0})=-X_{i}(y_{0})$ and
$\frac{\partial^{2}\Phi}{\partial x_{\text{c}}\partial x_{i}}(y_{0}%
)=-\frac{\partial X_{i}}{\partial x_{\text{c}}}(y_{0}).$

Therefore,

\qquad$\frac{\partial}{\partial x_{\text{c}}}[\frac{\text{L}\Psi}{\Psi}%
](y_{0})=$ $\frac{1}{2}\frac{\partial}{\partial x_{\text{c}}}[\Delta
\Phi](y_{0})-[\frac{\partial X_{i}}{\partial x_{\text{c}}}X_{i}+X_{i}%
\frac{\partial X_{i}}{\partial x_{\text{c}}}](y_{0})+$ $\frac{\partial
\text{V}}{\partial x_{\text{c}}}(y_{0})\qquad\qquad\left(  133\right)
\qquad\qquad\qquad\qquad\qquad\qquad\qquad\qquad\qquad\ \ \ $

We see from $\left(  127\right)  $ given in $\left(  B_{108}\right)  $ or in
(iii) of \textbf{Table B}$_{5}$ above that,

$\qquad\frac{1}{2}\frac{\partial}{\partial x_{\text{c}}}[\Delta\Phi
](y_{0})=2[X_{j}\frac{\partial X_{j}}{\partial x_{\text{c}}}](y_{0})-\frac
{1}{2}\left(  \frac{\partial^{2}X_{j}}{\partial x_{\text{c}}\partial x_{j}%
}\right)  (y_{0})+\frac{1}{2}[T_{\text{aa}i}\frac{\partial X_{i}}{\partial
x_{\text{c}}}](y_{0})$

\qquad\qquad\qquad\qquad$\ \ =2[X_{j}\frac{\partial X_{j}}{\partial
x_{\text{c}}}](y_{0})-\frac{1}{2}\left(  \frac{\partial^{2}X_{j}}{\partial
x_{\text{c}}\partial x_{j}}\right)  (y_{0})+\frac{1}{2}<H,i>(y_{0}%
)\frac{\partial X_{i}}{\partial x_{\text{c}}}(y_{0})$

Therefore,

$\left(  B_{110}\right)  \qquad\frac{\partial}{\partial x_{\text{c}}}%
[\frac{\text{L}\Psi}{\Psi}](y_{0})=2\left(  X_{i}\frac{\partial X_{i}%
}{\partial x_{\text{c}}}\right)  (y_{0})-\frac{1}{2}\left(  \frac{\partial
^{2}X_{i}}{\partial x_{\text{c}}\partial x_{i}}\right)  (y_{0})\qquad\left(
134\right)  $

$\qquad\qquad\qquad\qquad+\frac{1}{2}<H,i>(y_{0})\frac{\partial X_{i}%
}{\partial x_{\text{c}}}(y_{0})-2[X_{i}\frac{\partial X_{i}}{\partial
x_{\text{c}}}](y_{0})+$ $\frac{\partial\text{V}}{\partial x_{\text{c}}}%
(y_{0})$

We simplify and have:

$\left(  B_{111}\right)  \qquad\frac{\partial}{\partial x_{\text{c}}}%
[\frac{\text{L}\Psi}{\Psi}](y_{0})=-\frac{1}{2}\left(  \frac{\partial^{2}%
X_{i}}{\partial x_{\text{c}}\partial x_{i}}\right)  (y_{0})+\frac{1}%
{2}[<H,i>\frac{\partial X_{i}}{\partial x_{\text{c}}}](y_{0})+$ $\frac
{\partial\text{V}}{\partial x_{\text{c}}}(y_{0})\qquad\left(  135\right)
\qquad\qquad\qquad\qquad\qquad\qquad\qquad\qquad\qquad$

\qquad\qquad\qquad\qquad\qquad\qquad\qquad\qquad\qquad\qquad\qquad
$\qquad\qquad\qquad\qquad\qquad\qquad\qquad\blacksquare$

(vi)\qquad We next compute $\frac{\partial^{2}}{\partial x_{\text{c}}^{2}%
}[\frac{\text{L}\Psi}{\Psi}](y_{0}).$ We have from $\left(  132\right)  :$

$\ \frac{\text{L}\Psi}{\Psi}=\frac{1}{2}\theta^{\frac{1}{2}}\Delta
\theta^{-\frac{1}{2}}+$ $\frac{1}{2}\Phi^{-1}\Delta\Phi+$ $<\nabla\log
\theta^{-\frac{1}{2}},\nabla\log\Phi+X>+$ $\Phi^{-1}$ $<\nabla\Phi,X>+$ V

Therefore,

\qquad$\frac{\partial^{2}}{\partial x_{\text{c}}^{2}}[\frac{\text{L}\Psi}%
{\Psi}](y_{0})=\frac{1}{2}\frac{\partial^{2}}{\partial x_{\text{c}}^{2}%
}[\theta^{\frac{1}{2}}\Delta\theta^{-\frac{1}{2}}](y_{0})+$ $\frac{1}{2}%
\frac{\partial^{2}}{\partial x_{\text{c}}^{2}}[\Phi^{-1}\Delta\Phi](y_{0})$

$\qquad+$ $\frac{\partial^{2}}{\partial x_{\text{c}}^{2}}[(\frac{\partial
}{\partial x_{i}}\log\theta^{-\frac{1}{2}})((\nabla\log\Phi)_{i}+X_{i}%
)](y_{0})+\frac{\partial^{2}}{\partial x_{\text{c}}^{2}}[\Phi^{-1}%
\frac{\partial\Phi}{\partial x_{i}}X_{i}](y_{0})+$ $\frac{\partial^{2}%
\text{V}}{\partial x_{\text{c}}^{2}}(y_{0})$

\qquad Tangential derivatives of $\theta^{\frac{1}{2}}\Delta\theta^{-\frac
{1}{2}}$ and of $\frac{\partial}{\partial x_{i}}\log\theta^{-\frac{1}{2}}$
vanish at $y_{0}$. Pure tangential derivatives of $\Phi$ also vanish at
$y_{0}$. Consquently, the last equation above simplifies to:

$\frac{\partial^{2}}{\partial x_{\text{c}}^{2}}[\frac{\text{L}\Psi}{\Psi
}](y_{0})=$ $\frac{1}{2}\frac{\partial^{2}}{\partial x_{\text{c}}^{2}}%
[\Delta\Phi](y_{0})+\frac{\partial}{\partial x_{i}}\log\theta^{-\frac{1}{2}%
}(y_{0})[\frac{\partial^{2}}{\partial x_{\text{c}}^{2}}(\nabla\log\Phi
)_{i}+\frac{\partial^{2}X_{i}}{\partial x_{\text{c}}^{2}})](y_{0})$

$\qquad\qquad\qquad\qquad+\frac{\partial^{2}}{\partial x_{\text{c}}^{2}}%
[\frac{\partial\Phi}{\partial x_{i}}X_{i}](y_{0})+$ $\frac{\partial
^{2}\text{V}}{\partial x_{\text{c}}^{2}}(y_{0})$

Since $\frac{\partial}{\partial x_{i}}\log\theta^{-\frac{1}{2}}(y_{0})=0$ or,
alternatively, $\frac{\partial^{2}}{\partial x_{\text{c}}^{2}}(\nabla$%
log$\Phi_{P})_{i}(y_{0})=-\frac{\partial^{2}X_{i}}{\partial x_{\text{c}}^{2}%
}(y_{0})$ by (xii) of \textbf{Table B}$_{4},$ we have:

$\frac{\partial^{2}}{\partial x_{\text{c}}^{2}}[\frac{\text{L}\Psi}{\Psi
}](y_{0})=$ $\frac{1}{2}\frac{\partial^{2}}{\partial x_{\text{c}}^{2}}%
[\Delta\Phi](y_{0})+[\frac{\partial^{3}\Phi}{\partial x_{\text{c}}^{2}\partial
x_{i}}X_{i}+\frac{\partial\Phi}{\partial x_{i}}\frac{\partial^{2}X_{i}%
}{\partial x_{\text{c}}^{2}}+2\frac{\partial^{2}\Phi}{\partial x_{i}\partial
x_{\text{c}}}\frac{\partial X_{i}}{\partial x_{\text{c}}}](y_{0})+$
$\frac{\partial^{2}\text{V}}{\partial x_{\text{c}}^{2}}(y_{0})$

$\frac{\partial\Phi}{\partial x_{i}}(y_{0})=-X_{i}(y_{0}),$ $\frac
{\partial^{2}\Phi}{\partial x_{\text{c}}\partial x_{i}}(y_{0})=-\frac{\partial
X_{i}}{\partial x_{\text{c}}}(y_{0})$ by (xi) of \textbf{Table B}$_{4}$ and,

$\frac{\partial^{3}\Phi_{P}}{\partial x_{\text{c}}^{2}\partial x_{i}}%
(y_{0})=-\frac{\partial^{2}X_{i}}{\partial x_{\text{a}}^{2}}(y_{0})$ by (xiii)
of \textbf{Table B}$_{4}.$ Therefore,

\qquad$\frac{\partial^{2}}{\partial x_{\text{c}}^{2}}[\frac{\text{L}\Psi}%
{\Psi}](y_{0})$

$\qquad=$ $\frac{1}{2}\frac{\partial^{2}}{\partial x_{\text{c}}^{2}}%
[\Delta\Phi](y_{0})+[(-\frac{\partial^{2}X_{i}}{\partial x_{\text{c}}^{2}%
})X_{i}+(-X_{i})\frac{\partial^{2}X_{i}}{\partial x_{\text{c}}^{2}}%
+2(-\frac{\partial X_{i}}{\partial x_{\text{c}}})\frac{\partial X_{i}%
}{\partial x_{\text{c}}}](y_{0})+$ $\frac{\partial^{2}\text{V}}{\partial
x_{\text{c}}^{2}}(y_{0})\qquad\left(  136\right)  $

\qquad$=$ $\frac{1}{2}\frac{\partial^{2}}{\partial x_{\text{c}}^{2}}%
[\Delta\Phi](y_{0})-2[X_{i}\frac{\partial^{2}X_{i}}{\partial x_{\text{c}}^{2}%
}+(\frac{\partial X_{i}}{\partial x_{\text{c}}})^{2}](y_{0})+$ $\frac
{\partial^{2}\text{V}}{\partial x_{\text{c}}^{2}}(y_{0})$

From $\left(  B_{107}\right)  $ in $\left(  122\right)  ,$ we have:

$\frac{1}{2}\frac{\partial^{2}}{\partial x_{\text{c}}^{2}}[\Delta\Phi
](y_{0})=$ $[(\frac{\partial X_{i}}{\partial x_{\text{c}}})^{2}+X_{i}%
\frac{\partial^{2}X_{i}}{\partial x_{\text{c}}^{2}}](y_{0})$ $-\frac{1}%
{2}\frac{\partial^{3}X_{i}}{\partial x_{\text{c}}^{2}\partial x_{i}}%
(y_{0})+\frac{1}{2}\underset{\text{a}=1}{\overset{q}{\sum}}[T_{\text{aa}%
j}\frac{\partial^{2}X_{j}}{\partial x_{\text{c}}^{2}}](y_{0})$

Consequently, we have:

$\frac{\partial^{2}}{\partial x_{\text{c}}^{2}}[\frac{\text{L}\Psi}{\Psi
}](y_{0})=$ $[(\frac{\partial X_{i}}{\partial x_{\text{c}}})^{2}+X_{i}%
\frac{\partial^{2}X_{i}}{\partial x_{\text{c}}^{2}}](y_{0})$ $-\frac{1}%
{2}\frac{\partial^{3}X_{i}}{\partial x_{\text{c}}^{2}\partial x_{i}}%
(y_{0})+\frac{1}{2}[<H,j>\frac{\partial^{2}X_{j}}{\partial x_{\text{c}}^{2}%
}](y_{0})$

$-2[X_{i}\frac{\partial^{2}X_{i}}{\partial x_{\text{c}}^{2}}+(\frac{\partial
X_{i}}{\partial x_{\text{c}}})^{2}](y_{0})+$ $\frac{\partial^{2}\text{V}%
}{\partial x_{\text{c}}^{2}}(y_{0})$

The last expression above simplifies to:

$\left(  B_{112}\right)  \qquad\frac{\partial^{2}}{\partial x_{\text{c}}^{2}%
}[\frac{\text{L}\Psi}{\Psi}](y_{0})=$ $-$ $[\frac{1}{2}\frac{\partial^{3}%
X_{i}}{\partial x_{\text{c}}^{2}\partial x_{i}}+(\frac{\partial X_{i}%
}{\partial x_{\text{c}}})^{2}+X_{i}\frac{\partial^{2}X_{i}}{\partial
x_{\text{c}}^{2}}](y_{0})\qquad\left(  137\right)  $

\qquad\qquad\qquad\qquad\qquad\qquad$+\frac{1}{2}[<H,j>\frac{\partial^{2}%
X_{j}}{\partial x_{\text{c}}^{2}}](y_{0})+$ $\frac{\partial^{2}\text{V}%
}{\partial x_{\text{c}}^{2}}(y_{0})$

\subsubsection{Further Normal Derivatives}

(vii) We next compute:

$R=\frac{1}{12}$ $\frac{\partial}{\partial x_{i}}\left(  \frac{\text{L}\Psi
}{\Psi}\right)  (y_{0})\phi(y_{0})$ for $i=q+1,...,n:$

$\frac{\text{L}\Psi}{\Psi}=\frac{1}{2}\theta^{\frac{1}{2}}\Delta\theta
^{-\frac{1}{2}}+$ $\frac{1}{2}\Phi^{-1}\Delta\Phi+$ $<\nabla\log\theta
^{-\frac{1}{2}},\nabla\log\Phi+X>+$ $<\nabla\log\Phi,X>+$ V

$=R_{1}+R_{2}+R_{3}+R_{4}+R_{4}$

where,

$\qquad R_{1}=\frac{1}{24}\frac{\partial}{\partial x_{i}}(\theta^{\frac{1}{2}%
}\Delta\theta^{-\frac{1}{2}})(y_{0})\phi(y_{0})$

$\qquad R_{2}=\frac{1}{24}\frac{\partial}{\partial x_{i}}(\Phi^{-1}\Delta
\Phi)(y_{0})\phi(y_{0})$

$\qquad R_{3}=\frac{1}{12}\frac{\partial}{\partial x_{i}}[<\nabla\log
\theta^{-\frac{1}{2}},\nabla\log\Phi+X>](y_{0})\phi(y_{0})$

$\qquad R_{4}=\frac{1}{12}\frac{\partial}{\partial x_{i}}<\nabla\log
\Phi,X>\phi(y_{0})$

$\qquad R_{5}=\frac{1}{12}\frac{\partial}{\partial x_{i}}[$V$](y_{0}%
)\phi(y_{0})$

We compute each of these terms:

$\qquad R_{1}=\frac{1}{24}\frac{\partial}{\partial x_{i}}(\theta^{\frac{1}{2}%
}\Delta\theta^{-\frac{1}{2}})(y_{0})\phi(y_{0})$

\qquad\qquad$=\frac{1}{24}\frac{\partial\theta^{\frac{1}{2}}}{\partial x_{i}%
}(y_{0})(\Delta\theta^{-\frac{1}{2}})(y_{0})\phi(y_{0})+\frac{1}{24}%
\theta^{\frac{1}{2}}(y_{0})\frac{\partial}{\partial x_{i}}(\Delta
\theta^{-\frac{1}{2}})(y_{0})\phi(y_{0})$

\qquad\qquad$=\frac{1}{24}\frac{\partial\theta^{\frac{1}{2}}}{\partial x_{i}%
}(y_{0})(\Delta\theta^{-\frac{1}{2}})(y_{0})\phi(y_{0})+\frac{1}{24}%
\frac{\partial}{\partial x_{i}}(\Delta\theta^{-\frac{1}{2}})(y_{0})\phi
(y_{0})$

$\frac{\partial\theta^{\frac{1}{2}}}{\partial x_{i}}(y_{0})=-\frac{1}%
{2}<H,i>(y_{0})$ and so, we have:

$R_{1}=-\frac{1}{48}<H,i>(y_{0})(\Delta\theta^{-\frac{1}{2}})(y_{0})\phi
(y_{0})+\frac{1}{24}\frac{\partial}{\partial x_{i}}(\Delta\theta^{-\frac{1}%
{2}})(y_{0})\phi(y_{0})$

The expression for $\frac{1}{2}\Delta\theta^{-\frac{1}{2}}(y_{0})$ is given by
(ii) of \textbf{Table A}$_{10}$ and that of $\frac{1}{24}\frac{\partial
}{\partial x_{i}}(\Delta\theta^{-\frac{1}{2}})(y_{0})$ is given by (viii) of
\textbf{Table A}$_{10}.$ Since $\theta^{\frac{1}{2}}(y_{0})=1;$ we have:

$\left(  B_{112}\right)  $\qquad$R_{1}\ =-\frac{1}{576}<H,i>(y_{0}%
)\qquad\qquad\qquad\qquad\qquad\qquad\left(  138\right)  $

$\qquad\qquad\times\lbrack3<H,i>^{2}+2(\tau^{M}-3\tau^{P}%
\ +\overset{q}{\underset{\text{a=1}}{\sum}}\varrho_{\text{aa}}^{M}%
+\overset{q}{\underset{\text{a,b}=1}{\sum}}R_{\text{abab}}^{M})](y_{0}%
)\phi(y_{0})\qquad$

$\qquad-\frac{1}{24}\underset{j=q+1}{\overset{n}{%
{\textstyle\sum}
}}<H,j>(y_{0})\underset{\text{a,b=}1}{\overset{q}{%
{\textstyle\sum}
}}T_{\text{ab}j}^{2}(y_{0})\phi(y_{0})\qquad$

\qquad$\ +\frac{5}{64}\underset{j=q+1}{\overset{n}{%
{\textstyle\sum}
}}<H,i><H,j>^{2}\phi(y_{0})\qquad$

$\qquad+\frac{1}{96}\underset{j=q+1}{\overset{n}{%
{\textstyle\sum}
}}<H,j>(y_{0})$

$\qquad\times\lbrack(2\varrho_{ij}+4\overset{q}{\underset{\text{a}=1}{\sum}%
}R_{i\text{a}j\text{a}}-3\overset{q}{\underset{\text{a},b=1}{\sum}%
}T_{\text{aa}j}T_{\text{bb}i}-T_{\text{ab}j}T_{\text{ab}i}%
-3\overset{q}{\underset{\text{a},\text{b}=1}{\sum}}T_{\text{aa}i}%
T_{\text{bb}j}-T_{\text{ab}i}T_{\text{ab}j})](y_{0})\phi(y_{0})$

$\qquad+\frac{1}{96}<H,i>[$ $\tau^{M}-3\tau^{P}+\overset{q}{\underset{\text{a}%
=1}{\sum}}\varrho_{\text{aa}}+\overset{q}{\underset{\text{a,b}=1}{\sum}%
}R_{\text{abab}}](y_{0})\phi(y_{0})$

$\qquad+\frac{1}{288}[\nabla_{i}\varrho_{jj}-2\varrho_{ij}%
<H,j>+\overset{q}{\underset{\text{a}=1}{\sum}}(\nabla_{i}R_{\text{a}%
j\text{a}j}-4R_{i\text{a}j\text{a}}<H,j>)+4\overset{q}{\underset{\text{a,b}%
=1}{\sum}}R_{i\text{a}j\text{b}}T_{\text{ab}j}$

\qquad$+2\overset{q}{\underset{\text{a,b,c}=1}{\sum}}%
\underset{j=q+1}{\overset{n}{%
{\textstyle\sum}
}}(T_{\text{aa}i}T_{\text{bb}j}T_{\text{cc}j}-3T_{\text{aa}i}T_{\text{bc}%
j}T_{\text{bc}j}+2T_{\text{ab}i}T_{\text{bc}j}T_{\text{ca}j})](y_{0}%
)\phi(y_{0})$\qquad\qquad\qquad\qquad\qquad\ \ 

$\qquad+\frac{1}{288}\underset{j=q+1}{\overset{n}{%
{\textstyle\sum}
}}[\nabla_{j}\varrho_{ij}-2\varrho_{ij}<H,j>+\overset{q}{\underset{\text{a}%
=1}{\sum}}(\nabla_{j}R_{\text{a}i\text{a}j}-4R_{j\text{a}i\text{a}%
}<H,j>)+4\overset{q}{\underset{\text{a,b}=1}{\sum}}R_{j\text{a}i\text{b}%
}T_{\text{ab}j}$

\qquad$+2\overset{q}{\underset{\text{a,b,c}=1}{\sum}}(T_{\text{aa}%
j}T_{\text{bb}i}T_{\text{cc}j}-3T_{\text{aa}j}T_{\text{bc}i}T_{\text{bc}%
j}+2T_{\text{ab}j}T_{\text{bc}i}T_{\text{ca}j})](y_{0})\phi(y_{0})$

$\qquad+\frac{1}{288}\underset{j=q+1}{\overset{n}{%
{\textstyle\sum}
}}[\nabla_{j}\varrho_{ij}-2\varrho_{jj}<H,i>+\overset{q}{\underset{\text{a}%
=1}{\sum}}(\nabla_{j}R_{\text{a}i\text{a}j}-4R_{j\text{a}j\text{a}%
}<H,i>)+4\overset{q}{\underset{\text{a},b=1}{\sum}}R_{j\text{a}j\text{b}%
}T_{\text{ab}i}$

\qquad$+2\overset{q}{\underset{\text{a,b,c}=1}{\sum}}(T_{\text{aa}%
j}T_{\text{bb}j}T_{\text{cc}i}-3T_{\text{aa}j}T_{\text{bc}j}T_{\text{bc}%
i}+2T_{\text{ab}j}T_{\text{bc}j}T_{\text{ca}i})](y_{0})\phi(y_{0})$

$\qquad-\frac{1}{48}\underset{k=q+1}{\overset{n}{%
{\textstyle\sum}
}}<H,k>(y_{0})[R_{\text{a}i\text{a}k}-\underset{\text{c=1}}{\overset{\text{q}%
}{\sum}}T_{\text{ac}i}.T_{\text{ac}k}+\frac{2}{3}\underset{j=q+1}{\overset{n}{%
{\textstyle\sum}
}}R_{ijkj}](y_{0})\phi(y_{0})$

$\qquad-\frac{1}{24}\underset{j=q+1}{\overset{n}{%
{\textstyle\sum}
}}<H,j>(y_{0})[\frac{3}{4}<H,i><H,j>$\ $\qquad\ \ \ +\frac{1}{12}%
(2\varrho_{ij}+4\overset{q}{\underset{\text{a}=1}{\sum}}R_{i\text{a}k\text{a}%
}-3\overset{q}{\underset{\text{a},b=1}{\sum}}T_{\text{aa}j}T_{\text{bb}%
i}-T_{\text{ab}j}T_{\text{ab}i}-3\overset{q}{\underset{\text{a},\text{b}%
=1}{\sum}}T_{\text{aa}i}T_{\text{bb}j}-T_{\text{ab}i}T_{\text{ab}j}%
)]\phi(y_{0})$

Next,

\qquad$\ R_{2}=\frac{1}{24}\frac{\partial}{\partial x_{i}}(\Phi^{-1}\Delta
\Phi)(y_{0})$

\qquad$\ \ \ \ \ =\frac{1}{24}\frac{\partial\Phi^{-1}}{\partial x_{i}}%
(y_{0})\Delta\Phi(y_{0})+\frac{1}{24}\Phi^{-1}(y_{0})\frac{\partial}{\partial
x_{i}}(\Delta\Phi)(y_{0})$

We have: $\Phi^{-1}(y_{0})=1$ and $\frac{\partial\Phi^{-1}}{\partial x_{i}%
}(y_{0})=X_{i}(y_{0})$ is given in (iii) of \textbf{Table B}$_{4}$

$R_{2}=\frac{1}{24}X_{i}(y_{0})\Delta\Phi(y_{0})+\frac{1}{24}\frac{\partial
}{\partial x_{i}}(\Delta\Phi)(y_{0})$

The expression for $\Delta\Phi_{P}(y_{0})$ is in (iii) of \textbf{Table
B}$_{3}$ and that of $\frac{1}{24}\frac{\partial}{\partial x_{i}}[\Delta
\Phi](y_{0})$ in (i) of \textbf{Table B}$_{5}.$

$\qquad\Delta\Phi_{P}(y_{0})=$ $\left\Vert \text{X}\right\Vert ^{2}(y_{0})-$
divX$(y_{0})-\underset{\text{a}=1}{\overset{q}{\sum}}$X$_{\text{a}}^{2}%
(y_{0})+$ $\underset{\text{a}=1}{\overset{q}{\sum}}\frac{\partial X_{\text{a}%
}}{\partial x_{\text{a}}}(y_{0})$

\qquad\qquad$\qquad=$ $\left\Vert \text{X}\right\Vert _{M}^{2}(y_{0})-$
divX$_{M}(y_{0})-\left\Vert \text{X}\right\Vert _{P}^{2}(y_{0})+$
divX$_{P}(y_{0})$

$\qquad\frac{1}{24}\frac{\partial}{\partial x_{i}}[\Delta\Phi](y_{0})$ is
given by (i) of \textbf{Table B}$_{5}$

Consequently, the equation $R_{2}=\frac{1}{24}X_{i}(y_{0})\Delta\Phi
(y_{0})+\frac{1}{24}\frac{\partial}{\partial x_{i}}(\Delta\Phi)(y_{0})$ becomes:

$R_{2}=\frac{1}{24}X_{i}(y_{0})[\left\Vert \text{X}\right\Vert _{M}^{2}-$
divX$_{M}-$ $\left\Vert \text{X}\right\Vert _{P}^{2}$ $+$ divX$_{P}%
](y_{0})\qquad\qquad\qquad\qquad\qquad\qquad\qquad\qquad\left(  139\right)
\qquad\qquad\qquad\qquad\qquad$

\qquad$+\frac{1}{12}[T_{\text{ab}i}T_{\text{ab}j}X_{j}+2\perp_{\text{a}%
ij}(\frac{\partial X_{j}}{\partial x_{\text{c}}}-\perp_{\text{a}jk}%
X_{k})](y_{0})\qquad\qquad Q_{1}$

\qquad$+$ $\frac{1}{6}X_{j}(y_{0})[\frac{\partial X_{j}}{\partial x_{\text{a}%
}}+\frac{\partial X_{j}}{\partial x_{\text{a}}}](y_{0})-\frac{1}{24}\left(
\frac{\partial^{2}X_{j}}{\partial x_{\text{a}}\partial x_{j}}+\frac
{\partial^{2}X_{j}}{\partial x_{\text{a}}\partial x_{j}}\right)  (y_{0})\qquad
Q_{2}$

$\qquad-\frac{1}{24}\frac{\partial^{2}X_{i}}{\partial x_{\text{a}}^{2}}%
(y_{0})+\frac{1}{24}X_{i}(y_{0})[\operatorname{div}X_{M}-\left\Vert
X\right\Vert _{M}^{2}+\left\Vert X\right\Vert _{P}^{2}-\operatorname{div}%
X_{P}-\underset{j=q+1}{\overset{n}{\sum}}<H,j>(y_{0})X_{j}](y_{0})$

$\qquad+\frac{1}{12}[X_{j}\frac{\partial X_{i}}{\partial x_{j}}+$ $\frac{2}%
{3}X_{k}R_{jijk}-\frac{1}{2}\frac{\partial^{2}X_{i}}{\partial x_{j}^{2}%
}](y_{0})+\frac{1}{24}[R_{\text{a}i\text{a}k}-\underset{\text{c=1}%
}{\overset{\text{q}}{\sum}}T_{\text{ac}i}T_{\text{ac}k}](y_{0})X_{k}(y_{0})$

\qquad$+\frac{1}{24}<H,j>(y_{0})[X_{i}X_{j}-\frac{1}{2}\left(  \frac{\partial
X_{j}}{\partial x_{i}}+\frac{\partial X_{i}}{\partial x_{j}}\right)  ](y_{0}%
)$\qquad\qquad

We next compute:

\qquad$R_{3}=\frac{1}{12}\frac{\partial}{\partial x_{i}}[<\nabla\log
\theta^{-\frac{1}{2}},\nabla\log\Phi+X>](y_{0})$

It is immediate that for vector fields X and Yon the Riemannian manifold M, we have:

\qquad$<X,Y>$ $=$ $<X_{j}\frac{\partial}{\partial x_{j}},Y_{k}\frac{\partial
}{\partial x_{k}}>$ $=$ $<\frac{\partial}{\partial x_{j}},\frac{\partial
}{\partial x_{k}}>X_{j}Y_{k}$ $=g_{jk}X_{j}Y_{k}$

\qquad$R_{3}=\frac{1}{12}\frac{\partial}{\partial x_{i}}[g_{jk}(\nabla
\log\theta^{-\frac{1}{2}})_{j}((\nabla\log\Phi)_{k}+X_{k})](y_{0})$

\qquad\qquad$=\frac{1}{12}\frac{\partial g_{jk}}{\partial x_{i}}%
(y_{0})[(\nabla\log\theta^{-\frac{1}{2}})_{j}((\nabla\log\Phi)_{k}%
+X_{k})](y_{0})$

\qquad\qquad$+\frac{1}{12}g_{jk}(y_{0})\frac{\partial}{\partial x_{i}}%
[(\nabla\log\theta^{-\frac{1}{2}})_{j}((\nabla\log\Phi)_{k}+X_{k})](y_{0})$

$\qquad(\nabla\log\Phi)_{k}(y_{0})+X_{k}(y_{0})=0$ by (i) of \textbf{Table
B}$_{1}$ and since $g_{jk}(y_{0})=\delta_{jk},$ we have:

$\qquad R_{3}=\frac{1}{12}\frac{\partial}{\partial x_{i}}[(\nabla\log
\theta^{-\frac{1}{2}})_{j}((\nabla\log\Phi)_{j}+X_{j})](y_{0})$

$\qquad=\frac{1}{12}\frac{\partial}{\partial x_{i}}(\nabla\log\theta
^{-\frac{1}{2}})_{j}(y_{0})[(\nabla\log\Phi)_{j}+X_{j})](y_{0})$

$\qquad+\frac{1}{12}(\nabla\log\theta^{-\frac{1}{2}})_{j}(y_{0})[\frac
{\partial}{\partial x_{i}}(\nabla\log\Phi)_{j}+\frac{\partial X_{j}}{\partial
x_{i}}](y_{0})$

Again, since $(\nabla\log\Phi)_{j}(y_{0})+X_{j}(y_{0})=0,$ we have:

\qquad$R_{3}=\frac{1}{12}(\nabla\log\theta^{-\frac{1}{2}})_{j}(y_{0}%
)[\frac{\partial}{\partial x_{i}}(\nabla\log\Phi)_{j}+\frac{\partial X_{j}%
}{\partial x_{i}}](y_{0})$

From $B_{39}$ above we have:

$\qquad\frac{\partial}{\partial x_{i}}(\nabla\log\Phi)_{j}(y_{0})=-\frac{1}%
{2}\left(  \frac{\partial X_{j}}{\partial x_{i}}+\frac{\partial X_{i}%
}{\partial x_{j}}\right)  (y_{0})$

Further by (ix)$^{\ast\ast}$ of \textbf{Table A}$_{10},$ we have for
$i,j=q+1,...n,$

\qquad$(\nabla\log\theta^{-\frac{1}{2}})_{j}(y_{0})$ $=\frac{1}{2}%
<H,j>(y_{0})$

Therefore,

$\qquad R_{3}=\frac{1}{24}<H,j>(y_{0})[-\frac{1}{2}\left(  \frac{\partial
X_{j}}{\partial x_{i}}+\frac{\partial X_{i}}{\partial x_{j}}\right)
+\frac{\partial X_{j}}{\partial x_{i}}](y_{0})$

$\qquad\qquad=\frac{1}{48}<H,j>(y_{0})[\frac{\partial X_{j}}{\partial x_{i}%
}-\frac{\partial X_{i}}{\partial x_{j}}](y_{0})$

$\qquad R_{3}=\frac{1}{48}<H,j>(y_{0})[\frac{\partial X_{j}}{\partial x_{i}%
}-\frac{\partial X_{i}}{\partial x_{j}}](y_{0})$

We next compute:

\qquad$R_{4}=\frac{1}{12}\frac{\partial}{\partial x_{i}}[<\nabla\log
\Phi,X>](y_{0})=\frac{1}{12}\frac{\partial}{\partial x_{i}}[g_{jk}(\nabla
\log\Phi)_{j}X_{k}](y_{0})$

$\qquad\ \ \ \ =\frac{1}{12}\frac{\partial g_{jk}}{\partial x_{i}}%
(y_{0})[(\nabla\log\Phi)_{j}X_{k}](y_{0})+\frac{1}{12}g_{jk}(y_{0}%
)\frac{\partial}{\partial x_{i}}[(\nabla\log\Phi)_{j}X_{k}](y_{0})$

\qquad$\ \ \ =\frac{1}{12}\frac{\partial g_{jk}}{\partial x_{i}}%
(y_{0})[(\nabla\log\Phi)_{j}X_{k}](y_{0})+\frac{1}{12}g_{jk}(y_{0}%
)[\frac{\partial}{\partial x_{i}}(\nabla\log\Phi)_{j}X_{k}+(\nabla\log
\Phi)_{j}\frac{\partial X_{k}}{\partial x_{i}}](y_{0})$

Since $g_{jk}(y_{0})=\delta^{jk};\frac{\partial g_{\text{ab}}}{\partial x_{i}%
}(y_{0})=-2T_{\text{ab}i}(y_{0});$ $\frac{\partial g_{j\text{a}}}{\partial
x_{i}}(y_{0})=\perp_{\text{a}ij}(y_{0})$ for a,b = 1,...,q and $i,j=q+1,...,n$ and

$\frac{\partial g_{jk}}{\partial x_{i}}(y_{0})=0$ for all $i,j,k=q+1,...,n,$
we have:\qquad$\qquad R_{4}=-\frac{1}{6}T_{\text{ab}i}(y_{0})[(\nabla\log
\Phi)_{\text{a}}X_{\text{b}}](y_{0})+\frac{1}{12}\frac{\partial g_{jk}%
}{\partial x_{i}}(y_{0})[(\nabla\log\Phi)_{j}X_{k}](y_{0})$

$\qquad+\frac{1}{12}[\frac{\partial}{\partial x_{i}}(\nabla\log\Phi)_{j}%
X_{j}+(\nabla\log\Phi)_{j}\frac{\partial X_{j}}{\partial x_{i}}](y_{0})$

Since $(\nabla\log\Phi)_{\text{a}}(y_{0})=0$ by (xi) of \textbf{Table B}%
$_{1},$

\qquad$R_{4}=\frac{1}{12}\frac{\partial g_{j\text{a}}}{\partial x_{i}}%
(y_{0})[(\nabla\log\Phi)_{j}X_{\text{a}}](y_{0})+\frac{1}{12}[\frac{\partial
}{\partial x_{i}}(\nabla\log\Phi)_{j}X_{j}+(\nabla\log\Phi)_{j}\frac{\partial
X_{j}}{\partial x_{i}}](y_{0})$

Since from $B_{39}$ above we have:

$\qquad\frac{\partial}{\partial x_{i}}(\nabla\log\Phi)_{j}(y_{0})=-\frac{1}%
{2}\left(  \frac{\partial X_{j}}{\partial x_{i}}+\frac{\partial X_{i}%
}{\partial x_{j}}\right)  (y_{0})$ and $(\nabla\log\Phi)_{j}(y_{0})$

$\qquad\qquad\qquad\qquad\qquad=-X_{j}(y_{0})$ by (i) of \textbf{Table B}%
$_{4},$

\qquad$R_{4}=\frac{1}{12}\perp_{\text{a}ij}(y_{0})(-X_{j})(y_{0})X_{\text{a}%
}(y_{0})+\frac{1}{12}[-\frac{1}{2}\left(  \frac{\partial X_{j}}{\partial
x_{i}}+\frac{\partial X_{i}}{\partial x_{j}}\right)  X_{j}+(-X_{j}%
)\frac{\partial X_{j}}{\partial x_{i}}](y_{0})$

\qquad$R_{4}=-\frac{1}{12}\perp_{\text{a}ij}(y_{0})X_{\text{a}}(y_{0}%
)X_{j}(y_{0})-\frac{1}{24}\left(  3\frac{\partial X_{j}}{\partial x_{i}}%
+\frac{\partial X_{i}}{\partial x_{j}}\right)  (y_{0})X_{j}(y_{0})$

$\qquad R_{5}=\frac{1}{12}\frac{\partial\text{V}}{\partial x_{i}}(y_{0})$

In conclusion we have for $i=q+1,...,n,$

$\left(  B_{113}\right)  $\qquad$R=\frac{1}{12}\frac{\partial}{\partial x_{i}%
}\left(  \Psi^{-1}L\Psi\right)  (y_{0})=R_{1}+R_{2}+R_{3}+R_{4}+R_{5}$

$\qquad\qquad\ =\frac{1}{96}$%
$<$%
H,$i$%
$>$%
$[$ $\tau^{M}-3\tau^{P}+\overset{q}{\underset{\text{a}=1}{\sum}}%
\varrho_{\text{aa}}+\overset{q}{\underset{\text{a,b}=1}{\sum}}R_{\text{abab}%
}](y_{0})\phi(y_{0})\qquad\qquad\qquad R_{1}$

$\qquad+\frac{1}{288}[\nabla_{i}\varrho_{jj}-2\varrho_{ij}%
<H,j>+\overset{q}{\underset{\text{a}=1}{\sum}}(\nabla_{i}R_{\text{a}%
j\text{a}j}-4R_{i\text{a}j\text{a}}<H,j>)+4\overset{q}{\underset{\text{a,b}%
=1}{\sum}}R_{i\text{a}j\text{b}}T_{\text{ab}j}$

\qquad$+2\overset{q}{\underset{\text{a,b,c}=1}{\sum}}(T_{\text{aa}%
i}T_{\text{bb}j}T_{\text{cc}j}-3T_{\text{aa}i}T_{\text{bc}j}T_{\text{bc}%
j}+2T_{\text{ab}i}T_{\text{bc}j}T_{\text{ca}j})](y_{0})\phi(y_{0})$%
\qquad\qquad\qquad\qquad\qquad\ \ 

$\qquad+\frac{1}{288}[\nabla_{j}\varrho_{ij}-2\varrho_{ij}%
<H,j>+\overset{q}{\underset{\text{a}=1}{\sum}}(\nabla_{j}R_{\text{a}%
i\text{a}j}-4R_{j\text{a}i\text{a}}<H,j>)+4\overset{q}{\underset{\text{a,b}%
=1}{\sum}}R_{j\text{a}i\text{b}}T_{\text{ab}j}$

\qquad$+2\overset{q}{\underset{\text{a,b,c}=1}{\sum}}(T_{\text{aa}%
j}T_{\text{bb}i}T_{\text{cc}j}-3T_{\text{aa}j}T_{\text{bc}i}T_{\text{bc}%
j}+2T_{\text{ab}j}T_{\text{bc}i}T_{\text{ca}j})](y_{0})\phi(y_{0})$

$\qquad+\frac{1}{288}[\nabla_{j}\varrho_{ij}-2\varrho_{jj}%
<H,i>+\overset{q}{\underset{\text{a}=1}{\sum}}(\nabla_{j}R_{\text{a}%
i\text{a}j}-4R_{j\text{a}j\text{a}}<H,i>)+4\overset{q}{\underset{\text{a}%
,b=1}{\sum}}R_{j\text{a}j\text{b}}T_{\text{ab}i}$

\qquad$+2\overset{q}{\underset{\text{a,b,c}=1}{\sum}}(T_{\text{aa}%
j}T_{\text{bb}j}T_{\text{cc}i}-3T_{\text{aa}j}T_{\text{bc}j}T_{\text{bc}%
i}+2T_{\text{ab}j}T_{\text{bc}j}T_{\text{ca}i})](y_{0})\phi(y_{0})$

$\qquad-\frac{1}{48}<H,k>(y_{0})[R_{\text{a}i\text{a}k}-\underset{\text{c=1}%
}{\overset{\text{q}}{\sum}}T_{\text{ac}i}T_{\text{ac}k}+\frac{2}{3}%
R_{ijkj}](y_{0})\phi(y_{0})$

$\qquad-\frac{1}{24}<H,k>(y_{0})[\frac{3}{4}<H,i><H,k>$\ $+\frac{1}%
{12}(2\varrho_{ik}+4\overset{q}{\underset{\text{a}=1}{\sum}}R_{i\text{a}%
k\text{a}}-6\overset{q}{\underset{\text{a},b=1}{\sum}}T_{\text{aa}%
i}T_{\text{bb}k}-T_{\text{ab}i}T_{\text{ab}k})]\phi(y_{0})$

\qquad$+\frac{1}{24}X_{i}(y_{0})[\left\Vert \text{X}\right\Vert _{M}^{2}-$
divX$_{M}-$ $\left\Vert \text{X}\right\Vert _{P}^{2}$ $+$ divX$_{P}%
](y_{0})\qquad\qquad\qquad\qquad\qquad R_{2}\qquad$From $\left(  139\right)
\qquad\qquad\qquad\qquad\qquad$

\qquad$+\frac{1}{12}[T_{\text{ab}i}T_{\text{ab}j}X_{j}+2\perp_{\text{a}%
ij}(\frac{\partial X_{j}}{\partial x_{\text{c}}}-\perp_{\text{a}jk}%
X_{k})](y_{0})\qquad\qquad Q_{1}$

\qquad$+$ $\frac{1}{6}X_{j}(y_{0})[\frac{\partial X_{j}}{\partial x_{\text{a}%
}}+\frac{\partial X_{j}}{\partial x_{\text{a}}}](y_{0})-\frac{1}{24}\left(
\frac{\partial^{2}X_{j}}{\partial x_{\text{a}}\partial x_{j}}+\frac
{\partial^{2}X_{j}}{\partial x_{\text{a}}\partial x_{j}}\right)  (y_{0})\qquad
Q_{2}$

$\qquad-\frac{1}{24}\frac{\partial^{2}X_{i}}{\partial x_{\text{a}}^{2}}%
(y_{0})+\frac{1}{24}X_{i}(y_{0})[\operatorname{div}X_{M}-\left\Vert
X\right\Vert _{M}^{2}+\left\Vert X\right\Vert _{P}^{2}-\operatorname{div}%
X_{P}-\underset{j=q+1}{\overset{n}{\sum}}<H,j>(y_{0})X_{j}](y_{0})$

$\qquad+\frac{1}{12}[X_{j}\frac{\partial X_{i}}{\partial x_{j}}+$ $\frac{2}%
{3}X_{k}R_{jijk}-\frac{1}{2}\frac{\partial^{2}X_{i}}{\partial x_{j}^{2}%
}](y_{0})+\frac{1}{24}[R_{\text{a}i\text{a}k}-\underset{\text{c=1}%
}{\overset{\text{q}}{\sum}}T_{\text{ac}i}T_{\text{ac}k}](y_{0})X_{k}(y_{0})$

\qquad$+\frac{1}{24}<H,j>(y_{0})[X_{i}X_{j}-\frac{1}{2}\left(  \frac{\partial
X_{j}}{\partial x_{i}}+\frac{\partial X_{i}}{\partial x_{j}}\right)  ](y_{0})$

\qquad$+\frac{1}{48}<H,j>(y_{0})[\frac{\partial X_{j}}{\partial x_{i}}%
-\frac{\partial X_{i}}{\partial x_{j}}](y_{0})\qquad\qquad\qquad\qquad
\qquad\qquad R_{3}$

$\qquad-\frac{1}{12}\perp_{\text{a}ij}(y_{0})X_{\text{a}}(y_{0})X_{j}%
(y_{0})-\frac{1}{24}\left(  3\frac{\partial X_{j}}{\partial x_{i}}%
+\frac{\partial X_{i}}{\partial x_{j}}\right)  (y_{0})X_{j}(y_{0})\qquad\qquad
R_{4}$

$\qquad+\frac{1}{12}\frac{\partial\text{V}}{\partial x_{i}}(y_{0})\qquad
\qquad\qquad\qquad\qquad\qquad\qquad\qquad\qquad\qquad\qquad\qquad R_{5}%
\qquad\qquad\qquad$

(viii)\qquad I$_{321}\mathbf{=}\frac{1}{12}\frac{\partial^{2}}{\partial
x_{i}^{2}}\left(  \frac{\text{L}\Psi}{\Psi}\right)  (y_{0})\phi(y_{0})$ for
$i=q+1,...,n.$

We compute it here in the general case. We recall that by definition, the
scalar differential operator L is given by:

L $=$ $\frac{1}{2}\Delta+X+V$ where $\Delta$ is now the scalar Laplacian and
$\Psi(x)=\theta_{P}(x)^{-\frac{1}{2}}\ \ \Phi_{P}(x).$ Then,

L$\Psi=$ L$(\theta_{P}{}^{-\frac{1}{2}}\ \ \Phi_{P})=\frac{1}{2}\Phi
\Delta\theta_{P}{}^{-\frac{1}{2}}\ +\frac{1}{2}\theta_{P}{}^{-\frac{1}{2}%
}\ \Delta\Phi+$ $<\nabla\theta^{-\frac{1}{2}},\nabla\Phi>$

$+$ $\theta^{-\frac{1}{2}}$ $\nabla_{X}\Phi+\Phi\nabla_{X}\theta_{P}{}%
^{-\frac{1}{2}}$ $+$ $V(\theta_{P}(x)^{-\frac{1}{2}}\ \ \Phi_{P})$

Simplifying, we have:

$\frac{\text{L}\Psi}{\Psi}=\frac{1}{2}\theta^{\frac{1}{2}}\Delta\theta
^{-\frac{1}{2}}+$ $\frac{1}{2}\Phi^{-1}\Delta\Phi+$ $<\nabla\log\theta
^{-\frac{1}{2}},\nabla\log\Phi+X>+$ $<\nabla\log\Phi,X>+$ $V$

We set:

$\qquad$I$_{321}=\frac{1}{24}$ $\frac{\partial^{2}}{\partial x_{i}^{2}}%
(\frac{\text{L}\Psi}{\Psi})(y_{0})\phi(y_{0})=$ I$_{3211}+$ I$_{3212}+$
I$_{3213}+$ I$_{3214}+$ I$_{3215}$where,

\qquad I$_{3211}=\frac{1}{24}$ $\frac{\partial^{2}}{\partial x_{i}^{2}}%
(\theta^{\frac{1}{2}}\Delta\theta^{-\frac{1}{2}})(y_{0})\phi(y_{0})$

\ \ \ \ \ I$_{3212}$\ $=\frac{1}{24}$ $\frac{\partial^{2}}{\partial x_{i}^{2}%
}(\Phi^{-1}\Delta\Phi)(y_{0})\phi(y_{0})$

\ \ \ \ I$_{3213}=\frac{1}{12}\frac{\partial^{2}}{\partial x_{i}^{2}}%
(<\nabla\log\theta^{-\frac{1}{2}},\nabla\log\Phi+X>)\phi(y_{0})$

\ \ \ \ I$_{3214}=\frac{1}{12}\frac{\partial^{2}}{\partial x_{i}^{2}}%
[<\nabla\log\Phi,X>](y_{0})\phi(y_{0})$

\ \ \ \ I$_{3215}=\frac{1}{12}\frac{\partial^{2}\text{V}}{\partial x_{i}^{2}%
}(y_{0})\phi(y_{0})$

We compute the above items in geometric invariants:

$\left(  B_{114}\right)  $\qquad I$_{3211}=\frac{1}{24}\frac{\partial^{2}%
}{\partial\text{x}_{i}^{2}}(\theta^{\frac{1}{2}}\Delta\theta^{-\frac{1}{2}%
})(y_{0})\phi(y_{0})$ is A$_{321}$ in \textbf{Appendix A}$_{10}$

We then compute:

\qquad\ I$_{3212}$\ $=\frac{1}{24}$ $\frac{\partial^{2}}{\partial x_{i}^{2}%
}(\Phi^{-1}\Delta\Phi)(y_{0})\phi(y_{0})$

$\qquad\ \ =\frac{1}{24}$ $\frac{\partial^{2}\Phi^{-1}}{\partial x_{i}^{2}%
}(y_{0})(\Delta\Phi)(y_{0})\phi(y_{0})+\frac{1}{24}(\Phi^{-1})(y_{0}%
)\frac{\partial^{2}}{\partial x_{i}^{2}}(\Delta\Phi)(y_{0})\phi(y_{0})$

$\qquad\ \ \ +\frac{1}{12}$ $\frac{\partial\Phi^{-1}}{\partial x_{i}}%
(y_{0})\frac{\partial}{\partial x_{i}}(\Delta\Phi)(y_{0}).\phi(y_{0})$

\qquad$=\frac{1}{24}$ $\frac{\partial^{2}\Phi^{-1}}{\partial x_{i}^{2}}%
(y_{0})(\Delta\Phi)(y_{0})\phi(y_{0})+\frac{1}{24}\frac{\partial^{2}}{\partial
x_{i}^{2}}(\Delta\Phi)(y_{0})\phi(y_{0})$

$\qquad\ \ \ +\frac{1}{12}$ $\frac{\partial\Phi^{-1}}{\partial x_{i}}%
(y_{0})\frac{\partial}{\partial x_{i}}(\Delta\Phi)(y_{0})\phi(y_{0}%
)=$\ I$_{32121}+$\ I$_{32122}$\ $+$ I$_{3213}$\ where,

\qquad\ I$_{32121}=\frac{1}{24}\frac{\partial^{2}\Phi^{-1}}{\partial x_{i}%
^{2}}(y_{0})(\Delta\Phi)(y_{0})\phi(y_{0})$

\ \ \ \ \ \ I$_{32122}=\frac{1}{12}$ $\frac{\partial\Phi^{-1}}{\partial x_{i}%
}(y_{0})\frac{\partial}{\partial x_{i}}(\Delta\Phi)(y_{0})\phi(y_{0})$

$\qquad$\ I$_{32123}=\frac{1}{24}\frac{\partial^{2}}{\partial x_{i}^{2}%
}(\Delta\Phi)(y_{0})\phi(y_{0})$

\ \ 

All of the above have already been computed:

\qquad$\frac{\partial\Phi^{-1}}{\partial x_{i}}(y_{0})=X_{i}(y_{0})$ is given
in (iii) of \textbf{Table B}$_{4}$

\qquad\ $\frac{\partial^{2}\Phi^{-1}}{\partial x_{i}^{2}}(y_{0})=X_{i}%
^{2}(y_{0})+\frac{\partial X_{i}}{\partial x_{i}}(y_{0})$ is given in (iv) of
\textbf{Table B}$_{4}$ and so,

\qquad$\underset{i=1}{\overset{n}{\sum}}\frac{\partial^{2}\Phi^{-1}}{\partial
x_{i}^{2}}(y_{0})=$ $\left\Vert \text{X}\right\Vert ^{2}(y_{0})+$
divX$(y_{0})-\underset{\text{a}=1}{\overset{q}{\sum}}$X$_{\text{a}}^{2}%
(y_{0})-$ $\underset{\text{a}=1}{\overset{q}{\sum}}\frac{\partial X_{\text{a}%
}}{\partial x_{\text{a}}}(y_{0})$

\qquad\qquad\qquad\qquad$=$ $\left\Vert \text{X}\right\Vert _{M}^{2}(y_{0})+$
divX$_{M}(y_{0})-\left\Vert \text{X}\right\Vert _{P}^{2}(y_{0})-$
divX$_{P}(y_{0})$

Then (iii) of \textbf{Table B}$_{3}$ and (xvii) of Table B$_{4}$ gives:

$\qquad\Delta\Phi_{P}(y_{0})=$ $\left\Vert \text{X}\right\Vert ^{2}(y_{0})-$
divX$(y_{0})-\underset{\text{a}=1}{\overset{q}{\sum}}$X$_{\text{a}}^{2}%
(y_{0})+$ $\underset{\text{a}=1}{\overset{q}{\sum}}\frac{\partial X_{\text{a}%
}}{\partial x_{\text{a}}}(y_{0})$

\qquad\qquad\qquad$=$ $\left\Vert \text{X}\right\Vert _{M}^{2}(y_{0})-$
divX$_{M}(y_{0})-\left\Vert \text{X}\right\Vert _{P}^{2}(y_{0})+$
divX$_{P}(y_{0}).$Therefore,

\ I$_{32121}=\frac{1}{24}[\left\Vert \text{X}\right\Vert _{M}^{2}%
+\operatorname{div}$X$_{M}-\left\Vert \text{X}\right\Vert _{P}^{2}%
-\operatorname{div}$X$_{P}](y_{0})$

$\qquad\qquad\times\lbrack\left\Vert \text{X}\right\Vert _{M}^{2}%
-\operatorname{div}X_{M}-\left\Vert \text{X}\right\Vert _{P}^{2}%
+\operatorname{div}X_{P}](y_{0})\phi(y_{0})$

\ I$_{32122}=\frac{1}{12}$ $\frac{\partial\Phi^{-1}}{\partial x_{i}}%
(y_{0})\frac{\partial}{\partial x_{i}}(\Delta\Phi)(y_{0})\phi(y_{0})=\frac
{1}{12}$ $X_{i}(y_{0})\frac{\partial}{\partial x_{i}}(\Delta\Phi)(y_{0}%
)\phi(y_{0})$

where $\frac{1}{2}\frac{\partial}{\partial x_{i}}[\Delta\Phi](y_{0})$ is the
expression in $\left(  B_{102}\right)  $

\ I$_{32123}=\frac{1}{24}\frac{\partial^{2}}{\partial x_{i}^{2}}(\Delta
\Phi)(y_{0})\phi(y_{0})$ is given in $\left(  B_{106}\right)  .$\ \ \ \ \ \ \ \ \ \ \ \ \ 

Therefore, for $i=q+1,...,n,$ we have:

$\left(  B_{114}\right)  ^{\ast}\qquad$\ I$_{3212}$\ $=\frac{1}{24}$
$\frac{\partial^{2}}{\partial x_{i}^{2}}(\Phi^{-1}\Delta\Phi)(y_{0})\phi
(y_{0})$

\qquad$=\frac{1}{24}[\left\Vert \text{X}\right\Vert _{M}^{2}%
+\operatorname{div}$X$_{M}-\left\Vert \text{X}\right\Vert _{P}^{2}%
-\operatorname{div}X_{P}](y_{0})$

$\qquad\times\left\Vert \text{X}\right\Vert _{M}^{2}-\operatorname{div}$%
X$_{M}-\left\Vert \text{X}\right\Vert _{P}^{2}+\operatorname{div}$%
X$_{P}](y_{0})\phi(y_{0})\qquad$I$_{32121}$

\qquad$+\frac{1}{6}X_{i}(y_{0})T_{\text{ab}i}(y_{0})T_{\text{ab}j}(y_{0}%
)X_{j}(y_{0})+\frac{1}{3}\perp_{\text{a}ij}(y_{0})X_{i}(y_{0})[\frac{\partial
X_{j}}{\partial x_{\text{a}}}-\perp_{\text{a}jk}X_{k}](y_{0})\qquad$%
I$_{32122}\qquad Q_{1}$

\qquad$+$ $\frac{2}{3}X_{i}(y_{0})X_{j}(y_{0})\frac{\partial X_{j}}{\partial
x_{\text{a}}}(y_{0})-\frac{1}{6}X_{i}(y_{0})\frac{\partial^{2}X_{j}}{\partial
x_{\text{a}}\partial x_{j}}(y_{0})\qquad\qquad Q_{2}$

$\qquad-\frac{1}{12}X_{i}(y_{0})\frac{\partial^{2}X_{i}}{\partial x_{\text{a}%
}^{2}}(y_{0})+\frac{1}{12}X_{i}^{2}(y_{0})[\operatorname{div}X_{M}-\left\Vert
X\right\Vert _{M}^{2}+\left\Vert X\right\Vert _{P}^{2}-\operatorname{div}%
X_{P}-$ $<H,j>X_{j}](y_{0})$

$\qquad+\frac{1}{6}X_{i}(y_{0})X_{j}(y_{0})\frac{\partial X_{i}}{\partial
x_{j}}(y_{0})+$ $\frac{1}{18}X_{i}(y_{0})X_{k}(y_{0})R_{jijk}(y_{0})-\frac
{1}{12}X_{i}(y_{0})\frac{\partial^{2}X_{i}}{\partial x_{j}^{2}}(y_{0})$

\qquad$+\frac{1}{12}[R_{\text{a}i\text{a}k}-\underset{\text{c=1}%
}{\overset{\text{q}}{\sum}}T_{\text{ac}i}T_{\text{ac}k}-\perp_{\text{a}%
ik}\perp_{\text{a}jk}](y_{0})X_{k}(y_{0})+\frac{1}{18}R_{ijkj}(y_{0}%
)X_{i}(y_{0})X_{k}(y_{0})$

\qquad$+\frac{1}{12}<H,j>(y_{0})X_{i}(y_{0})[X_{i}X_{j}-\frac{1}{2}\left(
\frac{\partial X_{j}}{\partial x_{i}}+\frac{\partial X_{i}}{\partial x_{j}%
}\right)  ](y_{0})\qquad\qquad\qquad\qquad\qquad\qquad\qquad\qquad\qquad
\qquad\ \qquad\qquad\qquad\qquad\qquad\qquad\qquad\qquad\qquad\qquad\qquad$

$-\frac{1}{6}[-R_{\text{a}i\text{b}i}+5\overset{q}{\underset{\text{c}=1}{\sum
}}T_{\text{ac}i}T_{\text{bc}i}+2\overset{n}{\underset{j=q+1}{\sum}}%
\perp_{\text{a}ij}\perp_{\text{b}ij}](y_{0})\underset{k=q+1}{\overset{n}{\sum
}}T_{\text{ab}k}(y_{0})X_{k}(y_{0})\qquad\qquad$I$_{32123}\qquad S_{1}\qquad$

$-\frac{2}{9}\underset{j=q+1}{\overset{n}{\sum}}R_{i\text{a}ij}(y_{0}%
)[\frac{\partial X_{j}}{\partial x_{\text{a}}}%
-\underset{k=q+1}{\overset{n}{\sum}}\perp_{\text{a}jk}X_{k}](y_{0})$

$+\frac{1}{12}\times\frac{2}{3}\underset{j,k=q+1}{\overset{n}{\sum}}%
R_{ijik}(y_{0})[X_{j}X_{k}-\frac{1}{2}(\frac{\partial X_{j}}{\partial x_{k}%
}+\frac{\partial X_{k}}{\partial x_{j}})](y_{0})$

$-\frac{1}{6}T_{\text{ab}i}(y_{0})\frac{\partial^{2}X_{i}}{\partial
x_{\text{a}}\partial x_{\text{b}}}(y_{0})\qquad\qquad\qquad S_{2}\qquad\qquad
S_{21}\qquad\qquad\qquad\qquad\qquad\qquad$

$+$ $\frac{1}{12}T_{\text{ab}i}(y_{0})[$ $(R_{\text{a}i\text{b}j}%
+R_{\text{a}j\text{b}i})$ $-\underset{\text{c=1}}{\overset{\text{q}}{\sum}%
}(T_{\text{ac}i}T_{\text{bc}j}+T_{\text{ac}j}T_{\text{bc}i})$

$-\overset{n}{\underset{k=q+1}{\sum}}(\perp_{\text{a}ik}\perp_{\text{b}jk}+$
$\perp_{\text{a}jk}\perp_{\text{b}ik})](y_{0})X_{j}(y_{0})$

$-$ $\frac{1}{6}T_{\text{ab}i}(y_{0})T_{\text{ab}j}(y_{0})[X_{i}X_{j}-\frac
{1}{2}\left(  \frac{\partial X_{i}}{\partial x_{j}}+\frac{\partial X_{j}%
}{\partial x_{i}}\right)  ](y_{0})$

$-\frac{1}{3}\perp_{\text{a}ij}(y_{0})[(X_{i}\frac{\partial X_{j}}{\partial
x_{\text{a}}}+X_{j}\frac{\partial X_{i}}{\partial x_{\text{a}}})-\frac{1}%
{4}\left(  \frac{\partial^{2}X_{i}}{\partial x_{\text{a}}\partial x_{j}}%
+\frac{\partial^{2}X_{j}}{\partial x_{\text{a}}\partial x_{i}}\right)
](y_{0})\qquad S_{22}$

$-\frac{1}{6}\perp_{\text{a}ij}(y_{0})[T_{\text{ab}j}\frac{\partial X_{i}%
}{\partial x_{\text{b}}}](y_{0})$

$+\frac{1}{6}\perp_{\text{a}ij}(y_{0})[(\perp_{\text{b}ik}T_{\text{ab}%
j})+\frac{2}{3}(2R_{\text{a}ijk}+R_{\text{a}jik}+R_{\text{a}kji})](y_{0}%
)X_{k}(y_{0})$

$-\frac{1}{6}\perp_{\text{a}ij}(y_{0})\perp_{\text{a}jk}(y_{0})[X_{i}%
X_{k}-\frac{1}{2}\left(  \frac{\partial X_{i}}{\partial x_{k}}+\frac{\partial
X_{k}}{\partial x_{i}}\right)  ](y_{0})\qquad\qquad\qquad\qquad\qquad
\qquad\qquad\qquad$

$+\frac{1}{12}[(\frac{\partial X_{j}}{\partial x_{\text{a}}})^{2}+X_{j}%
\frac{\partial^{2}X_{j}}{\partial x_{\text{a}}^{2}}-\frac{1}{2}\frac
{\partial^{3}X_{j}}{\partial x_{\text{a}}^{2}\partial x_{j}}](y_{0})-\frac
{1}{6}\overset{n}{\underset{k=q+1}{\sum}}[\perp_{\text{b}ik}$T$_{\text{aa}%
k}\frac{\partial X_{i}}{\partial x_{\text{b}}^{2}}](y_{0})\qquad\qquad
S_{3}\qquad S_{31}$

$+\frac{1}{144}[\{4\nabla_{i}R_{i\text{a}j\text{a}}+2\nabla_{j}R_{i\text{a}%
i\text{a}}+$ $8(\overset{q}{\underset{\text{c=1}}{%
{\textstyle\sum}
}}R_{\text{a}i\text{c}i}^{{}}T_{\text{ac}j}+\;\overset{n}{\underset{k=q+1}{%
{\textstyle\sum}
}}R_{\text{a}iik}\perp_{\text{a}jk})$

$+8(\overset{q}{\underset{\text{c=1}}{%
{\textstyle\sum}
}}R_{\text{a}i\text{c}j}^{{}}T_{\text{ac}i}+\;\overset{n}{\underset{k=q+1}{%
{\textstyle\sum}
}}R_{\text{a}ijk}\perp_{\text{a}ik})+8(\overset{q}{\underset{\text{c=1}}{%
{\textstyle\sum}
}}R_{\text{a}j\text{c}i}^{{}}T_{\text{ac}i}+\;\overset{n}{\underset{k=q+1}{%
{\textstyle\sum}
}}R_{\text{a}jik}\perp_{\text{a}ik})\}$\ 

$\ +\frac{2}{3}\underset{k=q+1}{\overset{n}{\sum}}\{T_{\text{aa}k}%
(R_{ijik}+3\overset{q}{\underset{\text{c}=1}{\sum}}\perp_{\text{c}ij}%
\perp_{\text{c}ik})\}](y_{0})X_{k}(y_{0})$

$-\frac{1}{12}[$ R$_{\text{a}i\text{a}k}$ $-\underset{\text{c=1}%
}{\overset{\text{q}}{\sum}}T_{\text{ac}i}T_{\text{ac}k}%
-\overset{n}{\underset{l=q+1}{\sum}}(\perp_{\text{a}il}\perp_{\text{a}%
kl}](y_{0})\times\lbrack X_{i}X_{k}-\frac{1}{2}\left(  \frac{\partial X_{i}%
}{\partial x_{k}}+\frac{\partial X_{k}}{\partial x_{i}}\right)  ](y_{0})$

$-\frac{1}{24}T_{\text{aa}k}(y_{0})[-X_{i}^{2}X_{k}+X_{k}\frac{\partial X_{i}%
}{\partial x_{i}}\ +X_{i}\left(  \frac{\partial X_{k}}{\partial x_{i}}%
+\frac{\partial X_{i}}{\partial x_{k}}\right)  -\frac{1}{3}\left(
\frac{\partial^{2}X_{k}}{\partial x_{i}^{2}}+2\frac{\partial^{2}X_{i}%
}{\partial x_{i}\partial x_{k}}\right)  ](y_{0})$

$+\frac{1}{18}[R_{\text{a}jij}\frac{\partial X_{i}}{\partial x_{\text{a}}^{2}%
}](y_{0})\qquad\qquad\qquad\qquad\qquad\qquad\qquad\qquad\qquad S_{32}$

$+\frac{1}{24}[\frac{4}{3}\overset{q}{\underset{\text{a}=1}{\sum}}%
\perp_{\text{a}ki}R_{ij\text{a}j}-\frac{1}{3}(\nabla_{i}R_{kjij}+\nabla
_{j}R_{ijik}+\nabla_{k}R_{ijij})](y_{0})X_{k}(y_{0})$

$-\frac{1}{18}R_{ijkj}(y_{0})[X_{i}X_{k}-\frac{1}{2}\left(  \frac{\partial
X_{i}}{\partial x_{k}}+\frac{\partial X_{k}}{\partial x_{i}}\right)  ](y_{0})$

$+\frac{1}{24}[X_{i}^{2}X_{j}^{2}-2X_{i}X_{j}\left(  \frac{\partial X_{j}%
}{\partial x_{i}}+\frac{\partial X_{i}}{\partial x_{j}}\right)  -X_{i}%
^{2}\frac{\partial X_{j}}{\partial x_{j}}-X_{j}^{2}\frac{\partial X_{i}%
}{\partial x_{i}}](y_{0})$

$+\frac{1}{48}\left(  \frac{\partial X_{j}}{\partial x_{i}}+\frac{\partial
X_{i}}{\partial x_{j}}\right)  ^{2}(y_{0})+\frac{1}{24}\left(  \frac{\partial
X_{i}}{\partial x_{i}}\frac{\partial X_{j}}{\partial x_{j}}\right)
(y_{0})\qquad$

$+\frac{1}{36}X_{i}(y_{0})\left(  2\frac{\partial^{2}X_{j}}{\partial
x_{i}\partial x_{j}}+\frac{\partial^{2}X_{i}}{\partial x_{j}^{2}}\right)
(y_{0})+\frac{1}{36}X_{j}(y_{0})\left(  \frac{\partial^{2}X_{j}}{\partial
x_{i}^{2}}+2\frac{\partial^{2}X_{i}}{\partial x_{i}\partial x_{j}}\right)
(y_{0})$

$-\frac{1}{48}\left(  \frac{\partial^{3}X_{i}}{\partial x_{i}\partial
x_{j}^{2}}+\frac{\partial^{3}X_{j}}{\partial x_{i}^{2}\partial x_{j}}\right)
(y_{0})$

\qquad\qquad\qquad\qquad\qquad\qquad\qquad\qquad\qquad\qquad\qquad\qquad
\qquad\qquad\qquad$\blacksquare\qquad\qquad\qquad\qquad\qquad\qquad\qquad$

Next we recall that for any smooth function $f:M\longrightarrow R,$

$<X,\nabla^{0}f>$ $=$ $<X_{i}\frac{\partial}{\partial x_{i}},(\nabla^{0}%
f)_{k}\frac{\partial}{\partial x_{k}}>$ $=X_{i}(\nabla^{0}f)_{k}%
<\frac{\partial}{\partial x_{i}},\frac{\partial}{\partial x_{k}}>$
$=X_{i}(\nabla^{0}f)_{k}g_{ik}$

$\qquad\qquad=X_{i}(g^{jk}\frac{\partial f}{\partial x_{j}})g_{ik}=\delta
^{ij}X_{i}\frac{\partial f}{\partial x_{j}}=X_{j}\frac{\partial f}{\partial
x_{j}}$

\qquad$\ $I$_{3213}=\frac{1}{12}\frac{\partial^{2}}{\partial x_{i}^{2}%
}[<\nabla\log\theta^{-\frac{1}{2}},\nabla\log\Phi+X>](y_{0})\phi(y_{0})$

\qquad$\ \ \ =\frac{1}{12}\frac{\partial^{2}}{\partial x_{i}^{2}}%
[g_{jk}(\nabla\log\theta^{-\frac{1}{2}})_{j}((\nabla\log\Phi)_{k}%
+X_{k})](y_{0})\phi(y_{0})$

\qquad$\ \ \ =\frac{1}{12}\frac{\partial^{2}g_{jk}}{\partial x_{i}^{2}}%
(y_{0})[(\nabla\log\theta^{-\frac{1}{2}})_{j}((\nabla\log\Phi)_{k}%
+X_{k})](y_{0})\phi(y_{0})$

$\qquad+\frac{1}{12}g_{jk}(y_{0})\frac{\partial^{2}}{\partial x_{i}^{2}%
}[(\nabla\log\theta^{-\frac{1}{2}})_{j}((\nabla\log\Phi)_{k}+X_{k}%
)](y_{0})\phi(y_{0})$

\qquad$+\frac{1}{6}\frac{\partial g_{jk}}{\partial x_{i}}(y_{0})\frac
{\partial}{\partial x_{i}}[(\nabla\log\theta^{-\frac{1}{2}})_{j}((\nabla
\log\Phi)_{k}+X_{k})](y_{0})\phi(y_{0})$

$\qquad=$ I$_{32131}+$ I$_{32132}+$ I$_{32133}$

where,

\qquad\ I$_{32131}=\frac{1}{12}\frac{\partial^{2}g_{jk}}{\partial x_{i}^{2}%
}[(\nabla\log\theta^{-\frac{1}{2}})_{j}((\nabla\log\Phi)_{k}+X_{k}%
)](y_{0})\phi(y_{0})$

\qquad\ I$_{32132}=\frac{1}{12}g_{jk}(y_{0})\frac{\partial^{2}}{\partial
x_{i}^{2}}[(\nabla\log\theta^{-\frac{1}{2}})_{j}((\nabla\log\Phi)_{k}%
+X_{k})](y_{0})\phi(y_{0})$

\qquad\ I$_{32133}=\frac{1}{6}\frac{\partial g_{jk}}{\partial x_{i}}%
(y_{0})\frac{\partial}{\partial x_{i}}[(\nabla\log\theta^{-\frac{1}{2}}%
)_{j}((\nabla\log\Phi)_{k}+X_{k})](y_{0})\phi(y_{0})$

Since $(\nabla\log\Phi)_{k}(y_{0})+X_{k})(y_{0})=0,$

\qquad\ I$_{32131}=0$

Since $g_{jk}(y_{0})=\delta_{jk},$

\qquad\ I$_{32132}=\frac{1}{12}\frac{\partial^{2}}{\partial x_{i}^{2}}%
[(\nabla\log\theta^{-\frac{1}{2}})_{j}((\nabla\log\Phi)_{j}+X_{j})](y_{0}%
)\phi(y_{0})$

\qquad$\qquad=\frac{1}{12}\frac{\partial^{2}}{\partial x_{i}^{2}}(\nabla
\log\theta^{-\frac{1}{2}})_{j}(y_{0})[(\nabla\log\Phi)_{j}+X_{j})](y_{0}%
)\phi(y_{0})$

$\qquad\qquad+\frac{1}{12}[(\nabla\log\theta^{-\frac{1}{2}})_{j}(y_{0}%
)[\frac{\partial^{2}}{\partial x_{i}^{2}}(\nabla\log\Phi)_{j}+\frac
{\partial^{2}X_{j}}{\partial x_{i}^{2}}](y_{0})\phi(y_{0})$

\qquad$\qquad+\frac{1}{6}[\frac{\partial}{\partial x_{i}}(\nabla\log
\theta^{-\frac{1}{2}})_{j}(y_{0})[\frac{\partial}{\partial x_{i}}(\nabla
\log\Phi)_{j}+\frac{\partial X_{j}}{\partial x_{i}}](y_{0})\phi(y_{0})$

Since $[(\nabla\log\Phi)_{j}+X_{j})](y_{0})=0$ and $(\nabla\log\theta
)_{j}(y_{0})=-$ $<H,j>(y_{0})$

By (ix)$^{\ast\ast}$ of \textbf{Table A}$_{10},$ we have for $i,j=q+1,...n,$

\qquad\ I$_{32132}=-\frac{1}{12}$ $<H,j>(y_{0})[\frac{\partial^{2}}{\partial
x_{i}^{2}}(\nabla\log\Phi)_{j}+\frac{\partial^{2}X_{j}}{\partial x_{i}^{2}%
}](y_{0})\phi(y_{0})$

\qquad$+\frac{1}{6}[\frac{1}{2}<H,i><H,j>\ +\frac{1}{12}(2\varrho
_{ij}+4\overset{q}{\underset{\text{a}=1}{\sum}}R_{i\text{a}j\text{a}%
}-6\overset{q}{\underset{\text{a,b}=1}{\sum}}T_{\text{aa}i}T_{\text{bb}%
j}-T_{\text{ab}i}T_{\text{ab}j})](y_{0})$

$\qquad\times\lbrack\frac{\partial}{\partial x_{i}}(\nabla\log\Phi)_{j}%
+\frac{\partial X_{j}}{\partial x_{i}}](y_{0})\phi(y_{0})$

From $B_{39}$ above we have:

$\qquad\qquad\frac{\partial}{\partial x_{i}}(\nabla\log\Phi)_{j}(y_{0}%
)=-\frac{1}{2}\left(  \frac{\partial X_{j}}{\partial x_{i}}+\frac{\partial
X_{i}}{\partial x_{j}}\right)  (y_{0}),$

and from $B_{63}$ we have:

$[\frac{\partial^{2}}{\partial x_{j}\partial x_{i}}(\nabla\log\Phi_{P}%
)_{k}](y_{0})=-\frac{1}{3}\left(  \frac{\partial^{2}X_{j}}{\partial
x_{i}\partial x_{k}}+\frac{\partial^{2}X_{i}}{\partial x_{j}\partial x_{k}%
}+\frac{\partial^{2}X_{k}}{\partial x_{i}\partial x_{j}}\right)  (y_{0})$

$-\frac{1}{3}[R_{ikjl}+R_{jkil}](y_{0})X_{l}(y_{0})$

$\frac{\partial^{2}}{\partial x_{i}^{2}}(\nabla\log\Phi)_{j}(y_{0})=-\frac
{1}{3}\left(  \frac{\partial^{2}X_{i}}{\partial x_{i}\partial x_{j}}%
+\frac{\partial^{2}X_{i}}{\partial x_{i}\partial x_{j}}+\frac{\partial
^{2}X_{j}}{\partial x_{i}^{2}}\right)  (y_{0})-\frac{1}{3}[R_{ijil}%
+R_{ijil}](y_{0})X_{l}(y_{0})$

Therefore,

$\frac{\partial^{2}}{\partial x_{i}^{2}}(\nabla\log\Phi)_{j}(y_{0})=-\frac
{1}{3}\left(  2\frac{\partial^{2}X_{i}}{\partial x_{i}\partial x_{j}}%
+\frac{\partial^{2}X_{j}}{\partial x_{i}^{2}}\right)  (y_{0})-\frac{2}%
{3}R_{ijil}(y_{0})X_{l}(y_{0})$

Consequently,

\ I$_{32132}=-\frac{1}{12}$ $<H,j>(y_{0})[-\frac{1}{3}\left(  2\frac
{\partial^{2}X_{i}}{\partial x_{i}\partial x_{j}}+\frac{\partial^{2}X_{j}%
}{\partial x_{i}^{2}}\right)  +\frac{\partial^{2}X_{j}}{\partial x_{i}^{2}%
}](y_{0})\phi(y_{0})$

$\qquad\qquad+$ $\frac{1}{18}<H,j>(y_{0})R_{ijil}(y_{0})X_{l}(y_{0})\phi
(y_{0})$

$+\frac{1}{6}[\frac{1}{2}<H,i><H,j>\ +\frac{1}{12}(2\varrho_{ij}%
+4\overset{q}{\underset{\text{a}=1}{\sum}}R_{i\text{a}j\text{a}}%
-6\overset{q}{\underset{\text{a,b}=1}{\sum}}T_{\text{aa}i}T_{\text{bb}%
j}-T_{\text{ab}i}T_{\text{ab}j})](y_{0})\phi(y_{0})$

$\qquad\times\lbrack-\frac{1}{2}\left(  \frac{\partial X_{j}}{\partial x_{i}%
}+\frac{\partial X_{i}}{\partial x_{j}}\right)  +\frac{\partial X_{j}%
}{\partial x_{i}}](y_{0})\phi(y_{0})$

Simplifying, we have:

I$_{32132}=$ $\frac{1}{18}<H,j>(y_{0})[\frac{\partial^{2}X_{i}}{\partial
x_{i}\partial x_{j}}-\frac{\partial^{2}X_{j}}{\partial x_{i}^{2}}](y_{0}%
)\phi(y_{0})+$ $\frac{1}{18}<H,j>(y_{0})R_{ijik}(y_{0})X_{k}(y_{0})\phi
(y_{0})$

\qquad$+\frac{1}{24}[<H,i><H,j>\ +\frac{1}{6}(2\varrho_{ij}%
+4\overset{q}{\underset{\text{a}=1}{\sum}}R_{i\text{a}j\text{a}}%
-6\overset{q}{\underset{\text{a,b}=1}{\sum}}T_{\text{aa}i}T_{\text{bb}%
j}-T_{\text{ab}i}T_{\text{ab}j})](y_{0})\phi(y_{0})$

$\qquad\times\lbrack\frac{\partial X_{j}}{\partial x_{i}}-\frac{\partial
X_{i}}{\partial x_{j}}](y_{0})\phi(y_{0})$

Next we have:

\qquad I$_{32133}=\frac{1}{6}\frac{\partial g_{jk}}{\partial x_{i}}%
(y_{0})\frac{\partial}{\partial x_{i}}[(\nabla\log\theta^{-\frac{1}{2}}%
)_{j}((\nabla\log\Phi)_{k}+X_{k})](y_{0})\phi(y_{0})$

\qquad$=\frac{1}{6}\frac{\partial g_{jk}}{\partial x_{i}}(y_{0})\frac
{\partial}{\partial x_{i}}[(\nabla\log\theta^{-\frac{1}{2}})_{j}%
](y_{0})[(\nabla\log\Phi)_{k}+X_{k})](y_{0})\phi(y_{0})$

\qquad$\ +\frac{1}{6}\frac{\partial g_{jk}}{\partial x_{i}}(y_{0})[(\nabla
\log\theta^{-\frac{1}{2}})_{j}](y_{0})[\frac{\partial}{\partial x_{i}}%
(\nabla\log\Phi)_{k}+\frac{\partial X_{k}}{\partial x_{i}}](y_{0})\phi(y_{0})$

Since $(\nabla\log\Phi)_{k}(y_{0})=-X_{k}(y_{0}),$ we have: $(\nabla\log
\Phi)_{k}+X_{k})](y_{0})=0$

\qquad I$_{32133}=\frac{1}{6}\frac{\partial g_{jk}}{\partial x_{i}}%
(y_{0})[(\nabla\log\theta^{-\frac{1}{2}})_{j}](y_{0})[\frac{\partial}{\partial
x_{i}}(\nabla\log\Phi)_{k}+\frac{\partial X_{k}}{\partial x_{i}}](y_{0}%
)\phi(y_{0})$

Since $(\nabla\log\theta^{-\frac{1}{2}})_{\text{a}}(y_{0})=0$ by (iii)$^{\ast
}$ of \textbf{Table A}$_{9},$ we have,

$\frac{\partial g_{jk}}{\partial x_{\text{a}}}(y_{0})=0=(\nabla\log
\theta^{-\frac{1}{2}})_{\text{b}}](y_{0})$ for a,b =1,...,q.

We must therefore take $i,j=q+1,...,n.$

Now, $\frac{\partial g_{jk}}{\partial x_{i}}(y_{0})=0$ for $i,j,k=q+1,...,n.$
Therefore the only valid indices are:

$i,j=q+1,...,n$ and $k=$ a = 1,...,q. \ We have therefore,

I$_{32133}=\frac{1}{6}\frac{\partial g_{j\text{a}}}{\partial x_{i}}%
(y_{0})[(\nabla\log\theta^{-\frac{1}{2}})_{j}](y_{0})[\frac{\partial}{\partial
x_{i}}(\nabla\log\Phi)_{\text{a}}+\frac{\partial X_{\text{a}}}{\partial x_{i}%
}](y_{0})\phi(y_{0})$

for a = 1,...,q and $i,j=q+1,...,n.$

$\frac{\partial g_{j\text{a}}}{\partial x_{i}}(y_{0})=-\perp_{\text{a}ij}$ and
$(\nabla\log\theta^{-\frac{1}{2}})_{i}(y_{0})=\frac{1}{2}<H,i>(y_{0})$ by
(vi)$^{\ast}$ of \textbf{Table A}$_{9}$

By (xv) of \textbf{Table B}$_{1},$ $\frac{\partial}{\partial x_{i}}(\nabla
$log$\Phi_{P})_{\text{a}}(y)=$ $\left(  X_{j}\perp_{\text{a}ij}-\frac{\partial
X_{i}}{\partial x_{\text{a}}}\right)  (y_{0})$

Therefore, we have:

I$_{32133}=-\frac{1}{12}\perp_{\text{a}ij}(y_{0})<H,i>(y_{0})[X_{j}%
\perp_{\text{a}ij}-\frac{\partial X_{i}}{\partial x_{\text{a}}}+\frac{\partial
X_{\text{a}}}{\partial x_{i}}](y_{0})\phi(y_{0})$

Since I$_{31}=0$

$\left(  B_{115}\right)  $\qquad\ I$_{3213}=$ I$_{32131}+$ I$_{32132}+$
I$_{32133}$

$=$ $\frac{2}{3}<H,j>(y_{0})\left(  \frac{\partial^{2}X_{i}}{\partial
x_{i}\partial x_{j}}+2\frac{\partial^{2}X_{j}}{\partial x_{i}^{2}}\right)
(y_{0})\phi(y_{0})+$ $\frac{2}{3}<H,j>(y_{0})R_{ijik}(y_{0})X_{k}(y_{0}%
)\phi(y_{0})\qquad$I$_{3213}$

\qquad$+\frac{1}{12}[<H,i><H,j>\ +\frac{1}{6}(2\varrho_{ij}%
+4\overset{q}{\underset{\text{a}=1}{\sum}}R_{i\text{a}j\text{a}}%
-6\overset{q}{\underset{\text{a,b}=1}{\sum}}T_{\text{aa}i}T_{\text{bb}%
j}-T_{\text{ab}i}T_{\text{ab}j})](y_{0})\phi(y_{0})$

$\qquad\times\frac{1}{2}[\left(  \frac{\partial X_{j}}{\partial x_{i}}%
-\frac{\partial X_{i}}{\partial x_{j}}\right)  ](y_{0})\phi(y_{0})$

\qquad$-\frac{1}{12}\perp_{\text{a}ij}(y_{0})<H,i>(y_{0})[(X_{j}%
\perp_{\text{a}ij}-\frac{\partial X_{i}}{\partial x_{\text{a}}})+\frac
{\partial X_{\text{a}}}{\partial x_{i}}](y_{0})\phi(y_{0})$

We next compute:

\qquad I$_{3214}=\frac{1}{12}\frac{\partial^{2}}{\partial x_{i}^{2}}%
[<\nabla\log\Phi,X>](y_{0})\phi(y_{0})$

\qquad$=\frac{1}{12}\frac{\partial^{2}}{\partial x_{i}^{2}}[\Phi^{-1}%
<\nabla\Phi,X>](y_{0})\phi(y_{0})$

\qquad$=\frac{1}{12}\frac{\partial^{2}\Phi^{-1}}{\partial x_{i}^{2}}%
(y_{0})[<\nabla\Phi,X>](y_{0})\phi(y_{0})+\frac{1}{12}\Phi^{-1}(y_{0}%
)\frac{\partial^{2}}{\partial x_{i}^{2}}[<\nabla\Phi,X>](y_{0})\phi(y_{0})$

\qquad$+\frac{1}{6}\frac{\partial\Phi^{-1}}{\partial x_{i}}(y_{0}%
)\frac{\partial}{\partial x_{i}}[<\nabla\Phi,X>](y_{0})\phi(y_{0})=$
I$_{32141}+$ I$_{32142}+$ I$_{32143}$ where,

\qquad I$_{32141}=\frac{1}{12}\frac{\partial^{2}\Phi^{-1}}{\partial x_{i}^{2}%
}(y_{0})[<\nabla\Phi,X>](y_{0})\phi(y_{0})$

\qquad I$_{32142}=\frac{1}{12}\Phi^{-1}(y_{0})\frac{\partial^{2}}{\partial
x_{i}^{2}}[<\nabla\Phi,X>](y_{0})\phi(y_{0})$

\qquad I$_{32143}=\frac{1}{6}\frac{\partial\Phi^{-1}}{\partial x_{i}}%
(y_{0})\frac{\partial}{\partial x_{i}}[<\nabla\Phi,X>](y_{0})\phi(y_{0})$

\qquad\qquad\qquad\qquad\qquad\qquad\qquad\qquad\qquad\qquad\qquad
\qquad$\blacksquare$

Since $<X,\nabla f>$ $=X_{j}(\nabla f)_{k}g_{jk},$ we for $i=q+1,....n$ and
$j,k=1,....,q,q+1,...,n,$

I$_{3214}=\frac{1}{12}\frac{\partial^{2}}{\partial x_{i}^{2}}[<\nabla\log
\Phi,X>](y_{0})\phi(y_{0})=\frac{1}{12}\frac{\partial^{2}}{\partial x_{i}^{2}%
}[X_{j}(\nabla\log\Phi)_{k}g_{jk}](y_{0})\phi(y_{0})$

\qquad$=\frac{1}{12}\frac{\partial^{2}X_{j}}{\partial x_{i}^{2}}(y_{0}%
)[\nabla\log\Phi)_{k}g_{jk}](y_{0})\phi(y_{0})+\frac{1}{12}X_{j}(y_{0}%
)\frac{\partial^{2}}{\partial x_{i}^{2}}[\nabla\log\Phi)_{k}g_{jk}](y_{0}%
)\phi(y_{0})$

\qquad$+\frac{1}{6}\frac{\partial X_{j}}{\partial x_{i}}(y_{0})\frac{\partial
}{\partial x_{i}}[(\nabla\log\Phi)_{k}g_{jk}](y_{0})\phi(y_{0})=$ I$_{32141}+$
I$_{32142}+$ I$_{32143}$

where,

I$_{32141}=\frac{1}{12}\frac{\partial^{2}X_{j}}{\partial x_{i}^{2}}%
(y_{0})[\nabla\log\Phi)_{k}g_{jk}](y_{0})\phi(y_{0});$ I$_{32142}$

$=\frac{1}{12}X_{j}(y_{0})\frac{\partial^{2}}{\partial x_{i}^{2}}[\nabla
\log\Phi)_{k}g_{jk}](y_{0})\phi(y_{0})$

I$_{32143}=\frac{1}{6}\frac{\partial X_{j}}{\partial x_{i}}(y_{0}%
)\frac{\partial}{\partial x_{i}}[(\nabla\log\Phi)_{k}g_{jk}](y_{0})\phi
(y_{0})$

Since $g_{jk}(y_{0})=\delta_{jk},$ we have,

I$_{32141}=\frac{1}{12}\frac{\partial^{2}X_{j}}{\partial x_{i}^{2}}%
(y_{0})[\nabla\log\Phi)_{k}g_{jk}](y_{0})\phi(y_{0})=\frac{1}{12}%
\frac{\partial^{2}X_{j}}{\partial x_{i}^{2}}(y_{0})(\nabla\log\Phi)_{j}%
(y_{0})\phi(y_{0})$

Using the Einstein convention for repeated indices, we have

for a = 1,....,q and $i,j=q+1,...,n,$

I$_{32141}=\frac{1}{12}\frac{\partial^{2}X_{j}}{\partial x_{i}^{2}}%
(y_{0})(\nabla\log\Phi)_{j}(y_{0})\phi(y_{0})$

\qquad$\ =\frac{1}{12}\frac{\partial^{2}X_{\text{a}}}{\partial x_{i}^{2}%
}(y_{0})(\nabla\log\Phi)_{\text{a}}(y_{0})\phi(y_{0})+\frac{1}{12}%
\frac{\partial^{2}X_{j}}{\partial x_{i}^{2}}(y_{0})(\nabla\log\Phi)_{j}%
(y_{0})\phi(y_{0})$

Since $(\nabla\log\Phi)_{\text{a}}(y_{0})=0$ for a =1,...,q and $(\nabla
\log\Phi)_{j}(y_{0})=-X_{j}(y_{0})$ for $j=q+1,...,n,$

I$_{32141}=-\frac{1}{12}\frac{\partial^{2}X_{j}}{\partial x_{i}^{2}}%
(y_{0})X_{j}(y_{0})\phi(y_{0})=-\frac{1}{12}X_{j}(y_{0})\frac{\partial
^{2}X_{j}}{\partial x_{i}^{2}}(y_{0})\phi(y_{0})$ for $i,j=q+1,...,n.$

Next, since $g_{jk}(y_{0})=\delta_{jk},$ we have for $i=q+1,...,n$ and
$j,k=1,....,q,q+1,...,n,$

I$_{32142}=\frac{1}{12}X_{j}(y_{0})\frac{\partial^{2}}{\partial x_{i}^{2}%
}[\nabla\log\Phi)_{k}g_{jk}](y_{0})\phi(y_{0})$

\qquad$=\frac{1}{12}X_{j}(y_{0})\frac{\partial^{2}}{\partial x_{i}^{2}}%
(\nabla\log\Phi)_{j})(y_{0})\phi(y_{0})+\frac{1}{12}X_{j}(y_{0})(\nabla
\log\Phi)_{k}(y_{0})\frac{\partial^{2}g_{jk}}{\partial x_{i}^{2}}(y_{0}%
)\phi(y_{0})$

$\qquad\ \ \ \ +\frac{1}{6}X_{j}(y_{0})\frac{\partial}{\partial x_{i}}%
(\nabla\log\Phi)_{k})(y_{0})\frac{\partial g_{jk}}{\partial x_{i}}\phi
(y_{0})=$ I$_{321421}+$ I$_{321422}+$ I$_{321423}$

I$_{321421}=\frac{1}{12}X_{j}(y_{0})\frac{\partial^{2}}{\partial x_{i}^{2}%
}(\nabla\log\Phi)_{j})(y_{0})\phi(y_{0})$

Recalling that the Einstein convention for repeated indices apply, we have for:

a = 1,...,q; $i=q+1,...,n,$ and $j=1,....,q,q+1,...,n,$

I$_{321421}=\frac{1}{12}X_{\text{a}}(y_{0})\frac{\partial^{2}}{\partial
x_{i}^{2}}(\nabla\log\Phi)_{\text{a}})(y_{0})\phi(y_{0})+\frac{1}{12}%
X_{j}(y_{0})\frac{\partial^{2}}{\partial x_{i}^{2}}(\nabla\log\Phi)_{j}%
)(y_{0})\phi(y_{0})$

Then for a,b = 1,...,q and $i,j=q+1,...,n,$ we have:

\ I$_{321421}=\frac{1}{12}X_{\text{a}}(y_{0})\frac{\partial^{2}}{\partial
x_{i}^{2}}(\nabla\log\Phi)_{\text{a}})(y_{0})\phi(y_{0})+\frac{1}{12}%
X_{j}(y_{0})\frac{\partial^{2}}{\partial x_{i}^{2}}(\nabla\log\Phi)_{j}%
)(y_{0})\phi(y_{0})$

By (xvi) of \textbf{Table B}$_{1},$ and (v)$^{\ast\ast\ast}$ of \textbf{Table
B}$_{4},$ we have:

$\qquad\frac{\partial^{2}}{\partial x_{i}^{2}}(\nabla\log\Phi)_{\text{a}%
})(y_{0})=-4T_{\text{ab}i}(y_{0})\frac{\partial X_{i}}{\partial x_{\text{b}}%
}(y_{0})+$ $\perp_{\text{a}ij}(y_{0})\left[  \left(  \frac{\partial X_{i}%
}{\partial x_{j}}+\frac{\partial X_{j}}{\partial x_{i}}\right)  \right]
(y_{0})$

$\qquad\qquad\qquad\qquad\qquad+\frac{8}{3}R_{i\text{a}ij}(y)X_{j}%
(y_{0})+[2X_{i}\frac{\partial X_{i}}{\partial x_{\text{a}}}-\frac{\partial
^{2}X_{i}}{\partial x_{\text{a}}\partial x_{i}}](y_{0})$

By (v)$^{\ast\ast\ast}$ of \textbf{Table B}$_{4},$ we have:

$\qquad$

$\qquad\frac{\partial^{2}}{\partial x_{i}^{2}}(\nabla\log\Phi)_{j}%
)(y_{0})=-\frac{1}{3}\left(  \frac{\partial^{2}X_{j}}{\partial x_{i}^{2}%
}+2\frac{\partial^{2}X_{i}}{\partial x_{i}\partial x_{j}}\right)
(y_{0})-\frac{2}{3}\underset{l=q+1}{\overset{n}{\sum}}R_{ijik}(y_{0}%
)X_{k}(y_{0})\qquad$

\qquad\qquad\qquad\qquad\qquad$\qquad+[2\perp_{\text{a}ij}\frac{\partial
X_{i}}{\partial x_{\text{a}}}](y_{0})+[2\perp_{\text{a}ij}\perp_{\text{a}%
ik}X_{k}](y_{0})$

Therefore,

\ I$_{321421}=-\frac{1}{3}T_{\text{ab}i}(y_{0})X_{\text{a}}(y_{0}%
)\frac{\partial X_{i}}{\partial x_{\text{b}}}(y_{0})+$ $\frac{1}{12}%
\perp_{\text{a}ij}(y_{0})X_{\text{a}}(y_{0})\left[  \left(  \frac{\partial
X_{i}}{\partial x_{j}}+\frac{\partial X_{j}}{\partial x_{i}}\right)  \right]
(y_{0})$

$\qquad\qquad\qquad\qquad+\frac{2}{9}R_{i\text{a}ij}(y)X_{\text{a}}%
(y_{0})X_{j}(y_{0})+\frac{1}{12}X_{\text{a}}(y_{0})[2X_{i}\frac{\partial
X_{i}}{\partial x_{\text{a}}}-\frac{\partial^{2}X_{i}}{\partial x_{\text{a}%
}\partial x_{i}}](y_{0})$

$\qquad\qquad-\frac{1}{36}X_{j}(y_{0})\left(  \frac{\partial^{2}X_{j}%
}{\partial x_{i}^{2}}+2\frac{\partial^{2}X_{i}}{\partial x_{i}\partial x_{j}%
}\right)  (y_{0})-\frac{1}{18}\underset{l=q+1}{\overset{n}{\sum}}%
R_{ijik}(y_{0})X_{j}(y_{0})X_{k}(y_{0})$

$\qquad\qquad+\frac{1}{12}X_{j}(y_{0})[2\perp_{\text{a}ij}\frac{\partial
X_{i}}{\partial x_{\text{a}}}](y_{0})+\frac{1}{12}X_{j}(y_{0})[2\perp
_{\text{a}ij}\perp_{\text{a}ik}X_{k}](y_{0})$

Next we have for $i=q+1,...,n,$ and $j,k=1,....,q,q+1,...,n,$

I$_{321422}=\frac{1}{12}X_{j}(y_{0})(\nabla\log\Phi)_{k}(y_{0})\frac
{\partial^{2}g_{jk}}{\partial x_{i}^{2}}(y_{0})\phi(y_{0})$

Then, for a = 1,...,q; $i,j=q+1,...,n,$ and $k=1,....,q,q+1,...,n,$

I$_{321422}=\frac{1}{12}X_{\text{a}}(y_{0})(\nabla\log\Phi)_{k}(y_{0}%
)\frac{\partial^{2}g_{\text{a}k}}{\partial x_{i}^{2}}(y_{0})\phi(y_{0}%
)+\frac{1}{12}X_{j}(y_{0})(\nabla\log\Phi)_{k}(y_{0})\frac{\partial^{2}g_{jk}%
}{\partial x_{i}^{2}}(y_{0})\phi(y_{0})$

Then, finally for a,b = 1,...,q; $i,j,k=q+1,...,n,$ we have:

I$_{321422}=\frac{1}{12}X_{\text{a}}(y_{0})(\nabla\log\Phi)_{\text{b}}%
(y_{0})\frac{\partial^{2}g_{\text{ab}}}{\partial x_{i}^{2}}(y_{0})\phi
(y_{0})+\frac{1}{12}X_{\text{a}}(y_{0})(\nabla\log\Phi)_{k}(y_{0}%
)\frac{\partial^{2}g_{\text{a}k}}{\partial x_{i}^{2}}(y_{0})\phi(y_{0})$

$\qquad\ \ \ \ +\frac{1}{12}X_{j}(y_{0})(\nabla\log\Phi)_{\text{b}}%
(y_{0})\frac{\partial^{2}g_{j\text{b}}}{\partial x_{i}^{2}}(y_{0})\phi
(y_{0})+\frac{1}{12}X_{j}(y_{0})(\nabla\log\Phi)_{k}(y_{0})\frac{\partial
^{2}g_{jk}}{\partial x_{i}^{2}}(y_{0})\phi(y_{0})$

Since $(\nabla\log\Phi)_{\text{b}}(y_{0})=0,$ we have,

I$_{321422}=\frac{1}{12}X_{\text{a}}(y_{0})(\nabla\log\Phi)_{k}(y_{0}%
)\frac{\partial^{2}g_{\text{a}k}}{\partial x_{i}^{2}}(y_{0})\phi(y_{0}%
)+\frac{1}{12}X_{j}(y_{0})(\nabla\log\Phi)_{k}(y_{0})\frac{\partial^{2}g_{jk}%
}{\partial x_{i}^{2}}(y_{0})\phi(y_{0})$

$\frac{\partial^{2}g_{\text{a}k}}{\partial x_{i}^{2}}(y_{0})=-\frac{8}%
{3}R_{i\text{a}ik}(y_{0})$ by (iii)$^{\ast}$ of \textbf{Table A}$_{3}$ and
$\frac{\partial^{2}g_{jk}}{\partial x_{i}^{2}}(y_{0})\phi(y_{0})$

$=-\frac{2}{3}(R_{ijik})(y_{0})$ by (iii) of \textbf{Table A}$_{1}$

Consequently we have,

I$_{321422}=\frac{2}{9}X_{\text{a}}(y_{0})X_{k}(y_{0})R_{i\text{a}ik}%
(y_{0})\phi(y_{0})+\frac{1}{18}X_{j}(y_{0})X_{k}(y_{0})R_{ijik}(y_{0}%
)\phi(y_{0})$

Next, we have for $i=q+1,...,n,$ and $j,k=1,....,q,q+1,...,n,$

I$_{321423}=\frac{1}{6}X_{j}(y_{0})\frac{\partial}{\partial x_{i}}(\nabla
\log\Phi)_{k})(y_{0})\frac{\partial g_{jk}}{\partial x_{i}}\phi(y_{0})$

Then, for a = 1,...,q; $i,j=q+1,...,n,$ and $k=1,....,q,q+1,...,n,$

I$_{321423}=\frac{1}{6}X_{\text{a}}(y_{0})\frac{\partial}{\partial x_{i}%
}(\nabla\log\Phi)_{k})(y_{0})\frac{\partial g_{\text{a}k}}{\partial x_{i}}%
\phi(y_{0})+\frac{1}{6}X_{j}(y_{0})\frac{\partial}{\partial x_{i}}(\nabla
\log\Phi)_{k})(y_{0})\frac{\partial g_{jk}}{\partial x_{i}}\phi(y_{0})$

Finally here, we have for a,b = 1,...,q and $i,j,k=q+1,...,n,$ we have:

I$_{321423}=\frac{1}{6}X_{\text{a}}(y_{0})\frac{\partial}{\partial x_{i}%
}(\nabla\log\Phi)_{\text{b}})(y_{0})\frac{\partial g_{\text{ab}}}{\partial
x_{i}}\phi(y_{0})+\frac{1}{6}X_{\text{a}}(y_{0})\frac{\partial}{\partial
x_{i}}(\nabla\log\Phi)_{k})(y_{0})\frac{\partial g_{\text{a}k}}{\partial
x_{i}}\phi(y_{0})$

$+\frac{1}{6}X_{j}(y_{0})\frac{\partial}{\partial x_{i}}(\nabla\log
\Phi)_{\text{b}})(y_{0})\frac{\partial g_{j\text{b}}}{\partial x_{i}}%
\phi(y_{0})+\frac{1}{6}X_{j}(y_{0})\frac{\partial}{\partial x_{i}}(\nabla
\log\Phi)_{k})(y_{0})\frac{\partial g_{jk}}{\partial x_{i}}\phi(y_{0})$

Now, $\frac{\partial g_{jk}}{\partial x_{i}}(y_{0})=0$ for $i,j,k=q+1,...,n,$

(xv) $\frac{\partial}{\partial x_{j}}(\nabla$log$\Phi_{P})_{\text{a}}%
(y_{0})=\underset{k=q+1}{\overset{n}{\sum}}$ $X_{k}(y_{0})\perp_{\text{a}%
kj}(y_{0})-\frac{\partial X_{j}}{\partial x_{\text{a}}}(y_{0})$

$=-[$ $X_{k}(y_{0})\perp_{\text{a}jk}(y_{0})+\frac{\partial X_{j}}{\partial
x_{\text{a}}}](y_{0})$

\subsubsection{Tangential Derivatives}

(ix) For a $=$ 1,...,q and for $j=q+1,...,n,$

$\qquad\frac{\partial}{\partial x_{\text{a}}}(\nabla$log$\Phi_{P}%
)_{j}(y)=-\frac{\partial X_{j}}{\partial x_{\text{a}}}(y)\qquad$\qquad

(x) For a,b = 1,...,q,

\qquad$\frac{\partial^{2}}{\partial x_{\text{a}}\partial x_{\text{b}}}(\nabla
$log$\Phi_{P})_{j}(y)=-\frac{\partial^{2}X_{j}}{\partial x_{\text{a}}\partial
x_{\text{b}}}(y)$

Fomulae for higher derivatives follow.

(xi) For a = 1,...,q \ and $y\in U\subset P,$ we have:

\qquad$(\nabla$log$\Phi_{P})_{\text{a}}(y)=0$

(xii) For a, b = 1,...,q$\ $

$\qquad\frac{\partial}{\partial x_{\text{b}}}(\nabla$log$\Phi_{P})_{\text{a}%
}(y)\ =0$

(xiii) For a, b, c = 1,...,q,

\qquad$\ \frac{\partial^{2}}{\partial x_{\text{c}}\partial x_{\text{b}}%
}(\nabla$log$\Phi_{P})_{\text{a}}(y)\ =0\qquad$

\subsubsection{\textbf{Mixed Derivatives: }}

\textbf{For a =1,...,q and }$i,j,k=q+1,...,n:$

(xiv) $\frac{\partial^{2}}{\partial x_{\text{a}}\partial x_{k}}\nabla\log
\Phi_{P})_{j}(y)+\frac{\partial^{2}}{\partial x_{\text{a}}\partial x_{j}%
}\nabla\log\Phi_{P})_{k}(y)=-\frac{\partial^{2}X_{j}}{\partial x_{\text{a}%
}\partial x_{k}}(y)-$ $\frac{\partial^{2}X_{k}}{\partial x_{\text{a}}\partial
x_{j}}(y).$

In particular for $k=j,$

\qquad$\frac{\partial^{2}}{\partial x_{\text{a}}\partial x_{j}}\nabla\log
\Phi_{P})_{j}(y)=-\frac{\partial^{2}X_{j}}{\partial x_{\text{a}}\partial
x_{j}}(y)$

(xv) $\frac{\partial}{\partial x_{j}}(\nabla$log$\Phi_{P})_{\text{a}}(y)=$
$\underset{i=q+1}{\overset{n}{%
{\textstyle\sum}
}}X_{i}(y)\perp_{\text{a}ij}(y)-\frac{\partial X_{j}}{\partial x_{\text{a}}%
}(y)$

(xvi) $\ \frac{\partial^{2}}{\partial x_{i}\partial x_{j}}(\nabla$log$\Phi
_{P})_{\text{a}}(y)=-$ $2\underset{\text{b=1}}{\overset{\text{q}}{%
{\textstyle\sum}
}}$T$_{\text{ab}j}(y_{0})\frac{\partial X_{i}}{\partial x_{\text{b}}%
}(y)-2\underset{\text{b=1}}{\overset{\text{q}}{%
{\textstyle\sum}
}}T_{\text{ab}i}(y)\frac{\partial X_{j}}{\partial x_{\text{b}}}(y)$

$\qquad+\frac{1}{2}\underset{k=q+1}{\overset{n}{%
{\textstyle\sum}
}}\perp_{\text{a}jk}(y)\left[  \left(  \frac{\partial X_{i}}{\partial x_{k}%
}+\frac{\partial X_{k}}{\partial x_{i}}\right)  \right]  (y)+\frac{1}%
{2}\underset{k=q+1}{\overset{n}{%
{\textstyle\sum}
}}\perp_{\text{a}ik}(y)[\left(  \frac{\partial X_{k}}{\partial x_{j}}%
+\frac{\partial X_{j}}{\partial x_{k}}\right)  ](y)$\qquad

$\qquad+\frac{4}{3}\underset{k=q+1}{\overset{n}{%
{\textstyle\sum}
}}\left[  R_{i\text{a}jk}+R_{j\text{a}ik}\right]  (y)X_{k}(y)+[X_{i}%
\frac{\partial X_{j}}{\partial x_{\text{a}}}+X_{j}\frac{\partial X_{i}%
}{\partial x_{\text{a}}}-\frac{1}{2}\left(  \frac{\partial^{2}X_{i}}{\partial
x_{\text{a}}\partial x_{j}}+\frac{\partial^{2}X_{j}}{\partial x_{\text{a}%
}\partial x_{i}}\right)  ](y)$\qquad\qquad

In particular, taking $j=i,$ we have:

\qquad$\frac{\partial^{2}}{\partial x_{i}^{2}}(\nabla$log$\Phi_{P})_{\text{a}%
}(y)=-4\underset{\text{b=1}}{\overset{\text{q}}{%
{\textstyle\sum}
}}T_{\text{ab}i}(y)\frac{\partial X_{i}}{\partial x_{\text{b}}}%
(y)+\underset{k=q+1}{\overset{n}{%
{\textstyle\sum}
}}\perp_{\text{a}ik}(y)\left[  \left(  \frac{\partial X_{i}}{\partial x_{k}%
}+\frac{\partial X_{k}}{\partial x_{i}}\right)  \right]  (y)$

$\qquad\qquad\qquad\qquad\qquad\ \ +\frac{8}{3}\underset{k=q+1}{\overset{n}{%
{\textstyle\sum}
}}R_{i\text{a}ik}(y)X_{k}(y)+[2X_{i}\frac{\partial X_{i}}{\partial
x_{\text{a}}}-\frac{\partial^{2}X_{i}}{\partial x_{\text{a}}\partial x_{i}%
}](y)$

(v)$^{\ast\ast\ast}$ of \textbf{Table B}$_{4}$

$[\frac{\partial^{2}}{\partial x_{i}\partial x_{j}}(\nabla\log\Phi_{P}%
)_{k}](y)$

$\qquad=-\frac{1}{3}\left(  \frac{\partial^{2}X_{i}}{\partial x_{j}\partial
x_{k}}+\frac{\partial^{2}X_{j}}{\partial x_{i}\partial x_{k}}+\frac
{\partial^{2}X_{k}}{\partial x_{i}\partial x_{j}}\right)  (y)+\frac{1}%
{3}(R_{kjil}+R_{kijl})(y)X_{l}(y)$

$\qquad\qquad-[\perp_{\text{a}kj}\perp_{\text{a}il}X_{l}+\perp_{\text{a}%
kj}\frac{\partial X_{i}}{\partial x_{\text{a}}}](y)-[\perp_{\text{a}ki}%
\perp_{\text{a}jl}X_{l}+\perp_{\text{a}ki}\frac{\partial X_{j}}{\partial
x_{\text{a}}}](y)$

\qquad\qquad$=-\frac{1}{3}\left(  \frac{\partial^{2}X_{k}}{\partial
x_{i}\partial x_{j}}+\frac{\partial^{2}X_{j}}{\partial x_{i}\partial x_{k}%
}+\frac{\partial^{2}X_{i}}{\partial x_{j}\partial x_{k}}\right)  (y)-\frac
{1}{3}(R_{ikjl}+R_{jkil})(y)X_{l}(y)$

$\qquad\qquad+[\perp_{\text{a}jk}\perp_{\text{a}il}X_{l}+\perp_{\text{a}%
jk}\frac{\partial X_{i}}{\partial x_{\text{a}}}](y)+[\perp_{\text{a}ik}%
\perp_{\text{a}jl}X_{l}+\perp_{\text{a}ik}\frac{\partial X_{j}}{\partial
x_{\text{a}}}](y)$

\qquad\qquad\qquad\qquad\qquad\qquad\qquad\qquad\qquad\qquad\qquad\qquad
\qquad\qquad\qquad\qquad\qquad\qquad$\blacksquare$

We have thus proved the formula in (v)$^{\ast\ast\ast}$ of \textbf{Table
B}$_{4}$ in \textbf{Appendix B}.

In particular, we have at the centre of Fermi coordinates $y_{0}\in P:$

$[\frac{\partial^{2}}{\partial x_{i}\partial x_{j}}(\nabla\log\Phi_{P}%
)_{k}](y_{0})$

$=-\frac{1}{3}\left(  \frac{\partial^{2}X_{k}}{\partial x_{i}\partial x_{j}%
}+\frac{\partial^{2}X_{j}}{\partial x_{i}\partial x_{k}}+\frac{\partial
^{2}X_{i}}{\partial x_{j}\partial x_{k}}\right)  (y_{0})-\frac{1}%
{3}\underset{l=q+1}{\overset{n}{\sum}}(R_{ikjl}+R_{jkil})(y_{0})X_{l}%
(y_{0})\qquad\left(  140\right)  $

\qquad$+[\perp_{\text{a}ik}\frac{\partial X_{j}}{\partial x_{\text{a}}}%
+\perp_{\text{a}jk}\frac{\partial X_{i}}{\partial x_{\text{a}}}](y_{0}%
)+\underset{l=q+1}{\overset{n}{\sum}}[\perp_{\text{a}ik}\perp_{\text{a}%
jl}X_{l}+\perp_{\text{a}jk}\perp_{\text{a}il}X_{l}](y_{0})$

\qquad\qquad\qquad\qquad\qquad\qquad\qquad\qquad\qquad\qquad\qquad\qquad
\qquad\qquad\qquad\qquad\qquad\qquad\qquad\qquad$\blacksquare$

In particular,

$[\frac{\partial^{2}}{\partial x_{i}\partial x_{j}}(\nabla\log\Phi_{P}%
)_{j}](y_{0})$

$=-\frac{1}{3}\left(  \frac{\partial^{2}X_{j}}{\partial x_{i}\partial x_{j}%
}+\frac{\partial^{2}X_{j}}{\partial x_{i}\partial x_{j}}+\frac{\partial
^{2}X_{i}}{\partial x_{j}^{2}}\right)  (y_{0})-\frac{1}{3}%
\underset{l=q+1}{\overset{n}{\sum}}(R_{ijjl}+R_{jjil})(y_{0})X_{l}%
(y_{0})\qquad$

\qquad$+[\perp_{\text{a}ij}\frac{\partial X_{j}}{\partial x_{\text{a}}}%
+\perp_{\text{a}jj}\frac{\partial X_{i}}{\partial x_{\text{a}}}](y_{0}%
)+\underset{l=q+1}{\overset{n}{\sum}}[\perp_{\text{a}ij}\perp_{\text{a}%
jl}X_{l}+\perp_{\text{a}jj}\perp_{\text{a}il}X_{l}](y_{0})$

\qquad$=-\frac{1}{3}\left(  2\frac{\partial^{2}X_{j}}{\partial x_{i}\partial
x_{j}}+\frac{\partial^{2}X_{i}}{\partial x_{j}^{2}}\right)  (y_{0})-\frac
{1}{3}\underset{k=q+1}{\overset{n}{\sum}}R_{ijjk}(y_{0})X_{k}(y_{0})$

\qquad\qquad$+[\perp_{\text{a}ij}\frac{\partial X_{j}}{\partial x_{\text{a}}%
}](y_{0})+\underset{k=q+1}{\overset{n}{\sum}}[\perp_{\text{a}ij}%
\perp_{\text{a}jk}X_{k}](y_{0})$

$[\frac{\partial^{2}}{\partial x_{i}^{2}}(\nabla\log\Phi_{P})_{k}](y_{0})$

$=-\frac{1}{3}\left(  \frac{\partial^{2}X_{k}}{\partial x_{i}^{2}}%
+\frac{\partial^{2}X_{i}}{\partial x_{i}\partial x_{k}}+\frac{\partial
^{2}X_{i}}{\partial x_{i}\partial x_{k}}\right)  (y_{0})-\frac{1}%
{3}\underset{l=q+1}{\overset{n}{\sum}}(R_{ikil}+R_{ikil})(y_{0})X_{l}%
(y_{0})\qquad$

\qquad$+[\perp_{\text{a}ik}\frac{\partial X_{i}}{\partial x_{\text{a}}}%
+\perp_{\text{a}ik}\frac{\partial X_{i}}{\partial x_{\text{a}}}](y_{0}%
)+\underset{l=q+1}{\overset{n}{\sum}}[\perp_{\text{a}ik}\perp_{\text{a}%
il}X_{l}+\perp_{\text{a}ik}\perp_{\text{a}il}X_{l}](y_{0})$

$[\frac{\partial^{2}}{\partial x_{i}^{2}}(\nabla\log\Phi_{P})_{k}%
](y_{0})=-\frac{1}{3}\left(  \frac{\partial^{2}X_{k}}{\partial x_{i}^{2}%
}+2\frac{\partial^{2}X_{i}}{\partial x_{i}\partial x_{k}}\right)
(y_{0})-\frac{2}{3}\underset{l=q+1}{\overset{n}{\sum}}R_{ikil}(y_{0}%
)X_{l}(y_{0})\qquad$

\qquad\qquad\qquad\qquad\qquad$\qquad+[2\perp_{\text{a}ik}\frac{\partial
X_{i}}{\partial x_{\text{a}}}](y_{0})+\underset{l=q+1}{\overset{n}{\sum}%
}[2\perp_{\text{a}ik}\perp_{\text{a}il}X_{l}](y_{0})$

$[\frac{\partial^{2}}{\partial x_{i}^{2}}(\nabla\log\Phi_{P})_{j}%
](y_{0})=-\frac{1}{3}\left(  \frac{\partial^{2}X_{j}}{\partial x_{i}^{2}%
}+2\frac{\partial^{2}X_{i}}{\partial x_{i}\partial x_{j}}\right)
(y_{0})-\frac{2}{3}\underset{k=q+1}{\overset{n}{\sum}}R_{ijik}(y_{0}%
)X_{k}(y_{0})\qquad$

\qquad\qquad\qquad\qquad\qquad$\qquad+[2\perp_{\text{a}ij}\frac{\partial
X_{i}}{\partial x_{\text{a}}}](y_{0})+\underset{k=q+1}{\overset{n}{\sum}%
}[2\perp_{\text{a}ij}\perp_{\text{a}ik}X_{k}](y_{0})$

\textbf{Table B}$_{4}$

\textbf{Mixed Derivatives:}

\qquad(xi)\qquad$\frac{\partial^{2}\Phi_{P}}{\partial x_{\text{a}}\partial
x_{i}}(y_{0})=$ $\frac{\partial}{\partial x_{\text{a}}}(\nabla\log\Phi
_{P})_{i}(y_{0})=-\frac{\partial X_{i}}{\partial x_{\text{a}}}(y_{0})$

\qquad(xii)\qquad$\frac{\partial^{3}\Phi_{P}}{\partial x_{\text{a}}\partial
x_{\text{b}}\partial x_{i}}(x_{0})=\frac{\partial^{2}}{\partial x_{\text{a}%
}\partial x_{\text{b}}}(\nabla\log\Phi_{P})_{i}(y_{0})=-\frac{\partial
^{2}X_{i}}{\partial x_{\text{a}}\partial x_{\text{b}}}(y_{0})$

\qquad(xiii)$\qquad\frac{\partial^{3}\Phi_{P}}{\partial x_{\text{a}}%
^{2}\partial x_{i}}(y_{0})=\frac{\partial^{2}}{\partial x_{\text{a}}^{2}%
}(\nabla\log\Phi_{P})_{i}(y_{0})=-\frac{\partial^{2}X_{i}}{\partial
x_{\text{a}}^{2}}(y_{0})$

\qquad(xiv) From $\left(  B_{88}\right)  ,$

$\qquad\frac{\partial^{3}\Phi_{P}}{\partial x_{\text{c}}\partial x_{i}\partial
x_{j}}(y_{0})$

$\qquad\qquad=[X_{i}\frac{\partial X_{j}}{\partial x_{\text{c}}}+X_{j}%
\frac{\partial X_{i}}{\partial x_{\text{c}}}](y_{0})+[X_{j}\frac{\partial
X_{i}}{\partial x_{\text{c}}}+X_{i}\frac{\partial X_{j}}{\partial x_{\text{c}%
}}](y_{0})-\frac{1}{2}\left(  \frac{\partial^{2}X_{i}}{\partial x_{\text{c}%
}\partial x_{j}}+\frac{\partial^{2}X_{j}}{\partial x_{\text{c}}\partial x_{i}%
}\right)  (y_{0})$

$\qquad\qquad\qquad=2[X_{i}\frac{\partial X_{j}}{\partial x_{\text{c}}}%
+X_{j}\frac{\partial X_{i}}{\partial x_{\text{c}}}](y_{0})-\frac{1}{2}\left(
\frac{\partial^{2}X_{i}}{\partial x_{\text{c}}\partial x_{j}}+\frac
{\partial^{2}X_{j}}{\partial x_{\text{c}}\partial x_{i}}\right)  (y_{0}%
)\qquad\qquad$\qquad

\qquad(xv)\qquad$\frac{\partial^{4}\Phi_{P}}{\partial x_{\text{c}}^{2}\partial
x_{i}\partial x_{j}}(y_{0})$

$\qquad\qquad\qquad=$ $2\frac{\partial X_{i}}{\partial x_{\text{c}}}%
(y_{0})\frac{\partial X_{j}}{\partial x_{\text{c}}}(y_{0})+X_{i}(y_{0}%
)\frac{\partial^{2}X_{j}}{\partial x_{\text{c}}^{2}}(y_{0})+X_{j}(y_{0}%
)\frac{\partial^{2}X_{i}}{\partial x_{\text{c}}^{2}}(y_{0})$ $-\frac
{\partial^{3}X_{i}}{\partial x_{\text{c}}^{2}\partial x_{j}}(y_{0}%
)\qquad\qquad\qquad\qquad\qquad$

In particular,

\qquad(xvi)$\qquad\frac{\partial^{4}\Phi_{P}}{\partial x_{\text{c}}%
^{2}\partial x_{i}^{2}}(y_{0})=$ $2[(\frac{\partial X_{i}}{\partial
x_{\text{c}}})^{2}+X_{i}\frac{\partial^{2}X_{i}}{\partial x_{\text{c}}^{2}%
}](y_{0})$ $-\frac{\partial^{3}X_{i}}{\partial x_{\text{c}}^{2}\partial x_{i}%
}(y_{0})\qquad$from $\left(  B_{99}\right)  $

We recall that by (ii) of \textbf{Proposition }$\left(  3.2\right)  $ or by a
direct computation that,

\qquad$<\nabla\Phi,X>$ $=\frac{\partial\Phi}{\partial x_{j}}X_{j}$

\qquad I$_{32141}=\frac{1}{12}\frac{\partial^{2}\Phi^{-1}}{\partial x_{i}^{2}%
}(y_{0})[\frac{\partial\Phi}{\partial x_{j}}X_{j}](y_{0})\phi(y_{0})$

From (i) and (iv) of \textbf{Table B}$_{4},$ we have:

\qquad I$_{32141}=\frac{1}{12}[X_{i}^{2}(y_{0})+\frac{\partial X_{i}}{\partial
x_{i}}(y_{0})][-X_{j}X_{j}](y_{0})\phi(y_{0})$

$=-\frac{1}{12}[X_{i}^{2}X_{j}^{2}](y_{0})\phi(y_{0})-\frac{1}{12}[X_{j}%
^{2}\frac{\partial X_{i}}{\partial x_{i}}](y_{0})\phi(y_{0})$

Next we have for $i=q+1,...,n$ and $j=1,...,q,q+1,...,n:$

\qquad I$_{32142}=\frac{1}{12}\Phi^{-1}(y_{0})\frac{\partial^{2}}{\partial
x_{i}^{2}}[<\nabla\Phi,X>](y_{0})\phi(y_{0})=\frac{1}{12}\frac{\partial^{2}%
}{\partial x_{i}^{2}}[\frac{\partial\Phi}{\partial x_{j}}X_{j}](y_{0}%
)\phi(y_{0})\qquad\ \ \ \ $

\qquad$=\frac{1}{12}[\frac{\partial^{3}\Phi}{\partial x_{i}^{2}\partial x_{j}%
}X_{j}+\frac{\partial\Phi}{\partial x_{j}}\frac{\partial^{2}X_{j}}{\partial
x_{i}^{2}}+2\frac{\partial^{2}\Phi}{\partial x_{i}\partial x_{j}}%
\frac{\partial X_{j}}{\partial x_{i}}](y_{0})\phi(y_{0})$

Then for a = 1,...,q and $i,j=q+1,...,n,$ we have:

I$_{32142}=\frac{1}{12}[\frac{\partial^{3}\Phi}{\partial x_{i}^{2}\partial
x_{\text{a}}}X_{\text{a}}+\frac{\partial\Phi}{\partial x_{\text{a}}}%
\frac{\partial^{2}X_{\text{a}}}{\partial x_{i}^{2}}+2\frac{\partial^{2}\Phi
}{\partial x_{i}\partial x_{\text{a}}}\frac{\partial X_{\text{a}}}{\partial
x_{i}}](y_{0})\phi(y_{0})$

$+\frac{1}{12}[\frac{\partial^{3}\Phi}{\partial x_{i}^{2}\partial x_{j}}%
X_{j}+\frac{\partial\Phi}{\partial x_{j}}\frac{\partial^{2}X_{j}}{\partial
x_{i}^{2}}+2\frac{\partial^{2}\Phi}{\partial x_{i}\partial x_{j}}%
\frac{\partial X_{j}}{\partial x_{i}}](y_{0})\phi(y_{0})$

Since $\frac{\partial\Phi}{\partial x_{\text{a}}}(y_{0})=0$ and $\frac
{\partial\Phi}{\partial x_{j}}(y_{0})=-X_{j}(y_{0}),$ we have by (ix), (xiv), :

\qquad$\frac{\partial^{2}\Phi_{P}}{\partial x_{\text{a}}\partial x_{i}}%
(y_{0})=$ $\frac{\partial}{\partial x_{\text{a}}}(\nabla\log\Phi_{P}%
)_{i}(y_{0})=-\frac{\partial X_{i}}{\partial x_{\text{a}}}(y_{0})$ by (xi) of
\textbf{Table B}$_{4}$

$\qquad\frac{\partial^{3}\Phi_{P}}{\partial x_{\text{a}}\partial x_{i}^{2}%
}(y_{0})=4\left(  X_{i}\frac{\partial X_{i}}{\partial x_{\text{a}}}\right)
(y_{0})-\frac{\partial^{2}X_{i}}{\partial x_{\text{a}}\partial x_{i}}(y_{0})$
by (xiv) of \textbf{Table B}$_{4}$

\qquad$\frac{\partial^{2}\Phi}{\partial x_{i}\partial x_{j}}(y_{0})=\frac
{1}{12}[X_{i}X_{j}-\frac{1}{2}\left(  \frac{\partial X_{i}}{\partial x_{j}%
}+\frac{\partial X_{j}}{\partial x_{i}}\right)  ](y_{0})$ by (ii) of
\textbf{Table B}$_{4}$

$\qquad\frac{\partial^{3}\Phi_{P}}{\partial x_{i}^{2}\partial x_{j}}%
(y_{0})=[-X_{i}^{2}X_{j}+X_{j}\frac{\partial X_{i}}{\partial x_{i}}%
+X_{i}\left(  \frac{\partial X_{i}}{\partial x_{j}}+\frac{\partial X_{j}%
}{\partial x_{i}}\right)  -\frac{1}{3}\left(  \frac{\partial^{2}X_{j}%
}{\partial x_{i}^{2}}+2\frac{\partial^{2}X_{i}}{\partial x_{i}\partial x_{j}%
}\right)  ](y_{0})$ by (v) of \textbf{Table B}$_{4}.\qquad\qquad$

I$_{32142}=\frac{1}{12}[\frac{\partial^{3}\Phi}{\partial x_{i}^{2}\partial
x_{\text{a}}}X_{\text{a}}+\frac{\partial\Phi}{\partial x_{\text{a}}}%
\frac{\partial^{2}X_{\text{a}}}{\partial x_{i}^{2}}+2\frac{\partial^{2}\Phi
}{\partial x_{i}\partial x_{\text{a}}}\frac{\partial X_{\text{a}}}{\partial
x_{i}}](y_{0})\phi(y_{0})$

$\qquad\ \ +\frac{1}{12}[\frac{\partial^{3}\Phi}{\partial x_{i}^{2}\partial
x_{j}}X_{j}+\frac{\partial\Phi}{\partial x_{j}}\frac{\partial^{2}X_{j}%
}{\partial x_{i}^{2}}+2\frac{\partial^{2}\Phi}{\partial x_{i}\partial x_{j}%
}\frac{\partial X_{j}}{\partial x_{i}}](y_{0})\phi(y_{0})$

We use (i), (ii) and (v) of \textbf{Table B}$_{4},$ given by:

$\frac{\partial^{3}\Phi_{P}}{\partial x_{i}^{2}\partial x_{j}}(y_{0}%
)=[-X_{i}^{2}X_{j}+X_{j}\frac{\partial X_{i}}{\partial x_{i}}+X_{i}\left(
\frac{\partial X_{i}}{\partial x_{j}}+\frac{\partial X_{j}}{\partial x_{i}%
}\right)  -\frac{1}{3}\left(  \frac{\partial^{2}X_{j}}{\partial x_{i}^{2}%
}+2\frac{\partial^{2}X_{i}}{\partial x_{i}\partial x_{j}}\right)  ](y_{0})$

I$_{32142}=+\frac{1}{12}[-X_{i}^{2}X_{j}^{2}+X_{j}^{2}\frac{\partial X_{i}%
}{\partial x_{i}}\ +X_{i}X_{j}\left(  \frac{\partial X_{j}}{\partial x_{i}%
}+\frac{\partial X_{i}}{\partial x_{j}}\right)  -\frac{1}{3}X_{j}\left(
\frac{\partial^{2}X_{j}}{\partial x_{i}^{2}}+2\frac{\partial^{2}X_{i}%
}{\partial x_{i}\partial x_{j}}\right)  ](y_{0})\phi(y_{0})$

$\qquad+\frac{1}{12}[-X_{j}\frac{\partial^{2}X_{j}}{\partial x_{i}^{2}}%
](y_{0})\phi(y_{0})+\frac{1}{6}[X_{i}X_{j}-\frac{1}{2}\left(  \frac{\partial
X_{i}}{\partial x_{j}}+\frac{\partial X_{j}}{\partial x_{i}}\right)
]\frac{\partial X_{j}}{\partial x_{i}}(y_{0})\phi(y_{0})$

Next we have:

\qquad I$_{32143}=\frac{1}{6}\frac{\partial\Phi^{-1}}{\partial x_{i}}%
(y_{0})\frac{\partial}{\partial x_{i}}[<\nabla\Phi,X>](y_{0})\phi(y_{0})$

\qquad\qquad$=\frac{1}{6}X_{i}(y_{0})\frac{\partial}{\partial x_{i}}%
[\frac{\partial\Phi}{\partial x_{j}}X_{j}](y_{0})\phi(y_{0})=\frac{1}{6}%
X_{i}(y_{0})[\frac{\partial^{2}\Phi}{\partial x_{i}\partial x_{j}}X_{j}%
+\frac{\partial\Phi}{\partial x_{j}}\frac{\partial X_{j}}{\partial x_{i}%
}](y_{0})\phi(y_{0})$

\qquad\qquad$=\frac{1}{6}X_{i}(y_{0})[X_{i}X_{j}X_{j}-\frac{1}{2}X_{j}\left(
\frac{\partial X_{i}}{\partial x_{j}}+\frac{\partial X_{j}}{\partial x_{i}%
}\right)  +(-X_{j})\frac{\partial X_{j}}{\partial x_{i}}](y_{0})\phi(y_{0})$

\qquad$\ \ \ \ \ =+\frac{1}{6}[X_{i}^{2}X_{j}^{2}-\frac{1}{2}X_{i}X_{j}\left(
\frac{\partial X_{i}}{\partial x_{j}}+\frac{\partial X_{j}}{\partial x_{i}%
}\right)  -X_{i}X_{j}\frac{\partial X_{j}}{\partial x_{i}}](y_{0})\phi(y_{0})$

\qquad$\ \ \ \ =+\frac{1}{6}[X_{i}^{2}X_{j}^{2}-\frac{1}{2}X_{i}X_{j}%
\frac{\partial X_{i}}{\partial x_{j}}-\frac{3}{2}X_{i}X_{j}\frac{\partial
X_{j}}{\partial x_{i}}](y_{0})\phi(y_{0})$

\qquad I$_{32143}=+\frac{1}{12}[2X_{i}^{2}X_{j}^{2}-X_{i}X_{j}\frac{\partial
X_{i}}{\partial x_{j}}-3X_{i}X_{j}\frac{\partial X_{j}}{\partial x_{i}}%
](y_{0})\phi(y_{0})$

We see that,

\qquad I$_{3214}=$ I$_{32141}+$ I$_{32142}+$ I$_{32143}$

\qquad$-\frac{1}{12}[X_{i}^{2}X_{j}^{2}](y_{0})\phi(y_{0})-\frac{1}{12}%
[X_{j}^{2}\frac{\partial X_{i}}{\partial x_{i}}](y_{0})\phi(y_{0})\qquad$
I$_{32141}$

\qquad$+\frac{1}{12}[-X_{i}^{2}X_{j}^{2}+X_{j}^{2}\frac{\partial X_{i}%
}{\partial x_{i}}\ +X_{i}X_{j}\left(  \frac{\partial X_{j}}{\partial x_{i}%
}+\frac{\partial X_{i}}{\partial x_{j}}\right)  -\frac{1}{3}X_{j}\left(
\frac{\partial^{2}X_{j}}{\partial x_{i}^{2}}+2\frac{\partial^{2}X_{i}%
}{\partial x_{i}\partial x_{j}}\right)  ](y_{0})\phi(y_{0})\qquad$I$_{32142}$

$\qquad+\frac{1}{12}[-X_{j}\frac{\partial^{2}X_{j}}{\partial x_{i}^{2}}%
](y_{0})\phi(y_{0})+\frac{1}{6}[X_{i}X_{j}-\frac{1}{2}\left(  \frac{\partial
X_{i}}{\partial x_{j}}+\frac{\partial X_{j}}{\partial x_{i}}\right)
](y_{0})\frac{\partial X_{j}}{\partial x_{i}}(y_{0})\phi(y_{0})$

\qquad$+\frac{1}{12}[2X_{i}^{2}X_{j}^{2}-X_{i}X_{j}\frac{\partial X_{i}%
}{\partial x_{j}}-3X_{i}X_{j}\frac{\partial X_{j}}{\partial x_{i}}](y_{0}%
)\phi(y_{0})\qquad$I$_{32143}$

We simplify and have:

$\qquad\qquad\qquad$I$_{3214}=+\frac{1}{12}[-\frac{1}{3}X_{j}\left(
2\frac{\partial^{2}X_{i}}{\partial x_{i}\partial x_{j}}\right)  ](y_{0}%
)\phi(y_{0})\qquad$I$_{32142}$

$\qquad\qquad-\frac{4}{36}[X_{j}\frac{\partial^{2}X_{j}}{\partial x_{i}^{2}%
}](y_{0})\phi(y_{0})+\frac{1}{6}[-\frac{1}{2}\left(  \frac{\partial X_{i}%
}{\partial x_{j}}+\frac{\partial X_{j}}{\partial x_{i}}\right)  ](y_{0}%
)\frac{\partial X_{j}}{\partial x_{i}}(y_{0})\phi(y_{0})$

$\left(  B_{116}\right)  \qquad\ $I$_{3214}=-\frac{1}{18}[X_{j}\left(
2\frac{\partial^{2}X_{j}}{\partial x_{i}^{2}}+\frac{\partial^{2}X_{i}%
}{\partial x_{i}\partial x_{j}}\right)  ](y_{0})\phi(y_{0})$

$\qquad\qquad\qquad\qquad-\frac{1}{12}[\left(  \frac{\partial X_{i}}{\partial
x_{j}}+\frac{\partial X_{j}}{\partial x_{i}}\right)  ](y_{0})\frac{\partial
X_{j}}{\partial x_{i}}(y_{0})\phi(y_{0})\qquad$J$_{32142}$

$\qquad\qquad\qquad\qquad$\qquad\qquad\qquad\qquad\qquad\qquad\qquad
\qquad\qquad\qquad\qquad\qquad\qquad\qquad\qquad\qquad\qquad$\blacksquare$

$\left(  B_{117}\right)  \qquad$I$_{3215}=\frac{1}{12}\frac{\partial
^{2}\text{V}}{\partial x_{i}^{2}}(y_{0})\phi(y_{0})$

\qquad\qquad\qquad\qquad\qquad\qquad\qquad\qquad$\blacksquare$

\subsubsection{\textbf{Computation of I}$_{321}$\qquad\ }

Recall that that: I$_{321}$ $\mathbf{=}\frac{1}{12}$ $\frac{\partial^{2}%
}{\partial x_{i}^{2}}(\Psi^{-1}L\Psi)(y_{0})\phi(y_{0})=$ I$_{3211}$ $+$
I$_{3212}$ $+$ I$_{3213}$ $+$ I$_{3214}$ $+$ I$_{3215}$\qquad

Here, I$_{3211}=$ A$_{321}=\frac{1}{24}\frac{\partial^{2}}{\partial
\text{x}_{i}^{2}}(\theta^{\frac{1}{2}}\Delta\theta^{-\frac{1}{2}})$ is given
in $\left(  A_{31}\right)  $ above or by (xii) of \textbf{Table A}$_{10}.$ Then,

I$_{3212}=\frac{1}{24}$ $\frac{\partial^{2}}{\partial x_{i}^{2}}(\Phi
^{-1}\Delta\Phi)(y_{0})\phi(y_{0})$ is given in $\left(  B_{114}\right)  $ and
I$_{3213},$ I$_{3214}$ and I$_{3215}$ are given in $\left(  B_{115}\right)
,\left(  B_{116}\right)  $ and $\left(  B_{117}\right)  $ respectively.

\qquad\qquad\qquad\qquad\qquad\qquad\qquad\qquad\qquad\qquad\qquad\qquad
\qquad\qquad\qquad\qquad\qquad\qquad\qquad\qquad\qquad$\blacksquare$

$\left(  B_{118}\right)  \qquad$I$_{321}$ $\mathbf{=}\frac{1}{12}$
$\frac{\partial^{2}}{\partial x_{i}^{2}}(\Psi^{-1}L\Psi)(y_{0})\phi
(y_{0})\qquad$\ 

$\qquad=-\frac{1}{3456}[3<H,i>^{2}\ +2(\tau^{M}-3\tau^{P}%
\ +\overset{q}{\underset{\text{a=1}}{\sum}}\varrho_{\text{aa}}^{M}%
+\overset{q}{\underset{\text{a,b}=1}{\sum}}R_{\text{abab}}^{M})]^{2}%
(y_{0})\phi(y_{0})\phi(y_{0})$\ \qquad I$_{3211}$

$+\frac{1}{24}[2<H,i>^{2}(y_{0})+\frac{1}{3}(\tau^{M}-3\tau^{P}%
+\overset{q}{\underset{\text{a}=1}{\sum}}\varrho_{\text{aa}}%
+\overset{q}{\underset{\text{a,b}=1}{\sum}}R_{\text{abab}})](y_{0})\phi
(y_{0})\qquad$I$_{3212}=\frac{1}{24}(L_{1}+L_{2}+L_{3})$

$\times\lbrack\frac{1}{4}<H,j>^{2}(y_{0})+\frac{1}{6}(\tau^{M}-3\tau
^{P}+\overset{q}{\underset{\text{a}=1}{\sum}}\varrho_{\text{aa}}%
^{M}+\overset{q}{\underset{\text{a,b}=1}{\sum}}R_{\text{abab}}^{M}%
)](y_{0})\phi(y_{0})$

$-\frac{1}{4}\times\frac{1}{24}[2\varrho_{ij}+$
$\overset{q}{\underset{\text{a}=1}{4\sum}}R_{i\text{a}j\text{a}}%
-3\overset{q}{\underset{\text{a,b=1}}{\sum}}(T_{\text{aa}i}T_{\text{bb}%
j}-T_{\text{ab}i}T_{\text{ab}j})-3\overset{q}{\underset{\text{a,b=1}}{\sum}%
}(T_{\text{aa}j}T_{\text{bb}i}-T_{\text{ab}j}T_{\text{ab}i}](y_{0})\qquad
L_{2}\qquad L_{21}$

$\ \times\lbrack<H,i><H,j>](y_{0})\phi(y_{0})$

$-\frac{1}{24}\times\frac{1}{36}[2\varrho_{ij}+$
$\overset{q}{\underset{\text{a}=1}{4\sum}}R_{i\text{a}j\text{a}}%
-3\overset{q}{\underset{\text{a,b=1}}{\sum}}(T_{\text{aa}i}T_{\text{bb}%
j}-T_{\text{ab}i}T_{\text{ab}j})-3\overset{q}{\underset{\text{a,b=1}}{\sum}%
}(T_{\text{aa}j}T_{\text{bb}i}-T_{\text{ab}j}T_{\text{ab}i}]^{2}(y_{0}%
)\phi(y_{0})$

$-\frac{1}{24}\times\frac{1}{12}[<H,j>](y_{0})\times\lbrack\{\nabla_{i}%
\varrho_{ij}-2\varrho_{ij}<H,i>+\overset{q}{\underset{\text{a}=1}{\sum}%
}(\nabla_{i}R_{\text{a}i\text{a}j}-4R_{i\text{a}j\text{a}}<H,i>)\qquad
L_{212}$

$+4\overset{q}{\underset{\text{a,b=1}}{\sum}}R_{i\text{a}j\text{b}%
}T_{\text{ab}i}+2\overset{q}{\underset{\text{a,b,c=1}}{\sum}}(T_{\text{aa}%
i}T_{\text{bb}j}T_{\text{cc}i}-3T_{\text{aa}i}T_{\text{bc}j}T_{\text{bc}%
i}+2T_{\text{ab}i}T_{\text{bc}j}T_{\text{ac}i})](y_{0})\phi(y_{0})$%
\qquad\qquad\qquad\qquad\qquad\ \ 

$-\frac{1}{24}\times\frac{1}{12}[<H,j>](y_{0})\times\lbrack\nabla_{j}%
\varrho_{ii}-2\varrho_{ij}<H,i>+\overset{q}{\underset{\text{a}=1}{\sum}%
}(\nabla_{j}R_{\text{a}i\text{a}i}-4R_{i\text{a}j\text{a}}<H,i>)$

$+4\overset{q}{\underset{\text{a,b=1}}{\sum}}R_{j\text{a}i\text{b}%
}T_{\text{ab}i}+2\overset{q}{\underset{\text{a,b,c=1}}{\sum}}(T_{\text{aa}%
j}T_{\text{bb}i}T_{\text{cc}i}-3T_{\text{aa}j}T_{\text{bc}i}T_{\text{bc}%
i}+2T_{\text{ab}j}T_{\text{bc}i}T_{\text{ac}i})](y_{0})\phi(y_{0})$

$-\frac{1}{24}\times\frac{1}{12}[<H,j>](y_{0})\times\lbrack\nabla_{i}%
\varrho_{ij}-2\varrho_{ii}<H,j>+\overset{q}{\underset{\text{a}=1}{\sum}%
}(\nabla_{i}R_{\text{a}i\text{a}j}-4R_{i\text{a}i\text{a}}<H,j>)$

$+4\overset{q}{\underset{\text{a,b=1}}{\sum}}R_{i\text{a}i\text{b}%
}T_{\text{ab}j}+2\overset{q}{\underset{\text{a,b,c}=1}{\sum}}(T_{\text{aa}%
i}T_{\text{bb}i}T_{\text{cc}j}-3T_{\text{aa}i}T_{\text{bc}i}T_{\text{bc}%
j}+2T_{\text{ab}i}T_{\text{bc}i}T_{\text{ac}j})](y_{0})\phi(y_{0})$

$-\frac{1}{3}[<H,j><H,k>](y_{0})R_{ijik}(y_{0})-$ $\frac{1}{24}\times\frac
{15}{8}[<H,i>^{2}<H,j>^{2}](y_{0})\phi(y_{0})\qquad L_{213}$

$-\frac{1}{24}\times\frac{1}{4}<H,i><H,j>[2\varrho_{ij}%
+\overset{q}{\underset{\text{a}=1}{4\sum}}R_{i\text{a}j\text{a}}%
-3\overset{q}{\underset{\text{a,b=1}}{\sum}}(T_{\text{aa}i}T_{\text{bb}%
j}-T_{\text{ab}i}T_{\text{ab}j})$

$-3\overset{q}{\underset{\text{a,b=1}}{\sum}}(T_{\text{aa}j}T_{\text{bb}%
i}-T_{\text{ab}j}T_{\text{ab}i}](y_{0})\phi(y_{0})$

$-\frac{1}{24}\times\frac{1}{4}<H,j>^{2}[\tau^{M}\ -3\tau^{P}%
+\ \underset{\text{a}=1}{\overset{\text{q}}{\sum}}\varrho_{\text{aa}}^{M}+$
$\overset{q}{\underset{\text{a},\text{b}=1}{\sum}}R_{\text{abab}}^{M}$
$](y_{0})\phi(y_{0})$

$+\frac{1}{24}\times\frac{1}{12}<H,j>[\nabla_{i}\varrho_{ij}-2\varrho
_{ij}<H,i>+\overset{q}{\underset{\text{a}=1}{\sum}}(\nabla_{i}R_{\text{a}%
i\text{a}j}-4R_{i\text{a}j\text{a}}<H,i>)+4\overset{q}{\underset{\text{a,b=1}%
}{\sum}}R_{i\text{a}j\text{b}}T_{\text{ab}i}$

$\qquad+2\overset{q}{\underset{\text{a,b,c=1}}{\sum}}(T_{\text{aa}%
i}T_{\text{bb}j}T_{\text{cc}i}-3T_{\text{aa}i}T_{\text{bc}j}T_{\text{bc}%
i}+2T_{\text{ab}i}T_{\text{bc}j}T_{\text{ac}i})](y_{0})\phi(y_{0})$%
\qquad\qquad\qquad\qquad\qquad\ \ 

$+\frac{1}{12}<H,j>[\nabla_{j}\varrho_{ii}-2\varrho_{ij}%
<H,i>+\overset{q}{\underset{\text{a}=1}{\sum}}(\nabla_{j}R_{\text{a}%
i\text{a}i}-4R_{i\text{a}j\text{a}}<H,i>)+4\overset{q}{\underset{\text{a,b=1}%
}{\sum}}R_{j\text{a}i\text{b}}T_{\text{ab}i}$

$+2\overset{q}{\underset{\text{a,b,c=1}}{\sum}}(T_{\text{aa}j}T_{\text{bb}%
i}T_{\text{cc}i}-3T_{\text{aa}j}T_{\text{bc}i}T_{\text{bc}i}+2T_{\text{ab}%
j}T_{\text{bc}i}T_{\text{ac}i})](y_{0})\phi(y_{0})$

$+\frac{1}{24}\times\frac{1}{12}<H,j>[\nabla_{i}\varrho_{ij}-2\varrho
_{ii}<H,j>+\overset{q}{\underset{\text{a}=1}{\sum}}(\nabla_{i}R_{\text{a}%
i\text{a}j}-4R_{i\text{a}i\text{a}}<H,j>)+4\overset{q}{\underset{\text{a,b=1}%
}{\sum}}R_{i\text{a}i\text{b}}T_{\text{ab}j}$

$+2\overset{q}{\underset{\text{a,b,c}=1}{\sum}}(T_{\text{aa}i}T_{\text{bb}%
i}T_{\text{cc}j}-3T_{\text{aa}i}T_{\text{bc}i}T_{\text{bc}j}+2T_{\text{ab}%
i}T_{\text{bc}i}T_{\text{ac}j})](y_{0})\phi(y_{0})$

$-\frac{1}{24}\times\frac{1}{6}R_{jijk}(y_{0})$ \ $[<H,i><H,k>](y_{0}%
)\phi(y_{0})\qquad\qquad L_{22}$

$-\frac{1}{24}\times\frac{1}{18}R_{jijk}(y_{0})[2\varrho_{ik}+$
$\overset{q}{\underset{\text{a}=1}{4\sum}}R_{i\text{a}k\text{a}}%
-3\overset{q}{\underset{\text{a,b=1}}{\sum}}(T_{\text{aa}i}T_{\text{bb}%
k}-T_{\text{ab}i}T_{\text{ab}k})$

$-3\overset{q}{\underset{\text{a,b=1}}{\sum}}(T_{\text{aa}k}T_{\text{bb}%
i}-T_{\text{ab}k}T_{\text{ab}i}](y_{0})\phi(y_{0})$

$+\frac{1}{24}\times\frac{1}{6}<H,k>(y_{0})[\nabla_{j}$R$_{ijik}(y_{0}%
)-\nabla_{i}$R$_{jijk}](y_{0})\phi(y_{0})$

$-\frac{1}{24}\times\frac{15}{4}<H,i>^{2}(y_{0})<H,j>^{2}(y_{0})\phi
(y_{0})\qquad\qquad L_{23}\qquad L_{231}$

$-\frac{1}{24}\times\frac{1}{2}<H,i>(y_{0})<H,j>(y_{0})$

$\times\lbrack2\varrho_{ij}+$ $\overset{q}{\underset{\text{a}=1}{4\sum}%
}R_{i\text{a}j\text{a}}-3\overset{q}{\underset{\text{a,b=1}}{\sum}%
}(T_{\text{aa}i}T_{\text{bb}j}-T_{\text{ab}i}T_{\text{ab}j}%
)-3\overset{q}{\underset{\text{a,b=1}}{\sum}}(T_{\text{aa}j}T_{\text{bb}%
i}-T_{\text{ab}j}T_{\text{ab}i}](y_{0})\phi(y_{0})$

$-\frac{1}{24}\times\frac{1}{2}<H,i>^{2}(y_{0})[\varrho_{jj}+$
$\overset{q}{\underset{\text{a}=1}{2\sum}}R_{j\text{a}j\text{a}}%
-3\overset{q}{\underset{\text{a,b=1}}{\sum}}(T_{\text{aa}j}T_{\text{bb}%
j}-T_{\text{ab}j}T_{\text{ab}j})](y_{0})\phi(y_{0})$

$-\frac{1}{24}\times\frac{1}{6}<H,i>(y_{0})[\nabla_{i}\varrho_{jj}%
-2\varrho_{ij}<H,j>+\overset{q}{\underset{\text{a}=1}{\sum}}(\nabla
_{i}R_{\text{a}j\text{a}j}-4R_{i\text{a}j\text{a}}<H,j>)$

$+4\overset{q}{\underset{\text{a,b=1}}{\sum}}R_{i\text{a}j\text{b}%
}T_{\text{ab}j}+2\overset{q}{\underset{\text{a,b,c=1}}{\sum}}(T_{\text{aa}%
i}T_{\text{bb}j}T_{\text{cc}j}-3T_{\text{aa}i}T_{\text{bc}j}T_{\text{bc}%
j}+2T_{\text{ab}i}T_{\text{bc}j}T_{\text{ca}j})](y_{0})\phi(y_{0})$%
\qquad\qquad\qquad\qquad\qquad\ \ 

$-\frac{1}{24}\times\frac{1}{6}<H,i>(y_{0})[\nabla_{j}\varrho_{ij}%
-2\varrho_{ij}<H,j>+\overset{q}{\underset{\text{a}=1}{\sum}}(\nabla
_{j}R_{\text{a}i\text{a}j}-4R_{j\text{a}i\text{a}}<H,j>)$

$+4\overset{q}{\underset{\text{a,b=1}}{\sum}}R_{j\text{a}i\text{b}%
}T_{\text{ab}j}+2\overset{q}{\underset{\text{a,b,c=1}}{\sum}}(T_{\text{aa}%
j}T_{\text{bb}i}T_{\text{cc}j}-3T_{\text{aa}j}T_{\text{bc}i}T_{\text{bc}%
j}+2T_{\text{ab}j}T_{\text{bc}i}T_{\text{ac}j})](y_{0})$

$-\frac{1}{24}\times\frac{1}{6}<H,i>(y_{0})[\nabla_{j}\varrho_{ij}%
-2\varrho_{jj}<H,i>+\overset{q}{\underset{\text{a}=1}{\sum}}(\nabla
_{j}R_{\text{a}i\text{a}j}-4R_{j\text{a}j\text{a}}<H,i>)$

$+4\overset{q}{\underset{\text{a,b=1}}{\sum}}R_{j\text{a}j\text{b}%
}T_{\text{ab}i}+2\overset{q}{\underset{\text{a,b,c=1}}{\sum}}(T_{\text{aa}%
j}T_{\text{bb}j}T_{\text{cc}i}-3T_{\text{aa}j}T_{\text{bc}j}T_{\text{bc}%
i}+2T_{\text{ab}j}T_{\text{bc}j}T_{\text{ac}i})](y_{0})\phi(y_{0})$

$-\frac{1}{24}\times\frac{1}{4}<H,j>^{2}(y_{0})[\varrho_{ii}+$
$\overset{q}{\underset{\text{a}=1}{2\sum}}R_{i\text{a}i\text{a}}%
-3\overset{q}{\underset{\text{a,b=1}}{\sum}}(T_{\text{aa}i}T_{\text{bb}%
i}-T_{\text{ab}i}T_{\text{ab}i})](y_{0})\qquad L_{232}$

$-\frac{1}{24}\times\frac{1}{18}[\varrho_{ii}+$
$\overset{q}{\underset{\text{a}=1}{2\sum}}R_{i\text{a}i\text{a}}%
-3\overset{q}{\underset{\text{a,b=1}}{\sum}}(T_{\text{aa}i}T_{\text{bb}%
i}-T_{\text{ab}i}T_{\text{ab}i})](y_{0})\phi(y_{0})$

$\times\lbrack\varrho_{jj}+$ $\overset{q}{\underset{\text{a}=1}{2\sum}%
}R_{j\text{a}j\text{a}}-3\overset{q}{\underset{\text{a,b=1}}{\sum}%
}(T_{\text{aa}j}T_{\text{bb}j}-T_{\text{ab}j}T_{\text{ab}j})]\}(y_{0}%
)\phi(y_{0})$

$+\frac{1}{24}\times\frac{1}{2}R_{ijik}(y_{0})$ \ $[<H,j><H,k>](y_{0}%
)\phi(y_{0})\qquad\qquad\qquad\qquad$\ $L_{233}$

$+\frac{1}{24}\times\frac{1}{18}R_{ijik}(y_{0})[2\varrho_{jk}+$
$\overset{q}{\underset{\text{a}=1}{4\sum}}R_{j\text{a}k\text{a}}%
-3\overset{q}{\underset{\text{a,b=1}}{\sum}}(T_{\text{aa}j}T_{\text{bb}%
k}-T_{\text{ab}j}T_{\text{ab}k})$

$-3\overset{q}{\underset{\text{a,b=1}}{\sum}}(T_{\text{aa}k}T_{\text{bb}%
j}-T_{\text{ab}k}T_{\text{ab}j}](y_{0})\phi(y_{0})$

$+\overset{n}{\underset{i,j=q+1}{\sum}}\frac{35}{128}<H,i>^{2}(y_{0}%
)<H,j>^{2}(y_{0})\qquad\qquad\ \frac{1}{24}\frac{\partial^{4}\theta^{-\frac
{1}{2}}}{\partial x_{i}^{2}\partial x_{j}^{2}}(y_{0})$

$+\frac{5}{192}\overset{n}{\underset{j=q+1}{\sum}}<H,j>^{2}(y_{0})[\tau
^{M}\ -3\tau^{P}+\ \underset{\text{a}=1}{\overset{\text{q}}{\sum}}%
\varrho_{\text{aa}}^{M}+\overset{q}{\underset{\text{a},\text{b}=1}{\sum}%
}R_{\text{abab}}^{M}](y_{0})\qquad\ \ \ \ \ \ \ \ $

$+\frac{5}{192}\overset{n}{\underset{i=q+1}{\sum}}<H,i>^{2}(y_{0})[\tau
^{M}\ -3\tau^{P}+\ \underset{\text{a}=1}{\overset{\text{q}}{\sum}}%
\varrho_{\text{aa}}^{M}+\overset{q}{\underset{\text{a},\text{b}=1}{\sum}%
}R_{\text{abab}}^{M}](y_{0})\qquad\qquad$

$+\frac{5}{192}\overset{n}{\underset{i,j=q+1}{\sum}}[<H,i><H,j>](y_{0}%
)\qquad\qquad\qquad\qquad\qquad\qquad\qquad\qquad$

$\times\lbrack2\varrho_{ij}+4\overset{q}{\underset{\text{a}=1}{\sum}%
}R_{i\text{a}j\text{a}}-3\overset{q}{\underset{\text{a,b=1}}{\sum}%
}(T_{\text{aa}i}T_{\text{bb}j}-T_{\text{ab}i}T_{\text{ab}j}%
)-3\overset{q}{\underset{\text{a,b=1}}{\sum}}(T_{\text{aa}j}T_{\text{bb}%
i}-T_{\text{ab}j}T_{\text{ab}i})](y_{0})$

$+\frac{1}{96}\overset{n}{\underset{i,j=q+1}{\sum}}<H,j>(y_{0})[\{\nabla
_{i}\varrho_{ij}-2\varrho_{ij}<H,i>+\overset{q}{\underset{\text{a}=1}{\sum}%
}(\nabla_{i}R_{\text{a}i\text{a}j}-4R_{i\text{a}j\text{a}}<H,i>)\qquad$

$+4\overset{q}{\underset{\text{a,b=1}}{\sum}}R_{i\text{a}j\text{b}%
}T_{\text{ab}i}+2\overset{q}{\underset{\text{a,b,c=1}}{\sum}}(T_{\text{aa}%
i}T_{\text{bb}j}T_{\text{cc}i}-T_{\text{aa}i}T_{\text{bc}j}T_{\text{bc}%
i}-2T_{\text{bc}j}(T_{\text{aa}i}T_{\text{bc}i}-T_{\text{ab}i}T_{\text{ac}%
i}))\}$\qquad\qquad\qquad\ \ 

$+\{\nabla_{j}\varrho_{ii}-2\varrho_{ij}<H,i>+\overset{q}{\underset{\text{a}%
=1}{\sum}}(\nabla_{j}R_{\text{a}i\text{a}i}-4R_{i\text{a}j\text{a}}<H,i>)$

$+4\overset{q}{\underset{\text{a,b=1}}{\sum}}R_{j\text{a}i\text{b}%
}T_{\text{ab}i}+2\overset{q}{\underset{\text{a,b,c=1}}{\sum}}(T_{\text{aa}%
j}(T_{\text{bb}i}T_{\text{cc}i}-T_{\text{bc}i}T_{\text{bc}i})-2T_{\text{aa}%
j}T_{\text{bc}i}T_{\text{bc}i}+2T_{\text{ab}j}T_{\text{bc}i}T_{\text{ac}%
i})\}\qquad$

$+\{\nabla_{i}\varrho_{ij}-2\varrho_{ii}<H,j>+\overset{q}{\underset{\text{a}%
=1}{\sum}}(\nabla_{i}R_{\text{a}i\text{a}j}-4R_{i\text{a}i\text{a}%
}<H,j>)+4\overset{q}{\underset{\text{a,b=1}}{\sum}}R_{i\text{a}i\text{b}%
}T_{\text{ab}j}$

$+2\overset{q}{\underset{\text{a,b,c}=1}{\sum}}(T_{\text{aa}i}T_{\text{bb}%
i}T_{\text{cc}j}-3T_{\text{aa}i}T_{\text{bc}i}T_{\text{bc}j}+2T_{\text{ab}%
i}T_{\text{bc}i}T_{\text{ac}j})\}](y_{0})$

$+\frac{1}{96}\overset{n}{\underset{i,j=q+1}{\sum}}<H,i>(y_{0})[\{\nabla
_{i}\varrho_{jj}-2\varrho_{ij}<H,j>+\overset{q}{\underset{\text{a}=1}{\sum}%
}(\nabla_{i}R_{\text{a}j\text{a}j}-4R_{i\text{a}j\text{a}}<H,j>)\qquad$

$+4\overset{q}{\underset{\text{a,b=1}}{\sum}}R_{i\text{a}j\text{b}%
}T_{\text{ab}j}+2\overset{q}{\underset{\text{a,b,c=1}}{\sum}}T_{\text{aa}%
i}(T_{\text{bb}j}T_{\text{cc}j}-T_{\text{bc}j}T_{\text{bc}j})-2T_{\text{aa}%
i}T_{\text{bc}j}T_{\text{bc}j}+2T_{\text{ab}i}T_{\text{bc}j}T_{\text{ac}%
j})\}(y_{0})\qquad$\qquad\qquad\qquad\qquad\qquad\ \ 

$+\{\nabla_{j}\varrho_{ij}-2\varrho_{ij}<H,j>+\overset{q}{\underset{\text{a}%
=1}{\sum}}(\nabla_{j}R_{\text{a}i\text{a}j}-4R_{j\text{a}i\text{a}}<H,j>)$

$+4\overset{q}{\underset{\text{a,b=1}}{\sum}}R_{j\text{a}i\text{b}%
}T_{\text{ab}j}+2\overset{q}{\underset{\text{a,b,c=1}}{\sum}}(T_{\text{aa}%
j}T_{\text{bb}i}T_{\text{cc}j}-T_{\text{ab}j}T_{\text{bc}i}T_{\text{ac}%
j}-2T_{\text{bc}i}(T_{\text{aa}j}T_{\text{bc}j}-T_{\text{ab}j}T_{\text{ac}%
j}))\}(y_{0})$

$+\{\nabla_{j}\varrho_{ij}-2\varrho_{jj}<H,i>+\overset{q}{\underset{\text{a}%
=1}{\sum}}(\nabla_{j}R_{\text{a}i\text{a}j}-4R_{j\text{a}j\text{a}%
}<H,i>)+4\overset{q}{\underset{\text{a,b=1}}{\sum}}R_{j\text{a}j\text{b}%
}T_{\text{ab}i}$

$+2\overset{q}{\underset{\text{a,b,c=1}}{\sum}}(T_{\text{aa}j}T_{\text{bb}%
j}T_{\text{cc}i}-3T_{\text{aa}j}T_{\text{bc}j}T_{\text{bc}i}+2T_{\text{ab}%
j}T_{\text{bc}j}T_{\text{ac}i})\}](y_{0})$

$+\frac{1}{576}\overset{n}{\underset{i,j=q+1}{\sum}}[2\varrho_{ij}%
+4\overset{q}{\underset{\text{a}=1}{\sum}}R_{i\text{a}j\text{a}}%
-3\overset{q}{\underset{\text{a,b=1}}{\sum}}(T_{\text{aa}i}T_{\text{bb}%
j}-T_{\text{ab}i}T_{\text{ab}j})-3\overset{q}{\underset{\text{a,b=1}}{\sum}%
}(T_{\text{aa}j}T_{\text{bb}i}-T_{\text{ab}j}T_{\text{ab}i})]^{2}(y_{0})$

$+\frac{1}{288}[\tau^{M}\ -3\tau^{P}+\ \underset{\text{a}=1}{\overset{\text{q}%
}{\sum}}\varrho_{\text{aa}}^{M}+\overset{q}{\underset{\text{a},\text{b}%
=1}{\sum}}R_{\text{abab}}^{M}]^{2}(y_{0})$

$-\ \frac{1}{288}\overset{n}{\underset{i,j=q+1}{\sum}}[$
$\overset{q}{\underset{\text{a=1}}{\sum}}\{-(\nabla_{ii}^{2}R_{j\text{a}%
j\text{a}}+\nabla_{jj}^{2}R_{i\text{a}i\text{a}}+4\nabla_{ij}^{2}%
R_{i\text{a}j\text{a}}+2R_{ij}R_{i\text{a}j\text{a}})\qquad A$

$\qquad+\overset{n}{\underset{p=q+1}{\sum}}\overset{q}{\underset{\text{a=1}%
}{\sum}}(R_{\text{a}iip}R_{\text{a}jjp}+R_{\text{a}jjp}R_{\text{a}%
iip}+R_{\text{a}ijp}R_{\text{a}ijp}+R_{\text{a}ijp}R_{\text{a}jip}%
+R_{\text{a}jip}R_{\text{a}ijp}+R_{\text{a}jip}R_{\text{a}jip})$

$+2\overset{q}{\underset{\text{a,b=1}}{\sum}}\nabla_{i}(R)_{\text{a}%
i\text{b}j}T_{\text{ab}j}+2\overset{q}{\underset{\text{a,b=1}}{\sum}}%
\nabla_{j}(R)_{\text{a}j\text{b}i}T_{\text{ab}i}%
+2\overset{q}{\underset{\text{a,b=1}}{\sum}}\nabla_{i}(R)_{\text{a}j\text{b}%
i}T_{\text{ab}j}+2\overset{q}{\underset{\text{a,b=1}}{\sum}}\nabla
_{i}(R)_{\text{a}j\text{b}j}T_{\text{ab}i}$

$+2\overset{q}{\underset{\text{a,b=1}}{\sum}}\nabla_{j}(R)_{\text{a}%
i\text{b}i}T_{\text{ab}j}+2\overset{q}{\underset{\text{a,b=1}}{\sum}}%
\nabla_{j}(R)_{\text{a}i\text{b}j}T_{\text{ab}i}$

$+\overset{n}{\underset{p=q+1}{\sum}}(-\frac{3}{5}\nabla_{ii}^{2}%
(R)_{jpjp}+\overset{n}{\underset{p=q+1}{\sum}}(-\frac{3}{5}\nabla_{jj}%
^{2}(R)_{ipip}$

$+\overset{n}{\underset{p=q+1}{\sum}}(-\frac{3}{5}\nabla_{ij}^{2}%
(R)_{ipjp}+\overset{n}{\underset{p=q+1}{\sum}}(-\frac{3}{5}\nabla_{ij}%
^{2}(R)_{jpip}+\overset{n}{\underset{p=q+1}{\sum}}(-\frac{3}{5}\nabla_{ji}%
^{2}(R)_{ipjp}+\overset{n}{\underset{p=q+1}{\sum}}(-\frac{3}{5}\nabla_{ji}%
^{2}(R)_{jpip}$

$+\frac{1}{5}\overset{n}{\underset{m,p=q+1}{%
{\textstyle\sum}
}}R_{ipim}R_{jpjm}+\frac{1}{5}\overset{n}{\underset{m,p=q+1}{%
{\textstyle\sum}
}}R_{jpjm}R_{ipim}+\frac{1}{5}\overset{n}{\underset{m,p=q+1}{%
{\textstyle\sum}
}}R_{ipjm}R_{ipjm}+\frac{1}{5}\overset{n}{\underset{m,p=q+1}{%
{\textstyle\sum}
}}R_{ipjm}R_{jpim}$

$+\frac{1}{5}\overset{n}{\underset{m,p=q+1}{%
{\textstyle\sum}
}}R_{jpim}R_{ipjm}+\frac{1}{5}\overset{n}{\underset{m,p=q+1}{%
{\textstyle\sum}
}}R_{jpim}R_{jpim}\}(y_{0})$

$+4\overset{q}{\underset{\text{a,b=1}}{\sum}}\{(\nabla_{i}(R)_{i\text{a}%
j\text{a}}-\overset{q}{\underset{\text{c=1}}{%
{\textstyle\sum}
}}R_{\text{a}i\text{c}i}T_{\text{ac}j})$ $T_{\text{bb}j}+4(\nabla
_{j}(R)_{j\text{a}i\text{a}}-\overset{q}{\underset{\text{c=1}}{%
{\textstyle\sum}
}}R_{\text{a}j\text{c}j}T_{\text{ac}i})$ $T_{\text{bb}i}$

$+4(\nabla_{i}(R)_{j\text{a}i\text{a}}-\overset{q}{\underset{\text{c=1}}{%
{\textstyle\sum}
}}R_{\text{a}i\text{c}j}T_{\text{ac}i})$ $T_{\text{bb}j}$ $4B\ $

$+4(\nabla_{i}(R)_{j\text{a}j\text{a}}-\overset{q}{\underset{\text{c=1}}{%
{\textstyle\sum}
}}R_{\text{a}i\text{c}j}T_{\text{ac}j})$ $T_{\text{bb}i}+4(\nabla
_{j}(R)_{i\text{a}i\text{a}}-\overset{q}{\underset{\text{c=1}}{%
{\textstyle\sum}
}}R_{\text{a}j\text{c}i}T_{\text{ac}i})$ $T_{\text{bb}j}$

$+4(\nabla_{j}(R)_{i\text{a}j\text{a}}-\overset{q}{\underset{\text{c=1}}{%
{\textstyle\sum}
}}R_{\text{a}j\text{c}i}T_{\text{ac}j})$ $T_{\text{bb}i}$

$-4\overset{q}{\underset{\text{a,b=1}}{\sum}}(\nabla_{i}(R)_{i\text{a}%
j\text{b}}-\overset{q}{\underset{\text{c=1}}{%
{\textstyle\sum}
}}R_{\text{b}r\text{c}s}T_{\text{ac}t})T_{\text{ab}j}%
-4\overset{q}{\underset{\text{a,b=1}}{\sum}}(\nabla_{j}(R)_{j\text{a}%
i\text{b}}-\overset{q}{\underset{\text{c=1}}{%
{\textstyle\sum}
}}R_{\text{b}j\text{c}j}T_{\text{ac}i})T_{\text{ab}i}$

$-4\overset{q}{\underset{\text{a,b=1}}{\sum}}(\nabla_{i}(R)_{j\text{a}%
i\text{b}}-\overset{q}{\underset{\text{c=1}}{%
{\textstyle\sum}
}}R_{\text{b}i\text{c}j}T_{\text{ac}i})T_{\text{ab}j}%
-4\overset{q}{\underset{\text{a,b=1}}{\sum}}(\nabla_{i}(R)_{j\text{a}%
j\text{b}}-\overset{q}{\underset{\text{c=1}}{%
{\textstyle\sum}
}}R_{\text{b}i\text{c}j}T_{\text{ac}j})T_{\text{ab}i}$

$-4\overset{q}{\underset{\text{a,b=1}}{\sum}}(\nabla_{j}(R)_{i\text{a}%
i\text{b}}-\overset{q}{\underset{\text{c=1}}{%
{\textstyle\sum}
}}R_{\text{b}j\text{c}i}T_{\text{ac}i})T_{\text{ab}j}%
-4\overset{q}{\underset{\text{a,b=1}}{\sum}}(\nabla_{j}(R)_{i\text{a}%
j\text{b}}-\overset{q}{\underset{\text{c=1}}{%
{\textstyle\sum}
}}R_{\text{b}j\text{c}i}T_{\text{ac}j})T_{\text{ab}i}\}](y_{0})$

$-\frac{1}{48}$ $[\frac{4}{9}\overset{q}{\underset{\text{a,b=1}}{\sum}%
}(\varrho_{\text{aa}}-\overset{q}{\underset{\text{c}=1}{\sum}}R_{\text{acac}%
})(\varrho_{\text{bb}}-\overset{q}{\underset{\text{d}=1}{\sum}}R_{\text{bdbd}%
})+\frac{8}{9}\overset{n}{\underset{i,j=q+1}{\sum}}%
\overset{q}{\underset{\text{a,b}=1}{\sum}}(R_{i\text{a}j\text{a}}%
R_{i\text{b}j\text{b}})\qquad3C$

$+\frac{2}{9}\overset{q}{\underset{\text{a}=1}{\sum}}(\varrho_{\text{aa}}%
^{M}-\varrho_{\text{aa}}^{P})(\tau^{M}-\overset{q}{\underset{\text{c}=1}{\sum
}}\varrho_{\text{cc}}^{M})+\frac{4}{9}\overset{n}{\underset{i,j=q+1}{\sum}%
}\overset{q}{\underset{\text{a}=1}{\sum}}R_{i\text{a}j\text{a}}\varrho_{ij}\ $

$\ +\frac{2}{9}\overset{q}{\underset{\text{b}=1}{\sum}}(\varrho_{\text{bb}%
}^{M}-\varrho_{\text{bb}}^{P})(\tau^{M}-\overset{q}{\underset{\text{c}%
=1}{\sum}}\varrho_{\text{cc}}^{M})+\frac{4}{9}%
\overset{n}{\underset{i,j=q+1}{\sum}}\overset{q}{\underset{\text{b}=1}{\sum}%
}R_{i\text{b}j\text{b}}\varrho_{ij}\ $

$+\frac{1}{9}(\tau^{M}-\overset{q}{\underset{\text{a=1}}{\sum}}\varrho
_{\text{aa}})(\tau^{M}-\overset{q}{\underset{\text{b=1}}{\sum}}\varrho
_{\text{bb}})+\frac{2}{9}(\left\Vert \varrho^{M}\right\Vert ^{2}%
-\overset{q}{\underset{\text{a,b}=1}{\sum}}\varrho_{\text{ab}})$

$-\overset{n}{\underset{i,j=q+1}{\sum}}\overset{q}{\underset{\text{a,b}%
=1}{\sum}}R_{i\text{a}i\text{b}}R_{j\text{a}j\text{b}}\ -\frac{1}%
{2}\overset{n}{\underset{i,j=q+1}{\sum}}\overset{q}{\underset{\text{a,b}%
=1}{\sum}}R_{i\text{a}j\text{b}}^{2}-\overset{n}{\underset{i,j=q+1}{\sum}%
}\overset{q}{\underset{\text{a,b}=1}{\sum}}R_{i\text{a}j\text{b}}%
R_{j\text{a}i\text{b}}-\frac{1}{2}\overset{n}{\underset{i,j=q+1}{\sum}%
}\overset{q}{\underset{\text{a,b}=1}{\sum}}R_{j\text{a}i\text{b}}^{2}$

$-\frac{1}{9}\overset{n}{\underset{i,j,p,m=q+1}{\sum}}R_{ipim}R_{jpjm}%
\ -\frac{1}{18}\overset{n}{\underset{i,j,p,m=q+1}{\sum}}R_{ipjm}^{2}-\frac
{1}{9}\overset{n}{\underset{i,j,p,m=q+1}{\sum}}R_{ipjm}R_{jpim}-\frac{1}%
{18}\overset{n}{\underset{i,j,p,m=q+1}{\sum}}R_{jpim}^{2}$

$-\frac{1}{3}\overset{q}{\underset{\text{a}=1}{\sum}}%
\overset{n}{\underset{i,j,p=q+1}{\sum}}R_{i\text{a}ip}R_{j\text{a}jp}-\frac
{1}{6}\overset{q}{\underset{\text{a}=1}{\sum}}%
\overset{n}{\underset{i,j,p=q+1}{\sum}}R_{i\text{a}jp}^{2}-\frac{1}%
{3}\overset{q}{\underset{\text{a}=1i,j,}{\sum}}%
\overset{n}{\underset{p=q+1}{\sum}}R_{i\text{a}jp}R_{j\text{a}ip}$

$-\frac{1}{6}\overset{q}{\underset{\text{a}=1}{\sum}}%
\overset{n}{\underset{i,j,p=q+1}{\sum}}R_{j\text{a}ip}^{2}-\frac{1}%
{3}\overset{q}{\underset{\text{b}=1i,j,}{\sum}}%
\overset{n}{\underset{p=q+1}{\sum}}R_{i\text{b}ip}R_{j\text{b}jp}-\frac{1}%
{6}\overset{q}{\underset{\text{b}=1}{\sum}}%
\overset{n}{\underset{i,j,p=q+1}{\sum}}R_{i\text{b}jp}^{2}-\frac{1}%
{3}\overset{q}{\underset{\text{b}=1}{\sum}}%
\overset{n}{\underset{i.j,p=q+1}{\sum}}R_{i\text{b}jp}R_{j\text{b}ip}$

$-\frac{1}{6}\overset{q}{\underset{\text{b}=1}{\sum}}%
\overset{n}{\underset{i,j,p=q+1}{\sum}}R_{j\text{b}ip}^{2}](y_{0})$

$-\frac{1}{48}$ $\overset{q}{\underset{\text{a,b,c=1}}{\sum}}[$
$-\overset{n}{\underset{i=q+1}{\sum}}R_{i\text{a}i\text{a}}(R_{\text{bcbc}%
}^{P}-R_{\text{bcbc}}^{M})$ $-\overset{n}{\underset{j=q+1}{\sum}}%
R_{j\text{a}j\text{a}}(R_{\text{bcbc}}^{P}-R_{\text{bcbc}}^{M})\qquad\qquad6D$

\qquad\qquad\qquad\ $+\overset{n}{\underset{i=q+1}{\sum}}R_{i\text{a}%
i\text{b}}(R_{\text{acbc}}^{P}-R_{\text{acbc}}^{M}%
)\ -\overset{n}{\underset{i=q+1}{\sum}}R_{i\text{a}i\text{c}}(R_{\text{abbc}%
}^{P}-R_{\text{abbc}}^{M})$

$\qquad\qquad\qquad+\overset{n}{\underset{j=q+1}{\sum}}R_{j\text{a}j\text{b}%
}(R_{\text{acbc}}^{P}-R_{\text{acbc}}^{M})$%
\ $-\overset{n}{\underset{j=q+1}{\sum}}R_{j\text{a}j\text{c}}(R_{\text{abbc}%
}^{P}-R_{\text{abbc}}^{M})$

$\qquad\qquad+\underset{i,j=q+1}{\overset{n}{\sum}}$ $-R_{i\text{a}j\text{a}%
}(T_{\text{bb}i}T_{\text{cc}j}$ $-T_{\text{bc}i}T_{\text{bc}j})$
$-\underset{i,j=q+1}{\overset{n}{\sum}}R_{i\text{a}j\text{a}}(T_{\text{bb}%
j}T_{\text{cc}i}$ $-T_{\text{bc}j}T_{\text{bc}i})$

$\qquad\qquad+$ $\underset{i,j=q+1}{\overset{n}{\sum}}$ $-R_{j\text{a}%
i\text{a}}(T_{\text{bb}i}T_{\text{cc}j}$ $-T_{\text{bc}i}T_{\text{bc}j})$
$-\underset{i,j=q+1}{\overset{n}{\sum}}R_{j\text{a}i\text{a}}(T_{\text{bb}%
j}T_{\text{cc}i}$ $-T_{\text{bc}j}T_{\text{bc}i})$

$\qquad\qquad\ +\underset{i,j=q+1}{\overset{n}{\sum}}\ R_{i\text{a}j\text{b}%
}(T_{\text{ab}i}T_{\text{cc}j}-T_{\text{bc}i}T_{\text{ac}j}%
)\ +\underset{i,j=q+1}{\overset{n}{\sum}}\ R_{i\text{a}j\text{b}}%
(T_{\text{ab}j}T_{\text{cc}i}-T_{\text{bc}j}T_{\text{ac}i})$

$\qquad\qquad+\underset{i,j=q+1}{\overset{n}{\sum}}\ R_{j\text{a}i\text{ib}%
}(T_{\text{ab}i}T_{\text{cc}j}-T_{\text{bc}i}T_{\text{ac}j}%
)\ +\underset{i,j=q+1}{\overset{n}{\sum}}\ R_{j\text{a}i\text{b}}%
(T_{\text{ab}j}T_{\text{cc}i}-T_{\text{bc}j}T_{\text{ac}i})\qquad$

$\qquad\qquad+\underset{i,j=q+1}{\overset{n}{\sum}}-R_{i\text{a}j\text{c}%
}(T_{\text{ab}i}T_{\text{bc}j}-T_{\text{ac}i}T_{\text{bb}j}%
)-\underset{i,j=q+1}{\overset{n}{\sum}}R_{i\text{a}j\text{c}}(T_{\text{ba}%
j}T_{\text{bc}i}-T_{\text{ac}j}T_{\text{bb}i})$

$\qquad\qquad+\underset{i,j=q+1}{\overset{n}{\sum}}-R_{j\text{a}i\text{c}%
}(T_{\text{ba}i}T_{\text{bc}j}-T_{\text{ac}i}T_{\text{bb}j}%
)-\underset{i,j=q+1}{\overset{n}{\sum}}R_{j\text{a}i\text{c}}(T_{\text{ba}%
j}T_{\text{bc}i}-T_{\text{ac}j}T_{\text{bb}i})](y_{0})$

$\qquad\qquad+\frac{1}{144}\underset{p=q+1}{\overset{n}{\sum}}%
[\underset{i=q+1}{\overset{n}{\sum}}\overset{q}{\underset{\text{b,c=1}}{\sum}%
}R_{ipip}(R_{\text{bcbc}}^{P}-R_{\text{bcbc}}^{M}%
)+\underset{j=q+1}{\overset{n}{\sum}}$ $\overset{q}{\underset{\text{b,c=1}%
}{\sum}}R_{jpjp}(R_{\text{bcbc}}^{P}-R_{\text{bcbc}}^{M})](y_{0})$

$\qquad\qquad+\frac{1}{72}\underset{i,j,p=q+1}{\overset{n}{\sum}%
}\overset{q}{\underset{\text{b,c=1}}{\sum}}[R_{ipjp}(T_{\text{bb}%
i}T_{\text{cc}j}-T_{\text{bc}i}T_{\text{bc}j})+R_{ipjp}(T_{\text{bb}%
j}T_{\text{cc}i}-T_{\text{bc}j}T_{\text{bc}i})](y_{0})\qquad$

$\qquad\qquad-\frac{1}{288}\underset{i,j=q+1}{\overset{n}{\sum}}%
[T_{\text{aa}i}T_{\text{bb}j}(T_{\text{cc}i}T_{\text{dd}j}-T_{\text{cd}%
i}T_{\text{dc}j})+T_{\text{aa}i}T_{\text{bb}j}(T_{\text{cc}j}T_{\text{dd}%
i}-T_{\text{cd}j}T_{\text{dc}i})\qquad E$

$\qquad\qquad+T_{\text{aa}j}T_{\text{bb}i}(T_{\text{cc}i}T_{\text{dd}%
j}-T_{\text{cd}i}T_{\text{dc}j})+T_{\text{aa}j}T_{\text{bb}i}(T_{\text{cc}%
j}T_{\text{dd}i}-T_{\text{cd}j}T_{\text{dc}i})](y_{0})$

$\qquad\qquad+\frac{1}{288}\underset{i,j=q+1}{\overset{n}{\sum}}%
[T_{\text{aa}i}T_{\text{bc}j}(T_{\text{bc}i}T_{\text{dd}j}-T_{\text{bd}%
i}T_{\text{cd}j})+T_{\text{aa}i}T_{\text{bc}j}(T_{\text{bc}j}T_{\text{dd}%
i}-T_{\text{bd}j}T_{\text{cd}i})$

$\qquad\qquad+T_{\text{aa}j}T_{\text{bc}i}(T_{\text{bc}i}T_{\text{dd}%
j}-T_{\text{bd}i}T_{\text{cd}j})+T_{\text{aa}j}T_{\text{bc}i}(T_{\text{bc}%
j}T_{\text{dd}i}-T_{\text{bd}j}T_{\text{cd}i})](y_{0})$

$\qquad\qquad-\frac{1}{288}\underset{i,j=q+1}{\overset{n}{\sum}}%
[T_{\text{aa}i}T_{\text{bd}j}(T_{\text{bc}i}T_{\text{cd}j}-T_{\text{bd}%
i}T_{\text{cc}j})+T_{\text{aa}i}T_{\text{bd}j}(T_{\text{bc}j}T_{\text{cd}%
i}-T_{\text{bd}j}T_{\text{cc}i})$

$\qquad\qquad+T_{\text{aa}j}T_{\text{bd}i}(T_{\text{bc}i}T_{\text{cd}%
j}-T_{\text{bd}i}T_{\text{cc}j})+T_{\text{aa}j}T_{\text{bd}i}(T_{\text{bc}%
j}T_{\text{cd}i}-T_{\text{bd}j}T_{\text{cc}i})](y_{0})\qquad$

$\qquad\qquad+\frac{1}{288}\underset{i,j=q+1}{\overset{n}{\sum}}%
[T_{\text{ab}i}T_{\text{ab}j}(T_{\text{cc}i}T_{\text{dd}j}-T_{\text{cd}%
i}T_{\text{dc}j})+T_{\text{ab}i}T_{\text{ab}j}(T_{\text{cc}j}T_{\text{dd}%
i}-T_{\text{cd}j}T_{\text{dc}i})$

$\qquad\qquad+T_{\text{ab}j}T_{\text{ab}i}(T_{\text{cc}i}T_{\text{dd}%
j}-T_{\text{cd}i}T_{\text{dc}j})+T_{\text{ab}j}T_{\text{ab}i}(T_{\text{cc}%
j}T_{\text{dd}i}-T_{\text{cd}j}T_{\text{dc}i})](y_{0})$

$\qquad\qquad-\frac{1}{288}\underset{i,j=q+1}{\overset{n}{\sum}}%
[T_{\text{ab}i}T_{\text{bc}j}(T_{\text{ac}i}T_{\text{dd}j}-T_{\text{ad}%
i}T_{\text{cd}j})+T_{\text{ab}i}T_{\text{bc}j}(T_{\text{ac}j}T_{\text{dd}%
i}-T_{\text{ad}j}T_{\text{cd}i})$

$\qquad\qquad+T_{\text{ab}j}T_{\text{bc}i}(T_{\text{ac}i}T_{\text{dd}%
j}-T_{\text{ad}i}T_{\text{cd}j})+T_{\text{ab}j}T_{\text{bc}i}(T_{\text{ac}%
j}T_{\text{dd}i}-T_{\text{ad}j}T_{\text{cd}i})](y_{0})$

$\qquad\qquad+\frac{1}{288}\underset{i,j=q+1}{\overset{n}{\sum}}%
[T_{\text{ab}i}T_{\text{bd}j}(T_{\text{ac}i}T_{\text{cd}j}-T_{\text{ad}%
i}T_{\text{cc}j})+T_{\text{ab}i}T_{\text{bd}j}(T_{\text{ac}j}T_{\text{cd}%
i}-T_{\text{ad}j}T_{\text{cc}i})$

$\qquad\qquad+T_{\text{ab}i}T_{\text{bd}j}(T_{\text{ac}j}T_{\text{cd}%
i}-T_{\text{ad}j}T_{\text{cc}i})+T_{\text{ab}j}T_{\text{bd}i}(T_{\text{ac}%
j}T_{\text{cd}i}-T_{\text{ad}j}T_{\text{cc}i})](y_{0})$

$\qquad\qquad-\ \frac{1}{288}\underset{i,j=q+1}{\overset{n}{\sum}%
}[T_{\text{ac}i}T_{\text{ab}j}(T_{\text{bc}i}T_{\text{dd}j}-T_{\text{bd}%
i}T_{\text{dc}j})+T_{\text{ac}i}T_{\text{ab}j}(T_{\text{bc}j}T_{\text{dd}%
i}-T_{\text{bd}j}T_{\text{dc}i})$

$\qquad\qquad+T_{\text{ac}j}T_{\text{ab}i}(T_{\text{bc}i}T_{\text{dd}%
j}-T_{\text{bd}i}T_{\text{dc}j})+T_{\text{ac}j}T_{\text{ab}i}(T_{\text{bc}%
j}T_{\text{dd}i}-T_{\text{bd}j}T_{\text{dc}i})](y_{0})$

$\qquad\qquad+\ \frac{1}{288}\underset{i,j=q+1}{\overset{n}{\sum}%
}[T_{\text{ac}i}T_{\text{bb}j}(T_{\text{ac}i}T_{\text{dd}j}-T_{\text{ad}%
i}T_{\text{cd}j})+T_{\text{ac}i}T_{\text{bb}j}(T_{\text{ac}j}T_{\text{dd}%
i}-T_{\text{ad}j}T_{\text{cd}i})$

$\qquad\qquad+T_{\text{ac}j}T_{\text{bb}i}(T_{\text{ac}i}T_{\text{dd}%
j}-T_{\text{ad}i}T_{\text{cd}i})+T_{\text{ac}j}T_{\text{bb}i}(T_{\text{ac}%
j}T_{\text{dd}i}-T_{\text{ad}j}T_{\text{cd}i})](y_{0})$

$\qquad\qquad-\ \frac{1}{288}\underset{i,j=q+1}{\overset{n}{\sum}%
}[T_{\text{ac}i}T_{\text{bd}j}(T_{\text{ac}i}T_{\text{bd}j}-T_{\text{ad}%
i}T_{\text{bc}j})+T_{\text{ac}i}T_{\text{bd}j}(T_{\text{ac}j}T_{\text{bd}%
i}-T_{\text{ad}j}T_{\text{bc}i})$

$\qquad\qquad+T_{\text{ac}j}T_{\text{bd}i}(T_{\text{ac}i}T_{\text{bd}%
j}-T_{\text{ad}i}T_{\text{bc}j})+T_{\text{ac}j}T_{\text{bd}i}(T_{\text{ac}%
j}T_{\text{bd}i}-T_{\text{ad}j}T_{\text{bc}i})](y_{0})$

$\qquad\qquad+\frac{1}{288}\underset{i,j=q+1}{\overset{n}{\sum}}%
[T_{\text{ad}i}T_{\text{ab}j}(T_{\text{bc}i}T_{\text{cd}j}-T_{\text{bd}%
i}T_{\text{cc}j})+T_{\text{ad}i}T_{\text{ab}j}(T_{\text{bc}j}T_{\text{cd}%
i}-T_{\text{bd}j}T_{\text{cc}i})$

$\qquad\qquad+T_{\text{ad}j}T_{\text{ab}i}(T_{\text{bc}i}T_{\text{cd}%
j}-T_{\text{bd}i}T_{\text{cc}j})+T_{\text{ad}j}T_{\text{ab}i}(T_{\text{bc}%
j}T_{\text{cd}i}-T_{\text{bd}j}T_{\text{cc}i})](y_{0})$

$\qquad\qquad-\ \frac{1}{288}\underset{i,j=q+1}{\overset{n}{\sum}%
}[T_{\text{ad}i}T_{\text{bb}j}(T_{\text{ac}i}T_{\text{cd}j}-T_{\text{ad}%
i}T_{\text{cc}j})+T_{\text{ad}i}T_{\text{bb}j}(T_{\text{ac}j}T_{\text{cd}%
i}-T_{\text{ad}j}T_{\text{cc}i})$

$\qquad\qquad+T_{\text{ad}j}T_{\text{bb}i}(T_{\text{ac}i}T_{\text{cd}%
j}-T_{\text{ad}i}T_{\text{cc}j})+T_{\text{ad}j}T_{\text{bb}i}(T_{\text{ac}%
j}T_{\text{cd}i}-T_{\text{ad}j}T_{\text{cc}i})](y_{0})$

$\qquad\qquad+\ \frac{1}{288}\underset{i,j=q+1}{\overset{n}{\sum}%
}[T_{\text{ad}i}T_{\text{bc}j}(T_{\text{ac}i}T_{\text{bd}j}-T_{\text{ad}%
i}T_{\text{bc}j})+T_{\text{ad}i}T_{\text{bc}j}(T_{\text{ac}j}T_{\text{bd}%
i}-T_{\text{ad}j}T_{\text{bc}i})$

$\qquad\qquad+T_{\text{ad}j}T_{\text{bc}i}(T_{\text{ac}i}T_{\text{bd}%
j}-T_{\text{ad}i}T_{\text{bc}j})+T_{\text{ad}j}T_{\text{bc}i}(T_{\text{ac}%
j}T_{\text{bd}i}-T_{\text{ad}j}T_{\text{bc}i})](y_{0})$

$\qquad\qquad-\ \frac{1}{144}[(R_{\text{cdcd}}^{P}-R_{\text{cdcd}}%
^{M})(R_{\text{abab}}^{P}-R_{\text{abab}}^{M})](y_{0})\qquad\left(  1\right)
$

$\qquad\qquad\ +\frac{1}{144}[(R_{\text{bdcd}}^{P}-R_{\text{bdcd}}%
^{M})(R_{\text{abac}}^{P}-R_{\text{abac}}^{M})](y_{0})\qquad(2)$

$\ \qquad\qquad+\ \frac{1}{144}[(R_{\text{bcdc}}^{P}-R_{\text{bcdc}}%
^{M})(R_{\text{abad}}^{P}-R_{\text{abad}}^{M})](y_{0})\qquad(3)$

$\qquad\qquad\ -\ \frac{1}{144}[(R_{\text{adcd}}^{P}-R_{\text{adcd}}%
^{M})(R_{\text{abbc}}^{P}-R_{\text{abbc}}^{M})](y_{0})\qquad(4)\qquad$

$\ \qquad\qquad+\ \frac{1}{144}[(R_{\text{acdc}}^{P}-R_{\text{acdc}}%
^{M})(R_{\text{abdb}}^{P}-R_{\text{abdb}}^{M})](y_{0})\qquad(5)$

$\ \qquad\qquad-\ \frac{1}{576}[(R_{\text{abcd}}^{P}-R_{\text{abcd}}^{M}%
)]^{2}(y_{0})\qquad(6)$

\qquad\qquad$-\frac{1}{24}\times\frac{1}{6}<H,i>(y_{0})<H,j>(y_{0}%
)\qquad\qquad\qquad\qquad\qquad\qquad\qquad L_{3}$

$\times\lbrack2\varrho_{ij}+$ $\overset{q}{\underset{\text{a}=1}{4\sum}%
}R_{i\text{a}j\text{a}}-3\overset{q}{\underset{\text{a,b=1}}{\sum}%
}(T_{\text{aa}i}T_{\text{bb}j}-T_{\text{ab}i}T_{\text{ab}j}%
)-3\overset{q}{\underset{\text{a,b=1}}{\sum}}(T_{\text{aa}j}T_{\text{bb}%
i}-T_{\text{ab}j}T_{\text{ab}i}](y_{0})\phi(y_{0})$

$\qquad\qquad-\frac{1}{24}\times\frac{3}{2}[<H,i>^{2}(y_{0})<H,j>^{2}%
](y_{0})\phi(y_{0})$

$\qquad\qquad-\frac{1}{24}\times\frac{1}{6}<H,i>(y_{0})<H,j>(y_{0})$

$\times\lbrack2\varrho_{ij}+$ $\overset{q}{\underset{\text{a}=1}{4\sum}%
}R_{i\text{a}j\text{a}}-3\overset{q}{\underset{\text{a,b=1}}{\sum}%
}(T_{\text{aa}i}T_{\text{bb}j}-T_{\text{ab}i}T_{\text{ab}j}%
)-3\overset{q}{\underset{\text{a,b=1}}{\sum}}(T_{\text{aa}j}T_{\text{bb}%
i}-T_{\text{ab}j}T_{\text{ab}i}](y_{0})\phi(y_{0})$

\qquad\qquad$\ -\frac{1}{24}\times\frac{1}{3}<H,i>(y_{0})<H,k>(y_{0}%
)R_{jijk}(y_{0})\phi(y_{0})$

\qquad\qquad$-\frac{1}{24}\times\frac{3}{2}<H,i>^{2}(y_{0})<H,j>^{2}%
(y_{0})\phi(y_{0})$

$-\frac{1}{24}\times\frac{1}{3}<H,i>^{2}(y_{0})[\varrho_{jj}+$
$\overset{q}{\underset{\text{a}=1}{2\sum}}R_{j\text{a}j\text{a}}%
-3\overset{q}{\underset{\text{a,b=1}}{\sum}}(T_{\text{aa}j}T_{\text{bb}%
j}-T_{\text{ab}j}T_{\text{ab}j})](y_{0})\phi(y_{0})$

$+\frac{1}{24}\times\frac{15}{4}[<H,i>^{2}<H,j>^{2}](y_{0})\phi(y_{0})$

$+\frac{1}{24}\times\frac{1}{2}<H,i>(y_{0})<H,j>$

$\times\lbrack2\varrho_{ij}+$ $\overset{q}{\underset{\text{a}=1}{4\sum}%
}R_{i\text{a}j\text{a}}-3\overset{q}{\underset{\text{a,b=1}}{\sum}%
}(T_{\text{aa}i}T_{\text{bb}j}-T_{\text{ab}i}T_{\text{ab}j}%
)-3\overset{q}{\underset{\text{a,b=1}}{\sum}}(T_{\text{aa}j}T_{\text{bb}%
i}-T_{\text{ab}j}T_{\text{ab}i}](y_{0})\phi(y_{0})$

$+\frac{1}{24}\times\frac{1}{2}<H,i>(y_{0})<H,i>(y_{0})[\tau^{M}\ -3\tau
^{P}+\ \underset{\text{a}=1}{\overset{\text{q}}{\sum}}\varrho_{\text{aa}}%
^{M}+$ $\overset{q}{\underset{\text{a},\text{b}=1}{\sum}}R_{\text{abab}}^{M}$
$](y_{0})\phi(y_{0})$

$+\frac{1}{24}\times\frac{1}{6}<H,i>(y_{0})[\nabla_{i}\varrho_{jj}%
-2\varrho_{ij}<H,j>+\overset{q}{\underset{\text{a}=1}{\sum}}(\nabla
_{i}R_{\text{a}j\text{a}j}-4R_{i\text{a}j\text{a}}<H,j>)$

$+4\overset{q}{\underset{\text{a,b=1}}{\sum}}R_{i\text{a}j\text{b}%
}T_{\text{ab}j}+2\overset{q}{\underset{\text{a,b,c=1}}{\sum}}(T_{\text{aa}%
i}T_{\text{bb}j}T_{\text{cc}j}-3T_{\text{aa}i}T_{\text{bc}j}T_{\text{bc}%
j}+2T_{\text{ab}i}T_{\text{bc}j}T_{\text{ca}j})](y_{0})\phi(y_{0})$%
\qquad\qquad\qquad\qquad\qquad\ \ 

$+\frac{1}{24}\times\frac{1}{6}<H,i>(y_{0})[\nabla_{j}\varrho_{ij}%
-2\varrho_{ij}<H,j>+\overset{q}{\underset{\text{a}=1}{\sum}}(\nabla
_{j}R_{\text{a}i\text{a}j}-4R_{j\text{a}i\text{a}}<H,j>)$

$+4\overset{q}{\underset{\text{a,b=1}}{\sum}}R_{j\text{a}i\text{b}%
}T_{\text{ab}j}+2\overset{q}{\underset{\text{a,b,c=1}}{\sum}}(T_{\text{aa}%
j}T_{\text{bb}i}T_{\text{cc}j}-3T_{\text{aa}j}T_{\text{bc}i}T_{\text{bc}%
j}+2T_{\text{ab}j}T_{\text{bc}i}T_{\text{ac}j})](y_{0})\phi(y_{0})$

$+\frac{1}{24}\times\frac{1}{6}<H,i>(y_{0})[\nabla_{j}\varrho_{ij}%
-2\varrho_{jj}<H,i>+\overset{q}{\underset{\text{a}=1}{\sum}}(\nabla
_{j}R_{\text{a}i\text{a}j}-4R_{j\text{a}j\text{a}}<H,i>)$

$+4\overset{q}{\underset{\text{a,b=1}}{\sum}}R_{j\text{a}j\text{b}%
}T_{\text{ab}i}+2\overset{q}{\underset{\text{a,b,c=1}}{\sum}}(T_{\text{aa}%
j}T_{\text{bb}j}T_{\text{cc}i}-3T_{\text{aa}j}T_{\text{bc}j}T_{\text{bc}%
i}+2T_{\text{ab}j}T_{\text{bc}j}T_{\text{ac}i})](y_{0})\phi(y_{0})$

$-\frac{1}{192}<H,i>^{2}<H,j>^{2}(y_{0})\phi(y_{0})\qquad\qquad\qquad
\qquad\qquad\qquad\qquad$I$_{3213}\ $

$-\frac{1}{288}<H,i>^{2}(y_{0})[\tau^{M}\ -3\tau^{P}+\ \underset{\text{a}%
=1}{\overset{\text{q}}{\sum}}\varrho_{\text{aa}}^{M}+$
$\overset{q}{\underset{\text{a},\text{b}=1}{\sum}}R_{\text{abab}}^{M}$
$](y_{0})\phi(y_{0})$

$-\frac{1}{288}<H,i>(y_{0})<H,j>(y_{0})$

$\times\lbrack2\varrho_{ij}+$ $\overset{q}{\underset{\text{a}=1}{4\sum}%
}R_{i\text{a}j\text{a}}-3\overset{q}{\underset{\text{a,b=1}}{\sum}%
}(T_{\text{aa}i}T_{\text{bb}j}-T_{\text{ab}i}T_{\text{ab}j}%
)-3\overset{q}{\underset{\text{a,b=1}}{\sum}}(T_{\text{aa}j}T_{\text{bb}%
i}-T_{\text{ab}j}T_{\text{ab}i}](y_{0})\phi(y_{0})$

$\ +\frac{1}{144}<H,i>(y_{0})<H,k>(y_{0})R_{jijk}(y_{0})$

$+\frac{1}{144}<H,i>^{2}(y_{0})[\varrho_{jj}+$ $\overset{q}{\underset{\text{a}%
=1}{2\sum}}R_{j\text{a}j\text{a}}-3\overset{q}{\underset{\text{a,b=1}}{\sum}%
}(T_{\text{aa}j}T_{\text{bb}j}-T_{\text{ab}j}T_{\text{ab}j})](y_{0})\phi
(y_{0})$

$-\frac{1}{288}<H,i>(y_{0})[\nabla_{i}\varrho_{jj}-2\varrho_{ij}%
<H,j>+\overset{q}{\underset{\text{a}=1}{\sum}}(\nabla_{i}R_{\text{a}%
j\text{a}j}-4R_{i\text{a}j\text{a}}<H,j>)$

$+4\overset{q}{\underset{\text{a,b=1}}{\sum}}R_{i\text{a}j\text{b}%
}T_{\text{ab}j}+2\overset{q}{\underset{\text{a,b,c=1}}{\sum}}(T_{\text{aa}%
i}T_{\text{bb}j}T_{\text{cc}j}-3T_{\text{aa}i}T_{\text{bc}j}T_{\text{bc}%
j}+2T_{\text{ab}i}T_{\text{bc}j}T_{\text{ca}j})](y_{0})\phi(y_{0})$%
\qquad\qquad\qquad\qquad\qquad\ \ 

$-\frac{1}{288}<H,i>(y_{0})[\nabla_{j}\varrho_{ij}-2\varrho_{ij}%
<H,j>+\overset{q}{\underset{\text{a}=1}{\sum}}(\nabla_{j}R_{\text{a}%
i\text{a}j}-4R_{j\text{a}i\text{a}}<H,j>)$

$+4\overset{q}{\underset{\text{a,b=1}}{\sum}}R_{j\text{a}i\text{b}%
}T_{\text{ab}j}+2\overset{q}{\underset{\text{a,b,c=1}}{\sum}}(T_{\text{aa}%
j}T_{\text{bb}i}T_{\text{cc}j}-3T_{\text{aa}j}T_{\text{bc}i}T_{\text{bc}%
j}+2T_{\text{ab}j}T_{\text{bc}i}T_{\text{ac}j})](y_{0})\phi(y_{0})$

$-\frac{1}{288}<H,i>(y_{0})[\nabla_{j}\varrho_{ij}-2\varrho_{jj}%
<H,i>+\overset{q}{\underset{\text{a}=1}{\sum}}(\nabla_{j}R_{\text{a}%
i\text{a}j}-4R_{j\text{a}j\text{a}}<H,i>)$

$+4\overset{q}{\underset{\text{a,b=1}}{\sum}}R_{j\text{a}j\text{b}%
}T_{\text{ab}i}+2\overset{q}{\underset{\text{a,b,c=1}}{\sum}}(T_{\text{aa}%
j}T_{\text{bb}j}T_{\text{cc}i}-3T_{\text{aa}j}T_{\text{bc}j}T_{\text{bc}%
i}+2T_{\text{ab}j}T_{\text{bc}j}T_{\text{ac}i})](y_{0})\phi(y_{0})$

$+\frac{1}{24}[\left\Vert \text{X}\right\Vert _{M}^{2}+\operatorname{div}%
$X$_{M}-\left\Vert \text{X}\right\Vert _{P}^{2}-\operatorname{div}X_{P}%
](y_{0})[\left\Vert \text{X}\right\Vert _{M}^{2}-\operatorname{div}$%
X$_{M}-\left\Vert \text{X}\right\Vert _{P}^{2}+\operatorname{div}$%
X$_{P}](y_{0})\phi(y_{0})\qquad$\ I$_{3212}$

$+\frac{1}{6}X_{i}(y_{0})T_{\text{ab}i}(y_{0})T_{\text{ab}j}(y_{0})X_{j}%
(y_{0})+\frac{1}{3}\perp_{\text{a}ij}(y_{0})X_{i}(y_{0})[\frac{\partial X_{j}%
}{\partial x_{\text{a}}}-\perp_{\text{a}jk}X_{k}](y_{0})\phi(y_{0})\qquad
$I$_{32122}\qquad Q_{1}$

$+$ $\frac{2}{3}X_{i}(y_{0})X_{j}(y_{0})\frac{\partial X_{j}}{\partial
x_{\text{a}}}(y_{0})-\frac{1}{6}X_{i}(y_{0})\frac{\partial^{2}X_{j}}{\partial
x_{\text{a}}\partial x_{j}}(y_{0})\qquad\qquad Q_{2}$

$-\frac{1}{12}X_{i}(y_{0})\frac{\partial^{2}X_{i}}{\partial x_{\text{a}}^{2}%
}(y_{0})+\frac{1}{12}X_{i}^{2}(y_{0})[\operatorname{div}X_{M}-\left\Vert
X\right\Vert _{M}^{2}+\left\Vert X\right\Vert _{P}^{2}-\operatorname{div}%
X_{P}-$ $<H,j>X_{j}](y_{0})\phi(y_{0})$

$+\frac{1}{6}X_{i}(y_{0})X_{j}(y_{0})\frac{\partial X_{i}}{\partial x_{j}%
}(y_{0})+$ $\frac{1}{18}X_{i}(y_{0})X_{k}(y_{0})[R_{jijk}](y_{0})\phi
(y_{0})-\frac{1}{12}X_{i}(y_{0})\frac{\partial^{2}X_{i}}{\partial x_{j}^{2}%
}(y_{0})\phi(y_{0})$

$+\frac{1}{12}[R_{\text{a}i\text{a}k}-\underset{\text{c=1}}{\overset{\text{q}%
}{\sum}}T_{\text{ac}i}T_{\text{ac}k}-\perp_{\text{a}ik}\perp_{\text{a}%
jk}](y_{0})X_{k}(y_{0})+\frac{1}{18}R_{ijkj}(y_{0})X_{i}(y_{0})X_{k}(y_{0})$

$+\frac{1}{12}<H,j>(y_{0})X_{i}(y_{0})[X_{i}X_{j}-\frac{1}{2}\left(
\frac{\partial X_{j}}{\partial x_{i}}+\frac{\partial X_{i}}{\partial x_{j}%
}\right)  ](y_{0})\phi(y_{0})\qquad\qquad\qquad\qquad\qquad\qquad\qquad
\qquad\qquad\qquad\ \qquad\qquad\qquad\qquad\qquad\qquad\qquad\qquad
\qquad\qquad\qquad$

$-\frac{1}{6}[-R_{\text{a}i\text{b}i}+5\overset{q}{\underset{\text{c}=1}{\sum
}}T_{\text{ac}i}T_{\text{bc}i}+2\overset{n}{\underset{j=q+1}{\sum}}%
\perp_{\text{a}ij}\perp_{\text{b}ij}](y_{0})\underset{k=q+1}{\overset{n}{\sum
}}T_{\text{ab}k}(y_{0})X_{k}(y_{0})\phi(y_{0})\qquad\qquad$I$_{32123}\qquad
S_{1}\qquad$

$-\frac{2}{9}\underset{j=q+1}{\overset{n}{\sum}}R_{i\text{a}ij}(y_{0}%
)[\frac{\partial X_{j}}{\partial x_{\text{a}}}%
-\underset{k=q+1}{\overset{n}{\sum}}\perp_{\text{a}jk}X_{k}](y_{0})$

$+\frac{1}{12}\times\frac{2}{3}\underset{j,k=q+1}{\overset{n}{\sum}}%
R_{ijik}(y_{0})[X_{j}X_{k}-\frac{1}{2}(\frac{\partial X_{j}}{\partial x_{k}%
}+\frac{\partial X_{k}}{\partial x_{j}})](y_{0})\phi(y_{0})$

$-\frac{1}{6}T_{\text{ab}i}(y_{0})\frac{\partial^{2}X_{i}}{\partial
x_{\text{a}}\partial x_{\text{b}}}(y_{0})\phi(y_{0})\qquad\qquad\qquad
S_{2}\qquad\qquad S_{21}\qquad\qquad\qquad\qquad\qquad\qquad$

$+$ $\frac{1}{12}T_{\text{ab}i}(y_{0})$

$\times\lbrack$ $(R_{\text{a}i\text{b}j}+R_{\text{a}j\text{b}i})$
$-\underset{\text{c=1}}{\overset{\text{q}}{\sum}}(T_{\text{ac}i}T_{\text{bc}%
j}+T_{\text{ac}j}T_{\text{bc}i})-\overset{n}{\underset{k=q+1}{\sum}}%
(\perp_{\text{a}ik}\perp_{\text{b}jk}+$ $\perp_{\text{a}jk}\perp_{\text{b}%
ik})](y_{0})X_{j}(y_{0})\phi(y_{0})$

$\qquad-$ $\frac{1}{6}T_{\text{ab}i}(y_{0})T_{\text{ab}j}(y_{0})[X_{i}%
X_{j}-\frac{1}{2}\left(  \frac{\partial X_{i}}{\partial x_{j}}+\frac{\partial
X_{j}}{\partial x_{i}}\right)  ](y_{0})\phi(y_{0})$

$\qquad-\frac{1}{3}\perp_{\text{a}ij}(y_{0})[(X_{i}\frac{\partial X_{j}%
}{\partial x_{\text{a}}}+X_{j}\frac{\partial X_{i}}{\partial x_{\text{a}}%
})-\frac{1}{4}\left(  \frac{\partial^{2}X_{i}}{\partial x_{\text{a}}\partial
x_{j}}+\frac{\partial^{2}X_{j}}{\partial x_{\text{a}}\partial x_{i}}\right)
](y_{0})\phi(y_{0})\qquad S_{22}$

$\qquad-\frac{1}{6}\perp_{\text{a}ij}(y_{0})[T_{\text{ab}j}\frac{\partial
X_{i}}{\partial x_{\text{b}}}](y_{0})$

$\qquad+\frac{1}{6}\perp_{\text{a}ij}(y_{0})[(\perp_{\text{b}ik}T_{\text{ab}%
j})+\frac{2}{3}(2R_{\text{a}ijk}+R_{\text{a}jik}+R_{\text{a}kji})](y_{0}%
)X_{k}(y_{0})\phi(y_{0})$

$\qquad-\frac{1}{6}\perp_{\text{a}ij}(y_{0})\perp_{\text{a}jk}(y_{0}%
)[X_{i}X_{k}-\frac{1}{2}\left(  \frac{\partial X_{i}}{\partial x_{k}}%
+\frac{\partial X_{k}}{\partial x_{i}}\right)  ](y_{0})\phi(y_{0})\qquad
\qquad\qquad\qquad\qquad\qquad\qquad\qquad$

$\qquad+\frac{1}{12}[(\frac{\partial X_{j}}{\partial x_{\text{a}}})^{2}%
+X_{j}\frac{\partial^{2}X_{j}}{\partial x_{\text{a}}^{2}}-\frac{1}{2}%
\frac{\partial^{3}X_{j}}{\partial x_{\text{a}}^{2}\partial x_{j}}%
](y_{0})-\frac{1}{6}\overset{n}{\underset{k=q+1}{\sum}}[\perp_{\text{b}ik}%
$T$_{\text{aa}k}\frac{\partial X_{i}}{\partial x_{\text{b}}^{2}}](y_{0}%
)\phi(y_{0})\qquad\qquad S_{3}\qquad S_{31}$

$\qquad+\frac{1}{144}[\{4\nabla_{i}R_{i\text{a}j\text{a}}+2\nabla
_{j}R_{i\text{a}i\text{a}}+$ $8(\overset{q}{\underset{\text{c=1}}{%
{\textstyle\sum}
}}R_{\text{a}i\text{c}i}^{{}}T_{\text{ac}j}+\;\overset{n}{\underset{k=q+1}{%
{\textstyle\sum}
}}R_{\text{a}iik}\perp_{\text{a}jk})$

$\qquad+8(\overset{q}{\underset{\text{c=1}}{%
{\textstyle\sum}
}}R_{\text{a}i\text{c}j}^{{}}T_{\text{ac}i}+\;\overset{n}{\underset{k=q+1}{%
{\textstyle\sum}
}}R_{\text{a}ijk}\perp_{\text{a}ik})+8(\overset{q}{\underset{\text{c=1}}{%
{\textstyle\sum}
}}R_{\text{a}j\text{c}i}^{{}}T_{\text{ac}i}+\;\overset{n}{\underset{k=q+1}{%
{\textstyle\sum}
}}R_{\text{a}jik}\perp_{\text{a}ik})\}$\ 

$\qquad\ +\frac{2}{3}\underset{k=q+1}{\overset{n}{\sum}}\{T_{\text{aa}%
k}(R_{ijik}+3\overset{q}{\underset{\text{c}=1}{\sum}}\perp_{\text{c}ij}%
\perp_{\text{c}ik})\}](y_{0})X_{k}(y_{0})\phi(y_{0})$

$-\frac{1}{12}[$ R$_{\text{a}i\text{a}k}$ $-\underset{\text{c=1}%
}{\overset{\text{q}}{\sum}}T_{\text{ac}i}T_{\text{ac}k}%
-\overset{n}{\underset{l=q+1}{\sum}}(\perp_{\text{a}il}\perp_{\text{a}%
kl}](y_{0})\times\lbrack X_{i}X_{k}-\frac{1}{2}\left(  \frac{\partial X_{i}%
}{\partial x_{k}}+\frac{\partial X_{k}}{\partial x_{i}}\right)  ](y_{0}%
)\phi(y_{0})$

$-\frac{1}{24}T_{\text{aa}k}(y_{0})[-X_{i}^{2}X_{k}+X_{k}\frac{\partial X_{i}%
}{\partial x_{i}}\ +X_{i}\left(  \frac{\partial X_{k}}{\partial x_{i}}%
+\frac{\partial X_{i}}{\partial x_{k}}\right)  -\frac{1}{3}\left(
\frac{\partial^{2}X_{k}}{\partial x_{i}^{2}}+2\frac{\partial^{2}X_{i}%
}{\partial x_{i}\partial x_{k}}\right)  ](y_{0})\phi(y_{0})$

$\qquad+\frac{1}{18}[R_{\text{a}jij}\frac{\partial X_{i}}{\partial
x_{\text{a}}^{2}}](y_{0})\qquad\qquad\qquad\qquad\qquad\qquad\qquad
\qquad\qquad S_{32}$

$\qquad+\frac{1}{24}[\frac{4}{3}\overset{q}{\underset{\text{a}=1}{\sum}}%
\perp_{\text{a}ki}R_{ij\text{a}j}-\frac{1}{3}(\nabla_{i}R_{kjij}+\nabla
_{j}R_{ijik}+\nabla_{k}R_{ijij})](y_{0})X_{k}(y_{0})$

\qquad\ $-\frac{1}{18}R_{ijkj}(y_{0})[X_{i}X_{k}-\frac{1}{2}\left(
\frac{\partial X_{i}}{\partial x_{k}}+\frac{\partial X_{k}}{\partial x_{i}%
}\right)  ](y_{0})\phi(y_{0})$

$\qquad+\frac{1}{24}[X_{i}^{2}X_{j}^{2}-2X_{i}X_{j}\left(  \frac{\partial
X_{j}}{\partial x_{i}}+\frac{\partial X_{i}}{\partial x_{j}}\right)
-X_{i}^{2}\frac{\partial X_{j}}{\partial x_{j}}-X_{j}^{2}\frac{\partial X_{i}%
}{\partial x_{i}}](y_{0})\phi(y_{0})$

$\qquad+\frac{1}{48}\left(  \frac{\partial X_{j}}{\partial x_{i}}%
+\frac{\partial X_{i}}{\partial x_{j}}\right)  ^{2}(y_{0})+\frac{1}{24}\left(
\frac{\partial X_{i}}{\partial x_{i}}\frac{\partial X_{j}}{\partial x_{j}%
}\right)  (y_{0})\phi(y_{0})\qquad$

$+\frac{1}{36}X_{i}(y_{0})\left(  2\frac{\partial^{2}X_{j}}{\partial
x_{i}\partial x_{j}}+\frac{\partial^{2}X_{i}}{\partial x_{j}^{2}}\right)
(y_{0})\phi(y_{0})+\frac{1}{36}X_{j}(y_{0})\left(  \frac{\partial^{2}X_{j}%
}{\partial x_{i}^{2}}+2\frac{\partial^{2}X_{i}}{\partial x_{i}\partial x_{j}%
}\right)  (y_{0})\phi(y_{0})$

$-\frac{1}{48}\left(  \frac{\partial^{3}X_{i}}{\partial x_{i}\partial
x_{j}^{2}}+\frac{\partial^{3}X_{j}}{\partial x_{i}^{2}\partial x_{j}}\right)
(y_{0})\phi(y_{0})$

$+$ $\frac{2}{3}<H,j>(y_{0})\left(  \frac{\partial^{2}X_{i}}{\partial
x_{i}\partial x_{j}}+2\frac{\partial^{2}X_{j}}{\partial x_{i}^{2}}\right)
(y_{0})\phi(y_{0})+$ $\frac{2}{3}<H,j>(y_{0})R_{ijik}(y_{0})X_{k}(y_{0}%
)\phi(y_{0})\qquad$I$_{3213}$

$+\frac{1}{12}[<H,i><H,j>\ +\frac{1}{6}(2\varrho_{ij}%
+4\overset{q}{\underset{\text{a}=1}{\sum}}R_{i\text{a}j\text{a}}%
-6\overset{q}{\underset{\text{a,b}=1}{\sum}}T_{\text{aa}i}T_{\text{bb}%
j}-T_{\text{ab}i}T_{\text{ab}j})](y_{0})\phi(y_{0})$

$\times\frac{1}{2}[\left(  \frac{\partial X_{j}}{\partial x_{i}}%
-\frac{\partial X_{i}}{\partial x_{j}}\right)  ](y_{0})\phi(y_{0})$

$-\frac{1}{12}\perp_{\text{a}ij}(y_{0})<H,i>(y_{0})[(X_{j}\perp_{\text{a}%
ij}-\frac{\partial X_{i}}{\partial x_{\text{a}}})+\frac{\partial X_{\text{a}}%
}{\partial x_{i}}](y_{0})\phi(y_{0})$

$-\frac{1}{18}[X_{j}\left(  2\frac{\partial^{2}X_{j}}{\partial x_{i}^{2}%
}+\frac{\partial^{2}X_{i}}{\partial x_{i}\partial x_{j}}\right)  ](y_{0}%
)\phi(y_{0})-\frac{1}{12}[\left(  \frac{\partial X_{i}}{\partial x_{j}}%
+\frac{\partial X_{j}}{\partial x_{i}}\right)  ]\frac{\partial X_{j}}{\partial
x_{i}}(y_{0})\phi(y_{0})\qquad\ $I$_{3214}\qquad\qquad\qquad$\qquad
\qquad\qquad\qquad\qquad\qquad\qquad\qquad\qquad\qquad\qquad\qquad\qquad
\qquad\qquad\qquad

$+\frac{1}{12}\frac{\partial^{2}\text{V}}{\partial x_{i}^{2}}(y_{0})\phi
(y_{0})\qquad\qquad$I$_{3215}$\qquad\qquad

\chapter{Computation of the Second Coefficient$\qquad\ $}

\section{The Computation of b$_{1}$(y$_{0}$,P,$\phi$)}

The second coefficient is given in (iii) of \textbf{Theorem 4} by:

$\left(  C_{1}\right)  $\qquad b$_{1}($x,P,$\phi$) = $\int_{0}^{1}%
$F(1,1-r$_{1}$)L$_{\Psi}[\phi\circ\pi_{\text{P}}$](x)dr$_{1}$

It is well known that for L =$\frac{1}{2}\Delta+X+V$ we have for smooth
function f,g :M$\rightarrow$R and $\phi\in\Gamma(E):$

$\qquad\frac{\partial^{2}}{\partial x^{2}}$(fg) =$\frac{\partial^{2}%
f}{\partial x^{2}}$g + f$\frac{\partial^{2}g}{\partial x^{2}}$ +
2$\frac{\partial f}{\partial x}\frac{\partial g}{\partial x}$

\qquad$\nabla$(f$\phi$) = ($\nabla$f)$\phi$ + f($\nabla\phi$)

\qquad$\frac{1}{2}\Delta$(f$\phi$) =$\frac{1}{2}$($\Delta$f)$\phi$ + $\frac
{1}{2}$f($\Delta\phi$) +
$<$%
$\nabla$f,$\nabla\phi$%
$>$%

\qquad L(f$\phi$) = (Lf)$\phi$ + f(L$\phi$) +
$<$%
$\nabla$f,$\nabla\phi$%
$>$
$-$ V$($f$\phi)$\qquad\qquad\qquad

Using the definition of F(1,1-r$_{1}$) we have:

$\left(  C_{2}\right)  \qquad$F(1,1-r$_{1}$)L$_{\Psi}[\phi\circ\pi_{\text{P}}%
$](x) = L$_{\Psi}$[$\phi\circ\pi_{\text{P}}$]$(\gamma($r$_{1}))$

where $\gamma:\left[  0,1\right]  \longrightarrow$M$_{0}$ is the unique
minimal geodesic from a point x$\in$M$_{0}$ and meeting P orthogonally at a
point y$\in$P in time 1: $\gamma(0)=x\in M_{0}$ and $\gamma(1)=y\in P.$

So in Fermi coodinates, we write:

$\left(  C_{3}\right)  \qquad\gamma(s)=y+(1-s)(x-y)$ (vector-wise).

\qquad\qquad$=$ $\left(  \text{x}_{1},...,\text{x}_{q},(1-s)\text{x}%
_{q+1},...,(1-s)\text{x}_{n}\right)  $ (in local coordinates).

\qquad\qquad$=$ $\left(  \text{y}_{1}(y),...,\text{y}_{q}(y),(1-s)\text{x}%
_{q+1}(x),...,(1-s)\text{x}_{n}(x)\right)  $

The last equality in $\left(  C_{3}\right)  $ is due to the definition of
Fermi coordinates. See for example $\left(  2.2\right)  -\left(  2.3\right)  $
of \textbf{Gray }$\left[  4\right]  .$

We see that y = $\pi_{\text{P}}$(x) where $\pi_{\text{P}}:$M$_{0}%
\longrightarrow$P is the projection viewed in Fermi coordinates.

By $\left(  C_{3}\right)  $ and the definition of the geodesic $\gamma:\left[
0,1\right]  \longrightarrow$M$_{0},$ we have:

$\left(  C_{4}\right)  \qquad\gamma(0)=x=\left(  \text{y}_{1}(y),...,\text{y}%
_{q}(y),\text{x}_{q+1}(x),...,\text{x}_{n}(x)\right)  $

\qquad$\qquad\gamma(1)=y=\left(  \text{y}_{1}(y),...,\text{y}_{q}%
(y),0,...,0\right)  =\pi_{\text{P}}(x)$

By the definition of the geodesic $\gamma$ in $\left(  C_{3}\right)  ,$ we
have for r$_{1}\in\left[  0,1\right]  ,$

$\left(  C_{5}\right)  \qquad z_{0}=\gamma($r$_{1})=\left(  \text{y}%
_{1}(y),...,\text{y}_{q}(y),(1-r_{1})\text{x}_{q+1}(x_{0}),...,(1-r_{1}%
)\text{x}_{n}(x_{0})\right)  $

Consequently,

$\left(  C_{6}\right)  \qquad\pi_{\text{P}}$(z$_{0}$) $=\left(  \text{y}%
_{1}(y),...,\text{y}_{q}(y),0,...,0\right)  =$ y $=\pi_{\text{P}}$(x)

From $\left(  C_{3}\right)  :$

$\left(  C_{7}\right)  $ (i)$\qquad\overset{.}{\gamma}$(s) = $\left(
0,...,0,-\text{x}_{q+1}(x_{0}),...,-\text{x}_{n}(x_{0})\right)  $

$\qquad\ $(ii)$\qquad\frac{\partial}{\partial\text{x}_{i}}\gamma(s)=\left\{
\begin{array}
[c]{c}%
1\text{ for }i=1,...,q\\
(1-\text{s})\text{ for }i=q+1,...,n
\end{array}
\right.  $

\qquad\ (iii)\qquad$\frac{\partial^{2}}{\partial\text{x}_{i}\partial
\text{x}_{j}}\gamma(s)=0$ for $i,j=1,...,q,q+1,...,n$

From $\left(  C_{6}\right)  :$

$\left(  C_{8}\right)  $ $(i)\qquad\frac{\partial}{\partial\text{x}_{i}}%
\pi_{\text{P}}$(z$_{0}$) $=\frac{\partial}{\partial\text{x}_{i}}\pi_{\text{P}%
}(x)=\left\{
\begin{array}
[c]{c}%
1\text{ for }i=1,...,q\\
0\text{ for }i=q+1,...,n
\end{array}
\right.  $

$\qquad(ii)\qquad\qquad\frac{\partial^{2}}{\partial\text{x}_{i}\partial
\text{x}_{j}}\pi_{\text{P}}$(z$_{0}$) $=0=\frac{\partial^{2}}{\partial
\text{x}_{j}\partial\text{x}_{i}}\pi_{\text{P}}$(x) for all
$i,j=1,...,q,q+1,...,n$

Using the definition of L$_{\Psi}$ in $\left(  5.31\right)  ,$ we have:

$\left(  C_{9}\right)  \qquad\qquad$L$_{\Psi}[\phi\circ\pi_{\text{P}}](z_{0}$)
$=\frac{\text{L[}\Psi\phi\circ\pi_{\text{P}}\text{]}}{\Psi}$(z$_{0}$)

$\qquad\qquad\qquad=$ $\frac{\text{L}\Psi}{\Psi}(z_{0}).[\phi\circ
\pi_{\text{P}}$](z$_{0}$) + L$[\phi\circ\pi_{\text{P}}$](z$_{0}$)

$\qquad\qquad\qquad+$
$<$%
$\nabla\log\Psi,\nabla\lbrack\phi\circ\pi_{\text{P}}$]%
$>$%
(z$_{0}$) $-$ V$(z_{0})\phi\circ\pi_{P}$(z$_{0}$)

We denote the \textbf{FOUR terms }in $\left(  C_{9}\right)  $ by:

\qquad\qquad\qquad T$_{1}$ = $\frac{\text{L}\Psi}{\Psi}(z_{0}).[\phi\circ
\pi_{\text{P}}$](z$_{0})$

\qquad\qquad\qquad T$_{2}$ \ = L[$\phi\circ\pi_{\text{P}}$](z$_{0})$

\qquad\qquad\qquad T$_{3}$ = \
$<$%
$\nabla\log\Psi,\nabla\lbrack\phi\circ\pi_{\text{P}}]$%
$>$%
(z$_{0})$

\qquad\qquad\qquad T$_{4}=$ $-$ V$(z_{0})\phi\circ\pi_{P}$(z$_{0}$)

\textbf{COMPUTATION OF T}$_{1}$\textbf{:}

$\qquad\left(  C_{10}\right)  \qquad$\textbf{T}$_{1}$\textbf{ }$\mathbf{=}$
$\frac{\text{L}\Psi}{\Psi}\left(  z_{0}\right)  .[\phi\circ\pi_{\text{P}%
}](z_{0})=\frac{\text{L}\Psi}{\Psi}\left(  z_{0}\right)  .\phi\circ
\pi_{\text{P}}(z_{0})=\frac{\text{L}\Psi}{\Psi}\left(  z_{0}\right)
.\phi(y_{0})$

\textbf{COMPUTATION OF T}$_{2}:$

\qquad T$_{2}$ = L[$\phi\circ\pi_{\text{P}}$](z$_{0}$) $=$ $\frac{1}{2}%
[\Delta\phi\circ\pi_{\text{P}}](z_{0})+$ $<$X$,\nabla\phi\circ\pi_{\text{P}%
}>(z_{0})+$ V$(z_{0})\phi\circ\pi_{\text{P}}(z_{0})$

$\qquad\qquad=$ T$_{21}+$ T$_{22}+$ T$_{23}$

The formula for the Laplace-Type operator $\Delta$ acting on sections of the
vector bundle $\Gamma(E)$ is given in (iv) of \textbf{Proposition 4} by:

$\qquad\frac{1}{2}\Delta=\frac{1}{2}\left\{  g^{ij}(\nabla_{\partial_{i}%
}\nabla_{\partial_{j}}-\Gamma_{ij}^{k}(\frac{\partial}{\partial\text{x}_{k}%
}+\Lambda_{k}))+W\right\}  $

\ \ \ \ \ \ \ \qquad$=\frac{1}{2}g^{ij}\left\{  \frac{\partial^{2}}%
{\partial\text{x}_{i}\partial\text{x}_{j}}+\frac{\partial\Lambda_{j}}%
{\partial\text{x}_{i}}+\Lambda_{j}\frac{\partial}{\partial\text{x}_{i}%
}+\Lambda_{i}\frac{\partial}{\partial\text{x}_{j}}+\Lambda_{i}\Lambda
_{j}-\Gamma_{ij}^{k}(\frac{\partial}{\partial\text{x}_{k}}+\Lambda
_{k})\right\}  $ $+$ $\frac{1}{2}W$

where $\Gamma_{ij}^{k}$ are the Christoffel symbols of the Levi-Cevita
connection on M and where W is the Weitzenb\H{o}ck term. By definition, a
Levi-Cevita connection is a torsion-free connection. We emphasize here that
the Levi-Cevita connection is different from the Alfred Gray connection
defined earlier here.

We have:

$\left(  C_{11}\right)  \qquad$T$_{21}=\frac{1}{2}[\Delta\phi\circ
\pi_{\text{P}}]$(z$_{0}$)

$=\frac{1}{2}$ g$^{ij}(z_{0})\left[  \left\{  \frac{\partial^{2}}%
{\partial\text{x}_{i}\partial\text{x}_{j}}+\frac{\partial\Lambda_{j}}%
{\partial\text{x}_{i}}+\Lambda_{j}\frac{\partial}{\partial\text{x}_{i}%
}+\Lambda_{i}\frac{\partial}{\partial\text{x}_{j}}+\Lambda_{i}\Lambda
_{j}-\Gamma_{ij}^{k}(\frac{\partial}{\partial\text{x}_{k}}+\Lambda
_{k}\right\}  \phi\circ\pi_{\text{P}}\right]  (z_{0})$

$\qquad+\frac{1}{2}$(W$\phi\circ\pi_{\text{P}}$)$($ z$_{0})$

$=$ g$^{ij}(z_{0})\left\{  \frac{\partial^{2}\phi}{\partial\text{x}%
_{i}\partial\text{x}_{j}}(\pi_{\text{P}}(z_{0}))\frac{\partial}{\partial
\text{x}_{i}}\pi_{\text{P}}(z_{0})\frac{\partial}{\partial\text{x}_{j}}%
\pi_{\text{P}}(z_{0})\text{ + }\frac{\partial\phi}{\partial\text{x}_{j}}%
(\pi_{\text{P}}(z_{0}))\frac{\partial^{2}}{\partial\text{x}_{i}\partial
\text{x}_{j}}\pi_{\text{P}}(z_{0})\right\}  $

$+\frac{1}{2}$g$^{ij}(z_{0})\left\{  \frac{\partial\Lambda_{j}}{\partial
\text{x}_{i}}(z_{0})\phi(\pi_{\text{P}}(z_{0}))\text{ + }\Lambda_{j}%
(z_{0})\frac{\partial\phi}{\partial\text{x}_{i}}(\pi_{\text{P}}(z_{0}%
))\frac{\partial}{\partial\text{x}_{i}}\pi_{\text{P}}(z_{0})\text{ }\right\}
$

$+\frac{1}{2}$g$^{ij}(z_{0})\left\{  \Lambda_{i}(z_{0})\frac{\partial\phi
}{\partial\text{x}_{j}}(\pi_{\text{P}}(z_{0}))\frac{\partial}{\partial
\text{x}_{j}}\pi_{\text{P}}(z_{0})\text{ + }\Lambda_{i}(z_{0})\Lambda
_{j}(z_{0})\phi(\pi_{\text{P}}(z_{0}))\right\}  $

$-\frac{1}{2}$g$^{ij}(z_{0})\left\{  \Gamma_{ij}^{k}(z_{0})\frac{\partial
}{\partial\text{x}_{k}}\phi\text{(}\pi_{\text{P}}\text{(z}_{0}\text{))}%
\frac{\partial}{\partial\text{x}_{k}}\pi_{\text{P}}(z_{0})\text{ + }%
\Gamma_{ij}^{k}(z_{0})\Lambda_{k}(z_{0})\phi(\pi_{\text{P}}(z_{0}))\right\}  $

$+\frac{1}{2}$W(z$_{0}$)($\phi(\pi_{\text{P}}(z_{0})$)

By $\left(  C_{8}\right)  ,$ we have:

$\left(  C_{12}\right)  \qquad$T$_{21}$ \ $=$ $\frac{1}{2}\Delta\lbrack
\phi\circ\pi_{\text{P}}]$(z$_{0}$)\qquad\qquad\qquad

\qquad$=$ \ $\frac{1}{2}\underset{\text{a,b=1}}{\overset{\text{q}}{\sum}}%
$g$^{\text{ab}}(z_{0})\left\{  \frac{\partial^{2}\phi}{\partial\text{x}%
_{\text{a}}\partial\text{x}_{\text{b}}}\text{(y) }\right\}  +\frac{1}%
{2}\underset{\text{i,j=1}}{\overset{\text{n}}{\sum}}$g$^{ij}(z_{0})\left\{
\frac{\partial\Lambda_{j}}{\partial\text{x}_{i}}(z_{0})\phi\text{(y)}\right\}
$

\qquad$+\frac{1}{2}\underset{\text{j=1}}{\overset{\text{n}}{\sum}%
}\underset{\text{a=1}}{\overset{\text{q}}{\sum}}$g$^{\text{a}j}(z_{0})\left\{
\text{ }\Lambda_{j}(z_{0})\frac{\partial\phi}{\partial\text{x}_{\text{a}}%
}\text{(y)}\right\}  +\frac{1}{2}\underset{\text{i=1}}{\overset{\text{n}%
}{\sum}}\underset{\text{b=1}}{\overset{\text{q}}{\sum}}$g$^{i\text{b}}%
(z_{0})\left\{  \Lambda_{i}\text{(z}_{0}\text{)}\frac{\partial\phi}%
{\partial\text{x}_{\text{b}}}\text{(y) }\right\}  $

$\qquad+\frac{1}{2}\underset{\text{i,j=1}}{\overset{\text{n}}{\sum}}$%
g$^{ij}(z_{0})\Lambda_{i}$(z$_{0}$)$\Lambda_{j}$(z$_{0}$)$\phi$(y)

$\qquad-\frac{1}{2}\underset{\text{i,j=1}}{\overset{\text{n}}{\sum}}$%
g$^{ij}(z_{0})\left\{  \underset{\text{c=1}}{\overset{\text{q}}{\sum}}%
\Gamma_{ij}^{\text{c}}(z_{0})\frac{\partial\phi}{\partial\text{x}_{\text{c}}%
}\text{(y) + }\underset{\text{k=1}}{\overset{\text{n}}{\sum}}\Gamma_{ij}%
^{k}(z_{0})\Lambda_{k}(z_{0})\phi\text{(y)}\right\}  $

$\qquad+\frac{1}{2}$W$(z_{0})\phi$(y)

\textbf{COMPUTATION OF T}$_{3}$

$\qquad$T$_{22}$ = \ $<$X,$\nabla$[$\phi\circ\pi_{\text{P}}]>(z_{0})$ for a
general smooth vector field X.

By (i) of \textbf{Proposition 3.2},

\ $\left(  C_{13}\right)  \qquad<$X$,\nabla$[$\phi\circ\pi_{\text{P}}%
]>(z_{0})$

\qquad\qquad\ $=$ $\underset{j=1}{\overset{n}{\sum}}$X$_{j}(z_{0}%
)\nabla_{\partial_{j}}\left(  \phi\circ\pi_{\text{P}}\right)  (z_{0})=$
$\underset{\text{j=1}}{\overset{\text{n}}{\sum}}$X$_{j}(z_{0})\left(
\frac{\partial}{\partial x_{j}}+\Lambda_{j}\right)  (z_{0})\phi\circ
\pi_{\text{P}}](z_{0})$

\qquad\qquad$=$ $\underset{j=1}{\overset{n}{\sum}}$X$_{j}(z_{0})\left\{
\frac{\partial}{\partial\text{x}_{j}}\phi\circ\pi_{\text{P}}(z_{0})+\text{
}\Lambda_{j}(z_{0})\phi\circ\pi_{\text{P}}(z_{0})\right\}  $

\qquad\qquad$=$ $\underset{\text{a=1}}{\overset{\text{q}}{\sum}}$X$_{\text{a}%
}(z_{0})\frac{\partial\phi}{\partial\text{x}_{\text{a}}}(y)$ $+$
$\underset{j=1}{\overset{n}{\sum}}$ X$_{j}(z_{0})\Lambda_{j}(z_{0})\phi(y)$

The last equality above is due to (i) of $\left(  C_{8}\right)  .$

$\left(  C_{14}\right)  \qquad$T$_{23}=$V$(z_{0})\phi\circ\pi_{\text{P}}%
(z_{0})=$ V$(z_{0})\phi(y)$

The last step here is to compute:

\qquad\qquad\ \ T$_{3}=$ \
$<$%
$\nabla\log\Psi,\nabla\lbrack\phi\circ\pi_{\text{P}}]$%
$>$%
(z$_{0})$

We use $\left(  C_{13}\right)  $ where we take X $=$ $\nabla\log\Psi$ and have:

$\left(  C_{15}\right)  \qquad$T$_{3}$ $=$\ \
$<$%
$\nabla\log\Psi,\nabla\lbrack\phi\circ\pi_{\text{P}}]$%
$>$%
(z$_{0})$

\qquad\qquad$=$ $\underset{\text{a=1}}{\overset{\text{q}}{\sum}}(\nabla
\log\Psi)_{\text{a}}(z_{0})\frac{\partial\phi}{\partial\text{x}_{\text{a}}%
}(y)$ + $\underset{j=1}{\overset{n}{\sum}}$ $(\nabla\log\Psi)_{j}%
(z_{0})\Lambda_{j}(z_{0})\phi(y)$

We set:

$\qquad\qquad\Theta(z_{0})=$ L$_{\Psi}[\phi\circ\pi_{\text{P}}](z_{0})$

and recall that:

$\qquad\Theta(z_{0})=$ L$_{\Psi}$[$\phi\circ\pi_{\text{P}}$]$(z_{0})$ $=$
T$_{1}$ + T$_{2}$ + T$_{3}+$ T$_{4}$

where, we recall that: T$_{4}=$ $-$ V$(z_{0})\phi\circ\pi_{P}(z_{0})=$ $-$
V$(z_{0})\phi(y)$

We use $\left(  C_{10}\right)  ,$ $\left(  C_{12}\right)  ,\left(
C_{13}\right)  ,\left(  C_{14}\right)  $, $\left(  C_{15}\right)  $ and
$\left(  C_{16}\right)  $ to collect all the terms of L$_{\Psi}[\phi\circ
\pi_{\text{P}}](z_{0})$ in $\left(  C_{9}\right)  $ and have:

$\left(  C_{16}\right)  \qquad\Theta(z_{0})=$ L$_{\Psi}[\phi\circ\pi
_{\text{P}}](z_{0})$

\qquad\qquad\textbf{ }$=\frac{\text{L}\Psi}{\Psi}\left(  z_{0}\right)
\phi(y)$\qquad\qquad

+ \ $\frac{1}{2}\underset{\text{a,b=1}}{\overset{\text{q}}{\sum}}%
$g$^{\text{ab}}(z_{0})\left\{  \frac{\partial^{2}\phi}{\partial\text{x}%
_{\text{a}}\partial\text{x}_{\text{b}}}(y)\text{ }\right\}  +\frac{1}%
{2}\underset{i,j=1}{\overset{n}{\sum}}$g$^{ij}(z_{0})\left\{  \frac
{\partial\Lambda_{j}}{\partial\text{x}_{i}}(z_{0})\phi(y)\right\}  $

$+\frac{1}{2}\underset{j=1}{\overset{n}{\sum}}\underset{\text{a=1}%
}{\overset{\text{q}}{\sum}}$g$^{\text{a}j}(z_{0})\left\{  \text{ }\Lambda
_{j}(z_{0})\frac{\partial\phi}{\partial\text{x}_{\text{a}}}(y)\right\}
+\frac{1}{2}\underset{i=1}{\overset{n}{\sum}}\underset{\text{b=1}%
}{\overset{\text{q}}{\sum}}$g$^{i\text{b}}(z_{0})\left\{  \Lambda_{i}%
\text{(z}_{0}\text{)}\frac{\partial\phi}{\partial\text{x}_{\text{b}}}(y)\text{
}\right\}  $

$+\frac{1}{2}\underset{i,j=1}{\overset{n}{\sum}}$g$^{ij}(z_{0})\Lambda_{i}%
$(z$_{0}$)$\Lambda_{j}$(z$_{0}$)$\phi(y)$

$-\frac{1}{2}\underset{i,j=1}{\overset{n}{\sum}}$g$^{ij}(z_{0})\left\{
\Gamma_{ij}^{\text{c}}(z_{0})\frac{\partial\phi}{\partial\text{x}_{\text{c}}%
}(y)\text{ + }\Gamma_{ij}^{k}(z_{0})\Lambda_{k}(z_{0})\phi(y)\right\}
\ +\frac{1}{2}$W$(z_{0})\phi$(y)

$\ +$ $\underset{\text{a=1}}{\overset{\text{q}}{\sum}}(\nabla^{0}\log
\Psi)_{\text{a}}(z_{0})\frac{\partial\phi}{\partial\text{x}_{\text{a}}}(y)$ +
$\underset{j=1}{\overset{n}{\sum}}$ $(\nabla^{0}\log\Psi)_{j}(z_{0}%
)\Lambda_{j}(z_{0})\phi(y)$

$+$ $\underset{\text{a=1}}{\overset{\text{q}}{\sum}}$X$_{\text{a}}(z_{0}%
)\frac{\partial\phi}{\partial\text{x}_{\text{a}}}(y)+$ $\underset{\text{j=1}%
}{\overset{\text{n}}{\sum}}$ X$_{j}(z_{0})\Lambda_{j}(z_{0})\phi(y)+$
V$(z_{0})\phi(y)-$ V(z$_{0}$)$\phi(y)$

The two last items of \textbf{potential terms} cancel out. We set:
$\Theta(z_{0})=$ L$_{\Psi}[\phi\circ\pi_{\text{P}}](z_{0})$ and we have from
the above decomposition of $\Theta(z_{0})=$ L$_{\Psi}[\phi\circ\pi_{\text{P}%
}](z_{0}):$

$\left(  C_{17}\right)  \qquad$b$_{1}$(x,P,$\phi)=\int_{0}^{1}$F(1,1-r$_{1}%
$)$[$L$_{\Psi}[\phi\circ\pi_{\text{P}}](x)$dr$_{1}$

$=\int_{0}^{1}$L$_{\Psi}[\phi\circ\pi_{\text{P}}](z_{0})$dr $=\int_{0}%
^{1}\Theta$(z$_{0}$)dr$_{1}$\qquad\qquad\qquad\qquad\qquad$\ \ \ $

$=\int_{0}^{1}[\frac{\text{L}\Psi}{\Psi}\left(  z_{0}\right)  .\phi(y)$%
\qquad\qquad

$+$ \ $\frac{1}{2}\underset{\text{a,b=1}}{\overset{\text{q}}{\sum}}%
$g$^{\text{ab}}(z_{0})\left\{  \frac{\partial^{2}\phi}{\partial\text{x}%
_{\text{a}}\partial\text{x}_{\text{b}}}(y)\text{ }\right\}  +\frac{1}%
{2}\underset{\text{i,j=1}}{\overset{\text{n}}{\sum}}$g$^{ij}(z_{0})\left\{
\frac{\partial\Lambda_{j}}{\partial\text{x}_{i}}(z_{0})\phi(y)\right\}  $

$\bigskip+\frac{1}{2}\underset{\text{j=1}}{\overset{\text{n}}{\sum}%
}\underset{\text{a=1}}{\overset{\text{q}}{\sum}}$g$^{\text{a}j}(z_{0})\left\{
\text{ }\Lambda_{j}(z_{0})\frac{\partial\phi}{\partial\text{x}_{\text{a}}%
}(y)\right\}  +\frac{1}{2}\underset{\text{i=1}}{\overset{\text{n}}{\sum}%
}\underset{\text{b=1}}{\overset{\text{q}}{\sum}}$g$^{i\text{b}}(z_{0})\left\{
\Lambda_{i}\text{(z}_{0}\text{)}\frac{\partial\phi}{\partial\text{x}%
_{\text{b}}}(y)\text{ }\right\}  $

$-\frac{1}{2}\underset{\text{i,j=1}}{\overset{\text{n}}{\sum}}$g$^{ij}%
(z_{0})\left\{  \Gamma_{ij}^{\text{c}}(z_{0})\frac{\partial\phi}%
{\partial\text{x}_{\text{c}}}(y)\text{ + }\Gamma_{ij}^{k}(z_{0})\Lambda
_{k}(z_{0})\phi(y)\right\}  \ +\frac{1}{2}$W$(z_{0})\phi$(y)

+ $\underset{\text{a=1}}{\overset{\text{q}}{\sum}}(\nabla^{0}\log
\Psi)_{\text{a}}(z_{0})\frac{\partial\phi}{\partial\text{x}_{\text{a}}}(y)$ +
$\underset{\text{j=1}}{\overset{\text{n}}{\sum}}$ $(\nabla^{0}\log\Psi
)_{j}(z_{0})\Lambda_{j}(z_{0})\phi(y)$

\ + $\underset{\text{a=1}}{\overset{\text{q}}{\sum}}$X$_{\text{a}}(z_{0}%
)\frac{\partial\phi}{\partial\text{x}_{\text{a}}}(y)$ + $\underset{\text{j=1}%
}{\overset{\text{n}}{\sum}}$ X$_{j}(z_{0})\Lambda_{j}(z_{0})\phi(y)\phi
(y)]$dr$_{1}$

\qquad\qquad\qquad\qquad\qquad\qquad\qquad\qquad\qquad\qquad\qquad\qquad
\qquad\qquad\qquad\qquad\qquad$\blacksquare$

The integrand L$_{\Psi}$[$\phi\circ\pi_{\text{P}}$]$(z_{0})$ is not
independent of r$_{1}$ since z$_{0}=\gamma(r_{1})$ where $\gamma$ is the
unique minimal geodesic from x to y in time 1 and r$_{1}\in\left[  0,1\right]
.$

The above is the "raw expression" of b$_{1}$(x,P,$\phi).$ It will be made more
explicit in \textbf{Chaptet 7} by expressing it in geometric terms at the
centre of Fermi coordinates $y_{0}$ in \textbf{Theorem 7.1.}

\chapter{Computation of the Third Coefficient}

The "raw" expression for the third coefficient below is taken from $\left(
11.43\right)  $\textbf{ }in \textbf{Chapter 11. }

There were cancellations in the expression of b$_{2}($y$_{0}$,P$,\phi)$. In
these cancellations, \textbf{I}$_{2}$\textbf{ and I}$_{4}$\textbf{ were
entirely wiped off} and the final expression given in $\left(  11.43\right)  $
is follows:

$\left(  D_{1}\right)  \qquad$b$_{2}($y$_{0}$,P$,\phi)=$ b$_{2}($y$_{0}%
$,P$)\phi(y_{0})=$ I$_{1}+$ I$_{3}$

$=\frac{1}{2}\frac{\text{L}\Psi}{\Psi}(y_{0})\Theta(y_{0})\phi(y_{0})\qquad
$I$_{1}$

$+\frac{1}{4}\underset{\text{a=1}}{\overset{\text{q}}{\sum}}\frac{\partial
^{2}\Theta}{\partial x_{\text{a}}^{2}}(y_{0})\phi(y_{0})+\frac{1}%
{12}\underset{i=q+1}{\overset{n}{\sum}}\frac{\partial^{2}\Theta}{\partial
x_{i}^{2}}(y_{0})\phi(y_{0})+\frac{1}{2}\underset{\text{a=1}%
}{\overset{\text{q}}{\sum}}\Lambda_{\text{a}}(y_{0})\frac{\partial\Theta
}{\partial x_{\text{a}}}(y_{0})\phi(y_{0})\qquad$I$_{3}$ starts

$+\frac{1}{4}\underset{i=q+1}{\overset{n}{\sum}}\Lambda_{i}(y_{0}%
)\frac{\partial\Theta}{\partial x_{i}}(y_{0})\phi(y_{0})+\frac{1}%
{4}\underset{\text{a=1}}{\overset{\text{q}}{\sum}}\Lambda_{\text{a}}^{2}%
(y_{0})\Theta(y_{0})\phi(y_{0})$

$-\frac{1}{4}\underset{j=q+1}{\overset{n}{\sum}}<H,j>(y_{0})\Lambda_{j}%
(y_{0})\Theta(y_{0})\phi(y_{0})+\frac{1}{4}\Theta(y_{0})$W$(y_{0})\phi(y_{0})$

$+\frac{1}{2}$ $\underset{\text{a=1}}{\overset{\text{q}}{\sum}}$X$_{\text{a}%
}(y_{0})\frac{\partial\Theta}{\partial x_{\text{a}}}(y_{0})\phi(y_{0}%
)+\frac{1}{2}$ $\underset{j=1}{\overset{n}{\sum}}$X$_{j}(y_{0})\Lambda
_{j}(y_{0})\Theta(y_{0})\phi(y_{0})\qquad\qquad$I$_{3}$ ends

\qquad\qquad\qquad\qquad\qquad\qquad\qquad\qquad\qquad\qquad\qquad\qquad
\qquad\qquad\qquad$\qquad\qquad\blacksquare$

\qquad The numbering below has been slightly different from that of $\left(
11.44\right)  :$\qquad

\qquad b$_{2}($y$_{0}$,P$)\phi(y_{0})=$ I$_{1}+$ I$_{3}=$ I$_{31}+$ I$_{32}+$
I$_{33}+$I$_{34}+$ I$_{35}+$ I$_{36}+$ I$_{37}+$ I$_{38}+$ I$_{39}$

\qquad I$_{1}=\frac{1}{2}\frac{\text{L}\Psi}{\Psi}(y_{0})\Theta(y_{0}%
)\phi(y_{0})$

\qquad I$_{31}=\frac{1}{12}\underset{\text{a=1}}{\overset{\text{q}}{\sum}%
}\frac{\partial^{2}\Theta}{\partial x_{\text{a}}^{2}}(y_{0})\phi(y_{0});$
\qquad I$_{32}=\frac{1}{12}\underset{i=q+1}{\overset{n}{\sum}}\frac
{\partial^{2}\Theta}{\partial x_{i}^{2}}(y_{0})\phi(y_{0})$

$\qquad$I$_{33}=\frac{1}{2}\underset{\text{a=1}}{\overset{\text{q}}{\sum}%
}\Lambda_{\text{a}}(y_{0})\frac{\partial\Theta}{\partial x_{\text{a}}}%
(y_{0})\phi(y_{0});\qquad$I$_{34}=$ $\frac{1}{2}$ $\underset{\text{a=1}%
}{\overset{\text{q}}{\sum}}X_{\text{a}}(y_{0})\frac{\partial\Theta}{\partial
x_{\text{a}}}(y_{0})\phi(y_{0})$

\qquad I$_{35}=\frac{1}{4}\underset{i=q+1}{\overset{n}{\sum}}\Lambda_{i}%
(y_{0})\frac{\partial\Theta}{\partial x_{i}}(y_{0})\phi(y_{0});\qquad$%
I$_{36}=\frac{1}{4}\underset{\text{a=1}}{\overset{\text{q}}{\sum}}%
\Lambda_{\text{a}}^{2}(y_{0})\Theta(y_{0})\phi(y_{0});\qquad$

\qquad I$_{37}=-\frac{1}{4}\underset{j=q+1}{\overset{n}{\sum}}<H,j>(y_{0}%
)\Lambda_{j}(y_{0})\Theta(y_{0})\phi(y_{0});\qquad\qquad$

\qquad I$_{38}=\frac{1}{2}$ $\underset{j=1}{\overset{n}{\sum}}$X$_{j}%
(y_{0})\Lambda_{j}(y_{0})\Theta(y_{0})\phi(y_{0});\qquad$I$_{39}=$ $\frac
{1}{4}$W$(y_{0})\Theta(y_{0})\phi(y_{0})$

We now come to one of the most important theorems of this work. We make
computations of each of the above items in \textbf{geometric invariants} of
the Riemannian manifold M, the submanifold P and the vector bundle E:

We recall that $\Theta(y_{0})=$ L$_{\Psi}[\phi\circ\pi_{P}](y_{0})=$
\textbf{b}$_{1}$(y$_{0}$,P)$\phi\left(  y_{0}\right)  $ by $\left(
10.31\right)  $ in \textbf{Chapter 10.}

It is given in geometric invariants as follows:

$\left(  D_{2}\right)  \qquad\Theta(y_{0})=$ \textbf{b}$_{1}$(y$_{0}$%
,P)$\phi\left(  y_{0}\right)  $

$\qquad\qquad=\frac{1}{24}[\underset{i=q+1}{\overset{n}{\sum}}3<H,i>^{2}%
+2(\tau^{M}-3\tau^{P}\ +\overset{q}{\underset{\text{a=1}}{\sum}}%
\varrho_{\text{aa}}^{M}+\overset{q}{\underset{\text{a,b}=1}{\sum}%
}R_{\text{abab}}^{M})](y_{0})\phi\left(  y_{0}\right)  $

\qquad\qquad$\ -\frac{1}{2}[$ $\left\Vert \text{X}\right\Vert _{M}^{2}+$
$\operatorname{div}X_{M}-$ $\left\Vert \text{X}\right\Vert _{P}^{2}$
$-\operatorname{div}X_{P}](y_{0})\phi\left(  y_{0}\right)  +$ V$(y_{0}%
)\phi\left(  y_{0}\right)  $

$\qquad\qquad+$ $\frac{1}{2}\underset{\text{a=1}}{\overset{\text{q}}{\sum}%
}\frac{\partial^{2}\phi}{\partial\text{x}_{\text{a}}^{2}}(y_{0})$ $+$
$\underset{\text{a=1}}{\overset{\text{q}}{\sum}}\Lambda_{\text{a}}(y_{0}%
)\frac{\partial\phi}{\partial x_{\text{a}}}\left(  y_{0}\right)  \ +\frac
{1}{2}$ $\underset{\text{a=1}}{\overset{\text{q}}{\sum}}\Lambda_{\text{a}%
}(y_{0})\Lambda_{\text{a}}(y_{0})\phi\left(  y_{0}\right)  $

\qquad\qquad$+$ $\underset{\text{a=1}}{\overset{\text{q}}{\sum}}X_{\text{a}%
}(y_{0})\frac{\partial\phi}{\partial\text{x}_{\text{a}}}(y_{0})+$
$\underset{\text{a=1}}{\overset{\text{q}}{\sum}}$ $X_{\text{a}}(y_{0}%
)\Lambda_{\text{a}}(y_{0})\phi(y_{0})+\frac{1}{2}W\left(  y_{0}\right)
\phi\left(  y_{0}\right)  $

On the other hand, by $\left(  10.30\right)  $ of \textbf{Chapter 10,}\qquad

$\left(  D_{3}\right)  $\qquad$\frac{\text{L}\Psi}{\Psi}(y_{0})\phi\left(
y_{0}\right)  $

$\qquad=\frac{1}{24}[\underset{i=q+1}{\overset{n}{\sum}}3<H,i>^{2}+2(\tau
^{M}-3\tau^{P}\ +\overset{q}{\underset{\text{a=1}}{\sum}}\varrho_{\text{aa}%
}^{M}+\overset{q}{\underset{\text{a,b}=1}{\sum}}R_{\text{abab}}^{M}%
)](y_{0})\phi\left(  y_{0}\right)  $

\qquad$-\frac{1}{2}[$ $\left\Vert \text{X}\right\Vert _{M}^{2}+$
$\operatorname{div}X_{M}-\left\Vert \text{X}\right\Vert _{P}^{2}$ $-$
$\operatorname{div}X_{P}](y_{0})\phi\left(  y_{0}\right)  +$ V$(y_{0}%
)\phi\left(  y_{0}\right)  $

\qquad\qquad\qquad\qquad\qquad\qquad\qquad\qquad\qquad\qquad\qquad\qquad
\qquad\qquad\qquad\qquad\qquad$\blacksquare$

Therefore by the expressions of $\Theta(y_{0})\phi\left(  y_{0}\right)  $ and
$\frac{\text{L}\Psi}{\Psi}(y_{0})\phi\left(  y_{0}\right)  $ in $\left(
D_{2}\right)  $ and $\left(  D_{3}\right)  $ above$,$ we have:

$\left(  D_{4}\right)  \qquad$I$_{1}=\frac{1}{2}\frac{\text{L}\Psi}{\Psi
}(y_{0})\Theta(y_{0})\phi\left(  y_{0}\right)  $

\qquad\qquad$=\frac{1}{2}[\frac{1}{24}(\underset{i=q+1}{\overset{n}{\sum}%
}3<H,i>^{2}+2(\tau^{M}-3\tau^{P}\ +\overset{q}{\underset{\text{a=1}}{\sum}%
}\varrho_{\text{aa}}^{M}+\overset{q}{\underset{\text{a,b}=1}{\sum}%
}R_{\text{abab}}^{M}))$

\qquad\qquad$-\frac{1}{2}($ $\left\Vert \text{X}\right\Vert _{M}^{2}+$
$\operatorname{div}X_{M}-\left\Vert \text{X}\right\Vert _{P}^{2}$ $-$
$\operatorname{div}X_{P})+$ $V](y_{0})$

$\qquad\qquad\times\lbrack\frac{1}{24}(\underset{i=q+1}{\overset{n}{\sum}%
}3<H,i>^{2}+2(\tau^{M}-3\tau^{P}\ +\overset{q}{\underset{\text{a=1}}{\sum}%
}\varrho_{\text{aa}}^{M}+\overset{q}{\underset{\text{a,b}=1}{\sum}%
}R_{\text{abab}}^{M})$

\qquad$\qquad-\frac{1}{2}($ $\left\Vert \text{X}\right\Vert _{M}^{2}+$
$\operatorname{div}X_{M}-\left\Vert \text{X}\right\Vert _{P}^{2}$
$-\operatorname{div}X_{P})+$ $V$

$\qquad\qquad+$ $(\frac{1}{2}\underset{\text{a=1}}{\overset{\text{q}}{\sum}%
}\frac{\partial^{2}\phi}{\partial\text{x}_{\text{a}}^{2}}$ $+$
$\underset{\text{a=1}}{\overset{\text{q}}{\sum}}\Lambda_{\text{a}}(y_{0}%
)\frac{\partial\phi}{\partial x_{\text{a}}}\ +\frac{1}{2}$
$\underset{\text{a=1}}{\overset{\text{q}}{\sum}}\Lambda_{\text{a}}%
\Lambda_{\text{a}})$

\qquad$\qquad+($ $\underset{\text{a=1}}{\overset{\text{q}}{\sum}}X_{\text{a}%
}(y_{0})\frac{\partial\phi}{\partial\text{x}_{\text{a}}}+$
$\underset{\text{a=1}}{\overset{\text{q}}{\sum}}$ $X_{\text{a}}\Lambda
_{\text{a}}+\frac{1}{2}W)]\phi\left(  y_{0}\right)  $\qquad

\textbf{Computation of I}$_{31}\qquad$

$\left(  D_{5}\right)  \qquad$I$_{31}=\frac{1}{4}\underset{\text{c=1}%
}{\overset{\text{q}}{\sum}}\frac{\partial^{2}\Theta}{\partial x_{\text{c}}%
^{2}}(y_{0})=\frac{1}{4}\underset{\text{c=1}}{\overset{\text{q}}{\sum}}%
\frac{\partial^{2}}{\partial x_{\text{c}}^{2}}[L_{\Psi}\phi\circ\pi_{P}%
](y_{0})$

We have from $\left(  10.31\right)  :$

$\left(  D_{6}\right)  \qquad\Theta=$ L$_{\Psi}[\phi\circ\pi_{P}]$

\qquad\qquad\ \ \ \textbf{ }$=\frac{\text{L}\Psi}{\Psi}\left(  z_{1}\right)
\phi\circ\pi_{\text{P}}$\qquad\qquad

\qquad$+$ \ $\frac{1}{2}\underset{\text{a,b=1}}{\overset{\text{q}}{\sum}}%
$g$^{\text{ab}}\left\{  \frac{\partial^{2}\phi}{\partial\text{x}_{\text{a}%
}\partial\text{x}_{\text{b}}}\circ\pi_{\text{P}}\right\}  +\frac{1}%
{2}\underset{\text{i,j=1}}{\overset{\text{n}}{\sum}}$g$^{ij}\left\{
\frac{\partial\Lambda_{j}}{\partial\text{x}_{i}}\phi\circ\pi_{\text{P}%
}\right\}  $

$\qquad+$ $\underset{j=1}{\overset{n}{\sum}}\underset{\text{a=1}%
}{\overset{\text{q}}{\sum}}$g$^{\text{a}j}\left\{  \text{ }\Lambda_{j}%
\frac{\partial\phi}{\partial\text{x}_{\text{a}}}\circ\pi_{\text{P}}\right\}
+\frac{1}{2}\underset{\text{i,j=1}}{\overset{\text{n}}{\sum}}$g$^{ij}%
\Lambda_{i}\Lambda_{j}\phi\circ\pi_{\text{P}}$

$\qquad-\frac{1}{2}\underset{i,j=1}{\overset{n}{\sum}}$g$^{ij}\left\{
\underset{\text{c=1}}{\overset{\text{q}}{\sum}}\Gamma_{ij}^{\text{c}}%
\frac{\partial\phi}{\partial\text{x}_{\text{c}}}\circ\pi_{\text{P}}\text{ +
}\underset{\text{k=1}}{\overset{\text{n}}{\sum}}\Gamma_{ij}^{k}(z_{1}%
)\Lambda_{k}\phi\circ\pi_{\text{P}}\right\}  $

$\qquad+\frac{1}{2}$W$\phi\circ\pi_{\text{P}}$

\qquad$+$ $\underset{\text{a=1}}{\overset{\text{q}}{\sum}}(\nabla\log
\Psi)_{\text{a}}\frac{\partial\phi}{\partial\text{x}_{\text{a}}}\circ
\pi_{\text{P}}$ + $\underset{\text{j=1}}{\overset{\text{n}}{\sum}}$
$(\nabla\log\Psi)_{j}\Lambda_{j}\phi\circ\pi_{\text{P}}$

\qquad$+$ $\underset{\text{a=1}}{\overset{\text{q}}{\sum}}$X$^{\text{a}}%
\frac{\partial\phi}{\partial\text{x}_{\text{a}}}\circ\pi_{\text{P}}$ +
$\underset{\text{j=1}}{\overset{\text{n}}{\sum}}$ X$^{j}\Lambda_{j}\phi
\circ\pi_{\text{P}}$

Then,

$\left(  D_{7}\right)  \qquad$I$_{31}=\frac{1}{4}\underset{\text{c=1}%
}{\overset{\text{q}}{\sum}}\frac{\partial^{2}}{\partial x_{\text{c}}^{2}}%
$[L$_{\Psi}\phi\circ\pi_{P}$]$(y_{0})$

$\qquad=$ I$_{311}+$ I$_{312}$+ I$_{313}+$ I$_{314}+$ I$_{315}+$ I$_{316}+$
I$_{317}+$ I$_{318}+$ I$_{319}+$ L$_{1}+$ L$_{2},$

where,

(i)\qquad I$_{311}=\frac{1}{4}\frac{\partial^{2}}{\partial x_{\text{c}}^{2}%
}[\frac{\text{L}\Psi}{\Psi}\phi\circ\pi_{P}](y_{0})$

(ii)\qquad I$_{312}=\frac{1}{8}\frac{\partial^{2}}{\partial x_{\text{c}}^{2}%
}[\underset{\text{a,b=1}}{\overset{\text{q}}{\sum}}$g$^{\text{ab}}\left\{
\frac{\partial^{2}\phi}{\partial\text{x}_{\text{a}}\partial\text{x}_{\text{b}%
}}\circ\pi_{P}\text{ }\right\}  ](y_{0})$

(iii)\qquad I$_{313}=\frac{1}{8}\frac{\partial^{2}}{\partial x_{\text{c}}^{2}%
}[\underset{\text{i,j=1}}{\overset{\text{n}}{\sum}}g^{ij}\left\{
\frac{\partial\Lambda_{j}}{\partial\text{x}_{i}}\phi\circ\pi_{P}\right\}
](y_{0})\qquad$

(iv)\qquad I$_{314}$ $=\frac{1}{4}\frac{\partial^{2}}{\partial x_{\text{c}%
}^{2}}\underset{\text{j=1}}{\overset{\text{n}}{\sum}}\underset{\text{a=1}%
}{\overset{\text{q}}{\sum}}$g$^{\text{a}j}\left\{  \text{ }\Lambda_{j}%
\frac{\partial\phi}{\partial\text{x}_{\text{a}}}\circ\pi_{P}\right\}  (y_{0})$

(v)\qquad I$_{315}=$ $\frac{1}{8}\frac{\partial^{2}}{\partial x_{\text{c}}%
^{2}}[\underset{\text{i,j=1}}{\overset{\text{n}}{\sum}}$g$^{ij}\left\{
\Lambda_{i}\Lambda_{j}\phi\circ\pi_{P}\right\}  ](y_{0})$

(vi)\qquad I$_{316}=-\frac{1}{8}\frac{\partial^{2}}{\partial x_{\text{c}}^{2}%
}[\underset{\text{i,j=1}}{\overset{\text{n}}{\sum}}g^{ij}\left\{  \Gamma
_{ij}^{\text{c}}\frac{\partial\phi}{\partial\text{x}_{\text{c}}}\circ\pi
_{P}\text{ + }\Gamma_{ij}^{k}\Lambda_{k}\phi\circ\pi_{P}\right\}  ](y_{0})$

(vii)\qquad I$_{317}=\frac{1}{4}\frac{\partial^{2}}{\partial x_{\text{c}}^{2}%
}[$ $\underset{\text{a=1}}{\overset{\text{q}}{\sum}}(\nabla\log\Psi
)_{\text{a}}\frac{\partial\phi}{\partial\text{x}_{\text{a}}}\circ\pi
_{\text{P}}](y_{0})$

(viii)\qquad I$_{318}=\frac{1}{4}\frac{\partial^{2}}{\partial x_{\text{a}}%
^{2}}[\underset{\text{j=1}}{\overset{\text{n}}{\sum}}(\nabla\log\Psi
)_{j}\Lambda_{j}\phi\circ\pi_{P}](y_{0})$

(ix)\qquad I$_{319}=\frac{1}{8}\frac{\partial^{2}}{\partial x_{\text{a}}^{2}%
}[$W$\phi\circ\pi_{P}](y_{0})$

(x)\qquad L$_{1}=\frac{1}{4}\frac{\partial^{2}}{\partial x_{\text{a}}^{2}}[$
$\underset{\text{a=1}}{\overset{\text{q}}{\sum}}$X$_{\text{a}}\frac
{\partial\phi}{\partial\text{x}_{\text{a}}}\circ\pi](y_{0})$

(xi)\qquad L$_{2}=\frac{1}{4}\frac{\partial^{2}}{\partial x_{\text{a}}^{2}}[$
$\underset{\text{j=1}}{\overset{\text{n}}{\sum}}$ X$_{j}\Lambda_{j}\phi
\circ\pi_{\text{P}}](y_{0})$

\underline{\textbf{Computation of the above in Geometric Invariants:}}

In all the computations that follow we will often use the following well known
elementary equations:

$\qquad\frac{\partial^{2}}{\partial x^{2}}(fg)=$ $\frac{\partial^{2}%
f}{\partial x^{2}}g+f\frac{\partial^{2}g}{\partial x^{2}}+$ $2\frac{\partial
f}{\partial x}\frac{\partial g}{\partial x}$

\qquad$\nabla(fg)=(\nabla f)g+f(\nabla g)$

\qquad$\frac{1}{2}\Delta$(fg) =$\frac{1}{2}$($\Delta$f)g + $\frac{1}{2}%
$f($\Delta$g) +
$<$%
$\nabla$f,$\nabla$g%
$>$%

\qquad L(f$\phi$) = (Lf)$\phi$ + f(L$\phi$) +
$<$%
$\nabla$f,$\nabla\phi$%
$>$
$-$ V$($f$\phi)$

We will also frequently use $\left(  C_{8}\right)  $ without explicitely
referring to it:

(i)\qquad I$_{311}=\frac{1}{4}\frac{\partial^{2}}{\partial x_{\text{c}}^{2}%
}[\frac{\text{L}\Psi}{\Psi}\phi\circ\pi_{P}](y_{0})$

\qquad$=\frac{1}{4}\frac{\partial^{2}}{\partial x_{\text{c}}^{2}}%
[\frac{\text{L}\Psi}{\Psi}](y_{0})[\phi\circ\pi_{P}](y_{0})+\frac{1}{4}%
\frac{\text{L}\Psi}{\Psi}(y_{0})\frac{\partial^{2}}{\partial x_{\text{c}}^{2}%
}[\phi\circ\pi_{P}](y_{0})$

$\qquad+\frac{1}{2}\frac{\partial}{\partial x_{\text{c}}}[\frac{\text{L}\Psi
}{\Psi}](y_{0})\frac{\partial}{\partial x_{\text{c}}}[\phi\circ\pi_{P}%
](y_{0})$

\qquad$=\frac{1}{4}\frac{\partial^{2}}{\partial x_{\text{c}}^{2}}%
[\frac{\text{L}\Psi}{\Psi}](y_{0})\phi(y_{0})+\frac{1}{4}\frac{\text{L}\Psi
}{\Psi}(y_{0})\frac{\partial^{2}\phi}{\partial x_{\text{c}}^{2}}(y_{0}%
)+\frac{1}{2}\frac{\partial}{\partial x_{\text{c}}}[\frac{\text{L}\Psi}{\Psi
}](y_{0})\frac{\partial\phi}{\partial x_{\text{c}}}(y_{0})$

We write this in the order:

\qquad I$_{311}=\frac{1}{4}\frac{\text{L}\Psi}{\Psi}(y_{0})\frac{\partial
^{2}\phi}{\partial x_{\text{c}}^{2}}(y_{0})+\frac{1}{2}\frac{\partial
}{\partial x_{\text{c}}}[\frac{\text{L}\Psi}{\Psi}](y_{0})\frac{\partial\phi
}{\partial x_{\text{c}}}(y_{0})+\frac{1}{4}\frac{\partial^{2}}{\partial
x_{\text{c}}^{2}}[\frac{\text{L}\Psi}{\Psi}](y_{0})\phi(y_{0})$

Now, $\frac{\text{L}\Psi}{\Psi}(y_{0})$ is given by $\left(  D_{3}\right)  $
above, $\frac{\partial}{\partial x_{\text{c}}}[\frac{\text{L}\Psi}{\Psi
}](y_{0})$ is give by (v) of \textbf{Table B}$_{5}$ or from Equation $\left(
B_{111}\right)  $

and $\frac{\partial^{2}}{\partial x_{\text{c}}^{2}}[\frac{\text{L}\Psi}{\Psi
}](y_{0})$ is given by (vi) of \textbf{Table B}$_{5}$ or from $\left(
B_{112}\right)  .$ We have:

$\left(  D_{7}\right)  ^{\ast}$ I$_{311}=\frac{1}{96}%
[\underset{i=q+1}{\overset{n}{\sum}}3<H,i>^{2}+2(\tau^{M}-3\tau^{P}%
\ +\overset{q}{\underset{\text{a=1}}{\sum}}\varrho_{\text{aa}}^{M}%
+\overset{q}{\underset{\text{a,b}=1}{\sum}}R_{\text{abab}}^{M})](y_{0}%
)\frac{\partial^{2}\phi}{\partial x_{\text{c}}^{2}}(y_{0})$

\qquad\qquad$-\frac{1}{8}[\left\Vert X\right\Vert _{M}^{2}+\operatorname{div}%
X_{M}-\left\Vert \text{X}\right\Vert _{P}^{2}-\operatorname{div}X_{P}%
]\frac{\partial^{2}\phi}{\partial x_{\text{c}}^{2}}(y_{0})+$ $\frac{1}%
{4}V(y_{0})\frac{\partial^{2}\phi}{\partial x_{\text{c}}^{2}}(y_{0})$

\qquad\qquad$+\frac{1}{2}\left[  -\frac{1}{2}\left(  \frac{\partial^{2}X_{j}%
}{\partial x_{\text{c}}\partial x_{j}}\right)  +\frac{1}{2}<H,i>\frac{\partial
X_{i}}{\partial x_{\text{c}}}+\frac{\partial\text{V}}{\partial x_{\text{c}}%
}\right]  (y_{0})\frac{\partial\phi}{\partial x_{\text{c}}}(y_{0})$

\qquad$\qquad+\frac{1}{4}[(\frac{\partial X_{j}}{\partial x_{\text{c}}}%
)^{2}+X_{j}\frac{\partial^{2}X_{j}}{\partial x_{\text{c}}^{2}}](y_{0}%
)\phi(y_{0})$ $-\frac{1}{8}\frac{\partial^{3}X_{j}}{\partial x_{\text{c}}%
^{2}\partial x_{j}}(y_{0})\phi(y_{0})$

$\qquad\qquad+\frac{1}{8}<H,i>(y_{0})\frac{\partial^{2}X_{i}}{\partial
x_{\text{c}}^{2}}(y_{0})\phi(y_{0})\qquad$

\qquad\qquad$-\frac{1}{2}[(\frac{\partial X_{i}}{\partial x_{\text{c}}}%
)^{2}+X_{i}\frac{\partial^{2}X_{i}}{\partial x_{\text{c}}^{2}})](y_{0}%
)\phi(y_{0})+$ $\frac{1}{4}\frac{\partial^{2}\text{V}}{\partial x_{\text{c}%
}^{2}}(y_{0})\phi(y_{0}).$

\qquad\qquad\qquad\qquad\qquad\qquad\qquad\qquad\qquad\qquad\qquad\qquad
\qquad\qquad\qquad\qquad$\blacksquare$

(ii)\qquad I$_{312}=\frac{1}{24}\frac{\partial^{2}}{\partial x_{\text{c}}^{2}%
}[\underset{\text{a,b=1}}{\overset{\text{q}}{\sum}}$g$^{\text{ab}}%
(\frac{\partial^{2}\phi}{\partial\text{x}_{\text{a}}\partial\text{x}%
_{\text{b}}}\circ\pi_{P}$ $)](y_{0})$

Since $\frac{\partial g^{\text{ab}}}{\partial x_{\text{c}}}(y_{0}%
)=0=\frac{\partial^{2}g^{\text{ab}}}{\partial x_{\text{c}}^{2}}(y_{0})$ and
g$^{\text{ab}}(y_{0})=\delta^{\text{ab}},$ we have:

\qquad I$_{312}=\frac{1}{8}\frac{\partial^{2}}{\partial x_{\text{c}}^{2}%
}[\underset{\text{a,b=1}}{\overset{\text{q}}{\sum}}$g$^{\text{ab}}%
(\frac{\partial^{2}\phi}{\partial\text{x}_{\text{a}}\partial\text{x}%
_{\text{b}}}\circ\pi_{P}$ $)](y_{0})=\frac{1}{8}[\underset{\text{a,b=1}%
}{\overset{\text{q}}{\sum}}$g$^{\text{ab}}(y_{0})\frac{\partial^{2}}{\partial
x_{\text{c}}^{2}}(\frac{\partial^{2}\phi}{\partial\text{x}_{\text{a}}%
\partial\text{x}_{\text{b}}}\circ\pi_{P}$ $)](y_{0})$

\qquad\qquad$=\frac{1}{8}[\underset{\text{a=1}}{\overset{\text{q}}{\sum}}%
\frac{\partial^{2}}{\partial x_{\text{c}}^{2}}(\frac{\partial^{2}\phi
}{\partial x_{\text{a}}^{2}}\circ\pi_{P}$ $)](y_{0})=\frac{1}{8}%
\underset{\text{a,c=1}}{\overset{\text{q}}{\sum}}\frac{\partial^{4}\phi
}{\partial x_{\text{a}}^{2}\partial x_{\text{c}}^{2}}(y_{0})=\frac{1}%
{8}\underset{\text{a=1}}{\overset{\text{q}}{\sum}}\frac{\partial^{4}\phi
}{\partial x_{\text{a}}^{2}\partial x_{\text{b}}^{2}}(y_{0})$

\qquad\qquad\qquad\qquad\qquad\qquad\qquad\qquad\qquad\qquad\qquad\qquad
\qquad\qquad\qquad\qquad\qquad$\blacksquare$

(iii)\qquad I$_{313}=\frac{1}{24}\frac{\partial^{2}}{\partial x_{\text{c}}%
^{2}}[\underset{\text{i,j=1}}{\overset{\text{n}}{\sum}}g^{ij}\left\{
\frac{\partial\Lambda_{j}}{\partial\text{x}_{i}}\phi\circ\pi_{P}\right\}
](y_{0})$

Since $g^{ij}$ and $\Lambda_{j}$ are expanded in normal Fermi coordinates,
differentiation in tangential Fermi coordinates vanish. Therefore,

\qquad I$_{313}=\frac{1}{8}\frac{\partial^{2}}{\partial x_{\text{c}}^{2}%
}[\underset{i,j=1}{\overset{n}{\sum}}g^{ij}\left\{  \frac{\partial\Lambda_{j}%
}{\partial\text{x}_{i}}\phi\circ\pi_{P}\right\}  ](y_{0})$

$\qquad\ \ \ \ =\frac{1}{8}\underset{i,j=1}{\overset{n}{\sum}}g^{ij}%
(y_{0})[\frac{\partial\Lambda_{j}}{\partial\text{x}_{i}}(y_{0})\frac
{\partial^{2}}{\partial x_{\text{c}}^{2}}(\phi\circ\pi_{P})](y_{0})$

Since $g^{ij}(y_{0})=\delta^{ij}$ and $\frac{\partial\Lambda_{i}}%
{\partial\text{x}_{i}}(y_{0})=0$ by the skew-symmetry of $\frac{1}{2}%
\Omega_{ij}(y_{0})=\frac{\partial\Lambda_{j}}{\partial\text{x}_{i}}(y_{0}).$

Recall that from $\left(  C_{5}\right)  $ and $\left(  C_{6}\right)  ,$

$z_{0}=\gamma($r$_{1})=\left(  \text{y}_{1}(y),...,\text{y}_{q}(y),(1-r_{1}%
)\text{x}_{q+1}(x_{0}),...,(1-r_{1})\text{x}_{n}(x_{0})\right)  $

$\pi_{\text{P}}$(z$_{0}$) $=\left(  \text{y}_{1}(y),...,\text{y}%
_{q}(y),0,...,0\right)  =$ y $=\pi_{\text{P}}$(x)

Consequently,

$\frac{\partial}{\partial\text{x}_{i}}\pi_{\text{P}}(z_{0})=\left\{
\begin{array}
[c]{c}%
1\text{ for }i=\text{c = 1,...,q}\\
0\text{ for }i=q+1,...,n
\end{array}
\right.  $

We have:

\qquad I$_{313}=\frac{1}{8}\underset{i,j=1}{\overset{n}{\sum}}[\frac
{\partial\Lambda_{i}}{\partial\text{x}_{i}}(y_{0})\frac{\partial^{2}}{\partial
x_{\text{c}}^{2}}(\phi\circ\pi_{P})](y_{0})=0$

(iv) The same arguement as in (iii) shows that:\qquad

I$_{314}$ $=\frac{1}{4}\frac{\partial^{2}}{\partial x_{\text{c}}^{2}%
}\underset{\text{j=1}}{\overset{\text{n}}{\sum}}\underset{\text{a=1}%
}{\overset{\text{q}}{\sum}}$g$^{\text{a}j}[$ $\Lambda_{j}\frac{\partial\phi
}{\partial\text{x}_{\text{a}}}\circ\pi_{P}](y_{0})=\frac{1}{4}%
\underset{\text{j=1}}{\overset{\text{n}}{\sum}}\underset{\text{a=1}%
}{\overset{\text{q}}{\sum}}$g$^{\text{a}j}(y_{0})[$ $\Lambda_{j}\frac
{\partial^{2}}{\partial x_{\text{c}}^{2}}(\frac{\partial\phi}{\partial
\text{x}_{\text{a}}}\circ\pi_{P})](y_{0})$

$=\frac{1}{4}\underset{\text{j=1}}{\overset{\text{n}}{\sum}}%
\underset{\text{a=1}}{\overset{\text{q}}{\sum}}\delta^{\text{a}j}[$
$\Lambda_{j}\frac{\partial^{2}}{\partial x_{\text{c}}^{2}}(\frac{\partial\phi
}{\partial\text{x}_{\text{a}}}\circ\pi_{P})](y_{0})=\frac{1}{4}%
\underset{\text{a=1}}{\overset{\text{q}}{\sum}}[$ $\Lambda_{\text{a}}%
\frac{\partial^{2}}{\partial x_{\text{c}}^{2}}(\frac{\partial\phi}%
{\partial\text{x}_{\text{a}}}\circ\pi_{P})](y_{0})$

$=\frac{1}{4}\underset{\text{a=1}}{\overset{\text{q}}{\sum}}[\Lambda
_{\text{a}}\frac{\partial^{3}\phi}{\partial\text{x}_{\text{a}}\partial
x_{\text{b}}^{2}}](y_{0})$

\qquad\qquad\qquad\qquad\qquad\qquad\qquad\qquad\qquad\qquad\qquad\qquad
\qquad\qquad\qquad\qquad$\blacksquare$

(v) For the same reasons as in (iii) and (iv), we have:\qquad

\qquad I$_{315}=$ $\frac{1}{8}\frac{\partial^{2}}{\partial x_{\text{c}}^{2}%
}[\underset{i,j=1}{\overset{n}{\sum}}$g$^{ij}\left\{  \Lambda_{i}\Lambda
_{j}\phi\circ\pi_{P}\right\}  ](y_{0})=$ $\frac{1}{8}%
\underset{i,j=1}{\overset{n}{\sum}}$g$^{ij}(y_{0})[\Lambda_{i}\Lambda_{j}%
\frac{\partial^{2}}{\partial x_{\text{c}}^{2}}\phi\circ\pi_{P}](y_{0})$

\qquad\qquad$=$ $\frac{1}{8}\underset{i=1}{\overset{n}{\sum}}[\Lambda_{i}%
^{2}\frac{\partial^{2}}{\partial x_{\text{c}}^{2}}\phi\circ\pi_{P}](y_{0})$

Since $\Lambda_{i}^{2}(y_{0})=0$ for $i=q+1,...,n$ by $\left(  6.13\right)  $
of Chapter 6,

\qquad I$_{315}=$ $\frac{1}{8}\underset{\text{a=1}}{\overset{\text{q}}{\sum}%
}[\Lambda_{\text{a}}^{2}(y_{0})\frac{\partial^{2}\phi}{\partial x_{\text{c}%
}^{2}}](y_{0})\qquad$

\qquad\qquad\qquad\qquad\qquad\qquad\qquad\qquad\qquad\qquad\qquad\qquad
\qquad\qquad\qquad$\blacksquare$

(vi) Similarly, since expansions of $g^{ij},\Gamma_{ij}^{k}$ and $\Lambda_{k}$
are in normal Fermi coordinates,

I$_{316}=-\frac{1}{8}\frac{\partial^{2}}{\partial x_{\text{c}}^{2}%
}[\underset{i,j=1}{\overset{n}{\sum}}g^{ij}\left\{  \Gamma_{ij}^{\text{c}%
}\frac{\partial\phi}{\partial\text{x}_{\text{c}}}\circ\pi_{P}\text{ }%
+\Gamma_{ij}^{k}\Lambda_{k}\phi\circ\pi_{P}\right\}  ](y_{0})$

\qquad$=-\frac{1}{8}[\underset{i,j=1}{\overset{n}{\sum}}g^{ij}(y_{0})\left\{
\Gamma_{ij}^{\text{c}}\frac{\partial^{2}}{\partial x_{\text{c}}^{2}}%
(\frac{\partial\phi}{\partial\text{x}_{\text{c}}}\circ\pi_{P})+\text{ }%
\Gamma_{ij}^{k}\Lambda_{k}\frac{\partial^{2}}{\partial x_{\text{c}}^{2}}%
(\phi\circ\pi_{P})\right\}  ](y_{0})$

\ \ $=-\frac{1}{8}[\underset{i=1}{\overset{n}{\sum}}\left\{  \Gamma
_{ii}^{\text{c}}\frac{\partial^{2}}{\partial x_{\text{c}}^{2}}(\frac
{\partial\phi}{\partial\text{x}_{\text{c}}}\circ\pi_{P})\text{ }+\text{
}\underset{k=q+\text{1}}{\overset{n}{\sum}}\Gamma_{ii}^{k}\Lambda_{k}%
\frac{\partial^{2}}{\partial x_{\text{c}}^{2}}(\phi\circ\pi_{P})\right\}
](y_{0})$

$=-\frac{1}{8}\underset{\text{a=1}}{\overset{\text{q}}{\sum}}\Gamma
_{\text{aa}}^{\text{c}}(y_{0})\frac{\partial^{2}}{\partial x_{\text{c}}^{2}%
}(\frac{\partial\phi}{\partial\text{x}_{\text{c}}}\circ\pi_{P})(y_{0})$
$-\frac{1}{8}$ $\underset{i\text{=}i+1}{\overset{\text{n}}{\sum}}\Gamma
_{ii}^{\text{c}}(y_{0})\frac{\partial^{2}}{\partial x_{\text{c}}^{2}}%
(\frac{\partial\phi}{\partial\text{x}_{\text{c}}}\circ\pi_{P})$ $(y_{0})$

$-\frac{1}{8}\underset{\text{a=1}}{\overset{\text{q}}{\sum}}%
\underset{k=q+1}{\overset{n}{\sum}}\Gamma_{\text{aa}}^{k}(y_{0})\Lambda
_{k}(y_{0})\frac{\partial^{2}}{\partial x_{\text{c}}^{2}}[\phi\circ\pi
_{P}](y_{0})-\frac{1}{8}\underset{i,k=q+1}{\overset{\text{n}}{\sum}}%
\Gamma_{ii}^{k}(y_{0})\Lambda_{k}(y_{0})\frac{\partial^{2}}{\partial
x_{\text{c}}^{2}}\phi\circ\pi_{P}](y_{0})$

We have: $\Gamma_{\text{aa}}^{\text{c}}(y_{0})=0=\Gamma_{ii}^{\text{c}}%
(y_{0})$ and $\Gamma_{ii}^{k}(y_{0})=0$ for a,c =1,...,q and $i,k=q+1,...n.$

I$_{316}=-\frac{1}{8}\underset{\text{a=1}}{\overset{\text{q}}{\sum}%
}\underset{k=q+1}{\overset{n}{\sum}}\Gamma_{\text{aa}}^{k}(y_{0})\Lambda
_{k}(y_{0})\frac{\partial^{2}\phi}{\partial x_{\text{c}}^{2}}(y_{0})=-\frac
{1}{8}\underset{k=q+1}{\overset{n}{\sum}}<H,k>(y_{0})\Lambda_{k}(y_{0}%
)\frac{\partial^{2}\phi}{\partial x_{\text{c}}^{2}}(y_{0})$

\qquad\qquad\qquad\qquad\qquad\qquad\qquad\qquad\qquad\qquad\qquad\qquad
\qquad\qquad\qquad\qquad\qquad\qquad\qquad$\blacksquare$

(vii)\qquad I$_{317}=\frac{1}{12}\frac{\partial^{2}}{\partial x_{\text{c}}%
^{2}}[$ $\underset{\text{a=1}}{\overset{\text{q}}{\sum}}(\nabla\log
\Psi)_{\text{a}}\frac{\partial\phi}{\partial\text{x}_{\text{a}}}\circ
\pi_{\text{P}}](y_{0})$

Since $\frac{\partial}{\partial x_{\text{c}}}(\nabla\log\Psi)_{\text{a}}%
(y_{0})=0=\frac{\partial^{2}}{\partial x_{\text{c}}^{2}}(\nabla\log
\Psi)_{\text{a}}(y_{0})$ by (xii) and (xiii) of \textbf{Table B}$_{1},$ we have:

\qquad I$_{317}=\frac{1}{4}$ $\underset{\text{a=1}}{\overset{\text{q}}{\sum}%
}(\nabla\log\Psi)_{\text{a}}(y_{0})\frac{\partial^{2}}{\partial x_{\text{c}%
}^{2}}[\frac{\partial\phi}{\partial\text{x}_{\text{a}}}\circ\pi_{\text{P}%
}](y_{0})$

$(\nabla\log\Psi)_{\text{a}}(y_{0})=(\nabla\log\theta^{-\frac{1}{2}%
})_{\text{a}}(y_{0})+(\nabla\log\Phi)_{\text{a}}(y_{0})=0$ by (iii)$^{\ast}$
of \textbf{Table A}$_{9}$ and (xi) of \textbf{Table B}$_{1}.$

Therefore,\qquad

\qquad\qquad I$_{317}=0.$

\qquad\qquad\qquad$\qquad\qquad\qquad\qquad\qquad\qquad\qquad\qquad
\qquad\qquad\qquad\blacksquare$

(viii)\qquad I$_{318}=\frac{1}{4}\frac{\partial^{2}}{\partial x_{\text{a}}%
^{2}}[\underset{j=1}{\overset{n}{\sum}}(\nabla\log\Psi)_{j}\Lambda_{j}%
\phi\circ\pi_{P}](y_{0})$

\qquad$=\frac{1}{4}\underset{j=1}{\overset{n}{\sum}}[\frac{\partial^{2}%
}{\partial x_{\text{a}}^{2}}(\nabla\log\Psi)_{j})\Lambda_{j}\phi\circ\pi
_{P}](y_{0})+\frac{1}{4}\underset{j=1}{\overset{n}{\sum}}[(\nabla\log\Psi
)_{j}\frac{\partial^{2}}{\partial x_{\text{a}}^{2}}(\Lambda_{j}\phi\circ
\pi_{P})](y_{0})$

\qquad$+\frac{1}{2}\underset{j=1}{\overset{n}{\sum}}[\frac{\partial}{\partial
x_{\text{a}}}(\nabla\log\Psi)_{j})\frac{\partial}{\partial x_{\text{a}}%
}(\Lambda_{j}\phi\circ\pi_{P})](y_{0})=$ A$_{1}+$ A$_{2}+$ A$_{3}$where,

\qquad\ A$_{1}=\frac{1}{4}\underset{j=1}{\overset{n}{\sum}}[\frac{\partial
^{2}}{\partial x_{\text{a}}^{2}}(\nabla\log\Psi)_{j})(y_{0})(\Lambda_{j}%
\phi\circ\pi_{P})(y_{0})]$

$\qquad$A$_{2}=\frac{1}{4}\underset{j=1}{\overset{n}{\sum}}[(\nabla\log
\Psi)_{j}(y_{0})\frac{\partial^{2}}{\partial x_{\text{a}}^{2}}(\Lambda_{j}%
\phi\circ\pi_{P})(y_{0})]$

$\qquad$A$_{3}=\frac{1}{2}\underset{j=1}{\overset{n}{\sum}}[\frac{\partial
}{\partial x_{\text{a}}}(\nabla\log\Psi)_{j})(y_{0})\frac{\partial}{\partial
x_{\text{a}}}(\Lambda_{j}\phi\circ\pi_{P})(y_{0})]$

We examine each of these:

\qquad\ A$_{1}=\frac{1}{12}\underset{j=1}{\overset{n}{\sum}}[\frac
{\partial^{2}}{\partial x_{\text{a}}^{2}}(\nabla\log\Psi)_{j})(y_{0}%
)(\Lambda_{j}\phi\circ\pi_{P})(y_{0})]$

\qquad$\qquad=\frac{1}{12}\underset{\text{b=1}}{\overset{\text{q}}{\sum}%
}[\frac{\partial^{2}}{\partial x_{\text{a}}^{2}}(\nabla\log\Psi)_{\text{b}%
})(y_{0})(\Lambda_{\text{b}}\phi\circ\pi_{P})(y_{0})]$

$\qquad\qquad+\frac{1}{12}\underset{j=q+1}{\overset{\text{n}}{\sum}}%
[\frac{\partial^{2}}{\partial x_{\text{a}}^{2}}(\nabla\log\Psi)_{j}%
)(y_{0})(\Lambda_{j}\phi\circ\pi_{P})(y_{0})]$

We have:

\qquad$\frac{\partial}{\partial x_{\text{a}}}(\nabla\log\theta^{-\frac{1}{2}%
})_{\text{b}}(y_{0})=0=\frac{\partial^{2}}{\partial x_{\text{a}}^{2}}%
(\nabla\log\theta^{-\frac{1}{2}})_{\text{b}}(y_{0})$ because the expansion of
$\theta$ is in \textbf{normal Fermi coordinates} and so all differentiation
with respect to tangential coordinates vanish. On the other hand, by (xii) and
(xiii) of \textbf{Table B}$_{1},$ we also have:

\qquad$\frac{\partial}{\partial x_{\text{a}}}(\nabla\log\Phi)_{\text{b}%
})(y_{0})=0=\frac{\partial^{2}}{\partial x_{\text{a}}^{2}}(\nabla\log
\Phi)_{\text{b}})(y_{0}).$ Therefore $\frac{\partial^{2}}{\partial
x_{\text{a}}^{2}}(\nabla\log\Psi)_{\text{b}})(y_{0})=0$ and so,

$\qquad$A$_{1}=\frac{1}{4}\underset{j=q+1}{\overset{n}{\sum}}[\frac
{\partial^{2}}{\partial x_{\text{a}}^{2}}(\nabla\log\Psi)_{j})(y_{0}%
)(\Lambda_{j}\phi\circ\pi_{P})(y_{0})]$

$=\frac{1}{4}\underset{j=q+1}{\overset{n}{\sum}}[\frac{\partial^{2}}{\partial
x_{\text{a}}^{2}}(\nabla\log\theta^{-\frac{1}{2}})_{j})\Lambda_{j}\phi\circ
\pi_{P}](y_{0})+\frac{1}{4}\underset{j=q+1}{\overset{n}{\sum}}[\frac
{\partial^{2}}{\partial x_{\text{a}}^{2}}(\nabla\log\Phi)_{j})\Lambda_{j}%
\phi\circ\pi_{P}](y_{0})$

Again, $\frac{\partial^{2}}{\partial x_{\text{a}}^{2}}(\nabla\log
\theta^{-\frac{1}{2}})_{j})(y_{0})=0$ because the expansion of $\theta$ is in
\textbf{normal Fermi coordinates} and so all differentiation with respect to
tangential coordinates vanish. Therefore, we have:

\qquad A$_{1}=\frac{1}{4}\underset{j=q+1}{\overset{n}{\sum}}[\frac
{\partial^{2}}{\partial x_{\text{a}}^{2}}(\nabla\log\Phi)_{j})(y_{0}%
)(\Lambda_{j}\phi\circ\pi_{P})(y_{0})]$

By (x) of \textbf{Table B}$_{1},$

A$_{1}=\frac{1}{4}\underset{j=q+1}{\overset{n}{\sum}}[\frac{\partial^{2}%
}{\partial x_{\text{a}}^{2}}(\nabla\log\Phi)_{j})(y_{0})(\Lambda_{j}\phi
\circ\pi_{P})(y_{0})]$

$=-\frac{1}{4}\underset{j=q+1}{\overset{n}{\sum}}[\frac{\partial^{2}X_{j}%
}{\partial x_{\text{a}}^{2}}(y_{0})(\Lambda_{j}\phi\circ\pi_{P})(y_{0}%
)]=-\frac{1}{4}\underset{j=q+1}{\overset{n}{\sum}}[\frac{\partial^{2}X_{j}%
}{\partial x_{\text{a}}^{2}}\Lambda_{j}\phi](y_{0})$

Next we examine:

A$_{2}=\frac{1}{4}\underset{j=1}{\overset{n}{\sum}}[(\nabla\log\Psi)_{j}%
\frac{\partial^{2}}{\partial x_{\text{a}}^{2}}(\Lambda_{j}\phi\circ\pi
_{P})](y_{0})$

A$_{2}=\frac{1}{4}\underset{\text{b=1}}{\overset{\text{q}}{\sum}}[(\nabla
\log\Psi)_{\text{b}}\frac{\partial^{2}}{\partial x_{\text{a}}^{2}}(\Lambda
_{j}\phi\circ\pi_{P})](y_{0})$

$\qquad+\frac{1}{4}\underset{j=q+1}{\overset{n}{\sum}}[(\nabla\log\Psi
)_{j}\frac{\partial^{2}}{\partial x_{\text{a}}^{2}}(\Lambda_{j}\phi\circ
\pi_{P})](y_{0})$

By (iii)$^{\ast}$ of \textbf{Table A}$_{9}$ and by (xi) of \textbf{Table
B}$_{1}$

\qquad$(\nabla\log\Psi)_{\text{b}}(y_{0})=(\nabla\log\theta^{-\frac{1}{2}%
})_{\text{b}}(y_{0})+(\nabla\log\Phi)_{\text{b}}(y_{0})=0$

Next we have by (iv)$^{\ast}$ of \textbf{Table 9 }and (vi) of \textbf{Table
B}$_{1},$

$(\nabla\log\Psi)_{j}(y_{0})=(\nabla\log\theta^{-\frac{1}{2}})_{j}%
(y_{0})+(\nabla\log\Phi)_{j}(y_{0})=\frac{1}{2}<H,j>(y_{0})-X_{j}(y_{0})$

Next we have:

\qquad$\frac{\partial^{2}}{\partial x_{\text{a}}^{2}}(\Lambda_{j}\phi\circ
\pi_{P})](y_{0})=\frac{\partial^{2}\Lambda_{j}}{\partial x_{\text{a}}^{2}%
}(y_{0})(\phi\circ\pi_{P})(y_{0})+\Lambda_{j}(y_{0})\frac{\partial^{2}%
}{\partial x_{\text{a}}^{2}}(\phi\circ\pi_{P})(y_{0})$

$\qquad+\frac{\partial\Lambda_{j}}{\partial x_{\text{a}}}(y_{0})\frac
{\partial}{\partial x_{\text{a}}}\phi\circ\pi_{P}(y_{0})$

Since $\frac{\partial^{2}\Lambda_{j}}{\partial x_{\text{a}}^{2}}%
(y_{0})=0=\frac{\partial\Lambda_{j}}{\partial x_{\text{a}}}(y_{0}),$ and
$\frac{\partial^{2}}{\partial x_{\text{a}}^{2}}(\phi\circ\pi_{P})(y_{0}%
)=\frac{\partial^{2}\phi}{\partial x_{\text{a}}^{2}}(y_{0}),$ we have:

\qquad A$_{2}=\frac{1}{4}\underset{j=q+1}{\overset{\text{n}}{\sum}}[\frac
{1}{2}<H,j>-X_{j}](y_{0})\Lambda_{j}(y_{0})\frac{\partial^{2}\phi}{\partial
x_{\text{a}}^{2}}(y_{0})$

Here we finally consider:

A$_{3}=\frac{1}{2}\underset{j=1}{\overset{n}{\sum}}[\frac{\partial}{\partial
x_{\text{a}}}(\nabla\log\Psi)_{j})(y_{0})\frac{\partial}{\partial x_{\text{a}%
}}(\Lambda_{j}\phi\circ\pi_{P})(y_{0})]$

$=\frac{1}{2}\underset{\text{b=1}}{\overset{\text{q}}{\sum}}[\frac{\partial
}{\partial x_{\text{a}}}(\nabla\log\Psi)_{\text{b}})\frac{\partial}{\partial
x_{\text{a}}}(\Lambda_{j}\phi\circ\pi_{P})](y_{0})$

$+\frac{1}{2}\underset{j=q+1}{\overset{n}{\sum}}[\frac{\partial}{\partial
x_{\text{a}}}(\nabla\log\Psi)_{j})\frac{\partial}{\partial x_{\text{a}}%
}(\Lambda_{j}\phi\circ\pi_{P})](y_{0})$

We saw earlier that,

$\frac{\partial}{\partial x_{\text{a}}}(\nabla\log\theta^{-\frac{1}{2}%
})_{\text{b}}(y_{0})=0=\frac{\partial}{\partial x_{\text{a}}}(\nabla\log
\Phi_{P})_{\text{b}})(y_{0})$

Therefore,

$\frac{\partial}{\partial x_{\text{a}}}(\nabla\log\Psi)_{\text{b}}%
)(y_{0})=\frac{\partial}{\partial x_{\text{a}}}(\nabla\log\theta^{-\frac{1}%
{2}})_{\text{b}}(y_{0})+\frac{\partial}{\partial x_{\text{a}}}(\nabla\log
\Phi_{P})_{\text{b}})(y_{0})=0$

Consequently,

A$_{3}=\frac{1}{2}\underset{j=q+1}{\overset{n}{\sum}}[\frac{\partial}{\partial
x_{\text{a}}}(\nabla\log\Psi)_{j})\frac{\partial}{\partial x_{\text{a}}%
}(\Lambda_{j}\phi\circ\pi_{P})](y_{0})$

$=\frac{1}{2}\underset{j=q+1}{\overset{n}{\sum}}\frac{\partial}{\partial
x_{\text{a}}}(\nabla\log\theta^{-\frac{1}{2}})_{j})(y_{0})\frac{\partial
}{\partial x_{\text{a}}}(\Lambda_{j}\phi\circ\pi_{P})(y_{0})$

$+\frac{1}{2}\underset{j=q+1}{\overset{n}{\sum}}\frac{\partial}{\partial
x_{\text{a}}}(\nabla\log\Phi)_{j})(y_{0})\frac{\partial}{\partial x_{\text{a}%
}}(\Lambda_{j}\phi\circ\pi_{P})](y_{0})$

We saw earlier that, $\frac{\partial}{\partial x_{\text{a}}}(\nabla\log
\theta^{-\frac{1}{2}})_{j})(y_{0})=0.$ Now by (ix) of \textbf{Table B}$_{1}%
$\textbf{ }we have:

$\frac{\partial}{\partial x_{\text{a}}}(\nabla$log$\Phi_{P})_{j}(y_{0}%
)=-\frac{\partial X_{j}}{\partial x_{\text{a}}}(y_{0})$ for a $=$ 1,...,q and
$j$ $=q+1,...,n.$

\ $\frac{\partial}{\partial x_{\text{a}}}(\Lambda_{j}\phi\circ\pi_{P}%
)(y_{0})=\frac{\partial\Lambda_{j}}{\partial x_{\text{a}}}(y_{0})(\phi\circ
\pi_{P})(y_{0})+\Lambda_{j}(y_{0})\frac{\partial}{\partial x_{\text{a}}}%
(\phi\circ\pi_{P})(y_{0})$

$=0+\Lambda_{j}(y_{0})\frac{\partial\phi}{\partial x_{\text{a}}}%
(y_{0})=\Lambda_{j}(y_{0})\frac{\partial\phi}{\partial x_{\text{a}}}(y_{0})$

Therefore,

A$_{3}=\frac{1}{2}\underset{j=q+1}{\overset{n}{\sum}}-\frac{\partial X_{j}%
}{\partial x_{\text{a}}}(y_{0})\Lambda_{j}(y_{0})\frac{\partial\phi}{\partial
x_{\text{a}}}(y_{0})$ \ \ \ \ \ \ \ \ \ \ \ \ \ \ \ 

We have finally,

I$_{318}=$ A$_{1}$ + A$_{2}$ + A$_{3}=-\frac{1}{4}%
\underset{j=q+1}{\overset{\text{n}}{\sum}}[\frac{\partial^{2}X_{j}}{\partial
x_{\text{a}}^{2}}\Lambda_{j}\phi](y_{0})+\frac{1}{2}%
\underset{j=q+1}{\overset{n}{\sum}}[-\frac{\partial X_{j}}{\partial
x_{\text{a}}}\Lambda_{j}\frac{\partial\phi}{\partial x_{\text{a}}}](y_{0})$ \ \ 

$+\frac{1}{8}\underset{j=q+1}{\overset{\text{n}}{\sum}}[<H,j>-X_{j}%
]\Lambda_{j}\frac{\partial^{2}\phi}{\partial x_{\text{a}}^{2}}](y_{0})$

\qquad\qquad\qquad\qquad\qquad\qquad\qquad\qquad\qquad\qquad\qquad\qquad
\qquad\qquad\qquad\qquad\qquad$\blacksquare$

(ix)\qquad We recall that W is the Weitzenb\H{o}ckian which is an element of
the vector bundle $\Gamma($End(E)$).$

We have easily:

I$_{319}=\frac{1}{8}\frac{\partial^{2}}{\partial x_{\text{a}}^{2}}[$%
W$\phi\circ\pi_{P}](y_{0})$

\qquad$=\frac{1}{8}\frac{\partial^{2}\text{W}}{\partial x_{\text{a}}^{2}%
}(y_{0})(\phi\circ\pi_{P})(y_{0})+\frac{1}{4}\frac{\partial\text{W}}{\partial
x_{\text{a}}}(y_{0})\frac{\partial}{\partial x_{\text{a}}}(\phi\circ\pi
_{P})(y_{0})+\frac{1}{8}$W(y$_{0}$)$\phi(y_{0})\frac{\partial^{2}}{\partial
x_{\text{a}}^{2}}[\phi\circ\pi_{P}](y_{0})$

$\qquad\qquad$

\qquad$=\frac{1}{8}\frac{\partial^{2}\text{W}}{\partial x_{\text{a}}^{2}%
}(y_{0})\phi(y_{0})+\frac{1}{4}\frac{\partial\text{W}}{\partial x_{\text{a}}%
}(y_{0})\frac{\partial\phi}{\partial x_{\text{a}}}(y_{0})+\frac{1}{8}$%
W$(y_{0})\frac{\partial^{2}\phi}{\partial x_{\text{a}}^{2}}(y_{0})$

\qquad\qquad\qquad\qquad\qquad\qquad\qquad\qquad\qquad\qquad\qquad\qquad
\qquad\qquad\qquad\qquad$\blacksquare$

We next consider:

(x)\qquad L$_{1}=\frac{1}{4}\frac{\partial^{2}}{\partial x_{\text{a}}^{2}}[$
$\underset{\text{a=1}}{\overset{\text{q}}{\sum}}$X$_{\text{a}}\frac
{\partial\phi}{\partial\text{x}_{\text{a}}}\circ\pi](y_{0})$

$\qquad=\frac{1}{4}[$ $\underset{\text{a=1}}{\overset{\text{q}}{\sum}}%
\frac{\partial^{2}\text{X}_{\text{a}}}{\partial x_{\text{a}}^{2}}(y_{0}%
)\frac{\partial\phi}{\partial\text{x}_{\text{a}}}(y_{0})+\frac{1}{4}[$
$\underset{\text{a=1}}{\overset{\text{q}}{\sum}}$X$_{\text{a}}(y_{0}%
)\frac{\partial^{3}\phi}{\partial x_{\text{a}}^{3}}(y_{0})+\frac{1}{2}[$
$\underset{\text{a=1}}{\overset{\text{q}}{\sum}}\frac{\partial\text{X}%
_{\text{a}}}{\partial x_{\text{a}}}(y_{0})\frac{\partial^{2}\phi}%
{\partial\text{x}_{\text{a}}^{2}}(y_{0})$

\qquad\qquad\qquad\qquad\qquad\qquad\qquad\qquad\qquad\qquad\qquad\qquad
\qquad\qquad\qquad\qquad$\blacksquare$

(xi)\qquad L$_{2}=\frac{1}{4}\frac{\partial^{2}}{\partial x_{\text{a}}^{2}}[$
$\underset{j=1}{\overset{n}{\sum}}$ X$_{j}\Lambda_{j}\phi\circ\pi_{\text{P}%
}](y_{0})$

Since $\frac{\partial\Lambda_{j}}{\partial x_{\text{a}}}(y_{0})=0=\frac
{\partial^{2}\Lambda_{j}}{\partial x_{\text{a}}^{2}}(y_{0}),$ we have:

$\left(  D_{7}\right)  ^{\ast\ast}$\qquad L$_{2}=\frac{1}{4}[$
$\underset{j=1}{\overset{n}{\sum}}$ $\frac{\partial^{2}X_{j}}{\partial
x_{\text{a}}^{2}}(y_{0})\Lambda_{j}(y_{0})\phi(y_{0})+\frac{1}{2}%
[\underset{j=1}{\overset{n}{\sum}}\frac{\partial X_{j}}{\partial x_{\text{a}}%
}(y_{0})\Lambda_{j}(y_{0})\frac{\partial\phi}{\partial x_{\text{a}}}(y_{0})$

$\qquad+\frac{1}{4}[$ $\underset{j=1}{\overset{n}{\sum}}$ X$_{j}(y_{0}%
)\Lambda_{j}(y_{0})\frac{\partial^{2}\phi}{\partial x_{\text{a}}^{2}}(y_{0})$

\qquad\qquad\qquad\qquad\qquad\qquad\qquad\qquad\qquad\qquad\qquad\qquad
\qquad\qquad\qquad\qquad\qquad$\blacksquare$

We now gather all (sub)-expressions of I$_{31}$ and have:

$\left(  D_{8}\right)  \qquad$I$_{31}=\frac{1}{4}\underset{\text{c=1}%
}{\overset{\text{q}}{\sum}}\frac{\partial^{2}}{\partial x_{\text{c}}^{2}}%
$[L$_{\Psi}\phi\circ\pi_{P}$]$(y_{0})$

$\qquad=$ I$_{311}+$ I$_{312}$+ I$_{313}+$ I$_{314}+$ I$_{315}+$ I$_{316}+$
I$_{317}+$ I$_{318}+$ I$_{319}+$ L$_{1}+$ L$_{2}\qquad$

\qquad$=\frac{1}{96}[\underset{i=q+1}{\overset{n}{\sum}}3<H,i>^{2}+2(\tau
^{M}-3\tau^{P}\ +\overset{q}{\underset{\text{a=1}}{\sum}}\varrho_{\text{aa}%
}^{M}+\overset{q}{\underset{\text{a,b}=1}{\sum}}R_{\text{abab}}^{M}%
)](y_{0})\frac{\partial^{2}\phi}{\partial x_{\text{c}}^{2}}(y_{0})\qquad
$I$_{311}$

$-\frac{1}{4}[\left\Vert \text{X}(y_{0})\right\Vert ^{2}+$ divX$(y_{0})-$
$\underset{\text{a}=1}{\overset{q}{\sum}}($X$_{\text{a}})^{2}(y_{0})$ $-$
$\underset{\text{a}=1}{\overset{q}{\sum}}\frac{\partial X_{\text{a}}}{\partial
x_{\text{a}}}(y_{0})]\frac{\partial^{2}\phi}{\partial x_{\text{c}}^{2}}%
(y_{0})+$ $\frac{1}{4}$V(y$_{0}$)$\frac{\partial^{2}\phi}{\partial
x_{\text{c}}^{2}}(y_{0})$

$-\frac{1}{2}$ $[X_{j}\frac{\partial X_{j}}{\partial x_{\text{c}}}+\frac{1}%
{2}\frac{\partial^{2}X_{j}}{\partial x_{\text{c}}\partial x_{j}}](y_{0}%
).\frac{\partial\phi}{\partial x_{\text{c}}}(y_{0})+\frac{1}{4}[<H,j>\frac
{\partial X_{j}}{\partial x_{\text{c}}}](y_{0}).\frac{\partial\phi}{\partial
x_{\text{c}}}(y_{0})$

$+$ $\frac{1}{2}\frac{\partial\text{V}}{\partial x_{\text{c}}}(y_{0}%
).\frac{\partial\phi}{\partial x_{\text{c}}}(y_{0})+\frac{1}{4}[(\frac
{\partial X_{j}}{\partial x_{\text{c}}})^{2}-X_{j}\frac{\partial^{2}X_{j}%
}{\partial x_{\text{c}}^{2}}](y_{0})\phi(y_{0})$ $-\frac{1}{8}\frac
{\partial^{3}X_{j}}{\partial x_{\text{c}}^{2}\partial x_{j}}(y_{0})\phi
(y_{0})$

$-\frac{1}{2}\frac{\partial X_{i}}{\partial x_{\text{c}}}(y_{0}))\frac
{\partial X_{i}}{\partial x_{\text{c}}}(y_{0})\phi(y_{0})+\frac{1}{4}%
\frac{\partial^{2}\text{V}}{\partial x_{\text{c}}^{2}}(y_{0})\phi(y_{0})$

$+\frac{1}{8}\underset{\text{a=1}}{\overset{\text{q}}{\sum}}\frac{\partial
^{4}\phi}{\partial x_{\text{a}}^{2}\partial x_{\text{c}}^{2}}(y_{0}%
)\qquad\qquad\qquad$I$_{312}$

$+\frac{1}{4}\underset{\text{a=1}}{\overset{\text{q}}{\sum}}[\Lambda
_{\text{a}}\frac{\partial^{3}\phi}{\partial\text{x}_{\text{a}}\partial
x_{\text{c}}^{2}}](y_{0})\qquad\qquad$I$_{314}$

$+$ $\frac{1}{8}\underset{\text{a=1}}{\overset{\text{q}}{\sum}}[\Lambda
_{\text{a}}^{2}(y_{0})\frac{\partial^{2}\phi}{\partial x_{\text{c}}^{2}%
}](y_{0})\qquad\ \ \ \ $I$_{315}\qquad$

$+\frac{1}{8}\frac{\partial^{2}\text{W}}{\partial x_{\text{a}}^{2}}(y_{0}%
)\phi(y_{0})+\frac{1}{4}\frac{\partial\text{W}}{\partial x_{\text{a}}}%
(y_{0})\frac{\partial\phi}{\partial x_{\text{a}}}(y_{0})+\frac{1}{8}$%
W$(y_{0})\frac{\partial^{2}\phi}{\partial x_{\text{a}}^{2}}(y_{0})\qquad
$I$_{319}$

$+\frac{1}{4}$ $\underset{\text{a=1}}{\overset{\text{q}}{\sum}}\frac
{\partial^{2}\text{X}_{\text{a}}}{\partial x_{\text{a}}^{2}}(y_{0}%
)\frac{\partial\phi}{\partial\text{x}_{\text{a}}}(y_{0})+\frac{1}{4}$
$\underset{\text{a=1}}{\overset{\text{q}}{\sum}}$X$_{\text{a}}(y_{0}%
)\frac{\partial^{3}\phi}{\partial x_{\text{a}}^{3}}(y_{0})+\frac{1}{2}$
$\underset{\text{a=1}}{\overset{\text{q}}{\sum}}\frac{\partial\text{X}%
_{\text{a}}}{\partial x_{\text{a}}}(y_{0})\frac{\partial^{2}\phi}%
{\partial\text{x}_{\text{a}}^{2}}(y_{0})\qquad$L$_{1}$

$+\frac{1}{4}[$ $\underset{\text{b=1}}{\overset{\text{q}}{\sum}}$
$\frac{\partial^{2}X_{\text{b}}}{\partial x_{\text{a}}^{2}}\Lambda_{\text{b}%
}(y_{0})\phi(y_{0})+\frac{1}{4}[$ $\underset{\text{b=1}}{\overset{\text{q}%
}{\sum}}$ X$_{\text{b}}(y_{0})\Lambda_{\text{b}}(y_{0})\frac{\partial^{2}\phi
}{\partial x_{\text{a}}^{2}}(y_{0})\qquad$L$_{2}$

$+\frac{1}{2}[$ $\underset{\text{b=1}}{\overset{\text{q}}{\sum}}$
$\frac{\partial X_{\text{b}}}{\partial x_{\text{a}}}(y_{0})\Lambda_{\text{b}%
}(y_{0})\frac{\partial\phi}{\partial x_{\text{a}}}(y_{0})\qquad$

\qquad\qquad\qquad\qquad\qquad\qquad\qquad\qquad\qquad\qquad\qquad\qquad
\qquad\qquad\qquad\qquad$\blacksquare$

\section{\textbf{Computation of I}$_{32}$}

$\left(  D_{9}\right)  \qquad$\textbf{Computation of I}$_{32}=\frac{1}%
{12}\underset{i=q+1}{\overset{n}{\sum}}\frac{\partial^{2}\Theta}{\partial
x_{i}^{2}}(y_{0})$ for $i=q+1,...,n.$

$\qquad\ \ \ \Theta=\qquad$L$_{\Psi}$[$\phi\circ\pi_{\text{P}}$] as defined in
$\left(  10.10\right)  $ of \textbf{Chapter 10}.

\qquad\qquad\qquad\ \ \ See also $\left(  11.22\right)  $ of \textbf{Chapter
11} here.

\qquad\ \ \ $\Theta$\textbf{ }$=\frac{\text{L}\Psi}{\Psi}\phi\circ
\pi_{\text{P}}$ + \ $\frac{1}{2}\underset{\text{a,b=1}}{\overset{\text{q}%
}{\sum}}$g$^{\text{ab}}\left\{  \frac{\partial^{2}\phi}{\partial
\text{x}_{\text{a}}\partial\text{x}_{\text{b}}}\circ\pi_{\text{P}}\text{
}\right\}  +\frac{1}{2}\underset{i,j=1}{\overset{n}{\sum}}$g$^{jk}\left\{
\frac{\partial\Lambda_{k}}{\partial\text{x}_{j}}\phi\circ\pi_{\text{P}%
}\right\}  $

$\bigskip\qquad\ \ +$ $\underset{j=1}{\overset{n}{\sum}}\underset{\text{a=1}%
}{\overset{\text{q}}{\sum}}$g$^{\text{a}j}\left\{  \text{ }\Lambda_{j}%
\frac{\partial\phi}{\partial\text{x}_{\text{a}}}\circ\pi_{\text{P}}\right\}  $

$\qquad\ \ +\frac{1}{2}\underset{j,k=1}{\overset{n}{\sum}}$g$^{jk}\Lambda
_{j}\Lambda_{k}\phi\circ\pi_{\text{P}}-\frac{1}{2}%
\underset{j,k=1}{\overset{n}{\sum}}$g$^{jk}\left\{  \underset{\text{c=1}%
}{\overset{\text{q}}{\sum}}\Gamma_{jk}^{\text{c}}\frac{\partial\phi}%
{\partial\text{x}_{\text{c}}}\circ\pi_{\text{P}}\text{ + }%
\underset{l=1}{\overset{n}{\sum}}\Gamma_{jk}^{l}\Lambda_{l}\phi\circ
\pi_{\text{P}}\right\}  $

$\qquad\ +\frac{1}{2}$W$\phi\circ\pi_{\text{P}}$

\qquad\ \ + $\underset{\text{a=1}}{\overset{\text{q}}{\sum}}(\nabla\log
\Psi)_{\text{a}}\frac{\partial\phi}{\partial\text{x}_{\text{a}}}\circ
\pi_{\text{P}}$ + $\underset{j=1}{\overset{n}{\sum}}$ $(\nabla\log\Psi
)_{j}\Lambda_{j}\phi\circ\pi_{\text{P}}.$

\qquad\ \ + $\underset{\text{a=1}}{\overset{\text{q}}{\sum}}$X$_{\text{a}%
}\frac{\partial\phi}{\partial\text{x}_{\text{a}}}\circ\pi_{\text{P}}$ +
$\underset{j=1}{\overset{n}{\sum}}$ X$_{j}\Lambda_{j}\phi\circ\pi_{\text{P}}$

Therefore we have:

$\left(  D_{9}\right)  \qquad\frac{1}{12}\underset{i=q+1}{\overset{n}{\sum}%
}\frac{\partial^{2}\Theta}{\partial x_{i}^{2}}(y_{0})$

$\qquad\qquad=$ \textbf{I}$_{321}+$ \textbf{I}$_{322}+$ \textbf{I}$_{323}+$
\textbf{I}$_{324}+$ \textbf{I}$_{325}+$ \textbf{I}$_{326}+$ \textbf{I}%
$_{327}+$ \textbf{I}$_{328}+$ \textbf{I}$_{329}+$ L$_{1}+$ L$_{2}$ where,

(i)\qquad\textbf{I}$_{321}=\frac{1}{12}\underset{i=q+1}{\overset{n}{\sum}%
}\frac{\partial^{2}}{\partial x_{i}^{2}}[\frac{\text{L}\Psi}{\Psi}\phi\circ
\pi_{\text{P}}](y_{0})$

(ii)\qquad\textbf{I}$_{322}=\frac{1}{24}\underset{i=q+1}{\overset{n}{\sum}%
}\frac{\partial^{2}}{\partial x_{i}^{2}}[\underset{\text{a,b=1}%
}{\overset{\text{q}}{\sum}}$g$^{\text{ab}}\left\{  \frac{\partial^{2}\phi
}{\partial\text{x}_{\text{a}}\partial\text{x}_{\text{b}}}\circ\pi_{\text{P}%
}\text{ }\right\}  ](y_{0})$

(iii)\qquad\textbf{I}$_{323}=\frac{1}{24}\underset{i=q+1}{\overset{n}{\sum}%
}\frac{\partial^{2}}{\partial x_{i}^{2}}[\underset{j,k=1}{\overset{n}{\sum}}%
$g$^{jk}\left\{  \frac{\partial\Lambda_{k}}{\partial\text{x}_{j}}\phi\circ
\pi_{\text{P}}\right\}  ](y_{0})$

(iv)\qquad\textbf{I}$_{324}=\frac{1}{12}\underset{i=q+1}{\overset{n}{\sum}%
}\frac{\partial^{2}}{\partial x_{i}^{2}}[\underset{j=1}{\overset{n}{\sum}%
}\underset{\text{a=1}}{\overset{\text{q}}{\sum}}$g$^{\text{a}j}\left\{  \text{
}\Lambda_{j}\frac{\partial\phi}{\partial\text{x}_{\text{a}}}\circ\pi
_{\text{P}}\right\}  ](y_{0})$

(v)\qquad\textbf{I}$_{325}=\frac{1}{24}\underset{i=q+1}{\overset{n}{\sum}%
}\frac{\partial^{2}}{\partial x_{i}^{2}}[\underset{i,j=1}{\overset{n}{\sum}}%
$g$^{jk}\Lambda_{j}\Lambda_{k}\phi\circ\pi_{\text{P}}](y_{0})$

(vi)\qquad\textbf{I}$_{326}=-\frac{1}{24}\underset{i=q+1}{\overset{n}{\sum}%
}\frac{\partial^{2}}{\partial x_{i}^{2}}[\underset{j,k=1}{\overset{n}{\sum}}%
$g$^{jk}\left\{  \underset{\text{c=1}}{\overset{\text{q}}{\sum}}\Gamma
_{jk}^{\text{c}}\frac{\partial\phi}{\partial\text{x}_{\text{c}}}\circ
\pi_{\text{P}}\text{ }+\text{ }\underset{l=1}{\overset{n}{\sum}}\Gamma
_{jk}^{l}\Lambda_{l}\phi\circ\pi_{\text{P}}\right\}  ](y_{0})$

(vii)\qquad\textbf{I}$_{327}=\frac{1}{24}\underset{i=q+1}{\overset{n}{\sum}%
}\frac{\partial^{2}}{\partial x_{i}^{2}}[$W$\phi\circ\pi_{\text{P}}](y_{0})$

(viii)\qquad\textbf{I}$_{328}=\frac{1}{12}\underset{i=q+1}{\overset{n}{\sum}}$
$\underset{\text{a=1}}{\overset{\text{q}}{\sum}}\frac{\partial^{2}}{\partial
x_{i}^{2}}[(\nabla\log\Psi)_{\text{a}}\frac{\partial\phi}{\partial
\text{x}_{\text{a}}}\circ\pi_{\text{P}}](y_{0})$

(ix)\qquad\textbf{I}$_{329}=\frac{1}{12}\underset{i=q+1}{\overset{n}{\sum}}$
$\underset{j=1}{\overset{\text{n}}{\sum}}$ $\frac{\partial^{2}}{\partial
x_{i}^{2}}[(\nabla\log\Psi)_{j}\Lambda_{j}\phi\circ\pi_{\text{P}}](y_{0})$

(x)\qquad\textbf{L}$_{1}=\frac{1}{12}\underset{i=q+1}{\overset{n}{\sum}}$
$\underset{\text{a=1}}{\overset{\text{q}}{\sum}}\frac{\partial^{2}}{\partial
x_{i}^{2}}[$X$_{\text{a}}\frac{\partial\phi}{\partial\text{x}_{\text{a}}}%
\circ\pi_{\text{P}}](y_{0})$

(xi)\qquad\textbf{L}$_{2}=\frac{1}{12}\underset{i=q+1}{\overset{n}{\sum}}$
$\underset{\text{j=1}}{\overset{\text{n}}{\sum}}$ $\frac{\partial^{2}%
}{\partial x_{i}^{2}}[$X$_{j}\Lambda_{j}\phi\circ\pi_{\text{P}}](y_{0})$

\textbf{Computations:}

(i)\qquad$\frac{1}{12}\frac{\partial^{2}}{\partial x_{i}^{2}}[\frac
{\text{L}\Psi}{\Psi}\phi\circ\pi_{\text{P}}](y_{0})$ for $i=q+1,...,n$

\qquad$=\frac{1}{12}\frac{\partial^{2}}{\partial x_{i}^{2}}[\frac{\text{L}%
\Psi}{\Psi}](y_{0})[\phi\circ\pi_{\text{P}}](y_{0})+\frac{1}{12}\frac
{\text{L}\Psi}{\Psi}(y_{0})\frac{\partial^{2}}{\partial x_{i}^{2}}[\phi
\circ\pi_{\text{P}}](y_{0})$

$\qquad+\frac{1}{6}\frac{\partial}{\partial x_{i}}[\frac{\text{L}\Psi}{\Psi
}](y_{0})\frac{\partial}{\partial x_{i}}[\phi\circ\pi_{\text{P}}](y_{0})$

We note that by (i) and (ii) of $\left(  C_{8}\right)  :$ $\frac{\partial
}{\partial x_{i}}\pi_{\text{P}}(y_{0})=0=\frac{\partial^{2}}{\partial
x_{i}^{2}}\pi_{\text{P}}(y_{0})$ for $i=q+1,...,n.$

Hence,

\textbf{I}$_{321}=\frac{1}{12}\underset{i=q+1}{\overset{n}{\sum}}%
\frac{\partial^{2}}{\partial x_{i}^{2}}[\frac{\text{L}\Psi}{\Psi}](y_{0}%
)\phi(y_{0})$ which is given by (viii) of \textbf{Table B}$_{5}$ or from
$\left(  B_{118}\right)  .$

\qquad\qquad\qquad\qquad\qquad\qquad\qquad\qquad\qquad\qquad\qquad\qquad
\qquad\qquad\qquad\qquad\qquad\qquad$\blacksquare$

\textbf{I}$_{322}=\frac{1}{24}\underset{i=q+1}{\overset{n}{\sum}}%
\frac{\partial^{2}}{\partial x_{i}^{2}}[\underset{\text{a,b=1}%
}{\overset{\text{q}}{\sum}}$g$^{\text{ab}}\left\{  \frac{\partial^{2}\phi
}{\partial\text{x}_{\text{a}}\partial\text{x}_{\text{b}}}\circ\pi_{\text{P}%
}\text{ }\right\}  ](y_{0})$

\qquad\qquad$\mathbf{=}$ $\frac{1}{24}\underset{i=q+1}{\overset{n}{\sum}%
}\underset{\text{a,b=1}}{\overset{\text{q}}{\sum}}\frac{\partial^{2}%
\text{g}^{\text{ab}}}{\partial x_{i}^{2}}(y_{0})\left\{  \frac{\partial
^{2}\phi}{\partial\text{x}_{\text{a}}\partial\text{x}_{\text{b}}}\circ
\pi_{\text{P}}\right\}  (y_{0})=$ $\frac{1}{24}%
\underset{i=q+1}{\overset{n}{\sum}}\underset{\text{a,b=1}}{\overset{\text{q}%
}{\sum}}\frac{\partial^{2}\text{g}^{\text{ab}}}{\partial x_{i}^{2}}%
(y_{0})\frac{\partial^{2}\phi}{\partial\text{x}_{\text{a}}\partial
\text{x}_{\text{b}}}(y_{0})$

\qquad\qquad$\mathbf{+}\frac{1}{24}\underset{i=q+1}{\overset{n}{\sum}%
}\underset{\text{a,b=1}}{\overset{\text{q}}{\sum}}$g$^{\text{ab}}%
(y_{0})\left\{  \frac{\partial^{2}}{\partial x_{i}^{2}}[\frac{\partial^{2}%
\phi}{\partial\text{x}_{\text{a}}\partial\text{x}_{\text{b}}}\circ
\pi_{\text{P}}]\text{ }\right\}  (y_{0})=0$

$\qquad\qquad\mathbf{+}$ $\frac{1}{12}\underset{i=q+1}{\overset{n}{\sum}%
}\underset{\text{a,b=1}}{\overset{\text{q}}{\sum}}\frac{\partial
\text{g}^{\text{ab}}}{\partial x_{i}}(y_{0})\frac{\partial}{\partial x_{i}%
}[\frac{\partial^{2}\phi}{\partial\text{x}_{\text{a}}\partial\text{x}%
_{\text{b}}}\circ\pi_{\text{P}}]$ $(y_{0})=0$

$\qquad=$ $\frac{1}{24}\underset{i=q+1}{\overset{n}{\sum}}%
\underset{\text{a,b=1}}{\overset{\text{q}}{\sum}}\frac{\partial^{2}%
\text{g}^{\text{ab}}}{\partial x_{i}^{2}}(y_{0})\frac{\partial^{2}\phi
}{\partial\text{x}_{\text{a}}\partial\text{x}_{\text{b}}}(y_{0})$

The last two lines are zero because $\frac{\partial^{2}}{\partial x_{i}^{2}%
}\pi_{\text{P}}(y_{0})=0=\frac{\partial}{\partial x_{i}}\pi_{\text{P}}(y_{0})$
for $i=q+1,...,n$ by $\left(  C_{8}\right)  .$

We will use this property repeatedly without mentioning it explicitely.

By (iii) of \textbf{Table A}$_{6}$,

$\frac{\partial^{2}\text{g}^{\text{ab}}}{\partial\text{x}_{i}^{2}}%
(y_{0})=2[-R_{\text{a}i\text{b}i}+5\overset{q}{\underset{\text{c}=1}{\sum}%
}T_{\text{ac}i}T_{\text{bc}i}+2\overset{n}{\underset{j=q+1}{\sum}}%
\perp_{\text{a}ij}\perp_{\text{b}ij}](y_{0})$

Therefore we have:

$\left(  D_{10}\right)  \qquad$\textbf{I}$_{322}=\frac{1}{12}%
\underset{i=q+\text{1}}{\overset{\text{n}}{\sum}}\underset{\text{a,b=1}%
}{\overset{\text{q}}{\sum}}[-R_{\text{a}i\text{b}i}%
+5\overset{q}{\underset{\text{c}=1}{\sum}}T_{\text{ac}i}T_{\text{bc}%
i}+2\overset{n}{\underset{j=q+1}{\sum}}\perp_{\text{a}ij}\perp_{\text{b}%
ij}](y_{0})\times\frac{\partial^{2}\phi}{\partial\text{x}_{\text{a}}%
\partial\text{x}_{\text{b}}}(y_{0})$

\qquad\qquad\qquad\qquad\qquad\qquad\qquad\qquad\qquad\qquad\qquad\qquad
\qquad\qquad\qquad\qquad\qquad\qquad$\blacksquare$

(iii) Next we compute:

$\qquad$I$_{323}$ $\mathbf{=}$ $\frac{1}{24}\underset{i=q+1}{\overset{n}{\sum
}}\frac{\partial^{2}}{\partial x_{i}^{2}}[\underset{i,j=1}{\overset{n}{\sum}}%
$g$^{jk}\frac{\partial\Lambda_{k}}{\partial\text{x}_{j}}\phi\circ\pi
_{\text{P}}$ $](y_{0})$

$\qquad\mathbf{=}$ $\frac{1}{24}\underset{i=q+1}{\overset{n}{\sum}%
}\underset{j,k=1}{\overset{n}{\sum}}[\frac{\partial^{2}\text{g}^{jk}}{\partial
x_{i}^{2}}\frac{\partial\Lambda_{k}}{\partial\text{x}_{j}}\phi\circ
\pi_{\text{P}}$ $](y_{0})$ $\mathbf{+}$ $\frac{1}{24}%
\underset{i=q+1}{\overset{n}{\sum}}[\underset{j,k=1}{\overset{n}{\sum}}%
$g$^{jk}\frac{\partial^{2}}{\partial x_{i}^{2}}(\frac{\partial\Lambda_{k}%
}{\partial\text{x}_{j}}\phi\circ\pi_{\text{P}})$ ]$(y_{0})$

\qquad$\mathbf{+}\frac{1}{12}\underset{i=q+1}{\overset{n}{\sum}}%
\underset{j,k=1}{\overset{n}{\sum}}[\frac{\partial\text{g}^{jk}}{\partial
x_{i}}\frac{\partial}{\partial x_{i}}(\frac{\partial\Lambda_{k}}%
{\partial\text{x}_{j}}\phi\circ\pi_{\text{P}})$ $](y_{0})$

\qquad$=$ I$_{3231}+$ I$_{3232}$ + I$_{3233}$ where,

\ I$_{3231}\mathbf{=}$ $\frac{1}{24}\underset{i=q+1}{\overset{n}{\sum}%
}\underset{j,k=1}{\overset{n}{\sum}}[\frac{\partial^{2}\text{g}^{jk}}{\partial
x_{i}^{2}}\frac{\partial\Lambda_{k}}{\partial\text{x}_{j}}\phi\circ
\pi_{\text{P}}](y_{0});$\ \ \ I$_{3232}=\frac{1}{24}%
\underset{i=q+1}{\overset{n}{\sum}}[\underset{i,j=1}{\overset{n}{\sum}}%
$g$^{jk}\frac{\partial^{2}}{\partial x_{i}^{2}}(\frac{\partial\Lambda_{k}%
}{\partial\text{x}_{j}}\phi\circ\pi_{\text{P}})](y_{0})$

$\ \ \ \qquad$I$_{3233}=\frac{1}{12}\underset{i=q+1}{\overset{n}{\sum}%
}\underset{j,k=1}{\overset{n}{\sum}}[\frac{\partial\text{g}^{jk}}{\partial
x_{i}}\frac{\partial}{\partial x_{i}}(\frac{\partial\Lambda_{k}}%
{\partial\text{x}_{j}}\phi\circ\pi_{\text{P}})](y_{0})$

We will compute each of the above expressions in terms of invariants of the
manifold M, the submanifold P and the vector bundle E:

I$_{3231}\mathbf{=}$ $\frac{1}{24}\underset{j,k=1}{\overset{n}{\sum}}%
[\frac{\partial^{2}\text{g}^{jk}}{\partial x_{i}^{2}}\frac{\partial\Lambda
_{k}}{\partial\text{x}_{j}}\phi\circ\pi_{\text{P}}](y_{0})$

$\qquad=\frac{1}{24}\underset{i=1}{\overset{n}{\sum}}\underset{\text{a,b=1}%
}{\overset{\text{q}}{\sum}}\frac{\partial^{2}\text{g}^{\text{ab}}}{\partial
x_{i}^{2}}(y_{0})[\frac{\partial\Lambda_{\text{b}}}{\partial\text{x}%
_{\text{a}}}\phi\circ\pi_{\text{P}}](y_{0})+\frac{1}{24}%
\underset{i=q+1}{\overset{n}{\sum}}\underset{j,k=1}{\overset{n}{\sum}}%
\frac{\partial^{2}\text{g}^{jk}}{\partial x_{i}^{2}}(y_{0})[\frac
{\partial\Lambda_{k}}{\partial\text{x}_{j}}\phi\circ\pi_{\text{P}}](y_{0})$

$\qquad+\frac{1}{24}\underset{\text{a=1}}{\overset{\text{q}}{\sum}%
}\underset{i,k=q+1}{\overset{n}{\sum}}\frac{\partial^{2}\text{g}^{\text{a}k}%
}{\partial x_{i}^{2}}(y_{0})[\frac{\partial\Lambda_{k}}{\partial
\text{x}_{\text{a}}}\phi\circ\pi_{\text{P}}](y_{0})+\frac{1}{24}%
\underset{\text{b=1}}{\overset{\text{q}}{\sum}}%
\underset{i,j=q+1}{\overset{n}{\sum}}\frac{\partial^{2}\text{g}^{j\text{b}}%
}{\partial x_{i}^{2}}(y_{0})[\frac{\partial\Lambda_{\text{b}}}{\partial
\text{x}_{j}}\phi\circ\pi_{\text{P}}](y_{0})$

We have by $\left(  3.30\right)  ,$ $\frac{\partial\Lambda_{\text{b}}%
}{\partial\text{x}_{\text{a}}}(y_{0})=0=\frac{\partial\Lambda_{j}}%
{\partial\text{x}_{\text{a}}}(y_{0})\phi(y_{0})$.

Therefore we have:

I$_{3231}\mathbf{=}$ $\frac{1}{24}\underset{i,j,k=q+1}{\overset{n}{\sum}}%
\frac{\partial^{2}\text{g}^{jk}}{\partial x_{i}^{2}}(y_{0})[\frac
{\partial\Lambda_{k}}{\partial\text{x}_{j}}\phi\circ\pi_{\text{P}}%
](y_{0})+\frac{1}{24}\underset{\text{a=1}}{\overset{\text{q}}{\sum}%
}\underset{i,j\text{=q+1}}{\overset{\text{n}}{\sum}}\frac{\partial^{2}%
\text{g}^{j\text{a}}}{\partial x_{i}^{2}}(y_{0})[\frac{\partial\Lambda
_{\text{a}}}{\partial\text{x}_{j}}\phi\circ\pi_{\text{P}}](y_{0})$

$\frac{\partial^{2}\text{g}^{jk}}{\partial\text{x}_{i}^{2}}(y_{0})=\frac{2}%
{3}R_{jiki}(y_{0})=\frac{2}{3}R_{ijik}(y_{0})$ by (iii) of Table A$_{2},$

$\frac{\partial\Lambda_{k}}{\partial\text{x}_{j}}(y_{0})=\frac{1}{2}%
\Omega_{jk}(y_{0})$ by (x) of \textbf{Proposition 1.18, Berline, Getzler and
Vergne}

$\frac{\partial^{2}\text{g}^{\text{a}j}}{\partial\text{x}_{i}^{2}}%
(y_{0})=\frac{8}{3}R_{i\text{a}ij}(y_{0})+4\underset{\text{b=1}%
}{\overset{\text{q}}{\sum}}T_{\text{ab}i}(y_{0})\perp_{\text{b}ji}(y_{0})$
(iii) of \textbf{Table A}$_{4},$

$\frac{\partial\Lambda_{\text{a}}}{\partial x_{j}}=-\Omega_{\text{a}%
j}+[\Lambda_{\text{a}},\Lambda_{j}]$ for b = 1,...,q and $i=q+1,...,n$ by
(vii) of \textbf{Proposition 5 }above\textbf{.}

\textbf{ }Therefore we have:

$\left(  D_{11}\right)  \qquad$I$_{3231}=\frac{1}{72}%
\underset{i,j,k=q+1}{\overset{n}{\sum}}R_{ijik}(y_{0})\Omega_{jk}(y_{0}%
)\phi(y_{0})$

$+\frac{1}{24}\underset{\text{a=1}}{\overset{\text{q}}{\sum}}%
\underset{i,j\text{=q+1}}{\overset{\text{n}}{\sum}}[\frac{8}{3}R_{i\text{a}%
ij}+4\underset{\text{b=1}}{\overset{\text{q}}{\sum}}T_{\text{ab}i}%
\perp_{\text{b}ji}](y_{0})\left\{  -\Omega_{\text{a}j}+[\Lambda_{\text{a}%
},\Lambda_{j}]\right\}  (y_{0})\phi(y_{0})$

\qquad\qquad\qquad\qquad\qquad\qquad\qquad\qquad\qquad\qquad\qquad\qquad
\qquad\qquad\qquad\qquad\qquad\qquad$\blacksquare$

We next consider:

I$_{3232}=\frac{1}{24}\underset{i=q+1}{\overset{n}{\sum}}%
[\underset{j,k=1}{\overset{n}{\sum}}$g$^{jk}\frac{\partial^{2}}{\partial
x_{i}^{2}}(\frac{\partial\Lambda_{k}}{\partial\text{x}_{j}}\phi\circ
\pi_{\text{P}})](y_{0})=\frac{1}{24}\underset{i=q+1}{\overset{n}{\sum}%
}[\underset{j,k=1}{\overset{n}{\sum}}$g$^{jk}\frac{\partial^{3}\Lambda_{k}%
}{\partial x_{i}^{2}\partial\text{x}_{j}}\phi\circ\pi_{\text{P}}](y_{0})$

Since g$^{jk}(y_{0})=\delta^{jk},$

\ \ \ I$_{3232}=\frac{1}{24}\underset{i=q+1}{\overset{n}{\sum}}%
[\underset{j=1}{\overset{n}{\sum}}\frac{\partial^{3}\Lambda_{j}}{\partial
x_{i}^{2}\partial\text{x}_{j}}\phi](y_{0})=\frac{1}{24}%
\underset{i=q+1}{\overset{n}{\sum}}[\underset{\text{a}=1}{\overset{q}{\sum}%
}\frac{\partial^{3}\Lambda_{\text{a}}}{\partial x_{i}^{2}\partial
\text{x}_{\text{a}}}\phi](y_{0})+\frac{1}{24}%
[\underset{i,j=q+1}{\overset{n}{\sum}}\frac{\partial^{3}\Lambda_{j}}{\partial
x_{i}^{2}\partial\text{x}_{j}}\phi](y_{0})$

Since all differentiation with respect to tangential Fermi coordinates vanish,
we have:

I$_{3232}=\frac{1}{24}\underset{i=q+1}{\overset{n}{\sum}}%
[\underset{j=1}{\overset{n}{\sum}}\frac{\partial^{3}\Lambda_{j}}{\partial
x_{i}^{2}\partial\text{x}_{j}}\phi](y_{0})=\frac{1}{24}%
[\underset{i,j=q+1}{\overset{n}{\sum}}\frac{\partial^{3}\Lambda_{j}}{\partial
x_{i}^{2}\partial\text{x}_{j}}](y_{0})\phi(y_{0})$

$\frac{\partial^{3}\Lambda_{j}}{\partial\text{x}_{i}^{2}\partial\text{x}_{j}%
}(y_{0})=\frac{1}{4}\frac{\partial^{2}\Omega_{jj}}{\partial\text{x}_{i}^{2}%
}(y_{0})$

Since $\Omega_{ij}$ is skew-symmetric in the indices $\left(  i,j\right)  ,$
we have $\Omega_{jj}=0$ and so,

$\left(  D_{12}\right)  $\qquad I$_{3232}=\frac{1}{24}%
\underset{i,j=q+1}{\overset{n}{\sum}}\frac{\partial^{3}\Lambda_{j}}{\partial
x_{i}^{2}\partial\text{x}_{j}}(y_{0})\phi(y_{0})=\frac{1}{4}\frac{\partial
^{2}\Omega_{jj}}{\partial\text{x}_{i}^{2}}(y_{0})=0$

\qquad\qquad\qquad\qquad\qquad\qquad\qquad\qquad\qquad\qquad\qquad\qquad
\qquad\qquad$\blacksquare$

We next consider:

\qquad\qquad\ \ I$_{3233}=\frac{1}{12}\underset{i=q+1}{\overset{n}{\sum}%
}\underset{j,k=1}{\overset{n}{\sum}}\left\{  \frac{\partial\text{g}^{jk}%
}{\partial x_{i}}\frac{\partial}{\partial x_{i}}(\frac{\partial\Lambda_{k}%
}{\partial\text{x}_{j}}\phi\circ\pi_{\text{P}})\right\}  (y_{0})$

\qquad\qquad\qquad$\ \ \ =\frac{1}{12}\underset{i=q+1}{\overset{n}{\sum}%
}\underset{\text{a,b=1}}{\overset{\text{q}}{\sum}}\left\{  \frac
{\partial\text{g}^{\text{ab}}}{\partial x_{i}}\frac{\partial^{2}%
\Lambda_{\text{b}}}{\partial x_{i}\partial\text{x}_{\text{a}}}\phi\circ
\pi_{\text{P}}\right\}  (y_{0})$

\qquad\qquad$\qquad\ \ \ \ +\frac{1}{12}\underset{i=q+1}{\overset{n}{\sum}%
}\underset{j,k\text{=q+1}}{\overset{\text{n}}{\sum}}\left\{  \frac
{\partial\text{g}^{jk}}{\partial x_{i}}\frac{\partial^{2}\Lambda_{k}}{\partial
x_{i}\partial\text{x}_{j}}\phi\circ\pi_{\text{P}}\right\}  (y_{0})$

\qquad\qquad$\qquad\ \ \ \ +\frac{1}{6}\underset{i=q+1}{\overset{n}{\sum}%
}\underset{k=q+1}{\overset{n}{\sum}}\underset{\text{a=1}}{\overset{\text{q}%
}{\sum}}\left\{  \frac{\partial\text{g}^{\text{a}k}}{\partial x_{i}}%
\frac{\partial^{2}\Lambda_{k}}{\partial x_{i}\partial\text{x}_{\text{a}}}%
\phi\circ\pi_{\text{P}}\right\}  (y_{0})$

By (ii) of \textbf{Table A}$_{2}$\textbf{ }in \textbf{Appendix A},
$\frac{\partial\text{g}^{jk}}{\partial x_{i}}(y_{0})=0$ for $i,j,k=q+1,...,n$

Next since all differentiation with respect to tangential coordinate variables
vanish, we have:

$\qquad\frac{\partial^{2}\Lambda_{\text{b}}}{\partial x_{i}\partial
\text{x}_{\text{a}}}(y_{0})=0=\frac{\partial^{2}\Lambda_{k}}{\partial
x_{i}\partial\text{x}_{\text{a}}}(y_{0})$

We conclude that:

$\left(  D_{13}\right)  \qquad$I$_{3233}=0$

\qquad\qquad\qquad\qquad\qquad\qquad\qquad$\blacksquare$

We conclude by $\left(  D_{11}\right)  ,$ $\left(  D_{12}\right)  ,$ and
$\left(  D_{13}\right)  $ that:

\qquad\ \ I$_{323}=$ I$_{3231}+$ I$_{3232}+$ I$_{3233}$

$\left(  D_{14}\right)  \qquad$\ I$_{323}=\frac{1}{72}%
\underset{i,j,k=q+1}{\overset{n}{\sum}}R_{ijik}(y_{0})\Omega_{jk}(y_{0}%
)\phi(y_{0})$

$\qquad+\frac{1}{24}\underset{\text{a=1}}{\overset{\text{q}}{\sum}%
}\underset{i,j=q+1}{\overset{n}{\sum}}\left\{  \frac{8}{3}R_{i\text{a}%
ij}+4\underset{\text{b=1}}{\overset{\text{q}}{\sum}}T_{\text{ab}i}%
\perp_{\text{b}ji}\right\}  (y_{0})\left\{  -\Omega_{\text{a}j}+[\Lambda
_{\text{a}},\Lambda_{j}]\right\}  (y_{0})\phi(y_{0})$

\qquad\qquad\qquad\qquad\qquad\qquad\qquad\qquad\qquad\qquad\qquad\qquad
\qquad\qquad\qquad\qquad\qquad\qquad$\blacksquare$

We next consider:

(iv)\qquad\ \ I$_{324}$ $\mathbf{=}$ $\frac{1}{12}%
\underset{i=q+1}{\overset{n}{\sum}}\frac{\partial^{2}}{\partial x_{i}^{2}%
}\left\{  \text{ }\underset{j=1}{\overset{n}{\sum}}\underset{\text{a=1}%
}{\overset{\text{q}}{\sum}}\text{g}^{\text{a}j}\Lambda_{j}\frac{\partial\phi
}{\partial\text{x}_{\text{a}}}\circ\pi_{\text{P}}\right\}  (y_{0})$

$=$ $\frac{1}{12}\underset{i=q+1}{\overset{n}{\sum}}\left\{  \text{
}\underset{j=1}{\overset{n}{\sum}}\underset{\text{a=1}}{\overset{\text{q}%
}{\sum}}\frac{\partial^{2}\text{g}^{\text{a}j}}{\partial x_{i}^{2}}\Lambda
_{j}\frac{\partial\phi}{\partial\text{x}_{\text{a}}}\circ\pi_{\text{P}%
}\right\}  (y_{0})\mathbf{+}$ $\frac{1}{12}\underset{i=q+1}{\overset{n}{\sum}%
}\left\{  \text{ }\underset{j=1}{\overset{n}{\sum}}\underset{\text{a=1}%
}{\overset{\text{q}}{\sum}}\text{g}^{\text{a}j}\frac{\partial^{2}}{\partial
x_{i}^{2}}(\Lambda_{j}\frac{\partial\phi}{\partial\text{x}_{\text{a}}}\circ
\pi_{\text{P}})\right\}  (y_{0})$

$\mathbf{+}$ $\frac{1}{6}\underset{i=q+1}{\overset{n}{\sum}}\left\{  \text{
}\underset{j=1}{\overset{n}{\sum}}\underset{\text{a=1}}{\overset{\text{q}%
}{\sum}}\frac{\partial\text{g}^{\text{a}j}}{\partial x_{i}}\frac{\partial
}{\partial x_{i}}(\Lambda_{j}\frac{\partial\phi}{\partial\text{x}_{\text{a}}%
}\circ\pi_{\text{P}})\right\}  (y_{0})$

$=$ I$_{3241}+$ I$_{3242}+$ I$_{3243}$ where,

I$_{3241}=\frac{1}{12}\underset{i=q+1}{\overset{n}{\sum}}\left\{  \text{
}\underset{j=1}{\overset{n}{\sum}}\underset{\text{a=1}}{\overset{\text{q}%
}{\sum}}\frac{\partial^{2}\text{g}^{\text{a}j}}{\partial x_{i}^{2}}\Lambda
_{j}\frac{\partial\phi}{\partial\text{x}_{\text{a}}}\circ\pi_{\text{P}%
}\right\}  (y_{0})$

I$_{3242}$ $\ \mathbf{=}$ $\frac{1}{12}\underset{i=q+1}{\overset{n}{\sum}%
}\left\{  \text{ }\underset{j=1}{\overset{n}{\sum}}\underset{\text{a=1}%
}{\overset{\text{q}}{\sum}}\text{g}^{\text{a}j}\frac{\partial^{2}}{\partial
x_{i}^{2}}(\Lambda_{j}\frac{\partial\phi}{\partial\text{x}_{\text{a}}}\circ
\pi_{\text{P}})\right\}  (y_{0})$

I$_{3243}=\frac{1}{6}\underset{i=q+1}{\overset{n}{\sum}}\left\{  \text{
}\underset{j=1}{\overset{n}{\sum}}\underset{\text{a=1}}{\overset{\text{q}%
}{\sum}}\frac{\partial\text{g}^{\text{a}j}}{\partial x_{i}}\frac{\partial
}{\partial x_{i}}(\Lambda_{j}\frac{\partial\phi}{\partial\text{x}_{\text{a}}%
}\circ\pi_{\text{P}})\right\}  (y_{0})$

We express the above expressions in terms of geometric invariants:

\ \qquad I$_{3241}=\frac{1}{12}\underset{i=q+1}{\overset{n}{\sum}}\left\{
\text{ }\underset{j=1}{\overset{n}{\sum}}\underset{\text{a=1}%
}{\overset{\text{q}}{\sum}}\frac{\partial^{2}\text{g}^{\text{a}j}}{\partial
x_{i}^{2}}\Lambda_{j}\frac{\partial\phi}{\partial\text{x}_{\text{a}}}\circ
\pi_{\text{P}}\right\}  (y_{0})$

$=\frac{1}{12}\underset{i=q+1}{\overset{n}{\sum}}\left\{  \text{
}\underset{\text{b=1}}{\overset{\text{q}}{\sum}}\underset{\text{a=1}%
}{\overset{\text{q}}{\sum}}\frac{\partial^{2}\text{g}^{\text{ab}}}{\partial
x_{i}^{2}}\Lambda_{\text{b}}\frac{\partial\phi}{\partial\text{x}_{\text{a}}%
}\circ\pi_{\text{P}}\right\}  (y_{0})+\frac{1}{12}\left\{  \text{
}\underset{j=q+1}{\overset{n}{\sum}}\underset{\text{a=1}}{\overset{\text{q}%
}{\sum}}\frac{\partial^{2}\text{g}^{\text{a}j}}{\partial x_{i}^{2}}\Lambda
_{j}\frac{\partial\phi}{\partial\text{x}_{\text{a}}}\circ\pi_{\text{P}%
}\right\}  (y_{0})$

By (iii) of \textbf{Table A}$_{6},$

$\frac{1}{12}\underset{i=q+1}{\overset{n}{\sum}}\left\{  \text{ }%
\underset{\text{b=1}}{\overset{\text{q}}{\sum}}\underset{\text{a=1}%
}{\overset{\text{q}}{\sum}}\frac{\partial^{2}\text{g}^{\text{ab}}}{\partial
x_{i}^{2}}\Lambda_{\text{b}}\frac{\partial\phi}{\partial\text{x}_{\text{a}}%
}\circ\pi_{\text{P}}\right\}  (y_{0})$

$=\frac{1}{6}\underset{i=q+1}{\overset{n}{\sum}}$ $\underset{\text{a,b=1}%
}{\overset{\text{q}}{\sum}}\left\{  -R_{\text{a}i\text{b}i}%
+5\overset{q}{\underset{\text{c}=1}{\sum}}T_{\text{ac}i}T_{\text{bc}%
i}+2\overset{n}{\underset{\text{k}=q+1}{\sum}}\perp_{\text{a}i\text{k}}%
\perp_{\text{b}i\text{k}}\right\}  (y_{0})\times\Lambda_{\text{b}}(y_{0}%
)\frac{\partial\phi}{\partial\text{x}_{\text{a}}}(y_{0})$

By (iii) of Table A$_{4},$ $\frac{\partial^{2}\text{g}^{\text{a}j}}%
{\partial\text{x}_{i}^{2}}(y_{0})=\frac{8}{3}R_{i\text{a}ij}%
+4\underset{\text{b=1}}{\overset{\text{q}}{\sum}}T_{\text{ab}i}(y_{0}%
)\perp_{\text{b}ji}(y_{0})$ and so, we have:

$\frac{1}{12}\left\{  \text{ }\underset{j=q+1}{\overset{n}{\sum}%
}\underset{\text{a=1}}{\overset{\text{q}}{\sum}}\frac{\partial^{2}%
\text{g}^{\text{a}j}}{\partial x_{i}^{2}}\Lambda_{j}\frac{\partial\phi
}{\partial\text{x}_{\text{a}}}\circ\pi_{\text{P}}\right\}  (y_{0})$

$=\frac{1}{12}\underset{j=q+1}{\overset{n}{\sum}}\underset{\text{a=1}%
}{\overset{\text{q}}{\sum}}\left\{  \text{ }\frac{8}{3}R_{i\text{a}%
ij}+4\underset{\text{b=1}}{\overset{\text{q}}{\sum}}T_{\text{ab}i}(y_{0}%
)\perp_{\text{b}ji}(y_{0})\right\}  \Lambda_{j}(y_{0})\frac{\partial\phi
}{\partial\text{x}_{\text{a}}}(y_{0})$

We conclude here that:

$\left(  D_{15}\right)  \qquad$I$_{3241}$

$=\frac{1}{6}\underset{i=q+1}{\overset{n}{\sum}}$ $\underset{\text{a,b=1}%
}{\overset{\text{q}}{\sum}}\left\{  -R_{\text{a}i\text{b}i}%
+5\overset{q}{\underset{\text{c}=1}{\sum}}T_{\text{ac}i}T_{\text{bc}%
i}+2\overset{n}{\underset{\text{k}=q+1}{\sum}}\perp_{\text{a}i\text{k}}%
\perp_{\text{b}i\text{k}}\right\}  (y_{0})\times\Lambda_{\text{b}}(y_{0}%
)\frac{\partial\phi}{\partial\text{x}_{\text{a}}}(y_{0})$

$\qquad\qquad+\frac{1}{12}\underset{j=q+1}{\overset{n}{\sum}}%
\underset{\text{a=1}}{\overset{\text{q}}{\sum}}\left\{  \text{ }\frac{8}%
{3}R_{i\text{a}ij}+4\underset{\text{b=1}}{\overset{\text{q}}{\sum}%
}T_{\text{ab}i}(y_{0})\perp_{\text{b}ji}(y_{0})\right\}  \Lambda_{j}%
(y_{0})\frac{\partial\phi}{\partial\text{x}_{\text{a}}}(y_{0})$

\qquad\qquad\qquad\qquad\qquad\qquad\qquad\qquad\qquad\qquad\qquad\qquad
\qquad\qquad\qquad\qquad\qquad\qquad\qquad$\blacksquare$

Next we consider:

\qquad I$_{3242}$ $\ \mathbf{=}$ $\frac{1}{24}%
\underset{i=q+1}{\overset{n}{\sum}}\left\{  \text{ }%
\underset{j=1}{\overset{n}{\sum}}\underset{\text{a=1}}{\overset{\text{q}%
}{\sum}}\text{g}^{\text{a}j}\frac{\partial^{2}}{\partial x_{i}^{2}}%
(\Lambda_{j}\frac{\partial\phi}{\partial\text{x}_{\text{a}}}\circ\pi
_{\text{P}})\right\}  (y_{0})$

$\qquad\mathbf{=}$ $\frac{1}{24}\underset{i=q+1}{\overset{n}{\sum}}\left\{
\text{ }\underset{j=1}{\overset{n}{\sum}}\underset{\text{a=1}%
}{\overset{\text{q}}{\sum}}\text{g}^{\text{a}j}\frac{\partial^{2}\Lambda_{j}%
}{\partial x_{i}^{2}}\frac{\partial\phi}{\partial\text{x}_{\text{a}}}\circ
\pi_{\text{P}}\right\}  (y_{0})$

\qquad$\mathbf{=}$ $\frac{1}{24}\underset{i=q+1}{\overset{n}{\sum}}\left\{
\text{ }\underset{\text{b=1}}{\overset{\text{q}}{\sum}}\underset{\text{a=1}%
}{\overset{\text{q}}{\sum}}\text{g}^{\text{ab}}\frac{\partial^{2}%
\Lambda_{\text{b}}}{\partial x_{i}^{2}}\frac{\partial\phi}{\partial
\text{x}_{\text{a}}}\circ\pi_{\text{P}}\right\}  (y_{0})$

\qquad$\mathbf{+}$ $\frac{1}{24}\underset{i=q+\text{1}}{\overset{n}{\sum}%
}\left\{  \text{ }\underset{j=q+1}{\overset{n}{\sum}}\underset{\text{a=1}%
}{\overset{\text{q}}{\sum}}\text{g}^{\text{a}j}\frac{\partial^{2}\Lambda_{j}%
}{\partial x_{i}^{2}}\frac{\partial\phi}{\partial\text{x}_{\text{a}}}\circ
\pi_{\text{P}}\right\}  (y_{0})$

Since g$^{\text{ab}}(y_{0})=\delta^{\text{ab}}$ and g$^{\text{a}j}(y_{0})=0$
for a,b = 1,...,q and $j=q+1,...,n$

I$_{3242}$ $\ \mathbf{=}$ $\frac{1}{24}\underset{i=q+1}{\overset{n}{\sum}}[$
$\underset{\text{a=1}}{\overset{\text{q}}{\sum}}\frac{\partial^{2}%
\Lambda_{\text{a}}}{\partial x_{i}^{2}}\frac{\partial\phi}{\partial
\text{x}_{\text{a}}}](y_{0})$

Computations in (viii) of \textbf{Proposition 5} give:

$\frac{\partial^{2}\Lambda_{\text{a}}}{\partial x_{i}\partial x_{j}}%
(y_{0})=\frac{\partial\Omega_{j\text{a}}}{\partial x_{i}}(y_{0})+\left(
-\Omega_{\text{a}i}+\Lambda_{\text{a}}\Lambda_{i}-\Lambda_{i}\Lambda
_{\text{a}}\right)  (y_{0})\Lambda_{j}(y_{0})$\qquad$\qquad$

$-\Lambda_{j}(y_{0})\left(  -\Omega_{\text{a}i}+\Lambda_{\text{a}}\Lambda
_{i}-\Lambda_{i}\Lambda_{\text{a}}\right)  (y_{0})$

$+\frac{1}{2}\Lambda_{\text{a}}(y_{0})\Omega_{ij}(y_{0})-\frac{1}{2}%
\Omega_{ij}(y_{0})\Lambda_{\text{a}}(y_{0}).$

Therefore

$\frac{\partial^{2}\Lambda_{\text{a}}}{\partial x_{i}^{2}}(y_{0}%
)=\frac{\partial\Omega_{i\text{a}}}{\partial x_{i}}(y_{0})+\left(
-\Omega_{\text{a}i}+\Lambda_{\text{a}}\Lambda_{i}-\Lambda_{i}\Lambda
_{\text{a}}\right)  (y_{0})\Lambda_{i}(y_{0})$

$-\Lambda_{i}(y_{0})\left(  -\Omega_{\text{a}i}+\Lambda_{\text{a}}\Lambda
_{i}-\Lambda_{i}\Lambda_{\text{a}}\right)  (y_{0})+\frac{1}{2}\Lambda
_{\text{a}}(y_{0})\Omega_{ii}(y_{0})-\frac{1}{2}\Omega_{ii}(y_{0}%
)\Lambda_{\text{a}}(y_{0}).$

Since $\Omega_{ij}(y_{0})$ is skew-symmetric in the indices $\left(
i,j\right)  $ we have: $\Omega_{ii}(y_{0})=0$ and so,

$\frac{\partial^{2}\Lambda_{\text{a}}}{\partial x_{i}^{2}}(y_{0}%
)=\frac{\partial\Omega_{i\text{a}}}{\partial x_{i}}(y_{0})+\left(
-\Omega_{\text{a}i}+\Lambda_{\text{a}}\Lambda_{i}-\Lambda_{i}\Lambda
_{\text{a}}\right)  (y_{0})\Lambda_{i}(y_{0})-\Lambda_{i}(y_{0})\left(
-\Omega_{\text{a}i}+\Lambda_{\text{a}}\Lambda_{i}-\Lambda_{i}\Lambda
_{\text{a}}\right)  (y_{0})$

\qquad$=\frac{\partial\Omega_{i\text{a}}}{\partial x_{i}}(y_{0})+\left(
-\Omega_{\text{a}i}\Lambda_{i}+\Lambda_{\text{a}}\Lambda_{i}\Lambda
_{i}-\Lambda_{i}\Lambda_{\text{a}}\Lambda_{i}\right)  (y_{0})+\left(
\Lambda_{i}\Omega_{\text{a}i}-\Lambda_{i}\Lambda_{\text{a}}\Lambda_{i}%
+\Lambda_{i}\Lambda_{i}\Lambda_{\text{a}}\right)  (y_{0})$

Since $\Lambda_{i}^{2}(y_{0})=0,$ we have:

$\qquad\frac{\partial^{2}\Lambda_{\text{a}}}{\partial x_{i}^{2}}(y_{0}%
)=\frac{\partial\Omega_{i\text{a}}}{\partial x_{i}}(y_{0})+\left(  \Lambda
_{i}\Omega_{\text{a}i}-\Omega_{\text{a}i}\Lambda_{i}-2\Lambda_{i}%
\Lambda_{\text{a}}\Lambda_{i}\right)  (y_{0})$

$\qquad$I$_{3242}$ $\ \mathbf{=}$ $\frac{1}{24}%
\underset{i=q+1}{\overset{n}{\sum}}\underset{\text{a=1}}{\overset{\text{q}%
}{\sum}}\left[  \frac{\partial\Omega_{i\text{a}}}{\partial x_{i}}+\Lambda
_{i}\Omega_{\text{a}i}-\Omega_{\text{a}i}\Lambda_{i}-2\Lambda_{i}%
\Lambda_{\text{a}}\Lambda_{i}\right]  (y_{0})\frac{\partial\phi}%
{\partial\text{x}_{\text{a}}}(y_{0})$

\qquad\qquad\qquad\qquad\qquad\qquad\qquad\qquad\qquad\qquad\qquad\qquad
\qquad\qquad\qquad\qquad$\blacksquare$

We consider the last expression here:

I$_{3243}=\frac{1}{12}\underset{i=q+1}{\overset{n}{\sum}}\left\{  \text{
}\underset{j=1}{\overset{n}{\sum}}\underset{\text{a=1}}{\overset{\text{q}%
}{\sum}}\frac{\partial\text{g}^{\text{a}j}}{\partial x_{i}}\frac
{\partial\Lambda_{j}}{\partial x_{i}}\frac{\partial\phi}{\partial
\text{x}_{\text{a}}}\circ\pi_{\text{P}}\right\}  (y_{0})$

\qquad$=\frac{1}{12}\underset{i=q+1}{\overset{n}{\sum}}\left\{  \text{
}\underset{\text{b=1}}{\overset{\text{q}}{\sum}}\underset{\text{a=1}%
}{\overset{\text{q}}{\sum}}\frac{\partial\text{g}^{\text{ab}}}{\partial x_{i}%
}\frac{\partial\Lambda_{\text{b}}}{\partial x_{i}}\frac{\partial\phi}%
{\partial\text{x}_{\text{a}}}\circ\pi_{\text{P}}\right\}  (y_{0})$

\qquad$+\frac{1}{12}\underset{i=q+1}{\overset{n}{\sum}}\left\{  \text{
}\underset{j=q+1}{\overset{n}{\sum}}\underset{\text{a=1}}{\overset{\text{q}%
}{\sum}}\frac{\partial\text{g}^{\text{a}j}}{\partial x_{i}}\frac
{\partial\Lambda_{j}}{\partial x_{i}}\frac{\partial\phi}{\partial
\text{x}_{\text{a}}}\circ\pi_{\text{P}}\right\}  (y_{0})$

By (ii) of \textbf{Table A}$_{6}$, $\frac{\partial\text{g}^{\text{ab}}%
}{\partial x_{i}}(y_{0})=2$T$_{\text{ab}i}$ and (ii) of \textbf{Table A}$_{4}%
$, $\frac{\partial\text{g}^{\text{a}j}}{\partial x_{i}}(y_{0})=\perp
_{\text{a}ji}(y_{0})$

and from (vii) of \textbf{Proposition 5}, we have: $\ \frac{\partial
\Lambda_{\text{a}}}{\partial x_{i}}=-\Omega_{\text{a}i}+[\Lambda_{\text{a}%
},\Lambda_{i}]$

and since $\frac{\partial\Lambda_{j}}{\partial\text{x}_{i}}(y_{0})=\frac{1}%
{2}\Omega_{ij}(y_{0}),$ we have:

$\left(  D_{17}\right)  \qquad$I$_{3243}=\frac{1}{6}%
\underset{i,j=q+1}{\overset{n}{\sum}}\underset{\text{a=1}}{\overset{\text{q}%
}{\sum}}T_{\text{ab}i}(y_{0})\left(  -\Omega_{\text{a}i}+[\Lambda_{\text{a}%
},\Lambda_{i}]\right)  (y_{0})\frac{\partial\phi}{\partial\text{x}_{\text{a}}%
}(y_{0})$

$\qquad\qquad+\frac{1}{24}$ $\underset{i,j=q+1}{\overset{n}{\sum}%
}\underset{\text{a=1}}{\overset{\text{q}}{\sum}}\perp_{\text{a}ji}%
(y_{0})\Omega_{ij}(y_{0})\frac{\partial\phi}{\partial\text{x}_{\text{a}}%
}(y_{0})\qquad$

\qquad\qquad\qquad\qquad\qquad\qquad\qquad\qquad\qquad\qquad\qquad\qquad
\qquad\qquad\qquad$\blacksquare$

We now gather all the items that constitute I$_{324}=$ I$_{3241}+$ I$_{3242}+$
I$_{3243}$

from $\left(  D_{15}\right)  ,\left(  D_{16}\right)  ,\left(  D_{26}\right)  $
to give:

$\left(  D_{18}\right)  \qquad$I$_{324}$

$=\frac{1}{6}\underset{i=q+1}{\overset{n}{\sum}}$ $\underset{\text{a,b=1}%
}{\overset{\text{q}}{\sum}}\left\{  -R_{\text{a}i\text{b}i}%
+5\overset{q}{\underset{\text{c}=1}{\sum}}T_{\text{ac}i}T_{\text{bc}%
i}+2\overset{n}{\underset{\text{k}=q+1}{\sum}}\perp_{\text{a}i\text{k}}%
\perp_{\text{b}i\text{k}}\right\}  (y_{0})\times\Lambda_{\text{b}}(y_{0}%
)\frac{\partial\phi}{\partial\text{x}_{\text{a}}}(y_{0})\qquad$I$_{3241}$

$+\frac{1}{12}\underset{j=q+1}{\overset{n}{\sum}}\underset{\text{a=1}%
}{\overset{\text{q}}{\sum}}\left\{  \text{ }\frac{8}{3}R_{i\text{a}%
ij}+4\underset{\text{b=1}}{\overset{\text{q}}{\sum}}T_{\text{ab}i}(y_{0}%
)\perp_{\text{b}ji}(y_{0})\right\}  \Lambda_{j}(y_{0})\frac{\partial\phi
}{\partial\text{x}_{\text{a}}}(y_{0})$

$+$ $\frac{1}{24}\underset{i=q+1}{\overset{n}{\sum}}\underset{\text{a=1}%
}{\overset{\text{q}}{\sum}}\left[  \frac{\partial\Omega_{i\text{a}}}{\partial
x_{i}}+\Lambda_{i}\Omega_{\text{a}i}-\Omega_{\text{a}i}\Lambda_{i}%
-2\Lambda_{i}\Lambda_{\text{a}}\Lambda_{i}\right]  (y_{0})\frac{\partial\phi
}{\partial\text{x}_{\text{a}}}(y_{0})\qquad$I$_{3242}$

$+\frac{1}{6}\underset{i,j=q+1}{\overset{n}{\sum}}\underset{\text{a=1}%
}{\overset{\text{q}}{\sum}}T_{\text{ab}i}(y_{0})\left(  -\Omega_{\text{a}%
i}+[\Lambda_{\text{a}},\Lambda_{i}]\right)  (y_{0})\frac{\partial\phi
}{\partial\text{x}_{\text{a}}}(y_{0})\qquad\qquad\qquad$I$_{3243}$

$-\frac{1}{24}$ $\underset{i,j=q+1}{\overset{n}{\sum}}\underset{\text{a=1}%
}{\overset{\text{q}}{\sum}}\perp_{\text{a}ij}(y_{0})\Omega_{ij}(y_{0}%
)\frac{\partial\phi}{\partial\text{x}_{\text{a}}}(y_{0})\qquad$

\qquad\qquad\qquad\qquad\qquad\qquad\qquad\qquad\qquad\qquad\qquad\qquad
\qquad\qquad\qquad\qquad\qquad$\blacksquare$

(v) We next compute the expression for I$_{325}:$

\qquad I$_{325}$ $=\frac{1}{24}\frac{\partial^{2}}{\partial x_{i}^{2}%
}[\underset{\text{i,j=1}}{\overset{\text{n}}{\sum}}$g$^{jk}\Lambda_{j}%
\Lambda_{k}\phi\circ\pi_{\text{P}}](y_{0})$

Again here we will repeatedly use the fact given in (i) and (ii) of $\left(
C_{8}\right)  $ (without explicitely mentioning it) that:

$\frac{\partial}{\partial\text{x}_{i}}\pi_{\text{P}}$(z$_{0}$) $=\left\{
\begin{array}
[c]{c}%
1\text{ for }i=1,...,q\\
0\text{ for }i=q+1,...,n
\end{array}
\right.  $ and $\frac{\partial^{2}}{\partial\text{x}_{j}\partial\text{x}_{i}%
}\pi_{\text{P}}$(z$_{0}$) $=0$ for all $i,j=1,...,q,q+1,...,n.$

\qquad\ \ I$_{325}$ $\mathbf{=}$ $\frac{1}{24}%
\underset{i=q+1}{\overset{n}{\sum}}\frac{\partial^{2}}{\partial x_{i}^{2}%
}\left\{  \text{ }\underset{j,k\text{=1}}{\overset{\text{n}}{\sum}}%
\text{g}^{jk}\Lambda_{j}\Lambda_{k}\phi\circ\pi_{\text{P}}\right\}  (y_{0})$

\qquad$=$ $\frac{1}{24}[\underset{i=q+1}{\overset{n}{\sum}}%
\underset{j,k=1}{\overset{n}{\sum}}\frac{\partial^{2}\text{g}^{jk}}{\partial
x_{i}^{2}}\left\{  \text{ }\Lambda_{j}\Lambda_{k}\phi\circ\pi_{\text{P}%
}\right\}  ](y_{0})$

\qquad$+$ $\frac{1}{24}[\underset{j=q+1}{\overset{n}{\sum}}%
\underset{j,k=1}{\overset{n}{\sum}}$g$^{jk}\frac{\partial^{2}}{\partial
x_{i}^{2}}\left\{  \text{ }\Lambda_{j}\Lambda_{k}\phi\circ\pi_{\text{P}%
}\right\}  ](y_{0})$

\qquad$+$ $\frac{1}{12}[\underset{i=q+1}{\overset{n}{\sum}}%
\underset{j,k=1}{\overset{n}{\sum}}\frac{\partial\text{g}^{jk}}{\partial
x_{i}}\frac{\partial}{\partial x_{i}}\left\{  \text{ }\Lambda_{j}\Lambda
_{k}\phi\circ\pi_{\text{P}}\right\}  ](y_{0})$

\qquad$=$ I$_{3251}+$ I$_{3252}$ + I$_{3253}$ where,

I$_{3251}=\frac{1}{24}[\underset{i=q+1}{\overset{n}{\sum}}%
\underset{j,k=1}{\overset{n}{\sum}}\frac{\partial^{2}\text{g}^{jk}}{\partial
x_{i}^{2}}\left\{  \text{ }\Lambda_{j}\Lambda_{k}\phi\circ\pi_{\text{P}%
}\right\}  ](y_{0})$

I$_{3252}=$ $\frac{1}{24}[\underset{i=q+1}{\overset{n}{\sum}}%
\underset{j,k=1}{\overset{n}{\sum}}$g$^{jk}\frac{\partial^{2}}{\partial
x_{i}^{2}}\left\{  \text{ }\Lambda_{j}\Lambda_{k}\phi\circ\pi_{\text{P}%
}\right\}  ](y_{0})$

I$_{3253}=\frac{1}{12}[\underset{i=q+1}{\overset{n}{\sum}}%
\underset{j,k=1}{\overset{n}{\sum}}\frac{\partial\text{g}^{jk}}{\partial
x_{i}}\frac{\partial}{\partial x_{i}}\left\{  \text{ }\Lambda_{j}\Lambda
_{k}\phi\circ\pi_{\text{P}}\right\}  ](y_{0})$

We start with:

I$_{3251}=\frac{1}{24}[\underset{i=q+1}{\overset{n}{\sum}}%
\underset{j,k=1}{\overset{n}{\sum}}\frac{\partial^{2}\text{g}^{jk}}{\partial
x_{i}^{2}}\left\{  \Lambda_{j}\Lambda_{k}\phi\circ\pi_{\text{P}}\right\}
](y_{0})$

$\ \ \ \ \ \ =\frac{1}{24}[\underset{i=q+1}{\overset{n}{\sum}}%
\underset{\text{a,b=1}}{\overset{\text{q}}{\sum}}\frac{\partial^{2}%
\text{g}^{\text{ab}}}{\partial x_{i}^{2}}\left\{  \Lambda_{\text{a}}%
\Lambda_{\text{b}}\phi\circ\pi_{\text{P}}\right\}  ](y_{0})$

\qquad$\ $\ $+\frac{1}{24}[\underset{i,j,k=q+1}{\overset{n}{\sum}}%
\frac{\partial^{2}\text{g}^{jk}}{\partial x_{i}^{2}}\left\{  \Lambda
_{j}\Lambda_{k}\phi\circ\pi_{\text{P}}\right\}  ](y_{0})$

\qquad$\ +\frac{1}{12}[\underset{i,k=q+1}{\overset{n}{\sum}}%
\underset{\text{a=1}}{\overset{\text{q}}{\sum}}\frac{\partial^{2}%
\text{g}^{\text{a}k}}{\partial x_{i}^{2}}\left\{  \Lambda_{\text{a}}%
\Lambda_{k}\phi\circ\pi_{\text{P}}\right\}  ](y_{0})$

Since $\Lambda_{j}(y_{0})\Lambda_{k}(y_{0})=0$ for $j,k=q+1,...,n$ by $\left(
6.13\right)  ,$ we have:

I$_{3251}\ =\frac{1}{24}[\underset{i=q+1}{\overset{n}{\sum}}%
\underset{\text{a,b=1}}{\overset{\text{q}}{\sum}}\frac{\partial^{2}%
\text{g}^{\text{ab}}}{\partial x_{i}^{2}}\left\{  \Lambda_{\text{a}}%
\Lambda_{\text{b}}\phi\circ\pi_{\text{P}}\right\}  ](y_{0})$

$+\frac{1}{12}[\underset{i,j=q+1}{\overset{n}{\sum}}\underset{\text{a=1}%
}{\overset{\text{q}}{\sum}}\frac{\partial^{2}\text{g}^{\text{a}j}}{\partial
x_{i}^{2}}\left\{  \Lambda_{\text{a}}\Lambda_{j}\phi\circ\pi_{\text{P}%
}\right\}  ](y_{0})=$ I$_{32511}+$ I$_{32512}$

By (iii) of \textbf{Table A}$_{6}$,

I$_{32511}=\frac{1}{24}[\underset{i\text{=q+1}}{\overset{\text{n}}{\sum}%
}\underset{\text{a,b=1}}{\overset{\text{q}}{\sum}}\frac{\partial^{2}%
\text{g}^{\text{ab}}}{\partial x_{i}^{2}}\left\{  \Lambda_{\text{a}}%
\Lambda_{\text{b}}\phi\circ\pi_{\text{P}}\right\}  ](y_{0})$

$=\frac{1}{24}\underset{i\text{=}q+1}{\overset{\text{n}}{\sum}}%
\underset{\text{a,b=1}}{\overset{\text{q}}{\sum}}$ $2[-R_{\text{a}i\text{b}%
i}+5\overset{q}{\underset{\text{c}=1}{\sum}}T_{\text{ac}i}T_{\text{bc}%
i}+2\overset{n}{\underset{\text{k}=q+1}{\sum}}\perp_{\text{a}i\text{k}}%
\perp_{\text{b}i\text{k}}](y_{0})\times\lbrack\Lambda_{\text{a}}(y_{0}%
)\Lambda_{\text{b}}(y_{0})\phi(y_{0})]$

$=\frac{1}{12}\underset{i\text{=}q+1}{\overset{\text{n}}{\sum}}%
\underset{\text{a,b=1}}{\overset{\text{q}}{\sum}}$ $[-R_{\text{a}i\text{b}%
i}+5\overset{q}{\underset{\text{c}=1}{\sum}}T_{\text{ac}i}T_{\text{bc}%
i}+2\overset{n}{\underset{\text{k}=q+1}{\sum}}\perp_{\text{a}i\text{k}}%
\perp_{\text{b}i\text{k}}](y_{0})\times\lbrack\Lambda_{\text{a}}(y_{0}%
)\Lambda_{\text{b}}(y_{0})\phi(y_{0})]$

\qquad\qquad\qquad\qquad\qquad\qquad\qquad\qquad\qquad\qquad\qquad\qquad
\qquad\qquad\qquad\qquad\qquad\qquad$\blacksquare$

Since by (iii) of \textbf{Table A}$_{4},$ $\frac{\partial^{2}\text{g}%
^{\text{a}j}}{\partial\text{x}_{i}^{2}}(y_{0})=\frac{8}{3}R_{i\text{a}%
ij}+4\underset{\text{b=1}}{\overset{\text{q}}{\sum}}T_{\text{ab}i}(y_{0}%
)\perp_{\text{b}ji}(y_{0}),$ we next we have,

I$_{32512}=\frac{1}{12}[\underset{i,j=q+1}{\overset{n}{\sum}}%
\underset{\text{a=1}}{\overset{\text{q}}{\sum}}\frac{\partial^{2}%
\text{g}^{\text{a}j}}{\partial x_{i}^{2}}\Lambda_{\text{a}}\Lambda_{j}%
(\phi\circ\pi_{\text{P}})](y_{0})=\frac{1}{12}[\frac{8}{3}R_{i\text{a}%
ij}-4\underset{\text{b=1}}{\overset{\text{q}}{\sum}}T_{\text{ab}i}(y_{0}%
)\perp_{\text{b}ij}](y_{0})[\Lambda_{\text{a}}\Lambda_{j}\phi](y_{0})$

\qquad\qquad\qquad\qquad\qquad\qquad\qquad\qquad\qquad\qquad\qquad\qquad
\qquad\qquad\qquad\qquad\qquad\qquad$\blacksquare$

Therefore we have,

$\left(  D_{18}\right)  \qquad$I$_{3251}=$ I$_{32511}+$ I$_{32512}$

$=\frac{1}{12}\underset{i\text{=}q+1}{\overset{\text{n}}{\sum}}%
\underset{\text{a,b=1}}{\overset{\text{q}}{\sum}}$ $[-R_{\text{a}i\text{b}%
i}+5\overset{q}{\underset{\text{c}=1}{\sum}}T_{\text{ac}i}T_{\text{bc}%
i}+2\overset{n}{\underset{\text{k}=q+1}{\sum}}\perp_{\text{a}i\text{k}}%
\perp_{\text{b}i\text{k}}](y_{0})\times\lbrack\Lambda_{\text{a}}(y_{0}%
)\Lambda_{\text{b}}(y_{0})\phi(y_{0})]$

$\qquad+\frac{1}{12}[\frac{8}{3}R_{i\text{a}ij}-4\underset{\text{b=1}%
}{\overset{\text{q}}{\sum}}T_{\text{ab}i}(y_{0})\perp_{\text{b}ij}%
](y_{0})[\Lambda_{\text{a}}\Lambda_{j}\phi](y_{0})$

\qquad\qquad\qquad\qquad\qquad\qquad\qquad\qquad\qquad\qquad\qquad\qquad
\qquad\qquad\qquad\qquad\qquad$\blacksquare$

Since g$^{jk}(y_{0})=\delta^{jk},$

I$_{3252}=$ $\frac{1}{24}[\underset{i=q+1}{\overset{n}{\sum}}%
\underset{j,k=1}{\overset{n}{\sum}}$g$^{jk}$ $\frac{\partial^{2}}{\partial
x_{i}^{2}}\left\{  \Lambda_{j}\Lambda_{k}\phi\circ\pi_{\text{P}}\right\}
](y_{0})=$ $\frac{1}{24}[\underset{i=q+1}{\overset{n}{\sum}}%
\underset{i,j=1}{\overset{n}{\sum}}$g$\frac{\partial^{2}}{\partial x_{i}^{2}%
}(\Lambda_{j}^{2})\phi\circ\pi_{\text{P}}](y_{0})$

I$_{3252}=$ $\frac{1}{24}[\underset{i=q+1}{\overset{n}{\sum}}%
\underset{j=1}{\overset{n}{\sum}}\frac{\partial^{2}\Lambda_{j}^{2}}{\partial
x_{i}^{2}}(y_{0})\phi(y_{0})$

$=$ $\frac{1}{24}[\underset{i=q+1}{\overset{n}{\sum}}\underset{\text{a=1}%
}{\overset{\text{q}}{\sum}}\frac{\partial^{2}\Lambda_{\text{a}}^{2}}{\partial
x_{i}^{2}}(y_{0})\phi(y_{0})+$ $\frac{1}{24}%
[\underset{i,j=q+1}{\overset{n}{\sum}}\frac{\partial^{2}\Lambda_{j}^{2}%
}{\partial x_{i}^{2}}(y_{0})\phi(y_{0})=$ I$_{32521}+$ I$_{32522}$

$\left(  D_{18}\right)  ^{\ast}\qquad$I$_{32521}=\frac{1}{24}%
\underset{i=q+1}{\overset{n}{\sum}}\underset{\text{a=1}}{\overset{\text{q}%
}{\sum}}\frac{\partial^{2}\Lambda_{\text{a}}^{2}}{\partial x_{i}^{2}}%
(y_{0})\phi(y_{0})=\frac{1}{24}\underset{i=q+1}{\overset{n}{\sum}%
}\underset{\text{a=1}}{\overset{\text{q}}{\sum}}[\frac{\partial^{2}%
\Lambda_{\text{a}}}{\partial x_{i}^{2}}\Lambda_{\text{a}}+\Lambda_{\text{a}%
}\frac{\partial^{2}\Lambda_{\text{a}}}{\partial x_{i}^{2}}](y_{0})\phi(y_{0})$

$\qquad\qquad\qquad\qquad+\frac{1}{12}[$ $\frac{\partial\Lambda_{\text{a}}%
}{\partial x_{i}}$ $\frac{\partial\Lambda_{\text{a}}}{\partial x_{i}}%
](y_{0})\phi(y_{0})$

By (vii) \textbf{Proposition 5,}

$\left(  D_{18}\right)  ^{\ast\ast}$\qquad\qquad\ $\frac{\partial
\Lambda_{\text{a}}}{\partial x_{i}}(y_{0})=\Omega_{i\text{a}}(y_{0}%
)+[\Lambda_{\text{a}},\Lambda_{i}](y_{0})$

By (viii) of \textbf{Proposition 5, }we have\textbf{:}

$\ \frac{\partial^{2}\Lambda_{\text{a}}}{\partial x_{i}\partial x_{j}}%
(y_{0})=[\frac{\partial\Omega_{j\text{a}}}{\partial x_{i}}+\frac
{\partial\Lambda_{\text{a}}}{\partial x_{i}}\Lambda_{j}-\Lambda_{j}%
\frac{\partial\Lambda_{\text{a}}}{\partial x_{i}}+\Lambda_{\text{a}}%
\frac{\partial\Lambda_{j}}{\partial x_{i}}-\frac{\partial\Lambda_{j}}{\partial
x_{i}}\Lambda_{\text{a}}](y_{0})$

$\qquad=\frac{\partial\Omega_{j\text{a}}}{\partial x_{i}}(y_{0})+[\frac
{\partial\Lambda_{\text{a}}}{\partial x_{i}},\Lambda_{j}](y_{0})+[\Lambda
_{\text{a}},\frac{\partial\Lambda_{j}}{\partial x_{i}}](y_{0})$

Since \ $\frac{\partial\Lambda_{\text{a}}}{\partial x_{i}}(y_{0}%
)=\Omega_{i\text{a}}(y_{0})+[\Lambda_{\text{a}},\Lambda_{i}](y_{0}),$ and
$\frac{\partial\Lambda_{i}}{\partial x_{i}}](y_{0})=0,$ the last equation
above gives:

$\frac{\partial^{2}\Lambda_{\text{a}}}{\partial x_{i}^{2}}(y_{0}%
)=\frac{\partial\Omega_{i\text{a}}}{\partial x_{i}}(y_{0})+[\frac
{\partial\Lambda_{\text{a}}}{\partial x_{i}},\Lambda_{i}](y_{0})=\frac
{\partial\Omega_{i\text{a}}}{\partial x_{i}}(y_{0})+[\Omega_{i\text{a}%
}+[\Lambda_{\text{a}},\Lambda_{i}],\Lambda_{i}](y_{0})$

Therefore from $\left(  D_{18}\right)  ^{\ast}$ and $\left(  D_{18}\right)
^{\ast\ast}$ above,

I$_{32521}=\frac{1}{24}\underset{i=q+1}{\overset{n}{\sum}}\underset{\text{a=1}%
}{\overset{\text{q}}{\sum}}\frac{\partial^{2}\Lambda_{\text{a}}^{2}}{\partial
x_{i}^{2}}(y_{0})\phi(y_{0})$

$\qquad=\frac{1}{24}\underset{i=q+1}{\overset{n}{\sum}}\underset{\text{a=1}%
}{\overset{\text{q}}{\sum}}\left(  \frac{\partial\Omega_{i\text{a}}}{\partial
x_{i}}\Lambda_{\text{a}}+[\Omega_{i\text{a}}+[\Lambda_{\text{a}},\Lambda
_{i}],\Lambda_{i}]\right)  \Lambda_{\text{a}}(y_{0})\phi(y_{0})$

$\qquad+\frac{1}{24}\underset{i=q+1}{\overset{n}{\sum}}\underset{\text{a=1}%
}{\overset{\text{q}}{\sum}}\Lambda_{\text{a}}(y_{0})\left(  \frac
{\partial\Omega_{i\text{a}}}{\partial x_{i}}\Lambda_{\text{a}}+[\Omega
_{i\text{a}}+[\Lambda_{\text{a}},\Lambda_{i}],\Lambda_{i}]\right)  (y_{0}%
)\phi(y_{0})$

$\qquad+\frac{1}{12}\underset{i=q+1}{\overset{n}{\sum}}\underset{\text{a=1}%
}{\overset{\text{q}}{\sum}}2[\Omega_{i\text{a}}+[\Lambda_{\text{a}}%
,\Lambda_{i}]]^{2}(y_{0})\phi(y_{0})$

\qquad\qquad\qquad\qquad\qquad\qquad\qquad\qquad\qquad\qquad\qquad\qquad
\qquad\qquad$\qquad\blacksquare$

By (xiii) of \textbf{Proposition 5,} we have:

I$_{32522}\qquad\qquad\ \frac{\partial^{2}\Lambda_{j}^{2}}{\partial
\text{x}_{i}^{2}}(y_{0})=\frac{1}{2}\left(  \Omega_{ij}\Omega_{ij}\right)
(y_{0})+\frac{1}{3}[\frac{\partial\Omega_{ij}}{\partial\text{x}_{i}}%
\Lambda_{j}+\Lambda_{j}\frac{\partial\Omega_{ij}}{\partial\text{x}_{i}}%
](y_{0})$

\qquad\qquad\qquad\qquad\qquad\qquad\qquad\qquad\qquad\qquad\qquad\qquad
\qquad\qquad\qquad$\blacksquare$Therefore,

$\left(  D_{19}\right)  \qquad$I$_{3252}=\frac{1}{24}%
\underset{i=q+1}{\overset{n}{\sum}}\underset{\text{a=1}}{\overset{\text{q}%
}{\sum}}\left(  \frac{\partial\Omega_{i\text{a}}}{\partial x_{i}}%
\Lambda_{\text{a}}+[\Omega_{i\text{a}}+[\Lambda_{\text{a}},\Lambda
_{i}],\Lambda_{i}]\right)  \Lambda_{\text{a}}(y_{0})\phi(y_{0})$

$\qquad\qquad\qquad+\frac{1}{24}\underset{i=q+1}{\overset{n}{\sum}%
}\underset{\text{a=1}}{\overset{\text{q}}{\sum}}\Lambda_{\text{a}}%
(y_{0})\left(  \frac{\partial\Omega_{i\text{a}}}{\partial x_{i}}%
\Lambda_{\text{a}}+[\Omega_{i\text{a}}+[\Lambda_{\text{a}},\Lambda
_{i}],\Lambda_{i}]\right)  (y_{0})\phi(y_{0})$

$\qquad\qquad\qquad+\frac{1}{12}\underset{i=q+1}{\overset{n}{\sum}%
}\underset{\text{a=1}}{\overset{\text{q}}{\sum}}[\Omega_{i\text{a}}%
+[\Lambda_{\text{a}},\Lambda_{i}]]^{2}(y_{0})\phi(y_{0})$

$\qquad\qquad\qquad+\frac{1}{48}\underset{i,j=q+1}{\overset{n}{\sum}}\left(
\Omega_{ij}\Omega_{ij}\right)  (y_{0})\phi(y_{0})+\frac{1}{72}%
\underset{i,j=q+\text{1}}{\overset{n}{\sum}}\left(  \frac{\partial\Omega_{ij}%
}{\partial\text{x}_{i}}\Lambda_{j}+\Lambda_{j}\frac{\partial\Omega_{ij}%
}{\partial\text{x}_{i}}\right)  (y_{0})\phi(y_{0})$

\qquad\qquad\qquad\qquad\qquad\qquad\qquad\qquad\qquad\qquad\qquad\qquad
\qquad\qquad\qquad\qquad\qquad\qquad$\blacksquare$

We come to the last expression of I$_{325}:$

Since $\frac{\partial}{\partial x_{i}}(\phi\circ\pi_{\text{P}})(y_{0}%
)=\frac{\partial\phi}{\partial x_{i}}(\pi_{\text{P}}(y_{0}))\frac{\partial
\pi_{\text{P}}}{\partial x_{i}}(y_{0})=0$ for $i=q+1,...,n,$

I$_{3253}=\frac{1}{12}[\underset{i\text{=q+1}}{\overset{n}{\sum}%
}\underset{j,k\text{=1}}{\overset{n}{\sum}}\frac{\partial\text{g}^{jk}%
}{\partial x_{i}}\frac{\partial}{\partial x_{i}}\left\{  \text{ }\Lambda
_{j}\Lambda_{k}\phi\circ\pi_{\text{P}}\right\}  ](y_{0})$

$\qquad=\frac{1}{12}[\underset{i=q+\text{1}}{\overset{n}{\sum}}%
\underset{\text{a,b=1}}{\overset{\text{q}}{\sum}}\frac{\partial\text{g}%
^{\text{ab}}}{\partial x_{i}}\frac{\partial}{\partial x_{i}}(\Lambda
_{\text{a}}\Lambda_{\text{b}})\phi\circ\pi_{\text{P}}](y_{0})$

\qquad$\ +\frac{1}{12}[\underset{i=q+1}{\overset{n}{\sum}}%
\underset{j,k=q+1}{\overset{n}{\sum}}\frac{\partial\text{g}^{jk}}{\partial
x_{i}}\frac{\partial}{\partial x_{i}}(\Lambda_{j}\Lambda_{k})\phi\circ
\pi_{\text{P}}](y_{0})$

\qquad$+\frac{1}{12}[\underset{i,j=q+1}{\overset{n}{\sum}}\underset{\text{a=1}%
}{\overset{\text{q}}{\sum}}\frac{\partial\text{g}^{\text{a}j}}{\partial x_{i}%
}\frac{\partial}{\partial x_{i}}(\Lambda_{\text{a}}\Lambda_{j})\phi\circ
\pi_{\text{P}}](y_{0})$

\qquad$+\frac{1}{12}[\underset{i,j=q+1}{\overset{n}{\sum}}\underset{\text{b=1}%
}{\overset{\text{q}}{\sum}}\frac{\partial\text{g}^{j\text{b}}}{\partial x_{i}%
}\frac{\partial}{\partial x_{i}}(\Lambda_{j}\Lambda_{\text{b}})\phi\circ
\pi_{\text{P}}](y_{0})$

$\frac{\partial\text{g}^{\text{ab}}}{\partial x_{i}}(y_{0})=2T_{\text{ab}%
i}(y_{0})$ by (ii) of \textbf{Table A}$_{6}$.

$\frac{\partial\text{g}^{ij}}{\partial x_{\alpha}}(y_{0})=0$ for
$\alpha,i,j=q+1,...,n$ by (ii) of \textbf{Table A}$_{2}.$

$\frac{\partial\text{g}^{\text{a}j}}{\partial x_{i}}(y_{0})=\perp_{\text{a}%
ji}(y_{0})$ by (iii) of \textbf{Table A}$_{4}$:

$\frac{\partial\text{g}^{j\text{b}}}{\partial x_{i}}(y_{0})=\frac
{\partial\text{g}^{\text{b}j}}{\partial x_{i}}(y_{0})=\perp_{\text{b}ji}%
(y_{0})=-\perp_{\text{b}ij}(y_{0})$

I$_{3253}=\frac{1}{12}[\underset{i=q+1}{\overset{n}{\sum}}%
\underset{\text{a,b=1}}{\overset{\text{q}}{\sum}}2T_{\text{ab}i}%
(y_{0})\left\{  (\frac{\partial\Lambda_{\text{a}}}{\partial x_{i}}%
\Lambda_{\text{b}}+\Lambda_{\text{a}}\frac{\partial\Lambda_{\text{b}}%
}{\partial x_{i}})\phi\circ\pi_{\text{P}}\right\}  ](y_{0})$

\qquad$-\frac{1}{12}[\underset{i,j=q+1}{\overset{n}{\sum}}%
\underset{\text{a,b=1}}{\overset{\text{q}}{\sum}}\perp_{\text{b}ij}%
(y_{0})\left\{  (\frac{\partial\Lambda_{\text{a}}}{\partial x_{i}}\Lambda
_{j}+\Lambda_{\text{a}}\frac{\partial\Lambda_{j}}{\partial x_{i}})\phi\circ
\pi_{\text{P}}\right\}  ](y_{0})$

\qquad$-\frac{1}{12}[\underset{i,j=q+1}{\overset{n}{\sum}}\underset{\text{b=1}%
}{\overset{\text{q}}{\sum}}\perp_{\text{b}ij}(y_{0})\left\{  \frac
{\partial\Lambda_{j}}{\partial x_{i}}\Lambda_{\text{b}}+\Lambda_{j}%
\frac{\partial\Lambda_{\text{b}}}{\partial x_{i}})\phi\circ\pi_{\text{P}%
}\right\}  ](y_{0})$

\qquad\ $\frac{\partial\Lambda_{j}}{\partial\text{x}_{i}}(y_{0})=\frac{1}%
{2}\Omega_{ij}(y_{0})$ and $\frac{\partial\Lambda_{\text{a}}}{\partial x_{i}%
}=\Omega_{i\text{a}}+[\Lambda_{\text{a}},\Lambda_{i}]$ by (vii)
\textbf{Proposition 5. }

We then have:

$\left(  D_{20}\right)  $ I$_{3253}$

$\qquad=\frac{1}{12}[\underset{i=q+1}{\overset{n}{\sum}}\underset{\text{a,b=1}%
}{\overset{\text{q}}{\sum}}2T_{\text{ab}i}(y_{0})\left\{  (\Omega_{i\text{a}%
}+[\Lambda_{\text{a}},\Lambda_{i}])\Lambda_{\text{b}}+\Lambda_{\text{a}%
}(\Omega_{i\text{a}}+[\Lambda_{\text{a}},\Lambda_{i}])\right\}  (y_{0}%
)\phi(y_{0})$

\qquad$-\frac{1}{12}[\underset{i,j=q+1}{\overset{n}{\sum}}\underset{\text{a=1}%
}{\overset{\text{q}}{\sum}}\perp_{\text{a}ij}(y_{0})\left\{  (\Omega
_{i\text{a}}+[\Lambda_{\text{a}},\Lambda_{i}])\Lambda_{j}+\frac{1}{2}%
\Lambda_{\text{a}}\Omega_{ij}\right\}  ](y_{0})\phi(y_{0})$

\qquad$-\frac{1}{12}[\underset{i,j=q+1}{\overset{n}{\sum}}\underset{\text{b=1}%
}{\overset{\text{q}}{\sum}}\perp_{\text{b}ij}(y_{0})\left\{  \frac{1}{2}%
\Omega_{ij}\Lambda_{\text{b}}+\Lambda_{j}(\Omega_{i\text{b}}+[\Lambda
_{\text{b}},\Lambda_{i}])\right\}  ](y_{0})\phi(y_{0})$

The expression for I$_{325}=$ I$_{3251}+$I$_{3252}+$ I$_{3253}$ is given in
$\left(  D_{18}\right)  ,$ $\left(  D_{19}\right)  $ and $\left(
D_{20}\right)  :$

$\left(  D_{21}\right)  \qquad$I$_{325}$

$\qquad=\frac{1}{12}\underset{i\text{=}q+1}{\overset{\text{n}}{\sum}%
}\underset{\text{a,b=1}}{\overset{\text{q}}{\sum}}$ $[-R_{\text{a}i\text{b}%
i}+5\overset{q}{\underset{\text{c}=1}{\sum}}T_{\text{ac}i}T_{\text{bc}%
i}+2\overset{n}{\underset{\text{k}=q+1}{\sum}}\perp_{\text{a}i\text{k}}%
\perp_{\text{b}i\text{k}}](y_{0})\qquad$I$_{3251}$

\qquad\qquad$\times\lbrack\Lambda_{\text{a}}(y_{0})\Lambda_{\text{b}}%
(y_{0})\phi(y_{0})]$

$\qquad\qquad+\frac{1}{12}[\frac{8}{3}R_{i\text{a}ij}-4\underset{\text{b=1}%
}{\overset{\text{q}}{\sum}}T_{\text{ab}i}(y_{0})\perp_{\text{b}ij}%
](y_{0})[\Lambda_{\text{a}}\Lambda_{j}\phi](y_{0})$

\qquad\qquad$+\frac{1}{24}\underset{i=q+1}{\overset{n}{\sum}}%
\underset{\text{a=1}}{\overset{\text{q}}{\sum}}\left(  \frac{\partial
\Omega_{i\text{a}}}{\partial x_{i}}\Lambda_{\text{a}}+[\Omega_{i\text{a}%
}+[\Lambda_{\text{a}},\Lambda_{i}],\Lambda_{i}]\right)  \Lambda_{\text{a}%
}(y_{0})\phi(y_{0})\qquad$I$_{3252}$

$\qquad\qquad+\frac{1}{24}\underset{i=q+1}{\overset{n}{\sum}}%
\underset{\text{a=1}}{\overset{\text{q}}{\sum}}\Lambda_{\text{a}}%
(y_{0})\left(  \frac{\partial\Omega_{i\text{a}}}{\partial x_{i}}%
\Lambda_{\text{a}}+[\Omega_{i\text{a}}+[\Lambda_{\text{a}},\Lambda
_{i}],\Lambda_{i}]\right)  (y_{0})\phi(y_{0})$

$\qquad\qquad+\frac{1}{12}\underset{i=q+1}{\overset{n}{\sum}}%
\underset{\text{a=1}}{\overset{\text{q}}{\sum}}\left(  \Omega_{i\text{a}%
}+[\Lambda_{\text{a}},\Lambda_{i}]\right)  ^{2}(y_{0})\phi(y_{0})$

$\qquad\qquad+\frac{1}{48}\underset{i,j=q+1}{\overset{n}{\sum}}\left(
\Omega_{ij}\Omega_{ij}\right)  (y_{0})\phi(y_{0})+\frac{1}{72}%
\underset{i,j=q+\text{1}}{\overset{n}{\sum}}\left(  \frac{\partial\Omega_{ij}%
}{\partial\text{x}_{i}}\Lambda_{j}+\Lambda_{j}\frac{\partial\Omega_{ij}%
}{\partial\text{x}_{i}}\right)  (y_{0})\phi(y_{0})$

\qquad$+\frac{1}{12}[\underset{i=q+1}{\overset{n}{\sum}}\underset{\text{a,b=1}%
}{\overset{\text{q}}{\sum}}2T_{\text{ab}i}(y_{0})\left\{  (\Omega_{i\text{a}%
}+[\Lambda_{\text{a}},\Lambda_{i}])\Lambda_{\text{b}}+\Lambda_{\text{a}%
}(\Omega_{i\text{a}}+[\Lambda_{\text{a}},\Lambda_{i}])\right\}  (y_{0}%
)\phi(y_{0})\qquad$I$_{3253}$

\qquad$-\frac{1}{12}[\underset{i,j=q+1}{\overset{n}{\sum}}\underset{\text{a=1}%
}{\overset{\text{q}}{\sum}}\perp_{\text{a}ij}(y_{0})\left\{  (\Omega
_{i\text{a}}+[\Lambda_{\text{a}},\Lambda_{i}])\Lambda_{j}+\frac{1}{2}%
\Lambda_{\text{a}}\Omega_{ij}\right\}  ](y_{0})\phi(y_{0})$

\qquad$-\frac{1}{12}[\underset{i,j=q+1}{\overset{n}{\sum}}\underset{\text{b=1}%
}{\overset{\text{q}}{\sum}}\perp_{\text{b}ij}(y_{0})\left\{  \frac{1}{2}%
\Omega_{ij}\Lambda_{\text{b}}+\Lambda_{j}(\Omega_{i\text{b}}+[\Lambda
_{\text{b}},\Lambda_{i}])\right\}  ](y_{0})\phi(y_{0})$

\qquad\qquad\qquad\qquad\qquad\qquad\qquad\qquad\qquad\qquad\qquad\qquad
\qquad\qquad\qquad\qquad\qquad$\blacksquare$

Next we consider:

\textbf{I}$_{326}=-\frac{1}{24}\underset{i=q+1}{\overset{n}{\sum}}%
\frac{\partial^{2}}{\partial x_{i}^{2}}[\underset{j,k=1}{\overset{n}{\sum}}%
$g$^{jk}\left\{  \underset{\text{c=1}}{\overset{\text{q}}{\sum}}\Gamma
_{jk}^{\text{c}}\frac{\partial\phi}{\partial\text{x}_{\text{c}}}\circ
\pi_{\text{P}}+\text{ }\underset{l=1}{\overset{n}{\sum}}\Gamma_{jk}^{l}%
\Lambda_{l}\phi\circ\pi_{\text{P}}\right\}  ](y_{0})$

\qquad$=-\frac{1}{24}\underset{i=q+1}{\overset{n}{\sum}}%
[\underset{j,k=1}{\overset{n}{\sum}}\frac{\partial^{2}\text{g}^{jk}}{\partial
x_{i}^{2}}\left\{  \underset{\text{c=1}}{\overset{\text{q}}{\sum}}\Gamma
_{jk}^{\text{c}}\frac{\partial\phi}{\partial\text{x}_{\text{c}}}\circ
\pi_{\text{P}}+\text{ }\underset{l=1}{\overset{n}{\sum}}\Gamma_{jk}^{l}%
\Lambda_{l}\phi\circ\pi_{\text{P}}\right\}  ](y_{0})$

\qquad$-\frac{1}{24}\underset{i=q+1}{\overset{n}{\sum}}%
[\underset{j,k=1}{\overset{n}{\sum}}$g$^{jk}\left\{  \underset{\text{c=1}%
}{\overset{\text{q}}{\sum}}\frac{\partial^{2}}{\partial x_{i}^{2}}(\Gamma
_{jk}^{\text{c}}\frac{\partial\phi}{\partial\text{x}_{\text{c}}}\circ
\pi_{\text{P}})\text{ }+\text{ }\underset{l=1}{\overset{n}{\sum}}%
\frac{\partial^{2}}{\partial x_{i}^{2}}(\Gamma_{jk}^{l}\Lambda_{l}\phi\circ
\pi_{\text{P}})\right\}  ](y_{0})$

\qquad$-\frac{1}{12}\underset{i=q+1}{\overset{n}{\sum}}%
[\underset{j,k=1}{\overset{n}{\sum}}\frac{\partial\text{g}^{jk}}{\partial
x_{i}}\left\{  \underset{\text{c=1}}{\overset{\text{q}}{\sum}}\frac{\partial
}{\partial x_{i}}(\Gamma_{jk}^{\text{c}}\frac{\partial\phi}{\partial
\text{x}_{\text{c}}}\circ\pi_{\text{P}})\text{ }+\text{ }%
\underset{k=1}{\overset{n}{\sum}}\frac{\partial}{\partial x_{i}}(\Gamma
_{jk}^{l}\Lambda_{l}\phi\circ\pi_{\text{P}})\right\}  ](y_{0})$

$\qquad=$ \textbf{I}$_{3261}+$ \textbf{I}$_{3262}+$ \textbf{I}$_{3263},$ where,

I$_{3261}=$ $-\frac{1}{24}\underset{i=q+1}{\overset{n}{\sum}}%
[\underset{j,k=1}{\overset{n}{\sum}}\frac{\partial^{2}\text{g}^{jk}}{\partial
x_{i}^{2}}\left\{  \underset{\text{c=1}}{\overset{\text{q}}{\sum}}\Gamma
_{jk}^{\text{c}}\frac{\partial\phi}{\partial\text{x}_{\text{c}}}\circ
\pi_{\text{P}}+\text{ }\underset{l=1}{\overset{n}{\sum}}\Gamma_{jk}^{l}%
\Lambda_{l}\phi\circ\pi_{\text{P}}\right\}  ](y_{0})$

I$_{3262}=\ -\frac{1}{24}\underset{i=q+1}{\overset{n}{\sum}}%
[\underset{j,k=1}{\overset{n}{\sum}}$g$^{jk}\left\{  \underset{\text{c=1}%
}{\overset{\text{q}}{\sum}}\frac{\partial^{2}\Gamma_{jk}^{\text{c}}}{\partial
x_{i}^{2}}\frac{\partial\phi}{\partial\text{x}_{\text{c}}}\circ\pi_{\text{P}%
}\text{ }+\text{ }\underset{k=1}{\overset{n}{\sum}}\frac{\partial^{2}%
}{\partial x_{i}^{2}}(\Gamma_{jk}^{l}\Lambda_{l}\phi\circ\pi_{\text{P}%
})\right\}  ](y_{0})$

\textbf{I}$_{3263}=-\frac{1}{12}\underset{i=q+1}{\overset{n}{\sum}%
}[\underset{j,k=1}{\overset{n}{\sum}}\frac{\partial\text{g}^{jk}}{\partial
x_{i}}\left\{  \underset{\text{c=1}}{\overset{\text{q}}{\sum}}\frac
{\partial\Gamma_{jk}^{\text{c}}}{\partial x_{i}}\frac{\partial\phi}%
{\partial\text{x}_{\text{c}}}\circ\pi_{\text{P}}\text{ }+\text{ }%
\underset{l=1}{\overset{n}{\sum}}\frac{\partial}{\partial x_{i}}(\Gamma
_{jk}^{l}\Lambda_{l}\phi\circ\pi_{\text{P}})\right\}  ](y_{0})$

We compute:

I$_{3261}=-\frac{1}{24}\underset{i=q+1}{\overset{n}{\sum}}%
[\underset{j,k=1}{\overset{n}{\sum}}\frac{\partial^{2}\text{g}^{jk}}{\partial
x_{i}^{2}}\left\{  \underset{\text{c=1}}{\overset{\text{q}}{\sum}}\Gamma
_{jk}^{\text{c}}\frac{\partial\phi}{\partial\text{x}_{\text{c}}}\circ
\pi_{\text{P}}+\text{ }\underset{l=1}{\overset{n}{\sum}}\Gamma_{jk}^{l}%
\Lambda_{l}\phi\circ\pi_{\text{P}}\right\}  ](y_{0})$

$\qquad=\mathbf{-}$ $\frac{1}{24}\underset{i=q+\text{1}}{\overset{n}{\sum}%
}[\underset{\text{a,b=1}}{\overset{\text{q}}{\sum}}\frac{\partial^{2}%
\text{g}^{\text{ab}}}{\partial x_{i}^{2}}\left\{  \underset{\text{c=1}%
}{\overset{\text{q}}{\sum}}\Gamma_{\text{ab}}^{\text{c}}\frac{\partial\phi
}{\partial\text{x}_{\text{c}}}\circ\pi_{\text{P}}\text{ + }%
\underset{l\text{=1}}{\overset{\text{n}}{\sum}}\Gamma_{\text{ab}}^{l}%
\Lambda_{l}\phi\circ\pi_{\text{P}}\right\}  ](y_{0})$

\qquad$\ \ -\frac{1}{24}\underset{i=q+1}{\overset{n}{\sum}}%
[\underset{j,k=q+1}{\overset{n}{\sum}}\frac{\partial^{2}\text{g}^{jk}%
}{\partial x_{i}^{2}}\left\{  \underset{\text{c=1}}{\overset{\text{q}}{\sum}%
}\Gamma_{jk}^{\text{c}}\frac{\partial\phi}{\partial\text{x}_{\text{c}}}%
\circ\pi_{\text{P}}+\text{ }\underset{l=1}{\overset{n}{\sum}}\Gamma_{jk}%
^{l}\Lambda_{l}\phi\circ\pi_{\text{P}}\right\}  ](y_{0})$

$\qquad\ \ \ \ \ \mathbf{-}$ $\frac{1}{12}[\underset{i,j=q+1}{\overset{n}{\sum
}}\underset{\text{a=1}}{\overset{\text{q}}{\sum}}\frac{\partial^{2}%
\text{g}^{\text{a}j}}{\partial x_{i}^{2}}\left\{  \underset{\text{c=1}%
}{\overset{\text{q}}{\sum}}\Gamma_{\text{a}j}^{\text{c}}\frac{\partial\phi
}{\partial\text{x}_{\text{c}}}\circ\pi_{\text{P}}\text{ + }%
\underset{l\text{=1}}{\overset{\text{n}}{\sum}}\Gamma_{\text{a}j}^{l}%
\Lambda_{l}\phi\circ\pi_{\text{P}}\right\}  ](y_{0})$

\qquad$=$ I$_{32611}+$ I$_{32612}+$ I$_{32613}$ (numbered in descending order)

We are now ready to express the above expressions in terms of geometric invariants:

I$_{32611}=\mathbf{-}$ $\frac{1}{24}\underset{i\text{=q+1}}{\overset{\text{n}%
}{\sum}}[\underset{\text{a,b=1}}{\overset{\text{q}}{\sum}}\frac{\partial
^{2}\text{g}^{\text{ab}}}{\partial x_{i}^{2}}\left\{  \underset{\text{c=1}%
}{\overset{\text{q}}{\sum}}\Gamma_{\text{ab}}^{\text{c}}\frac{\partial\phi
}{\partial\text{x}_{\text{c}}}\circ\pi_{\text{P}}\text{ + }%
\underset{l\text{=1}}{\overset{\text{n}}{\sum}}\Gamma_{\text{ab}}^{l}%
\Lambda_{l}\phi\circ\pi_{\text{P}}\right\}  ](y_{0})$

$=\mathbf{-}$ $\frac{1}{24}\underset{i\text{=q+1}}{\overset{\text{n}}{\sum}%
}[\underset{\text{a,b=1}}{\overset{\text{q}}{\sum}}\frac{\partial^{2}%
\text{g}^{\text{ab}}}{\partial x_{i}^{2}}\left\{  \underset{\text{c=1}%
}{\overset{\text{q}}{\sum}}\Gamma_{\text{ab}}^{\text{c}}\frac{\partial\phi
}{\partial\text{x}_{\text{c}}}\circ\pi_{\text{P}}\text{ + }%
\underset{\text{d=1}}{\overset{\text{q}}{\sum}}\Gamma_{\text{ab}}^{\text{d}%
}\Lambda_{\text{d}}\phi\circ\pi_{\text{P}}+\underset{l=q+1}{\overset{n}{\sum}%
}\Gamma_{\text{ab}}^{l}\Lambda_{l}\phi\circ\pi_{\text{P}}\right\}  ](y_{0})$

$\bigskip\qquad\Gamma_{\text{ab}}^{\text{c}}(y_{0})=0=\Gamma_{\text{ab}%
}^{\text{d}}(y_{0})$ for a,b,c,d = 1,...,q by (ii) of \textbf{Table A}$_{7}.$ Therefore,

\qquad by (iii) of \textbf{Table A}$_{6}$ and by (i) of \textbf{Table A}$_{7}$

I$_{32611}=\mathbf{-}$ $\frac{1}{24}[\underset{\text{a,b=1}}{\overset{\text{q}%
}{\sum}}\underset{i,j=q+1}{\overset{n}{\sum}}\frac{\partial^{2}\text{g}%
^{\text{ab}}}{\partial x_{i}^{2}}\Gamma_{\text{ab}}^{j}\Lambda_{j}\phi\circ
\pi_{\text{P}}](y_{0})$

$\left(  D_{22}\right)  \qquad$I$_{32611}$

$=$ $\mathbf{-}$ $\frac{1}{12}\underset{\text{a,b=1}}{\overset{\text{q}}{\sum
}}\underset{i,j=q+1}{\overset{n}{\sum}}T_{\text{ab}j}(y_{0})\left\{
-R_{\text{a}i\text{b}i}+5\overset{q}{\underset{\text{c}=1}{\sum}}%
T_{\text{ac}i}T_{\text{bc}i}+2\overset{n}{\underset{k=q+1}{\sum}}%
\perp_{\text{a}ik}\perp_{\text{b}ik}\right\}  (y_{0})\Lambda_{j}(y_{0}%
)\phi(y_{0})$

\qquad\qquad\qquad\qquad\qquad\qquad\qquad\qquad\qquad\qquad\qquad\qquad
\qquad\qquad\qquad\qquad\qquad\qquad\qquad\qquad$\blacksquare$

We next examine:

\ I$_{32612}=\ -\frac{1}{24}\underset{i=q+1}{\overset{n}{\sum}}%
[\underset{j,k=q+1}{\overset{n}{\sum}}\frac{\partial^{2}\text{g}^{jk}%
}{\partial x_{i}^{2}}\left\{  \underset{\text{c=1}}{\overset{\text{q}}{\sum}%
}\Gamma_{jk}^{\text{c}}\frac{\partial\phi}{\partial\text{x}_{\text{c}}}%
\circ\pi_{\text{P}}+\text{ }\underset{l=1}{\overset{n}{\sum}}\Gamma_{jk}%
^{l}\Lambda_{l}\phi\circ\pi_{\text{P}}\right\}  ](y_{0})$

Since $\Gamma_{jk}^{l}(y_{0})=0=\Gamma_{jk}^{\text{a}}(y_{0})$ for a = 1,...,q
and $i,j,l=q+1,...,n$ by (i) and (ii) of \textbf{Table A}$_{8},$

$\left(  D_{23}\right)  \qquad$I$_{32612}=0$

We then consider:

I$_{32613}=\ \mathbf{-}$ $\frac{1}{12}\underset{i=q+1}{\overset{n}{\sum}%
}\underset{j=q+1}{\overset{n}{\sum}}\underset{\text{a=1}}{\overset{\text{q}%
}{\sum}}\frac{\partial^{2}\text{g}^{\text{a}j}}{\partial x_{i}^{2}%
}[\underset{\text{c=1}}{\overset{\text{q}}{\sum}}\Gamma_{\text{a}j}^{\text{c}%
}\frac{\partial\phi}{\partial\text{x}_{\text{c}}}\circ\pi_{\text{P}}$ $+$
$\Gamma_{\text{a}j}^{l}\Lambda_{l}\phi\circ\pi_{\text{P}}](y_{0})$

$=\ \mathbf{-}$ $\frac{1}{12}\underset{i,j=q+1}{\overset{n}{\sum}%
}\underset{\text{a=1}}{\overset{\text{q}}{\sum}}\frac{\partial^{2}%
\text{g}^{\text{a}j}}{\partial x_{i}^{2}}[\underset{\text{c=1}%
}{\overset{\text{q}}{\sum}}\Gamma_{\text{a}j}^{\text{c}}\frac{\partial\phi
}{\partial\text{x}_{\text{c}}}\circ\pi_{\text{P}}$ $+\underset{\text{a=1}%
}{\overset{\text{q}}{\sum}}\Gamma_{\text{a}j}^{\text{b}}\Lambda_{\text{b}}%
\phi\circ\pi_{\text{P}}+\underset{k=q+1}{\overset{n}{\sum}}\Gamma_{\text{a}%
j}^{k}\Lambda_{k}\phi\circ\pi_{\text{P}}](y_{0})$

We have $\frac{\partial^{2}\text{g}^{\text{a}j}}{\partial\text{x}_{i}^{2}%
}(y_{0})=$ $4\underset{c=1}{\overset{q}{\sum}}(T_{\text{ac}i})(\perp
_{j\text{c}i})+\frac{8}{3}R_{i\text{a}ij}$ by (iii) \textbf{Table A}$_{4}$ and
$\Gamma_{\text{a}j}^{\text{b}}=-$T$_{\text{ab}j}(y_{0})$

by (iii) of \textbf{Table A}$_{8}$; $\Gamma_{\text{a}j}^{k}=$ $\perp
_{\text{a}jk}(y_{0})$ by (iv) of \textbf{Table A}$_{8}$.Consequently,

$\left(  D_{24}\right)  \qquad$I$_{32613}=\ \mathbf{-}\frac{1}{12}%
\underset{i,j=q+1}{\overset{n}{\sum}}\underset{\text{a=1}}{\overset{\text{q}%
}{\sum}}[4\underset{\text{c=1}}{\overset{q}{\sum}}(T_{\text{ac}i}%
)(\perp_{j\text{c}i})+\frac{8}{3}R_{i\text{a}ij}](y_{0})$

$\qquad\qquad\qquad\qquad\times\lbrack\underset{\text{c=1}}{\overset{\text{q}%
}{\sum}}-T_{\text{ac}j}\frac{\partial\phi}{\partial\text{x}_{\text{c}}%
}-\underset{\text{b=1}}{\overset{\text{q}}{\sum}}T_{\text{ab}j}\Lambda
_{\text{b}}+\underset{k=q+1}{\overset{n}{\sum}}\perp_{\text{a}jk}\Lambda
_{k}](y_{0})\phi(y_{0})$

\qquad\qquad\qquad\qquad\qquad\qquad\qquad\qquad\qquad\qquad\qquad\qquad
\qquad\qquad\qquad$\blacksquare$

We gather all the terms of I$_{3261}.$ These are given in $\left(
D_{22}\right)  ,$ $\left(  D_{23}\right)  $ and $\left(  D_{24}\right)  :$

$\left(  D_{25}\right)  \qquad$ I$_{3261}=\mathbf{-}\frac{1}{12}%
\underset{\text{a,b=1}}{\overset{\text{q}}{\sum}}%
\underset{i,j=q+1}{\overset{n}{\sum}}T_{\text{ab}i}(y_{0})\qquad\qquad\qquad
$I$_{32611}$

$\qquad\times\lbrack-R_{\text{a}i\text{b}i}+5\overset{q}{\underset{\text{c}%
=1}{\sum}}T_{\text{ac}i}T_{\text{bc}i}+2\overset{n}{\underset{k=q+1}{\sum}%
}\perp_{\text{a}ik}\perp_{\text{b}ik}](y_{0})\Lambda_{j}(y_{0})\phi
(y_{0})\ \ $

\qquad$+\frac{1}{12}\underset{i\text{=q+1}}{\overset{\text{n}}{\sum}%
}\underset{j\text{=q+1}}{\overset{\text{n}}{\sum}}\underset{\text{a=1}%
}{\overset{\text{q}}{\sum}}[4\underset{\text{c=1}}{\overset{q}{\sum}%
}(T_{\text{ac}i})(\perp_{j\text{c}i})+\frac{8}{3}R_{i\text{a}ij}](y_{0}%
)\qquad\qquad$I$_{32613}$

$\qquad\times\lbrack\underset{\text{c=1}}{\overset{\text{q}}{\sum}%
}T_{\text{ac}j}\frac{\partial\phi}{\partial\text{x}_{\text{c}}}%
+\underset{\text{b=1}}{\overset{\text{q}}{\sum}}T_{\text{ab}j}\Lambda
_{\text{b}}-\underset{k=q+1}{\overset{n}{\sum}}\perp_{\text{a}jk}\Lambda
_{k}](y_{0})\phi(y_{0})$

\qquad\qquad\qquad\qquad\qquad\qquad\qquad\qquad\qquad\qquad\qquad\qquad
\qquad\qquad\qquad\qquad\qquad\qquad$\blacksquare$

We next compute:

I$_{3262}=-\frac{1}{24}\underset{i=q+1}{\overset{n}{\sum}}%
\underset{j,k=1}{\overset{n}{\sum}}$g$^{jk}[\underset{\text{b=1}%
}{\overset{\text{q}}{\sum}}\frac{\partial^{2}}{\partial x_{i}^{2}}(\Gamma
_{jk}^{\text{b}}\frac{\partial\phi}{\partial\text{x}_{\text{b}}}\circ
\pi_{\text{P}})$ $+$ $\underset{l=1}{\overset{n}{\sum}}\frac{\partial^{2}%
}{\partial x_{i}^{2}}(\Gamma_{jk}^{l}\Lambda_{l}\phi\circ\pi_{\text{P}%
})](y_{0})$

Since g$^{jk}(y_{0})=\delta^{jk},$

$=-\frac{1}{24}\underset{i=q+1}{\overset{n}{\sum}}%
\underset{j=1}{\overset{n}{\sum}}[\underset{\text{b=1}}{\overset{\text{q}%
}{\sum}}\frac{\partial^{2}\Gamma_{jj}^{\text{b}}}{\partial x_{i}^{2}}%
\frac{\partial\phi}{\partial\text{x}_{\text{b}}}\circ\pi_{\text{P}}$ $+$
$\underset{l=q+1}{\overset{n}{\sum}}\frac{\partial^{2}}{\partial x_{i}^{2}%
}(\Gamma_{jj}^{l}\Lambda_{l})\phi\circ\pi_{\text{P}})](y_{0})$

$=$ I$_{32621}+$ I$_{32622}:$

I$_{32621}=-\frac{1}{24}\underset{i=q+1}{\overset{n}{\sum}}%
\underset{j=1}{\overset{n}{\sum}}\underset{\text{b=1}}{\overset{\text{q}%
}{\sum}}[\frac{\partial^{2}\Gamma_{jj}^{\text{b}}}{\partial x_{i}^{2}}%
\frac{\partial\phi}{\partial\text{x}_{\text{b}}}$ $](y_{0})$

I$_{32622}=-\frac{1}{24}\underset{i=q+1}{\overset{n}{\sum}}%
\underset{j=1}{\overset{n}{\sum}}$ $\underset{l=q+1}{\overset{n}{\sum}}%
\frac{\partial^{2}}{\partial x_{i}^{2}}(\Gamma_{jj}^{l}\Lambda_{l})\phi
\circ\pi_{\text{P}})](y_{0})$

Then,

I$_{32621}=-\frac{1}{24}\underset{i=q+1}{\overset{n}{\sum}}%
\underset{\text{a,b=1}}{\overset{\text{q}}{\sum}}[\frac{\partial^{2}%
\Gamma_{\text{aa}}^{\text{b}}}{\partial x_{i}^{2}}\frac{\partial\phi}%
{\partial\text{x}_{\text{b}}}$ $](y_{0})-\frac{1}{24}%
\underset{i,j=q+1}{\overset{n}{\sum}}\underset{\text{b=1}}{\overset{\text{q}%
}{\sum}}[\frac{\partial^{2}\Gamma_{jj}^{\text{b}}}{\partial x_{i}^{2}}%
\frac{\partial\phi}{\partial\text{x}_{\text{b}}}$ $](y_{0})$

By (xi) of \textbf{Table A}$_{8}:$

$\frac{\partial^{2}\Gamma_{\text{aa}}^{\text{b}}}{\partial\text{x}_{i}^{2}%
}(y_{0})=$ $\overset{n}{\underset{k=q+1}{\sum}}$T$_{\text{aa}k}[\frac{8}%
{3}R_{i\text{c}ik}+4\underset{\text{d}=1}{\overset{q}{\sum}}(T_{\text{db}%
k})(\perp_{\text{d}ik})]$

$+$ $2\overset{n}{\underset{k,l=q+1}{\sum}}\perp_{\text{b}ik}(y_{0}%
)[-R_{\text{a}k\text{a}l}+\underset{\text{d=1}}{\overset{\text{q}}{\sum}%
}T_{\text{ad}k}T_{\text{ad}l}](y_{0})+2\overset{n}{\underset{k,l=q+1}{\sum}%
}\perp_{\text{b}ik}(y_{0})[\underset{r=q+1}{\overset{n}{\sum}}\perp
_{\text{a}kr}\perp_{\text{a}lr}](y_{0})$

By (x) of \textbf{Table A}$_{8}:$

$\frac{\partial^{2}\Gamma_{jj}^{\text{b}}}{\partial\text{x}_{i}^{2}}%
(y_{0})\frac{\partial\phi}{\partial\text{x}_{\text{b}}}(y_{0})=\frac{8}%
{3}\overset{q}{\underset{\text{c}=1}{\sum}}($T$_{\text{bc}i}R_{ij\text{c}%
j})(y_{0})\frac{\partial\phi}{\partial\text{x}_{\text{b}}}(y_{0})+\frac{2}%
{3}\overset{n}{\underset{k=q+1}{\sum}}(\perp_{\text{b}ik}R_{ijjk})(y_{0}%
)\frac{\partial\phi}{\partial\text{x}_{\text{b}}}(y_{0})$

$-\frac{1}{6}[4\overset{q}{\underset{\text{c=1}}{%
{\textstyle\sum}
}}$R$_{ij\text{c}i}^{{}}T_{\text{bc}j}+$ $4\overset{n}{\underset{k=q+1}{%
{\textstyle\sum}
}}$R$_{ijik}\perp_{\text{b}jk}+3\nabla_{i}$R$_{j\text{b}ij}%
+4\overset{q}{\underset{\text{c=1}}{%
{\textstyle\sum}
}}$R$_{ij\text{c}j}^{{}}T_{\text{bc}i}+$ $4$R$_{ijjk}\perp_{\text{b}ik}%
](y_{0})\frac{\partial\phi}{\partial\text{x}_{\text{b}}}(y_{0})$

$\left(  D_{26}\right)  \qquad$ I$_{32621}=-\frac{1}{24}%
\underset{i=q+1}{\overset{n}{\sum}}\underset{\text{a,b=1}}{\overset{\text{q}%
}{\sum}}$ $\overset{n}{\underset{k=q+1}{\sum}}$T$_{\text{aa}k}[\frac{8}%
{3}R_{i\text{c}ik}+4\underset{\text{d}=1}{\overset{q}{\sum}}(T_{\text{db}%
k})(\perp_{\text{d}ik})]\frac{\partial\phi}{\partial\text{x}_{\text{b}}}$
$(y_{0})$

$-\frac{1}{12}\underset{i=q+1}{\overset{n}{\sum}}\underset{\text{a,b=1}%
}{\overset{\text{q}}{\sum}}[$ $\overset{n}{\underset{k,l=q+1}{\sum}}%
\perp_{\text{b}ik}(-R_{\text{a}k\text{a}l}+\underset{\text{d=1}%
}{\overset{\text{q}}{\sum}}T_{\text{ad}k}T_{\text{ad}l}%
))-\overset{n}{\underset{k,l=q+1}{\sum}}\perp_{\text{b}ik}%
(\underset{r=q+1}{\overset{n}{\sum}}\perp_{\text{a}kr}\perp_{\text{a}%
lr})](y_{0})\frac{\partial\phi}{\partial\text{x}_{\text{b}}}$ $(y_{0})$

$-\frac{1}{12}\underset{i=q+1}{\overset{n}{\sum}}\underset{\text{b=1}%
}{\overset{\text{q}}{\sum}}[\frac{8}{3}\overset{q}{\underset{\text{c}=1}{\sum
}}($T$_{\text{bc}i}R_{ij\text{c}j})+\frac{2}{3}%
\overset{n}{\underset{k=q+1}{\sum}}(\perp_{\text{b}ik}R_{ijjk})](y_{0}%
)\frac{\partial\phi}{\partial\text{x}_{\text{b}}}(y_{0})$

$-\frac{1}{6}\underset{i=q+1}{\overset{n}{\sum}}\underset{\text{a,b=1}%
}{\overset{\text{q}}{\sum}}[4\overset{q}{\underset{\text{c=1}}{%
{\textstyle\sum}
}}$R$_{ij\text{c}i}^{{}}T_{\text{bc}j}+$ $4\overset{n}{\underset{k=q+1}{%
{\textstyle\sum}
}}$R$_{ijik}\perp_{\text{b}jk}+3\nabla_{i}$R$_{j\text{b}ij}$

$+4\overset{q}{\underset{\text{c=1}}{%
{\textstyle\sum}
}}$R$_{ij\text{c}j}^{{}}T_{\text{bc}i}+$ $4$R$_{ijjk}\perp_{\text{b}ik}%
](y_{0})\frac{\partial\phi}{\partial\text{x}_{\text{b}}}(y_{0})$

\qquad\qquad\qquad\qquad\qquad\qquad\qquad\qquad\qquad\qquad\qquad\qquad
\qquad\qquad\qquad\qquad\qquad$\blacksquare$

We next compute:

I$_{32622}=-\frac{1}{24}\underset{i=q+1}{\overset{n}{\sum}}%
\underset{j=1}{\overset{n}{\sum}}[$ $\underset{l=q+1}{\overset{n}{\sum}}%
\frac{\partial^{2}}{\partial x_{i}^{2}}(\Gamma_{jj}^{l}\Lambda_{l})\phi
\circ\pi_{\text{P}})](y_{0})$

$\qquad=\ \ \mathbf{-}$ $\frac{1}{24}\underset{i=q+1}{\overset{n}{\sum}%
}[\underset{j=1}{\overset{n}{\sum}}\left\{  \underset{l=q+1}{\overset{n}{\sum
}}(\frac{\partial^{2}\Gamma_{jj}^{l}}{\partial x_{i}^{2}}\Lambda_{l}\phi
\circ\pi_{\text{P}}+\Gamma_{jj}^{l}\frac{\partial^{2}\Lambda_{l}}{\partial
x_{i}^{2}}\phi\circ\pi_{\text{P}})\right\}  ](y_{0})$

$\qquad\qquad\mathbf{-}$ $\frac{1}{12}\underset{i=q+1}{\overset{n}{\sum}%
}[\underset{j=1}{\overset{n}{\sum}}\underset{l=q+1}{\overset{n}{\sum}}%
\frac{\partial\Gamma_{jj}^{l}}{\partial x_{i}}\frac{\partial\Lambda_{l}%
}{\partial x_{i}}\phi\circ\pi_{\text{P}}](y_{0})$

\ \ \qquad\ \ \ $=$ I$_{326221}+$ I$_{326222}+$ I$_{326223}$

where,

I$_{326221}=\ \ \mathbf{-}$ $\frac{1}{24}\underset{i,l=q+1}{\overset{n}{\sum}%
}[\underset{\text{a=1}}{\overset{q}{\sum}}\frac{\partial^{2}\Gamma_{\text{aa}%
}^{l}}{\partial x_{i}^{2}}\Lambda_{l}\phi](y_{0})$ $\mathbf{-}$ $\frac{1}%
{24}\underset{i,j,l=q+1}{\overset{n}{\sum}}[\frac{\partial^{2}\Gamma_{jj}^{l}%
}{\partial x_{i}^{2}}\Lambda_{l}\phi](y_{0})$

\ I$_{326222}=\ \ \mathbf{-}$ $\frac{1}{24}\underset{i,k=q+1}{\overset{n}{\sum
}}[\underset{\text{a=1}}{\overset{q}{\sum}}\Gamma_{\text{aa}}^{l}%
\frac{\partial^{2}\Lambda_{l}}{\partial x_{i}^{2}}\phi](y_{0})$ $\mathbf{-}$
$\frac{1}{24}\underset{i,j,k=q+1}{\overset{n}{\sum}}[\Gamma_{jj}^{l}%
\frac{\partial^{2}\Lambda_{l}}{\partial x_{i}^{2}}\phi](y_{0})$

I$_{326223}=\mathbf{-}$ $\frac{1}{12}\underset{i,l=q+1}{\overset{n}{\sum}%
}\underset{\text{a=1}}{\overset{q}{\sum}}\frac{\partial\Gamma_{\text{aa}}^{l}%
}{\partial x_{i}}\frac{\partial\Lambda_{l}}{\partial x_{i}}\phi](y_{0}%
)\mathbf{-}$ $\frac{1}{12}\underset{i,j,l=q+1}{\overset{n}{\sum}}%
[\frac{\partial\Gamma_{jj}^{l}}{\partial x_{i}}\frac{\partial\Lambda_{l}%
}{\partial x_{i}}\phi](y_{0})$

We have:

I$_{326221}=$ $\ \ \mathbf{-}$ $\frac{1}{24}%
[\underset{i,k=q+1}{\overset{n}{\sum}}\underset{\text{a}=1}{\overset{q}{\sum}%
}\frac{\partial^{2}\Gamma_{\text{aa}}^{k}}{\partial x_{i}^{2}}\Lambda
_{k}](y_{0})\phi(y_{0})\mathbf{-}$ $\frac{1}{24}%
[\overset{n}{\underset{i,j,k=q+1}{\sum}}\frac{\partial^{2}\Gamma_{jj}^{k}%
}{\partial x_{i}^{2}}\Lambda_{k}](y_{0})\phi(y_{0})$

$\frac{\partial^{2}\Gamma_{\text{aa}}^{k}}{\partial x_{i}^{2}}(y_{0})$ is from
(vii) of \textbf{Table A}$_{8}$ and $\frac{\partial^{2}\Gamma_{jj}^{k}%
}{\partial x_{i}^{2}}(y_{0})$ is from (ix) of \textbf{Table A}$_{8}.$

We thus have from these expressions:

I$_{326221}=-\frac{1}{144}[\{4\nabla_{i}$R$_{i\text{a}j\text{a}}$
$+2\nabla_{j}$R$_{i\text{a}i\text{a}}+$ $8$ $(\overset{q}{\underset{\text{c=1}%
}{%
{\textstyle\sum}
}}R_{\text{a}i\text{c}i}^{{}}T_{\text{ac}j}+\;\overset{n}{\underset{k=q+1}{%
{\textstyle\sum}
}}R_{\text{a}iik}\perp_{\text{a}jk})$

\ \ \ \ \ \ \ \ \ \ \ \ \qquad\ \ $+8(\overset{q}{\underset{\text{c=1}}{%
{\textstyle\sum}
}}R_{\text{a}i\text{c}j}^{{}}T_{\text{ac}i}+\overset{n}{\underset{l=q+1}{%
{\textstyle\sum}
}}R_{\text{a}ijl}\perp_{\text{a}il})+8(\overset{q}{\underset{\text{c=1}}{%
{\textstyle\sum}
}}R_{\text{a}j\text{c}i}^{{}}T_{\text{ac}i}+\overset{q}{\underset{\text{c=1}}{%
{\textstyle\sum}
}}R_{\text{a}j\text{c}i}^{{}}T_{\text{ac}i})\}$\ 

\qquad$\qquad\ \ \ \ \ \ +\frac{2}{3}\underset{k=q+1}{\overset{n}{\sum}%
}\{T_{\text{aa}k}(R_{ijik}+3\overset{q}{\underset{\text{c}=1}{\sum}}%
\perp_{\text{c}ij}\perp_{\text{c}ik})\}](y_{0})\Lambda_{k}(y_{0})\phi(y_{0})$

\qquad\qquad\ $+\frac{1}{24}[\frac{4}{3}\overset{q}{\underset{\text{a}%
=1}{\sum}}\perp_{\text{a}ik}R_{ij\text{a}j}+\frac{1}{3}(\nabla_{i}%
R_{kjij}+\nabla_{j}R_{ijik}+\nabla_{k}R_{ijij})](y_{0})\Lambda_{k}(y_{0}%
)\phi(y_{0})$

\qquad\qquad\qquad\qquad\qquad\qquad\qquad\qquad\qquad\qquad\qquad\qquad
\qquad\qquad\qquad\qquad\qquad$\blacksquare$

Next we have:

\qquad I$_{326222}=\mathbf{-}$ $\frac{1}{24}\underset{i=q+1}{\overset{n}{\sum
}}[\underset{j=1}{\overset{n}{\sum}}\underset{k=1}{\overset{n}{\sum}}%
\Gamma_{jj}^{k}\frac{\partial^{2}\Lambda_{k}}{\partial x_{i}^{2}}\phi\circ
\pi_{\text{P}}](y_{0})$

\qquad$=\mathbf{-}$ $\frac{1}{24}\underset{i=q+1}{\overset{n}{\sum}%
}[\underset{\text{a=1}}{\overset{q}{\sum}}\underset{k=1}{\overset{n}{\sum}%
}\Gamma_{\text{aa}}^{k}\frac{\partial^{2}\Lambda_{k}}{\partial x_{i}^{2}}%
\phi\circ\pi_{\text{P}}](y_{0})\mathbf{-}$ $\frac{1}{24}%
\underset{i=q+1}{\overset{n}{\sum}}[\underset{j=q+1}{\overset{n}{\sum}%
}\underset{k=1}{\overset{n}{\sum}}\Gamma_{jj}^{k}\frac{\partial^{2}\Lambda
_{k}}{\partial x_{i}^{2}}\phi\circ\pi_{\text{P}}](y_{0})$

\qquad$=\mathbf{-}$ $\frac{1}{24}\underset{i=q+1}{\overset{n}{\sum}%
}[\underset{\text{a,b=1}}{\overset{q}{\sum}}\Gamma_{\text{aa}}^{\text{b}}%
\frac{\partial^{2}\Lambda_{\text{b}}}{\partial x_{i}^{2}}\phi\circ
\pi_{\text{P}}](y_{0})\mathbf{-}$ $\frac{1}{24}%
\underset{i,k=q+1}{\overset{n}{\sum}}\underset{\text{a=1}}{\overset{q}{\sum}%
}[\Gamma_{\text{aa}}^{k}\frac{\partial^{2}\Lambda_{k}}{\partial x_{i}^{2}}%
\phi\circ\pi_{\text{P}}](y_{0})$

\qquad$\mathbf{-}$ $\frac{1}{24}\underset{i,j=q+1}{\overset{n}{\sum}%
}[\underset{\text{a=1}}{\overset{q}{\sum}}\Gamma_{jj}^{\text{a}}\frac
{\partial^{2}\Lambda_{\text{a}}}{\partial x_{i}^{2}}\phi\circ\pi_{\text{P}%
}](y_{0})\mathbf{-}$ $\frac{1}{24}\underset{i,j,k=q+1}{\overset{n}{\sum}%
}[\Gamma_{jj}^{k}\frac{\partial^{2}\Lambda_{k}}{\partial x_{i}^{2}}\phi
\circ\pi_{\text{P}}](y_{0})$

$\Gamma_{\text{aa}}^{\text{b}}(y_{0})=0$ by (ii) of \textbf{Table A}$_{7},$
$\Gamma_{jj}^{\text{a}}(y_{0})=0$ by (ii) of \textbf{Table A}$_{8}$ and
$\Gamma_{jj}^{k}(y_{0})=0$ by (i) of \textbf{Table A}$_{8}$.

Therefore the last equality above reduces to:

I$_{326222}=\mathbf{-}$ $\frac{1}{24}\underset{i,j=q+1}{\overset{n}{\sum}%
}[\underset{\text{a}=1}{\overset{q}{\sum}}\Gamma_{\text{aa}}^{j}\frac
{\partial^{2}\Lambda_{j}}{\partial x_{i}^{2}}]\phi(y_{0})$

$\Gamma_{\text{aa}}^{k}(y_{0})=T_{\text{aa}k}(y_{0})$ by (i) Table A$_{8}$ and
$\frac{\partial^{2}\Lambda_{j}}{\partial\text{x}_{i}^{2}}(y_{0})=\frac{1}%
{3}\frac{\partial\Omega_{ij}}{\partial\text{x}_{i}}(y_{0})$ by (xi) of
\textbf{Proposition 5. }Therefore,

I$_{326222}=$ $\mathbf{-}$ $\frac{1}{72}\underset{i,j=q+1}{\overset{n}{\sum}%
}\underset{\text{a}=1}{\overset{q}{\sum}}T_{\text{aa}j}(y_{0})\frac
{\partial\Omega_{ij}}{\partial\text{x}_{i}}(y_{0})\phi(y_{0}).$

Next we have:

I$_{326223}=\mathbf{-}$ $\frac{1}{12}\underset{i=q+1}{\overset{n}{\sum}%
}[\underset{j=1}{\overset{n}{\sum}}\underset{k=1}{\overset{n}{\sum}}%
\frac{\partial\Gamma_{jj}^{k}}{\partial x_{i}}\frac{\partial\Lambda_{k}%
}{\partial x_{i}}\phi\circ\pi_{\text{P}}](y_{0})$

$=\mathbf{-}$ $\frac{1}{12}\underset{i=q+1}{\overset{n}{\sum}}%
[\underset{j=1}{\overset{n}{\sum}}\underset{\text{b=1}}{\overset{q}{\sum}%
}\frac{\partial\Gamma_{jj}^{\text{b}}}{\partial x_{i}}\frac{\partial
\Lambda_{\text{b}}}{\partial x_{i}}\phi\circ\pi_{\text{P}}](y_{0})$
$\mathbf{-}$ $\frac{1}{12}\underset{i=q+1}{\overset{n}{\sum}}%
[\underset{j=1}{\overset{n}{\sum}}\underset{k=q+1}{\overset{n}{\sum}}%
\frac{\partial\Gamma_{jj}^{k}}{\partial x_{i}}\frac{\partial\Lambda_{k}%
}{\partial x_{i}}\phi\circ\pi_{\text{P}}](y_{0})$

$=\mathbf{-}$ $\frac{1}{12}\underset{i=q+1}{\overset{n}{\sum}}%
\underset{\text{a,b=1}}{\overset{q}{\sum}}\frac{\partial\Gamma_{\text{aa}%
}^{\text{b}}}{\partial x_{i}}\frac{\partial\Lambda_{\text{b}}}{\partial x_{i}%
}\phi\circ\pi_{\text{P}}](y_{0})$ $\mathbf{-}$ $\frac{1}{12}%
\underset{i,j=q+1}{\overset{n}{\sum}}[\underset{\text{b=1}}{\overset{q}{\sum}%
}\frac{\partial\Gamma_{jj}^{\text{b}}}{\partial x_{i}}\frac{\partial
\Lambda_{\text{b}}}{\partial x_{i}}\phi\circ\pi_{\text{P}}](y_{0})$

$\mathbf{-}$ $\frac{1}{12}\underset{i,j=q+1}{\overset{n}{\sum}}%
[\underset{\text{a=1}}{\overset{q}{\sum}}\frac{\partial\Gamma_{\text{aa}}^{j}%
}{\partial x_{i}}\frac{\partial\Lambda_{j}}{\partial x_{i}}\phi\circ
\pi_{\text{P}}](y_{0})$ $\mathbf{-}$ $\frac{1}{12}%
\underset{i,j,k=q+1}{\overset{n}{\sum}}[\frac{\partial\Gamma_{jj}^{k}%
}{\partial x_{i}}\frac{\partial\Lambda_{k}}{\partial x_{i}}\phi\circ
\pi_{\text{P}}](y_{0})$

$\frac{\partial\Gamma_{\text{aa}}^{\text{b}}}{\partial\text{x}_{i}}%
(y_{0})=-\overset{n}{\underset{j=q+1}{\sum}}(\perp_{\text{b}ij}$%
T$_{\text{aa}j})(y_{0})$ by (v) of \textbf{TableA}$_{7}.$

$\frac{\partial\Gamma_{jj}^{\text{b}}}{\partial\text{x}_{i}}(y_{0})=\frac
{2}{3}R_{\text{b}jij}(y_{0})$ by (v) of \textbf{Table A}$_{8}.$

\ $\frac{\partial\Gamma_{\text{aa}}^{j}}{\partial x_{i}}(y_{0})=[$
R$_{\text{a}i\text{a}j}$ $-\underset{\text{c=1}}{\overset{\text{q}}{\sum}%
}T_{\text{ac}i}T_{\text{ac}j}-\overset{n}{\underset{k=q+1}{\sum}}%
(\perp_{\text{a}ik}\perp_{\text{a}jk}](y_{0})$ by (iv) of \textbf{Table
A}$_{7}.$

$\frac{\partial\Gamma_{jj}^{k}}{\partial\text{x}_{i}}(y_{0})=$ $\frac{2}%
{3}R_{ijkj}(y_{0})$ by (viii) of Table A$_{8}$

$\frac{\partial\Lambda_{\text{a}}}{\partial x_{i}}(y_{0})=\Omega_{i\text{a}%
}(y_{0})+[\Lambda_{\text{a}},\Lambda_{i}](y_{0})$ by (vii) \textbf{Proposition
5.}

$\frac{\partial\Lambda_{j}}{\partial x_{i}}(y_{0})=\frac{1}{2}\Omega
_{ij}(y_{0})$ by (vi) of \textbf{Proposition 5.}

Consequently we have:

\qquad I$_{326223}=$ $\frac{1}{12}$ $\overset{n}{\underset{j=q+1}{\sum}}%
($T$_{\text{aa}j}\perp_{\text{b}ij})(y_{0})\left(  \Omega_{i\text{b}}%
(y_{0})+[\Lambda_{\text{b}},\Lambda_{i}]\right)  (y_{0})\phi(y_{0})$

\qquad$\mathbf{-}\frac{1}{18}$ $\underset{i,j=q+1}{\overset{n}{\sum}%
}\underset{\text{b=1}}{\overset{q}{\sum}}R_{\text{b}jij}(y_{0})\left(
\Omega_{i\text{a}}(y_{0})+[\Lambda_{\text{a}},\Lambda_{i}]\right)  (y_{0}%
)\phi(y_{0})$

\qquad\ $\mathbf{-}$ $\frac{1}{24}\underset{i,j=q+1}{\overset{n}{\sum}%
}\underset{\text{a=1}}{\overset{q}{\sum}}[$ R$_{\text{a}i\text{a}j}$
$-\underset{\text{c=1}}{\overset{\text{q}}{\sum}}T_{\text{ac}i}T_{\text{ac}%
j}-\overset{n}{\underset{k=q+1}{\sum}}(\perp_{\text{a}ik}\perp_{\text{a}%
jk}](y_{0})\Omega_{ij}(y_{0})\phi(y_{0})$

\qquad$\mathbf{-}$ $\frac{1}{36}\underset{i,j,k=q+1}{\overset{n}{\sum}%
}R_{ijkj}(y_{0})\Omega_{ik}(y_{0})(y_{0})\phi(y_{0})$

\qquad\qquad\qquad\qquad\qquad\qquad\qquad\qquad\qquad\qquad\qquad\qquad
\qquad\qquad\qquad\qquad\qquad\qquad$\blacksquare$

Therefore,

$\left(  D_{27}\right)  \qquad$I$_{32622}=$ I$_{326221}+$ I$_{326222}+$
I$_{326223}$

$\qquad\qquad=-\frac{1}{144}[\{4\nabla_{i}$R$_{i\text{a}j\text{a}}$
$+2\nabla_{j}$R$_{i\text{a}i\text{a}}+$ $8$ $(\overset{q}{\underset{\text{c=1}%
}{%
{\textstyle\sum}
}}R_{\text{a}i\text{c}i}^{{}}T_{\text{ac}j}+\;\overset{n}{\underset{k=q+1}{%
{\textstyle\sum}
}}R_{\text{a}iik}\perp_{\text{a}jk})\qquad$I$_{326221}$

\ \ \ \ \ \ \ \ \ \ \ \ \qquad\ \ $+8(\overset{q}{\underset{\text{c=1}}{%
{\textstyle\sum}
}}R_{\text{a}i\text{c}j}^{{}}T_{\text{ac}i}+\overset{n}{\underset{l=q+1}{%
{\textstyle\sum}
}}R_{\text{a}ijl}\perp_{\text{a}il})+8(\overset{q}{\underset{\text{c=1}}{%
{\textstyle\sum}
}}R_{\text{a}j\text{c}i}^{{}}T_{\text{ac}i}+\overset{q}{\underset{\text{c=1}}{%
{\textstyle\sum}
}}R_{\text{a}j\text{c}i}^{{}}T_{\text{ac}i})\}$\ 

\qquad$\qquad\ \ \ \ \ \ +\frac{2}{3}\underset{k=q+1}{\overset{n}{\sum}%
}\{T_{\text{aa}k}(R_{ijik}+3\overset{q}{\underset{\text{c}=1}{\sum}}%
\perp_{\text{c}ij}\perp_{\text{c}ik})\}](y_{0})\Lambda_{k}(y_{0})\phi(y_{0})$

\qquad\qquad$\ \mathbf{-}$ $\frac{1}{24}[\frac{8}{3}%
\overset{q}{\underset{\text{c}=1}{\sum}}($T$_{\text{bc}i}R_{ij\text{c}%
j})+\frac{2}{3}\overset{n}{\underset{k=q+1}{\sum}}(\perp_{\text{b}ik}%
R_{ijjk})](y_{0})\Lambda_{\text{b}}(y_{0})\phi(y_{0})$

\qquad\qquad\ $+\frac{1}{24}[\frac{4}{3}\overset{q}{\underset{\text{a}%
=1}{\sum}}\perp_{\text{a}ik}R_{ij\text{a}j}+\frac{1}{3}(\nabla_{i}%
R_{kjij}+\nabla_{j}R_{ijik}+\nabla_{k}R_{ijij})](y_{0})\Lambda_{k}(y_{0}%
)\phi(y_{0})$

\qquad\qquad\ $\mathbf{-}$ $\frac{1}{72}\underset{i,j=q+1}{\overset{n}{\sum}%
}\underset{\text{a}=1}{\overset{q}{\sum}}T_{\text{aa}j}(y_{0})\frac
{\partial\Omega_{ij}}{\partial\text{x}_{i}}(y_{0})\phi(y_{0})\qquad
$I$_{326222}$

\qquad$\qquad+$ $\frac{1}{12}$ $\overset{n}{\underset{j=q+1}{\sum}}%
(\perp_{\text{b}ij}$T$_{\text{aa}j})(y_{0})\left(  \Omega_{i\text{b}}%
(y_{0})+[\Lambda_{\text{b}},\Lambda_{i}]\right)  (y_{0})\phi(y_{0})$ \qquad
I$_{326223}$

\qquad\qquad\ \ \ $\mathbf{-}\frac{1}{18}$
$\underset{i,j=q+1}{\overset{n}{\sum}}\underset{\text{b=1}}{\overset{q}{\sum}%
}R_{\text{b}jij}(y_{0})\left(  \Omega_{i\text{a}}(y_{0})+[\Lambda_{\text{a}%
},\Lambda_{i}]\right)  (y_{0})\phi(y_{0})$

\qquad\ \ \ $\qquad\mathbf{-}$ $\frac{1}{24}%
\underset{i,j=q+1}{\overset{n}{\sum}}\underset{\text{a=1}}{\overset{q}{\sum}%
}[$ R$_{\text{a}i\text{a}j}$ $-\underset{\text{c=1}}{\overset{\text{q}}{\sum}%
}T_{\text{ac}i}T_{\text{ac}j}-\overset{n}{\underset{k=q+1}{\sum}}%
(\perp_{\text{a}ik}\perp_{\text{a}jk}](y_{0})\Omega_{ij}(y_{0})\phi(y_{0})$

\qquad\qquad\ $\mathbf{-}$ $\frac{1}{36}\underset{i,j,k=q+1}{\overset{n}{\sum
}}R_{ijkj}(y_{0})\Omega_{ik}(y_{0})(y_{0})\phi(y_{0})$

\qquad\qquad\qquad\qquad\qquad\qquad\qquad\qquad\qquad\qquad\qquad\qquad
\qquad\qquad\qquad\qquad\qquad$\qquad\blacksquare$

We have from $\left(  D_{26}\right)  $ and $\left(  D_{27}\right)  :$

$\left(  D_{28}\right)  \qquad$I$_{3262}=$ I$_{32621}+$ I$_{32622}$

$=-\frac{1}{24}\underset{i=q+1}{\overset{n}{\sum}}\underset{\text{a,b=1}%
}{\overset{\text{q}}{\sum}}$ $\overset{n}{\underset{k=q+1}{\sum}}%
$T$_{\text{aa}k}[\frac{8}{3}R_{i\text{c}ik}+4\underset{\text{d}%
=1}{\overset{q}{\sum}}(T_{\text{db}k})(\perp_{\text{d}ik})]\frac{\partial\phi
}{\partial\text{x}_{\text{b}}}$ $(y_{0})\qquad$ I$_{32621}$

$-\frac{1}{12}\underset{i=q+1}{\overset{n}{\sum}}\underset{\text{a,b=1}%
}{\overset{\text{q}}{\sum}}[$ $\overset{n}{\underset{k,l=q+1}{\sum}}%
\perp_{\text{b}ik}(-R_{\text{a}k\text{a}l}+\underset{\text{d=1}%
}{\overset{\text{q}}{\sum}}T_{\text{ad}k}T_{\text{ad}l}%
))-\overset{n}{\underset{k,l=q+1}{\sum}}\perp_{\text{b}ik}%
(\underset{r=q+1}{\overset{n}{\sum}}\perp_{\text{a}kr}\perp_{\text{a}%
lr})](y_{0})\frac{\partial\phi}{\partial\text{x}_{\text{b}}}$ $(y_{0})$

\qquad$-\frac{1}{12}\underset{i=q+1}{\overset{n}{\sum}}\underset{\text{b=1}%
}{\overset{\text{q}}{\sum}}[\frac{8}{3}\overset{q}{\underset{\text{c}=1}{\sum
}}($T$_{\text{bc}i}R_{ij\text{c}j})+\frac{2}{3}%
\overset{n}{\underset{k=q+1}{\sum}}(\perp_{\text{b}ik}R_{ijjk})](y_{0}%
)\frac{\partial\phi}{\partial\text{x}_{\text{b}}}(y_{0})$

$-\frac{1}{6}\underset{i=q+1}{\overset{n}{\sum}}\underset{\text{a,b=1}%
}{\overset{\text{q}}{\sum}}[4\overset{q}{\underset{\text{c=1}}{%
{\textstyle\sum}
}}$R$_{ij\text{c}i}^{{}}T_{\text{bc}j}+$ $4\overset{n}{\underset{k=q+1}{%
{\textstyle\sum}
}}$R$_{ijik}\perp_{\text{b}jk}+3\nabla_{i}$R$_{j\text{b}ij}%
+4\overset{q}{\underset{\text{c=1}}{%
{\textstyle\sum}
}}$R$_{ij\text{c}j}^{{}}T_{\text{bc}i}+$ $4$R$_{ijjk}\perp_{\text{b}ik}%
](y_{0})\frac{\partial\phi}{\partial\text{x}_{\text{b}}}(y_{0})$

$-\frac{1}{144}[\{4\nabla_{i}$R$_{i\text{a}j\text{a}}$ $+2\nabla_{j}%
$R$_{i\text{a}i\text{a}}+$ $8$ $(\overset{q}{\underset{\text{c=1}}{%
{\textstyle\sum}
}}R_{\text{a}i\text{c}i}^{{}}T_{\text{ac}j}+\;\overset{n}{\underset{k=q+1}{%
{\textstyle\sum}
}}R_{\text{a}iik}\perp_{\text{a}jk})\qquad$I$_{326221}\qquad$I$_{32622}$

\ $+8(\overset{q}{\underset{\text{c=1}}{%
{\textstyle\sum}
}}R_{\text{a}i\text{c}j}^{{}}T_{\text{ac}i}+\overset{n}{\underset{l=q+1}{%
{\textstyle\sum}
}}R_{\text{a}ijl}\perp_{\text{a}il})+8(\overset{q}{\underset{\text{c=1}}{%
{\textstyle\sum}
}}R_{\text{a}j\text{c}i}^{{}}T_{\text{ac}i}+\overset{q}{\underset{\text{c=1}}{%
{\textstyle\sum}
}}R_{\text{a}j\text{c}i}^{{}}T_{\text{ac}i})\}$\ 

$+\frac{2}{3}\underset{k=q+1}{\overset{n}{\sum}}\{T_{\text{aa}k}%
(R_{ijik}+3\overset{q}{\underset{\text{c}=1}{\sum}}\perp_{\text{c}ij}%
\perp_{\text{c}ik})\}](y_{0})\Lambda_{k}(y_{0})\phi(y_{0})$

$+\frac{1}{24}[\frac{4}{3}\overset{q}{\underset{\text{a}=1}{\sum}}%
\perp_{\text{a}ik}R_{ij\text{a}j}+\frac{1}{3}(\nabla_{i}R_{kjij}+\nabla
_{j}R_{ijik}+\nabla_{k}R_{ijij})](y_{0})\Lambda_{k}(y_{0})\phi(y_{0})$

$\mathbf{-}$ $\frac{1}{72}\underset{i,j=q+1}{\overset{n}{\sum}}%
\underset{\text{a}=1}{\overset{q}{\sum}}T_{\text{aa}j}(y_{0})\frac
{\partial\Omega_{ij}}{\partial\text{x}_{i}}(y_{0})\phi(y_{0})\qquad
$I$_{326222}$

$+$ $\frac{1}{12}$ $\overset{n}{\underset{j=q+1}{\sum}}(\perp_{\text{b}ij}%
$T$_{\text{aa}j})(y_{0})\left(  \Omega_{i\text{b}}(y_{0})+[\Lambda_{\text{b}%
},\Lambda_{i}]\right)  (y_{0})\phi(y_{0})$ \qquad I$_{326223}$

$\mathbf{-}\frac{1}{18}$ $\underset{i,j=q+1}{\overset{n}{\sum}}%
\underset{\text{b=1}}{\overset{q}{\sum}}R_{\text{b}jij}(y_{0})\left(
\Omega_{i\text{a}}(y_{0})+[\Lambda_{\text{a}},\Lambda_{i}]\right)  (y_{0}%
)\phi(y_{0})$

$\mathbf{-}$ $\frac{1}{24}\underset{i,j=q+1}{\overset{n}{\sum}}%
\underset{\text{a=1}}{\overset{q}{\sum}}[$ R$_{\text{a}i\text{a}j}$
$-\underset{\text{c=1}}{\overset{\text{q}}{\sum}}T_{\text{ac}i}T_{\text{ac}%
j}-\overset{n}{\underset{k=q+1}{\sum}}(\perp_{\text{a}ik}\perp_{\text{a}%
jk}](y_{0})\Omega_{ij}(y_{0})\phi(y_{0})$

$\mathbf{-}$ $\frac{1}{36}\underset{i,j,k=q+1}{\overset{n}{\sum}}%
R_{ijkj}(y_{0})\Omega_{ik}(y_{0})(y_{0})\phi(y_{0})$

\qquad\qquad\qquad\qquad\qquad\qquad\qquad\qquad\qquad\qquad\qquad\qquad
\qquad\qquad\qquad\qquad\qquad$\blacksquare$

We next compute \textbf{I}$_{3263}:$

\textbf{I}$_{3263}=-\frac{1}{12}\underset{i=q+1}{\overset{n}{\sum}%
}[\underset{j,k=1}{\overset{n}{\sum}}\frac{\partial\text{g}^{jk}}{\partial
x_{i}}\left\{  \underset{\text{c=1}}{\overset{\text{q}}{\sum}}\frac
{\partial\Gamma_{jk}^{\text{c}}}{\partial x_{i}}\frac{\partial\phi}%
{\partial\text{x}_{\text{c}}}\circ\pi_{\text{P}}\text{ }+\text{ }%
\underset{l=1}{\overset{n}{\sum}}\frac{\partial}{\partial x_{i}}(\Gamma
_{jk}^{l}\Lambda_{l}\phi\circ\pi_{\text{P}})\right\}  ](y_{0})$

\qquad$=-\frac{1}{12}\underset{i=q+1}{\overset{n}{\sum}}%
[\underset{k=1}{\overset{n}{\sum}}\underset{\text{a}=1}{\overset{q}{\sum}%
}\frac{\partial\text{g}^{\text{a}k}}{\partial x_{i}}\left\{
\underset{\text{c=1}}{\overset{\text{q}}{\sum}}\frac{\partial\Gamma
_{\text{a}k}^{\text{c}}}{\partial x_{i}}\frac{\partial\phi}{\partial
\text{x}_{\text{c}}}\circ\pi_{\text{P}}\text{ }+\text{ }%
\underset{l=1}{\overset{n}{\sum}}\frac{\partial}{\partial x_{i}}%
(\Gamma_{\text{a}k}^{l}\Lambda_{l}\phi\circ\pi_{\text{P}})\right\}  ](y_{0})$

\qquad$\ -\frac{1}{12}\underset{i,j=q+1}{\overset{n}{\sum}}%
[\underset{k=1}{\overset{n}{\sum}}\frac{\partial\text{g}^{jk}}{\partial x_{i}%
}\left\{  \underset{\text{c=1}}{\overset{\text{q}}{\sum}}\frac{\partial
\Gamma_{jk}^{\text{c}}}{\partial x_{i}}\frac{\partial\phi}{\partial
\text{x}_{\text{c}}}\circ\pi_{\text{P}}\text{ }+\text{ }%
\underset{l=1}{\overset{n}{\sum}}\frac{\partial}{\partial x_{i}}(\Gamma
_{jk}^{l}\Lambda_{l}\phi\circ\pi_{\text{P}})\right\}  ](y_{0})$

\qquad$=-\frac{1}{12}\underset{i=q+1}{\overset{n}{\sum}}[\underset{\text{a,b}%
=1}{\overset{q}{\sum}}\frac{\partial\text{g}^{\text{ab}}}{\partial x_{i}%
}\left\{  \underset{\text{c=1}}{\overset{\text{q}}{\sum}}\frac{\partial
\Gamma_{\text{ab}}^{\text{c}}}{\partial x_{i}}\frac{\partial\phi}%
{\partial\text{x}_{\text{c}}}\circ\pi_{\text{P}}\text{ }+\text{ }%
\underset{l=1}{\overset{n}{\sum}}\frac{\partial}{\partial x_{i}}%
(\Gamma_{\text{ab}}^{l}\Lambda_{l}\phi\circ\pi_{\text{P}})\right\}  ](y_{0})$

$\qquad-\frac{1}{6}\underset{i,j=q+1}{\overset{n}{\sum}}[\underset{\text{a}%
=1}{\overset{q}{\sum}}\frac{\partial\text{g}^{\text{a}j}}{\partial x_{i}%
}\left\{  \underset{\text{c=1}}{\overset{\text{q}}{\sum}}\frac{\partial
\Gamma_{\text{a}j}^{\text{c}}}{\partial x_{i}}\frac{\partial\phi}%
{\partial\text{x}_{\text{c}}}\circ\pi_{\text{P}}\text{ }+\text{ }%
\underset{l=1}{\overset{n}{\sum}}\frac{\partial}{\partial x_{i}}%
(\Gamma_{\text{a}j}^{l}\Lambda_{l}\phi\circ\pi_{\text{P}})\right\}  ](y_{0})$

$\ \ -\frac{1}{12}\underset{i,j,k=q+1}{\overset{n}{\sum}}[\frac{\partial
\text{g}^{jk}}{\partial x_{i}}\left\{  \underset{\text{c=1}}{\overset{\text{q}%
}{\sum}}\frac{\partial\Gamma_{jk}^{\text{c}}}{\partial x_{i}}\frac
{\partial\phi}{\partial\text{x}_{\text{c}}}\circ\pi_{\text{P}}\text{ }+\text{
}\underset{l=1}{\overset{n}{\sum}}\frac{\partial}{\partial x_{i}}(\Gamma
_{jk}^{l}\Lambda_{l}\phi\circ\pi_{\text{P}})\right\}  ](y_{0})$

Since $\frac{\partial\text{g}^{jk}}{\partial x_{i}}(y_{0})=0$ for
$i,j,k=q+1,...,n,$ we have,

\textbf{I}$_{3263}=-\frac{1}{12}\underset{i=q+1}{\overset{n}{\sum}%
}[\underset{\text{a,b}=1}{\overset{q}{\sum}}\frac{\partial\text{g}^{\text{ab}%
}}{\partial x_{i}}\left\{  \underset{\text{c=1}}{\overset{\text{q}}{\sum}%
}\frac{\partial\Gamma_{\text{ab}}^{\text{c}}}{\partial x_{i}}\frac
{\partial\phi}{\partial\text{x}_{\text{c}}}\circ\pi_{\text{P}}\text{ }+\text{
}\underset{l=1}{\overset{n}{\sum}}\frac{\partial}{\partial x_{i}}%
(\Gamma_{\text{ab}}^{l}\Lambda_{l}\phi\circ\pi_{\text{P}})\right\}  ](y_{0})$

$\qquad\qquad-\frac{1}{6}\underset{i,j=q+1}{\overset{n}{\sum}}%
\underset{\text{a}=1}{\overset{q}{\sum}}\frac{\partial\text{g}^{\text{a}j}%
}{\partial x_{i}}[\underset{\text{c=1}}{\overset{\text{q}}{\sum}}%
\frac{\partial\Gamma_{\text{a}j}^{\text{c}}}{\partial x_{i}}\frac{\partial
\phi}{\partial\text{x}_{\text{c}}}\circ\pi_{\text{P}}$ $+$
$\underset{l=1}{\overset{n}{\sum}}\frac{\partial}{\partial x_{i}}%
(\Gamma_{\text{a}j}^{l}\Lambda_{l}\phi\circ\pi_{\text{P}})](y_{0})$

$\qquad\ \ =$ \textbf{I}$_{32631}+$ \textbf{I}$_{32632}$

\textbf{I}$_{32631}=-\frac{1}{12}\underset{i=q+1}{\overset{n}{\sum}%
}\underset{\text{a,b}=1}{\overset{q}{\sum}}\frac{\partial\text{g}^{\text{ab}}%
}{\partial x_{i}}[\underset{\text{c=1}}{\overset{\text{q}}{\sum}}%
\frac{\partial\Gamma_{\text{ab}}^{\text{c}}}{\partial x_{i}}\frac{\partial
\phi}{\partial\text{x}_{\text{c}}}\circ\pi_{\text{P}}$ $+$
$\underset{l=1}{\overset{n}{\sum}}\frac{\partial}{\partial x_{i}}%
(\Gamma_{\text{ab}}^{l}\Lambda_{l}\phi\circ\pi_{\text{P}})](y_{0})$

\qquad$\ =-\frac{1}{12}\underset{i=q+1}{\overset{n}{\sum}}\underset{\text{a,b}%
=1}{\overset{q}{\sum}}\frac{\partial\text{g}^{\text{ab}}}{\partial x_{i}%
}[\underset{\text{c=1}}{\overset{\text{q}}{\sum}}\frac{\partial\Gamma
_{\text{ab}}^{\text{c}}}{\partial x_{i}}\frac{\partial\phi}{\partial
\text{x}_{\text{c}}}$ $+$ $\underset{\text{c}=1}{\overset{q}{\sum}}%
(\frac{\partial\Gamma_{\text{ab}}^{\text{c}}}{\partial x_{i}}\Lambda
_{\text{c}}\phi+\Gamma_{\text{ab}}^{\text{c}}\frac{\partial\Lambda_{\text{c}}%
}{\partial x_{i}}\phi)$

$\qquad\qquad+$ $\underset{j=q+1}{\overset{n}{\sum}}(\frac{\partial
\Gamma_{\text{ab}}^{j}}{\partial x_{i}}\Lambda_{j}\phi+\Gamma_{\text{ab}}%
^{j}\frac{\partial\Lambda_{j}}{\partial x_{i}}\phi)](y_{0})$

We have:

$\frac{\partial\text{g}^{\text{ab}}}{\partial x_{i}}(y_{0})=2T_{\text{ab}%
i}(y_{0})$ by (ii) of \textbf{Table A}$_{6};$ $\Gamma_{\text{ab}}^{\text{c}%
}(y_{0})=0$ by (ii) of \textbf{Table A}$_{7};$ $\Gamma_{\text{ab}}^{j}%
(y_{0})=T_{\text{ab}j}(y_{0})$ by \textbf{Table} \textbf{A}$_{7}.$

By (iii ) of \textbf{Table A}$_{7},$

$\frac{\partial\Gamma_{\text{ab}}^{j}}{\partial\text{x}_{i}}(y_{0})=\frac
{1}{2}[$ $(R_{\text{a}i\text{b}j}+R_{\text{a}j\text{b}i})$
$-\underset{\text{c=1}}{\overset{\text{q}}{\sum}}(T_{\text{ac}i}T_{\text{bc}%
j}+T_{\text{ac}j}T_{\text{bc}i})-\overset{n}{\underset{k=q+1}{\sum}}%
(\perp_{\text{a}ik}\perp_{\text{b}jk}+$ $\perp_{\text{a}jk}\perp_{\text{b}%
ik})](y_{0})$

$\frac{\partial\Gamma_{\text{ab}}^{\text{c}}}{\partial x_{i}}(y_{0}%
)=-\overset{n}{\underset{k=q+1}{\sum}}(\perp_{\text{c}ik}$T$_{\text{ab}%
k})(y_{0})$ by (v) of \textbf{Table A}$_{7};$Therefore,

\textbf{I}$_{32631}=\frac{1}{6}\underset{i,k=q+1}{\overset{n}{\sum}%
}\underset{\text{a,b,c}=1}{\overset{q}{\sum}}T_{\text{ab}i}(y_{0}%
)[(\perp_{\text{c}ik}T_{\text{ab}k})(\frac{\partial\phi}{\partial
\text{x}_{\text{c}}}$ $+$ $\Lambda_{\text{c}}\phi)](y_{0})$

$\qquad-\frac{1}{12}[$ $(R_{\text{a}i\text{b}j}+R_{\text{a}j\text{b}i})$
$-\underset{\text{c=1}}{\overset{\text{q}}{\sum}}(T_{\text{ac}i}T_{\text{bc}%
j}+T_{\text{ac}j}T_{\text{bc}i})$

$\qquad-\overset{n}{\underset{k=q+1}{\sum}}(\perp_{\text{a}ik}\perp
_{\text{b}jk}+$ $\perp_{\text{a}jk}\perp_{\text{b}ik})](y_{0})T_{\text{ab}%
i}(y_{0})\Lambda_{j}(y_{0})\phi(y_{0})$

$\qquad-\frac{1}{12}T_{\text{ab}i}^{2}(y_{0})\Omega_{ij}(y_{0})\phi(y_{0})$

\qquad\qquad\qquad\qquad\qquad\qquad\qquad\qquad\qquad\qquad\qquad\qquad
\qquad\qquad\qquad\qquad\qquad\qquad\qquad$\blacksquare$

Next we have:

\qquad\ \textbf{I}$_{32632}=-\frac{1}{6}\underset{i,j=q+1}{\overset{n}{\sum}%
}\underset{\text{a}=1}{\overset{q}{\sum}}[\frac{\partial\text{g}^{\text{a}j}%
}{\partial x_{i}}\underset{\text{c=1}}{\overset{\text{q}}{\sum}}\frac
{\partial\Gamma_{\text{a}j}^{\text{c}}}{\partial x_{i}}\frac{\partial\phi
}{\partial\text{x}_{\text{c}}}\circ\pi_{\text{P}}$ $+$
$\underset{k=1}{\overset{n}{\sum}}\frac{\partial}{\partial x_{i}}%
(\Gamma_{\text{a}j}^{k}\Lambda_{k}\phi\circ\pi_{\text{P}})](y_{0})$

\qquad$\qquad=-\frac{1}{6}\underset{i,j=q+1}{\overset{n}{\sum}}%
\underset{\text{a}=1}{\overset{q}{\sum}}[\frac{\partial\text{g}^{\text{a}j}%
}{\partial x_{i}}\underset{\text{c=1}}{\overset{\text{q}}{\sum}}\frac
{\partial\Gamma_{\text{a}j}^{\text{c}}}{\partial x_{i}}\frac{\partial\phi
}{\partial\text{x}_{\text{c}}}](y_{0})$

$\qquad\qquad-\frac{1}{6}\underset{i,j=q+1}{\overset{n}{\sum}}%
\underset{\text{a}=1}{\overset{q}{\sum}}\frac{\partial\text{g}^{\text{a}j}%
}{\partial x_{i}}[$ $\underset{k=1}{\overset{n}{\sum}}\frac{\partial
\Gamma_{\text{a}j}^{k}}{\partial x_{i}}\Lambda_{k}\phi+\Gamma_{\text{a}j}%
^{k}\frac{\partial\Lambda_{k}}{\partial x_{i}})](y_{0})\phi(y_{0})$

\qquad$\qquad=-\frac{1}{6}\underset{i,j=q+1}{\overset{n}{\sum}}%
\underset{\text{a,b=1}}{\overset{q}{\sum}}[\frac{\partial\text{g}^{\text{a}j}%
}{\partial x_{i}}\frac{\partial\Gamma_{\text{a}j}^{\text{b}}}{\partial x_{i}%
}\frac{\partial\phi}{\partial\text{x}_{\text{b}}}](y_{0})$

\qquad$\qquad\ -\frac{1}{6}\underset{i,j=q+1}{\overset{n}{\sum}}%
\underset{\text{a,b=1}}{\overset{q}{\sum}}\frac{\partial\text{g}^{\text{a}j}%
}{\partial x_{i}}[$ $\frac{\partial\Gamma_{\text{a}j}^{\text{b}}}{\partial
x_{i}}\Lambda_{\text{b}}+\Gamma_{\text{a}j}^{\text{b}}\frac{\partial
\Lambda_{\text{b}}}{\partial x_{i}}](y_{0})\phi(y_{0})$

$\qquad\qquad-\frac{1}{6}\underset{i,j,k=q+1}{\overset{n}{\sum}}%
\underset{\text{a}=1}{\overset{q}{\sum}}\frac{\partial\text{g}^{\text{a}j}%
}{\partial x_{i}}[$ $\frac{\partial\Gamma_{\text{a}j}^{k}}{\partial x_{i}%
}\Lambda_{k}+\Gamma_{\text{a}j}^{k}\frac{\partial\Lambda_{k}}{\partial x_{i}%
})](y_{0})\phi(y_{0})$

\qquad$\qquad=-\frac{1}{6}\underset{i,j=q+1}{\overset{n}{\sum}}%
\underset{\text{a,b=1}}{\overset{q}{\sum}}[\frac{\partial\text{g}^{\text{a}j}%
}{\partial x_{i}}$ $\frac{\partial\Gamma_{\text{a}j}^{\text{b}}}{\partial
x_{i}}(\frac{\partial\phi}{\partial\text{x}_{\text{b}}}+\Lambda_{\text{b}%
})](y_{0})\phi(y_{0})$

$\qquad\qquad-\frac{1}{6}\underset{i,j=q+1}{\overset{n}{\sum}}%
\underset{\text{a,b=1}}{\overset{q}{\sum}}[\frac{\partial\text{g}^{\text{a}j}%
}{\partial x_{i}}$ $\Gamma_{\text{a}j}^{\text{b}}\frac{\partial\Lambda
_{\text{b}}}{\partial x_{i}})](y_{0})\phi(y_{0})$

$\qquad\qquad-\frac{1}{6}\underset{i,j,k=q+1}{\overset{n}{\sum}}%
\underset{\text{a}=1}{\overset{q}{\sum}}[\frac{\partial\text{g}^{\text{a}j}%
}{\partial x_{i}}$ $\frac{\partial\Gamma_{\text{a}j}^{k}}{\partial x_{i}%
}\Lambda_{k}+\frac{\partial\text{g}^{\text{a}j}}{\partial x_{i}}%
\Gamma_{\text{a}j}^{k}\frac{\partial\Lambda_{k}}{\partial x_{i}})](y_{0}%
)\phi(y_{0})$

$\frac{\partial\text{g}^{\text{a}j}}{\partial x_{i}}(y_{0})=-$ $\perp
_{\text{a}ij}(y_{0})$ by (ii) of \textbf{Table A}$_{4};$ \ $\Gamma_{\text{a}%
j}^{\text{b}}(y_{0})=-\Gamma_{\text{ab}}^{j}(y_{0})=-$T$_{\text{ab}j}(y_{0})$
by (iii) of Table \textbf{Table A}$_{7}.$

$\Gamma_{\text{a}j}^{k}(y_{0})=$ $\perp_{\text{a}jk}(y_{0})$ by (x) of
\textbf{Table A}$_{7}.$

By (vi) of Table \textbf{Table A}$_{7},$

$\frac{\partial\Gamma_{\text{a}j}^{\text{b}}}{\partial x_{i}}(y_{0})$

$=\frac{1}{2}[-R_{\text{a}i\text{b}j}-R_{\text{a}j\text{b}i}%
+\underset{\text{c=1}}{\overset{\text{q}}{\sum}}T_{\text{ac}i}T_{\text{bc}%
j}-3\underset{\text{c=1}}{\overset{\text{q}}{\sum}}T_{\text{ac}j}%
T_{\text{bc}i}+\underset{\text{k=q+1}}{\overset{\text{n}}{\sum}}%
\perp_{\text{a}i\text{k}}\perp_{\text{b}j\text{k}}-\underset{\text{k=q+1}%
}{\overset{\text{n}}{\sum}}\perp_{\text{a}j\text{k}}\perp_{\text{b}i\text{k}%
}\ ](y_{0})$

$\frac{\partial\Gamma_{\text{a}j}^{k}}{\partial x_{i}}(y_{0}%
)=\ \overset{q}{\underset{\text{b}=1}{\sum(}}T_{\text{ab}j}\perp_{\text{b}%
ik})(y_{0})+\frac{2}{3}\left\{  2R_{\text{a}ijk}+R_{\text{a}jik}%
+R_{\text{a}kji}\right\}  (y_{0})$ by (xi) of \textbf{Table A}$_{7}.$

\ \textbf{I}$_{32632}=\frac{1}{12}\underset{i,j=q+1}{\overset{n}{\sum}%
}\underset{\text{a,b=1}}{\overset{q}{\sum}}[\perp_{\text{a}ij}(\frac
{\partial\phi}{\partial\text{x}_{\text{b}}}+\Lambda_{\text{b}})](y_{0})$

$\times\lbrack-R_{\text{a}i\text{b}j}-R_{\text{a}j\text{b}i}%
+\underset{\text{c=1}}{\overset{\text{q}}{\sum}}T_{\text{ac}i}T_{\text{bc}%
j}-3\underset{\text{c=1}}{\overset{\text{q}}{\sum}}T_{\text{ac}j}%
T_{\text{bc}i}+\underset{\text{k=q+1}}{\overset{\text{n}}{\sum}}%
\perp_{\text{a}i\text{k}}\perp_{\text{b}j\text{k}}-\underset{\text{k=q+1}%
}{\overset{\text{n}}{\sum}}\perp_{\text{a}j\text{k}}\perp_{\text{b}i\text{k}%
}\ ](y_{0})\phi(y_{0})$

\qquad$-\frac{1}{6}\underset{i,j=q+1}{\overset{n}{\sum}}\underset{\text{a,b=1}%
}{\overset{q}{\sum}}T_{\text{ab}j}(y_{0})\perp_{\text{a}ij}(y_{0}%
)\frac{\partial\Lambda_{\text{b}}}{\partial x_{i}}(y_{0})\phi(y_{0})$

\qquad$-\frac{1}{6}\underset{i,j,k=q+1}{\overset{n}{\sum}}\underset{\text{a}%
=1}{\overset{q}{\sum}}\perp_{\text{a}ij}(y_{0})[\overset{q}{\underset{\text{b}%
=1}{\sum(}}\perp_{\text{b}ik}T_{\text{ab}j})(y_{0})+\frac{2}{3}(2R_{\text{a}%
ijk}+R_{\text{a}jik}+R_{\text{a}kji})](y_{0})\Lambda_{k}(y_{0})\phi(y_{0})$

$\ \ \ \ +\frac{1}{6}\underset{i,j,k=q+1}{\overset{n}{\sum}}\underset{\text{a}%
=1}{\overset{q}{\sum}}\perp_{\text{a}ij}(y_{0})\perp_{\text{a}jk}(y_{0}%
)\Omega_{ik}(y_{0})\phi(y_{0})$

\qquad\qquad\qquad\qquad\qquad\qquad\qquad\qquad\qquad\qquad\qquad\qquad
\qquad\qquad\qquad\qquad\qquad\qquad$\blacksquare$

Therefore,

$\left(  D_{29}\right)  \qquad$ \textbf{I}$_{3263}=$ \textbf{I}$_{32631}+$
\textbf{I}$_{32632}$

$=\frac{1}{6}\underset{i,k=q+1}{\overset{n}{\sum}}\underset{\text{a,b,c}%
=1}{\overset{q}{\sum}}T_{\text{ab}i}(y_{0})[(\perp_{\text{c}ik}T_{\text{ab}%
k})(\frac{\partial\phi}{\partial\text{x}_{\text{c}}}$ $+$ $\Lambda_{\text{c}%
}\phi)](y_{0})\qquad$\textbf{I}$_{32631}$

$-\frac{1}{12}[$ $(R_{\text{a}i\text{b}j}+R_{\text{a}j\text{b}i})$
$-\underset{\text{c=1}}{\overset{\text{q}}{\sum}}(T_{\text{ac}i}T_{\text{bc}%
j}+T_{\text{ac}j}T_{\text{bc}i})$

$-\overset{n}{\underset{k=q+1}{\sum}}(\perp_{\text{a}ik}\perp_{\text{b}jk}+$
$\perp_{\text{a}jk}\perp_{\text{b}ik})](y_{0})T_{\text{ab}i}(y_{0})\Lambda
_{j}(y_{0})\phi(y_{0})-\frac{1}{12}T_{\text{ab}i}^{2}(y_{0})\Omega_{ij}%
(y_{0})\phi(y_{0}$

$+\frac{1}{12}\underset{i,j=q+1}{\overset{n}{\sum}}\underset{\text{a,b=1}%
}{\overset{q}{\sum}}[\perp_{\text{a}ij}(\frac{\partial\phi}{\partial
\text{x}_{\text{b}}}+\Lambda_{\text{b}})](y_{0})\qquad\qquad$\ \textbf{I}%
$_{32632}$

$\times\lbrack-R_{\text{a}i\text{b}j}-R_{\text{a}j\text{b}i}%
+\underset{\text{c=1}}{\overset{\text{q}}{\sum}}T_{\text{ac}i}T_{\text{bc}%
j}-3\underset{\text{c=1}}{\overset{\text{q}}{\sum}}T_{\text{ac}j}%
T_{\text{bc}i}+\underset{\text{k=q+1}}{\overset{\text{n}}{\sum}}%
\perp_{\text{a}i\text{k}}\perp_{\text{b}j\text{k}}-\underset{\text{k=q+1}%
}{\overset{\text{n}}{\sum}}\perp_{\text{a}j\text{k}}\perp_{\text{b}i\text{k}%
}\ ](y_{0})\phi(y_{0})$

$-\frac{1}{6}\underset{i,j=q+1}{\overset{n}{\sum}}\underset{\text{a,b=1}%
}{\overset{q}{\sum}}T_{\text{ab}j}(y_{0})\perp_{\text{a}ij}(y_{0}%
)\frac{\partial\Lambda_{\text{b}}}{\partial x_{i}}(y_{0})\phi(y_{0})$

$-\frac{1}{6}\underset{i,j,k=q+1}{\overset{n}{\sum}}\underset{\text{a}%
=1}{\overset{q}{\sum}}\perp_{\text{a}ij}(y_{0})[\overset{q}{\underset{\text{b}%
=1}{\sum(}}\perp_{\text{b}ik}T_{\text{ab}j})(y_{0})+\frac{2}{3}(2R_{\text{a}%
ijk}+R_{\text{a}jik}+R_{\text{a}kji})](y_{0})\Lambda_{k}(y_{0})\phi(y_{0})$

$\ +\frac{1}{6}\underset{i,j,k=q+1}{\overset{n}{\sum}}\underset{\text{a}%
=1}{\overset{q}{\sum}}\perp_{\text{a}ij}(y_{0})\perp_{\text{a}jk}(y_{0}%
)\Omega_{ik}(y_{0})\phi(y_{0})$

\qquad\qquad\qquad\qquad\qquad\qquad\qquad\qquad\qquad\qquad\qquad\qquad
\qquad\qquad$\blacksquare$

We finally conclude from $\left(  D_{25}\right)  ,$ $\left(  D_{28}\right)  $
and $\left(  D_{29}\right)  $ that,

$\left(  D_{30}\right)  $\qquad\ \textbf{I}$_{326}=$ \textbf{I}$_{3261}+$
\textbf{I}$_{3262}+$ \textbf{I}$_{3263}$

$=\mathbf{-}\frac{1}{12}\underset{\text{a,b=1}}{\overset{\text{q}}{\sum}%
}\underset{i,j=q+1}{\overset{n}{\sum}}T_{\text{ab}i}(y_{0})$

$\qquad\times\lbrack-R_{\text{a}i\text{b}i}+5\overset{q}{\underset{\text{c}%
=1}{\sum}}T_{\text{ac}i}T_{\text{bc}i}+2\overset{n}{\underset{k=q+1}{\sum}%
}\perp_{\text{a}ik}\perp_{\text{b}ik}](y_{0})\Lambda_{j}(y_{0})\phi
(y_{0})\ \ $I$_{32611}\qquad$\textbf{I}$_{3261}$

\qquad$+\frac{1}{12}\underset{i\text{=q+1}}{\overset{\text{n}}{\sum}%
}\underset{j\text{=q+1}}{\overset{\text{n}}{\sum}}\underset{\text{a=1}%
}{\overset{\text{q}}{\sum}}[4\underset{\text{c=1}}{\overset{q}{\sum}%
}(T_{\text{ac}i})(\perp_{j\text{c}i})+\frac{8}{3}R_{i\text{a}ij}](y_{0}%
)\qquad\qquad$I$_{32613}$

$\qquad\times\lbrack\underset{\text{c=1}}{\overset{\text{q}}{\sum}%
}T_{\text{ac}j}\frac{\partial\phi}{\partial\text{x}_{\text{c}}}%
+\underset{\text{b=1}}{\overset{\text{q}}{\sum}}T_{\text{ab}j}\Lambda
_{\text{b}}-\underset{k=q+1}{\overset{n}{\sum}}\perp_{\text{a}jk}\Lambda
_{k}](y_{0})\phi(y_{0})$

$\qquad-\frac{1}{24}\underset{i=q+1}{\overset{n}{\sum}}\underset{\text{a,b=1}%
}{\overset{\text{q}}{\sum}}$ $\overset{n}{\underset{k=q+1}{\sum}}%
$T$_{\text{aa}k}[\frac{8}{3}R_{i\text{c}ik}+4\underset{\text{d}%
=1}{\overset{q}{\sum}}(T_{\text{db}k})(\perp_{\text{d}ik})]\frac{\partial\phi
}{\partial\text{x}_{\text{b}}}$ $(y_{0})\qquad$ I$_{32621}\qquad$I$_{3262}$

$\qquad-\frac{1}{12}\underset{i=q+1}{\overset{n}{\sum}}\underset{\text{a,b=1}%
}{\overset{\text{q}}{\sum}}[$ $\overset{n}{\underset{k,l=q+1}{\sum}}%
\perp_{\text{b}ik}(-R_{\text{a}k\text{a}l}+\underset{\text{d=1}%
}{\overset{\text{q}}{\sum}}T_{\text{ad}k}T_{\text{ad}l}%
))-\overset{n}{\underset{k,l=q+1}{\sum}}\perp_{\text{b}ik}%
(\underset{r=q+1}{\overset{n}{\sum}}\perp_{\text{a}kr}\perp_{\text{a}%
lr})](y_{0})\frac{\partial\phi}{\partial\text{x}_{\text{b}}}$ $(y_{0})$

\qquad$-\frac{1}{12}\underset{i=q+1}{\overset{n}{\sum}}\underset{\text{b=1}%
}{\overset{\text{q}}{\sum}}[\frac{8}{3}\overset{q}{\underset{\text{c}=1}{\sum
}}($T$_{\text{bc}i}R_{ij\text{c}j})+\frac{2}{3}%
\overset{n}{\underset{k=q+1}{\sum}}(\perp_{\text{b}ik}R_{ijjk})](y_{0}%
)\frac{\partial\phi}{\partial\text{x}_{\text{b}}}(y_{0})$

$-\frac{1}{6}\underset{i=q+1}{\overset{n}{\sum}}\underset{\text{a,b=1}%
}{\overset{\text{q}}{\sum}}[4\overset{q}{\underset{\text{c=1}}{%
{\textstyle\sum}
}}$R$_{ij\text{c}i}^{{}}T_{\text{bc}j}+$ $4\overset{n}{\underset{k=q+1}{%
{\textstyle\sum}
}}$R$_{ijik}\perp_{\text{b}jk}+3\nabla_{i}$R$_{j\text{b}ij}%
+4\overset{q}{\underset{\text{c=1}}{%
{\textstyle\sum}
}}$R$_{ij\text{c}j}^{{}}T_{\text{bc}i}+$ $4$R$_{ijjk}\perp_{\text{b}ik}%
](y_{0})\frac{\partial\phi}{\partial\text{x}_{\text{b}}}(y_{0})$

$-\frac{1}{144}[\{4\nabla_{i}$R$_{i\text{a}j\text{a}}$ $+2\nabla_{j}%
$R$_{i\text{a}i\text{a}}+$ $8$ $(\overset{q}{\underset{\text{c=1}}{%
{\textstyle\sum}
}}R_{\text{a}i\text{c}i}^{{}}T_{\text{ac}j}+\;\overset{n}{\underset{k=q+1}{%
{\textstyle\sum}
}}R_{\text{a}iik}\perp_{\text{a}jk})\qquad$I$_{326221}\qquad$I$_{32622}$

\ \ $+8(\overset{q}{\underset{\text{c=1}}{%
{\textstyle\sum}
}}R_{\text{a}i\text{c}j}^{{}}T_{\text{ac}i}+\overset{n}{\underset{l=q+1}{%
{\textstyle\sum}
}}R_{\text{a}ijl}\perp_{\text{a}il})+8(\overset{q}{\underset{\text{c=1}}{%
{\textstyle\sum}
}}R_{\text{a}j\text{c}i}^{{}}T_{\text{ac}i}+\overset{q}{\underset{\text{c=1}}{%
{\textstyle\sum}
}}R_{\text{a}j\text{c}i}^{{}}T_{\text{ac}i})\}$\ 

$+\frac{2}{3}\underset{k=q+1}{\overset{n}{\sum}}\{T_{\text{aa}k}%
(R_{ijik}+3\overset{q}{\underset{\text{c}=1}{\sum}}\perp_{\text{c}ij}%
\perp_{\text{c}ik})\}](y_{0})\Lambda_{k}(y_{0})\phi(y_{0})$

\ $+\frac{1}{24}[\frac{4}{3}\overset{q}{\underset{\text{a}=1}{\sum}}%
\perp_{\text{a}ik}R_{ij\text{a}j}+\frac{1}{3}(\nabla_{i}R_{kjij}+\nabla
_{j}R_{ijik}+\nabla_{k}R_{ijij})](y_{0})\Lambda_{k}(y_{0})\phi(y_{0})$

\qquad\qquad\ $\mathbf{-}$ $\frac{1}{72}\underset{i,j=q+1}{\overset{n}{\sum}%
}\underset{\text{a}=1}{\overset{q}{\sum}}T_{\text{aa}j}(y_{0})\frac
{\partial\Omega_{ij}}{\partial\text{x}_{i}}(y_{0})\phi(y_{0})\qquad
$I$_{326222}$

$+$ $\frac{1}{12}$ $\overset{n}{\underset{j=q+1}{\sum}}(\perp_{\text{b}ij}%
$T$_{\text{aa}j})(y_{0})\left(  \Omega_{i\text{b}}(y_{0})+[\Lambda_{\text{b}%
},\Lambda_{i}]\right)  (y_{0})\phi(y_{0})$ \qquad I$_{326223}$

$\mathbf{-}\frac{1}{18}$ $\underset{i,j=q+1}{\overset{n}{\sum}}%
\underset{\text{b=1}}{\overset{q}{\sum}}R_{\text{b}jij}(y_{0})\left(
\Omega_{i\text{a}}(y_{0})+[\Lambda_{\text{a}},\Lambda_{i}]\right)  (y_{0}%
)\phi(y_{0})$

$\mathbf{-}$ $\frac{1}{24}\underset{i,j=q+1}{\overset{n}{\sum}}%
\underset{\text{a=1}}{\overset{q}{\sum}}[$ R$_{\text{a}i\text{a}j}$
$-\underset{\text{c=1}}{\overset{\text{q}}{\sum}}T_{\text{ac}i}T_{\text{ac}%
j}-\overset{n}{\underset{k=q+1}{\sum}}(\perp_{\text{a}ik}\perp_{\text{a}%
jk}](y_{0})\Omega_{ij}(y_{0})\phi(y_{0})$

\ $\mathbf{-}$ $\frac{1}{36}\underset{i,j,k=q+1}{\overset{n}{\sum}}%
R_{ijkj}(y_{0})\Omega_{ik}(y_{0})(y_{0})\phi(y_{0})$

$+\frac{1}{6}\underset{i,k=q+1}{\overset{n}{\sum}}\underset{\text{a,b,c}%
=1}{\overset{q}{\sum}}T_{\text{ab}i}(y_{0})[(\perp_{\text{c}ik}T_{\text{ab}%
k})(\frac{\partial\phi}{\partial\text{x}_{\text{c}}}$ $+$ $\Lambda_{\text{c}%
}\phi)](y_{0})\qquad$\textbf{I}$_{32631}\qquad$\textbf{I}$_{3263}$

$-\frac{1}{12}[$ $(R_{\text{a}i\text{b}j}+R_{\text{a}j\text{b}i})$
$-\underset{\text{c=1}}{\overset{\text{q}}{\sum}}(T_{\text{ac}i}T_{\text{bc}%
j}+T_{\text{ac}j}T_{\text{bc}i})$

$-\overset{n}{\underset{k=q+1}{\sum}}(\perp_{\text{a}ik}\perp_{\text{b}jk}+$
$\perp_{\text{a}jk}\perp_{\text{b}ik})](y_{0})T_{\text{ab}i}(y_{0})\Lambda
_{j}(y_{0})\phi(y_{0})-\frac{1}{12}T_{\text{ab}i}^{2}(y_{0})\Omega_{ij}%
(y_{0})\phi(y_{0}$

$+\frac{1}{12}\underset{i,j=q+1}{\overset{n}{\sum}}\underset{\text{a,b=1}%
}{\overset{q}{\sum}}[\perp_{\text{a}ij}(\frac{\partial\phi}{\partial
\text{x}_{\text{b}}}+\Lambda_{\text{b}})](y_{0})\qquad\qquad$\ \textbf{I}%
$_{32632}$

$\times\lbrack-R_{\text{a}i\text{b}j}-R_{\text{a}j\text{b}i}%
+\underset{\text{c=1}}{\overset{\text{q}}{\sum}}T_{\text{ac}i}T_{\text{bc}%
j}-3\underset{\text{c=1}}{\overset{\text{q}}{\sum}}T_{\text{ac}j}%
T_{\text{bc}i}+\underset{\text{k=q+1}}{\overset{\text{n}}{\sum}}%
\perp_{\text{a}i\text{k}}\perp_{\text{b}j\text{k}}-\underset{\text{k=q+1}%
}{\overset{\text{n}}{\sum}}\perp_{\text{a}j\text{k}}\perp_{\text{b}i\text{k}%
}\ ](y_{0})\phi(y_{0})$

$-\frac{1}{6}\underset{i,j=q+1}{\overset{n}{\sum}}\underset{\text{a,b=1}%
}{\overset{q}{\sum}}T_{\text{ab}j}(y_{0})\perp_{\text{a}ij}(y_{0}%
)\frac{\partial\Lambda_{\text{b}}}{\partial x_{i}}(y_{0})\phi(y_{0})$

$-\frac{1}{6}\underset{i,j,k=q+1}{\overset{n}{\sum}}\underset{\text{a}%
=1}{\overset{q}{\sum}}\perp_{\text{a}ij}(y_{0})[\overset{q}{\underset{\text{b}%
=1}{\sum(}}\perp_{\text{b}ik}T_{\text{ab}j})(y_{0})+\frac{2}{3}(2R_{\text{a}%
ijk}+R_{\text{a}jik}+R_{\text{a}kji})](y_{0})\Lambda_{k}(y_{0})\phi(y_{0})$

$\ \ \ +\frac{1}{6}\underset{i,j,k=q+1}{\overset{n}{\sum}}\underset{\text{a}%
=1}{\overset{q}{\sum}}\perp_{\text{a}ij}(y_{0})\perp_{\text{a}jk}(y_{0}%
)\Omega_{ik}(y_{0})\phi(y_{0})$

\qquad\qquad\qquad\qquad\qquad\qquad\qquad\qquad\qquad\qquad\qquad\qquad
\qquad\qquad\qquad\qquad\qquad\qquad$\blacksquare$

We next compute I$_{327}:$

We recall that $\frac{\partial}{\partial x_{i}}\pi_{P}(y_{0})=0=\frac
{\partial^{2}}{\partial x_{i}^{2}}\pi_{P}(y_{0}),$ and so we have:

(vii)\qquad\textbf{I}$_{327}=\frac{1}{24}\underset{i=q+1}{\overset{n}{\sum}%
}\frac{\partial^{2}}{\partial x_{i}^{2}}[$W$\phi\circ\pi_{\text{P}}](y_{0})$

$\left(  D_{31}\right)  \qquad$\ \textbf{I}$_{327}=\ \frac{1}{24}%
\underset{i\text{=q+1}}{\overset{\text{n}}{\sum}}\frac{\partial^{2}\text{W}%
}{\partial x_{i}^{2}}(y_{0})\phi(y_{0})$

\qquad\qquad\qquad\qquad\qquad\qquad\qquad\qquad\qquad\qquad\qquad
$\blacksquare$

(viii)\qquad\textbf{I}$_{328}=\frac{1}{12}\underset{i=q+1}{\overset{n}{\sum}}$
$\underset{\text{a=1}}{\overset{\text{q}}{\sum}}\frac{\partial^{2}}{\partial
x_{i}^{2}}[(\nabla\log\Psi)_{\text{a}}\frac{\partial\phi}{\partial
\text{x}_{\text{a}}}\circ\pi_{\text{P}}](y_{0})$

\qquad\qquad$=\frac{1}{12}$ $\underset{\text{a=1}}{\overset{\text{q}}{\sum}%
}\frac{\partial^{2}}{\partial x_{i}^{2}}(\nabla\log\theta^{-\frac{1}{2}%
})_{\text{a}}(y_{0})\frac{\partial\phi}{\partial\text{x}_{\text{a}}}%
(y_{0})+\frac{1}{12}$ $\underset{\text{a=1}}{\overset{\text{q}}{\sum}}%
\frac{\partial^{2}}{\partial x_{i}^{2}}(\nabla\log\Phi)_{\text{a}}(y_{0}%
)\frac{\partial\phi}{\partial\text{x}_{\text{a}}}(y_{0})$\qquad\qquad

By (xvii)$^{\ast}$ of \textbf{Table A}$_{9},$ we have:

$\qquad\frac{1}{12}\underset{i=q+1}{\overset{n}{\sum}}$ $\underset{\text{a=1}%
}{\overset{\text{q}}{\sum}}[\frac{\partial^{2}}{\partial x_{i}^{2}}(\nabla
\log\theta^{-\frac{1}{2}})_{\text{a}}](y_{0})$

$\qquad\qquad=\frac{1}{24}\underset{i,j=q+1}{\overset{n}{\sum}}%
\underset{\text{a=1}}{\overset{\text{q}}{\sum}}<H,j>[4\underset{\text{c}%
=1}{\overset{q}{\sum}}(T_{\text{ac}i})(\perp_{j\text{c}i})+\frac{8}%
{3}R_{i\text{a}ij}](y_{0})$

$\qquad\ \ \ \ \ \ \ \ -\frac{1}{24}\underset{i,j=q+1}{\overset{n}{\sum}%
}\underset{\text{a=1}}{\overset{\text{q}}{\sum}}\perp_{\text{a}ij}%
(y_{0})[<H,i><H,j>](y_{0})$\qquad$\ \ \ \ \qquad\qquad-\frac{1}{72}%
\underset{i,j=q+1}{\overset{n}{\sum}}\underset{\text{a=1}}{\overset{\text{q}%
}{\sum}}\perp_{\text{a}ij}(y_{0})[2\varrho_{ij}%
+4\overset{q}{\underset{\text{a}=1}{\sum}}R_{i\text{a}j\text{a}}%
-6\overset{q}{\underset{\text{b,c}=1}{\sum}}T_{\text{cc}i}T_{\text{bb}%
j}-T_{\text{bc}i}T_{\text{bc}j}](y_{0})$

By (xvi) of \textbf{Table B}$_{1},$

$\frac{\partial^{2}}{\partial x_{i}^{2}}(\nabla$log$\Phi_{P})_{\text{a}%
}(y)=-4\underset{\text{b=1}}{\overset{\text{q}}{%
{\textstyle\sum}
}}T_{\text{ab}i}(y)\frac{\partial X_{i}}{\partial x_{\text{b}}}%
(y)+\underset{k=q+1}{\overset{n}{%
{\textstyle\sum}
}}\perp_{\text{a}ik}(y)\left[  \left(  \frac{\partial X_{i}}{\partial x_{k}%
}+\frac{\partial X_{k}}{\partial x_{i}}\right)  \right]  (y)$

$\qquad\qquad\qquad\qquad+\frac{8}{3}\underset{k=q+1}{\overset{n}{%
{\textstyle\sum}
}}R_{i\text{a}ik}(y)X_{k}(y)+[2X_{i}\frac{\partial X_{i}}{\partial
x_{\text{a}}}-\frac{\partial^{2}X_{i}}{\partial x_{\text{a}}\partial x_{i}%
}](y)$

\qquad\qquad\qquad\qquad$=-4\underset{\text{b=1}}{\overset{\text{q}}{%
{\textstyle\sum}
}}T_{\text{ab}i}(y)\frac{\partial X_{i}}{\partial x_{\text{b}}}%
(y)+\underset{j=q+1}{\overset{n}{%
{\textstyle\sum}
}}\perp_{\text{a}ij}(y)\left(  \frac{\partial X_{i}}{\partial x_{j}}%
+\frac{\partial X_{j}}{\partial x_{i}}\right)  (y)$

$\qquad\qquad\qquad\qquad+\frac{8}{3}\underset{j=q+1}{\overset{n}{%
{\textstyle\sum}
}}R_{i\text{a}ij}(y)X_{j}(y)+\left(  2X_{i}\frac{\partial X_{i}}{\partial
x_{\text{a}}}-\frac{\partial^{2}X_{i}}{\partial x_{\text{a}}\partial x_{i}%
}\right)  (y)$

Therefore,

$\left(  D_{31}\right)  \qquad$\ \textbf{I}$_{327}=\ \frac{1}{24}%
\underset{i\text{=q+1}}{\overset{\text{n}}{\sum}}\frac{\partial^{2}\text{W}%
}{\partial x_{i}^{2}}(y_{0})\phi(y_{0})$

$=\frac{1}{12}$ $\underset{\text{a=1}}{\overset{\text{q}}{\sum}}\frac
{\partial^{2}}{\partial x_{i}^{2}}(\nabla\log\theta^{-\frac{1}{2}})_{\text{a}%
}(y_{0})\frac{\partial\phi}{\partial\text{x}_{\text{a}}}(y_{0})+\frac{1}{12}$
$\underset{\text{a=1}}{\overset{\text{q}}{\sum}}\frac{\partial^{2}}{\partial
x_{i}^{2}}(\nabla\log\Phi)_{\text{a}}(y_{0})\frac{\partial\phi}{\partial
\text{x}_{\text{a}}}(y_{0})$

$\left(  D_{32}\right)  $\qquad I$_{328}\qquad=\frac{1}{24}%
\underset{i,j\text{=q+1}}{\overset{\text{n}}{\sum}}\underset{\text{a=1}%
}{\overset{\text{q}}{\sum}}<H,j>[4\underset{\text{c}=1}{\overset{q}{\sum}%
}(T_{\text{ac}i})(\perp_{j\text{c}i})+\frac{8}{3}R_{i\text{a}ij}](y_{0}%
)\frac{\partial\phi}{\partial\text{x}_{\text{a}}}(y_{0})$

$\qquad\ \ \ \ \ \ +\frac{1}{24}\underset{i,j\text{=q+1}}{\overset{\text{n}%
}{\sum}}\underset{\text{a=1}}{\overset{\text{q}}{\sum}}\perp_{\text{a}%
ji}(y_{0})[<H,i><H,j>](y_{0})\frac{\partial\phi}{\partial\text{x}_{\text{a}}%
}(y_{0})$\qquad$\ \ \ \ \qquad\ \ \ \ \ +\frac{1}{72}%
\underset{i,j=q+1}{\overset{n}{\sum}}\underset{\text{a=1}}{\overset{\text{q}%
}{\sum}}\perp_{\text{a}ji}(y_{0})[2\varrho_{ij}%
+4\overset{q}{\underset{\text{a}=1}{\sum}}R_{i\text{a}j\text{a}}%
-6\overset{q}{\underset{\text{b,c}=1}{\sum}}T_{\text{cc}i}T_{\text{bb}%
j}-T_{\text{bc}i}T_{\text{bc}j}](y_{0})\frac{\partial\phi}{\partial
\text{x}_{\text{a}}}(y_{0})$

\qquad$\qquad+\frac{8}{3}R_{j\text{a}ji}(y_{0})X_{i}(y_{0})+[2X_{j}%
\frac{\partial X_{j}}{\partial x_{\text{a}}}-\frac{\partial^{2}X_{j}}{\partial
x_{\text{a}}\partial x_{j}}](y_{0})\frac{\partial\phi}{\partial\text{x}%
_{\text{a}}}(y_{0})$\qquad

\qquad\qquad\qquad\qquad\qquad\qquad\qquad\qquad\qquad\qquad\qquad\qquad
\qquad\qquad\qquad$\blacksquare$

(ix)\qquad We next compute I$_{329}:$

\qquad\textbf{I}$_{329}=\frac{1}{12}\underset{i=q+1}{\overset{n}{\sum}}$
$\underset{j=1}{\overset{n}{\sum}}$ $\frac{\partial^{2}}{\partial x_{i}^{2}%
}[(\nabla\log\Psi)_{j}\Lambda_{j}\phi\circ\pi_{\text{P}}](y_{0})$

$=\frac{1}{12}\underset{i=q+1}{\overset{n}{\sum}}$
$\underset{j=1}{\overset{n}{\sum}}$ $\frac{\partial^{2}}{\partial x_{i}^{2}%
}[(\nabla\log\theta^{-\frac{1}{2}})_{j}\Lambda_{j}\phi\circ\pi_{\text{P}%
}](y_{0})+\frac{1}{12}$ $\underset{j=1}{\overset{n}{\sum}}$ $\frac
{\partial^{2}}{\partial x_{i}^{2}}[(\nabla\log\Phi)_{j}\Lambda_{j}\phi\circ
\pi_{\text{P}}](y_{0})$

\qquad\qquad$=$ I$_{3291}+$ I$_{3292}$

where,

\qquad I$_{3291}=\frac{1}{12}$ $\underset{i=q+1}{\overset{n}{\sum}%
}\underset{j=1}{\overset{n}{\sum}}$ $\frac{\partial^{2}}{\partial x_{i}^{2}%
}[(\nabla\log\theta^{-\frac{1}{2}})_{j}\Lambda_{j}\phi\circ\pi_{\text{P}%
}](y_{0})$

\qquad I$_{3292}=\frac{1}{12}\underset{i=q+1}{\overset{n}{\sum}}$
$\underset{j=1}{\overset{n}{\sum}}$ $\frac{\partial^{2}}{\partial x_{i}^{2}%
}[(\nabla\log\Phi)_{j}\Lambda_{j}\phi\circ\pi_{\text{P}}](y_{0})$

We have:

\qquad I$_{3291}=\frac{1}{12}\underset{i=q+1}{\overset{n}{\sum}}$
$\underset{j=1}{\overset{n}{\sum}}$ $\frac{\partial^{2}}{\partial x_{i}^{2}%
}[(\nabla\log\theta^{-\frac{1}{2}})_{j}\Lambda_{j}\phi\circ\pi_{\text{P}%
}](y_{0})$

\qquad$\qquad=\frac{1}{12}$ $\underset{i=q+1}{\overset{n}{\sum}}%
\underset{j=\text{1}}{\overset{n}{\sum}}$ $\frac{\partial^{2}}{\partial
x_{i}^{2}}[(\nabla\log\theta^{-\frac{1}{2}})_{j}](y_{0})\Lambda_{j}(y_{0}%
)\phi\circ\pi_{P}(y_{0})$

\qquad\ \qquad$\ +$ $\frac{1}{12}\underset{i=q+1}{\overset{n}{\sum}}$
$\underset{j=1}{\overset{n}{\sum}}$ $(\nabla\log\theta^{-\frac{1}{2}}%
)_{j}(y_{0})\frac{\partial^{2}\Lambda_{j}}{\partial x_{i}^{2}}(y_{0})\phi
\circ\pi_{P}(y_{0})$

$\qquad\qquad+\frac{1}{6}\underset{i=q+1}{\overset{n}{\sum}}$
$\underset{j=1}{\overset{n}{\sum}}$ $\frac{\partial}{\partial x_{i}}%
(\nabla\log\theta^{-\frac{1}{2}})_{j}(y_{0})\frac{\partial\Lambda_{j}%
}{\partial x_{i}}(y_{0})\phi\circ\pi_{P}(y_{0})$

\qquad\qquad\qquad$\ =$ I$_{32911}+$ I$_{32912}$ $+$ I$_{32913}$ where,

\qquad\qquad I$_{32911}$ $=\frac{1}{12}$ $\underset{i=q+1}{\overset{n}{\sum}%
}\underset{j=1}{\overset{n}{\sum}}$ $\frac{\partial^{2}}{\partial x_{i}^{2}%
}(\nabla\log\theta^{-\frac{1}{2}})_{j}(y_{0})\Lambda_{j}(y_{0})\phi\circ
\pi_{P}(y_{0})$

\qquad\qquad I$_{32912}=\frac{1}{12}\underset{i=q+1}{\overset{n}{\sum}}$
$\underset{j=1}{\overset{n}{\sum}}$ $(\nabla\log\theta^{-\frac{1}{2}}%
)_{j}(y_{0})\frac{\partial^{2}\Lambda_{j}}{\partial x_{i}^{2}}(y_{0})\phi
\circ\pi_{P}(y_{0})$

\qquad\qquad I$_{32913}=\frac{1}{6}\underset{i=q+1}{\overset{n}{\sum}}$
$\underset{j=1}{\overset{n}{\sum}}$ $\frac{\partial}{\partial x_{i}}%
(\nabla\log\theta^{-\frac{1}{2}})_{j}(y_{0})\frac{\partial\Lambda_{j}%
}{\partial x_{i}}(y_{0})\phi\circ\pi(y_{0})$

We now compute each of the above expressions:

I$_{32911}$ $=\frac{1}{12}\underset{i=q+1}{\overset{n}{\sum}}$
$\underset{j=1}{\overset{n}{\sum}}$ $\frac{\partial^{2}}{\partial x_{i}^{2}%
}(\nabla\log\theta^{-\frac{1}{2}})_{j}(y_{0})\Lambda_{j}(y_{0})\phi\circ
\pi_{P}(y_{0})$

\qquad$=\frac{1}{12}\underset{i=q+1}{\overset{n}{\sum}}$ $\underset{\text{a=1}%
}{\overset{\text{q}}{\sum}}$ $\frac{\partial^{2}}{\partial x_{i}^{2}}%
(\nabla\log\theta^{-\frac{1}{2}})_{\text{a}}(y_{0})\Lambda_{\text{a}}%
(y_{0})\phi\circ\pi_{P}(y_{0})$

\qquad$+\frac{1}{12}$ $\underset{i,j=q+1}{\overset{n}{\sum}}$ $\frac
{\partial^{2}}{\partial x_{i}^{2}}(\nabla\log\theta^{-\frac{1}{2}})_{j}%
(y_{0})\Lambda_{j}(y_{0})\phi\circ\pi_{P}(y_{0})$

By (xvii)$^{\ast}$ and (xvii)$^{\ast\ast}$ of \textbf{Table A}$_{9}$ we have
for a = 1,...,q and $i,j=q+1,...,n:$

$\left(  D_{32}\right)  \qquad$I$_{32911}=\frac{1}{24}%
\underset{i,j=q+1}{\overset{n}{\sum}}\underset{\text{a=1}}{\overset{\text{q}%
}{\sum}}<H,j>[\frac{8}{3}R_{i\text{a}ij}-4\underset{\text{b=1}%
}{\overset{\text{q}}{\sum}}T_{\text{ab}i}\perp_{\text{b}ij}](y_{0}%
)\Lambda_{\text{a}}(y_{0})\phi(y_{0})$

$\qquad\ \ \ \ \ \ \ +\frac{1}{12}\underset{i,j=q+1}{\overset{n}{\sum}%
}\underset{\text{a=1}}{\overset{\text{q}}{\sum}}\perp_{\text{a}ij}%
(y_{0})[<H,i><H,j>](y_{0})\Lambda_{\text{a}}(y_{0})\phi(y_{0})$

\qquad$\ \ \ \ \ \ +\frac{1}{72}\underset{i,j=q+1}{\overset{n}{\sum}%
}\underset{\text{a=1}}{\overset{\text{q}}{\sum}}\perp_{\text{a}ji}%
(y_{0})[2\varrho_{ij}+4\overset{q}{\underset{\text{a}=1}{\sum}}R_{i\text{a}%
j\text{a}}-6\overset{q}{\underset{\text{b,c}=1}{\sum}}T_{\text{cc}%
i}T_{\text{bb}j}-T_{\text{bc}i}T_{\text{bc}j}](y_{0})\Lambda_{\text{a}}%
(y_{0})\phi(y_{0})$

$\qquad\qquad+\frac{1}{36}\underset{i,j,k=q+1}{\overset{n}{\sum}}%
<H,k>(y_{0})R_{ijik}(y_{0})\Lambda_{j}(y_{0})\phi(y_{0})$

$\qquad\ \ \ \ \ \ -\frac{1}{288}\underset{i,j=q+1}{\overset{n}{\sum}%
}<H,j>(y_{0})[3<H,i>^{2}+2(\tau^{M}-3\tau^{P}+\overset{q}{\underset{\text{a}%
=1}{\sum}}\varrho_{\text{aa}}+\overset{q}{\underset{\text{a,b}=1}{\sum}%
}R_{\text{abab}})](y_{0})\Lambda_{j}(y_{0})\phi(y_{0})$

$\qquad\qquad-\frac{1}{12}\underset{i,j=q+1}{\overset{n}{\sum}}<H,i>(y_{0})$

$\qquad\qquad\times\lbrack\frac{3}{4}<H,i><H,j>$\ $+\frac{1}{6}(\varrho
_{ij}+2\overset{q}{\underset{\text{a}=1}{\sum}}R_{i\text{a}j\text{a}%
}-3\overset{q}{\underset{\text{a,b=1}}{\sum}}T_{\text{aa}i}T_{\text{bb}%
j}-T_{\text{ab}i}T_{\text{ab}j})](y_{0})\Lambda_{j}(y_{0})\phi(y_{0})$

$\qquad\qquad+\frac{5}{32}\underset{i,j=q+1}{\overset{n}{\sum}}<H,i>^{2}%
<H,j>\Lambda_{j}(y_{0})\phi(y_{0})$

$\qquad+\frac{1}{48}\underset{i,j=q+1}{\overset{n}{\sum}}$%
$<$%
H,$i$%
$>$%
(y$_{0}$)$[$(2$\varrho_{ij}$+4$\overset{q}{\underset{\text{a}=1}{\sum}}%
$R$_{i\text{a}j\text{a}}$-3$\overset{q}{\underset{\text{a,b=1}}{\sum}}%
$T$_{\text{aa}i}$T$_{\text{bb}j}$-T$_{\text{ab}i}$T$_{\text{ab}j}%
$-3$\overset{q}{\underset{\text{a,b=1}}{\sum}}$T$_{\text{aa}j}$T$_{\text{bb}%
i}$-T$_{\text{ab}j}$T$_{\text{ab}i}$)$](y_{0})\Lambda_{j}(y_{0})\phi(y_{0})$

$\qquad+\frac{1}{48}\underset{i,j=q+1}{\overset{n}{\sum}}<H,j>[$ $\tau
^{M}-3\tau^{P}+\overset{q}{\underset{\text{a}=1}{\sum}}\varrho_{\text{aa}%
}+\overset{q}{\underset{\text{a,b}=1}{\sum}}R_{\text{abab}}](y_{0})\Lambda
_{j}(y_{0})\phi(y_{0})$

$\qquad+\frac{1}{144}\underset{i,j=q+1}{\overset{n}{\sum}}[\nabla_{i}%
\varrho_{ij}-2\varrho_{ij}<H,i>+\overset{q}{\underset{\text{a}=1}{\sum}%
}(\nabla_{i}R_{\text{a}i\text{a}j}-4R_{i\text{a}j\text{a}}<H,i>)$

$\qquad+4\overset{q}{\underset{\text{a,b=1}}{\sum}}R_{i\text{a}j\text{b}%
}T_{\text{ab}i}$+2$\overset{q}{\underset{\text{a,b,c=1}}{\sum}}(T_{\text{aa}%
i}T_{\text{bb}j}T_{\text{cc}i}-3T_{\text{aa}i}T_{\text{bc}j}T_{\text{bc}%
i}+2T_{\text{ab}i}T_{\text{bc}j}T_{\text{ca}i})](y_{0})\Lambda_{j}(y_{0}%
)\phi(y_{0})$\qquad\qquad\qquad\qquad\qquad\ \ 

$\qquad+\frac{1}{144}\underset{i,j=q+1}{\overset{n}{\sum}}[\nabla_{j}%
\varrho_{ii}-2\varrho_{ji}<H,i>+\overset{q}{\underset{\text{a}=1}{\sum}%
}(\nabla_{j}R_{\text{a}i\text{a}i}-4R_{j\text{a}i\text{a}}<H,i>)$

$\qquad+4\overset{q}{\underset{\text{a,b=1}}{\sum}}R_{j\text{a}i\text{b}%
}T_{\text{ab}i}$+2$\overset{q}{\underset{\text{a,b,c=1}}{\sum}}(T_{\text{aa}%
j}T_{\text{bb}i}T_{\text{cc}i}-3T_{\text{aa}j}T_{\text{bc}i}T_{\text{bc}%
i}+2T_{\text{ab}j}T_{\text{bc}i}T_{\text{ca}i})](y_{0})\Lambda_{j}(y_{0}%
)\phi(y_{0})$

$\qquad+\frac{1}{144}\underset{i,j=q+1}{\overset{n}{\sum}}[\nabla_{i}%
\varrho_{ij}-2\varrho_{ii}<H,j>+\overset{q}{\underset{\text{a}=1}{\sum}%
}(\nabla_{i}R_{\text{a}i\text{a}j}-4R_{i\text{a}i\text{a}}<H,j>)$

$\qquad+4\overset{q}{\underset{\text{a,b=1}}{\sum}}R_{i\text{a}i\text{b}%
}T_{\text{ab}j}$+2$\overset{q}{\underset{\text{a,b,c=1}}{\sum}}(T_{\text{aa}%
i}T_{\text{bb}i}T_{\text{cc}j}-3T_{\text{aa}i}T_{\text{bc}i}T_{\text{bc}%
j}+2T_{\text{ab}i}T_{\text{bc}i}T_{\text{ac}j})](y_{0})\Lambda_{j}(y_{0}%
)\phi(y_{0}).$

\qquad\qquad\qquad\qquad\qquad\qquad\qquad\qquad\qquad\qquad\qquad\qquad
\qquad\qquad\qquad\qquad\qquad\qquad\qquad$\qquad\blacksquare$

We next compute:

I$_{32912}=\frac{1}{12}\underset{i=q+1}{\overset{n}{\sum}}$
$\underset{j=1}{\overset{n}{\sum}}$ $(\nabla\log\theta^{-\frac{1}{2}}%
)_{j}(y_{0})\frac{\partial^{2}\Lambda_{j}}{\partial x_{i}^{2}}(y_{0})\phi
\circ\pi_{P}(y_{0})$

$=\frac{1}{12}\underset{i=q+1}{\overset{n}{\sum}}$ $\underset{\text{a=1}%
}{\overset{\text{q}}{\sum}}$ $(\nabla\log\theta^{-\frac{1}{2}})_{\text{a}%
}(y_{0})\frac{\partial^{2}\Lambda_{\text{a}}}{\partial x_{i}^{2}}(y_{0}%
)\phi\circ\pi_{P}(y_{0})$

$+\frac{1}{12}\underset{i=q+1}{\overset{n}{\sum}}$
$\underset{q+1=1}{\overset{n}{\sum}}$ $(\nabla\log\theta^{-\frac{1}{2}}%
)_{j}(y_{0})\frac{\partial^{2}\Lambda_{j}}{\partial x_{i}^{2}}(y_{0})\phi
\circ\pi_{P}(y_{0})$

$(\nabla\log\theta^{-\frac{1}{2}})_{\text{a}}(y_{0})=0$ by (iii)$^{\ast}$ of
\textbf{Table A}$_{9}$ and

$(\nabla\log\theta^{-\frac{1}{2}})_{j}(y_{0})=\frac{1}{2}<H,j>(y_{0})$ by
(iv)$^{\ast}$ of \textbf{Table A}$_{9}.$

$\frac{\partial^{2}\Lambda_{j}}{\partial x_{i}^{2}}(y_{0})=\frac{1}{3}%
\frac{\partial\Omega_{ij}}{\partial\text{x}_{i}}(y_{0})$ by (xi) of
\textbf{Proposition 5}

$\left(  D_{33}\right)  \qquad$I$_{32912}=\frac{1}{72}%
\underset{i,j=q+1}{\overset{n}{\sum}}$ $<H,j>(y_{0})\frac{\partial\Omega_{ij}%
}{\partial\text{x}_{i}}(y_{0})\phi(y_{0})$

\qquad\qquad\qquad\qquad\qquad\qquad\qquad\qquad\qquad\qquad\qquad
\qquad$\blacksquare$

We then compute the last expression here:

I$_{32913}=\frac{1}{6}\underset{i=q+1}{\overset{n}{\sum}}$
$\underset{j=1}{\overset{n}{\sum}}$ $\frac{\partial}{\partial x_{i}}%
(\nabla\log\theta^{-\frac{1}{2}})_{j}(y_{0})\frac{\partial\Lambda_{j}%
}{\partial x_{i}}(y_{0})\phi\circ\pi_{P}(y_{0})$

\qquad$=\frac{1}{6}\underset{i=q+1}{\overset{n}{\sum}}$ $\underset{\text{a=1}%
}{\overset{\text{q}}{\sum}}$ $\frac{\partial}{\partial x_{i}}(\nabla\log
\theta^{-\frac{1}{2}})_{\text{a}}(y_{0})\frac{\partial\Lambda_{\text{a}}%
}{\partial x_{i}}(y_{0})\phi\circ\pi_{P}(y_{0})$

\qquad$+\frac{1}{6}\underset{i=q+1}{\overset{n}{\sum}}$
$\underset{j=q+1}{\overset{n}{\sum}}$ $\frac{\partial}{\partial x_{i}}%
(\nabla\log\theta^{-\frac{1}{2}})_{j}(y_{0})\frac{\partial\Lambda_{j}%
}{\partial x_{i}}(y_{0})\phi\circ\pi_{P}(y_{0})$

\ For $i,j=q+1,...n,$ we have:

$\qquad\frac{\partial}{\partial x_{i}}(\nabla\log\theta^{-\frac{1}{2}%
})_{\text{a}}(y_{0})=-\frac{1}{2}\perp_{\text{a}ij}(y_{0})<H,j>(y_{0})$ by
(ix)$^{\ast}$ of \textbf{Table A}$_{9}$

\qquad$\frac{\partial}{\partial x_{i}}(\nabla\log\theta^{-\frac{1}{2}}%
)_{j}(y_{0})$ by (ix)$^{\ast\ast}$ of \textbf{Table A}$_{9}$

$\qquad=\frac{1}{2}<H,i><H,j>\ +\frac{1}{12}(2\varrho_{ij}%
+4\overset{q}{\underset{\text{a}=1}{\sum}}R_{i\text{a}j\text{a}}%
-6\overset{q}{\underset{\text{a,b}=1}{\sum}}T_{\text{aa}i}T_{\text{bb}%
j}-T_{\text{ab}i}T_{\text{ab}j})(y_{0})$

\qquad$=\frac{1}{6}[3<H,i><H,j>\ +(\varrho_{ij}%
+2\overset{q}{\underset{\text{a}=1}{\sum}}R_{i\text{a}j\text{a}}%
-3\overset{q}{\underset{\text{a,b}=1}{\sum}}T_{\text{aa}i}T_{\text{bb}%
j}-T_{\text{ab}i}T_{\text{ab}j})](y_{0})$

\qquad$\frac{\partial\Lambda_{\text{a}}}{\partial x_{i}}(y_{0})=\Omega
_{i\text{a}}(y_{0})+[\Lambda_{\text{a}},\Lambda_{i}](y_{0})$ by (vii)
\textbf{Proposition 5}

\qquad$\frac{\partial\Lambda_{j}}{\partial x_{i}}(y_{0})=\frac{1}{2}%
\Omega_{ij}(y_{0})$ by (x) of \textbf{Proposition 5}

\qquad$=\frac{1}{6}\underset{i,j=q+1}{\overset{\text{n}}{\sum}}$
$\frac{\partial}{\partial x_{i}}(\nabla\log\theta^{-\frac{1}{2}})_{j}%
(y_{0})\frac{\partial\Lambda_{j}}{\partial x_{i}}(y_{0})\phi\circ\pi_{P}%
(y_{0})$

$\qquad$I$_{32913}=\frac{1}{6}\underset{i=q+1}{\overset{n}{\sum}}$
$\underset{\text{a=1}}{\overset{\text{q}}{\sum}}$ $\frac{\partial}{\partial
x_{i}}(\nabla\log\theta^{-\frac{1}{2}})_{\text{a}}(y_{0})\frac{\partial
\Lambda_{\text{a}}}{\partial x_{i}}(y_{0})\phi(y_{0})$

\qquad$\qquad+\frac{1}{6}\underset{i=q+1}{\overset{n}{\sum}}$
$\underset{j=q+1}{\overset{n}{\sum}}$ $\frac{\partial}{\partial x_{i}}%
(\nabla\log\theta^{-\frac{1}{2}})_{j}(y_{0})\frac{\partial\Lambda_{j}%
}{\partial x_{i}}(y_{0})\phi(y_{0})$

$\left(  D_{34}\right)  \qquad$I$_{32913}=-\frac{1}{12}%
\underset{i=q+1}{\overset{n}{\sum}}\underset{\text{a=1}}{\overset{\text{q}%
}{\sum}}\perp_{\text{a}ij}(y_{0})<H,j>(y_{0})[\Omega_{i\text{a}}%
+[\Lambda_{\text{a}},\Lambda_{i}]](y_{0})\phi(y_{0})$

$\qquad\qquad\qquad+$ $\frac{1}{72}\underset{i,j=q+1}{\overset{n}{\sum}%
}[3<H,i><H,j>$

$\qquad\qquad\qquad\ +(\varrho_{ij}+2\overset{q}{\underset{\text{a}=1}{\sum}%
}R_{i\text{a}j\text{a}}-3\overset{q}{\underset{\text{a,b}=1}{\sum}%
}T_{\text{aa}i}T_{\text{bb}j}-T_{\text{ab}i}T_{\text{ab}j})](y_{0})\Omega
_{ij}(y_{0})\phi(y_{0})$

\qquad\qquad\qquad\qquad\qquad\qquad\qquad\qquad\qquad\qquad\qquad\qquad
\qquad\qquad\qquad\qquad\qquad$\blacksquare$

Collecting terms, we have by $\left(  D_{32}\right)  ,\left(  D_{33}\right)
,\left(  D_{34}\right)  :$

$\left(  D_{35}\right)  \qquad$I$_{3291}=$ I$_{32911}$ + I$_{32912}$ +
I$_{32913}\qquad$

$=\frac{1}{24}\underset{i,j=q+1}{\overset{n}{\sum}}\underset{\text{a=1}%
}{\overset{\text{q}}{\sum}}<H,j>[4\underset{\text{c}=1}{\overset{q}{\sum}%
}(T_{\text{ac}i})(\perp_{j\text{c}i})+\frac{8}{3}R_{i\text{a}ij}%
](y_{0})\Lambda_{\text{a}}(y_{0})\phi(y_{0})\qquad$I$_{32911}$

$\ +\frac{1}{12}\underset{i,j=q+1}{\overset{n}{\sum}}\underset{\text{a=1}%
}{\overset{\text{q}}{\sum}}\perp_{\text{a}ji}(y_{0})[<H,i><H,j>](y_{0}%
)\Lambda_{\text{a}}(y_{0})\phi(y_{0})$

$\ +\frac{1}{72}\underset{i,j=q+1}{\overset{n}{\sum}}\underset{\text{a=1}%
}{\overset{\text{q}}{\sum}}\perp_{\text{a}ji}(y_{0})[2\varrho_{ij}%
+4\overset{q}{\underset{\text{a}=1}{\sum}}R_{i\text{a}j\text{a}}%
-6\overset{q}{\underset{\text{b,c}=1}{\sum}}T_{\text{cc}i}T_{\text{bb}%
j}-T_{\text{bc}i}T_{\text{bc}j}](y_{0})\Lambda_{\text{a}}(y_{0})\phi(y_{0})$

$+\frac{1}{36}\underset{i,j=q+1}{\overset{n}{\sum}}%
\underset{k=q+1}{\overset{n}{\sum}}<H,k>(y_{0})$R$_{ijik}(y_{0})\Lambda
_{j}(y_{0})\phi(y_{0})$

$-\frac{1}{288}\underset{i,j=q+1}{\overset{n}{\sum}}<H,j>(y_{0})[3<H,i>^{2}$

$+2(\tau^{M}-3\tau^{P}+\overset{q}{\underset{\text{a}=1}{\sum}}\varrho
_{\text{aa}}+\overset{q}{\underset{\text{a,b}=1}{\sum}}R_{\text{abab}}%
)](y_{0})\Lambda_{j}(y_{0})\phi(y_{0})$

$-\frac{1}{12}\underset{i,j=q+1}{\overset{n}{\sum}}<H,i>(y_{0})$

$\times\lbrack\frac{3}{4}<H,i><H,j>$\ $+\frac{1}{6}(\varrho_{ij}%
+2\overset{q}{\underset{\text{a}=1}{\sum}}R_{i\text{a}j\text{a}}%
-3\overset{q}{\underset{\text{a,b=1}}{\sum}}T_{\text{aa}i}T_{\text{bb}%
j}-T_{\text{ab}i}T_{\text{ab}j})](y_{0})\Lambda_{j}(y_{0})\phi(y_{0})$

$+\frac{5}{32}\underset{i,j=q+1}{\overset{n}{\sum}}$%
$<$%
H,$i$%
$>$%
$^{2}$%
$<$%
H,$j$%
$>$%
$\Lambda_{j}(y_{0})\phi(y_{0})$

$+\frac{1}{48}\underset{i,j=q+1}{\overset{n}{\sum}}$%
$<$%
H,$i$%
$>$%
(y$_{0}$)$[$(2$\varrho_{ij}$+4$\overset{q}{\underset{\text{a}=1}{\sum}}%
$R$_{i\text{a}j\text{a}}$-3$\overset{q}{\underset{\text{a,b=1}}{\sum}}%
$T$_{\text{aa}i}$T$_{\text{bb}j}$-T$_{\text{ab}i}$T$_{\text{ab}j}%
$-3$\overset{q}{\underset{\text{a,b=1}}{\sum}}$T$_{\text{aa}j}$T$_{\text{bb}%
i}$-T$_{\text{ab}j}$T$_{\text{ab}i}$)$](y_{0})\Lambda_{j}(y_{0})\phi(y_{0})$

$+\frac{1}{48}\underset{i,j=q+1}{\overset{n}{\sum}}<H,j>[$ $\tau^{M}-3\tau
^{P}+\overset{q}{\underset{\text{a}=1}{\sum}}\varrho_{\text{aa}}%
+\overset{q}{\underset{\text{a,b}=1}{\sum}}R_{\text{abab}}](y_{0})\Lambda
_{j}(y_{0})\phi(y_{0})$

$+\frac{1}{144}\underset{i,j=q+1}{\overset{n}{\sum}}[\nabla_{i}\varrho
_{ij}-2\varrho_{ij}<H,i>+\overset{q}{\underset{\text{a}=1}{\sum}}(\nabla
_{i}R_{\text{a}i\text{a}j}-4R_{i\text{a}j\text{a}}<H,i>)$

$+4\overset{q}{\underset{\text{a,b=1}}{\sum}}R_{i\text{a}j\text{b}%
}T_{\text{ab}i}$+2$\overset{q}{\underset{\text{a,b,c=1}}{\sum}}(T_{\text{aa}%
i}T_{\text{bb}j}T_{\text{cc}i}-3T_{\text{aa}i}T_{\text{bc}j}T_{\text{bc}%
i}+2T_{\text{ab}i}T_{\text{bc}j}T_{\text{ca}i})](y_{0})\Lambda_{j}(y_{0}%
)\phi(y_{0})$\qquad\qquad\qquad\qquad\qquad\ \ 

$+\frac{1}{144}\underset{i,j=q+1}{\overset{n}{\sum}}[\nabla_{j}\varrho
_{ii}-2\varrho_{ji}<H,i>+\overset{q}{\underset{\text{a}=1}{\sum}}(\nabla
_{j}R_{\text{a}i\text{a}i}-4R_{j\text{a}i\text{a}}<H,i>)$

$+4\overset{q}{\underset{\text{a,b=1}}{\sum}}R_{j\text{a}i\text{b}%
}T_{\text{ab}i}$+2$\overset{q}{\underset{\text{a,b,c=1}}{\sum}}(T_{\text{aa}%
j}T_{\text{bb}i}T_{\text{cc}i}-3T_{\text{aa}j}T_{\text{bc}i}T_{\text{bc}%
i}+2T_{\text{ab}j}T_{\text{bc}i}T_{\text{ca}i})](y_{0})\Lambda_{j}(y_{0}%
)\phi(y_{0})$

$+\frac{1}{144}\underset{i,j=q+1}{\overset{n}{\sum}}[\nabla_{i}\varrho
_{ij}-2\varrho_{ii}<H,j>+\overset{q}{\underset{\text{a}=1}{\sum}}(\nabla
_{i}R_{\text{a}i\text{a}j}-4R_{i\text{a}i\text{a}}<H,j>)$

$+4\overset{q}{\underset{\text{a,b=1}}{\sum}}R_{i\text{a}i\text{b}%
}T_{\text{ab}j}$+2$\overset{q}{\underset{\text{a,b,c=1}}{\sum}}(T_{\text{aa}%
i}T_{\text{bb}i}T_{\text{cc}j}-3T_{\text{aa}i}T_{\text{bc}i}T_{\text{bc}%
j}+2T_{\text{ab}i}T_{\text{bc}i}T_{\text{ac}j})](y_{0})\Lambda_{j}(y_{0}%
)\phi(y_{0}).$

$+\frac{1}{72}\underset{i,j=q+1}{\overset{n}{\sum}}$ $<H,j>(y_{0}%
)\frac{\partial\Omega_{ij}}{\partial\text{x}_{i}}(y_{0})\phi(y_{0}%
)\qquad\qquad$I$_{32912}$

$-\frac{1}{12}\underset{i=q+1}{\overset{n}{\sum}}\underset{\text{a=1}%
}{\overset{\text{q}}{\sum}}\perp_{\text{a}ij}(y_{0})<H,j>(y_{0})[\Omega
_{i\text{a}}+[\Lambda_{\text{a}},\Lambda_{i}]](y_{0})\phi(y_{0})\qquad\qquad
$I$_{32913}$

$+$ $\frac{1}{72}\underset{i,j=q+1}{\overset{n}{\sum}}[3<H,i><H,j>$

$\ +(\varrho_{ij}+2\overset{q}{\underset{\text{a}=1}{\sum}}R_{i\text{a}%
j\text{a}}-3\overset{q}{\underset{\text{a,b}=1}{\sum}}T_{\text{aa}%
i}T_{\text{bb}j}-T_{\text{ab}i}T_{\text{ab}j})](y_{0})\Omega_{ij}(y_{0}%
)\phi(y_{0})$

\qquad\qquad\qquad\qquad\qquad\qquad\qquad\qquad\qquad\qquad\qquad\qquad
\qquad\qquad\qquad\qquad\qquad\qquad\qquad\qquad\qquad$\blacksquare$

Next we have:

I$_{3292}=\frac{1}{12}\underset{i=q+1}{\overset{n}{\sum}}$
$\underset{j=1}{\overset{n}{\sum}}$ $\frac{\partial^{2}}{\partial x_{i}^{2}%
}[(\nabla\log\Phi)_{j}\Lambda_{j}\phi\circ\pi_{\text{P}}](y_{0})$

$=\frac{1}{12}\underset{i=q+1}{\overset{n}{\sum}}$ $\underset{\text{a=1}%
}{\overset{\text{q}}{\sum}}$ $\frac{\partial^{2}}{\partial x_{i}^{2}}%
[(\nabla\log\Phi)_{\text{a}}\Lambda_{\text{a}}\phi\circ\pi_{\text{P}}%
](y_{0})+\frac{1}{12}\underset{i=q+1}{\overset{n}{\sum}}$
$\underset{j=q+1}{\overset{n}{\sum}}$ $\frac{\partial^{2}}{\partial x_{i}^{2}%
}[(\nabla\log\Phi)_{j}\Lambda_{j}\phi\circ\pi_{\text{P}}](y_{0})$

We set:

I$_{32921}=\frac{1}{12}\underset{i=q+1}{\overset{n}{\sum}}$
$\underset{\text{a=1}}{\overset{\text{q}}{\sum}}$ $\frac{\partial^{2}%
}{\partial x_{i}^{2}}[(\nabla\log\Phi)_{\text{a}}\Lambda_{\text{a}}\phi
\circ\pi_{\text{P}}](y_{0})$

I$_{32922}=\frac{1}{12}\underset{i=q+1}{\overset{n}{\sum}}$
$\underset{j=q+1}{\overset{n}{\sum}}$ $\frac{\partial^{2}}{\partial x_{i}^{2}%
}[(\nabla\log\Phi)_{j}\Lambda_{j}\phi\circ\pi_{\text{P}}](y_{0})$

We carry out the computations:

I$_{32921}=\frac{1}{12}\underset{i=q+1}{\overset{n}{\sum}}$
$\underset{\text{a=1}}{\overset{\text{q}}{\sum}}$ $\frac{\partial^{2}%
}{\partial x_{i}^{2}}[(\nabla\log\Phi)_{\text{a}}\Lambda_{\text{a}}\phi
\circ\pi_{\text{P}}](y_{0})$

\qquad$=\frac{1}{12}\underset{i=q+1}{\overset{n}{\sum}}$ $\underset{\text{a=1}%
}{\overset{\text{q}}{\sum}}$ $\frac{\partial^{2}}{\partial x_{i}^{2}}%
[(\nabla\log\Phi)_{\text{a}}](y_{0})[\Lambda_{\text{a}}\phi\circ\pi_{\text{P}%
}](y_{0})$

$\qquad+\frac{1}{12}$ $\underset{\text{a=1}}{\overset{\text{q}}{\sum}}$
$(\nabla\log\Phi)_{\text{a}}(y_{0})\frac{\partial^{2}}{\partial x_{i}^{2}%
}[\Lambda_{\text{a}}\phi\circ\pi_{\text{P}}](y_{0})$

\qquad$+\frac{1}{6}\underset{i=q+1}{\overset{n}{\sum}}$ $\underset{\text{a=1}%
}{\overset{\text{q}}{\sum}}$ $\frac{\partial}{\partial x_{i}}[(\nabla\log
\Phi)_{\text{a}}](y_{0})\frac{\partial}{\partial x_{i}}[\Lambda_{\text{a}}%
\phi\circ\pi_{\text{P}}](y_{0})$

Since $\frac{\partial^{2}}{\partial x_{i}^{2}}[\phi\circ\pi_{\text{P}}%
](y_{0})=0=\frac{\partial}{\partial x_{i}}[\phi\circ\pi_{\text{P}}](y_{0}),$
we have:

I$_{32921}=\frac{1}{12}\underset{i=q+1}{\overset{n}{\sum}}$
$\underset{\text{a=1}}{\overset{\text{q}}{\sum}}$ $\frac{\partial^{2}%
}{\partial x_{i}^{2}}(\nabla\log\Phi)_{\text{a}}(y_{0})\Lambda_{\text{a}%
}(y_{0})\phi(y_{0})$

$\qquad+\frac{1}{12}\underset{i=q+1}{\overset{n}{\sum}}$ $\underset{\text{a=1}%
}{\overset{\text{q}}{\sum}}$ $(\nabla\log\Phi)_{\text{a}}(y_{0})\frac
{\partial^{2}\Lambda_{\text{a}}}{\partial x_{i}^{2}}(y_{0})\phi(y_{0})$

\qquad$+\frac{1}{6}$ $\underset{i=q+1}{\overset{n}{\sum}}\underset{\text{a=1}%
}{\overset{\text{q}}{\sum}}$ $\frac{\partial}{\partial x_{i}}(\nabla\log
\Phi)_{\text{a}}(y_{0})\frac{\partial\Lambda_{\text{a}}}{\partial x_{i}}%
(y_{0})\phi(y_{0})$

$\qquad(\nabla\log\Phi)_{\text{a}}(y_{0})=0$ by (xi) of \textbf{Table B}%
$_{1}.$

$\qquad\frac{\partial}{\partial x_{i}}(\nabla$log$\Phi_{P})_{\text{a}}%
(y_{0})=$ $\underset{j=q+1}{\overset{n}{%
{\textstyle\sum}
}}X_{j}(y_{0})\perp_{\text{a}ij}(y_{0})-\frac{\partial X_{i}}{\partial
x_{\text{a}}}(y_{0})$ by (xv) of \textbf{Table B}$_{1}\qquad$

By (xvi) of \textbf{Table B}$_{1},$

$\qquad\frac{\partial^{2}}{\partial x_{i}^{2}}(\nabla$log$\Phi_{P})_{\text{a}%
}(y_{0})=-4\underset{\text{b=1}}{\overset{\text{q}}{%
{\textstyle\sum}
}}T_{\text{ab}i}(y)\frac{\partial X_{i}}{\partial x_{\text{b}}}(y_{0}%
)+\underset{k=q+1}{\overset{n}{%
{\textstyle\sum}
}}\perp_{\text{a}ik}(y)\left[  \left(  \frac{\partial X_{i}}{\partial x_{k}%
}+\frac{\partial X_{k}}{\partial x_{i}}\right)  \right]  (y_{0})$

$\qquad\qquad\qquad\qquad+\frac{8}{3}\underset{k=q+1}{\overset{n}{%
{\textstyle\sum}
}}R_{i\text{a}ik}(y)X_{k}(y_{0})+[2X_{i}\frac{\partial X_{i}}{\partial
x_{\text{a}}}-\frac{\partial^{2}X_{i}}{\partial x_{\text{a}}\partial x_{i}%
}](y_{0})$

\qquad\qquad\qquad\qquad$=-4\underset{\text{b=1}}{\overset{\text{q}}{%
{\textstyle\sum}
}}T_{\text{ab}i}(y)\frac{\partial X_{i}}{\partial x_{\text{b}}}(y_{0}%
)+\underset{j=q+1}{\overset{n}{%
{\textstyle\sum}
}}\perp_{\text{a}ij}(y)\left(  \frac{\partial X_{i}}{\partial x_{j}}%
+\frac{\partial X_{j}}{\partial x_{i}}\right)  (y_{0})$

$\qquad\qquad\qquad\qquad+\frac{8}{3}\underset{j=q+1}{\overset{n}{%
{\textstyle\sum}
}}R_{i\text{a}ij}(y_{0})X_{j}(y_{0})+\left(  2X_{i}\frac{\partial X_{i}%
}{\partial x_{\text{a}}}-\frac{\partial^{2}X_{i}}{\partial x_{\text{a}%
}\partial x_{i}}\right)  (y_{0})$

Therefore,

I$_{32921}=\frac{1}{12}\underset{i=q+1}{\overset{n}{\sum}}$
$\underset{\text{a=1}}{\overset{\text{q}}{\sum}}$ $[-4\underset{\text{b=1}%
}{\overset{\text{q}}{%
{\textstyle\sum}
}}T_{\text{ab}i}\frac{\partial X_{i}}{\partial x_{\text{b}}}%
+\underset{j=q+1}{\overset{n}{%
{\textstyle\sum}
}}\perp_{\text{a}ij}\left(  \frac{\partial X_{i}}{\partial x_{j}}%
+\frac{\partial X_{j}}{\partial x_{i}}\right)  $

$\qquad\qquad+\frac{8}{3}\underset{j=q+1}{\overset{n}{%
{\textstyle\sum}
}}R_{i\text{a}ij}X_{j}+\left(  2X_{i}\frac{\partial X_{i}}{\partial
x_{\text{a}}}-\frac{\partial^{2}X_{i}}{\partial x_{\text{a}}\partial x_{i}%
}\right)  ](y_{0})\Lambda_{\text{a}}(y_{0})\phi(y_{0})$

\qquad\qquad$+\frac{1}{6}$ $\underset{i=q+1}{\overset{n}{\sum}}%
\underset{\text{a=1}}{\overset{\text{q}}{\sum}}$
$[\underset{j=q+1}{\overset{n}{%
{\textstyle\sum}
}}X_{j}(y)\perp_{\text{a}ij}+\frac{\partial X_{i}}{\partial x_{\text{a}}}]$
$(y_{0})[\Omega_{i\text{a}}+[\Lambda_{\text{a}},\Lambda_{i}](y_{0})\phi
(y_{0})$

We then compute:

I$_{32922}=\frac{1}{12}$ $\underset{j=q+1}{\overset{n}{\sum}}$ $\frac
{\partial^{2}}{\partial x_{i}^{2}}[(\nabla\log\Phi)_{j}\Lambda_{j}\phi\circ
\pi_{\text{P}}](y_{0})$

$=\frac{1}{12}$ $\underset{j=q+1}{\overset{n}{\sum}}$ $\frac{\partial^{2}%
}{\partial x_{i}^{2}}[(\nabla\log\Phi)_{j}](y_{0})[\Lambda_{j}\phi\circ
\pi_{\text{P}}](y_{0})+\frac{1}{12}$ $\underset{j=q+1}{\overset{n}{\sum}}$
$(\nabla\log\Phi)_{j}(y_{0})\frac{\partial^{2}}{\partial x_{i}^{2}}%
[\Lambda_{j}\phi\circ\pi_{\text{P}}](y_{0})$

\qquad$+\frac{1}{6}$ $\underset{j=q+1}{\overset{n}{\sum}}$ $\frac{\partial
}{\partial x_{i}}(\nabla\log\Phi)_{j}(y_{0})\frac{\partial}{\partial x_{i}%
}[\Lambda_{j}\phi\circ\pi_{\text{P}}](y_{0})$

I$_{32922}=\frac{1}{12}\underset{j=q+1}{\overset{n}{\sum}}\frac{\partial^{2}%
}{\partial x_{i}^{2}}[(\nabla\log\Phi)_{j}](y_{0})\Lambda_{j}(y_{0})\phi
(y_{0})$

$+\frac{1}{12}$ $\underset{j=q+1}{\overset{n}{\sum}}$ $(\nabla\log\Phi
)_{j}(y_{0})\frac{\partial^{2}\Lambda_{j}}{\partial x_{i}^{2}}(y_{0}%
)\phi(y_{0})$

$+\frac{1}{6}$ $\underset{j=q+1}{\overset{n}{\sum}}$ $\frac{\partial}{\partial
x_{i}}(\nabla\log\Phi)_{j}(y_{0})\frac{\partial\Lambda_{j}}{\partial x_{i}%
}(y_{0})\phi(y_{0})$

$(\nabla\log\Phi)_{j}(y_{0})=-X_{j}(y_{0})$ by (i) of \textbf{Table B}%
$_{1};\frac{\partial}{\partial x_{i}}(\nabla\log\Phi)_{j}(y_{0})=-\frac
{\partial X_{j}}{\partial x_{i}}(y_{0})$ by \ (ii) of Table B$_{1};$

$\frac{\partial\Lambda_{j}}{\partial x_{i}}(y_{0})=\frac{1}{2}\Omega
_{ij}(y_{0})$ by (x) of \textbf{Proposition 5; }\ $\frac{\partial^{2}%
\Lambda_{j}}{\partial x_{i}^{2}}(y_{0})=\frac{1}{3}\frac{\partial\Omega_{ij}%
}{\partial\text{x}_{i}}(y_{0})$ by (xi) of \textbf{Proposition 5.}

By (viii)$^{\ast}$ of \textbf{Appendix B}$_{1}$ we have the formula:

$[\frac{\partial^{2}}{\partial x_{i}\partial x_{j}}(\nabla\log\Phi_{P}%
)_{k}](y_{0})$

$=-\frac{1}{3}\left(  \frac{\partial^{2}X_{i}}{\partial x_{j}\partial x_{k}%
}+\frac{\partial^{2}X_{j}}{\partial x_{i}\partial x_{k}}+\frac{\partial
^{2}X_{k}}{\partial x_{i}\partial x_{j}}\right)  (y_{0})-\frac{1}%
{3}\underset{l=q+1}{\overset{n}{\sum}}[R_{jkil}+R_{ikjl}](y_{0})X_{l}(y_{0})$

We deduce that:

$\frac{\partial^{2}}{\partial x_{i}^{2}}[(\nabla\log\Phi)_{j}](y_{0})$

$=-\frac{1}{3}\left(  \frac{\partial^{2}X_{i}}{\partial x_{i}\partial x_{j}%
}+\frac{\partial^{2}X_{i}}{\partial x_{i}\partial x_{j}}+\frac{\partial
^{2}X_{j}}{\partial x_{i}^{2}}\right)  (y_{0})-\frac{1}{3}%
\underset{k=q+1}{\overset{n}{\sum}}[R_{ijik}+R_{ijik}](y_{0})X_{l}(y_{0})$

We simplify and have:

$\frac{\partial^{2}}{\partial x_{i}^{2}}[(\nabla\log\Phi)_{j}](y_{0}%
)=-\frac{1}{3}\left(  2\frac{\partial^{2}X_{i}}{\partial x_{i}\partial x_{j}%
}+\frac{\partial^{2}X_{j}}{\partial x_{i}^{2}}\right)  (y_{0})-\frac{2}%
{3}\underset{k=q+1}{\overset{n}{\sum}}R_{ijik}(y_{0})X_{k}(y_{0})$\qquad

Therefore,

I$_{32922}=-\frac{1}{36}[\left(  2\frac{\partial^{2}X_{i}}{\partial
x_{i}\partial x_{j}}+\frac{\partial^{2}X_{j}}{\partial x_{i}^{2}}\right)
+2\underset{k=q+1}{\overset{n}{\sum}}R_{ijik}X_{k}](y_{0})\Lambda_{j}%
(y_{0})\phi(y_{0})$

$\qquad\qquad-\frac{1}{36}X_{j}(y_{0})$ $\frac{\partial\Omega_{ij}}{\partial
x_{i}}(y_{0})\phi(y_{0})-\frac{1}{12}\frac{\partial X_{j}}{\partial x_{i}%
}(y_{0})\Omega_{ij}(y_{0})$ $\phi(y_{0})$

Therefore,

$\left(  D_{36}\right)  $\qquad I$_{3292}=$ I$_{32921}+$ I$_{32922}$

$\qquad=\frac{1}{12}\underset{i=q+1}{\overset{n}{\sum}}$ $\underset{\text{a=1}%
}{\overset{\text{q}}{\sum}}$ $[-4\underset{\text{b=1}}{\overset{\text{q}}{%
{\textstyle\sum}
}}T_{\text{ab}i}\frac{\partial X_{i}}{\partial x_{\text{b}}}%
+\underset{j=q+1}{\overset{n}{%
{\textstyle\sum}
}}\perp_{\text{a}ij}\left(  \frac{\partial X_{i}}{\partial x_{j}}%
+\frac{\partial X_{j}}{\partial x_{i}}\right)  \qquad\qquad$I$_{32921}$

$\qquad+\frac{8}{3}\underset{j=q+1}{\overset{n}{%
{\textstyle\sum}
}}R_{i\text{a}ij}X_{j}+\left(  2X_{i}\frac{\partial X_{i}}{\partial
x_{\text{a}}}-\frac{\partial^{2}X_{i}}{\partial x_{\text{a}}\partial x_{i}%
}\right)  ](y_{0})\Lambda_{\text{a}}(y_{0})\phi(y_{0}$

\qquad$+\frac{1}{6}$ $\underset{i=q+1}{\overset{n}{\sum}}\underset{\text{a=1}%
}{\overset{\text{q}}{\sum}}$ $[\underset{j=q+1}{\overset{n}{%
{\textstyle\sum}
}}X_{j}\perp_{\text{a}ij}+\frac{\partial X_{i}}{\partial x_{\text{a}}}]$
$(y_{0})[\Omega_{i\text{a}}+[\Lambda_{\text{a}},\Lambda_{i}](y_{0})\phi
(y_{0})$

\qquad$-\frac{1}{36}\underset{i,j=q+1}{\overset{n}{\sum}}[\left(
2\frac{\partial^{2}X_{i}}{\partial x_{i}\partial x_{j}}+\frac{\partial
^{2}X_{j}}{\partial x_{i}^{2}}\right)  +2\underset{k=q+1}{\overset{n}{\sum}%
}R_{ijik}X_{k}](y_{0})\Lambda_{j}(y_{0})\phi(y_{0})\qquad$I$_{32922}$

$\qquad-\frac{1}{36}\underset{i,j=q+1}{\overset{n}{\sum}}X_{j}(y_{0})$
$\frac{\partial\Omega_{ij}}{\partial x_{i}}(y_{0})\phi(y_{0})-\frac{1}%
{12}\underset{i,j=q+1}{\overset{n}{\sum}}\frac{\partial X_{j}}{\partial x_{i}%
}(y_{0})\Omega_{ij}(y_{0})$ $\phi(y_{0})$

\qquad\qquad\qquad\qquad\qquad\qquad\qquad\qquad\qquad\qquad\qquad\qquad
\qquad\qquad\qquad\qquad\qquad$\blacksquare$

We conclude by $\left(  D_{35}\right)  $ and $\left(  D_{36}\right)  $ that:

$\qquad\left(  D_{37}\right)  $\qquad I$_{329}=$ I$_{3291}+$ I$_{3292}\qquad$

$=\frac{1}{24}\underset{i,j=q+1}{\overset{n}{\sum}}\underset{\text{a=1}%
}{\overset{\text{q}}{\sum}}<H,j>(y_{0})[4\underset{\text{c}%
=1}{\overset{q}{\sum}}(T_{\text{ac}i})(\perp_{j\text{c}i})+\frac{8}%
{3}R_{i\text{a}ij}](y_{0})\Lambda_{\text{a}}(y_{0})\phi(y_{0})\qquad$%
I$_{3291}\qquad$I$_{32911}\qquad$

$+\frac{1}{12}\underset{i,j=q+1}{\overset{n}{\sum}}\underset{\text{a=1}%
}{\overset{\text{q}}{\sum}}\perp_{\text{a}ji}(y_{0})[<H,i><H,j>](y_{0}%
)\Lambda_{\text{a}}(y_{0})\phi(y_{0})$

$+\frac{1}{72}\underset{i,j=q+1}{\overset{n}{\sum}}\underset{\text{a=1}%
}{\overset{\text{q}}{\sum}}\perp_{\text{a}ji}(y_{0})[2\varrho_{ij}%
+4\overset{q}{\underset{\text{a}=1}{\sum}}R_{i\text{a}j\text{a}}%
-6\overset{q}{\underset{\text{b,c}=1}{\sum}}T_{\text{cc}i}T_{\text{bb}%
j}-T_{\text{bc}i}T_{\text{bc}j}](y_{0})\Lambda_{\text{a}}(y_{0})\phi(y_{0})$

$+\frac{1}{36}\underset{i,j=q+1}{\overset{n}{\sum}}%
\underset{k=q+1}{\overset{n}{\sum}}<H,k>(y_{0})$R$_{ijik}(y_{0})\Lambda
_{j}(y_{0})\phi(y_{0})$

$-\frac{1}{288}\underset{i,j=q+1}{\overset{n}{\sum}}<H,j>(y_{0})$

$\times\lbrack3<H,i>^{2}+2(\tau^{M}-3\tau^{P}+\overset{q}{\underset{\text{a}%
=1}{\sum}}\varrho_{\text{aa}}+\overset{q}{\underset{\text{a,b}=1}{\sum}%
}R_{\text{abab}})](y_{0})\Lambda_{j}(y_{0})\phi(y_{0})$

$-\frac{1}{12}\underset{i,j=q+1}{\overset{n}{\sum}}<H,i>(y_{0})$

$\times\lbrack\frac{3}{4}<H,i><H,j>$\ $+\frac{1}{6}(\varrho_{ij}%
+2\overset{q}{\underset{\text{a}=1}{\sum}}R_{i\text{a}j\text{a}}%
-3\overset{q}{\underset{\text{a,b=1}}{\sum}}T_{\text{aa}i}T_{\text{bb}%
j}-T_{\text{ab}i}T_{\text{ab}j})](y_{0})\Lambda_{j}(y_{0})\phi(y_{0})$

$+\frac{5}{32}\underset{i,j=q+1}{\overset{n}{\sum}}$%
$<$%
H,$i$%
$>$%
$^{2}$%
$<$%
H,$j$%
$>$%
$\Lambda_{j}(y_{0})\phi(y_{0})$

$+\frac{1}{48}\underset{i,j=q+1}{\overset{n}{\sum}}$%
$<$%
H,$i$%
$>$%
(y$_{0}$)

$\times\lbrack$(2$\varrho_{ij}$+4$\overset{q}{\underset{\text{a}=1}{\sum}}%
$R$_{i\text{a}j\text{a}}$-3$\overset{q}{\underset{\text{a,b=1}}{\sum}}%
$T$_{\text{aa}i}$T$_{\text{bb}j}$-T$_{\text{ab}i}$T$_{\text{ab}j}%
$-3$\overset{q}{\underset{\text{a,b=1}}{\sum}}$T$_{\text{aa}j}$T$_{\text{bb}%
i}$-T$_{\text{ab}j}$T$_{\text{ab}i}$)$](y_{0})\Lambda_{j}(y_{0})\phi(y_{0})$

$+\frac{1}{48}\underset{i,j=q+1}{\overset{n}{\sum}}<H,j>(y_{0}\times\lbrack$
$\tau^{M}-3\tau^{P}+\overset{q}{\underset{\text{a}=1}{\sum}}\varrho
_{\text{aa}}+\overset{q}{\underset{\text{a,b}=1}{\sum}}R_{\text{abab}}%
](y_{0})\Lambda_{j}(y_{0})\phi(y_{0})$

$+\frac{1}{144}\underset{i,j=q+1}{\overset{n}{\sum}}[\nabla_{i}\varrho
_{ij}-2\varrho_{ij}<H,i>+\overset{q}{\underset{\text{a}=1}{\sum}}(\nabla
_{i}R_{\text{a}i\text{a}j}-4R_{i\text{a}j\text{a}}<H,i>)$

$+4\overset{q}{\underset{\text{a,b=1}}{\sum}}R_{i\text{a}j\text{b}%
}T_{\text{ab}i}+2\overset{q}{\underset{\text{a,b,c=1}}{\sum}}(T_{\text{aa}%
i}T_{\text{bb}j}T_{\text{cc}i}-3T_{\text{aa}i}T_{\text{bc}j}T_{\text{bc}%
i}+2T_{\text{ab}i}T_{\text{bc}j}T_{\text{ca}i})](y_{0})\Lambda_{j}(y_{0}%
)\phi(y_{0})$\qquad\qquad\qquad\qquad\qquad\ \ 

$+\frac{1}{144}\underset{i,j=q+1}{\overset{n}{\sum}}[\nabla_{j}\varrho
_{ii}-2\varrho_{ji}<H,i>+\overset{q}{\underset{\text{a}=1}{\sum}}(\nabla
_{j}R_{\text{a}i\text{a}i}-4R_{j\text{a}i\text{a}}<H,i>)$

$+4\overset{q}{\underset{\text{a,b=1}}{\sum}}R_{j\text{a}i\text{b}%
}T_{\text{ab}i}$+2$\overset{q}{\underset{\text{a,b,c=1}}{\sum}}(T_{\text{aa}%
j}T_{\text{bb}i}T_{\text{cc}i}-3T_{\text{aa}j}T_{\text{bc}i}T_{\text{bc}%
i}+2T_{\text{ab}j}T_{\text{bc}i}T_{\text{ca}i})](y_{0})\Lambda_{j}(y_{0}%
)\phi(y_{0})$

$+\frac{1}{144}\underset{i,j=q+1}{\overset{n}{\sum}}[\nabla_{i}\varrho
_{ij}-2\varrho_{ii}<H,j>+\overset{q}{\underset{\text{a}=1}{\sum}}(\nabla
_{i}R_{\text{a}i\text{a}j}-4R_{i\text{a}i\text{a}}<H,j>)$

$+4\overset{q}{\underset{\text{a,b=1}}{\sum}}R_{i\text{a}i\text{b}%
}T_{\text{ab}j}$+2$\overset{q}{\underset{\text{a,b,c=1}}{\sum}}(T_{\text{aa}%
i}T_{\text{bb}i}T_{\text{cc}j}-3T_{\text{aa}i}T_{\text{bc}i}T_{\text{bc}%
j}+2T_{\text{ab}i}T_{\text{bc}i}T_{\text{ac}j})](y_{0})\Lambda_{j}(y_{0}%
)\phi(y_{0})$

$+\frac{1}{72}\underset{i,j=q+1}{\overset{n}{\sum}}$ $<H,j>(y_{0}%
)\frac{\partial\Omega_{ij}}{\partial\text{x}_{i}}(y_{0})\phi(y_{0}%
)\qquad\qquad$I$_{32912}$

$-\frac{1}{12}\underset{i=q+1}{\overset{n}{\sum}}\underset{\text{a=1}%
}{\overset{\text{q}}{\sum}}\perp_{\text{a}ij}(y_{0})<H,j>(y_{0})[\Omega
_{i\text{a}}+[\Lambda_{\text{a}},\Lambda_{i}]](y_{0})\phi(y_{0})\qquad\qquad
$I$_{32913}$

$+$ $\frac{1}{72}\underset{i,j=q+1}{\overset{n}{\sum}}[3<H,i><H,j>\ +(\varrho
_{ij}+2\overset{q}{\underset{\text{a}=1}{\sum}}R_{i\text{a}j\text{a}%
}-3\overset{q}{\underset{\text{a,b}=1}{\sum}}T_{\text{aa}i}T_{\text{bb}%
j}-T_{\text{ab}i}T_{\text{ab}j})](y_{0})\Omega_{ij}(y_{0})\phi(y_{0})$

$+\frac{1}{12}\underset{i=q+1}{\overset{n}{\sum}}$ $\underset{\text{a=1}%
}{\overset{\text{q}}{\sum}}$ $[-4\underset{\text{b=1}}{\overset{\text{q}}{%
{\textstyle\sum}
}}T_{\text{ab}i}\frac{\partial X_{i}}{\partial x_{\text{b}}}%
+\underset{j=q+1}{\overset{n}{%
{\textstyle\sum}
}}\perp_{\text{a}ij}\left(  \frac{\partial X_{i}}{\partial x_{j}}%
+\frac{\partial X_{j}}{\partial x_{i}}\right)  \qquad\qquad$I$_{3292}\qquad
$I$_{32921}$

$+\frac{8}{3}\underset{j=q+1}{\overset{n}{%
{\textstyle\sum}
}}R_{i\text{a}ij}X_{j}+\left(  2X_{i}\frac{\partial X_{i}}{\partial
x_{\text{a}}}-\frac{\partial^{2}X_{i}}{\partial x_{\text{a}}\partial x_{i}%
}\right)  ](y_{0})\Lambda_{\text{a}}(y_{0})\phi(y_{0}$

$+\frac{1}{6}$ $\underset{i=q+1}{\overset{n}{\sum}}\underset{\text{a=1}%
}{\overset{\text{q}}{\sum}}$ $[\underset{j=q+1}{\overset{n}{%
{\textstyle\sum}
}}X_{j}\perp_{\text{a}ij}+\frac{\partial X_{i}}{\partial x_{\text{a}}}]$
$(y_{0})[\Omega_{i\text{a}}+[\Lambda_{\text{a}},\Lambda_{i}](y_{0})\phi
(y_{0})$

$-\frac{1}{36}[\left(  2\frac{\partial^{2}X_{i}}{\partial x_{i}\partial x_{j}%
}+\frac{\partial^{2}X_{j}}{\partial x_{i}^{2}}\right)
+2\underset{k=q+1}{\overset{n}{\sum}}R_{ijik}X_{k}](y_{0})\Lambda_{j}%
(y_{0})\phi(y_{0})\qquad$I$_{32922}$

$-\frac{1}{36}X_{j}(y_{0})$ $\frac{\partial\Omega_{ij}}{\partial x_{i}}%
(y_{0})\phi(y_{0})-\frac{1}{12}\frac{\partial X_{j}}{\partial x_{i}}%
(y_{0})\Omega_{ij}(y_{0})$ $\phi(y_{0})$

\qquad\qquad\qquad\qquad\qquad\qquad\qquad\qquad\qquad\qquad\qquad\qquad
\qquad\qquad\qquad$\blacksquare$

(x)\qquad\textbf{L}$_{1}=\frac{1}{12}$ $\underset{\text{a=1}%
}{\overset{\text{q}}{\sum}}\frac{\partial^{2}}{\partial x_{i}^{2}}%
[$X$_{\text{a}}\frac{\partial\phi}{\partial\text{x}_{\text{a}}}\circ
\pi_{\text{P}}](y_{0})$

We recall for the last time that $\frac{\partial^{2}}{\partial x_{i}^{2}}%
\pi_{\text{P}}(y_{0})=0=\frac{\partial}{\partial x_{i}}\pi_{\text{P}}(y_{0})$
for $i=q+1,...,n.$

$\left(  D_{38}\right)  \qquad$\textbf{L}$_{1}=\frac{1}{12}$
$\underset{\text{a=1}}{\overset{\text{q}}{\sum}}\frac{\partial^{2}X_{\text{a}%
}}{\partial x_{i}^{2}}(y_{0})[\frac{\partial\phi}{\partial\text{x}_{\text{a}}%
}\circ\pi_{\text{P}}](y_{0})=\frac{1}{12}$ $\underset{\text{a=1}%
}{\overset{\text{q}}{\sum}}\frac{\partial^{2}X_{\text{a}}}{\partial x_{i}^{2}%
}(y_{0})\frac{\partial\phi}{\partial\text{x}_{\text{a}}}(y_{0})$

\qquad\qquad\qquad\qquad\qquad\qquad\qquad\qquad\qquad\qquad\qquad\qquad
\qquad\qquad\qquad$\blacksquare$

(xi)\qquad\textbf{L}$_{2}=\frac{1}{12}$ $\underset{j=1}{\overset{\text{n}%
}{\sum}}$ $\frac{\partial^{2}}{\partial x_{i}^{2}}[$X$_{j}\Lambda_{j}\phi
\circ\pi_{\text{P}}](y_{0})$

\qquad\qquad$=\frac{1}{12}$ $\underset{j=1}{\overset{\text{n}}{\sum}}$
$\frac{\partial^{2}X_{j}}{\partial x_{i}^{2}}[\Lambda_{j}\phi\circ
\pi_{\text{P}}](y_{0})+\frac{1}{12}$ $\underset{j=1}{\overset{\text{n}}{\sum}%
}$ $[$X$_{j}\frac{\partial^{2}}{\partial x_{i}^{2}}[\Lambda_{j}\phi\circ
\pi_{\text{P}}](y_{0})$

$\qquad\qquad+\frac{1}{6}$ $\underset{j=1}{\overset{\text{n}}{\sum}}$
$\frac{\partial X_{j}}{\partial x_{i}}\frac{\partial}{\partial x_{i}}%
[\Lambda_{j}\phi\circ\pi_{\text{P}}](y_{0})=$ \textbf{L}$_{21}+$
\textbf{L}$_{22}+$ \textbf{L}$_{23}$

where,

\qquad\qquad\textbf{L}$_{21}=\frac{1}{12}$ $\underset{j=1}{\overset{\text{n}%
}{\sum}}$ $\frac{\partial^{2}X_{j}}{\partial x_{i}^{2}}[\Lambda_{j}\phi
\circ\pi_{\text{P}}](y_{0})$

\qquad\qquad\qquad$=\frac{1}{12}$ $\underset{\text{a}=1}{\overset{\text{q}%
}{\sum}}$ $\frac{\partial^{2}X_{\text{a}}}{\partial x_{i}^{2}}(y_{0}%
)\Lambda_{\text{a}}(y_{0})\phi(y_{0})+\frac{1}{12}$
$\underset{j=q+1}{\overset{\text{n}}{\sum}}$ $\frac{\partial^{2}X_{j}%
}{\partial x_{i}^{2}}(y_{0})\Lambda_{j}(y_{0})\phi(y_{0})$

\qquad\qquad\textbf{L}$_{22}=\frac{1}{12}$ $\underset{j=1}{\overset{\text{n}%
}{\sum}}$ $[$X$_{j}\frac{\partial^{2}}{\partial x_{i}^{2}}(\Lambda_{j}%
\phi\circ\pi_{\text{P}})](y_{0})$

\qquad\qquad\qquad$=\frac{1}{12}$ $\underset{j=1}{\overset{\text{n}}{\sum}}$
$[$X$_{j}\frac{\partial^{2}\Lambda_{j}}{\partial x_{i}^{2}}(\phi\circ
\pi_{\text{P}})](y_{0})$

$\qquad\frac{\partial^{2}\Lambda_{j}}{\partial x_{i}^{2}}(y_{0})=\frac{1}%
{3}\frac{\partial\Omega_{ij}}{\partial x_{i}}(y_{0})$ by (xi) of
\textbf{Proposition 5}

$\qquad\frac{\partial\Lambda_{j}}{\partial x_{i}}(y_{0})=\frac{1}{2}%
\Omega_{ij}(y_{0})$ by (x) of \textbf{Proposition 5}

\qquad\textbf{L}$_{22}=\frac{1}{36}$ $\underset{j=q+1}{\overset{n}{\sum}}$
X$_{j}(y_{0})\frac{\partial\Omega_{ij}}{\partial x_{i}}(y_{0})\phi(y_{0})$

\qquad\qquad\qquad\qquad\qquad\qquad\qquad\qquad\qquad$\qquad\qquad
\qquad\qquad\qquad\blacksquare$

\qquad\textbf{L}$_{23}=\frac{1}{6}$ $\underset{j=1}{\overset{n}{\sum}}$
$\frac{\partial X_{j}}{\partial x_{i}}\frac{\partial}{\partial x_{i}}%
[\Lambda_{j}\phi\circ\pi_{\text{P}}](y_{0})$

\qquad$=\frac{1}{6}$ $\underset{j=1}{\overset{n}{\sum}}$ $\frac{\partial
X_{j}}{\partial x_{i}}(y_{0})\frac{\partial\Lambda_{j}}{\partial x_{i}}%
(y_{0})\phi(y_{0})=\frac{1}{12}$ $\underset{j=q+1}{\overset{n}{\sum}}$
$\frac{\partial X_{j}}{\partial x_{i}}(y_{0})\Omega_{ij}(y_{0})\phi(y_{0})$

\qquad\qquad\qquad\qquad\qquad\qquad\qquad\qquad\qquad\qquad\qquad\qquad
\qquad\qquad$\blacksquare$\qquad\qquad\qquad\qquad\qquad\qquad\qquad
\qquad\qquad\qquad\qquad\qquad\qquad\qquad

Therefore,

$\left(  D_{39}\right)  $\qquad\textbf{L}$_{2}=$ \textbf{L}$_{21}+$
\textbf{L}$_{22}+$ \textbf{L}$_{23}$

\qquad$=\frac{1}{12}$ $\underset{\text{a}=1}{\overset{\text{q}}{\sum}}$
$\frac{\partial^{2}X_{\text{a}}}{\partial x_{i}^{2}}(y_{0})\Lambda_{\text{a}%
}(y_{0})\phi(y_{0})+\frac{1}{12}\underset{j=q+1}{\overset{n}{\sum}}%
\frac{\partial^{2}X_{j}}{\partial x_{i}^{2}}(y_{0})\Lambda_{j}(y_{0}%
)\phi(y_{0})\qquad$\textbf{L}$_{21}$

\qquad$+\frac{1}{36}$ $\underset{j=q+1}{\overset{n}{\sum}}$ X$_{j}(y_{0}%
)\frac{\partial\Omega_{ij}}{\partial x_{i}}(y_{0})\phi(y_{0})\qquad$%
\textbf{L}$_{22}$

\qquad$+\frac{1}{12}$ $\underset{j=q+1}{\overset{n}{\sum}}$ $\frac{\partial
X_{j}}{\partial x_{i}}(y_{0})\Omega_{ij}(y_{0})\phi(y_{0})\qquad$%
\textbf{L}$_{23}$

\qquad\qquad\qquad\qquad\qquad\qquad\qquad\qquad\qquad\qquad\qquad\qquad
\qquad$\qquad\qquad\qquad\qquad\blacksquare$\qquad

\subsection{Final Expression of I$_{32}$}

We come to the final expression of:

I$_{32}=\frac{1}{12}\underset{i=q+1}{\overset{n}{\sum}}\frac{\partial
^{2}\Theta}{\partial x_{i}^{2}}(y_{0})=$ I$_{321}+$ I$_{322}+$ I$_{323}+$
I$_{324}+$ I$_{325}+$ I$_{326}+$ I$_{327}+$ I$_{328}+$ I$_{329}$

The various expressions which define I$_{32}$ are given in $\left(
D_{9}\right)  ^{\ast}:$ I$_{321}$ $\mathbf{=}\frac{1}{12}$ $\frac{\partial
^{2}}{\partial x_{i}^{2}}\left\{  \Psi^{-1}\text{L}\Psi\right\}  (y_{0})$ is
from $\left(  B_{118}\right)  $ of \textbf{Table B}$_{5}.$ Then in, $\left(
D_{10}\right)  ,$ $\left(  D_{14}\right)  ,\left(  D_{18}\right)  ,$ $\left(
D_{21}\right)  ,$ $\left(  D_{30}\right)  ,$ $\left(  D_{31}\right)  ,\left(
D_{32}\right)  ,$ $\left(  D_{37}\right)  ,$ $\left(  D_{38}\right)  ,$
$\left(  D_{39}\right)  $ and the final expressions of L$_{1}$ and L$_{2}$
respectively to get:

$\left(  D_{40}\right)  $\qquad\textbf{I}$_{32}=\frac{1}{12}%
\underset{i=q+1}{\overset{n}{\sum}}\frac{\partial^{2}\Theta}{\partial
x_{i}^{2}}(y_{0})\qquad$\ 

$=-\frac{1}{3456}[3<H,i>^{2}\ +2(\tau^{M}-3\tau^{P}%
\ +\overset{q}{\underset{\text{a=1}}{\sum}}\varrho_{\text{aa}}^{M}%
+\overset{q}{\underset{\text{a,b}=1}{\sum}}R_{\text{abab}}^{M})]^{2}%
(y_{0})\phi(y_{0})\phi(y_{0})$\ \qquad I$_{321}$

\qquad$+\frac{1}{24}[2<H,i>^{2}(y_{0})+\frac{1}{3}(\tau^{M}-3\tau
^{P}+\overset{q}{\underset{\text{a}=1}{\sum}}\varrho_{\text{aa}}%
+\overset{q}{\underset{\text{a,b}=1}{\sum}}R_{\text{abab}})](y_{0})\phi
(y_{0})\qquad$I$_{3212}=\frac{1}{24}(L_{1}+L_{2}+L_{3})$

$\qquad\times\lbrack\frac{1}{4}<H,j>^{2}(y_{0})+\frac{1}{6}(\tau^{M}-3\tau
^{P}+\overset{q}{\underset{\text{a}=1}{\sum}}\varrho_{\text{aa}}%
^{M}+\overset{q}{\underset{\text{a,b}=1}{\sum}}R_{\text{abab}}^{M}%
)](y_{0})\phi(y_{0})$

$-\frac{1}{96}[<H,i><H,j>](y_{0})$

$\times\lbrack2\varrho_{ij}+$ $\overset{q}{\underset{\text{a}=1}{4\sum}%
}R_{i\text{a}j\text{a}}-3\overset{q}{\underset{\text{a,b=1}}{\sum}%
}(T_{\text{aa}i}T_{\text{bb}j}-T_{\text{ab}i}T_{\text{ab}j}%
)-3\overset{q}{\underset{\text{a,b=1}}{\sum}}(T_{\text{aa}j}T_{\text{bb}%
i}-T_{\text{ab}j}T_{\text{ab}i}](y_{0})\phi(y_{0})\qquad L_{2}\qquad
L_{21}\qquad\ \ \ \ \ $

$\qquad-\frac{1}{864}[2\varrho_{ij}+$ $\overset{q}{\underset{\text{a}%
=1}{4\sum}}R_{i\text{a}j\text{a}}-3\overset{q}{\underset{\text{a,b=1}}{\sum}%
}(T_{\text{aa}i}T_{\text{bb}j}-T_{\text{ab}i}T_{\text{ab}j}%
)-3\overset{q}{\underset{\text{a,b=1}}{\sum}}(T_{\text{aa}j}T_{\text{bb}%
i}-T_{\text{ab}j}T_{\text{ab}i}]^{2}(y_{0})\phi(y_{0})$

\qquad$-\frac{1}{288}[<H,j>](y_{0})\times\lbrack\{\nabla_{i}\varrho
_{ij}-2\varrho_{ij}<H,i>+\overset{q}{\underset{\text{a}=1}{\sum}}(\nabla
_{i}R_{\text{a}i\text{a}j}-4R_{i\text{a}j\text{a}}<H,i>)\qquad L_{212}$

$\qquad+4\overset{q}{\underset{\text{a,b=1}}{\sum}}R_{i\text{a}j\text{b}%
}T_{\text{ab}i}+2\overset{q}{\underset{\text{a,b,c=1}}{\sum}}(T_{\text{aa}%
i}T_{\text{bb}j}T_{\text{cc}i}-3T_{\text{aa}i}T_{\text{bc}j}T_{\text{bc}%
i}+2T_{\text{ab}i}T_{\text{bc}j}T_{\text{ac}i})](y_{0})\phi(y_{0})$%
\qquad\qquad\qquad\qquad\qquad\ \ 

$\qquad-\frac{1}{288}[<H,j>](y_{0})\times\lbrack\nabla_{j}\varrho
_{ii}-2\varrho_{ij}<H,i>+\overset{q}{\underset{\text{a}=1}{\sum}}(\nabla
_{j}R_{\text{a}i\text{a}i}-4R_{i\text{a}j\text{a}}<H,i>)$

$\qquad+4\overset{q}{\underset{\text{a,b=1}}{\sum}}R_{j\text{a}i\text{b}%
}T_{\text{ab}i}+2\overset{q}{\underset{\text{a,b,c=1}}{\sum}}(T_{\text{aa}%
j}T_{\text{bb}i}T_{\text{cc}i}-3T_{\text{aa}j}T_{\text{bc}i}T_{\text{bc}%
i}+2T_{\text{ab}j}T_{\text{bc}i}T_{\text{ac}i})](y_{0})\phi(y_{0})$

$\qquad-\frac{1}{288}[<H,j>](y_{0})\times\lbrack\nabla_{i}\varrho
_{ij}-2\varrho_{ii}<H,j>+\overset{q}{\underset{\text{a}=1}{\sum}}(\nabla
_{i}R_{\text{a}i\text{a}j}-4R_{i\text{a}i\text{a}}<H,j>)$

$\qquad+4\overset{q}{\underset{\text{a,b=1}}{\sum}}R_{i\text{a}i\text{b}%
}T_{\text{ab}j}+2\overset{q}{\underset{\text{a,b,c}=1}{\sum}}(T_{\text{aa}%
i}T_{\text{bb}i}T_{\text{cc}j}-3T_{\text{aa}i}T_{\text{bc}i}T_{\text{bc}%
j}+2T_{\text{ab}i}T_{\text{bc}i}T_{\text{ac}j})](y_{0})\phi(y_{0})$

$\qquad-\frac{1}{3}[<H,j><H,k>](y_{0})R_{ijik}(y_{0})-$ $\frac{5}%
{64}[<H,i>^{2}<H,j>^{2}](y_{0})\phi(y_{0})\qquad L_{213}$

$-\frac{1}{96}<H,i><H,j>[2\varrho_{ij}+\overset{q}{\underset{\text{a}%
=1}{4\sum}}R_{i\text{a}j\text{a}}-3\overset{q}{\underset{\text{a,b=1}}{\sum}%
}(T_{\text{aa}i}T_{\text{bb}j}-T_{\text{ab}i}T_{\text{ab}j}%
)-3\overset{q}{\underset{\text{a,b=1}}{\sum}}(T_{\text{aa}j}T_{\text{bb}%
i}-T_{\text{ab}j}T_{\text{ab}i}](y_{0})\phi(y_{0})$

$\qquad-\frac{1}{96}<H,j>^{2}[\tau^{M}\ -3\tau^{P}+\ \underset{\text{a}%
=1}{\overset{\text{q}}{\sum}}\varrho_{\text{aa}}^{M}+$
$\overset{q}{\underset{\text{a},\text{b}=1}{\sum}}R_{\text{abab}}^{M}$
$](y_{0})\phi(y_{0})$

$\qquad+\frac{1}{288}<H,j>[\nabla_{i}\varrho_{ij}-2\varrho_{ij}%
<H,i>+\overset{q}{\underset{\text{a}=1}{\sum}}(\nabla_{i}R_{\text{a}%
i\text{a}j}-4R_{i\text{a}j\text{a}}<H,i>)+4\overset{q}{\underset{\text{a,b=1}%
}{\sum}}R_{i\text{a}j\text{b}}T_{\text{ab}i}$

$\qquad+2\overset{q}{\underset{\text{a,b,c=1}}{\sum}}(T_{\text{aa}%
i}T_{\text{bb}j}T_{\text{cc}i}-3T_{\text{aa}i}T_{\text{bc}j}T_{\text{bc}%
i}+2T_{\text{ab}i}T_{\text{bc}j}T_{\text{ac}i})](y_{0})\phi(y_{0})$%
\qquad\qquad\qquad\qquad\qquad\ \ 

$\qquad+\frac{1}{12}<H,j>[\nabla_{j}\varrho_{ii}-2\varrho_{ij}%
<H,i>+\overset{q}{\underset{\text{a}=1}{\sum}}(\nabla_{j}R_{\text{a}%
i\text{a}i}-4R_{i\text{a}j\text{a}}<H,i>)+4\overset{q}{\underset{\text{a,b=1}%
}{\sum}}R_{j\text{a}i\text{b}}T_{\text{ab}i}$

$\qquad+2\overset{q}{\underset{\text{a,b,c=1}}{\sum}}(T_{\text{aa}%
j}T_{\text{bb}i}T_{\text{cc}i}-3T_{\text{aa}j}T_{\text{bc}i}T_{\text{bc}%
i}+2T_{\text{ab}j}T_{\text{bc}i}T_{\text{ac}i})](y_{0})\phi(y_{0})$

$\qquad+\frac{1}{288}<H,j>[\nabla_{i}\varrho_{ij}-2\varrho_{ii}%
<H,j>+\overset{q}{\underset{\text{a}=1}{\sum}}(\nabla_{i}R_{\text{a}%
i\text{a}j}-4R_{i\text{a}i\text{a}}<H,j>)+4\overset{q}{\underset{\text{a,b=1}%
}{\sum}}R_{i\text{a}i\text{b}}T_{\text{ab}j}$

$\qquad+2\overset{q}{\underset{\text{a,b,c}=1}{\sum}}(T_{\text{aa}%
i}T_{\text{bb}i}T_{\text{cc}j}-3T_{\text{aa}i}T_{\text{bc}i}T_{\text{bc}%
j}+2T_{\text{ab}i}T_{\text{bc}i}T_{\text{ac}j})](y_{0})\phi(y_{0})$

\qquad$-\frac{1}{144}R_{jijk}(y_{0})$ \ $[<H,i><H,k>](y_{0})\phi(y_{0}%
)\qquad\qquad L_{22}$

$-\frac{1}{432}R_{jijk}(y_{0})[2\varrho_{ik}+$ $\overset{q}{\underset{\text{a}%
=1}{4\sum}}R_{i\text{a}k\text{a}}-3\overset{q}{\underset{\text{a,b=1}}{\sum}%
}(T_{\text{aa}i}T_{\text{bb}k}-T_{\text{ab}i}T_{\text{ab}k}%
)-3\overset{q}{\underset{\text{a,b=1}}{\sum}}(T_{\text{aa}k}T_{\text{bb}%
i}-T_{\text{ab}k}T_{\text{ab}i}](y_{0})\phi(y_{0})$

$+\frac{1}{144}<H,k>(y_{0})[\nabla_{j}$R$_{ijik}(y_{0})-\nabla_{i}$%
R$_{jijk}](y_{0})\phi(y_{0})$

$-\frac{5}{32}<H,i>^{2}(y_{0})<H,j>^{2}(y_{0})\phi(y_{0})\qquad\qquad
L_{23}\qquad L_{231}$

$-\frac{1}{48}<H,i>(y_{0})<H,j>(y_{0})$

$\times\lbrack2\varrho_{ij}+$ $\overset{q}{\underset{\text{a}=1}{4\sum}%
}R_{i\text{a}j\text{a}}-3\overset{q}{\underset{\text{a,b=1}}{\sum}%
}(T_{\text{aa}i}T_{\text{bb}j}-T_{\text{ab}i}T_{\text{ab}j}%
)-3\overset{q}{\underset{\text{a,b=1}}{\sum}}(T_{\text{aa}j}T_{\text{bb}%
i}-T_{\text{ab}j}T_{\text{ab}i}](y_{0})\phi(y_{0})$

$-\frac{1}{48}<H,i>^{2}(y_{0})[\varrho_{jj}+$ $\overset{q}{\underset{\text{a}%
=1}{2\sum}}R_{j\text{a}j\text{a}}-3\overset{q}{\underset{\text{a,b=1}}{\sum}%
}(T_{\text{aa}j}T_{\text{bb}j}-T_{\text{ab}j}T_{\text{ab}j})](y_{0})\phi
(y_{0})$

$-\frac{1}{144}<H,i>(y_{0})[\nabla_{i}\varrho_{jj}-2\varrho_{ij}%
<H,j>+\overset{q}{\underset{\text{a}=1}{\sum}}(\nabla_{i}R_{\text{a}%
j\text{a}j}-4R_{i\text{a}j\text{a}}<H,j>)+4\overset{q}{\underset{\text{a,b=1}%
}{\sum}}R_{i\text{a}j\text{b}}T_{\text{ab}j}$

$+2\overset{q}{\underset{\text{a,b,c=1}}{\sum}}(T_{\text{aa}i}T_{\text{bb}%
j}T_{\text{cc}j}-3T_{\text{aa}i}T_{\text{bc}j}T_{\text{bc}j}+2T_{\text{ab}%
i}T_{\text{bc}j}T_{\text{ca}j})](y_{0})\phi(y_{0})$\qquad\qquad\qquad
\qquad\qquad\ \ 

$-\frac{1}{24}\times\frac{1}{6}<H,i>(y_{0})[\nabla_{j}\varrho_{ij}%
-2\varrho_{ij}<H,j>+\overset{q}{\underset{\text{a}=1}{\sum}}(\nabla
_{j}R_{\text{a}i\text{a}j}-4R_{j\text{a}i\text{a}}%
<H,j>)+4\overset{q}{\underset{\text{a,b=1}}{\sum}}R_{j\text{a}i\text{b}%
}T_{\text{ab}j}$

$+2\overset{q}{\underset{\text{a,b,c=1}}{\sum}}(T_{\text{aa}j}T_{\text{bb}%
i}T_{\text{cc}j}-3T_{\text{aa}j}T_{\text{bc}i}T_{\text{bc}j}+2T_{\text{ab}%
j}T_{\text{bc}i}T_{\text{ac}j})](y_{0})$

$-\frac{1}{144}<H,i>(y_{0})[\nabla_{j}\varrho_{ij}-2\varrho_{jj}%
<H,i>+\overset{q}{\underset{\text{a}=1}{\sum}}(\nabla_{j}R_{\text{a}%
i\text{a}j}-4R_{j\text{a}j\text{a}}<H,i>)+4\overset{q}{\underset{\text{a,b=1}%
}{\sum}}R_{j\text{a}j\text{b}}T_{\text{ab}i}$

$\qquad+2\overset{q}{\underset{\text{a,b,c=1}}{\sum}}(T_{\text{aa}%
j}T_{\text{bb}j}T_{\text{cc}i}-3T_{\text{aa}j}T_{\text{bc}j}T_{\text{bc}%
i}+2T_{\text{ab}j}T_{\text{bc}j}T_{\text{ac}i})](y_{0})\phi(y_{0})$

$\qquad-\frac{1}{96}<H,j>^{2}(y_{0})[\varrho_{ii}+$
$\overset{q}{\underset{\text{a}=1}{2\sum}}R_{i\text{a}i\text{a}}%
-3\overset{q}{\underset{\text{a,b=1}}{\sum}}(T_{\text{aa}i}T_{\text{bb}%
i}-T_{\text{ab}i}T_{\text{ab}i})](y_{0})\qquad L_{232}$

$\qquad-\frac{1}{432}[\varrho_{ii}+$ $\overset{q}{\underset{\text{a}=1}{2\sum
}}R_{i\text{a}i\text{a}}-3\overset{q}{\underset{\text{a,b=1}}{\sum}%
}(T_{\text{aa}i}T_{\text{bb}i}-T_{\text{ab}i}T_{\text{ab}i})](y_{0})\phi
(y_{0})$

\qquad$\times\lbrack\varrho_{jj}+$ $\overset{q}{\underset{\text{a}=1}{2\sum}%
}R_{j\text{a}j\text{a}}-3\overset{q}{\underset{\text{a,b=1}}{\sum}%
}(T_{\text{aa}j}T_{\text{bb}j}-T_{\text{ab}j}T_{\text{ab}j})]\}(y_{0}%
)\phi(y_{0})$

$\qquad+\frac{1}{48}R_{ijik}(y_{0})$ \ $[<H,j><H,k>](y_{0})\phi(y_{0}%
)\qquad\qquad\qquad\qquad$\ $L_{233}$

$+\frac{1}{432}R_{ijik}(y_{0})[2\varrho_{jk}+$ $\overset{q}{\underset{\text{a}%
=1}{4\sum}}R_{j\text{a}k\text{a}}-3\overset{q}{\underset{\text{a,b=1}}{\sum}%
}(T_{\text{aa}j}T_{\text{bb}k}-T_{\text{ab}j}T_{\text{ab}k}%
)-3\overset{q}{\underset{\text{a,b=1}}{\sum}}(T_{\text{aa}k}T_{\text{bb}%
j}-T_{\text{ab}k}T_{\text{ab}j}](y_{0})\phi(y_{0})$

$+\overset{n}{\underset{i,j=q+1}{\sum}}\frac{35}{128}<H,i>^{2}(y_{0}%
)<H,j>^{2}(y_{0})\qquad\qquad\ \frac{1}{24}\frac{\partial^{4}\theta^{-\frac
{1}{2}}}{\partial x_{i}^{2}\partial x_{j}^{2}}(y_{0})$

$+\frac{5}{192}\overset{n}{\underset{j=q+1}{\sum}}<H,j>^{2}(y_{0})[\tau
^{M}\ -3\tau^{P}+\ \underset{\text{a}=1}{\overset{\text{q}}{\sum}}%
\varrho_{\text{aa}}^{M}+\overset{q}{\underset{\text{a},\text{b}=1}{\sum}%
}R_{\text{abab}}^{M}](y_{0})\qquad\ \ \ \ \ \ \ \ $

$+\frac{5}{192}\overset{n}{\underset{i=q+1}{\sum}}<H,i>^{2}(y_{0})[\tau
^{M}\ -3\tau^{P}+\ \underset{\text{a}=1}{\overset{\text{q}}{\sum}}%
\varrho_{\text{aa}}^{M}+\overset{q}{\underset{\text{a},\text{b}=1}{\sum}%
}R_{\text{abab}}^{M}](y_{0})\qquad\qquad$

$+\frac{5}{192}\overset{n}{\underset{i,j=q+1}{\sum}}[<H,i><H,j>](y_{0}%
)\qquad\qquad\qquad\qquad\qquad\qquad\qquad\qquad$

$\times\lbrack2\varrho_{ij}+4\overset{q}{\underset{\text{a}=1}{\sum}%
}R_{i\text{a}j\text{a}}-3\overset{q}{\underset{\text{a,b=1}}{\sum}%
}(T_{\text{aa}i}T_{\text{bb}j}-T_{\text{ab}i}T_{\text{ab}j}%
)-3\overset{q}{\underset{\text{a,b=1}}{\sum}}(T_{\text{aa}j}T_{\text{bb}%
i}-T_{\text{ab}j}T_{\text{ab}i})](y_{0})$

$+\frac{1}{96}\overset{n}{\underset{i,j=q+1}{\sum}}<H,j>(y_{0})[\{\nabla
_{i}\varrho_{ij}-2\varrho_{ij}<H,i>+\overset{q}{\underset{\text{a}=1}{\sum}%
}(\nabla_{i}R_{\text{a}i\text{a}j}-4R_{i\text{a}j\text{a}}<H,i>)\qquad$

$+4\overset{q}{\underset{\text{a,b=1}}{\sum}}R_{i\text{a}j\text{b}%
}T_{\text{ab}i}+2\overset{q}{\underset{\text{a,b,c=1}}{\sum}}(T_{\text{aa}%
i}T_{\text{bb}j}T_{\text{cc}i}-T_{\text{aa}i}T_{\text{bc}j}T_{\text{bc}%
i}-2T_{\text{bc}j}(T_{\text{aa}i}T_{\text{bc}i}-T_{\text{ab}i}T_{\text{ac}%
i}))\}$\qquad\qquad\qquad\ \ 

$+\{\nabla_{j}\varrho_{ii}-2\varrho_{ij}<H,i>+\overset{q}{\underset{\text{a}%
=1}{\sum}}(\nabla_{j}R_{\text{a}i\text{a}i}-4R_{i\text{a}j\text{a}}<H,i>)$

$+4\overset{q}{\underset{\text{a,b=1}}{\sum}}R_{j\text{a}i\text{b}%
}T_{\text{ab}i}+2\overset{q}{\underset{\text{a,b,c=1}}{\sum}}(T_{\text{aa}%
j}(T_{\text{bb}i}T_{\text{cc}i}-T_{\text{bc}i}T_{\text{bc}i})-2T_{\text{aa}%
j}T_{\text{bc}i}T_{\text{bc}i}+2T_{\text{ab}j}T_{\text{bc}i}T_{\text{ac}%
i})\}\qquad$

$+\{\nabla_{i}\varrho_{ij}-2\varrho_{ii}<H,j>+\overset{q}{\underset{\text{a}%
=1}{\sum}}(\nabla_{i}R_{\text{a}i\text{a}j}-4R_{i\text{a}i\text{a}%
}<H,j>)+4\overset{q}{\underset{\text{a,b=1}}{\sum}}R_{i\text{a}i\text{b}%
}T_{\text{ab}j}$

$+2\overset{q}{\underset{\text{a,b,c}=1}{\sum}}(T_{\text{aa}i}T_{\text{bb}%
i}T_{\text{cc}j}-3T_{\text{aa}i}T_{\text{bc}i}T_{\text{bc}j}+2T_{\text{ab}%
i}T_{\text{bc}i}T_{\text{ac}j})\}](y_{0})$

$+\frac{1}{96}\overset{n}{\underset{i,j=q+1}{\sum}}<H,i>(y_{0})[\{\nabla
_{i}\varrho_{jj}-2\varrho_{ij}<H,j>+\overset{q}{\underset{\text{a}=1}{\sum}%
}(\nabla_{i}R_{\text{a}j\text{a}j}-4R_{i\text{a}j\text{a}}<H,j>)\qquad$

$+4\overset{q}{\underset{\text{a,b=1}}{\sum}}R_{i\text{a}j\text{b}%
}T_{\text{ab}j}+2\overset{q}{\underset{\text{a,b,c=1}}{\sum}}T_{\text{aa}%
i}(T_{\text{bb}j}T_{\text{cc}j}-T_{\text{bc}j}T_{\text{bc}j})-2T_{\text{aa}%
i}T_{\text{bc}j}T_{\text{bc}j}+2T_{\text{ab}i}T_{\text{bc}j}T_{\text{ac}%
j})\}(y_{0})\qquad$\qquad\qquad\qquad\qquad\qquad\ \ 

$+\{\nabla_{j}\varrho_{ij}-2\varrho_{ij}<H,j>+\overset{q}{\underset{\text{a}%
=1}{\sum}}(\nabla_{j}R_{\text{a}i\text{a}j}-4R_{j\text{a}i\text{a}}<H,j>)$

$+4\overset{q}{\underset{\text{a,b=1}}{\sum}}R_{j\text{a}i\text{b}%
}T_{\text{ab}j}+2\overset{q}{\underset{\text{a,b,c=1}}{\sum}}(T_{\text{aa}%
j}T_{\text{bb}i}T_{\text{cc}j}-T_{\text{ab}j}T_{\text{bc}i}T_{\text{ac}%
j}-2T_{\text{bc}i}(T_{\text{aa}j}T_{\text{bc}j}-T_{\text{ab}j}T_{\text{ac}%
j}))\}(y_{0})$

$+\{\nabla_{j}\varrho_{ij}-2\varrho_{jj}<H,i>+\overset{q}{\underset{\text{a}%
=1}{\sum}}(\nabla_{j}R_{\text{a}i\text{a}j}-4R_{j\text{a}j\text{a}%
}<H,i>)+4\overset{q}{\underset{\text{a,b=1}}{\sum}}R_{j\text{a}j\text{b}%
}T_{\text{ab}i}$

$+2\overset{q}{\underset{\text{a,b,c=1}}{\sum}}(T_{\text{aa}j}T_{\text{bb}%
j}T_{\text{cc}i}-3T_{\text{aa}j}T_{\text{bc}j}T_{\text{bc}i}+2T_{\text{ab}%
j}T_{\text{bc}j}T_{\text{ac}i})\}](y_{0})$

$+\frac{1}{576}\overset{n}{\underset{i,j=q+1}{\sum}}[2\varrho_{ij}%
+4\overset{q}{\underset{\text{a}=1}{\sum}}R_{i\text{a}j\text{a}}%
-3\overset{q}{\underset{\text{a,b=1}}{\sum}}(T_{\text{aa}i}T_{\text{bb}%
j}-T_{\text{ab}i}T_{\text{ab}j})-3\overset{q}{\underset{\text{a,b=1}}{\sum}%
}(T_{\text{aa}j}T_{\text{bb}i}-T_{\text{ab}j}T_{\text{ab}i})]^{2}(y_{0})$

$+\frac{1}{288}[\tau^{M}\ -3\tau^{P}+\ \underset{\text{a}=1}{\overset{\text{q}%
}{\sum}}\varrho_{\text{aa}}^{M}+\overset{q}{\underset{\text{a},\text{b}%
=1}{\sum}}R_{\text{abab}}^{M}]^{2}(y_{0})$

$-\ \frac{1}{288}\overset{n}{\underset{i,j=q+1}{\sum}}[$
$\overset{q}{\underset{\text{a=1}}{\sum}}\{-(\nabla_{ii}^{2}R_{j\text{a}%
j\text{a}}+\nabla_{jj}^{2}R_{i\text{a}i\text{a}}+4\nabla_{ij}^{2}%
R_{i\text{a}j\text{a}}+2R_{ij}R_{i\text{a}j\text{a}})\qquad A$

$+\overset{n}{\underset{p=q+1}{\sum}}\overset{q}{\underset{\text{a=1}}{\sum}%
}(R_{\text{a}iip}R_{\text{a}jjp}+R_{\text{a}jjp}R_{\text{a}iip}+R_{\text{a}%
ijp}R_{\text{a}ijp}+R_{\text{a}ijp}R_{\text{a}jip}+R_{\text{a}jip}%
R_{\text{a}ijp}+R_{\text{a}jip}R_{\text{a}jip})$

$+2\overset{q}{\underset{\text{a,b=1}}{\sum}}\nabla_{i}(R)_{\text{a}%
i\text{b}j}T_{\text{ab}j}+2\overset{q}{\underset{\text{a,b=1}}{\sum}}%
\nabla_{j}(R)_{\text{a}j\text{b}i}T_{\text{ab}i}%
+2\overset{q}{\underset{\text{a,b=1}}{\sum}}\nabla_{i}(R)_{\text{a}j\text{b}%
i}T_{\text{ab}j}+2\overset{q}{\underset{\text{a,b=1}}{\sum}}\nabla
_{i}(R)_{\text{a}j\text{b}j}T_{\text{ab}i}$

$+2\overset{q}{\underset{\text{a,b=1}}{\sum}}\nabla_{j}(R)_{\text{a}%
i\text{b}i}T_{\text{ab}j}+2\overset{q}{\underset{\text{a,b=1}}{\sum}}%
\nabla_{j}(R)_{\text{a}i\text{b}j}T_{\text{ab}i}$

$+\overset{n}{\underset{p=q+1}{\sum}}(-\frac{3}{5}\nabla_{ii}^{2}%
(R)_{jpjp}+\overset{n}{\underset{p=q+1}{\sum}}(-\frac{3}{5}\nabla_{jj}%
^{2}(R)_{ipip}+\overset{n}{\underset{p=q+1}{\sum}}(-\frac{3}{5}\nabla_{ij}%
^{2}(R)_{ipjp}+\overset{n}{\underset{p=q+1}{\sum}}(-\frac{3}{5}\nabla_{ij}%
^{2}(R)_{jpip}$

$+\overset{n}{\underset{p=q+1}{\sum}}(-\frac{3}{5}\nabla_{ji}^{2}%
(R)_{ipjp}+\overset{n}{\underset{p=q+1}{\sum}}(-\frac{3}{5}\nabla_{ji}%
^{2}(R)_{jpip}$

$+\frac{1}{5}\overset{n}{\underset{m,p=q+1}{%
{\textstyle\sum}
}}R_{ipim}R_{jpjm}+\frac{1}{5}\overset{n}{\underset{m,p=q+1}{%
{\textstyle\sum}
}}R_{jpjm}R_{ipim}+\frac{1}{5}\overset{n}{\underset{m,p=q+1}{%
{\textstyle\sum}
}}R_{ipjm}R_{ipjm}+\frac{1}{5}\overset{n}{\underset{m,p=q+1}{%
{\textstyle\sum}
}}R_{ipjm}R_{jpim}$

$+\frac{1}{5}\overset{n}{\underset{m,p=q+1}{%
{\textstyle\sum}
}}R_{jpim}R_{ipjm}+\frac{1}{5}\overset{n}{\underset{m,p=q+1}{%
{\textstyle\sum}
}}R_{jpim}R_{jpim}\}(y_{0})$

$+4\overset{q}{\underset{\text{a,b=1}}{\sum}}\{(\nabla_{i}(R)_{i\text{a}%
j\text{a}}-\overset{q}{\underset{\text{c=1}}{%
{\textstyle\sum}
}}R_{\text{a}i\text{c}i}T_{\text{ac}j})$ $T_{\text{bb}j}+4(\nabla
_{j}(R)_{j\text{a}i\text{a}}-\overset{q}{\underset{\text{c=1}}{%
{\textstyle\sum}
}}R_{\text{a}j\text{c}j}T_{\text{ac}i})$ $T_{\text{bb}i}+4(\nabla
_{i}(R)_{j\text{a}i\text{a}}-\overset{q}{\underset{\text{c=1}}{%
{\textstyle\sum}
}}R_{\text{a}i\text{c}j}T_{\text{ac}i})$ $T_{\text{bb}j}$ $4B\ $

$+4(\nabla_{i}(R)_{j\text{a}j\text{a}}-\overset{q}{\underset{\text{c=1}}{%
{\textstyle\sum}
}}R_{\text{a}i\text{c}j}T_{\text{ac}j})$ $T_{\text{bb}i}+4(\nabla
_{j}(R)_{i\text{a}i\text{a}}-\overset{q}{\underset{\text{c=1}}{%
{\textstyle\sum}
}}R_{\text{a}j\text{c}i}T_{\text{ac}i})$ $T_{\text{bb}j}+4(\nabla
_{j}(R)_{i\text{a}j\text{a}}-\overset{q}{\underset{\text{c=1}}{%
{\textstyle\sum}
}}R_{\text{a}j\text{c}i}T_{\text{ac}j})$ $T_{\text{bb}i}$

$-4\overset{q}{\underset{\text{a,b=1}}{\sum}}(\nabla_{i}(R)_{i\text{a}%
j\text{b}}-\overset{q}{\underset{\text{c=1}}{%
{\textstyle\sum}
}}R_{\text{b}r\text{c}s}T_{\text{ac}t})T_{\text{ab}j}%
-4\overset{q}{\underset{\text{a,b=1}}{\sum}}(\nabla_{j}(R)_{j\text{a}%
i\text{b}}-\overset{q}{\underset{\text{c=1}}{%
{\textstyle\sum}
}}R_{\text{b}j\text{c}j}T_{\text{ac}i})T_{\text{ab}i}$

$-4\overset{q}{\underset{\text{a,b=1}}{\sum}}(\nabla_{i}(R)_{j\text{a}%
i\text{b}}-\overset{q}{\underset{\text{c=1}}{%
{\textstyle\sum}
}}R_{\text{b}i\text{c}j}T_{\text{ac}i})T_{\text{ab}j}%
-4\overset{q}{\underset{\text{a,b=1}}{\sum}}(\nabla_{i}(R)_{j\text{a}%
j\text{b}}-\overset{q}{\underset{\text{c=1}}{%
{\textstyle\sum}
}}R_{\text{b}i\text{c}j}T_{\text{ac}j})T_{\text{ab}i}$

$-4\overset{q}{\underset{\text{a,b=1}}{\sum}}(\nabla_{j}(R)_{i\text{a}%
i\text{b}}-\overset{q}{\underset{\text{c=1}}{%
{\textstyle\sum}
}}R_{\text{b}j\text{c}i}T_{\text{ac}i})T_{\text{ab}j}%
-4\overset{q}{\underset{\text{a,b=1}}{\sum}}(\nabla_{j}(R)_{i\text{a}%
j\text{b}}-\overset{q}{\underset{\text{c=1}}{%
{\textstyle\sum}
}}R_{\text{b}j\text{c}i}T_{\text{ac}j})T_{\text{ab}i}\}](y_{0})$

$-\frac{1}{48}$ $[\frac{4}{9}\overset{q}{\underset{\text{a,b=1}}{\sum}%
}(\varrho_{\text{aa}}-\overset{q}{\underset{\text{c}=1}{\sum}}R_{\text{acac}%
})(\varrho_{\text{bb}}-\overset{q}{\underset{\text{d}=1}{\sum}}R_{\text{bdbd}%
})+\frac{8}{9}\overset{n}{\underset{i,j=q+1}{\sum}}%
\overset{q}{\underset{\text{a,b}=1}{\sum}}(R_{i\text{a}j\text{a}}%
R_{i\text{b}j\text{b}})\qquad3C$

$+\frac{2}{9}\overset{q}{\underset{\text{a}=1}{\sum}}(\varrho_{\text{aa}}%
^{M}-\varrho_{\text{aa}}^{P})(\tau^{M}-\overset{q}{\underset{\text{c}=1}{\sum
}}\varrho_{\text{cc}}^{M})+\frac{4}{9}\overset{n}{\underset{i,j=q+1}{\sum}%
}\overset{q}{\underset{\text{a}=1}{\sum}}R_{i\text{a}j\text{a}}\varrho_{ij}\ $

$\ +\frac{2}{9}\overset{q}{\underset{\text{b}=1}{\sum}}(\varrho_{\text{bb}%
}^{M}-\varrho_{\text{bb}}^{P})(\tau^{M}-\overset{q}{\underset{\text{c}%
=1}{\sum}}\varrho_{\text{cc}}^{M})+\frac{4}{9}%
\overset{n}{\underset{i,j=q+1}{\sum}}\overset{q}{\underset{\text{b}=1}{\sum}%
}R_{i\text{b}j\text{b}}\varrho_{ij}\ $

$+\frac{1}{9}(\tau^{M}-\overset{q}{\underset{\text{a=1}}{\sum}}\varrho
_{\text{aa}})(\tau^{M}-\overset{q}{\underset{\text{b=1}}{\sum}}\varrho
_{\text{bb}})+\frac{2}{9}(\left\Vert \varrho^{M}\right\Vert ^{2}%
-\overset{q}{\underset{\text{a,b}=1}{\sum}}\varrho_{\text{ab}})$

$-\overset{n}{\underset{i,j=q+1}{\sum}}\overset{q}{\underset{\text{a,b}%
=1}{\sum}}R_{i\text{a}i\text{b}}R_{j\text{a}j\text{b}}\ -\frac{1}%
{2}\overset{n}{\underset{i,j=q+1}{\sum}}\overset{q}{\underset{\text{a,b}%
=1}{\sum}}R_{i\text{a}j\text{b}}^{2}-\overset{n}{\underset{i,j=q+1}{\sum}%
}\overset{q}{\underset{\text{a,b}=1}{\sum}}R_{i\text{a}j\text{b}}%
R_{j\text{a}i\text{b}}-\frac{1}{2}\overset{n}{\underset{i,j=q+1}{\sum}%
}\overset{q}{\underset{\text{a,b}=1}{\sum}}R_{j\text{a}i\text{b}}^{2}$

$-\frac{1}{9}\overset{n}{\underset{i,j,p,m=q+1}{\sum}}R_{ipim}R_{jpjm}%
\ -\frac{1}{18}\overset{n}{\underset{i,j,p,m=q+1}{\sum}}R_{ipjm}^{2}-\frac
{1}{9}\overset{n}{\underset{i,j,p,m=q+1}{\sum}}R_{ipjm}R_{jpim}-\frac{1}%
{18}\overset{n}{\underset{i,j,p,m=q+1}{\sum}}R_{jpim}^{2}$

$-\frac{1}{3}\overset{q}{\underset{\text{a}=1}{\sum}}%
\overset{n}{\underset{i,j,p=q+1}{\sum}}R_{i\text{a}ip}R_{j\text{a}jp}-\frac
{1}{6}\overset{q}{\underset{\text{a}=1}{\sum}}%
\overset{n}{\underset{i,j,p=q+1}{\sum}}R_{i\text{a}jp}^{2}-\frac{1}%
{3}\overset{q}{\underset{\text{a}=1i,j,}{\sum}}%
\overset{n}{\underset{p=q+1}{\sum}}R_{i\text{a}jp}R_{j\text{a}ip}-\frac{1}%
{6}\overset{q}{\underset{\text{a}=1}{\sum}}%
\overset{n}{\underset{i,j,p=q+1}{\sum}}R_{j\text{a}ip}^{2}$

$-\frac{1}{3}\overset{q}{\underset{\text{b}=1i,j,}{\sum}}%
\overset{n}{\underset{p=q+1}{\sum}}R_{i\text{b}ip}R_{j\text{b}jp}-\frac{1}%
{6}\overset{q}{\underset{\text{b}=1}{\sum}}%
\overset{n}{\underset{i,j,p=q+1}{\sum}}R_{i\text{b}jp}^{2}-\frac{1}%
{3}\overset{q}{\underset{\text{b}=1}{\sum}}%
\overset{n}{\underset{i.j,p=q+1}{\sum}}R_{i\text{b}jp}R_{j\text{b}ip}-\frac
{1}{6}\overset{q}{\underset{\text{b}=1}{\sum}}%
\overset{n}{\underset{i,j,p=q+1}{\sum}}R_{j\text{b}ip}^{2}](y_{0})$

$-\frac{1}{48}$ $\overset{q}{\underset{\text{a,b,c=1}}{\sum}}[$
$-\overset{n}{\underset{i=q+1}{\sum}}R_{i\text{a}i\text{a}}(R_{\text{bcbc}%
}^{P}-R_{\text{bcbc}}^{M})$ $-\overset{n}{\underset{j=q+1}{\sum}}%
R_{j\text{a}j\text{a}}(R_{\text{bcbc}}^{P}-R_{\text{bcbc}}^{M})\qquad\qquad6D$

\ $+\overset{n}{\underset{i=q+1}{\sum}}R_{i\text{a}i\text{b}}(R_{\text{acbc}%
}^{P}-R_{\text{acbc}}^{M})\ -\overset{n}{\underset{i=q+1}{\sum}}%
R_{i\text{a}i\text{c}}(R_{\text{abbc}}^{P}-R_{\text{abbc}}^{M})$

$+\overset{n}{\underset{j=q+1}{\sum}}R_{j\text{a}j\text{b}}(R_{\text{acbc}%
}^{P}-R_{\text{acbc}}^{M})$\ $-\overset{n}{\underset{j=q+1}{\sum}}%
R_{j\text{a}j\text{c}}(R_{\text{abbc}}^{P}-R_{\text{abbc}}^{M})$

$+\underset{i,j=q+1}{\overset{n}{\sum}}$ $-R_{i\text{a}j\text{a}}%
(T_{\text{bb}i}T_{\text{cc}j}$ $-T_{\text{bc}i}T_{\text{bc}j})$
$-\underset{i,j=q+1}{\overset{n}{\sum}}R_{i\text{a}j\text{a}}(T_{\text{bb}%
j}T_{\text{cc}i}$ $-T_{\text{bc}j}T_{\text{bc}i})$

$+$ $\underset{i,j=q+1}{\overset{n}{\sum}}$ $-R_{j\text{a}i\text{a}%
}(T_{\text{bb}i}T_{\text{cc}j}$ $-T_{\text{bc}i}T_{\text{bc}j})$
$-\underset{i,j=q+1}{\overset{n}{\sum}}R_{j\text{a}i\text{a}}(T_{\text{bb}%
j}T_{\text{cc}i}$ $-T_{\text{bc}j}T_{\text{bc}i})$

$+\underset{i,j=q+1}{\overset{n}{\sum}}\ R_{i\text{a}j\text{b}}(T_{\text{ab}%
i}T_{\text{cc}j}-T_{\text{bc}i}T_{\text{ac}j}%
)\ +\underset{i,j=q+1}{\overset{n}{\sum}}\ R_{i\text{a}j\text{b}}%
(T_{\text{ab}j}T_{\text{cc}i}-T_{\text{bc}j}T_{\text{ac}i})$

$+\underset{i,j=q+1}{\overset{n}{\sum}}\ R_{j\text{a}i\text{ib}}%
(T_{\text{ab}i}T_{\text{cc}j}-T_{\text{bc}i}T_{\text{ac}j}%
)\ +\underset{i,j=q+1}{\overset{n}{\sum}}\ R_{j\text{a}i\text{b}}%
(T_{\text{ab}j}T_{\text{cc}i}-T_{\text{bc}j}T_{\text{ac}i})\qquad$

$+\underset{i,j=q+1}{\overset{n}{\sum}}-R_{i\text{a}j\text{c}}(T_{\text{ab}%
i}T_{\text{bc}j}-T_{\text{ac}i}T_{\text{bb}j}%
)-\underset{i,j=q+1}{\overset{n}{\sum}}R_{i\text{a}j\text{c}}(T_{\text{ba}%
j}T_{\text{bc}i}-T_{\text{ac}j}T_{\text{bb}i})$

$+\underset{i,j=q+1}{\overset{n}{\sum}}-R_{j\text{a}i\text{c}}(T_{\text{ba}%
i}T_{\text{bc}j}-T_{\text{ac}i}T_{\text{bb}j}%
)-\underset{i,j=q+1}{\overset{n}{\sum}}R_{j\text{a}i\text{c}}(T_{\text{ba}%
j}T_{\text{bc}i}-T_{\text{ac}j}T_{\text{bb}i})](y_{0})$

$+\frac{1}{144}\underset{p=q+1}{\overset{n}{\sum}}%
[\underset{i=q+1}{\overset{n}{\sum}}\overset{q}{\underset{\text{b,c=1}}{\sum}%
}R_{ipip}(R_{\text{bcbc}}^{P}-R_{\text{bcbc}}^{M}%
)+\underset{j=q+1}{\overset{n}{\sum}}$ $\overset{q}{\underset{\text{b,c=1}%
}{\sum}}R_{jpjp}(R_{\text{bcbc}}^{P}-R_{\text{bcbc}}^{M})](y_{0})$

$+\frac{1}{72}\underset{i,j,p=q+1}{\overset{n}{\sum}}%
\overset{q}{\underset{\text{b,c=1}}{\sum}}[R_{ipjp}(T_{\text{bb}i}%
T_{\text{cc}j}-T_{\text{bc}i}T_{\text{bc}j})+R_{ipjp}(T_{\text{bb}%
j}T_{\text{cc}i}-T_{\text{bc}j}T_{\text{bc}i})](y_{0})\qquad$

$-\frac{1}{288}\underset{i,j=q+1}{\overset{n}{\sum}}[T_{\text{aa}%
i}T_{\text{bb}j}(T_{\text{cc}i}T_{\text{dd}j}-T_{\text{cd}i}T_{\text{dc}%
j})+T_{\text{aa}i}T_{\text{bb}j}(T_{\text{cc}j}T_{\text{dd}i}-T_{\text{cd}%
j}T_{\text{dc}i})\qquad E$

$+T_{\text{aa}j}T_{\text{bb}i}(T_{\text{cc}i}T_{\text{dd}j}-T_{\text{cd}%
i}T_{\text{dc}j})+T_{\text{aa}j}T_{\text{bb}i}(T_{\text{cc}j}T_{\text{dd}%
i}-T_{\text{cd}j}T_{\text{dc}i})](y_{0})$

$+\frac{1}{288}\underset{i,j=q+1}{\overset{n}{\sum}}[T_{\text{aa}%
i}T_{\text{bc}j}(T_{\text{bc}i}T_{\text{dd}j}-T_{\text{bd}i}T_{\text{cd}%
j})+T_{\text{aa}i}T_{\text{bc}j}(T_{\text{bc}j}T_{\text{dd}i}-T_{\text{bd}%
j}T_{\text{cd}i})$

$+T_{\text{aa}j}T_{\text{bc}i}(T_{\text{bc}i}T_{\text{dd}j}-T_{\text{bd}%
i}T_{\text{cd}j})+T_{\text{aa}j}T_{\text{bc}i}(T_{\text{bc}j}T_{\text{dd}%
i}-T_{\text{bd}j}T_{\text{cd}i})](y_{0})$

$-\frac{1}{288}\underset{i,j=q+1}{\overset{n}{\sum}}[T_{\text{aa}%
i}T_{\text{bd}j}(T_{\text{bc}i}T_{\text{cd}j}-T_{\text{bd}i}T_{\text{cc}%
j})+T_{\text{aa}i}T_{\text{bd}j}(T_{\text{bc}j}T_{\text{cd}i}-T_{\text{bd}%
j}T_{\text{cc}i})$

$+T_{\text{aa}j}T_{\text{bd}i}(T_{\text{bc}i}T_{\text{cd}j}-T_{\text{bd}%
i}T_{\text{cc}j})+T_{\text{aa}j}T_{\text{bd}i}(T_{\text{bc}j}T_{\text{cd}%
i}-T_{\text{bd}j}T_{\text{cc}i})](y_{0})\qquad$

$+\frac{1}{288}\underset{i,j=q+1}{\overset{n}{\sum}}[T_{\text{ab}%
i}T_{\text{ab}j}(T_{\text{cc}i}T_{\text{dd}j}-T_{\text{cd}i}T_{\text{dc}%
j})+T_{\text{ab}i}T_{\text{ab}j}(T_{\text{cc}j}T_{\text{dd}i}-T_{\text{cd}%
j}T_{\text{dc}i})$

$+T_{\text{ab}j}T_{\text{ab}i}(T_{\text{cc}i}T_{\text{dd}j}-T_{\text{cd}%
i}T_{\text{dc}j})+T_{\text{ab}j}T_{\text{ab}i}(T_{\text{cc}j}T_{\text{dd}%
i}-T_{\text{cd}j}T_{\text{dc}i})](y_{0})$

$-\frac{1}{288}\underset{i,j=q+1}{\overset{n}{\sum}}[T_{\text{ab}%
i}T_{\text{bc}j}(T_{\text{ac}i}T_{\text{dd}j}-T_{\text{ad}i}T_{\text{cd}%
j})+T_{\text{ab}i}T_{\text{bc}j}(T_{\text{ac}j}T_{\text{dd}i}-T_{\text{ad}%
j}T_{\text{cd}i})$

$+T_{\text{ab}j}T_{\text{bc}i}(T_{\text{ac}i}T_{\text{dd}j}-T_{\text{ad}%
i}T_{\text{cd}j})+T_{\text{ab}j}T_{\text{bc}i}(T_{\text{ac}j}T_{\text{dd}%
i}-T_{\text{ad}j}T_{\text{cd}i})](y_{0})$

$+\frac{1}{288}\underset{i,j=q+1}{\overset{n}{\sum}}[T_{\text{ab}%
i}T_{\text{bd}j}(T_{\text{ac}i}T_{\text{cd}j}-T_{\text{ad}i}T_{\text{cc}%
j})+T_{\text{ab}i}T_{\text{bd}j}(T_{\text{ac}j}T_{\text{cd}i}-T_{\text{ad}%
j}T_{\text{cc}i})$

$+T_{\text{ab}i}T_{\text{bd}j}(T_{\text{ac}j}T_{\text{cd}i}-T_{\text{ad}%
j}T_{\text{cc}i})+T_{\text{ab}j}T_{\text{bd}i}(T_{\text{ac}j}T_{\text{cd}%
i}-T_{\text{ad}j}T_{\text{cc}i})](y_{0})$

$-\ \frac{1}{288}\underset{i,j=q+1}{\overset{n}{\sum}}[T_{\text{ac}%
i}T_{\text{ab}j}(T_{\text{bc}i}T_{\text{dd}j}-T_{\text{bd}i}T_{\text{dc}%
j})+T_{\text{ac}i}T_{\text{ab}j}(T_{\text{bc}j}T_{\text{dd}i}-T_{\text{bd}%
j}T_{\text{dc}i})$

$+T_{\text{ac}j}T_{\text{ab}i}(T_{\text{bc}i}T_{\text{dd}j}-T_{\text{bd}%
i}T_{\text{dc}j})+T_{\text{ac}j}T_{\text{ab}i}(T_{\text{bc}j}T_{\text{dd}%
i}-T_{\text{bd}j}T_{\text{dc}i})](y_{0})$

$+\ \frac{1}{288}\underset{i,j=q+1}{\overset{n}{\sum}}[T_{\text{ac}%
i}T_{\text{bb}j}(T_{\text{ac}i}T_{\text{dd}j}-T_{\text{ad}i}T_{\text{cd}%
j})+T_{\text{ac}i}T_{\text{bb}j}(T_{\text{ac}j}T_{\text{dd}i}-T_{\text{ad}%
j}T_{\text{cd}i})$

$+T_{\text{ac}j}T_{\text{bb}i}(T_{\text{ac}i}T_{\text{dd}j}-T_{\text{ad}%
i}T_{\text{cd}i})+T_{\text{ac}j}T_{\text{bb}i}(T_{\text{ac}j}T_{\text{dd}%
i}-T_{\text{ad}j}T_{\text{cd}i})](y_{0})$

$-\ \frac{1}{288}\underset{i,j=q+1}{\overset{n}{\sum}}[T_{\text{ac}%
i}T_{\text{bd}j}(T_{\text{ac}i}T_{\text{bd}j}-T_{\text{ad}i}T_{\text{bc}%
j})+T_{\text{ac}i}T_{\text{bd}j}(T_{\text{ac}j}T_{\text{bd}i}-T_{\text{ad}%
j}T_{\text{bc}i})$

$+T_{\text{ac}j}T_{\text{bd}i}(T_{\text{ac}i}T_{\text{bd}j}-T_{\text{ad}%
i}T_{\text{bc}j})+T_{\text{ac}j}T_{\text{bd}i}(T_{\text{ac}j}T_{\text{bd}%
i}-T_{\text{ad}j}T_{\text{bc}i})](y_{0})$

$+\frac{1}{288}\underset{i,j=q+1}{\overset{n}{\sum}}[T_{\text{ad}%
i}T_{\text{ab}j}(T_{\text{bc}i}T_{\text{cd}j}-T_{\text{bd}i}T_{\text{cc}%
j})+T_{\text{ad}i}T_{\text{ab}j}(T_{\text{bc}j}T_{\text{cd}i}-T_{\text{bd}%
j}T_{\text{cc}i})$

$+T_{\text{ad}j}T_{\text{ab}i}(T_{\text{bc}i}T_{\text{cd}j}-T_{\text{bd}%
i}T_{\text{cc}j})+T_{\text{ad}j}T_{\text{ab}i}(T_{\text{bc}j}T_{\text{cd}%
i}-T_{\text{bd}j}T_{\text{cc}i})](y_{0})$

$-\ \frac{1}{288}\underset{i,j=q+1}{\overset{n}{\sum}}[T_{\text{ad}%
i}T_{\text{bb}j}(T_{\text{ac}i}T_{\text{cd}j}-T_{\text{ad}i}T_{\text{cc}%
j})+T_{\text{ad}i}T_{\text{bb}j}(T_{\text{ac}j}T_{\text{cd}i}-T_{\text{ad}%
j}T_{\text{cc}i})$

$+T_{\text{ad}j}T_{\text{bb}i}(T_{\text{ac}i}T_{\text{cd}j}-T_{\text{ad}%
i}T_{\text{cc}j})+T_{\text{ad}j}T_{\text{bb}i}(T_{\text{ac}j}T_{\text{cd}%
i}-T_{\text{ad}j}T_{\text{cc}i})](y_{0})$

$+\ \frac{1}{288}\underset{i,j=q+1}{\overset{n}{\sum}}[T_{\text{ad}%
i}T_{\text{bc}j}(T_{\text{ac}i}T_{\text{bd}j}-T_{\text{ad}i}T_{\text{bc}%
j})+T_{\text{ad}i}T_{\text{bc}j}(T_{\text{ac}j}T_{\text{bd}i}-T_{\text{ad}%
j}T_{\text{bc}i})$

$+T_{\text{ad}j}T_{\text{bc}i}(T_{\text{ac}i}T_{\text{bd}j}-T_{\text{ad}%
i}T_{\text{bc}j})+T_{\text{ad}j}T_{\text{bc}i}(T_{\text{ac}j}T_{\text{bd}%
i}-T_{\text{ad}j}T_{\text{bc}i})](y_{0})$

$-\ \frac{1}{144}[(R_{\text{cdcd}}^{P}-R_{\text{cdcd}}^{M})(R_{\text{abab}%
}^{P}-R_{\text{abab}}^{M})](y_{0})$

$+\frac{1}{144}[(R_{\text{bdcd}}^{P}-R_{\text{bdcd}}^{M})(R_{\text{abac}}%
^{P}-R_{\text{abac}}^{M})](y_{0})$

$+\ \frac{1}{144}[(R_{\text{bcdc}}^{P}-R_{\text{bcdc}}^{M})(R_{\text{abad}%
}^{P}-R_{\text{abad}}^{M})](y_{0})$

$-\ \frac{1}{144}[(R_{\text{adcd}}^{P}-R_{\text{adcd}}^{M})(R_{\text{abbc}%
}^{P}-R_{\text{abbc}}^{M})](y_{0})\qquad$

$+\ \frac{1}{144}[(R_{\text{acdc}}^{P}-R_{\text{acdc}}^{M})(R_{\text{abdb}%
}^{P}-R_{\text{abdb}}^{M})](y_{0})$

$-\ \frac{1}{576}[(R_{\text{abcd}}^{P}-R_{\text{abcd}}^{M})]^{2}(y_{0})$

$-\frac{1}{144}<H,i>(y_{0})<H,j>(y_{0})\qquad\qquad\qquad\qquad\qquad
\qquad\qquad L_{3}$

$\times\lbrack2\varrho_{ij}+$ $\overset{q}{\underset{\text{a}=1}{4\sum}%
}R_{i\text{a}j\text{a}}-3\overset{q}{\underset{\text{a,b=1}}{\sum}%
}(T_{\text{aa}i}T_{\text{bb}j}-T_{\text{ab}i}T_{\text{ab}j}%
)-3\overset{q}{\underset{\text{a,b=1}}{\sum}}(T_{\text{aa}j}T_{\text{bb}%
i}-T_{\text{ab}j}T_{\text{ab}i}](y_{0})\phi(y_{0})$

$-\frac{1}{16}[<H,i>^{2}(y_{0})<H,j>^{2}](y_{0})\phi(y_{0})$

$-\frac{1}{144}<H,i>(y_{0})<H,j>(y_{0})$

$\times\lbrack2\varrho_{ij}+$ $\overset{q}{\underset{\text{a}=1}{4\sum}%
}R_{i\text{a}j\text{a}}-3\overset{q}{\underset{\text{a,b=1}}{\sum}%
}(T_{\text{aa}i}T_{\text{bb}j}-T_{\text{ab}i}T_{\text{ab}j}%
)-3\overset{q}{\underset{\text{a,b=1}}{\sum}}(T_{\text{aa}j}T_{\text{bb}%
i}-T_{\text{ab}j}T_{\text{ab}i}](y_{0})\phi(y_{0})$

$-\frac{1}{72}<H,i>(y_{0})<H,k>(y_{0})R_{jijk}(y_{0})\phi(y_{0})$

$-\frac{1}{16}<H,i>^{2}(y_{0})<H,j>^{2}(y_{0})\phi(y_{0})$

$-\frac{1}{72}<H,i>^{2}(y_{0})[\varrho_{jj}+$ $\overset{q}{\underset{\text{a}%
=1}{2\sum}}R_{j\text{a}j\text{a}}-3\overset{q}{\underset{\text{a,b=1}}{\sum}%
}(T_{\text{aa}j}T_{\text{bb}j}-T_{\text{ab}j}T_{\text{ab}j})](y_{0})\phi
(y_{0})$

$+\frac{5}{32}[<H,i>^{2}<H,j>^{2}](y_{0})\phi(y_{0})$

$+\frac{1}{48}<H,i>(y_{0})<H,j>$

$\times\lbrack2\varrho_{ij}+$ $\overset{q}{\underset{\text{a}=1}{4\sum}%
}R_{i\text{a}j\text{a}}-3\overset{q}{\underset{\text{a,b=1}}{\sum}%
}(T_{\text{aa}i}T_{\text{bb}j}-T_{\text{ab}i}T_{\text{ab}j}%
)-3\overset{q}{\underset{\text{a,b=1}}{\sum}}(T_{\text{aa}j}T_{\text{bb}%
i}-T_{\text{ab}j}T_{\text{ab}i}](y_{0})\phi(y_{0})$

$+\frac{1}{48}<H,i>(y_{0})<H,i>(y_{0})[\tau^{M}\ -3\tau^{P}%
+\ \underset{\text{a}=1}{\overset{\text{q}}{\sum}}\varrho_{\text{aa}}^{M}+$
$\overset{q}{\underset{\text{a},\text{b}=1}{\sum}}R_{\text{abab}}^{M}$
$](y_{0})\phi(y_{0})$

$+\frac{1}{144}<H,i>(y_{0})[\nabla_{i}\varrho_{jj}-2\varrho_{ij}%
<H,j>+\overset{q}{\underset{\text{a}=1}{\sum}}(\nabla_{i}R_{\text{a}%
j\text{a}j}-4R_{i\text{a}j\text{a}}<H,j>)$

$+4\overset{q}{\underset{\text{a,b=1}}{\sum}}R_{i\text{a}j\text{b}%
}T_{\text{ab}j}+2\overset{q}{\underset{\text{a,b,c=1}}{\sum}}(T_{\text{aa}%
i}T_{\text{bb}j}T_{\text{cc}j}-3T_{\text{aa}i}T_{\text{bc}j}T_{\text{bc}%
j}+2T_{\text{ab}i}T_{\text{bc}j}T_{\text{ca}j})](y_{0})\phi(y_{0})$%
\qquad\qquad\qquad\qquad\qquad\ \ 

$+\frac{1}{144}<H,i>(y_{0})[\nabla_{j}\varrho_{ij}-2\varrho_{ij}%
<H,j>+\overset{q}{\underset{\text{a}=1}{\sum}}(\nabla_{j}R_{\text{a}%
i\text{a}j}-4R_{j\text{a}i\text{a}}<H,j>)$

$+4\overset{q}{\underset{\text{a,b=1}}{\sum}}R_{j\text{a}i\text{b}%
}T_{\text{ab}j}+2\overset{q}{\underset{\text{a,b,c=1}}{\sum}}(T_{\text{aa}%
j}T_{\text{bb}i}T_{\text{cc}j}-3T_{\text{aa}j}T_{\text{bc}i}T_{\text{bc}%
j}+2T_{\text{ab}j}T_{\text{bc}i}T_{\text{ac}j})](y_{0})\phi(y_{0})$

$+\frac{1}{144}<H,i>(y_{0})[\nabla_{j}\varrho_{ij}-2\varrho_{jj}%
<H,i>+\overset{q}{\underset{\text{a}=1}{\sum}}(\nabla_{j}R_{\text{a}%
i\text{a}j}-4R_{j\text{a}j\text{a}}<H,i>)+4\overset{q}{\underset{\text{a,b=1}%
}{\sum}}R_{j\text{a}j\text{b}}T_{\text{ab}i}$

$+2\overset{q}{\underset{\text{a,b,c=1}}{\sum}}(T_{\text{aa}j}T_{\text{bb}%
j}T_{\text{cc}i}-3T_{\text{aa}j}T_{\text{bc}j}T_{\text{bc}i}+2T_{\text{ab}%
j}T_{\text{bc}j}T_{\text{ac}i})](y_{0})\phi(y_{0})$

$-\frac{1}{192}<H,i>^{2}<H,j>^{2}(y_{0})\phi(y_{0})\qquad\qquad\qquad
\qquad\qquad\qquad\qquad$I$_{3213}\ $

$-\frac{1}{288}<H,i>^{2}(y_{0})[\tau^{M}\ -3\tau^{P}+\ \underset{\text{a}%
=1}{\overset{\text{q}}{\sum}}\varrho_{\text{aa}}^{M}+$
$\overset{q}{\underset{\text{a},\text{b}=1}{\sum}}R_{\text{abab}}^{M}$
$](y_{0})\phi(y_{0})$

$-\frac{1}{288}<H,i>(y_{0})<H,j>(y_{0})$

$\times\lbrack2\varrho_{ij}+$ $\overset{q}{\underset{\text{a}=1}{4\sum}%
}R_{i\text{a}j\text{a}}-3\overset{q}{\underset{\text{a,b=1}}{\sum}%
}(T_{\text{aa}i}T_{\text{bb}j}-T_{\text{ab}i}T_{\text{ab}j}%
)-3\overset{q}{\underset{\text{a,b=1}}{\sum}}(T_{\text{aa}j}T_{\text{bb}%
i}-T_{\text{ab}j}T_{\text{ab}i}](y_{0})\phi(y_{0})$

$\ +\frac{1}{144}<H,i>(y_{0})<H,k>(y_{0})R_{jijk}(y_{0})$

$+\frac{1}{144}<H,i>^{2}(y_{0})[\varrho_{jj}+$ $\overset{q}{\underset{\text{a}%
=1}{2\sum}}R_{j\text{a}j\text{a}}-3\overset{q}{\underset{\text{a,b=1}}{\sum}%
}(T_{\text{aa}j}T_{\text{bb}j}-T_{\text{ab}j}T_{\text{ab}j})](y_{0})\phi
(y_{0})$

$-\frac{1}{288}<H,i>(y_{0})[\nabla_{i}\varrho_{jj}-2\varrho_{ij}%
<H,j>+\overset{q}{\underset{\text{a}=1}{\sum}}(\nabla_{i}R_{\text{a}%
j\text{a}j}-4R_{i\text{a}j\text{a}}<H,j>)$

$+4\overset{q}{\underset{\text{a,b=1}}{\sum}}R_{i\text{a}j\text{b}%
}T_{\text{ab}j}+2\overset{q}{\underset{\text{a,b,c=1}}{\sum}}(T_{\text{aa}%
i}T_{\text{bb}j}T_{\text{cc}j}-3T_{\text{aa}i}T_{\text{bc}j}T_{\text{bc}%
j}+2T_{\text{ab}i}T_{\text{bc}j}T_{\text{ca}j})](y_{0})\phi(y_{0})$%
\qquad\qquad\qquad\qquad\qquad\ \ 

$-\frac{1}{288}<H,i>(y_{0})[\nabla_{j}\varrho_{ij}-2\varrho_{ij}%
<H,j>+\overset{q}{\underset{\text{a}=1}{\sum}}(\nabla_{j}R_{\text{a}%
i\text{a}j}-4R_{j\text{a}i\text{a}}<H,j>)$

$+4\overset{q}{\underset{\text{a,b=1}}{\sum}}R_{j\text{a}i\text{b}%
}T_{\text{ab}j}+2\overset{q}{\underset{\text{a,b,c=1}}{\sum}}(T_{\text{aa}%
j}T_{\text{bb}i}T_{\text{cc}j}-3T_{\text{aa}j}T_{\text{bc}i}T_{\text{bc}%
j}+2T_{\text{ab}j}T_{\text{bc}i}T_{\text{ac}j})](y_{0})\phi(y_{0})$

$-\frac{1}{288}<H,i>(y_{0})[\nabla_{j}\varrho_{ij}-2\varrho_{jj}%
<H,i>+\overset{q}{\underset{\text{a}=1}{\sum}}(\nabla_{j}R_{\text{a}%
i\text{a}j}-4R_{j\text{a}j\text{a}}<H,i>)+4\overset{q}{\underset{\text{a,b=1}%
}{\sum}}R_{j\text{a}j\text{b}}T_{\text{ab}i}$

$+2\overset{q}{\underset{\text{a,b,c=1}}{\sum}}(T_{\text{aa}j}T_{\text{bb}%
j}T_{\text{cc}i}-3T_{\text{aa}j}T_{\text{bc}j}T_{\text{bc}i}+2T_{\text{ab}%
j}T_{\text{bc}j}T_{\text{ac}i})](y_{0})\phi(y_{0})$

$+\frac{1}{24}[\left\Vert \text{X}\right\Vert _{M}^{2}+\operatorname{div}%
$X$_{M}-\left\Vert \text{X}\right\Vert _{P}^{2}-\operatorname{div}X_{P}%
](y_{0})[\left\Vert \text{X}\right\Vert _{M}^{2}-\operatorname{div}$%
X$_{M}-\left\Vert \text{X}\right\Vert _{P}^{2}+\operatorname{div}$%
X$_{P}](y_{0})\phi(y_{0})\qquad$\ I$_{3212}$

$+\frac{1}{6}X_{i}(y_{0})T_{\text{ab}i}(y_{0})T_{\text{ab}j}(y_{0})X_{j}%
(y_{0})+\frac{1}{3}\perp_{\text{a}ij}(y_{0})X_{i}(y_{0})[\frac{\partial X_{j}%
}{\partial x_{\text{a}}}-\perp_{\text{a}jk}X_{k}](y_{0})\qquad$I$_{32122}%
\qquad Q_{1}$

$+$ $\frac{2}{3}X_{i}(y_{0})X_{j}(y_{0})\frac{\partial X_{j}}{\partial
x_{\text{a}}}(y_{0})-\frac{1}{6}X_{i}(y_{0})\frac{\partial^{2}X_{j}}{\partial
x_{\text{a}}\partial x_{j}}(y_{0})\qquad\qquad Q_{2}$

$-\frac{1}{12}X_{i}(y_{0})\frac{\partial^{2}X_{i}}{\partial x_{\text{a}}^{2}%
}(y_{0})+\frac{1}{12}X_{i}^{2}(y_{0})[\operatorname{div}X_{M}-\left\Vert
X\right\Vert _{M}^{2}+\left\Vert X\right\Vert _{P}^{2}-\operatorname{div}%
X_{P}-$ $<H,j>X_{j}](y_{0})$

$+\frac{1}{6}X_{i}(y_{0})X_{j}(y_{0})\frac{\partial X_{i}}{\partial x_{j}%
}(y_{0})+$ $\frac{1}{18}X_{i}(y_{0})X_{k}(y_{0})[R_{jijk}](y_{0})-\frac{1}%
{12}X_{i}(y_{0})\frac{\partial^{2}X_{i}}{\partial x_{j}^{2}}(y_{0})$

$+\frac{1}{12}[R_{\text{a}i\text{a}k}-\underset{\text{c=1}}{\overset{\text{q}%
}{\sum}}T_{\text{ac}i}T_{\text{ac}k}-\perp_{\text{a}ik}\perp_{\text{a}%
jk}](y_{0})X_{k}(y_{0})+\frac{1}{18}R_{ijkj}(y_{0})X_{i}(y_{0})X_{k}(y_{0})$

$+\frac{1}{12}<H,j>(y_{0})X_{i}(y_{0})[X_{i}X_{j}-\frac{1}{2}\left(
\frac{\partial X_{j}}{\partial x_{i}}+\frac{\partial X_{i}}{\partial x_{j}%
}\right)  ](y_{0})\qquad\qquad\qquad\qquad\qquad\qquad\qquad\qquad\qquad
\qquad\ \qquad\qquad\qquad\qquad\qquad\qquad\qquad\qquad\qquad\qquad\qquad$

$-\frac{1}{6}[-R_{\text{a}i\text{b}i}+5\overset{q}{\underset{\text{c}=1}{\sum
}}T_{\text{ac}i}T_{\text{bc}i}+2\overset{n}{\underset{j=q+1}{\sum}}%
\perp_{\text{a}ij}\perp_{\text{b}ij}](y_{0})\underset{k=q+1}{\overset{n}{\sum
}}T_{\text{ab}k}(y_{0})X_{k}(y_{0})\qquad\qquad$I$_{32123}\qquad S_{1}\qquad$

$-\frac{2}{9}\underset{j=q+1}{\overset{n}{\sum}}R_{i\text{a}ij}(y_{0}%
)[\frac{\partial X_{j}}{\partial x_{\text{a}}}%
-\underset{k=q+1}{\overset{n}{\sum}}\perp_{\text{a}jk}X_{k}](y_{0})+\frac
{1}{12}\times\frac{2}{3}\underset{j,k=q+1}{\overset{n}{\sum}}R_{ijik}%
(y_{0})[X_{j}X_{k}-\frac{1}{2}(\frac{\partial X_{j}}{\partial x_{k}}%
+\frac{\partial X_{k}}{\partial x_{j}})](y_{0})$

$-\frac{1}{6}T_{\text{ab}i}(y_{0})\frac{\partial^{2}X_{i}}{\partial
x_{\text{a}}\partial x_{\text{b}}}(y_{0})\qquad\qquad\qquad S_{2}\qquad\qquad
S_{21}\qquad\qquad\qquad\qquad\qquad\qquad$

$+$ $\frac{1}{12}T_{\text{ab}i}(y_{0})[$ $(R_{\text{a}i\text{b}j}%
+R_{\text{a}j\text{b}i})$ $-\underset{\text{c=1}}{\overset{\text{q}}{\sum}%
}(T_{\text{ac}i}T_{\text{bc}j}+T_{\text{ac}j}T_{\text{bc}i}%
)-\overset{n}{\underset{k=q+1}{\sum}}(\perp_{\text{a}ik}\perp_{\text{b}jk}+$
$\perp_{\text{a}jk}\perp_{\text{b}ik})](y_{0})X_{j}(y_{0})$

$-$ $\frac{1}{6}T_{\text{ab}i}(y_{0})T_{\text{ab}j}(y_{0})[X_{i}X_{j}-\frac
{1}{2}\left(  \frac{\partial X_{i}}{\partial x_{j}}+\frac{\partial X_{j}%
}{\partial x_{i}}\right)  ](y_{0})$

$-\frac{1}{3}\perp_{\text{a}ij}(y_{0})[(X_{i}\frac{\partial X_{j}}{\partial
x_{\text{a}}}+X_{j}\frac{\partial X_{i}}{\partial x_{\text{a}}})-\frac{1}%
{4}\left(  \frac{\partial^{2}X_{i}}{\partial x_{\text{a}}\partial x_{j}}%
+\frac{\partial^{2}X_{j}}{\partial x_{\text{a}}\partial x_{i}}\right)
](y_{0})\qquad S_{22}$

$-\frac{1}{6}\perp_{\text{a}ij}(y_{0})[T_{\text{ab}j}\frac{\partial X_{i}%
}{\partial x_{\text{b}}}](y_{0})$

$+\frac{1}{6}\perp_{\text{a}ij}(y_{0})[(\perp_{\text{b}ik}T_{\text{ab}%
j})+\frac{2}{3}(2R_{\text{a}ijk}+R_{\text{a}jik}+R_{\text{a}kji})](y_{0}%
)X_{k}(y_{0})$

$-\frac{1}{6}\perp_{\text{a}ij}(y_{0})\perp_{\text{a}jk}(y_{0})[X_{i}%
X_{k}-\frac{1}{2}\left(  \frac{\partial X_{i}}{\partial x_{k}}+\frac{\partial
X_{k}}{\partial x_{i}}\right)  ](y_{0})\qquad\qquad\qquad\qquad\qquad
\qquad\qquad\qquad$

$+\frac{1}{12}[(\frac{\partial X_{j}}{\partial x_{\text{a}}})^{2}+X_{j}%
\frac{\partial^{2}X_{j}}{\partial x_{\text{a}}^{2}}-\frac{1}{2}\frac
{\partial^{3}X_{j}}{\partial x_{\text{a}}^{2}\partial x_{j}}](y_{0})-\frac
{1}{6}\overset{n}{\underset{k=q+1}{\sum}}[\perp_{\text{b}ik}$T$_{\text{aa}%
k}\frac{\partial X_{i}}{\partial x_{\text{b}}^{2}}](y_{0})\qquad\qquad
S_{3}\qquad S_{31}$

$+\frac{1}{144}[\{4\nabla_{i}R_{i\text{a}j\text{a}}+2\nabla_{j}R_{i\text{a}%
i\text{a}}+$ $8(\overset{q}{\underset{\text{c=1}}{%
{\textstyle\sum}
}}R_{\text{a}i\text{c}i}^{{}}T_{\text{ac}j}+\;\overset{n}{\underset{k=q+1}{%
{\textstyle\sum}
}}R_{\text{a}iik}\perp_{\text{a}jk})$

$+8(\overset{q}{\underset{\text{c=1}}{%
{\textstyle\sum}
}}R_{\text{a}i\text{c}j}^{{}}T_{\text{ac}i}+\;\overset{n}{\underset{k=q+1}{%
{\textstyle\sum}
}}R_{\text{a}ijk}\perp_{\text{a}ik})+8(\overset{q}{\underset{\text{c=1}}{%
{\textstyle\sum}
}}R_{\text{a}j\text{c}i}^{{}}T_{\text{ac}i}+\;\overset{n}{\underset{k=q+1}{%
{\textstyle\sum}
}}R_{\text{a}jik}\perp_{\text{a}ik})\}$\ 

$+\frac{2}{3}\underset{k=q+1}{\overset{n}{\sum}}\{T_{\text{aa}k}%
(R_{ijik}+3\overset{q}{\underset{\text{c}=1}{\sum}}\perp_{\text{c}ij}%
\perp_{\text{c}ik})\}](y_{0})X_{k}(y_{0})$

$-\frac{1}{12}[$ R$_{\text{a}i\text{a}k}$ $-\underset{\text{c=1}%
}{\overset{\text{q}}{\sum}}T_{\text{ac}i}T_{\text{ac}k}%
-\overset{n}{\underset{l=q+1}{\sum}}(\perp_{\text{a}il}\perp_{\text{a}%
kl}](y_{0})\times\lbrack X_{i}X_{k}-\frac{1}{2}\left(  \frac{\partial X_{i}%
}{\partial x_{k}}+\frac{\partial X_{k}}{\partial x_{i}}\right)  ](y_{0})$

$-\frac{1}{24}T_{\text{aa}k}(y_{0})[-X_{i}^{2}X_{k}+X_{k}\frac{\partial X_{i}%
}{\partial x_{i}}\ +X_{i}\left(  \frac{\partial X_{k}}{\partial x_{i}}%
+\frac{\partial X_{i}}{\partial x_{k}}\right)  -\frac{1}{3}\left(
\frac{\partial^{2}X_{k}}{\partial x_{i}^{2}}+2\frac{\partial^{2}X_{i}%
}{\partial x_{i}\partial x_{k}}\right)  ](y_{0})$

$+\frac{1}{18}[R_{\text{a}jij}\frac{\partial X_{i}}{\partial x_{\text{a}}^{2}%
}](y_{0})\qquad\qquad\qquad\qquad\qquad\qquad\qquad\qquad\qquad S_{32}$

$+\frac{1}{24}[\frac{4}{3}\overset{q}{\underset{\text{a}=1}{\sum}}%
\perp_{\text{a}ki}R_{ij\text{a}j}-\frac{1}{3}(\nabla_{i}R_{kjij}+\nabla
_{j}R_{ijik}+\nabla_{k}R_{ijij})](y_{0})X_{k}(y_{0})$

\ $-\frac{1}{18}R_{ijkj}(y_{0})[X_{i}X_{k}-\frac{1}{2}\left(  \frac{\partial
X_{i}}{\partial x_{k}}+\frac{\partial X_{k}}{\partial x_{i}}\right)  ](y_{0})$

$+\frac{1}{24}[X_{i}^{2}X_{j}^{2}-2X_{i}X_{j}\left(  \frac{\partial X_{j}%
}{\partial x_{i}}+\frac{\partial X_{i}}{\partial x_{j}}\right)  -X_{i}%
^{2}\frac{\partial X_{j}}{\partial x_{j}}-X_{j}^{2}\frac{\partial X_{i}%
}{\partial x_{i}}](y_{0})$

$+\frac{1}{48}\left(  \frac{\partial X_{j}}{\partial x_{i}}+\frac{\partial
X_{i}}{\partial x_{j}}\right)  ^{2}(y_{0})+\frac{1}{24}\left(  \frac{\partial
X_{i}}{\partial x_{i}}\frac{\partial X_{j}}{\partial x_{j}}\right)
(y_{0})\qquad$

$+\frac{1}{36}X_{i}(y_{0})\left(  2\frac{\partial^{2}X_{j}}{\partial
x_{i}\partial x_{j}}+\frac{\partial^{2}X_{i}}{\partial x_{j}^{2}}\right)
(y_{0})+\frac{1}{36}X_{j}(y_{0})\left(  \frac{\partial^{2}X_{j}}{\partial
x_{i}^{2}}+2\frac{\partial^{2}X_{i}}{\partial x_{i}\partial x_{j}}\right)
(y_{0})$

$-\frac{1}{48}\left(  \frac{\partial^{3}X_{i}}{\partial x_{i}\partial
x_{j}^{2}}+\frac{\partial^{3}X_{j}}{\partial x_{i}^{2}\partial x_{j}}\right)
(y_{0})$

$+$ $\frac{2}{3}<H,j>(y_{0})\left(  \frac{\partial^{2}X_{i}}{\partial
x_{i}\partial x_{j}}+2\frac{\partial^{2}X_{j}}{\partial x_{i}^{2}}\right)
(y_{0})\phi(y_{0})+$ $\frac{2}{3}<H,j>(y_{0})R_{ijik}(y_{0})X_{k}(y_{0}%
)\phi(y_{0})\qquad$I$_{3213}$

$+\frac{1}{12}[<H,i><H,j>\ +\frac{1}{6}(2\varrho_{ij}%
+4\overset{q}{\underset{\text{a}=1}{\sum}}R_{i\text{a}j\text{a}}%
-6\overset{q}{\underset{\text{a,b}=1}{\sum}}T_{\text{aa}i}T_{\text{bb}%
j}-T_{\text{ab}i}T_{\text{ab}j})](y_{0})\phi(y_{0})$

$\times\frac{1}{2}[\left(  \frac{\partial X_{j}}{\partial x_{i}}%
-\frac{\partial X_{i}}{\partial x_{j}}\right)  ](y_{0})\phi(y_{0})$

$-\frac{1}{12}\perp_{\text{a}ij}(y_{0})<H,i>(y_{0})[(X_{j}\perp_{\text{a}%
ij}-\frac{\partial X_{i}}{\partial x_{\text{a}}})+\frac{\partial X_{\text{a}}%
}{\partial x_{i}}](y_{0})\phi(y_{0})$

$-\frac{1}{18}[X_{j}\left(  2\frac{\partial^{2}X_{j}}{\partial x_{i}^{2}%
}+\frac{\partial^{2}X_{i}}{\partial x_{i}\partial x_{j}}\right)  ](y_{0}%
)\phi(y_{0})-\frac{1}{12}[\left(  \frac{\partial X_{i}}{\partial x_{j}}%
+\frac{\partial X_{j}}{\partial x_{i}}\right)  ]\frac{\partial X_{j}}{\partial
x_{i}}(y_{0})\phi(y_{0})\qquad\ $I$_{3214}\qquad\qquad\qquad$\qquad
\qquad\qquad\qquad\qquad\qquad\qquad\qquad\qquad\qquad\qquad\qquad\qquad
\qquad\qquad\qquad

$+\frac{1}{12}\frac{\partial^{2}\text{V}}{\partial x_{i}^{2}}(y_{0})\phi
(y_{0})\qquad\qquad$I$_{3215}$

$+\frac{1}{12}\underset{i=q+\text{1}}{\overset{\text{n}}{\sum}}%
\underset{\text{a,b=1}}{\overset{\text{q}}{\sum}}[-R_{\text{a}i\text{b}%
i}+5\overset{q}{\underset{\text{c}=1}{\sum}}T_{\text{ac}i}T_{\text{bc}%
i}+2\overset{n}{\underset{j=q+1}{\sum}}\perp_{\text{a}ij}\perp_{\text{b}%
ij}](y_{0})\times\frac{\partial^{2}\phi}{\partial\text{x}_{\text{a}}%
\partial\text{x}_{\text{b}}}(y_{0})\qquad$\textbf{I}$_{322}$

$+\frac{1}{72}\underset{i,j,k=q+1}{\overset{n}{\sum}}R_{ijik}(y_{0}%
)\Omega_{jk}(y_{0})\phi(y_{0})\qquad\qquad\qquad\qquad$\ I$_{323}$

$+\frac{1}{24}\underset{\text{a=1}}{\overset{\text{q}}{\sum}}%
\underset{i,j=q+1}{\overset{n}{\sum}}\left\{  \frac{8}{3}R_{i\text{a}%
ij}+4\underset{\text{b=1}}{\overset{\text{q}}{\sum}}T_{\text{ab}i}%
\perp_{\text{b}ji}\right\}  (y_{0})\left\{  -\Omega_{\text{a}j}+[\Lambda
_{\text{a}},\Lambda_{j}]\right\}  (y_{0})\phi(y_{0})$

$+\frac{1}{12}\underset{i\text{=}q+1}{\overset{\text{n}}{\sum}}%
\underset{\text{a,b=1}}{\overset{\text{q}}{\sum}}$ $[-R_{\text{a}i\text{b}%
i}+5\overset{q}{\underset{\text{c}=1}{\sum}}T_{\text{ac}i}T_{\text{bc}%
i}+2\overset{n}{\underset{\text{k}=q+1}{\sum}}\perp_{\text{a}i\text{k}}%
\perp_{\text{b}i\text{k}}](y_{0})\times\lbrack\Lambda_{\text{a}}(y_{0}%
)\Lambda_{\text{b}}(y_{0})\phi(y_{0})]\qquad$\ I$_{324}$

$+\frac{1}{12}[\frac{8}{3}R_{i\text{a}ij}-4\underset{\text{b=1}%
}{\overset{\text{q}}{\sum}}T_{\text{ab}i}(y_{0})\perp_{\text{b}ij}%
](y_{0})[\Lambda_{\text{a}}\Lambda_{j}\phi](y_{0})$

$+\frac{1}{12}\underset{i\text{=}q+1}{\overset{\text{n}}{\sum}}%
\underset{\text{a,b=1}}{\overset{\text{q}}{\sum}}$ $[-R_{\text{a}i\text{b}%
i}+5\overset{q}{\underset{\text{c}=1}{\sum}}T_{\text{ac}i}T_{\text{bc}%
i}+2\overset{n}{\underset{\text{k}=q+1}{\sum}}\perp_{\text{a}i\text{k}}%
\perp_{\text{b}i\text{k}}](y_{0})\qquad$I$_{325}\qquad$I$_{3251}$

$\times\lbrack\Lambda_{\text{a}}(y_{0})\Lambda_{\text{b}}(y_{0})\phi(y_{0})]$

$+\frac{1}{12}[\frac{8}{3}R_{i\text{a}ij}-4\underset{\text{b=1}%
}{\overset{\text{q}}{\sum}}T_{\text{ab}i}(y_{0})\perp_{\text{b}ij}%
](y_{0})[\Lambda_{\text{a}}\Lambda_{j}\phi](y_{0})$

$+\frac{1}{24}\underset{i=q+1}{\overset{n}{\sum}}\underset{\text{a=1}%
}{\overset{\text{q}}{\sum}}\left(  \frac{\partial\Omega_{i\text{a}}}{\partial
x_{i}}\Lambda_{\text{a}}+[\Omega_{i\text{a}}+[\Lambda_{\text{a}},\Lambda
_{i}],\Lambda_{i}]\right)  \Lambda_{\text{a}}(y_{0})\phi(y_{0})\qquad
$I$_{3252}$

$+\frac{1}{24}\underset{i=q+1}{\overset{n}{\sum}}\underset{\text{a=1}%
}{\overset{\text{q}}{\sum}}\Lambda_{\text{a}}(y_{0})\left(  \frac
{\partial\Omega_{i\text{a}}}{\partial x_{i}}\Lambda_{\text{a}}+[\Omega
_{i\text{a}}+[\Lambda_{\text{a}},\Lambda_{i}],\Lambda_{i}]\right)  (y_{0}%
)\phi(y_{0})$

$+\frac{1}{12}\underset{i=q+1}{\overset{n}{\sum}}\underset{\text{a=1}%
}{\overset{\text{q}}{\sum}}\left(  \Omega_{i\text{a}}+[\Lambda_{\text{a}%
},\Lambda_{i}]\right)  ^{2}(y_{0})\phi(y_{0})$

$+\frac{1}{48}\underset{i,j=q+1}{\overset{n}{\sum}}\left(  \Omega_{ij}%
\Omega_{ij}\right)  (y_{0})\phi(y_{0})+\frac{1}{72}%
\underset{i,j=q+1}{\overset{n}{\sum}}\left(  \frac{\partial\Omega_{ij}%
}{\partial\text{x}_{i}}\Lambda_{j}+\Lambda_{j}\frac{\partial\Omega_{ij}%
}{\partial\text{x}_{i}}\right)  (y_{0})\phi(y_{0})$

$+\frac{1}{12}[\underset{i=q+1}{\overset{n}{\sum}}\underset{\text{a,b=1}%
}{\overset{\text{q}}{\sum}}2T_{\text{ab}i}(y_{0})\left\{  (\Omega_{i\text{a}%
}+[\Lambda_{\text{a}},\Lambda_{i}])\Lambda_{\text{b}}+\Lambda_{\text{a}%
}(\Omega_{i\text{a}}+[\Lambda_{\text{a}},\Lambda_{i}])\right\}  (y_{0}%
)\phi(y_{0})\qquad$I$_{3253}$

$-\frac{1}{12}[\underset{i,j=q+1}{\overset{n}{\sum}}\underset{\text{a=1}%
}{\overset{\text{q}}{\sum}}\perp_{\text{a}ij}(y_{0})\left\{  (\Omega
_{i\text{a}}+[\Lambda_{\text{a}},\Lambda_{i}])\Lambda_{j}+\frac{1}{2}%
\Lambda_{\text{a}}\Omega_{ij}\right\}  ](y_{0})\phi(y_{0})$

$-\frac{1}{12}[\underset{i,j=q+1}{\overset{n}{\sum}}\underset{\text{b=1}%
}{\overset{\text{q}}{\sum}}\perp_{\text{b}ij}(y_{0})\left\{  \frac{1}{2}%
\Omega_{ij}\Lambda_{\text{b}}+\Lambda_{j}(\Omega_{i\text{b}}+[\Lambda
_{\text{b}},\Lambda_{i}])\right\}  ](y_{0})\phi(y_{0})$

$\mathbf{-}\frac{1}{12}\underset{\text{a,b=1}}{\overset{\text{q}}{\sum}%
}\underset{i,j=q+1}{\overset{n}{\sum}}T_{\text{ab}i}(y_{0})[-R_{\text{a}%
i\text{b}i}+5\overset{q}{\underset{\text{c}=1}{\sum}}T_{\text{ac}%
i}T_{\text{bc}i}+2\overset{n}{\underset{k=q+1}{\sum}}\perp_{\text{a}ik}%
\perp_{\text{b}ik}](y_{0})\Lambda_{j}(y_{0})\phi(y_{0})\ \ \qquad$%
\textbf{I}$_{326}\qquad$\textbf{I}$_{3261}\qquad$

$+\frac{1}{12}\underset{i\text{=q+1}}{\overset{\text{n}}{\sum}}%
\underset{j\text{=q+1}}{\overset{\text{n}}{\sum}}\underset{\text{a=1}%
}{\overset{\text{q}}{\sum}}[4\underset{\text{c=1}}{\overset{q}{\sum}%
}(T_{\text{ac}i})(\perp_{j\text{c}i})+\frac{8}{3}R_{i\text{a}ij}](y_{0}%
)\qquad\qquad$I$_{32613}$

$\times\lbrack\underset{\text{c=1}}{\overset{\text{q}}{\sum}}T_{\text{ac}%
j}\frac{\partial\phi}{\partial\text{x}_{\text{c}}}+\underset{\text{b=1}%
}{\overset{\text{q}}{\sum}}T_{\text{ab}j}\Lambda_{\text{b}}%
-\underset{k=q+1}{\overset{n}{\sum}}\perp_{\text{a}jk}\Lambda_{k}](y_{0}%
)\phi(y_{0})$

$-\frac{1}{24}\underset{i=q+1}{\overset{n}{\sum}}\underset{\text{a,b=1}%
}{\overset{\text{q}}{\sum}}$ $\overset{n}{\underset{k=q+1}{\sum}}%
$T$_{\text{aa}k}[\frac{8}{3}R_{i\text{c}ik}+4\underset{\text{d}%
=1}{\overset{q}{\sum}}(T_{\text{db}k})(\perp_{\text{d}ik})]\frac{\partial\phi
}{\partial\text{x}_{\text{b}}}$ $(y_{0})\qquad$ I$_{32621}\qquad$I$_{3262}$

$-\frac{1}{12}\underset{i=q+1}{\overset{n}{\sum}}\underset{\text{a,b=1}%
}{\overset{\text{q}}{\sum}}[$ $\overset{n}{\underset{k,l=q+1}{\sum}}%
\perp_{\text{b}ik}(-R_{\text{a}k\text{a}l}+\underset{\text{d=1}%
}{\overset{\text{q}}{\sum}}T_{\text{ad}k}T_{\text{ad}l}%
))-\overset{n}{\underset{k,l=q+1}{\sum}}\perp_{\text{b}ik}%
(\underset{r=q+1}{\overset{n}{\sum}}\perp_{\text{a}kr}\perp_{\text{a}%
lr})](y_{0})\frac{\partial\phi}{\partial\text{x}_{\text{b}}}$ $(y_{0})$

$-\frac{1}{12}\underset{i=q+1}{\overset{n}{\sum}}\underset{\text{b=1}%
}{\overset{\text{q}}{\sum}}[\frac{8}{3}\overset{q}{\underset{\text{c}=1}{\sum
}}($T$_{\text{bc}i}R_{ij\text{c}j})+\frac{2}{3}%
\overset{n}{\underset{k=q+1}{\sum}}(\perp_{\text{b}ik}R_{ijjk})](y_{0}%
)\frac{\partial\phi}{\partial\text{x}_{\text{b}}}(y_{0})$

$-\frac{1}{6}\underset{i=q+1}{\overset{n}{\sum}}\underset{\text{a,b=1}%
}{\overset{\text{q}}{\sum}}[4\overset{q}{\underset{\text{c=1}}{%
{\textstyle\sum}
}}$R$_{ij\text{c}i}^{{}}T_{\text{bc}j}+$ $4\overset{n}{\underset{k=q+1}{%
{\textstyle\sum}
}}$R$_{ijik}\perp_{\text{b}jk}+3\nabla_{i}$R$_{j\text{b}ij}%
+4\overset{q}{\underset{\text{c=1}}{%
{\textstyle\sum}
}}$R$_{ij\text{c}j}^{{}}T_{\text{bc}i}+$ $4$R$_{ijjk}\perp_{\text{b}ik}%
](y_{0})\frac{\partial\phi}{\partial\text{x}_{\text{b}}}(y_{0})$

$-\frac{1}{144}[\{4\nabla_{i}$R$_{i\text{a}j\text{a}}$ $+2\nabla_{j}%
$R$_{i\text{a}i\text{a}}+$ $8$ $(\overset{q}{\underset{\text{c=1}}{%
{\textstyle\sum}
}}R_{\text{a}i\text{c}i}^{{}}T_{\text{ac}j}+\;\overset{n}{\underset{k=q+1}{%
{\textstyle\sum}
}}R_{\text{a}iik}\perp_{\text{a}jk})\qquad$I$_{326221}\qquad$I$_{32622}$

\ $+8(\overset{q}{\underset{\text{c=1}}{%
{\textstyle\sum}
}}R_{\text{a}i\text{c}j}^{{}}T_{\text{ac}i}+\overset{n}{\underset{l=q+1}{%
{\textstyle\sum}
}}R_{\text{a}ijl}\perp_{\text{a}il})+8(\overset{q}{\underset{\text{c=1}}{%
{\textstyle\sum}
}}R_{\text{a}j\text{c}i}^{{}}T_{\text{ac}i}+\overset{q}{\underset{\text{c=1}}{%
{\textstyle\sum}
}}R_{\text{a}j\text{c}i}^{{}}T_{\text{ac}i})\}$\ 

$+\frac{2}{3}\underset{k=q+1}{\overset{n}{\sum}}\{T_{\text{aa}k}%
(R_{ijik}+3\overset{q}{\underset{\text{c}=1}{\sum}}\perp_{\text{c}ij}%
\perp_{\text{c}ik})\}](y_{0})\Lambda_{k}(y_{0})\phi(y_{0})$

\ $+\frac{1}{24}[\frac{4}{3}\overset{q}{\underset{\text{a}=1}{\sum}}%
\perp_{\text{a}ik}R_{ij\text{a}j}+\frac{1}{3}(\nabla_{i}R_{kjij}+\nabla
_{j}R_{ijik}+\nabla_{k}R_{ijij})](y_{0})\Lambda_{k}(y_{0})\phi(y_{0})$

$\mathbf{-}$ $\frac{1}{72}\underset{i,j=q+1}{\overset{n}{\sum}}%
\underset{\text{a}=1}{\overset{q}{\sum}}T_{\text{aa}j}(y_{0})\frac
{\partial\Omega_{ij}}{\partial\text{x}_{i}}(y_{0})\phi(y_{0})\qquad
$I$_{326222}$

$+$ $\frac{1}{12}$ $\overset{n}{\underset{j=q+1}{\sum}}(\perp_{\text{b}ij}%
$T$_{\text{aa}j})(y_{0})\left(  \Omega_{i\text{b}}(y_{0})+[\Lambda_{\text{b}%
},\Lambda_{i}]\right)  (y_{0})\phi(y_{0})$ \qquad I$_{326223}$

\ $\mathbf{-}\frac{1}{18}$ $\underset{i,j=q+1}{\overset{n}{\sum}%
}\underset{\text{b=1}}{\overset{q}{\sum}}R_{\text{b}jij}(y_{0})\left(
\Omega_{i\text{a}}(y_{0})+[\Lambda_{\text{a}},\Lambda_{i}]\right)  (y_{0}%
)\phi(y_{0})$

$\mathbf{-}$ $\frac{1}{24}\underset{i,j=q+1}{\overset{n}{\sum}}%
\underset{\text{a=1}}{\overset{q}{\sum}}[$ R$_{\text{a}i\text{a}j}$
$-\underset{\text{c=1}}{\overset{\text{q}}{\sum}}T_{\text{ac}i}T_{\text{ac}%
j}-\overset{n}{\underset{k=q+1}{\sum}}(\perp_{\text{a}ik}\perp_{\text{a}%
jk}](y_{0})\Omega_{ij}(y_{0})\phi(y_{0})$

\ $\mathbf{-}$ $\frac{1}{36}\underset{i,j,k=q+1}{\overset{n}{\sum}}%
R_{ijkj}(y_{0})\Omega_{ik}(y_{0})(y_{0})\phi(y_{0})$

$+\frac{1}{6}\underset{i,k=q+1}{\overset{n}{\sum}}\underset{\text{a,b,c}%
=1}{\overset{q}{\sum}}T_{\text{ab}i}(y_{0})[(\perp_{\text{c}ik}T_{\text{ab}%
k})(\frac{\partial\phi}{\partial\text{x}_{\text{c}}}$ $+$ $\Lambda_{\text{c}%
}\phi)](y_{0})\qquad$\textbf{I}$_{32631}\qquad$\textbf{I}$_{3263}$

$-\frac{1}{12}[$ $(R_{\text{a}i\text{b}j}+R_{\text{a}j\text{b}i})$
$-\underset{\text{c=1}}{\overset{\text{q}}{\sum}}(T_{\text{ac}i}T_{\text{bc}%
j}+T_{\text{ac}j}T_{\text{bc}i})$

$-\overset{n}{\underset{k=q+1}{\sum}}(\perp_{\text{a}ik}\perp_{\text{b}jk}+$
$\perp_{\text{a}jk}\perp_{\text{b}ik})](y_{0})T_{\text{ab}i}(y_{0})\Lambda
_{j}(y_{0})\phi(y_{0})-\frac{1}{12}T_{\text{ab}i}^{2}(y_{0})\Omega_{ij}%
(y_{0})\phi(y_{0}$

$+\frac{1}{12}\underset{i,j=q+1}{\overset{n}{\sum}}\underset{\text{a,b=1}%
}{\overset{q}{\sum}}[\perp_{\text{a}ij}(\frac{\partial\phi}{\partial
\text{x}_{\text{b}}}+\Lambda_{\text{b}})](y_{0})\qquad\qquad$\ \textbf{I}%
$_{32632}$

$\times\lbrack-R_{\text{a}i\text{b}j}-R_{\text{a}j\text{b}i}%
+\underset{\text{c=1}}{\overset{\text{q}}{\sum}}T_{\text{ac}i}T_{\text{bc}%
j}-3\underset{\text{c=1}}{\overset{\text{q}}{\sum}}T_{\text{ac}j}%
T_{\text{bc}i}+\underset{\text{k=q+1}}{\overset{\text{n}}{\sum}}%
\perp_{\text{a}i\text{k}}\perp_{\text{b}j\text{k}}-\underset{\text{k=q+1}%
}{\overset{\text{n}}{\sum}}\perp_{\text{a}j\text{k}}\perp_{\text{b}i\text{k}%
}\ ](y_{0})\phi(y_{0})$

$-\frac{1}{6}\underset{i,j=q+1}{\overset{n}{\sum}}\underset{\text{a,b=1}%
}{\overset{q}{\sum}}T_{\text{ab}j}(y_{0})\perp_{\text{a}ij}(y_{0}%
)\frac{\partial\Lambda_{\text{b}}}{\partial x_{i}}(y_{0})\phi(y_{0})$

$-\frac{1}{6}\underset{i,j,k=q+1}{\overset{n}{\sum}}\underset{\text{a}%
=1}{\overset{q}{\sum}}\perp_{\text{a}ij}(y_{0})[\overset{q}{\underset{\text{b}%
=1}{\sum(}}\perp_{\text{b}ik}T_{\text{ab}j})(y_{0})+\frac{2}{3}(2R_{\text{a}%
ijk}+R_{\text{a}jik}+R_{\text{a}kji})](y_{0})\Lambda_{k}(y_{0})\phi(y_{0})$

$+\frac{1}{6}\underset{i,j,k=q+1}{\overset{n}{\sum}}\underset{\text{a}%
=1}{\overset{q}{\sum}}\perp_{\text{a}ij}(y_{0})\perp_{\text{a}jk}(y_{0}%
)\Omega_{ik}(y_{0})\phi(y_{0})$

$+\ \frac{1}{24}\underset{i\text{=q+1}}{\overset{\text{n}}{\sum}}%
\frac{\partial^{2}\text{W}}{\partial x_{i}^{2}}(y_{0})\phi(y_{0})\qquad\qquad
$\ \textbf{I}$_{327}$

$+\frac{1}{24}\underset{i,j\text{=q+1}}{\overset{\text{n}}{\sum}%
}\underset{\text{a=1}}{\overset{\text{q}}{\sum}}<H,j>[4\underset{\text{c}%
=1}{\overset{q}{\sum}}(T_{\text{ac}i})(\perp_{j\text{c}i})+\frac{8}%
{3}R_{i\text{a}ij}](y_{0})\frac{\partial\phi}{\partial\text{x}_{\text{a}}%
}(y_{0})\qquad$I$_{328}$

$+\frac{1}{24}\underset{i,j\text{=q+1}}{\overset{\text{n}}{\sum}%
}\underset{\text{a=1}}{\overset{\text{q}}{\sum}}\perp_{\text{a}ji}%
(y_{0})[<H,i><H,j>](y_{0})\frac{\partial\phi}{\partial\text{x}_{\text{a}}%
}(y_{0})$\qquad$\ \ \ \ +\frac{1}{72}\underset{i,j=q+1}{\overset{n}{\sum}%
}\underset{\text{a=1}}{\overset{\text{q}}{\sum}}\perp_{\text{a}ji}%
(y_{0})[2\varrho_{ij}+4\overset{q}{\underset{\text{a}=1}{\sum}}R_{i\text{a}%
j\text{a}}-6\overset{q}{\underset{\text{b,c}=1}{\sum}}T_{\text{cc}%
i}T_{\text{bb}j}-T_{\text{bc}i}T_{\text{bc}j}](y_{0})\frac{\partial\phi
}{\partial\text{x}_{\text{a}}}(y_{0})$

$+\frac{8}{3}R_{j\text{a}ji}(y_{0})X_{i}(y_{0})+[2X_{j}\frac{\partial X_{j}%
}{\partial x_{\text{a}}}-\frac{\partial^{2}X_{j}}{\partial x_{\text{a}%
}\partial x_{j}}](y_{0})\frac{\partial\phi}{\partial\text{x}_{\text{a}}}%
(y_{0})$

$+\frac{1}{24}\underset{i,j=q+1}{\overset{n}{\sum}}\underset{\text{a=1}%
}{\overset{\text{q}}{\sum}}<H,j>(y_{0})[4\underset{\text{c}%
=1}{\overset{q}{\sum}}(T_{\text{ac}i})(\perp_{j\text{c}i})+\frac{8}%
{3}R_{i\text{a}ij}](y_{0})\Lambda_{\text{a}}(y_{0})\phi(y_{0})\qquad$%
I$_{329}\qquad$I$_{3291}\qquad$I$_{32911}\qquad$

$+\frac{1}{12}\underset{i,j=q+1}{\overset{n}{\sum}}\underset{\text{a=1}%
}{\overset{\text{q}}{\sum}}\perp_{\text{a}ji}(y_{0})[<H,i><H,j>](y_{0}%
)\Lambda_{\text{a}}(y_{0})\phi(y_{0})$

$+\frac{1}{72}\underset{i,j=q+1}{\overset{n}{\sum}}\underset{\text{a=1}%
}{\overset{\text{q}}{\sum}}\perp_{\text{a}ji}(y_{0})[2\varrho_{ij}%
+4\overset{q}{\underset{\text{a}=1}{\sum}}R_{i\text{a}j\text{a}}%
-6\overset{q}{\underset{\text{b,c}=1}{\sum}}T_{\text{cc}i}T_{\text{bb}%
j}-T_{\text{bc}i}T_{\text{bc}j}](y_{0})\Lambda_{\text{a}}(y_{0})\phi(y_{0})$

$+\frac{1}{36}\underset{i,j=q+1}{\overset{n}{\sum}}%
\underset{k=q+1}{\overset{n}{\sum}}<H,k>(y_{0})$R$_{ijik}(y_{0})\Lambda
_{j}(y_{0})\phi(y_{0})$

$-\frac{1}{288}\underset{i,j=q+1}{\overset{n}{\sum}}<H,j>(y_{0})[3<H,i>^{2}%
+2(\tau^{M}-3\tau^{P}+\overset{q}{\underset{\text{a}=1}{\sum}}\varrho
_{\text{aa}}+\overset{q}{\underset{\text{a,b}=1}{\sum}}R_{\text{abab}}%
)](y_{0})\Lambda_{j}(y_{0})\phi(y_{0})$

$-\frac{1}{12}\underset{i,j=q+1}{\overset{n}{\sum}}<H,i>(y_{0})$

$\times\lbrack\frac{3}{4}<H,i><H,j>$\ $+\frac{1}{6}(\varrho_{ij}%
+2\overset{q}{\underset{\text{a}=1}{\sum}}R_{i\text{a}j\text{a}}%
-3\overset{q}{\underset{\text{a,b=1}}{\sum}}T_{\text{aa}i}T_{\text{bb}%
j}-T_{\text{ab}i}T_{\text{ab}j})](y_{0})\Lambda_{j}(y_{0})\phi(y_{0})$

$+\frac{5}{32}\underset{i,j=q+1}{\overset{n}{\sum}}$%
$<$%
H,$i$%
$>$%
$^{2}$%
$<$%
H,$j$%
$>$%
$\Lambda_{j}(y_{0})\phi(y_{0})$

$+\frac{1}{48}\underset{i,j=q+1}{\overset{n}{\sum}}$%
$<$%
H,$i$%
$>$%
(y$_{0}$)$[$(2$\varrho_{ij}$+4$\overset{q}{\underset{\text{a}=1}{\sum}}%
$R$_{i\text{a}j\text{a}}$-3$\overset{q}{\underset{\text{a,b=1}}{\sum}}%
$T$_{\text{aa}i}$T$_{\text{bb}j}$-T$_{\text{ab}i}$T$_{\text{ab}j}%
$-3$\overset{q}{\underset{\text{a,b=1}}{\sum}}$T$_{\text{aa}j}$T$_{\text{bb}%
i}$-T$_{\text{ab}j}$T$_{\text{ab}i}$)$](y_{0})\Lambda_{j}(y_{0})\phi(y_{0})$

$+\frac{1}{48}\underset{i,j=q+1}{\overset{n}{\sum}}<H,j>(y_{0})[$ $\tau
^{M}-3\tau^{P}+\overset{q}{\underset{\text{a}=1}{\sum}}\varrho_{\text{aa}%
}+\overset{q}{\underset{\text{a,b}=1}{\sum}}R_{\text{abab}}](y_{0})\Lambda
_{j}(y_{0})\phi(y_{0})$

$+\frac{1}{144}\underset{i,j=q+1}{\overset{n}{\sum}}[\nabla_{i}\varrho
_{ij}-2\varrho_{ij}<H,i>+\overset{q}{\underset{\text{a}=1}{\sum}}(\nabla
_{i}R_{\text{a}i\text{a}j}-4R_{i\text{a}j\text{a}}<H,i>)$

$+4\overset{q}{\underset{\text{a,b=1}}{\sum}}R_{i\text{a}j\text{b}%
}T_{\text{ab}i}+2\overset{q}{\underset{\text{a,b,c=1}}{\sum}}(T_{\text{aa}%
i}T_{\text{bb}j}T_{\text{cc}i}-3T_{\text{aa}i}T_{\text{bc}j}T_{\text{bc}%
i}+2T_{\text{ab}i}T_{\text{bc}j}T_{\text{ca}i})](y_{0})\Lambda_{j}(y_{0}%
)\phi(y_{0})$\qquad\qquad\qquad\qquad\qquad\ \ 

$+\frac{1}{144}\underset{i,j=q+1}{\overset{n}{\sum}}[\nabla_{j}\varrho
_{ii}-2\varrho_{ji}<H,i>+\overset{q}{\underset{\text{a}=1}{\sum}}(\nabla
_{j}R_{\text{a}i\text{a}i}-4R_{j\text{a}i\text{a}}<H,i>)$

$+4\overset{q}{\underset{\text{a,b=1}}{\sum}}R_{j\text{a}i\text{b}%
}T_{\text{ab}i}+2\overset{q}{\underset{\text{a,b,c=1}}{\sum}}(T_{\text{aa}%
j}T_{\text{bb}i}T_{\text{cc}i}-3T_{\text{aa}j}T_{\text{bc}i}T_{\text{bc}%
i}+2T_{\text{ab}j}T_{\text{bc}i}T_{\text{ca}i})](y_{0})\Lambda_{j}(y_{0}%
)\phi(y_{0})$

$+\frac{1}{144}\underset{i,j=q+1}{\overset{n}{\sum}}[\nabla_{i}\varrho
_{ij}-2\varrho_{ii}<H,j>+\overset{q}{\underset{\text{a}=1}{\sum}}(\nabla
_{i}R_{\text{a}i\text{a}j}-4R_{i\text{a}i\text{a}}<H,j>)$

$+4\overset{q}{\underset{\text{a,b=1}}{\sum}}R_{i\text{a}i\text{b}%
}T_{\text{ab}j}+2\overset{q}{\underset{\text{a,b,c=1}}{\sum}}(T_{\text{aa}%
i}T_{\text{bb}i}T_{\text{cc}j}-3T_{\text{aa}i}T_{\text{bc}i}T_{\text{bc}%
j}+2T_{\text{ab}i}T_{\text{bc}i}T_{\text{ac}j})](y_{0})\Lambda_{j}(y_{0}%
)\phi(y_{0})$

$+\frac{1}{72}\underset{i,j=q+1}{\overset{n}{\sum}}$ $<H,j>(y_{0}%
)\frac{\partial\Omega_{ij}}{\partial\text{x}_{i}}(y_{0})\phi(y_{0}%
)\qquad\qquad$I$_{32912}$

$-\frac{1}{12}\underset{i=q+1}{\overset{n}{\sum}}\underset{\text{a=1}%
}{\overset{\text{q}}{\sum}}\perp_{\text{a}ij}(y_{0})<H,j>(y_{0})[\Omega
_{i\text{a}}+[\Lambda_{\text{a}},\Lambda_{i}]](y_{0})\phi(y_{0})\qquad\qquad
$I$_{32913}$

$+$ $\frac{1}{72}\underset{i,j=q+1}{\overset{n}{\sum}}[3<H,i><H,j>\ +(\varrho
_{ij}+2\overset{q}{\underset{\text{a}=1}{\sum}}R_{i\text{a}j\text{a}%
}-3\overset{q}{\underset{\text{a,b}=1}{\sum}}T_{\text{aa}i}T_{\text{bb}%
j}-T_{\text{ab}i}T_{\text{ab}j})](y_{0})\Omega_{ij}(y_{0})\phi(y_{0})$

$+\frac{1}{12}\underset{i=q+1}{\overset{n}{\sum}}$ $\underset{\text{a=1}%
}{\overset{\text{q}}{\sum}}$ $[-4\underset{\text{b=1}}{\overset{\text{q}}{%
{\textstyle\sum}
}}T_{\text{ab}i}\frac{\partial X_{i}}{\partial x_{\text{b}}}%
+\underset{j=q+1}{\overset{n}{%
{\textstyle\sum}
}}\perp_{\text{a}ij}\left(  \frac{\partial X_{i}}{\partial x_{j}}%
+\frac{\partial X_{j}}{\partial x_{i}}\right)  \qquad\qquad$I$_{3292}\qquad
$I$_{32921}$

$+\frac{8}{3}\underset{j=q+1}{\overset{n}{%
{\textstyle\sum}
}}R_{i\text{a}ij}X_{j}+\left(  2X_{i}\frac{\partial X_{i}}{\partial
x_{\text{a}}}-\frac{\partial^{2}X_{i}}{\partial x_{\text{a}}\partial x_{i}%
}\right)  ](y_{0})\Lambda_{\text{a}}(y_{0})\phi(y_{0}$

$+\frac{1}{6}$ $\underset{i=q+1}{\overset{n}{\sum}}\underset{\text{a=1}%
}{\overset{\text{q}}{\sum}}$ $[\underset{j=q+1}{\overset{n}{%
{\textstyle\sum}
}}X_{j}\perp_{\text{a}ij}+\frac{\partial X_{i}}{\partial x_{\text{a}}}]$
$(y_{0})[\Omega_{i\text{a}}+[\Lambda_{\text{a}},\Lambda_{i}](y_{0})\phi
(y_{0})$

$-\frac{1}{36}[\left(  2\frac{\partial^{2}X_{i}}{\partial x_{i}\partial x_{j}%
}+\frac{\partial^{2}X_{j}}{\partial x_{i}^{2}}\right)
+2\underset{k=q+1}{\overset{n}{\sum}}R_{ijik}X_{k}](y_{0})\Lambda_{j}%
(y_{0})\phi(y_{0})\qquad$I$_{32922}$

$-\frac{1}{36}X_{j}(y_{0})$ $\frac{\partial\Omega_{ij}}{\partial x_{i}}%
(y_{0})\phi(y_{0})-\frac{1}{12}\frac{\partial X_{j}}{\partial x_{i}}%
(y_{0})\Omega_{ij}(y_{0})$ $\phi(y_{0})$

$+\frac{1}{12}$ $\underset{\text{a=1}}{\overset{\text{q}}{\sum}}\frac
{\partial^{2}X_{\text{a}}}{\partial x_{i}^{2}}(y_{0})\frac{\partial\phi
}{\partial\text{x}_{\text{a}}}(y_{0})\qquad$\ \ \textbf{L}$_{1}$\qquad\qquad

$+\frac{1}{12}$ $\underset{\text{a}=1}{\overset{\text{q}}{\sum}}$
$\frac{\partial^{2}X_{\text{a}}}{\partial x_{i}^{2}}(y_{0})\Lambda_{\text{a}%
}(y_{0})\phi(y_{0})+\frac{1}{12}\underset{j=q+1}{\overset{\text{n}}{\sum}%
}\frac{\partial^{2}X_{j}}{\partial x_{i}^{2}}(y_{0})\Lambda_{j}(y_{0}%
)\phi(y_{0})\qquad$\textbf{L}$_{2}\qquad$\textbf{L}$_{21}$

$+\frac{1}{36}$ $\underset{j=q+1}{\overset{n}{\sum}}$ $X_{j}(y_{0}%
)\frac{\partial\Omega_{ij}}{\partial x_{i}}(y_{0})\phi(y_{0})\qquad$%
\textbf{L}$_{22}$

$+\frac{1}{12}$ $\underset{j=q+1}{\overset{\text{n}}{\sum}}$ $\frac{\partial
X_{j}}{\partial x_{i}}(y_{0})\Omega_{ij}(y_{0})\phi(y_{0})\qquad$%
\textbf{L}$_{23}$\qquad\qquad\qquad\qquad\qquad\qquad\qquad\qquad\qquad
\qquad\qquad\qquad

\qquad\qquad\qquad\qquad\qquad\qquad\qquad\qquad\qquad\qquad\qquad\qquad
\qquad\qquad\qquad$\qquad\qquad\qquad\blacksquare$\qquad$\qquad$

\subsection{Computation of I$_{33}$}

I$_{33}=\frac{1}{2}\underset{\text{c=1}}{\overset{\text{q}}{\sum}}%
\Lambda_{\text{c}}(y_{0})\frac{\partial\Theta}{\partial x_{\text{c}}}(y_{0})$

The expression for $\Theta$ is given in $\left(  D_{6}\right)  :$

$\qquad\Theta=$ L$_{\Psi}[\phi\circ\pi_{P}]$\textbf{ }

$=\frac{\text{L}\Psi}{\Psi}\phi\circ\pi_{\text{P}}+$\ $\frac{1}{2}%
\underset{\text{a,b=1}}{\overset{\text{q}}{\sum}}$g$^{\text{ab}}\left\{
\frac{\partial^{2}\phi}{\partial\text{x}_{\text{a}}\partial\text{x}_{\text{b}%
}}\circ\pi_{\text{P}}\right\}  +\frac{1}{2}\underset{i,j=1}{\overset{n}{\sum}%
}$g$^{ij}\left\{  \frac{\partial\Lambda_{j}}{\partial\text{x}_{i}}\phi\circ
\pi_{\text{P}}\right\}  $

$\bigskip+$ $\underset{j=1}{\overset{n}{\sum}}\underset{\text{a=1}%
}{\overset{\text{q}}{\sum}}$g$^{\text{a}j}\left\{  \text{ }\Lambda_{j}%
\frac{\partial\phi}{\partial\text{x}_{\text{a}}}\circ\pi_{\text{P}}\right\}
+\frac{1}{2}\underset{i,j=1}{\overset{n}{\sum}}$g$^{ij}\Lambda_{i}\Lambda
_{j}\phi\circ\pi_{\text{P}}$

$-\frac{1}{2}\underset{i,j=1}{\overset{n}{\sum}}$g$^{ij}\left\{
\underset{\text{b=1}}{\overset{\text{q}}{\sum}}\Gamma_{ij}^{\text{b}}%
\frac{\partial\phi}{\partial\text{x}_{\text{b}}}\circ\pi_{\text{P}}\text{
}+\text{ }\underset{k=1}{\overset{n}{\sum}}\Gamma_{ij}^{k}(z_{1})\Lambda
_{k}\phi\circ\pi_{\text{P}}\right\}  $

$+\frac{1}{2}$W$\phi\circ\pi_{\text{P}}$

$+$ $\underset{\text{a=1}}{\overset{\text{q}}{\sum}}(\nabla\log\Psi
)_{\text{a}}\frac{\partial\phi}{\partial\text{x}_{\text{a}}}\circ\pi
_{\text{P}}$ + $\underset{\text{j=1}}{\overset{\text{n}}{\sum}}$ $(\nabla
\log\Psi)_{j}\Lambda_{j}\phi\circ\pi_{\text{P}}$

$+$ $\underset{\text{a=1}}{\overset{\text{q}}{\sum}}$X$^{\text{a}}%
\frac{\partial\phi}{\partial\text{x}_{\text{a}}}\circ\pi_{\text{P}}$ +
$\underset{\text{j=1}}{\overset{\text{n}}{\sum}}$ X$^{j}\Lambda_{j}\phi
\circ\pi_{\text{P}}$

Then we set:

I$_{33}=\frac{1}{2}\underset{\text{c=1}}{\overset{\text{q}}{\sum}}%
\Lambda_{\text{c}}(y_{0})\frac{\partial\Theta}{\partial x_{\text{c}}}%
(y_{0})=\frac{1}{2}\underset{\text{c=1}}{\overset{\text{q}}{\sum}}%
\Lambda_{\text{c}}(y_{0})\frac{\partial}{\partial x_{\text{c}}}[L_{\Psi}%
\phi\circ\pi_{P}](y_{0})$

$=$ I$_{331}+$ I$_{332}$+ I$_{333}+$ I$_{334}+$ I$_{335}+$ I$_{336}+$
I$_{337}+$ I$_{338}+$ I$_{338}+$ E$_{1}+$ E$_{2}$

where,

I$_{331}=\frac{1}{2}\underset{\text{c=1}}{\overset{\text{q}}{\sum}}%
\Lambda_{\text{c}}(y_{0})\frac{\partial}{\partial x_{\text{c}}}[\frac
{\text{L}\Psi}{\Psi}\phi\circ\pi_{P}](y_{0})$

I$_{332}=\frac{1}{4}\underset{\text{c=1}}{\overset{\text{q}}{\sum}}%
\Lambda_{\text{c}}(y_{0})\frac{\partial}{\partial x_{\text{c}}}%
[\underset{\text{a,b=1}}{\overset{\text{q}}{\sum}}$g$^{\text{ab}}\left\{
\frac{\partial^{2}\phi}{\partial\text{x}_{\text{a}}\partial\text{x}_{\text{b}%
}}\circ\pi_{P}\text{ }\right\}  ](y_{0})$

I$_{333}=\frac{1}{4}\underset{\text{c=1}}{\overset{\text{q}}{\sum}}%
\Lambda_{\text{c}}(y_{0})\frac{\partial}{\partial x_{\text{c}}}%
[\underset{i,j=1}{\overset{n}{\sum}}g^{ij}\left\{  \frac{\partial\Lambda_{j}%
}{\partial\text{x}_{i}}\phi\circ\pi_{P}\right\}  ](y_{0})\qquad$

I$_{334}$ $=\frac{1}{2}\underset{\text{c=1}}{\overset{\text{q}}{\sum}}%
\Lambda_{\text{c}}(y_{0})\frac{\partial}{\partial x_{\text{c}}}%
\underset{j=1}{\overset{n}{\sum}}\underset{\text{a=1}}{\overset{\text{q}%
}{\sum}}$g$^{\text{a}j}\left\{  \text{ }\Lambda_{j}\frac{\partial\phi
}{\partial\text{x}_{\text{a}}}\circ\pi_{P}\right\}  (y_{0})$

I$_{335}=$ $\frac{1}{4}\underset{\text{c=1}}{\overset{\text{q}}{\sum}}%
\Lambda_{\text{c}}(y_{0})\frac{\partial}{\partial x_{\text{c}}}%
[\underset{i,j=1}{\overset{n}{\sum}}$g$^{ij}\left\{  \Lambda_{i}\Lambda
_{j}\phi\circ\pi_{P}\right\}  ](y_{0})$

I$_{336}=-\frac{1}{4}\underset{\text{c=1}}{\overset{\text{q}}{\sum}}%
\Lambda_{\text{c}}(y_{0})\frac{\partial}{\partial x_{\text{c}}}%
[\underset{i,j=1}{\overset{n}{\sum}}g^{ij}\left\{  \Gamma_{ij}^{\text{b}}%
\frac{\partial\phi}{\partial\text{x}_{\text{b}}}\circ\pi_{P}+\text{
}\underset{k=1}{\overset{n}{\sum}}\Gamma_{ij}^{k}\Lambda_{k}\phi\circ\pi
_{P}\right\}  ](y_{0})$

I$_{337}=\frac{1}{2}\underset{\text{c=1}}{\overset{\text{q}}{\sum}}%
\Lambda_{\text{c}}(y_{0})\frac{\partial}{\partial x_{\text{c}}}[$
$\underset{\text{a=1}}{\overset{\text{q}}{\sum}}(\nabla\log\Psi)_{\text{a}%
}\frac{\partial\phi}{\partial\text{x}_{\text{a}}}\circ\pi_{P}](y_{0})$

I$_{338}=\frac{1}{2}\underset{\text{c=1}}{\overset{\text{q}}{\sum}}%
\Lambda_{\text{c}}(y_{0})\frac{\partial}{\partial x_{\text{a}}}%
[\underset{j=1}{\overset{n}{\sum}}(\nabla\log\Psi)_{j}\Lambda_{j}\phi\circ
\pi_{P}](y_{0})$

I$_{339}=\frac{1}{4}\underset{\text{c=1}}{\overset{\text{q}}{\sum}}%
\Lambda_{\text{c}}(y_{0})\frac{\partial}{\partial x_{\text{a}}}[$W$\phi
\circ\pi_{P}](y_{0})$

E$_{1}=$ $\underset{\text{a=1}}{\overset{\text{q}}{\sum}}\Lambda_{\text{c}%
}(y_{0})\frac{\partial}{\partial x_{\text{a}}}[$X$_{\text{a}}\frac
{\partial\phi}{\partial\text{x}_{\text{a}}}\circ\pi_{\text{P}}](y_{0})$

E$_{2}=$ $\underset{\text{j=1}}{\overset{\text{n}}{\sum}}$ $\Lambda_{\text{c}%
}(y_{0})\frac{\partial}{\partial x_{\text{a}}}[$X$_{j}\Lambda_{j}\phi\circ
\pi_{\text{P}}](y_{0})$

We start the computation of I$_{33}$ here:

I$_{331}=\frac{1}{2}\underset{\text{c=1}}{\overset{\text{q}}{\sum}}%
\Lambda_{\text{c}}(y_{0})\frac{\partial}{\partial x_{\text{c}}}[\frac
{\text{L}\Psi}{\Psi}\phi\circ\pi_{P}](y_{0})$

\qquad$=\frac{1}{2}\underset{\text{c=1}}{\overset{\text{q}}{\sum}}%
\Lambda_{\text{c}}(y_{0})\frac{\partial}{\partial x_{\text{c}}}\frac
{\text{L}\Psi}{\Psi}(y_{0})\phi(y_{0})+\frac{1}{2}\underset{\text{c=1}%
}{\overset{\text{q}}{\sum}}\Lambda_{\text{c}}(y_{0})\frac{\text{L}\Psi}{\Psi
}(y_{0})\frac{\partial\phi}{\partial x_{\text{c}}}(y_{0})$

\qquad$=\frac{1}{2}\underset{\text{c=1}}{\overset{\text{q}}{\sum}}%
\Lambda_{\text{c}}(y_{0})\frac{\text{L}\Psi}{\Psi}(y_{0})\frac{\partial\phi
}{\partial x_{\text{c}}}(y_{0})+\frac{1}{2}\underset{\text{c=1}%
}{\overset{\text{q}}{\sum}}\Lambda_{\text{c}}(y_{0})\frac{\partial}{\partial
x_{\text{c}}}\frac{\text{L}\Psi}{\Psi}(y_{0})\phi(y_{0})$

The expression of $\frac{\text{L}\Psi}{\Psi}(y_{0})$ is given in $\left(
10.30\right)  $ of \textbf{Chapter 10} and that of $\frac{\partial}{\partial
x_{\text{c}}}[\frac{\text{L}\Psi}{\Psi}](y_{0})$ is given in (v) of
\textbf{Table B}$_{5}:$

\qquad I$_{331}=\frac{1}{2}\underset{\text{c=1}}{\overset{\text{q}}{\sum}%
}\Lambda_{\text{c}}(y_{0})\frac{\partial}{\partial x_{\text{c}}}%
[\frac{\text{L}\Psi}{\Psi}\phi\circ\pi_{P}](y_{0})$

\qquad\qquad$=\frac{1}{2}\underset{\text{c=1}}{\overset{\text{q}}{\sum}%
}\Lambda_{\text{c}}(y_{0})\frac{\text{L}\Psi}{\Psi}(y_{0})\frac{\partial\phi
}{\partial x_{\text{c}}}(y_{0})+\frac{1}{2}\underset{\text{c=1}%
}{\overset{\text{q}}{\sum}}\Lambda_{\text{c}}(y_{0})\frac{\partial}{\partial
x_{\text{c}}}\frac{\text{L}\Psi}{\Psi}(y_{0})\phi(y_{0})$

$\left(  D_{41}\right)  $ \qquad I$_{331}=\frac{1}{48}$ $\underset{\text{c=1}%
}{\overset{\text{q}}{\sum}}\Lambda_{\text{c}}(y_{0}%
)[\underset{i=q+1}{\overset{n}{\sum}}3<H,i>^{2}+2(\tau^{M}-3\tau
^{P}\ +\overset{q}{\underset{\text{a=1}}{\sum}}\varrho_{\text{aa}}%
^{M}+\overset{q}{\underset{\text{a,b}=1}{\sum}}R_{\text{abab}}^{M}%
)](y_{0})\frac{\partial\phi}{\partial x_{\text{c}}}(y_{0})$

\qquad$-\frac{1}{4}\underset{\text{c=1}}{\overset{\text{q}}{\sum}}%
\Lambda_{\text{c}}(y_{0})[\left\Vert \text{X}\right\Vert ^{2}+$ divX $-$
$\underset{\text{a}=1}{\overset{q}{\sum}}($X$_{\text{a}})^{2}$ $-$
$\underset{\text{a}=1}{\overset{q}{\sum}}\frac{\partial X_{\text{a}}}{\partial
x_{\text{a}}}](y_{0})\frac{\partial\phi}{\partial x_{\text{c}}}(y_{0})+$
$\frac{1}{2}\underset{\text{c=1}}{\overset{\text{q}}{\sum}}\Lambda_{\text{c}%
}(y_{0})$V$(y_{0})\frac{\partial\phi}{\partial x_{\text{c}}}(y_{0})$

\qquad$+\frac{1}{2}\underset{\text{c=1}}{\overset{\text{q}}{\sum}}%
\Lambda_{\text{c}}(y_{0})[-(X_{j}\frac{\partial X_{j}}{\partial x_{\text{c}}%
}+\frac{1}{2}\frac{\partial^{2}X_{j}}{\partial x_{\text{c}}\partial x_{j}%
})(y_{0})+\frac{1}{2}(<H,j>\frac{\partial X_{j}}{\partial x_{\text{c}}}%
)(y_{0})+\frac{\partial\text{V}}{\partial x_{\text{c}}}(y_{0})]\phi(y_{0})$

Next we have:

I$_{332}=\frac{1}{4}\underset{\text{c=1}}{\overset{\text{q}}{\sum}}%
\Lambda_{\text{c}}(y_{0})\frac{\partial}{\partial x_{\text{c}}}%
[\underset{\text{a,b=1}}{\overset{\text{q}}{\sum}}$g$^{\text{ab}}\left\{
\frac{\partial^{2}\phi}{\partial\text{x}_{\text{a}}\partial\text{x}_{\text{b}%
}}\circ\pi_{P}\text{ }\right\}  ](y_{0})$

$\qquad=\frac{1}{4}\underset{\text{c=1}}{\overset{\text{q}}{\sum}}%
\Lambda_{\text{c}}(y_{0})\underset{\text{a,b=1}}{\overset{\text{q}}{\sum}}%
$g$^{\text{ab}}(y_{0})\left\{  \frac{\partial^{3}\phi}{\partial\text{x}%
_{\text{a}}\partial\text{x}_{\text{b}}\partial x_{\text{c}}}\text{ }\right\}
(y_{0})$\qquad

$\left(  D_{42}\right)  \qquad$I$_{332}=\frac{1}{4}\underset{\text{a,c=1}%
}{\overset{\text{q}}{\sum}}\Lambda_{\text{c}}(y_{0})\left\{  \frac
{\partial^{3}\phi}{\partial\text{x}_{\text{a}}^{2}\partial x_{\text{c}}}\text{
}\right\}  (y_{0})$

Next we have:

I$_{333}=\frac{1}{4}\underset{\text{c=1}}{\overset{\text{q}}{\sum}}%
\Lambda_{\text{c}}(y_{0})\frac{\partial}{\partial x_{\text{c}}}%
[\underset{i,j=1}{\overset{n}{\sum}}g^{ij}\left\{  \frac{\partial\Lambda_{j}%
}{\partial\text{x}_{i}}\phi\circ\pi_{P}\right\}  ](y_{0})$

\qquad$=\frac{1}{4}\underset{\text{c=1}}{\overset{\text{q}}{\sum}}%
\Lambda_{\text{c}}(y_{0})[\underset{i,j=1}{\overset{n}{\sum}}$g$^{ij}%
(y_{0})\left\{  \frac{\partial\Lambda_{j}}{\partial\text{x}_{i}}(y_{0}%
)\frac{\partial\phi}{\partial x_{\text{c}}}\circ\pi_{P}\right\}  ](y_{0})$

$\qquad=\frac{1}{4}\underset{\text{c=1}}{\overset{\text{q}}{\sum}}%
\Lambda_{\text{c}}(y_{0})[\underset{i=1}{\overset{n}{\sum}}\left\{
\frac{\partial\Lambda_{i}}{\partial\text{x}_{i}}(y_{0})\frac{\partial\phi
}{\partial x_{\text{c}}}\right\}  ](y_{0})$

Since $\frac{\partial\Lambda_{i}}{\partial\text{x}_{i}}(y_{0})=0$, we have:

$\left(  D_{43}\right)  \qquad$I$_{333}=0$

Next we have:

I$_{334}$ $=\frac{1}{2}\underset{\text{c=1}}{\overset{\text{q}}{\sum}}%
\Lambda_{\text{c}}(y_{0})\frac{\partial}{\partial x_{\text{c}}}%
\underset{j=1}{\overset{n}{\sum}}\underset{\text{a=1}}{\overset{\text{q}%
}{\sum}}$g$^{\text{a}j}\left\{  \text{ }\Lambda_{j}\frac{\partial\phi
}{\partial\text{x}_{\text{a}}}\circ\pi_{P}\right\}  (y_{0})$

\qquad$=\frac{1}{2}\underset{\text{c=1}}{\overset{\text{q}}{\sum}}%
\Lambda_{\text{c}}(y_{0})\underset{j=1}{\overset{n}{\sum}}\underset{\text{a=1}%
}{\overset{\text{q}}{\sum}}$g$^{\text{a}j}(y_{0})\left\{  \text{ }\Lambda
_{j}\frac{\partial^{2}\phi}{\partial\text{x}_{\text{a}}\partial x_{\text{c}}%
}\circ\pi_{P}\right\}  (y_{0})$

Since g$^{\text{a}j}(y_{0})=\delta^{\text{a}j}$ for a = 1,...,q and
$j=1,...,q,q+1,...,n,$ we have:

I$_{334}=\frac{1}{2}\underset{\text{c=1}}{\overset{\text{q}}{\sum}}%
\Lambda_{\text{c}}(y_{0})\underset{\text{a=1}}{\overset{\text{q}}{\sum}%
}\left\{  \Lambda_{\text{a}}(y_{0})\frac{\partial^{2}\phi}{\partial
\text{x}_{\text{a}}\partial x_{\text{c}}}\text{ }\right\}  (y_{0})=\frac{1}%
{2}\underset{\text{a,c=1}}{\overset{\text{q}}{\sum}}\left\{  \Lambda
_{\text{a}}(y_{0})\Lambda_{\text{c}}(y_{0})\frac{\partial^{2}\phi}%
{\partial\text{x}_{\text{a}}\partial x_{\text{c}}}\text{ }\right\}  (y_{0})$

$\left(  D_{44}\right)  \qquad$I$_{334}=\frac{1}{2}\underset{\text{a,b=1}%
}{\overset{\text{q}}{\sum}}\left\{  \Lambda_{\text{a}}(y_{0})\Lambda
_{\text{b}}(y_{0})\frac{\partial^{2}\phi}{\partial\text{x}_{\text{a}}\partial
x_{\text{b}}}\text{ }\right\}  (y_{0})$ \qquad

We remind that differentiation of g$^{ij}$ and $\Lambda_{i}$with respect to
tangential coordinates vanish, and so:

I$_{335}=\frac{1}{4}\underset{\text{c=1}}{\overset{\text{q}}{\sum}}%
\Lambda_{\text{c}}(y_{0})\frac{\partial}{\partial x_{\text{c}}}%
[\underset{i,j=1}{\overset{n}{\sum}}$g$^{ij}\Lambda_{i}\Lambda_{j}\phi\circ
\pi_{P}](y_{0})$

\qquad$=\frac{1}{4}\underset{\text{c=1}}{\overset{\text{q}}{\sum}}%
\Lambda_{\text{c}}(y_{0})[\underset{i,j=1}{\overset{n}{\sum}}$g$^{ij}%
(y_{0})\Lambda_{i}\Lambda_{j}\frac{\partial\phi}{\partial x_{\text{c}}}%
\circ\pi_{P}](y_{0})$

Since g$^{ij}(y_{0})=\delta^{ij},$ we have:

$\left(  D_{45}\right)  \qquad$I$_{335}=\frac{1}{4}\underset{\text{c=1}%
}{\overset{\text{q}}{\sum}}[\underset{i=1}{\overset{n}{\sum}}\Lambda
_{\text{c}}(y_{0})\Lambda_{i}^{2}(y_{0})\frac{\partial\phi}{\partial
x_{\text{c}}}(y_{0})]$

We next compute:

I$_{336}=-\frac{1}{4}\underset{\text{c=1}}{\overset{\text{q}}{\sum}}%
\Lambda_{\text{c}}(y_{0})\frac{\partial}{\partial x_{\text{c}}}%
[\underset{i,j=1}{\overset{n}{\sum}}g^{ij}\left\{  \Gamma_{ij}^{\text{b}}%
\frac{\partial\phi}{\partial\text{x}_{\text{b}}}\circ\pi_{P}\text{
}+\underset{k=1}{\overset{n}{\sum}}\Gamma_{ij}^{k}\Lambda_{k}\phi\circ\pi
_{P}\right\}  ](y_{0})$

$\qquad=-\frac{1}{4}\underset{\text{c=1}}{\overset{\text{q}}{\sum}}%
\Lambda_{\text{c}}(y_{0})[\underset{i,j=1}{\overset{n}{\sum}}$g$^{ij}\left\{
\Gamma_{ij}^{\text{b}}\frac{\partial^{2}\phi}{\partial x_{\text{b}}%
\partial\text{x}_{\text{c}}}\circ\pi_{P}+\text{ }%
\underset{k=1}{\overset{n}{\sum}}\Gamma_{ij}^{k}\Lambda_{k}\frac{\partial\phi
}{\partial x_{\text{c}}}\right\}  ](y_{0})$

I$_{336}=-\frac{1}{4}\underset{\text{c=1}}{\overset{\text{q}}{\sum}}%
\Lambda_{\text{c}}(y_{0})\underset{i=1}{\overset{n}{\sum}}\left\{  \Gamma
_{ii}^{\text{b}}\frac{\partial^{2}\phi}{\partial x_{\text{b}}\partial
\text{x}_{\text{c}}}\circ\pi_{P}+\text{ }\underset{k=1}{\overset{n}{\sum}%
}\Gamma_{ii}^{k}\Lambda_{k}\frac{\partial\phi}{\partial x_{\text{c}}}\right\}
(y_{0})$\qquad

Since $\Gamma_{\text{aa}}^{\text{b}}(y_{0})=0$ for a,b = 1,...,q and
$\Gamma_{ii}^{\text{b}}(y_{0})=0=\Gamma_{ii}^{j}(y_{0})$ for b = 1,...,q and
$i,j=q+1,...,n,$

I$_{336}=-\frac{1}{4}\underset{\text{c=1}}{\overset{\text{q}}{\sum}}%
\Lambda_{\text{c}}(y_{0})\underset{j=q+1}{\overset{n}{\sum}}%
\underset{\text{a=1}}{\overset{\text{q}}{\sum}}\left\{  \Gamma_{\text{aa}}%
^{j}\Lambda_{j}\frac{\partial\phi}{\partial x_{\text{c}}}\right\}  (y_{0})$

\qquad\qquad I$_{336}=-\frac{1}{4}\underset{j=q+1}{\overset{n}{\sum}%
}\underset{\text{a,c=1}}{\overset{\text{q}}{\sum}}(y_{0})\left\{
\Gamma_{\text{aa}}^{j}\Lambda_{\text{c}}\Lambda_{j}\frac{\partial\phi
}{\partial x_{\text{c}}}\right\}  (y_{0})$

Since $\Gamma_{\text{aa}}^{j}(y_{0})=T_{\text{aa}j}(y_{0})$ by (i) of
\textbf{Table A}$_{7},$

$\left(  D_{46}\right)  $\qquad I$_{336}=-\frac{1}{4}%
\underset{j=q+1}{\overset{n}{\sum}}\underset{\text{a,c=1}}{\overset{\text{q}%
}{\sum}}[\Lambda_{\text{c}}\Lambda_{j}T_{\text{aa}j}(y_{0})\frac{\partial\phi
}{\partial x_{\text{c}}}](y_{0})$

I$_{337}=\frac{1}{2}\underset{\text{c=1}}{\overset{\text{q}}{\sum}}%
\Lambda_{\text{c}}(y_{0})\frac{\partial}{\partial x_{\text{c}}}[$
$\underset{\text{a=1}}{\overset{\text{q}}{\sum}}(\nabla\log\Psi)_{\text{a}%
}\frac{\partial\phi}{\partial\text{x}_{\text{a}}}\circ\pi_{P}](y_{0})$

We have:

$(\nabla\log\theta^{-\frac{1}{2}})_{\text{a}}(y_{0})=0=(\nabla\log
\Phi)_{\text{a}}(y_{0})$ and $\frac{\partial}{\partial x_{\text{c}}}%
(\nabla\log\theta^{-\frac{1}{2}})_{\text{a}}(y_{0})=0=\frac{\partial}{\partial
x_{\text{c}}}(\nabla\log\Phi)_{\text{a}}(y_{0}).$Therefore,

$\left(  D_{47}\right)  \qquad$I$_{337}=0$

We next consider:

I$_{338}=$ $\frac{1}{2}\underset{\text{c=1}}{\overset{\text{q}}{\sum}}%
\Lambda_{\text{c}}(y_{0})\frac{\partial}{\partial x_{\text{c}}}[$
$\underset{j=1}{\overset{n}{\sum}}$ $(\nabla\log\Psi)_{j}\Lambda_{j}\phi
\circ\pi_{P}](y_{0})$

\qquad$=\frac{1}{2}\underset{\text{c=1}}{\overset{\text{q}}{\sum}}%
\Lambda_{\text{c}}(y_{0})[\underset{j=1}{\overset{n}{\sum}}\frac{\partial
}{\partial x_{\text{c}}}(\nabla\log\theta^{-\frac{1}{2}})_{j}(y_{0}%
)+\frac{\partial}{\partial x_{\text{c}}}(\nabla\log\Phi)_{j}](y_{0}%
)\Lambda_{j}(y_{0})\phi(y_{0})$

Since $\frac{\partial}{\partial x_{\text{c}}}(\nabla\log\theta^{-\frac{1}{2}%
})_{j}(y_{0})=0=\frac{\partial}{\partial x_{\text{c}}}(\nabla\log
\Phi)_{\text{a}}(y_{0})$ for a,c = 1,...,q and $j=1,...,q,q+1,...,n,$

\qquad I$_{338}=\frac{1}{2}\underset{\text{c=1}}{\overset{\text{q}}{\sum}%
}\Lambda_{\text{c}}(y_{0})[\underset{j=q+1}{\overset{n}{\sum}}\frac{\partial
}{\partial x_{\text{c}}}(\nabla\log\Phi)_{j}](y_{0})\Lambda_{j}(y_{0}%
)\phi(y_{0})$

Since $\frac{\partial}{\partial x_{\text{a}}}(\nabla\log\Phi)_{j}%
(y_{0})=-\frac{\partial X_{j}}{\partial x_{\text{a}}}(y_{0})$ for
$j=q+1,...,n,$ we have,

$\left(  D_{48}\right)  $\qquad I$_{338}=-$ $\frac{1}{2}%
\underset{j=q+1}{\overset{n}{\sum}}\underset{\text{c=1}}{\overset{\text{q}%
}{\sum}}[\frac{\partial X_{j}}{\partial x_{\text{a}}}\Lambda_{\text{c}}%
^{2}\Lambda_{j}](y_{0})\phi(y_{0})$

We finally consider:

I$_{339}=\frac{1}{4}\underset{\text{c=1}}{\overset{\text{q}}{\sum}}%
\Lambda_{\text{c}}(y_{0})\frac{\partial}{\partial x_{\text{c}}}[$W$\phi
\circ\pi_{P}](y_{0})$

$\left(  D_{49}\right)  \qquad$I$_{339}=\frac{1}{4}\underset{\text{c=1}%
}{\overset{\text{q}}{\sum}}\left[  \Lambda_{\text{c}}(y_{0})\frac
{\partial\text{W}}{\partial x_{\text{c}}}(y_{0})\phi(y_{0})+\Lambda_{\text{c}%
}(y_{0})\text{W}(y_{0})\frac{\partial\phi}{\partial x_{\text{c}}}%
(y_{0})\right]  $

Next have:

$\left(  D_{50}\right)  \qquad$E$_{1}=$ $\underset{\text{a,c=1}%
}{\overset{\text{q}}{\sum}}\Lambda_{\text{c}}(y_{0})\frac{\partial}{\partial
x_{\text{a}}}[$X$_{\text{a}}\frac{\partial\phi}{\partial\text{x}_{\text{a}}%
}\circ\pi_{\text{P}}](y_{0})=$ $\underset{\text{a=1}}{\overset{\text{q}}{\sum
}}\Lambda_{\text{c}}(y_{0})\frac{\partial\text{X}_{\text{a}}}{\partial
x_{\text{a}}}(y_{0})[\frac{\partial\phi}{\partial\text{x}_{\text{a}}}+$
X$_{\text{a}}\frac{\partial^{2}\phi}{\partial\text{x}_{\text{a}}^{2}}](y_{0})$

Finally here we have:

\qquad E$_{2}=$ $\underset{j=1}{\overset{n}{\sum}}$ $\Lambda_{\text{c}}%
(y_{0})\frac{\partial}{\partial x_{\text{a}}}[$X$_{j}\Lambda_{j}\phi\circ
\pi_{\text{P}}](y_{0})$

\qquad$=$ $\underset{\text{b=1}}{\overset{\text{q}}{\sum}}$ $\Lambda
_{\text{c}}(y_{0})\frac{\partial}{\partial x_{\text{a}}}[$X$_{\text{b}}%
\Lambda_{\text{b}}\phi\circ\pi_{\text{P}}](y_{0})+$
$\underset{j=q+1}{\overset{n}{\sum}}$ $\Lambda_{\text{c}}(y_{0})\frac
{\partial}{\partial x_{\text{a}}}[$X$_{j}\Lambda_{j}\phi\circ\pi_{\text{P}%
}](y_{0})$

$\left(  D_{51}\right)  $\qquad E$_{2}=$ $\underset{\text{b=1}%
}{\overset{\text{q}}{\sum}}$ $[\frac{\partial\text{X}_{\text{b}}}{\partial
x_{\text{a}}}\Lambda_{\text{c}}\Lambda_{\text{b}}](y_{0})\phi(y_{0})+$
$\underset{\text{b=1}}{\overset{\text{q}}{\sum}}$ $[$X$_{\text{b}}%
\Lambda_{\text{c}}\Lambda_{\text{b}}](y_{0})\frac{\partial\phi}{\partial
x_{\text{a}}}y_{0})$

$\qquad\qquad\qquad+$ $\underset{j=q+1}{\overset{\text{n}}{\sum}}$
$\frac{\partial X_{j}}{\partial x_{\text{a}}}[\Lambda_{\text{c}}\Lambda
_{j}](y_{0})\phi(y_{0})+$ $\underset{j=q+1}{\overset{\text{n}}{\sum}}$
$[X_{j}\Lambda_{\text{c}}\Lambda_{j}\frac{\partial\phi}{\partial x_{\text{a}}%
}](y_{0})$

We now gather all the terms of I$_{33}=\frac{1}{2}\underset{\text{c=1}%
}{\overset{\text{q}}{\sum}}\Lambda_{\text{c}}(y_{0})\frac{\partial\Theta
}{\partial x_{\text{c}}}(y_{0})$ in $\left(  D_{41}\right)  ,\left(
D_{42}\right)  ,\left(  D_{43}\right)  ,\left(  D_{44}\right)  ,\left(
D_{45}\right)  ,$

$\left(  D_{46}\right)  ,\left(  D_{47}\right)  ,\left(  D_{48}\right)
,\left(  D_{49}\right)  ,\left(  D_{50}\right)  ,\left(  D_{51}\right)  $ and have:

$\left(  D_{52}\right)  $ \qquad I$_{33}=\frac{1}{48}$ $\underset{\text{c=1}%
}{\overset{\text{q}}{\sum}}\Lambda_{\text{c}}(y_{0}%
)[\underset{i=q+1}{\overset{n}{\sum}}3<H,i>^{2}+2(\tau^{M}-3\tau
^{P}\ +\overset{q}{\underset{\text{a=1}}{\sum}}\varrho_{\text{aa}}%
^{M}+\overset{q}{\underset{\text{a,b}=1}{\sum}}R_{\text{abab}}^{M}%
)](y_{0})\frac{\partial\phi}{\partial x_{\text{c}}}(y_{0})\qquad$I$_{331}$

\qquad$-\frac{1}{4}\underset{\text{c=1}}{\overset{\text{q}}{\sum}}%
\Lambda_{\text{c}}(y_{0})[\left\Vert \text{X}\right\Vert ^{2}+$ divX $-$
$\underset{\text{a}=1}{\overset{q}{\sum}}($X$_{\text{a}})^{2}$ $-$
$\underset{\text{a}=1}{\overset{q}{\sum}}\frac{\partial X_{\text{a}}}{\partial
x_{\text{a}}}](y_{0})\frac{\partial\phi}{\partial x_{\text{c}}}(y_{0})+$
$\frac{1}{2}\underset{\text{c=1}}{\overset{\text{q}}{\sum}}\Lambda_{\text{c}%
}(y_{0})$V$(y_{0})\frac{\partial\phi}{\partial x_{\text{c}}}(y_{0})$

\qquad$+\frac{1}{2}\underset{\text{c=1}}{\overset{\text{q}}{\sum}}%
\Lambda_{\text{c}}(y_{0})[-(X_{j}\frac{\partial X_{j}}{\partial x_{\text{c}}%
}+\frac{1}{2}\frac{\partial^{2}X_{j}}{\partial x_{\text{c}}\partial x_{j}%
})(y_{0})+\frac{1}{2}(<H,j>\frac{\partial X_{j}}{\partial x_{\text{c}}}%
)(y_{0})+\frac{\partial\text{V}}{\partial x_{\text{c}}}(y_{0})]\phi(y_{0})$

$\qquad+\frac{1}{4}\underset{\text{a,c=1}}{\overset{\text{q}}{\sum}}%
\Lambda_{\text{c}}(y_{0})\left\{  \frac{\partial^{3}\phi}{\partial
\text{x}_{\text{a}}^{2}\partial x_{\text{c}}}\text{ }\right\}  (y_{0})\qquad
$I$_{332}$

$\qquad+\frac{1}{2}\underset{\text{a,b=1}}{\overset{\text{q}}{\sum}}\left\{
\Lambda_{\text{a}}(y_{0})\Lambda_{\text{b}}(y_{0})\frac{\partial^{2}\phi
}{\partial\text{x}_{\text{a}}\partial x_{\text{b}}}\text{ }\right\}
(y_{0})\qquad$I$_{334}$

\qquad$+\frac{1}{4}\underset{\text{c=1}}{\overset{\text{q}}{\sum}%
}[\underset{i=1}{\overset{n}{\sum}}\Lambda_{\text{c}}(y_{0})\Lambda_{i}%
^{2}(y_{0})\frac{\partial\phi}{\partial x_{\text{c}}}(y_{0})]\qquad$I$_{335}$

\qquad$-\frac{1}{4}\underset{j=q+1}{\overset{n}{\sum}}\underset{\text{a,c=1}%
}{\overset{\text{q}}{\sum}}[\Lambda_{\text{c}}\Lambda_{j}T_{\text{aa}j}%
\frac{\partial\phi}{\partial x_{\text{c}}}](y_{0})\qquad$I$_{336}$

\qquad$-$ $\frac{1}{2}\underset{j=q+1}{\overset{n}{\sum}}\underset{\text{c=1}%
}{\overset{\text{q}}{\sum}}[\frac{\partial X_{j}}{\partial x_{\text{a}}%
}\Lambda_{\text{c}}^{2}\Lambda_{j}](y_{0})\phi(y_{0})\qquad$I$_{338}$

\qquad$+\frac{1}{4}\underset{\text{c=1}}{\overset{\text{q}}{\sum}}\left[
\Lambda_{\text{c}}(y_{0})\frac{\partial\text{W}}{\partial x_{\text{c}}}%
(y_{0})\phi(y_{0})+\Lambda_{\text{c}}(y_{0})\text{W}(y_{0})\frac{\partial\phi
}{\partial x_{\text{c}}}(y_{0})\right]  \qquad$I$_{339}$

\qquad$+$ $\underset{\text{a=1}}{\overset{\text{q}}{\sum}}\Lambda_{\text{c}%
}(y_{0})\frac{\partial\text{X}_{\text{a}}}{\partial x_{\text{a}}}(y_{0}%
)[\frac{\partial\phi}{\partial\text{x}_{\text{a}}}+$ X$_{\text{a}}%
\frac{\partial^{2}\phi}{\partial\text{x}_{\text{a}}^{2}}](y_{0})\qquad$E$_{1}$

\qquad$+$ $\underset{\text{b=1}}{\overset{\text{q}}{\sum}}$ $[\Lambda
_{\text{b}}\Lambda_{\text{c}}\frac{\partial\text{X}_{\text{b}}}{\partial
x_{\text{a}}}](y_{0})\phi(y_{0})+$ $\underset{\text{b=1}}{\overset{\text{q}%
}{\sum}}$ $[\Lambda_{\text{c}}\Lambda_{\text{b}}$X$_{\text{b}}](y_{0}%
)\frac{\partial\phi}{\partial x_{\text{a}}}(y_{0})\qquad$E$_{2}$

$\qquad+$ $\underset{j=q+1}{\overset{n}{\sum}}$ $[\Lambda_{\text{c}}%
\Lambda_{j}\frac{\partial X_{j}}{\partial x_{\text{a}}}](y_{0})\phi(y_{0})+$
$\underset{j=q+1}{\overset{n}{\sum}}$ $[\Lambda_{\text{c}}\Lambda_{j}%
\frac{\partial\phi}{\partial x_{\text{a}}}](y_{0})$

\qquad\qquad\qquad\qquad\qquad\qquad\qquad\qquad\qquad\qquad\qquad\qquad
\qquad\qquad\qquad$\blacksquare$

\subsection{\textbf{Computation of I}$_{34}$\qquad}

\qquad I$_{34}=\frac{1}{4}\underset{\text{a=1}}{\overset{\text{q}}{\sum}%
}\Lambda_{\text{c}}^{2}(y_{0})\Theta(y_{0})\phi\left(  y_{0}\right)  $

From $\left(  10.31\right)  ,$ or from $\left(  D_{2}\right)  $ above, we have:

\qquad$\Theta(y_{0})\phi\left(  y_{0}\right)  =$ \textbf{b}$_{1}$(y$_{0}%
$,P)$\phi\left(  y_{0}\right)  $

$\qquad=\frac{1}{24}[\underset{i=q+1}{\overset{n}{\sum}}3<H,i>^{2}+2(\tau
^{M}-3\tau^{P}\ +\overset{q}{\underset{\text{a=1}}{\sum}}\varrho_{\text{aa}%
}^{M}+\overset{q}{\underset{\text{a,b}=1}{\sum}}R_{\text{abab}}^{M}%
)](y_{0})\phi\left(  y_{0}\right)  $

\qquad$-\frac{1}{2}[$ $\left\Vert \text{X}\right\Vert _{M}^{2}+\frac{1}{2}$
$\operatorname{div}X_{M}-\frac{1}{2}$ $\left\Vert \text{X}\right\Vert _{P}%
^{2}$ $-$ $\frac{1}{2}\operatorname{div}X_{P}](y_{0})\phi\left(  y_{0}\right)
$

$\qquad+$ $\frac{1}{2}\underset{\text{a=1}}{\overset{\text{q}}{\sum}}%
\frac{\partial^{2}\phi}{\partial\text{x}_{\text{a}}^{2}}(y_{0})$ $+$
$\underset{\text{a=1}}{\overset{\text{q}}{\sum}}\Lambda_{\text{a}}(y_{0}%
)\frac{\partial\phi}{\partial x_{\text{a}}}\left(  y_{0}\right)  \ +\frac
{1}{2}$ $\underset{\text{a=1}}{\overset{\text{q}}{\sum}}\Lambda_{\text{a}%
}(y_{0})\Lambda_{\text{a}}(y_{0})\phi\left(  y_{0}\right)  $

\qquad$+$ $\underset{\text{a=1}}{\overset{\text{q}}{\sum}}X_{\text{a}}%
(y_{0})\frac{\partial\phi}{\partial\text{x}_{\text{a}}}(y_{0})+$
$\underset{\text{a=1}}{\overset{\text{q}}{\sum}}$ $X_{\text{a}}(y_{0}%
)\Lambda_{\text{a}}(y_{0})\phi(y_{0})+\frac{1}{2}$W$\left(  y_{0}\right)
\phi\left(  y_{0}\right)  +$ V$(y_{0})\phi\left(  y_{0}\right)  $

Consequently we have,

$\left(  D_{53}\right)  \qquad$I$_{34}=\frac{1}{4}\underset{\text{c=1}%
}{\overset{\text{q}}{\sum}}\Lambda_{\text{c}}^{2}(y_{0})\Theta(y_{0})$

$=\frac{1}{96}\underset{\text{c=1}}{\overset{\text{q}}{\sum}}\Lambda
_{\text{c}}^{2}(y_{0})[\underset{i=q+1}{\overset{n}{\sum}}3<H,i>^{2}%
+2(\tau^{M}-3\tau^{P}\ +\overset{q}{\underset{\text{a=1}}{\sum}}%
\varrho_{\text{aa}}^{M}+\overset{q}{\underset{\text{a,b}=1}{\sum}%
}R_{\text{abab}}^{M})](y_{0})\phi\left(  y_{0}\right)  $

$-\frac{1}{8}\underset{\text{c=1}}{\overset{\text{q}}{\sum}}\Lambda_{\text{c}%
}^{2}(y_{0})[$ $\left\Vert \text{X}\right\Vert _{M}^{2}+\frac{1}{2}$
$\operatorname{div}X_{M}-\frac{1}{2}$ $\left\Vert \text{X}\right\Vert _{P}%
^{2}$ $-$ $\frac{1}{2}\operatorname{div}X_{P}](y_{0})\phi\left(  y_{0}\right)
$

$+$ $\frac{1}{8}\underset{\text{c=1}}{\overset{\text{q}}{\sum}}\Lambda
_{\text{c}}^{2}(y_{0})[\underset{\text{a=1}}{\overset{\text{q}}{\sum}}%
\frac{\partial^{2}\phi}{\partial\text{x}_{\text{a}}^{2}}$ $+$ $2$
$\underset{\text{a=1}}{\overset{\text{q}}{\sum}}\Lambda_{\text{a}}%
\frac{\partial\phi}{\partial x_{\text{a}}}\ +$ $\underset{\text{a=1}%
}{\overset{\text{q}}{\sum}}\Lambda_{\text{a}}^{2}](y_{0})\phi\left(
y_{0}\right)  $

$+\frac{1}{4}\underset{\text{c=1}}{\overset{\text{q}}{\sum}}\Lambda_{\text{c}%
}^{2}(y_{0})[$ $\underset{\text{a=1}}{\overset{\text{q}}{\sum}}X_{\text{a}%
}\frac{\partial\phi}{\partial\text{x}_{\text{a}}}+$ $\underset{\text{a=1}%
}{\overset{\text{q}}{\sum}}$ $X_{\text{a}}\Lambda_{\text{a}}+\frac{1}{2}$W $+$
V$](y_{0})\phi\left(  y_{0}\right)  $

\subsection{\textbf{Computation of I}$_{35}$}

The expression of $\Theta(y_{0})$ is given in $\left(  D_{2}\right)  $ and has
been repeatedly used above and so:

$\left(  D_{54}\right)  \qquad$I$_{35}=$ $\frac{1}{4}$W$(y_{0})\Theta(y_{0})=$
$\frac{1}{4}\Theta(y_{0})$W$(y_{0})$

$\qquad\qquad=\frac{1}{96}[\underset{i=q+1}{\overset{n}{\sum}}3<H,i>^{2}%
+2(\tau^{M}-3\tau^{P}\ +\overset{q}{\underset{\text{a=1}}{\sum}}%
\varrho_{\text{aa}}^{M}+\overset{q}{\underset{\text{a,b}=1}{\sum}%
}R_{\text{abab}}^{M})](y_{0})$W$(y_{0})\phi\left(  y_{0}\right)  $

\qquad\qquad$\ -\frac{1}{8}[$ $\left\Vert \text{X}\right\Vert _{M}^{2}+$
$\operatorname{div}X_{M}-$ $\left\Vert \text{X}\right\Vert _{P}^{2}$
$-\operatorname{div}X_{P}](y_{0})$W(y$_{0}$)$\phi\left(  y_{0}\right)  $

$\qquad\qquad+$ $\frac{1}{8}\underset{\text{a=1}}{\overset{\text{q}}{\sum}%
}\frac{\partial^{2}\phi}{\partial\text{x}_{\text{a}}^{2}}(y_{0})$W$(y_{0})$
$+\frac{1}{4}$ $\underset{\text{a=1}}{\overset{\text{q}}{\sum}}\Lambda
_{\text{a}}(y_{0})\frac{\partial\phi}{\partial x_{\text{a}}}\left(
y_{0}\right)  $W$(y_{0})\ +\frac{1}{8}$ $\underset{\text{a=1}%
}{\overset{\text{q}}{\sum}}\Lambda_{\text{a}}^{2}(y_{0})$W$(y_{0})\phi\left(
y_{0}\right)  $

\qquad$+\frac{1}{4}$ $\underset{\text{a=1}}{\overset{\text{q}}{\sum}}%
$X$_{\text{a}}(y_{0})\frac{\partial\phi}{\partial\text{x}_{\text{a}}}(y_{0}%
)$W$(y_{0})+\frac{1}{4}$ $\underset{\text{a=1}}{\overset{\text{q}}{\sum}}$
X$_{\text{a}}(y_{0})\Lambda_{\text{a}}(y_{0})$W$(y_{0})\phi(y_{0})+\frac{1}%
{8}$W$^{2}\left(  y_{0}\right)  \phi(y_{0})+\frac{1}{4}$ V$(y_{0})$%
W$(y_{0})\phi\left(  y_{0}\right)  $

\subsection{\protect\underline{\textbf{Computation of I}$_{36}$}$\qquad$}

I$_{36}=$ $\frac{1}{2}$ $\underset{\text{a=1}}{\overset{\text{q}}{\sum}}%
$X$_{\text{c}}(y_{0})\frac{\partial\Theta}{\partial x_{\text{c}}}(y_{0})$

Recall that I$_{33}=\frac{1}{2}\underset{\text{c=1}}{\overset{\text{q}}{\sum}%
}\Lambda_{\text{c}}(y_{0})\frac{\partial\Theta}{\partial x_{\text{c}}}%
(y_{0}).$ Therefore I$_{36}=$ $\frac{1}{2}$ $\underset{\text{a=1}%
}{\overset{\text{q}}{\sum}}$X$_{\text{c}}(y_{0})\frac{\partial\Theta}{\partial
x_{\text{c}}}(y_{0})$ here is similar to I$_{33}$ in $\left(  D_{52}\right)
:$

We replace $\Lambda_{\text{c}}(y_{0})$ there with X$_{\text{c}}(y_{0})$ here
and have:

$\left(  D_{55}\right)  \qquad$I$_{36}=$ $\frac{1}{2}$ $\underset{\text{c=1}%
}{\overset{\text{q}}{\sum}}$X$_{\text{c}}(y_{0})\frac{\partial\Theta}{\partial
x_{\text{c}}}(y_{0})$

$=\frac{1}{48}$ $\underset{\text{c=1}}{\overset{\text{q}}{\sum}}$X$_{\text{c}%
}(y_{0})[\underset{\alpha=q+1}{\overset{n}{\sum}}3<H,i>^{2}+2(\tau^{M}%
-3\tau^{P}\ +\overset{q}{\underset{\text{a=1}}{\sum}}\varrho_{\text{aa}}%
^{M}+\overset{q}{\underset{\text{a,b}=1}{\sum}}R_{\text{abab}}^{M}%
)](y_{0})\frac{\partial\phi}{\partial x_{\text{c}}}(y_{0})\qquad$I$_{361}$

$-\frac{1}{4}\underset{\text{c=1}}{\overset{\text{q}}{\sum}}$X$_{\text{c}%
}(y_{0})[\left\Vert \text{X}\right\Vert ^{2}+$ divX $-$ $\underset{\text{a}%
=1}{\overset{q}{\sum}}($X$_{\text{a}})^{2}$ $-$ $\underset{\text{a}%
=1}{\overset{q}{\sum}}\frac{\partial X_{\text{a}}}{\partial x_{\text{a}}%
}](y_{0})\frac{\partial\phi}{\partial x_{\text{c}}}(y_{0})+$ $\frac{1}%
{2}\underset{\text{c=1}}{\overset{\text{q}}{\sum}}\Lambda_{\text{c}}(y_{0}%
)$V$(y_{0})\frac{\partial\phi}{\partial x_{\text{c}}}(y_{0})$

$+\frac{1}{2}\underset{\text{c=1}}{\overset{\text{q}}{\sum}}$X$_{\text{c}%
}(y_{0})[-(X_{j}\frac{\partial X_{j}}{\partial x_{\text{c}}}+\frac{1}{2}%
\frac{\partial^{2}X_{j}}{\partial x_{\text{c}}\partial x_{j}})(y_{0})+\frac
{1}{2}(<H,j>\frac{\partial X_{j}}{\partial x_{\text{c}}})(y_{0})+\frac
{\partial\text{V}}{\partial x_{\text{c}}}(y_{0})]\phi(y_{0})$

$+\frac{1}{4}\underset{\text{a,c=1}}{\overset{\text{q}}{\sum}}$X$_{\text{c}%
}(y_{0})\left\{  \frac{\partial^{3}\phi}{\partial\text{x}_{\text{a}}%
^{2}\partial x_{\text{c}}}\text{ }\right\}  (y_{0})$ \qquad\qquad I$_{362}$

$+\frac{1}{2}\underset{\text{a,c=1}}{\overset{\text{q}}{\sum}}\left\{
\text{X}_{\text{c}}(y_{0})\Lambda_{\text{a}}(y_{0})\frac{\partial^{2}\phi
}{\partial\text{x}_{\text{a}}\partial x_{\text{c}}}\text{ }\right\}
(y_{0})\qquad$I$_{364}$

$+\frac{1}{4}\underset{\text{b,c=1}}{\overset{\text{q}}{\sum}}$X$_{\text{c}%
}(y_{0})\Lambda_{\text{b}}^{2}(y_{0})\frac{\partial\phi}{\partial x_{\text{c}%
}}(y_{0})\qquad$I$_{365}$

$+\frac{1}{4}\underset{\text{a,c=1}}{\overset{\text{q}}{\sum}}\left[
\text{X}_{\text{c}}(y_{0})\frac{\partial\text{W}}{\partial x_{\text{a}}}%
(y_{0})\phi(y_{0})+\text{X}_{\text{c}}(y_{0})\text{W}(y_{0})\frac{\partial
\phi}{\partial x_{\text{c}}}(y_{0})\right]  \qquad$I$_{369}\qquad$

$+$ $\underset{\text{a,c=1}}{\overset{\text{q}}{\sum}}$X$_{\text{c}}%
(y_{0})\frac{\partial\text{X}_{\text{a}}}{\partial x_{\text{a}}}(y_{0}%
)[\frac{\partial\phi}{\partial\text{x}_{\text{a}}}+$ X$_{\text{a}}%
\frac{\partial^{2}\phi}{\partial\text{x}_{\text{a}}^{2}}](y_{0})\qquad\qquad
$E$_{1}$

$+\underset{\text{a,b,c=1}}{\overset{\text{q}}{\sum}}$ X$_{\text{c}}%
(y_{0})\frac{\partial\text{X}_{\text{b}}}{\partial x_{\text{a}}}(y_{0}%
)\Lambda_{\text{b}}(y_{0})\phi(y_{0})+$ $\underset{\text{a,b,c=1}%
}{\overset{\text{q}}{\sum}}$ X$_{\text{c}}(y_{0})$X$_{\text{b}}(y_{0}%
)\Lambda_{\text{b}}(y_{0})\frac{\partial\phi}{\partial x_{\text{a}}}%
(y_{0})\qquad$E$_{2}$

\subsection{\textbf{Computation of I}$_{37}$}

$\left(  D_{56}\right)  \qquad$I$_{37}=\frac{1}{2}$
$\underset{j=1}{\overset{n}{\sum}}$X$_{j}(y_{0})\Lambda_{j}(y_{0})\Theta
(y_{0})$

\qquad$\qquad=[\frac{1}{2}$ $\underset{\text{a=1}}{\overset{\text{q}}{\sum}}%
$X$_{\text{a}}(y_{0})\Lambda_{\text{a}}(y_{0})+\frac{1}{2}$
$\underset{j=q+1}{\overset{n}{\sum}}$X$_{j}(y_{0})\Lambda_{j}(y_{0}%
)]\Theta(y_{0})$

I$_{37}=\frac{1}{48}\underset{\text{a=1}}{\overset{\text{q}}{\sum}}%
X_{\text{a}}(y_{0})\Lambda_{\text{a}}(y_{0})[\underset{i=q+1}{\overset{n}{\sum
}}3<H,i>^{2}+2(\tau^{M}-3\tau^{P}\ +\overset{q}{\underset{\text{a=1}}{\sum}%
}\varrho_{\text{aa}}^{M}+\overset{q}{\underset{\text{a,b}=1}{\sum}%
}R_{\text{abab}}^{M})](y_{0})\phi\left(  y_{0}\right)  $

\qquad$-\frac{1}{4}\underset{\text{a=1}}{\overset{\text{q}}{\sum}}X_{\text{a}%
}(y_{0})\Lambda_{\text{a}}(y_{0})[$ $\left\Vert \text{X}\right\Vert _{M}%
^{2}+\frac{1}{2}$ $\operatorname{div}X_{M}-\frac{1}{2}$ $\left\Vert
\text{X}\right\Vert _{P}^{2}$ $-$ $\frac{1}{2}\operatorname{div}X_{P}%
](y_{0})\phi\left(  y_{0}\right)  $

$\qquad+$ $\frac{1}{4}\underset{\text{a,b=1}}{\overset{\text{q}}{\sum}%
}X_{\text{a}}(y_{0})\Lambda_{\text{a}}(y_{0})[\frac{\partial^{2}\phi}%
{\partial\text{x}_{\text{b}}^{2}}$ $+$ $\Lambda_{\text{b}}\frac{\partial\phi
}{\partial x_{\text{b}}}]\left(  y_{0}\right)  \ +\frac{1}{4}$
$\underset{\text{a,b=1}}{\overset{\text{q}}{\sum}}X_{\text{a}}(y_{0}%
)\Lambda_{\text{a}}(y_{0})\Lambda_{\text{b}}^{2}(y_{0})\phi\left(
y_{0}\right)  $

\qquad$+\frac{1}{4}$ $\underset{\text{a,b=1}}{\overset{\text{q}}{\sum}}%
$X$_{\text{a}}(y_{0})\Lambda_{\text{a}}(y_{0})$ $X_{\text{b}}(y_{0}%
)\frac{\partial\phi}{\partial\text{x}_{\text{b}}}(y_{0})$

$\qquad+\frac{1}{2}$ $\underset{\text{a,b=1}}{\overset{\text{q}}{\sum}}$
$X_{\text{a}}(y_{0})\Lambda_{\text{a}}(y_{0})X_{\text{b}}(y_{0})\Lambda
_{\text{b}}(y_{0})\phi(y_{0})+\frac{1}{4}$ $\underset{\text{a,b=1}%
}{\overset{\text{q}}{\sum}}$ $X_{\text{a}}(y_{0})\Lambda_{\text{a}}(y_{0}%
)$W$\left(  y_{0}\right)  \phi\left(  y_{0}\right)  $

$\qquad+\frac{1}{2}$ $\underset{\text{a,b=1}}{\overset{\text{q}}{\sum}}$
$X_{\text{a}}(y_{0})\Lambda_{\text{a}}(y_{0})$V$(y_{0})\phi(y_{0})$

\qquad$+\frac{1}{48}$ $\underset{j=q+1}{\overset{n}{\sum}}$X$_{j}%
(y_{0})\Lambda_{j}(y_{0})[\underset{i=q+1}{\overset{n}{\sum}}3<H,i>^{2}%
+2(\tau^{M}-3\tau^{P}\ +\overset{q}{\underset{\text{a=1}}{\sum}}%
\varrho_{\text{aa}}^{M}+\overset{q}{\underset{\text{a,b}=1}{\sum}%
}R_{\text{abab}}^{M})](y_{0})\phi\left(  y_{0}\right)  $

\qquad$\ -\frac{1}{4}\underset{j=q+1}{\overset{n}{\sum}}X_{j}(y_{0}%
)\Lambda_{j}(y_{0})[$ $\left\Vert \text{X}\right\Vert _{M}^{2}+\frac{1}{2}$
$\operatorname{div}X_{M}-\frac{1}{2}$ $\left\Vert \text{X}\right\Vert _{P}%
^{2}$ $-$ $\frac{1}{2}\operatorname{div}X_{P}](y_{0})\phi\left(  y_{0}\right)
$

$\qquad+$ $\frac{1}{4}\underset{j=q+1}{\overset{n}{\sum}}\underset{\text{a=1}%
}{\overset{\text{q}}{\sum}}X_{j}(y_{0})\Lambda_{j}(y_{0})\frac{\partial
^{2}\phi}{\partial\text{x}_{\text{a}}^{2}}(y_{0})$ $+\frac{1}{2}%
\underset{j=q+1}{\overset{n}{\sum}}\underset{\text{a=1}}{\overset{\text{q}%
}{\sum}}X_{j}(y_{0})\Lambda_{j}(y_{0})\Lambda_{\text{a}}(y_{0})\frac
{\partial\phi}{\partial x_{\text{a}}}\left(  y_{0}\right)  $

$\qquad+\frac{1}{4}$ $\underset{j=q+1}{\overset{n}{\sum}}\underset{\text{a=1}%
}{\overset{\text{q}}{\sum}}X_{j}(y_{0})\Lambda_{j}(y_{0})\Lambda_{\text{a}%
}^{2}(y_{0})\phi\left(  y_{0}\right)  +$ $\frac{1}{4}$
$\underset{j=q+1}{\overset{n}{\sum}}\underset{\text{a=1}}{\overset{\text{q}%
}{\sum}}X_{j}(y_{0})\Lambda_{j}(y_{0})X_{\text{a}}(y_{0})\frac{\partial\phi
}{\partial\text{x}_{\text{a}}}(y_{0})$

$\qquad+$ $\frac{1}{2}$ $\underset{j=q+1}{\overset{n}{\sum}}%
\underset{\text{a=1}}{\overset{\text{q}}{\sum}}X_{\text{a}}(y_{0}%
)\Lambda_{\text{a}}(y_{0})X_{j}(y_{0})\Lambda_{j}(y_{0})$ $\phi(y_{0}%
)+\frac{1}{4}$ $\underset{j=q+1}{\overset{n}{\sum}}X_{j}(y_{0})\Lambda
_{j}(y_{0})$W$\left(  y_{0}\right)  \phi\left(  y_{0}\right)  $

$\qquad+$ $\frac{1}{2}$ $\underset{j=q+1}{\overset{n}{\sum}}X_{j}%
(y_{0})\Lambda_{j}(y_{0})$V$(y_{0})\phi\left(  y_{0}\right)  $

\qquad\qquad\qquad\qquad\qquad\qquad\qquad\qquad\qquad\qquad\qquad\qquad
\qquad\qquad\qquad\qquad\qquad\qquad\qquad\qquad$\blacksquare$

\section{\qquad EXPRESSION FOR b$_{2}($y$_{0}$,P$)\phi($y$_{0})$}

At a general point x$\in M_{0},$ the third coefficient b$_{2}(y_{0},P,\phi)$
is defined in \textbf{Theorem }$\left(  5.3\right)  $ by:

\begin{center}
b$_{2}($x,P$,\phi)=\int_{0}^{1}\int_{0}^{r_{1}}$F(1,1-r$_{2}$)[L$_{\Psi}%
$F(1-r$_{2}$,1-r$_{1}$)L$_{\Psi}[\phi\circ\pi_{P}$](x)dr$_{1}$dr$_{2}$
\end{center}

Computation of the third coefficient is impossible with the mathematical tools
presently available. Even the computation at the particular point y$_{0}\in P$
will be very long.

We now present the third coefficient, expressed in \textbf{geometric
invariants}. It is one of the most \textbf{significant achievements} of this work.

By $\left(  D_{4}\right)  ,\left(  D_{8}\right)  ,\left(  D_{40}\right)
,\left(  D_{52}\right)  ,\left(  D_{53}\right)  ,\left(  D_{54}\right)
,\left(  D_{55}\right)  ,\left(  D_{56}\right)  ,$ we have:

$\left(  D_{57}\right)  $\qquad b$_{2}($y$_{0}$,P$,\phi)=$ I$_{1}+$ I$_{31}+$
I$_{32}+$ I$_{33}+$ I$_{34}+$ I$_{35}+$ I$_{36}+$ I$_{37}$

$=\frac{1}{2}[\frac{1}{24}(\underset{i=q+1}{\overset{n}{\sum}}3<H,i>^{2}%
+2(\tau^{M}-3\tau^{P}\ +\overset{q}{\underset{\text{a=1}}{\sum}}%
\varrho_{\text{aa}}^{M}+\overset{q}{\underset{\text{a,b}=1}{\sum}%
}R_{\text{abab}}^{M}))\qquad$I$_{1}$

$-\frac{1}{2}($ $\left\Vert \text{X}\right\Vert _{M}^{2}+$ $\operatorname{div}%
X_{M}-\left\Vert \text{X}\right\Vert _{P}^{2}$ $-$ $\operatorname{div}X_{P})+$
$V](y_{0})$

$\times\lbrack\frac{1}{24}(\underset{i=q+1}{\overset{n}{\sum}}3<H,i>^{2}%
+2(\tau^{M}-3\tau^{P}\ +\overset{q}{\underset{\text{a=1}}{\sum}}%
\varrho_{\text{aa}}^{M}+\overset{q}{\underset{\text{a,b}=1}{\sum}%
}R_{\text{abab}}^{M})$

$-\frac{1}{2}($ $\left\Vert \text{X}\right\Vert _{M}^{2}+$ $\operatorname{div}%
X_{M}-\left\Vert \text{X}\right\Vert _{P}^{2}$ $-\operatorname{div}X_{P})+$
$V+\frac{1}{2}W$

$+$ $(\frac{1}{2}\underset{\text{a=1}}{\overset{\text{q}}{\sum}}\frac
{\partial^{2}\phi}{\partial\text{x}_{\text{a}}^{2}}$ $+$ $\underset{\text{a=1}%
}{\overset{\text{q}}{\sum}}\Lambda_{\text{a}}(y_{0})\frac{\partial\phi
}{\partial x_{\text{a}}}\ +\frac{1}{2}$ $\underset{\text{a=1}%
}{\overset{\text{q}}{\sum}}\Lambda_{\text{a}}\Lambda_{\text{a}})+$
$\underset{\text{a=1}}{\overset{\text{q}}{\sum}}X_{\text{a}}\frac{\partial
\phi}{\partial\text{x}_{\text{a}}}+$ $\underset{\text{a=1}}{\overset{\text{q}%
}{\sum}}$ $X_{\text{a}}\Lambda_{\text{a}})](y_{0})\phi\left(  y_{0}\right)  $

$+\frac{1}{96}[\underset{i=q+1}{\overset{n}{\sum}}3<H,i>^{2}+2(\tau^{M}%
-3\tau^{P}\ +\overset{q}{\underset{\text{a=1}}{\sum}}\varrho_{\text{aa}}%
^{M}+\overset{q}{\underset{\text{a,b}=1}{\sum}}R_{\text{abab}}^{M}%
)](y_{0})\frac{\partial^{2}\phi}{\partial x_{\text{c}}^{2}}(y_{0})\qquad
$I$_{31}\qquad$I$_{311}$

$-\frac{1}{4}[\left\Vert \text{X}(y_{0})\right\Vert ^{2}+$ divX$(y_{0})-$
$\underset{\text{a}=1}{\overset{q}{\sum}}($X$_{\text{a}})^{2}(y_{0})$ $-$
$\underset{\text{a}=1}{\overset{q}{\sum}}\frac{\partial X_{\text{a}}}{\partial
x_{\text{a}}}(y_{0})]\frac{\partial^{2}\phi}{\partial x_{\text{c}}^{2}}%
(y_{0})+$ $\frac{1}{4}$V(y$_{0}$)$\frac{\partial^{2}\phi}{\partial
x_{\text{c}}^{2}}(y_{0})$

$-\frac{1}{2}$ $[X_{j}\frac{\partial X_{j}}{\partial x_{\text{c}}}+\frac{1}%
{2}\frac{\partial^{2}X_{j}}{\partial x_{\text{c}}\partial x_{j}}](y_{0}%
).\frac{\partial\phi}{\partial x_{\text{c}}}(y_{0})+\frac{1}{4}[<H,j>\frac
{\partial X_{j}}{\partial x_{\text{c}}}](y_{0}).\frac{\partial\phi}{\partial
x_{\text{c}}}(y_{0})+$ $\frac{1}{2}\frac{\partial\text{V}}{\partial
x_{\text{c}}}(y_{0}).\frac{\partial\phi}{\partial x_{\text{c}}}(y_{0})$

$+\frac{1}{4}[(\frac{\partial X_{j}}{\partial x_{\text{c}}})^{2}-X_{j}%
\frac{\partial^{2}X_{j}}{\partial x_{\text{c}}^{2}}](y_{0})\phi(y_{0})$
$-\frac{1}{8}\frac{\partial^{3}X_{j}}{\partial x_{\text{c}}^{2}\partial x_{j}%
}(y_{0})\phi(y_{0})-\frac{1}{2}\frac{\partial X_{i}}{\partial x_{\text{c}}%
}(y_{0}))\frac{\partial X_{i}}{\partial x_{\text{c}}}(y_{0})\phi(y_{0})$

$+\frac{1}{4}\frac{\partial^{2}\text{V}}{\partial x_{\text{c}}^{2}}(y_{0}%
)\phi(y_{0})$

$+\frac{1}{8}\underset{\text{a=1}}{\overset{\text{q}}{\sum}}\frac{\partial
^{4}\phi}{\partial x_{\text{a}}^{2}\partial x_{\text{c}}^{2}}(y_{0}%
)\qquad\qquad\qquad$I$_{312}$

$+\frac{1}{4}\underset{\text{a=1}}{\overset{\text{q}}{\sum}}[\Lambda
_{\text{a}}\frac{\partial^{3}\phi}{\partial\text{x}_{\text{a}}\partial
x_{\text{c}}^{2}}](y_{0})\qquad\qquad$I$_{314}$

$+$ $\frac{1}{8}\underset{\text{a=1}}{\overset{\text{q}}{\sum}}[\Lambda
_{\text{a}}^{2}\frac{\partial^{2}\phi}{\partial x_{\text{c}}^{2}}%
](y_{0})\qquad\ \ \ \ $I$_{315}\qquad$

$+\frac{1}{8}\frac{\partial^{2}\text{W}}{\partial x_{\text{a}}^{2}}(y_{0}%
)\phi(y_{0})+\frac{1}{4}\frac{\partial\text{W}}{\partial x_{\text{a}}}%
(y_{0})\frac{\partial\phi}{\partial x_{\text{a}}}(y_{0})+\frac{1}{8}$%
W$(y_{0})\frac{\partial^{2}\phi}{\partial x_{\text{a}}^{2}}(y_{0})\qquad
$I$_{319}$

$+\frac{1}{4}$ $\underset{\text{a=1}}{\overset{\text{q}}{\sum}}\frac
{\partial^{2}\text{X}_{\text{a}}}{\partial x_{\text{a}}^{2}}(y_{0}%
)\frac{\partial\phi}{\partial\text{x}_{\text{a}}}(y_{0})+\frac{1}{4}$
$\underset{\text{a=1}}{\overset{\text{q}}{\sum}}$X$_{\text{a}}(y_{0}%
)\frac{\partial^{3}\phi}{\partial x_{\text{a}}^{3}}(y_{0})+\frac{1}{2}$
$\underset{\text{a=1}}{\overset{\text{q}}{\sum}}\frac{\partial\text{X}%
_{\text{a}}}{\partial x_{\text{a}}}(y_{0})\frac{\partial^{2}\phi}%
{\partial\text{x}_{\text{a}}^{2}}(y_{0})\qquad$L$_{1}$

$+\frac{1}{4}[$ $\underset{\text{b=1}}{\overset{\text{q}}{\sum}}$
$\frac{\partial^{2}X_{\text{b}}}{\partial x_{\text{a}}^{2}}\Lambda_{\text{b}%
}(y_{0})\phi(y_{0})+\frac{1}{4}[$ $\underset{\text{b=1}}{\overset{\text{q}%
}{\sum}}$ X$_{\text{b}}(y_{0})\Lambda_{\text{b}}(y_{0})\frac{\partial^{2}\phi
}{\partial x_{\text{a}}^{2}}(y_{0})+\frac{1}{2}[$ $\underset{\text{b=1}%
}{\overset{\text{q}}{\sum}}$ $\frac{\partial X_{\text{b}}}{\partial
x_{\text{a}}}(y_{0})\Lambda_{\text{b}}(y_{0})\frac{\partial\phi}{\partial
x_{\text{a}}}(y_{0})\qquad$L$_{2}$

$-\frac{1}{3456}[3<H,i>^{2}\ +2(\tau^{M}-3\tau^{P}%
\ +\overset{q}{\underset{\text{a=1}}{\sum}}\varrho_{\text{aa}}^{M}%
+\overset{q}{\underset{\text{a,b}=1}{\sum}}R_{\text{abab}}^{M})]^{2}%
(y_{0})\phi(y_{0})\phi(y_{0})$\ \qquad\textbf{I}$_{32}\qquad$I$_{321}$

$+\frac{1}{24}[2<H,i>^{2}(y_{0})+\frac{1}{3}(\tau^{M}-3\tau^{P}%
+\overset{q}{\underset{\text{a}=1}{\sum}}\varrho_{\text{aa}}%
+\overset{q}{\underset{\text{a,b}=1}{\sum}}R_{\text{abab}})](y_{0})\qquad
$I$_{3212}=\frac{1}{24}(L_{1}+L_{2}+L_{3})$

$\times\lbrack\frac{1}{4}<H,j>^{2}(y_{0})+\frac{1}{6}(\tau^{M}-3\tau
^{P}+\overset{q}{\underset{\text{a}=1}{\sum}}\varrho_{\text{aa}}%
^{M}+\overset{q}{\underset{\text{a,b}=1}{\sum}}R_{\text{abab}}^{M}%
)](y_{0})\phi(y_{0})$

$-\frac{1}{96}[<H,i><H,j>](y_{0})$

$\times\lbrack2\varrho_{ij}+$ $\overset{q}{\underset{\text{a}=1}{4\sum}%
}R_{i\text{a}j\text{a}}-3\overset{q}{\underset{\text{a,b=1}}{\sum}%
}(T_{\text{aa}i}T_{\text{bb}j}-T_{\text{ab}i}T_{\text{ab}j}%
)-3\overset{q}{\underset{\text{a,b=1}}{\sum}}(T_{\text{aa}j}T_{\text{bb}%
i}-T_{\text{ab}j}T_{\text{ab}i}](y_{0})\phi(y_{0})\qquad L_{2}\qquad
L_{21}\qquad\ \ \ \ \ $

$-\frac{1}{864}[2\varrho_{ij}+$ $\overset{q}{\underset{\text{a}=1}{4\sum}%
}R_{i\text{a}j\text{a}}-3\overset{q}{\underset{\text{a,b=1}}{\sum}%
}(T_{\text{aa}i}T_{\text{bb}j}-T_{\text{ab}i}T_{\text{ab}j}%
)-3\overset{q}{\underset{\text{a,b=1}}{\sum}}(T_{\text{aa}j}T_{\text{bb}%
i}-T_{\text{ab}j}T_{\text{ab}i}]^{2}(y_{0})\phi(y_{0})$

$-\frac{1}{288}[<H,j>](y_{0})\times\lbrack\{\nabla_{i}\varrho_{ij}%
-2\varrho_{ij}<H,i>+\overset{q}{\underset{\text{a}=1}{\sum}}(\nabla
_{i}R_{\text{a}i\text{a}j}-4R_{i\text{a}j\text{a}}<H,i>)\qquad L_{212}$

$+4\overset{q}{\underset{\text{a,b=1}}{\sum}}R_{i\text{a}j\text{b}%
}T_{\text{ab}i}+2\overset{q}{\underset{\text{a,b,c=1}}{\sum}}(T_{\text{aa}%
i}T_{\text{bb}j}T_{\text{cc}i}-3T_{\text{aa}i}T_{\text{bc}j}T_{\text{bc}%
i}+2T_{\text{ab}i}T_{\text{bc}j}T_{\text{ac}i})](y_{0})\phi(y_{0})$%
\qquad\qquad\qquad\qquad\qquad\ \ 

$-\frac{1}{288}[<H,j>](y_{0})\times\lbrack\nabla_{j}\varrho_{ii}-2\varrho
_{ij}<H,i>+\overset{q}{\underset{\text{a}=1}{\sum}}(\nabla_{j}R_{\text{a}%
i\text{a}i}-4R_{i\text{a}j\text{a}}<H,i>)$

$+4\overset{q}{\underset{\text{a,b=1}}{\sum}}R_{j\text{a}i\text{b}%
}T_{\text{ab}i}+2\overset{q}{\underset{\text{a,b,c=1}}{\sum}}(T_{\text{aa}%
j}T_{\text{bb}i}T_{\text{cc}i}-3T_{\text{aa}j}T_{\text{bc}i}T_{\text{bc}%
i}+2T_{\text{ab}j}T_{\text{bc}i}T_{\text{ac}i})](y_{0})\phi(y_{0})$

$-\frac{1}{288}[<H,j>](y_{0})\times\lbrack\nabla_{i}\varrho_{ij}-2\varrho
_{ii}<H,j>+\overset{q}{\underset{\text{a}=1}{\sum}}(\nabla_{i}R_{\text{a}%
i\text{a}j}-4R_{i\text{a}i\text{a}}<H,j>)$

$+4\overset{q}{\underset{\text{a,b=1}}{\sum}}R_{i\text{a}i\text{b}%
}T_{\text{ab}j}+2\overset{q}{\underset{\text{a,b,c}=1}{\sum}}(T_{\text{aa}%
i}T_{\text{bb}i}T_{\text{cc}j}-3T_{\text{aa}i}T_{\text{bc}i}T_{\text{bc}%
j}+2T_{\text{ab}i}T_{\text{bc}i}T_{\text{ac}j})](y_{0})\phi(y_{0})$

$-\frac{1}{3}[<H,j><H,k>](y_{0})R_{ijik}(y_{0})\phi(y_{0})-$ $\frac{5}%
{64}[<H,i>^{2}<H,j>^{2}](y_{0})\phi(y_{0})\qquad L_{213}$

$-\frac{1}{96}<H,i><H,j>$

$\times\lbrack2\varrho_{ij}+\overset{q}{\underset{\text{a}=1}{4\sum}%
}R_{i\text{a}j\text{a}}-3\overset{q}{\underset{\text{a,b=1}}{\sum}%
}(T_{\text{aa}i}T_{\text{bb}j}-T_{\text{ab}i}T_{\text{ab}j}%
)-3\overset{q}{\underset{\text{a,b=1}}{\sum}}(T_{\text{aa}j}T_{\text{bb}%
i}-T_{\text{ab}j}T_{\text{ab}i}](y_{0})\phi(y_{0})$

$-\frac{1}{96}<H,j>^{2}[\tau^{M}\ -3\tau^{P}+\ \underset{\text{a}%
=1}{\overset{\text{q}}{\sum}}\varrho_{\text{aa}}^{M}+$
$\overset{q}{\underset{\text{a},\text{b}=1}{\sum}}R_{\text{abab}}^{M}$
$](y_{0})\phi(y_{0})$

$+\frac{1}{288}<H,j>[\nabla_{i}\varrho_{ij}-2\varrho_{ij}%
<H,i>+\overset{q}{\underset{\text{a}=1}{\sum}}(\nabla_{i}R_{\text{a}%
i\text{a}j}-4R_{i\text{a}j\text{a}}<H,i>)$

$+4\overset{q}{\underset{\text{a,b=1}}{\sum}}R_{i\text{a}j\text{b}%
}T_{\text{ab}i}+2\overset{q}{\underset{\text{a,b,c=1}}{\sum}}(T_{\text{aa}%
i}T_{\text{bb}j}T_{\text{cc}i}-3T_{\text{aa}i}T_{\text{bc}j}T_{\text{bc}%
i}+2T_{\text{ab}i}T_{\text{bc}j}T_{\text{ac}i})](y_{0})\phi(y_{0})$%
\qquad\qquad\qquad\qquad\qquad\ \ 

$+\frac{1}{12}<H,j>[\nabla_{j}\varrho_{ii}-2\varrho_{ij}%
<H,i>+\overset{q}{\underset{\text{a}=1}{\sum}}(\nabla_{j}R_{\text{a}%
i\text{a}i}-4R_{i\text{a}j\text{a}}<H,i>)$

$+4\overset{q}{\underset{\text{a,b=1}}{\sum}}R_{j\text{a}i\text{b}%
}T_{\text{ab}i}+2\overset{q}{\underset{\text{a,b,c=1}}{\sum}}(T_{\text{aa}%
j}T_{\text{bb}i}T_{\text{cc}i}-3T_{\text{aa}j}T_{\text{bc}i}T_{\text{bc}%
i}+2T_{\text{ab}j}T_{\text{bc}i}T_{\text{ac}i})](y_{0})\phi(y_{0})$

$+\frac{1}{288}<H,j>[\nabla_{i}\varrho_{ij}-2\varrho_{ii}%
<H,j>+\overset{q}{\underset{\text{a}=1}{\sum}}(\nabla_{i}R_{\text{a}%
i\text{a}j}-4R_{i\text{a}i\text{a}}<H,j>)+4\overset{q}{\underset{\text{a,b=1}%
}{\sum}}R_{i\text{a}i\text{b}}T_{\text{ab}j}$

$+2\overset{q}{\underset{\text{a,b,c}=1}{\sum}}(T_{\text{aa}i}T_{\text{bb}%
i}T_{\text{cc}j}-3T_{\text{aa}i}T_{\text{bc}i}T_{\text{bc}j}+2T_{\text{ab}%
i}T_{\text{bc}i}T_{\text{ac}j})](y_{0})\phi(y_{0})$

$-\frac{1}{144}R_{jijk}(y_{0})$ \ $[<H,i><H,k>](y_{0})\phi(y_{0})\qquad\qquad
L_{22}$

$-\frac{1}{432}R_{jijk}(y_{0})[2\varrho_{ik}+$ $\overset{q}{\underset{\text{a}%
=1}{4\sum}}R_{i\text{a}k\text{a}}-3\overset{q}{\underset{\text{a,b=1}}{\sum}%
}(T_{\text{aa}i}T_{\text{bb}k}-T_{\text{ab}i}T_{\text{ab}k}%
)-3\overset{q}{\underset{\text{a,b=1}}{\sum}}(T_{\text{aa}k}T_{\text{bb}%
i}-T_{\text{ab}k}T_{\text{ab}i}](y_{0})\phi(y_{0})$

$+\frac{1}{144}<H,k>(y_{0})[\nabla_{j}$R$_{ijik}(y_{0})-\nabla_{i}$%
R$_{jijk}](y_{0})\phi(y_{0})$

$-\frac{5}{32}<H,i>^{2}(y_{0})<H,j>^{2}(y_{0})\phi(y_{0})\qquad\qquad
L_{23}\qquad L_{231}$

$-\frac{1}{48}<H,i>(y_{0})<H,j>(y_{0})$

$\times\lbrack2\varrho_{ij}+$ $\overset{q}{\underset{\text{a}=1}{4\sum}%
}R_{i\text{a}j\text{a}}-3\overset{q}{\underset{\text{a,b=1}}{\sum}%
}(T_{\text{aa}i}T_{\text{bb}j}-T_{\text{ab}i}T_{\text{ab}j}%
)-3\overset{q}{\underset{\text{a,b=1}}{\sum}}(T_{\text{aa}j}T_{\text{bb}%
i}-T_{\text{ab}j}T_{\text{ab}i}](y_{0})\phi(y_{0})$

$-\frac{1}{48}<H,i>^{2}(y_{0})[\varrho_{jj}+$ $\overset{q}{\underset{\text{a}%
=1}{2\sum}}R_{j\text{a}j\text{a}}-3\overset{q}{\underset{\text{a,b=1}}{\sum}%
}(T_{\text{aa}j}T_{\text{bb}j}-T_{\text{ab}j}T_{\text{ab}j})](y_{0})\phi
(y_{0})$

$-\frac{1}{144}<H,i>(y_{0})[\nabla_{i}\varrho_{jj}-2\varrho_{ij}%
<H,j>+\overset{q}{\underset{\text{a}=1}{\sum}}(\nabla_{i}R_{\text{a}%
j\text{a}j}-4R_{i\text{a}j\text{a}}<H,j>)+4\overset{q}{\underset{\text{a,b=1}%
}{\sum}}R_{i\text{a}j\text{b}}T_{\text{ab}j}$

$+2\overset{q}{\underset{\text{a,b,c=1}}{\sum}}(T_{\text{aa}i}T_{\text{bb}%
j}T_{\text{cc}j}-3T_{\text{aa}i}T_{\text{bc}j}T_{\text{bc}j}+2T_{\text{ab}%
i}T_{\text{bc}j}T_{\text{ca}j})](y_{0})\phi(y_{0})$\qquad\qquad\qquad
\qquad\qquad\ \ 

$-\frac{1}{24}\times\frac{1}{6}<H,i>(y_{0})[\nabla_{j}\varrho_{ij}%
-2\varrho_{ij}<H,j>+\overset{q}{\underset{\text{a}=1}{\sum}}(\nabla
_{j}R_{\text{a}i\text{a}j}-4R_{j\text{a}i\text{a}}%
<H,j>)+4\overset{q}{\underset{\text{a,b=1}}{\sum}}R_{j\text{a}i\text{b}%
}T_{\text{ab}j}$

$+2\overset{q}{\underset{\text{a,b,c=1}}{\sum}}(T_{\text{aa}j}T_{\text{bb}%
i}T_{\text{cc}j}-3T_{\text{aa}j}T_{\text{bc}i}T_{\text{bc}j}+2T_{\text{ab}%
j}T_{\text{bc}i}T_{\text{ac}j})](y_{0})\phi(y_{0})$

$-\frac{1}{144}<H,i>(y_{0})[\nabla_{j}\varrho_{ij}-2\varrho_{jj}%
<H,i>+\overset{q}{\underset{\text{a}=1}{\sum}}(\nabla_{j}R_{\text{a}%
i\text{a}j}-4R_{j\text{a}j\text{a}}<H,i>)+4\overset{q}{\underset{\text{a,b=1}%
}{\sum}}R_{j\text{a}j\text{b}}T_{\text{ab}i}$

$+2\overset{q}{\underset{\text{a,b,c=1}}{\sum}}(T_{\text{aa}j}T_{\text{bb}%
j}T_{\text{cc}i}-3T_{\text{aa}j}T_{\text{bc}j}T_{\text{bc}i}+2T_{\text{ab}%
j}T_{\text{bc}j}T_{\text{ac}i})](y_{0})\phi(y_{0})$

$-\frac{1}{96}<H,j>^{2}(y_{0})[\varrho_{ii}+$ $\overset{q}{\underset{\text{a}%
=1}{2\sum}}R_{i\text{a}i\text{a}}-3\overset{q}{\underset{\text{a,b=1}}{\sum}%
}(T_{\text{aa}i}T_{\text{bb}i}-T_{\text{ab}i}T_{\text{ab}i})](y_{0})\phi
(y_{0})\qquad L_{232}$

$-\frac{1}{432}[\varrho_{ii}+$ $\overset{q}{\underset{\text{a}=1}{2\sum}%
}R_{i\text{a}i\text{a}}-3\overset{q}{\underset{\text{a,b=1}}{\sum}%
}(T_{\text{aa}i}T_{\text{bb}i}-T_{\text{ab}i}T_{\text{ab}i})](y_{0})\phi
(y_{0})$

\ $\times\lbrack\varrho_{jj}+$ $\overset{q}{\underset{\text{a}=1}{2\sum}%
}R_{j\text{a}j\text{a}}-3\overset{q}{\underset{\text{a,b=1}}{\sum}%
}(T_{\text{aa}j}T_{\text{bb}j}-T_{\text{ab}j}T_{\text{ab}j})]\}(y_{0}%
)\phi(y_{0})$

$+\frac{1}{48}R_{ijik}(y_{0})$ \ $[<H,j><H,k>](y_{0})\phi(y_{0})\qquad
\qquad\qquad\qquad$\ $L_{233}$

$+\frac{1}{432}R_{ijik}(y_{0})[2\varrho_{jk}+$ $\overset{q}{\underset{\text{a}%
=1}{4\sum}}R_{j\text{a}k\text{a}}-3\overset{q}{\underset{\text{a,b=1}}{\sum}%
}(T_{\text{aa}j}T_{\text{bb}k}-T_{\text{ab}j}T_{\text{ab}k}%
)-3\overset{q}{\underset{\text{a,b=1}}{\sum}}(T_{\text{aa}k}T_{\text{bb}%
j}-T_{\text{ab}k}T_{\text{ab}j}](y_{0})\phi(y_{0})$

$+\overset{n}{\underset{i,j=q+1}{\sum}}\frac{35}{128}<H,i>^{2}(y_{0}%
)<H,j>^{2}(y_{0})\phi(y_{0})\qquad\qquad\ \frac{1}{24}\frac{\partial^{4}%
\theta^{-\frac{1}{2}}}{\partial x_{i}^{2}\partial x_{j}^{2}}(y_{0})$

$+\frac{5}{192}\overset{n}{\underset{j=q+1}{\sum}}<H,j>^{2}(y_{0})[\tau
^{M}\ -3\tau^{P}+\ \underset{\text{a}=1}{\overset{\text{q}}{\sum}}%
\varrho_{\text{aa}}^{M}+\overset{q}{\underset{\text{a},\text{b}=1}{\sum}%
}R_{\text{abab}}^{M}](y_{0})\phi(y_{0})\qquad\ \ \ \ \ \ \ \ $

$+\frac{5}{192}\overset{n}{\underset{i=q+1}{\sum}}<H,i>^{2}(y_{0})[\tau
^{M}\ -3\tau^{P}+\ \underset{\text{a}=1}{\overset{\text{q}}{\sum}}%
\varrho_{\text{aa}}^{M}+\overset{q}{\underset{\text{a},\text{b}=1}{\sum}%
}R_{\text{abab}}^{M}](y_{0})\phi(y_{0})\qquad\qquad$

$+\frac{5}{192}\overset{n}{\underset{i,j=q+1}{\sum}}[<H,i><H,j>](y_{0}%
)\qquad\qquad\qquad\qquad\qquad\qquad\qquad\qquad$

$\times\lbrack2\varrho_{ij}+4\overset{q}{\underset{\text{a}=1}{\sum}%
}R_{i\text{a}j\text{a}}-3\overset{q}{\underset{\text{a,b=1}}{\sum}%
}(T_{\text{aa}i}T_{\text{bb}j}-T_{\text{ab}i}T_{\text{ab}j}%
)-3\overset{q}{\underset{\text{a,b=1}}{\sum}}(T_{\text{aa}j}T_{\text{bb}%
i}-T_{\text{ab}j}T_{\text{ab}i})](y_{0})\phi(y_{0})$

$+\frac{1}{96}\overset{n}{\underset{i,j=q+1}{\sum}}<H,j>(y_{0})[\{\nabla
_{i}\varrho_{ij}-2\varrho_{ij}<H,i>+\overset{q}{\underset{\text{a}=1}{\sum}%
}(\nabla_{i}R_{\text{a}i\text{a}j}-4R_{i\text{a}j\text{a}}<H,i>)\qquad$

$+4\overset{q}{\underset{\text{a,b=1}}{\sum}}R_{i\text{a}j\text{b}%
}T_{\text{ab}i}+2\overset{q}{\underset{\text{a,b,c=1}}{\sum}}(T_{\text{aa}%
i}T_{\text{bb}j}T_{\text{cc}i}-T_{\text{aa}i}T_{\text{bc}j}T_{\text{bc}%
i}-2T_{\text{bc}j}(T_{\text{aa}i}T_{\text{bc}i}-T_{\text{ab}i}T_{\text{ac}%
i}))\}$\qquad\qquad\qquad\ \ 

$+\{\nabla_{j}\varrho_{ii}-2\varrho_{ij}<H,i>+\overset{q}{\underset{\text{a}%
=1}{\sum}}(\nabla_{j}R_{\text{a}i\text{a}i}-4R_{i\text{a}j\text{a}}<H,i>)$

$+4\overset{q}{\underset{\text{a,b=1}}{\sum}}R_{j\text{a}i\text{b}%
}T_{\text{ab}i}+2\overset{q}{\underset{\text{a,b,c=1}}{\sum}}(T_{\text{aa}%
j}(T_{\text{bb}i}T_{\text{cc}i}-T_{\text{bc}i}T_{\text{bc}i})-2T_{\text{aa}%
j}T_{\text{bc}i}T_{\text{bc}i}+2T_{\text{ab}j}T_{\text{bc}i}T_{\text{ac}%
i})\}\qquad$

$+\{\nabla_{i}\varrho_{ij}-2\varrho_{ii}<H,j>+\overset{q}{\underset{\text{a}%
=1}{\sum}}(\nabla_{i}R_{\text{a}i\text{a}j}-4R_{i\text{a}i\text{a}%
}<H,j>)+4\overset{q}{\underset{\text{a,b=1}}{\sum}}R_{i\text{a}i\text{b}%
}T_{\text{ab}j}$

$+2\overset{q}{\underset{\text{a,b,c}=1}{\sum}}(T_{\text{aa}i}T_{\text{bb}%
i}T_{\text{cc}j}-3T_{\text{aa}i}T_{\text{bc}i}T_{\text{bc}j}+2T_{\text{ab}%
i}T_{\text{bc}i}T_{\text{ac}j})\}](y_{0})\phi(y_{0})$

$+\frac{1}{96}\overset{n}{\underset{i,j=q+1}{\sum}}<H,i>(y_{0})[\{\nabla
_{i}\varrho_{jj}-2\varrho_{ij}<H,j>+\overset{q}{\underset{\text{a}=1}{\sum}%
}(\nabla_{i}R_{\text{a}j\text{a}j}-4R_{i\text{a}j\text{a}}<H,j>)\qquad$

$+4\overset{q}{\underset{\text{a,b=1}}{\sum}}R_{i\text{a}j\text{b}%
}T_{\text{ab}j}+2\overset{q}{\underset{\text{a,b,c=1}}{\sum}}T_{\text{aa}%
i}(T_{\text{bb}j}T_{\text{cc}j}-T_{\text{bc}j}T_{\text{bc}j})-2T_{\text{aa}%
i}T_{\text{bc}j}T_{\text{bc}j}+2T_{\text{ab}i}T_{\text{bc}j}T_{\text{ac}%
j})\}(y_{0})\qquad$\qquad\qquad\qquad\qquad\qquad\ \ 

$+\{\nabla_{j}\varrho_{ij}-2\varrho_{ij}<H,j>+\overset{q}{\underset{\text{a}%
=1}{\sum}}(\nabla_{j}R_{\text{a}i\text{a}j}-4R_{j\text{a}i\text{a}}<H,j>)$

$+4\overset{q}{\underset{\text{a,b=1}}{\sum}}R_{j\text{a}i\text{b}%
}T_{\text{ab}j}+2\overset{q}{\underset{\text{a,b,c=1}}{\sum}}(T_{\text{aa}%
j}T_{\text{bb}i}T_{\text{cc}j}-T_{\text{ab}j}T_{\text{bc}i}T_{\text{ac}%
j}-2T_{\text{bc}i}(T_{\text{aa}j}T_{\text{bc}j}-T_{\text{ab}j}T_{\text{ac}%
j}))\}(y_{0})$

$+\{\nabla_{j}\varrho_{ij}-2\varrho_{jj}<H,i>+\overset{q}{\underset{\text{a}%
=1}{\sum}}(\nabla_{j}R_{\text{a}i\text{a}j}-4R_{j\text{a}j\text{a}%
}<H,i>)+4\overset{q}{\underset{\text{a,b=1}}{\sum}}R_{j\text{a}j\text{b}%
}T_{\text{ab}i}$

$+2\overset{q}{\underset{\text{a,b,c=1}}{\sum}}(T_{\text{aa}j}T_{\text{bb}%
j}T_{\text{cc}i}-3T_{\text{aa}j}T_{\text{bc}j}T_{\text{bc}i}+2T_{\text{ab}%
j}T_{\text{bc}j}T_{\text{ac}i})\}](y_{0})\phi(y_{0})$

$+\frac{1}{576}\overset{n}{\underset{i,j=q+1}{\sum}}[2\varrho_{ij}%
+4\overset{q}{\underset{\text{a}=1}{\sum}}R_{i\text{a}j\text{a}}%
-3\overset{q}{\underset{\text{a,b=1}}{\sum}}(T_{\text{aa}i}T_{\text{bb}%
j}-T_{\text{ab}i}T_{\text{ab}j})-3\overset{q}{\underset{\text{a,b=1}}{\sum}%
}(T_{\text{aa}j}T_{\text{bb}i}-T_{\text{ab}j}T_{\text{ab}i})]^{2}(y_{0}%
)\phi(y_{0})$

$+\frac{1}{288}[\tau^{M}\ -3\tau^{P}+\ \underset{\text{a}=1}{\overset{\text{q}%
}{\sum}}\varrho_{\text{aa}}^{M}+\overset{q}{\underset{\text{a},\text{b}%
=1}{\sum}}R_{\text{abab}}^{M}]^{2}(y_{0})\phi(y_{0})$

$-\ \frac{1}{288}\overset{n}{\underset{i,j=q+1}{\sum}}[$
$\overset{q}{\underset{\text{a=1}}{\sum}}\{-(\nabla_{ii}^{2}R_{j\text{a}%
j\text{a}}+\nabla_{jj}^{2}R_{i\text{a}i\text{a}}+4\nabla_{ij}^{2}%
R_{i\text{a}j\text{a}}+2R_{ij}R_{i\text{a}j\text{a}})\qquad A$

$+\overset{n}{\underset{p=q+1}{\sum}}\overset{q}{\underset{\text{a=1}}{\sum}%
}(R_{\text{a}iip}R_{\text{a}jjp}+R_{\text{a}jjp}R_{\text{a}iip}+R_{\text{a}%
ijp}R_{\text{a}ijp}+R_{\text{a}ijp}R_{\text{a}jip}+R_{\text{a}jip}%
R_{\text{a}ijp}+R_{\text{a}jip}R_{\text{a}jip})$

$+2\overset{q}{\underset{\text{a,b=1}}{\sum}}\nabla_{i}(R)_{\text{a}%
i\text{b}j}T_{\text{ab}j}+2\overset{q}{\underset{\text{a,b=1}}{\sum}}%
\nabla_{j}(R)_{\text{a}j\text{b}i}T_{\text{ab}i}%
+2\overset{q}{\underset{\text{a,b=1}}{\sum}}\nabla_{i}(R)_{\text{a}j\text{b}%
i}T_{\text{ab}j}+2\overset{q}{\underset{\text{a,b=1}}{\sum}}\nabla
_{i}(R)_{\text{a}j\text{b}j}T_{\text{ab}i}$

$+2\overset{q}{\underset{\text{a,b=1}}{\sum}}\nabla_{j}(R)_{\text{a}%
i\text{b}i}T_{\text{ab}j}+2\overset{q}{\underset{\text{a,b=1}}{\sum}}%
\nabla_{j}(R)_{\text{a}i\text{b}j}T_{\text{ab}i}$

$+\overset{n}{\underset{p=q+1}{\sum}}(-\frac{3}{5}\nabla_{ii}^{2}%
(R)_{jpjp}+\overset{n}{\underset{p=q+1}{\sum}}(-\frac{3}{5}\nabla_{jj}%
^{2}(R)_{ipip}+\overset{n}{\underset{p=q+1}{\sum}}(-\frac{3}{5}\nabla_{ij}%
^{2}(R)_{ipjp}+\overset{n}{\underset{p=q+1}{\sum}}(-\frac{3}{5}\nabla_{ij}%
^{2}(R)_{jpip}$

$+\overset{n}{\underset{p=q+1}{\sum}}(-\frac{3}{5}\nabla_{ji}^{2}%
(R)_{ipjp}+\overset{n}{\underset{p=q+1}{\sum}}(-\frac{3}{5}\nabla_{ji}%
^{2}(R)_{jpip}$

$+\frac{1}{5}\overset{n}{\underset{m,p=q+1}{%
{\textstyle\sum}
}}R_{ipim}R_{jpjm}+\frac{1}{5}\overset{n}{\underset{m,p=q+1}{%
{\textstyle\sum}
}}R_{jpjm}R_{ipim}+\frac{1}{5}\overset{n}{\underset{m,p=q+1}{%
{\textstyle\sum}
}}R_{ipjm}R_{ipjm}+\frac{1}{5}\overset{n}{\underset{m,p=q+1}{%
{\textstyle\sum}
}}R_{ipjm}R_{jpim}$

$+\frac{1}{5}\overset{n}{\underset{m,p=q+1}{%
{\textstyle\sum}
}}R_{jpim}R_{ipjm}+\frac{1}{5}\overset{n}{\underset{m,p=q+1}{%
{\textstyle\sum}
}}R_{jpim}R_{jpim}\}(y_{0})$

$+4\overset{q}{\underset{\text{a,b=1}}{\sum}}\{(\nabla_{i}(R)_{i\text{a}%
j\text{a}}-\overset{q}{\underset{\text{c=1}}{%
{\textstyle\sum}
}}R_{\text{a}i\text{c}i}T_{\text{ac}j})$ $T_{\text{bb}j}+4(\nabla
_{j}(R)_{j\text{a}i\text{a}}-\overset{q}{\underset{\text{c=1}}{%
{\textstyle\sum}
}}R_{\text{a}j\text{c}j}T_{\text{ac}i})$ $T_{\text{bb}i}+4(\nabla
_{i}(R)_{j\text{a}i\text{a}}-\overset{q}{\underset{\text{c=1}}{%
{\textstyle\sum}
}}R_{\text{a}i\text{c}j}T_{\text{ac}i})$ $T_{\text{bb}j}$ $4B\ $

$+4(\nabla_{i}(R)_{j\text{a}j\text{a}}-\overset{q}{\underset{\text{c=1}}{%
{\textstyle\sum}
}}R_{\text{a}i\text{c}j}T_{\text{ac}j})$ $T_{\text{bb}i}+4(\nabla
_{j}(R)_{i\text{a}i\text{a}}-\overset{q}{\underset{\text{c=1}}{%
{\textstyle\sum}
}}R_{\text{a}j\text{c}i}T_{\text{ac}i})$ $T_{\text{bb}j}+4(\nabla
_{j}(R)_{i\text{a}j\text{a}}-\overset{q}{\underset{\text{c=1}}{%
{\textstyle\sum}
}}R_{\text{a}j\text{c}i}T_{\text{ac}j})$ $T_{\text{bb}i}$

$-4\overset{q}{\underset{\text{a,b=1}}{\sum}}(\nabla_{i}(R)_{i\text{a}%
j\text{b}}-\overset{q}{\underset{\text{c=1}}{%
{\textstyle\sum}
}}R_{\text{b}r\text{c}s}T_{\text{ac}t})T_{\text{ab}j}%
-4\overset{q}{\underset{\text{a,b=1}}{\sum}}(\nabla_{j}(R)_{j\text{a}%
i\text{b}}-\overset{q}{\underset{\text{c=1}}{%
{\textstyle\sum}
}}R_{\text{b}j\text{c}j}T_{\text{ac}i})T_{\text{ab}i}$

$-4\overset{q}{\underset{\text{a,b=1}}{\sum}}(\nabla_{i}(R)_{j\text{a}%
i\text{b}}-\overset{q}{\underset{\text{c=1}}{%
{\textstyle\sum}
}}R_{\text{b}i\text{c}j}T_{\text{ac}i})T_{\text{ab}j}%
-4\overset{q}{\underset{\text{a,b=1}}{\sum}}(\nabla_{i}(R)_{j\text{a}%
j\text{b}}-\overset{q}{\underset{\text{c=1}}{%
{\textstyle\sum}
}}R_{\text{b}i\text{c}j}T_{\text{ac}j})T_{\text{ab}i}$

$-4\overset{q}{\underset{\text{a,b=1}}{\sum}}(\nabla_{j}(R)_{i\text{a}%
i\text{b}}-\overset{q}{\underset{\text{c=1}}{%
{\textstyle\sum}
}}R_{\text{b}j\text{c}i}T_{\text{ac}i})T_{\text{ab}j}%
-4\overset{q}{\underset{\text{a,b=1}}{\sum}}(\nabla_{j}(R)_{i\text{a}%
j\text{b}}-\overset{q}{\underset{\text{c=1}}{%
{\textstyle\sum}
}}R_{\text{b}j\text{c}i}T_{\text{ac}j})T_{\text{ab}i}\}](y_{0})$

$-\frac{1}{48}$ $[\frac{4}{9}\overset{q}{\underset{\text{a,b=1}}{\sum}%
}(\varrho_{\text{aa}}-\overset{q}{\underset{\text{c}=1}{\sum}}R_{\text{acac}%
})(\varrho_{\text{bb}}-\overset{q}{\underset{\text{d}=1}{\sum}}R_{\text{bdbd}%
})+\frac{8}{9}\overset{n}{\underset{i,j=q+1}{\sum}}%
\overset{q}{\underset{\text{a,b}=1}{\sum}}(R_{i\text{a}j\text{a}}%
R_{i\text{b}j\text{b}})\qquad3C$

$+\frac{2}{9}\overset{q}{\underset{\text{a}=1}{\sum}}(\varrho_{\text{aa}}%
^{M}-\varrho_{\text{aa}}^{P})(\tau^{M}-\overset{q}{\underset{\text{c}=1}{\sum
}}\varrho_{\text{cc}}^{M})+\frac{4}{9}\overset{n}{\underset{i,j=q+1}{\sum}%
}\overset{q}{\underset{\text{a}=1}{\sum}}R_{i\text{a}j\text{a}}\varrho_{ij}\ $

$\ +\frac{2}{9}\overset{q}{\underset{\text{b}=1}{\sum}}(\varrho_{\text{bb}%
}^{M}-\varrho_{\text{bb}}^{P})(\tau^{M}-\overset{q}{\underset{\text{c}%
=1}{\sum}}\varrho_{\text{cc}}^{M})+\frac{4}{9}%
\overset{n}{\underset{i,j=q+1}{\sum}}\overset{q}{\underset{\text{b}=1}{\sum}%
}R_{i\text{b}j\text{b}}\varrho_{ij}\ $

$+\frac{1}{9}(\tau^{M}-\overset{q}{\underset{\text{a=1}}{\sum}}\varrho
_{\text{aa}})(\tau^{M}-\overset{q}{\underset{\text{b=1}}{\sum}}\varrho
_{\text{bb}})+\frac{2}{9}(\left\Vert \varrho^{M}\right\Vert ^{2}%
-\overset{q}{\underset{\text{a,b}=1}{\sum}}\varrho_{\text{ab}})$

$-\overset{n}{\underset{i,j=q+1}{\sum}}\overset{q}{\underset{\text{a,b}%
=1}{\sum}}R_{i\text{a}i\text{b}}R_{j\text{a}j\text{b}}\ -\frac{1}%
{2}\overset{n}{\underset{i,j=q+1}{\sum}}\overset{q}{\underset{\text{a,b}%
=1}{\sum}}R_{i\text{a}j\text{b}}^{2}-\overset{n}{\underset{i,j=q+1}{\sum}%
}\overset{q}{\underset{\text{a,b}=1}{\sum}}R_{i\text{a}j\text{b}}%
R_{j\text{a}i\text{b}}-\frac{1}{2}\overset{n}{\underset{i,j=q+1}{\sum}%
}\overset{q}{\underset{\text{a,b}=1}{\sum}}R_{j\text{a}i\text{b}}^{2}$

$-\frac{1}{9}\overset{n}{\underset{i,j,p,m=q+1}{\sum}}R_{ipim}R_{jpjm}%
\ -\frac{1}{18}\overset{n}{\underset{i,j,p,m=q+1}{\sum}}R_{ipjm}^{2}-\frac
{1}{9}\overset{n}{\underset{i,j,p,m=q+1}{\sum}}R_{ipjm}R_{jpim}-\frac{1}%
{18}\overset{n}{\underset{i,j,p,m=q+1}{\sum}}R_{jpim}^{2}$

$-\frac{1}{3}\overset{q}{\underset{\text{a}=1}{\sum}}%
\overset{n}{\underset{i,j,p=q+1}{\sum}}R_{i\text{a}ip}R_{j\text{a}jp}-\frac
{1}{6}\overset{q}{\underset{\text{a}=1}{\sum}}%
\overset{n}{\underset{i,j,p=q+1}{\sum}}R_{i\text{a}jp}^{2}-\frac{1}%
{3}\overset{q}{\underset{\text{a}=1i,j,}{\sum}}%
\overset{n}{\underset{p=q+1}{\sum}}R_{i\text{a}jp}R_{j\text{a}ip}-\frac{1}%
{6}\overset{q}{\underset{\text{a}=1}{\sum}}%
\overset{n}{\underset{i,j,p=q+1}{\sum}}R_{j\text{a}ip}^{2}$

$-\frac{1}{3}\overset{q}{\underset{\text{b}=1i,j,}{\sum}}%
\overset{n}{\underset{p=q+1}{\sum}}R_{i\text{b}ip}R_{j\text{b}jp}-\frac{1}%
{6}\overset{q}{\underset{\text{b}=1}{\sum}}%
\overset{n}{\underset{i,j,p=q+1}{\sum}}R_{i\text{b}jp}^{2}-\frac{1}%
{3}\overset{q}{\underset{\text{b}=1}{\sum}}%
\overset{n}{\underset{i.j,p=q+1}{\sum}}R_{i\text{b}jp}R_{j\text{b}ip}-\frac
{1}{6}\overset{q}{\underset{\text{b}=1}{\sum}}%
\overset{n}{\underset{i,j,p=q+1}{\sum}}R_{j\text{b}ip}^{2}](y_{0})\phi(y_{0})$

$-\frac{1}{48}$ $\overset{q}{\underset{\text{a,b,c=1}}{\sum}}[$
$-\overset{n}{\underset{i=q+1}{\sum}}R_{i\text{a}i\text{a}}(R_{\text{bcbc}%
}^{P}-R_{\text{bcbc}}^{M})$ $-\overset{n}{\underset{j=q+1}{\sum}}%
R_{j\text{a}j\text{a}}(R_{\text{bcbc}}^{P}-R_{\text{bcbc}}^{M})\qquad\qquad6D$

$+\overset{n}{\underset{i=q+1}{\sum}}R_{i\text{a}i\text{b}}(R_{\text{acbc}%
}^{P}-R_{\text{acbc}}^{M})\ -\overset{n}{\underset{i=q+1}{\sum}}%
R_{i\text{a}i\text{c}}(R_{\text{abbc}}^{P}-R_{\text{abbc}}^{M})$

$+\overset{n}{\underset{j=q+1}{\sum}}R_{j\text{a}j\text{b}}(R_{\text{acbc}%
}^{P}-R_{\text{acbc}}^{M})$\ $-\overset{n}{\underset{j=q+1}{\sum}}%
R_{j\text{a}j\text{c}}(R_{\text{abbc}}^{P}-R_{\text{abbc}}^{M})$

$+\underset{i,j=q+1}{\overset{n}{\sum}}$ $-R_{i\text{a}j\text{a}}%
(T_{\text{bb}i}T_{\text{cc}j}$ $-T_{\text{bc}i}T_{\text{bc}j})$
$-\underset{i,j=q+1}{\overset{n}{\sum}}R_{i\text{a}j\text{a}}(T_{\text{bb}%
j}T_{\text{cc}i}$ $-T_{\text{bc}j}T_{\text{bc}i})$

$+$ $\underset{i,j=q+1}{\overset{n}{\sum}}$ $-R_{j\text{a}i\text{a}%
}(T_{\text{bb}i}T_{\text{cc}j}$ $-T_{\text{bc}i}T_{\text{bc}j})$
$-\underset{i,j=q+1}{\overset{n}{\sum}}R_{j\text{a}i\text{a}}(T_{\text{bb}%
j}T_{\text{cc}i}$ $-T_{\text{bc}j}T_{\text{bc}i})$

$+\underset{i,j=q+1}{\overset{n}{\sum}}\ R_{i\text{a}j\text{b}}(T_{\text{ab}%
i}T_{\text{cc}j}-T_{\text{bc}i}T_{\text{ac}j}%
)\ +\underset{i,j=q+1}{\overset{n}{\sum}}\ R_{i\text{a}j\text{b}}%
(T_{\text{ab}j}T_{\text{cc}i}-T_{\text{bc}j}T_{\text{ac}i})$

$+\underset{i,j=q+1}{\overset{n}{\sum}}\ R_{j\text{a}i\text{ib}}%
(T_{\text{ab}i}T_{\text{cc}j}-T_{\text{bc}i}T_{\text{ac}j}%
)\ +\underset{i,j=q+1}{\overset{n}{\sum}}\ R_{j\text{a}i\text{b}}%
(T_{\text{ab}j}T_{\text{cc}i}-T_{\text{bc}j}T_{\text{ac}i})\qquad$

$+\underset{i,j=q+1}{\overset{n}{\sum}}-R_{i\text{a}j\text{c}}(T_{\text{ab}%
i}T_{\text{bc}j}-T_{\text{ac}i}T_{\text{bb}j}%
)-\underset{i,j=q+1}{\overset{n}{\sum}}R_{i\text{a}j\text{c}}(T_{\text{ba}%
j}T_{\text{bc}i}-T_{\text{ac}j}T_{\text{bb}i})$

$+\underset{i,j=q+1}{\overset{n}{\sum}}-R_{j\text{a}i\text{c}}(T_{\text{ba}%
i}T_{\text{bc}j}-T_{\text{ac}i}T_{\text{bb}j}%
)-\underset{i,j=q+1}{\overset{n}{\sum}}R_{j\text{a}i\text{c}}(T_{\text{ba}%
j}T_{\text{bc}i}-T_{\text{ac}j}T_{\text{bb}i})](y_{0})\phi(y_{0})$

$+\frac{1}{144}\underset{p=q+1}{\overset{n}{\sum}}%
[\underset{i=q+1}{\overset{n}{\sum}}\overset{q}{\underset{\text{b,c=1}}{\sum}%
}R_{ipip}(R_{\text{bcbc}}^{P}-R_{\text{bcbc}}^{M}%
)+\underset{j=q+1}{\overset{n}{\sum}}$ $\overset{q}{\underset{\text{b,c=1}%
}{\sum}}R_{jpjp}(R_{\text{bcbc}}^{P}-R_{\text{bcbc}}^{M})](y_{0})\phi(y_{0})$

$+\frac{1}{72}\underset{i,j,p=q+1}{\overset{n}{\sum}}%
\overset{q}{\underset{\text{b,c=1}}{\sum}}[R_{ipjp}(T_{\text{bb}i}%
T_{\text{cc}j}-T_{\text{bc}i}T_{\text{bc}j})+R_{ipjp}(T_{\text{bb}%
j}T_{\text{cc}i}-T_{\text{bc}j}T_{\text{bc}i})](y_{0})\phi(y_{0})\qquad$

$-\frac{1}{288}\underset{i,j=q+1}{\overset{n}{\sum}}[T_{\text{aa}%
i}T_{\text{bb}j}(T_{\text{cc}i}T_{\text{dd}j}-T_{\text{cd}i}T_{\text{dc}%
j})+T_{\text{aa}i}T_{\text{bb}j}(T_{\text{cc}j}T_{\text{dd}i}-T_{\text{cd}%
j}T_{\text{dc}i})\qquad E$

$+T_{\text{aa}j}T_{\text{bb}i}(T_{\text{cc}i}T_{\text{dd}j}-T_{\text{cd}%
i}T_{\text{dc}j})+T_{\text{aa}j}T_{\text{bb}i}(T_{\text{cc}j}T_{\text{dd}%
i}-T_{\text{cd}j}T_{\text{dc}i})](y_{0})\phi(y_{0})$

$+\frac{1}{288}\underset{i,j=q+1}{\overset{n}{\sum}}[T_{\text{aa}%
i}T_{\text{bc}j}(T_{\text{bc}i}T_{\text{dd}j}-T_{\text{bd}i}T_{\text{cd}%
j})+T_{\text{aa}i}T_{\text{bc}j}(T_{\text{bc}j}T_{\text{dd}i}-T_{\text{bd}%
j}T_{\text{cd}i})$

$+T_{\text{aa}j}T_{\text{bc}i}(T_{\text{bc}i}T_{\text{dd}j}-T_{\text{bd}%
i}T_{\text{cd}j})+T_{\text{aa}j}T_{\text{bc}i}(T_{\text{bc}j}T_{\text{dd}%
i}-T_{\text{bd}j}T_{\text{cd}i})](y_{0})\phi(y_{0})$

$-\frac{1}{288}\underset{i,j=q+1}{\overset{n}{\sum}}[T_{\text{aa}%
i}T_{\text{bd}j}(T_{\text{bc}i}T_{\text{cd}j}-T_{\text{bd}i}T_{\text{cc}%
j})+T_{\text{aa}i}T_{\text{bd}j}(T_{\text{bc}j}T_{\text{cd}i}-T_{\text{bd}%
j}T_{\text{cc}i})$

$+T_{\text{aa}j}T_{\text{bd}i}(T_{\text{bc}i}T_{\text{cd}j}-T_{\text{bd}%
i}T_{\text{cc}j})+T_{\text{aa}j}T_{\text{bd}i}(T_{\text{bc}j}T_{\text{cd}%
i}-T_{\text{bd}j}T_{\text{cc}i})](y_{0})\phi(y_{0})\qquad$

$+\frac{1}{288}\underset{i,j=q+1}{\overset{n}{\sum}}[T_{\text{ab}%
i}T_{\text{ab}j}(T_{\text{cc}i}T_{\text{dd}j}-T_{\text{cd}i}T_{\text{dc}%
j})+T_{\text{ab}i}T_{\text{ab}j}(T_{\text{cc}j}T_{\text{dd}i}-T_{\text{cd}%
j}T_{\text{dc}i})$

$+T_{\text{ab}j}T_{\text{ab}i}(T_{\text{cc}i}T_{\text{dd}j}-T_{\text{cd}%
i}T_{\text{dc}j})+T_{\text{ab}j}T_{\text{ab}i}(T_{\text{cc}j}T_{\text{dd}%
i}-T_{\text{cd}j}T_{\text{dc}i})](y_{0})\phi(y_{0})$

$-\frac{1}{288}\underset{i,j=q+1}{\overset{n}{\sum}}[T_{\text{ab}%
i}T_{\text{bc}j}(T_{\text{ac}i}T_{\text{dd}j}-T_{\text{ad}i}T_{\text{cd}%
j})+T_{\text{ab}i}T_{\text{bc}j}(T_{\text{ac}j}T_{\text{dd}i}-T_{\text{ad}%
j}T_{\text{cd}i})$

$+T_{\text{ab}j}T_{\text{bc}i}(T_{\text{ac}i}T_{\text{dd}j}-T_{\text{ad}%
i}T_{\text{cd}j})+T_{\text{ab}j}T_{\text{bc}i}(T_{\text{ac}j}T_{\text{dd}%
i}-T_{\text{ad}j}T_{\text{cd}i})](y_{0})\phi(y_{0})$

$+\frac{1}{288}\underset{i,j=q+1}{\overset{n}{\sum}}[T_{\text{ab}%
i}T_{\text{bd}j}(T_{\text{ac}i}T_{\text{cd}j}-T_{\text{ad}i}T_{\text{cc}%
j})+T_{\text{ab}i}T_{\text{bd}j}(T_{\text{ac}j}T_{\text{cd}i}-T_{\text{ad}%
j}T_{\text{cc}i})$

$+T_{\text{ab}i}T_{\text{bd}j}(T_{\text{ac}j}T_{\text{cd}i}-T_{\text{ad}%
j}T_{\text{cc}i})+T_{\text{ab}j}T_{\text{bd}i}(T_{\text{ac}j}T_{\text{cd}%
i}-T_{\text{ad}j}T_{\text{cc}i})](y_{0})\phi(y_{0})$

$-\ \frac{1}{288}\underset{i,j=q+1}{\overset{n}{\sum}}[T_{\text{ac}%
i}T_{\text{ab}j}(T_{\text{bc}i}T_{\text{dd}j}-T_{\text{bd}i}T_{\text{dc}%
j})+T_{\text{ac}i}T_{\text{ab}j}(T_{\text{bc}j}T_{\text{dd}i}-T_{\text{bd}%
j}T_{\text{dc}i})$

$+T_{\text{ac}j}T_{\text{ab}i}(T_{\text{bc}i}T_{\text{dd}j}-T_{\text{bd}%
i}T_{\text{dc}j})+T_{\text{ac}j}T_{\text{ab}i}(T_{\text{bc}j}T_{\text{dd}%
i}-T_{\text{bd}j}T_{\text{dc}i})](y_{0})\phi(y_{0})$

$+\ \frac{1}{288}\underset{i,j=q+1}{\overset{n}{\sum}}[T_{\text{ac}%
i}T_{\text{bb}j}(T_{\text{ac}i}T_{\text{dd}j}-T_{\text{ad}i}T_{\text{cd}%
j})+T_{\text{ac}i}T_{\text{bb}j}(T_{\text{ac}j}T_{\text{dd}i}-T_{\text{ad}%
j}T_{\text{cd}i})$

$+T_{\text{ac}j}T_{\text{bb}i}(T_{\text{ac}i}T_{\text{dd}j}-T_{\text{ad}%
i}T_{\text{cd}i})+T_{\text{ac}j}T_{\text{bb}i}(T_{\text{ac}j}T_{\text{dd}%
i}-T_{\text{ad}j}T_{\text{cd}i})](y_{0})\phi(y_{0})$

$-\ \frac{1}{288}\underset{i,j=q+1}{\overset{n}{\sum}}[T_{\text{ac}%
i}T_{\text{bd}j}(T_{\text{ac}i}T_{\text{bd}j}-T_{\text{ad}i}T_{\text{bc}%
j})+T_{\text{ac}i}T_{\text{bd}j}(T_{\text{ac}j}T_{\text{bd}i}-T_{\text{ad}%
j}T_{\text{bc}i})$

$+T_{\text{ac}j}T_{\text{bd}i}(T_{\text{ac}i}T_{\text{bd}j}-T_{\text{ad}%
i}T_{\text{bc}j})+T_{\text{ac}j}T_{\text{bd}i}(T_{\text{ac}j}T_{\text{bd}%
i}-T_{\text{ad}j}T_{\text{bc}i})](y_{0})\phi(y_{0})$

$+\frac{1}{288}\underset{i,j=q+1}{\overset{n}{\sum}}[T_{\text{ad}%
i}T_{\text{ab}j}(T_{\text{bc}i}T_{\text{cd}j}-T_{\text{bd}i}T_{\text{cc}%
j})+T_{\text{ad}i}T_{\text{ab}j}(T_{\text{bc}j}T_{\text{cd}i}-T_{\text{bd}%
j}T_{\text{cc}i})$

$+T_{\text{ad}j}T_{\text{ab}i}(T_{\text{bc}i}T_{\text{cd}j}-T_{\text{bd}%
i}T_{\text{cc}j})+T_{\text{ad}j}T_{\text{ab}i}(T_{\text{bc}j}T_{\text{cd}%
i}-T_{\text{bd}j}T_{\text{cc}i})](y_{0})\phi(y_{0})$

$-\ \frac{1}{288}\underset{i,j=q+1}{\overset{n}{\sum}}[T_{\text{ad}%
i}T_{\text{bb}j}(T_{\text{ac}i}T_{\text{cd}j}-T_{\text{ad}i}T_{\text{cc}%
j})+T_{\text{ad}i}T_{\text{bb}j}(T_{\text{ac}j}T_{\text{cd}i}-T_{\text{ad}%
j}T_{\text{cc}i})$

$+T_{\text{ad}j}T_{\text{bb}i}(T_{\text{ac}i}T_{\text{cd}j}-T_{\text{ad}%
i}T_{\text{cc}j})+T_{\text{ad}j}T_{\text{bb}i}(T_{\text{ac}j}T_{\text{cd}%
i}-T_{\text{ad}j}T_{\text{cc}i})](y_{0})\phi(y_{0})$

$+\ \frac{1}{288}\underset{i,j=q+1}{\overset{n}{\sum}}[T_{\text{ad}%
i}T_{\text{bc}j}(T_{\text{ac}i}T_{\text{bd}j}-T_{\text{ad}i}T_{\text{bc}%
j})+T_{\text{ad}i}T_{\text{bc}j}(T_{\text{ac}j}T_{\text{bd}i}-T_{\text{ad}%
j}T_{\text{bc}i})$

$+T_{\text{ad}j}T_{\text{bc}i}(T_{\text{ac}i}T_{\text{bd}j}-T_{\text{ad}%
i}T_{\text{bc}j})+T_{\text{ad}j}T_{\text{bc}i}(T_{\text{ac}j}T_{\text{bd}%
i}-T_{\text{ad}j}T_{\text{bc}i})](y_{0})\phi(y_{0})$

$-\ \frac{1}{144}[(R_{\text{cdcd}}^{P}-R_{\text{cdcd}}^{M})(R_{\text{abab}%
}^{P}-R_{\text{abab}}^{M})](y_{0})\phi(y_{0})$

$+\frac{1}{144}[(R_{\text{bdcd}}^{P}-R_{\text{bdcd}}^{M})(R_{\text{abac}}%
^{P}-R_{\text{abac}}^{M})](y_{0})\phi(y_{0})$

$\ +\ \frac{1}{144}[(R_{\text{bcdc}}^{P}-R_{\text{bcdc}}^{M})(R_{\text{abad}%
}^{P}-R_{\text{abad}}^{M})](y_{0})\phi(y_{0})$

$\ -\ \frac{1}{144}[(R_{\text{adcd}}^{P}-R_{\text{adcd}}^{M})(R_{\text{abbc}%
}^{P}-R_{\text{abbc}}^{M})](y_{0})\phi(y_{0})\qquad$

$\ +\ \frac{1}{144}[(R_{\text{acdc}}^{P}-R_{\text{acdc}}^{M})(R_{\text{abdb}%
}^{P}-R_{\text{abdb}}^{M})](y_{0})\phi(y_{0})$

$\ -\ \frac{1}{576}[(R_{\text{abcd}}^{P}-R_{\text{abcd}}^{M})]^{2}(y_{0}%
)\phi(y_{0})$

$-\frac{1}{144}<H,i>(y_{0})<H,j>(y_{0})\qquad\qquad\qquad\qquad\qquad
\qquad\qquad L_{3}$

$\times\lbrack2\varrho_{ij}+$ $\overset{q}{\underset{\text{a}=1}{4\sum}%
}R_{i\text{a}j\text{a}}-3\overset{q}{\underset{\text{a,b=1}}{\sum}%
}(T_{\text{aa}i}T_{\text{bb}j}-T_{\text{ab}i}T_{\text{ab}j}%
)-3\overset{q}{\underset{\text{a,b=1}}{\sum}}(T_{\text{aa}j}T_{\text{bb}%
i}-T_{\text{ab}j}T_{\text{ab}i}](y_{0})\phi(y_{0})$

$-\frac{1}{16}[<H,i>^{2}(y_{0})<H,j>^{2}](y_{0})\phi(y_{0})$

$-\frac{1}{144}<H,i>(y_{0})<H,j>(y_{0})$

$\times\lbrack2\varrho_{ij}+$ $\overset{q}{\underset{\text{a}=1}{4\sum}%
}R_{i\text{a}j\text{a}}-3\overset{q}{\underset{\text{a,b=1}}{\sum}%
}(T_{\text{aa}i}T_{\text{bb}j}-T_{\text{ab}i}T_{\text{ab}j}%
)-3\overset{q}{\underset{\text{a,b=1}}{\sum}}(T_{\text{aa}j}T_{\text{bb}%
i}-T_{\text{ab}j}T_{\text{ab}i}](y_{0})\phi(y_{0})$

$-\frac{1}{72}<H,i>(y_{0})<H,k>(y_{0})R_{jijk}(y_{0})\phi(y_{0})$

$-\frac{1}{16}<H,i>^{2}(y_{0})<H,j>^{2}(y_{0})\phi(y_{0})$

$-\frac{1}{72}<H,i>^{2}(y_{0})[\varrho_{jj}+$ $\overset{q}{\underset{\text{a}%
=1}{2\sum}}R_{j\text{a}j\text{a}}-3\overset{q}{\underset{\text{a,b=1}}{\sum}%
}(T_{\text{aa}j}T_{\text{bb}j}-T_{\text{ab}j}T_{\text{ab}j})](y_{0})\phi
(y_{0})$

$+\frac{5}{32}[<H,i>^{2}<H,j>^{2}](y_{0})\phi(y_{0})$

$+\frac{1}{48}<H,i>(y_{0})<H,j>$

$\times\lbrack2\varrho_{ij}+$ $\overset{q}{\underset{\text{a}=1}{4\sum}%
}R_{i\text{a}j\text{a}}-3\overset{q}{\underset{\text{a,b=1}}{\sum}%
}(T_{\text{aa}i}T_{\text{bb}j}-T_{\text{ab}i}T_{\text{ab}j}%
)-3\overset{q}{\underset{\text{a,b=1}}{\sum}}(T_{\text{aa}j}T_{\text{bb}%
i}-T_{\text{ab}j}T_{\text{ab}i}](y_{0})\phi(y_{0})$

$+\frac{1}{48}<H,i>^{2}(y_{0})[\tau^{M}\ -3\tau^{P}+\ \underset{\text{a}%
=1}{\overset{\text{q}}{\sum}}\varrho_{\text{aa}}^{M}+$
$\overset{q}{\underset{\text{a},\text{b}=1}{\sum}}R_{\text{abab}}^{M}$
$](y_{0})\phi(y_{0})$

$+\frac{1}{144}<H,i>(y_{0})[\nabla_{i}\varrho_{jj}-2\varrho_{ij}%
<H,j>+\overset{q}{\underset{\text{a}=1}{\sum}}(\nabla_{i}R_{\text{a}%
j\text{a}j}-4R_{i\text{a}j\text{a}}<H,j>)$

$+4\overset{q}{\underset{\text{a,b=1}}{\sum}}R_{i\text{a}j\text{b}%
}T_{\text{ab}j}+2\overset{q}{\underset{\text{a,b,c=1}}{\sum}}(T_{\text{aa}%
i}T_{\text{bb}j}T_{\text{cc}j}-3T_{\text{aa}i}T_{\text{bc}j}T_{\text{bc}%
j}+2T_{\text{ab}i}T_{\text{bc}j}T_{\text{ca}j})](y_{0})\phi(y_{0})$%
\qquad\qquad\qquad\qquad\qquad\ \ 

$+\frac{1}{144}<H,i>(y_{0})[\nabla_{j}\varrho_{ij}-2\varrho_{ij}%
<H,j>+\overset{q}{\underset{\text{a}=1}{\sum}}(\nabla_{j}R_{\text{a}%
i\text{a}j}-4R_{j\text{a}i\text{a}}<H,j>)$

$+4\overset{q}{\underset{\text{a,b=1}}{\sum}}R_{j\text{a}i\text{b}%
}T_{\text{ab}j}+2\overset{q}{\underset{\text{a,b,c=1}}{\sum}}(T_{\text{aa}%
j}T_{\text{bb}i}T_{\text{cc}j}-3T_{\text{aa}j}T_{\text{bc}i}T_{\text{bc}%
j}+2T_{\text{ab}j}T_{\text{bc}i}T_{\text{ac}j})](y_{0})\phi(y_{0})$

$+\frac{1}{144}<H,i>(y_{0})[\nabla_{j}\varrho_{ij}-2\varrho_{jj}%
<H,i>+\overset{q}{\underset{\text{a}=1}{\sum}}(\nabla_{j}R_{\text{a}%
i\text{a}j}-4R_{j\text{a}j\text{a}}<H,i>)$

$+4\overset{q}{\underset{\text{a,b=1}}{\sum}}R_{j\text{a}j\text{b}%
}T_{\text{ab}i}+2\overset{q}{\underset{\text{a,b,c=1}}{\sum}}(T_{\text{aa}%
j}T_{\text{bb}j}T_{\text{cc}i}-3T_{\text{aa}j}T_{\text{bc}j}T_{\text{bc}%
i}+2T_{\text{ab}j}T_{\text{bc}j}T_{\text{ac}i})](y_{0})\phi(y_{0})$

$-\frac{1}{192}<H,i>^{2}<H,j>^{2}(y_{0})\phi(y_{0})\qquad\qquad\qquad
\qquad\qquad\qquad\qquad$I$_{3213}\ $

$-\frac{1}{288}<H,i>^{2}(y_{0})[\tau^{M}\ -3\tau^{P}+\ \underset{\text{a}%
=1}{\overset{\text{q}}{\sum}}\varrho_{\text{aa}}^{M}+$
$\overset{q}{\underset{\text{a},\text{b}=1}{\sum}}R_{\text{abab}}^{M}$
$](y_{0})\phi(y_{0})$

$-\frac{1}{288}<H,i>(y_{0})<H,j>(y_{0})$

$\times\lbrack2\varrho_{ij}+$ $\overset{q}{\underset{\text{a}=1}{4\sum}%
}R_{i\text{a}j\text{a}}-3\overset{q}{\underset{\text{a,b=1}}{\sum}%
}(T_{\text{aa}i}T_{\text{bb}j}-T_{\text{ab}i}T_{\text{ab}j}%
)-3\overset{q}{\underset{\text{a,b=1}}{\sum}}(T_{\text{aa}j}T_{\text{bb}%
i}-T_{\text{ab}j}T_{\text{ab}i}](y_{0})\phi(y_{0})$

$+\frac{1}{144}<H,i>(y_{0})<H,k>(y_{0})R_{jijk}(y_{0})$

$+\frac{1}{144}<H,i>^{2}(y_{0})[\varrho_{jj}+$ $\overset{q}{\underset{\text{a}%
=1}{2\sum}}R_{j\text{a}j\text{a}}-3\overset{q}{\underset{\text{a,b=1}}{\sum}%
}(T_{\text{aa}j}T_{\text{bb}j}-T_{\text{ab}j}T_{\text{ab}j})](y_{0})\phi
(y_{0})$

$-\frac{1}{288}<H,i>(y_{0})[\nabla_{i}\varrho_{jj}-2\varrho_{ij}%
<H,j>+\overset{q}{\underset{\text{a}=1}{\sum}}(\nabla_{i}R_{\text{a}%
j\text{a}j}-4R_{i\text{a}j\text{a}}<H,j>)$

$+4\overset{q}{\underset{\text{a,b=1}}{\sum}}R_{i\text{a}j\text{b}%
}T_{\text{ab}j}+2\overset{q}{\underset{\text{a,b,c=1}}{\sum}}(T_{\text{aa}%
i}T_{\text{bb}j}T_{\text{cc}j}-3T_{\text{aa}i}T_{\text{bc}j}T_{\text{bc}%
j}+2T_{\text{ab}i}T_{\text{bc}j}T_{\text{ca}j})](y_{0})\phi(y_{0})$%
\qquad\qquad\qquad\qquad\qquad\ \ 

$-\frac{1}{288}<H,i>(y_{0})[\nabla_{j}\varrho_{ij}-2\varrho_{ij}%
<H,j>+\overset{q}{\underset{\text{a}=1}{\sum}}(\nabla_{j}R_{\text{a}%
i\text{a}j}-4R_{j\text{a}i\text{a}}<H,j>)$

$+4\overset{q}{\underset{\text{a,b=1}}{\sum}}R_{j\text{a}i\text{b}%
}T_{\text{ab}j}+2\overset{q}{\underset{\text{a,b,c=1}}{\sum}}(T_{\text{aa}%
j}T_{\text{bb}i}T_{\text{cc}j}-3T_{\text{aa}j}T_{\text{bc}i}T_{\text{bc}%
j}+2T_{\text{ab}j}T_{\text{bc}i}T_{\text{ac}j})](y_{0})\phi(y_{0})$

$-\frac{1}{288}<H,i>(y_{0})[\nabla_{j}\varrho_{ij}-2\varrho_{jj}%
<H,i>+\overset{q}{\underset{\text{a}=1}{\sum}}(\nabla_{j}R_{\text{a}%
i\text{a}j}-4R_{j\text{a}j\text{a}}<H,i>)+4\overset{q}{\underset{\text{a,b=1}%
}{\sum}}R_{j\text{a}j\text{b}}T_{\text{ab}i}$

$+2\overset{q}{\underset{\text{a,b,c=1}}{\sum}}(T_{\text{aa}j}T_{\text{bb}%
j}T_{\text{cc}i}-3T_{\text{aa}j}T_{\text{bc}j}T_{\text{bc}i}+2T_{\text{ab}%
j}T_{\text{bc}j}T_{\text{ac}i})](y_{0})\phi(y_{0})$

$+\frac{1}{24}[\left\Vert \text{X}\right\Vert _{M}^{2}+\operatorname{div}%
$X$_{M}-\left\Vert \text{X}\right\Vert _{P}^{2}-\operatorname{div}X_{P}%
](y_{0})[\left\Vert \text{X}\right\Vert _{M}^{2}-\operatorname{div}$%
X$_{M}-\left\Vert \text{X}\right\Vert _{P}^{2}+\operatorname{div}$%
X$_{P}](y_{0})\phi(y_{0})\qquad$\ I$_{3212}$

$+\frac{1}{6}X_{i}(y_{0})T_{\text{ab}i}(y_{0})T_{\text{ab}j}(y_{0})X_{j}%
(y_{0})+\frac{1}{3}\perp_{\text{a}ij}(y_{0})X_{i}(y_{0})[\frac{\partial X_{j}%
}{\partial x_{\text{a}}}-\perp_{\text{a}jk}X_{k}](y_{0})\phi(y_{0})\qquad
$I$_{32122}\qquad Q_{1}$

$+$ $\frac{2}{3}X_{i}(y_{0})X_{j}(y_{0})\frac{\partial X_{j}}{\partial
x_{\text{a}}}(y_{0})-\frac{1}{6}X_{i}(y_{0})\frac{\partial^{2}X_{j}}{\partial
x_{\text{a}}\partial x_{j}}(y_{0})\phi(y_{0})\qquad\qquad Q_{2}$

$-\frac{1}{12}X_{i}(y_{0})\frac{\partial^{2}X_{i}}{\partial x_{\text{a}}^{2}%
}(y_{0})+\frac{1}{12}X_{i}^{2}(y_{0})[\operatorname{div}X_{M}-\left\Vert
X\right\Vert _{M}^{2}+\left\Vert X\right\Vert _{P}^{2}-\operatorname{div}%
X_{P}-$ $<H,j>X_{j}](y_{0})\phi(y_{0})$

$+\frac{1}{6}X_{i}(y_{0})X_{j}(y_{0})\frac{\partial X_{i}}{\partial x_{j}%
}(y_{0})\phi(y_{0})+$ $\frac{1}{18}X_{i}(y_{0})X_{k}(y_{0})R_{jijk}(y_{0}%
)\phi(y_{0})-\frac{1}{12}X_{i}(y_{0})\frac{\partial^{2}X_{i}}{\partial
x_{j}^{2}}(y_{0})\phi(y_{0})$

$+\frac{1}{12}[R_{\text{a}i\text{a}k}-\underset{\text{c=1}}{\overset{\text{q}%
}{\sum}}T_{\text{ac}i}T_{\text{ac}k}-\perp_{\text{a}ik}\perp_{\text{a}%
jk}](y_{0})X_{k}(y_{0})\phi(y_{0})+\frac{1}{18}R_{ijkj}(y_{0})X_{i}%
(y_{0})X_{k}(y_{0})\phi(y_{0})$

$+\frac{1}{12}<H,j>(y_{0})X_{i}(y_{0})[X_{i}X_{j}-\frac{1}{2}\left(
\frac{\partial X_{j}}{\partial x_{i}}+\frac{\partial X_{i}}{\partial x_{j}%
}\right)  ](y_{0})\phi(y_{0})\qquad\qquad\qquad\qquad\qquad\qquad\qquad
\qquad\qquad\qquad\ \qquad\qquad\qquad\qquad\qquad\qquad\qquad\qquad
\qquad\qquad\qquad$

$-\frac{1}{6}[-R_{\text{a}i\text{b}i}+5\overset{q}{\underset{\text{c}=1}{\sum
}}T_{\text{ac}i}T_{\text{bc}i}+2\overset{n}{\underset{j=q+1}{\sum}}%
\perp_{\text{a}ij}\perp_{\text{b}ij}](y_{0})\underset{k=q+1}{\overset{n}{\sum
}}T_{\text{ab}k}(y_{0})X_{k}(y_{0})\phi(y_{0})\qquad\qquad$I$_{32123}\qquad
S_{1}\qquad$

$-\frac{2}{9}\underset{j=q+1}{\overset{n}{\sum}}R_{i\text{a}ij}(y_{0}%
)[\frac{\partial X_{j}}{\partial x_{\text{a}}}%
-\underset{k=q+1}{\overset{n}{\sum}}\perp_{\text{a}jk}X_{k}](y_{0})\phi
(y_{0})$

$+\frac{1}{12}\times\frac{2}{3}\underset{j,k=q+1}{\overset{n}{\sum}}%
R_{ijik}(y_{0})[X_{j}X_{k}-\frac{1}{2}(\frac{\partial X_{j}}{\partial x_{k}%
}+\frac{\partial X_{k}}{\partial x_{j}})](y_{0})\phi(y_{0})$

$-\frac{1}{6}T_{\text{ab}i}(y_{0})\frac{\partial^{2}X_{i}}{\partial
x_{\text{a}}\partial x_{\text{b}}}(y_{0})\phi(y_{0})\qquad\qquad\qquad
S_{2}\qquad\qquad S_{21}\qquad\qquad\qquad\qquad\qquad\qquad$

$+$ $\frac{1}{12}T_{\text{ab}i}(y_{0})$

$\times\lbrack$ $(R_{\text{a}i\text{b}j}+R_{\text{a}j\text{b}i})$
$-\underset{\text{c=1}}{\overset{\text{q}}{\sum}}(T_{\text{ac}i}T_{\text{bc}%
j}+T_{\text{ac}j}T_{\text{bc}i})-\overset{n}{\underset{k=q+1}{\sum}}%
(\perp_{\text{a}ik}\perp_{\text{b}jk}+$ $\perp_{\text{a}jk}\perp_{\text{b}%
ik})](y_{0})X_{j}(y_{0})\phi(y_{0})$

$-$ $\frac{1}{6}T_{\text{ab}i}(y_{0})T_{\text{ab}j}(y_{0})[X_{i}X_{j}-\frac
{1}{2}\left(  \frac{\partial X_{i}}{\partial x_{j}}+\frac{\partial X_{j}%
}{\partial x_{i}}\right)  ](y_{0})\phi(y_{0})$

$\qquad-\frac{1}{3}\perp_{\text{a}ij}(y_{0})[(X_{i}\frac{\partial X_{j}%
}{\partial x_{\text{a}}}+X_{j}\frac{\partial X_{i}}{\partial x_{\text{a}}%
})-\frac{1}{4}\left(  \frac{\partial^{2}X_{i}}{\partial x_{\text{a}}\partial
x_{j}}+\frac{\partial^{2}X_{j}}{\partial x_{\text{a}}\partial x_{i}}\right)
](y_{0})\qquad S_{22}$

$\qquad-\frac{1}{6}\perp_{\text{a}ij}(y_{0})[T_{\text{ab}j}\frac{\partial
X_{i}}{\partial x_{\text{b}}}](y_{0})$

$\qquad+\frac{1}{6}\perp_{\text{a}ij}(y_{0})[(\perp_{\text{b}ik}T_{\text{ab}%
j})+\frac{2}{3}(2R_{\text{a}ijk}+R_{\text{a}jik}+R_{\text{a}kji})](y_{0}%
)X_{k}(y_{0})$

$\qquad-\frac{1}{6}\perp_{\text{a}ij}(y_{0})\perp_{\text{a}jk}(y_{0}%
)[X_{i}X_{k}-\frac{1}{2}\left(  \frac{\partial X_{i}}{\partial x_{k}}%
+\frac{\partial X_{k}}{\partial x_{i}}\right)  ](y_{0})\qquad\qquad
\qquad\qquad\qquad\qquad\qquad\qquad$

$\qquad+\frac{1}{12}[(\frac{\partial X_{j}}{\partial x_{\text{a}}})^{2}%
+X_{j}\frac{\partial^{2}X_{j}}{\partial x_{\text{a}}^{2}}-\frac{1}{2}%
\frac{\partial^{3}X_{j}}{\partial x_{\text{a}}^{2}\partial x_{j}}](y_{0}%
)\phi(y_{0})-\frac{1}{6}\overset{n}{\underset{k=q+1}{\sum}}[\perp_{\text{b}%
ik}$T$_{\text{aa}k}\frac{\partial X_{i}}{\partial x_{\text{b}}^{2}}%
](y_{0})\phi(y_{0})\qquad\qquad S_{3}\qquad S_{31}$

$\qquad+\frac{1}{144}[\{4\nabla_{i}R_{i\text{a}j\text{a}}+2\nabla
_{j}R_{i\text{a}i\text{a}}+$ $8(\overset{q}{\underset{\text{c=1}}{%
{\textstyle\sum}
}}R_{\text{a}i\text{c}i}^{{}}T_{\text{ac}j}+\;\overset{n}{\underset{k=q+1}{%
{\textstyle\sum}
}}R_{\text{a}iik}\perp_{\text{a}jk})$

$\qquad+8(\overset{q}{\underset{\text{c=1}}{%
{\textstyle\sum}
}}R_{\text{a}i\text{c}j}^{{}}T_{\text{ac}i}+\;\overset{n}{\underset{k=q+1}{%
{\textstyle\sum}
}}R_{\text{a}ijk}\perp_{\text{a}ik})+8(\overset{q}{\underset{\text{c=1}}{%
{\textstyle\sum}
}}R_{\text{a}j\text{c}i}^{{}}T_{\text{ac}i}+\;\overset{n}{\underset{k=q+1}{%
{\textstyle\sum}
}}R_{\text{a}jik}\perp_{\text{a}ik})\}$\ 

$\qquad+\frac{2}{3}\underset{k=q+1}{\overset{n}{\sum}}\{T_{\text{aa}%
k}(R_{ijik}+3\overset{q}{\underset{\text{c}=1}{\sum}}\perp_{\text{c}ij}%
\perp_{\text{c}ik})\}](y_{0})X_{k}(y_{0})\phi(y_{0})$

$-\frac{1}{12}[$ R$_{\text{a}i\text{a}k}$ $-\underset{\text{c=1}%
}{\overset{\text{q}}{\sum}}T_{\text{ac}i}T_{\text{ac}k}%
-\overset{n}{\underset{l=q+1}{\sum}}(\perp_{\text{a}il}\perp_{\text{a}%
kl}](y_{0})\times\lbrack X_{i}X_{k}-\frac{1}{2}\left(  \frac{\partial X_{i}%
}{\partial x_{k}}+\frac{\partial X_{k}}{\partial x_{i}}\right)  ](y_{0}%
)\phi(y_{0})$

$-\frac{1}{24}T_{\text{aa}k}(y_{0})[-X_{i}^{2}X_{k}+X_{k}\frac{\partial X_{i}%
}{\partial x_{i}}\ +X_{i}\left(  \frac{\partial X_{k}}{\partial x_{i}}%
+\frac{\partial X_{i}}{\partial x_{k}}\right)  -\frac{1}{3}\left(
\frac{\partial^{2}X_{k}}{\partial x_{i}^{2}}+2\frac{\partial^{2}X_{i}%
}{\partial x_{i}\partial x_{k}}\right)  ](y_{0})\phi(y_{0})$

$+\frac{1}{18}[R_{\text{a}jij}\frac{\partial X_{i}}{\partial x_{\text{a}}^{2}%
}](y_{0})\phi(y_{0})\qquad\qquad\qquad\qquad\qquad\qquad\qquad\qquad\qquad
S_{32}$

$+\frac{1}{24}[\frac{4}{3}\overset{q}{\underset{\text{a}=1}{\sum}}%
\perp_{\text{a}ki}R_{ij\text{a}j}-\frac{1}{3}(\nabla_{i}R_{kjij}+\nabla
_{j}R_{ijik}+\nabla_{k}R_{ijij})](y_{0})X_{k}(y_{0})\phi(y_{0})$

\ $-\frac{1}{18}R_{ijkj}(y_{0})[X_{i}X_{k}-\frac{1}{2}\left(  \frac{\partial
X_{i}}{\partial x_{k}}+\frac{\partial X_{k}}{\partial x_{i}}\right)
](y_{0})\phi(y_{0})$

$+\frac{1}{24}[X_{i}^{2}X_{j}^{2}-2X_{i}X_{j}\left(  \frac{\partial X_{j}%
}{\partial x_{i}}+\frac{\partial X_{i}}{\partial x_{j}}\right)  -X_{i}%
^{2}\frac{\partial X_{j}}{\partial x_{j}}-X_{j}^{2}\frac{\partial X_{i}%
}{\partial x_{i}}](y_{0})\phi(y_{0})$

$+\frac{1}{48}\left(  \frac{\partial X_{j}}{\partial x_{i}}+\frac{\partial
X_{i}}{\partial x_{j}}\right)  ^{2}(y_{0})\phi(y_{0})+\frac{1}{24}\left(
\frac{\partial X_{i}}{\partial x_{i}}\frac{\partial X_{j}}{\partial x_{j}%
}\right)  (y_{0})\phi(y_{0})\qquad$

$+\frac{1}{36}X_{i}(y_{0})\left(  2\frac{\partial^{2}X_{j}}{\partial
x_{i}\partial x_{j}}+\frac{\partial^{2}X_{i}}{\partial x_{j}^{2}}\right)
(y_{0})\phi(y_{0})+\frac{1}{36}X_{j}(y_{0})\left(  \frac{\partial^{2}X_{j}%
}{\partial x_{i}^{2}}+2\frac{\partial^{2}X_{i}}{\partial x_{i}\partial x_{j}%
}\right)  (y_{0})\phi(y_{0})$

$-\frac{1}{48}\left(  \frac{\partial^{3}X_{i}}{\partial x_{i}\partial
x_{j}^{2}}+\frac{\partial^{3}X_{j}}{\partial x_{i}^{2}\partial x_{j}}\right)
(y_{0})\phi(y_{0})$

$+$ $\frac{2}{3}<H,j>(y_{0})\left(  \frac{\partial^{2}X_{i}}{\partial
x_{i}\partial x_{j}}+2\frac{\partial^{2}X_{j}}{\partial x_{i}^{2}}\right)
(y_{0})\phi(y_{0})+$ $\frac{2}{3}<H,j>(y_{0})R_{ijik}(y_{0})X_{k}(y_{0}%
)\phi(y_{0})\qquad$I$_{3213}$

$+\frac{1}{12}[<H,i><H,j>\ +\frac{1}{6}(2\varrho_{ij}%
+4\overset{q}{\underset{\text{a}=1}{\sum}}R_{i\text{a}j\text{a}}%
-6\overset{q}{\underset{\text{a,b}=1}{\sum}}T_{\text{aa}i}T_{\text{bb}%
j}-T_{\text{ab}i}T_{\text{ab}j})](y_{0})\phi(y_{0})$

$\times\frac{1}{2}[\left(  \frac{\partial X_{j}}{\partial x_{i}}%
-\frac{\partial X_{i}}{\partial x_{j}}\right)  ](y_{0})\phi(y_{0})$

$-\frac{1}{12}\perp_{\text{a}ij}(y_{0})<H,i>(y_{0})[(X_{j}\perp_{\text{a}%
ij}-\frac{\partial X_{i}}{\partial x_{\text{a}}})+\frac{\partial X_{\text{a}}%
}{\partial x_{i}}](y_{0})\phi(y_{0})$

$-\frac{1}{18}[X_{j}\left(  2\frac{\partial^{2}X_{j}}{\partial x_{i}^{2}%
}+\frac{\partial^{2}X_{i}}{\partial x_{i}\partial x_{j}}\right)  ](y_{0}%
)\phi(y_{0})-\frac{1}{12}[\left(  \frac{\partial X_{i}}{\partial x_{j}}%
+\frac{\partial X_{j}}{\partial x_{i}}\right)  ]\frac{\partial X_{j}}{\partial
x_{i}}(y_{0})\phi(y_{0})\qquad\ $I$_{3214}\qquad\qquad\qquad$\qquad
\qquad\qquad\qquad\qquad\qquad\qquad\qquad\qquad\qquad\qquad\qquad\qquad
\qquad\qquad\qquad

$+\frac{1}{12}\frac{\partial^{2}\text{V}}{\partial x_{i}^{2}}(y_{0})\phi
(y_{0})\qquad\qquad$I$_{3215}$

$+\frac{1}{12}\underset{i=q+\text{1}}{\overset{\text{n}}{\sum}}%
\underset{\text{a,b=1}}{\overset{\text{q}}{\sum}}[-R_{\text{a}i\text{b}%
i}+5\overset{q}{\underset{\text{c}=1}{\sum}}T_{\text{ac}i}T_{\text{bc}%
i}+2\overset{n}{\underset{j=q+1}{\sum}}\perp_{\text{a}ij}\perp_{\text{b}%
ij}](y_{0})\times\frac{\partial^{2}\phi}{\partial\text{x}_{\text{a}}%
\partial\text{x}_{\text{b}}}(y_{0})\qquad$\textbf{I}$_{322}$

$+\frac{1}{72}\underset{i,j,k=q+1}{\overset{n}{\sum}}R_{ijik}(y_{0}%
)\Omega_{jk}(y_{0})\phi(y_{0})\qquad\qquad\qquad\qquad$\ I$_{323}$

$+\frac{1}{24}\underset{\text{a=1}}{\overset{\text{q}}{\sum}}%
\underset{i,j=q+1}{\overset{n}{\sum}}\left\{  \frac{8}{3}R_{i\text{a}%
ij}+4\underset{\text{b=1}}{\overset{\text{q}}{\sum}}T_{\text{ab}i}%
\perp_{\text{b}ji}\right\}  (y_{0})\left\{  -\Omega_{\text{a}j}+[\Lambda
_{\text{a}},\Lambda_{j}]\right\}  (y_{0})\phi(y_{0})$

$+\frac{1}{12}\underset{i\text{=}q+1}{\overset{\text{n}}{\sum}}%
\underset{\text{a,b=1}}{\overset{\text{q}}{\sum}}$ $[-R_{\text{a}i\text{b}%
i}+5\overset{q}{\underset{\text{c}=1}{\sum}}T_{\text{ac}i}T_{\text{bc}%
i}+2\overset{n}{\underset{\text{k}=q+1}{\sum}}\perp_{\text{a}i\text{k}}%
\perp_{\text{b}i\text{k}}](y_{0})\times\lbrack\Lambda_{\text{a}}(y_{0}%
)\Lambda_{\text{b}}(y_{0})\phi(y_{0})]\qquad$\ I$_{324}$

$+\frac{1}{12}[\frac{8}{3}R_{i\text{a}ij}-4\underset{\text{b=1}%
}{\overset{\text{q}}{\sum}}T_{\text{ab}i}(y_{0})\perp_{\text{b}ij}%
](y_{0})[\Lambda_{\text{a}}\Lambda_{j}\phi](y_{0})$

$+\frac{1}{12}\underset{i\text{=}q+1}{\overset{\text{n}}{\sum}}%
\underset{\text{a,b=1}}{\overset{\text{q}}{\sum}}$ $[-R_{\text{a}i\text{b}%
i}+5\overset{q}{\underset{\text{c}=1}{\sum}}T_{\text{ac}i}T_{\text{bc}%
i}+2\overset{n}{\underset{\text{k}=q+1}{\sum}}\perp_{\text{a}i\text{k}}%
\perp_{\text{b}i\text{k}}](y_{0})\qquad$I$_{325}\qquad$I$_{3251}$

$\times\lbrack\Lambda_{\text{a}}(y_{0})\Lambda_{\text{b}}(y_{0})\phi(y_{0})]$

$+\frac{1}{12}[\frac{8}{3}R_{i\text{a}ij}-4\underset{\text{b=1}%
}{\overset{\text{q}}{\sum}}T_{\text{ab}i}(y_{0})\perp_{\text{b}ij}%
](y_{0})[\Lambda_{\text{a}}\Lambda_{j}\phi](y_{0})$

$+\frac{1}{24}\underset{i=q+1}{\overset{n}{\sum}}\underset{\text{a=1}%
}{\overset{\text{q}}{\sum}}\left(  \frac{\partial\Omega_{i\text{a}}}{\partial
x_{i}}\Lambda_{\text{a}}+[\Omega_{i\text{a}}+[\Lambda_{\text{a}},\Lambda
_{i}],\Lambda_{i}]\right)  \Lambda_{\text{a}}(y_{0})\phi(y_{0})\qquad
$I$_{3252}$

$+\frac{1}{24}\underset{i=q+1}{\overset{n}{\sum}}\underset{\text{a=1}%
}{\overset{\text{q}}{\sum}}\Lambda_{\text{a}}(y_{0})\left(  \frac
{\partial\Omega_{i\text{a}}}{\partial x_{i}}\Lambda_{\text{a}}+[\Omega
_{i\text{a}}+[\Lambda_{\text{a}},\Lambda_{i}],\Lambda_{i}]\right)  (y_{0}%
)\phi(y_{0})$

$+\frac{1}{12}\underset{i=q+1}{\overset{n}{\sum}}\underset{\text{a=1}%
}{\overset{\text{q}}{\sum}}\left(  \Omega_{i\text{a}}+[\Lambda_{\text{a}%
},\Lambda_{i}]\right)  ^{2}(y_{0})\phi(y_{0})$

$+\frac{1}{48}\underset{i,j=q+1}{\overset{n}{\sum}}\left(  \Omega_{ij}%
\Omega_{ij}\right)  (y_{0})\phi(y_{0})+\frac{1}{72}%
\underset{i,j=q+1}{\overset{n}{\sum}}\left(  \frac{\partial\Omega_{ij}%
}{\partial\text{x}_{i}}\Lambda_{j}+\Lambda_{j}\frac{\partial\Omega_{ij}%
}{\partial\text{x}_{i}}\right)  (y_{0})\phi(y_{0})$

$+\frac{1}{12}[\underset{i=q+1}{\overset{n}{\sum}}\underset{\text{a,b=1}%
}{\overset{\text{q}}{\sum}}2T_{\text{ab}i}(y_{0})\left\{  (\Omega_{i\text{a}%
}+[\Lambda_{\text{a}},\Lambda_{i}])\Lambda_{\text{b}}+\Lambda_{\text{a}%
}(\Omega_{i\text{a}}+[\Lambda_{\text{a}},\Lambda_{i}])\right\}  (y_{0}%
)\phi(y_{0})\qquad$I$_{3253}$

$-\frac{1}{12}[\underset{i,j=q+1}{\overset{n}{\sum}}\underset{\text{a=1}%
}{\overset{\text{q}}{\sum}}\perp_{\text{a}ij}(y_{0})\left\{  (\Omega
_{i\text{a}}+[\Lambda_{\text{a}},\Lambda_{i}])\Lambda_{j}+\frac{1}{2}%
\Lambda_{\text{a}}\Omega_{ij}\right\}  ](y_{0})\phi(y_{0})$

$-\frac{1}{12}[\underset{i,j=q+1}{\overset{n}{\sum}}\underset{\text{b=1}%
}{\overset{\text{q}}{\sum}}\perp_{\text{b}ij}(y_{0})\left\{  \frac{1}{2}%
\Omega_{ij}\Lambda_{\text{b}}+\Lambda_{j}(\Omega_{i\text{b}}+[\Lambda
_{\text{b}},\Lambda_{i}])\right\}  ](y_{0})\phi(y_{0})$

$\mathbf{-}\frac{1}{12}\underset{\text{a,b=1}}{\overset{\text{q}}{\sum}%
}\underset{i,j=q+1}{\overset{n}{\sum}}T_{\text{ab}i}(y_{0})[-R_{\text{a}%
i\text{b}i}+5\overset{q}{\underset{\text{c}=1}{\sum}}T_{\text{ac}%
i}T_{\text{bc}i}+2\overset{n}{\underset{k=q+1}{\sum}}\perp_{\text{a}ik}%
\perp_{\text{b}ik}](y_{0})\Lambda_{j}(y_{0})\phi(y_{0})\ \ \qquad$%
\textbf{I}$_{326}\qquad$\textbf{I}$_{3261}\qquad$

$+\frac{1}{12}\underset{i\text{=q+1}}{\overset{\text{n}}{\sum}}%
\underset{j\text{=q+1}}{\overset{\text{n}}{\sum}}\underset{\text{a=1}%
}{\overset{\text{q}}{\sum}}[4\underset{\text{c=1}}{\overset{q}{\sum}%
}(T_{\text{ac}i})(\perp_{j\text{c}i})+\frac{8}{3}R_{i\text{a}ij}](y_{0}%
)\phi(y_{0})\qquad\qquad$I$_{32613}$

$\times\lbrack\underset{\text{c=1}}{\overset{\text{q}}{\sum}}T_{\text{ac}%
j}\frac{\partial\phi}{\partial\text{x}_{\text{c}}}+\underset{\text{b=1}%
}{\overset{\text{q}}{\sum}}T_{\text{ab}j}\Lambda_{\text{b}}%
-\underset{k=q+1}{\overset{n}{\sum}}\perp_{\text{a}jk}\Lambda_{k}](y_{0}%
)\phi(y_{0})$

$-\frac{1}{24}\underset{i=q+1}{\overset{n}{\sum}}\underset{\text{a,b=1}%
}{\overset{\text{q}}{\sum}}$ $\overset{n}{\underset{k=q+1}{\sum}}%
$T$_{\text{aa}k}[\frac{8}{3}R_{i\text{c}ik}+4\underset{\text{d}%
=1}{\overset{q}{\sum}}(T_{\text{db}k})(\perp_{\text{d}ik})]\frac{\partial\phi
}{\partial\text{x}_{\text{b}}}$ $(y_{0})\qquad$ I$_{32621}\qquad$I$_{3262}$

$-\frac{1}{12}\underset{i=q+1}{\overset{n}{\sum}}\underset{\text{a,b=1}%
}{\overset{\text{q}}{\sum}}[$ $\overset{n}{\underset{k,l=q+1}{\sum}}%
\perp_{\text{b}ik}(-R_{\text{a}k\text{a}l}+\underset{\text{d=1}%
}{\overset{\text{q}}{\sum}}T_{\text{ad}k}T_{\text{ad}l}%
))-\overset{n}{\underset{k,l=q+1}{\sum}}\perp_{\text{b}ik}%
(\underset{r=q+1}{\overset{n}{\sum}}\perp_{\text{a}kr}\perp_{\text{a}%
lr})](y_{0})\frac{\partial\phi}{\partial\text{x}_{\text{b}}}$ $(y_{0})$

$-\frac{1}{12}\underset{i=q+1}{\overset{n}{\sum}}\underset{\text{b=1}%
}{\overset{\text{q}}{\sum}}[\frac{8}{3}\overset{q}{\underset{\text{c}=1}{\sum
}}($T$_{\text{bc}i}R_{ij\text{c}j})+\frac{2}{3}%
\overset{n}{\underset{k=q+1}{\sum}}(\perp_{\text{b}ik}R_{ijjk})](y_{0}%
)\frac{\partial\phi}{\partial\text{x}_{\text{b}}}(y_{0})$

$-\frac{1}{6}\underset{i=q+1}{\overset{n}{\sum}}\underset{\text{a,b=1}%
}{\overset{\text{q}}{\sum}}[4\overset{q}{\underset{\text{c=1}}{%
{\textstyle\sum}
}}$R$_{ij\text{c}i}^{{}}T_{\text{bc}j}+$ $4\overset{n}{\underset{k=q+1}{%
{\textstyle\sum}
}}$R$_{ijik}\perp_{\text{b}jk}+3\nabla_{i}$R$_{j\text{b}ij}%
+4\overset{q}{\underset{\text{c=1}}{%
{\textstyle\sum}
}}$R$_{ij\text{c}j}^{{}}T_{\text{bc}i}+$ $4$R$_{ijjk}\perp_{\text{b}ik}%
](y_{0})\frac{\partial\phi}{\partial\text{x}_{\text{b}}}(y_{0})$

$-\frac{1}{24}\overset{n}{\underset{k=q+1}{\sum}}$T$_{\text{aa}k}[\frac{8}%
{3}R_{i\text{c}ik}+4\underset{\text{d}=1}{\overset{q}{\sum}}(T_{\text{db}%
k})(\perp_{\text{d}ik})]\Lambda_{\text{b}}(y_{0})\phi(y_{0})\qquad$%
I$_{326221}\qquad$I$_{32622}$

$-\frac{1}{12}$ $\overset{n}{\underset{k,l=q+1}{\sum}}\perp_{\text{b}ik}%
(y_{0})[-R_{\text{a}k\text{a}l}+\underset{\text{d=1}}{\overset{\text{q}}{\sum
}}T_{\text{ad}k}T_{\text{ad}l}](y_{0})\Lambda_{\text{b}}(y_{0})\phi(y_{0})$

$-\frac{1}{12}\overset{n}{\underset{k,l=q+1}{\sum}}\perp_{\text{b}ik}%
(y_{0})[\underset{r=q+1}{\overset{n}{\sum}}\perp_{\text{a}kr}\perp
_{\text{a}lr}](y_{0})\Lambda_{\text{b}}(y_{0})\phi(y_{0})$

$-\frac{1}{144}[\{4\nabla_{i}$R$_{i\text{a}j\text{a}}$ $+2\nabla_{j}%
$R$_{i\text{a}i\text{a}}+$ $8$ $(\overset{q}{\underset{\text{c=1}}{%
{\textstyle\sum}
}}R_{\text{a}i\text{c}i}^{{}}T_{\text{ac}j}+\;\overset{n}{\underset{k=q+1}{%
{\textstyle\sum}
}}R_{\text{a}iik}\perp_{\text{a}jk})$

\ $+8(\overset{q}{\underset{\text{c=1}}{%
{\textstyle\sum}
}}R_{\text{a}i\text{c}j}^{{}}T_{\text{ac}i}+\overset{n}{\underset{l=q+1}{%
{\textstyle\sum}
}}R_{\text{a}ijl}\perp_{\text{a}il})+8(\overset{q}{\underset{\text{c=1}}{%
{\textstyle\sum}
}}R_{\text{a}j\text{c}i}^{{}}T_{\text{ac}i}+\overset{q}{\underset{\text{c=1}}{%
{\textstyle\sum}
}}R_{\text{a}j\text{c}i}^{{}}T_{\text{ac}i})\}$\ \qquad$\qquad\ \ \ \ \ $

$+\frac{2}{3}\underset{k=q+1}{\overset{n}{\sum}}\{T_{\text{aa}k}%
(R_{ijik}+3\overset{q}{\underset{\text{c}=1}{\sum}}\perp_{\text{c}ij}%
\perp_{\text{c}ik})\}](y_{0})\Lambda_{k}(y_{0})\phi(y_{0})$

\ $+\frac{1}{24}[\frac{4}{3}\overset{q}{\underset{\text{a}=1}{\sum}}%
\perp_{\text{a}ik}R_{ij\text{a}j}+\frac{1}{3}(\nabla_{i}R_{kjij}+\nabla
_{j}R_{ijik}+\nabla_{k}R_{ijij})](y_{0})\Lambda_{k}(y_{0})\phi(y_{0})$

\ $\mathbf{-}$ $\frac{1}{72}\underset{i,j=q+1}{\overset{n}{\sum}%
}\underset{\text{a}=1}{\overset{q}{\sum}}T_{\text{aa}j}(y_{0})\frac
{\partial\Omega_{ij}}{\partial\text{x}_{i}}(y_{0})\phi(y_{0})\qquad
$I$_{326222}$

$+$ $\frac{1}{12}$ $\overset{n}{\underset{j=q+1}{\sum}}(\perp_{\text{b}ij}%
$T$_{\text{aa}j})(y_{0})\left(  \Omega_{i\text{b}}(y_{0})+[\Lambda_{\text{b}%
},\Lambda_{i}]\right)  (y_{0})\phi(y_{0})$ \qquad I$_{326223}$

$\mathbf{-}\frac{1}{18}$ $\underset{i,j=q+1}{\overset{n}{\sum}}%
\underset{\text{b=1}}{\overset{q}{\sum}}R_{\text{b}jij}(y_{0})\left(
\Omega_{i\text{a}}(y_{0})+[\Lambda_{\text{a}},\Lambda_{i}]\right)  (y_{0}%
)\phi(y_{0})$

$\mathbf{-}$ $\frac{1}{24}\underset{i,j=q+1}{\overset{n}{\sum}}%
\underset{\text{a=1}}{\overset{q}{\sum}}[$ R$_{\text{a}i\text{a}j}$
$-\underset{\text{c=1}}{\overset{\text{q}}{\sum}}T_{\text{ac}i}T_{\text{ac}%
j}-\overset{n}{\underset{k=q+1}{\sum}}(\perp_{\text{a}ik}\perp_{\text{a}%
jk}](y_{0})\Omega_{ij}(y_{0})\phi(y_{0})$

\ $\mathbf{-}$ $\frac{1}{36}\underset{i,j,k=q+1}{\overset{n}{\sum}}%
R_{ijkj}(y_{0})\Omega_{ik}(y_{0})(y_{0})\phi(y_{0})$

$+\frac{1}{6}\underset{i,k=q+1}{\overset{n}{\sum}}\underset{\text{a,b,c}%
=1}{\overset{q}{\sum}}T_{\text{ab}i}(y_{0})[(\perp_{\text{c}ik}T_{\text{ab}%
k})(\frac{\partial\phi}{\partial\text{x}_{\text{c}}}$ $+$ $\Lambda_{\text{c}%
}\phi)](y_{0})\qquad$\textbf{I}$_{32631}\qquad$\textbf{I}$_{3263}$

$-\frac{1}{12}[$ $(R_{\text{a}i\text{b}j}+R_{\text{a}j\text{b}i})$
$-\underset{\text{c=1}}{\overset{\text{q}}{\sum}}(T_{\text{ac}i}T_{\text{bc}%
j}+T_{\text{ac}j}T_{\text{bc}i})$

$-\overset{n}{\underset{k=q+1}{\sum}}(\perp_{\text{a}ik}\perp_{\text{b}jk}+$
$\perp_{\text{a}jk}\perp_{\text{b}ik})](y_{0})T_{\text{ab}i}(y_{0})\Lambda
_{j}(y_{0})\phi(y_{0})-\frac{1}{12}T_{\text{ab}i}^{2}(y_{0})\Omega_{ij}%
(y_{0})\phi(y_{0}$

$+\frac{1}{12}\underset{i,j=q+1}{\overset{n}{\sum}}\underset{\text{a,b=1}%
}{\overset{q}{\sum}}[\perp_{\text{a}ij}(\frac{\partial\phi}{\partial
\text{x}_{\text{b}}}+\Lambda_{\text{b}})](y_{0})\qquad\qquad$\ \textbf{I}%
$_{32632}$

$\times\lbrack-R_{\text{a}i\text{b}j}-R_{\text{a}j\text{b}i}%
+\underset{\text{c=1}}{\overset{\text{q}}{\sum}}T_{\text{ac}i}T_{\text{bc}%
j}-3\underset{\text{c=1}}{\overset{\text{q}}{\sum}}T_{\text{ac}j}%
T_{\text{bc}i}+\underset{\text{k=q+1}}{\overset{\text{n}}{\sum}}%
\perp_{\text{a}i\text{k}}\perp_{\text{b}j\text{k}}-\underset{\text{k=q+1}%
}{\overset{\text{n}}{\sum}}\perp_{\text{a}j\text{k}}\perp_{\text{b}i\text{k}%
}\ ](y_{0})\phi(y_{0})$

$-\frac{1}{6}\underset{i,j=q+1}{\overset{n}{\sum}}\underset{\text{a,b=1}%
}{\overset{q}{\sum}}T_{\text{ab}j}(y_{0})\perp_{\text{a}ij}(y_{0}%
)\frac{\partial\Lambda_{\text{b}}}{\partial x_{i}}(y_{0})\phi(y_{0})$

$-\frac{1}{6}\underset{i,j,k=q+1}{\overset{n}{\sum}}\underset{\text{a}%
=1}{\overset{q}{\sum}}\perp_{\text{a}ij}(y_{0})[\overset{q}{\underset{\text{b}%
=1}{\sum(}}\perp_{\text{b}ik}T_{\text{ab}j})(y_{0})+\frac{2}{3}(2R_{\text{a}%
ijk}+R_{\text{a}jik}+R_{\text{a}kji})](y_{0})\Lambda_{k}(y_{0})\phi(y_{0})$

$+\frac{1}{6}\underset{i,j,k=q+1}{\overset{n}{\sum}}\underset{\text{a}%
=1}{\overset{q}{\sum}}\perp_{\text{a}ij}(y_{0})\perp_{\text{a}jk}(y_{0}%
)\Omega_{ik}(y_{0})\phi(y_{0})$

$+\ \frac{1}{24}\underset{i\text{=q+1}}{\overset{\text{n}}{\sum}}%
\frac{\partial^{2}\text{W}}{\partial x_{i}^{2}}(y_{0})\phi(y_{0})\qquad\qquad
$\ \textbf{I}$_{327}$

$+\frac{1}{24}\underset{i,j\text{=q+1}}{\overset{\text{n}}{\sum}%
}\underset{\text{a=1}}{\overset{\text{q}}{\sum}}<H,j>[4\underset{\text{c}%
=1}{\overset{q}{\sum}}(T_{\text{ac}i})(\perp_{j\text{c}i})+\frac{8}%
{3}R_{i\text{a}ij}](y_{0})\frac{\partial\phi}{\partial\text{x}_{\text{a}}%
}(y_{0})\qquad$I$_{328}$

$+\frac{1}{24}\underset{i,j\text{=q+1}}{\overset{\text{n}}{\sum}%
}\underset{\text{a=1}}{\overset{\text{q}}{\sum}}\perp_{\text{a}ji}%
(y_{0})[<H,i><H,j>](y_{0})\frac{\partial\phi}{\partial\text{x}_{\text{a}}%
}(y_{0})$\qquad$\ \ \ +\frac{1}{72}\underset{i,j=q+1}{\overset{n}{\sum}%
}\underset{\text{a=1}}{\overset{\text{q}}{\sum}}\perp_{\text{a}ji}%
(y_{0})[2\varrho_{ij}+4\overset{q}{\underset{\text{a}=1}{\sum}}R_{i\text{a}%
j\text{a}}-6\overset{q}{\underset{\text{b,c}=1}{\sum}}T_{\text{cc}%
i}T_{\text{bb}j}-T_{\text{bc}i}T_{\text{bc}j}](y_{0})\frac{\partial\phi
}{\partial\text{x}_{\text{a}}}(y_{0})$

$+\frac{8}{3}R_{j\text{a}ji}(y_{0})X_{i}(y_{0})+[2X_{j}\frac{\partial X_{j}%
}{\partial x_{\text{a}}}-\frac{\partial^{2}X_{j}}{\partial x_{\text{a}%
}\partial x_{j}}](y_{0})\frac{\partial\phi}{\partial\text{x}_{\text{a}}}%
(y_{0})$

$+\frac{1}{24}\underset{i,j=q+1}{\overset{n}{\sum}}\underset{\text{a=1}%
}{\overset{\text{q}}{\sum}}<H,j>(y_{0})[4\underset{\text{c}%
=1}{\overset{q}{\sum}}(T_{\text{ac}i})(\perp_{j\text{c}i})+\frac{8}%
{3}R_{i\text{a}ij}](y_{0})\Lambda_{\text{a}}(y_{0})\phi(y_{0})\qquad$%
I$_{329}\qquad$I$_{3291}\qquad$I$_{32911}\qquad$

$+\frac{1}{12}\underset{i,j=q+1}{\overset{n}{\sum}}\underset{\text{a=1}%
}{\overset{\text{q}}{\sum}}\perp_{\text{a}ji}(y_{0})[<H,i><H,j>](y_{0}%
)\Lambda_{\text{a}}(y_{0})\phi(y_{0})$

$+\frac{1}{72}\underset{i,j=q+1}{\overset{n}{\sum}}\underset{\text{a=1}%
}{\overset{\text{q}}{\sum}}\perp_{\text{a}ji}(y_{0})[2\varrho_{ij}%
+4\overset{q}{\underset{\text{a}=1}{\sum}}R_{i\text{a}j\text{a}}%
-6\overset{q}{\underset{\text{b,c}=1}{\sum}}T_{\text{cc}i}T_{\text{bb}%
j}-T_{\text{bc}i}T_{\text{bc}j}](y_{0})\Lambda_{\text{a}}(y_{0})\phi(y_{0})$

$+\frac{1}{36}\underset{i,j=q+1}{\overset{n}{\sum}}%
\underset{k=q+1}{\overset{n}{\sum}}<H,k>(y_{0})$R$_{ijik}(y_{0})\Lambda
_{j}(y_{0})\phi(y_{0})$

$-\frac{1}{288}\underset{i,j=q+1}{\overset{n}{\sum}}<H,j>(y_{0})[3<H,i>^{2}%
+2(\tau^{M}-3\tau^{P}+\overset{q}{\underset{\text{a}=1}{\sum}}\varrho
_{\text{aa}}+\overset{q}{\underset{\text{a,b}=1}{\sum}}R_{\text{abab}}%
)](y_{0})\Lambda_{j}(y_{0})\phi(y_{0})$

$-\frac{1}{12}\underset{i,j=q+1}{\overset{n}{\sum}}<H,i>(y_{0})$

$\times\lbrack\frac{3}{4}<H,i><H,j>$\ $+\frac{1}{6}(\varrho_{ij}%
+2\overset{q}{\underset{\text{a}=1}{\sum}}R_{i\text{a}j\text{a}}%
-3\overset{q}{\underset{\text{a,b=1}}{\sum}}T_{\text{aa}i}T_{\text{bb}%
j}-T_{\text{ab}i}T_{\text{ab}j})](y_{0})\Lambda_{j}(y_{0})\phi(y_{0})$

$+\frac{5}{32}\underset{i,j=q+1}{\overset{n}{\sum}}<H,i>^{2}<H,j>\Lambda
_{j}(y_{0})\phi(y_{0})$

$+\frac{1}{48}\underset{i,j=q+1}{\overset{n}{\sum}}$%
$<$%
H,$i$%
$>$%
(y$_{0}$)$[$(2$\varrho_{ij}$+4$\overset{q}{\underset{\text{a}=1}{\sum}}%
$R$_{i\text{a}j\text{a}}$-3$\overset{q}{\underset{\text{a,b=1}}{\sum}}%
$T$_{\text{aa}i}$T$_{\text{bb}j}$-T$_{\text{ab}i}$T$_{\text{ab}j}%
$-3$\overset{q}{\underset{\text{a,b=1}}{\sum}}$T$_{\text{aa}j}$T$_{\text{bb}%
i}$-T$_{\text{ab}j}$T$_{\text{ab}i}$)$](y_{0})\Lambda_{j}(y_{0})\phi(y_{0})$

$+\frac{1}{48}\underset{i,j=q+1}{\overset{n}{\sum}}<H,j>(y_{0})[$ $\tau
^{M}-3\tau^{P}+\overset{q}{\underset{\text{a}=1}{\sum}}\varrho_{\text{aa}%
}+\overset{q}{\underset{\text{a,b}=1}{\sum}}R_{\text{abab}}](y_{0})\Lambda
_{j}(y_{0})\phi(y_{0})$

$+\frac{1}{144}\underset{i,j=q+1}{\overset{n}{\sum}}[\nabla_{i}\varrho
_{ij}-2\varrho_{ij}<H,i>+\overset{q}{\underset{\text{a}=1}{\sum}}(\nabla
_{i}R_{\text{a}i\text{a}j}-4R_{i\text{a}j\text{a}}<H,i>)$

$+4\overset{q}{\underset{\text{a,b=1}}{\sum}}R_{i\text{a}j\text{b}%
}T_{\text{ab}i}+2\overset{q}{\underset{\text{a,b,c=1}}{\sum}}(T_{\text{aa}%
i}T_{\text{bb}j}T_{\text{cc}i}-3T_{\text{aa}i}T_{\text{bc}j}T_{\text{bc}%
i}+2T_{\text{ab}i}T_{\text{bc}j}T_{\text{ca}i})](y_{0})\Lambda_{j}(y_{0}%
)\phi(y_{0})$\qquad\qquad\qquad\qquad\qquad\ \ 

$+\frac{1}{144}\underset{i,j=q+1}{\overset{n}{\sum}}[\nabla_{j}\varrho
_{ii}-2\varrho_{ji}<H,i>+\overset{q}{\underset{\text{a}=1}{\sum}}(\nabla
_{j}R_{\text{a}i\text{a}i}-4R_{j\text{a}i\text{a}}<H,i>)$

$+4\overset{q}{\underset{\text{a,b=1}}{\sum}}R_{j\text{a}i\text{b}%
}T_{\text{ab}i}$+2$\overset{q}{\underset{\text{a,b,c=1}}{\sum}}(T_{\text{aa}%
j}T_{\text{bb}i}T_{\text{cc}i}-3T_{\text{aa}j}T_{\text{bc}i}T_{\text{bc}%
i}+2T_{\text{ab}j}T_{\text{bc}i}T_{\text{ca}i})](y_{0})\Lambda_{j}(y_{0}%
)\phi(y_{0})$

$+\frac{1}{144}\underset{i,j=q+1}{\overset{n}{\sum}}[\nabla_{i}\varrho
_{ij}-2\varrho_{ii}<H,j>+\overset{q}{\underset{\text{a}=1}{\sum}}(\nabla
_{i}R_{\text{a}i\text{a}j}-4R_{i\text{a}i\text{a}}<H,j>)$

$+4\overset{q}{\underset{\text{a,b=1}}{\sum}}R_{i\text{a}i\text{b}%
}T_{\text{ab}j}$+2$\overset{q}{\underset{\text{a,b,c=1}}{\sum}}(T_{\text{aa}%
i}T_{\text{bb}i}T_{\text{cc}j}-3T_{\text{aa}i}T_{\text{bc}i}T_{\text{bc}%
j}+2T_{\text{ab}i}T_{\text{bc}i}T_{\text{ac}j})](y_{0})\Lambda_{j}(y_{0}%
)\phi(y_{0})$

$+\frac{1}{72}\underset{i,j=q+1}{\overset{n}{\sum}}$ $<H,j>(y_{0}%
)\frac{\partial\Omega_{ij}}{\partial\text{x}_{i}}(y_{0})\phi(y_{0}%
)\qquad\qquad$I$_{32912}$

$-\frac{1}{12}\underset{i=q+1}{\overset{n}{\sum}}\underset{\text{a=1}%
}{\overset{\text{q}}{\sum}}\perp_{\text{a}ij}(y_{0})<H,j>(y_{0})[\Omega
_{i\text{a}}+[\Lambda_{\text{a}},\Lambda_{i}]](y_{0})\phi(y_{0})\qquad\qquad
$I$_{32913}$

$+$ $\frac{1}{72}\underset{i,j=q+1}{\overset{n}{\sum}}[3<H,i><H,j>\ +(\varrho
_{ij}+2\overset{q}{\underset{\text{a}=1}{\sum}}R_{i\text{a}j\text{a}%
}-3\overset{q}{\underset{\text{a,b}=1}{\sum}}T_{\text{aa}i}T_{\text{bb}%
j}-T_{\text{ab}i}T_{\text{ab}j})](y_{0})\Omega_{ij}(y_{0})\phi(y_{0})$

$+\frac{1}{12}\underset{i=q+1}{\overset{n}{\sum}}$ $\underset{\text{a=1}%
}{\overset{\text{q}}{\sum}}$ $[-4\underset{\text{b=1}}{\overset{\text{q}}{%
{\textstyle\sum}
}}T_{\text{ab}i}\frac{\partial X_{i}}{\partial x_{\text{b}}}%
+\underset{j=q+1}{\overset{n}{%
{\textstyle\sum}
}}\perp_{\text{a}ij}\left(  \frac{\partial X_{i}}{\partial x_{j}}%
+\frac{\partial X_{j}}{\partial x_{i}}\right)  \qquad\qquad$I$_{3292}\qquad
$I$_{32921}$

$+\frac{8}{3}\underset{j=q+1}{\overset{n}{%
{\textstyle\sum}
}}R_{i\text{a}ij}X_{j}+\left(  2X_{i}\frac{\partial X_{i}}{\partial
x_{\text{a}}}-\frac{\partial^{2}X_{i}}{\partial x_{\text{a}}\partial x_{i}%
}\right)  ](y_{0})\Lambda_{\text{a}}(y_{0})\phi(y_{0}$

$+\frac{1}{6}$ $\underset{i=q+1}{\overset{n}{\sum}}\underset{\text{a=1}%
}{\overset{\text{q}}{\sum}}$ $[\underset{j=q+1}{\overset{n}{%
{\textstyle\sum}
}}X_{j}\perp_{\text{a}ij}+\frac{\partial X_{i}}{\partial x_{\text{a}}}]$
$(y_{0})[\Omega_{i\text{a}}+[\Lambda_{\text{a}},\Lambda_{i}](y_{0})\phi
(y_{0})$

$-\frac{1}{36}[\left(  2\frac{\partial^{2}X_{i}}{\partial x_{i}\partial x_{j}%
}+\frac{\partial^{2}X_{j}}{\partial x_{i}^{2}}\right)
+2\underset{k=q+1}{\overset{n}{\sum}}R_{ijik}X_{k}](y_{0})\Lambda_{j}%
(y_{0})\phi(y_{0})\qquad$I$_{32922}$

$-\frac{1}{36}X_{j}(y_{0})$ $\frac{\partial\Omega_{ij}}{\partial x_{i}}%
(y_{0})\phi(y_{0})-\frac{1}{12}\frac{\partial X_{j}}{\partial x_{i}}%
(y_{0})\Omega_{ij}(y_{0})$ $\phi(y_{0})$

$+\frac{1}{12}$ $\underset{\text{a=1}}{\overset{\text{q}}{\sum}}\frac
{\partial^{2}X_{\text{a}}}{\partial x_{i}^{2}}(y_{0})\frac{\partial\phi
}{\partial\text{x}_{\text{a}}}(y_{0})\qquad$\ \ \textbf{L}$_{1}$\qquad\qquad

$+\frac{1}{12}$ $\underset{\text{a}=1}{\overset{\text{q}}{\sum}}$
$\frac{\partial^{2}X_{\text{a}}}{\partial x_{i}^{2}}(y_{0})\Lambda_{\text{a}%
}(y_{0})\phi(y_{0})+\frac{1}{12}\underset{j=q+1}{\overset{n}{\sum}}%
\frac{\partial^{2}X_{j}}{\partial x_{i}^{2}}(y_{0})\Lambda_{j}(y_{0}%
)\phi(y_{0})\qquad$\textbf{L}$_{2}\qquad$\textbf{L}$_{21}$

$+\frac{1}{36}$ $\underset{j=q+1}{\overset{n}{\sum}}$ $X_{j}(y_{0}%
)\frac{\partial\Omega_{ij}}{\partial x_{i}}(y_{0})\phi(y_{0})\qquad$%
\textbf{L}$_{22}$

$+\frac{1}{12}$ $\underset{j=q+1}{\overset{\text{n}}{\sum}}$ $\frac{\partial
X_{j}}{\partial x_{i}}(y_{0})\Omega_{ij}(y_{0})\phi(y_{0})\qquad$%
\textbf{L}$_{23}$\qquad$\qquad$

$+\frac{1}{48}$ $\underset{\text{c=1}}{\overset{\text{q}}{\sum}}%
\Lambda_{\text{c}}(y_{0})[\underset{\alpha=q+1}{\overset{n}{\sum}}%
3<H,i>^{2}+2(\tau^{M}-3\tau^{P}\ +\overset{q}{\underset{\text{a=1}}{\sum}%
}\varrho_{\text{aa}}^{M}+\overset{q}{\underset{\text{a,b}=1}{\sum}%
}R_{\text{abab}}^{M})](y_{0})\frac{\partial\phi}{\partial x_{\text{c}}}%
(y_{0})\qquad$I$_{33}\qquad$I$_{331}$

$-\frac{1}{4}\underset{\text{c=1}}{\overset{\text{q}}{\sum}}\Lambda_{\text{c}%
}(y_{0})[\left\Vert \text{X}\right\Vert ^{2}+$ divX $-$ $\underset{\text{a}%
=1}{\overset{q}{\sum}}($X$_{\text{a}})^{2}$ $-$ $\underset{\text{a}%
=1}{\overset{q}{\sum}}\frac{\partial X_{\text{a}}}{\partial x_{\text{a}}%
}](y_{0})\frac{\partial\phi}{\partial x_{\text{c}}}(y_{0})+$ $\frac{1}%
{2}\underset{\text{c=1}}{\overset{\text{q}}{\sum}}\Lambda_{\text{c}}(y_{0}%
)$V$(y_{0})\frac{\partial\phi}{\partial x_{\text{c}}}(y_{0})$

$+\frac{1}{2}\underset{\text{c=1}}{\overset{\text{q}}{\sum}}\Lambda_{\text{c}%
}(y_{0})[-(X_{j}\frac{\partial X_{j}}{\partial x_{\text{c}}}+\frac{1}{2}%
\frac{\partial^{2}X_{j}}{\partial x_{\text{c}}\partial x_{j}})(y_{0})+\frac
{1}{2}(<H,j>\frac{\partial X_{j}}{\partial x_{\text{c}}})(y_{0})+\frac
{\partial\text{V}}{\partial x_{\text{c}}}(y_{0})]\phi(y_{0})$

$+\frac{1}{4}\underset{\text{a,c=1}}{\overset{\text{q}}{\sum}}\Lambda
_{\text{c}}(y_{0})\left\{  \frac{\partial^{3}\phi}{\partial\text{x}_{\text{a}%
}^{2}\partial x_{\text{c}}}\text{ }\right\}  (y_{0})\qquad$I$_{332}$

$+\frac{1}{2}\underset{\text{a,b=1}}{\overset{\text{q}}{\sum}}\left\{
\Lambda_{\text{a}}(y_{0})\Lambda_{\text{b}}(y_{0})\frac{\partial^{2}\phi
}{\partial\text{x}_{\text{a}}\partial x_{\text{b}}}\text{ }\right\}
(y_{0})\qquad$I$_{334}$

$+\frac{1}{4}\underset{\text{c=1}}{\overset{\text{q}}{\sum}}%
[\underset{i=1}{\overset{n}{\sum}}\Lambda_{\text{c}}(y_{0})\Lambda_{i}%
^{2}(y_{0})\frac{\partial\phi}{\partial x_{\text{c}}}(y_{0})\qquad$I$_{335}$

$-\frac{1}{4}\underset{j=q+1}{\overset{n}{\sum}}\underset{\text{a,c=1}%
}{\overset{\text{q}}{\sum}}\Lambda_{\text{c}}(y_{0})T_{\text{aa}j}%
(y_{0})\left\{  \Lambda_{j}\frac{\partial\phi}{\partial x_{\text{c}}}\right\}
(y_{0})\qquad$I$_{336}$

$-$ $\frac{1}{2}\underset{j=q+1}{\overset{n}{\sum}}\underset{\text{c=1}%
}{\overset{\text{q}}{\sum}}\Lambda_{\text{c}}(y_{0})\frac{\partial X_{j}%
}{\partial x_{\text{a}}}(y_{0})\Lambda_{\text{c}}(y_{0})\Lambda_{j}(y_{0}%
)\phi(y_{0})\qquad$I$_{338}$

$+\frac{1}{4}\underset{\text{c=1}}{\overset{\text{q}}{\sum}}\left[
\Lambda_{\text{c}}\frac{\partial\text{W}}{\partial x_{\text{c}}}\phi
+\Lambda_{\text{c}}\text{W}\frac{\partial\phi}{\partial x_{\text{c}}}\right]
(y_{0})\qquad$I$_{339}$

$+$ $\underset{\text{a,c=1}}{\overset{\text{q}}{\sum}}[\Lambda_{\text{c}}%
\frac{\partial\text{X}_{\text{a}}}{\partial x_{\text{a}}}](y_{0}%
)[\frac{\partial\phi}{\partial\text{x}_{\text{a}}}+$ X$_{\text{a}}%
\frac{\partial^{2}\phi}{\partial\text{x}_{\text{a}}^{2}}](y_{0})\qquad$E$_{1}$

$+$ $\underset{\text{b=1}}{\overset{\text{q}}{\sum}}$ $[\frac{\partial
\text{X}_{\text{b}}}{\partial x_{\text{a}}}\Lambda_{\text{b}}\Lambda
_{\text{c}}](y_{0})\phi(y_{0})+$ $\underset{\text{b=1}}{\overset{\text{q}%
}{\sum}}$ $[$X$_{\text{b}}\Lambda_{\text{c}}\Lambda_{\text{b}}\frac
{\partial\phi}{\partial x_{\text{a}}}](y_{0})\qquad$E$_{2}$

$+$ $\underset{j=q+1}{\overset{\text{n}}{\sum}}$ $[\frac{\partial X_{j}%
}{\partial x_{\text{a}}}\Lambda_{\text{c}}\Lambda_{j}](y_{0})\phi(y_{0})+$
$\underset{j=q+1}{\overset{\text{n}}{\sum}}$ $\Lambda_{\text{c}}%
(y_{0})[\Lambda_{j}\frac{\partial\phi}{\partial x_{\text{a}}}](y_{0})$

$+\frac{1}{96}\underset{\text{c=1}}{\overset{\text{q}}{\sum}}\Lambda
_{\text{c}}^{2}(y_{0})[\underset{i=q+1}{\overset{n}{\sum}}3<H,i>^{2}%
+2(\tau^{M}-3\tau^{P}\ +\overset{q}{\underset{\text{a=1}}{\sum}}%
\varrho_{\text{aa}}^{M}+\overset{q}{\underset{\text{a,b}=1}{\sum}%
}R_{\text{abab}}^{M})](y_{0})\phi\left(  y_{0}\right)  \qquad$I$_{34}$

$-\frac{1}{8}\underset{\text{c=1}}{\overset{\text{q}}{\sum}}\Lambda_{\text{c}%
}^{2}(y_{0})[$ $\left\Vert \text{X}\right\Vert _{M}^{2}+\frac{1}{2}$
$\operatorname{div}X_{M}-\frac{1}{2}$ $\left\Vert \text{X}\right\Vert _{P}%
^{2}$ $-$ $\frac{1}{2}\operatorname{div}X_{P}](y_{0})\phi\left(  y_{0}\right)
$

$+$ $\frac{1}{8}\underset{\text{c=1}}{\overset{\text{q}}{\sum}}\Lambda
_{\text{c}}^{2}(y_{0})[\underset{\text{a=1}}{\overset{\text{q}}{\sum}}%
\frac{\partial^{2}\phi}{\partial\text{x}_{\text{a}}^{2}}$ $+$ $2$
$\underset{\text{a=1}}{\overset{\text{q}}{\sum}}\Lambda_{\text{a}}%
\frac{\partial\phi}{\partial x_{\text{a}}}\ +$ $\underset{\text{a=1}%
}{\overset{\text{q}}{\sum}}\Lambda_{\text{a}}^{2}](y_{0})\phi\left(
y_{0}\right)  $

$+\frac{1}{4}\underset{\text{c=1}}{\overset{\text{q}}{\sum}}\Lambda_{\text{c}%
}^{2}(y_{0})[$ $\underset{\text{a=1}}{\overset{\text{q}}{\sum}}X_{\text{a}%
}\frac{\partial\phi}{\partial\text{x}_{\text{a}}}+$ $\underset{\text{a=1}%
}{\overset{\text{q}}{\sum}}$ $X_{\text{a}}\Lambda_{\text{a}}+\frac{1}{2}$W $+$
V$](y_{0})\phi(y_{0})$

$+\frac{1}{96}[3\underset{j=q+1}{\overset{n}{\sum}}<H,j>^{2}+2(\tau^{M}%
-3\tau^{P}\ +\overset{q}{\underset{\text{a=1}}{\sum}}\varrho_{\text{aa}}%
^{M}+\overset{q}{\underset{\text{a,b}=1}{\sum}}R_{\text{abab}}^{M})](y_{0}%
)$W$(y_{0})\phi\left(  y_{0}\right)  \qquad$I$_{35}$

$-\frac{1}{8}[$ $\left\Vert \text{X}\right\Vert _{M}^{2}+$ $\operatorname{div}%
X_{M}-$ $\left\Vert \text{X}\right\Vert _{P}^{2}$ $-\operatorname{div}%
X_{P}](y_{0})$W$(y_{0})\phi\left(  y_{0}\right)  $

$+$ $\frac{1}{8}\underset{\text{a=1}}{\overset{\text{q}}{\sum}}\frac
{\partial^{2}\phi}{\partial\text{x}_{\text{a}}^{2}}(y_{0})$W$(y_{0})$
$+\frac{1}{4}$ $\underset{\text{a=1}}{\overset{\text{q}}{\sum}}\Lambda
_{\text{a}}(y_{0})\frac{\partial\phi}{\partial x_{\text{a}}}\left(
y_{0}\right)  $W$(y_{0})\ +\frac{1}{8}$ $\underset{\text{a=1}%
}{\overset{\text{q}}{\sum}}\Lambda_{\text{a}}^{2}(y_{0})$W$(y_{0})\phi\left(
y_{0}\right)  $

$+\frac{1}{4}$ $\underset{\text{a=1}}{\overset{\text{q}}{\sum}}$X$_{\text{a}%
}(y_{0})\frac{\partial\phi}{\partial\text{x}_{\text{a}}}(y_{0})$%
W$(y_{0})+\frac{1}{4}$ $\underset{\text{a=1}}{\overset{\text{q}}{\sum}}$
X$_{\text{a}}(y_{0})\Lambda_{\text{a}}(y_{0})$W$(y_{0})\phi(y_{0})+\frac{1}%
{8}$W$^{2}\left(  y_{0}\right)  \phi(y_{0})+\frac{1}{4}$ V$(y_{0})$%
W$(y_{0})\phi\left(  y_{0}\right)  $

$+\frac{1}{48}$ $\underset{\text{c=1}}{\overset{\text{q}}{\sum}}$X$_{\text{c}%
}(y_{0})[\underset{\alpha=q+1}{\overset{n}{\sum}}3<H,i>^{2}+2(\tau^{M}%
-3\tau^{P}\ +\overset{q}{\underset{\text{a=1}}{\sum}}\varrho_{\text{aa}}%
^{M}+\overset{q}{\underset{\text{a,b}=1}{\sum}}R_{\text{abab}}^{M}%
)](y_{0})\frac{\partial\phi}{\partial x_{\text{c}}}(y_{0})\qquad$I$_{36}%
\qquad$I$_{361}$

$-\frac{1}{4}\underset{\text{c=1}}{\overset{\text{q}}{\sum}}$X$_{\text{c}%
}(y_{0})[\left\Vert \text{X}\right\Vert ^{2}+$ divX $-$ $\underset{\text{a}%
=1}{\overset{q}{\sum}}($X$_{\text{a}})^{2}$ $-$ $\underset{\text{a}%
=1}{\overset{q}{\sum}}\frac{\partial X_{\text{a}}}{\partial x_{\text{a}}%
}](y_{0})\frac{\partial\phi}{\partial x_{\text{c}}}(y_{0})+$ $\frac{1}%
{2}\underset{\text{c=1}}{\overset{\text{q}}{\sum}}\Lambda_{\text{c}}(y_{0}%
)$V$(y_{0})\frac{\partial\phi}{\partial x_{\text{c}}}(y_{0})$

$+\frac{1}{2}\underset{\text{c=1}}{\overset{\text{q}}{\sum}}$X$_{\text{c}%
}(y_{0})[-(X_{j}\frac{\partial X_{j}}{\partial x_{\text{c}}}+\frac{1}{2}%
\frac{\partial^{2}X_{j}}{\partial x_{\text{c}}\partial x_{j}})(y_{0})+\frac
{1}{2}(<H,j>\frac{\partial X_{j}}{\partial x_{\text{c}}})(y_{0})+\frac
{\partial\text{V}}{\partial x_{\text{c}}}(y_{0})]\phi(y_{0})$

$+\frac{1}{4}\underset{\text{a,c=1}}{\overset{\text{q}}{\sum}}$X$_{\text{c}%
}(y_{0})\frac{\partial^{3}\phi}{\partial\text{x}_{\text{a}}^{2}\partial
x_{\text{c}}}$ $(y_{0})$ \qquad\qquad I$_{362}$

$+\frac{1}{2}\underset{\text{a,c=1}}{\overset{\text{q}}{\sum}}[$X$_{\text{c}%
}\Lambda_{\text{a}}\frac{\partial^{2}\phi}{\partial\text{x}_{\text{a}}\partial
x_{\text{c}}}$ $](y_{0})\qquad\qquad$I$_{364}$

$+\frac{1}{4}\underset{\text{b,c=1}}{\overset{\text{q}}{\sum}}[$X$_{\text{c}%
}\Lambda_{\text{b}}^{2}](y_{0})\frac{\partial\phi}{\partial x_{\text{c}}%
}(y_{0})\qquad$I$_{365}$

$+\frac{1}{4}\underset{\text{a,c=1}}{\overset{\text{q}}{\sum}}\left[
\text{X}_{\text{c}}\frac{\partial\text{W}}{\partial x_{\text{a}}}\phi
+\text{X}_{\text{c}}\text{W}\frac{\partial\phi}{\partial x_{\text{c}}}\right]
(y_{0})\qquad$I$_{369}\qquad$

$+$ $\underset{\text{a,c=1}}{\overset{\text{q}}{\sum}}[$X$_{\text{c}}%
\frac{\partial\text{X}_{\text{a}}}{\partial x_{\text{a}}}[\frac{\partial\phi
}{\partial\text{x}_{\text{a}}}+$ X$_{\text{a}}\frac{\partial^{2}\phi}%
{\partial\text{x}_{\text{a}}^{2}}](y_{0})\qquad\qquad$E$_{1}$

$+\underset{\text{a,b,c=1}}{\overset{\text{q}}{\sum}}$ $[$X$_{\text{c}}%
\frac{\partial\text{X}_{\text{b}}}{\partial x_{\text{a}}}\Lambda_{\text{b}%
}](y_{0})\phi(y_{0})+$ $\underset{\text{a,b,c=1}}{\overset{\text{q}}{\sum}}$
$[$X$_{\text{c}}$X$_{\text{b}}\Lambda_{\text{b}}](y_{0})\frac{\partial\phi
}{\partial x_{\text{a}}}(y_{0})\qquad$E$_{2}$

$+\frac{1}{48}\underset{\text{a=1}}{\overset{\text{q}}{\sum}}X_{\text{a}%
}(y_{0})\Lambda_{\text{a}}(y_{0})[\underset{i=q+1}{\overset{n}{\sum}%
}3<H,i>^{2}+2(\tau^{M}-3\tau^{P}\ +\overset{q}{\underset{\text{a=1}}{\sum}%
}\varrho_{\text{aa}}^{M}+\overset{q}{\underset{\text{a,b}=1}{\sum}%
}R_{\text{abab}}^{M})](y_{0})\phi\left(  y_{0}\right)  \qquad$I$_{37}$

$-\frac{1}{4}\underset{\text{a=1}}{\overset{\text{q}}{\sum}}X_{\text{a}}%
(y_{0})\Lambda_{\text{a}}(y_{0})[$ $\left\Vert \text{X}\right\Vert _{M}%
^{2}+\frac{1}{2}$ $\operatorname{div}X_{M}-\frac{1}{2}$ $\left\Vert
\text{X}\right\Vert _{P}^{2}$ $-$ $\frac{1}{2}\operatorname{div}X_{P}%
](y_{0})\phi\left(  y_{0}\right)  $

$+$ $\frac{1}{4}\underset{\text{a,b=1}}{\overset{\text{q}}{\sum}}X_{\text{a}%
}(y_{0})\Lambda_{\text{a}}(y_{0})[\frac{\partial^{2}\phi}{\partial
\text{x}_{\text{b}}^{2}}$ $+$ $\Lambda_{\text{b}}\frac{\partial\phi}{\partial
x_{\text{b}}}]\left(  y_{0}\right)  \ +\frac{1}{4}$ $\underset{\text{a,b=1}%
}{\overset{\text{q}}{\sum}}[X_{\text{a}}\Lambda_{\text{a}}\Lambda_{\text{b}%
}^{2}](y_{0})\phi\left(  y_{0}\right)  $

$+\frac{1}{4}$ $\underset{\text{a,b=1}}{\overset{\text{q}}{\sum}}%
[$X$_{\text{a}}\Lambda_{\text{a}}$ $X_{\text{b}}](y_{0})\frac{\partial\phi
}{\partial\text{x}_{\text{b}}}(y_{0})$

$+\frac{1}{2}$ $\underset{\text{a,b=1}}{\overset{\text{q}}{\sum}}$
$[X_{\text{a}}X_{\text{b}}\Lambda_{\text{a}}\Lambda_{\text{b}}](y_{0}%
)\phi(y_{0})+\frac{1}{4}$ $\underset{\text{a=1}}{\overset{\text{q}}{\sum}}$
$[X_{\text{a}}\Lambda_{\text{a}}$W$]\left(  y_{0}\right)  \phi\left(
y_{0}\right)  $

$+\frac{1}{2}$ $\underset{\text{a=1}}{\overset{\text{q}}{\sum}}$
$[X_{\text{a}}\Lambda_{\text{a}}$V$](y_{0})\phi(y_{0})$

$+\frac{1}{48}$ $\underset{j=q+1}{\overset{n}{\sum}}[$X$_{j}\Lambda_{j}%
](y_{0})[\underset{i=q+1}{\overset{n}{\sum}}3<H,i>^{2}+2(\tau^{M}-3\tau
^{P}\ +\overset{q}{\underset{\text{a=1}}{\sum}}\varrho_{\text{aa}}%
^{M}+\overset{q}{\underset{\text{a,b}=1}{\sum}}R_{\text{abab}}^{M}%
)](y_{0})\phi\left(  y_{0}\right)  $

$-\frac{1}{4}\underset{j=q+1}{\overset{n}{\sum}}[X_{j}\Lambda_{j}](y_{0})[$
$\left\Vert \text{X}\right\Vert _{M}^{2}+\frac{1}{2}$ $\operatorname{div}%
X_{M}-\frac{1}{2}$ $\left\Vert \text{X}\right\Vert _{P}^{2}$ $-$ $\frac{1}%
{2}\operatorname{div}X_{P}](y_{0})\phi\left(  y_{0}\right)  $

$+$ $\frac{1}{4}\underset{j=q+1}{\overset{n}{\sum}}\underset{\text{a=1}%
}{\overset{\text{q}}{\sum}}[X_{j}\Lambda_{j}](y_{0})\frac{\partial^{2}\phi
}{\partial\text{x}_{\text{a}}^{2}}(y_{0})$ $+\frac{1}{2}%
\underset{j=q+1}{\overset{n}{\sum}}\underset{\text{a=1}}{\overset{\text{q}%
}{\sum}}[X_{j}\Lambda_{j}\Lambda_{\text{a}}](y_{0})\frac{\partial\phi
}{\partial x_{\text{a}}}\left(  y_{0}\right)  $

$+\frac{1}{4}$ $\underset{j=q+1}{\overset{n}{\sum}}\underset{\text{a=1}%
}{\overset{\text{q}}{\sum}}[X_{j}\Lambda_{j}\Lambda_{\text{a}}^{2}](y_{0}%
)\phi\left(  y_{0}\right)  +$ $\frac{1}{4}$ $\underset{j=q+1}{\overset{n}{\sum
}}\underset{\text{a=1}}{\overset{\text{q}}{\sum}}[X_{j}\Lambda_{j}X_{\text{a}%
}](y_{0})\frac{\partial\phi}{\partial\text{x}_{\text{a}}}(y_{0})$

$+$ $\frac{1}{2}$ $\underset{j=q+1}{\overset{n}{\sum}}\underset{\text{a=1}%
}{\overset{\text{q}}{\sum}}[X_{\text{a}}\Lambda_{\text{a}}X_{j}\Lambda
_{j}](y_{0})$ $\phi(y_{0})+\frac{1}{4}$ $\underset{j=q+1}{\overset{n}{\sum}%
}[X_{j}\Lambda_{j}$W$]\left(  y_{0}\right)  \phi\left(  y_{0}\right)  $

$+$ $\frac{1}{2}$ $\underset{j=q+1}{\overset{n}{\sum}}[X_{j}\Lambda_{j}%
$V$](y_{0})\phi\left(  y_{0}\right)  $

\qquad\qquad\qquad\qquad\qquad\qquad\qquad\qquad\qquad\qquad\qquad\qquad
\qquad\qquad\qquad$\blacksquare$

\end{document}